\def\sphinxdocclass{report}
\documentclass[a4paper,10pt,english]{sphinxmanual}
\ifdefined\pdfpxdimen
   \let\sphinxpxdimen\pdfpxdimen\else\newdimen\sphinxpxdimen
\fi \sphinxpxdimen=.75bp\relax
\ifdefined\pdfimageresolution
    \pdfimageresolution= \numexpr \dimexpr1in\relax/\sphinxpxdimen\relax
\fi
\PassOptionsToPackage{bookmarksdepth=5}{hyperref}

\PassOptionsToPackage{warn}{textcomp}
\usepackage[utf8]{inputenc}
\ifdefined\DeclareUnicodeCharacter
  \ifdefined\DeclareUnicodeCharacterAsOptional
    \def\sphinxDUC#1{\DeclareUnicodeCharacter{"#1}}
  \else
    \let\sphinxDUC\DeclareUnicodeCharacter
  \fi
  \sphinxDUC{00A0}{\nobreakspace}
  \sphinxDUC{2500}{\sphinxunichar{2500}}
  \sphinxDUC{2502}{\sphinxunichar{2502}}
  \sphinxDUC{2514}{\sphinxunichar{2514}}
  \sphinxDUC{251C}{\sphinxunichar{251C}}
  \sphinxDUC{2572}{\textbackslash}
\fi
\usepackage{cmap}
\usepackage[T1]{fontenc}
\usepackage{amsmath,amssymb,amstext}
\usepackage{babel}

\usepackage{amsmath,amsfonts,amssymb,amsthm}

\usepackage{fncychap}
\usepackage[,numfigreset=1,mathnumfig]{sphinx}

\fvset{fontsize=auto}
\usepackage{geometry}

\usepackage{hyperref}
\usepackage{hypcap}
\addto\captionsenglish{}

\usepackage{sphinxmessages}
\title{Cardinal Optimizer (COPT) User~Guide}
\date{Jan 30, 2026}
\release{8.0.1}
\author{\hfill Cardinal Operations \\ \scriptsize \hfill D. Ge, Q. Huangfu, Z. Wang, J. Wu and Y. Ye}
\newcommand{\sphinxlogo}{\sphinxincludegraphics{shanshu_logo.pdf}\par}
\renewcommand{\releasename}{Ver}
\makeindex
\begin{document}

\ifdefined\shorthandoff
  \ifnum\catcode`\=\string=\active\shorthandoff{=}\fi
  \ifnum\catcode`\"=\active\shorthandoff{"}\fi
\fi

\pagestyle{empty}
\sphinxmaketitle
\pagestyle{plain}
\sphinxtableofcontents
\pagestyle{normal}
\phantomsection\label{\detokenize{index::doc}}

\sphinxstepscope

\chapter{Introduction to Cardinal Optimizer}
\label{\detokenize{intro:introduction-to-cardinal-optimizer}}\label{\detokenize{intro:chapintro}}\label{\detokenize{intro::doc}}
\sphinxAtStartPar
Cardinal Optimizer is a high\sphinxhyphen{}performance mathematical programming solver
for efficiently solving large\sphinxhyphen{}scale optimization problem. This documentation
provides basic introduction to the Cardinal Optimizer, including:
\begin{itemize}
\item {} 
\sphinxAtStartPar
{\hyperref[\detokenize{install:parcoptinstall}]{\sphinxcrossref{\DUrole{std,std-ref}{How to install Cardinal Optimizer}}}}

\item {} 
\sphinxAtStartPar
{\hyperref[\detokenize{install:parcoptgetlic}]{\sphinxcrossref{\DUrole{std,std-ref}{How to setup license files}}}}

\item {} 
\sphinxAtStartPar
{\hyperref[\detokenize{cmdline:chapcmdline}]{\sphinxcrossref{\DUrole{std,std-ref}{How to use Cardinal Optimizer in interactive shell}}}}

\end{itemize}

\sphinxAtStartPar
We suggest that all users to read the first two sections carefully before using
the Cardinal Optimizer.

\sphinxAtStartPar
Once the installation and license setup are done, we recommend users who
want to do a quick experiment on the Cardinal Optimizer to read the
{\hyperref[\detokenize{cmdline:chapcmdline}]{\sphinxcrossref{\DUrole{std,std-ref}{COPT Interactive Shell}}}} chapter for details.
Users who already have preferred programming language can select from available
\sphinxcode{\sphinxupquote{Application Programming Interfaces (APIs)}}, including:
\begin{itemize}
\item {} 
\sphinxAtStartPar
{\hyperref[\detokenize{cinterface:chapcinterface}]{\sphinxcrossref{\DUrole{std,std-ref}{C interface}}}}

\item {} 
\sphinxAtStartPar
{\hyperref[\detokenize{cppinterface:chapcppinterface}]{\sphinxcrossref{\DUrole{std,std-ref}{C++ interface}}}}

\item {} 
\sphinxAtStartPar
{\hyperref[\detokenize{csharpinterface:chapcsharpinterface}]{\sphinxcrossref{\DUrole{std,std-ref}{C\# interface}}}}

\item {} 
\sphinxAtStartPar
{\hyperref[\detokenize{javainterface:chapjavainterface}]{\sphinxcrossref{\DUrole{std,std-ref}{Java interface}}}}

\item {} 
\sphinxAtStartPar
{\hyperref[\detokenize{pythoninterface:chappythoninterface}]{\sphinxcrossref{\DUrole{std,std-ref}{Python interface}}}}

\item {} 
\sphinxAtStartPar
{\hyperref[\detokenize{amplinterface:chapamplinterface}]{\sphinxcrossref{\DUrole{std,std-ref}{AMPL interface}}}}

\item {} 
\sphinxAtStartPar
{\hyperref[\detokenize{pyomointerface:chappyomointerface}]{\sphinxcrossref{\DUrole{std,std-ref}{Pyomo interface}}}}

\item {} 
\sphinxAtStartPar
{\hyperref[\detokenize{pulpinterface:chappulpinterface}]{\sphinxcrossref{\DUrole{std,std-ref}{PuLP interface}}}}

\item {} 
\sphinxAtStartPar
{\hyperref[\detokenize{cvxpyinterface:chapcvxpyinterface}]{\sphinxcrossref{\DUrole{std,std-ref}{CVXPY interface}}}}

\end{itemize}

\section{Overview}
\label{\detokenize{intro:overview}}
\sphinxAtStartPar
Cardinal Optimizer supports solving Linear Programming (LP) problems, Second\sphinxhyphen{}Order\sphinxhyphen{}Cone
Programming (SOCP) problems, Quadratic Programming (QP) problems,
Quadratically Constrained Programming (QCP) problems, Exponential Cone Programming (ExpCone) problems,
Semidefinite Programming problems (SDP), General Nonlinear Programming (NLP) problems,
and Mixed Integer Programming (MIP) problems,
which include Mixed Integer Linear Programming (MILP), Mixed Integer Second\sphinxhyphen{}Order\sphinxhyphen{}Cone Programming (MISOCP),
Mixed Integer Convex Quadratic Programming (MIQP), Mixed Integer Convex Quadratically Constrained Programming (MIQCP).

\sphinxAtStartPar
We will support more problem types in the future.
The supported problem types and available algorithms are summarized in \hyperref[\detokenize{intro:copttab-probsalgs}]{Table \ref{\detokenize{intro:copttab-probsalgs}}}

\begin{savenotes}\sphinxattablestart
\sphinxthistablewithglobalstyle
\centering
\sphinxcapstartof{table}
\sphinxthecaptionisattop
\sphinxcaption{Supported problem types and available algorithms}\label{\detokenize{intro:copttab-probsalgs}}
\sphinxaftertopcaption
\begin{tabulary}{\linewidth}[t]{|T|T|}
\sphinxtoprule
\sphinxtableatstartofbodyhook
\sphinxAtStartPar
\sphinxstylestrong{Problem type}
&
\sphinxAtStartPar
Available algorithms
\\
\sphinxhline
\sphinxAtStartPar
Linear Programming (LP)
&
\sphinxAtStartPar
Simplex, Barrier (CPU/GPU),
First\sphinxhyphen{}order Method (PDLP)(CPU/GPU)
\\
\sphinxhline
\sphinxAtStartPar
Second\sphinxhyphen{}Order\sphinxhyphen{}Cone Programming (SOCP)
&
\sphinxAtStartPar
Barrier (CPU/GPU)
\\
\sphinxhline
\sphinxAtStartPar
Exponential Cone Programming (ExpCone)
&
\sphinxAtStartPar
Barrier (CPU/GPU)
\\
\sphinxhline
\sphinxAtStartPar
Convex Quadratic Programming (QP)
&
\sphinxAtStartPar
Barrier (CPU/GPU)
\\
\sphinxhline
\sphinxAtStartPar
Convex Quadratically Constrained Programming (QCP)
&
\sphinxAtStartPar
Barrier (CPU/GPU)
\\
\sphinxhline
\sphinxAtStartPar
Semidefinite Programming (SDP)
&
\sphinxAtStartPar
Barrier (CPU/GPU), ADMM
\\
\sphinxhline
\sphinxAtStartPar
Nonconvex Quadratic Programming (QP)
&
\sphinxAtStartPar
Spatial Branch\sphinxhyphen{}and\sphinxhyphen{}Bound (global optimum),
Barrier (local optimum)
\\
\sphinxhline
\sphinxAtStartPar
Nonconvex Quadratically Constrained Programming (QCP)
&
\sphinxAtStartPar
Spatial Branch\sphinxhyphen{}and\sphinxhyphen{}Bound (global optimum),
Barrier (local optimum)
\\
\sphinxhline
\sphinxAtStartPar
General Nonlinear Programming (NLP)
&
\sphinxAtStartPar
Barrier
\\
\sphinxhline
\sphinxAtStartPar
Mixed Integer Linear Programming (MILP)
&
\sphinxAtStartPar
Branch\sphinxhyphen{}and\sphinxhyphen{}Cut
\\
\sphinxhline
\sphinxAtStartPar
Mixed Integer Second\sphinxhyphen{}Order\sphinxhyphen{}Cone Programming (MISOCP)
&
\sphinxAtStartPar
Branch\sphinxhyphen{}and\sphinxhyphen{}Cut
\\
\sphinxhline
\sphinxAtStartPar
Mixed Integer Convex/Nonconvex Quadratic Programming (MIQP)
&
\sphinxAtStartPar
Branch\sphinxhyphen{}and\sphinxhyphen{}Cut
\\
\sphinxhline
\sphinxAtStartPar
Mixed Integer Convex/Nonconvex
Quadratically Constrained Programming (MIQCP)
&
\sphinxAtStartPar
Branch\sphinxhyphen{}and\sphinxhyphen{}Cut
\\
\sphinxbottomrule
\end{tabulary}
\sphinxtableafterendhook\par
\sphinxattableend\end{savenotes}

\sphinxAtStartPar
Cardinal Optimizer supports all major 64\sphinxhyphen{}bit operating systems including Windows,
Linux (including ARM64 platform) and MacOS (including ARM64 platform),
and currently provides programming interfaces shown below:
\begin{itemize}
\item {} 
\sphinxAtStartPar
C interface

\item {} 
\sphinxAtStartPar
C++ interface

\item {} 
\sphinxAtStartPar
C\# interface

\item {} 
\sphinxAtStartPar
Java interface

\item {} 
\sphinxAtStartPar
Python interface

\item {} 
\sphinxAtStartPar
AMPL interface

\item {} 
\sphinxAtStartPar
AIMMS interface

\item {} 
\sphinxAtStartPar
Pyomo interface

\item {} 
\sphinxAtStartPar
PuLP interface

\item {} 
\sphinxAtStartPar
CVXPY interface

\item {} 
\sphinxAtStartPar
GAMS interface

\item {} 
\sphinxAtStartPar
Julia interface

\end{itemize}

\sphinxAtStartPar
We are going to develop more programming interfaces to suit various needs of
users and situations.

\section{Licenses}
\label{\detokenize{intro:licenses}}
\sphinxAtStartPar
Now, we provides 4 types of license, which are Personal License, Server License,
Floating License, and Cluster License. They are listed in table below (\hyperref[\detokenize{intro:copttab-licensetype}]{Table \ref{\detokenize{intro:copttab-licensetype}}}) :

\begin{savenotes}\sphinxattablestart
\sphinxthistablewithglobalstyle
\centering
\sphinxcapstartof{table}
\sphinxthecaptionisattop
\sphinxcaption{License Type}\label{\detokenize{intro:copttab-licensetype}}
\sphinxaftertopcaption
\begin{tabular}[t]{|\X{8}{40}|\X{32}{40}|}
\sphinxtoprule
\sphinxstyletheadfamily 
\sphinxAtStartPar
License Type
&\sphinxstyletheadfamily 
\sphinxAtStartPar
Detail
\\
\sphinxmidrule
\sphinxtableatstartofbodyhook
\sphinxAtStartPar
Personal License
&
\sphinxAtStartPar
It is tied to personal computers by username. Only approved user can run COPT on his devices. No limitations on CPU cores and threads.
\\
\sphinxhline
\sphinxAtStartPar
Server License
&
\sphinxAtStartPar
It is tied to a single server computer by its hardware info (MAC and CPUID). An arbitrary number of users and programs can run COPT simultaneously. No limitations on CPU cores as well.
\\
\sphinxhline
\sphinxAtStartPar
Floating License
&
\sphinxAtStartPar
It is tied to a server machine running COPT floating token service, by its hardware info (MAC and CPUID). Any COPT floating client connected to server can borrow and use the floating license, thus run one process for optimization jobs simultaneously. The token number is max number of clients who can use floating licenses simultaneously.
\\
\sphinxhline
\sphinxAtStartPar
Cluster License
&
\sphinxAtStartPar
It is tied to a server machine running COPT compute cluster service, by its hardware info (MAC and CPUID). Any COPT compute cluster client connected to server can offload optimization computations. That is, clients are allowed to do modelling locally, execute optimization jobs remotely, and then obtain results interactively. Although server can have multiple clients connected, each connection must run optimization jobs sequentially. No limitations on CPU cores.
\\
\sphinxbottomrule
\end{tabular}
\sphinxtableafterendhook\par
\sphinxattableend\end{savenotes}

\section{How to Cite}
\label{\detokenize{intro:how-to-cite}}
\sphinxAtStartPar
If you used COPT in your research work, please mention us in your publication. For example:
\begin{itemize}
\item {} 
\sphinxAtStartPar
We used COPT {[}1{]} in our project.

\item {} 
\sphinxAtStartPar
To solve the integer problem, we used Cardinal Optimizer {[}1{]}.

\end{itemize}

\sphinxAtStartPar
with the following entry in the Reference section:

\begin{sphinxVerbatim}[commandchars=\\\{\}]
\PYG{p}{[}\PYG{l+m+mi}{1}\PYG{p}{]} \PYG{n}{D}\PYG{o}{.} \PYG{n}{Ge}\PYG{p}{,} \PYG{n}{Q}\PYG{o}{.} \PYG{n}{Huangfu}\PYG{p}{,} \PYG{n}{Z}\PYG{o}{.} \PYG{n}{Wang}\PYG{p}{,} \PYG{n}{J}\PYG{o}{.} \PYG{n}{Wu} \PYG{o+ow}{and} \PYG{n}{Y}\PYG{o}{.} \PYG{n}{Ye}\PYG{o}{.} \PYG{n}{Cardinal} \PYG{n}{Optimizer} \PYG{p}{(}\PYG{n}{COPT}\PYG{p}{)} \PYG{n}{user} \PYG{n}{guide}\PYG{o}{.} \PYG{n}{https}\PYG{p}{:}\PYG{o}{/}\PYG{o}{/}\PYG{n}{guide}\PYG{o}{.}\PYG{n}{coap}\PYG{o}{.}\PYG{n}{online}\PYG{o}{/}\PYG{n}{copt}\PYG{o}{/}\PYG{n}{en}\PYG{o}{\PYGZhy{}}\PYG{n}{doc}\PYG{p}{,} \PYG{l+m+mf}{2023.}
\end{sphinxVerbatim}

\sphinxAtStartPar
The corresponding BiBTeX citation is:

\begin{sphinxVerbatim}[commandchars=\\\{\}]
@misc\PYGZob{}copt,
  author=\PYGZob{}Dongdong Ge and Qi Huangfu and Zizhuo Wang and Jian Wu and Yinyu Ye\PYGZcb{},
  title=\PYGZob{}Cardinal \PYGZob{}O\PYGZcb{}ptimizer \PYGZob{}(COPT)\PYGZcb{} user guide\PYGZcb{},
  howpublished=\PYGZob{}https://guide.coap.online/copt/en\PYGZhy{}doc\PYGZcb{},
  year=2023
\PYGZcb{}
\end{sphinxVerbatim}

\section{Contact Information}
\label{\detokenize{intro:contact-information}}
\sphinxAtStartPar
Cardinal Optimizer is developed by \sphinxhref{https://www.shanshu.ai}{Cardinal Operations},
users who want any further help can contact us using information provided in
\hyperref[\detokenize{intro:copttab-contactinfo}]{Table \ref{\detokenize{intro:copttab-contactinfo}}}

\begin{savenotes}\sphinxattablestart
\sphinxthistablewithglobalstyle
\centering
\sphinxcapstartof{table}
\sphinxthecaptionisattop
\sphinxcaption{Contact information}\label{\detokenize{intro:copttab-contactinfo}}
\sphinxaftertopcaption
\begin{tabulary}{\linewidth}[t]{|T|T|T|}
\sphinxtoprule
\sphinxtableatstartofbodyhook
\sphinxAtStartPar
\sphinxstylestrong{Type}
&
\sphinxAtStartPar
\sphinxstylestrong{Information}
&
\sphinxAtStartPar
\sphinxstylestrong{Description}
\\
\sphinxhline
\sphinxAtStartPar
Website
&
\sphinxAtStartPar
\sphinxurl{https://www.shanshu.ai/}
&\\
\sphinxhline
\sphinxAtStartPar
Phone
&
\sphinxAtStartPar
400\sphinxhyphen{}680\sphinxhyphen{}5680
&\\
\sphinxhline
\sphinxAtStartPar
Email
&
\sphinxAtStartPar
\sphinxhref{mailto:coptsales@shanshu.ai}{coptsales@shanshu.ai}
&
\sphinxAtStartPar
business support
\\
\sphinxhline
\sphinxAtStartPar
Email
&
\sphinxAtStartPar
\sphinxhref{mailto:coptsupport@shanshu.ai}{coptsupport@shanshu.ai}
&
\sphinxAtStartPar
technical support
\\
\sphinxbottomrule
\end{tabulary}
\sphinxtableafterendhook\par
\sphinxattableend\end{savenotes}

\sphinxstepscope

\chapter{Installation Guide}
\label{\detokenize{install:installation-guide}}\label{\detokenize{install:chapinstall}}\label{\detokenize{install::doc}}
\sphinxAtStartPar
This chapter introduces how to install \sphinxstylestrong{Cardinal Optimizer} on all supported
operating systems, and how to obtain and setup license correctly. We recommend
all users ready this chapter carefully before using Cardinal Optimizer.

\section{Registration}
\label{\detokenize{install:registration}}
\sphinxAtStartPar
Before using Cardinal Optimizer, users need to registrate online and then install the
COPT package on your machine. If this is not done yet, please visit
official \sphinxhref{https://copt.shanshu.ai}{COPT home page} and
fill the registration form following the guidelines.

\sphinxAtStartPar
The online registration is for personal license application. Specifically, users
only need to provide username of machine, besides basic information.

\sphinxAtStartPar
Upon approval, you will receive a letter from \sphinxhref{mailto:coptsales@shanshu.ai}{coptsales@shanshu.ai}. It gives
both link to download COPT software package, a license key tied with registration
information, and also two attached license files. You may refer to {\hyperref[\detokenize{install:parcoptinstall}]{\sphinxcrossref{\DUrole{std,std-ref}{Software Installation}}}} below and {\hyperref[\detokenize{install:parcoptgetlic}]{\sphinxcrossref{\DUrole{std,std-ref}{Setting Up License}}}} for further steps.

\sphinxAtStartPar
If you encountered any problems, please contact \sphinxhref{mailto:coptsupport@shanshu.ai}{coptsupport@shanshu.ai} for help.

\section{Software Installation}
\label{\detokenize{install:software-installation}}\label{\detokenize{install:parcoptinstall}}

\subsection{Windows}
\label{\detokenize{install:windows}}
\sphinxAtStartPar
We provide two types of installation packages for Windows operating systems,
The user need to \sphinxstylestrong{choose one of the two options}.
One is an executable installer (Download link name contains \sphinxcode{\sphinxupquote{installer}}) for most of users and the other one is a zip\sphinxhyphen{}format archive specialized for expert users.

\sphinxAtStartPar
The executable installer provides a visual installation prompt window and automatically configures environment variables. Users only need to follow the instructions and click to complete the installation steps in sequence;
The ZIP format installation package requires the user to decompress the installation package first and manually configure the environment variables.
\begin{itemize}
\item {} 
\sphinxAtStartPar
Executable installer: CardinalOptimizer\sphinxhyphen{}8.0.1\sphinxhyphen{}win64\sphinxhyphen{}installer.zip

\item {} 
\sphinxAtStartPar
ZIP\sphinxhyphen{}format installation package: CardinalOptimizer\sphinxhyphen{}8.0.1\sphinxhyphen{}win64.zip

\end{itemize}

\sphinxAtStartPar
\sphinxstylestrong{We recommend users to download the excecutable installer.}

\subsubsection{Executable installer}
\label{\detokenize{install:executable-installer}}
\sphinxAtStartPar
If you download the executable installer for Windows from our website, e.g.
CardinalOptimizer\sphinxhyphen{}8.0.1\sphinxhyphen{}win64\sphinxhyphen{}installer.exe for 64\sphinxhyphen{}bit version of COPT 8.0.1,
just double\sphinxhyphen{}click it and follow the following guidance:
\begin{itemize}
\item {} 
\sphinxAtStartPar
Step 1: Click the installer and select the installation language. The default
installation language is \sphinxcode{\sphinxupquote{English}}, users can change it by select from the
drop\sphinxhyphen{}down menu, see \hyperref[\detokenize{install:coptfig-winexe-langsel}]{Fig.\@ \ref{\detokenize{install:coptfig-winexe-langsel}}}. Here we use the default
setting.

\begin{figure}[H]
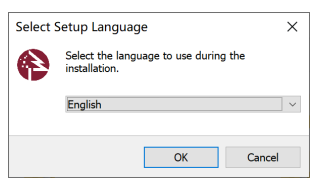

\centering
\capstart

\noindent\sphinxincludegraphics[scale=0.7]{{copt-winexe_langsel}.png}
\caption{Select installation language}\label{\detokenize{install:coptfig-winexe-langsel}}\end{figure}

\item {} 
\sphinxAtStartPar
Step 2: If you agree with the \sphinxcode{\sphinxupquote{\textquotesingle{}End\sphinxhyphen{}User License Agreement (EULA)\textquotesingle{}}},
just choose \sphinxcode{\sphinxupquote{\textquotesingle{}I accept the agreement\textquotesingle{}}} and then click \sphinxcode{\sphinxupquote{\textquotesingle{}Next\textquotesingle{}}}.
The software won’t install if you disagree with the \sphinxcode{\sphinxupquote{EULA}},
see \hyperref[\detokenize{install:coptfig-winexe-eula}]{Fig.\@ \ref{\detokenize{install:coptfig-winexe-eula}}}.

\begin{figure}[H]
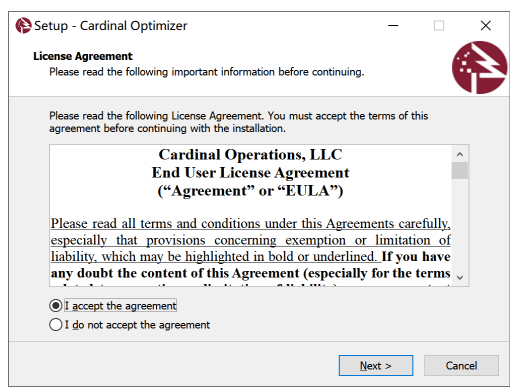

\centering
\capstart

\noindent\sphinxincludegraphics[scale=0.5]{{copt-winexe_eula}.png}
\caption{License agreement page}\label{\detokenize{install:coptfig-winexe-eula}}\end{figure}

\item {} 
\sphinxAtStartPar
Step 3: By default, the installer will place all files into directory
C:\textbackslash{}Program Files\textbackslash{}copt80, you may change it to any directories.
If you have decided the install directory, just click \sphinxcode{\sphinxupquote{\textquotesingle{}Next\textquotesingle{}}},
see \hyperref[\detokenize{install:coptfig-winexe-instdir}]{Fig.\@ \ref{\detokenize{install:coptfig-winexe-instdir}}}.

\begin{figure}[H]
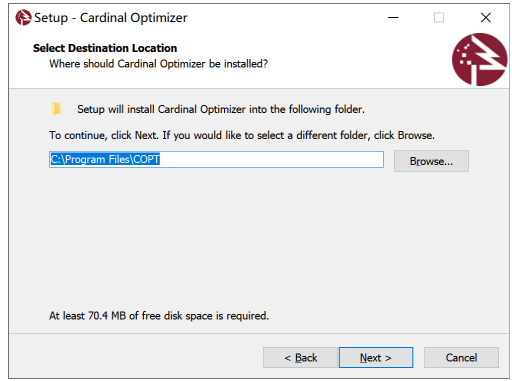

\centering
\capstart

\noindent\sphinxincludegraphics[scale=0.5]{{copt-winexe_instdir}.png}
\caption{Choose install directory}\label{\detokenize{install:coptfig-winexe-instdir}}\end{figure}

\item {} 
\sphinxAtStartPar
Step 4: To select the start menu folder, you can simply use the default
setting and click \sphinxcode{\sphinxupquote{\textquotesingle{}Next\textquotesingle{}}}, see \hyperref[\detokenize{install:coptfig-winexe-startmenu}]{Fig.\@ \ref{\detokenize{install:coptfig-winexe-startmenu}}}.

\begin{figure}[H]
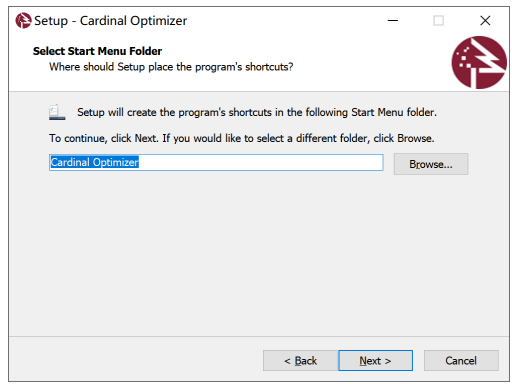

\centering
\capstart

\noindent\sphinxincludegraphics[scale=0.5]{{copt-winexe_startmenu}.png}
\caption{Select start menu folder}\label{\detokenize{install:coptfig-winexe-startmenu}}\end{figure}

\item {} 
\sphinxAtStartPar
Step 5: By now, the software is ready to install, just click \sphinxcode{\sphinxupquote{\textquotesingle{}Install\textquotesingle{}}},
see \hyperref[\detokenize{install:coptfig-winexe-install}]{Fig.\@ \ref{\detokenize{install:coptfig-winexe-install}}}.

\begin{figure}[H]
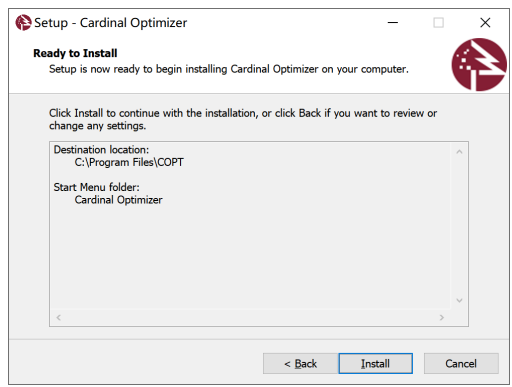

\centering
\capstart

\noindent\sphinxincludegraphics[scale=0.5]{{copt-winexe_install}.png}
\caption{Ready to install}\label{\detokenize{install:coptfig-winexe-install}}\end{figure}

\item {} 
\sphinxAtStartPar
Step 6: When installation completed, the software requires restart of your
machine since the installer has automatically made the required modifications
to environment variables for you. Be sure to save your working files and
close other running applications before restart, and then click \sphinxcode{\sphinxupquote{\textquotesingle{}Finish\textquotesingle{}}},
see \hyperref[\detokenize{install:coptfig-winexe-restart}]{Fig.\@ \ref{\detokenize{install:coptfig-winexe-restart}}}.

\begin{figure}[H]
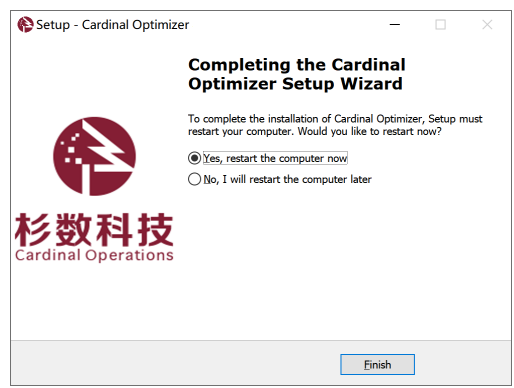

\centering
\capstart

\noindent\sphinxincludegraphics[scale=0.5]{{shanshu_winexe_restart}.png}
\caption{Installation completed and restart your machine}\label{\detokenize{install:coptfig-winexe-restart}}\end{figure}

\end{itemize}

\subsubsection{Zip\sphinxhyphen{}format installer archive}
\label{\detokenize{install:zip-format-installer-archive}}
\sphinxAtStartPar
If you downloaded the zip\sphinxhyphen{}format installer archive, just uncompress it to any
directories with any unarchiver and set up environment variables as follows.
Here we assumed the installation directory to be C:\textbackslash{}Program Files\textbackslash{}copt80:
\begin{itemize}
\item {} 
\sphinxAtStartPar
Step 1: Open \sphinxcode{\sphinxupquote{Command Prompt (cmd)}} with \sphinxstylestrong{administrator priviledge} and
execute the following command to pop\sphinxhyphen{}up the environment variables setting panel.

\begin{sphinxVerbatim}[commandchars=\\\{\}]
rundll32 sysdm.cpl,EditEnvironmentVariables
\end{sphinxVerbatim}

\item {} 
\sphinxAtStartPar
Step 2: Add directory C:\textbackslash{}Program Files\textbackslash{}copt80\textbackslash{}bin to system environment
variable \sphinxcode{\sphinxupquote{PATH}}.

\item {} 
\sphinxAtStartPar
Step 3: Create a new system environment variable named \sphinxcode{\sphinxupquote{COPT\_HOME}}, whose
value is C:\textbackslash{}Program Files\textbackslash{}copt80.

\item {} 
\sphinxAtStartPar
Step 4: Create a new system environment variable named \sphinxcode{\sphinxupquote{COPT\_LICENSE\_DIR}},
whose value is C:\textbackslash{}Program Files\textbackslash{}copt80.

\end{itemize}

\sphinxAtStartPar
Up to now, you have already setup the required modifications to the environment
variables. If you accept the \sphinxcode{\sphinxupquote{copt\sphinxhyphen{}eula\_en.pdf}} in installation directory,
then please go to {\hyperref[\detokenize{install:parcoptgetlic}]{\sphinxcrossref{\DUrole{std,std-ref}{Setting Up License}}}} for license issues.

\subsection{MacOS}
\label{\detokenize{install:macos}}
\sphinxAtStartPar
We provide two types of installation packages for MacOS, one is a DMG\sphinxhyphen{}format
installer for most of users and the other is a gzip\sphinxhyphen{}format archive for expert
users. The user need to \sphinxstylestrong{choose one of the two options}.

\sphinxAtStartPar
In addition, for MacOS systems, we provide the MacOS\sphinxhyphen{}Universal installation package,
which is universal to both Apple Silicon and old Intel chips (the installation package suffix is \sphinxcode{\sphinxupquote{universal\_mac}}).
\begin{itemize}
\item {} 
\sphinxAtStartPar
DMG\sphinxhyphen{}format installer: CardinalOptimizer\sphinxhyphen{}8.0.1\sphinxhyphen{}universal\_mac.dmg

\item {} 
\sphinxAtStartPar
Gzip\sphinxhyphen{}format Archive: CardinalOptimizer\sphinxhyphen{}8.0.1\sphinxhyphen{}universal\_mac.tar.gz

\end{itemize}

\sphinxAtStartPar
\sphinxstylestrong{We recommend users to download the DMG\sphinxhyphen{}format installer.}

\subsubsection{DMG installer}
\label{\detokenize{install:dmg-installer}}
\sphinxAtStartPar
If you download the DMG\sphinxhyphen{}format installer for MacOS\sphinxhyphen{}Universal from our website, e.g.
CardinalOptimizer\sphinxhyphen{}8.0.1\sphinxhyphen{}universal\_mac.dmg of COPT 8.0.1,
please follow the following guidance to install the software:
\begin{itemize}
\item {} 
\sphinxAtStartPar
Step 1: Double\sphinxhyphen{}click the DMG\sphinxhyphen{}format installer, waiting for the OS to mount
the DMG installer automatically.

\item {} 
\sphinxAtStartPar
Step 2: Simply drag the copt80 folder into \sphinxcode{\sphinxupquote{\textquotesingle{}Applications\textquotesingle{}}} folder,
see \hyperref[\detokenize{install:coptfig-macdmg-dragfolder}]{Fig.\@ \ref{\detokenize{install:coptfig-macdmg-dragfolder}}}:

\begin{figure}[H]
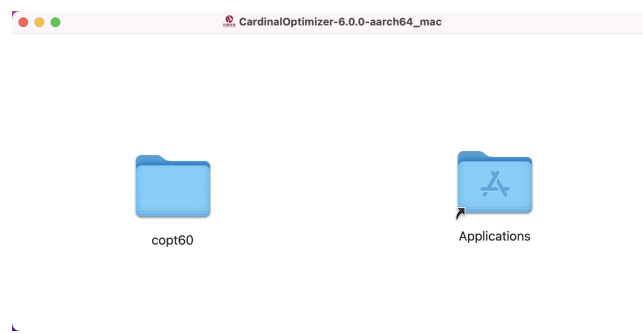

\centering
\capstart

\noindent\sphinxincludegraphics[scale=0.5]{{shanshu_macdmg_dragfolder}.png}
\caption{Drag copt80 into \sphinxcode{\sphinxupquote{\textquotesingle{}Applications\textquotesingle{}}}}\label{\detokenize{install:coptfig-macdmg-dragfolder}}\end{figure}

\end{itemize}

\sphinxAtStartPar
After the above steps are completed, the necessary environment variables need to be configured. Please continue to read: {\hyperref[\detokenize{install:chapinstall-macosenv}]{\sphinxcrossref{\DUrole{std,std-ref}{Environment variable configuration}}}} section.

\subsubsection{GZIP installer}
\label{\detokenize{install:gzip-installer}}
\sphinxAtStartPar
If you download the gzip\sphinxhyphen{}format installer archive, just uncompress it to any
directories with commands:
\begin{sphinxalltt}
tar \sphinxhyphen{}xzf CardinalOptimizer\sphinxhyphen{}8.0.1\sphinxhyphen{}universal\_mac.tar.gz
\end{sphinxalltt}

\sphinxAtStartPar
You will get a folder named copt80 in current directory, you can move it
to any other directories as you want. We recommend users move it to the
\sphinxcode{\sphinxupquote{\textquotesingle{}Applications\textquotesingle{}}} folder by executing commands below:
\begin{sphinxalltt}
mv copt80 /Applications/
\end{sphinxalltt}

\sphinxAtStartPar
After the above steps are completed, the necessary environment variables need to be configured. Please continue to read: {\hyperref[\detokenize{install:chapinstall-macosenv}]{\sphinxcrossref{\DUrole{std,std-ref}{Environment variable configuration}}}} section.

\subsubsection{Environment variable configuration}
\label{\detokenize{install:environment-variable-configuration}}\label{\detokenize{install:chapinstall-macosenv}}
\sphinxAtStartPar
The next step for both DMG\sphinxhyphen{}format and gzip\sphinxhyphen{}format installer is to setup the
environment variables.

\sphinxAtStartPar
First, please execute the following command on the terminal to determine the current Shell version:

\begin{sphinxVerbatim}[commandchars=\\\{\}]
echo \PYGZdl{}SHELL
\end{sphinxVerbatim}

\sphinxAtStartPar
If the output is: \sphinxcode{\sphinxupquote{\textquotesingle{}/bin/bash\textquotesingle{}}}, it means that the current terminal is \sphinxcode{\sphinxupquote{bash}};
if the output is: \sphinxcode{\sphinxupquote{\textquotesingle{}/bin/zsh\textquotesingle{}}}, It means that the current terminal is \sphinxcode{\sphinxupquote{zsh}}.

\sphinxAtStartPar
Next, please open the terminal and confirm that you are in the user directory (if not, you need to execute the following command to switch to the user directory):

\begin{sphinxVerbatim}[commandchars=\\\{\}]
\PYG{n}{cd} \PYG{o}{\PYGZti{}}\PYG{o}{/}
\end{sphinxVerbatim}

\sphinxAtStartPar
Please refer to the corresponding environment variable configuration guide according to the terminal.
\begin{itemize}
\item {} 
\sphinxAtStartPar
{\hyperref[\detokenize{install:parinstallenv-bash}]{\sphinxcrossref{\DUrole{std,std-ref}{BASH terminal}}}}

\item {} 
\sphinxAtStartPar
{\hyperref[\detokenize{install:parinstallenv-zsh}]{\sphinxcrossref{\DUrole{std,std-ref}{ZSH terminal}}}}

\end{itemize}

\sphinxAtStartPar
After the environment variable is configured, you can verify whether the configuration is successful, refer to {\hyperref[\detokenize{install:parinstallenv-valid}]{\sphinxcrossref{\DUrole{std,std-ref}{Verify the environment variable}}}}.

\phantomsection\label{\detokenize{install:parinstallenv-bash}}
\sphinxAtStartPar
\sphinxstylestrong{BASH terminal}

\sphinxAtStartPar
There are three main steps to configure environment variables under the BASH terminal:

\sphinxAtStartPar
\sphinxstylestrong{1. Check ‘bash\_profile’ file}

\sphinxAtStartPar
Enter the following command to output all files in the user directory
and check whether there is a \sphinxcode{\sphinxupquote{\textquotesingle{}.bash\_profile\textquotesingle{}}} hidden file:

\begin{sphinxVerbatim}[commandchars=\\\{\}]
\PYG{n}{ls} \PYG{o}{\PYGZhy{}}\PYG{n}{a}
\end{sphinxVerbatim}

\sphinxAtStartPar
If the file does not exist, execute the following command to create an empty \sphinxcode{\sphinxupquote{\textquotesingle{}.bash\_profile\textquotesingle{}}} file at your own (if the file already exists, please ignore this step):

\begin{sphinxVerbatim}[commandchars=\\\{\}]
\PYG{n}{touch} \PYG{o}{\PYGZti{}}\PYG{o}{/}\PYG{o}{.}\PYG{n}{bash\PYGZus{}profile}
\end{sphinxVerbatim}

\sphinxAtStartPar
\sphinxstylestrong{2. Add environment variables}

\sphinxAtStartPar
First, please execute the following command to open the \sphinxcode{\sphinxupquote{\textquotesingle{}.bash\_profile\textquotesingle{}}} file.

\begin{sphinxVerbatim}[commandchars=\\\{\}]
\PYG{n+nb}{open} \PYG{o}{\PYGZti{}}\PYG{o}{/}\PYG{o}{.}\PYG{n}{bash\PYGZus{}profile}
\end{sphinxVerbatim}

\sphinxAtStartPar
Next, please add the following content to the \sphinxcode{\sphinxupquote{\textquotesingle{}.bash\_profile\textquotesingle{}}} using any editors that you prefered:

\begin{figure}[H]
\centering

\noindent\sphinxincludegraphics[scale=0.6]{{copt-mac_installenv_bash}.png}
\end{figure}

\sphinxAtStartPar
\sphinxstylestrong{Note:} There can be no spaces on both sides of the equal sign.

\sphinxAtStartPar
\sphinxstylestrong{3. Check the environment variables to take effect}

\sphinxAtStartPar
Save the above modification and exit, the user can execute the following command in the terminal to view the modified \sphinxcode{\sphinxupquote{\textquotesingle{}.bash\_profile\textquotesingle{}}} file:

\begin{sphinxVerbatim}[commandchars=\\\{\}]
\PYG{n}{cat} \PYG{o}{\PYGZti{}}\PYG{o}{/}\PYG{o}{.}\PYG{n}{bash\PYGZus{}profile}
\end{sphinxVerbatim}

\sphinxAtStartPar
If the modification is successful, the content of the file output should include the newly added information above.

\sphinxAtStartPar
Then, the user needs to execute the following command on the terminal to make the above modification take effect.

\begin{sphinxVerbatim}[commandchars=\\\{\}]
\PYG{n}{source} \PYG{o}{\PYGZti{}}\PYG{o}{/}\PYG{o}{.}\PYG{n}{bash\PYGZus{}profile}
\end{sphinxVerbatim}
\phantomsection\label{\detokenize{install:parinstallenv-zsh}}
\sphinxAtStartPar
\sphinxstylestrong{ZSH terminal}

\sphinxAtStartPar
Likewise, there are three main steps to configure environment variables under the ZSH terminal:

\sphinxAtStartPar
\sphinxstylestrong{1. Check ‘zshrc’ file}

\sphinxAtStartPar
Enter the following command to output all files in the user directory
and check whether there is a \sphinxcode{\sphinxupquote{\textquotesingle{}.zshrc\textquotesingle{}}} hidden file:

\begin{sphinxVerbatim}[commandchars=\\\{\}]
\PYG{n}{ls} \PYG{o}{\PYGZhy{}}\PYG{n}{a}
\end{sphinxVerbatim}

\sphinxAtStartPar
If the file does not exist, execute the following command to create an empty \sphinxcode{\sphinxupquote{\textquotesingle{}.zshrc\textquotesingle{}}} file at your own (if the file already exists, please ignore this step):

\begin{sphinxVerbatim}[commandchars=\\\{\}]
\PYG{n}{touch} \PYG{o}{\PYGZti{}}\PYG{o}{/}\PYG{o}{.}\PYG{n}{zshrc}
\end{sphinxVerbatim}

\sphinxAtStartPar
\sphinxstylestrong{2. Add environment variables}

\sphinxAtStartPar
First, please execute the following command to open the \sphinxcode{\sphinxupquote{\textquotesingle{}.zshrc\textquotesingle{}}} file.

\begin{sphinxVerbatim}[commandchars=\\\{\}]
\PYG{n+nb}{open} \PYG{o}{\PYGZti{}}\PYG{o}{/}\PYG{o}{.}\PYG{n}{zshrc}
\end{sphinxVerbatim}

\sphinxAtStartPar
Next, please add the following content to the \sphinxcode{\sphinxupquote{\textquotesingle{}.zshrc\textquotesingle{}}} using any editors that you prefered:

\begin{figure}[H]
\centering

\noindent\sphinxincludegraphics[scale=0.6]{{copt-mac_installenv_zsh}.png}
\end{figure}

\sphinxAtStartPar
\sphinxstylestrong{3. Check the environment variables to take effect}

\sphinxAtStartPar
Save the above modification and exit, the user can execute the following command in the terminal to view the modified \sphinxcode{\sphinxupquote{\textquotesingle{}.zshrc\textquotesingle{}}} file:

\begin{sphinxVerbatim}[commandchars=\\\{\}]
\PYG{n}{cat} \PYG{o}{\PYGZti{}}\PYG{o}{/}\PYG{o}{.}\PYG{n}{zshrc}
\end{sphinxVerbatim}

\sphinxAtStartPar
If the modification is successful, the content of the file output should include the newly added information above.

\sphinxAtStartPar
Then, the user needs to execute the following command on the terminal to make the above modification take effect.

\begin{sphinxVerbatim}[commandchars=\\\{\}]
\PYG{n}{source} \PYG{o}{\PYGZti{}}\PYG{o}{/}\PYG{o}{.}\PYG{n}{zshrc}
\end{sphinxVerbatim}

\subsubsection{Verify environment variable configuration}
\label{\detokenize{install:verify-environment-variable-configuration}}\label{\detokenize{install:parinstallenv-valid}}
\sphinxAtStartPar
Now open a new terminal to check if the previous modifications work by
executing commands below respectively:

\begin{sphinxVerbatim}[commandchars=\\\{\}]
echo \PYGZdl{}COPT\PYGZus{}HOME
echo \PYGZdl{}COPT\PYGZus{}LICENSE\PYGZus{}DIR

echo \PYGZdl{}PATH
echo \PYGZdl{}DYLD\PYGZus{}LIBRARY\PYGZus{}PATH
\end{sphinxVerbatim}

\sphinxAtStartPar
If the terminal outputs respectively:
\begin{sphinxalltt}
/Applications/copt80
\end{sphinxalltt}
\begin{sphinxalltt}
/Applications/copt80
\end{sphinxalltt}
\begin{sphinxalltt}
/Applications/copt80/bin:\$PATH
\end{sphinxalltt}
\begin{sphinxalltt}
/Applications/copt80/lib:\$DYLD\_LIBRARY\_PATH
\end{sphinxalltt}

\sphinxAtStartPar
It means that the COPT\sphinxhyphen{}related environment variables are configured successfully.

\sphinxAtStartPar
\sphinxstylestrong{Notes:} The environment variables \sphinxcode{\sphinxupquote{\$PATH}} and \sphinxcode{\sphinxupquote{\$DYLD\_LIBRARY\_PATH}} may display different contents on different computers. However, the set COPT\sphinxhyphen{}related environment variables should be displayed normally to indicate that the configuration is successful.

\sphinxAtStartPar
If the user checks that the COPT\sphinxhyphen{}related environment variables have been successfully added, please carefully read the user agreement document \sphinxcode{\sphinxupquote{\textquotesingle{}copt\sphinxhyphen{}eula\_cn.pdf\textquotesingle{}}} in the installation directory. If you accept the the agreement, please go to {\hyperref[\detokenize{install:parcoptgetlic}]{\sphinxcrossref{\DUrole{std,std-ref}{Setting Up License}}}} for license issues.

\subsubsection{MacOS Security Checkup}
\label{\detokenize{install:macos-security-checkup}}
\sphinxAtStartPar
For MacOS 10.15 (Catalina), if users received error message below:

\begin{sphinxVerbatim}[commandchars=\\\{\}]
\PYGZdq{}libcopt.dylib\PYGZdq{} cannot be opened because the developer cannot be verified.
macOS cannot verify that this app is free from malware.
\end{sphinxVerbatim}

\sphinxAtStartPar
or

\begin{sphinxVerbatim}[commandchars=\\\{\}]
\PYGZdq{}libcopt\PYGZus{}cpp.dylib\PYGZdq{}: dlopen(libcopt\PYGZus{}cpp.dylib,6): no suitable image found.
Did find: libcopt\PYGZus{}cpp.dylib: code signature in (libcopt\PYGZus{}cpp.dylib) not valid
for use in process using Library Validation: library load disallowed by
system policy.
\end{sphinxVerbatim}

\sphinxAtStartPar
Then execute the following commands:
\begin{sphinxalltt}
xattr \sphinxhyphen{}d com.apple.quarantine CardinalOptimizer\sphinxhyphen{}8.0.1\sphinxhyphen{}universal\_mac.dmg
xattr \sphinxhyphen{}d com.apple.quarantine CardinalOptimizer\sphinxhyphen{}8.0.1\sphinxhyphen{}universal\_mac.tar.gz
\end{sphinxalltt}

\sphinxAtStartPar
or
\begin{sphinxalltt}
xattr \sphinxhyphen{}dr com.apple.quarantine /Applications/copt80
\end{sphinxalltt}

\sphinxAtStartPar
to disable security check of MacOS.

\subsection{Linux}
\label{\detokenize{install:linux}}
\sphinxAtStartPar
For Linux platform, we support X86 and ARM64 architectures and provide packages in \sphinxcode{\sphinxupquote{GZIP}} format.
The installation package names differ from different architectures:
\begin{itemize}
\item {} 
\sphinxAtStartPar
Linux\sphinxhyphen{}X86: CardinalOptimizer\sphinxhyphen{}8.0.1\sphinxhyphen{}lnx64.tar.gz

\item {} 
\sphinxAtStartPar
Linux\sphinxhyphen{}ARM64: CardinalOptimizer\sphinxhyphen{}8.0.1\sphinxhyphen{}aarch64\_lnx.tar.gz

\end{itemize}

\sphinxAtStartPar
If you download the software, e.g. CardinalOptimizer\sphinxhyphen{}8.0.1\sphinxhyphen{}lnx64.tar.gz for
64\sphinxhyphen{}bit version of COPT 8.0.1, just type commands below in shell to extract it
to any directories:
\begin{sphinxalltt}
tar \sphinxhyphen{}xzf CardinalOptimizer\sphinxhyphen{}8.0.1\sphinxhyphen{}lnx64.tar.gz
\end{sphinxalltt}

\sphinxAtStartPar
You will get a folder named copt80 in current directory, you can move it
to any other directories as you like. We recommend users move it to \sphinxcode{\sphinxupquote{\textquotesingle{}/opt\textquotesingle{}}}
directory by typing commands below in shell:
\begin{sphinxalltt}
sudo mv copt80 /opt/
\end{sphinxalltt}

\sphinxAtStartPar
Note that the above command requires \sphinxstylestrong{root privilege} to execute.

\sphinxAtStartPar
The next step users need to do is to set the required environment variables,
by adding the following commands to the \sphinxcode{\sphinxupquote{\textquotesingle{}.bashrc\textquotesingle{}}} file in your \sphinxcode{\sphinxupquote{\$HOME}}
directory using any editors that you prefered:

\begin{figure}[H]
\centering

\noindent\sphinxincludegraphics[scale=0.7]{{copt-linux_installenv}.png}
\end{figure}

\sphinxAtStartPar
Remember to save your modifications to the \sphinxcode{\sphinxupquote{\textquotesingle{}.bashrc\textquotesingle{}}} file, and open a new
terminal to check if it works by executing commands below respectively:

\begin{sphinxVerbatim}[commandchars=\\\{\}]
echo \PYGZdl{}COPT\PYGZus{}HOME
echo \PYGZdl{}COPT\PYGZus{}LICENSE\PYGZus{}DIR

echo \PYGZdl{}PATH
echo \PYGZdl{}LD\PYGZus{}LIBRARY\PYGZus{}PATH
\end{sphinxVerbatim}

\sphinxAtStartPar
If you accept the \sphinxcode{\sphinxupquote{copt\sphinxhyphen{}eula\_en.pdf}} in installation directory, then please
go to {\hyperref[\detokenize{install:parcoptgetlic}]{\sphinxcrossref{\DUrole{std,std-ref}{Setting Up License}}}} for license issues.

\section{Setting Up License}
\label{\detokenize{install:setting-up-license}}\label{\detokenize{install:parcoptgetlic}}
\sphinxAtStartPar
The Cardinal Optimizers requires a valid license to work properly.
We offer different types of licenses most suitable for user’s needs.
All users should read this section carefully. If you encounter any
problem about license, feel free to contact \sphinxhref{mailto:coptsupport@shanshu.ai}{coptsupport@shanshu.ai}.

\sphinxAtStartPar
The license specifically includes: \sphinxcode{\sphinxupquote{license.dat}} and \sphinxcode{\sphinxupquote{license.key}} two files.
Starting from COPT 6.5, if you apply for personal license applications via our official website,
we will send both files directly as attachments (no need to obtain it yourself).

\sphinxAtStartPar
Users can \sphinxstylestrong{directly download two license files to the local computer, skip the following steps of obtaining and verifying the license} , and can go straight to the {\hyperref[\detokenize{install:parcoptinslic}]{\sphinxcrossref{\DUrole{std,std-ref}{Install license}}}} step.

\sphinxAtStartPar
The configuration of the license mainly includes the following three steps:
\begin{itemize}
\item {} 
\sphinxAtStartPar
{\hyperref[\detokenize{install:parcoptgetlic-obtain}]{\sphinxcrossref{\DUrole{std,std-ref}{Obtaining License (skippable)}}}}

\item {} 
\sphinxAtStartPar
{\hyperref[\detokenize{install:parcoptgetlic-valid}]{\sphinxcrossref{\DUrole{std,std-ref}{Verifying License (skippable)}}}}

\item {} 
\sphinxAtStartPar
{\hyperref[\detokenize{install:parcoptinslic}]{\sphinxcrossref{\DUrole{std,std-ref}{Installing License}}}}

\end{itemize}

\subsection{Obtaining License}
\label{\detokenize{install:obtaining-license}}\label{\detokenize{install:parcoptgetlic-obtain}}
\sphinxAtStartPar
Once the registration is done, a license key is sent to
users. It is a unique token binding with user’s registration
information. Afterwards, users may
run the \sphinxcode{\sphinxupquote{copt\_licgen}} tool, shipped with Cardinal Optimizer, to
obtain license files from COPT licensing server (Internet connection required) .

\sphinxAtStartPar
\sphinxstylestrong{Note}: If the user has already obtained the two files \sphinxcode{\sphinxupquote{license.dat}} and \sphinxcode{\sphinxupquote{license.key}}, there is no need to repeat the following steps to obtain them again.

\sphinxAtStartPar
The following notes show you how to play with the \sphinxcode{\sphinxupquote{copt\_licgen}} tool.

\sphinxAtStartPar
For Windows system, open a new cmd window. The current path is the user directory, and the path is as follows:
\sphinxcode{\sphinxupquote{"C:\textbackslash{}Users\textbackslash{}copt80"}}.

\sphinxAtStartPar
For MacOS and Linux systems, open a new terminal. The current path is the user directory, represented by the symbol \sphinxcode{\sphinxupquote{\textasciitilde{}}}.

\sphinxAtStartPar
To obtain COPT license files, execute the following command using the option
\sphinxcode{\sphinxupquote{\textquotesingle{}\sphinxhyphen{}key\textquotesingle{}}} and the license key as argument. Below is an example, assuming
the license key is \sphinxcode{\sphinxupquote{\textquotesingle{}19200817f147gd9f60abc791def047fb\textquotesingle{}}}:

\begin{sphinxVerbatim}[commandchars=\\\{\}]
copt\PYGZus{}licgen \PYGZhy{}key 19200817f147gd9f60abc791def047fb
\end{sphinxVerbatim}

\sphinxAtStartPar
If the license key is saved to \sphinxcode{\sphinxupquote{key.txt}} file in format of \sphinxcode{\sphinxupquote{\textquotesingle{}KEY=xxx\textquotesingle{}}},
which resides in the same place
as \sphinxcode{\sphinxupquote{copt\_licgen}}, execute the following command using the option \sphinxcode{\sphinxupquote{\textquotesingle{}\sphinxhyphen{}file\textquotesingle{}}}
and \sphinxcode{\sphinxupquote{\textquotesingle{}key.txt\textquotesingle{}}} as arugment.

\begin{sphinxVerbatim}[commandchars=\\\{\}]
copt\PYGZus{}licgen \PYGZhy{}file key.txt
\end{sphinxVerbatim}

\sphinxAtStartPar
We recommend to use the \sphinxstylestrong{first way} to get your license files.

\sphinxAtStartPar
If the authorization server verifies the license, then it generates \sphinxcode{\sphinxupquote{license.dat}} and \sphinxcode{\sphinxupquote{license.key}}
and download it to the user’s computer. The default download directory is the current working directory.
\begin{sphinxalltt}
copt\_licgen \sphinxhyphen{}key 19200817f147gd9f60abc791def047fb
  {[}Info{]} Cardinal Optimizer   COPT v8.0.1 20240304
  {[}Info{]} Use specific key 19200817f147gd9f60abc791def047fb
  {[}Info{]} \sphinxstylestrong{*} get new COPT license from licensing server \sphinxstylestrong{*}
  {[}Info{]} Write to license.dat
  {[}Info{]} Write to license.key
  {[}Info{]} Received new license files from server
  {[}Info{]} Done !!!
\end{sphinxalltt}

\sphinxAtStartPar
\sphinxstylestrong{Note:} Users do not need to have internet connection to run COPT. However, obtaining license itself requires internet connection.
If you encounter any problem, please feel free to contact \sphinxhref{mailto:coptsupport@shanshu.ai}{coptsupport@shanshu.ai}.

\subsection{Verifying License}
\label{\detokenize{install:verifying-license}}\label{\detokenize{install:parcoptgetlic-valid}}
\sphinxAtStartPar
If the license key binding to registration information is verified by COPT
license server, two license files, \sphinxcode{\sphinxupquote{license.dat}} and \sphinxcode{\sphinxupquote{license.key}}, are
downloaded to the current working directory. To double check two
license files are valid in current version of COPT, execute the following command
with the option \sphinxcode{\sphinxupquote{\textquotesingle{}\sphinxhyphen{}v\textquotesingle{}}}. Note that this command requires both \sphinxcode{\sphinxupquote{license.dat}}
and \sphinxcode{\sphinxupquote{license.key}} existing in the current working directory.

\begin{sphinxVerbatim}[commandchars=\\\{\}]
copt\PYGZus{}licgen \PYGZhy{}v
\end{sphinxVerbatim}

\sphinxAtStartPar
If you see log information similar to the following, you have obtained and verified the license files sucessfully.
\begin{sphinxalltt}
copt\_licgen \sphinxhyphen{}v
  {[}Info{]} Cardinal Optimizer   COPT v8.0.1 20240304
  {[}Info{]} Run local validation
  {[}Info{]} Read license.dat
  {[}Info{]} Read license.key
  {[}Info{]} Expiry : Tue 2030\sphinxhyphen{}12\sphinxhyphen{}31 00:00:00 +0800
  {[}Info{]} Local validation result: Succeeded
  {[}Info{]} Done !!!
\end{sphinxalltt}

\subsection{Installing License}
\label{\detokenize{install:installing-license}}\label{\detokenize{install:parcoptinslic}}
\sphinxAtStartPar
Once you obtained the \sphinxcode{\sphinxupquote{license.dat}} and \sphinxcode{\sphinxupquote{license.key}} files from COPT
license server and verified they work as expected, setting up them is just
as simple as moving them to the same directory as the COPT dynamic libary
or {\hyperref[\detokenize{cmdline:chapcmdline}]{\sphinxcrossref{\DUrole{std,std-ref}{COPT Command\sphinxhyphen{}Line}}}}. Note this installation is applied only to current
COPT.

\sphinxAtStartPar
The following two ways of installation applied to all versions of COPT on your machine.

\subsubsection{\sphinxstylestrong{Via HOME directory method (recommended)}}
\label{\detokenize{install:via-home-directory-method-recommended}}\label{\detokenize{install:parcoptinslic-userlib}}
\sphinxAtStartPar
The simplest way to install COPT license is to create a folder of \sphinxcode{\sphinxupquote{copt}}
in your \sphinxcode{\sphinxupquote{HOME}} directory, and move the authenticated license files
\sphinxcode{\sphinxupquote{license.dat}} and \sphinxcode{\sphinxupquote{license.key}} to the new folder \sphinxcode{\sphinxupquote{copt}}.

\sphinxAtStartPar
\sphinxstylestrong{Note:} The folder name \sphinxcode{\sphinxupquote{"copt"}} is case sensitive and it must be lowercase.

\sphinxAtStartPar
The \sphinxcode{\sphinxupquote{HOME}} directory names are slightly different under different systems,
please refer to the corresponding license installation steps according to your operating system:
\begin{itemize}
\item {} 
\sphinxAtStartPar
{\hyperref[\detokenize{install:parcoptgetlic-win}]{\sphinxcrossref{\DUrole{std,std-ref}{Windows}}}}

\item {} 
\sphinxAtStartPar
{\hyperref[\detokenize{install:parcoptgetlic-mac}]{\sphinxcrossref{\DUrole{std,std-ref}{MacOS}}}}

\item {} 
\sphinxAtStartPar
{\hyperref[\detokenize{install:parcoptgetlic-lin}]{\sphinxcrossref{\DUrole{std,std-ref}{Linux}}}}

\end{itemize}
\phantomsection\label{\detokenize{install:parcoptgetlic-win}}
\sphinxAtStartPar
\sphinxstylestrong{Windows}

\sphinxAtStartPar
The \sphinxcode{\sphinxupquote{HOME}} directory looks like \sphinxcode{\sphinxupquote{"C:\textbackslash{}Users\textbackslash{}username\textbackslash{}"}} on Windows. User can manually move
the two license files \sphinxcode{\sphinxupquote{license.dat}} and \sphinxcode{\sphinxupquote{license.key}} to the directory \sphinxcode{\sphinxupquote{"C:\textbackslash{}Users\textbackslash{}username\textbackslash{}copt"}}.

\sphinxAtStartPar
Finally, please check whether the license files \sphinxcode{\sphinxupquote{license.dat}} and \sphinxcode{\sphinxupquote{license.key}} are located in this directory, as shown in the figure below, indicating that the license has been installed correctly:
\begin{quote}

\begin{figure}[H]
\centering

\noindent\sphinxincludegraphics[scale=0.4]{{copt-win_license}.png}
\end{figure}
\end{quote}
\phantomsection\label{\detokenize{install:parcoptgetlic-mac}}
\sphinxAtStartPar
\sphinxstylestrong{MacOS}

\sphinxAtStartPar
The \sphinxcode{\sphinxupquote{HOME}} directory looks like \sphinxcode{\sphinxupquote{"/home/your username/"}}
on MacOS. User can manually move
the two license files to \sphinxcode{\sphinxupquote{"/home/your username/copt/"}}.

\sphinxAtStartPar
Finally, please check whether the license files \sphinxcode{\sphinxupquote{license.dat}} and \sphinxcode{\sphinxupquote{license.key}} are located in this directory, as shown in the figure below, indicating that the license has been installed correctly:
\begin{quote}

\begin{figure}[H]
\centering

\noindent\sphinxincludegraphics[scale=0.4]{{copt-mac_license}.png}
\end{figure}
\end{quote}
\phantomsection\label{\detokenize{install:parcoptgetlic-lin}}
\sphinxAtStartPar
\sphinxstylestrong{Linux}

\sphinxAtStartPar
The \sphinxcode{\sphinxupquote{HOME}} directory looks like \sphinxcode{\sphinxupquote{"/home/your username/"}}
on Linux.

\sphinxAtStartPar
User can execute the following command in the terminal to move the two license files \sphinxcode{\sphinxupquote{license.dat}} and \sphinxcode{\sphinxupquote{license.key}} to the directory \sphinxcode{\sphinxupquote{"/home/your username/copt"}}.

\begin{sphinxVerbatim}[commandchars=\\\{\}]
\PYG{n}{mv} \PYG{n}{license}\PYG{o}{.}\PYG{o}{*}  \PYG{o}{\PYGZti{}}\PYG{o}{/}\PYG{n}{copt}\PYG{o}{/}
\end{sphinxVerbatim}

\sphinxAtStartPar
Next, please check that the license files \sphinxcode{\sphinxupquote{license.dat}} and \sphinxcode{\sphinxupquote{license.key}} are located in \sphinxcode{\sphinxupquote{"/home/username/copt"}}
directory, the command is:

\begin{sphinxVerbatim}[commandchars=\\\{\}]
\PYG{n}{ls} \PYG{o}{\PYGZti{}}\PYG{o}{/}\PYG{n}{copt}\PYG{o}{/}
\end{sphinxVerbatim}

\sphinxAtStartPar
If the terminal output shows that there are two files \sphinxcode{\sphinxupquote{license.dat}} and \sphinxcode{\sphinxupquote{license.key}}, it means that the license has been installed correctly, otherwise the installation fails.

\subsubsection{\sphinxstylestrong{Via environment variable method}}
\label{\detokenize{install:via-environment-variable-method}}
\sphinxAtStartPar
Alternatively, for users who prefer having licenses in a customized folder, they can set
environment variable \sphinxcode{\sphinxupquote{COPT\_LICENSE\_DIR}} to the customized folder. You may refer
to {\hyperref[\detokenize{install:parcoptinstall}]{\sphinxcrossref{\DUrole{std,std-ref}{Software Installation}}}} for how to set environment variable on Windows, Linux
and MacOS.

\sphinxAtStartPar
In addition, please double check that if license files \sphinxcode{\sphinxupquote{license.dat}} and \sphinxcode{\sphinxupquote{license.key}}
locate in path specified by environmental variable \sphinxcode{\sphinxupquote{COPT\_LICENSE\_DIR}}.

\sphinxAtStartPar
The viewing and operation of environment variables are slightly different under different systems, please refer to the following installation steps for your own operating system:
\begin{itemize}
\item {} 
\sphinxAtStartPar
{\hyperref[\detokenize{install:parcoptgetlic2-win}]{\sphinxcrossref{\DUrole{std,std-ref}{Windows}}}}

\item {} 
\sphinxAtStartPar
{\hyperref[\detokenize{install:parcoptgetlic2-mac}]{\sphinxcrossref{\DUrole{std,std-ref}{MacOS}}}}

\item {} 
\sphinxAtStartPar
{\hyperref[\detokenize{install:parcoptgetlic2-lin}]{\sphinxcrossref{\DUrole{std,std-ref}{Linux}}}}

\end{itemize}
\phantomsection\label{\detokenize{install:parcoptgetlic2-win}}
\sphinxAtStartPar
\sphinxstylestrong{Windows}

\sphinxAtStartPar
For Windows systems, execute the following command to view the path pointed to by the environment variable \sphinxcode{\sphinxupquote{COPT\_LICENSE\_DIR}} :

\begin{sphinxVerbatim}[commandchars=\\\{\}]
\PYG{n}{echo} \PYG{o}{\PYGZpc{}}\PYG{n}{COPT\PYGZus{}LICENSE\PYGZus{}DIR}\PYG{o}{\PYGZpc{}}
\end{sphinxVerbatim}

\sphinxAtStartPar
\sphinxstylestrong{Note:} If there is no output on the terminal, it means that the COPT\sphinxhyphen{}related environment variables have not been configured, please check {\hyperref[\detokenize{install:parcoptinstall}]{\sphinxcrossref{\DUrole{std,std-ref}{Software installation}}}}.

\sphinxAtStartPar
Then move the license files \sphinxcode{\sphinxupquote{license.dat}} and \sphinxcode{\sphinxupquote{license.key}} to the path \sphinxcode{\sphinxupquote{COPT\_LICENSE\_DIR}}
points to.

\sphinxAtStartPar
Here we assume that the environment variable \sphinxcode{\sphinxupquote{COPT\_LICENSE\_DIR}} points to the default installation directory of COPT. As shown in the figure, it means that the license has been installed correctly:
\begin{quote}

\begin{figure}[H]
\centering

\noindent\sphinxincludegraphics[scale=0.4]{{copt-win_license2}.png}
\end{figure}
\end{quote}
\phantomsection\label{\detokenize{install:parcoptgetlic2-mac}}
\sphinxAtStartPar
\sphinxstylestrong{MacOS}

\sphinxAtStartPar
For MacOS system, you can use the following command to view the path pointed to by the environment variable \sphinxcode{\sphinxupquote{COPT\_LICENSE\_DIR}}:

\begin{sphinxVerbatim}[commandchars=\\\{\}]
echo \PYGZdl{}COPT\PYGZus{}LICENSE\PYGZus{}DIR
\end{sphinxVerbatim}

\sphinxAtStartPar
Then move the license files \sphinxcode{\sphinxupquote{license.dat}} and \sphinxcode{\sphinxupquote{license.key}} to the path \sphinxcode{\sphinxupquote{COPT\_LICENSE\_DIR}}
points to.

\sphinxAtStartPar
Here we assume that the environment variable \sphinxcode{\sphinxupquote{COPT\_LICENSE\_DIR}} points to the default installation directory of COPT. As shown in the figure, it means that the license has been installed correctly:
\begin{quote}

\begin{figure}[H]
\centering

\noindent\sphinxincludegraphics[scale=0.4]{{copt-mac_license2}.png}
\end{figure}
\end{quote}
\phantomsection\label{\detokenize{install:parcoptgetlic2-lin}}
\sphinxAtStartPar
\sphinxstylestrong{Linux}

\sphinxAtStartPar
For Linux systems, user can use the following command to view the path pointed to by the environment variable \sphinxcode{\sphinxupquote{COPT\_LICENSE\_DIR}}:

\begin{sphinxVerbatim}[commandchars=\\\{\}]
echo \PYGZdl{}COPT\PYGZus{}LICENSE\PYGZus{}DIR
\end{sphinxVerbatim}

\sphinxAtStartPar
Then user can execute the following command to move the license files \sphinxcode{\sphinxupquote{license.dat}} and \sphinxcode{\sphinxupquote{license.key}} to to the path \sphinxcode{\sphinxupquote{COPT\_LICENSE\_DIR}} points to.

\begin{sphinxVerbatim}[commandchars=\\\{\}]
mv license.* \PYGZdl{}COPT\PYGZus{}LICENSE\PYGZus{}DIR/
\end{sphinxVerbatim}

\sphinxAtStartPar
\sphinxstylestrong{Note:} For Linux systems, if the path pointed to by the environment variable \sphinxcode{\sphinxupquote{COPT\_LICENSE\_DIR}} is /opt/copt80, user need to perform the move operation with \sphinxstylestrong{root authority} , and the command is:

\begin{sphinxVerbatim}[commandchars=\\\{\}]
sudo mv license.* \PYGZdl{}COPT\PYGZus{}LICENSE\PYGZus{}DIR/
\end{sphinxVerbatim}

\sphinxAtStartPar
In addition, please double check that if license files \sphinxcode{\sphinxupquote{license.dat}} and \sphinxcode{\sphinxupquote{license.key}}
locate in path specified by environmental variable \sphinxcode{\sphinxupquote{COPT\_LICENSE\_DIR}} , and the command is:

\begin{sphinxVerbatim}[commandchars=\\\{\}]
ls \PYGZdl{}COPT\PYGZus{}LICENSE\PYGZus{}DIR/
\end{sphinxVerbatim}

\sphinxAtStartPar
If the terminal displays the files \sphinxcode{\sphinxupquote{license.dat}} and \sphinxcode{\sphinxupquote{license.key}}, it means that the license has been installed correctly, otherwise the license installation fails.

\sphinxAtStartPar
\sphinxstylestrong{Note:} If there are license files \sphinxcode{\sphinxupquote{license.dat}} and \sphinxcode{\sphinxupquote{license.key}} in the \sphinxcode{\sphinxupquote{copt}} folder of the user directory and the directory pointed to by the environment variable, the former will be checked and used the former.

\subsection{Others}
\label{\detokenize{install:others}}
\sphinxAtStartPar
Basically, each type of licenses includes two license files: \sphinxcode{\sphinxupquote{license.dat}}
and \sphinxcode{\sphinxupquote{license.key}}, each of which has digital signature to protect its
content. When invoking COPT command\sphinxhyphen{}line or loading COPT dynamic library
to solve an optimization problem, the public RSA key stored in \sphinxcode{\sphinxupquote{\textquotesingle{}license.key\textquotesingle{}}}
is used to verify signature in \sphinxcode{\sphinxupquote{\textquotesingle{}license.dat\textquotesingle{}}}.

\sphinxAtStartPar
Afterwards, license data in format of key\sphinxhyphen{}value pair is parsed to verify
whether it is a legal license. Below is a sample of license data.
\begin{sphinxalltt}
\#\# COPT LICENSE FILE \#\#

USER = Trial User
MAC = 44:05:99:31:41:C2
CPUID = BFEBFBFF000706E5
EXPIRY = 2030\sphinxhyphen{}12\sphinxhyphen{}31
VERSION = 8.0.1
NOTE = Free For Trial
\end{sphinxalltt}

\sphinxAtStartPar
Note: Please make sure that the \sphinxcode{\sphinxupquote{VERSION}} field in the license file is the same major version as the software installed on the machine. Otherwise, it may cause error:
\sphinxcode{\sphinxupquote{Invalid signature in public key file}} .

\section{Verify the installation and configuration}
\label{\detokenize{install:verify-the-installation-and-configuration}}\label{\detokenize{install:parcoptvalid}}
\sphinxAtStartPar
After the user completes the above steps of software installation and license configuration, please enter the COPT command line tool to verify whether the installation is successful.

\sphinxAtStartPar
In MacOS and Linux systems, please open the terminal (in Windows system, please open the cmd window), and enter the \sphinxcode{\sphinxupquote{copt\_cmd}} command to use the COPT command line tool.

\sphinxAtStartPar
If the user sees the output as shown below (without any error message), the COPT command line tool can be entered normally, which means that the previous software and license configuration has been successfully completed.

\begin{sphinxVerbatim}[commandchars=\\\{\}]
Cardinal Optimizer v8.0.0. Build date Oct 17 2024
Copyright Cardinal Operations 2024. All Rights Reserved

COPT\PYGZgt{}
\end{sphinxVerbatim}

\sphinxAtStartPar
If you fail to enter the COPT command line tool normally, please refer to {\hyperref[\detokenize{faq:chapfaq}]{\sphinxcrossref{\DUrole{std,std-ref}{FAQ}}}} according to the error message and check whether the above steps of installation and license configuration are completed correctly.

\sphinxAtStartPar
After the installation is successfully completed, according to the interface of COPT as needed, please continue to refer to {\hyperref[\detokenize{quickstart:chapquickstart}]{\sphinxcrossref{\DUrole{std,std-ref}{Quick Start}}}}, and check the installation and applying methods of each interface.

\section{Upgrade}
\label{\detokenize{install:upgrade}}
\sphinxAtStartPar
The upgrade steps of COPT mainly include two parts: 1. Software upgrade; 2. License file update.

\sphinxAtStartPar
\sphinxstylestrong{1. Uninstall the old version of COPT}

\sphinxAtStartPar
Firstly, if the user does not need to run multiple versions at the same time, the old version of the software can be uninstalled first.

\sphinxAtStartPar
For different operating systems, the methods of uninstalling COPT are slightly different:

\sphinxAtStartPar
In Windows system, please open the old COPT installation package (named like copt80 ), please double\sphinxhyphen{}click to run \sphinxcode{\sphinxupquote{unins000.exe}} , and then the software can be uninstalled automatically. In MacOS and Linux systems, please delete the old version of the COPT installation package directly to complete the software uninstallation.

\sphinxAtStartPar
\sphinxstylestrong{2. Install the new version of COPT}

\sphinxAtStartPar
If it is an upgrade between major versions (for example: COPT 6.5 \sphinxhyphen{}\textgreater{} COPT 7.0), the user needs to re\sphinxhyphen{}apply for new version license.
Next, the user needs to reinstall the new version of COPT. For the same steps, please refer to {\hyperref[\detokenize{install:parcoptinstall}]{\sphinxcrossref{\DUrole{std,std-ref}{Software Installation}}}}.

\sphinxAtStartPar
\sphinxstylestrong{3. Update the license files}

\sphinxAtStartPar
In addition, after the software update is complete, the license file needs to be updated.
We recommend that the user configure the license file via {\hyperref[\detokenize{install:parcoptinslic-userlib}]{\sphinxcrossref{\DUrole{std,std-ref}{user directory}}}}. The user only needs to delete the license file of the old version under the copt folder and move the new license file into it to complete the license update.

\sphinxAtStartPar
After completing the above steps, users can refer to {\hyperref[\detokenize{install:parcoptvalid}]{\sphinxcrossref{\DUrole{std,std-ref}{Verify installation and configuration completion}}}} to observe whether the version number in the banner information output by the COPT command line tool is a new version.

\begin{sphinxadmonition}{note}{Notes}
\begin{enumerate}
\sphinxsetlistlabels{\arabic}{enumi}{enumii}{}{.}%
\item {} 
\sphinxAtStartPar
For the Python interface of COPT, users also need to \sphinxstylestrong{upgrade its
corresponding interface} \sphinxcode{\sphinxupquote{coptpy}} \sphinxstylestrong{at the same time}. If Python is an Anaconda distribution, it can be conducted via \sphinxcode{\sphinxupquote{pip install \sphinxhyphen{}\sphinxhyphen{}upgrade coptpy}} .

\item {} 
\sphinxAtStartPar
The above are the upgrade steps for major version iterations (for example: COPT 6.5 \sphinxhyphen{}\textgreater{} COPT 7.0).
In addition, on the basis of each major version, COPT will also provide patch release (eg: COPT 8.0.1 \sphinxhyphen{}\textgreater{} COPT 8.0.2), for some feature fixes and updates. For upgrades between minor versions, there is no need to update the license files.

\end{enumerate}
\end{sphinxadmonition}

\sphinxstepscope

\chapter{COPT Command\sphinxhyphen{}Line}
\label{\detokenize{cmdline:copt-command-line}}\label{\detokenize{cmdline:chapcmdline}}\label{\detokenize{cmdline::doc}}
\sphinxAtStartPar
The \sphinxstylestrong{Cardinal Optimizer} ships with \sphinxcode{\sphinxupquote{copt\_cmd}} executable on all supported
platforms, which let users solve optimization models in an interactive way.
Before running COPT command\sphinxhyphen{}line, please make sure that you have valid license
installed.

\section{Overview}
\label{\detokenize{cmdline:overview}}\label{\detokenize{cmdline:parcmdline-overview}}
\sphinxAtStartPar
COPT command\sphinxhyphen{}line is a COPT API interpreter that executes commands read from
the standard input or from a script file. COPT command\sphinxhyphen{}line interprets the
following options when it is invoked:
\begin{itemize}
\item {} 
\sphinxAtStartPar
\sphinxcode{\sphinxupquote{\sphinxhyphen{}c}}: If the \sphinxcode{\sphinxupquote{\textquotesingle{}\sphinxhyphen{}c\textquotesingle{}}} option is present, it reads from an inline scripts,
which is a quoted string and specified by the second argument.

\item {} 
\sphinxAtStartPar
\sphinxcode{\sphinxupquote{\sphinxhyphen{}i}}: If the \sphinxcode{\sphinxupquote{\textquotesingle{}\sphinxhyphen{}i\textquotesingle{}}} option is present, it reads from an input script file,
whose path is specified by the second argument.

\item {} 
\sphinxAtStartPar
\sphinxcode{\sphinxupquote{\sphinxhyphen{}o}}: If the \sphinxcode{\sphinxupquote{\textquotesingle{}\sphinxhyphen{}o\textquotesingle{}}} option is present, it reads from standard input, while
each \sphinxstylestrong{valid} command line is written to an output script file, whose path
is specified by the second argument.

\end{itemize}

\sphinxAtStartPar
Regardless of arguments, the tool is interactive. Besides wrapping COPT API
calls, it offers various features to help users move cursor around and
edit lines. We try to provide as much user experience as standard command prompt
(Windows console and Unix terminal).

\section{Edit mode}
\label{\detokenize{cmdline:edit-mode}}
\sphinxAtStartPar
This tool defines a number of commands to position the cursor, edit lines with
combination keys on a standard keyboard. The following notes show you how to use
the most important ones.
\begin{itemize}
\item {} 
\sphinxAtStartPar
\sphinxstylestrong{Basic commands}
\begin{enumerate}
\sphinxsetlistlabels{\arabic}{enumi}{enumii}{}{.}%
\item {} 
\sphinxAtStartPar
\sphinxcode{\sphinxupquote{\textless{}Insert\textgreater{}}}: Toggle between inserting characters and replacing the
existing ones.

\item {} 
\sphinxAtStartPar
\sphinxcode{\sphinxupquote{\textless{}Esc\textgreater{}}}: Discard inputs and move the cursor to the beginning of line.
Press \sphinxcode{\sphinxupquote{\textless{}ESC\textgreater{}}} twice on Linux/Mac platform to do the same thing.

\end{enumerate}

\item {} 
\sphinxAtStartPar
\sphinxstylestrong{Moving around}
\begin{enumerate}
\sphinxsetlistlabels{\arabic}{enumi}{enumii}{}{.}%
\item {} 
\sphinxAtStartPar
\sphinxcode{\sphinxupquote{\textless{}Home\textgreater{}/\textless{}End\textgreater{}}}: Jump to the beginning/end of line.

\item {} 
\sphinxAtStartPar
\sphinxcode{\sphinxupquote{\textless{}Left\textgreater{}/\textless{}Right\textgreater{}}} Arrow: Move the cursor one character to the left/right.

\item {} 
\sphinxAtStartPar
\sphinxcode{\sphinxupquote{\textless{}CTRL\textgreater{}+\textless{}Left\textgreater{}/\textless{}Right\textgreater{}}} Arrow: Move the cursor one word to the left/right.

\end{enumerate}

\item {} 
\sphinxAtStartPar
\sphinxstylestrong{Cut and Paste}

\sphinxAtStartPar
An internal paste buffer is available for the following cut operations.
\begin{enumerate}
\sphinxsetlistlabels{\arabic}{enumi}{enumii}{}{.}%
\item {} 
\sphinxAtStartPar
\sphinxcode{\sphinxupquote{\textless{}Delete\textgreater{}}}: Cut the character under the cursor.

\item {} 
\sphinxAtStartPar
\sphinxcode{\sphinxupquote{\textless{}Backspace\textgreater{}}}: Cut the character before the cursor.

\item {} 
\sphinxAtStartPar
\sphinxcode{\sphinxupquote{\textless{}CTRL\textgreater{}+\textless{}H\textgreater{}}}: Cut from the cursor to the beginning of line.

\item {} 
\sphinxAtStartPar
\sphinxcode{\sphinxupquote{\textless{}CTRL\textgreater{}+\textless{}E\textgreater{}}}: Cut from the cursor to the end of line.

\item {} 
\sphinxAtStartPar
\sphinxcode{\sphinxupquote{\textless{}CTRL\textgreater{}+\textless{}Y\textgreater{}}}: Paste text in paste buffer at the cursor position.

\end{enumerate}

\sphinxAtStartPar
Each of cut operations defines cut direction: cut forward or cut backward.
Obviously, \sphinxcode{\sphinxupquote{\textless{}Delete\textgreater{}}} and \sphinxcode{\sphinxupquote{\textless{}CTRL\textgreater{}+\textless{}E\textgreater{}}} cut forward; \sphinxcode{\sphinxupquote{\textless{}Backspace\textgreater{}}} and
\sphinxcode{\sphinxupquote{\textless{}CTRL\textgreater{}+\textless{}H\textgreater{}}} cut backward. When two consecutive cut operations have the
same cut direction, the cutting text is appended the paste buffer.
Otherwise, the paste buffer is overwritten by the latest chopped text.

\item {} 
\sphinxAtStartPar
\sphinxstylestrong{Command history}
\begin{enumerate}
\sphinxsetlistlabels{\arabic}{enumi}{enumii}{}{.}%
\item {} 
\sphinxAtStartPar
\sphinxcode{\sphinxupquote{\textless{}Up\textgreater{}/\textless{}Down\textgreater{}}}: Move through the history of command lines in the
older/newer direction. The tool remembers the history entry if the last
executed line is in history.

\item {} 
\sphinxAtStartPar
\sphinxcode{\sphinxupquote{\textless{}CTRL\textgreater{}+\textless{}R\textgreater{}}} or \sphinxcode{\sphinxupquote{\textless{}F8\textgreater{}}}: If you know what a previously executed line
starts with, and you want to run it again, type prefix characters and then
press \sphinxcode{\sphinxupquote{\textless{}CTRL\textgreater{}+\textless{}R\textgreater{}}}, or \sphinxcode{\sphinxupquote{\textless{}F8\textgreater{}}} on Windows platform, to iterate through
the history of commands with matching prefix.

\end{enumerate}

\item {} 
\sphinxAtStartPar
\sphinxstylestrong{Tab completion}

\sphinxAtStartPar
Use \sphinxcode{\sphinxupquote{\textless{}Tab\textgreater{}}} to complete shell commands, COPT paramters/attributes,
or files under specified path. To cycle though multiple matches, just repeat
pressing \sphinxcode{\sphinxupquote{\textless{}Tab\textgreater{}}}.
\begin{enumerate}
\sphinxsetlistlabels{\arabic}{enumi}{enumii}{}{.}%
\item {} 
\sphinxAtStartPar
If the cursor is over or right after the first word on the current
command line, press \sphinxcode{\sphinxupquote{\textless{}Tab\textgreater{}}} to complete available shell commands with
matching prefix (from the cursor to the first character of word).

\item {} 
\sphinxAtStartPar
Otherwise, press \sphinxcode{\sphinxupquote{\textless{}Tab\textgreater{}}} to complete COPT parameters/attributes, or file
names under path with matching prefix. Specifically, if the prefix
matches with COPT parameters as well as file under current working
directory, only COPT parameters will be listed. In this case, to iterate
file names, add relative path \sphinxcode{\sphinxupquote{\textquotesingle{}./\textquotesingle{}}} to start with.

\item {} 
\sphinxAtStartPar
For convenience, tab completion ignores case and support asterisk (\sphinxcode{\sphinxupquote{*}})
as wildcard to match file and directory pattern.

\item {} 
\sphinxAtStartPar
\sphinxcode{\sphinxupquote{\textless{}Shift\textgreater{}+\textless{}Tab\textgreater{}}}: Complete the next one in an opposite direction.

\end{enumerate}

\end{itemize}

\section{Script mode}
\label{\detokenize{cmdline:script-mode}}
\sphinxAtStartPar
There are two approaches to run scripts, a batch of commands, in COPT
command\sphinxhyphen{}line. One is to save scripts as a text file. The other is called
\sphinxstylestrong{inline scripts}, that is, a quoted string of commands separated by \sphinxcode{\sphinxupquote{\textquotesingle{};\textquotesingle{}}}.
Both of them can be loaded when COPT command\sphinxhyphen{}line is invoked, or loaded on fly
in the edit mode (see shell command \sphinxcode{\sphinxupquote{\textquotesingle{}load\textquotesingle{}}}). Below describes more details
about loading scripts as arguments.

\sphinxAtStartPar
This tool allows users to load a script file to do a batch job automatically.
As mentioned in {\hyperref[\detokenize{cmdline:parcmdline-overview}]{\sphinxcrossref{\DUrole{std,std-ref}{Overview}}}}, a script file is read
when its file path is specified as arguments of the option \sphinxcode{\sphinxupquote{\textquotesingle{}\sphinxhyphen{}i\textquotesingle{}}}.
\begin{itemize}
\item {} 
\sphinxAtStartPar
When reading a script file, COPT command\sphinxhyphen{}line double checks whether the
first non\sphinxhyphen{}blank line starts with special text: \sphinxcode{\sphinxupquote{\textquotesingle{}\#COPT script file\textquotesingle{}}}.
This is to make sure users do not load an invalid script file by mistake.
Indeed, only \sphinxcode{\sphinxupquote{\textquotesingle{}\#COPT\textquotesingle{}}} is verified. In addition, any line in scripts is
commented out if its first non\sphinxhyphen{}blank character is \sphinxcode{\sphinxupquote{\textquotesingle{}\#\textquotesingle{}}}.

\item {} 
\sphinxAtStartPar
After a script file is loaded, the tool keeps reading it as standard
inputs, until reaching end of file or a special character \sphinxcode{\sphinxupquote{\textquotesingle{}?\textquotesingle{}}}. Here,
we use question mark \sphinxcode{\sphinxupquote{\textquotesingle{}?\textquotesingle{}}} to pause scripts on purpose. To continue,
users can type \sphinxcode{\sphinxupquote{\textquotesingle{}load\textquotesingle{}}} in command line. Afterwards, the tool picks
whatever left in scripts and start to run from there, until reaching
end of file or another question mark \sphinxcode{\sphinxupquote{\textquotesingle{}?\textquotesingle{}}}. Once current scripts finish,
users can load any other script file on fly.

\end{itemize}

\sphinxAtStartPar
It also allows users to load special scripts, called \sphinxstylestrong{inline scripts}. The
only difference from a script file is that commands are separated by \sphinxcode{\sphinxupquote{\textquotesingle{};\textquotesingle{}}},
instead of \sphinxcode{\sphinxupquote{\textquotesingle{}\textbackslash{}n\textquotesingle{}}}. So inline scripts can be read by using arguments of the
option \sphinxcode{\sphinxupquote{\textquotesingle{}\sphinxhyphen{}c\textquotesingle{}}}, or loaded on fly by specifying a quoted string, and special
character \sphinxcode{\sphinxupquote{\textquotesingle{}?\textquotesingle{}}} works in the same way.

\sphinxAtStartPar
In addition, this tool provides a feature of recording \sphinxstylestrong{valid} command lines
sequentially to a script file, if users specify an output script file as
argument of the option \sphinxcode{\sphinxupquote{\textquotesingle{}\sphinxhyphen{}o\textquotesingle{}}}. Here, \sphinxstylestrong{valid} command must use known shell
commands and do not exceed number of allowed parameters.

\sphinxAtStartPar
In particular, if users load a script file or inline scripts on fly, all
commands in scripts are written to the output script file. Note that
command \sphinxcode{\sphinxupquote{\textquotesingle{}load\textquotesingle{}}} itself is not written to output script file on purpose.
Because we’ve expanded and written all commands in scripts. On the other hand,
it may trigger infinity loop if the script file loaded is actually itself.

\section{Shell commands}
\label{\detokenize{cmdline:shell-commands}}
\sphinxAtStartPar
COPT command\sphinxhyphen{}line supports the following shell commands for users to manipulate
optimization models. Moreover, shell commands are case\sphinxhyphen{}insensitive and
support tab\sphinxhyphen{}completion.

\subsection{General shell commands}
\label{\detokenize{cmdline:general-shell-commands}}
\sphinxAtStartPar
The shell commands below are in support of interactions.
\begin{itemize}
\item {} 
\sphinxAtStartPar
\sphinxcode{\sphinxupquote{cd}}:
This shell command works similar to DOS command \sphinxcode{\sphinxupquote{\textquotesingle{}cd\textquotesingle{}}}. That is,
it changes \sphinxstylestrong{current working directory}, if its argument is valid relative
or absolute path of a directory. Note that \sphinxstylestrong{current working directory} is
the base directory for relative path and tab completion. It is initialized
to current binary folder where \sphinxcode{\sphinxupquote{copt\_cmd}} exist. Users can change it by
shell command \sphinxcode{\sphinxupquote{\textquotesingle{}cd \textless{}dirpath\textgreater{}\textquotesingle{}}}. For example, if users change working
directory to a folder having mps files, reading model becomes much easier
because only filename is needed.

\item {} 
\sphinxAtStartPar
\sphinxcode{\sphinxupquote{dir/ls}}:
This shell command works similar to DOS command \sphinxcode{\sphinxupquote{\textquotesingle{}dir\textquotesingle{}}} or Bash
command \sphinxcode{\sphinxupquote{\textquotesingle{}ls\textquotesingle{}}}. That is, it lists all files and directories under given
relative or absolute path. To see files under current working directory, type
\sphinxcode{\sphinxupquote{\textquotesingle{}dir\textquotesingle{}}} or \sphinxcode{\sphinxupquote{\textquotesingle{}ls\textquotesingle{}}}; To see files under parent folder, type \sphinxcode{\sphinxupquote{\textquotesingle{}dir ../\textquotesingle{}}};
To see files under home path, type \sphinxcode{\sphinxupquote{\textquotesingle{}dir \textasciitilde{}/\textquotesingle{}}}, etc. In addition, wildcard
(\sphinxcode{\sphinxupquote{*}}) is supported as well. That is, \sphinxcode{\sphinxupquote{\textquotesingle{}dir net\textquotesingle{}}} lists all file names
starting with \sphinxcode{\sphinxupquote{\textquotesingle{}net\textquotesingle{}}} under current working directory;
\sphinxcode{\sphinxupquote{\textquotesingle{}dir /home/user/*.gz\textquotesingle{}}} lists all files of type of \sphinxcode{\sphinxupquote{\textquotesingle{}.gz\textquotesingle{}}} under path
of \sphinxcode{\sphinxupquote{\textquotesingle{}/home/user/\textquotesingle{}}}.

\item {} 
\sphinxAtStartPar
\sphinxcode{\sphinxupquote{exit/quit}}:
Leave COPT command\sphinxhyphen{}line.

\item {} 
\sphinxAtStartPar
\sphinxcode{\sphinxupquote{help}}:
It provides information on all shell commands. Typing \sphinxcode{\sphinxupquote{\textquotesingle{}help\textquotesingle{}}}
followed by a shell command name gives you more details on shell commands.
In particular, typing \sphinxcode{\sphinxupquote{\textquotesingle{}help\textquotesingle{}}} without arguments lists all shell commands with
short descriptions. Right after overview of shell commands, the text \sphinxcode{\sphinxupquote{\textquotesingle{}help\textquotesingle{}}}
with additional whitespace appear in the new prompt line. So users can
directly type, or possibly \sphinxcode{\sphinxupquote{\textless{}Tab\textgreater{}}} complete, actual shell command to read
more details.

\item {} 
\sphinxAtStartPar
\sphinxcode{\sphinxupquote{load}}:
Load a script file or inline scripts on fly and then execute a batch of commands.
The syntax of \sphinxcode{\sphinxupquote{\textquotesingle{}load\textquotesingle{}}} command should be followed by either relative/absolute
path of a script file, or quoted string of inline scripts. One special scenaio
is when current script is paused, that is, hit question mark (\sphinxcode{\sphinxupquote{\textquotesingle{}?\textquotesingle{}}})
during execution. In this case, type \sphinxcode{\sphinxupquote{\textquotesingle{}load\textquotesingle{}}} will continue the paused scripts.
If users forgot having scripts in progress and try to load another scripts, it
works as command \sphinxcode{\sphinxupquote{\textquotesingle{}load\textquotesingle{}}} and any additional argument is ignored.
This behavior is back to normal after reaching the end of current scripts.

\item {} 
\sphinxAtStartPar
\sphinxcode{\sphinxupquote{pwd}}:
This shell command works similar to Bash command \sphinxcode{\sphinxupquote{\textquotesingle{}pwd\textquotesingle{}}}. That is, display
current working directory to let users know where they are.

\item {} 
\sphinxAtStartPar
\sphinxcode{\sphinxupquote{status}}:
COPT command\sphinxhyphen{}line has a state machine on status of problem solving (see
\hyperref[\detokenize{cmdline:coptfig-cmdline-state}]{Fig.\@ \ref{\detokenize{cmdline:coptfig-cmdline-state}}}).
This is used to guide users through steps. Typing \sphinxcode{\sphinxupquote{\textquotesingle{}status\textquotesingle{}}} shows you
current interactive status. The status exposed to users are as follows:
\begin{itemize}
\item {} 
\sphinxAtStartPar
\sphinxcode{\sphinxupquote{Initial}}, initial status, either right after the tool is invoked,
or shell command \sphinxcode{\sphinxupquote{reset}} is called.

\item {} 
\sphinxAtStartPar
\sphinxcode{\sphinxupquote{Read}}, read an optimization model in format of mps successfully.

\item {} 
\sphinxAtStartPar
\sphinxcode{\sphinxupquote{SetParam}}, set value of any COPT parameter successfully.

\item {} 
\sphinxAtStartPar
\sphinxcode{\sphinxupquote{Optimize}}, shell command \sphinxcode{\sphinxupquote{opt}} is called to solve current optimization model.

\end{itemize}

\begin{figure}[H]
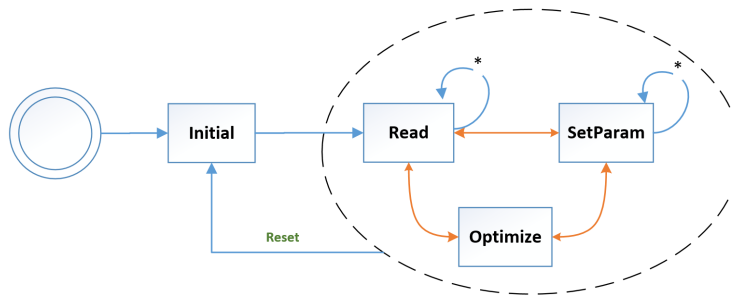

\centering
\capstart

\noindent\sphinxincludegraphics[scale=0.6]{{cmdline-state}.png}
\caption{Status of COPT command\sphinxhyphen{}line}\label{\detokenize{cmdline:coptfig-cmdline-state}}\end{figure}

\end{itemize}

\subsection{COPT shell commands}
\label{\detokenize{cmdline:copt-shell-commands}}
\sphinxAtStartPar
The shell commands below are related to COPT API calls.
\begin{itemize}
\item {} 
\sphinxAtStartPar
\sphinxcode{\sphinxupquote{display/get}}:
Display current setting of any COPT parameter or attribute. Typing \sphinxcode{\sphinxupquote{\textquotesingle{}display\textquotesingle{}}}
followed by COPT parameter or attribute name shows its current value.
Typing \sphinxcode{\sphinxupquote{\textquotesingle{}display\textquotesingle{}}} only lists all COPT parameters and attributes with
short descriptions. Right after overview of COPT parameters/attributes,
the text \sphinxcode{\sphinxupquote{\textquotesingle{}display\textquotesingle{}}} with additional whitespace appear in the new prompt line.
So users can directly type, or possibly \sphinxcode{\sphinxupquote{\textless{}Tab\textgreater{}}} complete, actual
parameter/attribute name to see its current value.

\item {} 
\sphinxAtStartPar
\sphinxcode{\sphinxupquote{opt/optimize}}:
Solve optimization model and display results on screen.
This command requires no parameter and its interactive status is set to
\sphinxcode{\sphinxupquote{\textquotesingle{}Optimize\textquotesingle{}}}.

\item {} 
\sphinxAtStartPar
\sphinxcode{\sphinxupquote{optlp/optimizelp}}:
Solve optimization model as LP model. This command requires no parameter and
its interactive status is set to \sphinxcode{\sphinxupquote{\textquotesingle{}Optimize\textquotesingle{}}}.

\item {} 
\sphinxAtStartPar
\sphinxcode{\sphinxupquote{iis}}:
Compute IIS for the loaded model.

\item {} 
\sphinxAtStartPar
\sphinxcode{\sphinxupquote{feasrelax}}:
Do feasible relaxation for infeasible problem. Note that an optimization
problem must exist before calling ‘feasrelax’. \sphinxcode{\sphinxupquote{\textquotesingle{}feasrelax\textquotesingle{}}} or \sphinxcode{\sphinxupquote{\textquotesingle{}feasrelax all\textquotesingle{}}}
means to relax all bounds of variables and constraints with penalty 1, while
\sphinxcode{\sphinxupquote{\textquotesingle{}feasrelax vars\textquotesingle{}}} means to only relax bounds of variables with penalty 1, and
\sphinxcode{\sphinxupquote{\textquotesingle{}feasrelax cons\textquotesingle{}}} means to only relax bounds of constraints with penalty 1.

\item {} 
\sphinxAtStartPar
\sphinxcode{\sphinxupquote{tune}}:
Tuning parameter automatically of the loaded model.

\item {} 
\sphinxAtStartPar
\sphinxcode{\sphinxupquote{loadtuneparam}}:
Loads the specified tuning results into the currently loaded model. The calling command is \sphinxcode{\sphinxupquote{\textquotesingle{}loadtuneparam idx\textquotesingle{}}}, where \sphinxcode{\sphinxupquote{\textquotesingle{}idx\textquotesingle{}}}
is the specified tuning result’s number. Use the command \sphinxcode{\sphinxupquote{\textquotesingle{}get TuneResults\textquotesingle{}}} to get the number of tuning results obtained after the current tuning.

\item {} 
\sphinxAtStartPar
\sphinxcode{\sphinxupquote{read}}:
Read optimization model, MIP solution, LP basis, MIP initial solution and
COPT parametes from file under given relative/absolute path. It supports
files including optimization problem file of type \sphinxcode{\sphinxupquote{\textquotesingle{}.mps\textquotesingle{}}} and
compressed format \sphinxcode{\sphinxupquote{\textquotesingle{}.mps.gz\textquotesingle{}}}, problem file of type \sphinxcode{\sphinxupquote{\textquotesingle{}.lp\textquotesingle{}}} and
compressed format \sphinxcode{\sphinxupquote{\textquotesingle{}.lp.gz\textquotesingle{}}}, problem file of type \sphinxcode{\sphinxupquote{\textquotesingle{}.dat\sphinxhyphen{}s\textquotesingle{}}} and
compressed format \sphinxcode{\sphinxupquote{\textquotesingle{}.dat\sphinxhyphen{}s.gz\textquotesingle{}}}, problem file of type \sphinxcode{\sphinxupquote{\textquotesingle{}.cbf\textquotesingle{}}} and
compressed format \sphinxcode{\sphinxupquote{\textquotesingle{}.cbf.gz\textquotesingle{}}}, problem file of type \sphinxcode{\sphinxupquote{\textquotesingle{}.bin\textquotesingle{}}} and
compressed format \sphinxcode{\sphinxupquote{\textquotesingle{}.bin.gz\textquotesingle{}}}, basis file of type \sphinxcode{\sphinxupquote{\textquotesingle{}.bas\textquotesingle{}}},
solution file of type \sphinxcode{\sphinxupquote{\textquotesingle{}.sol\textquotesingle{}}}, MIP initial solution file of type \sphinxcode{\sphinxupquote{\textquotesingle{}.mst\textquotesingle{}}},
parameter file of type \sphinxcode{\sphinxupquote{\textquotesingle{}.par\textquotesingle{}}} and tuning file of type \sphinxcode{\sphinxupquote{\textquotesingle{}.tune\textquotesingle{}}}.

\item {} 
\sphinxAtStartPar
\sphinxcode{\sphinxupquote{readmps}}:
Read optimization model in format of \sphinxcode{\sphinxupquote{\textquotesingle{}.mps\textquotesingle{}}} or \sphinxcode{\sphinxupquote{\textquotesingle{}.mps.gz\textquotesingle{}}} from file.
Note this command does not require file type \sphinxcode{\sphinxupquote{\textquotesingle{}.mps\textquotesingle{}}} or \sphinxcode{\sphinxupquote{\textquotesingle{}.mps.gz\textquotesingle{}}}.
That is, it is fine to have content in format of MPS, regardless of its file
name. In addtion, the interactive status is set to \sphinxcode{\sphinxupquote{\textquotesingle{}Read\textquotesingle{}}}.

\item {} 
\sphinxAtStartPar
\sphinxcode{\sphinxupquote{readlp}}:
Read optimization model in format of \sphinxcode{\sphinxupquote{\textquotesingle{}.lp\textquotesingle{}}} or \sphinxcode{\sphinxupquote{\textquotesingle{}.lp.gz\textquotesingle{}}} from file.
Note this command does not require file type \sphinxcode{\sphinxupquote{\textquotesingle{}.lp\textquotesingle{}}} or \sphinxcode{\sphinxupquote{\textquotesingle{}.lp.gz\textquotesingle{}}}.
That is, it is fine to have content in format of LP, regardless of its file
name. In addtion, the interactive status is set to \sphinxcode{\sphinxupquote{\textquotesingle{}Read\textquotesingle{}}}.

\item {} 
\sphinxAtStartPar
\sphinxcode{\sphinxupquote{readsdpa}}:
Read an optimization problem in format of \sphinxcode{\sphinxupquote{\textquotesingle{}.dat\sphinxhyphen{}s\textquotesingle{}}} or \sphinxcode{\sphinxupquote{\textquotesingle{}.dat\sphinxhyphen{}s.gz\textquotesingle{}}} from file.
Note that this command is similar to \sphinxcode{\sphinxupquote{\textquotesingle{}read\textquotesingle{}}} when the file type is \sphinxcode{\sphinxupquote{\textquotesingle{}.dat\sphinxhyphen{}s\textquotesingle{}}} .
However, \sphinxcode{\sphinxupquote{\textquotesingle{}readsdpa\textquotesingle{}}} does not require the file has type \sphinxcode{\sphinxupquote{\textquotesingle{}.dat\sphinxhyphen{}s\textquotesingle{}}}. It parses the
file as SDPA format, no matter what file type it is.

\item {} 
\sphinxAtStartPar
\sphinxcode{\sphinxupquote{readcbf}}:
Read an optimization problem in format of \sphinxcode{\sphinxupquote{\textquotesingle{}.cbf\textquotesingle{}}} or \sphinxcode{\sphinxupquote{\textquotesingle{}.cbf.gz\textquotesingle{}}} from a file.
Note that this command is similar to \sphinxcode{\sphinxupquote{\textquotesingle{}read\textquotesingle{}}} when the file type is \sphinxcode{\sphinxupquote{\textquotesingle{}.cbf\textquotesingle{}}} .
However, \sphinxcode{\sphinxupquote{\textquotesingle{}readcbf\textquotesingle{}}} does not require the file has type \sphinxcode{\sphinxupquote{\textquotesingle{}.cbf\textquotesingle{}}}. It parses the
file as CBF format, no matter what file type it is.

\item {} 
\sphinxAtStartPar
\sphinxcode{\sphinxupquote{readbin}}:
Read optimization model in format of \sphinxcode{\sphinxupquote{\textquotesingle{}.bin\textquotesingle{}}} or \sphinxcode{\sphinxupquote{\textquotesingle{}.bin.gz\textquotesingle{}}} from file.
Note this command does not require file type \sphinxcode{\sphinxupquote{\textquotesingle{}.bin\textquotesingle{}}} .
That is, it is fine to have content in format of COPT binary format,
regardless of its file name. In addtion, the interactive status is set to \sphinxcode{\sphinxupquote{\textquotesingle{}Read\textquotesingle{}}}.

\item {} 
\sphinxAtStartPar
\sphinxcode{\sphinxupquote{readsol}}:
Read MIP solution from file. Note this command does not require file type of
\sphinxcode{\sphinxupquote{\textquotesingle{}.sol\textquotesingle{}}}. That is, it is fine to have content of solution, regardless of its
file name.

\item {} 
\sphinxAtStartPar
\sphinxcode{\sphinxupquote{readbasis}}:
Read optimal basis from file. Note this command does not require file type of
\sphinxcode{\sphinxupquote{\textquotesingle{}.bas\textquotesingle{}}}. That is, it is fine to have content of basis, regardless of its
file name.

\item {} 
\sphinxAtStartPar
\sphinxcode{\sphinxupquote{readmst}}:
Read MIP initial solution from file. Note this command does not require file type of
\sphinxcode{\sphinxupquote{\textquotesingle{}.mst\textquotesingle{}}}. That is, it is fine to have content of solution, regardless of its
file name.

\item {} 
\sphinxAtStartPar
\sphinxcode{\sphinxupquote{readparam}}:
Read COPT parameters from file and set corresponding values. Note this command
does not require file type of \sphinxcode{\sphinxupquote{\textquotesingle{}.par\textquotesingle{}}}. That is, it is fine to have content
of COPT parameters, regardless of its file name.

\item {} 
\sphinxAtStartPar
\sphinxcode{\sphinxupquote{readtune}}:
Read the parameter file under the given relative or absolute path. Note this command
does not require file type of \sphinxcode{\sphinxupquote{\textquotesingle{}.tune\textquotesingle{}}} . That is, it is fine that the file conforms
to the COPT parameter tuning file format, regardless of its file name.

\item {} 
\sphinxAtStartPar
\sphinxcode{\sphinxupquote{reset}}:
Reset current optimization model and all parameters/attributes to defaults.
In addtion, the interactive status is set to \sphinxcode{\sphinxupquote{\textquotesingle{}Initial\textquotesingle{}}}.

\item {} 
\sphinxAtStartPar
\sphinxcode{\sphinxupquote{set}}:
Set value of any COPT parameter. The syntax of this command should
be \sphinxcode{\sphinxupquote{\textquotesingle{}set\textquotesingle{}}}, followed by COPT parameter name and an integer/double
number. Moreover, Typing \sphinxcode{\sphinxupquote{\textquotesingle{}set\textquotesingle{}}} only lists all COPT parameters
with short descriptions. Right after overview of COPT parameters,
the text \sphinxcode{\sphinxupquote{\textquotesingle{}set\textquotesingle{}}} with additional whitespace appears
in the new prompt line. So users can directly type, or possibly \sphinxcode{\sphinxupquote{\textless{}Tab\textgreater{}}}
complete, actual parameter name to complete as partial command. If so, its
current value, default value, min value, max value of given COPT parameter
are displayed on screen. Now users know how to add correct value to complete
the full command of \sphinxcode{\sphinxupquote{\textquotesingle{}set\textquotesingle{}}}. One tip of saving typing here is to get last
history entry by pressing \sphinxcode{\sphinxupquote{\textless{}Up\textgreater{}}}. At last, the interactive
status is set to \sphinxcode{\sphinxupquote{\textquotesingle{}SetParam\textquotesingle{}}}.

\item {} 
\sphinxAtStartPar
\sphinxcode{\sphinxupquote{set LogFile}}:

\sphinxAtStartPar
Set the log file for recording solver’s logging output, e.g., \sphinxcode{\sphinxupquote{\textquotesingle{}set LogFile filename.log\textquotesingle{}}}.

\sphinxAtStartPar
This must be specified before the solve command (\sphinxcode{\sphinxupquote{opt}}) is executed.
COPT will write log to the specified file during the solving process.

\item {} 
\sphinxAtStartPar
\sphinxcode{\sphinxupquote{set CpuAffinity}}:

\sphinxAtStartPar
Set the CPU affinity for the current solver process, e.g., \sphinxcode{\sphinxupquote{\textquotesingle{}set CpuAffinity 0f\textquotesingle{}}}.

\item {} 
\sphinxAtStartPar
\sphinxcode{\sphinxupquote{get CpuAffinity}}:

\sphinxAtStartPar
Get the CPU affinity configuration of the current solver process.

\item {} 
\sphinxAtStartPar
\sphinxcode{\sphinxupquote{get NumaNodeCnt}}:

\sphinxAtStartPar
Get the number of NUMA nodes available on the system.

\item {} 
\sphinxAtStartPar
\sphinxcode{\sphinxupquote{set CpuBind}}:

\sphinxAtStartPar
Bind the current solver process to the specified NUMA node, e.g., \sphinxcode{\sphinxupquote{\textquotesingle{}set CpuBind 0\textquotesingle{}}}.

\item {} 
\sphinxAtStartPar
\sphinxcode{\sphinxupquote{set MemBind}}:

\sphinxAtStartPar
Bind the memory allocations of the solver process to the specified NUMA node, e.g., \sphinxcode{\sphinxupquote{\textquotesingle{}set MemBind 0\textquotesingle{}}}.

\item {} 
\sphinxAtStartPar
\sphinxcode{\sphinxupquote{write}}:
Output MPS/LP/CBF format model, COPT binary format model, IIS model, FeasRelax model,
LP/MIP solution, optimal basis, settings of modified COPT parameters to file under
given relative/absolution path. This command detects file types. For an instance,
\sphinxcode{\sphinxupquote{\textquotesingle{}write diet.sol\textquotesingle{}}} outputs LP solution to file \sphinxcode{\sphinxupquote{\textquotesingle{}diet.sol\textquotesingle{}}}. An error message
will be shown to users if file type do not match. Supported file types are:
\sphinxcode{\sphinxupquote{\textquotesingle{}.mps\textquotesingle{}}}, \sphinxcode{\sphinxupquote{\textquotesingle{}.lp\textquotesingle{}}}, \sphinxcode{\sphinxupquote{\textquotesingle{}.bin\textquotesingle{}}}, \sphinxcode{\sphinxupquote{\textquotesingle{}cbf\textquotesingle{}}}, \sphinxcode{\sphinxupquote{\textquotesingle{}.iis\textquotesingle{}}}, \sphinxcode{\sphinxupquote{\textquotesingle{}.relax\textquotesingle{}}}, \sphinxcode{\sphinxupquote{\textquotesingle{}.sol\textquotesingle{}}},
\sphinxcode{\sphinxupquote{\textquotesingle{}.bas\textquotesingle{}}}, \sphinxcode{\sphinxupquote{\textquotesingle{}.mst\textquotesingle{}}} and \sphinxcode{\sphinxupquote{\textquotesingle{}.par\textquotesingle{}}}.

\item {} 
\sphinxAtStartPar
\sphinxcode{\sphinxupquote{writemps}}:
Output current optimization model to a file of type \sphinxcode{\sphinxupquote{\textquotesingle{}.mps\textquotesingle{}}}. Note
\sphinxcode{\sphinxupquote{\textquotesingle{}.mps\textquotesingle{}}} is appended to the file name, if users do not add it.

\item {} 
\sphinxAtStartPar
\sphinxcode{\sphinxupquote{writelp}}:
Output current optimization model to a file of type \sphinxcode{\sphinxupquote{\textquotesingle{}.lp\textquotesingle{}}}. Note
\sphinxcode{\sphinxupquote{\textquotesingle{}.lp\textquotesingle{}}} is appended to the file name, if users do not add it.

\item {} 
\sphinxAtStartPar
\sphinxcode{\sphinxupquote{writecbf}}:
Output problem to a file of type \sphinxcode{\sphinxupquote{\textquotesingle{}.cbf\textquotesingle{}}}. Note
\sphinxcode{\sphinxupquote{\textquotesingle{}.cbf\textquotesingle{}}} is appended to the file name, if users do not add it.

\item {} 
\sphinxAtStartPar
\sphinxcode{\sphinxupquote{writebin}}:
Output current optimization model to a file of type \sphinxcode{\sphinxupquote{\textquotesingle{}.bin\textquotesingle{}}}. Note
\sphinxcode{\sphinxupquote{\textquotesingle{}.bin\textquotesingle{}}} is appended to the file name, if users do not add it.

\item {} 
\sphinxAtStartPar
\sphinxcode{\sphinxupquote{writeiis}}:
Output computed IIS model to a file of type \sphinxcode{\sphinxupquote{\textquotesingle{}.iis\textquotesingle{}}}. Note
\sphinxcode{\sphinxupquote{\textquotesingle{}.iis\textquotesingle{}}} is appended to the file name, if users do not add it.

\item {} 
\sphinxAtStartPar
\sphinxcode{\sphinxupquote{writerelax}}:
Output feasibility relaxation problem to a file of type \sphinxcode{\sphinxupquote{\textquotesingle{}.relax\textquotesingle{}}}.
Note \sphinxcode{\sphinxupquote{\textquotesingle{}.relax\textquotesingle{}}} is appended to the file name, if users do not add it.

\item {} 
\sphinxAtStartPar
\sphinxcode{\sphinxupquote{writesol}}:
Output LP solution of problem to a file of type \sphinxcode{\sphinxupquote{\textquotesingle{}.sol\textquotesingle{}}}. Note \sphinxcode{\sphinxupquote{\textquotesingle{}.sol\textquotesingle{}}} is
appended to the file name, if users do not add it.

\item {} 
\sphinxAtStartPar
\sphinxcode{\sphinxupquote{writepoolsol}}:
For the current problem, save the solution with the specified number in the solution pool
to the result file under the given relative or absolute path. If the file extension is not\textasciigrave{}\textasciigrave{}’.sol’\textasciigrave{}\textasciigrave{} , the suffix \sphinxcode{\sphinxupquote{\textquotesingle{}.sol\textquotesingle{}}} is automatically added.
The calling command is \sphinxcode{\sphinxupquote{\textquotesingle{}writepoolsol idx pool\_idx.sol\textquotesingle{}}},
where \sphinxcode{\sphinxupquote{\textquotesingle{}idx\textquotesingle{}}} is the specified solution number in solution pool.
Use the command \sphinxcode{\sphinxupquote{\textquotesingle{}get PoolSols\textquotesingle{}}} to get the current number of solution pool results.

\item {} 
\sphinxAtStartPar
\sphinxcode{\sphinxupquote{writebasis}}:
Output optimal basis to a file of type \sphinxcode{\sphinxupquote{\textquotesingle{}.bas\textquotesingle{}}}. Note \sphinxcode{\sphinxupquote{\textquotesingle{}.bas\textquotesingle{}}}
is appended to the file name, if users do not add it.

\item {} 
\sphinxAtStartPar
\sphinxcode{\sphinxupquote{writemst}}:
Output current best MIP solution to a file of type \sphinxcode{\sphinxupquote{\textquotesingle{}.mst\textquotesingle{}}}. Note \sphinxcode{\sphinxupquote{\textquotesingle{}.mst\textquotesingle{}}}
is appended to the file name, if users do not add it.

\item {} 
\sphinxAtStartPar
\sphinxcode{\sphinxupquote{writeparam}}:
Output modified COPT parameters to a file of type \sphinxcode{\sphinxupquote{\textquotesingle{}.par\textquotesingle{}}}.
Note \sphinxcode{\sphinxupquote{\textquotesingle{}.par\textquotesingle{}}} is appended to the file name, if users do not add it.

\item {} 
\sphinxAtStartPar
\sphinxcode{\sphinxupquote{writetuneparam}}:
Save the specified tuning result to the parameter setting file under the given relative or absolute path.
The calling command is \sphinxcode{\sphinxupquote{\textquotesingle{}writetuneparam idx tune\_idx.par\textquotesingle{}}},
where \sphinxcode{\sphinxupquote{\textquotesingle{}idx\textquotesingle{}}} is the specified tuning result number.
Use the command \sphinxcode{\sphinxupquote{\textquotesingle{}get TuneResults\textquotesingle{}}} to get the number of tuning results obtained after the current tuning.

\end{itemize}

\section{Example Usage}
\label{\detokenize{cmdline:example-usage}}
\sphinxAtStartPar
This section shows how to use COPT command\sphinxhyphen{}line to interactively solve a
well\sphinxhyphen{}known problem, called “Diet Problem”. Please refer to
{\hyperref[\detokenize{amplinterface:examplediet}]{\sphinxcrossref{\DUrole{std,std-ref}{AMPL Interface \sphinxhyphen{} Example Usage}}}} for problem description
in more detail.

\sphinxAtStartPar
With valid license, COPT command\sphinxhyphen{}line should run as follows, after entering
\sphinxcode{\sphinxupquote{copt\_cmd}} in command prompt.

\begin{sphinxVerbatim}[commandchars=\\\{\}]
copt\PYGZus{}cmd
\end{sphinxVerbatim}

\sphinxAtStartPar
If users are new to COPT command\sphinxhyphen{}line, always start with shell command
\sphinxcode{\sphinxupquote{\textquotesingle{}help\textquotesingle{}}}.
\begin{sphinxalltt}
Cardinal Optimizer v8.0.1. Build date Mar 04 2024
Copyright Cardinal Operations 2025. All Rights Reserved

COPT\textgreater{}
\end{sphinxalltt}

\sphinxAtStartPar
Suppose diet model has mps format and exits in the current working directory.
In this way, we just type its file name to read,
without worrying about its path.

\begin{sphinxVerbatim}[commandchars=\\\{\}]
COPT\PYGZgt{} read diet.mps
Reading from \PYGZsq{}/home/username/copt/diet.mps\PYGZsq{}
Reading finished (0.00s)
\end{sphinxVerbatim}

\sphinxAtStartPar
Before solving it, users are free to tune any COPT parameter. Below
is an example to set value of double parameter \sphinxcode{\sphinxupquote{TimeLimit}} to 10s.

\sphinxAtStartPar
If users are not familiar with COPT parameters, just type \sphinxcode{\sphinxupquote{\textquotesingle{}set\textquotesingle{}}} to list
all public COPT paramters and attributes with short description. Furthermore,
\sphinxcode{\sphinxupquote{\textquotesingle{}set\textquotesingle{}}} with parameter name, for example \sphinxcode{\sphinxupquote{\textquotesingle{}set TimeLimit\textquotesingle{}}}, displays its
current value, default value, min value and max value of the given parameter.

\begin{sphinxVerbatim}[commandchars=\\\{\}]
COPT\PYGZgt{} set timelimit 10

Setting parameter \PYGZsq{}TimeLimit\PYGZsq{} to 10
\end{sphinxVerbatim}

\sphinxAtStartPar
After tuning parameters, it is time to solve the model. The messages during
problem solving are shown as follows.

\begin{sphinxVerbatim}[commandchars=\\\{\}]
COPT\PYGZgt{} opt
Model fingerprint: 129c032d

Hardware has 4 cores and 8 threads. Using instruction set X86\PYGZus{}NATIVE (1)
Minimizing an LP problem

The original problem has:
    4 rows, 8 columns and 31 non\PYGZhy{}zero elements
The presolved problem has:
    4 rows, 8 columns and 31 non\PYGZhy{}zero elements

Starting the simplex solver using up to 8 threads

Method   Iteration           Objective  Primal.NInf   Dual.NInf        Time
Dual             0    0.0000000000e+00            4           0       0.01s
Dual             1    8.8201991646e+01            0           0       0.01s
Postsolving
Dual             1    8.8200000000e+01            0           0       0.01s

Solving finished
Status: Optimal  Objective: 8.8200000000e+01  Iterations: 1  Time: 0.01s
\end{sphinxVerbatim}

\sphinxAtStartPar
After solving the model, users might check results by using shell command
\sphinxcode{\sphinxupquote{\textquotesingle{}get\textquotesingle{}}} with a parameter name. Note that, similar to \sphinxcode{\sphinxupquote{\textquotesingle{}set\textquotesingle{}}}, type
\sphinxcode{\sphinxupquote{\textquotesingle{}get\textquotesingle{}}} to list all public parameters and attributes. In particular,
\sphinxcode{\sphinxupquote{\textquotesingle{}get all\textquotesingle{}}} shows all parameters/attributes and their current value.

\begin{sphinxVerbatim}[commandchars=\\\{\}]
COPT\PYGZgt{} get TimeLimit
  DblParam: [Current]  10s
COPT\PYGZgt{} get LpObjval
  DblAttr:  [Current] 88.2
COPT\PYGZgt{} get LpStatus
  IntAttr:  [Current] 1 optimal
\end{sphinxVerbatim}

\sphinxAtStartPar
Before leaving COPT command\sphinxhyphen{}line, users might output the model in format of
mps, optimal basis, modified parameters, or LP solution to files. Below is an
example to write LP solution to current directory.

\begin{sphinxVerbatim}[commandchars=\\\{\}]
COPT\PYGZgt{} writesol diet
  Writing solutions to /home/username/copt/diet.sol
COPT\PYGZgt{} quit
  Leaving COPT...
\end{sphinxVerbatim}

\sphinxAtStartPar
Below is the script file \sphinxcode{\sphinxupquote{diet.in}} by putting everything together,
see \hyperref[\detokenize{cmdline:coptcode-cmddiet-bat}]{Listing \ref{\detokenize{cmdline:coptcode-cmddiet-bat}}}:
\sphinxSetupCaptionForVerbatim{\sphinxcode{\sphinxupquote{diet.in}}}
\def\sphinxLiteralBlockLabel{\label{\detokenize{cmdline:coptcode-cmddiet-bat}}}
\begin{sphinxVerbatim}[commandchars=\\\{\},numbers=left,firstnumber=1,stepnumber=1]
\PYGZsh{}COPT script\PYGZhy{}in file

read diet.mps
set timelimit 10
opt
writesol diet
quit
\end{sphinxVerbatim}

\sphinxAtStartPar
which is loaded by using the option \sphinxcode{\sphinxupquote{\textquotesingle{}\sphinxhyphen{}i\textquotesingle{}}} when starting \sphinxcode{\sphinxupquote{copt\_cmd}}:

\begin{sphinxVerbatim}[commandchars=\\\{\}]
copt\PYGZus{}cmd \PYGZhy{}i diet.in
\end{sphinxVerbatim}

\sphinxAtStartPar
or executing shell command \sphinxcode{\sphinxupquote{load}} on fly.

\begin{sphinxVerbatim}[commandchars=\\\{\}]
COPT\PYGZgt{} load diet.in
\end{sphinxVerbatim}

\sphinxstepscope

\chapter{COPT Floating Licensing service}
\label{\detokenize{floating:copt-floating-licensing-service}}\label{\detokenize{floating:chapfloating}}\label{\detokenize{floating::doc}}
\sphinxAtStartPar
The \sphinxstylestrong{Cardinal Optimizer} provides COPT Floating Token Server on all supported
platforms, who serve license tokens to COPT client applications over local network.

\sphinxAtStartPar
Once you have floating license properly installed, server owns a set of license tokens
up to number described in the license file. Any properly configured COPT client of
the same version can request a token from server and release it whenever quit.

\section{Server Setup}
\label{\detokenize{floating:server-setup}}
\sphinxAtStartPar
The application of COPT Floating Token server includes \sphinxcode{\sphinxupquote{copt\_flserver}} executable
and a configuration file \sphinxcode{\sphinxupquote{fls.ini}}. The very first thing to do when server starts
is to verify floating license locally, whose location is specified in \sphinxcode{\sphinxupquote{fls.ini}}.
If local validation passes, server connects to remote COPT licensing server for further
validation, including machine IP, which is supposed to match IP range that user provided
during registration. This means the machine running COPT Floating Token Server
should have internet access in specified area.
For details, please see descriptions below or refer to
{\hyperref[\detokenize{install:parcoptgetlic}]{\sphinxcrossref{\DUrole{std,std-ref}{How to obtain and setup license}}}}.

\subsection{Installation}
\label{\detokenize{floating:installation}}
\sphinxAtStartPar
The \sphinxstylestrong{Cardinal Optimizer} provides a separate package for remote services,
which include COPT floating token server. Users may apply for remote package from
customer service. Afterwards, unzip the remote package and move to any
folder on your computer. The software is portable and does not change
anything in the system it runs on. Below are details of installation.

\sphinxAtStartPar
\sphinxstylestrong{Windows}

\sphinxAtStartPar
Please unzip the remote package and move to any folder. Though, it is common
to move to folder under \sphinxcode{\sphinxupquote{C:\textbackslash{}Program Files}}.

\sphinxAtStartPar
\sphinxstylestrong{Linux}

\sphinxAtStartPar
To unzip the remote package, enter the following command in terminal:
\begin{sphinxalltt}
tar \sphinxhyphen{}xzf CardinalOptimizer\sphinxhyphen{}Remote\sphinxhyphen{}8.0.1\sphinxhyphen{}lnx64.tar.gz
\end{sphinxalltt}

\sphinxAtStartPar
Then, the following command moves folder copt\_remote80 in current
directory to other path. For an example, admin user may move it to folder
under \sphinxcode{\sphinxupquote{/opt}} and standard user may move it to \sphinxcode{\sphinxupquote{\$HOME}}.
\begin{sphinxalltt}
sudo mv copt\_remote80 /opt
\end{sphinxalltt}

\sphinxAtStartPar
Note that it requires \sphinxcode{\sphinxupquote{root}} privilege to execute this command.

\sphinxAtStartPar
\sphinxstylestrong{MacOS}

\sphinxAtStartPar
To unzip the remote package, enter the following command in terminal:
\begin{sphinxalltt}
tar \sphinxhyphen{}xzf CardinalOptimizer\sphinxhyphen{}Remote\sphinxhyphen{}8.0.1\sphinxhyphen{}universal\_mac.tar.gz
\end{sphinxalltt}

\sphinxAtStartPar
Then, the following command moves folder copt\_remote80 in current
directory to other path. For an example, admin user may move it to folder
under \sphinxcode{\sphinxupquote{/Applications}} and standard user may move it to \sphinxcode{\sphinxupquote{\$HOME}}.
\begin{sphinxalltt}
mv copt\_remote80 /Applications
\end{sphinxalltt}

\subsection{Floating License}
\label{\detokenize{floating:floating-license}}
\sphinxAtStartPar
After installing COPT remote package, it requirs floating license to run.
It is prefered to save floating license files, \sphinxcode{\sphinxupquote{license.dat}} and
\sphinxcode{\sphinxupquote{license.key}}, to \sphinxcode{\sphinxupquote{floating}} folder in path of remote package.

\sphinxAtStartPar
The following explains how to obtain the license file
via the \sphinxcode{\sphinxupquote{copt\_licgen}} tool and the license credential information \sphinxcode{\sphinxupquote{key}} under different systems.

\begin{sphinxadmonition}{note}{Note}

\sphinxAtStartPar
If the user has already obtained the two license files \sphinxcode{\sphinxupquote{license.dat}} and \sphinxcode{\sphinxupquote{license.key}}, there is no need to obtain them again.
You can skip the following steps to obtain the license file and refer to {\hyperref[\detokenize{floating:chapcluster-config}]{\sphinxcrossref{\DUrole{std,std-ref}{Configuration}}}} directly.
\end{sphinxadmonition}

\sphinxAtStartPar
\sphinxstylestrong{Windows}

\sphinxAtStartPar
If the COPT remote package is installed under \sphinxcode{\sphinxupquote{"C:\textbackslash{}Program Files"}},
execute the following command to enter \sphinxcode{\sphinxupquote{floating}} folder in path of
remote package.
\begin{sphinxalltt}
cd "C:\textbackslash{}Program Files\textbackslash{}copt\_remote80\textbackslash{}floating"
\end{sphinxalltt}

\sphinxAtStartPar
Note that the tool  \sphinxcode{\sphinxupquote{copt\_licgen}} creating license files exists
under \sphinxcode{\sphinxupquote{tools}} folder in path of remote package. The following command
creates floating license files in current directory, given a floating
license key, such as \sphinxcode{\sphinxupquote{7483dff0863ffdae9fff697d3573e8bc}} .

\begin{sphinxVerbatim}[commandchars=\\\{\}]
..\PYGZbs{}tools\PYGZbs{}copt\PYGZus{}licgen \PYGZhy{}key 7483dff0863ffdae9fff697d3573e8bc
\end{sphinxVerbatim}

\sphinxAtStartPar
\sphinxstylestrong{Linux and MacOS}

\sphinxAtStartPar
If the COPT remote package is installed under \sphinxcode{\sphinxupquote{"/Applications"}},
execute the following command to enter \sphinxcode{\sphinxupquote{floating}} folder in path of
remote package on MacOS system.
\begin{sphinxalltt}
cd /Applications/copt\_remote80/floating
\end{sphinxalltt}

\sphinxAtStartPar
The following command creates floating license files in current directory,
given a floating license key, such as \sphinxcode{\sphinxupquote{7483dff0863ffdae9fff697d3573e8bc}} .

\begin{sphinxVerbatim}[commandchars=\\\{\}]
../tools/copt\PYGZus{}licgen \PYGZhy{}key 7483dff0863ffdae9fff697d3573e8bc
\end{sphinxVerbatim}

\sphinxAtStartPar
In addition, if users run the above command when current directory is
different than \sphinxcode{\sphinxupquote{floating}} folder in path of remote package, it is prefered
to move them to \sphinxcode{\sphinxupquote{floating}}. The following command does so.
\begin{sphinxalltt}
mv license.* /Application/copt\_remote80/floating
\end{sphinxalltt}

\subsection{Configuration}
\label{\detokenize{floating:configuration}}\label{\detokenize{floating:chapcluster-config}}
\sphinxAtStartPar
Below is a typical configuration file, \sphinxcode{\sphinxupquote{fls.ini}}, of COPT Floating Token Server.

\begin{sphinxVerbatim}[commandchars=\\\{\}]
\PYG{o}{[}Main\PYG{o}{]}
\PYG{n+nv}{Port}\PYG{+w}{ }\PYG{o}{=}\PYG{+w}{ }\PYG{l+m}{7979}

\PYG{c+c1}{\PYGZsh{} password is case\PYGZhy{}sensitive and default is empty}
\PYG{c+c1}{\PYGZsh{} it applies to both copt clients and managing tool}
\PYG{n+nv}{PassWd}\PYG{+w}{ }\PYG{o}{=}

\PYG{o}{[}SSL\PYG{o}{]}
\PYG{c+c1}{\PYGZsh{} needed if connecting using SSL}
\PYG{n+nv}{CaFile}\PYG{+w}{ }\PYG{o}{=}
\PYG{n+nv}{CertFile}\PYG{+w}{ }\PYG{o}{=}
\PYG{n+nv}{CertkeyFile}\PYG{+w}{ }\PYG{o}{=}

\PYG{o}{[}Licensing\PYG{o}{]}
\PYG{c+c1}{\PYGZsh{} if empty or default license name, it is from binary folder}
\PYG{c+c1}{\PYGZsh{} to get license files from cwd, add prefix \PYGZdq{}./\PYGZdq{}}
\PYG{c+c1}{\PYGZsh{} full path is supported as well}
\PYG{n+nv}{LicenseFile}\PYG{+w}{ }\PYG{o}{=}\PYG{+w}{ }license.dat
\PYG{n+nv}{PubkeyFile}\PYG{+w}{ }\PYG{o}{=}\PYG{+w}{ }license.key

\PYG{o}{[}WLS\PYG{o}{]}
\PYG{c+c1}{\PYGZsh{} WebServer have a default host and no need to edit in most scenarios}
\PYG{c+c1}{\PYGZsh{} Must specify WebLicenseId and WebAccesskey to trigger web licensing}
\PYG{n+nv}{WebServer}\PYG{+w}{ }\PYG{o}{=}
\PYG{n+nv}{WebLicenseId}\PYG{+w}{ }\PYG{o}{=}
\PYG{n+nv}{WebAccessKey}\PYG{+w}{ }\PYG{o}{=}
\PYG{n+nv}{WebTokenDuration}\PYG{+w}{ }\PYG{o}{=}\PYG{+w}{ }\PYG{l+m}{300}

\PYG{o}{[}Filter\PYG{o}{]}
\PYG{c+c1}{\PYGZsh{} default policy 0 indicates accepting all connections, except for ones in blacklist}
\PYG{c+c1}{\PYGZsh{} otherwise, denying all connections except for ones in whitelist}
\PYG{n+nv}{DefaultPolicy}\PYG{+w}{ }\PYG{o}{=}\PYG{+w}{ }\PYG{l+m}{0}
\PYG{n+nv}{UseBlackList}\PYG{+w}{ }\PYG{o}{=}\PYG{+w}{ }\PYG{n+nb}{true}
\PYG{n+nv}{UseWhiteList}\PYG{+w}{ }\PYG{o}{=}\PYG{+w}{ }\PYG{n+nb}{true}
\PYG{n+nv}{FilterListFile}\PYG{+w}{ }\PYG{o}{=}\PYG{+w}{ }flsfilters.ini

\PYG{o}{[}Logs\PYG{o}{]}
\PYG{n+nv}{LogsFolder}\PYG{+w}{ }\PYG{o}{=}\PYG{+w}{ }./logs
\end{sphinxVerbatim}

\sphinxAtStartPar
The \sphinxcode{\sphinxupquote{Main}} section specifies port number, through which COPT clients connect to server
and then obtain the license token. The \sphinxcode{\sphinxupquote{Licensing}} section specifies location of
floating license. As described in comments above, if emtpy string or default license file name
is specified, floating license files are read from the binary folder where
the server executable reside.

\sphinxAtStartPar
It is possible to run COPT Floating Token server, even if
floating license files do not exist in the same folder as the server executive.
One solution is to set \sphinxcode{\sphinxupquote{LicenseFile = ./license.dat}} and \sphinxcode{\sphinxupquote{PubkeyFile = ./license.key}}.
By doing so, server read floating license from the current working directory.
That is, user could execute server application at the path where floating license files exist.

\sphinxAtStartPar
The other solution is to set full path of license files in configuration.
As mentioned before, Cardinal Optimizer allows users to set
environment variable \sphinxcode{\sphinxupquote{COPT\_LICENSE\_DIR}} for license files. For details,
please refer to {\hyperref[\detokenize{install:chapinstall}]{\sphinxcrossref{\DUrole{std,std-ref}{How to install Cardinal Optimizer}}}}.
If user prefers the way of environment variable, the configuration file should have
the full path to floating license.

\sphinxAtStartPar
In the \sphinxcode{\sphinxupquote{Filter}} section, \sphinxcode{\sphinxupquote{DefaultPolicy}} has default value 0, meaning
all connections are accepted except for those in black lists; if it is set
to non\sphinxhyphen{}zero value, then all connection are blocked except for those in white
lists. In addition, black list is enabled if \sphinxcode{\sphinxupquote{UseBlackList}} is true and white list
is enabled if \sphinxcode{\sphinxupquote{UseWhiteList}} is true. The filter configuration file is specified by
\sphinxcode{\sphinxupquote{FilterListFile}}. Below is an example of the filter configuration file.

\begin{sphinxVerbatim}[commandchars=\\\{\}]
\PYG{o}{[}BlackList\PYG{o}{]}
\PYG{c+c1}{\PYGZsh{} 127.0.*.* + user@machine*}

\PYG{o}{[}WhiteList\PYG{o}{]}
\PYG{c+c1}{\PYGZsh{} 127.0.1.2/16 \PYGZhy{} user@machine*}

\PYG{o}{[}ToolList\PYG{o}{]}
\PYG{c+c1}{\PYGZsh{} only tool client at server side can access by default}
\PYG{l+m}{127}.0.0.1/32
\end{sphinxVerbatim}

\sphinxAtStartPar
It has three sections and each section has its own rules. In section of \sphinxcode{\sphinxupquote{BlackList}},
one may add rules to block others from connection. In section of \sphinxcode{\sphinxupquote{WhiteList}},
one may add rules to grant others for connection, even if the default policy is to block
all connections. Only users listed in section of \sphinxcode{\sphinxupquote{ToolList}} are able to connect to
floating token server by Floating Token Server Managing Tool (see below for details).

\sphinxAtStartPar
Specifically, rules in filter configuration have format of starting with IP address.
To specify IP range, you may include wildcard (*) in IP address, or use CIDR notation,
that is, a IPv4 address and its associated network prefix.
In addtion, a rule may include (+) or exclude (\sphinxhyphen{}) given user at given machine, such as
\sphinxcode{\sphinxupquote{127.0.1.2/16 \sphinxhyphen{} user@machine}}. Here, \sphinxcode{\sphinxupquote{user}} refers to \sphinxcode{\sphinxupquote{username}}, which can be queried
by \sphinxcode{\sphinxupquote{whoami}} on Linux/MacOS platform; \sphinxcode{\sphinxupquote{machine}} refers to \sphinxcode{\sphinxupquote{computer name}}, which can
be queried by \sphinxcode{\sphinxupquote{hostname}} on Linux/MacOS platform.

\subsection{Web License for Floating Server}
\label{\detokenize{floating:web-license-for-floating-server}}
\sphinxAtStartPar
Besides local floating license above, users may use web license for floating
server to run floating service. This requires that the machine running
floating server must have internet access. However, hardware info are not
required any more. That is, users are free to deploy floating server to
any cloud machine or container, as long as they have internet access.
Please refer to \sphinxhref{https://copt.shanshu.ai/license}{COPT web license page} for details.

\sphinxAtStartPar
Below are brief steps:
\begin{itemize}
\item {} 
\sphinxAtStartPar
Follow steps to register an account and apply for trial of
web license for floating server.

\item {} 
\sphinxAtStartPar
Once approved, \sphinxtitleref{Web License ID} is generated for users

\item {} 
\sphinxAtStartPar
On page of \sphinxtitleref{API Keys}, create \sphinxtitleref{Web Access Key} using
given \sphinxtitleref{Web License ID}

\end{itemize}

\sphinxAtStartPar
Afterwards, users edit configuration file \sphinxtitleref{fls.ini} and add values of both
\sphinxtitleref{Web License ID} and \sphinxtitleref{Web Access Key} to related keywords in section of
\sphinxtitleref{WLS}. For instance,

\begin{sphinxVerbatim}[commandchars=\\\{\}]
\PYG{o}{[}WLS\PYG{o}{]}
\PYG{c+c1}{\PYGZsh{} WebServer have a default host and no need to edit in most scenarios}
\PYG{c+c1}{\PYGZsh{} Must specify WebLicenseId and WebAccesskey to trigger web licensing}
\PYG{n+nv}{WebServer}\PYG{+w}{ }\PYG{o}{=}
\PYG{n+nv}{WebLicenseId}\PYG{+w}{ }\PYG{o}{=}
\PYG{n+nv}{WebAccessKey}\PYG{+w}{ }\PYG{o}{=}
\PYG{n+nv}{WebTokenDuration}\PYG{+w}{ }\PYG{o}{=}\PYG{+w}{ }\PYG{l+m}{300}
\end{sphinxVerbatim}

\sphinxAtStartPar
As of now, floating server talks to \sphinxhref{https://copt.shanshu.ai/license}{COPT web license page}
for licensing. Users are able to monitor its token usage and other
informations online.

\subsection{Example Usage}
\label{\detokenize{floating:example-usage}}
\sphinxAtStartPar
Suppose that floating license exists in the same folder where the server
executable reside. To start the COPT Floating Token Server, just execute
the following command at any directory in Windows console, or Linux/Mac
terminal.

\begin{sphinxVerbatim}[commandchars=\\\{\}]
copt\PYGZus{}flserver
\end{sphinxVerbatim}

\sphinxAtStartPar
If you see log information as follows, the Floating Token Server has been
successfully started. Server monitors any connction from COPT clients,
manages approved clients as well as requests in queue. User can stop
Floating Token Server anytime when entering \sphinxcode{\sphinxupquote{q}} or \sphinxcode{\sphinxupquote{Q}}.
\begin{sphinxalltt}
\textgreater{} copt\_flserver
  {[} Info{]} Floating Token Server, COPT v8.0.1 20240304
  {[} Info{]} server started at port 7979
\end{sphinxalltt}

\sphinxAtStartPar
If failed to verify local floating license, or something is wrong on remote
COPT license server, you might see error logs as follows.
\begin{sphinxalltt}
\textgreater{} copt\_flserver
  {[} Info{]} Floating Token Server, COPT v8.0.1 20240304
  {[}Error{]} Invalid signature in public key file
  {[}Error{]} Fail to verify local license
\end{sphinxalltt}

\sphinxAtStartPar
and
\begin{sphinxalltt}
\textgreater{} copt\_flserver
  {[} Info{]} Floating Token Server, COPT v8.0.1 20240304
  {[}Error{]} Error to connect license server
  {[}Error{]} Fail to verify floating license by server
\end{sphinxalltt}

\section{Client Setup}
\label{\detokenize{floating:client-setup}}
\sphinxAtStartPar
COPT Clients can be COPT command\sphinxhyphen{}line tool, or any application which
solve problems using COPT api, such as COPT python interface. Floating licensing
is a better approach in terms of flexibility and efficiency.
Different from stand\sphinxhyphen{}alone license, any machine having properly configured COPT client
can legally run Cardingal Optimizer, as long as peak number of connected clients
does not exceed the token number.

\subsection{Configuration}
\label{\detokenize{floating:id1}}

\subsubsection{Via the configuration file}
\label{\detokenize{floating:via-the-configuration-file}}
\sphinxAtStartPar
Before running COPT as floating client, please make sure that you have
installed COPT locally. For details, please refer to
{\hyperref[\detokenize{install:chapinstall}]{\sphinxcrossref{\DUrole{std,std-ref}{How to install Cardinal Optimizer}}}}.
Users can skip obtaining local licenses by adding a floating configuration
file \sphinxcode{\sphinxupquote{client.ini}}.

\sphinxAtStartPar
Below is a typical configuration file, \sphinxcode{\sphinxupquote{client.ini}}, of COPT floating clients.

\begin{sphinxVerbatim}[commandchars=\\\{\}]
\PYG{n+nv}{Host}\PYG{+w}{ }\PYG{o}{=}\PYG{+w}{ }\PYG{l+m}{192}.168.1.11
\PYG{n+nv}{Port}\PYG{+w}{ }\PYG{o}{=}\PYG{+w}{ }\PYG{l+m}{7979}
\PYG{n+nv}{QueueTime}\PYG{+w}{ }\PYG{o}{=}\PYG{+w}{ }\PYG{l+m}{600}
\end{sphinxVerbatim}

\sphinxAtStartPar
As configured above, COPT floating client tries to connect to \sphinxcode{\sphinxupquote{192.168.1.11}} at port \sphinxcode{\sphinxupquote{7979}}
with wait time in queue up to 600 seconds. Here, \sphinxcode{\sphinxupquote{Host}} is set to \sphinxcode{\sphinxupquote{localhost}} if empty
or not specified; \sphinxcode{\sphinxupquote{QueueTime}} is set to 0 if empty or not specified. Specifically,
empty QueueTime means client does not wait and should quit immediately, if COPT Floating
Token Server have no tokens available. Port number must be great than zero and should be
the same as that specified in server configuration file.
Note that keywords in the client configuration file are case insensitive.

\sphinxAtStartPar
Without local license files, a COPT application still works if client configuration
file, \sphinxcode{\sphinxupquote{client.ini}}, exists in one of the following three locations in order, that is,
current working directory, environment directory by \sphinxcode{\sphinxupquote{COPT\_LICENSE\_DIR}} and binary
directory where COPT executable is located.

\sphinxAtStartPar
By design, COPT application reads local license files instead of client configuration file,
if they both exist in the same location. On the other hand, if local license files are under
the environment directory, to activate approach of floating licensing, user can simply
add a configuration file, \sphinxcode{\sphinxupquote{client.ini}}, under the current working directory (different from
the environment directory).

\sphinxAtStartPar
If a COPT application calls COPT API to solve problems, such as COPT python interface,
license is checked as soon as COPT environment object is created. If there exists proper client
configuration file, \sphinxcode{\sphinxupquote{client.ini}}, a license token is granted to COPT client. This license
token is released and sent back to token server, as soon as last COPT environment object in the
same process destroys.

\subsubsection{Via API Functions}
\label{\detokenize{floating:via-api-functions}}
\sphinxAtStartPar
In addition to the method using the \sphinxcode{\sphinxupquote{client.ini}} file mentioned above, users can also configure the client in their code through API functions. Taking the COPT Python interface as an example, the corresponding class is {\hyperref[\detokenize{pyapiref:chappyapi-envrconfig}]{\sphinxcrossref{\DUrole{std,std-ref}{EnvrConfig Class}}}}, and similar approaches apply to other programming languages. The configuration is as follows:

\begin{sphinxVerbatim}[commandchars=\\\{\}]
\PYG{c+c1}{\PYGZsh{} Set client configuration parameters}
\PYG{n}{envconfig}\PYG{o}{.}\PYG{n}{set}\PYG{p}{(}\PYG{n}{COPT}\PYG{o}{.}\PYG{n}{CLIENT\PYGZus{}FLOATING}\PYG{p}{,} \PYG{l+s+s2}{\PYGZdq{}}\PYG{l+s+s2}{192.168.1.11}\PYG{l+s+s2}{\PYGZdq{}}\PYG{p}{)}
\PYG{n}{envconfig}\PYG{o}{.}\PYG{n}{set}\PYG{p}{(}\PYG{n}{COPT}\PYG{o}{.}\PYG{n}{CLIENT\PYGZus{}PORT}\PYG{p}{,} \PYG{l+s+s2}{\PYGZdq{}}\PYG{l+s+s2}{7979}\PYG{l+s+s2}{\PYGZdq{}}\PYG{p}{)}
\PYG{n}{envconfig}\PYG{o}{.}\PYG{n}{set}\PYG{p}{(}\PYG{n}{COPT}\PYG{o}{.}\PYG{n}{CLIENT\PYGZus{}WAITTIME}\PYG{p}{,} \PYG{l+s+s2}{\PYGZdq{}}\PYG{l+s+s2}{600}\PYG{l+s+s2}{\PYGZdq{}}\PYG{p}{)}
\end{sphinxVerbatim}

\subsection{High Availability}
\label{\detokenize{floating:high-availability}}
\sphinxAtStartPar
In the case of multiple floating license servers, the client can achieve high availability by configuring multiple server IP addresses in the \sphinxcode{\sphinxupquote{Host}} field of the \sphinxcode{\sphinxupquote{client.ini}} file, as shown in the table below:

\begin{sphinxVerbatim}[commandchars=\\\{\}]
\PYG{n+nv}{Host}\PYG{+w}{ }\PYG{o}{=}\PYG{+w}{ }\PYG{l+m}{192}.168.1.11\PYG{p}{;}\PYG{+w}{ }\PYG{l+m}{192}.168.1.22\PYG{p}{;}\PYG{+w}{ }\PYG{l+m}{192}.168.1.33
\PYG{n+nv}{Port}\PYG{+w}{ }\PYG{o}{=}\PYG{+w}{ }\PYG{l+m}{7979}
\PYG{n+nv}{WaitTime}\PYG{+w}{ }\PYG{o}{=}\PYG{+w}{ }\PYG{l+m}{600}
\end{sphinxVerbatim}

\sphinxAtStartPar
Here, the \sphinxcode{\sphinxupquote{Host}} field contains the IP addresses of each floating license server. If the port number of a server is not \sphinxcode{\sphinxupquote{7979}}, the port can be specified after the IP address.

\sphinxAtStartPar
The above configuration indicates that the client will first attempt to connect to the server \sphinxcode{\sphinxupquote{192.168.1.11:7979}}. If that server is unavailable, the client will try connecting to \sphinxcode{\sphinxupquote{192.168.1.22:7979}}, and so on, until a successful connection to one of the floating servers is made.

\sphinxAtStartPar
If a floating server fails during token issuance or management, the client will reassign to another available floating license server as specified in the configuration file.
This setup ensures that the floating license service has disaster recovery capabilities, thereby improving the high availability of the COPT floating license service.

\subsection{Example Usage}
\label{\detokenize{floating:id2}}
\sphinxAtStartPar
Suppose that we’ve set client configuration file \sphinxcode{\sphinxupquote{client.ini}} properly and have no local
license, below is an example of obtaining a floating token by COPT command\sphinxhyphen{}line tool
\sphinxcode{\sphinxupquote{copt\_cmd}}. Execute the following command in Windows console, or Linux/Mac terminal.

\begin{sphinxVerbatim}[commandchars=\\\{\}]
copt\PYGZus{}cmd
\end{sphinxVerbatim}

\sphinxAtStartPar
If you see log information as follows, the COPT client, \sphinxcode{\sphinxupquote{copt\_cmd}}, has obtained the
floating token successfully. COPT command\sphinxhyphen{}line tool is ready to solve optimization problems.
\begin{sphinxalltt}
\textgreater{} copt\_cmd
  Cardinal Optimizer v8.0.1. Build date Mar 04 2024
  Copyright Cardinal Operations 2025. All Rights Reserved

  {[} Info{]} initialize floating client: ./client.ini

  {[} Info{]} connecting to server ...
  {[} Info{]} connection established
COPT\textgreater{}
\end{sphinxalltt}

\sphinxAtStartPar
If you see log information as follows, the COPT client, \sphinxcode{\sphinxupquote{copt\_cmd}}, has connected to COPT
Floating Token Server. But due to limited number of tokens, client waits in queue of size 1.
\begin{sphinxalltt}
\textgreater{} copt\_cmd
  Cardinal Optimizer v8.0.1. Build date Mar 04 2024
  Copyright Cardinal Operations 2025. All Rights Reserved

  {[} Info{]} Initialize floating client: ./client.ini

  {[} Info{]} connecting to server ...
  {[}Error{]} empty license and queue size 1
  {[} Info{]} Wait for license in  2 / 39 secs
  {[} Info{]} Wait for license in  4 / 39 secs
  {[} Info{]} Wait for license in  6 / 39 secs
  {[} Info{]} Wait for license in  8 / 39 secs
  {[} Info{]} Wait for license in 10 / 39 secs
  {[} Info{]} Wait for license in 20 / 39 secs
  {[} Info{]} Wait for license in 30 / 39 secs
\end{sphinxalltt}

\sphinxAtStartPar
If you see log information as follows, the COPT client, \sphinxcode{\sphinxupquote{copt\_cmd}}, has connected to COPT
Floating Token Server. But client refused to wait in queue, as Queuetime is 0.
\begin{sphinxalltt}
\textgreater{} copt\_cmd
  Cardinal Optimizer v8.0.1. Build date Mar 04 2024
  Copyright Cardinal Operations 2025. All Rights Reserved

  {[} Info{]} Initialize floating client: ./client.ini

  {[} Info{]} connecting to server ...
  {[}Error{]} Server error: "no more token available", code = 2
  {[}Error{]} Fail to open: ./license.dat

  {[}Error{]} Fail to initialize cmdline
\end{sphinxalltt}

\sphinxAtStartPar
If you see log information as follows, the COPT client, \sphinxcode{\sphinxupquote{copt\_cmd}}, fails to connect to COPT
Floating Token Server. Finally, client quits after time limit.
\begin{sphinxalltt}
\textgreater{} copt\_cmd
  Cardinal Optimizer v8.0.1. Build date Mar 04 2024
  Copyright Cardinal Operations 2025. All Rights Reserved

  {[} Info{]} initialize floating client: ./client.ini

  {[} Info{]} connecting to server ...
  {[} Info{]} wait for license in  2 / 10 secs
  {[} Info{]} wait for license in  4 / 10 secs
  {[} Info{]} wait for license in  6 / 10 secs
  {[} Info{]} wait for license in  8 / 10 secs
  {[} Info{]} wait for license in 10 / 10 secs
  {[}Error{]} timeout at waiting for license
  {[}Error{]} fail to open: ./license.dat

  {[}Error{]} Fail to initialize cmdline
\end{sphinxalltt}

\sphinxAtStartPar
In addition, users can start the COPT command\sphinxhyphen{}line tool on the client side to connect to the specified floating server IP address as follows:

\begin{sphinxVerbatim}[commandchars=\\\{\}]
\PYG{o}{\PYGZgt{}} \PYG{n}{copt\PYGZus{}cmd} \PYG{o}{\PYGZhy{}}\PYG{n}{floating} \PYG{o}{\PYGZlt{}}\PYG{n}{ip}\PYG{o}{\PYGZgt{}}
\end{sphinxVerbatim}

\section{Floating Token Server Managing Tool}
\label{\detokenize{floating:floating-token-server-managing-tool}}
\sphinxAtStartPar
COPT floating token service ships with a tool \sphinxcode{\sphinxupquote{copt\_flstool}}, for retrieving information
and tune parameters of floating token server on fly.

\subsection{Tool Usage}
\label{\detokenize{floating:tool-usage}}
\sphinxAtStartPar
Execute the following command in Windows console, Linux or MacOS terminal:

\begin{sphinxVerbatim}[commandchars=\\\{\}]
\PYGZgt{} ./copt\PYGZus{}flstool
\end{sphinxVerbatim}

\sphinxAtStartPar
Below displays help messages of this tool:

\begin{sphinxVerbatim}[commandchars=\\\{\}]
\PYGZgt{} ./copt\PYGZus{}flstool
  COPT Floating Token Server Managing Tool

  copt\PYGZus{}flstool [\PYGZhy{}s server ip] [\PYGZhy{}p port] [\PYGZhy{}x passwd] command \PYGZlt{}param\PYGZgt{}

  commands are:   addblackrule \PYGZlt{}127.0.0.1/20[\PYGZhy{}user@machine]\PYGZgt{}
                  addwhiterule \PYGZlt{}127.0.*.*[+user@machine]\PYGZgt{}
                  getfilters
                  getinfo
                  resetfilters
                  setpasswd \PYGZlt{}xxx\PYGZgt{}
                  toggleblackrule \PYGZlt{}n\PYGZhy{}th\PYGZgt{}
                  togglewhiterule \PYGZlt{}n\PYGZhy{}th\PYGZgt{}
                  writefilters
\end{sphinxVerbatim}

\sphinxAtStartPar
If the \sphinxcode{\sphinxupquote{\sphinxhyphen{}s}}  and \sphinxcode{\sphinxupquote{\sphinxhyphen{}p}} option are present, tool connects to floating token server
with given server IP and port. Otherwise, tool connections to localhost and
default port 7979. If floating token server sets a password, tool must provide
password string after the \sphinxcode{\sphinxupquote{\sphinxhyphen{}x}} option.

\sphinxAtStartPar
This tool provides the following commands:
\begin{itemize}
\item {} 
\sphinxAtStartPar
\sphinxcode{\sphinxupquote{AddBlackRule}}:
Add a new rule in black filters. each rule has format starting with non\sphinxhyphen{}empty
IP address, which may have wildcard to match IPs in the scope. In addition,
it is optional to be followed by including (+) or excluding (\sphinxhyphen{}) user name at machine name.

\item {} 
\sphinxAtStartPar
\sphinxcode{\sphinxupquote{AddWhiteRule}}:
Add a new rule in white filters. Note that a white rule has the same format as a black rule.

\item {} 
\sphinxAtStartPar
\sphinxcode{\sphinxupquote{GetFilters}}:
Get all rules of black filters, white filters and tool filters, along with relative sequence
numbers, which are parameters for command ToggleBlackRule and ToggleWhiteRule.

\item {} 
\sphinxAtStartPar
\sphinxcode{\sphinxupquote{GetInfo}}:
Get general information of floating token server, including token usage, connected clients,
and all COPT versions in support.

\item {} 
\sphinxAtStartPar
\sphinxcode{\sphinxupquote{ResetFilters}}:
Reset filter lists in memory to those on filter config file.

\item {} 
\sphinxAtStartPar
\sphinxcode{\sphinxupquote{SetPasswd}}:
Update password of target floating token server on fly.

\item {} 
\sphinxAtStartPar
\sphinxcode{\sphinxupquote{ToggleBlackRule}}:
Toggle between enabling and disabling a black rule, given its sequence number by GetFilters.

\item {} 
\sphinxAtStartPar
\sphinxcode{\sphinxupquote{ToggleWhiteRule}}:
Toggle between enabling and disabling a white rule, given its sequence number by GetFilters.

\item {} 
\sphinxAtStartPar
\sphinxcode{\sphinxupquote{WriteFilters}}:
Write filter lists in memory to filter config file.

\end{itemize}

\subsection{Example Usage}
\label{\detokenize{floating:id3}}
\sphinxAtStartPar
The following command lists general information on local machine.
\begin{sphinxalltt}
\textgreater{} ./copt\_flstool GetInfo

{[} Info{]} COPT Floating Token Server Managing Tool, COPT v8.0.1 20240304
{[} Info{]} connecting to localhost:7979
{[} Info{]} {[}command{]} wait for connecting to floating token server
{[} Info{]} {[}floating{]} general info
  \# of available tokens is 3 / 3, queue size is 0
  \# of active clients is 0
\end{sphinxalltt}

\sphinxAtStartPar
To run managing tool on other machine, its IP should be added to
a rule in \sphinxcode{\sphinxupquote{ToolList}} section in filter configuration file \sphinxcode{\sphinxupquote{flsfilters.ini}}.
The following command from other machine lists information of server 192.168.1.11.
\begin{sphinxalltt}
\textgreater{} ./copt\_flstool \sphinxhyphen{}s 192.168.1.11 GetInfo

{[} Info{]} COPT Floating Token Server Managing Tool, COPT v8.0.1 20240304
{[} Info{]} connecting to 192.168.1.11:7979
{[} Info{]} {[}command{]} wait for connecting to floating token server
{[} Info{]} {[}floating{]} general info
  \# of available tokens is 3 / 3, queue size is 0
  \# of active clients is 0
\end{sphinxalltt}

\sphinxAtStartPar
The following command shows all filter lists of server 192.168.1.11,
including those in BlackList section, WhiteList section and ToolList section.
\begin{sphinxalltt}
\textgreater{} ./copt\_flstool \sphinxhyphen{}s 192.168.1.11 GetFilters

{[} Info{]} COPT Floating Token Server Managing Tool, COPT v8.0.1 20240304
{[} Info{]} connecting to 192.168.1.11:7979
{[} Info{]} {[}command{]} wait for connecting to floating token server
{[} Info{]} {[}floating{]} filters info
{[}BlackList{]}

{[}WhiteList{]}

{[}ToolList{]}
  {[}1{]}  127.0.0.1
\end{sphinxalltt}

\sphinxAtStartPar
The following command added user of IP 192.168.3.13 to black list.
\begin{sphinxalltt}
\textgreater{} ./copt\_flstool \sphinxhyphen{}s 192.168.1.11 AddBlackRule 192.168.3.133

{[} Info{]} COPT Floating Token Server Managing Tool, COPT v8.0.1 20240304
{[} Info{]} connecting to 192.168.1.11:7979
{[} Info{]} {[}command{]} wait for connecting to floating token server
{[} Info{]} {[}floating{]} server added new black rule (succeeded)
\end{sphinxalltt}

\sphinxAtStartPar
The follwing command shows that a new rule in BlackList section is added.
\begin{sphinxalltt}
\textgreater{} ./copt\_flstool \sphinxhyphen{}s 192.168.1.11 GetFilters

{[} Info{]} COPT Floating Token Server Managing Tool, COPT v8.0.1 20240304
{[} Info{]} connecting to 192.168.1.11:7979
{[} Info{]} {[}command{]} wait for connecting to floating token server
{[} Info{]} {[}floating{]} filters info
{[}BlackList{]}
  {[}1{]} 192.168.3.133

{[}WhiteList{]}

{[}ToolList{]}
  {[}1{]}  127.0.0.1
\end{sphinxalltt}

\sphinxAtStartPar
The following command disable a rule in BlackList section.
\begin{sphinxalltt}
\textgreater{} ./copt\_flstool \sphinxhyphen{}s 192.168.1.11 ToggleBlackRule 1

{[} Info{]} COPT Floating Token Server Managing Tool, COPT v8.0.1 20240304
{[} Info{]} connecting to 192.168.1.11:7979
{[} Info{]} {[}command{]} wait for connecting to floating token server
{[} Info{]} {[}floating{]} server toggle black rule {[}1{]} (succeeded)
\end{sphinxalltt}

\section{Running as service}
\label{\detokenize{floating:running-as-service}}
\sphinxAtStartPar
To run COPT floating token server as a system service, follow steps described
in \sphinxcode{\sphinxupquote{readme.txt}} under \sphinxcode{\sphinxupquote{floating}} folder, and set config file
\sphinxcode{\sphinxupquote{copt\_flserver.service}} properly.

\sphinxAtStartPar
Below is \sphinxcode{\sphinxupquote{readme.txt}}, which lists installing steps in both Linux and MacOS
platforms.

\begin{sphinxVerbatim}[commandchars=\\\{\}]
[Linux] To run copt\PYGZus{}flserver as a service with systemd

Add a systemd file
    cp copt\PYGZus{}flserver.service to /lib/systemd/system/
    sudo systemctl daemon\PYGZhy{}reload

Enable new service
    sudo systemctl start copt\PYGZus{}flserver.service
    or
    sudo systemctl enable copt\PYGZus{}flserver.service

Restart service
    sudo systemctl restart copt\PYGZus{}flserver.service

Stop service
    sudo systemctl stop copt\PYGZus{}flserver.service
    or
    sudo systemctl disable copt\PYGZus{}flserver.service

Verify service is running
    sudo systemctl status copt\PYGZus{}flserver.service

[MacOS] To run copt\PYGZus{}flserver as a service with launchctrl

Add a plist file
    cp copt\PYGZus{}flserver.plist to /Library/LaunchAgents as current user
    or
    cp copt\PYGZus{}flserver.plist to /Library/LaunchDaemons with the key \PYGZsq{}UserName\PYGZsq{}

Enable new service
    sudo launchctl load \PYGZhy{}w /Library/LaunchAgents/copt\PYGZus{}flserver.plist
    or
    sudo launchctl load \PYGZhy{}w /Library/LaunchDaemons/copt\PYGZus{}flserver.plist

Stop service
    sudo launchctl unload \PYGZhy{}w /Library/LaunchAgents/copt\PYGZus{}flserver.plist
    or
    sudo launchctl unload \PYGZhy{}w /Library/LaunchDaemons/copt\PYGZus{}flserver.plist

Verify service is running
    sudo launchctl list solver.copt.flserver
\end{sphinxVerbatim}

\subsection{Linux}
\label{\detokenize{floating:linux}}
\sphinxAtStartPar
Below are steps in details of how to run COPT floating token server as a system service
in Linux platform.

\sphinxAtStartPar
For instance, assume that COPT remote service is installed under \sphinxcode{\sphinxupquote{\textquotesingle{}/home/eleven\textquotesingle{}}}.
In your terminal, type the following command to enter the root directory of
floating service.
\begin{sphinxalltt}
cd /home/eleven/copt\_remote80/floating
\end{sphinxalltt}

\sphinxAtStartPar
modify template of the service config file \sphinxcode{\sphinxupquote{copt\_flserver.service}} in text format:

\begin{sphinxVerbatim}[commandchars=\\\{\}]
[Unit]
Description=COPT Floating Token Server

[Service]
WorkingDirectory=/path/to/service
ExecStart=/path/to/service/copt\PYGZus{}flserver
Restart=always
RestartSec=1

[Install]
WantedBy=multi\PYGZhy{}user.target
\end{sphinxVerbatim}

\sphinxAtStartPar
That is, update template path in keyword \sphinxcode{\sphinxupquote{WorkingDirectory}} and \sphinxcode{\sphinxupquote{ExecStart}} to
actual path where the floating service exits.
\begin{sphinxalltt}
{[}Unit{]}
Description=COPT Floating Token Server

{[}Service{]}
WorkingDirectory=/home/eleven/copt\_remote80/floating
ExecStart=/home/eleven/copt\_remote80/floating/copt\_flserver
Restart=always
RestartSec=1

{[}Install{]}
WantedBy=multi\sphinxhyphen{}user.target
\end{sphinxalltt}

\sphinxAtStartPar
Afterwards, copy \sphinxcode{\sphinxupquote{copt\_flserver.service}} to system service folder
\sphinxcode{\sphinxupquote{/lib/systemd/system/}} (see below).

\begin{sphinxVerbatim}[commandchars=\\\{\}]
sudo cp copt\PYGZus{}flserver.service /lib/systemd/system/
\end{sphinxVerbatim}

\sphinxAtStartPar
The following command may be needed if you add or update service config file.
It is not needed if service unit has been loaded before.

\begin{sphinxVerbatim}[commandchars=\\\{\}]
sudo systemctl daemon\PYGZhy{}reload
\end{sphinxVerbatim}

\sphinxAtStartPar
The following command starts the new floating service.

\begin{sphinxVerbatim}[commandchars=\\\{\}]
sudo systemctl start copt\PYGZus{}flserver.service
\end{sphinxVerbatim}

\sphinxAtStartPar
To verify the floating service is actually running, type the following command

\begin{sphinxVerbatim}[commandchars=\\\{\}]
sudo systemctl status copt\PYGZus{}flserver.service
\end{sphinxVerbatim}

\sphinxAtStartPar
If you see logs similar to below, COPT floating server is running successfully as a system service.
\begin{sphinxalltt}
copt\_flserver.service \sphinxhyphen{} COPT Floating Token Server
Loaded: loaded (/lib/systemd/system/copt\_flserver.service; enabled; vendor preset: enabled)
Active: active (running) since Tue 2021\sphinxhyphen{}06\sphinxhyphen{}29 11:46:10 CST; 3s ago
Main PID: 3054 (copt\_flserver)
    Tasks: 6 (limit: 4915)
CGroup: /system.slice/copt\_flserver.service
        └─3054 /home/eleven/copt\_remote80/floating/copt\_flserver

eleven\sphinxhyphen{}ubuntu systemd{[}1{]}: Started COPT Floating Token Server.
eleven\sphinxhyphen{}ubuntu COPTCLS{[}3054{]}: LWS: 4.1.4\sphinxhyphen{}b2011a00, loglevel 1039
eleven\sphinxhyphen{}ubuntu COPTCLS{[}3054{]}: NET CLI SRV H1 H2 WS IPv6\sphinxhyphen{}absent
eleven\sphinxhyphen{}ubuntu COPTCLS{[}3054{]}: server started at port 7979
\end{sphinxalltt}

\sphinxAtStartPar
To stop the floating service, type the following command

\begin{sphinxVerbatim}[commandchars=\\\{\}]
sudo systemctl stop copt\PYGZus{}flserver.service
\end{sphinxVerbatim}

\subsection{MacOS}
\label{\detokenize{floating:macos}}
\sphinxAtStartPar
Below are steps in details of how to run COPT floating token server as a system service
in MacOS platform.

\sphinxAtStartPar
For instance, assume that COPT remote service is installed under \sphinxcode{\sphinxupquote{\textquotesingle{}/Applications\textquotesingle{}}}.
In your terminal, type the following command to enter the root directory of
floating service.
\begin{sphinxalltt}
cd /Applications/copt\_remote80/floating
\end{sphinxalltt}

\sphinxAtStartPar
modify template of the service config file \sphinxcode{\sphinxupquote{copt\_flserver.plist}} in xml format:

\begin{sphinxVerbatim}[commandchars=\\\{\}]
\PYG{c+cp}{\PYGZlt{}?xml version=\PYGZdq{}1.0\PYGZdq{} encoding=\PYGZdq{}UTF\PYGZhy{}8\PYGZdq{}?\PYGZgt{}}
\PYG{p}{\PYGZlt{}}\PYG{n+nt}{plist} \PYG{n+na}{version}\PYG{o}{=}\PYG{l+s}{\PYGZdq{}1.0\PYGZdq{}}\PYG{p}{\PYGZgt{}}
    \PYG{p}{\PYGZlt{}}\PYG{n+nt}{dict}\PYG{p}{\PYGZgt{}}
        \PYG{p}{\PYGZlt{}}\PYG{n+nt}{key}\PYG{p}{\PYGZgt{}}Label\PYG{p}{\PYGZlt{}}\PYG{p}{/}\PYG{n+nt}{key}\PYG{p}{\PYGZgt{}}
        \PYG{p}{\PYGZlt{}}\PYG{n+nt}{string}\PYG{p}{\PYGZgt{}}solver.copt.flserver\PYG{p}{\PYGZlt{}}\PYG{p}{/}\PYG{n+nt}{string}\PYG{p}{\PYGZgt{}}
        \PYG{p}{\PYGZlt{}}\PYG{n+nt}{key}\PYG{p}{\PYGZgt{}}Program\PYG{p}{\PYGZlt{}}\PYG{p}{/}\PYG{n+nt}{key}\PYG{p}{\PYGZgt{}}
        \PYG{p}{\PYGZlt{}}\PYG{n+nt}{string}\PYG{p}{\PYGZgt{}}/path/to/service/copt\PYGZus{}flserver\PYG{p}{\PYGZlt{}}\PYG{p}{/}\PYG{n+nt}{string}\PYG{p}{\PYGZgt{}}
        \PYG{p}{\PYGZlt{}}\PYG{n+nt}{key}\PYG{p}{\PYGZgt{}}RunAtLoad\PYG{p}{\PYGZlt{}}\PYG{p}{/}\PYG{n+nt}{key}\PYG{p}{\PYGZgt{}}
        \PYG{p}{\PYGZlt{}}\PYG{n+nt}{true}\PYG{p}{/}\PYG{p}{\PYGZgt{}}
        \PYG{p}{\PYGZlt{}}\PYG{n+nt}{key}\PYG{p}{\PYGZgt{}}KeepAlive\PYG{p}{\PYGZlt{}}\PYG{p}{/}\PYG{n+nt}{key}\PYG{p}{\PYGZgt{}}
        \PYG{p}{\PYGZlt{}}\PYG{n+nt}{true}\PYG{p}{/}\PYG{p}{\PYGZgt{}}
    \PYG{p}{\PYGZlt{}}\PYG{p}{/}\PYG{n+nt}{dict}\PYG{p}{\PYGZgt{}}
\PYG{p}{\PYGZlt{}}\PYG{p}{/}\PYG{n+nt}{plist}\PYG{p}{\PYGZgt{}}
\end{sphinxVerbatim}

\sphinxAtStartPar
That is, update template path in \sphinxcode{\sphinxupquote{Program}} tag to
actual path where the floating service exits.
\begin{sphinxalltt}
\textless{}?xml version="1.0" encoding="UTF\sphinxhyphen{}8"?\textgreater{}
\textless{}plist version="1.0"\textgreater{}
    \textless{}dict\textgreater{}
        \textless{}key\textgreater{}Label\textless{}/key\textgreater{}
        \textless{}string\textgreater{}solver.copt.flserver\textless{}/string\textgreater{}
        \textless{}key\textgreater{}Program\textless{}/key\textgreater{}
        \textless{}string\textgreater{}/Applications/copt\_remote80/floating/copt\_flserver\textless{}/string\textgreater{}
        \textless{}key\textgreater{}RunAtLoad\textless{}/key\textgreater{}
        \textless{}true/\textgreater{}
        \textless{}key\textgreater{}KeepAlive\textless{}/key\textgreater{}
        \textless{}true/\textgreater{}
    \textless{}/dict\textgreater{}
\textless{}/plist\textgreater{}
\end{sphinxalltt}

\sphinxAtStartPar
Afterwards, copy \sphinxcode{\sphinxupquote{copt\_flserver.plist}} to system service folder
\sphinxcode{\sphinxupquote{/Library/LaunchAgents}} (see below).

\begin{sphinxVerbatim}[commandchars=\\\{\}]
sudo cp copt\PYGZus{}flserver.plist /Library/LaunchAgents
\end{sphinxVerbatim}

\sphinxAtStartPar
The following command starts the new floating service.

\begin{sphinxVerbatim}[commandchars=\\\{\}]
sudo launchctl load \PYGZhy{}w /Library/LaunchAgents/copt\PYGZus{}flserver.plist
\end{sphinxVerbatim}

\sphinxAtStartPar
To verify the floating service is actually running, type the following command

\begin{sphinxVerbatim}[commandchars=\\\{\}]
sudo launchctl list solver.copt.flserver
\end{sphinxVerbatim}

\sphinxAtStartPar
If you see logs similar to below, COPT floating server is running successfully as a system service.
\begin{sphinxalltt}
\{
    "LimitLoadToSessionType" = "System";
    "Label" = "solver.copt.flserver";
    "OnDemand" = false;
    "LastExitStatus" = 0;
    "PID" = 16406;
    "Program" = "/Applications/copt\_remote80/floating/copt\_flserver";
\};
\end{sphinxalltt}

\sphinxAtStartPar
To stop the floating service, type the following command

\begin{sphinxVerbatim}[commandchars=\\\{\}]
sudo launchctl unload \PYGZhy{}w /Library/LaunchAgents/copt\PYGZus{}flserver.plist
\end{sphinxVerbatim}

\sphinxAtStartPar
If the floating service should be run by a specific user, add \sphinxcode{\sphinxupquote{UserName}} tag to conifg file.
Below adds a user \sphinxcode{\sphinxupquote{eleven}}, who has priviledge to run the floating service.
\begin{sphinxalltt}
\textless{}?xml version="1.0" encoding="UTF\sphinxhyphen{}8"?\textgreater{}
\textless{}plist version="1.0"\textgreater{}
    \textless{}dict\textgreater{}
        \textless{}key\textgreater{}Label\textless{}/key\textgreater{}
        \textless{}string\textgreater{}solver.copt.flserver\textless{}/string\textgreater{}
        \textless{}key\textgreater{}Program\textless{}/key\textgreater{}
        \textless{}string\textgreater{}/Applications/copt\_remote80/floating/copt\_flserver\textless{}/string\textgreater{}
        \textless{}key\textgreater{}UserName\textless{}/key\textgreater{}
        \textless{}string\textgreater{}eleven\textless{}/string\textgreater{}
        \textless{}key\textgreater{}RunAtLoad\textless{}/key\textgreater{}
        \textless{}true/\textgreater{}
        \textless{}key\textgreater{}KeepAlive\textless{}/key\textgreater{}
        \textless{}true/\textgreater{}
    \textless{}/dict\textgreater{}
\textless{}/plist\textgreater{}
\end{sphinxalltt}

\sphinxAtStartPar
Then copy new \sphinxcode{\sphinxupquote{copt\_flserver.plist}} to system service folder
\sphinxcode{\sphinxupquote{/Library/LaunchDaemons}} (see below).

\begin{sphinxVerbatim}[commandchars=\\\{\}]
sudo cp copt\PYGZus{}flserver.plist /Library/LaunchDaemons
\end{sphinxVerbatim}

\sphinxAtStartPar
The following command starts the new floating service.

\begin{sphinxVerbatim}[commandchars=\\\{\}]
sudo launchctl load \PYGZhy{}w /Library/LaunchDaemons/copt\PYGZus{}flserver.plist
\end{sphinxVerbatim}

\sphinxAtStartPar
To stop the floating service, type the following command

\begin{sphinxVerbatim}[commandchars=\\\{\}]
sudo launchctl unload \PYGZhy{}w /Library/LaunchDaemons/copt\PYGZus{}flserver.plist
\end{sphinxVerbatim}

\sphinxstepscope

\chapter{COPT Compute Cluster Service}
\label{\detokenize{cluster:copt-compute-cluster-service}}\label{\detokenize{cluster:chapcluster}}\label{\detokenize{cluster::doc}}
\sphinxAtStartPar
The \sphinxstylestrong{Cardinal Optimizer} provides COPT compute cluster service on all supported
platforms, which allows you to offload optimization computations
from COPT client applications over local network.

\sphinxAtStartPar
Once COPT compute cluster server runs at local network, any COPT client application with matching
COPT version can connect to server and offload optimization computations. That is, COPT compute cluster
clients are allowed to do modelling locally, execute optimization jobs remotely,
and then obtain results interactively.

\sphinxAtStartPar
Note that the more computing power server has, the more optimization jobs can run
simultaneously. Furthermore, COPT compute cluster service has functionality to cluster
multiple servers together and therefore serve more COPT compute cluster clients over local network.

\section{Server Setup}
\label{\detokenize{cluster:server-setup}}
\sphinxAtStartPar
The COPT compute cluster service includes \sphinxcode{\sphinxupquote{copt\_cluster}} executable
and a configuration file \sphinxcode{\sphinxupquote{cls.ini}}. The very first thing to do when cluster server starts
is to verify cluster license locally, whose path is specified in \sphinxcode{\sphinxupquote{cls.ini}}.
If local validation passes, cluster server might connect remotely to COPT licensing server
for further validation, such as verifying machine IP, which is supposed to match IP range
that user provided during registration. This means the server running COPT compute cluster service
should have internet access in specified area, such as campus network.
For details, please see descriptions below or refer to
{\hyperref[\detokenize{install:parcoptgetlic}]{\sphinxcrossref{\DUrole{std,std-ref}{How to obtain and setup license}}}}.

\subsection{Installation}
\label{\detokenize{cluster:installation}}
\sphinxAtStartPar
The \sphinxstylestrong{Cardinal Optimizer} provides a separate package for remote services,
which include COPT compute cluster. Users may apply for remote package from
customer service. Afterwards, unzip the remote package and move to any
folder on your computer. The software is portable and does not change
anything in the system it runs on. Below are details of installation.

\sphinxAtStartPar
\sphinxstylestrong{Windows}

\sphinxAtStartPar
Please unzip the remote package and move to any folder. Though, it is common
to move to folder under \sphinxcode{\sphinxupquote{C:\textbackslash{}Program Files}}.

\sphinxAtStartPar
\sphinxstylestrong{Linux}

\sphinxAtStartPar
To unzip the remote package, enter the following command in terminal:
\begin{sphinxalltt}
tar \sphinxhyphen{}xzf CardinalOptimizer\sphinxhyphen{}Remote\sphinxhyphen{}8.0.1\sphinxhyphen{}lnx64.tar.gz
\end{sphinxalltt}

\sphinxAtStartPar
Then, the following command moves folder copt\_remote80 in current
directory to other path. For an example, admin user may move it to folder
under \sphinxcode{\sphinxupquote{/opt}} and standard user may move it to \sphinxcode{\sphinxupquote{\$HOME}}.
\begin{sphinxalltt}
sudo mv copt\_remote80 /opt
\end{sphinxalltt}

\sphinxAtStartPar
Note that it requires \sphinxcode{\sphinxupquote{root}} privilege to execute this command.

\sphinxAtStartPar
\sphinxstylestrong{MacOS}

\sphinxAtStartPar
To unzip the remote package, enter the following command in terminal:
\begin{sphinxalltt}
tar \sphinxhyphen{}xzf CardinalOptimizer\sphinxhyphen{}Remote\sphinxhyphen{}8.0.1\sphinxhyphen{}universal\_mac.tar.gz
\end{sphinxalltt}

\sphinxAtStartPar
Then, the following command moves folder copt\_remote80 in current
directory to other path. For an example, admin user may move it to folder
under \sphinxcode{\sphinxupquote{/Applications}} and standard user may move it to \sphinxcode{\sphinxupquote{\$HOME}}.
\begin{sphinxalltt}
mv copt\_remote80 /Applications
\end{sphinxalltt}

\sphinxAtStartPar
If you see errors below or similar signature problem of COPT lib
during installation,

\begin{sphinxVerbatim}[commandchars=\\\{\}]
\PYGZdq{}libcopt.dylib\PYGZdq{} cannot be opened because the developer cannot be verified.
macOS cannot verify that this app is free from malware.
\end{sphinxVerbatim}

\sphinxAtStartPar
run the following command as root user, to bypass check of loading dynamic lib
on MacOS.
\begin{sphinxalltt}
xattr \sphinxhyphen{}d com.apple.quarantine CardinalOptimizer\sphinxhyphen{}Remote\sphinxhyphen{}8.0.1\sphinxhyphen{}universal\_mac.tar.gz
\end{sphinxalltt}

\sphinxAtStartPar
or
\begin{sphinxalltt}
xattr \sphinxhyphen{}dr com.apple.quarantine /Applications/copt\_remote80
\end{sphinxalltt}

\subsection{Cluster License}
\label{\detokenize{cluster:cluster-license}}
\sphinxAtStartPar
After installing COPT remote package, it requirs cluster license to run.
It is prefered to save cluster license files, \sphinxcode{\sphinxupquote{license.dat}} and
\sphinxcode{\sphinxupquote{license.key}}, to \sphinxcode{\sphinxupquote{cluster}} folder in path of remote package.

\sphinxAtStartPar
The following explains how to obtain the license file
via the \sphinxcode{\sphinxupquote{copt\_licgen}} tool and the license credential information \sphinxcode{\sphinxupquote{key}} under different systems.

\begin{sphinxadmonition}{note}{Note}

\sphinxAtStartPar
If the user has already obtained the two license files \sphinxcode{\sphinxupquote{license.dat}} and \sphinxcode{\sphinxupquote{license.key}}, there is no need to obtain them again.
You can skip the following steps to obtain the license file and refer to {\hyperref[\detokenize{floating:chapcluster-config}]{\sphinxcrossref{\DUrole{std,std-ref}{Configuration}}}} directly.
\end{sphinxadmonition}

\sphinxAtStartPar
\sphinxstylestrong{Windows}

\sphinxAtStartPar
If the COPT remote package is installed under \sphinxcode{\sphinxupquote{"C:\textbackslash{}Program Files"}},
execute the following command to enter \sphinxcode{\sphinxupquote{cluster}} folder in path of
remote package.
\begin{sphinxalltt}
cd "C:\textbackslash{}Program Files\textbackslash{}copt\_remote80\textbackslash{}cluster"
\end{sphinxalltt}

\sphinxAtStartPar
Note that the tool  \sphinxcode{\sphinxupquote{copt\_licgen}} creating license files exists
under \sphinxcode{\sphinxupquote{tools}} folder in path of remote package. The following command
creates cluster license files in current directory, given a cluster
license key, such as \sphinxcode{\sphinxupquote{7483dff0863ffdae9fff697d3573e8bc}} .

\begin{sphinxVerbatim}[commandchars=\\\{\}]
..\PYGZbs{}tools\PYGZbs{}copt\PYGZus{}licgen \PYGZhy{}key 7483dff0863ffdae9fff697d3573e8bc
\end{sphinxVerbatim}

\sphinxAtStartPar
\sphinxstylestrong{Linux and MacOS}

\sphinxAtStartPar
If the COPT remote package is installed under \sphinxcode{\sphinxupquote{"/Applications"}},
execute the following command to enter \sphinxcode{\sphinxupquote{cluster}} folder in path of
remote package on MacOS system.
\begin{sphinxalltt}
cd /Applications/copt\_remote80/cluster
\end{sphinxalltt}

\sphinxAtStartPar
The following command creates cluster license files in current directory,
given a cluster license key, such as \sphinxcode{\sphinxupquote{7483dff0863ffdae9fff697d3573e8bc}} .

\begin{sphinxVerbatim}[commandchars=\\\{\}]
../tools/copt\PYGZus{}licgen \PYGZhy{}key 7483dff0863ffdae9fff697d3573e8bc
\end{sphinxVerbatim}

\sphinxAtStartPar
In addition, if users run the above command when current directory is
different than \sphinxcode{\sphinxupquote{cluster}} folder in path of remote package, it is prefered
to move them to \sphinxcode{\sphinxupquote{cluster}}. The following command does so.
\begin{sphinxalltt}
mv license.* /Application/copt\_remote80/cluster
\end{sphinxalltt}

\subsection{Configuration}
\label{\detokenize{cluster:configuration}}
\sphinxAtStartPar
Below is a typical configuration file, \sphinxcode{\sphinxupquote{cls.ini}}, of COPT compute cluster.

\begin{sphinxVerbatim}[commandchars=\\\{\}]
\PYG{o}{[}Main\PYG{o}{]}
\PYG{n+nv}{Port}\PYG{+w}{ }\PYG{o}{=}\PYG{+w}{ }\PYG{l+m}{7878}
\PYG{c+c1}{\PYGZsh{} number of total tokens, which copt jobs can run simultaneously up to}
\PYG{n+nv}{NumToken}\PYG{+w}{ }\PYG{o}{=}\PYG{+w}{ }\PYG{l+m}{3}
\PYG{c+c1}{\PYGZsh{} password is case\PYGZhy{}sensitive and default is emtpy}
\PYG{c+c1}{\PYGZsh{} it applies to both copt clients and cluster nodes}
\PYG{n+nv}{PassWd}\PYG{+w}{ }\PYG{o}{=}
\PYG{c+c1}{\PYGZsh{} data folder of cluster relative to its binary folder,}
\PYG{c+c1}{\PYGZsh{} where various of copt libraries and temprary job files reside.}
\PYG{n+nv}{DataFolder}\PYG{+w}{ }\PYG{o}{=}\PYG{+w}{ }./data

\PYG{o}{[}SSL\PYG{o}{]}
\PYG{c+c1}{\PYGZsh{} needed if connecting using SSL}
\PYG{n+nv}{CaFile}\PYG{+w}{ }\PYG{o}{=}
\PYG{n+nv}{CertFile}\PYG{+w}{ }\PYG{o}{=}
\PYG{n+nv}{CertkeyFile}\PYG{+w}{ }\PYG{o}{=}

\PYG{o}{[}Licensing\PYG{o}{]}
\PYG{c+c1}{\PYGZsh{} if empty or default license name, it is from binary folder}
\PYG{c+c1}{\PYGZsh{} to get license files from cwd, add prefix \PYGZdq{}./\PYGZdq{}}
\PYG{c+c1}{\PYGZsh{} full path is supported as well}
\PYG{n+nv}{LicenseFile}\PYG{+w}{ }\PYG{o}{=}\PYG{+w}{ }license.dat
\PYG{n+nv}{PubkeyFile}\PYG{+w}{ }\PYG{o}{=}\PYG{+w}{ }license.key

\PYG{o}{[}WLS\PYG{o}{]}
\PYG{c+c1}{\PYGZsh{} WebServer have a default host and no need to edit in most scenarios}
\PYG{c+c1}{\PYGZsh{} Must specify WebLicenseId and WebAccesskey to trigger web licensing}
\PYG{c+c1}{\PYGZsh{} If specified, ignore settings in Licensing section}
\PYG{n+nv}{WebServer}\PYG{+w}{ }\PYG{o}{=}
\PYG{n+nv}{WebLicenseId}\PYG{+w}{ }\PYG{o}{=}
\PYG{n+nv}{WebAccessKey}\PYG{+w}{ }\PYG{o}{=}
\PYG{n+nv}{WebTokenDuration}\PYG{+w}{ }\PYG{o}{=}\PYG{+w}{ }\PYG{l+m}{300}

\PYG{o}{[}FLS\PYG{o}{]}
\PYG{c+c1}{\PYGZsh{} Lease token from floating server}
\PYG{c+c1}{\PYGZsh{} If specified, ignore settings in Licensing section and WLS section}
\PYG{n+nv}{TokenServer}\PYG{+w}{ }\PYG{o}{=}
\PYG{n+nv}{TPort}\PYG{+w}{ }\PYG{o}{=}\PYG{+w}{ }\PYG{l+m}{7979}

\PYG{o}{[}Cluster\PYG{o}{]}
\PYG{c+c1}{\PYGZsh{} host name and port of parent node in cluster}
\PYG{c+c1}{\PYGZsh{} If specified or changed on fly, connect to parent node}
\PYG{n+nv}{Parent}\PYG{+w}{ }\PYG{o}{=}
\PYG{n+nv}{PPort}\PYG{+w}{ }\PYG{o}{=}\PYG{+w}{ }\PYG{l+m}{7878}

\PYG{o}{[}Filter\PYG{o}{]}
\PYG{c+c1}{\PYGZsh{} default policy 0 indicates accepting all connections, except for ones in blacklist}
\PYG{c+c1}{\PYGZsh{} otherwise, denying all connections except for ones in whitelist}
\PYG{n+nv}{DefaultPolicy}\PYG{+w}{ }\PYG{o}{=}\PYG{+w}{ }\PYG{l+m}{0}
\PYG{n+nv}{UseBlackList}\PYG{+w}{ }\PYG{o}{=}\PYG{+w}{ }\PYG{n+nb}{true}
\PYG{n+nv}{UseWhiteList}\PYG{+w}{ }\PYG{o}{=}\PYG{+w}{ }\PYG{n+nb}{true}
\PYG{n+nv}{FilterListFile}\PYG{+w}{ }\PYG{o}{=}\PYG{+w}{ }clsfilters.ini

\PYG{o}{[}Logs\PYG{o}{]}
\PYG{n+nv}{LogsFolder}\PYG{+w}{ }\PYG{o}{=}\PYG{+w}{ }./logs
\end{sphinxVerbatim}

\sphinxAtStartPar
The \sphinxcode{\sphinxupquote{Main}} section specifies port number, through which COPT compute cluster clients
connect to server;
token number, the number of optimization jobs that
server can run simultaneously up to;
password string, if specified, cluster
clients should send the same password when requesting for service.

\sphinxAtStartPar
\sphinxstylestrong{Main:}

\sphinxAtStartPar
In the \sphinxcode{\sphinxupquote{{[}Main{]}}} section of the configuration file, users can set the following
connection information:
\begin{itemize}
\item {} 
\sphinxAtStartPar
\sphinxcode{\sphinxupquote{Port}} : Specifies the connection port for the cluster server, through which
COPT compute cluster clients connect to server;

\item {} 
\sphinxAtStartPar
\sphinxcode{\sphinxupquote{NumToken}} : Defines the number of optimization jobs that the COPT cluster server
can run simultaneously up to. This number can be estimated based on the server’s
hardware resource capacity and is user\sphinxhyphen{}configurable. The default value is 3, with no concurrency limit.

\item {} 
\sphinxAtStartPar
\sphinxcode{\sphinxupquote{PassWd}} : Sets the password for client connections to the cluster server.
The default is empty, meaning no password is required. If specified, clients
should send the same password when requesting for cluster service.
For a more secure connection, users can enable encrypted communication using RSA
certificates in the \sphinxcode{\sphinxupquote{{[}SSL{]}}} section.

\item {} 
\sphinxAtStartPar
\sphinxcode{\sphinxupquote{DataFolder}} : Specifies the path where the COPT solver is installed on the server.

\end{itemize}

\sphinxAtStartPar
The COPT compute cluster may install multiple versions of COPT to subfolder of \sphinxcode{\sphinxupquote{DataFolder}}.
Only clients with matching version (major and minor) will get approved and then offload optimization jobs at server side. Note that the COPT compute cluster pre\sphinxhyphen{}installs a COPT solver of the same version as server itself, which illustrate how to install other versions of COPT.

\sphinxAtStartPar
For instance, the COPT compute cluster has default COPT v8.0.1 installed and users plan to install
COPT of other version v6.0.7. Users may create a folder ./data/copt/6.0.7/ and
copy a COPT C lib of the same version to it. Specifically, on Linux platform,
copy C dynamic library \sphinxcode{\sphinxupquote{libcopt.so}} from the binary folder \sphinxcode{\sphinxupquote{\$COPT\_HOME/lib/}} of COPT v6.0.7
to subfolder ./data/copt/6.0.7/ of the COPT compute cluster.

\sphinxAtStartPar
Furthermore, users are allowed to install newer version of COPT than cluster server version,
such as COPT v9.0.0.
To do so, follow the same step of copying a C lib of COPT v9.0.0 to ./data/copt/9.0.0/.
In addition, users need a personal license of v9.0.0 to load C lib of COPT v9.0.0 at
server side. That is, copy valid personal license files to folder ./data/copt/9.0.0/ as well.
However, this simple procedure may break if the newer COPT solver has significant updates.
In this case, it is necessary to upgrade the COPT compute cluster to newer version, that is, v9.0.0.

\sphinxAtStartPar
Below is an example of directory structure of cluster server on Linux platform. It includes
pre\sphinxhyphen{}installed COPT v8.0.1, COPT of previous version v6.0.7,
and COPT of newer version v9.0.0.
\begin{sphinxalltt}
\textasciitilde{}/copt\_remote80/cluster
│  cls.ini
│  copt\_cluster
│  license.dat \sphinxhyphen{}\textgreater{} cluster license v8.0.1
│  license.key
│
└─data
    └─copt
        └─6.0.7
        │   libcopt.so
        └─8.0.1
        │   libcopt.so
        └─9.0.0
            libcopt.so
            license.dat \sphinxhyphen{}\textgreater{} license v9.0.0
            license.key
\end{sphinxalltt}

\sphinxAtStartPar
\sphinxstylestrong{Licensing:}

\sphinxAtStartPar
The \sphinxcode{\sphinxupquote{Licensing}} section specifies location of cluster license files.
As described in comments above, if emtpy string or default license file name
( \sphinxcode{\sphinxupquote{license.dat}} or \sphinxcode{\sphinxupquote{license.key}} ) is specified, cluster license files are
read from the binary folder where the cluster executive reside.

\sphinxAtStartPar
It is possible to run COPT compute cluster service, even if
cluster license files do not exist in the same folder as the cluster executive.
One option is to set \sphinxcode{\sphinxupquote{LicenseFile = ./license.dat}} and \sphinxcode{\sphinxupquote{PubkeyFile = ./license.key}}.
By doing so, the COPT compute cluster reads cluster license files from the current
working directory. That is, user could execute command at the path where cluster
license files exist to run service.

\sphinxAtStartPar
The other option is to set full path of license files in configuration.
As mentioned before, Cardinal Optimizer allows user to set
environment variable \sphinxcode{\sphinxupquote{COPT\_LICENSE\_DIR}} for license files. For details,
please refer to {\hyperref[\detokenize{install:chapinstall}]{\sphinxcrossref{\DUrole{std,std-ref}{How to install Cardinal Optimizer}}}}.
If user prefers the way of environment variable, \sphinxcode{\sphinxupquote{cls.ini}} should have
the full path to cluster license files.

\sphinxAtStartPar
\sphinxstylestrong{Cluster:}

\sphinxAtStartPar
In the \sphinxcode{\sphinxupquote{{[}Cluster{]}}} section of the configuration file, users can set the parent node (IP and port) for the current cluster server connection. By default, this is empty, indicating no connection to other nodes.
\begin{itemize}
\item {} 
\sphinxAtStartPar
\sphinxcode{\sphinxupquote{Parent}}: The IP address of the parent node within the local network.

\item {} 
\sphinxAtStartPar
\sphinxcode{\sphinxupquote{PPort}}: The port number of the parent node.

\end{itemize}

\sphinxAtStartPar
If there is only one cluster server or if the current server is the root node, this section does not need to be configured.

\sphinxAtStartPar
For multiple cluster servers, specifying a parent node allows users to form a tree\sphinxhyphen{}structured network topology for the COPT computing cluster. The current server can join the parent node’s cluster through this configuration.

\sphinxAtStartPar
\sphinxstylestrong{Filter:}

\sphinxAtStartPar
In the \sphinxcode{\sphinxupquote{{[}Filter{]}}} section of the configuration file, users can configure the filtering policy for the cluster server.
\begin{itemize}
\item {} 
\sphinxAtStartPar
\sphinxcode{\sphinxupquote{DefaultPolicy}}: Default is set to 0, meaning all clients are allowed to connect to   the cluster server, except those on the blacklist. If it is set to non\sphinxhyphen{}zero value, then all connection are blocked except for those in white lists.

\item {} 
\sphinxAtStartPar
\sphinxcode{\sphinxupquote{UseBlackList}}: If set to \sphinxcode{\sphinxupquote{True}}, the blacklist will be enabled.

\item {} 
\sphinxAtStartPar
\sphinxcode{\sphinxupquote{UseWhiteList}}: If set to \sphinxcode{\sphinxupquote{True}}, the whitelist will be enabled.

\item {} 
\sphinxAtStartPar
\sphinxcode{\sphinxupquote{FilterListFile}}: Specifies the name of the filtering configuration file for the cluster server, default is \sphinxcode{\sphinxupquote{clsfilters.ini}}. Below is an example of the filter configuration file:

\end{itemize}

\begin{sphinxVerbatim}[commandchars=\\\{\}]
\PYG{o}{[}BlackList\PYG{o}{]}
\PYG{c+c1}{\PYGZsh{} 127.0.*.* + user@machine*}

\PYG{o}{[}WhiteList\PYG{o}{]}
\PYG{c+c1}{\PYGZsh{} 127.0.1.2/16 \PYGZhy{} user@machine*}

\PYG{o}{[}ToolList\PYG{o}{]}
\PYG{c+c1}{\PYGZsh{} only tool client at server side can access by default}
\PYG{l+m}{127}.0.0.1/32
\end{sphinxVerbatim}

\sphinxAtStartPar
It has three sections and each section has its own rules. In section of \sphinxcode{\sphinxupquote{BlackList}},
one may add rules to block others from connection. In section of \sphinxcode{\sphinxupquote{WhiteList}},
one may add rules to grant others for connection, even if the default policy is to block
all connections. Only users listed in section of \sphinxcode{\sphinxupquote{ToolList}} are able to connect to
cluster server by Cluster Managing Tool (see below for details).

\sphinxAtStartPar
Specifically, rules in filter configuration have format of starting with IP address.
To specify IP range, users may include wildcard (*) in IP address, or use CIDR notation,
that is, a IPv4 address and its associated network prefix.
In addtion, a rule may include (+) or exclude (\sphinxhyphen{}) given user at given machine, such as
\sphinxcode{\sphinxupquote{127.0.1.2/16 \sphinxhyphen{} user@machine}}. Here, \sphinxcode{\sphinxupquote{user}} refers to \sphinxcode{\sphinxupquote{username}}, which can be queried
by \sphinxcode{\sphinxupquote{whoami}} on Linux/MacOS platform; \sphinxcode{\sphinxupquote{machine}} refers to \sphinxcode{\sphinxupquote{computer name}}, which can
be queried by \sphinxcode{\sphinxupquote{hostname}} on Linux/MacOS platform.

\sphinxAtStartPar
Note that after modifying the configuration file \sphinxcode{\sphinxupquote{clsfilters.ini}} , users can use the \sphinxcode{\sphinxupquote{ResetFilters}} command to reset the current rules to those in the filter configuration document. Users can also use the \sphinxcode{\sphinxupquote{WriteFilters}} command to output the current rules to the filter configuration document.

\sphinxAtStartPar
\sphinxstylestrong{Logs:}

\sphinxAtStartPar
In the \sphinxcode{\sphinxupquote{{[}Logs{]}}} section of the configuration file, users can set the log file path for the cluster server, which is by default stored in the \sphinxcode{\sphinxupquote{./logs}} path under the cluster server installation directory.

\subsection{Web License for Compute Cluster}
\label{\detokenize{cluster:web-license-for-compute-cluster}}
\sphinxAtStartPar
Besides local cluster license above, users may use web license for compute
cluster to run compute cluster service. This requires that machines running
compute cluster server must have internet access. However, hardware info are
not required any more. That is, users are free to deploy compute cluster
servers to any cloud machine or container, as long as they have internet
access. Please refer to \sphinxhref{https://copt.shanshu.ai/license}{COPT web license page} for details.

\sphinxAtStartPar
Below are brief steps:
\begin{itemize}
\item {} 
\sphinxAtStartPar
Follow steps to register an account and apply for trial of
web license for compute cluster.

\item {} 
\sphinxAtStartPar
Once approved, \sphinxcode{\sphinxupquote{Web License ID}} is generated for users

\item {} 
\sphinxAtStartPar
On page of \sphinxcode{\sphinxupquote{API Keys}} , create \sphinxcode{\sphinxupquote{Web Access Key}} using
given \sphinxcode{\sphinxupquote{Web License ID}}

\end{itemize}

\sphinxAtStartPar
Afterwards, users edit configuration file \sphinxcode{\sphinxupquote{cls.ini}} and add values of both
\sphinxcode{\sphinxupquote{Web License ID}} and \sphinxcode{\sphinxupquote{Web Access Key}} to related keywords in section of
\sphinxcode{\sphinxupquote{WLS}} . For instance,

\begin{sphinxVerbatim}[commandchars=\\\{\}]
\PYG{o}{[}WLS\PYG{o}{]}
\PYG{c+c1}{\PYGZsh{} WebServer have a default host and no need to edit in most scenarios}
\PYG{c+c1}{\PYGZsh{} Must specify WebLicenseId and WebAccesskey to trigger web licensing}
\PYG{n+nv}{WebServer}\PYG{+w}{ }\PYG{o}{=}
\PYG{n+nv}{WebLicenseId}\PYG{+w}{ }\PYG{o}{=}
\PYG{n+nv}{WebAccessKey}\PYG{+w}{ }\PYG{o}{=}
\PYG{n+nv}{WebTokenDuration}\PYG{+w}{ }\PYG{o}{=}\PYG{+w}{ }\PYG{l+m}{300}
\end{sphinxVerbatim}

\sphinxAtStartPar
As of now, compute cluster server talks to \sphinxhref{https://copt.shanshu.ai/license}{COPT web license page}
for licensing. Users are able to monitor its token usage and other
information online.

\subsection{Example Usage}
\label{\detokenize{cluster:example-usage}}
\sphinxAtStartPar
Suppose that cluster license files exist in the same folder where the
cluster executable reside. To start the COPT compute cluster, just
execute the following command at any directory in Windows console,
or Linux/Mac terminal.

\begin{sphinxVerbatim}[commandchars=\\\{\}]
./copt\PYGZus{}cluster
\end{sphinxVerbatim}

\sphinxAtStartPar
If you see log information as follows, the COPT compute cluster has been
successfully started. Server monitors any connection from COPT compute cluster clients, manages
client requests in queue as well as approved clients. User may stop cluster server
anytime when entering \sphinxcode{\sphinxupquote{q}} or \sphinxcode{\sphinxupquote{Q}}.
\begin{sphinxalltt}
\textgreater{} ./copt\_cluster
  {[} Info{]} start COPT Compute Cluster, COPT v8.0.1 20240304
  {[} Info{]} {[}NODE{]} node has been initialized
  {[} Info{]} server started at port 7878
\end{sphinxalltt}

\sphinxAtStartPar
If failed to verify local cluster license, or something is wrong on remote
COPT licensing server, you might see error logs as follows.
\begin{sphinxalltt}
\textgreater{} ./copt\_cluster
  {[} Info{]} start COPT Compute Cluster, COPT v8.0.1 20240304
  {[}Error{]} Invalid signature in public key file
  {[}Error{]} Fail to verify local license
\end{sphinxalltt}

\sphinxAtStartPar
and
\begin{sphinxalltt}
\textgreater{} ./copt\_cluster
  {[} Info{]} start COPT Compute Cluster, COPT v8.0.1 20240304
  {[}Error{]} Error to connect license server
  {[}Error{]} Fail to verify cluster license by server
\end{sphinxalltt}

\section{Client Setup}
\label{\detokenize{cluster:client-setup}}
\sphinxAtStartPar
The COPT compute cluster client can be COPT command\sphinxhyphen{}line, or any application which
solves problems by COPT API, such as COPT cpp/java/csharp/python interface.
The COPT compute cluster service
is a better approach in terms of flexibility and efficiency.
Any COPT compute cluster client can legally run Cardingal Optimizer without local license.

\subsection{Configuration}
\label{\detokenize{cluster:id1}}

\subsubsection{Via the configuration file}
\label{\detokenize{cluster:via-the-configuration-file}}
\sphinxAtStartPar
Before running COPT as cluster client, please make sure that you have
installed COPT locally. For details, please refer to
{\hyperref[\detokenize{install:chapinstall}]{\sphinxcrossref{\DUrole{std,std-ref}{How to install Cardinal Optimizer}}}}.
Users can skip obtaining local licenses by adding a cluster configuration
file \sphinxcode{\sphinxupquote{client.ini}} .

\sphinxAtStartPar
Below is a typical configuration file, \sphinxcode{\sphinxupquote{client.ini}}, of COPT compute cluster client.

\begin{sphinxVerbatim}[commandchars=\\\{\}]
\PYG{n+nv}{Cluster}\PYG{+w}{ }\PYG{o}{=}\PYG{+w}{ }\PYG{l+m}{192}.168.1.11
\PYG{n+nv}{Port}\PYG{+w}{ }\PYG{o}{=}\PYG{+w}{ }\PYG{l+m}{7878}
\PYG{n+nv}{WaitTime}\PYG{+w}{ }\PYG{o}{=}\PYG{+w}{ }\PYG{l+m}{600}
\PYG{n+nv}{Passwd}\PYG{+w}{ }\PYG{o}{=}
\end{sphinxVerbatim}

\sphinxAtStartPar
As configured above, COPT compute cluster client tries to connect to \sphinxcode{\sphinxupquote{192.168.1.11}} at port \sphinxcode{\sphinxupquote{7878}} with waiting time in queue up to 600 seconds.
Here, the default value of \sphinxcode{\sphinxupquote{Cluster}} is localhost.
\sphinxcode{\sphinxupquote{WaitTime}} (or \sphinxcode{\sphinxupquote{QueueTime}} ) is set to 0 if empty or not specified. Specifically, empty \sphinxcode{\sphinxupquote{WaitTime}} means client does not wait and should quit immediately, if the COPT compute cluster have no more token available.
\sphinxcode{\sphinxupquote{Port}} default number is 7878. It must be great than zero if specified
and should be the same as that specified in cluster configuration file \sphinxcode{\sphinxupquote{cls.ini}}.
Note that keywords in the configuration file are case insensitive.

\sphinxAtStartPar
In addition, users can set the password for connecting the remote server through \sphinxcode{\sphinxupquote{Passwd}}. The \sphinxcode{\sphinxupquote{Priority}} can be used to set the priority of cluster optimization jobs. Possible values range from 0 to 99, with higher values indicating higher priority.If jobs are queued, the priority setting will ensure that the next task is processed first, but it will not affect tasks that are already running.

\sphinxAtStartPar
To run as a COPT compute cluster client, an application must have configuration
file, \sphinxcode{\sphinxupquote{client.ini}}, in one of the following three locations, that is,
current working directory, environment directory by \sphinxcode{\sphinxupquote{COPT\_LICENSE\_DIR}} and binary
directory where COPT executable resides.

\sphinxAtStartPar
By design, COPT application reads local license files instead of \sphinxcode{\sphinxupquote{client.ini}},
if they both exist in the same location. However, if local license files are under
the environment directory, to connect to cluster server, user could simply
add a configuration file, \sphinxcode{\sphinxupquote{client.ini}}, under the current working directory (different from
the environment directory).

\sphinxAtStartPar
If a COPT application calls COPT API to solve problems, such as COPT python interface,
license is checked as soon as COPT environment object is created. If there only exists proper
configuration file, \sphinxcode{\sphinxupquote{client.ini}}, the application works as a COPT compute cluster client and obtains token to offload optimization jobs. As soon as COPT environment object is destroyed, the COPT compute cluster server is notified to release token and thus to approve more requests waiting in queue.

\subsubsection{Via API Functions}
\label{\detokenize{cluster:via-api-functions}}
\sphinxAtStartPar
In addition to the method using the \sphinxcode{\sphinxupquote{client.ini}} file mentioned above, users can also configure the client in their code through API functions. Taking the COPT Python interface as an example, the corresponding class is {\hyperref[\detokenize{pyapiref:chappyapi-envrconfig}]{\sphinxcrossref{\DUrole{std,std-ref}{EnvrConfig Class}}}}, and similar approaches apply to other programming languages. The configuration is as follows:

\begin{sphinxVerbatim}[commandchars=\\\{\}]
\PYG{c+c1}{\PYGZsh{} Set client configuration parameters}
\PYG{n}{envconfig}\PYG{o}{.}\PYG{n}{set}\PYG{p}{(}\PYG{n}{COPT}\PYG{o}{.}\PYG{n}{CLIENT\PYGZus{}CLUSTER}\PYG{p}{,} \PYG{l+s+s2}{\PYGZdq{}}\PYG{l+s+s2}{192.168.1.11}\PYG{l+s+s2}{\PYGZdq{}}\PYG{p}{)}
\PYG{n}{envconfig}\PYG{o}{.}\PYG{n}{set}\PYG{p}{(}\PYG{n}{COPT}\PYG{o}{.}\PYG{n}{CLIENT\PYGZus{}PORT}\PYG{p}{,} \PYG{l+s+s2}{\PYGZdq{}}\PYG{l+s+s2}{7878}\PYG{l+s+s2}{\PYGZdq{}}\PYG{p}{)}
\PYG{n}{envconfig}\PYG{o}{.}\PYG{n}{set}\PYG{p}{(}\PYG{n}{COPT}\PYG{o}{.}\PYG{n}{CLIENT\PYGZus{}WAITTIME}\PYG{p}{,} \PYG{l+s+s2}{\PYGZdq{}}\PYG{l+s+s2}{600}\PYG{l+s+s2}{\PYGZdq{}}\PYG{p}{)}
\end{sphinxVerbatim}

\subsection{High Availability}
\label{\detokenize{cluster:high-availability}}
\sphinxAtStartPar
When there are multiple cluster servers, the client can achieve high availability by configuring multiple Cluster server IP addresses in the \sphinxcode{\sphinxupquote{Cluster}} field of \sphinxcode{\sphinxupquote{client.ini}}, as shown in the following table:

\begin{sphinxVerbatim}[commandchars=\\\{\}]
\PYG{n+nv}{Cluster}\PYG{+w}{ }\PYG{o}{=}\PYG{+w}{ }\PYG{l+m}{192}.168.1.11\PYG{p}{;}\PYG{+w}{ }\PYG{l+m}{192}.168.1.22\PYG{p}{;}\PYG{+w}{ }\PYG{l+m}{192}.168.1.33
\PYG{n+nv}{Port}\PYG{+w}{ }\PYG{o}{=}\PYG{+w}{ }\PYG{l+m}{7878}
\PYG{n+nv}{WaitTime}\PYG{+w}{ }\PYG{o}{=}\PYG{+w}{ }\PYG{l+m}{600}
\PYG{n+nv}{Passwd}\PYG{+w}{ }\PYG{o}{=}
\end{sphinxVerbatim}

\sphinxAtStartPar
According to the above configuration, the first IP in the \sphinxcode{\sphinxupquote{Cluster}} field is the root node IP of the cluster. If the port of other child nodes of the cluster is not \sphinxcode{\sphinxupquote{7878}}, it could be added after the IP address.
The above configuration file indicates that the client will first try to connect to the root node \sphinxcode{\sphinxupquote{192.168.1.11:7878}} .
If the root node is not available, it will try to connect to the child node of \sphinxcode{\sphinxupquote{192.168.1.22:7878}} , and so on, until the connection with the cluster server is successful.

\sphinxAtStartPar
If a cluster node fails when performing a computing task, the client will reallocate the available cluster nodes according to the backup nodes in the configuration file. This can enable the cluster service to have a certain disaster recovery capability, thereby improving the high availability of the COPT computing cluster service.

\subsection{Example Usage}
\label{\detokenize{cluster:id2}}
\sphinxAtStartPar
Suppose that we’ve set configuration file \sphinxcode{\sphinxupquote{client.ini}} properly and have no local
license, below is an example of connecting to cluster server by COPT command\sphinxhyphen{}line tool
\sphinxcode{\sphinxupquote{copt\_cmd}}. Execute the following command in Windows console, or Linux/Mac terminal.

\begin{sphinxVerbatim}[commandchars=\\\{\}]
copt\PYGZus{}cmd
\end{sphinxVerbatim}

\sphinxAtStartPar
If you see log information as follows, the COPT compute cluster client, \sphinxcode{\sphinxupquote{copt\_cmd}}, has
connected to cluster server successfully. COPT command\sphinxhyphen{}line tool is ready to do
modelling locally and then offload optimization jobs to server.
\begin{sphinxalltt}
\textgreater{} copt\_cmd
  Cardinal Optimizer v8.0.1. Build date Mar 04 2024
  Copyright Cardinal Operations 2025. All Rights Reserved

  {[} Info{]} initialize cluster client with ./client.ini

  {[} Info{]} wait for server in 0 / 39 secs
  {[} Info{]} connecting to cluster server 192.168.1.11:7878
COPT\textgreater{}
\end{sphinxalltt}

\sphinxAtStartPar
If you see log information as follows, the COPT compute cluster client, \sphinxcode{\sphinxupquote{copt\_cmd}}, has connected to
cluster server. However, due to limited number of tokens, it waits in queue of size 5, until timeout.
\begin{sphinxalltt}
\textgreater{} copt\_cmd
  Cardinal Optimizer v8.0.1. Build date Mar 04 2024
  Copyright Cardinal Operations 2025. All Rights Reserved

  {[} Info{]} initialize cluster client with ./client.ini

  {[} Info{]} wait for server in 0 / 39 secs
  {[} Info{]} connecting to cluster server 192.168.1.11:7878

  {[} Warn{]} wait in queue of size 5
  {[} Info{]} wait for license in  2 / 39 secs
  {[} Info{]} wait for license in  4 / 39 secs
  {[} Info{]} wait for license in  6 / 39 secs
  {[} Info{]} wait for license in  8 / 39 secs
  {[} Info{]} wait for license in 10 / 39 secs
  {[} Info{]} wait for license in 20 / 39 secs
  {[} Info{]} wait for license in 30 / 39 secs
  {[}Error{]} timeout at waiting for server approval
  {[}Error{]} Fail to initialize copt command\sphinxhyphen{}line tool
\end{sphinxalltt}

\sphinxAtStartPar
If you see log information as follows, the COPT compute cluster client, \sphinxcode{\sphinxupquote{copt\_cmd}}, has connected to
cluster server. But it refused to wait in queue, as Queuetime was set to 0. Therefore,
client quits with error immediately.
\begin{sphinxalltt}
\textgreater{} copt\_cmd
  Cardinal Optimizer v8.0.1. Build date Mar 04 2024
  Copyright Cardinal Operations 2025. All Rights Reserved

  {[} Info{]} initialize cluster client with ./client.ini

  {[} Info{]} wait for server in 0 / 9 secs
  {[} Info{]} connecting to cluster server 192.168.1.11:7878
  {[} Warn{]} server error: "no more token available", code = 129
  {[}Error{]} Fail to initialize copt command\sphinxhyphen{}line tool
\end{sphinxalltt}

\sphinxAtStartPar
If you see log information as follows, the COPT compute cluster client, \sphinxcode{\sphinxupquote{copt\_cmd}}, fails to connect to
cluster server. Finally, client quits after timeout.
\begin{sphinxalltt}
\textgreater{} copt\_cmd
  Cardinal Optimizer v8.0.1. Build date Mar 04 2024
  Copyright Cardinal Operations 2025. All Rights Reserved

  {[} Info{]} initialize cluster client with ./client.ini

  {[} Info{]} wait for server in 0 / 39 secs
  {[} Info{]} connecting to cluster server 192.168.1.11:7878
  {[} Info{]} wait for license in  2 / 39 secs
  {[} Info{]} wait for license in  4 / 39 secs
  {[} Info{]} wait for license in  6 / 39 secs
  {[} Info{]} wait for license in  8 / 39 secs
  {[} Info{]} wait for license in 10 / 39 secs
  {[} Info{]} wait for license in 20 / 39 secs
  {[} Info{]} wait for license in 30 / 39 secs
  {[}Error{]} timeout at waiting for server approval
  {[}Error{]} Fail to initialize copt command\sphinxhyphen{}line tool
\end{sphinxalltt}

\sphinxAtStartPar
In addition, users can start the COPT command\sphinxhyphen{}line tool on the client side to connect to the specified cluster server IP address as follows:

\begin{sphinxVerbatim}[commandchars=\\\{\}]
\PYG{o}{\PYGZgt{}} \PYG{n}{copt\PYGZus{}cmd} \PYG{o}{\PYGZhy{}}\PYG{n}{cluster} \PYG{o}{\PYGZlt{}}\PYG{n}{ip}\PYG{o}{\PYGZgt{}}
\end{sphinxVerbatim}

\section{COPT Cluster Managing Tool}
\label{\detokenize{cluster:copt-cluster-managing-tool}}
\sphinxAtStartPar
COPT cluster service ships with a tool \sphinxcode{\sphinxupquote{copt\_clstool}}, for retrieving information
and tune parameters of cluster servers on fly.

\subsection{Tool Usage}
\label{\detokenize{cluster:tool-usage}}
\sphinxAtStartPar
Execute the following command in Windows console, Linux or MacOS terminal:

\begin{sphinxVerbatim}[commandchars=\\\{\}]
\PYGZgt{} ./copt\PYGZus{}clstool
\end{sphinxVerbatim}

\sphinxAtStartPar
Below displays help messages of this tool:

\begin{sphinxVerbatim}[commandchars=\\\{\}]
\PYGZgt{} ./copt\PYGZus{}clstool
  COPT Cluster Managing Tool

  copt\PYGZus{}clstool [\PYGZhy{}s server ip] [\PYGZhy{}p port] [\PYGZhy{}x passwd] command \PYGZlt{}param\PYGZgt{}

  commands are:   addblackrule \PYGZlt{}127.0.0.1/20[\PYGZhy{}user@machine]\PYGZgt{}
                  addwhiterule \PYGZlt{}127.0.*.*[+user@machine]\PYGZgt{}
                  getfilters
                  getinfo
                  getnodes
                  getjobs
                  interrupt
                  reload
                  resetfilters
                  setparent \PYGZlt{}xxx:7878\PYGZgt{}
                  setpasswd \PYGZlt{}xxx\PYGZgt{}
                  settoken  \PYGZlt{}num\PYGZgt{}
                  toggleblackrule \PYGZlt{}n\PYGZhy{}th\PYGZgt{}
                  togglewhiterule \PYGZlt{}n\PYGZhy{}th\PYGZgt{}
                  writefilters
\end{sphinxVerbatim}

\sphinxAtStartPar
If the \sphinxcode{\sphinxupquote{\sphinxhyphen{}s}}  and \sphinxcode{\sphinxupquote{\sphinxhyphen{}p}} option are present, tool connects to cluster server
with given server IP and port. Otherwise, tool connections to localhost and
default port 7878. If cluster server sets a password, tool must provide
password string after the \sphinxcode{\sphinxupquote{\sphinxhyphen{}x}} option.

\sphinxAtStartPar
This tool provides the following commands:
\begin{itemize}
\item {} 
\sphinxAtStartPar
\sphinxcode{\sphinxupquote{AddBlackRule}}:
Add a new rule in black filters. each rule has format starting with non\sphinxhyphen{}empty
IP address, which may have wildcard to match IPs in the scope. In addition,
it is optional to be followed by including (+) or excluding (\sphinxhyphen{}) user name at machine name.

\item {} 
\sphinxAtStartPar
\sphinxcode{\sphinxupquote{AddWhiteRule}}:
Add a new rule in white filters. Note that a white rule has the same format as a black rule.

\item {} 
\sphinxAtStartPar
\sphinxcode{\sphinxupquote{GetFilters}}:
Get all rules of black filters, white filters and tool filters, along with relative sequence
numbers, which are parameters for command ToggleBlackRule and ToggleWhiteRule.

\item {} 
\sphinxAtStartPar
\sphinxcode{\sphinxupquote{GetInfo}}:
Get general information of cluster server, including token usage, connected clients,
and all COPT versions in support.

\item {} 
\sphinxAtStartPar
\sphinxcode{\sphinxupquote{GetNodes}}:
Get information of nodes in cluster, including parent address and status, all children nodes.

\item {} 
\sphinxAtStartPar
\sphinxcode{\sphinxupquote{Reload}}:
Reload available token information of all child nodes, in case it is not consistent for various reasons.

\item {} 
\sphinxAtStartPar
\sphinxcode{\sphinxupquote{GetJobs}}:
Retrieves all currently running tasks on the server, including the task ID (TID), runtime (in seconds), and client ID.

\item {} 
\sphinxAtStartPar
\sphinxcode{\sphinxupquote{Interrupt}}:
Terminates the specified task (TID) on the current server (cluster server).
After executing this command, the task running on the client will be stopped,
and the optimization status will return \sphinxcode{\sphinxupquote{stopped (user interrupt)}}.

\item {} 
\sphinxAtStartPar
\sphinxcode{\sphinxupquote{ResetFilters}}:
Reset filter lists in memory to those on filter config file.

\item {} 
\sphinxAtStartPar
\sphinxcode{\sphinxupquote{SetParent}}:
Change parent node address on fly and then connecting to new parent. In this way,
it avoids \sphinxcode{\sphinxupquote{draining}} operation when stopping a node for maintenance purpose.

\item {} 
\sphinxAtStartPar
\sphinxcode{\sphinxupquote{SetPasswd}}:
Update password of target cluster server on fly.

\item {} 
\sphinxAtStartPar
\sphinxcode{\sphinxupquote{SetToken}}:
Change token number of target cluster server on fly.

\item {} 
\sphinxAtStartPar
\sphinxcode{\sphinxupquote{ToggleBlackRule}}:
Toggle between enabling and disabling a black rule, given its sequence number by GetFilters.

\item {} 
\sphinxAtStartPar
\sphinxcode{\sphinxupquote{ToggleWhiteRule}}:
Toggle between enabling and disabling a white rule, given its sequence number by GetFilters.

\item {} 
\sphinxAtStartPar
\sphinxcode{\sphinxupquote{WriteFilters}}:
Write filter lists in memory to filter config file.

\end{itemize}

\subsection{Example Usage}
\label{\detokenize{cluster:id3}}
\sphinxAtStartPar
The following command lists general information on local machine.
\begin{sphinxalltt}
\textgreater{} ./copt\_clstool GetInfo

{[} Info{]} COPT Cluster Managing Tool, COPT v8.0.1 20240304
{[} Info{]} connecting to localhost:7878
{[} Info{]} {[}command{]} wait for connecting to cluster
{[} Info{]} {[}cluster{]} general info
  \# of available tokens is 3 / 3, queue size is 0
  \# of active clients is 0
  \# of installed COPT versions is 1
    COPT v8.0.1
\end{sphinxalltt}

\sphinxAtStartPar
To run managing tool on other machine, its IP should be added to
a rule in \sphinxcode{\sphinxupquote{ToolList}} section in filter configuration file \sphinxcode{\sphinxupquote{clsfilters.ini}}.
The following command from other machine lists cluster information of server 192.168.1.11.
\begin{sphinxalltt}
\textgreater{} ./copt\_clstool \sphinxhyphen{}s 192.168.1.11 GetNodes

{[} Info{]} COPT Cluster Managing Tool, COPT v8.0.1 20240304
{[} Info{]} connecting to 192.168.1.11:7878
{[} Info{]} {[}command{]} wait for connecting to cluster
{[} Info{]} {[}cluster{]} node info
  {[}Parent{]} (null):7878 (Lost)
  {[}Child{]} Node\_192.168.1.12:7878\_N0001, v2.0=3
  Total num of child nodes is 1
\end{sphinxalltt}

\sphinxAtStartPar
The following command changes token number of server 192.168.1.11 from 3 to 0.
\begin{sphinxalltt}
\textgreater{} ./copt\_clstool \sphinxhyphen{}s 192.168.1.11 SetToken 0

{[} Info{]} COPT Cluster Managing Tool, COPT v8.0.1 20240304
{[} Info{]} connecting to 192.168.1.11:7878
{[} Info{]} {[}command{]} wait for connecting to cluster
{[} Info{]} {[}cluster{]} total token was 3 and now set to 0
\end{sphinxalltt}

\sphinxAtStartPar
The following command shows all filter lists of server 192.168.1.11,
including those in BlackList section, WhiteList section and ToolList section.
\begin{sphinxalltt}
\textgreater{} ./copt\_clstool \sphinxhyphen{}s 192.168.1.11 GetFilters

{[} Info{]} COPT Cluster Managing Tool, COPT v8.0.1 20240304
{[} Info{]} connecting to 192.168.1.11:7979
{[} Info{]} {[}command{]} wait for connecting to cluster
{[} Info{]} {[}cluster{]} filters info
{[}BlackList{]}

{[}WhiteList{]}

{[}ToolList{]}
  {[}1{]}  127.0.0.1
\end{sphinxalltt}

\sphinxAtStartPar
The following command added user of IP 192.168.3.13 to black list.
\begin{sphinxalltt}
\textgreater{} ./copt\_clstool \sphinxhyphen{}s 192.168.1.11 AddBlackRule 192.168.3.133

{[} Info{]} COPT Cluster Managing Tool, COPT v8.0.1 20240304
{[} Info{]} connecting to 192.168.1.11:7979
{[} Info{]} {[}command{]} wait for connecting to cluster
{[} Info{]} {[}cluster{]} server added new black rule (succeeded)
\end{sphinxalltt}

\sphinxAtStartPar
The follwing command shows that a new rule in BlackList section is added.
\begin{sphinxalltt}
\textgreater{} ./copt\_clstool \sphinxhyphen{}s 192.168.1.11 GetFilters

{[} Info{]} COPT Cluster Managing Tool, COPT v8.0.1 20240304
{[} Info{]} connecting to 192.168.1.11:7979
{[} Info{]} {[}command{]} wait for connecting to cluster
{[} Info{]} {[}cluster{]} filters info
{[}BlackList{]}
  {[}1{]} 192.168.3.133

{[}WhiteList{]}

{[}ToolList{]}
  {[}1{]}  127.0.0.1
\end{sphinxalltt}

\sphinxAtStartPar
The following command disable a rule in BlackList section.
\begin{sphinxalltt}
\textgreater{} ./copt\_clstool \sphinxhyphen{}s 192.168.1.11 ToggleBlackRule 1

{[} Info{]} COPT Cluster Managing Tool, COPT v8.0.1 20240304
{[} Info{]} connecting to 192.168.1.11:7979
{[} Info{]} {[}command{]} wait for connecting to cluster
{[} Info{]} {[}cluster{]} server toggle black rule {[}1{]} (succeeded)
\end{sphinxalltt}

\section{Running as service}
\label{\detokenize{cluster:running-as-service}}
\sphinxAtStartPar
To run COPT compute cluster server as a system service, follow steps described
in \sphinxcode{\sphinxupquote{readme.txt}} under \sphinxcode{\sphinxupquote{cluster}} folder, and set config file
\sphinxcode{\sphinxupquote{copt\_cluster.service}} properly.

\sphinxAtStartPar
Below is \sphinxcode{\sphinxupquote{readme.txt}}, which lists installing steps in both Linux and MacOS
platforms.

\begin{sphinxVerbatim}[commandchars=\\\{\}]
[Linux] To run copt\PYGZus{}cluster as a service with systemd

Add a systemd file
    cp copt\PYGZus{}cluster.service to /lib/systemd/system/
    sudo systemctl daemon\PYGZhy{}reload

Enable new service
    sudo systemctl start copt\PYGZus{}cluster.service
    or
    sudo systemctl enable copt\PYGZus{}cluster.service

Restart service
    sudo systemctl restart copt\PYGZus{}cluster.service

Stop service
    sudo systemctl stop copt\PYGZus{}cluster.service
    or
    sudo systemctl disable copt\PYGZus{}cluster.service

Verify service is running
    sudo systemctl status copt\PYGZus{}cluster.service

[MacOS] To run copt\PYGZus{}cluster as a service with launchctrl

Add a plist file
    cp copt\PYGZus{}cluster.plist to /Library/LaunchAgents as current user
    or
    cp copt\PYGZus{}cluster.plist to /Library/LaunchDaemons with the key \PYGZsq{}UserName\PYGZsq{}

Enable new service
    sudo launchctl load \PYGZhy{}w /Library/LaunchAgents/copt\PYGZus{}cluster.plist
    or
    sudo launchctl load \PYGZhy{}w /Library/LaunchDaemons/copt\PYGZus{}cluster.plist

Stop service
    sudo launchctl unload \PYGZhy{}w /Library/LaunchAgents/copt\PYGZus{}cluster.plist
    or
    sudo launchctl unload \PYGZhy{}w /Library/LaunchDaemons/copt\PYGZus{}cluster.plist

Verify service is running
    sudo launchctl list solver.copt.cluster
\end{sphinxVerbatim}

\subsection{Linux}
\label{\detokenize{cluster:linux}}
\sphinxAtStartPar
Below are steps in details of how to run COPT compute cluster server as a system service
in Linux platform.

\sphinxAtStartPar
For instance, assume that COPT remote service is installed under \sphinxcode{\sphinxupquote{\textquotesingle{}/home/eleven\textquotesingle{}}}.
In your terminal, type the following command to enter the root directory of
cluster service.
\begin{sphinxalltt}
cd /home/eleven/copt\_remote80/cluster
\end{sphinxalltt}

\sphinxAtStartPar
modify template of the service config file \sphinxcode{\sphinxupquote{copt\_cluster.service}} in text format:

\begin{sphinxVerbatim}[commandchars=\\\{\}]
[Unit]
Description=COPT Compute Cluster Server

[Service]
WorkingDirectory=/path/to/service
ExecStart=/path/to/service/copt\PYGZus{}cluster
Restart=always
RestartSec=1

[Install]
WantedBy=multi\PYGZhy{}user.target
\end{sphinxVerbatim}

\sphinxAtStartPar
That is, update template path in keyword \sphinxcode{\sphinxupquote{WorkingDirectory}} and \sphinxcode{\sphinxupquote{ExecStart}} to
actual path where the cluster service exits.
\begin{sphinxalltt}
{[}Unit{]}
Description=COPT Compute Cluster Server

{[}Service{]}
WorkingDirectory=/home/eleven/copt\_remote80/cluster
ExecStart=/home/eleven/copt\_remote80/cluster/copt\_cluster
Restart=always
RestartSec=1

{[}Install{]}
WantedBy=multi\sphinxhyphen{}user.target
\end{sphinxalltt}

\sphinxAtStartPar
Afterwards, copy \sphinxcode{\sphinxupquote{copt\_cluster.service}} to system service folder
\sphinxcode{\sphinxupquote{/lib/systemd/system/}} (see below).

\begin{sphinxVerbatim}[commandchars=\\\{\}]
sudo cp copt\PYGZus{}cluster.service /lib/systemd/system/
\end{sphinxVerbatim}

\sphinxAtStartPar
The following command may be needed if you add or update service config file.
It is not needed if service unit has been loaded before.

\begin{sphinxVerbatim}[commandchars=\\\{\}]
sudo systemctl daemon\PYGZhy{}reload
\end{sphinxVerbatim}

\sphinxAtStartPar
The following command starts the new cluster service.

\begin{sphinxVerbatim}[commandchars=\\\{\}]
sudo systemctl start copt\PYGZus{}cluster.service
\end{sphinxVerbatim}

\sphinxAtStartPar
To verify the cluster service is actually running, type the following command

\begin{sphinxVerbatim}[commandchars=\\\{\}]
sudo systemctl status copt\PYGZus{}cluster.service
\end{sphinxVerbatim}

\sphinxAtStartPar
If you see logs similar to below, COPT compute cluster server is running successfully as a system service.
\begin{sphinxalltt}
copt\_cluster.service \sphinxhyphen{} COPT Cluster Server
Loaded: loaded (/lib/systemd/system/copt\_cluster.service; enabled; vendor preset: enabled)
Active: active (running) since Sat 2021\sphinxhyphen{}08\sphinxhyphen{}28 11:46:10 CST; 3s ago
Main PID: 3054 (copt\_cluster)
    Tasks: 6 (limit: 4915)
CGroup: /system.slice/copt\_cluster.service
          └─3054 /home/eleven/copt\_remote80/cluster/copt\_cluster

eleven\sphinxhyphen{}ubuntu systemd{[}1{]}: Started COPT Cluster Server.
eleven\sphinxhyphen{}ubuntu COPTCLS{[}3054{]}: LWS: 4.1.4\sphinxhyphen{}b2011a00, loglevel 1039
eleven\sphinxhyphen{}ubuntu COPTCLS{[}3054{]}: NET CLI SRV H1 H2 WS IPv6\sphinxhyphen{}absent
eleven\sphinxhyphen{}ubuntu COPTCLS{[}3054{]}: server started at port 7878
eleven\sphinxhyphen{}ubuntu COPTCLS{[}3054{]}: LWS: 4.1.4\sphinxhyphen{}b2011a00, loglevel 1039
eleven\sphinxhyphen{}ubuntu COPTCLS{[}3054{]}: NET CLI SRV H1 H2 WS IPv6\sphinxhyphen{}absent
eleven\sphinxhyphen{}ubuntu COPTCLS{[}3054{]}: {[}NODE{]} node has been initialized
\end{sphinxalltt}

\sphinxAtStartPar
To stop the cluster service, type the following command

\begin{sphinxVerbatim}[commandchars=\\\{\}]
sudo systemctl stop copt\PYGZus{}cluster.service
\end{sphinxVerbatim}

\subsection{MacOS}
\label{\detokenize{cluster:macos}}
\sphinxAtStartPar
Below are steps in details of how to run COPT Compute Cluster server as a system service
in MacOS platform.

\sphinxAtStartPar
For instance, assume that COPT remote service is installed under \sphinxcode{\sphinxupquote{"/Applications"}}.
In your terminal, type the following command to enter the root directory of
cluster service.
\begin{sphinxalltt}
cd /Applications/copt\_remote80/cluster
\end{sphinxalltt}

\sphinxAtStartPar
modify template of the service config file \sphinxcode{\sphinxupquote{copt\_cluster.plist}} in xml format:

\begin{sphinxVerbatim}[commandchars=\\\{\}]
\PYG{c+cp}{\PYGZlt{}?xml version=\PYGZdq{}1.0\PYGZdq{} encoding=\PYGZdq{}UTF\PYGZhy{}8\PYGZdq{}?\PYGZgt{}}
\PYG{p}{\PYGZlt{}}\PYG{n+nt}{plist} \PYG{n+na}{version}\PYG{o}{=}\PYG{l+s}{\PYGZdq{}1.0\PYGZdq{}}\PYG{p}{\PYGZgt{}}
    \PYG{p}{\PYGZlt{}}\PYG{n+nt}{dict}\PYG{p}{\PYGZgt{}}
        \PYG{p}{\PYGZlt{}}\PYG{n+nt}{key}\PYG{p}{\PYGZgt{}}Label\PYG{p}{\PYGZlt{}}\PYG{p}{/}\PYG{n+nt}{key}\PYG{p}{\PYGZgt{}}
        \PYG{p}{\PYGZlt{}}\PYG{n+nt}{string}\PYG{p}{\PYGZgt{}}solver.copt.cluster\PYG{p}{\PYGZlt{}}\PYG{p}{/}\PYG{n+nt}{string}\PYG{p}{\PYGZgt{}}
        \PYG{p}{\PYGZlt{}}\PYG{n+nt}{key}\PYG{p}{\PYGZgt{}}Program\PYG{p}{\PYGZlt{}}\PYG{p}{/}\PYG{n+nt}{key}\PYG{p}{\PYGZgt{}}
        \PYG{p}{\PYGZlt{}}\PYG{n+nt}{string}\PYG{p}{\PYGZgt{}}/path/to/service/copt\PYGZus{}cluster\PYG{p}{\PYGZlt{}}\PYG{p}{/}\PYG{n+nt}{string}\PYG{p}{\PYGZgt{}}
        \PYG{p}{\PYGZlt{}}\PYG{n+nt}{key}\PYG{p}{\PYGZgt{}}RunAtLoad\PYG{p}{\PYGZlt{}}\PYG{p}{/}\PYG{n+nt}{key}\PYG{p}{\PYGZgt{}}
        \PYG{p}{\PYGZlt{}}\PYG{n+nt}{true}\PYG{p}{/}\PYG{p}{\PYGZgt{}}
        \PYG{p}{\PYGZlt{}}\PYG{n+nt}{key}\PYG{p}{\PYGZgt{}}KeepAlive\PYG{p}{\PYGZlt{}}\PYG{p}{/}\PYG{n+nt}{key}\PYG{p}{\PYGZgt{}}
        \PYG{p}{\PYGZlt{}}\PYG{n+nt}{true}\PYG{p}{/}\PYG{p}{\PYGZgt{}}
    \PYG{p}{\PYGZlt{}}\PYG{p}{/}\PYG{n+nt}{dict}\PYG{p}{\PYGZgt{}}
\PYG{p}{\PYGZlt{}}\PYG{p}{/}\PYG{n+nt}{plist}\PYG{p}{\PYGZgt{}}
\end{sphinxVerbatim}

\sphinxAtStartPar
That is, update template path in \sphinxcode{\sphinxupquote{Program}} tag to
actual path where the cluster service exits.
\begin{sphinxalltt}
\textless{}?xml version="1.0" encoding="UTF\sphinxhyphen{}8"?\textgreater{}
\textless{}plist version="1.0"\textgreater{}
    \textless{}dict\textgreater{}
        \textless{}key\textgreater{}Label\textless{}/key\textgreater{}
        \textless{}string\textgreater{}solver.copt.cluster\textless{}/string\textgreater{}
        \textless{}key\textgreater{}Program\textless{}/key\textgreater{}
        \textless{}string\textgreater{}/Applications/copt\_remote80/cluster/copt\_cluster\textless{}/string\textgreater{}
        \textless{}key\textgreater{}RunAtLoad\textless{}/key\textgreater{}
        \textless{}true/\textgreater{}
        \textless{}key\textgreater{}KeepAlive\textless{}/key\textgreater{}
        \textless{}true/\textgreater{}
    \textless{}/dict\textgreater{}
\textless{}/plist\textgreater{}
\end{sphinxalltt}

\sphinxAtStartPar
Afterwards, copy \sphinxcode{\sphinxupquote{copt\_cluster.plist}} to system service folder
\sphinxcode{\sphinxupquote{/Library/LaunchAgents}} (see below).

\begin{sphinxVerbatim}[commandchars=\\\{\}]
sudo cp copt\PYGZus{}cluster.plist /Library/LaunchAgents
\end{sphinxVerbatim}

\sphinxAtStartPar
The following command starts the new cluster service.

\begin{sphinxVerbatim}[commandchars=\\\{\}]
sudo launchctl load \PYGZhy{}w /Library/LaunchAgents/copt\PYGZus{}cluster.plist
\end{sphinxVerbatim}

\sphinxAtStartPar
To verify the cluster service is actually running, type the following command

\begin{sphinxVerbatim}[commandchars=\\\{\}]
sudo launchctl list solver.copt.cluster
\end{sphinxVerbatim}

\sphinxAtStartPar
If you see logs similar to below, COPT compute cluster server is running successfully as a system service.
\begin{sphinxalltt}
\{
    "LimitLoadToSessionType" = "System";
    "Label" = "solver.copt.cluster";
    "OnDemand" = false;
    "LastExitStatus" = 0;
    "PID" = 16406;
    "Program" = "/Applications/copt\_remote80/cluster/copt\_cluster";
\};
\end{sphinxalltt}

\sphinxAtStartPar
To stop the cluster service, type the following command

\begin{sphinxVerbatim}[commandchars=\\\{\}]
sudo launchctl unload \PYGZhy{}w /Library/LaunchAgents/copt\PYGZus{}cluster.plist
\end{sphinxVerbatim}

\sphinxAtStartPar
If the cluster service should be run by a specific user, add \sphinxcode{\sphinxupquote{UserName}} tag to conifg file.
Below adds a user \sphinxcode{\sphinxupquote{eleven}}, who has priviledge to run the cluster service.
\begin{sphinxalltt}
\textless{}?xml version="1.0" encoding="UTF\sphinxhyphen{}8"?\textgreater{}
\textless{}plist version="1.0"\textgreater{}
    \textless{}dict\textgreater{}
        \textless{}key\textgreater{}Label\textless{}/key\textgreater{}
        \textless{}string\textgreater{}solver.copt.cluster\textless{}/string\textgreater{}
        \textless{}key\textgreater{}Program\textless{}/key\textgreater{}
        \textless{}string\textgreater{}/Applications/copt\_remote80/cluster/copt\_cluster\textless{}/string\textgreater{}
        \textless{}key\textgreater{}UserName\textless{}/key\textgreater{}
        \textless{}string\textgreater{}eleven\textless{}/string\textgreater{}
        \textless{}key\textgreater{}RunAtLoad\textless{}/key\textgreater{}
        \textless{}true/\textgreater{}
        \textless{}key\textgreater{}KeepAlive\textless{}/key\textgreater{}
        \textless{}true/\textgreater{}
    \textless{}/dict\textgreater{}
\textless{}/plist\textgreater{}
\end{sphinxalltt}

\sphinxAtStartPar
Then copy new \sphinxcode{\sphinxupquote{copt\_cluster.plist}} to system service folder
\sphinxcode{\sphinxupquote{/Library/LaunchDaemons}} (see below).

\begin{sphinxVerbatim}[commandchars=\\\{\}]
sudo cp copt\PYGZus{}cluster.plist /Library/LaunchDaemons
\end{sphinxVerbatim}

\sphinxAtStartPar
The following command starts the new cluster service.

\begin{sphinxVerbatim}[commandchars=\\\{\}]
sudo launchctl load \PYGZhy{}w /Library/LaunchDaemons/copt\PYGZus{}cluster.plist
\end{sphinxVerbatim}

\sphinxAtStartPar
To stop the cluster service, type the following command

\begin{sphinxVerbatim}[commandchars=\\\{\}]
sudo launchctl unload \PYGZhy{}w /Library/LaunchDaemons/copt\PYGZus{}cluster.plist
\end{sphinxVerbatim}

\sphinxstepscope

\chapter{COPT Web Licensing Service}
\label{\detokenize{weblicense:copt-web-licensing-service}}\label{\detokenize{weblicense:chapweblicense}}\label{\detokenize{weblicense::doc}}
\sphinxAtStartPar
COPT’s Web Licensing Service provides users with remote licensing services.
Regardless of whether the client is located in the cloud or in a container,
as long as it can access the Internet through the HTTPS protocol, it  can
obtain the Token from COPT’s Web License server to run COPT without binding
any hardware information.

\sphinxAtStartPar
Therefore, compared with traditional authorization methods, Web License is
not limited to fixed hardware environments, flexibly supporting multiple
users, and is suitable for container deployment (such as Docker) and cloud
deployment (Internet environment is required). Its characteristics are
summarized as follows:
\begin{enumerate}
\sphinxsetlistlabels{\arabic}{enumi}{enumii}{}{.}%
\item {} 
\sphinxAtStartPar
It does not bind to any hardware information, flexibly supporting
scenarios such as cloud and container deployment (hardware do not need
to be fixed);

\item {} 
\sphinxAtStartPar
No version restrictions, support cross\sphinxhyphen{}version usage;

\item {} 
\sphinxAtStartPar
At least one machine needs to be connected to the Internet
(remote communication with the web server);

\item {} 
\sphinxAtStartPar
Provide a \sphinxhref{https://copt.shanshu.ai/license}{COPT web license page}
for users to obtain and manage licenses by themselves, which is
user\sphinxhyphen{}friendly and convenient.

\end{enumerate}

\sphinxAtStartPar
At the same time, corresponding to the traditional authorization method,
Web License also includes three subcategories: Web License,
Web License for Floating Server, and Web License for Compute Cluster.

\sphinxAtStartPar
\sphinxstylestrong{Web License}

\sphinxAtStartPar
It supports the deployment of servers running COPT in the cloud development
environment (without binding any machine hardware information), and supports
multiple modeling and solving tasks on the server at the same time.

\sphinxAtStartPar
\sphinxstylestrong{Web License for Floating Server}

\sphinxAtStartPar
It supports deploying a floating token server in a cloud development
environment (the server needs to connect to the Internet and obtain
remote authorization via Web License), and then use the floating token
server to authorize other machines (clients) within the local network
to run COPT.

\sphinxAtStartPar
\sphinxstylestrong{Web License for Compute Cluster}

\sphinxAtStartPar
It supports setting up one or more computing cluster servers in the cloud
development environment. Modeling can be performed on the local machine
(client) within the local network while optimization can be conducted on
the remote cluster server machine (server), so that the powerful computing
resources of the server can be efficiently utilized.

\sphinxAtStartPar
Users need to register on the \sphinxhref{https://copt.shanshu.ai/license}{COPT web license page} first.
After logging in, they can \sphinxstylestrong{apply directly on the page} to obtain the
above three types of Web licenses, create new API Keys, obtain license
files, and manage token occupancy and usage.

\sphinxAtStartPar
For more information on Web License or installation and usage tutorials,
please refer to the \sphinxhref{https://copt.shanshu.ai/license/help}{Web license help documentation}.
If you have any further questions about web licensing, please contact us
as follows:

\begin{savenotes}\sphinxattablestart
\sphinxthistablewithglobalstyle
\centering
\sphinxcapstartof{table}
\sphinxthecaptionisattop
\sphinxcaption{Contact information}\label{\detokenize{weblicense:coptweblicense-contactinfo}}
\sphinxaftertopcaption
\begin{tabulary}{\linewidth}[t]{|T|T|}
\sphinxtoprule
\sphinxtableatstartofbodyhook
\sphinxAtStartPar
\sphinxstylestrong{Email}
&
\sphinxAtStartPar
\sphinxstylestrong{Description}
\\
\sphinxhline
\sphinxAtStartPar
\sphinxhref{mailto:coptsales@shanshu.ai}{coptsales@shanshu.ai}
&
\sphinxAtStartPar
Business Support
\\
\sphinxhline
\sphinxAtStartPar
\sphinxhref{mailto:coptsupport@shanshu.ai}{coptsupport@shanshu.ai}
&
\sphinxAtStartPar
Technical Support
\\
\sphinxbottomrule
\end{tabulary}
\sphinxtableafterendhook\par
\sphinxattableend\end{savenotes}

\sphinxstepscope

\chapter{COPT Quick Start}
\label{\detokenize{quickstart:copt-quick-start}}\label{\detokenize{quickstart:chapquickstart}}\label{\detokenize{quickstart::doc}}
\sphinxstepscope

\section{C Interface}
\label{\detokenize{cinterface:c-interface}}\label{\detokenize{cinterface:chapcinterface}}\label{\detokenize{cinterface::doc}}
\sphinxAtStartPar
This chapter illustrate the use of C interface of Cardinal Optimizer through
a simple C example. The problem to solve is shown in Eq. \ref{equation:cinterface:coptEq_capilp1}:
\begin{equation}\label{equation:cinterface:coptEq_capilp1}
\begin{split}\text{Maximize: } & \\
                  & 1.2 x + 1.8 y + 2.1 z \\
\text{Subject to: } & \\
                    & 1.5 x + 1.2 y + 1.8 z \leq 2.6 \\
                    & 0.8 x + 0.6 y + 0.9 z \geq 1.2 \\
\text{Bounds: } & \\
               & 0.1 \leq x \leq 0.6 \\
               & 0.2 \leq y \leq 1.5 \\
               & 0.3 \leq z \leq 2.8\end{split}
\end{equation}

\subsection{Example details}
\label{\detokenize{cinterface:example-details}}
\sphinxAtStartPar
The source code for solving the above problem using C API of Cardinal Optimizer
is shown in \hyperref[\detokenize{cinterface:coptcode-capilp1}]{Listing \ref{\detokenize{cinterface:coptcode-capilp1}}}:
\sphinxSetupCaptionForVerbatim{\sphinxcode{\sphinxupquote{lp\_ex1.c}}}
\def\sphinxLiteralBlockLabel{\label{\detokenize{cinterface:coptcode-capilp1}}}
\begin{sphinxVerbatim}[commandchars=\\\{\},numbers=left,firstnumber=1,stepnumber=1]
\PYG{c+cm}{/*}
\PYG{c+cm}{ * This file is part of the Cardinal Optimizer, all rights reserved.}
\PYG{c+cm}{ */}

\PYG{c+cm}{/*}
\PYG{c+cm}{ * The problem to solve:}
\PYG{c+cm}{ *}
\PYG{c+cm}{ *  Maximize:}
\PYG{c+cm}{ *    1.2 x + 1.8 y + 2.1 z}
\PYG{c+cm}{ *}
\PYG{c+cm}{ *  Subject to:}
\PYG{c+cm}{ *    1.5 x + 1.2 y + 1.8 z \PYGZlt{}= 2.6}
\PYG{c+cm}{ *    0.8 x + 0.6 y + 0.9 z \PYGZgt{}= 1.2}
\PYG{c+cm}{ *}
\PYG{c+cm}{ *  where:}
\PYG{c+cm}{ *    0.1 \PYGZlt{}= x \PYGZlt{}= 0.6}
\PYG{c+cm}{ *    0.2 \PYGZlt{}= y \PYGZlt{}= 1.5}
\PYG{c+cm}{ *    0.3 \PYGZlt{}= z \PYGZlt{}= 2.8}
\PYG{c+cm}{ */}

\PYG{c+cp}{\PYGZsh{}}\PYG{c+cp}{include}\PYG{+w}{ }\PYG{c+cpf}{\PYGZdq{}copt.h\PYGZdq{}}

\PYG{c+cp}{\PYGZsh{}}\PYG{c+cp}{include}\PYG{+w}{ }\PYG{c+cpf}{\PYGZlt{}stdio.h\PYGZgt{}}
\PYG{c+cp}{\PYGZsh{}}\PYG{c+cp}{include}\PYG{+w}{ }\PYG{c+cpf}{\PYGZlt{}stdlib.h\PYGZgt{}}

\PYG{k+kt}{int}\PYG{+w}{ }\PYG{n+nf}{main}\PYG{p}{(}\PYG{k+kt}{int}\PYG{+w}{ }\PYG{n}{argc}\PYG{p}{,}\PYG{+w}{ }\PYG{k+kt}{char}\PYG{o}{*}\PYG{+w}{ }\PYG{n}{argv}\PYG{p}{[}\PYG{p}{]}\PYG{p}{)}
\PYG{p}{\PYGZob{}}
\PYG{+w}{  }\PYG{k+kt}{int}\PYG{+w}{ }\PYG{n}{errcode}\PYG{+w}{ }\PYG{o}{=}\PYG{+w}{ }\PYG{l+m+mi}{0}\PYG{p}{;}

\PYG{+w}{  }\PYG{n}{copt\PYGZus{}env}\PYG{o}{*}\PYG{+w}{ }\PYG{n}{env}\PYG{+w}{ }\PYG{o}{=}\PYG{+w}{ }\PYG{n+nb}{NULL}\PYG{p}{;}
\PYG{+w}{  }\PYG{n}{copt\PYGZus{}prob}\PYG{o}{*}\PYG{+w}{ }\PYG{n}{prob}\PYG{+w}{ }\PYG{o}{=}\PYG{+w}{ }\PYG{n+nb}{NULL}\PYG{p}{;}

\PYG{+w}{  }\PYG{c+c1}{// Create COPT environment}
\PYG{+w}{  }\PYG{n}{errcode}\PYG{+w}{ }\PYG{o}{=}\PYG{+w}{ }\PYG{n}{COPT\PYGZus{}CreateEnv}\PYG{p}{(}\PYG{o}{\PYGZam{}}\PYG{n}{env}\PYG{p}{)}\PYG{p}{;}
\PYG{+w}{  }\PYG{k}{if}\PYG{+w}{ }\PYG{p}{(}\PYG{n}{errcode}\PYG{p}{)}
\PYG{+w}{    }\PYG{k}{goto}\PYG{+w}{ }\PYG{n}{COPT\PYGZus{}EXIT}\PYG{p}{;}

\PYG{+w}{  }\PYG{c+c1}{// Create COPT problem}
\PYG{+w}{  }\PYG{n}{errcode}\PYG{+w}{ }\PYG{o}{=}\PYG{+w}{ }\PYG{n}{COPT\PYGZus{}CreateProb}\PYG{p}{(}\PYG{n}{env}\PYG{p}{,}\PYG{+w}{ }\PYG{o}{\PYGZam{}}\PYG{n}{prob}\PYG{p}{)}\PYG{p}{;}
\PYG{+w}{  }\PYG{k}{if}\PYG{+w}{ }\PYG{p}{(}\PYG{n}{errcode}\PYG{p}{)}
\PYG{+w}{    }\PYG{k}{goto}\PYG{+w}{ }\PYG{n}{COPT\PYGZus{}EXIT}\PYG{p}{;}
\PYG{+w}{  }
\PYG{+w}{  }\PYG{c+cm}{/*}
\PYG{c+cm}{   * Add variables}
\PYG{c+cm}{   *}
\PYG{c+cm}{   *   obj: 1.2 C0 + 1.8 C1 + 2.1 C2}
\PYG{c+cm}{   *}
\PYG{c+cm}{   *   var:}
\PYG{c+cm}{   *     0.1 \PYGZlt{}= C0 \PYGZlt{}= 0.6}
\PYG{c+cm}{   *     0.2 \PYGZlt{}= C1 \PYGZlt{}= 1.5}
\PYG{c+cm}{   *     0.3 \PYGZlt{}= C2 \PYGZlt{}= 2.8}
\PYG{c+cm}{   *}
\PYG{c+cm}{   */}
\PYG{+w}{  }\PYG{k+kt}{int}\PYG{+w}{ }\PYG{n}{ncol}\PYG{+w}{ }\PYG{o}{=}\PYG{+w}{ }\PYG{l+m+mi}{3}\PYG{p}{;}
\PYG{+w}{  }\PYG{k+kt}{double}\PYG{+w}{ }\PYG{n}{colcost}\PYG{p}{[}\PYG{p}{]}\PYG{+w}{ }\PYG{o}{=}\PYG{+w}{ }\PYG{p}{\PYGZob{}}\PYG{l+m+mf}{1.2}\PYG{p}{,}\PYG{+w}{ }\PYG{l+m+mf}{1.8}\PYG{p}{,}\PYG{+w}{ }\PYG{l+m+mf}{2.1}\PYG{p}{\PYGZcb{}}\PYG{p}{;}
\PYG{+w}{  }\PYG{k+kt}{double}\PYG{+w}{ }\PYG{n}{collb}\PYG{p}{[}\PYG{p}{]}\PYG{+w}{ }\PYG{o}{=}\PYG{+w}{ }\PYG{p}{\PYGZob{}}\PYG{l+m+mf}{0.1}\PYG{p}{,}\PYG{+w}{ }\PYG{l+m+mf}{0.2}\PYG{p}{,}\PYG{+w}{ }\PYG{l+m+mf}{0.3}\PYG{p}{\PYGZcb{}}\PYG{p}{;}
\PYG{+w}{  }\PYG{k+kt}{double}\PYG{+w}{ }\PYG{n}{colub}\PYG{p}{[}\PYG{p}{]}\PYG{+w}{ }\PYG{o}{=}\PYG{+w}{ }\PYG{p}{\PYGZob{}}\PYG{l+m+mf}{0.6}\PYG{p}{,}\PYG{+w}{ }\PYG{l+m+mf}{1.5}\PYG{p}{,}\PYG{+w}{ }\PYG{l+m+mf}{1.8}\PYG{p}{\PYGZcb{}}\PYG{p}{;}

\PYG{+w}{  }\PYG{n}{errcode}\PYG{+w}{ }\PYG{o}{=}\PYG{+w}{ }\PYG{n}{COPT\PYGZus{}AddCols}\PYG{p}{(}\PYG{n}{prob}\PYG{p}{,}\PYG{+w}{ }\PYG{n}{ncol}\PYG{p}{,}\PYG{+w}{ }\PYG{n}{colcost}\PYG{p}{,}\PYG{+w}{ }\PYG{n+nb}{NULL}\PYG{p}{,}\PYG{+w}{ }\PYG{n+nb}{NULL}\PYG{p}{,}\PYG{+w}{ }\PYG{n+nb}{NULL}\PYG{p}{,}\PYG{+w}{ }\PYG{n+nb}{NULL}\PYG{p}{,}\PYG{+w}{ }\PYG{n+nb}{NULL}\PYG{p}{,}\PYG{+w}{ }\PYG{n}{collb}\PYG{p}{,}\PYG{+w}{ }\PYG{n}{colub}\PYG{p}{,}\PYG{+w}{ }\PYG{n+nb}{NULL}\PYG{p}{)}\PYG{p}{;}
\PYG{+w}{  }\PYG{k}{if}\PYG{+w}{ }\PYG{p}{(}\PYG{n}{errcode}\PYG{p}{)}
\PYG{+w}{    }\PYG{k}{goto}\PYG{+w}{ }\PYG{n}{COPT\PYGZus{}EXIT}\PYG{p}{;}

\PYG{+w}{  }\PYG{c+cm}{/*}
\PYG{c+cm}{   * Add constraints}
\PYG{c+cm}{   *}
\PYG{c+cm}{   *   r0: 1.5 C0 + 1.2 C1 + 1.8 C2 \PYGZlt{}= 2.6}
\PYG{c+cm}{   *   r1: 0.8 C0 + 0.6 C1 + 0.9 C2 \PYGZgt{}= 1.2}
\PYG{c+cm}{   */}
\PYG{+w}{  }\PYG{k+kt}{int}\PYG{+w}{ }\PYG{n}{nrow}\PYG{+w}{ }\PYG{o}{=}\PYG{+w}{ }\PYG{l+m+mi}{2}\PYG{p}{;}
\PYG{+w}{  }\PYG{k+kt}{int}\PYG{+w}{ }\PYG{n}{rowbeg}\PYG{p}{[}\PYG{p}{]}\PYG{+w}{ }\PYG{o}{=}\PYG{+w}{ }\PYG{p}{\PYGZob{}}\PYG{l+m+mi}{0}\PYG{p}{,}\PYG{+w}{ }\PYG{l+m+mi}{3}\PYG{p}{\PYGZcb{}}\PYG{p}{;}
\PYG{+w}{  }\PYG{k+kt}{int}\PYG{+w}{ }\PYG{n}{rowcnt}\PYG{p}{[}\PYG{p}{]}\PYG{+w}{ }\PYG{o}{=}\PYG{+w}{ }\PYG{p}{\PYGZob{}}\PYG{l+m+mi}{3}\PYG{p}{,}\PYG{+w}{ }\PYG{l+m+mi}{3}\PYG{p}{\PYGZcb{}}\PYG{p}{;}
\PYG{+w}{  }\PYG{k+kt}{int}\PYG{+w}{ }\PYG{n}{rowind}\PYG{p}{[}\PYG{p}{]}\PYG{+w}{ }\PYG{o}{=}\PYG{+w}{ }\PYG{p}{\PYGZob{}}\PYG{l+m+mi}{0}\PYG{p}{,}\PYG{+w}{ }\PYG{l+m+mi}{1}\PYG{p}{,}\PYG{+w}{ }\PYG{l+m+mi}{2}\PYG{p}{,}\PYG{+w}{ }\PYG{l+m+mi}{0}\PYG{p}{,}\PYG{+w}{ }\PYG{l+m+mi}{1}\PYG{p}{,}\PYG{+w}{ }\PYG{l+m+mi}{2}\PYG{p}{\PYGZcb{}}\PYG{p}{;}
\PYG{+w}{  }\PYG{k+kt}{double}\PYG{+w}{ }\PYG{n}{rowelem}\PYG{p}{[}\PYG{p}{]}\PYG{+w}{ }\PYG{o}{=}\PYG{+w}{ }\PYG{p}{\PYGZob{}}\PYG{l+m+mf}{1.5}\PYG{p}{,}\PYG{+w}{ }\PYG{l+m+mf}{1.2}\PYG{p}{,}\PYG{+w}{ }\PYG{l+m+mf}{1.8}\PYG{p}{,}\PYG{+w}{ }\PYG{l+m+mf}{0.8}\PYG{p}{,}\PYG{+w}{ }\PYG{l+m+mf}{0.6}\PYG{p}{,}\PYG{+w}{ }\PYG{l+m+mf}{0.9}\PYG{p}{\PYGZcb{}}\PYG{p}{;}
\PYG{+w}{  }\PYG{k+kt}{char}\PYG{+w}{ }\PYG{n}{rowsen}\PYG{p}{[}\PYG{p}{]}\PYG{+w}{ }\PYG{o}{=}\PYG{+w}{ }\PYG{p}{\PYGZob{}}\PYG{n}{COPT\PYGZus{}LESS\PYGZus{}EQUAL}\PYG{p}{,}\PYG{+w}{ }\PYG{n}{COPT\PYGZus{}GREATER\PYGZus{}EQUAL}\PYG{p}{\PYGZcb{}}\PYG{p}{;}
\PYG{+w}{  }\PYG{k+kt}{double}\PYG{+w}{ }\PYG{n}{rowrhs}\PYG{p}{[}\PYG{p}{]}\PYG{+w}{ }\PYG{o}{=}\PYG{+w}{ }\PYG{p}{\PYGZob{}}\PYG{l+m+mf}{2.6}\PYG{p}{,}\PYG{+w}{ }\PYG{l+m+mf}{1.2}\PYG{p}{\PYGZcb{}}\PYG{p}{;}

\PYG{+w}{  }\PYG{n}{errcode}\PYG{+w}{ }\PYG{o}{=}\PYG{+w}{ }\PYG{n}{COPT\PYGZus{}AddRows}\PYG{p}{(}\PYG{n}{prob}\PYG{p}{,}\PYG{+w}{ }\PYG{n}{nrow}\PYG{p}{,}\PYG{+w}{ }\PYG{n}{rowbeg}\PYG{p}{,}\PYG{+w}{ }\PYG{n}{rowcnt}\PYG{p}{,}\PYG{+w}{ }\PYG{n}{rowind}\PYG{p}{,}\PYG{+w}{ }\PYG{n}{rowelem}\PYG{p}{,}\PYG{+w}{ }\PYG{n}{rowsen}\PYG{p}{,}\PYG{+w}{ }\PYG{n}{rowrhs}\PYG{p}{,}\PYG{+w}{ }\PYG{n+nb}{NULL}\PYG{p}{,}\PYG{+w}{ }\PYG{n+nb}{NULL}\PYG{p}{)}\PYG{p}{;}
\PYG{+w}{  }\PYG{k}{if}\PYG{+w}{ }\PYG{p}{(}\PYG{n}{errcode}\PYG{p}{)}
\PYG{+w}{    }\PYG{k}{goto}\PYG{+w}{ }\PYG{n}{COPT\PYGZus{}EXIT}\PYG{p}{;}

\PYG{+w}{  }\PYG{c+c1}{// Set parameters and attributes}
\PYG{+w}{  }\PYG{n}{errcode}\PYG{+w}{ }\PYG{o}{=}\PYG{+w}{ }\PYG{n}{COPT\PYGZus{}SetDblParam}\PYG{p}{(}\PYG{n}{prob}\PYG{p}{,}\PYG{+w}{ }\PYG{n}{COPT\PYGZus{}DBLPARAM\PYGZus{}TIMELIMIT}\PYG{p}{,}\PYG{+w}{ }\PYG{l+m+mi}{10}\PYG{p}{)}\PYG{p}{;}
\PYG{+w}{  }\PYG{k}{if}\PYG{+w}{ }\PYG{p}{(}\PYG{n}{errcode}\PYG{p}{)}
\PYG{+w}{    }\PYG{k}{goto}\PYG{+w}{ }\PYG{n}{COPT\PYGZus{}EXIT}\PYG{p}{;}
\PYG{+w}{  }\PYG{n}{errcode}\PYG{+w}{ }\PYG{o}{=}\PYG{+w}{ }\PYG{n}{COPT\PYGZus{}SetObjSense}\PYG{p}{(}\PYG{n}{prob}\PYG{p}{,}\PYG{+w}{ }\PYG{n}{COPT\PYGZus{}MAXIMIZE}\PYG{p}{)}\PYG{p}{;}
\PYG{+w}{  }\PYG{k}{if}\PYG{+w}{ }\PYG{p}{(}\PYG{n}{errcode}\PYG{p}{)}
\PYG{+w}{    }\PYG{k}{goto}\PYG{+w}{ }\PYG{n}{COPT\PYGZus{}EXIT}\PYG{p}{;}

\PYG{+w}{  }\PYG{c+c1}{// Solve problem}
\PYG{+w}{  }\PYG{n}{errcode}\PYG{+w}{ }\PYG{o}{=}\PYG{+w}{ }\PYG{n}{COPT\PYGZus{}SolveLp}\PYG{p}{(}\PYG{n}{prob}\PYG{p}{)}\PYG{p}{;}
\PYG{+w}{  }\PYG{k}{if}\PYG{+w}{ }\PYG{p}{(}\PYG{n}{errcode}\PYG{p}{)}
\PYG{+w}{    }\PYG{k}{goto}\PYG{+w}{ }\PYG{n}{COPT\PYGZus{}EXIT}\PYG{p}{;}

\PYG{+w}{  }\PYG{c+c1}{// Analyze solution}
\PYG{+w}{  }\PYG{k+kt}{int}\PYG{+w}{ }\PYG{n}{lpstat}\PYG{+w}{ }\PYG{o}{=}\PYG{+w}{ }\PYG{n}{COPT\PYGZus{}LPSTATUS\PYGZus{}UNSTARTED}\PYG{p}{;}
\PYG{+w}{  }\PYG{k+kt}{double}\PYG{+w}{ }\PYG{n}{lpobjval}\PYG{p}{;}
\PYG{+w}{  }\PYG{k+kt}{double}\PYG{o}{*}\PYG{+w}{ }\PYG{n}{lpsol}\PYG{+w}{ }\PYG{o}{=}\PYG{+w}{ }\PYG{n+nb}{NULL}\PYG{p}{;}
\PYG{+w}{  }\PYG{k+kt}{int}\PYG{o}{*}\PYG{+w}{ }\PYG{n}{colstat}\PYG{+w}{ }\PYG{o}{=}\PYG{+w}{ }\PYG{n+nb}{NULL}\PYG{p}{;}

\PYG{+w}{  }\PYG{n}{errcode}\PYG{+w}{ }\PYG{o}{=}\PYG{+w}{ }\PYG{n}{COPT\PYGZus{}GetIntAttr}\PYG{p}{(}\PYG{n}{prob}\PYG{p}{,}\PYG{+w}{ }\PYG{n}{COPT\PYGZus{}INTATTR\PYGZus{}LPSTATUS}\PYG{p}{,}\PYG{+w}{ }\PYG{o}{\PYGZam{}}\PYG{n}{lpstat}\PYG{p}{)}\PYG{p}{;}
\PYG{+w}{  }\PYG{k}{if}\PYG{+w}{ }\PYG{p}{(}\PYG{n}{errcode}\PYG{p}{)}
\PYG{+w}{    }\PYG{k}{goto}\PYG{+w}{ }\PYG{n}{COPT\PYGZus{}EXIT}\PYG{p}{;}

\PYG{+w}{  }\PYG{k}{if}\PYG{+w}{ }\PYG{p}{(}\PYG{n}{lpstat}\PYG{+w}{ }\PYG{o}{=}\PYG{o}{=}\PYG{+w}{ }\PYG{n}{COPT\PYGZus{}LPSTATUS\PYGZus{}OPTIMAL}\PYG{p}{)}
\PYG{+w}{  }\PYG{p}{\PYGZob{}}
\PYG{+w}{    }\PYG{n}{lpsol}\PYG{+w}{ }\PYG{o}{=}\PYG{+w}{ }\PYG{p}{(}\PYG{k+kt}{double}\PYG{o}{*}\PYG{p}{)}\PYG{n}{malloc}\PYG{p}{(}\PYG{n}{ncol}\PYG{+w}{ }\PYG{o}{*}\PYG{+w}{ }\PYG{k}{sizeof}\PYG{p}{(}\PYG{k+kt}{double}\PYG{p}{)}\PYG{p}{)}\PYG{p}{;}
\PYG{+w}{    }\PYG{n}{colstat}\PYG{+w}{ }\PYG{o}{=}\PYG{+w}{ }\PYG{p}{(}\PYG{k+kt}{int}\PYG{o}{*}\PYG{p}{)}\PYG{n}{malloc}\PYG{p}{(}\PYG{n}{ncol}\PYG{+w}{ }\PYG{o}{*}\PYG{+w}{ }\PYG{k}{sizeof}\PYG{p}{(}\PYG{k+kt}{int}\PYG{p}{)}\PYG{p}{)}\PYG{p}{;}

\PYG{+w}{    }\PYG{n}{errcode}\PYG{+w}{ }\PYG{o}{=}\PYG{+w}{ }\PYG{n}{COPT\PYGZus{}GetLpSolution}\PYG{p}{(}\PYG{n}{prob}\PYG{p}{,}\PYG{+w}{ }\PYG{n}{lpsol}\PYG{p}{,}\PYG{+w}{ }\PYG{n+nb}{NULL}\PYG{p}{,}\PYG{+w}{ }\PYG{n+nb}{NULL}\PYG{p}{,}\PYG{+w}{ }\PYG{n+nb}{NULL}\PYG{p}{)}\PYG{p}{;}
\PYG{+w}{    }\PYG{k}{if}\PYG{+w}{ }\PYG{p}{(}\PYG{n}{errcode}\PYG{p}{)}
\PYG{+w}{      }\PYG{k}{goto}\PYG{+w}{ }\PYG{n}{COPT\PYGZus{}EXIT}\PYG{p}{;}

\PYG{+w}{    }\PYG{n}{errcode}\PYG{+w}{ }\PYG{o}{=}\PYG{+w}{ }\PYG{n}{COPT\PYGZus{}GetBasis}\PYG{p}{(}\PYG{n}{prob}\PYG{p}{,}\PYG{+w}{ }\PYG{n}{colstat}\PYG{p}{,}\PYG{+w}{ }\PYG{n+nb}{NULL}\PYG{p}{)}\PYG{p}{;}
\PYG{+w}{    }\PYG{k}{if}\PYG{+w}{ }\PYG{p}{(}\PYG{n}{errcode}\PYG{p}{)}
\PYG{+w}{      }\PYG{k}{goto}\PYG{+w}{ }\PYG{n}{COPT\PYGZus{}EXIT}\PYG{p}{;}

\PYG{+w}{    }\PYG{n}{errcode}\PYG{+w}{ }\PYG{o}{=}\PYG{+w}{ }\PYG{n}{COPT\PYGZus{}GetDblAttr}\PYG{p}{(}\PYG{n}{prob}\PYG{p}{,}\PYG{+w}{ }\PYG{n}{COPT\PYGZus{}DBLATTR\PYGZus{}LPOBJVAL}\PYG{p}{,}\PYG{+w}{ }\PYG{o}{\PYGZam{}}\PYG{n}{lpobjval}\PYG{p}{)}\PYG{p}{;}
\PYG{+w}{    }\PYG{k}{if}\PYG{+w}{ }\PYG{p}{(}\PYG{n}{errcode}\PYG{p}{)}
\PYG{+w}{      }\PYG{k}{goto}\PYG{+w}{ }\PYG{n}{COPT\PYGZus{}EXIT}\PYG{p}{;}

\PYG{+w}{    }\PYG{n}{printf}\PYG{p}{(}\PYG{l+s}{\PYGZdq{}}\PYG{l+s+se}{\PYGZbs{}n}\PYG{l+s}{Objective value: \PYGZpc{}.6f}\PYG{l+s+se}{\PYGZbs{}n}\PYG{l+s}{\PYGZdq{}}\PYG{p}{,}\PYG{+w}{ }\PYG{n}{lpobjval}\PYG{p}{)}\PYG{p}{;}

\PYG{+w}{    }\PYG{n}{printf}\PYG{p}{(}\PYG{l+s}{\PYGZdq{}}\PYG{l+s}{Variable solution: }\PYG{l+s+se}{\PYGZbs{}n}\PYG{l+s}{\PYGZdq{}}\PYG{p}{)}\PYG{p}{;}
\PYG{+w}{    }\PYG{k}{for}\PYG{+w}{ }\PYG{p}{(}\PYG{k+kt}{int}\PYG{+w}{ }\PYG{n}{i}\PYG{+w}{ }\PYG{o}{=}\PYG{+w}{ }\PYG{l+m+mi}{0}\PYG{p}{;}\PYG{+w}{ }\PYG{n}{i}\PYG{+w}{ }\PYG{o}{\PYGZlt{}}\PYG{+w}{ }\PYG{n}{ncol}\PYG{p}{;}\PYG{+w}{ }\PYG{o}{+}\PYG{o}{+}\PYG{n}{i}\PYG{p}{)}
\PYG{+w}{      }\PYG{n}{printf}\PYG{p}{(}\PYG{l+s}{\PYGZdq{}}\PYG{l+s}{  x[\PYGZpc{}d] = \PYGZpc{}.6f}\PYG{l+s+se}{\PYGZbs{}n}\PYG{l+s}{\PYGZdq{}}\PYG{p}{,}\PYG{+w}{ }\PYG{n}{i}\PYG{p}{,}\PYG{+w}{ }\PYG{n}{lpsol}\PYG{p}{[}\PYG{n}{i}\PYG{p}{]}\PYG{p}{)}\PYG{p}{;}

\PYG{+w}{    }\PYG{n}{printf}\PYG{p}{(}\PYG{l+s}{\PYGZdq{}}\PYG{l+s}{Variable basis status: }\PYG{l+s+se}{\PYGZbs{}n}\PYG{l+s}{\PYGZdq{}}\PYG{p}{)}\PYG{p}{;}
\PYG{+w}{    }\PYG{k}{for}\PYG{+w}{ }\PYG{p}{(}\PYG{k+kt}{int}\PYG{+w}{ }\PYG{n}{i}\PYG{+w}{ }\PYG{o}{=}\PYG{+w}{ }\PYG{l+m+mi}{0}\PYG{p}{;}\PYG{+w}{ }\PYG{n}{i}\PYG{+w}{ }\PYG{o}{\PYGZlt{}}\PYG{+w}{ }\PYG{n}{ncol}\PYG{p}{;}\PYG{+w}{ }\PYG{o}{+}\PYG{o}{+}\PYG{n}{i}\PYG{p}{)}
\PYG{+w}{      }\PYG{n}{printf}\PYG{p}{(}\PYG{l+s}{\PYGZdq{}}\PYG{l+s}{  x[\PYGZpc{}d]: \PYGZpc{}d}\PYG{l+s+se}{\PYGZbs{}n}\PYG{l+s}{\PYGZdq{}}\PYG{p}{,}\PYG{+w}{ }\PYG{n}{i}\PYG{p}{,}\PYG{+w}{ }\PYG{n}{colstat}\PYG{p}{[}\PYG{n}{i}\PYG{p}{]}\PYG{p}{)}\PYG{p}{;}

\PYG{+w}{    }\PYG{n}{free}\PYG{p}{(}\PYG{n}{lpsol}\PYG{p}{)}\PYG{p}{;}
\PYG{+w}{    }\PYG{n}{free}\PYG{p}{(}\PYG{n}{colstat}\PYG{p}{)}\PYG{p}{;}
\PYG{+w}{  }\PYG{p}{\PYGZcb{}}

\PYG{+w}{  }\PYG{c+c1}{// Write problem, solution and modified parameters to files}
\PYG{+w}{  }\PYG{n}{errcode}\PYG{+w}{ }\PYG{o}{=}\PYG{+w}{ }\PYG{n}{COPT\PYGZus{}WriteMps}\PYG{p}{(}\PYG{n}{prob}\PYG{p}{,}\PYG{+w}{ }\PYG{l+s}{\PYGZdq{}}\PYG{l+s}{lp\PYGZus{}ex1.mps}\PYG{l+s}{\PYGZdq{}}\PYG{p}{)}\PYG{p}{;}
\PYG{+w}{  }\PYG{k}{if}\PYG{+w}{ }\PYG{p}{(}\PYG{n}{errcode}\PYG{p}{)}
\PYG{+w}{    }\PYG{k}{goto}\PYG{+w}{ }\PYG{n}{COPT\PYGZus{}EXIT}\PYG{p}{;}
\PYG{+w}{  }\PYG{n}{errcode}\PYG{+w}{ }\PYG{o}{=}\PYG{+w}{ }\PYG{n}{COPT\PYGZus{}WriteBasis}\PYG{p}{(}\PYG{n}{prob}\PYG{p}{,}\PYG{+w}{ }\PYG{l+s}{\PYGZdq{}}\PYG{l+s}{lp\PYGZus{}ex1.bas}\PYG{l+s}{\PYGZdq{}}\PYG{p}{)}\PYG{p}{;}
\PYG{+w}{  }\PYG{k}{if}\PYG{+w}{ }\PYG{p}{(}\PYG{n}{errcode}\PYG{p}{)}
\PYG{+w}{    }\PYG{k}{goto}\PYG{+w}{ }\PYG{n}{COPT\PYGZus{}EXIT}\PYG{p}{;}
\PYG{+w}{  }\PYG{n}{errcode}\PYG{+w}{ }\PYG{o}{=}\PYG{+w}{ }\PYG{n}{COPT\PYGZus{}WriteSol}\PYG{p}{(}\PYG{n}{prob}\PYG{p}{,}\PYG{+w}{ }\PYG{l+s}{\PYGZdq{}}\PYG{l+s}{lp\PYGZus{}ex1.sol}\PYG{l+s}{\PYGZdq{}}\PYG{p}{)}\PYG{p}{;}
\PYG{+w}{  }\PYG{k}{if}\PYG{+w}{ }\PYG{p}{(}\PYG{n}{errcode}\PYG{p}{)}
\PYG{+w}{    }\PYG{k}{goto}\PYG{+w}{ }\PYG{n}{COPT\PYGZus{}EXIT}\PYG{p}{;}
\PYG{+w}{  }\PYG{n}{errcode}\PYG{+w}{ }\PYG{o}{=}\PYG{+w}{ }\PYG{n}{COPT\PYGZus{}WriteParam}\PYG{p}{(}\PYG{n}{prob}\PYG{p}{,}\PYG{+w}{ }\PYG{l+s}{\PYGZdq{}}\PYG{l+s}{lp\PYGZus{}ex1.par}\PYG{l+s}{\PYGZdq{}}\PYG{p}{)}\PYG{p}{;}
\PYG{+w}{  }\PYG{k}{if}\PYG{+w}{ }\PYG{p}{(}\PYG{n}{errcode}\PYG{p}{)}
\PYG{+w}{    }\PYG{k}{goto}\PYG{+w}{ }\PYG{n}{COPT\PYGZus{}EXIT}\PYG{p}{;}

\PYG{+w}{  }\PYG{c+c1}{// Error handling}
\PYG{n+nl}{COPT\PYGZus{}EXIT}\PYG{p}{:}
\PYG{+w}{  }\PYG{k}{if}\PYG{+w}{ }\PYG{p}{(}\PYG{n}{errcode}\PYG{p}{)}
\PYG{+w}{  }\PYG{p}{\PYGZob{}}
\PYG{+w}{    }\PYG{k+kt}{char}\PYG{+w}{ }\PYG{n}{errmsg}\PYG{p}{[}\PYG{n}{COPT\PYGZus{}BUFFSIZE}\PYG{p}{]}\PYG{p}{;}

\PYG{+w}{    }\PYG{n}{COPT\PYGZus{}GetRetcodeMsg}\PYG{p}{(}\PYG{n}{errcode}\PYG{p}{,}\PYG{+w}{ }\PYG{n}{errmsg}\PYG{p}{,}\PYG{+w}{ }\PYG{n}{COPT\PYGZus{}BUFFSIZE}\PYG{p}{)}\PYG{p}{;}
\PYG{+w}{    }\PYG{n}{printf}\PYG{p}{(}\PYG{l+s}{\PYGZdq{}}\PYG{l+s}{ERROR \PYGZpc{}d: \PYGZpc{}s}\PYG{l+s+se}{\PYGZbs{}n}\PYG{l+s}{\PYGZdq{}}\PYG{p}{,}\PYG{+w}{ }\PYG{n}{errcode}\PYG{p}{,}\PYG{+w}{ }\PYG{n}{errmsg}\PYG{p}{)}\PYG{p}{;}

\PYG{+w}{    }\PYG{k}{return}\PYG{+w}{ }\PYG{l+m+mi}{0}\PYG{p}{;}
\PYG{+w}{  }\PYG{p}{\PYGZcb{}}

\PYG{+w}{  }\PYG{c+c1}{// Delete problem and environment}
\PYG{+w}{  }\PYG{n}{COPT\PYGZus{}DeleteProb}\PYG{p}{(}\PYG{o}{\PYGZam{}}\PYG{n}{prob}\PYG{p}{)}\PYG{p}{;}

\PYG{+w}{  }\PYG{n}{COPT\PYGZus{}DeleteEnv}\PYG{p}{(}\PYG{o}{\PYGZam{}}\PYG{n}{env}\PYG{p}{)}\PYG{p}{;}

\PYG{+w}{  }\PYG{k}{return}\PYG{+w}{ }\PYG{l+m+mi}{0}\PYG{p}{;}
\PYG{p}{\PYGZcb{}}
\end{sphinxVerbatim}

\sphinxAtStartPar
We will explain how to use the C API step by step based on code above, please
refer to {\hyperref[\detokenize{capiref:chapapi}]{\sphinxcrossref{\DUrole{std,std-ref}{C API Reference}}}} for detailed usage of C API.

\subsubsection{Creating the environment}
\label{\detokenize{cinterface:creating-the-environment}}
\sphinxAtStartPar
To solve any problem with Cardinal Optimizer, users are required to create
optimization environment first, and check if it was created successfully by
checking the return value:

\begin{sphinxVerbatim}[commandchars=\\\{\}]
\PYG{+w}{  }\PYG{c+c1}{// Create COPT environment}
\PYG{+w}{  }\PYG{n}{errcode}\PYG{+w}{ }\PYG{o}{=}\PYG{+w}{ }\PYG{n}{COPT\PYGZus{}CreateEnv}\PYG{p}{(}\PYG{o}{\PYGZam{}}\PYG{n}{env}\PYG{p}{)}\PYG{p}{;}
\PYG{+w}{  }\PYG{k}{if}\PYG{+w}{ }\PYG{p}{(}\PYG{n}{errcode}\PYG{p}{)}
\PYG{+w}{    }\PYG{k}{goto}\PYG{+w}{ }\PYG{n}{COPT\PYGZus{}EXIT}\PYG{p}{;}
\end{sphinxVerbatim}

\sphinxAtStartPar
If non\sphinxhyphen{}zero value was returned, it will jump to error reporting code block for
detailed information and exit.

\subsubsection{Creating the problem}
\label{\detokenize{cinterface:creating-the-problem}}
\sphinxAtStartPar
Once the optimization environment was successfully created, users will need to
create problem then, the problem is the main structure that consists of variables,
constraints etc. Users need to check the return value too.

\begin{sphinxVerbatim}[commandchars=\\\{\}]
\PYG{+w}{  }\PYG{c+c1}{// Create COPT problem}
\PYG{+w}{  }\PYG{n}{errcode}\PYG{+w}{ }\PYG{o}{=}\PYG{+w}{ }\PYG{n}{COPT\PYGZus{}CreateProb}\PYG{p}{(}\PYG{n}{env}\PYG{p}{,}\PYG{+w}{ }\PYG{o}{\PYGZam{}}\PYG{n}{prob}\PYG{p}{)}\PYG{p}{;}
\PYG{+w}{  }\PYG{k}{if}\PYG{+w}{ }\PYG{p}{(}\PYG{n}{errcode}\PYG{p}{)}
\PYG{+w}{    }\PYG{k}{goto}\PYG{+w}{ }\PYG{n}{COPT\PYGZus{}EXIT}\PYG{p}{;}
\end{sphinxVerbatim}

\sphinxAtStartPar
If non\sphinxhyphen{}zero value was returned, it will jump to error reporting code block for
detailed information and exit.

\subsubsection{Adding variables}
\label{\detokenize{cinterface:adding-variables}}
\sphinxAtStartPar
For linear problem, C API allows users to specify costs of variables in objective,
and lower and upper bound simultaneously. For the problem above, we use code below
to create variables:

\begin{sphinxVerbatim}[commandchars=\\\{\}]
\PYG{+w}{  }\PYG{c+cm}{/*}
\PYG{c+cm}{   * Add variables}
\PYG{c+cm}{   *}
\PYG{c+cm}{   *   obj: 1.2 C0 + 1.8 C1 + 2.1 C2}
\PYG{c+cm}{   *}
\PYG{c+cm}{   *   var:}
\PYG{c+cm}{   *     0.1 \PYGZlt{}= C0 \PYGZlt{}= 0.6}
\PYG{c+cm}{   *     0.2 \PYGZlt{}= C1 \PYGZlt{}= 1.5}
\PYG{c+cm}{   *     0.3 \PYGZlt{}= C2 \PYGZlt{}= 2.8}
\PYG{c+cm}{   *}
\PYG{c+cm}{   */}
\PYG{+w}{  }\PYG{k+kt}{int}\PYG{+w}{ }\PYG{n}{ncol}\PYG{+w}{ }\PYG{o}{=}\PYG{+w}{ }\PYG{l+m+mi}{3}\PYG{p}{;}
\PYG{+w}{  }\PYG{k+kt}{double}\PYG{+w}{ }\PYG{n}{colcost}\PYG{p}{[}\PYG{p}{]}\PYG{+w}{ }\PYG{o}{=}\PYG{+w}{ }\PYG{p}{\PYGZob{}}\PYG{l+m+mf}{1.2}\PYG{p}{,}\PYG{+w}{ }\PYG{l+m+mf}{1.8}\PYG{p}{,}\PYG{+w}{ }\PYG{l+m+mf}{2.1}\PYG{p}{\PYGZcb{}}\PYG{p}{;}
\PYG{+w}{  }\PYG{k+kt}{double}\PYG{+w}{ }\PYG{n}{collb}\PYG{p}{[}\PYG{p}{]}\PYG{+w}{ }\PYG{o}{=}\PYG{+w}{ }\PYG{p}{\PYGZob{}}\PYG{l+m+mf}{0.1}\PYG{p}{,}\PYG{+w}{ }\PYG{l+m+mf}{0.2}\PYG{p}{,}\PYG{+w}{ }\PYG{l+m+mf}{0.3}\PYG{p}{\PYGZcb{}}\PYG{p}{;}
\PYG{+w}{  }\PYG{k+kt}{double}\PYG{+w}{ }\PYG{n}{colub}\PYG{p}{[}\PYG{p}{]}\PYG{+w}{ }\PYG{o}{=}\PYG{+w}{ }\PYG{p}{\PYGZob{}}\PYG{l+m+mf}{0.6}\PYG{p}{,}\PYG{+w}{ }\PYG{l+m+mf}{1.5}\PYG{p}{,}\PYG{+w}{ }\PYG{l+m+mf}{1.8}\PYG{p}{\PYGZcb{}}\PYG{p}{;}

\PYG{+w}{  }\PYG{n}{errcode}\PYG{+w}{ }\PYG{o}{=}\PYG{+w}{ }\PYG{n}{COPT\PYGZus{}AddCols}\PYG{p}{(}\PYG{n}{prob}\PYG{p}{,}\PYG{+w}{ }\PYG{n}{ncol}\PYG{p}{,}\PYG{+w}{ }\PYG{n}{colcost}\PYG{p}{,}\PYG{+w}{ }\PYG{n+nb}{NULL}\PYG{p}{,}\PYG{+w}{ }\PYG{n+nb}{NULL}\PYG{p}{,}\PYG{+w}{ }\PYG{n+nb}{NULL}\PYG{p}{,}\PYG{+w}{ }\PYG{n+nb}{NULL}\PYG{p}{,}\PYG{+w}{ }\PYG{n+nb}{NULL}\PYG{p}{,}\PYG{+w}{ }\PYG{n}{collb}\PYG{p}{,}\PYG{+w}{ }\PYG{n}{colub}\PYG{p}{,}\PYG{+w}{ }\PYG{n+nb}{NULL}\PYG{p}{)}\PYG{p}{;}
\PYG{+w}{  }\PYG{k}{if}\PYG{+w}{ }\PYG{p}{(}\PYG{n}{errcode}\PYG{p}{)}
\PYG{+w}{    }\PYG{k}{goto}\PYG{+w}{ }\PYG{n}{COPT\PYGZus{}EXIT}\PYG{p}{;}
\end{sphinxVerbatim}

\sphinxAtStartPar
The argument \sphinxcode{\sphinxupquote{ncol}} specify that the number of variables to create is 3, while
the argument \sphinxcode{\sphinxupquote{colcost}}, \sphinxcode{\sphinxupquote{collb}} and \sphinxcode{\sphinxupquote{colub}} specify the costs in objective,
lower and upper bound respectively. Regarding other arguments of \sphinxcode{\sphinxupquote{COPT\_AddCols}}
for specifing variables types and names, we just pass \sphinxcode{\sphinxupquote{NULL}} to them, which
means all variables are continuous and names are automatically generated by
the Cardinal Optimizer. For the remaining arguments, we passed \sphinxcode{\sphinxupquote{NULL}} too
for further action.

\sphinxAtStartPar
Similarly, if non\sphinxhyphen{}zero value was returned, it will jump to error reporting code
block for detailed information and exit.

\subsubsection{Adding constraints}
\label{\detokenize{cinterface:adding-constraints}}
\sphinxAtStartPar
The next step to do after adding variables successfully is to add constraints to
problem. For the problem above, the implementation is shown below:

\begin{sphinxVerbatim}[commandchars=\\\{\}]
\PYG{+w}{  }\PYG{c+cm}{/*}
\PYG{c+cm}{   * Add constraints}
\PYG{c+cm}{   *}
\PYG{c+cm}{   *   r0: 1.5 C0 + 1.2 C1 + 1.8 C2 \PYGZlt{}= 2.6}
\PYG{c+cm}{   *   r1: 0.8 C0 + 0.6 C1 + 0.9 C2 \PYGZgt{}= 1.2}
\PYG{c+cm}{   */}
\PYG{+w}{  }\PYG{k+kt}{int}\PYG{+w}{ }\PYG{n}{nrow}\PYG{+w}{ }\PYG{o}{=}\PYG{+w}{ }\PYG{l+m+mi}{2}\PYG{p}{;}
\PYG{+w}{  }\PYG{k+kt}{int}\PYG{+w}{ }\PYG{n}{rowbeg}\PYG{p}{[}\PYG{p}{]}\PYG{+w}{ }\PYG{o}{=}\PYG{+w}{ }\PYG{p}{\PYGZob{}}\PYG{l+m+mi}{0}\PYG{p}{,}\PYG{+w}{ }\PYG{l+m+mi}{3}\PYG{p}{\PYGZcb{}}\PYG{p}{;}
\PYG{+w}{  }\PYG{k+kt}{int}\PYG{+w}{ }\PYG{n}{rowcnt}\PYG{p}{[}\PYG{p}{]}\PYG{+w}{ }\PYG{o}{=}\PYG{+w}{ }\PYG{p}{\PYGZob{}}\PYG{l+m+mi}{3}\PYG{p}{,}\PYG{+w}{ }\PYG{l+m+mi}{3}\PYG{p}{\PYGZcb{}}\PYG{p}{;}
\PYG{+w}{  }\PYG{k+kt}{int}\PYG{+w}{ }\PYG{n}{rowind}\PYG{p}{[}\PYG{p}{]}\PYG{+w}{ }\PYG{o}{=}\PYG{+w}{ }\PYG{p}{\PYGZob{}}\PYG{l+m+mi}{0}\PYG{p}{,}\PYG{+w}{ }\PYG{l+m+mi}{1}\PYG{p}{,}\PYG{+w}{ }\PYG{l+m+mi}{2}\PYG{p}{,}\PYG{+w}{ }\PYG{l+m+mi}{0}\PYG{p}{,}\PYG{+w}{ }\PYG{l+m+mi}{1}\PYG{p}{,}\PYG{+w}{ }\PYG{l+m+mi}{2}\PYG{p}{\PYGZcb{}}\PYG{p}{;}
\PYG{+w}{  }\PYG{k+kt}{double}\PYG{+w}{ }\PYG{n}{rowelem}\PYG{p}{[}\PYG{p}{]}\PYG{+w}{ }\PYG{o}{=}\PYG{+w}{ }\PYG{p}{\PYGZob{}}\PYG{l+m+mf}{1.5}\PYG{p}{,}\PYG{+w}{ }\PYG{l+m+mf}{1.2}\PYG{p}{,}\PYG{+w}{ }\PYG{l+m+mf}{1.8}\PYG{p}{,}\PYG{+w}{ }\PYG{l+m+mf}{0.8}\PYG{p}{,}\PYG{+w}{ }\PYG{l+m+mf}{0.6}\PYG{p}{,}\PYG{+w}{ }\PYG{l+m+mf}{0.9}\PYG{p}{\PYGZcb{}}\PYG{p}{;}
\PYG{+w}{  }\PYG{k+kt}{char}\PYG{+w}{ }\PYG{n}{rowsen}\PYG{p}{[}\PYG{p}{]}\PYG{+w}{ }\PYG{o}{=}\PYG{+w}{ }\PYG{p}{\PYGZob{}}\PYG{n}{COPT\PYGZus{}LESS\PYGZus{}EQUAL}\PYG{p}{,}\PYG{+w}{ }\PYG{n}{COPT\PYGZus{}GREATER\PYGZus{}EQUAL}\PYG{p}{\PYGZcb{}}\PYG{p}{;}
\PYG{+w}{  }\PYG{k+kt}{double}\PYG{+w}{ }\PYG{n}{rowrhs}\PYG{p}{[}\PYG{p}{]}\PYG{+w}{ }\PYG{o}{=}\PYG{+w}{ }\PYG{p}{\PYGZob{}}\PYG{l+m+mf}{2.6}\PYG{p}{,}\PYG{+w}{ }\PYG{l+m+mf}{1.2}\PYG{p}{\PYGZcb{}}\PYG{p}{;}

\PYG{+w}{  }\PYG{n}{errcode}\PYG{+w}{ }\PYG{o}{=}\PYG{+w}{ }\PYG{n}{COPT\PYGZus{}AddRows}\PYG{p}{(}\PYG{n}{prob}\PYG{p}{,}\PYG{+w}{ }\PYG{n}{nrow}\PYG{p}{,}\PYG{+w}{ }\PYG{n}{rowbeg}\PYG{p}{,}\PYG{+w}{ }\PYG{n}{rowcnt}\PYG{p}{,}\PYG{+w}{ }\PYG{n}{rowind}\PYG{p}{,}\PYG{+w}{ }\PYG{n}{rowelem}\PYG{p}{,}\PYG{+w}{ }\PYG{n}{rowsen}\PYG{p}{,}\PYG{+w}{ }\PYG{n}{rowrhs}\PYG{p}{,}\PYG{+w}{ }\PYG{n+nb}{NULL}\PYG{p}{,}\PYG{+w}{ }\PYG{n+nb}{NULL}\PYG{p}{)}\PYG{p}{;}
\PYG{+w}{  }\PYG{k}{if}\PYG{+w}{ }\PYG{p}{(}\PYG{n}{errcode}\PYG{p}{)}
\PYG{+w}{    }\PYG{k}{goto}\PYG{+w}{ }\PYG{n}{COPT\PYGZus{}EXIT}\PYG{p}{;}
\end{sphinxVerbatim}

\sphinxAtStartPar
The argument \sphinxcode{\sphinxupquote{nrow}} specifies that the number of constraints to create is 2,
while argument \sphinxcode{\sphinxupquote{rowbeg}}, \sphinxcode{\sphinxupquote{rowcnt}}, \sphinxcode{\sphinxupquote{rowind}} and \sphinxcode{\sphinxupquote{rowelem}} define the
coefficient matrix in CSR format. The argument \sphinxcode{\sphinxupquote{rowsen}} represents the sense
of constraints, while argument \sphinxcode{\sphinxupquote{rowrhs}} specifies the right hand side of
constraints. For remaining arguments in \sphinxcode{\sphinxupquote{COPT\_AddRows}}, we simply pass \sphinxcode{\sphinxupquote{NULL}}
to them.

\sphinxAtStartPar
If the return value is non\sphinxhyphen{}zero, then it jump to error reporting code block for
detailed information and exit.

\subsubsection{Setting parameters and attributes}
\label{\detokenize{cinterface:setting-parameters-and-attributes}}
\sphinxAtStartPar
Users are allowed to set parameters and attributes of problem before solving.
For example, to set the time limit to 10 seconds, and to set the optimization
direction to maximization, the code is shown below:

\begin{sphinxVerbatim}[commandchars=\\\{\}]
\PYG{+w}{  }\PYG{c+c1}{// Set parameters and attributes}
\PYG{+w}{  }\PYG{n}{errcode}\PYG{+w}{ }\PYG{o}{=}\PYG{+w}{ }\PYG{n}{COPT\PYGZus{}SetDblParam}\PYG{p}{(}\PYG{n}{prob}\PYG{p}{,}\PYG{+w}{ }\PYG{n}{COPT\PYGZus{}DBLPARAM\PYGZus{}TIMELIMIT}\PYG{p}{,}\PYG{+w}{ }\PYG{l+m+mi}{10}\PYG{p}{)}\PYG{p}{;}
\PYG{+w}{  }\PYG{k}{if}\PYG{+w}{ }\PYG{p}{(}\PYG{n}{errcode}\PYG{p}{)}
\PYG{+w}{    }\PYG{k}{goto}\PYG{+w}{ }\PYG{n}{COPT\PYGZus{}EXIT}\PYG{p}{;}
\PYG{+w}{  }\PYG{n}{errcode}\PYG{+w}{ }\PYG{o}{=}\PYG{+w}{ }\PYG{n}{COPT\PYGZus{}SetObjSense}\PYG{p}{(}\PYG{n}{prob}\PYG{p}{,}\PYG{+w}{ }\PYG{n}{COPT\PYGZus{}MAXIMIZE}\PYG{p}{)}\PYG{p}{;}
\PYG{+w}{  }\PYG{k}{if}\PYG{+w}{ }\PYG{p}{(}\PYG{n}{errcode}\PYG{p}{)}
\PYG{+w}{    }\PYG{k}{goto}\PYG{+w}{ }\PYG{n}{COPT\PYGZus{}EXIT}\PYG{p}{;}
\end{sphinxVerbatim}

\sphinxAtStartPar
If non\sphinxhyphen{}zero value was returned, then it will jump to error reporting code block
for detailed information and exit.

\subsubsection{Solve the problem}
\label{\detokenize{cinterface:solve-the-problem}}
\sphinxAtStartPar
The next step to do is to solve the problem using code beblow:

\begin{sphinxVerbatim}[commandchars=\\\{\}]
\PYG{+w}{  }\PYG{c+c1}{// Solve problem}
\PYG{+w}{  }\PYG{n}{errcode}\PYG{+w}{ }\PYG{o}{=}\PYG{+w}{ }\PYG{n}{COPT\PYGZus{}SolveLp}\PYG{p}{(}\PYG{n}{prob}\PYG{p}{)}\PYG{p}{;}
\PYG{+w}{  }\PYG{k}{if}\PYG{+w}{ }\PYG{p}{(}\PYG{n}{errcode}\PYG{p}{)}
\PYG{+w}{    }\PYG{k}{goto}\PYG{+w}{ }\PYG{n}{COPT\PYGZus{}EXIT}\PYG{p}{;}
\end{sphinxVerbatim}

\sphinxAtStartPar
Non\sphinxhyphen{}zero return value indicates unsuccessful solve and jump to error reporting
code block for detailed information and exit.

\subsubsection{Analyze the solution}
\label{\detokenize{cinterface:analyze-the-solution}}
\sphinxAtStartPar
Once the solving process was finished, check the solution status first. If it
claimed to have found the optimal solution, then use code below to obtain objective
value, variables’ solution and basis status:

\begin{sphinxVerbatim}[commandchars=\\\{\}]
\PYG{+w}{  }\PYG{c+c1}{// Analyze solution}
\PYG{+w}{  }\PYG{k+kt}{int}\PYG{+w}{ }\PYG{n}{lpstat}\PYG{+w}{ }\PYG{o}{=}\PYG{+w}{ }\PYG{n}{COPT\PYGZus{}LPSTATUS\PYGZus{}UNSTARTED}\PYG{p}{;}
\PYG{+w}{  }\PYG{k+kt}{double}\PYG{+w}{ }\PYG{n}{lpobjval}\PYG{p}{;}
\PYG{+w}{  }\PYG{k+kt}{double}\PYG{o}{*}\PYG{+w}{ }\PYG{n}{lpsol}\PYG{+w}{ }\PYG{o}{=}\PYG{+w}{ }\PYG{n+nb}{NULL}\PYG{p}{;}
\PYG{+w}{  }\PYG{k+kt}{int}\PYG{o}{*}\PYG{+w}{ }\PYG{n}{colstat}\PYG{+w}{ }\PYG{o}{=}\PYG{+w}{ }\PYG{n+nb}{NULL}\PYG{p}{;}

\PYG{+w}{  }\PYG{n}{errcode}\PYG{+w}{ }\PYG{o}{=}\PYG{+w}{ }\PYG{n}{COPT\PYGZus{}GetIntAttr}\PYG{p}{(}\PYG{n}{prob}\PYG{p}{,}\PYG{+w}{ }\PYG{n}{COPT\PYGZus{}INTATTR\PYGZus{}LPSTATUS}\PYG{p}{,}\PYG{+w}{ }\PYG{o}{\PYGZam{}}\PYG{n}{lpstat}\PYG{p}{)}\PYG{p}{;}
\PYG{+w}{  }\PYG{k}{if}\PYG{+w}{ }\PYG{p}{(}\PYG{n}{errcode}\PYG{p}{)}
\PYG{+w}{    }\PYG{k}{goto}\PYG{+w}{ }\PYG{n}{COPT\PYGZus{}EXIT}\PYG{p}{;}

\PYG{+w}{  }\PYG{k}{if}\PYG{+w}{ }\PYG{p}{(}\PYG{n}{lpstat}\PYG{+w}{ }\PYG{o}{=}\PYG{o}{=}\PYG{+w}{ }\PYG{n}{COPT\PYGZus{}LPSTATUS\PYGZus{}OPTIMAL}\PYG{p}{)}
\PYG{+w}{  }\PYG{p}{\PYGZob{}}
\PYG{+w}{    }\PYG{n}{lpsol}\PYG{+w}{ }\PYG{o}{=}\PYG{+w}{ }\PYG{p}{(}\PYG{k+kt}{double}\PYG{o}{*}\PYG{p}{)}\PYG{n}{malloc}\PYG{p}{(}\PYG{n}{ncol}\PYG{+w}{ }\PYG{o}{*}\PYG{+w}{ }\PYG{k}{sizeof}\PYG{p}{(}\PYG{k+kt}{double}\PYG{p}{)}\PYG{p}{)}\PYG{p}{;}
\PYG{+w}{    }\PYG{n}{colstat}\PYG{+w}{ }\PYG{o}{=}\PYG{+w}{ }\PYG{p}{(}\PYG{k+kt}{int}\PYG{o}{*}\PYG{p}{)}\PYG{n}{malloc}\PYG{p}{(}\PYG{n}{ncol}\PYG{+w}{ }\PYG{o}{*}\PYG{+w}{ }\PYG{k}{sizeof}\PYG{p}{(}\PYG{k+kt}{int}\PYG{p}{)}\PYG{p}{)}\PYG{p}{;}

\PYG{+w}{    }\PYG{n}{errcode}\PYG{+w}{ }\PYG{o}{=}\PYG{+w}{ }\PYG{n}{COPT\PYGZus{}GetLpSolution}\PYG{p}{(}\PYG{n}{prob}\PYG{p}{,}\PYG{+w}{ }\PYG{n}{lpsol}\PYG{p}{,}\PYG{+w}{ }\PYG{n+nb}{NULL}\PYG{p}{,}\PYG{+w}{ }\PYG{n+nb}{NULL}\PYG{p}{,}\PYG{+w}{ }\PYG{n+nb}{NULL}\PYG{p}{)}\PYG{p}{;}
\PYG{+w}{    }\PYG{k}{if}\PYG{+w}{ }\PYG{p}{(}\PYG{n}{errcode}\PYG{p}{)}
\PYG{+w}{      }\PYG{k}{goto}\PYG{+w}{ }\PYG{n}{COPT\PYGZus{}EXIT}\PYG{p}{;}

\PYG{+w}{    }\PYG{n}{errcode}\PYG{+w}{ }\PYG{o}{=}\PYG{+w}{ }\PYG{n}{COPT\PYGZus{}GetBasis}\PYG{p}{(}\PYG{n}{prob}\PYG{p}{,}\PYG{+w}{ }\PYG{n}{colstat}\PYG{p}{,}\PYG{+w}{ }\PYG{n+nb}{NULL}\PYG{p}{)}\PYG{p}{;}
\PYG{+w}{    }\PYG{k}{if}\PYG{+w}{ }\PYG{p}{(}\PYG{n}{errcode}\PYG{p}{)}
\PYG{+w}{      }\PYG{k}{goto}\PYG{+w}{ }\PYG{n}{COPT\PYGZus{}EXIT}\PYG{p}{;}

\PYG{+w}{    }\PYG{n}{errcode}\PYG{+w}{ }\PYG{o}{=}\PYG{+w}{ }\PYG{n}{COPT\PYGZus{}GetDblAttr}\PYG{p}{(}\PYG{n}{prob}\PYG{p}{,}\PYG{+w}{ }\PYG{n}{COPT\PYGZus{}DBLATTR\PYGZus{}LPOBJVAL}\PYG{p}{,}\PYG{+w}{ }\PYG{o}{\PYGZam{}}\PYG{n}{lpobjval}\PYG{p}{)}\PYG{p}{;}
\PYG{+w}{    }\PYG{k}{if}\PYG{+w}{ }\PYG{p}{(}\PYG{n}{errcode}\PYG{p}{)}
\PYG{+w}{      }\PYG{k}{goto}\PYG{+w}{ }\PYG{n}{COPT\PYGZus{}EXIT}\PYG{p}{;}

\PYG{+w}{    }\PYG{n}{printf}\PYG{p}{(}\PYG{l+s}{\PYGZdq{}}\PYG{l+s+se}{\PYGZbs{}n}\PYG{l+s}{Objective value: \PYGZpc{}.6f}\PYG{l+s+se}{\PYGZbs{}n}\PYG{l+s}{\PYGZdq{}}\PYG{p}{,}\PYG{+w}{ }\PYG{n}{lpobjval}\PYG{p}{)}\PYG{p}{;}

\PYG{+w}{    }\PYG{n}{printf}\PYG{p}{(}\PYG{l+s}{\PYGZdq{}}\PYG{l+s}{Variable solution: }\PYG{l+s+se}{\PYGZbs{}n}\PYG{l+s}{\PYGZdq{}}\PYG{p}{)}\PYG{p}{;}
\PYG{+w}{    }\PYG{k}{for}\PYG{+w}{ }\PYG{p}{(}\PYG{k+kt}{int}\PYG{+w}{ }\PYG{n}{i}\PYG{+w}{ }\PYG{o}{=}\PYG{+w}{ }\PYG{l+m+mi}{0}\PYG{p}{;}\PYG{+w}{ }\PYG{n}{i}\PYG{+w}{ }\PYG{o}{\PYGZlt{}}\PYG{+w}{ }\PYG{n}{ncol}\PYG{p}{;}\PYG{+w}{ }\PYG{o}{+}\PYG{o}{+}\PYG{n}{i}\PYG{p}{)}
\PYG{+w}{      }\PYG{n}{printf}\PYG{p}{(}\PYG{l+s}{\PYGZdq{}}\PYG{l+s}{  x[\PYGZpc{}d] = \PYGZpc{}.6f}\PYG{l+s+se}{\PYGZbs{}n}\PYG{l+s}{\PYGZdq{}}\PYG{p}{,}\PYG{+w}{ }\PYG{n}{i}\PYG{p}{,}\PYG{+w}{ }\PYG{n}{lpsol}\PYG{p}{[}\PYG{n}{i}\PYG{p}{]}\PYG{p}{)}\PYG{p}{;}

\PYG{+w}{    }\PYG{n}{printf}\PYG{p}{(}\PYG{l+s}{\PYGZdq{}}\PYG{l+s}{Variable basis status: }\PYG{l+s+se}{\PYGZbs{}n}\PYG{l+s}{\PYGZdq{}}\PYG{p}{)}\PYG{p}{;}
\PYG{+w}{    }\PYG{k}{for}\PYG{+w}{ }\PYG{p}{(}\PYG{k+kt}{int}\PYG{+w}{ }\PYG{n}{i}\PYG{+w}{ }\PYG{o}{=}\PYG{+w}{ }\PYG{l+m+mi}{0}\PYG{p}{;}\PYG{+w}{ }\PYG{n}{i}\PYG{+w}{ }\PYG{o}{\PYGZlt{}}\PYG{+w}{ }\PYG{n}{ncol}\PYG{p}{;}\PYG{+w}{ }\PYG{o}{+}\PYG{o}{+}\PYG{n}{i}\PYG{p}{)}
\PYG{+w}{      }\PYG{n}{printf}\PYG{p}{(}\PYG{l+s}{\PYGZdq{}}\PYG{l+s}{  x[\PYGZpc{}d]: \PYGZpc{}d}\PYG{l+s+se}{\PYGZbs{}n}\PYG{l+s}{\PYGZdq{}}\PYG{p}{,}\PYG{+w}{ }\PYG{n}{i}\PYG{p}{,}\PYG{+w}{ }\PYG{n}{colstat}\PYG{p}{[}\PYG{n}{i}\PYG{p}{]}\PYG{p}{)}\PYG{p}{;}

\PYG{+w}{    }\PYG{n}{free}\PYG{p}{(}\PYG{n}{lpsol}\PYG{p}{)}\PYG{p}{;}
\PYG{+w}{    }\PYG{n}{free}\PYG{p}{(}\PYG{n}{colstat}\PYG{p}{)}\PYG{p}{;}
\PYG{+w}{  }\PYG{p}{\PYGZcb{}}
\end{sphinxVerbatim}

\subsubsection{Write problem and solution}
\label{\detokenize{cinterface:write-problem-and-solution}}
\sphinxAtStartPar
Users are allowed not only to save the problem to solve to standard MPS file,
but also the solution, basis status and modified parameters to files:

\begin{sphinxVerbatim}[commandchars=\\\{\}]
\PYG{+w}{  }\PYG{c+c1}{// Write problem, solution and modified parameters to files}
\PYG{+w}{  }\PYG{n}{errcode}\PYG{+w}{ }\PYG{o}{=}\PYG{+w}{ }\PYG{n}{COPT\PYGZus{}WriteMps}\PYG{p}{(}\PYG{n}{prob}\PYG{p}{,}\PYG{+w}{ }\PYG{l+s}{\PYGZdq{}}\PYG{l+s}{lp\PYGZus{}ex1.mps}\PYG{l+s}{\PYGZdq{}}\PYG{p}{)}\PYG{p}{;}
\PYG{+w}{  }\PYG{k}{if}\PYG{+w}{ }\PYG{p}{(}\PYG{n}{errcode}\PYG{p}{)}
\PYG{+w}{    }\PYG{k}{goto}\PYG{+w}{ }\PYG{n}{COPT\PYGZus{}EXIT}\PYG{p}{;}
\PYG{+w}{  }\PYG{n}{errcode}\PYG{+w}{ }\PYG{o}{=}\PYG{+w}{ }\PYG{n}{COPT\PYGZus{}WriteBasis}\PYG{p}{(}\PYG{n}{prob}\PYG{p}{,}\PYG{+w}{ }\PYG{l+s}{\PYGZdq{}}\PYG{l+s}{lp\PYGZus{}ex1.bas}\PYG{l+s}{\PYGZdq{}}\PYG{p}{)}\PYG{p}{;}
\PYG{+w}{  }\PYG{k}{if}\PYG{+w}{ }\PYG{p}{(}\PYG{n}{errcode}\PYG{p}{)}
\PYG{+w}{    }\PYG{k}{goto}\PYG{+w}{ }\PYG{n}{COPT\PYGZus{}EXIT}\PYG{p}{;}
\PYG{+w}{  }\PYG{n}{errcode}\PYG{+w}{ }\PYG{o}{=}\PYG{+w}{ }\PYG{n}{COPT\PYGZus{}WriteSol}\PYG{p}{(}\PYG{n}{prob}\PYG{p}{,}\PYG{+w}{ }\PYG{l+s}{\PYGZdq{}}\PYG{l+s}{lp\PYGZus{}ex1.sol}\PYG{l+s}{\PYGZdq{}}\PYG{p}{)}\PYG{p}{;}
\PYG{+w}{  }\PYG{k}{if}\PYG{+w}{ }\PYG{p}{(}\PYG{n}{errcode}\PYG{p}{)}
\PYG{+w}{    }\PYG{k}{goto}\PYG{+w}{ }\PYG{n}{COPT\PYGZus{}EXIT}\PYG{p}{;}
\PYG{+w}{  }\PYG{n}{errcode}\PYG{+w}{ }\PYG{o}{=}\PYG{+w}{ }\PYG{n}{COPT\PYGZus{}WriteParam}\PYG{p}{(}\PYG{n}{prob}\PYG{p}{,}\PYG{+w}{ }\PYG{l+s}{\PYGZdq{}}\PYG{l+s}{lp\PYGZus{}ex1.par}\PYG{l+s}{\PYGZdq{}}\PYG{p}{)}\PYG{p}{;}
\PYG{+w}{  }\PYG{k}{if}\PYG{+w}{ }\PYG{p}{(}\PYG{n}{errcode}\PYG{p}{)}
\PYG{+w}{    }\PYG{k}{goto}\PYG{+w}{ }\PYG{n}{COPT\PYGZus{}EXIT}\PYG{p}{;}
\end{sphinxVerbatim}

\subsubsection{Error handling}
\label{\detokenize{cinterface:error-handling}}
\sphinxAtStartPar
The error handling block report error code and message by checking if the return
value was non\sphinxhyphen{}zero:

\begin{sphinxVerbatim}[commandchars=\\\{\}]
\PYG{+w}{  }\PYG{c+c1}{// Error handling}
\PYG{n+nl}{COPT\PYGZus{}EXIT}\PYG{p}{:}
\PYG{+w}{  }\PYG{k}{if}\PYG{+w}{ }\PYG{p}{(}\PYG{n}{errcode}\PYG{p}{)}
\PYG{+w}{  }\PYG{p}{\PYGZob{}}
\PYG{+w}{    }\PYG{k+kt}{char}\PYG{+w}{ }\PYG{n}{errmsg}\PYG{p}{[}\PYG{n}{COPT\PYGZus{}BUFFSIZE}\PYG{p}{]}\PYG{p}{;}

\PYG{+w}{    }\PYG{n}{COPT\PYGZus{}GetRetcodeMsg}\PYG{p}{(}\PYG{n}{errcode}\PYG{p}{,}\PYG{+w}{ }\PYG{n}{errmsg}\PYG{p}{,}\PYG{+w}{ }\PYG{n}{COPT\PYGZus{}BUFFSIZE}\PYG{p}{)}\PYG{p}{;}
\PYG{+w}{    }\PYG{n}{printf}\PYG{p}{(}\PYG{l+s}{\PYGZdq{}}\PYG{l+s}{ERROR \PYGZpc{}d: \PYGZpc{}s}\PYG{l+s+se}{\PYGZbs{}n}\PYG{l+s}{\PYGZdq{}}\PYG{p}{,}\PYG{+w}{ }\PYG{n}{errcode}\PYG{p}{,}\PYG{+w}{ }\PYG{n}{errmsg}\PYG{p}{)}\PYG{p}{;}

\PYG{+w}{    }\PYG{k}{return}\PYG{+w}{ }\PYG{l+m+mi}{0}\PYG{p}{;}
\PYG{+w}{  }\PYG{p}{\PYGZcb{}}
\end{sphinxVerbatim}

\subsubsection{Delete environment and problem}
\label{\detokenize{cinterface:delete-environment-and-problem}}
\sphinxAtStartPar
Before exiting, delete problem and environment respectively:

\begin{sphinxVerbatim}[commandchars=\\\{\}]
\PYG{+w}{  }\PYG{c+c1}{// Delete problem and environment}
\PYG{+w}{  }\PYG{n}{COPT\PYGZus{}DeleteProb}\PYG{p}{(}\PYG{o}{\PYGZam{}}\PYG{n}{prob}\PYG{p}{)}\PYG{p}{;}

\PYG{+w}{  }\PYG{n}{COPT\PYGZus{}DeleteEnv}\PYG{p}{(}\PYG{o}{\PYGZam{}}\PYG{n}{env}\PYG{p}{)}\PYG{p}{;}
\end{sphinxVerbatim}

\subsection{Build and run}
\label{\detokenize{cinterface:build-and-run}}
\sphinxAtStartPar
To ease the work for running the example for users on different operating systems,
we provide Visual Studio project and Makefile for Windows, Linux and MacOS
respectively, details are shown below.

\subsubsection{Windows}
\label{\detokenize{cinterface:windows}}
\sphinxAtStartPar
For users on Windows platform, we provide Visual Studio project, all users are
required to install Visual Studio 2017 beforehand. Assume that the installation
directory is: \sphinxcode{\sphinxupquote{\textquotesingle{}\textless{}instdir\textgreater{}\textquotesingle{}}}, users that install the Cardinal Optimizer with
executable installer can change directory to \sphinxcode{\sphinxupquote{\textquotesingle{}\textless{}instdir\textgreater{}\textbackslash{}examples\textbackslash{}c\textbackslash{}vsprojects\textquotesingle{}}}
and open the Visual Studio project \sphinxcode{\sphinxupquote{lp\_ex1.vcxproj}} to build the solution.
Users that install the Cardinal Optimizer using ZIP\sphinxhyphen{}format archive should make
sure that all required environment variables are set correctly, see
{\hyperref[\detokenize{install:chapinstall}]{\sphinxcrossref{\DUrole{std,std-ref}{Install Guide for Cardinal Optimizer}}}} for details.

\sphinxAtStartPar
In addition, for Windows systems, COPT also supports the MinGW\sphinxhyphen{}w64 toolchain. Please refer to the corresponding instructions
in the Makefile under the path \sphinxcode{\sphinxupquote{\textquotesingle{}\textless{}instdir\textgreater{}\textbackslash{}examples\textbackslash{}c\textquotesingle{}}} for usage.

\subsubsection{Linux and MacOS}
\label{\detokenize{cinterface:linux-and-macos}}
\sphinxAtStartPar
For users on Linux or MacOS, we provide Makefile to build the example. Please
install GCC toolchain for Linux and Clang toolchain for MacOS, together with
the \sphinxcode{\sphinxupquote{make}} utility beforehand. What’s more, users should make sure also that
all required environment variables are set correctly, see
{\hyperref[\detokenize{install:chapinstall}]{\sphinxcrossref{\DUrole{std,std-ref}{Install Guide for Cardinal Optimizer}}}} for details. Let’s
assume that the installation directory of Cardinal Optimizer is \sphinxcode{\sphinxupquote{\textquotesingle{}\textless{}instdir\textgreater{}\textquotesingle{}}},
then users need to change directory to \sphinxcode{\sphinxupquote{\textquotesingle{}\textless{}instdir\textgreater{}\textbackslash{}examples\textbackslash{}c\textquotesingle{}}} and execute
command \sphinxcode{\sphinxupquote{make}} in terminal.

\sphinxstepscope

\section{C++ Interface}
\label{\detokenize{cppinterface:c-interface}}\label{\detokenize{cppinterface:chapcppinterface}}\label{\detokenize{cppinterface::doc}}
\sphinxAtStartPar
This chapter walks through a simple C++ example to illustrate the use of the
COPT C++ interface. In short words, the example creates an environment, builds
a model, add variables and constraints, optimizes it, and then outputs the
optimal objective value.

\sphinxAtStartPar
The example solves the following linear problem:
\begin{equation}\label{equation:cppinterface:coptcppEq_capilp1}
\begin{split}\text{Maximize: } & \\
                  & 1.2 x + 1.8 y + 2.1 z \\
\text{Subject to: } & \\
                    & 1.5 x + 1.2 y + 1.8 z \leq 2.6 \\
                    & 0.8 x + 0.6 y + 0.9 z \geq 1.2 \\
\text{Bounds: } & \\
               & 0.1 \leq x \leq 0.6 \\
               & 0.2 \leq y \leq 1.5 \\
               & 0.3 \leq z \leq 2.8 \\\end{split}
\end{equation}
\sphinxAtStartPar
Note that this is the same problem that was modelled and optimized in
chapter of {\hyperref[\detokenize{cinterface:chapcinterface}]{\sphinxcrossref{\DUrole{std,std-ref}{C Interface}}}}.

\subsection{Example details}
\label{\detokenize{cppinterface:example-details}}
\sphinxAtStartPar
Below is the source code solving the above problem using COPT C++ interface.
\sphinxSetupCaptionForVerbatim{\sphinxcode{\sphinxupquote{lp\_ex1.cpp}}}
\def\sphinxLiteralBlockLabel{\label{\detokenize{cppinterface:id1}}}
\begin{sphinxVerbatim}[commandchars=\\\{\},numbers=left,firstnumber=1,stepnumber=1]
\PYG{c+cm}{/*}
\PYG{c+cm}{ * This file is part of the Cardinal Optimizer, all rights reserved.}
\PYG{c+cm}{ */}
\PYG{c+cp}{\PYGZsh{}}\PYG{c+cp}{include}\PYG{+w}{ }\PYG{c+cpf}{\PYGZdq{}coptcpp\PYGZus{}pch.h\PYGZdq{}}

\PYG{k}{using}\PYG{+w}{ }\PYG{k}{namespace}\PYG{+w}{ }\PYG{n+nn}{std}\PYG{p}{;}

\PYG{c+cm}{/*}
\PYG{c+cm}{ * This example solves the following LP model:}
\PYG{c+cm}{ *}
\PYG{c+cm}{ *  Maximize:}
\PYG{c+cm}{ *    1.2 x + 1.8 y + 2.1 z}
\PYG{c+cm}{ *}
\PYG{c+cm}{ *  Subject to:}
\PYG{c+cm}{ *    R0: 1.5 x + 1.2 y + 1.8 z \PYGZlt{}= 2.6}
\PYG{c+cm}{ *    R1: 0.8 x + 0.6 y + 0.9 z \PYGZgt{}= 1.2}
\PYG{c+cm}{ *}
\PYG{c+cm}{ *  Where:}
\PYG{c+cm}{ *    0.1 \PYGZlt{}= x \PYGZlt{}= 0.6}
\PYG{c+cm}{ *    0.2 \PYGZlt{}= y \PYGZlt{}= 1.5}
\PYG{c+cm}{ *    0.3 \PYGZlt{}= z \PYGZlt{}= 2.8}
\PYG{c+cm}{ */}
\PYG{k+kt}{int}\PYG{+w}{ }\PYG{n+nf}{main}\PYG{p}{(}\PYG{k+kt}{int}\PYG{+w}{ }\PYG{n}{argc}\PYG{p}{,}\PYG{+w}{ }\PYG{k+kt}{char}\PYG{o}{*}\PYG{+w}{ }\PYG{n}{argv}\PYG{p}{[}\PYG{p}{]}\PYG{p}{)}
\PYG{p}{\PYGZob{}}
\PYG{+w}{  }\PYG{k}{try}
\PYG{+w}{  }\PYG{p}{\PYGZob{}}
\PYG{+w}{    }\PYG{n}{Envr}\PYG{+w}{ }\PYG{n}{env}\PYG{p}{;}
\PYG{+w}{    }\PYG{n}{Model}\PYG{+w}{ }\PYG{n}{model}\PYG{+w}{ }\PYG{o}{=}\PYG{+w}{ }\PYG{n}{env}\PYG{p}{.}\PYG{n}{CreateModel}\PYG{p}{(}\PYG{l+s}{\PYGZdq{}}\PYG{l+s}{lp\PYGZus{}ex1}\PYG{l+s}{\PYGZdq{}}\PYG{p}{)}\PYG{p}{;}

\PYG{+w}{    }\PYG{c+c1}{// Add variables}
\PYG{+w}{    }\PYG{n}{Var}\PYG{+w}{ }\PYG{n}{x}\PYG{+w}{ }\PYG{o}{=}\PYG{+w}{ }\PYG{n}{model}\PYG{p}{.}\PYG{n}{AddVar}\PYG{p}{(}\PYG{l+m+mf}{0.1}\PYG{p}{,}\PYG{+w}{ }\PYG{l+m+mf}{0.6}\PYG{p}{,}\PYG{+w}{ }\PYG{l+m+mf}{0.0}\PYG{p}{,}\PYG{+w}{ }\PYG{n}{COPT\PYGZus{}CONTINUOUS}\PYG{p}{,}\PYG{+w}{ }\PYG{l+s}{\PYGZdq{}}\PYG{l+s}{x}\PYG{l+s}{\PYGZdq{}}\PYG{p}{)}\PYG{p}{;}
\PYG{+w}{    }\PYG{n}{Var}\PYG{+w}{ }\PYG{n}{y}\PYG{+w}{ }\PYG{o}{=}\PYG{+w}{ }\PYG{n}{model}\PYG{p}{.}\PYG{n}{AddVar}\PYG{p}{(}\PYG{l+m+mf}{0.2}\PYG{p}{,}\PYG{+w}{ }\PYG{l+m+mf}{1.5}\PYG{p}{,}\PYG{+w}{ }\PYG{l+m+mf}{0.0}\PYG{p}{,}\PYG{+w}{ }\PYG{n}{COPT\PYGZus{}CONTINUOUS}\PYG{p}{,}\PYG{+w}{ }\PYG{l+s}{\PYGZdq{}}\PYG{l+s}{y}\PYG{l+s}{\PYGZdq{}}\PYG{p}{)}\PYG{p}{;}
\PYG{+w}{    }\PYG{n}{Var}\PYG{+w}{ }\PYG{n}{z}\PYG{+w}{ }\PYG{o}{=}\PYG{+w}{ }\PYG{n}{model}\PYG{p}{.}\PYG{n}{AddVar}\PYG{p}{(}\PYG{l+m+mf}{0.3}\PYG{p}{,}\PYG{+w}{ }\PYG{l+m+mf}{2.8}\PYG{p}{,}\PYG{+w}{ }\PYG{l+m+mf}{0.0}\PYG{p}{,}\PYG{+w}{ }\PYG{n}{COPT\PYGZus{}CONTINUOUS}\PYG{p}{,}\PYG{+w}{ }\PYG{l+s}{\PYGZdq{}}\PYG{l+s}{z}\PYG{l+s}{\PYGZdq{}}\PYG{p}{)}\PYG{p}{;}

\PYG{+w}{    }\PYG{c+c1}{// Set objective}
\PYG{+w}{    }\PYG{n}{model}\PYG{p}{.}\PYG{n}{SetObjective}\PYG{p}{(}\PYG{l+m+mf}{1.2}\PYG{+w}{ }\PYG{o}{*}\PYG{+w}{ }\PYG{n}{x}\PYG{+w}{ }\PYG{o}{+}\PYG{+w}{ }\PYG{l+m+mf}{1.8}\PYG{+w}{ }\PYG{o}{*}\PYG{+w}{ }\PYG{n}{y}\PYG{+w}{ }\PYG{o}{+}\PYG{+w}{ }\PYG{l+m+mf}{2.1}\PYG{+w}{ }\PYG{o}{*}\PYG{+w}{ }\PYG{n}{z}\PYG{p}{,}\PYG{+w}{ }\PYG{n}{COPT\PYGZus{}MAXIMIZE}\PYG{p}{)}\PYG{p}{;}

\PYG{+w}{    }\PYG{c+c1}{// Add linear constraints using linear expression}
\PYG{+w}{    }\PYG{n}{model}\PYG{p}{.}\PYG{n}{AddConstr}\PYG{p}{(}\PYG{l+m+mf}{1.5}\PYG{+w}{ }\PYG{o}{*}\PYG{+w}{ }\PYG{n}{x}\PYG{+w}{ }\PYG{o}{+}\PYG{+w}{ }\PYG{l+m+mf}{1.2}\PYG{+w}{ }\PYG{o}{*}\PYG{+w}{ }\PYG{n}{y}\PYG{+w}{ }\PYG{o}{+}\PYG{+w}{ }\PYG{l+m+mf}{1.8}\PYG{+w}{ }\PYG{o}{*}\PYG{+w}{ }\PYG{n}{z}\PYG{+w}{ }\PYG{o}{\PYGZlt{}}\PYG{o}{=}\PYG{+w}{ }\PYG{l+m+mf}{2.6}\PYG{p}{,}\PYG{+w}{ }\PYG{l+s}{\PYGZdq{}}\PYG{l+s}{R0}\PYG{l+s}{\PYGZdq{}}\PYG{p}{)}\PYG{p}{;}

\PYG{+w}{    }\PYG{n}{Expr}\PYG{+w}{ }\PYG{n}{expr}\PYG{p}{(}\PYG{n}{x}\PYG{p}{,}\PYG{+w}{ }\PYG{l+m+mf}{0.8}\PYG{p}{)}\PYG{p}{;}
\PYG{+w}{    }\PYG{n}{expr}\PYG{p}{.}\PYG{n}{AddTerm}\PYG{p}{(}\PYG{n}{y}\PYG{p}{,}\PYG{+w}{ }\PYG{l+m+mf}{0.6}\PYG{p}{)}\PYG{p}{;}
\PYG{+w}{    }\PYG{n}{expr}\PYG{+w}{ }\PYG{o}{+}\PYG{o}{=}\PYG{+w}{ }\PYG{l+m+mf}{0.9}\PYG{+w}{ }\PYG{o}{*}\PYG{+w}{ }\PYG{n}{z}\PYG{p}{;}
\PYG{+w}{    }\PYG{n}{model}\PYG{p}{.}\PYG{n}{AddConstr}\PYG{p}{(}\PYG{n}{expr}\PYG{+w}{ }\PYG{o}{\PYGZgt{}}\PYG{o}{=}\PYG{+w}{ }\PYG{l+m+mf}{1.2}\PYG{p}{,}\PYG{+w}{ }\PYG{l+s}{\PYGZdq{}}\PYG{l+s}{R1}\PYG{l+s}{\PYGZdq{}}\PYG{p}{)}\PYG{p}{;}

\PYG{+w}{    }\PYG{c+c1}{// Set parameters}
\PYG{+w}{    }\PYG{n}{model}\PYG{p}{.}\PYG{n}{SetDblParam}\PYG{p}{(}\PYG{n}{COPT\PYGZus{}DBLPARAM\PYGZus{}TIMELIMIT}\PYG{p}{,}\PYG{+w}{ }\PYG{l+m+mi}{10}\PYG{p}{)}\PYG{p}{;}

\PYG{+w}{    }\PYG{c+c1}{// Solve problem}
\PYG{+w}{    }\PYG{n}{model}\PYG{p}{.}\PYG{n}{Solve}\PYG{p}{(}\PYG{p}{)}\PYG{p}{;}

\PYG{+w}{    }\PYG{c+c1}{// Output solution}
\PYG{+w}{    }\PYG{k}{if}\PYG{+w}{ }\PYG{p}{(}\PYG{n}{model}\PYG{p}{.}\PYG{n}{GetIntAttr}\PYG{p}{(}\PYG{n}{COPT\PYGZus{}INTATTR\PYGZus{}HASLPSOL}\PYG{p}{)}\PYG{+w}{ }\PYG{o}{!}\PYG{o}{=}\PYG{+w}{ }\PYG{l+m+mi}{0}\PYG{p}{)}
\PYG{+w}{    }\PYG{p}{\PYGZob{}}
\PYG{+w}{      }\PYG{n}{cout}\PYG{+w}{ }\PYG{o}{\PYGZlt{}}\PYG{o}{\PYGZlt{}}\PYG{+w}{ }\PYG{l+s}{\PYGZdq{}}\PYG{l+s+se}{\PYGZbs{}n}\PYG{l+s}{Found optimal solution:}\PYG{l+s}{\PYGZdq{}}\PYG{+w}{ }\PYG{o}{\PYGZlt{}}\PYG{o}{\PYGZlt{}}\PYG{+w}{ }\PYG{n}{endl}\PYG{p}{;}
\PYG{+w}{      }\PYG{n}{VarArray}\PYG{+w}{ }\PYG{n}{vars}\PYG{+w}{ }\PYG{o}{=}\PYG{+w}{ }\PYG{n}{model}\PYG{p}{.}\PYG{n}{GetVars}\PYG{p}{(}\PYG{p}{)}\PYG{p}{;}
\PYG{+w}{      }\PYG{k}{for}\PYG{+w}{ }\PYG{p}{(}\PYG{k+kt}{int}\PYG{+w}{ }\PYG{n}{i}\PYG{+w}{ }\PYG{o}{=}\PYG{+w}{ }\PYG{l+m+mi}{0}\PYG{p}{;}\PYG{+w}{ }\PYG{n}{i}\PYG{+w}{ }\PYG{o}{\PYGZlt{}}\PYG{+w}{ }\PYG{n}{vars}\PYG{p}{.}\PYG{n}{Size}\PYG{p}{(}\PYG{p}{)}\PYG{p}{;}\PYG{+w}{ }\PYG{n}{i}\PYG{o}{+}\PYG{o}{+}\PYG{p}{)}
\PYG{+w}{      }\PYG{p}{\PYGZob{}}
\PYG{+w}{        }\PYG{n}{Var}\PYG{+w}{ }\PYG{n}{var}\PYG{+w}{ }\PYG{o}{=}\PYG{+w}{ }\PYG{n}{vars}\PYG{p}{.}\PYG{n}{GetVar}\PYG{p}{(}\PYG{n}{i}\PYG{p}{)}\PYG{p}{;}
\PYG{+w}{        }\PYG{n}{cout}\PYG{+w}{ }\PYG{o}{\PYGZlt{}}\PYG{o}{\PYGZlt{}}\PYG{+w}{ }\PYG{l+s}{\PYGZdq{}}\PYG{l+s}{  }\PYG{l+s}{\PYGZdq{}}\PYG{+w}{ }\PYG{o}{\PYGZlt{}}\PYG{o}{\PYGZlt{}}\PYG{+w}{ }\PYG{n}{var}\PYG{p}{.}\PYG{n}{GetName}\PYG{p}{(}\PYG{p}{)}\PYG{+w}{ }\PYG{o}{\PYGZlt{}}\PYG{o}{\PYGZlt{}}\PYG{+w}{ }\PYG{l+s}{\PYGZdq{}}\PYG{l+s}{ = }\PYG{l+s}{\PYGZdq{}}\PYG{+w}{ }\PYG{o}{\PYGZlt{}}\PYG{o}{\PYGZlt{}}\PYG{+w}{ }\PYG{n}{var}\PYG{p}{.}\PYG{n}{Get}\PYG{p}{(}\PYG{n}{COPT\PYGZus{}DBLINFO\PYGZus{}VALUE}\PYG{p}{)}\PYG{+w}{ }\PYG{o}{\PYGZlt{}}\PYG{o}{\PYGZlt{}}\PYG{+w}{ }\PYG{n}{endl}\PYG{p}{;}
\PYG{+w}{      }\PYG{p}{\PYGZcb{}}
\PYG{+w}{      }\PYG{n}{cout}\PYG{+w}{ }\PYG{o}{\PYGZlt{}}\PYG{o}{\PYGZlt{}}\PYG{+w}{ }\PYG{l+s}{\PYGZdq{}}\PYG{l+s}{Obj = }\PYG{l+s}{\PYGZdq{}}\PYG{+w}{ }\PYG{o}{\PYGZlt{}}\PYG{o}{\PYGZlt{}}\PYG{+w}{ }\PYG{n}{model}\PYG{p}{.}\PYG{n}{GetDblAttr}\PYG{p}{(}\PYG{n}{COPT\PYGZus{}DBLATTR\PYGZus{}LPOBJVAL}\PYG{p}{)}\PYG{+w}{ }\PYG{o}{\PYGZlt{}}\PYG{o}{\PYGZlt{}}\PYG{+w}{ }\PYG{n}{endl}\PYG{p}{;}
\PYG{+w}{    }\PYG{p}{\PYGZcb{}}
\PYG{+w}{  }\PYG{p}{\PYGZcb{}}
\PYG{+w}{  }\PYG{k}{catch}\PYG{+w}{ }\PYG{p}{(}\PYG{n}{CoptException}\PYG{+w}{ }\PYG{n}{e}\PYG{p}{)}
\PYG{+w}{  }\PYG{p}{\PYGZob{}}
\PYG{+w}{    }\PYG{n}{cout}\PYG{+w}{ }\PYG{o}{\PYGZlt{}}\PYG{o}{\PYGZlt{}}\PYG{+w}{ }\PYG{l+s}{\PYGZdq{}}\PYG{l+s}{Error Code = }\PYG{l+s}{\PYGZdq{}}\PYG{+w}{ }\PYG{o}{\PYGZlt{}}\PYG{o}{\PYGZlt{}}\PYG{+w}{ }\PYG{n}{e}\PYG{p}{.}\PYG{n}{GetCode}\PYG{p}{(}\PYG{p}{)}\PYG{+w}{ }\PYG{o}{\PYGZlt{}}\PYG{o}{\PYGZlt{}}\PYG{+w}{ }\PYG{n}{endl}\PYG{p}{;}
\PYG{+w}{    }\PYG{n}{cout}\PYG{+w}{ }\PYG{o}{\PYGZlt{}}\PYG{o}{\PYGZlt{}}\PYG{+w}{ }\PYG{n}{e}\PYG{p}{.}\PYG{n}{what}\PYG{p}{(}\PYG{p}{)}\PYG{+w}{ }\PYG{o}{\PYGZlt{}}\PYG{o}{\PYGZlt{}}\PYG{+w}{ }\PYG{n}{endl}\PYG{p}{;}
\PYG{+w}{  }\PYG{p}{\PYGZcb{}}
\PYG{+w}{  }\PYG{k}{catch}\PYG{+w}{ }\PYG{p}{(}\PYG{p}{.}\PYG{p}{.}\PYG{p}{.}\PYG{p}{)}
\PYG{+w}{  }\PYG{p}{\PYGZob{}}
\PYG{+w}{    }\PYG{n}{cout}\PYG{+w}{ }\PYG{o}{\PYGZlt{}}\PYG{o}{\PYGZlt{}}\PYG{+w}{ }\PYG{l+s}{\PYGZdq{}}\PYG{l+s}{Unknown exception occurs!}\PYG{l+s}{\PYGZdq{}}\PYG{p}{;}
\PYG{+w}{  }\PYG{p}{\PYGZcb{}}
\PYG{p}{\PYGZcb{}}
\end{sphinxVerbatim}

\sphinxAtStartPar
Let’s now walk through the example, line by line, to understand how it
achieves the desired result of optimizing the model. Note that the example must
include header coptcpp\_pch.h.

\subsubsection{Creating environment and model}
\label{\detokenize{cppinterface:creating-environment-and-model}}
\sphinxAtStartPar
Essentially, any C++ application using Cardinal Optimizer should start with a
COPT environment, where user could add one or more models. Note that each model
encapsulates a problem and corresponding data.

\sphinxAtStartPar
Furthermore, to create multiple problems, one
can load them one by one in the same model, besides the naive option of creating
multiple models in the environment.

\begin{sphinxVerbatim}[commandchars=\\\{\}]
\PYG{+w}{    }\PYG{n}{Envr}\PYG{+w}{ }\PYG{n}{env}\PYG{p}{;}
\PYG{+w}{    }\PYG{n}{Model}\PYG{+w}{ }\PYG{n}{model}\PYG{+w}{ }\PYG{o}{=}\PYG{+w}{ }\PYG{n}{env}\PYG{p}{.}\PYG{n}{CreateModel}\PYG{p}{(}\PYG{l+s}{\PYGZdq{}}\PYG{l+s}{lp\PYGZus{}ex1}\PYG{l+s}{\PYGZdq{}}\PYG{p}{)}\PYG{p}{;}
\end{sphinxVerbatim}

\sphinxAtStartPar
The above call instantiates a COPT environment and a model with name “lp\_ex1”.

\subsubsection{Adding variables}
\label{\detokenize{cppinterface:adding-variables}}
\sphinxAtStartPar
The next step in our example is to add variables to the model.
Variables are added through AddVar() or AddVars() method on the model object.
A variable is always associated with a particular model.

\begin{sphinxVerbatim}[commandchars=\\\{\}]
\PYG{+w}{    }\PYG{c+c1}{// Add variables}
\PYG{+w}{    }\PYG{n}{Var}\PYG{+w}{ }\PYG{n}{x}\PYG{+w}{ }\PYG{o}{=}\PYG{+w}{ }\PYG{n}{model}\PYG{p}{.}\PYG{n}{AddVar}\PYG{p}{(}\PYG{l+m+mf}{0.1}\PYG{p}{,}\PYG{+w}{ }\PYG{l+m+mf}{0.6}\PYG{p}{,}\PYG{+w}{ }\PYG{l+m+mf}{0.0}\PYG{p}{,}\PYG{+w}{ }\PYG{n}{COPT\PYGZus{}CONTINUOUS}\PYG{p}{,}\PYG{+w}{ }\PYG{l+s}{\PYGZdq{}}\PYG{l+s}{x}\PYG{l+s}{\PYGZdq{}}\PYG{p}{)}\PYG{p}{;}
\PYG{+w}{    }\PYG{n}{Var}\PYG{+w}{ }\PYG{n}{y}\PYG{+w}{ }\PYG{o}{=}\PYG{+w}{ }\PYG{n}{model}\PYG{p}{.}\PYG{n}{AddVar}\PYG{p}{(}\PYG{l+m+mf}{0.2}\PYG{p}{,}\PYG{+w}{ }\PYG{l+m+mf}{1.5}\PYG{p}{,}\PYG{+w}{ }\PYG{l+m+mf}{0.0}\PYG{p}{,}\PYG{+w}{ }\PYG{n}{COPT\PYGZus{}CONTINUOUS}\PYG{p}{,}\PYG{+w}{ }\PYG{l+s}{\PYGZdq{}}\PYG{l+s}{y}\PYG{l+s}{\PYGZdq{}}\PYG{p}{)}\PYG{p}{;}
\PYG{+w}{    }\PYG{n}{Var}\PYG{+w}{ }\PYG{n}{z}\PYG{+w}{ }\PYG{o}{=}\PYG{+w}{ }\PYG{n}{model}\PYG{p}{.}\PYG{n}{AddVar}\PYG{p}{(}\PYG{l+m+mf}{0.3}\PYG{p}{,}\PYG{+w}{ }\PYG{l+m+mf}{2.8}\PYG{p}{,}\PYG{+w}{ }\PYG{l+m+mf}{0.0}\PYG{p}{,}\PYG{+w}{ }\PYG{n}{COPT\PYGZus{}CONTINUOUS}\PYG{p}{,}\PYG{+w}{ }\PYG{l+s}{\PYGZdq{}}\PYG{l+s}{z}\PYG{l+s}{\PYGZdq{}}\PYG{p}{)}\PYG{p}{;}
\end{sphinxVerbatim}

\sphinxAtStartPar
The first and second arguments to the AddVar() call are the variable lower and
upper bounds, respectively. The third argument is the linear objective coefficient
(zero here \sphinxhyphen{} we’ll set the objective later). The fourth argument is the variable
type. Our variables are all continuous in this example. The final argument is the name
of the variable.

\sphinxAtStartPar
The AddVar() method has been overloaded to accept several different argument lists.
Please refer to {\hyperref[\detokenize{cppapiref:chapcppapiref}]{\sphinxcrossref{\DUrole{std,std-ref}{C++ API Reference}}}} for further details.

\sphinxAtStartPar
The objective is built here using overloaded operators. The C++ API overloads the
arithmetic operators to allow you to build linear expressions by COPT variables.
The second argument indicates that the sense is maximization.

\subsubsection{Adding constraints}
\label{\detokenize{cppinterface:adding-constraints}}
\sphinxAtStartPar
The next step in the example is to add the linear constraints. As with variables,
constraints are always associated with a specific model. They are created using
AddConstr() or AddConstrs() methods on the model object.

\begin{sphinxVerbatim}[commandchars=\\\{\}]
\PYG{+w}{    }\PYG{c+c1}{// Add linear constraints using linear expression}
\PYG{+w}{    }\PYG{n}{model}\PYG{p}{.}\PYG{n}{AddConstr}\PYG{p}{(}\PYG{l+m+mf}{1.5}\PYG{+w}{ }\PYG{o}{*}\PYG{+w}{ }\PYG{n}{x}\PYG{+w}{ }\PYG{o}{+}\PYG{+w}{ }\PYG{l+m+mf}{1.2}\PYG{+w}{ }\PYG{o}{*}\PYG{+w}{ }\PYG{n}{y}\PYG{+w}{ }\PYG{o}{+}\PYG{+w}{ }\PYG{l+m+mf}{1.8}\PYG{+w}{ }\PYG{o}{*}\PYG{+w}{ }\PYG{n}{z}\PYG{+w}{ }\PYG{o}{\PYGZlt{}}\PYG{o}{=}\PYG{+w}{ }\PYG{l+m+mf}{2.6}\PYG{p}{,}\PYG{+w}{ }\PYG{l+s}{\PYGZdq{}}\PYG{l+s}{R0}\PYG{l+s}{\PYGZdq{}}\PYG{p}{)}\PYG{p}{;}

\PYG{+w}{    }\PYG{n}{Expr}\PYG{+w}{ }\PYG{n+nf}{expr}\PYG{p}{(}\PYG{n}{x}\PYG{p}{,}\PYG{+w}{ }\PYG{l+m+mf}{0.8}\PYG{p}{)}\PYG{p}{;}
\PYG{+w}{    }\PYG{n}{expr}\PYG{p}{.}\PYG{n}{AddTerm}\PYG{p}{(}\PYG{n}{y}\PYG{p}{,}\PYG{+w}{ }\PYG{l+m+mf}{0.6}\PYG{p}{)}\PYG{p}{;}
\PYG{+w}{    }\PYG{n}{expr}\PYG{+w}{ }\PYG{o}{+}\PYG{o}{=}\PYG{+w}{ }\PYG{l+m+mf}{0.9}\PYG{+w}{ }\PYG{o}{*}\PYG{+w}{ }\PYG{n}{z}\PYG{p}{;}
\PYG{+w}{    }\PYG{n}{model}\PYG{p}{.}\PYG{n}{AddConstr}\PYG{p}{(}\PYG{n}{expr}\PYG{+w}{ }\PYG{o}{\PYGZgt{}}\PYG{o}{=}\PYG{+w}{ }\PYG{l+m+mf}{1.2}\PYG{p}{,}\PYG{+w}{ }\PYG{l+s}{\PYGZdq{}}\PYG{l+s}{R1}\PYG{l+s}{\PYGZdq{}}\PYG{p}{)}\PYG{p}{;}
\end{sphinxVerbatim}

\sphinxAtStartPar
The first constraint is to use overloaded arithmetic operators to build the
linear expression. The comparison operators are also overloaded to make it easier
to build constraints.

\sphinxAtStartPar
The second constraint is created by building a linear expression incrementally.
That is, an expression can be built by constructor of a variable and its coefficient,
by AddTerm() method, and by overloaded operators.

\subsubsection{Setting parameters and attributes}
\label{\detokenize{cppinterface:setting-parameters-and-attributes}}
\sphinxAtStartPar
The next step in the example is to set parameters and attributes of
the problem before optimization.

\begin{sphinxVerbatim}[commandchars=\\\{\}]
\PYG{+w}{    }\PYG{c+c1}{// Set parameters}
\PYG{+w}{    }\PYG{n}{model}\PYG{p}{.}\PYG{n}{SetDblParam}\PYG{p}{(}\PYG{n}{COPT\PYGZus{}DBLPARAM\PYGZus{}TIMELIMIT}\PYG{p}{,}\PYG{+w}{ }\PYG{l+m+mi}{10}\PYG{p}{)}\PYG{p}{;}
\end{sphinxVerbatim}

\sphinxAtStartPar
The SetDblParam() call here with COPT\_DBLPARAM\_TIMELIMIT argument sets solver to
optimize up to 10 seconds.

\subsubsection{Solving problem}
\label{\detokenize{cppinterface:solving-problem}}
\sphinxAtStartPar
Now that the model has been built, the next step is to optimize it:

\begin{sphinxVerbatim}[commandchars=\\\{\}]
\PYG{+w}{    }\PYG{c+c1}{// Solve problem}
\PYG{+w}{    }\PYG{n}{model}\PYG{p}{.}\PYG{n}{Solve}\PYG{p}{(}\PYG{p}{)}\PYG{p}{;}
\end{sphinxVerbatim}

\sphinxAtStartPar
This routine performs the optimization and populates several internal model attributes
(including the status of the optimization, the solution, etc.).

\subsubsection{Outputting solution}
\label{\detokenize{cppinterface:outputting-solution}}
\sphinxAtStartPar
After solving the problem, one can query the values of the attributes for
various of purposes.

\begin{sphinxVerbatim}[commandchars=\\\{\}]
\PYG{+w}{    }\PYG{c+c1}{// Output solution}
\PYG{+w}{    }\PYG{k}{if}\PYG{+w}{ }\PYG{p}{(}\PYG{n}{model}\PYG{p}{.}\PYG{n}{GetIntAttr}\PYG{p}{(}\PYG{n}{COPT\PYGZus{}INTATTR\PYGZus{}HASLPSOL}\PYG{p}{)}\PYG{+w}{ }\PYG{o}{!}\PYG{o}{=}\PYG{+w}{ }\PYG{l+m+mi}{0}\PYG{p}{)}
\PYG{+w}{    }\PYG{p}{\PYGZob{}}
\PYG{+w}{      }\PYG{n}{cout}\PYG{+w}{ }\PYG{o}{\PYGZlt{}}\PYG{o}{\PYGZlt{}}\PYG{+w}{ }\PYG{l+s}{\PYGZdq{}}\PYG{l+s+se}{\PYGZbs{}n}\PYG{l+s}{Found optimal solution:}\PYG{l+s}{\PYGZdq{}}\PYG{+w}{ }\PYG{o}{\PYGZlt{}}\PYG{o}{\PYGZlt{}}\PYG{+w}{ }\PYG{n}{endl}\PYG{p}{;}
\PYG{+w}{      }\PYG{n}{VarArray}\PYG{+w}{ }\PYG{n}{vars}\PYG{+w}{ }\PYG{o}{=}\PYG{+w}{ }\PYG{n}{model}\PYG{p}{.}\PYG{n}{GetVars}\PYG{p}{(}\PYG{p}{)}\PYG{p}{;}
\PYG{+w}{      }\PYG{k}{for}\PYG{+w}{ }\PYG{p}{(}\PYG{k+kt}{int}\PYG{+w}{ }\PYG{n}{i}\PYG{+w}{ }\PYG{o}{=}\PYG{+w}{ }\PYG{l+m+mi}{0}\PYG{p}{;}\PYG{+w}{ }\PYG{n}{i}\PYG{+w}{ }\PYG{o}{\PYGZlt{}}\PYG{+w}{ }\PYG{n}{vars}\PYG{p}{.}\PYG{n}{Size}\PYG{p}{(}\PYG{p}{)}\PYG{p}{;}\PYG{+w}{ }\PYG{n}{i}\PYG{o}{+}\PYG{o}{+}\PYG{p}{)}
\PYG{+w}{      }\PYG{p}{\PYGZob{}}
\PYG{+w}{        }\PYG{n}{Var}\PYG{+w}{ }\PYG{n}{var}\PYG{+w}{ }\PYG{o}{=}\PYG{+w}{ }\PYG{n}{vars}\PYG{p}{.}\PYG{n}{GetVar}\PYG{p}{(}\PYG{n}{i}\PYG{p}{)}\PYG{p}{;}
\PYG{+w}{        }\PYG{n}{cout}\PYG{+w}{ }\PYG{o}{\PYGZlt{}}\PYG{o}{\PYGZlt{}}\PYG{+w}{ }\PYG{l+s}{\PYGZdq{}}\PYG{l+s}{  }\PYG{l+s}{\PYGZdq{}}\PYG{+w}{ }\PYG{o}{\PYGZlt{}}\PYG{o}{\PYGZlt{}}\PYG{+w}{ }\PYG{n}{var}\PYG{p}{.}\PYG{n}{GetName}\PYG{p}{(}\PYG{p}{)}\PYG{+w}{ }\PYG{o}{\PYGZlt{}}\PYG{o}{\PYGZlt{}}\PYG{+w}{ }\PYG{l+s}{\PYGZdq{}}\PYG{l+s}{ = }\PYG{l+s}{\PYGZdq{}}\PYG{+w}{ }\PYG{o}{\PYGZlt{}}\PYG{o}{\PYGZlt{}}\PYG{+w}{ }\PYG{n}{var}\PYG{p}{.}\PYG{n}{Get}\PYG{p}{(}\PYG{n}{COPT\PYGZus{}DBLINFO\PYGZus{}VALUE}\PYG{p}{)}\PYG{+w}{ }\PYG{o}{\PYGZlt{}}\PYG{o}{\PYGZlt{}}\PYG{+w}{ }\PYG{n}{endl}\PYG{p}{;}
\PYG{+w}{      }\PYG{p}{\PYGZcb{}}
\PYG{+w}{      }\PYG{n}{cout}\PYG{+w}{ }\PYG{o}{\PYGZlt{}}\PYG{o}{\PYGZlt{}}\PYG{+w}{ }\PYG{l+s}{\PYGZdq{}}\PYG{l+s}{Obj = }\PYG{l+s}{\PYGZdq{}}\PYG{+w}{ }\PYG{o}{\PYGZlt{}}\PYG{o}{\PYGZlt{}}\PYG{+w}{ }\PYG{n}{model}\PYG{p}{.}\PYG{n}{GetDblAttr}\PYG{p}{(}\PYG{n}{COPT\PYGZus{}DBLATTR\PYGZus{}LPOBJVAL}\PYG{p}{)}\PYG{+w}{ }\PYG{o}{\PYGZlt{}}\PYG{o}{\PYGZlt{}}\PYG{+w}{ }\PYG{n}{endl}\PYG{p}{;}
\PYG{+w}{    }\PYG{p}{\PYGZcb{}}
\end{sphinxVerbatim}

\sphinxAtStartPar
Spcifically, one can query the COPT\_INTATTR\_HASLPSOL attribute on the model to know whether
we have optimal LP solution; query the COPT\_DBLINFO\_VALUE attribute of a variable to obtain its
solution value; query the COPT\_DBLATTR\_LPOBJVAL attribute on the model to obtain the
objective value for the current solution.

\sphinxAtStartPar
The names and types of all model, variable, and constraint attributes can be found in
{\hyperref[\detokenize{capiref:chapapi-attrs}]{\sphinxcrossref{\DUrole{std,std-ref}{Attributes}}}} of C API reference.

\subsubsection{Error handling}
\label{\detokenize{cppinterface:error-handling}}
\sphinxAtStartPar
Errors in the COPT C++ interface are handled through the C++ exception mechanism.
In the example, all COPT statements are enclosed inside a try block, and any associated
errors would be caught by the catch block.

\begin{sphinxVerbatim}[commandchars=\\\{\}]
\PYG{+w}{  }\PYG{k}{catch}\PYG{+w}{ }\PYG{p}{(}\PYG{n}{CoptException}\PYG{+w}{ }\PYG{n}{e}\PYG{p}{)}
\PYG{+w}{  }\PYG{p}{\PYGZob{}}
\PYG{+w}{    }\PYG{n}{cout}\PYG{+w}{ }\PYG{o}{\PYGZlt{}}\PYG{o}{\PYGZlt{}}\PYG{+w}{ }\PYG{l+s}{\PYGZdq{}}\PYG{l+s}{Error Code = }\PYG{l+s}{\PYGZdq{}}\PYG{+w}{ }\PYG{o}{\PYGZlt{}}\PYG{o}{\PYGZlt{}}\PYG{+w}{ }\PYG{n}{e}\PYG{p}{.}\PYG{n}{GetCode}\PYG{p}{(}\PYG{p}{)}\PYG{+w}{ }\PYG{o}{\PYGZlt{}}\PYG{o}{\PYGZlt{}}\PYG{+w}{ }\PYG{n}{endl}\PYG{p}{;}
\PYG{+w}{    }\PYG{n}{cout}\PYG{+w}{ }\PYG{o}{\PYGZlt{}}\PYG{o}{\PYGZlt{}}\PYG{+w}{ }\PYG{n}{e}\PYG{p}{.}\PYG{n}{what}\PYG{p}{(}\PYG{p}{)}\PYG{+w}{ }\PYG{o}{\PYGZlt{}}\PYG{o}{\PYGZlt{}}\PYG{+w}{ }\PYG{n}{endl}\PYG{p}{;}
\PYG{+w}{  }\PYG{p}{\PYGZcb{}}
\PYG{+w}{  }\PYG{k}{catch}\PYG{+w}{ }\PYG{p}{(}\PYG{p}{.}\PYG{p}{.}\PYG{p}{.}\PYG{p}{)}
\PYG{+w}{  }\PYG{p}{\PYGZob{}}
\PYG{+w}{    }\PYG{n}{cout}\PYG{+w}{ }\PYG{o}{\PYGZlt{}}\PYG{o}{\PYGZlt{}}\PYG{+w}{ }\PYG{l+s}{\PYGZdq{}}\PYG{l+s}{Unknown exception occurs!}\PYG{l+s}{\PYGZdq{}}\PYG{p}{;}
\PYG{+w}{  }\PYG{p}{\PYGZcb{}}
\end{sphinxVerbatim}

\subsection{Build and Run}
\label{\detokenize{cppinterface:build-and-run}}
\sphinxAtStartPar
To build and run the example, users may refer to files under
\sphinxcode{\sphinxupquote{\$COPT\_HOME/examples/cpp}}. Specifically, We provide visual studio
project on Windows, as well as makefile project on Linux and Mac platforms.

\subsubsection{Windows}
\label{\detokenize{cppinterface:windows}}
\sphinxAtStartPar
\sphinxstylestrong{1. Visual Studio project}

\sphinxAtStartPar
For Windows platform, Visual Studio project is located at
\sphinxcode{\sphinxupquote{\$COPT\_HOME/examples/cpp/vsprojects}}. Double\sphinxhyphen{}clicking the project file \sphinxcode{\sphinxupquote{lp\_ex1.vcxproj}}
will bring Visual Studio. Note that it requires Visual studio 2017 or 2019 installed
on Windows 10 to build.

\sphinxAtStartPar
The Visual Studio project has dependency on COPT cpp shared libary \sphinxcode{\sphinxupquote{copt\_cpp.dll}}, which is
refered in project file, along with its import library \sphinxcode{\sphinxupquote{copt\_cpp.lib}}.
The required headers are \sphinxcode{\sphinxupquote{copt.h}}, \sphinxcode{\sphinxupquote{coptcpp.h}} and \sphinxcode{\sphinxupquote{coptcpp.idl.h}}, which declare COPT
constants, interfaces and methods exported from \sphinxcode{\sphinxupquote{copt\_cpp.dll}}. In addition,
the example provides class header files under \sphinxcode{\sphinxupquote{\$COPT\_HOME/include/coptcpp\_inc}},
which wraps COPT cpp interfaces and redefines
overloaded operators.

\sphinxAtStartPar
In simple terms, users only need to include the header file \sphinxcode{\sphinxupquote{coptcpp\_pch.h}} as shown in the example, configure additional dependencies as \sphinxcode{\sphinxupquote{copt\_cpp.lib}}, and set the directory of additional link libraries as \sphinxcode{\sphinxupquote{\$(COPT\_HOME)/lib}}, and make sure that the dynamic library \sphinxcode{\sphinxupquote{copt\_cpp.dll}} has been installed in the appropriate path, and can be loaded at runtime.

\sphinxAtStartPar
To run the example, users should have COPT installed. Specifically, it requires
COPT cpp libary, \sphinxcode{\sphinxupquote{copt\_cpp.dll}}, and valid license files to run. Please refer to
{\hyperref[\detokenize{install:chapinstall}]{\sphinxcrossref{\DUrole{std,std-ref}{Install Guide for Cardinal Optimizer}}}} for further details.

\sphinxAtStartPar
In addition, on Windows systems, the Clang and Intel compilers, like the aforementioned MSVC compiler,
all support compiling the C++ examples of COPT.

\sphinxAtStartPar
\sphinxstylestrong{2. MinGW\sphinxhyphen{}w64 Toolchain}

\sphinxAtStartPar
For Windows systems, COPT also supports the MinGW\sphinxhyphen{}w64 toolchain. For the two variants of the MinGW\sphinxhyphen{}w64, we provide
the dynamic library \sphinxcode{\sphinxupquote{copt\_mmingwcpp.dll}} supporting MSVCRT and the dynamic library \sphinxcode{\sphinxupquote{copt\_umingwcpp.dll}} supporting UCRT.
Please refer to the corresponding instructions in the Makefile under the path \sphinxcode{\sphinxupquote{\$COPT\_HOME/examples/cpp}} for usage.

\sphinxAtStartPar
We recommend using \sphinxcode{\sphinxupquote{copt\_mmingwcpp.dll}} supporting MSVCRT by default, and trying the UCRT if any issues arise.

\subsubsection{Linux and MacOS}
\label{\detokenize{cppinterface:linux-and-macos}}
\sphinxAtStartPar
We provide Makefile to build the example for Linux and MacOS platforms. Please
make sure tools, \sphinxcode{\sphinxupquote{gcc}} and \sphinxcode{\sphinxupquote{make}}, are already installed on the platforms.
To build the example, change directory to \sphinxcode{\sphinxupquote{\$COPT\_HOME/examples/cpp}} and
execute command \sphinxcode{\sphinxupquote{\textquotesingle{}make\textquotesingle{}}} on unix terminal.

\sphinxAtStartPar
The project has dependency on COPT cpp shared libary, that is, \sphinxcode{\sphinxupquote{libcopt\_cpp.so}}
on Linux platform and \sphinxcode{\sphinxupquote{libcopt\_cpp.dylib}} on MacOS. Similar to Windows VS project,
user should refer to \sphinxcode{\sphinxupquote{coptcpp\_pch.h}} under \sphinxcode{\sphinxupquote{\$COPT\_HOME/include/coptcpp\_inc\textquotesingle{}}} to
include all necessary headers, as shown in the example.

\sphinxAtStartPar
To run the example, users should have COPT installed. Specifically, it requires
COPT cpp libary, \sphinxcode{\sphinxupquote{libcopt\_cpp.so}} on Linux, \sphinxcode{\sphinxupquote{libcopt\_cpp.dylib}} on MacOS,
and valid license files to run. Alternatively, user might set \sphinxcode{\sphinxupquote{LD\_LIBRARY\_PATH}}
and \sphinxcode{\sphinxupquote{COPT\_LICENSE\_DIR}} properly to work around. Please refer to
{\hyperref[\detokenize{install:chapinstall}]{\sphinxcrossref{\DUrole{std,std-ref}{Install Guide for Cardinal Optimizer}}}} for further details.

\sphinxstepscope

\section{C\# Interface}
\label{\detokenize{csharpinterface:c-interface}}\label{\detokenize{csharpinterface:chapcsharpinterface}}\label{\detokenize{csharpinterface::doc}}
\sphinxAtStartPar
This chapter walks through a simple C\# example to illustrate the use of the
COPT C\# interface. In short words, the example creates an environment, builds
a model, add variables and constraints, optimizes it, and then outputs the
optimal objective value.

\sphinxAtStartPar
The example solves the following linear problem:
\begin{equation}\label{equation:csharpinterface:coptcsharpEq_capilp1}
\begin{split}\text{Maximize: } & \\
                  & 1.2 x + 1.8 y + 2.1 z \\
\text{Subject to: } & \\
                    & 1.5 x + 1.2 y + 1.8 z \leq 2.6 \\
                    & 0.8 x + 0.6 y + 0.9 z \geq 1.2 \\
\text{Bounds: } & \\
               & 0.1 \leq x \leq 0.6 \\
               & 0.2 \leq y \leq 1.5 \\
               & 0.3 \leq z \leq 2.8 \\\end{split}
\end{equation}
\sphinxAtStartPar
Note that this is the same problem that was modelled and optimized in
chapter of {\hyperref[\detokenize{cinterface:chapcinterface}]{\sphinxcrossref{\DUrole{std,std-ref}{C Interface}}}}.

\subsection{Example details}
\label{\detokenize{csharpinterface:example-details}}
\sphinxAtStartPar
Below is the source code solving the above problem using COPT C\# interface.
\sphinxSetupCaptionForVerbatim{\sphinxcode{\sphinxupquote{lp\_ex1.cs}}}
\def\sphinxLiteralBlockLabel{\label{\detokenize{csharpinterface:id1}}}
\begin{sphinxVerbatim}[commandchars=\\\{\},numbers=left,firstnumber=1,stepnumber=1]
\PYG{c+cm}{/*}
\PYG{c+cm}{ * This file is part of the Cardinal Optimizer, all rights reserved.}
\PYG{c+cm}{ */}
\PYG{k}{using}\PYG{+w}{ }\PYG{n+nn}{Copt}\PYG{p}{;}
\PYG{k}{using}\PYG{+w}{ }\PYG{n+nn}{System}\PYG{p}{;}

\PYG{c+cm}{/*}
\PYG{c+cm}{ * This C\PYGZsh{} example solves the following LP model:}
\PYG{c+cm}{ *}
\PYG{c+cm}{ * }
\PYG{c+cm}{ * Maximize:}
\PYG{c+cm}{ *  1.2 x + 1.8 y + 2.1 z}
\PYG{c+cm}{ *}
\PYG{c+cm}{ * Subject to:}
\PYG{c+cm}{ *  1.5 x + 1.2 y + 1.8 z \PYGZlt{}= 2.6}
\PYG{c+cm}{ *  0.8 x + 0.6 y + 0.9 z \PYGZgt{}= 1.2}

\PYG{c+cm}{ * where:}
\PYG{c+cm}{ *  0.1 \PYGZlt{}= x \PYGZlt{}= 0.6}
\PYG{c+cm}{ *  0.2 \PYGZlt{}= y \PYGZlt{}= 1.5}
\PYG{c+cm}{ *  0.3 \PYGZlt{}= z \PYGZlt{}= 2.8}
\PYG{c+cm}{ */}
\PYG{k}{public}\PYG{+w}{ }\PYG{k}{class}\PYG{+w}{ }\PYG{n+nc}{lp\PYGZus{}ex1}
\PYG{p}{\PYGZob{}}
\PYG{+w}{  }\PYG{k}{public}\PYG{+w}{ }\PYG{k}{static}\PYG{+w}{ }\PYG{k}{void}\PYG{+w}{ }\PYG{n+nf}{Main}\PYG{p}{(}\PYG{p}{)}
\PYG{+w}{  }\PYG{p}{\PYGZob{}}
\PYG{+w}{    }\PYG{k}{try}
\PYG{+w}{    }\PYG{p}{\PYGZob{}}
\PYG{+w}{      }\PYG{n}{Envr}\PYG{+w}{ }\PYG{n}{env}\PYG{+w}{ }\PYG{o}{=}\PYG{+w}{ }\PYG{k}{new}\PYG{+w}{ }\PYG{n}{Envr}\PYG{p}{(}\PYG{p}{)}\PYG{p}{;}
\PYG{+w}{      }\PYG{n}{Model}\PYG{+w}{ }\PYG{n}{model}\PYG{+w}{ }\PYG{o}{=}\PYG{+w}{ }\PYG{n}{env}\PYG{p}{.}\PYG{n}{CreateModel}\PYG{p}{(}\PYG{l+s}{\PYGZdq{}lp\PYGZus{}ex1\PYGZdq{}}\PYG{p}{)}\PYG{p}{;}

\PYG{+w}{      }\PYG{c+cm}{/* }
\PYG{c+cm}{       * Add variables x, y, z}
\PYG{c+cm}{       *}
\PYG{c+cm}{       * obj: 1.2 x + 1.8 y + 2.1 z}
\PYG{c+cm}{       *}
\PYG{c+cm}{       * var:}
\PYG{c+cm}{       *  0.1 \PYGZlt{}= x \PYGZlt{}= 0.6}
\PYG{c+cm}{       *  0.2 \PYGZlt{}= y \PYGZlt{}= 1.5}
\PYG{c+cm}{       *  0.3 \PYGZlt{}= z \PYGZlt{}= 2.8}
\PYG{c+cm}{       */}
\PYG{+w}{      }\PYG{n}{Var}\PYG{+w}{ }\PYG{n}{x}\PYG{+w}{ }\PYG{o}{=}\PYG{+w}{ }\PYG{n}{model}\PYG{p}{.}\PYG{n}{AddVar}\PYG{p}{(}\PYG{l+m+mf}{0.1}\PYG{p}{,}\PYG{+w}{ }\PYG{l+m+mf}{0.6}\PYG{p}{,}\PYG{+w}{ }\PYG{l+m+mf}{0.0}\PYG{p}{,}\PYG{+w}{ }\PYG{n}{Copt}\PYG{p}{.}\PYG{n}{Consts}\PYG{p}{.}\PYG{n}{CONTINUOUS}\PYG{p}{,}\PYG{+w}{ }\PYG{l+s}{\PYGZdq{}x\PYGZdq{}}\PYG{p}{)}\PYG{p}{;}
\PYG{+w}{      }\PYG{n}{Var}\PYG{+w}{ }\PYG{n}{y}\PYG{+w}{ }\PYG{o}{=}\PYG{+w}{ }\PYG{n}{model}\PYG{p}{.}\PYG{n}{AddVar}\PYG{p}{(}\PYG{l+m+mf}{0.2}\PYG{p}{,}\PYG{+w}{ }\PYG{l+m+mf}{1.5}\PYG{p}{,}\PYG{+w}{ }\PYG{l+m+mf}{0.0}\PYG{p}{,}\PYG{+w}{ }\PYG{n}{Copt}\PYG{p}{.}\PYG{n}{Consts}\PYG{p}{.}\PYG{n}{CONTINUOUS}\PYG{p}{,}\PYG{+w}{ }\PYG{l+s}{\PYGZdq{}y\PYGZdq{}}\PYG{p}{)}\PYG{p}{;}
\PYG{+w}{      }\PYG{n}{Var}\PYG{+w}{ }\PYG{n}{z}\PYG{+w}{ }\PYG{o}{=}\PYG{+w}{ }\PYG{n}{model}\PYG{p}{.}\PYG{n}{AddVar}\PYG{p}{(}\PYG{l+m+mf}{0.3}\PYG{p}{,}\PYG{+w}{ }\PYG{l+m+mf}{2.8}\PYG{p}{,}\PYG{+w}{ }\PYG{l+m+mf}{0.0}\PYG{p}{,}\PYG{+w}{ }\PYG{n}{Copt}\PYG{p}{.}\PYG{n}{Consts}\PYG{p}{.}\PYG{n}{CONTINUOUS}\PYG{p}{,}\PYG{+w}{ }\PYG{l+s}{\PYGZdq{}z\PYGZdq{}}\PYG{p}{)}\PYG{p}{;}

\PYG{+w}{      }\PYG{n}{model}\PYG{p}{.}\PYG{n}{SetObjective}\PYG{p}{(}\PYG{l+m+mf}{1.2}\PYG{+w}{ }\PYG{o}{*}\PYG{+w}{ }\PYG{n}{x}\PYG{+w}{ }\PYG{o}{+}\PYG{+w}{ }\PYG{l+m+mf}{1.8}\PYG{+w}{ }\PYG{o}{*}\PYG{+w}{ }\PYG{n}{y}\PYG{+w}{ }\PYG{o}{+}\PYG{+w}{ }\PYG{l+m+mf}{2.1}\PYG{+w}{ }\PYG{o}{*}\PYG{+w}{ }\PYG{n}{z}\PYG{p}{,}\PYG{+w}{ }\PYG{n}{Copt}\PYG{p}{.}\PYG{n}{Consts}\PYG{p}{.}\PYG{n}{MAXIMIZE}\PYG{p}{)}\PYG{p}{;}

\PYG{+w}{      }\PYG{c+cm}{/*}
\PYG{c+cm}{       * Add two constraints using linear expression}
\PYG{c+cm}{       *}
\PYG{c+cm}{       * r0: 1.5 x + 1.2 y + 1.8 z \PYGZlt{}= 2.6}
\PYG{c+cm}{       * r1: 0.8 x + 0.6 y + 0.9 z \PYGZgt{}= 1.2}
\PYG{c+cm}{       */}
\PYG{+w}{      }\PYG{n}{model}\PYG{p}{.}\PYG{n}{AddConstr}\PYG{p}{(}\PYG{l+m+mf}{1.5}\PYG{+w}{ }\PYG{o}{*}\PYG{+w}{ }\PYG{n}{x}\PYG{+w}{ }\PYG{o}{+}\PYG{+w}{ }\PYG{l+m+mf}{1.2}\PYG{+w}{ }\PYG{o}{*}\PYG{+w}{ }\PYG{n}{y}\PYG{+w}{ }\PYG{o}{+}\PYG{+w}{ }\PYG{l+m+mf}{1.8}\PYG{+w}{ }\PYG{o}{*}\PYG{+w}{ }\PYG{n}{z}\PYG{+w}{ }\PYG{o}{\PYGZlt{}=}\PYG{+w}{ }\PYG{l+m+mf}{2.6}\PYG{p}{,}\PYG{+w}{ }\PYG{l+s}{\PYGZdq{}r0\PYGZdq{}}\PYG{p}{)}\PYG{p}{;}

\PYG{+w}{      }\PYG{n}{Expr}\PYG{+w}{ }\PYG{n}{expr}\PYG{+w}{ }\PYG{o}{=}\PYG{+w}{ }\PYG{k}{new}\PYG{+w}{ }\PYG{n}{Expr}\PYG{p}{(}\PYG{n}{x}\PYG{p}{,}\PYG{+w}{ }\PYG{l+m+mf}{0.8}\PYG{p}{)}\PYG{p}{;}
\PYG{+w}{      }\PYG{n}{expr}\PYG{p}{.}\PYG{n}{AddTerm}\PYG{p}{(}\PYG{n}{y}\PYG{p}{,}\PYG{+w}{ }\PYG{l+m+mf}{0.6}\PYG{p}{)}\PYG{p}{;}
\PYG{+w}{      }\PYG{n}{expr}\PYG{+w}{ }\PYG{o}{+=}\PYG{+w}{ }\PYG{l+m+mf}{0.9}\PYG{+w}{ }\PYG{o}{*}\PYG{+w}{ }\PYG{n}{z}\PYG{p}{;}
\PYG{+w}{      }\PYG{n}{model}\PYG{p}{.}\PYG{n}{AddConstr}\PYG{p}{(}\PYG{n}{expr}\PYG{+w}{ }\PYG{o}{\PYGZgt{}=}\PYG{+w}{ }\PYG{l+m+mf}{1.2}\PYG{p}{,}\PYG{+w}{ }\PYG{l+s}{\PYGZdq{}r1\PYGZdq{}}\PYG{p}{)}\PYG{p}{;}

\PYG{+w}{      }\PYG{c+c1}{// Set parameters}
\PYG{+w}{      }\PYG{n}{model}\PYG{p}{.}\PYG{n}{SetDblParam}\PYG{p}{(}\PYG{n}{Copt}\PYG{p}{.}\PYG{n}{DblParam}\PYG{p}{.}\PYG{n}{TimeLimit}\PYG{p}{,}\PYG{+w}{ }\PYG{l+m+mi}{10}\PYG{p}{)}\PYG{p}{;}

\PYG{+w}{      }\PYG{c+c1}{// Solve problem}
\PYG{+w}{      }\PYG{n}{model}\PYG{p}{.}\PYG{n}{Solve}\PYG{p}{(}\PYG{p}{)}\PYG{p}{;}

\PYG{+w}{      }\PYG{c+c1}{// Output solution}
\PYG{+w}{      }\PYG{k}{if}\PYG{+w}{ }\PYG{p}{(}\PYG{n}{model}\PYG{p}{.}\PYG{n}{GetIntAttr}\PYG{p}{(}\PYG{n}{Copt}\PYG{p}{.}\PYG{n}{IntAttr}\PYG{p}{.}\PYG{n}{LpStatus}\PYG{p}{)}\PYG{+w}{ }\PYG{o}{==}\PYG{+w}{ }\PYG{n}{Copt}\PYG{p}{.}\PYG{n}{Status}\PYG{p}{.}\PYG{n}{OPTIMAL}\PYG{p}{)}
\PYG{+w}{      }\PYG{p}{\PYGZob{}}
\PYG{+w}{        }\PYG{n}{Console}\PYG{p}{.}\PYG{n}{WriteLine}\PYG{p}{(}\PYG{l+s}{\PYGZdq{}\PYGZbs{}nFound optimal solution:\PYGZdq{}}\PYG{p}{)}\PYG{p}{;}
\PYG{+w}{        }\PYG{n}{VarArray}\PYG{+w}{ }\PYG{n}{vars}\PYG{+w}{ }\PYG{o}{=}\PYG{+w}{ }\PYG{n}{model}\PYG{p}{.}\PYG{n}{GetVars}\PYG{p}{(}\PYG{p}{)}\PYG{p}{;}
\PYG{+w}{        }\PYG{k}{for}\PYG{+w}{ }\PYG{p}{(}\PYG{k+kt}{int}\PYG{+w}{ }\PYG{n}{i}\PYG{+w}{ }\PYG{o}{=}\PYG{+w}{ }\PYG{l+m+mi}{0}\PYG{p}{;}\PYG{+w}{ }\PYG{n}{i}\PYG{+w}{ }\PYG{o}{\PYGZlt{}}\PYG{+w}{ }\PYG{n}{vars}\PYG{p}{.}\PYG{n}{Size}\PYG{p}{(}\PYG{p}{)}\PYG{p}{;}\PYG{+w}{ }\PYG{n}{i}\PYG{o}{++}\PYG{p}{)}
\PYG{+w}{        }\PYG{p}{\PYGZob{}}
\PYG{+w}{          }\PYG{n}{Var}\PYG{+w}{ }\PYG{n}{x}\PYG{+w}{ }\PYG{o}{=}\PYG{+w}{ }\PYG{n}{vars}\PYG{p}{.}\PYG{n}{GetVar}\PYG{p}{(}\PYG{n}{i}\PYG{p}{)}\PYG{p}{;}
\PYG{+w}{          }\PYG{n}{Console}\PYG{p}{.}\PYG{n}{WriteLine}\PYG{p}{(}\PYG{l+s}{\PYGZdq{}  \PYGZob{}0\PYGZcb{} = \PYGZob{}1\PYGZcb{}\PYGZdq{}}\PYG{p}{,}\PYG{+w}{ }\PYG{n}{x}\PYG{p}{.}\PYG{n}{GetName}\PYG{p}{(}\PYG{p}{)}\PYG{p}{,}\PYG{+w}{ }\PYG{n}{x}\PYG{p}{.}\PYG{n}{Get}\PYG{p}{(}\PYG{n}{Copt}\PYG{p}{.}\PYG{n}{DblInfo}\PYG{p}{.}\PYG{n}{Value}\PYG{p}{)}\PYG{p}{)}\PYG{p}{;}
\PYG{+w}{        }\PYG{p}{\PYGZcb{}}
\PYG{+w}{        }\PYG{n}{Console}\PYG{p}{.}\PYG{n}{WriteLine}\PYG{p}{(}\PYG{l+s}{\PYGZdq{}Obj = \PYGZob{}0\PYGZcb{}\PYGZdq{}}\PYG{p}{,}\PYG{+w}{ }\PYG{n}{model}\PYG{p}{.}\PYG{n}{GetDblAttr}\PYG{p}{(}\PYG{n}{Copt}\PYG{p}{.}\PYG{n}{DblAttr}\PYG{p}{.}\PYG{n}{LpObjVal}\PYG{p}{)}\PYG{p}{)}\PYG{p}{;}
\PYG{+w}{      }\PYG{p}{\PYGZcb{}}

\PYG{+w}{      }\PYG{n}{Console}\PYG{p}{.}\PYG{n}{WriteLine}\PYG{p}{(}\PYG{l+s}{\PYGZdq{}\PYGZbs{}nDone\PYGZdq{}}\PYG{p}{)}\PYG{p}{;}
\PYG{+w}{    }\PYG{p}{\PYGZcb{}}
\PYG{+w}{    }\PYG{k}{catch}\PYG{+w}{ }\PYG{p}{(}\PYG{n}{CoptException}\PYG{+w}{ }\PYG{n}{e}\PYG{p}{)}
\PYG{+w}{    }\PYG{p}{\PYGZob{}}
\PYG{+w}{      }\PYG{n}{Console}\PYG{p}{.}\PYG{n}{WriteLine}\PYG{p}{(}\PYG{l+s}{\PYGZdq{}Error Code = \PYGZob{}0\PYGZcb{}\PYGZdq{}}\PYG{p}{,}\PYG{+w}{ }\PYG{n}{e}\PYG{p}{.}\PYG{n}{GetCode}\PYG{p}{(}\PYG{p}{)}\PYG{p}{)}\PYG{p}{;}
\PYG{+w}{      }\PYG{n}{Console}\PYG{p}{.}\PYG{n}{WriteLine}\PYG{p}{(}\PYG{n}{e}\PYG{p}{.}\PYG{n}{Message}\PYG{p}{)}\PYG{p}{;}
\PYG{+w}{    }\PYG{p}{\PYGZcb{}}
\PYG{+w}{  }\PYG{p}{\PYGZcb{}}
\PYG{p}{\PYGZcb{}}
\end{sphinxVerbatim}

\sphinxAtStartPar
Let’s now walk through the example, line by line, to understand how it
achieves the desired result of optimizing the model.

\subsubsection{Creating environment and model}
\label{\detokenize{csharpinterface:creating-environment-and-model}}
\sphinxAtStartPar
Essentially, any C\# application using Cardinal Optimizer should start with a
COPT environment, where user could add one or more models. Note that each model
encapsulates a problem and corresponding data.

\sphinxAtStartPar
Furthermore, to create multiple problems, one
can load them one by one in the same model, besides the naive option of creating
multiple models in the environment.

\begin{sphinxVerbatim}[commandchars=\\\{\}]
\PYG{+w}{      }\PYG{n}{Envr}\PYG{+w}{ }\PYG{n}{env}\PYG{+w}{ }\PYG{o}{=}\PYG{+w}{ }\PYG{k}{new}\PYG{+w}{ }\PYG{n}{Envr}\PYG{p}{(}\PYG{p}{)}\PYG{p}{;}
\PYG{+w}{      }\PYG{n}{Model}\PYG{+w}{ }\PYG{n}{model}\PYG{+w}{ }\PYG{o}{=}\PYG{+w}{ }\PYG{n}{env}\PYG{p}{.}\PYG{n}{CreateModel}\PYG{p}{(}\PYG{l+s}{\PYGZdq{}lp\PYGZus{}ex1\PYGZdq{}}\PYG{p}{)}\PYG{p}{;}
\end{sphinxVerbatim}

\sphinxAtStartPar
The above call instantiates a COPT environment and a model with name “lp\_ex1”.

\subsubsection{Adding variables}
\label{\detokenize{csharpinterface:adding-variables}}
\sphinxAtStartPar
The next step in our example is to add variables to the model.
Variables are added through AddVar() or AddVars() method on the model object.
A variable is always associated with a particular model.

\begin{sphinxVerbatim}[commandchars=\\\{\}]
\PYG{+w}{      }\PYG{c+cm}{/* }
\PYG{c+cm}{       * Add variables x, y, z}
\PYG{c+cm}{       *}
\PYG{c+cm}{       * obj: 1.2 x + 1.8 y + 2.1 z}
\PYG{c+cm}{       *}
\PYG{c+cm}{       * var:}
\PYG{c+cm}{       *  0.1 \PYGZlt{}= x \PYGZlt{}= 0.6}
\PYG{c+cm}{       *  0.2 \PYGZlt{}= y \PYGZlt{}= 1.5}
\PYG{c+cm}{       *  0.3 \PYGZlt{}= z \PYGZlt{}= 2.8}
\PYG{c+cm}{       */}
\PYG{+w}{      }\PYG{n}{Var}\PYG{+w}{ }\PYG{n}{x}\PYG{+w}{ }\PYG{o}{=}\PYG{+w}{ }\PYG{n}{model}\PYG{p}{.}\PYG{n}{AddVar}\PYG{p}{(}\PYG{l+m+mf}{0.1}\PYG{p}{,}\PYG{+w}{ }\PYG{l+m+mf}{0.6}\PYG{p}{,}\PYG{+w}{ }\PYG{l+m+mf}{0.0}\PYG{p}{,}\PYG{+w}{ }\PYG{n}{Copt}\PYG{p}{.}\PYG{n}{Consts}\PYG{p}{.}\PYG{n}{CONTINUOUS}\PYG{p}{,}\PYG{+w}{ }\PYG{l+s}{\PYGZdq{}x\PYGZdq{}}\PYG{p}{)}\PYG{p}{;}
\PYG{+w}{      }\PYG{n}{Var}\PYG{+w}{ }\PYG{n}{y}\PYG{+w}{ }\PYG{o}{=}\PYG{+w}{ }\PYG{n}{model}\PYG{p}{.}\PYG{n}{AddVar}\PYG{p}{(}\PYG{l+m+mf}{0.2}\PYG{p}{,}\PYG{+w}{ }\PYG{l+m+mf}{1.5}\PYG{p}{,}\PYG{+w}{ }\PYG{l+m+mf}{0.0}\PYG{p}{,}\PYG{+w}{ }\PYG{n}{Copt}\PYG{p}{.}\PYG{n}{Consts}\PYG{p}{.}\PYG{n}{CONTINUOUS}\PYG{p}{,}\PYG{+w}{ }\PYG{l+s}{\PYGZdq{}y\PYGZdq{}}\PYG{p}{)}\PYG{p}{;}
\PYG{+w}{      }\PYG{n}{Var}\PYG{+w}{ }\PYG{n}{z}\PYG{+w}{ }\PYG{o}{=}\PYG{+w}{ }\PYG{n}{model}\PYG{p}{.}\PYG{n}{AddVar}\PYG{p}{(}\PYG{l+m+mf}{0.3}\PYG{p}{,}\PYG{+w}{ }\PYG{l+m+mf}{2.8}\PYG{p}{,}\PYG{+w}{ }\PYG{l+m+mf}{0.0}\PYG{p}{,}\PYG{+w}{ }\PYG{n}{Copt}\PYG{p}{.}\PYG{n}{Consts}\PYG{p}{.}\PYG{n}{CONTINUOUS}\PYG{p}{,}\PYG{+w}{ }\PYG{l+s}{\PYGZdq{}z\PYGZdq{}}\PYG{p}{)}\PYG{p}{;}

\PYG{+w}{      }\PYG{n}{model}\PYG{p}{.}\PYG{n}{SetObjective}\PYG{p}{(}\PYG{l+m+mf}{1.2}\PYG{+w}{ }\PYG{o}{*}\PYG{+w}{ }\PYG{n}{x}\PYG{+w}{ }\PYG{o}{+}\PYG{+w}{ }\PYG{l+m+mf}{1.8}\PYG{+w}{ }\PYG{o}{*}\PYG{+w}{ }\PYG{n}{y}\PYG{+w}{ }\PYG{o}{+}\PYG{+w}{ }\PYG{l+m+mf}{2.1}\PYG{+w}{ }\PYG{o}{*}\PYG{+w}{ }\PYG{n}{z}\PYG{p}{,}\PYG{+w}{ }\PYG{n}{Copt}\PYG{p}{.}\PYG{n}{Consts}\PYG{p}{.}\PYG{n}{MAXIMIZE}\PYG{p}{)}\PYG{p}{;}
\end{sphinxVerbatim}

\sphinxAtStartPar
The first and second arguments to the AddVar() call are the variable lower and
upper bounds, respectively. The third argument is the linear objective coefficient
(zero here \sphinxhyphen{} we’ll set the objective later). The fourth argument is the variable
type. Our variables are all continuous in this example. The final argument is the name
of the variable.

\sphinxAtStartPar
The AddVar() method has been overloaded to accept several different argument lists.
Please refer to {\hyperref[\detokenize{csharpapiref:chapcsharpapiref-model}]{\sphinxcrossref{\DUrole{std,std-ref}{Model}}}} of C\# API reference for further details.

\sphinxAtStartPar
The objective is built here using overloaded operators. The C\# API overloads the
arithmetic operators to allow you to build linear expressions by COPT variables.
The second argument indicates that the sense is maximization.

\subsubsection{Adding constraints}
\label{\detokenize{csharpinterface:adding-constraints}}
\sphinxAtStartPar
The next step in the example is to add the linear constraints. As with variables,
constraints are always associated with a specific model. They are created using
AddConstr() or AddConstrs() methods on the model object.

\begin{sphinxVerbatim}[commandchars=\\\{\}]
\PYG{+w}{      }\PYG{c+cm}{/*}
\PYG{c+cm}{       * Add two constraints using linear expression}
\PYG{c+cm}{       *}
\PYG{c+cm}{       * r0: 1.5 x + 1.2 y + 1.8 z \PYGZlt{}= 2.6}
\PYG{c+cm}{       * r1: 0.8 x + 0.6 y + 0.9 z \PYGZgt{}= 1.2}
\PYG{c+cm}{       */}
\PYG{+w}{      }\PYG{n}{model}\PYG{p}{.}\PYG{n}{AddConstr}\PYG{p}{(}\PYG{l+m+mf}{1.5}\PYG{+w}{ }\PYG{o}{*}\PYG{+w}{ }\PYG{n}{x}\PYG{+w}{ }\PYG{o}{+}\PYG{+w}{ }\PYG{l+m+mf}{1.2}\PYG{+w}{ }\PYG{o}{*}\PYG{+w}{ }\PYG{n}{y}\PYG{+w}{ }\PYG{o}{+}\PYG{+w}{ }\PYG{l+m+mf}{1.8}\PYG{+w}{ }\PYG{o}{*}\PYG{+w}{ }\PYG{n}{z}\PYG{+w}{ }\PYG{o}{\PYGZlt{}=}\PYG{+w}{ }\PYG{l+m+mf}{2.6}\PYG{p}{,}\PYG{+w}{ }\PYG{l+s}{\PYGZdq{}r0\PYGZdq{}}\PYG{p}{)}\PYG{p}{;}

\PYG{+w}{      }\PYG{n}{Expr}\PYG{+w}{ }\PYG{n}{expr}\PYG{+w}{ }\PYG{o}{=}\PYG{+w}{ }\PYG{k}{new}\PYG{+w}{ }\PYG{n}{Expr}\PYG{p}{(}\PYG{n}{x}\PYG{p}{,}\PYG{+w}{ }\PYG{l+m+mf}{0.8}\PYG{p}{)}\PYG{p}{;}
\PYG{+w}{      }\PYG{n}{expr}\PYG{p}{.}\PYG{n}{AddTerm}\PYG{p}{(}\PYG{n}{y}\PYG{p}{,}\PYG{+w}{ }\PYG{l+m+mf}{0.6}\PYG{p}{)}\PYG{p}{;}
\PYG{+w}{      }\PYG{n}{expr}\PYG{+w}{ }\PYG{o}{+=}\PYG{+w}{ }\PYG{l+m+mf}{0.9}\PYG{+w}{ }\PYG{o}{*}\PYG{+w}{ }\PYG{n}{z}\PYG{p}{;}
\PYG{+w}{      }\PYG{n}{model}\PYG{p}{.}\PYG{n}{AddConstr}\PYG{p}{(}\PYG{n}{expr}\PYG{+w}{ }\PYG{o}{\PYGZgt{}=}\PYG{+w}{ }\PYG{l+m+mf}{1.2}\PYG{p}{,}\PYG{+w}{ }\PYG{l+s}{\PYGZdq{}r1\PYGZdq{}}\PYG{p}{)}\PYG{p}{;}
\end{sphinxVerbatim}

\sphinxAtStartPar
The first constraint is to use overloaded arithmetic operators to build the
linear expression. The comparison operators are also overloaded to make it easier
to build constraints.

\sphinxAtStartPar
The second constraint is created by building a linear expression incrementally.
That is, an expression can be built by constructor of a variable and its coefficient,
by AddTerm(), and by overloaded operators.

\subsubsection{Setting parameters and attributes}
\label{\detokenize{csharpinterface:setting-parameters-and-attributes}}
\sphinxAtStartPar
The next step in the example is to set parameters and attributes of
the problem before optimization.

\begin{sphinxVerbatim}[commandchars=\\\{\}]
\PYG{+w}{      }\PYG{c+c1}{// Set parameters}
\PYG{+w}{      }\PYG{n}{model}\PYG{p}{.}\PYG{n}{SetDblParam}\PYG{p}{(}\PYG{n}{Copt}\PYG{p}{.}\PYG{n}{DblParam}\PYG{p}{.}\PYG{n}{TimeLimit}\PYG{p}{,}\PYG{+w}{ }\PYG{l+m+mi}{10}\PYG{p}{)}\PYG{p}{;}
\end{sphinxVerbatim}

\sphinxAtStartPar
The SetDblParam() call here with Copt.DblParam.TimeLimit argument sets solver to
optimize up to 10 seconds.

\subsubsection{Solving problem}
\label{\detokenize{csharpinterface:solving-problem}}
\sphinxAtStartPar
Now that the model has been built, the next step is to optimize it:

\begin{sphinxVerbatim}[commandchars=\\\{\}]
\PYG{+w}{      }\PYG{c+c1}{// Solve problem}
\PYG{+w}{      }\PYG{n}{model}\PYG{p}{.}\PYG{n}{Solve}\PYG{p}{(}\PYG{p}{)}\PYG{p}{;}
\end{sphinxVerbatim}

\sphinxAtStartPar
This routine performs the optimization and populates several internal model attributes
(including the status of the optimization, the solution, etc.).

\subsubsection{Outputting solution}
\label{\detokenize{csharpinterface:outputting-solution}}
\sphinxAtStartPar
After solving the problem, one can query the values of the attributes for
various of purposes.

\begin{sphinxVerbatim}[commandchars=\\\{\}]
\PYG{+w}{      }\PYG{c+c1}{// Output solution}
\PYG{+w}{      }\PYG{k}{if}\PYG{+w}{ }\PYG{p}{(}\PYG{n}{model}\PYG{p}{.}\PYG{n}{GetIntAttr}\PYG{p}{(}\PYG{n}{Copt}\PYG{p}{.}\PYG{n}{IntAttr}\PYG{p}{.}\PYG{n}{LpStatus}\PYG{p}{)}\PYG{+w}{ }\PYG{o}{==}\PYG{+w}{ }\PYG{n}{Copt}\PYG{p}{.}\PYG{n}{Status}\PYG{p}{.}\PYG{n}{OPTIMAL}\PYG{p}{)}
\PYG{+w}{      }\PYG{p}{\PYGZob{}}
\PYG{+w}{        }\PYG{n}{Console}\PYG{p}{.}\PYG{n}{WriteLine}\PYG{p}{(}\PYG{l+s}{\PYGZdq{}\PYGZbs{}nFound optimal solution:\PYGZdq{}}\PYG{p}{)}\PYG{p}{;}
\PYG{+w}{        }\PYG{n}{VarArray}\PYG{+w}{ }\PYG{n}{vars}\PYG{+w}{ }\PYG{o}{=}\PYG{+w}{ }\PYG{n}{model}\PYG{p}{.}\PYG{n}{GetVars}\PYG{p}{(}\PYG{p}{)}\PYG{p}{;}
\PYG{+w}{        }\PYG{k}{for}\PYG{+w}{ }\PYG{p}{(}\PYG{k+kt}{int}\PYG{+w}{ }\PYG{n}{i}\PYG{+w}{ }\PYG{o}{=}\PYG{+w}{ }\PYG{l+m+mi}{0}\PYG{p}{;}\PYG{+w}{ }\PYG{n}{i}\PYG{+w}{ }\PYG{o}{\PYGZlt{}}\PYG{+w}{ }\PYG{n}{vars}\PYG{p}{.}\PYG{n}{Size}\PYG{p}{(}\PYG{p}{)}\PYG{p}{;}\PYG{+w}{ }\PYG{n}{i}\PYG{o}{++}\PYG{p}{)}
\PYG{+w}{        }\PYG{p}{\PYGZob{}}
\PYG{+w}{          }\PYG{n}{Var}\PYG{+w}{ }\PYG{n}{x}\PYG{+w}{ }\PYG{o}{=}\PYG{+w}{ }\PYG{n}{vars}\PYG{p}{.}\PYG{n}{GetVar}\PYG{p}{(}\PYG{n}{i}\PYG{p}{)}\PYG{p}{;}
\PYG{+w}{          }\PYG{n}{Console}\PYG{p}{.}\PYG{n}{WriteLine}\PYG{p}{(}\PYG{l+s}{\PYGZdq{}  \PYGZob{}0\PYGZcb{} = \PYGZob{}1\PYGZcb{}\PYGZdq{}}\PYG{p}{,}\PYG{+w}{ }\PYG{n}{x}\PYG{p}{.}\PYG{n}{GetName}\PYG{p}{(}\PYG{p}{)}\PYG{p}{,}\PYG{+w}{ }\PYG{n}{x}\PYG{p}{.}\PYG{n}{Get}\PYG{p}{(}\PYG{n}{Copt}\PYG{p}{.}\PYG{n}{DblInfo}\PYG{p}{.}\PYG{n}{Value}\PYG{p}{)}\PYG{p}{)}\PYG{p}{;}
\PYG{+w}{        }\PYG{p}{\PYGZcb{}}
\PYG{+w}{        }\PYG{n}{Console}\PYG{p}{.}\PYG{n}{WriteLine}\PYG{p}{(}\PYG{l+s}{\PYGZdq{}Obj = \PYGZob{}0\PYGZcb{}\PYGZdq{}}\PYG{p}{,}\PYG{+w}{ }\PYG{n}{model}\PYG{p}{.}\PYG{n}{GetDblAttr}\PYG{p}{(}\PYG{n}{Copt}\PYG{p}{.}\PYG{n}{DblAttr}\PYG{p}{.}\PYG{n}{LpObjVal}\PYG{p}{)}\PYG{p}{)}\PYG{p}{;}
\PYG{+w}{      }\PYG{p}{\PYGZcb{}}

\PYG{+w}{      }\PYG{n}{Console}\PYG{p}{.}\PYG{n}{WriteLine}\PYG{p}{(}\PYG{l+s}{\PYGZdq{}\PYGZbs{}nDone\PYGZdq{}}\PYG{p}{)}\PYG{p}{;}
\end{sphinxVerbatim}

\sphinxAtStartPar
Spcifically, one can query the Copt.IntAttr.LpStatus attribute of the model to determine whether
we have found optimal LP solution; query the Copt.DblInfo.Value attribute of a variable to obtain its
solution value; query the Copt.DblAttr.LpObjVal attribute on the model to obtain the
objective value for the current solution.

\sphinxAtStartPar
The names and types of all model, variable, and constraint attributes can be found in
{\hyperref[\detokenize{csharpapiref:chapcsharpapiref-const}]{\sphinxcrossref{\DUrole{std,std-ref}{Constants}}}} of C\# API reference.

\subsubsection{Error handling}
\label{\detokenize{csharpinterface:error-handling}}
\sphinxAtStartPar
Errors in the COPT C\# interface are handled through the C\# exception mechanism.
In the example, all COPT statements are enclosed inside a try block, and any associated
errors would be caught by the catch block.

\begin{sphinxVerbatim}[commandchars=\\\{\}]
\PYG{+w}{    }\PYG{k}{catch}\PYG{+w}{ }\PYG{p}{(}\PYG{n}{CoptException}\PYG{+w}{ }\PYG{n}{e}\PYG{p}{)}
\PYG{+w}{    }\PYG{p}{\PYGZob{}}
\PYG{+w}{      }\PYG{n}{Console}\PYG{p}{.}\PYG{n}{WriteLine}\PYG{p}{(}\PYG{l+s}{\PYGZdq{}Error Code = \PYGZob{}0\PYGZcb{}\PYGZdq{}}\PYG{p}{,}\PYG{+w}{ }\PYG{n}{e}\PYG{p}{.}\PYG{n}{GetCode}\PYG{p}{(}\PYG{p}{)}\PYG{p}{)}\PYG{p}{;}
\PYG{+w}{      }\PYG{n}{Console}\PYG{p}{.}\PYG{n}{WriteLine}\PYG{p}{(}\PYG{n}{e}\PYG{p}{.}\PYG{n}{Message}\PYG{p}{)}\PYG{p}{;}
\PYG{+w}{    }\PYG{p}{\PYGZcb{}}
\end{sphinxVerbatim}

\subsection{Build and Run}
\label{\detokenize{csharpinterface:build-and-run}}
\sphinxAtStartPar
To build and run csharp example, users may refer to project under
\sphinxcode{\sphinxupquote{\$COPT\_HOME/examples/csharp}}. Specifically, We provide a csharp project file
in cross\sphinxhyphen{}platform framework of dotnet core 2.0. This example shows a single project
working on Windows, as well as Linux and Mac platforms.

\sphinxAtStartPar
First of all, download and install \sphinxhref{https://dotnet.microsoft.com/download/dotnet-core/2.0}{dotnet core 2.0} on your platform.
To get started, follow \sphinxhref{https://docs.microsoft.com/en-us/aspnet/core/getting-started/?view=aspnetcore-3.1\&tabs=windows}{instructions} in the dotnet core docs.

\subsubsection{Dotnet core 2.0 project}
\label{\detokenize{csharpinterface:dotnet-core-2-0-project}}
\sphinxAtStartPar
The dotnet core 2.0 project file \sphinxcode{\sphinxupquote{example.csproj}} example locates in folder
\sphinxcode{\sphinxupquote{\$COPT\_HOME/examples/csharp/dotnetprojects}}. Copy example file \sphinxcode{\sphinxupquote{lp\_ex1.cs}} to
this folder and change directory to there by Windows command line prompt, then run
with command \sphinxcode{\sphinxupquote{\textquotesingle{}dotnet run \sphinxhyphen{}\sphinxhyphen{}framework netcoreapp2.0\textquotesingle{}}}. For users of dotnet core 3.0,
just run with \sphinxcode{\sphinxupquote{\textquotesingle{}dotnet run \sphinxhyphen{}\sphinxhyphen{}framework netcoreapp3.0\textquotesingle{}}} instead will work too.

\sphinxAtStartPar
This csharp project has dependency on COPT dotnet 2.0 shared libary \sphinxcode{\sphinxupquote{copt\_dotnet20.dll}}, which
is refered in the project file and defines all managed classes of COPT solver.
In addition, \sphinxcode{\sphinxupquote{copt\_dotnet20.dll}} loads two shared libraries, that is, \sphinxcode{\sphinxupquote{coptcswrap.dll}}
and \sphinxcode{\sphinxupquote{copt\_cpp.dll}} on Windows,  \sphinxcode{\sphinxupquote{libcoptcswrap.so}} and \sphinxcode{\sphinxupquote{libcopt\_cpp.so}} on Linux,
\sphinxcode{\sphinxupquote{libcoptcswrap.dylib}} and \sphinxcode{\sphinxupquote{libcopt\_cpp.dylib}} on Mac respectively. Note that \sphinxcode{\sphinxupquote{coptcswrap}}
library is a bridge between managed COPT library and native library \sphinxcode{\sphinxupquote{copt\_cpp.dll}},
which declares and implements COPT constants, interfaces and methods.
So users should make sure they are installed properly on runtime search paths.

\sphinxAtStartPar
In summary, to run csharp example, users should have COPT installed properly. Specifically, it requires
three related COPT shared libaries existing on runtime search paths, and valid license files to run. Please refer to
{\hyperref[\detokenize{install:chapinstall}]{\sphinxcrossref{\DUrole{std,std-ref}{Install Guide for Cardinal Optimizer}}}} for further details.

\sphinxstepscope

\section{Java Interface}
\label{\detokenize{javainterface:java-interface}}\label{\detokenize{javainterface:chapjavainterface}}\label{\detokenize{javainterface::doc}}
\sphinxAtStartPar
This chapter walks through a simple Java example to illustrate the use of the
COPT Java interface. In short words, the example creates an environment, builds
a model, add variables and constraints, optimizes it, and then outputs the
optimal objective value.

\sphinxAtStartPar
The example solves the following linear problem:
\begin{equation}\label{equation:javainterface:coptjavaEq_capilp1}
\begin{split}\text{Maximize: } & \\
                  & 1.2 x + 1.8 y + 2.1 z \\
\text{Subject to: } & \\
                    & 1.5 x + 1.2 y + 1.8 z \leq 2.6 \\
                    & 0.8 x + 0.6 y + 0.9 z \geq 1.2 \\
\text{Bounds: } & \\
               & 0.1 \leq x \leq 0.6 \\
               & 0.2 \leq y \leq 1.5 \\
               & 0.3 \leq z \leq 2.8 \\\end{split}
\end{equation}
\sphinxAtStartPar
Note that this is the same problem that was modelled and optimized in
chapter of {\hyperref[\detokenize{cinterface:chapcinterface}]{\sphinxcrossref{\DUrole{std,std-ref}{C Interface}}}}.

\subsection{Example details}
\label{\detokenize{javainterface:example-details}}
\sphinxAtStartPar
Below is the source code solving the above problem using COPT Java interface.
\sphinxSetupCaptionForVerbatim{\sphinxcode{\sphinxupquote{Lp\_ex1.java}}}
\def\sphinxLiteralBlockLabel{\label{\detokenize{javainterface:id1}}}
\begin{sphinxVerbatim}[commandchars=\\\{\},numbers=left,firstnumber=1,stepnumber=1]
\PYG{c+cm}{/*}
\PYG{c+cm}{ * This file is part of the Cardinal Optimizer, all rights reserved.}
\PYG{c+cm}{ */}
\PYG{k+kn}{import}\PYG{+w}{ }\PYG{n+nn}{copt.*}\PYG{p}{;}

\PYG{c+cm}{/*}
\PYG{c+cm}{ * This Java example solves the following LP model:}
\PYG{c+cm}{ *}
\PYG{c+cm}{ * Maximize:}
\PYG{c+cm}{ *  1.2 x + 1.8 y + 2.1 z}
\PYG{c+cm}{ *}
\PYG{c+cm}{ * Subject to:}
\PYG{c+cm}{ *  1.5 x + 1.2 y + 1.8 z \PYGZlt{}= 2.6}
\PYG{c+cm}{ *  0.8 x + 0.6 y + 0.9 z \PYGZgt{}= 1.2}
\PYG{c+cm}{ *}
\PYG{c+cm}{ * where:}
\PYG{c+cm}{ *  0.1 \PYGZlt{}= x \PYGZlt{}= 0.6}
\PYG{c+cm}{ *  0.2 \PYGZlt{}= y \PYGZlt{}= 1.5}
\PYG{c+cm}{ *  0.3 \PYGZlt{}= z \PYGZlt{}= 2.8}
\PYG{c+cm}{ */}
\PYG{k+kd}{public}\PYG{+w}{ }\PYG{k+kd}{class} \PYG{n+nc}{Lp\PYGZus{}ex1}\PYG{+w}{ }\PYG{p}{\PYGZob{}}
\PYG{+w}{  }\PYG{k+kd}{public}\PYG{+w}{ }\PYG{k+kd}{static}\PYG{+w}{ }\PYG{k+kt}{void}\PYG{+w}{ }\PYG{n+nf}{main}\PYG{p}{(}\PYG{k+kd}{final}\PYG{+w}{ }\PYG{n}{String}\PYG{+w}{ }\PYG{n}{argv}\PYG{o}{[}\PYG{o}{]}\PYG{p}{)}\PYG{+w}{ }\PYG{p}{\PYGZob{}}
\PYG{+w}{    }\PYG{k}{try}\PYG{+w}{ }\PYG{p}{\PYGZob{}}
\PYG{+w}{      }\PYG{n}{Envr}\PYG{+w}{ }\PYG{n}{env}\PYG{+w}{ }\PYG{o}{=}\PYG{+w}{ }\PYG{k}{new}\PYG{+w}{ }\PYG{n}{Envr}\PYG{p}{(}\PYG{p}{)}\PYG{p}{;}
\PYG{+w}{      }\PYG{n}{Model}\PYG{+w}{ }\PYG{n}{model}\PYG{+w}{ }\PYG{o}{=}\PYG{+w}{ }\PYG{n}{env}\PYG{p}{.}\PYG{n+na}{createModel}\PYG{p}{(}\PYG{l+s}{\PYGZdq{}}\PYG{l+s}{lp\PYGZus{}ex1}\PYG{l+s}{\PYGZdq{}}\PYG{p}{)}\PYG{p}{;}

\PYG{+w}{      }\PYG{c+cm}{/* }
\PYG{c+cm}{       * Add variables x, y, z}
\PYG{c+cm}{       *}
\PYG{c+cm}{       * obj: 1.2 x + 1.8 y + 2.1 z}
\PYG{c+cm}{       *}
\PYG{c+cm}{       * var:}
\PYG{c+cm}{       *  0.1 \PYGZlt{}= x \PYGZlt{}= 0.6}
\PYG{c+cm}{       *  0.2 \PYGZlt{}= y \PYGZlt{}= 1.5}
\PYG{c+cm}{       *  0.3 \PYGZlt{}= z \PYGZlt{}= 2.8}
\PYG{c+cm}{       */}
\PYG{+w}{      }\PYG{n}{Var}\PYG{+w}{ }\PYG{n}{x}\PYG{+w}{ }\PYG{o}{=}\PYG{+w}{ }\PYG{n}{model}\PYG{p}{.}\PYG{n+na}{addVar}\PYG{p}{(}\PYG{l+m+mf}{0.1}\PYG{p}{,}\PYG{+w}{ }\PYG{l+m+mf}{0.6}\PYG{p}{,}\PYG{+w}{ }\PYG{l+m+mf}{1.2}\PYG{p}{,}\PYG{+w}{ }\PYG{n}{copt}\PYG{p}{.}\PYG{n+na}{Consts}\PYG{p}{.}\PYG{n+na}{CONTINUOUS}\PYG{p}{,}\PYG{+w}{ }\PYG{l+s}{\PYGZdq{}}\PYG{l+s}{x}\PYG{l+s}{\PYGZdq{}}\PYG{p}{)}\PYG{p}{;}
\PYG{+w}{      }\PYG{n}{Var}\PYG{+w}{ }\PYG{n}{y}\PYG{+w}{ }\PYG{o}{=}\PYG{+w}{ }\PYG{n}{model}\PYG{p}{.}\PYG{n+na}{addVar}\PYG{p}{(}\PYG{l+m+mf}{0.2}\PYG{p}{,}\PYG{+w}{ }\PYG{l+m+mf}{1.5}\PYG{p}{,}\PYG{+w}{ }\PYG{l+m+mf}{1.8}\PYG{p}{,}\PYG{+w}{ }\PYG{n}{copt}\PYG{p}{.}\PYG{n+na}{Consts}\PYG{p}{.}\PYG{n+na}{CONTINUOUS}\PYG{p}{,}\PYG{+w}{ }\PYG{l+s}{\PYGZdq{}}\PYG{l+s}{y}\PYG{l+s}{\PYGZdq{}}\PYG{p}{)}\PYG{p}{;}
\PYG{+w}{      }\PYG{n}{Var}\PYG{+w}{ }\PYG{n}{z}\PYG{+w}{ }\PYG{o}{=}\PYG{+w}{ }\PYG{n}{model}\PYG{p}{.}\PYG{n+na}{addVar}\PYG{p}{(}\PYG{l+m+mf}{0.3}\PYG{p}{,}\PYG{+w}{ }\PYG{l+m+mf}{2.8}\PYG{p}{,}\PYG{+w}{ }\PYG{l+m+mf}{2.1}\PYG{p}{,}\PYG{+w}{ }\PYG{n}{copt}\PYG{p}{.}\PYG{n+na}{Consts}\PYG{p}{.}\PYG{n+na}{CONTINUOUS}\PYG{p}{,}\PYG{+w}{ }\PYG{l+s}{\PYGZdq{}}\PYG{l+s}{z}\PYG{l+s}{\PYGZdq{}}\PYG{p}{)}\PYG{p}{;}

\PYG{+w}{      }\PYG{c+cm}{/*}
\PYG{c+cm}{       * Add two constraints using linear expression}
\PYG{c+cm}{       *}
\PYG{c+cm}{       * r0: 1.5 x + 1.2 y + 1.8 z \PYGZlt{}= 2.6}
\PYG{c+cm}{       * r1: 0.8 x + 0.6 y + 0.9 z \PYGZgt{}= 1.2}
\PYG{c+cm}{       */}
\PYG{+w}{      }\PYG{n}{Expr}\PYG{+w}{ }\PYG{n}{e0}\PYG{+w}{ }\PYG{o}{=}\PYG{+w}{ }\PYG{k}{new}\PYG{+w}{ }\PYG{n}{Expr}\PYG{p}{(}\PYG{n}{x}\PYG{p}{,}\PYG{+w}{ }\PYG{l+m+mf}{1.5}\PYG{p}{)}\PYG{p}{;}
\PYG{+w}{      }\PYG{n}{e0}\PYG{p}{.}\PYG{n+na}{addTerm}\PYG{p}{(}\PYG{n}{y}\PYG{p}{,}\PYG{+w}{ }\PYG{l+m+mf}{1.2}\PYG{p}{)}\PYG{p}{;}
\PYG{+w}{      }\PYG{n}{e0}\PYG{p}{.}\PYG{n+na}{addTerm}\PYG{p}{(}\PYG{n}{z}\PYG{p}{,}\PYG{+w}{ }\PYG{l+m+mf}{1.8}\PYG{p}{)}\PYG{p}{;}
\PYG{+w}{      }\PYG{n}{model}\PYG{p}{.}\PYG{n+na}{addConstr}\PYG{p}{(}\PYG{n}{e0}\PYG{p}{,}\PYG{+w}{ }\PYG{n}{copt}\PYG{p}{.}\PYG{n+na}{Consts}\PYG{p}{.}\PYG{n+na}{LESS\PYGZus{}EQUAL}\PYG{p}{,}\PYG{+w}{ }\PYG{l+m+mf}{2.6}\PYG{p}{,}\PYG{+w}{ }\PYG{l+s}{\PYGZdq{}}\PYG{l+s}{r0}\PYG{l+s}{\PYGZdq{}}\PYG{p}{)}\PYG{p}{;}

\PYG{+w}{      }\PYG{n}{Expr}\PYG{+w}{ }\PYG{n}{e1}\PYG{+w}{ }\PYG{o}{=}\PYG{+w}{ }\PYG{k}{new}\PYG{+w}{ }\PYG{n}{Expr}\PYG{p}{(}\PYG{n}{x}\PYG{p}{,}\PYG{+w}{ }\PYG{l+m+mf}{0.8}\PYG{p}{)}\PYG{p}{;}
\PYG{+w}{      }\PYG{n}{e1}\PYG{p}{.}\PYG{n+na}{addTerm}\PYG{p}{(}\PYG{n}{y}\PYG{p}{,}\PYG{+w}{ }\PYG{l+m+mf}{0.6}\PYG{p}{)}\PYG{p}{;}
\PYG{+w}{      }\PYG{n}{e1}\PYG{p}{.}\PYG{n+na}{addTerm}\PYG{p}{(}\PYG{n}{z}\PYG{p}{,}\PYG{+w}{ }\PYG{l+m+mf}{0.9}\PYG{p}{)}\PYG{p}{;}
\PYG{+w}{      }\PYG{n}{model}\PYG{p}{.}\PYG{n+na}{addConstr}\PYG{p}{(}\PYG{n}{e1}\PYG{p}{,}\PYG{+w}{ }\PYG{n}{copt}\PYG{p}{.}\PYG{n+na}{Consts}\PYG{p}{.}\PYG{n+na}{GREATER\PYGZus{}EQUAL}\PYG{p}{,}\PYG{+w}{ }\PYG{l+m+mf}{1.2}\PYG{p}{,}\PYG{+w}{ }\PYG{l+s}{\PYGZdq{}}\PYG{l+s}{r1}\PYG{l+s}{\PYGZdq{}}\PYG{p}{)}\PYG{p}{;}

\PYG{+w}{      }\PYG{c+c1}{// Set parameters and attributes}
\PYG{+w}{      }\PYG{n}{model}\PYG{p}{.}\PYG{n+na}{setDblParam}\PYG{p}{(}\PYG{n}{copt}\PYG{p}{.}\PYG{n+na}{DblParam}\PYG{p}{.}\PYG{n+na}{TimeLimit}\PYG{p}{,}\PYG{+w}{ }\PYG{l+m+mi}{10}\PYG{p}{)}\PYG{p}{;}
\PYG{+w}{      }\PYG{n}{model}\PYG{p}{.}\PYG{n+na}{setObjSense}\PYG{p}{(}\PYG{n}{copt}\PYG{p}{.}\PYG{n+na}{Consts}\PYG{p}{.}\PYG{n+na}{MAXIMIZE}\PYG{p}{)}\PYG{p}{;}

\PYG{+w}{      }\PYG{c+c1}{// Solve problem}
\PYG{+w}{      }\PYG{n}{model}\PYG{p}{.}\PYG{n+na}{solve}\PYG{p}{(}\PYG{p}{)}\PYG{p}{;}

\PYG{+w}{      }\PYG{c+c1}{// Output solution}
\PYG{+w}{      }\PYG{k}{if}\PYG{+w}{ }\PYG{p}{(}\PYG{n}{model}\PYG{p}{.}\PYG{n+na}{getIntAttr}\PYG{p}{(}\PYG{n}{copt}\PYG{p}{.}\PYG{n+na}{IntAttr}\PYG{p}{.}\PYG{n+na}{HasLpSol}\PYG{p}{)}\PYG{+w}{ }\PYG{o}{!}\PYG{o}{=}\PYG{+w}{ }\PYG{l+m+mi}{0}\PYG{p}{)}\PYG{+w}{ }\PYG{p}{\PYGZob{}}
\PYG{+w}{        }\PYG{n}{System}\PYG{p}{.}\PYG{n+na}{out}\PYG{p}{.}\PYG{n+na}{println}\PYG{p}{(}\PYG{l+s}{\PYGZdq{}}\PYG{l+s}{\PYGZbs{}}\PYG{l+s}{nFound optimal solution:}\PYG{l+s}{\PYGZdq{}}\PYG{p}{)}\PYG{p}{;}
\PYG{+w}{        }\PYG{n}{VarArray}\PYG{+w}{ }\PYG{n}{vars}\PYG{+w}{ }\PYG{o}{=}\PYG{+w}{ }\PYG{n}{model}\PYG{p}{.}\PYG{n+na}{getVars}\PYG{p}{(}\PYG{p}{)}\PYG{p}{;}
\PYG{+w}{        }\PYG{k}{for}\PYG{+w}{ }\PYG{p}{(}\PYG{k+kt}{int}\PYG{+w}{ }\PYG{n}{i}\PYG{+w}{ }\PYG{o}{=}\PYG{+w}{ }\PYG{l+m+mi}{0}\PYG{p}{;}\PYG{+w}{ }\PYG{n}{i}\PYG{+w}{ }\PYG{o}{\PYGZlt{}}\PYG{+w}{ }\PYG{n}{vars}\PYG{p}{.}\PYG{n+na}{size}\PYG{p}{(}\PYG{p}{)}\PYG{p}{;}\PYG{+w}{ }\PYG{n}{i}\PYG{o}{+}\PYG{o}{+}\PYG{p}{)}\PYG{+w}{ }\PYG{p}{\PYGZob{}}
\PYG{+w}{          }\PYG{n}{Var}\PYG{+w}{ }\PYG{n}{x}\PYG{+w}{ }\PYG{o}{=}\PYG{+w}{ }\PYG{n}{vars}\PYG{p}{.}\PYG{n+na}{getVar}\PYG{p}{(}\PYG{n}{i}\PYG{p}{)}\PYG{p}{;}
\PYG{+w}{          }\PYG{n}{System}\PYG{p}{.}\PYG{n+na}{out}\PYG{p}{.}\PYG{n+na}{println}\PYG{p}{(}\PYG{l+s}{\PYGZdq{}}\PYG{l+s}{  }\PYG{l+s}{\PYGZdq{}}\PYG{+w}{ }\PYG{o}{+}\PYG{+w}{ }\PYG{n}{x}\PYG{p}{.}\PYG{n+na}{getName}\PYG{p}{(}\PYG{p}{)}\PYG{+w}{ }\PYG{o}{+}\PYG{+w}{ }\PYG{l+s}{\PYGZdq{}}\PYG{l+s}{ = }\PYG{l+s}{\PYGZdq{}}\PYG{+w}{ }\PYG{o}{+}\PYG{+w}{ }\PYG{n}{x}\PYG{p}{.}\PYG{n+na}{get}\PYG{p}{(}\PYG{n}{copt}\PYG{p}{.}\PYG{n+na}{DblInfo}\PYG{p}{.}\PYG{n+na}{Value}\PYG{p}{)}\PYG{p}{)}\PYG{p}{;}
\PYG{+w}{        }\PYG{p}{\PYGZcb{}}
\PYG{+w}{        }\PYG{n}{System}\PYG{p}{.}\PYG{n+na}{out}\PYG{p}{.}\PYG{n+na}{println}\PYG{p}{(}\PYG{l+s}{\PYGZdq{}}\PYG{l+s}{Obj = }\PYG{l+s}{\PYGZdq{}}\PYG{+w}{ }\PYG{o}{+}\PYG{+w}{ }\PYG{n}{model}\PYG{p}{.}\PYG{n+na}{getDblAttr}\PYG{p}{(}\PYG{n}{copt}\PYG{p}{.}\PYG{n+na}{DblAttr}\PYG{p}{.}\PYG{n+na}{LpObjVal}\PYG{p}{)}\PYG{p}{)}\PYG{p}{;}
\PYG{+w}{      }\PYG{p}{\PYGZcb{}}

\PYG{+w}{      }\PYG{n}{System}\PYG{p}{.}\PYG{n+na}{out}\PYG{p}{.}\PYG{n+na}{println}\PYG{p}{(}\PYG{l+s}{\PYGZdq{}}\PYG{l+s}{\PYGZbs{}}\PYG{l+s}{nDone}\PYG{l+s}{\PYGZdq{}}\PYG{p}{)}\PYG{p}{;}
\PYG{+w}{    }\PYG{p}{\PYGZcb{}}\PYG{+w}{ }\PYG{k}{catch}\PYG{+w}{ }\PYG{p}{(}\PYG{n}{CoptException}\PYG{+w}{ }\PYG{n}{e}\PYG{p}{)}\PYG{+w}{ }\PYG{p}{\PYGZob{}}
\PYG{+w}{      }\PYG{n}{System}\PYG{p}{.}\PYG{n+na}{out}\PYG{p}{.}\PYG{n+na}{println}\PYG{p}{(}\PYG{l+s}{\PYGZdq{}}\PYG{l+s}{Error Code = }\PYG{l+s}{\PYGZdq{}}\PYG{+w}{ }\PYG{o}{+}\PYG{+w}{ }\PYG{n}{e}\PYG{p}{.}\PYG{n+na}{getCode}\PYG{p}{(}\PYG{p}{)}\PYG{p}{)}\PYG{p}{;}
\PYG{+w}{      }\PYG{n}{System}\PYG{p}{.}\PYG{n+na}{out}\PYG{p}{.}\PYG{n+na}{println}\PYG{p}{(}\PYG{n}{e}\PYG{p}{.}\PYG{n+na}{getMessage}\PYG{p}{(}\PYG{p}{)}\PYG{p}{)}\PYG{p}{;}
\PYG{+w}{    }\PYG{p}{\PYGZcb{}}
\PYG{+w}{  }\PYG{p}{\PYGZcb{}}
\PYG{p}{\PYGZcb{}}
\end{sphinxVerbatim}

\sphinxAtStartPar
Let’s now walk through the example, line by line, to understand how it
achieves the desired result of optimizing the model.

\subsubsection{Import COPT class}
\label{\detokenize{javainterface:import-copt-class}}
\sphinxAtStartPar
To use the Java interface of COPT, users need to import the Java interface class
of COPT first.

\begin{sphinxVerbatim}[commandchars=\\\{\}]
\PYG{k+kn}{import}\PYG{+w}{ }\PYG{n+nn}{copt.*}\PYG{p}{;}
\end{sphinxVerbatim}

\subsubsection{Creating environment and model}
\label{\detokenize{javainterface:creating-environment-and-model}}
\sphinxAtStartPar
Essentially, any Java application using Cardinal Optimizer should start with a
COPT environment, where user could add one or more models. Note that each model
encapsulates a problem and corresponding data.

\sphinxAtStartPar
Furthermore, to create multiple problems, one
can load them one by one in the same model, besides the naive option of creating
multiple models in the environment.

\begin{sphinxVerbatim}[commandchars=\\\{\}]
\PYG{+w}{      }\PYG{n}{Envr}\PYG{+w}{ }\PYG{n}{env}\PYG{+w}{ }\PYG{o}{=}\PYG{+w}{ }\PYG{k}{new}\PYG{+w}{ }\PYG{n}{Envr}\PYG{p}{(}\PYG{p}{)}\PYG{p}{;}
\PYG{+w}{      }\PYG{n}{Model}\PYG{+w}{ }\PYG{n}{model}\PYG{+w}{ }\PYG{o}{=}\PYG{+w}{ }\PYG{n}{env}\PYG{p}{.}\PYG{n+na}{createModel}\PYG{p}{(}\PYG{l+s}{\PYGZdq{}}\PYG{l+s}{lp\PYGZus{}ex1}\PYG{l+s}{\PYGZdq{}}\PYG{p}{)}\PYG{p}{;}
\end{sphinxVerbatim}

\sphinxAtStartPar
The above call instantiates a COPT environment and a model with name “lp\_ex1”.

\subsubsection{Adding variables}
\label{\detokenize{javainterface:adding-variables}}
\sphinxAtStartPar
The next step in our example is to add variables to the model.
Variables are added through addVar() or addVars() method on the model object.
A variable is always associated with a particular model.

\begin{sphinxVerbatim}[commandchars=\\\{\}]
\PYG{+w}{      }\PYG{c+cm}{/* }
\PYG{c+cm}{       * Add variables x, y, z}
\PYG{c+cm}{       *}
\PYG{c+cm}{       * obj: 1.2 x + 1.8 y + 2.1 z}
\PYG{c+cm}{       *}
\PYG{c+cm}{       * var:}
\PYG{c+cm}{       *  0.1 \PYGZlt{}= x \PYGZlt{}= 0.6}
\PYG{c+cm}{       *  0.2 \PYGZlt{}= y \PYGZlt{}= 1.5}
\PYG{c+cm}{       *  0.3 \PYGZlt{}= z \PYGZlt{}= 2.8}
\PYG{c+cm}{       */}
\PYG{+w}{      }\PYG{n}{Var}\PYG{+w}{ }\PYG{n}{x}\PYG{+w}{ }\PYG{o}{=}\PYG{+w}{ }\PYG{n}{model}\PYG{p}{.}\PYG{n+na}{addVar}\PYG{p}{(}\PYG{l+m+mf}{0.1}\PYG{p}{,}\PYG{+w}{ }\PYG{l+m+mf}{0.6}\PYG{p}{,}\PYG{+w}{ }\PYG{l+m+mf}{1.2}\PYG{p}{,}\PYG{+w}{ }\PYG{n}{copt}\PYG{p}{.}\PYG{n+na}{Consts}\PYG{p}{.}\PYG{n+na}{CONTINUOUS}\PYG{p}{,}\PYG{+w}{ }\PYG{l+s}{\PYGZdq{}}\PYG{l+s}{x}\PYG{l+s}{\PYGZdq{}}\PYG{p}{)}\PYG{p}{;}
\PYG{+w}{      }\PYG{n}{Var}\PYG{+w}{ }\PYG{n}{y}\PYG{+w}{ }\PYG{o}{=}\PYG{+w}{ }\PYG{n}{model}\PYG{p}{.}\PYG{n+na}{addVar}\PYG{p}{(}\PYG{l+m+mf}{0.2}\PYG{p}{,}\PYG{+w}{ }\PYG{l+m+mf}{1.5}\PYG{p}{,}\PYG{+w}{ }\PYG{l+m+mf}{1.8}\PYG{p}{,}\PYG{+w}{ }\PYG{n}{copt}\PYG{p}{.}\PYG{n+na}{Consts}\PYG{p}{.}\PYG{n+na}{CONTINUOUS}\PYG{p}{,}\PYG{+w}{ }\PYG{l+s}{\PYGZdq{}}\PYG{l+s}{y}\PYG{l+s}{\PYGZdq{}}\PYG{p}{)}\PYG{p}{;}
\PYG{+w}{      }\PYG{n}{Var}\PYG{+w}{ }\PYG{n}{z}\PYG{+w}{ }\PYG{o}{=}\PYG{+w}{ }\PYG{n}{model}\PYG{p}{.}\PYG{n+na}{addVar}\PYG{p}{(}\PYG{l+m+mf}{0.3}\PYG{p}{,}\PYG{+w}{ }\PYG{l+m+mf}{2.8}\PYG{p}{,}\PYG{+w}{ }\PYG{l+m+mf}{2.1}\PYG{p}{,}\PYG{+w}{ }\PYG{n}{copt}\PYG{p}{.}\PYG{n+na}{Consts}\PYG{p}{.}\PYG{n+na}{CONTINUOUS}\PYG{p}{,}\PYG{+w}{ }\PYG{l+s}{\PYGZdq{}}\PYG{l+s}{z}\PYG{l+s}{\PYGZdq{}}\PYG{p}{)}\PYG{p}{;}
\end{sphinxVerbatim}

\sphinxAtStartPar
The first and second arguments to the addVar() call are the variable lower and
upper bounds, respectively. The third argument is the linear objective coefficient.
The fourth argument is the variable type. Our variables are all continuous in this
example. The final argument is the name of the variable.

\sphinxAtStartPar
The addVar() method has been overloaded to accept several different argument lists.
Please refer to {\hyperref[\detokenize{javaapiref:chapjavaapiref-class}]{\sphinxcrossref{\DUrole{std,std-ref}{Java Modeling Classes}}}} of Java API reference for further details.

\subsubsection{Adding constraints}
\label{\detokenize{javainterface:adding-constraints}}
\sphinxAtStartPar
The next step in the example is to add the linear constraints. As with variables,
constraints are always associated with a specific model. They are created using
addConstr() or addConstrs() methods on the model object.

\begin{sphinxVerbatim}[commandchars=\\\{\}]
\PYG{+w}{      }\PYG{c+cm}{/*}
\PYG{c+cm}{       * Add two constraints using linear expression}
\PYG{c+cm}{       *}
\PYG{c+cm}{       * r0: 1.5 x + 1.2 y + 1.8 z \PYGZlt{}= 2.6}
\PYG{c+cm}{       * r1: 0.8 x + 0.6 y + 0.9 z \PYGZgt{}= 1.2}
\PYG{c+cm}{       */}
\PYG{+w}{      }\PYG{n}{Expr}\PYG{+w}{ }\PYG{n}{e0}\PYG{+w}{ }\PYG{o}{=}\PYG{+w}{ }\PYG{k}{new}\PYG{+w}{ }\PYG{n}{Expr}\PYG{p}{(}\PYG{n}{x}\PYG{p}{,}\PYG{+w}{ }\PYG{l+m+mf}{1.5}\PYG{p}{)}\PYG{p}{;}
\PYG{+w}{      }\PYG{n}{e0}\PYG{p}{.}\PYG{n+na}{addTerm}\PYG{p}{(}\PYG{n}{y}\PYG{p}{,}\PYG{+w}{ }\PYG{l+m+mf}{1.2}\PYG{p}{)}\PYG{p}{;}
\PYG{+w}{      }\PYG{n}{e0}\PYG{p}{.}\PYG{n+na}{addTerm}\PYG{p}{(}\PYG{n}{z}\PYG{p}{,}\PYG{+w}{ }\PYG{l+m+mf}{1.8}\PYG{p}{)}\PYG{p}{;}
\PYG{+w}{      }\PYG{n}{model}\PYG{p}{.}\PYG{n+na}{addConstr}\PYG{p}{(}\PYG{n}{e0}\PYG{p}{,}\PYG{+w}{ }\PYG{n}{copt}\PYG{p}{.}\PYG{n+na}{Consts}\PYG{p}{.}\PYG{n+na}{LESS\PYGZus{}EQUAL}\PYG{p}{,}\PYG{+w}{ }\PYG{l+m+mf}{2.6}\PYG{p}{,}\PYG{+w}{ }\PYG{l+s}{\PYGZdq{}}\PYG{l+s}{r0}\PYG{l+s}{\PYGZdq{}}\PYG{p}{)}\PYG{p}{;}

\PYG{+w}{      }\PYG{n}{Expr}\PYG{+w}{ }\PYG{n}{e1}\PYG{+w}{ }\PYG{o}{=}\PYG{+w}{ }\PYG{k}{new}\PYG{+w}{ }\PYG{n}{Expr}\PYG{p}{(}\PYG{n}{x}\PYG{p}{,}\PYG{+w}{ }\PYG{l+m+mf}{0.8}\PYG{p}{)}\PYG{p}{;}
\PYG{+w}{      }\PYG{n}{e1}\PYG{p}{.}\PYG{n+na}{addTerm}\PYG{p}{(}\PYG{n}{y}\PYG{p}{,}\PYG{+w}{ }\PYG{l+m+mf}{0.6}\PYG{p}{)}\PYG{p}{;}
\PYG{+w}{      }\PYG{n}{e1}\PYG{p}{.}\PYG{n+na}{addTerm}\PYG{p}{(}\PYG{n}{z}\PYG{p}{,}\PYG{+w}{ }\PYG{l+m+mf}{0.9}\PYG{p}{)}\PYG{p}{;}
\PYG{+w}{      }\PYG{n}{model}\PYG{p}{.}\PYG{n+na}{addConstr}\PYG{p}{(}\PYG{n}{e1}\PYG{p}{,}\PYG{+w}{ }\PYG{n}{copt}\PYG{p}{.}\PYG{n+na}{Consts}\PYG{p}{.}\PYG{n+na}{GREATER\PYGZus{}EQUAL}\PYG{p}{,}\PYG{+w}{ }\PYG{l+m+mf}{1.2}\PYG{p}{,}\PYG{+w}{ }\PYG{l+s}{\PYGZdq{}}\PYG{l+s}{r1}\PYG{l+s}{\PYGZdq{}}\PYG{p}{)}\PYG{p}{;}
\end{sphinxVerbatim}

\sphinxAtStartPar
Two constraints here are created by building linear expressions incrementally.
That is, an expression can be built by constructor of a variable and its coefficient,
and then by addTerm() method.

\subsubsection{Setting parameters and attributes}
\label{\detokenize{javainterface:setting-parameters-and-attributes}}
\sphinxAtStartPar
The next step in the example is to set parameters and attributes of
the problem before optimization.

\begin{sphinxVerbatim}[commandchars=\\\{\}]
\PYG{+w}{      }\PYG{c+c1}{// Set parameters and attributes}
\PYG{+w}{      }\PYG{n}{model}\PYG{p}{.}\PYG{n+na}{setDblParam}\PYG{p}{(}\PYG{n}{copt}\PYG{p}{.}\PYG{n+na}{DblParam}\PYG{p}{.}\PYG{n+na}{TimeLimit}\PYG{p}{,}\PYG{+w}{ }\PYG{l+m+mi}{10}\PYG{p}{)}\PYG{p}{;}
\PYG{+w}{      }\PYG{n}{model}\PYG{p}{.}\PYG{n+na}{setObjSense}\PYG{p}{(}\PYG{n}{copt}\PYG{p}{.}\PYG{n+na}{Consts}\PYG{p}{.}\PYG{n+na}{MAXIMIZE}\PYG{p}{)}\PYG{p}{;}
\end{sphinxVerbatim}

\sphinxAtStartPar
The setDblParam() call here with copt.DblParam.TimeLimit argument sets solver to
optimize up to 10 seconds. The setObjSense() call with copt.Consts.MAXIMIZE argument
sets objective sense as maximization.

\subsubsection{Solving problem}
\label{\detokenize{javainterface:solving-problem}}
\sphinxAtStartPar
Now that the model has been built, the next step is to optimize it:

\begin{sphinxVerbatim}[commandchars=\\\{\}]
\PYG{+w}{      }\PYG{c+c1}{// Solve problem}
\PYG{+w}{      }\PYG{n}{model}\PYG{p}{.}\PYG{n+na}{solve}\PYG{p}{(}\PYG{p}{)}\PYG{p}{;}
\end{sphinxVerbatim}

\sphinxAtStartPar
This routine performs the optimization and populates several internal model attributes
(including the status of the optimization, the solution, etc.).

\subsubsection{Outputting solution}
\label{\detokenize{javainterface:outputting-solution}}
\sphinxAtStartPar
After solving the problem, one can query the values of the attributes for
various of purposes.

\begin{sphinxVerbatim}[commandchars=\\\{\}]
\PYG{+w}{      }\PYG{c+c1}{// Output solution}
\PYG{+w}{      }\PYG{k}{if}\PYG{+w}{ }\PYG{p}{(}\PYG{n}{model}\PYG{p}{.}\PYG{n+na}{getIntAttr}\PYG{p}{(}\PYG{n}{copt}\PYG{p}{.}\PYG{n+na}{IntAttr}\PYG{p}{.}\PYG{n+na}{HasLpSol}\PYG{p}{)}\PYG{+w}{ }\PYG{o}{!}\PYG{o}{=}\PYG{+w}{ }\PYG{l+m+mi}{0}\PYG{p}{)}\PYG{+w}{ }\PYG{p}{\PYGZob{}}
\PYG{+w}{        }\PYG{n}{System}\PYG{p}{.}\PYG{n+na}{out}\PYG{p}{.}\PYG{n+na}{println}\PYG{p}{(}\PYG{l+s}{\PYGZdq{}}\PYG{l+s}{\PYGZbs{}}\PYG{l+s}{nFound optimal solution:}\PYG{l+s}{\PYGZdq{}}\PYG{p}{)}\PYG{p}{;}
\PYG{+w}{        }\PYG{n}{VarArray}\PYG{+w}{ }\PYG{n}{vars}\PYG{+w}{ }\PYG{o}{=}\PYG{+w}{ }\PYG{n}{model}\PYG{p}{.}\PYG{n+na}{getVars}\PYG{p}{(}\PYG{p}{)}\PYG{p}{;}
\PYG{+w}{        }\PYG{k}{for}\PYG{+w}{ }\PYG{p}{(}\PYG{k+kt}{int}\PYG{+w}{ }\PYG{n}{i}\PYG{+w}{ }\PYG{o}{=}\PYG{+w}{ }\PYG{l+m+mi}{0}\PYG{p}{;}\PYG{+w}{ }\PYG{n}{i}\PYG{+w}{ }\PYG{o}{\PYGZlt{}}\PYG{+w}{ }\PYG{n}{vars}\PYG{p}{.}\PYG{n+na}{size}\PYG{p}{(}\PYG{p}{)}\PYG{p}{;}\PYG{+w}{ }\PYG{n}{i}\PYG{o}{+}\PYG{o}{+}\PYG{p}{)}\PYG{+w}{ }\PYG{p}{\PYGZob{}}
\PYG{+w}{          }\PYG{n}{Var}\PYG{+w}{ }\PYG{n}{x}\PYG{+w}{ }\PYG{o}{=}\PYG{+w}{ }\PYG{n}{vars}\PYG{p}{.}\PYG{n+na}{getVar}\PYG{p}{(}\PYG{n}{i}\PYG{p}{)}\PYG{p}{;}
\PYG{+w}{          }\PYG{n}{System}\PYG{p}{.}\PYG{n+na}{out}\PYG{p}{.}\PYG{n+na}{println}\PYG{p}{(}\PYG{l+s}{\PYGZdq{}}\PYG{l+s}{  }\PYG{l+s}{\PYGZdq{}}\PYG{+w}{ }\PYG{o}{+}\PYG{+w}{ }\PYG{n}{x}\PYG{p}{.}\PYG{n+na}{getName}\PYG{p}{(}\PYG{p}{)}\PYG{+w}{ }\PYG{o}{+}\PYG{+w}{ }\PYG{l+s}{\PYGZdq{}}\PYG{l+s}{ = }\PYG{l+s}{\PYGZdq{}}\PYG{+w}{ }\PYG{o}{+}\PYG{+w}{ }\PYG{n}{x}\PYG{p}{.}\PYG{n+na}{get}\PYG{p}{(}\PYG{n}{copt}\PYG{p}{.}\PYG{n+na}{DblInfo}\PYG{p}{.}\PYG{n+na}{Value}\PYG{p}{)}\PYG{p}{)}\PYG{p}{;}
\PYG{+w}{        }\PYG{p}{\PYGZcb{}}
\PYG{+w}{        }\PYG{n}{System}\PYG{p}{.}\PYG{n+na}{out}\PYG{p}{.}\PYG{n+na}{println}\PYG{p}{(}\PYG{l+s}{\PYGZdq{}}\PYG{l+s}{Obj = }\PYG{l+s}{\PYGZdq{}}\PYG{+w}{ }\PYG{o}{+}\PYG{+w}{ }\PYG{n}{model}\PYG{p}{.}\PYG{n+na}{getDblAttr}\PYG{p}{(}\PYG{n}{copt}\PYG{p}{.}\PYG{n+na}{DblAttr}\PYG{p}{.}\PYG{n+na}{LpObjVal}\PYG{p}{)}\PYG{p}{)}\PYG{p}{;}
\PYG{+w}{      }\PYG{p}{\PYGZcb{}}

\PYG{+w}{      }\PYG{n}{System}\PYG{p}{.}\PYG{n+na}{out}\PYG{p}{.}\PYG{n+na}{println}\PYG{p}{(}\PYG{l+s}{\PYGZdq{}}\PYG{l+s}{\PYGZbs{}}\PYG{l+s}{nDone}\PYG{l+s}{\PYGZdq{}}\PYG{p}{)}\PYG{p}{;}
\end{sphinxVerbatim}

\sphinxAtStartPar
Spcifically, one can query the copt.IntAttr.HasLpSol attribute on the model to know whether
we have optimal LP solution; query the copt.DblInfo.Value attribute of a variable to obtain its
solution value; query the copt.DblAttr.LpObjVal attribute on the model to obtain the
objective value for the current solution.

\sphinxAtStartPar
The names and types of all model, variable, and constraint attributes can be found in
{\hyperref[\detokenize{javaapiref:chapjavaapiref-const}]{\sphinxcrossref{\DUrole{std,std-ref}{Constants}}}} of Java API reference.

\subsubsection{Error Handling}
\label{\detokenize{javainterface:error-handling}}
\sphinxAtStartPar
Errors in the COPT Java interface are handled through the Java exception mechanism.
In the example, all COPT statements are enclosed inside a try block, and any associated
errors would be caught by the catch block.

\begin{sphinxVerbatim}[commandchars=\\\{\}]
\PYG{+w}{    }\PYG{p}{\PYGZcb{}}\PYG{+w}{ }\PYG{k}{catch}\PYG{+w}{ }\PYG{p}{(}\PYG{n}{CoptException}\PYG{+w}{ }\PYG{n}{e}\PYG{p}{)}\PYG{+w}{ }\PYG{p}{\PYGZob{}}
\PYG{+w}{      }\PYG{n}{System}\PYG{p}{.}\PYG{n+na}{out}\PYG{p}{.}\PYG{n+na}{println}\PYG{p}{(}\PYG{l+s}{\PYGZdq{}}\PYG{l+s}{Error Code = }\PYG{l+s}{\PYGZdq{}}\PYG{+w}{ }\PYG{o}{+}\PYG{+w}{ }\PYG{n}{e}\PYG{p}{.}\PYG{n+na}{getCode}\PYG{p}{(}\PYG{p}{)}\PYG{p}{)}\PYG{p}{;}
\PYG{+w}{      }\PYG{n}{System}\PYG{p}{.}\PYG{n+na}{out}\PYG{p}{.}\PYG{n+na}{println}\PYG{p}{(}\PYG{n}{e}\PYG{p}{.}\PYG{n+na}{getMessage}\PYG{p}{(}\PYG{p}{)}\PYG{p}{)}\PYG{p}{;}
\PYG{+w}{    }\PYG{p}{\PYGZcb{}}
\end{sphinxVerbatim}

\subsection{Build and Run}
\label{\detokenize{javainterface:build-and-run}}
\sphinxAtStartPar
To build and run java example, users may refer to files under
\sphinxcode{\sphinxupquote{\$COPT\_HOME/examples/java}}. Specifically, We provide an example file in java and
a script file to build. This single example runs on all platforms that support Java.

\sphinxAtStartPar
First of all, download and install \sphinxtitleref{Java 8 or above} on your platform.

\subsubsection{Java example detail}
\label{\detokenize{javainterface:java-example-detail}}
\sphinxAtStartPar
In the java example folder \sphinxcode{\sphinxupquote{\$COPT\_HOME/examples/java}}, the easiest way to
run the example is to enter the java example folder in console or terminal
and then execute command \sphinxcode{\sphinxupquote{\textquotesingle{}sh run.sh\textquotesingle{}}}.

\sphinxAtStartPar
This java project has dependency on COPT java package \sphinxcode{\sphinxupquote{copt\_java.jar}}, which
defines all java classes of COPT solver.
In addition, \sphinxcode{\sphinxupquote{copt\_java.jar}} loads two shared libraries, that is,
\sphinxcode{\sphinxupquote{coptjniwrap.dll}} and \sphinxcode{\sphinxupquote{copt\_cpp.dll}} on Windows,  \sphinxcode{\sphinxupquote{libcoptjniwrap.so}}
and \sphinxcode{\sphinxupquote{libcopt\_cpp.so}} on Linux,  \sphinxcode{\sphinxupquote{libcoptjniwrap.dylib}} and
\sphinxcode{\sphinxupquote{libcopt\_cpp.dylib}} on Mac respectively. Note that \sphinxcode{\sphinxupquote{coptjniwrap}}
library is a JNI swig wrapper and acts as a bridge between COPT java package
and native library \sphinxcode{\sphinxupquote{copt\_cpp}}, which declares and implements COPT constants,
interfaces and methods. So users should make sure they are installed properly
on runtime search paths.

\sphinxAtStartPar
In summary, to run java example, users should have COPT installed properly.
Specifically, it requires two related COPT shared libaries existing on runtime
search paths, and valid license files to run. Please refer to
{\hyperref[\detokenize{install:chapinstall}]{\sphinxcrossref{\DUrole{std,std-ref}{Install Guide for Cardinal Optimizer}}}} for further details.

\sphinxstepscope

\section{Python Interface}
\label{\detokenize{pythoninterface:python-interface}}\label{\detokenize{pythoninterface:chappythoninterface}}\label{\detokenize{pythoninterface::doc}}

\subsection{Installation guide}
\label{\detokenize{pythoninterface:installation-guide}}
\sphinxAtStartPar
Currently, the Python interface of Cardinal Optimizer supports Python 3.7\sphinxhyphen{}3.13.
Among them, for Python 3.8\sphinxhyphen{}3.13 versions, COPT’s MacOS\sphinxhyphen{}Universal can provide support. For Python 3.7, only MacOS\sphinxhyphen{}X86 is available.

\sphinxAtStartPar
Before using the Python interface, please ensure that
COPT has been installed and configured correctly. For details, please refer to
{\hyperref[\detokenize{install:chapinstall}]{\sphinxcrossref{\DUrole{std,std-ref}{How to install the Cardinal Optimizer}}}}.
Users can download Python from \sphinxhref{https://www.anaconda.com/distribution/}{Anaconda distribution} or
\sphinxhref{https://www.python.org/}{Python official distribution} . We recommend users
install the Anaconda distribution, because it is more user\sphinxhyphen{}friendly and
convenient for Python novices (For Windows, please don’t install Python via
Microsoft Store). If you use official Python distribution or Python
shipped with system, then make sure you have installed the \sphinxcode{\sphinxupquote{pip}} and \sphinxcode{\sphinxupquote{setuptools}}
Python packages beforehand.

\begin{sphinxadmonition}{note}{Note}

\sphinxAtStartPar
We recommend using versions 3.8\sphinxhyphen{}3.13 because the minimum version requirement for Matrix Modeling function of the COPT\sphinxhyphen{}Python interface is 3.8.
\end{sphinxadmonition}

\subsubsection{Windows}
\label{\detokenize{pythoninterface:windows}}
\sphinxAtStartPar
\sphinxstylestrong{Method 1: via pip install (recommended)}

\sphinxAtStartPar
Open the cmd command window (if Python is an Anaconda distribution, open the Anaconda command line window) and enter the following command:

\begin{sphinxVerbatim}[commandchars=\\\{\}]
pip\PYG{+w}{ }install\PYG{+w}{ }coptpy
\end{sphinxVerbatim}

\sphinxAtStartPar
If an older version of the \sphinxcode{\sphinxupquote{coptpy}} package has been installed, please open the cmd command window (if Python is an Anaconda distribution, open the Anaconda command line window) and
enter the following command to upgrade to the latest version of the \sphinxcode{\sphinxupquote{coptpy}} package:

\begin{sphinxVerbatim}[commandchars=\\\{\}]
pip\PYG{+w}{ }install\PYG{+w}{ }\PYGZhy{}\PYGZhy{}upgrade\PYG{+w}{ }coptpy
\end{sphinxVerbatim}

\sphinxAtStartPar
\sphinxstylestrong{Method 2: via COPT installation package}

\sphinxAtStartPar
For Windows, assuming the installation path of COPT is: \sphinxcode{\sphinxupquote{"C:\textbackslash{}Program Files\textbackslash{}COPT"}},
please switch to the directory \sphinxcode{\sphinxupquote{"C:\textbackslash{}Program Files\textbackslash{}COPT\textbackslash{}lib\textbackslash{}python"}} and
execute the following commands on command line:

\begin{sphinxVerbatim}[commandchars=\\\{\}]
python\PYG{+w}{ }setup.py\PYG{+w}{ }install
\end{sphinxVerbatim}

\sphinxAtStartPar
Note that if COPT is installed on the system disk, you need to
\sphinxstylestrong{execute with administrator privileges} to open the command prompt.
To test whether the Python interface is installed correctly, users can switch to
the directory \sphinxcode{\sphinxupquote{"C:\textbackslash{}Program Files\textbackslash{}COPT\textbackslash{}examples\textbackslash{}python"}} and execute
the following commands on the command line:

\begin{sphinxVerbatim}[commandchars=\\\{\}]
python\PYG{+w}{ }lp\PYGZus{}ex1.py
\end{sphinxVerbatim}

\sphinxAtStartPar
If the model is solved correctly, it means that the Python interface of COPT
has been installed correctly.

\sphinxAtStartPar
\sphinxstylestrong{Note} If you use the official release of Python 3.8, assume that its installation path is:
\sphinxcode{\sphinxupquote{"C:\textbackslash{}Program Files\textbackslash{}Python38"}}, you need to copy the \sphinxcode{\sphinxupquote{copt\_cpp.dll}} file in the
\sphinxcode{\sphinxupquote{"bin"}} subdirectory of the COPT installation path to \sphinxcode{\sphinxupquote{"C:\textbackslash{}Program Files\textbackslash{}Python38\textbackslash{}Lib\textbackslash{}site\sphinxhyphen{}packages\textbackslash{}coptpy"}}
to solve the problem of dynamic library dependency.

\sphinxAtStartPar
Currently, \sphinxcode{\sphinxupquote{coptpy}} supports type hints. After \sphinxcode{\sphinxupquote{coptpy}} is successfully installed, when you code in Python IDE,
variable name completion and function prototype information will be automatically prompted.

\subsubsection{Linux}
\label{\detokenize{pythoninterface:linux}}
\sphinxAtStartPar
\sphinxstylestrong{Method 1: via pip install (recommended)}

\sphinxAtStartPar
Open the terminal and enter the following command:

\begin{sphinxVerbatim}[commandchars=\\\{\}]
pip\PYG{+w}{ }install\PYG{+w}{ }coptpy
\end{sphinxVerbatim}

\sphinxAtStartPar
If an older version of the \sphinxcode{\sphinxupquote{coptpy}} package has been installed, please open the terminal and
enter the following command to upgrade to the latest version of the \sphinxcode{\sphinxupquote{coptpy}} package:

\begin{sphinxVerbatim}[commandchars=\\\{\}]
pip\PYG{+w}{ }install\PYG{+w}{ }\PYGZhy{}\PYGZhy{}upgrade\PYG{+w}{ }coptpy
\end{sphinxVerbatim}

\sphinxAtStartPar
\sphinxstylestrong{Method 2: via COPT installation package}

\sphinxAtStartPar
For Linux, suppose the installation path of COPT is: /opt/copt80,
please switch to the directory /opt/copt80/lib/python  and execute
the following commands on terminal:

\begin{sphinxVerbatim}[commandchars=\\\{\}]
sudo\PYG{+w}{ }python\PYG{+w}{ }setup.py\PYG{+w}{ }install
\end{sphinxVerbatim}

\sphinxAtStartPar
For users using Python from Anaconda distribution, if above commands fails,
assuming the installation path of Anaconda is: \sphinxcode{\sphinxupquote{"/opt/anaconda3"}},
please execute the following commands instead on terminal to install
the Python interface of COPT:

\begin{sphinxVerbatim}[commandchars=\\\{\}]
sudo\PYG{+w}{ }/opt/anaconda3/bin/python\PYG{+w}{ }setup.py\PYG{+w}{ }install
\end{sphinxVerbatim}

\sphinxAtStartPar
To test whether the Python interface is installed correctly, users can switch
to the directory /opt/copt80/examples/python and execute the following
commands on terminal:

\begin{sphinxVerbatim}[commandchars=\\\{\}]
python\PYG{+w}{ }lp\PYGZus{}ex1.py
\end{sphinxVerbatim}

\sphinxAtStartPar
If the model solved correctly, it means that the Python interface of COPT
has been installed correctly.

\sphinxAtStartPar
Currently, \sphinxcode{\sphinxupquote{coptpy}} supports type hints. After \sphinxcode{\sphinxupquote{coptpy}} is successfully installed, when you code in Python IDE,
variable name completion and function prototype information will be automatically prompted.

\subsubsection{MacOS}
\label{\detokenize{pythoninterface:macos}}
\sphinxAtStartPar
\sphinxstylestrong{Method 1: via pip install (recommended)}

\sphinxAtStartPar
Open the terminal and enter the following command:

\begin{sphinxVerbatim}[commandchars=\\\{\}]
pip\PYG{+w}{ }install\PYG{+w}{ }coptpy
\end{sphinxVerbatim}

\sphinxAtStartPar
If an older version of the \sphinxcode{\sphinxupquote{coptpy}} package has been installed, please open the terminal and
enter the following command to upgrade to the latest version of the \sphinxcode{\sphinxupquote{coptpy}} package:

\begin{sphinxVerbatim}[commandchars=\\\{\}]
pip\PYG{+w}{ }install\PYG{+w}{ }\PYGZhy{}\PYGZhy{}upgrade\PYG{+w}{ }coptpy
\end{sphinxVerbatim}

\sphinxAtStartPar
\sphinxstylestrong{Method 2: via COPT installation package}

\sphinxAtStartPar
For MacOS, assuming that the installation path of COPT is: /Applications/copt80,
please switch the directory to /Applications/copt80/lib/python and execute
the following commands on terminal:

\begin{sphinxVerbatim}[commandchars=\\\{\}]
sudo\PYG{+w}{ }python\PYG{+w}{ }setup.py\PYG{+w}{ }install
\end{sphinxVerbatim}

\sphinxAtStartPar
To test whether the Python interface is installed correctly, users can switch
to the directory /Applications/copt80/examples/python and execute
the following commands on terminal:

\begin{sphinxVerbatim}[commandchars=\\\{\}]
python\PYG{+w}{ }lp\PYGZus{}ex1.py
\end{sphinxVerbatim}

\sphinxAtStartPar
If the model solved correctly, it means that the Python interface of COPT
has been installed correctly.

\sphinxAtStartPar
Currently, \sphinxcode{\sphinxupquote{coptpy}} supports type hints. After \sphinxcode{\sphinxupquote{coptpy}} is successfully installed, when you code in Python IDE,
variable name completion and function prototype information will be automatically prompted.

\subsection{Example details}
\label{\detokenize{pythoninterface:example-details}}
\sphinxAtStartPar
This chapter illustrate the use of C interface of Cardinal Optimizer through
a simple Python example. The problem to solve is shown in Eq. \ref{equation:pythoninterface:coptEq_pyapilp1}:
\begin{equation}\label{equation:pythoninterface:coptEq_pyapilp1}
\begin{split}\text{Maximize: } & \\
                    & 1.2 x + 1.8 y + 2.1 z \\
\text{Subject to: } & \\
                      & 1.5 x + 1.2 y + 1.8 z \leq 2.6 \\
                      & 0.8 x + 0.6 y + 0.9 z \geq 1.2 \\
\text{Bounds: } & \\
             & 0.1 \leq x \leq 0.6 \\
             & 0.2 \leq y \leq 1.5 \\
             & 0.3 \leq z \leq 2.8\end{split}
\end{equation}
\sphinxAtStartPar
The source code for solving the above problem using Python API of Cardinal Optimizer
is shown in \hyperref[\detokenize{pythoninterface:coptcode-pyapilp1}]{Listing \ref{\detokenize{pythoninterface:coptcode-pyapilp1}}}:
\sphinxSetupCaptionForVerbatim{\sphinxcode{\sphinxupquote{lp\_ex1.py}}}
\def\sphinxLiteralBlockLabel{\label{\detokenize{pythoninterface:coptcode-pyapilp1}}}
\begin{sphinxVerbatim}[commandchars=\\\{\},numbers=left,firstnumber=1,stepnumber=1]
\PYG{c+c1}{\PYGZsh{}}
\PYG{c+c1}{\PYGZsh{} This file is part of the Cardinal Optimizer, all rights reserved.}
\PYG{c+c1}{\PYGZsh{}}

\PYG{l+s+sd}{\PYGZdq{}\PYGZdq{}\PYGZdq{}}
\PYG{l+s+sd}{The problem to solve:}

\PYG{l+s+sd}{  Maximize:}
\PYG{l+s+sd}{    1.2 x + 1.8 y + 2.1 z}

\PYG{l+s+sd}{  Subject to:}
\PYG{l+s+sd}{    1.5 x + 1.2 y + 1.8 z \PYGZlt{}= 2.6}
\PYG{l+s+sd}{    0.8 x + 0.6 y + 0.9 z \PYGZgt{}= 1.2}

\PYG{l+s+sd}{  where:}
\PYG{l+s+sd}{    0.1 \PYGZlt{}= x \PYGZlt{}= 0.6}
\PYG{l+s+sd}{    0.2 \PYGZlt{}= y \PYGZlt{}= 1.5}
\PYG{l+s+sd}{    0.3 \PYGZlt{}= z \PYGZlt{}= 2.8}
\PYG{l+s+sd}{\PYGZdq{}\PYGZdq{}\PYGZdq{}}

\PYG{k+kn}{import}\PYG{+w}{ }\PYG{n+nn}{coptpy}\PYG{+w}{ }\PYG{k}{as}\PYG{+w}{ }\PYG{n+nn}{cp}
\PYG{k+kn}{from}\PYG{+w}{ }\PYG{n+nn}{coptpy}\PYG{+w}{ }\PYG{k+kn}{import} \PYG{n}{COPT}

\PYG{c+c1}{\PYGZsh{} Create COPT environment}
\PYG{n}{env} \PYG{o}{=} \PYG{n}{cp}\PYG{o}{.}\PYG{n}{Envr}\PYG{p}{(}\PYG{p}{)}

\PYG{c+c1}{\PYGZsh{} Create COPT model}
\PYG{n}{model} \PYG{o}{=} \PYG{n}{env}\PYG{o}{.}\PYG{n}{createModel}\PYG{p}{(}\PYG{l+s+s2}{\PYGZdq{}}\PYG{l+s+s2}{lp\PYGZus{}ex1}\PYG{l+s+s2}{\PYGZdq{}}\PYG{p}{)}

\PYG{c+c1}{\PYGZsh{} Add variables: x, y, z}
\PYG{n}{x} \PYG{o}{=} \PYG{n}{model}\PYG{o}{.}\PYG{n}{addVar}\PYG{p}{(}\PYG{n}{lb}\PYG{o}{=}\PYG{l+m+mf}{0.1}\PYG{p}{,} \PYG{n}{ub}\PYG{o}{=}\PYG{l+m+mf}{0.6}\PYG{p}{,} \PYG{n}{name}\PYG{o}{=}\PYG{l+s+s2}{\PYGZdq{}}\PYG{l+s+s2}{x}\PYG{l+s+s2}{\PYGZdq{}}\PYG{p}{)}
\PYG{n}{y} \PYG{o}{=} \PYG{n}{model}\PYG{o}{.}\PYG{n}{addVar}\PYG{p}{(}\PYG{n}{lb}\PYG{o}{=}\PYG{l+m+mf}{0.2}\PYG{p}{,} \PYG{n}{ub}\PYG{o}{=}\PYG{l+m+mf}{1.5}\PYG{p}{,} \PYG{n}{name}\PYG{o}{=}\PYG{l+s+s2}{\PYGZdq{}}\PYG{l+s+s2}{y}\PYG{l+s+s2}{\PYGZdq{}}\PYG{p}{)}
\PYG{n}{z} \PYG{o}{=} \PYG{n}{model}\PYG{o}{.}\PYG{n}{addVar}\PYG{p}{(}\PYG{n}{lb}\PYG{o}{=}\PYG{l+m+mf}{0.3}\PYG{p}{,} \PYG{n}{ub}\PYG{o}{=}\PYG{l+m+mf}{2.8}\PYG{p}{,} \PYG{n}{name}\PYG{o}{=}\PYG{l+s+s2}{\PYGZdq{}}\PYG{l+s+s2}{z}\PYG{l+s+s2}{\PYGZdq{}}\PYG{p}{)}

\PYG{c+c1}{\PYGZsh{} Add constraints}
\PYG{n}{model}\PYG{o}{.}\PYG{n}{addConstr}\PYG{p}{(}\PYG{l+m+mf}{1.5}\PYG{o}{*}\PYG{n}{x} \PYG{o}{+} \PYG{l+m+mf}{1.2}\PYG{o}{*}\PYG{n}{y} \PYG{o}{+} \PYG{l+m+mf}{1.8}\PYG{o}{*}\PYG{n}{z} \PYG{o}{\PYGZlt{}}\PYG{o}{=} \PYG{l+m+mf}{2.6}\PYG{p}{)}
\PYG{n}{model}\PYG{o}{.}\PYG{n}{addConstr}\PYG{p}{(}\PYG{l+m+mf}{0.8}\PYG{o}{*}\PYG{n}{x} \PYG{o}{+} \PYG{l+m+mf}{0.6}\PYG{o}{*}\PYG{n}{y} \PYG{o}{+} \PYG{l+m+mf}{0.9}\PYG{o}{*}\PYG{n}{z} \PYG{o}{\PYGZgt{}}\PYG{o}{=} \PYG{l+m+mf}{1.2}\PYG{p}{)}

\PYG{c+c1}{\PYGZsh{} Set objective function}
\PYG{n}{model}\PYG{o}{.}\PYG{n}{setObjective}\PYG{p}{(}\PYG{l+m+mf}{1.2}\PYG{o}{*}\PYG{n}{x} \PYG{o}{+} \PYG{l+m+mf}{1.8}\PYG{o}{*}\PYG{n}{y} \PYG{o}{+} \PYG{l+m+mf}{2.1}\PYG{o}{*}\PYG{n}{z}\PYG{p}{,} \PYG{n}{sense}\PYG{o}{=}\PYG{n}{COPT}\PYG{o}{.}\PYG{n}{MAXIMIZE}\PYG{p}{)}

\PYG{c+c1}{\PYGZsh{} Set parameter}
\PYG{n}{model}\PYG{o}{.}\PYG{n}{setParam}\PYG{p}{(}\PYG{n}{COPT}\PYG{o}{.}\PYG{n}{Param}\PYG{o}{.}\PYG{n}{TimeLimit}\PYG{p}{,} \PYG{l+m+mf}{10.0}\PYG{p}{)}

\PYG{c+c1}{\PYGZsh{} Solve the model}
\PYG{n}{model}\PYG{o}{.}\PYG{n}{solve}\PYG{p}{(}\PYG{p}{)}

\PYG{c+c1}{\PYGZsh{} Analyze solution}
\PYG{k}{if} \PYG{n}{model}\PYG{o}{.}\PYG{n}{status} \PYG{o}{==} \PYG{n}{COPT}\PYG{o}{.}\PYG{n}{OPTIMAL}\PYG{p}{:}
  \PYG{n+nb}{print}\PYG{p}{(}\PYG{l+s+s2}{\PYGZdq{}}\PYG{l+s+s2}{Objective value: }\PYG{l+s+si}{\PYGZob{}\PYGZcb{}}\PYG{l+s+s2}{\PYGZdq{}}\PYG{o}{.}\PYG{n}{format}\PYG{p}{(}\PYG{n}{model}\PYG{o}{.}\PYG{n}{objval}\PYG{p}{)}\PYG{p}{)}
  \PYG{n}{allvars} \PYG{o}{=} \PYG{n}{model}\PYG{o}{.}\PYG{n}{getVars}\PYG{p}{(}\PYG{p}{)}

  \PYG{n+nb}{print}\PYG{p}{(}\PYG{l+s+s2}{\PYGZdq{}}\PYG{l+s+s2}{Variable solution:}\PYG{l+s+s2}{\PYGZdq{}}\PYG{p}{)}
  \PYG{k}{for} \PYG{n}{var} \PYG{o+ow}{in} \PYG{n}{allvars}\PYG{p}{:}
    \PYG{n+nb}{print}\PYG{p}{(}\PYG{l+s+s2}{\PYGZdq{}}\PYG{l+s+s2}{ x[}\PYG{l+s+si}{\PYGZob{}0\PYGZcb{}}\PYG{l+s+s2}{]: }\PYG{l+s+si}{\PYGZob{}1\PYGZcb{}}\PYG{l+s+s2}{\PYGZdq{}}\PYG{o}{.}\PYG{n}{format}\PYG{p}{(}\PYG{n}{var}\PYG{o}{.}\PYG{n}{index}\PYG{p}{,} \PYG{n}{var}\PYG{o}{.}\PYG{n}{x}\PYG{p}{)}\PYG{p}{)}

  \PYG{n+nb}{print}\PYG{p}{(}\PYG{l+s+s2}{\PYGZdq{}}\PYG{l+s+s2}{Variable basis status:}\PYG{l+s+s2}{\PYGZdq{}}\PYG{p}{)}
  \PYG{k}{for} \PYG{n}{var} \PYG{o+ow}{in} \PYG{n}{allvars}\PYG{p}{:}
    \PYG{n+nb}{print}\PYG{p}{(}\PYG{l+s+s2}{\PYGZdq{}}\PYG{l+s+s2}{ x[}\PYG{l+s+si}{\PYGZob{}0\PYGZcb{}}\PYG{l+s+s2}{]: }\PYG{l+s+si}{\PYGZob{}1\PYGZcb{}}\PYG{l+s+s2}{\PYGZdq{}}\PYG{o}{.}\PYG{n}{format}\PYG{p}{(}\PYG{n}{var}\PYG{o}{.}\PYG{n}{index}\PYG{p}{,} \PYG{n}{var}\PYG{o}{.}\PYG{n}{basis}\PYG{p}{)}\PYG{p}{)}

  \PYG{c+c1}{\PYGZsh{} Write model, solution and modified parameters to file}
  \PYG{n}{model}\PYG{o}{.}\PYG{n}{write}\PYG{p}{(}\PYG{l+s+s2}{\PYGZdq{}}\PYG{l+s+s2}{lp\PYGZus{}ex1.mps}\PYG{l+s+s2}{\PYGZdq{}}\PYG{p}{)}
  \PYG{n}{model}\PYG{o}{.}\PYG{n}{write}\PYG{p}{(}\PYG{l+s+s2}{\PYGZdq{}}\PYG{l+s+s2}{lp\PYGZus{}ex1.bas}\PYG{l+s+s2}{\PYGZdq{}}\PYG{p}{)}
  \PYG{n}{model}\PYG{o}{.}\PYG{n}{write}\PYG{p}{(}\PYG{l+s+s2}{\PYGZdq{}}\PYG{l+s+s2}{lp\PYGZus{}ex1.sol}\PYG{l+s+s2}{\PYGZdq{}}\PYG{p}{)}
  \PYG{n}{model}\PYG{o}{.}\PYG{n}{write}\PYG{p}{(}\PYG{l+s+s2}{\PYGZdq{}}\PYG{l+s+s2}{lp\PYGZus{}ex1.par}\PYG{l+s+s2}{\PYGZdq{}}\PYG{p}{)}
\end{sphinxVerbatim}

\sphinxAtStartPar
We will explain how to use the Python API step by step based on code above, please
refer to {\hyperref[\detokenize{capiref:chapapi}]{\sphinxcrossref{\DUrole{std,std-ref}{C API Reference}}}} for detailed usage of Python API.

\subsubsection{Import Python interface}
\label{\detokenize{pythoninterface:import-python-interface}}
\sphinxAtStartPar
To use the Python interface of COPT, users need to import the Python interface
library first.

\begin{sphinxVerbatim}[commandchars=\\\{\}]
\PYG{k+kn}{import}\PYG{+w}{ }\PYG{n+nn}{coptpy}\PYG{+w}{ }\PYG{k}{as}\PYG{+w}{ }\PYG{n+nn}{cp}
\PYG{k+kn}{from}\PYG{+w}{ }\PYG{n+nn}{coptpy}\PYG{+w}{ }\PYG{k+kn}{import} \PYG{n}{COPT}
\end{sphinxVerbatim}

\subsubsection{Create environment}
\label{\detokenize{pythoninterface:create-environment}}
\sphinxAtStartPar
To solve any problem with COPT, users need to create optimization environment before
creating any model.

\begin{sphinxVerbatim}[commandchars=\\\{\}]
\PYG{c+c1}{\PYGZsh{} Create COPT environment}
\PYG{n}{env} \PYG{o}{=} \PYG{n}{cp}\PYG{o}{.}\PYG{n}{Envr}\PYG{p}{(}\PYG{p}{)}
\end{sphinxVerbatim}

\subsubsection{Create model}
\label{\detokenize{pythoninterface:create-model}}
\sphinxAtStartPar
If the optimization environment was created successfully, users need to create
the model to solve, which includes variables and constraints information.

\begin{sphinxVerbatim}[commandchars=\\\{\}]
\PYG{c+c1}{\PYGZsh{} Create COPT model}
\PYG{n}{model} \PYG{o}{=} \PYG{n}{env}\PYG{o}{.}\PYG{n}{createModel}\PYG{p}{(}\PYG{l+s+s2}{\PYGZdq{}}\PYG{l+s+s2}{lp\PYGZus{}ex1}\PYG{l+s+s2}{\PYGZdq{}}\PYG{p}{)}
\end{sphinxVerbatim}

\subsubsection{Add variables}
\label{\detokenize{pythoninterface:add-variables}}
\sphinxAtStartPar
Users can specify information such as objective costs, lower and upper bounds
of variables when creating them. In this example, we just set the lower and upper
bounds of variables and their names.

\begin{sphinxVerbatim}[commandchars=\\\{\}]
\PYG{c+c1}{\PYGZsh{} Add variables: x, y, z}
\PYG{n}{x} \PYG{o}{=} \PYG{n}{model}\PYG{o}{.}\PYG{n}{addVar}\PYG{p}{(}\PYG{n}{lb}\PYG{o}{=}\PYG{l+m+mf}{0.1}\PYG{p}{,} \PYG{n}{ub}\PYG{o}{=}\PYG{l+m+mf}{0.6}\PYG{p}{,} \PYG{n}{name}\PYG{o}{=}\PYG{l+s+s2}{\PYGZdq{}}\PYG{l+s+s2}{x}\PYG{l+s+s2}{\PYGZdq{}}\PYG{p}{)}
\PYG{n}{y} \PYG{o}{=} \PYG{n}{model}\PYG{o}{.}\PYG{n}{addVar}\PYG{p}{(}\PYG{n}{lb}\PYG{o}{=}\PYG{l+m+mf}{0.2}\PYG{p}{,} \PYG{n}{ub}\PYG{o}{=}\PYG{l+m+mf}{1.5}\PYG{p}{,} \PYG{n}{name}\PYG{o}{=}\PYG{l+s+s2}{\PYGZdq{}}\PYG{l+s+s2}{y}\PYG{l+s+s2}{\PYGZdq{}}\PYG{p}{)}
\PYG{n}{z} \PYG{o}{=} \PYG{n}{model}\PYG{o}{.}\PYG{n}{addVar}\PYG{p}{(}\PYG{n}{lb}\PYG{o}{=}\PYG{l+m+mf}{0.3}\PYG{p}{,} \PYG{n}{ub}\PYG{o}{=}\PYG{l+m+mf}{2.8}\PYG{p}{,} \PYG{n}{name}\PYG{o}{=}\PYG{l+s+s2}{\PYGZdq{}}\PYG{l+s+s2}{z}\PYG{l+s+s2}{\PYGZdq{}}\PYG{p}{)}
\end{sphinxVerbatim}

\subsubsection{Add constraints}
\label{\detokenize{pythoninterface:add-constraints}}
\sphinxAtStartPar
After adding variables, we can then add constraints to the model.

\begin{sphinxVerbatim}[commandchars=\\\{\}]
\PYG{c+c1}{\PYGZsh{} Add constraints}
\PYG{n}{model}\PYG{o}{.}\PYG{n}{addConstr}\PYG{p}{(}\PYG{l+m+mf}{1.5}\PYG{o}{*}\PYG{n}{x} \PYG{o}{+} \PYG{l+m+mf}{1.2}\PYG{o}{*}\PYG{n}{y} \PYG{o}{+} \PYG{l+m+mf}{1.8}\PYG{o}{*}\PYG{n}{z} \PYG{o}{\PYGZlt{}}\PYG{o}{=} \PYG{l+m+mf}{2.6}\PYG{p}{)}
\PYG{n}{model}\PYG{o}{.}\PYG{n}{addConstr}\PYG{p}{(}\PYG{l+m+mf}{0.8}\PYG{o}{*}\PYG{n}{x} \PYG{o}{+} \PYG{l+m+mf}{0.6}\PYG{o}{*}\PYG{n}{y} \PYG{o}{+} \PYG{l+m+mf}{0.9}\PYG{o}{*}\PYG{n}{z} \PYG{o}{\PYGZgt{}}\PYG{o}{=} \PYG{l+m+mf}{1.2}\PYG{p}{)}
\end{sphinxVerbatim}

\subsubsection{Set objective function}
\label{\detokenize{pythoninterface:set-objective-function}}
\sphinxAtStartPar
After adding variables and constraints, we can further specify objective function
for the model.

\begin{sphinxVerbatim}[commandchars=\\\{\}]
\PYG{c+c1}{\PYGZsh{} Set objective function}
\PYG{n}{model}\PYG{o}{.}\PYG{n}{setObjective}\PYG{p}{(}\PYG{l+m+mf}{1.2}\PYG{o}{*}\PYG{n}{x} \PYG{o}{+} \PYG{l+m+mf}{1.8}\PYG{o}{*}\PYG{n}{y} \PYG{o}{+} \PYG{l+m+mf}{2.1}\PYG{o}{*}\PYG{n}{z}\PYG{p}{,} \PYG{n}{sense}\PYG{o}{=}\PYG{n}{COPT}\PYG{o}{.}\PYG{n}{MAXIMIZE}\PYG{p}{)}
\end{sphinxVerbatim}

\subsubsection{Set parameters}
\label{\detokenize{pythoninterface:set-parameters}}
\sphinxAtStartPar
Users can set optimization parameters before solving the model, e.g. set optimization
time limit to 10 seconds.

\begin{sphinxVerbatim}[commandchars=\\\{\}]
\PYG{c+c1}{\PYGZsh{} Set parameter}
\PYG{n}{model}\PYG{o}{.}\PYG{n}{setParam}\PYG{p}{(}\PYG{n}{COPT}\PYG{o}{.}\PYG{n}{Param}\PYG{o}{.}\PYG{n}{TimeLimit}\PYG{p}{,} \PYG{l+m+mf}{10.0}\PYG{p}{)}
\end{sphinxVerbatim}

\subsubsection{Solve model}
\label{\detokenize{pythoninterface:solve-model}}
\sphinxAtStartPar
Solve the model via \sphinxcode{\sphinxupquote{solve}} method.

\begin{sphinxVerbatim}[commandchars=\\\{\}]
\PYG{c+c1}{\PYGZsh{} Solve the model}
\PYG{n}{model}\PYG{o}{.}\PYG{n}{solve}\PYG{p}{(}\PYG{p}{)}
\end{sphinxVerbatim}

\subsubsection{Analyze solution}
\label{\detokenize{pythoninterface:analyze-solution}}
\sphinxAtStartPar
When solving finished, we should query the solution status first. If the solution
status is optimal, then we can retrieve objective value, solution and basis status
of variables.

\begin{sphinxVerbatim}[commandchars=\\\{\}]
\PYG{c+c1}{\PYGZsh{} Analyze solution}
\PYG{k}{if} \PYG{n}{model}\PYG{o}{.}\PYG{n}{status} \PYG{o}{==} \PYG{n}{COPT}\PYG{o}{.}\PYG{n}{OPTIMAL}\PYG{p}{:}
  \PYG{n+nb}{print}\PYG{p}{(}\PYG{l+s+s2}{\PYGZdq{}}\PYG{l+s+s2}{Objective value: }\PYG{l+s+si}{\PYGZob{}\PYGZcb{}}\PYG{l+s+s2}{\PYGZdq{}}\PYG{o}{.}\PYG{n}{format}\PYG{p}{(}\PYG{n}{model}\PYG{o}{.}\PYG{n}{objval}\PYG{p}{)}\PYG{p}{)}
  \PYG{n}{allvars} \PYG{o}{=} \PYG{n}{model}\PYG{o}{.}\PYG{n}{getVars}\PYG{p}{(}\PYG{p}{)}

  \PYG{n+nb}{print}\PYG{p}{(}\PYG{l+s+s2}{\PYGZdq{}}\PYG{l+s+s2}{Variable solution:}\PYG{l+s+s2}{\PYGZdq{}}\PYG{p}{)}
  \PYG{k}{for} \PYG{n}{var} \PYG{o+ow}{in} \PYG{n}{allvars}\PYG{p}{:}
    \PYG{n+nb}{print}\PYG{p}{(}\PYG{l+s+s2}{\PYGZdq{}}\PYG{l+s+s2}{ x[}\PYG{l+s+si}{\PYGZob{}0\PYGZcb{}}\PYG{l+s+s2}{]: }\PYG{l+s+si}{\PYGZob{}1\PYGZcb{}}\PYG{l+s+s2}{\PYGZdq{}}\PYG{o}{.}\PYG{n}{format}\PYG{p}{(}\PYG{n}{var}\PYG{o}{.}\PYG{n}{index}\PYG{p}{,} \PYG{n}{var}\PYG{o}{.}\PYG{n}{x}\PYG{p}{)}\PYG{p}{)}

  \PYG{n+nb}{print}\PYG{p}{(}\PYG{l+s+s2}{\PYGZdq{}}\PYG{l+s+s2}{Variable basis status:}\PYG{l+s+s2}{\PYGZdq{}}\PYG{p}{)}
  \PYG{k}{for} \PYG{n}{var} \PYG{o+ow}{in} \PYG{n}{allvars}\PYG{p}{:}
    \PYG{n+nb}{print}\PYG{p}{(}\PYG{l+s+s2}{\PYGZdq{}}\PYG{l+s+s2}{ x[}\PYG{l+s+si}{\PYGZob{}0\PYGZcb{}}\PYG{l+s+s2}{]: }\PYG{l+s+si}{\PYGZob{}1\PYGZcb{}}\PYG{l+s+s2}{\PYGZdq{}}\PYG{o}{.}\PYG{n}{format}\PYG{p}{(}\PYG{n}{var}\PYG{o}{.}\PYG{n}{index}\PYG{p}{,} \PYG{n}{var}\PYG{o}{.}\PYG{n}{basis}\PYG{p}{)}\PYG{p}{)}
\end{sphinxVerbatim}

\subsubsection{Write files}
\label{\detokenize{pythoninterface:write-files}}
\sphinxAtStartPar
Users can write current model to MPS format file, and write solution, basis status
and modified parameters to file.

\begin{sphinxVerbatim}[commandchars=\\\{\}]
  \PYG{c+c1}{\PYGZsh{} Write model, solution and modified parameters to file}
  \PYG{n}{model}\PYG{o}{.}\PYG{n}{write}\PYG{p}{(}\PYG{l+s+s2}{\PYGZdq{}}\PYG{l+s+s2}{lp\PYGZus{}ex1.mps}\PYG{l+s+s2}{\PYGZdq{}}\PYG{p}{)}
  \PYG{n}{model}\PYG{o}{.}\PYG{n}{write}\PYG{p}{(}\PYG{l+s+s2}{\PYGZdq{}}\PYG{l+s+s2}{lp\PYGZus{}ex1.bas}\PYG{l+s+s2}{\PYGZdq{}}\PYG{p}{)}
  \PYG{n}{model}\PYG{o}{.}\PYG{n}{write}\PYG{p}{(}\PYG{l+s+s2}{\PYGZdq{}}\PYG{l+s+s2}{lp\PYGZus{}ex1.sol}\PYG{l+s+s2}{\PYGZdq{}}\PYG{p}{)}
  \PYG{n}{model}\PYG{o}{.}\PYG{n}{write}\PYG{p}{(}\PYG{l+s+s2}{\PYGZdq{}}\PYG{l+s+s2}{lp\PYGZus{}ex1.par}\PYG{l+s+s2}{\PYGZdq{}}\PYG{p}{)}
\end{sphinxVerbatim}

\subsection{Best Practice}
\label{\detokenize{pythoninterface:best-practice}}

\subsubsection{Upgrade to the newer version}
\label{\detokenize{pythoninterface:upgrade-to-the-newer-version}}\label{\detokenize{pythoninterface:chappythoninterface-practiceupgrade}}
\sphinxAtStartPar
If users have COPT Python Interface \sphinxcode{\sphinxupquote{coptpy}} installed and need to upgrade latest version, it is recommended
to remove previous version before installing new version. To remove previous version, it
is as simple as deleting the folder \sphinxcode{\sphinxupquote{coptpy}} at \sphinxcode{\sphinxupquote{site\sphinxhyphen{}package}}.

\sphinxAtStartPar
If the user is accessible to the Internet, we recommend using \sphinxcode{\sphinxupquote{pip}} to uninstall or upgrade \sphinxcode{\sphinxupquote{coptpy}}.

\sphinxAtStartPar
\sphinxstylestrong{Uninstall the old package}

\begin{sphinxVerbatim}[commandchars=\\\{\}]
pip\PYG{+w}{ }uninstall\PYG{+w}{ }coptpy
\end{sphinxVerbatim}

\sphinxAtStartPar
\sphinxstylestrong{Upgrade to the new package}

\begin{sphinxVerbatim}[commandchars=\\\{\}]
pip\PYG{+w}{ }install\PYG{+w}{ }\PYGZhy{}\PYGZhy{}upgrade\PYG{+w}{ }coptpy
\end{sphinxVerbatim}

\sphinxAtStartPar
If the user installed \sphinxcode{\sphinxupquote{coptpy}} from the offline COPT installation package by \sphinxcode{\sphinxupquote{python setup.py install}},
please refer to the following steps to uninstall and upgrade:

\sphinxAtStartPar
1.Firstly, please confirm the installation location of \sphinxcode{\sphinxupquote{coptpy}}, which can be viewed through \sphinxcode{\sphinxupquote{pip show coptpy}},
and the output is as follows:

\begin{sphinxVerbatim}[commandchars=\\\{\}]
Name: coptpy
Version: 8.0.1
Summary: The Python interface of the Cardinal Optimizer
Home\PYGZhy{}page: www.shanshu.ai
Author: Cardinal Operations, LLC
Author\PYGZhy{}email: coptsales@shanshu.ai
License: Cardinal Operations, LLC
End User License Agreement
Location: /Users/user/anaconda3/lib/python3.11/site\PYGZhy{}packages
Requires:
Required\PYGZhy{}by:
\end{sphinxVerbatim}

\sphinxAtStartPar
2.Secondly, please manually delete the \sphinxcode{\sphinxupquote{coptpy\sphinxhyphen{}*.egg\sphinxhyphen{}info}} and \sphinxcode{\sphinxupquote{coptpy}} under the \sphinxcode{\sphinxupquote{site\sphinxhyphen{}packages}} directory of the \sphinxcode{\sphinxupquote{Location}} above.

\sphinxAtStartPar
3.Finally, please install the new version of \sphinxcode{\sphinxupquote{coptpy}}: enter the COPT installation package directory (such as copt80/lib/python), and execute:

\begin{sphinxVerbatim}[commandchars=\\\{\}]
python\PYG{+w}{ }setup.py\PYG{+w}{ }install
\end{sphinxVerbatim}

\subsubsection{Multi\sphinxhyphen{}Thread Programming}
\label{\detokenize{pythoninterface:multi-thread-programming}}
\sphinxAtStartPar
COPT does not guarrantee thread safe and modelling APIs are not reentrant in general.
It is safe to share COPT Envr objects among threads. However, it is not recommended to
share Model objects among threads, unless you understand what you are doing. For instance,
if you share the same model between two threads. One thread is responsible for modelling
and solving. The other thread is used to monitor the progress and may interrupt at some
circumstances, such as running out of time.

\subsubsection{Dictionary order guarranteed after Python 3.7}
\label{\detokenize{pythoninterface:dictionary-order-guarranteed-after-python-3-7}}
\sphinxAtStartPar
As you know, Python dictionaries did not preserve the order in which items were added to them.
For instance, you might type \sphinxcode{\sphinxupquote{\{\textquotesingle{}fruits\textquotesingle{}: {[}\textquotesingle{}apple\textquotesingle{}, \textquotesingle{}orange\textquotesingle{}{]}, \textquotesingle{}veggies\textquotesingle{}: {[}\textquotesingle{}carrot\textquotesingle{}, \textquotesingle{}pea\textquotesingle{}{]}\}}}
and get back \sphinxcode{\sphinxupquote{\{\textquotesingle{}veggies\textquotesingle{}: {[}\textquotesingle{}carrot\textquotesingle{}, \textquotesingle{}pea\textquotesingle{}{]}, \textquotesingle{}fruits\textquotesingle{}: {[}\textquotesingle{}apple\textquotesingle{}, \textquotesingle{}orange\textquotesingle{}{]}\}}}.

\sphinxAtStartPar
However, the situation is changed. Standard dict objects preserve order in the implementation
of Python 3.6. This order\sphinxhyphen{}preserving property is becoming a language feature in Python 3.7.

\sphinxAtStartPar
If your program has dependency on dictionary orders, plese install COPT Python Interface \sphinxcode{\sphinxupquote{coptpy}} for Python 3.7 or later version.
For instance, if your model is implemented in Python 2.7 as follows:

\begin{sphinxVerbatim}[commandchars=\\\{\}]
\PYG{n}{m} \PYG{o}{=} \PYG{n}{Envr}\PYG{p}{(}\PYG{p}{)}\PYG{o}{.}\PYG{n}{createModel}\PYG{p}{(}\PYG{l+s+s2}{\PYGZdq{}}\PYG{l+s+s2}{customized model}\PYG{l+s+s2}{\PYGZdq{}}\PYG{p}{)}
\PYG{n}{vx} \PYG{o}{=} \PYG{n}{m}\PYG{o}{.}\PYG{n}{addVars}\PYG{p}{(}\PYG{p}{[}\PYG{l+s+s1}{\PYGZsq{}}\PYG{l+s+s1}{hello}\PYG{l+s+s1}{\PYGZsq{}}\PYG{p}{,} \PYG{l+s+s1}{\PYGZsq{}}\PYG{l+s+s1}{world}\PYG{l+s+s1}{\PYGZsq{}}\PYG{p}{]}\PYG{p}{,} \PYG{p}{[}\PYG{l+m+mi}{0}\PYG{p}{,} \PYG{l+m+mi}{1}\PYG{p}{,} \PYG{l+m+mi}{2}\PYG{p}{]}\PYG{p}{,} \PYG{n}{nameprefix} \PYG{o}{=} \PYG{l+s+s2}{\PYGZdq{}}\PYG{l+s+s2}{X}\PYG{l+s+s2}{\PYGZdq{}}\PYG{p}{)}
\PYG{c+c1}{\PYGZsh{} add a constraint for each var in tupledict \PYGZsq{}vx\PYGZsq{}}
\PYG{n}{m}\PYG{o}{.}\PYG{n}{addConstrs}\PYG{p}{(}\PYG{n}{vx}\PYG{p}{[}\PYG{n}{key}\PYG{p}{]} \PYG{o}{\PYGZgt{}}\PYG{o}{=} \PYG{l+m+mf}{1.0} \PYG{k}{for} \PYG{n}{key} \PYG{o+ow}{in} \PYG{n}{vx}\PYG{p}{)}
\end{sphinxVerbatim}

\sphinxAtStartPar
Your model might end up with rows \sphinxcode{\sphinxupquote{\{R(hello,1), R(hello,0), R(world,1), R(world,0), R(hello,2), R(world,2)\}}} .

\subsubsection{Use quicksum and psdquicksum when possible}
\label{\detokenize{pythoninterface:use-quicksum-and-psdquicksum-when-possible}}
\sphinxAtStartPar
The Python interface of COPT supports building linear expression, quadratic expression
and PSD expression in natural way.
For linear and quadratic expression, it is recommended to use quicksum() to build expression objects.
For linear and PSD expression, it is recommended to use psdquicksum() to build expression objects.
Both of them implement inplace summation, which is much faster than standard plus operator.

\subsubsection{Operate the model in batch}
\label{\detokenize{pythoninterface:operate-the-model-in-batch}}
\sphinxAtStartPar
The Python interface of COPT supports performing batch operations on the models, such as:
\begin{itemize}
\item {} 
\sphinxAtStartPar
Add multiple variables or constraints: \sphinxcode{\sphinxupquote{Model.addVars()/Model.addConstrs()}} .

\item {} 
\sphinxAtStartPar
Set the coefficients of multiple variables in the linear constraints: \sphinxcode{\sphinxupquote{Model.setCoeffs()}}
(\sphinxstylestrong{Note:} The index pairs of variables and constraints cannot appear repeatedly).

\item {} 
\sphinxAtStartPar
Set the names of multiple variables or constraintss: \sphinxcode{\sphinxupquote{Model.setNames()}} .

\end{itemize}

\sphinxAtStartPar
Please refer to {\hyperref[\detokenize{pyapiref:chappyapi-model}]{\sphinxcrossref{\DUrole{std,std-ref}{Python API Reference: Model Class}}}} for function descriptions.

\sphinxstepscope

\section{AMPL Interface}
\label{\detokenize{amplinterface:ampl-interface}}\label{\detokenize{amplinterface:chapamplinterface}}\label{\detokenize{amplinterface::doc}}
\sphinxAtStartPar
\sphinxhref{https://ampl.com}{AMPL} is an algebraic modeling language for describing
large\sphinxhyphen{}scale complex mathematical problems, it was hooked to many commercial
and open\sphinxhyphen{}source mathematical optimizers, with various data interfaces and
extensions, and received high popularity among both industries and institutes,
see \sphinxhref{https://ampl.com/about-us/customers/}{Who uses AMPL?} for more information.
The solver \sphinxcode{\sphinxupquote{coptampl}} (located in the \sphinxcode{\sphinxupquote{"\textbackslash{}bin"}} directory of the COPT installation package)
uses \sphinxstylestrong{Cardinal Optimizer} to solve linear programming, convex quadratic programming, convex quadratic constrainted programming and mixed integer programming problems. Normally \sphinxcode{\sphinxupquote{coptampl}} is invoked
by AMPL’s solve command, which gives the invocation:

\begin{sphinxVerbatim}[commandchars=\\\{\}]
coptampl\PYG{+w}{ }stub\PYG{+w}{ }\PYGZhy{}AMPL
\end{sphinxVerbatim}

\sphinxAtStartPar
in which \sphinxcode{\sphinxupquote{stub.nl}} is an AMPL generic output file (possibly written by
\sphinxcode{\sphinxupquote{\textquotesingle{}ampl \sphinxhyphen{}obstub\textquotesingle{}}} or \sphinxcode{\sphinxupquote{\textquotesingle{}ampl \sphinxhyphen{}ogstub\textquotesingle{}}}). After solving the problem,
\sphinxcode{\sphinxupquote{coptampl}} writes a \sphinxcode{\sphinxupquote{stub.sol}} file for use by AMPL’s solve and solution
commands. When you run AMPL, this all happens automatically if you
give the AMPL commands:

\begin{sphinxVerbatim}[commandchars=\\\{\}]
ampl:\PYG{+w}{ }option\PYG{+w}{ }solver\PYG{+w}{ }coptampl\PYG{p}{;}
ampl:\PYG{+w}{ }solve\PYG{p}{;}
\end{sphinxVerbatim}

\subsection{Installation Guide}
\label{\detokenize{amplinterface:installation-guide}}
\sphinxAtStartPar
To use \sphinxcode{\sphinxupquote{coptampl}} in AMPL, you must have a valid AMPL license and make sure
that you have installed Cardinal Optimizer and setup its license properly,
see {\hyperref[\detokenize{install:chapinstall}]{\sphinxcrossref{\DUrole{std,std-ref}{How to install Cardinal Optimizer}}}} for details.
Be sure to check if it satisfies the following requirements for different
operating systems.

\subsubsection{Windows}
\label{\detokenize{amplinterface:windows}}
\sphinxAtStartPar
On Windows platform, the \sphinxcode{\sphinxupquote{coptampl.exe}} utility and the \sphinxcode{\sphinxupquote{copt.dll}} dynamic
library contained in the Cardinal Optimizer must appear somewhere in your
user or system \sphinxcode{\sphinxupquote{PATH}} environment variable (or in the current directory).

\sphinxAtStartPar
To test if your setting meets the above requirements, you can check it by
executing commands below in command prompt:

\begin{sphinxVerbatim}[commandchars=\\\{\}]
coptampl\PYG{+w}{ }\PYGZhy{}v
\end{sphinxVerbatim}

\sphinxAtStartPar
And you are expected to see output similar to the following on screen:
\begin{sphinxalltt}
AMPL/x\sphinxhyphen{}COPT Optimizer {[}8.0.1{]} (windows\sphinxhyphen{}x86), driver(20220526), MP(20220526)
\end{sphinxalltt}

\sphinxAtStartPar
If the commands failed, then you should recheck your settings.

\subsubsection{Linux}
\label{\detokenize{amplinterface:linux}}
\sphinxAtStartPar
On Linux platform, the \sphinxcode{\sphinxupquote{coptampl}} utility must appears somewhere in your
\sphinxcode{\sphinxupquote{\$PATH}} environment variable, while the \sphinxcode{\sphinxupquote{libcopt.so}} shared library must
appears somewhere in your \sphinxcode{\sphinxupquote{\$LD\_LIBRARY\_PATH}} environment variable.

\sphinxAtStartPar
Similarly, to test if your setting meets the above requirements, just execute
commands below in shell:

\begin{sphinxVerbatim}[commandchars=\\\{\}]
coptampl\PYG{+w}{ }\PYGZhy{}v
\end{sphinxVerbatim}

\sphinxAtStartPar
And you are expected to see output similar to the following on screen:
\begin{sphinxalltt}
AMPL/x\sphinxhyphen{}COPT Optimizer {[}8.0.1{]} (linux\sphinxhyphen{}x86), driver(20220526), MP(20220526)
\end{sphinxalltt}

\sphinxAtStartPar
If the commands failed, please recheck your settings.

\subsubsection{MacOS}
\label{\detokenize{amplinterface:macos}}
\sphinxAtStartPar
On MacOS platform, the \sphinxcode{\sphinxupquote{coptampl}} utility must appears somewhere in your
\sphinxcode{\sphinxupquote{\$PATH}} environment variable, while the \sphinxcode{\sphinxupquote{libcopt.dylib}} dynamic library
must appears somewhere in your \sphinxcode{\sphinxupquote{\$DYLD\_LIBRARY\_PATH}} environment variable.

\sphinxAtStartPar
You can execute commands below in shell to see if your settings meets the
above requirements:

\begin{sphinxVerbatim}[commandchars=\\\{\}]
coptampl\PYG{+w}{ }\PYGZhy{}v
\end{sphinxVerbatim}

\sphinxAtStartPar
And you are expected to see output similar to the following on screen:
\begin{sphinxalltt}
AMPL/x\sphinxhyphen{}COPT Optimizer {[}8.0.1{]} (macos\sphinxhyphen{}x86), driver(20220526), MP(20220526)
\end{sphinxalltt}

\sphinxAtStartPar
If the commands failed, then please recheck your settings.

\subsection{Solver Options and Exit Codes}
\label{\detokenize{amplinterface:solver-options-and-exit-codes}}
\sphinxAtStartPar
The \sphinxcode{\sphinxupquote{coptampl}} utility offers some options to customize its behavior. Uers can
control it by setting the environment variable \sphinxcode{\sphinxupquote{copt\_options}} or use AMPL’s
\sphinxcode{\sphinxupquote{option}} command. To see all available options, please invoke:

\begin{sphinxVerbatim}[commandchars=\\\{\}]
coptampl\PYG{+w}{ }\PYGZhy{}\PYG{o}{=}
\end{sphinxVerbatim}

\sphinxAtStartPar
The supported parameters and their interpretation for current version are shown
in \hyperref[\detokenize{amplinterface:copttab-amplparam}]{Table \ref{\detokenize{amplinterface:copttab-amplparam}}}:

\begin{savenotes}
\sphinxatlongtablestart
\sphinxthistablewithglobalstyle
\begin{longtable}[c]{|l|l|}
\sphinxthelongtablecaptionisattop
\caption{Parameters of \sphinxstyleliteralintitle{\sphinxupquote{coptampl}}\strut}\label{\detokenize{amplinterface:copttab-amplparam}}\\*[\sphinxlongtablecapskipadjust]
\sphinxtoprule
\endfirsthead

\multicolumn{2}{c}{\sphinxnorowcolor
    \makebox[0pt]{\sphinxtablecontinued{\tablename\ \thetable{} \textendash{} continued from previous page}}%
}\\
\sphinxtoprule
\endhead

\sphinxbottomrule
\multicolumn{2}{r}{\sphinxnorowcolor
    \makebox[0pt][r]{\sphinxtablecontinued{continues on next page}}%
}\\
\endfoot

\endlastfoot
\sphinxtableatstartofbodyhook

\sphinxAtStartPar
\sphinxstylestrong{Parameter}
&
\sphinxAtStartPar
\sphinxstylestrong{Interpretation}
\\
\sphinxhline
\sphinxAtStartPar
barhomogeneous
&
\sphinxAtStartPar
whether to use homogeneous self\sphinxhyphen{}dual form in barrier
\\
\sphinxhline
\sphinxAtStartPar
bariterlimit
&
\sphinxAtStartPar
iteration limit of barrier method
\\
\sphinxhline
\sphinxAtStartPar
barthreads
&
\sphinxAtStartPar
number of threads used by barrier
\\
\sphinxhline
\sphinxAtStartPar
basis
&
\sphinxAtStartPar
whether to use or return basis status
\\
\sphinxhline
\sphinxAtStartPar
bestbound
&
\sphinxAtStartPar
whether to return best bound by suffix
\\
\sphinxhline
\sphinxAtStartPar
conflictanalysis
&
\sphinxAtStartPar
whether to perform conflict analysis
\\
\sphinxhline
\sphinxAtStartPar
crossoverthreads
&
\sphinxAtStartPar
number of threads used by crossover
\\
\sphinxhline
\sphinxAtStartPar
cutlevel
&
\sphinxAtStartPar
level of cutting\sphinxhyphen{}planes generation
\\
\sphinxhline
\sphinxAtStartPar
divingheurlevel
&
\sphinxAtStartPar
level of diving heuristics
\\
\sphinxhline
\sphinxAtStartPar
dualize
&
\sphinxAtStartPar
whether to dualize a problem before solving it
\\
\sphinxhline
\sphinxAtStartPar
dualperturb
&
\sphinxAtStartPar
whether to allow the objective function perturbation
\\
\sphinxhline
\sphinxAtStartPar
dualprice
&
\sphinxAtStartPar
specifies the dual simplex pricing algorithm
\\
\sphinxhline
\sphinxAtStartPar
dualtol
&
\sphinxAtStartPar
the tolerance for dual solutions and reduced cost
\\
\sphinxhline
\sphinxAtStartPar
feastol
&
\sphinxAtStartPar
the feasibility tolerance
\\
\sphinxhline
\sphinxAtStartPar
heurlevel
&
\sphinxAtStartPar
level of heuristics
\\
\sphinxhline
\sphinxAtStartPar
iisfind
&
\sphinxAtStartPar
whether to compute IIS and return result
\\
\sphinxhline
\sphinxAtStartPar
iismethod
&
\sphinxAtStartPar
specify the IIS method
\\
\sphinxhline
\sphinxAtStartPar
inttol
&
\sphinxAtStartPar
the integrality tolerance for variables
\\
\sphinxhline
\sphinxAtStartPar
logging
&
\sphinxAtStartPar
whether to print solving logs
\\
\sphinxhline
\sphinxAtStartPar
logfile
&
\sphinxAtStartPar
name of log file
\\
\sphinxhline
\sphinxAtStartPar
exportfile
&
\sphinxAtStartPar
name of model file to be exported
\\
\sphinxhline
\sphinxAtStartPar
lpmethod
&
\sphinxAtStartPar
method to solve the LP problem
\\
\sphinxhline
\sphinxAtStartPar
matrixtol
&
\sphinxAtStartPar
input matrix coefficient tolerance
\\
\sphinxhline
\sphinxAtStartPar
mipstart
&
\sphinxAtStartPar
whether to use initial values for MIP problem
\\
\sphinxhline
\sphinxAtStartPar
miptasks
&
\sphinxAtStartPar
number of MIP tasks in parallel
\\
\sphinxhline
\sphinxAtStartPar
nodecutrounds
&
\sphinxAtStartPar
rounds of cutting\sphinxhyphen{}planes generation of tree node
\\
\sphinxhline
\sphinxAtStartPar
nodelimit
&
\sphinxAtStartPar
node limit of the optimization
\\
\sphinxhline
\sphinxAtStartPar
objno
&
\sphinxAtStartPar
objective number to solve
\\
\sphinxhline
\sphinxAtStartPar
count
&
\sphinxAtStartPar
whether to count the number of solutions
\\
\sphinxhline
\sphinxAtStartPar
stub
&
\sphinxAtStartPar
name prefix for alternative MIP solutions written
\\
\sphinxhline
\sphinxAtStartPar
presolve
&
\sphinxAtStartPar
level of presolving before solving a problem
\\
\sphinxhline
\sphinxAtStartPar
relgap
&
\sphinxAtStartPar
the relative gap of optimization
\\
\sphinxhline
\sphinxAtStartPar
absgap
&
\sphinxAtStartPar
the absolute gap of optimization
\\
\sphinxhline
\sphinxAtStartPar
return\_mipgap
&
\sphinxAtStartPar
whether to return absolute/relative gap by suffix
\\
\sphinxhline
\sphinxAtStartPar
rootcutlevel
&
\sphinxAtStartPar
level of cutting\sphinxhyphen{}planes generation of root node
\\
\sphinxhline
\sphinxAtStartPar
rootcutrounds
&
\sphinxAtStartPar
rounds of cutting\sphinxhyphen{}planes generation of root node
\\
\sphinxhline
\sphinxAtStartPar
roundingheurlevel
&
\sphinxAtStartPar
level of rounding heuristics
\\
\sphinxhline
\sphinxAtStartPar
scaling
&
\sphinxAtStartPar
whether to perform scaling before solving a problem
\\
\sphinxhline
\sphinxAtStartPar
simplexthreads
&
\sphinxAtStartPar
number of threads used by dual simplex
\\
\sphinxhline
\sphinxAtStartPar
sos
&
\sphinxAtStartPar
whether to use ‘.sosno’ and ‘.ref’ suffix
\\
\sphinxhline
\sphinxAtStartPar
sos2
&
\sphinxAtStartPar
whether to use SOS2 to represent piecewise linear terms
\\
\sphinxhline
\sphinxAtStartPar
strongbranching
&
\sphinxAtStartPar
level of strong branching
\\
\sphinxhline
\sphinxAtStartPar
submipheurlevel
&
\sphinxAtStartPar
level of Sub\sphinxhyphen{}MIP heuristics
\\
\sphinxhline
\sphinxAtStartPar
threads
&
\sphinxAtStartPar
number of threads to use
\\
\sphinxhline
\sphinxAtStartPar
timelimit
&
\sphinxAtStartPar
time limit of the optimization
\\
\sphinxhline
\sphinxAtStartPar
treecutlevel
&
\sphinxAtStartPar
level of cutting\sphinxhyphen{}planes generation of search tree
\\
\sphinxhline
\sphinxAtStartPar
wantsol
&
\sphinxAtStartPar
whether to generate ‘.sol’ file
\\
\sphinxbottomrule
\end{longtable}
\sphinxtableafterendhook
\sphinxatlongtableend
\end{savenotes}

\sphinxAtStartPar
Please refer to {\hyperref[\detokenize{capiref:chapapi-param}]{\sphinxcrossref{\DUrole{std,std-ref}{COPT Parameters}}}} for details.

\sphinxAtStartPar
AMPL uses suffix to store or pass model and solution information, and also some
extension features, such as support for SOS constraints. Currently, \sphinxcode{\sphinxupquote{coptampl}}
support suffix information as shown in {\hyperref[\detokenize{amplinterface:copttab-amplsuffix}]{\sphinxcrossref{\DUrole{std,std-ref}{Suffix supported by coptampl}}}} :

\begin{savenotes}\sphinxattablestart
\sphinxthistablewithglobalstyle
\centering
\sphinxcapstartof{table}
\sphinxthecaptionisattop
\sphinxcaption{Suffix supported by \sphinxstyleliteralintitle{\sphinxupquote{coptampl}}}\label{\detokenize{amplinterface:copttab-amplsuffix}}
\sphinxaftertopcaption
\begin{tabulary}{\linewidth}[t]{|T|T|}
\sphinxtoprule
\sphinxtableatstartofbodyhook
\sphinxAtStartPar
\sphinxstylestrong{Suffix}
&
\sphinxAtStartPar
\sphinxstylestrong{Interpretation}
\\
\sphinxhline
\sphinxAtStartPar
absmipgap
&
\sphinxAtStartPar
absolute gap for MIP problem
\\
\sphinxhline
\sphinxAtStartPar
bestbound
&
\sphinxAtStartPar
best bound for MIP problem
\\
\sphinxhline
\sphinxAtStartPar
iis
&
\sphinxAtStartPar
store IIS status of variables or constraints
\\
\sphinxhline
\sphinxAtStartPar
nsol
&
\sphinxAtStartPar
number of pool solutions written
\\
\sphinxhline
\sphinxAtStartPar
ref
&
\sphinxAtStartPar
weight of variable in SOS constraint
\\
\sphinxhline
\sphinxAtStartPar
relmipgap
&
\sphinxAtStartPar
relative gap for MIP problem
\\
\sphinxhline
\sphinxAtStartPar
sos
&
\sphinxAtStartPar
store type of SOS constraint
\\
\sphinxhline
\sphinxAtStartPar
sosno
&
\sphinxAtStartPar
type of SOS constraint
\\
\sphinxhline
\sphinxAtStartPar
sosref
&
\sphinxAtStartPar
store variable weight in SOS constraint
\\
\sphinxhline
\sphinxAtStartPar
sstatus
&
\sphinxAtStartPar
basis status of variables and constraints
\\
\sphinxbottomrule
\end{tabulary}
\sphinxtableafterendhook\par
\sphinxattableend\end{savenotes}

\sphinxAtStartPar
Users who want to know how to use SOS constraints in AMPL, please refer to
resources in AMPL’s website:
\sphinxhref{https://ampl.com/faqs/how-can-i-use-the-solvers-special-ordered-sets-feature/}{How to use SOS constraints in AMPL} .

\sphinxAtStartPar
When solving finished, \sphinxcode{\sphinxupquote{coptampl}} will display a status message and return
exit code to AMPL. The exit code can be displayed by:

\begin{sphinxVerbatim}[commandchars=\\\{\}]
ampl:\PYG{+w}{ }display\PYG{+w}{ }solve\PYGZus{}result\PYGZus{}num\PYG{p}{;}
\end{sphinxVerbatim}

\sphinxAtStartPar
If no solution was found or something unexpected happened, \sphinxcode{\sphinxupquote{coptampl}} will
return non\sphinxhyphen{}zero code to AMPL from \hyperref[\detokenize{amplinterface:copttab-amplexitcode}]{Table \ref{\detokenize{amplinterface:copttab-amplexitcode}}}:

\begin{savenotes}\sphinxattablestart
\sphinxthistablewithglobalstyle
\centering
\sphinxcapstartof{table}
\sphinxthecaptionisattop
\sphinxcaption{Exit codes of \sphinxstyleliteralintitle{\sphinxupquote{coptampl}}}\label{\detokenize{amplinterface:copttab-amplexitcode}}
\sphinxaftertopcaption
\begin{tabulary}{\linewidth}[t]{|T|T|}
\sphinxtoprule
\sphinxtableatstartofbodyhook
\sphinxAtStartPar
\sphinxstylestrong{Exit Code}
&
\sphinxAtStartPar
\sphinxstylestrong{Interpretation}
\\
\sphinxhline
\sphinxAtStartPar
0
&
\sphinxAtStartPar
optimal solution
\\
\sphinxhline
\sphinxAtStartPar
200
&
\sphinxAtStartPar
infeasible
\\
\sphinxhline
\sphinxAtStartPar
300
&
\sphinxAtStartPar
unbounded
\\
\sphinxhline
\sphinxAtStartPar
301
&
\sphinxAtStartPar
infeasible or unbounded
\\
\sphinxhline
\sphinxAtStartPar
600
&
\sphinxAtStartPar
user interrupted
\\
\sphinxbottomrule
\end{tabulary}
\sphinxtableafterendhook\par
\sphinxattableend\end{savenotes}

\subsection{Example Usage}
\label{\detokenize{amplinterface:example-usage}}\label{\detokenize{amplinterface:examplediet}}
\sphinxAtStartPar
The following section will illustrate the use of AMPL by a well\sphinxhyphen{}known example
called “Diet problem”, which finds a mix of foods that satisfies requirements
on the amounts of various vitamins, see
\sphinxhref{https://ampl.com/BOOK/CHAPTERS/05-tut2.pdf}{AMPL book} for details.

\sphinxAtStartPar
Suppose the following kinds of foods are available for the following prices
per unit, see \hyperref[\detokenize{amplinterface:copttab-ampldiet-1}]{Table \ref{\detokenize{amplinterface:copttab-ampldiet-1}}}:

\begin{savenotes}\sphinxattablestart
\sphinxthistablewithglobalstyle
\centering
\sphinxcapstartof{table}
\sphinxthecaptionisattop
\sphinxcaption{Prices of foods}\label{\detokenize{amplinterface:copttab-ampldiet-1}}
\sphinxaftertopcaption
\begin{tabulary}{\linewidth}[t]{|T|T|}
\sphinxtoprule
\sphinxtableatstartofbodyhook
\sphinxAtStartPar
\sphinxstylestrong{Food}
&
\sphinxAtStartPar
\sphinxstylestrong{Price}
\\
\sphinxhline
\sphinxAtStartPar
BEEF
&
\sphinxAtStartPar
3.19
\\
\sphinxhline
\sphinxAtStartPar
CHK
&
\sphinxAtStartPar
2.59
\\
\sphinxhline
\sphinxAtStartPar
FISH
&
\sphinxAtStartPar
2.29
\\
\sphinxhline
\sphinxAtStartPar
HAM
&
\sphinxAtStartPar
2.89
\\
\sphinxhline
\sphinxAtStartPar
MCH
&
\sphinxAtStartPar
1.89
\\
\sphinxhline
\sphinxAtStartPar
MTL
&
\sphinxAtStartPar
1.99
\\
\sphinxhline
\sphinxAtStartPar
SPG
&
\sphinxAtStartPar
1.99
\\
\sphinxhline
\sphinxAtStartPar
TUR
&
\sphinxAtStartPar
2.49
\\
\sphinxbottomrule
\end{tabulary}
\sphinxtableafterendhook\par
\sphinxattableend\end{savenotes}

\sphinxAtStartPar
These foods provide the following percentages, per unit, of the minimum daily
requirements for vitamins A, C, B1 and B2, see \hyperref[\detokenize{amplinterface:copttab-ampldiet-2}]{Table \ref{\detokenize{amplinterface:copttab-ampldiet-2}}}:

\begin{savenotes}\sphinxattablestart
\sphinxthistablewithglobalstyle
\centering
\sphinxcapstartof{table}
\sphinxthecaptionisattop
\sphinxcaption{Nutrition of foods (\%)}\label{\detokenize{amplinterface:copttab-ampldiet-2}}
\sphinxaftertopcaption
\begin{tabulary}{\linewidth}[t]{|T|T|T|T|T|}
\sphinxtoprule
\sphinxtableatstartofbodyhook&
\sphinxAtStartPar
A
&
\sphinxAtStartPar
C
&
\sphinxAtStartPar
B1
&
\sphinxAtStartPar
B2
\\
\sphinxhline
\sphinxAtStartPar
BEEF
&
\sphinxAtStartPar
60\%
&
\sphinxAtStartPar
20\%
&
\sphinxAtStartPar
10\%
&
\sphinxAtStartPar
15\%
\\
\sphinxhline
\sphinxAtStartPar
CHK
&
\sphinxAtStartPar
8
&
\sphinxAtStartPar
0
&
\sphinxAtStartPar
20
&
\sphinxAtStartPar
20
\\
\sphinxhline
\sphinxAtStartPar
FISH
&
\sphinxAtStartPar
8
&
\sphinxAtStartPar
10
&
\sphinxAtStartPar
15
&
\sphinxAtStartPar
10
\\
\sphinxhline
\sphinxAtStartPar
HAM
&
\sphinxAtStartPar
40
&
\sphinxAtStartPar
40
&
\sphinxAtStartPar
35
&
\sphinxAtStartPar
10
\\
\sphinxhline
\sphinxAtStartPar
MCH
&
\sphinxAtStartPar
15
&
\sphinxAtStartPar
35
&
\sphinxAtStartPar
15
&
\sphinxAtStartPar
15
\\
\sphinxhline
\sphinxAtStartPar
MTL
&
\sphinxAtStartPar
70
&
\sphinxAtStartPar
30
&
\sphinxAtStartPar
15
&
\sphinxAtStartPar
15
\\
\sphinxhline
\sphinxAtStartPar
SPG
&
\sphinxAtStartPar
25
&
\sphinxAtStartPar
50
&
\sphinxAtStartPar
25
&
\sphinxAtStartPar
15
\\
\sphinxhline
\sphinxAtStartPar
TUR
&
\sphinxAtStartPar
60
&
\sphinxAtStartPar
20
&
\sphinxAtStartPar
15
&
\sphinxAtStartPar
10
\\
\sphinxbottomrule
\end{tabulary}
\sphinxtableafterendhook\par
\sphinxattableend\end{savenotes}

\sphinxAtStartPar
The problem is to find the cheapest combination that meets a week’s
requirements, that is, at least 700\% of the daily requirements for
each nutrient.

\sphinxAtStartPar
To summarize, the mathematical form for the above problem can be modeled
as shown in Eq. \ref{equation:amplinterface:coptEq_ampldiet_math}:
\begin{equation}\label{equation:amplinterface:coptEq_ampldiet_math}
\begin{split}\textrm{Minimize: } & \\
& \sum_{j \in J} cost_j \cdot buy_j \\
\textrm{Subject to: } & \\
& n\_min_i \leq \sum_{j \in J} amt_{i, j} \cdot buy_j \leq n\_max_i \,\,\, \forall i \in I \\
& f\_min_j \leq buy_j \leq f\_max_j \,\,\, \forall j \in J\end{split}
\end{equation}
\sphinxAtStartPar
The AMPL model for above problem is shown in \sphinxcode{\sphinxupquote{diet.mod}},
see \hyperref[\detokenize{amplinterface:coptcode-ampldiet-model}]{Listing \ref{\detokenize{amplinterface:coptcode-ampldiet-model}}}:
\sphinxSetupCaptionForVerbatim{\sphinxcode{\sphinxupquote{diet.mod}}}
\def\sphinxLiteralBlockLabel{\label{\detokenize{amplinterface:coptcode-ampldiet-model}}}
\begin{sphinxVerbatim}[commandchars=\\\{\},numbers=left,firstnumber=1,stepnumber=1]
\PYG{c+c1}{\PYGZsh{} The code is adopted from:}
\PYG{c+c1}{\PYGZsh{}}
\PYG{c+c1}{\PYGZsh{} https://github.com/Pyomo/pyomo/blob/master/examples/pyomo/amplbook2/diet.mod}
\PYG{c+c1}{\PYGZsh{}}
\PYG{c+c1}{\PYGZsh{} with some modification by developer of the Cardinal Optimizer}

\PYG{k+kd}{set}\PYG{+w}{ }\PYG{n+nv}{NUTR}\PYG{p}{;}
\PYG{k+kd}{set}\PYG{+w}{ }\PYG{n+nv}{FOOD}\PYG{p}{;}

\PYG{k+kd}{param}\PYG{+w}{ }\PYG{n+nv}{cost}\PYG{+w}{ }\PYG{p}{\PYGZob{}}FOOD\PYG{p}{\PYGZcb{}}\PYG{+w}{ }\PYG{o}{\PYGZgt{}}\PYG{+w}{ }\PYG{l+m+mi}{0}\PYG{p}{;}
\PYG{k+kd}{param}\PYG{+w}{ }\PYG{n+nv}{f\PYGZus{}min}\PYG{+w}{ }\PYG{p}{\PYGZob{}}FOOD\PYG{p}{\PYGZcb{}}\PYG{+w}{ }\PYG{o}{\PYGZgt{}=}\PYG{+w}{ }\PYG{l+m+mi}{0}\PYG{p}{;}
\PYG{k+kd}{param}\PYG{+w}{ }\PYG{n+nv}{f\PYGZus{}max}\PYG{+w}{ }\PYG{p}{\PYGZob{}}j\PYG{+w}{ }\PYG{k+kr}{in}\PYG{+w}{ }FOOD\PYG{p}{\PYGZcb{}}\PYG{+w}{ }\PYG{o}{\PYGZgt{}=}\PYG{+w}{ }f\PYGZus{}min\PYG{p}{[}j\PYG{p}{];}

\PYG{k+kd}{param}\PYG{+w}{ }\PYG{n+nv}{n\PYGZus{}min}\PYG{+w}{ }\PYG{p}{\PYGZob{}}NUTR\PYG{p}{\PYGZcb{}}\PYG{+w}{ }\PYG{o}{\PYGZgt{}=}\PYG{+w}{ }\PYG{l+m+mi}{0}\PYG{p}{;}
\PYG{k+kd}{param}\PYG{+w}{ }\PYG{n+nv}{n\PYGZus{}max}\PYG{+w}{ }\PYG{p}{\PYGZob{}}i\PYG{+w}{ }\PYG{k+kr}{in}\PYG{+w}{ }NUTR\PYG{p}{\PYGZcb{}}\PYG{+w}{ }\PYG{o}{\PYGZgt{}=}\PYG{+w}{ }n\PYGZus{}min\PYG{p}{[}i\PYG{p}{];}

\PYG{k+kd}{param}\PYG{+w}{ }\PYG{n+nv}{amt}\PYG{+w}{ }\PYG{p}{\PYGZob{}}NUTR\PYG{p}{,}\PYG{+w}{ }FOOD\PYG{p}{\PYGZcb{}}\PYG{+w}{ }\PYG{o}{\PYGZgt{}=}\PYG{+w}{ }\PYG{l+m+mi}{0}\PYG{p}{;}

\PYG{k+kd}{var}\PYG{+w}{ }\PYG{n+nv}{Buy}\PYG{+w}{ }\PYG{p}{\PYGZob{}}j\PYG{+w}{ }\PYG{k+kr}{in}\PYG{+w}{ }FOOD\PYG{p}{\PYGZcb{}}\PYG{+w}{ }\PYG{o}{\PYGZgt{}=}\PYG{+w}{ }f\PYGZus{}min\PYG{p}{[}j\PYG{p}{],}\PYG{+w}{ }\PYG{o}{\PYGZlt{}=}\PYG{+w}{ }f\PYGZus{}max\PYG{p}{[}j\PYG{p}{];}

\PYG{k+kd}{minimize}\PYG{+w}{ }\PYG{n+nv}{Total\PYGZus{}Cost}\PYG{p}{:}
\PYG{+w}{  }\PYG{k+kr}{sum}\PYG{+w}{ }\PYG{p}{\PYGZob{}}j\PYG{+w}{ }\PYG{k+kr}{in}\PYG{+w}{ }FOOD\PYG{p}{\PYGZcb{}}\PYG{+w}{ }cost\PYG{p}{[}j\PYG{p}{]}\PYG{+w}{ }\PYG{o}{*}\PYG{+w}{ }Buy\PYG{p}{[}j\PYG{p}{];}

\PYG{k+kd}{subject to}\PYG{+w}{ }\PYG{n+nv}{Diet}\PYG{+w}{ }\PYG{p}{\PYGZob{}}i\PYG{+w}{ }\PYG{k+kr}{in}\PYG{+w}{ }NUTR\PYG{p}{\PYGZcb{}:}
\PYG{+w}{  }n\PYGZus{}min\PYG{p}{[}i\PYG{p}{]}\PYG{+w}{ }\PYG{o}{\PYGZlt{}=}\PYG{+w}{ }\PYG{k+kr}{sum}\PYG{+w}{ }\PYG{p}{\PYGZob{}}j\PYG{+w}{ }\PYG{k+kr}{in}\PYG{+w}{ }FOOD\PYG{p}{\PYGZcb{}}\PYG{+w}{ }amt\PYG{p}{[}i\PYG{p}{,}\PYG{+w}{ }j\PYG{p}{]}\PYG{+w}{ }\PYG{o}{*}\PYG{+w}{ }Buy\PYG{p}{[}j\PYG{p}{]}\PYG{+w}{ }\PYG{o}{\PYGZlt{}=}\PYG{+w}{ }n\PYGZus{}max\PYG{p}{[}i\PYG{p}{];}
\end{sphinxVerbatim}

\sphinxAtStartPar
The data file for above problem is shown in \sphinxcode{\sphinxupquote{diet.dat}},
see \hyperref[\detokenize{amplinterface:coptcode-ampldiet-data}]{Listing \ref{\detokenize{amplinterface:coptcode-ampldiet-data}}}:
\sphinxSetupCaptionForVerbatim{\sphinxcode{\sphinxupquote{diet.dat}}}
\def\sphinxLiteralBlockLabel{\label{\detokenize{amplinterface:coptcode-ampldiet-data}}}
\begin{sphinxVerbatim}[commandchars=\\\{\},numbers=left,firstnumber=1,stepnumber=1]
\PYG{c+c1}{\PYGZsh{} The data is adopted from:}
\PYG{c+c1}{\PYGZsh{} }
\PYG{c+c1}{\PYGZsh{} https://github.com/Pyomo/pyomo/blob/master/examples/pyomo/amplbook2/diet.dat}
\PYG{c+c1}{\PYGZsh{}}
\PYG{c+c1}{\PYGZsh{} with some modification by developer of the Cardinal Optimizer}

\PYG{k+kr}{data}\PYG{p}{;}

\PYG{k+kd}{set}\PYG{+w}{ }\PYG{n+nv}{NUTR}\PYG{+w}{ }\PYG{p}{:}\PYG{o}{=}\PYG{+w}{ }A\PYG{+w}{ }B1\PYG{+w}{ }B2\PYG{+w}{ }C\PYG{+w}{ }\PYG{p}{;}
\PYG{k+kd}{set}\PYG{+w}{ }\PYG{n+nv}{FOOD}\PYG{+w}{ }\PYG{p}{:}\PYG{o}{=}\PYG{+w}{ }BEEF\PYG{+w}{ }CHK\PYG{+w}{ }FISH\PYG{+w}{ }HAM\PYG{+w}{ }MCH\PYG{+w}{ }MTL\PYG{+w}{ }SPG\PYG{+w}{ }TUR\PYG{+w}{ }\PYG{p}{;}

param\PYG{p}{:}\PYG{+w}{   }cost\PYG{+w}{  }f\PYGZus{}min\PYG{+w}{  }f\PYGZus{}max\PYG{+w}{ }\PYG{p}{:}\PYG{o}{=}
\PYG{+w}{  }BEEF\PYG{+w}{   }\PYG{l+m+mf}{3.19}\PYG{+w}{    }\PYG{l+m+mi}{0}\PYG{+w}{     }\PYG{l+m+mi}{100}
\PYG{+w}{  }CHK\PYG{+w}{    }\PYG{l+m+mf}{2.59}\PYG{+w}{    }\PYG{l+m+mi}{0}\PYG{+w}{     }\PYG{l+m+mi}{100}
\PYG{+w}{  }FISH\PYG{+w}{   }\PYG{l+m+mf}{2.29}\PYG{+w}{    }\PYG{l+m+mi}{0}\PYG{+w}{     }\PYG{l+m+mi}{100}
\PYG{+w}{  }HAM\PYG{+w}{    }\PYG{l+m+mf}{2.89}\PYG{+w}{    }\PYG{l+m+mi}{0}\PYG{+w}{     }\PYG{l+m+mi}{100}
\PYG{+w}{  }MCH\PYG{+w}{    }\PYG{l+m+mf}{1.89}\PYG{+w}{    }\PYG{l+m+mi}{0}\PYG{+w}{     }\PYG{l+m+mi}{100}
\PYG{+w}{  }MTL\PYG{+w}{    }\PYG{l+m+mf}{1.99}\PYG{+w}{    }\PYG{l+m+mi}{0}\PYG{+w}{     }\PYG{l+m+mi}{100}
\PYG{+w}{  }SPG\PYG{+w}{    }\PYG{l+m+mf}{1.99}\PYG{+w}{    }\PYG{l+m+mi}{0}\PYG{+w}{     }\PYG{l+m+mi}{100}
\PYG{+w}{  }TUR\PYG{+w}{    }\PYG{l+m+mf}{2.49}\PYG{+w}{    }\PYG{l+m+mi}{0}\PYG{+w}{     }\PYG{l+m+mi}{100}\PYG{+w}{ }\PYG{p}{;}

param\PYG{p}{:}\PYG{+w}{   }n\PYGZus{}min\PYG{+w}{  }n\PYGZus{}max\PYG{+w}{ }\PYG{p}{:}\PYG{o}{=}
\PYG{+w}{   }A\PYG{+w}{      }\PYG{l+m+mi}{700}\PYG{+w}{   }\PYG{l+m+mi}{10000}
\PYG{+w}{   }C\PYG{+w}{      }\PYG{l+m+mi}{700}\PYG{+w}{   }\PYG{l+m+mi}{10000}
\PYG{+w}{   }B1\PYG{+w}{     }\PYG{l+m+mi}{700}\PYG{+w}{   }\PYG{l+m+mi}{10000}
\PYG{+w}{   }B2\PYG{+w}{     }\PYG{l+m+mi}{700}\PYG{+w}{   }\PYG{l+m+mi}{10000}\PYG{+w}{ }\PYG{p}{;}

\PYG{k+kd}{param}\PYG{+w}{ }\PYG{n+nv}{amt}\PYG{+w}{ }\PYG{p}{(}tr\PYG{p}{):}
\PYG{+w}{           }A\PYG{+w}{    }C\PYG{+w}{   }B1\PYG{+w}{   }B2\PYG{+w}{ }\PYG{p}{:}\PYG{o}{=}
\PYG{+w}{   }BEEF\PYG{+w}{   }\PYG{l+m+mi}{60}\PYG{+w}{   }\PYG{l+m+mi}{20}\PYG{+w}{   }\PYG{l+m+mi}{10}\PYG{+w}{   }\PYG{l+m+mi}{15}
\PYG{+w}{   }CHK\PYG{+w}{     }\PYG{l+m+mi}{8}\PYG{+w}{    }\PYG{l+m+mi}{0}\PYG{+w}{   }\PYG{l+m+mi}{20}\PYG{+w}{   }\PYG{l+m+mi}{20}
\PYG{+w}{   }FISH\PYG{+w}{    }\PYG{l+m+mi}{8}\PYG{+w}{   }\PYG{l+m+mi}{10}\PYG{+w}{   }\PYG{l+m+mi}{15}\PYG{+w}{   }\PYG{l+m+mi}{10}
\PYG{+w}{   }HAM\PYG{+w}{    }\PYG{l+m+mi}{40}\PYG{+w}{   }\PYG{l+m+mi}{40}\PYG{+w}{   }\PYG{l+m+mi}{35}\PYG{+w}{   }\PYG{l+m+mi}{10}
\PYG{+w}{   }MCH\PYG{+w}{    }\PYG{l+m+mi}{15}\PYG{+w}{   }\PYG{l+m+mi}{35}\PYG{+w}{   }\PYG{l+m+mi}{15}\PYG{+w}{   }\PYG{l+m+mi}{15}
\PYG{+w}{   }MTL\PYG{+w}{    }\PYG{l+m+mi}{70}\PYG{+w}{   }\PYG{l+m+mi}{30}\PYG{+w}{   }\PYG{l+m+mi}{15}\PYG{+w}{   }\PYG{l+m+mi}{15}
\PYG{+w}{   }SPG\PYG{+w}{    }\PYG{l+m+mi}{25}\PYG{+w}{   }\PYG{l+m+mi}{50}\PYG{+w}{   }\PYG{l+m+mi}{25}\PYG{+w}{   }\PYG{l+m+mi}{15}
\PYG{+w}{   }TUR\PYG{+w}{    }\PYG{l+m+mi}{60}\PYG{+w}{   }\PYG{l+m+mi}{20}\PYG{+w}{   }\PYG{l+m+mi}{15}\PYG{+w}{   }\PYG{l+m+mi}{10}\PYG{+w}{ }\PYG{p}{;}
\end{sphinxVerbatim}

\sphinxAtStartPar
To solve the problem with \sphinxcode{\sphinxupquote{coptampl}} in AMPL, just type commands in command
prompt on Windows or shell on Linux and MacOS:

\begin{sphinxVerbatim}[commandchars=\\\{\}]
ampl\PYG{p}{:}\PYG{+w}{ }\PYG{k+kr}{model}\PYG{+w}{ }diet.\PYG{k+kr}{mod}\PYG{p}{;}
ampl\PYG{p}{:}\PYG{+w}{ }\PYG{k+kr}{data}\PYG{+w}{ }diet.dat\PYG{p}{;}
ampl\PYG{p}{:}\PYG{+w}{ }\PYG{k+kr}{option}\PYG{+w}{ }solver\PYG{+w}{ }coptampl\PYG{p}{;}
ampl\PYG{p}{:}\PYG{+w}{ }\PYG{k+kr}{option}\PYG{+w}{ }copt\PYGZus{}options\PYG{+w}{ }\PYG{l+s+s1}{\PYGZsq{}logging 1\PYGZsq{}}\PYG{p}{;}
ampl\PYG{p}{:}\PYG{+w}{ }\PYG{k+kr}{solve}\PYG{p}{;}
\end{sphinxVerbatim}

\sphinxAtStartPar
\sphinxcode{\sphinxupquote{coptampl}} solve it quickly and display solving log and status message on screen:

\begin{sphinxVerbatim}[commandchars=\\\{\}]
x\PYGZhy{}COPT 5.0.1: optimal solution; objective 88.2
1 simplex iterations
\end{sphinxVerbatim}

\sphinxAtStartPar
So \sphinxcode{\sphinxupquote{coptampl}} claimed it found the optimal solution, and the minimal cost is
88.2 units. You can further check the solution by:

\begin{sphinxVerbatim}[commandchars=\\\{\}]
ampl\PYG{p}{:}\PYG{+w}{ }\PYG{k+kr}{display}\PYG{+w}{ }Buy\PYG{p}{;}
\end{sphinxVerbatim}

\sphinxAtStartPar
And you will get:

\begin{sphinxVerbatim}[commandchars=\\\{\}]
Buy\PYG{+w}{ }\PYG{p}{[}\PYG{o}{*}\PYG{p}{]}\PYG{+w}{ }\PYG{p}{:}\PYG{o}{=}
BEEF\PYG{+w}{   }\PYG{l+m+mi}{0}
\PYG{+w}{ }CHK\PYG{+w}{   }\PYG{l+m+mi}{0}
FISH\PYG{+w}{   }\PYG{l+m+mi}{0}
\PYG{+w}{ }HAM\PYG{+w}{   }\PYG{l+m+mi}{0}
\PYG{+w}{ }MCH\PYG{+w}{  }\PYG{l+m+mf}{46.6667}
\PYG{+w}{ }MTL\PYG{+w}{   }\PYG{l+m+mi}{0}
\PYG{+w}{ }SPG\PYG{+w}{   }\PYG{l+m+mi}{0}
\PYG{+w}{ }TUR\PYG{+w}{   }\PYG{l+m+mi}{0}
\PYG{p}{;}
\end{sphinxVerbatim}

\sphinxAtStartPar
So if we buy 46.667 units of MCH, we will have a minimal cost of 88.2 units.

\sphinxstepscope

\section{Pyomo Interface}
\label{\detokenize{pyomointerface:pyomo-interface}}\label{\detokenize{pyomointerface:chappyomointerface}}\label{\detokenize{pyomointerface::doc}}
\sphinxAtStartPar
\sphinxhref{https://www.pyomo.org/about}{Pyomo} is a Python based, open source
optimization modeling language with a diverse set of optimization capabilities.
It is used by researchers to solve complex real\sphinxhyphen{}world applications,
see \sphinxhref{https://www.pyomo.org/impact}{Who uses Pyomo?} for more introduction.
The following documentation explains how to use the \sphinxstylestrong{Cardinal Optimizer}.

\subsection{Installation Guide}
\label{\detokenize{pyomointerface:installation-guide}}
\sphinxAtStartPar
To use the Cardinal Optimizer in Pyomo, you should setup Pyomo and
Cardinal Optimizer correctly first. Pyomo currently supports Python 2.7, 3.6\sphinxhyphen{}3.9,
you can install Python from
\sphinxhref{https://www.anaconda.com/distribution/}{Anaconda Distribution of Python} or
from \sphinxhref{https://www.python.org/}{Official Python}. We recommend install Python
from Anaconda since it is much more friendly and convenient for fresh users.

\subsubsection{Using conda}
\label{\detokenize{pyomointerface:using-conda}}
\sphinxAtStartPar
The recommended way to install Pyomo in Anaconda Distribution of Python is
to use \sphinxcode{\sphinxupquote{conda}} which is built\sphinxhyphen{}in supported. Just execute the following
commands in command prompt on Windows or shell on Linux and MacOS:

\begin{sphinxVerbatim}[commandchars=\\\{\}]
conda\PYG{+w}{ }install\PYG{+w}{ }\PYGZhy{}c\PYG{+w}{ }conda\PYGZhy{}forge\PYG{+w}{ }pyomo
\end{sphinxVerbatim}

\sphinxAtStartPar
Pyomo also has conditional dependencies on a variety of third\sphinxhyphen{}party Python
packages, they can be installed using \sphinxcode{\sphinxupquote{conda}} with commands:

\begin{sphinxVerbatim}[commandchars=\\\{\}]
conda\PYG{+w}{ }install\PYG{+w}{ }\PYGZhy{}c\PYG{+w}{ }conda\PYGZhy{}forge\PYG{+w}{ }pyomo.extras
\end{sphinxVerbatim}

\subsubsection{Using pip}
\label{\detokenize{pyomointerface:using-pip}}
\sphinxAtStartPar
The alternative way to install Pyomo is to use the standard \sphinxcode{\sphinxupquote{pip}} utitility,
just execute the following commands in command prompt on Windows or shell
on Linux and MacOS:

\begin{sphinxVerbatim}[commandchars=\\\{\}]
pip\PYG{+w}{ }install\PYG{+w}{ }pyomo
\end{sphinxVerbatim}

\sphinxAtStartPar
If you encounter any problems when installing Pyomo, please refer to
\sphinxhref{https://www.pyomo.org/installation}{How to install Pyomo} for details.
To install Cardinal Optimizer and setup its license properly, please refer to
{\hyperref[\detokenize{install:chapinstall}]{\sphinxcrossref{\DUrole{std,std-ref}{How to install Cardinal Optimizer}}}} for details.

\subsection{Example Usage}
\label{\detokenize{pyomointerface:example-usage}}
\sphinxAtStartPar
We are going to make a simple introduction on how to use the Cardinal Optimizer
in Pyomo by solving the example described in
{\hyperref[\detokenize{amplinterface:examplediet}]{\sphinxcrossref{\DUrole{std,std-ref}{AMPL Interface \sphinxhyphen{} Example Usage}}}}. Users who want to
learn more information about Pyomo may refer to
\sphinxhref{https://pyomo.readthedocs.io/en/stable/}{Pyomo documentation} for details.

\subsubsection{Abstract Model}
\label{\detokenize{pyomointerface:abstract-model}}
\sphinxAtStartPar
Pyomo provides two major approaches to construct any supported model types,
here we show the \sphinxstylestrong{Abstract Model} approach to solve the above problem.

\sphinxAtStartPar
The source code \sphinxcode{\sphinxupquote{pydiet\_abstract.py}} is shown below, see
\hyperref[\detokenize{pyomointerface:coptcode-pydietabs-code}]{Listing \ref{\detokenize{pyomointerface:coptcode-pydietabs-code}}}:
\sphinxSetupCaptionForVerbatim{\sphinxcode{\sphinxupquote{pydiet\_abstract.py}}}
\def\sphinxLiteralBlockLabel{\label{\detokenize{pyomointerface:coptcode-pydietabs-code}}}
\begin{sphinxVerbatim}[commandchars=\\\{\},numbers=left,firstnumber=1,stepnumber=1]
\PYG{c+c1}{\PYGZsh{} The code is adopted from:}
\PYG{c+c1}{\PYGZsh{}}
\PYG{c+c1}{\PYGZsh{} https://github.com/Pyomo/pyomo/blob/master/examples/pyomo/amplbook2/diet.py}
\PYG{c+c1}{\PYGZsh{}}
\PYG{c+c1}{\PYGZsh{} with some modification by developer of the Cardinal Optimizer}

\PYG{k+kn}{from}\PYG{+w}{ }\PYG{n+nn}{pyomo}\PYG{n+nn}{.}\PYG{n+nn}{core}\PYG{+w}{ }\PYG{k+kn}{import} \PYG{o}{*}

\PYG{n}{model} \PYG{o}{=} \PYG{n}{AbstractModel}\PYG{p}{(}\PYG{p}{)}

\PYG{n}{model}\PYG{o}{.}\PYG{n}{NUTR} \PYG{o}{=} \PYG{n}{Set}\PYG{p}{(}\PYG{p}{)}
\PYG{n}{model}\PYG{o}{.}\PYG{n}{FOOD} \PYG{o}{=} \PYG{n}{Set}\PYG{p}{(}\PYG{p}{)}

\PYG{n}{model}\PYG{o}{.}\PYG{n}{cost}  \PYG{o}{=} \PYG{n}{Param}\PYG{p}{(}\PYG{n}{model}\PYG{o}{.}\PYG{n}{FOOD}\PYG{p}{,} \PYG{n}{within}\PYG{o}{=}\PYG{n}{NonNegativeReals}\PYG{p}{)}
\PYG{n}{model}\PYG{o}{.}\PYG{n}{f\PYGZus{}min} \PYG{o}{=} \PYG{n}{Param}\PYG{p}{(}\PYG{n}{model}\PYG{o}{.}\PYG{n}{FOOD}\PYG{p}{,} \PYG{n}{within}\PYG{o}{=}\PYG{n}{NonNegativeReals}\PYG{p}{)}

\PYG{n}{model}\PYG{o}{.}\PYG{n}{f\PYGZus{}max} \PYG{o}{=} \PYG{n}{Param}\PYG{p}{(}\PYG{n}{model}\PYG{o}{.}\PYG{n}{FOOD}\PYG{p}{)}
\PYG{n}{model}\PYG{o}{.}\PYG{n}{n\PYGZus{}min} \PYG{o}{=} \PYG{n}{Param}\PYG{p}{(}\PYG{n}{model}\PYG{o}{.}\PYG{n}{NUTR}\PYG{p}{,} \PYG{n}{within}\PYG{o}{=}\PYG{n}{NonNegativeReals}\PYG{p}{)}
\PYG{n}{model}\PYG{o}{.}\PYG{n}{n\PYGZus{}max} \PYG{o}{=} \PYG{n}{Param}\PYG{p}{(}\PYG{n}{model}\PYG{o}{.}\PYG{n}{NUTR}\PYG{p}{)}
\PYG{n}{model}\PYG{o}{.}\PYG{n}{amt}   \PYG{o}{=} \PYG{n}{Param}\PYG{p}{(}\PYG{n}{model}\PYG{o}{.}\PYG{n}{NUTR}\PYG{p}{,} \PYG{n}{model}\PYG{o}{.}\PYG{n}{FOOD}\PYG{p}{,} \PYG{n}{within}\PYG{o}{=}\PYG{n}{NonNegativeReals}\PYG{p}{)}

\PYG{k}{def}\PYG{+w}{ }\PYG{n+nf}{Buy\PYGZus{}bounds}\PYG{p}{(}\PYG{n}{model}\PYG{p}{,} \PYG{n}{i}\PYG{p}{)}\PYG{p}{:}
  \PYG{k}{return} \PYG{p}{(}\PYG{n}{model}\PYG{o}{.}\PYG{n}{f\PYGZus{}min}\PYG{p}{[}\PYG{n}{i}\PYG{p}{]}\PYG{p}{,} \PYG{n}{model}\PYG{o}{.}\PYG{n}{f\PYGZus{}max}\PYG{p}{[}\PYG{n}{i}\PYG{p}{]}\PYG{p}{)}
\PYG{n}{model}\PYG{o}{.}\PYG{n}{Buy} \PYG{o}{=} \PYG{n}{Var}\PYG{p}{(}\PYG{n}{model}\PYG{o}{.}\PYG{n}{FOOD}\PYG{p}{,} \PYG{n}{bounds}\PYG{o}{=}\PYG{n}{Buy\PYGZus{}bounds}\PYG{p}{)}

\PYG{k}{def}\PYG{+w}{ }\PYG{n+nf}{Objective\PYGZus{}rule}\PYG{p}{(}\PYG{n}{model}\PYG{p}{)}\PYG{p}{:}
  \PYG{k}{return} \PYG{n}{sum\PYGZus{}product}\PYG{p}{(}\PYG{n}{model}\PYG{o}{.}\PYG{n}{cost}\PYG{p}{,} \PYG{n}{model}\PYG{o}{.}\PYG{n}{Buy}\PYG{p}{)}
\PYG{n}{model}\PYG{o}{.}\PYG{n}{totalcost} \PYG{o}{=} \PYG{n}{Objective}\PYG{p}{(}\PYG{n}{rule}\PYG{o}{=}\PYG{n}{Objective\PYGZus{}rule}\PYG{p}{,} \PYG{n}{sense}\PYG{o}{=}\PYG{n}{minimize}\PYG{p}{)}

\PYG{k}{def}\PYG{+w}{ }\PYG{n+nf}{Diet\PYGZus{}rule}\PYG{p}{(}\PYG{n}{model}\PYG{p}{,} \PYG{n}{i}\PYG{p}{)}\PYG{p}{:}
  \PYG{n}{expr} \PYG{o}{=} \PYG{l+m+mi}{0}

  \PYG{k}{for} \PYG{n}{j} \PYG{o+ow}{in} \PYG{n}{model}\PYG{o}{.}\PYG{n}{FOOD}\PYG{p}{:}
    \PYG{n}{expr} \PYG{o}{=} \PYG{n}{expr} \PYG{o}{+} \PYG{n}{model}\PYG{o}{.}\PYG{n}{amt}\PYG{p}{[}\PYG{n}{i}\PYG{p}{,} \PYG{n}{j}\PYG{p}{]} \PYG{o}{*} \PYG{n}{model}\PYG{o}{.}\PYG{n}{Buy}\PYG{p}{[}\PYG{n}{j}\PYG{p}{]}

  \PYG{k}{return} \PYG{p}{(}\PYG{n}{model}\PYG{o}{.}\PYG{n}{n\PYGZus{}min}\PYG{p}{[}\PYG{n}{i}\PYG{p}{]}\PYG{p}{,} \PYG{n}{expr}\PYG{p}{,} \PYG{n}{model}\PYG{o}{.}\PYG{n}{n\PYGZus{}max}\PYG{p}{[}\PYG{n}{i}\PYG{p}{]}\PYG{p}{)}
\PYG{n}{model}\PYG{o}{.}\PYG{n}{Diet} \PYG{o}{=} \PYG{n}{Constraint}\PYG{p}{(}\PYG{n}{model}\PYG{o}{.}\PYG{n}{NUTR}\PYG{p}{,} \PYG{n}{rule}\PYG{o}{=}\PYG{n}{Diet\PYGZus{}rule}\PYG{p}{)}
\end{sphinxVerbatim}

\sphinxAtStartPar
And the data file \sphinxcode{\sphinxupquote{pydiet\_abstract.dat}} in \hyperref[\detokenize{pyomointerface:coptcode-pydietabs-data}]{Listing \ref{\detokenize{pyomointerface:coptcode-pydietabs-data}}}:
\sphinxSetupCaptionForVerbatim{\sphinxcode{\sphinxupquote{pydiet\_abstract.dat}}}
\def\sphinxLiteralBlockLabel{\label{\detokenize{pyomointerface:coptcode-pydietabs-data}}}
\begin{sphinxVerbatim}[commandchars=\\\{\},numbers=left,firstnumber=1,stepnumber=1]
\PYGZsh{} The data is adopted from:
\PYGZsh{} 
\PYGZsh{} https://github.com/Pyomo/pyomo/blob/master/examples/pyomo/amplbook2/diet.dat
\PYGZsh{}
\PYGZsh{} with some modification by developer of the Cardinal Optimizer

data;

set NUTR := A B1 B2 C ;
set FOOD := BEEF CHK FISH HAM MCH MTL SPG TUR ;

param:   cost  f\PYGZus{}min  f\PYGZus{}max :=
  BEEF   3.19    0     100
  CHK    2.59    0     100
  FISH   2.29    0     100
  HAM    2.89    0     100
  MCH    1.89    0     100
  MTL    1.99    0     100
  SPG    1.99    0     100
  TUR    2.49    0     100 ;

param:   n\PYGZus{}min  n\PYGZus{}max :=
   A      700   10000
   C      700   10000
   B1     700   10000
   B2     700   10000 ;

param amt (tr):
           A    C   B1   B2 :=
   BEEF   60   20   10   15
   CHK     8    0   20   20
   FISH    8   10   15   10
   HAM    40   40   35   10
   MCH    15   35   15   15
   MTL    70   30   15   15
   SPG    25   50   25   15
   TUR    60   20   15   10 ;
\end{sphinxVerbatim}

\sphinxAtStartPar
To solve the problem using Pyomo and the Cardinal Optimizer, just type commands
below in command prompt on Windows or Bash shell on Linux and MacOS.

\begin{sphinxVerbatim}[commandchars=\\\{\}]
pyomo\PYG{+w}{ }solve\PYG{+w}{ }\PYGZhy{}\PYGZhy{}solver\PYG{o}{=}coptampl\PYG{+w}{ }pydiet\PYGZus{}abstract.py\PYG{+w}{ }pydiet\PYGZus{}abstract.dat
\end{sphinxVerbatim}

\sphinxAtStartPar
When solving the problem, Pyomo write log information to the screen:

\begin{sphinxVerbatim}[commandchars=\\\{\}]
[    0.00] Setting up Pyomo environment
[    0.00] Applying Pyomo preprocessing actions
[    0.00] Creating model
[    0.01] Applying solver
[    0.05] Processing results
    Number of solutions: 1
    Solution Information
      Gap: None
      Status: optimal
      Function Value: 88.19999999999999
    Solver results file: results.yml
[    0.05] Applying Pyomo postprocessing actions
[    0.05] Pyomo Finished
\end{sphinxVerbatim}

\sphinxAtStartPar
Upon completion, you can check the solution summary in \sphinxcode{\sphinxupquote{results.yml}}:

\begin{sphinxVerbatim}[commandchars=\\\{\}]
\PYGZsh{} ==========================================================
\PYGZsh{} = Solver Results                                         =
\PYGZsh{} ==========================================================
\PYGZsh{} \PYGZhy{}\PYGZhy{}\PYGZhy{}\PYGZhy{}\PYGZhy{}\PYGZhy{}\PYGZhy{}\PYGZhy{}\PYGZhy{}\PYGZhy{}\PYGZhy{}\PYGZhy{}\PYGZhy{}\PYGZhy{}\PYGZhy{}\PYGZhy{}\PYGZhy{}\PYGZhy{}\PYGZhy{}\PYGZhy{}\PYGZhy{}\PYGZhy{}\PYGZhy{}\PYGZhy{}\PYGZhy{}\PYGZhy{}\PYGZhy{}\PYGZhy{}\PYGZhy{}\PYGZhy{}\PYGZhy{}\PYGZhy{}\PYGZhy{}\PYGZhy{}\PYGZhy{}\PYGZhy{}\PYGZhy{}\PYGZhy{}\PYGZhy{}\PYGZhy{}\PYGZhy{}\PYGZhy{}\PYGZhy{}\PYGZhy{}\PYGZhy{}\PYGZhy{}\PYGZhy{}\PYGZhy{}\PYGZhy{}\PYGZhy{}\PYGZhy{}\PYGZhy{}\PYGZhy{}\PYGZhy{}\PYGZhy{}\PYGZhy{}\PYGZhy{}\PYGZhy{}
\PYGZsh{}   Problem Information
\PYGZsh{} \PYGZhy{}\PYGZhy{}\PYGZhy{}\PYGZhy{}\PYGZhy{}\PYGZhy{}\PYGZhy{}\PYGZhy{}\PYGZhy{}\PYGZhy{}\PYGZhy{}\PYGZhy{}\PYGZhy{}\PYGZhy{}\PYGZhy{}\PYGZhy{}\PYGZhy{}\PYGZhy{}\PYGZhy{}\PYGZhy{}\PYGZhy{}\PYGZhy{}\PYGZhy{}\PYGZhy{}\PYGZhy{}\PYGZhy{}\PYGZhy{}\PYGZhy{}\PYGZhy{}\PYGZhy{}\PYGZhy{}\PYGZhy{}\PYGZhy{}\PYGZhy{}\PYGZhy{}\PYGZhy{}\PYGZhy{}\PYGZhy{}\PYGZhy{}\PYGZhy{}\PYGZhy{}\PYGZhy{}\PYGZhy{}\PYGZhy{}\PYGZhy{}\PYGZhy{}\PYGZhy{}\PYGZhy{}\PYGZhy{}\PYGZhy{}\PYGZhy{}\PYGZhy{}\PYGZhy{}\PYGZhy{}\PYGZhy{}\PYGZhy{}\PYGZhy{}\PYGZhy{}
Problem: 
\PYGZhy{} Lower bound: \PYGZhy{}inf
  Upper bound: inf
  Number of objectives: 1
  Number of constraints: 4
  Number of variables: 8
  Sense: unknown
\PYGZsh{} \PYGZhy{}\PYGZhy{}\PYGZhy{}\PYGZhy{}\PYGZhy{}\PYGZhy{}\PYGZhy{}\PYGZhy{}\PYGZhy{}\PYGZhy{}\PYGZhy{}\PYGZhy{}\PYGZhy{}\PYGZhy{}\PYGZhy{}\PYGZhy{}\PYGZhy{}\PYGZhy{}\PYGZhy{}\PYGZhy{}\PYGZhy{}\PYGZhy{}\PYGZhy{}\PYGZhy{}\PYGZhy{}\PYGZhy{}\PYGZhy{}\PYGZhy{}\PYGZhy{}\PYGZhy{}\PYGZhy{}\PYGZhy{}\PYGZhy{}\PYGZhy{}\PYGZhy{}\PYGZhy{}\PYGZhy{}\PYGZhy{}\PYGZhy{}\PYGZhy{}\PYGZhy{}\PYGZhy{}\PYGZhy{}\PYGZhy{}\PYGZhy{}\PYGZhy{}\PYGZhy{}\PYGZhy{}\PYGZhy{}\PYGZhy{}\PYGZhy{}\PYGZhy{}\PYGZhy{}\PYGZhy{}\PYGZhy{}\PYGZhy{}\PYGZhy{}\PYGZhy{}
\PYGZsh{}   Solver Information
\PYGZsh{} \PYGZhy{}\PYGZhy{}\PYGZhy{}\PYGZhy{}\PYGZhy{}\PYGZhy{}\PYGZhy{}\PYGZhy{}\PYGZhy{}\PYGZhy{}\PYGZhy{}\PYGZhy{}\PYGZhy{}\PYGZhy{}\PYGZhy{}\PYGZhy{}\PYGZhy{}\PYGZhy{}\PYGZhy{}\PYGZhy{}\PYGZhy{}\PYGZhy{}\PYGZhy{}\PYGZhy{}\PYGZhy{}\PYGZhy{}\PYGZhy{}\PYGZhy{}\PYGZhy{}\PYGZhy{}\PYGZhy{}\PYGZhy{}\PYGZhy{}\PYGZhy{}\PYGZhy{}\PYGZhy{}\PYGZhy{}\PYGZhy{}\PYGZhy{}\PYGZhy{}\PYGZhy{}\PYGZhy{}\PYGZhy{}\PYGZhy{}\PYGZhy{}\PYGZhy{}\PYGZhy{}\PYGZhy{}\PYGZhy{}\PYGZhy{}\PYGZhy{}\PYGZhy{}\PYGZhy{}\PYGZhy{}\PYGZhy{}\PYGZhy{}\PYGZhy{}\PYGZhy{}
Solver: 
\PYGZhy{} Status: ok
  Message: COPT\PYGZhy{}AMPL\PYGZbs{}x3a optimal solution; objective 88.2, iterations 1
  Termination condition: optimal
  Id: 0
  Error rc: 0
  Time: 0.03171110153198242
\PYGZsh{} \PYGZhy{}\PYGZhy{}\PYGZhy{}\PYGZhy{}\PYGZhy{}\PYGZhy{}\PYGZhy{}\PYGZhy{}\PYGZhy{}\PYGZhy{}\PYGZhy{}\PYGZhy{}\PYGZhy{}\PYGZhy{}\PYGZhy{}\PYGZhy{}\PYGZhy{}\PYGZhy{}\PYGZhy{}\PYGZhy{}\PYGZhy{}\PYGZhy{}\PYGZhy{}\PYGZhy{}\PYGZhy{}\PYGZhy{}\PYGZhy{}\PYGZhy{}\PYGZhy{}\PYGZhy{}\PYGZhy{}\PYGZhy{}\PYGZhy{}\PYGZhy{}\PYGZhy{}\PYGZhy{}\PYGZhy{}\PYGZhy{}\PYGZhy{}\PYGZhy{}\PYGZhy{}\PYGZhy{}\PYGZhy{}\PYGZhy{}\PYGZhy{}\PYGZhy{}\PYGZhy{}\PYGZhy{}\PYGZhy{}\PYGZhy{}\PYGZhy{}\PYGZhy{}\PYGZhy{}\PYGZhy{}\PYGZhy{}\PYGZhy{}\PYGZhy{}\PYGZhy{}
\PYGZsh{}   Solution Information
\PYGZsh{} \PYGZhy{}\PYGZhy{}\PYGZhy{}\PYGZhy{}\PYGZhy{}\PYGZhy{}\PYGZhy{}\PYGZhy{}\PYGZhy{}\PYGZhy{}\PYGZhy{}\PYGZhy{}\PYGZhy{}\PYGZhy{}\PYGZhy{}\PYGZhy{}\PYGZhy{}\PYGZhy{}\PYGZhy{}\PYGZhy{}\PYGZhy{}\PYGZhy{}\PYGZhy{}\PYGZhy{}\PYGZhy{}\PYGZhy{}\PYGZhy{}\PYGZhy{}\PYGZhy{}\PYGZhy{}\PYGZhy{}\PYGZhy{}\PYGZhy{}\PYGZhy{}\PYGZhy{}\PYGZhy{}\PYGZhy{}\PYGZhy{}\PYGZhy{}\PYGZhy{}\PYGZhy{}\PYGZhy{}\PYGZhy{}\PYGZhy{}\PYGZhy{}\PYGZhy{}\PYGZhy{}\PYGZhy{}\PYGZhy{}\PYGZhy{}\PYGZhy{}\PYGZhy{}\PYGZhy{}\PYGZhy{}\PYGZhy{}\PYGZhy{}\PYGZhy{}\PYGZhy{}
Solution: 
\PYGZhy{} number of solutions: 1
  number of solutions displayed: 1
\PYGZhy{} Gap: None
  Status: optimal
  Message: COPT\PYGZhy{}AMPL\PYGZbs{}x3a optimal solution; objective 88.2, iterations 1
  Objective:
    totalcost:
      Value: 88.19999999999999
  Variable:
    Buy[MCH]:
      Value: 46.666666666666664
  Constraint: No values
\end{sphinxVerbatim}

\sphinxAtStartPar
So the minimal total cost is about 88.2 units when buying 46.67 units of MCH.

\subsubsection{Concrete Model}
\label{\detokenize{pyomointerface:concrete-model}}
\sphinxAtStartPar
The other approach to construct model in Pyomo is to use \sphinxstylestrong{Concrete Model},
we will show how to model and solve the above problem in this way.

\sphinxAtStartPar
Concrete models can be solved using the \sphinxcode{\sphinxupquote{"Direct"}} and \sphinxcode{\sphinxupquote{"Persistent"}} interface methods. This method relies on the Pyomo plugin file \sphinxcode{\sphinxupquote{"copt\_pyomo.py"}} of COPT, which is located in the \sphinxcode{\sphinxupquote{"lib/pyomo"}} subfolder of the installation package.

\sphinxAtStartPar
To use this plugin, you need to copy the \sphinxcode{\sphinxupquote{"copt\_pyomo.py"}} file to the same directory of your program, and have correctly installed the corresponding version of the Python interface of COPT(\sphinxcode{\sphinxupquote{coptpy}}).

\sphinxAtStartPar
The source code \sphinxcode{\sphinxupquote{pydiet\_concrete.py}} is shown in
\hyperref[\detokenize{pyomointerface:coptcode-pydietcon-code}]{Listing \ref{\detokenize{pyomointerface:coptcode-pydietcon-code}}}:
\sphinxSetupCaptionForVerbatim{\sphinxcode{\sphinxupquote{pydiet\_concrete.py}}}
\def\sphinxLiteralBlockLabel{\label{\detokenize{pyomointerface:coptcode-pydietcon-code}}}
\begin{sphinxVerbatim}[commandchars=\\\{\},numbers=left,firstnumber=1,stepnumber=1]
\PYG{c+c1}{\PYGZsh{} The code is adopted from:}
\PYG{c+c1}{\PYGZsh{}}
\PYG{c+c1}{\PYGZsh{} https://github.com/Pyomo/pyomo/blob/master/examples/pyomo/amplbook2/diet.py}
\PYG{c+c1}{\PYGZsh{}}
\PYG{c+c1}{\PYGZsh{} with some modification by developer of the Cardinal Optimizer}

\PYG{k+kn}{from}\PYG{+w}{ }\PYG{n+nn}{\PYGZus{}\PYGZus{}future\PYGZus{}\PYGZus{}}\PYG{+w}{ }\PYG{k+kn}{import} \PYG{n}{print\PYGZus{}function}\PYG{p}{,} \PYG{n}{division}

\PYG{k+kn}{import}\PYG{+w}{ }\PYG{n+nn}{pyomo}\PYG{n+nn}{.}\PYG{n+nn}{environ}\PYG{+w}{ }\PYG{k}{as}\PYG{+w}{ }\PYG{n+nn}{pyo}
\PYG{k+kn}{import}\PYG{+w}{ }\PYG{n+nn}{pyomo}\PYG{n+nn}{.}\PYG{n+nn}{opt}\PYG{+w}{ }\PYG{k}{as}\PYG{+w}{ }\PYG{n+nn}{pyopt}

\PYG{k+kn}{from}\PYG{+w}{ }\PYG{n+nn}{copt\PYGZus{}pyomo}\PYG{+w}{ }\PYG{k+kn}{import} \PYG{o}{*}

\PYG{c+c1}{\PYGZsh{} Nutrition set}
\PYG{n}{NUTR} \PYG{o}{=} \PYG{p}{[}\PYG{l+s+s2}{\PYGZdq{}}\PYG{l+s+s2}{A}\PYG{l+s+s2}{\PYGZdq{}}\PYG{p}{,} \PYG{l+s+s2}{\PYGZdq{}}\PYG{l+s+s2}{C}\PYG{l+s+s2}{\PYGZdq{}}\PYG{p}{,} \PYG{l+s+s2}{\PYGZdq{}}\PYG{l+s+s2}{B1}\PYG{l+s+s2}{\PYGZdq{}}\PYG{p}{,} \PYG{l+s+s2}{\PYGZdq{}}\PYG{l+s+s2}{B2}\PYG{l+s+s2}{\PYGZdq{}}\PYG{p}{]}
\PYG{c+c1}{\PYGZsh{} Food set}
\PYG{n}{FOOD} \PYG{o}{=} \PYG{p}{[}\PYG{l+s+s2}{\PYGZdq{}}\PYG{l+s+s2}{BEEF}\PYG{l+s+s2}{\PYGZdq{}}\PYG{p}{,} \PYG{l+s+s2}{\PYGZdq{}}\PYG{l+s+s2}{CHK}\PYG{l+s+s2}{\PYGZdq{}}\PYG{p}{,} \PYG{l+s+s2}{\PYGZdq{}}\PYG{l+s+s2}{FISH}\PYG{l+s+s2}{\PYGZdq{}}\PYG{p}{,} \PYG{l+s+s2}{\PYGZdq{}}\PYG{l+s+s2}{HAM}\PYG{l+s+s2}{\PYGZdq{}}\PYG{p}{,} \PYG{l+s+s2}{\PYGZdq{}}\PYG{l+s+s2}{MCH}\PYG{l+s+s2}{\PYGZdq{}}\PYG{p}{,} \PYG{l+s+s2}{\PYGZdq{}}\PYG{l+s+s2}{MTL}\PYG{l+s+s2}{\PYGZdq{}}\PYG{p}{,} \PYG{l+s+s2}{\PYGZdq{}}\PYG{l+s+s2}{SPG}\PYG{l+s+s2}{\PYGZdq{}}\PYG{p}{,} \PYG{l+s+s2}{\PYGZdq{}}\PYG{l+s+s2}{TUR}\PYG{l+s+s2}{\PYGZdq{}}\PYG{p}{]}

\PYG{c+c1}{\PYGZsh{} Price of foods}
\PYG{n}{cost} \PYG{o}{=} \PYG{p}{\PYGZob{}}\PYG{l+s+s2}{\PYGZdq{}}\PYG{l+s+s2}{BEEF}\PYG{l+s+s2}{\PYGZdq{}}\PYG{p}{:} \PYG{l+m+mf}{3.19}\PYG{p}{,} \PYG{l+s+s2}{\PYGZdq{}}\PYG{l+s+s2}{CHK}\PYG{l+s+s2}{\PYGZdq{}}\PYG{p}{:} \PYG{l+m+mf}{2.59}\PYG{p}{,} \PYG{l+s+s2}{\PYGZdq{}}\PYG{l+s+s2}{FISH}\PYG{l+s+s2}{\PYGZdq{}}\PYG{p}{:} \PYG{l+m+mf}{2.29}\PYG{p}{,} \PYG{l+s+s2}{\PYGZdq{}}\PYG{l+s+s2}{HAM}\PYG{l+s+s2}{\PYGZdq{}}\PYG{p}{:} \PYG{l+m+mf}{2.89}\PYG{p}{,} \PYG{l+s+s2}{\PYGZdq{}}\PYG{l+s+s2}{MCH}\PYG{l+s+s2}{\PYGZdq{}}\PYG{p}{:} \PYG{l+m+mf}{1.89}\PYG{p}{,}
        \PYG{l+s+s2}{\PYGZdq{}}\PYG{l+s+s2}{MTL}\PYG{l+s+s2}{\PYGZdq{}}\PYG{p}{:}  \PYG{l+m+mf}{1.99}\PYG{p}{,} \PYG{l+s+s2}{\PYGZdq{}}\PYG{l+s+s2}{SPG}\PYG{l+s+s2}{\PYGZdq{}}\PYG{p}{:} \PYG{l+m+mf}{1.99}\PYG{p}{,} \PYG{l+s+s2}{\PYGZdq{}}\PYG{l+s+s2}{TUR}\PYG{l+s+s2}{\PYGZdq{}}\PYG{p}{:}  \PYG{l+m+mf}{2.49}\PYG{p}{\PYGZcb{}}
\PYG{c+c1}{\PYGZsh{} Nutrition of foods}
\PYG{n}{amt} \PYG{o}{=} \PYG{p}{\PYGZob{}}\PYG{l+s+s2}{\PYGZdq{}}\PYG{l+s+s2}{BEEF}\PYG{l+s+s2}{\PYGZdq{}}\PYG{p}{:} \PYG{p}{\PYGZob{}}\PYG{l+s+s2}{\PYGZdq{}}\PYG{l+s+s2}{A}\PYG{l+s+s2}{\PYGZdq{}}\PYG{p}{:} \PYG{l+m+mi}{60}\PYG{p}{,} \PYG{l+s+s2}{\PYGZdq{}}\PYG{l+s+s2}{C}\PYG{l+s+s2}{\PYGZdq{}}\PYG{p}{:} \PYG{l+m+mi}{20}\PYG{p}{,} \PYG{l+s+s2}{\PYGZdq{}}\PYG{l+s+s2}{B1}\PYG{l+s+s2}{\PYGZdq{}}\PYG{p}{:} \PYG{l+m+mi}{10}\PYG{p}{,} \PYG{l+s+s2}{\PYGZdq{}}\PYG{l+s+s2}{B2}\PYG{l+s+s2}{\PYGZdq{}}\PYG{p}{:} \PYG{l+m+mi}{15}\PYG{p}{\PYGZcb{}}\PYG{p}{,}
       \PYG{l+s+s2}{\PYGZdq{}}\PYG{l+s+s2}{CHK}\PYG{l+s+s2}{\PYGZdq{}}\PYG{p}{:}  \PYG{p}{\PYGZob{}}\PYG{l+s+s2}{\PYGZdq{}}\PYG{l+s+s2}{A}\PYG{l+s+s2}{\PYGZdq{}}\PYG{p}{:} \PYG{l+m+mi}{8}\PYG{p}{,}  \PYG{l+s+s2}{\PYGZdq{}}\PYG{l+s+s2}{C}\PYG{l+s+s2}{\PYGZdq{}}\PYG{p}{:} \PYG{l+m+mi}{0}\PYG{p}{,}  \PYG{l+s+s2}{\PYGZdq{}}\PYG{l+s+s2}{B1}\PYG{l+s+s2}{\PYGZdq{}}\PYG{p}{:} \PYG{l+m+mi}{20}\PYG{p}{,} \PYG{l+s+s2}{\PYGZdq{}}\PYG{l+s+s2}{B2}\PYG{l+s+s2}{\PYGZdq{}}\PYG{p}{:} \PYG{l+m+mi}{20}\PYG{p}{\PYGZcb{}}\PYG{p}{,}
       \PYG{l+s+s2}{\PYGZdq{}}\PYG{l+s+s2}{FISH}\PYG{l+s+s2}{\PYGZdq{}}\PYG{p}{:} \PYG{p}{\PYGZob{}}\PYG{l+s+s2}{\PYGZdq{}}\PYG{l+s+s2}{A}\PYG{l+s+s2}{\PYGZdq{}}\PYG{p}{:} \PYG{l+m+mi}{8}\PYG{p}{,}  \PYG{l+s+s2}{\PYGZdq{}}\PYG{l+s+s2}{C}\PYG{l+s+s2}{\PYGZdq{}}\PYG{p}{:} \PYG{l+m+mi}{10}\PYG{p}{,} \PYG{l+s+s2}{\PYGZdq{}}\PYG{l+s+s2}{B1}\PYG{l+s+s2}{\PYGZdq{}}\PYG{p}{:} \PYG{l+m+mi}{15}\PYG{p}{,} \PYG{l+s+s2}{\PYGZdq{}}\PYG{l+s+s2}{B2}\PYG{l+s+s2}{\PYGZdq{}}\PYG{p}{:} \PYG{l+m+mi}{10}\PYG{p}{\PYGZcb{}}\PYG{p}{,}
       \PYG{l+s+s2}{\PYGZdq{}}\PYG{l+s+s2}{HAM}\PYG{l+s+s2}{\PYGZdq{}}\PYG{p}{:}  \PYG{p}{\PYGZob{}}\PYG{l+s+s2}{\PYGZdq{}}\PYG{l+s+s2}{A}\PYG{l+s+s2}{\PYGZdq{}}\PYG{p}{:} \PYG{l+m+mi}{40}\PYG{p}{,} \PYG{l+s+s2}{\PYGZdq{}}\PYG{l+s+s2}{C}\PYG{l+s+s2}{\PYGZdq{}}\PYG{p}{:} \PYG{l+m+mi}{40}\PYG{p}{,} \PYG{l+s+s2}{\PYGZdq{}}\PYG{l+s+s2}{B1}\PYG{l+s+s2}{\PYGZdq{}}\PYG{p}{:} \PYG{l+m+mi}{35}\PYG{p}{,} \PYG{l+s+s2}{\PYGZdq{}}\PYG{l+s+s2}{B2}\PYG{l+s+s2}{\PYGZdq{}}\PYG{p}{:} \PYG{l+m+mi}{10}\PYG{p}{\PYGZcb{}}\PYG{p}{,}
       \PYG{l+s+s2}{\PYGZdq{}}\PYG{l+s+s2}{MCH}\PYG{l+s+s2}{\PYGZdq{}}\PYG{p}{:}  \PYG{p}{\PYGZob{}}\PYG{l+s+s2}{\PYGZdq{}}\PYG{l+s+s2}{A}\PYG{l+s+s2}{\PYGZdq{}}\PYG{p}{:} \PYG{l+m+mi}{15}\PYG{p}{,} \PYG{l+s+s2}{\PYGZdq{}}\PYG{l+s+s2}{C}\PYG{l+s+s2}{\PYGZdq{}}\PYG{p}{:} \PYG{l+m+mi}{35}\PYG{p}{,} \PYG{l+s+s2}{\PYGZdq{}}\PYG{l+s+s2}{B1}\PYG{l+s+s2}{\PYGZdq{}}\PYG{p}{:} \PYG{l+m+mi}{15}\PYG{p}{,} \PYG{l+s+s2}{\PYGZdq{}}\PYG{l+s+s2}{B2}\PYG{l+s+s2}{\PYGZdq{}}\PYG{p}{:} \PYG{l+m+mi}{15}\PYG{p}{\PYGZcb{}}\PYG{p}{,}
       \PYG{l+s+s2}{\PYGZdq{}}\PYG{l+s+s2}{MTL}\PYG{l+s+s2}{\PYGZdq{}}\PYG{p}{:}  \PYG{p}{\PYGZob{}}\PYG{l+s+s2}{\PYGZdq{}}\PYG{l+s+s2}{A}\PYG{l+s+s2}{\PYGZdq{}}\PYG{p}{:} \PYG{l+m+mi}{70}\PYG{p}{,} \PYG{l+s+s2}{\PYGZdq{}}\PYG{l+s+s2}{C}\PYG{l+s+s2}{\PYGZdq{}}\PYG{p}{:} \PYG{l+m+mi}{30}\PYG{p}{,} \PYG{l+s+s2}{\PYGZdq{}}\PYG{l+s+s2}{B1}\PYG{l+s+s2}{\PYGZdq{}}\PYG{p}{:} \PYG{l+m+mi}{15}\PYG{p}{,} \PYG{l+s+s2}{\PYGZdq{}}\PYG{l+s+s2}{B2}\PYG{l+s+s2}{\PYGZdq{}}\PYG{p}{:} \PYG{l+m+mi}{15}\PYG{p}{\PYGZcb{}}\PYG{p}{,}
       \PYG{l+s+s2}{\PYGZdq{}}\PYG{l+s+s2}{SPG}\PYG{l+s+s2}{\PYGZdq{}}\PYG{p}{:}  \PYG{p}{\PYGZob{}}\PYG{l+s+s2}{\PYGZdq{}}\PYG{l+s+s2}{A}\PYG{l+s+s2}{\PYGZdq{}}\PYG{p}{:} \PYG{l+m+mi}{25}\PYG{p}{,} \PYG{l+s+s2}{\PYGZdq{}}\PYG{l+s+s2}{C}\PYG{l+s+s2}{\PYGZdq{}}\PYG{p}{:} \PYG{l+m+mi}{50}\PYG{p}{,} \PYG{l+s+s2}{\PYGZdq{}}\PYG{l+s+s2}{B1}\PYG{l+s+s2}{\PYGZdq{}}\PYG{p}{:} \PYG{l+m+mi}{25}\PYG{p}{,} \PYG{l+s+s2}{\PYGZdq{}}\PYG{l+s+s2}{B2}\PYG{l+s+s2}{\PYGZdq{}}\PYG{p}{:} \PYG{l+m+mi}{15}\PYG{p}{\PYGZcb{}}\PYG{p}{,}
       \PYG{l+s+s2}{\PYGZdq{}}\PYG{l+s+s2}{TUR}\PYG{l+s+s2}{\PYGZdq{}}\PYG{p}{:}  \PYG{p}{\PYGZob{}}\PYG{l+s+s2}{\PYGZdq{}}\PYG{l+s+s2}{A}\PYG{l+s+s2}{\PYGZdq{}}\PYG{p}{:} \PYG{l+m+mi}{60}\PYG{p}{,} \PYG{l+s+s2}{\PYGZdq{}}\PYG{l+s+s2}{C}\PYG{l+s+s2}{\PYGZdq{}}\PYG{p}{:} \PYG{l+m+mi}{20}\PYG{p}{,} \PYG{l+s+s2}{\PYGZdq{}}\PYG{l+s+s2}{B1}\PYG{l+s+s2}{\PYGZdq{}}\PYG{p}{:} \PYG{l+m+mi}{15}\PYG{p}{,} \PYG{l+s+s2}{\PYGZdq{}}\PYG{l+s+s2}{B2}\PYG{l+s+s2}{\PYGZdq{}}\PYG{p}{:} \PYG{l+m+mi}{10}\PYG{p}{\PYGZcb{}}\PYG{p}{\PYGZcb{}}

\PYG{c+c1}{\PYGZsh{} The \PYGZdq{}diet problem\PYGZdq{} using ConcreteModel}
\PYG{n}{model} \PYG{o}{=} \PYG{n}{pyo}\PYG{o}{.}\PYG{n}{ConcreteModel}\PYG{p}{(}\PYG{p}{)}

\PYG{n}{model}\PYG{o}{.}\PYG{n}{NUTR} \PYG{o}{=} \PYG{n}{pyo}\PYG{o}{.}\PYG{n}{Set}\PYG{p}{(}\PYG{n}{initialize}\PYG{o}{=}\PYG{n}{NUTR}\PYG{p}{)}
\PYG{n}{model}\PYG{o}{.}\PYG{n}{FOOD} \PYG{o}{=} \PYG{n}{pyo}\PYG{o}{.}\PYG{n}{Set}\PYG{p}{(}\PYG{n}{initialize}\PYG{o}{=}\PYG{n}{FOOD}\PYG{p}{)}

\PYG{n}{model}\PYG{o}{.}\PYG{n}{cost} \PYG{o}{=} \PYG{n}{pyo}\PYG{o}{.}\PYG{n}{Param}\PYG{p}{(}\PYG{n}{model}\PYG{o}{.}\PYG{n}{FOOD}\PYG{p}{,} \PYG{n}{initialize}\PYG{o}{=}\PYG{n}{cost}\PYG{p}{)}

\PYG{k}{def}\PYG{+w}{ }\PYG{n+nf}{amt\PYGZus{}rule}\PYG{p}{(}\PYG{n}{model}\PYG{p}{,} \PYG{n}{i}\PYG{p}{,} \PYG{n}{j}\PYG{p}{)}\PYG{p}{:}
  \PYG{k}{return} \PYG{n}{amt}\PYG{p}{[}\PYG{n}{i}\PYG{p}{]}\PYG{p}{[}\PYG{n}{j}\PYG{p}{]}
\PYG{n}{model}\PYG{o}{.}\PYG{n}{amt}  \PYG{o}{=} \PYG{n}{pyo}\PYG{o}{.}\PYG{n}{Param}\PYG{p}{(}\PYG{n}{model}\PYG{o}{.}\PYG{n}{FOOD}\PYG{p}{,} \PYG{n}{model}\PYG{o}{.}\PYG{n}{NUTR}\PYG{p}{,} \PYG{n}{initialize}\PYG{o}{=}\PYG{n}{amt\PYGZus{}rule}\PYG{p}{)}

\PYG{n}{model}\PYG{o}{.}\PYG{n}{f\PYGZus{}min} \PYG{o}{=} \PYG{n}{pyo}\PYG{o}{.}\PYG{n}{Param}\PYG{p}{(}\PYG{n}{model}\PYG{o}{.}\PYG{n}{FOOD}\PYG{p}{,} \PYG{n}{default}\PYG{o}{=}\PYG{l+m+mi}{0}\PYG{p}{)}
\PYG{n}{model}\PYG{o}{.}\PYG{n}{f\PYGZus{}max} \PYG{o}{=} \PYG{n}{pyo}\PYG{o}{.}\PYG{n}{Param}\PYG{p}{(}\PYG{n}{model}\PYG{o}{.}\PYG{n}{FOOD}\PYG{p}{,} \PYG{n}{default}\PYG{o}{=}\PYG{l+m+mi}{100}\PYG{p}{)}

\PYG{n}{model}\PYG{o}{.}\PYG{n}{n\PYGZus{}min} \PYG{o}{=} \PYG{n}{pyo}\PYG{o}{.}\PYG{n}{Param}\PYG{p}{(}\PYG{n}{model}\PYG{o}{.}\PYG{n}{NUTR}\PYG{p}{,} \PYG{n}{default}\PYG{o}{=}\PYG{l+m+mi}{700}\PYG{p}{)}
\PYG{n}{model}\PYG{o}{.}\PYG{n}{n\PYGZus{}max} \PYG{o}{=} \PYG{n}{pyo}\PYG{o}{.}\PYG{n}{Param}\PYG{p}{(}\PYG{n}{model}\PYG{o}{.}\PYG{n}{NUTR}\PYG{p}{,} \PYG{n}{default}\PYG{o}{=}\PYG{l+m+mi}{10000}\PYG{p}{)}

\PYG{k}{def}\PYG{+w}{ }\PYG{n+nf}{Buy\PYGZus{}bounds}\PYG{p}{(}\PYG{n}{model}\PYG{p}{,} \PYG{n}{i}\PYG{p}{)}\PYG{p}{:}
  \PYG{k}{return} \PYG{p}{(}\PYG{n}{model}\PYG{o}{.}\PYG{n}{f\PYGZus{}min}\PYG{p}{[}\PYG{n}{i}\PYG{p}{]}\PYG{p}{,} \PYG{n}{model}\PYG{o}{.}\PYG{n}{f\PYGZus{}max}\PYG{p}{[}\PYG{n}{i}\PYG{p}{]}\PYG{p}{)}
\PYG{n}{model}\PYG{o}{.}\PYG{n}{buy} \PYG{o}{=} \PYG{n}{pyo}\PYG{o}{.}\PYG{n}{Var}\PYG{p}{(}\PYG{n}{model}\PYG{o}{.}\PYG{n}{FOOD}\PYG{p}{,} \PYG{n}{bounds}\PYG{o}{=}\PYG{n}{Buy\PYGZus{}bounds}\PYG{p}{)}

\PYG{k}{def}\PYG{+w}{ }\PYG{n+nf}{Objective\PYGZus{}rule}\PYG{p}{(}\PYG{n}{model}\PYG{p}{)}\PYG{p}{:}
  \PYG{k}{return} \PYG{n}{pyo}\PYG{o}{.}\PYG{n}{sum\PYGZus{}product}\PYG{p}{(}\PYG{n}{model}\PYG{o}{.}\PYG{n}{cost}\PYG{p}{,} \PYG{n}{model}\PYG{o}{.}\PYG{n}{buy}\PYG{p}{)}
\PYG{n}{model}\PYG{o}{.}\PYG{n}{totalcost} \PYG{o}{=} \PYG{n}{pyo}\PYG{o}{.}\PYG{n}{Objective}\PYG{p}{(}\PYG{n}{rule}\PYG{o}{=}\PYG{n}{Objective\PYGZus{}rule}\PYG{p}{,} \PYG{n}{sense}\PYG{o}{=}\PYG{n}{pyo}\PYG{o}{.}\PYG{n}{minimize}\PYG{p}{)}

\PYG{k}{def}\PYG{+w}{ }\PYG{n+nf}{Diet\PYGZus{}rule}\PYG{p}{(}\PYG{n}{model}\PYG{p}{,} \PYG{n}{j}\PYG{p}{)}\PYG{p}{:}
  \PYG{n}{expr} \PYG{o}{=} \PYG{l+m+mi}{0}

  \PYG{k}{for} \PYG{n}{i} \PYG{o+ow}{in} \PYG{n}{model}\PYG{o}{.}\PYG{n}{FOOD}\PYG{p}{:}
    \PYG{n}{expr} \PYG{o}{=} \PYG{n}{expr} \PYG{o}{+} \PYG{n}{model}\PYG{o}{.}\PYG{n}{amt}\PYG{p}{[}\PYG{n}{i}\PYG{p}{,} \PYG{n}{j}\PYG{p}{]} \PYG{o}{*} \PYG{n}{model}\PYG{o}{.}\PYG{n}{buy}\PYG{p}{[}\PYG{n}{i}\PYG{p}{]}

  \PYG{k}{return} \PYG{p}{(}\PYG{n}{model}\PYG{o}{.}\PYG{n}{n\PYGZus{}min}\PYG{p}{[}\PYG{n}{j}\PYG{p}{]}\PYG{p}{,} \PYG{n}{expr}\PYG{p}{,} \PYG{n}{model}\PYG{o}{.}\PYG{n}{n\PYGZus{}max}\PYG{p}{[}\PYG{n}{j}\PYG{p}{]}\PYG{p}{)}
\PYG{n}{model}\PYG{o}{.}\PYG{n}{Diet} \PYG{o}{=} \PYG{n}{pyo}\PYG{o}{.}\PYG{n}{Constraint}\PYG{p}{(}\PYG{n}{model}\PYG{o}{.}\PYG{n}{NUTR}\PYG{p}{,} \PYG{n}{rule}\PYG{o}{=}\PYG{n}{Diet\PYGZus{}rule}\PYG{p}{)}

\PYG{c+c1}{\PYGZsh{} Reduced costs of variables}
\PYG{n}{model}\PYG{o}{.}\PYG{n}{rc} \PYG{o}{=} \PYG{n}{pyo}\PYG{o}{.}\PYG{n}{Suffix}\PYG{p}{(}\PYG{n}{direction}\PYG{o}{=}\PYG{n}{pyo}\PYG{o}{.}\PYG{n}{Suffix}\PYG{o}{.}\PYG{n}{IMPORT}\PYG{p}{)}

\PYG{c+c1}{\PYGZsh{} Activities and duals of constraints}
\PYG{n}{model}\PYG{o}{.}\PYG{n}{slack} \PYG{o}{=} \PYG{n}{pyo}\PYG{o}{.}\PYG{n}{Suffix}\PYG{p}{(}\PYG{n}{direction}\PYG{o}{=}\PYG{n}{pyo}\PYG{o}{.}\PYG{n}{Suffix}\PYG{o}{.}\PYG{n}{IMPORT}\PYG{p}{)}
\PYG{n}{model}\PYG{o}{.}\PYG{n}{dual} \PYG{o}{=} \PYG{n}{pyo}\PYG{o}{.}\PYG{n}{Suffix}\PYG{p}{(}\PYG{n}{direction}\PYG{o}{=}\PYG{n}{pyo}\PYG{o}{.}\PYG{n}{Suffix}\PYG{o}{.}\PYG{n}{IMPORT}\PYG{p}{)}

\PYG{c+c1}{\PYGZsh{} Use \PYGZsq{}copt\PYGZus{}direct\PYGZsq{} solver to solve the problem}
\PYG{n}{solver} \PYG{o}{=} \PYG{n}{pyopt}\PYG{o}{.}\PYG{n}{SolverFactory}\PYG{p}{(}\PYG{l+s+s1}{\PYGZsq{}}\PYG{l+s+s1}{copt\PYGZus{}direct}\PYG{l+s+s1}{\PYGZsq{}}\PYG{p}{)}

\PYG{c+c1}{\PYGZsh{} Use \PYGZsq{}copt\PYGZus{}persistent\PYGZsq{} solver to solve the problem}
\PYG{c+c1}{\PYGZsh{} solver = pyopt.SolverFactory(\PYGZsq{}copt\PYGZus{}persistent\PYGZsq{})}
\PYG{c+c1}{\PYGZsh{} solver.set\PYGZus{}instance(model)}

\PYG{n}{results} \PYG{o}{=} \PYG{n}{solver}\PYG{o}{.}\PYG{n}{solve}\PYG{p}{(}\PYG{n}{model}\PYG{p}{,} \PYG{n}{tee}\PYG{o}{=}\PYG{k+kc}{True}\PYG{p}{)}

\PYG{c+c1}{\PYGZsh{} Check result}
\PYG{n+nb}{print}\PYG{p}{(}\PYG{l+s+s2}{\PYGZdq{}}\PYG{l+s+s2}{\PYGZdq{}}\PYG{p}{)}
\PYG{k}{if} \PYG{n}{results}\PYG{o}{.}\PYG{n}{solver}\PYG{o}{.}\PYG{n}{status} \PYG{o}{==} \PYG{n}{pyopt}\PYG{o}{.}\PYG{n}{SolverStatus}\PYG{o}{.}\PYG{n}{ok} \PYG{o+ow}{and} \PYGZbs{}
   \PYG{n}{results}\PYG{o}{.}\PYG{n}{solver}\PYG{o}{.}\PYG{n}{termination\PYGZus{}condition} \PYG{o}{==} \PYG{n}{pyopt}\PYG{o}{.}\PYG{n}{TerminationCondition}\PYG{o}{.}\PYG{n}{optimal}\PYG{p}{:}
  \PYG{n+nb}{print}\PYG{p}{(}\PYG{l+s+s2}{\PYGZdq{}}\PYG{l+s+s2}{Optimal solution found}\PYG{l+s+s2}{\PYGZdq{}}\PYG{p}{)}
\PYG{k}{else}\PYG{p}{:}
  \PYG{n+nb}{print}\PYG{p}{(}\PYG{l+s+s2}{\PYGZdq{}}\PYG{l+s+s2}{Something unexpected happened: }\PYG{l+s+s2}{\PYGZdq{}}\PYG{p}{,} \PYG{n+nb}{str}\PYG{p}{(}\PYG{n}{results}\PYG{o}{.}\PYG{n}{solver}\PYG{p}{)}\PYG{p}{)}

\PYG{n+nb}{print}\PYG{p}{(}\PYG{l+s+s2}{\PYGZdq{}}\PYG{l+s+s2}{\PYGZdq{}}\PYG{p}{)}
\PYG{n+nb}{print}\PYG{p}{(}\PYG{l+s+s2}{\PYGZdq{}}\PYG{l+s+s2}{Optimal objective value:}\PYG{l+s+s2}{\PYGZdq{}}\PYG{p}{)}
\PYG{n+nb}{print}\PYG{p}{(}\PYG{l+s+s2}{\PYGZdq{}}\PYG{l+s+s2}{  totalcost: }\PYG{l+s+si}{\PYGZob{}0:6f\PYGZcb{}}\PYG{l+s+s2}{\PYGZdq{}}\PYG{o}{.}\PYG{n}{format}\PYG{p}{(}\PYG{n}{pyo}\PYG{o}{.}\PYG{n}{value}\PYG{p}{(}\PYG{n}{model}\PYG{o}{.}\PYG{n}{totalcost}\PYG{p}{)}\PYG{p}{)}\PYG{p}{)}

\PYG{n+nb}{print}\PYG{p}{(}\PYG{l+s+s2}{\PYGZdq{}}\PYG{l+s+s2}{\PYGZdq{}}\PYG{p}{)}
\PYG{n+nb}{print}\PYG{p}{(}\PYG{l+s+s2}{\PYGZdq{}}\PYG{l+s+s2}{Variables solution:}\PYG{l+s+s2}{\PYGZdq{}}\PYG{p}{)}
\PYG{k}{for} \PYG{n}{i} \PYG{o+ow}{in} \PYG{n}{FOOD}\PYG{p}{:}
  \PYG{n+nb}{print}\PYG{p}{(}\PYG{l+s+s2}{\PYGZdq{}}\PYG{l+s+s2}{  buy[}\PYG{l+s+si}{\PYGZob{}0:4s\PYGZcb{}}\PYG{l+s+s2}{] = }\PYG{l+s+si}{\PYGZob{}1:9.6f\PYGZcb{}}\PYG{l+s+s2}{ (rc: }\PYG{l+s+si}{\PYGZob{}2:9.6f\PYGZcb{}}\PYG{l+s+s2}{)}\PYG{l+s+s2}{\PYGZdq{}}\PYG{o}{.}\PYG{n}{format}\PYG{p}{(}\PYG{n}{i}\PYG{p}{,} \PYGZbs{}
                                                  \PYG{n}{pyo}\PYG{o}{.}\PYG{n}{value}\PYG{p}{(}\PYG{n}{model}\PYG{o}{.}\PYG{n}{buy}\PYG{p}{[}\PYG{n}{i}\PYG{p}{]}\PYG{p}{)}\PYG{p}{,} \PYGZbs{}
                                                  \PYG{n}{model}\PYG{o}{.}\PYG{n}{rc}\PYG{p}{[}\PYG{n}{model}\PYG{o}{.}\PYG{n}{buy}\PYG{p}{[}\PYG{n}{i}\PYG{p}{]}\PYG{p}{]}\PYG{p}{)}\PYG{p}{)}

\PYG{n+nb}{print}\PYG{p}{(}\PYG{l+s+s2}{\PYGZdq{}}\PYG{l+s+s2}{\PYGZdq{}}\PYG{p}{)}
\PYG{n+nb}{print}\PYG{p}{(}\PYG{l+s+s2}{\PYGZdq{}}\PYG{l+s+s2}{Constraint solution:}\PYG{l+s+s2}{\PYGZdq{}}\PYG{p}{)}
\PYG{k}{for} \PYG{n}{i} \PYG{o+ow}{in} \PYG{n}{NUTR}\PYG{p}{:}
  \PYG{n+nb}{print}\PYG{p}{(}\PYG{l+s+s2}{\PYGZdq{}}\PYG{l+s+s2}{  diet[}\PYG{l+s+si}{\PYGZob{}0:2s\PYGZcb{}}\PYG{l+s+s2}{] = }\PYG{l+s+si}{\PYGZob{}1:12.6f\PYGZcb{}}\PYG{l+s+s2}{ (dual: }\PYG{l+s+si}{\PYGZob{}2:9.6f\PYGZcb{}}\PYG{l+s+s2}{)}\PYG{l+s+s2}{\PYGZdq{}}\PYG{o}{.}\PYG{n}{format}\PYG{p}{(}\PYG{n}{i}\PYG{p}{,} \PYGZbs{}
                                                      \PYG{n}{model}\PYG{o}{.}\PYG{n}{slack}\PYG{p}{[}\PYG{n}{model}\PYG{o}{.}\PYG{n}{Diet}\PYG{p}{[}\PYG{n}{i}\PYG{p}{]}\PYG{p}{]}\PYG{p}{,} \PYGZbs{}
                                                      \PYG{n}{model}\PYG{o}{.}\PYG{n}{dual}\PYG{p}{[}\PYG{n}{model}\PYG{o}{.}\PYG{n}{Diet}\PYG{p}{[}\PYG{n}{i}\PYG{p}{]}\PYG{p}{]}\PYG{p}{)}\PYG{p}{)}
\end{sphinxVerbatim}

\sphinxAtStartPar
To solve the problem using Pyomo and the Cardinal Optimizer, just execute
commands below:

\begin{sphinxVerbatim}[commandchars=\\\{\}]
python\PYG{+w}{ }pydiet\PYGZus{}concrete.py
\end{sphinxVerbatim}

\sphinxAtStartPar
Up completion, you should see solution summary on screen as below:

\begin{sphinxVerbatim}[commandchars=\\\{\}]
Optimal solution found
Objective:
  totalcost: 88.200000
Variables:
  buy[BEEF] = 0.000000
  buy[CHK ] = 0.000000
  buy[FISH] = 0.000000
  buy[HAM ] = 0.000000
  buy[MCH ] = 46.666667
  buy[MTL ] = 0.000000
  buy[SPG ] = 0.000000
  buy[TUR ] = 0.000000
\end{sphinxVerbatim}

\sphinxAtStartPar
So the Cardinal Optimizer found the optimal solution, which is about 88.2 units
when buying about 46.67 units of MCH.

\sphinxstepscope

\section{PuLP Interface}
\label{\detokenize{pulpinterface:pulp-interface}}\label{\detokenize{pulpinterface:chappulpinterface}}\label{\detokenize{pulpinterface::doc}}
\sphinxAtStartPar
\sphinxhref{http://coin-or.github.io/pulp/}{PuLP} is an open source modeling tool based on Python,
it is mainly used for modeling integer programming problems.
This chapter introduces how to use the Cardinal Optimizer (COPT) in PuLP.

\subsection{Installation guide}
\label{\detokenize{pulpinterface:installation-guide}}
\sphinxAtStartPar
Before calling COPT in PuLP to solve problem, users need to setup PuLP and COPT correctly.
PuLP currently supports Python 2.7 and later versions of Python. Users can download Python from
\sphinxhref{https://www.anaconda.com/distribution/}{Anaconda distribution} or
\sphinxhref{https://www.python.org/}{Python official distribution} .
We recommend users install the Anaconda distribution, because it is more user\sphinxhyphen{}friendly
and convenient for Python novices.

\subsubsection{Install via conda}
\label{\detokenize{pulpinterface:install-via-conda}}
\sphinxAtStartPar
We recommend that users who have installed the Anaconda distribution of Python
use its own \sphinxcode{\sphinxupquote{conda}} tool to install PuLP. Execute the following commands in
Windows command prompt or terminal on Linux and MacOS:

\begin{sphinxVerbatim}[commandchars=\\\{\}]
conda\PYG{+w}{ }install\PYG{+w}{ }\PYGZhy{}c\PYG{+w}{ }conda\PYGZhy{}forge\PYG{+w}{ }pulp
\end{sphinxVerbatim}

\subsubsection{Install via pip}
\label{\detokenize{pulpinterface:install-via-pip}}
\sphinxAtStartPar
Users can also install PuLP through the standard \sphinxcode{\sphinxupquote{pip}} tool, execute the following
command in Windows command prompt or Linux and MacOS terminal:

\begin{sphinxVerbatim}[commandchars=\\\{\}]
pip\PYG{+w}{ }install\PYG{+w}{ }pulp
\end{sphinxVerbatim}

\subsubsection{Setup PuLP interface}
\label{\detokenize{pulpinterface:setup-pulp-interface}}
\sphinxAtStartPar
For PuLP V2.8.0 and above, COPT can be applied directly. After installing and configuring the COPT, users can proceed with the following steps:

\begin{sphinxVerbatim}[commandchars=\\\{\}]
\PYG{k+kn}{from}\PYG{+w}{ }\PYG{n+nn}{pulp}\PYG{+w}{ }\PYG{k+kn}{import} \PYG{o}{*}
\end{sphinxVerbatim}

\sphinxAtStartPar
In the solving function \sphinxcode{\sphinxupquote{solve}}, specify the solver to COPT to solve:

\begin{sphinxVerbatim}[commandchars=\\\{\}]
\PYG{n}{solver} \PYG{o}{=} \PYG{n}{COPT}\PYG{p}{(}\PYG{p}{)}
\PYG{n}{result} \PYG{o}{=} \PYG{n}{prob}\PYG{o}{.}\PYG{n}{solve}\PYG{p}{(}\PYG{n}{solver}\PYG{p}{)}
\end{sphinxVerbatim}

\begin{sphinxadmonition}{note}{Notes:}
\begin{enumerate}
\sphinxsetlistlabels{\arabic}{enumi}{enumii}{}{.}%
\item {} 
\sphinxAtStartPar
Calling COPT in PuLP via \sphinxcode{\sphinxupquote{solver = COPT()}} relies on the Python interface of COPT so \sphinxcode{\sphinxupquote{coptpy}} needs to be installed first.

\item {} 
\sphinxAtStartPar
By specifying: \sphinxcode{\sphinxupquote{solver = COPT\_DLL()}} , \sphinxcode{\sphinxupquote{solver = COPT\_CMD()}} , you can also call COPT, depending on COPT installation package.

\end{enumerate}
\end{sphinxadmonition}

\subsection{Introduction of features}
\label{\detokenize{pulpinterface:introduction-of-features}}
\sphinxAtStartPar
The PuLP interface of COPT provides two methods: command line and dynamic library,
which are introduced as follows:

\subsubsection{Command line}
\label{\detokenize{pulpinterface:command-line}}
\sphinxAtStartPar
The command\sphinxhyphen{}line method actually calls the interactive tool \sphinxcode{\sphinxupquote{copt\_cmd}} of COPT
to solve problems. In this way, PuLP generates the MPS format file
corresponding to the model, and combines the parameter settings passed by the user
to generate the solving commands. Upon finish of solving, COPT writes and reads
the result file, and assigns values to the corresponding variables and return them to PuLP.

\sphinxAtStartPar
Functions of the command line method are encapsulated as class \sphinxcode{\sphinxupquote{COPT\_CMD}}.
Users can set parameters when creating the object of the class and
provides the following parameters:
\begin{itemize}
\item {} 
\sphinxAtStartPar
\sphinxcode{\sphinxupquote{keepFiles}}

\sphinxAtStartPar
This option controls whether to keep the generated temporary files.
The default value is \sphinxcode{\sphinxupquote{0}}, which means no temporary files are kept.

\item {} 
\sphinxAtStartPar
\sphinxcode{\sphinxupquote{mip}}

\sphinxAtStartPar
This option controls whether to support solving integer programming models.
The default value is \sphinxcode{\sphinxupquote{True}}, which means support solving integer programming models.

\item {} 
\sphinxAtStartPar
\sphinxcode{\sphinxupquote{msg}}

\sphinxAtStartPar
This option controls whether to print log information on the screen.
The default value is \sphinxcode{\sphinxupquote{True}}, that is, print log information.

\item {} 
\sphinxAtStartPar
\sphinxcode{\sphinxupquote{mip\_start}}

\sphinxAtStartPar
This option controls whether to use initial solution information for integer programming models.
The default value is \sphinxcode{\sphinxupquote{False}}, that is, the initial solution information will not be used.

\item {} 
\sphinxAtStartPar
\sphinxcode{\sphinxupquote{logfile}}

\sphinxAtStartPar
This option specifies the solver log. The default value is \sphinxcode{\sphinxupquote{None}},
which means no solver log will be generated.

\item {} 
\sphinxAtStartPar
\sphinxcode{\sphinxupquote{params}}

\sphinxAtStartPar
This option sets optimization parameters in the form of \sphinxcode{\sphinxupquote{key=value}}.
Please refer to the chapter {\hyperref[\detokenize{capiref:chapapi-param}]{\sphinxcrossref{\DUrole{std,std-ref}{Parameters}}}} for currently supported parameters.

\end{itemize}

\subsubsection{Dynamic library}
\label{\detokenize{pulpinterface:dynamic-library}}
\sphinxAtStartPar
The dynamic library method directly calls COPT C APIs to solve problems.
In this way, PuLP generates problem data and call COPT APIs to load the problem
and parameters set by the user. When optimization finishes, the solution is
obtained by calling COPT APIs, and then assigned to the corresponding variables
and constraints, and passed back to PuLP.

\sphinxAtStartPar
Functions of the dynamic library method are encapsulated as class \sphinxcode{\sphinxupquote{COPT\_DLL}}.
Users can set parameters when creating the object of the class and
provides the following parameters:
\begin{itemize}
\item {} 
\sphinxAtStartPar
\sphinxcode{\sphinxupquote{mip}}

\sphinxAtStartPar
This option controls whether to support solving integer programming models.
The default value is \sphinxcode{\sphinxupquote{True}}, which means support solving integer programming models.

\item {} 
\sphinxAtStartPar
\sphinxcode{\sphinxupquote{msg}}

\sphinxAtStartPar
This option controls whether to print log information on the screen.
The default value is \sphinxcode{\sphinxupquote{True}}, that is, print log information.

\item {} 
\sphinxAtStartPar
\sphinxcode{\sphinxupquote{mip\_start}}

\sphinxAtStartPar
This option controls whether to use initial solution information for integer programming models.
The default value is \sphinxcode{\sphinxupquote{False}}, that is, the initial solution information will not used.

\item {} 
\sphinxAtStartPar
\sphinxcode{\sphinxupquote{logfile}}

\sphinxAtStartPar
This option specifies the solver log. The default value is \sphinxcode{\sphinxupquote{None}},
which means no solver log will be generated.

\item {} 
\sphinxAtStartPar
\sphinxcode{\sphinxupquote{params}}

\sphinxAtStartPar
This option sets optimization parameters in the form of \sphinxcode{\sphinxupquote{key=value}}.
Please refer to the chapter {\hyperref[\detokenize{capiref:chapapi-param}]{\sphinxcrossref{\DUrole{std,std-ref}{Parameters}}}} for currently supported parameters.

\end{itemize}

\sphinxAtStartPar
In addition, the following methods are provided:
\begin{itemize}
\item {} 
\sphinxAtStartPar
\sphinxcode{\sphinxupquote{setParam(self, name, val)}}

\sphinxAtStartPar
Set optimization solution parameters.

\item {} 
\sphinxAtStartPar
\sphinxcode{\sphinxupquote{getParam(self, name)}}

\sphinxAtStartPar
Obtain optimized solution parameters.

\item {} 
\sphinxAtStartPar
\sphinxcode{\sphinxupquote{getAttr(self, name)}}

\sphinxAtStartPar
Get the attribute information of the model.

\item {} 
\sphinxAtStartPar
\sphinxcode{\sphinxupquote{write(self, filename)}}

\sphinxAtStartPar
Output MPS/LP format model file, COPT binary format model file, result file,
basic solution file, initial solution file and parameter setting file.

\end{itemize}

\sphinxstepscope

\section{CVXPY Interface}
\label{\detokenize{cvxpyinterface:cvxpy-interface}}\label{\detokenize{cvxpyinterface:chapcvxpyinterface}}\label{\detokenize{cvxpyinterface::doc}}
\sphinxAtStartPar
\sphinxhref{https://www.cvxpy.org/}{CVXPY} is an open source Python\sphinxhyphen{}embedded modeling language
for convex optimization problems. It lets you express your problem in a natural way
that follows the math, which is quite flexible and efficient. This chapter introduces
how to use the Cardinal Optimizer (COPT) in CVXPY.

\subsection{Installation guide}
\label{\detokenize{cvxpyinterface:installation-guide}}
\sphinxAtStartPar
Before calling COPT in CVXPY to solve problem, users need to setup CVXPY and COPT correctly.
CVXPY currently supports Python 3.7 and later versions of Python. Users can download Python from
\sphinxhref{https://www.anaconda.com/distribution/}{Anaconda distribution} or
\sphinxhref{https://www.python.org/}{Python official distribution} .
We recommend users install the Anaconda distribution, because it is more user\sphinxhyphen{}friendly
and convenient for Python novices.

\subsubsection{Install via conda}
\label{\detokenize{cvxpyinterface:install-via-conda}}
\sphinxAtStartPar
We recommend that users who have installed the Anaconda distribution of Python
use its own \sphinxcode{\sphinxupquote{conda}} tool to install CVXPY. Execute the following commands in
Windows command prompt or terminal on Linux and MacOS:

\begin{sphinxVerbatim}[commandchars=\\\{\}]
conda\PYG{+w}{ }install\PYG{+w}{ }\PYGZhy{}c\PYG{+w}{ }conda\PYGZhy{}forge\PYG{+w}{ }cvxpy
\end{sphinxVerbatim}

\subsubsection{Install via pip}
\label{\detokenize{cvxpyinterface:install-via-pip}}
\sphinxAtStartPar
Users can also install CVXPY through the standard \sphinxcode{\sphinxupquote{pip}} tool, execute the following
command in Windows command prompt or Linux and MacOS terminal:

\begin{sphinxVerbatim}[commandchars=\\\{\}]
pip\PYG{+w}{ }install\PYG{+w}{ }cvxpy
\end{sphinxVerbatim}

\subsubsection{Setup CVXPY interface}
\label{\detokenize{cvxpyinterface:setup-cvxpy-interface}}
\sphinxAtStartPar
\sphinxstylestrong{CVXPY V1.3} and its above versions support calling COPT directly. users need to install and configure COPT in advance and then:

\begin{sphinxVerbatim}[commandchars=\\\{\}]
\PYG{k+kn}{import}\PYG{+w}{ }\PYG{n+nn}{cvxpy}\PYG{+w}{ }\PYG{k}{as}\PYG{+w}{ }\PYG{n+nn}{cp}
\end{sphinxVerbatim}

\sphinxAtStartPar
In CVXPY’s solving function \sphinxcode{\sphinxupquote{solve}}, specify the parameter \sphinxcode{\sphinxupquote{solver="COPT"}} to use the COPT solver to solve:

\begin{sphinxVerbatim}[commandchars=\\\{\}]
\PYG{n}{prob}\PYG{o}{.}\PYG{n}{solve}\PYG{p}{(}\PYG{n}{solver}\PYG{o}{=}\PYG{l+s+s2}{\PYGZdq{}}\PYG{l+s+s2}{COPT}\PYG{l+s+s2}{\PYGZdq{}}\PYG{p}{)}
\end{sphinxVerbatim}

\subsection{Introduction of features}
\label{\detokenize{cvxpyinterface:introduction-of-features}}
\sphinxAtStartPar
The CVXPY interface of COPT supports Linear Programming (LP), Mixed Integer Programming (MIP),
Convex Quadratic Programming (QP), Second\sphinxhyphen{}Order\sphinxhyphen{}Cone Programming (SOCP), Semi\sphinxhyphen{}definite Programming  (SDP), Mixed Integer Convex Quadratic Programming (MIQP) and Mixed Integer Second\sphinxhyphen{}Order\sphinxhyphen{}Cone Programming (MISOCP), common used parameters are:
\begin{itemize}
\item {} 
\sphinxAtStartPar
\sphinxcode{\sphinxupquote{verbose}}

\sphinxAtStartPar
CVXPY builtin parameter, which controls whether to display solving log to the screen.
The default value is \sphinxcode{\sphinxupquote{False}}, which means no log to be displayed.

\item {} 
\sphinxAtStartPar
\sphinxcode{\sphinxupquote{params}}

\sphinxAtStartPar
This option sets optimization parameters in the form of \sphinxcode{\sphinxupquote{key=value}}.
Please refer to the chapter {\hyperref[\detokenize{capiref:chapapi-param}]{\sphinxcrossref{\DUrole{std,std-ref}{Parameters}}}} for currently supported parameters.

\end{itemize}

\sphinxstepscope

\chapter{General Constants}
\label{\detokenize{constant:general-constants}}\label{\detokenize{constant:chapconst}}\label{\detokenize{constant::doc}}
\sphinxAtStartPar
There are three types of constants.
\begin{enumerate}
\sphinxsetlistlabels{\arabic}{enumi}{enumii}{}{.}%
\item {} 
\sphinxAtStartPar
Constructing models, such as optimization directions, constraint senses or variable types.

\item {} 
\sphinxAtStartPar
Accessing solution results, such as API return code, basis status and LP status.

\item {} 
\sphinxAtStartPar
Monitoring optimization progress, such as callback context.

\end{enumerate}

\section{Version information}
\label{\detokenize{constant:version-information}}\begin{itemize}
\item {} 
\sphinxAtStartPar
\sphinxcode{\sphinxupquote{VERSION\_MAJOR}}
\begin{quote}

\sphinxAtStartPar
The major version number.
\end{quote}

\item {} 
\sphinxAtStartPar
\sphinxcode{\sphinxupquote{VERSION\_MINOR}}
\begin{quote}

\sphinxAtStartPar
The minor version number.
\end{quote}

\item {} 
\sphinxAtStartPar
\sphinxcode{\sphinxupquote{VERSION\_TECHNICAL}}
\begin{quote}

\sphinxAtStartPar
The technical version number.
\end{quote}

\end{itemize}

\section{Optimization directions}
\label{\detokenize{constant:optimization-directions}}\label{\detokenize{constant:chapconst-sense}}
\sphinxAtStartPar
For different optimization scenarios, it may be required
to either maximize or minimize the objective function.
There are two optimization directions:
\begin{itemize}
\item {} 
\sphinxAtStartPar
\sphinxcode{\sphinxupquote{MINIMIZE}}
\begin{quote}

\sphinxAtStartPar
For minimizing the objective function.
\end{quote}

\item {} 
\sphinxAtStartPar
\sphinxcode{\sphinxupquote{MAXIMIZE}}
\begin{quote}

\sphinxAtStartPar
For maximizing the objective function.
\end{quote}

\end{itemize}

\sphinxAtStartPar
The optimization direction is automatically
set when reading in a problem from file.
Besides, COPT provides relevant functions, allowing user to explicitly set.
Functions for different APIs are listed below:

\begin{savenotes}\sphinxattablestart
\sphinxthistablewithglobalstyle
\centering
\sphinxcapstartof{table}
\sphinxthecaptionisattop
\sphinxcaption{Functions for setting optimization directions}\label{\detokenize{constant:copttab-const-objsense}}
\sphinxaftertopcaption
\begin{tabular}[t]{|\X{5}{15}|\X{10}{15}|}
\sphinxtoprule
\sphinxstyletheadfamily 
\sphinxAtStartPar
Programming API
&\sphinxstyletheadfamily 
\sphinxAtStartPar
Function
\\
\sphinxmidrule
\sphinxtableatstartofbodyhook
\sphinxAtStartPar
C
&
\sphinxAtStartPar
\sphinxcode{\sphinxupquote{COPT\_SetObjSense}}
\\
\sphinxhline
\sphinxAtStartPar
C++
&
\sphinxAtStartPar
\sphinxcode{\sphinxupquote{Model::SetObjSense()}}
\\
\sphinxhline
\sphinxAtStartPar
C\#
&
\sphinxAtStartPar
\sphinxcode{\sphinxupquote{Model.SetObjSense()}}
\\
\sphinxhline
\sphinxAtStartPar
Java
&
\sphinxAtStartPar
\sphinxcode{\sphinxupquote{Model.setObjSense()}}
\\
\sphinxhline
\sphinxAtStartPar
Python
&
\sphinxAtStartPar
\sphinxcode{\sphinxupquote{Model.setObjSense()}}
\\
\sphinxbottomrule
\end{tabular}
\sphinxtableafterendhook\par
\sphinxattableend\end{savenotes}

\sphinxAtStartPar
\sphinxstylestrong{NOTE:} The function names, calling methods and parameter names in different programming interfaces are slightly different. For details, please refer to the API parameters of each programming language.

\section{Infinity and undefined value}
\label{\detokenize{constant:infinity-and-undefined-value}}
\sphinxAtStartPar
\sphinxstylestrong{Infinity}

\sphinxAtStartPar
In COPT, the infinite bound is represented by a large value,
whose default value is also available as a constant:
\begin{itemize}
\item {} 
\sphinxAtStartPar
\sphinxcode{\sphinxupquote{INFINITY}}
\begin{quote}

\sphinxAtStartPar
The default value (\sphinxcode{\sphinxupquote{1e30}}) of the infinite bound.
\end{quote}

\end{itemize}

\sphinxAtStartPar
\sphinxstylestrong{Undefined Value}

\sphinxAtStartPar
In COPT, the undefined value is represented by another large value.
For example, the default solution value of MIP start is set to
a constant:
\begin{itemize}
\item {} 
\sphinxAtStartPar
\sphinxcode{\sphinxupquote{UNDEFINED}}
\begin{quote}

\sphinxAtStartPar
Undefined value(\sphinxcode{\sphinxupquote{1e40\textasciigrave{}}}).
\end{quote}

\end{itemize}

\section{Constraint senses}
\label{\detokenize{constant:constraint-senses}}\label{\detokenize{constant:chapconst-constrtype}}
\sphinxAtStartPar
Traditionally, for optimization models, constraints are defined using \sphinxstylestrong{senses}.
The most common constraint senses are:
\begin{itemize}
\item {} 
\sphinxAtStartPar
\sphinxcode{\sphinxupquote{LESS\_EQUAL}}
\begin{quote}

\sphinxAtStartPar
For constraint in the form of \(g(x) \leq b\)
\end{quote}

\item {} 
\sphinxAtStartPar
\sphinxcode{\sphinxupquote{GREATER\_EQUAL}}
\begin{quote}

\sphinxAtStartPar
For constraint in the form of \(g(x) \geq b\)
\end{quote}

\item {} 
\sphinxAtStartPar
\sphinxcode{\sphinxupquote{EQUAL}}
\begin{quote}

\sphinxAtStartPar
For constraint in the form of \(g(x) = b\)
\end{quote}

\end{itemize}

\sphinxAtStartPar
In additional, there are two less used constraint senses:
\begin{itemize}
\item {} 
\sphinxAtStartPar
\sphinxcode{\sphinxupquote{FREE}}
\begin{quote}

\sphinxAtStartPar
For unconstrained expression
\end{quote}

\item {} 
\sphinxAtStartPar
\sphinxcode{\sphinxupquote{RANGE}}
\begin{quote}

\sphinxAtStartPar
For constraints with both lower and upper bounds in the form of
\(l \leq g(x) \leq u\).
\end{quote}

\end{itemize}

\sphinxAtStartPar
\sphinxstylestrong{NOTE:} Using constraint senses is supported by COPT but not recommended.
We recommend defining constraints using explicit lower and upper bounds.

\section{Variable types}
\label{\detokenize{constant:variable-types}}\label{\detokenize{constant:chapconst-vartype}}
\sphinxAtStartPar
Variable types are used for defining whether
a variable is continuous or integral.
\begin{itemize}
\item {} 
\sphinxAtStartPar
\sphinxcode{\sphinxupquote{CONTINUOUS}}
\begin{quote}

\sphinxAtStartPar
Non\sphinxhyphen{}integer continuous variables
\end{quote}

\item {} 
\sphinxAtStartPar
\sphinxcode{\sphinxupquote{BINARY}}
\begin{quote}

\sphinxAtStartPar
Binary variables
\end{quote}

\item {} 
\sphinxAtStartPar
\sphinxcode{\sphinxupquote{INTEGER}}
\begin{quote}

\sphinxAtStartPar
Integer variables
\end{quote}

\end{itemize}

\section{SOS constraint types}
\label{\detokenize{constant:sos-constraint-types}}\label{\detokenize{constant:chapconst-sostype}}
\sphinxAtStartPar
SOS constraint (Special Ordered Set) is a kind of special constraint that
places restrictions on the values that a set of variables can take. COPT currently support two types of SOS constraints:
\begin{itemize}
\item {} 
\sphinxAtStartPar
\sphinxcode{\sphinxupquote{SOS\_TYPE1}}
\begin{quote}

\sphinxAtStartPar
SOS1 constraint

\sphinxAtStartPar
At most one variable in the constraint is allowed to take a non\sphinxhyphen{}zero value.
\end{quote}

\item {} 
\sphinxAtStartPar
\sphinxcode{\sphinxupquote{SOS\_TYPE2}}
\begin{quote}

\sphinxAtStartPar
SOS2 constraint

\sphinxAtStartPar
At most two variables in the constraint are allowed to take non\sphinxhyphen{}zero value, and those non\sphinxhyphen{}zero variables must be contiguous.
\end{quote}

\end{itemize}

\sphinxAtStartPar
\sphinxstylestrong{NOTE:} Variables in SOS constraints are allowed to be continuous, binary and integer.

\section{Indicator constraint types}
\label{\detokenize{constant:indicator-constraint-types}}\label{\detokenize{constant:chapconst-indicatortype}}
\sphinxAtStartPar
Indicator constraint is a kind of logical constraints in COPT, used to describe the relationship between the value of
the binary variable \(y\) and whether the linear constraint \(a^{T}x \leq b\) is satisfied.
Currently, COPT supports three types of indicator constraints:
\begin{itemize}
\item {} 
\sphinxAtStartPar
\sphinxcode{\sphinxupquote{INDICATOR\_IF}}

\sphinxAtStartPar
If\sphinxhyphen{}Then:
\begin{quote}

\sphinxAtStartPar
If \(y=f\) , then the linear constraint \(a^{T}x \leq b\) is satisfied.

\sphinxAtStartPar
If \(y\ne f\) , then the linear constraint \(a^{T}x \leq b\) is invalid (may be violated).
\end{quote}

\end{itemize}
\begin{equation}\label{equation:constant:constant:0}
\begin{split}y &= f \rightarrow a^{T}x \leq b\\
f &\in \{0, 1\}\end{split}
\end{equation}\begin{itemize}
\item {} 
\sphinxAtStartPar
\sphinxcode{\sphinxupquote{INDICATOR\_ONLYIF}}

\sphinxAtStartPar
Only\sphinxhyphen{}If:
\begin{quote}

\sphinxAtStartPar
If the linear constraint \(a^{T}x \leq b\) is satisfied, then \(y=f\) .

\sphinxAtStartPar
If the linear constraint \(a^{T}x \leq b\) is not satisfied, then \(y\) can be 0 or 1.
\end{quote}

\end{itemize}
\begin{equation}\label{equation:constant:constant:1}
\begin{split}a^{T}x &\leq b \rightarrow y = f\\
f &\in \{0, 1\}\end{split}
\end{equation}\begin{itemize}
\item {} 
\sphinxAtStartPar
\sphinxcode{\sphinxupquote{INDICATOR\_IFANDONLYIF}}

\sphinxAtStartPar
If\sphinxhyphen{}and\sphinxhyphen{}Only\sphinxhyphen{}If:
\begin{quote}

\sphinxAtStartPar
The linear constraint \(a^{T}x \leq b\) and \(y=f\) must be satisfied simultaneously or not satisfied simultaneously.
\end{quote}

\end{itemize}
\begin{equation}\label{equation:constant:constant:2}
\begin{split}a^{T}x &\leq b \leftrightarrow y = f\\
f &\in \{0, 1\}\end{split}
\end{equation}

\section{SOC constraint types}
\label{\detokenize{constant:soc-constraint-types}}\label{\detokenize{constant:chapconst-conetype}}
\sphinxAtStartPar
The Second\sphinxhyphen{}Order\sphinxhyphen{}Cone (SOC) constraint is a special type of quadratic constraints.
COPT supports two types of SOC constraints:
\begin{itemize}
\item {} 
\sphinxAtStartPar
\sphinxcode{\sphinxupquote{CONE\_QUAD}} : Regular Second\sphinxhyphen{}Order\sphinxhyphen{}Cone

\end{itemize}
\begin{equation}\label{equation:constant:constant:3}
\begin{split}Q^n= \left\{x\in \mathbb{R}^n\ \left|\ x_0\geq\sqrt{\sum_{i=1}^{n-1} x_i^2}, x_0\geq0 \right. \right\}\end{split}
\end{equation}\begin{itemize}
\item {} 
\sphinxAtStartPar
\sphinxcode{\sphinxupquote{CONE\_RQUAD}} : Rotated Second\sphinxhyphen{}Order\sphinxhyphen{}Cone

\end{itemize}
\begin{equation}\label{equation:constant:constant:4}
\begin{split}Q^n_r= \left\{x\in \mathbb{R}^n\ \left|\ 2x_0x_1\geq\sum_{i=2}^{n-1} x_i^2, x_0\geq0, x_1\geq 0 \right. \right\}\end{split}
\end{equation}

\section{Exponential Cone type}
\label{\detokenize{constant:exponential-cone-type}}\label{\detokenize{constant:chapconst-expconetype}}
\sphinxAtStartPar
COPT supports two types of exponential cone contraints:
\begin{itemize}
\item {} 
\sphinxAtStartPar
\sphinxcode{\sphinxupquote{EXPCONE\_PRIMAL}} : Primal exponential cone

\end{itemize}
\begin{equation}\label{equation:constant:constant:5}
\begin{split}\mathrm{cl}(S_1) = S_1 \cup S_2\end{split}
\end{equation}\begin{equation}\label{equation:constant:constant:6}
\begin{split}\begin{aligned}
S_1 &= \left\{\begin{pmatrix} t \\ s \\ r \end{pmatrix}\in \mathbb{R}^3\ |\ s > 0,\ t \geq s\ \mathrm{exp}\left(\frac{r}{s} \right) \right\}, \\
S_2 &= \left\{\begin{pmatrix} t \\ s \\ r \end{pmatrix}\in \mathbb{R}^3\ |\ s=0,\ t\geq 0,\ r\leq 0 \right\}
\end{aligned}\end{split}
\end{equation}\begin{itemize}
\item {} 
\sphinxAtStartPar
\sphinxcode{\sphinxupquote{EXPCONE\_DUAL}} : Dual exponential cone

\end{itemize}
\begin{equation}\label{equation:constant:constant:7}
\begin{split}\mathrm{cl}(S_1) = S_1 \cup S_2\end{split}
\end{equation}\begin{equation}\label{equation:constant:constant:8}
\begin{split}\begin{aligned}
S_1 &= \left\{\begin{pmatrix} t \\ s \\ r \end{pmatrix}\in \mathbb{R}^3\ |\ r < 0,\ t \geq -r\ \mathrm{exp}\left(\frac{s}{r}-1\right) \right\}, \\
S_2 &= \left\{\begin{pmatrix} t \\ s \\ r \end{pmatrix}\in \mathbb{R}^3\ |\ r = 0,\ t\geq 0,\ s\geq 0 \right\}
\end{aligned}\end{split}
\end{equation}

\section{Quadratic objective function}
\label{\detokenize{constant:quadratic-objective-function}}
\sphinxAtStartPar
Besides linear objective function, COPT also supports general convex quadratic objective function.

\sphinxAtStartPar
The mathematical form is:
\begin{equation}\label{equation:constant:coptTab_Eq_qobj}
\begin{split}x^{T}Qx + c^{T}x\end{split}
\end{equation}
\sphinxAtStartPar
Where, \(x\) is an array of variables, \(Q\) is the quadratic part of the quadratic
objective funtion and \(c\) is the linear part.

\section{Quadratic constraint}
\label{\detokenize{constant:quadratic-constraint}}
\sphinxAtStartPar
Besides the special type of quadratic constraint, Second\sphinxhyphen{}Order\sphinxhyphen{}Cone (SOC) constraint, COPT
also supports general convex quadratic constraint.

\sphinxAtStartPar
The mathematical form is:
\begin{equation}\label{equation:constant:coptTab_Eq_qconstr}
\begin{split}x^{T}Qx + q^{T}x \leq b\end{split}
\end{equation}
\sphinxAtStartPar
Where, \(x\) is an array of variables, \(Q\) is the quadratic part of the quadratic
constraint and \(c\) is the linear part.

\section{Basis status}
\label{\detokenize{constant:basis-status}}
\sphinxAtStartPar
For an LP problem with \(n\) variables and \(m\) constraints,
the constraints are treated as slack variables \sphinxstyleemphasis{internally},
resulting in \(n+m\) variables.
When solving an LP problem using the simplex method,
the simplex method fixes \(n\) variables at one of their bounds,
and then computes solutions for the other \(m\) variables.
The \(m\) variables with computed solution are called \sphinxstyleemphasis{basic} variables,
and the other \(n\) variables are called \sphinxstyleemphasis{non\sphinxhyphen{}basic} variables.
The simplex progress and its final solution can be defined using
the basis status of all the variables and constraints.
The basis status supported by COPT are listed below:

\begin{savenotes}\sphinxattablestart
\sphinxthistablewithglobalstyle
\centering
\sphinxcapstartof{table}
\sphinxthecaptionisattop
\sphinxcaption{Basis status values and descriptions}\label{\detokenize{constant:copttab-basiscodes}}
\sphinxaftertopcaption
\begin{tabular}[t]{|\X{10}{37}|\X{5}{37}|\X{22}{37}|}
\sphinxtoprule
\sphinxstyletheadfamily 
\sphinxAtStartPar
Basis status codes
&\sphinxstyletheadfamily 
\sphinxAtStartPar
Value
&\sphinxstyletheadfamily 
\sphinxAtStartPar
Description
\\
\sphinxmidrule
\sphinxtableatstartofbodyhook
\sphinxAtStartPar
\sphinxcode{\sphinxupquote{BASIS\_LOWER}}
&
\sphinxAtStartPar
0
&
\sphinxAtStartPar
The variable is non\sphinxhyphen{}basic at its lower bound
\\
\sphinxhline
\sphinxAtStartPar
\sphinxcode{\sphinxupquote{BASIS\_BASIC}}
&
\sphinxAtStartPar
1
&
\sphinxAtStartPar
The variable is basic
\\
\sphinxhline
\sphinxAtStartPar
\sphinxcode{\sphinxupquote{BASIS\_UPPER}}
&
\sphinxAtStartPar
2
&
\sphinxAtStartPar
The variable is non\sphinxhyphen{}basic at its upper bound
\\
\sphinxhline
\sphinxAtStartPar
\sphinxcode{\sphinxupquote{BASIS\_SUPERBASIC}}
&
\sphinxAtStartPar
3
&
\sphinxAtStartPar
The variable is non\sphinxhyphen{}basic but not any of its bounds
\\
\sphinxhline
\sphinxAtStartPar
\sphinxcode{\sphinxupquote{BASIS\_FIXED}}
&
\sphinxAtStartPar
4
&
\sphinxAtStartPar
The variable is non\sphinxhyphen{}basic and fixed at its bound
\\
\sphinxbottomrule
\end{tabular}
\sphinxtableafterendhook\par
\sphinxattableend\end{savenotes}

\section{Solution status}
\label{\detokenize{constant:solution-status}}\label{\detokenize{constant:chapconst-solstatus}}
\sphinxAtStartPar
Possible solution status values are listed below:

\begin{savenotes}\sphinxattablestart
\sphinxthistablewithglobalstyle
\centering
\sphinxcapstartof{table}
\sphinxthecaptionisattop
\sphinxcaption{Solution Status}\label{\detokenize{constant:copttab-statuscodes}}
\sphinxaftertopcaption
\begin{tabular}[t]{|\X{15}{65}|\X{5}{65}|\X{45}{65}|}
\sphinxtoprule
\sphinxstyletheadfamily 
\sphinxAtStartPar
Status Codes
&\sphinxstyletheadfamily 
\sphinxAtStartPar
Value
&\sphinxstyletheadfamily 
\sphinxAtStartPar
Description
\\
\sphinxmidrule
\sphinxtableatstartofbodyhook
\sphinxAtStartPar
\sphinxcode{\sphinxupquote{UNSTARTED}}
&
\sphinxAtStartPar
0
&
\sphinxAtStartPar
The solving process has not been started yet
\\
\sphinxhline
\sphinxAtStartPar
\sphinxcode{\sphinxupquote{OPTIMAL}}
&
\sphinxAtStartPar
1
&
\sphinxAtStartPar
The optimal solutions are found
\\
\sphinxhline
\sphinxAtStartPar
\sphinxcode{\sphinxupquote{INFEASIBLE}}
&
\sphinxAtStartPar
2
&
\sphinxAtStartPar
The model is infeasible
\\
\sphinxhline
\sphinxAtStartPar
\sphinxcode{\sphinxupquote{UNBOUNDED}}
&
\sphinxAtStartPar
3
&
\sphinxAtStartPar
The objective is unbounded
\\
\sphinxhline
\sphinxAtStartPar
\sphinxcode{\sphinxupquote{INF\_OR\_UNB}}
&
\sphinxAtStartPar
4
&
\sphinxAtStartPar
The model is infeasible or unbounded
\\
\sphinxhline
\sphinxAtStartPar
\sphinxcode{\sphinxupquote{NUMERICAL}}
&
\sphinxAtStartPar
5
&
\sphinxAtStartPar
The solver encountered numerical difficulties
\\
\sphinxhline
\sphinxAtStartPar
\sphinxcode{\sphinxupquote{NODELIMIT}}
&
\sphinxAtStartPar
6
&
\sphinxAtStartPar
The solving process was not finished within node limit
\\
\sphinxhline
\sphinxAtStartPar
\sphinxcode{\sphinxupquote{IMPRECISE}}
&
\sphinxAtStartPar
7
&
\sphinxAtStartPar
The solution is imprecise
\\
\sphinxhline
\sphinxAtStartPar
\sphinxcode{\sphinxupquote{TIMEOUT}}
&
\sphinxAtStartPar
8
&
\sphinxAtStartPar
The solving process was not finished within time limit
\\
\sphinxhline
\sphinxAtStartPar
\sphinxcode{\sphinxupquote{UNFINISHED}}
&
\sphinxAtStartPar
9
&
\sphinxAtStartPar
The solving process was not finished due to an internal error
\\
\sphinxhline
\sphinxAtStartPar
\sphinxcode{\sphinxupquote{INTERRUPTED}}
&
\sphinxAtStartPar
10
&
\sphinxAtStartPar
The solving process was stopped by user interruption
\\
\sphinxhline
\sphinxAtStartPar
\sphinxcode{\sphinxupquote{ITERLIMIT}}
&
\sphinxAtStartPar
11
&
\sphinxAtStartPar
The solving process was not finished within iteration limit
\\
\sphinxbottomrule
\end{tabular}
\sphinxtableafterendhook\par
\sphinxattableend\end{savenotes}

\begin{sphinxadmonition}{note}{Notes}
\begin{itemize}
\item {} 
\sphinxAtStartPar
In the Python API, solution status are defined in COPT’s General Constants Class. They can be accessed via the \sphinxcode{\sphinxupquote{"COPT."}} prefix or \sphinxcode{\sphinxupquote{Model.status}} ;

\item {} 
\sphinxAtStartPar
In the Constant class of the Java API and C\# API, the constants about the solution status are defined in the \sphinxcode{\sphinxupquote{Status}} class;

\item {} 
\sphinxAtStartPar
The linear programming solution status can be obtained through the attribute \sphinxcode{\sphinxupquote{"LpStatus"}} , and the integer programming solution status can be obtained through the attribute \sphinxcode{\sphinxupquote{"MipStatus"}} .

\item {} 
\sphinxAtStartPar
The LP\sphinxhyphen{}relaxation status of the current node can be obtained through the Callback information
\sphinxcode{\sphinxupquote{"NodeStatus"}} . The return value is as above, except for \sphinxcode{\sphinxupquote{NODELIMIT}}, \sphinxcode{\sphinxupquote{UNSTARTED}},   \sphinxcode{\sphinxupquote{INF\_OR\_UNB}} .

\end{itemize}
\end{sphinxadmonition}

\section{Client configuration}
\label{\detokenize{constant:client-configuration}}\label{\detokenize{constant:chapconst-client}}
\sphinxAtStartPar
For floating, cluster and Web License, the user can set client configuration parameters by calling API functions.
The available parameters are:
\begin{itemize}
\item {} 
\sphinxAtStartPar
\sphinxcode{\sphinxupquote{CLIENT\_CAFILE}}
\begin{quote}

\sphinxAtStartPar
Path to the CA public certificate file for configuring SSL proxy.
\end{quote}

\item {} 
\sphinxAtStartPar
\sphinxcode{\sphinxupquote{CLIENT\_CERTFILE}}
\begin{quote}

\sphinxAtStartPar
Path to the client certificate file for configuring SSL proxy.
\end{quote}

\item {} 
\sphinxAtStartPar
\sphinxcode{\sphinxupquote{CLIENT\_CERTKEYFILE}}
\begin{quote}

\sphinxAtStartPar
Path to the client private key file for configuring SSL proxy.
\end{quote}

\item {} 
\sphinxAtStartPar
\sphinxcode{\sphinxupquote{CLIENT\_CLUSTER}}
\begin{quote}

\sphinxAtStartPar
IP address of the remote cluster server.
\end{quote}

\item {} 
\sphinxAtStartPar
\sphinxcode{\sphinxupquote{CLIENT\_FLOATING}}
\begin{quote}

\sphinxAtStartPar
IP address of the remote floating token server.
\end{quote}

\item {} 
\sphinxAtStartPar
\sphinxcode{\sphinxupquote{CLIENT\_PASSWORD}}
\begin{quote}

\sphinxAtStartPar
Password for accessing the remote server.
\end{quote}

\item {} 
\sphinxAtStartPar
\sphinxcode{\sphinxupquote{CLIENT\_PORT}}
\begin{quote}

\sphinxAtStartPar
Connection port of the remote server.
\end{quote}

\item {} 
\sphinxAtStartPar
\sphinxcode{\sphinxupquote{CLIENT\_PRIORITY}}
\begin{quote}

\sphinxAtStartPar
Task priority for the remote cluster server.
Possible values range from 0 to 99, with higher values indicating higher priority.
\end{quote}

\item {} 
\sphinxAtStartPar
\sphinxcode{\sphinxupquote{CLIENT\_WAITTIME}}
\begin{quote}

\sphinxAtStartPar
Wait time of the client connection .
\end{quote}

\item {} 
\sphinxAtStartPar
\sphinxcode{\sphinxupquote{CLIENT\_WEBSERVER}}
\begin{quote}

\sphinxAtStartPar
Domain name of the COPT Web License server.
\end{quote}

\item {} 
\sphinxAtStartPar
\sphinxcode{\sphinxupquote{CLIENT\_WEBLICENSEID}}
\begin{quote}

\sphinxAtStartPar
Unique identifier for the Web License.
\end{quote}

\item {} 
\sphinxAtStartPar
\sphinxcode{\sphinxupquote{CLIENT\_WEBACCESSKEY}}
\begin{quote}

\sphinxAtStartPar
Each Web License can have multiple API keys, with each API key having a unique \sphinxcode{\sphinxupquote{WebAccessKey}}.
\end{quote}

\item {} 
\sphinxAtStartPar
\sphinxcode{\sphinxupquote{CLIENT\_WEBTOKENDURATION}}
\begin{quote}

\sphinxAtStartPar
Validity period for each authorization from the Web License.
\end{quote}

\end{itemize}

\begin{sphinxadmonition}{note}{Notes:}

\sphinxAtStartPar
The above client configuration related parameters should be set within the API function \sphinxcode{\sphinxupquote{EnvrConfig}} .
\end{sphinxadmonition}

\section{Callback context}
\label{\detokenize{constant:callback-context}}\label{\detokenize{constant:chapconst-cbc}}\begin{itemize}
\item {} 
\sphinxAtStartPar
\sphinxcode{\sphinxupquote{CBCONTEXT\_INCUMBENT}}
\begin{quote}

\sphinxAtStartPar
Invokes the callback after a new incumbent was found.
\end{quote}

\item {} 
\sphinxAtStartPar
\sphinxcode{\sphinxupquote{CBCONTEXT\_MIPRELAX}}
\begin{quote}

\sphinxAtStartPar
Invokes the callback when a new LP\sphinxhyphen{}relaxation solution is found.
\end{quote}

\item {} 
\sphinxAtStartPar
\sphinxcode{\sphinxupquote{CBCONTEXT\_MIPSOL}}
\begin{quote}

\sphinxAtStartPar
Invokes the callback when a new MIP candidate solution is found.
\end{quote}

\item {} 
\sphinxAtStartPar
\sphinxcode{\sphinxupquote{CBCONTEXT\_MIPNODE}}
\begin{quote}

\sphinxAtStartPar
Invokes the callback when a MIP node is finished and LP\sphinxhyphen{}relaxation has been solved.
\end{quote}

\end{itemize}

\section{Nonlinear Model Callback Evaluation Types}
\label{\detokenize{constant:nonlinear-model-callback-evaluation-types}}\label{\detokenize{constant:chapconst-nlpevaltype}}\begin{itemize}
\item {} 
\sphinxAtStartPar
\sphinxcode{\sphinxupquote{EVALTYPE\_OBJVAL}}
\begin{quote}

\sphinxAtStartPar
Indicates that the callback object can provide the objective function value.
\end{quote}

\item {} 
\sphinxAtStartPar
\sphinxcode{\sphinxupquote{EVALTYPE\_CONSTRVAL}}
\begin{quote}

\sphinxAtStartPar
Indicates that the callback object can provide constraint function values.
\end{quote}

\item {} 
\sphinxAtStartPar
\sphinxcode{\sphinxupquote{EVALTYPE\_GRADIENT}}
\begin{quote}

\sphinxAtStartPar
Indicates that the callback object can provide first\sphinxhyphen{}order derivative information
(gradient) of the objective function.
\end{quote}

\item {} 
\sphinxAtStartPar
\sphinxcode{\sphinxupquote{EVALTYPE\_JACOBIAN}}
\begin{quote}

\sphinxAtStartPar
Indicates that the callback object can provide first\sphinxhyphen{}order derivative information
(Jacobian matrix) of the constraint functions.
\end{quote}

\item {} 
\sphinxAtStartPar
\sphinxcode{\sphinxupquote{EVALTYPE\_HESSIAN}}
\begin{quote}

\sphinxAtStartPar
Indicates that the callback object can provide second\sphinxhyphen{}order derivative information
(Hessian matrix) of the nonlinear model Lagrangian function.
\end{quote}

\end{itemize}

\section{Dense Storage Types for Nonlinear Matrices}
\label{\detokenize{constant:dense-storage-types-for-nonlinear-matrices}}\label{\detokenize{constant:chapconst-nlpdensetype}}\begin{itemize}
\item {} 
\sphinxAtStartPar
\sphinxcode{\sphinxupquote{DENSETYPE\_ROWMAJOR}}
\begin{quote}

\sphinxAtStartPar
Indicates that the corresponding gradient, Jacobian, or Hessian matrix is provided
in dense format and stored in row\sphinxhyphen{}major order.
\end{quote}

\item {} 
\sphinxAtStartPar
\sphinxcode{\sphinxupquote{DENSETYPE\_COLMAJOR}}
\begin{quote}

\sphinxAtStartPar
Indicates that the corresponding gradient, Jacobian, or Hessian matrix is provided
in dense format and stored in column\sphinxhyphen{}major order.
\end{quote}

\end{itemize}

\begin{sphinxadmonition}{note}{Notes:}

\sphinxAtStartPar
When these two flag is used,
the related index array parameters in \sphinxcode{\sphinxupquote{Model.loadNlData}} may be set to \sphinxcode{\sphinxupquote{None}}.
\end{sphinxadmonition}

\section{API function return code}
\label{\detokenize{constant:api-function-return-code}}\label{\detokenize{constant:chapconst-return}}
\sphinxAtStartPar
When an API function finishes, it returns an integer \sphinxstylestrong{return code},
which indicates whether the API call was finished okay or failed.
In case of failure, it specifies the reason of the failure.

\sphinxAtStartPar
Possible COPT API function return codes are listed below:

\begin{savenotes}\sphinxattablestart
\sphinxthistablewithglobalstyle
\centering
\sphinxcapstartof{table}
\sphinxthecaptionisattop
\sphinxcaption{COPT API Function Return Codes}\label{\detokenize{constant:copttab-returncodes}}
\sphinxaftertopcaption
\begin{tabular}[t]{|\X{15}{65}|\X{5}{65}|\X{45}{65}|}
\sphinxtoprule
\sphinxstyletheadfamily 
\sphinxAtStartPar
Return Codes
&\sphinxstyletheadfamily 
\sphinxAtStartPar
Value
&\sphinxstyletheadfamily 
\sphinxAtStartPar
Description
\\
\sphinxmidrule
\sphinxtableatstartofbodyhook
\sphinxAtStartPar
\sphinxcode{\sphinxupquote{OK}}
&
\sphinxAtStartPar
0
&
\sphinxAtStartPar
The API call finished successfully
\\
\sphinxhline
\sphinxAtStartPar
\sphinxcode{\sphinxupquote{MEMORY}}
&
\sphinxAtStartPar
1
&
\sphinxAtStartPar
The API call failed because of memory allocation failure
\\
\sphinxhline
\sphinxAtStartPar
\sphinxcode{\sphinxupquote{FILE}}
&
\sphinxAtStartPar
2
&
\sphinxAtStartPar
The API call failed because of file input or output failure
\\
\sphinxhline
\sphinxAtStartPar
\sphinxcode{\sphinxupquote{INVALID}}
&
\sphinxAtStartPar
3
&
\sphinxAtStartPar
The API call failed because of invalid data
\\
\sphinxhline
\sphinxAtStartPar
\sphinxcode{\sphinxupquote{LICENSE}}
&
\sphinxAtStartPar
4
&
\sphinxAtStartPar
The API call failed because of license validation failure
\\
\sphinxhline
\sphinxAtStartPar
\sphinxcode{\sphinxupquote{INTERNAL}}
&
\sphinxAtStartPar
5
&
\sphinxAtStartPar
The API call failed because of internal error
\\
\sphinxhline
\sphinxAtStartPar
\sphinxcode{\sphinxupquote{THREAD}}
&
\sphinxAtStartPar
6
&
\sphinxAtStartPar
The API call failed because of thread error
\\
\sphinxhline
\sphinxAtStartPar
\sphinxcode{\sphinxupquote{SERVER}}
&
\sphinxAtStartPar
7
&
\sphinxAtStartPar
The API call failed because of remote server error
\\
\sphinxhline
\sphinxAtStartPar
\sphinxcode{\sphinxupquote{NONCONVEX}}
&
\sphinxAtStartPar
8
&
\sphinxAtStartPar
The API call failed because of problem is nonconvex”
\\
\sphinxhline
\sphinxAtStartPar
\sphinxcode{\sphinxupquote{MEMORY\_GPU}}
&
\sphinxAtStartPar
9
&
\sphinxAtStartPar
The API call failed because of GPU memory allocation failure”
\\
\sphinxbottomrule
\end{tabular}
\sphinxtableafterendhook\par
\sphinxattableend\end{savenotes}

\section{Client configuration}
\label{\detokenize{constant:id1}}
\sphinxAtStartPar
For floating and cluster clients, users are allowed to set client configuration parameters,
currently available settings are:
\begin{itemize}
\item {} 
\sphinxAtStartPar
\sphinxcode{\sphinxupquote{COPT\_CLIENT\_CLUSTER}}
\begin{quote}

\sphinxAtStartPar
IP address of cluster server.
\end{quote}

\item {} 
\sphinxAtStartPar
\sphinxcode{\sphinxupquote{COPT\_CLIENT\_FLOATING}}
\begin{quote}

\sphinxAtStartPar
IP address of token server.
\end{quote}

\item {} 
\sphinxAtStartPar
\sphinxcode{\sphinxupquote{COPT\_CLIENT\_PASSWORD}}
\begin{quote}

\sphinxAtStartPar
Password of cluster server.
\end{quote}

\item {} 
\sphinxAtStartPar
\sphinxcode{\sphinxupquote{CLIENT\_PORT}}
\begin{quote}

\sphinxAtStartPar
Connection port of token server.
\end{quote}

\item {} 
\sphinxAtStartPar
\sphinxcode{\sphinxupquote{COPT\_CLIENT\_WAITTIME}}
\begin{quote}

\sphinxAtStartPar
Wait time of client.
\end{quote}

\end{itemize}

\section{Methods for accessing constants}
\label{\detokenize{constant:methods-for-accessing-constants}}\label{\detokenize{constant:chapconst-method}}
\sphinxAtStartPar
In different programming interfaces, the ways of accessing constants are slightly different. In the C language interface,
the constant name is prefixed with \sphinxcode{\sphinxupquote{"COPT\_"}} (like \sphinxcode{\sphinxupquote{COPT\_MAXIMIZE}}). For details, please refer to the corresponding
chapters for each programming language API:
\begin{itemize}
\item {} 
\sphinxAtStartPar
C API: {\hyperref[\detokenize{capiref:chapapi-const}]{\sphinxcrossref{\DUrole{std,std-ref}{C API Reference: Constants}}}}

\item {} 
\sphinxAtStartPar
C++ API: {\hyperref[\detokenize{cppapiref:chapcppapiref-const}]{\sphinxcrossref{\DUrole{std,std-ref}{C++ API Reference: Constants}}}}

\item {} 
\sphinxAtStartPar
C\# API: {\hyperref[\detokenize{csharpapiref:chapcsharpapiref-general}]{\sphinxcrossref{\DUrole{std,std-ref}{C\# API Reference: General Constants}}}}

\item {} 
\sphinxAtStartPar
Java API: {\hyperref[\detokenize{javaapiref:chapjavaapiref-const-general}]{\sphinxcrossref{\DUrole{std,std-ref}{Java API Reference: General Constants}}}}

\item {} 
\sphinxAtStartPar
Python API: {\hyperref[\detokenize{pyapiref:chappyapi-const-general}]{\sphinxcrossref{\DUrole{std,std-ref}{Python API Reference: General Constants}}}}

\end{itemize}

\sphinxstepscope

\chapter{Attributes}
\label{\detokenize{attribute:attributes}}\label{\detokenize{attribute:chapattrs}}\label{\detokenize{attribute::doc}}
\sphinxAtStartPar
To query and modify properties of a COPT model is through the
attribute interface. A variety of different attributes are available,
and they can be associated with solutions, or the model.

\section{Problem related}
\label{\detokenize{attribute:problem-related}}
\sphinxAtStartPar
Problem related attributes provide the relevant information of the model composition and description.
The names and descriptions of these attributes are summarized below.

\begin{savenotes}\sphinxattablestart
\sphinxthistablewithglobalstyle
\centering
\sphinxcapstartof{table}
\sphinxthecaptionisattop
\sphinxcaption{Problem related attributes}\label{\detokenize{attribute:id1}}
\sphinxaftertopcaption
\begin{tabular}[t]{|\X{15}{59}|\X{9}{59}|\X{35}{59}|}
\sphinxtoprule
\sphinxstyletheadfamily 
\sphinxAtStartPar
Name
&\sphinxstyletheadfamily 
\sphinxAtStartPar
Type
&\sphinxstyletheadfamily 
\sphinxAtStartPar
Description
\\
\sphinxmidrule
\sphinxtableatstartofbodyhook
\sphinxAtStartPar
{\hyperref[\detokenize{attribute:cols}]{\sphinxcrossref{\DUrole{std,std-ref}{Cols}}}}
&
\sphinxAtStartPar
Integer
&
\sphinxAtStartPar
Number of variables (columns) in the problem
\\
\sphinxhline
\sphinxAtStartPar
{\hyperref[\detokenize{attribute:psdcols}]{\sphinxcrossref{\DUrole{std,std-ref}{PSDCols}}}}
&
\sphinxAtStartPar
Integer
&
\sphinxAtStartPar
Number of PSD variables in the problem
\\
\sphinxhline
\sphinxAtStartPar
{\hyperref[\detokenize{attribute:rows}]{\sphinxcrossref{\DUrole{std,std-ref}{Rows}}}}
&
\sphinxAtStartPar
Integer
&
\sphinxAtStartPar
Number of constraints (rows) in the problem
\\
\sphinxhline
\sphinxAtStartPar
{\hyperref[\detokenize{attribute:elems}]{\sphinxcrossref{\DUrole{std,std-ref}{Elems}}}}
&
\sphinxAtStartPar
Integer
&
\sphinxAtStartPar
Number of non\sphinxhyphen{}zero elements in the coefficient matrix
\\
\sphinxhline
\sphinxAtStartPar
{\hyperref[\detokenize{attribute:qelems}]{\sphinxcrossref{\DUrole{std,std-ref}{QElems}}}}
&
\sphinxAtStartPar
Integer
&
\sphinxAtStartPar
Number of non\sphinxhyphen{}zero quadratic elements in the quadratic objective function
\\
\sphinxhline
\sphinxAtStartPar
{\hyperref[\detokenize{attribute:psdelems}]{\sphinxcrossref{\DUrole{std,std-ref}{PSDElems}}}}
&
\sphinxAtStartPar
Integer
&
\sphinxAtStartPar
Number of PSD terms in objective function
\\
\sphinxhline
\sphinxAtStartPar
{\hyperref[\detokenize{attribute:symmats}]{\sphinxcrossref{\DUrole{std,std-ref}{SymMats}}}}
&
\sphinxAtStartPar
Integer
&
\sphinxAtStartPar
Number of symmetric matrices in the problem
\\
\sphinxhline
\sphinxAtStartPar
{\hyperref[\detokenize{attribute:bins}]{\sphinxcrossref{\DUrole{std,std-ref}{Bins}}}}
&
\sphinxAtStartPar
Integer
&
\sphinxAtStartPar
Number of binary variables
\\
\sphinxhline
\sphinxAtStartPar
{\hyperref[\detokenize{attribute:ints}]{\sphinxcrossref{\DUrole{std,std-ref}{Ints}}}}
&
\sphinxAtStartPar
Integer
&
\sphinxAtStartPar
Number of integer variables
\\
\sphinxhline
\sphinxAtStartPar
{\hyperref[\detokenize{attribute:soss}]{\sphinxcrossref{\DUrole{std,std-ref}{Soss}}}}
&
\sphinxAtStartPar
Integer
&
\sphinxAtStartPar
Number of SOS constraints
\\
\sphinxhline
\sphinxAtStartPar
{\hyperref[\detokenize{attribute:cones}]{\sphinxcrossref{\DUrole{std,std-ref}{Cones}}}}
&
\sphinxAtStartPar
Integer
&
\sphinxAtStartPar
Number of Second\sphinxhyphen{}Order\sphinxhyphen{}Cone constraints
\\
\sphinxhline
\sphinxAtStartPar
{\hyperref[\detokenize{attribute:expcones}]{\sphinxcrossref{\DUrole{std,std-ref}{ExpCones}}}}
&
\sphinxAtStartPar
Integer
&
\sphinxAtStartPar
Number of exponential cone constraints
\\
\sphinxhline
\sphinxAtStartPar
{\hyperref[\detokenize{attribute:affinecones}]{\sphinxcrossref{\DUrole{std,std-ref}{AffineCones}}}}
&
\sphinxAtStartPar
Integer
&
\sphinxAtStartPar
Number of affine cone constraints
\\
\sphinxhline
\sphinxAtStartPar
{\hyperref[\detokenize{attribute:qconstrs}]{\sphinxcrossref{\DUrole{std,std-ref}{QConstrs}}}}
&
\sphinxAtStartPar
Integer
&
\sphinxAtStartPar
Number of quadratic constraints
\\
\sphinxhline
\sphinxAtStartPar
{\hyperref[\detokenize{attribute:psdconstrs}]{\sphinxcrossref{\DUrole{std,std-ref}{PSDConstrs}}}}
&
\sphinxAtStartPar
Integer
&
\sphinxAtStartPar
Number of PSD constraints
\\
\sphinxhline
\sphinxAtStartPar
{\hyperref[\detokenize{attribute:lmiconstrs}]{\sphinxcrossref{\DUrole{std,std-ref}{LMIConstrs}}}}
&
\sphinxAtStartPar
Integer
&
\sphinxAtStartPar
Number of LMI constraints
\\
\sphinxhline
\sphinxAtStartPar
{\hyperref[\detokenize{attribute:indicators}]{\sphinxcrossref{\DUrole{std,std-ref}{Indicators}}}}
&
\sphinxAtStartPar
Integer
&
\sphinxAtStartPar
Number of indicator constraints
\\
\sphinxhline
\sphinxAtStartPar
{\hyperref[\detokenize{attribute:multiobjs}]{\sphinxcrossref{\DUrole{std,std-ref}{MultiObjs}}}}
&
\sphinxAtStartPar
Integer
&
\sphinxAtStartPar
Number of objectives in a multi\sphinxhyphen{}objective model
\\
\sphinxhline
\sphinxAtStartPar
{\hyperref[\detokenize{attribute:objsense}]{\sphinxcrossref{\DUrole{std,std-ref}{ObjSense}}}}
&
\sphinxAtStartPar
Integer
&
\sphinxAtStartPar
The optimization direction
\\
\sphinxhline
\sphinxAtStartPar
{\hyperref[\detokenize{attribute:objconst}]{\sphinxcrossref{\DUrole{std,std-ref}{ObjConst}}}}
&
\sphinxAtStartPar
Double
&
\sphinxAtStartPar
The constant part of the objective function
\\
\sphinxhline
\sphinxAtStartPar
{\hyperref[\detokenize{attribute:hasqobj}]{\sphinxcrossref{\DUrole{std,std-ref}{HasQObj}}}}
&
\sphinxAtStartPar
Integer
&
\sphinxAtStartPar
Whether the problem has quadratic objective function
\\
\sphinxhline
\sphinxAtStartPar
{\hyperref[\detokenize{attribute:haspsdobj}]{\sphinxcrossref{\DUrole{std,std-ref}{HasPSDObj}}}}
&
\sphinxAtStartPar
Integer
&
\sphinxAtStartPar
Whether the problem has PSD terms in objective function
\\
\sphinxhline
\sphinxAtStartPar
{\hyperref[\detokenize{attribute:hasnlobj}]{\sphinxcrossref{\DUrole{std,std-ref}{HasNLObj}}}}
&
\sphinxAtStartPar
Integer
&
\sphinxAtStartPar
Whether the problem has nonlinear terms in objective function
\\
\sphinxhline
\sphinxAtStartPar
{\hyperref[\detokenize{attribute:ismip}]{\sphinxcrossref{\DUrole{std,std-ref}{IsMIP}}}}
&
\sphinxAtStartPar
Integer
&
\sphinxAtStartPar
Whether the problem is a MIP
\\
\sphinxbottomrule
\end{tabular}
\sphinxtableafterendhook\par
\sphinxattableend\end{savenotes}
\phantomsection\label{\detokenize{attribute:cols}}\begin{itemize}
\item {} 
\sphinxAtStartPar
\sphinxcode{\sphinxupquote{Cols}}
\begin{quote}

\sphinxAtStartPar
Integer attribute.

\sphinxAtStartPar
Number of variables (columns) in the problem.
\end{quote}

\end{itemize}
\phantomsection\label{\detokenize{attribute:psdcols}}\begin{itemize}
\item {} 
\sphinxAtStartPar
\sphinxcode{\sphinxupquote{PSDCols}}
\begin{quote}

\sphinxAtStartPar
Integer attribute.

\sphinxAtStartPar
Number of PSD variables in the problem.
\end{quote}

\end{itemize}
\phantomsection\label{\detokenize{attribute:rows}}\begin{itemize}
\item {} 
\sphinxAtStartPar
\sphinxcode{\sphinxupquote{Rows}}
\begin{quote}

\sphinxAtStartPar
Integer attribute.

\sphinxAtStartPar
Number of constraints (rows) in the problem.
\end{quote}

\end{itemize}
\phantomsection\label{\detokenize{attribute:elems}}\begin{itemize}
\item {} 
\sphinxAtStartPar
\sphinxcode{\sphinxupquote{Elems}}
\begin{quote}

\sphinxAtStartPar
Integer attribute.

\sphinxAtStartPar
Number of non\sphinxhyphen{}zero elements in the coefficient matrix.
\end{quote}

\end{itemize}
\phantomsection\label{\detokenize{attribute:qelems}}\begin{itemize}
\item {} 
\sphinxAtStartPar
\sphinxcode{\sphinxupquote{QElems}}
\begin{quote}

\sphinxAtStartPar
Integer attribute.

\sphinxAtStartPar
Number of non\sphinxhyphen{}zero quadratic elements in the quadratic objective function.
\end{quote}

\end{itemize}
\phantomsection\label{\detokenize{attribute:psdelems}}\begin{itemize}
\item {} 
\sphinxAtStartPar
\sphinxcode{\sphinxupquote{PSDElems}}
\begin{quote}

\sphinxAtStartPar
Integer attribute.

\sphinxAtStartPar
Number of PSD terms in objective function.
\end{quote}

\end{itemize}
\phantomsection\label{\detokenize{attribute:symmats}}\begin{itemize}
\item {} 
\sphinxAtStartPar
\sphinxcode{\sphinxupquote{SymMats}}
\begin{quote}

\sphinxAtStartPar
Integer attribute.

\sphinxAtStartPar
Number of symmetric matrices in the problem.
\end{quote}

\end{itemize}
\phantomsection\label{\detokenize{attribute:bins}}\begin{itemize}
\item {} 
\sphinxAtStartPar
\sphinxcode{\sphinxupquote{Bins}}
\begin{quote}

\sphinxAtStartPar
Integer attribute.

\sphinxAtStartPar
Number of binary variables.
\end{quote}

\end{itemize}
\phantomsection\label{\detokenize{attribute:ints}}\begin{itemize}
\item {} 
\sphinxAtStartPar
\sphinxcode{\sphinxupquote{Ints}}
\begin{quote}

\sphinxAtStartPar
Integer attribute.

\sphinxAtStartPar
Number of integer variables.
\end{quote}

\end{itemize}
\phantomsection\label{\detokenize{attribute:soss}}\begin{itemize}
\item {} 
\sphinxAtStartPar
\sphinxcode{\sphinxupquote{Soss}}
\begin{quote}

\sphinxAtStartPar
Integer attribute.

\sphinxAtStartPar
Number of SOS constraints.
\end{quote}

\end{itemize}
\phantomsection\label{\detokenize{attribute:cones}}\begin{itemize}
\item {} 
\sphinxAtStartPar
\sphinxcode{\sphinxupquote{Cones}}
\begin{quote}

\sphinxAtStartPar
Integer attribute.

\sphinxAtStartPar
Number of Second\sphinxhyphen{}Order\sphinxhyphen{}Cone constraints.
\end{quote}

\end{itemize}
\phantomsection\label{\detokenize{attribute:expcones}}\begin{itemize}
\item {} 
\sphinxAtStartPar
\sphinxcode{\sphinxupquote{ExpCones}}
\begin{quote}

\sphinxAtStartPar
Integer attribute.

\sphinxAtStartPar
Number of exponential cone constraints.
\end{quote}

\end{itemize}
\phantomsection\label{\detokenize{attribute:affinecones}}\begin{itemize}
\item {} 
\sphinxAtStartPar
\sphinxcode{\sphinxupquote{AffineCones}}
\begin{quote}

\sphinxAtStartPar
Integer attribute.

\sphinxAtStartPar
Number of affine cone constraints.
\end{quote}

\end{itemize}
\phantomsection\label{\detokenize{attribute:qconstrs}}\begin{itemize}
\item {} 
\sphinxAtStartPar
\sphinxcode{\sphinxupquote{QConstrs}}
\begin{quote}

\sphinxAtStartPar
Integer attribute.

\sphinxAtStartPar
Number of quadratic constraints.
\end{quote}

\end{itemize}
\phantomsection\label{\detokenize{attribute:psdconstrs}}\begin{itemize}
\item {} 
\sphinxAtStartPar
\sphinxcode{\sphinxupquote{PSDConstrs}}
\begin{quote}

\sphinxAtStartPar
Integer attribute.

\sphinxAtStartPar
Number of PSD constraints.
\end{quote}

\end{itemize}
\phantomsection\label{\detokenize{attribute:lmiconstrs}}\begin{itemize}
\item {} 
\sphinxAtStartPar
\sphinxcode{\sphinxupquote{LMIConstrs}}
\begin{quote}

\sphinxAtStartPar
Integer attribute.

\sphinxAtStartPar
Number of LMI (Linear Matrix Inequalities) constraints.
\end{quote}

\end{itemize}
\phantomsection\label{\detokenize{attribute:indicators}}\begin{itemize}
\item {} 
\sphinxAtStartPar
\sphinxcode{\sphinxupquote{Indicators}}
\begin{quote}

\sphinxAtStartPar
Integer attribute.

\sphinxAtStartPar
Number of indicator constraints.
\end{quote}

\end{itemize}
\phantomsection\label{\detokenize{attribute:multiobjs}}\begin{itemize}
\item {} 
\sphinxAtStartPar
\sphinxcode{\sphinxupquote{MultiObjs}}
\begin{quote}

\sphinxAtStartPar
Integer attribute.

\sphinxAtStartPar
Number of objectives in a multi\sphinxhyphen{}objective model.
\end{quote}

\end{itemize}
\phantomsection\label{\detokenize{attribute:objsense}}\begin{itemize}
\item {} 
\sphinxAtStartPar
\sphinxcode{\sphinxupquote{ObjSense}}
\begin{quote}

\sphinxAtStartPar
Integer attribute.

\sphinxAtStartPar
The optimization direction.
\end{quote}

\end{itemize}
\phantomsection\label{\detokenize{attribute:objconst}}\begin{itemize}
\item {} 
\sphinxAtStartPar
\sphinxcode{\sphinxupquote{ObjConst}}
\begin{quote}

\sphinxAtStartPar
Double attribute.

\sphinxAtStartPar
The constant part of the objective function.
\end{quote}

\end{itemize}
\phantomsection\label{\detokenize{attribute:hasqobj}}\begin{itemize}
\item {} 
\sphinxAtStartPar
\sphinxcode{\sphinxupquote{HasQObj}}
\begin{quote}

\sphinxAtStartPar
Integer attribute.

\sphinxAtStartPar
Whether the problem has quadratic objective function.
\end{quote}

\end{itemize}
\phantomsection\label{\detokenize{attribute:haspsdobj}}\begin{itemize}
\item {} 
\sphinxAtStartPar
\sphinxcode{\sphinxupquote{HasPSDObj}}
\begin{quote}

\sphinxAtStartPar
Integer attribute.

\sphinxAtStartPar
Whether the problem has PSD terms in objective function.
\end{quote}

\end{itemize}
\phantomsection\label{\detokenize{attribute:hasnlobj}}\begin{itemize}
\item {} 
\sphinxAtStartPar
\sphinxcode{\sphinxupquote{HasNLObj}}
\begin{quote}

\sphinxAtStartPar
Integer attribute.

\sphinxAtStartPar
Whether the problem has nonlinear terms in objective function.
\end{quote}

\end{itemize}
\phantomsection\label{\detokenize{attribute:ismip}}\begin{itemize}
\item {} 
\sphinxAtStartPar
\sphinxcode{\sphinxupquote{IsMIP}}
\begin{quote}

\sphinxAtStartPar
Integer attribute.

\sphinxAtStartPar
Whether the problem is a MIP.
\end{quote}

\end{itemize}

\section{Solution related}
\label{\detokenize{attribute:solution-related}}
\sphinxAtStartPar
Solution related attributes provide the relevant information of the solution composition and description.
The names and descriptions of these attributes are summarized below.

\begin{savenotes}\sphinxattablestart
\sphinxthistablewithglobalstyle
\centering
\sphinxcapstartof{table}
\sphinxthecaptionisattop
\sphinxcaption{Solution related attributes}\label{\detokenize{attribute:id2}}
\sphinxaftertopcaption
\begin{tabular}[t]{|\X{15}{69}|\X{9}{69}|\X{45}{69}|}
\sphinxtoprule
\sphinxstyletheadfamily 
\sphinxAtStartPar
Name
&\sphinxstyletheadfamily 
\sphinxAtStartPar
Type
&\sphinxstyletheadfamily 
\sphinxAtStartPar
Description
\\
\sphinxmidrule
\sphinxtableatstartofbodyhook
\sphinxAtStartPar
{\hyperref[\detokenize{attribute:lpstatus}]{\sphinxcrossref{\DUrole{std,std-ref}{LpStatus}}}}
&
\sphinxAtStartPar
Integer
&
\sphinxAtStartPar
The LP status
\\
\sphinxhline
\sphinxAtStartPar
{\hyperref[\detokenize{attribute:mipstatus}]{\sphinxcrossref{\DUrole{std,std-ref}{MipStatus}}}}
&
\sphinxAtStartPar
Integer
&
\sphinxAtStartPar
The MIP status
\\
\sphinxhline
\sphinxAtStartPar
{\hyperref[\detokenize{attribute:simplexiter}]{\sphinxcrossref{\DUrole{std,std-ref}{SimplexIter}}}}
&
\sphinxAtStartPar
Integer
&
\sphinxAtStartPar
Number of simplex iterations performed
\\
\sphinxhline
\sphinxAtStartPar
{\hyperref[\detokenize{attribute:barrieriter}]{\sphinxcrossref{\DUrole{std,std-ref}{BarrierIter}}}}
&
\sphinxAtStartPar
Integer
&
\sphinxAtStartPar
Number of barrier iterations performed
\\
\sphinxhline
\sphinxAtStartPar
{\hyperref[\detokenize{attribute:nodecnt}]{\sphinxcrossref{\DUrole{std,std-ref}{NodeCnt}}}}
&
\sphinxAtStartPar
Integer
&
\sphinxAtStartPar
Number of explored nodes
\\
\sphinxhline
\sphinxAtStartPar
{\hyperref[\detokenize{solutionpool:poolsols}]{\sphinxcrossref{\DUrole{std,std-ref}{PoolSols}}}}
&
\sphinxAtStartPar
Integer
&
\sphinxAtStartPar
Number of solutions in solution pool
\\
\sphinxhline
\sphinxAtStartPar
{\hyperref[\detokenize{attribute:tuneresults}]{\sphinxcrossref{\DUrole{std,std-ref}{TuneResults}}}}
&
\sphinxAtStartPar
Integer
&
\sphinxAtStartPar
Number of parameter tuning results
\\
\sphinxhline
\sphinxAtStartPar
{\hyperref[\detokenize{attribute:haslpsol}]{\sphinxcrossref{\DUrole{std,std-ref}{HasLpSol}}}}
&
\sphinxAtStartPar
Integer
&
\sphinxAtStartPar
Whether LP solution is available
\\
\sphinxhline
\sphinxAtStartPar
{\hyperref[\detokenize{attribute:hasbasis}]{\sphinxcrossref{\DUrole{std,std-ref}{HasBasis}}}}
&
\sphinxAtStartPar
Integer
&
\sphinxAtStartPar
Whether LP basis is available
\\
\sphinxhline
\sphinxAtStartPar
{\hyperref[\detokenize{attribute:hasdualfarkas}]{\sphinxcrossref{\DUrole{std,std-ref}{HasDualFarkas}}}}
&
\sphinxAtStartPar
Integer
&
\sphinxAtStartPar
Whether the dual Farkas of an infeasible LP problem is available
\\
\sphinxhline
\sphinxAtStartPar
{\hyperref[\detokenize{attribute:hasprimalray}]{\sphinxcrossref{\DUrole{std,std-ref}{HasPrimalRay}}}}
&
\sphinxAtStartPar
Integer
&
\sphinxAtStartPar
Whether the primal ray of an unbounded LP problem is available
\\
\sphinxhline
\sphinxAtStartPar
{\hyperref[\detokenize{attribute:hasmipsol}]{\sphinxcrossref{\DUrole{std,std-ref}{HasMipSol}}}}
&
\sphinxAtStartPar
Integer
&
\sphinxAtStartPar
Whether MIP solution is available
\\
\sphinxhline
\sphinxAtStartPar
{\hyperref[\detokenize{attribute:iiscols}]{\sphinxcrossref{\DUrole{std,std-ref}{IISCols}}}}
&
\sphinxAtStartPar
Integer
&
\sphinxAtStartPar
Number of bounds of columns in IIS
\\
\sphinxhline
\sphinxAtStartPar
{\hyperref[\detokenize{attribute:iisrows}]{\sphinxcrossref{\DUrole{std,std-ref}{IISRows}}}}
&
\sphinxAtStartPar
Integer
&
\sphinxAtStartPar
Number of rows in IIS
\\
\sphinxhline
\sphinxAtStartPar
{\hyperref[\detokenize{attribute:iissoss}]{\sphinxcrossref{\DUrole{std,std-ref}{IISSOSs}}}}
&
\sphinxAtStartPar
Integer
&
\sphinxAtStartPar
Number of SOS constraints in IIS
\\
\sphinxhline
\sphinxAtStartPar
{\hyperref[\detokenize{attribute:iisindicators}]{\sphinxcrossref{\DUrole{std,std-ref}{IISIndicators}}}}
&
\sphinxAtStartPar
Integer
&
\sphinxAtStartPar
Number of indicator constraints in IIS
\\
\sphinxhline
\sphinxAtStartPar
{\hyperref[\detokenize{attribute:hasiis}]{\sphinxcrossref{\DUrole{std,std-ref}{HasIIS}}}}
&
\sphinxAtStartPar
Double
&
\sphinxAtStartPar
Whether IIS is available
\\
\sphinxhline
\sphinxAtStartPar
{\hyperref[\detokenize{attribute:hasfeasrelaxsol}]{\sphinxcrossref{\DUrole{std,std-ref}{HasFeasRelaxSol}}}}
&
\sphinxAtStartPar
Integer
&
\sphinxAtStartPar
Whether feasibility LP\sphinxhyphen{}relaxation solution is available
\\
\sphinxhline
\sphinxAtStartPar
{\hyperref[\detokenize{attribute:isminiis}]{\sphinxcrossref{\DUrole{std,std-ref}{IsMinIIS}}}}
&
\sphinxAtStartPar
Integer
&
\sphinxAtStartPar
Whether the computed IIS is minimal
\\
\sphinxhline
\sphinxAtStartPar
{\hyperref[\detokenize{attribute:hassensitivity}]{\sphinxcrossref{\DUrole{std,std-ref}{HasSensitivity}}}}
&
\sphinxAtStartPar
Integer
&
\sphinxAtStartPar
Whether sensitivity analysis results are available for LP problem
\\
\sphinxhline
\sphinxAtStartPar
{\hyperref[\detokenize{attribute:lpobjval}]{\sphinxcrossref{\DUrole{std,std-ref}{LpObjval}}}}
&
\sphinxAtStartPar
Double
&
\sphinxAtStartPar
The LP objective value
\\
\sphinxhline
\sphinxAtStartPar
{\hyperref[\detokenize{information:bestobj}]{\sphinxcrossref{\DUrole{std,std-ref}{BestObj}}}}
&
\sphinxAtStartPar
Double
&
\sphinxAtStartPar
Best integer objective value for MIP
\\
\sphinxhline
\sphinxAtStartPar
{\hyperref[\detokenize{information:bestbnd}]{\sphinxcrossref{\DUrole{std,std-ref}{BestBnd}}}}
&
\sphinxAtStartPar
Double
&
\sphinxAtStartPar
Best bound for MIP
\\
\sphinxhline
\sphinxAtStartPar
{\hyperref[\detokenize{attribute:bestgap}]{\sphinxcrossref{\DUrole{std,std-ref}{BestGap}}}}
&
\sphinxAtStartPar
Double
&
\sphinxAtStartPar
Best relative gap for MIP
\\
\sphinxhline
\sphinxAtStartPar
{\hyperref[\detokenize{attribute:feasrelaxobj}]{\sphinxcrossref{\DUrole{std,std-ref}{FeasRelaxObj}}}}
&
\sphinxAtStartPar
Double
&
\sphinxAtStartPar
Feasibility relaxation objective value
\\
\sphinxhline
\sphinxAtStartPar
{\hyperref[\detokenize{attribute:solvingtime}]{\sphinxcrossref{\DUrole{std,std-ref}{SolvingTime}}}}
&
\sphinxAtStartPar
Double
&
\sphinxAtStartPar
The time spent for the optimization (in seconds)
\\
\sphinxbottomrule
\end{tabular}
\sphinxtableafterendhook\par
\sphinxattableend\end{savenotes}
\phantomsection\label{\detokenize{attribute:lpstatus}}\begin{itemize}
\item {} 
\sphinxAtStartPar
\sphinxcode{\sphinxupquote{LpStatus}}
\begin{quote}

\sphinxAtStartPar
Integer attribute.

\sphinxAtStartPar
The LP status. Please refer to {\hyperref[\detokenize{constant:chapconst-solstatus}]{\sphinxcrossref{\DUrole{std,std-ref}{General Constants: Solution Status}}}} for possible values.
\end{quote}

\end{itemize}
\phantomsection\label{\detokenize{attribute:mipstatus}}\begin{itemize}
\item {} 
\sphinxAtStartPar
\sphinxcode{\sphinxupquote{MipStatus}}
\begin{quote}

\sphinxAtStartPar
Integer attribute.

\sphinxAtStartPar
The MIP status. Please refer to {\hyperref[\detokenize{constant:chapconst-solstatus}]{\sphinxcrossref{\DUrole{std,std-ref}{General Constants: Solution Status}}}} for possible values.
\end{quote}

\end{itemize}
\phantomsection\label{\detokenize{attribute:simplexiter}}\begin{itemize}
\item {} 
\sphinxAtStartPar
\sphinxcode{\sphinxupquote{SimplexIter}}
\begin{quote}

\sphinxAtStartPar
Integer attribute.

\sphinxAtStartPar
Number of simplex iterations performed.
\end{quote}

\end{itemize}
\phantomsection\label{\detokenize{attribute:barrieriter}}\begin{itemize}
\item {} 
\sphinxAtStartPar
\sphinxcode{\sphinxupquote{BarrierIter}}
\begin{quote}

\sphinxAtStartPar
Integer attribute.

\sphinxAtStartPar
Number of barrier iterations performed.
\end{quote}

\end{itemize}
\phantomsection\label{\detokenize{attribute:nodecnt}}\begin{itemize}
\item {} 
\sphinxAtStartPar
\sphinxcode{\sphinxupquote{NodeCnt}}
\begin{quote}

\sphinxAtStartPar
Integer attribute.

\sphinxAtStartPar
Number of explored nodes.
\end{quote}

\end{itemize}
\phantomsection\label{\detokenize{attribute:poolsols}}\begin{itemize}
\item {} 
\sphinxAtStartPar
\sphinxcode{\sphinxupquote{PoolSols}}
\begin{quote}

\sphinxAtStartPar
Integer attribute.

\sphinxAtStartPar
Number of solutions in solution pool.
\end{quote}

\end{itemize}
\phantomsection\label{\detokenize{attribute:tuneresults}}\begin{itemize}
\item {} 
\sphinxAtStartPar
\sphinxcode{\sphinxupquote{TuneResults}}
\begin{quote}

\sphinxAtStartPar
Integer attribute.

\sphinxAtStartPar
Number of parameter tuning results
\end{quote}

\end{itemize}
\phantomsection\label{\detokenize{attribute:haslpsol}}\begin{itemize}
\item {} 
\sphinxAtStartPar
\sphinxcode{\sphinxupquote{HasLpSol}}
\begin{quote}

\sphinxAtStartPar
Integer attribute.

\sphinxAtStartPar
Whether LP solution is available.
\end{quote}

\end{itemize}
\phantomsection\label{\detokenize{attribute:hasbasis}}\begin{itemize}
\item {} 
\sphinxAtStartPar
\sphinxcode{\sphinxupquote{HasBasis}}
\begin{quote}

\sphinxAtStartPar
Integer attribute.

\sphinxAtStartPar
Whether LP basis is available.
\end{quote}

\end{itemize}
\phantomsection\label{\detokenize{attribute:hasdualfarkas}}\begin{itemize}
\item {} 
\sphinxAtStartPar
\sphinxcode{\sphinxupquote{HasDualFarkas}}
\begin{quote}

\sphinxAtStartPar
Integer attribute.

\sphinxAtStartPar
Whether the dual Farkas of an infeasible LP problem is available.
\end{quote}

\end{itemize}
\phantomsection\label{\detokenize{attribute:hasprimalray}}\begin{itemize}
\item {} 
\sphinxAtStartPar
\sphinxcode{\sphinxupquote{HasPrimalRay}}
\begin{quote}

\sphinxAtStartPar
Integer attribute.

\sphinxAtStartPar
Whether the primal ray of an unbounded LP problem is available.
\end{quote}

\end{itemize}
\phantomsection\label{\detokenize{attribute:hasmipsol}}\begin{itemize}
\item {} 
\sphinxAtStartPar
\sphinxcode{\sphinxupquote{HasMipSol}}
\begin{quote}

\sphinxAtStartPar
Integer attribute.

\sphinxAtStartPar
Whether MIP solution is available.
\end{quote}

\end{itemize}
\phantomsection\label{\detokenize{attribute:iiscols}}\begin{itemize}
\item {} 
\sphinxAtStartPar
\sphinxcode{\sphinxupquote{IISCols}}
\begin{quote}

\sphinxAtStartPar
Integer attribute.

\sphinxAtStartPar
Number of bounds of columns in IIS.
\end{quote}

\end{itemize}
\phantomsection\label{\detokenize{attribute:iisrows}}\begin{itemize}
\item {} 
\sphinxAtStartPar
\sphinxcode{\sphinxupquote{IISRows}}
\begin{quote}

\sphinxAtStartPar
Integer attribute.

\sphinxAtStartPar
Number of rows in IIS.
\end{quote}

\end{itemize}
\phantomsection\label{\detokenize{attribute:iissoss}}\begin{itemize}
\item {} 
\sphinxAtStartPar
\sphinxcode{\sphinxupquote{IISSOSs}}
\begin{quote}

\sphinxAtStartPar
Integer attribute.

\sphinxAtStartPar
Number of SOS constraints in IIS.
\end{quote}

\end{itemize}
\phantomsection\label{\detokenize{attribute:iisindicators}}\begin{itemize}
\item {} 
\sphinxAtStartPar
\sphinxcode{\sphinxupquote{IISIndicators}}
\begin{quote}

\sphinxAtStartPar
Integer attribute.

\sphinxAtStartPar
Number of indicator constraints in IIS.
\end{quote}

\end{itemize}
\phantomsection\label{\detokenize{attribute:hasiis}}\begin{itemize}
\item {} 
\sphinxAtStartPar
\sphinxcode{\sphinxupquote{HasIIS}}
\begin{quote}

\sphinxAtStartPar
Integer attribute.

\sphinxAtStartPar
Whether IIS is available.
\end{quote}

\end{itemize}
\phantomsection\label{\detokenize{attribute:hasfeasrelaxsol}}\begin{itemize}
\item {} 
\sphinxAtStartPar
\sphinxcode{\sphinxupquote{HasFeasRelaxSol}}
\begin{quote}

\sphinxAtStartPar
Integer attribute.

\sphinxAtStartPar
Whether feasibility LP\sphinxhyphen{}relaxation solution is available.
\end{quote}

\end{itemize}
\phantomsection\label{\detokenize{attribute:isminiis}}\begin{itemize}
\item {} 
\sphinxAtStartPar
\sphinxcode{\sphinxupquote{IsMinIIS}}
\begin{quote}

\sphinxAtStartPar
Integer attribute.

\sphinxAtStartPar
Whether the computed IIS is minimal.
\end{quote}

\end{itemize}
\phantomsection\label{\detokenize{attribute:hassensitivity}}\begin{itemize}
\item {} 
\sphinxAtStartPar
\sphinxcode{\sphinxupquote{HasSensitivity}}
\begin{quote}

\sphinxAtStartPar
Integer attribute.

\sphinxAtStartPar
Whether sensitivity analysis results are available for LP problem.
\end{quote}

\end{itemize}
\phantomsection\label{\detokenize{attribute:lpobjval}}\begin{itemize}
\item {} 
\sphinxAtStartPar
\sphinxcode{\sphinxupquote{LpObjval}}
\begin{quote}

\sphinxAtStartPar
Double attribute.

\sphinxAtStartPar
The LP objective value.
\end{quote}

\end{itemize}
\phantomsection\label{\detokenize{attribute:bestobj}}\begin{itemize}
\item {} 
\sphinxAtStartPar
\sphinxcode{\sphinxupquote{BestObj}}
\begin{quote}

\sphinxAtStartPar
Double attribute.

\sphinxAtStartPar
Best integer objective value for MIP.
\end{quote}

\end{itemize}
\phantomsection\label{\detokenize{attribute:bestbnd}}\begin{itemize}
\item {} 
\sphinxAtStartPar
\sphinxcode{\sphinxupquote{BestBnd}}
\begin{quote}

\sphinxAtStartPar
Double attribute.

\sphinxAtStartPar
Best bound for MIP.
\end{quote}

\end{itemize}
\phantomsection\label{\detokenize{attribute:bestgap}}\begin{itemize}
\item {} 
\sphinxAtStartPar
\sphinxcode{\sphinxupquote{BestGap}}
\begin{quote}

\sphinxAtStartPar
Double attribute.

\sphinxAtStartPar
Best relative gap for MIP.
\end{quote}

\end{itemize}
\phantomsection\label{\detokenize{attribute:feasrelaxobj}}\begin{itemize}
\item {} 
\sphinxAtStartPar
\sphinxcode{\sphinxupquote{FeasRelaxObj}}
\begin{quote}

\sphinxAtStartPar
Double attribute.

\sphinxAtStartPar
Feasibility relaxation objective value.
\end{quote}

\end{itemize}
\phantomsection\label{\detokenize{attribute:solvingtime}}\begin{itemize}
\item {} 
\sphinxAtStartPar
\sphinxcode{\sphinxupquote{SolvingTime}}
\begin{quote}

\sphinxAtStartPar
Double attribute.

\sphinxAtStartPar
The time spent for the optimization (in seconds).
\end{quote}

\end{itemize}

\section{Methods for accessing attributes}
\label{\detokenize{attribute:methods-for-accessing-attributes}}\label{\detokenize{attribute:chapattrs-method}}
\sphinxAtStartPar
In different programming interfaces, the ways to access attributes are slightly different. For details, please refer to
the corresponding chapters for each programming language API:
\begin{itemize}
\item {} 
\sphinxAtStartPar
C API: {\hyperref[\detokenize{capiref:chapapi-getattr}]{\sphinxcrossref{\DUrole{std,std-ref}{C API Functions: Accessing attributes}}}}

\item {} 
\sphinxAtStartPar
C++ API: {\hyperref[\detokenize{cppapiref:chapcppapiref-attrs}]{\sphinxcrossref{\DUrole{std,std-ref}{C++ API Reference: Attributes}}}}

\item {} 
\sphinxAtStartPar
C\# API: {\hyperref[\detokenize{csharpapiref:chapcsharpapiref-attrs}]{\sphinxcrossref{\DUrole{std,std-ref}{C\# API Reference: Attributes}}}}

\item {} 
\sphinxAtStartPar
Java API: {\hyperref[\detokenize{javaapiref:chapjavaapiref-const-attrs}]{\sphinxcrossref{\DUrole{std,std-ref}{Java API Reference: Attributes}}}}

\item {} 
\sphinxAtStartPar
Python API: {\hyperref[\detokenize{pyapiref:chappyapi-const-attrs}]{\sphinxcrossref{\DUrole{std,std-ref}{Python API Reference: Attributes}}}}

\end{itemize}

\sphinxstepscope

\chapter{Information}
\label{\detokenize{information:information}}\label{\detokenize{information:chapinfo}}\label{\detokenize{information::doc}}
\sphinxAtStartPar
The information constants describe the relevant information of the model components (objective function, constraints and variables), solution results, and feasibility relaxation calculation results. This chapter will introduce the information constants provided by COPT and their meanings.

\begin{savenotes}\sphinxattablestart
\sphinxthistablewithglobalstyle
\centering
\sphinxcapstartof{table}
\sphinxthecaptionisattop
\sphinxcaption{COPT Information}\label{\detokenize{information:id1}}
\sphinxaftertopcaption
\begin{tabular}[t]{|\X{15}{69}|\X{9}{69}|\X{45}{69}|}
\sphinxtoprule
\sphinxstyletheadfamily 
\sphinxAtStartPar
Name
&\sphinxstyletheadfamily 
\sphinxAtStartPar
Type
&\sphinxstyletheadfamily 
\sphinxAtStartPar
Description
\\
\sphinxmidrule
\sphinxtableatstartofbodyhook
\sphinxAtStartPar
{\hyperref[\detokenize{information:obj}]{\sphinxcrossref{\DUrole{std,std-ref}{Obj}}}}
&
\sphinxAtStartPar
Double
&
\sphinxAtStartPar
Objective cost of columns
\\
\sphinxhline
\sphinxAtStartPar
{\hyperref[\detokenize{information:lb}]{\sphinxcrossref{\DUrole{std,std-ref}{LB}}}}
&
\sphinxAtStartPar
Double
&
\sphinxAtStartPar
Lower bounds of columns or rows
\\
\sphinxhline
\sphinxAtStartPar
{\hyperref[\detokenize{information:ub}]{\sphinxcrossref{\DUrole{std,std-ref}{UB}}}}
&
\sphinxAtStartPar
Double
&
\sphinxAtStartPar
Upper bounds of columns or rows
\\
\sphinxhline
\sphinxAtStartPar
{\hyperref[\detokenize{information:value}]{\sphinxcrossref{\DUrole{std,std-ref}{Value}}}}
&
\sphinxAtStartPar
Double
&
\sphinxAtStartPar
Solution of columns
\\
\sphinxhline
\sphinxAtStartPar
{\hyperref[\detokenize{information:slack}]{\sphinxcrossref{\DUrole{std,std-ref}{Slack}}}}
&
\sphinxAtStartPar
Double
&
\sphinxAtStartPar
Solution of slack variables, also known as activities of constraints. Only available for LP problem
\\
\sphinxhline
\sphinxAtStartPar
{\hyperref[\detokenize{information:dual}]{\sphinxcrossref{\DUrole{std,std-ref}{Dual}}}}
&
\sphinxAtStartPar
Double
&
\sphinxAtStartPar
Solution of dual variables. Only available for LP problem
\\
\sphinxhline
\sphinxAtStartPar
{\hyperref[\detokenize{information:redcost}]{\sphinxcrossref{\DUrole{std,std-ref}{RedCost}}}}
&
\sphinxAtStartPar
Double
&
\sphinxAtStartPar
Reduced cost of columns. Only available for LP problem
\\
\sphinxhline
\sphinxAtStartPar
{\hyperref[\detokenize{information:saobjlow}]{\sphinxcrossref{\DUrole{std,std-ref}{SAObjLow}}}}
&
\sphinxAtStartPar
Double
&
\sphinxAtStartPar
Indicates the minimum value to which the objective coefficient of a variable can be reduced while keeping the current basis optimal
\\
\sphinxhline
\sphinxAtStartPar
{\hyperref[\detokenize{information:saobjup}]{\sphinxcrossref{\DUrole{std,std-ref}{SAObjUp}}}}
&
\sphinxAtStartPar
Double
&
\sphinxAtStartPar
Indicates the maximum value to which the objective coefficient of a variable can be increased while keeping the current basis optimal
\\
\sphinxhline
\sphinxAtStartPar
{\hyperref[\detokenize{information:salblow}]{\sphinxcrossref{\DUrole{std,std-ref}{SALBLow}}}}
&
\sphinxAtStartPar
Double
&
\sphinxAtStartPar
Indicates the minimum value to which the lower bound of the variable/constraint can be reduced while keeping the current basis optimal
\\
\sphinxhline
\sphinxAtStartPar
{\hyperref[\detokenize{information:salbup}]{\sphinxcrossref{\DUrole{std,std-ref}{SALBUp}}}}
&
\sphinxAtStartPar
Double
&
\sphinxAtStartPar
Indicates the maximum value to which the lower bound of the variable/constraint can be increased while keeping the current basis optimal
\\
\sphinxhline
\sphinxAtStartPar
{\hyperref[\detokenize{information:saublow}]{\sphinxcrossref{\DUrole{std,std-ref}{SAUBLow}}}}
&
\sphinxAtStartPar
Double
&
\sphinxAtStartPar
Indicates the minimum value to which the upper bound of the variable/constraint can be reduced while keeping the current basis optimal
\\
\sphinxhline
\sphinxAtStartPar
{\hyperref[\detokenize{information:saubup}]{\sphinxcrossref{\DUrole{std,std-ref}{SAUBUp}}}}
&
\sphinxAtStartPar
Double
&
\sphinxAtStartPar
Indicates the maximum value to which the upper bound of the variable/constraint can be increased while keeping the current basis optimal
\\
\sphinxhline
\sphinxAtStartPar
{\hyperref[\detokenize{information:dualfarkas}]{\sphinxcrossref{\DUrole{std,std-ref}{DualFarkas}}}}
&
\sphinxAtStartPar
Double
&
\sphinxAtStartPar
The dual Farkas for constraints of an infeasible LP problem
\\
\sphinxhline
\sphinxAtStartPar
{\hyperref[\detokenize{information:primalray}]{\sphinxcrossref{\DUrole{std,std-ref}{PrimalRay}}}}
&
\sphinxAtStartPar
Double
&
\sphinxAtStartPar
The primal ray for variables of an unbounded LP problem
\\
\sphinxhline
\sphinxAtStartPar
{\hyperref[\detokenize{information:relaxlb}]{\sphinxcrossref{\DUrole{std,std-ref}{RelaxLB}}}}
&
\sphinxAtStartPar
Double
&
\sphinxAtStartPar
Feasibility relaxation values for lower bounds of columns or rows
\\
\sphinxhline
\sphinxAtStartPar
{\hyperref[\detokenize{information:relaxub}]{\sphinxcrossref{\DUrole{std,std-ref}{RelaxUB}}}}
&
\sphinxAtStartPar
Double
&
\sphinxAtStartPar
Feasibility relaxation values for upper bounds of columns or rows
\\
\sphinxhline
\sphinxAtStartPar
{\hyperref[\detokenize{information:relaxvalue}]{\sphinxcrossref{\DUrole{std,std-ref}{RelaxValue}}}}
&
\sphinxAtStartPar
Double
&
\sphinxAtStartPar
Solutions for the original model variables (columns) in the feasibility relaxation model
\\
\sphinxbottomrule
\end{tabular}
\sphinxtableafterendhook\par
\sphinxattableend\end{savenotes}

\begin{savenotes}\sphinxattablestart
\sphinxthistablewithglobalstyle
\centering
\sphinxcapstartof{table}
\sphinxthecaptionisattop
\sphinxcaption{COPT Callback Information}\label{\detokenize{information:id2}}
\sphinxaftertopcaption
\begin{tabular}[t]{|\X{15}{69}|\X{9}{69}|\X{45}{69}|}
\sphinxtoprule
\sphinxstyletheadfamily 
\sphinxAtStartPar
Name
&\sphinxstyletheadfamily 
\sphinxAtStartPar
Type
&\sphinxstyletheadfamily 
\sphinxAtStartPar
Description
\\
\sphinxmidrule
\sphinxtableatstartofbodyhook
\sphinxAtStartPar
{\hyperref[\detokenize{information:bestobj}]{\sphinxcrossref{\DUrole{std,std-ref}{BestObj}}}}
&
\sphinxAtStartPar
Double
&
\sphinxAtStartPar
Current best objective
\\
\sphinxhline
\sphinxAtStartPar
{\hyperref[\detokenize{information:bestbnd}]{\sphinxcrossref{\DUrole{std,std-ref}{BestBnd}}}}
&
\sphinxAtStartPar
Double
&
\sphinxAtStartPar
Current best objective bound
\\
\sphinxhline
\sphinxAtStartPar
{\hyperref[\detokenize{information:hasincumbent}]{\sphinxcrossref{\DUrole{std,std-ref}{HasIncumbent}}}}
&
\sphinxAtStartPar
Integer
&
\sphinxAtStartPar
Whether an incumbent is available
\\
\sphinxhline
\sphinxAtStartPar
{\hyperref[\detokenize{information:incumbent}]{\sphinxcrossref{\DUrole{std,std-ref}{Incumbent}}}}
&
\sphinxAtStartPar
Double
&
\sphinxAtStartPar
Current best feasible solution
\\
\sphinxhline
\sphinxAtStartPar
{\hyperref[\detokenize{information:mipcandobj}]{\sphinxcrossref{\DUrole{std,std-ref}{MipCandObj}}}}
&
\sphinxAtStartPar
Double
&
\sphinxAtStartPar
Objective value for current feasible solution candidate
\\
\sphinxhline
\sphinxAtStartPar
{\hyperref[\detokenize{information:mipcandidate}]{\sphinxcrossref{\DUrole{std,std-ref}{MipCandidate}}}}
&
\sphinxAtStartPar
Double
&
\sphinxAtStartPar
Current feasible solution candidate
\\
\sphinxhline
\sphinxAtStartPar
{\hyperref[\detokenize{information:relaxsolobj}]{\sphinxcrossref{\DUrole{std,std-ref}{RelaxSolObj}}}}
&
\sphinxAtStartPar
Double
&
\sphinxAtStartPar
Current Objective of LP\sphinxhyphen{}relaxation
\\
\sphinxhline
\sphinxAtStartPar
{\hyperref[\detokenize{information:relaxsolution}]{\sphinxcrossref{\DUrole{std,std-ref}{RelaxSolution}}}}
&
\sphinxAtStartPar
Double
&
\sphinxAtStartPar
Current solution of LP\sphinxhyphen{}relaxation
\\
\sphinxhline
\sphinxAtStartPar
{\hyperref[\detokenize{information:nodestatus}]{\sphinxcrossref{\DUrole{std,std-ref}{NodeStatus}}}}
&
\sphinxAtStartPar
Integer
&
\sphinxAtStartPar
The solution status of the LP\sphinxhyphen{}relaxation problem at the current node
\\
\sphinxbottomrule
\end{tabular}
\sphinxtableafterendhook\par
\sphinxattableend\end{savenotes}

\section{Problem information}
\label{\detokenize{information:problem-information}}\phantomsection\label{\detokenize{information:obj}}\begin{itemize}
\item {} 
\sphinxAtStartPar
\sphinxcode{\sphinxupquote{Obj}}
\begin{quote}

\sphinxAtStartPar
Double information.

\sphinxAtStartPar
Objective cost of columns.
\end{quote}

\end{itemize}
\phantomsection\label{\detokenize{information:lb}}\begin{itemize}
\item {} 
\sphinxAtStartPar
\sphinxcode{\sphinxupquote{LB}}
\begin{quote}

\sphinxAtStartPar
Double information.

\sphinxAtStartPar
Lower bounds of columns or rows.
\end{quote}

\end{itemize}
\phantomsection\label{\detokenize{information:ub}}\begin{itemize}
\item {} 
\sphinxAtStartPar
\sphinxcode{\sphinxupquote{UB}}
\begin{quote}

\sphinxAtStartPar
Double information.

\sphinxAtStartPar
Upper bounds of columns or rows.
\end{quote}

\end{itemize}

\section{Solution and sensitivity analysis information}
\label{\detokenize{information:solution-and-sensitivity-analysis-information}}\phantomsection\label{\detokenize{information:value}}\begin{itemize}
\item {} 
\sphinxAtStartPar
\sphinxcode{\sphinxupquote{Value}}
\begin{quote}

\sphinxAtStartPar
Double information.

\sphinxAtStartPar
Solution of columns.
\end{quote}

\end{itemize}
\phantomsection\label{\detokenize{information:slack}}\begin{itemize}
\item {} 
\sphinxAtStartPar
\sphinxcode{\sphinxupquote{Slack}}
\begin{quote}

\sphinxAtStartPar
Double information.

\sphinxAtStartPar
Solution of slack variables, also known as activities of constraints.
Only available for LP problem.
\end{quote}

\end{itemize}
\phantomsection\label{\detokenize{information:dual}}\begin{itemize}
\item {} 
\sphinxAtStartPar
\sphinxcode{\sphinxupquote{Dual}}
\begin{quote}

\sphinxAtStartPar
Double information.

\sphinxAtStartPar
Solution of dual variables. Only available for LP problem.
\end{quote}

\end{itemize}
\phantomsection\label{\detokenize{information:redcost}}\begin{itemize}
\item {} 
\sphinxAtStartPar
\sphinxcode{\sphinxupquote{RedCost}}
\begin{quote}

\sphinxAtStartPar
Double information.

\sphinxAtStartPar
Reduced cost of columns. Only available for LP problem.
\end{quote}

\end{itemize}
\phantomsection\label{\detokenize{information:saobjlow}}\begin{itemize}
\item {} 
\sphinxAtStartPar
\sphinxcode{\sphinxupquote{SAObjLow}}
\begin{quote}

\sphinxAtStartPar
Double information.

\sphinxAtStartPar
Sensitivity analysis information of the objective coefficient.

\sphinxAtStartPar
Indicates the minimum value to which the objective coefficient of a variable
can be reduced while keeping the current basis optimal.
\end{quote}

\end{itemize}
\phantomsection\label{\detokenize{information:saobjup}}\begin{itemize}
\item {} 
\sphinxAtStartPar
\sphinxcode{\sphinxupquote{SAObjUp}}
\begin{quote}

\sphinxAtStartPar
Double information.

\sphinxAtStartPar
Sensitivity analysis information of the objective coefficient.

\sphinxAtStartPar
Indicates the maximum value to which the objective coefficient of a variable
can be increased while keeping the current basis optimal.
\end{quote}

\end{itemize}
\phantomsection\label{\detokenize{information:salblow}}\begin{itemize}
\item {} 
\sphinxAtStartPar
\sphinxcode{\sphinxupquote{SALBLow}}
\begin{quote}

\sphinxAtStartPar
Double information.

\sphinxAtStartPar
Sensitivity analysis information of the lower bound of the variable/constraint.

\sphinxAtStartPar
Indicates the minimum value to which the lower bound of the variable/constraint
can be reduced while keeping the current basis optimal.
\end{quote}

\end{itemize}
\phantomsection\label{\detokenize{information:salbup}}\begin{itemize}
\item {} 
\sphinxAtStartPar
\sphinxcode{\sphinxupquote{SALBUp}}
\begin{quote}

\sphinxAtStartPar
Double information.

\sphinxAtStartPar
Sensitivity analysis information of the lower bound of the variable/constraint.

\sphinxAtStartPar
Indicates the maximum value to which the lower bound of the variable/constraint
can be increased while keeping the current basis optimal.
\end{quote}

\end{itemize}
\phantomsection\label{\detokenize{information:saublow}}\begin{itemize}
\item {} 
\sphinxAtStartPar
\sphinxcode{\sphinxupquote{SAUBLow}}
\begin{quote}

\sphinxAtStartPar
Double information.

\sphinxAtStartPar
Sensitivity analysis information of the upper bound of the variable/constraint.

\sphinxAtStartPar
Indicates the minimum value to which the upper bound of the variable/constraint
can be reduced while keeping the current basis optimal.
\end{quote}

\end{itemize}
\phantomsection\label{\detokenize{information:saubup}}\begin{itemize}
\item {} 
\sphinxAtStartPar
\sphinxcode{\sphinxupquote{SAUBUp}}
\begin{quote}

\sphinxAtStartPar
Double information.

\sphinxAtStartPar
Sensitivity analysis information of the upper bound of the variable/constraint.

\sphinxAtStartPar
Indicates the maximum value to which the upper bound of the variable/constraint
can be increased while keeping the current basis optimal.
\end{quote}

\end{itemize}

\section{Dual Farkas and primal ray}
\label{\detokenize{information:dual-farkas-and-primal-ray}}\begin{quote}

\sphinxAtStartPar
Advanced topic. When an LP is infeasible or unbounded,
the solver can return the dual Farkas or primal ray to prove it.
\end{quote}
\phantomsection\label{\detokenize{information:dualfarkas}}\begin{itemize}
\item {} 
\sphinxAtStartPar
\sphinxcode{\sphinxupquote{DualFarkas}}
\begin{quote}

\sphinxAtStartPar
Double information.

\sphinxAtStartPar
The dual Farkas for constraints of an infeasible LP problem.
Please enable the parameter \sphinxcode{\sphinxupquote{"ReqFarkasRay"}} to ensure that
the dual Farkas is available when the LP is infeasible.

\sphinxAtStartPar
Without loss of generality, the concept of the dual Farkas can be
conveniently demonstrated using an LP problem
with general variable bounds and equality constraints:
\(Ax = 0 \text{ and } l \leq x \leq u\).
When the LP is infeasible, a dual Farkas vector \(y\)
can prove that the system has conflict that \(\max y^TAx < y^T b = 0\).
Computing \(\max y^TAx\): with the vector \(\hat{a} = y^TA\),
choosing variable bound
\(x_i = l_i\) when \(\hat{a}_i < 0\) and
\(x_i = u_i\) when \(\hat{a}_i > 0\)
gives the maximal possible value of \(y^TAx\) for any \(x\) within their bounds.

\sphinxAtStartPar
Some application relies on the alternate conflict \(\min \bar{y}^TAx > \bar{y}^T b = 0\).
This can be achieved by negating the dual Farkas, i.e. \(\bar{y}=-y\) returned by the solver.

\sphinxAtStartPar
In very rare cases, the solver may fail to return a valid dual Farkas.
For example when the LP problem slightly infeasible by tiny amount, which
We recommend to study and to repair the infeasibility using FeasRelax instead.
\end{quote}

\end{itemize}
\phantomsection\label{\detokenize{information:primalray}}\begin{itemize}
\item {} 
\sphinxAtStartPar
\sphinxcode{\sphinxupquote{PrimalRay}}
\begin{quote}

\sphinxAtStartPar
Double information.

\sphinxAtStartPar
The primal ray for variables of an unbounded LP problem.
Please enable the parameter \sphinxcode{\sphinxupquote{"ReqFarkasRay"}} to ensure that the primal
ray is available when an LP is unbounded.

\sphinxAtStartPar
For a minimization LP problem in the standard form:
\(\min c^T x, Ax = b  \text{ and } x \geq 0\),
a primal ray vector \(r\) satisfies that \(r \geq 0, Ar = 0  \text{ and } c^T r < 0\).
\end{quote}

\end{itemize}

\section{Feasibility relaxation information}
\label{\detokenize{information:feasibility-relaxation-information}}\phantomsection\label{\detokenize{information:relaxlb}}\begin{itemize}
\item {} 
\sphinxAtStartPar
\sphinxcode{\sphinxupquote{RelaxLB}}
\begin{quote}

\sphinxAtStartPar
Double information.

\sphinxAtStartPar
Feasibility relaxation values for lower bounds of columns or rows.
\end{quote}

\end{itemize}
\phantomsection\label{\detokenize{information:relaxub}}\begin{itemize}
\item {} 
\sphinxAtStartPar
\sphinxcode{\sphinxupquote{RelaxUB}}
\begin{quote}

\sphinxAtStartPar
Double information.

\sphinxAtStartPar
Feasibility relaxation values for upper bounds of columns or rows.
\end{quote}

\end{itemize}
\phantomsection\label{\detokenize{information:relaxvalue}}\begin{itemize}
\item {} 
\sphinxAtStartPar
\sphinxcode{\sphinxupquote{RelaxValue}}
\begin{quote}

\sphinxAtStartPar
Double information.

\sphinxAtStartPar
Solutions for the original model variables (columns) in the feasibility relaxation model.
\end{quote}

\end{itemize}

\section{Callback information}
\label{\detokenize{information:callback-information}}\label{\detokenize{information:chapinfo-cbc}}\phantomsection\label{\detokenize{information:bestobj}}\begin{itemize}
\item {} 
\sphinxAtStartPar
\sphinxcode{\sphinxupquote{BestObj}}
\begin{quote}

\sphinxAtStartPar
Double information.

\sphinxAtStartPar
Current best objective.
\end{quote}

\end{itemize}
\phantomsection\label{\detokenize{information:bestbnd}}\begin{itemize}
\item {} 
\sphinxAtStartPar
\sphinxcode{\sphinxupquote{BestBnd}}
\begin{quote}

\sphinxAtStartPar
Double information.

\sphinxAtStartPar
Current best objective bound.
\end{quote}

\end{itemize}
\phantomsection\label{\detokenize{information:hasincumbent}}\begin{itemize}
\item {} 
\sphinxAtStartPar
\sphinxcode{\sphinxupquote{HasIncumbent}}
\begin{quote}

\sphinxAtStartPar
Integer information.

\sphinxAtStartPar
Whether an incumbent is available.
\end{quote}

\end{itemize}
\phantomsection\label{\detokenize{information:incumbent}}\begin{itemize}
\item {} 
\sphinxAtStartPar
\sphinxcode{\sphinxupquote{Incumbent}}
\begin{quote}

\sphinxAtStartPar
Double information.

\sphinxAtStartPar
Current best feasible solution.
\end{quote}

\end{itemize}
\phantomsection\label{\detokenize{information:mipcandobj}}\begin{itemize}
\item {} 
\sphinxAtStartPar
\sphinxcode{\sphinxupquote{MipCandObj}}
\begin{quote}

\sphinxAtStartPar
Double information.

\sphinxAtStartPar
Objective value for current feasible solution candidate.
\end{quote}

\end{itemize}
\phantomsection\label{\detokenize{information:mipcandidate}}\begin{itemize}
\item {} 
\sphinxAtStartPar
\sphinxcode{\sphinxupquote{MipCandidate}}
\begin{quote}

\sphinxAtStartPar
Double information.

\sphinxAtStartPar
Current feasible solution candidate.
\end{quote}

\end{itemize}
\phantomsection\label{\detokenize{information:relaxsolobj}}\begin{itemize}
\item {} 
\sphinxAtStartPar
\sphinxcode{\sphinxupquote{RelaxSolObj}}
\begin{quote}

\sphinxAtStartPar
Double information.

\sphinxAtStartPar
Current Objective of LP\sphinxhyphen{}relaxation.
\end{quote}

\end{itemize}
\phantomsection\label{\detokenize{information:relaxsolution}}\begin{itemize}
\item {} 
\sphinxAtStartPar
\sphinxcode{\sphinxupquote{RelaxSolution}}
\begin{quote}

\sphinxAtStartPar
Double information.

\sphinxAtStartPar
Current solution of LP\sphinxhyphen{}relaxation.
\end{quote}

\end{itemize}
\phantomsection\label{\detokenize{information:nodestatus}}\begin{itemize}
\item {} 
\sphinxAtStartPar
\sphinxcode{\sphinxupquote{NodeStatus}}
\begin{quote}

\sphinxAtStartPar
Integer information.

\sphinxAtStartPar
The solution status of the LP\sphinxhyphen{}relaxation problem at the current node. For possible values, please refer to: {\hyperref[\detokenize{constant:copttab-statuscodes}]{\sphinxcrossref{\DUrole{std,std-ref}{General Constants Chapter: Solution Status (Part)}}}}, except for \sphinxcode{\sphinxupquote{NODELIMIT}}, \sphinxcode{\sphinxupquote{UNSTARTED}}, \sphinxcode{\sphinxupquote{INF\_OR\_UNB}} .
\end{quote}

\end{itemize}

\section{Methods for accessing information}
\label{\detokenize{information:methods-for-accessing-information}}
\sphinxAtStartPar
In different programming interfaces, the ways to access and set information are slightly different. For details,
please refer to the corresponding chapters for each programming language API:
\begin{itemize}
\item {} 
\sphinxAtStartPar
C API: {\hyperref[\detokenize{capiref:chapapi-getinfo}]{\sphinxcrossref{\DUrole{std,std-ref}{C API Functions: Accessing information of problem}}}}

\item {} 
\sphinxAtStartPar
C\# API: {\hyperref[\detokenize{csharpapiref:chapcsharpapiref-info}]{\sphinxcrossref{\DUrole{std,std-ref}{C\# API Reference: Information}}}}

\item {} 
\sphinxAtStartPar
Java API: {\hyperref[\detokenize{javaapiref:chapjavaapiref-const-info}]{\sphinxcrossref{\DUrole{std,std-ref}{Java API Reference: Information}}}}

\item {} 
\sphinxAtStartPar
Python API: {\hyperref[\detokenize{pyapiref:chappyapi-const-info}]{\sphinxcrossref{\DUrole{std,std-ref}{Python API Reference: Information}}}}

\end{itemize}

\sphinxstepscope

\chapter{Parameters}
\label{\detokenize{parameter:parameters}}\label{\detokenize{parameter:chapparams}}\label{\detokenize{parameter::doc}}
\sphinxAtStartPar
Parameters control the operation of the \sphinxstylestrong{Cardinal Optimizer}. They can be modified before the optimization begins.
Each parameter has its own dfault value and value range. Before starting the solution, user can set the parameters
to different values, so as to put forward specific requirements on the solution algorithm and solution process.
Obviously, the default settings can also be maintained.

\sphinxAtStartPar
According to the task performed by the solver COPT and the optimization problem solved, it can be divided into
different types of parameters.

\sphinxAtStartPar
This chapter will introduce the parameters constants provided by COPT and their meanings.

\section{Limits and tolerances}
\label{\detokenize{parameter:limits-and-tolerances}}\label{\detokenize{parameter:chapparam-limittol}}

\begin{savenotes}\sphinxattablestart
\sphinxthistablewithglobalstyle
\centering
\sphinxcapstartof{table}
\sphinxthecaptionisattop
\sphinxcaption{Limits and tolerances parameters}\label{\detokenize{parameter:id1}}
\sphinxaftertopcaption
\begin{tabular}[t]{|\X{15}{59}|\X{9}{59}|\X{35}{59}|}
\sphinxtoprule
\sphinxstyletheadfamily 
\sphinxAtStartPar
Name
&\sphinxstyletheadfamily 
\sphinxAtStartPar
Type
&\sphinxstyletheadfamily 
\sphinxAtStartPar
Description
\\
\sphinxmidrule
\sphinxtableatstartofbodyhook
\sphinxAtStartPar
{\hyperref[\detokenize{parameter:timelimit}]{\sphinxcrossref{\DUrole{std,std-ref}{TimeLimit}}}}
&
\sphinxAtStartPar
Double
&
\sphinxAtStartPar
Time limit of the optimization
\\
\sphinxhline
\sphinxAtStartPar
{\hyperref[\detokenize{parameter:soltimelimit}]{\sphinxcrossref{\DUrole{std,std-ref}{SolTimeLimit}}}}
&
\sphinxAtStartPar
Double
&
\sphinxAtStartPar
Time limit if a primal feasible solution has been found
\\
\sphinxhline
\sphinxAtStartPar
{\hyperref[\detokenize{parameter:nodelimit}]{\sphinxcrossref{\DUrole{std,std-ref}{NodeLimit}}}}
&
\sphinxAtStartPar
Integer
&
\sphinxAtStartPar
Node limit of the optimization
\\
\sphinxhline
\sphinxAtStartPar
{\hyperref[\detokenize{parameter:bariterlimit}]{\sphinxcrossref{\DUrole{std,std-ref}{BarIterLimit}}}}
&
\sphinxAtStartPar
Integer
&
\sphinxAtStartPar
Iteration limit of barrier method
\\
\sphinxhline
\sphinxAtStartPar
{\hyperref[\detokenize{parameter:nlpiterlimit}]{\sphinxcrossref{\DUrole{std,std-ref}{NLPIterLimit}}}}
&
\sphinxAtStartPar
Integer
&
\sphinxAtStartPar
Iteration limit for the nonlinear solver
\\
\sphinxhline
\sphinxAtStartPar
{\hyperref[\detokenize{parameter:mipnlpiterlimit}]{\sphinxcrossref{\DUrole{std,std-ref}{MipNLPIterLimit}}}}
&
\sphinxAtStartPar
Integer
&
\sphinxAtStartPar
Iteration limit for solving NLP problem(s) within the MIP solver
\\
\sphinxhline
\sphinxAtStartPar
{\hyperref[\detokenize{parameter:matrixtol}]{\sphinxcrossref{\DUrole{std,std-ref}{MatrixTol}}}}
&
\sphinxAtStartPar
Double
&
\sphinxAtStartPar
Input matrix coefficient tolerance
\\
\sphinxhline
\sphinxAtStartPar
{\hyperref[\detokenize{parameter:feastol}]{\sphinxcrossref{\DUrole{std,std-ref}{FeasTol}}}}
&
\sphinxAtStartPar
Double
&
\sphinxAtStartPar
The feasibility tolerance
\\
\sphinxhline
\sphinxAtStartPar
{\hyperref[\detokenize{parameter:dualtol}]{\sphinxcrossref{\DUrole{std,std-ref}{DualTol}}}}
&
\sphinxAtStartPar
Double
&
\sphinxAtStartPar
The tolerance for dual solutions and reduced cost
\\
\sphinxhline
\sphinxAtStartPar
{\hyperref[\detokenize{parameter:inttol}]{\sphinxcrossref{\DUrole{std,std-ref}{IntTol}}}}
&
\sphinxAtStartPar
Double
&
\sphinxAtStartPar
The integrality tolerance for variables
\\
\sphinxhline
\sphinxAtStartPar
{\hyperref[\detokenize{parameter:pdlptol}]{\sphinxcrossref{\DUrole{std,std-ref}{PDLPTol}}}}
&
\sphinxAtStartPar
Double
&
\sphinxAtStartPar
The PDLP tolerance.
\\
\sphinxhline
\sphinxAtStartPar
{\hyperref[\detokenize{parameter:nlptol}]{\sphinxcrossref{\DUrole{std,std-ref}{NLPTol}}}}
&
\sphinxAtStartPar
Double
&
\sphinxAtStartPar
The NLP tolerance
\\
\sphinxhline
\sphinxAtStartPar
{\hyperref[\detokenize{parameter:relgap}]{\sphinxcrossref{\DUrole{std,std-ref}{RelGap}}}}
&
\sphinxAtStartPar
Double
&
\sphinxAtStartPar
The relative gap of optimization
\\
\sphinxhline
\sphinxAtStartPar
{\hyperref[\detokenize{parameter:absgap}]{\sphinxcrossref{\DUrole{std,std-ref}{AbsGap}}}}
&
\sphinxAtStartPar
Double
&
\sphinxAtStartPar
The absolute gap of optimization
\\
\sphinxbottomrule
\end{tabular}
\sphinxtableafterendhook\par
\sphinxattableend\end{savenotes}
\phantomsection\label{\detokenize{parameter:timelimit}}\begin{itemize}
\item {} 
\sphinxAtStartPar
\sphinxcode{\sphinxupquote{TimeLimit}}
\begin{quote}

\sphinxAtStartPar
Double parameter.

\sphinxAtStartPar
Time limit of the optimization.

\sphinxAtStartPar
\sphinxstylestrong{Default:} 1e20

\sphinxAtStartPar
\sphinxstylestrong{Minimal:} 0

\sphinxAtStartPar
\sphinxstylestrong{Maximal:} 1e20
\end{quote}

\end{itemize}
\phantomsection\label{\detokenize{parameter:soltimelimit}}\begin{itemize}
\item {} 
\sphinxAtStartPar
\sphinxcode{\sphinxupquote{SolTimeLimit}}
\begin{quote}

\sphinxAtStartPar
Double parameter.

\sphinxAtStartPar
Time limit if a primal feasible solution has been found.

\sphinxAtStartPar
\sphinxstylestrong{Default:} 1e20

\sphinxAtStartPar
\sphinxstylestrong{Minimal:} 0

\sphinxAtStartPar
\sphinxstylestrong{Maximal:} 1e20
\end{quote}

\end{itemize}
\phantomsection\label{\detokenize{parameter:nodelimit}}\begin{itemize}
\item {} 
\sphinxAtStartPar
\sphinxcode{\sphinxupquote{NodeLimit}}
\begin{quote}

\sphinxAtStartPar
Integer parameter.

\sphinxAtStartPar
Node limit of the optimization.

\sphinxAtStartPar
\sphinxstylestrong{Default:} \sphinxhyphen{}1 (Choose automatically)

\sphinxAtStartPar
\sphinxstylestrong{Minimal:} \sphinxhyphen{}1 (Choose automatically)

\sphinxAtStartPar
\sphinxstylestrong{Maximal:} \sphinxcode{\sphinxupquote{INT\_MAX}}
\end{quote}

\end{itemize}
\phantomsection\label{\detokenize{parameter:bariterlimit}}\begin{itemize}
\item {} 
\sphinxAtStartPar
\sphinxcode{\sphinxupquote{BarIterLimit}}
\begin{quote}

\sphinxAtStartPar
Integer parameter.

\sphinxAtStartPar
Iteration limit of barrier method.

\sphinxAtStartPar
\sphinxstylestrong{Default:} 500

\sphinxAtStartPar
\sphinxstylestrong{Minimal:} 0

\sphinxAtStartPar
\sphinxstylestrong{Maximal:} \sphinxcode{\sphinxupquote{INT\_MAX}}
\end{quote}

\end{itemize}
\phantomsection\label{\detokenize{parameter:nlpiterlimit}}\begin{itemize}
\item {} 
\sphinxAtStartPar
\sphinxcode{\sphinxupquote{NLPIterLimit}}
\begin{quote}

\sphinxAtStartPar
Integer parameter.

\sphinxAtStartPar
Iteration limit for the nonlinear solver.

\sphinxAtStartPar
\sphinxstylestrong{Default:} 1e4

\sphinxAtStartPar
\sphinxstylestrong{Minimum:} 0

\sphinxAtStartPar
\sphinxstylestrong{Maximum:} \sphinxcode{\sphinxupquote{INT\_MAX}}
\end{quote}

\end{itemize}
\phantomsection\label{\detokenize{parameter:mipnlpiterlimit}}\begin{itemize}
\item {} 
\sphinxAtStartPar
\sphinxcode{\sphinxupquote{MipNLPIterLimit}}
\begin{quote}

\sphinxAtStartPar
Integer parameter.

\sphinxAtStartPar
Iteration limit for solving NLP problem(s) within the MIP solver.

\sphinxAtStartPar
\sphinxstylestrong{Default:} 100

\sphinxAtStartPar
\sphinxstylestrong{Minimal:} \sphinxhyphen{}1 (no limit)

\sphinxAtStartPar
\sphinxstylestrong{Maximal:} \sphinxcode{\sphinxupquote{INT\_MAX}}
\end{quote}

\end{itemize}
\phantomsection\label{\detokenize{parameter:matrixtol}}\begin{itemize}
\item {} 
\sphinxAtStartPar
\sphinxcode{\sphinxupquote{MatrixTol}}
\begin{quote}

\sphinxAtStartPar
Double parameter.

\sphinxAtStartPar
Input matrix coefficient tolerance.

\sphinxAtStartPar
\sphinxstylestrong{Default:} 1e\sphinxhyphen{}10

\sphinxAtStartPar
\sphinxstylestrong{Minimal:} 0

\sphinxAtStartPar
\sphinxstylestrong{Maximal:} 1e\sphinxhyphen{}7
\end{quote}

\end{itemize}
\phantomsection\label{\detokenize{parameter:feastol}}\begin{itemize}
\item {} 
\sphinxAtStartPar
\sphinxcode{\sphinxupquote{FeasTol}}
\begin{quote}

\sphinxAtStartPar
Double parameter.

\sphinxAtStartPar
The feasibility tolerance.

\sphinxAtStartPar
\sphinxstylestrong{Default:} 1e\sphinxhyphen{}6

\sphinxAtStartPar
\sphinxstylestrong{Minimal:} 1e\sphinxhyphen{}9.

\sphinxAtStartPar
\sphinxstylestrong{Maximal:} 1e\sphinxhyphen{}4
\end{quote}

\end{itemize}
\phantomsection\label{\detokenize{parameter:dualtol}}\begin{itemize}
\item {} 
\sphinxAtStartPar
\sphinxcode{\sphinxupquote{DualTol}}
\begin{quote}

\sphinxAtStartPar
Double parameter.

\sphinxAtStartPar
The tolerance for dual solutions and reduced cost.

\sphinxAtStartPar
\sphinxstylestrong{Default:} 1e\sphinxhyphen{}6

\sphinxAtStartPar
\sphinxstylestrong{Minimal:} 1e\sphinxhyphen{}9

\sphinxAtStartPar
\sphinxstylestrong{Maximal:} 1e\sphinxhyphen{}4
\end{quote}

\end{itemize}
\phantomsection\label{\detokenize{parameter:inttol}}\begin{itemize}
\item {} 
\sphinxAtStartPar
\sphinxcode{\sphinxupquote{IntTol}}
\begin{quote}

\sphinxAtStartPar
Double parameter.

\sphinxAtStartPar
The integrality tolerance for variables.

\sphinxAtStartPar
\sphinxstylestrong{Default:} 1e\sphinxhyphen{}6

\sphinxAtStartPar
\sphinxstylestrong{Minimal:} 1e\sphinxhyphen{}9

\sphinxAtStartPar
\sphinxstylestrong{Maximal:} 1e\sphinxhyphen{}1
\end{quote}

\end{itemize}
\phantomsection\label{\detokenize{parameter:pdlptol}}\begin{itemize}
\item {} 
\sphinxAtStartPar
\sphinxcode{\sphinxupquote{PDLPTol}}
\begin{quote}

\sphinxAtStartPar
Double parameter.

\sphinxAtStartPar
The PDLP tolerance.

\sphinxAtStartPar
\sphinxstylestrong{Default:} 1e\sphinxhyphen{}6

\sphinxAtStartPar
\sphinxstylestrong{Minimal:} 1e\sphinxhyphen{}12

\sphinxAtStartPar
\sphinxstylestrong{Maximal:} 1e\sphinxhyphen{}4
\end{quote}

\end{itemize}
\phantomsection\label{\detokenize{parameter:nlptol}}\begin{itemize}
\item {} 
\sphinxAtStartPar
\sphinxcode{\sphinxupquote{NLPTol}}
\begin{quote}

\sphinxAtStartPar
Double parameter.

\sphinxAtStartPar
The NLP tolerance.

\sphinxAtStartPar
\sphinxstylestrong{Default:} 1e\sphinxhyphen{}8

\sphinxAtStartPar
\sphinxstylestrong{Minimum:} 1e\sphinxhyphen{}13

\sphinxAtStartPar
\sphinxstylestrong{Maximum:} 1e\sphinxhyphen{}3
\end{quote}

\end{itemize}
\phantomsection\label{\detokenize{parameter:relgap}}\begin{itemize}
\item {} 
\sphinxAtStartPar
\sphinxcode{\sphinxupquote{RelGap}}
\begin{quote}

\sphinxAtStartPar
Double parameter.

\sphinxAtStartPar
The relative gap of optimization.

\sphinxAtStartPar
\sphinxstylestrong{Default:} 1e\sphinxhyphen{}4

\sphinxAtStartPar
\sphinxstylestrong{Minimal:} 0

\sphinxAtStartPar
\sphinxstylestrong{Maximal:} \sphinxcode{\sphinxupquote{DBL\_MAX}}
\end{quote}

\end{itemize}
\phantomsection\label{\detokenize{parameter:absgap}}\begin{itemize}
\item {} 
\sphinxAtStartPar
\sphinxcode{\sphinxupquote{AbsGap}}
\begin{quote}

\sphinxAtStartPar
Double parameter.

\sphinxAtStartPar
The absolute gap of optimization.

\sphinxAtStartPar
\sphinxstylestrong{Default:} 1e\sphinxhyphen{}6

\sphinxAtStartPar
\sphinxstylestrong{Minimal:} 0

\sphinxAtStartPar
\sphinxstylestrong{Maximal:} \sphinxcode{\sphinxupquote{DBL\_MAX}}
\end{quote}

\end{itemize}

\section{Presolving and scaling}
\label{\detokenize{parameter:presolving-and-scaling}}\label{\detokenize{parameter:chapparam-pre}}

\begin{savenotes}\sphinxattablestart
\sphinxthistablewithglobalstyle
\centering
\sphinxcapstartof{table}
\sphinxthecaptionisattop
\sphinxcaption{Presolving and scaling parameters}\label{\detokenize{parameter:id2}}
\sphinxaftertopcaption
\begin{tabular}[t]{|\X{15}{59}|\X{9}{59}|\X{35}{59}|}
\sphinxtoprule
\sphinxstyletheadfamily 
\sphinxAtStartPar
Name
&\sphinxstyletheadfamily 
\sphinxAtStartPar
Type
&\sphinxstyletheadfamily 
\sphinxAtStartPar
Description
\\
\sphinxmidrule
\sphinxtableatstartofbodyhook
\sphinxAtStartPar
{\hyperref[\detokenize{parameter:presolve}]{\sphinxcrossref{\DUrole{std,std-ref}{Presolve}}}}
&
\sphinxAtStartPar
Integer
&
\sphinxAtStartPar
Level of presolving before solving a model
\\
\sphinxhline
\sphinxAtStartPar
{\hyperref[\detokenize{parameter:scaling}]{\sphinxcrossref{\DUrole{std,std-ref}{Scaling}}}}
&
\sphinxAtStartPar
Integer
&
\sphinxAtStartPar
Whether to perform scaling before solving a problem
\\
\sphinxhline
\sphinxAtStartPar
{\hyperref[\detokenize{parameter:dualize}]{\sphinxcrossref{\DUrole{std,std-ref}{Dualize}}}}
&
\sphinxAtStartPar
Integer
&
\sphinxAtStartPar
Whether to dualize a problem before solving it
\\
\sphinxbottomrule
\end{tabular}
\sphinxtableafterendhook\par
\sphinxattableend\end{savenotes}
\phantomsection\label{\detokenize{parameter:presolve}}\begin{itemize}
\item {} 
\sphinxAtStartPar
\sphinxcode{\sphinxupquote{Presolve}}
\begin{quote}

\sphinxAtStartPar
Integer parameter.

\sphinxAtStartPar
Level of presolving before solving a model.

\sphinxAtStartPar
\sphinxstylestrong{Default:} \sphinxhyphen{}1

\sphinxAtStartPar
\sphinxstylestrong{Possible values:}
\begin{quote}

\sphinxAtStartPar
\sphinxhyphen{}1: Choose automatically.

\sphinxAtStartPar
0: Off.

\sphinxAtStartPar
1: Fast.

\sphinxAtStartPar
2: Normal.

\sphinxAtStartPar
3: Aggressive.

\sphinxAtStartPar
4: No Limitations, continues until the model cannot be modified (may be very time\sphinxhyphen{}consuming).
\end{quote}
\end{quote}

\end{itemize}
\phantomsection\label{\detokenize{parameter:scaling}}\begin{itemize}
\item {} 
\sphinxAtStartPar
\sphinxcode{\sphinxupquote{Scaling}}
\begin{quote}

\sphinxAtStartPar
Integer parameter.

\sphinxAtStartPar
Whether to perform scaling before solving a problem.

\sphinxAtStartPar
\sphinxstylestrong{Default:} \sphinxhyphen{}1

\sphinxAtStartPar
\sphinxstylestrong{Possible values:}
\begin{quote}

\sphinxAtStartPar
\sphinxhyphen{}1: Choose automatically.

\sphinxAtStartPar
0: No scaling.

\sphinxAtStartPar
1: Apply scaling.
\end{quote}
\end{quote}

\end{itemize}
\phantomsection\label{\detokenize{parameter:dualize}}\begin{itemize}
\item {} 
\sphinxAtStartPar
\sphinxcode{\sphinxupquote{Dualize}}
\begin{quote}

\sphinxAtStartPar
Integer parameter.

\sphinxAtStartPar
Whether to dualize a problem before solving it.

\sphinxAtStartPar
\sphinxstylestrong{Default:} \sphinxhyphen{}1

\sphinxAtStartPar
\sphinxstylestrong{Possible values:}
\begin{quote}

\sphinxAtStartPar
\sphinxhyphen{}1: Choose automatically.

\sphinxAtStartPar
0: No dualizing.

\sphinxAtStartPar
1: Dualizing the problem.
\end{quote}
\end{quote}

\end{itemize}

\section{Linear programming related}
\label{\detokenize{parameter:linear-programming-related}}\label{\detokenize{parameter:chapparam-lp}}

\begin{savenotes}\sphinxattablestart
\sphinxthistablewithglobalstyle
\centering
\sphinxcapstartof{table}
\sphinxthecaptionisattop
\sphinxcaption{Linear programming related parameters}\label{\detokenize{parameter:id3}}
\sphinxaftertopcaption
\begin{tabular}[t]{|\X{15}{59}|\X{9}{59}|\X{35}{59}|}
\sphinxtoprule
\sphinxstyletheadfamily 
\sphinxAtStartPar
Name
&\sphinxstyletheadfamily 
\sphinxAtStartPar
Type
&\sphinxstyletheadfamily 
\sphinxAtStartPar
Description
\\
\sphinxmidrule
\sphinxtableatstartofbodyhook
\sphinxAtStartPar
{\hyperref[\detokenize{parameter:lpmethod}]{\sphinxcrossref{\DUrole{std,std-ref}{LpMethod}}}}
&
\sphinxAtStartPar
Integer
&
\sphinxAtStartPar
Method to solve the LP problem
\\
\sphinxhline
\sphinxAtStartPar
{\hyperref[\detokenize{parameter:dualprice}]{\sphinxcrossref{\DUrole{std,std-ref}{DualPrice}}}}
&
\sphinxAtStartPar
Integer
&
\sphinxAtStartPar
Specifies the dual simplex pricing algorithm
\\
\sphinxhline
\sphinxAtStartPar
{\hyperref[\detokenize{parameter:dualperturb}]{\sphinxcrossref{\DUrole{std,std-ref}{DualPerturb}}}}
&
\sphinxAtStartPar
Integer
&
\sphinxAtStartPar
Whether to allow the objective function perturbation when using the dual simplex method
\\
\sphinxhline
\sphinxAtStartPar
{\hyperref[\detokenize{parameter:barhomogeneous}]{\sphinxcrossref{\DUrole{std,std-ref}{BarHomogeneous}}}}
&
\sphinxAtStartPar
Integer
&
\sphinxAtStartPar
Whether to use homogeneous self\sphinxhyphen{}dual form in barrier
\\
\sphinxhline
\sphinxAtStartPar
{\hyperref[\detokenize{parameter:barorder}]{\sphinxcrossref{\DUrole{std,std-ref}{BarOrder}}}}
&
\sphinxAtStartPar
Integer
&
\sphinxAtStartPar
Ordering algorithm in barrier method
\\
\sphinxhline
\sphinxAtStartPar
{\hyperref[\detokenize{parameter:barstart}]{\sphinxcrossref{\DUrole{std,std-ref}{BarStart}}}}
&
\sphinxAtStartPar
Integer
&
\sphinxAtStartPar
Algorithm for finding initial points in barrier method
\\
\sphinxhline
\sphinxAtStartPar
{\hyperref[\detokenize{parameter:crossover}]{\sphinxcrossref{\DUrole{std,std-ref}{Crossover}}}}
&
\sphinxAtStartPar
Integer
&
\sphinxAtStartPar
Whether to use crossover
\\
\sphinxhline
\sphinxAtStartPar
{\hyperref[\detokenize{parameter:reqfarkasray}]{\sphinxcrossref{\DUrole{std,std-ref}{ReqFarkasRay}}}}
&
\sphinxAtStartPar
Integer
&
\sphinxAtStartPar
Advanced topic. Whether to compute the dual Farkas or primal ray when the LP is infeasible or unbounded
\\
\sphinxhline
\sphinxAtStartPar
{\hyperref[\detokenize{parameter:reqsensitivity}]{\sphinxcrossref{\DUrole{std,std-ref}{ReqSensitivity}}}}
&
\sphinxAtStartPar
Integer
&
\sphinxAtStartPar
Whether to compute sensitivity analysis when an optimal basis is available for an LP problem (when solved by the simplex method or by other methods followed by crossover)
\\
\sphinxbottomrule
\end{tabular}
\sphinxtableafterendhook\par
\sphinxattableend\end{savenotes}
\phantomsection\label{\detokenize{parameter:lpmethod}}\begin{itemize}
\item {} 
\sphinxAtStartPar
\sphinxcode{\sphinxupquote{LpMethod}}
\begin{quote}

\sphinxAtStartPar
Integer parameter.

\sphinxAtStartPar
Method to solve the LP problem.

\sphinxAtStartPar
\sphinxstylestrong{Default:} \sphinxhyphen{}1

\sphinxAtStartPar
\sphinxstylestrong{Possible values:}
\begin{quote}

\sphinxAtStartPar
\sphinxhyphen{}1: Choose automatically.
\begin{quote}

\sphinxAtStartPar
For Linear Programming, choose dual simplex method;

\sphinxAtStartPar
For Mixed Integer Linear Programming, choose dual simplex or barrier method.
\end{quote}

\sphinxAtStartPar
1: Dual simplex.

\sphinxAtStartPar
2: Barrier.

\sphinxAtStartPar
3: Crossover.

\sphinxAtStartPar
4: Concurrent (Use multiple algorithms simultaneously).

\sphinxAtStartPar
5: Choose between simplex and barrier automatically (Based on features such as sparsity and/or coefficients ranges).

\sphinxAtStartPar
6: First\sphinxhyphen{}order method (PDLP).
\end{quote}
\end{quote}

\end{itemize}

\begin{sphinxadmonition}{note}{Note:}

\sphinxAtStartPar
Currently, COPT GPU solver supports solving Linear Programming problems and the root relaxation of Integer Programming problems.
To enable it, either the first\sphinxhyphen{}order method (PDLP) or barrier method need to be used, and {\hyperref[\detokenize{parameter:gpumode}]{\sphinxcrossref{\DUrole{std,std-ref}{GPUMode}}}} need to be set accordingly.
\end{sphinxadmonition}
\phantomsection\label{\detokenize{parameter:dualprice}}\begin{itemize}
\item {} 
\sphinxAtStartPar
\sphinxcode{\sphinxupquote{DualPrice}}
\begin{quote}

\sphinxAtStartPar
Integer parameter.

\sphinxAtStartPar
Specifies the dual simplex pricing algorithm.

\sphinxAtStartPar
\sphinxstylestrong{Default:} \sphinxhyphen{}1

\sphinxAtStartPar
\sphinxstylestrong{Possible values:}
\begin{quote}

\sphinxAtStartPar
\sphinxhyphen{}1: Choose automatically.

\sphinxAtStartPar
0: Using Devex pricing algorithm.

\sphinxAtStartPar
1: Using dual steepest\sphinxhyphen{}edge pricing algorithm.
\end{quote}
\end{quote}

\end{itemize}
\phantomsection\label{\detokenize{parameter:dualperturb}}\begin{itemize}
\item {} 
\sphinxAtStartPar
\sphinxcode{\sphinxupquote{DualPerturb}}
\begin{quote}

\sphinxAtStartPar
Integer parameter.

\sphinxAtStartPar
Whether to allow the objective function perturbation
when using the dual simplex method.

\sphinxAtStartPar
\sphinxstylestrong{Default:} \sphinxhyphen{}1

\sphinxAtStartPar
\sphinxstylestrong{Possible values:}
\begin{quote}

\sphinxAtStartPar
\sphinxhyphen{}1: Choose automatically.

\sphinxAtStartPar
0: No perturbation.

\sphinxAtStartPar
1: Allow objective function perturbation.
\end{quote}
\end{quote}

\end{itemize}
\phantomsection\label{\detokenize{parameter:barhomogeneous}}\begin{itemize}
\item {} 
\sphinxAtStartPar
\sphinxcode{\sphinxupquote{BarHomogeneous}}
\begin{quote}

\sphinxAtStartPar
Integer parameter.

\sphinxAtStartPar
Whether to use homogeneous self\sphinxhyphen{}dual form in barrier.

\sphinxAtStartPar
\sphinxstylestrong{Default:} \sphinxhyphen{}1

\sphinxAtStartPar
\sphinxstylestrong{Possible values:}
\begin{quote}

\sphinxAtStartPar
\sphinxhyphen{}1: Choose automatically.

\sphinxAtStartPar
0: No.

\sphinxAtStartPar
1: Yes.
\end{quote}
\end{quote}

\end{itemize}
\phantomsection\label{\detokenize{parameter:barorder}}\begin{itemize}
\item {} 
\sphinxAtStartPar
\sphinxcode{\sphinxupquote{BarOrder}}
\begin{quote}

\sphinxAtStartPar
Integer parameter.

\sphinxAtStartPar
Ordering algorithm in barrier method.

\sphinxAtStartPar
\sphinxstylestrong{Default:} \sphinxhyphen{}1

\sphinxAtStartPar
\sphinxstylestrong{Possible values:}
\begin{quote}

\sphinxAtStartPar
\sphinxhyphen{}1: Choose automatically.

\sphinxAtStartPar
0: Approximate Minimum Degree (AMD).

\sphinxAtStartPar
1: Nested Dissection (ND1).

\sphinxAtStartPar
2: Modified Nested Dissection (ND2).
\end{quote}
\end{quote}

\end{itemize}
\phantomsection\label{\detokenize{parameter:barstart}}\begin{itemize}
\item {} 
\sphinxAtStartPar
\sphinxcode{\sphinxupquote{BarStart}}
\begin{quote}

\sphinxAtStartPar
Integer parameter.

\sphinxAtStartPar
Algorithm for finding initial points in barrier method.

\sphinxAtStartPar
\sphinxstylestrong{Default:} \sphinxhyphen{}1

\sphinxAtStartPar
\sphinxstylestrong{Possible values:}
\begin{quote}

\sphinxAtStartPar
\sphinxhyphen{}1: Choose automatically.

\sphinxAtStartPar
0: Simple.

\sphinxAtStartPar
1: Mehrotra.

\sphinxAtStartPar
2: Modified Mehrotra.
\end{quote}
\end{quote}

\end{itemize}
\phantomsection\label{\detokenize{parameter:crossover}}\begin{itemize}
\item {} 
\sphinxAtStartPar
\sphinxcode{\sphinxupquote{Crossover}}
\begin{quote}

\sphinxAtStartPar
Integer parameter.

\sphinxAtStartPar
Whether to use crossover.

\sphinxAtStartPar
\sphinxstylestrong{Default:} 1

\sphinxAtStartPar
\sphinxstylestrong{Possible values:}
\begin{quote}

\sphinxAtStartPar
\sphinxhyphen{}1: Choose automatically.
\begin{quote}

\sphinxAtStartPar
Only run crossover when the LP solution is not primal\sphinxhyphen{}dual feasible.
\end{quote}

\sphinxAtStartPar
0: No.

\sphinxAtStartPar
1: Yes.
\end{quote}
\end{quote}

\end{itemize}
\phantomsection\label{\detokenize{parameter:reqfarkasray}}\begin{itemize}
\item {} 
\sphinxAtStartPar
\sphinxcode{\sphinxupquote{ReqFarkasRay}}
\begin{quote}

\sphinxAtStartPar
Integer parameter.

\sphinxAtStartPar
Advanced topic. Whether to compute the dual Farkas or primal ray when the LP is infeasible or unbounded.

\sphinxAtStartPar
\sphinxstylestrong{Default:} 0

\sphinxAtStartPar
\sphinxstylestrong{Possible values:}
\begin{quote}

\sphinxAtStartPar
0: No.

\sphinxAtStartPar
1: Yes.
\end{quote}
\end{quote}

\end{itemize}
\phantomsection\label{\detokenize{parameter:reqsensitivity}}\begin{itemize}
\item {} 
\sphinxAtStartPar
\sphinxcode{\sphinxupquote{ReqSensitivity}}
\begin{quote}

\sphinxAtStartPar
Integer parameter.

\sphinxAtStartPar
Whether to compute sensitivity analysis when an optimal basis is available for an LP problem
(when solved by the simplex method or by other methods followed by crossover).

\sphinxAtStartPar
\sphinxstylestrong{Default:} 0

\sphinxAtStartPar
\sphinxstylestrong{Possible values:}
\begin{quote}

\sphinxAtStartPar
0: No.

\sphinxAtStartPar
1: Yes.
\end{quote}
\end{quote}

\end{itemize}

\section{Integer Programming related}
\label{\detokenize{parameter:integer-programming-related}}\label{\detokenize{parameter:chapparam-mip}}

\begin{savenotes}\sphinxattablestart
\sphinxthistablewithglobalstyle
\centering
\sphinxcapstartof{table}
\sphinxthecaptionisattop
\sphinxcaption{Integer programming related parameters}\label{\detokenize{parameter:id4}}
\sphinxaftertopcaption
\begin{tabular}[t]{|\X{15}{59}|\X{9}{59}|\X{35}{59}|}
\sphinxtoprule
\sphinxstyletheadfamily 
\sphinxAtStartPar
Name
&\sphinxstyletheadfamily 
\sphinxAtStartPar
Type
&\sphinxstyletheadfamily 
\sphinxAtStartPar
Description
\\
\sphinxmidrule
\sphinxtableatstartofbodyhook
\sphinxAtStartPar
{\hyperref[\detokenize{parameter:cutlevel}]{\sphinxcrossref{\DUrole{std,std-ref}{CutLevel}}}}
&
\sphinxAtStartPar
Integer
&
\sphinxAtStartPar
Level of cutting\sphinxhyphen{}planes generation
\\
\sphinxhline
\sphinxAtStartPar
{\hyperref[\detokenize{parameter:rootcutlevel}]{\sphinxcrossref{\DUrole{std,std-ref}{RootCutLevel}}}}
&
\sphinxAtStartPar
Integer
&
\sphinxAtStartPar
Level of cutting\sphinxhyphen{}planes generation of root node
\\
\sphinxhline
\sphinxAtStartPar
{\hyperref[\detokenize{parameter:treecutlevel}]{\sphinxcrossref{\DUrole{std,std-ref}{TreeCutLevel}}}}
&
\sphinxAtStartPar
Integer
&
\sphinxAtStartPar
Level of cutting\sphinxhyphen{}planes generation of search tree
\\
\sphinxhline
\sphinxAtStartPar
{\hyperref[\detokenize{parameter:rootcutrounds}]{\sphinxcrossref{\DUrole{std,std-ref}{RootCutRounds}}}}
&
\sphinxAtStartPar
Integer
&
\sphinxAtStartPar
Rounds of cutting\sphinxhyphen{}planes generation of root node
\\
\sphinxhline
\sphinxAtStartPar
{\hyperref[\detokenize{parameter:nodecutrounds}]{\sphinxcrossref{\DUrole{std,std-ref}{NodeCutRounds}}}}
&
\sphinxAtStartPar
Integer
&
\sphinxAtStartPar
Rounds of cutting\sphinxhyphen{}planes generation of search tree node
\\
\sphinxhline
\sphinxAtStartPar
{\hyperref[\detokenize{parameter:heurlevel}]{\sphinxcrossref{\DUrole{std,std-ref}{HeurLevel}}}}
&
\sphinxAtStartPar
Integer
&
\sphinxAtStartPar
Level of heuristics
\\
\sphinxhline
\sphinxAtStartPar
{\hyperref[\detokenize{parameter:prerootheurlevel}]{\sphinxcrossref{\DUrole{std,std-ref}{PreRootHeurLevel}}}}
&
\sphinxAtStartPar
Integer
&
\sphinxAtStartPar
Level of pre\sphinxhyphen{}root heuristics
\\
\sphinxhline
\sphinxAtStartPar
{\hyperref[\detokenize{parameter:roundingheurlevel}]{\sphinxcrossref{\DUrole{std,std-ref}{RoundingHeurLevel}}}}
&
\sphinxAtStartPar
Integer
&
\sphinxAtStartPar
Level of rounding heuristics
\\
\sphinxhline
\sphinxAtStartPar
{\hyperref[\detokenize{parameter:divingheurlevel}]{\sphinxcrossref{\DUrole{std,std-ref}{DivingHeurLevel}}}}
&
\sphinxAtStartPar
Integer
&
\sphinxAtStartPar
Level of diving heuristics
\\
\sphinxhline
\sphinxAtStartPar
{\hyperref[\detokenize{parameter:submipheurlevel}]{\sphinxcrossref{\DUrole{std,std-ref}{SubMipHeurLevel}}}}
&
\sphinxAtStartPar
Integer
&
\sphinxAtStartPar
Level of Sub\sphinxhyphen{}MIP heuristics
\\
\sphinxhline
\sphinxAtStartPar
{\hyperref[\detokenize{parameter:fapheurlevel}]{\sphinxcrossref{\DUrole{std,std-ref}{FAPHeurLevel}}}}
&
\sphinxAtStartPar
Integer
&
\sphinxAtStartPar
Level of Fix\sphinxhyphen{}and\sphinxhyphen{}propagate heuristics
\\
\sphinxhline
\sphinxAtStartPar
{\hyperref[\detokenize{parameter:strongbranching}]{\sphinxcrossref{\DUrole{std,std-ref}{StrongBranching}}}}
&
\sphinxAtStartPar
Integer
&
\sphinxAtStartPar
Level of strong branching
\\
\sphinxhline
\sphinxAtStartPar
{\hyperref[\detokenize{parameter:conflictanalysis}]{\sphinxcrossref{\DUrole{std,std-ref}{ConflictAnalysis}}}}
&
\sphinxAtStartPar
Integer
&
\sphinxAtStartPar
Whether to perform conflict analysis
\\
\sphinxhline
\sphinxAtStartPar
{\hyperref[\detokenize{parameter:mipstartmode}]{\sphinxcrossref{\DUrole{std,std-ref}{MipStartMode}}}}
&
\sphinxAtStartPar
Integer
&
\sphinxAtStartPar
Mode of MIP starts
\\
\sphinxhline
\sphinxAtStartPar
{\hyperref[\detokenize{parameter:mipstartnodelimit}]{\sphinxcrossref{\DUrole{std,std-ref}{MipStartNodeLimit}}}}
&
\sphinxAtStartPar
Integer
&
\sphinxAtStartPar
Limit of nodes for MIP start sub\sphinxhyphen{}MIPs
\\
\sphinxhline
\sphinxAtStartPar
{\hyperref[\detokenize{parameter:linearizeindicators}]{\sphinxcrossref{\DUrole{std,std-ref}{LinearizeIndicators}}}}
&
\sphinxAtStartPar
Integer
&
\sphinxAtStartPar
Controls whether to force the linearization of indicator constraints
\\
\sphinxhline
\sphinxAtStartPar
{\hyperref[\detokenize{parameter:linearizesos}]{\sphinxcrossref{\DUrole{std,std-ref}{LinearizeSos}}}}
&
\sphinxAtStartPar
Integer
&
\sphinxAtStartPar
Controls whether to force the linearization of SOS constraints
\\
\sphinxhline
\sphinxAtStartPar
{\hyperref[\detokenize{parameter:miprepair}]{\sphinxcrossref{\DUrole{std,std-ref}{MipRepair}}}}
&
\sphinxAtStartPar
Integer
&
\sphinxAtStartPar
Level for repairing the MIP solution in case of numerical issues
\\
\sphinxbottomrule
\end{tabular}
\sphinxtableafterendhook\par
\sphinxattableend\end{savenotes}
\phantomsection\label{\detokenize{parameter:cutlevel}}\begin{itemize}
\item {} 
\sphinxAtStartPar
\sphinxcode{\sphinxupquote{CutLevel}}
\begin{quote}

\sphinxAtStartPar
Integer parameter.

\sphinxAtStartPar
Level of cutting\sphinxhyphen{}planes generation.

\sphinxAtStartPar
\sphinxstylestrong{Default:} \sphinxhyphen{}1

\sphinxAtStartPar
\sphinxstylestrong{Possible values:}
\begin{quote}

\sphinxAtStartPar
\sphinxhyphen{}1: Choose automatically.

\sphinxAtStartPar
0: Off

\sphinxAtStartPar
1: Fast

\sphinxAtStartPar
2: Normal

\sphinxAtStartPar
3: Aggressive
\end{quote}
\end{quote}

\end{itemize}
\phantomsection\label{\detokenize{parameter:rootcutlevel}}\begin{itemize}
\item {} 
\sphinxAtStartPar
\sphinxcode{\sphinxupquote{RootCutLevel}}
\begin{quote}

\sphinxAtStartPar
Integer parameter.

\sphinxAtStartPar
Level of cutting\sphinxhyphen{}planes generation of root node.

\sphinxAtStartPar
\sphinxstylestrong{Default:} \sphinxhyphen{}1

\sphinxAtStartPar
\sphinxstylestrong{Possible values:}
\begin{quote}

\sphinxAtStartPar
\sphinxhyphen{}1: Choose automatically.

\sphinxAtStartPar
0: Off

\sphinxAtStartPar
1: Fast

\sphinxAtStartPar
2: Normal

\sphinxAtStartPar
3: Aggressive
\end{quote}
\end{quote}

\end{itemize}
\phantomsection\label{\detokenize{parameter:treecutlevel}}\begin{itemize}
\item {} 
\sphinxAtStartPar
\sphinxcode{\sphinxupquote{TreeCutLevel}}
\begin{quote}

\sphinxAtStartPar
Integer parameter.

\sphinxAtStartPar
Level of cutting\sphinxhyphen{}planes generation of search tree.

\sphinxAtStartPar
\sphinxstylestrong{Default:} \sphinxhyphen{}1

\sphinxAtStartPar
\sphinxstylestrong{Possible values:}
\begin{quote}

\sphinxAtStartPar
\sphinxhyphen{}1: Choose automatically.

\sphinxAtStartPar
0: Off

\sphinxAtStartPar
1: Fast

\sphinxAtStartPar
2: Normal

\sphinxAtStartPar
3: Aggressive
\end{quote}
\end{quote}

\end{itemize}
\phantomsection\label{\detokenize{parameter:rootcutrounds}}\begin{itemize}
\item {} 
\sphinxAtStartPar
\sphinxcode{\sphinxupquote{RootCutRounds}}
\begin{quote}

\sphinxAtStartPar
Integer parameter.

\sphinxAtStartPar
Rounds of cutting\sphinxhyphen{}planes generation of root node.

\sphinxAtStartPar
\sphinxstylestrong{Default:} \sphinxhyphen{}1 (Choose automatically)

\sphinxAtStartPar
\sphinxstylestrong{Minimal:} \sphinxhyphen{}1 (Choose automatically)

\sphinxAtStartPar
\sphinxstylestrong{Maximal:} \sphinxcode{\sphinxupquote{INT\_MAX}}
\end{quote}

\end{itemize}
\phantomsection\label{\detokenize{parameter:nodecutrounds}}\begin{itemize}
\item {} 
\sphinxAtStartPar
\sphinxcode{\sphinxupquote{NodeCutRounds}}
\begin{quote}

\sphinxAtStartPar
Integer parameter.

\sphinxAtStartPar
Rounds of cutting\sphinxhyphen{}planes generation of search tree node.

\sphinxAtStartPar
\sphinxstylestrong{Default:} \sphinxhyphen{}1 (Choose automatically)

\sphinxAtStartPar
\sphinxstylestrong{Minimal:} \sphinxhyphen{}1 (Choose automatically)

\sphinxAtStartPar
\sphinxstylestrong{Maximal:} \sphinxcode{\sphinxupquote{INT\_MAX}}
\end{quote}

\end{itemize}
\phantomsection\label{\detokenize{parameter:heurlevel}}\begin{itemize}
\item {} 
\sphinxAtStartPar
\sphinxcode{\sphinxupquote{HeurLevel}}
\begin{quote}

\sphinxAtStartPar
Integer parameter.

\sphinxAtStartPar
Level of heuristics.

\sphinxAtStartPar
\sphinxstylestrong{Default:} \sphinxhyphen{}1

\sphinxAtStartPar
\sphinxstylestrong{Possible values:}
\begin{quote}

\sphinxAtStartPar
\sphinxhyphen{}1: Choose automatically.

\sphinxAtStartPar
0: Off

\sphinxAtStartPar
1: Fast

\sphinxAtStartPar
2: Normal

\sphinxAtStartPar
3: Aggressive
\end{quote}
\end{quote}

\end{itemize}
\phantomsection\label{\detokenize{parameter:prerootheurlevel}}\begin{itemize}
\item {} 
\sphinxAtStartPar
\sphinxcode{\sphinxupquote{PreRootHeurLevel}}
\begin{quote}

\sphinxAtStartPar
Integer parameter.

\sphinxAtStartPar
Level of pre\sphinxhyphen{}root heuristics.

\sphinxAtStartPar
\sphinxstylestrong{Default:} \sphinxhyphen{}1

\sphinxAtStartPar
\sphinxstylestrong{Possible values:}
\begin{quote}

\sphinxAtStartPar
\sphinxhyphen{}1: Choose automatically.

\sphinxAtStartPar
0: Off

\sphinxAtStartPar
1: Fast

\sphinxAtStartPar
2: Normal

\sphinxAtStartPar
3: Aggressive
\end{quote}
\end{quote}

\end{itemize}
\phantomsection\label{\detokenize{parameter:roundingheurlevel}}\begin{itemize}
\item {} 
\sphinxAtStartPar
\sphinxcode{\sphinxupquote{RoundingHeurLevel}}
\begin{quote}

\sphinxAtStartPar
Integer parameter.

\sphinxAtStartPar
Level of rounding heuristics.

\sphinxAtStartPar
\sphinxstylestrong{Default:} \sphinxhyphen{}1

\sphinxAtStartPar
\sphinxstylestrong{Possible values:}
\begin{quote}

\sphinxAtStartPar
\sphinxhyphen{}1: Choose automatically.

\sphinxAtStartPar
0: Off

\sphinxAtStartPar
1: Fast

\sphinxAtStartPar
2: Normal

\sphinxAtStartPar
3: Aggressive
\end{quote}
\end{quote}

\end{itemize}
\phantomsection\label{\detokenize{parameter:divingheurlevel}}\begin{itemize}
\item {} 
\sphinxAtStartPar
\sphinxcode{\sphinxupquote{DivingHeurLevel}}
\begin{quote}

\sphinxAtStartPar
Integer parameter.

\sphinxAtStartPar
Level of diving heuristics.

\sphinxAtStartPar
\sphinxstylestrong{Default:} \sphinxhyphen{}1

\sphinxAtStartPar
\sphinxstylestrong{Possible values:}
\begin{quote}

\sphinxAtStartPar
\sphinxhyphen{}1: Choose automatically.

\sphinxAtStartPar
0: Off

\sphinxAtStartPar
1: Fast

\sphinxAtStartPar
2: Normal

\sphinxAtStartPar
3: Aggressive
\end{quote}
\end{quote}

\end{itemize}
\phantomsection\label{\detokenize{parameter:submipheurlevel}}\begin{itemize}
\item {} 
\sphinxAtStartPar
\sphinxcode{\sphinxupquote{SubMipHeurLevel}}
\begin{quote}

\sphinxAtStartPar
Integer parameter.

\sphinxAtStartPar
Level of Sub\sphinxhyphen{}MIP heuristics.

\sphinxAtStartPar
\sphinxstylestrong{Default:} \sphinxhyphen{}1

\sphinxAtStartPar
\sphinxstylestrong{Possible values:}
\begin{quote}

\sphinxAtStartPar
\sphinxhyphen{}1: Choose automatically.

\sphinxAtStartPar
0: Off

\sphinxAtStartPar
1: Fast

\sphinxAtStartPar
2: Normal

\sphinxAtStartPar
3: Aggressive
\end{quote}
\end{quote}

\end{itemize}
\phantomsection\label{\detokenize{parameter:fapheurlevel}}\begin{itemize}
\item {} 
\sphinxAtStartPar
\sphinxcode{\sphinxupquote{FAPHeurLevel}}
\begin{quote}

\sphinxAtStartPar
Integer parameter.

\sphinxAtStartPar
Level of Fix\sphinxhyphen{}and\sphinxhyphen{}propagate heuristics.

\sphinxAtStartPar
\sphinxstylestrong{Default:} \sphinxhyphen{}1

\sphinxAtStartPar
\sphinxstylestrong{Possible values:}
\begin{quote}

\sphinxAtStartPar
\sphinxhyphen{}1: Choose automatically.

\sphinxAtStartPar
0: Off

\sphinxAtStartPar
1: Fast

\sphinxAtStartPar
2: Normal

\sphinxAtStartPar
3: Aggressive

\sphinxAtStartPar
4: Heavy
\end{quote}
\end{quote}

\end{itemize}
\phantomsection\label{\detokenize{parameter:strongbranching}}\begin{itemize}
\item {} 
\sphinxAtStartPar
\sphinxcode{\sphinxupquote{StrongBranching}}
\begin{quote}

\sphinxAtStartPar
Integer parameter.

\sphinxAtStartPar
Level of strong branching.

\sphinxAtStartPar
\sphinxstylestrong{Default:} \sphinxhyphen{}1

\sphinxAtStartPar
\sphinxstylestrong{Possible values:}
\begin{quote}

\sphinxAtStartPar
\sphinxhyphen{}1: Choose automatically.

\sphinxAtStartPar
0: Off

\sphinxAtStartPar
1: Fast

\sphinxAtStartPar
2: Normal

\sphinxAtStartPar
3: Aggressive
\end{quote}
\end{quote}

\end{itemize}
\phantomsection\label{\detokenize{parameter:conflictanalysis}}\begin{itemize}
\item {} 
\sphinxAtStartPar
\sphinxcode{\sphinxupquote{ConflictAnalysis}}
\begin{quote}

\sphinxAtStartPar
Integer parameter.

\sphinxAtStartPar
Whether to perform conflict analysis.

\sphinxAtStartPar
\sphinxstylestrong{Default:} \sphinxhyphen{}1

\sphinxAtStartPar
\sphinxstylestrong{Possible values:}
\begin{quote}

\sphinxAtStartPar
\sphinxhyphen{}1: Choose automatically.

\sphinxAtStartPar
0: No

\sphinxAtStartPar
1: Yes
\end{quote}
\end{quote}

\end{itemize}
\phantomsection\label{\detokenize{parameter:mipstartmode}}\begin{itemize}
\item {} 
\sphinxAtStartPar
\sphinxcode{\sphinxupquote{MipStartMode}}
\begin{quote}

\sphinxAtStartPar
Integer parameter.

\sphinxAtStartPar
Mode of MIP starts.

\sphinxAtStartPar
\sphinxstylestrong{Default:} \sphinxhyphen{}1

\sphinxAtStartPar
\sphinxstylestrong{Possible values:}
\begin{quote}

\sphinxAtStartPar
\sphinxhyphen{}1: Choose automatically.

\sphinxAtStartPar
0: Do not use any MIP starts.

\sphinxAtStartPar
1: Only load full and feasible MIP starts.

\sphinxAtStartPar
2: Only load feasible ones (complete partial solutions by solving subMIPs).
\end{quote}
\end{quote}

\end{itemize}
\phantomsection\label{\detokenize{parameter:mipstartnodelimit}}\begin{itemize}
\item {} 
\sphinxAtStartPar
\sphinxcode{\sphinxupquote{MipStartNodeLimit}}
\begin{quote}

\sphinxAtStartPar
Integer parameter.

\sphinxAtStartPar
Limit of nodes for MIP start sub\sphinxhyphen{}MIPs.

\sphinxAtStartPar
\sphinxstylestrong{Default:} \sphinxhyphen{}1 (Choose automatically)

\sphinxAtStartPar
\sphinxstylestrong{Minimal:} \sphinxhyphen{}1 (Choose automatically)

\sphinxAtStartPar
\sphinxstylestrong{Maximal:} \sphinxcode{\sphinxupquote{INT\_MAX}}
\end{quote}

\end{itemize}
\phantomsection\label{\detokenize{parameter:linearizeindicators}}\begin{itemize}
\item {} 
\sphinxAtStartPar
\sphinxcode{\sphinxupquote{LinearizeIndicators}}
\begin{quote}

\sphinxAtStartPar
Integer parameter.

\sphinxAtStartPar
Controls whether to force the linearization of Indicator constraints.

\sphinxAtStartPar
\sphinxstylestrong{Default:} 0

\sphinxAtStartPar
\sphinxstylestrong{Possible values}:
\begin{quote}

\sphinxAtStartPar
0: Do not force linearization of indicator constraints.

\sphinxAtStartPar
1: Force linearization of all indicator constraints.

\sphinxAtStartPar
(All indicator constraints will be transformed into equivalent linear constraints if enabled.)
\end{quote}
\end{quote}

\end{itemize}
\phantomsection\label{\detokenize{parameter:linearizesos}}\begin{itemize}
\item {} 
\sphinxAtStartPar
\sphinxcode{\sphinxupquote{LinearizeSos}}
\begin{quote}

\sphinxAtStartPar
Integer parameter.

\sphinxAtStartPar
Controls whether to force the linearization of SOS constraints.

\sphinxAtStartPar
\sphinxstylestrong{Default:} 0

\sphinxAtStartPar
\sphinxstylestrong{Possible values}:
\begin{quote}

\sphinxAtStartPar
0: Do not force linearization of SOS constraints.

\sphinxAtStartPar
1: Force linearization of all SOS constraints.

\sphinxAtStartPar
(All SOS constraints will be transformed into equivalent linear constraints if enabled.)
\end{quote}
\end{quote}

\end{itemize}
\phantomsection\label{\detokenize{parameter:miprepair}}\begin{itemize}
\item {} 
\sphinxAtStartPar
\sphinxcode{\sphinxupquote{MipRepair}}
\begin{quote}

\sphinxAtStartPar
Integer parameter.

\sphinxAtStartPar
Level for repairing the MIP solution in case of numerical issues.

\sphinxAtStartPar
\sphinxstylestrong{Default:} \sphinxhyphen{}1

\sphinxAtStartPar
\sphinxstylestrong{Possible values:}
\begin{quote}

\sphinxAtStartPar
\sphinxhyphen{}1: Only when time left.

\sphinxAtStartPar
0: Off.

\sphinxAtStartPar
1: Extend time limit and attempt repair (Fast).

\sphinxAtStartPar
2: Extend time limit and attempt repair (Normal).

\sphinxAtStartPar
3: Extend time limit and attempt repair (Aggressive).
\end{quote}
\end{quote}

\end{itemize}

\section{Semidefinite Programming related}
\label{\detokenize{parameter:semidefinite-programming-related}}\label{\detokenize{parameter:chapparam-sdp}}\begin{itemize}
\item {} 
\sphinxAtStartPar
\sphinxcode{\sphinxupquote{SDPMethod}}
\begin{quote}

\sphinxAtStartPar
Integer parameter.

\sphinxAtStartPar
Method to solve the SDP problem.

\sphinxAtStartPar
\sphinxstylestrong{Default:} \sphinxhyphen{}1

\sphinxAtStartPar
\sphinxstylestrong{Possible values:}
\begin{quote}

\sphinxAtStartPar
\sphinxhyphen{}1: Choose automatically.

\sphinxAtStartPar
0: Primal\sphinxhyphen{}Dual method.

\sphinxAtStartPar
1: ADMM method.

\sphinxAtStartPar
2: Dual method.
\end{quote}
\end{quote}

\end{itemize}

\section{Nonlinear Programming related}
\label{\detokenize{parameter:nonlinear-programming-related}}\label{\detokenize{parameter:chapparam-nlp}}
\sphinxAtStartPar
Nonlinear programming parameters control the workflow of nonlinear programming solvers.

\begin{savenotes}\sphinxattablestart
\sphinxthistablewithglobalstyle
\centering
\sphinxcapstartof{table}
\sphinxthecaptionisattop
\sphinxcaption{Overview of Nonlinear Programming Parameters}\label{\detokenize{parameter:copttab-paramnlp}}
\sphinxaftertopcaption
\begin{tabular}[t]{|\X{15}{65}|\X{10}{65}|\X{40}{65}|}
\sphinxtoprule
\sphinxstyletheadfamily 
\sphinxAtStartPar
Parameter
&\sphinxstyletheadfamily 
\sphinxAtStartPar
Type
&\sphinxstyletheadfamily 
\sphinxAtStartPar
Description
\\
\sphinxmidrule
\sphinxtableatstartofbodyhook
\sphinxAtStartPar
{\hyperref[\detokenize{parameter:nonconvex}]{\sphinxcrossref{\DUrole{std,std-ref}{NonConvex}}}}
&
\sphinxAtStartPar
Integer
&
\sphinxAtStartPar
Handling strategy for nonconvex models
\\
\sphinxhline
\sphinxAtStartPar
{\hyperref[\detokenize{parameter:nlpmuupdate}]{\sphinxcrossref{\DUrole{std,std-ref}{NLPMuUpdate}}}}
&
\sphinxAtStartPar
Integer
&
\sphinxAtStartPar
Barrier parameter update strategy of the nonlinear solver
\\
\sphinxhline
\sphinxAtStartPar
{\hyperref[\detokenize{parameter:nlplinscale}]{\sphinxcrossref{\DUrole{std,std-ref}{NLPLinScale}}}}
&
\sphinxAtStartPar
Integer
&
\sphinxAtStartPar
Linear system scaling strategy of the nonlinear solver
\\
\sphinxbottomrule
\end{tabular}
\sphinxtableafterendhook\par
\sphinxattableend\end{savenotes}
\phantomsection\label{\detokenize{parameter:nonconvex}}\begin{itemize}
\item {} 
\sphinxAtStartPar
\sphinxcode{\sphinxupquote{NonConvex}}
\begin{quote}

\sphinxAtStartPar
Integer parameter.

\sphinxAtStartPar
Handling strategy for continuous nonconvex models.

\sphinxAtStartPar
\sphinxstylestrong{Default:} \sphinxhyphen{}1

\sphinxAtStartPar
\sphinxstylestrong{Possible values:}
\begin{quote}

\sphinxAtStartPar
\sphinxhyphen{}1: Choose automatically.

\sphinxAtStartPar
0: Report nonconvexity and terminate.

\sphinxAtStartPar
1: Search for a local optimal solution.

\sphinxAtStartPar
2: Search for a global optimal solution.
\end{quote}
\end{quote}

\end{itemize}
\phantomsection\label{\detokenize{parameter:nlpmuupdate}}\begin{itemize}
\item {} 
\sphinxAtStartPar
\sphinxcode{\sphinxupquote{NLPMuUpdate}}
\begin{quote}

\sphinxAtStartPar
Integer parameter.

\sphinxAtStartPar
Barrier parameter update strategy of the nonlinear solver.

\sphinxAtStartPar
\sphinxstylestrong{Default:} \sphinxhyphen{}1

\sphinxAtStartPar
\sphinxstylestrong{Possible values:}
\begin{quote}

\sphinxAtStartPar
\sphinxhyphen{}1: Choose automatically.

\sphinxAtStartPar
0: Monotonically decreasing.

\sphinxAtStartPar
1: Adaptive adjustment.
\end{quote}
\end{quote}

\end{itemize}
\phantomsection\label{\detokenize{parameter:nlplinscale}}\begin{itemize}
\item {} 
\sphinxAtStartPar
\sphinxcode{\sphinxupquote{NLPLinScale}}
\begin{quote}

\sphinxAtStartPar
Integer parameter.

\sphinxAtStartPar
Linear system scaling strategy of the nonlinear solver.

\sphinxAtStartPar
\sphinxstylestrong{Default:} \sphinxhyphen{}1

\sphinxAtStartPar
\sphinxstylestrong{Possible values:}
\begin{quote}

\sphinxAtStartPar
\sphinxhyphen{}1: Choose automatically.

\sphinxAtStartPar
0: No scaling.

\sphinxAtStartPar
1: Always scale.
\end{quote}
\end{quote}

\end{itemize}

\section{Multi\sphinxhyphen{}objective Optimization related}
\label{\detokenize{parameter:multi-objective-optimization-related}}\label{\detokenize{parameter:chapparam-mobj}}\phantomsection\label{\detokenize{parameter:multiobjtimelimit}}\begin{itemize}
\item {} 
\sphinxAtStartPar
\sphinxcode{\sphinxupquote{MultiObjTimeLimit}}
\begin{quote}

\sphinxAtStartPar
Double parameter.

\sphinxAtStartPar
Time limit (in seconds) for solving the whole multi\sphinxhyphen{}objective model.

\sphinxAtStartPar
\sphinxstylestrong{Default:} 1e20

\sphinxAtStartPar
\sphinxstylestrong{Minimum:} 0

\sphinxAtStartPar
\sphinxstylestrong{Maximum:} 1e20
\end{quote}

\end{itemize}
\phantomsection\label{\detokenize{parameter:multiobjparammode}}\begin{itemize}
\item {} 
\sphinxAtStartPar
\sphinxcode{\sphinxupquote{MultiObjParamMode}}
\begin{quote}

\sphinxAtStartPar
Integer parameter.

\sphinxAtStartPar
Specifies how solver parameters are applied to models in multi\sphinxhyphen{}objective optimization.

\sphinxAtStartPar
\sphinxstylestrong{Default:} 0

\sphinxAtStartPar
\sphinxstylestrong{Possible values:}
\begin{quote}

\sphinxAtStartPar
0: Use the same parameters for all objectives.

\sphinxAtStartPar
1: Use individual solver parameters for each objective.
\end{quote}
\end{quote}

\end{itemize}

\section{Parallel computing related}
\label{\detokenize{parameter:parallel-computing-related}}\label{\detokenize{parameter:chapparam-threads}}

\begin{savenotes}\sphinxattablestart
\sphinxthistablewithglobalstyle
\centering
\sphinxcapstartof{table}
\sphinxthecaptionisattop
\sphinxcaption{Parallel computing related parameters}\label{\detokenize{parameter:id5}}
\sphinxaftertopcaption
\begin{tabular}[t]{|\X{15}{59}|\X{9}{59}|\X{35}{59}|}
\sphinxtoprule
\sphinxstyletheadfamily 
\sphinxAtStartPar
Name
&\sphinxstyletheadfamily 
\sphinxAtStartPar
Type
&\sphinxstyletheadfamily 
\sphinxAtStartPar
Description
\\
\sphinxmidrule
\sphinxtableatstartofbodyhook
\sphinxAtStartPar
{\hyperref[\detokenize{parameter:threads}]{\sphinxcrossref{\DUrole{std,std-ref}{Threads}}}}
&
\sphinxAtStartPar
Integer
&
\sphinxAtStartPar
Number of threads to use
\\
\sphinxhline
\sphinxAtStartPar
{\hyperref[\detokenize{parameter:barthreads}]{\sphinxcrossref{\DUrole{std,std-ref}{BarThreads}}}}
&
\sphinxAtStartPar
Integer
&
\sphinxAtStartPar
Number of threads used by barrier. If value is \sphinxhyphen{}1, the thread count is determined by parameter \sphinxcode{\sphinxupquote{Threads}}
\\
\sphinxhline
\sphinxAtStartPar
{\hyperref[\detokenize{parameter:simplexthreads}]{\sphinxcrossref{\DUrole{std,std-ref}{SimplexThreads}}}}
&
\sphinxAtStartPar
Integer
&
\sphinxAtStartPar
Number of threads used by dual simplex. If value is \sphinxhyphen{}1, the thread count is determined by parameter \sphinxcode{\sphinxupquote{Threads}}
\\
\sphinxhline
\sphinxAtStartPar
{\hyperref[\detokenize{parameter:crossoverthreads}]{\sphinxcrossref{\DUrole{std,std-ref}{CrossoverThreads}}}}
&
\sphinxAtStartPar
Integer
&
\sphinxAtStartPar
Number of threads used by crossover. If value is \sphinxhyphen{}1, the thread count is determined by parameter \sphinxcode{\sphinxupquote{Threads}}
\\
\sphinxhline
\sphinxAtStartPar
{\hyperref[\detokenize{parameter:miptasks}]{\sphinxcrossref{\DUrole{std,std-ref}{MipTasks}}}}
&
\sphinxAtStartPar
Integer
&
\sphinxAtStartPar
Number of MIP tasks in parallel
\\
\sphinxhline
\sphinxAtStartPar
{\hyperref[\detokenize{parameter:concurrentlpmode}]{\sphinxcrossref{\DUrole{std,std-ref}{ConcurrentLpMode}}}}
&
\sphinxAtStartPar
Integer
&
\sphinxAtStartPar
The LP concurrent solving mode
\\
\sphinxbottomrule
\end{tabular}
\sphinxtableafterendhook\par
\sphinxattableend\end{savenotes}
\phantomsection\label{\detokenize{parameter:threads}}\begin{itemize}
\item {} 
\sphinxAtStartPar
\sphinxcode{\sphinxupquote{Threads}}
\begin{quote}

\sphinxAtStartPar
Integer parameter.

\sphinxAtStartPar
Number of threads to use.

\sphinxAtStartPar
\sphinxstylestrong{Default:} \sphinxhyphen{}1 (Choose automatically)

\sphinxAtStartPar
\sphinxstylestrong{Minimal:} \sphinxhyphen{}1 (Choose automatically)

\sphinxAtStartPar
\sphinxstylestrong{Maximal:} 128
\end{quote}

\end{itemize}
\phantomsection\label{\detokenize{parameter:barthreads}}\begin{itemize}
\item {} 
\sphinxAtStartPar
\sphinxcode{\sphinxupquote{BarThreads}}
\begin{quote}

\sphinxAtStartPar
Integer parameter.

\sphinxAtStartPar
Number of threads used by barrier. If value is \sphinxhyphen{}1, the thread count is determined
by parameter \sphinxcode{\sphinxupquote{Threads}}.

\sphinxAtStartPar
\sphinxstylestrong{Default:} \sphinxhyphen{}1

\sphinxAtStartPar
\sphinxstylestrong{Minimal:} \sphinxhyphen{}1

\sphinxAtStartPar
\sphinxstylestrong{Maximal:} 128
\end{quote}

\end{itemize}
\phantomsection\label{\detokenize{parameter:simplexthreads}}\begin{itemize}
\item {} 
\sphinxAtStartPar
\sphinxcode{\sphinxupquote{SimplexThreads}}
\begin{quote}

\sphinxAtStartPar
Integer parameter.

\sphinxAtStartPar
Number of threads used by dual simplex. If value is \sphinxhyphen{}1, the thread count is determined
by parameter \sphinxcode{\sphinxupquote{Threads}}.

\sphinxAtStartPar
\sphinxstylestrong{Default:} \sphinxhyphen{}1

\sphinxAtStartPar
\sphinxstylestrong{Minimal:} \sphinxhyphen{}1

\sphinxAtStartPar
\sphinxstylestrong{Maximal:} 128
\end{quote}

\end{itemize}
\phantomsection\label{\detokenize{parameter:crossoverthreads}}\begin{itemize}
\item {} 
\sphinxAtStartPar
\sphinxcode{\sphinxupquote{CrossoverThreads}}
\begin{quote}

\sphinxAtStartPar
Integer parameter.

\sphinxAtStartPar
Number of threads used by crossover. If value is \sphinxhyphen{}1, the thread count is determined
by parameter \sphinxcode{\sphinxupquote{Threads}}.

\sphinxAtStartPar
\sphinxstylestrong{Default:} \sphinxhyphen{}1

\sphinxAtStartPar
\sphinxstylestrong{Minimal:} \sphinxhyphen{}1

\sphinxAtStartPar
\sphinxstylestrong{Maximal:} 128
\end{quote}

\end{itemize}
\phantomsection\label{\detokenize{parameter:miptasks}}\begin{itemize}
\item {} 
\sphinxAtStartPar
\sphinxcode{\sphinxupquote{MipTasks}}
\begin{quote}

\sphinxAtStartPar
Integer parameter.

\sphinxAtStartPar
Number of MIP tasks in parallel.

\sphinxAtStartPar
\sphinxstylestrong{Default:} \sphinxhyphen{}1 (Choose automatically)

\sphinxAtStartPar
\sphinxstylestrong{Minimal:} \sphinxhyphen{}1 (Choose automatically)

\sphinxAtStartPar
\sphinxstylestrong{Maximal:} 256
\end{quote}

\end{itemize}
\phantomsection\label{\detokenize{parameter:concurrentlpmode}}\begin{itemize}
\item {} 
\sphinxAtStartPar
\sphinxcode{\sphinxupquote{ConcurrentLpMode}}
\begin{quote}

\sphinxAtStartPar
Integer parameter.

\sphinxAtStartPar
The LP concurrent solving mode.

\sphinxAtStartPar
Only effective when \sphinxcode{\sphinxupquote{LpMethod = 4}}.

\sphinxAtStartPar
When \sphinxcode{\sphinxupquote{LpMethod = 4}} is enabled, the parameters \sphinxcode{\sphinxupquote{GPUMode}} and
\sphinxcode{\sphinxupquote{GPUDevice}} are ignored. GPU usage and device selection are
fully controlled by \sphinxcode{\sphinxupquote{ConcurrentLpMode}}.

\sphinxAtStartPar
\sphinxstylestrong{Default:} \sphinxhyphen{}1 (Choose automatically)

\sphinxAtStartPar
\sphinxstylestrong{Possible values:}
\begin{quote}

\sphinxAtStartPar
\sphinxhyphen{}1: Choose automatically.

\sphinxAtStartPar
0: CPU only (simplex + barrier).

\sphinxAtStartPar
1: CPU (simplex + barrier) + GPU first\sphinxhyphen{}order method PDLP
(uses GPU device 0).

\sphinxAtStartPar
2: CPU (simplex + barrier) + GPU barrier
(uses GPU device 0).

\sphinxAtStartPar
3: CPU (simplex + barrier) + GPU first\sphinxhyphen{}order method PDLP + GPU barrier
(uses GPU device 0 and GPU device 1).
\end{quote}
\end{quote}

\end{itemize}

\begin{sphinxadmonition}{note}{Note:}

\sphinxAtStartPar
If \sphinxcode{\sphinxupquote{ConcurrentLpMode}} is set to 1 or 2 and no GPU device is detected,
or if it is set to 3 and fewer than two GPU devices are detected,
the solver reports an error and terminates.
\end{sphinxadmonition}

\section{GPU computing related}
\label{\detokenize{parameter:gpu-computing-related}}\label{\detokenize{parameter:chapparam-gpu}}\phantomsection\label{\detokenize{parameter:gpumode}}\begin{itemize}
\item {} 
\sphinxAtStartPar
\sphinxcode{\sphinxupquote{GPUMode}}
\begin{quote}

\sphinxAtStartPar
Integer parameter.

\sphinxAtStartPar
Specifies GPU mode.

\sphinxAtStartPar
\sphinxstylestrong{Default:} \sphinxhyphen{}1

\sphinxAtStartPar
\sphinxstylestrong{Default:} \sphinxhyphen{}1

\sphinxAtStartPar
\sphinxstylestrong{Possible values:}
\begin{quote}

\sphinxAtStartPar
\sphinxhyphen{}1: Choose automatically. The first\sphinxhyphen{}order method (PDLP) will attempt to use the GPU by default,
while the barrier method will use the CPU by default.

\sphinxAtStartPar
0: Force CPU mode.

\sphinxAtStartPar
1: Attempt to use the standard GPU mode.

\sphinxAtStartPar
2: For the barrier method, attempt to use the high\sphinxhyphen{}performance GPU mode when solving LP problems,
which may lead to higher memory usage. For the first\sphinxhyphen{}order method (PDLP),
this is equivalent to \sphinxcode{\sphinxupquote{GPUMode=1}} (standard GPU mode).
\end{quote}
\end{quote}

\end{itemize}

\begin{sphinxadmonition}{note}{Notes}
\begin{enumerate}
\sphinxsetlistlabels{\arabic}{enumi}{enumii}{}{.}%
\item {} 
\sphinxAtStartPar
COPT’s GPU solver supports solving LP and the root relaxation of MILP problems.
It is only effective when using either the first\sphinxhyphen{}order method (PDLP) or the barrier method.

\item {} 
\sphinxAtStartPar
COPT’s GPU solver also supports solving SOCP, Q(C)P, Exponential Cone Programming, and SDP problems. It is only effective when using the barrier method.

\item {} 
\sphinxAtStartPar
On platforms such as Windows, Linux\sphinxhyphen{}x86, and Linux\sphinxhyphen{}aarch64, when \sphinxcode{\sphinxupquote{GPUMode}} is set to an appropriate value,
COPT will attempt to detect whether the required CUDA libraries are available and whether a supported GPU is present.
If both conditions are met, COPT will run in GPU mode; otherwise, it will fall back to CPU mode.
The startup messages differ slightly between the two modes.

\item {} 
\sphinxAtStartPar
On platforms such as macOS, only the CPU version of COPT is currently available. Even if \sphinxcode{\sphinxupquote{GPUMode=1}} or \sphinxcode{\sphinxupquote{2}} is set,
COPT will still run in CPU mode. The startup messages differ slightly between the two modes.

\item {} 
\sphinxAtStartPar
Installing CUDA libraries is not required for running COPT. If GPU mode is not enabled, COPT will function normally
regardless of whether the CUDA libraries are installed. For more information about installing CUDA libraries,
see {\hyperref[\detokenize{faq:chapfaq-gpu}]{\sphinxcrossref{\DUrole{std,std-ref}{FAQ \sphinxhyphen{} GPU Usage}}}}.

\end{enumerate}
\end{sphinxadmonition}
\phantomsection\label{\detokenize{parameter:gpudevice}}\begin{itemize}
\item {} 
\sphinxAtStartPar
\sphinxcode{\sphinxupquote{GPUDevice}}
\begin{quote}

\sphinxAtStartPar
Integer parameter.

\sphinxAtStartPar
Specifies GPU device to use (in cases where the running machine has multiple GPUs).

\sphinxAtStartPar
\sphinxstylestrong{Default:} \sphinxhyphen{}1 (Choose automatically)

\sphinxAtStartPar
\sphinxstylestrong{Minimal:} \sphinxhyphen{}1 (Choose automatically)

\sphinxAtStartPar
\sphinxstylestrong{Maximal:} \sphinxcode{\sphinxupquote{INT\_MAX}}
\end{quote}

\end{itemize}

\section{IIS computation related}
\label{\detokenize{parameter:iis-computation-related}}\label{\detokenize{parameter:chapparam-iis}}\begin{itemize}
\item {} 
\sphinxAtStartPar
\sphinxcode{\sphinxupquote{IISMethod}}
\begin{quote}

\sphinxAtStartPar
Integer parameter.

\sphinxAtStartPar
Method for IIS computation.

\sphinxAtStartPar
\sphinxstylestrong{Default:} \sphinxhyphen{}1

\sphinxAtStartPar
\sphinxstylestrong{Possible values:}
\begin{quote}

\sphinxAtStartPar
\sphinxhyphen{}1: Choose automatically.

\sphinxAtStartPar
0: Find smaller IIS.

\sphinxAtStartPar
1: Find IIS quickly.
\end{quote}
\end{quote}

\end{itemize}

\section{Feasibility relaxation related}
\label{\detokenize{parameter:feasibility-relaxation-related}}\label{\detokenize{parameter:chapparam-feas}}\begin{itemize}
\item {} 
\sphinxAtStartPar
\sphinxcode{\sphinxupquote{FeasRelaxMode}}
\begin{quote}

\sphinxAtStartPar
Integer parameter.

\sphinxAtStartPar
Method for feasibility relaxation.

\sphinxAtStartPar
\sphinxstylestrong{Default:} 0

\sphinxAtStartPar
\sphinxstylestrong{Possible values:}
\begin{quote}

\sphinxAtStartPar
0: Minimize sum of violations.

\sphinxAtStartPar
1: Optimize original objective function under minimal sum of violations.

\sphinxAtStartPar
2: Minimize number of violations.

\sphinxAtStartPar
3: Optimize original objective function under minimal number of violations.

\sphinxAtStartPar
4: Minimize sum of squared violations.

\sphinxAtStartPar
5: Optimize original objective function under minimal sum of squared violations.
\end{quote}
\end{quote}

\end{itemize}

\begin{sphinxadmonition}{note}{Notices:}
\begin{itemize}
\item {} 
\sphinxAtStartPar
Either \sphinxcode{\sphinxupquote{FeasRelaxMode=4}} or \sphinxcode{\sphinxupquote{FeasRelaxMode=5}}  does not support MILP but only LP. Other methods support both
two types of problems.

\end{itemize}
\end{sphinxadmonition}

\section{Parameter Tuning related}
\label{\detokenize{parameter:parameter-tuning-related}}\label{\detokenize{parameter:chapparam-tune}}

\begin{savenotes}\sphinxattablestart
\sphinxthistablewithglobalstyle
\centering
\sphinxcapstartof{table}
\sphinxthecaptionisattop
\sphinxcaption{Tuner related parameters}\label{\detokenize{parameter:id6}}
\sphinxaftertopcaption
\begin{tabular}[t]{|\X{15}{59}|\X{9}{59}|\X{35}{59}|}
\sphinxtoprule
\sphinxstyletheadfamily 
\sphinxAtStartPar
Name
&\sphinxstyletheadfamily 
\sphinxAtStartPar
Type
&\sphinxstyletheadfamily 
\sphinxAtStartPar
Description
\\
\sphinxmidrule
\sphinxtableatstartofbodyhook
\sphinxAtStartPar
{\hyperref[\detokenize{parameter:tunetimelimit}]{\sphinxcrossref{\DUrole{std,std-ref}{TuneTimeLimit}}}}
&
\sphinxAtStartPar
Double
&
\sphinxAtStartPar
Time limit for parameter tuning
\\
\sphinxhline
\sphinxAtStartPar
{\hyperref[\detokenize{parameter:tunetargettime}]{\sphinxcrossref{\DUrole{std,std-ref}{TuneTargetTime}}}}
&
\sphinxAtStartPar
Double
&
\sphinxAtStartPar
Time target for parameter tuning
\\
\sphinxhline
\sphinxAtStartPar
{\hyperref[\detokenize{parameter:tunetargetrelgap}]{\sphinxcrossref{\DUrole{std,std-ref}{TuneTargetRelGap}}}}
&
\sphinxAtStartPar
Double
&
\sphinxAtStartPar
Optimal relative tolerance target for parameter tuning
\\
\sphinxhline
\sphinxAtStartPar
{\hyperref[\detokenize{parameter:tunemethod}]{\sphinxcrossref{\DUrole{std,std-ref}{TuneMethod}}}}
&
\sphinxAtStartPar
Integer
&
\sphinxAtStartPar
Method for parameter tuning
\\
\sphinxhline
\sphinxAtStartPar
{\hyperref[\detokenize{parameter:tunemode}]{\sphinxcrossref{\DUrole{std,std-ref}{TuneMode}}}}
&
\sphinxAtStartPar
Integer
&
\sphinxAtStartPar
Mode for parameter tuning
\\
\sphinxhline
\sphinxAtStartPar
{\hyperref[\detokenize{parameter:tunemeasure}]{\sphinxcrossref{\DUrole{std,std-ref}{TuneMeasure}}}}
&
\sphinxAtStartPar
Integer
&
\sphinxAtStartPar
Parameter tuning result calculation method
\\
\sphinxhline
\sphinxAtStartPar
{\hyperref[\detokenize{parameter:tunepermutes}]{\sphinxcrossref{\DUrole{std,std-ref}{TunePermutes}}}}
&
\sphinxAtStartPar
Integer
&
\sphinxAtStartPar
Permutations for each trial parameter set
\\
\sphinxhline
\sphinxAtStartPar
{\hyperref[\detokenize{parameter:tuneoutputlevel}]{\sphinxcrossref{\DUrole{std,std-ref}{TuneOutputLevel}}}}
&
\sphinxAtStartPar
Integer
&
\sphinxAtStartPar
Parameter tuning log output intensity
\\
\sphinxbottomrule
\end{tabular}
\sphinxtableafterendhook\par
\sphinxattableend\end{savenotes}
\phantomsection\label{\detokenize{parameter:tunetimelimit}}\begin{itemize}
\item {} 
\sphinxAtStartPar
\sphinxcode{\sphinxupquote{TuneTimeLimit}}
\begin{quote}

\sphinxAtStartPar
Double parameter.

\sphinxAtStartPar
Time limit for parameter tuning. If the parameter value is 0, it will automatically set by the solver.

\sphinxAtStartPar
\sphinxstylestrong{Default:} 0

\sphinxAtStartPar
\sphinxstylestrong{Minimal:} 0

\sphinxAtStartPar
\sphinxstylestrong{Maximal:} 1e20
\end{quote}

\end{itemize}
\phantomsection\label{\detokenize{parameter:tunetargettime}}\begin{itemize}
\item {} 
\sphinxAtStartPar
\sphinxcode{\sphinxupquote{TuneTargetTime}}
\begin{quote}

\sphinxAtStartPar
Double parameter.

\sphinxAtStartPar
Time target for parameter tuning.

\sphinxAtStartPar
\sphinxstylestrong{Default:} 1e\sphinxhyphen{}2

\sphinxAtStartPar
\sphinxstylestrong{Minimal:} 0

\sphinxAtStartPar
\sphinxstylestrong{Maximal:} \sphinxcode{\sphinxupquote{DBL\_MAX}}
\end{quote}

\end{itemize}
\phantomsection\label{\detokenize{parameter:tunetargetrelgap}}\begin{itemize}
\item {} 
\sphinxAtStartPar
\sphinxcode{\sphinxupquote{TuneTargetRelGap}}
\begin{quote}

\sphinxAtStartPar
Double parameter.

\sphinxAtStartPar
Optimal relative tolerance target for parameter tuning.

\sphinxAtStartPar
\sphinxstylestrong{Default:} 1e\sphinxhyphen{}4

\sphinxAtStartPar
\sphinxstylestrong{Minimal:} 0

\sphinxAtStartPar
\sphinxstylestrong{Maximal:} \sphinxcode{\sphinxupquote{DBL\_MAX}}
\end{quote}

\end{itemize}
\phantomsection\label{\detokenize{parameter:tunemethod}}\begin{itemize}
\item {} 
\sphinxAtStartPar
\sphinxcode{\sphinxupquote{TuneMethod}}
\begin{quote}

\sphinxAtStartPar
Integer parameter.

\sphinxAtStartPar
Method for parameter tuning.

\sphinxAtStartPar
\sphinxstylestrong{Default:} \sphinxhyphen{}1

\sphinxAtStartPar
\sphinxstylestrong{Possible values:}
\begin{quote}

\sphinxAtStartPar
\sphinxhyphen{}1: Choose automatically.

\sphinxAtStartPar
0: Greedy search strategy.

\sphinxAtStartPar
1: Broader search strategy.
\end{quote}
\end{quote}

\end{itemize}
\phantomsection\label{\detokenize{parameter:tunemode}}\begin{itemize}
\item {} 
\sphinxAtStartPar
\sphinxcode{\sphinxupquote{TuneMode}}
\begin{quote}

\sphinxAtStartPar
Integer parameter.

\sphinxAtStartPar
Mode for parameter tuning.

\sphinxAtStartPar
\sphinxstylestrong{Default:} \sphinxhyphen{}1

\sphinxAtStartPar
\sphinxstylestrong{Possible values:}
\begin{quote}

\sphinxAtStartPar
\sphinxhyphen{}1: Choose automatically.

\sphinxAtStartPar
0: Solving time.

\sphinxAtStartPar
1: Optimal relative tolerance.

\sphinxAtStartPar
2: Objective function value.

\sphinxAtStartPar
3: The lower bound of the objective function value.
\end{quote}
\end{quote}

\end{itemize}
\phantomsection\label{\detokenize{parameter:tunemeasure}}\begin{itemize}
\item {} 
\sphinxAtStartPar
\sphinxcode{\sphinxupquote{TuneMeasure}}
\begin{quote}

\sphinxAtStartPar
Integer parameter.

\sphinxAtStartPar
Parameter tuning result calculation method.

\sphinxAtStartPar
\sphinxstylestrong{Default:} \sphinxhyphen{}1

\sphinxAtStartPar
\sphinxstylestrong{Possible values:}
\begin{quote}

\sphinxAtStartPar
\sphinxhyphen{}1: Choose automatically.

\sphinxAtStartPar
0: Calculate the average value.

\sphinxAtStartPar
1: Calculate the maximum value.
\end{quote}
\end{quote}

\end{itemize}
\phantomsection\label{\detokenize{parameter:tunepermutes}}\begin{itemize}
\item {} 
\sphinxAtStartPar
\sphinxcode{\sphinxupquote{TunePermutes}}
\begin{quote}

\sphinxAtStartPar
Integer parameter.

\sphinxAtStartPar
Permutations for each trial parameter set. If the parameter value is 0, it will automatically set by the solver.

\sphinxAtStartPar
\sphinxstylestrong{Default:} 0

\sphinxAtStartPar
\sphinxstylestrong{Minimal:} 0

\sphinxAtStartPar
\sphinxstylestrong{Maximal:} \sphinxcode{\sphinxupquote{INT\_MAX}}
\end{quote}

\end{itemize}
\phantomsection\label{\detokenize{parameter:tuneoutputlevel}}\begin{itemize}
\item {} 
\sphinxAtStartPar
\sphinxcode{\sphinxupquote{TuneOutputLevel}}
\begin{quote}

\sphinxAtStartPar
Integer parameter.

\sphinxAtStartPar
Parameter tuning log output level.

\sphinxAtStartPar
\sphinxstylestrong{Default:} 2

\sphinxAtStartPar
\sphinxstylestrong{Possible values:}
\begin{quote}

\sphinxAtStartPar
0: No output of tuning log.

\sphinxAtStartPar
1: Output only a summary of the improved parameters.

\sphinxAtStartPar
2: Output a summary of each tuning attempt.

\sphinxAtStartPar
3: Output a detailed log of each tuning attempt.
\end{quote}
\end{quote}

\end{itemize}

\section{Callback related}
\label{\detokenize{parameter:callback-related}}\label{\detokenize{parameter:chapparam-cbc}}\phantomsection\label{\detokenize{parameter:lazyconstraints}}\begin{itemize}
\item {} 
\sphinxAtStartPar
\sphinxcode{\sphinxupquote{LazyConstraints}}
\begin{quote}

\sphinxAtStartPar
Integer parameter.

\sphinxAtStartPar
Whether lazy constraints are part of the model.

\sphinxAtStartPar
\sphinxstylestrong{Default:} \sphinxhyphen{}1

\sphinxAtStartPar
\sphinxstylestrong{Possible values:}
\begin{quote}

\sphinxAtStartPar
\sphinxhyphen{}1: Choose automatically.

\sphinxAtStartPar
0: No.

\sphinxAtStartPar
1: Yes.
\end{quote}
\end{quote}

\end{itemize}

\begin{sphinxadmonition}{note}{Notes}
\begin{itemize}
\item {} 
\sphinxAtStartPar
This parameter only affects MIP.

\end{itemize}
\end{sphinxadmonition}

\section{Other parameters}
\label{\detokenize{parameter:other-parameters}}\label{\detokenize{parameter:chapparam-other}}\begin{itemize}
\item {} 
\sphinxAtStartPar
\sphinxcode{\sphinxupquote{Logging}}
\begin{quote}

\sphinxAtStartPar
Integer parameter.

\sphinxAtStartPar
Whether to print optimization logs.

\sphinxAtStartPar
\sphinxstylestrong{Default:} 1

\sphinxAtStartPar
\sphinxstylestrong{Possible values:}
\begin{quote}

\sphinxAtStartPar
0: No optimization logs.

\sphinxAtStartPar
1: Print optimization logs.
\end{quote}
\end{quote}

\item {} 
\sphinxAtStartPar
\sphinxcode{\sphinxupquote{LogLevel}}
\begin{quote}

\sphinxAtStartPar
Integer parameter.

\sphinxAtStartPar
Controls the level of detail in the optimization logs.

\sphinxAtStartPar
\sphinxstylestrong{Default:} 2

\sphinxAtStartPar
\sphinxstylestrong{Possible values:}
\begin{quote}

\sphinxAtStartPar
2: Print basic optimization logs.

\sphinxAtStartPar
3: Print memory usage information in addition to basic optimization logs (for MIP problems).
\end{quote}
\end{quote}

\item {} 
\sphinxAtStartPar
\sphinxcode{\sphinxupquote{LogToConsole}}
\begin{quote}

\sphinxAtStartPar
Integer parameter.

\sphinxAtStartPar
Whether to print optimization logs to console.

\sphinxAtStartPar
\sphinxstylestrong{Default:} 1

\sphinxAtStartPar
\sphinxstylestrong{Possible values:}
\begin{quote}

\sphinxAtStartPar
0: No optimization logs to console.

\sphinxAtStartPar
1: Print optimization logs to console.
\end{quote}
\end{quote}

\end{itemize}

\section{Methods for accessing and setting parameters}
\label{\detokenize{parameter:methods-for-accessing-and-setting-parameters}}\label{\detokenize{parameter:chapparam-method}}
\sphinxAtStartPar
In different programming interfaces, the ways to access and set parameters are slightly different.
For details, please refer to the corresponding chapters for each programming language API:
\begin{itemize}
\item {} 
\sphinxAtStartPar
C API: {\hyperref[\detokenize{capiref:chapapi-getparam}]{\sphinxcrossref{\DUrole{std,std-ref}{C API Functions: Accessing and setting parameters}}}}

\item {} 
\sphinxAtStartPar
C++ API: {\hyperref[\detokenize{cppapiref:chapcppapiref-params}]{\sphinxcrossref{\DUrole{std,std-ref}{C++ API Reference: Parameters}}}}

\item {} 
\sphinxAtStartPar
C\# API: {\hyperref[\detokenize{csharpapiref:chapcsharpapiref-params}]{\sphinxcrossref{\DUrole{std,std-ref}{C\# API Reference: Parameters}}}}

\item {} 
\sphinxAtStartPar
Java API: {\hyperref[\detokenize{javaapiref:chapjavaapiref-params}]{\sphinxcrossref{\DUrole{std,std-ref}{Java API Reference: Parameters}}}}

\item {} 
\sphinxAtStartPar
Python API: {\hyperref[\detokenize{pyapiref:chappyapi-const-param}]{\sphinxcrossref{\DUrole{std,std-ref}{Python API Reference: Parameters}}}}

\end{itemize}

\sphinxstepscope

\chapter{Modeling and Solving Optimization Problems}
\label{\detokenize{modeling:modeling-and-solving-optimization-problems}}\label{\detokenize{modeling:chapmodeling}}\label{\detokenize{modeling::doc}}
\sphinxAtStartPar
This chapter will introduce how to model and solve different types of optimization problems using the COPT solver. The content of the chapter is structured as follows:
\begin{itemize}
\item {} 
\sphinxAtStartPar
{\hyperref[\detokenize{modeling:chapmodeling-lp}]{\sphinxcrossref{\DUrole{std,std-ref}{Linear Programming (LP)}}}}

\item {} 
\sphinxAtStartPar
{\hyperref[\detokenize{modeling:chapmodeling-socp}]{\sphinxcrossref{\DUrole{std,std-ref}{Second\sphinxhyphen{}Order Cone Programming (SOCP)}}}}

\item {} 
\sphinxAtStartPar
{\hyperref[\detokenize{modeling:chapmodeling-expcone}]{\sphinxcrossref{\DUrole{std,std-ref}{Exponential Cone Programming (ExpCone Programming)}}}}

\item {} 
\sphinxAtStartPar
{\hyperref[\detokenize{modeling:chapmodeling-sdp}]{\sphinxcrossref{\DUrole{std,std-ref}{Semidefinite Programming (SDP)}}}}

\item {} 
\sphinxAtStartPar
{\hyperref[\detokenize{modeling:chapmodeling-qp}]{\sphinxcrossref{\DUrole{std,std-ref}{Quadratic Programming (QP)}}}}

\item {} 
\sphinxAtStartPar
{\hyperref[\detokenize{modeling:chapmodeling-qcp}]{\sphinxcrossref{\DUrole{std,std-ref}{Quadratically Constrained Programming (QCP)}}}}

\item {} 
\sphinxAtStartPar
{\hyperref[\detokenize{modeling:chapmodeling-nlp}]{\sphinxcrossref{\DUrole{std,std-ref}{General Nonlinear Programming (NLP)}}}}

\item {} 
\sphinxAtStartPar
{\hyperref[\detokenize{modeling:chapmodeling-mip}]{\sphinxcrossref{\DUrole{std,std-ref}{Mixed\sphinxhyphen{}Integer Programming (MIP)}}}}

\item {} 
\sphinxAtStartPar
{\hyperref[\detokenize{modeling:chapmodeling-spec}]{\sphinxcrossref{\DUrole{std,std-ref}{Special Constraints}}}}

\end{itemize}

\sphinxAtStartPar
The types of problems and corresponding optimization algorithms currently supported by the COPT solver are shown in \hyperref[\detokenize{intro:copttab-probsalgs}]{Table \ref{\detokenize{intro:copttab-probsalgs}}}:

\begin{savenotes}\sphinxattablestart
\sphinxthistablewithglobalstyle
\centering
\sphinxcapstartof{table}
\sphinxthecaptionisattop
\sphinxcaption{Supported problem types and available algorithms}\label{\detokenize{modeling:coptmodel-probsalgs}}
\sphinxaftertopcaption
\begin{tabulary}{\linewidth}[t]{|T|T|}
\sphinxtoprule
\sphinxtableatstartofbodyhook
\sphinxAtStartPar
\sphinxstylestrong{Problem type}
&
\sphinxAtStartPar
Available algorithms
\\
\sphinxhline
\sphinxAtStartPar
Linear Programming (LP)
&
\sphinxAtStartPar
Simplex, Barrier (CPU/GPU),
First\sphinxhyphen{}order Method (PDLP)(CPU/GPU)
\\
\sphinxhline
\sphinxAtStartPar
Second\sphinxhyphen{}Order\sphinxhyphen{}Cone Programming (SOCP)
&
\sphinxAtStartPar
Barrier (CPU/GPU)
\\
\sphinxhline
\sphinxAtStartPar
Exponential Cone Programming (ExpCone)
&
\sphinxAtStartPar
Barrier (CPU/GPU)
\\
\sphinxhline
\sphinxAtStartPar
Convex Quadratic Programming (QP)
&
\sphinxAtStartPar
Barrier (CPU/GPU)
\\
\sphinxhline
\sphinxAtStartPar
Convex Quadratically Constrained Programming (QCP)
&
\sphinxAtStartPar
Barrier (CPU/GPU)
\\
\sphinxhline
\sphinxAtStartPar
Semidefinite Programming (SDP)
&
\sphinxAtStartPar
Barrier (CPU/GPU), ADMM
\\
\sphinxhline
\sphinxAtStartPar
Nonconvex Quadratic Programming (QP)
&
\sphinxAtStartPar
Spatial Branch\sphinxhyphen{}and\sphinxhyphen{}Bound (global optimum),
Barrier (local optimum)
\\
\sphinxhline
\sphinxAtStartPar
Nonconvex Quadratically Constrained Programming (QCP)
&
\sphinxAtStartPar
Spatial Branch\sphinxhyphen{}and\sphinxhyphen{}Bound (global optimum),
Barrier (local optimum)
\\
\sphinxhline
\sphinxAtStartPar
General Nonlinear Programming (NLP)
&
\sphinxAtStartPar
Barrier
\\
\sphinxhline
\sphinxAtStartPar
Mixed Integer Linear Programming (MILP)
&
\sphinxAtStartPar
Branch\sphinxhyphen{}and\sphinxhyphen{}Cut
\\
\sphinxhline
\sphinxAtStartPar
Mixed Integer Second\sphinxhyphen{}Order\sphinxhyphen{}Cone Programming (MISOCP)
&
\sphinxAtStartPar
Branch\sphinxhyphen{}and\sphinxhyphen{}Cut
\\
\sphinxhline
\sphinxAtStartPar
Mixed Integer Convex/Nonconvex Quadratic Programming (MIQP)
&
\sphinxAtStartPar
Branch\sphinxhyphen{}and\sphinxhyphen{}Cut
\\
\sphinxhline
\sphinxAtStartPar
Mixed Integer Convex/Nonconvex
Quadratically Constrained Programming (MIQCP)
&
\sphinxAtStartPar
Branch\sphinxhyphen{}and\sphinxhyphen{}Cut
\\
\sphinxbottomrule
\end{tabulary}
\sphinxtableafterendhook\par
\sphinxattableend\end{savenotes}

\sphinxAtStartPar
COPT supports three methods of inputting optimization problems: model file reader, programming interface,
and third\sphinxhyphen{}party tools.
\begin{enumerate}
\sphinxsetlistlabels{\arabic}{enumi}{enumii}{}{.}%
\item {} 
\sphinxAtStartPar
For model file input, please refer to {\hyperref[\detokenize{fileformats:chapfileformat}]{\sphinxcrossref{\DUrole{std,std-ref}{File Formats}}}}.

\item {} 
\sphinxAtStartPar
The programming language interfaces supported by COPT include:
\begin{itemize}
\item {} 
\sphinxAtStartPar
C

\item {} 
\sphinxAtStartPar
Python

\item {} 
\sphinxAtStartPar
C++

\item {} 
\sphinxAtStartPar
C\#

\item {} 
\sphinxAtStartPar
Java

\item {} 
\sphinxAtStartPar
Fortran

\end{itemize}

\item {} 
\sphinxAtStartPar
The main third\sphinxhyphen{}party tools interfaces currently supported by COPT include:
\begin{itemize}
\item {} 
\sphinxAtStartPar
AMPL

\item {} 
\sphinxAtStartPar
AIMMS

\item {} 
\sphinxAtStartPar
GAMS

\item {} 
\sphinxAtStartPar
Julia

\item {} 
\sphinxAtStartPar
Pyomo

\item {} 
\sphinxAtStartPar
PuLP

\item {} 
\sphinxAtStartPar
CVXPY

\end{itemize}

\end{enumerate}

\section{Linear Programming (LP)}
\label{\detokenize{modeling:linear-programming-lp}}\label{\detokenize{modeling:chapmodeling-lp}}
\sphinxAtStartPar
Linear Programming, as the most fundamental and important branch of operations research, has a wide range of applications.
The objective function and constraints of a linear programming problem are both linear.

\subsection{Mathematical Model}
\label{\detokenize{modeling:mathematical-model}}
\sphinxAtStartPar
The mathematical formulation is as follows:
\begin{equation}\label{equation:modeling:modeling:0}
\begin{split}\min\quad &\sum_{i=0}^{m-1}c_jx_j +c^f\\
\mathrm{s.t.}\quad &l_i^c\leq \sum_{j=0}^{n-1}a_{ij}x_j\leq u_i^c,\qquad i=1,2,\cdots,m-1 \\
&l_j^v\leq x_j\leq u_j^v,\qquad j=1,2,\cdots,n-1\end{split}
\end{equation}
\sphinxAtStartPar
Alternatively, the model can be concisely expressed using vectors and matrices:
\begin{equation}\label{equation:modeling:modeling:1}
\begin{split}\min\quad &c^Tx +c^f\\
\mathrm{s.t.}\quad &l^c\leq Ax\leq u^c \\
&l^v\leq x\leq u^v\end{split}
\end{equation}
\sphinxAtStartPar
The variables and parameters in the model have the following meanings:
\begin{itemize}
\item {} 
\sphinxAtStartPar
Problem size: \(m\) represents the number of constraints, and \(n\) represents the number of decision variables.

\item {} 
\sphinxAtStartPar
Decision variables: \(x=(x_j)_{j=0}^{n-1}\in\mathbb{R}^n\)

\item {} 
\sphinxAtStartPar
Variable bounds: \(l^v, u^v\in\mathbb{R}^n\), where \(l^v\) represents the lower bounds and \(u^v\)
represents the upper bounds of the variables.

\item {} 
\sphinxAtStartPar
Constraint bounds: \(l^c, u^c\in\mathbb{R}^m\), where \(l^c\) represents the lower bounds and \(u^c\)
represents the upper bounds of the constraints.

\item {} 
\sphinxAtStartPar
Coefficient matrix for linear constraints: \(A=(a_{ij})_{m\times n}\in\mathbb{R}^{m\times n}\)

\item {} 
\sphinxAtStartPar
Coefficients of variables in the objective function: \(c\in\mathbb{R}^n\) represents the coefficients of the variables
in the objective function, and \(c^f\) represents the constant term in the objective function.

\end{itemize}

\subsection{Modeling}
\label{\detokenize{modeling:modeling}}
\sphinxAtStartPar
The basic steps to build and solve a Linear Programming (LP) model in COPT are as follows:
\begin{enumerate}
\sphinxsetlistlabels{\arabic}{enumi}{enumii}{}{.}%
\item {} 
\sphinxAtStartPar
Create the COPT environment and model.

\item {} 
\sphinxAtStartPar
Add required data.

\item {} 
\sphinxAtStartPar
Construct the linear programming model:
\begin{itemize}
\item {} 
\sphinxAtStartPar
Add decision variables.

\item {} 
\sphinxAtStartPar
Add linear constraints.

\item {} 
\sphinxAtStartPar
Set the linear objective function.

\end{itemize}

\item {} 
\sphinxAtStartPar
Set optimization parameters and solve the model.

\item {} 
\sphinxAtStartPar
Retrieve the solution results.

\end{enumerate}

\sphinxAtStartPar

\sphinxAtStartPar
\sphinxstylestrong{Modeling: Adding Linear Constraints}

\sphinxAtStartPar
The COPT provides three methods to add linear constraints:
\begin{enumerate}
\sphinxsetlistlabels{\arabic}{enumi}{enumii}{}{.}%
\item {} 
\sphinxAtStartPar
Add a single linear constraint to the model.

\item {} 
\sphinxAtStartPar
Batch add a group of linear constraints.

\item {} 
\sphinxAtStartPar
Add a single linear constraint with lower and upper bounds.

\end{enumerate}

\sphinxAtStartPar
When adding linear constraints to the model, the main parameters that can be specified are:
\begin{itemize}
\item {} 
\sphinxAtStartPar
\sphinxcode{\sphinxupquote{expr}}/ \sphinxcode{\sphinxupquote{builder}} : Linear constraint expression or linear constraint builder.

\item {} 
\sphinxAtStartPar
\sphinxcode{\sphinxupquote{sense}}: Type of constraint. For possible values, please refer to {\hyperref[\detokenize{constant:chapconst-constrtype}]{\sphinxcrossref{\DUrole{std,std-ref}{Constants Section: Constraint Type}}}}.

\item {} 
\sphinxAtStartPar
\sphinxcode{\sphinxupquote{name}}: Name of the linear constraint.

\end{itemize}

\sphinxAtStartPar
If adding a linear constraint with bounds, the following must also be specified:
\begin{itemize}
\item {} 
\sphinxAtStartPar
\sphinxcode{\sphinxupquote{lb}}: Lower bound of the linear constraint.

\item {} 
\sphinxAtStartPar
\sphinxcode{\sphinxupquote{ub}}: Upper bound of the linear constraint.

\end{itemize}

\sphinxAtStartPar
The implementation methods in different programming interfaces are shown in \hyperref[\detokenize{modeling:copttab-modelinglpmodel-constr}]{Table \ref{\detokenize{modeling:copttab-modelinglpmodel-constr}}}:

\begin{savenotes}\sphinxattablestart
\sphinxthistablewithglobalstyle
\centering
\sphinxcapstartof{table}
\sphinxthecaptionisattop
\sphinxcaption{Functions for Adding Linear Constraints}\label{\detokenize{modeling:copttab-modelinglpmodel-constr}}
\sphinxaftertopcaption
\begin{tabular}[t]{|\X{5}{25}|\X{10}{25}|\X{10}{25}|}
\sphinxtoprule
\sphinxstyletheadfamily 
\sphinxAtStartPar
API
&\sphinxstyletheadfamily 
\sphinxAtStartPar
Add a single constraint
&\sphinxstyletheadfamily 
\sphinxAtStartPar
Add a group of constraints in batch
\\
\sphinxmidrule
\sphinxtableatstartofbodyhook
\sphinxAtStartPar
C
&
\sphinxAtStartPar
\sphinxcode{\sphinxupquote{COPT\_AddRow}}
&
\sphinxAtStartPar
\sphinxcode{\sphinxupquote{COPT\_AddRows}}
\\
\sphinxhline
\sphinxAtStartPar
C++
&
\sphinxAtStartPar
\sphinxcode{\sphinxupquote{Model::AddConstr()}}
&
\sphinxAtStartPar
\sphinxcode{\sphinxupquote{Model::AddConstrs()}}
\\
\sphinxhline
\sphinxAtStartPar
C\#
&
\sphinxAtStartPar
\sphinxcode{\sphinxupquote{Model.AddConstr()}}
&
\sphinxAtStartPar
\sphinxcode{\sphinxupquote{Model.AddConstrs()}}
\\
\sphinxhline
\sphinxAtStartPar
Java
&
\sphinxAtStartPar
\sphinxcode{\sphinxupquote{Model.addConstr()}}
&
\sphinxAtStartPar
\sphinxcode{\sphinxupquote{Model.addConstrs()}}
\\
\sphinxhline
\sphinxAtStartPar
Python
&
\sphinxAtStartPar
\sphinxcode{\sphinxupquote{Model.addConstr()}}
&
\sphinxAtStartPar
\sphinxcode{\sphinxupquote{Model.addConstrs()}}
\\
\sphinxbottomrule
\end{tabular}
\sphinxtableafterendhook\par
\sphinxattableend\end{savenotes}

\begin{sphinxadmonition}{note}{Notes}
\begin{itemize}
\item {} 
\sphinxAtStartPar
For linear constraint operations, the function names, calling methods, and \sphinxstylestrong{argument names may vary slightly}
across different programming interfaces, but the functionality and argument meanings are consistent.

\item {} 
\sphinxAtStartPar
COPT supports defining constraints using constraint types, but we \sphinxstylestrong{recommend defining constraints directly using bounds}.

\item {} 
\sphinxAtStartPar
In the C API, linear constraints are added using the coefficient matrix as an input argument.

\item {} 
\sphinxAtStartPar
In the Python API, an additional method for adding linear constraints with bounds is provided: \sphinxcode{\sphinxupquote{Model.addBoundConstr()}}.

\end{itemize}
\end{sphinxadmonition}

\sphinxAtStartPar

\sphinxAtStartPar
\sphinxstylestrong{Linear Programming Problem Example}
\begin{equation}\label{equation:modeling:modeling:2}
\begin{split}\text{Maximize:} & \\
                    & 1.2 x + 1.8 y + 2.1 z \\
\text{Subject to:} & \\
                      & 1.5 x + 1.2 y + 1.8 z \leq 2.6 \\
                      & 0.8 x + 0.6 y + 0.9 z \geq 1.2 \\
\text{Bounds:} & \\
             & 0.1 \leq x \leq 0.6 \\
             & 0.2 \leq y \leq 1.5 \\
             & 0.3 \leq z \leq 2.8\end{split}
\end{equation}
\sphinxAtStartPar
For corresponding code implementations in different programming interfaces, please refer to: {\hyperref[\detokenize{quickstart:chapquickstart}]{\sphinxcrossref{\DUrole{std,std-ref}{COPT Quick Start Section}}}}.

\sphinxAtStartPar

\sphinxAtStartPar
In the programming interfaces provided by COPT, except for the C language, the other object\sphinxhyphen{}oriented programming interfaces
(C\#, C++, Java, Python) offer classes related to linear constraints:
\begin{itemize}
\item {} 
\sphinxAtStartPar
Encapsulation of operations related to linear constraints:
\begin{enumerate}
\sphinxsetlistlabels{\arabic}{enumi}{enumii}{}{.}%
\item {} 
\sphinxAtStartPar
\sphinxcode{\sphinxupquote{Constraint}} class: Encapsulation of operations related to linear constraints in COPT.

\item {} 
\sphinxAtStartPar
\sphinxcode{\sphinxupquote{ConstrArray}} class: Facilitates operations on a group of \sphinxcode{\sphinxupquote{Constraint}} class objects.

\end{enumerate}

\item {} 
\sphinxAtStartPar
Encapsulation of linear constraint builders:
\begin{enumerate}
\sphinxsetlistlabels{\arabic}{enumi}{enumii}{}{.}%
\item {} 
\sphinxAtStartPar
\sphinxcode{\sphinxupquote{ConstrBuilder}} class: Encapsulation of linear constraint builders in COPT.

\item {} 
\sphinxAtStartPar
\sphinxcode{\sphinxupquote{ConstrBuilderArray}} class: Facilitates operations on a group of \sphinxcode{\sphinxupquote{ConstrBuilder}} class objects.

\end{enumerate}

\item {} 
\sphinxAtStartPar
C++ API: {\hyperref[\detokenize{cppapiref:chapcppapiref-constraint}]{\sphinxcrossref{\DUrole{std,std-ref}{Constraint Class}}}} , {\hyperref[\detokenize{cppapiref:chapcppapiref-constrarray}]{\sphinxcrossref{\DUrole{std,std-ref}{ConstrArray Class}}}} , {\hyperref[\detokenize{cppapiref:chapcppapiref-constrbuilder}]{\sphinxcrossref{\DUrole{std,std-ref}{ConstrBuilder Class}}}} , {\hyperref[\detokenize{cppapiref:chapcppapiref-constrbuilderarray}]{\sphinxcrossref{\DUrole{std,std-ref}{ConstrBuilderArray Class}}}}

\item {} 
\sphinxAtStartPar
C\# API: {\hyperref[\detokenize{csharpapiref:chapcsharpapiref-constraint}]{\sphinxcrossref{\DUrole{std,std-ref}{Constraint Class}}}} , {\hyperref[\detokenize{csharpapiref:chapcsharpapiref-constrarray}]{\sphinxcrossref{\DUrole{std,std-ref}{ConstrArray Class}}}} , {\hyperref[\detokenize{csharpapiref:chapcsharpapiref-constrbuilder}]{\sphinxcrossref{\DUrole{std,std-ref}{ConstrBuilder Class}}}} , {\hyperref[\detokenize{csharpapiref:chapcsharpapiref-constrbuilderarray}]{\sphinxcrossref{\DUrole{std,std-ref}{ConstrBuilderArray Class}}}}

\item {} 
\sphinxAtStartPar
Java API: {\hyperref[\detokenize{javaapiref:chapjavaapiref-constraint}]{\sphinxcrossref{\DUrole{std,std-ref}{Constraint Class}}}} , {\hyperref[\detokenize{javaapiref:chapjavaapiref-constrarray}]{\sphinxcrossref{\DUrole{std,std-ref}{ConstrArray Class}}}} , {\hyperref[\detokenize{javaapiref:chapjavaapiref-constrbuilder}]{\sphinxcrossref{\DUrole{std,std-ref}{ConstrBuilder Class}}}} , {\hyperref[\detokenize{javaapiref:chapjavaapiref-constrbuilderarray}]{\sphinxcrossref{\DUrole{std,std-ref}{ConstrBuilderArray Class}}}}

\item {} 
\sphinxAtStartPar
Python API: {\hyperref[\detokenize{pyapiref:chappyapi-constraint}]{\sphinxcrossref{\DUrole{std,std-ref}{Constraint Class}}}} , {\hyperref[\detokenize{pyapiref:chappyapi-constrarray}]{\sphinxcrossref{\DUrole{std,std-ref}{ConstrArray Class}}}} , {\hyperref[\detokenize{pyapiref:chappyapi-constrbuilder}]{\sphinxcrossref{\DUrole{std,std-ref}{ConstrBuilder Class}}}} , {\hyperref[\detokenize{pyapiref:chappyapi-constrbuilderarray}]{\sphinxcrossref{\DUrole{std,std-ref}{ConstrBuilderArray Class}}}}

\end{itemize}

\subsection{Solving}
\label{\detokenize{modeling:solving}}
\sphinxAtStartPar
For linear programming problems, the COPT solver provides the Simplex method and Barrier method.
The specific method can be selected by setting the optimization parameter \sphinxcode{\sphinxupquote{"LpMethod"}}.
By configuring other related optimization parameters for linear programming, you can control the detailed workflow of the solving algorithm.
For more details, please refer to {\hyperref[\detokenize{parameter:chapparam-lp}]{\sphinxcrossref{\DUrole{std,std-ref}{Parameter: Linear Programming Related}}}}.

\sphinxAtStartPar
For the solving logs of linear programming problems, please refer to {\hyperref[\detokenize{logging:chaplogging-simplex}]{\sphinxcrossref{\DUrole{std,std-ref}{Logging Section: Simplex Method}}}}
and {\hyperref[\detokenize{logging:chaplogging-barrier}]{\sphinxcrossref{\DUrole{std,std-ref}{Interior Point Method}}}}.

\subsection{Related Attributes and Information}
\label{\detokenize{modeling:related-attributes-and-information}}
\sphinxAtStartPar
\sphinxstylestrong{Linear Programming Related Attributes}

\sphinxAtStartPar
Attributes for linear programming are shown in \hyperref[\detokenize{modeling:copttab-modelinglp-attr}]{Table \ref{\detokenize{modeling:copttab-modelinglp-attr}}}:

\begin{savenotes}\sphinxattablestart
\sphinxthistablewithglobalstyle
\centering
\sphinxcapstartof{table}
\sphinxthecaptionisattop
\sphinxcaption{Attributes for Linear Programming}\label{\detokenize{modeling:copttab-modelinglp-attr}}
\sphinxaftertopcaption
\begin{tabular}[t]{|\X{15}{70}|\X{10}{70}|\X{45}{70}|}
\sphinxtoprule
\sphinxstyletheadfamily 
\sphinxAtStartPar
Name
&\sphinxstyletheadfamily 
\sphinxAtStartPar
Type
&\sphinxstyletheadfamily 
\sphinxAtStartPar
Description
\\
\sphinxmidrule
\sphinxtableatstartofbodyhook
\sphinxAtStartPar
\sphinxcode{\sphinxupquote{Cols}}
&
\sphinxAtStartPar
Integer
&
\sphinxAtStartPar
Number of variables (columns) in the problem
\\
\sphinxhline
\sphinxAtStartPar
\sphinxcode{\sphinxupquote{Rows}}
&
\sphinxAtStartPar
Integer
&
\sphinxAtStartPar
Number of constraints (rows) in the problem
\\
\sphinxhline
\sphinxAtStartPar
\sphinxcode{\sphinxupquote{Elems}}
&
\sphinxAtStartPar
Integer
&
\sphinxAtStartPar
Number of non\sphinxhyphen{}zero elements in the coefficient matrix
\\
\sphinxhline
\sphinxAtStartPar
\sphinxcode{\sphinxupquote{LpObjval}}
&
\sphinxAtStartPar
Double
&
\sphinxAtStartPar
The LP objective value
\\
\sphinxhline
\sphinxAtStartPar
\sphinxcode{\sphinxupquote{SimplexIter}}
&
\sphinxAtStartPar
Integer
&
\sphinxAtStartPar
Number of simplex iterations performed
\\
\sphinxhline
\sphinxAtStartPar
\sphinxcode{\sphinxupquote{BarrierIter}}
&
\sphinxAtStartPar
Integer
&
\sphinxAtStartPar
Number of barrier iterations performed
\\
\sphinxbottomrule
\end{tabular}
\sphinxtableafterendhook\par
\sphinxattableend\end{savenotes}

\sphinxAtStartPar
Attributes for results of a linear programming problem are shown in \hyperref[\detokenize{modeling:copttab-modelinglpsol-attr}]{Table \ref{\detokenize{modeling:copttab-modelinglpsol-attr}}}:

\begin{savenotes}\sphinxattablestart
\sphinxthistablewithglobalstyle
\centering
\sphinxcapstartof{table}
\sphinxthecaptionisattop
\sphinxcaption{Attributes for Linear Programming Results}\label{\detokenize{modeling:copttab-modelinglpsol-attr}}
\sphinxaftertopcaption
\begin{tabular}[t]{|\X{15}{70}|\X{10}{70}|\X{45}{70}|}
\sphinxtoprule
\sphinxstyletheadfamily 
\sphinxAtStartPar
Name
&\sphinxstyletheadfamily 
\sphinxAtStartPar
Type
&\sphinxstyletheadfamily 
\sphinxAtStartPar
Description
\\
\sphinxmidrule
\sphinxtableatstartofbodyhook
\sphinxAtStartPar
\sphinxcode{\sphinxupquote{LpStatus}}
&
\sphinxAtStartPar
Integer
&
\sphinxAtStartPar
The LP status
\\
\sphinxhline
\sphinxAtStartPar
\sphinxcode{\sphinxupquote{HasLpSol}}
&
\sphinxAtStartPar
Integer
&
\sphinxAtStartPar
Whether LP solution is available
\\
\sphinxhline
\sphinxAtStartPar
\sphinxcode{\sphinxupquote{HasBasis}}
&
\sphinxAtStartPar
Integer
&
\sphinxAtStartPar
Whether LP basis is available
\\
\sphinxhline
\sphinxAtStartPar
\sphinxcode{\sphinxupquote{LpObjval}}
&
\sphinxAtStartPar
Double
&
\sphinxAtStartPar
The LP objective value
\\
\sphinxbottomrule
\end{tabular}
\sphinxtableafterendhook\par
\sphinxattableend\end{savenotes}

\sphinxAtStartPar
For the results of a linear programming problem, COPT also provides relevant information constants,
as shown in \hyperref[\detokenize{modeling:copttab-modelinglpsol-info}]{Table \ref{\detokenize{modeling:copttab-modelinglpsol-info}}}:

\begin{savenotes}\sphinxattablestart
\sphinxthistablewithglobalstyle
\centering
\sphinxcapstartof{table}
\sphinxthecaptionisattop
\sphinxcaption{Information for Linear Programming}\label{\detokenize{modeling:copttab-modelinglpsol-info}}
\sphinxaftertopcaption
\begin{tabular}[t]{|\X{15}{70}|\X{10}{70}|\X{45}{70}|}
\sphinxtoprule
\sphinxstyletheadfamily 
\sphinxAtStartPar
Name
&\sphinxstyletheadfamily 
\sphinxAtStartPar
Type
&\sphinxstyletheadfamily 
\sphinxAtStartPar
Description
\\
\sphinxmidrule
\sphinxtableatstartofbodyhook
\sphinxAtStartPar
\sphinxcode{\sphinxupquote{Slack}}
&
\sphinxAtStartPar
Double
&
\sphinxAtStartPar
Solution of slack variables, also known as activities of constraints. Only available for LP problem
\\
\sphinxhline
\sphinxAtStartPar
\sphinxcode{\sphinxupquote{Dual}}
&
\sphinxAtStartPar
Double
&
\sphinxAtStartPar
Solution of dual variables. Only available for LP problem
\\
\sphinxhline
\sphinxAtStartPar
\sphinxcode{\sphinxupquote{RedCost}}
&
\sphinxAtStartPar
Double
&
\sphinxAtStartPar
Reduced cost of columns. Only available for LP problem
\\
\sphinxbottomrule
\end{tabular}
\sphinxtableafterendhook\par
\sphinxattableend\end{savenotes}

\sphinxAtStartPar
For different programming interfaces, please refer to {\hyperref[\detokenize{attribute:chapattrs}]{\sphinxcrossref{\DUrole{std,std-ref}{Attributes}}}} and {\hyperref[\detokenize{information:chapinfo}]{\sphinxcrossref{\DUrole{std,std-ref}{Information}}}} sections
to see how to access these attributes and information.

\section{Second\sphinxhyphen{}Order Cone Programming (SOCP)}
\label{\detokenize{modeling:second-order-cone-programming-socp}}\label{\detokenize{modeling:chapmodeling-socp}}
\sphinxAtStartPar
Second\sphinxhyphen{}Order Cone Programming (SOCP) is an optimization problem where the objective function is linear,
and the constraints include second\sphinxhyphen{}order cone.

\subsection{Mathematical Model}
\label{\detokenize{modeling:id1}}
\sphinxAtStartPar
Second\sphinxhyphen{}Order Cone (SOC):
\begin{equation}\label{equation:modeling:modeling:3}
\begin{split}\mathcal{Q}^{n+1} = \left\{(t,x)\in \mathbb{R}\times\mathbb{R}^n\ |\ t\geq \|x\|_2 \right\}\end{split}
\end{equation}
\sphinxAtStartPar
\sphinxstylestrong{Second\sphinxhyphen{}Order Cone Constraint:}

\sphinxAtStartPar
When \(t\in \mathbb{R}\text{ and } x\in \mathbb{R}^n\) are decision variables, a constraint of the form \(t\geq \|x\|_2\)
is called a second\sphinxhyphen{}order cone constraint.

\sphinxAtStartPar
The mathematical formulation is as follows:
\begin{equation}\label{equation:modeling:modeling:4}
\begin{split}\min\quad &c^Tx +c^f\\
\mathrm{s.t.}\quad &l^c\leq Ax\leq u^c \\
&l^v\leq x\leq u^v\\
&Fx+g\in \mathcal{Q}\end{split}
\end{equation}
\sphinxAtStartPar
The variables and arguments in the model have the following meanings:
\begin{itemize}
\item {} 
\sphinxAtStartPar
Decision variables: \(x=(x_j)_{j=0}^{n-1}\in\mathbb{R}^n\)

\item {} 
\sphinxAtStartPar
Decision variable bounds: \(l^v, u^v\in\mathbb{R}^n\), where \(l^v\) represents the lower bounds, and \(u^v\)
represents the upper bounds of the variables.

\item {} 
\sphinxAtStartPar
Constraint boundaries: \(l^c, u^c\in\mathbb{R}^m\), where \(l^c\) represents the lower bounds, and \(u^c\)
represents the upper bounds of the linear constraints.

\item {} 
\sphinxAtStartPar
Coefficient matrix of linear constraints: \(A=(a_{ij})_{m\times n}\in\mathbb{R}^{m\times n}\)

\item {} 
\sphinxAtStartPar
Coefficients in the objective function: \(c\in\mathbb{R}^n\) represents the coefficients of the variables in the objective function,
and \(c^f\) represents the constant term in the objective function.

\item {} 
\sphinxAtStartPar
Problem size: \(m\) represents the number of linear constraints, \(n\) represents the number of decision variables,
and \(k\) represents the number of second\sphinxhyphen{}order cone constraints.

\end{itemize}

\sphinxAtStartPar
The following are the meanings of arguments related to the cone constraints:
\begin{itemize}
\item {} 
\sphinxAtStartPar
\(F\): \(F\in\mathbb{R}^{k\times n}\) is the coefficient matrix of the cone constraints.

\item {} 
\sphinxAtStartPar
\(g\): \(g\in\mathbb{R}^k\) is the constant vector in the cone constraints.

\item {} 
\sphinxAtStartPar
\(\mathcal{Q}\): Represents the Cartesian product of \(k\) sets, \(\mathcal{Q} = \mathcal{Q}_1\times \mathcal{Q}_2\times \cdots \times \mathcal{Q}_p\), where \(p\) represents the number of second\sphinxhyphen{}order cone constraints, and each \(\mathcal{Q}_i,i\in\{1,2,\cdots,p\}\) represents a second\sphinxhyphen{}order cone.

\end{itemize}

\subsection{Modeling}
\label{\detokenize{modeling:id2}}
\sphinxAtStartPar
The basic steps to construct and solve a Second\sphinxhyphen{}Order Cone Programming (SOCP) model in COPT are as follows:
\begin{enumerate}
\sphinxsetlistlabels{\arabic}{enumi}{enumii}{}{.}%
\item {} 
\sphinxAtStartPar
Create the COPT environment and model.

\item {} 
\sphinxAtStartPar
Add required data.

\item {} 
\sphinxAtStartPar
Construct the SOCP model:
\begin{itemize}
\item {} 
\sphinxAtStartPar
Add decision variables.

\item {} 
\sphinxAtStartPar
Add constraints (second\sphinxhyphen{}order cone constraints, linear constraints).

\item {} 
\sphinxAtStartPar
Set the objective function.

\end{itemize}

\item {} 
\sphinxAtStartPar
Set solver parameters and solve.

\item {} 
\sphinxAtStartPar
Retrieve the solution results.

\end{enumerate}

\sphinxAtStartPar
\sphinxstylestrong{Modeling: Adding Second\sphinxhyphen{}Order Cone Constraints}

\sphinxAtStartPar
COPT supports modeling the following two types of second\sphinxhyphen{}order cone constraints:

\sphinxAtStartPar
\sphinxstylestrong{Standard Second\sphinxhyphen{}Order Cone}
\begin{equation}\label{equation:modeling:modeling:5}
\begin{split}Q^n= \left\{x\in \mathbb{R}^n\ \left|\ x_0\geq\sqrt{\sum_{i=1}^{n-1} x_i^2}, x_0\geq0 \right. \right\}\end{split}
\end{equation}
\sphinxAtStartPar
Constant representation: \sphinxcode{\sphinxupquote{CONE\_QUAD}}

\sphinxAtStartPar
\sphinxstylestrong{Rotated Second\sphinxhyphen{}Order Cone}
\begin{equation}\label{equation:modeling:modeling:6}
\begin{split}Q^n_r= \left\{x\in \mathbb{R}^n\ \left|\ 2x_0x_1\geq\sum_{i=2}^{n-1} x_i^2, x_0\geq0, x_1\geq 0 \right. \right\}\end{split}
\end{equation}
\sphinxAtStartPar
Constant representation: \sphinxcode{\sphinxupquote{CONE\_RQUAD}}

\sphinxAtStartPar
When adding second\sphinxhyphen{}order cone constraints to the model, the main function arguments that can be specified by the user are:
\begin{itemize}
\item {} 
\sphinxAtStartPar
\sphinxcode{\sphinxupquote{ctype}}: The type of second\sphinxhyphen{}order cone constraint.
\begin{itemize}
\item {} 
\sphinxAtStartPar
\sphinxcode{\sphinxupquote{CONE\_QUAD}} : Standard second\sphinxhyphen{}order cone.

\item {} 
\sphinxAtStartPar
\sphinxcode{\sphinxupquote{CONE\_RQUAD}}: Rotated second\sphinxhyphen{}order cone.

\end{itemize}

\item {} 
\sphinxAtStartPar
\sphinxcode{\sphinxupquote{cvars}}: Variables that form the second\sphinxhyphen{}order cone constraint.

\item {} 
\sphinxAtStartPar
\sphinxcode{\sphinxupquote{dim}}: Dimension of the second\sphinxhyphen{}order cone constraint.

\end{itemize}

\sphinxAtStartPar
The implementation in different programming interfaces is shown in the table below:

\begin{savenotes}\sphinxattablestart
\sphinxthistablewithglobalstyle
\centering
\sphinxcapstartof{table}
\sphinxthecaptionisattop
\sphinxcaption{Functions for Adding Second\sphinxhyphen{}Order Cone}\label{\detokenize{modeling:copttab-modelingconemodel-constr}}
\sphinxaftertopcaption
\begin{tabular}[t]{|\X{5}{15}|\X{10}{15}|}
\sphinxtoprule
\sphinxstyletheadfamily 
\sphinxAtStartPar
API
&\sphinxstyletheadfamily 
\sphinxAtStartPar
Function
\\
\sphinxmidrule
\sphinxtableatstartofbodyhook
\sphinxAtStartPar
C
&
\sphinxAtStartPar
\sphinxcode{\sphinxupquote{COPT\_AddCones}}
\\
\sphinxhline
\sphinxAtStartPar
C++
&
\sphinxAtStartPar
\sphinxcode{\sphinxupquote{Model::AddCone()}}
\\
\sphinxhline
\sphinxAtStartPar
C\#
&
\sphinxAtStartPar
\sphinxcode{\sphinxupquote{Model.addCone()}}
\\
\sphinxhline
\sphinxAtStartPar
Java
&
\sphinxAtStartPar
\sphinxcode{\sphinxupquote{Model.addCone()}}
\\
\sphinxhline
\sphinxAtStartPar
Python
&
\sphinxAtStartPar
\sphinxcode{\sphinxupquote{Model.addCone()}}
\\
\sphinxbottomrule
\end{tabular}
\sphinxtableafterendhook\par
\sphinxattableend\end{savenotes}

\begin{sphinxadmonition}{note}{Note}
\begin{itemize}
\item {} 
\sphinxAtStartPar
The function names, calling methods, and \sphinxstylestrong{argument names differ slightly} in different programming interfaces
for operations related to modeling second\sphinxhyphen{}order cone constraints, but the functionality and argument meanings are consistent.

\item {} 
\sphinxAtStartPar
When providing variables that form the second\sphinxhyphen{}order cone constraint, please input them sequentially according to
their order in the constraint expression.

\item {} 
\sphinxAtStartPar
In the C API, the members of the second\sphinxhyphen{}order cone constraints are provided in compressed row storage format.
For more information on sparse matrix compressed storage in the C API, please refer to the specific example in the relevant section.

\item {} 
\sphinxAtStartPar
In the Python API, an additional member method \sphinxcode{\sphinxupquote{Model.addConeByDim()}} is provided for specifying the dimension of the second\sphinxhyphen{}order cone.

\end{itemize}
\end{sphinxadmonition}

\sphinxAtStartPar
\sphinxstylestrong{Example of Second\sphinxhyphen{}Order Cone Constraints}

\sphinxAtStartPar
A standard second\sphinxhyphen{}order cone formed by \(x_4,x_1,x_2,x_3\):
\begin{equation}\label{equation:modeling:modeling:7}
\begin{split}x_4 \geq \sqrt{x_1^2 + x_2^2 + x_3^2}\end{split}
\end{equation}
\sphinxAtStartPar
A rotated second\sphinxhyphen{}order cone formed by \(x_3,x_4,x_1,x_2\):
\begin{equation}\label{equation:modeling:modeling:8}
\begin{split}2 x_3 x_4 \geq x_1^2 + x_2^2\end{split}
\end{equation}
\sphinxAtStartPar
Taking the standard second\sphinxhyphen{}order cone formed by \(x_4,x_1,x_2,x_3\) as an example, the code implementation in different programming interfaces is as follows:

\sphinxAtStartPar
\sphinxstylestrong{C Interface:}

\begin{sphinxVerbatim}[commandchars=\\\{\}]
\PYG{k+kt}{int}\PYG{+w}{ }\PYG{n}{ncone}\PYG{+w}{ }\PYG{o}{=}\PYG{+w}{ }\PYG{l+m+mi}{1}\PYG{p}{;}
\PYG{k+kt}{int}\PYG{+w}{ }\PYG{n}{conetype}\PYG{p}{[}\PYG{p}{]}\PYG{+w}{ }\PYG{o}{=}\PYG{+w}{ }\PYG{p}{\PYGZob{}}\PYG{n}{COPT\PYGZus{}CONE\PYGZus{}QUAD}\PYG{p}{\PYGZcb{}}\PYG{p}{;}
\PYG{k+kt}{int}\PYG{+w}{ }\PYG{n}{conebeg}\PYG{p}{[}\PYG{p}{]}\PYG{+w}{ }\PYG{o}{=}\PYG{+w}{ }\PYG{p}{\PYGZob{}}\PYG{l+m+mi}{0}\PYG{p}{\PYGZcb{}}\PYG{p}{;}
\PYG{k+kt}{int}\PYG{+w}{ }\PYG{n}{conecnt}\PYG{p}{[}\PYG{p}{]}\PYG{+w}{ }\PYG{o}{=}\PYG{+w}{ }\PYG{p}{\PYGZob{}}\PYG{l+m+mi}{4}\PYG{p}{\PYGZcb{}}\PYG{p}{;}
\PYG{k+kt}{int}\PYG{+w}{ }\PYG{n}{coneind}\PYG{p}{[}\PYG{p}{]}\PYG{+w}{ }\PYG{o}{=}\PYG{+w}{ }\PYG{p}{\PYGZob{}}\PYG{l+m+mi}{3}\PYG{p}{,}\PYG{+w}{ }\PYG{l+m+mi}{0}\PYG{p}{,}\PYG{+w}{ }\PYG{l+m+mi}{1}\PYG{p}{,}\PYG{+w}{ }\PYG{l+m+mi}{2}\PYG{p}{\PYGZcb{}}\PYG{p}{;}
\PYG{n}{COPT\PYGZus{}AddCones}\PYG{p}{(}\PYG{n}{prob}\PYG{p}{,}\PYG{+w}{ }\PYG{n}{ncone}\PYG{p}{,}\PYG{+w}{ }\PYG{n}{conetype}\PYG{p}{,}\PYG{+w}{ }\PYG{n}{conebeg}\PYG{p}{,}\PYG{+w}{ }\PYG{n}{conecnt}\PYG{p}{,}\PYG{+w}{ }\PYG{n}{coneind}\PYG{p}{)}
\end{sphinxVerbatim}

\sphinxAtStartPar
\sphinxstylestrong{C++ Interface:}

\begin{sphinxVerbatim}[commandchars=\\\{\}]
\PYG{n}{VarArray}\PYG{+w}{ }\PYG{n}{cvars}\PYG{p}{;}
\PYG{n}{cvars}\PYG{p}{.}\PYG{n}{PushBack}\PYG{p}{(}\PYG{n}{x4}\PYG{p}{)}\PYG{p}{;}
\PYG{n}{cvars}\PYG{p}{.}\PYG{n}{PushBack}\PYG{p}{(}\PYG{n}{x1}\PYG{p}{)}\PYG{p}{;}
\PYG{n}{cvars}\PYG{p}{.}\PYG{n}{PushBack}\PYG{p}{(}\PYG{n}{x2}\PYG{p}{)}\PYG{p}{;}
\PYG{n}{cvars}\PYG{p}{.}\PYG{n}{PushBack}\PYG{p}{(}\PYG{n}{x3}\PYG{p}{)}\PYG{p}{;}
\PYG{n}{model}\PYG{p}{.}\PYG{n}{AddCone}\PYG{p}{(}\PYG{n}{cvars}\PYG{p}{,}\PYG{+w}{ }\PYG{n}{COPT\PYGZus{}CONE\PYGZus{}QUAD}\PYG{p}{)}\PYG{p}{;}
\end{sphinxVerbatim}

\sphinxAtStartPar
\sphinxstylestrong{C\# Interface:}

\begin{sphinxVerbatim}[commandchars=\\\{\}]
\PYG{n}{VarArray}\PYG{+w}{ }\PYG{n}{cvars}\PYG{+w}{ }\PYG{o}{=}\PYG{+w}{ }\PYG{k}{new}\PYG{+w}{ }\PYG{n}{VarArray}\PYG{p}{(}\PYG{p}{)}\PYG{p}{;}
\PYG{n}{cvars}\PYG{p}{.}\PYG{n}{PushBack}\PYG{p}{(}\PYG{n}{x4}\PYG{p}{)}\PYG{p}{;}
\PYG{n}{cvars}\PYG{p}{.}\PYG{n}{PushBack}\PYG{p}{(}\PYG{n}{x1}\PYG{p}{)}\PYG{p}{;}
\PYG{n}{cvars}\PYG{p}{.}\PYG{n}{PushBack}\PYG{p}{(}\PYG{n}{x2}\PYG{p}{)}\PYG{p}{;}
\PYG{n}{cvars}\PYG{p}{.}\PYG{n}{PushBack}\PYG{p}{(}\PYG{n}{x3}\PYG{p}{)}\PYG{p}{;}
\PYG{n}{model}\PYG{p}{.}\PYG{n}{AddCone}\PYG{p}{(}\PYG{n}{cvars}\PYG{p}{,}\PYG{+w}{ }\PYG{n}{Copt}\PYG{p}{.}\PYG{n}{Consts}\PYG{p}{.}\PYG{n}{CONE\PYGZus{}QUAD}\PYG{p}{)}\PYG{p}{;}
\end{sphinxVerbatim}

\sphinxAtStartPar
\sphinxstylestrong{Java Interface:}

\begin{sphinxVerbatim}[commandchars=\\\{\}]
\PYG{n}{VarArray}\PYG{+w}{ }\PYG{n}{cvars}\PYG{+w}{ }\PYG{o}{=}\PYG{+w}{ }\PYG{k}{new}\PYG{+w}{ }\PYG{n}{VarArray}\PYG{p}{(}\PYG{p}{)}\PYG{p}{;}
\PYG{n}{cvars}\PYG{p}{.}\PYG{n+na}{PushBack}\PYG{p}{(}\PYG{n}{x4}\PYG{p}{)}\PYG{p}{;}
\PYG{n}{cvars}\PYG{p}{.}\PYG{n+na}{PushBack}\PYG{p}{(}\PYG{n}{x1}\PYG{p}{)}\PYG{p}{;}
\PYG{n}{cvars}\PYG{p}{.}\PYG{n+na}{PushBack}\PYG{p}{(}\PYG{n}{x2}\PYG{p}{)}\PYG{p}{;}
\PYG{n}{cvars}\PYG{p}{.}\PYG{n+na}{PushBack}\PYG{p}{(}\PYG{n}{x3}\PYG{p}{)}\PYG{p}{;}
\PYG{n}{model}\PYG{p}{.}\PYG{n+na}{addCone}\PYG{p}{(}\PYG{n}{cvars}\PYG{p}{,}\PYG{+w}{ }\PYG{n}{copt}\PYG{p}{.}\PYG{n+na}{Consts}\PYG{p}{.}\PYG{n+na}{CONE\PYGZus{}QUAD}\PYG{p}{)}\PYG{p}{;}
\end{sphinxVerbatim}

\sphinxAtStartPar
\sphinxstylestrong{Python Interface:}

\begin{sphinxVerbatim}[commandchars=\\\{\}]
\PYG{n}{model}\PYG{o}{.}\PYG{n}{addCone}\PYG{p}{(}\PYG{p}{[}\PYG{n}{x4}\PYG{p}{,} \PYG{n}{x1}\PYG{p}{,} \PYG{n}{x2}\PYG{p}{,} \PYG{n}{x3}\PYG{p}{]}\PYG{p}{,} \PYG{n}{COPT}\PYG{o}{.}\PYG{n}{CONE\PYGZus{}QUAD}\PYG{p}{)}
\end{sphinxVerbatim}

\sphinxAtStartPar
In the programming interfaces provided by COPT, with the exception of the C language,
object\sphinxhyphen{}oriented programming interfaces (C\#, C++, Java, Python) provide classes related to second\sphinxhyphen{}order cone constraints:
\begin{itemize}
\item {} 
\sphinxAtStartPar
\sphinxstylestrong{Encapsulation of operations related to second\sphinxhyphen{}order cone constraints:}
\begin{enumerate}
\sphinxsetlistlabels{\arabic}{enumi}{enumii}{}{.}%
\item {} 
\sphinxAtStartPar
\sphinxcode{\sphinxupquote{Cone}} class: Encapsulation of operations related to second\sphinxhyphen{}order cone constraints in COPT.

\item {} 
\sphinxAtStartPar
\sphinxcode{\sphinxupquote{ConeArray}} class: Conveniently allows users to operate on a group of \sphinxcode{\sphinxupquote{Cone}} objects.

\end{enumerate}

\item {} 
\sphinxAtStartPar
\sphinxstylestrong{Encapsulation of second\sphinxhyphen{}order cone constraint builders:}
\begin{enumerate}
\sphinxsetlistlabels{\arabic}{enumi}{enumii}{}{.}%
\item {} 
\sphinxAtStartPar
\sphinxcode{\sphinxupquote{ConeBuilder}} class: Encapsulation of builders for constructing second\sphinxhyphen{}order cone constraints in COPT.

\item {} 
\sphinxAtStartPar
\sphinxcode{\sphinxupquote{ConeBuilderArray}} class: Conveniently allows users to operate on a group of \sphinxcode{\sphinxupquote{ConeBuilder}} objects.

\end{enumerate}

\item {} 
\sphinxAtStartPar
C++ API: {\hyperref[\detokenize{cppapiref:chapcppapiref-cone}]{\sphinxcrossref{\DUrole{std,std-ref}{Cone Class}}}} , {\hyperref[\detokenize{cppapiref:chapcppapiref-conearray}]{\sphinxcrossref{\DUrole{std,std-ref}{ConeArray Class}}}} , {\hyperref[\detokenize{cppapiref:chapcppapiref-conebuilder}]{\sphinxcrossref{\DUrole{std,std-ref}{ConeBuilder Class}}}} , {\hyperref[\detokenize{cppapiref:chapcppapiref-conebuilderarray}]{\sphinxcrossref{\DUrole{std,std-ref}{ConeBuilderArray Class}}}}

\item {} 
\sphinxAtStartPar
C\# API: {\hyperref[\detokenize{csharpapiref:chapcsharpapiref-cone}]{\sphinxcrossref{\DUrole{std,std-ref}{Cone Class}}}} , {\hyperref[\detokenize{csharpapiref:chapcsharpapiref-conearray}]{\sphinxcrossref{\DUrole{std,std-ref}{ConeArray Class}}}} , {\hyperref[\detokenize{csharpapiref:chapcsharpapiref-conebuilder}]{\sphinxcrossref{\DUrole{std,std-ref}{ConeBuilder Class}}}} , {\hyperref[\detokenize{csharpapiref:chapcsharpapiref-conebuilderarray}]{\sphinxcrossref{\DUrole{std,std-ref}{ConeBuilderArray Class}}}}

\item {} 
\sphinxAtStartPar
Java API: {\hyperref[\detokenize{javaapiref:chapjavaapiref-cone}]{\sphinxcrossref{\DUrole{std,std-ref}{Cone Class}}}} , {\hyperref[\detokenize{javaapiref:chapjavaapiref-conearray}]{\sphinxcrossref{\DUrole{std,std-ref}{ConeArray Class}}}} , {\hyperref[\detokenize{javaapiref:chapjavaapiref-conebuilder}]{\sphinxcrossref{\DUrole{std,std-ref}{ConeBuilder Class}}}} , {\hyperref[\detokenize{javaapiref:chapjavaapiref-conebuilderarray}]{\sphinxcrossref{\DUrole{std,std-ref}{ConeBuilderArray Class}}}}

\item {} 
\sphinxAtStartPar
Python API: {\hyperref[\detokenize{pyapiref:chappyapi-cone}]{\sphinxcrossref{\DUrole{std,std-ref}{Cone Class}}}} , {\hyperref[\detokenize{pyapiref:chappyapi-conearray}]{\sphinxcrossref{\DUrole{std,std-ref}{ConeArray Class}}}} , {\hyperref[\detokenize{pyapiref:chappyapi-conebuilder}]{\sphinxcrossref{\DUrole{std,std-ref}{ConeBuilder Class}}}} , {\hyperref[\detokenize{pyapiref:chappyapi-conebuilderarray}]{\sphinxcrossref{\DUrole{std,std-ref}{ConeBuilderArray Class}}}}

\end{itemize}

\sphinxAtStartPar
\sphinxstylestrong{Attributes Related to Second\sphinxhyphen{}Order Cone Constraints}
\begin{itemize}
\item {} 
\sphinxAtStartPar
\sphinxcode{\sphinxupquote{COPT\_INTATTR\_CONES}} or \sphinxcode{\sphinxupquote{"Cones"}}
\begin{quote}

\sphinxAtStartPar
Integer attribute.
The number of second\sphinxhyphen{}order cone constraints in the model.
\end{quote}

\end{itemize}

\sphinxAtStartPar
This attribute provides the count of second\sphinxhyphen{}order cone constraints that have been added to the model.
It can be useful for monitoring or validating the structure of the model during or after the modeling process.

\section{ExpCone Programming}
\label{\detokenize{modeling:expcone-programming}}\label{\detokenize{modeling:chapmodeling-expcone}}

\subsection{Mathematical Formulation}
\label{\detokenize{modeling:mathematical-formulation}}
\sphinxAtStartPar
COPT supports two types of exponential cone contraints:
\begin{itemize}
\item {} 
\sphinxAtStartPar
\sphinxcode{\sphinxupquote{EXPCONE\_PRIMAL}} : Primal exponential cone

\end{itemize}
\begin{equation}\label{equation:modeling:modeling:9}
\begin{split}\mathrm{cl}(S_1) = S_1 \cup S_2\end{split}
\end{equation}\begin{equation}\label{equation:modeling:modeling:10}
\begin{split}\begin{aligned}
S_1 &= \left\{\begin{pmatrix} t \\ s \\ r \end{pmatrix}\in \mathbb{R}^3\ |\ s > 0,\ t \geq s\ \mathrm{exp}\left(\frac{r}{s} \right) \right\}, \\
S_2 &= \left\{\begin{pmatrix} t \\ s \\ r \end{pmatrix}\in \mathbb{R}^3\ |\ s=0,\ t\geq 0,\ r\leq 0 \right\}
\end{aligned}\end{split}
\end{equation}\begin{itemize}
\item {} 
\sphinxAtStartPar
\sphinxcode{\sphinxupquote{EXPCONE\_DUAL}} : Dual exponential cone

\end{itemize}
\begin{equation}\label{equation:modeling:modeling:11}
\begin{split}\mathrm{cl}(S_1) = S_1 \cup S_2\end{split}
\end{equation}\begin{equation}\label{equation:modeling:modeling:12}
\begin{split}\begin{aligned}
S_1 &= \left\{\begin{pmatrix} t \\ s \\ r \end{pmatrix}\in \mathbb{R}^3\ |\ r < 0,\ t \geq -r\ \mathrm{exp}\left(\frac{s}{r}-1\right) \right\}, \\
S_2 &= \left\{\begin{pmatrix} t \\ s \\ r \end{pmatrix}\in \mathbb{R}^3\ |\ r = 0,\ t\geq 0,\ s\geq 0 \right\}
\end{aligned}\end{split}
\end{equation}

\subsection{Exponential Cone Example}
\label{\detokenize{modeling:exponential-cone-example}}
\sphinxAtStartPar
Take \((u_1, 1, u_3) \in K_{\text{exp}}\) for example, which mathematical form is \(u_1 \geq e^{u_3}\).
The code implementation is as follows:

\sphinxAtStartPar
\sphinxstylestrong{C Interface:}

\sphinxAtStartPar
Please refer to the example code in the installation package \sphinxcode{\sphinxupquote{expcone\_gp.c}}.

\sphinxAtStartPar
\sphinxstylestrong{Python Interface:}

\begin{sphinxVerbatim}[commandchars=\\\{\}]
\PYG{n}{u1} \PYG{o}{=} \PYG{n}{model}\PYG{o}{.}\PYG{n}{addVar}\PYG{p}{(}\PYG{n}{lb}\PYG{o}{=}\PYG{o}{\PYGZhy{}}\PYG{n}{COPT}\PYG{o}{.}\PYG{n}{INFINITY}\PYG{p}{)}
\PYG{n}{u2} \PYG{o}{=} \PYG{n}{model}\PYG{o}{.}\PYG{n}{addVar}\PYG{p}{(}\PYG{n}{lb}\PYG{o}{=}\PYG{l+m+mf}{1.0}\PYG{p}{,} \PYG{n}{ub}\PYG{o}{=}\PYG{l+m+mf}{1.0}\PYG{p}{)}
\PYG{n}{u3} \PYG{o}{=} \PYG{n}{model}\PYG{o}{.}\PYG{n}{addVar}\PYG{p}{(}\PYG{n}{lb}\PYG{o}{=}\PYG{o}{\PYGZhy{}}\PYG{n}{COPT}\PYG{o}{.}\PYG{n}{INFINITY}\PYG{p}{)}
\PYG{n}{model}\PYG{o}{.}\PYG{n}{addExpCone}\PYG{p}{(}\PYG{p}{[}\PYG{n}{u1}\PYG{p}{,} \PYG{n}{u2}\PYG{p}{,} \PYG{n}{u3}\PYG{p}{]}\PYG{p}{,} \PYG{n}{COPT}\PYG{o}{.}\PYG{n}{EXPCONE\PYGZus{}PRIMAL}\PYG{p}{)}
\end{sphinxVerbatim}

\sphinxAtStartPar
\sphinxstylestrong{C++ Interface:}

\begin{sphinxVerbatim}[commandchars=\\\{\}]
\PYG{n}{Var}\PYG{+w}{ }\PYG{n}{u1}\PYG{+w}{ }\PYG{o}{=}\PYG{+w}{ }\PYG{n}{model}\PYG{p}{.}\PYG{n}{AddVar}\PYG{p}{(}\PYG{o}{\PYGZhy{}}\PYG{n}{COPT\PYGZus{}INFINITY}\PYG{p}{,}\PYG{+w}{ }\PYG{o}{+}\PYG{n}{COPT\PYGZus{}INFINITY}\PYG{p}{,}\PYG{+w}{ }\PYG{l+m+mf}{0.0}\PYG{p}{,}\PYG{+w}{ }\PYG{n}{COPT\PYGZus{}CONTINUOUS}\PYG{p}{,}\PYG{+w}{ }\PYG{l+s}{\PYGZdq{}}\PYG{l+s}{u1}\PYG{l+s}{\PYGZdq{}}\PYG{p}{)}\PYG{p}{;}
\PYG{n}{Var}\PYG{+w}{ }\PYG{n}{u2}\PYG{+w}{ }\PYG{o}{=}\PYG{+w}{ }\PYG{n}{model}\PYG{p}{.}\PYG{n}{AddVar}\PYG{p}{(}\PYG{l+m+mf}{1.0}\PYG{p}{,}\PYG{+w}{ }\PYG{l+m+mf}{1.0}\PYG{p}{,}\PYG{+w}{ }\PYG{l+m+mf}{0.0}\PYG{p}{,}\PYG{+w}{ }\PYG{n}{COPT\PYGZus{}CONTINUOUS}\PYG{p}{,}\PYG{+w}{ }\PYG{l+s}{\PYGZdq{}}\PYG{l+s}{u2}\PYG{l+s}{\PYGZdq{}}\PYG{p}{)}\PYG{p}{;}
\PYG{n}{Var}\PYG{+w}{ }\PYG{n}{u3}\PYG{+w}{ }\PYG{o}{=}\PYG{+w}{ }\PYG{n}{model}\PYG{p}{.}\PYG{n}{AddVar}\PYG{p}{(}\PYG{o}{\PYGZhy{}}\PYG{n}{COPT\PYGZus{}INFINITY}\PYG{p}{,}\PYG{+w}{ }\PYG{o}{+}\PYG{n}{COPT\PYGZus{}INFINITY}\PYG{p}{,}\PYG{+w}{ }\PYG{l+m+mf}{0.0}\PYG{p}{,}\PYG{+w}{ }\PYG{n}{COPT\PYGZus{}CONTINUOUS}\PYG{p}{,}\PYG{+w}{ }\PYG{l+s}{\PYGZdq{}}\PYG{l+s}{u3}\PYG{l+s}{\PYGZdq{}}\PYG{p}{)}\PYG{p}{;}

\PYG{n}{VarArray}\PYG{+w}{ }\PYG{n}{uconevars}\PYG{p}{;}
\PYG{n}{uconevars}\PYG{p}{.}\PYG{n}{PushBack}\PYG{p}{(}\PYG{n}{u1}\PYG{p}{)}\PYG{p}{;}
\PYG{n}{uconevars}\PYG{p}{.}\PYG{n}{PushBack}\PYG{p}{(}\PYG{n}{u2}\PYG{p}{)}\PYG{p}{;}
\PYG{n}{uconevars}\PYG{p}{.}\PYG{n}{PushBack}\PYG{p}{(}\PYG{n}{u3}\PYG{p}{)}\PYG{p}{;}
\PYG{n}{model}\PYG{p}{.}\PYG{n}{AddExpCone}\PYG{p}{(}\PYG{n}{uconevars}\PYG{p}{,}\PYG{+w}{ }\PYG{n}{COPT\PYGZus{}EXPCONE\PYGZus{}PRIMAL}\PYG{p}{)}\PYG{p}{;}
\end{sphinxVerbatim}

\sphinxAtStartPar
\sphinxstylestrong{C\# Interface:}

\begin{sphinxVerbatim}[commandchars=\\\{\}]
\PYG{n}{VarArray}\PYG{+w}{ }\PYG{n}{uconevars}\PYG{+w}{ }\PYG{o}{=}\PYG{+w}{ }\PYG{k}{new}\PYG{+w}{ }\PYG{n}{VarArray}\PYG{p}{(}\PYG{p}{)}\PYG{p}{;}
\PYG{n}{uconevars}\PYG{p}{.}\PYG{n}{PushBack}\PYG{p}{(}\PYG{n}{u1}\PYG{p}{)}\PYG{p}{;}
\PYG{n}{uconevars}\PYG{p}{.}\PYG{n}{PushBack}\PYG{p}{(}\PYG{n}{u2}\PYG{p}{)}\PYG{p}{;}
\PYG{n}{uconevars}\PYG{p}{.}\PYG{n}{PushBack}\PYG{p}{(}\PYG{n}{u3}\PYG{p}{)}\PYG{p}{;}
\PYG{n}{model}\PYG{p}{.}\PYG{n}{AddExpCone}\PYG{p}{(}\PYG{n}{uconevars}\PYG{p}{,}\PYG{+w}{ }\PYG{n}{Copt}\PYG{p}{.}\PYG{n}{Consts}\PYG{p}{.}\PYG{n}{EXPCONE\PYGZus{}PRIMAL}\PYG{p}{)}\PYG{p}{;}

\PYG{n}{Var}\PYG{+w}{ }\PYG{n}{u1}\PYG{+w}{ }\PYG{o}{=}\PYG{+w}{ }\PYG{n}{model}\PYG{p}{.}\PYG{n}{AddVar}\PYG{p}{(}\PYG{o}{\PYGZhy{}}\PYG{n}{Copt}\PYG{p}{.}\PYG{n}{Consts}\PYG{p}{.}\PYG{n}{INFINITY}\PYG{p}{,}\PYG{+w}{ }\PYG{o}{+}\PYG{n}{Copt}\PYG{p}{.}\PYG{n}{Consts}\PYG{p}{.}\PYG{n}{INFINITY}\PYG{p}{,}\PYG{+w}{ }\PYG{l+m+mf}{0.0}\PYG{p}{,}\PYG{+w}{ }\PYG{n}{Copt}\PYG{p}{.}\PYG{n}{Consts}\PYG{p}{.}\PYG{n}{CONTINUOUS}\PYG{p}{,}\PYG{+w}{ }\PYG{l+s}{\PYGZdq{}u1\PYGZdq{}}\PYG{p}{)}\PYG{p}{;}
\PYG{n}{Var}\PYG{+w}{ }\PYG{n}{u2}\PYG{+w}{ }\PYG{o}{=}\PYG{+w}{ }\PYG{n}{model}\PYG{p}{.}\PYG{n}{AddVar}\PYG{p}{(}\PYG{l+m+mf}{1.0}\PYG{p}{,}\PYG{+w}{ }\PYG{l+m+mf}{1.0}\PYG{p}{,}\PYG{+w}{ }\PYG{l+m+mf}{0.0}\PYG{p}{,}\PYG{+w}{ }\PYG{n}{Copt}\PYG{p}{.}\PYG{n}{Consts}\PYG{p}{.}\PYG{n}{CONTINUOUS}\PYG{p}{,}\PYG{+w}{ }\PYG{l+s}{\PYGZdq{}u2\PYGZdq{}}\PYG{p}{)}\PYG{p}{;}
\PYG{n}{Var}\PYG{+w}{ }\PYG{n}{u3}\PYG{+w}{ }\PYG{o}{=}\PYG{+w}{ }\PYG{n}{model}\PYG{p}{.}\PYG{n}{AddVar}\PYG{p}{(}\PYG{o}{\PYGZhy{}}\PYG{n}{Copt}\PYG{p}{.}\PYG{n}{Consts}\PYG{p}{.}\PYG{n}{INFINITY}\PYG{p}{,}\PYG{+w}{ }\PYG{o}{+}\PYG{n}{Copt}\PYG{p}{.}\PYG{n}{Consts}\PYG{p}{.}\PYG{n}{INFINITY}\PYG{p}{,}\PYG{+w}{ }\PYG{l+m+mf}{0.0}\PYG{p}{,}\PYG{+w}{ }\PYG{n}{Copt}\PYG{p}{.}\PYG{n}{Consts}\PYG{p}{.}\PYG{n}{CONTINUOUS}\PYG{p}{,}\PYG{+w}{ }\PYG{l+s}{\PYGZdq{}u3\PYGZdq{}}\PYG{p}{)}\PYG{p}{;}
\end{sphinxVerbatim}

\sphinxAtStartPar
\sphinxstylestrong{Java Interface:}

\begin{sphinxVerbatim}[commandchars=\\\{\}]
\PYG{n}{Var}\PYG{+w}{ }\PYG{n}{u1}\PYG{+w}{ }\PYG{o}{=}\PYG{+w}{ }\PYG{n}{model}\PYG{p}{.}\PYG{n+na}{addVar}\PYG{p}{(}\PYG{o}{\PYGZhy{}}\PYG{n}{copt}\PYG{p}{.}\PYG{n+na}{Consts}\PYG{p}{.}\PYG{n+na}{INFINITY}\PYG{p}{,}\PYG{+w}{ }\PYG{o}{+}\PYG{n}{copt}\PYG{p}{.}\PYG{n+na}{Consts}\PYG{p}{.}\PYG{n+na}{INFINITY}\PYG{p}{,}\PYG{+w}{ }\PYG{l+m+mf}{0.0}\PYG{p}{,}\PYG{+w}{ }\PYG{n}{copt}\PYG{p}{.}\PYG{n+na}{Consts}\PYG{p}{.}\PYG{n+na}{CONTINUOUS}\PYG{p}{,}\PYG{+w}{ }\PYG{l+s}{\PYGZdq{}}\PYG{l+s}{u1}\PYG{l+s}{\PYGZdq{}}\PYG{p}{)}\PYG{p}{;}
\PYG{n}{Var}\PYG{+w}{ }\PYG{n}{u2}\PYG{+w}{ }\PYG{o}{=}\PYG{+w}{ }\PYG{n}{model}\PYG{p}{.}\PYG{n+na}{addVar}\PYG{p}{(}\PYG{l+m+mf}{1.0}\PYG{p}{,}\PYG{+w}{ }\PYG{l+m+mf}{1.0}\PYG{p}{,}\PYG{+w}{ }\PYG{l+m+mf}{0.0}\PYG{p}{,}\PYG{+w}{ }\PYG{n}{copt}\PYG{p}{.}\PYG{n+na}{Consts}\PYG{p}{.}\PYG{n+na}{CONTINUOUS}\PYG{p}{,}\PYG{+w}{ }\PYG{l+s}{\PYGZdq{}}\PYG{l+s}{u2}\PYG{l+s}{\PYGZdq{}}\PYG{p}{)}\PYG{p}{;}
\PYG{n}{Var}\PYG{+w}{ }\PYG{n}{u3}\PYG{+w}{ }\PYG{o}{=}\PYG{+w}{ }\PYG{n}{model}\PYG{p}{.}\PYG{n+na}{addVar}\PYG{p}{(}\PYG{o}{\PYGZhy{}}\PYG{n}{copt}\PYG{p}{.}\PYG{n+na}{Consts}\PYG{p}{.}\PYG{n+na}{INFINITY}\PYG{p}{,}\PYG{+w}{ }\PYG{o}{+}\PYG{n}{copt}\PYG{p}{.}\PYG{n+na}{Consts}\PYG{p}{.}\PYG{n+na}{INFINITY}\PYG{p}{,}\PYG{+w}{ }\PYG{l+m+mf}{0.0}\PYG{p}{,}\PYG{+w}{ }\PYG{n}{copt}\PYG{p}{.}\PYG{n+na}{Consts}\PYG{p}{.}\PYG{n+na}{CONTINUOUS}\PYG{p}{,}\PYG{+w}{ }\PYG{l+s}{\PYGZdq{}}\PYG{l+s}{u3}\PYG{l+s}{\PYGZdq{}}\PYG{p}{)}\PYG{p}{;}

\PYG{n}{VarArray}\PYG{+w}{ }\PYG{n}{uconevars}\PYG{+w}{ }\PYG{o}{=}\PYG{+w}{ }\PYG{k}{new}\PYG{+w}{ }\PYG{n}{VarArray}\PYG{p}{(}\PYG{p}{)}\PYG{p}{;}
\PYG{n}{uconevars}\PYG{p}{.}\PYG{n+na}{pushBack}\PYG{p}{(}\PYG{n}{u1}\PYG{p}{)}\PYG{p}{;}
\PYG{n}{uconevars}\PYG{p}{.}\PYG{n+na}{pushBack}\PYG{p}{(}\PYG{n}{u2}\PYG{p}{)}\PYG{p}{;}
\PYG{n}{uconevars}\PYG{p}{.}\PYG{n+na}{pushBack}\PYG{p}{(}\PYG{n}{u3}\PYG{p}{)}\PYG{p}{;}
\PYG{n}{model}\PYG{p}{.}\PYG{n+na}{addExpCone}\PYG{p}{(}\PYG{n}{uconevars}\PYG{p}{,}\PYG{+w}{ }\PYG{n}{copt}\PYG{p}{.}\PYG{n+na}{Consts}\PYG{p}{.}\PYG{n+na}{EXPCONE\PYGZus{}PRIMAL}\PYG{p}{)}\PYG{p}{;}
\end{sphinxVerbatim}

\section{Semidefinite Programming (SDP)}
\label{\detokenize{modeling:semidefinite-programming-sdp}}\label{\detokenize{modeling:chapmodeling-sdp}}
\sphinxAtStartPar
Semidefinite Programming (SDP) consists semidefinite variables and cone constraints.

\subsection{Mathematical Model}
\label{\detokenize{modeling:id3}}
\sphinxAtStartPar
\sphinxstylestrong{Positive Semidefinite Cone:}

\sphinxAtStartPar
Let \(\mathcal{S}^n\) denote the set of \(n\) \sphinxhyphen{}dimensional symmetric matrices, the positive semidefinite cone is defined as:
\begin{equation}\label{equation:modeling:modeling:13}
\begin{split}\mathcal{S}_+^n = \left\{X\in \mathcal{S}^n\ |\ u^TXu \geq 0, \forall u\in \mathbb{R}^n\right\}\end{split}
\end{equation}
\sphinxAtStartPar
\sphinxstylestrong{Positive Semidefinite (PSD) Variable:}

\sphinxAtStartPar
In an SDP model, the decision variable \(X\in \mathcal{S}_+^n\), can also be denoted as \(X\succeq 0\).
This decision variable is referred to as a semidefinite variable. Besides, the linear constraints of the model include semidefinite variables.

\sphinxAtStartPar
\sphinxstylestrong{Semidefinite Cone Constraint}

\sphinxAtStartPar
In an SDP model, \(X\succeq 0\) is commonly referred to as the \sphinxstylestrong{semidefinite cone constraint},
while a constraint of the form \(A\bullet X = b\), which involves a semidefinite variable,
is referred to as a linear constraint consisting of PSD variables. For simplicity, in the following chapter,
we will refer to it as the \sphinxstylestrong{semidefinite constraint} and distinguish it from the semidefinite cone constraint (i.e., \(X\succeq 0\)).

\sphinxAtStartPar
The mathematical form of semidefinite programming is as follows:
\begin{equation}\label{equation:modeling:modeling:14}
\begin{split}\min\quad &\sum_{j=0}^{n-1} c_jx_j + \sum_{j=0}^{p-1} C_j \bullet X_j + c^f \\
\text{s.t.}\quad &l_i^c \leq \sum_{j=0}^{n-1} a_{ij} x_j + \sum_{j=0}^{p-1} A_j \bullet X_j \leq u_i^c,\qquad i=0,1, ..., m-1 \\
&l_j^v\leq\qquad\qquad\quad x_j\qquad\qquad \leq u_j^v,\qquad j=0,1, ..., n-1\\
&\qquad\qquad\qquad X_j \succeq 0,\qquad j = 0,1, ..., p-1\end{split}
\end{equation}
\sphinxAtStartPar
Here, the operator \(\bullet\) denotes the matrix inner product operation:

\sphinxAtStartPar
Given two matrices \(A = \{a_{ij}\} \in \mathbb{R}^{m\times n}\) and \(B = \{b_{ij}\}\in \mathbb{R}^{m\times n}\),
the inner product of matrix \(A\) and matrix \(B\) is defined as:
\begin{equation}\label{equation:modeling:modeling:15}
\begin{split}A\bullet B = \sum_{i=0}^{m-1} \sum_{j=0}^{n-1} a_{ij}b_{ij}\end{split}
\end{equation}
\sphinxAtStartPar
The variables and arguments in the model have the following meanings:
\begin{itemize}
\item {} 
\sphinxAtStartPar
Problem size: \(m\) denotes the number of constraints, \(n\) denotes the number of non\sphinxhyphen{}semidefinite variables,
and \(p\) denotes the number of semidefinite variables.

\item {} 
\sphinxAtStartPar
Decision variables: non\sphinxhyphen{}semidefinite variables \(x=(x_j)_{j=0}^{n-1}\in\mathbb{R}^n\),
semidefinite variables \(X_j\in \mathcal{S}_+^{r_j}\) (for \(j=0, ..., p-1\)).

\item {} 
\sphinxAtStartPar
Non\sphinxhyphen{}semidefinite variable range: \(l^v, u^v\in\mathbb{R}^n\), where \(l^v\) denotes the lower bounds and \(u^v\)
denotes the upper bounds of the non\sphinxhyphen{}semidefinite variables.

\item {} 
\sphinxAtStartPar
Constraint boundaries: \(l^c, u^c\in\mathbb{R}^m\), where \(l^c\) denotes the lower bounds and \(u^c\) denotes
the upper bounds of the constraints.

\item {} 
\sphinxAtStartPar
Coefficient matrix of the linear constraints: \(A=(a_{ij})_{m\times n}\in\mathbb{R}^{m\times n}\),
\(A_j\in \mathbb{R}^{r_j\times r_j}\).

\item {} 
\sphinxAtStartPar
Coefficients in the objective function: \(c\in\mathbb{R}^n\) represents the coefficients of the non\sphinxhyphen{}semidefinite variables,
\(C_j\in\mathbb{R}^{r_j\times r_j}\) represents the coefficients of the semidefinite variables, and \(c^f\) represents the constant term in the objective function.

\end{itemize}

\subsection{Modeling}
\label{\detokenize{modeling:id4}}
\sphinxAtStartPar
The basic steps to construct and solve a semidefinite programming (SDP) model in COPT are as follows:
\begin{enumerate}
\sphinxsetlistlabels{\arabic}{enumi}{enumii}{}{.}%
\item {} 
\sphinxAtStartPar
Create the COPT environment and model.

\item {} 
\sphinxAtStartPar
Add required data.

\item {} 
\sphinxAtStartPar
Construct the semidefinite programming model:
\begin{itemize}
\item {} 
\sphinxAtStartPar
Add decision variables (semidefinite and non\sphinxhyphen{}semidefinite variables).

\item {} 
\sphinxAtStartPar
Add constraints (including semidefinite constraints).

\item {} 
\sphinxAtStartPar
Set the objective function.

\end{itemize}

\item {} 
\sphinxAtStartPar
Set optimization parameters and solve.

\item {} 
\sphinxAtStartPar
Retrieve the solution.

\end{enumerate}

\sphinxAtStartPar
\sphinxstylestrong{Modeling: Adding Semidefinite Variables and Semidefinite Constraints}

\sphinxAtStartPar
In COPT, semidefinite variables are added by specifying the dimension ( \sphinxcode{\sphinxupquote{dim}} ) of the semidefinite variables:
\begin{enumerate}
\sphinxsetlistlabels{\arabic}{enumi}{enumii}{}{.}%
\item {} 
\sphinxAtStartPar
Add a single semidefinite variable.

\item {} 
\sphinxAtStartPar
Add multiple semidefinite variables.

\end{enumerate}

\sphinxAtStartPar
COPT provides two ways to add semidefinite constraints:
\begin{enumerate}
\sphinxsetlistlabels{\arabic}{enumi}{enumii}{}{.}%
\item {} 
\sphinxAtStartPar
First construct the semidefinite constraint expression by combining the semidefinite terms
(semidefinite variables and their corresponding coefficient matrices) and the linear terms, and then add them to the model.

\item {} 
\sphinxAtStartPar
Directly provide the semidefinite constraint expression (or semidefinite constraint builder) as a argument input, and add it to the model.

\end{enumerate}

\sphinxAtStartPar
When adding semidefinite constraints to the model, the main arguments that can be specified are:
\begin{itemize}
\item {} 
\sphinxAtStartPar
\sphinxcode{\sphinxupquote{expr}} / \sphinxcode{\sphinxupquote{builder}} : The semidefinite constraint expression or semidefinite constraint builder.

\item {} 
\sphinxAtStartPar
\sphinxcode{\sphinxupquote{rhs}} : The right\sphinxhyphen{}hand side of the semidefinite constraint.

\item {} 
\sphinxAtStartPar
\sphinxcode{\sphinxupquote{sense}} : The constraint type. For possible values, refer to {\hyperref[\detokenize{constant:chapconst-constrtype}]{\sphinxcrossref{\DUrole{std,std-ref}{Constants Chapter: Constraint Type}}}}.

\item {} 
\sphinxAtStartPar
\sphinxcode{\sphinxupquote{name}} : The name of the semidefinite constraint.

\end{itemize}

\sphinxAtStartPar
The implementation in different programming interfaces is shown in \hyperref[\detokenize{modeling:copttab-modelingsdp}]{Table \ref{\detokenize{modeling:copttab-modelingsdp}}}:

\begin{savenotes}\sphinxattablestart
\sphinxthistablewithglobalstyle
\centering
\sphinxcapstartof{table}
\sphinxthecaptionisattop
\sphinxcaption{Functions for Adding Semidefinite Variables/Constraints}\label{\detokenize{modeling:copttab-modelingsdp}}
\sphinxaftertopcaption
\begin{tabular}[t]{|\X{7}{50}|\X{25}{50}|\X{18}{50}|}
\sphinxtoprule
\sphinxstyletheadfamily 
\sphinxAtStartPar
API
&\sphinxstyletheadfamily 
\sphinxAtStartPar
Add Semidefinite Variable
&\sphinxstyletheadfamily 
\sphinxAtStartPar
Add Semidefinite Constraint
\\
\sphinxmidrule
\sphinxtableatstartofbodyhook
\sphinxAtStartPar
C
&
\sphinxAtStartPar
\sphinxcode{\sphinxupquote{COPT\_AddPSDCol}} / \sphinxcode{\sphinxupquote{COPT\_AddPSDCols}}
&
\sphinxAtStartPar
\sphinxcode{\sphinxupquote{COPT\_AddPSDConstr}}
\\
\sphinxhline
\sphinxAtStartPar
C++
&
\sphinxAtStartPar
\sphinxcode{\sphinxupquote{Model::AddPsdVar()}} / \sphinxcode{\sphinxupquote{Model.AddPsdVars()}}
&
\sphinxAtStartPar
\sphinxcode{\sphinxupquote{Model::AddPsdConstr()}}
\\
\sphinxhline
\sphinxAtStartPar
C\#
&
\sphinxAtStartPar
\sphinxcode{\sphinxupquote{Model.AddPsdVar()}} / \sphinxcode{\sphinxupquote{Model.AddPsdVars()}}
&
\sphinxAtStartPar
\sphinxcode{\sphinxupquote{Model.AddPsdConstr()}}
\\
\sphinxhline
\sphinxAtStartPar
Java
&
\sphinxAtStartPar
\sphinxcode{\sphinxupquote{Model.addPsdVar()}} / \sphinxcode{\sphinxupquote{Model.addPsdVars()}}
&
\sphinxAtStartPar
\sphinxcode{\sphinxupquote{Model.addPsdConstr()}}
\\
\sphinxhline
\sphinxAtStartPar
Python
&
\sphinxAtStartPar
\sphinxcode{\sphinxupquote{Model.addPsdVar()}} / \sphinxcode{\sphinxupquote{Model.addPsdVars()}}
&
\sphinxAtStartPar
\sphinxcode{\sphinxupquote{Model.addPsdConstr()}}
\\
\sphinxbottomrule
\end{tabular}
\sphinxtableafterendhook\par
\sphinxattableend\end{savenotes}

\begin{sphinxadmonition}{note}{Notes}
\begin{itemize}
\item {} 
\sphinxAtStartPar
Regarding the operations for modeling semidefinite variables and semidefinite constraints, the function names, calling methods,
and \sphinxstylestrong{argument names differ slightly} among different programming interfaces, but the functionality and argument meanings are consistent.

\item {} 
\sphinxAtStartPar
In the C API, when adding semidefinite constraints, arguments such as the indices of the non\sphinxhyphen{}zero linear term coefficients
and the semidefinite variable indices need to be provided. For more details, refer to function \sphinxcode{\sphinxupquote{COPT\_AddPSDConstr}} .

\item {} 
\sphinxAtStartPar
In the Python API, \sphinxcode{\sphinxupquote{Model.addConstr()}} can be used to add a linear constraint, a semidefinite constraint, or an indicator constraint to the model. For more details, please refer to {\hyperref[\detokenize{pyapiref:chappyapi-model}]{\sphinxcrossref{\DUrole{std,std-ref}{Python API Functions: Model Class}}}}.

\end{itemize}
\end{sphinxadmonition}

\sphinxAtStartPar
\sphinxstylestrong{Examples of PSD varibales and PSD constraints}
\begin{equation}\label{equation:modeling:modeling:16}
\begin{split}&A\bullet X + x_1 + x_2 = 0.6 \\
&x_1 \geq 0, x_2 \geq 0, X\succeq 0\end{split}
\end{equation}
\sphinxAtStartPar
where
\begin{equation}\label{equation:modeling:modeling:17}
\begin{split}&A=\begin{pmatrix}
1 & 1 & 1 \\
1 & 1 & 1 \\
1 & 1 & 1
\end{pmatrix}\end{split}
\end{equation}
\sphinxAtStartPar
The code implementation in different programming interfaces is as follows:

\sphinxAtStartPar
C API:

\begin{sphinxVerbatim}[commandchars=\\\{\}]
\PYG{p}{\PYGZob{}}
\PYG{+w}{  }\PYG{c+cm}{/* Matrix A */}
\PYG{+w}{  }\PYG{k+kt}{int}\PYG{+w}{ }\PYG{n}{ndim}\PYG{+w}{ }\PYG{o}{=}\PYG{+w}{ }\PYG{l+m+mi}{3}\PYG{p}{;}
\PYG{+w}{  }\PYG{k+kt}{int}\PYG{+w}{ }\PYG{n}{nelem}\PYG{+w}{ }\PYG{o}{=}\PYG{+w}{ }\PYG{l+m+mi}{6}\PYG{p}{;}
\PYG{+w}{  }\PYG{k+kt}{int}\PYG{+w}{ }\PYG{n}{rows}\PYG{p}{[}\PYG{p}{]}\PYG{+w}{ }\PYG{o}{=}\PYG{+w}{ }\PYG{p}{\PYGZob{}}\PYG{l+m+mi}{0}\PYG{p}{,}\PYG{+w}{ }\PYG{l+m+mi}{1}\PYG{p}{,}\PYG{+w}{ }\PYG{l+m+mi}{2}\PYG{p}{,}\PYG{+w}{ }\PYG{l+m+mi}{1}\PYG{p}{,}\PYG{+w}{ }\PYG{l+m+mi}{2}\PYG{p}{,}\PYG{+w}{ }\PYG{l+m+mi}{2}\PYG{p}{\PYGZcb{}}\PYG{p}{;}
\PYG{+w}{  }\PYG{k+kt}{int}\PYG{+w}{ }\PYG{n}{cols}\PYG{p}{[}\PYG{p}{]}\PYG{+w}{ }\PYG{o}{=}\PYG{+w}{ }\PYG{p}{\PYGZob{}}\PYG{l+m+mi}{0}\PYG{p}{,}\PYG{+w}{ }\PYG{l+m+mi}{0}\PYG{p}{,}\PYG{+w}{ }\PYG{l+m+mi}{0}\PYG{p}{,}\PYG{+w}{ }\PYG{l+m+mi}{1}\PYG{p}{,}\PYG{+w}{ }\PYG{l+m+mi}{1}\PYG{p}{,}\PYG{+w}{ }\PYG{l+m+mi}{2}\PYG{p}{\PYGZcb{}}\PYG{p}{;}
\PYG{+w}{  }\PYG{k+kt}{double}\PYG{+w}{ }\PYG{n}{elems}\PYG{p}{[}\PYG{p}{]}\PYG{+w}{ }\PYG{o}{=}\PYG{+w}{ }\PYG{p}{\PYGZob{}}\PYG{l+m+mf}{1.0}\PYG{p}{,}\PYG{+w}{ }\PYG{l+m+mf}{1.0}\PYG{p}{,}\PYG{+w}{ }\PYG{l+m+mf}{1.0}\PYG{p}{,}\PYG{+w}{ }\PYG{l+m+mf}{1.0}\PYG{p}{,}\PYG{+w}{ }\PYG{l+m+mf}{1.0}\PYG{p}{,}\PYG{+w}{ }\PYG{l+m+mf}{1.0}\PYG{p}{\PYGZcb{}}\PYG{p}{;}
\PYG{+w}{  }\PYG{n}{retcode}\PYG{+w}{ }\PYG{o}{=}\PYG{+w}{ }\PYG{n}{COPT\PYGZus{}AddSymMat}\PYG{p}{(}\PYG{n}{prob}\PYG{p}{,}\PYG{+w}{ }\PYG{n}{ndim}\PYG{p}{,}\PYG{+w}{ }\PYG{n}{nelem}\PYG{p}{,}\PYG{+w}{ }\PYG{n}{rows}\PYG{p}{,}\PYG{+w}{ }\PYG{n}{cols}\PYG{p}{,}\PYG{+w}{ }\PYG{n}{elems}\PYG{p}{)}\PYG{p}{;}
\PYG{+w}{  }\PYG{k}{if}\PYG{+w}{ }\PYG{p}{(}\PYG{n}{retcode}\PYG{p}{)}\PYG{+w}{ }\PYG{k}{goto}\PYG{+w}{ }\PYG{n}{exit\PYGZus{}cleanup}\PYG{p}{;}
\PYG{p}{\PYGZcb{}}
\PYG{c+cm}{/* Add PSD columns */}
\PYG{k+kt}{int}\PYG{+w}{ }\PYG{n}{nPSDCol}\PYG{+w}{ }\PYG{o}{=}\PYG{+w}{ }\PYG{l+m+mi}{1}\PYG{p}{;}
\PYG{k+kt}{int}\PYG{+w}{ }\PYG{n}{colDims}\PYG{p}{[}\PYG{p}{]}\PYG{+w}{ }\PYG{o}{=}\PYG{+w}{ }\PYG{p}{\PYGZob{}}\PYG{l+m+mi}{3}\PYG{p}{\PYGZcb{}}\PYG{p}{;}
\PYG{n}{retcode}\PYG{+w}{ }\PYG{o}{=}\PYG{+w}{ }\PYG{n}{COPT\PYGZus{}AddPSDCols}\PYG{p}{(}\PYG{n}{prob}\PYG{p}{,}\PYG{+w}{ }\PYG{n}{nPSDCol}\PYG{p}{,}\PYG{+w}{ }\PYG{n}{colDims}\PYG{p}{,}\PYG{+w}{ }\PYG{n+nb}{NULL}\PYG{p}{)}\PYG{p}{;}
\PYG{k}{if}\PYG{+w}{ }\PYG{p}{(}\PYG{n}{retcode}\PYG{p}{)}\PYG{+w}{ }\PYG{k}{goto}\PYG{+w}{ }\PYG{n}{exit\PYGZus{}cleanup}\PYG{p}{;}
\PYG{c+cm}{/* Add columns */}
\PYG{k+kt}{int}\PYG{+w}{ }\PYG{n}{nCol}\PYG{+w}{ }\PYG{o}{=}\PYG{+w}{ }\PYG{l+m+mi}{2}\PYG{p}{;}
\PYG{n}{retcode}\PYG{+w}{ }\PYG{o}{=}\PYG{+w}{ }\PYG{n}{COPT\PYGZus{}AddCols}\PYG{p}{(}\PYG{n}{prob}\PYG{p}{,}\PYG{+w}{ }\PYG{n}{nCol}\PYG{p}{,}\PYG{+w}{ }\PYG{n+nb}{NULL}\PYG{p}{,}\PYG{+w}{ }\PYG{n+nb}{NULL}\PYG{p}{,}\PYG{+w}{ }\PYG{n+nb}{NULL}\PYG{p}{,}\PYG{+w}{ }\PYG{n+nb}{NULL}\PYG{p}{,}\PYG{+w}{ }\PYG{n+nb}{NULL}\PYG{p}{,}\PYG{+w}{ }\PYG{n+nb}{NULL}\PYG{p}{,}
\PYG{n+nb}{NULL}\PYG{p}{,}\PYG{+w}{ }\PYG{n+nb}{NULL}\PYG{p}{,}\PYG{+w}{ }\PYG{n+nb}{NULL}\PYG{p}{)}\PYG{p}{;}
\PYG{k}{if}\PYG{+w}{ }\PYG{p}{(}\PYG{n}{retcode}\PYG{p}{)}\PYG{+w}{ }\PYG{k}{goto}\PYG{+w}{ }\PYG{n}{exit\PYGZus{}cleanup}\PYG{p}{;}
\PYG{c+cm}{/* Add PSD constraints */}
\PYG{p}{\PYGZob{}}
\PYG{+w}{  }\PYG{k+kt}{int}\PYG{+w}{ }\PYG{n}{nRowMatCnt}\PYG{+w}{ }\PYG{o}{=}\PYG{+w}{ }\PYG{l+m+mi}{2}\PYG{p}{;}
\PYG{+w}{  }\PYG{k+kt}{int}\PYG{+w}{ }\PYG{n}{rowMatIdx}\PYG{p}{[}\PYG{p}{]}\PYG{+w}{ }\PYG{o}{=}\PYG{+w}{ }\PYG{p}{\PYGZob{}}\PYG{l+m+mi}{0}\PYG{p}{,}\PYG{+w}{ }\PYG{l+m+mi}{1}\PYG{p}{\PYGZcb{}}\PYG{p}{;}
\PYG{+w}{  }\PYG{k+kt}{double}\PYG{+w}{ }\PYG{n}{rowMatElem}\PYG{p}{[}\PYG{p}{]}\PYG{+w}{ }\PYG{o}{=}\PYG{+w}{ }\PYG{p}{\PYGZob{}}\PYG{l+m+mf}{1.0}\PYG{p}{,}\PYG{+w}{ }\PYG{l+m+mf}{1.0}\PYG{p}{\PYGZcb{}}\PYG{p}{;}
\PYG{+w}{  }\PYG{k+kt}{int}\PYG{+w}{ }\PYG{n}{nColCnt}\PYG{+w}{ }\PYG{o}{=}\PYG{+w}{ }\PYG{l+m+mi}{1}\PYG{p}{;}
\PYG{+w}{  }\PYG{k+kt}{int}\PYG{+w}{ }\PYG{n}{psdColIdx}\PYG{p}{[}\PYG{p}{]}\PYG{+w}{ }\PYG{o}{=}\PYG{+w}{ }\PYG{p}{\PYGZob{}}\PYG{l+m+mi}{0}\PYG{p}{\PYGZcb{}}\PYG{p}{;}
\PYG{+w}{  }\PYG{k+kt}{int}\PYG{+w}{ }\PYG{n}{symMatIdx}\PYG{p}{[}\PYG{p}{]}\PYG{+w}{ }\PYG{o}{=}\PYG{+w}{ }\PYG{p}{\PYGZob{}}\PYG{l+m+mi}{2}\PYG{p}{\PYGZcb{}}\PYG{p}{;}
\PYG{+w}{  }\PYG{k+kt}{char}\PYG{+w}{ }\PYG{n}{cRowSense}\PYG{+w}{ }\PYG{o}{=}\PYG{+w}{ }\PYG{n}{COPT\PYGZus{}EQUAL}\PYG{p}{;}
\PYG{+w}{  }\PYG{k+kt}{double}\PYG{+w}{ }\PYG{n}{dRowBound}\PYG{+w}{ }\PYG{o}{=}\PYG{+w}{ }\PYG{l+m+mf}{0.6}\PYG{p}{;}
\PYG{+w}{  }\PYG{n}{retcode}\PYG{+w}{ }\PYG{o}{=}\PYG{+w}{ }\PYG{n}{COPT\PYGZus{}AddPSDConstr}\PYG{p}{(}\PYG{n}{prob}\PYG{p}{,}\PYG{+w}{ }\PYG{n}{nRowMatCnt}\PYG{p}{,}\PYG{+w}{ }\PYG{n}{rowMatIdx}\PYG{p}{,}\PYG{+w}{ }\PYG{n}{rowMatElem}\PYG{p}{,}\PYG{+w}{ }\PYG{n}{nColCnt}\PYG{p}{,}
\PYG{+w}{    }\PYG{n}{psdColIdx}\PYG{p}{,}\PYG{+w}{ }\PYG{n}{symMatIdx}\PYG{p}{,}\PYG{+w}{ }\PYG{n}{cRowSense}\PYG{p}{,}\PYG{+w}{ }\PYG{n}{dRowBound}\PYG{p}{,}\PYG{+w}{ }\PYG{l+m+mf}{0.0}\PYG{p}{,}\PYG{+w}{ }\PYG{n+nb}{NULL}\PYG{p}{)}\PYG{p}{;}
\PYG{+w}{  }\PYG{k}{if}\PYG{+w}{ }\PYG{p}{(}\PYG{n}{retcode}\PYG{p}{)}\PYG{+w}{ }\PYG{k}{goto}\PYG{+w}{ }\PYG{n}{exit\PYGZus{}cleanup}\PYG{p}{;}
\PYG{p}{\PYGZcb{}}
\end{sphinxVerbatim}

\sphinxAtStartPar
C++ API:

\begin{sphinxVerbatim}[commandchars=\\\{\}]
\PYG{n}{SymMatrix}\PYG{+w}{ }\PYG{n}{A}\PYG{+w}{ }\PYG{o}{=}\PYG{+w}{ }\PYG{n}{model}\PYG{p}{.}\PYG{n}{AddOnesMat}\PYG{p}{(}\PYG{l+m+mi}{3}\PYG{p}{)}\PYG{p}{;}
\PYG{n}{Var}\PYG{+w}{ }\PYG{n}{x1}\PYG{+w}{ }\PYG{o}{=}\PYG{+w}{ }\PYG{n}{model}\PYG{p}{.}\PYG{n}{AddVar}\PYG{p}{(}\PYG{l+m+mf}{0.0}\PYG{p}{,}\PYG{+w}{ }\PYG{n}{COPT\PYGZus{}INFINITY}\PYG{p}{,}\PYG{+w}{ }\PYG{l+m+mi}{0}\PYG{p}{,}\PYG{+w}{ }\PYG{n}{COPT\PYGZus{}CONTINUOUS}\PYG{p}{,}\PYG{+w}{ }\PYG{l+s}{\PYGZdq{}}\PYG{l+s}{x1}\PYG{l+s}{\PYGZdq{}}\PYG{p}{)}\PYG{p}{;}
\PYG{n}{Var}\PYG{+w}{ }\PYG{n}{x2}\PYG{+w}{ }\PYG{o}{=}\PYG{+w}{ }\PYG{n}{model}\PYG{p}{.}\PYG{n}{AddVar}\PYG{p}{(}\PYG{l+m+mf}{0.0}\PYG{p}{,}\PYG{+w}{ }\PYG{n}{COPT\PYGZus{}INFINITY}\PYG{p}{,}\PYG{+w}{ }\PYG{l+m+mi}{0}\PYG{p}{,}\PYG{+w}{ }\PYG{n}{COPT\PYGZus{}CONTINUOUS}\PYG{p}{,}\PYG{+w}{ }\PYG{l+s}{\PYGZdq{}}\PYG{l+s}{x2}\PYG{l+s}{\PYGZdq{}}\PYG{p}{)}\PYG{p}{;}
\PYG{n}{PsdVar}\PYG{+w}{ }\PYG{n}{barX}\PYG{+w}{ }\PYG{o}{=}\PYG{+w}{ }\PYG{n}{model}\PYG{p}{.}\PYG{n}{AddPsdVar}\PYG{p}{(}\PYG{l+m+mi}{3}\PYG{p}{,}\PYG{+w}{ }\PYG{l+s}{\PYGZdq{}}\PYG{l+s}{X}\PYG{l+s}{\PYGZdq{}}\PYG{p}{)}\PYG{p}{;}
\PYG{n}{model}\PYG{p}{.}\PYG{n}{AddPsdConstr}\PYG{p}{(}\PYG{n}{A}\PYG{+w}{ }\PYG{o}{*}\PYG{+w}{ }\PYG{n}{barX}\PYG{+w}{ }\PYG{o}{+}\PYG{+w}{ }\PYG{n}{x1}\PYG{+w}{ }\PYG{o}{+}\PYG{+w}{ }\PYG{n}{x2}\PYG{+w}{ }\PYG{o}{=}\PYG{o}{=}\PYG{+w}{ }\PYG{l+m+mf}{0.6}\PYG{p}{,}\PYG{+w}{ }\PYG{l+s}{\PYGZdq{}}\PYG{l+s}{PSD\PYGZus{}R}\PYG{l+s}{\PYGZdq{}}\PYG{p}{)}\PYG{p}{;}
\end{sphinxVerbatim}

\sphinxAtStartPar
C\# API:

\begin{sphinxVerbatim}[commandchars=\\\{\}]
\PYG{n}{SymMatrix}\PYG{+w}{ }\PYG{n}{A}\PYG{+w}{ }\PYG{o}{=}\PYG{+w}{ }\PYG{n}{model}\PYG{p}{.}\PYG{n}{AddOnesMat}\PYG{p}{(}\PYG{l+m+mi}{3}\PYG{p}{)}\PYG{p}{;}
\PYG{n}{Var}\PYG{+w}{ }\PYG{n}{x1}\PYG{+w}{ }\PYG{o}{=}\PYG{+w}{ }\PYG{n}{model}\PYG{p}{.}\PYG{n}{AddVar}\PYG{p}{(}\PYG{l+m+mf}{0.0}\PYG{p}{,}\PYG{+w}{ }\PYG{n}{Copt}\PYG{p}{.}\PYG{n}{Consts}\PYG{p}{.}\PYG{n}{INFINITY}\PYG{p}{,}\PYG{+w}{ }\PYG{l+m+mi}{0}\PYG{p}{,}\PYG{+w}{ }\PYG{n}{Copt}\PYG{p}{.}\PYG{n}{Consts}\PYG{p}{.}\PYG{n}{CONTINUOUS}\PYG{p}{,}\PYG{+w}{ }\PYG{l+s}{\PYGZdq{}x1\PYGZdq{}}\PYG{p}{)}\PYG{p}{;}
\PYG{n}{Var}\PYG{+w}{ }\PYG{n}{x2}\PYG{+w}{ }\PYG{o}{=}\PYG{+w}{ }\PYG{n}{model}\PYG{p}{.}\PYG{n}{AddVar}\PYG{p}{(}\PYG{l+m+mf}{0.0}\PYG{p}{,}\PYG{+w}{ }\PYG{n}{Copt}\PYG{p}{.}\PYG{n}{Consts}\PYG{p}{.}\PYG{n}{INFINITY}\PYG{p}{,}\PYG{+w}{ }\PYG{l+m+mi}{0}\PYG{p}{,}\PYG{+w}{ }\PYG{n}{Copt}\PYG{p}{.}\PYG{n}{Consts}\PYG{p}{.}\PYG{n}{CONTINUOUS}\PYG{p}{,}\PYG{+w}{ }\PYG{l+s}{\PYGZdq{}x2\PYGZdq{}}\PYG{p}{)}\PYG{p}{;}
\PYG{n}{PsdVar}\PYG{+w}{ }\PYG{n}{barX}\PYG{+w}{ }\PYG{o}{=}\PYG{+w}{ }\PYG{n}{model}\PYG{p}{.}\PYG{n}{AddPsdVar}\PYG{p}{(}\PYG{l+m+mi}{3}\PYG{p}{,}\PYG{+w}{ }\PYG{l+s}{\PYGZdq{}X\PYGZdq{}}\PYG{p}{)}\PYG{p}{;}
\PYG{n}{model}\PYG{p}{.}\PYG{n}{AddPsdConstr}\PYG{p}{(}\PYG{n}{A}\PYG{+w}{ }\PYG{o}{*}\PYG{+w}{ }\PYG{n}{barX}\PYG{+w}{ }\PYG{o}{+}\PYG{+w}{ }\PYG{n}{x1}\PYG{+w}{ }\PYG{o}{+}\PYG{+w}{ }\PYG{n}{x2}\PYG{+w}{ }\PYG{o}{==}\PYG{+w}{ }\PYG{l+m+mf}{0.6}\PYG{p}{,}\PYG{+w}{ }\PYG{l+s}{\PYGZdq{}PSD\PYGZus{}R\PYGZdq{}}\PYG{p}{)}\PYG{p}{;}
\end{sphinxVerbatim}

\sphinxAtStartPar
Java API:

\begin{sphinxVerbatim}[commandchars=\\\{\}]
\PYG{n}{SymMatrix}\PYG{+w}{ }\PYG{n}{A}\PYG{+w}{ }\PYG{o}{=}\PYG{+w}{ }\PYG{n}{model}\PYG{p}{.}\PYG{n+na}{addOnesMat}\PYG{p}{(}\PYG{l+m+mi}{3}\PYG{p}{)}\PYG{p}{;}
\PYG{n}{Var}\PYG{+w}{ }\PYG{n}{x1}\PYG{+w}{ }\PYG{o}{=}\PYG{+w}{ }\PYG{n}{model}\PYG{p}{.}\PYG{n+na}{addVar}\PYG{p}{(}\PYG{l+m+mf}{0.0}\PYG{p}{,}\PYG{+w}{ }\PYG{n}{copt}\PYG{p}{.}\PYG{n+na}{Consts}\PYG{p}{.}\PYG{n+na}{INFINITY}\PYG{p}{,}\PYG{+w}{ }\PYG{l+m+mi}{0}\PYG{p}{,}\PYG{+w}{ }\PYG{n}{copt}\PYG{p}{.}\PYG{n+na}{Consts}\PYG{p}{.}\PYG{n+na}{CONTINUOUS}\PYG{p}{,}\PYG{+w}{ }\PYG{l+s}{\PYGZdq{}}\PYG{l+s}{x1}\PYG{l+s}{\PYGZdq{}}\PYG{p}{)}\PYG{p}{;}
\PYG{n}{Var}\PYG{+w}{ }\PYG{n}{x2}\PYG{+w}{ }\PYG{o}{=}\PYG{+w}{ }\PYG{n}{model}\PYG{p}{.}\PYG{n+na}{addVar}\PYG{p}{(}\PYG{l+m+mf}{0.0}\PYG{p}{,}\PYG{+w}{ }\PYG{n}{copt}\PYG{p}{.}\PYG{n+na}{Consts}\PYG{p}{.}\PYG{n+na}{INFINITY}\PYG{p}{,}\PYG{+w}{ }\PYG{l+m+mi}{0}\PYG{p}{,}\PYG{+w}{ }\PYG{n}{copt}\PYG{p}{.}\PYG{n+na}{Consts}\PYG{p}{.}\PYG{n+na}{CONTINUOUS}\PYG{p}{,}\PYG{+w}{ }\PYG{l+s}{\PYGZdq{}}\PYG{l+s}{x2}\PYG{l+s}{\PYGZdq{}}\PYG{p}{)}\PYG{p}{;}
\PYG{n}{PsdVar}\PYG{+w}{ }\PYG{n}{barX}\PYG{+w}{ }\PYG{o}{=}\PYG{+w}{ }\PYG{n}{model}\PYG{p}{.}\PYG{n+na}{addPsdVar}\PYG{p}{(}\PYG{l+m+mi}{3}\PYG{p}{,}\PYG{+w}{ }\PYG{l+s}{\PYGZdq{}}\PYG{l+s}{X}\PYG{l+s}{\PYGZdq{}}\PYG{p}{)}\PYG{p}{;}
\PYG{n}{PsdExpr}\PYG{+w}{ }\PYG{n}{pexpr}\PYG{+w}{ }\PYG{o}{=}\PYG{+w}{ }\PYG{k}{new}\PYG{+w}{ }\PYG{n}{PsdExpr}\PYG{p}{(}\PYG{n}{x1}\PYG{p}{,}\PYG{+w}{ }\PYG{l+m+mf}{1.0}\PYG{p}{)}\PYG{p}{;}
\PYG{n}{pexpr}\PYG{p}{.}\PYG{n+na}{addTerm}\PYG{p}{(}\PYG{n}{x2}\PYG{p}{,}\PYG{+w}{ }\PYG{l+m+mf}{1.0}\PYG{p}{)}\PYG{p}{;}
\PYG{n}{pexpr}\PYG{p}{.}\PYG{n+na}{addTerm}\PYG{p}{(}\PYG{n}{barX}\PYG{p}{,}\PYG{+w}{ }\PYG{n}{A}\PYG{p}{)}\PYG{p}{;}
\PYG{n}{model}\PYG{p}{.}\PYG{n+na}{addPsdConstr}\PYG{p}{(}\PYG{n}{pexpr}\PYG{p}{,}\PYG{+w}{ }\PYG{n}{copt}\PYG{p}{.}\PYG{n+na}{Consts}\PYG{p}{.}\PYG{n+na}{EQUAL}\PYG{p}{,}\PYG{+w}{ }\PYG{l+m+mf}{0.6}\PYG{p}{,}\PYG{+w}{ }\PYG{l+s}{\PYGZdq{}}\PYG{l+s}{PSD\PYGZus{}R}\PYG{l+s}{\PYGZdq{}}\PYG{p}{)}\PYG{p}{;}
\end{sphinxVerbatim}

\sphinxAtStartPar
Python API:

\begin{sphinxVerbatim}[commandchars=\\\{\}]
\PYG{n}{A} \PYG{o}{=} \PYG{n}{m}\PYG{o}{.}\PYG{n}{addOnesMat}\PYG{p}{(}\PYG{l+m+mi}{3}\PYG{p}{)}
\PYG{n}{x1} \PYG{o}{=} \PYG{n}{m}\PYG{o}{.}\PYG{n}{addVar}\PYG{p}{(}\PYG{n}{lb}\PYG{o}{=}\PYG{l+m+mf}{0.0}\PYG{p}{,} \PYG{n}{ub}\PYG{o}{=}\PYG{n}{COPT}\PYG{o}{.}\PYG{n}{INFINITY}\PYG{p}{,} \PYG{n}{name}\PYG{o}{=}\PYG{l+s+s2}{\PYGZdq{}}\PYG{l+s+s2}{x1}\PYG{l+s+s2}{\PYGZdq{}}\PYG{p}{)}
\PYG{n}{x2} \PYG{o}{=} \PYG{n}{m}\PYG{o}{.}\PYG{n}{addVar}\PYG{p}{(}\PYG{n}{lb}\PYG{o}{=}\PYG{l+m+mf}{0.0}\PYG{p}{,} \PYG{n}{ub}\PYG{o}{=}\PYG{n}{COPT}\PYG{o}{.}\PYG{n}{INFINITY}\PYG{p}{,} \PYG{n}{name}\PYG{o}{=}\PYG{l+s+s2}{\PYGZdq{}}\PYG{l+s+s2}{x2}\PYG{l+s+s2}{\PYGZdq{}}\PYG{p}{)}
\PYG{n}{X} \PYG{o}{=} \PYG{n}{m}\PYG{o}{.}\PYG{n}{addPsdVars}\PYG{p}{(}\PYG{l+m+mi}{3}\PYG{p}{,} \PYG{l+s+s2}{\PYGZdq{}}\PYG{l+s+s2}{BAR\PYGZus{}X}\PYG{l+s+s2}{\PYGZdq{}}\PYG{p}{)}
\PYG{n}{psdc} \PYG{o}{=} \PYG{n}{model}\PYG{o}{.}\PYG{n}{addConstr}\PYG{p}{(}\PYG{n}{A} \PYG{o}{*} \PYG{n}{X} \PYG{o}{+} \PYG{n}{x1} \PYG{o}{+} \PYG{n}{x2} \PYG{o}{==} \PYG{l+m+mf}{0.6}\PYG{p}{,} \PYG{n}{name}\PYG{o}{=}\PYG{l+s+s2}{\PYGZdq{}}\PYG{l+s+s2}{PSD\PYGZus{}C}\PYG{l+s+s2}{\PYGZdq{}}\PYG{p}{)}
\end{sphinxVerbatim}

\sphinxAtStartPar
In the programming interfaces provided by COPT, except for the C language, other object\sphinxhyphen{}oriented programming interfaces (C\#, C++, Java, Python)
offer classes related to semidefinite constraints:
\begin{itemize}
\item {} 
\sphinxAtStartPar
Encapsulation of operations for constructing semidefinite expressions
\begin{enumerate}
\sphinxsetlistlabels{\arabic}{enumi}{enumii}{}{.}%
\item {} 
\sphinxAtStartPar
\sphinxcode{\sphinxupquote{PsdExpr}} Class: Encapsulation of operations related to the combination of semidefinite variables when constructing semidefinite expressions in COPT.

\end{enumerate}

\item {} 
\sphinxAtStartPar
Encapsulation of operations related to semidefinite constraints
\begin{enumerate}
\sphinxsetlistlabels{\arabic}{enumi}{enumii}{}{.}%
\item {} 
\sphinxAtStartPar
\sphinxcode{\sphinxupquote{PsdConstraint}} Class: Encapsulation of operations related to semidefinite constraints in COPT.

\item {} 
\sphinxAtStartPar
\sphinxcode{\sphinxupquote{PsdConstrArray}} Class: Facilitates user operations on a group of \sphinxcode{\sphinxupquote{PsdConstraint}} class objects.

\end{enumerate}

\item {} 
\sphinxAtStartPar
Encapsulation of semidefinite constraint builders
\begin{enumerate}
\sphinxsetlistlabels{\arabic}{enumi}{enumii}{}{.}%
\item {} 
\sphinxAtStartPar
\sphinxcode{\sphinxupquote{PsdConstrBuilder}} Class: Encapsulation for constructing semidefinite constraint builders in COPT.

\item {} 
\sphinxAtStartPar
\sphinxcode{\sphinxupquote{PsdConstrBuilderArray}} Class: Facilitates user operations on a group of \sphinxcode{\sphinxupquote{PsdConstrBuilder}} class objects.

\end{enumerate}

\item {} 
\sphinxAtStartPar
C++ API: {\hyperref[\detokenize{cppapiref:chapcppapiref-psdexpr}]{\sphinxcrossref{\DUrole{std,std-ref}{PsdExpr Class}}}} , {\hyperref[\detokenize{cppapiref:chapcppapiref-psdconstrarray}]{\sphinxcrossref{\DUrole{std,std-ref}{PsdConstrArray Class}}}} , {\hyperref[\detokenize{cppapiref:chapcppapiref-psdconstrbuilder}]{\sphinxcrossref{\DUrole{std,std-ref}{PsdConstrBuilder Class}}}} , {\hyperref[\detokenize{cppapiref:chapcppapiref-psdconstrbuilderarray}]{\sphinxcrossref{\DUrole{std,std-ref}{PsdConstrBuilderArray Class}}}}

\item {} 
\sphinxAtStartPar
C\# API:  {\hyperref[\detokenize{csharpapiref:chapcsharpapiref-psdexpr}]{\sphinxcrossref{\DUrole{std,std-ref}{PsdExpr Class}}}} , {\hyperref[\detokenize{csharpapiref:chapcsharpapiref-psdconstraint}]{\sphinxcrossref{\DUrole{std,std-ref}{PsdConstraint Class}}}} , {\hyperref[\detokenize{csharpapiref:chapcsharpapiref-psdconstrarray}]{\sphinxcrossref{\DUrole{std,std-ref}{PsdConstrArray Class}}}} , {\hyperref[\detokenize{csharpapiref:chapcsharpapiref-psdconstrbuilder}]{\sphinxcrossref{\DUrole{std,std-ref}{PsdConstrBuilder Class}}}} , {\hyperref[\detokenize{csharpapiref:chapcsharpapiref-psdconstrbuilderarray}]{\sphinxcrossref{\DUrole{std,std-ref}{PsdConstrBuilderArray Class}}}}

\item {} 
\sphinxAtStartPar
Java API:  {\hyperref[\detokenize{javaapiref:chapjavaapiref-psdexpr}]{\sphinxcrossref{\DUrole{std,std-ref}{PsdExpr Class}}}} , {\hyperref[\detokenize{javaapiref:chapjavaapiref-psdconstraint}]{\sphinxcrossref{\DUrole{std,std-ref}{PsdConstraint Class}}}} , {\hyperref[\detokenize{javaapiref:chapjavaapiref-psdconstrarray}]{\sphinxcrossref{\DUrole{std,std-ref}{PsdConstrArray Class}}}} , {\hyperref[\detokenize{javaapiref:chapjavaapiref-psdconstrbuilder}]{\sphinxcrossref{\DUrole{std,std-ref}{PsdConstrBuilder Class}}}} , {\hyperref[\detokenize{javaapiref:chapjavaapiref-psdconstrbuilderarray}]{\sphinxcrossref{\DUrole{std,std-ref}{PsdConstrBuilderArray Class}}}}

\item {} 
\sphinxAtStartPar
Python API:  {\hyperref[\detokenize{pyapiref:chappyapi-psdexpr}]{\sphinxcrossref{\DUrole{std,std-ref}{PsdExpr Class}}}} , {\hyperref[\detokenize{pyapiref:chappyapi-psdconstraint}]{\sphinxcrossref{\DUrole{std,std-ref}{PsdConstraint Class}}}} , {\hyperref[\detokenize{pyapiref:chappyapi-psdconstrarray}]{\sphinxcrossref{\DUrole{std,std-ref}{PsdConstrArray Class}}}} , {\hyperref[\detokenize{pyapiref:chappyapi-psdconstrbuilder}]{\sphinxcrossref{\DUrole{std,std-ref}{PsdConstrBuilder Class}}}} , {\hyperref[\detokenize{pyapiref:chappyapi-psdconstrbuilderarray}]{\sphinxcrossref{\DUrole{std,std-ref}{PsdConstrBuilderArray Class}}}}

\end{itemize}

\subsection{Solving}
\label{\detokenize{modeling:id5}}
\sphinxAtStartPar
For SDP problems, COPT provides the Barrier method and the ADMM method.
You can specify the solving method by setting the optimization parameter \sphinxcode{\sphinxupquote{SDPMethod}}.

\sphinxAtStartPar
For more details, please refer to {\hyperref[\detokenize{parameter:chapparam-sdp}]{\sphinxcrossref{\DUrole{std,std-ref}{Parameter Section: Semidefinite Programming Related}}}}.

\subsection{Related Attributes}
\label{\detokenize{modeling:related-attributes}}
\sphinxAtStartPar
Attributes related to SDP model are as shown in \hyperref[\detokenize{modeling:copttab-modelingsdp-attr}]{Table \ref{\detokenize{modeling:copttab-modelingsdp-attr}}}:

\begin{savenotes}\sphinxattablestart
\sphinxthistablewithglobalstyle
\centering
\sphinxcapstartof{table}
\sphinxthecaptionisattop
\sphinxcaption{Attributes for Semidefinite Variables and Constraints}\label{\detokenize{modeling:copttab-modelingsdp-attr}}
\sphinxaftertopcaption
\begin{tabular}[t]{|\X{15}{59}|\X{9}{59}|\X{35}{59}|}
\sphinxtoprule
\sphinxstyletheadfamily 
\sphinxAtStartPar
Name
&\sphinxstyletheadfamily 
\sphinxAtStartPar
Type
&\sphinxstyletheadfamily 
\sphinxAtStartPar
Description
\\
\sphinxmidrule
\sphinxtableatstartofbodyhook
\sphinxAtStartPar
\sphinxcode{\sphinxupquote{PSDCols}}
&
\sphinxAtStartPar
Integer
&
\sphinxAtStartPar
Number of PSD variables in the problem
\\
\sphinxhline
\sphinxAtStartPar
\sphinxcode{\sphinxupquote{PSDElems}}
&
\sphinxAtStartPar
Integer
&
\sphinxAtStartPar
Number of PSD terms in objective function
\\
\sphinxhline
\sphinxAtStartPar
\sphinxcode{\sphinxupquote{SymMats}}
&
\sphinxAtStartPar
Integer
&
\sphinxAtStartPar
Number of symmetric matrices in the problem
\\
\sphinxhline
\sphinxAtStartPar
\sphinxcode{\sphinxupquote{PSDConstrs}}
&
\sphinxAtStartPar
Integer
&
\sphinxAtStartPar
Number of PSD constraints
\\
\sphinxhline
\sphinxAtStartPar
\sphinxcode{\sphinxupquote{HasPSDObj}}
&
\sphinxAtStartPar
Integer
&
\sphinxAtStartPar
Whether the problem has PSD terms in objective function
\\
\sphinxbottomrule
\end{tabular}
\sphinxtableafterendhook\par
\sphinxattableend\end{savenotes}

\section{Quadratic Programming (QP)}
\label{\detokenize{modeling:quadratic-programming-qp}}\label{\detokenize{modeling:chapmodeling-qp}}

\subsection{Mathematical Formulation}
\label{\detokenize{modeling:id6}}
\sphinxAtStartPar
Convex Quadratic Programming (QP) has a convex quadratic objective function and linear constraints.

\sphinxAtStartPar
The mathematical formulation is as follows:
\begin{equation}\label{equation:modeling:modeling:18}
\begin{split}\min\quad x^TQx + c^Tx + c^f\\
\text{s.t.}\quad &l_i^c\leq\sum_{j=0}^{n-1}a_{ij}x_j\leq u_i^c,\qquad i=0,\cdots,m-1\\
&l_j^v\leq x_j\leq u_j^v,\qquad j=0,1,\cdots,n-1\end{split}
\end{equation}
\sphinxAtStartPar
The variables and arguments in the model are defined as follows:
\begin{itemize}
\item {} 
\sphinxAtStartPar
Decision variables: \(x=(x_j)_{j=0}^{n-1}\in\mathbb{R}^n\);

\item {} 
\sphinxAtStartPar
Variable bounds: \(l^v, u^v\in\mathbb{R}^n\); where \(l^v\) denotes the lower bounds and \(u^v\) denotes the upper bounds of the variables;

\item {} 
\sphinxAtStartPar
Constraint bounds: \(l^c, u^c\in\mathbb{R}^m\); where \(l^c\) denotes the lower bounds and \(u^c\) denotes the upper bounds of the constraints;

\item {} 
\sphinxAtStartPar
Coefficient matrix of linear constraints: \(A=(a_{ij})_{m\times n}\in\mathbb{R}^{m\times n}\)

\item {} 
\sphinxAtStartPar
Coefficients in the objective function:
\begin{itemize}
\item {} 
\sphinxAtStartPar
\(Q\in\mathbb{R}^{n\times n}\) represents the coefficients of the quadratic term in the objective function

\item {} 
\sphinxAtStartPar
\(c\in\mathbb{R}^n\) represents the coefficients of the linear term in the objective function, and \(c^f\) represents the constant term in the objective function

\end{itemize}

\item {} 
\sphinxAtStartPar
Problem size: \(m\) denotes the number of constraints, and \(n\) denotes the number of decision variables.

\end{itemize}

\section{Quadratically Constrained Programming (QCP)}
\label{\detokenize{modeling:quadratically-constrained-programming-qcp}}\label{\detokenize{modeling:chapmodeling-qcp}}
\sphinxAtStartPar
Convex Quadratically Constrained Programming (QCP) consists of convex quadratic constraints.

\subsection{Mathematical Model}
\label{\detokenize{modeling:id7}}
\sphinxAtStartPar
The mathematical formulation is as follows:
\begin{equation}\label{equation:modeling:modeling:19}
\begin{split}\min\quad &\frac{1}{2} x^TQx + c^Tx + c^f\\
\text{s.t.}\quad&\frac{1}{2}x^TP_ix+a^T_ix_j\leq u_i^c,\qquad i=0,\cdots,m-1\\
&l^v\leq x_j\leq u^v,\qquad j=0,1,\cdots,n-1\end{split}
\end{equation}
\sphinxAtStartPar
The variables and parameters in the model are defined as follows:
\begin{itemize}
\item {} 
\sphinxAtStartPar
Decision variables: \(x=(x_j)_{j=0}^{n-1}\in\mathbb{R}^n\);

\item {} 
\sphinxAtStartPar
Variable bounds: \(l^v, u^v\in\mathbb{R}^n\), where \(l^v\) denotes the lower bounds and \(u^v\) denotes the upper bounds of the variables;

\item {} 
\sphinxAtStartPar
Constraint bounds: \(u^c\in\mathbb{R}^m\), where \(u^c\) denotes the upper bounds of the constraints;

\item {} 
\sphinxAtStartPar
Coefficients in the quadratic constraints:
\begin{itemize}
\item {} 
\sphinxAtStartPar
\(P_i\in \mathbb{R}^{n\times n}\) (for \(i=0,...,m-1\)),

\item {} 
\sphinxAtStartPar
\(a_i = (a_{ij})_{n}\in\mathbb{R}^{n}\);

\end{itemize}

\item {} 
\sphinxAtStartPar
Coefficients in the objective function:
\begin{itemize}
\item {} 
\sphinxAtStartPar
\(Q\in\mathbb{R}^{n\times n}\) represents the coefficients of the quadratic term in the objective function;

\item {} 
\sphinxAtStartPar
\(c\in\mathbb{R}^n\) represents the coefficients of the linear term in the objective function, and \(c^f\)
represents the constant term in the objective function.

\end{itemize}

\item {} 
\sphinxAtStartPar
Problem size: \(m\) denotes the number of constraints, and \(n\) denotes the number of decision variables.

\end{itemize}

\begin{sphinxadmonition}{note}{Notes}
\begin{itemize}
\item {} 
\sphinxAtStartPar
The COPT solver currently supports solving convex quadratic programming and convex quadratically constrained programming problems.
The matrices \(Q\) and \(P_i\) (for \(i=0,1,\cdots,n\)) must be \sphinxstylestrong{positive semidefinite}.

\item {} 
\sphinxAtStartPar
When the constraint type is \sphinxcode{\sphinxupquote{\textgreater{}=}} (of the form \(\quad\frac{1}{2}x^TP_ix+\sum_{j=0}^{n-1}a_{ij}x_j\geq l_i^c\)),
the matrix \(Q\) must be \sphinxstylestrong{negative semidefinite}.

\item {} 
\sphinxAtStartPar
As seen from the mathematical formulation above, the quadratic constraint expression includes quadratic terms, linear terms, and constant terms.

\end{itemize}
\end{sphinxadmonition}

\subsection{Modeling}
\label{\detokenize{modeling:id8}}
\sphinxAtStartPar
The basic steps to construct and solve a Quadratically Constrained Programming (QCP) model in COPT are as follows:
\begin{enumerate}
\sphinxsetlistlabels{\arabic}{enumi}{enumii}{}{.}%
\item {} 
\sphinxAtStartPar
Create a COPT environment and model.

\item {} 
\sphinxAtStartPar
Add model parameters.

\item {} 
\sphinxAtStartPar
Construct the QCP model:
\begin{itemize}
\item {} 
\sphinxAtStartPar
Add variables.

\item {} 
\sphinxAtStartPar
Add quadratic constraints.

\item {} 
\sphinxAtStartPar
Set the quadratic objective function.

\end{itemize}

\item {} 
\sphinxAtStartPar
Set solving parameters and solve.

\item {} 
\sphinxAtStartPar
Retrieve the solution.

\end{enumerate}

\sphinxAtStartPar
\sphinxstylestrong{Modeling: Adding Quadratic Constraints}

\sphinxAtStartPar
COPT provides two ways to add quadratic constraints:
\begin{enumerate}
\sphinxsetlistlabels{\arabic}{enumi}{enumii}{}{.}%
\item {} 
\sphinxAtStartPar
First, construct the quadratic constraint expression by combining the quadratic and linear terms, then add it to the model.

\item {} 
\sphinxAtStartPar
Directly provide the quadratic constraint expression (or a quadratic constraint builder) as a parameter input and add it to the model.

\end{enumerate}

\sphinxAtStartPar
When adding quadratic constraints to the model, the following parameters can be specified by the user:
\begin{itemize}
\item {} 
\sphinxAtStartPar
\sphinxcode{\sphinxupquote{expr}} / \sphinxcode{\sphinxupquote{builder}}: The quadratic constraint expression or quadratic constraint builder.

\item {} 
\sphinxAtStartPar
\sphinxcode{\sphinxupquote{rhs}}: The right\sphinxhyphen{}hand side of the quadratic constraint.

\item {} 
\sphinxAtStartPar
\sphinxcode{\sphinxupquote{sense}}: The type of constraint, which can be \sphinxcode{\sphinxupquote{COPT\_LESS\_EQUAL}} or \sphinxcode{\sphinxupquote{COPT\_GREATER\_EQUAL}}.

\item {} 
\sphinxAtStartPar
\sphinxcode{\sphinxupquote{name}}: The name of the quadratic constraint.

\end{itemize}

\sphinxAtStartPar
The implementation in different programming interfaces is shown in \hyperref[\detokenize{modeling:copttab-modelingqmodel-constr}]{Table \ref{\detokenize{modeling:copttab-modelingqmodel-constr}}}:

\begin{savenotes}\sphinxattablestart
\sphinxthistablewithglobalstyle
\centering
\sphinxcapstartof{table}
\sphinxthecaptionisattop
\sphinxcaption{Functions for Adding Quadratic Constraints}\label{\detokenize{modeling:copttab-modelingqmodel-constr}}
\sphinxaftertopcaption
\begin{tabular}[t]{|\X{5}{15}|\X{10}{15}|}
\sphinxtoprule
\sphinxstyletheadfamily 
\sphinxAtStartPar
Programming Interface
&\sphinxstyletheadfamily 
\sphinxAtStartPar
Function
\\
\sphinxmidrule
\sphinxtableatstartofbodyhook
\sphinxAtStartPar
C
&
\sphinxAtStartPar
\sphinxcode{\sphinxupquote{COPT\_AddQConstr}}
\\
\sphinxhline
\sphinxAtStartPar
C++
&
\sphinxAtStartPar
\sphinxcode{\sphinxupquote{Model::AddQConstr()}}
\\
\sphinxhline
\sphinxAtStartPar
C\#
&
\sphinxAtStartPar
\sphinxcode{\sphinxupquote{Model.AddQConstr()}}
\\
\sphinxhline
\sphinxAtStartPar
Java
&
\sphinxAtStartPar
\sphinxcode{\sphinxupquote{Model.addQConstr()}}
\\
\sphinxhline
\sphinxAtStartPar
Python
&
\sphinxAtStartPar
\sphinxcode{\sphinxupquote{Model.addQConstr()}}
\\
\sphinxbottomrule
\end{tabular}
\sphinxtableafterendhook\par
\sphinxattableend\end{savenotes}

\begin{sphinxadmonition}{note}{Note:}
\begin{itemize}
\item {} 
\sphinxAtStartPar
The operations related to quadratic constraint modeling may vary slightly in terms of function names, calling methods, and \sphinxstylestrong{parameter names} in different programming interfaces, but the functionality and parameter meanings are consistent.

\item {} 
\sphinxAtStartPar
In the C API, when adding a quadratic constraint, non\sphinxhyphen{}zero linear term coefficient indices, quadratic term indices, and other parameters must be provided. For more details, please refer to {\hyperref[\detokenize{capiref:chapapi-model}]{\sphinxcrossref{\DUrole{std,std-ref}{C API Functions: Constructing and Modifying the Model}}}} under the function \sphinxcode{\sphinxupquote{COPT\_AddQConstr}}.

\item {} 
\sphinxAtStartPar
In the Python API, \sphinxcode{\sphinxupquote{Model.addQConstr()}} can be used to add a quadratic constraint to the model. For more details, please refer to {\hyperref[\detokenize{pyapiref:chappyapi-model}]{\sphinxcrossref{\DUrole{std,std-ref}{Python API Functions: Model Class}}}}.

\end{itemize}
\end{sphinxadmonition}

\sphinxAtStartPar
\sphinxstylestrong{Examples of quadratic constraints}
\begin{equation}\label{equation:modeling:modeling:20}
\begin{split}x_1^2 + x_2^2 + x_1 + 2x_2 <= 0\end{split}
\end{equation}
\sphinxAtStartPar
The code implementation in different programming interfaces is as follows:

\sphinxAtStartPar
C API:

\begin{sphinxVerbatim}[commandchars=\\\{\}]
\PYG{k+kt}{int}\PYG{+w}{ }\PYG{n}{nRowMatCnt}\PYG{+w}{ }\PYG{o}{=}\PYG{+w}{ }\PYG{l+m+mi}{2}\PYG{p}{;}
\PYG{k+kt}{int}\PYG{+w}{ }\PYG{n}{rowMatIdx}\PYG{p}{[}\PYG{p}{]}\PYG{+w}{ }\PYG{o}{=}\PYG{+w}{ }\PYG{p}{\PYGZob{}}\PYG{l+m+mi}{0}\PYG{p}{,}\PYG{+w}{ }\PYG{l+m+mi}{1}\PYG{p}{\PYGZcb{}}\PYG{p}{;}
\PYG{k+kt}{double}\PYG{+w}{ }\PYG{n}{rowMatElem}\PYG{p}{[}\PYG{p}{]}\PYG{+w}{ }\PYG{o}{=}\PYG{+w}{ }\PYG{p}{\PYGZob{}}\PYG{l+m+mf}{1.0}\PYG{p}{,}\PYG{+w}{ }\PYG{l+m+mf}{2.0}\PYG{p}{\PYGZcb{}}\PYG{p}{;}

\PYG{k+kt}{int}\PYG{+w}{ }\PYG{n}{nQMatCnt}\PYG{+w}{ }\PYG{o}{=}\PYG{+w}{ }\PYG{l+m+mi}{2}\PYG{p}{;}
\PYG{k+kt}{int}\PYG{+w}{ }\PYG{n}{qMatRow}\PYG{p}{[}\PYG{p}{]}\PYG{+w}{ }\PYG{o}{=}\PYG{+w}{ }\PYG{p}{\PYGZob{}}\PYG{l+m+mi}{0}\PYG{p}{,}\PYG{+w}{ }\PYG{l+m+mi}{1}\PYG{p}{\PYGZcb{}}\PYG{p}{;}
\PYG{k+kt}{int}\PYG{+w}{ }\PYG{n}{qMatCol}\PYG{p}{[}\PYG{p}{]}\PYG{+w}{ }\PYG{o}{=}\PYG{+w}{ }\PYG{p}{\PYGZob{}}\PYG{l+m+mi}{0}\PYG{p}{,}\PYG{+w}{ }\PYG{l+m+mi}{1}\PYG{p}{\PYGZcb{}}\PYG{p}{;}
\PYG{k+kt}{double}\PYG{+w}{ }\PYG{n}{qMatElem}\PYG{p}{[}\PYG{p}{]}\PYG{+w}{ }\PYG{o}{=}\PYG{+w}{ }\PYG{p}{\PYGZob{}}\PYG{l+m+mf}{1.0}\PYG{p}{,}\PYG{+w}{ }\PYG{l+m+mf}{1.0}\PYG{p}{\PYGZcb{}}\PYG{p}{;}
\PYG{k+kt}{char}\PYG{+w}{ }\PYG{n}{cRowSense}\PYG{+w}{ }\PYG{o}{=}\PYG{+w}{ }\PYG{n}{COPT\PYGZus{}LESS\PYGZus{}EQUAL}\PYG{p}{;}
\PYG{k+kt}{double}\PYG{+w}{ }\PYG{n}{dRowBound}\PYG{+w}{ }\PYG{o}{=}\PYG{+w}{ }\PYG{l+m+mf}{0.0}\PYG{p}{;}
\PYG{k+kt}{char}\PYG{+w}{ }\PYG{o}{*}\PYG{n}{name}\PYG{+w}{ }\PYG{o}{=}\PYG{+w}{ }\PYG{l+s}{\PYGZdq{}}\PYG{l+s}{q1}\PYG{l+s}{\PYGZdq{}}\PYG{p}{;}
\PYG{n}{errcode}\PYG{+w}{ }\PYG{o}{=}\PYG{+w}{ }\PYG{n}{COPT\PYGZus{}AddQConstr}\PYG{p}{(}\PYG{n}{prob}\PYG{p}{,}\PYG{+w}{ }\PYG{n}{nRowMatCnt}\PYG{p}{,}\PYG{+w}{ }\PYG{n}{rowMatIdx}\PYG{p}{,}\PYG{+w}{ }\PYG{n}{rowMatElem}\PYG{p}{,}
\PYG{+w}{                          }\PYG{n}{nQMatCnt}\PYG{p}{,}\PYG{+w}{ }\PYG{n}{qMatRow}\PYG{p}{,}\PYG{+w}{ }\PYG{n}{qMatCol}\PYG{p}{,}\PYG{+w}{ }\PYG{n}{qMatElem}\PYG{p}{,}
\PYG{+w}{                          }\PYG{n}{cRowSense}\PYG{p}{,}\PYG{+w}{ }\PYG{n}{dRowBound}\PYG{p}{,}\PYG{+w}{ }\PYG{n}{name}\PYG{p}{)}\PYG{p}{;}
\end{sphinxVerbatim}

\sphinxAtStartPar
C++ API:

\begin{sphinxVerbatim}[commandchars=\\\{\}]
\PYG{n}{model}\PYG{p}{.}\PYG{n}{AddQConstr}\PYG{p}{(}\PYG{n}{x1}\PYG{o}{*}\PYG{n}{x1}\PYG{+w}{ }\PYG{o}{+}\PYG{+w}{ }\PYG{n}{x2}\PYG{o}{*}\PYG{n}{x2}\PYG{+w}{ }\PYG{o}{+}\PYG{+w}{ }\PYG{n}{x1}\PYG{+w}{ }\PYG{o}{+}\PYG{+w}{ }\PYG{l+m+mi}{2}\PYG{o}{*}\PYG{n}{x2}\PYG{+w}{ }\PYG{o}{\PYGZlt{}}\PYG{o}{=}\PYG{+w}{ }\PYG{l+m+mi}{0}\PYG{p}{,}\PYG{+w}{ }\PYG{l+s}{\PYGZdq{}}\PYG{l+s}{q1}\PYG{l+s}{\PYGZdq{}}\PYG{p}{)}\PYG{p}{;}
\end{sphinxVerbatim}

\sphinxAtStartPar
C\# API:

\begin{sphinxVerbatim}[commandchars=\\\{\}]
\PYG{n}{model}\PYG{p}{.}\PYG{n}{AddQConstr}\PYG{p}{(}\PYG{n}{x1}\PYG{o}{*}\PYG{n}{x1}\PYG{+w}{ }\PYG{o}{+}\PYG{+w}{ }\PYG{n}{x2}\PYG{o}{*}\PYG{n}{x2}\PYG{+w}{ }\PYG{o}{+}\PYG{+w}{ }\PYG{n}{x1}\PYG{+w}{ }\PYG{o}{+}\PYG{+w}{ }\PYG{l+m+mi}{2}\PYG{o}{*}\PYG{n}{x2}\PYG{+w}{ }\PYG{o}{\PYGZlt{}=}\PYG{+w}{ }\PYG{l+m+mi}{0}\PYG{p}{,}\PYG{+w}{ }\PYG{l+s}{\PYGZdq{}q1\PYGZdq{}}\PYG{p}{)}\PYG{p}{;}
\end{sphinxVerbatim}

\sphinxAtStartPar
Java API:

\begin{sphinxVerbatim}[commandchars=\\\{\}]
\PYG{n}{QuadExpr}\PYG{+w}{ }\PYG{n}{q1}\PYG{+w}{ }\PYG{o}{=}\PYG{+w}{ }\PYG{k}{new}\PYG{+w}{ }\PYG{n}{QuadExpr}\PYG{p}{(}\PYG{l+m+mf}{0.0}\PYG{p}{)}\PYG{p}{;}
\PYG{n}{q1}\PYG{p}{.}\PYG{n+na}{addTerm}\PYG{p}{(}\PYG{n}{x1}\PYG{p}{,}\PYG{+w}{ }\PYG{n}{x1}\PYG{p}{,}\PYG{+w}{ }\PYG{l+m+mi}{1}\PYG{p}{)}\PYG{p}{;}
\PYG{n}{q1}\PYG{p}{.}\PYG{n+na}{addTerm}\PYG{p}{(}\PYG{n}{x2}\PYG{p}{,}\PYG{+w}{ }\PYG{n}{x2}\PYG{p}{,}\PYG{+w}{ }\PYG{l+m+mi}{1}\PYG{p}{)}\PYG{p}{;}
\PYG{n}{q1}\PYG{p}{.}\PYG{n+na}{addTerm}\PYG{p}{(}\PYG{n}{x1}\PYG{p}{,}\PYG{+w}{ }\PYG{l+m+mi}{1}\PYG{p}{)}\PYG{p}{;}
\PYG{n}{q1}\PYG{p}{.}\PYG{n+na}{addTerm}\PYG{p}{(}\PYG{n}{x2}\PYG{p}{,}\PYG{+w}{ }\PYG{l+m+mi}{2}\PYG{p}{)}\PYG{p}{;}
\PYG{n}{model}\PYG{p}{.}\PYG{n+na}{addQConstr}\PYG{p}{(}\PYG{n}{q1}\PYG{p}{,}\PYG{+w}{ }\PYG{n}{copt}\PYG{p}{.}\PYG{n+na}{Consts}\PYG{p}{.}\PYG{n+na}{LESS\PYGZus{}EQUAL}\PYG{p}{,}\PYG{+w}{ }\PYG{l+m+mi}{0}\PYG{p}{,}\PYG{+w}{ }\PYG{l+s}{\PYGZdq{}}\PYG{l+s}{q1}\PYG{l+s}{\PYGZdq{}}\PYG{p}{)}\PYG{p}{;}
\end{sphinxVerbatim}

\sphinxAtStartPar
Python API:

\begin{sphinxVerbatim}[commandchars=\\\{\}]
\PYG{n}{model}\PYG{o}{.}\PYG{n}{addQConstr}\PYG{p}{(}\PYG{n}{x1}\PYG{o}{*}\PYG{n}{x1} \PYG{o}{+} \PYG{n}{x2}\PYG{o}{*}\PYG{n}{x2} \PYG{o}{+} \PYG{n}{x1} \PYG{o}{+} \PYG{l+m+mi}{2}\PYG{o}{*}\PYG{n}{x2} \PYG{o}{\PYGZlt{}}\PYG{o}{=} \PYG{l+m+mi}{0}\PYG{p}{,} \PYG{n}{name}\PYG{o}{=}\PYG{l+s+s2}{\PYGZdq{}}\PYG{l+s+s2}{q1}\PYG{l+s+s2}{\PYGZdq{}}\PYG{p}{)}
\end{sphinxVerbatim}

\sphinxAtStartPar
In the programming interfaces provided by COPT, except for C language, other object\sphinxhyphen{}oriented programming interfaces (C\#, C++, Java, Python)
offer classes related to quadratic constraints:
\begin{itemize}
\item {} 
\sphinxAtStartPar
Encapsulation of operations for constructing quadratic expressions
\begin{enumerate}
\sphinxsetlistlabels{\arabic}{enumi}{enumii}{}{.}%
\item {} 
\sphinxAtStartPar
\sphinxcode{\sphinxupquote{QuadExpr}} Class: Encapsulation of operations related to combining variables when constructing quadratic expressions in the COPT.

\end{enumerate}

\item {} 
\sphinxAtStartPar
Encapsulation of operations related to quadratic constraints
\begin{enumerate}
\sphinxsetlistlabels{\arabic}{enumi}{enumii}{}{.}%
\item {} 
\sphinxAtStartPar
\sphinxcode{\sphinxupquote{QConstraint}} Class: Encapsulation of operations related to quadratic constraints in the COPT.

\item {} 
\sphinxAtStartPar
\sphinxcode{\sphinxupquote{QConstrArray}} Class: Facilitates user operations on a group of \sphinxcode{\sphinxupquote{QConstraint}} class objects.

\end{enumerate}

\item {} 
\sphinxAtStartPar
Encapsulation of quadratic constraint builders
\begin{enumerate}
\sphinxsetlistlabels{\arabic}{enumi}{enumii}{}{.}%
\item {} 
\sphinxAtStartPar
\sphinxcode{\sphinxupquote{QConstrBuilder}} Class: Encapsulation of builders for constructing quadratic constraints in the COPT.

\item {} 
\sphinxAtStartPar
\sphinxcode{\sphinxupquote{QConstrBuilderArray}} Class: Facilitates user operations on a group of \sphinxcode{\sphinxupquote{QConstrBuilder}} class objects.

\end{enumerate}

\item {} 
\sphinxAtStartPar
C++ API: {\hyperref[\detokenize{cppapiref:chapcppapiref-quadexpr}]{\sphinxcrossref{\DUrole{std,std-ref}{QuadExpr Class}}}} , {\hyperref[\detokenize{cppapiref:chapcppapiref-qconstr}]{\sphinxcrossref{\DUrole{std,std-ref}{QConstraint Class}}}} , {\hyperref[\detokenize{cppapiref:chapcppapiref-qconstrarray}]{\sphinxcrossref{\DUrole{std,std-ref}{QConstrArray Class}}}} , {\hyperref[\detokenize{cppapiref:chapcppapiref-qconstrbuilder}]{\sphinxcrossref{\DUrole{std,std-ref}{QConstrBuilder Class}}}} , {\hyperref[\detokenize{cppapiref:chapcppapiref-qconstrbuilderarray}]{\sphinxcrossref{\DUrole{std,std-ref}{QConstrBuilderArray Class}}}}

\item {} 
\sphinxAtStartPar
C\# API: {\hyperref[\detokenize{csharpapiref:chapcsharpapiref-quadexpr}]{\sphinxcrossref{\DUrole{std,std-ref}{QuadExpr Class}}}} , {\hyperref[\detokenize{csharpapiref:chapcsharpapiref-qconstr}]{\sphinxcrossref{\DUrole{std,std-ref}{QConstraint Class}}}} , {\hyperref[\detokenize{csharpapiref:chapcsharpapiref-qconstrarray}]{\sphinxcrossref{\DUrole{std,std-ref}{QConstrArray Class}}}} , {\hyperref[\detokenize{csharpapiref:chapcsharpapiref-qconstrbuilder}]{\sphinxcrossref{\DUrole{std,std-ref}{QConstrBuilder Class}}}} , {\hyperref[\detokenize{csharpapiref:chapcsharpapiref-qconstrbuilderarray}]{\sphinxcrossref{\DUrole{std,std-ref}{QConstrBuilderArray Class}}}}

\item {} 
\sphinxAtStartPar
Java API: {\hyperref[\detokenize{javaapiref:chapjavaapiref-quadexpr}]{\sphinxcrossref{\DUrole{std,std-ref}{QuadExpr Class}}}} , {\hyperref[\detokenize{javaapiref:chapjavaapiref-qconstr}]{\sphinxcrossref{\DUrole{std,std-ref}{QConstraint Class}}}} , {\hyperref[\detokenize{javaapiref:chapjavaapiref-qconstrarray}]{\sphinxcrossref{\DUrole{std,std-ref}{QConstrArray Class}}}} , {\hyperref[\detokenize{javaapiref:chapjavaapiref-qconstrbuilder}]{\sphinxcrossref{\DUrole{std,std-ref}{QConstrBuilder Class}}}} , {\hyperref[\detokenize{javaapiref:chapjavaapiref-qconstrbuilderarray}]{\sphinxcrossref{\DUrole{std,std-ref}{QConstrBuilderArray Class}}}}

\item {} 
\sphinxAtStartPar
Python API: {\hyperref[\detokenize{pyapiref:chappyapi-quadexpr}]{\sphinxcrossref{\DUrole{std,std-ref}{QuadExpr Class}}}} , {\hyperref[\detokenize{pyapiref:chappyapi-qconstraint}]{\sphinxcrossref{\DUrole{std,std-ref}{QConstraint Class}}}} , {\hyperref[\detokenize{pyapiref:chappyapi-qconstrarray}]{\sphinxcrossref{\DUrole{std,std-ref}{QConstrArray Class}}}} , {\hyperref[\detokenize{pyapiref:chappyapi-qconstrbuilder}]{\sphinxcrossref{\DUrole{std,std-ref}{QConstrBuilder Class}}}} , {\hyperref[\detokenize{pyapiref:chappyapi-qconstrbuilderarray}]{\sphinxcrossref{\DUrole{std,std-ref}{QConstrBuilderArray Class}}}}

\end{itemize}

\subsection{Related Attributes}
\label{\detokenize{modeling:id9}}
\sphinxAtStartPar
Attributes for Quadratic Programming (QP) and Quadratically Constrained Programming (QCP) models are shown in \hyperref[\detokenize{modeling:copttab-modelingqcp-attr}]{Table \ref{\detokenize{modeling:copttab-modelingqcp-attr}}}:

\begin{savenotes}\sphinxattablestart
\sphinxthistablewithglobalstyle
\centering
\sphinxcapstartof{table}
\sphinxthecaptionisattop
\sphinxcaption{Attributes for QP and QCP}\label{\detokenize{modeling:copttab-modelingqcp-attr}}
\sphinxaftertopcaption
\begin{tabular}[t]{|\X{15}{54}|\X{9}{54}|\X{30}{54}|}
\sphinxtoprule
\sphinxstyletheadfamily 
\sphinxAtStartPar
Name
&\sphinxstyletheadfamily 
\sphinxAtStartPar
Type
&\sphinxstyletheadfamily 
\sphinxAtStartPar
Description
\\
\sphinxmidrule
\sphinxtableatstartofbodyhook
\sphinxAtStartPar
\sphinxcode{\sphinxupquote{QConstrs}}
&
\sphinxAtStartPar
Integer
&
\sphinxAtStartPar
Number of quadratic constraints
\\
\sphinxhline
\sphinxAtStartPar
\sphinxcode{\sphinxupquote{QElems}}
&
\sphinxAtStartPar
Integer
&
\sphinxAtStartPar
Number of non\sphinxhyphen{}zero quadratic elements in the quadratic objective function
\\
\sphinxhline
\sphinxAtStartPar
\sphinxcode{\sphinxupquote{HasQObj}}
&
\sphinxAtStartPar
Integer
&
\sphinxAtStartPar
Whether the problem has quadratic objective function
\\
\sphinxbottomrule
\end{tabular}
\sphinxtableafterendhook\par
\sphinxattableend\end{savenotes}

\section{General Nonlinear Programming (NLP)}
\label{\detokenize{modeling:general-nonlinear-programming-nlp}}\label{\detokenize{modeling:chapmodeling-nlp}}
\sphinxAtStartPar
Nonlinear Programming (NLP) refers to optimization problems in which the objective
function or constraint functions contain nonlinear expressions.
The Cardinal Optimizer (COPT) provides support for solving a wide range of nonlinear
programming problems and is suitable for applications such as engineering optimization,
numerical computation, and complex system modeling.

\sphinxAtStartPar
Currently, COPT supports nonlinear programming problems through the following approaches:
\begin{itemize}
\item {} 
\sphinxAtStartPar
For nonconvex quadratic programming (Nonconvex QP) and nonconvex quadratically
constrained programming (Nonconvex QCQP), the problems can be solved by setting
the optimization parameter \sphinxcode{\sphinxupquote{"NonConvex"}};

\item {} 
\sphinxAtStartPar
For more general nonlinear programming problems, COPT supports solving the model
either by explicit nonlinear expression modeling or by using a callback\sphinxhyphen{}based
nonlinear modeling approach;

\item {} 
\sphinxAtStartPar
In addition, models in \sphinxcode{\sphinxupquote{NL}} format can be solved by reading the model files
via the \sphinxcode{\sphinxupquote{coptampl}} tool.

\end{itemize}

\sphinxAtStartPar
The following sections introduce the two main nonlinear modeling approaches supported by COPT.

\subsection{Explicit Expression Nonlinear  Modeling}
\label{\detokenize{modeling:explicit-expression-nonlinear-modeling}}
\sphinxAtStartPar
In the explicit expression modeling approach, users construct nonlinear operators
and function expressions to directly describe the mathematical form of the objective
function and constraint functions.
This approach is suitable for nonlinear problems with clear structures that can be
represented in closed\sphinxhyphen{}form analytic expressions.

\sphinxAtStartPar
COPT supports the following types of nonlinear expressions:
\begin{itemize}
\item {} 
\sphinxAtStartPar
Basic arithmetic operations (add, minus, mult, div)

\item {} 
\sphinxAtStartPar
Trigonometric functions (sin, cos, tan, etc.)

\item {} 
\sphinxAtStartPar
Exponential functions (exp), power functions (pow), logarithmic functions (log),
and root functions (sqrt)

\end{itemize}

\sphinxAtStartPar
In addition, COPT also supports commonly used nonlinear functions such as ceiling (ceil),
floor (floor), and absolute value (abs).

\sphinxAtStartPar
For the C interface of COPT, the nonlinear expression operators are described in
{\hyperref[\detokenize{capiref:chapapi-const-nlp}]{\sphinxcrossref{\DUrole{std,std-ref}{Nonlinear Expression Operator Constants}}}};

\sphinxAtStartPar
For object\sphinxhyphen{}oriented interfaces, nonlinear expression operations are encapsulated
in the \sphinxcode{\sphinxupquote{NL}} namespace:
\begin{itemize}
\item {} 
\sphinxAtStartPar
C++ API: {\hyperref[\detokenize{cppapiref:chapcppapiref-nl}]{\sphinxcrossref{\DUrole{std,std-ref}{NL Namespace}}}}

\item {} 
\sphinxAtStartPar
C\# API: {\hyperref[\detokenize{csharpapiref:chapcsharpapiref-nl}]{\sphinxcrossref{\DUrole{std,std-ref}{NL Namespace}}}}

\item {} 
\sphinxAtStartPar
Java API: {\hyperref[\detokenize{javaapiref:chapjavaapiref-nl}]{\sphinxcrossref{\DUrole{std,std-ref}{NL Namespace}}}}

\item {} 
\sphinxAtStartPar
Python API: {\hyperref[\detokenize{pyapiref:chappyapi-nl}]{\sphinxcrossref{\DUrole{std,std-ref}{nl Namespace}}}}

\end{itemize}

\sphinxAtStartPar
The following example illustrates a nonlinear programming model that includes
nonlinear expressions in both the objective function and the constraint functions.

\sphinxAtStartPar
\sphinxstylestrong{Objective Function}
\begin{equation}\label{equation:modeling:modeling:21}
\begin{split}\min \quad x_1 x_4 \left(\sin(x_1 + x_2) + \cos(x_2 x_3) + \tan\left(\frac{x_3}{x_4}\right)\right) + x_3\end{split}
\end{equation}
\sphinxAtStartPar
\sphinxstylestrong{Constraints}
\begin{align}\label{equation:modeling:modeling:22}\!\begin{aligned}
x_1 x_2 x_3 x_4 + x_1 + x_2 \geq 35\\
\log(x_1) + 2 \log(x_2) + 3 \log(x_3) + 4 \log(x_4) + x_3 + x_4 \geq 15\\
x_1^2 + x_2^2 + x_3^2 + x_4^2 + x_1 + x_3 \geq 50\\
\end{aligned}\end{align}
\sphinxAtStartPar
\sphinxstylestrong{Variable Bounds}
\begin{equation}\label{equation:modeling:modeling:23}
\begin{split}1 \leq x_1, x_2, x_3, x_4 \leq 5\end{split}
\end{equation}
\sphinxAtStartPar
The complete implementation of the above model can be found in the example
\sphinxcode{\sphinxupquote{nlp\_ex1}} under the \sphinxcode{\sphinxupquote{examples}} directory of the COPT installation package.

\subsection{Callback\sphinxhyphen{}Based Nonlinear Modeling}
\label{\detokenize{modeling:callback-based-nonlinear-modeling}}
\sphinxAtStartPar
In addition to expression\sphinxhyphen{}based nonlinear modeling, COPT also provides a
callback\sphinxhyphen{}based nonlinear modeling mechanism, which is designed for nonlinear
optimization problems that are difficult to describe using explicit symbolic
expressions.

\sphinxAtStartPar
Under this mode, users do not need to explicitly construct analytical
expressions for the objective function and constraints. Instead, by implementing
a nonlinear callback class, users provide numerical values of the objective
function, constraint functions, and their first\sphinxhyphen{}order and second\sphinxhyphen{}order derivatives
to the solver at given variable points.

\sphinxAtStartPar
The structural information of the nonlinear model (such as the number of
variables, the number of constraints, and the sparsity patterns of derivatives)
is loaded once via the \sphinxcode{\sphinxupquote{Model.loadNlData}} interface. The numerical
computation logic is encapsulated in the callback object and invoked by the
solver during the solution process.
During the iterations of the nonlinear optimization algorithm, COPT
automatically calls the corresponding callback methods at different variable
points to obtain the required numerical information and advance the solution
process.

\sphinxAtStartPar
Currently, the callback\sphinxhyphen{}based nonlinear modeling mechanism is supported in the
Python and C++ interfaces of COPT:
\begin{itemize}
\item {} 
\sphinxAtStartPar
C++ API: {\hyperref[\detokenize{cppapiref:chapcppapiref-nlpcallbackbase}]{\sphinxcrossref{\DUrole{std,std-ref}{NlpCallbackBase class}}}}

\item {} 
\sphinxAtStartPar
Python API: {\hyperref[\detokenize{pyapiref:chappyapi-nlpcbcbase}]{\sphinxcrossref{\DUrole{std,std-ref}{NlpCallbackBase class}}}}

\end{itemize}

\sphinxAtStartPar
Taking object\sphinxhyphen{}oriented programming interfaces as an example, the callback\sphinxhyphen{}based
nonlinear modeling workflow typically consists of the following steps:
\begin{enumerate}
\sphinxsetlistlabels{\arabic}{enumi}{enumii}{}{.}%
\item {} 
\sphinxAtStartPar
Define a custom nonlinear callback class and inherit from
\sphinxcode{\sphinxupquote{NlpCallbackBase}}.

\sphinxAtStartPar
The callback class is responsible for implementing the numerical computation
interfaces of the nonlinear model, including the objective function, the
constraint functions, and their first\sphinxhyphen{} and second\sphinxhyphen{}order derivatives.

\item {} 
\sphinxAtStartPar
Implement the required numerical evaluation methods in the callback class.

\sphinxAtStartPar
The computation rules of the objective function, constraint functions, and
their first\sphinxhyphen{} and second\sphinxhyphen{}order derivatives are provided through the callback
interfaces. Users may implement some or all of the following methods to
define the numerical computation logic of the nonlinear model:
\begin{itemize}
\item {} 
\sphinxAtStartPar
\sphinxcode{\sphinxupquote{NlpCallbackBase.EvalObj}}: compute the objective function value;

\item {} 
\sphinxAtStartPar
\sphinxcode{\sphinxupquote{NlpCallbackBase.EvalCon}}: compute the constraint function values;

\item {} 
\sphinxAtStartPar
\sphinxcode{\sphinxupquote{NlpCallbackBase.EvalGrad}}: compute the first\sphinxhyphen{}order derivative
(gradient) of the objective function;

\item {} 
\sphinxAtStartPar
\sphinxcode{\sphinxupquote{NlpCallbackBase.EvalJac}}: compute the first\sphinxhyphen{}order derivatives
(Jacobian matrix) of the constraint functions;

\item {} 
\sphinxAtStartPar
\sphinxcode{\sphinxupquote{NlpCallbackBase.EvalHess}}: compute the second\sphinxhyphen{}order derivatives
(Hessian matrix) of the Lagrangian function.

\end{itemize}

\sphinxAtStartPar
Each method takes the current variable values as input and provides the
corresponding numerical results through the designated output arrays.

\item {} 
\sphinxAtStartPar
Create an instance of the custom nonlinear callback class.

\sphinxAtStartPar
This instance encapsulates the numerical computation logic of the nonlinear
model and is scheduled uniformly by the solver during the nonlinear solution
process.

\item {} 
\sphinxAtStartPar
Register the nonlinear model information and the callback object by calling
the \sphinxcode{\sphinxupquote{Model.loadNlData}} interface.

\sphinxAtStartPar
In this step, users provide the solver with the model scale, variable and
constraint bounds, and sparsity structures of derivatives, and pass the
custom callback instance as an argument.

\sphinxAtStartPar
Meanwhile, the \sphinxcode{\sphinxupquote{evalType}} parameter can be used to declare the types of
numerical information that the callback object is able to provide, allowing
the solver to invoke the corresponding callback methods as required during
the iterations of the nonlinear optimization algorithm.

\end{enumerate}

\sphinxAtStartPar
After completing the above steps, the nonlinear optimization process can be
started by calling \sphinxcode{\sphinxupquote{Model.solve}}. During the solution process, the solver
automatically invokes the relevant callback methods at different iteration
points to obtain the required numerical information and advance the nonlinear
optimization algorithm.

\section{Mixed\sphinxhyphen{}Integer Programming (MIP)}
\label{\detokenize{modeling:mixed-integer-programming-mip}}\label{\detokenize{modeling:chapmodeling-mip}}

\subsection{Modeling}
\label{\detokenize{modeling:id10}}
\sphinxAtStartPar
Mixed\sphinxhyphen{}Integer Programming (MIP) refers to optimization problems where some of the decision variables are restricted to integer values.
Currently, COPT supports integer variables combined with linear programming, second\sphinxhyphen{}order cone programming, quadratic programming, and quadratically constrained programming.
\begin{quote}

\begin{savenotes}\sphinxattablestart
\sphinxthistablewithglobalstyle
\centering
\begin{tabulary}{\linewidth}[t]{|T|T|}
\sphinxtoprule
\sphinxtableatstartofbodyhook
\sphinxAtStartPar
\sphinxstylestrong{MIP Problem Type}
&
\sphinxAtStartPar
\sphinxstylestrong{Solving Algorithm}
\\
\sphinxhline
\sphinxAtStartPar
MILP
&
\sphinxAtStartPar
Branch\sphinxhyphen{}and\sphinxhyphen{}Cut
\\
\sphinxhline
\sphinxAtStartPar
MISOCP
&
\sphinxAtStartPar
Branch\sphinxhyphen{}and\sphinxhyphen{}Cut
\\
\sphinxhline
\sphinxAtStartPar
MIQP
&
\sphinxAtStartPar
Branch\sphinxhyphen{}and\sphinxhyphen{}Cut
\\
\sphinxhline
\sphinxAtStartPar
MIQCP
&
\sphinxAtStartPar
Branch\sphinxhyphen{}and\sphinxhyphen{}Cut
\\
\sphinxbottomrule
\end{tabulary}
\sphinxtableafterendhook\par
\sphinxattableend\end{savenotes}
\end{quote}

\sphinxAtStartPar
In COPT, the supported integer variable types and their corresponding constants are as follows:
\begin{itemize}
\item {} 
\sphinxAtStartPar
\sphinxcode{\sphinxupquote{BINARY}}
\begin{quote}

\sphinxAtStartPar
Binary variable
\end{quote}

\item {} 
\sphinxAtStartPar
\sphinxcode{\sphinxupquote{INTEGER}}
\begin{quote}

\sphinxAtStartPar
Integer variable
\end{quote}

\end{itemize}

\sphinxAtStartPar
When adding decision variables to the model:
\begin{itemize}
\item {} 
\sphinxAtStartPar
Specify the parameter \sphinxcode{\sphinxupquote{vtype}} as \sphinxcode{\sphinxupquote{BINARY}} to add binary variables;

\item {} 
\sphinxAtStartPar
Specify the parameter \sphinxcode{\sphinxupquote{vtype}} as \sphinxcode{\sphinxupquote{INTEGER}} to add integer variables.

\end{itemize}

\subsection{Solving}
\label{\detokenize{modeling:id11}}
\sphinxAtStartPar
For integer programming problems, COPT provides the Branch\sphinxhyphen{}and\sphinxhyphen{}Cut algorithm, which can be specified by setting the optimization parameter \sphinxcode{\sphinxupquote{"MipMethod"}}.
By configuring other related optimization parameters, you can control the specific workflow of the Branch\sphinxhyphen{}and\sphinxhyphen{}Cut algorithm. For more details, please refer to {\hyperref[\detokenize{parameter:chapparam-mip}]{\sphinxcrossref{\DUrole{std,std-ref}{Parameter Section: Integer Programming Related}}}}.

\sphinxAtStartPar
For the solving logs of integer programming problems, please refer to {\hyperref[\detokenize{logging:chaplogging-branch}]{\sphinxcrossref{\DUrole{std,std-ref}{Logging Section: Branch\sphinxhyphen{}and\sphinxhyphen{}Cut}}}}.

\subsection{Related Attributes}
\label{\detokenize{modeling:id12}}

\begin{savenotes}\sphinxattablestart
\sphinxthistablewithglobalstyle
\centering
\sphinxcapstartof{table}
\sphinxthecaptionisattop
\sphinxcaption{Overview of attributes for MIP}\label{\detokenize{modeling:id13}}
\sphinxaftertopcaption
\begin{tabular}[t]{|\X{15}{54}|\X{9}{54}|\X{30}{54}|}
\sphinxtoprule
\sphinxstyletheadfamily 
\sphinxAtStartPar
Name
&\sphinxstyletheadfamily 
\sphinxAtStartPar
Type
&\sphinxstyletheadfamily 
\sphinxAtStartPar
Description
\\
\sphinxmidrule
\sphinxtableatstartofbodyhook
\sphinxAtStartPar
\sphinxcode{\sphinxupquote{Bins}}
&
\sphinxAtStartPar
Integer
&
\sphinxAtStartPar
Number of binary variables
\\
\sphinxhline
\sphinxAtStartPar
\sphinxcode{\sphinxupquote{Ints}}
&
\sphinxAtStartPar
Integer
&
\sphinxAtStartPar
Number of integer variables
\\
\sphinxhline
\sphinxAtStartPar
\sphinxcode{\sphinxupquote{Indicators}}
&
\sphinxAtStartPar
Integer
&
\sphinxAtStartPar
Number of indicator constraints
\\
\sphinxhline
\sphinxAtStartPar
\sphinxcode{\sphinxupquote{IsMIP}}
&
\sphinxAtStartPar
Integer
&
\sphinxAtStartPar
Whether the problem is a MIP
\\
\sphinxhline
\sphinxAtStartPar
\sphinxcode{\sphinxupquote{NodeCnt}}
&
\sphinxAtStartPar
Integer
&
\sphinxAtStartPar
Number of explored nodes
\\
\sphinxhline
\sphinxAtStartPar
\sphinxcode{\sphinxupquote{HasMipSol}}
&
\sphinxAtStartPar
Integer
&
\sphinxAtStartPar
Whether MIP solution is available
\\
\sphinxhline
\sphinxAtStartPar
\sphinxcode{\sphinxupquote{BestObj}}
&
\sphinxAtStartPar
Double
&
\sphinxAtStartPar
Best integer objective value for MIP
\\
\sphinxhline
\sphinxAtStartPar
\sphinxcode{\sphinxupquote{BestBnd}}
&
\sphinxAtStartPar
Double
&
\sphinxAtStartPar
Best bound for MIP
\\
\sphinxhline
\sphinxAtStartPar
\sphinxcode{\sphinxupquote{BestGap}}
&
\sphinxAtStartPar
Double
&
\sphinxAtStartPar
Best relative gap for MIP
\\
\sphinxbottomrule
\end{tabular}
\sphinxtableafterendhook\par
\sphinxattableend\end{savenotes}

\section{Special Constraints}
\label{\detokenize{modeling:special-constraints}}\label{\detokenize{modeling:chapmodeling-spec}}
\sphinxAtStartPar
COPT supports the construction of two types of special constraints: SOS constraints and Indicator constraints.

\subsection{SOS Constraints}
\label{\detokenize{modeling:sos-constraints}}\label{\detokenize{modeling:chapmodeling-sos}}
\sphinxAtStartPar
SOS constraints (Special Ordered Set) are a special type of constraint that restricts the values of a group of variables.
Currently, COPT supports two types of SOS constraints:
\begin{enumerate}
\sphinxsetlistlabels{\arabic}{enumi}{enumii}{}{.}%
\item {} 
\sphinxAtStartPar
\sphinxstylestrong{SOS1 Constraint}: In this type of constraint, at most one variable in the specified group can take a non\sphinxhyphen{}zero value.

\item {} 
\sphinxAtStartPar
\sphinxstylestrong{SOS2 Constraint}: In this type of constraint, at most two variables in the specified group can take non\sphinxhyphen{}zero values,
and the variables with non\sphinxhyphen{}zero values must be adjacent in the order.

\end{enumerate}

\sphinxAtStartPar
These two types of SOS constraints correspond to the following constants in COPT, which can be specified when adding SOS constraints to the model:
\begin{itemize}
\item {} 
\sphinxAtStartPar
\sphinxcode{\sphinxupquote{SOS\_TYPE1}}
\begin{quote}

\sphinxAtStartPar
SOS1 Constraint
\end{quote}

\item {} 
\sphinxAtStartPar
\sphinxcode{\sphinxupquote{SOS\_TYPE2}}
\begin{quote}

\sphinxAtStartPar
SOS2 Constraint
\end{quote}

\end{itemize}

\sphinxAtStartPar
When adding SOS constraints to the model, the following arguments can be specified by the user:
\begin{itemize}
\item {} 
\sphinxAtStartPar
\sphinxcode{\sphinxupquote{sostype}}: Specifies the type of the SOS constraint.

\item {} 
\sphinxAtStartPar
\sphinxcode{\sphinxupquote{vars}}: The variables involved in the SOS constraint.

\item {} 
\sphinxAtStartPar
\sphinxcode{\sphinxupquote{weights}}: The weights of the variables involved in the SOS constraint; an optional argument, default is \sphinxcode{\sphinxupquote{None}}.

\end{itemize}

\begin{sphinxadmonition}{note}{Note}
\begin{itemize}
\item {} 
\sphinxAtStartPar
The variables involved in the SOS constraint can be continuous variables, binary variables, or integer variables.

\item {} 
\sphinxAtStartPar
If the model includes SOS constraints, the model is an integer programming model.

\item {} 
\sphinxAtStartPar
The specific operations and usage of SOS constraints, including the function names, calling methods, and \sphinxstylestrong{argument names,
may vary slightly} in different programming interfaces, but the functionality and argument meanings are consistent.

\end{itemize}
\end{sphinxadmonition}

\sphinxAtStartPar
COPT provides related functions to support operations on SOS constraints. Below are the corresponding functions for adding
and retrieving SOS constraints in different programming interfaces:

\begin{savenotes}\sphinxattablestart
\sphinxthistablewithglobalstyle
\centering
\sphinxcapstartof{table}
\sphinxthecaptionisattop
\sphinxcaption{Adding and Retrieving SOS Constraints}\label{\detokenize{modeling:copttab-addgetsos}}
\sphinxaftertopcaption
\begin{tabular}[t]{|\X{10}{60}|\X{25}{60}|\X{25}{60}|}
\sphinxtoprule
\sphinxstyletheadfamily 
\sphinxAtStartPar
API
&\sphinxstyletheadfamily 
\sphinxAtStartPar
Add SOS Constraints
&\sphinxstyletheadfamily 
\sphinxAtStartPar
Retrieve All SOS Constraints in the Model
\\
\sphinxmidrule
\sphinxtableatstartofbodyhook
\sphinxAtStartPar
C
&
\sphinxAtStartPar
\sphinxcode{\sphinxupquote{COPT\_AddSOSs}}
&
\sphinxAtStartPar
\sphinxcode{\sphinxupquote{COPT\_GetSOSs}}
\\
\sphinxhline
\sphinxAtStartPar
C++
&
\sphinxAtStartPar
\sphinxcode{\sphinxupquote{Model::AddSos()}}
&
\sphinxAtStartPar
\sphinxcode{\sphinxupquote{Model::GetSoss()}}
\\
\sphinxhline
\sphinxAtStartPar
C\#
&
\sphinxAtStartPar
\sphinxcode{\sphinxupquote{Model.AddSos()}}
&
\sphinxAtStartPar
\sphinxcode{\sphinxupquote{Model.GetSoss()}}
\\
\sphinxhline
\sphinxAtStartPar
Java
&
\sphinxAtStartPar
\sphinxcode{\sphinxupquote{Model.addSos()}}
&
\sphinxAtStartPar
\sphinxcode{\sphinxupquote{Model.getSoss()}}
\\
\sphinxhline
\sphinxAtStartPar
Python
&
\sphinxAtStartPar
\sphinxcode{\sphinxupquote{Model.addSOS()}}
&
\sphinxAtStartPar
\sphinxcode{\sphinxupquote{Model.getSOSs()}}
\\
\sphinxbottomrule
\end{tabular}
\sphinxtableafterendhook\par
\sphinxattableend\end{savenotes}

\sphinxAtStartPar
For operations related to SOS constraints, the function names, calling methods, and argument names may vary slightly in different programming interfaces,
but the functionality and argument meanings are consistent. Please refer to the corresponding sections of each programming language’s API reference manual
for specific details:
\begin{itemize}
\item {} 
\sphinxAtStartPar
C API: {\hyperref[\detokenize{capiref:chapapi-model}]{\sphinxcrossref{\DUrole{std,std-ref}{Building and modifying a problem}}}}

\item {} 
\sphinxAtStartPar
C++ API: {\hyperref[\detokenize{cppapiref:chapcppapiref-model}]{\sphinxcrossref{\DUrole{std,std-ref}{Model Class}}}}

\item {} 
\sphinxAtStartPar
C\# API: {\hyperref[\detokenize{csharpapiref:chapcsharpapiref-model}]{\sphinxcrossref{\DUrole{std,std-ref}{Model Class}}}}

\item {} 
\sphinxAtStartPar
Java API: {\hyperref[\detokenize{javaapiref:chapjavaapiref-model}]{\sphinxcrossref{\DUrole{std,std-ref}{Model Class}}}}

\item {} 
\sphinxAtStartPar
Python API: {\hyperref[\detokenize{pyapiref:chappyapi-model}]{\sphinxcrossref{\DUrole{std,std-ref}{Model Class}}}}

\end{itemize}

\sphinxAtStartPar
In the programming interfaces provided by COPT, except for the C language, other object\sphinxhyphen{}oriented programming interfaces
(C\#, C++, Java, Python) offer classes related to SOS constraints:
\begin{itemize}
\item {} 
\sphinxAtStartPar
Encapsulation of operations related to SOS constraints:
\begin{enumerate}
\sphinxsetlistlabels{\arabic}{enumi}{enumii}{}{.}%
\item {} 
\sphinxAtStartPar
\sphinxcode{\sphinxupquote{Sos}} Class: Encapsulation of operations related to SOS constraints in COPT.

\item {} 
\sphinxAtStartPar
\sphinxcode{\sphinxupquote{SosArray}} Class: Facilitates user operations on a group of \sphinxcode{\sphinxupquote{Sos}} class objects.

\end{enumerate}

\item {} 
\sphinxAtStartPar
Encapsulation of SOS constraint builders:
\begin{enumerate}
\sphinxsetlistlabels{\arabic}{enumi}{enumii}{}{.}%
\item {} 
\sphinxAtStartPar
\sphinxcode{\sphinxupquote{SosBuilder}} Class: Encapsulation of builders for constructing SOS constraints in COPT, providing the following member methods:

\item {} 
\sphinxAtStartPar
\sphinxcode{\sphinxupquote{SosBuilderArray}} Class: Facilitates user operations on a group of \sphinxcode{\sphinxupquote{SosBuilder}} class objects.

\end{enumerate}

\end{itemize}

\sphinxAtStartPar
For more details on the member methods and specific descriptions of the above SOS constraint\sphinxhyphen{}related classes, please refer to the corresponding sections in the API reference manuals of each programming language.
\begin{itemize}
\item {} 
\sphinxAtStartPar
C++ API: {\hyperref[\detokenize{cppapiref:chapcppapiref-sos}]{\sphinxcrossref{\DUrole{std,std-ref}{Sos Class}}}} , {\hyperref[\detokenize{cppapiref:chapcppapiref-sosarray}]{\sphinxcrossref{\DUrole{std,std-ref}{SosArray Class}}}} , {\hyperref[\detokenize{cppapiref:chapcppapiref-sosbuilder}]{\sphinxcrossref{\DUrole{std,std-ref}{SosBuilder Class}}}} , {\hyperref[\detokenize{cppapiref:chapcppapiref-sosbuilderarray}]{\sphinxcrossref{\DUrole{std,std-ref}{SosBuilderArray Class}}}}

\item {} 
\sphinxAtStartPar
C\# API: {\hyperref[\detokenize{csharpapiref:chapcsharpapiref-sos}]{\sphinxcrossref{\DUrole{std,std-ref}{Sos Class}}}} , {\hyperref[\detokenize{csharpapiref:chapcsharpapiref-sosarray}]{\sphinxcrossref{\DUrole{std,std-ref}{SosArray Class}}}} , {\hyperref[\detokenize{csharpapiref:chapcsharpapiref-sosbuilder}]{\sphinxcrossref{\DUrole{std,std-ref}{SosBuilder Class}}}} , {\hyperref[\detokenize{csharpapiref:chapcsharpapiref-sosbuilderarray}]{\sphinxcrossref{\DUrole{std,std-ref}{SosBuilderArray Class}}}}

\item {} 
\sphinxAtStartPar
Java API: {\hyperref[\detokenize{javaapiref:chapjavaapiref-sos}]{\sphinxcrossref{\DUrole{std,std-ref}{Sos Class}}}} , {\hyperref[\detokenize{javaapiref:chapjavaapiref-sosarray}]{\sphinxcrossref{\DUrole{std,std-ref}{SosArray Class}}}} , {\hyperref[\detokenize{javaapiref:chapjavaapiref-sosbuilder}]{\sphinxcrossref{\DUrole{std,std-ref}{SosBuilder Class}}}} , {\hyperref[\detokenize{javaapiref:chapjavaapiref-sosbuilderarray}]{\sphinxcrossref{\DUrole{std,std-ref}{SosBuilderArray Class}}}}

\item {} 
\sphinxAtStartPar
Python API: {\hyperref[\detokenize{pyapiref:chappyapi-sos}]{\sphinxcrossref{\DUrole{std,std-ref}{SOS Class}}}} , {\hyperref[\detokenize{pyapiref:chappyapi-sosarray}]{\sphinxcrossref{\DUrole{std,std-ref}{SOSArray Class}}}} , {\hyperref[\detokenize{pyapiref:chappyapi-sosbuilder}]{\sphinxcrossref{\DUrole{std,std-ref}{SOSBuilder Class}}}} , {\hyperref[\detokenize{pyapiref:chappyapi-sosbuilderarray}]{\sphinxcrossref{\DUrole{std,std-ref}{SOSBuilderArray Class}}}}

\end{itemize}

\subsection{Indicator Constraints}
\label{\detokenize{modeling:indicator-constraints}}\label{\detokenize{modeling:chapmodeling-indicator}}
\sphinxAtStartPar
Indicator constraint is a type of logical constraint that uses a binary variable \(y\) as the indicator variable to
determine the logical relationship between the value of \(y\) and whether the linear constraint \(a^{T}x \leq b\)
is satisfied. Currently, COPT supports three types of Indicator constraints: If\sphinxhyphen{}Then, Only\sphinxhyphen{}If, and If\sphinxhyphen{}and\sphinxhyphen{}Only\sphinxhyphen{}If.
\begin{itemize}
\item {} 
\sphinxAtStartPar
\sphinxcode{\sphinxupquote{INDICATOR\_IF}}

\sphinxAtStartPar
If\sphinxhyphen{}Then:
\begin{quote}

\sphinxAtStartPar
If \(y = f\), then the linear constraint is satisfied;

\sphinxAtStartPar
If \(y \ne f\), then the linear constraint can be violated.
\end{quote}

\end{itemize}
\begin{equation}\label{equation:modeling:modeling:24}
\begin{split}y &= f \rightarrow a^{T}x \leq b\\
f &\in \{0, 1\}\end{split}
\end{equation}\begin{itemize}
\item {} 
\sphinxAtStartPar
\sphinxcode{\sphinxupquote{INDICATOR\_ONLYIF}}

\sphinxAtStartPar
Only\sphinxhyphen{}If:
\begin{quote}

\sphinxAtStartPar
If the linear constraint \(a^{T}x \leq b\) is satisfied, then \(y = f\);

\sphinxAtStartPar
If the linear constraint \(a^{T}x \leq b\) is not satisfied, then \(y\) can take the value of 0 or 1.
\end{quote}

\end{itemize}
\begin{equation}\label{equation:modeling:modeling:25}
\begin{split}a^{T}x &\leq b \rightarrow y = f\\
f &\in \{0, 1\}\end{split}
\end{equation}\begin{itemize}
\item {} 
\sphinxAtStartPar
\sphinxcode{\sphinxupquote{INDICATOR\_IFANDONLYIF}}

\sphinxAtStartPar
If\sphinxhyphen{}and\sphinxhyphen{}Only\sphinxhyphen{}If:
\begin{quote}

\sphinxAtStartPar
The linear constraint \(a^{T}x \leq b\) and \(y = f\) are equivalent.
They must both be satisfied or both be unsatisfied.
\end{quote}

\end{itemize}
\begin{equation}\label{equation:modeling:modeling:26}
\begin{split}a^{T}x &\leq b \leftrightarrow y = f\\
f &\in \{0, 1\}\end{split}
\end{equation}
\sphinxAtStartPar
COPT provides two methods to add Indicator constraints:
\begin{enumerate}
\sphinxsetlistlabels{\arabic}{enumi}{enumii}{}{.}%
\item {} 
\sphinxAtStartPar
By calling an API function (e.g., in Python: \sphinxcode{\sphinxupquote{Model.addGenConstrIndicator()}}). The key parameters for constructing an Indicator constraint are:
\begin{itemize}
\item {} 
\sphinxAtStartPar
\sphinxcode{\sphinxupquote{binVar}}: The binary indicator variable.

\item {} 
\sphinxAtStartPar
\sphinxcode{\sphinxupquote{binval}}: The value (\sphinxcode{\sphinxupquote{True}} or \sphinxcode{\sphinxupquote{False}}) of the indicator variable that is conditionally related to the satisfaction of the linear constraint.

\item {} 
\sphinxAtStartPar
\sphinxcode{\sphinxupquote{builder}}: The linear constraint builder.

\item {} 
\sphinxAtStartPar
\sphinxcode{\sphinxupquote{type}}: The type of Indicator constraint (possible values are listed in {\hyperref[\detokenize{constant:chapconst-indicatortype}]{\sphinxcrossref{\DUrole{std,std-ref}{Indicator constraint types}}}}).

\end{itemize}

\item {} 
\sphinxAtStartPar
By overloading operators (for If\sphinxhyphen{}Then and Only\sphinxhyphen{}If constraints):
\begin{itemize}
\item {} 
\sphinxAtStartPar
\sphinxcode{\sphinxupquote{\textgreater{}\textgreater{}}}: Represents the If\sphinxhyphen{}Then logical relationship, corresponding to \sphinxcode{\sphinxupquote{INDICATOR\_IF}};

\item {} 
\sphinxAtStartPar
\sphinxcode{\sphinxupquote{\textless{}\textless{}}}: Represents the Only\sphinxhyphen{}If logical relationship, corresponding to \sphinxcode{\sphinxupquote{INDICATOR\_ONLYIF}}.

\end{itemize}

\end{enumerate}

\sphinxAtStartPar
Here are examples of how to implement these methods in Python API:
\begin{enumerate}
\sphinxsetlistlabels{\arabic}{enumi}{enumii}{}{.}%
\item {} 
\sphinxAtStartPar
Adding an If\sphinxhyphen{}Then type Indicator constraint where the linear constraint \(y + 2z >= 3\) is satisfied if \(x\) is true.
\begin{quote}
\begin{equation}\label{equation:modeling:modeling:27}
\begin{split}x = 1 \rightarrow y+2z \geq 3\end{split}
\end{equation}
\begin{sphinxVerbatim}[commandchars=\\\{\}]
\PYG{n}{model}\PYG{o}{.}\PYG{n}{addGenConstrIndicator}\PYG{p}{(}\PYG{n}{x}\PYG{p}{,} \PYG{k+kc}{True}\PYG{p}{,} \PYG{n}{y} \PYG{o}{+} \PYG{l+m+mi}{2}\PYG{o}{*}\PYG{n}{z} \PYG{o}{\PYGZgt{}}\PYG{o}{=} \PYG{l+m+mi}{3}\PYG{p}{)}
\end{sphinxVerbatim}

\begin{sphinxVerbatim}[commandchars=\\\{\}]
\PYG{n}{model}\PYG{o}{.}\PYG{n}{addConstr}\PYG{p}{(}\PYG{p}{(}\PYG{n}{x}\PYG{o}{==}\PYG{l+m+mi}{1}\PYG{p}{)} \PYG{o}{\PYGZgt{}\PYGZgt{}} \PYG{p}{(}\PYG{n}{y} \PYG{o}{+} \PYG{l+m+mi}{2}\PYG{o}{*}\PYG{n}{z} \PYG{o}{\PYGZgt{}}\PYG{o}{=} \PYG{l+m+mi}{3}\PYG{p}{)}\PYG{p}{)}
\end{sphinxVerbatim}
\end{quote}

\item {} 
\sphinxAtStartPar
Adding an Only\sphinxhyphen{}If type Indicator constraint where \(x\) is false if the linear constraint \(y + 2z <= 3\) is satisfied.
\begin{quote}
\begin{equation}\label{equation:modeling:modeling:28}
\begin{split}x = 0 \leftarrow y+2z \leq 3\end{split}
\end{equation}
\begin{sphinxVerbatim}[commandchars=\\\{\}]
\PYG{n}{model}\PYG{o}{.}\PYG{n}{addGenConstrIndicator}\PYG{p}{(}\PYG{n}{x}\PYG{p}{,} \PYG{k+kc}{False}\PYG{p}{,} \PYG{n}{y} \PYG{o}{+} \PYG{l+m+mi}{2}\PYG{o}{*}\PYG{n}{z} \PYG{o}{\PYGZlt{}}\PYG{o}{=} \PYG{l+m+mi}{3}\PYG{p}{,} \PYG{n+nb}{type}\PYG{o}{=}\PYG{n}{COPT}\PYG{o}{.}\PYG{n}{INDICATOR\PYGZus{}ONLYIF}\PYG{p}{)}
\end{sphinxVerbatim}

\begin{sphinxVerbatim}[commandchars=\\\{\}]
\PYG{n}{model}\PYG{o}{.}\PYG{n}{addConstr}\PYG{p}{(}\PYG{p}{(}\PYG{n}{x}\PYG{o}{==}\PYG{l+m+mi}{0}\PYG{p}{)} \PYG{o}{\PYGZlt{}\PYGZlt{}} \PYG{p}{(}\PYG{n}{y} \PYG{o}{+} \PYG{l+m+mi}{2}\PYG{o}{*}\PYG{n}{z} \PYG{o}{\PYGZlt{}}\PYG{o}{=} \PYG{l+m+mi}{3}\PYG{p}{)}\PYG{p}{)}
\end{sphinxVerbatim}
\end{quote}

\item {} 
\sphinxAtStartPar
Adding an If\sphinxhyphen{}and\sphinxhyphen{}Only\sphinxhyphen{}If type Indicator constraint where the binary variable \(x\) being true is equivalent to the
linear constraint \(y + 2z = 3\) being satisfied.
\begin{quote}
\begin{equation}\label{equation:modeling:modeling:29}
\begin{split}x = 1 \leftrightarrow y+2z = 3\end{split}
\end{equation}
\begin{sphinxVerbatim}[commandchars=\\\{\}]
\PYG{n}{model}\PYG{o}{.}\PYG{n}{addGenConstrIndicator}\PYG{p}{(}\PYG{n}{x}\PYG{p}{,} \PYG{k+kc}{True}\PYG{p}{,} \PYG{n}{y} \PYG{o}{+} \PYG{l+m+mi}{2}\PYG{o}{*}\PYG{n}{z} \PYG{o}{==} \PYG{l+m+mi}{3}\PYG{p}{,} \PYG{n+nb}{type}\PYG{o}{=}\PYG{n}{COPT}\PYG{o}{.}\PYG{n}{INDICATOR\PYGZus{}IFANDONLYIF}\PYG{p}{)}
\end{sphinxVerbatim}
\end{quote}

\end{enumerate}

\begin{sphinxadmonition}{note}{Note}
\begin{enumerate}
\sphinxsetlistlabels{\arabic}{enumi}{enumii}{}{.}%
\item {} 
\sphinxAtStartPar
The general expression of the linear constraint given above, \(a^{T}x \leq b\), can actually take the forms \(\leq\), \(\geq\), or \(=\).

\item {} 
\sphinxAtStartPar
If the model includes Indicator constraints, it is considered an integer programming model.

\item {} 
\sphinxAtStartPar
COPT supports adding a batch of Indicator constraints to the model by calling an API function. In Python, the corresponding function is: \sphinxcode{\sphinxupquote{Model.addGenConstrIndicators()}}.

\item {} 
\sphinxAtStartPar
The method of adding Indicator constraints by overloading operators only supports If\sphinxhyphen{}Then and Only\sphinxhyphen{}If constraints. If you need to add an If\sphinxhyphen{}and\sphinxhyphen{}Only\sphinxhyphen{}If constraint, you must use the API function and specify the \sphinxcode{\sphinxupquote{type}} as \sphinxcode{\sphinxupquote{INDICATOR\_IFANDONLYIF}}.

\item {} 
\sphinxAtStartPar
The specific operations and usage of Indicator constraints, including function names, calling methods, and \sphinxstylestrong{parameter names, may vary slightly} in different programming interfaces, but the functionality and parameter meanings are consistent.

\end{enumerate}
\end{sphinxadmonition}

\sphinxAtStartPar
COPT provides related functions to support adding Indicator constraints and retrieving the corresponding Indicator constraint builders, as listed below:

\begin{savenotes}\sphinxattablestart
\sphinxthistablewithglobalstyle
\centering
\sphinxcapstartof{table}
\sphinxthecaptionisattop
\sphinxcaption{Adding and Retrieving Indicator Constraints in Different APIs}\label{\detokenize{modeling:copttab-addgetindicator}}
\sphinxaftertopcaption
\begin{tabular}[t]{|\X{10}{60}|\X{25}{60}|\X{25}{60}|}
\sphinxtoprule
\sphinxstyletheadfamily 
\sphinxAtStartPar
API
&\sphinxstyletheadfamily 
\sphinxAtStartPar
Add Indicator
&\sphinxstyletheadfamily 
\sphinxAtStartPar
Retrieve Indicator
\\
\sphinxmidrule
\sphinxtableatstartofbodyhook
\sphinxAtStartPar
C
&
\sphinxAtStartPar
\sphinxcode{\sphinxupquote{COPT\_AddIndicator}}
&
\sphinxAtStartPar
\sphinxcode{\sphinxupquote{COPT\_GetIndicator}}
\\
\sphinxhline
\sphinxAtStartPar
C++
&
\sphinxAtStartPar
\sphinxcode{\sphinxupquote{Model::AddGenConstrIndicator()}}
&
\sphinxAtStartPar
\sphinxcode{\sphinxupquote{Model::GetGenConstrIndicator()}}
\\
\sphinxhline
\sphinxAtStartPar
C\#
&
\sphinxAtStartPar
\sphinxcode{\sphinxupquote{Model.AddSos()}}
&
\sphinxAtStartPar
\sphinxcode{\sphinxupquote{Model.GetGenConstrIndicator()}}
\\
\sphinxhline
\sphinxAtStartPar
Java
&
\sphinxAtStartPar
\sphinxcode{\sphinxupquote{Model.addSos()}}
&
\sphinxAtStartPar
\sphinxcode{\sphinxupquote{Model.getGenConstrIndicator()}}
\\
\sphinxhline
\sphinxAtStartPar
Python
&
\sphinxAtStartPar
\sphinxcode{\sphinxupquote{Model.addGenConstrIndicator()}}
&
\sphinxAtStartPar
\sphinxcode{\sphinxupquote{Model.getGenConstrIndicator()}}
\\
\sphinxbottomrule
\end{tabular}
\sphinxtableafterendhook\par
\sphinxattableend\end{savenotes}

\sphinxAtStartPar
For operations such as building and adding Indicator constraints, the function names, calling methods, and argument names may vary slightly in different programming interfaces, but the functionality and argument meanings are consistent. Please refer to the corresponding sections of each programming language’s API reference manual for specific details:
\begin{itemize}
\item {} 
\sphinxAtStartPar
C API: {\hyperref[\detokenize{capiref:chapapi-model}]{\sphinxcrossref{\DUrole{std,std-ref}{Constructing and Modifying the Model}}}}

\item {} 
\sphinxAtStartPar
C++ API: {\hyperref[\detokenize{cppapiref:chapcppapiref-model}]{\sphinxcrossref{\DUrole{std,std-ref}{Model Class}}}}

\item {} 
\sphinxAtStartPar
C\# API: {\hyperref[\detokenize{csharpapiref:chapcsharpapiref-model}]{\sphinxcrossref{\DUrole{std,std-ref}{Model Class}}}}

\item {} 
\sphinxAtStartPar
Java API: {\hyperref[\detokenize{javaapiref:chapjavaapiref-model}]{\sphinxcrossref{\DUrole{std,std-ref}{Model Class}}}}

\item {} 
\sphinxAtStartPar
Python API: {\hyperref[\detokenize{pyapiref:chappyapi-model}]{\sphinxcrossref{\DUrole{std,std-ref}{Model Class}}}}

\end{itemize}

\sphinxAtStartPar
In the programming interfaces supported by COPT, except for the C language, other object\sphinxhyphen{}oriented programming interfaces (C\#, C++, Java, Python) offer classes related to Indicator constraints:
\begin{itemize}
\item {} 
\sphinxAtStartPar
Encapsulation of operations related to Indicator constraints:
\begin{enumerate}
\sphinxsetlistlabels{\arabic}{enumi}{enumii}{}{.}%
\item {} 
\sphinxAtStartPar
\sphinxcode{\sphinxupquote{GenConstr}} Class: Encapsulation of operations related to Indicator constraints in COPT.

\item {} 
\sphinxAtStartPar
\sphinxcode{\sphinxupquote{GenConstrArray}} Class: Facilitates user operations on a group of \sphinxcode{\sphinxupquote{GenConstr}} class objects.

\end{enumerate}

\item {} 
\sphinxAtStartPar
Encapsulation of Indicator constraint builders:
\begin{enumerate}
\sphinxsetlistlabels{\arabic}{enumi}{enumii}{}{.}%
\item {} 
\sphinxAtStartPar
\sphinxcode{\sphinxupquote{GenConstrBuilder}} Class: Encapsulation of builders for constructing Indicator constraints in COPT.

\item {} 
\sphinxAtStartPar
\sphinxcode{\sphinxupquote{GenConstrBuilderArray}} Class: Facilitates user operations on a group of \sphinxcode{\sphinxupquote{GenConstrBuilder}} class objects.

\end{enumerate}

\end{itemize}

\sphinxAtStartPar
For more details on the member methods and specific descriptions of the above Indicator constraint\sphinxhyphen{}related classes, please refer to the corresponding sections in the API reference manuals of each programming language.
\begin{itemize}
\item {} 
\sphinxAtStartPar
C++ API: {\hyperref[\detokenize{cppapiref:chapcppapiref-genconstr}]{\sphinxcrossref{\DUrole{std,std-ref}{GenConstr Class}}}} , {\hyperref[\detokenize{cppapiref:chapcppapiref-genconstrarray}]{\sphinxcrossref{\DUrole{std,std-ref}{GenConstrArray Class}}}} , {\hyperref[\detokenize{cppapiref:chapcppapiref-genconstrbuilder}]{\sphinxcrossref{\DUrole{std,std-ref}{GenConstrBuilder Class}}}} , {\hyperref[\detokenize{cppapiref:chapcppapiref-genconstrbuilderarray}]{\sphinxcrossref{\DUrole{std,std-ref}{GenConstrBuilderArray Class}}}}

\item {} 
\sphinxAtStartPar
C\# API: {\hyperref[\detokenize{csharpapiref:chapcsharpapiref-genconstr}]{\sphinxcrossref{\DUrole{std,std-ref}{GenConstr Class}}}} , {\hyperref[\detokenize{csharpapiref:chapcsharpapiref-genconstrarray}]{\sphinxcrossref{\DUrole{std,std-ref}{GenConstrArray Class}}}} , {\hyperref[\detokenize{csharpapiref:chapcsharpapiref-genconstrbuilder}]{\sphinxcrossref{\DUrole{std,std-ref}{GenConstrBuilder Class}}}} , {\hyperref[\detokenize{csharpapiref:chapcsharpapiref-genconstrbuilderarray}]{\sphinxcrossref{\DUrole{std,std-ref}{GenConstrBuilderArray Class}}}}

\item {} 
\sphinxAtStartPar
Java API: {\hyperref[\detokenize{javaapiref:chapjavaapiref-genconstr}]{\sphinxcrossref{\DUrole{std,std-ref}{GenConstr Class}}}} , {\hyperref[\detokenize{javaapiref:chapjavaapiref-genconstrarray}]{\sphinxcrossref{\DUrole{std,std-ref}{GenConstrArray Class}}}} , {\hyperref[\detokenize{javaapiref:chapjavaapiref-genconstrbuilder}]{\sphinxcrossref{\DUrole{std,std-ref}{GenConstrBuilder Class}}}} , {\hyperref[\detokenize{javaapiref:chapjavaapiref-genconstrbuilderarray}]{\sphinxcrossref{\DUrole{std,std-ref}{GenConstrBuilderArray Class}}}}

\item {} 
\sphinxAtStartPar
Python API: {\hyperref[\detokenize{pyapiref:chappyapi-genconstr}]{\sphinxcrossref{\DUrole{std,std-ref}{GenConstr Class}}}} , {\hyperref[\detokenize{pyapiref:chappyapi-genconstrarray}]{\sphinxcrossref{\DUrole{std,std-ref}{GenConstrArray Class}}}} , {\hyperref[\detokenize{pyapiref:chappyapi-sosbuilder}]{\sphinxcrossref{\DUrole{std,std-ref}{GenConstrBuilder Class}}}} , {\hyperref[\detokenize{pyapiref:chappyapi-genconstrbuilderarray}]{\sphinxcrossref{\DUrole{std,std-ref}{GenConstrBuilderArray Class}}}}

\end{itemize}

\subsection{Attributes for Special Constraints}
\label{\detokenize{modeling:attributes-for-special-constraints}}
\sphinxAtStartPar
COPT provides the following attributes to describe the number of special constraints in the model, as shown in \hyperref[\detokenize{modeling:copttab-special-attr}]{Table \ref{\detokenize{modeling:copttab-special-attr}}}.

\sphinxAtStartPar
For methods of retrieving these attributes in different programming interfaces, please refer to: {\hyperref[\detokenize{attribute:chapattrs}]{\sphinxcrossref{\DUrole{std,std-ref}{Attributes Section}}}}.

\begin{savenotes}\sphinxattablestart
\sphinxthistablewithglobalstyle
\centering
\sphinxcapstartof{table}
\sphinxthecaptionisattop
\sphinxcaption{Overview of Special Constraints\sphinxhyphen{}Related Attributes}\label{\detokenize{modeling:copttab-special-attr}}
\sphinxaftertopcaption
\begin{tabular}[t]{|\X{15}{54}|\X{9}{54}|\X{30}{54}|}
\sphinxtoprule
\sphinxstyletheadfamily 
\sphinxAtStartPar
Attribute
&\sphinxstyletheadfamily 
\sphinxAtStartPar
Type
&\sphinxstyletheadfamily 
\sphinxAtStartPar
Description
\\
\sphinxmidrule
\sphinxtableatstartofbodyhook
\sphinxAtStartPar
\sphinxcode{\sphinxupquote{Soss}}
&
\sphinxAtStartPar
Integer
&
\sphinxAtStartPar
Number of SOS constraints
\\
\sphinxhline
\sphinxAtStartPar
\sphinxcode{\sphinxupquote{Indicators}}
&
\sphinxAtStartPar
Integer
&
\sphinxAtStartPar
Number of indicator constraints
\\
\sphinxhline
\sphinxAtStartPar
\sphinxcode{\sphinxupquote{IISSOSs}}
&
\sphinxAtStartPar
Integer
&
\sphinxAtStartPar
Number of SOS constraints in IIS
\\
\sphinxhline
\sphinxAtStartPar
\sphinxcode{\sphinxupquote{IISIndicators}}
&
\sphinxAtStartPar
Integer
&
\sphinxAtStartPar
Number of indicator constraints in IIS
\\
\sphinxbottomrule
\end{tabular}
\sphinxtableafterendhook\par
\sphinxattableend\end{savenotes}

\subsection{IIS Status of Special Constraints}
\label{\detokenize{modeling:iis-status-of-special-constraints}}
\sphinxAtStartPar
Regarding the IIS (Irreducible Infeasible Set) calculation results for infeasible models,
COPT provides related functions to obtain the IIS status of SOS constraints.
Please refer to {\hyperref[\detokenize{infeasible:chapinfeas-getconstriis-special}]{\sphinxcrossref{\DUrole{std,std-ref}{Handling Infeasible Models Section: Retrieving IIS Status of Special Constraints}}}}.

\sphinxstepscope

\chapter{Handling Infeasible Models}
\label{\detokenize{infeasible:handling-infeasible-models}}\label{\detokenize{infeasible:chapinfeas}}\label{\detokenize{infeasible::doc}}
\sphinxAtStartPar
This chapter introduces two approaches supported by COPT for handling infeasible problems:
\begin{itemize}
\item {} 
\sphinxAtStartPar
{\hyperref[\detokenize{infeasible:chapinfeas-iis}]{\sphinxcrossref{\DUrole{std,std-ref}{IIS for Infeasible Models}}}}

\item {} 
\sphinxAtStartPar
{\hyperref[\detokenize{infeasible:chapinfeas-feasrelax}]{\sphinxcrossref{\DUrole{std,std-ref}{Feasibility Relaxation}}}}

\end{itemize}

\sphinxAtStartPar
In real\sphinxhyphen{}world problems, it is common to encounter infeasible models, which correspond to the solution status
code \sphinxcode{\sphinxupquote{COPT.INFEASIBLE}}. The main reasons for infeasibility are usually:
\begin{enumerate}
\sphinxsetlistlabels{\arabic}{enumi}{enumii}{}{.}%
\item {} 
\sphinxAtStartPar
Making some mistakes when modeling or inputting data (e.g., an empty left\sphinxhyphen{}hand side in a constraint).

\item {} 
\sphinxAtStartPar
The problem itself is infeasible, meaning some constraints or variable bounds are conflicting.

\end{enumerate}

\sphinxAtStartPar
COPT provides two methods for analyzing and handling infeasible models, which are supported for both
Linear Programming (LP) and Mixed\sphinxhyphen{}Integer linear programming (MILP):
\begin{enumerate}
\sphinxsetlistlabels{\arabic}{enumi}{enumii}{}{.}%
\item {} 
\sphinxAtStartPar
\sphinxstylestrong{Compute IIS}: Identify the key constraints and variable bounds causing infeasibility.

\item {} 
\sphinxAtStartPar
\sphinxstylestrong{Feasibility Relaxation (FeasRelax)}: Quantitatively compute the conflicts in constraints or variable bounds
(violations) that lead to infeasibility.

\end{enumerate}

\section{IIS for Infeasible Models}
\label{\detokenize{infeasible:iis-for-infeasible-models}}\label{\detokenize{infeasible:chapinfeas-iis}}
\sphinxAtStartPar
IIS (Irreducible Inconsistent Subsystem) refers to a minimal conflicting set in
the model that causes infeasibility, and has the following properties:
\begin{enumerate}
\sphinxsetlistlabels{\arabic}{enumi}{enumii}{}{.}%
\item {} 
\sphinxAtStartPar
The subsystem is still infeasible.

\item {} 
\sphinxAtStartPar
Removing any single constraint or variable bound from the IIS will make the subsystem feasible.

\end{enumerate}

\sphinxAtStartPar
\sphinxstylestrong{Note:} The IIS computed by COPT may not be minimal or unique. It may require
several iterations of modifying constraints and recomputing the IIS before the model becomes feasible.

\sphinxAtStartPar
Below is an example of an infeasible Linear Programming model:
\begin{equation}\label{equation:infeasible:infeasible:0}
\begin{split}\max\quad &z=12x_{1}+8x_{2}\\
\text{s.t.}\quad
&5 x_{1} + 2 x_{2} \geq 140\\
&2 x_{1} + 3 x_{2} \leq 90\\
&4 x_{1} + 2 x_{2} \leq 100\\
&x_{1}, x_{2}\geq0\end{split}
\end{equation}
\sphinxAtStartPar
The feasible region of the model is shown below:
\begin{quote}

\begin{figure}[H]
\centering

\noindent\sphinxincludegraphics[scale=0.8]{{copt-iisexample}.png}
\end{figure}
\end{quote}

\sphinxAtStartPar
From the figure, the conflicting constraints can be clearly seen:
\begin{equation}\label{equation:infeasible:infeasible:1}
\begin{split}&c1:5 x_{1} + 2x_{2} \geq 140\\
&c3:4 x_{1} + 2x_{2} \leq 100\end{split}
\end{equation}
\sphinxAtStartPar
After computing the IIS for the above model and writing it to a file (\sphinxcode{\sphinxupquote{.iis}} format),
the file content is as follows, consistent with the figure:

\begin{sphinxVerbatim}[commandchars=\\\{\}]
\PYGZbs{}Generated by Cardinal Operations

Maximize
   12 x[1] + 8 x[2]
Subject To
 c1: 5 x[1] + 2 x[2] \PYGZgt{}= 140
 c3: 4 x[1] + 2 x[2] \PYGZlt{}= 100
END
\end{sphinxVerbatim}

\subsection{Computing IIS}
\label{\detokenize{infeasible:computing-iis}}
\sphinxAtStartPar
COPT provides functions in various Programming API to compute the IIS of an infeasible model, returning a set of conflicting
constraints and variable bounds. The IIS can also be written to a file by specifying the file name suffix as \sphinxcode{\sphinxupquote{.iis}} (e.g., \sphinxcode{\sphinxupquote{example.iis}}).
Related function names are shown in \hyperref[\detokenize{infeasible:copttab-iiscompute}]{Table \ref{\detokenize{infeasible:copttab-iiscompute}}}:

\begin{savenotes}\sphinxattablestart
\sphinxthistablewithglobalstyle
\centering
\sphinxcapstartof{table}
\sphinxthecaptionisattop
\sphinxcaption{Functions for Computing IIS of Infeasible Models}\label{\detokenize{infeasible:copttab-iiscompute}}
\sphinxaftertopcaption
\begin{tabular}[t]{|\X{10}{65}|\X{25}{65}|\X{30}{65}|}
\sphinxtoprule
\sphinxstyletheadfamily 
\sphinxAtStartPar
API
&\sphinxstyletheadfamily 
\sphinxAtStartPar
Compute IIS
&\sphinxstyletheadfamily 
\sphinxAtStartPar
Write IIS to File
\\
\sphinxmidrule
\sphinxtableatstartofbodyhook
\sphinxAtStartPar
C
&
\sphinxAtStartPar
\sphinxcode{\sphinxupquote{COPT\_ComputeIIS}}
&
\sphinxAtStartPar
\sphinxcode{\sphinxupquote{COPT\_WriteIIS}}
\\
\sphinxhline
\sphinxAtStartPar
C++
&
\sphinxAtStartPar
\sphinxcode{\sphinxupquote{Model::ComputeIIS()}}
&
\sphinxAtStartPar
\sphinxcode{\sphinxupquote{Model::WriteIIS()}}
\\
\sphinxhline
\sphinxAtStartPar
C\#
&
\sphinxAtStartPar
\sphinxcode{\sphinxupquote{Model.ComputeIIS()}}
&
\sphinxAtStartPar
\sphinxcode{\sphinxupquote{Model.WriteIIS()}}
\\
\sphinxhline
\sphinxAtStartPar
Java
&
\sphinxAtStartPar
\sphinxcode{\sphinxupquote{Model.computeIIS()}}
&
\sphinxAtStartPar
\sphinxcode{\sphinxupquote{Model.writeIIS()}}
\\
\sphinxhline
\sphinxAtStartPar
Python
&
\sphinxAtStartPar
\sphinxcode{\sphinxupquote{Model.computeIIS()}}
&
\sphinxAtStartPar
\sphinxcode{\sphinxupquote{Model.writeIIS()}}
\\
\sphinxbottomrule
\end{tabular}
\sphinxtableafterendhook\par
\sphinxattableend\end{savenotes}

\subsection{Getting IIS status of variables and constraints}
\label{\detokenize{infeasible:getting-iis-status-of-variables-and-constraints}}

\subsubsection{Functions for variable’s IIS status}
\label{\detokenize{infeasible:functions-for-variable-s-iis-status}}
\sphinxAtStartPar
After computing the IIS, the IIS status of variables (lower/upper bound) can be obtained. The status indicates whether the bound is in the IIS:
\begin{itemize}
\item {} 
\sphinxAtStartPar
1: The specified variable bound (lower/upper) is in the IIS

\item {} 
\sphinxAtStartPar
0: The specified variable bound (lower/upper) is not in the IIS

\end{itemize}

\sphinxAtStartPar
Supported functions for different interfaces are shown in \hyperref[\detokenize{infeasible:copttab-getvariis}]{Table \ref{\detokenize{infeasible:copttab-getvariis}}}:

\begin{savenotes}\sphinxattablestart
\sphinxthistablewithglobalstyle
\centering
\sphinxcapstartof{table}
\sphinxthecaptionisattop
\sphinxcaption{Functions for Getting Variable IIS Status}\label{\detokenize{infeasible:copttab-getvariis}}
\sphinxaftertopcaption
\begin{tabular}[t]{|\X{10}{70}|\X{30}{70}|\X{30}{70}|}
\sphinxtoprule
\sphinxstyletheadfamily 
\sphinxAtStartPar
API
&\sphinxstyletheadfamily 
\sphinxAtStartPar
Lower Bound IIS Status
&\sphinxstyletheadfamily 
\sphinxAtStartPar
Upper Bound IIS Status
\\
\sphinxmidrule
\sphinxtableatstartofbodyhook
\sphinxAtStartPar
C
&
\sphinxAtStartPar
\sphinxcode{\sphinxupquote{COPT\_GetColLowerIIS}}
&
\sphinxAtStartPar
\sphinxcode{\sphinxupquote{COPT\_GetColUpperIIS}}
\\
\sphinxhline
\sphinxAtStartPar
C++
&
\sphinxAtStartPar
\sphinxcode{\sphinxupquote{Model::GetVarLowerIIS()}}
&
\sphinxAtStartPar
\sphinxcode{\sphinxupquote{Model::GetVarUpperIIS()}}
\\
\sphinxhline
\sphinxAtStartPar
C\#
&
\sphinxAtStartPar
\sphinxcode{\sphinxupquote{Model.GetVarLowerIIS()}}
&
\sphinxAtStartPar
\sphinxcode{\sphinxupquote{Model.GetVarUpperIIS()}}
\\
\sphinxhline
\sphinxAtStartPar
Java
&
\sphinxAtStartPar
\sphinxcode{\sphinxupquote{Model.getVarLowerIIS()}}
&
\sphinxAtStartPar
\sphinxcode{\sphinxupquote{Model.getVarUpperIIS()}}
\\
\sphinxhline
\sphinxAtStartPar
Python
&
\sphinxAtStartPar
\sphinxcode{\sphinxupquote{Model.getVarLowerIIS()}}
&
\sphinxAtStartPar
\sphinxcode{\sphinxupquote{Model.getVarUpperIIS()}}
\\
\sphinxbottomrule
\end{tabular}
\sphinxtableafterendhook\par
\sphinxattableend\end{savenotes}

\begin{sphinxadmonition}{note}{Notes:}
\begin{itemize}
\item {} 
\sphinxAtStartPar
The above functions accept either a single variable (\sphinxcode{\sphinxupquote{Var}} object) or a set of variables (\sphinxcode{\sphinxupquote{VarArray}} or \sphinxcode{\sphinxupquote{tupledict}} objects).

\item {} 
\sphinxAtStartPar
Except for the C API, other interfaces provide member functions in the \sphinxcode{\sphinxupquote{Var}} class to get the IIS status for a single variable.
For example, in Python: \sphinxcode{\sphinxupquote{Var.getLowerIIS()}} and \sphinxcode{\sphinxupquote{Var.getUpperIIS()}}.

\item {} 
\sphinxAtStartPar
In the Python API, you can also directly access the IIS status as member attributes of the \sphinxcode{\sphinxupquote{Var}} object:
\sphinxcode{\sphinxupquote{Var.iislb}} for the lower bound IIS status, \sphinxcode{\sphinxupquote{Var.iisub}} for the upper bound IIS status.

\end{itemize}
\end{sphinxadmonition}

\subsubsection{Functions for constraint’s IIS status}
\label{\detokenize{infeasible:functions-for-constraint-s-iis-status}}
\sphinxAtStartPar
After computing the IIS, COPT also supports querying the IIS status of constraint bounds (lower/upper).
The status indicates whether the bound is in the IIS:
\begin{itemize}
\item {} 
\sphinxAtStartPar
1: The specified constraint bound (lower/upper) is in the IIS

\item {} 
\sphinxAtStartPar
0: The specified constraint bound (lower/upper) is not in the IIS

\end{itemize}

\sphinxAtStartPar
Supported functions for different interfaces are shown in \hyperref[\detokenize{infeasible:copttab-getconstriis}]{Table \ref{\detokenize{infeasible:copttab-getconstriis}}}:

\begin{savenotes}\sphinxattablestart
\sphinxthistablewithglobalstyle
\centering
\sphinxcapstartof{table}
\sphinxthecaptionisattop
\sphinxcaption{Functions for Getting Constraint IIS Status}\label{\detokenize{infeasible:copttab-getconstriis}}
\sphinxaftertopcaption
\begin{tabular}[t]{|\X{10}{70}|\X{30}{70}|\X{30}{70}|}
\sphinxtoprule
\sphinxstyletheadfamily 
\sphinxAtStartPar
API
&\sphinxstyletheadfamily 
\sphinxAtStartPar
Lower Bound IIS Status
&\sphinxstyletheadfamily 
\sphinxAtStartPar
Upper Bound IIS Status
\\
\sphinxmidrule
\sphinxtableatstartofbodyhook
\sphinxAtStartPar
C
&
\sphinxAtStartPar
\sphinxcode{\sphinxupquote{COPT\_GetRowLowerIIS}}
&
\sphinxAtStartPar
\sphinxcode{\sphinxupquote{COPT\_GetRowUpperIIS}}
\\
\sphinxhline
\sphinxAtStartPar
C++
&
\sphinxAtStartPar
\sphinxcode{\sphinxupquote{Model::GetConstrLowerIIS()}}
&
\sphinxAtStartPar
\sphinxcode{\sphinxupquote{Model::GetConstrUpperIIS()}}
\\
\sphinxhline
\sphinxAtStartPar
C\#
&
\sphinxAtStartPar
\sphinxcode{\sphinxupquote{Model.GetConstrLowerIIS()}}
&
\sphinxAtStartPar
\sphinxcode{\sphinxupquote{Model.GetConstrUpperIIS()}}
\\
\sphinxhline
\sphinxAtStartPar
Java
&
\sphinxAtStartPar
\sphinxcode{\sphinxupquote{Model.getConstrLowerIIS()}}
&
\sphinxAtStartPar
\sphinxcode{\sphinxupquote{Model.getConstrUpperIIS()}}
\\
\sphinxhline
\sphinxAtStartPar
Python
&
\sphinxAtStartPar
\sphinxcode{\sphinxupquote{Model.getConstrLowerIIS()}}
&
\sphinxAtStartPar
\sphinxcode{\sphinxupquote{Model.getConstrUpperIIS()}}
\\
\sphinxbottomrule
\end{tabular}
\sphinxtableafterendhook\par
\sphinxattableend\end{savenotes}

\begin{sphinxadmonition}{note}{Notes:}
\begin{itemize}
\item {} 
\sphinxAtStartPar
The above functions accept either a single constraint (\sphinxcode{\sphinxupquote{Constraint}} object) or a set of constraints (\sphinxcode{\sphinxupquote{Constraint}}
or \sphinxcode{\sphinxupquote{ConstrArray}} objects).

\item {} 
\sphinxAtStartPar
Except for the C API, other interfaces provide member functions in the \sphinxcode{\sphinxupquote{Constraint}} class to get the IIS status
for a single constraint. For example, in Python: \sphinxcode{\sphinxupquote{Constraint.getLowerIIS()}} and \sphinxcode{\sphinxupquote{Constraint.getUpperIIS()}}.

\item {} 
\sphinxAtStartPar
In the Python API, you can also directly access the IIS status as member attributes of the \sphinxcode{\sphinxupquote{Constraint}} object:
\sphinxcode{\sphinxupquote{Constraint.iislb}} for the lower bound IIS status, \sphinxcode{\sphinxupquote{Constraint.iisub}} for the upper bound IIS status.

\end{itemize}
\end{sphinxadmonition}

\subsubsection{Functions for special constraints}
\label{\detokenize{infeasible:functions-for-special-constraints}}\label{\detokenize{infeasible:chapinfeas-getconstriis-special}}
\sphinxAtStartPar
COPT also provides functions, as shown in \hyperref[\detokenize{infeasible:copttab-getconstriis-special}]{Table \ref{\detokenize{infeasible:copttab-getconstriis-special}}}, for getting the IIS status of SOS and Indicator constraints:

\begin{savenotes}\sphinxattablestart
\sphinxthistablewithglobalstyle
\centering
\sphinxcapstartof{table}
\sphinxthecaptionisattop
\sphinxcaption{Functions for Getting IIS Status of Special Constraints}\label{\detokenize{infeasible:copttab-getconstriis-special}}
\sphinxaftertopcaption
\begin{tabular}[t]{|\X{10}{70}|\X{30}{70}|\X{30}{70}|}
\sphinxtoprule
\sphinxstyletheadfamily 
\sphinxAtStartPar
API
&\sphinxstyletheadfamily 
\sphinxAtStartPar
SOS Constraint
&\sphinxstyletheadfamily 
\sphinxAtStartPar
Indicator Constraint
\\
\sphinxmidrule
\sphinxtableatstartofbodyhook
\sphinxAtStartPar
C
&
\sphinxAtStartPar
\sphinxcode{\sphinxupquote{COPT\_GetSOSIIS}}
&
\sphinxAtStartPar
\sphinxcode{\sphinxupquote{COPT\_GetIndicatorIIS}}
\\
\sphinxhline
\sphinxAtStartPar
C++
&
\sphinxAtStartPar
\sphinxcode{\sphinxupquote{Model::GetSOSIIS()}}
&
\sphinxAtStartPar
\sphinxcode{\sphinxupquote{Model::GetIndicatorIIS()}}
\\
\sphinxhline
\sphinxAtStartPar
C\#
&
\sphinxAtStartPar
\sphinxcode{\sphinxupquote{Model.GetSOSIIS()}}
&
\sphinxAtStartPar
\sphinxcode{\sphinxupquote{Model.GetIndicatorIIS()}}
\\
\sphinxhline
\sphinxAtStartPar
Java
&
\sphinxAtStartPar
\sphinxcode{\sphinxupquote{Model.getSOSIIS()}}
&
\sphinxAtStartPar
\sphinxcode{\sphinxupquote{Model.getIndicatorIIS()}}
\\
\sphinxhline
\sphinxAtStartPar
Python
&
\sphinxAtStartPar
\sphinxcode{\sphinxupquote{Model.getSOSIIS()}}
&
\sphinxAtStartPar
\sphinxcode{\sphinxupquote{Model.getIndicatorIIS()}}
\\
\sphinxbottomrule
\end{tabular}
\sphinxtableafterendhook\par
\sphinxattableend\end{savenotes}

\sphinxAtStartPar
For details on the usage of the above functions, please refer to the API reference manual for each programming interface:
\begin{itemize}
\item {} 
\sphinxAtStartPar
C API: {\hyperref[\detokenize{capiref:chapapi-iiscompute}]{\sphinxcrossref{\DUrole{std,std-ref}{IIS computation functions}}}}

\item {} 
\sphinxAtStartPar
C++ API: {\hyperref[\detokenize{cppapiref:chapcppapiref-model}]{\sphinxcrossref{\DUrole{std,std-ref}{Model class}}}}, {\hyperref[\detokenize{cppapiref:chapcppapiref-model}]{\sphinxcrossref{\DUrole{std,std-ref}{Var class}}}}, {\hyperref[\detokenize{cppapiref:chapcppapiref-model}]{\sphinxcrossref{\DUrole{std,std-ref}{Constraint class}}}}

\item {} 
\sphinxAtStartPar
C\# API: {\hyperref[\detokenize{csharpapiref:chapcsharpapiref-model}]{\sphinxcrossref{\DUrole{std,std-ref}{Model class}}}}, {\hyperref[\detokenize{cppapiref:chapcppapiref-model}]{\sphinxcrossref{\DUrole{std,std-ref}{Var class}}}}, {\hyperref[\detokenize{cppapiref:chapcppapiref-model}]{\sphinxcrossref{\DUrole{std,std-ref}{Constraint class}}}}

\item {} 
\sphinxAtStartPar
Java API: {\hyperref[\detokenize{javaapiref:chapjavaapiref-model}]{\sphinxcrossref{\DUrole{std,std-ref}{Model class}}}}, {\hyperref[\detokenize{cppapiref:chapcppapiref-model}]{\sphinxcrossref{\DUrole{std,std-ref}{Var class}}}}, {\hyperref[\detokenize{cppapiref:chapcppapiref-model}]{\sphinxcrossref{\DUrole{std,std-ref}{Constraint class}}}}

\item {} 
\sphinxAtStartPar
Python API: {\hyperref[\detokenize{pyapiref:chappyapi-model}]{\sphinxcrossref{\DUrole{std,std-ref}{Model class}}}}, {\hyperref[\detokenize{cppapiref:chapcppapiref-model}]{\sphinxcrossref{\DUrole{std,std-ref}{Var class}}}}, {\hyperref[\detokenize{cppapiref:chapcppapiref-model}]{\sphinxcrossref{\DUrole{std,std-ref}{Constraint class}}}}

\end{itemize}

\subsection{IIS\sphinxhyphen{}related parameters, attributes and information}
\label{\detokenize{infeasible:iis-related-parameters-attributes-and-information}}

\subsubsection{Parameters}
\label{\detokenize{infeasible:parameters}}
\sphinxAtStartPar
Users can select the method for computing IIS by setting the \sphinxcode{\sphinxupquote{"IISMethod"}} parameter.
For parameter setting in each interface, see {\hyperref[\detokenize{parameter:chapparams}]{\sphinxcrossref{\DUrole{std,std-ref}{Parameter section}}}}.
\begin{itemize}
\item {} 
\sphinxAtStartPar
\sphinxcode{\sphinxupquote{IISMethod}}
\begin{quote}

\sphinxAtStartPar
Integer parameter.

\sphinxAtStartPar
The method for computing IIS.

\sphinxAtStartPar
\sphinxstylestrong{Default value:} \sphinxhyphen{}1

\sphinxAtStartPar
\sphinxstylestrong{Possible values:}
\begin{quote}

\sphinxAtStartPar
\sphinxhyphen{}1: Automatically selected.

\sphinxAtStartPar
0: Prioritize quality of the IIS result.

\sphinxAtStartPar
1: Prioritize computation efficiency.
\end{quote}
\end{quote}

\end{itemize}

\subsubsection{Attributes}
\label{\detokenize{infeasible:attributes}}
\sphinxAtStartPar
COPT provides related attributes to describe the IIS result, mainly indicating the existence of an IIS and the number
of variables/constraints in the IIS. See {\hyperref[\detokenize{attribute:chapattrs}]{\sphinxcrossref{\DUrole{std,std-ref}{Attributes section}}}} for details.

\begin{savenotes}\sphinxattablestart
\sphinxthistablewithglobalstyle
\centering
\sphinxcapstartof{table}
\sphinxthecaptionisattop
\sphinxcaption{Overview of IIS\sphinxhyphen{}Related Attributes}\label{\detokenize{infeasible:id4}}
\sphinxaftertopcaption
\begin{tabular}[t]{|\X{15}{69}|\X{9}{69}|\X{45}{69}|}
\sphinxtoprule
\sphinxstyletheadfamily 
\sphinxAtStartPar
Attribute Name
&\sphinxstyletheadfamily 
\sphinxAtStartPar
Type
&\sphinxstyletheadfamily 
\sphinxAtStartPar
Description
\\
\sphinxmidrule
\sphinxtableatstartofbodyhook
\sphinxAtStartPar
{\hyperref[\detokenize{attribute:hasiis}]{\sphinxcrossref{\DUrole{std,std-ref}{HasIIS}}}}
&
\sphinxAtStartPar
Integer
&
\sphinxAtStartPar
Indicates whether an IIS exists.
\\
\sphinxhline
\sphinxAtStartPar
{\hyperref[\detokenize{attribute:isminiis}]{\sphinxcrossref{\DUrole{std,std-ref}{IsMinIIS}}}}
&
\sphinxAtStartPar
Integer
&
\sphinxAtStartPar
Indicates whether the computed IIS is minimal.
\\
\sphinxhline
\sphinxAtStartPar
{\hyperref[\detokenize{attribute:iiscols}]{\sphinxcrossref{\DUrole{std,std-ref}{IISCols}}}}
&
\sphinxAtStartPar
Integer
&
\sphinxAtStartPar
Number of variable bounds in the IIS.
\\
\sphinxhline
\sphinxAtStartPar
{\hyperref[\detokenize{attribute:iisrows}]{\sphinxcrossref{\DUrole{std,std-ref}{IISRows}}}}
&
\sphinxAtStartPar
Integer
&
\sphinxAtStartPar
Number of constraints in the IIS.
\\
\sphinxhline
\sphinxAtStartPar
{\hyperref[\detokenize{attribute:iissoss}]{\sphinxcrossref{\DUrole{std,std-ref}{IISSOSs}}}}
&
\sphinxAtStartPar
Integer
&
\sphinxAtStartPar
Number of SOS constraints in the IIS.
\\
\sphinxhline
\sphinxAtStartPar
{\hyperref[\detokenize{attribute:iisindicators}]{\sphinxcrossref{\DUrole{std,std-ref}{IISIndicators}}}}
&
\sphinxAtStartPar
Integer
&
\sphinxAtStartPar
Number of Indicator constraints in the IIS.
\\
\sphinxbottomrule
\end{tabular}
\sphinxtableafterendhook\par
\sphinxattableend\end{savenotes}

\subsubsection{Example Code}
\label{\detokenize{infeasible:example-code}}
\sphinxAtStartPar
Users can write the infeasible model to a file and directly compute the IIS after loading it in COPT. For example, in Python:

\begin{sphinxVerbatim}[commandchars=\\\{\}]
\PYG{k+kn}{from}\PYG{+w}{ }\PYG{n+nn}{coptpy}\PYG{+w}{ }\PYG{k+kn}{import} \PYG{o}{*}

\PYG{n}{env} \PYG{o}{=} \PYG{n}{Envr}\PYG{p}{(}\PYG{p}{)}
\PYG{n}{model} \PYG{o}{=} \PYG{n}{env}\PYG{o}{.}\PYG{n}{createModel}\PYG{p}{(}\PYG{l+s+s2}{\PYGZdq{}}\PYG{l+s+s2}{example}\PYG{l+s+s2}{\PYGZdq{}}\PYG{p}{)}
\PYG{n}{model}\PYG{o}{.}\PYG{n}{read}\PYG{p}{(}\PYG{l+s+s2}{\PYGZdq{}}\PYG{l+s+s2}{example.lp}\PYG{l+s+s2}{\PYGZdq{}}\PYG{p}{)}
\PYG{n}{model}\PYG{o}{.}\PYG{n}{computeIIS}\PYG{p}{(}\PYG{p}{)}
\PYG{n}{model}\PYG{o}{.}\PYG{n}{writeIIS}\PYG{p}{(}\PYG{l+s+s2}{\PYGZdq{}}\PYG{l+s+s2}{example.iis}\PYG{l+s+s2}{\PYGZdq{}}\PYG{p}{)}
\end{sphinxVerbatim}

\sphinxAtStartPar
Sample code for computing IIS in different APIs can be found in the \sphinxcode{\sphinxupquote{"examples"}} folder of the COPT installation package.
For COPT Python Interface, the sample path is: \sphinxcode{\sphinxupquote{"/examples/python/iis\_ex1.py"}}.

\section{Feasibility Relaxation for infeasible models}
\label{\detokenize{infeasible:feasibility-relaxation-for-infeasible-models}}\label{\detokenize{infeasible:chapinfeas-feasrelax}}
\sphinxAtStartPar
Feasibility relaxation is the process of minimizing the conflicts in the bounds of variables and constraints in the original infeasible model.
Users can use the quantitative results from feasibility relaxation to adjust constraints or variable bounds, thus making the model feasible.

\subsection{Computing Feasibility Relaxation}
\label{\detokenize{infeasible:computing-feasibility-relaxation}}
\sphinxAtStartPar
COPT provides functions in different interfaces for computing feasibility relaxation, and also allows writing the relaxed
model to a file with the suffix \sphinxcode{\sphinxupquote{.relax}} (e.g., \sphinxcode{\sphinxupquote{example.relax}}). Related function names are listed in \hyperref[\detokenize{infeasible:copttab-feasrelax}]{Table \ref{\detokenize{infeasible:copttab-feasrelax}}}:

\begin{savenotes}\sphinxattablestart
\sphinxthistablewithglobalstyle
\centering
\sphinxcapstartof{table}
\sphinxthecaptionisattop
\sphinxcaption{Functions for Computing Feasibility Relaxation}\label{\detokenize{infeasible:copttab-feasrelax}}
\sphinxaftertopcaption
\begin{tabular}[t]{|\X{10}{65}|\X{25}{65}|\X{30}{65}|}
\sphinxtoprule
\sphinxstyletheadfamily 
\sphinxAtStartPar
API
&\sphinxstyletheadfamily 
\sphinxAtStartPar
Compute FeasRelax
&\sphinxstyletheadfamily 
\sphinxAtStartPar
Write Relaxed Model to File
\\
\sphinxmidrule
\sphinxtableatstartofbodyhook
\sphinxAtStartPar
C
&
\sphinxAtStartPar
\sphinxcode{\sphinxupquote{COPT\_FeasRelax}}
&
\sphinxAtStartPar
\sphinxcode{\sphinxupquote{COPT\_WriteRelax}}
\\
\sphinxhline
\sphinxAtStartPar
C++
&
\sphinxAtStartPar
\sphinxcode{\sphinxupquote{Model::FeasRelax()}}
&
\sphinxAtStartPar
\sphinxcode{\sphinxupquote{Model::WriteRelax()}}
\\
\sphinxhline
\sphinxAtStartPar
C\#
&
\sphinxAtStartPar
\sphinxcode{\sphinxupquote{Model.FeasRelax()}}
&
\sphinxAtStartPar
\sphinxcode{\sphinxupquote{Model.WriteRelax()}}
\\
\sphinxhline
\sphinxAtStartPar
Java
&
\sphinxAtStartPar
\sphinxcode{\sphinxupquote{Model.feasRelax()}}
&
\sphinxAtStartPar
\sphinxcode{\sphinxupquote{Model.writeRelax()}}
\\
\sphinxhline
\sphinxAtStartPar
Python
&
\sphinxAtStartPar
\sphinxcode{\sphinxupquote{Model.feasRelax()}}
&
\sphinxAtStartPar
\sphinxcode{\sphinxupquote{Model.writeRelax()}}
\\
\sphinxbottomrule
\end{tabular}
\sphinxtableafterendhook\par
\sphinxattableend\end{savenotes}

\sphinxAtStartPar
COPT supports two approaches for computing feasibility relaxation.
In interfaces other than Python, users can use different function arguments:
\begin{enumerate}
\sphinxsetlistlabels{\arabic}{enumi}{enumii}{}{.}%
\item {} 
\sphinxAtStartPar
Simplified version: Relax all variables and/or constraints in the model, only requiring two parameters to specify
whether to relax all variables or all constraints.
\begin{itemize}
\item {} 
\sphinxAtStartPar
\sphinxcode{\sphinxupquote{ifRelaxVars}} : whether to relax variables (default: \sphinxcode{\sphinxupquote{True}})

\item {} 
\sphinxAtStartPar
\sphinxcode{\sphinxupquote{ifRelaxCons}} : whether to relax constraints (default: \sphinxcode{\sphinxupquote{True}})

\end{itemize}

\item {} 
\sphinxAtStartPar
Full version: Accepts more parameters (specify variables/constraints, set penalty factors for bounds).
\begin{itemize}
\item {} 
\sphinxAtStartPar
\sphinxcode{\sphinxupquote{vars}}: variables to relax

\item {} 
\sphinxAtStartPar
\sphinxcode{\sphinxupquote{cons}}: constraints to relax

\item {} 
\sphinxAtStartPar
\sphinxcode{\sphinxupquote{colLowPen}}: penalty factor for variable lower bounds

\item {} 
\sphinxAtStartPar
\sphinxcode{\sphinxupquote{colUppPen}}: penalty factor for variable upper bounds

\item {} 
\sphinxAtStartPar
\sphinxcode{\sphinxupquote{rowBndPen}}: penalty factor for constraint bounds

\item {} 
\sphinxAtStartPar
\sphinxcode{\sphinxupquote{rowUppPen}}: penalty factor for constraint upper bounds (for double\sphinxhyphen{}sided constraints)

\end{itemize}

\sphinxAtStartPar
\sphinxstylestrong{Note:} In most cases, set \sphinxcode{\sphinxupquote{rowUppPen}} as \sphinxcode{\sphinxupquote{NULL}}.

\end{enumerate}

\sphinxAtStartPar
The function names and argument conventions differ slightly across APIs,
but the functionality and meaning are consistent. See each interface’s API manual for details:
\begin{itemize}
\item {} 
\sphinxAtStartPar
C API: {\hyperref[\detokenize{capiref:chapapi-feasrelax}]{\sphinxcrossref{\DUrole{std,std-ref}{Feasibility relaxation functions}}}}

\item {} 
\sphinxAtStartPar
C++ API: {\hyperref[\detokenize{cppapiref:chapcppapiref-model}]{\sphinxcrossref{\DUrole{std,std-ref}{Model class}}}}

\item {} 
\sphinxAtStartPar
C\# API: {\hyperref[\detokenize{csharpapiref:chapcsharpapiref-model}]{\sphinxcrossref{\DUrole{std,std-ref}{Model class}}}}

\item {} 
\sphinxAtStartPar
Java API: {\hyperref[\detokenize{javaapiref:chapjavaapiref-model}]{\sphinxcrossref{\DUrole{std,std-ref}{Model class}}}}

\item {} 
\sphinxAtStartPar
Python API: {\hyperref[\detokenize{pyapiref:chappyapi-model}]{\sphinxcrossref{\DUrole{std,std-ref}{Model class}}}}

\end{itemize}

\begin{sphinxadmonition}{note}{Notes:}
\begin{itemize}
\item {} 
\sphinxAtStartPar
For the simplified feasibility relaxation approach, the Python API additionally
provides the function \sphinxcode{\sphinxupquote{Model.feasrelaxS(vrelax, crelax)}}, requiring only two arguments:
\begin{itemize}
\item {} 
\sphinxAtStartPar
\sphinxcode{\sphinxupquote{vrelax}}: whether to relax variables (default: \sphinxcode{\sphinxupquote{True}})

\item {} 
\sphinxAtStartPar
\sphinxcode{\sphinxupquote{crelax}}: whether to relax constraints (default: \sphinxcode{\sphinxupquote{True}})

\end{itemize}

\end{itemize}
\end{sphinxadmonition}

\subsection{Feasrelax\sphinxhyphen{}related parameters, attributes and information}
\label{\detokenize{infeasible:feasrelax-related-parameters-attributes-and-information}}

\subsubsection{Parameters}
\label{\detokenize{infeasible:id1}}
\sphinxAtStartPar
Users can select the method for feasibility relaxation by setting the \sphinxcode{\sphinxupquote{"FeasRelaxMode"}} parameter.
See {\hyperref[\detokenize{parameter:chapparams}]{\sphinxcrossref{\DUrole{std,std-ref}{Parameter section}}}} for details.
\begin{itemize}
\item {} 
\sphinxAtStartPar
\sphinxcode{\sphinxupquote{FeasRelaxMode}}
\begin{quote}

\sphinxAtStartPar
Integer parameter.

\sphinxAtStartPar
The method for computing feasibility relaxation.

\sphinxAtStartPar
\sphinxstylestrong{Default value:} 0

\sphinxAtStartPar
\sphinxstylestrong{Possible values:}
\begin{quote}

\sphinxAtStartPar
0: Minimize the weighted\sphinxhyphen{}sum of violations.

\sphinxAtStartPar
1: Minimize the original model objective under the weighted sum of violations.

\sphinxAtStartPar
2: Minimize the number of violations.

\sphinxAtStartPar
3: Minimize the original model objective under the minimum number of violations.

\sphinxAtStartPar
4: Minimize the weighted\sphinxhyphen{}sum of squared violations.

\sphinxAtStartPar
5: Minimize the original model objective under the weighted\sphinxhyphen{}sum of squared violations.
\end{quote}
\end{quote}

\end{itemize}

\subsubsection{Attributes}
\label{\detokenize{infeasible:id2}}
\sphinxAtStartPar
COPT provides attributes to describe the result of feasibility relaxation, as shown in \hyperref[\detokenize{infeasible:copttab-feasrelax-attr}]{Table \ref{\detokenize{infeasible:copttab-feasrelax-attr}}}.
See {\hyperref[\detokenize{attribute:chapattrs}]{\sphinxcrossref{\DUrole{std,std-ref}{Attributes section}}}} for details.

\begin{savenotes}\sphinxattablestart
\sphinxthistablewithglobalstyle
\centering
\sphinxcapstartof{table}
\sphinxthecaptionisattop
\sphinxcaption{Overview of Feasibility Relaxation Attributes}\label{\detokenize{infeasible:copttab-feasrelax-attr}}
\sphinxaftertopcaption
\begin{tabular}[t]{|\X{15}{60}|\X{15}{60}|\X{30}{60}|}
\sphinxtoprule
\sphinxstyletheadfamily 
\sphinxAtStartPar
Attribute Name
&\sphinxstyletheadfamily 
\sphinxAtStartPar
Type
&\sphinxstyletheadfamily 
\sphinxAtStartPar
Description
\\
\sphinxmidrule
\sphinxtableatstartofbodyhook
\sphinxAtStartPar
{\hyperref[\detokenize{attribute:hasfeasrelaxsol}]{\sphinxcrossref{\DUrole{std,std-ref}{HasFeasRelaxSol}}}}
&
\sphinxAtStartPar
Integer
&
\sphinxAtStartPar
Indicates whether a feasible relaxation solution exists.
\\
\sphinxhline
\sphinxAtStartPar
{\hyperref[\detokenize{attribute:feasrelaxobj}]{\sphinxcrossref{\DUrole{std,std-ref}{FeasRelaxObj}}}}
&
\sphinxAtStartPar
Double
&
\sphinxAtStartPar
Feasibility relaxation value.
\\
\sphinxbottomrule
\end{tabular}
\sphinxtableafterendhook\par
\sphinxattableend\end{savenotes}

\subsubsection{Information}
\label{\detokenize{infeasible:information}}
\sphinxAtStartPar
COPT provides the following information, representing the amount of relaxation for lower and upper bounds of variables
(or constraints). See {\hyperref[\detokenize{information:chapinfo}]{\sphinxcrossref{\DUrole{std,std-ref}{Information section}}}}.

\phantomsection\label{\detokenize{infeasible:relaxlb}}\begin{itemize}
\item {} 
\sphinxAtStartPar
\sphinxcode{\sphinxupquote{RelaxLB}}
\begin{quote}

\sphinxAtStartPar
Double information.

\sphinxAtStartPar
The relaxation amount for the lower bound of a variable (column) or constraint (row).
\end{quote}

\end{itemize}
\phantomsection\label{\detokenize{infeasible:relaxub}}\begin{itemize}
\item {} 
\sphinxAtStartPar
\sphinxcode{\sphinxupquote{RelaxUB}}
\begin{quote}

\sphinxAtStartPar
Double information.

\sphinxAtStartPar
The relaxation amount for the upper bound of a variable (column) or constraint (row).
\end{quote}

\end{itemize}

\subsubsection{Example Code}
\label{\detokenize{infeasible:id3}}
\sphinxAtStartPar
Sample code for feasibility relaxation can be found in the \sphinxcode{\sphinxupquote{"examples"}} folder of the COPT installation package, with the file name \sphinxcode{\sphinxupquote{"feasrelax\_ex1.py"}}.
For COPT Python API, the sample path is: \sphinxcode{\sphinxupquote{"/examples/python/feasrelax\_ex1.py"}}.

\sphinxstepscope

\chapter{MIP Starts}
\label{\detokenize{mipstart:mip-starts}}\label{\detokenize{mipstart:chapmipstart}}\label{\detokenize{mipstart::doc}}

\section{Utilities of MIP Starts}
\label{\detokenize{mipstart:utilities-of-mip-starts}}

\subsection{Set and Load MIP Starts}
\label{\detokenize{mipstart:set-and-load-mip-starts}}
\sphinxAtStartPar
For MIP problems, COPT provides methods to specify initial solution value(s) for a single variable or set of variables and load it/them into model. The parameters that can be specified are:
\begin{itemize}
\item {} 
\sphinxAtStartPar
\sphinxcode{\sphinxupquote{vars}} :variables

\item {} 
\sphinxAtStartPar
\sphinxcode{\sphinxupquote{startvals}} :variables’ solution values

\end{itemize}

\sphinxAtStartPar
The functions in different APIs are shown in \hyperref[\detokenize{mipstart:copttab-mipstart-set}]{Table \ref{\detokenize{mipstart:copttab-mipstart-set}}}:

\begin{savenotes}\sphinxattablestart
\sphinxthistablewithglobalstyle
\centering
\sphinxcapstartof{table}
\sphinxthecaptionisattop
\sphinxcaption{Functions for setting MIP starts}\label{\detokenize{mipstart:copttab-mipstart-set}}
\sphinxaftertopcaption
\begin{tabular}[t]{|\X{5}{15}|\X{10}{15}|}
\sphinxtoprule
\sphinxstyletheadfamily 
\sphinxAtStartPar
API
&\sphinxstyletheadfamily 
\sphinxAtStartPar
Function
\\
\sphinxmidrule
\sphinxtableatstartofbodyhook
\sphinxAtStartPar
C
&
\sphinxAtStartPar
\sphinxcode{\sphinxupquote{COPT\_AddMipStart}}
\\
\sphinxhline
\sphinxAtStartPar
C++
&
\sphinxAtStartPar
\sphinxcode{\sphinxupquote{Model::SetMipStart()}}
\\
\sphinxhline
\sphinxAtStartPar
C\#
&
\sphinxAtStartPar
\sphinxcode{\sphinxupquote{Model.SetMipStart()}}
\\
\sphinxhline
\sphinxAtStartPar
Java
&
\sphinxAtStartPar
\sphinxcode{\sphinxupquote{Model.setMipStart()}}
\\
\sphinxhline
\sphinxAtStartPar
Python
&
\sphinxAtStartPar
\sphinxcode{\sphinxupquote{Model.setMipStart(vars, startvals)}}
\\
\sphinxbottomrule
\end{tabular}
\sphinxtableafterendhook\par
\sphinxattableend\end{savenotes}

\begin{sphinxadmonition}{note}{Note}
\begin{itemize}
\item {} 
\sphinxAtStartPar
Regarding the operations of MIP starts, their \sphinxstylestrong{function names}, \sphinxstylestrong{calling methods}, and \sphinxstylestrong{parameter names} are slightly different in different programming interfaces, but the implementation of functions and meanings of parameter are the same.

\item {} 
\sphinxAtStartPar
For more details on setting initial solutions in C, please refer to funtion \sphinxcode{\sphinxupquote{COPT\_AddMipStart}} in chapter {\hyperref[\detokenize{capiref:chapapi-mipstart}]{\sphinxcrossref{\DUrole{std,std-ref}{C API Function: MIP start utilities}}}}.

\item {} 
\sphinxAtStartPar
You may want to call this method several times to input the MIP start. Please call \sphinxcode{\sphinxupquote{loadMipStart()}} once when the input is done.

\end{itemize}
\end{sphinxadmonition}

\subsection{Read and Write MIP Starts}
\label{\detokenize{mipstart:read-and-write-mip-starts}}
\sphinxAtStartPar
COPT provides functions for file read/write. It can read variable values from a MIP start file (\sphinxcode{\sphinxupquote{".mst"}}) as the initial solution values of variables, and write the solving results or existing initial solution to a MIP start file (\sphinxcode{\sphinxupquote{".mst"}}). The functions for reading and writing MIP start file in different APIs are shown in \hyperref[\detokenize{mipstart:copttab-mipstart-read}]{Table \ref{\detokenize{mipstart:copttab-mipstart-read}}}:

\begin{savenotes}\sphinxattablestart
\sphinxthistablewithglobalstyle
\centering
\sphinxcapstartof{table}
\sphinxthecaptionisattop
\sphinxcaption{Functions for reading and writing MIP starts}\label{\detokenize{mipstart:copttab-mipstart-read}}
\sphinxaftertopcaption
\begin{tabular}[t]{|\X{5}{25}|\X{10}{25}|\X{10}{25}|}
\sphinxtoprule
\sphinxstyletheadfamily 
\sphinxAtStartPar
API
&\sphinxstyletheadfamily 
\sphinxAtStartPar
read MIP starts
&\sphinxstyletheadfamily 
\sphinxAtStartPar
write MIP starts
\\
\sphinxmidrule
\sphinxtableatstartofbodyhook
\sphinxAtStartPar
C
&
\sphinxAtStartPar
\sphinxcode{\sphinxupquote{COPT\_ReadMst}}
&
\sphinxAtStartPar
\sphinxcode{\sphinxupquote{COPT\_WriteMst}}
\\
\sphinxhline
\sphinxAtStartPar
C++
&
\sphinxAtStartPar
\sphinxcode{\sphinxupquote{Model::ReadMst}}
&
\sphinxAtStartPar
\sphinxcode{\sphinxupquote{Model::WriteMst}}
\\
\sphinxhline
\sphinxAtStartPar
C\#
&
\sphinxAtStartPar
\sphinxcode{\sphinxupquote{Model.ReadMst()}}
&
\sphinxAtStartPar
\sphinxcode{\sphinxupquote{Model.WriteMst()}}
\\
\sphinxhline
\sphinxAtStartPar
Java
&
\sphinxAtStartPar
\sphinxcode{\sphinxupquote{Model.readMst()}}
&
\sphinxAtStartPar
\sphinxcode{\sphinxupquote{Model.writeMst()}}
\\
\sphinxhline
\sphinxAtStartPar
Python
&
\sphinxAtStartPar
\sphinxcode{\sphinxupquote{Model.readMst()}}
&
\sphinxAtStartPar
\sphinxcode{\sphinxupquote{Model.writeMst()}}
\\
\sphinxbottomrule
\end{tabular}
\sphinxtableafterendhook\par
\sphinxattableend\end{savenotes}

\section{Related Parameters}
\label{\detokenize{mipstart:related-parameters}}
\sphinxAtStartPar
COPT provides the following parameters that control how MIP starts (initial solutions) are handled inside.

\phantomsection\label{\detokenize{mipstart:mipstartmode}}\begin{itemize}
\item {} 
\sphinxAtStartPar
\sphinxcode{\sphinxupquote{MipStartMode}}
\begin{quote}

\sphinxAtStartPar
Integer parameter.

\sphinxAtStartPar
Mode of MIP starts, i.e. how MIP starts are handled.

\sphinxAtStartPar
\sphinxstylestrong{Default} \sphinxhyphen{}1

\sphinxAtStartPar
\sphinxstylestrong{Possible values:}
\begin{quote}

\sphinxAtStartPar
\sphinxhyphen{}1: Automatic.

\sphinxAtStartPar
0: Do not use any MIP starts.

\sphinxAtStartPar
1: Only load full and feasible MIP starts.

\sphinxAtStartPar
2: Only load feasible ones (complete partial solutions by solving subMIPs).
\end{quote}
\end{quote}

\end{itemize}

\sphinxAtStartPar
\sphinxstylestrong{Note:} If the provided initial solution is incomplete (partial), \sphinxcode{\sphinxupquote{MipStartMode=2}} needs to be set, otherwise the initial solution will be rejected.

\phantomsection\label{\detokenize{mipstart:mipstartnodelimit}}\begin{itemize}
\item {} 
\sphinxAtStartPar
\sphinxcode{\sphinxupquote{MipStartNodeLimit}}
\begin{quote}

\sphinxAtStartPar
Integer parameter.

\sphinxAtStartPar
Limit of nodes for MIP start sub\sphinxhyphen{}MIPs.

\sphinxAtStartPar
\sphinxstylestrong{Default:} \sphinxhyphen{}1

\sphinxAtStartPar
\sphinxstylestrong{Minimal:} \sphinxhyphen{}1

\sphinxAtStartPar
\sphinxstylestrong{Maximal:} \sphinxcode{\sphinxupquote{INT\_MAX}}
\end{quote}

\end{itemize}

\section{Log of MIP starts}
\label{\detokenize{mipstart:log-of-mip-starts}}

\subsection{MIP starts are accepted}
\label{\detokenize{mipstart:mip-starts-are-accepted}}
\sphinxAtStartPar
1.A (better) initial solution was provided

\begin{sphinxVerbatim}[commandchars=\\\{\}]
Initial MIP solution \PYGZsh{} 1 with objective value 9.73987 was accepted
\end{sphinxVerbatim}

\sphinxAtStartPar
2.A partial initial solution was provided, and set \sphinxcode{\sphinxupquote{MipStartMode=2}}. (complete it by solving a subMIP)

\begin{sphinxVerbatim}[commandchars=\\\{\}]
Loading 1 initial MIP solution
Extending partial MIP solution \PYGZsh{} 1
Extending partial MIP solution \PYGZsh{} 1 succeed (0.2s)
Initial MIP solution \PYGZsh{} 1 with objective value 9.66566 was accepted
\end{sphinxVerbatim}

\subsection{MIP starts are rejected}
\label{\detokenize{mipstart:mip-starts-are-rejected}}
\sphinxAtStartPar
\sphinxstylestrong{1.The provided initial solution was infeasible}

\begin{sphinxVerbatim}[commandchars=\\\{\}]
Initial MIP solution \PYGZsh{} 1 was rejected: Primal Inf 1.00e+00  Int Inf 1.78e\PYGZhy{}15
\end{sphinxVerbatim}

\sphinxAtStartPar
\sphinxstylestrong{2.The provided initial solution was  not better  than the current best one}

\begin{sphinxVerbatim}[commandchars=\\\{\}]
Initial MIP solution \PYGZsh{} 2 with objective value 10.3312 was rejected (not better)
\end{sphinxVerbatim}

\sphinxAtStartPar
\sphinxstylestrong{3.The provided initial solution was incomplete (partial), and not set MipStartMode=2}

\begin{sphinxVerbatim}[commandchars=\\\{\}]
Loading 1 initial MIP solution
Initial MIP solution \PYGZsh{} 1 was rejected: partial
\end{sphinxVerbatim}

\sphinxAtStartPar
\sphinxstylestrong{4.The provided initial solution was incomplete (partial), and COPT failed to find a feasible solution by solving subMIP}

\begin{sphinxVerbatim}[commandchars=\\\{\}]
Loading 1 initial MIP solution
Extending partial MIP solution \PYGZsh{} 1
Extending partial MIP solution \PYGZsh{} 1 failed (infeasible)
Initial MIP solution \PYGZsh{} 1 was rejected: partial
\end{sphinxVerbatim}

\sphinxstepscope

\chapter{MIP Solution Pool}
\label{\detokenize{solutionpool:mip-solution-pool}}\label{\detokenize{solutionpool:chapsolutionpool}}\label{\detokenize{solutionpool::doc}}
\sphinxAtStartPar
In general, the solver finds multiple feasible solutions in the process of solving a MIP problem with Branch\sphinxhyphen{}and\sphinxhyphen{}cut method. COPT provides a solution pool for MIP problem, from which users can obtain solutions and the corresponding objective function values. Supported optimization problem types are: MILP, MISOCP, MIQ(C)P.

\sphinxAtStartPar
COPT provides functions that users can get the \sphinxcode{\sphinxupquote{iSol}} th solution’s objective function value and solution values (of specified variables) by specifying the following parameters.
\begin{itemize}
\item {} 
\sphinxAtStartPar
\sphinxcode{\sphinxupquote{iSol}}: Index of the solution to obtain. (0\sphinxhyphen{}based)

\item {} 
\sphinxAtStartPar
\sphinxcode{\sphinxupquote{vars}}: Variables

\end{itemize}

\sphinxAtStartPar
The functions in different APIs are shown in \hyperref[\detokenize{solutionpool:copttab-solutionpool}]{Table \ref{\detokenize{solutionpool:copttab-solutionpool}}}:

\begin{savenotes}\sphinxattablestart
\sphinxthistablewithglobalstyle
\centering
\sphinxcapstartof{table}
\sphinxthecaptionisattop
\sphinxcaption{Get solutions and objective values from solution pool}\label{\detokenize{solutionpool:copttab-solutionpool}}
\sphinxaftertopcaption
\begin{tabular}[t]{|\X{10}{70}|\X{30}{70}|\X{30}{70}|}
\sphinxtoprule
\sphinxstyletheadfamily 
\sphinxAtStartPar
API
&\sphinxstyletheadfamily 
\sphinxAtStartPar
Get solution
&\sphinxstyletheadfamily 
\sphinxAtStartPar
Get objective value
\\
\sphinxmidrule
\sphinxtableatstartofbodyhook
\sphinxAtStartPar
C
&
\sphinxAtStartPar
\sphinxcode{\sphinxupquote{COPT\_GetSolution}}
&
\sphinxAtStartPar
\sphinxcode{\sphinxupquote{COPT\_GetPoolObjVal}}
\\
\sphinxhline
\sphinxAtStartPar
C++
&
\sphinxAtStartPar
\sphinxcode{\sphinxupquote{Model::GetPoolSolution()}}
&
\sphinxAtStartPar
\sphinxcode{\sphinxupquote{Model::GetPoolObjVal()}}
\\
\sphinxhline
\sphinxAtStartPar
C\#
&
\sphinxAtStartPar
\sphinxcode{\sphinxupquote{Model.GetPoolSolution()}}
&
\sphinxAtStartPar
\sphinxcode{\sphinxupquote{Model.GetPoolObjVal()}}
\\
\sphinxhline
\sphinxAtStartPar
Java
&
\sphinxAtStartPar
\sphinxcode{\sphinxupquote{Model.getPoolSolution()}}
&
\sphinxAtStartPar
\sphinxcode{\sphinxupquote{Model.getPoolObjVal()}}
\\
\sphinxhline
\sphinxAtStartPar
Python
&
\sphinxAtStartPar
\sphinxcode{\sphinxupquote{Model.getPoolSolution()}}
&
\sphinxAtStartPar
\sphinxcode{\sphinxupquote{Model.getPoolObjVal()}}
\\
\sphinxbottomrule
\end{tabular}
\sphinxtableafterendhook\par
\sphinxattableend\end{savenotes}

\sphinxAtStartPar
\sphinxstylestrong{Note:} Regarding the operations of solution pool, their \sphinxstylestrong{function names}, \sphinxstylestrong{calling methods}, and \sphinxstylestrong{parameter names} are slightly different in different programming interfaces, but the implementation of functions and meanings of parameter are the same.

\sphinxAtStartPar
\sphinxstylestrong{Attributes of Solution Pool}

\phantomsection\label{\detokenize{solutionpool:poolsols}}\begin{itemize}
\item {} 
\sphinxAtStartPar
\sphinxcode{\sphinxupquote{PoolSols}}
\begin{quote}

\sphinxAtStartPar
Integer attribute

\sphinxAtStartPar
Number of solutions in the solution pool.
\end{quote}

\end{itemize}

\sphinxstepscope

\chapter{COPT Tuner}
\label{\detokenize{tuner:copt-tuner}}\label{\detokenize{tuner:chaptuner}}\label{\detokenize{tuner::doc}}

\section{Introduction}
\label{\detokenize{tuner:introduction}}
\sphinxAtStartPar
The COPT Tuner is a tool designed for tuning performance automatically for all supported problem types.
\begin{itemize}
\item {} 
\sphinxAtStartPar
For MIP problems, it supports tuning for solving time, relative gap, best objective value and objective bound;

\item {} 
\sphinxAtStartPar
For non\sphinxhyphen{}MIP problems, only solving time supported.

\end{itemize}

\sphinxAtStartPar
The workflow of the COPT tuning tool is as follows:
\begin{enumerate}
\sphinxsetlistlabels{\arabic}{enumi}{enumii}{}{.}%
\item {} 
\sphinxAtStartPar
First, perform benchmark calculation and allow users to customize benchmark calculation parameters;

\item {} 
\sphinxAtStartPar
Next, generate tuning parameters one by one, and find parameter combinations that improve the solution performance through parameter tuning calculations.

\end{enumerate}

\section{Related parameters}
\label{\detokenize{tuner:related-parameters}}
\sphinxAtStartPar
Firstly, the tuner will do a baseline run, possibly with fixed parameters from users, then move to the improvement run, where the tuner will generate trial parameter sets and search for parameters that can improve the performance. To summarize, the COPT Tuner provides the following capabilities:

\begin{savenotes}\sphinxattablestart
\sphinxthistablewithglobalstyle
\centering
\sphinxcapstartof{table}
\sphinxthecaptionisattop
\sphinxcaption{Tuner related parameters}\label{\detokenize{tuner:id1}}
\sphinxaftertopcaption
\begin{tabular}[t]{|\X{15}{59}|\X{9}{59}|\X{35}{59}|}
\sphinxtoprule
\sphinxstyletheadfamily 
\sphinxAtStartPar
Name
&\sphinxstyletheadfamily 
\sphinxAtStartPar
Type
&\sphinxstyletheadfamily 
\sphinxAtStartPar
Description
\\
\sphinxmidrule
\sphinxtableatstartofbodyhook
\sphinxAtStartPar
{\hyperref[\detokenize{parameter:tunetimelimit}]{\sphinxcrossref{\DUrole{std,std-ref}{TuneTimeLimit}}}}
&
\sphinxAtStartPar
Double
&
\sphinxAtStartPar
Time limit for parameter tuning
\\
\sphinxhline
\sphinxAtStartPar
{\hyperref[\detokenize{parameter:tunetargettime}]{\sphinxcrossref{\DUrole{std,std-ref}{TuneTargetTime}}}}
&
\sphinxAtStartPar
Double
&
\sphinxAtStartPar
Time target for parameter tuning
\\
\sphinxhline
\sphinxAtStartPar
{\hyperref[\detokenize{parameter:tunetargetrelgap}]{\sphinxcrossref{\DUrole{std,std-ref}{TuneTargetRelGap}}}}
&
\sphinxAtStartPar
Double
&
\sphinxAtStartPar
Optimal relative tolerance target for parameter tuning
\\
\sphinxhline
\sphinxAtStartPar
{\hyperref[\detokenize{parameter:tunemethod}]{\sphinxcrossref{\DUrole{std,std-ref}{TuneMethod}}}}
&
\sphinxAtStartPar
Integer
&
\sphinxAtStartPar
Method for parameter tuning
\\
\sphinxhline
\sphinxAtStartPar
{\hyperref[\detokenize{parameter:tunemode}]{\sphinxcrossref{\DUrole{std,std-ref}{TuneMode}}}}
&
\sphinxAtStartPar
Integer
&
\sphinxAtStartPar
Mode for parameter tuning
\\
\sphinxhline
\sphinxAtStartPar
{\hyperref[\detokenize{parameter:tunemeasure}]{\sphinxcrossref{\DUrole{std,std-ref}{TuneMeasure}}}}
&
\sphinxAtStartPar
Integer
&
\sphinxAtStartPar
Parameter tuning result calculation method
\\
\sphinxhline
\sphinxAtStartPar
{\hyperref[\detokenize{parameter:tunepermutes}]{\sphinxcrossref{\DUrole{std,std-ref}{TunePermutes}}}}
&
\sphinxAtStartPar
Integer
&
\sphinxAtStartPar
Permutations for each trial parameter set
\\
\sphinxhline
\sphinxAtStartPar
{\hyperref[\detokenize{parameter:tuneoutputlevel}]{\sphinxcrossref{\DUrole{std,std-ref}{TuneOutputLevel}}}}
&
\sphinxAtStartPar
Integer
&
\sphinxAtStartPar
Parameter tuning log output intensity
\\
\sphinxbottomrule
\end{tabular}
\sphinxtableafterendhook\par
\sphinxattableend\end{savenotes}

\section{Provided capabilities}
\label{\detokenize{tuner:provided-capabilities}}
\sphinxAtStartPar
The COPT tuning tool provides the following capabilities:

\subsection{Tuning method}
\label{\detokenize{tuner:tuning-method}}
\sphinxAtStartPar
Controlled by the parameter \sphinxcode{\sphinxupquote{TuneMethod}}, options are: greedy search and aggressive search. The greedy method tries to find better parameter settings within limited number of trials, while the aggressive method search for more combinations and has much larger search space than the greedy one, and can potentially find even better parameter settings at the expense of more elapsed tuning time. Default setting is to choose automatically.
\begin{itemize}
\item {} 
\sphinxAtStartPar
Greedy search strategy: It is expected to optimize the calculation with fewer parameters and find better parameter settings;

\item {} 
\sphinxAtStartPar
Broader search strategy: try more parameter combinations, have a larger search space, and are more likely to find better parameter settings, but also consume more tuning time.

\end{itemize}

\sphinxAtStartPar
The possible values and corresponding meanings of the parameter \sphinxcode{\sphinxupquote{TuneMethod}} are as follows. By setting it to a different value, the search method can be selected. The default setting is automatic selection.
\begin{itemize}
\item {} 
\sphinxAtStartPar
\sphinxhyphen{}1: Automatic selection

\item {} 
\sphinxAtStartPar
0: Greedy search strategy

\item {} 
\sphinxAtStartPar
1: Broader search strategy

\end{itemize}

\subsection{Tuning mode}
\label{\detokenize{tuner:tuning-mode}}
\sphinxAtStartPar
Controlled by the parameter \sphinxcode{\sphinxupquote{TuneMode}}, options are: solving time, relative gap, objective value and objective bound. For MIP problem, by default, if the baseline run is not solved to optimality within specified time limit, tuner will change tuning mode to relative gap automatically. Default setting is to choose automatically.
\begin{itemize}
\item {} 
\sphinxAtStartPar
0: Solving time

\item {} 
\sphinxAtStartPar
1: Optimal relative tolerance

\item {} 
\sphinxAtStartPar
2: Objective function value

\item {} 
\sphinxAtStartPar
3: Lower bound of objective function value

\end{itemize}

\sphinxAtStartPar
\sphinxstylestrong{Note:} For integer programming problems, by default, if the benchmark calculation does not optimize the model within the given time limit, the tuning tool will automatically switch the tuning mode to the optimal relative tolerance.

\subsection{Tuning permutations}
\label{\detokenize{tuner:tuning-permutations}}
\sphinxAtStartPar
Controlled by the parameter \sphinxcode{\sphinxupquote{TunePermutes}}. Tuner allow users to run more permutations for each trial parameter set to evaluate performance variability. Default setting is to choose automatically.

\subsection{Tuning measure}
\label{\detokenize{tuner:tuning-measure}}
\sphinxAtStartPar
Controlled by the parameter \sphinxcode{\sphinxupquote{TuneMeasure}}, options are: by average or maximum. When users run more permutations for each trial, tuner will compute the aggregated tuning value by this measure. Default setting is to choose automatically.
\begin{itemize}
\item {} 
\sphinxAtStartPar
0: Calculate the average

\item {} 
\sphinxAtStartPar
1: Calculate the maximum value

\end{itemize}

\subsection{Tuning targets}
\label{\detokenize{tuner:tuning-targets}}
\sphinxAtStartPar
Controlled by the parameter \sphinxcode{\sphinxupquote{TuneTargetTime}} and \sphinxcode{\sphinxupquote{TuneTargetRelGap}}. Tuner enables users to specify target solving time or relative gap for tuning, when tuner finds out parameters that satisfy the specified target, it will stop tuning. For solving time, default value is 0.01 seconds, while for relative gap, default value is 1e\sphinxhyphen{}4.

\subsection{Tuning output}
\label{\detokenize{tuner:tuning-output}}
\sphinxAtStartPar
Controlled by the parameter \sphinxcode{\sphinxupquote{TuneOutputLevel}}, options are: no output, show summary for improved trials, show summary for each trial and show detailed log for each trial. Default setting is to show summary for each trial.
\begin{itemize}
\item {} 
\sphinxAtStartPar
0: Do not output tuning log

\item {} 
\sphinxAtStartPar
1: Output only a summary of the improved parameters

\item {} 
\sphinxAtStartPar
2: Output a summary of each tuning attempt

\item {} 
\sphinxAtStartPar
3: Output a detailed log of each tuning attempt

\end{itemize}

\subsection{TuneTimeLimit}
\label{\detokenize{tuner:tunetimelimit}}
\sphinxAtStartPar
Controlled by the parameter \sphinxcode{\sphinxupquote{TuneTimeLimit}}. This parameter is used to control the overall time limit for the improvement run of tuning. Default setting is to choose automatically.

\subsection{User defined parts}
\label{\detokenize{tuner:user-defined-parts}}\begin{itemize}
\item {} 
\sphinxAtStartPar
User defined parameters

\sphinxAtStartPar
The tool enables users to set parameters for the baseline run, which will also be used as fixed parameters for each trial run. Tuner will not tune parameters in the fixed parameters.

\item {} 
\sphinxAtStartPar
User defined MIP start

\sphinxAtStartPar
The COPT tuner enables users to set MIP start for the baseline run, which will also be used for each trial run also.

\item {} 
\sphinxAtStartPar
User defined tuning file

\sphinxAtStartPar
The COPT tuner enables users to read tuning parameter sets from tuning file, if so, the tuner will try to tune from the given parameter sets, otherwise, tuner will generate tuning parameter sets automatically. The tuning file is similar to parameter file, with the difference that it allow multiple values for each parameter name.

\sphinxAtStartPar
\sphinxstylestrong{Note:} The COPT tuning file has a similar format to the COPT parameter file, except that the tuning file allows multiple values to be specified for a single parameter.

\end{itemize}

\subsection{Load or writing tuning parameter}
\label{\detokenize{tuner:load-or-writing-tuning-parameter}}
\sphinxAtStartPar
After the parameter tuning is completed, the number of parameter tuning results can be obtained through the attribute \sphinxcode{\sphinxupquote{TuneResults}}, and the tuning results of the specified number can also be loaded into the model or written into the parameter file.

\sphinxAtStartPar
COPT can output the parameter tuning results of the specified number to the parameter file (\sphinxcode{\sphinxupquote{".par"}}). The parameters that need to be specified are:
\begin{itemize}
\item {} 
\sphinxAtStartPar
\sphinxcode{\sphinxupquote{idx}}: parameter tuning result number

\item {} 
\sphinxAtStartPar
\sphinxcode{\sphinxupquote{filename}}: file name

\end{itemize}

\sphinxAtStartPar
The corresponding functions in different programming interfaces are as follows:

\begin{savenotes}\sphinxattablestart
\sphinxthistablewithglobalstyle
\centering
\sphinxcapstartof{table}
\sphinxthecaptionisattop
\sphinxcaption{functions for writing parameter tuning results in different interfaces}\label{\detokenize{tuner:copttab-tuner}}
\sphinxaftertopcaption
\begin{tabular}[t]{|\X{5}{15}|\X{10}{15}|}
\sphinxtoprule
\sphinxstyletheadfamily 
\sphinxAtStartPar
API
&\sphinxstyletheadfamily 
\sphinxAtStartPar
function
\\
\sphinxmidrule
\sphinxtableatstartofbodyhook
\sphinxAtStartPar
C
&
\sphinxAtStartPar
\sphinxcode{\sphinxupquote{COPT\_WriteTuneParam}}
\\
\sphinxhline
\sphinxAtStartPar
C++
&
\sphinxAtStartPar
\sphinxcode{\sphinxupquote{Model::WriteTuneParam()}}
\\
\sphinxhline
\sphinxAtStartPar
C\#
&
\sphinxAtStartPar
\sphinxcode{\sphinxupquote{Model.WriteTuneParam()}}
\\
\sphinxhline
\sphinxAtStartPar
Java
&
\sphinxAtStartPar
\sphinxcode{\sphinxupquote{Model.writeTuneParam()}}
\\
\sphinxhline
\sphinxAtStartPar
Python
&
\sphinxAtStartPar
\sphinxcode{\sphinxupquote{Model.writeTuneParam()}}
\\
\sphinxbottomrule
\end{tabular}
\sphinxtableafterendhook\par
\sphinxattableend\end{savenotes}

\section{Example}
\label{\detokenize{tuner:example}}
\sphinxAtStartPar
For example, to tune model \sphinxcode{\sphinxupquote{"foo.mps"}} from command line for solving time with COPT command line tool, the commands are:

\begin{sphinxVerbatim}[commandchars=\\\{\}]
copt\PYGZus{}cmd\PYG{+w}{ }\PYGZhy{}c\PYG{+w}{ }\PYG{l+s+s2}{\PYGZdq{}read foo.mps; tune; exit\PYGZdq{}}
\end{sphinxVerbatim}

\sphinxAtStartPar
To use the tuner in API such as Python, the codes are:

\begin{sphinxVerbatim}[commandchars=\\\{\}]
\PYG{n}{env} \PYG{o}{=} \PYG{n}{Envr}\PYG{p}{(}\PYG{p}{)}
\PYG{n}{m} \PYG{o}{=} \PYG{n}{env}\PYG{o}{.}\PYG{n}{createModel}\PYG{p}{(}\PYG{p}{)}
\PYG{n}{m}\PYG{o}{.}\PYG{n}{read}\PYG{p}{(}\PYG{l+s+s2}{\PYGZdq{}}\PYG{l+s+s2}{foo.mps}\PYG{l+s+s2}{\PYGZdq{}}\PYG{p}{)}
\PYG{n}{m}\PYG{o}{.}\PYG{n}{tune}\PYG{p}{(}\PYG{p}{)}
\end{sphinxVerbatim}

\sphinxstepscope

\chapter{Callbacks}
\label{\detokenize{callback:callbacks}}\label{\detokenize{callback:chapcallback}}\label{\detokenize{callback::doc}}
\sphinxAtStartPar
COPT provides the callbacks utility, which supports users in obtaining information during the MIP solving process, e.g., the current best bound, the current optimal objective value, etc.; or controlling the solving process, e.g., by adding lazy constraints and cutting planes, or terminating the solving process. The problem types supporting the use of callbacks are MILP, MISOCP, MIQ(C)P.

\sphinxAtStartPar
A callback function is a user\sphinxhyphen{}provided function called by COPT during the solving process. The user can register one custom callback function via their preferred API for one or multiple callback contexts. Section {\hyperref[\detokenize{callback:chapcallback-api}]{\sphinxcrossref{\DUrole{std,std-ref}{Using the callback utilities in different APIs}}}} gives a detailed introduction to how to setup a callback function. The callback function will be invoked at certain moments during the solving process, depending on the callback contexts. When invoked, the user can {\hyperref[\detokenize{callback:chapcallback-getinfo}]{\sphinxcrossref{\DUrole{std,std-ref}{access information}}}} and {\hyperref[\detokenize{callback:chapcallback-control}]{\sphinxcrossref{\DUrole{std,std-ref}{control the solving process}}}}, respectively. The available information and operations depend on the context. Currently, COPT supports four callback contexts:
\begin{itemize}
\item {} 
\sphinxAtStartPar
\sphinxcode{\sphinxupquote{CBCONTEXT\_INCUMBENT}}: Invokes the callback after a new incumbent was found.

\item {} 
\sphinxAtStartPar
\sphinxcode{\sphinxupquote{CBCONTEXT\_MIPNODE}}: Invokes the callback after a MIP node was processed.

\item {} 
\sphinxAtStartPar
\sphinxcode{\sphinxupquote{CBCONTEXT\_MIPRELAX}} : Invokes the callback when an LP\sphinxhyphen{}relaxation was solved.

\item {} 
\sphinxAtStartPar
\sphinxcode{\sphinxupquote{CBCONTEXT\_MIPSOL}}: Invokes the callback when a new MIP candidate solution is found.

\end{itemize}

\sphinxAtStartPar
The content of this chapter is organized as follows:
\begin{itemize}
\item {} 
\sphinxAtStartPar
{\hyperref[\detokenize{callback:chapcallback-getinfo}]{\sphinxcrossref{\DUrole{std,std-ref}{Obtaining information during the solving process}}}}

\item {} 
\sphinxAtStartPar
{\hyperref[\detokenize{callback:chapcallback-control}]{\sphinxcrossref{\DUrole{std,std-ref}{Controlling the MIP solving process}}}}

\item {} 
\sphinxAtStartPar
{\hyperref[\detokenize{callback:chapcallback-api}]{\sphinxcrossref{\DUrole{std,std-ref}{Using the callback utilities in different APIs}}}}

\end{itemize}

\begin{sphinxadmonition}{note}{Notes}

\sphinxAtStartPar
Only \sphinxstylestrong{one} callback function can be registered in COPT at a time. But one callback can be registered for multiple contexts. If a user wants to call different operations for different contexts (such as adding lazy constraints under \sphinxcode{\sphinxupquote{CBCONTEXT\_MIPSOL}} and adding user cuts under \sphinxcode{\sphinxupquote{CBCONTEXT\_MIPRELAX}} ), they need to register one callback for all relevant contexts and have this callback call the respective operations based on the context in which it was called.
\end{sphinxadmonition}

\section{Obtaining information during the solving process}
\label{\detokenize{callback:obtaining-information-during-the-solving-process}}\label{\detokenize{callback:chapcallback-getinfo}}
\sphinxAtStartPar
The information that can be obtained during the MIP solving process depends on the context the callback is invoked in, see the table below. Information is usually obtained by calling (an API dependent version of) \sphinxcode{\sphinxupquote{getInfo}} / \sphinxcode{\sphinxupquote{GetCallbackInfo}} from within the callback function, specifying the desired information via a string supplied as the function argument. For a detailed description of the available callback information arguments, please refer to {\hyperref[\detokenize{information:chapinfo-cbc}]{\sphinxcrossref{\DUrole{std,std-ref}{Callback information}}}}.

\sphinxAtStartPar
The following table lists information that can be obtained in different contexts:

\begin{savenotes}\sphinxattablestart
\sphinxthistablewithglobalstyle
\centering
\begin{tabular}[t]{|\X{20}{60}|\X{40}{60}|}
\sphinxtoprule
\sphinxstyletheadfamily 
\sphinxAtStartPar
Context
&\sphinxstyletheadfamily 
\sphinxAtStartPar
Callback Information
\\
\sphinxmidrule
\sphinxtableatstartofbodyhook
\sphinxAtStartPar
\sphinxcode{\sphinxupquote{CBCONTEXT\_MIPNODE}}
&
\sphinxAtStartPar
\sphinxcode{\sphinxupquote{NodeStatus}}, \sphinxcode{\sphinxupquote{RelaxSolution}}, \sphinxcode{\sphinxupquote{RelaxSolObj}}, \sphinxcode{\sphinxupquote{MipCandObj}}, \sphinxcode{\sphinxupquote{MipCandidate}}
\\
\sphinxhline
\sphinxAtStartPar
\sphinxcode{\sphinxupquote{CBCONTEXT\_MIPRELAX}}
&
\sphinxAtStartPar
\sphinxcode{\sphinxupquote{RelaxSolution}}, \sphinxcode{\sphinxupquote{RelaxSolObj}}
\\
\sphinxhline
\sphinxAtStartPar
\sphinxcode{\sphinxupquote{CBCONTEXT\_MIPSOL}}
&
\sphinxAtStartPar
\sphinxcode{\sphinxupquote{MipCandObj}}, \sphinxcode{\sphinxupquote{MipCandidate}}
\\
\sphinxhline
\sphinxAtStartPar
\sphinxcode{\sphinxupquote{CBCONTEXT\_INCUMBENT}}
&\\
\sphinxbottomrule
\end{tabular}
\sphinxtableafterendhook\par
\sphinxattableend\end{savenotes}

\sphinxAtStartPar
In addition to the corresponding Callback Context and Information listed above, \sphinxcode{\sphinxupquote{BestObj}}, \sphinxcode{\sphinxupquote{BestBnd}}, \sphinxcode{\sphinxupquote{HasIncumbent}}, and \sphinxcode{\sphinxupquote{Incumbent}} can be obtained in any context.

\begin{sphinxadmonition}{note}{Notes}
\begin{enumerate}
\sphinxsetlistlabels{\arabic}{enumi}{enumii}{}{.}%
\item {} 
\sphinxAtStartPar
If \sphinxcode{\sphinxupquote{HasIncumbent == False}}, then \sphinxcode{\sphinxupquote{Incumbent}} cannot be obtained.

\item {} 
\sphinxAtStartPar
The return value of the “NodeStatus” information is constant, representing the solving status of the
current node’s LP relaxation. For possible values, please refer to {\hyperref[\detokenize{constant:copttab-statuscodes}]{\sphinxcrossref{\DUrole{std,std-ref}{General Constants Section:
Solution Status (Partial)}}}}.

\item {} 
\sphinxAtStartPar
\sphinxcode{\sphinxupquote{Incumbent}}, \sphinxcode{\sphinxupquote{RelaxSolution}} , and  \sphinxcode{\sphinxupquote{MipCandidate}} are obtained through different methods in different interfaces:
\begin{itemize}
\item {} 
\sphinxAtStartPar
C API: through the function \sphinxcode{\sphinxupquote{COPT\_GetCallbackInfo}}, the name of the intermediate information to
be obtained is provided as arguments of the function;

\item {} 
\sphinxAtStartPar
In object\sphinxhyphen{}oriented programming languages (C++/C\#/Java/Python), the \sphinxcode{\sphinxupquote{CallbackBase}} class provides
specialized functions to obtain the corresponding intermediate information. E.g., in Python/C++
\sphinxcode{\sphinxupquote{CallbackBase}} provides \sphinxcode{\sphinxupquote{GetIncumbent}}, \sphinxcode{\sphinxupquote{GetRelaxSol}}, and \sphinxcode{\sphinxupquote{GetSolution}}. Other programming
language interfaces are similar, please refer to the \sphinxcode{\sphinxupquote{CallbackBase}} class of each API.

\end{itemize}

\end{enumerate}
\end{sphinxadmonition}

\section{Controlling the MIP solving process}
\label{\detokenize{callback:controlling-the-mip-solving-process}}\label{\detokenize{callback:chapcallback-control}}
\sphinxAtStartPar
COPT provides functions to allow the user to interactively add lazy constraints or cutting planes during the solving process of the MIP branch\sphinxhyphen{}and\sphinxhyphen{}cut to control the MIP solving process. There are three main types of operations:
\begin{enumerate}
\sphinxsetlistlabels{\arabic}{enumi}{enumii}{}{.}%
\item {} 
\sphinxAtStartPar
Adding lazy constraints

\item {} 
\sphinxAtStartPar
Adding cutting planes

\item {} 
\sphinxAtStartPar
Adding feasible solutions

\end{enumerate}

\subsection{Adding lazy constraints}
\label{\detokenize{callback:adding-lazy-constraints}}
\sphinxAtStartPar
Lazy constraints are constraints that are added to the model only when they are violated. For some models with a large number of constraints, adding lazy constraints only when violated can effectively reduce the size of the model during the solution process and improve the efficiency of the solving process. A popular example of this is the TSP model, see \sphinxcode{\sphinxupquote{"cb\_ex1"}} in the examples directory in the installation package.

\sphinxAtStartPar
COPT supports two ways of adding lazy constraints. One is to explicitly add lazy constraints to the model \sphinxstylestrong{before} starting the solution process. The other is to add lazy constraints \sphinxstylestrong{during} the solving process through a user callback. For this purpose, each API provides two sets of methods, one for adding lazy constraints to the initial model and one for adding them from a callback. In the C API, the methods can be distinguished according to whether or not the function name contains \sphinxcode{\sphinxupquote{"Callback"}}, e.g., \sphinxcode{\sphinxupquote{COPT\_AddLazyConstr}} and \sphinxcode{\sphinxupquote{COPT\_AddCallbackLazyConstr}}. In object\sphinxhyphen{}oriented APIs, the two sets of functions correspond to the \sphinxcode{\sphinxupquote{Model}} class and the \sphinxcode{\sphinxupquote{CallbackBase}} class respectively. Taking python as an example:
\begin{itemize}
\item {} 
\sphinxAtStartPar
Before solving, a user can directly add lazy constraints to the model by calling \sphinxcode{\sphinxupquote{Model.addLazyConstr()}} or \sphinxcode{\sphinxupquote{Model.addLazyConstrs()}} respectively.

\item {} 
\sphinxAtStartPar
During the solving process, a user can dynamically add lazy constraints (if supported by the current context, see below) from within the callback function via \sphinxcode{\sphinxupquote{CallbackBase.addLazyConstr()}} or \sphinxcode{\sphinxupquote{CallbackBase.addLazyConstrs()}}.

\end{itemize}

\sphinxAtStartPar
In both cases, the added lazy constraints will be stored by COPT separately from the actual model and are only added to the model, when they are violated by a solution found during the solving process.

\sphinxAtStartPar
In order to ensure correctness, COPT will check whether any lazy constraints added so far are violated by any solution found during the solving process. This will increase the solving time, especially when many non\sphinxhyphen{}violated lazy constraints have been added. It is recommended that the user only adds lazy constraints when necessary, e.g. when they are violated by a solution.

\sphinxAtStartPar
It should be noted that although lazy constraints can be added in the
\sphinxcode{\sphinxupquote{CBCONTEXT\_MIPRELAX}} context, this is supported only when the
\sphinxcode{\sphinxupquote{CBCONTEXT\_MIPSOL}} context is also enabled. Otherwise, attempts to add
lazy constraints in \sphinxcode{\sphinxupquote{CBCONTEXT\_MIPRELAX}} will result in an error.
Therefore, while it is not strictly necessary to check every LP
relaxation solution for violations of lazy constraints, users must
verify that every solution provided in \sphinxcode{\sphinxupquote{CBCONTEXT\_MIPSOL}} satisfies
all lazy constraints. Failing to do so may result in incorrect solutions
being accepted.

\sphinxAtStartPar
To avoid adding unnecessarily many lazy constraints, COPT has some simple redundancy checks for lazy constraints in place. Exact duplicates will be discarded. However, adding many, very similar but redundant lazy constraints will negatively affect COPT’s performance. This should be avoided by the user.

\begin{sphinxadmonition}{note}{Notes}
\begin{itemize}
\item {} 
\sphinxAtStartPar
Registering a callback function for the \sphinxcode{\sphinxupquote{CBCONTEXT\_MIPSOL}} will make COPT believe that the user wants to add lazy constraints. As lazy constraints are not actually part of the model, this will lead to the deactivation of dual reductions during COPT’s presolve, as dual arguments rely on the knowledge of all model rows. If the user does not intend to add lazy constraints but still wants to use the \sphinxcode{\sphinxupquote{CBCONTEXT\_MIPSOL}}, COPT provides the {\hyperref[\detokenize{parameter:lazyconstraints}]{\sphinxcrossref{\DUrole{std,std-ref}{LazyConstaints}}}} parameter which enables the user to explicitly tell COPT whether or not lazy constraints will be added to the model. By default, this parameter is set to \sphinxcode{\sphinxupquote{\sphinxhyphen{}1}} meaning COPT will turn off dual presolve reductions if either lazy constraints are part of the model or a callback for context \sphinxcode{\sphinxupquote{CBCONTEXT\_MIPSOL}} has been installed. Explicitly setting the parameter to \sphinxcode{\sphinxupquote{0}} will allow dual reductions during COPT’s presolve even if lazy constraints or a callback for context \sphinxcode{\sphinxupquote{CBCONTEXT\_MIPSOL}} are present. This is useful only in very rare cases, e.g., if the callback only prints information about solution candidates but never adds lazy constraints. As soon as lazy constraints are added, this might lead to wrong results, however. For printing information about solutions, consider using the \sphinxcode{\sphinxupquote{CBCONTEXT\_INCUMBENT}} context instead.

\item {} 
\sphinxAtStartPar
If a user invokes a function to add lazy constraints from a callback in the \sphinxcode{\sphinxupquote{CBCONTEXT\_MIPSOL}} context, the current MIP candidate solution will be rejected, no matter whether the added lazy constraint(s) are actually violated or not. This enables the user to reject arbitrary solutions by adding empty lazy constraints when an undesirable solution is found. Note however, that COPT might find the same solution multiple times if no lazy constraint is provided. The LP relaxation solution will not necessarily be rejected if lazy constraints are added in \sphinxcode{\sphinxupquote{CBCONTEXT\_MIPRELAX}}, only when these are actually violated.

\item {} 
\sphinxAtStartPar
It is invalid to call the any functions of the \sphinxcode{\sphinxupquote{Model}} class for object\sphinxhyphen{}oriented languages (or their C equivalent) to add lazy constraints in a callback. More generally, the model cannot be changed during the solving process, except by adding lazy constraints or cutting planes.

\end{itemize}
\end{sphinxadmonition}

\subsection{Adding cutting planes}
\label{\detokenize{callback:adding-cutting-planes}}
\sphinxAtStartPar
Cutting planes are added to the model during the solving process to strengthen the LP relaxation, e.g., cut off fractional LP solutions and improve the lower bound of the MIP problem.

\sphinxAtStartPar
COPT supports the addition of custom cutting planes to the model during the solving process. Similar to lazy constraints, cutting planes can be added to the model \sphinxstylestrong{before} and, via the callback, \sphinxstylestrong{during} the solving process. Each API provides two sets of methods, one for adding cutting planes to the initial model and one for adding them from a callback. In the C API, the methods can be distinguished according to whether or not the function name contains \sphinxcode{\sphinxupquote{"Callback"}}, e.g., \sphinxcode{\sphinxupquote{COPT\_AddUserCut}} and \sphinxcode{\sphinxupquote{COPT\_AddCallbackUserCut}}. In object\sphinxhyphen{}oriented APIs, the two sets of functions correspond to the \sphinxcode{\sphinxupquote{Model}} class and the \sphinxcode{\sphinxupquote{CallbackBase}} class respectively. Taking python as an example:
\begin{itemize}
\item {} 
\sphinxAtStartPar
Before solving, user can directly add cutting planes to the model by calling  \sphinxcode{\sphinxupquote{Model.addUserCut()}} or \sphinxcode{\sphinxupquote{Model.addUserCuts()}}.

\item {} 
\sphinxAtStartPar
During the solving process, a user can dynamically add cutting planes (if supported by the current context, see below) from within the callback function via \sphinxcode{\sphinxupquote{CallbackBase.addUserCut()}} or \sphinxcode{\sphinxupquote{CallbackBase.addUserCuts()}}.

\end{itemize}

\sphinxAtStartPar
Cutting planes can \sphinxstylestrong{only} be added in the \sphinxcode{\sphinxupquote{CBCONTEXT\_MIPRELAX}} context. Here, the user is provided with the current LP relaxation solution to separate their own cutting planes.

\begin{sphinxadmonition}{note}{Notes}
\begin{itemize}
\item {} 
\sphinxAtStartPar
Cutting planes that do not violate the current LP relaxation solution are discarded by COPT.

\item {} 
\sphinxAtStartPar
It is invalid to call any functions of the \sphinxcode{\sphinxupquote{Model}} class (or their C equivalent) to add cutting planes in a callback. More generally, the model cannot be changed during  the solution process, except by adding lazy constraints or cutting planes.

\end{itemize}
\end{sphinxadmonition}

\subsection{Adding feasible solutions}
\label{\detokenize{callback:adding-feasible-solutions}}
\sphinxAtStartPar
COPT supports adding feasible solutions during the MIP solving process. This enables the user to provide any feasible solution they found in parallel to the COPT solution process, e.g., in a self\sphinxhyphen{}implemented heuristic. Known solutions can either be supplied as starting solutions by calling \sphinxcode{\sphinxupquote{COPT\_AddMipStart}} (see {\hyperref[\detokenize{mipstart:chapmipstart}]{\sphinxcrossref{\DUrole{std,std-ref}{MIP Starts}}}}) or from within a callback function. If a solution is know beforehand, supplying it as a MIP starting solution is preferred. For solutions found during the solving process, in the C API, a solution can be added by calling \sphinxcode{\sphinxupquote{COPT\_AddCallbackSolution}} inside the callback. In object\sphinxhyphen{}oriented APIs, the functions needed to add a solution are provided by the \sphinxcode{\sphinxupquote{CallbackBase}} class and the workflow consists of calling two functions:
\begin{itemize}
\item {} 
\sphinxAtStartPar
Set the feasible solution: \sphinxcode{\sphinxupquote{CallbackBase.setSolution(vars, val)}}

\item {} 
\sphinxAtStartPar
Load a custom solution into the model: \sphinxcode{\sphinxupquote{CallbackBase.loadSolution()}}

\end{itemize}

\sphinxAtStartPar
COPT will check any provided solution for feasibility and compute its objective value. The computed objective is return by \sphinxcode{\sphinxupquote{loadSolution}} / \sphinxcode{\sphinxupquote{COPT\_AddCallbackSolution}}. If a solution is infeasible or worse than the current incumbent, it is discarded and the objective value returned by COPT is set to \sphinxcode{\sphinxupquote{1.0e+30}}.

\sphinxAtStartPar
Solutions can be added in any callback context.

\begin{sphinxadmonition}{note}{Notes}

\sphinxAtStartPar
Currently, COPT only supports complete feasible solutions within callbacks.
\end{sphinxadmonition}

\section{Using the callback utilities in different APIs}
\label{\detokenize{callback:using-the-callback-utilities-in-different-apis}}\label{\detokenize{callback:chapcallback-api}}
\sphinxAtStartPar
For object oriented programming languages the basic steps of setting up a callback function are:
\begin{enumerate}
\sphinxsetlistlabels{\arabic}{enumi}{enumii}{}{.}%
\item {} 
\sphinxAtStartPar
Implement a custom callback class, inheriting from the \sphinxcode{\sphinxupquote{CallbackBase}} class.

\item {} 
\sphinxAtStartPar
Implement the \sphinxcode{\sphinxupquote{CallbackBase.callback()}} function. This is the callback function that will be invoked by COPT. Here, the user can invoke the callback\sphinxhyphen{}specific functions for {\hyperref[\detokenize{callback:chapcallback-getinfo}]{\sphinxcrossref{\DUrole{std,std-ref}{obtaining information}}}} or {\hyperref[\detokenize{callback:chapcallback-control}]{\sphinxcrossref{\DUrole{std,std-ref}{controlling the solution process}}}}.

\item {} 
\sphinxAtStartPar
Create an object of the custom callback class.

\item {} 
\sphinxAtStartPar
Register the callback in COPT through \sphinxcode{\sphinxupquote{Model.setCallback()}}, and input the Callback Context as a parameter. For registering the callback function for multiple contexts, one can bitwise\sphinxhyphen{}or the desired contexts, e.g., \sphinxcode{\sphinxupquote{COPT.CBCONTEXT\_MIPSOL | COPT.CBCONTEXT\_MIPNODE}}.

\end{enumerate}

\sphinxAtStartPar
In the subsequent solving process, the user\sphinxhyphen{}supplied \sphinxcode{\sphinxupquote{CallbackBase.callback()}} function will be called in each context registered. The currently invoked context can be obtained by calling the \sphinxcode{\sphinxupquote{CallbackBase}} class method \sphinxcode{\sphinxupquote{where()}}.

\sphinxAtStartPar
As already mentioned in the sections above, functions in \sphinxcode{\sphinxupquote{CallbackBase}} class (or their corresponding C functions) can only be called in certain contexts. The following table lists for each callback context the allowed callback operations, using the Python API:

\begin{savenotes}\sphinxattablestart
\sphinxthistablewithglobalstyle
\centering
\begin{tabular}[t]{|\X{20}{55}|\X{35}{55}|}
\sphinxtoprule
\sphinxstyletheadfamily 
\sphinxAtStartPar
Context
&\sphinxstyletheadfamily 
\sphinxAtStartPar
Function
\\
\sphinxmidrule
\sphinxtableatstartofbodyhook
\sphinxAtStartPar
CBCONTEXT\_INCUMBENT
&
\sphinxAtStartPar
\sphinxcode{\sphinxupquote{getInfo}}, \sphinxcode{\sphinxupquote{getIncumbent}}, \sphinxcode{\sphinxupquote{load/setSolution}}
\\
\sphinxhline
\sphinxAtStartPar
CBCONTEXT\_MIPNODE
&
\sphinxAtStartPar
\sphinxcode{\sphinxupquote{getInfo}}, \sphinxcode{\sphinxupquote{getIncumbent}}, \sphinxcode{\sphinxupquote{getRelaxSol}}, \sphinxcode{\sphinxupquote{load/setSolution}}
\\
\sphinxhline
\sphinxAtStartPar
CBCONTEXT\_MIPRELAX
&
\sphinxAtStartPar
\sphinxcode{\sphinxupquote{addUserCut(s)}}, \sphinxcode{\sphinxupquote{getIncumbent}}, \sphinxcode{\sphinxupquote{getInfo}}, \sphinxcode{\sphinxupquote{getRelaxSol}}, \sphinxcode{\sphinxupquote{load/setSolution}}
\\
\sphinxhline
\sphinxAtStartPar
CBCONTEXT\_MIPSOL
&
\sphinxAtStartPar
\sphinxcode{\sphinxupquote{addLazyConstr(s)}}, \sphinxcode{\sphinxupquote{getIncumbent}}, \sphinxcode{\sphinxupquote{getInfo}}, \sphinxcode{\sphinxupquote{getSolution}}, \sphinxcode{\sphinxupquote{load/setSolution}}
\\
\sphinxbottomrule
\end{tabular}
\sphinxtableafterendhook\par
\sphinxattableend\end{savenotes}

\begin{sphinxadmonition}{note}{Notes}
\begin{itemize}
\item {} 
\sphinxAtStartPar
While \sphinxcode{\sphinxupquote{getInfo}} can be called in all contexts, the information available depends on the context. See {\hyperref[\detokenize{callback:chapcallback-getinfo}]{\sphinxcrossref{\DUrole{std,std-ref}{Obtaining information during the solving process}}}} for details.

\item {} 
\sphinxAtStartPar
In other APIs \sphinxcode{\sphinxupquote{getInfo}} is often split into \sphinxcode{\sphinxupquote{getIntInfo}} and \sphinxcode{\sphinxupquote{getDblInfo}}.

\item {} 
\sphinxAtStartPar
The functions above may look slightly different for a certain API, but the presented relationships are the same.

\end{itemize}
\end{sphinxadmonition}

\sphinxAtStartPar
While the above steps use Python as the reference API, the implementation in each object oriented programming language API is similar, and the user can refer to the provided sample code. For Python, \sphinxcode{\sphinxupquote{"cb\_ex1.py"}} is available in the examples directory in the installation package. For C, the main differences when implementing and registering a callback are as follows:
\begin{itemize}
\item {} 
\sphinxAtStartPar
The custom callback function can be any function using the signature \sphinxcode{\sphinxupquote{int COPT\_CALL \textless{}function\textgreater{}(copt\_prob* prob, void* cbdata, int cbctx, void* usrdata)}} where \sphinxcode{\sphinxupquote{\textless{}function\textgreater{}}} is arbitrary.

\item {} 
\sphinxAtStartPar
Instead of using \sphinxcode{\sphinxupquote{Where()}} for obtaining the current context, the callback context is supplied in \sphinxcode{\sphinxupquote{cbctx}}.

\item {} 
\sphinxAtStartPar
Callback\sphinxhyphen{}relevant information can be passed by defining a custom \sphinxcode{\sphinxupquote{struct}} and passing it as the \sphinxcode{\sphinxupquote{usrdata}} argument.

\end{itemize}

\sphinxAtStartPar
See \sphinxcode{\sphinxupquote{cb\_ex1.c}} in the C examples folder for a reference implementation.

\sphinxAtStartPar
The calling method and function name of the callback function in different programming interfaces are \sphinxstylestrong{slightly different}, but the supported functions and function meanings are the same. Please refer to the corresponding chapters of different programming interface API reference manuals for specific introductions:
\begin{itemize}
\item {} 
\sphinxAtStartPar
C API: {\hyperref[\detokenize{capiref:chapapi-cbc}]{\sphinxcrossref{\DUrole{std,std-ref}{Callback utilities}}}}

\item {} 
\sphinxAtStartPar
C++ API: {\hyperref[\detokenize{cppapiref:chapcppapiref-callbackbase}]{\sphinxcrossref{\DUrole{std,std-ref}{CallbackBase Class}}}}

\item {} 
\sphinxAtStartPar
C\# API: {\hyperref[\detokenize{csharpapiref:chapcsharpapiref-callback}]{\sphinxcrossref{\DUrole{std,std-ref}{CallbackBase Class}}}}

\item {} 
\sphinxAtStartPar
Java API: {\hyperref[\detokenize{javaapiref:chapjavaapiref-callback}]{\sphinxcrossref{\DUrole{std,std-ref}{CallbackBase Class}}}}

\item {} 
\sphinxAtStartPar
Python API: {\hyperref[\detokenize{pyapiref:chappyapi-cbcbase}]{\sphinxcrossref{\DUrole{std,std-ref}{CallbackBase Class}}}}

\end{itemize}

\sphinxstepscope

\chapter{Matrix Modeling}
\label{\detokenize{matrix:matrix-modeling}}\label{\detokenize{matrix:chapmatrix}}\label{\detokenize{matrix::doc}}
\sphinxAtStartPar
The COPT Python API provides matrix modeling, supports \sphinxcode{\sphinxupquote{NumPy}} multi\sphinxhyphen{}dimensional array, a two\sphinxhyphen{}dimensional NumPy matrix, SciPy compressed sparse column matrix ( \sphinxcode{\sphinxupquote{csc\_matrix}} ) and compressed sparse row matrix ( \sphinxcode{\sphinxupquote{csr\_matrix}} ) operations and can be combined with ordinary (scalar) variables and constraints.

\sphinxAtStartPar
\sphinxcode{\sphinxupquote{NumPy}} version should be 1.23 or above and \sphinxcode{\sphinxupquote{Python}} minimum version requirement is 3.8. COPT mainly provides the following utilities:
\begin{enumerate}
\sphinxsetlistlabels{\arabic}{enumi}{enumii}{}{.}%
\item {} 
\sphinxAtStartPar
Add multi\sphinxhyphen{}dimensional variables ( \sphinxcode{\sphinxupquote{MVar}} ) and other related operations;

\item {} 
\sphinxAtStartPar
Construct multi\sphinxhyphen{}dimensional linear expressions ( \sphinxcode{\sphinxupquote{MLinExpr}} ), add multi\sphinxhyphen{}dimensional linear constraints ( \sphinxcode{\sphinxupquote{MConstr}} ) and other related operations;

\item {} 
\sphinxAtStartPar
Construct multi\sphinxhyphen{}dimensional quadratic expression ( \sphinxcode{\sphinxupquote{MQuadExpr}} ), add multi\sphinxhyphen{}dimensional convex quadratic constraint ( \sphinxcode{\sphinxupquote{QConstraint}} ) and other related operations.

\end{enumerate}

\section{Two different matrix modeling modes}
\label{\detokenize{matrix:two-different-matrix-modeling-modes}}
\sphinxAtStartPar
The matrix modeling function provided by COPT currently supports two modes: original legacy mode and experimental mode.
Among them, the legacy mode relies on Python’s \sphinxcode{\sphinxupquote{NumPy}} library.
The experimental mode (by default) is based on the built\sphinxhyphen{}in matrix modeling class implementation of the COPT C++ interface.
The modeling speed could be improved to a certain extent compared with the former.

\sphinxAtStartPar
In the Python interface, the two modes can be controlled and switched through \sphinxcode{\sphinxupquote{Model.matrixmodelmode}} , and the default is the experimental mode.
Please set according to your usage scenarios:
\begin{itemize}
\item {} 
\sphinxAtStartPar
\sphinxcode{\sphinxupquote{Model.matrixmodelmode = "legacy"}} , depends on the \sphinxcode{\sphinxupquote{NumPy}} library.

\item {} 
\sphinxAtStartPar
\sphinxcode{\sphinxupquote{Model.matrixmodelmode = "experimental"}} (by default), C++ built\sphinxhyphen{}in matrix modeling class based on COPT ( \sphinxcode{\sphinxupquote{NdArray}} ) , no external dependencies.

\end{itemize}

\sphinxAtStartPar
The functions and operation types supported by the two modes are slightly different.

\sphinxAtStartPar
1.The order (storage) method of \sphinxcode{\sphinxupquote{MConstr.reshape()/MVar.reshape()}} , that is, the value of the function parameter \sphinxcode{\sphinxupquote{order}}:
\begin{itemize}
\item {} 
\sphinxAtStartPar
Legacy mode: supports both \sphinxcode{\sphinxupquote{\textquotesingle{}C\textquotesingle{}}} (by row) and \sphinxcode{\sphinxupquote{\textquotesingle{}F\textquotesingle{}}} (by column)

\item {} 
\sphinxAtStartPar
Experimental mode: only supports \sphinxcode{\sphinxupquote{\textquotesingle{}C\textquotesingle{}}} (by row)

\end{itemize}

\sphinxAtStartPar
2.Advanced indexing support:

\sphinxAtStartPar
Take the nqueen problem for example, only one queen can appear on each diagonal:
\begin{itemize}
\item {} 
\sphinxAtStartPar
Legacy mode: Compatible with \sphinxcode{\sphinxupquote{NumPy}}’s advanced index, which can be implemented as follows:

\begin{sphinxVerbatim}[commandchars=\\\{\}]
\PYG{k+kn}{import}\PYG{+w}{ }\PYG{n+nn}{coptpy}\PYG{+w}{ }\PYG{k}{as}\PYG{+w}{ }\PYG{n+nn}{cp}
\PYG{k+kn}{import}\PYG{+w}{ }\PYG{n+nn}{numpy}\PYG{+w}{ }\PYG{k}{as}\PYG{+w}{ }\PYG{n+nn}{np}

\PYG{n}{env} \PYG{o}{=} \PYG{n}{cp}\PYG{o}{.}\PYG{n}{Envr}\PYG{p}{(}\PYG{p}{)}
\PYG{n}{model} \PYG{o}{=} \PYG{n}{env}\PYG{o}{.}\PYG{n}{createModel}\PYG{p}{(}\PYG{p}{)}
\PYG{n}{model}\PYG{o}{.}\PYG{n}{matrixmodelmode} \PYG{o}{=} \PYG{l+s+s2}{\PYGZdq{}}\PYG{l+s+s2}{legacy}\PYG{l+s+s2}{\PYGZdq{}}
\PYG{k}{for} \PYG{n}{i} \PYG{o+ow}{in} \PYG{n+nb}{range}\PYG{p}{(}\PYG{l+m+mi}{1}\PYG{p}{,} \PYG{l+m+mi}{2}\PYG{o}{*}\PYG{n}{n}\PYG{p}{)}\PYG{p}{:}
    \PYG{c+c1}{\PYGZsh{} At most one queen per diagonal}
    \PYG{n}{diagn} \PYG{o}{=} \PYG{p}{(}\PYG{n+nb}{range}\PYG{p}{(}\PYG{n+nb}{max}\PYG{p}{(}\PYG{l+m+mi}{0}\PYG{p}{,} \PYG{n}{i}\PYG{o}{\PYGZhy{}}\PYG{n}{n}\PYG{p}{)}\PYG{p}{,} \PYG{n+nb}{min}\PYG{p}{(}\PYG{n}{n}\PYG{p}{,} \PYG{n}{i}\PYG{p}{)}\PYG{p}{)}\PYG{p}{,} \PYG{n+nb}{range}\PYG{p}{(}\PYG{n+nb}{min}\PYG{p}{(}\PYG{n}{n}\PYG{p}{,} \PYG{n}{i}\PYG{p}{)}\PYG{o}{\PYGZhy{}}\PYG{l+m+mi}{1}\PYG{p}{,} \PYG{n+nb}{max}\PYG{p}{(}\PYG{l+m+mi}{0}\PYG{p}{,} \PYG{n}{i}\PYG{o}{\PYGZhy{}}\PYG{n}{n}\PYG{p}{)}\PYG{o}{\PYGZhy{}}\PYG{l+m+mi}{1}\PYG{p}{,} \PYG{o}{\PYGZhy{}}\PYG{l+m+mi}{1}\PYG{p}{)}\PYG{p}{)}
    \PYG{n}{model}\PYG{o}{.}\PYG{n}{addConstrs}\PYG{p}{(}\PYG{n}{x}\PYG{p}{[}\PYG{n}{diagn}\PYG{p}{]}\PYG{o}{.}\PYG{n}{sum}\PYG{p}{(}\PYG{p}{)} \PYG{o}{\PYGZlt{}}\PYG{o}{=} \PYG{l+m+mi}{1}\PYG{p}{,} \PYG{n}{nameprefix}\PYG{o}{=}\PYG{l+s+s2}{\PYGZdq{}}\PYG{l+s+s2}{diag}\PYG{l+s+s2}{\PYGZdq{}}\PYG{o}{+}\PYG{n+nb}{str}\PYG{p}{(}\PYG{n}{i}\PYG{p}{)}\PYG{p}{)}
    \PYG{c+c1}{\PYGZsh{} At most one queen per anti\PYGZhy{}diagonal}
    \PYG{n}{adiagn} \PYG{o}{=} \PYG{p}{(}\PYG{n+nb}{range}\PYG{p}{(}\PYG{n+nb}{max}\PYG{p}{(}\PYG{l+m+mi}{0}\PYG{p}{,} \PYG{n}{i}\PYG{o}{\PYGZhy{}}\PYG{n}{n}\PYG{p}{)}\PYG{p}{,} \PYG{n+nb}{min}\PYG{p}{(}\PYG{n}{n}\PYG{p}{,} \PYG{n}{i}\PYG{p}{)}\PYG{p}{)}\PYG{p}{,} \PYG{n+nb}{range}\PYG{p}{(}\PYG{n+nb}{max}\PYG{p}{(}\PYG{l+m+mi}{0}\PYG{p}{,} \PYG{n}{n}\PYG{o}{\PYGZhy{}}\PYG{n}{i}\PYG{p}{)}\PYG{p}{,} \PYG{n+nb}{min}\PYG{p}{(}\PYG{n}{n}\PYG{p}{,} \PYG{l+m+mi}{2}\PYG{o}{*}\PYG{n}{n}\PYG{o}{\PYGZhy{}}\PYG{n}{i}\PYG{p}{)}\PYG{p}{)}\PYG{p}{)}
    \PYG{n}{model}\PYG{o}{.}\PYG{n}{addConstrs}\PYG{p}{(}\PYG{n}{x}\PYG{p}{[}\PYG{n}{adiagn}\PYG{p}{]}\PYG{o}{.}\PYG{n}{sum}\PYG{p}{(}\PYG{p}{)} \PYG{o}{\PYGZlt{}}\PYG{o}{=} \PYG{l+m+mi}{1}\PYG{p}{,} \PYG{n}{nameprefix}\PYG{o}{=}\PYG{l+s+s2}{\PYGZdq{}}\PYG{l+s+s2}{adiag}\PYG{l+s+s2}{\PYGZdq{}}\PYG{o}{+}\PYG{n+nb}{str}\PYG{p}{(}\PYG{n}{i}\PYG{p}{)}\PYG{p}{)}
\end{sphinxVerbatim}

\item {} 
\sphinxAtStartPar
Experimental mode: The above advanced indexing method is not supported and requires the built\sphinxhyphen{}in Python indexing method,
which can be implemented as follows:

\begin{sphinxVerbatim}[commandchars=\\\{\}]
\PYG{k+kn}{import}\PYG{+w}{ }\PYG{n+nn}{coptpy}\PYG{+w}{ }\PYG{k}{as}\PYG{+w}{ }\PYG{n+nn}{cp}
\PYG{k+kn}{import}\PYG{+w}{ }\PYG{n+nn}{numpy}\PYG{+w}{ }\PYG{k}{as}\PYG{+w}{ }\PYG{n+nn}{np}

\PYG{n}{env} \PYG{o}{=} \PYG{n}{cp}\PYG{o}{.}\PYG{n}{Envr}\PYG{p}{(}\PYG{p}{)}
\PYG{n}{model} \PYG{o}{=} \PYG{n}{env}\PYG{o}{.}\PYG{n}{createModel}\PYG{p}{(}\PYG{p}{)}
\PYG{n}{model}\PYG{o}{.}\PYG{n}{matrixmodelmode} \PYG{o}{=} \PYG{l+s+s2}{\PYGZdq{}}\PYG{l+s+s2}{experimental}\PYG{l+s+s2}{\PYGZdq{}}
\PYG{k}{for} \PYG{n}{i} \PYG{n}{i} \PYG{o+ow}{in} \PYG{n+nb}{range}\PYG{p}{(}\PYG{o}{\PYGZhy{}}\PYG{n}{n}\PYG{o}{+}\PYG{l+m+mi}{1}\PYG{p}{,} \PYG{n}{n}\PYG{p}{)}\PYG{p}{:}
    \PYG{c+c1}{\PYGZsh{} At most one queen per diagonal}
    \PYG{n}{model}\PYG{o}{.}\PYG{n}{addConstrs}\PYG{p}{(}\PYG{n}{x}\PYG{o}{.}\PYG{n}{diagonal}\PYG{p}{(}\PYG{n}{i}\PYG{p}{)}\PYG{o}{.}\PYG{n}{sum}\PYG{p}{(}\PYG{p}{)} \PYG{o}{\PYGZlt{}}\PYG{o}{=} \PYG{l+m+mi}{1}\PYG{p}{,} \PYG{n}{nameprefix}\PYG{o}{=}\PYG{l+s+s2}{\PYGZdq{}}\PYG{l+s+s2}{diag}\PYG{l+s+s2}{\PYGZdq{}}\PYG{o}{+}\PYG{n+nb}{str}\PYG{p}{(}\PYG{n}{i}\PYG{p}{)}\PYG{p}{)}
    \PYG{c+c1}{\PYGZsh{} At most one queen per anti\PYGZhy{}diagonal}
    \PYG{n}{model}\PYG{o}{.}\PYG{n}{addConstrs}\PYG{p}{(}\PYG{n}{x}\PYG{p}{[}\PYG{p}{:}\PYG{p}{,} \PYG{p}{:}\PYG{p}{:}\PYG{o}{\PYGZhy{}}\PYG{l+m+mi}{1}\PYG{p}{]}\PYG{o}{.}\PYG{n}{diagonal}\PYG{p}{(}\PYG{n}{i}\PYG{p}{)}\PYG{o}{.}\PYG{n}{sum}\PYG{p}{(}\PYG{p}{)} \PYG{o}{\PYGZlt{}}\PYG{o}{=} \PYG{l+m+mi}{1}\PYG{p}{,} \PYG{n}{nameprefix}\PYG{o}{=}\PYG{l+s+s2}{\PYGZdq{}}\PYG{l+s+s2}{adiag}\PYG{l+s+s2}{\PYGZdq{}}\PYG{o}{+}\PYG{n+nb}{str}\PYG{p}{(}\PYG{n}{i}\PYG{p}{)}\PYG{p}{)}
\end{sphinxVerbatim}

\end{itemize}

\sphinxAtStartPar
3.The relevant information obtained from matrix objects (objects such as \sphinxcode{\sphinxupquote{MConstr}} and \sphinxcode{\sphinxupquote{MVar}}), and the return value types are different:
\begin{itemize}
\item {} 
\sphinxAtStartPar
Legacy mode: \sphinxcode{\sphinxupquote{numpy.ndarray}}

\item {} 
\sphinxAtStartPar
Experimental mode: \sphinxcode{\sphinxupquote{coptcore.NdArray}}

\end{itemize}

\sphinxAtStartPar
4.The supported dimensions are different:
\begin{itemize}
\item {} 
\sphinxAtStartPar
Legacy mode: There is no restriction on dimensions and can support higher\sphinxhyphen{}dimensional matrices.

\item {} 
\sphinxAtStartPar
Experimental mode: The highest supported dimensions are three.

\end{itemize}

\sphinxAtStartPar
5.When adding multi\sphinxhyphen{}dimensional quadratic constraints, the support for variable dimensions participating in constraint formation
and the type of return value are different:
\begin{itemize}
\item {} 
\sphinxAtStartPar
Legacy mode: The dimension of the variables involved in forming the quadratic constraint needs to be 1 (i.e. vector),
higher dimensions are not supported, and the return value type is {\hyperref[\detokenize{pyapiref:chappyapi-qconstraint}]{\sphinxcrossref{\DUrole{std,std-ref}{QConstraint Class}}}} object.

\begin{sphinxVerbatim}[commandchars=\\\{\}]
\PYG{k+kn}{import}\PYG{+w}{ }\PYG{n+nn}{coptpy}\PYG{+w}{ }\PYG{k}{as}\PYG{+w}{ }\PYG{n+nn}{cp}
\PYG{k+kn}{import}\PYG{+w}{ }\PYG{n+nn}{numpy}\PYG{+w}{ }\PYG{k}{as}\PYG{+w}{ }\PYG{n+nn}{np}

\PYG{n}{env} \PYG{o}{=} \PYG{n}{cp}\PYG{o}{.}\PYG{n}{Envr}\PYG{p}{(}\PYG{p}{)}
\PYG{n}{model} \PYG{o}{=} \PYG{n}{env}\PYG{o}{.}\PYG{n}{createModel}\PYG{p}{(}\PYG{p}{)}
\PYG{n}{model}\PYG{o}{.}\PYG{n}{matrixmodelmode} \PYG{o}{=} \PYG{l+s+s2}{\PYGZdq{}}\PYG{l+s+s2}{legacy}\PYG{l+s+s2}{\PYGZdq{}}
\PYG{n}{Q} \PYG{o}{=} \PYG{n}{np}\PYG{o}{.}\PYG{n}{full}\PYG{p}{(}\PYG{p}{(}\PYG{l+m+mi}{3}\PYG{p}{,} \PYG{l+m+mi}{3}\PYG{p}{)}\PYG{p}{,} \PYG{l+m+mi}{1}\PYG{p}{)}
\PYG{n}{mx} \PYG{o}{=} \PYG{n}{model}\PYG{o}{.}\PYG{n}{addMVar}\PYG{p}{(}\PYG{l+m+mi}{3}\PYG{p}{,} \PYG{n}{nameprefix}\PYG{o}{=}\PYG{l+s+s2}{\PYGZdq{}}\PYG{l+s+s2}{mx}\PYG{l+s+s2}{\PYGZdq{}}\PYG{p}{)}
\PYG{c+c1}{\PYGZsh{} mqc \PYGZlt{}coptpy.QConstraint: \PYGZgt{}}
\PYG{n}{mqc} \PYG{o}{=} \PYG{n}{model}\PYG{o}{.}\PYG{n}{addQConstr}\PYG{p}{(}\PYG{n}{mx}\PYG{n+nd}{@Q}\PYG{n+nd}{@mx}\PYG{o}{\PYGZlt{}}\PYG{o}{=}\PYG{l+m+mi}{1}\PYG{p}{)}
\end{sphinxVerbatim}

\item {} 
\sphinxAtStartPar
Experimental mode: The dimensions of variables participating in the quadratic constraints can be one or two\sphinxhyphen{}dimensional,
and the return value type is {\hyperref[\detokenize{pyapiref:chappyapi-mqconstr}]{\sphinxcrossref{\DUrole{std,std-ref}{MQConstr Class}}}} object.

\begin{sphinxVerbatim}[commandchars=\\\{\}]
\PYG{k+kn}{import}\PYG{+w}{ }\PYG{n+nn}{coptpy}\PYG{+w}{ }\PYG{k}{as}\PYG{+w}{ }\PYG{n+nn}{cp}
\PYG{k+kn}{import}\PYG{+w}{ }\PYG{n+nn}{numpy}\PYG{+w}{ }\PYG{k}{as}\PYG{+w}{ }\PYG{n+nn}{np}

\PYG{n}{env} \PYG{o}{=} \PYG{n}{cp}\PYG{o}{.}\PYG{n}{Envr}\PYG{p}{(}\PYG{p}{)}
\PYG{n}{model} \PYG{o}{=} \PYG{n}{env}\PYG{o}{.}\PYG{n}{createModel}\PYG{p}{(}\PYG{p}{)}
\PYG{n}{model}\PYG{o}{.}\PYG{n}{matrixmodelmode} \PYG{o}{=} \PYG{l+s+s2}{\PYGZdq{}}\PYG{l+s+s2}{experimental}\PYG{l+s+s2}{\PYGZdq{}}
\PYG{n}{Q} \PYG{o}{=} \PYG{n}{np}\PYG{o}{.}\PYG{n}{full}\PYG{p}{(}\PYG{p}{(}\PYG{l+m+mi}{3}\PYG{p}{,} \PYG{l+m+mi}{3}\PYG{p}{)}\PYG{p}{,} \PYG{l+m+mi}{1}\PYG{p}{)}
\PYG{n}{mx} \PYG{o}{=} \PYG{n}{model}\PYG{o}{.}\PYG{n}{addMVar}\PYG{p}{(}\PYG{p}{(}\PYG{l+m+mi}{3}\PYG{p}{,} \PYG{l+m+mi}{3}\PYG{p}{)}\PYG{p}{,} \PYG{n}{nameprefix}\PYG{o}{=}\PYG{l+s+s2}{\PYGZdq{}}\PYG{l+s+s2}{mx}\PYG{l+s+s2}{\PYGZdq{}}\PYG{p}{)}
\PYG{c+c1}{\PYGZsh{} mqc \PYGZlt{}coptpy.MQConstr: shape=(3, 3)\PYGZgt{}}
\PYG{n}{mqc} \PYG{o}{=} \PYG{n}{model}\PYG{o}{.}\PYG{n}{addQConstr}\PYG{p}{(}\PYG{n}{mx}\PYG{n+nd}{@Q}\PYG{n+nd}{@mx} \PYG{o}{\PYGZlt{}}\PYG{o}{=} \PYG{l+m+mi}{1}\PYG{p}{)}
\end{sphinxVerbatim}

\end{itemize}

\section{Multi\sphinxhyphen{}dimensional Variables}
\label{\detokenize{matrix:multi-dimensional-variables}}\begin{enumerate}
\sphinxsetlistlabels{\arabic}{enumi}{enumii}{}{.}%
\item {} 
\sphinxAtStartPar
Add multi\sphinxhyphen{}dimensional variable \sphinxcode{\sphinxupquote{MVar}}

\end{enumerate}
\begin{quote}

\sphinxAtStartPar
\sphinxcode{\sphinxupquote{MVar}} constains operations related to multi\sphinxhyphen{}dimensional variables. Users can use \sphinxcode{\sphinxupquote{Model.addMVar()}} to add a matmulti\sphinxhyphen{}dimensionalrix variable \sphinxcode{\sphinxupquote{MVar}} of any dimension and shape to the model. In addition to the need of specifying the argument \sphinxcode{\sphinxupquote{shape}} (matrix shape), the rest of the arguments are consistent with ordinary variables, including: \sphinxcode{\sphinxupquote{lb}} , \sphinxcode{\sphinxupquote{ub}} , \sphinxcode{\sphinxupquote{vtype}} , \sphinxcode{\sphinxupquote{nameprefix}} .
\begin{itemize}
\item {} 
\sphinxAtStartPar
Add one\sphinxhyphen{}dimensional continuous multi\sphinxhyphen{}dimensional variables: \sphinxcode{\sphinxupquote{x = Model.addMVar(3)}}

\item {} 
\sphinxAtStartPar
Add two\sphinxhyphen{}dimensional 3x3 binary multi\sphinxhyphen{}dimensional variables: \sphinxcode{\sphinxupquote{y = Model.addMVar(shape(3,3), vtype=COPT.BINARY)}}

\end{itemize}

\sphinxAtStartPar
In addition, the \sphinxcode{\sphinxupquote{MVar}} multi\sphinxhyphen{}dimensional variables can also be sliced, such as: \sphinxcode{\sphinxupquote{y1 = y{[}:,0:2{]}}}
\end{quote}
\begin{enumerate}
\sphinxsetlistlabels{\arabic}{enumi}{enumii}{}{.}%
\setcounter{enumi}{1}
\item {} 
\sphinxAtStartPar
Get multi\sphinxhyphen{}dimensional variable related attributes:
\begin{itemize}
\item {} 
\sphinxAtStartPar
Number of dimensions: \sphinxcode{\sphinxupquote{MVar.ndim}}

\item {} 
\sphinxAtStartPar
The shape of the multi\sphinxhyphen{}dimensional variable: \sphinxcode{\sphinxupquote{MVar.shape}}

\item {} 
\sphinxAtStartPar
Number of elements in the multi\sphinxhyphen{}dimensional variable: \sphinxcode{\sphinxupquote{MVar.size}}

\end{itemize}

\end{enumerate}

\section{Multi\sphinxhyphen{}dimensional array operations and expressions}
\label{\detokenize{matrix:multi-dimensional-array-operations-and-expressions}}

\subsection{Multi\sphinxhyphen{}dimensional Linear Expressions}
\label{\detokenize{matrix:multi-dimensional-linear-expressions}}
\sphinxAtStartPar
Multi\sphinxhyphen{}dimensional variables and their coefficients (can be \sphinxcode{\sphinxupquote{ndarray}} ) form a Multi\sphinxhyphen{}dimensional linear expression (\sphinxcode{\sphinxupquote{MLinExpr}}), and the supported operations mainly include:
\begin{enumerate}
\sphinxsetlistlabels{\arabic}{enumi}{enumii}{}{.}%
\item {} 
\sphinxAtStartPar
Matrix multiplication: A @ x

\end{enumerate}
\begin{quote}

\begin{sphinxVerbatim}[commandchars=\\\{\}]
\PYG{n}{x} \PYG{o}{=} \PYG{n}{model}\PYG{o}{.}\PYG{n}{addMVar}\PYG{p}{(}\PYG{l+m+mi}{3}\PYG{p}{)}
\PYG{n}{A} \PYG{o}{=} \PYG{n}{np}\PYG{o}{.}\PYG{n}{array}\PYG{p}{(}\PYG{p}{[}\PYG{p}{[}\PYG{l+m+mi}{1}\PYG{p}{,} \PYG{l+m+mi}{0}\PYG{p}{,} \PYG{l+m+mi}{1}\PYG{p}{]}\PYG{p}{,}\PYG{p}{[}\PYG{l+m+mi}{0}\PYG{p}{,} \PYG{l+m+mi}{0}\PYG{p}{,} \PYG{l+m+mi}{1}\PYG{p}{]}\PYG{p}{]}\PYG{p}{)}
\PYG{n}{expr1} \PYG{o}{=} \PYG{n}{A} \PYG{o}{@} \PYG{n}{x}
\end{sphinxVerbatim}
\end{quote}
\begin{enumerate}
\sphinxsetlistlabels{\arabic}{enumi}{enumii}{}{.}%
\setcounter{enumi}{1}
\item {} 
\sphinxAtStartPar
Vector inner product

\end{enumerate}
\begin{quote}

\begin{sphinxVerbatim}[commandchars=\\\{\}]
\PYG{n}{x} \PYG{o}{=} \PYG{n}{model}\PYG{o}{.}\PYG{n}{addMVar}\PYG{p}{(}\PYG{l+m+mi}{3}\PYG{p}{)}
\PYG{n}{c} \PYG{o}{=} \PYG{n}{np}\PYG{o}{.}\PYG{n}{array}\PYG{p}{(}\PYG{p}{[}\PYG{l+m+mi}{1}\PYG{p}{,} \PYG{l+m+mi}{2}\PYG{p}{,} \PYG{l+m+mi}{3}\PYG{p}{]}\PYG{p}{)}
\PYG{n}{expr2} \PYG{o}{=} \PYG{n}{c} \PYG{o}{@} \PYG{n}{x}
\end{sphinxVerbatim}
\end{quote}

\subsection{Multi\sphinxhyphen{}dimensional Quadratic Expression}
\label{\detokenize{matrix:multi-dimensional-quadratic-expression}}
\sphinxAtStartPar
Common multi\sphinxhyphen{}dimensional quadratic expressions and their corresponding mathematical forms are as follows:
\begin{itemize}
\item {} 
\sphinxAtStartPar
x @ Q @ x: \(x^TQx\)

\item {} 
\sphinxAtStartPar
x @ x: \(x^Tx\)

\item {} 
\sphinxAtStartPar
x @ Q @ x + c @ x + b: \(x^TQx+c^Tx+b\)

\end{itemize}

\subsection{Other multi\sphinxhyphen{}dimensional array operations}
\label{\detokenize{matrix:other-multi-dimensional-array-operations}}\begin{enumerate}
\sphinxsetlistlabels{\arabic}{enumi}{enumii}{}{.}%
\item {} 
\sphinxAtStartPar
Combine with regular linear variables, regular linear expressions, and constants:

\end{enumerate}
\begin{quote}

\begin{sphinxVerbatim}[commandchars=\\\{\}]
\PYG{n}{x} \PYG{o}{=} \PYG{n}{model}\PYG{o}{.}\PYG{n}{addMVar}\PYG{p}{(}\PYG{l+m+mi}{3}\PYG{p}{)}
\PYG{n}{y} \PYG{o}{=} \PYG{n}{model}\PYG{o}{.}\PYG{n}{addVar}\PYG{p}{(}\PYG{p}{)}
\PYG{n}{c} \PYG{o}{=} \PYG{n}{np}\PYG{o}{.}\PYG{n}{array}\PYG{p}{(}\PYG{p}{[}\PYG{l+m+mi}{1}\PYG{p}{,} \PYG{l+m+mi}{2}\PYG{p}{,} \PYG{l+m+mi}{3}\PYG{p}{]}\PYG{p}{)}
\PYG{n}{Q} \PYG{o}{=} \PYG{n}{np}\PYG{o}{.}\PYG{n}{full}\PYG{p}{(}\PYG{p}{(}\PYG{l+m+mi}{3}\PYG{p}{,} \PYG{l+m+mi}{3}\PYG{p}{)}\PYG{p}{,} \PYG{l+m+mi}{1}\PYG{p}{)}
\PYG{n}{expr3} \PYG{o}{=} \PYG{l+m+mi}{2} \PYG{o}{*} \PYG{n}{x} \PYG{o}{@} \PYG{n}{Q} \PYG{o}{@} \PYG{n}{x} \PYG{o}{+} \PYG{n}{c} \PYG{o}{@} \PYG{n}{x} \PYG{o}{+} \PYG{l+m+mi}{2} \PYG{o}{*} \PYG{n}{y} \PYG{o}{+} \PYG{l+m+mi}{1}
\end{sphinxVerbatim}
\end{quote}
\begin{enumerate}
\sphinxsetlistlabels{\arabic}{enumi}{enumii}{}{.}%
\setcounter{enumi}{1}
\item {} 
\sphinxAtStartPar
Self\sphinxhyphen{}increment/self\sphinxhyphen{}subtraction/self\sphinxhyphen{}multiplication operations:

\end{enumerate}
\begin{quote}

\begin{sphinxVerbatim}[commandchars=\\\{\}]
\PYG{n}{mx} \PYG{o}{=} \PYG{n}{m}\PYG{o}{.}\PYG{n}{addMVar}\PYG{p}{(}\PYG{p}{(}\PYG{l+m+mi}{3}\PYG{p}{,} \PYG{l+m+mi}{3}\PYG{p}{)}\PYG{p}{)}
\PYG{n}{B} \PYG{o}{=} \PYG{n}{np}\PYG{o}{.}\PYG{n}{array}\PYG{p}{(}\PYG{p}{[}\PYG{p}{[}\PYG{l+m+mi}{1}\PYG{p}{,} \PYG{l+m+mi}{0}\PYG{p}{,} \PYG{l+m+mi}{1}\PYG{p}{]}\PYG{p}{,} \PYG{p}{[}\PYG{l+m+mi}{0}\PYG{p}{,} \PYG{l+m+mi}{1}\PYG{p}{,} \PYG{l+m+mi}{1}\PYG{p}{]}\PYG{p}{]}\PYG{p}{)}
\PYG{n}{expr\PYGZus{}add} \PYG{o}{=} \PYG{n}{B} \PYG{o}{@} \PYG{n}{mx}
\PYG{n}{expr\PYGZus{}add} \PYG{o}{+}\PYG{o}{=} \PYG{l+m+mi}{1}
\PYG{n}{expr\PYGZus{}add} \PYG{o}{*}\PYG{o}{=} \PYG{l+m+mi}{2}
\end{sphinxVerbatim}
\end{quote}

\begin{sphinxadmonition}{note}{Notes}
\begin{itemize}
\item {} 
\sphinxAtStartPar
When we directly print the multi\sphinxhyphen{}dimensional expression with \sphinxcode{\sphinxupquote{print(MLinExpr)/print(MQuadExpr)}} , the \sphinxcode{\sphinxupquote{shape}} of the expression will be output at the same time. When \sphinxcode{\sphinxupquote{shape=()}} , it means that the expression is a scalar (single linear/quadratic expression), corresponding to \sphinxcode{\sphinxupquote{ndim=0}} , \sphinxcode{\sphinxupquote{size=1}} . The same is true for multi\sphinxhyphen{}dimensional variables \sphinxcode{\sphinxupquote{MVar}} ;

\item {} 
\sphinxAtStartPar
When performing matrix multiplication (A@x), the matrix multiplication algorithm needs to be satisfied, and the number of columns of A and the number of rows of X need to be the same;

\item {} 
\sphinxAtStartPar
COPT supports the combination of \sphinxcode{\sphinxupquote{MLinExpr}} and \sphinxcode{\sphinxupquote{LinExpr}} , but it should be noted that the \sphinxcode{\sphinxupquote{MLinExpr}} needs \sphinxcode{\sphinxupquote{shape=()}} at this time, and the final returned expression is \sphinxcode{\sphinxupquote{MLinExpr}} with \sphinxcode{\sphinxupquote{shape=()}} .

\end{itemize}
\end{sphinxadmonition}

\section{Matrix Constraints}
\label{\detokenize{matrix:matrix-constraints}}

\subsection{Matrix linear Constraints}
\label{\detokenize{matrix:matrix-linear-constraints}}
\sphinxAtStartPar
COPT supports two ways of adding multi\sphinxhyphen{}dimensional linear constraints, and the format provided by the function arguments is different:
\begin{enumerate}
\sphinxsetlistlabels{\arabic}{enumi}{enumii}{}{.}%
\item {} 
\sphinxAtStartPar
\sphinxcode{\sphinxupquote{Model.addMConstr()}} that specifically adds multi\sphinxhyphen{}dimensional linear constraints, the arguments that can be specified are:
\begin{itemize}
\item {} 
\sphinxAtStartPar
\sphinxcode{\sphinxupquote{A}} : coefficient matrix for linear constraints

\item {} 
\sphinxAtStartPar
\sphinxcode{\sphinxupquote{x}} : decision variables ( \sphinxcode{\sphinxupquote{MVar}} )

\item {} 
\sphinxAtStartPar
\sphinxcode{\sphinxupquote{sense}} : type of linear constraint, the possible values are: \sphinxcode{\sphinxupquote{\textquotesingle{}L\textquotesingle{}}} (\textless{}=), \sphinxcode{\sphinxupquote{\textquotesingle{}G\textquotesingle{}}} (\textgreater{}=), \sphinxcode{\sphinxupquote{\textquotesingle{}E\textquotesingle{}}} (=)

\item {} 
\sphinxAtStartPar
\sphinxcode{\sphinxupquote{b}} : right\sphinxhyphen{}hand\sphinxhyphen{}side of linear constraints (vector with dimensions equal to the number of rows of matrix \sphinxcode{\sphinxupquote{A}})

\item {} 
\sphinxAtStartPar
\sphinxcode{\sphinxupquote{name}} : name prefix for linear constraints

\end{itemize}

\end{enumerate}
\begin{quote}

\begin{sphinxVerbatim}[commandchars=\\\{\}]
\PYG{n}{x} \PYG{o}{=} \PYG{n}{model}\PYG{o}{.}\PYG{n}{addMVar}\PYG{p}{(}\PYG{n}{shape}\PYG{o}{=}\PYG{l+m+mi}{3}\PYG{p}{,} \PYG{n}{vtype}\PYG{o}{=}\PYG{n}{COPT}\PYG{o}{.}\PYG{n}{BINARY}\PYG{p}{,} \PYG{n}{nameprefix}\PYG{o}{=}\PYG{l+s+s1}{\PYGZsq{}}\PYG{l+s+s1}{x}\PYG{l+s+s1}{\PYGZsq{}}\PYG{p}{)}
\PYG{n}{A} \PYG{o}{=} \PYG{n}{np}\PYG{o}{.}\PYG{n}{array}\PYG{p}{(}\PYG{p}{[}\PYG{p}{[}\PYG{l+m+mi}{1}\PYG{p}{,} \PYG{l+m+mi}{2}\PYG{p}{,} \PYG{l+m+mi}{3}\PYG{p}{]}\PYG{p}{,} \PYG{p}{[}\PYG{l+m+mi}{3}\PYG{p}{,} \PYG{l+m+mi}{2}\PYG{p}{,} \PYG{l+m+mi}{1}\PYG{p}{]}\PYG{p}{]}\PYG{p}{)}
\PYG{n}{b} \PYG{o}{=} \PYG{n}{np}\PYG{o}{.}\PYG{n}{array}\PYG{p}{(}\PYG{p}{[}\PYG{l+m+mi}{2}\PYG{p}{,} \PYG{l+m+mi}{5}\PYG{p}{]}\PYG{p}{)}
\PYG{n}{mconstrs} \PYG{o}{=} \PYG{n}{model}\PYG{o}{.}\PYG{n}{addMConstr}\PYG{p}{(}\PYG{n}{A}\PYG{p}{,} \PYG{n}{x}\PYG{p}{,} \PYG{l+s+s1}{\PYGZsq{}}\PYG{l+s+s1}{L}\PYG{l+s+s1}{\PYGZsq{}}\PYG{p}{,} \PYG{n}{b}\PYG{p}{,} \PYG{n}{nameprefix}\PYG{o}{=}\PYG{l+s+s1}{\PYGZsq{}}\PYG{l+s+s1}{c}\PYG{l+s+s1}{\PYGZsq{}}\PYG{p}{)}
\PYG{n}{obj} \PYG{o}{=} \PYG{n}{np}\PYG{o}{.}\PYG{n}{array}\PYG{p}{(}\PYG{p}{[}\PYG{l+m+mi}{1}\PYG{p}{,} \PYG{l+m+mi}{2}\PYG{p}{,} \PYG{l+m+mi}{1}\PYG{p}{]}\PYG{p}{)}
\PYG{n}{model}\PYG{o}{.}\PYG{n}{setObjective}\PYG{p}{(}\PYG{n}{obj} \PYG{o}{@} \PYG{n}{x}\PYG{p}{,} \PYG{n}{COPT}\PYG{o}{.}\PYG{n}{MINIMIZE}\PYG{p}{)}
\end{sphinxVerbatim}
\end{quote}
\begin{enumerate}
\sphinxsetlistlabels{\arabic}{enumi}{enumii}{}{.}%
\setcounter{enumi}{1}
\item {} 
\sphinxAtStartPar
Matrix linear constraints can be regarded as a set of linear constraints, so \sphinxcode{\sphinxupquote{Model.addConstrs()}} can also add multi\sphinxhyphen{}dimensional linear constraints:

\end{enumerate}
\begin{quote}

\begin{sphinxVerbatim}[commandchars=\\\{\}]
\PYG{n}{x} \PYG{o}{=} \PYG{n}{model}\PYG{o}{.}\PYG{n}{addMVar}\PYG{p}{(}\PYG{n}{shape}\PYG{o}{=}\PYG{l+m+mi}{3}\PYG{p}{,} \PYG{n}{vtype}\PYG{o}{=}\PYG{n}{COPT}\PYG{o}{.}\PYG{n}{BINARY}\PYG{p}{,} \PYG{n}{nameprefix}\PYG{o}{=}\PYG{l+s+s1}{\PYGZsq{}}\PYG{l+s+s1}{x}\PYG{l+s+s1}{\PYGZsq{}}\PYG{p}{)}
\PYG{n}{A} \PYG{o}{=} \PYG{n}{np}\PYG{o}{.}\PYG{n}{array}\PYG{p}{(}\PYG{p}{[}\PYG{p}{[}\PYG{l+m+mi}{1}\PYG{p}{,} \PYG{l+m+mi}{2}\PYG{p}{,} \PYG{l+m+mi}{3}\PYG{p}{]}\PYG{p}{,} \PYG{p}{[}\PYG{l+m+mi}{3}\PYG{p}{,} \PYG{l+m+mi}{2}\PYG{p}{,} \PYG{l+m+mi}{1}\PYG{p}{]}\PYG{p}{]}\PYG{p}{)}
\PYG{n}{b} \PYG{o}{=} \PYG{n}{np}\PYG{o}{.}\PYG{n}{array}\PYG{p}{(}\PYG{p}{[}\PYG{l+m+mi}{2}\PYG{p}{,} \PYG{l+m+mi}{5}\PYG{p}{]}\PYG{p}{)}
\PYG{n}{mconstrs} \PYG{o}{=} \PYG{n}{model}\PYG{o}{.}\PYG{n}{addConstrs}\PYG{p}{(}\PYG{n}{A} \PYG{o}{@} \PYG{n}{x} \PYG{o}{\PYGZlt{}}\PYG{o}{=} \PYG{n}{b}\PYG{p}{,} \PYG{n}{nameprefix}\PYG{o}{=}\PYG{l+s+s1}{\PYGZsq{}}\PYG{l+s+s1}{c}\PYG{l+s+s1}{\PYGZsq{}}\PYG{p}{)}
\PYG{n}{obj} \PYG{o}{=} \PYG{n}{np}\PYG{o}{.}\PYG{n}{array}\PYG{p}{(}\PYG{p}{[}\PYG{l+m+mi}{1}\PYG{p}{,} \PYG{l+m+mi}{2}\PYG{p}{,} \PYG{l+m+mi}{1}\PYG{p}{]}\PYG{p}{)}
\PYG{n}{model}\PYG{o}{.}\PYG{n}{setObjective}\PYG{p}{(}\PYG{n}{obj} \PYG{o}{@} \PYG{n}{x}\PYG{p}{,} \PYG{n}{COPT}\PYG{o}{.}\PYG{n}{MINIMIZE}\PYG{p}{)}
\end{sphinxVerbatim}
\end{quote}

\subsection{Quadratic Constraints}
\label{\detokenize{matrix:quadratic-constraints}}
\sphinxAtStartPar
COPT supports two ways of constructing multi\sphinxhyphen{}dimensional quadratic constraints, and the format provided by the function arguments is different:
\begin{enumerate}
\sphinxsetlistlabels{\arabic}{enumi}{enumii}{}{.}%
\item {} 
\sphinxAtStartPar
\sphinxcode{\sphinxupquote{Model.addMQConstr()}} that specifically adds multi\sphinxhyphen{}dimensional quadratic constraints, the arguments that can be specified are:
\begin{itemize}
\item {} 
\sphinxAtStartPar
\sphinxcode{\sphinxupquote{Q}} : quadratic coefficient matrix

\item {} 
\sphinxAtStartPar
\sphinxcode{\sphinxupquote{c}} : vector of linear term coefficients, or \sphinxcode{\sphinxupquote{None}} if there is no linear term

\item {} 
\sphinxAtStartPar
\sphinxcode{\sphinxupquote{sense}} : type of quadratic constraint, the possible values are: \sphinxcode{\sphinxupquote{\textquotesingle{}L\textquotesingle{}}} (\textless{}=), \sphinxcode{\sphinxupquote{\textquotesingle{}G\textquotesingle{}}} (\textgreater{}=), \sphinxcode{\sphinxupquote{\textquotesingle{}E\textquotesingle{}}} (=)

\item {} 
\sphinxAtStartPar
\sphinxcode{\sphinxupquote{rhs}} : right\sphinxhyphen{}hand\sphinxhyphen{}side of quadratic constraints

\item {} 
\sphinxAtStartPar
\sphinxcode{\sphinxupquote{xQ\_L}} : the left\sphinxhyphen{}hand variable of the quadratic coefficient matrix Q (vector whose length is consistent with the number of rows of the matrix \sphinxcode{\sphinxupquote{Q}} )

\item {} 
\sphinxAtStartPar
\sphinxcode{\sphinxupquote{xQ\_R}} : right\sphinxhyphen{}hand variable of the quadratic coefficient matrix Q (vector whose length is consistent with the number of columns of the matrix \sphinxcode{\sphinxupquote{Q}} )

\item {} 
\sphinxAtStartPar
\sphinxcode{\sphinxupquote{xc}} : the variables for the linear term, or \sphinxcode{\sphinxupquote{None}} if there is no linear term

\item {} 
\sphinxAtStartPar
\sphinxcode{\sphinxupquote{name}} : name prefix for quadratic constraints

\end{itemize}

\end{enumerate}
\begin{quote}

\begin{sphinxVerbatim}[commandchars=\\\{\}]
\PYG{n}{Q} \PYG{o}{=} \PYG{n}{np}\PYG{o}{.}\PYG{n}{diag}\PYG{p}{(}\PYG{p}{[}\PYG{l+m+mi}{3}\PYG{p}{,} \PYG{l+m+mi}{2}\PYG{p}{,} \PYG{l+m+mi}{1}\PYG{p}{]}\PYG{p}{)}
\PYG{n}{x} \PYG{o}{=} \PYG{n}{model}\PYG{o}{.}\PYG{n}{addMVar}\PYG{p}{(}\PYG{l+m+mi}{3}\PYG{p}{)}
\PYG{n}{c1} \PYG{o}{=} \PYG{n}{model}\PYG{o}{.}\PYG{n}{addMQConstr}\PYG{p}{(}\PYG{n}{Q}\PYG{p}{,} \PYG{k+kc}{None}\PYG{p}{,} \PYG{l+s+s1}{\PYGZsq{}}\PYG{l+s+s1}{L}\PYG{l+s+s1}{\PYGZsq{}}\PYG{p}{,} \PYG{l+m+mf}{1.0}\PYG{p}{,} \PYG{n}{x}\PYG{p}{,} \PYG{n}{x}\PYG{p}{)}
\end{sphinxVerbatim}
\end{quote}
\begin{enumerate}
\sphinxsetlistlabels{\arabic}{enumi}{enumii}{}{.}%
\setcounter{enumi}{1}
\item {} 
\sphinxAtStartPar
\sphinxcode{\sphinxupquote{Model.addQConstr()}}, directly gives the multi\sphinxhyphen{}dimensional quadratic expression

\end{enumerate}
\begin{quote}
\begin{itemize}
\item {} 
\sphinxAtStartPar
\sphinxcode{\sphinxupquote{lhs}} : multi\sphinxhyphen{}dimensional quadratic expression

\item {} 
\sphinxAtStartPar
\sphinxcode{\sphinxupquote{sense}} : constraint type

\item {} 
\sphinxAtStartPar
\sphinxcode{\sphinxupquote{rhs}} : right\sphinxhyphen{}hand\sphinxhyphen{}side of quadratic constraints

\end{itemize}

\begin{sphinxVerbatim}[commandchars=\\\{\}]
\PYG{n}{Q} \PYG{o}{=} \PYG{n}{np}\PYG{o}{.}\PYG{n}{diag}\PYG{p}{(}\PYG{p}{[}\PYG{l+m+mi}{3}\PYG{p}{,} \PYG{l+m+mi}{2}\PYG{p}{,} \PYG{l+m+mi}{1}\PYG{p}{]}\PYG{p}{)}
\PYG{n}{x} \PYG{o}{=} \PYG{n}{model}\PYG{o}{.}\PYG{n}{addMVar}\PYG{p}{(}\PYG{l+m+mi}{3}\PYG{p}{)}
\PYG{n}{c2} \PYG{o}{=} \PYG{n}{model}\PYG{o}{.}\PYG{n}{addQConstr}\PYG{p}{(}\PYG{n}{x}\PYG{n+nd}{@Q}\PYG{n+nd}{@x}\PYG{o}{\PYGZlt{}}\PYG{o}{=}\PYG{l+m+mf}{1.0}\PYG{p}{)}
\end{sphinxVerbatim}
\end{quote}

\section{Objective function composed of multi\sphinxhyphen{}dimensional variables}
\label{\detokenize{matrix:objective-function-composed-of-multi-dimensional-variables}}
\sphinxAtStartPar
COPT supports setting linear and quadratic objective functions, and provides two ways to set objective functions. The format of function arguments is different:
\begin{enumerate}
\sphinxsetlistlabels{\arabic}{enumi}{enumii}{}{.}%
\item {} 
\sphinxAtStartPar
\sphinxcode{\sphinxupquote{Model.setMObjective()}} that specifically sets the objective function composed of multi\sphinxhyphen{}dimensional variables, the arguments that can be specified are:
\begin{itemize}
\item {} 
\sphinxAtStartPar
\sphinxcode{\sphinxupquote{Q}} : quadratic coefficient matrix, or \sphinxcode{\sphinxupquote{None}} if the objective function is linear

\item {} 
\sphinxAtStartPar
\sphinxcode{\sphinxupquote{c}} : vector of linear term coefficients, or \sphinxcode{\sphinxupquote{None}} if there is no linear term

\item {} 
\sphinxAtStartPar
\sphinxcode{\sphinxupquote{constant}} : the constant term of the objective function

\item {} 
\sphinxAtStartPar
\sphinxcode{\sphinxupquote{xQ\_L}}: the left\sphinxhyphen{}hand variable of the quadratic term coefficient matrix Q (vector whose length is consistent with the number of rows of the matrix \sphinxcode{\sphinxupquote{Q}}), or \sphinxcode{\sphinxupquote{None}} if the objective function is linear

\item {} 
\sphinxAtStartPar
\sphinxcode{\sphinxupquote{xQ\_R}}: the right\sphinxhyphen{}hand variable of the quadratic coefficient matrix Q (vector, whose length is consistent with the number of columns in the matrix \sphinxcode{\sphinxupquote{Q}}), or \sphinxcode{\sphinxupquote{None}} if the objective function is linear

\item {} 
\sphinxAtStartPar
\sphinxcode{\sphinxupquote{xc}}: the variable for the linear term, or \sphinxcode{\sphinxupquote{None}} if there is no linear term

\item {} 
\sphinxAtStartPar
\sphinxcode{\sphinxupquote{sense}}: direction of optimization, possible values are: \sphinxcode{\sphinxupquote{COPT.MINIMIZE}} or \sphinxcode{\sphinxupquote{COPT.MAXIMIZE}}

\end{itemize}

\item {} 
\sphinxAtStartPar
\sphinxcode{\sphinxupquote{Model.setObjective()}} that directly gives the expression of the objective function
\begin{itemize}
\item {} 
\sphinxAtStartPar
\sphinxcode{\sphinxupquote{expr}}: Objective function expression, which can be linear or quadratic

\item {} 
\sphinxAtStartPar
\sphinxcode{\sphinxupquote{sense}}: optimization direction, possible values are: \sphinxcode{\sphinxupquote{COPT.MINIMIZE}} or \sphinxcode{\sphinxupquote{COPT.MINIMIZE}}

\end{itemize}

\end{enumerate}
\begin{quote}

\begin{sphinxVerbatim}[commandchars=\\\{\}]
\PYG{n}{x} \PYG{o}{=} \PYG{n}{model}\PYG{o}{.}\PYG{n}{addMVar}\PYG{p}{(}\PYG{n}{shape}\PYG{o}{=}\PYG{l+m+mi}{3}\PYG{p}{,} \PYG{n}{vtype}\PYG{o}{=}\PYG{n}{COPT}\PYG{o}{.}\PYG{n}{BINARY}\PYG{p}{,} \PYG{n}{nameprefix}\PYG{o}{=}\PYG{l+s+s2}{\PYGZdq{}}\PYG{l+s+s2}{x}\PYG{l+s+s2}{\PYGZdq{}}\PYG{p}{)}
\PYG{n}{obj} \PYG{o}{=} \PYG{n}{np}\PYG{o}{.} \PYG{n}{array}\PYG{p}{(}\PYG{p}{[}\PYG{l+m+mi}{1}\PYG{p}{,} \PYG{l+m+mi}{2}\PYG{p}{,} \PYG{l+m+mi}{1}\PYG{p}{]}\PYG{p}{)}
\PYG{n}{model}\PYG{o}{.}\PYG{n}{setObjective}\PYG{p}{(}\PYG{n}{obj} \PYG{o}{@} \PYG{n}{x}\PYG{p}{,} \PYG{n}{COPT}\PYG{o}{.}\PYG{n}{MINIMIZE}\PYG{p}{)}
\end{sphinxVerbatim}
\end{quote}

\section{Indexing and Slicing}
\label{\detokenize{matrix:indexing-and-slicing}}\label{\detokenize{matrix:chapmatrixindex}}
\sphinxAtStartPar
Multidimensional array objects in COPT (such as \sphinxcode{\sphinxupquote{NdArray}}, \sphinxcode{\sphinxupquote{MVar}},
\sphinxcode{\sphinxupquote{MConstr}}, \sphinxcode{\sphinxupquote{MLinExpr}}, etc.) support a unified set of indexing and slicing
rules for accessing, selecting, or updating elements or sub\sphinxhyphen{}objects.
Supported indexing forms include integer indexing, slice indexing,
and multidimensional tuple indexing.

\subsection{Integer Indexing}
\label{\detokenize{matrix:integer-indexing}}\begin{quote}

\sphinxAtStartPar
\sphinxstylestrong{Synopsis}
\begin{quote}

\sphinxAtStartPar
\sphinxcode{\sphinxupquote{mobj{[}i{]}}}
\end{quote}

\sphinxAtStartPar
\sphinxstylestrong{Description}
\begin{quote}

\sphinxAtStartPar
Integer indexing is used to access a single element in a multidimensional
array object, where \sphinxcode{\sphinxupquote{i}} is an integer index.
\end{quote}

\sphinxAtStartPar
\sphinxstylestrong{Return Value}
\begin{quote}

\sphinxAtStartPar
Returns a scalar value representing the element at the specified index.
The return type corresponds to the underlying element type of the array.
\end{quote}

\sphinxAtStartPar
\sphinxstylestrong{Example}

\begin{sphinxVerbatim}[commandchars=\\\{\}]
\PYG{n}{x} \PYG{o}{=} \PYG{n}{NdArray}\PYG{p}{(}\PYG{n}{args}\PYG{o}{=}\PYG{p}{[}\PYG{l+m+mf}{1.0}\PYG{p}{,} \PYG{l+m+mf}{2.0}\PYG{p}{,} \PYG{l+m+mf}{3.0}\PYG{p}{]}\PYG{p}{,} \PYG{n}{shape}\PYG{o}{=}\PYG{p}{(}\PYG{l+m+mi}{3}\PYG{p}{,}\PYG{p}{)}\PYG{p}{)}
\PYG{n}{v1} \PYG{o}{=} \PYG{n}{x}\PYG{p}{[}\PYG{l+m+mi}{0}\PYG{p}{]}    \PYG{c+c1}{\PYGZsh{} v1 is a scalar with value 1.0}
\PYG{n}{v2} \PYG{o}{=} \PYG{n}{x}\PYG{p}{[}\PYG{o}{\PYGZhy{}}\PYG{l+m+mi}{1}\PYG{p}{]}   \PYG{c+c1}{\PYGZsh{} v2 is a scalar with value 3.0}
\end{sphinxVerbatim}
\end{quote}

\subsection{Slice Indexing}
\label{\detokenize{matrix:slice-indexing}}\begin{quote}

\sphinxAtStartPar
\sphinxstylestrong{Synopsis}
\begin{quote}

\sphinxAtStartPar
\sphinxcode{\sphinxupquote{mobj{[}start:stop:step{]}}}

\sphinxAtStartPar
\sphinxcode{\sphinxupquote{mobj{[}start:stop{]}}}

\sphinxAtStartPar
\sphinxcode{\sphinxupquote{mobj{[}:{]}}}
\end{quote}

\sphinxAtStartPar
\sphinxstylestrong{Description}
\begin{quote}

\sphinxAtStartPar
Slice indexing is used to select a contiguous range of elements along
a given dimension. The parameters \sphinxcode{\sphinxupquote{start}}, \sphinxcode{\sphinxupquote{stop}}, and \sphinxcode{\sphinxupquote{step}}
are all optional and follow the semantics of Python’s built\sphinxhyphen{}in slicing:
\begin{itemize}
\item {} 
\sphinxAtStartPar
\sphinxcode{\sphinxupquote{start}}: starting index (inclusive)

\item {} 
\sphinxAtStartPar
\sphinxcode{\sphinxupquote{stop}}: ending index (exclusive)

\item {} 
\sphinxAtStartPar
\sphinxcode{\sphinxupquote{step}}: step size

\end{itemize}
\end{quote}

\sphinxAtStartPar
\sphinxstylestrong{Return Value}
\begin{quote}

\sphinxAtStartPar
Returns a new multidimensional array object representing the sliced
sub\sphinxhyphen{}array. The element type is the same as that of the original object.
\end{quote}

\sphinxAtStartPar
\sphinxstylestrong{Example}

\begin{sphinxVerbatim}[commandchars=\\\{\}]
\PYG{n}{x} \PYG{o}{=} \PYG{n}{NdArray}\PYG{p}{(}\PYG{n}{args}\PYG{o}{=}\PYG{p}{[}\PYG{l+m+mf}{1.0}\PYG{p}{,} \PYG{l+m+mf}{2.0}\PYG{p}{,} \PYG{l+m+mf}{3.0}\PYG{p}{,} \PYG{l+m+mf}{4.0}\PYG{p}{]}\PYG{p}{,} \PYG{n}{shape}\PYG{o}{=}\PYG{p}{(}\PYG{l+m+mi}{4}\PYG{p}{,}\PYG{p}{)}\PYG{p}{)}
\PYG{n}{y} \PYG{o}{=} \PYG{n}{x}\PYG{p}{[}\PYG{l+m+mi}{1}\PYG{p}{:}\PYG{l+m+mi}{3}\PYG{p}{]}      \PYG{c+c1}{\PYGZsh{} NdArray, [2.0, 3.0]}
\PYG{n}{z} \PYG{o}{=} \PYG{n}{x}\PYG{p}{[}\PYG{o}{\PYGZhy{}}\PYG{l+m+mi}{3}\PYG{p}{:}\PYG{o}{\PYGZhy{}}\PYG{l+m+mi}{1}\PYG{p}{]}    \PYG{c+c1}{\PYGZsh{} NdArray, [2.0, 3.0]}
\PYG{n}{w} \PYG{o}{=} \PYG{n}{x}\PYG{p}{[}\PYG{p}{:}\PYG{p}{]}        \PYG{c+c1}{\PYGZsh{} NdArray, same contents as x}
\end{sphinxVerbatim}
\end{quote}

\subsection{Tuple Indexing}
\label{\detokenize{matrix:tuple-indexing}}\begin{quote}

\sphinxAtStartPar
\sphinxstylestrong{Synopsis}
\begin{quote}

\sphinxAtStartPar
\sphinxcode{\sphinxupquote{mobj{[}i, j{]}}}

\sphinxAtStartPar
\sphinxcode{\sphinxupquote{mobj{[}i, :{]}}}

\sphinxAtStartPar
\sphinxcode{\sphinxupquote{mobj{[}:, j{]}}}

\sphinxAtStartPar
\sphinxcode{\sphinxupquote{mobj{[}i, j:k{]}}}
\end{quote}

\sphinxAtStartPar
\sphinxstylestrong{Description}
\begin{quote}

\sphinxAtStartPar
For multidimensional array objects, integer indexing and slice indexing
can be applied simultaneously across multiple dimensions using a tuple
indexing expression. Each element in the tuple corresponds to one
dimension of the array and may be either an integer index or a slice.
\end{quote}

\sphinxAtStartPar
\sphinxstylestrong{Return Value}
\begin{itemize}
\item {} 
\sphinxAtStartPar
If all dimensions are indexed using integers, a scalar value is returned.

\item {} 
\sphinxAtStartPar
If any dimension uses slice indexing, a multidimensional array object
is returned.

\end{itemize}

\sphinxAtStartPar
\sphinxstylestrong{Example}

\begin{sphinxVerbatim}[commandchars=\\\{\}]
\PYG{n}{x} \PYG{o}{=} \PYG{n}{NdArray}\PYG{p}{(}\PYG{n}{shape}\PYG{o}{=}\PYG{p}{(}\PYG{l+m+mi}{3}\PYG{p}{,} \PYG{l+m+mi}{3}\PYG{p}{)}\PYG{p}{)}

\PYG{n}{v} \PYG{o}{=} \PYG{n}{x}\PYG{p}{[}\PYG{l+m+mi}{0}\PYG{p}{,} \PYG{l+m+mi}{1}\PYG{p}{]}    \PYG{c+c1}{\PYGZsh{} scalar}
\PYG{n}{r} \PYG{o}{=} \PYG{n}{x}\PYG{p}{[}\PYG{l+m+mi}{0}\PYG{p}{,} \PYG{p}{:}\PYG{p}{]}    \PYG{c+c1}{\PYGZsh{} NdArray with shape (3,)}
\PYG{n}{c} \PYG{o}{=} \PYG{n}{x}\PYG{p}{[}\PYG{p}{:}\PYG{p}{,} \PYG{l+m+mi}{1}\PYG{p}{]}    \PYG{c+c1}{\PYGZsh{} NdArray with shape (3,)}
\end{sphinxVerbatim}
\end{quote}

\sphinxAtStartPar
Regarding the matrix modeling method, the COPT Python interface provides multi\sphinxhyphen{}dimensional variables, (linear and convex quadratic) expressions, and matrix constraint classes respectively, and contains related operations. For the methods and specific introductions included, please refer to the corresponding part of the Python API :
\begin{itemize}
\item {} 
\sphinxAtStartPar
Multi\sphinxhyphen{}dimensional variable: {\hyperref[\detokenize{pyapiref:chappyapi-mvar}]{\sphinxcrossref{\DUrole{std,std-ref}{MVar}}}}

\item {} 
\sphinxAtStartPar
Multi\sphinxhyphen{}dimensional linear expression: {\hyperref[\detokenize{pyapiref:chappyapi-mlinexpr}]{\sphinxcrossref{\DUrole{std,std-ref}{MLinExpr}}}}

\item {} 
\sphinxAtStartPar
Multi\sphinxhyphen{}dimensional quadratic expressions: {\hyperref[\detokenize{pyapiref:chappyapi-mquadexpr}]{\sphinxcrossref{\DUrole{std,std-ref}{MQuadExpr}}}}

\item {} 
\sphinxAtStartPar
Matrix constraints: {\hyperref[\detokenize{pyapiref:chappyapi-mconstr}]{\sphinxcrossref{\DUrole{std,std-ref}{MConstr}}}}

\end{itemize}

\sphinxstepscope

\chapter{Multi\sphinxhyphen{}objective Optimization}
\label{\detokenize{multiobj:multi-objective-optimization}}\label{\detokenize{multiobj:chapmultiobj}}\label{\detokenize{multiobj::doc}}
\sphinxAtStartPar
In real\sphinxhyphen{}world applications and decision\sphinxhyphen{}making systems, there is often more than
one objective to optimize. For example, in supply chain management,
one may aim to minimize inventory cost while maximizing order fulfillment rate.
COPT provides multi\sphinxhyphen{}objective optimization functionality
to properly balance the priorities or weights of multiple objectives,
using either a hierarchy method or a weighted\sphinxhyphen{}sum method, to achieve optimal decisions
under multi\sphinxhyphen{}objective scenarios.

\section{Modeling Multiple objectives}
\label{\detokenize{multiobj:modeling-multiple-objectives}}
\sphinxAtStartPar
COPT currently supports linear objective functions for multi\sphinxhyphen{}objective optimization.
The API function \sphinxcode{\sphinxupquote{Model.setObjectiveN()}} is used to define them.

\sphinxAtStartPar
When constructing multi\sphinxhyphen{}objective functions,
the following key objective parameters must be specified:

\phantomsection\label{\detokenize{multiobj:multiobjpriority}}\begin{itemize}
\item {} 
\sphinxAtStartPar
\sphinxcode{\sphinxupquote{MultiObjPriority}}
\begin{quote}

\sphinxAtStartPar
The priority of the objective in multi\sphinxhyphen{}objective optimization.
Higher values indicate higher priority.
This determines the order of optimization in hierarchy method.

\sphinxAtStartPar
\sphinxstylestrong{Default value} 0.0
\end{quote}

\end{itemize}
\phantomsection\label{\detokenize{multiobj:multiobjweight}}\begin{itemize}
\item {} 
\sphinxAtStartPar
\sphinxcode{\sphinxupquote{MultiObjWeight}}
\begin{quote}

\sphinxAtStartPar
The weight of the objective in multi\sphinxhyphen{}objective optimization.
This defines the coefficient of the objective in the weighted\sphinxhyphen{}sum method.

\sphinxAtStartPar
\sphinxstylestrong{Default value} 1.0
\end{quote}

\end{itemize}
\phantomsection\label{\detokenize{multiobj:multiobjabstol}}\begin{itemize}
\item {} 
\sphinxAtStartPar
\sphinxcode{\sphinxupquote{MultiObjAbsTol}}
\begin{quote}

\sphinxAtStartPar
The absolute tolerance allowed for degradation for the specified objective
in multi\sphinxhyphen{}objective optimization.

\sphinxAtStartPar
COPT optimizes objectives sequentially based on their priority:
In multi\sphinxhyphen{}objective MIP, if the current objective’s optimal value is \sphinxcode{\sphinxupquote{z}},
the solution is allowed to degrade \sphinxcode{\sphinxupquote{z}} by at most \sphinxcode{\sphinxupquote{MultiObjAbsTol}} in subsequent groups.
In multi\sphinxhyphen{}objective (continuous) LP, this parameter applies to variables.
To maintain solution quality of higher\sphinxhyphen{}priority objectives,
COPT fixes some variable values from the previous optimal solution based on reduced cost.
More precisely, all variables with absolute reduced costs larger than \sphinxcode{\sphinxupquote{MultiObjAbsTol}} will be fixed to their value in the previous optimal solution.

\sphinxAtStartPar
\sphinxstylestrong{Default value} 1e\sphinxhyphen{}6
\end{quote}

\end{itemize}
\phantomsection\label{\detokenize{multiobj:multiobjreltol}}\begin{itemize}
\item {} 
\sphinxAtStartPar
\sphinxcode{\sphinxupquote{MultiObjRelTol}}
\begin{quote}

\sphinxAtStartPar
The relative tolerance allowed for degradation for the specified objective
in multi\sphinxhyphen{}objective optimization.

\sphinxAtStartPar
COPT optimizes objectives sequentially based on their priority:
In multi\sphinxhyphen{}objective MIP, if the current objective’s optimal value is \sphinxcode{\sphinxupquote{z}},
subsequent optimization phases allow its objective value to deviate within
\sphinxcode{\sphinxupquote{MultiObjRelTol * |z|}}.
In multi\sphinxhyphen{}objective (continuous) linear programming, the degradation
is controlled by \sphinxcode{\sphinxupquote{MultiObjAbsTol}}, and \sphinxcode{\sphinxupquote{MultiObjRelTol}} is ignored.

\sphinxAtStartPar
\sphinxstylestrong{Default value} 0.0
\end{quote}

\end{itemize}

\begin{sphinxadmonition}{note}{Notes}
\begin{itemize}
\item {} 
\sphinxAtStartPar
The above parameters apply only to multiple objectives.

\item {} 
\sphinxAtStartPar
When specifying a multi\sphinxhyphen{}objective function with \sphinxcode{\sphinxupquote{setObjectiveN}},
in addition to specifying the above parameters,
the objective index, objective expression, and optimization sense can also
be specified.

\end{itemize}
\end{sphinxadmonition}

\sphinxAtStartPar
Following are operations on multi\sphinxhyphen{}objective parameters (e.g. COPT Python API):
\begin{enumerate}
\sphinxsetlistlabels{\arabic}{enumi}{enumii}{}{.}%
\item {} 
\sphinxAtStartPar
Use \sphinxcode{\sphinxupquote{Model.setObjParamN(idx, paramname, newval)}}
to set the parameter of a specific objective function by index.

\item {} 
\sphinxAtStartPar
Use \sphinxcode{\sphinxupquote{Model.getObjParamN(idx, paramname)}}
to retrieve the value of a parameter for a specific objective function by index.

\item {} 
\sphinxAtStartPar
Use \sphinxcode{\sphinxupquote{Model.resetObjParamN(idx)}}
to reset the parameters of a specific objective function to its default value.

\end{enumerate}

\section{Solving Multiple objectives}
\label{\detokenize{multiobj:solving-multiple-objectives}}
\sphinxAtStartPar
When solving multi\sphinxhyphen{}objective models, COPT supports two methods: hierarchy method
and weighted\sphinxhyphen{}sum method.

\subsection{Hierarchy Method}
\label{\detokenize{multiobj:hierarchy-method}}
\sphinxAtStartPar
The hierarchy method optimizes objectives in descending order of priority.
COPT solves each objective while ensuring that the quality of higher\sphinxhyphen{}priority
objectives is not degraded.

\sphinxAtStartPar
Users can specify priorities while calling \sphinxcode{\sphinxupquote{Model.setObjectiveN}} or using \sphinxcode{\sphinxupquote{Model.setObjParamN}}
to assign a priority to the specified objective.
The default priority of each objective is 0.0.

\sphinxAtStartPar
Assuming a model with two objectives, \(obj_1\) and \(obj_2\),
with priorities 3 and 2 respectively, COPT will first solve the problem for \(obj_1\),
and then solve for \(obj_2\) within the allowable degradation defined for \(obj_1\).

\subsection{Weighted\sphinxhyphen{}sum Method}
\label{\detokenize{multiobj:weighted-sum-method}}
\sphinxAtStartPar
For objectives with the same priority, COPT uses the weighted\sphinxhyphen{}sum method to combine
them into a single objective using the specified weights.
Weights can be set while calling \sphinxcode{\sphinxupquote{Model.setObjectiveN}} or using \sphinxcode{\sphinxupquote{Model.setObjParamN}}.
The default weight of each objective is 1.0.

\sphinxAtStartPar
Assuming a model includes two objectives at the same priority level:
\(\text{min}\ f_1(x)\) and \(\text{min}\ f_2(x)\),
with weights \(w_1\) and \(w_2\) respectively.
Then, COPT constructs the combined objective below:
\begin{equation}\label{equation:multiobj:multiobj:0}
\begin{split}\min\ w_1 \cdot f_1(x) + w_2 \cdot f_2(x)\end{split}
\end{equation}
\begin{sphinxadmonition}{note}{Note:}
\begin{itemize}
\item {} 
\sphinxAtStartPar
Please be cautious when assigning a negative weight, as it effectively
changes the optimization direction.

\item {} 
\sphinxAtStartPar
Avoid assigning excessively large or small weights. A large weight may cause numeric instability,
while a small weight may result in an insignificant contribution,
potentially ignored due to objective degradation tolerances.

\end{itemize}
\end{sphinxadmonition}

\sphinxAtStartPar
COPT also supports combining the hierarchy method and weighted\sphinxhyphen{}sum method together.
Each objective can have different/same priority and weight.
Objectives are grouped and solved sequentially by priority.
Within the same priority group, objectives are combined using the weighted\sphinxhyphen{}sum method.

\sphinxAtStartPar
Consider a multi\sphinxhyphen{}objective model with four objectives:
\begin{equation}\label{equation:multiobj:multiobj:1}
\begin{split}\begin{aligned}
&\min\ \text{obj}_1 = f_1(x) && \text{(priority = 2,\ weight = 1.0)} \\
&\max\ \text{obj}_2 = f_2(x) && \text{(priority = 1,\ weight = 2.5)} \\
&\max\ \text{obj}_3 = f_3(x) && \text{(priority = 1,\ weight = 1.5)} \\
&\min\ \text{obj}_4 = f_4(x) && \text{(priority = 0,\ weight = 1.0)} \\
\end{aligned}\end{split}
\end{equation}
\sphinxAtStartPar
COPT will first group the objectives by priority (2, 1, 0),
and solve them in descending order.
For objectives with the same priority, a weighted sum will be used.

\sphinxAtStartPar
Priority = 2:
\begin{equation}\label{equation:multiobj:multiobj:2}
\begin{split}\min\ f_1(x)\end{split}
\end{equation}
\sphinxAtStartPar
Priority = 1:
\begin{equation}\label{equation:multiobj:multiobj:3}
\begin{split}\max\ 2.5 \cdot f_2(x) + 1.5 \cdot f_3(x)\end{split}
\end{equation}
\sphinxAtStartPar
Priority = 0:
\begin{equation}\label{equation:multiobj:multiobj:4}
\begin{split}\min\ f_4(x)\end{split}
\end{equation}

\subsection{Objective degradation tolerance}
\label{\detokenize{multiobj:objective-degradation-tolerance}}
\sphinxAtStartPar
When using the hierarchy method for multi\sphinxhyphen{}objective optimization,
the degree to which a lower\sphinxhyphen{}priority objective is allowed to degrade
for higher ones can be adjusted via two parameters:
\sphinxcode{\sphinxupquote{MultiObjRelTol}} (relative degradation tolerance) and \sphinxcode{\sphinxupquote{MultiObjAbsTol}} (absolute degradation tolerance).

\sphinxAtStartPar
1.MIP Model

\sphinxAtStartPar
For MIP models, we suppose an objective \(\min\ z_1\) has an optimal value \(v\),
and the corresponding \sphinxcode{\sphinxupquote{MultiObjAbsTol = e}}.
Then, when optimizing the next lower\sphinxhyphen{}priority objective \(\min\ z_2\),
COPT searches for the optimal solution of \(z_2\) within the feasible region
defined by \(z_2 \leq v + e\).

\sphinxAtStartPar
For MIP models, we suppose an objective \(\min\ z_1\) has an optimal value \(v\),
and the corresponding \sphinxcode{\sphinxupquote{MultiObjRelTol = q}}.
Then, when optimizing the next lower\sphinxhyphen{}priority objective \(\min\ z_2\),
COPT searches for the optimal solution of \(z_2\) within the feasible region
defined by \(z_2 \leq v + q*|v|\).

\sphinxAtStartPar
Note: If both relative and absolute tolerances are set for an objective,
COPT will use the larger one of the two allowed degradation tolerance.

\sphinxAtStartPar
2.LP Model

\sphinxAtStartPar
For continuous LP models, in order to maintain the quality of higher\sphinxhyphen{}priority objectives,
COPT fixes certain variables to their optimal values obtained in the previous priority group.
Whether a variable is fixed will be determined by whether its reduced cost
is zero within the \sphinxcode{\sphinxupquote{MultiObjAbsTol}} tolerance.

\sphinxAtStartPar
Note: \sphinxcode{\sphinxupquote{MultiObjRelTol}} is ignored in LP models.

\section{Retrieving the results of multiple objectives}
\label{\detokenize{multiobj:retrieving-the-results-of-multiple-objectives}}
\sphinxAtStartPar
The optimization results of the multi\sphinxhyphen{}objective model can be queried
via \sphinxcode{\sphinxupquote{Model.getAttrN(idx, attrname)}}
(taking COPT\sphinxhyphen{}Python API as an example, other interfaces work the same way).
\begin{itemize}
\item {} 
\sphinxAtStartPar
\sphinxcode{\sphinxupquote{idx}}: The index of the objective.

\item {} 
\sphinxAtStartPar
\sphinxcode{\sphinxupquote{attrname}}: The attribute name of the objective.

\end{itemize}

\sphinxAtStartPar
Supported attributes for multiple objectives include:
{\hyperref[\detokenize{attribute:hasqobj}]{\sphinxcrossref{\DUrole{std,std-ref}{HasQObj}}}}, {\hyperref[\detokenize{attribute:hasnlobj}]{\sphinxcrossref{\DUrole{std,std-ref}{HasNLObj}}}}, {\hyperref[\detokenize{attribute:lpobjval}]{\sphinxcrossref{\DUrole{std,std-ref}{LpObjval}}}},
{\hyperref[\detokenize{information:bestobj}]{\sphinxcrossref{\DUrole{std,std-ref}{BestObj}}}}, {\hyperref[\detokenize{attribute:objsense}]{\sphinxcrossref{\DUrole{std,std-ref}{ObjSense}}}}, {\hyperref[\detokenize{attribute:objconst}]{\sphinxcrossref{\DUrole{std,std-ref}{ObjConst}}}}.

\sphinxstepscope

\chapter{Logging}
\label{\detokenize{logging:logging}}\label{\detokenize{logging::doc}}
\sphinxAtStartPar
Logging related parameters and functions are essential for users to control the display of solving logs of COPT.
This chapter provides an interpretation of loggings for different algorithms, including the following sections:
\begin{itemize}
\item {} 
\sphinxAtStartPar
{\hyperref[\detokenize{logging:chaplogging-param}]{\sphinxcrossref{\DUrole{std,std-ref}{Parameters and Functions for Logging}}}}

\item {} 
\sphinxAtStartPar
{\hyperref[\detokenize{logging:chaplogging-basic}]{\sphinxcrossref{\DUrole{std,std-ref}{Basic Information Section}}}}

\item {} 
\sphinxAtStartPar
{\hyperref[\detokenize{logging:chaplogging-simplex}]{\sphinxcrossref{\DUrole{std,std-ref}{Simplex Logging}}}}

\item {} 
\sphinxAtStartPar
{\hyperref[\detokenize{logging:chaplogging-barrier}]{\sphinxcrossref{\DUrole{std,std-ref}{Barrier Logging}}}}

\item {} 
\sphinxAtStartPar
{\hyperref[\detokenize{logging:chaplogging-branch}]{\sphinxcrossref{\DUrole{std,std-ref}{Branch\sphinxhyphen{}and\sphinxhyphen{}Cut Logging}}}}

\item {} 
\sphinxAtStartPar
{\hyperref[\detokenize{logging:chaplogging-gpu}]{\sphinxcrossref{\DUrole{std,std-ref}{First\sphinxhyphen{}order Method (PDLP) Logging for GPU Solver}}}}

\end{itemize}

\section{Parameters and Functions for Logging}
\label{\detokenize{logging:parameters-and-functions-for-logging}}\label{\detokenize{logging:chaplogging-param}}
\sphinxAtStartPar
Users can control whether to show logs and how detailed they are by setting logging\sphinxhyphen{}related arguments.
\begin{itemize}
\item {} 
\sphinxAtStartPar
\sphinxcode{\sphinxupquote{Logging}}
\begin{quote}

\sphinxAtStartPar
Integer argument.

\sphinxAtStartPar
Whether to print optimization logs.

\sphinxAtStartPar
\sphinxstylestrong{Default:} 1

\sphinxAtStartPar
\sphinxstylestrong{Possible values:}
\begin{quote}

\sphinxAtStartPar
0: Do not print optimization logs.

\sphinxAtStartPar
1: Print optimization logs.
\end{quote}
\end{quote}

\item {} 
\sphinxAtStartPar
\sphinxcode{\sphinxupquote{LogLevel}}
\begin{quote}

\sphinxAtStartPar
Integer argument.

\sphinxAtStartPar
Controls the level of detail in the optimization logs.

\sphinxAtStartPar
\sphinxstylestrong{Default:} 2

\sphinxAtStartPar
\sphinxstylestrong{Possible values:}
\begin{quote}

\sphinxAtStartPar
2: Print basic optimization logs.

\sphinxAtStartPar
3: In addition to the basic logs, print memory usage information (for MIP problems).
\end{quote}
\end{quote}

\item {} 
\sphinxAtStartPar
\sphinxcode{\sphinxupquote{LogToConsole}}
\begin{quote}

\sphinxAtStartPar
Integer argument.

\sphinxAtStartPar
Whether to print optimization logs to the console.

\sphinxAtStartPar
\sphinxstylestrong{Default:} 1

\sphinxAtStartPar
\sphinxstylestrong{Possible values:}
\begin{quote}

\sphinxAtStartPar
0: Do not print logs to the console.

\sphinxAtStartPar
1: Print logs to the console.
\end{quote}
\end{quote}

\end{itemize}

\sphinxAtStartPar
COPT provides operations related to logs, such as setting the logging file. COPT provides functions to set the logging file,
write the optimization logs into the specified file (with the file name suffix \sphinxcode{\sphinxupquote{.log}} ), allowing users to save and review the logs.
The functions for different programming interfaces are shown in \hyperref[\detokenize{logging:copttab-setlogfile}]{Table \ref{\detokenize{logging:copttab-setlogfile}}}:

\begin{savenotes}\sphinxattablestart
\sphinxthistablewithglobalstyle
\centering
\sphinxcapstartof{table}
\sphinxthecaptionisattop
\sphinxcaption{Functions for setting logging files in different programming interfaces}\label{\detokenize{logging:copttab-setlogfile}}
\sphinxaftertopcaption
\begin{tabular}[t]{|\X{10}{35}|\X{25}{35}|}
\sphinxtoprule
\sphinxstyletheadfamily 
\sphinxAtStartPar
Programming Interface
&\sphinxstyletheadfamily 
\sphinxAtStartPar
Functions for setting logging files
\\
\sphinxmidrule
\sphinxtableatstartofbodyhook
\sphinxAtStartPar
C
&
\sphinxAtStartPar
\sphinxcode{\sphinxupquote{COPT\_SetLogFile}}
\\
\sphinxhline
\sphinxAtStartPar
C++
&
\sphinxAtStartPar
\sphinxcode{\sphinxupquote{Model::SetLogFile()}}
\\
\sphinxhline
\sphinxAtStartPar
C\#
&
\sphinxAtStartPar
\sphinxcode{\sphinxupquote{Model.SetLogFile()}}
\\
\sphinxhline
\sphinxAtStartPar
Java
&
\sphinxAtStartPar
\sphinxcode{\sphinxupquote{Model.setLogFile()}}
\\
\sphinxhline
\sphinxAtStartPar
Python
&
\sphinxAtStartPar
\sphinxcode{\sphinxupquote{Model.setLogFile()}}
\\
\sphinxbottomrule
\end{tabular}
\sphinxtableafterendhook\par
\sphinxattableend\end{savenotes}

\sphinxAtStartPar
\sphinxstylestrong{Note:} When calling these functions, users should specify the file name for saving logs using the \sphinxcode{\sphinxupquote{logfile}} parameter.

\section{Basic Information Section}
\label{\detokenize{logging:basic-information-section}}\label{\detokenize{logging:chaplogging-basic}}
\sphinxAtStartPar
COPT outputs basic information before starting the solving process,
depending on the problem types. The following information is typically displayed:

\sphinxAtStartPar
This section prints the overview information of the model and optimization task,
such as model size and presolve results, and the numerical characteristics,
helping users quickly understand the problem structure and numerical status.

\subsection{Solver Environment and Model Overview}
\label{\detokenize{logging:solver-environment-and-model-overview}}
\sphinxAtStartPar
At the beginning of the log,
COPT prints an overview of the solver environment and the model,
including the model fingerprint (\sphinxcode{\sphinxupquote{fingerprint}}), COPT version, platform, and problem type.

\sphinxAtStartPar
Example:
\begin{sphinxalltt}
Using Cardinal Optimizer v8.0.1 on Windows
Model fingerprint: 2c27ab28
Hardware has 64 cores and 128 threads. Using instruction set X86\_AVX512\_E1 (14)
Minimizing a MIP problem
\end{sphinxalltt}

\sphinxAtStartPar
Here, \sphinxcode{\sphinxupquote{Model fingerprint}} is a unique identifier of the current model. Hardware information includes CPU cores (\sphinxcode{\sphinxupquote{cores}}) and threads (\sphinxcode{\sphinxupquote{threads}}).
COPT automatically detects the problem type and optimization sense, and then selects appropriate algorithms to solve it, for example:

\sphinxAtStartPar
\sphinxcode{\sphinxupquote{Minimizing an LP problem}}, \sphinxcode{\sphinxupquote{Minimizing a MIP problem}}, \sphinxcode{\sphinxupquote{Minimizing an SDP problem}}, etc.

\subsection{Model Size and Presolve}
\label{\detokenize{logging:model-size-and-presolve}}
\sphinxAtStartPar
First, COPT prints the size of the original model, for example:

\begin{sphinxVerbatim}[commandchars=\\\{\}]
The original problem has:
    404 rows, 1200 columns and 2598 non\PYGZhy{}zero elements
    200 binaries and 1000 integers
\end{sphinxVerbatim}

\sphinxAtStartPar
For LP problems, the size includes:
\begin{itemize}
\item {} 
\sphinxAtStartPar
Constraints (\sphinxcode{\sphinxupquote{rows}})

\item {} 
\sphinxAtStartPar
Variables (\sphinxcode{\sphinxupquote{columns}})

\item {} 
\sphinxAtStartPar
Nonzeros in the constraint matrix (\sphinxcode{\sphinxupquote{non\sphinxhyphen{}zero elements}})

\end{itemize}

\sphinxAtStartPar
For MIP problems, it additionally includes integer\sphinxhyphen{}variable statistics:
\begin{itemize}
\item {} 
\sphinxAtStartPar
Binary variables (\sphinxcode{\sphinxupquote{binaries}})

\item {} 
\sphinxAtStartPar
General integer variables (\sphinxcode{\sphinxupquote{integers}})

\end{itemize}

\begin{sphinxadmonition}{note}{Note}
\begin{itemize}
\item {} 
\sphinxAtStartPar
For MIP problems, the variable count (\sphinxcode{\sphinxupquote{columns}}) includes all variables in the model: continuous, integer, and binary.

\item {} 
\sphinxAtStartPar
In theory, binaries are a special case of integer variables, but here they are reported separately.

\end{itemize}
\end{sphinxadmonition}

\sphinxAtStartPar
For nonlinear problems, the log also reports corresponding nonlinear structures, for example:
\begin{itemize}
\item {} 
\sphinxAtStartPar
Semidefinite variables (\sphinxcode{\sphinxupquote{PSD columns}})

\item {} 
\sphinxAtStartPar
Quadratic terms in the objective (\sphinxcode{\sphinxupquote{quadratic objective elements}})

\item {} 
\sphinxAtStartPar
Quadratic constraints (\sphinxcode{\sphinxupquote{quadratic constraints}})

\item {} 
\sphinxAtStartPar
Second\sphinxhyphen{}Order cones (\sphinxcode{\sphinxupquote{SOC rows}})

\item {} 
\sphinxAtStartPar
Exponential cones (\sphinxcode{\sphinxupquote{exponential cones}})

\end{itemize}

\sphinxAtStartPar
With default settings, before the optimization starts,
COPT presolves to transform or reduce the original model to improve quality; the algorithm then works on the presolved model.
The presolve part prints the model sizes. Example:

\begin{sphinxVerbatim}[commandchars=\\\{\}]
Presolving the problem
The presolved problem has:
    373 rows, 1169 columns and 2505 non\PYGZhy{}zero elements
    369 binaries and 800 integers
\end{sphinxVerbatim}

\subsection{Problem Numerical Characteristics}
\label{\detokenize{logging:problem-numerical-characteristics}}
\sphinxAtStartPar
This part prints key numerical characteristics of the model to help assess numerical issue
(e.g., whether coefficient ranges span too many orders of magnitude). The items are:
\begin{itemize}
\item {} 
\sphinxAtStartPar
\sphinxcode{\sphinxupquote{Range of matrix coefficients}}: Range of the constraint matrix coefficients

\item {} 
\sphinxAtStartPar
\sphinxcode{\sphinxupquote{Range of rhs coefficients}}: Range of right\sphinxhyphen{}hand\sphinxhyphen{}side coefficients

\item {} 
\sphinxAtStartPar
\sphinxcode{\sphinxupquote{Range of bound coefficients}}: Range of variable bounds

\item {} 
\sphinxAtStartPar
\sphinxcode{\sphinxupquote{Range of cost coefficients}}: Range of objective function coefficients

\item {} 
\sphinxAtStartPar
\sphinxcode{\sphinxupquote{Density of cost}}: Density of the objective coefficients (ratio of nonzeros)

\end{itemize}

\sphinxAtStartPar
Example:

\begin{sphinxVerbatim}[commandchars=\\\{\}]
Problem info:
    Range of matrix coefficients:    [5e\PYGZhy{}01,4e+00]
    Range of rhs coefficients:       [7e+02,2e+04]
    Range of bound coefficients:     [8e+02,3e+03]
    Range of cost coefficients:      [6e\PYGZhy{}02,3e\PYGZhy{}01]
\end{sphinxVerbatim}

\section{Simplex Logging}
\label{\detokenize{logging:simplex-logging}}\label{\detokenize{logging:chaplogging-simplex}}
\sphinxAtStartPar
Based on different stages of the optimization process, the logs for Simplex Method can be divided into two parts:
\begin{enumerate}
\sphinxsetlistlabels{\arabic}{enumi}{enumii}{}{.}%
\item {} 
\sphinxAtStartPar
Simplex iteration process

\item {} 
\sphinxAtStartPar
Summary of optimization results

\end{enumerate}

\sphinxAtStartPar
This section uses the logs for example problem \sphinxhref{https://netlib.org/lp/data/afiro}{afiro.mps} to
interpret the information in the Simplex Method logging.

\subsection{Simplex Iteration Process}
\label{\detokenize{logging:simplex-iteration-process}}
\sphinxAtStartPar
This part of the logs provides relevant information about the iteration process using the Simplex method.

\begin{sphinxVerbatim}[commandchars=\\\{\}]
Starting the simplex solver using up to 8 threads

Method   Iteration           Objective  Primal.NInf   Dual.NInf        Time
Dual             0   \PYGZhy{}4.8553789460e+02            3           0       0.00s
Dual             3   \PYGZhy{}4.6476735494e+02            0           0       0.00s
Postsolving
Dual             3   \PYGZhy{}4.6475314286e+02            0           0       0.00s
\end{sphinxVerbatim}

\sphinxAtStartPar
Here, the first line indicates that the current optimization algorithm is the Simplex method, and it uses 8 threads (\sphinxcode{\sphinxupquote{threads}}) for computation.

\sphinxAtStartPar
The subsequent lines represent the simplex iteration process with 6 columns:
\begin{itemize}
\item {} 
\sphinxAtStartPar
\sphinxcode{\sphinxupquote{Method}}: The optimization algorithm used, where \sphinxcode{\sphinxupquote{"Dual"}} represents the dual simplex method.

\item {} 
\sphinxAtStartPar
\sphinxcode{\sphinxupquote{Iteration}}: The number of iterations.

\item {} 
\sphinxAtStartPar
\sphinxcode{\sphinxupquote{Objective}}: The objective function value.

\item {} 
\sphinxAtStartPar
\sphinxcode{\sphinxupquote{Primal.NInf}}: The number of primal infeasibilities in the primal problem.

\item {} 
\sphinxAtStartPar
\sphinxcode{\sphinxupquote{Dual.NInf}}: The number of dual infeasibilities in the dual problem.

\item {} 
\sphinxAtStartPar
\sphinxcode{\sphinxupquote{Time}}: The time taken for solving (in seconds).

\end{itemize}

\subsection{Solution Summary}
\label{\detokenize{logging:solution-summary}}
\sphinxAtStartPar
This part of the logs summarizes the results and the iteration process of the Simplex method after completing the solving process.

\begin{sphinxVerbatim}[commandchars=\\\{\}]
Solving finished
Status: Optimal  Objective: \PYGZhy{}4.6475314286e+02  Iterations: 3  Time: 0.00s
\end{sphinxVerbatim}

\sphinxAtStartPar
The included information consists of:
\begin{itemize}
\item {} 
\sphinxAtStartPar
Solving status (\sphinxcode{\sphinxupquote{Status}}): If the model has an optimal solution, it is \sphinxcode{\sphinxupquote{Optimal}}.

\item {} 
\sphinxAtStartPar
Objective function value (\sphinxcode{\sphinxupquote{Objective}}): If the model has an optimal solution, \sphinxcode{\sphinxupquote{Objective}} displays the optimal objective function value.

\item {} 
\sphinxAtStartPar
Total number of iterations (\sphinxcode{\sphinxupquote{Iterations}}).

\item {} 
\sphinxAtStartPar
Total solving time (\sphinxcode{\sphinxupquote{Time}}).

\end{itemize}

\sphinxAtStartPar
If the model is infeasible, the corresponding log output is as follows:

\begin{sphinxVerbatim}[commandchars=\\\{\}]
Solving finished
Status: Infeasible  Objective: \PYGZhy{}  Iterations: 2  Time: 0.00s
\end{sphinxVerbatim}

\section{Barrier Logging}
\label{\detokenize{logging:barrier-logging}}\label{\detokenize{logging:chaplogging-barrier}}
\sphinxAtStartPar
By solution phase, the Barrier log consists of two parts:
\begin{enumerate}
\sphinxsetlistlabels{\arabic}{enumi}{enumii}{}{.}%
\item {} 
\sphinxAtStartPar
Barrier solution Process

\item {} 
\sphinxAtStartPar
Solution summary

\end{enumerate}

\sphinxAtStartPar
\sphinxstylestrong{Note:} By setting the optimization parameter \sphinxcode{\sphinxupquote{"LpMethod = 2"}}, you can choose the Barrier method as the algorithm
for solving linear programming problems.

\sphinxAtStartPar
Again, we take \sphinxhref{https://netlib.org/lp/data/afiro}{afiro.mps} as an example to explain the Barrier log for LP.

\subsection{Barrier Iteration Process}
\label{\detokenize{logging:barrier-iteration-process}}
\sphinxAtStartPar
First, the log prints numerical information related to Barrier:

\begin{sphinxVerbatim}[commandchars=\\\{\}]
Starting barrier solver using 64 threads

Problem info:
Dualized in presolve:            No
Range of matrix coefficients:    [4e\PYGZhy{}01,4e+00]
Range of rhs coefficients:       [8e+01,3e+02]
Range of bound coefficients:     [4e+01,1e+02]
Range of cost coefficients:      [2e\PYGZhy{}01,2e+00]

Factor info:
Number of free columns:          0
Number of dense columns:         0
Number of matrix entries:        2.800e+01
Number of factor entries:        2.800e+01
Number of factor flops:          1.140e+02
\end{sphinxVerbatim}

\sphinxAtStartPar
The first line shows the algorithm (Barrier) and thread count (\sphinxcode{\sphinxupquote{threads}}). Then:
\begin{itemize}
\item {} 
\sphinxAtStartPar
\sphinxcode{\sphinxupquote{Problem info}} includes whether the model is dualized in presolve, and the ranges of matrix, RHS, bounds, and costs.

\item {} 
\sphinxAtStartPar
\sphinxcode{\sphinxupquote{Factor info}} includes the number of free columns, dense columns, and factorization\sphinxhyphen{}related statistics.

\end{itemize}

\sphinxAtStartPar
Next comes the Barrier iteration table, which reports iteration (\sphinxcode{\sphinxupquote{Iter}}), objectives, and time:
\begin{itemize}
\item {} 
\sphinxAtStartPar
\sphinxcode{\sphinxupquote{Iter}}: Iteration number

\item {} 
\sphinxAtStartPar
\sphinxcode{\sphinxupquote{Primal.Obj}}: Primal objective

\item {} 
\sphinxAtStartPar
\sphinxcode{\sphinxupquote{Dual.Obj}}: Dual objective

\item {} 
\sphinxAtStartPar
\sphinxcode{\sphinxupquote{Compl}}: Complementarity violation

\item {} 
\sphinxAtStartPar
\sphinxcode{\sphinxupquote{Primal.Inf}}: Primal infeasibility

\item {} 
\sphinxAtStartPar
\sphinxcode{\sphinxupquote{Dual.Inf}}: Dual infeasibility

\item {} 
\sphinxAtStartPar
\sphinxcode{\sphinxupquote{Time}}: Elapsed time (s)

\end{itemize}

\begin{sphinxVerbatim}[commandchars=\\\{\}]
Iter       Primal.Obj         Dual.Obj      Compl  Primal.Inf  Dual.Inf    Time
   0  +2.07010046e+01  \PYGZhy{}4.97632246e+02   5.89e+03    4.50e+02  2.65e+00   0.02s
   1  \PYGZhy{}1.18912241e+02  \PYGZhy{}5.91808560e+02   7.58e+02    3.36e+01  1.61e\PYGZhy{}01   0.02s
   2  \PYGZhy{}3.98096520e+02  \PYGZhy{}4.77597371e+02   2.28e+02    9.32e+00  7.45e\PYGZhy{}02   0.02s
   3  \PYGZhy{}4.55223227e+02  \PYGZhy{}4.68222895e+02   1.86e+01    3.60e\PYGZhy{}01  2.63e\PYGZhy{}03   0.02s
   4  \PYGZhy{}4.64587726e+02  \PYGZhy{}4.64803786e+02   2.52e\PYGZhy{}01    7.80e\PYGZhy{}03  7.93e\PYGZhy{}06   0.02s
   5  \PYGZhy{}4.64753143e+02  \PYGZhy{}4.64753143e+02   3.11e\PYGZhy{}07    7.80e\PYGZhy{}09  1.56e\PYGZhy{}11   0.02s
\end{sphinxVerbatim}

\subsection{Barrier Summary}
\label{\detokenize{logging:barrier-summary}}
\sphinxAtStartPar
Key items include solution status, primal/dual optimal objectives, etc.

\begin{sphinxVerbatim}[commandchars=\\\{\}]
Barrier status:                  OPTIMAL
Primal objective:                \PYGZhy{}4.64753143e+02
Dual objective:                  \PYGZhy{}4.64753143e+02
Duality gap (abs/rel):           2.61e\PYGZhy{}07 / 5.63e\PYGZhy{}10
Primal infeasibility (abs/rel):  7.80e\PYGZhy{}09 / 2.60e\PYGZhy{}11
Dual infeasibility (abs/rel):    1.56e\PYGZhy{}11 / 6.99e\PYGZhy{}12
\end{sphinxVerbatim}

\sphinxAtStartPar
\sphinxstylestrong{Crossover}

\begin{sphinxVerbatim}[commandchars=\\\{\}]
Starting crossover using up to 8 threads

    1 primal pushes remaining      0.03s
    0 primal pushes remaining      0.03s
    1 dual pushes remaining        0.03s
    0 dual pushes remaining        0.03s

Method   Iteration           Objective  Primal.NInf   Dual.NInf        Time
Dual             1   \PYGZhy{}4.6475314286e+02            0           0       0.03s
Postsolving
Dual             1   \PYGZhy{}4.6475314286e+02            0           0       0.03s
\end{sphinxVerbatim}

\subsection{Solution Summary}
\label{\detokenize{logging:id2}}
\sphinxAtStartPar
Printed after solving finishes; it summarizes the model solution and the (final) crossover cleanup.

\begin{sphinxVerbatim}[commandchars=\\\{\}]
Solving finished
Status: Optimal  Objective: \PYGZhy{}4.6475314286e+02  Iterations: 1  Time: 0.03s
\end{sphinxVerbatim}

\sphinxAtStartPar
Same as the Simplex log, it includes:
\begin{itemize}
\item {} 
\sphinxAtStartPar
\sphinxcode{\sphinxupquote{Status}}: If the model is solved to optimality, the status will be \sphinxcode{\sphinxupquote{Optimal}}

\item {} 
\sphinxAtStartPar
\sphinxcode{\sphinxupquote{Objective}}: If solved to optimality, the optimal objective value

\item {} 
\sphinxAtStartPar
\sphinxcode{\sphinxupquote{Iterations}}: Total iterations

\item {} 
\sphinxAtStartPar
\sphinxcode{\sphinxupquote{Time}}: Total solving time

\end{itemize}

\section{Branch\sphinxhyphen{}and\sphinxhyphen{}Cut Logging}
\label{\detokenize{logging:branch-and-cut-logging}}\label{\detokenize{logging:chaplogging-branch}}
\sphinxAtStartPar
By solution phase, the Branch\sphinxhyphen{}and\sphinxhyphen{}Cut log consists of two parts:
\begin{enumerate}
\sphinxsetlistlabels{\arabic}{enumi}{enumii}{}{.}%
\item {} 
\sphinxAtStartPar
Search Process

\item {} 
\sphinxAtStartPar
Solution Summary

\end{enumerate}

\sphinxAtStartPar
We use the \sphinxcode{\sphinxupquote{cutstock.mps.gz}} example (available under \sphinxcode{\sphinxupquote{/examples/data}} in the COPT installation) to explain the MIP log.

\subsection{Branch\sphinxhyphen{}and\sphinxhyphen{}Cut Search Process}
\label{\detokenize{logging:branch-and-cut-search-process}}
\sphinxAtStartPar
This section of the log provides information about the Branch\sphinxhyphen{}and\sphinxhyphen{}Cut search process.

\begin{sphinxVerbatim}[commandchars=\\\{\}]
Starting the MIP solver with 8 threads and 32 tasks

    Nodes    Active  LPit/n  IntInf     BestBound  BestSolution    Gap   Time
        0         1      \PYGZhy{}\PYGZhy{}       0  3.100000e+01            \PYGZhy{}\PYGZhy{}    Inf  0.05s
H       0         1      \PYGZhy{}\PYGZhy{}       0  3.100000e+01  6.800000e+01  54.4\PYGZpc{}  0.70s
H       0         1      \PYGZhy{}\PYGZhy{}       0  3.100000e+01  6.600000e+01  53.0\PYGZpc{}  0.70s
H       0         1      \PYGZhy{}\PYGZhy{}       0  3.100000e+01  6.500000e+01  52.3\PYGZpc{}  0.71s
        0         1      \PYGZhy{}\PYGZhy{}      86  5.591304e+01  6.500000e+01  14.0\PYGZpc{}  0.72s
H       0         1      \PYGZhy{}\PYGZhy{}      86  5.591304e+01  6.200000e+01  9.82\PYGZpc{}  0.74s
        1         2     0.0      86  5.591304e+01  6.200000e+01  9.82\PYGZpc{}  0.76s
H       1         1    2129       6  6.000000e+01  6.000000e+01  0.00\PYGZpc{}  0.87s
        2         0    1064       6  6.000000e+01  6.000000e+01  0.00\PYGZpc{}  0.87s
        2         0    1064       6  6.000000e+01  6.000000e+01  0.00\PYGZpc{}  0.87s
\end{sphinxVerbatim}

\sphinxAtStartPar
\sphinxstylestrong{Note:} For brevity, only a portion of the log is shown for explanation.

\sphinxAtStartPar
By meaning, the columns can be grouped as follows; we explain each group below:
\begin{itemize}
\item {} 
\sphinxAtStartPar
Node search information (columns 1\sphinxhyphen{}4)

\item {} 
\sphinxAtStartPar
Feasible solution interval information (columns 5\sphinxhyphen{}7)

\item {} 
\sphinxAtStartPar
Solving time information (column 8)

\end{itemize}

\sphinxAtStartPar
\sphinxstylestrong{Node search information (columns 1\sphinxhyphen{}4):}
\begin{itemize}
\item {} 
\sphinxAtStartPar
\sphinxcode{\sphinxupquote{Nodes}}: Number of nodes searched

\item {} 
\sphinxAtStartPar
\sphinxcode{\sphinxupquote{Active}}: Number of leaf nodes yet to be searched

\item {} 
\sphinxAtStartPar
\sphinxcode{\sphinxupquote{LPit/n}}: Average simplex iterations per node

\item {} 
\sphinxAtStartPar
\sphinxcode{\sphinxupquote{IntInf}}: Number of integer variables that are fractional in the current LP relaxation

\end{itemize}

\sphinxAtStartPar
\sphinxstylestrong{Feasible solution interval information (columns 5\sphinxhyphen{}7):}
\begin{itemize}
\item {} 
\sphinxAtStartPar
\sphinxcode{\sphinxupquote{BestBound}}: Current best bound on the objective

\item {} 
\sphinxAtStartPar
\sphinxcode{\sphinxupquote{BestSolution}}: Current incumbent objective value

\item {} 
\sphinxAtStartPar
\sphinxcode{\sphinxupquote{Gap}}: Relative gap between bound and incumbent. If it drops below the \sphinxcode{\sphinxupquote{RelGap}} threshold, solving process stops

\end{itemize}

\sphinxAtStartPar
\sphinxstylestrong{Solving time information (column 8):}
\begin{itemize}
\item {} 
\sphinxAtStartPar
\sphinxcode{\sphinxupquote{Time}}: Total solving time

\end{itemize}

\begin{sphinxadmonition}{note}{Notes}
\begin{itemize}
\item {} 
\sphinxAtStartPar
The prefix before the first column (\sphinxcode{\sphinxupquote{H/*}}) indicates a new incumbent solution:
\begin{itemize}
\item {} 
\sphinxAtStartPar
\sphinxcode{\sphinxupquote{H}}: Found by a heuristic

\item {} 
\sphinxAtStartPar
\sphinxcode{\sphinxupquote{*}}: Found by branching (subproblem) search

\end{itemize}

\item {} 
\sphinxAtStartPar
Sometimes \sphinxcode{\sphinxupquote{Nodes}} remains 0 for a while, which means COPT is still processing the root node—typically generating cuts or running heuristics to obtain a good incumbent, in order to reduce subsequent search.

\end{itemize}
\end{sphinxadmonition}

\subsection{Solution Summary}
\label{\detokenize{logging:id3}}
\sphinxAtStartPar
Printed after solving finishes, it summarizes the final MIP status and the search effort, including the model solution and search workload.

\begin{sphinxVerbatim}[commandchars=\\\{\}]
Best solution   : 60.000000000
Best bound      : 60.000000000
Best gap        : 0.0000\PYGZpc{}
Solve time      : 0.87
Solve node      : 2
MIP status      : solved
Solution status : integer optimal (relative gap limit 0.0001)

Violations    :     absolute     relative
bounds        :            0            0
rows          :            0            0
integrality   :            0
\end{sphinxVerbatim}

\sphinxAtStartPar
Solution summary includes:
\begin{itemize}
\item {} 
\sphinxAtStartPar
Best objective (\sphinxcode{\sphinxupquote{Best solution}})

\item {} 
\sphinxAtStartPar
Best bound (\sphinxcode{\sphinxupquote{Best bound}})

\item {} 
\sphinxAtStartPar
Best gap (\sphinxcode{\sphinxupquote{Best gap}})

\item {} 
\sphinxAtStartPar
Solution status (\sphinxcode{\sphinxupquote{Solution status}})

\end{itemize}

\sphinxAtStartPar
Search workload includes:
\begin{itemize}
\item {} 
\sphinxAtStartPar
Solve time (\sphinxcode{\sphinxupquote{Solve time}})

\item {} 
\sphinxAtStartPar
Number of explored nodes (\sphinxcode{\sphinxupquote{Solve node}})

\end{itemize}

\sphinxAtStartPar
The \sphinxcode{\sphinxupquote{Violations}} block reports the satisfaction of bounds and constraints at the solution, including:
\begin{itemize}
\item {} 
\sphinxAtStartPar
Violations in variable bounds (\sphinxcode{\sphinxupquote{bounds}}) and constraints (\sphinxcode{\sphinxupquote{rows}})

\item {} 
\sphinxAtStartPar
Integrality violations (\sphinxcode{\sphinxupquote{integrality}})

\end{itemize}

\section{First\sphinxhyphen{}order Method (PDLP) Logging for GPU Solver}
\label{\detokenize{logging:first-order-method-pdlp-logging-for-gpu-solver}}\label{\detokenize{logging:chaplogging-gpu}}
\sphinxAtStartPar
For Linear Programming Problems, if the PDLP method is selected (setting the parameter \sphinxcode{\sphinxupquote{LpMethod=6}} ),
GPU solver can be enabled.  (The machine needs the compatible GPU and the necessary CUDA library needs to be configured).

\sphinxAtStartPar
The logging for GPU solver can be divided into the following sections, which are slightly different from the CPU solver,
with the main distinction in the second section:
\begin{enumerate}
\sphinxsetlistlabels{\arabic}{enumi}{enumii}{}{.}%
\item {} 
\sphinxAtStartPar
GPU hardware information of the machine

\item {} 
\sphinxAtStartPar
First\sphinxhyphen{}order method PDLP iteration process

\item {} 
\sphinxAtStartPar
Crossover section

\item {} 
\sphinxAtStartPar
Summary of optimization results

\end{enumerate}

\sphinxAtStartPar
Taking the \sphinxcode{\sphinxupquote{"thk\_63"}} instance from public LP benchmark as an example, the following is the solving log for the GPU solver.

\subsection{GPU Hardware Information}
\label{\detokenize{logging:gpu-hardware-information}}
\sphinxAtStartPar
This section outputs information about the current machine’s GPU information and CUDA version.

\begin{sphinxVerbatim}[commandchars=\\\{\}]
Hardware has 1 supported GPU device with CUDA 12.3
  GPU 0: NVIDIA GeForce RTX 4090 (CUDA capability 8.9)
\end{sphinxVerbatim}

\begin{sphinxadmonition}{note}{Notes}
\begin{enumerate}
\sphinxsetlistlabels{\arabic}{enumi}{enumii}{}{.}%
\item {} 
\sphinxAtStartPar
The \sphinxcode{\sphinxupquote{"CUDA 12.3"}} mentioned in the log refers to the highest version of the CUDA Toolkit supported by the currently installed CUDA driver.

\item {} 
\sphinxAtStartPar
COPT’s GPU solver requires a minimum CUDA Toolkit version of 12.0.1.

\end{enumerate}
\end{sphinxadmonition}

\subsection{First\sphinxhyphen{}order Method(PDLP) Iteration Process}
\label{\detokenize{logging:first-order-method-pdlp-iteration-process}}
\sphinxAtStartPar
This section outputs the solving iteration process of the PDLP method,
and the summary after the completion, including the iteration count of PDLP,
and the optimal objective values, dual gaps for primal and dual problems.

\begin{sphinxVerbatim}[commandchars=\\\{\}]
Starting PDLP solver on GPU 0

Iterations      Primal.Obj        Dual.Obj         Gap  Primal.Inf  Dual.Inf    Time
        0  +6.00000000e+00  +6.00000000e+00  +0.00e+00    7.87e+00  0.00e+00  21.63s
     4000  +1.95436674e+03  \PYGZhy{}1.77166004e+03  +3.73e+03    2.09e\PYGZhy{}02  0.00e+00  33.37s
     8000  +1.90433201e+03  +1.55817851e+03  +3.46e+02    1.51e\PYGZhy{}02  0.00e+00  44.75s
    12000  +1.87801607e+03  +1.85689627e+03  +2.11e+01    1.74e\PYGZhy{}02  0.00e+00  56.27s
    16000  +1.86810632e+03  +1.86897715e+03  +8.71e\PYGZhy{}01    4.92e\PYGZhy{}03  0.00e+00  67.72s
    20000  +1.87022842e+03  +1.86994685e+03  +2.82e\PYGZhy{}01    3.42e\PYGZhy{}03  0.00e+00  79.18s
    23640  +1.87099459e+03  +1.87099144e+03  +3.15e\PYGZhy{}03    4.69e\PYGZhy{}05  0.00e+00  89.68s

PDLP status:                     OPTIMAL
PDLP iterations:                 23640
Primal objective:                1.87099459e+03
Dual objective:                  1.87099144e+03
Primal infeasibility (abs/rel):  4.69e\PYGZhy{}05 / 6.20e\PYGZhy{}07
Dual infeasibility (abs/rel):    0.00e+00 / 0.00e+00
Duality gap (abs/rel):           3.15e\PYGZhy{}03 / 8.43e\PYGZhy{}07

Experimental: using crossover to find a basic solution after PDLP
Please set parameter Crossover to 0 if the basic solution is not required
Please set parameter PDLPTol to a smaller value if the crossover cleanup takes too long

Starting crossover using up to 8 threads

50320 primal pushes remaining     92.09s
12495 primal pushes remaining     97.71s
 4124 primal pushes remaining       102s
  202 primal pushes remaining       104s
    0 primal pushes remaining       104s
1480858 dual pushes remaining       104s
 589582 dual pushes remaining       106s
      0 dual pushes remaining       107s

Method   Iteration           Objective  Primal.NInf   Dual.NInf        Time
Dual        986011    1.8710000000e+03            0           0     107.97s
Postsolving
Dual        986011    1.8710000000e+03            0           0     110.56s
Unfolding
Dual       1036742    1.8710000000e+03            0           0     157.25s

Solving finished
Status: Optimal  Objective: 1.8710000000e+03  Iterations: 1036742  Time: 157.69s
\end{sphinxVerbatim}

\sphinxAtStartPar
\sphinxstylestrong{Note:} After solving with a First\sphinxhyphen{}order method (PDLP), if the solution
reaches optimality (\sphinxcode{\sphinxupquote{Status: Optimal}}), the default behavior is to perform
the Crossover process to the basic solution. This process can also be disabled
by setting the parameter {\hyperref[\detokenize{parameter:crossover}]{\sphinxcrossref{\DUrole{std,std-ref}{Crossover}}}} to 0.

\sphinxstepscope

\chapter{File formats}
\label{\detokenize{fileformats:file-formats}}\label{\detokenize{fileformats:chapfileformat}}\label{\detokenize{fileformats::doc}}

\section{File format list}
\label{\detokenize{fileformats:file-format-list}}
\sphinxAtStartPar
The file formats currently supported by COPT are listed below:

\begin{savenotes}\sphinxattablestart
\sphinxthistablewithglobalstyle
\centering
\sphinxcapstartof{table}
\sphinxthecaptionisattop
\sphinxcaption{Supported file formats}\label{\detokenize{fileformats:copttab-file}}
\sphinxaftertopcaption
\begin{tabular}[t]{|\X{15}{35}|\X{20}{35}|}
\sphinxtoprule
\sphinxstyletheadfamily 
\sphinxAtStartPar
Fileformat
&\sphinxstyletheadfamily 
\sphinxAtStartPar
File extension
\\
\sphinxmidrule
\sphinxtableatstartofbodyhook
\sphinxAtStartPar
MPS model file
&
\sphinxAtStartPar
\sphinxcode{\sphinxupquote{.mps}}, \sphinxcode{\sphinxupquote{.mps.gz}}
\\
\sphinxhline
\sphinxAtStartPar
LP model file
&
\sphinxAtStartPar
\sphinxcode{\sphinxupquote{.lp}}, \sphinxcode{\sphinxupquote{.lp.gz}}
\\
\sphinxhline
\sphinxAtStartPar
SDPA model file
&
\sphinxAtStartPar
\sphinxcode{\sphinxupquote{.dat\sphinxhyphen{}s}}, \sphinxcode{\sphinxupquote{.dat\sphinxhyphen{}s.gz}}
\\
\sphinxhline
\sphinxAtStartPar
CBF model file
&
\sphinxAtStartPar
\sphinxcode{\sphinxupquote{.cbf}}, \sphinxcode{\sphinxupquote{.cbf.gz}}
\\
\sphinxhline
\sphinxAtStartPar
model file in COPT binary format
&
\sphinxAtStartPar
\sphinxcode{\sphinxupquote{.bin}}, \sphinxcode{\sphinxupquote{.bin.gz}}
\\
\sphinxhline
\sphinxAtStartPar
basis file
&
\sphinxAtStartPar
\sphinxcode{\sphinxupquote{.bas}}
\\
\sphinxhline
\sphinxAtStartPar
solution file
&
\sphinxAtStartPar
\sphinxcode{\sphinxupquote{.sol}}
\\
\sphinxhline
\sphinxAtStartPar
IIS file
&
\sphinxAtStartPar
\sphinxcode{\sphinxupquote{.iis}}
\\
\sphinxhline
\sphinxAtStartPar
FeasRelax file
&
\sphinxAtStartPar
\sphinxcode{\sphinxupquote{.relax}}
\\
\sphinxhline
\sphinxAtStartPar
MIP initial solution file
&
\sphinxAtStartPar
\sphinxcode{\sphinxupquote{.mst}}
\\
\sphinxhline
\sphinxAtStartPar
parameter file
&
\sphinxAtStartPar
\sphinxcode{\sphinxupquote{.par}}
\\
\sphinxhline
\sphinxAtStartPar
parameter tuning file
&
\sphinxAtStartPar
\sphinxcode{\sphinxupquote{.tune}}
\\
\sphinxhline
\sphinxAtStartPar
branching\sphinxhyphen{}order file
&
\sphinxAtStartPar
\sphinxcode{\sphinxupquote{.ord}}
\\
\sphinxbottomrule
\end{tabular}
\sphinxtableafterendhook\par
\sphinxattableend\end{savenotes}

\section{File I/O operations}
\label{\detokenize{fileformats:file-i-o-operations}}
\sphinxAtStartPar
By calling the relevant functions, users can input external model files to COPT for reading. At the same time, they can also save the output of the modeling and optimization results of COPT and output the files.

\sphinxAtStartPar
Let’s take reading/writing the model in MPS format in the current directory as an example (similar operations are performed for other files). The implementation methods in different interfaces are as follows:

\begin{savenotes}\sphinxattablestart
\sphinxthistablewithglobalstyle
\centering
\sphinxcapstartof{table}
\sphinxthecaptionisattop
\sphinxcaption{Functions for input and output files}\label{\detokenize{fileformats:copttab-fileoi}}
\sphinxaftertopcaption
\begin{tabular}[t]{|\X{10}{60}|\X{25}{60}|\X{25}{60}|}
\sphinxtoprule
\sphinxstyletheadfamily 
\sphinxAtStartPar
API
&\sphinxstyletheadfamily 
\sphinxAtStartPar
Input
&\sphinxstyletheadfamily 
\sphinxAtStartPar
Output
\\
\sphinxmidrule
\sphinxtableatstartofbodyhook
\sphinxAtStartPar
COPT cmd
&
\sphinxAtStartPar
\sphinxcode{\sphinxupquote{read example.mps}}
&
\sphinxAtStartPar
\sphinxcode{\sphinxupquote{write example.mps}}
\\
\sphinxhline
\sphinxAtStartPar
C
&
\sphinxAtStartPar
\sphinxcode{\sphinxupquote{COPT\_ReadMps}}
&
\sphinxAtStartPar
\sphinxcode{\sphinxupquote{COPT\_WriteMps}}
\\
\sphinxhline
\sphinxAtStartPar
Python
&
\sphinxAtStartPar
\sphinxcode{\sphinxupquote{Model.read()}} / \sphinxcode{\sphinxupquote{Model.readMps()}}
&
\sphinxAtStartPar
\sphinxcode{\sphinxupquote{Model.write()}} / \sphinxcode{\sphinxupquote{Model.writeMps()}}
\\
\sphinxhline
\sphinxAtStartPar
C++
&
\sphinxAtStartPar
\sphinxcode{\sphinxupquote{Model::read()}} / \sphinxcode{\sphinxupquote{Model::readMps()}}
&
\sphinxAtStartPar
\sphinxcode{\sphinxupquote{Model::write()}} / \sphinxcode{\sphinxupquote{Model::writeMps()}}
\\
\sphinxhline
\sphinxAtStartPar
C\#
&
\sphinxAtStartPar
\sphinxcode{\sphinxupquote{Model.Read()}} / \sphinxcode{\sphinxupquote{Model.ReadMps()}}
&
\sphinxAtStartPar
\sphinxcode{\sphinxupquote{Model.Write()}} / \sphinxcode{\sphinxupquote{Model.WriteMps()}}
\\
\sphinxhline
\sphinxAtStartPar
Java
&
\sphinxAtStartPar
\sphinxcode{\sphinxupquote{Model.read()}} / \sphinxcode{\sphinxupquote{Model.readMps()}}
&
\sphinxAtStartPar
\sphinxcode{\sphinxupquote{Model.write()}} / \sphinxcode{\sphinxupquote{Model.writeMps()}}
\\
\sphinxbottomrule
\end{tabular}
\sphinxtableafterendhook\par
\sphinxattableend\end{savenotes}

\section{Model file introduction}
\label{\detokenize{fileformats:model-file-introduction}}
\sphinxAtStartPar
Users can find the model file examples that come with COPT in the \sphinxcode{\sphinxupquote{"examples/data"}} directory of the installation package. Here we introduce the specific contents of two common model file formats: MPS and LP.
\begin{quote}

\begin{figure}[H]
\centering

\noindent\sphinxincludegraphics[scale=0.5]{{copt-file_example}.png}
\end{figure}
\end{quote}

\sphinxAtStartPar
\sphinxstylestrong{MPS format}

\sphinxAtStartPar
MPS is a universal model file standard format. Different types of optimization problems can be output and stored in mps format, which is widely used in optimization softwares.

\sphinxAtStartPar
The following is an example of a model file in MPS format:

\begin{sphinxVerbatim}[commandchars=\\\{\}]
NAME          COPTPROB
OBJSENSE
    MAX
ROWS
N  \PYGZus{}\PYGZus{}OBJ\PYGZus{}\PYGZus{}\PYGZus{}
L  R0000000
G  R0000001
COLUMNS
    x         \PYGZus{}\PYGZus{}OBJ\PYGZus{}\PYGZus{}\PYGZus{}  1.2
    x         R0000000  1.5
    x         R0000001  0.80000000000000004
    y         \PYGZus{}\PYGZus{}OBJ\PYGZus{}\PYGZus{}\PYGZus{}  1.8
    y         R0000000  1.2
    y         R0000001  0.59999999999999998
    z         \PYGZus{}\PYGZus{}OBJ\PYGZus{}\PYGZus{}\PYGZus{}  2.1000000000000001
    z         R0000000  1.8
    z         R0000001  0.90000000000000002
RHS
    RHS       R0000000  2.6000000000000001
    RHS       R0000001  1.2
BOUNDS
LO BOUND     x         0.10000000000000001
UP BOUND     x         0.59999999999999998
LO BOUND     y         0.20000000000000001
UP BOUND     y         1.5
LO BOUND     z         0.29999999999999999
UP BOUND     z         2.7999999999999998
ENDATA
\end{sphinxVerbatim}

\sphinxAtStartPar
This MPS format example mainly includes several parts: NAME, OBJSENSE, ROWS, COLUMNS, RHS, and BOUNDS.
\begin{enumerate}
\sphinxsetlistlabels{\arabic}{enumi}{enumii}{}{.}%
\item {} 
\sphinxAtStartPar
NAME: The name of the model

\item {} 
\sphinxAtStartPar
OBJSENSE: Optimization direction of the objective function

\item {} 
\sphinxAtStartPar
ROWS: Constraints and their directions in the model (L means \textless{}= constraints, G means \textgreater{}= constraints, N means no boundaries)

\item {} 
\sphinxAtStartPar
COLUMNS: Variables and their coefficients in the model

\item {} 
\sphinxAtStartPar
RHS: The value of the right\sphinxhyphen{}hand term of the constraint

\item {} 
\sphinxAtStartPar
BOUNDS: Bounds of variables (LO means lower bound, UP means upper bound, FR means no bounds)

\end{enumerate}

\begin{sphinxadmonition}{note}{Notes}
\begin{enumerate}
\sphinxsetlistlabels{\arabic}{enumi}{enumii}{}{.}%
\item {} 
\sphinxAtStartPar
In the ROWS section, the first line \sphinxcode{\sphinxupquote{\_\_OBJ\_\_\_}} represents the objective function.

\item {} 
\sphinxAtStartPar
In the COLS part, the form \sphinxcode{\sphinxupquote{x \_\_OBJ\_\_\_ 1.2}} indicates that the coefficient of variable x in the objective function is 1.2.

\item {} 
\sphinxAtStartPar
In MPS format, integer variables will be identified by the following fields:
\begin{itemize}
\item {} 
\sphinxAtStartPar
First integer variable: \sphinxcode{\sphinxupquote{MARKER \textquotesingle{}MARKER\textquotesingle{} \textquotesingle{}INTORG\textquotesingle{}}}

\item {} 
\sphinxAtStartPar
Last integer variable: \sphinxcode{\sphinxupquote{MARKER \textquotesingle{}MARKER\textquotesingle{} \textquotesingle{}INTEND\textquotesingle{}}}

\end{itemize}

\end{enumerate}
\end{sphinxadmonition}

\sphinxAtStartPar
\sphinxstylestrong{LP format}

\sphinxAtStartPar
The LP format is closer to the algebraic form. It is more readable than MPS, and can easily correspond to its original mathematical model.

\sphinxAtStartPar
The following is an example of a model file in LP format:

\begin{sphinxVerbatim}[commandchars=\\\{\}]
\PYGZbs{}Generated by Cardinal Operations

Maximize
1.2 x + 1.8 y + 2.1 z
Subject To
1.5 x + 1.2 y + 1.8 z \PYGZlt{}= 2.6
0.8 x + 0.6 y + 0.9 z \PYGZgt{}= 1.2
Bounds
0.1 \PYGZlt{}= x \PYGZlt{}= 0.6
0.2 \PYGZlt{}= y \PYGZlt{}= 1.5
0.3 \PYGZlt{}= z \PYGZlt{}= 2.8
END
\end{sphinxVerbatim}

\sphinxAtStartPar
This LP format mainly includes several parts: objective function (Maximize), constraints (Subject To), and variable scope (Bounds).

\begin{sphinxadmonition}{note}{Notes}
\begin{enumerate}
\sphinxsetlistlabels{\arabic}{enumi}{enumii}{}{.}%
\item {} 
\sphinxAtStartPar
In some LP format file, we can see variable names in the form \sphinxcode{\sphinxupquote{x\#1}}, which marks x as the first variable. When the user does not specify a variable name, this is the name automatically generated when outputting the lp model file.

\item {} 
\sphinxAtStartPar
If there is a binary type in the variable, it will be identified by the \sphinxcode{\sphinxupquote{Binaries}} field.

\end{enumerate}
\end{sphinxadmonition}

\sphinxstepscope

\chapter{FAQs}
\label{\detokenize{faq:faqs}}\label{\detokenize{faq:chapfaq}}\label{\detokenize{faq::doc}}

\section{Installation and Licensing Configuration Related}
\label{\detokenize{faq:installation-and-licensing-configuration-related}}\begin{itemize}
\item {} 
\sphinxAtStartPar
\sphinxstylestrong{Q:} What is the reason for the error \sphinxcode{\sphinxupquote{invalid username}} when configuring the license?

\sphinxAtStartPar
\sphinxstylestrong{A:} This error indicates that the username was incorrectly filled in when applying, and you
can re\sphinxhyphen{}fill the form with the correct username to apply. For information on how to obtain
username under different operating systems, please refer to the \sphinxhref{https://www.shanshu.ai/copt}{COPT application page}, please remark in the
application reason with \sphinxstylestrong{“The username is incorrectly filled in, reapply”} . we will issue new
license for the correct username.

\end{itemize}

\sphinxAtStartPar

\begin{itemize}
\item {} 
\sphinxAtStartPar
\sphinxstylestrong{Q:} After downloading COPT, an antivirus software installed on the computer reports a virus
and automatically isolates it.

\sphinxAtStartPar
\sphinxstylestrong{A:} The COPT software downloaded from the COPT official download link is the official version,
which has not been developed with any suspicious virus behavior. It can be determined that the
anti\sphinxhyphen{}virus software is falsely reporting. Please temporarily close the anti\sphinxhyphen{}virus software before
downloading it.

\end{itemize}

\sphinxAtStartPar

\begin{itemize}
\item {} 
\sphinxAtStartPar
\sphinxstylestrong{Q:} After validating the license (executing \sphinxcode{\sphinxupquote{copt\_licgen \sphinxhyphen{}v}} ), it reports an error:
\sphinxcode{\sphinxupquote{Missing Files}} or \sphinxcode{\sphinxupquote{Invalid Signature}} .

\sphinxAtStartPar
\sphinxstylestrong{A:} This type of error indicates that the license file configuration fails. Please refer to
{\hyperref[\detokenize{install:parcoptgetlic}]{\sphinxcrossref{\DUrole{std,std-ref}{Installation Instructions: Configuring License File}}}} to check whether the
steps for configuring the license file are correctly followed. Common reasons are as follows:
\begin{enumerate}
\sphinxsetlistlabels{\arabic}{enumi}{enumii}{}{.}%
\item {} 
\sphinxAtStartPar
The license file in the current working directory is not compatible with the version of COPT
(for example: the license is version 4.0, while the COPT is version 5.0), please check the
\sphinxcode{\sphinxupquote{VERSION}} in \sphinxcode{\sphinxupquote{"license.dat"}} to confirm whether the major version is the same, if not,
please go to the \sphinxhref{https://www.shanshu.ai/copt}{COPT application page} to re\sphinxhyphen{}apply, and we will issue you the latest license.

\item {} 
\sphinxAtStartPar
In Windows system, if the COPT software is installed in the system disk (usually \sphinxcode{\sphinxupquote{C}}) in a
non\sphinxhyphen{}user directory (eg: the default installation path \sphinxcode{\sphinxupquote{"C:\textbackslash{}Program Files\textbackslash{}copt70"}}) , you
need to \sphinxstylestrong{open the command line window with administrator privileges} and execute the license
acquisition command \sphinxcode{\sphinxupquote{copt\_licgen}} again.

\end{enumerate}

\end{itemize}

\sphinxAtStartPar

\begin{itemize}
\item {} 
\sphinxAtStartPar
\sphinxstylestrong{Q:} I have already installed an old version of the COPT Python interface ( \sphinxcode{\sphinxupquote{coptpy}} ), how
do I upgrade to the new version?

\sphinxAtStartPar
\sphinxstylestrong{A:} Please refer to {\hyperref[\detokenize{pythoninterface:chappythoninterface-practiceupgrade}]{\sphinxcrossref{\DUrole{std,std-ref}{Python Interface Quick Start: Upgrade to the newer version}}}} for detailed steps.

\end{itemize}

\subsection{MacOS System}
\label{\detokenize{faq:macos-system}}\begin{itemize}
\item {} 
\sphinxAtStartPar
\sphinxstylestrong{Q:} When calling \sphinxcode{\sphinxupquote{coptpy}} on MacOS, an error is reported: \sphinxcode{\sphinxupquote{ImportError: from .coptpywrap
import * symbol not found in flat namespace}}.

\sphinxAtStartPar
\sphinxstylestrong{A:} This error may occur before COPT 6.5.12, since the architecture of \sphinxcode{\sphinxupquote{coptpy}} and Anaconda
do not match. For example, \sphinxcode{\sphinxupquote{coptpy}} is the M chip version (arm64), and Anaconda is x86 version,
you could install Anaconda that supports arm64 architecture to solve this problem.

\sphinxAtStartPar
However, from COPT 6.5.12, for MacOS systems, we provide universal package, so this problem could
be resolved by upgrading COPT to the latest version.

\end{itemize}

\sphinxAtStartPar

\begin{itemize}
\item {} 
\sphinxAtStartPar
\sphinxstylestrong{Q:} When the license is configured in the MacOS system, the \sphinxcode{\sphinxupquote{copt\_licgen}} command is
executed in the terminal, and an error is reported: \sphinxcode{\sphinxupquote{command not found: copt\_licgen}} .

\sphinxAtStartPar
\sphinxstylestrong{A:} This error is because the relevant environment variables of COPT are not configured. In
the MacOS system, you need to configure the environment variables after installing COPT. Please
refer to {\hyperref[\detokenize{install:chapinstall-macosenv}]{\sphinxcrossref{\DUrole{std,std-ref}{Installation Instructions: MacOS System}}}} chapter to obtain
detailed installation instructions.

\end{itemize}

\sphinxAtStartPar

\begin{itemize}
\item {} 
\sphinxAtStartPar
\sphinxstylestrong{Q:} When manually configuring environment variables, I copy content directly from the document
to \sphinxcode{\sphinxupquote{.zshrc}} file or \sphinxcode{\sphinxupquote{.bash\_profile}} file, causing the configuration to fail.

\sphinxAtStartPar
\sphinxstylestrong{A:} Due to the document encoding problem, the above environment variables cannot be directly
copied to the corresponding file, and the contents of the environment variables need to be
manually entered.

\end{itemize}

\subsection{Windows System}
\label{\detokenize{faq:windows-system}}\begin{itemize}
\item {} 
\sphinxAtStartPar
\sphinxstylestrong{Q:} In Windows system, when executing \sphinxcode{\sphinxupquote{copt\_licgen}} to generate the license file, an error
is reported that the license file cannot be written to the hard disk, and the error message is:
\sphinxcode{\sphinxupquote{error opening file}} .

\sphinxAtStartPar
\sphinxstylestrong{A:} If the COPT software is installed in the system disk (usually \sphinxcode{\sphinxupquote{C}}) in a non\sphinxhyphen{}user
directory (eg: the default installation path \sphinxcode{\sphinxupquote{"C:\textbackslash{}Program Files\textbackslash{}copt80"}}), y you need to \sphinxstylestrong{open
the command line window with administrator privileges} and execute the license acquisition
command \sphinxcode{\sphinxupquote{copt\_licgen}}, in order to normally write the license file to the \sphinxcode{\sphinxupquote{C}} drive.
Administrator privileges are not required to execute permission acquisition commands under user
directories such as \sphinxcode{\sphinxupquote{"C:\textbackslash{}Users\textbackslash{}copt80"}}.

\end{itemize}

\sphinxAtStartPar

\begin{itemize}
\item {} 
\sphinxAtStartPar
\sphinxstylestrong{Q:} In Windows system, when installing COPT Python interface via command \sphinxcode{\sphinxupquote{pip install
coptpy}}, an error is displayed: \sphinxcode{\sphinxupquote{could not find a version, no matching distribution}} , what is
the reason?

\sphinxAtStartPar
\sphinxstylestrong{A:} Please do not use Python installed through Microsoft Store, it is recommended to download
from \sphinxhref{https://www.anaconda.com/distribution/}{Anaconda distribution} or \sphinxhref{https://www.python.org/}{Python official
distribution} Download Python.

\end{itemize}

\sphinxAtStartPar

\begin{itemize}
\item {} 
\sphinxAtStartPar
\sphinxstylestrong{Q:} In Windows system, when installing COPT Python interface through COPT installation package
( \sphinxcode{\sphinxupquote{python setup.py install}} ), an error \sphinxcode{\sphinxupquote{could not create build: access denied}} is reported.

\sphinxAtStartPar
\sphinxstylestrong{A:} If COPT is installed in the system disk (usually \sphinxcode{\sphinxupquote{C}}) in a non\sphinxhyphen{}user directory (eg: the
default installation path \sphinxcode{\sphinxupquote{"C:\textbackslash{}Program Files\textbackslash{}copt70"}}), you need to first \sphinxstylestrong{Open the command
line window with administrator privileges} and execute the command \sphinxcode{\sphinxupquote{python setup.py install}}.

\end{itemize}

\section{Modeling and Solving Functions Related}
\label{\detokenize{faq:modeling-and-solving-functions-related}}\begin{itemize}
\item {} 
\sphinxAtStartPar
\sphinxstylestrong{Q:} When creating a COPT solution environment, two lines of version information will be
output. If I want to turn off this information, how should I do it?

\item {} 
\sphinxAtStartPar
\sphinxstylestrong{A:} You can turn it off by setting \sphinxcode{\sphinxupquote{"nobanner"}} to \sphinxcode{\sphinxupquote{"1"}} in \sphinxcode{\sphinxupquote{EnvrConfig}} before creating
the solution environment. Taking the Python API as an example, the specific operations are as
follows:

\begin{sphinxVerbatim}[commandchars=\\\{\}]
\PYG{n}{envconfig} \PYG{o}{=} \PYG{n}{coptpy}\PYG{o}{.}\PYG{n}{EnvrConfig}\PYG{p}{(}\PYG{p}{)}
\PYG{n}{envconfig}\PYG{o}{.}\PYG{n}{set}\PYG{p}{(}\PYG{l+s+s2}{\PYGZdq{}}\PYG{l+s+s2}{nobanner}\PYG{l+s+s2}{\PYGZdq{}}\PYG{p}{,} \PYG{l+s+s2}{\PYGZdq{}}\PYG{l+s+s2}{1}\PYG{l+s+s2}{\PYGZdq{}}\PYG{p}{)}
\PYG{n}{env} \PYG{o}{=} \PYG{n}{coptpy}\PYG{o}{.}\PYG{n}{Envr}\PYG{p}{(}\PYG{n}{envconfig}\PYG{p}{)}
\PYG{n}{model} \PYG{o}{=} \PYG{n}{env}\PYG{o}{.}\PYG{n}{createModel}\PYG{p}{(}\PYG{p}{)}
\end{sphinxVerbatim}

\item {} 
\sphinxAtStartPar
\sphinxstylestrong{Q:} How to deal with the situation where the model is infeasible?

\sphinxAtStartPar
\sphinxstylestrong{A:} COPT provides functions to calculate IIS and feasible relaxation to analyze the reasons
for model infeasibility: Computing IIS will obtain the minimum set of infeasible constraints and
variables, and feasibility relaxation attempts to make the model feasible with minimal changes.

\item {} 
\sphinxAtStartPar
\sphinxstylestrong{Q:} Are there default value ranges for variables created in COPT?

\sphinxAtStartPar
\sphinxstylestrong{A:} Yes, the default lower bound of variables in COPT is 0, and the upper bound is
\sphinxcode{\sphinxupquote{COPT.INFINITY}}. Users can specify the lower bound of variables through the function argument
\sphinxtitleref{lb}, and the upper bound of variables through \sphinxtitleref{ub}.

\item {} 
\sphinxAtStartPar
\sphinxstylestrong{Q:} What is the reason for the error \sphinxcode{\sphinxupquote{ValueError: cannot create object arrays from
iterator.}} when adding matrix variables using the matrix modeling method of COPT?

\sphinxAtStartPar
\sphinxstylestrong{A:} The matrix modeling function supported by COPT Python has a minimum version requirement,
the minimum version requirement for \sphinxcode{\sphinxupquote{NumPy}} is 1.23, and the minimum version requirement for
\sphinxcode{\sphinxupquote{Python}} is 3.8 . \sphinxcode{\sphinxupquote{NumPy}} can be upgraded to the latest version by \sphinxcode{\sphinxupquote{pip install \sphinxhyphen{}\sphinxhyphen{}upgrade
numpy}}.

\item {} 
\sphinxAtStartPar
\sphinxstylestrong{Q:} When adding constraints using Python interface, if the modeling efficiency is slow, are
there any ways to improve the modeling process?

\sphinxAtStartPar
\sphinxstylestrong{A:} The Python interface of COPT supports building linear expression, quadratic expression and
PSD expression in natural way. For linear and quadratic expression, it is recommended to use
quicksum() to build expression objects. For linear and PSD expression, it is recommended to use
psdquicksum() to build expression objects. Both of them implement inplace summation, which is
much faster than standard plus operator.

\end{itemize}

\section{GPU Usage Related}
\label{\detokenize{faq:gpu-usage-related}}\label{\detokenize{faq:chapfaq-gpu}}\begin{itemize}
\item {} 
\sphinxAtStartPar
\sphinxstylestrong{Q:} Are there any requirements for the CUDA library version when enabling COPT’s GPU solver?

\sphinxAtStartPar
\sphinxstylestrong{A:} COPT requires a minimum version of 12.0.1 for the CUDA library.

\item {} 
\sphinxAtStartPar
\sphinxstylestrong{Q:} Are there any requirements for the GPU architecture when enabling GPU solver?

\sphinxAtStartPar
\sphinxstylestrong{A:} The GPU architecture must be at least Maxwell or a more recent version. (Maxwell is a GPU
architecture introduced by NVIDIA in 2014 as an upgrade to the earlier Kepler architecture.)
Besides, COPT supports the Blackwell architecture for Windows and Linux platform.

\item {} 
\sphinxAtStartPar
\sphinxstylestrong{Q:} What are the common error messages and possible reasons when the machine cannot use COPT
GPU solver?

\sphinxAtStartPar
\sphinxstylestrong{A:} Common error messages and possible reasons are as follows:

\sphinxAtStartPar
1.Solving log indicates \sphinxcode{\sphinxupquote{"CUDA libraries could not be loaded: cuBLAS cuSPARSE"}} which suggests
that the listed necessary CUDA libraries are missing. You could try checking and setting the
environment variable \sphinxcode{\sphinxupquote{"LD\_LIBRARY\_PATH"}} to point to the directory where CUDA is installed.
(Please follow the instructions provided after CUDA installation. Environment variables are
automatically configured during installation on Windows systems. For Linux systems, manual
configuration of environment variables is typically required. The directory will be like:
\sphinxcode{\sphinxupquote{"/usr/local/cuda/lib64"}} ).

\sphinxAtStartPar
2.Solving error message \sphinxcode{\sphinxupquote{"Fail to solve problem"}} is usually due to the lower version of CUDA
Driver. Please upgrade the COPT Driver (for Linux systems: 525.60.13 or above; for Windows
systems: 527.41 or above) to resolve this issue.

\sphinxAtStartPar
3.Solving error message \sphinxcode{\sphinxupquote{"sparse matrix format CUSPARSE\_FORMAT\_CSC is not supported"}} is
usually caused by the lower version of CUDA Toolkit (typically occurring between CUDA V11.2 and
V11.6). Please upgrade CUDA to version 12.0.1 or above to resolve this issue.

\item {} 
\sphinxAtStartPar
\sphinxstylestrong{Q:} On a client machine with multiple GPUs, when setting the parameter \sphinxcode{\sphinxupquote{GPUDevice}} to use a
specific GPU number, why does it still detect only GPU with number 0 during solving?

\sphinxAtStartPar
\sphinxstylestrong{A:} Please check if the environment variable \sphinxcode{\sphinxupquote{"CUDA\_VISIBLE\_DEVICES"}} has been manually set
to specify the visible GPU devices for CUDA. Try not to set this environment variable so that
COPT can detect all available GPUs on the current machine.

\item {} 
\sphinxAtStartPar
\sphinxstylestrong{Q:} Why do I encounter errors when using COPT’s GPU solver via Windows Subsystem for Linux
(WSL) despite installing CUDA libraries that meet the version requirements (V12.0.1 or higher)?

\sphinxAtStartPar
\sphinxstylestrong{A:} Please check if the CUDA Driver version meets the requirements. WSL usually skips the
Driver installation when installing CUDA and directly uses the Driver already installed in
Windows. Please manually upgrade the CUDA Driver version and then restart WSL to resolve this
issue.

\item {} 
\sphinxAtStartPar
\sphinxstylestrong{Q:} When solving with the GPU Barrier method, an error message \sphinxcode{\sphinxupquote{"GPU memory issue"}} is
reported. How should I address this?

\sphinxAtStartPar
\sphinxstylestrong{A:} This error is typically caused by insufficient GPU memory. Please try setting the
parameter \sphinxcode{\sphinxupquote{GPUMode=1}}. If the issue persists, consider solving the model using the CPU
(\sphinxcode{\sphinxupquote{GPUMode=0}}).

\item {} 
\sphinxAtStartPar
\sphinxstylestrong{Q:} When solving with the GPU Barrier method, a message says \sphinxcode{\sphinxupquote{"Performance may degenerate for
GPUMode=2 on this problem, consider trying GPUMode=1 instead"}}. What does this mean?

\sphinxAtStartPar
\sphinxstylestrong{A:} This message indicates that the problem is still solvable, but using \sphinxcode{\sphinxupquote{GPUMode=2}}
(high\sphinxhyphen{}performance mode) may result in longer solving time. You are advised to set \sphinxcode{\sphinxupquote{GPUMode=1}}
to solve the problem in the standard mode.

\end{itemize}

\section{CPU and Memory Binding Related}
\label{\detokenize{faq:cpu-and-memory-binding-related}}\label{\detokenize{faq:chapfaq-cpumem}}\begin{itemize}
\item {} 
\sphinxAtStartPar
\sphinxstylestrong{Q:} What is a NUMA node? How can I obtain the number of NUMA nodes in the operating system?

\sphinxAtStartPar
\sphinxstylestrong{A:} NUMA (Non\sphinxhyphen{}Uniform Memory Access) is a memory architecture designed for multi\sphinxhyphen{}processor
systems. A NUMA node is the basic unit in this architecture, consisting of a group of CPU cores
and the local memory directly attached to them. If the current machine has \sphinxstyleemphasis{N} physical CPUs,
the return value of \sphinxcode{\sphinxupquote{env.getNumaNodeCount()}} will be \sphinxstyleemphasis{N}.

\item {} 
\sphinxAtStartPar
\sphinxstylestrong{Q:} How can I solve a model using a specified NUMA node (CPU)?

\sphinxAtStartPar
\sphinxstylestrong{A:} If \sphinxcode{\sphinxupquote{env.getNumaNodeCnt()}} returns a value greater than 1, the machine contains multiple
NUMA nodes. Since cross\sphinxhyphen{}node memory access is typically more expensive than access within the
same node, you may call \sphinxcode{\sphinxupquote{env.bindNumaCpu(numaNode=n)}} to bind the solver to NUMA node \sphinxstyleemphasis{n} in
order to reduce cross\sphinxhyphen{}node access overhead. The \sphinxcode{\sphinxupquote{numaNode}} argument starts from 0.

\item {} 
\sphinxAtStartPar
\sphinxstylestrong{Q:} How can I bind memory to a specified NUMA node?

\sphinxAtStartPar
\sphinxstylestrong{A:} Similar to CPU binding, you may call \sphinxcode{\sphinxupquote{env.bindNumaMem(numaNode=0)}} to restrict memory
allocation to the local memory of the specified NUMA node, thereby avoiding cross\sphinxhyphen{}node memory
access.

\item {} 
\sphinxAtStartPar
\sphinxstylestrong{Q:} How does \sphinxcode{\sphinxupquote{Envr.setCpuAffinity()}} specify the CPU cores used for solving? How should the
argument \sphinxcode{\sphinxupquote{hexMask}} be interpreted?

\sphinxAtStartPar
\sphinxstylestrong{A:} \sphinxcode{\sphinxupquote{hexMask}} is a CPU mask represented as a hexadecimal string and used to specify the set
of CPU cores available to the current process. After converting the mask to binary, bits are
interpreted from right to left, corresponding to CPU core indices 0, 1, 2, and so on. A bit set
to 1 indicates that the corresponding core is enabled for solving. For example, on a machine with
32 CPU cores, the call \sphinxcode{\sphinxupquote{env.setCpuAffinity("F1")}} uses the mask \sphinxcode{\sphinxupquote{"F1"}}, which corresponds to
the binary \sphinxcode{\sphinxupquote{1111 0001}}. Thus, CPU cores 0 and 4\textendash{}7 will participate in solving.

\item {} 
\sphinxAtStartPar
\sphinxstylestrong{Q:} How can I check which CPU cores are used by the current solving process?

\sphinxAtStartPar
\sphinxstylestrong{A:} You may call \sphinxcode{\sphinxupquote{env.getCpuAffinity()}} to retrieve the CPU affinity setting of the current
solving process. It returns a list of core indices, such as \sphinxcode{\sphinxupquote{{[}0, 4, 5, 6, 7{]}}}. If no CPU
affinity is set, an empty list is returned, meaning the operating system determines the CPU
scheduling.

\item {} 
\sphinxAtStartPar
\sphinxstylestrong{Q:} What is the effective scope of CPU affinity set via \sphinxcode{\sphinxupquote{Envr.setCpuAffinity()}}?

\sphinxAtStartPar
\sphinxstylestrong{A:} This function sets the CPU affinity for the current solving process. Once set, the
affinity remains effective for the entire lifetime of the process and will not revert when
\sphinxcode{\sphinxupquote{env.close()}} is called.

\end{itemize}

\sphinxstepscope

\chapter{C API Reference}
\label{\detokenize{capiref:c-api-reference}}\label{\detokenize{capiref:chapapi}}\label{\detokenize{capiref::doc}}
\sphinxAtStartPar
The \sphinxstylestrong{Cardinal Optimizer} provides a C API library for advanced usage.
This section documents all the COPT constants,
API functions, parameters and attributes listed in \sphinxcode{\sphinxupquote{copt.h}}.

\section{Constants}
\label{\detokenize{capiref:constants}}\label{\detokenize{capiref:chapapi-const}}
\sphinxAtStartPar
There are three types of constants.
\begin{enumerate}
\sphinxsetlistlabels{\arabic}{enumi}{enumii}{}{.}%
\item {} 
\sphinxAtStartPar
Constructing models, such as optimization directions, constraint senses or variable types.

\item {} 
\sphinxAtStartPar
Accessing solution results, such as API return code, basis status and LP status.

\item {} 
\sphinxAtStartPar
Monitoring optimization progress, such as callback context.

\end{enumerate}

\subsection{Optimization directions}
\label{\detokenize{capiref:optimization-directions}}
\sphinxAtStartPar
For different optimization scenarios, it may be required
to either maximize or minimize the objective function.
There are two optimization directions:
\begin{itemize}
\item {} 
\sphinxAtStartPar
\sphinxcode{\sphinxupquote{COPT\_MINIMIZE}}
\begin{quote}

\sphinxAtStartPar
For minimizing the objective function.
\end{quote}

\item {} 
\sphinxAtStartPar
\sphinxcode{\sphinxupquote{COPT\_MAXIMIZE}}
\begin{quote}

\sphinxAtStartPar
For maximizing the objective function.
\end{quote}

\end{itemize}

\sphinxAtStartPar
The optimization direction is automatically
set when reading in a problem from file.
It can also be set explicitly using \sphinxcode{\sphinxupquote{COPT\_SetObjSense}}.

\subsection{Infinity}
\label{\detokenize{capiref:infinity}}
\sphinxAtStartPar
In COPT, the infinite bound is represented by a large value,
which can be set using the double parameter
\sphinxcode{\sphinxupquote{COPT\_DBLPARAM\_INFBOUND}},
whose default value is also available as a constant:
\begin{itemize}
\item {} 
\sphinxAtStartPar
\sphinxcode{\sphinxupquote{COPT\_INFINITY}}
\begin{quote}

\sphinxAtStartPar
The default value (\sphinxcode{\sphinxupquote{1e30}}) of the infinite bound.
\end{quote}

\end{itemize}

\subsection{Undefined Value}
\label{\detokenize{capiref:undefined-value}}
\sphinxAtStartPar
In COPT, the undefined value is represented by another large value.
For example, the default solution value of MIP start is set to
a constant:
\begin{itemize}
\item {} 
\sphinxAtStartPar
\sphinxcode{\sphinxupquote{COPT\_UNDEFINED}}
\begin{quote}

\sphinxAtStartPar
Undefined value(\sphinxcode{\sphinxupquote{1e40\textasciigrave{}}}).
\end{quote}

\end{itemize}

\subsection{Constraint senses}
\label{\detokenize{capiref:constraint-senses}}
\sphinxAtStartPar
\sphinxstylestrong{NOTE:} Using constraint senses is supported by COPT but not recommended.
We recommend defining constraints using explicit lower and upper bounds.

\sphinxAtStartPar
Traditionally, for optimization models, constraints are defined using \sphinxstylestrong{senses}.
The most common constraint senses are:
\begin{itemize}
\item {} 
\sphinxAtStartPar
\sphinxcode{\sphinxupquote{COPT\_LESS\_EQUAL}}
\begin{quote}

\sphinxAtStartPar
For constraint in the form of \(g(x) \leq b\)
\end{quote}

\item {} 
\sphinxAtStartPar
\sphinxcode{\sphinxupquote{COPT\_GREATER\_EQUAL}}
\begin{quote}

\sphinxAtStartPar
For constraint in the form of \(g(x) \geq b\)
\end{quote}

\item {} 
\sphinxAtStartPar
\sphinxcode{\sphinxupquote{COPT\_EQUAL}}
\begin{quote}

\sphinxAtStartPar
For constraint in the form of \(g(x) = b\)
\end{quote}

\end{itemize}

\sphinxAtStartPar
In additional, there are two less used constraint senses:
\begin{itemize}
\item {} 
\sphinxAtStartPar
\sphinxcode{\sphinxupquote{COPT\_FREE}}
\begin{quote}

\sphinxAtStartPar
For unconstrained expression
\end{quote}

\item {} 
\sphinxAtStartPar
\sphinxcode{\sphinxupquote{COPT\_RANGE}}
\begin{quote}

\sphinxAtStartPar
For constraints with both lower and upper bounds in the form of
\(l \leq g(x) \leq u\).

\sphinxAtStartPar
Please refer to documentation of \sphinxcode{\sphinxupquote{COPT\_LoadProb}} regarding how to
use \sphinxcode{\sphinxupquote{COPT\_RANGE}} to define a constraints with both lower and upper bounds.
\end{quote}

\end{itemize}

\subsection{Variable types}
\label{\detokenize{capiref:variable-types}}
\sphinxAtStartPar
Variable types are used for defining whether
a variable is continuous or integral.
\begin{itemize}
\item {} 
\sphinxAtStartPar
\sphinxcode{\sphinxupquote{COPT\_CONTINUOUS}}
\begin{quote}

\sphinxAtStartPar
Non\sphinxhyphen{}integer continuous variables
\end{quote}

\item {} 
\sphinxAtStartPar
\sphinxcode{\sphinxupquote{COPT\_BINARY}}
\begin{quote}

\sphinxAtStartPar
Binary variables
\end{quote}

\item {} 
\sphinxAtStartPar
\sphinxcode{\sphinxupquote{COPT\_INTEGER}}
\begin{quote}

\sphinxAtStartPar
Integer variables
\end{quote}

\end{itemize}

\subsection{SOS\sphinxhyphen{}constraint types}
\label{\detokenize{capiref:sos-constraint-types}}
\sphinxAtStartPar
SOS constraint (Special Ordered Set) is a kind of special constraint that
places restrictions on the values that a set of variables can take.

\sphinxAtStartPar
COPT currently support two types of SOS constraints, one is SOS1 constraint,
where at most one variable in the constraint is allowed to take a non\sphinxhyphen{}zero value,
the other is SOS2 constraint, where at most two variables in the constraint are
allowed to take non\sphinxhyphen{}zero value, and those non\sphinxhyphen{}zero variables must be contiguous.
Variables in SOS constraints are allowed to be continuous, binary and integer.
\begin{itemize}
\item {} 
\sphinxAtStartPar
\sphinxcode{\sphinxupquote{COPT\_SOS\_TYPE1}}
\begin{quote}

\sphinxAtStartPar
SOS1 constraint
\end{quote}

\item {} 
\sphinxAtStartPar
\sphinxcode{\sphinxupquote{COPT\_SOS\_TYPE2}}
\begin{quote}

\sphinxAtStartPar
SOS2 constraint
\end{quote}

\end{itemize}

\subsection{Indicator constraint}
\label{\detokenize{capiref:indicator-constraint}}
\sphinxAtStartPar
Indicator constraint is a kind of logical constraints in COPT, used to describe the relationship between the value of
the binary variable \(y\) and whether the linear constraint \(a^{T}x \leq b\) is satisfied.
Currently, COPT supports three types of indicator constraints:
\begin{itemize}
\item {} 
\sphinxAtStartPar
\sphinxcode{\sphinxupquote{COPT\_INDICATOR\_IF}}

\sphinxAtStartPar
If\sphinxhyphen{}Then:
\begin{quote}

\sphinxAtStartPar
If \(y=f\) , then the linear constraint \(a^{T}x \leq b\) is satisfied.

\sphinxAtStartPar
If \(y\ne f\) , then the linear constraint \(a^{T}x \leq b\) is invalid (may be violated).
\end{quote}

\end{itemize}
\begin{equation}\label{equation:capiref:capiref:0}
\begin{split}y &= f \rightarrow a^{T}x \leq b\\
f &\in \{0, 1\}\end{split}
\end{equation}\begin{itemize}
\item {} 
\sphinxAtStartPar
\sphinxcode{\sphinxupquote{COPT\_INDICATOR\_ONLYIF}}

\sphinxAtStartPar
Only\sphinxhyphen{}If:
\begin{quote}

\sphinxAtStartPar
If the linear constraint \(a^{T}x \leq b\) is satisfied, then \(y=f\) .

\sphinxAtStartPar
If the linear constraint \(a^{T}x \leq b\) is not satisfied, then \(y\) can be 0 or 1.
\end{quote}

\end{itemize}
\begin{equation}\label{equation:capiref:capiref:1}
\begin{split}a^{T}x &\leq b \rightarrow y = f\\
f &\in \{0, 1\}\end{split}
\end{equation}\begin{itemize}
\item {} 
\sphinxAtStartPar
\sphinxcode{\sphinxupquote{COPT\_INDICATOR\_IFANDONLYIF}}

\sphinxAtStartPar
If\sphinxhyphen{}and\sphinxhyphen{}Only\sphinxhyphen{}If:
\begin{quote}

\sphinxAtStartPar
The linear constraint \(a^{T}x \leq b\) and \(y=f\) must be satisfied simultaneously or not satisfied simultaneously.
\end{quote}

\end{itemize}
\begin{equation}\label{equation:capiref:capiref:2}
\begin{split}a^{T}x &\leq b \leftrightarrow y = f\\
f &\in \{0, 1\}\end{split}
\end{equation}

\subsection{SOC constraint type}
\label{\detokenize{capiref:soc-constraint-type}}
\sphinxAtStartPar
The Second\sphinxhyphen{}Order\sphinxhyphen{}Cone (SOC) constraint is a special type of quadratic constraints.
COPT supports two types of SOC constraints:
\begin{itemize}
\item {} 
\sphinxAtStartPar
\sphinxcode{\sphinxupquote{COPT\_CONE\_QUAD}} : Regular Second\sphinxhyphen{}Order\sphinxhyphen{}Cone

\end{itemize}
\begin{equation}\label{equation:capiref:capiref:3}
\begin{split}Q^n= \left\{x\in \mathbb{R}^n\ \left|\ x_0\geq\sqrt{\sum_{i=1}^{n-1} x_i^2}, x_0\geq0 \right. \right\}\end{split}
\end{equation}\begin{itemize}
\item {} 
\sphinxAtStartPar
\sphinxcode{\sphinxupquote{COPT\_CONE\_RQUAD}} : Rotated Second\sphinxhyphen{}Order\sphinxhyphen{}Cone

\end{itemize}
\begin{equation}\label{equation:capiref:capiref:4}
\begin{split}Q^n_r= \left\{x\in \mathbb{R}^n\ \left|\ 2x_0x_1\geq\sum_{i=2}^{n-1} x_i^2, x_0\geq0, x_1\geq 0 \right. \right\}\end{split}
\end{equation}

\subsection{Exponential cone constraint type}
\label{\detokenize{capiref:exponential-cone-constraint-type}}
\sphinxAtStartPar
COPT supports two types of exponential cone contraints:
\begin{itemize}
\item {} 
\sphinxAtStartPar
\sphinxcode{\sphinxupquote{COPT\_EXPCONE\_PRIMAL}} : Primal exponential cone

\end{itemize}
\begin{equation}\label{equation:capiref:capiref:5}
\begin{split}\mathrm{cl}(S_1) = S_1 \cup S_2\end{split}
\end{equation}\begin{equation}\label{equation:capiref:capiref:6}
\begin{split}\begin{aligned}
S_1 &= \left\{\begin{pmatrix} t \\ s \\ r \end{pmatrix}\in \mathbb{R}^3\ |\ s > 0,\ t \geq s\ \mathrm{exp}\left(\frac{r}{s} \right) \right\}, \\
S_2 &= \left\{\begin{pmatrix} t \\ s \\ r \end{pmatrix}\in \mathbb{R}^3\ |\ s=0,\ t\geq 0,\ r\leq 0 \right\}
\end{aligned}\end{split}
\end{equation}\begin{itemize}
\item {} 
\sphinxAtStartPar
\sphinxcode{\sphinxupquote{COPT\_EXPCONE\_DUAL}} : Dual exponential cone

\end{itemize}
\begin{equation}\label{equation:capiref:capiref:7}
\begin{split}\mathrm{cl}(S_1) = S_1 \cup S_2\end{split}
\end{equation}\begin{equation}\label{equation:capiref:capiref:8}
\begin{split}\begin{aligned}
S_1 &= \left\{\begin{pmatrix} t \\ s \\ r \end{pmatrix}\in \mathbb{R}^3\ |\ r < 0,\ t \geq -r\ \mathrm{exp}\left(\frac{s}{r}-1\right) \right\}, \\
S_2 &= \left\{\begin{pmatrix} t \\ s \\ r \end{pmatrix}\in \mathbb{R}^3\ |\ r = 0,\ t\geq 0,\ s\geq 0 \right\}
\end{aligned}\end{split}
\end{equation}

\subsection{Quadratic objective function}
\label{\detokenize{capiref:quadratic-objective-function}}
\sphinxAtStartPar
Besides linear objective function, COPT also supports general convex quadratic objective function.

\sphinxAtStartPar
The mathematical form is:
\begin{equation}\label{equation:capiref:coptEq_capiref_qobj}
\begin{split}x^{T}Qx + c^{T}x\end{split}
\end{equation}
\sphinxAtStartPar
Where, \(x\) is an array of variables, \(Q\) is the quadratic part of the quadratic
objective funtion and \(c\) is the linear part.

\subsection{Quadratic constraint}
\label{\detokenize{capiref:quadratic-constraint}}
\sphinxAtStartPar
Besides the special type of quadratic constraint, Second\sphinxhyphen{}Order\sphinxhyphen{}Cone (SOC) constraint, COPT
also supports general convex quadratic constraint.

\sphinxAtStartPar
The mathematical form is:
\begin{equation}\label{equation:capiref:coptEq_capiref_qconstr}
\begin{split}x^{T}Qx + q^{T}x \leq b\end{split}
\end{equation}
\sphinxAtStartPar
Where, \(x\) is an array of variables, \(Q\) is the quadratic part of the quadratic
constraint and \(c\) is the linear part.

\subsection{Basis status}
\label{\detokenize{capiref:basis-status}}
\sphinxAtStartPar
For an LP problem with \(n\) variables and \(m\) constraints,
the constraints are treated as slack variables \sphinxstyleemphasis{internally},
resulting in \(n+m\) variables.
When solving an LP problem using the simplex method,
the simplex method fixes \(n\) variables at one of their bounds,
and then computes solutions for the other \(m\) variables.
The \(m\) variables with computed solution are called \sphinxstyleemphasis{basic} variables,
and the other \(n\) variables are called \sphinxstyleemphasis{non\sphinxhyphen{}basic} variables.
The simplex progress and its final solution can be defined using
the basis status of all the variables and constraints.

\sphinxAtStartPar
The basis status supported by COPT are:
\begin{itemize}
\item {} 
\sphinxAtStartPar
\sphinxcode{\sphinxupquote{COPT\_BASIS\_LOWER}}
\begin{quote}

\sphinxAtStartPar
The variable is non\sphinxhyphen{}basic at its lower bound.
\end{quote}

\item {} 
\sphinxAtStartPar
\sphinxcode{\sphinxupquote{COPT\_BASIS\_BASIC}}
\begin{quote}

\sphinxAtStartPar
The variable is basic.
\end{quote}

\item {} 
\sphinxAtStartPar
\sphinxcode{\sphinxupquote{COPT\_BASIS\_UPPER}}
\begin{quote}

\sphinxAtStartPar
The variable is non\sphinxhyphen{}basic at its upper bound.
\end{quote}

\item {} 
\sphinxAtStartPar
\sphinxcode{\sphinxupquote{COPT\_BASIS\_SUPERBASIC}}
\begin{quote}

\sphinxAtStartPar
The variable is non\sphinxhyphen{}basic but not any of its bounds.
\end{quote}

\item {} 
\sphinxAtStartPar
\sphinxcode{\sphinxupquote{COPT\_BASIS\_FIXED}}
\begin{quote}

\sphinxAtStartPar
The variable is non\sphinxhyphen{}basic and fixed at its bound.
\end{quote}

\end{itemize}

\subsection{LP solution status}
\label{\detokenize{capiref:lp-solution-status}}\label{\detokenize{capiref:chapapi-const-lpstatus}}
\sphinxAtStartPar
The solution status of an LP problem is called LP status,
which can be obtained using integer attribute
\sphinxcode{\sphinxupquote{COPT\_INTATTR\_LPSTATUS}}.

\sphinxAtStartPar
Possible LP status values are:
\begin{itemize}
\item {} 
\sphinxAtStartPar
\sphinxcode{\sphinxupquote{COPT\_LPSTATUS\_UNSTARTED}}
\begin{quote}

\sphinxAtStartPar
The LP optimization is not started yet.
\end{quote}

\item {} 
\sphinxAtStartPar
\sphinxcode{\sphinxupquote{COPT\_LPSTATUS\_OPTIMAL}}
\begin{quote}

\sphinxAtStartPar
The LP problem is solved to optimality.
\end{quote}

\item {} 
\sphinxAtStartPar
\sphinxcode{\sphinxupquote{COPT\_LPSTATUS\_INFEASIBLE}}
\begin{quote}

\sphinxAtStartPar
The LP problem is infeasible.
\end{quote}

\item {} 
\sphinxAtStartPar
\sphinxcode{\sphinxupquote{COPT\_LPSTATUS\_UNBOUNDED}}
\begin{quote}

\sphinxAtStartPar
The LP problem is unbounded.
\end{quote}

\item {} 
\sphinxAtStartPar
\sphinxcode{\sphinxupquote{COPT\_LPSTATUS\_NUMERICAL}}
\begin{quote}

\sphinxAtStartPar
Numerical trouble encountered.
\end{quote}

\item {} 
\sphinxAtStartPar
\sphinxcode{\sphinxupquote{COPT\_LPSTATUS\_TIMEOUT}}
\begin{quote}

\sphinxAtStartPar
The LP optimization is stopped because of time limit.
\end{quote}

\item {} 
\sphinxAtStartPar
\sphinxcode{\sphinxupquote{COPT\_LPSTATUS\_UNFINISHED}}
\begin{quote}

\sphinxAtStartPar
The LP optimization is stopped but the solver cannot provide
a solution because of numerical difficulties.
\end{quote}

\item {} 
\sphinxAtStartPar
\sphinxcode{\sphinxupquote{COPT\_LPSTATUS\_IMPRECISE}}
\begin{quote}

\sphinxAtStartPar
The solution is imprecise.
\end{quote}

\item {} 
\sphinxAtStartPar
\sphinxcode{\sphinxupquote{COPT\_LPSTATUS\_INTERRUPTED}}
\begin{quote}

\sphinxAtStartPar
The LP optimization is stopped by user interrupt.
\end{quote}

\item {} 
\sphinxAtStartPar
\sphinxcode{\sphinxupquote{COPT\_LPSTATUS\_ITERLIMIT}}
\begin{quote}

\sphinxAtStartPar
The LP optimization is stopped because of iteration limit.
\end{quote}

\end{itemize}

\subsection{MIP solution status}
\label{\detokenize{capiref:mip-solution-status}}\label{\detokenize{capiref:chapapi-const-mipstatus}}
\sphinxAtStartPar
The solution status of an MIP problem is called MIP status,
which can be obtained using integer attribute
\sphinxcode{\sphinxupquote{COPT\_INTATTR\_MIPSTATUS}}.

\sphinxAtStartPar
Possible MIP status values are:
\begin{itemize}
\item {} 
\sphinxAtStartPar
\sphinxcode{\sphinxupquote{COPT\_MIPSTATUS\_UNSTARTED}}
\begin{quote}

\sphinxAtStartPar
The MIP optimization is not started yet.
\end{quote}

\item {} 
\sphinxAtStartPar
\sphinxcode{\sphinxupquote{COPT\_MIPSTATUS\_OPTIMAL}}
\begin{quote}

\sphinxAtStartPar
The MIP problem is solved to optimality.
\end{quote}

\item {} 
\sphinxAtStartPar
\sphinxcode{\sphinxupquote{COPT\_MIPSTATUS\_INFEASIBLE}}
\begin{quote}

\sphinxAtStartPar
The MIP problem is infeasible.
\end{quote}

\item {} 
\sphinxAtStartPar
\sphinxcode{\sphinxupquote{COPT\_MIPSTATUS\_UNBOUNDED}}
\begin{quote}

\sphinxAtStartPar
The MIP problem is unbounded.
\end{quote}

\item {} 
\sphinxAtStartPar
\sphinxcode{\sphinxupquote{COPT\_MIPSTATUS\_INF\_OR\_UNB}}
\begin{quote}

\sphinxAtStartPar
The MIP problem is infeasible or unbounded.
\end{quote}

\item {} 
\sphinxAtStartPar
\sphinxcode{\sphinxupquote{COPT\_MIPSTATUS\_NODELIMIT}}
\begin{quote}

\sphinxAtStartPar
The MIP optimization is stopped because of node limit.
\end{quote}

\item {} 
\sphinxAtStartPar
\sphinxcode{\sphinxupquote{COPT\_MIPSTATUS\_TIMEOUT}}
\begin{quote}

\sphinxAtStartPar
The MIP optimization is stopped because of time limit.
\end{quote}

\item {} 
\sphinxAtStartPar
\sphinxcode{\sphinxupquote{COPT\_MIPSTATUS\_UNFINISHED}}
\begin{quote}

\sphinxAtStartPar
The MIP optimization is stopped but the solver cannot provide
a solution because of numerical difficulties.
\end{quote}

\item {} 
\sphinxAtStartPar
\sphinxcode{\sphinxupquote{COPT\_MIPSTATUS\_INTERRUPTED}}
\begin{quote}

\sphinxAtStartPar
The MIP optimization is stopped by user interrupt.
\end{quote}

\end{itemize}

\subsection{Callback context}
\label{\detokenize{capiref:callback-context}}\label{\detokenize{capiref:chapapi-const-cbc}}\begin{itemize}
\item {} 
\sphinxAtStartPar
\sphinxcode{\sphinxupquote{CBCONTEXT\_INCUMBENT}}
\begin{quote}

\sphinxAtStartPar
Invokes the callback after a new incumbent was found.
\end{quote}

\item {} 
\sphinxAtStartPar
\sphinxcode{\sphinxupquote{COPT\_CBCONTEXT\_MIPNODE}}
\begin{quote}

\sphinxAtStartPar
Invokes the callback after a MIP node was processed.
\end{quote}

\item {} 
\sphinxAtStartPar
\sphinxcode{\sphinxupquote{COPT\_CBCONTEXT\_MIPRELAX}}
\begin{quote}

\sphinxAtStartPar
Invokes the callback when an LP\sphinxhyphen{}relaxation was solved.
\end{quote}

\item {} 
\sphinxAtStartPar
\sphinxcode{\sphinxupquote{COPT\_CBCONTEXT\_MIPSOL}}
\begin{quote}

\sphinxAtStartPar
Invokes the callback when a new MIP candidate solution is found.
\end{quote}

\end{itemize}

\subsection{Nonlinear Expression Operators}
\label{\detokenize{capiref:nonlinear-expression-operators}}\label{\detokenize{capiref:chapapi-const-nlp}}\begin{itemize}
\item {} 
\sphinxAtStartPar
\sphinxcode{\sphinxupquote{COPT\_NL\_PLUS}}
\begin{quote}

\sphinxAtStartPar
Addition operator.
\end{quote}

\item {} 
\sphinxAtStartPar
\sphinxcode{\sphinxupquote{COPT\_NL\_MINUS}}
\begin{quote}

\sphinxAtStartPar
Subtraction operator.
\end{quote}

\item {} 
\sphinxAtStartPar
\sphinxcode{\sphinxupquote{COPT\_NL\_MULT}}
\begin{quote}

\sphinxAtStartPar
Multiplication operator.
\end{quote}

\item {} 
\sphinxAtStartPar
\sphinxcode{\sphinxupquote{COPT\_NL\_DIV}}
\begin{quote}

\sphinxAtStartPar
Division operator.
\end{quote}

\item {} 
\sphinxAtStartPar
\sphinxcode{\sphinxupquote{COPT\_NL\_POW}}
\begin{quote}

\sphinxAtStartPar
Power operator.
\end{quote}

\item {} 
\sphinxAtStartPar
\sphinxcode{\sphinxupquote{COPT\_NL\_SQRT}}
\begin{quote}

\sphinxAtStartPar
Square root operator.
\end{quote}

\item {} 
\sphinxAtStartPar
\sphinxcode{\sphinxupquote{COPT\_NL\_EXP}}
\begin{quote}

\sphinxAtStartPar
Exponential function operator.
\end{quote}

\item {} 
\sphinxAtStartPar
\sphinxcode{\sphinxupquote{COPT\_NL\_LOG}}
\begin{quote}

\sphinxAtStartPar
Natural logarithm operator.
\end{quote}

\item {} 
\sphinxAtStartPar
\sphinxcode{\sphinxupquote{COPT\_NL\_LOG10}}
\begin{quote}

\sphinxAtStartPar
Base\sphinxhyphen{}10 logarithm operator.
\end{quote}

\item {} 
\sphinxAtStartPar
\sphinxcode{\sphinxupquote{COPT\_NL\_NEG}}
\begin{quote}

\sphinxAtStartPar
Unary negation operator.
\end{quote}

\item {} 
\sphinxAtStartPar
\sphinxcode{\sphinxupquote{COPT\_NL\_ABS}}
\begin{quote}

\sphinxAtStartPar
Absolute value operator.
\end{quote}

\item {} 
\sphinxAtStartPar
\sphinxcode{\sphinxupquote{COPT\_NL\_FLOOR}}
\begin{quote}

\sphinxAtStartPar
Floor function operator.
\end{quote}

\item {} 
\sphinxAtStartPar
\sphinxcode{\sphinxupquote{COPT\_NL\_CEIL}}
\begin{quote}

\sphinxAtStartPar
Ceiling function operator.
\end{quote}

\item {} 
\sphinxAtStartPar
\sphinxcode{\sphinxupquote{COPT\_NL\_SIN}}
\begin{quote}

\sphinxAtStartPar
Sine function operator.
\end{quote}

\item {} 
\sphinxAtStartPar
\sphinxcode{\sphinxupquote{COPT\_NL\_COS}}
\begin{quote}

\sphinxAtStartPar
Cosine function operator.
\end{quote}

\item {} 
\sphinxAtStartPar
\sphinxcode{\sphinxupquote{COPT\_NL\_TAN}}
\begin{quote}

\sphinxAtStartPar
Tangent function operator.
\end{quote}

\item {} 
\sphinxAtStartPar
\sphinxcode{\sphinxupquote{COPT\_NL\_SINH}}
\begin{quote}

\sphinxAtStartPar
Hyperbolic sine function operator.
\end{quote}

\item {} 
\sphinxAtStartPar
\sphinxcode{\sphinxupquote{COPT\_NL\_COSH}}
\begin{quote}

\sphinxAtStartPar
Hyperbolic cosine function operator.
\end{quote}

\item {} 
\sphinxAtStartPar
\sphinxcode{\sphinxupquote{COPT\_NL\_TANH}}
\begin{quote}

\sphinxAtStartPar
Hyperbolic tangent function operator.
\end{quote}

\item {} 
\sphinxAtStartPar
\sphinxcode{\sphinxupquote{COPT\_NL\_ASIN}}
\begin{quote}

\sphinxAtStartPar
Inverse sine (arcsin) function operator.
\end{quote}

\item {} 
\sphinxAtStartPar
\sphinxcode{\sphinxupquote{COPT\_NL\_ACOS}}
\begin{quote}

\sphinxAtStartPar
Inverse cosine (arccos) function operator.
\end{quote}

\item {} 
\sphinxAtStartPar
\sphinxcode{\sphinxupquote{COPT\_NL\_ATAN}}
\begin{quote}

\sphinxAtStartPar
Inverse tangent (arctan) function operator.
\end{quote}

\item {} 
\sphinxAtStartPar
\sphinxcode{\sphinxupquote{COPT\_NL\_ASINH}}
\begin{quote}

\sphinxAtStartPar
Inverse hyperbolic sine (arsinh) function operator.
\end{quote}

\item {} 
\sphinxAtStartPar
\sphinxcode{\sphinxupquote{COPT\_NL\_ACOSH}}
\begin{quote}

\sphinxAtStartPar
Inverse hyperbolic cosine (arcosh) function operator.
\end{quote}

\item {} 
\sphinxAtStartPar
\sphinxcode{\sphinxupquote{COPT\_NL\_ATANH}}
\begin{quote}

\sphinxAtStartPar
Inverse hyperbolic tangent (artanh) function operator.
\end{quote}

\item {} 
\sphinxAtStartPar
\sphinxcode{\sphinxupquote{COPT\_NL\_ATAN2}}
\begin{quote}

\sphinxAtStartPar
Two\sphinxhyphen{}argument inverse tangent function operator.
\end{quote}

\item {} 
\sphinxAtStartPar
\sphinxcode{\sphinxupquote{COPT\_NL\_SUM}}
\begin{quote}

\sphinxAtStartPar
Summation operator.
\end{quote}

\item {} 
\sphinxAtStartPar
\sphinxcode{\sphinxupquote{COPT\_NL\_GET}}
\begin{quote}

\sphinxAtStartPar
Constant retrieval operator.
\end{quote}

\end{itemize}

\subsection{API function return code}
\label{\detokenize{capiref:api-function-return-code}}
\sphinxAtStartPar
When an API function finishes, it returns an integer \sphinxstylestrong{return code},
which indicates whether the API call was finished okay or failed.
In case of failure, it specifies the reason of the failure.

\sphinxAtStartPar
Possible COPT API function return codes are:
\begin{itemize}
\item {} 
\sphinxAtStartPar
\sphinxcode{\sphinxupquote{COPT\_RETCODE\_OK}}
\begin{quote}

\sphinxAtStartPar
The API call finished successfully.
\end{quote}

\item {} 
\sphinxAtStartPar
\sphinxcode{\sphinxupquote{COPT\_RETCODE\_MEMORY}}
\begin{quote}

\sphinxAtStartPar
The API call failed because of memory allocation failure.
\end{quote}

\item {} 
\sphinxAtStartPar
\sphinxcode{\sphinxupquote{COPT\_RETCODE\_FILE}}
\begin{quote}

\sphinxAtStartPar
The API call failed because of file input or output failure.
\end{quote}

\item {} 
\sphinxAtStartPar
\sphinxcode{\sphinxupquote{COPT\_RETCODE\_INVALID}}
\begin{quote}

\sphinxAtStartPar
The API call failed because of invalid data.
\end{quote}

\item {} 
\sphinxAtStartPar
\sphinxcode{\sphinxupquote{COPT\_RETCODE\_LICENSE}}
\begin{quote}

\sphinxAtStartPar
The API call failed because of license validation failure.
In this case, further information can be obtained by calling
\sphinxcode{\sphinxupquote{COPT\_GetLicenseMsg}}.
\end{quote}

\item {} 
\sphinxAtStartPar
\sphinxcode{\sphinxupquote{COPT\_RETCODE\_INTERNAL}}
\begin{quote}

\sphinxAtStartPar
The API call failed because of internal error.
\end{quote}

\item {} 
\sphinxAtStartPar
\sphinxcode{\sphinxupquote{COPT\_RETCODE\_THREAD}}
\begin{quote}

\sphinxAtStartPar
The API call failed because of thread error.
\end{quote}

\item {} 
\sphinxAtStartPar
\sphinxcode{\sphinxupquote{COPT\_RETCODE\_SERVER}}
\begin{quote}

\sphinxAtStartPar
The API call failed because of remote server error.
\end{quote}

\item {} 
\sphinxAtStartPar
\sphinxcode{\sphinxupquote{COPT\_RETCODE\_NONCONVEX}}
\begin{quote}

\sphinxAtStartPar
The API call failed because of problem is nonconvex.
\end{quote}

\item {} 
\sphinxAtStartPar
\sphinxcode{\sphinxupquote{COPT\_RETCODE\_MEMORY\_GPU}}
\begin{quote}

\sphinxAtStartPar
The API call failed because of GPU memory allocation failure.
\end{quote}

\end{itemize}

\subsection{Client configuration}
\label{\detokenize{capiref:client-configuration}}
\sphinxAtStartPar
For floating and cluster clients, users are allowed to set client configuration parameters,
currently available settings are:
\begin{itemize}
\item {} 
\sphinxAtStartPar
\sphinxcode{\sphinxupquote{COPT\_CLIENT\_CLUSTER}}
\begin{quote}

\sphinxAtStartPar
IP address of cluster server.
\end{quote}

\item {} 
\sphinxAtStartPar
\sphinxcode{\sphinxupquote{COPT\_CLIENT\_FLOATING}}
\begin{quote}

\sphinxAtStartPar
IP address of token server.
\end{quote}

\item {} 
\sphinxAtStartPar
\sphinxcode{\sphinxupquote{COPT\_CLIENT\_PASSWORD}}
\begin{quote}

\sphinxAtStartPar
Password of cluster server.
\end{quote}

\item {} 
\sphinxAtStartPar
\sphinxcode{\sphinxupquote{COPT\_CLIENT\_PORT}}
\begin{quote}

\sphinxAtStartPar
Connection port of token server.
\end{quote}

\item {} 
\sphinxAtStartPar
\sphinxcode{\sphinxupquote{COPT\_CLIENT\_WAITTIME}}
\begin{quote}

\sphinxAtStartPar
Wait time of client.
\end{quote}

\end{itemize}

\subsection{Other constants}
\label{\detokenize{capiref:other-constants}}\begin{itemize}
\item {} 
\sphinxAtStartPar
\sphinxcode{\sphinxupquote{COPT\_BUFFSIZE}}
\begin{quote}

\sphinxAtStartPar
Defines the recommended buffer size when obtaining a
C\sphinxhyphen{}style string message from COPT library.
This can be used with, for example,
\sphinxcode{\sphinxupquote{COPT\_GetBanner}}, \sphinxcode{\sphinxupquote{COPT\_GetRetcodeMsg}} etc.
\end{quote}

\end{itemize}

\section{Attributes}
\label{\detokenize{capiref:attributes}}\label{\detokenize{capiref:chapapi-attrs}}

\subsection{Problem information}
\label{\detokenize{capiref:problem-information}}\begin{itemize}
\item {} 
\sphinxAtStartPar
\sphinxcode{\sphinxupquote{COPT\_INTATTR\_COLS}} or \sphinxcode{\sphinxupquote{"Cols"}}
\begin{quote}

\sphinxAtStartPar
Integer attribute.

\sphinxAtStartPar
Number of variables (columns) in the problem.
\end{quote}

\item {} 
\sphinxAtStartPar
\sphinxcode{\sphinxupquote{COPT\_INTATTR\_PSDCOLS}} or \sphinxcode{\sphinxupquote{"PSDCols"}}
\begin{quote}

\sphinxAtStartPar
Integer attribute.

\sphinxAtStartPar
Number of PSD variables in the problem.
\end{quote}

\item {} 
\sphinxAtStartPar
\sphinxcode{\sphinxupquote{COPT\_INTATTR\_ROWS}} or \sphinxcode{\sphinxupquote{"Rows"}}
\begin{quote}

\sphinxAtStartPar
Integer attribute.

\sphinxAtStartPar
Number of constraints (rows) in the problem.
\end{quote}

\item {} 
\sphinxAtStartPar
\sphinxcode{\sphinxupquote{COPT\_INTATTR\_ELEMS}} or \sphinxcode{\sphinxupquote{"Elems"}}
\begin{quote}

\sphinxAtStartPar
Integer attribute.

\sphinxAtStartPar
Number of non\sphinxhyphen{}zero elements in the coefficient matrix.
\end{quote}

\item {} 
\sphinxAtStartPar
\sphinxcode{\sphinxupquote{COPT\_INTATTR\_QELEMS}} or \sphinxcode{\sphinxupquote{"QElems"}}
\begin{quote}

\sphinxAtStartPar
Integer attribute.

\sphinxAtStartPar
Number of non\sphinxhyphen{}zero quadratic elements in the quadratic objective function.
\end{quote}

\item {} 
\sphinxAtStartPar
\sphinxcode{\sphinxupquote{COPT\_INTATTR\_PSDELEMS}} or \sphinxcode{\sphinxupquote{"PSDElems"}}
\begin{quote}

\sphinxAtStartPar
Integer attribute.

\sphinxAtStartPar
Number of PSD terms in objective function.
\end{quote}

\item {} 
\sphinxAtStartPar
\sphinxcode{\sphinxupquote{COPT\_INTATTR\_SYMMATS}} or \sphinxcode{\sphinxupquote{"SymMats"}}
\begin{quote}

\sphinxAtStartPar
Integer attribute.

\sphinxAtStartPar
Number of symmetric matrices in the problem.
\end{quote}

\item {} 
\sphinxAtStartPar
\sphinxcode{\sphinxupquote{COPT\_INTATTR\_BINS}} or \sphinxcode{\sphinxupquote{"Bins"}}
\begin{quote}

\sphinxAtStartPar
Integer attribute.

\sphinxAtStartPar
Number of binary variables.
\end{quote}

\item {} 
\sphinxAtStartPar
\sphinxcode{\sphinxupquote{COPT\_INTATTR\_INTS}} or \sphinxcode{\sphinxupquote{"Ints"}}
\begin{quote}

\sphinxAtStartPar
Integer attribute.

\sphinxAtStartPar
Number of integer variables.
\end{quote}

\item {} 
\sphinxAtStartPar
\sphinxcode{\sphinxupquote{COPT\_INTATTR\_SOSS}} or \sphinxcode{\sphinxupquote{"Soss"}}
\begin{quote}

\sphinxAtStartPar
Integer attribute.

\sphinxAtStartPar
Number of SOS constraints.
\end{quote}

\item {} 
\sphinxAtStartPar
\sphinxcode{\sphinxupquote{COPT\_INTATTR\_CONES}} or \sphinxcode{\sphinxupquote{"Cones"}}
\begin{quote}

\sphinxAtStartPar
Integer attribute.

\sphinxAtStartPar
Number of Second\sphinxhyphen{}Order\sphinxhyphen{}Cone constraints.
\end{quote}

\item {} 
\sphinxAtStartPar
\sphinxcode{\sphinxupquote{COPT\_INTATTR\_EXPCONES}} or \sphinxcode{\sphinxupquote{"ExpCones"}}
\begin{quote}

\sphinxAtStartPar
Integer attribute.

\sphinxAtStartPar
Number of exponential cone constraints.
\end{quote}

\item {} 
\sphinxAtStartPar
\sphinxcode{\sphinxupquote{COPT\_INTATTR\_AFFINECONES}} or \sphinxcode{\sphinxupquote{"AffineCones"}}
\begin{quote}

\sphinxAtStartPar
Integer attribute.

\sphinxAtStartPar
Number of affine cone constraints.
\end{quote}

\item {} 
\sphinxAtStartPar
\sphinxcode{\sphinxupquote{COPT\_INTATTR\_QCONSTRS}} or \sphinxcode{\sphinxupquote{"QConstrs"}}
\begin{quote}

\sphinxAtStartPar
Integer attribute.

\sphinxAtStartPar
Number of quadratic constraints.
\end{quote}

\item {} 
\sphinxAtStartPar
\sphinxcode{\sphinxupquote{COPT\_INTATTR\_PSDCONSTRS}} or \sphinxcode{\sphinxupquote{"PSDConstrs"}}
\begin{quote}

\sphinxAtStartPar
Integer attribute.

\sphinxAtStartPar
Number of PSD constraints.
\end{quote}

\item {} 
\sphinxAtStartPar
\sphinxcode{\sphinxupquote{COPT\_INTATTR\_LMICONSTRS}} or \sphinxcode{\sphinxupquote{"LMIConstrs"}}
\begin{quote}

\sphinxAtStartPar
Integer attribute.

\sphinxAtStartPar
Number of LMI constraints.
\end{quote}

\item {} 
\sphinxAtStartPar
\sphinxcode{\sphinxupquote{COPT\_INTATTR\_INDICATORS}} or \sphinxcode{\sphinxupquote{"Indicators"}}
\begin{quote}

\sphinxAtStartPar
Integer attribute.

\sphinxAtStartPar
Number of indicator constraints.
\end{quote}

\item {} 
\sphinxAtStartPar
\sphinxcode{\sphinxupquote{COPT\_INTATTR\_OBJSENSE}} or \sphinxcode{\sphinxupquote{"ObjSense"}}
\begin{quote}

\sphinxAtStartPar
Integer attribute.

\sphinxAtStartPar
The optimization direction.
\end{quote}

\item {} 
\sphinxAtStartPar
\sphinxcode{\sphinxupquote{COPT\_DBLATTR\_OBJCONST}} or \sphinxcode{\sphinxupquote{"ObjConst"}}
\begin{quote}

\sphinxAtStartPar
Double attribute.

\sphinxAtStartPar
The constant part of the objective function.
\end{quote}

\item {} 
\sphinxAtStartPar
\sphinxcode{\sphinxupquote{COPT\_INTATTR\_HASQOBJ}} or \sphinxcode{\sphinxupquote{"HasQObj"}}
\begin{quote}

\sphinxAtStartPar
Integer attribute.

\sphinxAtStartPar
Whether the problem has quadratic objective function.
\end{quote}

\item {} 
\sphinxAtStartPar
\sphinxcode{\sphinxupquote{COPT\_INTATTR\_HASPSDOBJ}} or \sphinxcode{\sphinxupquote{"HasPSDObj"}}
\begin{quote}

\sphinxAtStartPar
Integer attribute.

\sphinxAtStartPar
Whether the problem has PSD terms in objective function.
\end{quote}

\item {} 
\sphinxAtStartPar
\sphinxcode{\sphinxupquote{COPT\_INTATTR\_ISMIP}} or \sphinxcode{\sphinxupquote{"IsMIP"}}
\begin{quote}

\sphinxAtStartPar
Integer attribute.

\sphinxAtStartPar
Whether the problem is a MIP.
\end{quote}

\item {} 
\sphinxAtStartPar
\sphinxcode{\sphinxupquote{COPT\_INTATTR\_NLELEMS}} or \sphinxcode{\sphinxupquote{"NLElems"}}
\begin{quote}

\sphinxAtStartPar
Integer attribute.

\sphinxAtStartPar
The number of nonlinear expression terms in the objective of the model.
\end{quote}

\item {} 
\sphinxAtStartPar
\sphinxcode{\sphinxupquote{COPT\_INTATTR\_NLCONSTRS}} or \sphinxcode{\sphinxupquote{"NLConstrs"}}
\begin{quote}

\sphinxAtStartPar
Integer attribute.

\sphinxAtStartPar
The number of nonlinear expression constraints in the model.
\end{quote}

\item {} 
\sphinxAtStartPar
\sphinxcode{\sphinxupquote{COPT\_INTATTR\_HASNLOBJ}} or \sphinxcode{\sphinxupquote{"HasNLObj"}}
\begin{quote}

\sphinxAtStartPar
Integer attribute.

\sphinxAtStartPar
Indicates whether the model has nonlinear expressions in the objective.
\end{quote}

\item {} 
\sphinxAtStartPar
\sphinxcode{\sphinxupquote{COPT\_INTATTR\_MULTIOBJS}} or \sphinxcode{\sphinxupquote{"MultiObjs"}}
\begin{quote}

\sphinxAtStartPar
Integer attribute.

\sphinxAtStartPar
The number of objectives in a multi\sphinxhyphen{}objective model.
\end{quote}

\end{itemize}

\subsection{Solution information}
\label{\detokenize{capiref:solution-information}}\begin{itemize}
\item {} 
\sphinxAtStartPar
\sphinxcode{\sphinxupquote{COPT\_INTATTR\_LPSTATUS}} or \sphinxcode{\sphinxupquote{"LpStatus"}}
\begin{quote}

\sphinxAtStartPar
Integer attribute.

\sphinxAtStartPar
The LP status. Please refer to {\hyperref[\detokenize{capiref:chapapi-const-lpstatus}]{\sphinxcrossref{\DUrole{std,std-ref}{Constants: LP solution status}}}} for possible values.
\end{quote}

\item {} 
\sphinxAtStartPar
\sphinxcode{\sphinxupquote{COPT\_INTATTR\_MIPSTATUS}} or \sphinxcode{\sphinxupquote{"MipStatus"}}
\begin{quote}

\sphinxAtStartPar
Integer attribute.

\sphinxAtStartPar
The MIP status. Please refer to {\hyperref[\detokenize{capiref:chapapi-const-mipstatus}]{\sphinxcrossref{\DUrole{std,std-ref}{Constants: MIP solution status}}}} for possible values.
\end{quote}

\item {} 
\sphinxAtStartPar
\sphinxcode{\sphinxupquote{COPT\_INTATTR\_SIMPLEXITER}} or \sphinxcode{\sphinxupquote{"SimplexIter"}}
\begin{quote}

\sphinxAtStartPar
Integer attribute.

\sphinxAtStartPar
Number of simplex iterations performed.
\end{quote}

\item {} 
\sphinxAtStartPar
\sphinxcode{\sphinxupquote{COPT\_INTATTR\_BARRIERITER}} or \sphinxcode{\sphinxupquote{"BarrierIter"}}
\begin{quote}

\sphinxAtStartPar
Integer attribute.

\sphinxAtStartPar
Number of barrier iterations performed.
\end{quote}

\item {} 
\sphinxAtStartPar
\sphinxcode{\sphinxupquote{COPT\_INTATTR\_NODECNT}} or \sphinxcode{\sphinxupquote{"NodeCnt"}}
\begin{quote}

\sphinxAtStartPar
Integer attribute.

\sphinxAtStartPar
Number of explored nodes.
\end{quote}

\item {} 
\sphinxAtStartPar
\sphinxcode{\sphinxupquote{COPT\_INTATTR\_POOLSOLS}} or \sphinxcode{\sphinxupquote{"PoolSols"}}
\begin{quote}

\sphinxAtStartPar
Integer attribute.

\sphinxAtStartPar
Number of solutions in solution pool.
\end{quote}

\item {} 
\sphinxAtStartPar
\sphinxcode{\sphinxupquote{COPT\_INTATTR\_TUNERESULTS}} or \sphinxcode{\sphinxupquote{"TuneResults"}}
\begin{quote}

\sphinxAtStartPar
Integer attribute.

\sphinxAtStartPar
Number of parameter tuning results
\end{quote}

\item {} 
\sphinxAtStartPar
\sphinxcode{\sphinxupquote{COPT\_INTATTR\_HASLPSOL}} or \sphinxcode{\sphinxupquote{"HasLpSol"}}
\begin{quote}

\sphinxAtStartPar
Integer attribute.

\sphinxAtStartPar
Whether LP solution is available.
\end{quote}

\item {} 
\sphinxAtStartPar
\sphinxcode{\sphinxupquote{COPT\_INTATTR\_HASBASIS}} or \sphinxcode{\sphinxupquote{"HasBasis"}}
\begin{quote}

\sphinxAtStartPar
Integer attribute.

\sphinxAtStartPar
Whether LP basis is available.
\end{quote}

\item {} 
\sphinxAtStartPar
\sphinxcode{\sphinxupquote{COPT\_INTATTR\_HASDUALFARKAS}} or \sphinxcode{\sphinxupquote{"HasDualFarkas"}}
\begin{quote}

\sphinxAtStartPar
Integer attribute.

\sphinxAtStartPar
Whether the dual Farkas of an infeasible LP problem is available.
\end{quote}

\item {} 
\sphinxAtStartPar
\sphinxcode{\sphinxupquote{COPT\_INTATTR\_HASPRIMALRAY}} or \sphinxcode{\sphinxupquote{"HasPrimalRay"}}
\begin{quote}

\sphinxAtStartPar
Integer attribute.

\sphinxAtStartPar
Whether the primal ray of an unbounded LP problem is available.
\end{quote}

\item {} 
\sphinxAtStartPar
\sphinxcode{\sphinxupquote{COPT\_INTATTR\_HASMIPSOL}} or \sphinxcode{\sphinxupquote{"HasMipSol"}}
\begin{quote}

\sphinxAtStartPar
Integer attribute.

\sphinxAtStartPar
Whether MIP solution is available.
\end{quote}

\item {} 
\sphinxAtStartPar
\sphinxcode{\sphinxupquote{COPT\_INTATTR\_IISCOLS}} or \sphinxcode{\sphinxupquote{"IISCols"}}
\begin{quote}

\sphinxAtStartPar
Integer attribute.

\sphinxAtStartPar
Number of bounds of columns in IIS.
\end{quote}

\item {} 
\sphinxAtStartPar
\sphinxcode{\sphinxupquote{COPT\_INTATTR\_IISROWS}} or \sphinxcode{\sphinxupquote{"IISRows"}}
\begin{quote}

\sphinxAtStartPar
Integer attribute.

\sphinxAtStartPar
Number of rows in IIS.
\end{quote}

\item {} 
\sphinxAtStartPar
\sphinxcode{\sphinxupquote{COPT\_INTATTR\_IISSOSS}} or \sphinxcode{\sphinxupquote{"IISSOSs"}}
\begin{quote}

\sphinxAtStartPar
Integer attribute.

\sphinxAtStartPar
Number of SOS constraints in IIS.
\end{quote}

\item {} 
\sphinxAtStartPar
\sphinxcode{\sphinxupquote{COPT\_INTATTR\_IISINDICATORS}} or \sphinxcode{\sphinxupquote{"IISIndicators"}}
\begin{quote}

\sphinxAtStartPar
Integer attribute.

\sphinxAtStartPar
Number of indicator constraints in IIS.
\end{quote}

\item {} 
\sphinxAtStartPar
\sphinxcode{\sphinxupquote{COPT\_INTATTR\_HASIIS}} or \sphinxcode{\sphinxupquote{"HasIIS"}}
\begin{quote}

\sphinxAtStartPar
Integer attribute.

\sphinxAtStartPar
Whether IIS is available.
\end{quote}

\item {} 
\sphinxAtStartPar
\sphinxcode{\sphinxupquote{COPT\_INTATTR\_HASFEASRELAXSOL}} or \sphinxcode{\sphinxupquote{"HasFeasRelaxSol"}}
\begin{quote}

\sphinxAtStartPar
Integer attribute.

\sphinxAtStartPar
Whether feasibility LP\sphinxhyphen{}relaxation solution is available.
\end{quote}

\item {} 
\sphinxAtStartPar
\sphinxcode{\sphinxupquote{COPT\_INTATTR\_ISMINIIS}} or \sphinxcode{\sphinxupquote{"IsMinIIS"}}
\begin{quote}

\sphinxAtStartPar
Integer attribute.

\sphinxAtStartPar
Whether the computed IIS is minimal.
\end{quote}

\item {} 
\sphinxAtStartPar
\sphinxcode{\sphinxupquote{COPT\_INTATTR\_HASSENSITIVITY}} or \sphinxcode{\sphinxupquote{"HasSensitivity"}}
\begin{quote}

\sphinxAtStartPar
Integer attribute.

\sphinxAtStartPar
Whether sensitivity analysis results are available for LP problem.
\end{quote}

\item {} 
\sphinxAtStartPar
\sphinxcode{\sphinxupquote{COPT\_DBLATTR\_LPOBJVAL}} or \sphinxcode{\sphinxupquote{"LpObjval"}}
\begin{quote}

\sphinxAtStartPar
Double attribute.

\sphinxAtStartPar
The LP objective value.
\end{quote}

\item {} 
\sphinxAtStartPar
\sphinxcode{\sphinxupquote{COPT\_DBLATTR\_BESTOBJ}} or \sphinxcode{\sphinxupquote{"BestObj"}}
\begin{quote}

\sphinxAtStartPar
Double attribute.

\sphinxAtStartPar
Best integer objective value for MIP.
\end{quote}

\item {} 
\sphinxAtStartPar
\sphinxcode{\sphinxupquote{COPT\_DBLATTR\_BESTBND}} or \sphinxcode{\sphinxupquote{"BestBnd"}}
\begin{quote}

\sphinxAtStartPar
Double attribute.

\sphinxAtStartPar
Best bound for MIP.
\end{quote}

\item {} 
\sphinxAtStartPar
\sphinxcode{\sphinxupquote{COPT\_DBLATTR\_BESTGAP}} or \sphinxcode{\sphinxupquote{"BestGap"}}
\begin{quote}

\sphinxAtStartPar
Double attribute.

\sphinxAtStartPar
Best relative gap for MIP.
\end{quote}

\item {} 
\sphinxAtStartPar
\sphinxcode{\sphinxupquote{COPT\_DBLATTR\_FEASRELAXOBJ}} or \sphinxcode{\sphinxupquote{FeasRelaxObj}}
\begin{quote}

\sphinxAtStartPar
Double attribute.

\sphinxAtStartPar
Feasibility relaxation objective value.
\end{quote}

\item {} 
\sphinxAtStartPar
\sphinxcode{\sphinxupquote{COPT\_DBLATTR\_SOLVINGTIME}} or \sphinxcode{\sphinxupquote{"SolvingTime"}}
\begin{quote}

\sphinxAtStartPar
Double attribute.

\sphinxAtStartPar
The time spent for the optimization (in seconds).
\end{quote}

\end{itemize}

\section{Information}
\label{\detokenize{capiref:information}}\label{\detokenize{capiref:chapapi-info}}

\subsection{Problem information}
\label{\detokenize{capiref:id1}}\begin{itemize}
\item {} 
\sphinxAtStartPar
\sphinxcode{\sphinxupquote{COPT\_DBLINFO\_OBJ}} or \sphinxcode{\sphinxupquote{"Obj"}}
\begin{quote}

\sphinxAtStartPar
Double information.

\sphinxAtStartPar
Objective cost of columns.
\end{quote}

\item {} 
\sphinxAtStartPar
\sphinxcode{\sphinxupquote{COPT\_DBLINFO\_LB}} or \sphinxcode{\sphinxupquote{"LB"}}
\begin{quote}

\sphinxAtStartPar
Double information.

\sphinxAtStartPar
Lower bounds of columns or rows.
\end{quote}

\item {} 
\sphinxAtStartPar
\sphinxcode{\sphinxupquote{COPT\_DBLINFO\_UB}} or \sphinxcode{\sphinxupquote{"UB"}}
\begin{quote}

\sphinxAtStartPar
Double information.

\sphinxAtStartPar
Upper bounds of columns or rows.
\end{quote}

\end{itemize}

\subsection{Solution and sensitivity analysis information}
\label{\detokenize{capiref:solution-and-sensitivity-analysis-information}}\begin{itemize}
\item {} 
\sphinxAtStartPar
\sphinxcode{\sphinxupquote{COPT\_DBLINFO\_VALUE}} or \sphinxcode{\sphinxupquote{"Value"}}
\begin{quote}

\sphinxAtStartPar
Double information.

\sphinxAtStartPar
Solution of columns.
\end{quote}

\item {} 
\sphinxAtStartPar
\sphinxcode{\sphinxupquote{COPT\_DBLINFO\_SLACK}} or \sphinxcode{\sphinxupquote{"Slack"}}
\begin{quote}

\sphinxAtStartPar
Double information.

\sphinxAtStartPar
Solution of slack variables, also known as activities of constraints.
Only available for LP problem.
\end{quote}

\item {} 
\sphinxAtStartPar
\sphinxcode{\sphinxupquote{COPT\_DBLINFO\_DUAL}} or \sphinxcode{\sphinxupquote{"Dual"}}
\begin{quote}

\sphinxAtStartPar
Double information.

\sphinxAtStartPar
Solution of dual variables. Only available for LP problem.
\end{quote}

\item {} 
\sphinxAtStartPar
\sphinxcode{\sphinxupquote{COPT\_DBLINFO\_REDCOST}} or \sphinxcode{\sphinxupquote{"RedCost"}}
\begin{quote}

\sphinxAtStartPar
Double information.

\sphinxAtStartPar
Reduced cost of columns. Only available for LP problem.
\end{quote}

\item {} 
\sphinxAtStartPar
\sphinxcode{\sphinxupquote{COPT\_DBLINFO\_SAOBJLOW}} or \sphinxcode{\sphinxupquote{"SAObjLow"}}
\begin{quote}

\sphinxAtStartPar
Double information.

\sphinxAtStartPar
Sensitivity analysis information of the objective coefficient.

\sphinxAtStartPar
Indicates the minimum value to which the objective coefficient of a variable
can be reduced while keeping the current basis optimal.
\end{quote}

\item {} 
\sphinxAtStartPar
\sphinxcode{\sphinxupquote{COPT\_DBLINFO\_SAOBJUP}} or \sphinxcode{\sphinxupquote{"SAObjUp"}}
\begin{quote}

\sphinxAtStartPar
Double information.

\sphinxAtStartPar
Sensitivity analysis information of the objective coefficient.

\sphinxAtStartPar
Indicates the maximum value to which the objective coefficient of a variable
can be increased while keeping the current basis optimal.
\end{quote}

\item {} 
\sphinxAtStartPar
\sphinxcode{\sphinxupquote{COPT\_DBLINFO\_SALBLOW}} or \sphinxcode{\sphinxupquote{"SALBLow"}}
\begin{quote}

\sphinxAtStartPar
Double information.

\sphinxAtStartPar
Sensitivity analysis information of the lower bound of the variable/constraint.

\sphinxAtStartPar
Indicates the minimum value to which the lower bound of the variable/constraint
can be reduced while keeping the current basis optimal.
\end{quote}

\item {} 
\sphinxAtStartPar
\sphinxcode{\sphinxupquote{COPT\_DBLINFO\_SALBUP}} or \sphinxcode{\sphinxupquote{"SALBUp"}}
\begin{quote}

\sphinxAtStartPar
Double information.

\sphinxAtStartPar
Sensitivity analysis information of the lower bound of the variable/constraint.

\sphinxAtStartPar
Indicates the maximum value to which the lower bound of the variable/constraint
can be increased while keeping the current basis optimal.
\end{quote}

\item {} 
\sphinxAtStartPar
\sphinxcode{\sphinxupquote{COPT\_DBLINFO\_SAUBLOW}} or \sphinxcode{\sphinxupquote{"SAUBLow"}}
\begin{quote}

\sphinxAtStartPar
Double information.

\sphinxAtStartPar
Sensitivity analysis information of the upper bound of the variable/constraint.

\sphinxAtStartPar
Indicates the minimum value to which the upper bound of the variable/constraint
can be reduced while keeping the current basis optimal.
\end{quote}

\item {} 
\sphinxAtStartPar
\sphinxcode{\sphinxupquote{COPT\_DBLINFO\_SAUBUP}} or \sphinxcode{\sphinxupquote{"SAUBUp"}}
\begin{quote}

\sphinxAtStartPar
Double information.

\sphinxAtStartPar
Sensitivity analysis information of the upper bound of the variable/constraint.

\sphinxAtStartPar
Indicates the maximum value to which the upper bound of the variable/constraint
can be increased while keeping the current basis optimal.
\end{quote}

\end{itemize}

\subsection{Dual Farkas and primal ray}
\label{\detokenize{capiref:dual-farkas-and-primal-ray}}\begin{quote}

\sphinxAtStartPar
Advanced topic. When an LP is infeasible or unbounded,
the solver can return the dual Farkas or primal ray to prove it.
\end{quote}
\begin{itemize}
\item {} 
\sphinxAtStartPar
\sphinxcode{\sphinxupquote{COPT\_DBLINFO\_DUALFARKAS}} or \sphinxcode{\sphinxupquote{"DualFarkas"}}
\begin{quote}

\sphinxAtStartPar
Double information.

\sphinxAtStartPar
The dual Farkas for constraints of an infeasible LP problem.
Please enable the parameter \sphinxcode{\sphinxupquote{"ReqFarkasRay"}} to ensure that
the dual Farkas is available when the LP is infeasible.

\sphinxAtStartPar
Without loss of generality, the concept of the dual Farkas can be
conveniently demonstrated using an LP problem
with general variable bounds and equality constraints:
\(Ax = 0 \text{ and } l \leq x \leq u\).
When the LP is infeasible, a dual Farkas vector \(y\)
can prove that the system has conflict that \(\max y^TAx < y^T b = 0\).
Computing \(\max y^TAx\): with the vector \(\hat{a} = y^TA\),
choosing variable bound
\(x_i = l_i\) when \(\hat{a}_i < 0\) and
\(x_i = u_i\) when \(\hat{a}_i > 0\)
gives the maximal possible value of \(y^TAx\) for any \(x\) within their bounds.

\sphinxAtStartPar
Some application relies on the alternate conflict \(\min \bar{y}^TAx > \bar{y}^T b = 0\).
This can be achieved by negating the dual Farkas, i.e. \(\bar{y}=-y\) returned by the solver.

\sphinxAtStartPar
In very rare cases, the solver may fail to return a valid dual Farkas.
For example when the LP problem slightly infeasible by tiny amount, which
We recommend to study and to repair the infeasibility using FeasRelax instead.
\end{quote}

\item {} 
\sphinxAtStartPar
\sphinxcode{\sphinxupquote{COPT\_DBLINFO\_PRIMALRAY}} or \sphinxcode{\sphinxupquote{"PrimalRay"}}
\begin{quote}

\sphinxAtStartPar
Double information.

\sphinxAtStartPar
The primal ray for variables of an unbounded LP problem.
Please enable the parameter \sphinxcode{\sphinxupquote{"ReqFarkasRay"}} to ensure that the primal
ray is available when an LP is unbounded.

\sphinxAtStartPar
For a minimization LP problem in the standard form:
\(\min c^T x, Ax = b  \text{ and } x \geq 0\),
a primal ray vector \(r\) satisfies that \(r \geq 0, Ar = 0  \text{ and } c^T r < 0\).
\end{quote}

\end{itemize}

\subsection{Feasibility relaxation information}
\label{\detokenize{capiref:feasibility-relaxation-information}}\begin{itemize}
\item {} 
\sphinxAtStartPar
\sphinxcode{\sphinxupquote{COPT\_DBLINFO\_RELAXLB}} or \sphinxcode{\sphinxupquote{"RelaxLB"}}
\begin{quote}

\sphinxAtStartPar
Double information.

\sphinxAtStartPar
Feasibility relaxation values for lower bounds of columns or rows.
\end{quote}

\item {} 
\sphinxAtStartPar
\sphinxcode{\sphinxupquote{COPT\_DBLINFO\_RELAXUB}} or \sphinxcode{\sphinxupquote{"RelaxUB"}}
\begin{quote}

\sphinxAtStartPar
Double information.

\sphinxAtStartPar
Feasibility relaxation values for upper bounds of columns or rows.
\end{quote}

\item {} 
\sphinxAtStartPar
\sphinxcode{\sphinxupquote{COPT\_DBLINFO\_RELAXVALUE}} or \sphinxcode{\sphinxupquote{"RelaxValue"}}
\begin{quote}

\sphinxAtStartPar
Double information.

\sphinxAtStartPar
Solutions for the original model variables (columns) in the feasibility relaxation model.
\end{quote}

\end{itemize}

\section{Callback information}
\label{\detokenize{capiref:callback-information}}\label{\detokenize{capiref:chapapi-cbcinfo}}\phantomsection\label{\detokenize{capiref:bestobj}}\begin{itemize}
\item {} 
\sphinxAtStartPar
\sphinxcode{\sphinxupquote{COPT\_CBINFO\_BESTOBJ}} or  \sphinxcode{\sphinxupquote{"BestObj"}}
\begin{quote}

\sphinxAtStartPar
Double information.

\sphinxAtStartPar
Current best objective.
\end{quote}

\end{itemize}
\phantomsection\label{\detokenize{capiref:bestbnd}}\begin{itemize}
\item {} 
\sphinxAtStartPar
\sphinxcode{\sphinxupquote{COPT\_CBINFO\_BESTBND}} or  \sphinxcode{\sphinxupquote{"BestBnd"}}
\begin{quote}

\sphinxAtStartPar
Double information.

\sphinxAtStartPar
Current best objective bound.
\end{quote}

\end{itemize}
\phantomsection\label{\detokenize{capiref:hasincumbent}}\begin{itemize}
\item {} 
\sphinxAtStartPar
\sphinxcode{\sphinxupquote{COPT\_CBINFO\_HASINCUMBENT}} or  \sphinxcode{\sphinxupquote{"HasIncumbent"}}
\begin{quote}

\sphinxAtStartPar
Integer information.

\sphinxAtStartPar
Whether an incumbent is available.
\end{quote}

\end{itemize}
\phantomsection\label{\detokenize{capiref:incumbent}}\begin{itemize}
\item {} 
\sphinxAtStartPar
\sphinxcode{\sphinxupquote{COPT\_CBINFO\_INCUMBENT}} or  \sphinxcode{\sphinxupquote{"Incumbent"}}
\begin{quote}

\sphinxAtStartPar
Double information.

\sphinxAtStartPar
Current best feasible solution.
\end{quote}

\end{itemize}
\phantomsection\label{\detokenize{capiref:mipcandidate}}\begin{itemize}
\item {} 
\sphinxAtStartPar
\sphinxcode{\sphinxupquote{COPT\_CBINFO\_MIPCANDIDATE}} or  \sphinxcode{\sphinxupquote{"MipCandidate"}}
\begin{quote}

\sphinxAtStartPar
Double information.

\sphinxAtStartPar
Current feasible solution candidate.
\end{quote}

\end{itemize}
\phantomsection\label{\detokenize{capiref:mipcandobj}}\begin{itemize}
\item {} 
\sphinxAtStartPar
\sphinxcode{\sphinxupquote{COPT\_CBINFO\_MIPCANDOBJ}} or  \sphinxcode{\sphinxupquote{"MipCandObj"}}
\begin{quote}

\sphinxAtStartPar
Double information.

\sphinxAtStartPar
Objective value for current feasible solution candidate.
\end{quote}

\end{itemize}
\phantomsection\label{\detokenize{capiref:relaxsolution}}\begin{itemize}
\item {} 
\sphinxAtStartPar
\sphinxcode{\sphinxupquote{COPT\_CBINFO\_RELAXSOLUTION}} or  \sphinxcode{\sphinxupquote{"RelaxSolution"}}
\begin{quote}

\sphinxAtStartPar
Double information.

\sphinxAtStartPar
Current solution of LP\sphinxhyphen{}relaxation.
\end{quote}

\end{itemize}
\phantomsection\label{\detokenize{capiref:relaxsolobj}}\begin{itemize}
\item {} 
\sphinxAtStartPar
\sphinxcode{\sphinxupquote{COPT\_CBINFO\_RELAXSOLOBJ}} or  \sphinxcode{\sphinxupquote{"RelaxSolObj"}}
\begin{quote}

\sphinxAtStartPar
Double information.

\sphinxAtStartPar
Current objective of LP\sphinxhyphen{}relaxation.
\end{quote}

\end{itemize}
\phantomsection\label{\detokenize{capiref:nodestatus}}\begin{itemize}
\item {} 
\sphinxAtStartPar
\sphinxcode{\sphinxupquote{COPT\_CBINFO\_NODESTATUS}} or  \sphinxcode{\sphinxupquote{"NodeStatus"}}
\begin{quote}

\sphinxAtStartPar
Integer information.

\sphinxAtStartPar
The solution status of the LP\sphinxhyphen{}relaxation problem at the current node.

\sphinxAtStartPar
For possible values, please refer to: {\hyperref[\detokenize{constant:copttab-statuscodes}]{\sphinxcrossref{\DUrole{std,std-ref}{General Constants Chapter: Solution Status (Part)}}}}, except for \sphinxcode{\sphinxupquote{NODELIMIT}}, \sphinxcode{\sphinxupquote{UNSTARTED}}, \sphinxcode{\sphinxupquote{INF\_OR\_UNB}} .
\end{quote}

\end{itemize}

\section{Parameters}
\label{\detokenize{capiref:parameters}}\label{\detokenize{capiref:chapapi-param}}

\subsection{Limits and tolerances}
\label{\detokenize{capiref:limits-and-tolerances}}\begin{itemize}
\item {} 
\sphinxAtStartPar
\sphinxcode{\sphinxupquote{COPT\_DBLPARAM\_TIMELIMIT}} or \sphinxcode{\sphinxupquote{"TimeLimit"}}
\begin{quote}

\sphinxAtStartPar
Double parameter.

\sphinxAtStartPar
Time limit of the optimization.

\sphinxAtStartPar
\sphinxstylestrong{Default:} 1e20

\sphinxAtStartPar
\sphinxstylestrong{Minimal:} 0

\sphinxAtStartPar
\sphinxstylestrong{Maximal:} 1e20
\end{quote}

\item {} 
\sphinxAtStartPar
\sphinxcode{\sphinxupquote{COPT\_DBLPARAM\_SOLTIMELIMIT}} or \sphinxcode{\sphinxupquote{"SolTimeLimit"}}
\begin{quote}

\sphinxAtStartPar
Double parameter.

\sphinxAtStartPar
Time limit if a primal feasible solution has been found.

\sphinxAtStartPar
\sphinxstylestrong{Default:} 1e20

\sphinxAtStartPar
\sphinxstylestrong{Minimal:} 0

\sphinxAtStartPar
\sphinxstylestrong{Maximal:} 1e20
\end{quote}

\item {} 
\sphinxAtStartPar
\sphinxcode{\sphinxupquote{COPT\_INTPARAM\_NODELIMIT}} or \sphinxcode{\sphinxupquote{"NodeLimit"}}
\begin{quote}

\sphinxAtStartPar
Integer parameter.

\sphinxAtStartPar
Node limit of the optimization.

\sphinxAtStartPar
\sphinxstylestrong{Default:} \sphinxhyphen{}1

\sphinxAtStartPar
\sphinxstylestrong{Minimal:} \sphinxhyphen{}1

\sphinxAtStartPar
\sphinxstylestrong{Maximal:} \sphinxcode{\sphinxupquote{INT\_MAX}}
\end{quote}

\item {} 
\sphinxAtStartPar
\sphinxcode{\sphinxupquote{COPT\_INTPARAM\_BARITERLIMIT}} or \sphinxcode{\sphinxupquote{"BarIterLimit"}}
\begin{quote}

\sphinxAtStartPar
Integer parameter.

\sphinxAtStartPar
Iteration limit of barrier method.

\sphinxAtStartPar
\sphinxstylestrong{Default:} 500

\sphinxAtStartPar
\sphinxstylestrong{Minimal:} 0

\sphinxAtStartPar
\sphinxstylestrong{Maximal:} \sphinxcode{\sphinxupquote{INT\_MAX}}
\end{quote}

\item {} 
\sphinxAtStartPar
\sphinxcode{\sphinxupquote{COPT\_INTPARAM\_NLPITERLIMIT}} or \sphinxcode{\sphinxupquote{"NLPIterLimit"}}
\begin{quote}

\sphinxAtStartPar
Integer parameter.

\sphinxAtStartPar
Iteration limit for the nonlinear solver.

\sphinxAtStartPar
\sphinxstylestrong{Default:} 1e4

\sphinxAtStartPar
\sphinxstylestrong{Minimum value}: 0

\sphinxAtStartPar
\sphinxstylestrong{Maximum value}: \sphinxcode{\sphinxupquote{INT\_MAX}}
\end{quote}

\item {} 
\sphinxAtStartPar
\sphinxcode{\sphinxupquote{COPT\_INTPARAM\_MIPNLPITERLIMIT}} or \sphinxcode{\sphinxupquote{"MipNLPIterLimit"}}
\begin{quote}

\sphinxAtStartPar
Integer parameter.

\sphinxAtStartPar
Iteration limit for solving NLP problem(s) within the MIP solver.

\sphinxAtStartPar
\sphinxstylestrong{Default:} 100

\sphinxAtStartPar
\sphinxstylestrong{Minimal:} \sphinxhyphen{}1 (no limit)

\sphinxAtStartPar
\sphinxstylestrong{Maximal:} \sphinxcode{\sphinxupquote{INT\_MAX}}
\end{quote}

\item {} 
\sphinxAtStartPar
\sphinxcode{\sphinxupquote{COPT\_DBLPARAM\_MATRIXTOL}} or \sphinxcode{\sphinxupquote{"MatrixTol"}}
\begin{quote}

\sphinxAtStartPar
Double parameter.

\sphinxAtStartPar
Input matrix coefficient tolerance.

\sphinxAtStartPar
\sphinxstylestrong{Default:} 1e\sphinxhyphen{}10

\sphinxAtStartPar
\sphinxstylestrong{Minimal:} 0

\sphinxAtStartPar
\sphinxstylestrong{Maximal:} 1e\sphinxhyphen{}7
\end{quote}

\item {} 
\sphinxAtStartPar
\sphinxcode{\sphinxupquote{COPT\_DBLPARAM\_FEASTOL}} or \sphinxcode{\sphinxupquote{"FeasTol"}}
\begin{quote}

\sphinxAtStartPar
Double parameter.

\sphinxAtStartPar
The feasibility tolerance.

\sphinxAtStartPar
\sphinxstylestrong{Default:} 1e\sphinxhyphen{}6

\sphinxAtStartPar
\sphinxstylestrong{Minimal:} 1e\sphinxhyphen{}9.

\sphinxAtStartPar
\sphinxstylestrong{Maximal:} 1e\sphinxhyphen{}4
\end{quote}

\item {} 
\sphinxAtStartPar
\sphinxcode{\sphinxupquote{COPT\_DBLPARAM\_DUALTOL}} or \sphinxcode{\sphinxupquote{"DualTol"}}
\begin{quote}

\sphinxAtStartPar
Double parameter.

\sphinxAtStartPar
The tolerance for dual solutions and reduced cost.

\sphinxAtStartPar
\sphinxstylestrong{Default:} 1e\sphinxhyphen{}6

\sphinxAtStartPar
\sphinxstylestrong{Minimal:} 1e\sphinxhyphen{}9

\sphinxAtStartPar
\sphinxstylestrong{Maximal:} 1e\sphinxhyphen{}4
\end{quote}

\item {} 
\sphinxAtStartPar
\sphinxcode{\sphinxupquote{COPT\_DBLPARAM\_INTTOL}} or \sphinxcode{\sphinxupquote{"IntTol"}}
\begin{quote}

\sphinxAtStartPar
Double parameter.

\sphinxAtStartPar
The integrality tolerance for variables.

\sphinxAtStartPar
\sphinxstylestrong{Default:} 1e\sphinxhyphen{}6

\sphinxAtStartPar
\sphinxstylestrong{Minimal:} 1e\sphinxhyphen{}9

\sphinxAtStartPar
\sphinxstylestrong{Maximal:} 1e\sphinxhyphen{}1
\end{quote}

\item {} 
\sphinxAtStartPar
\sphinxcode{\sphinxupquote{COPT\_DBLPARAM\_PDLPTOL}} or \sphinxcode{\sphinxupquote{"PDLPTol"}}
\begin{quote}

\sphinxAtStartPar
Double parameter.

\sphinxAtStartPar
The PDLP tolerance.

\sphinxAtStartPar
\sphinxstylestrong{Default:} 1e\sphinxhyphen{}6

\sphinxAtStartPar
\sphinxstylestrong{Minimal:} 1e\sphinxhyphen{}12

\sphinxAtStartPar
\sphinxstylestrong{Maximal:} 1e\sphinxhyphen{}4
\end{quote}

\item {} 
\sphinxAtStartPar
\sphinxcode{\sphinxupquote{COPT\_DBLPARAM\_NLPTOL}} or \sphinxcode{\sphinxupquote{"NLPTol"}}
\begin{quote}

\sphinxAtStartPar
Double parameter.

\sphinxAtStartPar
The NLP tolerance.

\sphinxAtStartPar
\sphinxstylestrong{Default:} 1e\sphinxhyphen{}8

\sphinxAtStartPar
\sphinxstylestrong{Minimum value}: 1e\sphinxhyphen{}13

\sphinxAtStartPar
\sphinxstylestrong{Maximum value}: 1e\sphinxhyphen{}3
\end{quote}

\item {} 
\sphinxAtStartPar
\sphinxcode{\sphinxupquote{COPT\_DBLPARAM\_RELGAP}} or \sphinxcode{\sphinxupquote{"RelGap"}}
\begin{quote}

\sphinxAtStartPar
Double parameter.

\sphinxAtStartPar
The relative gap of optimization.

\sphinxAtStartPar
\sphinxstylestrong{Default:} 1e\sphinxhyphen{}4

\sphinxAtStartPar
\sphinxstylestrong{Minimal:} 0

\sphinxAtStartPar
\sphinxstylestrong{Maximal:} \sphinxcode{\sphinxupquote{DBL\_MAX}}
\end{quote}

\item {} 
\sphinxAtStartPar
\sphinxcode{\sphinxupquote{COPT\_DBLPARAM\_ABSGAP}} or \sphinxcode{\sphinxupquote{"AbsGap"}}
\begin{quote}

\sphinxAtStartPar
Double parameter.

\sphinxAtStartPar
The absolute gap of optimization.

\sphinxAtStartPar
\sphinxstylestrong{Default:} 1e\sphinxhyphen{}6

\sphinxAtStartPar
\sphinxstylestrong{Minimal:} 0

\sphinxAtStartPar
\sphinxstylestrong{Maximal:} \sphinxcode{\sphinxupquote{DBL\_MAX}}
\end{quote}

\end{itemize}

\subsection{Presolving and scaling}
\label{\detokenize{capiref:presolving-and-scaling}}\begin{itemize}
\item {} 
\sphinxAtStartPar
\sphinxcode{\sphinxupquote{COPT\_INTPARAM\_PRESOLVE}} or \sphinxcode{\sphinxupquote{"Presolve"}}
\begin{quote}

\sphinxAtStartPar
Integer parameter.

\sphinxAtStartPar
Level of presolving before solving a model.

\sphinxAtStartPar
\sphinxstylestrong{Default:} \sphinxhyphen{}1

\sphinxAtStartPar
\sphinxstylestrong{Possible values:}
\begin{quote}

\sphinxAtStartPar
\sphinxhyphen{}1: Automatic

\sphinxAtStartPar
0: Off

\sphinxAtStartPar
1: Fast

\sphinxAtStartPar
2: Normal

\sphinxAtStartPar
3: Aggressive

\sphinxAtStartPar
4: No Limitations, continues until the model cannot be modified (may be very time\sphinxhyphen{}consuming).
\end{quote}
\end{quote}

\item {} 
\sphinxAtStartPar
\sphinxcode{\sphinxupquote{COPT\_INTPARAM\_SCALING}} or \sphinxcode{\sphinxupquote{"Scaling"}}
\begin{quote}

\sphinxAtStartPar
Integer parameter.

\sphinxAtStartPar
Whether to perform scaling before solving a problem.

\sphinxAtStartPar
\sphinxstylestrong{Default:} \sphinxhyphen{}1

\sphinxAtStartPar
\sphinxstylestrong{Possible values:}
\begin{quote}

\sphinxAtStartPar
\sphinxhyphen{}1: Choose automatically.

\sphinxAtStartPar
0: No scaling.

\sphinxAtStartPar
1: Apply scaling.
\end{quote}
\end{quote}

\item {} 
\sphinxAtStartPar
\sphinxcode{\sphinxupquote{COPT\_INTPARAM\_DUALIZE}} or \sphinxcode{\sphinxupquote{"Dualize"}}
\begin{quote}

\sphinxAtStartPar
Integer parameter.

\sphinxAtStartPar
Whether to dualize a problem before solving it.

\sphinxAtStartPar
\sphinxstylestrong{Default:} \sphinxhyphen{}1

\sphinxAtStartPar
\sphinxstylestrong{Possible values:}
\begin{quote}

\sphinxAtStartPar
\sphinxhyphen{}1: Choose automatically.

\sphinxAtStartPar
0: No dualizing.

\sphinxAtStartPar
1: Dualizing the problem.
\end{quote}
\end{quote}

\end{itemize}

\subsection{Linear Programming related}
\label{\detokenize{capiref:linear-programming-related}}\begin{itemize}
\item {} 
\sphinxAtStartPar
\sphinxcode{\sphinxupquote{COPT\_INTPARAM\_LPMETHOD}} or \sphinxcode{\sphinxupquote{"LpMethod"}}
\begin{quote}

\sphinxAtStartPar
Integer parameter.

\sphinxAtStartPar
Method to solve the LP problem.

\sphinxAtStartPar
\sphinxstylestrong{Default:} \sphinxhyphen{}1

\sphinxAtStartPar
\sphinxstylestrong{Possible values:}
\begin{quote}

\sphinxAtStartPar
\sphinxhyphen{}1: Choose automatically.
\begin{quote}

\sphinxAtStartPar
For Linear Programming, choose dual simplex method;

\sphinxAtStartPar
For Mixed Integer Linear Programming, choose dual simplex or barrier method.
\end{quote}

\sphinxAtStartPar
1: Dual simplex.

\sphinxAtStartPar
2: Barrier.

\sphinxAtStartPar
3: Crossover.

\sphinxAtStartPar
4: Concurrent (Use multiple algorithms simultaneously).

\sphinxAtStartPar
5: Choose between simplex and barrier automatically (Based on features such as sparsity and/or coefficients ranges).

\sphinxAtStartPar
6: First\sphinxhyphen{}order method (PDLP).
\end{quote}
\end{quote}

\end{itemize}

\begin{sphinxadmonition}{note}{Note:}

\sphinxAtStartPar
Currently, COPT’s GPU mode only supports solving Linear Programming problems using the first\sphinxhyphen{}order method (PDLP).
To enable it, you need to set \sphinxcode{\sphinxupquote{LpMethod=6}} first.
\end{sphinxadmonition}
\begin{itemize}
\item {} 
\sphinxAtStartPar
\sphinxcode{\sphinxupquote{COPT\_INTPARAM\_DUALPRICE}} or \sphinxcode{\sphinxupquote{"DualPrice"}}
\begin{quote}

\sphinxAtStartPar
Integer parameter.

\sphinxAtStartPar
Specifies the dual simplex pricing algorithm.

\sphinxAtStartPar
\sphinxstylestrong{Default:} \sphinxhyphen{}1

\sphinxAtStartPar
\sphinxstylestrong{Possible values:}
\begin{quote}

\sphinxAtStartPar
\sphinxhyphen{}1: Choose automatically.

\sphinxAtStartPar
0: Using Devex pricing algorithm.

\sphinxAtStartPar
1: Using dual steepest\sphinxhyphen{}edge pricing algorithm.
\end{quote}
\end{quote}

\item {} 
\sphinxAtStartPar
\sphinxcode{\sphinxupquote{COPT\_INTPARAM\_DUALPERTURB}} or \sphinxcode{\sphinxupquote{"DualPerturb"}}
\begin{quote}

\sphinxAtStartPar
Integer parameter.

\sphinxAtStartPar
Whether to allow the objective function perturbation
when using the dual simplex method.

\sphinxAtStartPar
\sphinxstylestrong{Default:} \sphinxhyphen{}1

\sphinxAtStartPar
\sphinxstylestrong{Possible values:}
\begin{quote}

\sphinxAtStartPar
\sphinxhyphen{}1: Choose automatically.

\sphinxAtStartPar
0: No perturbation.

\sphinxAtStartPar
1: Allow objective function perturbation.
\end{quote}
\end{quote}

\item {} 
\sphinxAtStartPar
\sphinxcode{\sphinxupquote{COPT\_INTPARAM\_BARHOMOGENEOUS}} or \sphinxcode{\sphinxupquote{"BarHomogeneous"}}
\begin{quote}

\sphinxAtStartPar
Integer parameter.

\sphinxAtStartPar
Whether to use homogeneous self\sphinxhyphen{}dual form in barrier.

\sphinxAtStartPar
\sphinxstylestrong{Default:} \sphinxhyphen{}1

\sphinxAtStartPar
\sphinxstylestrong{Possible values:}
\begin{quote}

\sphinxAtStartPar
\sphinxhyphen{}1: Choose automatically.

\sphinxAtStartPar
0: No.

\sphinxAtStartPar
1: Yes.
\end{quote}
\end{quote}

\item {} 
\sphinxAtStartPar
\sphinxcode{\sphinxupquote{COPT\_INTPARAM\_BARORDER}} or \sphinxcode{\sphinxupquote{"BarOrder"}}
\begin{quote}

\sphinxAtStartPar
Integer parameter.

\sphinxAtStartPar
Barrier ordering algorithm.

\sphinxAtStartPar
\sphinxstylestrong{Default:} \sphinxhyphen{}1

\sphinxAtStartPar
\sphinxstylestrong{Possible values:}
\begin{quote}

\sphinxAtStartPar
\sphinxhyphen{}1: Choose automatically.

\sphinxAtStartPar
0: Approximate Minimum Degree (AMD).

\sphinxAtStartPar
1: Nested Dissection (ND1).

\sphinxAtStartPar
2: Modified Nested Dissection (ND2).
\end{quote}
\end{quote}

\item {} 
\sphinxAtStartPar
\sphinxcode{\sphinxupquote{COPT\_INTPARAM\_BARSTART}} or \sphinxcode{\sphinxupquote{"BarStart"}}
\begin{quote}

\sphinxAtStartPar
Integer parameter.

\sphinxAtStartPar
Algorithm for finding initial points in barrier method.

\sphinxAtStartPar
\sphinxstylestrong{Default:} \sphinxhyphen{}1

\sphinxAtStartPar
\sphinxstylestrong{Possible values:}
\begin{quote}

\sphinxAtStartPar
\sphinxhyphen{}1: Choose automatically.

\sphinxAtStartPar
0: Simple.

\sphinxAtStartPar
1: Mehrotra.

\sphinxAtStartPar
2: Modified Mehrotra.
\end{quote}
\end{quote}

\item {} 
\sphinxAtStartPar
\sphinxcode{\sphinxupquote{COPT\_INTPARAM\_CROSSOVER}} or \sphinxcode{\sphinxupquote{"Crossover"}}
\begin{quote}

\sphinxAtStartPar
Integer parameter.

\sphinxAtStartPar
Whether to use crossover.

\sphinxAtStartPar
\sphinxstylestrong{Default:} 1

\sphinxAtStartPar
\sphinxstylestrong{Possible values:}
\begin{quote}

\sphinxAtStartPar
\sphinxhyphen{}1: Choose automatically.
\begin{quote}

\sphinxAtStartPar
Only run crossover when the LP solution is not primal\sphinxhyphen{}dual feasible.
\end{quote}

\sphinxAtStartPar
0: No.

\sphinxAtStartPar
1: Yes.
\end{quote}
\end{quote}

\item {} 
\sphinxAtStartPar
\sphinxcode{\sphinxupquote{COPT\_INTPARAM\_REQFARKASRAY}} or \sphinxcode{\sphinxupquote{"ReqFarkasRay"}}
\begin{quote}

\sphinxAtStartPar
Integer parameter.

\sphinxAtStartPar
Advanced topic. Whether to compute the dual Farkas or primal ray when the LP is infeasible or unbounded.

\sphinxAtStartPar
\sphinxstylestrong{Default:} 0

\sphinxAtStartPar
\sphinxstylestrong{Possible values:}
\begin{quote}

\sphinxAtStartPar
0: No.

\sphinxAtStartPar
1: Yes.
\end{quote}
\end{quote}

\end{itemize}
\phantomsection\label{\detokenize{capiref:reqsensitivity}}\begin{itemize}
\item {} 
\sphinxAtStartPar
\sphinxcode{\sphinxupquote{ReqSensitivity}}
\begin{quote}

\sphinxAtStartPar
Integer parameter.

\sphinxAtStartPar
Whether to compute sensitivity analysis when an optimal basis is available for an LP problem
(when solved by the simplex method or by other methods followed by crossover).

\sphinxAtStartPar
\sphinxstylestrong{Default:} 0

\sphinxAtStartPar
\sphinxstylestrong{Possible values:}
\begin{quote}

\sphinxAtStartPar
0: No.

\sphinxAtStartPar
1: Yes.
\end{quote}
\end{quote}

\end{itemize}

\subsection{Semidefinite Programming related}
\label{\detokenize{capiref:semidefinite-programming-related}}\begin{itemize}
\item {} 
\sphinxAtStartPar
\sphinxcode{\sphinxupquote{COPT\_INTPARAM\_SDPMETHOD}} or \sphinxcode{\sphinxupquote{"SDPMethod"}}
\begin{quote}

\sphinxAtStartPar
Integer parameter.

\sphinxAtStartPar
Method to solve the SDP problem.

\sphinxAtStartPar
\sphinxstylestrong{Default:} \sphinxhyphen{}1

\sphinxAtStartPar
\sphinxstylestrong{Possible values:}
\begin{quote}

\sphinxAtStartPar
\sphinxhyphen{}1: Choose automatically.

\sphinxAtStartPar
0: Primal\sphinxhyphen{}Dual method.

\sphinxAtStartPar
1: Alternating direction method of multipliers (ADMM).

\sphinxAtStartPar
2: Dual method.
\end{quote}
\end{quote}

\end{itemize}

\subsection{Integer Programming related}
\label{\detokenize{capiref:integer-programming-related}}\begin{itemize}
\item {} 
\sphinxAtStartPar
\sphinxcode{\sphinxupquote{COPT\_INTPARAM\_CUTLEVEL}} or \sphinxcode{\sphinxupquote{"CutLevel"}}
\begin{quote}

\sphinxAtStartPar
Integer parameter.

\sphinxAtStartPar
Level of cutting\sphinxhyphen{}planes generation.

\sphinxAtStartPar
\sphinxstylestrong{Default:} \sphinxhyphen{}1

\sphinxAtStartPar
\sphinxstylestrong{Possible values:}
\begin{quote}

\sphinxAtStartPar
\sphinxhyphen{}1: Choose automatically.

\sphinxAtStartPar
0: Off

\sphinxAtStartPar
1: Fast

\sphinxAtStartPar
2: Normal

\sphinxAtStartPar
3: Aggressive
\end{quote}
\end{quote}

\item {} 
\sphinxAtStartPar
\sphinxcode{\sphinxupquote{COPT\_INTPARAM\_ROOTCUTLEVEL}} or \sphinxcode{\sphinxupquote{"RootCutLevel"}}
\begin{quote}

\sphinxAtStartPar
Integer parameter.

\sphinxAtStartPar
Level of cutting\sphinxhyphen{}planes generation of root node.

\sphinxAtStartPar
\sphinxstylestrong{Default:} \sphinxhyphen{}1

\sphinxAtStartPar
\sphinxstylestrong{Possible values:}
\begin{quote}

\sphinxAtStartPar
\sphinxhyphen{}1: Choose automatically.

\sphinxAtStartPar
0: Off

\sphinxAtStartPar
1: Fast

\sphinxAtStartPar
2: Normal

\sphinxAtStartPar
3: Aggressive
\end{quote}
\end{quote}

\item {} 
\sphinxAtStartPar
\sphinxcode{\sphinxupquote{COPT\_INTPARAM\_TREECUTLEVEL}} or \sphinxcode{\sphinxupquote{"TreeCutLevel"}}
\begin{quote}

\sphinxAtStartPar
Integer parameter.

\sphinxAtStartPar
Level of cutting\sphinxhyphen{}planes generation of search tree.

\sphinxAtStartPar
\sphinxstylestrong{Default:} \sphinxhyphen{}1

\sphinxAtStartPar
\sphinxstylestrong{Possible values:}
\begin{quote}

\sphinxAtStartPar
\sphinxhyphen{}1: Choose automatically.

\sphinxAtStartPar
0: Off

\sphinxAtStartPar
1: Fast

\sphinxAtStartPar
2: Normal

\sphinxAtStartPar
3: Aggressive
\end{quote}
\end{quote}

\item {} 
\sphinxAtStartPar
\sphinxcode{\sphinxupquote{COPT\_INTPARAM\_ROOTCUTROUNDS}} or \sphinxcode{\sphinxupquote{"RootCutRounds"}}
\begin{quote}

\sphinxAtStartPar
Integer parameter.

\sphinxAtStartPar
Rounds of cutting\sphinxhyphen{}planes generation of root node.

\sphinxAtStartPar
\sphinxstylestrong{Default:} \sphinxhyphen{}1 (Choose automatically)

\sphinxAtStartPar
\sphinxstylestrong{Minimal:} \sphinxhyphen{}1

\sphinxAtStartPar
\sphinxstylestrong{Maximal:} \sphinxcode{\sphinxupquote{INT\_MAX}}
\end{quote}

\item {} 
\sphinxAtStartPar
\sphinxcode{\sphinxupquote{COPT\_INTPARAM\_NODECUTROUNDS}} or \sphinxcode{\sphinxupquote{"NodeCutRounds"}}
\begin{quote}

\sphinxAtStartPar
Integer parameter.

\sphinxAtStartPar
Rounds of cutting\sphinxhyphen{}planes generation of search tree node.

\sphinxAtStartPar
\sphinxstylestrong{Default:} \sphinxhyphen{}1 (Choose automatically)

\sphinxAtStartPar
\sphinxstylestrong{Minimal:} \sphinxhyphen{}1

\sphinxAtStartPar
\sphinxstylestrong{Maximal:} \sphinxcode{\sphinxupquote{INT\_MAX}}
\end{quote}

\item {} 
\sphinxAtStartPar
\sphinxcode{\sphinxupquote{COPT\_INTPARAM\_HEURLEVEL}} or \sphinxcode{\sphinxupquote{"HeurLevel"}}
\begin{quote}

\sphinxAtStartPar
Integer parameter.

\sphinxAtStartPar
Level of heuristics.

\sphinxAtStartPar
\sphinxstylestrong{Default:} \sphinxhyphen{}1

\sphinxAtStartPar
\sphinxstylestrong{Possible values:}
\begin{quote}

\sphinxAtStartPar
\sphinxhyphen{}1: Choose automatically

\sphinxAtStartPar
0: Off

\sphinxAtStartPar
1: Fast

\sphinxAtStartPar
2: Normal

\sphinxAtStartPar
3: Aggressive
\end{quote}
\end{quote}

\end{itemize}
\phantomsection\label{\detokenize{capiref:prerootheurlevel}}\begin{itemize}
\item {} 
\sphinxAtStartPar
\sphinxcode{\sphinxupquote{COPT\_INTPARAM\_PREROOTHEURLEVEL}} or \sphinxcode{\sphinxupquote{"PreRootHeurLevel"}}
\begin{quote}

\sphinxAtStartPar
Integer parameter.

\sphinxAtStartPar
Level of pre\sphinxhyphen{}root heuristics.

\sphinxAtStartPar
\sphinxstylestrong{Default:} \sphinxhyphen{}1

\sphinxAtStartPar
\sphinxstylestrong{Possible value:}
\begin{quote}

\sphinxAtStartPar
\sphinxhyphen{}1: Choose automatically

\sphinxAtStartPar
0: Off

\sphinxAtStartPar
1: Fast

\sphinxAtStartPar
2: Normal

\sphinxAtStartPar
3: Aggressive
\end{quote}
\end{quote}

\item {} 
\sphinxAtStartPar
\sphinxcode{\sphinxupquote{COPT\_INTPARAM\_ROUNDINGHEURLEVEL}} or \sphinxcode{\sphinxupquote{"RoundingHeurLevel"}}
\begin{quote}

\sphinxAtStartPar
Integer parameter.

\sphinxAtStartPar
Level of rounding heuristics.

\sphinxAtStartPar
\sphinxstylestrong{Default:} \sphinxhyphen{}1

\sphinxAtStartPar
\sphinxstylestrong{Possible values:}
\begin{quote}

\sphinxAtStartPar
\sphinxhyphen{}1: Choose automatically

\sphinxAtStartPar
0: Off

\sphinxAtStartPar
1: Fast

\sphinxAtStartPar
2: Normal

\sphinxAtStartPar
3: Aggressive
\end{quote}
\end{quote}

\item {} 
\sphinxAtStartPar
\sphinxcode{\sphinxupquote{COPT\_INTPARAM\_DIVINGHEURLEVEL}} or \sphinxcode{\sphinxupquote{"DivingHeurLevel"}}
\begin{quote}

\sphinxAtStartPar
Integer parameter.

\sphinxAtStartPar
Level of diving heuristics.

\sphinxAtStartPar
\sphinxstylestrong{Default:} \sphinxhyphen{}1

\sphinxAtStartPar
\sphinxstylestrong{Possible values:}
\begin{quote}

\sphinxAtStartPar
\sphinxhyphen{}1: Choose automatically

\sphinxAtStartPar
0: Off

\sphinxAtStartPar
1: Fast

\sphinxAtStartPar
2: Normal

\sphinxAtStartPar
3: Aggressive
\end{quote}
\end{quote}

\item {} 
\sphinxAtStartPar
\sphinxcode{\sphinxupquote{COPT\_INTPARAM\_TREECUTLEVEL}} or \sphinxcode{\sphinxupquote{"SubMipHeurLevel"}}
\begin{quote}

\sphinxAtStartPar
Integer parameter.

\sphinxAtStartPar
Level of Sub\sphinxhyphen{}MIP heuristics.

\sphinxAtStartPar
\sphinxstylestrong{Default:} \sphinxhyphen{}1

\sphinxAtStartPar
\sphinxstylestrong{Possible values:}
\begin{quote}

\sphinxAtStartPar
\sphinxhyphen{}1: Choose automatically

\sphinxAtStartPar
0: Off

\sphinxAtStartPar
1: Fast

\sphinxAtStartPar
2: Normal

\sphinxAtStartPar
3: Aggressive
\end{quote}
\end{quote}

\item {} 
\sphinxAtStartPar
\sphinxcode{\sphinxupquote{COPT\_INTPARAM\_FAPHEURLEVEL}} or \sphinxcode{\sphinxupquote{"FAPHeurLevel"}}
\begin{quote}

\sphinxAtStartPar
Integer parameter.

\sphinxAtStartPar
Level of Fix\sphinxhyphen{}and\sphinxhyphen{}propagate heuristics.

\sphinxAtStartPar
\sphinxstylestrong{Default:} \sphinxhyphen{}1

\sphinxAtStartPar
\sphinxstylestrong{Possible values:}
\begin{quote}

\sphinxAtStartPar
\sphinxhyphen{}1: Choose automatically

\sphinxAtStartPar
0: Off

\sphinxAtStartPar
1: Fast

\sphinxAtStartPar
2: Normal

\sphinxAtStartPar
3: Aggressive

\sphinxAtStartPar
4: Heavy
\end{quote}
\end{quote}

\item {} 
\sphinxAtStartPar
\sphinxcode{\sphinxupquote{COPT\_INTPARAM\_STRONGBRANCHING}} or \sphinxcode{\sphinxupquote{"StrongBranching"}}
\begin{quote}

\sphinxAtStartPar
Integer parameter.

\sphinxAtStartPar
Level of strong branching.

\sphinxAtStartPar
\sphinxstylestrong{Default:} \sphinxhyphen{}1

\sphinxAtStartPar
\sphinxstylestrong{Possible values:}
\begin{quote}

\sphinxAtStartPar
\sphinxhyphen{}1: Choose automatically

\sphinxAtStartPar
0: Off

\sphinxAtStartPar
1: Fast

\sphinxAtStartPar
2: Normal

\sphinxAtStartPar
3: Aggressive
\end{quote}
\end{quote}

\item {} 
\sphinxAtStartPar
\sphinxcode{\sphinxupquote{COPT\_INTPARAM\_CONFLICTANALYSIS}} or \sphinxcode{\sphinxupquote{"ConflictAnalysis"}}
\begin{quote}

\sphinxAtStartPar
Integer parameter.

\sphinxAtStartPar
Whether to perform conflict analysis.

\sphinxAtStartPar
\sphinxstylestrong{Default:} \sphinxhyphen{}1

\sphinxAtStartPar
\sphinxstylestrong{Possible values:}
\begin{quote}

\sphinxAtStartPar
\sphinxhyphen{}1: Choose automatically.

\sphinxAtStartPar
0: No

\sphinxAtStartPar
1: Yes
\end{quote}
\end{quote}

\item {} 
\sphinxAtStartPar
\sphinxcode{\sphinxupquote{COPT\_INTPARAM\_MIPSTARTMODE}} or \sphinxcode{\sphinxupquote{"MipStartMode"}}
\begin{quote}

\sphinxAtStartPar
Integer parameter.

\sphinxAtStartPar
Mode of MIP starts.

\sphinxAtStartPar
\sphinxstylestrong{Default:} \sphinxhyphen{}1

\sphinxAtStartPar
\sphinxstylestrong{Possible values:}
\begin{quote}

\sphinxAtStartPar
\sphinxhyphen{}1: Choose automatically.

\sphinxAtStartPar
0: Do not use any MIP starts.

\sphinxAtStartPar
1: Only load full and feasible MIP starts.

\sphinxAtStartPar
2: Only load feasible ones (complete partial solutions by solving subMIPs).
\end{quote}
\end{quote}

\item {} 
\sphinxAtStartPar
\sphinxcode{\sphinxupquote{COPT\_INTPARAM\_MIPSTARTNODELIMIT}} or \sphinxcode{\sphinxupquote{"MipStartNodeLimit"}}
\begin{quote}

\sphinxAtStartPar
Integer parameter.

\sphinxAtStartPar
Limit of nodes for MIP start sub\sphinxhyphen{}MIPs.

\sphinxAtStartPar
\sphinxstylestrong{Default:} \sphinxhyphen{}1 (Choose automatically)

\sphinxAtStartPar
\sphinxstylestrong{Minimal:} \sphinxhyphen{}1 (Choose automatically)

\sphinxAtStartPar
\sphinxstylestrong{Maximal:} \sphinxcode{\sphinxupquote{INT\_MAX}}
\end{quote}

\item {} 
\sphinxAtStartPar
\sphinxcode{\sphinxupquote{COPT\_INTPARAM\_LINEARIZEINDICATORS}} or \sphinxcode{\sphinxupquote{"LinearizeIndicators"}}
\begin{quote}

\sphinxAtStartPar
Integer parameter.

\sphinxAtStartPar
Controls whether to force the linearization of Indicator constraints.

\sphinxAtStartPar
\sphinxstylestrong{Default:} 0

\sphinxAtStartPar
\sphinxstylestrong{Possible values}
\begin{quote}

\sphinxAtStartPar
0: Do not force linearization of indicator constraints.

\sphinxAtStartPar
1: Force linearization of all indicator constraints.

\sphinxAtStartPar
(All indicator constraints will be transformed into equivalent linear constraints if enabled.)
\end{quote}
\end{quote}

\item {} 
\sphinxAtStartPar
\sphinxcode{\sphinxupquote{COPT\_INTPARAM\_LINEARIZESOS}} or \sphinxcode{\sphinxupquote{"LinearizeSos"}}
\begin{quote}

\sphinxAtStartPar
Integer parameter.

\sphinxAtStartPar
Controls whether to force the linearization of SOS constraints.

\sphinxAtStartPar
\sphinxstylestrong{Default:} 0

\sphinxAtStartPar
\sphinxstylestrong{Possible values}
\begin{quote}

\sphinxAtStartPar
0: Do not force linearization of SOS constraints.

\sphinxAtStartPar
1: Force linearization of all SOS constraints.

\sphinxAtStartPar
(All SOS constraints will be transformed into equivalent linear constraints if enabled.)
\end{quote}
\end{quote}

\item {} 
\sphinxAtStartPar
\sphinxcode{\sphinxupquote{COPT\_INTPARAM\_MIPREPAIR}} or \sphinxcode{\sphinxupquote{"MipRepair"}}
\begin{quote}

\sphinxAtStartPar
Integer parameter.

\sphinxAtStartPar
Level for repairing the MIP solution in case of numerical issues.

\sphinxAtStartPar
\sphinxstylestrong{Default:} \sphinxhyphen{}1

\sphinxAtStartPar
\sphinxstylestrong{Possible values:}
\begin{quote}

\sphinxAtStartPar
\sphinxhyphen{}1: Only when time left.

\sphinxAtStartPar
0: Off.

\sphinxAtStartPar
1: Extend time limit and attempt repair (Fast).

\sphinxAtStartPar
2: Extend time limit and attempt repair (Normal).

\sphinxAtStartPar
3: Extend time limit and attempt repair (Aggressive).
\end{quote}
\end{quote}

\end{itemize}

\subsection{Nonlinear Programming related}
\label{\detokenize{capiref:nonlinear-programming-related}}\begin{itemize}
\item {} 
\sphinxAtStartPar
\sphinxcode{\sphinxupquote{COPT\_INTPARAM\_NONCONVEX}} or \sphinxcode{\sphinxupquote{"NonConvex"}}
\begin{quote}

\sphinxAtStartPar
Integer parameter.

\sphinxAtStartPar
Handling strategy for nonconvex models.

\sphinxAtStartPar
\sphinxstylestrong{Default:} \sphinxhyphen{}1

\sphinxAtStartPar
\sphinxstylestrong{Possible values}:
\begin{quote}

\sphinxAtStartPar
\sphinxhyphen{}1: Choose automatically.

\sphinxAtStartPar
0: Report model as nonconvex and terminate.

\sphinxAtStartPar
1: Search for a local optimal solution.
\end{quote}
\end{quote}

\item {} 
\sphinxAtStartPar
\sphinxcode{\sphinxupquote{COPT\_INTPARAM\_NLPMUUPDATE}} or \sphinxcode{\sphinxupquote{"NLPMuUpdate"}}
\begin{quote}

\sphinxAtStartPar
Integer parameter.

\sphinxAtStartPar
Barrier parameter update strategy for the nonlinear solver.

\sphinxAtStartPar
\sphinxstylestrong{Default:} \sphinxhyphen{}1

\sphinxAtStartPar
\sphinxstylestrong{Possible values}:
\begin{quote}

\sphinxAtStartPar
\sphinxhyphen{}1: Choose automatically.

\sphinxAtStartPar
0: Monotonic decrease.

\sphinxAtStartPar
1: Adaptive adjustment.
\end{quote}
\end{quote}

\item {} 
\sphinxAtStartPar
\sphinxcode{\sphinxupquote{COPT\_INTPARAM\_NLPLINSCALE}} or \sphinxcode{\sphinxupquote{"NLPLinScale"}}
\begin{quote}

\sphinxAtStartPar
Integer parameter.

\sphinxAtStartPar
Scaling strategy for the linear system in the nonlinear solver.

\sphinxAtStartPar
\sphinxstylestrong{Default:} \sphinxhyphen{}1

\sphinxAtStartPar
\sphinxstylestrong{Possible values}:
\begin{quote}

\sphinxAtStartPar
\sphinxhyphen{}1: Choose automatically.

\sphinxAtStartPar
0: No scaling.

\sphinxAtStartPar
1: Always scale.
\end{quote}
\end{quote}

\end{itemize}

\subsection{Parallel computing related}
\label{\detokenize{capiref:parallel-computing-related}}\begin{itemize}
\item {} 
\sphinxAtStartPar
\sphinxcode{\sphinxupquote{COPT\_INTPARAM\_THREADS}} or \sphinxcode{\sphinxupquote{"Threads"}}
\begin{quote}

\sphinxAtStartPar
Integer parameter.

\sphinxAtStartPar
Number of threads to use.

\sphinxAtStartPar
\sphinxstylestrong{Default:} \sphinxhyphen{}1 (Choose automatically)

\sphinxAtStartPar
\sphinxstylestrong{Minimal:} \sphinxhyphen{}1 (Choose automatically)

\sphinxAtStartPar
\sphinxstylestrong{Maximal:} 128
\end{quote}

\item {} 
\sphinxAtStartPar
\sphinxcode{\sphinxupquote{COPT\_INTPARAM\_BARTHREADS}} or \sphinxcode{\sphinxupquote{"BarThreads"}}
\begin{quote}

\sphinxAtStartPar
Integer parameter.

\sphinxAtStartPar
Number of threads used by barrier. If value is \sphinxhyphen{}1, the thread count is determined
by parameter \sphinxcode{\sphinxupquote{Threads}}.

\sphinxAtStartPar
\sphinxstylestrong{Default:} \sphinxhyphen{}1 (Choose automatically)

\sphinxAtStartPar
\sphinxstylestrong{Minimal:} \sphinxhyphen{}1 (Choose automatically)

\sphinxAtStartPar
\sphinxstylestrong{Maximal:} 128
\end{quote}

\item {} 
\sphinxAtStartPar
\sphinxcode{\sphinxupquote{COPT\_INTPARAM\_SIMPLEXTHREADS}} or \sphinxcode{\sphinxupquote{"SimplexThreads"}}
\begin{quote}

\sphinxAtStartPar
Integer parameter.

\sphinxAtStartPar
Number of threads used by dual simplex. If value is \sphinxhyphen{}1, the thread count is determined
by parameter \sphinxcode{\sphinxupquote{Threads}}.

\sphinxAtStartPar
\sphinxstylestrong{Default:} \sphinxhyphen{}1 (Choose automatically)

\sphinxAtStartPar
\sphinxstylestrong{Minimal:} \sphinxhyphen{}1 (Choose automatically)

\sphinxAtStartPar
\sphinxstylestrong{Maximal:} 128
\end{quote}

\item {} 
\sphinxAtStartPar
\sphinxcode{\sphinxupquote{COPT\_INTPARAM\_CROSSOVERTHREADS}} or \sphinxcode{\sphinxupquote{"CrossoverThreads"}}
\begin{quote}

\sphinxAtStartPar
Integer parameter.

\sphinxAtStartPar
Number of threads used by crossover. If value is \sphinxhyphen{}1, the thread count is determined
by parameter \sphinxcode{\sphinxupquote{Threads}}.

\sphinxAtStartPar
\sphinxstylestrong{Default:} \sphinxhyphen{}1 (Choose automatically)

\sphinxAtStartPar
\sphinxstylestrong{Minimal:} \sphinxhyphen{}1 (Choose automatically)

\sphinxAtStartPar
\sphinxstylestrong{Maximal:} 128
\end{quote}

\item {} 
\sphinxAtStartPar
\sphinxcode{\sphinxupquote{COPT\_INTPARAM\_MIPTASKS}} or \sphinxcode{\sphinxupquote{"MipTasks"}}
\begin{quote}

\sphinxAtStartPar
Integer parameter.

\sphinxAtStartPar
Number of MIP tasks in parallel.

\sphinxAtStartPar
\sphinxstylestrong{Default:} \sphinxhyphen{}1 (Choose automatically)

\sphinxAtStartPar
\sphinxstylestrong{Minimal:} \sphinxhyphen{}1 (Choose automatically)

\sphinxAtStartPar
\sphinxstylestrong{Maximal:} 256
\end{quote}

\item {} 
\sphinxAtStartPar
\sphinxcode{\sphinxupquote{COPT\_INTPARAM\_ConcurrentLpMode}} or \sphinxcode{\sphinxupquote{"ConcurrentLpMode"}}
\begin{quote}

\sphinxAtStartPar
Integer parameter.

\sphinxAtStartPar
The LP concurrent solving mode.
Only effective when \sphinxcode{\sphinxupquote{LpMethod = 4}}.

\sphinxAtStartPar
When \sphinxcode{\sphinxupquote{LpMethod = 4}} is enabled, the parameters \sphinxcode{\sphinxupquote{GPUMode}} and
\sphinxcode{\sphinxupquote{GPUDevice}} are ignored. GPU usage and device selection are
fully controlled by \sphinxcode{\sphinxupquote{ConcurrentLpMode}}.

\sphinxAtStartPar
\sphinxstylestrong{Default:} \sphinxhyphen{}1 (Choose automatically)

\sphinxAtStartPar
\sphinxstylestrong{Possible values:}
\begin{quote}

\sphinxAtStartPar
\sphinxhyphen{}1: Choose automatically.

\sphinxAtStartPar
0: CPU only (simplex + barrier).

\sphinxAtStartPar
1: CPU (simplex + barrier) + GPU first\sphinxhyphen{}order method PDLP
(uses GPU device 0).

\sphinxAtStartPar
2: CPU (simplex + barrier) + GPU barrier
(uses GPU device 0).

\sphinxAtStartPar
3: CPU (simplex + barrier) + GPU first\sphinxhyphen{}order method PDLP + GPU barrier
(uses GPU device 0 and GPU device 1).
\end{quote}
\end{quote}

\end{itemize}

\begin{sphinxadmonition}{note}{Note:}

\sphinxAtStartPar
If \sphinxcode{\sphinxupquote{ConcurrentLpMode}} is set to 1 or 2 and no GPU device is detected,
or if it is set to 3 and fewer than two GPU devices are detected,
the solver reports an error and terminates.
\end{sphinxadmonition}

\subsection{IIS computation related}
\label{\detokenize{capiref:iis-computation-related}}\begin{itemize}
\item {} 
\sphinxAtStartPar
\sphinxcode{\sphinxupquote{COPT\_INTPARAM\_IISMETHOD}} or \sphinxcode{\sphinxupquote{"IISMethod"}}
\begin{quote}

\sphinxAtStartPar
Integer parameter.

\sphinxAtStartPar
Method for IIS computation.

\sphinxAtStartPar
\sphinxstylestrong{Default:} \sphinxhyphen{}1

\sphinxAtStartPar
\sphinxstylestrong{Possible values:}
\begin{quote}

\sphinxAtStartPar
\sphinxhyphen{}1: Choose automatically.

\sphinxAtStartPar
0: Find smaller IIS.

\sphinxAtStartPar
1: Find IIS quickly.
\end{quote}
\end{quote}

\end{itemize}

\subsection{Feasibility relaxation related}
\label{\detokenize{capiref:feasibility-relaxation-related}}\begin{itemize}
\item {} 
\sphinxAtStartPar
\sphinxcode{\sphinxupquote{COPT\_INTPARAM\_FEASRELAXMODE}} or \sphinxcode{\sphinxupquote{"FeasRelaxMode"}}
\begin{quote}

\sphinxAtStartPar
Integer parameter.

\sphinxAtStartPar
Method for feasibility relaxation.

\sphinxAtStartPar
\sphinxstylestrong{Default:} 0

\sphinxAtStartPar
\sphinxstylestrong{Possible values:}
\begin{quote}

\sphinxAtStartPar
0: Minimize sum of violations.

\sphinxAtStartPar
1: Optimize original objective function under minimal sum of violations.

\sphinxAtStartPar
2: Minimize number of violations.

\sphinxAtStartPar
3: Optimize original objective function under minimal number of violations.

\sphinxAtStartPar
4: Minimize sum of squared violations.

\sphinxAtStartPar
5: Optimize original objective function under minimal sum of squared violations.
\end{quote}
\end{quote}

\end{itemize}

\subsection{Tuner related}
\label{\detokenize{capiref:tuner-related}}\phantomsection\label{\detokenize{capiref:tunetimelimit}}\begin{itemize}
\item {} 
\sphinxAtStartPar
\sphinxcode{\sphinxupquote{COPT\_DBLPARAM\_TUNETIMELIMIT}} or \sphinxcode{\sphinxupquote{"TuneTimeLimit"}}
\begin{quote}

\sphinxAtStartPar
Double parameter.

\sphinxAtStartPar
Time limit for parameter tuning. If the parameter value is 0, it will automatically set by the solver.

\sphinxAtStartPar
\sphinxstylestrong{Default:} 0

\sphinxAtStartPar
\sphinxstylestrong{Minimal:} 0

\sphinxAtStartPar
\sphinxstylestrong{Maximal:} 1e20
\end{quote}

\end{itemize}
\phantomsection\label{\detokenize{capiref:tunetargettime}}\begin{itemize}
\item {} 
\sphinxAtStartPar
\sphinxcode{\sphinxupquote{COPT\_DBLPARAM\_TUNETARGETTIME}} or \sphinxcode{\sphinxupquote{"TuneTargetTime"}}
\begin{quote}

\sphinxAtStartPar
Double parameter.

\sphinxAtStartPar
Time target for parameter tuning.

\sphinxAtStartPar
\sphinxstylestrong{Default:} 1e\sphinxhyphen{}2

\sphinxAtStartPar
\sphinxstylestrong{Minimal:} 0

\sphinxAtStartPar
\sphinxstylestrong{Maximal:} \sphinxcode{\sphinxupquote{DBL\_MAX}}
\end{quote}

\end{itemize}
\phantomsection\label{\detokenize{capiref:tunetargetrelgap}}\begin{itemize}
\item {} 
\sphinxAtStartPar
\sphinxcode{\sphinxupquote{COPT\_DBLPARAM\_TUNETARGETRELGAP}} or \sphinxcode{\sphinxupquote{"TuneTargetRelGap"}}
\begin{quote}

\sphinxAtStartPar
Double parameter.

\sphinxAtStartPar
Optimal relative tolerance target for parameter tuning.

\sphinxAtStartPar
\sphinxstylestrong{Default:} 1e\sphinxhyphen{}4

\sphinxAtStartPar
\sphinxstylestrong{Minimal:} 0

\sphinxAtStartPar
\sphinxstylestrong{Maximal:} \sphinxcode{\sphinxupquote{DBL\_MAX}}
\end{quote}

\end{itemize}
\phantomsection\label{\detokenize{capiref:tunemethod}}\begin{itemize}
\item {} 
\sphinxAtStartPar
\sphinxcode{\sphinxupquote{COPT\_INTPARAM\_TUNEMETHOD}} or \sphinxcode{\sphinxupquote{"TuneMethod"}}
\begin{quote}

\sphinxAtStartPar
Integer parameter.

\sphinxAtStartPar
Method for parameter tuning.

\sphinxAtStartPar
\sphinxstylestrong{Default:} \sphinxhyphen{}1

\sphinxAtStartPar
\sphinxstylestrong{Possible values:}
\begin{quote}

\sphinxAtStartPar
\sphinxhyphen{}1: Choose automatically.

\sphinxAtStartPar
0: Greedy search strategy.

\sphinxAtStartPar
1: Broader search strategy.
\end{quote}
\end{quote}

\end{itemize}
\phantomsection\label{\detokenize{capiref:tunemode}}\begin{itemize}
\item {} 
\sphinxAtStartPar
\sphinxcode{\sphinxupquote{COPT\_INTPARAM\_TUNEMODE}} or \sphinxcode{\sphinxupquote{"TuneMode"}}
\begin{quote}

\sphinxAtStartPar
Integer parameter.

\sphinxAtStartPar
Mode for parameter tuning.

\sphinxAtStartPar
\sphinxstylestrong{Default:} \sphinxhyphen{}1

\sphinxAtStartPar
\sphinxstylestrong{Possible values:}
\begin{quote}

\sphinxAtStartPar
\sphinxhyphen{}1: Choose automatically.

\sphinxAtStartPar
0: Solving time.

\sphinxAtStartPar
1: Optimal relative tolerance.

\sphinxAtStartPar
2: Objective function value.

\sphinxAtStartPar
3: The lower bound of the objective function value.
\end{quote}
\end{quote}

\end{itemize}
\phantomsection\label{\detokenize{capiref:tunemeasure}}\begin{itemize}
\item {} 
\sphinxAtStartPar
\sphinxcode{\sphinxupquote{COPT\_INTPARAM\_TUNEMEASURE}} or \sphinxcode{\sphinxupquote{"TuneMeasure"}}
\begin{quote}

\sphinxAtStartPar
Integer parameter.

\sphinxAtStartPar
Parameter tuning result calculation method.

\sphinxAtStartPar
\sphinxstylestrong{Default:} \sphinxhyphen{}1

\sphinxAtStartPar
\sphinxstylestrong{Possible values:}
\begin{quote}

\sphinxAtStartPar
\sphinxhyphen{}1: Choose automatically.

\sphinxAtStartPar
0: Calculate the average value.

\sphinxAtStartPar
1: Calculate the maximum value.
\end{quote}
\end{quote}

\end{itemize}
\phantomsection\label{\detokenize{capiref:tunepermutes}}\begin{itemize}
\item {} 
\sphinxAtStartPar
\sphinxcode{\sphinxupquote{COPT\_INTPARAM\_TUNEPERMUTES}} or \sphinxcode{\sphinxupquote{"TunePermutes"}}
\begin{quote}

\sphinxAtStartPar
Integer parameter.

\sphinxAtStartPar
Permutations for each trial parameter set. If the parameter value is 0, it will automatically set by the solver.

\sphinxAtStartPar
\sphinxstylestrong{Default:} 0

\sphinxAtStartPar
\sphinxstylestrong{Minimal:} 0

\sphinxAtStartPar
\sphinxstylestrong{Maximal:} \sphinxcode{\sphinxupquote{INT\_MAX}}
\end{quote}

\end{itemize}
\phantomsection\label{\detokenize{capiref:tuneoutputlevel}}\begin{itemize}
\item {} 
\sphinxAtStartPar
\sphinxcode{\sphinxupquote{COPT\_INTPARAM\_TUNEOUTPUTLEVEL}} or \sphinxcode{\sphinxupquote{"TuneOutputLevel"}}
\begin{quote}

\sphinxAtStartPar
Integer parameter.

\sphinxAtStartPar
Parameter tuning log output level.

\sphinxAtStartPar
\sphinxstylestrong{Default:} 2

\sphinxAtStartPar
\sphinxstylestrong{Possible values:}
\begin{quote}

\sphinxAtStartPar
0: Not displayed.

\sphinxAtStartPar
1: Display only improved parameter results.

\sphinxAtStartPar
2: Displays a summary of the results for each set of parameters.

\sphinxAtStartPar
3: Display each group of parameter results in detail.
\end{quote}
\end{quote}

\end{itemize}

\subsection{Callback related}
\label{\detokenize{capiref:callback-related}}\begin{itemize}
\item {} 
\sphinxAtStartPar
\sphinxcode{\sphinxupquote{COPT\_INTPARAM\_LAZYCONSTRAINTS}} or \sphinxcode{\sphinxupquote{"LazyConstraints"}}
\begin{quote}

\sphinxAtStartPar
Integer parameter.

\sphinxAtStartPar
Whether lazy constraints are part of the model.

\sphinxAtStartPar
\sphinxstylestrong{Default:} \sphinxhyphen{}1

\sphinxAtStartPar
\sphinxstylestrong{Possible values:}
\begin{quote}

\sphinxAtStartPar
\sphinxhyphen{}1: Choose automatically.

\sphinxAtStartPar
0: No.

\sphinxAtStartPar
1: Yes.
\end{quote}
\end{quote}

\end{itemize}

\begin{sphinxadmonition}{note}{Notes}
\begin{itemize}
\item {} 
\sphinxAtStartPar
This parameter only affects MIP.

\end{itemize}
\end{sphinxadmonition}

\subsection{GPU computing related}
\label{\detokenize{capiref:gpu-computing-related}}\begin{itemize}
\item {} 
\sphinxAtStartPar
\sphinxcode{\sphinxupquote{COPT\_INTPARAM\_GPUMODE}} or \sphinxcode{\sphinxupquote{"GPUMode"}}
\begin{quote}

\sphinxAtStartPar
Integer parameter.

\sphinxAtStartPar
Specifies GPU mode.

\sphinxAtStartPar
\sphinxstylestrong{Default:} \sphinxhyphen{}1

\sphinxAtStartPar
\sphinxstylestrong{Possible values:}
\begin{quote}

\sphinxAtStartPar
\sphinxhyphen{}1: Choose automatically. The first\sphinxhyphen{}order method (PDLP) will attempt to use the GPU by default,
while the barrier method will use the CPU by default.

\sphinxAtStartPar
0: Force CPU mode.

\sphinxAtStartPar
1: Attempt to use the standard GPU mode.

\sphinxAtStartPar
2: For the barrier method, attempt to use the high\sphinxhyphen{}performance GPU mode when solving LP problems,
which may lead to higher memory usage. For the first\sphinxhyphen{}order method (PDLP),
this is equivalent to \sphinxcode{\sphinxupquote{GPUMode=1}} (standard GPU mode).
\end{quote}
\end{quote}

\item {} 
\sphinxAtStartPar
\sphinxcode{\sphinxupquote{COPT\_INTPARAM\_GPUDEVICE}} or \sphinxcode{\sphinxupquote{"GPUDevice"}}
\begin{quote}

\sphinxAtStartPar
Integer parameter.

\sphinxAtStartPar
Specifies GPU device to use (in cases where the running machine has multiple GPUs).

\sphinxAtStartPar
\sphinxstylestrong{Default:} \sphinxhyphen{}1 (Choose automatically)

\sphinxAtStartPar
\sphinxstylestrong{Minimal:} \sphinxhyphen{}1 (Choose automatically)

\sphinxAtStartPar
\sphinxstylestrong{Maximal:} \sphinxcode{\sphinxupquote{INT\_MAX}}
\end{quote}

\end{itemize}

\subsection{Multi\sphinxhyphen{}objective Optimization}
\label{\detokenize{capiref:multi-objective-optimization}}\begin{itemize}
\item {} 
\sphinxAtStartPar
\sphinxcode{\sphinxupquote{COPT\_DBLPARAM\_MULTIOBJTIMELIMIT}} or \sphinxcode{\sphinxupquote{"MultiObjTimeLimit"}}
\begin{quote}

\sphinxAtStartPar
Double parameter.

\sphinxAtStartPar
Time limit (in seconds) for solving a multi\sphinxhyphen{}objective model.

\sphinxAtStartPar
\sphinxstylestrong{Default} 1e20

\sphinxAtStartPar
\sphinxstylestrong{Minimum} 0

\sphinxAtStartPar
\sphinxstylestrong{Maximum} 1e20
\end{quote}

\item {} 
\sphinxAtStartPar
\sphinxcode{\sphinxupquote{COPT\_INTPARAM\_MULTIOBJPARAMMODE}} or \sphinxcode{\sphinxupquote{"MultiObjParamMode"}}
\begin{quote}

\sphinxAtStartPar
Integer parameter.

\sphinxAtStartPar
Parameter usage mode for the model of each objective in multi\sphinxhyphen{}objective optimization.

\sphinxAtStartPar
\sphinxstylestrong{Default} 0

\sphinxAtStartPar
\sphinxstylestrong{Possible values}
\begin{quote}

\sphinxAtStartPar
0: Use the global solver parameters for all multi\sphinxhyphen{}objective models

\sphinxAtStartPar
1: Use the solver parameters specific to each multi\sphinxhyphen{}objective model.
\end{quote}
\end{quote}

\end{itemize}

\subsection{Other parameters}
\label{\detokenize{capiref:other-parameters}}\begin{itemize}
\item {} 
\sphinxAtStartPar
\sphinxcode{\sphinxupquote{COPT\_INTPARAM\_LOGGING}} or \sphinxcode{\sphinxupquote{"Logging"}}
\begin{quote}

\sphinxAtStartPar
Integer parameter.

\sphinxAtStartPar
Whether to print optimization logs.

\sphinxAtStartPar
\sphinxstylestrong{Default:} 1

\sphinxAtStartPar
\sphinxstylestrong{Possible values:}
\begin{quote}

\sphinxAtStartPar
0: No optimization logs.

\sphinxAtStartPar
1: Print optimization logs.
\end{quote}
\end{quote}

\item {} 
\sphinxAtStartPar
\sphinxcode{\sphinxupquote{COPT\_INTPARAM\_LOGLEVEL}} or \sphinxcode{\sphinxupquote{"LogLevel"}}
\begin{quote}

\sphinxAtStartPar
Integer parameter.

\sphinxAtStartPar
Controls the level of detail in the optimization logs.

\sphinxAtStartPar
\sphinxstylestrong{Default:} 2

\sphinxAtStartPar
\sphinxstylestrong{Possible values:}
\begin{quote}

\sphinxAtStartPar
2: Print basic optimization logs.

\sphinxAtStartPar
3: Print memory usage information in addition to basic optimization logs (for MIP problems).
\end{quote}
\end{quote}

\item {} 
\sphinxAtStartPar
\sphinxcode{\sphinxupquote{COPT\_INTPARAM\_LOGTOCONSOLE}} or \sphinxcode{\sphinxupquote{"LogToConsole"}}
\begin{quote}

\sphinxAtStartPar
Integer parameter.

\sphinxAtStartPar
Whether to print optimization logs to console.

\sphinxAtStartPar
\sphinxstylestrong{Default:} 1

\sphinxAtStartPar
\sphinxstylestrong{Possible values:}
\begin{quote}

\sphinxAtStartPar
0: No optimization logs to console.

\sphinxAtStartPar
1: Print optimization logs to console.
\end{quote}
\end{quote}

\end{itemize}

\section{API Functions}
\label{\detokenize{capiref:api-functions}}\label{\detokenize{capiref:chapapi-funcs}}
\sphinxAtStartPar
The documentations for COPT API functions are
grouped by their purposes.

\sphinxAtStartPar
All the return values of COPT API functions are integers,
and possible return values are documented in the constants section.

\subsection{Creating the environment and problem}
\label{\detokenize{capiref:creating-the-environment-and-problem}}

\subsubsection{COPT\_CreateEnvConfig}
\label{\detokenize{capiref:copt-createenvconfig}}\begin{quote}

\sphinxAtStartPar
\sphinxstylestrong{Synopsis}
\begin{quote}

\sphinxAtStartPar
\sphinxcode{\sphinxupquote{int COPT\_CreateEnvConfig(copt\_env\_config **p\_config)}}
\end{quote}

\sphinxAtStartPar
\sphinxstylestrong{Description}
\begin{quote}

\sphinxAtStartPar
Create a COPT client configuration.
\end{quote}

\sphinxAtStartPar
\sphinxstylestrong{Arguments}
\begin{quote}

\sphinxAtStartPar
\sphinxcode{\sphinxupquote{p\_config}}
\begin{quote}

\sphinxAtStartPar
Output pointer to COPT client configuration.
\end{quote}
\end{quote}
\end{quote}

\subsubsection{COPT\_DeleteEnvConfig}
\label{\detokenize{capiref:copt-deleteenvconfig}}\begin{quote}

\sphinxAtStartPar
\sphinxstylestrong{Synopsis}
\begin{quote}

\sphinxAtStartPar
\sphinxcode{\sphinxupquote{int COPT\_DeleteEnvConfig(copt\_env\_config **p\_config)}}
\end{quote}

\sphinxAtStartPar
\sphinxstylestrong{Description}
\begin{quote}

\sphinxAtStartPar
Delete COPT client configuration.
\end{quote}

\sphinxAtStartPar
\sphinxstylestrong{Arguments}
\begin{quote}

\sphinxAtStartPar
\sphinxcode{\sphinxupquote{p\_config}}
\begin{quote}

\sphinxAtStartPar
Input pointer to COPT client configuration.
\end{quote}
\end{quote}
\end{quote}

\subsubsection{COPT\_SetEnvConfig}
\label{\detokenize{capiref:copt-setenvconfig}}\begin{quote}

\sphinxAtStartPar
\sphinxstylestrong{Synopsis}
\begin{quote}

\sphinxAtStartPar
\sphinxcode{\sphinxupquote{int COPT\_SetEnvConfig(copt\_env\_config *config, const char *name, const char *value)}}
\end{quote}

\sphinxAtStartPar
\sphinxstylestrong{Description}
\begin{quote}

\sphinxAtStartPar
Set COPT client configuration.
\end{quote}

\sphinxAtStartPar
\sphinxstylestrong{Arguments}
\begin{quote}

\sphinxAtStartPar
\sphinxcode{\sphinxupquote{config}}
\begin{quote}

\sphinxAtStartPar
COPT client configuration.
\end{quote}

\sphinxAtStartPar
\sphinxcode{\sphinxupquote{name}}
\begin{quote}

\sphinxAtStartPar
Name of configuration parameter.
\end{quote}

\sphinxAtStartPar
\sphinxcode{\sphinxupquote{value}}
\begin{quote}

\sphinxAtStartPar
Value of configuration parameter.
\end{quote}
\end{quote}
\end{quote}

\subsubsection{COPT\_CreateEnv}
\label{\detokenize{capiref:copt-createenv}}\begin{quote}

\sphinxAtStartPar
\sphinxstylestrong{Synopsis}
\begin{quote}

\sphinxAtStartPar
\sphinxcode{\sphinxupquote{int COPT\_CreateEnv(copt\_env **p\_env)}}
\end{quote}

\sphinxAtStartPar
\sphinxstylestrong{Description}
\begin{quote}

\sphinxAtStartPar
Creates a COPT environment.

\sphinxAtStartPar
Calling this function is the first step when using the COPT library.
It validates the license, and when the license is okay,
the resulting environment variable will allow for
creating COPT problems.
When the license validation fails, more information can be
obtained using \sphinxcode{\sphinxupquote{COPT\_GetLicenseMsg}} to help identify the issue.
\end{quote}

\sphinxAtStartPar
\sphinxstylestrong{Arguments}
\begin{quote}

\sphinxAtStartPar
\sphinxcode{\sphinxupquote{p\_env}}
\begin{quote}

\sphinxAtStartPar
The output pointer to a variable holding COPT environment.
\end{quote}
\end{quote}
\end{quote}

\subsubsection{COPT\_CreateEnvWithPath}
\label{\detokenize{capiref:copt-createenvwithpath}}\begin{quote}

\sphinxAtStartPar
\sphinxstylestrong{Synopsis}
\begin{quote}

\sphinxAtStartPar
\sphinxcode{\sphinxupquote{int COPT\_CreateEnvWithPath(const char *licDir, copt\_env **p\_env)}}
\end{quote}

\sphinxAtStartPar
\sphinxstylestrong{Description}
\begin{quote}

\sphinxAtStartPar
Creates a COPT environment, directory of license files is specified by
argument \sphinxcode{\sphinxupquote{licDir}}.

\sphinxAtStartPar
Calling this function is the first step when using the COPT library.
It validates the license, and when the license is okay,
the resulting environment variable will allow for
creating COPT problems.
When the license validation fails, more information can be
obtained using \sphinxcode{\sphinxupquote{COPT\_GetLicenseMsg}} to help identify the issue.
\end{quote}

\sphinxAtStartPar
\sphinxstylestrong{Arguments}
\begin{quote}

\sphinxAtStartPar
\sphinxcode{\sphinxupquote{licDir}}
\begin{quote}

\sphinxAtStartPar
Directory of license files.
\end{quote}

\sphinxAtStartPar
\sphinxcode{\sphinxupquote{p\_env}}
\begin{quote}

\sphinxAtStartPar
Output pointer to a variable holding COPT environment.
\end{quote}
\end{quote}
\end{quote}

\subsubsection{COPT\_CreateEnvWithConfig}
\label{\detokenize{capiref:copt-createenvwithconfig}}\begin{quote}

\sphinxAtStartPar
\sphinxstylestrong{Synopsis}
\begin{quote}

\sphinxAtStartPar
\sphinxcode{\sphinxupquote{int COPT\_CreateEnvWithConfig(copt\_env\_config *config, copt\_env **p\_env)}}
\end{quote}

\sphinxAtStartPar
\sphinxstylestrong{Description}
\begin{quote}

\sphinxAtStartPar
Creates a COPT environment, client configuration is specified by
argument \sphinxcode{\sphinxupquote{config}}.

\sphinxAtStartPar
Calling this function is the first step when using the COPT library.
It validates the client configuration, and when the license is okay,
the resulting environment variable will allow for
creating COPT problems.
When the license validation fails, more information can be
obtained using \sphinxcode{\sphinxupquote{COPT\_GetLicenseMsg}} to help identify the issue.
\end{quote}

\sphinxAtStartPar
\sphinxstylestrong{Arguments}
\begin{quote}

\sphinxAtStartPar
\sphinxcode{\sphinxupquote{config}}
\begin{quote}

\sphinxAtStartPar
Client configuration.
\end{quote}

\sphinxAtStartPar
\sphinxcode{\sphinxupquote{p\_env}}
\begin{quote}

\sphinxAtStartPar
Output pointer to a variable holding COPT environment.
\end{quote}
\end{quote}
\end{quote}

\subsubsection{COPT\_DeleteEnv}
\label{\detokenize{capiref:copt-deleteenv}}\begin{quote}

\sphinxAtStartPar
\sphinxstylestrong{Synopsis}
\begin{quote}

\sphinxAtStartPar
\sphinxcode{\sphinxupquote{int COPT\_DeleteEnv(copt\_env **p\_env)}}
\end{quote}

\sphinxAtStartPar
\sphinxstylestrong{Description}
\begin{quote}

\sphinxAtStartPar
Deletes the COPT environment created by \sphinxcode{\sphinxupquote{COPT\_CreateEnv}}.
\end{quote}

\sphinxAtStartPar
\sphinxstylestrong{Arguments}
\begin{quote}

\sphinxAtStartPar
\sphinxcode{\sphinxupquote{p\_env}}
\begin{quote}

\sphinxAtStartPar
Input pointer to a variable holding COPT environment.
\end{quote}
\end{quote}
\end{quote}

\subsubsection{COPT\_GetLicenseMsg}
\label{\detokenize{capiref:copt-getlicensemsg}}\begin{quote}

\sphinxAtStartPar
\sphinxstylestrong{Synopsis}
\begin{quote}

\sphinxAtStartPar
\sphinxcode{\sphinxupquote{int COPT\_GetLicenseMsg(copt\_env *env, char *buff, int buffSize)}}
\end{quote}

\sphinxAtStartPar
\sphinxstylestrong{Description}
\begin{quote}

\sphinxAtStartPar
Returns a C\sphinxhyphen{}style string regarding the license validation information.
Please refer to this function when \sphinxcode{\sphinxupquote{COPT\_CreateEnv}} fails.
\end{quote}

\sphinxAtStartPar
\sphinxstylestrong{Arguments}
\begin{quote}

\sphinxAtStartPar
\sphinxcode{\sphinxupquote{env}}
\begin{quote}

\sphinxAtStartPar
The COPT environment.
\end{quote}

\sphinxAtStartPar
\sphinxcode{\sphinxupquote{buff}}
\begin{quote}

\sphinxAtStartPar
A buffer for holding the resulting string.
\end{quote}

\sphinxAtStartPar
\sphinxcode{\sphinxupquote{buffSize}}
\begin{quote}

\sphinxAtStartPar
The size of the above buffer.
\end{quote}
\end{quote}
\end{quote}

\subsubsection{COPT\_CreateProb}
\label{\detokenize{capiref:copt-createprob}}\begin{quote}

\sphinxAtStartPar
\sphinxstylestrong{Synopsis}
\begin{quote}

\sphinxAtStartPar
\sphinxcode{\sphinxupquote{int COPT\_CreateProb(copt\_env *env, copt\_prob **p\_prob)}}
\end{quote}

\sphinxAtStartPar
\sphinxstylestrong{Description}
\begin{quote}

\sphinxAtStartPar
Creates an empty COPT problem.
\end{quote}

\sphinxAtStartPar
\sphinxstylestrong{Arguments}
\begin{quote}

\sphinxAtStartPar
\sphinxcode{\sphinxupquote{env}}
\begin{quote}

\sphinxAtStartPar
The COPT environment.
\end{quote}

\sphinxAtStartPar
\sphinxcode{\sphinxupquote{p\_prob}}
\begin{quote}

\sphinxAtStartPar
Output pointer to a variable holding the COPT problem.
\end{quote}
\end{quote}
\end{quote}

\subsubsection{COPT\_CreateCopy}
\label{\detokenize{capiref:copt-createcopy}}\begin{quote}

\sphinxAtStartPar
\sphinxstylestrong{Synopsis}
\begin{quote}

\sphinxAtStartPar
\sphinxcode{\sphinxupquote{int COPT\_CreateCopy(copt\_prob *src\_prob, copt\_prob **p\_dst\_prob)}}
\end{quote}

\sphinxAtStartPar
\sphinxstylestrong{Description}
\begin{quote}

\sphinxAtStartPar
Create a deep\sphinxhyphen{}copy of an COPT problem.

\sphinxAtStartPar
\sphinxstylestrong{Note:} The parameter settings will be copied too. To solve the
copied problem with different parameters, users should reset
its parameters to default by calling \sphinxcode{\sphinxupquote{COPT\_ResetParam}} first,
and then set parameters as needed.
\end{quote}

\sphinxAtStartPar
\sphinxstylestrong{Arguments}
\begin{quote}

\sphinxAtStartPar
\sphinxcode{\sphinxupquote{src\_prob}}
\begin{quote}

\sphinxAtStartPar
The pointer to a varialbe hoding the COPT problem to be copied.
\end{quote}

\sphinxAtStartPar
\sphinxcode{\sphinxupquote{p\_dst\_prob}}
\begin{quote}

\sphinxAtStartPar
Output pointer to a variable hodling the copied COPT problem.
\end{quote}
\end{quote}
\end{quote}

\subsubsection{COPT\_ClearProb}
\label{\detokenize{capiref:copt-clearprob}}\begin{quote}

\sphinxAtStartPar
\sphinxstylestrong{Synopsis}
\begin{quote}

\sphinxAtStartPar
\sphinxcode{\sphinxupquote{int COPT\_ClearProb(copt\_prob *prob)}}
\end{quote}

\sphinxAtStartPar
\sphinxstylestrong{Description}
\begin{quote}

\sphinxAtStartPar
Clear COPT problem data (excluding callback function).
\end{quote}

\sphinxAtStartPar
\sphinxstylestrong{Arguments}
\begin{quote}

\sphinxAtStartPar
\sphinxcode{\sphinxupquote{prob}}
\begin{quote}

\sphinxAtStartPar
COPT problem.
\end{quote}
\end{quote}
\end{quote}

\subsubsection{COPT\_DeleteProb}
\label{\detokenize{capiref:copt-deleteprob}}\begin{quote}

\sphinxAtStartPar
\sphinxstylestrong{Synopsis}
\begin{quote}

\sphinxAtStartPar
\sphinxcode{\sphinxupquote{int COPT\_DeleteProb(copt\_prob **p\_prob)}}
\end{quote}

\sphinxAtStartPar
\sphinxstylestrong{Description}
\begin{quote}

\sphinxAtStartPar
Deletes the COPT problem created using \sphinxcode{\sphinxupquote{COPT\_CreateProb}}
\end{quote}

\sphinxAtStartPar
\sphinxstylestrong{Arguments}
\begin{quote}

\sphinxAtStartPar
\sphinxcode{\sphinxupquote{p\_prob}}
\begin{quote}

\sphinxAtStartPar
Input pointer to a variable holding the COPT problem.
\end{quote}
\end{quote}
\end{quote}

\subsection{Building and modifying a problem}
\label{\detokenize{capiref:building-and-modifying-a-problem}}\label{\detokenize{capiref:chapapi-model}}

\subsubsection{COPT\_LoadProb}
\label{\detokenize{capiref:copt-loadprob}}\begin{quote}

\sphinxAtStartPar
\sphinxstylestrong{Synopsis}
\begin{quote}

\sphinxAtStartPar
\sphinxcode{\sphinxupquote{int COPT\_LoadProb(}}
\begin{quote}

\sphinxAtStartPar
\sphinxcode{\sphinxupquote{copt\_prob *prob,}}

\sphinxAtStartPar
\sphinxcode{\sphinxupquote{int nCol,}}

\sphinxAtStartPar
\sphinxcode{\sphinxupquote{int nRow,}}

\sphinxAtStartPar
\sphinxcode{\sphinxupquote{int iObjSense,}}

\sphinxAtStartPar
\sphinxcode{\sphinxupquote{double dObjConst,}}

\sphinxAtStartPar
\sphinxcode{\sphinxupquote{const double *obj,}}

\sphinxAtStartPar
\sphinxcode{\sphinxupquote{const int *colMatBeg,}}

\sphinxAtStartPar
\sphinxcode{\sphinxupquote{const int *colMatCnt,}}

\sphinxAtStartPar
\sphinxcode{\sphinxupquote{const int *colMatIdx,}}

\sphinxAtStartPar
\sphinxcode{\sphinxupquote{const double *colMatElem,}}

\sphinxAtStartPar
\sphinxcode{\sphinxupquote{const char *colType,}}

\sphinxAtStartPar
\sphinxcode{\sphinxupquote{const double *colLower,}}

\sphinxAtStartPar
\sphinxcode{\sphinxupquote{const double *colUpper,}}

\sphinxAtStartPar
\sphinxcode{\sphinxupquote{const char *rowSense,}}

\sphinxAtStartPar
\sphinxcode{\sphinxupquote{const double *rowBound,}}

\sphinxAtStartPar
\sphinxcode{\sphinxupquote{const double *rowUpper,}}

\sphinxAtStartPar
\sphinxcode{\sphinxupquote{char const *const *colNames,}}

\sphinxAtStartPar
\sphinxcode{\sphinxupquote{char const *const *rowNames)}}
\end{quote}
\end{quote}

\sphinxAtStartPar
\sphinxstylestrong{Description}
\begin{quote}

\sphinxAtStartPar
Loads a problem defined by arrays.
\end{quote}

\sphinxAtStartPar
\sphinxstylestrong{Arguments}
\begin{quote}

\sphinxAtStartPar
\sphinxcode{\sphinxupquote{prob}}
\begin{quote}

\sphinxAtStartPar
The COPT problem.
\end{quote}

\sphinxAtStartPar
\sphinxcode{\sphinxupquote{nCol}}
\begin{quote}

\sphinxAtStartPar
Number of variables (coefficient matrix columns).
\end{quote}

\sphinxAtStartPar
\sphinxcode{\sphinxupquote{nRow}}
\begin{quote}

\sphinxAtStartPar
Number of constraints (coefficient matrix rows).
\end{quote}

\sphinxAtStartPar
\sphinxcode{\sphinxupquote{iObjSense}}
\begin{quote}

\sphinxAtStartPar
The optimization sense, either \sphinxcode{\sphinxupquote{COPT\_MAXIMIZE}} or \sphinxcode{\sphinxupquote{COPT\_MINIMIZE}}.
\end{quote}

\sphinxAtStartPar
\sphinxcode{\sphinxupquote{dObjConst}}
\begin{quote}

\sphinxAtStartPar
The constant part of the objective function.
\end{quote}

\sphinxAtStartPar
\sphinxcode{\sphinxupquote{obj}}
\begin{quote}

\sphinxAtStartPar
Objective coefficients of variables.
\end{quote}

\sphinxAtStartPar
\sphinxcode{\sphinxupquote{colMatBeg, colMatCnt, colMatIdx}} and \sphinxcode{\sphinxupquote{colMatElem}}
\begin{quote}

\sphinxAtStartPar
Defines the coefficient matrix in compressed column storage (CCS) format.
Please see \sphinxstylestrong{other information} for an example of the CCS format.

\sphinxAtStartPar
If \sphinxcode{\sphinxupquote{colMatCnt}} is \sphinxcode{\sphinxupquote{NULL}}, \sphinxcode{\sphinxupquote{colMatBeg}} will need to have
length of \sphinxcode{\sphinxupquote{nCol+1}}, and the begin and end pointers to the
i\sphinxhyphen{}th matrix column coefficients are defined using
\sphinxcode{\sphinxupquote{colMatBeg{[}i{]}}} and \sphinxcode{\sphinxupquote{colMatBeg{[}i+1{]}}}.

\sphinxAtStartPar
If \sphinxcode{\sphinxupquote{colMatCnt}} is provided, the begin and end pointers
to the i\sphinxhyphen{}th column coefficients are defined using
\sphinxcode{\sphinxupquote{colMatBeg{[}i{]}}} and \sphinxcode{\sphinxupquote{colMatBeg{[}i{]} + colMatCnt{[}i{]}}}.
\end{quote}

\sphinxAtStartPar
\sphinxcode{\sphinxupquote{colType}}
\begin{quote}

\sphinxAtStartPar
Types of variables.

\sphinxAtStartPar
If \sphinxcode{\sphinxupquote{colType}} is \sphinxcode{\sphinxupquote{NULL}}, all variables will be continuous.
\end{quote}

\sphinxAtStartPar
\sphinxcode{\sphinxupquote{colLower}} and \sphinxcode{\sphinxupquote{colUpper}}
\begin{quote}

\sphinxAtStartPar
Lower and upper bounds of variables.

\sphinxAtStartPar
If \sphinxcode{\sphinxupquote{colLower}} is \sphinxcode{\sphinxupquote{NULL}}, lower bounds will be 0.

\sphinxAtStartPar
If \sphinxcode{\sphinxupquote{colUpper}} is \sphinxcode{\sphinxupquote{NULL}}, upper bounds will be infinity.
\end{quote}

\sphinxAtStartPar
\sphinxcode{\sphinxupquote{rowSense}}
\begin{quote}

\sphinxAtStartPar
Senses of constraint.

\sphinxAtStartPar
Please refer to the list of all senses constants for
all the supported types.

\sphinxAtStartPar
If \sphinxcode{\sphinxupquote{rowSense}} is \sphinxcode{\sphinxupquote{NULL}}, then \sphinxcode{\sphinxupquote{rowBound}} and \sphinxcode{\sphinxupquote{rowUpper}}
will be treated as lower and upper bounds for constraints.
This is the recommended method for defining constraints.

\sphinxAtStartPar
If \sphinxcode{\sphinxupquote{rowSense}} is provided, then \sphinxcode{\sphinxupquote{rowBound}} and \sphinxcode{\sphinxupquote{rowUpper}}
will be treated as RHS and \sphinxstylestrong{range} for constraints.
In this case, \sphinxcode{\sphinxupquote{rowUpper}} is only required when there
are \sphinxcode{\sphinxupquote{COPT\_RANGE}} constraints, where the
\begin{quote}

\sphinxAtStartPar
lower bound is \sphinxcode{\sphinxupquote{rowBound{[}i{]} \sphinxhyphen{} fabs(rowUpper{[}i{]})}}

\sphinxAtStartPar
upper bound is \sphinxcode{\sphinxupquote{rowBound{[}i{]}}}
\end{quote}
\end{quote}

\sphinxAtStartPar
\sphinxcode{\sphinxupquote{rowBound}}
\begin{quote}

\sphinxAtStartPar
Lower bounds or RHS of constraints.
\end{quote}

\sphinxAtStartPar
\sphinxcode{\sphinxupquote{rowUpper}}
\begin{quote}

\sphinxAtStartPar
Upper bounds or \sphinxstylestrong{range} of constraints.
\end{quote}

\sphinxAtStartPar
\sphinxcode{\sphinxupquote{colNames}} and \sphinxcode{\sphinxupquote{rowNames}}
\begin{quote}

\sphinxAtStartPar
Names of variables and constraints. Can be \sphinxcode{\sphinxupquote{NULL}}.
\end{quote}
\end{quote}

\sphinxAtStartPar
\sphinxstylestrong{Other information}
\begin{quote}

\sphinxAtStartPar
The compressed column storage (CCS) is a common format for
storing sparse matrix. We demonstrate how to store the example
matrix with 4 columns and 3 rows in the CCS format.
\end{quote}
\end{quote}
\begin{equation}\label{equation:capiref:capiref:9}
\begin{split}A = \begin{bmatrix}
        1.1 & 1.2 &     &     \\
            & 2.2 & 2.3 &     \\
            &     & 3.3 & 3.4 \\
    \end{bmatrix}\end{split}
\end{equation}
\begin{sphinxVerbatim}[commandchars=\\\{\}]
\PYG{c+c1}{// Compressed column storage using colMatBeg}
\PYG{n}{colMatBeg}\PYG{p}{[}\PYG{l+m+mi}{5}\PYG{p}{]}\PYG{+w}{  }\PYG{o}{=}\PYG{+w}{ }\PYG{p}{\PYGZob{}}\PYG{+w}{  }\PYG{l+m+mi}{0}\PYG{p}{,}\PYG{+w}{   }\PYG{l+m+mi}{1}\PYG{p}{,}\PYG{+w}{   }\PYG{l+m+mi}{3}\PYG{p}{,}\PYG{+w}{   }\PYG{l+m+mi}{5}\PYG{p}{,}\PYG{+w}{   }\PYG{l+m+mi}{6}\PYG{p}{\PYGZcb{}}\PYG{p}{;}
\PYG{n}{colMatIdx}\PYG{p}{[}\PYG{l+m+mi}{6}\PYG{p}{]}\PYG{+w}{  }\PYG{o}{=}\PYG{+w}{ }\PYG{p}{\PYGZob{}}\PYG{+w}{  }\PYG{l+m+mi}{0}\PYG{p}{,}\PYG{+w}{   }\PYG{l+m+mi}{0}\PYG{p}{,}\PYG{+w}{   }\PYG{l+m+mi}{1}\PYG{p}{,}\PYG{+w}{   }\PYG{l+m+mi}{1}\PYG{p}{,}\PYG{+w}{   }\PYG{l+m+mi}{2}\PYG{p}{,}\PYG{+w}{   }\PYG{l+m+mi}{2}\PYG{p}{\PYGZcb{}}\PYG{p}{;}
\PYG{n}{colMatElem}\PYG{p}{[}\PYG{l+m+mi}{6}\PYG{p}{]}\PYG{+w}{ }\PYG{o}{=}\PYG{+w}{ }\PYG{p}{\PYGZob{}}\PYG{l+m+mf}{1.1}\PYG{p}{,}\PYG{+w}{ }\PYG{l+m+mf}{1.2}\PYG{p}{,}\PYG{+w}{ }\PYG{l+m+mf}{2.2}\PYG{p}{,}\PYG{+w}{ }\PYG{l+m+mf}{2.3}\PYG{p}{,}\PYG{+w}{ }\PYG{l+m+mf}{3.3}\PYG{p}{,}\PYG{+w}{ }\PYG{l+m+mf}{3.4}\PYG{p}{\PYGZcb{}}\PYG{p}{;}

\PYG{c+c1}{// Compressed column storage using both colMatBeg and colMatCnt.}
\PYG{c+c1}{// The * in the example represents unused spaces.}
\PYG{n}{colMatBeg}\PYG{p}{[}\PYG{l+m+mi}{4}\PYG{p}{]}\PYG{+w}{  }\PYG{o}{=}\PYG{+w}{ }\PYG{p}{\PYGZob{}}\PYG{+w}{  }\PYG{l+m+mi}{0}\PYG{p}{,}\PYG{+w}{   }\PYG{l+m+mi}{1}\PYG{p}{,}\PYG{+w}{   }\PYG{l+m+mi}{5}\PYG{p}{,}\PYG{+w}{   }\PYG{l+m+mi}{7}\PYG{p}{\PYGZcb{}}\PYG{p}{;}
\PYG{n}{colMatCnt}\PYG{p}{[}\PYG{l+m+mi}{4}\PYG{p}{]}\PYG{+w}{  }\PYG{o}{=}\PYG{+w}{ }\PYG{p}{\PYGZob{}}\PYG{+w}{  }\PYG{l+m+mi}{1}\PYG{p}{,}\PYG{+w}{   }\PYG{l+m+mi}{2}\PYG{p}{,}\PYG{+w}{   }\PYG{l+m+mi}{2}\PYG{p}{,}\PYG{+w}{   }\PYG{l+m+mi}{1}\PYG{p}{\PYGZcb{}}\PYG{p}{;}
\PYG{n}{colMatIdx}\PYG{p}{[}\PYG{l+m+mi}{6}\PYG{p}{]}\PYG{+w}{  }\PYG{o}{=}\PYG{+w}{ }\PYG{p}{\PYGZob{}}\PYG{+w}{  }\PYG{l+m+mi}{0}\PYG{p}{,}\PYG{+w}{   }\PYG{l+m+mi}{0}\PYG{p}{,}\PYG{+w}{   }\PYG{l+m+mi}{1}\PYG{p}{,}\PYG{+w}{   }\PYG{l+m+mi}{1}\PYG{p}{,}\PYG{+w}{   }\PYG{l+m+mi}{2}\PYG{p}{,}\PYG{+w}{   }\PYG{o}{*}\PYG{p}{,}\PYG{+w}{   }\PYG{o}{*}\PYG{p}{,}\PYG{+w}{   }\PYG{l+m+mi}{2}\PYG{p}{\PYGZcb{}}\PYG{p}{;}
\PYG{n}{colMatElem}\PYG{p}{[}\PYG{l+m+mi}{6}\PYG{p}{]}\PYG{+w}{ }\PYG{o}{=}\PYG{+w}{ }\PYG{p}{\PYGZob{}}\PYG{l+m+mf}{1.1}\PYG{p}{,}\PYG{+w}{ }\PYG{l+m+mf}{1.2}\PYG{p}{,}\PYG{+w}{ }\PYG{l+m+mf}{2.2}\PYG{p}{,}\PYG{+w}{ }\PYG{l+m+mf}{2.3}\PYG{p}{,}\PYG{+w}{ }\PYG{l+m+mf}{3.3}\PYG{p}{,}\PYG{+w}{   }\PYG{o}{*}\PYG{p}{,}\PYG{+w}{   }\PYG{o}{*}\PYG{p}{,}\PYG{+w}{ }\PYG{l+m+mf}{3.4}\PYG{p}{\PYGZcb{}}\PYG{p}{;}
\end{sphinxVerbatim}

\subsubsection{COPT\_AddCol}
\label{\detokenize{capiref:copt-addcol}}\begin{quote}

\sphinxAtStartPar
\sphinxstylestrong{Synopsis}
\begin{quote}

\sphinxAtStartPar
\sphinxcode{\sphinxupquote{int COPT\_AddCol(}}
\begin{quote}

\sphinxAtStartPar
\sphinxcode{\sphinxupquote{copt\_prob *prob,}}

\sphinxAtStartPar
\sphinxcode{\sphinxupquote{double dColObj,}}

\sphinxAtStartPar
\sphinxcode{\sphinxupquote{int nColMatCnt,}}

\sphinxAtStartPar
\sphinxcode{\sphinxupquote{const int *colMatIdx,}}

\sphinxAtStartPar
\sphinxcode{\sphinxupquote{const double *colMatElem,}}

\sphinxAtStartPar
\sphinxcode{\sphinxupquote{char cColType,}}

\sphinxAtStartPar
\sphinxcode{\sphinxupquote{double dColLower,}}

\sphinxAtStartPar
\sphinxcode{\sphinxupquote{double dColUpper,}}

\sphinxAtStartPar
\sphinxcode{\sphinxupquote{const char *colName)}}
\end{quote}
\end{quote}

\sphinxAtStartPar
\sphinxstylestrong{Description}
\begin{quote}

\sphinxAtStartPar
Adds one variable (column) to the problem.
\end{quote}

\sphinxAtStartPar
\sphinxstylestrong{Arguments}
\begin{quote}

\sphinxAtStartPar
\sphinxcode{\sphinxupquote{prob}}
\begin{quote}

\sphinxAtStartPar
The COPT problem.
\end{quote}

\sphinxAtStartPar
\sphinxcode{\sphinxupquote{dColObj}}
\begin{quote}

\sphinxAtStartPar
The objective coefficient of the variable.
\end{quote}

\sphinxAtStartPar
\sphinxcode{\sphinxupquote{nColMatCnt}}
\begin{quote}

\sphinxAtStartPar
Number of non\sphinxhyphen{}zero elements in the column.
\end{quote}

\sphinxAtStartPar
\sphinxcode{\sphinxupquote{colMatIdx}}
\begin{quote}

\sphinxAtStartPar
Row index of non\sphinxhyphen{}zero elements in the column.
\end{quote}

\sphinxAtStartPar
\sphinxcode{\sphinxupquote{colMatElem}}
\begin{quote}

\sphinxAtStartPar
Values of non\sphinxhyphen{}zero elements in the column.
\end{quote}

\sphinxAtStartPar
\sphinxcode{\sphinxupquote{cColType}}
\begin{quote}

\sphinxAtStartPar
The type of the variable.
\end{quote}

\sphinxAtStartPar
\sphinxcode{\sphinxupquote{dColLower}} and \sphinxcode{\sphinxupquote{dColUpper}}
\begin{quote}

\sphinxAtStartPar
The lower and upper bounds of the variable.
\end{quote}

\sphinxAtStartPar
\sphinxcode{\sphinxupquote{colName}}
\begin{quote}

\sphinxAtStartPar
The name of the variable. Can be \sphinxcode{\sphinxupquote{NULL}}.
\end{quote}
\end{quote}
\end{quote}

\subsubsection{COPT\_AddPSDCol}
\label{\detokenize{capiref:copt-addpsdcol}}\begin{quote}

\sphinxAtStartPar
\sphinxstylestrong{Synopsis}
\begin{quote}

\sphinxAtStartPar
\sphinxcode{\sphinxupquote{int COPT\_AddPSDCol(copt\_prob *prob, int colDim, const char *name)}}
\end{quote}

\sphinxAtStartPar
\sphinxstylestrong{Description}
\begin{quote}

\sphinxAtStartPar
Add a PSD variable to the problem.
\end{quote}

\sphinxAtStartPar
\sphinxstylestrong{Arguments}
\begin{quote}

\sphinxAtStartPar
\sphinxcode{\sphinxupquote{prob}}
\begin{quote}

\sphinxAtStartPar
The COPT problem.
\end{quote}

\sphinxAtStartPar
\sphinxcode{\sphinxupquote{colDim}}
\begin{quote}

\sphinxAtStartPar
Dimension of new PSD variable.
\end{quote}

\sphinxAtStartPar
\sphinxcode{\sphinxupquote{name}}
\begin{quote}

\sphinxAtStartPar
Name of new PSD variable. Can be \sphinxcode{\sphinxupquote{NULL}}.
\end{quote}
\end{quote}
\end{quote}

\subsubsection{COPT\_AddRow}
\label{\detokenize{capiref:copt-addrow}}\begin{quote}

\sphinxAtStartPar
\sphinxstylestrong{Synopsis}
\begin{quote}

\sphinxAtStartPar
\sphinxcode{\sphinxupquote{int COPT\_AddRow(}}
\begin{quote}

\sphinxAtStartPar
\sphinxcode{\sphinxupquote{copt\_prob *prob,}}

\sphinxAtStartPar
\sphinxcode{\sphinxupquote{int nRowMatCnt,}}

\sphinxAtStartPar
\sphinxcode{\sphinxupquote{const int *rowMatIdx,}}

\sphinxAtStartPar
\sphinxcode{\sphinxupquote{const double *rowMatElem,}}

\sphinxAtStartPar
\sphinxcode{\sphinxupquote{char cRowSense,}}

\sphinxAtStartPar
\sphinxcode{\sphinxupquote{double dRowBound,}}

\sphinxAtStartPar
\sphinxcode{\sphinxupquote{double dRowUpper,}}

\sphinxAtStartPar
\sphinxcode{\sphinxupquote{const char *rowName)}}
\end{quote}
\end{quote}

\sphinxAtStartPar
\sphinxstylestrong{Description}
\begin{quote}

\sphinxAtStartPar
Adds one constraint (row) to the problem.
\end{quote}

\sphinxAtStartPar
\sphinxstylestrong{Arguments}
\begin{quote}

\sphinxAtStartPar
\sphinxcode{\sphinxupquote{prob}}
\begin{quote}

\sphinxAtStartPar
The COPT problem.
\end{quote}

\sphinxAtStartPar
\sphinxcode{\sphinxupquote{nRowMatCnt}}
\begin{quote}

\sphinxAtStartPar
Number of non\sphinxhyphen{}zero elements in the row.
\end{quote}

\sphinxAtStartPar
\sphinxcode{\sphinxupquote{rowMatIdx}}
\begin{quote}

\sphinxAtStartPar
Column index of non\sphinxhyphen{}zero elements in the row.
\end{quote}

\sphinxAtStartPar
\sphinxcode{\sphinxupquote{rowMatElem}}
\begin{quote}

\sphinxAtStartPar
Values of non\sphinxhyphen{}zero elements in the row.
\end{quote}

\sphinxAtStartPar
\sphinxcode{\sphinxupquote{cRowSense}}
\begin{quote}

\sphinxAtStartPar
The sense of the row.

\sphinxAtStartPar
Please refer to the list of all senses constants for all the supported types.

\sphinxAtStartPar
If \sphinxcode{\sphinxupquote{cRowSense}} is 0, then \sphinxcode{\sphinxupquote{dRowBound}} and \sphinxcode{\sphinxupquote{dRowUpper}}
will be treated as lower and upper bounds for the constraint.
This is the recommended method for defining constraints.

\sphinxAtStartPar
If \sphinxcode{\sphinxupquote{cRowSense}} is provided, then \sphinxcode{\sphinxupquote{dRowBound}} and \sphinxcode{\sphinxupquote{dRowUpper}}
will be treated as RHS and \sphinxstylestrong{range} for the constraint.
In this case, \sphinxcode{\sphinxupquote{dRowUpper}}  is only required when
\sphinxcode{\sphinxupquote{cRowSense = COPT\_RANGE}}, where
\begin{quote}

\sphinxAtStartPar
lower bound is \sphinxcode{\sphinxupquote{dRowBound \sphinxhyphen{} dRowUpper}}

\sphinxAtStartPar
upper bound is \sphinxcode{\sphinxupquote{dRowBound}}
\end{quote}
\end{quote}

\sphinxAtStartPar
\sphinxcode{\sphinxupquote{dRowBound}}
\begin{quote}

\sphinxAtStartPar
Lower bound or RHS of the constraint.
\end{quote}

\sphinxAtStartPar
\sphinxcode{\sphinxupquote{dRowUpper}}
\begin{quote}

\sphinxAtStartPar
Upper bound or \sphinxstylestrong{range} of the constraint.
\end{quote}

\sphinxAtStartPar
\sphinxcode{\sphinxupquote{rowName}}
\begin{quote}

\sphinxAtStartPar
The name of the constraint. Can be \sphinxcode{\sphinxupquote{NULL}}.
\end{quote}
\end{quote}
\end{quote}

\subsubsection{COPT\_AddCols}
\label{\detokenize{capiref:copt-addcols}}\begin{quote}

\sphinxAtStartPar
\sphinxstylestrong{Synopsis}
\begin{quote}

\sphinxAtStartPar
\sphinxcode{\sphinxupquote{int COPT\_AddCols(}}
\begin{quote}

\sphinxAtStartPar
\sphinxcode{\sphinxupquote{copt\_prob *prob,}}

\sphinxAtStartPar
\sphinxcode{\sphinxupquote{int nAddCol,}}

\sphinxAtStartPar
\sphinxcode{\sphinxupquote{const double *colObj,}}

\sphinxAtStartPar
\sphinxcode{\sphinxupquote{const int *colMatBeg,}}

\sphinxAtStartPar
\sphinxcode{\sphinxupquote{const int *colMatCnt,}}

\sphinxAtStartPar
\sphinxcode{\sphinxupquote{const int *colMatIdx,}}

\sphinxAtStartPar
\sphinxcode{\sphinxupquote{const double *colMatElem,}}

\sphinxAtStartPar
\sphinxcode{\sphinxupquote{const char *colType,}}

\sphinxAtStartPar
\sphinxcode{\sphinxupquote{const double *colLower,}}

\sphinxAtStartPar
\sphinxcode{\sphinxupquote{const double *colUpper,}}

\sphinxAtStartPar
\sphinxcode{\sphinxupquote{char const *const *colNames)}}
\end{quote}
\end{quote}

\sphinxAtStartPar
\sphinxstylestrong{Description}
\begin{quote}

\sphinxAtStartPar
Adds \sphinxcode{\sphinxupquote{nAddCol}} variables (columns) to the problem.
\end{quote}

\sphinxAtStartPar
\sphinxstylestrong{Arguments}
\begin{quote}

\sphinxAtStartPar
\sphinxcode{\sphinxupquote{prob}}
\begin{quote}

\sphinxAtStartPar
The COPT problem.
\end{quote}

\sphinxAtStartPar
\sphinxcode{\sphinxupquote{nAddCol}}
\begin{quote}

\sphinxAtStartPar
Number of new variables.
\end{quote}

\sphinxAtStartPar
\sphinxcode{\sphinxupquote{colObj}}
\begin{quote}

\sphinxAtStartPar
Objective coefficients of new variables.
\end{quote}

\sphinxAtStartPar
\sphinxcode{\sphinxupquote{colMatBeg, colMatCnt, colMatIdx}} and \sphinxcode{\sphinxupquote{colMatElem}}
\begin{quote}

\sphinxAtStartPar
Defines the coefficient matrix in compressed column storage (CCS) format.
Please see \sphinxstylestrong{other information} of \sphinxcode{\sphinxupquote{COPT\_LoadProb}} for
an example of the CCS format.
\end{quote}

\sphinxAtStartPar
\sphinxcode{\sphinxupquote{colType}}
\begin{quote}

\sphinxAtStartPar
Types of new variables.
\end{quote}

\sphinxAtStartPar
\sphinxcode{\sphinxupquote{colLower}} and \sphinxcode{\sphinxupquote{colUpper}}
\begin{quote}

\sphinxAtStartPar
Lower and upper bounds of new variables.

\sphinxAtStartPar
If \sphinxcode{\sphinxupquote{colLower}} is \sphinxcode{\sphinxupquote{NULL}}, lower bounds will be 0.

\sphinxAtStartPar
If \sphinxcode{\sphinxupquote{colUpper}} is \sphinxcode{\sphinxupquote{NULL}}, upper bounds will be \sphinxcode{\sphinxupquote{COPT\_INFINITY}}.
\end{quote}

\sphinxAtStartPar
\sphinxcode{\sphinxupquote{colNames}}
\begin{quote}

\sphinxAtStartPar
Names of new variables. Can be \sphinxcode{\sphinxupquote{NULL}}.
\end{quote}
\end{quote}
\end{quote}

\subsubsection{COPT\_AddPSDCols}
\label{\detokenize{capiref:copt-addpsdcols}}\begin{quote}

\sphinxAtStartPar
\sphinxstylestrong{Synopsis}
\begin{quote}

\sphinxAtStartPar
\sphinxcode{\sphinxupquote{int COPT\_AddPSDCols(copt\_prob *prob, int nAddCol, const int* colDims, char const *const *names)}}
\end{quote}

\sphinxAtStartPar
\sphinxstylestrong{Description}
\begin{quote}

\sphinxAtStartPar
Add \sphinxcode{\sphinxupquote{nAddCol}} PSD variables to the problem.
\end{quote}

\sphinxAtStartPar
\sphinxstylestrong{Arguments}
\begin{quote}

\sphinxAtStartPar
\sphinxcode{\sphinxupquote{prob}}
\begin{quote}

\sphinxAtStartPar
The COPT problem.
\end{quote}

\sphinxAtStartPar
\sphinxcode{\sphinxupquote{nAddCol}}
\begin{quote}

\sphinxAtStartPar
Number of new PSD variables.
\end{quote}

\sphinxAtStartPar
\sphinxcode{\sphinxupquote{colDims}}
\begin{quote}

\sphinxAtStartPar
Dimensions of new PSD variables.
\end{quote}

\sphinxAtStartPar
\sphinxcode{\sphinxupquote{names}}
\begin{quote}

\sphinxAtStartPar
Names of new PSD variables. Can be \sphinxcode{\sphinxupquote{NULL}}.
\end{quote}
\end{quote}
\end{quote}

\subsubsection{COPT\_AddRows}
\label{\detokenize{capiref:copt-addrows}}\begin{quote}

\sphinxAtStartPar
\sphinxstylestrong{Synopsis}
\begin{quote}

\sphinxAtStartPar
\sphinxcode{\sphinxupquote{int COPT\_AddRows(}}
\begin{quote}

\sphinxAtStartPar
\sphinxcode{\sphinxupquote{copt\_prob *prob,}}

\sphinxAtStartPar
\sphinxcode{\sphinxupquote{int nAddRow,}}

\sphinxAtStartPar
\sphinxcode{\sphinxupquote{const int *rowMatBeg,}}

\sphinxAtStartPar
\sphinxcode{\sphinxupquote{const int *rowMatCnt,}}

\sphinxAtStartPar
\sphinxcode{\sphinxupquote{const int *rowMatIdx,}}

\sphinxAtStartPar
\sphinxcode{\sphinxupquote{const double *rowMatElem,}}

\sphinxAtStartPar
\sphinxcode{\sphinxupquote{const char *rowSense,}}

\sphinxAtStartPar
\sphinxcode{\sphinxupquote{const double *rowBound,}}

\sphinxAtStartPar
\sphinxcode{\sphinxupquote{const double *rowUpper,}}

\sphinxAtStartPar
\sphinxcode{\sphinxupquote{char const *const *rowNames)}}
\end{quote}
\end{quote}

\sphinxAtStartPar
\sphinxstylestrong{Description}
\begin{quote}

\sphinxAtStartPar
Adds \sphinxcode{\sphinxupquote{nAddRow}} constraints (rows) to the problem.
\end{quote}

\sphinxAtStartPar
\sphinxstylestrong{Arguments}
\begin{quote}

\sphinxAtStartPar
\sphinxcode{\sphinxupquote{prob}}
\begin{quote}

\sphinxAtStartPar
The COPT problem.
\end{quote}

\sphinxAtStartPar
\sphinxcode{\sphinxupquote{nAddRow}}
\begin{quote}

\sphinxAtStartPar
Number of new constraints.
\end{quote}

\sphinxAtStartPar
\sphinxcode{\sphinxupquote{rowMatBeg, rowMatCnt, rowMatIdx}} and \sphinxcode{\sphinxupquote{rowMatElem}}
\begin{quote}

\sphinxAtStartPar
Defines the coefficient matrix in compressed row storage (CRS) format.
The CRS format is similar to the CCS format described
in the \sphinxstylestrong{other information} of \sphinxcode{\sphinxupquote{COPT\_LoadProb}}.
\end{quote}

\sphinxAtStartPar
\sphinxcode{\sphinxupquote{rowSense}}
\begin{quote}

\sphinxAtStartPar
Senses of new constraints.

\sphinxAtStartPar
Please refer to the list of all senses constants for all the supported types.

\sphinxAtStartPar
If \sphinxcode{\sphinxupquote{rowSense}} is \sphinxcode{\sphinxupquote{NULL}}, then \sphinxcode{\sphinxupquote{rowBound}} and \sphinxcode{\sphinxupquote{rowUpper}}
will be treated as lower and upper bounds for constraints.
This is the recommended method for defining constraints.

\sphinxAtStartPar
If \sphinxcode{\sphinxupquote{rowSense}} is provided, then \sphinxcode{\sphinxupquote{rowBound}} and \sphinxcode{\sphinxupquote{rowUpper}}
will be treated as RHS and \sphinxstylestrong{range} for constraints.
In this case, \sphinxcode{\sphinxupquote{rowUpper}} is only required when there
are \sphinxcode{\sphinxupquote{COPT\_RANGE}} constraints, where the
\begin{quote}

\sphinxAtStartPar
lower bound is \sphinxcode{\sphinxupquote{rowBound{[}i{]} \sphinxhyphen{} fabs(rowUpper{[}i{]})}}

\sphinxAtStartPar
upper bound is \sphinxcode{\sphinxupquote{rowBound{[}i{]}}}
\end{quote}
\end{quote}

\sphinxAtStartPar
\sphinxcode{\sphinxupquote{rowBound}}
\begin{quote}

\sphinxAtStartPar
Lower bounds or RHS of new constraints.
\end{quote}

\sphinxAtStartPar
\sphinxcode{\sphinxupquote{rowUpper}}
\begin{quote}

\sphinxAtStartPar
Upper bounds or \sphinxstylestrong{range} of new constraints.
\end{quote}

\sphinxAtStartPar
\sphinxcode{\sphinxupquote{rowNames}}
\begin{quote}

\sphinxAtStartPar
Names of new constraints. Can be \sphinxcode{\sphinxupquote{NULL}}.
\end{quote}
\end{quote}
\end{quote}

\subsubsection{COPT\_AddSOSs}
\label{\detokenize{capiref:copt-addsoss}}\begin{quote}

\sphinxAtStartPar
\sphinxstylestrong{Synopsis}
\begin{quote}

\sphinxAtStartPar
\sphinxcode{\sphinxupquote{int COPT\_AddSOSs(}}
\begin{quote}

\sphinxAtStartPar
\sphinxcode{\sphinxupquote{copt\_prob *prob,}}

\sphinxAtStartPar
\sphinxcode{\sphinxupquote{int nAddSOS,}}

\sphinxAtStartPar
\sphinxcode{\sphinxupquote{const int *sosType,}}

\sphinxAtStartPar
\sphinxcode{\sphinxupquote{const int *sosMatBeg,}}

\sphinxAtStartPar
\sphinxcode{\sphinxupquote{const int *sosMatCnt,}}

\sphinxAtStartPar
\sphinxcode{\sphinxupquote{const int *sosMatIdx,}}

\sphinxAtStartPar
\sphinxcode{\sphinxupquote{const double *sosMatWt)}}
\end{quote}
\end{quote}

\sphinxAtStartPar
\sphinxstylestrong{Description}
\begin{quote}

\sphinxAtStartPar
Add \sphinxcode{\sphinxupquote{nAddSOS}} SOS constraints to the problem.
If \sphinxcode{\sphinxupquote{sosMatWt}} is \sphinxcode{\sphinxupquote{NULL}}, then COPT will generate it internally.

\sphinxAtStartPar
\sphinxstylestrong{Note:} if a problem contains SOS constraints, the problem is a MIP.
\end{quote}

\sphinxAtStartPar
\sphinxstylestrong{Arguments}
\begin{quote}

\sphinxAtStartPar
\sphinxcode{\sphinxupquote{prob}}
\begin{quote}

\sphinxAtStartPar
The COPT problem.
\end{quote}

\sphinxAtStartPar
\sphinxcode{\sphinxupquote{nAddSOS}}
\begin{quote}

\sphinxAtStartPar
Number of new SOS constraints.
\end{quote}

\sphinxAtStartPar
\sphinxcode{\sphinxupquote{sosType}}
\begin{quote}

\sphinxAtStartPar
Types of SOS constraints.
\end{quote}

\sphinxAtStartPar
\sphinxcode{\sphinxupquote{sosMatBeg, sosMatCnt, sosMatIdx}} and \sphinxcode{\sphinxupquote{sosMatWt}}
\begin{quote}

\sphinxAtStartPar
Defines the coefficient matrix in compressed row storage (CRS) format.
The CRS format is similar to the CCS format described
in the \sphinxstylestrong{other information} of \sphinxcode{\sphinxupquote{COPT\_LoadProb}}.
\end{quote}

\sphinxAtStartPar
\sphinxcode{\sphinxupquote{sosMatWt}}
\begin{quote}

\sphinxAtStartPar
Weights of variables in SOS constraints. Can be \sphinxcode{\sphinxupquote{NULL}}.
\end{quote}
\end{quote}
\end{quote}

\subsubsection{COPT\_AddCones}
\label{\detokenize{capiref:copt-addcones}}\begin{quote}

\sphinxAtStartPar
\sphinxstylestrong{Synopsis}
\begin{quote}

\sphinxAtStartPar
\sphinxcode{\sphinxupquote{int COPT\_AddCones(}}
\begin{quote}

\sphinxAtStartPar
\sphinxcode{\sphinxupquote{copt\_prob *prob,}}

\sphinxAtStartPar
\sphinxcode{\sphinxupquote{int nAddCone,}}

\sphinxAtStartPar
\sphinxcode{\sphinxupquote{const int *coneType,}}

\sphinxAtStartPar
\sphinxcode{\sphinxupquote{const int *coneBeg,}}

\sphinxAtStartPar
\sphinxcode{\sphinxupquote{const int *coneCnt,}}

\sphinxAtStartPar
\sphinxcode{\sphinxupquote{const int *coneIdx)}}
\end{quote}
\end{quote}

\sphinxAtStartPar
\sphinxstylestrong{Description}
\begin{quote}

\sphinxAtStartPar
Add \sphinxcode{\sphinxupquote{nAddCone}} Second\sphinxhyphen{}Order\sphinxhyphen{}Cone (SOC) constraints.
\end{quote}

\sphinxAtStartPar
\sphinxstylestrong{Arguments}
\begin{quote}

\sphinxAtStartPar
\sphinxcode{\sphinxupquote{prob}}
\begin{quote}

\sphinxAtStartPar
The COPT problem.
\end{quote}

\sphinxAtStartPar
\sphinxcode{\sphinxupquote{nAddCone}}
\begin{quote}

\sphinxAtStartPar
Number of new SOC constraints.
\end{quote}

\sphinxAtStartPar
\sphinxcode{\sphinxupquote{coneType}}
\begin{quote}

\sphinxAtStartPar
Types of SOC constraints.
\end{quote}

\sphinxAtStartPar
\sphinxcode{\sphinxupquote{coneBeg, coneCnt, coneIdx}}
\begin{quote}

\sphinxAtStartPar
Defines the coefficient matrix in compressed row storage (CRS) format.
The CRS format is similar to the CCS format described
in the \sphinxstylestrong{other information} of \sphinxcode{\sphinxupquote{COPT\_LoadProb}}.
\end{quote}
\end{quote}
\end{quote}

\subsubsection{COPT\_AddExpCones}
\label{\detokenize{capiref:copt-addexpcones}}\begin{quote}

\sphinxAtStartPar
\sphinxstylestrong{Synopsis}
\begin{quote}

\sphinxAtStartPar
\sphinxcode{\sphinxupquote{int COPT\_AddExpCones(}}
\begin{quote}

\sphinxAtStartPar
\sphinxcode{\sphinxupquote{copt\_prob *prob,}}

\sphinxAtStartPar
\sphinxcode{\sphinxupquote{int nAddCone,}}

\sphinxAtStartPar
\sphinxcode{\sphinxupquote{const int *coneType,}}

\sphinxAtStartPar
\sphinxcode{\sphinxupquote{const int *coneIdx)}}
\end{quote}
\end{quote}

\sphinxAtStartPar
\sphinxstylestrong{Description}
\begin{quote}

\sphinxAtStartPar
Add \sphinxcode{\sphinxupquote{nAddCone}} exponential cone constraints.
\end{quote}

\sphinxAtStartPar
\sphinxstylestrong{Arguments}
\begin{quote}

\sphinxAtStartPar
\sphinxcode{\sphinxupquote{prob}}
\begin{quote}

\sphinxAtStartPar
The COPT problem.
\end{quote}

\sphinxAtStartPar
\sphinxcode{\sphinxupquote{nAddCone}}
\begin{quote}

\sphinxAtStartPar
Number of new exponential cone constraints.
\end{quote}

\sphinxAtStartPar
\sphinxcode{\sphinxupquote{coneType}}
\begin{quote}

\sphinxAtStartPar
Types of exponential cone constraints.
\end{quote}

\sphinxAtStartPar
\sphinxcode{\sphinxupquote{coneIdx}}
\begin{quote}

\sphinxAtStartPar
Array of subscripts for the variables that constitute the exponential cone constraints.
\end{quote}
\end{quote}
\end{quote}

\subsubsection{COPT\_AddAffineCone}
\label{\detokenize{capiref:copt-addaffinecone}}\begin{quote}

\sphinxAtStartPar
\sphinxstylestrong{Synopsis}
\begin{quote}

\sphinxAtStartPar
\sphinxcode{\sphinxupquote{int COPT\_AddAffineCone(}}
\begin{quote}

\sphinxAtStartPar
\sphinxcode{\sphinxupquote{copt\_prob *prob,}}

\sphinxAtStartPar
\sphinxcode{\sphinxupquote{int coneType,}}

\sphinxAtStartPar
\sphinxcode{\sphinxupquote{int nConeDim,}}

\sphinxAtStartPar
\sphinxcode{\sphinxupquote{int nAlphaDim,}}

\sphinxAtStartPar
\sphinxcode{\sphinxupquote{const double *alphaElem,}}

\sphinxAtStartPar
\sphinxcode{\sphinxupquote{const int *psdBeg,}}

\sphinxAtStartPar
\sphinxcode{\sphinxupquote{const int *psdCnt,}}

\sphinxAtStartPar
\sphinxcode{\sphinxupquote{const int *psdColIdx,}}

\sphinxAtStartPar
\sphinxcode{\sphinxupquote{const int *psdMatIdx,}}

\sphinxAtStartPar
\sphinxcode{\sphinxupquote{const int *rowMatBeg,}}

\sphinxAtStartPar
\sphinxcode{\sphinxupquote{const int *rowMatCnt,}}

\sphinxAtStartPar
\sphinxcode{\sphinxupquote{const int *rowMatIdx,}}

\sphinxAtStartPar
\sphinxcode{\sphinxupquote{const double *rowMatElem,}}

\sphinxAtStartPar
\sphinxcode{\sphinxupquote{const double *rowConst,}}

\sphinxAtStartPar
\sphinxcode{\sphinxupquote{const char *name)}}
\end{quote}
\end{quote}

\sphinxAtStartPar
\sphinxstylestrong{Description}
\begin{quote}

\sphinxAtStartPar
Add the affine cone constraints.
\end{quote}

\sphinxAtStartPar
\sphinxstylestrong{Arguments}
\begin{quote}

\sphinxAtStartPar
\sphinxcode{\sphinxupquote{prob}}
\begin{quote}

\sphinxAtStartPar
The COPT model.
\end{quote}

\sphinxAtStartPar
\sphinxcode{\sphinxupquote{coneType}}
\begin{quote}

\sphinxAtStartPar
The type of the affine cone. Please refer to {\hyperref[\detokenize{constant:chapconst-conetype}]{\sphinxcrossref{\DUrole{std,std-ref}{SOC constraint types}}}} and {\hyperref[\detokenize{constant:chapconst-expconetype}]{\sphinxcrossref{\DUrole{std,std-ref}{Exponential Cone type}}}}
for possible values.
\end{quote}

\sphinxAtStartPar
\sphinxcode{\sphinxupquote{nConeDim}}
\begin{quote}

\sphinxAtStartPar
The dimension of the affine cone.
\end{quote}

\sphinxAtStartPar
\sphinxcode{\sphinxupquote{nAlphaDim}}
\begin{quote}

\sphinxAtStartPar
Reserved parameter, currently not in use.
\end{quote}

\sphinxAtStartPar
\sphinxcode{\sphinxupquote{alphaElem}}
\begin{quote}

\sphinxAtStartPar
Reserved parameter, currently not in use.
\end{quote}

\sphinxAtStartPar
\sphinxcode{\sphinxupquote{psdBeg, psdCnt, psdColIdx, psdMatIdx}}
\begin{quote}

\sphinxAtStartPar
Represents the PSD terms in the affine cone.

\sphinxAtStartPar
\sphinxcode{\sphinxupquote{psdBeg}} indicates the starting position of the PSD terms in each affine cone term,
\sphinxcode{\sphinxupquote{psdCnt}} specifies the number of PSD terms, \sphinxcode{\sphinxupquote{psdColIdx}} indicates the index of the
PSD variable, and \sphinxcode{\sphinxupquote{psdMatIdx}} refers to the index of the symmetric matrix.
\end{quote}

\sphinxAtStartPar
\sphinxcode{\sphinxupquote{rowMatBeg, rowMatCnt, rowMatIdx, rowMatElem}}
\begin{quote}

\sphinxAtStartPar
Represents the linear terms in the affine cone.

\sphinxAtStartPar
The coefficient matrix is provided in CRS\sphinxhyphen{}format. For detailed examples, please refer to \sphinxstylestrong{Additional Information} in \sphinxcode{\sphinxupquote{COPT\_LoadProb}}.
\end{quote}

\sphinxAtStartPar
\sphinxcode{\sphinxupquote{rowConst}}
\begin{quote}

\sphinxAtStartPar
The constant terms in the affine cone.
\end{quote}

\sphinxAtStartPar
\sphinxcode{\sphinxupquote{name}}
\begin{quote}

\sphinxAtStartPar
The name of the affine cone.
\end{quote}
\end{quote}
\end{quote}

\subsubsection{COPT\_AddQConstr}
\label{\detokenize{capiref:copt-addqconstr}}\begin{quote}

\sphinxAtStartPar
\sphinxstylestrong{Synopsis}
\begin{quote}

\sphinxAtStartPar
\sphinxcode{\sphinxupquote{int COPT\_AddQConstr(}}
\begin{quote}

\sphinxAtStartPar
\sphinxcode{\sphinxupquote{copt\_prob *prob,}}

\sphinxAtStartPar
\sphinxcode{\sphinxupquote{int nRowMatCnt,}}

\sphinxAtStartPar
\sphinxcode{\sphinxupquote{const int *rowMatIdx,}}

\sphinxAtStartPar
\sphinxcode{\sphinxupquote{const int *rowMatElem,}}

\sphinxAtStartPar
\sphinxcode{\sphinxupquote{int nQMatCnt,}}

\sphinxAtStartPar
\sphinxcode{\sphinxupquote{const int *qMatRow,}}

\sphinxAtStartPar
\sphinxcode{\sphinxupquote{const int *qMatCol,}}

\sphinxAtStartPar
\sphinxcode{\sphinxupquote{const double *qMatElem,}}

\sphinxAtStartPar
\sphinxcode{\sphinxupquote{char cRowSense,}}

\sphinxAtStartPar
\sphinxcode{\sphinxupquote{double dRowBound, const char *name)}}
\end{quote}
\end{quote}

\sphinxAtStartPar
\sphinxstylestrong{Description}
\begin{quote}

\sphinxAtStartPar
Add a general quadratic constraint.

\sphinxAtStartPar
\sphinxstylestrong{Note} Only convex quadratic constraint is currently supported.
\end{quote}

\sphinxAtStartPar
\sphinxstylestrong{Arguments}
\begin{quote}

\sphinxAtStartPar
\sphinxcode{\sphinxupquote{prob}}
\begin{quote}

\sphinxAtStartPar
The COPT problem.
\end{quote}

\sphinxAtStartPar
\sphinxcode{\sphinxupquote{nRowMatCnt}}
\begin{quote}

\sphinxAtStartPar
Number of non\sphinxhyphen{}zero linear terms of the quadratic constraint (row).
\end{quote}

\sphinxAtStartPar
\sphinxcode{\sphinxupquote{rowMatIdx}}
\begin{quote}

\sphinxAtStartPar
Column index of non\sphinxhyphen{}zero linear terms of the quadratic constraint (row).
\end{quote}

\sphinxAtStartPar
\sphinxcode{\sphinxupquote{rowMatElem}}
\begin{quote}

\sphinxAtStartPar
Values of non\sphinxhyphen{}zero linear terms of the quadratic constraint (row).
\end{quote}

\sphinxAtStartPar
\sphinxcode{\sphinxupquote{nQMatCnt}}
\begin{quote}

\sphinxAtStartPar
Number of non\sphinxhyphen{}zero quadratic terms of the quadratic constraint (row).
\end{quote}

\sphinxAtStartPar
\sphinxcode{\sphinxupquote{qMatRow}}
\begin{quote}

\sphinxAtStartPar
Row index of non\sphinxhyphen{}zero quadratic terms of the quadratic constraint (row).
\end{quote}

\sphinxAtStartPar
\sphinxcode{\sphinxupquote{qMatCol}}
\begin{quote}

\sphinxAtStartPar
Column index of non\sphinxhyphen{}zero quadratic terms of the quadratic constraint (row).
\end{quote}

\sphinxAtStartPar
\sphinxcode{\sphinxupquote{qMatElem}}
\begin{quote}

\sphinxAtStartPar
Values of non\sphinxhyphen{}zero quadratic terms of the quadratic constraint (row).
\end{quote}

\sphinxAtStartPar
\sphinxcode{\sphinxupquote{cRowSense}}
\begin{quote}

\sphinxAtStartPar
The sense of the quadratic constraint (row).
\end{quote}

\sphinxAtStartPar
\sphinxcode{\sphinxupquote{dRowBound}}
\begin{quote}

\sphinxAtStartPar
Right hand side of the quadratic constraint (row).
\end{quote}

\sphinxAtStartPar
\sphinxcode{\sphinxupquote{name}}
\begin{quote}

\sphinxAtStartPar
Name of the quadratic constraint (row).
\end{quote}
\end{quote}
\end{quote}

\subsubsection{COPT\_AddPSDConstr}
\label{\detokenize{capiref:copt-addpsdconstr}}\begin{quote}

\sphinxAtStartPar
\sphinxstylestrong{Synopsis}
\begin{quote}

\sphinxAtStartPar
\sphinxcode{\sphinxupquote{int COPT\_AddPSDConstr(}}
\begin{quote}

\sphinxAtStartPar
\sphinxcode{\sphinxupquote{copt\_prob *prob,}}

\sphinxAtStartPar
\sphinxcode{\sphinxupquote{int nRowMatCnt,}}

\sphinxAtStartPar
\sphinxcode{\sphinxupquote{const int *rowMatIdx,}}

\sphinxAtStartPar
\sphinxcode{\sphinxupquote{const int *rowMatElem,}}

\sphinxAtStartPar
\sphinxcode{\sphinxupquote{int nColCnt,}}

\sphinxAtStartPar
\sphinxcode{\sphinxupquote{const int *psdColIdx,}}

\sphinxAtStartPar
\sphinxcode{\sphinxupquote{const int *symMatIdx,}}

\sphinxAtStartPar
\sphinxcode{\sphinxupquote{char cRowSense,}}

\sphinxAtStartPar
\sphinxcode{\sphinxupquote{double dRowBound,}}

\sphinxAtStartPar
\sphinxcode{\sphinxupquote{double dRowUpper,}}

\sphinxAtStartPar
\sphinxcode{\sphinxupquote{const char *name)}}
\end{quote}
\end{quote}

\sphinxAtStartPar
\sphinxstylestrong{Description}
\begin{quote}

\sphinxAtStartPar
Add a PSD constraint.
\end{quote}

\sphinxAtStartPar
\sphinxstylestrong{Arguments}
\begin{quote}

\sphinxAtStartPar
\sphinxcode{\sphinxupquote{prob}}
\begin{quote}

\sphinxAtStartPar
The COPT problem.
\end{quote}

\sphinxAtStartPar
\sphinxcode{\sphinxupquote{nRowMatCnt}}
\begin{quote}

\sphinxAtStartPar
Number of non\sphinxhyphen{}zero linear terms of the PSD constraint.
\end{quote}

\sphinxAtStartPar
\sphinxcode{\sphinxupquote{rowMatIdx}}
\begin{quote}

\sphinxAtStartPar
Column index of non\sphinxhyphen{}zero linear terms of the PSD constraint.
\end{quote}

\sphinxAtStartPar
\sphinxcode{\sphinxupquote{rowMatElem}}
\begin{quote}

\sphinxAtStartPar
Values of non\sphinxhyphen{}zero linear terms of the PSD constraint.
\end{quote}

\sphinxAtStartPar
\sphinxcode{\sphinxupquote{nColCnt}}
\begin{quote}

\sphinxAtStartPar
Number of PSD terms of the PSD constraint.
\end{quote}

\sphinxAtStartPar
\sphinxcode{\sphinxupquote{psdColIdx}}
\begin{quote}

\sphinxAtStartPar
PSD variable index of PSD terms of the PSD constraint.
\end{quote}

\sphinxAtStartPar
\sphinxcode{\sphinxupquote{symMatIdx}}
\begin{quote}

\sphinxAtStartPar
Symmetric matrix index of PSD terms of the PSD constraint.
\end{quote}

\sphinxAtStartPar
\sphinxcode{\sphinxupquote{cRowSense}}
\begin{quote}

\sphinxAtStartPar
Senses of new PSD constraint.

\sphinxAtStartPar
Please refer to the list of all senses constants for all the supported types.

\sphinxAtStartPar
If \sphinxcode{\sphinxupquote{cRowSense}} is 0, then \sphinxcode{\sphinxupquote{dRowBound}} and \sphinxcode{\sphinxupquote{dRowUpper}}
will be treated as lower and upper bounds for the constraint.
This is the recommended method for defining constraints.

\sphinxAtStartPar
If \sphinxcode{\sphinxupquote{cRowSense}} is provided, then \sphinxcode{\sphinxupquote{dRowBound}} and \sphinxcode{\sphinxupquote{dRowUpper}}
will be treated as RHS and \sphinxstylestrong{range} for the constraint.
In this case, \sphinxcode{\sphinxupquote{dRowUpper}} is only required when
\sphinxcode{\sphinxupquote{cRowSense = COPT\_RANGE}}, where
\begin{quote}

\sphinxAtStartPar
lower bound is \sphinxcode{\sphinxupquote{dRowBound \sphinxhyphen{} dRowUpper}}

\sphinxAtStartPar
upper bound is \sphinxcode{\sphinxupquote{dRowBound}}
\end{quote}
\end{quote}

\sphinxAtStartPar
\sphinxcode{\sphinxupquote{dRowBound}}
\begin{quote}

\sphinxAtStartPar
Lower bound or RHS of the PSD constraint.
\end{quote}

\sphinxAtStartPar
\sphinxcode{\sphinxupquote{dRowUpper}}
\begin{quote}

\sphinxAtStartPar
Upper bound or \sphinxstylestrong{range} of the PSD constraint.
\end{quote}

\sphinxAtStartPar
\sphinxcode{\sphinxupquote{name}}
\begin{quote}

\sphinxAtStartPar
Name of the PSD constraint. Can be \sphinxcode{\sphinxupquote{NULL}}.
\end{quote}
\end{quote}
\end{quote}

\subsubsection{COPT\_AddLMIConstr}
\label{\detokenize{capiref:copt-addlmiconstr}}\begin{quote}

\sphinxAtStartPar
\sphinxstylestrong{Synopsis}
\begin{quote}

\sphinxAtStartPar
\sphinxcode{\sphinxupquote{int COPT\_AddLMIConstr(}}
\begin{quote}

\sphinxAtStartPar
\sphinxcode{\sphinxupquote{copt\_prob *prob,}}

\sphinxAtStartPar
\sphinxcode{\sphinxupquote{int nDim,}}

\sphinxAtStartPar
\sphinxcode{\sphinxupquote{int nLMIMatCnt,}}

\sphinxAtStartPar
\sphinxcode{\sphinxupquote{const int *colIdx,}}

\sphinxAtStartPar
\sphinxcode{\sphinxupquote{const int *symMatIdx,}}

\sphinxAtStartPar
\sphinxcode{\sphinxupquote{int constMatIdx,}}

\sphinxAtStartPar
\sphinxcode{\sphinxupquote{const char *name)}}
\end{quote}
\end{quote}

\sphinxAtStartPar
\sphinxstylestrong{Description}
\begin{quote}

\sphinxAtStartPar
Add a LMI constraint to the problem.
\end{quote}

\sphinxAtStartPar
\sphinxstylestrong{Arguments}
\begin{quote}

\sphinxAtStartPar
\sphinxcode{\sphinxupquote{prob}}
\begin{quote}

\sphinxAtStartPar
The COPT problem.
\end{quote}

\sphinxAtStartPar
\sphinxcode{\sphinxupquote{nDim}}
\begin{quote}

\sphinxAtStartPar
Dimension of symmetric matrix in the LMI constraint.
\end{quote}

\sphinxAtStartPar
\sphinxcode{\sphinxupquote{nLMIMatCnt}}
\begin{quote}

\sphinxAtStartPar
Number of coefficient matrix entries in the LMI constraint.
\end{quote}

\sphinxAtStartPar
\sphinxcode{\sphinxupquote{colIdx}}
\begin{quote}

\sphinxAtStartPar
Index of scalar variable in the LMI constraint.
\end{quote}

\sphinxAtStartPar
\sphinxcode{\sphinxupquote{symMatIdx}}
\begin{quote}

\sphinxAtStartPar
Index of symmetric coefficient matrix in the LMI constraint.
\end{quote}

\sphinxAtStartPar
\sphinxcode{\sphinxupquote{constMatIdx}}
\begin{quote}

\sphinxAtStartPar
Index of constant\sphinxhyphen{}term symmetric matrix in the LMI constraint.
\end{quote}

\sphinxAtStartPar
\sphinxcode{\sphinxupquote{name}}
\begin{quote}

\sphinxAtStartPar
Name of LMI constraint. Can be \sphinxcode{\sphinxupquote{NULL}}.
\end{quote}
\end{quote}
\end{quote}

\subsubsection{COPT\_AddIndicator}
\label{\detokenize{capiref:copt-addindicator}}\begin{quote}

\sphinxAtStartPar
\sphinxstylestrong{Synopsis}
\begin{quote}

\sphinxAtStartPar
\sphinxcode{\sphinxupquote{int COPT\_AddIndicator(}}
\begin{quote}

\sphinxAtStartPar
\sphinxcode{\sphinxupquote{copt\_prob *prob,}}

\sphinxAtStartPar
\sphinxcode{\sphinxupquote{int binColIdx,}}

\sphinxAtStartPar
\sphinxcode{\sphinxupquote{int binColVal,}}

\sphinxAtStartPar
\sphinxcode{\sphinxupquote{int nRowMatCnt,}}

\sphinxAtStartPar
\sphinxcode{\sphinxupquote{const int *rowMatIdx,}}

\sphinxAtStartPar
\sphinxcode{\sphinxupquote{const double *rowMatElem,}}

\sphinxAtStartPar
\sphinxcode{\sphinxupquote{char cRowSense, double dRowBound)}}
\end{quote}
\end{quote}

\sphinxAtStartPar
\sphinxstylestrong{Description}
\begin{quote}

\sphinxAtStartPar
Add an indicator constraint to the problem.
\end{quote}

\sphinxAtStartPar
\sphinxstylestrong{Arguments}
\begin{quote}

\sphinxAtStartPar
\sphinxcode{\sphinxupquote{prob}}
\begin{quote}

\sphinxAtStartPar
The COPT problem.
\end{quote}

\sphinxAtStartPar
\sphinxcode{\sphinxupquote{binColIdx}}
\begin{quote}

\sphinxAtStartPar
Index of indicator variable (column).
\end{quote}

\sphinxAtStartPar
\sphinxcode{\sphinxupquote{binColVal}}
\begin{quote}

\sphinxAtStartPar
Value of indicator variable (column).
\end{quote}

\sphinxAtStartPar
\sphinxcode{\sphinxupquote{nRowMatCnt}}
\begin{quote}

\sphinxAtStartPar
Number of non\sphinxhyphen{}zero elements in the linear constraint (row).
\end{quote}

\sphinxAtStartPar
\sphinxcode{\sphinxupquote{rowMatIdx}}
\begin{quote}

\sphinxAtStartPar
Column index of non\sphinxhyphen{}zero elements in the linear constraint (row).
\end{quote}

\sphinxAtStartPar
\sphinxcode{\sphinxupquote{rowMatElem}}
\begin{quote}

\sphinxAtStartPar
Values of non\sphinxhyphen{}zero elements in the linear constraint (row).
\end{quote}

\sphinxAtStartPar
\sphinxcode{\sphinxupquote{cRowSense}}
\begin{quote}

\sphinxAtStartPar
The sense of the linear constraint (row). Options are: \sphinxcode{\sphinxupquote{COPT\_EQUAL}} ,
\sphinxcode{\sphinxupquote{COPT\_LESS\_EQUAL}} and \sphinxcode{\sphinxupquote{COPT\_GREATER\_EQUAL}} .
\end{quote}

\sphinxAtStartPar
\sphinxcode{\sphinxupquote{dRowBound}}
\begin{quote}

\sphinxAtStartPar
Right hand side of the linear constraint (row).
\end{quote}
\end{quote}
\end{quote}

\subsubsection{COPT\_AddIndicators}
\label{\detokenize{capiref:copt-addindicators}}\begin{quote}

\sphinxAtStartPar
\sphinxstylestrong{Synopsis}
\begin{quote}

\sphinxAtStartPar
\sphinxcode{\sphinxupquote{int COPT\_AddIndicators(}}
\begin{quote}

\sphinxAtStartPar
\sphinxcode{\sphinxupquote{copt\_prob *prob,}}

\sphinxAtStartPar
\sphinxcode{\sphinxupquote{int nInd,}}

\sphinxAtStartPar
\sphinxcode{\sphinxupquote{int *indType,}}

\sphinxAtStartPar
\sphinxcode{\sphinxupquote{int *binColIdx,}}

\sphinxAtStartPar
\sphinxcode{\sphinxupquote{int *binColVal,}}

\sphinxAtStartPar
\sphinxcode{\sphinxupquote{const int *rowMatBeg,}}

\sphinxAtStartPar
\sphinxcode{\sphinxupquote{const int *RowMatCnt,}}

\sphinxAtStartPar
\sphinxcode{\sphinxupquote{const int *rowMatIdx,}}

\sphinxAtStartPar
\sphinxcode{\sphinxupquote{const double *rowMatElem,}}

\sphinxAtStartPar
\sphinxcode{\sphinxupquote{char cRowSense,}}

\sphinxAtStartPar
\sphinxcode{\sphinxupquote{double *dRowBound,}}

\sphinxAtStartPar
\sphinxcode{\sphinxupquote{char const *const indNames);}}
\end{quote}
\end{quote}

\sphinxAtStartPar
\sphinxstylestrong{Description}
\begin{quote}

\sphinxAtStartPar
Add  \sphinxcode{\sphinxupquote{nInd}} indicator constraints to the problem.
\end{quote}

\sphinxAtStartPar
\sphinxstylestrong{Arguments}
\begin{quote}

\sphinxAtStartPar
\sphinxcode{\sphinxupquote{prob}}
\begin{quote}

\sphinxAtStartPar
The COPT problem.
\end{quote}

\sphinxAtStartPar
\sphinxcode{\sphinxupquote{nInd}}
\begin{quote}

\sphinxAtStartPar
Number of indicator constraints.
\end{quote}

\sphinxAtStartPar
\sphinxcode{\sphinxupquote{indType}}
\begin{quote}

\sphinxAtStartPar
Type for indcator constraints. Please refer to {\hyperref[\detokenize{constant:chapconst-indicatortype}]{\sphinxcrossref{\DUrole{std,std-ref}{Indicator constraint types}}}}
for possible values.
\end{quote}

\sphinxAtStartPar
\sphinxcode{\sphinxupquote{binColIdx}}
\begin{quote}

\sphinxAtStartPar
Index of indicator variable (column).
\end{quote}

\sphinxAtStartPar
\sphinxcode{\sphinxupquote{binColVal}}
\begin{quote}

\sphinxAtStartPar
Value of indicator variable (column).
\end{quote}

\sphinxAtStartPar
\sphinxcode{\sphinxupquote{rowMatBeg, rowMatCnt, rowMatIdx}} and \sphinxcode{\sphinxupquote{rowMatElem}}
\begin{quote}

\sphinxAtStartPar
Defines the coefficient matrix in compressed row storage (CRS) format
for linear constraint in the indcator constrsints.

\sphinxAtStartPar
Please see \sphinxstylestrong{other information} of \sphinxcode{\sphinxupquote{COPT\_LoadProb}} for
an example of the CRS format.
\end{quote}

\sphinxAtStartPar
\sphinxcode{\sphinxupquote{cRowSense}}
\begin{quote}

\sphinxAtStartPar
The sense of the linear constraint (row). Options are: \sphinxcode{\sphinxupquote{COPT\_EQUAL}} ,
\sphinxcode{\sphinxupquote{COPT\_LESS\_EQUAL}} and \sphinxcode{\sphinxupquote{COPT\_GREATER\_EQUAL}} .
\end{quote}

\sphinxAtStartPar
\sphinxcode{\sphinxupquote{dRowBound}}
\begin{quote}

\sphinxAtStartPar
Right hand side of the linear constraint (row).
\end{quote}

\sphinxAtStartPar
\sphinxcode{\sphinxupquote{indNames}}
\begin{quote}

\sphinxAtStartPar
Names for indicator constraints.
\end{quote}
\end{quote}
\end{quote}

\subsubsection{COPT\_AddSymMat}
\label{\detokenize{capiref:copt-addsymmat}}\begin{quote}

\sphinxAtStartPar
\sphinxstylestrong{Synopsis}
\begin{quote}

\sphinxAtStartPar
\sphinxcode{\sphinxupquote{int COPT\_AddSymMat(copt\_prob *prob, int ndim, int nelem, int *rows, int *cols, double *elems)}}
\end{quote}

\sphinxAtStartPar
\sphinxstylestrong{Description}
\begin{quote}

\sphinxAtStartPar
Add a symmetric matrix to the problem. (Expect lower triangle part)
\end{quote}

\sphinxAtStartPar
\sphinxstylestrong{Arguments}
\begin{quote}

\sphinxAtStartPar
\sphinxcode{\sphinxupquote{prob}}
\begin{quote}

\sphinxAtStartPar
The COPT problem.
\end{quote}

\sphinxAtStartPar
\sphinxcode{\sphinxupquote{ndim}}
\begin{quote}

\sphinxAtStartPar
Dimension of symmetric matrix.
\end{quote}

\sphinxAtStartPar
\sphinxcode{\sphinxupquote{nelem}}
\begin{quote}

\sphinxAtStartPar
Number of non\sphinxhyphen{}zeros of symmetric matrix.
\end{quote}

\sphinxAtStartPar
\sphinxcode{\sphinxupquote{rows}}
\begin{quote}

\sphinxAtStartPar
Row index of symmetric matrix.
\end{quote}

\sphinxAtStartPar
\sphinxcode{\sphinxupquote{cols}}
\begin{quote}

\sphinxAtStartPar
Column index of symmetric matrix.
\end{quote}

\sphinxAtStartPar
\sphinxcode{\sphinxupquote{elems}}
\begin{quote}

\sphinxAtStartPar
Nonzero elements of symmetric matrix.
\end{quote}
\end{quote}
\end{quote}

\subsubsection{COPT\_AddNLConstr}
\label{\detokenize{capiref:copt-addnlconstr}}\begin{quote}

\sphinxAtStartPar
\sphinxstylestrong{Synopsis}
\begin{quote}

\sphinxAtStartPar
\sphinxcode{\sphinxupquote{int COPT\_AddNLConstr(}}
\begin{quote}

\sphinxAtStartPar
\sphinxcode{\sphinxupquote{copt\_prob *prob,}}

\sphinxAtStartPar
\sphinxcode{\sphinxupquote{int nToken,}}

\sphinxAtStartPar
\sphinxcode{\sphinxupquote{int nTokenElem,}}

\sphinxAtStartPar
\sphinxcode{\sphinxupquote{const int *token,}}

\sphinxAtStartPar
\sphinxcode{\sphinxupquote{const double *tokenElem,}}

\sphinxAtStartPar
\sphinxcode{\sphinxupquote{int nRowMatCnt,}}

\sphinxAtStartPar
\sphinxcode{\sphinxupquote{const int *rowMatIdx,}}

\sphinxAtStartPar
\sphinxcode{\sphinxupquote{const double *rowMatElem,}}

\sphinxAtStartPar
\sphinxcode{\sphinxupquote{char cRowSense,}}

\sphinxAtStartPar
\sphinxcode{\sphinxupquote{double dRowBound,}}

\sphinxAtStartPar
\sphinxcode{\sphinxupquote{double dRowUpper,}}

\sphinxAtStartPar
\sphinxcode{\sphinxupquote{const char *name)}}
\end{quote}
\end{quote}

\sphinxAtStartPar
\sphinxstylestrong{Description}
\begin{quote}

\sphinxAtStartPar
Add a nonlinear expression constraint.
Currently, only prefix notation is supported.
\end{quote}

\sphinxAtStartPar
\sphinxstylestrong{Arguments}
\begin{quote}

\sphinxAtStartPar
\sphinxcode{\sphinxupquote{prob}}
\begin{quote}

\sphinxAtStartPar
The COPT problem.
\end{quote}

\sphinxAtStartPar
\sphinxcode{\sphinxupquote{nToken}}
\begin{quote}

\sphinxAtStartPar
Number of tokens in the expression.
\end{quote}

\sphinxAtStartPar
\sphinxcode{\sphinxupquote{nTokenElem}}
\begin{quote}

\sphinxAtStartPar
Number of constants in the expression.
\end{quote}

\sphinxAtStartPar
\sphinxcode{\sphinxupquote{token}}
\begin{quote}

\sphinxAtStartPar
Array of expression tokens.
\end{quote}

\sphinxAtStartPar
\sphinxcode{\sphinxupquote{tokenElem}}
\begin{quote}

\sphinxAtStartPar
Array of constants in the expression.
\end{quote}

\sphinxAtStartPar
\sphinxcode{\sphinxupquote{nRowMatCnt}}
\begin{quote}

\sphinxAtStartPar
Number of linear terms in the constraint.
\end{quote}

\sphinxAtStartPar
\sphinxcode{\sphinxupquote{rowMatIdx}}
\begin{quote}

\sphinxAtStartPar
Indices of linear terms in the constraint.
\end{quote}

\sphinxAtStartPar
\sphinxcode{\sphinxupquote{rowMatElem}}
\begin{quote}

\sphinxAtStartPar
Coefficients of linear terms in the constraint.
\end{quote}

\sphinxAtStartPar
\sphinxcode{\sphinxupquote{cRowSense}}
\begin{quote}

\sphinxAtStartPar
Constraint type.

\sphinxAtStartPar
Refer to the constants section for the supported constraint types in COPT.

\sphinxAtStartPar
If \sphinxcode{\sphinxupquote{cRowSense}} is 0, then \sphinxcode{\sphinxupquote{dRowBound}} and \sphinxcode{\sphinxupquote{dRowUpper}}
are treated as the lower and upper bounds of the constraint.
This is the recommended way to define constraints.

\sphinxAtStartPar
If \sphinxcode{\sphinxupquote{cRowSense}} is a meaningful value, then \sphinxcode{\sphinxupquote{dRowBound}} and \sphinxcode{\sphinxupquote{dRowUpper}}
are treated as the RHS and \sphinxstylestrong{range}.
In this case, \sphinxcode{\sphinxupquote{dRowUpper}} is used only when the constraint type is \sphinxcode{\sphinxupquote{COPT\_RANGE}}.
The constraint bounds in this scenario are:
\begin{quote}

\sphinxAtStartPar
Lower bound: \sphinxcode{\sphinxupquote{dRowBound \sphinxhyphen{} dRowUpper}}

\sphinxAtStartPar
Upper bound: \sphinxcode{\sphinxupquote{dRowBound}}
\end{quote}
\end{quote}

\sphinxAtStartPar
\sphinxcode{\sphinxupquote{dRowBound}}
\begin{quote}

\sphinxAtStartPar
Lower bound or RHS of the constraint.
\end{quote}

\sphinxAtStartPar
\sphinxcode{\sphinxupquote{dRowUpper}}
\begin{quote}

\sphinxAtStartPar
Upper bound or \sphinxstylestrong{range} of the constraint.
\end{quote}

\sphinxAtStartPar
\sphinxcode{\sphinxupquote{name}}
\begin{quote}

\sphinxAtStartPar
Name of the constraint. Can be set to \sphinxcode{\sphinxupquote{NULL}}.
\end{quote}
\end{quote}
\end{quote}

\subsubsection{COPT\_AddNLConstrs}
\label{\detokenize{capiref:copt-addnlconstrs}}\begin{quote}

\sphinxAtStartPar
\sphinxstylestrong{Synopsis}
\begin{quote}

\sphinxAtStartPar
\sphinxcode{\sphinxupquote{int COPT\_AddNLConstrs(}}
\begin{quote}

\sphinxAtStartPar
\sphinxcode{\sphinxupquote{copt\_prob *prob,}}

\sphinxAtStartPar
\sphinxcode{\sphinxupquote{int nConstrs,}}

\sphinxAtStartPar
\sphinxcode{\sphinxupquote{const int *tokenBeg,}}

\sphinxAtStartPar
\sphinxcode{\sphinxupquote{const int *tokenCnt,}}

\sphinxAtStartPar
\sphinxcode{\sphinxupquote{const int *tokenElemBeg,}}

\sphinxAtStartPar
\sphinxcode{\sphinxupquote{const int *tokenElemCnt,}}

\sphinxAtStartPar
\sphinxcode{\sphinxupquote{const int *token,}}

\sphinxAtStartPar
\sphinxcode{\sphinxupquote{const double *tokenElem,}}

\sphinxAtStartPar
\sphinxcode{\sphinxupquote{const int *rowMatBeg,}}

\sphinxAtStartPar
\sphinxcode{\sphinxupquote{const int *rowMatCnt,}}

\sphinxAtStartPar
\sphinxcode{\sphinxupquote{const int *rowMatIdx,}}

\sphinxAtStartPar
\sphinxcode{\sphinxupquote{const double *rowMatElem,}}

\sphinxAtStartPar
\sphinxcode{\sphinxupquote{const char *rowSense,}}

\sphinxAtStartPar
\sphinxcode{\sphinxupquote{const double *rowBound,}}

\sphinxAtStartPar
\sphinxcode{\sphinxupquote{const double *rowUpper,}}

\sphinxAtStartPar
\sphinxcode{\sphinxupquote{char const *const *rowNames)}}
\end{quote}
\end{quote}

\sphinxAtStartPar
\sphinxstylestrong{Description}
\begin{quote}

\sphinxAtStartPar
Add a set of nonlinear expression constraints.
Currently, only prefix notation is supported.
\end{quote}

\sphinxAtStartPar
\sphinxstylestrong{Arguments}
\begin{quote}

\sphinxAtStartPar
\sphinxcode{\sphinxupquote{prob}}
\begin{quote}

\sphinxAtStartPar
The COPT problem.
\end{quote}

\sphinxAtStartPar
\sphinxcode{\sphinxupquote{nConstrs}}
\begin{quote}

\sphinxAtStartPar
Number of constraints.
\end{quote}

\sphinxAtStartPar
\sphinxcode{\sphinxupquote{tokenBeg}}
\begin{quote}

\sphinxAtStartPar
Array of starting positions of tokens in expressions.
\end{quote}

\sphinxAtStartPar
\sphinxcode{\sphinxupquote{tokenCnt}}
\begin{quote}

\sphinxAtStartPar
Array of counts of tokens in expressions.
\end{quote}

\sphinxAtStartPar
\sphinxcode{\sphinxupquote{tokenElemBeg}}
\begin{quote}

\sphinxAtStartPar
Array of starting positions of constants in expressions.
\end{quote}

\sphinxAtStartPar
\sphinxcode{\sphinxupquote{tokenElemCnt}}
\begin{quote}

\sphinxAtStartPar
Array of counts of constants in expressions.
\end{quote}

\sphinxAtStartPar
\sphinxcode{\sphinxupquote{token}}
\begin{quote}

\sphinxAtStartPar
Array of tokens in expressions.
\end{quote}

\sphinxAtStartPar
\sphinxcode{\sphinxupquote{tokenElem}}
\begin{quote}

\sphinxAtStartPar
Array of constants in expressions.
\end{quote}

\sphinxAtStartPar
\sphinxcode{\sphinxupquote{rowMatBeg, rowMatCnt, rowMatIdx}} and \sphinxcode{\sphinxupquote{rowMatElem}}
\begin{quote}

\sphinxAtStartPar
Provide the coefficient matrix in a row\sphinxhyphen{}compressed storage format.
For specific examples of sparse matrix compression storage format,
refer to the \sphinxstylestrong{Additional Information} section in \sphinxcode{\sphinxupquote{COPT\_LoadProb}}.
\end{quote}

\sphinxAtStartPar
\sphinxcode{\sphinxupquote{rowSense}}
\begin{quote}

\sphinxAtStartPar
Constraint types.

\sphinxAtStartPar
Refer to the constants section for supported constraint types in COPT.

\sphinxAtStartPar
If \sphinxcode{\sphinxupquote{rowSense}} is not provided, \sphinxcode{\sphinxupquote{rowBound}} and \sphinxcode{\sphinxupquote{rowUpper}}
are treated as the lower and upper bounds of constraints.
This is the recommended way to define constraints.

\sphinxAtStartPar
If \sphinxcode{\sphinxupquote{rowSense}} is provided, \sphinxcode{\sphinxupquote{rowBound}} and \sphinxcode{\sphinxupquote{rowUpper}}
are treated as the RHS and \sphinxstylestrong{range}.
In this case, the \sphinxcode{\sphinxupquote{rowUpper}} array is only used when
there is a \sphinxcode{\sphinxupquote{COPT\_RANGE}} constraint type.
The constraint bounds in this scenario are:
\begin{quote}

\sphinxAtStartPar
Lower bound: \sphinxcode{\sphinxupquote{rowBound{[}i{]} \sphinxhyphen{} fabs(rowUpper{[}i{]})}}

\sphinxAtStartPar
Upper bound: \sphinxcode{\sphinxupquote{rowBound{[}i{]}}}
\end{quote}
\end{quote}

\sphinxAtStartPar
\sphinxcode{\sphinxupquote{rowBound}}
\begin{quote}

\sphinxAtStartPar
Lower bound or RHS of the new constraints.
\end{quote}

\sphinxAtStartPar
\sphinxcode{\sphinxupquote{rowUpper}}
\begin{quote}

\sphinxAtStartPar
Upper bound or \sphinxstylestrong{range} of the new constraints.
\end{quote}

\sphinxAtStartPar
\sphinxcode{\sphinxupquote{rowNames}}
\begin{quote}

\sphinxAtStartPar
Names of the new constraints. Can be set to \sphinxcode{\sphinxupquote{NULL}}.
\end{quote}
\end{quote}
\end{quote}

\subsubsection{COPT\_DelCols}
\label{\detokenize{capiref:copt-delcols}}\begin{quote}

\sphinxAtStartPar
\sphinxstylestrong{Synopsis}
\begin{quote}

\sphinxAtStartPar
\sphinxcode{\sphinxupquote{int COPT\_DelCols(copt\_prob *prob, int num, const int *list)}}
\end{quote}

\sphinxAtStartPar
\sphinxstylestrong{Description}
\begin{quote}

\sphinxAtStartPar
Deletes \sphinxcode{\sphinxupquote{num}} variables (columns) from the problem.
\end{quote}

\sphinxAtStartPar
\sphinxstylestrong{Arguments}
\begin{quote}

\sphinxAtStartPar
\sphinxcode{\sphinxupquote{prob}}
\begin{quote}

\sphinxAtStartPar
The COPT problem.
\end{quote}

\sphinxAtStartPar
\sphinxcode{\sphinxupquote{num}}
\begin{quote}

\sphinxAtStartPar
Number of variables to be deleted.
\end{quote}

\sphinxAtStartPar
\sphinxcode{\sphinxupquote{list}}
\begin{quote}

\sphinxAtStartPar
A list of index of variables to be deleted.
\end{quote}
\end{quote}
\end{quote}

\subsubsection{COPT\_DelPSDCols}
\label{\detokenize{capiref:copt-delpsdcols}}\begin{quote}

\sphinxAtStartPar
\sphinxstylestrong{Synopsis}
\begin{quote}

\sphinxAtStartPar
\sphinxcode{\sphinxupquote{int COPT\_DelPSDCols(copt\_prob *prob, int num, const int *list)}}
\end{quote}

\sphinxAtStartPar
\sphinxstylestrong{Description}
\begin{quote}

\sphinxAtStartPar
Deletes \sphinxcode{\sphinxupquote{num}} PSD variables from the problem.
\end{quote}

\sphinxAtStartPar
\sphinxstylestrong{Arguments}
\begin{quote}

\sphinxAtStartPar
\sphinxcode{\sphinxupquote{prob}}
\begin{quote}

\sphinxAtStartPar
The COPT problem.
\end{quote}

\sphinxAtStartPar
\sphinxcode{\sphinxupquote{num}}
\begin{quote}

\sphinxAtStartPar
Number of PSD variables to be deleted.
\end{quote}

\sphinxAtStartPar
\sphinxcode{\sphinxupquote{list}}
\begin{quote}

\sphinxAtStartPar
A list of index of PSD variables to be deleted.
\end{quote}
\end{quote}
\end{quote}

\subsubsection{COPT\_DelRows}
\label{\detokenize{capiref:copt-delrows}}\begin{quote}

\sphinxAtStartPar
\sphinxstylestrong{Synopsis}
\begin{quote}

\sphinxAtStartPar
\sphinxcode{\sphinxupquote{int COPT\_DelRows(copt\_prob *prob, int num, const int *list)}}
\end{quote}

\sphinxAtStartPar
\sphinxstylestrong{Description}
\begin{quote}

\sphinxAtStartPar
Deletes \sphinxcode{\sphinxupquote{num}} constraints (rows) from the problem.
\end{quote}

\sphinxAtStartPar
\sphinxstylestrong{Arguments}
\begin{quote}

\sphinxAtStartPar
\sphinxcode{\sphinxupquote{prob}}
\begin{quote}

\sphinxAtStartPar
The COPT problem.
\end{quote}

\sphinxAtStartPar
\sphinxcode{\sphinxupquote{num}}
\begin{quote}

\sphinxAtStartPar
The number of constraints to be deleted.
\end{quote}

\sphinxAtStartPar
\sphinxcode{\sphinxupquote{list}}
\begin{quote}

\sphinxAtStartPar
The list of index of constraints to be deleted.
\end{quote}
\end{quote}
\end{quote}

\subsubsection{COPT\_DelSOSs}
\label{\detokenize{capiref:copt-delsoss}}\begin{quote}

\sphinxAtStartPar
\sphinxstylestrong{Synopsis}
\begin{quote}

\sphinxAtStartPar
\sphinxcode{\sphinxupquote{int COPT\_DelSOSs(copt\_prob *prob, int num, const int *list)}}
\end{quote}

\sphinxAtStartPar
\sphinxstylestrong{Description}
\begin{quote}

\sphinxAtStartPar
Deletes \sphinxcode{\sphinxupquote{num}} SOS constraints from the problem.
\end{quote}

\sphinxAtStartPar
\sphinxstylestrong{Arguments}
\begin{quote}

\sphinxAtStartPar
\sphinxcode{\sphinxupquote{prob}}
\begin{quote}

\sphinxAtStartPar
The COPT problem.
\end{quote}

\sphinxAtStartPar
\sphinxcode{\sphinxupquote{num}}
\begin{quote}

\sphinxAtStartPar
The number of SOS constraints to be deleted.
\end{quote}

\sphinxAtStartPar
\sphinxcode{\sphinxupquote{list}}
\begin{quote}

\sphinxAtStartPar
The list of index of SOS constraints to be deleted.
\end{quote}
\end{quote}
\end{quote}

\subsubsection{COPT\_DelCones}
\label{\detokenize{capiref:copt-delcones}}\begin{quote}

\sphinxAtStartPar
\sphinxstylestrong{Synopsis}
\begin{quote}

\sphinxAtStartPar
\sphinxcode{\sphinxupquote{int COPT\_DelCones(copt\_prob *prob, int num, const int *list)}}
\end{quote}

\sphinxAtStartPar
\sphinxstylestrong{Description}
\begin{quote}

\sphinxAtStartPar
Deletes \sphinxcode{\sphinxupquote{num}} Second\sphinxhyphen{}Order\sphinxhyphen{}Cone (SOC) constraints from the problem.
\end{quote}

\sphinxAtStartPar
\sphinxstylestrong{Arguments}
\begin{quote}

\sphinxAtStartPar
\sphinxcode{\sphinxupquote{prob}}
\begin{quote}

\sphinxAtStartPar
The COPT problem.
\end{quote}

\sphinxAtStartPar
\sphinxcode{\sphinxupquote{num}}
\begin{quote}

\sphinxAtStartPar
The number of SOC constraints to be deleted.
\end{quote}

\sphinxAtStartPar
\sphinxcode{\sphinxupquote{list}}
\begin{quote}

\sphinxAtStartPar
The list of index of SOC constraints to be deleted.
\end{quote}
\end{quote}
\end{quote}

\subsubsection{COPT\_DelExpCones}
\label{\detokenize{capiref:copt-delexpcones}}\begin{quote}

\sphinxAtStartPar
\sphinxstylestrong{Synopsis}
\begin{quote}

\sphinxAtStartPar
\sphinxcode{\sphinxupquote{int COPT\_DelExpCones(copt\_prob *prob, int num, const int *list)}}
\end{quote}

\sphinxAtStartPar
\sphinxstylestrong{Description}
\begin{quote}

\sphinxAtStartPar
Deletes \sphinxcode{\sphinxupquote{num}} exponential cone constraints from the problem.
\end{quote}

\sphinxAtStartPar
\sphinxstylestrong{Arguments}
\begin{quote}

\sphinxAtStartPar
\sphinxcode{\sphinxupquote{prob}}
\begin{quote}

\sphinxAtStartPar
The COPT problem.
\end{quote}

\sphinxAtStartPar
\sphinxcode{\sphinxupquote{num}}
\begin{quote}

\sphinxAtStartPar
The number of exponential cone constraints to be deleted.
\end{quote}

\sphinxAtStartPar
\sphinxcode{\sphinxupquote{list}}
\begin{quote}

\sphinxAtStartPar
The list of index of exponential cone constraints to be deleted.
\end{quote}
\end{quote}
\end{quote}

\subsubsection{COPT\_DelAffineCones}
\label{\detokenize{capiref:copt-delaffinecones}}\begin{quote}

\sphinxAtStartPar
\sphinxstylestrong{Synopsis}
\begin{quote}

\sphinxAtStartPar
\sphinxcode{\sphinxupquote{int COPT\_DelAffineCones(copt\_prob *prob, int num, const int *list)}}
\end{quote}

\sphinxAtStartPar
\sphinxstylestrong{Description}
\begin{quote}

\sphinxAtStartPar
Deletes \sphinxcode{\sphinxupquote{num}} affine cone constraints from the problem.
\end{quote}

\sphinxAtStartPar
\sphinxstylestrong{Arguments}
\begin{quote}

\sphinxAtStartPar
\sphinxcode{\sphinxupquote{prob}}
\begin{quote}

\sphinxAtStartPar
The COPT problem.
\end{quote}

\sphinxAtStartPar
\sphinxcode{\sphinxupquote{num}}
\begin{quote}

\sphinxAtStartPar
The number of affine cone constraints to be deleted.
\end{quote}

\sphinxAtStartPar
\sphinxcode{\sphinxupquote{list}}
\begin{quote}

\sphinxAtStartPar
The list of index of affine cone constraints to be deleted.
\end{quote}
\end{quote}
\end{quote}

\subsubsection{COPT\_DelQConstrs}
\label{\detokenize{capiref:copt-delqconstrs}}\begin{quote}

\sphinxAtStartPar
\sphinxstylestrong{Synopsis}
\begin{quote}

\sphinxAtStartPar
\sphinxcode{\sphinxupquote{int COPT\_DelQConstrs(copt\_prob *prob, int num, const int *list)}}
\end{quote}

\sphinxAtStartPar
\sphinxstylestrong{Description}
\begin{quote}

\sphinxAtStartPar
Deletes \sphinxcode{\sphinxupquote{num}} quadratic constraints from the problem.
\end{quote}

\sphinxAtStartPar
\sphinxstylestrong{Arguments}
\begin{quote}

\sphinxAtStartPar
\sphinxcode{\sphinxupquote{prob}}
\begin{quote}

\sphinxAtStartPar
The COPT problem.
\end{quote}

\sphinxAtStartPar
\sphinxcode{\sphinxupquote{num}}
\begin{quote}

\sphinxAtStartPar
The number of quadratic constraints to be deleted.
\end{quote}

\sphinxAtStartPar
\sphinxcode{\sphinxupquote{list}}
\begin{quote}

\sphinxAtStartPar
The list of index of quadratic constraints to be deleted.
\end{quote}
\end{quote}
\end{quote}

\subsubsection{COPT\_DelPSDConstrs}
\label{\detokenize{capiref:copt-delpsdconstrs}}\begin{quote}

\sphinxAtStartPar
\sphinxstylestrong{Synopsis}
\begin{quote}

\sphinxAtStartPar
\sphinxcode{\sphinxupquote{int COPT\_DelPSDConstrs(copt\_prob *prob, int num, const int *list)}}
\end{quote}

\sphinxAtStartPar
\sphinxstylestrong{Description}
\begin{quote}

\sphinxAtStartPar
Deletes \sphinxcode{\sphinxupquote{num}} PSD constraints from the problem.
\end{quote}

\sphinxAtStartPar
\sphinxstylestrong{Arguments}
\begin{quote}

\sphinxAtStartPar
\sphinxcode{\sphinxupquote{prob}}
\begin{quote}

\sphinxAtStartPar
The COPT problem.
\end{quote}

\sphinxAtStartPar
\sphinxcode{\sphinxupquote{num}}
\begin{quote}

\sphinxAtStartPar
The number of PSD constraints to be deleted.
\end{quote}

\sphinxAtStartPar
\sphinxcode{\sphinxupquote{list}}
\begin{quote}

\sphinxAtStartPar
The list of index of PSD constraints to be deleted.
\end{quote}
\end{quote}
\end{quote}

\subsubsection{COPT\_DelLMIConstrs}
\label{\detokenize{capiref:copt-dellmiconstrs}}\begin{quote}

\sphinxAtStartPar
\sphinxstylestrong{Synopsis}
\begin{quote}

\sphinxAtStartPar
\sphinxcode{\sphinxupquote{int COPT\_DelLMIConstrs(copt\_prob *prob, int num, const int *list)}}
\end{quote}

\sphinxAtStartPar
\sphinxstylestrong{Description}
\begin{quote}

\sphinxAtStartPar
Deletes \sphinxcode{\sphinxupquote{num}} LMI constraints from the problem.
\end{quote}

\sphinxAtStartPar
\sphinxstylestrong{Arguments}
\begin{quote}

\sphinxAtStartPar
\sphinxcode{\sphinxupquote{prob}}
\begin{quote}

\sphinxAtStartPar
The COPT problem.
\end{quote}

\sphinxAtStartPar
\sphinxcode{\sphinxupquote{num}}
\begin{quote}

\sphinxAtStartPar
The number of LMI constraints to be deleted.
\end{quote}

\sphinxAtStartPar
\sphinxcode{\sphinxupquote{list}}
\begin{quote}

\sphinxAtStartPar
The list of index of LMI constraints to be deleted.
\end{quote}
\end{quote}
\end{quote}

\subsubsection{COPT\_DelIndicators}
\label{\detokenize{capiref:copt-delindicators}}\begin{quote}

\sphinxAtStartPar
\sphinxstylestrong{Synopsis}
\begin{quote}

\sphinxAtStartPar
\sphinxcode{\sphinxupquote{int COPT\_DelIndicators(copt\_prob *prob, int num, const int *list)}}
\end{quote}

\sphinxAtStartPar
\sphinxstylestrong{Description}
\begin{quote}

\sphinxAtStartPar
Deletes \sphinxcode{\sphinxupquote{num}} indicator constraints from the problem.
\end{quote}

\sphinxAtStartPar
\sphinxstylestrong{Arguments}
\begin{quote}

\sphinxAtStartPar
\sphinxcode{\sphinxupquote{prob}}
\begin{quote}

\sphinxAtStartPar
The COPT problem.
\end{quote}

\sphinxAtStartPar
\sphinxcode{\sphinxupquote{num}}
\begin{quote}

\sphinxAtStartPar
The number of indicator constraints to be deleted.
\end{quote}

\sphinxAtStartPar
\sphinxcode{\sphinxupquote{list}}
\begin{quote}

\sphinxAtStartPar
The list of index of indicator constraints to be deleted.
\end{quote}
\end{quote}
\end{quote}

\subsubsection{COPT\_DelNLConstrs}
\label{\detokenize{capiref:copt-delnlconstrs}}\begin{quote}

\sphinxAtStartPar
\sphinxstylestrong{Synopsis}
\begin{quote}

\sphinxAtStartPar
\sphinxcode{\sphinxupquote{int COPT\_DelNLConstrs(copt\_prob *prob, int num, const int *list)}}
\end{quote}

\sphinxAtStartPar
\sphinxstylestrong{Description}
\begin{quote}

\sphinxAtStartPar
Delete \sphinxcode{\sphinxupquote{num}} nonlinear expression constraints.
\end{quote}

\sphinxAtStartPar
\sphinxstylestrong{Arguments}
\begin{quote}

\sphinxAtStartPar
\sphinxcode{\sphinxupquote{prob}}
\begin{quote}

\sphinxAtStartPar
The COPT problem.
\end{quote}

\sphinxAtStartPar
\sphinxcode{\sphinxupquote{num}}
\begin{quote}

\sphinxAtStartPar
Number of nonlinear expression constraints to delete.
\end{quote}

\sphinxAtStartPar
\sphinxcode{\sphinxupquote{list}}
\begin{quote}

\sphinxAtStartPar
List of indices of nonlinear expression constraints to delete.
\end{quote}
\end{quote}
\end{quote}

\subsubsection{COPT\_DelQuadObj}
\label{\detokenize{capiref:copt-delquadobj}}\begin{quote}

\sphinxAtStartPar
\sphinxstylestrong{Synopsis}
\begin{quote}

\sphinxAtStartPar
\sphinxcode{\sphinxupquote{int COPT\_DelQuadObj(copt\_prob *prob)}}
\end{quote}

\sphinxAtStartPar
\sphinxstylestrong{Description}
\begin{quote}

\sphinxAtStartPar
Deletes the quadratic terms from the quadratic objective function.
\end{quote}

\sphinxAtStartPar
\sphinxstylestrong{Arguments}
\begin{quote}

\sphinxAtStartPar
\sphinxcode{\sphinxupquote{prob}}
\begin{quote}

\sphinxAtStartPar
The COPT problem.
\end{quote}
\end{quote}
\end{quote}

\subsubsection{COPT\_DelNLObj}
\label{\detokenize{capiref:copt-delnlobj}}\begin{quote}

\sphinxAtStartPar
\sphinxstylestrong{Synopsis}
\begin{quote}

\sphinxAtStartPar
\sphinxcode{\sphinxupquote{int COPT\_DelNLObj(copt\_prob *prob)}}
\end{quote}

\sphinxAtStartPar
\sphinxstylestrong{Description}
\begin{quote}

\sphinxAtStartPar
Delete the nonlinear expression terms in the objective function.
\end{quote}

\sphinxAtStartPar
\sphinxstylestrong{Arguments}
\begin{quote}

\sphinxAtStartPar
\sphinxcode{\sphinxupquote{prob}}
\begin{quote}

\sphinxAtStartPar
The COPT problem.
\end{quote}
\end{quote}
\end{quote}

\subsubsection{COPT\_DelPSDObj}
\label{\detokenize{capiref:copt-delpsdobj}}\begin{quote}

\sphinxAtStartPar
\sphinxstylestrong{Synopsis}
\begin{quote}

\sphinxAtStartPar
\sphinxcode{\sphinxupquote{int COPT\_DelPSDObj(copt\_prob *prob)}}
\end{quote}

\sphinxAtStartPar
\sphinxstylestrong{Description}
\begin{quote}

\sphinxAtStartPar
Deletes the PSD terms from objective function.
\end{quote}

\sphinxAtStartPar
\sphinxstylestrong{Arguments}
\begin{quote}

\sphinxAtStartPar
\sphinxcode{\sphinxupquote{prob}}
\begin{quote}

\sphinxAtStartPar
The COPT problem.
\end{quote}
\end{quote}
\end{quote}

\subsubsection{COPT\_SetElem}
\label{\detokenize{capiref:copt-setelem}}\begin{quote}

\sphinxAtStartPar
\sphinxstylestrong{Synopsis}
\begin{quote}

\sphinxAtStartPar
\sphinxcode{\sphinxupquote{int COPT\_SetElem(copt\_prob *prob, int iCol, int iRow, double newElem)}}
\end{quote}

\sphinxAtStartPar
\sphinxstylestrong{Description}
\begin{quote}

\sphinxAtStartPar
Set coefficient of specified row and column.

\sphinxAtStartPar
\sphinxstylestrong{Note:} If \sphinxcode{\sphinxupquote{newElem}} is smaller than or equal to parameter \sphinxcode{\sphinxupquote{MatrixTol}},
the coefficient will be set as zero.
\end{quote}

\sphinxAtStartPar
\sphinxstylestrong{Arguments}
\begin{quote}

\sphinxAtStartPar
\sphinxcode{\sphinxupquote{prob}}
\begin{quote}

\sphinxAtStartPar
The COPT problem.
\end{quote}

\sphinxAtStartPar
\sphinxcode{\sphinxupquote{iCol}}
\begin{quote}

\sphinxAtStartPar
Column index.
\end{quote}

\sphinxAtStartPar
\sphinxcode{\sphinxupquote{iRow}}
\begin{quote}

\sphinxAtStartPar
Row index.
\end{quote}

\sphinxAtStartPar
\sphinxcode{\sphinxupquote{newElem}}
\begin{quote}

\sphinxAtStartPar
New coefficient.
\end{quote}
\end{quote}
\end{quote}

\subsubsection{COPT\_SetElems}
\label{\detokenize{capiref:copt-setelems}}\begin{quote}

\sphinxAtStartPar
\sphinxstylestrong{Synopsis}
\begin{quote}

\sphinxAtStartPar
\sphinxcode{\sphinxupquote{int COPT\_SetElems(copt\_prob *prob, int nelem, const int *cols, const int *rows, const double *elems)}}
\end{quote}

\sphinxAtStartPar
\sphinxstylestrong{Description}
\begin{quote}

\sphinxAtStartPar
Set the coefficients of the specified columns and rows in batches.

\sphinxAtStartPar
\sphinxstylestrong{Note} The index pairs of columns and rows cannot appear repeatedly.
\end{quote}

\sphinxAtStartPar
\sphinxstylestrong{Arguments}
\begin{quote}

\sphinxAtStartPar
\sphinxcode{\sphinxupquote{prob}}
\begin{quote}

\sphinxAtStartPar
The COPT problem.
\end{quote}

\sphinxAtStartPar
\sphinxcode{\sphinxupquote{nelem}}
\begin{quote}

\sphinxAtStartPar
The number of new coefficients to be set.
\end{quote}

\sphinxAtStartPar
\sphinxcode{\sphinxupquote{cols}}
\begin{quote}

\sphinxAtStartPar
Column indexes.
\end{quote}

\sphinxAtStartPar
\sphinxcode{\sphinxupquote{rows}}
\begin{quote}

\sphinxAtStartPar
Row indexes.
\end{quote}

\sphinxAtStartPar
\sphinxcode{\sphinxupquote{elems}}
\begin{quote}

\sphinxAtStartPar
The values of the new coefficients to be set.
\end{quote}
\end{quote}
\end{quote}

\subsubsection{COPT\_SetPSDElem}
\label{\detokenize{capiref:copt-setpsdelem}}\begin{quote}

\sphinxAtStartPar
\sphinxstylestrong{Synopsis}
\begin{quote}

\sphinxAtStartPar
\sphinxcode{\sphinxupquote{int COPT\_SetPSDElem(copt\_prob *prob, int iCol, int iRow, int newIdx)}}
\end{quote}

\sphinxAtStartPar
\sphinxstylestrong{Description}
\begin{quote}

\sphinxAtStartPar
Set symmetric matrix index for given PSD term of PSD constraint.
\end{quote}

\sphinxAtStartPar
\sphinxstylestrong{Arguments}
\begin{quote}

\sphinxAtStartPar
\sphinxcode{\sphinxupquote{prob}}
\begin{quote}

\sphinxAtStartPar
The COPT problem.
\end{quote}

\sphinxAtStartPar
\sphinxcode{\sphinxupquote{iCol}}
\begin{quote}

\sphinxAtStartPar
PSD variable index.
\end{quote}

\sphinxAtStartPar
\sphinxcode{\sphinxupquote{iRow}}
\begin{quote}

\sphinxAtStartPar
PSD constraint index.
\end{quote}

\sphinxAtStartPar
\sphinxcode{\sphinxupquote{newIdx}}
\begin{quote}

\sphinxAtStartPar
New symmetric matrix index.
\end{quote}
\end{quote}
\end{quote}

\subsubsection{COPT\_SetLMIElem}
\label{\detokenize{capiref:copt-setlmielem}}\begin{quote}

\sphinxAtStartPar
\sphinxstylestrong{Synopsis}
\begin{quote}

\sphinxAtStartPar
\sphinxcode{\sphinxupquote{int COPT\_SetLMIElem(copt\_prob *prob, int iCol, int iRow, int newIdx)}}
\end{quote}

\sphinxAtStartPar
\sphinxstylestrong{Description}
\begin{quote}

\sphinxAtStartPar
Set symmetric matrix index for given term of LMI constraint.
\end{quote}

\sphinxAtStartPar
\sphinxstylestrong{Arguments}
\begin{quote}

\sphinxAtStartPar
\sphinxcode{\sphinxupquote{prob}}
\begin{quote}

\sphinxAtStartPar
The COPT problem.
\end{quote}

\sphinxAtStartPar
\sphinxcode{\sphinxupquote{iCol}}
\begin{quote}

\sphinxAtStartPar
Scalar variable index.
\end{quote}

\sphinxAtStartPar
\sphinxcode{\sphinxupquote{iRow}}
\begin{quote}

\sphinxAtStartPar
LMI constraint index.
\end{quote}

\sphinxAtStartPar
\sphinxcode{\sphinxupquote{newIdx}}
\begin{quote}

\sphinxAtStartPar
New coefficient symmetric matrix index.
\end{quote}
\end{quote}
\end{quote}

\subsubsection{COPT\_SetObjSense}
\label{\detokenize{capiref:copt-setobjsense}}\begin{quote}

\sphinxAtStartPar
\sphinxstylestrong{Synopsis}
\begin{quote}

\sphinxAtStartPar
\sphinxcode{\sphinxupquote{int COPT\_SetObjSense(copt\_prob *prob, int iObjSense)}}
\end{quote}

\sphinxAtStartPar
\sphinxstylestrong{Description}
\begin{quote}

\sphinxAtStartPar
Change the objective function sense.
\end{quote}

\sphinxAtStartPar
\sphinxstylestrong{Arguments}
\begin{quote}

\sphinxAtStartPar
\sphinxcode{\sphinxupquote{prob}}
\begin{quote}

\sphinxAtStartPar
The COPT problem.
\end{quote}

\sphinxAtStartPar
\sphinxcode{\sphinxupquote{iObjSense}}
\begin{quote}

\sphinxAtStartPar
The optimization sense, either \sphinxcode{\sphinxupquote{COPT\_MAXIMIZE}} or \sphinxcode{\sphinxupquote{COPT\_MINIMIZE}}.
\end{quote}
\end{quote}
\end{quote}

\subsubsection{COPT\_SetObjConst}
\label{\detokenize{capiref:copt-setobjconst}}\begin{quote}

\sphinxAtStartPar
\sphinxstylestrong{Synopsis}
\begin{quote}

\sphinxAtStartPar
\sphinxcode{\sphinxupquote{int COPT\_SetObjConst(copt\_prob *prob, double dObjConst)}}
\end{quote}

\sphinxAtStartPar
\sphinxstylestrong{Description}
\begin{quote}

\sphinxAtStartPar
Set the constant term of objective function.
\end{quote}

\sphinxAtStartPar
\sphinxstylestrong{Arguments}
\begin{quote}

\sphinxAtStartPar
\sphinxcode{\sphinxupquote{prob}}
\begin{quote}

\sphinxAtStartPar
The COPT problem.
\end{quote}

\sphinxAtStartPar
\sphinxcode{\sphinxupquote{dObjConst}}
\begin{quote}

\sphinxAtStartPar
The constant term of objective function.
\end{quote}
\end{quote}
\end{quote}

\subsubsection{COPT\_SetColObj/Type/Lower/Upper/Names}
\label{\detokenize{capiref:copt-setcolobj-type-lower-upper-names}}\begin{quote}

\sphinxAtStartPar
\sphinxstylestrong{Synopsis}
\begin{quote}

\sphinxAtStartPar
\sphinxcode{\sphinxupquote{int COPT\_SetColObj(copt\_prob *prob, int num, const int *list, const double *obj)}}

\sphinxAtStartPar
\sphinxcode{\sphinxupquote{int COPT\_SetColType(copt\_prob *prob, int num, const int *list, const char *type)}}

\sphinxAtStartPar
\sphinxcode{\sphinxupquote{int COPT\_SetColLower(copt\_prob *prob, int num, const int *list, const double *lower)}}

\sphinxAtStartPar
\sphinxcode{\sphinxupquote{int COPT\_SetColUpper(copt\_prob *prob, int num, const int *list, const double *upper)}}

\sphinxAtStartPar
\sphinxcode{\sphinxupquote{int COPT\_SetColNames(copt\_prob *prob, int num, const int *list, char const *const *names)}}
\end{quote}

\sphinxAtStartPar
\sphinxstylestrong{Description}
\begin{quote}

\sphinxAtStartPar
These five API functions each modifies
\begin{quote}

\sphinxAtStartPar
objective coefficients

\sphinxAtStartPar
variable types

\sphinxAtStartPar
lower bounds

\sphinxAtStartPar
upper bounds

\sphinxAtStartPar
names
\end{quote}

\sphinxAtStartPar
of \sphinxcode{\sphinxupquote{num}} variables (columns) in the problem.
\end{quote}

\sphinxAtStartPar
\sphinxstylestrong{Arguments}
\begin{quote}

\sphinxAtStartPar
\sphinxcode{\sphinxupquote{prob}}
\begin{quote}

\sphinxAtStartPar
The COPT problem.
\end{quote}

\sphinxAtStartPar
\sphinxcode{\sphinxupquote{num}}
\begin{quote}

\sphinxAtStartPar
Number of variables to modify.
\end{quote}

\sphinxAtStartPar
\sphinxcode{\sphinxupquote{list}}
\begin{quote}

\sphinxAtStartPar
A list of index of variables to modify.
\end{quote}

\sphinxAtStartPar
\sphinxcode{\sphinxupquote{obj}}
\begin{quote}

\sphinxAtStartPar
New objective coefficients for each variable in the \sphinxcode{\sphinxupquote{list}}.
\end{quote}

\sphinxAtStartPar
\sphinxcode{\sphinxupquote{types}}
\begin{quote}

\sphinxAtStartPar
New types for each variable in the \sphinxcode{\sphinxupquote{list}}.
\end{quote}

\sphinxAtStartPar
\sphinxcode{\sphinxupquote{lower}}
\begin{quote}

\sphinxAtStartPar
New lower bounds for each variable in the \sphinxcode{\sphinxupquote{list}}.
\end{quote}

\sphinxAtStartPar
\sphinxcode{\sphinxupquote{upper}}
\begin{quote}

\sphinxAtStartPar
New upper bounds for each variable in the \sphinxcode{\sphinxupquote{list}}.
\end{quote}

\sphinxAtStartPar
\sphinxcode{\sphinxupquote{names}}
\begin{quote}

\sphinxAtStartPar
New names for each variable in the \sphinxcode{\sphinxupquote{list}}.
\end{quote}
\end{quote}
\end{quote}

\subsubsection{COPT\_SetPSDColNames}
\label{\detokenize{capiref:copt-setpsdcolnames}}\begin{quote}

\sphinxAtStartPar
\sphinxstylestrong{Synopsis}
\begin{quote}

\sphinxAtStartPar
\sphinxcode{\sphinxupquote{int COPT\_SetPSDColNames(copt\_prob *prob, int num, const int *list, char const *const *names)}}
\end{quote}

\sphinxAtStartPar
\sphinxstylestrong{Description}
\begin{quote}

\sphinxAtStartPar
Modify names of \sphinxcode{\sphinxupquote{num}} PSD variables.
\end{quote}

\sphinxAtStartPar
\sphinxstylestrong{Arguments}
\begin{quote}

\sphinxAtStartPar
\sphinxcode{\sphinxupquote{prob}}
\begin{quote}

\sphinxAtStartPar
The COPT problem.
\end{quote}

\sphinxAtStartPar
\sphinxcode{\sphinxupquote{num}}
\begin{quote}

\sphinxAtStartPar
Number of PSD variables to modify.
\end{quote}

\sphinxAtStartPar
\sphinxcode{\sphinxupquote{list}}
\begin{quote}

\sphinxAtStartPar
A list of index of PSD variables to modify.
\end{quote}

\sphinxAtStartPar
\sphinxcode{\sphinxupquote{names}}
\begin{quote}

\sphinxAtStartPar
New names for each PSD variable in the \sphinxcode{\sphinxupquote{list}}.
\end{quote}
\end{quote}
\end{quote}

\subsubsection{COPT\_SetRowLower/Upper/Names}
\label{\detokenize{capiref:copt-setrowlower-upper-names}}\begin{quote}

\sphinxAtStartPar
\sphinxstylestrong{Synopsis}
\begin{quote}

\sphinxAtStartPar
\sphinxcode{\sphinxupquote{int COPT\_SetRowLower(copt\_prob *prob, int num, const int *list, const double *lower)}}

\sphinxAtStartPar
\sphinxcode{\sphinxupquote{int COPT\_SetRowUpper(copt\_prob *prob, int num, const int *list, const double *upper)}}

\sphinxAtStartPar
\sphinxcode{\sphinxupquote{int COPT\_SetRowNames(copt\_prob *prob, int num, const int *list, char const *const *names)}}
\end{quote}

\sphinxAtStartPar
\sphinxstylestrong{Description}
\begin{quote}

\sphinxAtStartPar
These three API functions each modifies
\begin{quote}

\sphinxAtStartPar
lower bounds

\sphinxAtStartPar
upper bounds

\sphinxAtStartPar
names
\end{quote}

\sphinxAtStartPar
of \sphinxcode{\sphinxupquote{num}} constraints (rows) in the problem.
\end{quote}

\sphinxAtStartPar
\sphinxstylestrong{Arguments}
\begin{quote}

\sphinxAtStartPar
\sphinxcode{\sphinxupquote{prob}}
\begin{quote}

\sphinxAtStartPar
The COPT problem.
\end{quote}

\sphinxAtStartPar
\sphinxcode{\sphinxupquote{num}}
\begin{quote}

\sphinxAtStartPar
Number of constraints to modify.
\end{quote}

\sphinxAtStartPar
\sphinxcode{\sphinxupquote{list}}
\begin{quote}

\sphinxAtStartPar
A list of index of constraints to modify.
\end{quote}

\sphinxAtStartPar
\sphinxcode{\sphinxupquote{lower}}
\begin{quote}

\sphinxAtStartPar
New lower bounds for each constraint in the \sphinxcode{\sphinxupquote{list}}.
\end{quote}

\sphinxAtStartPar
\sphinxcode{\sphinxupquote{upper}}
\begin{quote}

\sphinxAtStartPar
New upper bounds for each constraint in the \sphinxcode{\sphinxupquote{list}}.
\end{quote}

\sphinxAtStartPar
\sphinxcode{\sphinxupquote{names}}
\begin{quote}

\sphinxAtStartPar
New names for each constraint in the \sphinxcode{\sphinxupquote{list}}.
\end{quote}
\end{quote}
\end{quote}

\subsubsection{COPT\_SetQConstrSense/Rhs/Names}
\label{\detokenize{capiref:copt-setqconstrsense-rhs-names}}\begin{quote}

\sphinxAtStartPar
\sphinxstylestrong{Synopsis}
\begin{quote}

\sphinxAtStartPar
\sphinxcode{\sphinxupquote{int COPT\_SetQConstrSense(copt\_prob *prob, int num, const int *list, const char *sense)}}

\sphinxAtStartPar
\sphinxcode{\sphinxupquote{int COPT\_SetQConstrRhs(copt\_prob *prob, int num, const int *list, const double *rhs)}}

\sphinxAtStartPar
\sphinxcode{\sphinxupquote{int COPT\_SetQConstrNames(copt\_prob *prob, int num, const int *list, char const *const *names)}}
\end{quote}

\sphinxAtStartPar
\sphinxstylestrong{Description}
\begin{quote}

\sphinxAtStartPar
These three API functions each modifies
\begin{quote}

\sphinxAtStartPar
senses

\sphinxAtStartPar
RHS

\sphinxAtStartPar
names
\end{quote}

\sphinxAtStartPar
of \sphinxcode{\sphinxupquote{num}} quadratic constraints (rows) in the problem.
\end{quote}

\sphinxAtStartPar
\sphinxstylestrong{Arguments}
\begin{quote}

\sphinxAtStartPar
\sphinxcode{\sphinxupquote{prob}}
\begin{quote}

\sphinxAtStartPar
The COPT problem.
\end{quote}

\sphinxAtStartPar
\sphinxcode{\sphinxupquote{num}}
\begin{quote}

\sphinxAtStartPar
Number of quadratic constraints to modify.
\end{quote}

\sphinxAtStartPar
\sphinxcode{\sphinxupquote{list}}
\begin{quote}

\sphinxAtStartPar
A list of index of quadratic constraints to modify.
\end{quote}

\sphinxAtStartPar
\sphinxcode{\sphinxupquote{sense}}
\begin{quote}

\sphinxAtStartPar
New senses for each quadratic constraint in the \sphinxcode{\sphinxupquote{list}}.
\end{quote}

\sphinxAtStartPar
\sphinxcode{\sphinxupquote{rhs}}
\begin{quote}

\sphinxAtStartPar
New RHS for each quadratic constraint in the \sphinxcode{\sphinxupquote{list}}.
\end{quote}

\sphinxAtStartPar
\sphinxcode{\sphinxupquote{names}}
\begin{quote}

\sphinxAtStartPar
New names for each quadratic constraint in the \sphinxcode{\sphinxupquote{list}}.
\end{quote}
\end{quote}
\end{quote}

\subsubsection{COPT\_SetPSDConstrLower/Upper/Names}
\label{\detokenize{capiref:copt-setpsdconstrlower-upper-names}}\begin{quote}

\sphinxAtStartPar
\sphinxstylestrong{Synopsis}
\begin{quote}

\sphinxAtStartPar
\sphinxcode{\sphinxupquote{int COPT\_SetPSDConstrLower(copt\_prob *prob, int num, const int *list, const double *lower)}}

\sphinxAtStartPar
\sphinxcode{\sphinxupquote{int COPT\_SetPSDConstrUpper(copt\_prob *prob, int num, const int *list, const double *upper)}}

\sphinxAtStartPar
\sphinxcode{\sphinxupquote{int COPT\_SetPSDConstrNames(copt\_prob *prob, int num, const int *list, char const *const *names)}}
\end{quote}

\sphinxAtStartPar
\sphinxstylestrong{Description}
\begin{quote}

\sphinxAtStartPar
These three API functions each modifies
\begin{quote}

\sphinxAtStartPar
lower bounds

\sphinxAtStartPar
upper bounds

\sphinxAtStartPar
names
\end{quote}

\sphinxAtStartPar
of \sphinxcode{\sphinxupquote{num}} PSD constraints in the problem.
\end{quote}

\sphinxAtStartPar
\sphinxstylestrong{Arguments}
\begin{quote}

\sphinxAtStartPar
\sphinxcode{\sphinxupquote{prob}}
\begin{quote}

\sphinxAtStartPar
The COPT problem.
\end{quote}

\sphinxAtStartPar
\sphinxcode{\sphinxupquote{num}}
\begin{quote}

\sphinxAtStartPar
Number of PSD constraints to modify.
\end{quote}

\sphinxAtStartPar
\sphinxcode{\sphinxupquote{list}}
\begin{quote}

\sphinxAtStartPar
A list of index of PSD constraints to modify.
\end{quote}

\sphinxAtStartPar
\sphinxcode{\sphinxupquote{lower}}
\begin{quote}

\sphinxAtStartPar
New lower bounds for each PSD constraint in the \sphinxcode{\sphinxupquote{list}}.
\end{quote}

\sphinxAtStartPar
\sphinxcode{\sphinxupquote{upper}}
\begin{quote}

\sphinxAtStartPar
New upper bounds for each PSD constraint in the \sphinxcode{\sphinxupquote{list}}.
\end{quote}

\sphinxAtStartPar
\sphinxcode{\sphinxupquote{names}}
\begin{quote}

\sphinxAtStartPar
New names for each PSD constraint in the \sphinxcode{\sphinxupquote{list}}.
\end{quote}
\end{quote}
\end{quote}

\subsubsection{COPT\_SetLMIConstrRhs}
\label{\detokenize{capiref:copt-setlmiconstrrhs}}\begin{quote}

\sphinxAtStartPar
\sphinxstylestrong{Synopsis}
\begin{quote}

\sphinxAtStartPar
\sphinxcode{\sphinxupquote{int COPT\_SetLMIConstrRhs(copt\_prob *prob, int num, const int *list, const int *newIdx)}}
\end{quote}

\sphinxAtStartPar
\sphinxstylestrong{Description}
\begin{quote}

\sphinxAtStartPar
Modify the constant\sphinxhyphen{}term symmetric matrix of \sphinxcode{\sphinxupquote{num}} LMI constraints.
\end{quote}

\sphinxAtStartPar
\sphinxstylestrong{Arguments}
\begin{quote}

\sphinxAtStartPar
\sphinxcode{\sphinxupquote{prob}}
\begin{quote}

\sphinxAtStartPar
The COPT problem.
\end{quote}

\sphinxAtStartPar
\sphinxcode{\sphinxupquote{num}}
\begin{quote}

\sphinxAtStartPar
Number of LMI constraints to modify.
\end{quote}

\sphinxAtStartPar
\sphinxcode{\sphinxupquote{list}}
\begin{quote}

\sphinxAtStartPar
A list of index of LMI constraints to modify.
\end{quote}

\sphinxAtStartPar
\sphinxcode{\sphinxupquote{newIdx}}
\begin{quote}

\sphinxAtStartPar
The new index of the constant\sphinxhyphen{}term symmetric matrix to be set.
\end{quote}
\end{quote}
\end{quote}

\subsubsection{COPT\_SetLMIConstrNames}
\label{\detokenize{capiref:copt-setlmiconstrnames}}\begin{quote}

\sphinxAtStartPar
\sphinxstylestrong{Synopsis}
\begin{quote}

\sphinxAtStartPar
\sphinxcode{\sphinxupquote{int COPT\_SetLMIConstrNames(copt\_prob *prob, int num, const int *list, char const *const *names)}}
\end{quote}

\sphinxAtStartPar
\sphinxstylestrong{Description}
\begin{quote}

\sphinxAtStartPar
Modify the names of \sphinxcode{\sphinxupquote{num}} LMI constraints.
\end{quote}

\sphinxAtStartPar
\sphinxstylestrong{Arguments}
\begin{quote}

\sphinxAtStartPar
\sphinxcode{\sphinxupquote{prob}}
\begin{quote}

\sphinxAtStartPar
The COPT problem.
\end{quote}

\sphinxAtStartPar
\sphinxcode{\sphinxupquote{num}}
\begin{quote}

\sphinxAtStartPar
Number of LMI constraints to modify.
\end{quote}

\sphinxAtStartPar
\sphinxcode{\sphinxupquote{list}}
\begin{quote}

\sphinxAtStartPar
A list of index of LMI constraints to modify.
\end{quote}

\sphinxAtStartPar
\sphinxcode{\sphinxupquote{names}}
\begin{quote}

\sphinxAtStartPar
New names for each LMI constraint in the \sphinxcode{\sphinxupquote{list}}.
\end{quote}
\end{quote}
\end{quote}

\subsubsection{COPT\_SetIndicatorNames}
\label{\detokenize{capiref:copt-setindicatornames}}\begin{quote}

\sphinxAtStartPar
\sphinxstylestrong{Synopsis}
\begin{quote}

\sphinxAtStartPar
\sphinxcode{\sphinxupquote{int COPT\_SetIndicatorNames(copt\_prob *prob, int num, const int *list, char const *const *names)}}
\end{quote}

\sphinxAtStartPar
\sphinxstylestrong{Description}
\begin{quote}

\sphinxAtStartPar
Modifies names of \sphinxcode{\sphinxupquote{num}} indicator constraints in the problem.
\end{quote}

\sphinxAtStartPar
\sphinxstylestrong{Arguments}
\begin{quote}

\sphinxAtStartPar
\sphinxcode{\sphinxupquote{prob}}
\begin{quote}

\sphinxAtStartPar
The COPT problem.
\end{quote}

\sphinxAtStartPar
\sphinxcode{\sphinxupquote{num}}
\begin{quote}

\sphinxAtStartPar
Number of indicator constraints to modify.
\end{quote}

\sphinxAtStartPar
\sphinxcode{\sphinxupquote{list}}
\begin{quote}

\sphinxAtStartPar
A list of indexes of indicator constraints to modify.
\end{quote}

\sphinxAtStartPar
\sphinxcode{\sphinxupquote{names}}
\begin{quote}

\sphinxAtStartPar
New names for each indicator constraint in the \sphinxcode{\sphinxupquote{list}}.
\end{quote}
\end{quote}
\end{quote}

\subsubsection{COPT\_SetAffineConeNames}
\label{\detokenize{capiref:copt-setaffineconenames}}\begin{quote}

\sphinxAtStartPar
\sphinxstylestrong{Synopsis}
\begin{quote}

\sphinxAtStartPar
\sphinxcode{\sphinxupquote{int COPT\_SetAffineConeNames(copt\_prob *prob, int num, const int *list, char const *const *names)}}
\end{quote}

\sphinxAtStartPar
\sphinxstylestrong{Description}
\begin{quote}

\sphinxAtStartPar
Modifies names of \sphinxcode{\sphinxupquote{num}} affine cones in the problem.
\end{quote}

\sphinxAtStartPar
\sphinxstylestrong{Arguments}
\begin{quote}

\sphinxAtStartPar
\sphinxcode{\sphinxupquote{prob}}
\begin{quote}

\sphinxAtStartPar
The COPT problem.
\end{quote}

\sphinxAtStartPar
\sphinxcode{\sphinxupquote{num}}
\begin{quote}

\sphinxAtStartPar
Number of affine cones to modify.
\end{quote}

\sphinxAtStartPar
\sphinxcode{\sphinxupquote{list}}
\begin{quote}

\sphinxAtStartPar
A list of indexes of affine cones to modify.
\end{quote}

\sphinxAtStartPar
\sphinxcode{\sphinxupquote{names}}
\begin{quote}

\sphinxAtStartPar
New names for each affine cone in the \sphinxcode{\sphinxupquote{list}}.
\end{quote}
\end{quote}
\end{quote}

\subsubsection{COPT\_SetNLConstrLower/Upper/Names}
\label{\detokenize{capiref:copt-setnlconstrlower-upper-names}}\begin{quote}

\sphinxAtStartPar
\sphinxstylestrong{Synopsis}
\begin{quote}

\sphinxAtStartPar
\sphinxcode{\sphinxupquote{int COPT\_SetNLConstrLower(copt\_prob *prob, int num, const int *list,}}
\sphinxcode{\sphinxupquote{const double *lower)}}

\sphinxAtStartPar
\sphinxcode{\sphinxupquote{int COPT\_SetNLConstrUpper(copt\_prob *prob, int num, const int *list,}}
\sphinxcode{\sphinxupquote{const double *upper)}}

\sphinxAtStartPar
\sphinxcode{\sphinxupquote{int COPT\_SetNLConstrNames(copt\_prob *prob, int num, const int *list,}}
\sphinxcode{\sphinxupquote{char const *const *names)}}
\end{quote}

\sphinxAtStartPar
\sphinxstylestrong{Description}
\begin{quote}

\sphinxAtStartPar
These three functions respectively modify the following properties
of \sphinxcode{\sphinxupquote{num}} nonlinear expression constraints:
\begin{quote}

\sphinxAtStartPar
Lower bound

\sphinxAtStartPar
Upper bound

\sphinxAtStartPar
Name
\end{quote}
\end{quote}

\sphinxAtStartPar
\sphinxstylestrong{Arguments}
\begin{quote}

\sphinxAtStartPar
\sphinxcode{\sphinxupquote{prob}}
\begin{quote}

\sphinxAtStartPar
The COPT problem.
\end{quote}

\sphinxAtStartPar
\sphinxcode{\sphinxupquote{num}}
\begin{quote}

\sphinxAtStartPar
Number of nonlinear expression constraints to modify.
\end{quote}

\sphinxAtStartPar
\sphinxcode{\sphinxupquote{list}}
\begin{quote}

\sphinxAtStartPar
List of nonlinear expression constraints to modify.
\end{quote}

\sphinxAtStartPar
\sphinxcode{\sphinxupquote{lower}}
\begin{quote}

\sphinxAtStartPar
New lower bounds for the listed nonlinear expression constraints.
\end{quote}

\sphinxAtStartPar
\sphinxcode{\sphinxupquote{upper}}
\begin{quote}

\sphinxAtStartPar
New upper bounds for the listed nonlinear expression constraints.
\end{quote}

\sphinxAtStartPar
\sphinxcode{\sphinxupquote{names}}
\begin{quote}

\sphinxAtStartPar
New names for the listed nonlinear expression constraints.
\end{quote}
\end{quote}
\end{quote}

\subsubsection{COPT\_ReplaceColObj}
\label{\detokenize{capiref:copt-replacecolobj}}\begin{quote}

\sphinxAtStartPar
\sphinxstylestrong{Synopsis}
\begin{quote}

\sphinxAtStartPar
\sphinxcode{\sphinxupquote{int COPT\_ReplaceColObj(copt\_prob *prob, int num, const int *list, const double *obj)}}
\end{quote}

\sphinxAtStartPar
\sphinxstylestrong{Description}
\begin{quote}

\sphinxAtStartPar
Replace objective function with new objective function represented by
specified objective costs.
\end{quote}

\sphinxAtStartPar
\sphinxstylestrong{Arguments}
\begin{quote}

\sphinxAtStartPar
\sphinxcode{\sphinxupquote{prob}}
\begin{quote}

\sphinxAtStartPar
The COPT problem.
\end{quote}

\sphinxAtStartPar
\sphinxcode{\sphinxupquote{num}}
\begin{quote}

\sphinxAtStartPar
Number of variables to be modified.
\end{quote}

\sphinxAtStartPar
\sphinxcode{\sphinxupquote{list}}
\begin{quote}

\sphinxAtStartPar
Index of variables to be modified.
\end{quote}

\sphinxAtStartPar
\sphinxcode{\sphinxupquote{obj}}
\begin{quote}

\sphinxAtStartPar
Objective costs of modified variables.
\end{quote}
\end{quote}
\end{quote}

\subsubsection{COPT\_ReplacePSDObj}
\label{\detokenize{capiref:copt-replacepsdobj}}\begin{quote}

\sphinxAtStartPar
\sphinxstylestrong{Synopsis}
\begin{quote}

\sphinxAtStartPar
\sphinxcode{\sphinxupquote{int COPT\_ReplacePSDObj(copt\_prob *prob, int num, const int *list, const int *idx)}}
\end{quote}

\sphinxAtStartPar
\sphinxstylestrong{Description}
\begin{quote}

\sphinxAtStartPar
Replace PSD terms in objective function with specified PSD terms.
\end{quote}

\sphinxAtStartPar
\sphinxstylestrong{Arguments}
\begin{quote}

\sphinxAtStartPar
\sphinxcode{\sphinxupquote{prob}}
\begin{quote}

\sphinxAtStartPar
The COPT problem.
\end{quote}

\sphinxAtStartPar
\sphinxcode{\sphinxupquote{num}}
\begin{quote}

\sphinxAtStartPar
Number of PSD variables to be modified.
\end{quote}

\sphinxAtStartPar
\sphinxcode{\sphinxupquote{list}}
\begin{quote}

\sphinxAtStartPar
Index of PSD variables to be modified.
\end{quote}

\sphinxAtStartPar
\sphinxcode{\sphinxupquote{idx}}
\begin{quote}

\sphinxAtStartPar
Symmetric matrix index of modified PSD variables.
\end{quote}
\end{quote}
\end{quote}

\subsubsection{COPT\_SetQuadObj}
\label{\detokenize{capiref:copt-setquadobj}}\begin{quote}

\sphinxAtStartPar
\sphinxstylestrong{Synopsis}
\begin{quote}

\sphinxAtStartPar
\sphinxcode{\sphinxupquote{int COPT\_SetQuadObj(copt\_prob *prob, int num, int *qRow, int *qCol, double *qElem)}}
\end{quote}

\sphinxAtStartPar
\sphinxstylestrong{Description}
\begin{quote}

\sphinxAtStartPar
Set the quadratic terms of the quadratic objective function.
\end{quote}

\sphinxAtStartPar
\sphinxstylestrong{Arguments}
\begin{quote}

\sphinxAtStartPar
\sphinxcode{\sphinxupquote{prob}}
\begin{quote}

\sphinxAtStartPar
The COPT problem.
\end{quote}

\sphinxAtStartPar
\sphinxcode{\sphinxupquote{num}}
\begin{quote}

\sphinxAtStartPar
Number of non\sphinxhyphen{}zero quadratic terms of the quadratic objective function.
\end{quote}

\sphinxAtStartPar
\sphinxcode{\sphinxupquote{qRow}}
\begin{quote}

\sphinxAtStartPar
Row index of non\sphinxhyphen{}zero quadratic terms of the quadratic objective function.
\end{quote}

\sphinxAtStartPar
\sphinxcode{\sphinxupquote{qCol}}
\begin{quote}

\sphinxAtStartPar
Column index of non\sphinxhyphen{}zero quadratic terms of the quadratic objective function.
\end{quote}

\sphinxAtStartPar
\sphinxcode{\sphinxupquote{qElem}}
\begin{quote}

\sphinxAtStartPar
Values of non\sphinxhyphen{}zero quadratic terms of the quadratic objective function.
\end{quote}
\end{quote}
\end{quote}

\subsubsection{COPT\_SetNLObj}
\label{\detokenize{capiref:copt-setnlobj}}\begin{quote}

\sphinxAtStartPar
\sphinxstylestrong{Synopsis}
\begin{quote}

\sphinxAtStartPar
\sphinxcode{\sphinxupquote{int COPT\_SetNLObj(copt\_prob *prob, int nToken, int nTokenElem,}}
\sphinxcode{\sphinxupquote{const int *token, const double *tokenElem)}}
\end{quote}

\sphinxAtStartPar
\sphinxstylestrong{Description}
\begin{quote}

\sphinxAtStartPar
Set the nonlinear expression terms in the nonlinear objective function.
\end{quote}

\sphinxAtStartPar
\sphinxstylestrong{Arguments}
\begin{quote}

\sphinxAtStartPar
\sphinxcode{\sphinxupquote{prob}}
\begin{quote}

\sphinxAtStartPar
The COPT problem.
\end{quote}

\sphinxAtStartPar
\sphinxcode{\sphinxupquote{nToken}}
\begin{quote}

\sphinxAtStartPar
Number of tokens in the expression.
\end{quote}

\sphinxAtStartPar
\sphinxcode{\sphinxupquote{nTokenElem}}
\begin{quote}

\sphinxAtStartPar
Number of constants in the expression.
\end{quote}

\sphinxAtStartPar
\sphinxcode{\sphinxupquote{token}}
\begin{quote}

\sphinxAtStartPar
Array of tokens in the expression.
\end{quote}

\sphinxAtStartPar
\sphinxcode{\sphinxupquote{tokenElem}}
\begin{quote}

\sphinxAtStartPar
Array of constants in the expression.
\end{quote}
\end{quote}
\end{quote}

\subsubsection{COPT\_SetPSDObj}
\label{\detokenize{capiref:copt-setpsdobj}}\begin{quote}

\sphinxAtStartPar
\sphinxstylestrong{Synopsis}
\begin{quote}

\sphinxAtStartPar
\sphinxcode{\sphinxupquote{int COPT\_SetPSDObj(copt\_prob *prob, int iCol, int newIdx)}}
\end{quote}

\sphinxAtStartPar
\sphinxstylestrong{Description}
\begin{quote}

\sphinxAtStartPar
Set PSD terms of objective function.
\end{quote}

\sphinxAtStartPar
\sphinxstylestrong{Arguments}
\begin{quote}

\sphinxAtStartPar
\sphinxcode{\sphinxupquote{prob}}
\begin{quote}

\sphinxAtStartPar
The COPT problem.
\end{quote}

\sphinxAtStartPar
\sphinxcode{\sphinxupquote{iCol}}
\begin{quote}

\sphinxAtStartPar
PSD variable index.
\end{quote}

\sphinxAtStartPar
\sphinxcode{\sphinxupquote{newIdx}}
\begin{quote}

\sphinxAtStartPar
Symmetric matrix index.
\end{quote}
\end{quote}
\end{quote}

\subsection{Reading and writing the problem}
\label{\detokenize{capiref:reading-and-writing-the-problem}}

\subsubsection{COPT\_ReadMps}
\label{\detokenize{capiref:copt-readmps}}\begin{quote}

\sphinxAtStartPar
\sphinxstylestrong{Synopsis}
\begin{quote}

\sphinxAtStartPar
\sphinxcode{\sphinxupquote{int COPT\_ReadMps(copt\_prob *prob, const char *mpsfilename)}}
\end{quote}

\sphinxAtStartPar
\sphinxstylestrong{Description}
\begin{quote}

\sphinxAtStartPar
Reads a problem from a MPS file.
\end{quote}

\sphinxAtStartPar
\sphinxstylestrong{Arguments}
\begin{quote}

\sphinxAtStartPar
\sphinxcode{\sphinxupquote{prob}}
\begin{quote}

\sphinxAtStartPar
The COPT problem.
\end{quote}

\sphinxAtStartPar
\sphinxcode{\sphinxupquote{mpsfilename}}
\begin{quote}

\sphinxAtStartPar
The path to the MPS file.
\end{quote}
\end{quote}
\end{quote}

\subsubsection{COPT\_ReadLp}
\label{\detokenize{capiref:copt-readlp}}\begin{quote}

\sphinxAtStartPar
\sphinxstylestrong{Synopsis}
\begin{quote}

\sphinxAtStartPar
\sphinxcode{\sphinxupquote{int COPT\_ReadLp(copt\_prob *prob, const char *lpfilename)}}
\end{quote}

\sphinxAtStartPar
\sphinxstylestrong{Description}
\begin{quote}

\sphinxAtStartPar
Read a problem from a LP file.
\end{quote}

\sphinxAtStartPar
\sphinxstylestrong{Arguments}
\begin{quote}

\sphinxAtStartPar
\sphinxcode{\sphinxupquote{prob}}
\begin{quote}

\sphinxAtStartPar
The COPT problem.
\end{quote}

\sphinxAtStartPar
\sphinxcode{\sphinxupquote{lpfilename}}
\begin{quote}

\sphinxAtStartPar
The path to the LP file.
\end{quote}
\end{quote}
\end{quote}

\subsubsection{COPT\_ReadSDPA}
\label{\detokenize{capiref:copt-readsdpa}}\begin{quote}

\sphinxAtStartPar
\sphinxstylestrong{Synopsis}
\begin{quote}

\sphinxAtStartPar
\sphinxcode{\sphinxupquote{int COPT\_ReadSDPA(copt\_prob *prob, const char *sdpafilename)}}
\end{quote}

\sphinxAtStartPar
\sphinxstylestrong{Description}
\begin{quote}

\sphinxAtStartPar
Reads a problem from SDPA format file.
\end{quote}

\sphinxAtStartPar
\sphinxstylestrong{Arguments}
\begin{quote}

\sphinxAtStartPar
\sphinxcode{\sphinxupquote{prob}}
\begin{quote}

\sphinxAtStartPar
The COPT problem.
\end{quote}

\sphinxAtStartPar
\sphinxcode{\sphinxupquote{sdpafilename}}
\begin{quote}

\sphinxAtStartPar
The path to the SDPA format file.
\end{quote}
\end{quote}
\end{quote}

\subsubsection{COPT\_ReadCbf}
\label{\detokenize{capiref:copt-readcbf}}\begin{quote}

\sphinxAtStartPar
\sphinxstylestrong{Synopsis}
\begin{quote}

\sphinxAtStartPar
\sphinxcode{\sphinxupquote{int COPT\_ReadCbf(copt\_prob *prob, const char *cbffilename)}}
\end{quote}

\sphinxAtStartPar
\sphinxstylestrong{Description}
\begin{quote}

\sphinxAtStartPar
Reads a problem from CBF format file.
\end{quote}

\sphinxAtStartPar
\sphinxstylestrong{Arguments}
\begin{quote}

\sphinxAtStartPar
\sphinxcode{\sphinxupquote{prob}}
\begin{quote}

\sphinxAtStartPar
The COPT problem.
\end{quote}

\sphinxAtStartPar
\sphinxcode{\sphinxupquote{cbffilename}}
\begin{quote}

\sphinxAtStartPar
The path to the CBF format file.
\end{quote}
\end{quote}
\end{quote}

\subsubsection{COPT\_ReadBin}
\label{\detokenize{capiref:copt-readbin}}\begin{quote}

\sphinxAtStartPar
\sphinxstylestrong{Synopsis}
\begin{quote}

\sphinxAtStartPar
\sphinxcode{\sphinxupquote{int COPT\_ReadBin(copt\_prob *prob, const char *binfilename)}}
\end{quote}

\sphinxAtStartPar
\sphinxstylestrong{Description}
\begin{quote}

\sphinxAtStartPar
Reads a problem from a COPT binary format file.
\end{quote}

\sphinxAtStartPar
\sphinxstylestrong{Arguments}
\begin{quote}

\sphinxAtStartPar
\sphinxcode{\sphinxupquote{prob}}
\begin{quote}

\sphinxAtStartPar
The COPT problem.
\end{quote}

\sphinxAtStartPar
\sphinxcode{\sphinxupquote{binfilename}}
\begin{quote}

\sphinxAtStartPar
The path to the COPT binary format file.
\end{quote}
\end{quote}
\end{quote}

\subsubsection{COPT\_ReadBlob}
\label{\detokenize{capiref:copt-readblob}}\begin{quote}

\sphinxAtStartPar
\sphinxstylestrong{Synopsis}
\begin{quote}

\sphinxAtStartPar
\sphinxcode{\sphinxupquote{int COPT\_ReadBlob(copt\_prob *prob, void *blob, COPT\_INT64 len)}}
\end{quote}

\sphinxAtStartPar
\sphinxstylestrong{Description}
\begin{quote}

\sphinxAtStartPar
Reads a problem from COPT serialized data.
\end{quote}

\sphinxAtStartPar
\sphinxstylestrong{Arguments}
\begin{quote}

\sphinxAtStartPar
\sphinxcode{\sphinxupquote{prob}}
\begin{quote}

\sphinxAtStartPar
The COPT problem.
\end{quote}

\sphinxAtStartPar
\sphinxcode{\sphinxupquote{blob}}
\begin{quote}

\sphinxAtStartPar
Serialized data.
\end{quote}

\sphinxAtStartPar
\sphinxcode{\sphinxupquote{len}}
\begin{quote}

\sphinxAtStartPar
Length of serialized data.
\end{quote}
\end{quote}
\end{quote}

\subsubsection{COPT\_WriteMps}
\label{\detokenize{capiref:copt-writemps}}\begin{quote}

\sphinxAtStartPar
\sphinxstylestrong{Synopsis}
\begin{quote}

\sphinxAtStartPar
\sphinxcode{\sphinxupquote{int COPT\_WriteMps(copt\_prob *prob, const char *mpsfilename)}}
\end{quote}

\sphinxAtStartPar
\sphinxstylestrong{Description}
\begin{quote}

\sphinxAtStartPar
Writes the problem to a MPS file.
\end{quote}

\sphinxAtStartPar
\sphinxstylestrong{Arguments}
\begin{quote}

\sphinxAtStartPar
\sphinxcode{\sphinxupquote{prob}}
\begin{quote}

\sphinxAtStartPar
The COPT problem.
\end{quote}

\sphinxAtStartPar
\sphinxcode{\sphinxupquote{mpsfilename}}
\begin{quote}

\sphinxAtStartPar
The path to the MPS file.
\end{quote}
\end{quote}
\end{quote}

\subsubsection{COPT\_WriteMpsStr}
\label{\detokenize{capiref:copt-writempsstr}}\begin{quote}

\sphinxAtStartPar
\sphinxstylestrong{Synopsis}
\begin{quote}

\sphinxAtStartPar
\sphinxcode{\sphinxupquote{int COPT\_WriteMpsStr(copt\_prob *prob, char *str, int nStrSize, int *pReqSize)}}
\end{quote}

\sphinxAtStartPar
\sphinxstylestrong{Description}
\begin{quote}

\sphinxAtStartPar
Writes the problem to a string buffer as MPS format.
\end{quote}

\sphinxAtStartPar
\sphinxstylestrong{Arguments}
\begin{quote}

\sphinxAtStartPar
\sphinxcode{\sphinxupquote{prob}}
\begin{quote}

\sphinxAtStartPar
The COPT problem.
\end{quote}

\sphinxAtStartPar
\sphinxcode{\sphinxupquote{str}}
\begin{quote}

\sphinxAtStartPar
String buffer of MPS\sphinxhyphen{}format problem.
\end{quote}

\sphinxAtStartPar
\sphinxcode{\sphinxupquote{nStrSize}}
\begin{quote}

\sphinxAtStartPar
The size of string buffer.
\end{quote}

\sphinxAtStartPar
\sphinxcode{\sphinxupquote{pReqSize}}
\begin{quote}

\sphinxAtStartPar
Minimum space requirement of string buffer for problem.
\end{quote}
\end{quote}
\end{quote}

\subsubsection{COPT\_WriteLp}
\label{\detokenize{capiref:copt-writelp}}\begin{quote}

\sphinxAtStartPar
\sphinxstylestrong{Synopsis}
\begin{quote}

\sphinxAtStartPar
\sphinxcode{\sphinxupquote{int COPT\_WriteLp(copt\_prob *prob, const char *lpfilename)}}
\end{quote}

\sphinxAtStartPar
\sphinxstylestrong{Description}
\begin{quote}

\sphinxAtStartPar
Writes the problem to a LP file.
\end{quote}

\sphinxAtStartPar
\sphinxstylestrong{Arguments}
\begin{quote}

\sphinxAtStartPar
\sphinxcode{\sphinxupquote{prob}}
\begin{quote}

\sphinxAtStartPar
The COPT problem.
\end{quote}

\sphinxAtStartPar
\sphinxcode{\sphinxupquote{lpfilename}}
\begin{quote}

\sphinxAtStartPar
The path to the LP file.
\end{quote}
\end{quote}
\end{quote}

\subsubsection{COPT\_WriteCbf}
\label{\detokenize{capiref:copt-writecbf}}\begin{quote}

\sphinxAtStartPar
\sphinxstylestrong{Synopsis}
\begin{quote}

\sphinxAtStartPar
\sphinxcode{\sphinxupquote{int COPT\_WriteCbf(copt\_prob *prob, const char *cbffilename)}}
\end{quote}

\sphinxAtStartPar
\sphinxstylestrong{Description}
\begin{quote}

\sphinxAtStartPar
Writes the problem to a CBF format file.
\end{quote}

\sphinxAtStartPar
\sphinxstylestrong{Arguments}
\begin{quote}

\sphinxAtStartPar
\sphinxcode{\sphinxupquote{prob}}
\begin{quote}

\sphinxAtStartPar
The COPT problem.
\end{quote}

\sphinxAtStartPar
\sphinxcode{\sphinxupquote{cbffilename}}
\begin{quote}

\sphinxAtStartPar
The path to the CBF format file.
\end{quote}
\end{quote}
\end{quote}

\subsubsection{COPT\_WriteBin}
\label{\detokenize{capiref:copt-writebin}}\begin{quote}

\sphinxAtStartPar
\sphinxstylestrong{Synopsis}
\begin{quote}

\sphinxAtStartPar
\sphinxcode{\sphinxupquote{int COPT\_WriteBin(copt\_prob *prob, const char *binfilename)}}
\end{quote}

\sphinxAtStartPar
\sphinxstylestrong{Description}
\begin{quote}

\sphinxAtStartPar
Writes the problem to a COPT binary format file.
\end{quote}

\sphinxAtStartPar
\sphinxstylestrong{Arguments}
\begin{quote}

\sphinxAtStartPar
\sphinxcode{\sphinxupquote{prob}}
\begin{quote}

\sphinxAtStartPar
The COPT problem.
\end{quote}

\sphinxAtStartPar
\sphinxcode{\sphinxupquote{binfilename}}
\begin{quote}

\sphinxAtStartPar
The path to the COPT binary format file.
\end{quote}
\end{quote}
\end{quote}

\subsubsection{COPT\_WriteNL}
\label{\detokenize{capiref:copt-writenl}}\begin{quote}

\sphinxAtStartPar
\sphinxstylestrong{Synopsis}
\begin{quote}

\sphinxAtStartPar
\sphinxcode{\sphinxupquote{int COPT\_WriteNL(copt\_prob *prob, const char *nlfilename)}}
\end{quote}

\sphinxAtStartPar
\sphinxstylestrong{Description}
\begin{quote}

\sphinxAtStartPar
Export the internal COPT model to an NL format file.
\end{quote}

\sphinxAtStartPar
\sphinxstylestrong{Arguments}
\begin{quote}

\sphinxAtStartPar
\sphinxcode{\sphinxupquote{prob}}
\begin{quote}

\sphinxAtStartPar
The COPT problem.
\end{quote}

\sphinxAtStartPar
\sphinxcode{\sphinxupquote{nlfilename}}
\begin{quote}

\sphinxAtStartPar
Path to the NL format file.
\end{quote}
\end{quote}
\end{quote}

\subsubsection{COPT\_WriteBlob}
\label{\detokenize{capiref:copt-writeblob}}\begin{quote}

\sphinxAtStartPar
\sphinxstylestrong{Synopsis}
\begin{quote}

\sphinxAtStartPar
\sphinxcode{\sphinxupquote{int COPT\_WriteBlob(copt\_prob *prob, int tryCompress, void **p\_blob, COPT\_INT64 *pLen)}}
\end{quote}

\sphinxAtStartPar
\sphinxstylestrong{Description}
\begin{quote}

\sphinxAtStartPar
Writes the problem to COPT serialized data.
\end{quote}

\sphinxAtStartPar
\sphinxstylestrong{Arguments}
\begin{quote}

\sphinxAtStartPar
\sphinxcode{\sphinxupquote{prob}}
\begin{quote}

\sphinxAtStartPar
The COPT problem.
\end{quote}

\sphinxAtStartPar
\sphinxcode{\sphinxupquote{tryCompress}}
\begin{quote}

\sphinxAtStartPar
Whether to compress data.
\end{quote}

\sphinxAtStartPar
\sphinxcode{\sphinxupquote{p\_blob}}
\begin{quote}

\sphinxAtStartPar
Output pointer of serialized data.
\end{quote}

\sphinxAtStartPar
\sphinxcode{\sphinxupquote{pLen}}
\begin{quote}

\sphinxAtStartPar
Pointer to length of serialized data.
\end{quote}
\end{quote}
\end{quote}

\subsection{Solving the problem and accessing solutions}
\label{\detokenize{capiref:solving-the-problem-and-accessing-solutions}}

\subsubsection{COPT\_SolveLp}
\label{\detokenize{capiref:copt-solvelp}}\begin{quote}

\sphinxAtStartPar
\sphinxstylestrong{Synopsis}
\begin{quote}

\sphinxAtStartPar
\sphinxcode{\sphinxupquote{int COPT\_SolveLp(copt\_prob *prob)}}
\end{quote}

\sphinxAtStartPar
\sphinxstylestrong{Description}
\begin{quote}

\sphinxAtStartPar
Solves the LP, QP, QCP, SOCP and SDP problem.

\sphinxAtStartPar
If problem is a MIP, then integer restrictions on variables will be ignored,
and SOS constraints, indicator constraints will be discarded,
and the problem will be solved as a LP.
\end{quote}

\sphinxAtStartPar
\sphinxstylestrong{Arguments}
\begin{quote}

\sphinxAtStartPar
\sphinxcode{\sphinxupquote{prob}}
\begin{quote}

\sphinxAtStartPar
The COPT problem.
\end{quote}
\end{quote}
\end{quote}

\subsubsection{COPT\_Solve}
\label{\detokenize{capiref:copt-solve}}\begin{quote}

\sphinxAtStartPar
\sphinxstylestrong{Synopsis}
\begin{quote}

\sphinxAtStartPar
\sphinxcode{\sphinxupquote{int COPT\_Solve(copt\_prob *prob)}}
\end{quote}

\sphinxAtStartPar
\sphinxstylestrong{Description}
\begin{quote}

\sphinxAtStartPar
Solves the problem.
\end{quote}

\sphinxAtStartPar
\sphinxstylestrong{Arguments}
\begin{quote}

\sphinxAtStartPar
\sphinxcode{\sphinxupquote{prob}}
\begin{quote}

\sphinxAtStartPar
The COPT problem.
\end{quote}
\end{quote}
\end{quote}

\subsubsection{COPT\_GetSolution}
\label{\detokenize{capiref:copt-getsolution}}\begin{quote}

\sphinxAtStartPar
\sphinxstylestrong{Synopsis}
\begin{quote}

\sphinxAtStartPar
\sphinxcode{\sphinxupquote{int COPT\_GetSolution(copt\_prob *prob, double *colVal)}}
\end{quote}

\sphinxAtStartPar
\sphinxstylestrong{Description}
\begin{quote}

\sphinxAtStartPar
Obtains MIP solution.
\end{quote}

\sphinxAtStartPar
\sphinxstylestrong{Arguments}
\begin{quote}

\sphinxAtStartPar
\sphinxcode{\sphinxupquote{prob}}
\begin{quote}

\sphinxAtStartPar
The COPT problem.
\end{quote}

\sphinxAtStartPar
\sphinxcode{\sphinxupquote{colVal}}
\begin{quote}

\sphinxAtStartPar
Solution values of variables.
\end{quote}
\end{quote}
\end{quote}

\subsubsection{COPT\_GetPoolObjVal}
\label{\detokenize{capiref:copt-getpoolobjval}}\begin{quote}

\sphinxAtStartPar
\sphinxstylestrong{Synopsis}
\begin{quote}

\sphinxAtStartPar
\sphinxcode{\sphinxupquote{int COPT\_GetPoolObjVal(copt\_prob *prob, int iSol, double *p\_objVal)}}
\end{quote}

\sphinxAtStartPar
\sphinxstylestrong{Description}
\begin{quote}

\sphinxAtStartPar
Obtains the \sphinxcode{\sphinxupquote{iSol}} \sphinxhyphen{}th objective value in solution pool.
\end{quote}

\sphinxAtStartPar
\sphinxstylestrong{Arguments}
\begin{quote}

\sphinxAtStartPar
\sphinxcode{\sphinxupquote{prob}}
\begin{quote}

\sphinxAtStartPar
The COPT problem.
\end{quote}

\sphinxAtStartPar
\sphinxcode{\sphinxupquote{iSol}}
\begin{quote}

\sphinxAtStartPar
Index of solution.
\end{quote}

\sphinxAtStartPar
\sphinxcode{\sphinxupquote{p\_objVal}}
\begin{quote}

\sphinxAtStartPar
Pointer to objective value.
\end{quote}
\end{quote}
\end{quote}

\subsubsection{COPT\_GetPoolSolution}
\label{\detokenize{capiref:copt-getpoolsolution}}\begin{quote}

\sphinxAtStartPar
\sphinxstylestrong{Synopsis}
\begin{quote}

\sphinxAtStartPar
\sphinxcode{\sphinxupquote{int COPT\_GetPoolSolution(copt\_prob *prob, int iSol, int num, const int *list, double *colVal)}}
\end{quote}

\sphinxAtStartPar
\sphinxstylestrong{Description}
\begin{quote}

\sphinxAtStartPar
Obtains the \sphinxcode{\sphinxupquote{iSol}} \sphinxhyphen{}th solution.
\end{quote}

\sphinxAtStartPar
\sphinxstylestrong{Arguments}
\begin{quote}

\sphinxAtStartPar
\sphinxcode{\sphinxupquote{prob}}
\begin{quote}

\sphinxAtStartPar
The COPT problem.
\end{quote}

\sphinxAtStartPar
\sphinxcode{\sphinxupquote{iSol}}
\begin{quote}

\sphinxAtStartPar
Index of solution.
\end{quote}

\sphinxAtStartPar
\sphinxcode{\sphinxupquote{num}}
\begin{quote}

\sphinxAtStartPar
Number of columns.
\end{quote}

\sphinxAtStartPar
\sphinxcode{\sphinxupquote{list}}
\begin{quote}

\sphinxAtStartPar
Index of columns. Can be \sphinxcode{\sphinxupquote{NULL}}.
\end{quote}

\sphinxAtStartPar
\sphinxcode{\sphinxupquote{colVal}}
\begin{quote}

\sphinxAtStartPar
Array of solution.
\end{quote}
\end{quote}
\end{quote}

\subsubsection{COPT\_GetLpSolution}
\label{\detokenize{capiref:copt-getlpsolution}}\begin{quote}

\sphinxAtStartPar
\sphinxstylestrong{Synopsis}
\begin{quote}

\sphinxAtStartPar
\sphinxcode{\sphinxupquote{int COPT\_GetLpSolution(copt\_prob *prob, double *value, double *slack, double *rowDual, double *redCost)}}
\end{quote}

\sphinxAtStartPar
\sphinxstylestrong{Description}
\begin{quote}

\sphinxAtStartPar
Obtains LP, QP, QCP, SOCP and SDP solutions.

\sphinxAtStartPar
\sphinxstylestrong{Note:} For SDP, please use \sphinxcode{\sphinxupquote{COPT\_GetPSDColInfo}} to obtain primal/dual
solution of PSD variable.
\end{quote}

\sphinxAtStartPar
\sphinxstylestrong{Arguments}
\begin{quote}

\sphinxAtStartPar
\sphinxcode{\sphinxupquote{prob}}
\begin{quote}

\sphinxAtStartPar
The COPT problem.
\end{quote}

\sphinxAtStartPar
\sphinxcode{\sphinxupquote{value}}
\begin{quote}

\sphinxAtStartPar
Solution values of variables. Can be \sphinxcode{\sphinxupquote{NULL}}.
\end{quote}

\sphinxAtStartPar
\sphinxcode{\sphinxupquote{slack}}
\begin{quote}

\sphinxAtStartPar
Solution values of slack variables.
They are also known as activities of constraints.
Can be \sphinxcode{\sphinxupquote{NULL}}.
\end{quote}

\sphinxAtStartPar
\sphinxcode{\sphinxupquote{rowDual}}
\begin{quote}

\sphinxAtStartPar
Dual values of constraints. Can be \sphinxcode{\sphinxupquote{NULL}}.
\end{quote}

\sphinxAtStartPar
\sphinxcode{\sphinxupquote{redCost}}
\begin{quote}

\sphinxAtStartPar
Reduced cost of variables. Can be \sphinxcode{\sphinxupquote{NULL}}.
\end{quote}
\end{quote}
\end{quote}

\subsubsection{COPT\_SetLpSolution}
\label{\detokenize{capiref:copt-setlpsolution}}\begin{quote}

\sphinxAtStartPar
\sphinxstylestrong{Synopsis}
\begin{quote}

\sphinxAtStartPar
\sphinxcode{\sphinxupquote{int COPT\_SetLpSolution(copt\_prob *prob, double *value, double *slack, double *rowDual, double *redCost)}}
\end{quote}

\sphinxAtStartPar
\sphinxstylestrong{Description}
\begin{quote}

\sphinxAtStartPar
Set LP solution.
\end{quote}

\sphinxAtStartPar
\sphinxstylestrong{Arguments}
\begin{quote}

\sphinxAtStartPar
\sphinxcode{\sphinxupquote{prob}}
\begin{quote}

\sphinxAtStartPar
The COPT problem.
\end{quote}

\sphinxAtStartPar
\sphinxcode{\sphinxupquote{value}}
\begin{quote}

\sphinxAtStartPar
Solution values of variables.
\end{quote}

\sphinxAtStartPar
\sphinxcode{\sphinxupquote{slack}}
\begin{quote}

\sphinxAtStartPar
Solution values of slack variables.
\end{quote}

\sphinxAtStartPar
\sphinxcode{\sphinxupquote{rowDual}}
\begin{quote}

\sphinxAtStartPar
Dual values of constraints.
\end{quote}

\sphinxAtStartPar
\sphinxcode{\sphinxupquote{redCost}}
\begin{quote}

\sphinxAtStartPar
Reduced cost of variables.
\end{quote}
\end{quote}
\end{quote}

\subsubsection{COPT\_GetBasis}
\label{\detokenize{capiref:copt-getbasis}}\begin{quote}

\sphinxAtStartPar
\sphinxstylestrong{Synopsis}
\begin{quote}

\sphinxAtStartPar
\sphinxcode{\sphinxupquote{int COPT\_GetBasis(copt\_prob *prob, int *colBasis, int *rowBasis)}}
\end{quote}

\sphinxAtStartPar
\sphinxstylestrong{Description}
\begin{quote}

\sphinxAtStartPar
Obtains LP basis.
\end{quote}

\sphinxAtStartPar
\sphinxstylestrong{Arguments}
\begin{quote}

\sphinxAtStartPar
\sphinxcode{\sphinxupquote{prob}}
\begin{quote}

\sphinxAtStartPar
The COPT problem.
\end{quote}

\sphinxAtStartPar
\sphinxcode{\sphinxupquote{colBasis}} and \sphinxcode{\sphinxupquote{rowBasis}}
\begin{quote}

\sphinxAtStartPar
The basis status of variables and constraints.
Please refer to basis constants for possible
values and their meanings.
\end{quote}
\end{quote}
\end{quote}

\subsubsection{COPT\_SetBasis}
\label{\detokenize{capiref:copt-setbasis}}\begin{quote}

\sphinxAtStartPar
\sphinxstylestrong{Synopsis}
\begin{quote}

\sphinxAtStartPar
\sphinxcode{\sphinxupquote{int COPT\_SetBasis(copt\_prob *prob, const int *colBasis, const int *rowBasis)}}
\end{quote}

\sphinxAtStartPar
\sphinxstylestrong{Description}
\begin{quote}

\sphinxAtStartPar
Sets LP basis. It can be used to warm\sphinxhyphen{}start an LP optimization.
\end{quote}

\sphinxAtStartPar
\sphinxstylestrong{Arguments}
\begin{quote}

\sphinxAtStartPar
\sphinxcode{\sphinxupquote{prob}}
\begin{quote}

\sphinxAtStartPar
The COPT problem.
\end{quote}

\sphinxAtStartPar
\sphinxcode{\sphinxupquote{colBasis}} and \sphinxcode{\sphinxupquote{rowBasis}}
\begin{quote}

\sphinxAtStartPar
The basis status of variables and constraints.
Please refer to basis constants for possible
values and their meanings.
\end{quote}
\end{quote}
\end{quote}

\subsubsection{COPT\_SetSlackBasis}
\label{\detokenize{capiref:copt-setslackbasis}}\begin{quote}

\sphinxAtStartPar
\sphinxstylestrong{Synopsis}
\begin{quote}

\sphinxAtStartPar
\sphinxcode{\sphinxupquote{int COPT\_SetSlackBasis(copt\_prob *prob)}}
\end{quote}

\sphinxAtStartPar
\sphinxstylestrong{Description}
\begin{quote}

\sphinxAtStartPar
Sets a slack basis for LP. The slack basis is the default starting basis
for an LP problem. This API function can be used to
restore an LP problem to its starting basis.
\end{quote}

\sphinxAtStartPar
\sphinxstylestrong{Arguments}
\begin{quote}

\sphinxAtStartPar
\sphinxcode{\sphinxupquote{prob}}
\begin{quote}

\sphinxAtStartPar
The COPT problem.
\end{quote}
\end{quote}
\end{quote}

\subsubsection{COPT\_Reset}
\label{\detokenize{capiref:copt-reset}}\begin{quote}

\sphinxAtStartPar
\sphinxstylestrong{Synopsis}
\begin{quote}

\sphinxAtStartPar
\sphinxcode{\sphinxupquote{int COPT\_Reset(copt\_prob *prob, int iClearAll)}}
\end{quote}

\sphinxAtStartPar
\sphinxstylestrong{Description}
\begin{quote}

\sphinxAtStartPar
Reset basis and LP/MIP solution in problem, which forces next solve
start from scratch. If \sphinxcode{\sphinxupquote{iClearAll}} is \sphinxcode{\sphinxupquote{1}}, then clear additional
information such as MIP start as well.
\end{quote}

\sphinxAtStartPar
\sphinxstylestrong{Arguments}
\begin{quote}

\sphinxAtStartPar
\sphinxcode{\sphinxupquote{prob}}
\begin{quote}

\sphinxAtStartPar
The COPT problem.
\end{quote}

\sphinxAtStartPar
\sphinxcode{\sphinxupquote{iClearAll}}
\begin{quote}

\sphinxAtStartPar
Whether to clear additional information.
\end{quote}
\end{quote}
\end{quote}

\subsubsection{COPT\_ReadSol}
\label{\detokenize{capiref:copt-readsol}}\begin{quote}

\sphinxAtStartPar
\sphinxstylestrong{Synopsis}
\begin{quote}

\sphinxAtStartPar
\sphinxcode{\sphinxupquote{int COPT\_ReadSol(copt\_prob *prob, const char *solfilename)}}
\end{quote}

\sphinxAtStartPar
\sphinxstylestrong{Description}
\begin{quote}

\sphinxAtStartPar
Reads a MIP solution from file as MIP start information.

\sphinxAtStartPar
\sphinxstylestrong{Note:} The default solution value is 0, i.e. a partial solution
will be automatically filled in with zeros.
\end{quote}

\sphinxAtStartPar
\sphinxstylestrong{Arguments}
\begin{quote}

\sphinxAtStartPar
\sphinxcode{\sphinxupquote{prob}}
\begin{quote}

\sphinxAtStartPar
The COPT problem.
\end{quote}

\sphinxAtStartPar
\sphinxcode{\sphinxupquote{solfilename}}
\begin{quote}

\sphinxAtStartPar
The path to the solution file.
\end{quote}
\end{quote}
\end{quote}

\subsubsection{COPT\_ReadJsonSol}
\label{\detokenize{capiref:copt-readjsonsol}}\begin{quote}

\sphinxAtStartPar
\sphinxstylestrong{Synopsis}
\begin{quote}

\sphinxAtStartPar
\sphinxcode{\sphinxupquote{int COPT\_ReadJsonSol(copt\_prob *prob, const char *solfilename)}}
\end{quote}

\sphinxAtStartPar
\sphinxstylestrong{Description}
\begin{quote}

\sphinxAtStartPar
Read solution in format of JSON from file.
\end{quote}

\sphinxAtStartPar
\sphinxstylestrong{Arguments}
\begin{quote}

\sphinxAtStartPar
\sphinxcode{\sphinxupquote{prob}}
\begin{quote}

\sphinxAtStartPar
The COPT problem.
\end{quote}

\sphinxAtStartPar
\sphinxcode{\sphinxupquote{solfilename}}
\begin{quote}

\sphinxAtStartPar
The path to the solution file.
\end{quote}
\end{quote}
\end{quote}

\subsubsection{COPT\_WriteSol}
\label{\detokenize{capiref:copt-writesol}}\begin{quote}

\sphinxAtStartPar
\sphinxstylestrong{Synopsis}
\begin{quote}

\sphinxAtStartPar
\sphinxcode{\sphinxupquote{int COPT\_WriteSol(copt\_prob *prob, const char *solfilename)}}
\end{quote}

\sphinxAtStartPar
\sphinxstylestrong{Description}
\begin{quote}

\sphinxAtStartPar
Writes a LP/MIP solution to a file.
\end{quote}

\sphinxAtStartPar
\sphinxstylestrong{Arguments}
\begin{quote}

\sphinxAtStartPar
\sphinxcode{\sphinxupquote{prob}}
\begin{quote}

\sphinxAtStartPar
The COPT problem.
\end{quote}

\sphinxAtStartPar
\sphinxcode{\sphinxupquote{solfilename}}
\begin{quote}

\sphinxAtStartPar
The path to the solution file.
\end{quote}
\end{quote}
\end{quote}

\subsubsection{COPT\_WriteJsonSol}
\label{\detokenize{capiref:copt-writejsonsol}}\begin{quote}

\sphinxAtStartPar
\sphinxstylestrong{Synopsis}
\begin{quote}

\sphinxAtStartPar
\sphinxcode{\sphinxupquote{int COPT\_WriteJsonSol(copt\_prob *prob, const char *solfilename)}}
\end{quote}

\sphinxAtStartPar
\sphinxstylestrong{Description}
\begin{quote}

\sphinxAtStartPar
Writes solution to a file of type \sphinxcode{\sphinxupquote{".json"}}.
\end{quote}

\sphinxAtStartPar
\sphinxstylestrong{Arguments}
\begin{quote}

\sphinxAtStartPar
\sphinxcode{\sphinxupquote{prob}}
\begin{quote}

\sphinxAtStartPar
The COPT problem.
\end{quote}

\sphinxAtStartPar
\sphinxcode{\sphinxupquote{solfilename}}
\begin{quote}

\sphinxAtStartPar
The path to the solution file.
\end{quote}
\end{quote}
\end{quote}

\subsubsection{COPT\_WritePoolSol}
\label{\detokenize{capiref:copt-writepoolsol}}\begin{quote}

\sphinxAtStartPar
\sphinxstylestrong{Synopsis}
\begin{quote}

\sphinxAtStartPar
\sphinxcode{\sphinxupquote{int COPT\_WritePoolSol(copt\_prob *prob, int iSol, const char *solfilename)}}
\end{quote}

\sphinxAtStartPar
\sphinxstylestrong{Description}
\begin{quote}

\sphinxAtStartPar
Writes selected pool solution to a file.
\end{quote}

\sphinxAtStartPar
\sphinxstylestrong{Arguments}
\begin{quote}

\sphinxAtStartPar
\sphinxcode{\sphinxupquote{prob}}
\begin{quote}

\sphinxAtStartPar
The COPT problem.
\end{quote}

\sphinxAtStartPar
\sphinxcode{\sphinxupquote{iSol}}
\begin{quote}

\sphinxAtStartPar
Index of pool solution.
\end{quote}

\sphinxAtStartPar
\sphinxcode{\sphinxupquote{solfilename}}
\begin{quote}

\sphinxAtStartPar
The path to the solution file.
\end{quote}
\end{quote}
\end{quote}

\subsubsection{COPT\_WriteBasis}
\label{\detokenize{capiref:copt-writebasis}}\begin{quote}

\sphinxAtStartPar
\sphinxstylestrong{Synopsis}
\begin{quote}

\sphinxAtStartPar
\sphinxcode{\sphinxupquote{int COPT\_WriteBasis(copt\_prob *prob, const char *basfilename)}}
\end{quote}

\sphinxAtStartPar
\sphinxstylestrong{Description}
\begin{quote}

\sphinxAtStartPar
Writes the internal LP basis to a file.
\end{quote}

\sphinxAtStartPar
\sphinxstylestrong{Arguments}
\begin{quote}

\sphinxAtStartPar
\sphinxcode{\sphinxupquote{prob}}
\begin{quote}

\sphinxAtStartPar
The COPT problem.
\end{quote}

\sphinxAtStartPar
\sphinxcode{\sphinxupquote{basfilename}}
\begin{quote}

\sphinxAtStartPar
The path to the basis file.
\end{quote}
\end{quote}
\end{quote}

\subsubsection{COPT\_ReadBasis}
\label{\detokenize{capiref:copt-readbasis}}\begin{quote}

\sphinxAtStartPar
\sphinxstylestrong{Synopsis}
\begin{quote}

\sphinxAtStartPar
\sphinxcode{\sphinxupquote{int COPT\_ReadBasis(copt\_prob *prob, const char *basfilename)}}
\end{quote}

\sphinxAtStartPar
\sphinxstylestrong{Description}
\begin{quote}

\sphinxAtStartPar
Reads the LP basis from a file.
It can be used to warm\sphinxhyphen{}start an LP optimization.
\end{quote}

\sphinxAtStartPar
\sphinxstylestrong{Arguments}
\begin{quote}

\sphinxAtStartPar
\sphinxcode{\sphinxupquote{prob}}
\begin{quote}

\sphinxAtStartPar
The COPT problem.
\end{quote}

\sphinxAtStartPar
\sphinxcode{\sphinxupquote{basfilename}}
\begin{quote}

\sphinxAtStartPar
The path to the basis file.
\end{quote}
\end{quote}
\end{quote}

\subsection{Accessing information of problem}
\label{\detokenize{capiref:accessing-information-of-problem}}\label{\detokenize{capiref:chapapi-getinfo}}

\subsubsection{COPT\_GetCols}
\label{\detokenize{capiref:copt-getcols}}\begin{quote}

\sphinxAtStartPar
\sphinxstylestrong{Synopsis}
\begin{quote}

\sphinxAtStartPar
\sphinxcode{\sphinxupquote{int COPT\_GetCols(}}
\begin{quote}

\sphinxAtStartPar
\sphinxcode{\sphinxupquote{copt\_prob *prob,}}

\sphinxAtStartPar
\sphinxcode{\sphinxupquote{int nCol,}}

\sphinxAtStartPar
\sphinxcode{\sphinxupquote{const int *list,}}

\sphinxAtStartPar
\sphinxcode{\sphinxupquote{int *colMatBeg,}}

\sphinxAtStartPar
\sphinxcode{\sphinxupquote{int *colMatCnt,}}

\sphinxAtStartPar
\sphinxcode{\sphinxupquote{int *colMatIdx,}}

\sphinxAtStartPar
\sphinxcode{\sphinxupquote{double *colMatElem,}}

\sphinxAtStartPar
\sphinxcode{\sphinxupquote{int nElemSize,}}

\sphinxAtStartPar
\sphinxcode{\sphinxupquote{int *pReqSize)}}
\end{quote}
\end{quote}

\sphinxAtStartPar
\sphinxstylestrong{Description}
\begin{quote}

\sphinxAtStartPar
Extract coefficient matrix by columns.

\sphinxAtStartPar
In general, users need to call this function twice to accomplish the task.
Firstly, by passing \sphinxcode{\sphinxupquote{NULL}} to arguments \sphinxcode{\sphinxupquote{colMatBeg}}, \sphinxcode{\sphinxupquote{colMatCnt}},
\sphinxcode{\sphinxupquote{colMatIdx}} and \sphinxcode{\sphinxupquote{colMatElem}}, we get number of non\sphinxhyphen{}zeros elements by
\sphinxcode{\sphinxupquote{pReqSize}} specified by \sphinxcode{\sphinxupquote{nCol}} and \sphinxcode{\sphinxupquote{list}}. Secondly, allocate
sufficient memory for CCS\sphinxhyphen{}format matrix and call this function again to
extract coefficient matrix. If the memory of coefficient matrix passed to
function is not sufficient, then return the first \sphinxcode{\sphinxupquote{nElemSize}} non\sphinxhyphen{}zero elements,
and the minimal required length of non\sphinxhyphen{}zero elements by \sphinxcode{\sphinxupquote{pReqSize}}.
If \sphinxcode{\sphinxupquote{list}} is \sphinxcode{\sphinxupquote{NULL}}, then the first \sphinxcode{\sphinxupquote{nCol}} columns will be returned.
\end{quote}

\sphinxAtStartPar
\sphinxstylestrong{Arguments}
\begin{quote}

\sphinxAtStartPar
\sphinxcode{\sphinxupquote{prob}}
\begin{quote}

\sphinxAtStartPar
The COPT problem.
\end{quote}

\sphinxAtStartPar
\sphinxcode{\sphinxupquote{nCol}}
\begin{quote}

\sphinxAtStartPar
Number of columns.
\end{quote}

\sphinxAtStartPar
\sphinxcode{\sphinxupquote{list}}
\begin{quote}

\sphinxAtStartPar
Index of columns. Can be \sphinxcode{\sphinxupquote{NULL}}.
\end{quote}

\sphinxAtStartPar
\sphinxcode{\sphinxupquote{colMatBeg, colMatCnt, colMatIdx}} and \sphinxcode{\sphinxupquote{colMatElem}}
\begin{quote}

\sphinxAtStartPar
Defines the coefficient matrix in compressed column storage (CCS) format.
Please see \sphinxstylestrong{other information} of \sphinxcode{\sphinxupquote{COPT\_LoadProb}} for
an example of the CCS format.
\end{quote}

\sphinxAtStartPar
\sphinxcode{\sphinxupquote{nElemSize}}
\begin{quote}

\sphinxAtStartPar
Length of array for non\sphinxhyphen{}zero coefficients.
\end{quote}

\sphinxAtStartPar
\sphinxcode{\sphinxupquote{pReqSize}}
\begin{quote}

\sphinxAtStartPar
Pointer to minimal length of array for non\sphinxhyphen{}zero coefficients.
Can be \sphinxcode{\sphinxupquote{NULL}}.
\end{quote}
\end{quote}
\end{quote}

\subsubsection{COPT\_GetPSDCols}
\label{\detokenize{capiref:copt-getpsdcols}}\begin{quote}

\sphinxAtStartPar
\sphinxstylestrong{Synopsis}
\begin{quote}

\sphinxAtStartPar
\sphinxcode{\sphinxupquote{int COPT\_GetPSDCols(copt\_prob *prob, int nCol, int *list, int* colDims, int *colLens)}}
\end{quote}

\sphinxAtStartPar
\sphinxstylestrong{Description}
\begin{quote}

\sphinxAtStartPar
Get dimension and flattened length of \sphinxcode{\sphinxupquote{nCol}} PSD variables.
\end{quote}

\sphinxAtStartPar
\sphinxstylestrong{Arguments}
\begin{quote}

\sphinxAtStartPar
\sphinxcode{\sphinxupquote{prob}}
\begin{quote}

\sphinxAtStartPar
The COPT problem.
\end{quote}

\sphinxAtStartPar
\sphinxcode{\sphinxupquote{nCol}}
\begin{quote}

\sphinxAtStartPar
Number of PSD variables.
\end{quote}

\sphinxAtStartPar
\sphinxcode{\sphinxupquote{list}}
\begin{quote}

\sphinxAtStartPar
Index of PSD variables. Can be \sphinxcode{\sphinxupquote{NULL}}.
\end{quote}

\sphinxAtStartPar
\sphinxcode{\sphinxupquote{colDims}}
\begin{quote}

\sphinxAtStartPar
Dimension of PSD variables.
\end{quote}

\sphinxAtStartPar
\sphinxcode{\sphinxupquote{colLens}}
\begin{quote}

\sphinxAtStartPar
Flattened length of PSD variables.
\end{quote}
\end{quote}
\end{quote}

\subsubsection{COPT\_GetRows}
\label{\detokenize{capiref:copt-getrows}}\begin{quote}

\sphinxAtStartPar
\sphinxstylestrong{Synopsis}
\begin{quote}

\sphinxAtStartPar
\sphinxcode{\sphinxupquote{int COPT\_GetRows(}}
\begin{quote}

\sphinxAtStartPar
\sphinxcode{\sphinxupquote{copt\_prob *prob,}}

\sphinxAtStartPar
\sphinxcode{\sphinxupquote{int nRow,}}

\sphinxAtStartPar
\sphinxcode{\sphinxupquote{const int *list,}}

\sphinxAtStartPar
\sphinxcode{\sphinxupquote{int *rowMatBeg,}}

\sphinxAtStartPar
\sphinxcode{\sphinxupquote{int *rowMatCnt,}}

\sphinxAtStartPar
\sphinxcode{\sphinxupquote{int *rowMatIdx,}}

\sphinxAtStartPar
\sphinxcode{\sphinxupquote{double *rowMatElem,}}

\sphinxAtStartPar
\sphinxcode{\sphinxupquote{int nElemSize,}}

\sphinxAtStartPar
\sphinxcode{\sphinxupquote{int *pReqSize)}}
\end{quote}
\end{quote}

\sphinxAtStartPar
\sphinxstylestrong{Description}
\begin{quote}

\sphinxAtStartPar
Extract coefficient matrix by rows.

\sphinxAtStartPar
In general, users need to call this function twice to accomplish the task.
Firstly, by passing \sphinxcode{\sphinxupquote{NULL}} to arguments \sphinxcode{\sphinxupquote{rowMatBeg}}, \sphinxcode{\sphinxupquote{rowMatCnt}},
\sphinxcode{\sphinxupquote{rowMatIdx}} and \sphinxcode{\sphinxupquote{rowMatElem}}, we get number of non\sphinxhyphen{}zeros elements by
\sphinxcode{\sphinxupquote{pReqSize}} specified by \sphinxcode{\sphinxupquote{nRow}} and \sphinxcode{\sphinxupquote{list}}. Secondly, allocate
sufficient memory for CRS\sphinxhyphen{}format matrix and call this function again to
extract coefficient matrix. If the memory of coefficient matrix passed to
function is not sufficient, then return the first \sphinxcode{\sphinxupquote{nElemSize}} non\sphinxhyphen{}zero elements,
and the minimal required length of non\sphinxhyphen{}zero elements by \sphinxcode{\sphinxupquote{pReqSize}}.
If \sphinxcode{\sphinxupquote{list}} is \sphinxcode{\sphinxupquote{NULL}}, then the first \sphinxcode{\sphinxupquote{nRow}} rows will be returned.
\end{quote}

\sphinxAtStartPar
\sphinxstylestrong{Arguments}
\begin{quote}

\sphinxAtStartPar
\sphinxcode{\sphinxupquote{prob}}
\begin{quote}

\sphinxAtStartPar
The COPT problem.
\end{quote}

\sphinxAtStartPar
\sphinxcode{\sphinxupquote{nRow}}
\begin{quote}

\sphinxAtStartPar
Number of rows.
\end{quote}

\sphinxAtStartPar
\sphinxcode{\sphinxupquote{list}}
\begin{quote}

\sphinxAtStartPar
Index of rows. Can be \sphinxcode{\sphinxupquote{NULL}}.
\end{quote}

\sphinxAtStartPar
\sphinxcode{\sphinxupquote{rowMatBeg, rowMatCnt, rowMatIdx}} and \sphinxcode{\sphinxupquote{rowMatElem}}
\begin{quote}

\sphinxAtStartPar
Defines the coefficient matrix in compressed row storage (CRS) format.
Please see \sphinxstylestrong{other information} of \sphinxcode{\sphinxupquote{COPT\_LoadProb}} for
an example of the CRS format.
\end{quote}

\sphinxAtStartPar
\sphinxcode{\sphinxupquote{nElemSize}}
\begin{quote}

\sphinxAtStartPar
Length of array for non\sphinxhyphen{}zero coefficients.
\end{quote}

\sphinxAtStartPar
\sphinxcode{\sphinxupquote{pReqSize}}
\begin{quote}

\sphinxAtStartPar
Pointer to minimal length of array for non\sphinxhyphen{}zero coefficients.
Can be \sphinxcode{\sphinxupquote{NULL}}.
\end{quote}
\end{quote}
\end{quote}

\subsubsection{COPT\_GetElem}
\label{\detokenize{capiref:copt-getelem}}\begin{quote}

\sphinxAtStartPar
\sphinxstylestrong{Synopsis}
\begin{quote}

\sphinxAtStartPar
\sphinxcode{\sphinxupquote{int COPT\_GetElem(copt\_prob *prob, int iCol, int iRow, double *p\_elem)}}
\end{quote}

\sphinxAtStartPar
\sphinxstylestrong{Description}
\begin{quote}

\sphinxAtStartPar
Get coefficient of specified row and column.
\end{quote}

\sphinxAtStartPar
\sphinxstylestrong{Arguments}
\begin{quote}

\sphinxAtStartPar
\sphinxcode{\sphinxupquote{prob}}
\begin{quote}

\sphinxAtStartPar
The COPT problem.
\end{quote}

\sphinxAtStartPar
\sphinxcode{\sphinxupquote{iCol}}
\begin{quote}

\sphinxAtStartPar
Column index.
\end{quote}

\sphinxAtStartPar
\sphinxcode{\sphinxupquote{iRow}}
\begin{quote}

\sphinxAtStartPar
Row index.
\end{quote}

\sphinxAtStartPar
\sphinxcode{\sphinxupquote{p\_elem}}
\begin{quote}

\sphinxAtStartPar
Pointer to requested coefficient.
\end{quote}
\end{quote}
\end{quote}

\subsubsection{COPT\_GetPSDElem}
\label{\detokenize{capiref:copt-getpsdelem}}\begin{quote}

\sphinxAtStartPar
\sphinxstylestrong{Synopsis}
\begin{quote}

\sphinxAtStartPar
\sphinxcode{\sphinxupquote{int COPT\_GetPSDElem(copt\_prob *prob, int iCol, int iRow, int *p\_idx)}}
\end{quote}

\sphinxAtStartPar
\sphinxstylestrong{Description}
\begin{quote}

\sphinxAtStartPar
Get symmetric matrix index of specified PSD constraint and PSD variable.
\end{quote}

\sphinxAtStartPar
\sphinxstylestrong{Arguments}
\begin{quote}

\sphinxAtStartPar
\sphinxcode{\sphinxupquote{prob}}
\begin{quote}

\sphinxAtStartPar
The COPT problem.
\end{quote}

\sphinxAtStartPar
\sphinxcode{\sphinxupquote{iCol}}
\begin{quote}

\sphinxAtStartPar
PSD variable index.
\end{quote}

\sphinxAtStartPar
\sphinxcode{\sphinxupquote{iRow}}
\begin{quote}

\sphinxAtStartPar
PSD constraint index.
\end{quote}

\sphinxAtStartPar
\sphinxcode{\sphinxupquote{p\_idx}}
\begin{quote}

\sphinxAtStartPar
Pointer to requested symmetric matrix index.
\end{quote}
\end{quote}
\end{quote}

\subsubsection{COPT\_GetLMIElem}
\label{\detokenize{capiref:copt-getlmielem}}\begin{quote}

\sphinxAtStartPar
\sphinxstylestrong{Synopsis}
\begin{quote}

\sphinxAtStartPar
\sphinxcode{\sphinxupquote{int COPT\_GetLMIElem(copt\_prob *prob, int iCol, int iRow, int *p\_idx)}}
\end{quote}

\sphinxAtStartPar
\sphinxstylestrong{Description}
\begin{quote}

\sphinxAtStartPar
Get symmetric matrix index of specified LMI constraint and scalar variable.
\end{quote}

\sphinxAtStartPar
\sphinxstylestrong{Arguments}
\begin{quote}

\sphinxAtStartPar
\sphinxcode{\sphinxupquote{prob}}
\begin{quote}

\sphinxAtStartPar
The COPT problem.
\end{quote}

\sphinxAtStartPar
\sphinxcode{\sphinxupquote{iCol}}
\begin{quote}

\sphinxAtStartPar
Scalar variable index.
\end{quote}

\sphinxAtStartPar
\sphinxcode{\sphinxupquote{iRow}}
\begin{quote}

\sphinxAtStartPar
LMI constraint index.
\end{quote}

\sphinxAtStartPar
\sphinxcode{\sphinxupquote{p\_idx}}
\begin{quote}

\sphinxAtStartPar
Pointer to requested coefficient matrix index.
\end{quote}
\end{quote}
\end{quote}

\subsubsection{COPT\_GetSymMat}
\label{\detokenize{capiref:copt-getsymmat}}\begin{quote}

\sphinxAtStartPar
\sphinxstylestrong{Synopsis}
\begin{quote}

\sphinxAtStartPar
\sphinxcode{\sphinxupquote{int COPT\_GetSymMat(}}
\begin{quote}

\sphinxAtStartPar
\sphinxcode{\sphinxupquote{copt\_prob *prob,}}

\sphinxAtStartPar
\sphinxcode{\sphinxupquote{int iMat,}}

\sphinxAtStartPar
\sphinxcode{\sphinxupquote{int *p\_nDim,}}

\sphinxAtStartPar
\sphinxcode{\sphinxupquote{int *p\_nElem,}}

\sphinxAtStartPar
\sphinxcode{\sphinxupquote{int *rows,}}

\sphinxAtStartPar
\sphinxcode{\sphinxupquote{int *cols,}}

\sphinxAtStartPar
\sphinxcode{\sphinxupquote{double *elems)}}
\end{quote}
\end{quote}

\sphinxAtStartPar
\sphinxstylestrong{Description}
\begin{quote}

\sphinxAtStartPar
Get specified symmetric matrix.

\sphinxAtStartPar
In general, users need to call this function twice to accomplish the task.
Firstly, by passing \sphinxcode{\sphinxupquote{NULL}} to arguments \sphinxcode{\sphinxupquote{rows}}, \sphinxcode{\sphinxupquote{cols}} and \sphinxcode{\sphinxupquote{elems}},
we get dimension and number of non\sphinxhyphen{}zeros of symmetric matrix by \sphinxcode{\sphinxupquote{p\_nDim}}
and \sphinxcode{\sphinxupquote{p\_nElem}}, then allocate enough memory for \sphinxcode{\sphinxupquote{rows}}, \sphinxcode{\sphinxupquote{cols}} and
\sphinxcode{\sphinxupquote{elems}} and call this function to get the data of symmetric matrix.
\end{quote}

\sphinxAtStartPar
\sphinxstylestrong{Arguments}
\begin{quote}

\sphinxAtStartPar
\sphinxcode{\sphinxupquote{prob}}
\begin{quote}

\sphinxAtStartPar
The COPT problem.
\end{quote}

\sphinxAtStartPar
\sphinxcode{\sphinxupquote{iMat}}
\begin{quote}

\sphinxAtStartPar
Symmetric matrix index.
\end{quote}

\sphinxAtStartPar
\sphinxcode{\sphinxupquote{p\_nDim}}
\begin{quote}

\sphinxAtStartPar
Pointer to dimension of symmetric matrix.
\end{quote}

\sphinxAtStartPar
\sphinxcode{\sphinxupquote{p\_nElem}}
\begin{quote}

\sphinxAtStartPar
Pointer to number of nonzeros of symmetric matrix.
\end{quote}

\sphinxAtStartPar
\sphinxcode{\sphinxupquote{rows}}
\begin{quote}

\sphinxAtStartPar
Row index of symmetric matrix.
\end{quote}

\sphinxAtStartPar
\sphinxcode{\sphinxupquote{cols}}
\begin{quote}

\sphinxAtStartPar
Column index of symmetric matrix.
\end{quote}

\sphinxAtStartPar
\sphinxcode{\sphinxupquote{elems}}
\begin{quote}

\sphinxAtStartPar
Nonzero elements of symmetric matrix.
\end{quote}
\end{quote}
\end{quote}

\subsubsection{COPT\_GetQuadObj}
\label{\detokenize{capiref:copt-getquadobj}}\begin{quote}

\sphinxAtStartPar
\sphinxstylestrong{Synopsis}
\begin{quote}

\sphinxAtStartPar
\sphinxcode{\sphinxupquote{int COPT\_GetQuadObj(copt\_prob* prob, int* p\_nQElem, int* qRow, int* qCol, double* qElem)}}
\end{quote}

\sphinxAtStartPar
\sphinxstylestrong{Description}
\begin{quote}

\sphinxAtStartPar
Get the quadratic terms of the quadratic objective function.
\end{quote}

\sphinxAtStartPar
\sphinxstylestrong{Arguments}
\begin{quote}

\sphinxAtStartPar
\sphinxcode{\sphinxupquote{prob}}
\begin{quote}

\sphinxAtStartPar
The COPT problem.
\end{quote}

\sphinxAtStartPar
\sphinxcode{\sphinxupquote{p\_nQElem}}
\begin{quote}

\sphinxAtStartPar
Pointer to number of non\sphinxhyphen{}zero quadratic terms .
Can be \sphinxcode{\sphinxupquote{NULL}}.
\end{quote}

\sphinxAtStartPar
\sphinxcode{\sphinxupquote{qRow}}
\begin{quote}

\sphinxAtStartPar
Row index of non\sphinxhyphen{}zero quadratic terms of the quadratic objective function.
\end{quote}

\sphinxAtStartPar
\sphinxcode{\sphinxupquote{qCol}}
\begin{quote}

\sphinxAtStartPar
Column index of non\sphinxhyphen{}zero quadratic terms of the quadratic objective function.
\end{quote}

\sphinxAtStartPar
\sphinxcode{\sphinxupquote{qElem}}
\begin{quote}

\sphinxAtStartPar
Values of non\sphinxhyphen{}zero quadratic terms of the quadratic objective function.
\end{quote}
\end{quote}
\end{quote}

\subsubsection{COPT\_GetPSDObj}
\label{\detokenize{capiref:copt-getpsdobj}}\begin{quote}

\sphinxAtStartPar
\sphinxstylestrong{Synopsis}
\begin{quote}

\sphinxAtStartPar
\sphinxcode{\sphinxupquote{int COPT\_GetPSDObj(copt\_prob *prob, int iCol, int *p\_idx)}}
\end{quote}

\sphinxAtStartPar
\sphinxstylestrong{Description}
\begin{quote}

\sphinxAtStartPar
Get the specified PSD term of objective function.
\end{quote}

\sphinxAtStartPar
\sphinxstylestrong{Arguments}
\begin{quote}

\sphinxAtStartPar
\sphinxcode{\sphinxupquote{prob}}
\begin{quote}

\sphinxAtStartPar
The COPT problem.
\end{quote}

\sphinxAtStartPar
\sphinxcode{\sphinxupquote{iCol}}
\begin{quote}

\sphinxAtStartPar
PSD variable index.
\end{quote}

\sphinxAtStartPar
\sphinxcode{\sphinxupquote{p\_idx}}
\begin{quote}

\sphinxAtStartPar
Pointer to symmetric matrix index.
\end{quote}
\end{quote}
\end{quote}

\subsubsection{COPT\_GetNLObj}
\label{\detokenize{capiref:copt-getnlobj}}\begin{quote}

\sphinxAtStartPar
\sphinxstylestrong{Synopsis}
\begin{quote}

\sphinxAtStartPar
\sphinxcode{\sphinxupquote{int COPT\_GetNLObj(}}
\begin{quote}

\sphinxAtStartPar
\sphinxcode{\sphinxupquote{copt\_prob *prob,}}

\sphinxAtStartPar
\sphinxcode{\sphinxupquote{int *p\_nToken,}}

\sphinxAtStartPar
\sphinxcode{\sphinxupquote{int *p\_nTokenElem,}}

\sphinxAtStartPar
\sphinxcode{\sphinxupquote{int *token,}}

\sphinxAtStartPar
\sphinxcode{\sphinxupquote{double *tokenElem)}}
\end{quote}
\end{quote}

\sphinxAtStartPar
\sphinxstylestrong{Description}
\begin{quote}

\sphinxAtStartPar
Retrieve the nonlinear expression terms in the objective function.
\end{quote}

\sphinxAtStartPar
\sphinxstylestrong{Arguments}
\begin{quote}

\sphinxAtStartPar
\sphinxcode{\sphinxupquote{prob}}
\begin{quote}

\sphinxAtStartPar
The COPT problem.
\end{quote}

\sphinxAtStartPar
\sphinxcode{\sphinxupquote{p\_nToken}}
\begin{quote}

\sphinxAtStartPar
Pointer to the number of tokens in the expression.
\end{quote}

\sphinxAtStartPar
\sphinxcode{\sphinxupquote{p\_nTokenElem}}
\begin{quote}

\sphinxAtStartPar
Pointer to the number of constants in the expression.
\end{quote}

\sphinxAtStartPar
\sphinxcode{\sphinxupquote{token}}
\begin{quote}

\sphinxAtStartPar
Array of tokens in the expression.
\end{quote}

\sphinxAtStartPar
\sphinxcode{\sphinxupquote{tokenElem}}
\begin{quote}

\sphinxAtStartPar
Array of constants in the expression.
\end{quote}
\end{quote}
\end{quote}

\subsubsection{COPT\_GetSOSs}
\label{\detokenize{capiref:copt-getsoss}}\begin{quote}

\sphinxAtStartPar
\sphinxstylestrong{Synopsis}
\begin{quote}

\sphinxAtStartPar
\sphinxcode{\sphinxupquote{int COPT\_GetSOSs(}}
\begin{quote}

\sphinxAtStartPar
\sphinxcode{\sphinxupquote{copt\_prob *prob,}}

\sphinxAtStartPar
\sphinxcode{\sphinxupquote{int nSos,}}

\sphinxAtStartPar
\sphinxcode{\sphinxupquote{const int *list,}}

\sphinxAtStartPar
\sphinxcode{\sphinxupquote{int *sosMatBeg,}}

\sphinxAtStartPar
\sphinxcode{\sphinxupquote{int *sosMatCnt,}}

\sphinxAtStartPar
\sphinxcode{\sphinxupquote{int *sosMatIdx,}}

\sphinxAtStartPar
\sphinxcode{\sphinxupquote{double *sosMatWt,}}

\sphinxAtStartPar
\sphinxcode{\sphinxupquote{int nElemSize,}}

\sphinxAtStartPar
\sphinxcode{\sphinxupquote{int *pReqSize)}}
\end{quote}
\end{quote}

\sphinxAtStartPar
\sphinxstylestrong{Description}
\begin{quote}

\sphinxAtStartPar
Get the weight matrix of SOS constraints.

\sphinxAtStartPar
In general, users need to call this function twice to accomplish the task.
Firstly, by passing \sphinxcode{\sphinxupquote{NULL}} to arguments \sphinxcode{\sphinxupquote{sosMatBeg}}, \sphinxcode{\sphinxupquote{sosMatCnt}},
\sphinxcode{\sphinxupquote{sosMatIdx}} and \sphinxcode{\sphinxupquote{sosMatWt}}, we get number of non\sphinxhyphen{}zeros elements by
\sphinxcode{\sphinxupquote{pReqSize}} specified by \sphinxcode{\sphinxupquote{nSos}} and \sphinxcode{\sphinxupquote{list}}. Secondly, allocate
sufficient memory for CRS\sphinxhyphen{}format matrix and call this function again to
extract weight matrix. If the memory of weight matrix passed to function
is not sufficient, then return the first \sphinxcode{\sphinxupquote{nElemSize}} non\sphinxhyphen{}zero elements,
and the minimal required length of non\sphinxhyphen{}zero elements by \sphinxcode{\sphinxupquote{pReqSize}}.
If \sphinxcode{\sphinxupquote{list}} is \sphinxcode{\sphinxupquote{NULL}}, then the first \sphinxcode{\sphinxupquote{nSos}} rows will be returned.
\end{quote}

\sphinxAtStartPar
\sphinxstylestrong{Arguments}
\begin{quote}

\sphinxAtStartPar
\sphinxcode{\sphinxupquote{prob}}
\begin{quote}

\sphinxAtStartPar
The COPT problem.
\end{quote}

\sphinxAtStartPar
\sphinxcode{\sphinxupquote{nSos}}
\begin{quote}

\sphinxAtStartPar
Number of SOS constraints.
\end{quote}

\sphinxAtStartPar
\sphinxcode{\sphinxupquote{list}}
\begin{quote}

\sphinxAtStartPar
Index of SOS constraints. Can be \sphinxcode{\sphinxupquote{NULL}}.
\end{quote}

\sphinxAtStartPar
\sphinxcode{\sphinxupquote{sosMatBeg, sosMatCnt, sosMatIdx}} and \sphinxcode{\sphinxupquote{sosMatWt}}
\begin{quote}

\sphinxAtStartPar
Defines the weight matrix of SOS constraints in compressed row
storage (CRS) format. Please see \sphinxstylestrong{other information} of
\sphinxcode{\sphinxupquote{COPT\_LoadProb}} for an example of the CRS format.
\end{quote}

\sphinxAtStartPar
\sphinxcode{\sphinxupquote{nElemSize}}
\begin{quote}

\sphinxAtStartPar
Length of array for non\sphinxhyphen{}zero weights.
\end{quote}

\sphinxAtStartPar
\sphinxcode{\sphinxupquote{pReqSize}}
\begin{quote}

\sphinxAtStartPar
Pointer to minimal length of array for non\sphinxhyphen{}zero weights.
Can be \sphinxcode{\sphinxupquote{NULL}}.
\end{quote}
\end{quote}
\end{quote}

\subsubsection{COPT\_GetCones}
\label{\detokenize{capiref:copt-getcones}}\begin{quote}

\sphinxAtStartPar
\sphinxstylestrong{Synopsis}
\begin{quote}

\sphinxAtStartPar
\sphinxcode{\sphinxupquote{int COPT\_GetCones(}}
\begin{quote}

\sphinxAtStartPar
\sphinxcode{\sphinxupquote{copt\_prob *prob,}}

\sphinxAtStartPar
\sphinxcode{\sphinxupquote{int nCone,}}

\sphinxAtStartPar
\sphinxcode{\sphinxupquote{const int *list,}}

\sphinxAtStartPar
\sphinxcode{\sphinxupquote{int *coneBeg,}}

\sphinxAtStartPar
\sphinxcode{\sphinxupquote{int *coneCnt,}}

\sphinxAtStartPar
\sphinxcode{\sphinxupquote{int *coneIdx,}}

\sphinxAtStartPar
\sphinxcode{\sphinxupquote{int nElemSize,}}

\sphinxAtStartPar
\sphinxcode{\sphinxupquote{int *pReqSize)}}
\end{quote}
\end{quote}

\sphinxAtStartPar
\sphinxstylestrong{Description}
\begin{quote}

\sphinxAtStartPar
Get the matrix of Second\sphinxhyphen{}Order\sphinxhyphen{}Cone (SOC) constraints.

\sphinxAtStartPar
In general, users need to call this function twice to accomplish the task.
Firstly, by passing \sphinxcode{\sphinxupquote{NULL}} to arguments \sphinxcode{\sphinxupquote{coneBeg}}, \sphinxcode{\sphinxupquote{coneCnt}} and \sphinxcode{\sphinxupquote{coneIdx}},
we get number of subscripts of variables by \sphinxcode{\sphinxupquote{pReqSize}} specified by \sphinxcode{\sphinxupquote{nCone}}
and \sphinxcode{\sphinxupquote{list}}. Secondly, allocate sufficient memory for CRS\sphinxhyphen{}format matrix
and call this function again to extract weight matrix. If the memory
of weight matrix passed to function is not sufficient, then return
the first \sphinxcode{\sphinxupquote{nElemSize}} subscripts of variables, and the minimal required length
of non\sphinxhyphen{}zero elements by \sphinxcode{\sphinxupquote{pReqSize}}. If \sphinxcode{\sphinxupquote{list}} is \sphinxcode{\sphinxupquote{NULL}}, then
the first \sphinxcode{\sphinxupquote{nCone}} rows will be returned.
\end{quote}

\sphinxAtStartPar
\sphinxstylestrong{Arguments}
\begin{quote}

\sphinxAtStartPar
\sphinxcode{\sphinxupquote{prob}}
\begin{quote}

\sphinxAtStartPar
The COPT problem.
\end{quote}

\sphinxAtStartPar
\sphinxcode{\sphinxupquote{nCone}}
\begin{quote}

\sphinxAtStartPar
Number of SOC constraints.
\end{quote}

\sphinxAtStartPar
\sphinxcode{\sphinxupquote{list}}
\begin{quote}

\sphinxAtStartPar
Index of SOC constraints. Can be \sphinxcode{\sphinxupquote{NULL}}.
\end{quote}

\sphinxAtStartPar
\sphinxcode{\sphinxupquote{coneBeg, coneCnt, coneIdx}}
\begin{quote}

\sphinxAtStartPar
Defines the matrix of SOC constraints in compressed row
storage (CRS) format. Please see \sphinxstylestrong{other information} of
\sphinxcode{\sphinxupquote{COPT\_LoadProb}} for an example of the CRS format.
\end{quote}

\sphinxAtStartPar
\sphinxcode{\sphinxupquote{nElemSize}}
\begin{quote}

\sphinxAtStartPar
Length of array for non\sphinxhyphen{}zero weights.
\end{quote}

\sphinxAtStartPar
\sphinxcode{\sphinxupquote{pReqSize}}
\begin{quote}

\sphinxAtStartPar
Pointer to minimal length of array for non\sphinxhyphen{}zero weights.
Can be \sphinxcode{\sphinxupquote{NULL}}.
\end{quote}
\end{quote}
\end{quote}

\subsubsection{COPT\_GetExpCones}
\label{\detokenize{capiref:copt-getexpcones}}\begin{quote}

\sphinxAtStartPar
\sphinxstylestrong{Synopsis}
\begin{quote}

\sphinxAtStartPar
\sphinxcode{\sphinxupquote{int COPT\_GetExpCones(}}
\begin{quote}

\sphinxAtStartPar
\sphinxcode{\sphinxupquote{copt\_prob *prob,}}

\sphinxAtStartPar
\sphinxcode{\sphinxupquote{int nCone,}}

\sphinxAtStartPar
\sphinxcode{\sphinxupquote{const int *list,}}

\sphinxAtStartPar
\sphinxcode{\sphinxupquote{int *coneType,}}

\sphinxAtStartPar
\sphinxcode{\sphinxupquote{int *coneIdx,}}

\sphinxAtStartPar
\sphinxcode{\sphinxupquote{int nElemSize,}}

\sphinxAtStartPar
\sphinxcode{\sphinxupquote{int *pReqSize)}}
\end{quote}
\end{quote}

\sphinxAtStartPar
\sphinxstylestrong{Description}
\begin{quote}

\sphinxAtStartPar
Get the array of exponential cone constraints.

\sphinxAtStartPar
In general, users need to call this function twice to accomplish the task.
Firstly, by passing \sphinxcode{\sphinxupquote{NULL}} to arguments \sphinxcode{\sphinxupquote{coneIdx}},
we get number of subscripts of variables by \sphinxcode{\sphinxupquote{nElemSizes}} specified by \sphinxcode{\sphinxupquote{nCone}}
and \sphinxcode{\sphinxupquote{list}}. Secondly, allocate sufficient memory for array
and call this function again to extract weight array. If the memory
of weight array passed to function is not sufficient, then return
the first \sphinxcode{\sphinxupquote{nElemSize}} subscripts, and the minimal required length
of subscripts by \sphinxcode{\sphinxupquote{pReqSize}}. If \sphinxcode{\sphinxupquote{list}} is \sphinxcode{\sphinxupquote{NULL}}, then
the first \sphinxcode{\sphinxupquote{nCone}} rows will be returned.
\end{quote}

\sphinxAtStartPar
\sphinxstylestrong{Arguments}
\begin{quote}

\sphinxAtStartPar
\sphinxcode{\sphinxupquote{prob}}
\begin{quote}

\sphinxAtStartPar
The COPT problem.
\end{quote}

\sphinxAtStartPar
\sphinxcode{\sphinxupquote{nCone}}
\begin{quote}

\sphinxAtStartPar
Number of exponential cone constraints.
\end{quote}

\sphinxAtStartPar
\sphinxcode{\sphinxupquote{list}}
\begin{quote}

\sphinxAtStartPar
Index of exponential cone constraints. Can be \sphinxcode{\sphinxupquote{NULL}}.
\end{quote}

\sphinxAtStartPar
\sphinxcode{\sphinxupquote{coneType}}
\begin{quote}

\sphinxAtStartPar
Type of exponential cone constraints. Please refer to {\hyperref[\detokenize{constant:chapconst-expconetype}]{\sphinxcrossref{\DUrole{std,std-ref}{Exponential Cone type}}}} for possible values.
\end{quote}

\sphinxAtStartPar
\sphinxcode{\sphinxupquote{coneIdx}}
\begin{quote}

\sphinxAtStartPar
Array of subscripts of variables constituting the exponential cone constraints.
\end{quote}

\sphinxAtStartPar
\sphinxcode{\sphinxupquote{nElemSize}}
\begin{quote}

\sphinxAtStartPar
Length of array for subscripts of variables.
\end{quote}

\sphinxAtStartPar
\sphinxcode{\sphinxupquote{pReqSize}}
\begin{quote}

\sphinxAtStartPar
Pointer to minimal length of array for subscripts of variables.
Can be \sphinxcode{\sphinxupquote{NULL}}.
\end{quote}
\end{quote}
\end{quote}

\subsubsection{COPT\_GetAffineCone}
\label{\detokenize{capiref:copt-getaffinecone}}\begin{quote}

\sphinxAtStartPar
\sphinxstylestrong{Synopsis}
\begin{quote}

\sphinxAtStartPar
\sphinxcode{\sphinxupquote{int COPT\_GetAffineCone(}}
\begin{quote}

\sphinxAtStartPar
\sphinxcode{\sphinxupquote{copt\_prob *prob,}}

\sphinxAtStartPar
\sphinxcode{\sphinxupquote{int affConeIdx,}}

\sphinxAtStartPar
\sphinxcode{\sphinxupquote{int *coneType,}}

\sphinxAtStartPar
\sphinxcode{\sphinxupquote{int *nConeDim,}}

\sphinxAtStartPar
\sphinxcode{\sphinxupquote{int *nAlphaDim,}}

\sphinxAtStartPar
\sphinxcode{\sphinxupquote{double *alphaElem,}}

\sphinxAtStartPar
\sphinxcode{\sphinxupquote{int *psdBeg,}}

\sphinxAtStartPar
\sphinxcode{\sphinxupquote{int *psdCnt,}}

\sphinxAtStartPar
\sphinxcode{\sphinxupquote{int *psdColIdx,}}

\sphinxAtStartPar
\sphinxcode{\sphinxupquote{int *psdMatIdx,}}

\sphinxAtStartPar
\sphinxcode{\sphinxupquote{int nPsdElemSize,}}

\sphinxAtStartPar
\sphinxcode{\sphinxupquote{int *pPsdReqSize,}}

\sphinxAtStartPar
\sphinxcode{\sphinxupquote{int *rowMatBeg,}}

\sphinxAtStartPar
\sphinxcode{\sphinxupquote{int *rowMatCnt,}}

\sphinxAtStartPar
\sphinxcode{\sphinxupquote{int *rowMatIdx,}}

\sphinxAtStartPar
\sphinxcode{\sphinxupquote{double *rowMatElem,}}

\sphinxAtStartPar
\sphinxcode{\sphinxupquote{double *rowConst,}}

\sphinxAtStartPar
\sphinxcode{\sphinxupquote{int nElemSize,}}

\sphinxAtStartPar
\sphinxcode{\sphinxupquote{int *pReqSize)}}
\end{quote}
\end{quote}

\sphinxAtStartPar
\sphinxstylestrong{Description}
\begin{quote}

\sphinxAtStartPar
Retrieve the array of affine cone constraints from the model.

\sphinxAtStartPar
Generally, the user needs to call this function twice to complete the extraction of the affine cone constraint.

\sphinxAtStartPar
First, pass \sphinxcode{\sphinxupquote{NULL}} to \sphinxcode{\sphinxupquote{affConeIdx}} to obtain the number of indices for the affine cone constraint array specified by \sphinxcode{\sphinxupquote{nConeDim}} and \sphinxcode{\sphinxupquote{list}} through \sphinxcode{\sphinxupquote{pReqSize}}. Then, allocate appropriate space for the array parameters, and call the function again to retrieve the specified affine cone constraint array.

\sphinxAtStartPar
If the length of the affine cone constraint array passed in is insufficient, the function will return the first \sphinxcode{\sphinxupquote{nElemSize}} variable index arrays and return the minimum required array length through \sphinxcode{\sphinxupquote{pReqSize}}.
If \sphinxcode{\sphinxupquote{list}} is \sphinxcode{\sphinxupquote{NULL}}, the function returns the affine cone constraint array corresponding to the first \sphinxcode{\sphinxupquote{nConeDim}} rows.
\end{quote}

\sphinxAtStartPar
\sphinxstylestrong{Arguments}
\begin{quote}

\sphinxAtStartPar
\sphinxcode{\sphinxupquote{prob}}
\begin{quote}

\sphinxAtStartPar
The COPT problem.
\end{quote}

\sphinxAtStartPar
\sphinxcode{\sphinxupquote{affConeIdx}}
\begin{quote}

\sphinxAtStartPar
The index of the affine cone.
\end{quote}

\sphinxAtStartPar
\sphinxcode{\sphinxupquote{coneType}}
\begin{quote}

\sphinxAtStartPar
The type of the affine cone. Please refer to {\hyperref[\detokenize{constant:chapconst-conetype}]{\sphinxcrossref{\DUrole{std,std-ref}{SOC constraint types}}}} and {\hyperref[\detokenize{constant:chapconst-expconetype}]{\sphinxcrossref{\DUrole{std,std-ref}{Exponential Cone type}}}}
for possible values.
\end{quote}

\sphinxAtStartPar
\sphinxcode{\sphinxupquote{nConeDim}}
\begin{quote}

\sphinxAtStartPar
The dimension of the affine cone.
\end{quote}

\sphinxAtStartPar
\sphinxcode{\sphinxupquote{nAlphaDim}}
\begin{quote}

\sphinxAtStartPar
Reserved parameter, currently not in use.
\end{quote}

\sphinxAtStartPar
\sphinxcode{\sphinxupquote{alphaElem}}
\begin{quote}

\sphinxAtStartPar
Reserved parameter, currently not in use.
\end{quote}

\sphinxAtStartPar
\sphinxcode{\sphinxupquote{psdBeg, psdCnt, psdColIdx, psdMatIdx}}
\begin{quote}

\sphinxAtStartPar
Represents the PSD terms in the affine cone.

\sphinxAtStartPar
\sphinxcode{\sphinxupquote{psdBeg}} indicates the starting position of the PSD terms in each affine cone term,
\sphinxcode{\sphinxupquote{psdCnt}} specifies the number of PSD terms, \sphinxcode{\sphinxupquote{psdColIdx}} indicates the index of the
PSD variable, and \sphinxcode{\sphinxupquote{psdMatIdx}} refers to the index of the symmetric matrix.

\sphinxAtStartPar
First, the number of PSD terms is returned via \sphinxcode{\sphinxupquote{pPsdReqSize}}, and \sphinxcode{\sphinxupquote{nPsdElemSize}} represents the number of PSD terms to be retrieved.
\end{quote}

\sphinxAtStartPar
\sphinxcode{\sphinxupquote{rowMatBeg, rowMatCnt, rowMatIdx, rowMatElem}}
\begin{quote}

\sphinxAtStartPar
Represents the linear terms in the affine cone.

\sphinxAtStartPar
The coefficient matrix is provided in CRS\sphinxhyphen{}format. For detailed examples, please refer to \sphinxstylestrong{Additional Information} in \sphinxcode{\sphinxupquote{COPT\_LoadProb}}.
\end{quote}

\sphinxAtStartPar
\sphinxcode{\sphinxupquote{rowConst}}
\begin{quote}

\sphinxAtStartPar
The constant term in the affine cone.
\end{quote}

\sphinxAtStartPar
\sphinxcode{\sphinxupquote{nElemSize}}
\begin{quote}

\sphinxAtStartPar
Length of array for subscripts of variables.
\end{quote}

\sphinxAtStartPar
\sphinxcode{\sphinxupquote{pReqSize}}
\begin{quote}

\sphinxAtStartPar
Pointer to minimal length of array for subscripts of variables of the affine cone.
Can be \sphinxcode{\sphinxupquote{NULL}}.
\end{quote}
\end{quote}
\end{quote}

\subsubsection{COPT\_GetQConstr}
\label{\detokenize{capiref:copt-getqconstr}}\begin{quote}

\sphinxAtStartPar
\sphinxstylestrong{Synopsis}
\begin{quote}

\sphinxAtStartPar
\sphinxcode{\sphinxupquote{int COPT\_GetQConstr(}}
\begin{quote}

\sphinxAtStartPar
\sphinxcode{\sphinxupquote{copt\_prob *prob,}}

\sphinxAtStartPar
\sphinxcode{\sphinxupquote{int qConstrIdx,}}

\sphinxAtStartPar
\sphinxcode{\sphinxupquote{int *qMatRow,}}

\sphinxAtStartPar
\sphinxcode{\sphinxupquote{int *qMatCol,}}

\sphinxAtStartPar
\sphinxcode{\sphinxupquote{double *qMatElem,}}

\sphinxAtStartPar
\sphinxcode{\sphinxupquote{int nQElemSize,}}

\sphinxAtStartPar
\sphinxcode{\sphinxupquote{int *pQReqSize,}}

\sphinxAtStartPar
\sphinxcode{\sphinxupquote{int *rowMatIdx,}}

\sphinxAtStartPar
\sphinxcode{\sphinxupquote{double *rowMatElem,}}

\sphinxAtStartPar
\sphinxcode{\sphinxupquote{char *cRowSense,}}

\sphinxAtStartPar
\sphinxcode{\sphinxupquote{double *dRowBound,}}

\sphinxAtStartPar
\sphinxcode{\sphinxupquote{int nElemSize,}}

\sphinxAtStartPar
\sphinxcode{\sphinxupquote{int *pReqSize)}}
\end{quote}
\end{quote}

\sphinxAtStartPar
\sphinxstylestrong{Description}
\begin{quote}

\sphinxAtStartPar
Get quadratic constraint.

\sphinxAtStartPar
In general, users need to call this function twice to accomplish the task.
Firstly, by passing \sphinxcode{\sphinxupquote{NULL}} to arguments \sphinxcode{\sphinxupquote{qMatRow}}, \sphinxcode{\sphinxupquote{qMatCol}},
\sphinxcode{\sphinxupquote{qMatElem}}, \sphinxcode{\sphinxupquote{rowMatIdx}} and \sphinxcode{\sphinxupquote{rowMatElem}}, we get number of non\sphinxhyphen{}zero
quadratic terms by \sphinxcode{\sphinxupquote{pQReqSize}} and number of non\sphinxhyphen{}zero linear terms by
\sphinxcode{\sphinxupquote{pReqSize}} specified by \sphinxcode{\sphinxupquote{qConstrIdx}}. Secondly, allocate sufficient
memory for the quadratic terms and the linear terms, and call this function
again to extract the quadratic constraint. If the memory of the array of
the quadratic terms passed to function is not sufficient, then return
the first \sphinxcode{\sphinxupquote{nQElemSize}} quadratic terms, and the minimal required length
of quadratic terms by \sphinxcode{\sphinxupquote{pQReqSize}}. If  the memory of the array of the
linear terms passed to function is not sufficient, then return the first
\sphinxcode{\sphinxupquote{nElemSize}} linear terms, and the minimal required length of linear terms
by \sphinxcode{\sphinxupquote{pReqSize}}.
\end{quote}

\sphinxAtStartPar
\sphinxstylestrong{Arguments}
\begin{quote}

\sphinxAtStartPar
\sphinxcode{\sphinxupquote{prob}}
\begin{quote}

\sphinxAtStartPar
The COPT problem.
\end{quote}

\sphinxAtStartPar
\sphinxcode{\sphinxupquote{qConstrIdx}}
\begin{quote}

\sphinxAtStartPar
Index of the quadratic constraint.
\end{quote}

\sphinxAtStartPar
\sphinxcode{\sphinxupquote{qMatRow}}
\begin{quote}

\sphinxAtStartPar
Row index of non\sphinxhyphen{}zero quadratic terms of the quadratic constraint (row).
\end{quote}

\sphinxAtStartPar
\sphinxcode{\sphinxupquote{qMatCol}}
\begin{quote}

\sphinxAtStartPar
Column index of non\sphinxhyphen{}zero quadratic terms of the quadratic constraint (row).
\end{quote}

\sphinxAtStartPar
\sphinxcode{\sphinxupquote{qMatElem}}
\begin{quote}

\sphinxAtStartPar
Values of non\sphinxhyphen{}zero quadratic terms of the quadratic constraint (row).
\end{quote}

\sphinxAtStartPar
\sphinxcode{\sphinxupquote{nQElemSize}}
\begin{quote}

\sphinxAtStartPar
Length of array for non\sphinxhyphen{}zero quadratic terms of the quadratic constraint (row).
\end{quote}

\sphinxAtStartPar
\sphinxcode{\sphinxupquote{pQReqSize}}
\begin{quote}

\sphinxAtStartPar
Pointer to minimal length of array for non\sphinxhyphen{}zero quadratic terms of the quadratic constraint (row).
Can be \sphinxcode{\sphinxupquote{NULL}}.
\end{quote}

\sphinxAtStartPar
\sphinxcode{\sphinxupquote{rowMatIdx}}
\begin{quote}

\sphinxAtStartPar
Column index of non\sphinxhyphen{}zero linear terms of the quadratic constraint (row).
\end{quote}

\sphinxAtStartPar
\sphinxcode{\sphinxupquote{rowMatElem}}
\begin{quote}

\sphinxAtStartPar
Values of non\sphinxhyphen{}zero linear terms of the quadratic constraint (row).
\end{quote}

\sphinxAtStartPar
\sphinxcode{\sphinxupquote{cRowSense}}
\begin{quote}

\sphinxAtStartPar
The sense of the quadratic constraint (row).
\end{quote}

\sphinxAtStartPar
\sphinxcode{\sphinxupquote{dRowBound}}
\begin{quote}

\sphinxAtStartPar
Right hand side of the quadratic constraint (row).
\end{quote}

\sphinxAtStartPar
\sphinxcode{\sphinxupquote{nElemSize}}
\begin{quote}

\sphinxAtStartPar
Length of array for non\sphinxhyphen{}zero linear terms of the quadratic constraint (row).
\end{quote}

\sphinxAtStartPar
\sphinxcode{\sphinxupquote{pReqSize}}
\begin{quote}

\sphinxAtStartPar
Pointer to minimal length of array for non\sphinxhyphen{}zero linear terms of the quadratic constraint (row).
Can be \sphinxcode{\sphinxupquote{NULL}}.
\end{quote}
\end{quote}
\end{quote}

\subsubsection{COPT\_GetPSDConstr}
\label{\detokenize{capiref:copt-getpsdconstr}}\begin{quote}

\sphinxAtStartPar
\sphinxstylestrong{Synopsis}
\begin{quote}

\sphinxAtStartPar
\sphinxcode{\sphinxupquote{int COPT\_GetPSDConstr(}}
\begin{quote}

\sphinxAtStartPar
\sphinxcode{\sphinxupquote{copt\_prob *prob,}}

\sphinxAtStartPar
\sphinxcode{\sphinxupquote{int psdConstrIdx,}}

\sphinxAtStartPar
\sphinxcode{\sphinxupquote{int *psdColIdx,}}

\sphinxAtStartPar
\sphinxcode{\sphinxupquote{int *symMatIdx,}}

\sphinxAtStartPar
\sphinxcode{\sphinxupquote{int nColSize,}}

\sphinxAtStartPar
\sphinxcode{\sphinxupquote{int *pColReqSize,}}

\sphinxAtStartPar
\sphinxcode{\sphinxupquote{int *rowMatIdx,}}

\sphinxAtStartPar
\sphinxcode{\sphinxupquote{double *rowMatElem,}}

\sphinxAtStartPar
\sphinxcode{\sphinxupquote{double *dRowLower,}}

\sphinxAtStartPar
\sphinxcode{\sphinxupquote{double *dRowUpper,}}

\sphinxAtStartPar
\sphinxcode{\sphinxupquote{int nElemSize,}}

\sphinxAtStartPar
\sphinxcode{\sphinxupquote{int *pReqSize)}}
\end{quote}
\end{quote}

\sphinxAtStartPar
\sphinxstylestrong{Description}
\begin{quote}

\sphinxAtStartPar
Get PSD constraint.

\sphinxAtStartPar
In general, users need to call this function twice to accomplish the task.
Fisrtly, by passing \sphinxcode{\sphinxupquote{NULL}} to arguments \sphinxcode{\sphinxupquote{psdColIdx}} and \sphinxcode{\sphinxupquote{symMatIdx}},
we get number of PSD terms by \sphinxcode{\sphinxupquote{pColReqSize}} specified by \sphinxcode{\sphinxupquote{psdConstrIdx}},
by passing \sphinxcode{\sphinxupquote{NULL}} to arguments \sphinxcode{\sphinxupquote{rowMatIdx}} and \sphinxcode{\sphinxupquote{rowMatElem}}, we get
number of linear terms by \sphinxcode{\sphinxupquote{pReqSize}} specified by \sphinxcode{\sphinxupquote{qConstrIdx}}.
Secondly, allocate sufficient memory for the PSD terms and the linear terms,
and call this function again to extract the PSD constraint. If the memory
of the array of the PSD terms passed to function is not sufficient, then return
the first \sphinxcode{\sphinxupquote{nColSize}} PSD terms, and the minimal required length of PSD terms
by \sphinxcode{\sphinxupquote{pColReqSize}}. If the memory of the array of the linear terms passed to
function is not sufficient, then return the first \sphinxcode{\sphinxupquote{nElemSize}} linear terms,
and the minimal required length of linear terms by \sphinxcode{\sphinxupquote{pReqSize}}.
\end{quote}

\sphinxAtStartPar
\sphinxstylestrong{Arguments}
\begin{quote}

\sphinxAtStartPar
\sphinxcode{\sphinxupquote{prob}}
\begin{quote}

\sphinxAtStartPar
The COPT problem.
\end{quote}

\sphinxAtStartPar
\sphinxcode{\sphinxupquote{psdConstrIdx}}
\begin{quote}

\sphinxAtStartPar
PSD constraint index.
\end{quote}

\sphinxAtStartPar
\sphinxcode{\sphinxupquote{psdColIdx}}
\begin{quote}

\sphinxAtStartPar
PSD variable index.
\end{quote}

\sphinxAtStartPar
\sphinxcode{\sphinxupquote{symMatIdx}}
\begin{quote}

\sphinxAtStartPar
Symmetric matrix index.
\end{quote}

\sphinxAtStartPar
\sphinxcode{\sphinxupquote{nColSize}}
\begin{quote}

\sphinxAtStartPar
Length of array for PSD terms of the PSD constraint.
\end{quote}

\sphinxAtStartPar
\sphinxcode{\sphinxupquote{pColReqSize}}
\begin{quote}

\sphinxAtStartPar
Pointer to minimal length of array for PSD terms of the PSD constraint.
Can be \sphinxcode{\sphinxupquote{NULL}}.
\end{quote}

\sphinxAtStartPar
\sphinxcode{\sphinxupquote{rowMatIdx}}
\begin{quote}

\sphinxAtStartPar
Column index of non\sphinxhyphen{}zero linear terms of the PSD constraint.
\end{quote}

\sphinxAtStartPar
\sphinxcode{\sphinxupquote{rowMatElem}}
\begin{quote}

\sphinxAtStartPar
Values of non\sphinxhyphen{}zero linear terms of the PSD constraint.
\end{quote}

\sphinxAtStartPar
\sphinxcode{\sphinxupquote{dRowLower}}
\begin{quote}

\sphinxAtStartPar
Pointer to lower bound of the PSD constraint.
\end{quote}

\sphinxAtStartPar
\sphinxcode{\sphinxupquote{dRowUpper}}
\begin{quote}

\sphinxAtStartPar
Pointer to upper bound of the PSD constraint.
\end{quote}

\sphinxAtStartPar
\sphinxcode{\sphinxupquote{nElemSize}}
\begin{quote}

\sphinxAtStartPar
Length of array for non\sphinxhyphen{}zero linear terms of the PSD constraint.
\end{quote}

\sphinxAtStartPar
\sphinxcode{\sphinxupquote{pReqSize}}
\begin{quote}

\sphinxAtStartPar
Pointer to minimal length of array for non\sphinxhyphen{}zero linear terms of the
PSD constraint (row). Can be \sphinxcode{\sphinxupquote{NULL}}.
\end{quote}
\end{quote}
\end{quote}

\subsubsection{COPT\_GetLMIConstr}
\label{\detokenize{capiref:copt-getlmiconstr}}\begin{quote}

\sphinxAtStartPar
\sphinxstylestrong{Synopsis}
\begin{quote}

\sphinxAtStartPar
\sphinxcode{\sphinxupquote{int COPT\_GetLMIConstr(}}
\begin{quote}

\sphinxAtStartPar
\sphinxcode{\sphinxupquote{copt\_prob *prob,}}

\sphinxAtStartPar
\sphinxcode{\sphinxupquote{int lmiConstrIdx,}}

\sphinxAtStartPar
\sphinxcode{\sphinxupquote{int *nDim,}}

\sphinxAtStartPar
\sphinxcode{\sphinxupquote{int *nLMILen,}}

\sphinxAtStartPar
\sphinxcode{\sphinxupquote{int *colIdx,}}

\sphinxAtStartPar
\sphinxcode{\sphinxupquote{int *symMatIdx,}}

\sphinxAtStartPar
\sphinxcode{\sphinxupquote{int *constMatIdx,}}

\sphinxAtStartPar
\sphinxcode{\sphinxupquote{int nElemSize,}}

\sphinxAtStartPar
\sphinxcode{\sphinxupquote{int *pReqSize)}}
\end{quote}
\end{quote}

\sphinxAtStartPar
\sphinxstylestrong{Description}
\begin{quote}

\sphinxAtStartPar
Gets the LMI constraint with the specified index in the model.
\end{quote}

\sphinxAtStartPar
\sphinxstylestrong{Arguments}
\begin{quote}

\sphinxAtStartPar
\sphinxcode{\sphinxupquote{prob}}
\begin{quote}

\sphinxAtStartPar
The COPT problem.
\end{quote}

\sphinxAtStartPar
\sphinxcode{\sphinxupquote{lmiConstrIdx}}
\begin{quote}

\sphinxAtStartPar
LMI constraint index.
\end{quote}

\sphinxAtStartPar
\sphinxcode{\sphinxupquote{nDim}}
\begin{quote}

\sphinxAtStartPar
Pointer to the dimension of symmetric matrix in the LMI constraint.
\end{quote}

\sphinxAtStartPar
\sphinxcode{\sphinxupquote{nLMILen}}
\begin{quote}

\sphinxAtStartPar
Pointer to the flattened length of the LMI constraint.
\end{quote}

\sphinxAtStartPar
\sphinxcode{\sphinxupquote{colIdx}}
\begin{quote}

\sphinxAtStartPar
Index of scalar variable in the LMI constraint.
\end{quote}

\sphinxAtStartPar
\sphinxcode{\sphinxupquote{symMatIdx}}
\begin{quote}

\sphinxAtStartPar
Index of symmetric coefficient matrix in the LMI constraint.
\end{quote}

\sphinxAtStartPar
\sphinxcode{\sphinxupquote{constMatIdx}}
\begin{quote}

\sphinxAtStartPar
Pointer to the index of symmetric constant\sphinxhyphen{}term matrix in the LMI constraint.
\end{quote}

\sphinxAtStartPar
\sphinxcode{\sphinxupquote{nElemSize}}
\begin{quote}

\sphinxAtStartPar
Length of array for non\sphinxhyphen{}zero linear terms of the LMI constraint.
\end{quote}

\sphinxAtStartPar
\sphinxcode{\sphinxupquote{pReqSize}}
\begin{quote}

\sphinxAtStartPar
Pointer to minimal length of array for non\sphinxhyphen{}zero linear terms of the
LMI constraint (row). Can be \sphinxcode{\sphinxupquote{NULL}}.
\end{quote}
\end{quote}
\end{quote}

\subsubsection{COPT\_GetIndicator}
\label{\detokenize{capiref:copt-getindicator}}\begin{quote}

\sphinxAtStartPar
\sphinxstylestrong{Synopsis}
\begin{quote}

\sphinxAtStartPar
\sphinxcode{\sphinxupquote{int COPT\_GetIndicator(}}
\begin{quote}

\sphinxAtStartPar
\sphinxcode{\sphinxupquote{copt\_prob *prob,}}

\sphinxAtStartPar
\sphinxcode{\sphinxupquote{int rowIdx,}}

\sphinxAtStartPar
\sphinxcode{\sphinxupquote{int *binColIdx,}}

\sphinxAtStartPar
\sphinxcode{\sphinxupquote{int *binColVal,}}

\sphinxAtStartPar
\sphinxcode{\sphinxupquote{int *nRowMatCnt,}}

\sphinxAtStartPar
\sphinxcode{\sphinxupquote{int *rowMatIdx,}}

\sphinxAtStartPar
\sphinxcode{\sphinxupquote{double *rowMatElem,}}

\sphinxAtStartPar
\sphinxcode{\sphinxupquote{char *cRowSense,}}

\sphinxAtStartPar
\sphinxcode{\sphinxupquote{double *dRowBound,}}

\sphinxAtStartPar
\sphinxcode{\sphinxupquote{int nElemSize,}}

\sphinxAtStartPar
\sphinxcode{\sphinxupquote{int *pReqSize)}}
\end{quote}
\end{quote}

\sphinxAtStartPar
\sphinxstylestrong{Description}
\begin{quote}

\sphinxAtStartPar
Get the data of an indicator constraint.

\sphinxAtStartPar
In general, users need to call this function twice to accomplish the task.
Firstly, by passing \sphinxcode{\sphinxupquote{NULL}} to arguments \sphinxcode{\sphinxupquote{nRowMatCnt}}, \sphinxcode{\sphinxupquote{rowMatIdx}} and \sphinxcode{\sphinxupquote{rowMatElem}},
we get number of non\sphinxhyphen{}zeros elements by \sphinxcode{\sphinxupquote{pReqSize}} specified by \sphinxcode{\sphinxupquote{rowIdx}}.
Secondly, allocate sufficient memory for sparse row vector and call
this function again to extract data.

\sphinxAtStartPar
If the memory of sparse row vector
passed to function is not sufficient, then return the first \sphinxcode{\sphinxupquote{nElemSize}}
non\sphinxhyphen{}zero elements, and the minimal required length of non\sphinxhyphen{}zero elements
by \sphinxcode{\sphinxupquote{pReqSize}}.
\end{quote}

\sphinxAtStartPar
\sphinxstylestrong{Arguments}
\begin{quote}

\sphinxAtStartPar
\sphinxcode{\sphinxupquote{prob}}
\begin{quote}

\sphinxAtStartPar
The COPT problem.
\end{quote}

\sphinxAtStartPar
\sphinxcode{\sphinxupquote{rowIdx}}
\begin{quote}

\sphinxAtStartPar
Index of the indicator constraint.
\end{quote}

\sphinxAtStartPar
\sphinxcode{\sphinxupquote{binColIdx}}
\begin{quote}

\sphinxAtStartPar
Index of the indicator variable (column).
\end{quote}

\sphinxAtStartPar
\sphinxcode{\sphinxupquote{binColVal}}
\begin{quote}

\sphinxAtStartPar
Value of the indicator variable (column).
\end{quote}

\sphinxAtStartPar
\sphinxcode{\sphinxupquote{nRowMatCnt}}
\begin{quote}

\sphinxAtStartPar
Number of non\sphinxhyphen{}zeros elements in the linear constraint (row).
\end{quote}

\sphinxAtStartPar
\sphinxcode{\sphinxupquote{rowMatIdx}}
\begin{quote}

\sphinxAtStartPar
Column index of non\sphinxhyphen{}zeros elements in the linear constraint (row).
\end{quote}

\sphinxAtStartPar
\sphinxcode{\sphinxupquote{rowMatElem}}
\begin{quote}

\sphinxAtStartPar
Values of non\sphinxhyphen{}zero elements in the linear constraint (row).
\end{quote}

\sphinxAtStartPar
\sphinxcode{\sphinxupquote{cRowSense}}
\begin{quote}

\sphinxAtStartPar
The sense of the linear constraint (row).
\end{quote}

\sphinxAtStartPar
\sphinxcode{\sphinxupquote{dRowBound}}
\begin{quote}

\sphinxAtStartPar
Right hand side of the linear constraint (row).
\end{quote}

\sphinxAtStartPar
\sphinxcode{\sphinxupquote{nElemSize}}
\begin{quote}

\sphinxAtStartPar
Length of array for non\sphinxhyphen{}zero coefficients.
\end{quote}

\sphinxAtStartPar
\sphinxcode{\sphinxupquote{pReqSize}}
\begin{quote}

\sphinxAtStartPar
Pointer to minimal length of array for non\sphinxhyphen{}zero coefficients.
Can be \sphinxcode{\sphinxupquote{NULL}}.
\end{quote}
\end{quote}
\end{quote}

\subsubsection{COPT\_GetIndicators}
\label{\detokenize{capiref:copt-getindicators}}\begin{quote}

\sphinxAtStartPar
\sphinxstylestrong{Synopsis}
\begin{quote}

\sphinxAtStartPar
\sphinxcode{\sphinxupquote{int COPT\_GetIndicators(}}
\begin{quote}

\sphinxAtStartPar
\sphinxcode{\sphinxupquote{copt\_prob *prob,}}

\sphinxAtStartPar
\sphinxcode{\sphinxupquote{int nInd,}}

\sphinxAtStartPar
\sphinxcode{\sphinxupquote{int *list,}}

\sphinxAtStartPar
\sphinxcode{\sphinxupquote{int *indType,}}

\sphinxAtStartPar
\sphinxcode{\sphinxupquote{int *binColIdx,}}

\sphinxAtStartPar
\sphinxcode{\sphinxupquote{int *binColVal,}}

\sphinxAtStartPar
\sphinxcode{\sphinxupquote{int *rowMatBeg,}}

\sphinxAtStartPar
\sphinxcode{\sphinxupquote{int *rowMatCnt,}}

\sphinxAtStartPar
\sphinxcode{\sphinxupquote{int *rowMatIdx,}}

\sphinxAtStartPar
\sphinxcode{\sphinxupquote{double *rowMatElem,}}

\sphinxAtStartPar
\sphinxcode{\sphinxupquote{char *cRowSense,}}

\sphinxAtStartPar
\sphinxcode{\sphinxupquote{double *dRowBound,}}

\sphinxAtStartPar
\sphinxcode{\sphinxupquote{int nElemSize,}}

\sphinxAtStartPar
\sphinxcode{\sphinxupquote{int *pReqSize)}}
\end{quote}
\end{quote}

\sphinxAtStartPar
\sphinxstylestrong{Description}
\begin{quote}

\sphinxAtStartPar
Get the data of a set of indicator constraints.

\sphinxAtStartPar
In general, users need to call this function twice to accomplish the task.
Firstly, by passing \sphinxcode{\sphinxupquote{NULL}} to arguments \sphinxcode{\sphinxupquote{rowMatBeg}}, \sphinxcode{\sphinxupquote{rowMatCnt}},
\sphinxcode{\sphinxupquote{rowMatIdx}} , we get number of non\sphinxhyphen{}zeros elements by
\sphinxcode{\sphinxupquote{pReqSize}} specified by \sphinxcode{\sphinxupquote{nInd}} and \sphinxcode{\sphinxupquote{list}}. Secondly, allocate
sufficient memory for CRS\sphinxhyphen{}format matrix and call this function again to
extract coefficient matrix.

\sphinxAtStartPar
If the memory of coefficient matrix passed to function
is not sufficient, then return the first \sphinxcode{\sphinxupquote{nElemSize}} non\sphinxhyphen{}zero elements,
and the minimal required length of non\sphinxhyphen{}zero elements by \sphinxcode{\sphinxupquote{pReqSize}}.
If \sphinxcode{\sphinxupquote{list}} is \sphinxcode{\sphinxupquote{NULL}}, then the first \sphinxcode{\sphinxupquote{nInd}} rows will be returned.
\end{quote}

\sphinxAtStartPar
\sphinxstylestrong{Arguments}
\begin{quote}

\sphinxAtStartPar
\sphinxcode{\sphinxupquote{prob}}
\begin{quote}

\sphinxAtStartPar
The COPT problem.
\end{quote}

\sphinxAtStartPar
\sphinxcode{\sphinxupquote{nInd}}
\begin{quote}

\sphinxAtStartPar
Number of indicator constraints (rows).
\end{quote}

\sphinxAtStartPar
\sphinxcode{\sphinxupquote{list}}
\begin{quote}

\sphinxAtStartPar
Index of indicator constraints. Can be \sphinxcode{\sphinxupquote{NULL}}.
\end{quote}

\sphinxAtStartPar
\sphinxcode{\sphinxupquote{indType}}
\begin{quote}

\sphinxAtStartPar
Type of indicator constraints. Please refer to {\hyperref[\detokenize{constant:chapconst-indicatortype}]{\sphinxcrossref{\DUrole{std,std-ref}{Indicator constraint types}}}} for possible values.
\end{quote}

\sphinxAtStartPar
\sphinxcode{\sphinxupquote{binColIdx}}
\begin{quote}

\sphinxAtStartPar
Index of the indicator variable (column).
\end{quote}

\sphinxAtStartPar
\sphinxcode{\sphinxupquote{binColVal}}
\begin{quote}

\sphinxAtStartPar
Value of the indicator variable (column).
\end{quote}

\sphinxAtStartPar
\sphinxcode{\sphinxupquote{rowMatBeg, rowMatCnt, rowMatIdx}} and \sphinxcode{\sphinxupquote{rowMatElem}}
\begin{quote}

\sphinxAtStartPar
Defines the coefficient matrix of indcator constrants in compressed row storage (CRS) format.
Please see \sphinxstylestrong{other information} of \sphinxcode{\sphinxupquote{COPT\_LoadProb}} for
an example of the CRS format.
\end{quote}

\sphinxAtStartPar
\sphinxcode{\sphinxupquote{cRowSense}}
\begin{quote}

\sphinxAtStartPar
The sense of the linear constraint (row).
\end{quote}

\sphinxAtStartPar
\sphinxcode{\sphinxupquote{dRowBound}}
\begin{quote}

\sphinxAtStartPar
Right hand side of the linear constraint (row).
\end{quote}

\sphinxAtStartPar
\sphinxcode{\sphinxupquote{nElemSize}}
\begin{quote}

\sphinxAtStartPar
Length of array for non\sphinxhyphen{}zero coefficients.
\end{quote}

\sphinxAtStartPar
\sphinxcode{\sphinxupquote{pReqSize}}
\begin{quote}

\sphinxAtStartPar
Pointer to minimal length of array for non\sphinxhyphen{}zero coefficients.
Can be \sphinxcode{\sphinxupquote{NULL}}.
\end{quote}
\end{quote}
\end{quote}

\subsubsection{COPT\_GetNLConstr}
\label{\detokenize{capiref:copt-getnlconstr}}\begin{quote}

\sphinxAtStartPar
\sphinxstylestrong{Synopsis}
\begin{quote}

\sphinxAtStartPar
\sphinxcode{\sphinxupquote{int COPT\_GetNLConstr(}}
\begin{quote}

\sphinxAtStartPar
\sphinxcode{\sphinxupquote{copt\_prob *prob,}}

\sphinxAtStartPar
\sphinxcode{\sphinxupquote{int nlConstrIdx,}}

\sphinxAtStartPar
\sphinxcode{\sphinxupquote{int *token,}}

\sphinxAtStartPar
\sphinxcode{\sphinxupquote{double *tokenElem,}}

\sphinxAtStartPar
\sphinxcode{\sphinxupquote{int nToken,}}

\sphinxAtStartPar
\sphinxcode{\sphinxupquote{int nTokenElem,}}

\sphinxAtStartPar
\sphinxcode{\sphinxupquote{int *pReqToken,}}

\sphinxAtStartPar
\sphinxcode{\sphinxupquote{int *pReqTokenElem,}}

\sphinxAtStartPar
\sphinxcode{\sphinxupquote{int *rowMatIdx,}}

\sphinxAtStartPar
\sphinxcode{\sphinxupquote{double *rowMatElem,}}

\sphinxAtStartPar
\sphinxcode{\sphinxupquote{double *dRowLower,}}

\sphinxAtStartPar
\sphinxcode{\sphinxupquote{double *dRowUpper,}}

\sphinxAtStartPar
\sphinxcode{\sphinxupquote{int nElemSize,}}

\sphinxAtStartPar
\sphinxcode{\sphinxupquote{int *pReqElemSize)}}
\end{quote}
\end{quote}

\sphinxAtStartPar
\sphinxstylestrong{Description}
\begin{quote}

\sphinxAtStartPar
Retrieve the nonlinear expression constraint at the specified index.
\end{quote}

\sphinxAtStartPar
\sphinxstylestrong{Arguments}
\begin{quote}

\sphinxAtStartPar
\sphinxcode{\sphinxupquote{prob}}
\begin{quote}

\sphinxAtStartPar
The COPT problem.
\end{quote}

\sphinxAtStartPar
\sphinxcode{\sphinxupquote{nlConstrIdx}}
\begin{quote}

\sphinxAtStartPar
Index of the constraint.
\end{quote}

\sphinxAtStartPar
\sphinxcode{\sphinxupquote{token}}
\begin{quote}

\sphinxAtStartPar
Array of tokens in the expression.
\end{quote}

\sphinxAtStartPar
\sphinxcode{\sphinxupquote{tokenElem}}
\begin{quote}

\sphinxAtStartPar
Array of constants in the expression.
\end{quote}

\sphinxAtStartPar
\sphinxcode{\sphinxupquote{nToken}}
\begin{quote}

\sphinxAtStartPar
Number of tokens in the expression.
\end{quote}

\sphinxAtStartPar
\sphinxcode{\sphinxupquote{nTokenElem}}
\begin{quote}

\sphinxAtStartPar
Number of constants in the expression.
\end{quote}

\sphinxAtStartPar
\sphinxcode{\sphinxupquote{pReqToken}}
\begin{quote}

\sphinxAtStartPar
Pointer to the number of tokens in the expression.
\end{quote}

\sphinxAtStartPar
\sphinxcode{\sphinxupquote{pReqTokenElem}}
\begin{quote}

\sphinxAtStartPar
Pointer to the number of constants in the expression.
\end{quote}

\sphinxAtStartPar
\sphinxcode{\sphinxupquote{rowMatIdx}}
\begin{quote}

\sphinxAtStartPar
Indices of linear terms in the constraint.
\end{quote}

\sphinxAtStartPar
\sphinxcode{\sphinxupquote{rowMatElem}}
\begin{quote}

\sphinxAtStartPar
Coefficients of linear terms in the constraint.
\end{quote}

\sphinxAtStartPar
\sphinxcode{\sphinxupquote{dRowLower}}
\begin{quote}

\sphinxAtStartPar
Pointer to the lower bound of the constraint.
\end{quote}

\sphinxAtStartPar
\sphinxcode{\sphinxupquote{dRowUpper}}
\begin{quote}

\sphinxAtStartPar
Pointer to the upper bound of the constraint.
\end{quote}

\sphinxAtStartPar
\sphinxcode{\sphinxupquote{nElemSize}}
\begin{quote}

\sphinxAtStartPar
Number of linear terms in the constraint.
\end{quote}

\sphinxAtStartPar
\sphinxcode{\sphinxupquote{pReqElemSize}}
\begin{quote}

\sphinxAtStartPar
Pointer to the number of linear terms.
\end{quote}
\end{quote}
\end{quote}

\subsubsection{COPT\_GetColIdx}
\label{\detokenize{capiref:copt-getcolidx}}\begin{quote}

\sphinxAtStartPar
\sphinxstylestrong{Synopsis}
\begin{quote}

\sphinxAtStartPar
\sphinxcode{\sphinxupquote{int COPT\_GetColIdx(copt\_prob *prob, const char *colName, int *p\_iCol)}}
\end{quote}

\sphinxAtStartPar
\sphinxstylestrong{Description}
\begin{quote}

\sphinxAtStartPar
Get index of column by name.
\end{quote}

\sphinxAtStartPar
\sphinxstylestrong{Arguments}
\begin{quote}

\sphinxAtStartPar
\sphinxcode{\sphinxupquote{prob}}
\begin{quote}

\sphinxAtStartPar
The COPT problem.
\end{quote}

\sphinxAtStartPar
\sphinxcode{\sphinxupquote{colName}}
\begin{quote}

\sphinxAtStartPar
Name of column.
\end{quote}

\sphinxAtStartPar
\sphinxcode{\sphinxupquote{p\_iCol}}
\begin{quote}

\sphinxAtStartPar
Pointer to index of column.
\end{quote}
\end{quote}
\end{quote}

\subsubsection{COPT\_GetPSDColIdx}
\label{\detokenize{capiref:copt-getpsdcolidx}}\begin{quote}

\sphinxAtStartPar
\sphinxstylestrong{Synopsis}
\begin{quote}

\sphinxAtStartPar
\sphinxcode{\sphinxupquote{int COPT\_GetPSDColIdx(copt\_prob *prob, const char *psdColName, int *p\_iPSDCol)}}
\end{quote}

\sphinxAtStartPar
\sphinxstylestrong{Description}
\begin{quote}

\sphinxAtStartPar
Get index of PSD variable by name.
\end{quote}

\sphinxAtStartPar
\sphinxstylestrong{Arguments}
\begin{quote}

\sphinxAtStartPar
\sphinxcode{\sphinxupquote{prob}}
\begin{quote}

\sphinxAtStartPar
The COPT problem.
\end{quote}

\sphinxAtStartPar
\sphinxcode{\sphinxupquote{psdColName}}
\begin{quote}

\sphinxAtStartPar
Name of PSD variable.
\end{quote}

\sphinxAtStartPar
\sphinxcode{\sphinxupquote{p\_iPSDCol}}
\begin{quote}

\sphinxAtStartPar
Pointer to index of PSD variable.
\end{quote}
\end{quote}
\end{quote}

\subsubsection{COPT\_GetRowIdx}
\label{\detokenize{capiref:copt-getrowidx}}\begin{quote}

\sphinxAtStartPar
\sphinxstylestrong{Synopsis}
\begin{quote}

\sphinxAtStartPar
\sphinxcode{\sphinxupquote{int COPT\_GetRowIdx(copt\_prob *prob, const char *rowName, int *p\_iRow)}}
\end{quote}

\sphinxAtStartPar
\sphinxstylestrong{Description}
\begin{quote}

\sphinxAtStartPar
Get index of row by name.
\end{quote}

\sphinxAtStartPar
\sphinxstylestrong{Arguments}
\begin{quote}

\sphinxAtStartPar
\sphinxcode{\sphinxupquote{prob}}
\begin{quote}

\sphinxAtStartPar
The COPT problem.
\end{quote}

\sphinxAtStartPar
\sphinxcode{\sphinxupquote{rowName}}
\begin{quote}

\sphinxAtStartPar
Name of row.
\end{quote}

\sphinxAtStartPar
\sphinxcode{\sphinxupquote{p\_iRow}}
\begin{quote}

\sphinxAtStartPar
Pointer to index of row.
\end{quote}
\end{quote}
\end{quote}

\subsubsection{COPT\_GetQConstrIdx}
\label{\detokenize{capiref:copt-getqconstridx}}\begin{quote}

\sphinxAtStartPar
\sphinxstylestrong{Synopsis}
\begin{quote}

\sphinxAtStartPar
\sphinxcode{\sphinxupquote{int COPT\_GetQConstrIdx(copt\_prob *prob, const char *qConstrName, int *p\_iQConstr)}}
\end{quote}

\sphinxAtStartPar
\sphinxstylestrong{Description}
\begin{quote}

\sphinxAtStartPar
Get index of quadratic constraint by name.
\end{quote}

\sphinxAtStartPar
\sphinxstylestrong{Arguments}
\begin{quote}

\sphinxAtStartPar
\sphinxcode{\sphinxupquote{prob}}
\begin{quote}

\sphinxAtStartPar
The COPT problem.
\end{quote}

\sphinxAtStartPar
\sphinxcode{\sphinxupquote{qConstrName}}
\begin{quote}

\sphinxAtStartPar
Name of quadratic constraint.
\end{quote}

\sphinxAtStartPar
\sphinxcode{\sphinxupquote{p\_iQConstr}}
\begin{quote}

\sphinxAtStartPar
Pointer to index of quadratic constraint.
\end{quote}
\end{quote}
\end{quote}

\subsubsection{COPT\_GetPSDConstrIdx}
\label{\detokenize{capiref:copt-getpsdconstridx}}\begin{quote}

\sphinxAtStartPar
\sphinxstylestrong{Synopsis}
\begin{quote}

\sphinxAtStartPar
\sphinxcode{\sphinxupquote{int COPT\_GetPSDConstrIdx(copt\_prob *prob, const char *psdConstrName, int *p\_iPSDConstr)}}
\end{quote}

\sphinxAtStartPar
\sphinxstylestrong{Description}
\begin{quote}

\sphinxAtStartPar
Get index of PSD constraint by name.
\end{quote}

\sphinxAtStartPar
\sphinxstylestrong{Arguments}
\begin{quote}

\sphinxAtStartPar
\sphinxcode{\sphinxupquote{prob}}
\begin{quote}

\sphinxAtStartPar
The COPT problem.
\end{quote}

\sphinxAtStartPar
\sphinxcode{\sphinxupquote{psdConstrName}}
\begin{quote}

\sphinxAtStartPar
Name of PSD constraint.
\end{quote}

\sphinxAtStartPar
\sphinxcode{\sphinxupquote{p\_iPSDConstr}}
\begin{quote}

\sphinxAtStartPar
Pointer to index of PSD constraint.
\end{quote}
\end{quote}
\end{quote}

\subsubsection{COPT\_GetLMIConstrIdx}
\label{\detokenize{capiref:copt-getlmiconstridx}}\begin{quote}

\sphinxAtStartPar
\sphinxstylestrong{Synopsis}
\begin{quote}

\sphinxAtStartPar
\sphinxcode{\sphinxupquote{int COPT\_GetLMIConstrIdx(copt\_prob *prob, const char *lmiConstrName, int *p\_iLMIConstr)}}
\end{quote}

\sphinxAtStartPar
\sphinxstylestrong{Description}
\begin{quote}

\sphinxAtStartPar
Get index of LMI constraint by name.
\end{quote}

\sphinxAtStartPar
\sphinxstylestrong{Arguments}
\begin{quote}

\sphinxAtStartPar
\sphinxcode{\sphinxupquote{prob}}
\begin{quote}

\sphinxAtStartPar
The COPT problem.
\end{quote}

\sphinxAtStartPar
\sphinxcode{\sphinxupquote{lmiConstrName}}
\begin{quote}

\sphinxAtStartPar
Name of LMI constraint.
\end{quote}

\sphinxAtStartPar
\sphinxcode{\sphinxupquote{p\_iLMIConstr}}
\begin{quote}

\sphinxAtStartPar
Pointer to index of LMI constraint.
\end{quote}
\end{quote}
\end{quote}

\subsubsection{COPT\_GetIndicatorIdx}
\label{\detokenize{capiref:copt-getindicatoridx}}\begin{quote}

\sphinxAtStartPar
\sphinxstylestrong{Synopsis}
\begin{quote}

\sphinxAtStartPar
\sphinxcode{\sphinxupquote{int COPT\_GetIndicatorIdx(copt\_prob *prob, const char *indicatorName, int *p\_iIndicator)}}
\end{quote}

\sphinxAtStartPar
\sphinxstylestrong{Description}
\begin{quote}

\sphinxAtStartPar
Get index of indicator constraint by name.
\end{quote}

\sphinxAtStartPar
\sphinxstylestrong{Arguments}
\begin{quote}

\sphinxAtStartPar
\sphinxcode{\sphinxupquote{prob}}
\begin{quote}

\sphinxAtStartPar
The COPT problem.
\end{quote}

\sphinxAtStartPar
\sphinxcode{\sphinxupquote{indicatorName}}
\begin{quote}

\sphinxAtStartPar
Name of indicator constraint.
\end{quote}

\sphinxAtStartPar
\sphinxcode{\sphinxupquote{p\_iIndicator}}
\begin{quote}

\sphinxAtStartPar
Pointer to index of indicator constraint.
\end{quote}
\end{quote}
\end{quote}

\subsubsection{COPT\_GetAffineConeIdx}
\label{\detokenize{capiref:copt-getaffineconeidx}}\begin{quote}

\sphinxAtStartPar
\sphinxstylestrong{Synopsis}
\begin{quote}

\sphinxAtStartPar
\sphinxcode{\sphinxupquote{int COPT\_GetAffineConeIdx(copt\_prob *prob, const char *affConeName, int *p\_iAffCone)}}
\end{quote}

\sphinxAtStartPar
\sphinxstylestrong{Description}
\begin{quote}

\sphinxAtStartPar
Get index of affine cone by name.
\end{quote}

\sphinxAtStartPar
\sphinxstylestrong{Arguments}
\begin{quote}

\sphinxAtStartPar
\sphinxcode{\sphinxupquote{prob}}
\begin{quote}

\sphinxAtStartPar
The COPT problem.
\end{quote}

\sphinxAtStartPar
\sphinxcode{\sphinxupquote{affConeName}}
\begin{quote}

\sphinxAtStartPar
Name of affine cone.
\end{quote}

\sphinxAtStartPar
\sphinxcode{\sphinxupquote{p\_iAffCone}}
\begin{quote}

\sphinxAtStartPar
Pointer to index of affine cone.
\end{quote}
\end{quote}
\end{quote}

\subsubsection{COPT\_GetNLConstrIdx}
\label{\detokenize{capiref:copt-getnlconstridx}}\begin{quote}

\sphinxAtStartPar
\sphinxstylestrong{Synopsis}
\begin{quote}

\sphinxAtStartPar
\sphinxcode{\sphinxupquote{int COPT\_GetNLConstrIdx(copt\_prob *prob, const char *nlConstrName,
int *p\_iNLConstr)}}
\end{quote}

\sphinxAtStartPar
\sphinxstylestrong{Description}
\begin{quote}

\sphinxAtStartPar
Retrieve the index of a nonlinear expression constraint by name.
\end{quote}

\sphinxAtStartPar
\sphinxstylestrong{Arguments}
\begin{quote}

\sphinxAtStartPar
\sphinxcode{\sphinxupquote{prob}}
\begin{quote}

\sphinxAtStartPar
The COPT problem.
\end{quote}

\sphinxAtStartPar
\sphinxcode{\sphinxupquote{nlConstrName}}
\begin{quote}

\sphinxAtStartPar
Name of the nonlinear expression constraint.
\end{quote}

\sphinxAtStartPar
\sphinxcode{\sphinxupquote{p\_iNLConstr}}
\begin{quote}

\sphinxAtStartPar
Pointer to the internal index of the nonlinear expression constraint in COPT.
\end{quote}
\end{quote}
\end{quote}

\subsubsection{COPT\_GetColInfo}
\label{\detokenize{capiref:copt-getcolinfo}}\begin{quote}

\sphinxAtStartPar
\sphinxstylestrong{Synopsis}
\begin{quote}

\sphinxAtStartPar
\sphinxcode{\sphinxupquote{int COPT\_GetColInfo(copt\_prob *prob, const char *infoName, int num, const int *list, double *info)}}
\end{quote}

\sphinxAtStartPar
\sphinxstylestrong{Description}
\begin{quote}

\sphinxAtStartPar
Get information of column. If \sphinxcode{\sphinxupquote{list}} is \sphinxcode{\sphinxupquote{NULL}}, then information
of the first \sphinxcode{\sphinxupquote{num}} columns will be returned.
\end{quote}

\sphinxAtStartPar
\sphinxstylestrong{Arguments}
\begin{quote}

\sphinxAtStartPar
\sphinxcode{\sphinxupquote{prob}}
\begin{quote}

\sphinxAtStartPar
The COPT problem.
\end{quote}

\sphinxAtStartPar
\sphinxcode{\sphinxupquote{infoName}}
\begin{quote}

\sphinxAtStartPar
Name of information. Please refer to {\hyperref[\detokenize{capiref:chapapi-info}]{\sphinxcrossref{\DUrole{std,std-ref}{Information}}}} for supported information.
\end{quote}

\sphinxAtStartPar
\sphinxcode{\sphinxupquote{num}}
\begin{quote}

\sphinxAtStartPar
Number of columns.
\end{quote}

\sphinxAtStartPar
\sphinxcode{\sphinxupquote{list}}
\begin{quote}

\sphinxAtStartPar
Index of columns. Can be \sphinxcode{\sphinxupquote{NULL}}.
\end{quote}

\sphinxAtStartPar
\sphinxcode{\sphinxupquote{info}}
\begin{quote}

\sphinxAtStartPar
Array of information.
\end{quote}
\end{quote}
\end{quote}

\subsubsection{COPT\_GetPSDColInfo}
\label{\detokenize{capiref:copt-getpsdcolinfo}}\begin{quote}

\sphinxAtStartPar
\sphinxstylestrong{Synopsis}
\begin{quote}

\sphinxAtStartPar
\sphinxcode{\sphinxupquote{int COPT\_GetPSDColInfo(copt\_prob *prob, const char *infoName, int iCol, double *info)}}
\end{quote}

\sphinxAtStartPar
\sphinxstylestrong{Description}
\begin{quote}

\sphinxAtStartPar
Get information of PSD variable.
\end{quote}

\sphinxAtStartPar
\sphinxstylestrong{Arguments}
\begin{quote}

\sphinxAtStartPar
\sphinxcode{\sphinxupquote{prob}}
\begin{quote}

\sphinxAtStartPar
The COPT problem.
\end{quote}

\sphinxAtStartPar
\sphinxcode{\sphinxupquote{infoName}}
\begin{quote}

\sphinxAtStartPar
Name of information. Please refer to {\hyperref[\detokenize{capiref:chapapi-info}]{\sphinxcrossref{\DUrole{std,std-ref}{Information}}}} for supported information.
\end{quote}

\sphinxAtStartPar
\sphinxcode{\sphinxupquote{iCol}}
\begin{quote}

\sphinxAtStartPar
Index of PSD variable.
\end{quote}

\sphinxAtStartPar
\sphinxcode{\sphinxupquote{info}}
\begin{quote}

\sphinxAtStartPar
Array of information.
\end{quote}
\end{quote}
\end{quote}

\subsubsection{COPT\_GetRowInfo}
\label{\detokenize{capiref:copt-getrowinfo}}\begin{quote}

\sphinxAtStartPar
\sphinxstylestrong{Synopsis}
\begin{quote}

\sphinxAtStartPar
\sphinxcode{\sphinxupquote{int COPT\_GetRowInfo(copt\_prob *prob, const char *infoName, int num, const int *list, double *info)}}
\end{quote}

\sphinxAtStartPar
\sphinxstylestrong{Description}
\begin{quote}

\sphinxAtStartPar
Get information of row. If \sphinxcode{\sphinxupquote{list}} is \sphinxcode{\sphinxupquote{NULL}}, then information
of the first \sphinxcode{\sphinxupquote{num}} rows will be returned.
\end{quote}

\sphinxAtStartPar
\sphinxstylestrong{Arguments}
\begin{quote}

\sphinxAtStartPar
\sphinxcode{\sphinxupquote{prob}}
\begin{quote}

\sphinxAtStartPar
The COPT problem.
\end{quote}

\sphinxAtStartPar
\sphinxcode{\sphinxupquote{infoName}}
\begin{quote}

\sphinxAtStartPar
Name of information. Please refer to {\hyperref[\detokenize{capiref:chapapi-info}]{\sphinxcrossref{\DUrole{std,std-ref}{Information}}}} for supported information.
\end{quote}

\sphinxAtStartPar
\sphinxcode{\sphinxupquote{num}}
\begin{quote}

\sphinxAtStartPar
Number of rows.
\end{quote}

\sphinxAtStartPar
\sphinxcode{\sphinxupquote{list}}
\begin{quote}

\sphinxAtStartPar
Index of rows. Can be \sphinxcode{\sphinxupquote{NULL}}.
\end{quote}

\sphinxAtStartPar
\sphinxcode{\sphinxupquote{info}}
\begin{quote}

\sphinxAtStartPar
Array of information.
\end{quote}
\end{quote}
\end{quote}

\subsubsection{COPT\_GetQConstrInfo}
\label{\detokenize{capiref:copt-getqconstrinfo}}\begin{quote}

\sphinxAtStartPar
\sphinxstylestrong{Synopsis}
\begin{quote}

\sphinxAtStartPar
\sphinxcode{\sphinxupquote{int COPT\_GetQConstrInfo(copt\_prob *prob, const char *infoName, int num, const int *list, double *info)}}
\end{quote}

\sphinxAtStartPar
\sphinxstylestrong{Description}
\begin{quote}

\sphinxAtStartPar
Get information of quadratic constraints. If \sphinxcode{\sphinxupquote{list}} is \sphinxcode{\sphinxupquote{NULL}}, then information
of the first \sphinxcode{\sphinxupquote{num}} quadratic constraints will be returned.
\end{quote}

\sphinxAtStartPar
\sphinxstylestrong{Arguments}
\begin{quote}

\sphinxAtStartPar
\sphinxcode{\sphinxupquote{prob}}
\begin{quote}

\sphinxAtStartPar
The COPT problem.
\end{quote}

\sphinxAtStartPar
\sphinxcode{\sphinxupquote{infoName}}
\begin{quote}

\sphinxAtStartPar
Name of information. Please refer to {\hyperref[\detokenize{capiref:chapapi-info}]{\sphinxcrossref{\DUrole{std,std-ref}{Information}}}} for supported information.
\end{quote}

\sphinxAtStartPar
\sphinxcode{\sphinxupquote{num}}
\begin{quote}

\sphinxAtStartPar
Number of quadratic constraints.
\end{quote}

\sphinxAtStartPar
\sphinxcode{\sphinxupquote{list}}
\begin{quote}

\sphinxAtStartPar
Index of quadratic constraints. Can be \sphinxcode{\sphinxupquote{NULL}}.
\end{quote}

\sphinxAtStartPar
\sphinxcode{\sphinxupquote{info}}
\begin{quote}

\sphinxAtStartPar
Array of information.
\end{quote}
\end{quote}
\end{quote}

\subsubsection{COPT\_GetPSDConstrInfo}
\label{\detokenize{capiref:copt-getpsdconstrinfo}}\begin{quote}

\sphinxAtStartPar
\sphinxstylestrong{Synopsis}
\begin{quote}

\sphinxAtStartPar
\sphinxcode{\sphinxupquote{int COPT\_GetPSDConstrInfo(copt\_prob *prob, const char *infoName, int num, const int* list, double *info)}}
\end{quote}

\sphinxAtStartPar
\sphinxstylestrong{Description}
\begin{quote}

\sphinxAtStartPar
Get information of PSD constraints. If \sphinxcode{\sphinxupquote{list}} is \sphinxcode{\sphinxupquote{NULL}}, then information
of the first \sphinxcode{\sphinxupquote{num}} PSD constraints will be returned.
\end{quote}

\sphinxAtStartPar
\sphinxstylestrong{Arguments}
\begin{quote}

\sphinxAtStartPar
\sphinxcode{\sphinxupquote{prob}}
\begin{quote}

\sphinxAtStartPar
The COPT problem.
\end{quote}

\sphinxAtStartPar
\sphinxcode{\sphinxupquote{infoName}}
\begin{quote}

\sphinxAtStartPar
Name of information. Please refer to {\hyperref[\detokenize{capiref:chapapi-info}]{\sphinxcrossref{\DUrole{std,std-ref}{Information}}}} for supported information.
\end{quote}

\sphinxAtStartPar
\sphinxcode{\sphinxupquote{num}}
\begin{quote}

\sphinxAtStartPar
Number of PSD constraints.
\end{quote}

\sphinxAtStartPar
\sphinxcode{\sphinxupquote{list}}
\begin{quote}

\sphinxAtStartPar
Index of PSD constraints. Can be \sphinxcode{\sphinxupquote{NULL}}.
\end{quote}

\sphinxAtStartPar
\sphinxcode{\sphinxupquote{info}}
\begin{quote}

\sphinxAtStartPar
Array of information.
\end{quote}
\end{quote}
\end{quote}

\subsubsection{COPT\_GetLMIConstrInfo}
\label{\detokenize{capiref:copt-getlmiconstrinfo}}\begin{quote}

\sphinxAtStartPar
\sphinxstylestrong{Synopsis}
\begin{quote}

\sphinxAtStartPar
\sphinxcode{\sphinxupquote{int COPT\_GetLMIConstrInfo(copt\_prob *prob, const char *infoName, int iLMI, double *info)}}
\end{quote}

\sphinxAtStartPar
\sphinxstylestrong{Description}
\begin{quote}

\sphinxAtStartPar
Get a set of information about LMI constraints.
\end{quote}

\sphinxAtStartPar
\sphinxstylestrong{Arguments}
\begin{quote}

\sphinxAtStartPar
\sphinxcode{\sphinxupquote{prob}}
\begin{quote}

\sphinxAtStartPar
The COPT problem.
\end{quote}

\sphinxAtStartPar
\sphinxcode{\sphinxupquote{infoName}}
\begin{quote}

\sphinxAtStartPar
Name of information. Possible values are: \sphinxcode{\sphinxupquote{COPT\_DBLINFO\_SLACK}} and \sphinxcode{\sphinxupquote{COPT\_DBLINFO\_DUAL}}.
\end{quote}

\sphinxAtStartPar
\sphinxcode{\sphinxupquote{iLMI}}
\begin{quote}

\sphinxAtStartPar
The index of the LMI constraint whose information is to be retrieved.
\end{quote}

\sphinxAtStartPar
\sphinxcode{\sphinxupquote{info}}
\begin{quote}

\sphinxAtStartPar
Array of information.
\end{quote}
\end{quote}
\end{quote}

\subsubsection{COPT\_GetNLConstrInfo}
\label{\detokenize{capiref:copt-getnlconstrinfo}}\begin{quote}

\sphinxAtStartPar
\sphinxstylestrong{Synopsis}
\begin{quote}

\sphinxAtStartPar
\sphinxcode{\sphinxupquote{int COPT\_GetNLConstrInfo(copt\_prob *prob, const char *infoName,}}
\sphinxcode{\sphinxupquote{int num, const int *list, double *info)}}
\end{quote}

\sphinxAtStartPar
\sphinxstylestrong{Description}
\begin{quote}

\sphinxAtStartPar
Retrieve information about nonlinear expression constraints.
If \sphinxcode{\sphinxupquote{list}} is \sphinxcode{\sphinxupquote{NULL}}, the function returns information for
the first \sphinxcode{\sphinxupquote{num}} nonlinear expression constraints.
\end{quote}

\sphinxAtStartPar
\sphinxstylestrong{Arguments}
\begin{quote}

\sphinxAtStartPar
\sphinxcode{\sphinxupquote{prob}}
\begin{quote}

\sphinxAtStartPar
The COPT problem.
\end{quote}

\sphinxAtStartPar
\sphinxcode{\sphinxupquote{infoName}}
\begin{quote}

\sphinxAtStartPar
Name of the information to retrieve.
Currently supported nonlinear expression constraint
information is: \sphinxcode{\sphinxupquote{"LB"}} , \sphinxcode{\sphinxupquote{"UB"}} and \sphinxcode{\sphinxupquote{"Slack"}} .
\end{quote}

\sphinxAtStartPar
\sphinxcode{\sphinxupquote{num}}
\begin{quote}

\sphinxAtStartPar
Number of nonlinear expression constraints to retrieve information for.
\end{quote}

\sphinxAtStartPar
\sphinxcode{\sphinxupquote{list}}
\begin{quote}

\sphinxAtStartPar
Index list of nonlinear expression constraints to retrieve information for.
Can be \sphinxcode{\sphinxupquote{NULL}}.
\end{quote}

\sphinxAtStartPar
\sphinxcode{\sphinxupquote{info}}
\begin{quote}

\sphinxAtStartPar
Array to store the retrieved information.
\end{quote}
\end{quote}
\end{quote}

\subsubsection{COPT\_GetColType}
\label{\detokenize{capiref:copt-getcoltype}}\begin{quote}

\sphinxAtStartPar
\sphinxstylestrong{Synopsis}
\begin{quote}

\sphinxAtStartPar
\sphinxcode{\sphinxupquote{int COPT\_GetColType(copt\_prob *prob, int num, const int *list, char *type)}}
\end{quote}

\sphinxAtStartPar
\sphinxstylestrong{Description}
\begin{quote}

\sphinxAtStartPar
Get types of columns. If \sphinxcode{\sphinxupquote{list}} is \sphinxcode{\sphinxupquote{NULL}}, then types
of the first \sphinxcode{\sphinxupquote{num}} columns will be returned.
\end{quote}

\sphinxAtStartPar
\sphinxstylestrong{Arguments}
\begin{quote}

\sphinxAtStartPar
\sphinxcode{\sphinxupquote{prob}}
\begin{quote}

\sphinxAtStartPar
The COPT problem.
\end{quote}

\sphinxAtStartPar
\sphinxcode{\sphinxupquote{num}}
\begin{quote}

\sphinxAtStartPar
Number of columns.
\end{quote}

\sphinxAtStartPar
\sphinxcode{\sphinxupquote{list}}
\begin{quote}

\sphinxAtStartPar
Index of columns. Can be \sphinxcode{\sphinxupquote{NULL}}.
\end{quote}

\sphinxAtStartPar
\sphinxcode{\sphinxupquote{type}}
\begin{quote}

\sphinxAtStartPar
Types of columns.
\end{quote}
\end{quote}
\end{quote}

\subsubsection{COPT\_GetColBasis}
\label{\detokenize{capiref:copt-getcolbasis}}\begin{quote}

\sphinxAtStartPar
\sphinxstylestrong{Synopsis}
\begin{quote}

\sphinxAtStartPar
\sphinxcode{\sphinxupquote{int COPT\_GetColBasis(copt\_prob *prob, int num, const int *list, int *colBasis)}}
\end{quote}

\sphinxAtStartPar
\sphinxstylestrong{Description}
\begin{quote}

\sphinxAtStartPar
Get basis status of columns. If \sphinxcode{\sphinxupquote{list}} is \sphinxcode{\sphinxupquote{NULL}}, then basis status
of the first \sphinxcode{\sphinxupquote{num}} columns will be returned.
\end{quote}

\sphinxAtStartPar
\sphinxstylestrong{Arguments}
\begin{quote}

\sphinxAtStartPar
\sphinxcode{\sphinxupquote{prob}}
\begin{quote}

\sphinxAtStartPar
The COPT problem.
\end{quote}

\sphinxAtStartPar
\sphinxcode{\sphinxupquote{num}}
\begin{quote}

\sphinxAtStartPar
Number of columns.
\end{quote}

\sphinxAtStartPar
\sphinxcode{\sphinxupquote{list}}
\begin{quote}

\sphinxAtStartPar
Index of columns. Can be \sphinxcode{\sphinxupquote{NULL}}.
\end{quote}

\sphinxAtStartPar
\sphinxcode{\sphinxupquote{colBasis}}
\begin{quote}

\sphinxAtStartPar
Basis status of columns.
\end{quote}
\end{quote}
\end{quote}

\subsubsection{COPT\_GetRowBasis}
\label{\detokenize{capiref:copt-getrowbasis}}\begin{quote}

\sphinxAtStartPar
\sphinxstylestrong{Synopsis}
\begin{quote}

\sphinxAtStartPar
\sphinxcode{\sphinxupquote{int COPT\_GetRowBasis(copt\_prob *prob, int num, const int *list, int *rowBasis)}}
\end{quote}

\sphinxAtStartPar
\sphinxstylestrong{Description}
\begin{quote}

\sphinxAtStartPar
Get basis status of rows. If \sphinxcode{\sphinxupquote{list}} is \sphinxcode{\sphinxupquote{NULL}}, then basis status
of the first \sphinxcode{\sphinxupquote{num}} rows will be returned.
\end{quote}

\sphinxAtStartPar
\sphinxstylestrong{Arguments}
\begin{quote}

\sphinxAtStartPar
\sphinxcode{\sphinxupquote{prob}}
\begin{quote}

\sphinxAtStartPar
The COPT problem.
\end{quote}

\sphinxAtStartPar
\sphinxcode{\sphinxupquote{num}}
\begin{quote}

\sphinxAtStartPar
Number of rows.
\end{quote}

\sphinxAtStartPar
\sphinxcode{\sphinxupquote{list}}
\begin{quote}

\sphinxAtStartPar
Index of rows. Can be \sphinxcode{\sphinxupquote{NULL}}.
\end{quote}

\sphinxAtStartPar
\sphinxcode{\sphinxupquote{rowBasis}}
\begin{quote}

\sphinxAtStartPar
Basis status of rows.
\end{quote}
\end{quote}
\end{quote}

\subsubsection{COPT\_GetQConstrSense}
\label{\detokenize{capiref:copt-getqconstrsense}}\begin{quote}

\sphinxAtStartPar
\sphinxstylestrong{Synopsis}
\begin{quote}

\sphinxAtStartPar
\sphinxcode{\sphinxupquote{int COPT\_GetQConstrSense(copt\_prob *prob, int num, const int *list, char *sense)}}
\end{quote}

\sphinxAtStartPar
\sphinxstylestrong{Description}
\begin{quote}

\sphinxAtStartPar
Get senses of quadratic constraints. If \sphinxcode{\sphinxupquote{list}} is \sphinxcode{\sphinxupquote{NULL}}, then types
of the first \sphinxcode{\sphinxupquote{num}} quadratic constraints will be returned.
\end{quote}

\sphinxAtStartPar
\sphinxstylestrong{Arguments}
\begin{quote}

\sphinxAtStartPar
\sphinxcode{\sphinxupquote{prob}}
\begin{quote}

\sphinxAtStartPar
The COPT problem.
\end{quote}

\sphinxAtStartPar
\sphinxcode{\sphinxupquote{num}}
\begin{quote}

\sphinxAtStartPar
Number of quadratic constraints.
\end{quote}

\sphinxAtStartPar
\sphinxcode{\sphinxupquote{list}}
\begin{quote}

\sphinxAtStartPar
Index of quadratic constraints. Can be \sphinxcode{\sphinxupquote{NULL}}.
\end{quote}

\sphinxAtStartPar
\sphinxcode{\sphinxupquote{sense}}
\begin{quote}

\sphinxAtStartPar
Array of senses.
\end{quote}
\end{quote}
\end{quote}

\subsubsection{COPT\_GetQConstrRhs}
\label{\detokenize{capiref:copt-getqconstrrhs}}\begin{quote}

\sphinxAtStartPar
\sphinxstylestrong{Synopsis}
\begin{quote}

\sphinxAtStartPar
\sphinxcode{\sphinxupquote{int COPT\_GetQConstrRhs(copt\_prob *prob, int num, const int *list, double *rhs)}}
\end{quote}

\sphinxAtStartPar
\sphinxstylestrong{Description}
\begin{quote}

\sphinxAtStartPar
Get RHS of quadratic constraints. If \sphinxcode{\sphinxupquote{list}} is \sphinxcode{\sphinxupquote{NULL}}, then types
of the first \sphinxcode{\sphinxupquote{num}} quadratic constraints will be returned.
\end{quote}

\sphinxAtStartPar
\sphinxstylestrong{Arguments}
\begin{quote}

\sphinxAtStartPar
\sphinxcode{\sphinxupquote{prob}}
\begin{quote}

\sphinxAtStartPar
The COPT problem.
\end{quote}

\sphinxAtStartPar
\sphinxcode{\sphinxupquote{num}}
\begin{quote}

\sphinxAtStartPar
Number of quadratic constraints.
\end{quote}

\sphinxAtStartPar
\sphinxcode{\sphinxupquote{list}}
\begin{quote}

\sphinxAtStartPar
Index of quadratic constraints. Can be \sphinxcode{\sphinxupquote{NULL}}.
\end{quote}

\sphinxAtStartPar
\sphinxcode{\sphinxupquote{rhs}}
\begin{quote}

\sphinxAtStartPar
Array of RHS.
\end{quote}
\end{quote}
\end{quote}

\subsubsection{COPT\_GetColName}
\label{\detokenize{capiref:copt-getcolname}}\begin{quote}

\sphinxAtStartPar
\sphinxstylestrong{Synopsis}
\begin{quote}

\sphinxAtStartPar
\sphinxcode{\sphinxupquote{int COPT\_GetColName(copt\_prob *prob, int iCol, char *buff, int buffSize, int *pReqSize)}}
\end{quote}

\sphinxAtStartPar
\sphinxstylestrong{Description}
\begin{quote}

\sphinxAtStartPar
Get name of column by index. If memory of \sphinxcode{\sphinxupquote{buff}} is not sufficient,
then return the first \sphinxcode{\sphinxupquote{buffSize}} length of sub\sphinxhyphen{}string, and the length
of name requested by \sphinxcode{\sphinxupquote{pReqSize}}. If \sphinxcode{\sphinxupquote{buff}} is \sphinxcode{\sphinxupquote{NULL}}, then we can
get the length of name requested by \sphinxcode{\sphinxupquote{pReqSize}}.
\end{quote}

\sphinxAtStartPar
\sphinxstylestrong{Arguments}
\begin{quote}

\sphinxAtStartPar
\sphinxcode{\sphinxupquote{prob}}
\begin{quote}

\sphinxAtStartPar
The COPT problem.
\end{quote}

\sphinxAtStartPar
\sphinxcode{\sphinxupquote{iCol}}
\begin{quote}

\sphinxAtStartPar
Index of column.
\end{quote}

\sphinxAtStartPar
\sphinxcode{\sphinxupquote{buff}}
\begin{quote}

\sphinxAtStartPar
Buffer for storing name.
\end{quote}

\sphinxAtStartPar
\sphinxcode{\sphinxupquote{buffSize}}
\begin{quote}

\sphinxAtStartPar
Length of the buffer.
\end{quote}

\sphinxAtStartPar
\sphinxcode{\sphinxupquote{pReqSize}}
\begin{quote}

\sphinxAtStartPar
Length of the requested name. Can be \sphinxcode{\sphinxupquote{NULL}}.
\end{quote}
\end{quote}
\end{quote}

\subsubsection{COPT\_GetPSDColName}
\label{\detokenize{capiref:copt-getpsdcolname}}\begin{quote}

\sphinxAtStartPar
\sphinxstylestrong{Synopsis}
\begin{quote}

\sphinxAtStartPar
\sphinxcode{\sphinxupquote{int COPT\_GetPSDColName(copt\_prob *prob, int iPSDCol, char *buff, int buffSize, int *pReqSize)}}
\end{quote}

\sphinxAtStartPar
\sphinxstylestrong{Description}
\begin{quote}

\sphinxAtStartPar
Get name of PSD variable by index. If memory of \sphinxcode{\sphinxupquote{buff}} is not sufficient,
then return the first \sphinxcode{\sphinxupquote{buffSize}} length of sub\sphinxhyphen{}string, and the length
of name requested by \sphinxcode{\sphinxupquote{pReqSize}}. If \sphinxcode{\sphinxupquote{buff}} is \sphinxcode{\sphinxupquote{NULL}}, then we can
get the length of name requested by \sphinxcode{\sphinxupquote{pReqSize}}.
\end{quote}

\sphinxAtStartPar
\sphinxstylestrong{Arguments}
\begin{quote}

\sphinxAtStartPar
\sphinxcode{\sphinxupquote{prob}}
\begin{quote}

\sphinxAtStartPar
The COPT problem.
\end{quote}

\sphinxAtStartPar
\sphinxcode{\sphinxupquote{iPSDCol}}
\begin{quote}

\sphinxAtStartPar
Index of PSD variable.
\end{quote}

\sphinxAtStartPar
\sphinxcode{\sphinxupquote{buff}}
\begin{quote}

\sphinxAtStartPar
Buffer for storing name.
\end{quote}

\sphinxAtStartPar
\sphinxcode{\sphinxupquote{buffSize}}
\begin{quote}

\sphinxAtStartPar
Length of the buffer.
\end{quote}

\sphinxAtStartPar
\sphinxcode{\sphinxupquote{pReqSize}}
\begin{quote}

\sphinxAtStartPar
Length of the requested name. Can be \sphinxcode{\sphinxupquote{NULL}}.
\end{quote}
\end{quote}
\end{quote}

\subsubsection{COPT\_GetRowName}
\label{\detokenize{capiref:copt-getrowname}}\begin{quote}

\sphinxAtStartPar
\sphinxstylestrong{Synopsis}
\begin{quote}

\sphinxAtStartPar
\sphinxcode{\sphinxupquote{int COPT\_GetRowName(copt\_prob *prob, int iRow, char *buff, int buffSize, int *pReqSize)}}
\end{quote}

\sphinxAtStartPar
\sphinxstylestrong{Description}
\begin{quote}

\sphinxAtStartPar
Get name of row by index. If memory of \sphinxcode{\sphinxupquote{buff}} is not sufficient,
then return the first \sphinxcode{\sphinxupquote{buffSize}} length of sub\sphinxhyphen{}string, and the length
of name requested by \sphinxcode{\sphinxupquote{pReqSize}}. If \sphinxcode{\sphinxupquote{buff}} is \sphinxcode{\sphinxupquote{NULL}}, then we can
get the length of name requested by \sphinxcode{\sphinxupquote{pReqSize}}.
\end{quote}

\sphinxAtStartPar
\sphinxstylestrong{Arguments}
\begin{quote}

\sphinxAtStartPar
\sphinxcode{\sphinxupquote{prob}}
\begin{quote}

\sphinxAtStartPar
The COPT problem.
\end{quote}

\sphinxAtStartPar
\sphinxcode{\sphinxupquote{iRow}}
\begin{quote}

\sphinxAtStartPar
Index of row.
\end{quote}

\sphinxAtStartPar
\sphinxcode{\sphinxupquote{buff}}
\begin{quote}

\sphinxAtStartPar
Buffer for storing name.
\end{quote}

\sphinxAtStartPar
\sphinxcode{\sphinxupquote{buffSize}}
\begin{quote}

\sphinxAtStartPar
Length of the buffer.
\end{quote}

\sphinxAtStartPar
\sphinxcode{\sphinxupquote{pReqSize}}
\begin{quote}

\sphinxAtStartPar
Length of the requested name. Can be \sphinxcode{\sphinxupquote{NULL}}.
\end{quote}
\end{quote}
\end{quote}

\subsubsection{COPT\_GetQConstrName}
\label{\detokenize{capiref:copt-getqconstrname}}\begin{quote}

\sphinxAtStartPar
\sphinxstylestrong{Synopsis}
\begin{quote}

\sphinxAtStartPar
\sphinxcode{\sphinxupquote{int COPT\_GetQConstrName(copt\_prob *prob, int iQConstr, char *buff, int buffSize, int *pReqSize)}}
\end{quote}

\sphinxAtStartPar
\sphinxstylestrong{Description}
\begin{quote}

\sphinxAtStartPar
Get name of quadratic constraint by index. If memory of \sphinxcode{\sphinxupquote{buff}} is not sufficient,
then return the first \sphinxcode{\sphinxupquote{buffSize}} length of sub\sphinxhyphen{}string, and the length
of name requested by \sphinxcode{\sphinxupquote{pReqSize}}. If \sphinxcode{\sphinxupquote{buff}} is \sphinxcode{\sphinxupquote{NULL}}, then we can
get the length of name requested by \sphinxcode{\sphinxupquote{pReqSize}}.
\end{quote}

\sphinxAtStartPar
\sphinxstylestrong{Arguments}
\begin{quote}

\sphinxAtStartPar
\sphinxcode{\sphinxupquote{prob}}
\begin{quote}

\sphinxAtStartPar
The COPT problem.
\end{quote}

\sphinxAtStartPar
\sphinxcode{\sphinxupquote{iQConstr\textasciigrave{}}}
\begin{quote}

\sphinxAtStartPar
Index of quadratic constraint.
\end{quote}

\sphinxAtStartPar
\sphinxcode{\sphinxupquote{buff}}
\begin{quote}

\sphinxAtStartPar
Buffer for storing name.
\end{quote}

\sphinxAtStartPar
\sphinxcode{\sphinxupquote{buffSize}}
\begin{quote}

\sphinxAtStartPar
Length of the buffer.
\end{quote}

\sphinxAtStartPar
\sphinxcode{\sphinxupquote{pReqSize}}
\begin{quote}

\sphinxAtStartPar
Length of the requested name. Can be \sphinxcode{\sphinxupquote{NULL}}.
\end{quote}
\end{quote}
\end{quote}

\subsubsection{COPT\_GetPSDConstrName}
\label{\detokenize{capiref:copt-getpsdconstrname}}\begin{quote}

\sphinxAtStartPar
\sphinxstylestrong{Synopsis}
\begin{quote}

\sphinxAtStartPar
\sphinxcode{\sphinxupquote{int COPT\_GetPSDConstrName(copt\_prob *prob, int iPSDConstr, char *buff, int buffSize, int *pReqSize)}}
\end{quote}

\sphinxAtStartPar
\sphinxstylestrong{Description}
\begin{quote}

\sphinxAtStartPar
Get name of PSD constraint by index. If memory of \sphinxcode{\sphinxupquote{buff}} is not sufficient,
then return the first \sphinxcode{\sphinxupquote{buffSize}} length of sub\sphinxhyphen{}string, and the length
of name requested by \sphinxcode{\sphinxupquote{pReqSize}}. If \sphinxcode{\sphinxupquote{buff}} is \sphinxcode{\sphinxupquote{NULL}}, then we can
get the length of name requested by \sphinxcode{\sphinxupquote{pReqSize}}.
\end{quote}

\sphinxAtStartPar
\sphinxstylestrong{Arguments}
\begin{quote}

\sphinxAtStartPar
\sphinxcode{\sphinxupquote{prob}}
\begin{quote}

\sphinxAtStartPar
The COPT problem.
\end{quote}

\sphinxAtStartPar
\sphinxcode{\sphinxupquote{iPSDConstr\textasciigrave{}}}
\begin{quote}

\sphinxAtStartPar
Index of PSD constraint.
\end{quote}

\sphinxAtStartPar
\sphinxcode{\sphinxupquote{buff}}
\begin{quote}

\sphinxAtStartPar
Buffer for storing name.
\end{quote}

\sphinxAtStartPar
\sphinxcode{\sphinxupquote{buffSize}}
\begin{quote}

\sphinxAtStartPar
Length of the buffer.
\end{quote}

\sphinxAtStartPar
\sphinxcode{\sphinxupquote{pReqSize}}
\begin{quote}

\sphinxAtStartPar
Length of the requested name. Can be \sphinxcode{\sphinxupquote{NULL}}.
\end{quote}
\end{quote}
\end{quote}

\subsubsection{COPT\_GetLMIConstrName}
\label{\detokenize{capiref:copt-getlmiconstrname}}\begin{quote}

\sphinxAtStartPar
\sphinxstylestrong{Synopsis}
\begin{quote}

\sphinxAtStartPar
\sphinxcode{\sphinxupquote{int COPT\_CALL COPT\_GetLMIConstrName(copt\_prob *prob, int iLMIConstr, char *buff, int buffSize, int *pReqSize)}}
\end{quote}

\sphinxAtStartPar
\sphinxstylestrong{Description}
\begin{quote}

\sphinxAtStartPar
Get name of LMI constraint by index. If memory of \sphinxcode{\sphinxupquote{buff}} is not sufficient,
then return the first \sphinxcode{\sphinxupquote{buffSize}} length of sub\sphinxhyphen{}string, and the length
of name requested by \sphinxcode{\sphinxupquote{pReqSize}}. If \sphinxcode{\sphinxupquote{buff}} is \sphinxcode{\sphinxupquote{NULL}}, then we can
get the length of name requested by \sphinxcode{\sphinxupquote{pReqSize}}.
\end{quote}

\sphinxAtStartPar
\sphinxstylestrong{Arguments}
\begin{quote}

\sphinxAtStartPar
\sphinxcode{\sphinxupquote{prob}}
\begin{quote}

\sphinxAtStartPar
The COPT problem.
\end{quote}

\sphinxAtStartPar
\sphinxcode{\sphinxupquote{iLMIConstr}}
\begin{quote}

\sphinxAtStartPar
Index of LMI constraint.
\end{quote}

\sphinxAtStartPar
\sphinxcode{\sphinxupquote{buff}}
\begin{quote}

\sphinxAtStartPar
Buffer for storing name.
\end{quote}

\sphinxAtStartPar
\sphinxcode{\sphinxupquote{buffSize}}
\begin{quote}

\sphinxAtStartPar
Length of the buffer.
\end{quote}

\sphinxAtStartPar
\sphinxcode{\sphinxupquote{pReqSize}}
\begin{quote}

\sphinxAtStartPar
Length of the requested name. Can be \sphinxcode{\sphinxupquote{NULL}}.
\end{quote}
\end{quote}
\end{quote}

\subsubsection{COPT\_GetIndicatorName}
\label{\detokenize{capiref:copt-getindicatorname}}\begin{quote}

\sphinxAtStartPar
\sphinxstylestrong{Synopsis}
\begin{quote}

\sphinxAtStartPar
\sphinxcode{\sphinxupquote{int COPT\_GetIndicatorName(copt\_prob *prob, int iIndicator, char *buff, int buffSize, int *pReqSize)}}
\end{quote}

\sphinxAtStartPar
\sphinxstylestrong{Description}
\begin{quote}

\sphinxAtStartPar
Get name of indicator constraints by index. If memory of \sphinxcode{\sphinxupquote{buff}} is not sufficient,
then return the first \sphinxcode{\sphinxupquote{buffSize}} length of sub\sphinxhyphen{}string, and the length
of name requested by \sphinxcode{\sphinxupquote{pReqSize}}. If \sphinxcode{\sphinxupquote{buff}} is \sphinxcode{\sphinxupquote{NULL}}, then we can
get the length of name requested by \sphinxcode{\sphinxupquote{pReqSize}}.
\end{quote}

\sphinxAtStartPar
\sphinxstylestrong{Arguments}
\begin{quote}

\sphinxAtStartPar
\sphinxcode{\sphinxupquote{prob}}
\begin{quote}

\sphinxAtStartPar
The COPT problem.
\end{quote}

\sphinxAtStartPar
\sphinxcode{\sphinxupquote{iIndicator}}
\begin{quote}

\sphinxAtStartPar
Index of indicator constraint.
\end{quote}

\sphinxAtStartPar
\sphinxcode{\sphinxupquote{buff}}
\begin{quote}

\sphinxAtStartPar
Buffer for storing name.
\end{quote}

\sphinxAtStartPar
\sphinxcode{\sphinxupquote{buffSize}}
\begin{quote}

\sphinxAtStartPar
Length of the buffer.
\end{quote}

\sphinxAtStartPar
\sphinxcode{\sphinxupquote{pReqSize}}
\begin{quote}

\sphinxAtStartPar
Length of the requested name. Can be \sphinxcode{\sphinxupquote{NULL}}.
\end{quote}
\end{quote}
\end{quote}

\subsubsection{COPT\_GetAffineConeName}
\label{\detokenize{capiref:copt-getaffineconename}}\begin{quote}

\sphinxAtStartPar
\sphinxstylestrong{Synopsis}
\begin{quote}

\sphinxAtStartPar
\sphinxcode{\sphinxupquote{int COPT\_GetAffineConeName(copt\_prob *prob, int iAffCone, char *buff, int buffSize, int *pReqSize)}}
\end{quote}

\sphinxAtStartPar
\sphinxstylestrong{Description}
\begin{quote}

\sphinxAtStartPar
Get name of affine cone by index. If memory of \sphinxcode{\sphinxupquote{buff}} is not sufficient,
then return the first \sphinxcode{\sphinxupquote{buffSize}} length of sub\sphinxhyphen{}string, and the length
of name requested by \sphinxcode{\sphinxupquote{pReqSize}}. If \sphinxcode{\sphinxupquote{buff}} is \sphinxcode{\sphinxupquote{NULL}}, then we can
get the length of name requested by \sphinxcode{\sphinxupquote{pReqSize}}.
\end{quote}

\sphinxAtStartPar
\sphinxstylestrong{Arguments}
\begin{quote}

\sphinxAtStartPar
\sphinxcode{\sphinxupquote{prob}}
\begin{quote}

\sphinxAtStartPar
The COPT problem.
\end{quote}

\sphinxAtStartPar
\sphinxcode{\sphinxupquote{iAffCone}}
\begin{quote}

\sphinxAtStartPar
Index of affine cone.
\end{quote}

\sphinxAtStartPar
\sphinxcode{\sphinxupquote{buff}}
\begin{quote}

\sphinxAtStartPar
Buffer for storing name.
\end{quote}

\sphinxAtStartPar
\sphinxcode{\sphinxupquote{buffSize}}
\begin{quote}

\sphinxAtStartPar
Length of the buffer.
\end{quote}

\sphinxAtStartPar
\sphinxcode{\sphinxupquote{pReqSize}}
\begin{quote}

\sphinxAtStartPar
Length of the requested name. Can be \sphinxcode{\sphinxupquote{NULL}}.
\end{quote}
\end{quote}
\end{quote}

\subsubsection{COPT\_GetNLConstrName}
\label{\detokenize{capiref:copt-getnlconstrname}}\begin{quote}

\sphinxAtStartPar
\sphinxstylestrong{Synopsis}
\begin{quote}

\sphinxAtStartPar
\sphinxcode{\sphinxupquote{int COPT\_GetNLConstrName(copt\_prob *prob, int iNLConstr, char *buff,}}
\sphinxcode{\sphinxupquote{int buffSize, int *pReqSize)}}
\end{quote}

\sphinxAtStartPar
\sphinxstylestrong{Description}
\begin{quote}

\sphinxAtStartPar
Retrieve the name of nonlinear constraint by index.
If the provided \sphinxcode{\sphinxupquote{buff}} length is insufficient, the function
returns a substring of length \sphinxcode{\sphinxupquote{buffSize}}, and the full required
string length is returned via \sphinxcode{\sphinxupquote{pReqSize}}.

\sphinxAtStartPar
If \sphinxcode{\sphinxupquote{buff}} is \sphinxcode{\sphinxupquote{NULL}}, the function returns the required string
length via \sphinxcode{\sphinxupquote{pReqSize}}.
\end{quote}

\sphinxAtStartPar
\sphinxstylestrong{Arguments}
\begin{quote}

\sphinxAtStartPar
\sphinxcode{\sphinxupquote{prob}}
\begin{quote}

\sphinxAtStartPar
The COPT problem.
\end{quote}

\sphinxAtStartPar
\sphinxcode{\sphinxupquote{iNLConstr}}
\begin{quote}

\sphinxAtStartPar
Index of the nonlinear expression constraint.
\end{quote}

\sphinxAtStartPar
\sphinxcode{\sphinxupquote{buff}}
\begin{quote}

\sphinxAtStartPar
Array to store the retrieved name string.
\end{quote}

\sphinxAtStartPar
\sphinxcode{\sphinxupquote{buffSize}}
\begin{quote}

\sphinxAtStartPar
Size of the provided array.
\end{quote}

\sphinxAtStartPar
\sphinxcode{\sphinxupquote{pReqSize}}
\begin{quote}

\sphinxAtStartPar
Pointer to store the minimum required string length
for the full name. Can be \sphinxcode{\sphinxupquote{NULL}}.
\end{quote}
\end{quote}
\end{quote}

\subsubsection{COPT\_GetLMIConstrRhs}
\label{\detokenize{capiref:copt-getlmiconstrrhs}}\begin{quote}

\sphinxAtStartPar
\sphinxstylestrong{Synopsis}
\begin{quote}

\sphinxAtStartPar
\sphinxcode{\sphinxupquote{int COPT\_GetLMIConstrRhs(copt\_prob *prob, int num, const int *list, int *constMatIdx)}}
\end{quote}

\sphinxAtStartPar
\sphinxstylestrong{Description}
\begin{quote}

\sphinxAtStartPar
Get the constant\sphinxhyphen{}term symmetric matrix of \sphinxcode{\sphinxupquote{num}} LMI constraints.
\end{quote}

\sphinxAtStartPar
\sphinxstylestrong{Arguments}
\begin{quote}

\sphinxAtStartPar
\sphinxcode{\sphinxupquote{prob}}
\begin{quote}

\sphinxAtStartPar
The COPT problem.
\end{quote}

\sphinxAtStartPar
\sphinxcode{\sphinxupquote{num}}
\begin{quote}

\sphinxAtStartPar
Number of LMI constraints.
\end{quote}

\sphinxAtStartPar
\sphinxcode{\sphinxupquote{list}}
\begin{quote}

\sphinxAtStartPar
Index of LMI constraints.
\end{quote}

\sphinxAtStartPar
\sphinxcode{\sphinxupquote{constMatIdx}}
\begin{quote}

\sphinxAtStartPar
Index of constant\sphinxhyphen{}term symmetric in the LMI constraints.
\end{quote}
\end{quote}
\end{quote}

\subsection{Accessing and setting parameters}
\label{\detokenize{capiref:accessing-and-setting-parameters}}\label{\detokenize{capiref:chapapi-getparam}}

\subsubsection{COPT\_SetIntParam}
\label{\detokenize{capiref:copt-setintparam}}\begin{quote}

\sphinxAtStartPar
\sphinxstylestrong{Synopsis}
\begin{quote}

\sphinxAtStartPar
\sphinxcode{\sphinxupquote{int COPT\_SetIntParam(copt\_prob *prob, const char *paramName, int intParam)}}
\end{quote}

\sphinxAtStartPar
\sphinxstylestrong{Description}
\begin{quote}

\sphinxAtStartPar
Sets an integer parameter.
\end{quote}

\sphinxAtStartPar
\sphinxstylestrong{Arguments}
\begin{quote}

\sphinxAtStartPar
\sphinxcode{\sphinxupquote{prob}}
\begin{quote}

\sphinxAtStartPar
The COPT problem.
\end{quote}

\sphinxAtStartPar
\sphinxcode{\sphinxupquote{paramName}}
\begin{quote}

\sphinxAtStartPar
The name of the integer parameter.
\end{quote}

\sphinxAtStartPar
\sphinxcode{\sphinxupquote{intParam}}
\begin{quote}

\sphinxAtStartPar
The value of the integer parameter.
\end{quote}
\end{quote}
\end{quote}

\subsubsection{COPT\_GetIntParam, COPT\_GetIntParamDef/Min/Max}
\label{\detokenize{capiref:copt-getintparam-copt-getintparamdef-min-max}}\begin{quote}

\sphinxAtStartPar
\sphinxstylestrong{Synopsis}
\begin{quote}

\sphinxAtStartPar
\sphinxcode{\sphinxupquote{int COPT\_GetIntParam(copt\_prob *prob, const char *paramName, int *p\_intParam)}}

\sphinxAtStartPar
\sphinxcode{\sphinxupquote{int COPT\_GetIntParamDef(copt\_prob *prob, const char *paramName, int *p\_intParam)}}

\sphinxAtStartPar
\sphinxcode{\sphinxupquote{int COPT\_GetIntParamMin(copt\_prob *prob, const char *paramName, int *p\_intParam)}}

\sphinxAtStartPar
\sphinxcode{\sphinxupquote{int COPT\_GetIntParamMax(copt\_prob *prob, const char *paramName, int *p\_intParam)}}
\end{quote}

\sphinxAtStartPar
\sphinxstylestrong{Description}
\begin{quote}

\sphinxAtStartPar
Gets the
\begin{quote}

\sphinxAtStartPar
current

\sphinxAtStartPar
default

\sphinxAtStartPar
minimal

\sphinxAtStartPar
maximal
\end{quote}

\sphinxAtStartPar
value of an integer parameter.
\end{quote}

\sphinxAtStartPar
\sphinxstylestrong{Arguments}
\begin{quote}

\sphinxAtStartPar
\sphinxcode{\sphinxupquote{prob}}
\begin{quote}

\sphinxAtStartPar
The COPT problem.
\end{quote}

\sphinxAtStartPar
\sphinxcode{\sphinxupquote{paramName}}
\begin{quote}

\sphinxAtStartPar
The name of the integer parameter.
\end{quote}

\sphinxAtStartPar
\sphinxcode{\sphinxupquote{p\_intParam}}
\begin{quote}

\sphinxAtStartPar
Pointer to the value of the integer parameter.
\end{quote}
\end{quote}
\end{quote}

\subsubsection{COPT\_SetDblParam}
\label{\detokenize{capiref:copt-setdblparam}}\begin{quote}

\sphinxAtStartPar
\sphinxstylestrong{Synopsis}
\begin{quote}

\sphinxAtStartPar
\sphinxcode{\sphinxupquote{int COPT\_SetDblParam(copt\_prob *prob, const char *paramName, double dblParam)}}
\end{quote}

\sphinxAtStartPar
\sphinxstylestrong{Description}
\begin{quote}

\sphinxAtStartPar
Sets a double parameter.
\end{quote}

\sphinxAtStartPar
\sphinxstylestrong{Arguments}
\begin{quote}

\sphinxAtStartPar
\sphinxcode{\sphinxupquote{prob}}
\begin{quote}

\sphinxAtStartPar
The COPT problem.
\end{quote}

\sphinxAtStartPar
\sphinxcode{\sphinxupquote{paramName}}
\begin{quote}

\sphinxAtStartPar
The name of the double parameter.
\end{quote}

\sphinxAtStartPar
\sphinxcode{\sphinxupquote{dblParam}}
\begin{quote}

\sphinxAtStartPar
The value of the double parameter.
\end{quote}
\end{quote}
\end{quote}

\subsubsection{COPT\_GetDblParam, COPT\_GetDblParamDef/Min/Max}
\label{\detokenize{capiref:copt-getdblparam-copt-getdblparamdef-min-max}}\begin{quote}

\sphinxAtStartPar
\sphinxstylestrong{Synopsis}
\begin{quote}

\sphinxAtStartPar
\sphinxcode{\sphinxupquote{int COPT\_GetDblParam(copt\_prob *prob, const char *paramName, double *p\_dblParam)}}

\sphinxAtStartPar
\sphinxcode{\sphinxupquote{int COPT\_GetDblParamDef(copt\_prob *prob, const char *paramName, double *p\_dblParam)}}

\sphinxAtStartPar
\sphinxcode{\sphinxupquote{int COPT\_GetDblParamMin(copt\_prob *prob, const char *paramName, double *p\_dblParam)}}

\sphinxAtStartPar
\sphinxcode{\sphinxupquote{int COPT\_GetDblParamMax(copt\_prob *prob, const char *paramName, double *p\_dblParam)}}
\end{quote}

\sphinxAtStartPar
\sphinxstylestrong{Description}
\begin{quote}

\sphinxAtStartPar
Gets the
\begin{quote}

\sphinxAtStartPar
current

\sphinxAtStartPar
default

\sphinxAtStartPar
minimal

\sphinxAtStartPar
maximal
\end{quote}

\sphinxAtStartPar
value of a double parameter.
\end{quote}

\sphinxAtStartPar
\sphinxstylestrong{Arguments}
\begin{quote}

\sphinxAtStartPar
\sphinxcode{\sphinxupquote{prob}}
\begin{quote}

\sphinxAtStartPar
The COPT problem.
\end{quote}

\sphinxAtStartPar
\sphinxcode{\sphinxupquote{paramName}}
\begin{quote}

\sphinxAtStartPar
The name of the double parameter.
\end{quote}

\sphinxAtStartPar
\sphinxcode{\sphinxupquote{p\_dblParam}}
\begin{quote}

\sphinxAtStartPar
Pointer to the value of the double parameter.
\end{quote}
\end{quote}
\end{quote}

\subsubsection{COPT\_ResetParam}
\label{\detokenize{capiref:copt-resetparam}}\begin{quote}

\sphinxAtStartPar
\sphinxstylestrong{Synopsis}
\begin{quote}

\sphinxAtStartPar
\sphinxcode{\sphinxupquote{int COPT\_ResetParam(copt\_prob *prob)}}
\end{quote}

\sphinxAtStartPar
\sphinxstylestrong{Description}
\begin{quote}

\sphinxAtStartPar
Reset parameters to default settings.
\end{quote}

\sphinxAtStartPar
\sphinxstylestrong{Arguments}
\begin{quote}

\sphinxAtStartPar
\sphinxcode{\sphinxupquote{prob}}
\begin{quote}

\sphinxAtStartPar
The COPT problem.
\end{quote}
\end{quote}
\end{quote}

\subsubsection{COPT\_WriteParam}
\label{\detokenize{capiref:copt-writeparam}}\begin{quote}

\sphinxAtStartPar
\sphinxstylestrong{Synopsis}
\begin{quote}

\sphinxAtStartPar
\sphinxcode{\sphinxupquote{int COPT\_WriteParam(copt\_prob *prob, const char *parfilename)}}
\end{quote}

\sphinxAtStartPar
\sphinxstylestrong{Description}
\begin{quote}

\sphinxAtStartPar
Writes user defined parameters to a file. This API function
will write out all the parameters that are different from their
default values.
\end{quote}

\sphinxAtStartPar
\sphinxstylestrong{Arguments}
\begin{quote}

\sphinxAtStartPar
\sphinxcode{\sphinxupquote{prob}}
\begin{quote}

\sphinxAtStartPar
The COPT problem.
\end{quote}

\sphinxAtStartPar
\sphinxcode{\sphinxupquote{parfilename}}
\begin{quote}

\sphinxAtStartPar
The path to the parameter file.
\end{quote}
\end{quote}
\end{quote}

\subsubsection{COPT\_WriteParamStr}
\label{\detokenize{capiref:copt-writeparamstr}}\begin{quote}

\sphinxAtStartPar
\sphinxstylestrong{Synopsis}
\begin{quote}

\sphinxAtStartPar
\sphinxcode{\sphinxupquote{int COPT\_WriteParamStr(copt\_prob *prob, char *str, int nStrSize, int *pReqSize)}}
\end{quote}

\sphinxAtStartPar
\sphinxstylestrong{Description}
\begin{quote}

\sphinxAtStartPar
Writes the modified parameters to a string buffer.
\end{quote}

\sphinxAtStartPar
\sphinxstylestrong{Arguments}
\begin{quote}

\sphinxAtStartPar
\sphinxcode{\sphinxupquote{prob}}
\begin{quote}

\sphinxAtStartPar
The COPT problem.
\end{quote}

\sphinxAtStartPar
\sphinxcode{\sphinxupquote{str}}
\begin{quote}

\sphinxAtStartPar
String buffer of modified parameters.
\end{quote}

\sphinxAtStartPar
\sphinxcode{\sphinxupquote{nStrSize}}
\begin{quote}

\sphinxAtStartPar
The size of string buffer.
\end{quote}

\sphinxAtStartPar
\sphinxcode{\sphinxupquote{pReqSize}}
\begin{quote}

\sphinxAtStartPar
Minimum space requirement of string buffer for modified parameters.
\end{quote}
\end{quote}
\end{quote}

\subsubsection{COPT\_ReadParam}
\label{\detokenize{capiref:copt-readparam}}\begin{quote}

\sphinxAtStartPar
\sphinxstylestrong{Synopsis}
\begin{quote}

\sphinxAtStartPar
\sphinxcode{\sphinxupquote{int COPT\_ReadParam(copt\_prob *prob, const char *parfilename)}}
\end{quote}

\sphinxAtStartPar
\sphinxstylestrong{Description}
\begin{quote}

\sphinxAtStartPar
Reads and applies parameters settings as
defined in the parameter file.
\end{quote}

\sphinxAtStartPar
\sphinxstylestrong{Arguments}
\begin{quote}

\sphinxAtStartPar
\sphinxcode{\sphinxupquote{prob}}
\begin{quote}

\sphinxAtStartPar
The COPT problem.
\end{quote}

\sphinxAtStartPar
\sphinxcode{\sphinxupquote{parfilename}}
\begin{quote}

\sphinxAtStartPar
The path to the parameter file.
\end{quote}
\end{quote}
\end{quote}

\subsubsection{COPT\_ReadParamStr}
\label{\detokenize{capiref:copt-readparamstr}}\begin{quote}

\sphinxAtStartPar
\sphinxstylestrong{Synopsis}
\begin{quote}

\sphinxAtStartPar
\sphinxcode{\sphinxupquote{int COPT\_ReadParamStr(copt\_prob *prob, const char *strParam)}}
\end{quote}

\sphinxAtStartPar
\sphinxstylestrong{Description}
\begin{quote}

\sphinxAtStartPar
Read parameter settings from string buffer, and set parameters in COPT.
\end{quote}

\sphinxAtStartPar
\sphinxstylestrong{Arguments}
\begin{quote}

\sphinxAtStartPar
\sphinxcode{\sphinxupquote{prob}}
\begin{quote}

\sphinxAtStartPar
The COPT problem.
\end{quote}

\sphinxAtStartPar
\sphinxcode{\sphinxupquote{strParam}}
\begin{quote}

\sphinxAtStartPar
String buffer of parameter settings.
\end{quote}
\end{quote}
\end{quote}

\subsection{Accessing attributes}
\label{\detokenize{capiref:accessing-attributes}}\label{\detokenize{capiref:chapapi-getattr}}

\subsubsection{COPT\_GetIntAttr}
\label{\detokenize{capiref:copt-getintattr}}\begin{quote}

\sphinxAtStartPar
\sphinxstylestrong{Synopsis}
\begin{quote}

\sphinxAtStartPar
\sphinxcode{\sphinxupquote{int COPT\_GetIntAttr(copt\_prob *prob, const char *attrName, int *p\_intAttr)}}
\end{quote}

\sphinxAtStartPar
\sphinxstylestrong{Description}
\begin{quote}

\sphinxAtStartPar
Gets the value of an integer attribute.
\end{quote}

\sphinxAtStartPar
\sphinxstylestrong{Arguments}
\begin{quote}

\sphinxAtStartPar
\sphinxcode{\sphinxupquote{prob}}
\begin{quote}

\sphinxAtStartPar
The COPT problem.
\end{quote}

\sphinxAtStartPar
\sphinxcode{\sphinxupquote{attrName}}
\begin{quote}

\sphinxAtStartPar
The name of the integer attribute.
\end{quote}

\sphinxAtStartPar
\sphinxcode{\sphinxupquote{p\_intAttr}}
\begin{quote}

\sphinxAtStartPar
Pointer to the value of the integer attribute.
\end{quote}
\end{quote}
\end{quote}

\subsubsection{COPT\_GetDblAttr}
\label{\detokenize{capiref:copt-getdblattr}}\begin{quote}

\sphinxAtStartPar
\sphinxstylestrong{Synopsis}
\begin{quote}

\sphinxAtStartPar
\sphinxcode{\sphinxupquote{int COPT\_GetDblAttr(copt\_prob *prob, const char *attrName, int *p\_dblAttr)}}
\end{quote}

\sphinxAtStartPar
\sphinxstylestrong{Description}
\begin{quote}

\sphinxAtStartPar
Gets the value of an double attribute.
\end{quote}

\sphinxAtStartPar
\sphinxstylestrong{Arguments}
\begin{quote}

\sphinxAtStartPar
\sphinxcode{\sphinxupquote{prob}}
\begin{quote}

\sphinxAtStartPar
The COPT problem.
\end{quote}

\sphinxAtStartPar
\sphinxcode{\sphinxupquote{attrName}}
\begin{quote}

\sphinxAtStartPar
The name of the double attribute.
\end{quote}

\sphinxAtStartPar
\sphinxcode{\sphinxupquote{p\_dblAttr}}
\begin{quote}

\sphinxAtStartPar
Pointer to the value of the double attribute.
\end{quote}
\end{quote}
\end{quote}

\subsection{Logging utilities}
\label{\detokenize{capiref:logging-utilities}}

\subsubsection{COPT\_SetLogFile}
\label{\detokenize{capiref:copt-setlogfile}}\begin{quote}

\sphinxAtStartPar
\sphinxstylestrong{Synopsis}
\begin{quote}

\sphinxAtStartPar
\sphinxcode{\sphinxupquote{int COPT\_SetLogFile(copt\_prob *prob, char *logfilename)}}
\end{quote}

\sphinxAtStartPar
\sphinxstylestrong{Description}
\begin{quote}

\sphinxAtStartPar
Set log file for the problem.
\end{quote}

\sphinxAtStartPar
\sphinxstylestrong{Arguments}
\begin{quote}

\sphinxAtStartPar
\sphinxcode{\sphinxupquote{prob}}
\begin{quote}

\sphinxAtStartPar
The COPT problem.
\end{quote}

\sphinxAtStartPar
\sphinxcode{\sphinxupquote{logfilename}}
\begin{quote}

\sphinxAtStartPar
The path to the log file.
\end{quote}
\end{quote}
\end{quote}

\subsubsection{COPT\_SetLogCallback}
\label{\detokenize{capiref:copt-setlogcallback}}\begin{quote}

\sphinxAtStartPar
\sphinxstylestrong{Synopsis}
\begin{quote}

\sphinxAtStartPar
\sphinxcode{\sphinxupquote{int COPT\_SetLogCallback(copt\_prob *prob, void (*logcb)(char *msg, void *userdata), void *userdata)}}
\end{quote}

\sphinxAtStartPar
\sphinxstylestrong{Description}
\begin{quote}

\sphinxAtStartPar
Set message callback for the problem.
\end{quote}

\sphinxAtStartPar
\sphinxstylestrong{Arguments}
\begin{quote}

\sphinxAtStartPar
\sphinxcode{\sphinxupquote{prob}}
\begin{quote}

\sphinxAtStartPar
The COPT problem.
\end{quote}

\sphinxAtStartPar
\sphinxcode{\sphinxupquote{logcb}}
\begin{quote}

\sphinxAtStartPar
Callback function for message.
\end{quote}

\sphinxAtStartPar
\sphinxcode{\sphinxupquote{userdata}}
\begin{quote}

\sphinxAtStartPar
User defined data. The data will be passed to the solver without
modification.
\end{quote}
\end{quote}
\end{quote}

\subsection{MIP start utilities}
\label{\detokenize{capiref:mip-start-utilities}}\label{\detokenize{capiref:chapapi-mipstart}}

\subsubsection{COPT\_AddMipStart}
\label{\detokenize{capiref:copt-addmipstart}}\begin{quote}

\sphinxAtStartPar
\sphinxstylestrong{Synopsis}
\begin{quote}

\sphinxAtStartPar
\sphinxcode{\sphinxupquote{int COPT\_AddMipStart(copt\_prob *prob, int num, const int *list, double *colVal)}}
\end{quote}

\sphinxAtStartPar
\sphinxstylestrong{Description}
\begin{quote}

\sphinxAtStartPar
Add MIP start information for the problem. If \sphinxcode{\sphinxupquote{list}} is \sphinxcode{\sphinxupquote{NULL}},
then information of the first \sphinxcode{\sphinxupquote{num}} columns will be added.

\sphinxAtStartPar
One MIP start information will be added for each call to this function.
\end{quote}

\sphinxAtStartPar
\sphinxstylestrong{Arguments}
\begin{quote}

\sphinxAtStartPar
\sphinxcode{\sphinxupquote{prob}}
\begin{quote}

\sphinxAtStartPar
The COPT problem.
\end{quote}

\sphinxAtStartPar
\sphinxcode{\sphinxupquote{num}}
\begin{quote}

\sphinxAtStartPar
Number of variables (columns).
\end{quote}

\sphinxAtStartPar
\sphinxcode{\sphinxupquote{list}}
\begin{quote}

\sphinxAtStartPar
Index of variables (columns). Can be \sphinxcode{\sphinxupquote{NULL}}.
\end{quote}

\sphinxAtStartPar
\sphinxcode{\sphinxupquote{colVal}}
\begin{quote}

\sphinxAtStartPar
MIP start information.
\end{quote}
\end{quote}
\end{quote}

\subsubsection{COPT\_ReadMst}
\label{\detokenize{capiref:copt-readmst}}\begin{quote}

\sphinxAtStartPar
\sphinxstylestrong{Synopsis}
\begin{quote}

\sphinxAtStartPar
\sphinxcode{\sphinxupquote{int COPT\_ReadMst(copt\_prob *prob, const char *mstfilename)}}
\end{quote}

\sphinxAtStartPar
\sphinxstylestrong{Description}
\begin{quote}

\sphinxAtStartPar
Read MIP start information from file, and used as initial solution for
the problem.
\end{quote}

\sphinxAtStartPar
\sphinxstylestrong{Arguments}
\begin{quote}

\sphinxAtStartPar
\sphinxcode{\sphinxupquote{prob}}
\begin{quote}

\sphinxAtStartPar
The COPT problem.
\end{quote}

\sphinxAtStartPar
\sphinxcode{\sphinxupquote{mstfilename}}
\begin{quote}

\sphinxAtStartPar
The path to the MIP start file.
\end{quote}
\end{quote}
\end{quote}

\subsubsection{COPT\_WriteMst}
\label{\detokenize{capiref:copt-writemst}}\begin{quote}

\sphinxAtStartPar
\sphinxstylestrong{Synopsis}
\begin{quote}

\sphinxAtStartPar
\sphinxcode{\sphinxupquote{int COPT\_WriteMst(copt\_prob *prob, const char *mstfilename)}}
\end{quote}

\sphinxAtStartPar
\sphinxstylestrong{Description}
\begin{quote}

\sphinxAtStartPar
Write solution or existed MIP start information in problem to file.
\end{quote}

\sphinxAtStartPar
\sphinxstylestrong{Arguments}
\begin{quote}

\sphinxAtStartPar
\sphinxcode{\sphinxupquote{prob}}
\begin{quote}

\sphinxAtStartPar
The COPT problem.
\end{quote}

\sphinxAtStartPar
\sphinxcode{\sphinxupquote{mstfilename}}
\begin{quote}

\sphinxAtStartPar
The path to the MIP start file.
\end{quote}
\end{quote}
\end{quote}

\subsection{Nonlinear start point utilities}
\label{\detokenize{capiref:nonlinear-start-point-utilities}}\label{\detokenize{capiref:chapapi-nlprimstart}}

\subsubsection{COPT\_SetNLPrimalStart}
\label{\detokenize{capiref:copt-setnlprimalstart}}\begin{quote}

\sphinxAtStartPar
\sphinxstylestrong{Synopsis}
\begin{quote}

\sphinxAtStartPar
\sphinxcode{\sphinxupquote{int COPT\_SetNLPrimalStart(copt\_prob *prob, int num, const int *list, double *colVal)}}
\end{quote}

\sphinxAtStartPar
\sphinxstylestrong{Description}
\begin{quote}

\sphinxAtStartPar
Set nonlinear primal start information for the problem. If \sphinxcode{\sphinxupquote{list}} is \sphinxcode{\sphinxupquote{NULL}},
then information of the first \sphinxcode{\sphinxupquote{num}} columns will be set.
\end{quote}

\sphinxAtStartPar
\sphinxstylestrong{Arguments}
\begin{quote}

\sphinxAtStartPar
\sphinxcode{\sphinxupquote{prob}}
\begin{quote}

\sphinxAtStartPar
The COPT problem.
\end{quote}

\sphinxAtStartPar
\sphinxcode{\sphinxupquote{num}}
\begin{quote}

\sphinxAtStartPar
Number of variables (columns).
\end{quote}

\sphinxAtStartPar
\sphinxcode{\sphinxupquote{list}}
\begin{quote}

\sphinxAtStartPar
Index of variables (columns). Can be \sphinxcode{\sphinxupquote{NULL}}.
\end{quote}

\sphinxAtStartPar
\sphinxcode{\sphinxupquote{colVal}}
\begin{quote}

\sphinxAtStartPar
Nonlinear primal start information.
\end{quote}
\end{quote}
\end{quote}

\subsection{IIS utilities}
\label{\detokenize{capiref:iis-utilities}}\label{\detokenize{capiref:chapapi-iiscompute}}

\subsubsection{COPT\_ComputeIIS}
\label{\detokenize{capiref:copt-computeiis}}\begin{quote}

\sphinxAtStartPar
\sphinxstylestrong{Synopsis}
\begin{quote}

\sphinxAtStartPar
\sphinxcode{\sphinxupquote{int COPT\_ComputeIIS(copt\_prob *prob)}}
\end{quote}

\sphinxAtStartPar
\sphinxstylestrong{Description}
\begin{quote}

\sphinxAtStartPar
Compute IIS (Irreducible Inconsistent Subsystem) for infeasible problem.
\end{quote}

\sphinxAtStartPar
\sphinxstylestrong{Arguments}
\begin{quote}

\sphinxAtStartPar
\sphinxcode{\sphinxupquote{prob}}
\begin{quote}

\sphinxAtStartPar
The COPT problem.
\end{quote}
\end{quote}
\end{quote}

\subsubsection{COPT\_GetColLowerIIS}
\label{\detokenize{capiref:copt-getcolloweriis}}\begin{quote}

\sphinxAtStartPar
\sphinxstylestrong{Synopsis}
\begin{quote}

\sphinxAtStartPar
\sphinxcode{\sphinxupquote{int COPT\_GetColLowerIIS(copt\_prob *prob, int num, const int *list, int *colLowerIIS)}}
\end{quote}

\sphinxAtStartPar
\sphinxstylestrong{Description}
\begin{quote}

\sphinxAtStartPar
Get IIS status of lower bounds of columns. If \sphinxcode{\sphinxupquote{list}} is \sphinxcode{\sphinxupquote{NULL}}, then IIS status
of the first \sphinxcode{\sphinxupquote{num}} columns will be returned.
\end{quote}

\sphinxAtStartPar
\sphinxstylestrong{Arguments}
\begin{quote}

\sphinxAtStartPar
\sphinxcode{\sphinxupquote{prob}}
\begin{quote}

\sphinxAtStartPar
The COPT problem.
\end{quote}

\sphinxAtStartPar
\sphinxcode{\sphinxupquote{num}}
\begin{quote}

\sphinxAtStartPar
Number of columns.
\end{quote}

\sphinxAtStartPar
\sphinxcode{\sphinxupquote{list}}
\begin{quote}

\sphinxAtStartPar
Index of columns. Can be \sphinxcode{\sphinxupquote{NULL}}.
\end{quote}

\sphinxAtStartPar
\sphinxcode{\sphinxupquote{colLowerIIS}}
\begin{quote}

\sphinxAtStartPar
IIS status of lower bounds of columns.
\end{quote}
\end{quote}
\end{quote}

\subsubsection{COPT\_GetColUpperIIS}
\label{\detokenize{capiref:copt-getcolupperiis}}\begin{quote}

\sphinxAtStartPar
\sphinxstylestrong{Synopsis}
\begin{quote}

\sphinxAtStartPar
\sphinxcode{\sphinxupquote{int COPT\_GetColUpperIIS(copt\_prob *prob, int num, const int *list, int *colUpperIIS)}}
\end{quote}

\sphinxAtStartPar
\sphinxstylestrong{Description}
\begin{quote}

\sphinxAtStartPar
Get IIS status of upper bounds of columns. If \sphinxcode{\sphinxupquote{list}} is \sphinxcode{\sphinxupquote{NULL}}, then IIS status
of the first \sphinxcode{\sphinxupquote{num}} columns will be returned.
\end{quote}

\sphinxAtStartPar
\sphinxstylestrong{Arguments}
\begin{quote}

\sphinxAtStartPar
\sphinxcode{\sphinxupquote{prob}}
\begin{quote}

\sphinxAtStartPar
The COPT problem.
\end{quote}

\sphinxAtStartPar
\sphinxcode{\sphinxupquote{num}}
\begin{quote}

\sphinxAtStartPar
Number of columns.
\end{quote}

\sphinxAtStartPar
\sphinxcode{\sphinxupquote{list}}
\begin{quote}

\sphinxAtStartPar
Index of columns. Can be \sphinxcode{\sphinxupquote{NULL}}.
\end{quote}

\sphinxAtStartPar
\sphinxcode{\sphinxupquote{colUpperIIS}}
\begin{quote}

\sphinxAtStartPar
IIS status of upper bounds of columns.
\end{quote}
\end{quote}
\end{quote}

\subsubsection{COPT\_GetRowLowerIIS}
\label{\detokenize{capiref:copt-getrowloweriis}}\begin{quote}

\sphinxAtStartPar
\sphinxstylestrong{Synopsis}
\begin{quote}

\sphinxAtStartPar
\sphinxcode{\sphinxupquote{int COPT\_GetRowLowerIIS(copt\_prob *prob, int num, const int *list, int *rowLowerIIS)}}
\end{quote}

\sphinxAtStartPar
\sphinxstylestrong{Description}
\begin{quote}

\sphinxAtStartPar
Get IIS status of lower bounds of rows. If \sphinxcode{\sphinxupquote{list}} is \sphinxcode{\sphinxupquote{NULL}}, then IIS status
of the first \sphinxcode{\sphinxupquote{num}} rows will be returned.
\end{quote}

\sphinxAtStartPar
\sphinxstylestrong{Arguments}
\begin{quote}

\sphinxAtStartPar
\sphinxcode{\sphinxupquote{prob}}
\begin{quote}

\sphinxAtStartPar
The COPT problem.
\end{quote}

\sphinxAtStartPar
\sphinxcode{\sphinxupquote{num}}
\begin{quote}

\sphinxAtStartPar
Number of rows.
\end{quote}

\sphinxAtStartPar
\sphinxcode{\sphinxupquote{list}}
\begin{quote}

\sphinxAtStartPar
Index of rows. Can be \sphinxcode{\sphinxupquote{NULL}}.
\end{quote}

\sphinxAtStartPar
\sphinxcode{\sphinxupquote{rowLowerIIS}}
\begin{quote}

\sphinxAtStartPar
IIS status of lower bounds of rows.
\end{quote}
\end{quote}
\end{quote}

\subsubsection{COPT\_GetRowUpperIIS}
\label{\detokenize{capiref:copt-getrowupperiis}}\begin{quote}

\sphinxAtStartPar
\sphinxstylestrong{Synopsis}
\begin{quote}

\sphinxAtStartPar
\sphinxcode{\sphinxupquote{int COPT\_GetRowUpperIIS(copt\_prob *prob, int num, const int *list, int *rowUpperIIS)}}
\end{quote}

\sphinxAtStartPar
\sphinxstylestrong{Description}
\begin{quote}

\sphinxAtStartPar
Get IIS status of upper bounds of rows. If \sphinxcode{\sphinxupquote{list}} is \sphinxcode{\sphinxupquote{NULL}}, then IIS status
of the first \sphinxcode{\sphinxupquote{num}} rows will be returned.
\end{quote}

\sphinxAtStartPar
\sphinxstylestrong{Arguments}
\begin{quote}

\sphinxAtStartPar
\sphinxcode{\sphinxupquote{prob}}
\begin{quote}

\sphinxAtStartPar
The COPT problem.
\end{quote}

\sphinxAtStartPar
\sphinxcode{\sphinxupquote{num}}
\begin{quote}

\sphinxAtStartPar
Number of rows.
\end{quote}

\sphinxAtStartPar
\sphinxcode{\sphinxupquote{list}}
\begin{quote}

\sphinxAtStartPar
Index of rows. Can be \sphinxcode{\sphinxupquote{NULL}}.
\end{quote}

\sphinxAtStartPar
\sphinxcode{\sphinxupquote{rowUpperIIS}}
\begin{quote}

\sphinxAtStartPar
IIS status of upper bounds of rows.
\end{quote}
\end{quote}
\end{quote}

\subsubsection{COPT\_GetSOSIIS}
\label{\detokenize{capiref:copt-getsosiis}}\begin{quote}

\sphinxAtStartPar
\sphinxstylestrong{Synopsis}
\begin{quote}

\sphinxAtStartPar
\sphinxcode{\sphinxupquote{int COPT\_GetSOSIIS(copt\_prob *prob, int num, const int *list, int *sosIIS)}}
\end{quote}

\sphinxAtStartPar
\sphinxstylestrong{Description}
\begin{quote}

\sphinxAtStartPar
Get IIS status of SOS constraints. If \sphinxcode{\sphinxupquote{list}} is \sphinxcode{\sphinxupquote{NULL}}, then IIS status
of the first \sphinxcode{\sphinxupquote{num}} SOS constraints will be returned.
\end{quote}

\sphinxAtStartPar
\sphinxstylestrong{Arguments}
\begin{quote}

\sphinxAtStartPar
\sphinxcode{\sphinxupquote{prob}}
\begin{quote}

\sphinxAtStartPar
The COPT problem.
\end{quote}

\sphinxAtStartPar
\sphinxcode{\sphinxupquote{num}}
\begin{quote}

\sphinxAtStartPar
Number of SOS constraints.
\end{quote}

\sphinxAtStartPar
\sphinxcode{\sphinxupquote{list}}
\begin{quote}

\sphinxAtStartPar
Index of SOS constraints. Can be \sphinxcode{\sphinxupquote{NULL}}.
\end{quote}

\sphinxAtStartPar
\sphinxcode{\sphinxupquote{sosIIS}}
\begin{quote}

\sphinxAtStartPar
IIS status of SOS constraints.
\end{quote}
\end{quote}
\end{quote}

\subsubsection{COPT\_GetIndicatorIIS}
\label{\detokenize{capiref:copt-getindicatoriis}}\begin{quote}

\sphinxAtStartPar
\sphinxstylestrong{Synopsis}
\begin{quote}

\sphinxAtStartPar
\sphinxcode{\sphinxupquote{int COPT\_GetIndicatorIIS(copt\_prob *prob, int num, const int *list, int *indicatorIIS)}}
\end{quote}

\sphinxAtStartPar
\sphinxstylestrong{Description}
\begin{quote}

\sphinxAtStartPar
Get IIS status of indicator constraints. If \sphinxcode{\sphinxupquote{list}} is \sphinxcode{\sphinxupquote{NULL}}, then IIS status
of the first \sphinxcode{\sphinxupquote{num}} indicator constraints will be returned.
\end{quote}

\sphinxAtStartPar
\sphinxstylestrong{Arguments}
\begin{quote}

\sphinxAtStartPar
\sphinxcode{\sphinxupquote{prob}}
\begin{quote}

\sphinxAtStartPar
The COPT problem.
\end{quote}

\sphinxAtStartPar
\sphinxcode{\sphinxupquote{num}}
\begin{quote}

\sphinxAtStartPar
Number of indicator constraints.
\end{quote}

\sphinxAtStartPar
\sphinxcode{\sphinxupquote{list}}
\begin{quote}

\sphinxAtStartPar
Index of indicator constraints. Can be \sphinxcode{\sphinxupquote{NULL}}.
\end{quote}

\sphinxAtStartPar
\sphinxcode{\sphinxupquote{indicatorIIS}}
\begin{quote}

\sphinxAtStartPar
IIS status of indicator constraints.
\end{quote}
\end{quote}
\end{quote}

\subsection{Feasibility relaxation utilities}
\label{\detokenize{capiref:feasibility-relaxation-utilities}}\label{\detokenize{capiref:chapapi-feasrelax}}

\subsubsection{COPT\_FeasRelax}
\label{\detokenize{capiref:copt-feasrelax}}\begin{quote}

\sphinxAtStartPar
\sphinxstylestrong{Synopsis}
\begin{quote}

\sphinxAtStartPar
\sphinxcode{\sphinxupquote{int COPT\_FeasRelax(copt\_prob *prob, double *colLowPen, double *colUppPen, double *rowBndPen, double *rowUppPen)}}
\end{quote}

\sphinxAtStartPar
\sphinxstylestrong{Description}
\begin{quote}

\sphinxAtStartPar
Compute feasibility relaxation to infeasible problem.
\end{quote}

\sphinxAtStartPar
\sphinxstylestrong{Arguments}
\begin{quote}

\sphinxAtStartPar
\sphinxcode{\sphinxupquote{prob}}
\begin{quote}

\sphinxAtStartPar
The COPT problem.
\end{quote}

\sphinxAtStartPar
\sphinxcode{\sphinxupquote{colLowPen}}
\begin{quote}

\sphinxAtStartPar
Penalty for lower bounds of columns. If \sphinxcode{\sphinxupquote{NULL}}, then no relaxation for lower bounds of columns are allowed;
If penalty in \sphinxcode{\sphinxupquote{colLowPen}} is \sphinxcode{\sphinxupquote{COPT\_INFINITY}}, then no relaxation is allowed for corresponding bound.
\end{quote}

\sphinxAtStartPar
\sphinxcode{\sphinxupquote{colUppPen}}
\begin{quote}

\sphinxAtStartPar
Penalty for upper bounds of columns. If \sphinxcode{\sphinxupquote{NULL}}, then no relaxation for upper bounds of columns are allowed;
If penalty in \sphinxcode{\sphinxupquote{colUppen}} is \sphinxcode{\sphinxupquote{COPT\_INFINITY}}, then no relaxation is allowed for corresponding bound.
\end{quote}

\sphinxAtStartPar
\sphinxcode{\sphinxupquote{rowBndPen}}
\begin{quote}

\sphinxAtStartPar
Penalty for bounds of rows. If \sphinxcode{\sphinxupquote{NULL}}, then no relaxation for bounds of rows are allowed;
If penalty in \sphinxcode{\sphinxupquote{rowBndPen}} is \sphinxcode{\sphinxupquote{COPT\_INFINITY}}, then no relaxation is allowed for corresponding bound.
\end{quote}

\sphinxAtStartPar
\sphinxcode{\sphinxupquote{rowUppPen}}
\begin{quote}

\sphinxAtStartPar
Penalty for upper bounds of rows. If problem has two\sphinxhyphen{}sided rows, and \sphinxcode{\sphinxupquote{rowUppPen}} is not \sphinxcode{\sphinxupquote{NULL}}, then
\sphinxcode{\sphinxupquote{rowUppPen}} is penalty for upper bounds of rows; If penalty in \sphinxcode{\sphinxupquote{rowUppPen}} is \sphinxcode{\sphinxupquote{COPT\_INFINITY}},
then no relaxation is allowed for corresponding bound.

\sphinxAtStartPar
\sphinxstylestrong{Note:} Normally, just set \sphinxcode{\sphinxupquote{rowUppPen}} to \sphinxcode{\sphinxupquote{NULL}}.
\end{quote}
\end{quote}
\end{quote}

\subsubsection{COPT\_WriteRelax}
\label{\detokenize{capiref:copt-writerelax}}\begin{quote}

\sphinxAtStartPar
\sphinxstylestrong{Synopsis}
\begin{quote}

\sphinxAtStartPar
\sphinxcode{\sphinxupquote{int COPT\_WriteRelax(copt\_prob *prob, const char *relaxfilename)}}
\end{quote}

\sphinxAtStartPar
\sphinxstylestrong{Description}
\begin{quote}

\sphinxAtStartPar
Write feasrelax problem to file.
\end{quote}

\sphinxAtStartPar
\sphinxstylestrong{Arguments}
\begin{quote}

\sphinxAtStartPar
\sphinxcode{\sphinxupquote{prob}}
\begin{quote}

\sphinxAtStartPar
The COPT problem.
\end{quote}

\sphinxAtStartPar
\sphinxcode{\sphinxupquote{relaxfilename}}
\begin{quote}

\sphinxAtStartPar
Name of feasrelax problem file.
\end{quote}
\end{quote}
\end{quote}

\subsection{Parameter tuning utilities}
\label{\detokenize{capiref:parameter-tuning-utilities}}\label{\detokenize{capiref:chapapi-tune}}

\subsubsection{COPT\_Tune}
\label{\detokenize{capiref:copt-tune}}\begin{quote}

\sphinxAtStartPar
\sphinxstylestrong{Synopsis}
\begin{quote}

\sphinxAtStartPar
\sphinxcode{\sphinxupquote{int COPT\_Tune(copt\_prob *prob)}}
\end{quote}

\sphinxAtStartPar
\sphinxstylestrong{Description}
\begin{quote}

\sphinxAtStartPar
Parameter tuning of the model.
\end{quote}

\sphinxAtStartPar
\sphinxstylestrong{Arguments}
\begin{quote}

\sphinxAtStartPar
\sphinxcode{\sphinxupquote{prob}}
\begin{quote}

\sphinxAtStartPar
COPT model.
\end{quote}
\end{quote}
\end{quote}

\subsubsection{COPT\_LoadTuneParam}
\label{\detokenize{capiref:copt-loadtuneparam}}\begin{quote}

\sphinxAtStartPar
\sphinxstylestrong{Synopsis}
\begin{quote}

\sphinxAtStartPar
\sphinxcode{\sphinxupquote{int COPT\_LoadTuneParam(copt\_prob *prob, int idx)}}
\end{quote}

\sphinxAtStartPar
\sphinxstylestrong{Description}
\begin{quote}

\sphinxAtStartPar
Load the parameter tuning results of the specified number into the model.
\end{quote}

\sphinxAtStartPar
\sphinxstylestrong{Arguments}
\begin{quote}

\sphinxAtStartPar
\sphinxcode{\sphinxupquote{prob}}
\begin{quote}

\sphinxAtStartPar
COPT model.
\end{quote}

\sphinxAtStartPar
\sphinxcode{\sphinxupquote{idx}}
\begin{quote}

\sphinxAtStartPar
The number of the parameter tuning result.
\end{quote}
\end{quote}
\end{quote}

\subsubsection{COPT\_ReadTune}
\label{\detokenize{capiref:copt-readtune}}\begin{quote}

\sphinxAtStartPar
\sphinxstylestrong{Synopsis}
\begin{quote}

\sphinxAtStartPar
\sphinxcode{\sphinxupquote{int COPT\_ReadTune(copt\_prob *prob, const char *tunefilename)}}
\end{quote}

\sphinxAtStartPar
\sphinxstylestrong{Description}
\begin{quote}

\sphinxAtStartPar
Read the parameter combination to be tuned from the tuning file to the model.
\end{quote}

\sphinxAtStartPar
\sphinxstylestrong{Arguments}
\begin{quote}

\sphinxAtStartPar
\sphinxcode{\sphinxupquote{prob}}
\begin{quote}

\sphinxAtStartPar
COPT model.
\end{quote}

\sphinxAtStartPar
\sphinxcode{\sphinxupquote{tunefilename}}
\begin{quote}

\sphinxAtStartPar
Tuning file names.
\end{quote}
\end{quote}
\end{quote}

\subsubsection{COPT\_WriteTuneParam}
\label{\detokenize{capiref:copt-writetuneparam}}\begin{quote}

\sphinxAtStartPar
\sphinxstylestrong{Synopsis}
\begin{quote}

\sphinxAtStartPar
\sphinxcode{\sphinxupquote{int COPT\_WriteTuneParam(copt\_prob *prob, int idx, const char *parfilename)}}
\end{quote}

\sphinxAtStartPar
\sphinxstylestrong{Description}
\begin{quote}

\sphinxAtStartPar
Output the parameter tuning result of the specified number to the parameter file.
\end{quote}

\sphinxAtStartPar
\sphinxstylestrong{Arguments}
\begin{quote}

\sphinxAtStartPar
\sphinxcode{\sphinxupquote{prob}}
\begin{quote}

\sphinxAtStartPar
COPT model.
\end{quote}

\sphinxAtStartPar
\sphinxcode{\sphinxupquote{idx}}
\begin{quote}

\sphinxAtStartPar
The number of the parameter tuning result.
\end{quote}

\sphinxAtStartPar
\sphinxcode{\sphinxupquote{parfilename}}
\begin{quote}

\sphinxAtStartPar
parameter file name.
\end{quote}
\end{quote}
\end{quote}

\subsubsection{COPT\_ReadOrd}
\label{\detokenize{capiref:copt-readord}}\begin{quote}

\sphinxAtStartPar
\sphinxstylestrong{Synopsis}
\begin{quote}

\sphinxAtStartPar
\sphinxcode{\sphinxupquote{int COPT\_ReadOrd(copt\_prob *prob, const char *ordfilename)}}
\end{quote}

\sphinxAtStartPar
\sphinxstylestrong{Description}
\begin{quote}

\sphinxAtStartPar
Read a branching order file (ORD\sphinxhyphen{}format) into the current model
to reuse a predefined branching strategy.
\end{quote}

\sphinxAtStartPar
\sphinxstylestrong{Arguments}
\begin{quote}

\sphinxAtStartPar
\sphinxcode{\sphinxupquote{prob}}
\begin{quote}

\sphinxAtStartPar
The COPT problem.
\end{quote}

\sphinxAtStartPar
\sphinxcode{\sphinxupquote{ordfilename}}
\begin{quote}

\sphinxAtStartPar
Name of the branching order file.
\end{quote}
\end{quote}
\end{quote}

\subsubsection{COPT\_WriteOrd}
\label{\detokenize{capiref:copt-writeord}}\begin{quote}

\sphinxAtStartPar
\sphinxstylestrong{Synopsis}
\begin{quote}

\sphinxAtStartPar
\sphinxcode{\sphinxupquote{int COPT\_WriteOrd(copt\_prob *prob, const char *ordfilename)}}
\end{quote}

\sphinxAtStartPar
\sphinxstylestrong{Description}
\begin{quote}

\sphinxAtStartPar
Write the branching order information of the current model into an ORD\sphinxhyphen{}format file,
so that the branching strategy can be saved and reused.
\end{quote}

\sphinxAtStartPar
\sphinxstylestrong{Arguments}
\begin{quote}

\sphinxAtStartPar
\sphinxcode{\sphinxupquote{prob}}
\begin{quote}

\sphinxAtStartPar
The COPT problem.
\end{quote}

\sphinxAtStartPar
\sphinxcode{\sphinxupquote{ordfilename}}
\begin{quote}

\sphinxAtStartPar
Path of the ORD\sphinxhyphen{}format file.
\end{quote}
\end{quote}
\end{quote}

\subsection{Callback utilities}
\label{\detokenize{capiref:callback-utilities}}\label{\detokenize{capiref:chapapi-cbc}}
\sphinxAtStartPar
Certain callback utilities methods can only be called in certain contexts, which are listed below:

\begin{savenotes}\sphinxattablestart
\sphinxthistablewithglobalstyle
\centering
\sphinxcapstartof{table}
\sphinxthecaptionisattop
\sphinxcaption{Callback utilities}\label{\detokenize{capiref:id2}}
\sphinxaftertopcaption
\begin{tabular}[t]{|\X{25}{65}|\X{40}{65}|}
\sphinxtoprule
\sphinxstyletheadfamily 
\sphinxAtStartPar
Context
&\sphinxstyletheadfamily 
\sphinxAtStartPar
Methods
\\
\sphinxmidrule
\sphinxtableatstartofbodyhook
\sphinxAtStartPar
COPT\_CBCONTEXT\_INCUMBENT
&
\sphinxAtStartPar
{\hyperref[\detokenize{capiref:chapapi-cbc-addcallbacksolution}]{\sphinxcrossref{\DUrole{std,std-ref}{COPT\_AddCallbackSolution}}}}, {\hyperref[\detokenize{capiref:chapapi-cbc-getcallbackinfo}]{\sphinxcrossref{\DUrole{std,std-ref}{COPT\_GetCallbackInfo}}}}
\\
\sphinxhline
\sphinxAtStartPar
COPT\_CBCONTEXT\_MIPSOL
&
\sphinxAtStartPar
{\hyperref[\detokenize{capiref:chapapi-cbc-addcallbacklazyconstr}]{\sphinxcrossref{\DUrole{std,std-ref}{COPT\_AddCallbackLazyConstr}}}}, {\hyperref[\detokenize{capiref:chapapi-cbc-addcallbacklazyconstrs}]{\sphinxcrossref{\DUrole{std,std-ref}{COPT\_AddCallbackLazyConstrs}}}}, {\hyperref[\detokenize{capiref:chapapi-cbc-addcallbacksolution}]{\sphinxcrossref{\DUrole{std,std-ref}{COPT\_AddCallbackSolution}}}}, {\hyperref[\detokenize{capiref:chapapi-cbc-getcallbackinfo}]{\sphinxcrossref{\DUrole{std,std-ref}{COPT\_GetCallbackInfo}}}}
\\
\sphinxhline
\sphinxAtStartPar
COPT\_CBCONTEXT\_MIPRELAX
&
\sphinxAtStartPar
{\hyperref[\detokenize{capiref:chapapi-cbc-addcallbackusercut}]{\sphinxcrossref{\DUrole{std,std-ref}{COPT\_AddCallbackUserCut}}}}, {\hyperref[\detokenize{capiref:chapapi-cbc-addcallbackusercuts}]{\sphinxcrossref{\DUrole{std,std-ref}{COPT\_AddCallbackUserCuts}}}}, {\hyperref[\detokenize{capiref:chapapi-cbc-addcallbacksolution}]{\sphinxcrossref{\DUrole{std,std-ref}{COPT\_AddCallbackSolution}}}}, {\hyperref[\detokenize{capiref:chapapi-cbc-getcallbackinfo}]{\sphinxcrossref{\DUrole{std,std-ref}{COPT\_GetCallbackInfo}}}}
\\
\sphinxhline
\sphinxAtStartPar
COPT\_CBCONTEXT\_MIPNODE
&
\sphinxAtStartPar
{\hyperref[\detokenize{capiref:chapapi-cbc-addcallbacksolution}]{\sphinxcrossref{\DUrole{std,std-ref}{COPT\_AddCallbackSolution}}}}, {\hyperref[\detokenize{capiref:chapapi-cbc-getcallbackinfo}]{\sphinxcrossref{\DUrole{std,std-ref}{COPT\_GetCallbackInfo}}}}
\\
\sphinxbottomrule
\end{tabular}
\sphinxtableafterendhook\par
\sphinxattableend\end{savenotes}

\subsubsection{COPT\_AddLazyConstr}
\label{\detokenize{capiref:copt-addlazyconstr}}\label{\detokenize{capiref:chapapi-cbc-addlazyconstr}}\begin{quote}

\sphinxAtStartPar
\sphinxstylestrong{Synopsis}
\begin{quote}

\sphinxAtStartPar
\sphinxcode{\sphinxupquote{int COPT\_AddLazyConstr(}}
\begin{quote}

\sphinxAtStartPar
\sphinxcode{\sphinxupquote{copt\_prob *prob,}}

\sphinxAtStartPar
\sphinxcode{\sphinxupquote{int nRowMatCnt,}}

\sphinxAtStartPar
\sphinxcode{\sphinxupquote{const int *rowMatIdx,}}

\sphinxAtStartPar
\sphinxcode{\sphinxupquote{const double *rowMatElem,}}

\sphinxAtStartPar
\sphinxcode{\sphinxupquote{char cRowSense,}}

\sphinxAtStartPar
\sphinxcode{\sphinxupquote{double dRowBound,}}

\sphinxAtStartPar
\sphinxcode{\sphinxupquote{double dRowUpper,}}

\sphinxAtStartPar
\sphinxcode{\sphinxupquote{const char *rowName)}}
\end{quote}
\end{quote}

\sphinxAtStartPar
\sphinxstylestrong{Description}
\begin{quote}

\sphinxAtStartPar
Add a lazy constraint to the MIP model.
\end{quote}

\sphinxAtStartPar
\sphinxstylestrong{Arguments}
\begin{quote}

\sphinxAtStartPar
\sphinxcode{\sphinxupquote{prob}}
\begin{quote}

\sphinxAtStartPar
The COPT problem.
\end{quote}

\sphinxAtStartPar
\sphinxcode{\sphinxupquote{nRowMatCnt}}
\begin{quote}

\sphinxAtStartPar
Number of non\sphinxhyphen{}zero elements in the lazy constraint.
\end{quote}

\sphinxAtStartPar
\sphinxcode{\sphinxupquote{rowMatIdx}}
\begin{quote}

\sphinxAtStartPar
Column index of non\sphinxhyphen{}zero elements in the lazy constraint.
\end{quote}

\sphinxAtStartPar
\sphinxcode{\sphinxupquote{rowMatElem}}
\begin{quote}

\sphinxAtStartPar
Values of non\sphinxhyphen{}zero elements in the lazy constraint.
\end{quote}

\sphinxAtStartPar
\sphinxcode{\sphinxupquote{cRowSense}}
\begin{quote}

\sphinxAtStartPar
The sense of the new lazy constraint.

\sphinxAtStartPar
Please refer to {\hyperref[\detokenize{constant:chapconst-constrtype}]{\sphinxcrossref{\DUrole{std,std-ref}{Constraint senses}}}} for all the supported types.

\sphinxAtStartPar
If \sphinxcode{\sphinxupquote{cRowSense}} is 0, then \sphinxcode{\sphinxupquote{dRowBound}} and \sphinxcode{\sphinxupquote{dRowUpper}}
will be treated as lower and upper bounds for the constraint.
This is the recommended method for defining constraints.

\sphinxAtStartPar
If \sphinxcode{\sphinxupquote{cRowSense}} is provided, then \sphinxcode{\sphinxupquote{dRowBound}} and \sphinxcode{\sphinxupquote{dRowUpper}}
will be treated as RHS and \sphinxstylestrong{range} for the constraint.
In this case, \sphinxcode{\sphinxupquote{dRowUpper}}  is only required when
\sphinxcode{\sphinxupquote{cRowSense = COPT\_RANGE}}, where
\begin{quote}

\sphinxAtStartPar
lower bound is \sphinxcode{\sphinxupquote{dRowBound \sphinxhyphen{} dRowUpper}}

\sphinxAtStartPar
upper bound is \sphinxcode{\sphinxupquote{dRowBound}}
\end{quote}
\end{quote}

\sphinxAtStartPar
\sphinxcode{\sphinxupquote{dRowBound}}
\begin{quote}

\sphinxAtStartPar
Lower bound or RHS of the lazy constraint.
\end{quote}

\sphinxAtStartPar
\sphinxcode{\sphinxupquote{dRowUpper}}
\begin{quote}

\sphinxAtStartPar
Upper bound or \sphinxstylestrong{range} of the lazy constraint.
\end{quote}

\sphinxAtStartPar
\sphinxcode{\sphinxupquote{rowName}}
\begin{quote}

\sphinxAtStartPar
The name of the lazy constraint. Can be \sphinxcode{\sphinxupquote{NULL}}.
\end{quote}
\end{quote}
\end{quote}

\subsubsection{COPT\_AddLazyConstrs}
\label{\detokenize{capiref:copt-addlazyconstrs}}\label{\detokenize{capiref:chapapi-cbc-addlazyconstrs}}\begin{quote}

\sphinxAtStartPar
\sphinxstylestrong{Synopsis}
\begin{quote}

\sphinxAtStartPar
\sphinxcode{\sphinxupquote{int COPT\_AddLazyConstrs(}}
\begin{quote}

\sphinxAtStartPar
\sphinxcode{\sphinxupquote{copt\_prob *prob,}}

\sphinxAtStartPar
\sphinxcode{\sphinxupquote{int nAddRow,}}

\sphinxAtStartPar
\sphinxcode{\sphinxupquote{const int *rowMatBeg,}}

\sphinxAtStartPar
\sphinxcode{\sphinxupquote{const int *rowMatCnt,}}

\sphinxAtStartPar
\sphinxcode{\sphinxupquote{const int *rowMatIdx,}}

\sphinxAtStartPar
\sphinxcode{\sphinxupquote{const double *rowMatElem,}}

\sphinxAtStartPar
\sphinxcode{\sphinxupquote{const char *rowSense,}}

\sphinxAtStartPar
\sphinxcode{\sphinxupquote{const double *rowBound,}}

\sphinxAtStartPar
\sphinxcode{\sphinxupquote{const double *rowUpper,}}

\sphinxAtStartPar
\sphinxcode{\sphinxupquote{char const *const *rowNames)}}
\end{quote}
\end{quote}

\sphinxAtStartPar
\sphinxstylestrong{Description}
\begin{quote}

\sphinxAtStartPar
Add a set of lazy constraints to the MIP model.
\end{quote}

\sphinxAtStartPar
\sphinxstylestrong{Arguments}
\begin{quote}

\sphinxAtStartPar
\sphinxcode{\sphinxupquote{prob}}
\begin{quote}

\sphinxAtStartPar
The COPT problem.
\end{quote}

\sphinxAtStartPar
\sphinxcode{\sphinxupquote{nAddRow}}
\begin{quote}

\sphinxAtStartPar
Number of new lazy constraints.
\end{quote}

\sphinxAtStartPar
\sphinxcode{\sphinxupquote{rowMatBeg, rowMatCnt, rowMatIdx}} and \sphinxcode{\sphinxupquote{rowMatElem}}
\begin{quote}

\sphinxAtStartPar
Defines the coefficient matrix in compressed row storage (CRS) format.
The CRS format is similar to the CCS format described
in the \sphinxstylestrong{other information} of \sphinxcode{\sphinxupquote{COPT\_LoadProb}}.
\end{quote}

\sphinxAtStartPar
\sphinxcode{\sphinxupquote{rowSense}}
\begin{quote}

\sphinxAtStartPar
Senses of new lazy constraints.

\sphinxAtStartPar
Please refer to {\hyperref[\detokenize{constant:chapconst-constrtype}]{\sphinxcrossref{\DUrole{std,std-ref}{Constraint senses}}}} for all the supported types.

\sphinxAtStartPar
If \sphinxcode{\sphinxupquote{rowSense}} is \sphinxcode{\sphinxupquote{NULL}}, then \sphinxcode{\sphinxupquote{rowBound}} and \sphinxcode{\sphinxupquote{rowUpper}}
will be treated as lower and upper bounds for constraints.
This is the recommended method for defining constraints.

\sphinxAtStartPar
If \sphinxcode{\sphinxupquote{rowSense}} is provided, then \sphinxcode{\sphinxupquote{rowBound}} and \sphinxcode{\sphinxupquote{rowUpper}}
will be treated as RHS and \sphinxstylestrong{range} for constraints.
In this case, \sphinxcode{\sphinxupquote{rowUpper}} is only required when there
are \sphinxcode{\sphinxupquote{COPT\_RANGE}} constraints, where the
\begin{quote}

\sphinxAtStartPar
lower bound is \sphinxcode{\sphinxupquote{rowBound{[}i{]} \sphinxhyphen{} fabs(rowUpper{[}i{]})}}

\sphinxAtStartPar
upper bound is \sphinxcode{\sphinxupquote{rowBound{[}i{]}}}
\end{quote}
\end{quote}

\sphinxAtStartPar
\sphinxcode{\sphinxupquote{rowBound}}
\begin{quote}

\sphinxAtStartPar
Lower bounds or RHS of new lazy constraints.
\end{quote}

\sphinxAtStartPar
\sphinxcode{\sphinxupquote{rowUpper}}
\begin{quote}

\sphinxAtStartPar
Upper bounds or \sphinxstylestrong{range} of new lazy constraints.
\end{quote}

\sphinxAtStartPar
\sphinxcode{\sphinxupquote{rowNames}}
\begin{quote}

\sphinxAtStartPar
Names of new lazy constraints. Can be \sphinxcode{\sphinxupquote{NULL}}.
\end{quote}
\end{quote}
\end{quote}

\subsubsection{COPT\_AddUserCut}
\label{\detokenize{capiref:copt-addusercut}}\label{\detokenize{capiref:chapapi-cbc-addusercut}}\begin{quote}

\sphinxAtStartPar
\sphinxstylestrong{Synopsis}
\begin{quote}

\sphinxAtStartPar
\sphinxcode{\sphinxupquote{int COPT\_AddUserCut(}}
\begin{quote}

\sphinxAtStartPar
\sphinxcode{\sphinxupquote{copt\_prob* prob,}}

\sphinxAtStartPar
\sphinxcode{\sphinxupquote{int nRowMatCnt,}}

\sphinxAtStartPar
\sphinxcode{\sphinxupquote{const int* rowMatIdx,}}

\sphinxAtStartPar
\sphinxcode{\sphinxupquote{const double* rowMmatElem,}}

\sphinxAtStartPar
\sphinxcode{\sphinxupquote{char cRowSense,}}

\sphinxAtStartPar
\sphinxcode{\sphinxupquote{double dRowBound,}}

\sphinxAtStartPar
\sphinxcode{\sphinxupquote{double dRowUpper,}}

\sphinxAtStartPar
\sphinxcode{\sphinxupquote{const char* rowName)}}
\end{quote}
\end{quote}

\sphinxAtStartPar
\sphinxstylestrong{Description}
\begin{quote}

\sphinxAtStartPar
Add a user cut to the MIP model.
\end{quote}

\sphinxAtStartPar
\sphinxstylestrong{Arguments}
\begin{quote}

\sphinxAtStartPar
\sphinxcode{\sphinxupquote{prob}}
\begin{quote}

\sphinxAtStartPar
The COPT problem.
\end{quote}

\sphinxAtStartPar
\sphinxcode{\sphinxupquote{nRowMatCnt}}
\begin{quote}

\sphinxAtStartPar
Number of non\sphinxhyphen{}zero elements in the user cut.
\end{quote}

\sphinxAtStartPar
\sphinxcode{\sphinxupquote{rowMatIdx}}
\begin{quote}

\sphinxAtStartPar
Column index of non\sphinxhyphen{}zero elements in the user cut.
\end{quote}

\sphinxAtStartPar
\sphinxcode{\sphinxupquote{rowMatElem}}
\begin{quote}

\sphinxAtStartPar
Values of non\sphinxhyphen{}zero elements in the user cut.
\end{quote}

\sphinxAtStartPar
\sphinxcode{\sphinxupquote{cRowSense}}
\begin{quote}

\sphinxAtStartPar
The sense of the user cut.

\sphinxAtStartPar
Please refer to {\hyperref[\detokenize{constant:chapconst-constrtype}]{\sphinxcrossref{\DUrole{std,std-ref}{Constraint senses}}}} for all the supported types.

\sphinxAtStartPar
If \sphinxcode{\sphinxupquote{cRowSense}} is 0, then \sphinxcode{\sphinxupquote{dRowBound}} and \sphinxcode{\sphinxupquote{dRowUpper}}
will be treated as lower and upper bounds for the constraint.
This is the recommended method for defining constraints.

\sphinxAtStartPar
If \sphinxcode{\sphinxupquote{cRowSense}} is provided, then \sphinxcode{\sphinxupquote{dRowBound}} and \sphinxcode{\sphinxupquote{dRowUpper}}
will be treated as RHS and \sphinxstylestrong{range} for the constraint.
In this case, \sphinxcode{\sphinxupquote{dRowUpper}}  is only required when
\sphinxcode{\sphinxupquote{cRowSense = COPT\_RANGE}}, where
\begin{quote}

\sphinxAtStartPar
lower bound is \sphinxcode{\sphinxupquote{dRowBound \sphinxhyphen{} dRowUpper}}

\sphinxAtStartPar
upper bound is \sphinxcode{\sphinxupquote{dRowBound}}
\end{quote}
\end{quote}

\sphinxAtStartPar
\sphinxcode{\sphinxupquote{dRowBound}}
\begin{quote}

\sphinxAtStartPar
Lower bound or RHS of the user cut.
\end{quote}

\sphinxAtStartPar
\sphinxcode{\sphinxupquote{dRowUpper}}
\begin{quote}

\sphinxAtStartPar
Upper bound or \sphinxstylestrong{range} of the user cut.
\end{quote}

\sphinxAtStartPar
\sphinxcode{\sphinxupquote{rowName}}
\begin{quote}

\sphinxAtStartPar
The name of the user cut. Can be \sphinxcode{\sphinxupquote{NULL}}.
\end{quote}
\end{quote}
\end{quote}

\subsubsection{COPT\_AddUserCuts}
\label{\detokenize{capiref:copt-addusercuts}}\label{\detokenize{capiref:chapapi-cbc-addusercuts}}\begin{quote}

\sphinxAtStartPar
\sphinxstylestrong{Synopsis}
\begin{quote}

\sphinxAtStartPar
\sphinxcode{\sphinxupquote{int COPT\_CALL COPT\_AddUserCuts(}}
\begin{quote}

\sphinxAtStartPar
\sphinxcode{\sphinxupquote{copt\_prob *prob,}}

\sphinxAtStartPar
\sphinxcode{\sphinxupquote{int nAddRow,}}

\sphinxAtStartPar
\sphinxcode{\sphinxupquote{const int *rowMatBeg,}}

\sphinxAtStartPar
\sphinxcode{\sphinxupquote{const int *rowMatCnt,}}

\sphinxAtStartPar
\sphinxcode{\sphinxupquote{const int *rowMatIdx,}}

\sphinxAtStartPar
\sphinxcode{\sphinxupquote{const double *rowMatElem,}}

\sphinxAtStartPar
\sphinxcode{\sphinxupquote{const char *rowSense,}}

\sphinxAtStartPar
\sphinxcode{\sphinxupquote{const double *rowBound,}}

\sphinxAtStartPar
\sphinxcode{\sphinxupquote{const double *rowUpper,}}

\sphinxAtStartPar
\sphinxcode{\sphinxupquote{char const *const *rowNames)}}
\end{quote}
\end{quote}

\sphinxAtStartPar
\sphinxstylestrong{Description}
\begin{quote}

\sphinxAtStartPar
Add a set of user cuts to the MIP model.
\end{quote}

\sphinxAtStartPar
\sphinxstylestrong{Arguments}
\begin{quote}

\sphinxAtStartPar
\sphinxcode{\sphinxupquote{prob}}
\begin{quote}

\sphinxAtStartPar
The COPT problem.
\end{quote}

\sphinxAtStartPar
\sphinxcode{\sphinxupquote{nAddRow}}
\begin{quote}

\sphinxAtStartPar
Number of new user cuts.
\end{quote}

\sphinxAtStartPar
\sphinxcode{\sphinxupquote{rowMatBeg, rowMatCnt, rowMatIdx}} and \sphinxcode{\sphinxupquote{rowMatElem}}
\begin{quote}

\sphinxAtStartPar
Defines the coefficient matrix in compressed row storage (CRS) format.
The CRS format is similar to the CCS format described
in the \sphinxstylestrong{other information} of \sphinxcode{\sphinxupquote{COPT\_LoadProb}}.
\end{quote}

\sphinxAtStartPar
\sphinxcode{\sphinxupquote{rowSense}}
\begin{quote}

\sphinxAtStartPar
Senses of new user cuts.

\sphinxAtStartPar
Please refer to {\hyperref[\detokenize{constant:chapconst-constrtype}]{\sphinxcrossref{\DUrole{std,std-ref}{Constraint senses}}}} for all the supported types.

\sphinxAtStartPar
If \sphinxcode{\sphinxupquote{rowSense}} is \sphinxcode{\sphinxupquote{NULL}}, then \sphinxcode{\sphinxupquote{rowBound}} and \sphinxcode{\sphinxupquote{rowUpper}}
will be treated as lower and upper bounds for constraints.
This is the recommended method for defining constraints.

\sphinxAtStartPar
If \sphinxcode{\sphinxupquote{rowSense}} is provided, then \sphinxcode{\sphinxupquote{rowBound}} and \sphinxcode{\sphinxupquote{rowUpper}}
will be treated as RHS and \sphinxstylestrong{range} for constraints.
In this case, \sphinxcode{\sphinxupquote{rowUpper}} is only required when there
are \sphinxcode{\sphinxupquote{COPT\_RANGE}} constraints, where the
\begin{quote}

\sphinxAtStartPar
lower bound is \sphinxcode{\sphinxupquote{rowBound{[}i{]} \sphinxhyphen{} fabs(rowUpper{[}i{]})}}

\sphinxAtStartPar
upper bound is \sphinxcode{\sphinxupquote{rowBound{[}i{]}}}
\end{quote}
\end{quote}

\sphinxAtStartPar
\sphinxcode{\sphinxupquote{rowBound}}
\begin{quote}

\sphinxAtStartPar
Lower bounds or RHS of new user cuts.
\end{quote}

\sphinxAtStartPar
\sphinxcode{\sphinxupquote{rowUpper}}
\begin{quote}

\sphinxAtStartPar
Upper bounds or \sphinxstylestrong{range} of new user cuts.
\end{quote}

\sphinxAtStartPar
\sphinxcode{\sphinxupquote{rowNames}}
\begin{quote}

\sphinxAtStartPar
Names of new user cuts. Can be \sphinxcode{\sphinxupquote{NULL}}.
\end{quote}
\end{quote}
\end{quote}

\subsubsection{COPT\_SetCallback}
\label{\detokenize{capiref:copt-setcallback}}\label{\detokenize{capiref:chapapi-cbc-setcallback}}\begin{quote}

\sphinxAtStartPar
\sphinxstylestrong{Synopsis}
\begin{quote}

\sphinxAtStartPar
\sphinxcode{\sphinxupquote{int COPT\_SetCallback(}}
\begin{quote}

\sphinxAtStartPar
\sphinxcode{\sphinxupquote{copt\_prob *prob,}}

\sphinxAtStartPar
\sphinxcode{\sphinxupquote{int (COPT\_CALL *cb)(copt\_prob *prob, void *cbdata, int cbctx, void *userdata),}}

\sphinxAtStartPar
\sphinxcode{\sphinxupquote{int cbctx,}}

\sphinxAtStartPar
\sphinxcode{\sphinxupquote{void *userdata)}}
\end{quote}
\end{quote}

\sphinxAtStartPar
\sphinxstylestrong{Description}
\begin{quote}

\sphinxAtStartPar
Set the callback function of the model.
\end{quote}

\sphinxAtStartPar
\sphinxstylestrong{Arguments}
\begin{quote}

\sphinxAtStartPar
\sphinxcode{\sphinxupquote{prob}}
\begin{quote}

\sphinxAtStartPar
The COPT problem.
\end{quote}

\sphinxAtStartPar
\sphinxcode{\sphinxupquote{cb}}
\begin{quote}

\sphinxAtStartPar
Callback function.
\end{quote}

\sphinxAtStartPar
\sphinxcode{\sphinxupquote{cbctx}}
\begin{quote}

\sphinxAtStartPar
Callback context. Please refer to {\hyperref[\detokenize{constant:chapconst-cbc}]{\sphinxcrossref{\DUrole{std,std-ref}{Callback context}}}} .
\end{quote}

\sphinxAtStartPar
\sphinxcode{\sphinxupquote{userdata}}
\begin{quote}

\sphinxAtStartPar
User defined data. The data will be passed to the solver without modification.
\end{quote}
\end{quote}
\end{quote}

\subsubsection{COPT\_GetCallbackInfo}
\label{\detokenize{capiref:copt-getcallbackinfo}}\label{\detokenize{capiref:chapapi-cbc-getcallbackinfo}}\begin{quote}

\sphinxAtStartPar
\sphinxstylestrong{Synopsis}
\begin{quote}

\sphinxAtStartPar
\sphinxcode{\sphinxupquote{int COPT\_GetCallbackInfo(void *cbdata, const char* cbinfo, void *p\_val)}}
\end{quote}

\sphinxAtStartPar
\sphinxstylestrong{Description}
\begin{quote}

\sphinxAtStartPar
Retrieve the value of the specified callback information.
\end{quote}

\sphinxAtStartPar
\sphinxstylestrong{Arguments}
\begin{quote}

\sphinxAtStartPar
\sphinxcode{\sphinxupquote{cbdata}}
\begin{quote}

\sphinxAtStartPar
The cbdata argument that was passed into the user callback by COPT. This argument must be passed unmodified from the user callback to COPT\_GetCallbackInfo().
\end{quote}

\sphinxAtStartPar
\sphinxcode{\sphinxupquote{cbinfo}}
\begin{quote}

\sphinxAtStartPar
The name of the callback information. Please refer to {\hyperref[\detokenize{capiref:chapapi-cbcinfo}]{\sphinxcrossref{\DUrole{std,std-ref}{Callback information}}}} for possible values.
\end{quote}

\sphinxAtStartPar
\sphinxcode{\sphinxupquote{p\_val}}
\begin{quote}

\sphinxAtStartPar
Pointer to the value of the callback information.
\end{quote}
\end{quote}
\end{quote}

\subsubsection{COPT\_AddCallbackSolution}
\label{\detokenize{capiref:copt-addcallbacksolution}}\label{\detokenize{capiref:chapapi-cbc-addcallbacksolution}}\begin{quote}

\sphinxAtStartPar
\sphinxstylestrong{Synopsis}
\begin{quote}

\sphinxAtStartPar
\sphinxcode{\sphinxupquote{int COPT\_AddCallbackSolution(void *cbdata, const double *sol, double* p\_objval)}}
\end{quote}

\sphinxAtStartPar
\sphinxstylestrong{Description}
\begin{quote}

\sphinxAtStartPar
Set feasible solutions for the specified variables.
\end{quote}

\sphinxAtStartPar
\sphinxstylestrong{Arguments}
\begin{quote}

\sphinxAtStartPar
\sphinxcode{\sphinxupquote{cbdata}}
\begin{quote}

\sphinxAtStartPar
The cbdata argument that was passed into the user callback by COPT. This argument must be passed unmodified from the user callback to COPT\_AddCallbackSolution().
\end{quote}

\sphinxAtStartPar
\sphinxcode{\sphinxupquote{sol}}
\begin{quote}

\sphinxAtStartPar
The solution vector.
\end{quote}

\sphinxAtStartPar
\sphinxcode{\sphinxupquote{p\_objval}}
\begin{quote}

\sphinxAtStartPar
Pointer to the objective value for solution.
\end{quote}
\end{quote}
\end{quote}

\subsubsection{COPT\_AddCallbackUserCut}
\label{\detokenize{capiref:copt-addcallbackusercut}}\label{\detokenize{capiref:chapapi-cbc-addcallbackusercut}}\begin{quote}

\sphinxAtStartPar
\sphinxstylestrong{Synopsis}
\begin{quote}

\sphinxAtStartPar
\sphinxcode{\sphinxupquote{int COPT\_AddCallbackUserCut(}}
\begin{quote}

\sphinxAtStartPar
\sphinxcode{\sphinxupquote{void *cbdata,}}

\sphinxAtStartPar
\sphinxcode{\sphinxupquote{int nRowMatCnt,}}

\sphinxAtStartPar
\sphinxcode{\sphinxupquote{const int *rowMatIdx,}}

\sphinxAtStartPar
\sphinxcode{\sphinxupquote{const double *rowMatElem,}}

\sphinxAtStartPar
\sphinxcode{\sphinxupquote{char cRowSense,}}

\sphinxAtStartPar
\sphinxcode{\sphinxupquote{double dRowRhs)}}
\end{quote}
\end{quote}

\sphinxAtStartPar
\sphinxstylestrong{Description}
\begin{quote}

\sphinxAtStartPar
Add a user cut to the MIP model from within the callback function.
\end{quote}

\sphinxAtStartPar
\sphinxstylestrong{Arguments}
\begin{quote}

\sphinxAtStartPar
\sphinxcode{\sphinxupquote{cbdata}}
\begin{quote}

\sphinxAtStartPar
The cbdata argument that was passed into the user callback by COPT. This argument must be passed unmodified from the user callback to COPT\_AddCallbackUserCut().
\end{quote}

\sphinxAtStartPar
\sphinxcode{\sphinxupquote{nRowMatCnt}}
\begin{quote}

\sphinxAtStartPar
Number of non\sphinxhyphen{}zero elements in the user cut.
\end{quote}

\sphinxAtStartPar
\sphinxcode{\sphinxupquote{rowMatIdx}}
\begin{quote}

\sphinxAtStartPar
Column index of non\sphinxhyphen{}zero elements in the user cut.
\end{quote}

\sphinxAtStartPar
\sphinxcode{\sphinxupquote{rowMatElem}}
\begin{quote}

\sphinxAtStartPar
Values of non\sphinxhyphen{}zero elements in the user cut.
\end{quote}

\sphinxAtStartPar
\sphinxcode{\sphinxupquote{cRowSense}}
\begin{quote}

\sphinxAtStartPar
The sense of the new user cut. It supports for \sphinxcode{\sphinxupquote{LESS\_EQUAL}}, \sphinxcode{\sphinxupquote{GREATER\_EQUAL}}, \sphinxcode{\sphinxupquote{EQUAL}} and \sphinxcode{\sphinxupquote{FREE}} .

\sphinxAtStartPar
The user cut added from within callback can only have a single bound.
\end{quote}

\sphinxAtStartPar
\sphinxcode{\sphinxupquote{dRowRhs}}
\begin{quote}

\sphinxAtStartPar
RHS of the user cut.
\end{quote}
\end{quote}
\end{quote}

\subsubsection{COPT\_AddCallbackUserCuts}
\label{\detokenize{capiref:copt-addcallbackusercuts}}\label{\detokenize{capiref:chapapi-cbc-addcallbackusercuts}}\begin{quote}

\sphinxAtStartPar
\sphinxstylestrong{Synopsis}
\begin{quote}

\sphinxAtStartPar
\sphinxcode{\sphinxupquote{int COPT\_AddCallbackUserCuts(}}
\begin{quote}

\sphinxAtStartPar
\sphinxcode{\sphinxupquote{void *cbdata,}}

\sphinxAtStartPar
\sphinxcode{\sphinxupquote{int nAddRow,}}

\sphinxAtStartPar
\sphinxcode{\sphinxupquote{const int *rowMatBeg,}}

\sphinxAtStartPar
\sphinxcode{\sphinxupquote{const int *rowMatCnt,}}

\sphinxAtStartPar
\sphinxcode{\sphinxupquote{const int *rowMatIdx,}}

\sphinxAtStartPar
\sphinxcode{\sphinxupquote{const double *rowMatElem,}}

\sphinxAtStartPar
\sphinxcode{\sphinxupquote{const char *rowSense,}}

\sphinxAtStartPar
\sphinxcode{\sphinxupquote{const double *rowRhs)}}
\end{quote}
\end{quote}

\sphinxAtStartPar
\sphinxstylestrong{Description}
\begin{quote}

\sphinxAtStartPar
Add a set of user cuts to the MIP model from within the callback function.
\end{quote}

\sphinxAtStartPar
\sphinxstylestrong{Arguments}
\begin{quote}

\sphinxAtStartPar
\sphinxcode{\sphinxupquote{cbdata}}
\begin{quote}

\sphinxAtStartPar
The cbdata argument that was passed into the user callback by COPT. This argument must be passed unmodified from the user callback to COPT\_AddCallbackUserCuts().
\end{quote}

\sphinxAtStartPar
\sphinxcode{\sphinxupquote{nAddRow}}
\begin{quote}

\sphinxAtStartPar
Number of new user cuts.
\end{quote}

\sphinxAtStartPar
\sphinxcode{\sphinxupquote{rowMatBeg, rowMatCnt, rowMatIdx}} and \sphinxcode{\sphinxupquote{rowMatElem}}
\begin{quote}

\sphinxAtStartPar
Defines the coefficient matrix in compressed row storage (CRS) format.
The CRS format is similar to the CCS format described
in the \sphinxstylestrong{other information} of \sphinxcode{\sphinxupquote{COPT\_LoadProb}}.
\end{quote}

\sphinxAtStartPar
\sphinxcode{\sphinxupquote{rowSense}}
\begin{quote}

\sphinxAtStartPar
Senses of new user cuts. It supports for \sphinxcode{\sphinxupquote{LESS\_EQUAL}}, \sphinxcode{\sphinxupquote{GREATER\_EQUAL}}, \sphinxcode{\sphinxupquote{EQUAL}} and \sphinxcode{\sphinxupquote{FREE}} .

\sphinxAtStartPar
The user cuts added from within callback can only have single bounds.
\end{quote}

\sphinxAtStartPar
\sphinxcode{\sphinxupquote{rowRhs}}
\begin{quote}

\sphinxAtStartPar
RHS of new user cuts.
\end{quote}
\end{quote}
\end{quote}

\subsubsection{COPT\_AddCallbackLazyConstr}
\label{\detokenize{capiref:copt-addcallbacklazyconstr}}\label{\detokenize{capiref:chapapi-cbc-addcallbacklazyconstr}}\begin{quote}

\sphinxAtStartPar
\sphinxstylestrong{Synopsis}
\begin{quote}

\sphinxAtStartPar
\sphinxcode{\sphinxupquote{int COPT\_AddCallbackLazyConstr(}}
\begin{quote}

\sphinxAtStartPar
\sphinxcode{\sphinxupquote{void *cbdata,}}

\sphinxAtStartPar
\sphinxcode{\sphinxupquote{int nRowMatCnt,}}

\sphinxAtStartPar
\sphinxcode{\sphinxupquote{const int *rowMatIdx,}}

\sphinxAtStartPar
\sphinxcode{\sphinxupquote{const double *rowMatElem,}}

\sphinxAtStartPar
\sphinxcode{\sphinxupquote{char cRowSense,}}

\sphinxAtStartPar
\sphinxcode{\sphinxupquote{double dRowRhs)}}
\end{quote}
\end{quote}

\sphinxAtStartPar
\sphinxstylestrong{Description}
\begin{quote}

\sphinxAtStartPar
Add a lazy constraint to the MIP model from within the callback function.
\end{quote}

\sphinxAtStartPar
\sphinxstylestrong{Arguments}
\begin{quote}

\sphinxAtStartPar
\sphinxcode{\sphinxupquote{cbdata}}
\begin{quote}

\sphinxAtStartPar
The cbdata argument that was passed into the user callback by COPT. This argument must be passed unmodified from the user callback to COPT\_AddCallbackLazyConstr().
\end{quote}

\sphinxAtStartPar
\sphinxcode{\sphinxupquote{nRowMatCnt}}
\begin{quote}

\sphinxAtStartPar
Number of non\sphinxhyphen{}zero elements in the lazy constraint.

\sphinxAtStartPar
When \sphinxcode{\sphinxupquote{nRowMatCnt\textless{}=0}}, the MIP candidate solution will be directly rejected without adding a lazy constraint.
And \sphinxcode{\sphinxupquote{rowMatIdx}}, \sphinxcode{\sphinxupquote{rowMatElem}}, \sphinxcode{\sphinxupquote{cRowSense}} and \sphinxcode{\sphinxupquote{dRowRhs}} will be ignored.
\end{quote}

\sphinxAtStartPar
\sphinxcode{\sphinxupquote{rowMatIdx}}
\begin{quote}

\sphinxAtStartPar
Column index of non\sphinxhyphen{}zero elements in the lazy constraint.
\end{quote}

\sphinxAtStartPar
\sphinxcode{\sphinxupquote{rowMatElem}}
\begin{quote}

\sphinxAtStartPar
Values of non\sphinxhyphen{}zero elements in the lazy constraint.
\end{quote}

\sphinxAtStartPar
\sphinxcode{\sphinxupquote{cRowSense}}
\begin{quote}

\sphinxAtStartPar
The sense of new lazy constraint. It supports for \sphinxcode{\sphinxupquote{LESS\_EQUAL}}, \sphinxcode{\sphinxupquote{GREATER\_EQUAL}}, \sphinxcode{\sphinxupquote{EQUAL}} and \sphinxcode{\sphinxupquote{FREE}} .

\sphinxAtStartPar
The lazy constraint added from within callback can only have a single bound.
\end{quote}

\sphinxAtStartPar
\sphinxcode{\sphinxupquote{dRowRhs}}
\begin{quote}

\sphinxAtStartPar
RHS of the lazy constraint.
\end{quote}
\end{quote}
\end{quote}

\subsubsection{COPT\_AddCallbackLazyConstrs}
\label{\detokenize{capiref:copt-addcallbacklazyconstrs}}\label{\detokenize{capiref:chapapi-cbc-addcallbacklazyconstrs}}\begin{quote}

\sphinxAtStartPar
\sphinxstylestrong{Synopsis}
\begin{quote}

\sphinxAtStartPar
\sphinxcode{\sphinxupquote{int COPT\_AddCallbackLazyConstrs(}}
\begin{quote}

\sphinxAtStartPar
\sphinxcode{\sphinxupquote{void *cbdata,}}

\sphinxAtStartPar
\sphinxcode{\sphinxupquote{int nAddRow,}}

\sphinxAtStartPar
\sphinxcode{\sphinxupquote{const int *rowMatBeg,}}

\sphinxAtStartPar
\sphinxcode{\sphinxupquote{const int *rowMatCnt,}}

\sphinxAtStartPar
\sphinxcode{\sphinxupquote{const int *rowMatIdx,}}

\sphinxAtStartPar
\sphinxcode{\sphinxupquote{const double *rowMatElem,}}

\sphinxAtStartPar
\sphinxcode{\sphinxupquote{const char *rowSense,}}

\sphinxAtStartPar
\sphinxcode{\sphinxupquote{const double *rowRhs)}}
\end{quote}
\end{quote}

\sphinxAtStartPar
\sphinxstylestrong{Description}
\begin{quote}

\sphinxAtStartPar
Add a set of lazy constraints to the MIP model from within the callback function.
\end{quote}

\sphinxAtStartPar
\sphinxstylestrong{Arguments}
\begin{quote}

\sphinxAtStartPar
\sphinxcode{\sphinxupquote{cbdata}}
\begin{quote}

\sphinxAtStartPar
The cbdata argument that was passed into the user callback by COPT. This argument must be passed unmodified from the user callback to COPT\_AddCallbackLazyConstrs().
\end{quote}

\sphinxAtStartPar
\sphinxcode{\sphinxupquote{nAddRow}}
\begin{quote}

\sphinxAtStartPar
Number of new lazy constraints.

\sphinxAtStartPar
When \sphinxcode{\sphinxupquote{nAddRow\textless{}=0}}, the MIP candidate solution will be directly rejected without adding lazy constraints.

\sphinxAtStartPar
And other parameters(apart from \sphinxcode{\sphinxupquote{cbdata}}) will be ignored.
\end{quote}

\sphinxAtStartPar
\sphinxcode{\sphinxupquote{rowMatBeg, rowMatCnt, rowMatIdx}} and \sphinxcode{\sphinxupquote{rowMatElem}}
\begin{quote}

\sphinxAtStartPar
Defines the coefficient matrix in compressed row storage (CRS) format.
The CRS format is similar to the CCS format described
in the \sphinxstylestrong{other information} of \sphinxcode{\sphinxupquote{COPT\_LoadProb}}.
\end{quote}

\sphinxAtStartPar
\sphinxcode{\sphinxupquote{rowSense}}
\begin{quote}

\sphinxAtStartPar
Senses of new lazy constraints. It supports for \sphinxcode{\sphinxupquote{LESS\_EQUAL}}, \sphinxcode{\sphinxupquote{GREATER\_EQUAL}}, \sphinxcode{\sphinxupquote{EQUAL}} and \sphinxcode{\sphinxupquote{FREE}} .

\sphinxAtStartPar
The lazy constraints added from within callback can only have single bounds.
\end{quote}

\sphinxAtStartPar
\sphinxcode{\sphinxupquote{rowRhs}}
\begin{quote}

\sphinxAtStartPar
RHS of new lazy constraints.
\end{quote}
\end{quote}
\end{quote}

\sphinxAtStartPar
Note: If the user gives a solution to COPT in a callback, the user should ensure that this fulfills the lazy constraint callback since that won’t be called for user\sphinxhyphen{}provided solutions.

\subsection{Other API functions}
\label{\detokenize{capiref:other-api-functions}}

\subsubsection{COPT\_GetBanner}
\label{\detokenize{capiref:copt-getbanner}}\begin{quote}

\sphinxAtStartPar
\sphinxstylestrong{Synopsis}
\begin{quote}

\sphinxAtStartPar
\sphinxcode{\sphinxupquote{int COPT\_GetBanner(char *buff, int buffSize)}}
\end{quote}

\sphinxAtStartPar
\sphinxstylestrong{Description}
\begin{quote}

\sphinxAtStartPar
Obtains a C\sphinxhyphen{}style string as banner, which describes
the COPT version information.
\end{quote}

\sphinxAtStartPar
\sphinxstylestrong{Arguments}
\begin{quote}

\sphinxAtStartPar
\sphinxcode{\sphinxupquote{buff}}
\begin{quote}

\sphinxAtStartPar
A buffer for holding the resulting string.
\end{quote}

\sphinxAtStartPar
\sphinxcode{\sphinxupquote{buffSize}}
\begin{quote}

\sphinxAtStartPar
The size of the above buffer.
\end{quote}
\end{quote}
\end{quote}

\subsubsection{COPT\_GetRetcodeMsg}
\label{\detokenize{capiref:copt-getretcodemsg}}\begin{quote}

\sphinxAtStartPar
\sphinxstylestrong{Synopsis}
\begin{quote}

\sphinxAtStartPar
\sphinxcode{\sphinxupquote{int COPT\_GetRetcodeMsg(int code, char *buff, int buffSize)}}
\end{quote}

\sphinxAtStartPar
\sphinxstylestrong{Description}
\begin{quote}

\sphinxAtStartPar
Obtains a C\sphinxhyphen{}style string which explains a
return code value in plain text.
\end{quote}

\sphinxAtStartPar
\sphinxstylestrong{Arguments}
\begin{quote}

\sphinxAtStartPar
\sphinxcode{\sphinxupquote{code}}
\begin{quote}

\sphinxAtStartPar
The return code from a COPT API function.
\end{quote}

\sphinxAtStartPar
\sphinxcode{\sphinxupquote{buff}}
\begin{quote}

\sphinxAtStartPar
A buffer for holding the resulting string.
\end{quote}

\sphinxAtStartPar
\sphinxcode{\sphinxupquote{buffSize}}
\begin{quote}

\sphinxAtStartPar
The size of the above buffer.
\end{quote}
\end{quote}
\end{quote}

\subsubsection{COPT\_Interrupt}
\label{\detokenize{capiref:copt-interrupt}}\begin{quote}

\sphinxAtStartPar
\sphinxstylestrong{Synopsis}
\begin{quote}

\sphinxAtStartPar
\sphinxcode{\sphinxupquote{int COPT\_Interrupt(copt\_prob *prob)}}
\end{quote}

\sphinxAtStartPar
\sphinxstylestrong{Description}
\begin{quote}

\sphinxAtStartPar
Interrupt solving process of current problem.
\end{quote}

\sphinxAtStartPar
\sphinxstylestrong{Arguments}
\begin{quote}

\sphinxAtStartPar
\sphinxcode{\sphinxupquote{prob}}
\begin{quote}

\sphinxAtStartPar
The COPT problem.
\end{quote}
\end{quote}
\end{quote}

\subsection{Multi\sphinxhyphen{}objective Model API Functions}
\label{\detokenize{capiref:multi-objective-model-api-functions}}\label{\detokenize{capiref:chapapi-mobj}}

\subsubsection{COPT\_MultiObjSetColObj}
\label{\detokenize{capiref:copt-multiobjsetcolobj}}\begin{quote}

\sphinxAtStartPar
\sphinxstylestrong{Synopsis}
\begin{quote}

\sphinxAtStartPar
\sphinxcode{\sphinxupquote{int COPT\_MultiObjSetColObj(}}
\begin{quote}

\sphinxAtStartPar
\sphinxcode{\sphinxupquote{copt\_prob *prob,}}

\sphinxAtStartPar
\sphinxcode{\sphinxupquote{int iObj,}}

\sphinxAtStartPar
\sphinxcode{\sphinxupquote{int num,}}

\sphinxAtStartPar
\sphinxcode{\sphinxupquote{const int *list,}}

\sphinxAtStartPar
\sphinxcode{\sphinxupquote{const double *colObj)}}
\end{quote}
\end{quote}

\sphinxAtStartPar
\sphinxstylestrong{Description}
\begin{quote}

\sphinxAtStartPar
Set the specified linear objective in a multi\sphinxhyphen{}objective model.

\sphinxAtStartPar
Note: this operation overwrites the existing objective definition.
\end{quote}

\sphinxAtStartPar
\sphinxstylestrong{Arguments}
\begin{quote}

\sphinxAtStartPar
\sphinxcode{\sphinxupquote{prob}}
\begin{quote}

\sphinxAtStartPar
COPT problem.
\end{quote}

\sphinxAtStartPar
\sphinxcode{\sphinxupquote{iObj}}
\begin{quote}

\sphinxAtStartPar
Index of the objective.
\end{quote}

\sphinxAtStartPar
\sphinxcode{\sphinxupquote{num}}
\begin{quote}

\sphinxAtStartPar
Number of linear objective coefficients.
\end{quote}

\sphinxAtStartPar
\sphinxcode{\sphinxupquote{list}}
\begin{quote}

\sphinxAtStartPar
Index list of linear objective coefficients.
\end{quote}

\sphinxAtStartPar
\sphinxcode{\sphinxupquote{colObj}}
\begin{quote}

\sphinxAtStartPar
Value list of linear objective coefficients.
\end{quote}
\end{quote}
\end{quote}

\subsubsection{COPT\_MultiObjGetColObj}
\label{\detokenize{capiref:copt-multiobjgetcolobj}}\begin{quote}

\sphinxAtStartPar
\sphinxstylestrong{Synopsis}
\begin{quote}

\sphinxAtStartPar
\sphinxcode{\sphinxupquote{int COPT\_MultiObjGetColObj(}}
\begin{quote}

\sphinxAtStartPar
\sphinxcode{\sphinxupquote{copt\_prob *prob,}}

\sphinxAtStartPar
\sphinxcode{\sphinxupquote{int iObj,}}

\sphinxAtStartPar
\sphinxcode{\sphinxupquote{int num,}}

\sphinxAtStartPar
\sphinxcode{\sphinxupquote{const int *list,}}

\sphinxAtStartPar
\sphinxcode{\sphinxupquote{double *colObj)}}
\end{quote}
\end{quote}

\sphinxAtStartPar
\sphinxstylestrong{Description}
\begin{quote}

\sphinxAtStartPar
Retrieve the specified linear objective in a multi\sphinxhyphen{}objective model.
\end{quote}

\sphinxAtStartPar
\sphinxstylestrong{Arguments}
\begin{quote}

\sphinxAtStartPar
\sphinxcode{\sphinxupquote{prob}}
\begin{quote}

\sphinxAtStartPar
COPT problem.
\end{quote}

\sphinxAtStartPar
\sphinxcode{\sphinxupquote{iObj}}
\begin{quote}

\sphinxAtStartPar
Index of the objective.
\end{quote}

\sphinxAtStartPar
\sphinxcode{\sphinxupquote{num}}
\begin{quote}

\sphinxAtStartPar
Number of nonzero coefficients in the linear objective.
\end{quote}

\sphinxAtStartPar
\sphinxcode{\sphinxupquote{list}}
\begin{quote}

\sphinxAtStartPar
Indices of the nonzero coefficients in the linear objective.
\end{quote}

\sphinxAtStartPar
\sphinxcode{\sphinxupquote{colObj}}
\begin{quote}

\sphinxAtStartPar
Values of the nonzero coefficients in the linear objective.
\end{quote}
\end{quote}
\end{quote}

\subsubsection{COPT\_MultiObjDelColObj}
\label{\detokenize{capiref:copt-multiobjdelcolobj}}\begin{quote}

\sphinxAtStartPar
\sphinxstylestrong{Synopsis}
\begin{quote}

\sphinxAtStartPar
\sphinxcode{\sphinxupquote{int COPT\_MultiObjDelColObj(}}
\begin{quote}

\sphinxAtStartPar
\sphinxcode{\sphinxupquote{copt\_prob *prob,}}

\sphinxAtStartPar
\sphinxcode{\sphinxupquote{int iObj)}}
\end{quote}
\end{quote}

\sphinxAtStartPar
\sphinxstylestrong{Description}
\begin{quote}

\sphinxAtStartPar
Delete the specified linear objective in a multi\sphinxhyphen{}objective model.
\end{quote}

\sphinxAtStartPar
\sphinxstylestrong{Arguments}
\begin{quote}

\sphinxAtStartPar
\sphinxcode{\sphinxupquote{prob}}
\begin{quote}

\sphinxAtStartPar
COPT problem.
\end{quote}

\sphinxAtStartPar
\sphinxcode{\sphinxupquote{iObj}}
\begin{quote}

\sphinxAtStartPar
Index of the objective.
\end{quote}
\end{quote}
\end{quote}

\subsubsection{COPT\_MultiObjSetObjSense}
\label{\detokenize{capiref:copt-multiobjsetobjsense}}\begin{quote}

\sphinxAtStartPar
\sphinxstylestrong{Synopsis}
\begin{quote}

\sphinxAtStartPar
\sphinxcode{\sphinxupquote{int COPT\_MultiObjSetObjSense(}}
\begin{quote}

\sphinxAtStartPar
\sphinxcode{\sphinxupquote{copt\_prob *prob,}}

\sphinxAtStartPar
\sphinxcode{\sphinxupquote{int iObj,}}

\sphinxAtStartPar
\sphinxcode{\sphinxupquote{int iObjSense)}}
\end{quote}
\end{quote}

\sphinxAtStartPar
\sphinxstylestrong{Description}
\begin{quote}

\sphinxAtStartPar
Set the objective sense for the specified objective in a multi\sphinxhyphen{}objective model.
\end{quote}

\sphinxAtStartPar
\sphinxstylestrong{Arguments}
\begin{quote}

\sphinxAtStartPar
\sphinxcode{\sphinxupquote{prob}}
\begin{quote}

\sphinxAtStartPar
COPT problem.
\end{quote}

\sphinxAtStartPar
\sphinxcode{\sphinxupquote{iObj}}
\begin{quote}

\sphinxAtStartPar
Index of the objective.
\end{quote}

\sphinxAtStartPar
\sphinxcode{\sphinxupquote{iObjSense}}
\begin{quote}

\sphinxAtStartPar
The optimization sense, either \sphinxcode{\sphinxupquote{COPT\_MAXIMIZE}} or \sphinxcode{\sphinxupquote{COPT\_MINIMIZE}}.
\end{quote}
\end{quote}
\end{quote}

\subsubsection{COPT\_MultiObjSetObjConst}
\label{\detokenize{capiref:copt-multiobjsetobjconst}}\begin{quote}

\sphinxAtStartPar
\sphinxstylestrong{Synopsis}
\begin{quote}

\sphinxAtStartPar
\sphinxcode{\sphinxupquote{int COPT\_MultiObjSetObjConst(}}
\begin{quote}

\sphinxAtStartPar
\sphinxcode{\sphinxupquote{copt\_prob *prob,}}

\sphinxAtStartPar
\sphinxcode{\sphinxupquote{int iObj,}}

\sphinxAtStartPar
\sphinxcode{\sphinxupquote{double dObjConst)}}
\end{quote}
\end{quote}

\sphinxAtStartPar
\sphinxstylestrong{Description}
\begin{quote}

\sphinxAtStartPar
Set the objective constant for the specified objective in a multi\sphinxhyphen{}objective model.
\end{quote}

\sphinxAtStartPar
\sphinxstylestrong{Arguments}
\begin{quote}

\sphinxAtStartPar
\sphinxcode{\sphinxupquote{prob}}
\begin{quote}

\sphinxAtStartPar
COPT problem.
\end{quote}

\sphinxAtStartPar
\sphinxcode{\sphinxupquote{iObj}}
\begin{quote}

\sphinxAtStartPar
Index of the objective.
\end{quote}

\sphinxAtStartPar
\sphinxcode{\sphinxupquote{dObjConst}}
\begin{quote}

\sphinxAtStartPar
Constant term of the objective function.
\end{quote}
\end{quote}
\end{quote}

\subsubsection{COPT\_MultiObjSetObjParam}
\label{\detokenize{capiref:copt-multiobjsetobjparam}}\begin{quote}

\sphinxAtStartPar
\sphinxstylestrong{Synopsis}
\begin{quote}

\sphinxAtStartPar
\sphinxcode{\sphinxupquote{int COPT\_MultiObjSetObjParam(}}
\begin{quote}

\sphinxAtStartPar
\sphinxcode{\sphinxupquote{copt\_prob *prob,}}

\sphinxAtStartPar
\sphinxcode{\sphinxupquote{int iObj,}}

\sphinxAtStartPar
\sphinxcode{\sphinxupquote{const char *paramName,}}

\sphinxAtStartPar
\sphinxcode{\sphinxupquote{double dblParam)}}
\end{quote}
\end{quote}

\sphinxAtStartPar
\sphinxstylestrong{Description}
\begin{quote}

\sphinxAtStartPar
Set the objective parameter value for the specified objective
in a multi\sphinxhyphen{}objective model.
\end{quote}

\sphinxAtStartPar
\sphinxstylestrong{Arguments}
\begin{quote}

\sphinxAtStartPar
\sphinxcode{\sphinxupquote{prob}}
\begin{quote}

\sphinxAtStartPar
COPT problem.
\end{quote}

\sphinxAtStartPar
\sphinxcode{\sphinxupquote{iObj}}
\begin{quote}

\sphinxAtStartPar
Index of the objective.
\end{quote}

\sphinxAtStartPar
\sphinxcode{\sphinxupquote{paramName}}
\begin{quote}

\sphinxAtStartPar
Name of the objective parameter.

\sphinxAtStartPar
Possible values:
\sphinxcode{\sphinxupquote{COPT\_MULTIOBJ\_PRIORITY}} , \sphinxcode{\sphinxupquote{COPT\_MULTIOBJ\_WEIGHT}} ,
\sphinxcode{\sphinxupquote{COPT\_MULTIOBJ\_ABSTOL}} , \sphinxcode{\sphinxupquote{COPT\_MULTIOBJ\_RELTOL}} .
\end{quote}

\sphinxAtStartPar
\sphinxcode{\sphinxupquote{dblParam}}
\begin{quote}

\sphinxAtStartPar
Value to set for the objective parameter.
\end{quote}
\end{quote}
\end{quote}

\subsubsection{COPT\_MultiObjGetObjParam}
\label{\detokenize{capiref:copt-multiobjgetobjparam}}\begin{quote}

\sphinxAtStartPar
\sphinxstylestrong{Synopsis}
\begin{quote}

\sphinxAtStartPar
\sphinxcode{\sphinxupquote{int COPT\_MultiObjGetObjParam(}}
\begin{quote}

\sphinxAtStartPar
\sphinxcode{\sphinxupquote{copt\_prob *prob,}}

\sphinxAtStartPar
\sphinxcode{\sphinxupquote{int iObj,}}

\sphinxAtStartPar
\sphinxcode{\sphinxupquote{const char *paramName,}}

\sphinxAtStartPar
\sphinxcode{\sphinxupquote{double *p\_dblParam)}}
\end{quote}
\end{quote}

\sphinxAtStartPar
\sphinxstylestrong{Description}
\begin{quote}

\sphinxAtStartPar
Retrieve the objective parameter value for the specified objective
in a multi\sphinxhyphen{}objective model.
\end{quote}

\sphinxAtStartPar
\sphinxstylestrong{Arguments}
\begin{quote}

\sphinxAtStartPar
\sphinxcode{\sphinxupquote{prob}}
\begin{quote}

\sphinxAtStartPar
COPT problem.
\end{quote}

\sphinxAtStartPar
\sphinxcode{\sphinxupquote{iObj}}
\begin{quote}

\sphinxAtStartPar
Index of the objective.
\end{quote}

\sphinxAtStartPar
\sphinxcode{\sphinxupquote{paramName}}
\begin{quote}

\sphinxAtStartPar
Name of the objective parameter.

\sphinxAtStartPar
Possible values are:
\sphinxcode{\sphinxupquote{COPT\_MULTIOBJ\_PRIORITY}} , \sphinxcode{\sphinxupquote{COPT\_MULTIOBJ\_WEIGHT}} ,
\sphinxcode{\sphinxupquote{COPT\_MULTIOBJ\_ABSTOL}} , \sphinxcode{\sphinxupquote{COPT\_MULTIOBJ\_RELTOL}} .
\end{quote}

\sphinxAtStartPar
\sphinxcode{\sphinxupquote{p\_dblParam}}
\begin{quote}

\sphinxAtStartPar
Pointer to store the retrieved parameter value.
\end{quote}
\end{quote}
\end{quote}

\subsubsection{COPT\_MultiObjResetObjParam}
\label{\detokenize{capiref:copt-multiobjresetobjparam}}\begin{quote}

\sphinxAtStartPar
\sphinxstylestrong{Synopsis}
\begin{quote}

\sphinxAtStartPar
\sphinxcode{\sphinxupquote{int COPT\_MultiObjResetObjParam(}}
\begin{quote}

\sphinxAtStartPar
\sphinxcode{\sphinxupquote{copt\_prob *prob,}}

\sphinxAtStartPar
\sphinxcode{\sphinxupquote{int iObj)}}
\end{quote}
\end{quote}

\sphinxAtStartPar
\sphinxstylestrong{Description}
\begin{quote}

\sphinxAtStartPar
Reset all objective parameters of the specified objective
in a multi\sphinxhyphen{}objective model to their default values.
\end{quote}

\sphinxAtStartPar
\sphinxstylestrong{Arguments}
\begin{quote}

\sphinxAtStartPar
\sphinxcode{\sphinxupquote{prob}}
\begin{quote}

\sphinxAtStartPar
COPT problem.
\end{quote}

\sphinxAtStartPar
\sphinxcode{\sphinxupquote{iObj}}
\begin{quote}

\sphinxAtStartPar
Index of the objective.
\end{quote}
\end{quote}
\end{quote}

\subsubsection{COPT\_MultiObjSetIntParam}
\label{\detokenize{capiref:copt-multiobjsetintparam}}\begin{quote}

\sphinxAtStartPar
\sphinxstylestrong{Synopsis}
\begin{quote}

\sphinxAtStartPar
\sphinxcode{\sphinxupquote{int COPT\_MultiObjSetIntParam(}}
\begin{quote}

\sphinxAtStartPar
\sphinxcode{\sphinxupquote{copt\_prob *prob,}}

\sphinxAtStartPar
\sphinxcode{\sphinxupquote{int iObj,}}

\sphinxAtStartPar
\sphinxcode{\sphinxupquote{const char *paramName,}}

\sphinxAtStartPar
\sphinxcode{\sphinxupquote{int intParam)}}
\end{quote}
\end{quote}

\sphinxAtStartPar
\sphinxstylestrong{Description}
\begin{quote}

\sphinxAtStartPar
Set an integer solver parameter for the model
corresponding to the specified objective in a multi\sphinxhyphen{}objective problem.
\end{quote}

\sphinxAtStartPar
\sphinxstylestrong{Arguments}
\begin{quote}

\sphinxAtStartPar
\sphinxcode{\sphinxupquote{prob}}
\begin{quote}

\sphinxAtStartPar
COPT problem.
\end{quote}

\sphinxAtStartPar
\sphinxcode{\sphinxupquote{iObj}}
\begin{quote}

\sphinxAtStartPar
Index of the objective.
\end{quote}

\sphinxAtStartPar
\sphinxcode{\sphinxupquote{paramName}}
\begin{quote}

\sphinxAtStartPar
Name of the parameter.
\end{quote}

\sphinxAtStartPar
\sphinxcode{\sphinxupquote{intParam}}
\begin{quote}

\sphinxAtStartPar
Integer value of the parameter to be set.
\end{quote}
\end{quote}
\end{quote}

\subsubsection{COPT\_MultiObjGetIntParam}
\label{\detokenize{capiref:copt-multiobjgetintparam}}\begin{quote}

\sphinxAtStartPar
\sphinxstylestrong{Synopsis}
\begin{quote}

\sphinxAtStartPar
\sphinxcode{\sphinxupquote{int COPT\_MultiObjGetIntParam(}}
\begin{quote}

\sphinxAtStartPar
\sphinxcode{\sphinxupquote{copt\_prob *prob,}}

\sphinxAtStartPar
\sphinxcode{\sphinxupquote{int iObj,}}

\sphinxAtStartPar
\sphinxcode{\sphinxupquote{const char *paramName,}}

\sphinxAtStartPar
\sphinxcode{\sphinxupquote{int *p\_intParam)}}
\end{quote}
\end{quote}

\sphinxAtStartPar
\sphinxstylestrong{Description}
\begin{quote}

\sphinxAtStartPar
Get the integer solver parameter for the model
corresponding to the specified objective in a multi\sphinxhyphen{}objective problem.
\end{quote}

\sphinxAtStartPar
\sphinxstylestrong{Arguments}
\begin{quote}

\sphinxAtStartPar
\sphinxcode{\sphinxupquote{prob}}
\begin{quote}

\sphinxAtStartPar
COPT problem.
\end{quote}

\sphinxAtStartPar
\sphinxcode{\sphinxupquote{iObj}}
\begin{quote}

\sphinxAtStartPar
Index of the objective.
\end{quote}

\sphinxAtStartPar
\sphinxcode{\sphinxupquote{paramName}}
\begin{quote}

\sphinxAtStartPar
Name of the parameter.
\end{quote}

\sphinxAtStartPar
\sphinxcode{\sphinxupquote{p\_intParam}}
\begin{quote}

\sphinxAtStartPar
Pointer to the retrieved integer parameter value.
\end{quote}
\end{quote}
\end{quote}

\subsubsection{COPT\_MultiObjSetDblParam}
\label{\detokenize{capiref:copt-multiobjsetdblparam}}\begin{quote}

\sphinxAtStartPar
\sphinxstylestrong{Synopsis}
\begin{quote}

\sphinxAtStartPar
\sphinxcode{\sphinxupquote{int COPT\_MultiObjSetDblParam(}}
\begin{quote}

\sphinxAtStartPar
\sphinxcode{\sphinxupquote{copt\_prob *prob,}}

\sphinxAtStartPar
\sphinxcode{\sphinxupquote{int iObj,}}

\sphinxAtStartPar
\sphinxcode{\sphinxupquote{const char *paramName,}}

\sphinxAtStartPar
\sphinxcode{\sphinxupquote{int dblParam)}}
\end{quote}
\end{quote}

\sphinxAtStartPar
\sphinxstylestrong{Description}
\begin{quote}

\sphinxAtStartPar
Set the double solver parameter for the model corresponding to the specified
objective in a multi\sphinxhyphen{}objective problem.
\end{quote}

\sphinxAtStartPar
\sphinxstylestrong{Arguments}
\begin{quote}

\sphinxAtStartPar
\sphinxcode{\sphinxupquote{prob}}
\begin{quote}

\sphinxAtStartPar
COPT problem.
\end{quote}

\sphinxAtStartPar
\sphinxcode{\sphinxupquote{iObj}}
\begin{quote}

\sphinxAtStartPar
Index of the objective.
\end{quote}

\sphinxAtStartPar
\sphinxcode{\sphinxupquote{paramName}}
\begin{quote}

\sphinxAtStartPar
Name of the parameter.
\end{quote}

\sphinxAtStartPar
\sphinxcode{\sphinxupquote{dblParam}}
\begin{quote}

\sphinxAtStartPar
Value of the double parameter to set.
\end{quote}
\end{quote}
\end{quote}

\subsubsection{COPT\_MultiObjGetDblParam}
\label{\detokenize{capiref:copt-multiobjgetdblparam}}\begin{quote}

\sphinxAtStartPar
\sphinxstylestrong{Synopsis}
\begin{quote}

\sphinxAtStartPar
\sphinxcode{\sphinxupquote{int COPT\_MultiObjGetDblParam(}}
\begin{quote}

\sphinxAtStartPar
\sphinxcode{\sphinxupquote{copt\_prob *prob,}}

\sphinxAtStartPar
\sphinxcode{\sphinxupquote{int iObj,}}

\sphinxAtStartPar
\sphinxcode{\sphinxupquote{const char *paramName,}}

\sphinxAtStartPar
\sphinxcode{\sphinxupquote{int *p\_dblParam)}}
\end{quote}
\end{quote}

\sphinxAtStartPar
\sphinxstylestrong{Description}
\begin{quote}

\sphinxAtStartPar
Get the double solver parameter of the model corresponding to the specified
objective in a multi\sphinxhyphen{}objective problem.
\end{quote}

\sphinxAtStartPar
\sphinxstylestrong{Arguments}
\begin{quote}

\sphinxAtStartPar
\sphinxcode{\sphinxupquote{prob}}
\begin{quote}

\sphinxAtStartPar
COPT problem.
\end{quote}

\sphinxAtStartPar
\sphinxcode{\sphinxupquote{iObj}}
\begin{quote}

\sphinxAtStartPar
Index of the objective.
\end{quote}

\sphinxAtStartPar
\sphinxcode{\sphinxupquote{paramName}}
\begin{quote}

\sphinxAtStartPar
Name of the parameter.
\end{quote}

\sphinxAtStartPar
\sphinxcode{\sphinxupquote{p\_dblParam}}
\begin{quote}

\sphinxAtStartPar
Pointer to store the value of the double parameter.
\end{quote}
\end{quote}
\end{quote}

\subsubsection{COPT\_MultiObjResetParam}
\label{\detokenize{capiref:copt-multiobjresetparam}}\begin{quote}

\sphinxAtStartPar
\sphinxstylestrong{Synopsis}
\begin{quote}

\sphinxAtStartPar
\sphinxcode{\sphinxupquote{int COPT\_MultiObjResetParam(}}
\begin{quote}

\sphinxAtStartPar
\sphinxcode{\sphinxupquote{copt\_prob *prob,}}

\sphinxAtStartPar
\sphinxcode{\sphinxupquote{int iObj)}}
\end{quote}
\end{quote}

\sphinxAtStartPar
\sphinxstylestrong{Description}
\begin{quote}

\sphinxAtStartPar
Reset all solver parameters of the model corresponding to the specified
objective in a multi\sphinxhyphen{}objective problem to their default settings.
\end{quote}

\sphinxAtStartPar
\sphinxstylestrong{Arguments}
\begin{quote}

\sphinxAtStartPar
\sphinxcode{\sphinxupquote{prob}}
\begin{quote}

\sphinxAtStartPar
COPT problem.
\end{quote}

\sphinxAtStartPar
\sphinxcode{\sphinxupquote{iObj}}
\begin{quote}

\sphinxAtStartPar
Index of the objective.
\end{quote}
\end{quote}
\end{quote}

\subsubsection{COPT\_MultiObjGetIntAttr}
\label{\detokenize{capiref:copt-multiobjgetintattr}}\begin{quote}

\sphinxAtStartPar
\sphinxstylestrong{Synopsis}
\begin{quote}

\sphinxAtStartPar
\sphinxcode{\sphinxupquote{int COPT\_MultiObjGetIntAttr(}}
\begin{quote}

\sphinxAtStartPar
\sphinxcode{\sphinxupquote{copt\_prob *prob,}}

\sphinxAtStartPar
\sphinxcode{\sphinxupquote{int iObj,}}

\sphinxAtStartPar
\sphinxcode{\sphinxupquote{const char *attrName,}}

\sphinxAtStartPar
\sphinxcode{\sphinxupquote{int *p\_intAttr)}}
\end{quote}
\end{quote}

\sphinxAtStartPar
\sphinxstylestrong{Description}
\begin{quote}

\sphinxAtStartPar
Retrieve an integer attribute of the specified
objective in a multi\sphinxhyphen{}objective problem.
\end{quote}

\sphinxAtStartPar
\sphinxstylestrong{Arguments}
\begin{quote}

\sphinxAtStartPar
\sphinxcode{\sphinxupquote{prob}}
\begin{quote}

\sphinxAtStartPar
COPT problem.
\end{quote}

\sphinxAtStartPar
\sphinxcode{\sphinxupquote{iObj}}
\begin{quote}

\sphinxAtStartPar
Index of the objective.
\end{quote}

\sphinxAtStartPar
\sphinxcode{\sphinxupquote{attrName}}
\begin{quote}

\sphinxAtStartPar
Name of the attribute.

\sphinxAtStartPar
Possible values are:
\sphinxcode{\sphinxupquote{COPT\_DBLATTR\_LPOBJVAL}}, \sphinxcode{\sphinxupquote{COPT\_DBLATTR\_BESTOBJ}},
\sphinxcode{\sphinxupquote{COPT\_DBLATTR\_OBJCONST}}.
\end{quote}

\sphinxAtStartPar
\sphinxcode{\sphinxupquote{p\_intAttr}}
\begin{quote}

\sphinxAtStartPar
Pointer to the integer value of the attribute.
\end{quote}
\end{quote}
\end{quote}

\subsubsection{COPT\_MultiObjGetDblAttr}
\label{\detokenize{capiref:copt-multiobjgetdblattr}}\begin{quote}

\sphinxAtStartPar
\sphinxstylestrong{Synopsis}
\begin{quote}

\sphinxAtStartPar
\sphinxcode{\sphinxupquote{int COPT\_MultiObjGetDblAttr(}}
\begin{quote}

\sphinxAtStartPar
\sphinxcode{\sphinxupquote{copt\_prob *prob,}}

\sphinxAtStartPar
\sphinxcode{\sphinxupquote{int iObj,}}

\sphinxAtStartPar
\sphinxcode{\sphinxupquote{const char *attrName,}}

\sphinxAtStartPar
\sphinxcode{\sphinxupquote{double *p\_dblAttr)}}
\end{quote}
\end{quote}

\sphinxAtStartPar
\sphinxstylestrong{Description}
\begin{quote}

\sphinxAtStartPar
Retrieve a double attribute of the specified
objective in a multi\sphinxhyphen{}objective problem.
\end{quote}

\sphinxAtStartPar
\sphinxstylestrong{Arguments}
\begin{quote}

\sphinxAtStartPar
\sphinxcode{\sphinxupquote{prob}}
\begin{quote}

\sphinxAtStartPar
COPT problem.
\end{quote}

\sphinxAtStartPar
\sphinxcode{\sphinxupquote{iObj}}
\begin{quote}

\sphinxAtStartPar
Index of the objective.
\end{quote}

\sphinxAtStartPar
\sphinxcode{\sphinxupquote{attrName}}
\begin{quote}

\sphinxAtStartPar
Name of the attribute.

\sphinxAtStartPar
Possible values are:
\sphinxcode{\sphinxupquote{COPT\_DBLATTR\_LPOBJVAL}}, \sphinxcode{\sphinxupquote{COPT\_DBLATTR\_BESTOBJ}},
\sphinxcode{\sphinxupquote{COPT\_DBLATTR\_OBJCONST}}.
\end{quote}

\sphinxAtStartPar
\sphinxcode{\sphinxupquote{p\_dblAttr}}
\begin{quote}

\sphinxAtStartPar
Pointer to the double value of the attribute.
\end{quote}
\end{quote}
\end{quote}

\subsubsection{COPT\_MultiObjGetPoolObjVal}
\label{\detokenize{capiref:copt-multiobjgetpoolobjval}}\begin{quote}

\sphinxAtStartPar
\sphinxstylestrong{Synopsis}
\begin{quote}

\sphinxAtStartPar
\sphinxcode{\sphinxupquote{int COPT\_MultiObjGetPoolObjVal(}}
\begin{quote}

\sphinxAtStartPar
\sphinxcode{\sphinxupquote{copt\_prob *prob,}}

\sphinxAtStartPar
\sphinxcode{\sphinxupquote{int iObj,}}

\sphinxAtStartPar
\sphinxcode{\sphinxupquote{int iSol,}}

\sphinxAtStartPar
\sphinxcode{\sphinxupquote{double *p\_objVal)}}
\end{quote}
\end{quote}

\sphinxAtStartPar
\sphinxstylestrong{Description}
\begin{quote}

\sphinxAtStartPar
Retrieve the value of the specified objective at the solution
with index \sphinxcode{\sphinxupquote{iSol}} in the solution pool.
\end{quote}

\sphinxAtStartPar
\sphinxstylestrong{Arguments}
\begin{quote}

\sphinxAtStartPar
\sphinxcode{\sphinxupquote{prob}}
\begin{quote}

\sphinxAtStartPar
COPT problem.
\end{quote}

\sphinxAtStartPar
\sphinxcode{\sphinxupquote{iObj}}
\begin{quote}

\sphinxAtStartPar
Index of the objective.
\end{quote}

\sphinxAtStartPar
\sphinxcode{\sphinxupquote{iSol}}
\begin{quote}

\sphinxAtStartPar
Index of the solution in the solution pool.
\end{quote}

\sphinxAtStartPar
\sphinxcode{\sphinxupquote{p\_objVal}}
\begin{quote}

\sphinxAtStartPar
Pointer to the value of the objective.
\end{quote}
\end{quote}
\end{quote}

\subsection{CPU and Memory Binding Functions}
\label{\detokenize{capiref:cpu-and-memory-binding-functions}}\label{\detokenize{capiref:chapapi-cpumem}}

\subsubsection{COPT\_BindNumaCpu}
\label{\detokenize{capiref:copt-bindnumacpu}}\begin{quote}

\sphinxAtStartPar
\sphinxstylestrong{Synopsis}
\begin{quote}

\sphinxAtStartPar
\sphinxcode{\sphinxupquote{int COPT\_BindNumaCpu(}}
\begin{quote}

\sphinxAtStartPar
\sphinxcode{\sphinxupquote{copt\_env *env,}}

\sphinxAtStartPar
\sphinxcode{\sphinxupquote{int numaNode)}}
\end{quote}
\end{quote}

\sphinxAtStartPar
\sphinxstylestrong{Description}
\begin{quote}

\sphinxAtStartPar
Bind the CPUs of the current process to a NUMA node.
\end{quote}

\sphinxAtStartPar
\sphinxstylestrong{Arguments}
\begin{quote}

\sphinxAtStartPar
\sphinxcode{\sphinxupquote{env}}
\begin{quote}

\sphinxAtStartPar
The COPT environment.
\end{quote}

\sphinxAtStartPar
\sphinxcode{\sphinxupquote{numaNode}}
\begin{quote}

\sphinxAtStartPar
ID of the NUMA node.
\end{quote}
\end{quote}
\end{quote}

\subsubsection{COPT\_BindNumaMem}
\label{\detokenize{capiref:copt-bindnumamem}}\begin{quote}

\sphinxAtStartPar
\sphinxstylestrong{Synopsis}
\begin{quote}

\sphinxAtStartPar
\sphinxcode{\sphinxupquote{int COPT\_BindNumaMem(}}
\begin{quote}

\sphinxAtStartPar
\sphinxcode{\sphinxupquote{copt\_env *env,}}

\sphinxAtStartPar
\sphinxcode{\sphinxupquote{int numaNode)}}
\end{quote}
\end{quote}

\sphinxAtStartPar
\sphinxstylestrong{Description}
\begin{quote}

\sphinxAtStartPar
Bind memory for the current process to a NUMA node (Linux only).
\end{quote}

\sphinxAtStartPar
\sphinxstylestrong{Arguments}
\begin{quote}

\sphinxAtStartPar
\sphinxcode{\sphinxupquote{env}}
\begin{quote}

\sphinxAtStartPar
The COPT environment.
\end{quote}

\sphinxAtStartPar
\sphinxcode{\sphinxupquote{numaNode}}
\begin{quote}

\sphinxAtStartPar
ID of the NUMA node.
\end{quote}
\end{quote}
\end{quote}

\subsubsection{COPT\_GetCpuAffinity}
\label{\detokenize{capiref:copt-getcpuaffinity}}\begin{quote}

\sphinxAtStartPar
\sphinxstylestrong{Synopsis}
\begin{quote}

\sphinxAtStartPar
\sphinxcode{\sphinxupquote{int COPT\_GetCpuAffinity(}}
\begin{quote}

\sphinxAtStartPar
\sphinxcode{\sphinxupquote{copt\_env *env,}}

\sphinxAtStartPar
\sphinxcode{\sphinxupquote{int *cpuList,}}

\sphinxAtStartPar
\sphinxcode{\sphinxupquote{int len,}}

\sphinxAtStartPar
\sphinxcode{\sphinxupquote{int *pReqSize)}}
\end{quote}
\end{quote}

\sphinxAtStartPar
\sphinxstylestrong{Description}
\begin{quote}

\sphinxAtStartPar
Get CPU affinity for the current process, which is saved in an integer array.
\end{quote}

\sphinxAtStartPar
\sphinxstylestrong{Arguments}
\begin{quote}

\sphinxAtStartPar
\sphinxcode{\sphinxupquote{env}}
\begin{quote}

\sphinxAtStartPar
The COPT environment.
\end{quote}

\sphinxAtStartPar
\sphinxcode{\sphinxupquote{cpuList}}
\begin{quote}

\sphinxAtStartPar
Integer array to store the list of CPU IDs.
\end{quote}

\sphinxAtStartPar
\sphinxcode{\sphinxupquote{len}}
\begin{quote}

\sphinxAtStartPar
Length of the \sphinxcode{\sphinxupquote{cpuList}} array.
\end{quote}

\sphinxAtStartPar
\sphinxcode{\sphinxupquote{pReqSize}}
\begin{quote}

\sphinxAtStartPar
Pointer to store the actual number of bound CPUs.
\end{quote}
\end{quote}
\end{quote}

\subsubsection{COPT\_GetNumaNodeCount}
\label{\detokenize{capiref:copt-getnumanodecount}}\begin{quote}

\sphinxAtStartPar
\sphinxstylestrong{Synopsis}
\begin{quote}

\sphinxAtStartPar
\sphinxcode{\sphinxupquote{int COPT\_GetNumaNodeCount(}}
\begin{quote}

\sphinxAtStartPar
\sphinxcode{\sphinxupquote{copt\_env *env,}}

\sphinxAtStartPar
\sphinxcode{\sphinxupquote{int *pCnt)}}
\end{quote}
\end{quote}

\sphinxAtStartPar
\sphinxstylestrong{Description}
\begin{quote}

\sphinxAtStartPar
Get count of NUMA nodes.
\end{quote}

\sphinxAtStartPar
\sphinxstylestrong{Arguments}
\begin{quote}

\sphinxAtStartPar
\sphinxcode{\sphinxupquote{env}}
\begin{quote}

\sphinxAtStartPar
The COPT environment.
\end{quote}

\sphinxAtStartPar
\sphinxcode{\sphinxupquote{pCnt}}
\begin{quote}

\sphinxAtStartPar
Pointer to store the count of NUMA nodes.
\end{quote}
\end{quote}
\end{quote}

\subsubsection{COPT\_SetCpuAffinity}
\label{\detokenize{capiref:copt-setcpuaffinity}}\begin{quote}

\sphinxAtStartPar
\sphinxstylestrong{Synopsis}
\begin{quote}

\sphinxAtStartPar
\sphinxcode{\sphinxupquote{int COPT\_SetCpuAffinity(}}
\begin{quote}

\sphinxAtStartPar
\sphinxcode{\sphinxupquote{copt\_env *env,}}

\sphinxAtStartPar
\sphinxcode{\sphinxupquote{const char *hexMask)}}
\end{quote}
\end{quote}

\sphinxAtStartPar
\sphinxstylestrong{Description}
\begin{quote}

\sphinxAtStartPar
Set CPU affinity with given mask string.
\end{quote}

\sphinxAtStartPar
\sphinxstylestrong{Arguments}
\begin{quote}

\sphinxAtStartPar
\sphinxcode{\sphinxupquote{env}}
\begin{quote}

\sphinxAtStartPar
The COPT environment.
\end{quote}

\sphinxAtStartPar
\sphinxcode{\sphinxupquote{hexMask}}
\begin{quote}

\sphinxAtStartPar
CPU mask string of hexadecimal characters.
\end{quote}
\end{quote}
\end{quote}

\sphinxstepscope

\chapter{Python API Reference}
\label{\detokenize{pyapiref:python-api-reference}}\label{\detokenize{pyapiref:chappyapi}}\label{\detokenize{pyapiref::doc}}
\sphinxAtStartPar
The \sphinxstylestrong{Cardinal Optimizer} provides a Python API library. This chapter documents all
COPT Python constants and API functions for python applications.

\section{Constants}
\label{\detokenize{pyapiref:constants}}\label{\detokenize{pyapiref:chappyapi-const}}
\sphinxAtStartPar
Python Constants are necessary to solve a problem using the Python interface.
There are four types of constants defined in COPT Python API library. They are general
constants, information constants, parameters and attributes.

\subsection{General Constants}
\label{\detokenize{pyapiref:general-constants}}\label{\detokenize{pyapiref:chappyapi-const-general}}
\sphinxAtStartPar
For the contents of Python general constants, see {\hyperref[\detokenize{constant:chapconst}]{\sphinxcrossref{\DUrole{std,std-ref}{General Constants}}}}.

\sphinxAtStartPar
General constants are those commonly used in modeling, such as optimization directions, variable types,
and solving status, etc. Users may refer to general constants with a \sphinxcode{\sphinxupquote{\textquotesingle{}COPT\textquotesingle{}}} prefix. For instance,
\sphinxcode{\sphinxupquote{COPT.VERSION\_MAJOR}} is the major version number of the \sphinxstylestrong{Cardinal Optimizer}.

\subsection{Attributes}
\label{\detokenize{pyapiref:attributes}}\label{\detokenize{pyapiref:chappyapi-const-attrs}}
\sphinxAtStartPar
For the contents of Python attribute constants, see {\hyperref[\detokenize{attribute:chapattrs}]{\sphinxcrossref{\DUrole{std,std-ref}{Attributes}}}}.

\sphinxAtStartPar
In the Python API, users may refer to attributes using a \sphinxcode{\sphinxupquote{\textquotesingle{}COPT.Attr\textquotesingle{}}} prefix.
For instance, \sphinxcode{\sphinxupquote{COPT.Attr.Cols}} is the number of variables or columns in the model.

\sphinxAtStartPar
In the Python API, user can get the attribute value by specifying the attribute name.
Attributes are mostly used in \sphinxcode{\sphinxupquote{Model.getAttr()}} method to query properties of the model, please refer to {\hyperref[\detokenize{pyapiref:chappyapi-model}]{\sphinxcrossref{\DUrole{std,std-ref}{Python API: Model Class}}}} for details. Here is an example:
\begin{itemize}
\item {} 
\sphinxAtStartPar
\sphinxcode{\sphinxupquote{Model.getAttr()}}: \sphinxcode{\sphinxupquote{Model.getAttr("Cols")}} obtain the number of variables or columns in the model.

\end{itemize}

\subsection{Information}
\label{\detokenize{pyapiref:information}}\label{\detokenize{pyapiref:chappyapi-const-info}}
\sphinxAtStartPar
For the contents of Python API information class constants, see {\hyperref[\detokenize{information:chapinfo}]{\sphinxcrossref{\DUrole{std,std-ref}{Information}}}}.

\sphinxAtStartPar
In the Python API, user can access the information through the \sphinxcode{\sphinxupquote{COPT.Info}} prefix. For instance, \sphinxcode{\sphinxupquote{COPT.Info.Obj}} is the objective function coefficients for variables (columns).

\sphinxAtStartPar
In the Python API, user can get or set the information value of the object by specifying the information name:
\begin{itemize}
\item {} 
\sphinxAtStartPar
Get the value of variable or constraint information: \sphinxcode{\sphinxupquote{Model.getInfo()}} / \sphinxcode{\sphinxupquote{Var.getInfo()}} / \sphinxcode{\sphinxupquote{Constraint.getInfo()}}

\item {} 
\sphinxAtStartPar
Set the value of variable or constraint information: \sphinxcode{\sphinxupquote{Model.setInfo()}} / \sphinxcode{\sphinxupquote{Var.setInfo()}} / \sphinxcode{\sphinxupquote{Constraint.setInfo()}}

\end{itemize}

\subsection{Callback Information}
\label{\detokenize{pyapiref:callback-information}}\label{\detokenize{pyapiref:chappyapi-const-cbcinfo}}
\sphinxAtStartPar
For the content of Python API callback information class constants, see {\hyperref[\detokenize{information:chapinfo-cbc}]{\sphinxcrossref{\DUrole{std,std-ref}{Callback Information}}}}.

\sphinxAtStartPar
In the Python API, callback\sphinxhyphen{}related information constants are defined in the \sphinxcode{\sphinxupquote{CbInfo}} class.
User can access the callback information via \sphinxcode{\sphinxupquote{COPT.CbInfo.}} prefix.

\sphinxAtStartPar
For instance, \sphinxcode{\sphinxupquote{COPT.CbInfo.BestObj}} is the current best objective.

\sphinxAtStartPar
In the Python API, user can get the value of callback information by specifying the information name.

\sphinxAtStartPar
For instance, \sphinxcode{\sphinxupquote{CallbackBase.getInfo(COPT.CbInfo.BestObj)}} : get the value of the current best objective.

\subsection{Parameters}
\label{\detokenize{pyapiref:parameters}}\label{\detokenize{pyapiref:chappyapi-const-param}}
\sphinxAtStartPar
For the contents of Python API Parameters class constants, see {\hyperref[\detokenize{parameter:chapparams}]{\sphinxcrossref{\DUrole{std,std-ref}{Parameters}}}}.

\sphinxAtStartPar
Parameters control the operation of the \sphinxstylestrong{Cardinal Optimizer}. They can be modified before the
optimization begins.

\sphinxAtStartPar
In the Python API, user can access parameters through the \sphinxcode{\sphinxupquote{COPT.Param}} prefix. For instance, \sphinxcode{\sphinxupquote{COPT.Param.TimeLimit}} is time limit in seconds of the optimization.

\sphinxAtStartPar
In the Python API, user can get and set the parameter value by specifying the parameter name.
The provided functions are as follows, please refer to {\hyperref[\detokenize{pyapiref:chappyapi-model}]{\sphinxcrossref{\DUrole{std,std-ref}{Python API: Model Class}}}} for details.
\begin{itemize}
\item {} 
\sphinxAtStartPar
Get detailed information of the specified parameter (current value/max/min): \sphinxcode{\sphinxupquote{Model.getParamInfo()}}

\item {} 
\sphinxAtStartPar
Get the current value of the specified parameter: \sphinxcode{\sphinxupquote{Model.getParam()}}

\item {} 
\sphinxAtStartPar
Set the specified parameter value: \sphinxcode{\sphinxupquote{Model.setParam()}}

\end{itemize}

\section{Python Modeling Classes}
\label{\detokenize{pyapiref:python-modeling-classes}}
\sphinxAtStartPar
Python modeling classes are essential for the Python interface of Cardinal Optimizer.
It provides plentiful easy\sphinxhyphen{}to\sphinxhyphen{}use methods to quickly build optimization models in complex
practical scenarios. This section will explains these functions and their usage.

\subsection{EnvrConfig Class}
\label{\detokenize{pyapiref:envrconfig-class}}\label{\detokenize{pyapiref:chappyapi-envrconfig}}
\sphinxAtStartPar
EnvrConfig object contains operations related to client configuration, and provides the following methods:

\subsubsection{EnvrConfig()}
\label{\detokenize{pyapiref:envrconfig}}\begin{quote}

\sphinxAtStartPar
\sphinxstylestrong{Synopsis}
\begin{quote}

\sphinxAtStartPar
\sphinxcode{\sphinxupquote{EnvrConfig()}}
\end{quote}

\sphinxAtStartPar
\sphinxstylestrong{Description}
\begin{quote}

\sphinxAtStartPar
Constructor of EnvrConfig class. This method creates and returns an {\hyperref[\detokenize{pyapiref:chappyapi-envrconfig}]{\sphinxcrossref{\DUrole{std,std-ref}{EnvrConfig Class}}}} object.
\end{quote}

\sphinxAtStartPar
\sphinxstylestrong{Example}
\end{quote}

\begin{sphinxVerbatim}[commandchars=\\\{\}]
\PYG{c+c1}{\PYGZsh{} Create client configuration}
\PYG{n}{envconfig} \PYG{o}{=} \PYG{n}{EnvrConfig}\PYG{p}{(}\PYG{p}{)}
\end{sphinxVerbatim}

\subsubsection{EnvrConfig.set()}
\label{\detokenize{pyapiref:envrconfig-set}}\begin{quote}

\sphinxAtStartPar
\sphinxstylestrong{Synopsis}
\begin{quote}

\sphinxAtStartPar
\sphinxcode{\sphinxupquote{set(name, value)}}
\end{quote}

\sphinxAtStartPar
\sphinxstylestrong{Description}
\begin{quote}

\sphinxAtStartPar
Set client configuration.
\end{quote}

\sphinxAtStartPar
\sphinxstylestrong{Arguments}
\begin{quote}

\sphinxAtStartPar
\sphinxcode{\sphinxupquote{name}}
\begin{quote}

\sphinxAtStartPar
Name of configuration parameter. Please refer to {\hyperref[\detokenize{constant:chapconst-client}]{\sphinxcrossref{\DUrole{std,std-ref}{Client configuration}}}} for possible values.
\end{quote}

\sphinxAtStartPar
\sphinxcode{\sphinxupquote{value}}
\begin{quote}

\sphinxAtStartPar
Value of configuration parameter.
\end{quote}
\end{quote}

\sphinxAtStartPar
\sphinxstylestrong{Example}
\end{quote}

\begin{sphinxVerbatim}[commandchars=\\\{\}]
\PYG{c+c1}{\PYGZsh{} Set client configuration}
\PYG{n}{envconfig}\PYG{o}{.}\PYG{n}{set}\PYG{p}{(}\PYG{n}{COPT}\PYG{o}{.}\PYG{n}{CLIENT\PYGZus{}WAITTIME}\PYG{p}{,} \PYG{l+m+mi}{600}\PYG{p}{)}
\PYG{n}{envconfig}\PYG{o}{.}\PYG{n}{set}\PYG{p}{(}\PYG{n}{COPT}\PYG{o}{.}\PYG{n}{CLIENT\PYGZus{}CLUSTER}\PYG{p}{,} \PYG{l+s+s2}{\PYGZdq{}}\PYG{l+s+s2}{127.0.0.1}\PYG{l+s+s2}{\PYGZdq{}}\PYG{p}{)}
\PYG{c+c1}{\PYGZsh{} Turn off the banner output when creating COPT environment (such as version, etc.)}
\PYG{n}{envconfig}\PYG{o}{.}\PYG{n}{set}\PYG{p}{(}\PYG{l+s+s2}{\PYGZdq{}}\PYG{l+s+s2}{nobanner}\PYG{l+s+s2}{\PYGZdq{}}\PYG{p}{,} \PYG{l+s+s2}{\PYGZdq{}}\PYG{l+s+s2}{1}\PYG{l+s+s2}{\PYGZdq{}}\PYG{p}{)}
\end{sphinxVerbatim}

\subsection{Envr Class}
\label{\detokenize{pyapiref:envr-class}}\label{\detokenize{pyapiref:chappyapi-envr}}
\sphinxAtStartPar
Envr object contains operations related to COPT optimization environment, and provides the following methods:

\subsubsection{Envr()}
\label{\detokenize{pyapiref:envr}}\begin{quote}

\sphinxAtStartPar
\sphinxstylestrong{Synopsis}
\begin{quote}

\sphinxAtStartPar
\sphinxcode{\sphinxupquote{Envr(arg=None)}}
\end{quote}

\sphinxAtStartPar
\sphinxstylestrong{Description}
\begin{quote}

\sphinxAtStartPar
Function for constructing Envr object. This method creates and returns an {\hyperref[\detokenize{pyapiref:chappyapi-envr}]{\sphinxcrossref{\DUrole{std,std-ref}{Envr Class}}}} object.
\end{quote}

\sphinxAtStartPar
\sphinxstylestrong{Arguments}
\begin{quote}

\sphinxAtStartPar
\sphinxcode{\sphinxupquote{arg}}
\begin{quote}

\sphinxAtStartPar
Path of license file or client configuration. Optional argument, defaults to \sphinxcode{\sphinxupquote{None}}.
\end{quote}
\end{quote}

\sphinxAtStartPar
\sphinxstylestrong{Example}
\end{quote}

\begin{sphinxVerbatim}[commandchars=\\\{\}]
\PYG{c+c1}{\PYGZsh{} Create solving environment}
\PYG{n}{env} \PYG{o}{=} \PYG{n}{Envr}\PYG{p}{(}\PYG{p}{)}
\end{sphinxVerbatim}

\subsubsection{Envr.createModel()}
\label{\detokenize{pyapiref:envr-createmodel}}\begin{quote}

\sphinxAtStartPar
\sphinxstylestrong{Synopsis}
\begin{quote}

\sphinxAtStartPar
\sphinxcode{\sphinxupquote{createModel(name="")}}
\end{quote}

\sphinxAtStartPar
\sphinxstylestrong{Description}
\begin{quote}

\sphinxAtStartPar
Create optimization model and return a {\hyperref[\detokenize{pyapiref:chappyapi-model}]{\sphinxcrossref{\DUrole{std,std-ref}{Model Class}}}} object.
\end{quote}

\sphinxAtStartPar
\sphinxstylestrong{Arguments}
\begin{quote}

\sphinxAtStartPar
\sphinxcode{\sphinxupquote{name}}
\begin{quote}

\sphinxAtStartPar
The name of the Model object to be created. Optional,  \sphinxcode{\sphinxupquote{""}} by default.
\end{quote}
\end{quote}

\sphinxAtStartPar
\sphinxstylestrong{Example}
\end{quote}

\begin{sphinxVerbatim}[commandchars=\\\{\}]
\PYG{c+c1}{\PYGZsh{} Create optimization model}
\PYG{n}{model} \PYG{o}{=} \PYG{n}{env}\PYG{o}{.}\PYG{n}{createModel}\PYG{p}{(}\PYG{l+s+s2}{\PYGZdq{}}\PYG{l+s+s2}{coptprob}\PYG{l+s+s2}{\PYGZdq{}}\PYG{p}{)}
\end{sphinxVerbatim}

\subsubsection{Envr.close()}
\label{\detokenize{pyapiref:envr-close}}\begin{quote}

\sphinxAtStartPar
\sphinxstylestrong{Synopsis}
\begin{quote}

\sphinxAtStartPar
\sphinxcode{\sphinxupquote{close()}}
\end{quote}

\sphinxAtStartPar
\sphinxstylestrong{Description}
\begin{quote}

\sphinxAtStartPar
Close connection to remote server. (For floating or cluster license)
\end{quote}

\sphinxAtStartPar
\sphinxstylestrong{Example}
\end{quote}

\begin{sphinxVerbatim}[commandchars=\\\{\}]
\PYG{c+c1}{\PYGZsh{} Close connection to remote server}
\PYG{n}{env}\PYG{o}{.}\PYG{n}{close}\PYG{p}{(}\PYG{p}{)}
\end{sphinxVerbatim}

\subsubsection{Envr.bindNumaCpu()}
\label{\detokenize{pyapiref:envr-bindnumacpu}}\begin{quote}

\sphinxAtStartPar
\sphinxstylestrong{Synopsis}
\begin{quote}

\sphinxAtStartPar
\sphinxcode{\sphinxupquote{bindNumaCpu(numaNode)}}
\end{quote}

\sphinxAtStartPar
\sphinxstylestrong{Description}
\begin{quote}

\sphinxAtStartPar
Bind the CPUs for the current process to a NUMA node.
\end{quote}

\sphinxAtStartPar
\sphinxstylestrong{Arguments}
\begin{quote}

\sphinxAtStartPar
\sphinxcode{\sphinxupquote{numaNode}}
\begin{quote}

\sphinxAtStartPar
ID of a NUMA node.
\end{quote}
\end{quote}

\sphinxAtStartPar
\sphinxstylestrong{Example}
\end{quote}

\begin{sphinxVerbatim}[commandchars=\\\{\}]
\PYG{n}{env}\PYG{o}{.}\PYG{n}{bindNumaCpu}\PYG{p}{(}\PYG{l+m+mi}{0}\PYG{p}{)}
\end{sphinxVerbatim}

\subsubsection{Envr.bindNumaMem()}
\label{\detokenize{pyapiref:envr-bindnumamem}}\begin{quote}

\sphinxAtStartPar
\sphinxstylestrong{Synopsis}
\begin{quote}

\sphinxAtStartPar
\sphinxcode{\sphinxupquote{bindNumaMem(numaNode)}}
\end{quote}

\sphinxAtStartPar
\sphinxstylestrong{Description}
\begin{quote}

\sphinxAtStartPar
Bind memory for the current process to a NUMA node (Linux only).
\end{quote}

\sphinxAtStartPar
\sphinxstylestrong{Arguments}
\begin{quote}

\sphinxAtStartPar
\sphinxcode{\sphinxupquote{numaNode}}
\begin{quote}

\sphinxAtStartPar
ID of a NUMA node.
\end{quote}
\end{quote}

\sphinxAtStartPar
\sphinxstylestrong{Example}
\end{quote}

\begin{sphinxVerbatim}[commandchars=\\\{\}]
\PYG{n}{env}\PYG{o}{.}\PYG{n}{bindNumaMem}\PYG{p}{(}\PYG{l+m+mi}{0}\PYG{p}{)}
\end{sphinxVerbatim}

\subsubsection{Envr.getCpuAffinity()}
\label{\detokenize{pyapiref:envr-getcpuaffinity}}\begin{quote}

\sphinxAtStartPar
\sphinxstylestrong{Synopsis}
\begin{quote}

\sphinxAtStartPar
\sphinxcode{\sphinxupquote{getCpuAffinity()}}
\end{quote}

\sphinxAtStartPar
\sphinxstylestrong{Description}
\begin{quote}

\sphinxAtStartPar
Get CPU affinity for the current process, which is saved in an integer array.
\end{quote}

\sphinxAtStartPar
\sphinxstylestrong{Example}
\end{quote}

\begin{sphinxVerbatim}[commandchars=\\\{\}]
\PYG{n}{cpu\PYGZus{}list} \PYG{o}{=} \PYG{n}{env}\PYG{o}{.}\PYG{n}{getCpuAffinity}\PYG{p}{(}\PYG{p}{)}
\end{sphinxVerbatim}

\subsubsection{Envr.getNumaNodeCount()}
\label{\detokenize{pyapiref:envr-getnumanodecount}}\begin{quote}

\sphinxAtStartPar
\sphinxstylestrong{Synopsis}
\begin{quote}

\sphinxAtStartPar
\sphinxcode{\sphinxupquote{getNumaNodeCount()}}
\end{quote}

\sphinxAtStartPar
\sphinxstylestrong{Description}
\begin{quote}

\sphinxAtStartPar
Get count of NUMA nodes.
\end{quote}

\sphinxAtStartPar
\sphinxstylestrong{Example}
\end{quote}

\begin{sphinxVerbatim}[commandchars=\\\{\}]
\PYG{n}{count} \PYG{o}{=} \PYG{n}{env}\PYG{o}{.}\PYG{n}{getNumaNodeCount}\PYG{p}{(}\PYG{p}{)}
\end{sphinxVerbatim}

\subsubsection{Envr.setCpuAffinity()}
\label{\detokenize{pyapiref:envr-setcpuaffinity}}\begin{quote}

\sphinxAtStartPar
\sphinxstylestrong{Synopsis}
\begin{quote}

\sphinxAtStartPar
\sphinxcode{\sphinxupquote{setCpuAffinity(hexMask)}}
\end{quote}

\sphinxAtStartPar
\sphinxstylestrong{Description}
\begin{quote}

\sphinxAtStartPar
Set CPU affinity with given mask string.
\end{quote}

\sphinxAtStartPar
\sphinxstylestrong{Arguments}
\begin{quote}

\sphinxAtStartPar
\sphinxcode{\sphinxupquote{hexMask}}
\begin{quote}

\sphinxAtStartPar
CPU mask string of hexadecimal characters.
\end{quote}
\end{quote}

\sphinxAtStartPar
\sphinxstylestrong{Example}
\end{quote}

\begin{sphinxVerbatim}[commandchars=\\\{\}]
\PYG{n}{env}\PYG{o}{.}\PYG{n}{setCpuAffinity}\PYG{p}{(}\PYG{l+s+s2}{\PYGZdq{}}\PYG{l+s+s2}{0f}\PYG{l+s+s2}{\PYGZdq{}}\PYG{p}{)}
\end{sphinxVerbatim}

\subsection{Model Class}
\label{\detokenize{pyapiref:model-class}}\label{\detokenize{pyapiref:chappyapi-model}}
\sphinxAtStartPar
For easy access to model’s attributes and optimization parameters, Model object provides methods such as \sphinxcode{\sphinxupquote{Model.Rows}}.
The full list of attributes can be found in  {\hyperref[\detokenize{attribute:chapattrs}]{\sphinxcrossref{\DUrole{std,std-ref}{Attributes}}}}  section.
For convenience, attributes can be accessed by their names in capital or lower case.

\sphinxAtStartPar
Note that for LP or MIP, both the objective value and the solution status can be accessed through \sphinxcode{\sphinxupquote{Model.objval}} and \sphinxcode{\sphinxupquote{Model.status}}.

\sphinxAtStartPar
For optimization parameters, they can be set in the form \sphinxcode{\sphinxupquote{"Model.Param.TimeLimit = 10"}}.
For details of the parameter names supported, please refer to {\hyperref[\detokenize{parameter:chapparams}]{\sphinxcrossref{\DUrole{std,std-ref}{Parameters}}}} section.

\sphinxAtStartPar
Class Model contains COPT model\sphinxhyphen{}related operations and provides the following methods:

\subsubsection{Model.addVar()}
\label{\detokenize{pyapiref:model-addvar}}\begin{quote}

\sphinxAtStartPar
\sphinxstylestrong{Synopsis}
\begin{quote}

\sphinxAtStartPar
\sphinxcode{\sphinxupquote{addVar(lb=0.0, ub=COPT.INFINITY, obj=0.0, vtype=COPT.CONTINUOUS, name="", column=None)}}
\end{quote}

\sphinxAtStartPar
\sphinxstylestrong{Description}
\begin{quote}

\sphinxAtStartPar
Add a decision variable to model and return the added {\hyperref[\detokenize{pyapiref:chappyapi-var}]{\sphinxcrossref{\DUrole{std,std-ref}{Var Class}}}} object.
\end{quote}

\sphinxAtStartPar
\sphinxstylestrong{Arguments}
\begin{quote}

\sphinxAtStartPar
\sphinxcode{\sphinxupquote{lb}}
\begin{quote}

\sphinxAtStartPar
Lower bound for new variable. Optional, 0.0 by default.
\end{quote}

\sphinxAtStartPar
\sphinxcode{\sphinxupquote{ub}}
\begin{quote}

\sphinxAtStartPar
Upper bound for new variable. Optional, \sphinxcode{\sphinxupquote{COPT.INFINITY}} by default.
\end{quote}

\sphinxAtStartPar
\sphinxcode{\sphinxupquote{obj}}
\begin{quote}

\sphinxAtStartPar
Objective parameter for new variable. Optional, 0.0 by default.
\end{quote}

\sphinxAtStartPar
\sphinxcode{\sphinxupquote{vtype}}
\begin{quote}

\sphinxAtStartPar
Variable type. Optional, \sphinxcode{\sphinxupquote{COPT.CONTINUOUS}} by default.
Please refer to {\hyperref[\detokenize{constant:chapconst-vartype}]{\sphinxcrossref{\DUrole{std,std-ref}{Variable types}}}} for possible types.
\end{quote}

\sphinxAtStartPar
\sphinxcode{\sphinxupquote{name}}
\begin{quote}

\sphinxAtStartPar
Name for new variable. Optional, \sphinxcode{\sphinxupquote{""}} by default, which is automatically generated by solver.
\end{quote}

\sphinxAtStartPar
\sphinxcode{\sphinxupquote{column}}
\begin{quote}

\sphinxAtStartPar
Column corresponds to the variable. Optional, \sphinxcode{\sphinxupquote{None}} by default.
\end{quote}
\end{quote}

\sphinxAtStartPar
\sphinxstylestrong{Example}
\end{quote}

\begin{sphinxVerbatim}[commandchars=\\\{\}]
\PYG{c+c1}{\PYGZsh{} Add a continuous variable}
\PYG{n}{x} \PYG{o}{=} \PYG{n}{m}\PYG{o}{.}\PYG{n}{addVar}\PYG{p}{(}\PYG{p}{)}
\PYG{c+c1}{\PYGZsh{} Add a binary variable}
\PYG{n}{y} \PYG{o}{=} \PYG{n}{m}\PYG{o}{.}\PYG{n}{addVar}\PYG{p}{(}\PYG{n}{vtype}\PYG{o}{=}\PYG{n}{COPT}\PYG{o}{.}\PYG{n}{BINARY}\PYG{p}{)}
\PYG{c+c1}{\PYGZsh{} Add an integer variable with lowerbound \PYGZhy{}1.0, upperbound 1.0, objective coefficient 1.0 and variable name \PYGZdq{}z\PYGZdq{}}
\PYG{n}{z} \PYG{o}{=} \PYG{n}{m}\PYG{o}{.}\PYG{n}{addVar}\PYG{p}{(}\PYG{o}{\PYGZhy{}}\PYG{l+m+mf}{1.0}\PYG{p}{,} \PYG{l+m+mf}{1.0}\PYG{p}{,} \PYG{l+m+mf}{1.0}\PYG{p}{,} \PYG{n}{COPT}\PYG{o}{.}\PYG{n}{INTEGER}\PYG{p}{,} \PYG{l+s+s2}{\PYGZdq{}}\PYG{l+s+s2}{z}\PYG{l+s+s2}{\PYGZdq{}}\PYG{p}{)}
\end{sphinxVerbatim}

\subsubsection{Model.addVars()}
\label{\detokenize{pyapiref:model-addvars}}\begin{quote}

\sphinxAtStartPar
\sphinxstylestrong{Synopsis}
\begin{quote}

\sphinxAtStartPar
\sphinxcode{\sphinxupquote{addVars(*indices, lb=0.0, ub=COPT.INFINITY, obj=0.0, vtype=COPT.CONTINUOUS, nameprefix="C")}}
\end{quote}

\sphinxAtStartPar
\sphinxstylestrong{Description}
\begin{quote}

\sphinxAtStartPar
Add multiple new variables to a model.
Return a {\hyperref[\detokenize{pyapiref:chappyapi-util-tupledict}]{\sphinxcrossref{\DUrole{std,std-ref}{tupledict Class}}}},
whose key is indice of the variable and value is the {\hyperref[\detokenize{pyapiref:chappyapi-var}]{\sphinxcrossref{\DUrole{std,std-ref}{Var Class}}}} object.
\end{quote}

\sphinxAtStartPar
\sphinxstylestrong{Arguments}
\begin{quote}

\sphinxAtStartPar
\sphinxcode{\sphinxupquote{*indices}}
\begin{quote}

\sphinxAtStartPar
Indices for accessing the new variables.
\end{quote}

\sphinxAtStartPar
\sphinxcode{\sphinxupquote{lb}}
\begin{quote}

\sphinxAtStartPar
Lower bounds for new variables. Optional, 0.0 by default.
\end{quote}

\sphinxAtStartPar
\sphinxcode{\sphinxupquote{ub}}
\begin{quote}

\sphinxAtStartPar
Upper bounds for new variables. Optional,  \sphinxcode{\sphinxupquote{COPT.INFINITY}} by default.
\end{quote}

\sphinxAtStartPar
\sphinxcode{\sphinxupquote{obj}}
\begin{quote}

\sphinxAtStartPar
Objective costs for new variables. Optional, 0.0 by default.
\end{quote}

\sphinxAtStartPar
\sphinxcode{\sphinxupquote{vtype}}
\begin{quote}

\sphinxAtStartPar
Variable types. Optional, \sphinxcode{\sphinxupquote{COPT.CONTINUOUS}} by default.
Please refer to {\hyperref[\detokenize{constant:chapconst-vartype}]{\sphinxcrossref{\DUrole{std,std-ref}{Variable types}}}} for possible types.
\end{quote}

\sphinxAtStartPar
\sphinxcode{\sphinxupquote{nameprefix}}
\begin{quote}

\sphinxAtStartPar
Name prefix for new variables. Optional, \sphinxcode{\sphinxupquote{"C"}} by default.
The actual name and the index of the variables are automatically generated by COPT.
\end{quote}
\end{quote}

\sphinxAtStartPar
\sphinxstylestrong{Example}
\end{quote}

\begin{sphinxVerbatim}[commandchars=\\\{\}]
\PYG{c+c1}{\PYGZsh{} Add three\PYGZhy{}dimensional integer variable, 6 variables in total}
\PYG{n}{x} \PYG{o}{=} \PYG{n}{m}\PYG{o}{.}\PYG{n}{addVars}\PYG{p}{(}\PYG{l+m+mi}{1}\PYG{p}{,} \PYG{l+m+mi}{2}\PYG{p}{,} \PYG{l+m+mi}{3}\PYG{p}{,} \PYG{n}{vtype}\PYG{o}{=}\PYG{n}{COPT}\PYG{o}{.}\PYG{n}{INTEGER}\PYG{p}{)}
\PYG{c+c1}{\PYGZsh{} Add two continuous variable y, whose indice is designated by elements in tuplelist and prefix is \PYGZdq{}tl\PYGZdq{}}
\PYG{n}{tl} \PYG{o}{=} \PYG{n}{tuplelist}\PYG{p}{(}\PYG{p}{[}\PYG{p}{(}\PYG{l+m+mi}{0}\PYG{p}{,} \PYG{l+m+mi}{1}\PYG{p}{)}\PYG{p}{,} \PYG{p}{(}\PYG{l+m+mi}{1}\PYG{p}{,} \PYG{l+m+mi}{2}\PYG{p}{)}\PYG{p}{]}\PYG{p}{)}
\PYG{n}{y}  \PYG{o}{=} \PYG{n}{m}\PYG{o}{.}\PYG{n}{addVars}\PYG{p}{(}\PYG{n}{tl}\PYG{p}{,} \PYG{n}{nameprefix}\PYG{o}{=}\PYG{l+s+s2}{\PYGZdq{}}\PYG{l+s+s2}{tl}\PYG{l+s+s2}{\PYGZdq{}}\PYG{p}{)}
\end{sphinxVerbatim}

\subsubsection{Model.addMVar()}
\label{\detokenize{pyapiref:model-addmvar}}\begin{quote}

\sphinxAtStartPar
\sphinxstylestrong{Synopsis}
\begin{quote}

\sphinxAtStartPar
\sphinxcode{\sphinxupquote{addMVar(shape, lb=0.0, ub=COPT.INFINITY, obj=0.0, vtype=COPT.CONTINUOUS, nameprefix="")}}
\end{quote}

\sphinxAtStartPar
\sphinxstylestrong{Description}
\begin{quote}

\sphinxAtStartPar
Add {\hyperref[\detokenize{pyapiref:chappyapi-mvar}]{\sphinxcrossref{\DUrole{std,std-ref}{MVar Class}}}} object to the model. It is used in matrix modeling and can be operated like a multidimensional array in NumPy, its shape and dimensions are similarly defined.
\end{quote}

\sphinxAtStartPar
\sphinxstylestrong{Arguments}
\begin{quote}

\sphinxAtStartPar
\sphinxcode{\sphinxupquote{shape}}
\begin{quote}

\sphinxAtStartPar
The value is an integer, or a tuple of integers. which represents the shape of a {\hyperref[\detokenize{pyapiref:chappyapi-mvar}]{\sphinxcrossref{\DUrole{std,std-ref}{MVar Class}}}} object.
\end{quote}

\sphinxAtStartPar
\sphinxcode{\sphinxupquote{lb}}
\begin{quote}

\sphinxAtStartPar
The lower bound of the variable. Optional parameter, defaults to 0.0.
\end{quote}

\sphinxAtStartPar
\sphinxcode{\sphinxupquote{ub}}
\begin{quote}

\sphinxAtStartPar
The upper bound of the variable. Optional parameter, defaults to \sphinxcode{\sphinxupquote{COPT.INFINITY}}.
\end{quote}

\sphinxAtStartPar
\sphinxcode{\sphinxupquote{obj}}
\begin{quote}

\sphinxAtStartPar
The objective function coefficients for the variables. Optional parameter, defaults to 0.0.
\end{quote}

\sphinxAtStartPar
\sphinxcode{\sphinxupquote{vtype}}
\begin{quote}

\sphinxAtStartPar
The type of the variable. Optional parameter, the default is \sphinxcode{\sphinxupquote{COPT.CONTINUOUS}}, see the possible values in
{\hyperref[\detokenize{constant:chapconst-vartype}]{\sphinxcrossref{\DUrole{std,std-ref}{Variable types}}}}.
\end{quote}

\sphinxAtStartPar
\sphinxcode{\sphinxupquote{nameprefix}}
\begin{quote}

\sphinxAtStartPar
Variable name prefix. Optional parameter, the default is \sphinxcode{\sphinxupquote{""}}, its actual name is automatically generated by combining the subscript of the variable.
\end{quote}
\end{quote}

\sphinxAtStartPar
\sphinxstylestrong{Return value}
\begin{quote}

\sphinxAtStartPar
Returns a {\hyperref[\detokenize{pyapiref:chappyapi-mvar}]{\sphinxcrossref{\DUrole{std,std-ref}{MVar Class}}}} object
\end{quote}

\sphinxAtStartPar
\sphinxstylestrong{Example}

\begin{sphinxVerbatim}[commandchars=\\\{\}]
\PYG{n}{model}\PYG{o}{.}\PYG{n}{addMVar}\PYG{p}{(}\PYG{p}{(}\PYG{l+m+mi}{2}\PYG{p}{,} \PYG{l+m+mi}{3}\PYG{p}{)}\PYG{p}{,} \PYG{n}{lb}\PYG{o}{=}\PYG{l+m+mf}{0.0}\PYG{p}{,} \PYG{n}{nameprefix}\PYG{o}{=}\PYG{l+s+s2}{\PYGZdq{}}\PYG{l+s+s2}{mx}\PYG{l+s+s2}{\PYGZdq{}}\PYG{p}{)}
\end{sphinxVerbatim}
\end{quote}

\subsubsection{Model.addConstr()}
\label{\detokenize{pyapiref:model-addconstr}}\begin{quote}

\sphinxAtStartPar
\sphinxstylestrong{Synopsis}
\begin{quote}

\sphinxAtStartPar
\sphinxcode{\sphinxupquote{addConstr(lhs, sense=None, rhs=None, name="")}}
\end{quote}

\sphinxAtStartPar
\sphinxstylestrong{Description}
\begin{quote}

\sphinxAtStartPar
Add a linear constraint to the model, return {\hyperref[\detokenize{pyapiref:chappyapi-constraint}]{\sphinxcrossref{\DUrole{std,std-ref}{Constraint Class}}}} object
or {\hyperref[\detokenize{pyapiref:chappyapi-mconstr}]{\sphinxcrossref{\DUrole{std,std-ref}{MConstr Class}}}} object;

\sphinxAtStartPar
Add a semidefinite constraint to the model, return {\hyperref[\detokenize{pyapiref:chappyapi-psdconstraint}]{\sphinxcrossref{\DUrole{std,std-ref}{PsdConstraint Class}}}} object
or {\hyperref[\detokenize{pyapiref:chappyapi-mpsdconstr}]{\sphinxcrossref{\DUrole{std,std-ref}{MPsdConstr Class}}}} object;

\sphinxAtStartPar
Add an indicator constraint to the model and return the {\hyperref[\detokenize{pyapiref:chappyapi-genconstr}]{\sphinxcrossref{\DUrole{std,std-ref}{GenConstr Class}}}} object;

\sphinxAtStartPar
Adds a LMI constraint to the model, returning a {\hyperref[\detokenize{pyapiref:chappyapi-lmiconstraint}]{\sphinxcrossref{\DUrole{std,std-ref}{LmiConstraint Class}}}} object.

\sphinxAtStartPar
If a linear constraint added, then the parameter \sphinxcode{\sphinxupquote{lhs}} can take the value of {\hyperref[\detokenize{pyapiref:chappyapi-var}]{\sphinxcrossref{\DUrole{std,std-ref}{Var Class}}}} object,
{\hyperref[\detokenize{pyapiref:chappyapi-linexpr}]{\sphinxcrossref{\DUrole{std,std-ref}{LinExpr Class}}}} object, {\hyperref[\detokenize{pyapiref:chappyapi-constrbuilder}]{\sphinxcrossref{\DUrole{std,std-ref}{ConstrBuilder Class}}}} object or {\hyperref[\detokenize{pyapiref:chappyapi-mconstrbuilder}]{\sphinxcrossref{\DUrole{std,std-ref}{MConstrBuilder Class}}}} .

\sphinxAtStartPar
If a positive semi\sphinxhyphen{}definite constraint added, then the parameter \sphinxcode{\sphinxupquote{lhs}} can take the value of
{\hyperref[\detokenize{pyapiref:chappyapi-psdexpr}]{\sphinxcrossref{\DUrole{std,std-ref}{PsdExpr Class}}}} object, {\hyperref[\detokenize{pyapiref:chappyapi-psdconstrbuilder}]{\sphinxcrossref{\DUrole{std,std-ref}{PsdConstrBuilder Class}}}} object or {\hyperref[\detokenize{pyapiref:chappyapi-mpsdconstrbuilder}]{\sphinxcrossref{\DUrole{std,std-ref}{MPsdConstrBuilder Class}}}} .

\sphinxAtStartPar
If an indicator constraint added, then the parameter \sphinxcode{\sphinxupquote{lhs}} is {\hyperref[\detokenize{pyapiref:chappyapi-genconstrbuilder}]{\sphinxcrossref{\DUrole{std,std-ref}{GenConstrBuilder Class}}}} object, ignoring other parameters.

\sphinxAtStartPar
If a LMI constraint added, then the parameter \sphinxcode{\sphinxupquote{lhs}} can take the value of
{\hyperref[\detokenize{pyapiref:chappyapi-lmiexpr}]{\sphinxcrossref{\DUrole{std,std-ref}{LmiExpr Class}}}} object.
\end{quote}

\sphinxAtStartPar
\sphinxstylestrong{Arguments}
\begin{quote}

\sphinxAtStartPar
\sphinxcode{\sphinxupquote{lhs}}
\begin{quote}

\sphinxAtStartPar
Left\sphinxhyphen{}hand side expression for new linear constraint or constraint builder.
\end{quote}

\sphinxAtStartPar
\sphinxcode{\sphinxupquote{sense}}
\begin{quote}

\sphinxAtStartPar
Sense for the new constraint. Optional, None by default.
Please refer to  {\hyperref[\detokenize{constant:chapconst-constrtype}]{\sphinxcrossref{\DUrole{std,std-ref}{Constraint type}}}}  for possible values.
\end{quote}

\sphinxAtStartPar
\sphinxcode{\sphinxupquote{rhs}}
\begin{quote}

\sphinxAtStartPar
Right\sphinxhyphen{}hand side expression for the new constraint. Optional, None by default.
It can be a constant, or {\hyperref[\detokenize{pyapiref:chappyapi-var}]{\sphinxcrossref{\DUrole{std,std-ref}{Var Class}}}} object, or {\hyperref[\detokenize{pyapiref:chappyapi-linexpr}]{\sphinxcrossref{\DUrole{std,std-ref}{LinExpr Class}}}} object.
\end{quote}

\sphinxAtStartPar
\sphinxcode{\sphinxupquote{name}}
\begin{quote}

\sphinxAtStartPar
Name for new constraint. Optional, \sphinxcode{\sphinxupquote{""}} by default, generated by solver automatically.
\end{quote}
\end{quote}

\sphinxAtStartPar
\sphinxstylestrong{Example}
\end{quote}

\begin{sphinxVerbatim}[commandchars=\\\{\}]
\PYG{c+c1}{\PYGZsh{} Add a linear constraint: x + y == 2}
\PYG{n}{m}\PYG{o}{.}\PYG{n}{addConstr}\PYG{p}{(}\PYG{n}{x} \PYG{o}{+} \PYG{n}{y}\PYG{p}{,} \PYG{n}{COPT}\PYG{o}{.}\PYG{n}{EQUAL}\PYG{p}{,} \PYG{l+m+mi}{2}\PYG{p}{)}
\PYG{c+c1}{\PYGZsh{} Add a linear constraint: x + 2*y \PYGZgt{}= 3}
\PYG{n}{m}\PYG{o}{.}\PYG{n}{addConstr}\PYG{p}{(}\PYG{n}{x} \PYG{o}{+} \PYG{l+m+mi}{2}\PYG{o}{*}\PYG{n}{y} \PYG{o}{\PYGZgt{}}\PYG{o}{=} \PYG{l+m+mf}{3.0}\PYG{p}{)}
\PYG{c+c1}{\PYGZsh{} Add an indicator constraint}
\PYG{n}{m}\PYG{o}{.}\PYG{n}{addConstr}\PYG{p}{(}\PYG{p}{(}\PYG{n}{x} \PYG{o}{==} \PYG{l+m+mi}{1}\PYG{p}{)} \PYG{o}{\PYGZgt{}\PYGZgt{}} \PYG{p}{(}\PYG{l+m+mi}{2}\PYG{o}{*}\PYG{n}{y} \PYG{o}{+} \PYG{l+m+mi}{3}\PYG{o}{*}\PYG{n}{z} \PYG{o}{\PYGZlt{}}\PYG{o}{=} \PYG{l+m+mi}{4}\PYG{p}{)}\PYG{p}{)}
\end{sphinxVerbatim}

\subsubsection{Model.addBoundConstr()}
\label{\detokenize{pyapiref:model-addboundconstr}}\begin{quote}

\sphinxAtStartPar
\sphinxstylestrong{Synopsis}
\begin{quote}

\sphinxAtStartPar
\sphinxcode{\sphinxupquote{addBoundConstr(expr, lb=\sphinxhyphen{}COPT.INFINITY, ub=COPT.INFINITY, name="")}}
\end{quote}

\sphinxAtStartPar
\sphinxstylestrong{Description}
\begin{quote}

\sphinxAtStartPar
Add a constraint with a lower bound and an upper bound to a model and return the added {\hyperref[\detokenize{pyapiref:chappyapi-constraint}]{\sphinxcrossref{\DUrole{std,std-ref}{Constraint Class}}}}  object.
\end{quote}

\sphinxAtStartPar
\sphinxstylestrong{Arguments}
\begin{quote}

\sphinxAtStartPar
\sphinxcode{\sphinxupquote{expr}}
\begin{quote}

\sphinxAtStartPar
Expression for the new constraint, which can be {\hyperref[\detokenize{pyapiref:chappyapi-var}]{\sphinxcrossref{\DUrole{std,std-ref}{Var Class}}}} object or {\hyperref[\detokenize{pyapiref:chappyapi-linexpr}]{\sphinxcrossref{\DUrole{std,std-ref}{LinExpr Class}}}} object.
\end{quote}

\sphinxAtStartPar
\sphinxcode{\sphinxupquote{lb}}
\begin{quote}

\sphinxAtStartPar
Lower bound for the new constraint. Optional,  \sphinxcode{\sphinxupquote{\sphinxhyphen{}COPT.INFINITY}} by default.
\end{quote}

\sphinxAtStartPar
\sphinxcode{\sphinxupquote{ub}}
\begin{quote}

\sphinxAtStartPar
Upper bound for the new constraint. Optional, \sphinxcode{\sphinxupquote{COPT.INFINITY}} by default.
\end{quote}

\sphinxAtStartPar
\sphinxcode{\sphinxupquote{name}}
\begin{quote}

\sphinxAtStartPar
Name for new constraint. Optional, \sphinxcode{\sphinxupquote{""}} by default, automatically generated by solver.
\end{quote}
\end{quote}

\sphinxAtStartPar
\sphinxstylestrong{Example}

\begin{sphinxVerbatim}[commandchars=\\\{\}]
\PYG{c+c1}{\PYGZsh{} Add linear bilateral constraint: \PYGZhy{}1 \PYGZlt{}= x + y \PYGZlt{}= 1}
\PYG{n}{m}\PYG{o}{.}\PYG{n}{addBoundConstr}\PYG{p}{(}\PYG{n}{x} \PYG{o}{+} \PYG{n}{y}\PYG{p}{,} \PYG{o}{\PYGZhy{}}\PYG{l+m+mf}{1.0}\PYG{p}{,} \PYG{l+m+mf}{1.0}\PYG{p}{)}
\end{sphinxVerbatim}
\end{quote}

\subsubsection{Model.addConstrs()}
\label{\detokenize{pyapiref:model-addconstrs}}\begin{quote}

\sphinxAtStartPar
\sphinxstylestrong{Synopsis}
\begin{quote}

\sphinxAtStartPar
\sphinxcode{\sphinxupquote{addConstrs(generator, nameprefix="R")}}
\end{quote}

\sphinxAtStartPar
\sphinxstylestrong{Description}
\begin{quote}

\sphinxAtStartPar
Add a set of linear constraints, semidefinite constraints, or indicator constraints to the model.

\sphinxAtStartPar
If paramter \sphinxcode{\sphinxupquote{generator}} is integer, the return a {\hyperref[\detokenize{pyapiref:chappyapi-constrarray}]{\sphinxcrossref{\DUrole{std,std-ref}{ConstrArray Class}}}} object consisting of
\sphinxcode{\sphinxupquote{generator}} number of empty {\hyperref[\detokenize{pyapiref:chappyapi-constraint}]{\sphinxcrossref{\DUrole{std,std-ref}{Constraint Class}}}} objects, and users need to specify these constraints.

\sphinxAtStartPar
If parameter \sphinxcode{\sphinxupquote{generator}} is expression generator, then return a {\hyperref[\detokenize{pyapiref:chappyapi-util-tupledict}]{\sphinxcrossref{\DUrole{std,std-ref}{tupledict Class}}}} object whose key
is the indice of linear constraint and value is the corresponding {\hyperref[\detokenize{pyapiref:chappyapi-constraint}]{\sphinxcrossref{\DUrole{std,std-ref}{Constraint Class}}}} object. Every iteration generates a {\hyperref[\detokenize{pyapiref:chappyapi-constraint}]{\sphinxcrossref{\DUrole{std,std-ref}{Constraint Class}}}} object.

\sphinxAtStartPar
If the parameter \sphinxcode{\sphinxupquote{generator}} is a matrix expression generator,
return a {\hyperref[\detokenize{pyapiref:chappyapi-mconstr}]{\sphinxcrossref{\DUrole{std,std-ref}{MConstr Class}}}} object or a {\hyperref[\detokenize{pyapiref:chappyapi-mpsdconstr}]{\sphinxcrossref{\DUrole{std,std-ref}{MPsdConstr Class}}}} object.

\sphinxAtStartPar
If the parameter \sphinxcode{\sphinxupquote{generator}} is an indicator expression generator, return a {\hyperref[\detokenize{pyapiref:chappyapi-genconstrarray}]{\sphinxcrossref{\DUrole{std,std-ref}{GenConstrArray Class}}}} object.
\end{quote}

\sphinxAtStartPar
\sphinxstylestrong{Arguments}
\begin{quote}

\sphinxAtStartPar
\sphinxcode{\sphinxupquote{generator}}
\begin{quote}

\sphinxAtStartPar
A generator expression, where each iteration produces a  {\hyperref[\detokenize{pyapiref:chappyapi-constraint}]{\sphinxcrossref{\DUrole{std,std-ref}{Constraint Class}}}} object,
or a matrix expression builder, or an indicator expression generator.
\end{quote}

\sphinxAtStartPar
\sphinxcode{\sphinxupquote{nameprefix}}
\begin{quote}

\sphinxAtStartPar
Name prefix for new constraints. Optional, \sphinxcode{\sphinxupquote{"R"}} by default.
The actual name and the index of the constraints are automatically generated by COPT.
\end{quote}
\end{quote}

\sphinxAtStartPar
\sphinxstylestrong{Example}

\begin{sphinxVerbatim}[commandchars=\\\{\}]
\PYG{c+c1}{\PYGZsh{} Add 10 linear constraints, each constraint shaped like: x[0] + y[0] \PYGZgt{}= 2.0}
\PYG{n}{m}\PYG{o}{.}\PYG{n}{addConstrs}\PYG{p}{(}\PYG{n}{x}\PYG{p}{[}\PYG{n}{i}\PYG{p}{]} \PYG{o}{+} \PYG{n}{y}\PYG{p}{[}\PYG{n}{i}\PYG{p}{]} \PYG{o}{\PYGZgt{}}\PYG{o}{=} \PYG{l+m+mf}{2.0} \PYG{k}{for} \PYG{n}{i} \PYG{o+ow}{in} \PYG{n+nb}{range}\PYG{p}{(}\PYG{l+m+mi}{10}\PYG{p}{)}\PYG{p}{)}
\end{sphinxVerbatim}
\end{quote}

\subsubsection{Model.addMConstr()}
\label{\detokenize{pyapiref:model-addmconstr}}\begin{quote}

\sphinxAtStartPar
\sphinxstylestrong{Synopsis}
\begin{quote}

\sphinxAtStartPar
\sphinxcode{\sphinxupquote{addMConstr(A, x, sense, b, nameprefix="")}}
\end{quote}

\sphinxAtStartPar
\sphinxstylestrong{Description}
\begin{quote}

\sphinxAtStartPar
By means of matrix modeling, a set of linear constraints are added to the model. If the value of \sphinxcode{\sphinxupquote{sense}} here is \sphinxcode{\sphinxupquote{COPT.LESS\_EQUAL}} , the added constraint is \(Ax <= b\).

\sphinxAtStartPar
It is more convenient to generate {\hyperref[\detokenize{pyapiref:chappyapi-mlinexpr}]{\sphinxcrossref{\DUrole{std,std-ref}{MLinExpr Class}}}} objects by matrix multiplication, and then use overloaded comparison operators
to generate {\hyperref[\detokenize{pyapiref:chappyapi-mconstrbuilder}]{\sphinxcrossref{\DUrole{std,std-ref}{MConstrBuilder Class}}}} object, which can be used as input of \sphinxcode{\sphinxupquote{Model.addConstrs()}}
to generate a set of linear constraints.
\end{quote}

\sphinxAtStartPar
\sphinxstylestrong{Arguments}
\begin{quote}

\sphinxAtStartPar
\sphinxcode{\sphinxupquote{A}}
\begin{quote}

\sphinxAtStartPar
Parameter A is a two\sphinxhyphen{}dimensional NumPy matrix, SciPy compressed sparse column matrix ( \sphinxcode{\sphinxupquote{csc\_matrix}} ) or compressed sparse row matrix ( \sphinxcode{\sphinxupquote{csr\_matrix}} ).
\end{quote}

\sphinxAtStartPar
\sphinxcode{\sphinxupquote{x}}
\begin{quote}

\sphinxAtStartPar
The variable corresponding to the linear term can be a {\hyperref[\detokenize{pyapiref:chappyapi-mvar}]{\sphinxcrossref{\DUrole{std,std-ref}{MVar Class}}}} object, {\hyperref[\detokenize{pyapiref:chappyapi-vararray}]{\sphinxcrossref{\DUrole{std,std-ref}{VarArray Class}}}} object,
list, dictionary or {\hyperref[\detokenize{pyapiref:chappyapi-util-tupledict}]{\sphinxcrossref{\DUrole{std,std-ref}{tupledict Class}}}} object.
If it is empty, but the parameter c is not empty, all variables in the model are taken.
\end{quote}

\sphinxAtStartPar
\sphinxcode{\sphinxupquote{sense}}
\begin{quote}

\sphinxAtStartPar
The type of constraint. Possible values refer to {\hyperref[\detokenize{constant:chapconst-constrtype}]{\sphinxcrossref{\DUrole{std,std-ref}{Constraint type}}}}.
\end{quote}

\sphinxAtStartPar
\sphinxcode{\sphinxupquote{b}}
\begin{quote}

\sphinxAtStartPar
The value on the right side of the constraint, usually a floating point number, can also be a set of numbers, or a one\sphinxhyphen{}dimensional array of NumPy.
\end{quote}

\sphinxAtStartPar
\sphinxcode{\sphinxupquote{nameprefix}}
\begin{quote}

\sphinxAtStartPar
Constraint name prefix.
\end{quote}
\end{quote}

\sphinxAtStartPar
\sphinxstylestrong{Return value}
\begin{quote}

\sphinxAtStartPar
Returns a {\hyperref[\detokenize{pyapiref:chappyapi-mconstr}]{\sphinxcrossref{\DUrole{std,std-ref}{MConstr Class}}}} object
\end{quote}

\sphinxAtStartPar
\sphinxstylestrong{Example}

\begin{sphinxVerbatim}[commandchars=\\\{\}]
\PYG{n}{A} \PYG{o}{=} \PYG{n}{np}\PYG{o}{.}\PYG{n}{full}\PYG{p}{(}\PYG{p}{(}\PYG{l+m+mi}{2}\PYG{p}{,} \PYG{l+m+mi}{3}\PYG{p}{)}\PYG{p}{,} \PYG{l+m+mi}{1}\PYG{p}{)}
\PYG{n}{mx} \PYG{o}{=} \PYG{n}{model}\PYG{o}{.}\PYG{n}{addMVar}\PYG{p}{(}\PYG{l+m+mi}{3}\PYG{p}{,} \PYG{n}{nameprefix}\PYG{o}{=}\PYG{l+s+s2}{\PYGZdq{}}\PYG{l+s+s2}{mx}\PYG{l+s+s2}{\PYGZdq{}}\PYG{p}{)}
\PYG{n}{mc} \PYG{o}{=} \PYG{n}{model}\PYG{o}{.}\PYG{n}{addMConstr}\PYG{p}{(}\PYG{n}{A}\PYG{p}{,} \PYG{n}{mx}\PYG{p}{,} \PYG{l+s+s1}{\PYGZsq{}}\PYG{l+s+s1}{L}\PYG{l+s+s1}{\PYGZsq{}}\PYG{p}{,} \PYG{l+m+mf}{1.0}\PYG{p}{,} \PYG{n}{nameprefix}\PYG{o}{=}\PYG{l+s+s2}{\PYGZdq{}}\PYG{l+s+s2}{mc}\PYG{l+s+s2}{\PYGZdq{}}\PYG{p}{)}
\end{sphinxVerbatim}
\end{quote}

\subsubsection{Model.addSOS()}
\label{\detokenize{pyapiref:model-addsos}}\begin{quote}

\sphinxAtStartPar
\sphinxstylestrong{Synopsis}
\begin{quote}

\sphinxAtStartPar
\sphinxcode{\sphinxupquote{addSOS(sostype, vars, weights=None)}}
\end{quote}

\sphinxAtStartPar
\sphinxstylestrong{Description}
\begin{quote}

\sphinxAtStartPar
Add a SOS constraint to model and return the added {\hyperref[\detokenize{pyapiref:chappyapi-sos}]{\sphinxcrossref{\DUrole{std,std-ref}{SOS Class}}}} object.

\sphinxAtStartPar
If param \sphinxcode{\sphinxupquote{sostype}} is {\hyperref[\detokenize{pyapiref:chappyapi-sosbuilder}]{\sphinxcrossref{\DUrole{std,std-ref}{SOSBuilder Class}}}} object, then the values of param \sphinxcode{\sphinxupquote{vars}} and param\textasciigrave{}\textasciigrave{}weights\textasciigrave{}\textasciigrave{} will be ignored;

\sphinxAtStartPar
If param \sphinxcode{\sphinxupquote{sostype}} is SOS constraint type, it can be valued as {\hyperref[\detokenize{constant:chapconst-sostype}]{\sphinxcrossref{\DUrole{std,std-ref}{SOS\sphinxhyphen{}constraint types}}}},
then param \sphinxcode{\sphinxupquote{vars}} represents variables of SOS constraint, taking the value of {\hyperref[\detokenize{pyapiref:chappyapi-vararray}]{\sphinxcrossref{\DUrole{std,std-ref}{VarArray Class}}}} object,
list, dict or {\hyperref[\detokenize{pyapiref:chappyapi-util-tupledict}]{\sphinxcrossref{\DUrole{std,std-ref}{tupledict Class}}}} object;

\sphinxAtStartPar
If param \sphinxcode{\sphinxupquote{weights}} is \sphinxcode{\sphinxupquote{None}}, then variable weights of SOS constraints will be automatically generated by solver.
Otherwise take the input from user as weights, possible values can be list, dictionary or {\hyperref[\detokenize{pyapiref:chappyapi-util-tupledict}]{\sphinxcrossref{\DUrole{std,std-ref}{tupledict Class}}}} object.
\end{quote}

\sphinxAtStartPar
\sphinxstylestrong{Arguments}
\begin{quote}

\sphinxAtStartPar
\sphinxcode{\sphinxupquote{sostype}}
\begin{quote}

\sphinxAtStartPar
SOS contraint type or SOS constraint builder.
\end{quote}

\sphinxAtStartPar
\sphinxcode{\sphinxupquote{vars}}
\begin{quote}

\sphinxAtStartPar
Variables of SOS constraints.
\end{quote}

\sphinxAtStartPar
\sphinxcode{\sphinxupquote{weights}}
\begin{quote}

\sphinxAtStartPar
Weights of variables in SOS constraints, optional, \sphinxcode{\sphinxupquote{None}} by default.
\end{quote}
\end{quote}

\sphinxAtStartPar
\sphinxstylestrong{Example}
\end{quote}

\begin{sphinxVerbatim}[commandchars=\\\{\}]
\PYG{c+c1}{\PYGZsh{} Add an SOS1 constraint, including variable x and y, weights of 1 and 2.}
\PYG{n}{m}\PYG{o}{.}\PYG{n}{addSOS}\PYG{p}{(}\PYG{n}{COPT}\PYG{o}{.}\PYG{n}{SOS\PYGZus{}TYPE1}\PYG{p}{,} \PYG{p}{[}\PYG{n}{x}\PYG{p}{,} \PYG{n}{y}\PYG{p}{]}\PYG{p}{,} \PYG{p}{[}\PYG{l+m+mi}{1}\PYG{p}{,} \PYG{l+m+mi}{2}\PYG{p}{]}\PYG{p}{)}
\end{sphinxVerbatim}

\subsubsection{Model.addGenConstrIndicator()}
\label{\detokenize{pyapiref:model-addgenconstrindicator}}\begin{quote}

\sphinxAtStartPar
\sphinxstylestrong{Synopsis}
\begin{quote}

\sphinxAtStartPar
\sphinxcode{\sphinxupquote{addGenConstrIndicator(binvar, binval, lhs, sense=None, rhs=None, type=COPT.INDICATOR\_IF, name="")}}
\end{quote}

\sphinxAtStartPar
\sphinxstylestrong{Description}
\begin{quote}

\sphinxAtStartPar
Add an indicator constraint with the specified type to a model and return the added {\hyperref[\detokenize{pyapiref:chappyapi-genconstr}]{\sphinxcrossref{\DUrole{std,std-ref}{GenConstr Class}}}} object.

\sphinxAtStartPar
If the parameter \sphinxcode{\sphinxupquote{lhs}} is {\hyperref[\detokenize{pyapiref:chappyapi-constrbuilder}]{\sphinxcrossref{\DUrole{std,std-ref}{ConstrBuilder Class}}}} object, then the values of parameter \sphinxcode{\sphinxupquote{sense}} and
parameter \sphinxcode{\sphinxupquote{rhs}} will be ignored.

\sphinxAtStartPar
If parameter \sphinxcode{\sphinxupquote{lhs}} represents Left\sphinxhyphen{}hand side expression, it can take value of {\hyperref[\detokenize{pyapiref:chappyapi-var}]{\sphinxcrossref{\DUrole{std,std-ref}{Var Class}}}} object or
{\hyperref[\detokenize{pyapiref:chappyapi-linexpr}]{\sphinxcrossref{\DUrole{std,std-ref}{LinExpr Class}}}} object.
\end{quote}

\sphinxAtStartPar
\sphinxstylestrong{Arguments}
\begin{quote}

\sphinxAtStartPar
\sphinxcode{\sphinxupquote{binvar}}
\begin{quote}

\sphinxAtStartPar
Indicator variable.
\end{quote}

\sphinxAtStartPar
\sphinxcode{\sphinxupquote{binval}}
\begin{quote}

\sphinxAtStartPar
Value of indicator variable, can be \sphinxcode{\sphinxupquote{True}} or \sphinxcode{\sphinxupquote{False}}.
\end{quote}

\sphinxAtStartPar
\sphinxcode{\sphinxupquote{lhs}}
\begin{quote}

\sphinxAtStartPar
Left\sphinxhyphen{}hand side expression for the linear constraint triggered by the indicator or linear constraint builder.
\end{quote}

\sphinxAtStartPar
\sphinxcode{\sphinxupquote{sense}}
\begin{quote}

\sphinxAtStartPar
Sense for the linear constraint. Optional,  \sphinxcode{\sphinxupquote{None}} by default. Please refer to
{\hyperref[\detokenize{constant:chapconst-constrtype}]{\sphinxcrossref{\DUrole{std,std-ref}{Constraint type}}}} for possible values.
\end{quote}

\sphinxAtStartPar
\sphinxcode{\sphinxupquote{rhs}}
\begin{quote}

\sphinxAtStartPar
Right\sphinxhyphen{}hand\sphinxhyphen{}side value for the linear constraint triggered by the indicator. Optional,  \sphinxcode{\sphinxupquote{None}} by default,
value type is constant.
\end{quote}

\sphinxAtStartPar
\sphinxcode{\sphinxupquote{type}}
\begin{quote}

\sphinxAtStartPar
Type of the indicator constraint. Optional,  \sphinxcode{\sphinxupquote{COPT.INDICATOR\_IF}} by default.
Please refer to {\hyperref[\detokenize{constant:chapconst-indicatortype}]{\sphinxcrossref{\DUrole{std,std-ref}{Indicator Constraint type}}}} for possible values.
\end{quote}

\sphinxAtStartPar
\sphinxcode{\sphinxupquote{name}}
\begin{quote}

\sphinxAtStartPar
Name for new indicator constraint. Optional, \sphinxcode{\sphinxupquote{""}} by default, generated by solver automatically.
\end{quote}
\end{quote}

\sphinxAtStartPar
\sphinxstylestrong{Example}
\end{quote}

\begin{sphinxVerbatim}[commandchars=\\\{\}]
\PYG{c+c1}{\PYGZsh{} Add an indicator constraint, if x is True, then the linear constraint y + 2*z \PYGZgt{}= 3 should hold}
\PYG{n}{m}\PYG{o}{.}\PYG{n}{addGenConstrIndicator}\PYG{p}{(}\PYG{n}{x}\PYG{p}{,} \PYG{k+kc}{True}\PYG{p}{,} \PYG{n}{y} \PYG{o}{+} \PYG{l+m+mi}{2}\PYG{o}{*}\PYG{n}{z} \PYG{o}{\PYGZgt{}}\PYG{o}{=} \PYG{l+m+mi}{3}\PYG{p}{)}
\PYG{c+c1}{\PYGZsh{} Add an indicator constraint, if the linear constraint y + 2*z \PYGZgt{}= 3 hold, then x should be True}
\PYG{n}{m}\PYG{o}{.}\PYG{n}{addGenConstrIndicator}\PYG{p}{(}\PYG{n}{x}\PYG{p}{,} \PYG{k+kc}{True}\PYG{p}{,} \PYG{n}{y} \PYG{o}{+} \PYG{l+m+mi}{2}\PYG{o}{*}\PYG{n}{z} \PYG{o}{\PYGZgt{}}\PYG{o}{=} \PYG{l+m+mi}{3}\PYG{p}{,} \PYG{n+nb}{type}\PYG{o}{=}\PYG{n}{COPT}\PYG{o}{.}\PYG{n}{INDICATOR\PYGZus{}ONLYIF}\PYG{p}{)}
\end{sphinxVerbatim}

\subsubsection{Model.addGenConstrIndicators()}
\label{\detokenize{pyapiref:model-addgenconstrindicators}}\begin{quote}

\sphinxAtStartPar
\sphinxstylestrong{Synopsis}
\begin{quote}

\sphinxAtStartPar
\sphinxcode{\sphinxupquote{addGenConstrIndicators(builders, nameprefix="")}}
\end{quote}

\sphinxAtStartPar
\sphinxstylestrong{Description}
\begin{quote}

\sphinxAtStartPar
Add a set of indicator constraints with the specified type to a model and return the added {\hyperref[\detokenize{pyapiref:chappyapi-genconstrarray}]{\sphinxcrossref{\DUrole{std,std-ref}{GenConstrArray Class}}}} object.
\end{quote}

\sphinxAtStartPar
\sphinxstylestrong{Arguments}
\begin{quote}

\sphinxAtStartPar
\sphinxcode{\sphinxupquote{builders}}
\begin{quote}

\sphinxAtStartPar
A set of indicator constraint builders. The possible values could be {\hyperref[\detokenize{pyapiref:chappyapi-genconstrbuilderarray}]{\sphinxcrossref{\DUrole{std,std-ref}{GenConstrBuilderArray Class}}}} object,
a list or dictionary of {\hyperref[\detokenize{pyapiref:chappyapi-genconstrbuilder}]{\sphinxcrossref{\DUrole{std,std-ref}{GenConstrBuilder Class}}}} object.
\end{quote}

\sphinxAtStartPar
\sphinxcode{\sphinxupquote{nameprefix}}
\begin{quote}

\sphinxAtStartPar
Nameprefix of indicator variable. Optional,  \sphinxcode{\sphinxupquote{""}} by default, generated by solver automatically.
\end{quote}
\end{quote}
\end{quote}

\subsubsection{Model.addGenConstrMin()}
\label{\detokenize{pyapiref:model-addgenconstrmin}}\begin{quote}

\sphinxAtStartPar
\sphinxstylestrong{Synopsis}
\begin{quote}

\sphinxAtStartPar
\sphinxcode{\sphinxupquote{addGenConstrMin(resvar, vars, constant=None, name="")}}
\end{quote}

\sphinxAtStartPar
\sphinxstylestrong{Description}
\begin{quote}

\sphinxAtStartPar
Add a constraint of the form \(y=\min\{x_1, x_2, \cdots, x_n, c\}\) to the model.
\end{quote}

\sphinxAtStartPar
\sphinxstylestrong{Arguments}
\begin{quote}

\sphinxAtStartPar
\sphinxcode{\sphinxupquote{resvar}}
\begin{quote}

\sphinxAtStartPar
The term \sphinxcode{\sphinxupquote{y}} on the left side of the equation and can be an object of class \sphinxcode{\sphinxupquote{Var}} or \sphinxcode{\sphinxupquote{MVar}} .
\end{quote}

\sphinxAtStartPar
\sphinxcode{\sphinxupquote{vars}}
\begin{quote}

\sphinxAtStartPar
The variable of the \(\min\{\}\) function on the right side of the equation,

\sphinxAtStartPar
Possible values are \sphinxcode{\sphinxupquote{list}} class objects.
\end{quote}

\sphinxAtStartPar
\sphinxcode{\sphinxupquote{constant}}
\begin{quote}

\sphinxAtStartPar
The constant term in the \(\min\{\}\) function on the right side of the equation,

\sphinxAtStartPar
Optional parameter, the possible value is a floating number,

\sphinxAtStartPar
The default value is \sphinxcode{\sphinxupquote{None}}.
\end{quote}

\sphinxAtStartPar
\sphinxcode{\sphinxupquote{name}}
\begin{quote}

\sphinxAtStartPar
Constraint name, optional parameter, default value is \sphinxcode{\sphinxupquote{""}} .
\end{quote}
\end{quote}

\sphinxAtStartPar
\sphinxstylestrong{Return value}
\begin{quote}

\sphinxAtStartPar
It returns a \sphinxcode{\sphinxupquote{GenConstrX}} Class object.
\end{quote}
\end{quote}

\subsubsection{Model.addGenConstrMax()}
\label{\detokenize{pyapiref:model-addgenconstrmax}}\begin{quote}

\sphinxAtStartPar
\sphinxstylestrong{Synopsis}
\begin{quote}

\sphinxAtStartPar
\sphinxcode{\sphinxupquote{addGenConstrMax(resvar, vars, constant=None, name="")}}
\end{quote}

\sphinxAtStartPar
\sphinxstylestrong{Description}
\begin{quote}

\sphinxAtStartPar
Add a constraint of the form \(y=\max\{x_1, x_2, \cdots, x_n, c\}\) to the model.
\end{quote}

\sphinxAtStartPar
\sphinxstylestrong{Arguments}
\begin{quote}

\sphinxAtStartPar
\sphinxcode{\sphinxupquote{resvar}}
\begin{quote}

\sphinxAtStartPar
The term \sphinxcode{\sphinxupquote{y}} on the left side of the equation and can be an object of class \sphinxcode{\sphinxupquote{Var}} or \sphinxcode{\sphinxupquote{MVar}} .
\end{quote}

\sphinxAtStartPar
\sphinxcode{\sphinxupquote{vars}}
\begin{quote}

\sphinxAtStartPar
The variable of the \(\max\{\}\) function on the right side of the equation.

\sphinxAtStartPar
Possible values are \sphinxcode{\sphinxupquote{list}} class objects.
\end{quote}

\sphinxAtStartPar
\sphinxcode{\sphinxupquote{constant}}
\begin{quote}

\sphinxAtStartPar
The constant term in the \(\max\{\}\) function on the right side of the equation.

\sphinxAtStartPar
Optional parameter, the possible value is a floating number.

\sphinxAtStartPar
The default value is \sphinxcode{\sphinxupquote{None}} .
\end{quote}

\sphinxAtStartPar
\sphinxcode{\sphinxupquote{name}}
\begin{quote}

\sphinxAtStartPar
Constraint name, optional parameter, default value is \sphinxcode{\sphinxupquote{""}} .
\end{quote}
\end{quote}

\sphinxAtStartPar
\sphinxstylestrong{Return value}
\begin{quote}

\sphinxAtStartPar
It returns a \sphinxcode{\sphinxupquote{GenConstrX}} Class object.
\end{quote}
\end{quote}

\subsubsection{Model.addGenConstrAbs()}
\label{\detokenize{pyapiref:model-addgenconstrabs}}\begin{quote}

\sphinxAtStartPar
\sphinxstylestrong{Synopsis}
\begin{quote}

\sphinxAtStartPar
\sphinxcode{\sphinxupquote{addGenConstrAbs(resvar, argvar, name="")}}
\end{quote}

\sphinxAtStartPar
\sphinxstylestrong{Description}
\begin{quote}

\sphinxAtStartPar
Add a constraint of the form \(cy+d=|ax+b|\) to the model.
\end{quote}

\sphinxAtStartPar
\sphinxstylestrong{Arguments}
\begin{quote}

\sphinxAtStartPar
\sphinxcode{\sphinxupquote{resvar}}
\begin{quote}

\sphinxAtStartPar
\(cy+d\) , possible values are objects of class \sphinxcode{\sphinxupquote{Var}} / \sphinxcode{\sphinxupquote{MVar}} or class \sphinxcode{\sphinxupquote{LinExpr}} / \sphinxcode{\sphinxupquote{MLinExpr}} .
\end{quote}

\sphinxAtStartPar
\sphinxcode{\sphinxupquote{argvar}}
\begin{quote}

\sphinxAtStartPar
\(ax+b\) , the possible value is object of class \sphinxcode{\sphinxupquote{Var}} / \sphinxcode{\sphinxupquote{MVar}} or class \sphinxcode{\sphinxupquote{LinExpr}} / \sphinxcode{\sphinxupquote{MLinExpr}} .
\end{quote}

\sphinxAtStartPar
\sphinxcode{\sphinxupquote{name}}
\begin{quote}

\sphinxAtStartPar
Constraint name, optional parameter, default value is \sphinxcode{\sphinxupquote{""}} .
\end{quote}
\end{quote}

\sphinxAtStartPar
\sphinxstylestrong{Return value}
\begin{quote}

\sphinxAtStartPar
It returns a \sphinxcode{\sphinxupquote{GenConstrX}} Class object.
\end{quote}
\end{quote}

\subsubsection{Model.addGenConstrAnd()}
\label{\detokenize{pyapiref:model-addgenconstrand}}\begin{quote}

\sphinxAtStartPar
\sphinxstylestrong{Synopsis}
\begin{quote}

\sphinxAtStartPar
\sphinxcode{\sphinxupquote{addGenConstrAnd(resvar, vars, name="")}}
\end{quote}

\sphinxAtStartPar
\sphinxstylestrong{Description}
\begin{quote}

\sphinxAtStartPar
Add a logical \sphinxcode{\sphinxupquote{and}} constraint of the form \(y = x_1 \text{ and } x_2 \cdots \text{ and } x_n\) to the model.
\end{quote}

\sphinxAtStartPar
\sphinxstylestrong{Arguments}
\begin{quote}

\sphinxAtStartPar
\sphinxcode{\sphinxupquote{resvar}}
\begin{quote}

\sphinxAtStartPar
The term \sphinxcode{\sphinxupquote{y}} on the left side of the equation and can be an object of class \sphinxcode{\sphinxupquote{Var}} or \sphinxcode{\sphinxupquote{MVar}} .
\end{quote}

\sphinxAtStartPar
\sphinxcode{\sphinxupquote{vars}}
\begin{quote}

\sphinxAtStartPar
Elements connected by logical operator \sphinxcode{\sphinxupquote{and}} \(x_i, \text{for } i \in \{1,2,\cdots,n\}\)

\sphinxAtStartPar
Possible values are \sphinxcode{\sphinxupquote{List}} class (where the elements are binary \sphinxcode{\sphinxupquote{Var}} class or \sphinxcode{\sphinxupquote{MVar}} class objects).
\end{quote}

\sphinxAtStartPar
\sphinxcode{\sphinxupquote{name}}
\begin{quote}

\sphinxAtStartPar
Constraint name, optional parameter, default value is \sphinxcode{\sphinxupquote{""}} .
\end{quote}
\end{quote}

\sphinxAtStartPar
\sphinxstylestrong{Return value}
\begin{quote}

\sphinxAtStartPar
It returns a \sphinxcode{\sphinxupquote{GenConstrX}} Class object.
\end{quote}
\end{quote}

\subsubsection{Model.addGenConstrOr()}
\label{\detokenize{pyapiref:model-addgenconstror}}\begin{quote}

\sphinxAtStartPar
\sphinxstylestrong{Synopsis}
\begin{quote}

\sphinxAtStartPar
\sphinxcode{\sphinxupquote{addGenConstrOr(resvar, vars, name="")}}
\end{quote}

\sphinxAtStartPar
\sphinxstylestrong{Description}
\begin{quote}

\sphinxAtStartPar
Add a logical \sphinxcode{\sphinxupquote{or}} constraint of the form \(y = x_1 \text{ or } x_2 \cdots \text{ or } x_n\) to the model.
\end{quote}

\sphinxAtStartPar
\sphinxstylestrong{Arguments}
\begin{quote}

\sphinxAtStartPar
\sphinxcode{\sphinxupquote{resvar}}
\begin{quote}

\sphinxAtStartPar
The term \sphinxcode{\sphinxupquote{y}} on the left side of the equation and can be an object of class \sphinxcode{\sphinxupquote{Var}} or \sphinxcode{\sphinxupquote{MVar}} .
\end{quote}

\sphinxAtStartPar
\sphinxcode{\sphinxupquote{vars}}
\begin{quote}

\sphinxAtStartPar
Elements connected by logical operator \sphinxcode{\sphinxupquote{or}} \(x_i, \text{for } i \in \{1,2,\cdots,n\}\)

\sphinxAtStartPar
Possible values are \sphinxcode{\sphinxupquote{List}} class (where the elements are binary \sphinxcode{\sphinxupquote{Var}} class or \sphinxcode{\sphinxupquote{MVar}} objects).
\end{quote}

\sphinxAtStartPar
\sphinxcode{\sphinxupquote{name}}
\begin{quote}

\sphinxAtStartPar
Constraint name, optional parameter, default value is \sphinxcode{\sphinxupquote{""}} .
\end{quote}
\end{quote}

\sphinxAtStartPar
\sphinxstylestrong{Return value}
\begin{quote}

\sphinxAtStartPar
It returns a \sphinxcode{\sphinxupquote{GenConstrX}} Class object.
\end{quote}
\end{quote}

\subsubsection{Model.addGenConstrPWL()}
\label{\detokenize{pyapiref:model-addgenconstrpwl}}\begin{quote}

\sphinxAtStartPar
\sphinxstylestrong{Synopsis}
\begin{quote}

\sphinxAtStartPar
\sphinxcode{\sphinxupquote{addGenConstrPWL(xvar, yvar, xpts, ypts, name="")}}
\end{quote}

\sphinxAtStartPar
\sphinxstylestrong{Description}
\begin{quote}

\sphinxAtStartPar
Add a constraint of the form \(y=f(x)\), where a piecewise linear function is defined as:
\begin{equation}\label{equation:pyapiref:pyapiref:0}
\begin{split}f(v) = \begin{cases}
\tilde{y}_1 + \frac{\tilde{y}_2-\tilde{y}_1}{\tilde{x}_2-\tilde{x}_1} (v-\tilde{x}_1),\quad &\text{if } v\leq x_1 \\
\tilde{y}_i + \frac{\tilde{y}_{i+1} - \tilde{y}_i}{\tilde{x}_{i+1} - \tilde{x}_i} (v-\tilde{x}_i),\quad &\text{if } \tilde{x}_i\leq v\leq \tilde{x}_{i+1} \\
\tilde{y}_n + \frac{\tilde{y}_n - \tilde{y}_{n-1}}{\tilde{x}_n - \tilde{x}_{n-1}}(v-\tilde{x}_n),\quad &\text{if } v\geq \tilde{x}_n
\end{cases}\notag\end{split}
\end{equation}\end{quote}

\sphinxAtStartPar
\sphinxstylestrong{Arguments}
\begin{quote}

\sphinxAtStartPar
\sphinxcode{\sphinxupquote{xvar}}
\begin{quote}

\sphinxAtStartPar
\sphinxcode{\sphinxupquote{x}}, which can be an object of class \sphinxcode{\sphinxupquote{Var}} or class \sphinxcode{\sphinxupquote{MVar}} .
\end{quote}

\sphinxAtStartPar
\sphinxcode{\sphinxupquote{yvar}}
\begin{quote}

\sphinxAtStartPar
The term \sphinxcode{\sphinxupquote{y}} on the left side of the equation,

\sphinxAtStartPar
Possible values are objects of class \sphinxcode{\sphinxupquote{Var}} / \sphinxcode{\sphinxupquote{MVar}} or class \sphinxcode{\sphinxupquote{LinExpr}} / \sphinxcode{\sphinxupquote{MLinExpr}} .
\end{quote}

\sphinxAtStartPar
\sphinxcode{\sphinxupquote{xpts}}
\begin{quote}

\sphinxAtStartPar
\(\tilde{\boldsymbol{x}}\), the abscissa of the segmentation point.

\sphinxAtStartPar
It should be arranged in ascending order of values, possible values are \sphinxcode{\sphinxupquote{List}} class.
\end{quote}

\sphinxAtStartPar
\sphinxcode{\sphinxupquote{ypts}}
\begin{quote}

\sphinxAtStartPar
\(\tilde{\boldsymbol{y}}\) , the vertical coordinate of the segmentation point,

\sphinxAtStartPar
Possible values are \sphinxcode{\sphinxupquote{List}} class.
\end{quote}

\sphinxAtStartPar
\sphinxcode{\sphinxupquote{name}}
\begin{quote}

\sphinxAtStartPar
Constraint name, optional parameter, default value is \sphinxcode{\sphinxupquote{""}} .
\end{quote}
\end{quote}

\sphinxAtStartPar
\sphinxstylestrong{Return value}
\begin{quote}

\sphinxAtStartPar
It returns a \sphinxcode{\sphinxupquote{GenConstrX}} Class object.
\end{quote}
\end{quote}

\subsubsection{Model.addConeByDim()}
\label{\detokenize{pyapiref:model-addconebydim}}\begin{quote}

\sphinxAtStartPar
\sphinxstylestrong{Synopsis}
\begin{quote}

\sphinxAtStartPar
\sphinxcode{\sphinxupquote{addConeByDim(dim, ctype, vtype, nameprefix="ConeV")}}
\end{quote}

\sphinxAtStartPar
\sphinxstylestrong{Description}
\begin{quote}

\sphinxAtStartPar
Add a Second\sphinxhyphen{}Order\sphinxhyphen{}Cone (SOC) constraint with given dimension, and return the added
{\hyperref[\detokenize{pyapiref:chappyapi-cone}]{\sphinxcrossref{\DUrole{std,std-ref}{Cone Class}}}} object.
\end{quote}

\sphinxAtStartPar
\sphinxstylestrong{Arguments}
\begin{quote}

\sphinxAtStartPar
\sphinxcode{\sphinxupquote{dim}}
\begin{quote}

\sphinxAtStartPar
Dimension of SOC constraint.
\end{quote}

\sphinxAtStartPar
\sphinxcode{\sphinxupquote{ctype}}
\begin{quote}

\sphinxAtStartPar
Type of SOC constraint.
\end{quote}

\sphinxAtStartPar
\sphinxcode{\sphinxupquote{vtype}}
\begin{quote}

\sphinxAtStartPar
Variable types of SOC constraint.
\end{quote}

\sphinxAtStartPar
\sphinxcode{\sphinxupquote{nameprefix}}
\begin{quote}

\sphinxAtStartPar
Name prefix of variables in SOC constraint. Optional, default to \sphinxcode{\sphinxupquote{"ConeV"}}.
\end{quote}
\end{quote}

\sphinxAtStartPar
\sphinxstylestrong{Example}
\end{quote}

\begin{sphinxVerbatim}[commandchars=\\\{\}]
\PYG{c+c1}{\PYGZsh{} Add a 5 dimension rotated SOC constraint}
\PYG{n}{m}\PYG{o}{.}\PYG{n}{addConeByDim}\PYG{p}{(}\PYG{l+m+mi}{5}\PYG{p}{,} \PYG{n}{COPT}\PYG{o}{.}\PYG{n}{CONE\PYGZus{}RQUAD}\PYG{p}{,} \PYG{k+kc}{None}\PYG{p}{)}
\end{sphinxVerbatim}

\subsubsection{Model.addExpConeByDim()}
\label{\detokenize{pyapiref:model-addexpconebydim}}\begin{quote}

\sphinxAtStartPar
\sphinxstylestrong{Synopsis}
\begin{quote}

\sphinxAtStartPar
\sphinxcode{\sphinxupquote{addExpConeByDim(ctype, vtype, nameprefix="ExpConeV")}}
\end{quote}

\sphinxAtStartPar
\sphinxstylestrong{Description}
\begin{quote}

\sphinxAtStartPar
Add an exponential cone constraint and return the added
{\hyperref[\detokenize{pyapiref:chappyapi-expcone}]{\sphinxcrossref{\DUrole{std,std-ref}{ExpCone Class}}}} object.
\end{quote}

\sphinxAtStartPar
\sphinxstylestrong{Arguments}
\begin{quote}

\sphinxAtStartPar
\sphinxcode{\sphinxupquote{ctype}}
\begin{quote}

\sphinxAtStartPar
Type of exponential cone constraint.
\end{quote}

\sphinxAtStartPar
\sphinxcode{\sphinxupquote{vtype}}
\begin{quote}

\sphinxAtStartPar
Variable types of exponential cone constraint.
\end{quote}

\sphinxAtStartPar
\sphinxcode{\sphinxupquote{nameprefix}}
\begin{quote}

\sphinxAtStartPar
Name prefix of variables in exponential cone constraint. Optional, default to \sphinxcode{\sphinxupquote{"ExpConeV"}} .
\end{quote}
\end{quote}

\sphinxAtStartPar
\sphinxstylestrong{Example}
\end{quote}

\begin{sphinxVerbatim}[commandchars=\\\{\}]
\PYG{c+c1}{\PYGZsh{} Add a primal exponential cone}
\PYG{n}{m}\PYG{o}{.}\PYG{n}{addExpConeByDim}\PYG{p}{(}\PYG{n}{COPT}\PYG{o}{.}\PYG{n}{EXPCONE\PYGZus{}PRIMAL}\PYG{p}{,} \PYG{k+kc}{None}\PYG{p}{)}
\end{sphinxVerbatim}

\subsubsection{Model.addCone()}
\label{\detokenize{pyapiref:model-addcone}}\begin{quote}

\sphinxAtStartPar
\sphinxstylestrong{Synopsis}
\begin{quote}

\sphinxAtStartPar
\sphinxcode{\sphinxupquote{addCone(vars, ctype)}}
\end{quote}

\sphinxAtStartPar
\sphinxstylestrong{Description}
\begin{quote}

\sphinxAtStartPar
Add a Second\sphinxhyphen{}Order\sphinxhyphen{}Cone (SOC) constraint with given variables.

\sphinxAtStartPar
If argument \sphinxcode{\sphinxupquote{vars}} is a {\hyperref[\detokenize{pyapiref:chappyapi-conebuilder}]{\sphinxcrossref{\DUrole{std,std-ref}{ConeBuilder Class}}}} object, then the value of
argument \sphinxcode{\sphinxupquote{ctype}} will be ignored;
If argument \sphinxcode{\sphinxupquote{vars}} are variables, the optional values are {\hyperref[\detokenize{pyapiref:chappyapi-vararray}]{\sphinxcrossref{\DUrole{std,std-ref}{VarArray Class}}}} objects,
Python list, Python dictionary or {\hyperref[\detokenize{pyapiref:chappyapi-util-tupledict}]{\sphinxcrossref{\DUrole{std,std-ref}{tupledict Class}}}} objects,
argument \sphinxcode{\sphinxupquote{ctype}} is the type of SOC constraint.
\end{quote}

\sphinxAtStartPar
\sphinxstylestrong{Arguments}
\begin{quote}

\sphinxAtStartPar
\sphinxcode{\sphinxupquote{vars}}
\begin{quote}

\sphinxAtStartPar
Variables of SOC constraint.
\end{quote}

\sphinxAtStartPar
\sphinxcode{\sphinxupquote{ctype}}
\begin{quote}

\sphinxAtStartPar
Type of SOC constraint. Please refer {\hyperref[\detokenize{constant:chapconst-conetype}]{\sphinxcrossref{\DUrole{std,std-ref}{SOC constraint types}}}} for possible values.
\end{quote}
\end{quote}

\sphinxAtStartPar
\sphinxstylestrong{Example}
\end{quote}

\begin{sphinxVerbatim}[commandchars=\\\{\}]
\PYG{c+c1}{\PYGZsh{} Add a SOC constraint with [z, x, y] as variables}
\PYG{n}{m}\PYG{o}{.}\PYG{n}{addCone}\PYG{p}{(}\PYG{p}{[}\PYG{n}{z}\PYG{p}{,} \PYG{n}{x}\PYG{p}{,} \PYG{n}{y}\PYG{p}{]}\PYG{p}{,} \PYG{n}{COPT}\PYG{o}{.}\PYG{n}{CONE\PYGZus{}QUAD}\PYG{p}{)}
\end{sphinxVerbatim}

\subsubsection{Model.addCones()}
\label{\detokenize{pyapiref:model-addcones}}\begin{quote}

\sphinxAtStartPar
\sphinxstylestrong{Synopsis}
\begin{quote}

\sphinxAtStartPar
\sphinxcode{\sphinxupquote{addCones(vars, ctype)}}
\end{quote}

\sphinxAtStartPar
\sphinxstylestrong{Description}
\begin{quote}

\sphinxAtStartPar
Add a set of Second\sphinxhyphen{}Order\sphinxhyphen{}Cone (SOC) constraint with given variables.

\sphinxAtStartPar
If argument \sphinxcode{\sphinxupquote{vars}} is a {\hyperref[\detokenize{pyapiref:chappyapi-conebuilder}]{\sphinxcrossref{\DUrole{std,std-ref}{ConeBuilder Class}}}} object or {\hyperref[\detokenize{pyapiref:chappyapi-conebuilderarray}]{\sphinxcrossref{\DUrole{std,std-ref}{ConeBuilderArray Class}}}}
object, then the value of argument \sphinxcode{\sphinxupquote{ctype}} will be ignored;
If argument \sphinxcode{\sphinxupquote{vars}} are {\hyperref[\detokenize{pyapiref:chappyapi-mvar}]{\sphinxcrossref{\DUrole{std,std-ref}{MVar Class}}}} , the argument \sphinxcode{\sphinxupquote{ctype}} is the type of SOC constraint,
and it must be specified.
\end{quote}

\sphinxAtStartPar
\sphinxstylestrong{Arguments}
\begin{quote}

\sphinxAtStartPar
\sphinxcode{\sphinxupquote{vars}}
\begin{quote}

\sphinxAtStartPar
Variables of SOC constraint.
\end{quote}

\sphinxAtStartPar
\sphinxcode{\sphinxupquote{ctype}}
\begin{quote}

\sphinxAtStartPar
Type of SOC constraint. Please refer {\hyperref[\detokenize{constant:chapconst-conetype}]{\sphinxcrossref{\DUrole{std,std-ref}{SOC constraint types}}}} for possible values.
\end{quote}
\end{quote}
\end{quote}

\subsubsection{Model.addExpCone()}
\label{\detokenize{pyapiref:model-addexpcone}}\begin{quote}

\sphinxAtStartPar
\sphinxstylestrong{Synopsis}
\begin{quote}

\sphinxAtStartPar
\sphinxcode{\sphinxupquote{addExpCone(vars, ctype)}}
\end{quote}

\sphinxAtStartPar
\sphinxstylestrong{Description}
\begin{quote}

\sphinxAtStartPar
Add an exponential cone to the model.

\sphinxAtStartPar
For example, the primal exponential cone defined by
the \(\boldsymbol{x}\in \mathbb{R}^3\) :
\begin{equation}\label{equation:pyapiref:pyapiref:1}
\begin{split}\begin{align*}
K_{exp}= \mathrm{cl}\left\{ \boldsymbol{x}\in \mathbb{R}^3 \mid x_0 \geq x_1\ \mathrm{exp}\left(\frac{x_2}{x_1} \right),\ x_1 > 0 \right\}
\end{align*}\end{split}
\end{equation}
\sphinxAtStartPar
Please refer to {\hyperref[\detokenize{constant:chapconst-expconetype}]{\sphinxcrossref{\DUrole{std,std-ref}{Exponential Cone type}}}} for more details.

\sphinxAtStartPar
If argument \sphinxcode{\sphinxupquote{vars}} is a {\hyperref[\detokenize{pyapiref:chappyapi-expconebuilder}]{\sphinxcrossref{\DUrole{std,std-ref}{ExpConeBuilder Class}}}} object, then the value of
argument \sphinxcode{\sphinxupquote{ctype}} will be ignored.

\sphinxAtStartPar
If argument \sphinxcode{\sphinxupquote{vars}} are variables, the optional values are
{\hyperref[\detokenize{pyapiref:chappyapi-vararray}]{\sphinxcrossref{\DUrole{std,std-ref}{VarArray Class}}}} objects, a one\sphinxhyphen{}dimensional {\hyperref[\detokenize{pyapiref:chappyapi-mvar}]{\sphinxcrossref{\DUrole{std,std-ref}{MVar Class}}}} object,
{\hyperref[\detokenize{pyapiref:chappyapi-util-tupledict}]{\sphinxcrossref{\DUrole{std,std-ref}{tupledict Class}}}} object, dictionary or list,
the argument \sphinxcode{\sphinxupquote{ctype}}  specifies the type of the exponential cone and must be explicitly provided.
\end{quote}

\sphinxAtStartPar
\sphinxstylestrong{Arguments}
\begin{quote}

\sphinxAtStartPar
\sphinxcode{\sphinxupquote{vars}}
\begin{quote}

\sphinxAtStartPar
The variables forming the exponential cone.

\sphinxAtStartPar
Possible values are {\hyperref[\detokenize{pyapiref:chappyapi-expconebuilder}]{\sphinxcrossref{\DUrole{std,std-ref}{ExpConeBuilder Class}}}} object,
{\hyperref[\detokenize{pyapiref:chappyapi-vararray}]{\sphinxcrossref{\DUrole{std,std-ref}{VarArray Class}}}} object, one\sphinxhyphen{}dimensional {\hyperref[\detokenize{pyapiref:chappyapi-mvar}]{\sphinxcrossref{\DUrole{std,std-ref}{MVar Class}}}} object,
{\hyperref[\detokenize{pyapiref:chappyapi-util-tupledict}]{\sphinxcrossref{\DUrole{std,std-ref}{tupledict Class}}}} object, dictionary or list.
\end{quote}

\sphinxAtStartPar
\sphinxcode{\sphinxupquote{ctype}}
\begin{quote}

\sphinxAtStartPar
Type of the exponential cone.
Please refer to {\hyperref[\detokenize{constant:chapconst-expconetype}]{\sphinxcrossref{\DUrole{std,std-ref}{Exponential Cone type}}}} for possible values.
\end{quote}
\end{quote}

\sphinxAtStartPar
\sphinxstylestrong{Example}
\end{quote}

\begin{sphinxVerbatim}[commandchars=\\\{\}]
\PYG{c+c1}{\PYGZsh{} Add a primal exponential cone formed by mx}
\PYG{n}{mx} \PYG{o}{=} \PYG{n}{m}\PYG{o}{.}\PYG{n}{addMVar}\PYG{p}{(}\PYG{l+m+mi}{3}\PYG{p}{)}
\PYG{n}{model}\PYG{o}{.}\PYG{n}{addExpCone}\PYG{p}{(}\PYG{n}{mx}\PYG{p}{,} \PYG{n}{ctype}\PYG{o}{=}\PYG{n}{COPT}\PYG{o}{.}\PYG{n}{EXPCONE\PYGZus{}PRIMAL}\PYG{p}{)}
\end{sphinxVerbatim}

\subsubsection{Model.addExpCones()}
\label{\detokenize{pyapiref:model-addexpcones}}\begin{quote}

\sphinxAtStartPar
\sphinxstylestrong{Synopsis}
\begin{quote}

\sphinxAtStartPar
\sphinxcode{\sphinxupquote{addExpCones(vars, ctype)}}
\end{quote}

\sphinxAtStartPar
\sphinxstylestrong{Description}
\begin{quote}

\sphinxAtStartPar
Add a batch of exponential cones to the model.
Return a {\hyperref[\detokenize{pyapiref:chappyapi-expconearray}]{\sphinxcrossref{\DUrole{std,std-ref}{ExpConeArray Class}}}} object.

\sphinxAtStartPar
If argument \sphinxcode{\sphinxupquote{vars}} is a {\hyperref[\detokenize{pyapiref:chappyapi-expconebuilder}]{\sphinxcrossref{\DUrole{std,std-ref}{ExpConeBuilder Class}}}} object,
then the value of argument \sphinxcode{\sphinxupquote{ctype}} will be ignored.

\sphinxAtStartPar
If argument \sphinxcode{\sphinxupquote{vars}} is a {\hyperref[\detokenize{pyapiref:chappyapi-mvar}]{\sphinxcrossref{\DUrole{std,std-ref}{MVar Class}}}} object,
the \sphinxcode{\sphinxupquote{ctype}} argument specifies the type of the exponential cones
and must be explicitly provided.
\end{quote}

\sphinxAtStartPar
\sphinxstylestrong{Arguments}
\begin{quote}

\sphinxAtStartPar
\sphinxcode{\sphinxupquote{vars}}
\begin{quote}

\sphinxAtStartPar
Variables forming the exponential cones. Possible values are
{\hyperref[\detokenize{pyapiref:chappyapi-expconebuilderarray}]{\sphinxcrossref{\DUrole{std,std-ref}{ExpConeBuilderArray Class}}}} object, 2\sphinxhyphen{}dimensional {\hyperref[\detokenize{pyapiref:chappyapi-mvar}]{\sphinxcrossref{\DUrole{std,std-ref}{MVar Class}}}} object.
\end{quote}

\sphinxAtStartPar
\sphinxcode{\sphinxupquote{ctype}}
\begin{quote}

\sphinxAtStartPar
Type of the exponential cone.
Please refer to {\hyperref[\detokenize{constant:chapconst-expconetype}]{\sphinxcrossref{\DUrole{std,std-ref}{Exponential Cone type}}}} for possible values.
\end{quote}
\end{quote}

\sphinxAtStartPar
\sphinxstylestrong{Example}
\end{quote}

\begin{sphinxVerbatim}[commandchars=\\\{\}]
\PYG{c+c1}{\PYGZsh{} Add two primal exponential cones formed by my}
\PYG{n}{model}\PYG{o}{.}\PYG{n}{matrixmodelmode} \PYG{o}{=} \PYG{l+s+s1}{\PYGZsq{}}\PYG{l+s+s1}{experimental}\PYG{l+s+s1}{\PYGZsq{}}
\PYG{n}{my} \PYG{o}{=} \PYG{n}{model}\PYG{o}{.}\PYG{n}{addMVar}\PYG{p}{(}\PYG{n}{shape}\PYG{o}{=}\PYG{p}{(}\PYG{l+m+mi}{2}\PYG{p}{,} \PYG{l+m+mi}{3}\PYG{p}{)}\PYG{p}{)}
\PYG{n}{model}\PYG{o}{.}\PYG{n}{addExpCones}\PYG{p}{(}\PYG{n}{my}\PYG{p}{,} \PYG{n}{COPT}\PYG{o}{.}\PYG{n}{EXPCONE\PYGZus{}PRIMAL}\PYG{p}{)}
\end{sphinxVerbatim}

\subsubsection{Model.addAffineCone()}
\label{\detokenize{pyapiref:model-addaffinecone}}\begin{quote}

\sphinxAtStartPar
\sphinxstylestrong{Synopsis}
\begin{quote}

\sphinxAtStartPar
\sphinxcode{\sphinxupquote{addAffineCone(exprs, ctype=None, name="")}}
\end{quote}

\sphinxAtStartPar
\sphinxstylestrong{Description}
\begin{quote}

\sphinxAtStartPar
Adds an affine cone to the model.

\sphinxAtStartPar
If the argument \sphinxcode{\sphinxupquote{exprs}} is a {\hyperref[\detokenize{pyapiref:chappyapi-affineconebuilder}]{\sphinxcrossref{\DUrole{std,std-ref}{AffineConeBuilder Class}}}} object, the value of the argument \sphinxcode{\sphinxupquote{ctype}} will be ignored.

\sphinxAtStartPar
If the argument \sphinxcode{\sphinxupquote{exprs}} is an {\hyperref[\detokenize{pyapiref:chappyapi-mlinexpr}]{\sphinxcrossref{\DUrole{std,std-ref}{MLinExpr Class}}}} object or an {\hyperref[\detokenize{pyapiref:chappyapi-mpsdexpr}]{\sphinxcrossref{\DUrole{std,std-ref}{MPsdExpr Class}}}} object, the argument \sphinxcode{\sphinxupquote{ctype}} specifies the type of the affine cone and must be explicitly provided.
\end{quote}

\sphinxAtStartPar
\sphinxstylestrong{Arguments}
\begin{quote}

\sphinxAtStartPar
\sphinxcode{\sphinxupquote{exprs}}
\begin{quote}

\sphinxAtStartPar
An affine cone generator or a multi\sphinxhyphen{}dimensional array expression that forms the affine cone.
\end{quote}

\sphinxAtStartPar
\sphinxcode{\sphinxupquote{ctype}}
\begin{quote}

\sphinxAtStartPar
The type of the affine cone. Possible values are detailed in {\hyperref[\detokenize{constant:chapconst-conetype}]{\sphinxcrossref{\DUrole{std,std-ref}{SOC constraint types}}}} or {\hyperref[\detokenize{constant:chapconst-expconetype}]{\sphinxcrossref{\DUrole{std,std-ref}{Exponential Cone type}}}}.
\end{quote}
\end{quote}

\sphinxAtStartPar
\sphinxstylestrong{Example}
\end{quote}

\begin{sphinxVerbatim}[commandchars=\\\{\}]
\PYG{c+c1}{\PYGZsh{} Add a standard second\PYGZhy{}order affine cone formed by Ax+b and c\PYGZca{}Tx+d}
\PYG{n}{A} \PYG{o}{=} \PYG{n}{np}\PYG{o}{.}\PYG{n}{array}\PYG{p}{(}\PYG{p}{[}\PYG{p}{[}\PYG{l+m+mi}{1}\PYG{p}{,} \PYG{l+m+mi}{2}\PYG{p}{,} \PYG{l+m+mi}{3}\PYG{p}{]}\PYG{p}{,}
             \PYG{p}{[}\PYG{l+m+mi}{4}\PYG{p}{,} \PYG{l+m+mi}{5}\PYG{p}{,} \PYG{l+m+mi}{6}\PYG{p}{]}\PYG{p}{,}
             \PYG{p}{[}\PYG{l+m+mi}{7}\PYG{p}{,} \PYG{l+m+mi}{8}\PYG{p}{,} \PYG{l+m+mi}{9}\PYG{p}{]}\PYG{p}{]}\PYG{p}{)}
\PYG{n}{b} \PYG{o}{=} \PYG{l+m+mi}{10}
\PYG{n}{c} \PYG{o}{=} \PYG{n}{np}\PYG{o}{.}\PYG{n}{array}\PYG{p}{(}\PYG{p}{[}\PYG{l+m+mi}{1}\PYG{p}{,} \PYG{o}{\PYGZhy{}}\PYG{l+m+mi}{1}\PYG{p}{,} \PYG{l+m+mi}{2}\PYG{p}{]}\PYG{p}{)}
\PYG{n}{d} \PYG{o}{=} \PYG{l+m+mi}{5}
\PYG{n}{x} \PYG{o}{=} \PYG{n}{model}\PYG{o}{.}\PYG{n}{addMVar}\PYG{p}{(}\PYG{l+m+mi}{3}\PYG{p}{)}
\PYG{n}{model}\PYG{o}{.}\PYG{n}{addAffineCone}\PYG{p}{(}\PYG{n}{cp}\PYG{o}{.}\PYG{n}{vstack}\PYG{p}{(}\PYG{n}{c} \PYG{o}{@} \PYG{n}{x} \PYG{o}{+} \PYG{n}{d}\PYG{p}{,} \PYG{n}{A} \PYG{o}{@} \PYG{n}{x} \PYG{o}{+} \PYG{n}{b}\PYG{p}{)}\PYG{p}{,} \PYG{n}{ctype}\PYG{o}{=}\PYG{n}{COPT}\PYG{o}{.}\PYG{n}{CONE\PYGZus{}QUAD}\PYG{p}{)}
\end{sphinxVerbatim}

\subsubsection{Model.addAffineCones()}
\label{\detokenize{pyapiref:model-addaffinecones}}\begin{quote}

\sphinxAtStartPar
\sphinxstylestrong{Synopsis}
\begin{quote}

\sphinxAtStartPar
\sphinxcode{\sphinxupquote{addAffineCones(exprs, ctype=None, nameprefix="AffineConV")}}
\end{quote}

\sphinxAtStartPar
\sphinxstylestrong{Description}
\begin{quote}

\sphinxAtStartPar
Adds a set of affine cones to the model.

\sphinxAtStartPar
If the argument \sphinxcode{\sphinxupquote{exprs}} is a {\hyperref[\detokenize{pyapiref:chappyapi-affineconebuilder}]{\sphinxcrossref{\DUrole{std,std-ref}{AffineConeBuilder Class}}}} object or a {\hyperref[\detokenize{pyapiref:chappyapi-affineconebuilderarray}]{\sphinxcrossref{\DUrole{std,std-ref}{AffineConeBuilderArray Class}}}} object, the value of the argument \sphinxcode{\sphinxupquote{ctype}} is ignored.

\sphinxAtStartPar
If the argument \sphinxcode{\sphinxupquote{exprs}} is an {\hyperref[\detokenize{pyapiref:chappyapi-mlinexpr}]{\sphinxcrossref{\DUrole{std,std-ref}{MLinExpr Class}}}} object or an {\hyperref[\detokenize{pyapiref:chappyapi-mpsdexpr}]{\sphinxcrossref{\DUrole{std,std-ref}{MPsdExpr Class}}}} object, the argument \sphinxcode{\sphinxupquote{ctype}} specifies the type of the affine cones and must be explicitly provided.
\end{quote}

\sphinxAtStartPar
\sphinxstylestrong{Arguments}
\begin{quote}

\sphinxAtStartPar
\sphinxcode{\sphinxupquote{exprs}}
\begin{quote}

\sphinxAtStartPar
Affine cone generators or multi\sphinxhyphen{}dimensional array expressions that form the affine cones.
\end{quote}

\sphinxAtStartPar
\sphinxcode{\sphinxupquote{ctype}}
\begin{quote}

\sphinxAtStartPar
The type of the affine cone. Possible values are detailed in {\hyperref[\detokenize{constant:chapconst-conetype}]{\sphinxcrossref{\DUrole{std,std-ref}{SOC constraint types}}}} or {\hyperref[\detokenize{constant:chapconst-expconetype}]{\sphinxcrossref{\DUrole{std,std-ref}{Exponential Cone type}}}}.
\end{quote}
\end{quote}
\end{quote}

\subsubsection{Model.addQConstr()}
\label{\detokenize{pyapiref:model-addqconstr}}\begin{quote}

\sphinxAtStartPar
\sphinxstylestrong{Synopsis}
\begin{quote}

\sphinxAtStartPar
\sphinxcode{\sphinxupquote{addQConstr(lhs, sense=None, rhs=None, name="")}}
\end{quote}

\sphinxAtStartPar
\sphinxstylestrong{Description}
\begin{quote}

\sphinxAtStartPar
Add a linear or quadratic constraint, and return the added {\hyperref[\detokenize{pyapiref:chappyapi-constraint}]{\sphinxcrossref{\DUrole{std,std-ref}{Constraint Class}}}} object
or {\hyperref[\detokenize{pyapiref:chappyapi-qconstraint}]{\sphinxcrossref{\DUrole{std,std-ref}{QConstraint Class}}}} object.

\sphinxAtStartPar
If the constraint is linear, then value of parameter \sphinxcode{\sphinxupquote{lhs}} can be taken {\hyperref[\detokenize{pyapiref:chappyapi-var}]{\sphinxcrossref{\DUrole{std,std-ref}{Var Class}}}} object,
{\hyperref[\detokenize{pyapiref:chappyapi-linexpr}]{\sphinxcrossref{\DUrole{std,std-ref}{LinExpr Class}}}} object or {\hyperref[\detokenize{pyapiref:chappyapi-constrbuilder}]{\sphinxcrossref{\DUrole{std,std-ref}{ConstrBuilder Class}}}} object;
If the constraint is quadratic, then value of parameter \sphinxcode{\sphinxupquote{lhs}} can be taken {\hyperref[\detokenize{pyapiref:chappyapi-qconstrbuilder}]{\sphinxcrossref{\DUrole{std,std-ref}{QConstrBuilder Class}}}} object, or {\hyperref[\detokenize{pyapiref:chappyapi-mqconstrbuilder}]{\sphinxcrossref{\DUrole{std,std-ref}{MQConstrBuilder Class}}}} object
and other parameters will be ignored.
\end{quote}

\sphinxAtStartPar
\sphinxstylestrong{Arguments}
\begin{quote}

\sphinxAtStartPar
\sphinxcode{\sphinxupquote{lhs}}
\begin{quote}

\sphinxAtStartPar
Left\sphinxhyphen{}hand side expression for new constraint or constraint builder.
\end{quote}

\sphinxAtStartPar
\sphinxcode{\sphinxupquote{sense}}
\begin{quote}

\sphinxAtStartPar
Sense for the new constraint. Optional, None by default.
Please refer to  {\hyperref[\detokenize{constant:chapconst-constrtype}]{\sphinxcrossref{\DUrole{std,std-ref}{Constraint type}}}}  for possible values.
\end{quote}

\sphinxAtStartPar
\sphinxcode{\sphinxupquote{rhs}}
\begin{quote}

\sphinxAtStartPar
Right\sphinxhyphen{}hand side expression for the new constraint. Optional, None by default.
It can be a constant, {\hyperref[\detokenize{pyapiref:chappyapi-var}]{\sphinxcrossref{\DUrole{std,std-ref}{Var Class}}}} object, {\hyperref[\detokenize{pyapiref:chappyapi-linexpr}]{\sphinxcrossref{\DUrole{std,std-ref}{LinExpr Class}}}} object
or {\hyperref[\detokenize{pyapiref:chappyapi-quadexpr}]{\sphinxcrossref{\DUrole{std,std-ref}{QuadExpr Class}}}} object.
\end{quote}

\sphinxAtStartPar
\sphinxcode{\sphinxupquote{name}}
\begin{quote}

\sphinxAtStartPar
Name for new constraint. Optional, \sphinxcode{\sphinxupquote{""}} by default, generated by solver automatically.
\end{quote}
\end{quote}

\sphinxAtStartPar
\sphinxstylestrong{Example}
\end{quote}

\begin{sphinxVerbatim}[commandchars=\\\{\}]
\PYG{c+c1}{\PYGZsh{} add a linear equality: x + y == 2}
\PYG{n}{m}\PYG{o}{.}\PYG{n}{addQConstr}\PYG{p}{(}\PYG{n}{x} \PYG{o}{+} \PYG{n}{y}\PYG{p}{,} \PYG{n}{COPT}\PYG{o}{.}\PYG{n}{EQUAL}\PYG{p}{,} \PYG{l+m+mi}{2}\PYG{p}{)}
\PYG{c+c1}{\PYGZsh{} add a quadratic inequality: x*x + y*y \PYGZlt{}= 3}
\PYG{n}{m}\PYG{o}{.}\PYG{n}{addQConstr}\PYG{p}{(}\PYG{n}{x}\PYG{o}{*}\PYG{n}{x} \PYG{o}{+} \PYG{n}{y}\PYG{o}{*}\PYG{n}{y} \PYG{o}{\PYGZlt{}}\PYG{o}{=} \PYG{l+m+mf}{3.0}\PYG{p}{)}
\end{sphinxVerbatim}

\subsubsection{Model.addMQConstr()}
\label{\detokenize{pyapiref:model-addmqconstr}}\begin{quote}

\sphinxAtStartPar
\sphinxstylestrong{Synopsis}
\begin{quote}

\sphinxAtStartPar
\sphinxcode{\sphinxupquote{addMQConstr(Q, c, sense, rhs, xQ\_L=None, xQ\_R=None, xc=None, name="")}}
\end{quote}

\sphinxAtStartPar
\sphinxstylestrong{Description}
\begin{quote}

\sphinxAtStartPar
By means of matrix modeling, a quadratic constraint is added to the model. If the value of \sphinxcode{\sphinxupquote{sense}} here is \sphinxcode{\sphinxupquote{COPT.LESS\_EQUAL}} , the added constraint is \(x_{Q_L} Q x_{Q_R} + c x_c <= rhs\).

\sphinxAtStartPar
It is more convenient to generate {\hyperref[\detokenize{pyapiref:chappyapi-mquadexpr}]{\sphinxcrossref{\DUrole{std,std-ref}{MQuadExpr Class}}}} objects by matrix multiplication, and then use overloaded comparison operators to generate {\hyperref[\detokenize{pyapiref:chappyapi-mqconstrbuilder}]{\sphinxcrossref{\DUrole{std,std-ref}{MQConstrBuilder Class}}}} object, which can be used as input of \sphinxcode{\sphinxupquote{Model.addQConstr()}} to generate quadratic constraints.
\end{quote}

\sphinxAtStartPar
\sphinxstylestrong{Arguments}
\begin{quote}

\sphinxAtStartPar
\sphinxcode{\sphinxupquote{Q}}
\begin{quote}

\sphinxAtStartPar
If the quadratic term is not empty, the parameter Q needs to be provided, which is a two\sphinxhyphen{}dimensional NumPy matrix, SciPy compressed sparse column matrix ( \sphinxcode{\sphinxupquote{csc\_matrix}} ) or compressed sparse row matrix ( \sphinxcode{\sphinxupquote{csr\_matrix}} ).
\end{quote}

\sphinxAtStartPar
\sphinxcode{\sphinxupquote{c}}
\begin{quote}

\sphinxAtStartPar
If the item is non\sphinxhyphen{}empty, you need to provide the parameter c, which is a one\sphinxhyphen{}dimensional NumPy array, or a Python list.
\end{quote}

\sphinxAtStartPar
\sphinxcode{\sphinxupquote{sense}}
\begin{quote}

\sphinxAtStartPar
The type of constraint. Possible values refer to {\hyperref[\detokenize{constant:chapconst-constrtype}]{\sphinxcrossref{\DUrole{std,std-ref}{Constraint type}}}}.
\end{quote}

\sphinxAtStartPar
\sphinxcode{\sphinxupquote{rhs}}
\begin{quote}

\sphinxAtStartPar
The value on the right side of the constraint, usually a floating point number.
\end{quote}

\sphinxAtStartPar
\sphinxcode{\sphinxupquote{xQ\_L}}
\begin{quote}

\sphinxAtStartPar
The variable on the left side of the quadratic term, can be a {\hyperref[\detokenize{pyapiref:chappyapi-mvar}]{\sphinxcrossref{\DUrole{std,std-ref}{MVar Class}}}} object, {\hyperref[\detokenize{pyapiref:chappyapi-vararray}]{\sphinxcrossref{\DUrole{std,std-ref}{VarArray Class}}}} object, list, dictionary or {\hyperref[\detokenize{pyapiref:chappyapi-util-tupledict}]{\sphinxcrossref{\DUrole{std,std-ref}{tupledict Class}}}} object.

\sphinxAtStartPar
If empty, all variables in the model are taken.
\end{quote}

\sphinxAtStartPar
\sphinxcode{\sphinxupquote{xQ\_R}}
\begin{quote}

\sphinxAtStartPar
The variable on the right side of the quadratic term, can be a {\hyperref[\detokenize{pyapiref:chappyapi-mvar}]{\sphinxcrossref{\DUrole{std,std-ref}{MVar Class}}}} object, {\hyperref[\detokenize{pyapiref:chappyapi-vararray}]{\sphinxcrossref{\DUrole{std,std-ref}{VarArray Class}}}} object,
list, dictionary or {\hyperref[\detokenize{pyapiref:chappyapi-util-tupledict}]{\sphinxcrossref{\DUrole{std,std-ref}{tupledict Class}}}} object.

\sphinxAtStartPar
If empty, all variables in the model are taken.
\end{quote}

\sphinxAtStartPar
\sphinxcode{\sphinxupquote{xc}}
\begin{quote}

\sphinxAtStartPar
The variable corresponding to the linear term can be a {\hyperref[\detokenize{pyapiref:chappyapi-mvar}]{\sphinxcrossref{\DUrole{std,std-ref}{MVar Class}}}} object, {\hyperref[\detokenize{pyapiref:chappyapi-vararray}]{\sphinxcrossref{\DUrole{std,std-ref}{VarArray Class}}}} object, list, dictionary or {\hyperref[\detokenize{pyapiref:chappyapi-util-tupledict}]{\sphinxcrossref{\DUrole{std,std-ref}{tupledict Class}}}} object.

\sphinxAtStartPar
If it is empty, but the parameter c is not empty, all variables in the model are taken.
\end{quote}

\sphinxAtStartPar
\sphinxcode{\sphinxupquote{name}}
\begin{quote}

\sphinxAtStartPar
constraint name.
\end{quote}
\end{quote}

\sphinxAtStartPar
\sphinxstylestrong{Return value}
\begin{quote}

\sphinxAtStartPar
Returns a {\hyperref[\detokenize{pyapiref:chappyapi-qconstraint}]{\sphinxcrossref{\DUrole{std,std-ref}{QConstraint Class}}}} object
\end{quote}

\sphinxAtStartPar
\sphinxstylestrong{Example}

\begin{sphinxVerbatim}[commandchars=\\\{\}]
\PYG{n}{Q} \PYG{o}{=} \PYG{n}{np}\PYG{o}{.}\PYG{n}{full}\PYG{p}{(}\PYG{p}{(}\PYG{l+m+mi}{3}\PYG{p}{,} \PYG{l+m+mi}{3}\PYG{p}{)}\PYG{p}{,} \PYG{l+m+mi}{1}\PYG{p}{)}
\PYG{n}{mx} \PYG{o}{=} \PYG{n}{model}\PYG{o}{.}\PYG{n}{addMVar}\PYG{p}{(}\PYG{l+m+mi}{3}\PYG{p}{,} \PYG{n}{nameprefix}\PYG{o}{=}\PYG{l+s+s2}{\PYGZdq{}}\PYG{l+s+s2}{mx}\PYG{l+s+s2}{\PYGZdq{}}\PYG{p}{)}
\PYG{n}{mqc} \PYG{o}{=} \PYG{n}{model}\PYG{o}{.}\PYG{n}{addMQConstr}\PYG{p}{(}\PYG{n}{Q}\PYG{p}{,} \PYG{k+kc}{None}\PYG{p}{,} \PYG{l+s+s1}{\PYGZsq{}}\PYG{l+s+s1}{L}\PYG{l+s+s1}{\PYGZsq{}}\PYG{p}{,} \PYG{l+m+mf}{1.0}\PYG{p}{,} \PYG{n}{mx}\PYG{p}{,} \PYG{n}{mx}\PYG{p}{,} \PYG{k+kc}{None}\PYG{p}{,} \PYG{n}{name}\PYG{o}{=}\PYG{l+s+s2}{\PYGZdq{}}\PYG{l+s+s2}{mqc}\PYG{l+s+s2}{\PYGZdq{}}\PYG{p}{)}
\end{sphinxVerbatim}
\end{quote}

\subsubsection{Model.addNlConstr()}
\label{\detokenize{pyapiref:model-addnlconstr}}\begin{quote}

\sphinxAtStartPar
\sphinxstylestrong{Synopsis}
\begin{quote}

\sphinxAtStartPar
\sphinxcode{\sphinxupquote{addNlConstr(lhs, sense=None, rhs=None, name="")}}
\end{quote}

\sphinxAtStartPar
\sphinxstylestrong{Description}
\begin{quote}

\sphinxAtStartPar
Add a nonlinear constraint to the model.

\sphinxAtStartPar
If \sphinxcode{\sphinxupquote{lhs}} is a nonlinear expression type,
both \sphinxcode{\sphinxupquote{sense}} and \sphinxcode{\sphinxupquote{rhs}} must be provided to specify
the constraint direction and right\sphinxhyphen{}hand side.

\sphinxAtStartPar
If \sphinxcode{\sphinxupquote{lhs}} is a {\hyperref[\detokenize{pyapiref:chappyapi-nlconstrbuilder}]{\sphinxcrossref{\DUrole{std,std-ref}{NlConstrBuilder Class}}}} object,
the constraint is created using the builder content.
Other arguments are ignored.

\sphinxAtStartPar
If \sphinxcode{\sphinxupquote{lhs}} is a {\hyperref[\detokenize{pyapiref:chappyapi-constrbuilder}]{\sphinxcrossref{\DUrole{std,std-ref}{ConstrBuilder Class}}}} object,
a linear constraint is added.
The return value is {\hyperref[\detokenize{pyapiref:chappyapi-constraint}]{\sphinxcrossref{\DUrole{std,std-ref}{Constraint Class}}}}.

\sphinxAtStartPar
If \sphinxcode{\sphinxupquote{lhs}} is a {\hyperref[\detokenize{pyapiref:chappyapi-qconstrbuilder}]{\sphinxcrossref{\DUrole{std,std-ref}{QConstrBuilder Class}}}} object,
a quadratic constraint is added.
The return value is {\hyperref[\detokenize{pyapiref:chappyapi-qconstraint}]{\sphinxcrossref{\DUrole{std,std-ref}{QConstraint Class}}}}.
\end{quote}

\sphinxAtStartPar
\sphinxstylestrong{Arguments}
\begin{quote}

\sphinxAtStartPar
\sphinxcode{\sphinxupquote{lhs}}
\begin{quote}

\sphinxAtStartPar
The left\sphinxhyphen{}hand side expression or constraint builder.
\end{quote}

\sphinxAtStartPar
\sphinxcode{\sphinxupquote{sense}}
\begin{quote}

\sphinxAtStartPar
The constraint type. Optional. Defaults to \sphinxcode{\sphinxupquote{None}}.

\sphinxAtStartPar
See {\hyperref[\detokenize{constant:chapconst-constrtype}]{\sphinxcrossref{\DUrole{std,std-ref}{Constraint Types}}}} for valid values.
\end{quote}

\sphinxAtStartPar
\sphinxcode{\sphinxupquote{rhs}}
\begin{quote}

\sphinxAtStartPar
The right\sphinxhyphen{}hand side of the constraint. Optional. Defaults to \sphinxcode{\sphinxupquote{None}}.

\sphinxAtStartPar
Can be a constant or a {\hyperref[\detokenize{pyapiref:chappyapi-nlexpr}]{\sphinxcrossref{\DUrole{std,std-ref}{NlExpr Class}}}} object.
\end{quote}

\sphinxAtStartPar
\sphinxcode{\sphinxupquote{name}}
\begin{quote}

\sphinxAtStartPar
The name of the constraint. Optional.

\sphinxAtStartPar
Defaults to \sphinxcode{\sphinxupquote{""}}, and a name will be generated automatically by the solver.
\end{quote}
\end{quote}
\end{quote}

\subsubsection{Model.addNlConstrs()}
\label{\detokenize{pyapiref:model-addnlconstrs}}\begin{quote}

\sphinxAtStartPar
\sphinxstylestrong{Synopsis}
\begin{quote}

\sphinxAtStartPar
\sphinxcode{\sphinxupquote{addNlConstrs(generator, nameprefix="NR")}}
\end{quote}

\sphinxAtStartPar
\sphinxstylestrong{Description}
\begin{quote}

\sphinxAtStartPar
Add a set of nonlinear constraints to the model.
\end{quote}

\sphinxAtStartPar
\sphinxstylestrong{Arguments}
\begin{quote}

\sphinxAtStartPar
\sphinxcode{\sphinxupquote{generator}}
\begin{quote}

\sphinxAtStartPar
A constraint builder or a sequence of nonlinear constraint builders.

\sphinxAtStartPar
Possible values include:
a single {\hyperref[\detokenize{pyapiref:chappyapi-nlconstrbuilder}]{\sphinxcrossref{\DUrole{std,std-ref}{NlConstrBuilder Class}}}}, {\hyperref[\detokenize{pyapiref:chappyapi-qconstrbuilder}]{\sphinxcrossref{\DUrole{std,std-ref}{QConstrBuilder Class}}}},
or {\hyperref[\detokenize{pyapiref:chappyapi-constrbuilder}]{\sphinxcrossref{\DUrole{std,std-ref}{ConstrBuilder Class}}}} object,
or an iterable of {\hyperref[\detokenize{pyapiref:chappyapi-nlconstrbuilder}]{\sphinxcrossref{\DUrole{std,std-ref}{NlConstrBuilder Class}}}} objects.
\end{quote}

\sphinxAtStartPar
\sphinxcode{\sphinxupquote{nameprefix}}
\begin{quote}

\sphinxAtStartPar
The prefix used for naming constraints.

\sphinxAtStartPar
Optional. Defaults to \sphinxcode{\sphinxupquote{"NR"}}.
\end{quote}
\end{quote}
\end{quote}

\subsubsection{Model.addPsdVar()}
\label{\detokenize{pyapiref:model-addpsdvar}}\begin{quote}

\sphinxAtStartPar
\sphinxstylestrong{Synopsis}
\begin{quote}

\sphinxAtStartPar
\sphinxcode{\sphinxupquote{addPsdVar(dim, name="")}}
\end{quote}

\sphinxAtStartPar
\sphinxstylestrong{Description}
\begin{quote}

\sphinxAtStartPar
Add a positive semi\sphinxhyphen{}definite variable.
\end{quote}

\sphinxAtStartPar
\sphinxstylestrong{Arguments}
\begin{quote}

\sphinxAtStartPar
\sphinxcode{\sphinxupquote{dim}}
\begin{quote}

\sphinxAtStartPar
Dimension for the positive semi\sphinxhyphen{}definite variable.
\end{quote}

\sphinxAtStartPar
\sphinxcode{\sphinxupquote{name}}
\begin{quote}

\sphinxAtStartPar
Name for the positive semi\sphinxhyphen{}definite variable.
\end{quote}
\end{quote}

\sphinxAtStartPar
\sphinxstylestrong{Example}
\end{quote}

\begin{sphinxVerbatim}[commandchars=\\\{\}]
\PYG{c+c1}{\PYGZsh{} Add a three\PYGZhy{}dimensional positive semi\PYGZhy{}definite variable, \PYGZdq{}X\PYGZdq{}}
\PYG{n}{m}\PYG{o}{.}\PYG{n}{addPsdVar}\PYG{p}{(}\PYG{l+m+mi}{3}\PYG{p}{,} \PYG{l+s+s2}{\PYGZdq{}}\PYG{l+s+s2}{X}\PYG{l+s+s2}{\PYGZdq{}}\PYG{p}{)}
\end{sphinxVerbatim}

\subsubsection{Model.addPsdVars()}
\label{\detokenize{pyapiref:model-addpsdvars}}\begin{quote}

\sphinxAtStartPar
\sphinxstylestrong{Synopsis}
\begin{quote}

\sphinxAtStartPar
\sphinxcode{\sphinxupquote{addPsdVars(dims, nameprefix="PSDV")}}
\end{quote}

\sphinxAtStartPar
\sphinxstylestrong{Description}
\begin{quote}

\sphinxAtStartPar
Add multiple new positive semi\sphinxhyphen{}definite variables to a model.
\end{quote}

\sphinxAtStartPar
\sphinxstylestrong{Arguments}
\begin{quote}

\sphinxAtStartPar
\sphinxcode{\sphinxupquote{dim}}
\begin{quote}

\sphinxAtStartPar
Dimensions for new positive semi\sphinxhyphen{}definite variables.
\end{quote}

\sphinxAtStartPar
\sphinxcode{\sphinxupquote{nameprefix}}
\begin{quote}

\sphinxAtStartPar
Name prefix for new positive semi\sphinxhyphen{}definite variables.
\end{quote}
\end{quote}

\sphinxAtStartPar
\sphinxstylestrong{Example}
\end{quote}

\begin{sphinxVerbatim}[commandchars=\\\{\}]
\PYG{c+c1}{\PYGZsh{} Add two three\PYGZhy{}dimensional positive semi\PYGZhy{}definite variables}
\PYG{n}{m}\PYG{o}{.}\PYG{n}{addPsdVars}\PYG{p}{(}\PYG{p}{[}\PYG{l+m+mi}{3}\PYG{p}{,} \PYG{l+m+mi}{3}\PYG{p}{]}\PYG{p}{)}
\end{sphinxVerbatim}

\subsubsection{Model.addUserCut()}
\label{\detokenize{pyapiref:model-addusercut}}\begin{quote}

\sphinxAtStartPar
\sphinxstylestrong{Synopsis}
\begin{quote}

\sphinxAtStartPar
\sphinxcode{\sphinxupquote{addUserCut(lhs, sense = None, rhs = None, name="")}}
\end{quote}

\sphinxAtStartPar
\sphinxstylestrong{Description}
\begin{quote}

\sphinxAtStartPar
Add a user cut to the MIP model.
\end{quote}

\sphinxAtStartPar
\sphinxstylestrong{Arguments}
\begin{quote}

\sphinxAtStartPar
\sphinxcode{\sphinxupquote{lhs}}
\begin{quote}

\sphinxAtStartPar
Left\sphinxhyphen{}hand side expression for the new user cut. It can take the value of {\hyperref[\detokenize{pyapiref:chappyapi-var}]{\sphinxcrossref{\DUrole{std,std-ref}{Var Class}}}} object, {\hyperref[\detokenize{pyapiref:chappyapi-linexpr}]{\sphinxcrossref{\DUrole{std,std-ref}{LinExpr Class}}}} object or {\hyperref[\detokenize{pyapiref:chappyapi-constrbuilder}]{\sphinxcrossref{\DUrole{std,std-ref}{ConstrBuilder Class}}}} object.
\end{quote}

\sphinxAtStartPar
\sphinxcode{\sphinxupquote{sense}}
\begin{quote}

\sphinxAtStartPar
The sense of the new user cut. Please refer to  {\hyperref[\detokenize{constant:chapconst-constrtype}]{\sphinxcrossref{\DUrole{std,std-ref}{Constraint type}}}}  for possible values.

\sphinxAtStartPar
Optional. None by default.
\end{quote}

\sphinxAtStartPar
\sphinxcode{\sphinxupquote{rhs}}
\begin{quote}

\sphinxAtStartPar
Right\sphinxhyphen{}hand side expression for the new user cut.

\sphinxAtStartPar
Optional. None by default.

\sphinxAtStartPar
It can be a constant, or {\hyperref[\detokenize{pyapiref:chappyapi-var}]{\sphinxcrossref{\DUrole{std,std-ref}{Var Class}}}} object, or {\hyperref[\detokenize{pyapiref:chappyapi-linexpr}]{\sphinxcrossref{\DUrole{std,std-ref}{LinExpr Class}}}} object.
\end{quote}

\sphinxAtStartPar
\sphinxcode{\sphinxupquote{name}}
\begin{quote}

\sphinxAtStartPar
Name for the new user cut. Optional, \sphinxcode{\sphinxupquote{""}} by default, automatically generated by solver.
\end{quote}
\end{quote}

\sphinxAtStartPar
\sphinxstylestrong{Example}

\begin{sphinxVerbatim}[commandchars=\\\{\}]
\PYG{n}{model}\PYG{o}{.}\PYG{n}{addUserCut}\PYG{p}{(}\PYG{n}{x}\PYG{o}{+}\PYG{n}{y} \PYG{o}{\PYGZlt{}}\PYG{o}{=} \PYG{l+m+mi}{1}\PYG{p}{)}

\PYG{n}{model}\PYG{o}{.}\PYG{n}{addUserCut}\PYG{p}{(}\PYG{n}{x}\PYG{o}{+}\PYG{n}{y} \PYG{o}{==} \PYG{p}{[}\PYG{l+m+mi}{0}\PYG{p}{,} \PYG{l+m+mi}{1}\PYG{p}{]}\PYG{p}{)}
\end{sphinxVerbatim}
\end{quote}

\subsubsection{Model.addUserCuts()}
\label{\detokenize{pyapiref:model-addusercuts}}\begin{quote}

\sphinxAtStartPar
\sphinxstylestrong{Synopsis}
\begin{quote}

\sphinxAtStartPar
\sphinxcode{\sphinxupquote{addUserCuts(generator, nameprefix="U")}}
\end{quote}

\sphinxAtStartPar
\sphinxstylestrong{Description}
\begin{quote}

\sphinxAtStartPar
Add a set of user cuts to the MIP model.
\end{quote}

\sphinxAtStartPar
\sphinxstylestrong{Arguments}
\begin{quote}

\sphinxAtStartPar
\sphinxcode{\sphinxupquote{generator}}
\begin{quote}

\sphinxAtStartPar
A generator expression, where each iteration produces a  {\hyperref[\detokenize{pyapiref:chappyapi-constraint}]{\sphinxcrossref{\DUrole{std,std-ref}{Constraint Class}}}} object,
or {\hyperref[\detokenize{pyapiref:chappyapi-mconstrbuilder}]{\sphinxcrossref{\DUrole{std,std-ref}{MConstrBuilder Class}}}} object.
\end{quote}

\sphinxAtStartPar
\sphinxcode{\sphinxupquote{nameprefix}}
\begin{quote}

\sphinxAtStartPar
Name prefix for new user cuts. Optional, \sphinxcode{\sphinxupquote{"U"}} by default.
The actual name and the index of the constraints are automatically generated by COPT.
\end{quote}
\end{quote}

\sphinxAtStartPar
\sphinxstylestrong{Example}

\begin{sphinxVerbatim}[commandchars=\\\{\}]
\PYG{n}{model}\PYG{o}{.}\PYG{n}{addUserCuts}\PYG{p}{(}\PYG{n}{x}\PYG{p}{[}\PYG{n}{i}\PYG{p}{]}\PYG{o}{+}\PYG{n}{y}\PYG{p}{[}\PYG{n}{i}\PYG{p}{]} \PYG{o}{\PYGZlt{}}\PYG{o}{=} \PYG{l+m+mi}{1} \PYG{k}{for} \PYG{n}{i} \PYG{o+ow}{in} \PYG{n+nb}{range}\PYG{p}{(}\PYG{l+m+mi}{10}\PYG{p}{)}\PYG{p}{)}
\end{sphinxVerbatim}
\end{quote}

\subsubsection{Model.addLazyConstr()}
\label{\detokenize{pyapiref:model-addlazyconstr}}\begin{quote}

\sphinxAtStartPar
\sphinxstylestrong{Synopsis}
\begin{quote}

\sphinxAtStartPar
\sphinxcode{\sphinxupquote{addLazyConstr(lhs, sense = None, rhs = None, name="")}}
\end{quote}

\sphinxAtStartPar
\sphinxstylestrong{Description}
\begin{quote}

\sphinxAtStartPar
Add a lazy constraint to the MIP model.
\end{quote}

\sphinxAtStartPar
\sphinxstylestrong{Arguments}
\begin{quote}

\sphinxAtStartPar
\sphinxcode{\sphinxupquote{lhs}}
\begin{quote}

\sphinxAtStartPar
Left\sphinxhyphen{}hand side expression for the new lazy constraint. It can take the value of {\hyperref[\detokenize{pyapiref:chappyapi-var}]{\sphinxcrossref{\DUrole{std,std-ref}{Var Class}}}} object, {\hyperref[\detokenize{pyapiref:chappyapi-linexpr}]{\sphinxcrossref{\DUrole{std,std-ref}{LinExpr Class}}}} object or {\hyperref[\detokenize{pyapiref:chappyapi-constrbuilder}]{\sphinxcrossref{\DUrole{std,std-ref}{ConstrBuilder Class}}}} object.
\end{quote}

\sphinxAtStartPar
\sphinxcode{\sphinxupquote{sense}}
\begin{quote}

\sphinxAtStartPar
The sense of the lazy constraint. Please refer to  {\hyperref[\detokenize{constant:chapconst-constrtype}]{\sphinxcrossref{\DUrole{std,std-ref}{Constraint type}}}}  for possible values.

\sphinxAtStartPar
Optional. None by default.
\end{quote}

\sphinxAtStartPar
\sphinxcode{\sphinxupquote{rhs}}
\begin{quote}

\sphinxAtStartPar
Right\sphinxhyphen{}hand side expression for the new lazy constraint.

\sphinxAtStartPar
Optional. None by default.

\sphinxAtStartPar
It can be a constant, or {\hyperref[\detokenize{pyapiref:chappyapi-var}]{\sphinxcrossref{\DUrole{std,std-ref}{Var Class}}}} object, or {\hyperref[\detokenize{pyapiref:chappyapi-linexpr}]{\sphinxcrossref{\DUrole{std,std-ref}{LinExpr Class}}}} object.
\end{quote}

\sphinxAtStartPar
\sphinxcode{\sphinxupquote{name}}
\begin{quote}

\sphinxAtStartPar
Name for the new lazy constraint. Optional, \sphinxcode{\sphinxupquote{""}} by default, automatically generated by solver.
\end{quote}
\end{quote}

\sphinxAtStartPar
\sphinxstylestrong{Example}

\begin{sphinxVerbatim}[commandchars=\\\{\}]
\PYG{n}{model}\PYG{o}{.}\PYG{n}{addLazyConstr}\PYG{p}{(}\PYG{n}{x}\PYG{o}{+}\PYG{n}{y} \PYG{o}{\PYGZlt{}}\PYG{o}{=} \PYG{l+m+mi}{1}\PYG{p}{)}
\PYG{n}{model}\PYG{o}{.}\PYG{n}{addLazyConstr}\PYG{p}{(}\PYG{n}{x}\PYG{o}{+}\PYG{n}{y} \PYG{o}{==} \PYG{p}{[}\PYG{l+m+mi}{0}\PYG{p}{,} \PYG{l+m+mi}{1}\PYG{p}{]}\PYG{p}{)}
\end{sphinxVerbatim}
\end{quote}

\subsubsection{Model.addLazyConstrs()}
\label{\detokenize{pyapiref:model-addlazyconstrs}}\begin{quote}

\sphinxAtStartPar
\sphinxstylestrong{Synopsis}
\begin{quote}

\sphinxAtStartPar
\sphinxcode{\sphinxupquote{addLazyConstrs(generator, nameprefix="L")}}
\end{quote}

\sphinxAtStartPar
\sphinxstylestrong{Description}
\begin{quote}

\sphinxAtStartPar
Add a set of lazy constraints to the MIP model.
\end{quote}

\sphinxAtStartPar
\sphinxstylestrong{Arguments}
\begin{quote}

\sphinxAtStartPar
\sphinxcode{\sphinxupquote{generator}}
\begin{quote}

\sphinxAtStartPar
A generator expression, where each iteration produces a {\hyperref[\detokenize{pyapiref:chappyapi-constraint}]{\sphinxcrossref{\DUrole{std,std-ref}{Constraint Class}}}} object,
or {\hyperref[\detokenize{pyapiref:chappyapi-mconstrbuilder}]{\sphinxcrossref{\DUrole{std,std-ref}{MConstrBuilder Class}}}} object.
\end{quote}

\sphinxAtStartPar
\sphinxcode{\sphinxupquote{nameprefix}}
\begin{quote}

\sphinxAtStartPar
Name prefix for new lazy constraints. Optional, \sphinxcode{\sphinxupquote{"L"}} by default.
The actual name and the index of the constraints are automatically generated by COPT.
\end{quote}
\end{quote}

\sphinxAtStartPar
\sphinxstylestrong{Example}

\begin{sphinxVerbatim}[commandchars=\\\{\}]
\PYG{n}{model}\PYG{o}{.}\PYG{n}{addLazyConstrs}\PYG{p}{(}\PYG{n}{x}\PYG{p}{[}\PYG{n}{i}\PYG{p}{]}\PYG{o}{+}\PYG{n}{y}\PYG{p}{[}\PYG{n}{i}\PYG{p}{]} \PYG{o}{\PYGZlt{}}\PYG{o}{=} \PYG{l+m+mi}{1} \PYG{k}{for} \PYG{n}{i} \PYG{o+ow}{in} \PYG{n+nb}{range}\PYG{p}{(}\PYG{l+m+mi}{10}\PYG{p}{)}\PYG{p}{)}
\end{sphinxVerbatim}
\end{quote}

\subsubsection{Model.addSparseMat()}
\label{\detokenize{pyapiref:model-addsparsemat}}\begin{quote}

\sphinxAtStartPar
\sphinxstylestrong{Synopsis}
\begin{quote}

\sphinxAtStartPar
\sphinxcode{\sphinxupquote{addSparseMat(dim, rows, cols=None, vals=None)}}
\end{quote}

\sphinxAtStartPar
\sphinxstylestrong{Description}
\begin{quote}

\sphinxAtStartPar
Add a sparse symmetric matrix in triplet format
\end{quote}

\sphinxAtStartPar
\sphinxstylestrong{Arguments}
\begin{quote}

\sphinxAtStartPar
\sphinxcode{\sphinxupquote{dim}}
\begin{quote}

\sphinxAtStartPar
Dimension for the matrix.
\end{quote}

\sphinxAtStartPar
\sphinxcode{\sphinxupquote{rows}}
\begin{quote}

\sphinxAtStartPar
Row indices for accessing rows of non\sphinxhyphen{}zero elements.
\end{quote}

\sphinxAtStartPar
\sphinxcode{\sphinxupquote{cols}}
\begin{quote}

\sphinxAtStartPar
Column indices for accessing columns of non\sphinxhyphen{}zero elements.
\end{quote}

\sphinxAtStartPar
\sphinxcode{\sphinxupquote{vals}}
\begin{quote}

\sphinxAtStartPar
Coefficient values for non\sphinxhyphen{}zero elements
\end{quote}
\end{quote}

\sphinxAtStartPar
\sphinxstylestrong{Example}
\end{quote}

\begin{sphinxVerbatim}[commandchars=\\\{\}]
\PYG{c+c1}{\PYGZsh{} Add a three\PYGZhy{}dimentional symmetric matrix}
\PYG{n}{m}\PYG{o}{.}\PYG{n}{addSparseMat}\PYG{p}{(}\PYG{l+m+mi}{3}\PYG{p}{,} \PYG{p}{[}\PYG{l+m+mi}{0}\PYG{p}{,} \PYG{l+m+mi}{1}\PYG{p}{,} \PYG{l+m+mi}{2}\PYG{p}{]}\PYG{p}{,} \PYG{p}{[}\PYG{l+m+mi}{0}\PYG{p}{,} \PYG{l+m+mi}{1}\PYG{p}{,} \PYG{l+m+mi}{2}\PYG{p}{]}\PYG{p}{,} \PYG{p}{[}\PYG{l+m+mf}{2.0}\PYG{p}{,} \PYG{l+m+mf}{5.0}\PYG{p}{,} \PYG{l+m+mf}{8.0}\PYG{p}{]}\PYG{p}{)}
\PYG{c+c1}{\PYGZsh{} Add a two\PYGZhy{}dimentional symmetric matrix}
\PYG{n}{m}\PYG{o}{.}\PYG{n}{addSparseMat}\PYG{p}{(}\PYG{l+m+mi}{2}\PYG{p}{,} \PYG{p}{[}\PYG{p}{(}\PYG{l+m+mi}{0}\PYG{p}{,} \PYG{l+m+mi}{0}\PYG{p}{,} \PYG{l+m+mf}{3.0}\PYG{p}{)}\PYG{p}{,} \PYG{p}{(}\PYG{l+m+mi}{1}\PYG{p}{,} \PYG{l+m+mi}{0}\PYG{p}{,} \PYG{l+m+mf}{1.0}\PYG{p}{)}\PYG{p}{]}\PYG{p}{)}
\end{sphinxVerbatim}

\subsubsection{Model.addDenseMat()}
\label{\detokenize{pyapiref:model-adddensemat}}\begin{quote}

\sphinxAtStartPar
\sphinxstylestrong{Synopsis}
\begin{quote}

\sphinxAtStartPar
\sphinxcode{\sphinxupquote{addDenseMat(dim, vals)}}
\end{quote}

\sphinxAtStartPar
\sphinxstylestrong{Description}
\begin{quote}

\sphinxAtStartPar
Add a dense symmetric matrix
\end{quote}

\sphinxAtStartPar
\sphinxstylestrong{Arguments}
\begin{quote}

\sphinxAtStartPar
\sphinxcode{\sphinxupquote{dim}}
\begin{quote}

\sphinxAtStartPar
Dimension for the matrix.
\end{quote}

\sphinxAtStartPar
\sphinxcode{\sphinxupquote{vals}}
\begin{quote}

\sphinxAtStartPar
Coefficient values, which can be a constant or a list.
\end{quote}
\end{quote}

\sphinxAtStartPar
\sphinxstylestrong{Example}
\end{quote}

\begin{sphinxVerbatim}[commandchars=\\\{\}]
\PYG{c+c1}{\PYGZsh{} Add a tree\PYGZhy{}dimentional matrix (filled with ones)}
\PYG{n}{m}\PYG{o}{.}\PYG{n}{addDenseMat}\PYG{p}{(}\PYG{l+m+mi}{3}\PYG{p}{,} \PYG{l+m+mf}{1.0}\PYG{p}{)}
\end{sphinxVerbatim}

\subsubsection{Model.addDiagMat()}
\label{\detokenize{pyapiref:model-adddiagmat}}\begin{quote}

\sphinxAtStartPar
\sphinxstylestrong{Synopsis}
\begin{quote}

\sphinxAtStartPar
\sphinxcode{\sphinxupquote{addDiagMat(dim, vals, offset=None)}}
\end{quote}

\sphinxAtStartPar
\sphinxstylestrong{Description}
\begin{quote}

\sphinxAtStartPar
Add a diagnal symmetric matrix
\end{quote}

\sphinxAtStartPar
\sphinxstylestrong{Arguments}
\begin{quote}

\sphinxAtStartPar
\sphinxcode{\sphinxupquote{dim}}
\begin{quote}

\sphinxAtStartPar
Dimension for the matrix.
\end{quote}

\sphinxAtStartPar
\sphinxcode{\sphinxupquote{vals}}
\begin{quote}

\sphinxAtStartPar
Coefficient values, which can be a constant or a list.
\end{quote}

\sphinxAtStartPar
\sphinxcode{\sphinxupquote{offset}}
\begin{quote}

\sphinxAtStartPar
Offset of diagnal elements. Positive: above the diagnal; Negative: below the diagnal
\end{quote}
\end{quote}

\sphinxAtStartPar
\sphinxstylestrong{Example}
\end{quote}

\begin{sphinxVerbatim}[commandchars=\\\{\}]
\PYG{c+c1}{\PYGZsh{} Add a tree\PYGZhy{}dimentional identity matrix}
\PYG{n}{m}\PYG{o}{.}\PYG{n}{addDiagMat}\PYG{p}{(}\PYG{l+m+mi}{3}\PYG{p}{,} \PYG{l+m+mf}{1.0}\PYG{p}{)}
\end{sphinxVerbatim}

\subsubsection{Model.addOnesMat()}
\label{\detokenize{pyapiref:model-addonesmat}}\begin{quote}

\sphinxAtStartPar
\sphinxstylestrong{Synopsis}
\begin{quote}

\sphinxAtStartPar
\sphinxcode{\sphinxupquote{addOnesMat(dim)}}
\end{quote}

\sphinxAtStartPar
\sphinxstylestrong{Description}
\begin{quote}

\sphinxAtStartPar
Add a matrix filled with ones.
\end{quote}

\sphinxAtStartPar
\sphinxstylestrong{Arguments}
\begin{quote}

\sphinxAtStartPar
\sphinxcode{\sphinxupquote{dim}}
\begin{quote}

\sphinxAtStartPar
Dimension for the matrix.
\end{quote}
\end{quote}

\sphinxAtStartPar
\sphinxstylestrong{Example}
\end{quote}

\begin{sphinxVerbatim}[commandchars=\\\{\}]
\PYG{c+c1}{\PYGZsh{} Add a tree\PYGZhy{}dimentional matrix (filled with ones)}
\PYG{n}{m}\PYG{o}{.}\PYG{n}{addOnesMat}\PYG{p}{(}\PYG{l+m+mi}{3}\PYG{p}{)}
\end{sphinxVerbatim}

\subsubsection{Model.addEyeMat()}
\label{\detokenize{pyapiref:model-addeyemat}}\begin{quote}

\sphinxAtStartPar
\sphinxstylestrong{Synopsis}
\begin{quote}

\sphinxAtStartPar
\sphinxcode{\sphinxupquote{addEyeMat(dim)}}
\end{quote}

\sphinxAtStartPar
\sphinxstylestrong{Description}
\begin{quote}

\sphinxAtStartPar
Add an identity matrix
\end{quote}

\sphinxAtStartPar
\sphinxstylestrong{Arguments}
\begin{quote}

\sphinxAtStartPar
\sphinxcode{\sphinxupquote{dim}}
\begin{quote}

\sphinxAtStartPar
Dimension for the matrix.
\end{quote}
\end{quote}

\sphinxAtStartPar
\sphinxstylestrong{Example}
\end{quote}

\begin{sphinxVerbatim}[commandchars=\\\{\}]
\PYG{c+c1}{\PYGZsh{} Add a tree\PYGZhy{}dimentional identity matrix}
\PYG{n}{m}\PYG{o}{.}\PYG{n}{addEyeMat}\PYG{p}{(}\PYG{l+m+mi}{3}\PYG{p}{)}
\end{sphinxVerbatim}

\subsubsection{Model.setObjective()}
\label{\detokenize{pyapiref:model-setobjective}}\begin{quote}

\sphinxAtStartPar
\sphinxstylestrong{Synopsis}
\begin{quote}

\sphinxAtStartPar
\sphinxcode{\sphinxupquote{setObjective(expr, sense=None)}}
\end{quote}

\sphinxAtStartPar
\sphinxstylestrong{Description}
\begin{quote}

\sphinxAtStartPar
Set the model objective.
\end{quote}

\sphinxAtStartPar
\sphinxstylestrong{Arguments}
\begin{quote}

\sphinxAtStartPar
\sphinxcode{\sphinxupquote{expr}}
\begin{quote}

\sphinxAtStartPar
Objective expression. Argument can be a constant, {\hyperref[\detokenize{pyapiref:chappyapi-var}]{\sphinxcrossref{\DUrole{std,std-ref}{Var Class}}}} object,
{\hyperref[\detokenize{pyapiref:chappyapi-linexpr}]{\sphinxcrossref{\DUrole{std,std-ref}{LinExpr Class}}}} object, {\hyperref[\detokenize{pyapiref:chappyapi-quadexpr}]{\sphinxcrossref{\DUrole{std,std-ref}{QuadExpr Class}}}} object,
{\hyperref[\detokenize{pyapiref:chappyapi-mlinexpr}]{\sphinxcrossref{\DUrole{std,std-ref}{MLinExpr Class}}}} object, or {\hyperref[\detokenize{pyapiref:chappyapi-mquadexpr}]{\sphinxcrossref{\DUrole{std,std-ref}{MQuadExpr Class}}}} object,
or {\hyperref[\detokenize{pyapiref:chappyapi-nlexpr}]{\sphinxcrossref{\DUrole{std,std-ref}{NlExpr Class}}}} object.

\sphinxAtStartPar
Note: If \sphinxcode{\sphinxupquote{expr}} is a {\hyperref[\detokenize{pyapiref:chappyapi-linexpr}]{\sphinxcrossref{\DUrole{std,std-ref}{LinExpr Class}}}} object, the linear term in the objective will be updated;
If it is a {\hyperref[\detokenize{pyapiref:chappyapi-quadexpr}]{\sphinxcrossref{\DUrole{std,std-ref}{QuadExpr Class}}}} object, both the quadratic and linear terms in the objective will be updated.
\end{quote}

\sphinxAtStartPar
\sphinxcode{\sphinxupquote{sense}}
\begin{quote}

\sphinxAtStartPar
Optimization sense. Optional, \sphinxcode{\sphinxupquote{None}} by default, which means no change to objective sense.
Users can get access to current objective sense by attribute \sphinxcode{\sphinxupquote{ObjSense}} .
Please refer to {\hyperref[\detokenize{constant:chapconst-sense}]{\sphinxcrossref{\DUrole{std,std-ref}{Optimization directions}}}} for possible values.
\end{quote}
\end{quote}

\sphinxAtStartPar
\sphinxstylestrong{Example}
\end{quote}

\begin{sphinxVerbatim}[commandchars=\\\{\}]
\PYG{c+c1}{\PYGZsh{} Set objective function = x+y, optimization sense is maximization.}
\PYG{n}{m}\PYG{o}{.}\PYG{n}{setObjective}\PYG{p}{(}\PYG{n}{x} \PYG{o}{+} \PYG{n}{y}\PYG{p}{,} \PYG{n}{COPT}\PYG{o}{.}\PYG{n}{MAXIMIZE}\PYG{p}{)}
\end{sphinxVerbatim}

\subsubsection{Model.setObjectiveN()}
\label{\detokenize{pyapiref:model-setobjectiven}}\begin{quote}

\sphinxAtStartPar
\sphinxstylestrong{Synopsis}
\begin{quote}

\sphinxAtStartPar
\sphinxcode{\sphinxupquote{setObjectiveN(idx, expr, sense=None, priority=0.0, weight=1.0, abstol=1e\sphinxhyphen{}6, reltol=0.0)}}
\end{quote}

\sphinxAtStartPar
\sphinxstylestrong{Description}
\begin{quote}

\sphinxAtStartPar
Specify the objective function for multi\sphinxhyphen{}objective optimization.
\end{quote}

\sphinxAtStartPar
\sphinxstylestrong{Arguments}
\begin{quote}

\sphinxAtStartPar
\sphinxcode{\sphinxupquote{idx}}
\begin{quote}

\sphinxAtStartPar
The index of the objective function. Argument can be an integer constant.
\end{quote}

\sphinxAtStartPar
\sphinxcode{\sphinxupquote{expr}}
\begin{quote}

\sphinxAtStartPar
The expression of the objective function. Argument can be a linear expression.
\end{quote}

\sphinxAtStartPar
\sphinxcode{\sphinxupquote{sense}}
\begin{quote}

\sphinxAtStartPar
Optimization sense.

\sphinxAtStartPar
Optional, \sphinxcode{\sphinxupquote{None}} by default, which means no change to objective sense.
Please refer to {\hyperref[\detokenize{constant:chapconst-sense}]{\sphinxcrossref{\DUrole{std,std-ref}{Optimization directions}}}} for possible values.
\end{quote}

\sphinxAtStartPar
\sphinxcode{\sphinxupquote{priority}}
\begin{quote}

\sphinxAtStartPar
The priority of the objective function.

\sphinxAtStartPar
Optional, default is 0.0.
\end{quote}

\sphinxAtStartPar
\sphinxcode{\sphinxupquote{weight}}
\begin{quote}

\sphinxAtStartPar
The weight of the objective function.

\sphinxAtStartPar
Optional, default is 1.0.
\end{quote}

\sphinxAtStartPar
\sphinxcode{\sphinxupquote{abstol}}
\begin{quote}

\sphinxAtStartPar
The absolute tolerance for degradation.

\sphinxAtStartPar
Optional, default is 1e\sphinxhyphen{}6.
\end{quote}

\sphinxAtStartPar
\sphinxcode{\sphinxupquote{reltol}}
\begin{quote}

\sphinxAtStartPar
The relative tolerance for degradation.

\sphinxAtStartPar
Optional, default is 0.0.
\end{quote}
\end{quote}
\end{quote}

\subsubsection{Model.setMObjective()}
\label{\detokenize{pyapiref:model-setmobjective}}\begin{quote}

\sphinxAtStartPar
\sphinxstylestrong{Synopsis}
\begin{quote}

\sphinxAtStartPar
\sphinxcode{\sphinxupquote{setMObjective(Q, c, constant, xQ\_L=None, xQ\_R=None, xc=None, sense=None)}}
\end{quote}

\sphinxAtStartPar
\sphinxstylestrong{Description}
\begin{quote}

\sphinxAtStartPar
Set the secondary objective of the model by matrix modeling. Objective functions of the form \(x_{Q_L} Q x_{Q_R} + c x_c + constant\) can be added.

\sphinxAtStartPar
Even more convenient is to generate a {\hyperref[\detokenize{pyapiref:chappyapi-mquadexpr}]{\sphinxcrossref{\DUrole{std,std-ref}{MQuadExpr Class}}}} object via matrix multiplication, available as the input to \sphinxcode{\sphinxupquote{setObjective()}} to set the objective function.
\end{quote}

\sphinxAtStartPar
\sphinxstylestrong{Arguments}
\begin{quote}

\sphinxAtStartPar
\sphinxcode{\sphinxupquote{Q}}
\begin{quote}

\sphinxAtStartPar
If the quadratic term is not empty, the parameter Q needs to be provided, which is a two\sphinxhyphen{}dimensional NumPy matrix, SciPy compressed sparse column matrix ( \sphinxcode{\sphinxupquote{csc\_matrix}} ) or compressed sparse row matrix ( \sphinxcode{\sphinxupquote{csr\_matrix}} ).
\end{quote}

\sphinxAtStartPar
\sphinxcode{\sphinxupquote{c}}
\begin{quote}

\sphinxAtStartPar
If the item is non\sphinxhyphen{}empty, you need to provide the parameter c, which is a one\sphinxhyphen{}dimensional NumPy array, or a Python list.
\end{quote}

\sphinxAtStartPar
\sphinxcode{\sphinxupquote{constant}}
\begin{quote}

\sphinxAtStartPar
Constant term, usually a floating point number.
\end{quote}

\sphinxAtStartPar
\sphinxcode{\sphinxupquote{xQ\_L}}
\begin{quote}

\sphinxAtStartPar
The variable on the left side of the quadratic term, can be a {\hyperref[\detokenize{pyapiref:chappyapi-mvar}]{\sphinxcrossref{\DUrole{std,std-ref}{MVar Class}}}} object, {\hyperref[\detokenize{pyapiref:chappyapi-vararray}]{\sphinxcrossref{\DUrole{std,std-ref}{VarArray Class}}}} object,
list, dictionary or {\hyperref[\detokenize{pyapiref:chappyapi-util-tupledict}]{\sphinxcrossref{\DUrole{std,std-ref}{tupledict Class}}}} object.

\sphinxAtStartPar
If empty, all variables in the model are taken.
\end{quote}

\sphinxAtStartPar
\sphinxcode{\sphinxupquote{xQ\_R}}
\begin{quote}

\sphinxAtStartPar
The variable on the right side of the quadratic term, can be a {\hyperref[\detokenize{pyapiref:chappyapi-mvar}]{\sphinxcrossref{\DUrole{std,std-ref}{MVar Class}}}} object, {\hyperref[\detokenize{pyapiref:chappyapi-vararray}]{\sphinxcrossref{\DUrole{std,std-ref}{VarArray Class}}}} object,
list, dictionary or {\hyperref[\detokenize{pyapiref:chappyapi-util-tupledict}]{\sphinxcrossref{\DUrole{std,std-ref}{tupledict Class}}}} object.

\sphinxAtStartPar
If empty, all variables in the model are taken.
\end{quote}

\sphinxAtStartPar
\sphinxcode{\sphinxupquote{xc}}
\begin{quote}

\sphinxAtStartPar
The variable corresponding to the linear term can be a {\hyperref[\detokenize{pyapiref:chappyapi-mvar}]{\sphinxcrossref{\DUrole{std,std-ref}{MVar Class}}}} object, {\hyperref[\detokenize{pyapiref:chappyapi-vararray}]{\sphinxcrossref{\DUrole{std,std-ref}{VarArray Class}}}} object,
list, dictionary or {\hyperref[\detokenize{pyapiref:chappyapi-util-tupledict}]{\sphinxcrossref{\DUrole{std,std-ref}{tupledict Class}}}} object.

\sphinxAtStartPar
If it is empty, but the parameter c is not empty, all variables in the model are taken.
\end{quote}

\sphinxAtStartPar
\sphinxcode{\sphinxupquote{sense}}
\begin{quote}

\sphinxAtStartPar
The optimization direction of the objective function. Optional parameter, the default is \sphinxcode{\sphinxupquote{None}}, which means that the optimization direction of the model will not be changed.
The current optimization direction of the model is viewed through the property \sphinxcode{\sphinxupquote{ObjSense}}.
See {\hyperref[\detokenize{constant:chapconst-sense}]{\sphinxcrossref{\DUrole{std,std-ref}{Optimization Direction}}}} for possible values.
\end{quote}
\end{quote}

\sphinxAtStartPar
\sphinxstylestrong{Example}

\begin{sphinxVerbatim}[commandchars=\\\{\}]
\PYG{n}{Q} \PYG{o}{=} \PYG{n}{np}\PYG{o}{.}\PYG{n}{full}\PYG{p}{(}\PYG{p}{(}\PYG{l+m+mi}{3}\PYG{p}{,} \PYG{l+m+mi}{3}\PYG{p}{)}\PYG{p}{,} \PYG{l+m+mi}{1}\PYG{p}{)}
\PYG{n}{mx} \PYG{o}{=} \PYG{n}{model}\PYG{o}{.}\PYG{n}{addMVar}\PYG{p}{(}\PYG{l+m+mi}{3}\PYG{p}{,} \PYG{n}{nameprefix}\PYG{o}{=}\PYG{l+s+s2}{\PYGZdq{}}\PYG{l+s+s2}{mx}\PYG{l+s+s2}{\PYGZdq{}}\PYG{p}{)}
\PYG{n}{my} \PYG{o}{=} \PYG{n}{model}\PYG{o}{.}\PYG{n}{addVars}\PYG{p}{(}\PYG{l+m+mi}{3}\PYG{p}{,} \PYG{n}{nameprefix}\PYG{o}{=}\PYG{l+s+s2}{\PYGZdq{}}\PYG{l+s+s2}{my}\PYG{l+s+s2}{\PYGZdq{}}\PYG{p}{)}
\PYG{n}{mqc} \PYG{o}{=} \PYG{n}{model}\PYG{o}{.}\PYG{n}{setMObjective}\PYG{p}{(}\PYG{n}{Q}\PYG{p}{,} \PYG{k+kc}{None}\PYG{p}{,} \PYG{l+m+mf}{0.0}\PYG{p}{,} \PYG{n}{mx}\PYG{p}{,} \PYG{n}{my}\PYG{p}{,} \PYG{k+kc}{None}\PYG{p}{,} \PYG{n}{sense}\PYG{o}{=}\PYG{n}{COPT}\PYG{o}{.}\PYG{n}{MINIMIZE}\PYG{p}{)}
\end{sphinxVerbatim}
\end{quote}

\subsubsection{Model.setObjSense()}
\label{\detokenize{pyapiref:model-setobjsense}}\begin{quote}

\sphinxAtStartPar
\sphinxstylestrong{Synopsis}
\begin{quote}

\sphinxAtStartPar
\sphinxcode{\sphinxupquote{setObjSense(sense)}}
\end{quote}

\sphinxAtStartPar
\sphinxstylestrong{Description}
\begin{quote}

\sphinxAtStartPar
Set optimization sense.
\end{quote}

\sphinxAtStartPar
\sphinxstylestrong{Arguments}
\begin{quote}

\sphinxAtStartPar
\sphinxcode{\sphinxupquote{sense}}
\begin{quote}

\sphinxAtStartPar
Optimization sense. Please refer to {\hyperref[\detokenize{constant:chapconst-sense}]{\sphinxcrossref{\DUrole{std,std-ref}{Optimization directions}}}} for possible values.
\end{quote}
\end{quote}

\sphinxAtStartPar
\sphinxstylestrong{Example}
\end{quote}

\begin{sphinxVerbatim}[commandchars=\\\{\}]
\PYG{c+c1}{\PYGZsh{} Set optimization sense as maximization}
\PYG{n}{m}\PYG{o}{.}\PYG{n}{setObjSense}\PYG{p}{(}\PYG{n}{COPT}\PYG{o}{.}\PYG{n}{MAXIMIZE}\PYG{p}{)}
\end{sphinxVerbatim}

\subsubsection{Model.setObjConst()}
\label{\detokenize{pyapiref:model-setobjconst}}\begin{quote}

\sphinxAtStartPar
\sphinxstylestrong{Synopsis}
\begin{quote}

\sphinxAtStartPar
\sphinxcode{\sphinxupquote{setObjConst(const)}}
\end{quote}

\sphinxAtStartPar
\sphinxstylestrong{Description}
\begin{quote}

\sphinxAtStartPar
Set constant objective offset.
\end{quote}

\sphinxAtStartPar
\sphinxstylestrong{Arguments}
\begin{quote}

\sphinxAtStartPar
\sphinxcode{\sphinxupquote{const}}
\begin{quote}

\sphinxAtStartPar
Constant objective offset.
\end{quote}
\end{quote}

\sphinxAtStartPar
\sphinxstylestrong{Example}
\end{quote}

\begin{sphinxVerbatim}[commandchars=\\\{\}]
\PYG{c+c1}{\PYGZsh{} Set constant objective offset 1.0}
\PYG{n}{m}\PYG{o}{.}\PYG{n}{setObjConst}\PYG{p}{(}\PYG{l+m+mf}{1.0}\PYG{p}{)}
\end{sphinxVerbatim}

\subsubsection{Model.getObjective()}
\label{\detokenize{pyapiref:model-getobjective}}\begin{quote}

\sphinxAtStartPar
\sphinxstylestrong{Synopsis}
\begin{quote}

\sphinxAtStartPar
\sphinxcode{\sphinxupquote{getObjective()}}
\end{quote}

\sphinxAtStartPar
\sphinxstylestrong{Description}
\begin{quote}

\sphinxAtStartPar
Retrieve current model objective. Return a {\hyperref[\detokenize{pyapiref:chappyapi-linexpr}]{\sphinxcrossref{\DUrole{std,std-ref}{LinExpr Class}}}} object.
\end{quote}

\sphinxAtStartPar
\sphinxstylestrong{Example}
\end{quote}

\begin{sphinxVerbatim}[commandchars=\\\{\}]
\PYG{c+c1}{\PYGZsh{} Retrieve the optimization objective.}
\PYG{n}{obj} \PYG{o}{=} \PYG{n}{m}\PYG{o}{.}\PYG{n}{getObjective}\PYG{p}{(}\PYG{p}{)}
\end{sphinxVerbatim}

\subsubsection{Model.getObjectiveN()}
\label{\detokenize{pyapiref:model-getobjectiven}}\begin{quote}

\sphinxAtStartPar
\sphinxstylestrong{Synopsis}
\begin{quote}

\sphinxAtStartPar
\sphinxcode{\sphinxupquote{getObjectiveN(idx)}}
\end{quote}

\sphinxAtStartPar
\sphinxstylestrong{Description}
\begin{quote}

\sphinxAtStartPar
Retrieve the expression of the objective function with the specified index in multi\sphinxhyphen{}objective optimization.
\end{quote}

\sphinxAtStartPar
\sphinxstylestrong{Arguments}
\begin{quote}

\sphinxAtStartPar
\sphinxcode{\sphinxupquote{idx}}
\begin{quote}

\sphinxAtStartPar
The index of the objective function.
\end{quote}
\end{quote}

\sphinxAtStartPar
\sphinxstylestrong{Return Value}
\begin{quote}

\sphinxAtStartPar
Returns a {\hyperref[\detokenize{pyapiref:chappyapi-linexpr}]{\sphinxcrossref{\DUrole{std,std-ref}{LinExpr Class}}}} object.
\end{quote}
\end{quote}

\subsubsection{Model.delObjN()}
\label{\detokenize{pyapiref:model-delobjn}}\begin{quote}

\sphinxAtStartPar
\sphinxstylestrong{Synopsis}
\begin{quote}

\sphinxAtStartPar
\sphinxcode{\sphinxupquote{delObjN(idx)}}
\end{quote}

\sphinxAtStartPar
\sphinxstylestrong{Description}
\begin{quote}

\sphinxAtStartPar
Delete the linear part of the objective
with the specified index from the multi\sphinxhyphen{}objective model.
\end{quote}

\sphinxAtStartPar
\sphinxstylestrong{Arguments}
\begin{quote}

\sphinxAtStartPar
\sphinxcode{\sphinxupquote{idx}}
\begin{quote}

\sphinxAtStartPar
The index of the objective function.
\end{quote}
\end{quote}
\end{quote}

\subsubsection{Model.delQuadObj()}
\label{\detokenize{pyapiref:model-delquadobj}}\begin{quote}

\sphinxAtStartPar
\sphinxstylestrong{Synopsis}
\begin{quote}

\sphinxAtStartPar
\sphinxcode{\sphinxupquote{delQuadObj()}}
\end{quote}

\sphinxAtStartPar
\sphinxstylestrong{Description}
\begin{quote}

\sphinxAtStartPar
Deletes the quadratic terms from the quadratic objective function.
\end{quote}

\sphinxAtStartPar
\sphinxstylestrong{Example}
\end{quote}

\begin{sphinxVerbatim}[commandchars=\\\{\}]
\PYG{c+c1}{\PYGZsh{} Deletes the quadratic terms from the quadratic objective function}
\PYG{n}{m}\PYG{o}{.}\PYG{n}{delQuadObj}\PYG{p}{(}\PYG{p}{)}
\end{sphinxVerbatim}

\subsubsection{Model.delNlObj()}
\label{\detokenize{pyapiref:model-delnlobj}}\begin{quote}

\sphinxAtStartPar
\sphinxstylestrong{Synopsis}
\begin{quote}

\sphinxAtStartPar
\sphinxcode{\sphinxupquote{delNlObj()}}
\end{quote}

\sphinxAtStartPar
\sphinxstylestrong{Description}
\begin{quote}

\sphinxAtStartPar
Delete all nonlinear terms from the nonlinear objective function.
\end{quote}

\sphinxAtStartPar
\sphinxstylestrong{Example}
\end{quote}

\begin{sphinxVerbatim}[commandchars=\\\{\}]
\PYG{c+c1}{\PYGZsh{} Delete nonlinear terms from the objective}
\PYG{n}{m}\PYG{o}{.}\PYG{n}{delNlObj}\PYG{p}{(}\PYG{p}{)}
\end{sphinxVerbatim}

\subsubsection{Model.delPsdObj()}
\label{\detokenize{pyapiref:model-delpsdobj}}\begin{quote}

\sphinxAtStartPar
\sphinxstylestrong{Synopsis}
\begin{quote}

\sphinxAtStartPar
\sphinxcode{\sphinxupquote{delPsdObj()}}
\end{quote}

\sphinxAtStartPar
\sphinxstylestrong{Description}
\begin{quote}

\sphinxAtStartPar
Delete the positive semi\sphinxhyphen{}definite terms from the objective function
\end{quote}

\sphinxAtStartPar
\sphinxstylestrong{Example}
\end{quote}

\begin{sphinxVerbatim}[commandchars=\\\{\}]
\PYG{c+c1}{\PYGZsh{} Delete the positive semi\PYGZhy{}definite terms from the objective function}
\PYG{n}{m}\PYG{o}{.}\PYG{n}{delPsdObj}\PYG{p}{(}\PYG{p}{)}
\end{sphinxVerbatim}

\subsubsection{Model.getCol()}
\label{\detokenize{pyapiref:model-getcol}}\begin{quote}

\sphinxAtStartPar
\sphinxstylestrong{Synopsis}
\begin{quote}

\sphinxAtStartPar
\sphinxcode{\sphinxupquote{getCol(var)}}
\end{quote}

\sphinxAtStartPar
\sphinxstylestrong{Description}
\begin{quote}

\sphinxAtStartPar
Retrieve the list of constraints in which a variable participates,
Return value is a {\hyperref[\detokenize{pyapiref:chappyapi-column}]{\sphinxcrossref{\DUrole{std,std-ref}{Column Class}}}}  object that captures the set of constraints in which the variable participates.
\end{quote}

\sphinxAtStartPar
\sphinxstylestrong{Example}
\end{quote}

\begin{sphinxVerbatim}[commandchars=\\\{\}]
\PYG{c+c1}{\PYGZsh{} Get column that captures the set of constraints in which x participates.}
\PYG{n}{col} \PYG{o}{=} \PYG{n}{m}\PYG{o}{.}\PYG{n}{getCol}\PYG{p}{(}\PYG{n}{x}\PYG{p}{)}
\end{sphinxVerbatim}

\subsubsection{Model.getRow()}
\label{\detokenize{pyapiref:model-getrow}}\begin{quote}

\sphinxAtStartPar
\sphinxstylestrong{Synopsis}
\begin{quote}

\sphinxAtStartPar
\sphinxcode{\sphinxupquote{getRow(constr)}}
\end{quote}
\begin{description}
\sphinxlineitem{\sphinxstylestrong{Description}}
\sphinxAtStartPar
Retrieve the list of variables that participate in a specific constraint and return a {\hyperref[\detokenize{pyapiref:chappyapi-linexpr}]{\sphinxcrossref{\DUrole{std,std-ref}{LinExpr Class}}}} object.

\end{description}

\sphinxAtStartPar
\sphinxstylestrong{Example}
\end{quote}

\begin{sphinxVerbatim}[commandchars=\\\{\}]
\PYG{c+c1}{\PYGZsh{} Return variables that participate in conx.}
\PYG{n}{linexpr} \PYG{o}{=} \PYG{n}{m}\PYG{o}{.}\PYG{n}{getRow}\PYG{p}{(}\PYG{n}{conx}\PYG{p}{)}
\end{sphinxVerbatim}

\subsubsection{Model.getQuadRow()}
\label{\detokenize{pyapiref:model-getquadrow}}\begin{quote}

\sphinxAtStartPar
\sphinxstylestrong{Synopsis}
\begin{quote}

\sphinxAtStartPar
\sphinxcode{\sphinxupquote{getQuadRow(qconstr)}}
\end{quote}

\sphinxAtStartPar
\sphinxstylestrong{Description}
\begin{quote}

\sphinxAtStartPar
Retrieve the list of variables that participate in a specific quadratic constraint and return a {\hyperref[\detokenize{pyapiref:chappyapi-quadexpr}]{\sphinxcrossref{\DUrole{std,std-ref}{QuadExpr Class}}}} object.
\end{quote}

\sphinxAtStartPar
\sphinxstylestrong{Example}
\end{quote}

\begin{sphinxVerbatim}[commandchars=\\\{\}]
\PYG{c+c1}{\PYGZsh{} Return variables that participate in qconx}
\PYG{n}{quadexpr} \PYG{o}{=} \PYG{n}{m}\PYG{o}{.}\PYG{n}{getQuadRow}\PYG{p}{(}\PYG{n}{qconx}\PYG{p}{)}
\end{sphinxVerbatim}

\subsubsection{Model.getPsdRow()}
\label{\detokenize{pyapiref:model-getpsdrow}}\begin{quote}

\sphinxAtStartPar
\sphinxstylestrong{Synopsis}
\begin{quote}

\sphinxAtStartPar
\sphinxcode{\sphinxupquote{getPsdRow(constr)}}
\end{quote}

\sphinxAtStartPar
\sphinxstylestrong{Description}
\begin{quote}

\sphinxAtStartPar
Retrieve the list of variables that participate in a specific positive semi\sphinxhyphen{}definite constraint and return a {\hyperref[\detokenize{pyapiref:chappyapi-psdexpr}]{\sphinxcrossref{\DUrole{std,std-ref}{PsdExpr Class}}}} object.
\end{quote}

\sphinxAtStartPar
\sphinxstylestrong{Example}
\end{quote}

\begin{sphinxVerbatim}[commandchars=\\\{\}]
\PYG{c+c1}{\PYGZsh{} Retrieve the row corresponding to a specific positive semi\PYGZhy{}definite constraint}
\PYG{n}{psdexpr} \PYG{o}{=} \PYG{n}{m}\PYG{o}{.}\PYG{n}{getPsdRow}\PYG{p}{(}\PYG{n}{psdcon}\PYG{p}{)}
\end{sphinxVerbatim}

\subsubsection{Model.getNlRow()}
\label{\detokenize{pyapiref:model-getnlrow}}\begin{quote}

\sphinxAtStartPar
\sphinxstylestrong{Synopsis}
\begin{quote}

\sphinxAtStartPar
\sphinxcode{\sphinxupquote{getNlRow(constr)}}
\end{quote}

\sphinxAtStartPar
\sphinxstylestrong{Description}
\begin{quote}

\sphinxAtStartPar
Retrieve the row expression associated with the specified nonlinear constraint.
Returns a {\hyperref[\detokenize{pyapiref:chappyapi-nlexpr}]{\sphinxcrossref{\DUrole{std,std-ref}{NlExpr Class}}}} object.
\end{quote}

\sphinxAtStartPar
\sphinxstylestrong{Example}
\end{quote}

\begin{sphinxVerbatim}[commandchars=\\\{\}]
\PYG{c+c1}{\PYGZsh{} Retrieve the row expression of nonlinear constraint}
\PYG{n}{nlexpr} \PYG{o}{=} \PYG{n}{m}\PYG{o}{.}\PYG{n}{getNlRow}\PYG{p}{(}\PYG{n}{nlcon}\PYG{p}{)}
\end{sphinxVerbatim}

\subsubsection{Model.getVar()}
\label{\detokenize{pyapiref:model-getvar}}\begin{quote}

\sphinxAtStartPar
\sphinxstylestrong{Synopsis}
\begin{quote}

\sphinxAtStartPar
\sphinxcode{\sphinxupquote{getVar(idx)}}
\end{quote}

\sphinxAtStartPar
\sphinxstylestrong{Description}
\begin{quote}

\sphinxAtStartPar
Retrieve a variable according to its index in the coefficient matrix. Return a {\hyperref[\detokenize{pyapiref:chappyapi-var}]{\sphinxcrossref{\DUrole{std,std-ref}{Var Class}}}} object.
\end{quote}

\sphinxAtStartPar
\sphinxstylestrong{Arguments}
\begin{quote}

\sphinxAtStartPar
\sphinxcode{\sphinxupquote{idx}}
\begin{quote}

\sphinxAtStartPar
Index of the desired variable in the coefficient matrix, starting with \sphinxcode{\sphinxupquote{0}}.
\end{quote}
\end{quote}

\sphinxAtStartPar
\sphinxstylestrong{Example}
\end{quote}

\begin{sphinxVerbatim}[commandchars=\\\{\}]
\PYG{c+c1}{\PYGZsh{} Retrive variable with indice of 1.}
 \PYG{n}{x} \PYG{o}{=} \PYG{n}{m}\PYG{o}{.}\PYG{n}{getVar}\PYG{p}{(}\PYG{l+m+mi}{1}\PYG{p}{)}
\end{sphinxVerbatim}

\subsubsection{Model.getVarByName()}
\label{\detokenize{pyapiref:model-getvarbyname}}\begin{quote}

\sphinxAtStartPar
\sphinxstylestrong{Synopsis}
\begin{quote}

\sphinxAtStartPar
\sphinxcode{\sphinxupquote{getVarByName(name)}}
\end{quote}

\sphinxAtStartPar
\sphinxstylestrong{Description}
\begin{quote}

\sphinxAtStartPar
Retrieves a variable by name. Return a {\hyperref[\detokenize{pyapiref:chappyapi-var}]{\sphinxcrossref{\DUrole{std,std-ref}{Var Class}}}} object.
\end{quote}

\sphinxAtStartPar
\sphinxstylestrong{Arguments}
\begin{quote}

\sphinxAtStartPar
\sphinxcode{\sphinxupquote{name}}
\begin{quote}

\sphinxAtStartPar
Name of the desired variable.
\end{quote}
\end{quote}

\sphinxAtStartPar
\sphinxstylestrong{Example}
\end{quote}

\begin{sphinxVerbatim}[commandchars=\\\{\}]
\PYG{c+c1}{\PYGZsh{}  Retrive variable with name \PYGZdq{}x\PYGZdq{}.}
\PYG{n}{x} \PYG{o}{=} \PYG{n}{m}\PYG{o}{.}\PYG{n}{getVarByName}\PYG{p}{(}\PYG{l+s+s2}{\PYGZdq{}}\PYG{l+s+s2}{x}\PYG{l+s+s2}{\PYGZdq{}}\PYG{p}{)}
\end{sphinxVerbatim}

\subsubsection{Model.getVars()}
\label{\detokenize{pyapiref:model-getvars}}\begin{quote}

\sphinxAtStartPar
\sphinxstylestrong{Synopsis}
\begin{quote}

\sphinxAtStartPar
\sphinxcode{\sphinxupquote{getVars()}}
\end{quote}

\sphinxAtStartPar
\sphinxstylestrong{Description}
\begin{quote}

\sphinxAtStartPar
Retrieve all variables in the model. Return a  {\hyperref[\detokenize{pyapiref:chappyapi-vararray}]{\sphinxcrossref{\DUrole{std,std-ref}{VarArray Class}}}} object.
\end{quote}

\sphinxAtStartPar
\sphinxstylestrong{Example}
\end{quote}

\begin{sphinxVerbatim}[commandchars=\\\{\}]
\PYG{c+c1}{\PYGZsh{} Retrieve all variables in the model}
\PYG{n+nb}{vars} \PYG{o}{=} \PYG{n}{m}\PYG{o}{.}\PYG{n}{getVars}\PYG{p}{(}\PYG{p}{)}
\end{sphinxVerbatim}

\subsubsection{Model.getConstr()}
\label{\detokenize{pyapiref:model-getconstr}}\begin{quote}

\sphinxAtStartPar
\sphinxstylestrong{Synopsis}
\begin{quote}

\sphinxAtStartPar
\sphinxcode{\sphinxupquote{getConstr(idx)}}
\end{quote}

\sphinxAtStartPar
\sphinxstylestrong{Description}
\begin{quote}

\sphinxAtStartPar
Retrieve a constraint by its indice in the coefficient matrix. Return a {\hyperref[\detokenize{pyapiref:chappyapi-constraint}]{\sphinxcrossref{\DUrole{std,std-ref}{Constraint Class}}}} object.
\end{quote}

\sphinxAtStartPar
\sphinxstylestrong{Arguments}
\begin{quote}

\sphinxAtStartPar
\sphinxcode{\sphinxupquote{idx}}
\begin{quote}

\sphinxAtStartPar
Index of the desired constraint in the coefficient matrix, starting with \sphinxcode{\sphinxupquote{0}}.
\end{quote}
\end{quote}

\sphinxAtStartPar
\sphinxstylestrong{Example}
\end{quote}

\begin{sphinxVerbatim}[commandchars=\\\{\}]
\PYG{c+c1}{\PYGZsh{} Retrive linear constraint with indice of 1.}
\PYG{n}{r} \PYG{o}{=} \PYG{n}{m}\PYG{o}{.}\PYG{n}{getConstr}\PYG{p}{(}\PYG{l+m+mi}{1}\PYG{p}{)}
\end{sphinxVerbatim}

\subsubsection{Model.getConstrByName()}
\label{\detokenize{pyapiref:model-getconstrbyname}}\begin{quote}

\sphinxAtStartPar
\sphinxstylestrong{Synopsis}
\begin{quote}

\sphinxAtStartPar
\sphinxcode{\sphinxupquote{getConstrByName(name)}}
\end{quote}

\sphinxAtStartPar
\sphinxstylestrong{Description}
\begin{quote}

\sphinxAtStartPar
Retrieves a linear constraint by name. Return a {\hyperref[\detokenize{pyapiref:chappyapi-constraint}]{\sphinxcrossref{\DUrole{std,std-ref}{Constraint Class}}}} object.
\end{quote}

\sphinxAtStartPar
\sphinxstylestrong{Arguments}
\begin{quote}

\sphinxAtStartPar
\sphinxcode{\sphinxupquote{name}}
\begin{quote}

\sphinxAtStartPar
The name of the constraint.
\end{quote}
\end{quote}

\sphinxAtStartPar
\sphinxstylestrong{Example}
\end{quote}

\begin{sphinxVerbatim}[commandchars=\\\{\}]
\PYG{c+c1}{\PYGZsh{} Retrieve linear constraint with name \PYGZdq{}r\PYGZdq{}.}
\PYG{n}{r} \PYG{o}{=} \PYG{n}{m}\PYG{o}{.}\PYG{n}{getConstrByName}\PYG{p}{(}\PYG{l+s+s2}{\PYGZdq{}}\PYG{l+s+s2}{r}\PYG{l+s+s2}{\PYGZdq{}}\PYG{p}{)}
\end{sphinxVerbatim}

\subsubsection{Model.getConstrs()}
\label{\detokenize{pyapiref:model-getconstrs}}\begin{quote}

\sphinxAtStartPar
\sphinxstylestrong{Synopsis}
\begin{quote}

\sphinxAtStartPar
\sphinxcode{\sphinxupquote{getConstrs()}}
\end{quote}

\sphinxAtStartPar
\sphinxstylestrong{Description}
\begin{quote}

\sphinxAtStartPar
Retrieve all constraints in the model. Return a {\hyperref[\detokenize{pyapiref:chappyapi-constrarray}]{\sphinxcrossref{\DUrole{std,std-ref}{ConstrArray Class}}}} object.
\end{quote}

\sphinxAtStartPar
\sphinxstylestrong{Example}
\end{quote}

\begin{sphinxVerbatim}[commandchars=\\\{\}]
\PYG{c+c1}{\PYGZsh{} Retrieve all constraints in the model}
\PYG{n}{cons} \PYG{o}{=} \PYG{n}{m}\PYG{o}{.}\PYG{n}{getConstrs}\PYG{p}{(}\PYG{p}{)}
\end{sphinxVerbatim}

\subsubsection{Model.getConstrBuilders()}
\label{\detokenize{pyapiref:model-getconstrbuilders}}\begin{quote}

\sphinxAtStartPar
\sphinxstylestrong{Synopsis}
\begin{quote}

\sphinxAtStartPar
\sphinxcode{\sphinxupquote{getConstrBuilders(constrs=None)}}
\end{quote}

\sphinxAtStartPar
\sphinxstylestrong{Description}
\begin{quote}

\sphinxAtStartPar
Retrieve linear constraint builders in current model.

\sphinxAtStartPar
If parameter \sphinxcode{\sphinxupquote{constrs}} is \sphinxcode{\sphinxupquote{None}}, then return a {\hyperref[\detokenize{pyapiref:chappyapi-constrbuilderarray}]{\sphinxcrossref{\DUrole{std,std-ref}{ConstrBuilderArray Class}}}} object composed of all
linear constraint builders.

\sphinxAtStartPar
If parameter \sphinxcode{\sphinxupquote{constrs}} is {\hyperref[\detokenize{pyapiref:chappyapi-constraint}]{\sphinxcrossref{\DUrole{std,std-ref}{Constraint Class}}}} object, then return the {\hyperref[\detokenize{pyapiref:chappyapi-constrbuilder}]{\sphinxcrossref{\DUrole{std,std-ref}{ConstrBuilder Class}}}}
object corresponding to the specific constraint.

\sphinxAtStartPar
If parameter \sphinxcode{\sphinxupquote{constrs}} is a list or a {\hyperref[\detokenize{pyapiref:chappyapi-constrarray}]{\sphinxcrossref{\DUrole{std,std-ref}{ConstrArray Class}}}} object, then return a {\hyperref[\detokenize{pyapiref:chappyapi-constrbuilderarray}]{\sphinxcrossref{\DUrole{std,std-ref}{ConstrBuilderArray Class}}}}
object composed of specified constraints’ builders.

\sphinxAtStartPar
If parameter \sphinxcode{\sphinxupquote{constrs}} is dictionary or {\hyperref[\detokenize{pyapiref:chappyapi-util-tupledict}]{\sphinxcrossref{\DUrole{std,std-ref}{tupledict Class}}}} object, then the indice of the specified
constraint is returned as key, the value is a {\hyperref[\detokenize{pyapiref:chappyapi-util-tupledict}]{\sphinxcrossref{\DUrole{std,std-ref}{tupledict Class}}}} object composed of the specified constraints’ builders.
\end{quote}

\sphinxAtStartPar
\sphinxstylestrong{Arguments}
\begin{quote}

\sphinxAtStartPar
\sphinxcode{\sphinxupquote{constrs}}
\begin{quote}

\sphinxAtStartPar
The specified linear constraint. Optional, \sphinxcode{\sphinxupquote{None}} by default.
\end{quote}
\end{quote}

\sphinxAtStartPar
\sphinxstylestrong{Example}
\end{quote}

\begin{sphinxVerbatim}[commandchars=\\\{\}]
\PYG{c+c1}{\PYGZsh{} Retrieve all of linear constraint builders.}
\PYG{n}{conbuilders} \PYG{o}{=} \PYG{n}{m}\PYG{o}{.}\PYG{n}{getConstrBuilders}\PYG{p}{(}\PYG{p}{)}
\PYG{c+c1}{\PYGZsh{} Retrive the builder corresponding to linear contstraint x.}
\PYG{n}{conbuilders} \PYG{o}{=} \PYG{n}{m}\PYG{o}{.}\PYG{n}{getConstrBuilders}\PYG{p}{(}\PYG{n}{x}\PYG{p}{)}
\PYG{c+c1}{\PYGZsh{} Retrieve builders corresponding to linear constraint x and y.}
\PYG{n}{conbuilders} \PYG{o}{=} \PYG{n}{m}\PYG{o}{.}\PYG{n}{getConstrBuilders}\PYG{p}{(}\PYG{p}{[}\PYG{n}{x}\PYG{p}{,} \PYG{n}{y}\PYG{p}{]}\PYG{p}{)}
\PYG{c+c1}{\PYGZsh{} Retrieve builders corresponding to linear constraint in tupledict object xx.}
\PYG{n}{conbuilders} \PYG{o}{=} \PYG{n}{m}\PYG{o}{.}\PYG{n}{getConstrBuilders}\PYG{p}{(}\PYG{n}{xx}\PYG{p}{)}
\end{sphinxVerbatim}

\subsubsection{Model.getNlConstr()}
\label{\detokenize{pyapiref:model-getnlconstr}}\begin{quote}

\sphinxAtStartPar
\sphinxstylestrong{Synopsis}
\begin{quote}

\sphinxAtStartPar
\sphinxcode{\sphinxupquote{getNlConstr(idx)}}
\end{quote}

\sphinxAtStartPar
\sphinxstylestrong{Description}
\begin{quote}

\sphinxAtStartPar
Retrieve a nonlinear constraint by its index in the model.
Return a {\hyperref[\detokenize{pyapiref:chappyapi-nlconstraint}]{\sphinxcrossref{\DUrole{std,std-ref}{NlConstraint Class}}}} object.
\end{quote}

\sphinxAtStartPar
\sphinxstylestrong{Arguments}
\begin{quote}

\sphinxAtStartPar
\sphinxcode{\sphinxupquote{idx}}
\begin{quote}

\sphinxAtStartPar
Index of the nonlinear constraint in the model, starting with \sphinxcode{\sphinxupquote{0}}.
\end{quote}
\end{quote}

\sphinxAtStartPar
\sphinxstylestrong{Example}
\end{quote}

\begin{sphinxVerbatim}[commandchars=\\\{\}]
\PYG{c+c1}{\PYGZsh{} Retrieve nonlinear constraint with index 1}
\PYG{n}{r} \PYG{o}{=} \PYG{n}{m}\PYG{o}{.}\PYG{n}{getNlConstr}\PYG{p}{(}\PYG{l+m+mi}{1}\PYG{p}{)}
\end{sphinxVerbatim}

\subsubsection{Model.getNlConstrByName()}
\label{\detokenize{pyapiref:model-getnlconstrbyname}}\begin{quote}

\sphinxAtStartPar
\sphinxstylestrong{Synopsis}
\begin{quote}

\sphinxAtStartPar
\sphinxcode{\sphinxupquote{getNlConstrByName(name)}}
\end{quote}

\sphinxAtStartPar
\sphinxstylestrong{Description}
\begin{quote}

\sphinxAtStartPar
Retrieve a nonlinear constraint by its name.
Return a {\hyperref[\detokenize{pyapiref:chappyapi-nlconstraint}]{\sphinxcrossref{\DUrole{std,std-ref}{NlConstraint Class}}}} object.
\end{quote}

\sphinxAtStartPar
\sphinxstylestrong{Arguments}
\begin{quote}

\sphinxAtStartPar
\sphinxcode{\sphinxupquote{name}}
\begin{quote}

\sphinxAtStartPar
Name of the nonlinear constraint.
\end{quote}
\end{quote}

\sphinxAtStartPar
\sphinxstylestrong{Example}
\end{quote}

\begin{sphinxVerbatim}[commandchars=\\\{\}]
\PYG{c+c1}{\PYGZsh{} Retrieve nonlinear constraint with name \PYGZdq{}r\PYGZdq{}}
\PYG{n}{r} \PYG{o}{=} \PYG{n}{m}\PYG{o}{.}\PYG{n}{getNlConstrByName}\PYG{p}{(}\PYG{l+s+s2}{\PYGZdq{}}\PYG{l+s+s2}{r}\PYG{l+s+s2}{\PYGZdq{}}\PYG{p}{)}
\end{sphinxVerbatim}

\subsubsection{Model.getNlConstrs()}
\label{\detokenize{pyapiref:model-getnlconstrs}}\begin{quote}

\sphinxAtStartPar
\sphinxstylestrong{Synopsis}
\begin{quote}

\sphinxAtStartPar
\sphinxcode{\sphinxupquote{getNlConstrs()}}
\end{quote}

\sphinxAtStartPar
\sphinxstylestrong{Description}
\begin{quote}

\sphinxAtStartPar
Retrieve all nonlinear constraints in the model.
Return a {\hyperref[\detokenize{pyapiref:chappyapi-nlconstrarray}]{\sphinxcrossref{\DUrole{std,std-ref}{NlConstrArray Class}}}} object.
\end{quote}

\sphinxAtStartPar
\sphinxstylestrong{Example}
\end{quote}

\begin{sphinxVerbatim}[commandchars=\\\{\}]
\PYG{c+c1}{\PYGZsh{} Retrieve all nonlinear constraints in the model}
\PYG{n}{cons} \PYG{o}{=} \PYG{n}{m}\PYG{o}{.}\PYG{n}{getNlConstrs}\PYG{p}{(}\PYG{p}{)}
\end{sphinxVerbatim}

\subsubsection{Model.getNlConstrBuilders()}
\label{\detokenize{pyapiref:model-getnlconstrbuilders}}\begin{quote}

\sphinxAtStartPar
\sphinxstylestrong{Synopsis}
\begin{quote}

\sphinxAtStartPar
\sphinxcode{\sphinxupquote{getNlConstrBuilders(constrs=None)}}
\end{quote}

\sphinxAtStartPar
\sphinxstylestrong{Description}
\begin{quote}

\sphinxAtStartPar
Retrieve nonlinear constraint builders in the current model.

\sphinxAtStartPar
If parameter \sphinxcode{\sphinxupquote{constrs}} is \sphinxcode{\sphinxupquote{None}},
return a {\hyperref[\detokenize{pyapiref:chappyapi-nlconstrbuilderarray}]{\sphinxcrossref{\DUrole{std,std-ref}{NlConstrBuilderArray Class}}}} object composed of all
nonlinear constraint builders.

\sphinxAtStartPar
If parameter \sphinxcode{\sphinxupquote{constrs}} is a {\hyperref[\detokenize{pyapiref:chappyapi-nlconstraint}]{\sphinxcrossref{\DUrole{std,std-ref}{NlConstraint Class}}}} object,
return the {\hyperref[\detokenize{pyapiref:chappyapi-nlconstrbuilder}]{\sphinxcrossref{\DUrole{std,std-ref}{NlConstrBuilder Class}}}} object corresponding to
the specified constraint.

\sphinxAtStartPar
If parameter \sphinxcode{\sphinxupquote{constrs}} is a list or a {\hyperref[\detokenize{pyapiref:chappyapi-nlconstrarray}]{\sphinxcrossref{\DUrole{std,std-ref}{NlConstrArray Class}}}} object,
return a {\hyperref[\detokenize{pyapiref:chappyapi-nlconstrbuilderarray}]{\sphinxcrossref{\DUrole{std,std-ref}{NlConstrBuilderArray Class}}}} object composed of the specified constraints’ builders.

\sphinxAtStartPar
If parameter \sphinxcode{\sphinxupquote{constrs}} is a dictionary or {\hyperref[\detokenize{pyapiref:chappyapi-util-tupledict}]{\sphinxcrossref{\DUrole{std,std-ref}{tupledict Class}}}} object,
return a {\hyperref[\detokenize{pyapiref:chappyapi-util-tupledict}]{\sphinxcrossref{\DUrole{std,std-ref}{tupledict Class}}}} where keys are the same as the input, and values are the corresponding nonlinear constraint builders.
\end{quote}

\sphinxAtStartPar
\sphinxstylestrong{Arguments}
\begin{quote}

\sphinxAtStartPar
\sphinxcode{\sphinxupquote{constrs}}
\begin{quote}

\sphinxAtStartPar
The specified nonlinear constraint(s).

\sphinxAtStartPar
Optional, \sphinxcode{\sphinxupquote{None}} by default.
\end{quote}
\end{quote}

\sphinxAtStartPar
\sphinxstylestrong{Example}
\end{quote}

\begin{sphinxVerbatim}[commandchars=\\\{\}]
\PYG{c+c1}{\PYGZsh{} Retrieve all nonlinear constraint builders}
\PYG{n}{conbuilders} \PYG{o}{=} \PYG{n}{m}\PYG{o}{.}\PYG{n}{getNlConstrBuilders}\PYG{p}{(}\PYG{p}{)}
\PYG{c+c1}{\PYGZsh{} Retrieve the builder corresponding to nonlinear constraint nl}
\PYG{n}{conbuilder} \PYG{o}{=} \PYG{n}{m}\PYG{o}{.}\PYG{n}{getNlConstrBuilders}\PYG{p}{(}\PYG{n}{nl}\PYG{p}{)}
\PYG{c+c1}{\PYGZsh{} Retrieve builders corresponding to nonlinear constraints nl1 and nl2}
\PYG{n}{conbuilders} \PYG{o}{=} \PYG{n}{m}\PYG{o}{.}\PYG{n}{getNlConstrBuilders}\PYG{p}{(}\PYG{p}{[}\PYG{n}{nl1}\PYG{p}{,} \PYG{n}{nl2}\PYG{p}{]}\PYG{p}{)}
\PYG{c+c1}{\PYGZsh{} Retrieve builders corresponding to nonlinear constraints in a tupledict}
\PYG{n}{conbuilders} \PYG{o}{=} \PYG{n}{m}\PYG{o}{.}\PYG{n}{getNlConstrBuilders}\PYG{p}{(}\PYG{n}{td}\PYG{p}{)}
\end{sphinxVerbatim}

\subsubsection{Model.getSOS(sos)}
\label{\detokenize{pyapiref:model-getsos-sos}}\begin{quote}

\sphinxAtStartPar
\sphinxstylestrong{Synopsis}
\begin{quote}

\sphinxAtStartPar
\sphinxcode{\sphinxupquote{getSOS(sos)}}
\end{quote}

\sphinxAtStartPar
\sphinxstylestrong{Description}
\begin{quote}

\sphinxAtStartPar
Retrieve the SOS constraint builder corresponding to specific SOS constraint. Return a {\hyperref[\detokenize{pyapiref:chappyapi-sosbuilder}]{\sphinxcrossref{\DUrole{std,std-ref}{SOSBuilder Class}}}} object
\end{quote}

\sphinxAtStartPar
\sphinxstylestrong{Arguments}
\begin{quote}

\sphinxAtStartPar
\sphinxcode{\sphinxupquote{sos}}
\begin{quote}

\sphinxAtStartPar
The specified SOS constraint.
\end{quote}
\end{quote}

\sphinxAtStartPar
\sphinxstylestrong{Example}
\end{quote}

\begin{sphinxVerbatim}[commandchars=\\\{\}]
\PYG{c+c1}{\PYGZsh{} Retrieve the builder corresponding to SOS constraint sosx.}
\PYG{n}{sosbuilder} \PYG{o}{=} \PYG{n}{m}\PYG{o}{.}\PYG{n}{getSOS}\PYG{p}{(}\PYG{n}{sosx}\PYG{p}{)}
\end{sphinxVerbatim}

\subsubsection{Model.getSOSs()}
\label{\detokenize{pyapiref:model-getsoss}}\begin{quote}

\sphinxAtStartPar
\sphinxstylestrong{Synopsis}
\begin{quote}

\sphinxAtStartPar
\sphinxcode{\sphinxupquote{getSOSs()}}
\end{quote}

\sphinxAtStartPar
\sphinxstylestrong{Description}
\begin{quote}

\sphinxAtStartPar
Retrieve all SOS constraints in model and return a {\hyperref[\detokenize{pyapiref:chappyapi-sosarray}]{\sphinxcrossref{\DUrole{std,std-ref}{SOSArray Class}}}} object.
\end{quote}

\sphinxAtStartPar
\sphinxstylestrong{Example}
\end{quote}

\begin{sphinxVerbatim}[commandchars=\\\{\}]
\PYG{c+c1}{\PYGZsh{}  Retrieve all SOS constraints in model.}
\PYG{n}{soss} \PYG{o}{=} \PYG{n}{m}\PYG{o}{.}\PYG{n}{getSOSs}\PYG{p}{(}\PYG{p}{)}
\end{sphinxVerbatim}

\subsubsection{Model.getSOSBuilders()}
\label{\detokenize{pyapiref:model-getsosbuilders}}\begin{quote}

\sphinxAtStartPar
\sphinxstylestrong{Synopsis}
\begin{quote}

\sphinxAtStartPar
\sphinxcode{\sphinxupquote{getSOSBuilders(soss=None)}}
\end{quote}

\sphinxAtStartPar
\sphinxstylestrong{Description}
\begin{quote}

\sphinxAtStartPar
Retrieve the SOS constraint builder corresponding to the specified SOS constraint.

\sphinxAtStartPar
If parameter \sphinxcode{\sphinxupquote{soss}} is \sphinxcode{\sphinxupquote{None}}, then return a {\hyperref[\detokenize{pyapiref:chappyapi-sosbuilderarray}]{\sphinxcrossref{\DUrole{std,std-ref}{SOSBuilderArray Class}}}} object consisting of builders
corresponding to all SOS constraints.

\sphinxAtStartPar
If parameter \sphinxcode{\sphinxupquote{soss}} is {\hyperref[\detokenize{pyapiref:chappyapi-sos}]{\sphinxcrossref{\DUrole{std,std-ref}{SOS Class}}}} object, then return a {\hyperref[\detokenize{pyapiref:chappyapi-sosbuilder}]{\sphinxcrossref{\DUrole{std,std-ref}{SOSBuilder Class}}}}
corresponding to the specified SOS constraint.

\sphinxAtStartPar
If parameter \sphinxcode{\sphinxupquote{soss}} is list or {\hyperref[\detokenize{pyapiref:chappyapi-sosarray}]{\sphinxcrossref{\DUrole{std,std-ref}{SOSArray Class}}}} object, then return a {\hyperref[\detokenize{pyapiref:chappyapi-sosbuilderarray}]{\sphinxcrossref{\DUrole{std,std-ref}{SOSBuilderArray Class}}}}
object consisting of builders corresponding to the specific SOS constraints.
\end{quote}

\sphinxAtStartPar
\sphinxstylestrong{Arguments}
\begin{quote}
\begin{description}
\sphinxlineitem{\sphinxcode{\sphinxupquote{soss}}}
\sphinxAtStartPar
The specific SOS constraint. Optional, \sphinxcode{\sphinxupquote{None}} by default.

\end{description}
\end{quote}

\sphinxAtStartPar
\sphinxstylestrong{Example}
\end{quote}

\begin{sphinxVerbatim}[commandchars=\\\{\}]
\PYG{c+c1}{\PYGZsh{} Retrieve builders corresponding to all SOS constraints in the model.}
\PYG{n}{soss} \PYG{o}{=} \PYG{n}{m}\PYG{o}{.}\PYG{n}{getSOSBuilders}\PYG{p}{(}\PYG{p}{)}
\end{sphinxVerbatim}

\subsubsection{Model.getGenConstrIndicator()}
\label{\detokenize{pyapiref:model-getgenconstrindicator}}\begin{quote}

\sphinxAtStartPar
\sphinxstylestrong{Synopsis}
\begin{quote}

\sphinxAtStartPar
\sphinxcode{\sphinxupquote{getGenConstrIndicator(genconstr)}}
\end{quote}

\sphinxAtStartPar
\sphinxstylestrong{Description}
\begin{quote}

\sphinxAtStartPar
Retrieve the builder corresponding to specific indicator constraint. Return a {\hyperref[\detokenize{pyapiref:chappyapi-genconstrbuilder}]{\sphinxcrossref{\DUrole{std,std-ref}{GenConstrBuilder Class}}}} object.
\end{quote}

\sphinxAtStartPar
\sphinxstylestrong{Arguments}
\begin{quote}

\sphinxAtStartPar
\sphinxcode{\sphinxupquote{genconstr}}
\begin{quote}

\sphinxAtStartPar
The specified indicator constraint.
\end{quote}
\end{quote}

\sphinxAtStartPar
\sphinxstylestrong{Example}
\end{quote}

\begin{sphinxVerbatim}[commandchars=\\\{\}]
\PYG{c+c1}{\PYGZsh{} Retrieve the builder corresponding to indicator constraint genx.}
\PYG{n}{indic} \PYG{o}{=} \PYG{n}{m}\PYG{o}{.}\PYG{n}{getGenConstrIndicator}\PYG{p}{(}\PYG{n}{genx}\PYG{p}{)}
\end{sphinxVerbatim}

\subsubsection{Model.getGenConstr()}
\label{\detokenize{pyapiref:model-getgenconstr}}\begin{quote}

\sphinxAtStartPar
\sphinxstylestrong{Synopsis}
\begin{quote}

\sphinxAtStartPar
\sphinxcode{\sphinxupquote{getGenConstr(idx)}}
\end{quote}

\sphinxAtStartPar
\sphinxstylestrong{Description}
\begin{quote}

\sphinxAtStartPar
Retrieve an indicator constraint by its indice in the model. Return a {\hyperref[\detokenize{pyapiref:chappyapi-genconstr}]{\sphinxcrossref{\DUrole{std,std-ref}{GenConstr Class}}}} object.
\end{quote}

\sphinxAtStartPar
\sphinxstylestrong{Arguments}
\begin{quote}

\sphinxAtStartPar
\sphinxcode{\sphinxupquote{idx}}
\begin{quote}

\sphinxAtStartPar
Index of the desired constraint in the model, starting with \sphinxcode{\sphinxupquote{0}}.
\end{quote}
\end{quote}

\sphinxAtStartPar
\sphinxstylestrong{Example}
\end{quote}

\begin{sphinxVerbatim}[commandchars=\\\{\}]
\PYG{c+c1}{\PYGZsh{} Retrive indicator constraint with indice of 0}
\PYG{n}{genx} \PYG{o}{=} \PYG{n}{m}\PYG{o}{.}\PYG{n}{getGenConstr}\PYG{p}{(}\PYG{l+m+mi}{0}\PYG{p}{)}
\end{sphinxVerbatim}

\subsubsection{Model.getGenConstrs()}
\label{\detokenize{pyapiref:model-getgenconstrs}}\begin{quote}

\sphinxAtStartPar
\sphinxstylestrong{Synopsis}
\begin{quote}

\sphinxAtStartPar
\sphinxcode{\sphinxupquote{getGenConstrs()}}
\end{quote}

\sphinxAtStartPar
\sphinxstylestrong{Description}
\begin{quote}

\sphinxAtStartPar
Retrieve all indicator constraints in the model. Return a {\hyperref[\detokenize{pyapiref:chappyapi-genconstrarray}]{\sphinxcrossref{\DUrole{std,std-ref}{GenConstrArray Class}}}} object.
\end{quote}

\sphinxAtStartPar
\sphinxstylestrong{Example}
\end{quote}

\begin{sphinxVerbatim}[commandchars=\\\{\}]
\PYG{c+c1}{\PYGZsh{} Retrieve all indicator constraints in the model}
\PYG{n}{cons} \PYG{o}{=} \PYG{n}{m}\PYG{o}{.}\PYG{n}{getGenConstrs}\PYG{p}{(}\PYG{p}{)}
\end{sphinxVerbatim}

\subsubsection{Model.getGenConstrByName()}
\label{\detokenize{pyapiref:model-getgenconstrbyname}}\begin{quote}

\sphinxAtStartPar
\sphinxstylestrong{Synopsis}
\begin{quote}

\sphinxAtStartPar
\sphinxcode{\sphinxupquote{getGenConstrByName(name)}}
\end{quote}

\sphinxAtStartPar
\sphinxstylestrong{Description}
\begin{quote}

\sphinxAtStartPar
Retrieves an indicator constraint by the specified name. Return a {\hyperref[\detokenize{pyapiref:chappyapi-genconstr}]{\sphinxcrossref{\DUrole{std,std-ref}{GenConstr Class}}}} object.
\end{quote}

\sphinxAtStartPar
\sphinxstylestrong{Arguments}
\begin{quote}

\sphinxAtStartPar
\sphinxcode{\sphinxupquote{name}}
\begin{quote}

\sphinxAtStartPar
The name of the indicator constraint.
\end{quote}
\end{quote}

\sphinxAtStartPar
\sphinxstylestrong{Example}
\end{quote}

\begin{sphinxVerbatim}[commandchars=\\\{\}]
\PYG{c+c1}{\PYGZsh{} Retrieve indicator constraint with name \PYGZdq{}r\PYGZdq{}}
\PYG{n}{r} \PYG{o}{=} \PYG{n}{m}\PYG{o}{.}\PYG{n}{getGenConstrByName}\PYG{p}{(}\PYG{l+s+s2}{\PYGZdq{}}\PYG{l+s+s2}{r}\PYG{l+s+s2}{\PYGZdq{}}\PYG{p}{)}
\end{sphinxVerbatim}

\subsubsection{Model.getGenConstrIndicators()}
\label{\detokenize{pyapiref:model-getgenconstrindicators}}\begin{quote}

\sphinxAtStartPar
\sphinxstylestrong{Synopsis}
\begin{quote}

\sphinxAtStartPar
\sphinxcode{\sphinxupquote{getGenConstrIndicators(genconstrs=None)}}
\end{quote}

\sphinxAtStartPar
\sphinxstylestrong{Description}
\begin{quote}

\sphinxAtStartPar
Retrieve the specified indicator constraints’ builders in the model.
All indicator constraints will be retrieved in default.
Return a {\hyperref[\detokenize{pyapiref:chappyapi-genconstrbuilder}]{\sphinxcrossref{\DUrole{std,std-ref}{GenConstrBuilder Class}}}} or {\hyperref[\detokenize{pyapiref:chappyapi-genconstrbuilderarray}]{\sphinxcrossref{\DUrole{std,std-ref}{GenConstrBuilderArray Class}}}} object.
\end{quote}

\sphinxAtStartPar
\sphinxstylestrong{Example}
\end{quote}

\begin{sphinxVerbatim}[commandchars=\\\{\}]
\PYG{c+c1}{\PYGZsh{} Retrieve all indicator constraints\PYGZsq{} builders in the model}
\PYG{n}{cons} \PYG{o}{=} \PYG{n}{m}\PYG{o}{.}\PYG{n}{getGenConstrIndicators}\PYG{p}{(}\PYG{p}{)}
\end{sphinxVerbatim}

\subsubsection{Model.getCone()}
\label{\detokenize{pyapiref:model-getcone}}\begin{quote}

\sphinxAtStartPar
\sphinxstylestrong{Synopsis}
\begin{quote}

\sphinxAtStartPar
\sphinxcode{\sphinxupquote{getCone(idx)}}
\end{quote}

\sphinxAtStartPar
\sphinxstylestrong{Description}
\begin{quote}

\sphinxAtStartPar
Retrieves the second\sphinxhyphen{}order cone at the specified index in the model.
Returns a {\hyperref[\detokenize{pyapiref:chappyapi-cone}]{\sphinxcrossref{\DUrole{std,std-ref}{Cone Class}}}} object.
\end{quote}

\sphinxAtStartPar
\sphinxstylestrong{Arguments}
\begin{quote}

\sphinxAtStartPar
\sphinxcode{\sphinxupquote{idx}}
\begin{quote}

\sphinxAtStartPar
The specified index. Indexing starts at 0.
\end{quote}
\end{quote}

\sphinxAtStartPar
\sphinxstylestrong{Example}
\end{quote}

\begin{sphinxVerbatim}[commandchars=\\\{\}]
\PYG{c+c1}{\PYGZsh{} Retrieve the second\PYGZhy{}order cone at index 1 in the model}
\PYG{n}{cones} \PYG{o}{=} \PYG{n}{m}\PYG{o}{.}\PYG{n}{getCone}\PYG{p}{(}\PYG{l+m+mi}{1}\PYG{p}{)}
\end{sphinxVerbatim}

\subsubsection{Model.getExpCone()}
\label{\detokenize{pyapiref:model-getexpcone}}\begin{quote}

\sphinxAtStartPar
\sphinxstylestrong{Synopsis}
\begin{quote}

\sphinxAtStartPar
\sphinxcode{\sphinxupquote{getExpCone(idx)}}
\end{quote}

\sphinxAtStartPar
\sphinxstylestrong{Description}
\begin{quote}

\sphinxAtStartPar
Retrieves the exponential cone at the specified index in the model.
Returns an {\hyperref[\detokenize{pyapiref:chappyapi-expcone}]{\sphinxcrossref{\DUrole{std,std-ref}{ExpCone Class}}}} object.
\end{quote}

\sphinxAtStartPar
\sphinxstylestrong{Arguments}
\begin{quote}

\sphinxAtStartPar
\sphinxcode{\sphinxupquote{idx}}
\begin{quote}

\sphinxAtStartPar
The specified index. Indexing starts at 0.
\end{quote}
\end{quote}

\sphinxAtStartPar
\sphinxstylestrong{Example}
\end{quote}

\begin{sphinxVerbatim}[commandchars=\\\{\}]
\PYG{c+c1}{\PYGZsh{} Retrieve the exponential cone at index 1 in the model}
\PYG{n}{cones} \PYG{o}{=} \PYG{n}{m}\PYG{o}{.}\PYG{n}{getExpCone}\PYG{p}{(}\PYG{l+m+mi}{1}\PYG{p}{)}
\end{sphinxVerbatim}

\subsubsection{Model.getAffineCone()}
\label{\detokenize{pyapiref:model-getaffinecone}}\begin{quote}

\sphinxAtStartPar
\sphinxstylestrong{Synopsis}
\begin{quote}

\sphinxAtStartPar
\sphinxcode{\sphinxupquote{getAffineCone(idx)}}
\end{quote}

\sphinxAtStartPar
\sphinxstylestrong{Description}
\begin{quote}

\sphinxAtStartPar
Retrieves the affine cone at the specified index in the model.
Returns an {\hyperref[\detokenize{pyapiref:chappyapi-affinecone}]{\sphinxcrossref{\DUrole{std,std-ref}{AffineCone Class}}}} object.
\end{quote}

\sphinxAtStartPar
\sphinxstylestrong{Arguments}
\begin{quote}

\sphinxAtStartPar
\sphinxcode{\sphinxupquote{idx}}
\begin{quote}

\sphinxAtStartPar
The specified index. Indexing starts at 0.
\end{quote}
\end{quote}

\sphinxAtStartPar
\sphinxstylestrong{Example}
\end{quote}

\begin{sphinxVerbatim}[commandchars=\\\{\}]
\PYG{c+c1}{\PYGZsh{} Retrieve the affine cone at index 1 in the model}
\PYG{n}{cones} \PYG{o}{=} \PYG{n}{m}\PYG{o}{.}\PYG{n}{getAffineCone}\PYG{p}{(}\PYG{l+m+mi}{1}\PYG{p}{)}
\end{sphinxVerbatim}

\subsubsection{Model.getAffineConeByName()}
\label{\detokenize{pyapiref:model-getaffineconebyname}}\begin{quote}

\sphinxAtStartPar
\sphinxstylestrong{Synopsis}
\begin{quote}

\sphinxAtStartPar
\sphinxcode{\sphinxupquote{getAffineConeByName(name)}}
\end{quote}

\sphinxAtStartPar
\sphinxstylestrong{Description}
\begin{quote}

\sphinxAtStartPar
Retrieves the affine cone with the specified name in the model.
Returns a {\hyperref[\detokenize{pyapiref:chappyapi-affinecone}]{\sphinxcrossref{\DUrole{std,std-ref}{AffineCone Class}}}} object.
\end{quote}

\sphinxAtStartPar
\sphinxstylestrong{Arguments}
\begin{quote}

\sphinxAtStartPar
\sphinxcode{\sphinxupquote{name}}
\begin{quote}

\sphinxAtStartPar
The specified name.
\end{quote}
\end{quote}

\sphinxAtStartPar
\sphinxstylestrong{Example}
\end{quote}

\begin{sphinxVerbatim}[commandchars=\\\{\}]
\PYG{c+c1}{\PYGZsh{} Retrieve the affine cone with the name \PYGZdq{}afcone\PYGZdq{} in the model}
\PYG{n}{cones} \PYG{o}{=} \PYG{n}{m}\PYG{o}{.}\PYG{n}{getAffineConeByName}\PYG{p}{(}\PYG{l+s+s2}{\PYGZdq{}}\PYG{l+s+s2}{afcone}\PYG{l+s+s2}{\PYGZdq{}}\PYG{p}{)}
\end{sphinxVerbatim}

\subsubsection{Model.getCones()}
\label{\detokenize{pyapiref:model-getcones}}\begin{quote}

\sphinxAtStartPar
\sphinxstylestrong{Synopsis}
\begin{quote}

\sphinxAtStartPar
\sphinxcode{\sphinxupquote{getCones()}}
\end{quote}

\sphinxAtStartPar
\sphinxstylestrong{Description}
\begin{quote}

\sphinxAtStartPar
Retrieve all Second\sphinxhyphen{}Order\sphinxhyphen{}Cone (SOC) constraints in model, and return a
{\hyperref[\detokenize{pyapiref:chappyapi-conearray}]{\sphinxcrossref{\DUrole{std,std-ref}{ConeArray Class}}}} object.
\end{quote}

\sphinxAtStartPar
\sphinxstylestrong{Example}
\end{quote}

\begin{sphinxVerbatim}[commandchars=\\\{\}]
\PYG{c+c1}{\PYGZsh{} Retrieve all SOC constraints}
\PYG{n}{cones} \PYG{o}{=} \PYG{n}{m}\PYG{o}{.}\PYG{n}{getCones}\PYG{p}{(}\PYG{p}{)}
\end{sphinxVerbatim}

\subsubsection{Model.getExpCones()}
\label{\detokenize{pyapiref:model-getexpcones}}\begin{quote}

\sphinxAtStartPar
\sphinxstylestrong{Synopsis}
\begin{quote}

\sphinxAtStartPar
\sphinxcode{\sphinxupquote{getExpCones()}}
\end{quote}

\sphinxAtStartPar
\sphinxstylestrong{Description}
\begin{quote}

\sphinxAtStartPar
Retrieve all exponential cone constraints in model, and return a
{\hyperref[\detokenize{pyapiref:chappyapi-expconearray}]{\sphinxcrossref{\DUrole{std,std-ref}{ExpConeArray Class}}}} object.
\end{quote}

\sphinxAtStartPar
\sphinxstylestrong{Example}
\end{quote}

\begin{sphinxVerbatim}[commandchars=\\\{\}]
\PYG{c+c1}{\PYGZsh{} Retrieve all exponential cone constraints}
\PYG{n}{cones} \PYG{o}{=} \PYG{n}{m}\PYG{o}{.}\PYG{n}{getExpCones}\PYG{p}{(}\PYG{p}{)}
\end{sphinxVerbatim}

\subsubsection{Model.getAffineCones()}
\label{\detokenize{pyapiref:model-getaffinecones}}\begin{quote}

\sphinxAtStartPar
\sphinxstylestrong{Synopsis}
\begin{quote}

\sphinxAtStartPar
\sphinxcode{\sphinxupquote{getAffineCones()}}
\end{quote}

\sphinxAtStartPar
\sphinxstylestrong{Description}
\begin{quote}

\sphinxAtStartPar
Retrieves all affine cones in the model.
Returns a {\hyperref[\detokenize{pyapiref:chappyapi-affineconearray}]{\sphinxcrossref{\DUrole{std,std-ref}{AffineConeArray Class}}}} object.
\end{quote}

\sphinxAtStartPar
\sphinxstylestrong{Example}
\end{quote}

\begin{sphinxVerbatim}[commandchars=\\\{\}]
\PYG{c+c1}{\PYGZsh{} Retrieve all affine cones in the model}
\PYG{n}{cones} \PYG{o}{=} \PYG{n}{m}\PYG{o}{.}\PYG{n}{getAffineCones}\PYG{p}{(}\PYG{p}{)}
\end{sphinxVerbatim}

\subsubsection{Model.getConeBuilders()}
\label{\detokenize{pyapiref:model-getconebuilders}}\begin{quote}

\sphinxAtStartPar
\sphinxstylestrong{Synopsis}
\begin{quote}

\sphinxAtStartPar
\sphinxcode{\sphinxupquote{getConeBuilders(cones=None)}}
\end{quote}

\sphinxAtStartPar
\sphinxstylestrong{Description}
\begin{quote}

\sphinxAtStartPar
Retrieve Second\sphinxhyphen{}Order\sphinxhyphen{}Cone (SOC) constraint builders for given SOC constraints.

\sphinxAtStartPar
If argument \sphinxcode{\sphinxupquote{cones}} is \sphinxcode{\sphinxupquote{None}}, then return a {\hyperref[\detokenize{pyapiref:chappyapi-conebuilderarray}]{\sphinxcrossref{\DUrole{std,std-ref}{ConeBuilderArray Class}}}} object
consists of all SOC constraints’ builders;
If argument \sphinxcode{\sphinxupquote{cones}} is {\hyperref[\detokenize{pyapiref:chappyapi-cone}]{\sphinxcrossref{\DUrole{std,std-ref}{Cone Class}}}} object, then return a {\hyperref[\detokenize{pyapiref:chappyapi-conebuilder}]{\sphinxcrossref{\DUrole{std,std-ref}{ConeBuilder Class}}}} object
of given SOC constraint;
If \sphinxcode{\sphinxupquote{cones}} is Python list or {\hyperref[\detokenize{pyapiref:chappyapi-conearray}]{\sphinxcrossref{\DUrole{std,std-ref}{ConeArray Class}}}} object,
then return a {\hyperref[\detokenize{pyapiref:chappyapi-conebuilderarray}]{\sphinxcrossref{\DUrole{std,std-ref}{ConeBuilderArray Class}}}} object consists of builders of given SOC constraints.
\end{quote}

\sphinxAtStartPar
\sphinxstylestrong{Arguments}
\begin{quote}

\sphinxAtStartPar
\sphinxcode{\sphinxupquote{cones}}
\begin{quote}

\sphinxAtStartPar
Given SOC constraints. Optional, default to \sphinxcode{\sphinxupquote{None}}.
\end{quote}
\end{quote}

\sphinxAtStartPar
\sphinxstylestrong{Example}
\end{quote}

\begin{sphinxVerbatim}[commandchars=\\\{\}]
\PYG{c+c1}{\PYGZsh{} Retrieve all SOC constraints\PYGZsq{} builders}
\PYG{n}{cones} \PYG{o}{=} \PYG{n}{m}\PYG{o}{.}\PYG{n}{getConeBuilders}\PYG{p}{(}\PYG{p}{)}
\end{sphinxVerbatim}

\subsubsection{Model.getExpConeBuilders()}
\label{\detokenize{pyapiref:model-getexpconebuilders}}\begin{quote}

\sphinxAtStartPar
\sphinxstylestrong{Synopsis}
\begin{quote}

\sphinxAtStartPar
\sphinxcode{\sphinxupquote{getExpConeBuilders(cones=None)}}
\end{quote}

\sphinxAtStartPar
\sphinxstylestrong{Description}
\begin{quote}

\sphinxAtStartPar
Retrieve the exponential cone constraint builders corresponding to the specified exponential cone constraint.

\sphinxAtStartPar
If argument \sphinxcode{\sphinxupquote{cones}} is \sphinxcode{\sphinxupquote{None}}, then return a {\hyperref[\detokenize{pyapiref:chappyapi-expconebuilderarray}]{\sphinxcrossref{\DUrole{std,std-ref}{ExpConeBuilderArray Class}}}} object
consists of all exponential cone constraints’ builders;
If argument \sphinxcode{\sphinxupquote{cones}} is {\hyperref[\detokenize{pyapiref:chappyapi-expcone}]{\sphinxcrossref{\DUrole{std,std-ref}{ExpCone Class}}}} object, then return a {\hyperref[\detokenize{pyapiref:chappyapi-expconebuilder}]{\sphinxcrossref{\DUrole{std,std-ref}{ExpConeBuilder Class}}}} object
of given exponential cone constraints;
If \sphinxcode{\sphinxupquote{cones}} is Python list or {\hyperref[\detokenize{pyapiref:chappyapi-expconearray}]{\sphinxcrossref{\DUrole{std,std-ref}{ExpConeArray Class}}}} object,
then return a {\hyperref[\detokenize{pyapiref:chappyapi-expconebuilderarray}]{\sphinxcrossref{\DUrole{std,std-ref}{ExpConeBuilderArray Class}}}} object consists of builders of given exponential cone constraints.
\end{quote}

\sphinxAtStartPar
\sphinxstylestrong{Arguments}
\begin{quote}

\sphinxAtStartPar
\sphinxcode{\sphinxupquote{cones}}
\begin{quote}

\sphinxAtStartPar
Given exponential cone constraints. Optional, default to \sphinxcode{\sphinxupquote{None}}.
\end{quote}
\end{quote}

\sphinxAtStartPar
\sphinxstylestrong{Example}
\end{quote}

\begin{sphinxVerbatim}[commandchars=\\\{\}]
\PYG{c+c1}{\PYGZsh{} Retrieve all exponential cone constraints\PYGZsq{} builders}
\PYG{n}{cones} \PYG{o}{=} \PYG{n}{m}\PYG{o}{.}\PYG{n}{getExpConeBuilders}\PYG{p}{(}\PYG{p}{)}
\end{sphinxVerbatim}

\subsubsection{Model.getAffineConeBuilders()}
\label{\detokenize{pyapiref:model-getaffineconebuilders}}\begin{quote}

\sphinxAtStartPar
\sphinxstylestrong{Synopsis}
\begin{quote}

\sphinxAtStartPar
\sphinxcode{\sphinxupquote{getAffineConeBuilders(cones=None)}}
\end{quote}

\sphinxAtStartPar
\sphinxstylestrong{Description}
\begin{quote}

\sphinxAtStartPar
Retrieves the affine cone builders corresponding to the specified affine cones.

\sphinxAtStartPar
If the argument \sphinxcode{\sphinxupquote{cones}} is \sphinxcode{\sphinxupquote{None}}, returns a {\hyperref[\detokenize{pyapiref:chappyapi-affineconebuilderarray}]{\sphinxcrossref{\DUrole{std,std-ref}{AffineConeBuilderArray Class}}}} object containing the builders for all affine cones in the model.

\sphinxAtStartPar
If the argument \sphinxcode{\sphinxupquote{cones}} is a {\hyperref[\detokenize{pyapiref:chappyapi-affinecone}]{\sphinxcrossref{\DUrole{std,std-ref}{AffineCone Class}}}} object, returns the corresponding {\hyperref[\detokenize{pyapiref:chappyapi-affineconebuilder}]{\sphinxcrossref{\DUrole{std,std-ref}{AffineConeBuilder Class}}}} object.

\sphinxAtStartPar
If the argument \sphinxcode{\sphinxupquote{cones}} is a list or a {\hyperref[\detokenize{pyapiref:chappyapi-affineconearray}]{\sphinxcrossref{\DUrole{std,std-ref}{AffineConeArray Class}}}} object, returns a {\hyperref[\detokenize{pyapiref:chappyapi-affineconebuilderarray}]{\sphinxcrossref{\DUrole{std,std-ref}{AffineConeBuilderArray Class}}}} object containing the builders for the specified affine cones.
\end{quote}

\sphinxAtStartPar
\sphinxstylestrong{Arguments}
\begin{quote}

\sphinxAtStartPar
\sphinxcode{\sphinxupquote{cones}}
\begin{quote}

\sphinxAtStartPar
The specified affine cones.
Optional, defaulting to \sphinxcode{\sphinxupquote{None}}.
\end{quote}
\end{quote}

\sphinxAtStartPar
\sphinxstylestrong{Example}
\end{quote}

\begin{sphinxVerbatim}[commandchars=\\\{\}]
\PYG{c+c1}{\PYGZsh{} Retrieve the builders for all affine cones in the model}
\PYG{n}{cones} \PYG{o}{=} \PYG{n}{m}\PYG{o}{.}\PYG{n}{getAffineConeBuilders}\PYG{p}{(}\PYG{p}{)}
\end{sphinxVerbatim}

\subsubsection{Model.getQConstr()}
\label{\detokenize{pyapiref:model-getqconstr}}\begin{quote}

\sphinxAtStartPar
\sphinxstylestrong{Synopsis}
\begin{quote}

\sphinxAtStartPar
\sphinxcode{\sphinxupquote{getQConstr(idx)}}
\end{quote}

\sphinxAtStartPar
\sphinxstylestrong{Description}
\begin{quote}

\sphinxAtStartPar
Retrieve a quadratic constraint by its indice, and return a {\hyperref[\detokenize{pyapiref:chappyapi-qconstraint}]{\sphinxcrossref{\DUrole{std,std-ref}{QConstraint Class}}}} object.
\end{quote}

\sphinxAtStartPar
\sphinxstylestrong{Arguments}
\begin{quote}

\sphinxAtStartPar
\sphinxcode{\sphinxupquote{idx}}
\begin{quote}

\sphinxAtStartPar
Index of the desired quadratic constraint, starting with \sphinxcode{\sphinxupquote{0}}.
\end{quote}
\end{quote}

\sphinxAtStartPar
\sphinxstylestrong{Example}
\end{quote}

\begin{sphinxVerbatim}[commandchars=\\\{\}]
\PYG{c+c1}{\PYGZsh{} Retrieve a quadratic constraint with indice of 1}
\PYG{n}{qr} \PYG{o}{=} \PYG{n}{m}\PYG{o}{.}\PYG{n}{getQConstr}\PYG{p}{(}\PYG{l+m+mi}{1}\PYG{p}{)}
\end{sphinxVerbatim}

\subsubsection{Model.getQConstrByName()}
\label{\detokenize{pyapiref:model-getqconstrbyname}}\begin{quote}

\sphinxAtStartPar
\sphinxstylestrong{Synopsis}
\begin{quote}

\sphinxAtStartPar
\sphinxcode{\sphinxupquote{getQConstrByName(name)}}
\end{quote}

\sphinxAtStartPar
\sphinxstylestrong{Description}
\begin{quote}

\sphinxAtStartPar
Retrieve a quadratic constraint by its name, and return a {\hyperref[\detokenize{pyapiref:chappyapi-qconstraint}]{\sphinxcrossref{\DUrole{std,std-ref}{QConstraint Class}}}} object.
\end{quote}

\sphinxAtStartPar
\sphinxstylestrong{Arguments}
\begin{quote}

\sphinxAtStartPar
\sphinxcode{\sphinxupquote{name}}
\begin{quote}

\sphinxAtStartPar
Name of the desired quadratic constraint.
\end{quote}
\end{quote}

\sphinxAtStartPar
\sphinxstylestrong{Example}
\end{quote}

\begin{sphinxVerbatim}[commandchars=\\\{\}]
\PYG{c+c1}{\PYGZsh{} Retrieve a quadratic constraint with name \PYGZdq{}qr\PYGZdq{}}
\PYG{n}{qr} \PYG{o}{=} \PYG{n}{m}\PYG{o}{.}\PYG{n}{getQConstrByName}\PYG{p}{(}\PYG{l+s+s2}{\PYGZdq{}}\PYG{l+s+s2}{qr}\PYG{l+s+s2}{\PYGZdq{}}\PYG{p}{)}
\end{sphinxVerbatim}

\subsubsection{Model.getQConstrs()}
\label{\detokenize{pyapiref:model-getqconstrs}}\begin{quote}

\sphinxAtStartPar
\sphinxstylestrong{Synopsis}
\begin{quote}

\sphinxAtStartPar
\sphinxcode{\sphinxupquote{getQConstrs()}}
\end{quote}

\sphinxAtStartPar
\sphinxstylestrong{Description}
\begin{quote}

\sphinxAtStartPar
Retrieve all quadratic constraints in the model. Return a {\hyperref[\detokenize{pyapiref:chappyapi-qconstrarray}]{\sphinxcrossref{\DUrole{std,std-ref}{QConstrArray Class}}}} object.
\end{quote}

\sphinxAtStartPar
\sphinxstylestrong{Example}
\end{quote}

\begin{sphinxVerbatim}[commandchars=\\\{\}]
\PYG{c+c1}{\PYGZsh{} Retrieve all quadratic constraints in the model}
\PYG{n}{qcons} \PYG{o}{=} \PYG{n}{m}\PYG{o}{.}\PYG{n}{getQConstrs}\PYG{p}{(}\PYG{p}{)}
\end{sphinxVerbatim}

\subsubsection{Model.getQConstrBuilders()}
\label{\detokenize{pyapiref:model-getqconstrbuilders}}\begin{quote}

\sphinxAtStartPar
\sphinxstylestrong{Synopsis}
\begin{quote}

\sphinxAtStartPar
\sphinxcode{\sphinxupquote{getQConstrBuilders(qconstrs=None)}}
\end{quote}

\sphinxAtStartPar
\sphinxstylestrong{Description}
\begin{quote}

\sphinxAtStartPar
Retrieve quadratic constraint builders in current model.

\sphinxAtStartPar
If parameter \sphinxcode{\sphinxupquote{qconstrs}} is \sphinxcode{\sphinxupquote{None}}, then return a {\hyperref[\detokenize{pyapiref:chappyapi-qconstrbuilderarray}]{\sphinxcrossref{\DUrole{std,std-ref}{QConstrBuilderArray Class}}}} object composed of all
quadratic constraint builders.

\sphinxAtStartPar
If parameter \sphinxcode{\sphinxupquote{qconstrs}} is {\hyperref[\detokenize{pyapiref:chappyapi-qconstraint}]{\sphinxcrossref{\DUrole{std,std-ref}{QConstraint Class}}}} object, then return the {\hyperref[\detokenize{pyapiref:chappyapi-qconstrbuilder}]{\sphinxcrossref{\DUrole{std,std-ref}{QConstrBuilder Class}}}}
object corresponding to the specific quadratic constraint.

\sphinxAtStartPar
If parameter \sphinxcode{\sphinxupquote{qconstrs}} is a list or a {\hyperref[\detokenize{pyapiref:chappyapi-qconstrarray}]{\sphinxcrossref{\DUrole{std,std-ref}{QConstrArray Class}}}} object, then return a {\hyperref[\detokenize{pyapiref:chappyapi-qconstrbuilderarray}]{\sphinxcrossref{\DUrole{std,std-ref}{QConstrBuilderArray Class}}}}
object composed of specified quadratic constraints’ builders.

\sphinxAtStartPar
If parameter \sphinxcode{\sphinxupquote{qconstrs}} is dictionary or {\hyperref[\detokenize{pyapiref:chappyapi-util-tupledict}]{\sphinxcrossref{\DUrole{std,std-ref}{tupledict Class}}}} object, then the indice of the specified quadratic
constraint is returned as key, the value is a {\hyperref[\detokenize{pyapiref:chappyapi-util-tupledict}]{\sphinxcrossref{\DUrole{std,std-ref}{tupledict Class}}}} object composed of the specified quadratic constraints’ builders.
\end{quote}

\sphinxAtStartPar
\sphinxstylestrong{Arguments}
\begin{quote}

\sphinxAtStartPar
\sphinxcode{\sphinxupquote{qconstrs}}
\begin{quote}

\sphinxAtStartPar
The specified quadratic constraint. Optional, \sphinxcode{\sphinxupquote{None}} by default.
\end{quote}
\end{quote}

\sphinxAtStartPar
\sphinxstylestrong{Example}
\end{quote}

\begin{sphinxVerbatim}[commandchars=\\\{\}]
\PYG{c+c1}{\PYGZsh{} Retrieve all of quadratic constraint builders.}
\PYG{n}{qconbuilders} \PYG{o}{=} \PYG{n}{m}\PYG{o}{.}\PYG{n}{getQConstrBuilders}\PYG{p}{(}\PYG{p}{)}
\PYG{c+c1}{\PYGZsh{} Retrive the builder corresponding to quadratic contstraint qx.}
\PYG{n}{qconbuilders} \PYG{o}{=} \PYG{n}{m}\PYG{o}{.}\PYG{n}{getQConstrBuilders}\PYG{p}{(}\PYG{n}{qx}\PYG{p}{)}
\PYG{c+c1}{\PYGZsh{} Retrieve builders corresponding to quadratic constraint qx and qy.}
\PYG{n}{qconbuilders} \PYG{o}{=} \PYG{n}{m}\PYG{o}{.}\PYG{n}{getQConstrBuilders}\PYG{p}{(}\PYG{p}{[}\PYG{n}{qx}\PYG{p}{,} \PYG{n}{qy}\PYG{p}{]}\PYG{p}{)}
\PYG{c+c1}{\PYGZsh{} Retrieve builders corresponding to quadratic constraint in tupledict object qxx.}
\PYG{n}{qconbuilders} \PYG{o}{=} \PYG{n}{m}\PYG{o}{.}\PYG{n}{getQConstrBuilders}\PYG{p}{(}\PYG{n}{qxx}\PYG{p}{)}
\end{sphinxVerbatim}

\subsubsection{Model.getPsdVar()}
\label{\detokenize{pyapiref:model-getpsdvar}}\begin{quote}

\sphinxAtStartPar
\sphinxstylestrong{Synopsis}
\begin{quote}

\sphinxAtStartPar
\sphinxcode{\sphinxupquote{getPsdVar(idx)}}
\end{quote}

\sphinxAtStartPar
\sphinxstylestrong{Description}
\begin{quote}

\sphinxAtStartPar
Retrieve a positive semi\sphinxhyphen{}definite variable according to its index in the model. Return a {\hyperref[\detokenize{pyapiref:chappyapi-psdvar}]{\sphinxcrossref{\DUrole{std,std-ref}{PsdVar Class}}}} object.
\end{quote}

\sphinxAtStartPar
\sphinxstylestrong{Arguments}
\begin{quote}

\sphinxAtStartPar
\sphinxcode{\sphinxupquote{idx}}
\begin{quote}

\sphinxAtStartPar
Index of the desired positive semi\sphinxhyphen{}definite variable in the model, starting with \sphinxcode{\sphinxupquote{0}}.
\end{quote}
\end{quote}

\sphinxAtStartPar
\sphinxstylestrong{Example}
\end{quote}

\begin{sphinxVerbatim}[commandchars=\\\{\}]
\PYG{c+c1}{\PYGZsh{} Retrieve a positive semi\PYGZhy{}definite variable with index of 1}
\PYG{n}{x} \PYG{o}{=} \PYG{n}{m}\PYG{o}{.}\PYG{n}{getPsdVar}\PYG{p}{(}\PYG{l+m+mi}{1}\PYG{p}{)}
\end{sphinxVerbatim}

\subsubsection{Model.getPsdVarByName()}
\label{\detokenize{pyapiref:model-getpsdvarbyname}}\begin{quote}

\sphinxAtStartPar
\sphinxstylestrong{Synopsis}
\begin{quote}

\sphinxAtStartPar
\sphinxcode{\sphinxupquote{getPsdVarByName(name)}}
\end{quote}

\sphinxAtStartPar
\sphinxstylestrong{Description}
\begin{quote}

\sphinxAtStartPar
Retrieve a positive semi\sphinxhyphen{}definite variable by name. Return a {\hyperref[\detokenize{pyapiref:chappyapi-psdvar}]{\sphinxcrossref{\DUrole{std,std-ref}{PsdVar Class}}}} object.
\end{quote}

\sphinxAtStartPar
\sphinxstylestrong{Arguments}
\begin{quote}

\sphinxAtStartPar
\sphinxcode{\sphinxupquote{name}}
\begin{quote}

\sphinxAtStartPar
The name of the positive semi\sphinxhyphen{}definite variable.
\end{quote}
\end{quote}

\sphinxAtStartPar
\sphinxstylestrong{Example}
\end{quote}

\begin{sphinxVerbatim}[commandchars=\\\{\}]
\PYG{c+c1}{\PYGZsh{} Retrieve a positive semi\PYGZhy{}definite variable with name \PYGZdq{}x\PYGZdq{}.}
\PYG{n}{x} \PYG{o}{=} \PYG{n}{m}\PYG{o}{.}\PYG{n}{getPsdVarByName}\PYG{p}{(}\PYG{l+s+s2}{\PYGZdq{}}\PYG{l+s+s2}{x}\PYG{l+s+s2}{\PYGZdq{}}\PYG{p}{)}
\end{sphinxVerbatim}

\subsubsection{Model.getPsdVars()}
\label{\detokenize{pyapiref:model-getpsdvars}}\begin{quote}

\sphinxAtStartPar
\sphinxstylestrong{Synopsis}
\begin{quote}

\sphinxAtStartPar
\sphinxcode{\sphinxupquote{getPsdVars()}}
\end{quote}

\sphinxAtStartPar
\sphinxstylestrong{Description}
\begin{quote}

\sphinxAtStartPar
Retrieve all positive semi\sphinxhyphen{}definite variables in the model, and return a {\hyperref[\detokenize{pyapiref:chappyapi-psdvararray}]{\sphinxcrossref{\DUrole{std,std-ref}{PsdVarArray Class}}}} object.
\end{quote}

\sphinxAtStartPar
\sphinxstylestrong{Example}
\end{quote}

\begin{sphinxVerbatim}[commandchars=\\\{\}]
\PYG{c+c1}{\PYGZsh{} Retrieve all positive semi\PYGZhy{}definite variables in the model.}
\PYG{n+nb}{vars} \PYG{o}{=} \PYG{n}{m}\PYG{o}{.}\PYG{n}{getPsdVars}\PYG{p}{(}\PYG{p}{)}
\end{sphinxVerbatim}

\subsubsection{Model.getPsdConstr()}
\label{\detokenize{pyapiref:model-getpsdconstr}}\begin{quote}

\sphinxAtStartPar
\sphinxstylestrong{Synopsis}
\begin{quote}

\sphinxAtStartPar
\sphinxcode{\sphinxupquote{getPsdConstr(idx)}}
\end{quote}

\sphinxAtStartPar
\sphinxstylestrong{Description}
\begin{quote}

\sphinxAtStartPar
Retrieve the positive semi\sphinxhyphen{}definite constraint according to its index in the model. Return a {\hyperref[\detokenize{pyapiref:chappyapi-psdconstraint}]{\sphinxcrossref{\DUrole{std,std-ref}{PsdConstraint Class}}}} object.
\end{quote}

\sphinxAtStartPar
\sphinxstylestrong{Arguments}
\begin{quote}

\sphinxAtStartPar
\sphinxcode{\sphinxupquote{idx}}
\begin{quote}

\sphinxAtStartPar
Index for the positive semi\sphinxhyphen{}definite constraint, starting with 0.
\end{quote}
\end{quote}

\sphinxAtStartPar
\sphinxstylestrong{Example}
\end{quote}

\begin{sphinxVerbatim}[commandchars=\\\{\}]
\PYG{c+c1}{\PYGZsh{} Retrieve the positive semi\PYGZhy{}definite constraint with index of 1}
\PYG{n}{r} \PYG{o}{=} \PYG{n}{m}\PYG{o}{.}\PYG{n}{getPsdConstr}\PYG{p}{(}\PYG{l+m+mi}{1}\PYG{p}{)}
\end{sphinxVerbatim}

\subsubsection{Model.getPsdConstrByName()}
\label{\detokenize{pyapiref:model-getpsdconstrbyname}}\begin{quote}

\sphinxAtStartPar
\sphinxstylestrong{Synopsis}
\begin{quote}

\sphinxAtStartPar
\sphinxcode{\sphinxupquote{getPsdConstrByName(name)}}
\end{quote}

\sphinxAtStartPar
\sphinxstylestrong{Description}
\begin{quote}

\sphinxAtStartPar
Retrieve a positive semi\sphinxhyphen{}definite constraint by name. Return a {\hyperref[\detokenize{pyapiref:chappyapi-psdconstraint}]{\sphinxcrossref{\DUrole{std,std-ref}{PsdConstraint Class}}}} object.
\end{quote}

\sphinxAtStartPar
\sphinxstylestrong{Arguments}
\begin{quote}

\sphinxAtStartPar
\sphinxcode{\sphinxupquote{name}}
\begin{quote}

\sphinxAtStartPar
The name of the positive semi\sphinxhyphen{}definite constraint.
\end{quote}
\end{quote}

\sphinxAtStartPar
\sphinxstylestrong{Example}
\end{quote}

\begin{sphinxVerbatim}[commandchars=\\\{\}]
\PYG{c+c1}{\PYGZsh{} Retrieve the positive semi\PYGZhy{}definite constraint with name \PYGZdq{}r\PYGZdq{}.}
\PYG{n}{r} \PYG{o}{=} \PYG{n}{m}\PYG{o}{.}\PYG{n}{getPsdConstrByName}\PYG{p}{(}\PYG{l+s+s2}{\PYGZdq{}}\PYG{l+s+s2}{r}\PYG{l+s+s2}{\PYGZdq{}}\PYG{p}{)}
\end{sphinxVerbatim}

\subsubsection{Model.getPsdConstrs()}
\label{\detokenize{pyapiref:model-getpsdconstrs}}\begin{quote}

\sphinxAtStartPar
\sphinxstylestrong{Synopsis}
\begin{quote}

\sphinxAtStartPar
\sphinxcode{\sphinxupquote{getPsdConstrs()}}
\end{quote}

\sphinxAtStartPar
\sphinxstylestrong{Description}
\begin{quote}

\sphinxAtStartPar
Retrieve all positive semi\sphinxhyphen{}definite constraints in the model. Return a {\hyperref[\detokenize{pyapiref:chappyapi-psdconstrarray}]{\sphinxcrossref{\DUrole{std,std-ref}{PsdConstrArray Class}}}} object.
\end{quote}

\sphinxAtStartPar
\sphinxstylestrong{Example}
\end{quote}

\begin{sphinxVerbatim}[commandchars=\\\{\}]
\PYG{c+c1}{\PYGZsh{} Retrieve all positive semi\PYGZhy{}definite constraints in the model}
\PYG{n}{cons} \PYG{o}{=} \PYG{n}{m}\PYG{o}{.}\PYG{n}{getPsdConstrs}\PYG{p}{(}\PYG{p}{)}
\end{sphinxVerbatim}

\subsubsection{Model.getPsdConstrBuilders()}
\label{\detokenize{pyapiref:model-getpsdconstrbuilders}}\begin{quote}

\sphinxAtStartPar
\sphinxstylestrong{Synopsis}
\begin{quote}

\sphinxAtStartPar
\sphinxcode{\sphinxupquote{getPsdConstrBuilders(constrs=None)}}
\end{quote}

\sphinxAtStartPar
\sphinxstylestrong{Description}
\begin{quote}

\sphinxAtStartPar
Retrieve positive semi\sphinxhyphen{}definite constraint builders in current model.

\sphinxAtStartPar
If parameter \sphinxcode{\sphinxupquote{constrs}} is \sphinxcode{\sphinxupquote{None}}, then return a {\hyperref[\detokenize{pyapiref:chappyapi-psdconstrbuilderarray}]{\sphinxcrossref{\DUrole{std,std-ref}{PsdConstrBuilderArray Class}}}} object composed of all
positive semi\sphinxhyphen{}definite constraint builders.

\sphinxAtStartPar
If parameter \sphinxcode{\sphinxupquote{constrs}} is {\hyperref[\detokenize{pyapiref:chappyapi-psdconstraint}]{\sphinxcrossref{\DUrole{std,std-ref}{PsdConstraint Class}}}} object, then return the {\hyperref[\detokenize{pyapiref:chappyapi-psdconstrbuilder}]{\sphinxcrossref{\DUrole{std,std-ref}{PsdConstrBuilder Class}}}}
object corresponding to the specific positive semi\sphinxhyphen{}definite constraint.

\sphinxAtStartPar
If parameter \sphinxcode{\sphinxupquote{constrs}} is a list or a {\hyperref[\detokenize{pyapiref:chappyapi-psdconstrarray}]{\sphinxcrossref{\DUrole{std,std-ref}{PsdConstrArray Class}}}} object, then return a {\hyperref[\detokenize{pyapiref:chappyapi-psdconstrbuilderarray}]{\sphinxcrossref{\DUrole{std,std-ref}{PsdConstrBuilderArray Class}}}}
object composed of specified positive semi\sphinxhyphen{}definite constraints’ builders.

\sphinxAtStartPar
If parameter \sphinxcode{\sphinxupquote{constrs}} is dictionary or {\hyperref[\detokenize{pyapiref:chappyapi-util-tupledict}]{\sphinxcrossref{\DUrole{std,std-ref}{tupledict Class}}}} object, then the indice of the specified
positive semi\sphinxhyphen{}definite constraint is returned as key, the value is a {\hyperref[\detokenize{pyapiref:chappyapi-util-tupledict}]{\sphinxcrossref{\DUrole{std,std-ref}{tupledict Class}}}} object
composed of the specified positive semi\sphinxhyphen{}definite constraints’ builders.
\end{quote}

\sphinxAtStartPar
\sphinxstylestrong{Arguments}
\begin{quote}

\sphinxAtStartPar
\sphinxcode{\sphinxupquote{constrs}}
\begin{quote}

\sphinxAtStartPar
The specified positive semi\sphinxhyphen{}definite constraint. Optional, \sphinxcode{\sphinxupquote{None}} by default.
\end{quote}
\end{quote}

\sphinxAtStartPar
\sphinxstylestrong{Example}
\end{quote}

\begin{sphinxVerbatim}[commandchars=\\\{\}]
\PYG{c+c1}{\PYGZsh{} Retrieve all of positive semi\PYGZhy{}definite constraint builders.}
\PYG{n}{conbuilders} \PYG{o}{=} \PYG{n}{m}\PYG{o}{.}\PYG{n}{getPsdConstrBuilders}\PYG{p}{(}\PYG{p}{)}
\PYG{c+c1}{\PYGZsh{} Retrive the builder corresponding to positive semi\PYGZhy{}definite contstraint x.}
\PYG{n}{conbuilders} \PYG{o}{=} \PYG{n}{m}\PYG{o}{.}\PYG{n}{getPsdConstrBuilders}\PYG{p}{(}\PYG{n}{x}\PYG{p}{)}
\PYG{c+c1}{\PYGZsh{} Retrieve builders corresponding to positive semi\PYGZhy{}definite constraint x and y.}
\PYG{n}{conbuilders} \PYG{o}{=} \PYG{n}{m}\PYG{o}{.}\PYG{n}{getPsdConstrBuilders}\PYG{p}{(}\PYG{p}{[}\PYG{n}{x}\PYG{p}{,} \PYG{n}{y}\PYG{p}{]}\PYG{p}{)}
\PYG{c+c1}{\PYGZsh{} Retrieve builders corresponding to positive semi\PYGZhy{}definite constraint in tupledict object xx.}
\PYG{n}{conbuilders} \PYG{o}{=} \PYG{n}{m}\PYG{o}{.}\PYG{n}{getPsdConstrBuilders}\PYG{p}{(}\PYG{n}{xx}\PYG{p}{)}
\end{sphinxVerbatim}

\subsubsection{Model.getLmiRow()}
\label{\detokenize{pyapiref:model-getlmirow}}\begin{quote}

\sphinxAtStartPar
\sphinxstylestrong{Synopsis}
\begin{quote}

\sphinxAtStartPar
\sphinxcode{\sphinxupquote{getLmiRow(constr)}}
\end{quote}

\sphinxAtStartPar
\sphinxstylestrong{Description}
\begin{quote}

\sphinxAtStartPar
Get the LMI expression involved in the specified LMI constraint, including variables and corresponding coefficient matrices.
\end{quote}

\sphinxAtStartPar
\sphinxstylestrong{Arguments}
\begin{quote}

\sphinxAtStartPar
\sphinxcode{\sphinxupquote{constr}}
\begin{quote}

\sphinxAtStartPar
The specified LMI constraint.
\end{quote}
\end{quote}

\sphinxAtStartPar
\sphinxstylestrong{Example}
\end{quote}

\begin{sphinxVerbatim}[commandchars=\\\{\}]
\PYG{c+c1}{\PYGZsh{} Get the expression in LMI constraint c}
\PYG{n}{expr} \PYG{o}{=} \PYG{n}{m}\PYG{o}{.}\PYG{n}{getLmiRow}\PYG{p}{(}\PYG{n}{c}\PYG{p}{)}
\end{sphinxVerbatim}

\subsubsection{Model.getLmiConstr()}
\label{\detokenize{pyapiref:model-getlmiconstr}}\begin{quote}

\sphinxAtStartPar
\sphinxstylestrong{Synopsis}
\begin{quote}

\sphinxAtStartPar
\sphinxcode{\sphinxupquote{getLmiConstr(idx)}}
\end{quote}

\sphinxAtStartPar
\sphinxstylestrong{Description}
\begin{quote}

\sphinxAtStartPar
Get the LMI constraint corresponding to the specified index in the model.
\end{quote}

\sphinxAtStartPar
\sphinxstylestrong{Arguments}
\begin{quote}

\sphinxAtStartPar
\sphinxcode{\sphinxupquote{idx}}
\begin{quote}

\sphinxAtStartPar
The index of the LMI constraint in the model. Starts with \sphinxcode{\sphinxupquote{0}}.
\end{quote}
\end{quote}

\sphinxAtStartPar
\sphinxstylestrong{Example}
\end{quote}

\begin{sphinxVerbatim}[commandchars=\\\{\}]
\PYG{c+c1}{\PYGZsh{} Get the 1st LMI constraint in the model}
\PYG{n}{coeff} \PYG{o}{=} \PYG{n}{m}\PYG{o}{.}\PYG{n}{getLmiConstr}\PYG{p}{(}\PYG{l+m+mi}{1}\PYG{p}{)}
\end{sphinxVerbatim}

\subsubsection{Model.getLmiConstrByName()}
\label{\detokenize{pyapiref:model-getlmiconstrbyname}}\begin{quote}

\sphinxAtStartPar
\sphinxstylestrong{Synopsis}
\begin{quote}

\sphinxAtStartPar
\sphinxcode{\sphinxupquote{getLmiConstrByName(name)}}
\end{quote}

\sphinxAtStartPar
\sphinxstylestrong{Description}
\begin{quote}

\sphinxAtStartPar
Get the LMI constraint of the specified name in the model.
\end{quote}

\sphinxAtStartPar
\sphinxstylestrong{Arguments}
\begin{quote}

\sphinxAtStartPar
\sphinxcode{\sphinxupquote{name}}
\begin{quote}

\sphinxAtStartPar
The specified name of the LMI constraint.
\end{quote}
\end{quote}

\sphinxAtStartPar
\sphinxstylestrong{Example}
\end{quote}

\begin{sphinxVerbatim}[commandchars=\\\{\}]
\PYG{c+c1}{\PYGZsh{} Get the LMI constraint named r1 in the model}
\PYG{n}{name} \PYG{o}{=} \PYG{n}{m}\PYG{o}{.}\PYG{n}{getLmiConstrByName}\PYG{p}{(}\PYG{l+s+s2}{\PYGZdq{}}\PYG{l+s+s2}{r1}\PYG{l+s+s2}{\PYGZdq{}}\PYG{p}{)}
\end{sphinxVerbatim}

\subsubsection{Model.getLmiConstrs()}
\label{\detokenize{pyapiref:model-getlmiconstrs}}\begin{quote}

\sphinxAtStartPar
\sphinxstylestrong{Synopsis}
\begin{quote}

\sphinxAtStartPar
\sphinxcode{\sphinxupquote{getLmiConstrs()}}
\end{quote}

\sphinxAtStartPar
\sphinxstylestrong{Description}
\begin{quote}

\sphinxAtStartPar
Get all LMI constraints in the model. Returns a {\hyperref[\detokenize{pyapiref:chappyapi-lmiconstrarray}]{\sphinxcrossref{\DUrole{std,std-ref}{LmiConstrArray Class}}}} object composed of LMI constraints.
\end{quote}
\end{quote}

\subsubsection{Model.getLmiRhs()}
\label{\detokenize{pyapiref:model-getlmirhs}}\begin{quote}

\sphinxAtStartPar
\sphinxstylestrong{Synopsis}
\begin{quote}

\sphinxAtStartPar
\sphinxcode{\sphinxupquote{getLmiRhs(constr)}}
\end{quote}

\sphinxAtStartPar
\sphinxstylestrong{Description}
\begin{quote}

\sphinxAtStartPar
Get the constant term of the specified LMI constraint. Returns a {\hyperref[\detokenize{pyapiref:chappyapi-symmatrix}]{\sphinxcrossref{\DUrole{std,std-ref}{SymMatrix Class}}}} object.
\end{quote}

\sphinxAtStartPar
\sphinxstylestrong{Arguments}
\begin{quote}

\sphinxAtStartPar
\sphinxcode{\sphinxupquote{constr}}
\begin{quote}

\sphinxAtStartPar
The specified LMI constraint.
\end{quote}
\end{quote}
\end{quote}

\subsubsection{Model.setLmiRhs()}
\label{\detokenize{pyapiref:model-setlmirhs}}\begin{quote}

\sphinxAtStartPar
\sphinxstylestrong{Synopsis}
\begin{quote}

\sphinxAtStartPar
\sphinxcode{\sphinxupquote{setLmiRhs(constr, mat)}}
\end{quote}

\sphinxAtStartPar
\sphinxstylestrong{Description}
\begin{quote}

\sphinxAtStartPar
Set the constant term of the specified LMI constraint.
\end{quote}

\sphinxAtStartPar
\sphinxstylestrong{Arguments}
\begin{quote}

\sphinxAtStartPar
\sphinxcode{\sphinxupquote{constr}}
\begin{quote}

\sphinxAtStartPar
The specified LMI constraint.
\end{quote}

\sphinxAtStartPar
\sphinxcode{\sphinxupquote{mat}}
\begin{quote}

\sphinxAtStartPar
The new constant\sphinxhyphen{}term symmetric to set.
\end{quote}
\end{quote}

\sphinxAtStartPar
\sphinxstylestrong{Example}
\end{quote}

\begin{sphinxVerbatim}[commandchars=\\\{\}]
\PYG{c+c1}{\PYGZsh{} Set the constant\PYGZhy{}term symmetric of the LMI constraint con to D}
\PYG{n}{m}\PYG{o}{.}\PYG{n}{setLmiRhs}\PYG{p}{(}\PYG{n}{con}\PYG{p}{,} \PYG{n}{D}\PYG{p}{)}
\end{sphinxVerbatim}

\subsubsection{Model.getLmiSolution()}
\label{\detokenize{pyapiref:model-getlmisolution}}\begin{quote}

\sphinxAtStartPar
\sphinxstylestrong{Synopsis}
\begin{quote}

\sphinxAtStartPar
\sphinxcode{\sphinxupquote{getLmiSolution()}}
\end{quote}

\sphinxAtStartPar
\sphinxstylestrong{Description}
\begin{quote}

\sphinxAtStartPar
Get the value and dual value of the LMI constraint.
\end{quote}
\end{quote}

\subsubsection{Model.getLmiSlacks()}
\label{\detokenize{pyapiref:model-getlmislacks}}\begin{quote}

\sphinxAtStartPar
\sphinxstylestrong{Synopsis}
\begin{quote}

\sphinxAtStartPar
\sphinxcode{\sphinxupquote{getLmiSlacks()}}
\end{quote}

\sphinxAtStartPar
\sphinxstylestrong{Description}
\begin{quote}

\sphinxAtStartPar
Get the values of all slack variables of LMI constraints. Returns a list object.
\end{quote}
\end{quote}

\subsubsection{Model.getLmiDuals()}
\label{\detokenize{pyapiref:model-getlmiduals}}\begin{quote}

\sphinxAtStartPar
\sphinxstylestrong{Synopsis}
\begin{quote}

\sphinxAtStartPar
\sphinxcode{\sphinxupquote{getLmiDuals()}}
\end{quote}

\sphinxAtStartPar
\sphinxstylestrong{Description}
\begin{quote}

\sphinxAtStartPar
Get the values of all dual variables of the LMI constraint. Returns a list object.
\end{quote}
\end{quote}

\subsubsection{Model.getCoeff()}
\label{\detokenize{pyapiref:model-getcoeff}}\begin{quote}

\sphinxAtStartPar
\sphinxstylestrong{Synopsis}
\begin{quote}

\sphinxAtStartPar
\sphinxcode{\sphinxupquote{getCoeff(constr, var)}}
\end{quote}

\sphinxAtStartPar
\sphinxstylestrong{Description}
\begin{quote}

\sphinxAtStartPar
Get the coefficient of variable in linear constraint, PSD constraint or LMI constraint.
\end{quote}

\sphinxAtStartPar
\sphinxstylestrong{Arguments}
\begin{quote}

\sphinxAtStartPar
\sphinxcode{\sphinxupquote{constr}}
\begin{quote}

\sphinxAtStartPar
The requested linear constraint, PSD constraint or LMI constraint.
\end{quote}

\sphinxAtStartPar
\sphinxcode{\sphinxupquote{var}}
\begin{quote}

\sphinxAtStartPar
The requested variable or PSD variable.
\end{quote}
\end{quote}

\sphinxAtStartPar
\sphinxstylestrong{Example}
\end{quote}

\begin{sphinxVerbatim}[commandchars=\\\{\}]
\PYG{c+c1}{\PYGZsh{} Get the coefficient of variable x in linear constraint c1}
\PYG{n}{coeff1} \PYG{o}{=} \PYG{n}{m}\PYG{o}{.}\PYG{n}{getCoeff}\PYG{p}{(}\PYG{n}{c1}\PYG{p}{,} \PYG{n}{x}\PYG{p}{)}
\PYG{c+c1}{\PYGZsh{} Get the coefficient of PSD variable X in PSD constraint c2}
\PYG{n}{coeff2} \PYG{o}{=} \PYG{n}{m}\PYG{o}{.}\PYG{n}{getCoeff}\PYG{p}{(}\PYG{n}{c2}\PYG{p}{,} \PYG{n}{X}\PYG{p}{)}
\end{sphinxVerbatim}

\subsubsection{Model.setCoeff()}
\label{\detokenize{pyapiref:model-setcoeff}}\begin{quote}

\sphinxAtStartPar
\sphinxstylestrong{Synopsis}
\begin{quote}

\sphinxAtStartPar
\sphinxcode{\sphinxupquote{setCoeff(constr, var, newval)}}
\end{quote}

\sphinxAtStartPar
\sphinxstylestrong{Description}
\begin{quote}

\sphinxAtStartPar
Set the coefficient of variable in linear constraint, PSD constraint or LMI constraint.
\end{quote}

\sphinxAtStartPar
\sphinxstylestrong{Arguments}
\begin{quote}

\sphinxAtStartPar
\sphinxcode{\sphinxupquote{constr}}
\begin{quote}

\sphinxAtStartPar
The requested linear constraint, PSD constraint or LMI constraint.
\end{quote}

\sphinxAtStartPar
\sphinxcode{\sphinxupquote{var}}
\begin{quote}

\sphinxAtStartPar
The requested variable or PSD variable.
\end{quote}

\sphinxAtStartPar
\sphinxcode{\sphinxupquote{newval}}
\begin{quote}

\sphinxAtStartPar
New coefficient or symmetric matrix coefficient.
\end{quote}
\end{quote}

\sphinxAtStartPar
\sphinxstylestrong{Example}
\end{quote}

\begin{sphinxVerbatim}[commandchars=\\\{\}]
\PYG{c+c1}{\PYGZsh{} Set the coefficient of variable x in linear constraint c to 1.0}
\PYG{n}{m}\PYG{o}{.}\PYG{n}{setCoeff}\PYG{p}{(}\PYG{n}{c}\PYG{p}{,} \PYG{n}{x}\PYG{p}{,} \PYG{l+m+mf}{1.0}\PYG{p}{)}
\end{sphinxVerbatim}

\subsubsection{Model.setCoeffs()}
\label{\detokenize{pyapiref:model-setcoeffs}}\begin{quote}

\sphinxAtStartPar
\sphinxstylestrong{Synopsis}
\begin{quote}

\sphinxAtStartPar
\sphinxcode{\sphinxupquote{setCoeffs(constrs, vars, vals)}}
\end{quote}

\sphinxAtStartPar
\sphinxstylestrong{Description}
\begin{quote}

\sphinxAtStartPar
Set the coefficients of variables in the linear constraints in batches.

\sphinxAtStartPar
\sphinxstylestrong{Note} The constraint and variable pair cannot repeat with the same or different coefficient to be set.
\end{quote}

\sphinxAtStartPar
\sphinxstylestrong{Arguments}
\begin{quote}

\sphinxAtStartPar
\sphinxcode{\sphinxupquote{constrs}}
\begin{quote}

\sphinxAtStartPar
Specifies the constraints related to the coefficients to be set, which could be a dictionary, {\hyperref[\detokenize{pyapiref:chappyapi-util-tupledict}]{\sphinxcrossref{\DUrole{std,std-ref}{tupledict Class}}}} , {\hyperref[\detokenize{pyapiref:chappyapi-constrarray}]{\sphinxcrossref{\DUrole{std,std-ref}{ConstrArray Class}}}} or a list of {\hyperref[\detokenize{pyapiref:chappyapi-constraint}]{\sphinxcrossref{\DUrole{std,std-ref}{Constraint Class}}}} objects.
\end{quote}

\sphinxAtStartPar
\sphinxcode{\sphinxupquote{vars}}
\begin{quote}

\sphinxAtStartPar
Specifies the variables related to the coefficients to be set, which could be a dictionary, {\hyperref[\detokenize{pyapiref:chappyapi-util-tupledict}]{\sphinxcrossref{\DUrole{std,std-ref}{tupledict Class}}}} , {\hyperref[\detokenize{pyapiref:chappyapi-vararray}]{\sphinxcrossref{\DUrole{std,std-ref}{VarArray Class}}}} or a list of {\hyperref[\detokenize{pyapiref:chappyapi-var}]{\sphinxcrossref{\DUrole{std,std-ref}{Var Class}}}} objects.
\end{quote}

\sphinxAtStartPar
\sphinxcode{\sphinxupquote{vals}}
\begin{quote}

\sphinxAtStartPar
The new coefficient values to be set, which could be a constant, or the list/dictionary corresponding to the \sphinxcode{\sphinxupquote{constrs}} .
\end{quote}
\end{quote}
\end{quote}

\subsubsection{Model.getA()}
\label{\detokenize{pyapiref:model-geta}}\begin{quote}

\sphinxAtStartPar
\sphinxstylestrong{Synopsis}
\begin{quote}

\sphinxAtStartPar
\sphinxcode{\sphinxupquote{getA()}}
\end{quote}

\sphinxAtStartPar
\sphinxstylestrong{Description}
\begin{quote}

\sphinxAtStartPar
Get the coefficient matrix of model, returns a \sphinxcode{\sphinxupquote{scipy.sparse.csc\_matrix}} object.
This method requires the \sphinxcode{\sphinxupquote{scipy}} package.
\end{quote}

\sphinxAtStartPar
\sphinxstylestrong{Example}
\end{quote}

\begin{sphinxVerbatim}[commandchars=\\\{\}]
\PYG{c+c1}{\PYGZsh{} Get the coefficient matrix}
\PYG{n}{A} \PYG{o}{=} \PYG{n}{model}\PYG{o}{.}\PYG{n}{getA}\PYG{p}{(}\PYG{p}{)}
\end{sphinxVerbatim}

\subsubsection{Model.loadMatrix()}
\label{\detokenize{pyapiref:model-loadmatrix}}\begin{quote}

\sphinxAtStartPar
\sphinxstylestrong{Synopsis}
\begin{quote}

\sphinxAtStartPar
\sphinxcode{\sphinxupquote{loadMatrix(c, A, lhs, rhs, lb, ub, vtype=None)}}
\end{quote}

\sphinxAtStartPar
\sphinxstylestrong{Description}
\begin{quote}

\sphinxAtStartPar
Load matrix and vector data to build model. This method requires the \sphinxcode{\sphinxupquote{scipy}} package.
\end{quote}

\sphinxAtStartPar
\sphinxstylestrong{Arguments}
\begin{quote}

\sphinxAtStartPar
\sphinxcode{\sphinxupquote{c}}
\begin{quote}

\sphinxAtStartPar
Objective costs. If \sphinxcode{\sphinxupquote{None}}, the objective costs are all zeros.
\end{quote}

\sphinxAtStartPar
\sphinxcode{\sphinxupquote{A}}
\begin{quote}

\sphinxAtStartPar
Coefficient matrix. Must be of type \sphinxcode{\sphinxupquote{scipy.sparse.csc\_matrix}}.
\end{quote}

\sphinxAtStartPar
\sphinxcode{\sphinxupquote{lhs}}
\begin{quote}

\sphinxAtStartPar
Lower bounds of constraints.
\end{quote}

\sphinxAtStartPar
\sphinxcode{\sphinxupquote{rhs}}
\begin{quote}

\sphinxAtStartPar
Upper bounds of constraints.
\end{quote}

\sphinxAtStartPar
\sphinxcode{\sphinxupquote{lb}}
\begin{quote}

\sphinxAtStartPar
Lower bounds of variables. If \sphinxcode{\sphinxupquote{None}}, the lower bounds are all zeros.
\end{quote}

\sphinxAtStartPar
\sphinxcode{\sphinxupquote{ub}}
\begin{quote}

\sphinxAtStartPar
Upper bounds of variables. If \sphinxcode{\sphinxupquote{None}}, the upper bounds are all \sphinxcode{\sphinxupquote{COPT.INFINITY}}.
\end{quote}

\sphinxAtStartPar
\sphinxcode{\sphinxupquote{vtype}}
\begin{quote}

\sphinxAtStartPar
Variable types. Default to \sphinxcode{\sphinxupquote{None}}, which means all variables are continuous.
\end{quote}
\end{quote}

\sphinxAtStartPar
\sphinxstylestrong{Example}
\end{quote}

\begin{sphinxVerbatim}[commandchars=\\\{\}]
\PYG{c+c1}{\PYGZsh{} Build model by problem matrix}
\PYG{n}{m}\PYG{o}{.}\PYG{n}{loadMatrix}\PYG{p}{(}\PYG{n}{c}\PYG{p}{,} \PYG{n}{A}\PYG{p}{,} \PYG{n}{lhs}\PYG{p}{,} \PYG{n}{rhs}\PYG{p}{,} \PYG{n}{lb}\PYG{p}{,} \PYG{n}{ub}\PYG{p}{)}
\end{sphinxVerbatim}

\subsubsection{Model.loadCone()}
\label{\detokenize{pyapiref:model-loadcone}}\begin{quote}

\sphinxAtStartPar
\sphinxstylestrong{Synopsis}
\begin{quote}

\sphinxAtStartPar
\sphinxcode{\sphinxupquote{loadCone(ncone, types, dims, indices)}}
\end{quote}

\sphinxAtStartPar
\sphinxstylestrong{Description}
\begin{quote}

\sphinxAtStartPar
Load Second\sphinxhyphen{}Order\sphinxhyphen{}Cones (SOC) to the model.
\end{quote}

\sphinxAtStartPar
\sphinxstylestrong{Arguments}
\begin{quote}

\sphinxAtStartPar
\sphinxcode{\sphinxupquote{ncone}}
\begin{quote}

\sphinxAtStartPar
Number of SOC.
\end{quote}

\sphinxAtStartPar
\sphinxcode{\sphinxupquote{types}}
\begin{quote}

\sphinxAtStartPar
Type of SOC, please refer to {\hyperref[\detokenize{constant:chapconst-conetype}]{\sphinxcrossref{\DUrole{std,std-ref}{SOC constraint types}}}} for possible values.
\end{quote}

\sphinxAtStartPar
\sphinxcode{\sphinxupquote{dims}}
\begin{quote}

\sphinxAtStartPar
Dimension of SOC.
\end{quote}

\sphinxAtStartPar
\sphinxcode{\sphinxupquote{indices}}
\begin{quote}

\sphinxAtStartPar
Array of subscripts for the variables that constitute the SOC.
\end{quote}
\end{quote}
\end{quote}

\subsubsection{Model.loadExpCone()}
\label{\detokenize{pyapiref:model-loadexpcone}}\begin{quote}

\sphinxAtStartPar
\sphinxstylestrong{Synopsis}
\begin{quote}

\sphinxAtStartPar
\sphinxcode{\sphinxupquote{loadExpCone(ncone, types, indices)}}
\end{quote}

\sphinxAtStartPar
\sphinxstylestrong{Description}
\begin{quote}

\sphinxAtStartPar
Load exponential cones to the model.
\end{quote}

\sphinxAtStartPar
\sphinxstylestrong{Arguments}
\begin{quote}

\sphinxAtStartPar
\sphinxcode{\sphinxupquote{ncone}}
\begin{quote}

\sphinxAtStartPar
Number of exponential cones.
\end{quote}

\sphinxAtStartPar
\sphinxcode{\sphinxupquote{types}}
\begin{quote}

\sphinxAtStartPar
Type of exponential cones, please refer to {\hyperref[\detokenize{constant:chapconst-expconetype}]{\sphinxcrossref{\DUrole{std,std-ref}{Exponential Cone type}}}} for possible values.
\end{quote}

\sphinxAtStartPar
\sphinxcode{\sphinxupquote{indices}}
\begin{quote}

\sphinxAtStartPar
Array of subscripts for the variables that constitute the exponential cones.
\end{quote}
\end{quote}
\end{quote}

\subsubsection{Model.loadNlData()}
\label{\detokenize{pyapiref:model-loadnldata}}\begin{quote}

\sphinxAtStartPar
\sphinxstylestrong{Synopsis}
\begin{quote}

\sphinxAtStartPar
\sphinxcode{\sphinxupquote{loadNlData(nCols, nRows, sense, nGrad, idxGrad, nJac, idxJacRow, idxJacCol,
nHess, idxHessRow, idxHessCol, colLower, colUpper, rowLower, rowUpper,
initX, evalType, cb)}}
\end{quote}

\sphinxAtStartPar
\sphinxstylestrong{Description}
\begin{quote}

\sphinxAtStartPar
Load the structural information and initial data of a nonlinear optimization model,
and register a user\sphinxhyphen{}implemented nonlinear callback object.

\sphinxAtStartPar
This interface is used to provide the solver with the model size information,
variable and constraint bounds, and the sparsity structure of derivatives
of the objective and constraint functions, while registering callback interfaces
for evaluation in nonlinear optimization.

\sphinxAtStartPar
During the nonlinear optimization process, the solver will invoke
the callback object \sphinxcode{\sphinxupquote{cb}} at different iteration points to dynamically evaluate
the values of the objective function, constraint functions,
and their first\sphinxhyphen{}order and second\sphinxhyphen{}order derivatives, to advance the
nonlinear optimization algorithm.

\sphinxAtStartPar
This interface adopts a callback\sphinxhyphen{}based nonlinear modeling approach,
decoupling model structure definition from numerical computation logic.
Users provide the required model and derivative information through
callback methods.
\end{quote}

\sphinxAtStartPar
\sphinxstylestrong{Arguments}
\begin{quote}

\sphinxAtStartPar
\sphinxcode{\sphinxupquote{nCols}}
\begin{quote}

\sphinxAtStartPar
Number of variables.

\sphinxAtStartPar
Possible values: positive integers.
\end{quote}

\sphinxAtStartPar
\sphinxcode{\sphinxupquote{nRows}}
\begin{quote}

\sphinxAtStartPar
Number of constraints.

\sphinxAtStartPar
Possible values: non\sphinxhyphen{}negative integers.
\end{quote}

\sphinxAtStartPar
\sphinxcode{\sphinxupquote{sense}}
\begin{quote}

\sphinxAtStartPar
Optimization sense.

\sphinxAtStartPar
Possible values: \sphinxcode{\sphinxupquote{COPT.MINIMIZE}} or \sphinxcode{\sphinxupquote{COPT.MAXIMIZE}}.
\end{quote}

\sphinxAtStartPar
\sphinxcode{\sphinxupquote{nGrad}}
\begin{quote}

\sphinxAtStartPar
Number of nonzero elements in the objective gradient, or an identifier
specifying the gradient storage format.

\sphinxAtStartPar
Possible values: a non\sphinxhyphen{}negative integer, or a dense storage type
such as \sphinxcode{\sphinxupquote{COPT.DENSETYPE\_ROWMAJOR}} or \sphinxcode{\sphinxupquote{COPT.DENSETYPE\_COLMAJOR}}.
\end{quote}

\sphinxAtStartPar
\sphinxcode{\sphinxupquote{idxGrad}}
\begin{quote}

\sphinxAtStartPar
Index sequence of variables corresponding to the nonzero elements
in the objective gradient.

\sphinxAtStartPar
Possible values: a sequence of integers, or \sphinxcode{\sphinxupquote{None}}.
\end{quote}

\sphinxAtStartPar
\sphinxcode{\sphinxupquote{nJac}}
\begin{quote}

\sphinxAtStartPar
Number of nonzero elements in the constraint Jacobian matrix,
or an identifier specifying the Jacobian storage format.

\sphinxAtStartPar
Possible values: a non\sphinxhyphen{}negative integer, or a dense storage type
such as \sphinxcode{\sphinxupquote{COPT.DENSETYPE\_ROWMAJOR}} or \sphinxcode{\sphinxupquote{COPT.DENSETYPE\_COLMAJOR}}.
\end{quote}

\sphinxAtStartPar
\sphinxcode{\sphinxupquote{idxJacRow}}
\begin{quote}

\sphinxAtStartPar
Row index sequence of nonzero elements in the Jacobian matrix.

\sphinxAtStartPar
Possible values: a sequence of integers, or \sphinxcode{\sphinxupquote{None}}.
\end{quote}

\sphinxAtStartPar
\sphinxcode{\sphinxupquote{idxJacCol}}
\begin{quote}

\sphinxAtStartPar
Column index sequence of nonzero elements in the Jacobian matrix.

\sphinxAtStartPar
Possible values: a sequence of integers, or \sphinxcode{\sphinxupquote{None}}.
\end{quote}

\sphinxAtStartPar
\sphinxcode{\sphinxupquote{nHess}}
\begin{quote}

\sphinxAtStartPar
Number of nonzero elements in the Hessian matrix,
or an identifier specifying the Hessian storage format.

\sphinxAtStartPar
Possible values: a non\sphinxhyphen{}negative integer, or a dense storage type
such as \sphinxcode{\sphinxupquote{COPT.DENSETYPE\_ROWMAJOR}} or \sphinxcode{\sphinxupquote{COPT.DENSETYPE\_COLMAJOR}}.
\end{quote}

\sphinxAtStartPar
\sphinxcode{\sphinxupquote{idxHessRow}}
\begin{quote}

\sphinxAtStartPar
Row index sequence of nonzero elements in the Hessian matrix.

\sphinxAtStartPar
Possible values: a sequence of integers, or \sphinxcode{\sphinxupquote{None}}.
\end{quote}

\sphinxAtStartPar
\sphinxcode{\sphinxupquote{idxHessCol}}
\begin{quote}

\sphinxAtStartPar
Column index sequence of nonzero elements in the Hessian matrix.

\sphinxAtStartPar
Possible values: a sequence of integers, or \sphinxcode{\sphinxupquote{None}}.
\end{quote}

\sphinxAtStartPar
\sphinxcode{\sphinxupquote{colLower}}
\begin{quote}

\sphinxAtStartPar
Lower bounds of variables.

\sphinxAtStartPar
Possible values: a sequence of \sphinxcode{\sphinxupquote{float}} numbers, or \sphinxcode{\sphinxupquote{None}}.

\sphinxAtStartPar
If \sphinxcode{\sphinxupquote{None}}, default lower bounds 0 are applied to all variables.
\end{quote}

\sphinxAtStartPar
\sphinxcode{\sphinxupquote{colUpper}}
\begin{quote}

\sphinxAtStartPar
Upper bounds of variables.

\sphinxAtStartPar
Possible values: a sequence of \sphinxcode{\sphinxupquote{float}} numbers, or \sphinxcode{\sphinxupquote{None}}.

\sphinxAtStartPar
If \sphinxcode{\sphinxupquote{None}}, default upper bounds \sphinxcode{\sphinxupquote{COPT.INFINITY}} are applied to all variables.
\end{quote}

\sphinxAtStartPar
\sphinxcode{\sphinxupquote{rowLower}}
\begin{quote}

\sphinxAtStartPar
Lower bounds of constraint functions.

\sphinxAtStartPar
Possible values: a sequence of \sphinxcode{\sphinxupquote{float}} numbers.

\sphinxAtStartPar
If \sphinxcode{\sphinxupquote{None}}, default lower bounds \sphinxcode{\sphinxupquote{\sphinxhyphen{}COPT.INFINITY}} are applied to all constraints.
\end{quote}

\sphinxAtStartPar
\sphinxcode{\sphinxupquote{rowUpper}}
\begin{quote}

\sphinxAtStartPar
Upper bounds of constraint functions.

\sphinxAtStartPar
Possible values: a sequence of \sphinxcode{\sphinxupquote{float}} numbers.

\sphinxAtStartPar
If \sphinxcode{\sphinxupquote{None}}, default upper bounds \sphinxcode{\sphinxupquote{COPT.INFINITY}} are applied to all constraints.
\end{quote}

\sphinxAtStartPar
\sphinxcode{\sphinxupquote{initX}}
\begin{quote}

\sphinxAtStartPar
Initial values of variables.

\sphinxAtStartPar
Possible values: a sequence of \sphinxcode{\sphinxupquote{float}} numbers, or \sphinxcode{\sphinxupquote{None}}.
\end{quote}

\sphinxAtStartPar
\sphinxcode{\sphinxupquote{evalType}}
\begin{quote}

\sphinxAtStartPar
Declares which callback methods are implemented.

\sphinxAtStartPar
This argument declares which information of the nonlinear model
(including the objective function, constraint functions, and their first\sphinxhyphen{}
and second\sphinxhyphen{}order derivatives) are available through the callback interface.
During the nonlinear optimization process, the solver invokes the
corresponding callback methods as required by the algorithm.

\sphinxAtStartPar
Possible values: see {\hyperref[\detokenize{constant:chapconst-nlpevaltype}]{\sphinxcrossref{\DUrole{std,std-ref}{Nonlinear Model Callback Evaluation Types}}}}.
Multiple flags can be combined using bitwise OR.
\end{quote}

\sphinxAtStartPar
\sphinxcode{\sphinxupquote{cb}}
\begin{quote}

\sphinxAtStartPar
User\sphinxhyphen{}implemented nonlinear callback object.

\sphinxAtStartPar
Possible values: an object derived from
{\hyperref[\detokenize{pyapiref:chappyapi-nlpcbcbase}]{\sphinxcrossref{\DUrole{std,std-ref}{NlpCallbackBase Class}}}}.
\end{quote}
\end{quote}

\sphinxAtStartPar
\sphinxstylestrong{Example}

\begin{sphinxVerbatim}[commandchars=\\\{\}]
\PYG{c+c1}{\PYGZsh{} Customized NLP callback}
\PYG{n}{nlpcb} \PYG{o}{=} \PYG{n}{MyNlpCallback}\PYG{p}{(}\PYG{p}{)}

\PYG{n}{model}\PYG{o}{.}\PYG{n}{loadNlData}\PYG{p}{(}
    \PYG{n}{nCols}\PYG{o}{=}\PYG{l+m+mi}{2}\PYG{p}{,}          \PYG{c+c1}{\PYGZsh{} number of variables}
    \PYG{n}{nRows}\PYG{o}{=}\PYG{l+m+mi}{1}\PYG{p}{,}          \PYG{c+c1}{\PYGZsh{} number of constraints}
    \PYG{n}{sense}\PYG{o}{=}\PYG{n}{COPT}\PYG{o}{.}\PYG{n}{MAXIMIZE}\PYG{p}{,}

    \PYG{n}{nGrad}\PYG{o}{=}\PYG{l+m+mi}{2}\PYG{p}{,}
    \PYG{n}{idxGrad}\PYG{o}{=}\PYG{p}{[}\PYG{l+m+mi}{0}\PYG{p}{,} \PYG{l+m+mi}{1}\PYG{p}{]}\PYG{p}{,}   \PYG{c+c1}{\PYGZsh{} objective gradient structure}

    \PYG{n}{nJac}\PYG{o}{=}\PYG{l+m+mi}{2}\PYG{p}{,}
    \PYG{n}{idxJacRow}\PYG{o}{=}\PYG{p}{[}\PYG{l+m+mi}{0}\PYG{p}{,} \PYG{l+m+mi}{0}\PYG{p}{]}\PYG{p}{,} \PYG{c+c1}{\PYGZsh{} Jacobian structure: row indices}
    \PYG{n}{idxJacCol}\PYG{o}{=}\PYG{p}{[}\PYG{l+m+mi}{0}\PYG{p}{,} \PYG{l+m+mi}{1}\PYG{p}{]}\PYG{p}{,} \PYG{c+c1}{\PYGZsh{} Jacobian structure: column indices}

    \PYG{n}{nHess}\PYG{o}{=}\PYG{l+m+mi}{2}\PYG{p}{,}
    \PYG{n}{idxHessRow}\PYG{o}{=}\PYG{p}{[}\PYG{l+m+mi}{0}\PYG{p}{,} \PYG{l+m+mi}{1}\PYG{p}{]}\PYG{p}{,}  \PYG{c+c1}{\PYGZsh{} Hessian (lower triangular)}
    \PYG{n}{idxHessCol}\PYG{o}{=}\PYG{p}{[}\PYG{l+m+mi}{0}\PYG{p}{,} \PYG{l+m+mi}{1}\PYG{p}{]}\PYG{p}{,}

    \PYG{n}{colLower}\PYG{o}{=}\PYG{p}{[}\PYG{o}{\PYGZhy{}}\PYG{l+m+mf}{1.0}\PYG{p}{,} \PYG{o}{\PYGZhy{}}\PYG{n}{COPT}\PYG{o}{.}\PYG{n}{INFINITY}\PYG{p}{]}\PYG{p}{,}
    \PYG{n}{colUpper}\PYG{o}{=}\PYG{p}{[}\PYG{l+m+mf}{1.0}\PYG{p}{,}  \PYG{n}{COPT}\PYG{o}{.}\PYG{n}{INFINITY}\PYG{p}{]}\PYG{p}{,}
    \PYG{n}{rowLower}\PYG{o}{=}\PYG{p}{[}\PYG{l+m+mf}{0.0}\PYG{p}{]}\PYG{p}{,}
    \PYG{n}{rowUpper}\PYG{o}{=}\PYG{p}{[}\PYG{l+m+mf}{0.0}\PYG{p}{]}\PYG{p}{,}
    \PYG{n}{initX}\PYG{o}{=}\PYG{p}{[}\PYG{l+m+mf}{0.5}\PYG{p}{,} \PYG{l+m+mf}{1.5}\PYG{p}{]}\PYG{p}{,}

    \PYG{n}{evalType}\PYG{o}{=}\PYG{o}{\PYGZhy{}}\PYG{l+m+mi}{1}\PYG{p}{,}      \PYG{c+c1}{\PYGZsh{} all evaluation callbacks are provided}
    \PYG{n}{cb}\PYG{o}{=}\PYG{n}{nlpcb}\PYG{p}{,}
\PYG{p}{)}
\end{sphinxVerbatim}
\end{quote}

\begin{sphinxadmonition}{note}{Notes:}
\begin{enumerate}
\sphinxsetlistlabels{\arabic}{enumi}{enumii}{}{.}%
\item {} 
\sphinxAtStartPar
Setting \sphinxcode{\sphinxupquote{evalType}}

\sphinxAtStartPar
If the nonlinear callback object implements all supported numerical
computation methods, the \sphinxcode{\sphinxupquote{evalType}} parameter may be set to a
negative value (for example, \sphinxcode{\sphinxupquote{\sphinxhyphen{}1}}).

\sphinxAtStartPar
In this case, it is not necessary to explicitly combine individual
\sphinxcode{\sphinxupquote{EVALTYPE\_*}} constants.
During the nonlinear optimization process, the solver will
automatically invoke the corresponding callback methods as required
by the algorithm.

\item {} 
\sphinxAtStartPar
Rules for derivative structure arguments
(\sphinxcode{\sphinxupquote{nGrad}}, \sphinxcode{\sphinxupquote{nJac}}, and \sphinxcode{\sphinxupquote{nHess}})

\sphinxAtStartPar
The parameters \sphinxcode{\sphinxupquote{nGrad}}, \sphinxcode{\sphinxupquote{nJac}}, and \sphinxcode{\sphinxupquote{nHess}} are used to describe
the structural information of the objective gradient, the
Jacobian matrix, and the Hessian matrix, respectively.
These parameters follow the same specification rules:
\begin{itemize}
\item {} 
\sphinxAtStartPar
If the value is a non\sphinxhyphen{}negative integer, it represents the number of
nonzero elements in the corresponding derivative matrix.
In this case, the locations of the nonzero elements must be specified
through the associated index arrays
(such as \sphinxcode{\sphinxupquote{idxGrad}}, \sphinxcode{\sphinxupquote{idxJacRow}} / \sphinxcode{\sphinxupquote{idxJacCol}},
and \sphinxcode{\sphinxupquote{idxHessRow}} / \sphinxcode{\sphinxupquote{idxHessCol}}).

\item {} 
\sphinxAtStartPar
If the value is \sphinxcode{\sphinxupquote{COPT.DENSETYPE\_ROWMAJOR}} or
\sphinxcode{\sphinxupquote{COPT.DENSETYPE\_COLMAJOR}}, the corresponding derivative matrix is
treated as dense and stored in row\sphinxhyphen{}major or column\sphinxhyphen{}major order,
respectively.
In this case, the associated index arrays may be set to \sphinxcode{\sphinxupquote{None}}.

\end{itemize}

\sphinxAtStartPar
The chosen sparse or dense structure must be consistent with the
ordering of numerical values returned by the callback methods.

\end{enumerate}
\end{sphinxadmonition}

\subsubsection{Model.getLpSolution()}
\label{\detokenize{pyapiref:model-getlpsolution}}\begin{quote}

\sphinxAtStartPar
\sphinxstylestrong{Synopsis}
\begin{quote}

\sphinxAtStartPar
\sphinxcode{\sphinxupquote{getLpSolution()}}
\end{quote}

\sphinxAtStartPar
\sphinxstylestrong{Description}
\begin{quote}

\sphinxAtStartPar
Retrieve the values of variables, slack variables, dual variables and reduced cost of variables.
Return a quad tuple object, in which each element is a list.
\end{quote}

\sphinxAtStartPar
\sphinxstylestrong{Example}
\end{quote}

\begin{sphinxVerbatim}[commandchars=\\\{\}]
\PYG{c+c1}{\PYGZsh{} Retrieve solutions of linear model.}
\PYG{n}{values}\PYG{p}{,} \PYG{n}{slacks}\PYG{p}{,} \PYG{n}{duals}\PYG{p}{,} \PYG{n}{redcosts} \PYG{o}{=} \PYG{n}{m}\PYG{o}{.}\PYG{n}{getLpSolution}\PYG{p}{(}\PYG{p}{)}
\end{sphinxVerbatim}

\subsubsection{Model.setLpSolution()}
\label{\detokenize{pyapiref:model-setlpsolution}}\begin{quote}

\sphinxAtStartPar
\sphinxstylestrong{Synopsis}
\begin{quote}

\sphinxAtStartPar
\sphinxcode{\sphinxupquote{setLpSolution(values, slack, duals, redcost)}}
\end{quote}

\sphinxAtStartPar
\sphinxstylestrong{Description}
\begin{quote}

\sphinxAtStartPar
Set LP solution.
\end{quote}

\sphinxAtStartPar
\sphinxstylestrong{Arguments}
\begin{quote}

\sphinxAtStartPar
\sphinxcode{\sphinxupquote{values}}
\begin{quote}

\sphinxAtStartPar
Solution of variables.
\end{quote}

\sphinxAtStartPar
\sphinxcode{\sphinxupquote{slack}}
\begin{quote}

\sphinxAtStartPar
Solution of slack variables.
\end{quote}

\sphinxAtStartPar
\sphinxcode{\sphinxupquote{duals}}
\begin{quote}

\sphinxAtStartPar
Solution of dual variables.
\end{quote}

\sphinxAtStartPar
\sphinxcode{\sphinxupquote{redcost}}
\begin{quote}

\sphinxAtStartPar
Reduced costs of variables.
\end{quote}
\end{quote}

\sphinxAtStartPar
\sphinxstylestrong{Example}
\end{quote}

\begin{sphinxVerbatim}[commandchars=\\\{\}]
\PYG{c+c1}{\PYGZsh{} Set LP solution}
\PYG{n}{m}\PYG{o}{.}\PYG{n}{setLpSolution}\PYG{p}{(}\PYG{n}{values}\PYG{p}{,} \PYG{n}{slack}\PYG{p}{,} \PYG{n}{duals}\PYG{p}{,} \PYG{n}{redcost}\PYG{p}{)}
\end{sphinxVerbatim}

\subsubsection{Model.getValues()}
\label{\detokenize{pyapiref:model-getvalues}}\begin{quote}

\sphinxAtStartPar
\sphinxstylestrong{Synopsis}
\begin{quote}

\sphinxAtStartPar
\sphinxcode{\sphinxupquote{getValues()}}
\end{quote}

\sphinxAtStartPar
\sphinxstylestrong{Description}
\begin{quote}

\sphinxAtStartPar
Retrieve solution values of all variables in a LP or MIP. Return a Python list.
\end{quote}

\sphinxAtStartPar
\sphinxstylestrong{Example}
\end{quote}

\begin{sphinxVerbatim}[commandchars=\\\{\}]
\PYG{n}{values} \PYG{o}{=} \PYG{n}{m}\PYG{o}{.}\PYG{n}{getValues}\PYG{p}{(}\PYG{p}{)}
\end{sphinxVerbatim}

\subsubsection{Model.getRedcosts()}
\label{\detokenize{pyapiref:model-getredcosts}}\begin{quote}

\sphinxAtStartPar
\sphinxstylestrong{Synopsis}
\begin{quote}

\sphinxAtStartPar
\sphinxcode{\sphinxupquote{getRedcosts()}}
\end{quote}

\sphinxAtStartPar
\sphinxstylestrong{Description}
\begin{quote}

\sphinxAtStartPar
Retrieve reduced costs of all variables in a LP. Return a list.
\end{quote}

\sphinxAtStartPar
\sphinxstylestrong{Example}
\end{quote}

\begin{sphinxVerbatim}[commandchars=\\\{\}]
\PYG{c+c1}{\PYGZsh{} Retrieve reduced cost of all variables in model.}
\PYG{n}{redcosts} \PYG{o}{=} \PYG{n}{m}\PYG{o}{.}\PYG{n}{getRedcosts}\PYG{p}{(}\PYG{p}{)}
\end{sphinxVerbatim}

\subsubsection{Model.getSlacks()}
\label{\detokenize{pyapiref:model-getslacks}}\begin{quote}

\sphinxAtStartPar
\sphinxstylestrong{Synopsis}
\begin{quote}

\sphinxAtStartPar
\sphinxcode{\sphinxupquote{getSlacks()}}
\end{quote}

\sphinxAtStartPar
\sphinxstylestrong{Description}
\begin{quote}

\sphinxAtStartPar
Retrieve values of all slack variables in a LP. Return a Python list.
\end{quote}

\sphinxAtStartPar
\sphinxstylestrong{Example}
\end{quote}

\begin{sphinxVerbatim}[commandchars=\\\{\}]
\PYG{c+c1}{\PYGZsh{} Retrieve value of all slack variables in model.}
 \PYG{n}{slacks} \PYG{o}{=} \PYG{n}{m}\PYG{o}{.}\PYG{n}{getSlacks}\PYG{p}{(}\PYG{p}{)}
\end{sphinxVerbatim}

\subsubsection{Model.getDuals()}
\label{\detokenize{pyapiref:model-getduals}}\begin{quote}

\sphinxAtStartPar
\sphinxstylestrong{Synopsis}
\begin{quote}

\sphinxAtStartPar
\sphinxcode{\sphinxupquote{getDuals()}}
\end{quote}

\sphinxAtStartPar
\sphinxstylestrong{Description}
\begin{quote}

\sphinxAtStartPar
Obtain values of all dual variables in a LP. Return a Python list.
\end{quote}

\sphinxAtStartPar
\sphinxstylestrong{Example}
\end{quote}

\begin{sphinxVerbatim}[commandchars=\\\{\}]
\PYG{c+c1}{\PYGZsh{} Retrieve value of all dual variables in model.}
\PYG{n}{duals} \PYG{o}{=} \PYG{n}{m}\PYG{o}{.}\PYG{n}{getDuals}\PYG{p}{(}\PYG{p}{)}
\end{sphinxVerbatim}

\subsubsection{Model.getVarBasis()}
\label{\detokenize{pyapiref:model-getvarbasis}}\begin{quote}

\sphinxAtStartPar
\sphinxstylestrong{Synopsis}
\begin{quote}

\sphinxAtStartPar
\sphinxcode{\sphinxupquote{getVarBasis(vars=None)}}
\end{quote}

\sphinxAtStartPar
\sphinxstylestrong{Description}
\begin{quote}

\sphinxAtStartPar
Obtain basis status of specified variables.

\sphinxAtStartPar
If parameter \sphinxcode{\sphinxupquote{vars}} is \sphinxcode{\sphinxupquote{None}}, then return a list object consistinf of all variables’ basis status.
If parameter \sphinxcode{\sphinxupquote{vars}} is {\hyperref[\detokenize{pyapiref:chappyapi-var}]{\sphinxcrossref{\DUrole{std,std-ref}{Var Class}}}} object, then return basis status of the specified variable.
If parameter \sphinxcode{\sphinxupquote{vars}} is list or {\hyperref[\detokenize{pyapiref:chappyapi-vararray}]{\sphinxcrossref{\DUrole{std,std-ref}{VarArray Class}}}} object, then return a list object consisting of
the specified variables’ basis status.
If parameter \sphinxcode{\sphinxupquote{vars}} is dictionary or {\hyperref[\detokenize{pyapiref:chappyapi-util-tupledict}]{\sphinxcrossref{\DUrole{std,std-ref}{tupledict Class}}}} object, then return indice of the
specified variable as key and {\hyperref[\detokenize{pyapiref:chappyapi-util-tupledict}]{\sphinxcrossref{\DUrole{std,std-ref}{tupledict Class}}}} object consisting of the specified
variables’ basis status as value.
\end{quote}

\sphinxAtStartPar
\sphinxstylestrong{Arguments}
\begin{quote}

\sphinxAtStartPar
\sphinxcode{\sphinxupquote{vars}}
\begin{quote}

\sphinxAtStartPar
The specified variables. Optional, \sphinxcode{\sphinxupquote{None}} by default,
\end{quote}
\end{quote}

\sphinxAtStartPar
\sphinxstylestrong{Example}
\end{quote}

\begin{sphinxVerbatim}[commandchars=\\\{\}]
\PYG{c+c1}{\PYGZsh{} Retrieve all variables\PYGZsq{} basis status in model.}
\PYG{n}{varbasis} \PYG{o}{=} \PYG{n}{m}\PYG{o}{.}\PYG{n}{getVarBasis}\PYG{p}{(}\PYG{p}{)}
\PYG{c+c1}{\PYGZsh{} Retrieve basis status of variable x and y.}
\PYG{n}{varbasis} \PYG{o}{=} \PYG{n}{m}\PYG{o}{.}\PYG{n}{getVarBasis}\PYG{p}{(}\PYG{p}{[}\PYG{n}{x}\PYG{p}{,} \PYG{n}{y}\PYG{p}{]}\PYG{p}{)}
\PYG{c+c1}{\PYGZsh{} Retrieve basis status of tupledict object xx.}
\PYG{n}{varbasis} \PYG{o}{=} \PYG{n}{m}\PYG{o}{.}\PYG{n}{getVarBasis}\PYG{p}{(}\PYG{n}{xx}\PYG{p}{)}
\end{sphinxVerbatim}

\subsubsection{Model.getConstrBasis()}
\label{\detokenize{pyapiref:model-getconstrbasis}}\begin{quote}

\sphinxAtStartPar
\sphinxstylestrong{Synopsis}
\begin{quote}

\sphinxAtStartPar
\sphinxcode{\sphinxupquote{getConstrBasis(constrs=None)}}
\end{quote}

\sphinxAtStartPar
\sphinxstylestrong{Description}
\begin{quote}

\sphinxAtStartPar
Obtain the basis status of linear constraints in LP.

\sphinxAtStartPar
If parameter \sphinxcode{\sphinxupquote{constrs}} is \sphinxcode{\sphinxupquote{None}}, then return a list object consisting of all linear constraints’ basis status.
If parameter \sphinxcode{\sphinxupquote{constrs}} is {\hyperref[\detokenize{pyapiref:chappyapi-constraint}]{\sphinxcrossref{\DUrole{std,std-ref}{Constraint Class}}}} object, then return basis status of the specified linear constraint.
If parameter \sphinxcode{\sphinxupquote{constrs}} is list or {\hyperref[\detokenize{pyapiref:chappyapi-constrarray}]{\sphinxcrossref{\DUrole{std,std-ref}{ConstrArray Class}}}} object, then return a list object consisting of
the specified linear constraints’ basis status.
If parameter \sphinxcode{\sphinxupquote{constrs}} is dictionary or {\hyperref[\detokenize{pyapiref:chappyapi-util-tupledict}]{\sphinxcrossref{\DUrole{std,std-ref}{tupledict Class}}}} object, then return the indice of the
specified linear constraint as key and return {\hyperref[\detokenize{pyapiref:chappyapi-util-tupledict}]{\sphinxcrossref{\DUrole{std,std-ref}{tupledict Class}}}} object consisting of the specified
linear constraints’ basis status as value.
\end{quote}

\sphinxAtStartPar
\sphinxstylestrong{Arguments}
\begin{quote}

\sphinxAtStartPar
\sphinxcode{\sphinxupquote{constrs}}
\begin{quote}

\sphinxAtStartPar
The specified linear constraint. Optional, \sphinxcode{\sphinxupquote{None}} by default.
\end{quote}
\end{quote}

\sphinxAtStartPar
\sphinxstylestrong{Example}
\end{quote}

\begin{sphinxVerbatim}[commandchars=\\\{\}]
\PYG{c+c1}{\PYGZsh{} Retrieve all linear constraints\PYGZsq{} basis status in model.}
\PYG{n}{conbasis} \PYG{o}{=} \PYG{n}{m}\PYG{o}{.}\PYG{n}{getConstrBasis}\PYG{p}{(}\PYG{p}{)}
\PYG{c+c1}{\PYGZsh{} Retrieve basis status corresponding to linear constraint r0 and r1 in model.}
\PYG{n}{conbasis} \PYG{o}{=} \PYG{n}{m}\PYG{o}{.}\PYG{n}{getConstrBasis}\PYG{p}{(}\PYG{p}{[}\PYG{n}{r0}\PYG{p}{,} \PYG{n}{r1}\PYG{p}{]}\PYG{p}{)}
\PYG{c+c1}{\PYGZsh{} Retrieve basis status of linear constraints in tupledict rr.}
\PYG{n}{conbasis} \PYG{o}{=} \PYG{n}{m}\PYG{o}{.}\PYG{n}{getConstrBasis}\PYG{p}{(}\PYG{n}{rr}\PYG{p}{)}
\end{sphinxVerbatim}

\subsubsection{Model.getPoolObjVal()}
\label{\detokenize{pyapiref:model-getpoolobjval}}\begin{quote}

\sphinxAtStartPar
\sphinxstylestrong{Synopsis}
\begin{quote}

\sphinxAtStartPar
\sphinxcode{\sphinxupquote{getPoolObjVal(isol)}}
\end{quote}

\sphinxAtStartPar
\sphinxstylestrong{Description}
\begin{quote}

\sphinxAtStartPar
Obtain the \sphinxcode{\sphinxupquote{isol}} \sphinxhyphen{}th objective value in solution pool, return a constant.
\end{quote}

\sphinxAtStartPar
\sphinxstylestrong{Arguments}
\begin{quote}

\sphinxAtStartPar
\sphinxcode{\sphinxupquote{isol}}
\begin{quote}

\sphinxAtStartPar
Index of solution.
\end{quote}
\end{quote}

\sphinxAtStartPar
\sphinxstylestrong{Example}
\end{quote}

\begin{sphinxVerbatim}[commandchars=\\\{\}]
\PYG{c+c1}{\PYGZsh{} Obtain the second objective value}
\PYG{n}{objval} \PYG{o}{=} \PYG{n}{m}\PYG{o}{.}\PYG{n}{getPoolObjVal}\PYG{p}{(}\PYG{l+m+mi}{2}\PYG{p}{)}
\end{sphinxVerbatim}

\subsubsection{Model.getPoolSolution()}
\label{\detokenize{pyapiref:model-getpoolsolution}}\begin{quote}

\sphinxAtStartPar
\sphinxstylestrong{Synopsis}
\begin{quote}

\sphinxAtStartPar
\sphinxcode{\sphinxupquote{getPoolSolution(isol, vars)}}
\end{quote}

\sphinxAtStartPar
\sphinxstylestrong{Description}
\begin{quote}

\sphinxAtStartPar
Obtain variable values in the \sphinxcode{\sphinxupquote{isol}} \sphinxhyphen{}th solution of solution pool.

\sphinxAtStartPar
If parameter \sphinxcode{\sphinxupquote{vars}} is {\hyperref[\detokenize{pyapiref:chappyapi-var}]{\sphinxcrossref{\DUrole{std,std-ref}{Var Class}}}} object, then return values of the specified variable.
If parameter \sphinxcode{\sphinxupquote{vars}} is list or {\hyperref[\detokenize{pyapiref:chappyapi-vararray}]{\sphinxcrossref{\DUrole{std,std-ref}{VarArray Class}}}} object, then return a list object consisting of
the specified variables’ values.
If parameter \sphinxcode{\sphinxupquote{vars}} is dictionary or {\hyperref[\detokenize{pyapiref:chappyapi-util-tupledict}]{\sphinxcrossref{\DUrole{std,std-ref}{tupledict Class}}}} object, then return indice of the
specified variable as key and {\hyperref[\detokenize{pyapiref:chappyapi-util-tupledict}]{\sphinxcrossref{\DUrole{std,std-ref}{tupledict Class}}}} object consisting of the specified
variables’ values as value.
\end{quote}

\sphinxAtStartPar
\sphinxstylestrong{Arguments}
\begin{quote}

\sphinxAtStartPar
\sphinxcode{\sphinxupquote{isol}}
\begin{quote}

\sphinxAtStartPar
Index of solution
\end{quote}

\sphinxAtStartPar
\sphinxcode{\sphinxupquote{vars}}
\begin{quote}

\sphinxAtStartPar
The specified variables.
\end{quote}
\end{quote}

\sphinxAtStartPar
\sphinxstylestrong{Example}
\end{quote}

\begin{sphinxVerbatim}[commandchars=\\\{\}]
\PYG{c+c1}{\PYGZsh{} Get value of x in the second solution}
\PYG{n}{xval} \PYG{o}{=} \PYG{n}{m}\PYG{o}{.}\PYG{n}{getPoolSolution}\PYG{p}{(}\PYG{l+m+mi}{2}\PYG{p}{,} \PYG{n}{x}\PYG{p}{)}
\end{sphinxVerbatim}

\subsubsection{Model.getPoolObjValN()}
\label{\detokenize{pyapiref:model-getpoolobjvaln}}\begin{quote}

\sphinxAtStartPar
\sphinxstylestrong{Synopsis}
\begin{quote}

\sphinxAtStartPar
\sphinxcode{\sphinxupquote{getPoolObjValN(idx, isol)}}
\end{quote}

\sphinxAtStartPar
\sphinxstylestrong{Description}
\begin{quote}

\sphinxAtStartPar
In a multi\sphinxhyphen{}objective model,
retrieve the objective value of the specified objective
at the given solution index in the solution pool.
\end{quote}

\sphinxAtStartPar
\sphinxstylestrong{Arguments}
\begin{quote}

\sphinxAtStartPar
\sphinxcode{\sphinxupquote{idx}}
\begin{quote}

\sphinxAtStartPar
Index of the target objective in the multi\sphinxhyphen{}objective model.
\end{quote}

\sphinxAtStartPar
\sphinxcode{\sphinxupquote{isol}}
\begin{quote}

\sphinxAtStartPar
Index of the solution in the solution pool.
\end{quote}
\end{quote}

\sphinxAtStartPar
\sphinxstylestrong{Return Value}
\begin{quote}

\sphinxAtStartPar
Returns a double value.
\end{quote}
\end{quote}

\subsubsection{Model.getVarLowerIIS()}
\label{\detokenize{pyapiref:model-getvarloweriis}}\begin{quote}

\sphinxAtStartPar
\sphinxstylestrong{Synopsis}
\begin{quote}

\sphinxAtStartPar
\sphinxcode{\sphinxupquote{getVarLowerIIS(vars)}}
\end{quote}

\sphinxAtStartPar
\sphinxstylestrong{Description}
\begin{quote}

\sphinxAtStartPar
Obtain IIS status of lower bounds of variables.

\sphinxAtStartPar
If parameter \sphinxcode{\sphinxupquote{vars}} is {\hyperref[\detokenize{pyapiref:chappyapi-var}]{\sphinxcrossref{\DUrole{std,std-ref}{Var Class}}}} object, then return IIS status of lower bound of variable.
If parameter \sphinxcode{\sphinxupquote{vars}} is list or {\hyperref[\detokenize{pyapiref:chappyapi-vararray}]{\sphinxcrossref{\DUrole{std,std-ref}{VarArray Class}}}} object, then return a list object consisting of
the IIS status of lower bounds of variables.
If parameter \sphinxcode{\sphinxupquote{vars}} is dictionary or {\hyperref[\detokenize{pyapiref:chappyapi-util-tupledict}]{\sphinxcrossref{\DUrole{std,std-ref}{tupledict Class}}}} object, then return indice of the
specified variable as key and {\hyperref[\detokenize{pyapiref:chappyapi-util-tupledict}]{\sphinxcrossref{\DUrole{std,std-ref}{tupledict Class}}}} object consisting of the IIS status of lower
bounds of variables as value.
\end{quote}

\sphinxAtStartPar
\sphinxstylestrong{Arguments}
\begin{quote}

\sphinxAtStartPar
\sphinxcode{\sphinxupquote{vars}}
\begin{quote}

\sphinxAtStartPar
The specified variables.
\end{quote}
\end{quote}

\sphinxAtStartPar
\sphinxstylestrong{Example}
\end{quote}

\begin{sphinxVerbatim}[commandchars=\\\{\}]
\PYG{c+c1}{\PYGZsh{} Retrieve IIS status of lower bounds of variable x and y.}
\PYG{n}{lowerIIS} \PYG{o}{=} \PYG{n}{m}\PYG{o}{.}\PYG{n}{getVarLowerIIS}\PYG{p}{(}\PYG{p}{[}\PYG{n}{x}\PYG{p}{,} \PYG{n}{y}\PYG{p}{]}\PYG{p}{)}
\PYG{c+c1}{\PYGZsh{} Retrieve IIS status of lower bounds of variables in tupledict object xx.}
\PYG{n}{lowerIIS} \PYG{o}{=} \PYG{n}{m}\PYG{o}{.}\PYG{n}{getVarLowerIIS}\PYG{p}{(}\PYG{n}{xx}\PYG{p}{)}
\end{sphinxVerbatim}

\subsubsection{Model.getVarUpperIIS()}
\label{\detokenize{pyapiref:model-getvarupperiis}}\begin{quote}

\sphinxAtStartPar
\sphinxstylestrong{Synopsis}
\begin{quote}

\sphinxAtStartPar
\sphinxcode{\sphinxupquote{getVarUpperIIS(vars)}}
\end{quote}

\sphinxAtStartPar
\sphinxstylestrong{Description}
\begin{quote}

\sphinxAtStartPar
Obtain IIS status of upper bounds of variables.

\sphinxAtStartPar
If parameter \sphinxcode{\sphinxupquote{vars}} is {\hyperref[\detokenize{pyapiref:chappyapi-var}]{\sphinxcrossref{\DUrole{std,std-ref}{Var Class}}}} object, then return IIS status of upper bound of variable.
If parameter \sphinxcode{\sphinxupquote{vars}} is list or {\hyperref[\detokenize{pyapiref:chappyapi-vararray}]{\sphinxcrossref{\DUrole{std,std-ref}{VarArray Class}}}} object, then return a list object consisting of
the IIS status of upper bounds of variables.
If parameter \sphinxcode{\sphinxupquote{vars}} is dictionary or {\hyperref[\detokenize{pyapiref:chappyapi-util-tupledict}]{\sphinxcrossref{\DUrole{std,std-ref}{tupledict Class}}}} object, then return indice of the
specified variable as key and {\hyperref[\detokenize{pyapiref:chappyapi-util-tupledict}]{\sphinxcrossref{\DUrole{std,std-ref}{tupledict Class}}}} object consisting of the IIS status of upper
bounds of variables as value.
\end{quote}

\sphinxAtStartPar
\sphinxstylestrong{Arguments}
\begin{quote}

\sphinxAtStartPar
\sphinxcode{\sphinxupquote{vars}}
\begin{quote}

\sphinxAtStartPar
The specified variables.
\end{quote}
\end{quote}

\sphinxAtStartPar
\sphinxstylestrong{Example}
\end{quote}

\begin{sphinxVerbatim}[commandchars=\\\{\}]
\PYG{c+c1}{\PYGZsh{} Retrieve IIS status of upper bounds of variable x and y.}
\PYG{n}{upperIIS} \PYG{o}{=} \PYG{n}{m}\PYG{o}{.}\PYG{n}{getVarUpperIIS}\PYG{p}{(}\PYG{p}{[}\PYG{n}{x}\PYG{p}{,} \PYG{n}{y}\PYG{p}{]}\PYG{p}{)}
\PYG{c+c1}{\PYGZsh{} Retrieve IIS status of upper bounds of variables in tupledict object xx.}
\PYG{n}{upperIIS} \PYG{o}{=} \PYG{n}{m}\PYG{o}{.}\PYG{n}{getVarUpperIIS}\PYG{p}{(}\PYG{n}{xx}\PYG{p}{)}
\end{sphinxVerbatim}

\subsubsection{Model.getConstrLowerIIS()}
\label{\detokenize{pyapiref:model-getconstrloweriis}}\begin{quote}

\sphinxAtStartPar
\sphinxstylestrong{Synopsis}
\begin{quote}

\sphinxAtStartPar
\sphinxcode{\sphinxupquote{getConstrLowerIIS(constrs)}}
\end{quote}

\sphinxAtStartPar
\sphinxstylestrong{Description}
\begin{quote}

\sphinxAtStartPar
Obtain the IIS status of lower bounds of constraints.

\sphinxAtStartPar
If parameter \sphinxcode{\sphinxupquote{constrs}} is {\hyperref[\detokenize{pyapiref:chappyapi-constraint}]{\sphinxcrossref{\DUrole{std,std-ref}{Constraint Class}}}} object, then return IIS status of lower bound of constraint.
If parameter \sphinxcode{\sphinxupquote{constrs}} is list or {\hyperref[\detokenize{pyapiref:chappyapi-constrarray}]{\sphinxcrossref{\DUrole{std,std-ref}{ConstrArray Class}}}} object, then return a list object consisting of
the IIS status of lower bounds of constraints.
If parameter \sphinxcode{\sphinxupquote{constrs}} is dictionary or {\hyperref[\detokenize{pyapiref:chappyapi-util-tupledict}]{\sphinxcrossref{\DUrole{std,std-ref}{tupledict Class}}}} object, then return the indice of the
specified linear constraint as key and return {\hyperref[\detokenize{pyapiref:chappyapi-util-tupledict}]{\sphinxcrossref{\DUrole{std,std-ref}{tupledict Class}}}} object consisting of the IIS status
of lower bounds of constraints.
\end{quote}

\sphinxAtStartPar
\sphinxstylestrong{Arguments}
\begin{quote}

\sphinxAtStartPar
\sphinxcode{\sphinxupquote{constrs}}
\begin{quote}

\sphinxAtStartPar
The specified linear constraint.
\end{quote}
\end{quote}

\sphinxAtStartPar
\sphinxstylestrong{Example}
\end{quote}

\begin{sphinxVerbatim}[commandchars=\\\{\}]
\PYG{c+c1}{\PYGZsh{} Retrieve IIS status corresponding to lower bounds of linear constraint r0 and r1 in model.}
\PYG{n}{lowerIIS} \PYG{o}{=} \PYG{n}{m}\PYG{o}{.}\PYG{n}{getConstrLowerIIS}\PYG{p}{(}\PYG{p}{[}\PYG{n}{r0}\PYG{p}{,} \PYG{n}{r1}\PYG{p}{]}\PYG{p}{)}
\PYG{c+c1}{\PYGZsh{} Retrieve IIS status of lower bounds of linear constraints in tupledict rr.}
\PYG{n}{lowerIIS} \PYG{o}{=} \PYG{n}{m}\PYG{o}{.}\PYG{n}{getConstrLowerIIS}\PYG{p}{(}\PYG{n}{rr}\PYG{p}{)}
\end{sphinxVerbatim}

\subsubsection{Model.getConstrUpperIIS()}
\label{\detokenize{pyapiref:model-getconstrupperiis}}\begin{quote}

\sphinxAtStartPar
\sphinxstylestrong{Synopsis}
\begin{quote}

\sphinxAtStartPar
\sphinxcode{\sphinxupquote{getConstrUpperIIS(constrs)}}
\end{quote}

\sphinxAtStartPar
\sphinxstylestrong{Description}
\begin{quote}

\sphinxAtStartPar
Obtain the IIS status of upper bounds of constraints.

\sphinxAtStartPar
If parameter \sphinxcode{\sphinxupquote{constrs}} is {\hyperref[\detokenize{pyapiref:chappyapi-constraint}]{\sphinxcrossref{\DUrole{std,std-ref}{Constraint Class}}}} object, then return IIS status of upper bound of constraint.
If parameter \sphinxcode{\sphinxupquote{constrs}} is list or {\hyperref[\detokenize{pyapiref:chappyapi-constrarray}]{\sphinxcrossref{\DUrole{std,std-ref}{ConstrArray Class}}}} object, then return a list object consisting of
the IIS status of upper bounds of constraints.
If parameter \sphinxcode{\sphinxupquote{constrs}} is dictionary or {\hyperref[\detokenize{pyapiref:chappyapi-util-tupledict}]{\sphinxcrossref{\DUrole{std,std-ref}{tupledict Class}}}} object, then return the indice of the
specified linear constraint as key and return {\hyperref[\detokenize{pyapiref:chappyapi-util-tupledict}]{\sphinxcrossref{\DUrole{std,std-ref}{tupledict Class}}}} object consisting of the IIS status
of upper bounds of constraints.
\end{quote}

\sphinxAtStartPar
\sphinxstylestrong{Arguments}
\begin{quote}

\sphinxAtStartPar
\sphinxcode{\sphinxupquote{constrs}}
\begin{quote}

\sphinxAtStartPar
The specified linear constraint.
\end{quote}
\end{quote}

\sphinxAtStartPar
\sphinxstylestrong{Example}
\end{quote}

\begin{sphinxVerbatim}[commandchars=\\\{\}]
\PYG{c+c1}{\PYGZsh{} Retrieve IIS status corresponding to upper bounds of linear constraint r0 and r1 in model.}
\PYG{n}{upperIIS} \PYG{o}{=} \PYG{n}{m}\PYG{o}{.}\PYG{n}{getConstrUpperIIS}\PYG{p}{(}\PYG{p}{[}\PYG{n}{r0}\PYG{p}{,} \PYG{n}{r1}\PYG{p}{]}\PYG{p}{)}
\PYG{c+c1}{\PYGZsh{} Retrieve IIS status of upper bounds of linear constraints in tupledict rr.}
\PYG{n}{upperIIS} \PYG{o}{=} \PYG{n}{m}\PYG{o}{.}\PYG{n}{getConstrUpperIIS}\PYG{p}{(}\PYG{n}{rr}\PYG{p}{)}
\end{sphinxVerbatim}

\subsubsection{Model.getSOSIIS()}
\label{\detokenize{pyapiref:model-getsosiis}}\begin{quote}

\sphinxAtStartPar
\sphinxstylestrong{Synopsis}
\begin{quote}

\sphinxAtStartPar
\sphinxcode{\sphinxupquote{getSOSIIS(soss)}}
\end{quote}

\sphinxAtStartPar
\sphinxstylestrong{Description}
\begin{quote}

\sphinxAtStartPar
Obtain the IIS status of SOS constraints.

\sphinxAtStartPar
If parameter \sphinxcode{\sphinxupquote{soss}} is {\hyperref[\detokenize{pyapiref:chappyapi-sos}]{\sphinxcrossref{\DUrole{std,std-ref}{SOS Class}}}} object, then return IIS status of SOS constraint.
If parameter \sphinxcode{\sphinxupquote{soss}} is list or {\hyperref[\detokenize{pyapiref:chappyapi-sosarray}]{\sphinxcrossref{\DUrole{std,std-ref}{SOSArray Class}}}} object, then return a list object consisting of
the IIS status of SOS constraints.
If parameter \sphinxcode{\sphinxupquote{soss}} is dictionary or {\hyperref[\detokenize{pyapiref:chappyapi-util-tupledict}]{\sphinxcrossref{\DUrole{std,std-ref}{tupledict Class}}}} object, then return the indice of the
specified SOS constraint as key and return {\hyperref[\detokenize{pyapiref:chappyapi-util-tupledict}]{\sphinxcrossref{\DUrole{std,std-ref}{tupledict Class}}}} object consisting of the IIS status
of SOS constraints.
\end{quote}

\sphinxAtStartPar
\sphinxstylestrong{Arguments}
\begin{quote}

\sphinxAtStartPar
\sphinxcode{\sphinxupquote{soss}}
\begin{quote}

\sphinxAtStartPar
The specified SOS constraint.
\end{quote}
\end{quote}

\sphinxAtStartPar
\sphinxstylestrong{Example}
\end{quote}

\begin{sphinxVerbatim}[commandchars=\\\{\}]
\PYG{c+c1}{\PYGZsh{} Retrieve IIS status corresponding to SOS constraint r0 and r1 in model.}
\PYG{n}{sosIIS} \PYG{o}{=} \PYG{n}{m}\PYG{o}{.}\PYG{n}{getSOSIIS}\PYG{p}{(}\PYG{p}{[}\PYG{n}{r0}\PYG{p}{,} \PYG{n}{r1}\PYG{p}{]}\PYG{p}{)}
\PYG{c+c1}{\PYGZsh{} Retrieve IIS status of SOS constraints in tupledict rr.}
\PYG{n}{sosIIS} \PYG{o}{=} \PYG{n}{m}\PYG{o}{.}\PYG{n}{getSOSIIS}\PYG{p}{(}\PYG{n}{rr}\PYG{p}{)}
\end{sphinxVerbatim}

\subsubsection{Model.getIndicatorIIS()}
\label{\detokenize{pyapiref:model-getindicatoriis}}\begin{quote}

\sphinxAtStartPar
\sphinxstylestrong{Synopsis}
\begin{quote}

\sphinxAtStartPar
\sphinxcode{\sphinxupquote{getIndicatorIIS(genconstrs)}}
\end{quote}

\sphinxAtStartPar
\sphinxstylestrong{Description}
\begin{quote}

\sphinxAtStartPar
Obtain the IIS status of indicator constraints.

\sphinxAtStartPar
If parameter \sphinxcode{\sphinxupquote{genconstrs}} is {\hyperref[\detokenize{pyapiref:chappyapi-genconstr}]{\sphinxcrossref{\DUrole{std,std-ref}{GenConstr Class}}}} object, then return IIS status of indicator constraint.
If parameter \sphinxcode{\sphinxupquote{genconstrs}} is list or {\hyperref[\detokenize{pyapiref:chappyapi-genconstrarray}]{\sphinxcrossref{\DUrole{std,std-ref}{GenConstrArray Class}}}} object, then return a list object consisting of
the IIS status of indicator constraints.
If parameter \sphinxcode{\sphinxupquote{genconstrs}} is dictionary or {\hyperref[\detokenize{pyapiref:chappyapi-util-tupledict}]{\sphinxcrossref{\DUrole{std,std-ref}{tupledict Class}}}} object, then return the indice of the
specified indicator constraint as key and return {\hyperref[\detokenize{pyapiref:chappyapi-util-tupledict}]{\sphinxcrossref{\DUrole{std,std-ref}{tupledict Class}}}} object consisting of the IIS status
of indicator constraints.
\end{quote}

\sphinxAtStartPar
\sphinxstylestrong{Arguments}
\begin{quote}

\sphinxAtStartPar
\sphinxcode{\sphinxupquote{genconstrs}}
\begin{quote}

\sphinxAtStartPar
The specified indicator constraint.
\end{quote}
\end{quote}

\sphinxAtStartPar
\sphinxstylestrong{Example}
\end{quote}

\begin{sphinxVerbatim}[commandchars=\\\{\}]
\PYG{c+c1}{\PYGZsh{} Retrieve IIS status corresponding to indicator constraint r0 and r1 in model.}
\PYG{n}{indicatorIIS} \PYG{o}{=} \PYG{n}{m}\PYG{o}{.}\PYG{n}{getIndicatorIIS}\PYG{p}{(}\PYG{p}{[}\PYG{n}{r0}\PYG{p}{,} \PYG{n}{r1}\PYG{p}{]}\PYG{p}{)}
\PYG{c+c1}{\PYGZsh{} Retrieve IIS status of indicator constraints in tupledict rr.}
\PYG{n}{indicatorIIS} \PYG{o}{=} \PYG{n}{m}\PYG{o}{.}\PYG{n}{getIndicatorIIS}\PYG{p}{(}\PYG{n}{rr}\PYG{p}{)}
\end{sphinxVerbatim}

\subsubsection{Model.getAttr()}
\label{\detokenize{pyapiref:model-getattr}}\begin{quote}

\sphinxAtStartPar
\sphinxstylestrong{Synopsis}
\begin{quote}

\sphinxAtStartPar
\sphinxcode{\sphinxupquote{getAttr(attrname)}}
\end{quote}

\sphinxAtStartPar
\sphinxstylestrong{Description}
\begin{quote}

\sphinxAtStartPar
Get the value of an attribute of model. Return a constant.
\end{quote}

\sphinxAtStartPar
\sphinxstylestrong{Arguments}
\begin{quote}

\sphinxAtStartPar
\sphinxcode{\sphinxupquote{attrname}}
\begin{quote}

\sphinxAtStartPar
The specified attribute name. The full list of available attributes can be found in {\hyperref[\detokenize{attribute:chapattrs}]{\sphinxcrossref{\DUrole{std,std-ref}{Attributes}}}} section.
\end{quote}
\end{quote}

\sphinxAtStartPar
\sphinxstylestrong{Example}
\end{quote}

\begin{sphinxVerbatim}[commandchars=\\\{\}]
\PYG{c+c1}{\PYGZsh{} Retrieve the constant terms of objective.}
\PYG{n}{objconst} \PYG{o}{=} \PYG{n}{m}\PYG{o}{.}\PYG{n}{getAttr}\PYG{p}{(}\PYG{n}{COPT}\PYG{o}{.}\PYG{n}{Attr}\PYG{o}{.}\PYG{n}{ObjConst}\PYG{p}{)}
\end{sphinxVerbatim}

\subsubsection{Model.getInfo()}
\label{\detokenize{pyapiref:model-getinfo}}\begin{quote}

\sphinxAtStartPar
\sphinxstylestrong{Synopsis}
\begin{quote}

\sphinxAtStartPar
\sphinxcode{\sphinxupquote{getInfo(infoname, args)}}
\end{quote}

\sphinxAtStartPar
\sphinxstylestrong{Description}
\begin{quote}

\sphinxAtStartPar
Retrieve specified information.

\sphinxAtStartPar
If parameter \sphinxcode{\sphinxupquote{args}} is {\hyperref[\detokenize{pyapiref:chappyapi-var}]{\sphinxcrossref{\DUrole{std,std-ref}{Var Class}}}} object,
{\hyperref[\detokenize{pyapiref:chappyapi-constraint}]{\sphinxcrossref{\DUrole{std,std-ref}{Constraint Class}}}} object, {\hyperref[\detokenize{pyapiref:chappyapi-qconstraint}]{\sphinxcrossref{\DUrole{std,std-ref}{QConstraint Class}}}} object,
or {\hyperref[\detokenize{pyapiref:chappyapi-nlconstraint}]{\sphinxcrossref{\DUrole{std,std-ref}{NlConstraint Class}}}} object,
then return info of the specified variable or constraint.

\sphinxAtStartPar
If parameter \sphinxcode{\sphinxupquote{args}} is list or {\hyperref[\detokenize{pyapiref:chappyapi-vararray}]{\sphinxcrossref{\DUrole{std,std-ref}{VarArray Class}}}} object,
{\hyperref[\detokenize{pyapiref:chappyapi-constrarray}]{\sphinxcrossref{\DUrole{std,std-ref}{ConstrArray Class}}}} object, iterable {\hyperref[\detokenize{pyapiref:chappyapi-qconstraint}]{\sphinxcrossref{\DUrole{std,std-ref}{QConstraint Class}}}} objects,
or {\hyperref[\detokenize{pyapiref:chappyapi-nlconstraint}]{\sphinxcrossref{\DUrole{std,std-ref}{NlConstraint Class}}}} objects,
then return a list consisting of the specified variables or constraints.

\sphinxAtStartPar
If parameter \sphinxcode{\sphinxupquote{args}} is dictionary or {\hyperref[\detokenize{pyapiref:chappyapi-util-tupledict}]{\sphinxcrossref{\DUrole{std,std-ref}{tupledict Class}}}} object, then return the indice of the
specified variables or constraints as key and return {\hyperref[\detokenize{pyapiref:chappyapi-util-tupledict}]{\sphinxcrossref{\DUrole{std,std-ref}{tupledict Class}}}} object consisting of info
corresponding to the specified variables or constraints as value.

\sphinxAtStartPar
If parameter \sphinxcode{\sphinxupquote{args}} is {\hyperref[\detokenize{pyapiref:chappyapi-mvar}]{\sphinxcrossref{\DUrole{std,std-ref}{MVar Class}}}} object, {\hyperref[\detokenize{pyapiref:chappyapi-mconstr}]{\sphinxcrossref{\DUrole{std,std-ref}{MConstr Class}}}} object,
{\hyperref[\detokenize{pyapiref:chappyapi-mqconstr}]{\sphinxcrossref{\DUrole{std,std-ref}{MQConstr Class}}}} object or {\hyperref[\detokenize{pyapiref:chappyapi-mpsdconstr}]{\sphinxcrossref{\DUrole{std,std-ref}{MPsdConstr Class}}}} object,
then return a \sphinxcode{\sphinxupquote{numpy.ndarray}} consisting of the specified variables or constraints.
\end{quote}

\sphinxAtStartPar
\sphinxstylestrong{Arguments}
\begin{quote}

\sphinxAtStartPar
\sphinxcode{\sphinxupquote{infoname}}
\begin{quote}

\sphinxAtStartPar
The specified information name. The full list of available attributes can be found in {\hyperref[\detokenize{information:chapinfo}]{\sphinxcrossref{\DUrole{std,std-ref}{Information}}}} section.
\end{quote}

\sphinxAtStartPar
\sphinxcode{\sphinxupquote{args}}
\begin{quote}

\sphinxAtStartPar
Variables and constraints to get information.
\end{quote}
\end{quote}

\sphinxAtStartPar
\sphinxstylestrong{Example}
\end{quote}

\begin{sphinxVerbatim}[commandchars=\\\{\}]
\PYG{c+c1}{\PYGZsh{} Retrieve lower bound information of all linear constraints in model.}
\PYG{n}{lb} \PYG{o}{=} \PYG{n}{m}\PYG{o}{.}\PYG{n}{getInfo}\PYG{p}{(}\PYG{n}{COPT}\PYG{o}{.}\PYG{n}{Info}\PYG{o}{.}\PYG{n}{LB}\PYG{p}{,} \PYG{n}{m}\PYG{o}{.}\PYG{n}{getConstrs}\PYG{p}{(}\PYG{p}{)}\PYG{p}{)}
\PYG{c+c1}{\PYGZsh{} Retrieve value information of variables x and y.}
\PYG{n}{sol} \PYG{o}{=} \PYG{n}{m}\PYG{o}{.}\PYG{n}{getInfo}\PYG{p}{(}\PYG{n}{COPT}\PYG{o}{.}\PYG{n}{Info}\PYG{o}{.}\PYG{n}{Value}\PYG{p}{,} \PYG{p}{[}\PYG{n}{x}\PYG{p}{,} \PYG{n}{y}\PYG{p}{]}\PYG{p}{)}
\PYG{c+c1}{\PYGZsh{} Retrieve the dual variable value corresponding to linear constraint in tupledict object shipconstr.}
\PYG{n}{dual} \PYG{o}{=} \PYG{n}{m}\PYG{o}{.}\PYG{n}{getInfo}\PYG{p}{(}\PYG{n}{COPT}\PYG{o}{.}\PYG{n}{Info}\PYG{o}{.}\PYG{n}{Dual}\PYG{p}{,} \PYG{n}{shipconstr}\PYG{p}{)}
\end{sphinxVerbatim}

\subsubsection{Model.getAttrN()}
\label{\detokenize{pyapiref:model-getattrn}}\begin{quote}

\sphinxAtStartPar
\sphinxstylestrong{Synopsis}
\begin{quote}

\sphinxAtStartPar
\sphinxcode{\sphinxupquote{getAttrN(idx, attrname)}}
\end{quote}

\sphinxAtStartPar
\sphinxstylestrong{Description}
\begin{quote}

\sphinxAtStartPar
Get the attribute value of the specified objective in multi\sphinxhyphen{}objective optimization.
\end{quote}

\sphinxAtStartPar
\sphinxstylestrong{Arguments}
\begin{quote}

\sphinxAtStartPar
\sphinxcode{\sphinxupquote{idx}}
\begin{quote}

\sphinxAtStartPar
Index of the objective function.
\end{quote}

\sphinxAtStartPar
\sphinxcode{\sphinxupquote{attrname}}
\begin{quote}

\sphinxAtStartPar
Name of the objective attribute.

\sphinxAtStartPar
Possible values:
{\hyperref[\detokenize{attribute:hasqobj}]{\sphinxcrossref{\DUrole{std,std-ref}{HasQObj}}}}, {\hyperref[\detokenize{attribute:hasnlobj}]{\sphinxcrossref{\DUrole{std,std-ref}{HasNLObj}}}}, {\hyperref[\detokenize{attribute:lpobjval}]{\sphinxcrossref{\DUrole{std,std-ref}{LpObjval}}}},
{\hyperref[\detokenize{information:bestobj}]{\sphinxcrossref{\DUrole{std,std-ref}{BestObj}}}}, {\hyperref[\detokenize{attribute:objsense}]{\sphinxcrossref{\DUrole{std,std-ref}{ObjSense}}}}, {\hyperref[\detokenize{attribute:objconst}]{\sphinxcrossref{\DUrole{std,std-ref}{ObjConst}}}}.
\end{quote}
\end{quote}

\sphinxAtStartPar
\sphinxstylestrong{Return Value}
\begin{quote}

\sphinxAtStartPar
A double or integer value.
\end{quote}
\end{quote}

\subsubsection{Model.getParam()}
\label{\detokenize{pyapiref:model-getparam}}\begin{quote}

\sphinxAtStartPar
\sphinxstylestrong{Synopsis}
\begin{quote}

\sphinxAtStartPar
\sphinxcode{\sphinxupquote{getParam(paramname)}}
\end{quote}

\sphinxAtStartPar
\sphinxstylestrong{Description}
\begin{quote}

\sphinxAtStartPar
Retrive the current value of the specified parameter. Return a constant.
\end{quote}

\sphinxAtStartPar
\sphinxstylestrong{Arguments}
\begin{quote}

\sphinxAtStartPar
\sphinxcode{\sphinxupquote{paramname}}
\begin{quote}

\sphinxAtStartPar
The name of the parameter to get access to. The full list of available attributes can be found in  {\hyperref[\detokenize{pyapiref:chappyapi-const-param}]{\sphinxcrossref{\DUrole{std,std-ref}{Parameters}}}}  section.
\end{quote}
\end{quote}

\sphinxAtStartPar
\sphinxstylestrong{Example}
\end{quote}

\begin{sphinxVerbatim}[commandchars=\\\{\}]
\PYG{c+c1}{\PYGZsh{} Retrieve current value of time limit.}
\PYG{n}{timelimit} \PYG{o}{=} \PYG{n}{m}\PYG{o}{.}\PYG{n}{getParam}\PYG{p}{(}\PYG{n}{COPT}\PYG{o}{.}\PYG{n}{Param}\PYG{o}{.}\PYG{n}{TimeLimit}\PYG{p}{)}
\end{sphinxVerbatim}

\subsubsection{Model.getParamInfo()}
\label{\detokenize{pyapiref:model-getparaminfo}}\begin{quote}

\sphinxAtStartPar
\sphinxstylestrong{Synopsis}
\begin{quote}

\sphinxAtStartPar
\sphinxcode{\sphinxupquote{getParamInfo(paramname)}}
\end{quote}

\sphinxAtStartPar
\sphinxstylestrong{Description}
\begin{quote}

\sphinxAtStartPar
Retrieve information of the specified optimization parameter. Return a tuple object, consisiting of param name,
current value, default value, minimum value, maximum value.
\end{quote}

\sphinxAtStartPar
\sphinxstylestrong{Arguments}
\begin{quote}

\sphinxAtStartPar
\sphinxcode{\sphinxupquote{paramname}}
\begin{quote}

\sphinxAtStartPar
Name of the specified parameter.
The full list of available values can be found in {\hyperref[\detokenize{pyapiref:chappyapi-const-param}]{\sphinxcrossref{\DUrole{std,std-ref}{Parameters}}}} section.
\end{quote}
\end{quote}

\sphinxAtStartPar
\sphinxstylestrong{Example}
\end{quote}

\begin{sphinxVerbatim}[commandchars=\\\{\}]
\PYG{c+c1}{\PYGZsh{} Retrieve information of time limit.}
\PYG{n}{pname}\PYG{p}{,} \PYG{n}{pcur}\PYG{p}{,} \PYG{n}{pdef}\PYG{p}{,} \PYG{n}{pmin}\PYG{p}{,} \PYG{n}{pmax} \PYG{o}{=} \PYG{n}{m}\PYG{o}{.}\PYG{n}{getParamInfo}\PYG{p}{(}\PYG{n}{COPT}\PYG{o}{.}\PYG{n}{Param}\PYG{o}{.}\PYG{n}{TimeLimit}\PYG{p}{)}
\end{sphinxVerbatim}

\subsubsection{Model.getObjParamN()}
\label{\detokenize{pyapiref:model-getobjparamn}}\begin{quote}

\sphinxAtStartPar
\sphinxstylestrong{Synopsis}
\begin{quote}

\sphinxAtStartPar
\sphinxcode{\sphinxupquote{getObjParamN(idx, paramname)}}
\end{quote}

\sphinxAtStartPar
\sphinxstylestrong{Description}
\begin{quote}

\sphinxAtStartPar
Get the value of a parameter for the specified objective
in multi\sphinxhyphen{}objective optimization.
\end{quote}

\sphinxAtStartPar
\sphinxstylestrong{Arguments}
\begin{quote}

\sphinxAtStartPar
\sphinxcode{\sphinxupquote{idx}}
\begin{quote}

\sphinxAtStartPar
Index of the objective function.
\end{quote}

\sphinxAtStartPar
\sphinxcode{\sphinxupquote{paramname}}
\begin{quote}

\sphinxAtStartPar
Name of the objective parameter.

\sphinxAtStartPar
Possible values are {\hyperref[\detokenize{multiobj:multiobjpriority}]{\sphinxcrossref{\DUrole{std,std-ref}{MultiObjPriority}}}},
{\hyperref[\detokenize{multiobj:multiobjweight}]{\sphinxcrossref{\DUrole{std,std-ref}{MultiObjWeight}}}},
{\hyperref[\detokenize{multiobj:multiobjreltol}]{\sphinxcrossref{\DUrole{std,std-ref}{MultiObjRelTol}}}},
{\hyperref[\detokenize{multiobj:multiobjabstol}]{\sphinxcrossref{\DUrole{std,std-ref}{MultiObjAbsTol}}}}.
\end{quote}
\end{quote}

\sphinxAtStartPar
\sphinxstylestrong{Return Value}
\begin{quote}

\sphinxAtStartPar
A double or integer value.
\end{quote}
\end{quote}

\subsubsection{Model.getParamN()}
\label{\detokenize{pyapiref:model-getparamn}}\begin{quote}

\sphinxAtStartPar
\sphinxstylestrong{Synopsis}
\begin{quote}

\sphinxAtStartPar
\sphinxcode{\sphinxupquote{getParamN(idx, paramname)}}
\end{quote}

\sphinxAtStartPar
\sphinxstylestrong{Description}
\begin{quote}

\sphinxAtStartPar
Get the value of a solver parameter for the model associated with
the specified objective in multi\sphinxhyphen{}objective optimization.
\end{quote}

\sphinxAtStartPar
\sphinxstylestrong{Arguments}
\begin{quote}

\sphinxAtStartPar
\sphinxcode{\sphinxupquote{idx}}
\begin{quote}

\sphinxAtStartPar
Index of the objective function.
\end{quote}

\sphinxAtStartPar
\sphinxcode{\sphinxupquote{paramname}}
\begin{quote}

\sphinxAtStartPar
Name of the solver parameter.

\sphinxAtStartPar
See {\hyperref[\detokenize{pyapiref:chappyapi-const-param}]{\sphinxcrossref{\DUrole{std,std-ref}{Parameters}}}} for possible values.
\end{quote}
\end{quote}

\sphinxAtStartPar
\sphinxstylestrong{Return Value}
\begin{quote}

\sphinxAtStartPar
A double or integer value.
\end{quote}
\end{quote}

\subsubsection{Model.setBasis()}
\label{\detokenize{pyapiref:model-setbasis}}\begin{quote}

\sphinxAtStartPar
\sphinxstylestrong{Synopsis}
\begin{quote}

\sphinxAtStartPar
\sphinxcode{\sphinxupquote{setBasis(varbasis, constrbasis)}}
\end{quote}

\sphinxAtStartPar
\sphinxstylestrong{Description}
\begin{quote}

\sphinxAtStartPar
Set basis status for all variables and linear constraints in LP.
The parameters \sphinxcode{\sphinxupquote{varbasis}} and \sphinxcode{\sphinxupquote{constrbasis}} are list objects whose number of elements is the total number of
variables or linear constraints.
\end{quote}

\sphinxAtStartPar
\sphinxstylestrong{Arguments}
\begin{quote}

\sphinxAtStartPar
\sphinxcode{\sphinxupquote{varbasis}}
\begin{quote}

\sphinxAtStartPar
The basis status of variables.
\end{quote}

\sphinxAtStartPar
\sphinxcode{\sphinxupquote{constrbasis}}
\begin{quote}

\sphinxAtStartPar
The basis status of constraints.
\end{quote}
\end{quote}

\sphinxAtStartPar
\sphinxstylestrong{Example}
\end{quote}

\begin{sphinxVerbatim}[commandchars=\\\{\}]
\PYG{c+c1}{\PYGZsh{} Set basis status for all variables and linear constraints in the model.}
\PYG{n}{m}\PYG{o}{.}\PYG{n}{setBasis}\PYG{p}{(}\PYG{n}{varbasis}\PYG{p}{,} \PYG{n}{constrbasis}\PYG{p}{)}
\end{sphinxVerbatim}

\subsubsection{Model.setSlackBasis()}
\label{\detokenize{pyapiref:model-setslackbasis}}\begin{quote}

\sphinxAtStartPar
\sphinxstylestrong{Synopsis}
\begin{quote}

\sphinxAtStartPar
\sphinxcode{\sphinxupquote{setSlackBasis()}}
\end{quote}

\sphinxAtStartPar
\sphinxstylestrong{Description}
\begin{quote}

\sphinxAtStartPar
Set LP basis to be slack.
\end{quote}

\sphinxAtStartPar
\sphinxstylestrong{Example}
\end{quote}

\begin{sphinxVerbatim}[commandchars=\\\{\}]
\PYG{c+c1}{\PYGZsh{}  Set LP basis to be slack.}
\PYG{n}{m}\PYG{o}{.}\PYG{n}{setSlackBasis}\PYG{p}{(}\PYG{p}{)}
\end{sphinxVerbatim}

\subsubsection{Model.setVarType()}
\label{\detokenize{pyapiref:model-setvartype}}\begin{quote}

\sphinxAtStartPar
\sphinxstylestrong{Synopsis}
\begin{quote}

\sphinxAtStartPar
\sphinxcode{\sphinxupquote{setVarType(vars, vartypes)}}
\end{quote}

\sphinxAtStartPar
\sphinxstylestrong{Description}
\begin{quote}

\sphinxAtStartPar
Set the type of specific variable.

\sphinxAtStartPar
If parameter \sphinxcode{\sphinxupquote{vars}} is {\hyperref[\detokenize{pyapiref:chappyapi-var}]{\sphinxcrossref{\DUrole{std,std-ref}{Var Class}}}} object,
then parameter \sphinxcode{\sphinxupquote{vartypes}} is {\hyperref[\detokenize{constant:chapconst-vartype}]{\sphinxcrossref{\DUrole{std,std-ref}{Variable types}}}} constant;

\sphinxAtStartPar
If parameter \sphinxcode{\sphinxupquote{vars}} is dictionary or {\hyperref[\detokenize{pyapiref:chappyapi-util-tupledict}]{\sphinxcrossref{\DUrole{std,std-ref}{tupledict Class}}}} object,
then parameter \sphinxcode{\sphinxupquote{vartypes}} can be {\hyperref[\detokenize{constant:chapconst-vartype}]{\sphinxcrossref{\DUrole{std,std-ref}{Variable types}}}} constart, dictionary or
{\hyperref[\detokenize{pyapiref:chappyapi-util-tupledict}]{\sphinxcrossref{\DUrole{std,std-ref}{tupledict Class}}}} object;

\sphinxAtStartPar
If parameter \sphinxcode{\sphinxupquote{vars}} is list or {\hyperref[\detokenize{pyapiref:chappyapi-vararray}]{\sphinxcrossref{\DUrole{std,std-ref}{VarArray Class}}}} object,
then parameter \sphinxcode{\sphinxupquote{vartypes}} can be {\hyperref[\detokenize{constant:chapconst-vartype}]{\sphinxcrossref{\DUrole{std,std-ref}{Variable types}}}} constart or list object.
\end{quote}

\sphinxAtStartPar
\sphinxstylestrong{Arguments}
\begin{quote}

\sphinxAtStartPar
\sphinxcode{\sphinxupquote{vars}}
\begin{quote}

\sphinxAtStartPar
The specified variable.
\end{quote}

\sphinxAtStartPar
\sphinxcode{\sphinxupquote{vartypes}}
\begin{quote}

\sphinxAtStartPar
Type of the specified variable.
\end{quote}
\end{quote}

\sphinxAtStartPar
\sphinxstylestrong{Example}
\end{quote}

\begin{sphinxVerbatim}[commandchars=\\\{\}]
\PYG{c+c1}{\PYGZsh{} Set variable x as integer variable.}
\PYG{n}{m}\PYG{o}{.}\PYG{n}{setVarType}\PYG{p}{(}\PYG{n}{x}\PYG{p}{,} \PYG{n}{COPT}\PYG{o}{.}\PYG{n}{INTEGER}\PYG{p}{)}
\PYG{c+c1}{\PYGZsh{} Set variables x and y as binary variables.}
\PYG{n}{m}\PYG{o}{.}\PYG{n}{setVarType}\PYG{p}{(}\PYG{p}{[}\PYG{n}{x}\PYG{p}{,} \PYG{n}{y}\PYG{p}{]}\PYG{p}{,} \PYG{n}{COPT}\PYG{o}{.}\PYG{n}{BINARY}\PYG{p}{)}
\PYG{c+c1}{\PYGZsh{} Set the variables in tupledict object xdict as continuous variables.}
\PYG{n}{m}\PYG{o}{.}\PYG{n}{setVarType}\PYG{p}{(}\PYG{n}{xdict}\PYG{p}{,} \PYG{n}{COPT}\PYG{o}{.}\PYG{n}{CONTINUOUS}\PYG{p}{)}
\end{sphinxVerbatim}

\subsubsection{Model.setNames()}
\label{\detokenize{pyapiref:model-setnames}}\begin{quote}

\sphinxAtStartPar
\sphinxstylestrong{Synopsis}
\begin{quote}

\sphinxAtStartPar
\sphinxcode{\sphinxupquote{setNames(args, names)}}
\end{quote}

\sphinxAtStartPar
\sphinxstylestrong{Description}
\begin{quote}

\sphinxAtStartPar
Sets the name(s) of the specified variable(s) or constraint(s).
\end{quote}

\sphinxAtStartPar
\sphinxstylestrong{Arguments}
\begin{quote}

\sphinxAtStartPar
\sphinxcode{\sphinxupquote{args}}
\begin{quote}

\sphinxAtStartPar
The specified variable(s) or constraint(s), which could be:
{\hyperref[\detokenize{pyapiref:chappyapi-var}]{\sphinxcrossref{\DUrole{std,std-ref}{Var Class}}}} , {\hyperref[\detokenize{pyapiref:chappyapi-constraint}]{\sphinxcrossref{\DUrole{std,std-ref}{Constraint Class}}}} , {\hyperref[\detokenize{pyapiref:chappyapi-qconstraint}]{\sphinxcrossref{\DUrole{std,std-ref}{QConstraint Class}}}} ,
{\hyperref[\detokenize{pyapiref:chappyapi-psdvar}]{\sphinxcrossref{\DUrole{std,std-ref}{PsdVar Class}}}} , {\hyperref[\detokenize{pyapiref:chappyapi-psdconstraint}]{\sphinxcrossref{\DUrole{std,std-ref}{PsdConstraint Class}}}} , {\hyperref[\detokenize{pyapiref:chappyapi-lmiconstraint}]{\sphinxcrossref{\DUrole{std,std-ref}{LmiConstraint Class}}}} ,
{\hyperref[\detokenize{pyapiref:chappyapi-genconstr}]{\sphinxcrossref{\DUrole{std,std-ref}{GenConstr Class}}}} object, {\hyperref[\detokenize{pyapiref:chappyapi-affinecone}]{\sphinxcrossref{\DUrole{std,std-ref}{AffineCone Class}}}} or the dictionary/list they constitute.
\end{quote}

\sphinxAtStartPar
\sphinxcode{\sphinxupquote{names}}
\begin{quote}

\sphinxAtStartPar
The name(s) of the specified variable(s) or constraint(s), which could be a single string or the list/dictionary
corresponding to the \sphinxcode{\sphinxupquote{args}} .
\end{quote}
\end{quote}

\sphinxAtStartPar
\sphinxstylestrong{Example}
\end{quote}

\begin{sphinxVerbatim}[commandchars=\\\{\}]
\PYG{c+c1}{\PYGZsh{} Set the name of x to \PYGZdq{}var\PYGZdq{}}
\PYG{n}{m}\PYG{o}{.}\PYG{n}{setNames}\PYG{p}{(}\PYG{n}{x}\PYG{p}{,} \PYG{l+s+s2}{\PYGZdq{}}\PYG{l+s+s2}{var}\PYG{l+s+s2}{\PYGZdq{}}\PYG{p}{)}
\PYG{c+c1}{\PYGZsh{}Set the name of constr1 to \PYGZdq{}c1\PYGZdq{} and the name of constr2 to \PYGZdq{}c2\PYGZdq{}}
\PYG{n}{m}\PYG{o}{.}\PYG{n}{setNames}\PYG{p}{(}\PYG{p}{[}\PYG{n}{constr1}\PYG{p}{,} \PYG{n}{constr2}\PYG{p}{]}\PYG{p}{,} \PYG{p}{[}\PYG{l+s+s2}{\PYGZdq{}}\PYG{l+s+s2}{c1}\PYG{l+s+s2}{\PYGZdq{}}\PYG{p}{,} \PYG{l+s+s2}{\PYGZdq{}}\PYG{l+s+s2}{c2}\PYG{l+s+s2}{\PYGZdq{}}\PYG{p}{]}\PYG{p}{)}
\end{sphinxVerbatim}

\subsubsection{Model.setMipStart()}
\label{\detokenize{pyapiref:model-setmipstart}}\begin{quote}

\sphinxAtStartPar
\sphinxstylestrong{Synopsis}
\begin{quote}

\sphinxAtStartPar
\sphinxcode{\sphinxupquote{setMipStart(vars, startvals)}}
\end{quote}

\sphinxAtStartPar
\sphinxstylestrong{Description}
\begin{quote}

\sphinxAtStartPar
Set initial value for specified variables, valid only for integer programming.

\sphinxAtStartPar
If parameter \sphinxcode{\sphinxupquote{vars}} is {\hyperref[\detokenize{pyapiref:chappyapi-var}]{\sphinxcrossref{\DUrole{std,std-ref}{Var Class}}}} object, then parameter \sphinxcode{\sphinxupquote{startvals}} is constant;
If parameter  \sphinxcode{\sphinxupquote{vars}} is dictionary or {\hyperref[\detokenize{pyapiref:chappyapi-util-tupledict}]{\sphinxcrossref{\DUrole{std,std-ref}{tupledict Class}}}} object,
then parameter \sphinxcode{\sphinxupquote{startvals}} can be constant, dictionary or {\hyperref[\detokenize{pyapiref:chappyapi-util-tupledict}]{\sphinxcrossref{\DUrole{std,std-ref}{tupledict Class}}}} object;
If parameter  \sphinxcode{\sphinxupquote{vars}} is list or {\hyperref[\detokenize{pyapiref:chappyapi-vararray}]{\sphinxcrossref{\DUrole{std,std-ref}{VarArray Class}}}} object,
then parameter \sphinxcode{\sphinxupquote{startvals}} can be constant or list object.

\sphinxAtStartPar
\sphinxstylestrong{Notice:} You may want to call this method several times to input the MIP start. Please call \sphinxcode{\sphinxupquote{loadMipStart()}} once when the input is done.
\end{quote}

\sphinxAtStartPar
\sphinxstylestrong{Arguments}
\begin{quote}

\sphinxAtStartPar
\sphinxcode{\sphinxupquote{vars}}
\begin{quote}

\sphinxAtStartPar
The specified variable.
\end{quote}

\sphinxAtStartPar
\sphinxcode{\sphinxupquote{startvals}}
\begin{quote}

\sphinxAtStartPar
Initial value of the specified variable
\end{quote}
\end{quote}

\sphinxAtStartPar
\sphinxstylestrong{Example}
\end{quote}

\begin{sphinxVerbatim}[commandchars=\\\{\}]
\PYG{c+c1}{\PYGZsh{} Set initial value of x as 1.}
\PYG{n}{m}\PYG{o}{.}\PYG{n}{setMipStart}\PYG{p}{(}\PYG{n}{x}\PYG{p}{,} \PYG{l+m+mi}{1}\PYG{p}{)}
\PYG{c+c1}{\PYGZsh{} Set initial value of x, y as 2, 3.}
\PYG{n}{m}\PYG{o}{.}\PYG{n}{setMipStart}\PYG{p}{(}\PYG{p}{[}\PYG{n}{x}\PYG{p}{,} \PYG{n}{y}\PYG{p}{]}\PYG{p}{,} \PYG{p}{[}\PYG{l+m+mi}{2}\PYG{p}{,} \PYG{l+m+mi}{3}\PYG{p}{]}\PYG{p}{)}
\PYG{c+c1}{\PYGZsh{} Set initial value of all variables in tupledict xdict as 1.}
\PYG{n}{m}\PYG{o}{.}\PYG{n}{setMipStart}\PYG{p}{(}\PYG{n}{xdict}\PYG{p}{,} \PYG{l+m+mi}{1}\PYG{p}{)}

\PYG{c+c1}{\PYGZsh{} Load initial solution to model}
\PYG{n}{m}\PYG{o}{.}\PYG{n}{loadMipStart}\PYG{p}{(}\PYG{p}{)}
\end{sphinxVerbatim}

\subsubsection{Model.loadMipStart()}
\label{\detokenize{pyapiref:model-loadmipstart}}\begin{quote}

\sphinxAtStartPar
\sphinxstylestrong{Synopsis}
\begin{quote}

\sphinxAtStartPar
\sphinxcode{\sphinxupquote{loadMipStart()}}
\end{quote}

\sphinxAtStartPar
\sphinxstylestrong{Description}
\begin{quote}

\sphinxAtStartPar
Load the currently specified initial values into model.

\sphinxAtStartPar
\sphinxstylestrong{Notice:} After calling this method, the previously initial values will be cleared, and users can continue to set
initial values for specified variables.
\end{quote}
\end{quote}

\subsubsection{Model.setNlPrimalStart()}
\label{\detokenize{pyapiref:model-setnlprimalstart}}\begin{quote}

\sphinxAtStartPar
\sphinxstylestrong{Synopsis}
\begin{quote}

\sphinxAtStartPar
\sphinxcode{\sphinxupquote{setNlPrimalStart(vars, startvals)}}
\end{quote}

\sphinxAtStartPar
\sphinxstylestrong{Description}
\begin{quote}

\sphinxAtStartPar
Set the initial values of variables in a nonlinear model.
\end{quote}

\sphinxAtStartPar
\sphinxstylestrong{Arguments}
\begin{quote}

\sphinxAtStartPar
\sphinxcode{\sphinxupquote{vars}}
\begin{quote}

\sphinxAtStartPar
The variables to be assigned initial values to.

\sphinxAtStartPar
Acceptable values include:
a single {\hyperref[\detokenize{pyapiref:chappyapi-var}]{\sphinxcrossref{\DUrole{std,std-ref}{Var Class}}}} object,
a {\hyperref[\detokenize{pyapiref:chappyapi-vararray}]{\sphinxcrossref{\DUrole{std,std-ref}{VarArray Class}}}} or any iterable of variables,
a dictionary mapping arbitrary keys to {\hyperref[\detokenize{pyapiref:chappyapi-var}]{\sphinxcrossref{\DUrole{std,std-ref}{Var Class}}}} objects, or
a {\hyperref[\detokenize{pyapiref:chappyapi-mvar}]{\sphinxcrossref{\DUrole{std,std-ref}{MVar Class}}}} multidimensional variable.
\end{quote}

\sphinxAtStartPar
\sphinxcode{\sphinxupquote{startvals}}
\begin{quote}

\sphinxAtStartPar
The corresponding start values.
Can be a double value, a list, array, or a mapping structure consistent
with the shape or keys of the variables.
\end{quote}
\end{quote}
\end{quote}

\subsubsection{Model.setInfo()}
\label{\detokenize{pyapiref:model-setinfo}}\begin{quote}

\sphinxAtStartPar
\sphinxstylestrong{Synopsis}
\begin{quote}

\sphinxAtStartPar
\sphinxcode{\sphinxupquote{setInfo(infoname, args, newvals)}}
\end{quote}

\sphinxAtStartPar
\sphinxstylestrong{Description}
\begin{quote}

\sphinxAtStartPar
Set new information value for specific variables or constraints.

\sphinxAtStartPar
If parameter \sphinxcode{\sphinxupquote{args}} is {\hyperref[\detokenize{pyapiref:chappyapi-var}]{\sphinxcrossref{\DUrole{std,std-ref}{Var Class}}}} object or {\hyperref[\detokenize{pyapiref:chappyapi-constraint}]{\sphinxcrossref{\DUrole{std,std-ref}{Constraint Class}}}} object,
then parameter \sphinxcode{\sphinxupquote{newvals}} is constant;
If parameter \sphinxcode{\sphinxupquote{args}} is dictionary or {\hyperref[\detokenize{pyapiref:chappyapi-util-tupledict}]{\sphinxcrossref{\DUrole{std,std-ref}{tupledict Class}}}} object,
then parameter \sphinxcode{\sphinxupquote{newvals}} can be constant, dictionary or {\hyperref[\detokenize{pyapiref:chappyapi-util-tupledict}]{\sphinxcrossref{\DUrole{std,std-ref}{tupledict Class}}}} object;
If parameter \sphinxcode{\sphinxupquote{args}} is list, {\hyperref[\detokenize{pyapiref:chappyapi-vararray}]{\sphinxcrossref{\DUrole{std,std-ref}{VarArray Class}}}} object or {\hyperref[\detokenize{pyapiref:chappyapi-constrarray}]{\sphinxcrossref{\DUrole{std,std-ref}{ConstrArray Class}}}} object,
then parameter \sphinxcode{\sphinxupquote{newvals}} can be constant or list;
If parameter \sphinxcode{\sphinxupquote{args}} is {\hyperref[\detokenize{pyapiref:chappyapi-mvar}]{\sphinxcrossref{\DUrole{std,std-ref}{MVar Class}}}} object, {\hyperref[\detokenize{pyapiref:chappyapi-mconstr}]{\sphinxcrossref{\DUrole{std,std-ref}{MConstr Class}}}} object,
{\hyperref[\detokenize{pyapiref:chappyapi-mqconstr}]{\sphinxcrossref{\DUrole{std,std-ref}{MQConstr Class}}}} object or {\hyperref[\detokenize{pyapiref:chappyapi-mpsdconstr}]{\sphinxcrossref{\DUrole{std,std-ref}{MPsdConstr Class}}}} object.
then parameter \sphinxcode{\sphinxupquote{newvals}} can be constant or \sphinxcode{\sphinxupquote{numpy.ndarray}}.
\end{quote}

\sphinxAtStartPar
\sphinxstylestrong{Arguments}
\begin{quote}

\sphinxAtStartPar
\sphinxcode{\sphinxupquote{infoname}}
\begin{quote}

\sphinxAtStartPar
The specified information name. The full list of available names can be found in {\hyperref[\detokenize{pyapiref:chappyapi-const-info}]{\sphinxcrossref{\DUrole{std,std-ref}{Information}}}} section.
\end{quote}

\sphinxAtStartPar
\sphinxcode{\sphinxupquote{args}}
\begin{quote}

\sphinxAtStartPar
The specified variables of constraints.
\end{quote}

\sphinxAtStartPar
\sphinxcode{\sphinxupquote{newvals}}
\begin{quote}

\sphinxAtStartPar
Value of the specified information.
\end{quote}
\end{quote}

\sphinxAtStartPar
\sphinxstylestrong{Example}
\end{quote}

\begin{sphinxVerbatim}[commandchars=\\\{\}]
\PYG{n}{m}\PYG{o}{.}\PYG{n}{setInfo}\PYG{p}{(}\PYG{n}{COPT}\PYG{o}{.}\PYG{n}{Info}\PYG{o}{.}\PYG{n}{LB}\PYG{p}{,} \PYG{p}{[}\PYG{n}{x}\PYG{p}{,} \PYG{n}{y}\PYG{p}{]}\PYG{p}{,} \PYG{p}{[}\PYG{l+m+mf}{1.0}\PYG{p}{,} \PYG{l+m+mf}{2.0}\PYG{p}{]}\PYG{p}{)}

\PYG{c+c1}{\PYGZsh{} Set the upperbound of variable x as 1.0}
\PYG{n}{m}\PYG{o}{.}\PYG{n}{setInfo}\PYG{p}{(}\PYG{n}{COPT}\PYG{o}{.}\PYG{n}{Info}\PYG{o}{.}\PYG{n}{UB}\PYG{p}{,} \PYG{n}{x}\PYG{p}{,} \PYG{l+m+mf}{1.0}\PYG{p}{)}
\PYG{c+c1}{\PYGZsh{} Set the lowerbound of variables x and y as 1.0, 2.0, respectively.}
\PYG{n}{m}\PYG{o}{.}\PYG{n}{setInfo}\PYG{p}{(}\PYG{n}{COPT}\PYG{o}{.}\PYG{n}{Info}\PYG{o}{.}\PYG{n}{LB}\PYG{p}{,} \PYG{p}{[}\PYG{n}{x}\PYG{p}{,} \PYG{n}{y}\PYG{p}{]}\PYG{p}{,} \PYG{p}{[}\PYG{l+m+mf}{1.0}\PYG{p}{,} \PYG{l+m+mf}{2.0}\PYG{p}{]}\PYG{p}{)}
\PYG{c+c1}{\PYGZsh{} Set the objective of all variables in tupledict xdict as 0.}
\PYG{n}{m}\PYG{o}{.}\PYG{n}{setInfo}\PYG{p}{(}\PYG{n}{COPT}\PYG{o}{.}\PYG{n}{Info}\PYG{o}{.}\PYG{n}{OBJ}\PYG{p}{,} \PYG{n}{xdict}\PYG{p}{,} \PYG{l+m+mf}{0.0}\PYG{p}{)}
\end{sphinxVerbatim}

\subsubsection{Model.setParam()}
\label{\detokenize{pyapiref:model-setparam}}\begin{quote}

\sphinxAtStartPar
\sphinxstylestrong{Synopsis}
\begin{quote}

\sphinxAtStartPar
\sphinxcode{\sphinxupquote{setParam(paramname, newval)}}
\end{quote}

\sphinxAtStartPar
\sphinxstylestrong{Description}
\begin{quote}

\sphinxAtStartPar
Set the value of a parameter to a specific value.
\end{quote}

\sphinxAtStartPar
\sphinxstylestrong{Arguments}
\begin{quote}

\sphinxAtStartPar
\sphinxcode{\sphinxupquote{paramname}}
\begin{quote}

\sphinxAtStartPar
The name of parameter to be set. The list of available names can be found in {\hyperref[\detokenize{pyapiref:chappyapi-const-param}]{\sphinxcrossref{\DUrole{std,std-ref}{Parameters}}}} section.
\end{quote}

\sphinxAtStartPar
\sphinxcode{\sphinxupquote{newval}}
\begin{quote}

\sphinxAtStartPar
New value of parameter.
\end{quote}
\end{quote}

\sphinxAtStartPar
\sphinxstylestrong{Example}
\end{quote}

\begin{sphinxVerbatim}[commandchars=\\\{\}]
\PYG{c+c1}{\PYGZsh{} Set time limit of solving to 1 hour.}
 \PYG{n}{m}\PYG{o}{.}\PYG{n}{setParam}\PYG{p}{(}\PYG{n}{COPT}\PYG{o}{.}\PYG{n}{Param}\PYG{o}{.}\PYG{n}{TimeLimit}\PYG{p}{,} \PYG{l+m+mi}{3600}\PYG{p}{)}
\end{sphinxVerbatim}

\subsubsection{Model.resetParam()}
\label{\detokenize{pyapiref:model-resetparam}}\begin{quote}

\sphinxAtStartPar
\sphinxstylestrong{Synopsis}
\begin{quote}

\sphinxAtStartPar
\sphinxcode{\sphinxupquote{resetParam()}}
\end{quote}

\sphinxAtStartPar
\sphinxstylestrong{Description}
\begin{quote}

\sphinxAtStartPar
Reset all parameters to their default values.
\end{quote}

\sphinxAtStartPar
\sphinxstylestrong{Example}
\end{quote}

\begin{sphinxVerbatim}[commandchars=\\\{\}]
\PYG{c+c1}{\PYGZsh{} Reset all parameters to their default values.}
\PYG{n}{m}\PYG{o}{.}\PYG{n}{resetParam}\PYG{p}{(}\PYG{p}{)}
\end{sphinxVerbatim}

\subsubsection{Model.setObjParamN()}
\label{\detokenize{pyapiref:model-setobjparamn}}\begin{quote}

\sphinxAtStartPar
\sphinxstylestrong{Synopsis}
\begin{quote}

\sphinxAtStartPar
\sphinxcode{\sphinxupquote{setObjParamN(idx, paramname, newval)}}
\end{quote}

\sphinxAtStartPar
\sphinxstylestrong{Description}
\begin{quote}

\sphinxAtStartPar
Set a parameter of the specified objective in multi\sphinxhyphen{}objective optimization.
\end{quote}

\sphinxAtStartPar
\sphinxstylestrong{Arguments}
\begin{quote}

\sphinxAtStartPar
\sphinxcode{\sphinxupquote{idx}}
\begin{quote}

\sphinxAtStartPar
Index of the objective function.
\end{quote}

\sphinxAtStartPar
\sphinxcode{\sphinxupquote{paramname}}
\begin{quote}

\sphinxAtStartPar
Name of the objective parameter.

\sphinxAtStartPar
Possible values are {\hyperref[\detokenize{multiobj:multiobjpriority}]{\sphinxcrossref{\DUrole{std,std-ref}{MultiObjPriority}}}},
{\hyperref[\detokenize{multiobj:multiobjweight}]{\sphinxcrossref{\DUrole{std,std-ref}{MultiObjWeight}}}},
{\hyperref[\detokenize{multiobj:multiobjreltol}]{\sphinxcrossref{\DUrole{std,std-ref}{MultiObjRelTol}}}},
{\hyperref[\detokenize{multiobj:multiobjabstol}]{\sphinxcrossref{\DUrole{std,std-ref}{MultiObjAbsTol}}}}.
\end{quote}

\sphinxAtStartPar
\sphinxcode{\sphinxupquote{newval}}
\begin{quote}

\sphinxAtStartPar
New value for the specified objective function parameter.
\end{quote}
\end{quote}
\end{quote}

\subsubsection{Model.setParamN()}
\label{\detokenize{pyapiref:model-setparamn}}\begin{quote}

\sphinxAtStartPar
\sphinxstylestrong{Synopsis}
\begin{quote}

\sphinxAtStartPar
\sphinxcode{\sphinxupquote{setParamN(idx, paramname, newval)}}
\end{quote}

\sphinxAtStartPar
\sphinxstylestrong{Description}
\begin{quote}

\sphinxAtStartPar
Set a solver parameter for the model associated with
the specified objective in multi\sphinxhyphen{}objective optimization.
\end{quote}

\sphinxAtStartPar
\sphinxstylestrong{Arguments}
\begin{quote}

\sphinxAtStartPar
\sphinxcode{\sphinxupquote{idx}}
\begin{quote}

\sphinxAtStartPar
Index of the objective function.
\end{quote}

\sphinxAtStartPar
\sphinxcode{\sphinxupquote{paramname}}
\begin{quote}

\sphinxAtStartPar
Name of the solver parameter.

\sphinxAtStartPar
See {\hyperref[\detokenize{pyapiref:chappyapi-const-param}]{\sphinxcrossref{\DUrole{std,std-ref}{Parameters}}}} for possible values.
\end{quote}

\sphinxAtStartPar
\sphinxcode{\sphinxupquote{newval}}
\begin{quote}

\sphinxAtStartPar
New value of the specified solver parameter.
\end{quote}
\end{quote}
\end{quote}

\subsubsection{Model.resetObjParamN()}
\label{\detokenize{pyapiref:model-resetobjparamn}}\begin{quote}

\sphinxAtStartPar
\sphinxstylestrong{Synopsis}
\begin{quote}

\sphinxAtStartPar
\sphinxcode{\sphinxupquote{resetObjParamN(idx)}}
\end{quote}

\sphinxAtStartPar
\sphinxstylestrong{Description}
\begin{quote}

\sphinxAtStartPar
Reset all parameters of the
specified objective in multi\sphinxhyphen{}objective optimization to their default values.
\end{quote}

\sphinxAtStartPar
\sphinxstylestrong{Arguments}
\begin{quote}

\sphinxAtStartPar
\sphinxcode{\sphinxupquote{idx}}
\begin{quote}

\sphinxAtStartPar
Index of the objective function.
\end{quote}
\end{quote}
\end{quote}

\subsubsection{Model.resetParamN()}
\label{\detokenize{pyapiref:model-resetparamn}}\begin{quote}

\sphinxAtStartPar
\sphinxstylestrong{Synopsis}
\begin{quote}

\sphinxAtStartPar
\sphinxcode{\sphinxupquote{resetParamN(idx)}}
\end{quote}

\sphinxAtStartPar
\sphinxstylestrong{Description}
\begin{quote}

\sphinxAtStartPar
Reset all solver parameters associated
with the specified objective in multi\sphinxhyphen{}objective
optimization to their default values.
\end{quote}

\sphinxAtStartPar
\sphinxstylestrong{Arguments}
\begin{quote}

\sphinxAtStartPar
\sphinxcode{\sphinxupquote{idx}}
\begin{quote}

\sphinxAtStartPar
Index of the objective function.
\end{quote}
\end{quote}
\end{quote}

\subsubsection{Model.read()}
\label{\detokenize{pyapiref:model-read}}\begin{quote}

\sphinxAtStartPar
\sphinxstylestrong{Synopsis}
\begin{quote}

\sphinxAtStartPar
\sphinxcode{\sphinxupquote{read(filename)}}
\end{quote}

\sphinxAtStartPar
\sphinxstylestrong{Description}
\begin{quote}

\sphinxAtStartPar
Determine the type of data by the file suffix and read it into a model.

\sphinxAtStartPar
Currently, it supports MPS files (suffix \sphinxcode{\sphinxupquote{\textquotesingle{}.mps\textquotesingle{}}} or \sphinxcode{\sphinxupquote{\textquotesingle{}.mps.gz\textquotesingle{}}}),
LP files (suffix \sphinxcode{\sphinxupquote{\textquotesingle{}.lp\textquotesingle{}}} or \sphinxcode{\sphinxupquote{\textquotesingle{}.lp.gz\textquotesingle{}}}),
SDPA files (suffix \sphinxcode{\sphinxupquote{\textquotesingle{}.dat\sphinxhyphen{}s\textquotesingle{}}} or \sphinxcode{\sphinxupquote{\textquotesingle{}.dat\sphinxhyphen{}s.gz\textquotesingle{}}}),
CBF files (suffix \sphinxcode{\sphinxupquote{\textquotesingle{}.cbf\textquotesingle{}}} or \sphinxcode{\sphinxupquote{\textquotesingle{}.cbf.gz\textquotesingle{}}}),
COPT binary format files (suffix \sphinxcode{\sphinxupquote{\textquotesingle{}.bin\textquotesingle{}}}),
basis files (suffix \sphinxcode{\sphinxupquote{\textquotesingle{}.bas\textquotesingle{}}}),
result files (suffix \sphinxcode{\sphinxupquote{\textquotesingle{}.sol\textquotesingle{}}}),
start files (suffix \sphinxcode{\sphinxupquote{\textquotesingle{}.mst\textquotesingle{}}}),
branching order files (suffix \sphinxcode{\sphinxupquote{\textquotesingle{}.ord\textquotesingle{}}}),
and parameter files (suffix \sphinxcode{\sphinxupquote{\textquotesingle{}.par\textquotesingle{}}}).
\end{quote}

\sphinxAtStartPar
\sphinxstylestrong{Arguments}
\begin{quote}

\sphinxAtStartPar
\sphinxcode{\sphinxupquote{filename}}
\begin{quote}

\sphinxAtStartPar
Name of the file to be read.
\end{quote}
\end{quote}

\sphinxAtStartPar
\sphinxstylestrong{Example}
\end{quote}

\begin{sphinxVerbatim}[commandchars=\\\{\}]
\PYG{c+c1}{\PYGZsh{} Read MPS format model file}
\PYG{n}{m}\PYG{o}{.}\PYG{n}{read}\PYG{p}{(}\PYG{l+s+s1}{\PYGZsq{}}\PYG{l+s+s1}{test.mps.gz}\PYG{l+s+s1}{\PYGZsq{}}\PYG{p}{)}
\PYG{c+c1}{\PYGZsh{} Read LP format model file}
\PYG{n}{m}\PYG{o}{.}\PYG{n}{read}\PYG{p}{(}\PYG{l+s+s1}{\PYGZsq{}}\PYG{l+s+s1}{test.lp.gz}\PYG{l+s+s1}{\PYGZsq{}}\PYG{p}{)}
\PYG{c+c1}{\PYGZsh{} Read COPT binary format model file}
\PYG{n}{m}\PYG{o}{.}\PYG{n}{read}\PYG{p}{(}\PYG{l+s+s1}{\PYGZsq{}}\PYG{l+s+s1}{test.bin}\PYG{l+s+s1}{\PYGZsq{}}\PYG{p}{)}
\PYG{c+c1}{\PYGZsh{} Read basis file}
\PYG{n}{m}\PYG{o}{.}\PYG{n}{read}\PYG{p}{(}\PYG{l+s+s1}{\PYGZsq{}}\PYG{l+s+s1}{testlp.bas}\PYG{l+s+s1}{\PYGZsq{}}\PYG{p}{)}
\PYG{c+c1}{\PYGZsh{} Read solution file}
\PYG{n}{m}\PYG{o}{.}\PYG{n}{read}\PYG{p}{(}\PYG{l+s+s1}{\PYGZsq{}}\PYG{l+s+s1}{testmip.sol}\PYG{l+s+s1}{\PYGZsq{}}\PYG{p}{)}
\PYG{c+c1}{\PYGZsh{} Read start file}
\PYG{n}{m}\PYG{o}{.}\PYG{n}{read}\PYG{p}{(}\PYG{l+s+s1}{\PYGZsq{}}\PYG{l+s+s1}{testmip.mst}\PYG{l+s+s1}{\PYGZsq{}}\PYG{p}{)}
\PYG{c+c1}{\PYGZsh{} Read paramter file}
\PYG{n}{m}\PYG{o}{.}\PYG{n}{read}\PYG{p}{(}\PYG{l+s+s1}{\PYGZsq{}}\PYG{l+s+s1}{test.par}\PYG{l+s+s1}{\PYGZsq{}}\PYG{p}{)}
\end{sphinxVerbatim}

\subsubsection{Model.readMps()}
\label{\detokenize{pyapiref:model-readmps}}\begin{quote}

\sphinxAtStartPar
\sphinxstylestrong{Synopsis}
\begin{quote}

\sphinxAtStartPar
\sphinxcode{\sphinxupquote{readMps(filename)}}
\end{quote}

\sphinxAtStartPar
\sphinxstylestrong{Description}
\begin{quote}

\sphinxAtStartPar
Read MPS file to model.
\end{quote}

\sphinxAtStartPar
\sphinxstylestrong{Arguments}
\begin{quote}

\sphinxAtStartPar
\sphinxcode{\sphinxupquote{filename}}
\begin{quote}

\sphinxAtStartPar
The name of the MPS file to be read.
\end{quote}
\end{quote}

\sphinxAtStartPar
\sphinxstylestrong{Example}
\end{quote}

\begin{sphinxVerbatim}[commandchars=\\\{\}]
\PYG{c+c1}{\PYGZsh{} Read file \PYGZdq{}test.mps.gz\PYGZdq{} according to mps file format}
\PYG{n}{m}\PYG{o}{.}\PYG{n}{readMps}\PYG{p}{(}\PYG{l+s+s1}{\PYGZsq{}}\PYG{l+s+s1}{test.mps.gz}\PYG{l+s+s1}{\PYGZsq{}}\PYG{p}{)}
\PYG{c+c1}{\PYGZsh{} Read file \PYGZdq{}test.lp.gz\PYGZdq{} according mps file format}
\PYG{n}{m}\PYG{o}{.}\PYG{n}{readMps}\PYG{p}{(}\PYG{l+s+s1}{\PYGZsq{}}\PYG{l+s+s1}{test.lp.gz}\PYG{l+s+s1}{\PYGZsq{}}\PYG{p}{)}
\end{sphinxVerbatim}

\subsubsection{Model.readLp()}
\label{\detokenize{pyapiref:model-readlp}}\begin{quote}

\sphinxAtStartPar
\sphinxstylestrong{Synopsis}
\begin{quote}

\sphinxAtStartPar
\sphinxcode{\sphinxupquote{readLp(filename)}}
\end{quote}

\sphinxAtStartPar
\sphinxstylestrong{Description}
\begin{quote}

\sphinxAtStartPar
Read a file to model according to LP file format.
\end{quote}

\sphinxAtStartPar
\sphinxstylestrong{Arguments}
\begin{quote}

\sphinxAtStartPar
\sphinxcode{\sphinxupquote{filename}}
\begin{quote}

\sphinxAtStartPar
Name of the LP file to be read.
\end{quote}
\end{quote}

\sphinxAtStartPar
\sphinxstylestrong{Example}
\end{quote}

\begin{sphinxVerbatim}[commandchars=\\\{\}]
\PYG{c+c1}{\PYGZsh{} Read file\PYGZdq{}test.mps.gz\PYGZdq{} according to LP file format}
\PYG{n}{m}\PYG{o}{.}\PYG{n}{readLp}\PYG{p}{(}\PYG{l+s+s1}{\PYGZsq{}}\PYG{l+s+s1}{test.mps.gz}\PYG{l+s+s1}{\PYGZsq{}}\PYG{p}{)}
\PYG{c+c1}{\PYGZsh{} Read file\PYGZdq{}test.lp.gz\PYGZdq{} according to LP file format}
\PYG{n}{m}\PYG{o}{.}\PYG{n}{readLp}\PYG{p}{(}\PYG{l+s+s1}{\PYGZsq{}}\PYG{l+s+s1}{test.lp.gz}\PYG{l+s+s1}{\PYGZsq{}}\PYG{p}{)}
\end{sphinxVerbatim}

\subsubsection{Model.readSdpa()}
\label{\detokenize{pyapiref:model-readsdpa}}\begin{quote}

\sphinxAtStartPar
\sphinxstylestrong{Synopsis}
\begin{quote}

\sphinxAtStartPar
\sphinxcode{\sphinxupquote{readSdpa(filename)}}
\end{quote}

\sphinxAtStartPar
\sphinxstylestrong{Description}
\begin{quote}

\sphinxAtStartPar
Read a file to model according to SDPA file format.
\end{quote}

\sphinxAtStartPar
\sphinxstylestrong{Arguments}
\begin{quote}

\sphinxAtStartPar
\sphinxcode{\sphinxupquote{filename}}
\begin{quote}

\sphinxAtStartPar
Name of the SDPA file to be read.
\end{quote}
\end{quote}

\sphinxAtStartPar
\sphinxstylestrong{Example}
\end{quote}

\begin{sphinxVerbatim}[commandchars=\\\{\}]
\PYG{c+c1}{\PYGZsh{} Read file\PYGZdq{}test.dat\PYGZhy{}s\PYGZdq{} according to SDPA file format}
\PYG{n}{m}\PYG{o}{.}\PYG{n}{readSdpa}\PYG{p}{(}\PYG{l+s+s1}{\PYGZsq{}}\PYG{l+s+s1}{test.dat\PYGZhy{}s}\PYG{l+s+s1}{\PYGZsq{}}\PYG{p}{)}
\end{sphinxVerbatim}

\subsubsection{Model.readCbf()}
\label{\detokenize{pyapiref:model-readcbf}}\begin{quote}

\sphinxAtStartPar
\sphinxstylestrong{Synopsis}
\begin{quote}

\sphinxAtStartPar
\sphinxcode{\sphinxupquote{readCbf(filename)}}
\end{quote}

\sphinxAtStartPar
\sphinxstylestrong{Description}
\begin{quote}

\sphinxAtStartPar
Read a file to model according to CBF file format.
\end{quote}

\sphinxAtStartPar
\sphinxstylestrong{Arguments}
\begin{quote}

\sphinxAtStartPar
\sphinxcode{\sphinxupquote{filename}}
\begin{quote}

\sphinxAtStartPar
Name of the CBF file to be read.
\end{quote}
\end{quote}

\sphinxAtStartPar
\sphinxstylestrong{Example}
\end{quote}

\begin{sphinxVerbatim}[commandchars=\\\{\}]
\PYG{c+c1}{\PYGZsh{} Read file\PYGZdq{}test.cbf\PYGZdq{} according to CBF file format}
\PYG{n}{m}\PYG{o}{.}\PYG{n}{readCbf}\PYG{p}{(}\PYG{l+s+s1}{\PYGZsq{}}\PYG{l+s+s1}{test.cbf}\PYG{l+s+s1}{\PYGZsq{}}\PYG{p}{)}
\end{sphinxVerbatim}

\subsubsection{Model.readBin()}
\label{\detokenize{pyapiref:model-readbin}}\begin{quote}

\sphinxAtStartPar
\sphinxstylestrong{Synopsis}
\begin{quote}

\sphinxAtStartPar
\sphinxcode{\sphinxupquote{readBin(filename)}}
\end{quote}

\sphinxAtStartPar
\sphinxstylestrong{Description}
\begin{quote}

\sphinxAtStartPar
Read COPT binary format file to model.
\end{quote}

\sphinxAtStartPar
\sphinxstylestrong{Arguments}
\begin{quote}

\sphinxAtStartPar
\sphinxcode{\sphinxupquote{filename}}
\begin{quote}

\sphinxAtStartPar
The name of the COPT binary format file to be read.
\end{quote}
\end{quote}

\sphinxAtStartPar
\sphinxstylestrong{Example}
\end{quote}

\begin{sphinxVerbatim}[commandchars=\\\{\}]
\PYG{c+c1}{\PYGZsh{} Read file \PYGZdq{}test.bin\PYGZdq{} according COPT binary file format}
\PYG{n}{m}\PYG{o}{.}\PYG{n}{readBin}\PYG{p}{(}\PYG{l+s+s1}{\PYGZsq{}}\PYG{l+s+s1}{test.bin}\PYG{l+s+s1}{\PYGZsq{}}\PYG{p}{)}
\end{sphinxVerbatim}

\subsubsection{Model.readSol()}
\label{\detokenize{pyapiref:model-readsol}}\begin{quote}

\sphinxAtStartPar
\sphinxstylestrong{Synopsis}
\begin{quote}

\sphinxAtStartPar
\sphinxcode{\sphinxupquote{readSol(filename)}}
\end{quote}

\sphinxAtStartPar
\sphinxstylestrong{Description}
\begin{quote}

\sphinxAtStartPar
Read a file to model according to solution file format.

\sphinxAtStartPar
\sphinxstylestrong{Notice:} The default solution value is 0, i.e. a partial solution
will be automatically filled in with zeros. If read successfully, then
the values read can be act as initial solution for integer programming.
\end{quote}

\sphinxAtStartPar
\sphinxstylestrong{Arguments}
\begin{quote}

\sphinxAtStartPar
\sphinxcode{\sphinxupquote{filename}}
\begin{quote}

\sphinxAtStartPar
Name of file to be read.
\end{quote}
\end{quote}

\sphinxAtStartPar
\sphinxstylestrong{Example}
\end{quote}

\begin{sphinxVerbatim}[commandchars=\\\{\}]
\PYG{c+c1}{\PYGZsh{} Read file \PYGZdq{}testmip.sol\PYGZdq{} according to solution file format.}
\PYG{n}{m}\PYG{o}{.}\PYG{n}{readSol}\PYG{p}{(}\PYG{l+s+s1}{\PYGZsq{}}\PYG{l+s+s1}{testmip.sol}\PYG{l+s+s1}{\PYGZsq{}}\PYG{p}{)}
\PYG{c+c1}{\PYGZsh{} Read file testmip.txt\PYGZdq{} according to solution file format.}
\PYG{n}{m}\PYG{o}{.}\PYG{n}{readSol}\PYG{p}{(}\PYG{l+s+s1}{\PYGZsq{}}\PYG{l+s+s1}{testmip.txt}\PYG{l+s+s1}{\PYGZsq{}}\PYG{p}{)}
\end{sphinxVerbatim}

\subsubsection{Model.readJsonSol()}
\label{\detokenize{pyapiref:model-readjsonsol}}\begin{quote}

\sphinxAtStartPar
\sphinxstylestrong{Synopsis}
\begin{quote}

\sphinxAtStartPar
\sphinxcode{\sphinxupquote{readJsonSol(filename)}}
\end{quote}

\sphinxAtStartPar
\sphinxstylestrong{Description}
\begin{quote}

\sphinxAtStartPar
Read the complete solution in JSON format from a file.
\end{quote}

\sphinxAtStartPar
\sphinxstylestrong{Arguments}
\begin{quote}

\sphinxAtStartPar
\sphinxcode{\sphinxupquote{filename}}
\begin{quote}

\sphinxAtStartPar
Name of file to be read.
\end{quote}
\end{quote}
\end{quote}

\subsubsection{Model.readBasis()}
\label{\detokenize{pyapiref:model-readbasis}}\begin{quote}

\sphinxAtStartPar
\sphinxstylestrong{Synopsis}
\begin{quote}

\sphinxAtStartPar
\sphinxcode{\sphinxupquote{readBasis(filename)}}
\end{quote}

\sphinxAtStartPar
\sphinxstylestrong{Description}
\begin{quote}

\sphinxAtStartPar
Read basis status of variables and linear constraints to model accoring to basis solution format,
valid only for linear programming.
\end{quote}

\sphinxAtStartPar
\sphinxstylestrong{Arguments}
\begin{quote}

\sphinxAtStartPar
\sphinxcode{\sphinxupquote{filename}}
\begin{quote}

\sphinxAtStartPar
The name of the basis file to be read.
\end{quote}
\end{quote}

\sphinxAtStartPar
\sphinxstylestrong{Example}
\end{quote}

\begin{sphinxVerbatim}[commandchars=\\\{\}]
\PYG{c+c1}{\PYGZsh{} Read file \PYGZdq{}testmip.bas\PYGZdq{} to basis solution format}
\PYG{n}{m}\PYG{o}{.}\PYG{n}{readBasis}\PYG{p}{(}\PYG{l+s+s1}{\PYGZsq{}}\PYG{l+s+s1}{testmip.bas}\PYG{l+s+s1}{\PYGZsq{}}\PYG{p}{)}
\PYG{c+c1}{\PYGZsh{} Read file \PYGZdq{}testmip.txt\PYGZdq{} to basis solution format}
\PYG{n}{m}\PYG{o}{.}\PYG{n}{readBasis}\PYG{p}{(}\PYG{l+s+s1}{\PYGZsq{}}\PYG{l+s+s1}{testmip.txt}\PYG{l+s+s1}{\PYGZsq{}}\PYG{p}{)}
\end{sphinxVerbatim}

\subsubsection{Model.readMst()}
\label{\detokenize{pyapiref:model-readmst}}\begin{quote}

\sphinxAtStartPar
\sphinxstylestrong{Synopsis}
\begin{quote}

\sphinxAtStartPar
\sphinxcode{\sphinxupquote{readMst(filename)}}
\end{quote}

\sphinxAtStartPar
\sphinxstylestrong{Description}
\begin{quote}

\sphinxAtStartPar
Read initial solution data to model according to initial solution file format.

\sphinxAtStartPar
\sphinxstylestrong{Notice:}
If read successfully, the read value will be act as initial solution for integer programming model.
Variable values may not be specified completely, if the value of variable is specified for multiple times,
the last specified value is used.
\end{quote}

\sphinxAtStartPar
\sphinxstylestrong{Arguments}
\begin{quote}

\sphinxAtStartPar
\sphinxcode{\sphinxupquote{filename}}
\begin{quote}

\sphinxAtStartPar
Name of the file to be read.
\end{quote}
\end{quote}

\sphinxAtStartPar
\sphinxstylestrong{Example}
\end{quote}

\begin{sphinxVerbatim}[commandchars=\\\{\}]
\PYG{c+c1}{\PYGZsh{} Read file \PYGZdq{}testmip.mst\PYGZdq{} according to initial solution file format.}
\PYG{n}{m}\PYG{o}{.}\PYG{n}{readMst}\PYG{p}{(}\PYG{l+s+s1}{\PYGZsq{}}\PYG{l+s+s1}{testmip.mst}\PYG{l+s+s1}{\PYGZsq{}}\PYG{p}{)}
\PYG{c+c1}{\PYGZsh{} Read file \PYGZdq{}testmip.txt\PYGZdq{} according to initial solution file format.}
\PYG{n}{m}\PYG{o}{.}\PYG{n}{readMst}\PYG{p}{(}\PYG{l+s+s1}{\PYGZsq{}}\PYG{l+s+s1}{testmip.txt}\PYG{l+s+s1}{\PYGZsq{}}\PYG{p}{)}
\end{sphinxVerbatim}

\subsubsection{Model.readParam()}
\label{\detokenize{pyapiref:model-readparam}}\begin{quote}

\sphinxAtStartPar
\sphinxstylestrong{Synopsis}
\begin{quote}

\sphinxAtStartPar
\sphinxcode{\sphinxupquote{readParam(filename)}}
\end{quote}

\sphinxAtStartPar
\sphinxstylestrong{Description}
\begin{quote}

\sphinxAtStartPar
Read optimization parameters to model according to parameter file format.

\sphinxAtStartPar
\sphinxstylestrong{Notice:} If any optimization parameter is specified multiple times, the last spcified value is used.
\end{quote}

\sphinxAtStartPar
\sphinxstylestrong{Arguments}
\begin{quote}

\sphinxAtStartPar
\sphinxcode{\sphinxupquote{filename}}
\begin{quote}

\sphinxAtStartPar
The name of the parameter file to be read.
\end{quote}
\end{quote}

\sphinxAtStartPar
\sphinxstylestrong{Example}
\end{quote}

\begin{sphinxVerbatim}[commandchars=\\\{\}]
\PYG{c+c1}{\PYGZsh{} Read file \PYGZdq{}testmip.par\PYGZdq{} according to parameter file format.}
\PYG{n}{m}\PYG{o}{.}\PYG{n}{readParam}\PYG{p}{(}\PYG{l+s+s1}{\PYGZsq{}}\PYG{l+s+s1}{testmip.par}\PYG{l+s+s1}{\PYGZsq{}}\PYG{p}{)}
\PYG{c+c1}{\PYGZsh{} Read file \PYGZdq{}testmip.txt\PYGZdq{} according to parameter file format.}
\PYG{n}{m}\PYG{o}{.}\PYG{n}{readParam}\PYG{p}{(}\PYG{l+s+s1}{\PYGZsq{}}\PYG{l+s+s1}{testmip.txt}\PYG{l+s+s1}{\PYGZsq{}}\PYG{p}{)}
\end{sphinxVerbatim}

\subsubsection{Model.readTune()}
\label{\detokenize{pyapiref:model-readtune}}\begin{quote}

\sphinxAtStartPar
\sphinxstylestrong{Synopsis}
\begin{quote}

\sphinxAtStartPar
\sphinxcode{\sphinxupquote{readTune(filename)}}
\end{quote}

\sphinxAtStartPar
\sphinxstylestrong{Description}
\begin{quote}

\sphinxAtStartPar
Read the tuning parameters and combine them into the model according to the tuning file format.
\end{quote}

\sphinxAtStartPar
\sphinxstylestrong{Arguments}
\begin{quote}

\sphinxAtStartPar
\sphinxcode{\sphinxupquote{filename}}
\begin{quote}

\sphinxAtStartPar
The name of the file to be read.
\end{quote}
\end{quote}

\sphinxAtStartPar
\sphinxstylestrong{Example}
\end{quote}

\subsubsection{Model.readOrd()}
\label{\detokenize{pyapiref:model-readord}}\begin{quote}

\sphinxAtStartPar
\sphinxstylestrong{Synopsis}
\begin{quote}

\sphinxAtStartPar
\sphinxcode{\sphinxupquote{readOrd(filename)}}
\end{quote}

\sphinxAtStartPar
\sphinxstylestrong{Description}
\begin{quote}

\sphinxAtStartPar
Read a branching order file (ORD format) into the current model
to reuse a predefined branching strategy.
\end{quote}

\sphinxAtStartPar
\sphinxstylestrong{Arguments}
\begin{quote}

\sphinxAtStartPar
\sphinxcode{\sphinxupquote{filename}}
\begin{quote}

\sphinxAtStartPar
Name of the branching order file to be read.
\end{quote}
\end{quote}

\sphinxAtStartPar
\sphinxstylestrong{Example}

\begin{sphinxVerbatim}[commandchars=\\\{\}]
\PYG{c+c1}{\PYGZsh{} Read branching order information from \PYGZdq{}testmip.ord\PYGZdq{}}
\PYG{n}{m}\PYG{o}{.}\PYG{n}{readOrd}\PYG{p}{(}\PYG{l+s+s2}{\PYGZdq{}}\PYG{l+s+s2}{testmip.ord}\PYG{l+s+s2}{\PYGZdq{}}\PYG{p}{)}
\end{sphinxVerbatim}
\end{quote}

\subsubsection{Model.write()}
\label{\detokenize{pyapiref:model-write}}\begin{quote}

\sphinxAtStartPar
\sphinxstylestrong{Synopsis}
\begin{quote}

\sphinxAtStartPar
\sphinxcode{\sphinxupquote{write(filename)}}
\end{quote}

\sphinxAtStartPar
\sphinxstylestrong{Description}
\begin{quote}

\sphinxAtStartPar
Currently, COPT supports writing of
MPS files (suffix \sphinxcode{\sphinxupquote{\textquotesingle{}.mps\textquotesingle{}}} ),
LP files (suffix \sphinxcode{\sphinxupquote{\textquotesingle{}.lp\textquotesingle{}}}),
CBF files (suffix \sphinxcode{\sphinxupquote{\textquotesingle{}.cbf\textquotesingle{}}}),
COPT binary format files (suffix \sphinxcode{\sphinxupquote{\textquotesingle{}.bin\textquotesingle{}}}),
basis files (suffix \sphinxcode{\sphinxupquote{\textquotesingle{}.bas\textquotesingle{}}}),
LP solution files (suffix \sphinxcode{\sphinxupquote{\textquotesingle{}.sol\textquotesingle{}}}),
initial solution files (suffix \sphinxcode{\sphinxupquote{\textquotesingle{}.mst\textquotesingle{}}}),
branching order files (suffix \sphinxcode{\sphinxupquote{\textquotesingle{}.ord\textquotesingle{}}}),
and parameter files (suffix \sphinxcode{\sphinxupquote{\textquotesingle{}.par\textquotesingle{}}}).
\end{quote}

\sphinxAtStartPar
\sphinxstylestrong{Arguments}
\begin{quote}

\sphinxAtStartPar
\sphinxcode{\sphinxupquote{filename}}
\begin{quote}

\sphinxAtStartPar
The file name to be written.
\end{quote}
\end{quote}

\sphinxAtStartPar
\sphinxstylestrong{Example}
\end{quote}

\begin{sphinxVerbatim}[commandchars=\\\{\}]
\PYG{c+c1}{\PYGZsh{} Write MPS file}
\PYG{n}{m}\PYG{o}{.}\PYG{n}{write}\PYG{p}{(}\PYG{l+s+s1}{\PYGZsq{}}\PYG{l+s+s1}{test.mps}\PYG{l+s+s1}{\PYGZsq{}}\PYG{p}{)}
\PYG{c+c1}{\PYGZsh{} Write LP file}
\PYG{n}{m}\PYG{o}{.}\PYG{n}{write}\PYG{p}{(}\PYG{l+s+s1}{\PYGZsq{}}\PYG{l+s+s1}{test.lp}\PYG{l+s+s1}{\PYGZsq{}}\PYG{p}{)}
\PYG{c+c1}{\PYGZsh{} Write COPT binary format file}
\PYG{n}{m}\PYG{o}{.}\PYG{n}{write}\PYG{p}{(}\PYG{l+s+s1}{\PYGZsq{}}\PYG{l+s+s1}{test.bin}\PYG{l+s+s1}{\PYGZsq{}}\PYG{p}{)}
\PYG{c+c1}{\PYGZsh{} Write basis file}
\PYG{n}{m}\PYG{o}{.}\PYG{n}{write}\PYG{p}{(}\PYG{l+s+s1}{\PYGZsq{}}\PYG{l+s+s1}{testlp.bas}\PYG{l+s+s1}{\PYGZsq{}}\PYG{p}{)}
\PYG{c+c1}{\PYGZsh{} Write solution file}
\PYG{n}{m}\PYG{o}{.}\PYG{n}{write}\PYG{p}{(}\PYG{l+s+s1}{\PYGZsq{}}\PYG{l+s+s1}{testmip.sol}\PYG{l+s+s1}{\PYGZsq{}}\PYG{p}{)}
\PYG{c+c1}{\PYGZsh{} Write initial solution file}
\PYG{n}{m}\PYG{o}{.}\PYG{n}{write}\PYG{p}{(}\PYG{l+s+s1}{\PYGZsq{}}\PYG{l+s+s1}{testmip.mst}\PYG{l+s+s1}{\PYGZsq{}}\PYG{p}{)}
\PYG{c+c1}{\PYGZsh{} Write parameter file}
\PYG{n}{m}\PYG{o}{.}\PYG{n}{write}\PYG{p}{(}\PYG{l+s+s1}{\PYGZsq{}}\PYG{l+s+s1}{test.par}\PYG{l+s+s1}{\PYGZsq{}}\PYG{p}{)}
\end{sphinxVerbatim}

\subsubsection{Model.writeMps()}
\label{\detokenize{pyapiref:model-writemps}}\begin{quote}

\sphinxAtStartPar
\sphinxstylestrong{Synopsis}
\begin{quote}

\sphinxAtStartPar
\sphinxcode{\sphinxupquote{writeMps(filename)}}
\end{quote}

\sphinxAtStartPar
\sphinxstylestrong{Description}
\begin{quote}

\sphinxAtStartPar
Write current model into an MPS file.
\end{quote}

\sphinxAtStartPar
\sphinxstylestrong{Arguments}
\begin{quote}

\sphinxAtStartPar
\sphinxcode{\sphinxupquote{filename}}
\begin{quote}

\sphinxAtStartPar
The name of the MPS file to be written.
\end{quote}
\end{quote}

\sphinxAtStartPar
\sphinxstylestrong{Example}
\end{quote}

\begin{sphinxVerbatim}[commandchars=\\\{\}]
\PYG{c+c1}{\PYGZsh{} Write MPS model file \PYGZdq{}test.mps\PYGZdq{}}
\PYG{n}{m}\PYG{o}{.}\PYG{n}{writeMps}\PYG{p}{(}\PYG{l+s+s1}{\PYGZsq{}}\PYG{l+s+s1}{test.mps}\PYG{l+s+s1}{\PYGZsq{}}\PYG{p}{)}
\end{sphinxVerbatim}

\subsubsection{Model.writeMpsStr()}
\label{\detokenize{pyapiref:model-writempsstr}}\begin{quote}

\sphinxAtStartPar
\sphinxstylestrong{Synopsis}
\begin{quote}

\sphinxAtStartPar
\sphinxcode{\sphinxupquote{writeMpsStr()}}
\end{quote}

\sphinxAtStartPar
\sphinxstylestrong{Description}
\begin{quote}

\sphinxAtStartPar
Write current model into a buffer as MPS format.
\end{quote}

\sphinxAtStartPar
\sphinxstylestrong{Example}
\end{quote}

\begin{sphinxVerbatim}[commandchars=\\\{\}]
\PYG{c+c1}{\PYGZsh{} Write model to buffer \PYGZsq{}buff\PYGZsq{} and print model}
\PYG{n}{buff} \PYG{o}{=} \PYG{n}{m}\PYG{o}{.}\PYG{n}{writeMpsStr}\PYG{p}{(}\PYG{p}{)}
\PYG{n+nb}{print}\PYG{p}{(}\PYG{n}{buff}\PYG{o}{.}\PYG{n}{getData}\PYG{p}{(}\PYG{p}{)}\PYG{p}{)}
\end{sphinxVerbatim}

\subsubsection{Model.writeLp()}
\label{\detokenize{pyapiref:model-writelp}}\begin{quote}

\sphinxAtStartPar
\sphinxstylestrong{Synopsis}
\begin{quote}

\sphinxAtStartPar
\sphinxcode{\sphinxupquote{writeLp(filename)}}
\end{quote}

\sphinxAtStartPar
\sphinxstylestrong{Description}
\begin{quote}

\sphinxAtStartPar
Write current optimization model to a LP file.
\end{quote}

\sphinxAtStartPar
\sphinxstylestrong{Arguments}
\begin{quote}

\sphinxAtStartPar
\sphinxcode{\sphinxupquote{filename}}
\begin{quote}

\sphinxAtStartPar
The name of the LP file to be written.
\end{quote}
\end{quote}

\sphinxAtStartPar
\sphinxstylestrong{Example}
\end{quote}

\begin{sphinxVerbatim}[commandchars=\\\{\}]
\PYG{c+c1}{\PYGZsh{} Write LP model file \PYGZdq{}test.lp\PYGZdq{}}
\PYG{n}{m}\PYG{o}{.}\PYG{n}{writeLp}\PYG{p}{(}\PYG{l+s+s1}{\PYGZsq{}}\PYG{l+s+s1}{test.lp}\PYG{l+s+s1}{\PYGZsq{}}\PYG{p}{)}
\end{sphinxVerbatim}

\subsubsection{Model.writeCbf()}
\label{\detokenize{pyapiref:model-writecbf}}\begin{quote}

\sphinxAtStartPar
\sphinxstylestrong{Synopsis}
\begin{quote}

\sphinxAtStartPar
\sphinxcode{\sphinxupquote{writeCbf(filename)}}
\end{quote}

\sphinxAtStartPar
\sphinxstylestrong{Description}
\begin{quote}

\sphinxAtStartPar
Write current optimization model to a CBF file.
\end{quote}

\sphinxAtStartPar
\sphinxstylestrong{Arguments}
\begin{quote}

\sphinxAtStartPar
\sphinxcode{\sphinxupquote{filename}}
\begin{quote}

\sphinxAtStartPar
The name of the CBF file to be written.
\end{quote}
\end{quote}

\sphinxAtStartPar
\sphinxstylestrong{Example}
\end{quote}

\begin{sphinxVerbatim}[commandchars=\\\{\}]
\PYG{c+c1}{\PYGZsh{} Write CBF model file \PYGZdq{}test.cbf\PYGZdq{}}
\PYG{n}{m}\PYG{o}{.}\PYG{n}{writeCbf}\PYG{p}{(}\PYG{l+s+s1}{\PYGZsq{}}\PYG{l+s+s1}{test.cbf}\PYG{l+s+s1}{\PYGZsq{}}\PYG{p}{)}
\end{sphinxVerbatim}

\subsubsection{Model.writeBin()}
\label{\detokenize{pyapiref:model-writebin}}\begin{quote}

\sphinxAtStartPar
\sphinxstylestrong{Synopsis}
\begin{quote}

\sphinxAtStartPar
\sphinxcode{\sphinxupquote{writeBin(filename)}}
\end{quote}

\sphinxAtStartPar
\sphinxstylestrong{Description}
\begin{quote}

\sphinxAtStartPar
Write current model into an COPT binary format file.
\end{quote}

\sphinxAtStartPar
\sphinxstylestrong{Arguments}
\begin{quote}

\sphinxAtStartPar
\sphinxcode{\sphinxupquote{filename}}
\begin{quote}

\sphinxAtStartPar
The name of the COPT binary format file to be written.
\end{quote}
\end{quote}

\sphinxAtStartPar
\sphinxstylestrong{Example}
\end{quote}

\begin{sphinxVerbatim}[commandchars=\\\{\}]
\PYG{c+c1}{\PYGZsh{} Write COPT binary format model file \PYGZdq{}test.bin\PYGZdq{}}
\PYG{n}{m}\PYG{o}{.}\PYG{n}{writeBin}\PYG{p}{(}\PYG{l+s+s1}{\PYGZsq{}}\PYG{l+s+s1}{test.bin}\PYG{l+s+s1}{\PYGZsq{}}\PYG{p}{)}
\end{sphinxVerbatim}

\subsubsection{Model.writeIIS()}
\label{\detokenize{pyapiref:model-writeiis}}\begin{quote}

\sphinxAtStartPar
\sphinxstylestrong{Synopsis}
\begin{quote}

\sphinxAtStartPar
\sphinxcode{\sphinxupquote{writeIIS(filename)}}
\end{quote}

\sphinxAtStartPar
\sphinxstylestrong{Description}
\begin{quote}

\sphinxAtStartPar
Write current irreducible inconsistent subsystem into an IIS file.
\end{quote}

\sphinxAtStartPar
\sphinxstylestrong{Arguments}
\begin{quote}

\sphinxAtStartPar
\sphinxcode{\sphinxupquote{filename}}
\begin{quote}

\sphinxAtStartPar
The name of the IIS file to be written.
\end{quote}
\end{quote}

\sphinxAtStartPar
\sphinxstylestrong{Example}
\end{quote}

\begin{sphinxVerbatim}[commandchars=\\\{\}]
\PYG{c+c1}{\PYGZsh{} Write IIS file \PYGZdq{}test.iis\PYGZdq{}}
\PYG{n}{m}\PYG{o}{.}\PYG{n}{writeIIS}\PYG{p}{(}\PYG{l+s+s1}{\PYGZsq{}}\PYG{l+s+s1}{test.iis}\PYG{l+s+s1}{\PYGZsq{}}\PYG{p}{)}
\end{sphinxVerbatim}

\subsubsection{Model.writeRelax()}
\label{\detokenize{pyapiref:model-writerelax}}\begin{quote}

\sphinxAtStartPar
\sphinxstylestrong{Synopsis}
\begin{quote}

\sphinxAtStartPar
\sphinxcode{\sphinxupquote{writeRelax(filename)}}
\end{quote}

\sphinxAtStartPar
\sphinxstylestrong{Description}
\begin{quote}

\sphinxAtStartPar
Write the feasibility relaxation model into a Relax file.
\end{quote}

\sphinxAtStartPar
\sphinxstylestrong{Arguments}
\begin{quote}

\sphinxAtStartPar
\sphinxcode{\sphinxupquote{filename}}
\begin{quote}

\sphinxAtStartPar
The name of the Relax file to be written.
\end{quote}
\end{quote}

\sphinxAtStartPar
\sphinxstylestrong{Example}
\end{quote}

\begin{sphinxVerbatim}[commandchars=\\\{\}]
\PYG{c+c1}{\PYGZsh{} Write Relax file \PYGZdq{}test.relax\PYGZdq{}}
\PYG{n}{m}\PYG{o}{.}\PYG{n}{writeRelax}\PYG{p}{(}\PYG{l+s+s1}{\PYGZsq{}}\PYG{l+s+s1}{test.relax}\PYG{l+s+s1}{\PYGZsq{}}\PYG{p}{)}
\end{sphinxVerbatim}

\subsubsection{Model.writeSol()}
\label{\detokenize{pyapiref:model-writesol}}\begin{quote}

\sphinxAtStartPar
\sphinxstylestrong{Synopsis}
\begin{quote}

\sphinxAtStartPar
\sphinxcode{\sphinxupquote{writeSol(filename)}}
\end{quote}

\sphinxAtStartPar
\sphinxstylestrong{Description}
\begin{quote}

\sphinxAtStartPar
Output the model solution to a solution file.
\end{quote}

\sphinxAtStartPar
\sphinxstylestrong{Arguments}
\begin{quote}

\sphinxAtStartPar
\sphinxcode{\sphinxupquote{filename}}
\begin{quote}

\sphinxAtStartPar
The name of the solution file to be written.
\end{quote}
\end{quote}

\sphinxAtStartPar
\sphinxstylestrong{Example}
\end{quote}

\begin{sphinxVerbatim}[commandchars=\\\{\}]
\PYG{c+c1}{\PYGZsh{} Write solution file \PYGZdq{}test.sol\PYGZdq{}}
\PYG{n}{m}\PYG{o}{.}\PYG{n}{writeSol}\PYG{p}{(}\PYG{l+s+s1}{\PYGZsq{}}\PYG{l+s+s1}{test.sol}\PYG{l+s+s1}{\PYGZsq{}}\PYG{p}{)}
\end{sphinxVerbatim}

\subsubsection{Model.writeJsonSol()}
\label{\detokenize{pyapiref:model-writejsonsol}}\begin{quote}

\sphinxAtStartPar
\sphinxstylestrong{Synopsis}
\begin{quote}

\sphinxAtStartPar
\sphinxcode{\sphinxupquote{writeJsonSol(filename)}}
\end{quote}

\sphinxAtStartPar
\sphinxstylestrong{Description}
\begin{quote}

\sphinxAtStartPar
Write the solution of the model to a \sphinxcode{\sphinxupquote{".json"}} format result file.
\end{quote}

\sphinxAtStartPar
\sphinxstylestrong{Arguments}
\begin{quote}

\sphinxAtStartPar
\sphinxcode{\sphinxupquote{filename}}
\begin{quote}

\sphinxAtStartPar
The name of the solution file to be written.
\end{quote}
\end{quote}
\end{quote}

\subsubsection{Model.writePoolSol()}
\label{\detokenize{pyapiref:model-writepoolsol}}\begin{quote}

\sphinxAtStartPar
\sphinxstylestrong{Synopsis}
\begin{quote}

\sphinxAtStartPar
\sphinxcode{\sphinxupquote{writePoolSol(isol, filename)}}
\end{quote}

\sphinxAtStartPar
\sphinxstylestrong{Description}
\begin{quote}

\sphinxAtStartPar
Output selected pool solution to a solution file.
\end{quote}

\sphinxAtStartPar
\sphinxstylestrong{Arguments}
\begin{quote}

\sphinxAtStartPar
\sphinxcode{\sphinxupquote{isol}}
\begin{quote}

\sphinxAtStartPar
Index of pool solution.
\end{quote}

\sphinxAtStartPar
\sphinxcode{\sphinxupquote{filename}}
\begin{quote}

\sphinxAtStartPar
The name of the solution file to be written.
\end{quote}
\end{quote}

\sphinxAtStartPar
\sphinxstylestrong{Example}
\end{quote}

\begin{sphinxVerbatim}[commandchars=\\\{\}]
\PYG{c+c1}{\PYGZsh{} Write 1\PYGZhy{}th pool solution to solution file \PYGZdq{}poolsol\PYGZus{}1.sol\PYGZdq{}}
\PYG{n}{m}\PYG{o}{.}\PYG{n}{writePoolSol}\PYG{p}{(}\PYG{l+s+s1}{\PYGZsq{}}\PYG{l+s+s1}{poolsol\PYGZus{}1.sol}\PYG{l+s+s1}{\PYGZsq{}}\PYG{p}{)}
\end{sphinxVerbatim}

\subsubsection{Model.writeBasis()}
\label{\detokenize{pyapiref:model-writebasis}}\begin{quote}

\sphinxAtStartPar
\sphinxstylestrong{Synopsis}
\begin{quote}

\sphinxAtStartPar
\sphinxcode{\sphinxupquote{writeBasis(filename)}}
\end{quote}

\sphinxAtStartPar
\sphinxstylestrong{Description}
\begin{quote}

\sphinxAtStartPar
Write the LP basic solution to a basis file.
\end{quote}

\sphinxAtStartPar
\sphinxstylestrong{Arguments}
\begin{quote}

\sphinxAtStartPar
\sphinxcode{\sphinxupquote{filename}}
\begin{quote}

\sphinxAtStartPar
The name of the basis file to be written.
\end{quote}
\end{quote}

\sphinxAtStartPar
\sphinxstylestrong{Example}
\end{quote}

\begin{sphinxVerbatim}[commandchars=\\\{\}]
\PYG{c+c1}{\PYGZsh{} Write the basis file \PYGZdq{}testlp.bas\PYGZdq{}}
\PYG{n}{m}\PYG{o}{.}\PYG{n}{writeBasis}\PYG{p}{(}\PYG{l+s+s1}{\PYGZsq{}}\PYG{l+s+s1}{testlp.bas}\PYG{l+s+s1}{\PYGZsq{}}\PYG{p}{)}
\end{sphinxVerbatim}

\subsubsection{Model.writeMst()}
\label{\detokenize{pyapiref:model-writemst}}\begin{quote}

\sphinxAtStartPar
\sphinxstylestrong{Synopsis}
\begin{quote}

\sphinxAtStartPar
\sphinxcode{\sphinxupquote{writeMst(filename)}}
\end{quote}
\begin{description}
\sphinxlineitem{\sphinxstylestrong{Description}}
\sphinxAtStartPar
For integer programming models, write the best integer solution currently to the initial solution file.
If there are no integer solutions, then the first set of initial solution stored in model is output.

\end{description}

\sphinxAtStartPar
\sphinxstylestrong{Arguments}
\begin{quote}

\sphinxAtStartPar
\sphinxcode{\sphinxupquote{filename}}
\begin{quote}

\sphinxAtStartPar
Name of the file to be written.
\end{quote}
\end{quote}

\sphinxAtStartPar
\sphinxstylestrong{Example}
\end{quote}

\begin{sphinxVerbatim}[commandchars=\\\{\}]
\PYG{c+c1}{\PYGZsh{} Output initial solution file \PYGZdq{}testmip.mst\PYGZdq{}}
\PYG{n}{m}\PYG{o}{.}\PYG{n}{writeMst}\PYG{p}{(}\PYG{l+s+s1}{\PYGZsq{}}\PYG{l+s+s1}{testmip.mst}\PYG{l+s+s1}{\PYGZsq{}}\PYG{p}{)}
\end{sphinxVerbatim}

\subsubsection{Model.writeParam()}
\label{\detokenize{pyapiref:model-writeparam}}\begin{quote}

\sphinxAtStartPar
\sphinxstylestrong{Synopsis}
\begin{quote}

\sphinxAtStartPar
\sphinxcode{\sphinxupquote{writeParam(filename)}}
\end{quote}

\sphinxAtStartPar
\sphinxstylestrong{Description}
\begin{quote}

\sphinxAtStartPar
Output modified parameters to a parameter file.
\end{quote}

\sphinxAtStartPar
\sphinxstylestrong{Arguments}
\begin{quote}

\sphinxAtStartPar
\sphinxcode{\sphinxupquote{filename}}
\begin{quote}

\sphinxAtStartPar
The name of the parameter file to be written.
\end{quote}
\end{quote}

\sphinxAtStartPar
\sphinxstylestrong{Example}
\end{quote}

\begin{sphinxVerbatim}[commandchars=\\\{\}]
\PYG{c+c1}{\PYGZsh{} Output parameter file \PYGZdq{}testmip.par\PYGZdq{}}
\PYG{n}{m}\PYG{o}{.}\PYG{n}{writeParam}\PYG{p}{(}\PYG{l+s+s1}{\PYGZsq{}}\PYG{l+s+s1}{testmip.par}\PYG{l+s+s1}{\PYGZsq{}}\PYG{p}{)}
\end{sphinxVerbatim}

\subsubsection{Model.writeTuneParam()}
\label{\detokenize{pyapiref:model-writetuneparam}}\begin{quote}

\sphinxAtStartPar
\sphinxstylestrong{Synopsis}
\begin{quote}

\sphinxAtStartPar
\sphinxcode{\sphinxupquote{writeTuneParam(idx, filename)}}
\end{quote}

\sphinxAtStartPar
\sphinxstylestrong{Description}
\begin{quote}

\sphinxAtStartPar
Output the parameter tuning result of the specified number to the parameter file.
\end{quote}

\sphinxAtStartPar
\sphinxstylestrong{Arguments}
\begin{quote}

\sphinxAtStartPar
\sphinxcode{\sphinxupquote{idx}}
\begin{quote}

\sphinxAtStartPar
Parameter tuning result number.
\end{quote}

\sphinxAtStartPar
\sphinxcode{\sphinxupquote{filename}}
\begin{quote}

\sphinxAtStartPar
The name of the parameter file to be output.
\end{quote}
\end{quote}

\sphinxAtStartPar
\sphinxstylestrong{Example}
\end{quote}

\subsubsection{Model.writeOrd()}
\label{\detokenize{pyapiref:model-writeord}}\begin{quote}

\sphinxAtStartPar
\sphinxstylestrong{Synopsis}
\begin{quote}

\sphinxAtStartPar
\sphinxcode{\sphinxupquote{writeOrd(filename)}}
\end{quote}

\sphinxAtStartPar
\sphinxstylestrong{Description}
\begin{quote}

\sphinxAtStartPar
Write the branching order information of the current model into an ORD\sphinxhyphen{}format file,
so that the branching strategy can be saved and reused.
\end{quote}

\sphinxAtStartPar
\sphinxstylestrong{Arguments}
\begin{quote}

\sphinxAtStartPar
\sphinxcode{\sphinxupquote{filename}}
\begin{quote}

\sphinxAtStartPar
Name of the branching order file to be written.
\end{quote}
\end{quote}

\sphinxAtStartPar
\sphinxstylestrong{Example}

\begin{sphinxVerbatim}[commandchars=\\\{\}]
\PYG{c+c1}{\PYGZsh{} Write branching order information of the current model into \PYGZdq{}testmip.ord\PYGZdq{}}
\PYG{n}{m}\PYG{o}{.}\PYG{n}{writeOrd}\PYG{p}{(}\PYG{l+s+s2}{\PYGZdq{}}\PYG{l+s+s2}{testmip.ord}\PYG{l+s+s2}{\PYGZdq{}}\PYG{p}{)}
\end{sphinxVerbatim}
\end{quote}

\subsubsection{Model.setLogFile()}
\label{\detokenize{pyapiref:model-setlogfile}}\begin{quote}

\sphinxAtStartPar
\sphinxstylestrong{Synopsis}
\begin{quote}

\sphinxAtStartPar
\sphinxcode{\sphinxupquote{setLogFile(logfile)}}
\end{quote}

\sphinxAtStartPar
\sphinxstylestrong{Description}
\begin{quote}

\sphinxAtStartPar
Set the optimizer log file.
\end{quote}

\sphinxAtStartPar
\sphinxstylestrong{Arguments}
\begin{quote}

\sphinxAtStartPar
\sphinxcode{\sphinxupquote{logfile}}
\begin{quote}

\sphinxAtStartPar
The log file.
\end{quote}
\end{quote}

\sphinxAtStartPar
\sphinxstylestrong{Example}
\end{quote}

\begin{sphinxVerbatim}[commandchars=\\\{\}]
\PYG{c+c1}{\PYGZsh{} Set the log file as \PYGZdq{}copt.log\PYGZdq{}}
\PYG{n}{m}\PYG{o}{.}\PYG{n}{setLogFile}\PYG{p}{(}\PYG{l+s+s1}{\PYGZsq{}}\PYG{l+s+s1}{copt.log}\PYG{l+s+s1}{\PYGZsq{}}\PYG{p}{)}
\end{sphinxVerbatim}

\subsubsection{Model.setLogCallback()}
\label{\detokenize{pyapiref:model-setlogcallback}}\begin{quote}

\sphinxAtStartPar
\sphinxstylestrong{Synopsis}
\begin{quote}

\sphinxAtStartPar
\sphinxcode{\sphinxupquote{setLogCallback(logcb)}}
\end{quote}

\sphinxAtStartPar
\sphinxstylestrong{Description}
\begin{quote}

\sphinxAtStartPar
Set the call back function of log.
\end{quote}

\sphinxAtStartPar
\sphinxstylestrong{Arguments}
\begin{quote}

\sphinxAtStartPar
\sphinxcode{\sphinxupquote{logcb}}
\begin{quote}

\sphinxAtStartPar
Call back function of log.
\end{quote}
\end{quote}

\sphinxAtStartPar
\sphinxstylestrong{Example}
\end{quote}

\begin{sphinxVerbatim}[commandchars=\\\{\}]
\PYG{c+c1}{\PYGZsh{} Set the call back function of log as a python function \PYGZsq{}logcbfun\PYGZsq{}.}
\PYG{n}{m}\PYG{o}{.}\PYG{n}{setLogCallback}\PYG{p}{(}\PYG{n}{logcbfun}\PYG{p}{)}
\end{sphinxVerbatim}

\subsubsection{Model.solve()}
\label{\detokenize{pyapiref:model-solve}}\begin{quote}

\sphinxAtStartPar
\sphinxstylestrong{Synopsis}
\begin{quote}

\sphinxAtStartPar
\sphinxcode{\sphinxupquote{solve()}}
\end{quote}

\sphinxAtStartPar
\sphinxstylestrong{Description}
\begin{quote}

\sphinxAtStartPar
Solve an optimization problem.
\end{quote}

\sphinxAtStartPar
\sphinxstylestrong{Example}
\end{quote}

\begin{sphinxVerbatim}[commandchars=\\\{\}]
\PYG{c+c1}{\PYGZsh{} Solve the model.}
\PYG{n}{m}\PYG{o}{.}\PYG{n}{solve}\PYG{p}{(}\PYG{p}{)}
\end{sphinxVerbatim}

\subsubsection{Model.solveLP()}
\label{\detokenize{pyapiref:model-solvelp}}\begin{quote}

\sphinxAtStartPar
\sphinxstylestrong{Synopsis}
\begin{quote}

\sphinxAtStartPar
\sphinxcode{\sphinxupquote{solveLP()}}
\end{quote}

\sphinxAtStartPar
\sphinxstylestrong{Description}
\begin{quote}

\sphinxAtStartPar
Solve LP model.
If the model is integer programming,  then the model will be solved as LP.
\end{quote}

\sphinxAtStartPar
\sphinxstylestrong{Example}
\end{quote}

\begin{sphinxVerbatim}[commandchars=\\\{\}]
\PYG{c+c1}{\PYGZsh{} Solve a model calling LP solver.}
\PYG{n}{m}\PYG{o}{.}\PYG{n}{solveLP}\PYG{p}{(}\PYG{p}{)}
\end{sphinxVerbatim}

\subsubsection{Model.computeIIS()}
\label{\detokenize{pyapiref:model-computeiis}}\begin{quote}

\sphinxAtStartPar
\sphinxstylestrong{Synopsis}
\begin{quote}

\sphinxAtStartPar
\sphinxcode{\sphinxupquote{computeIIS()}}
\end{quote}

\sphinxAtStartPar
\sphinxstylestrong{Description}
\begin{quote}

\sphinxAtStartPar
Compute IIS for infeasible model.
\end{quote}

\sphinxAtStartPar
\sphinxstylestrong{Example}
\end{quote}

\begin{sphinxVerbatim}[commandchars=\\\{\}]
\PYG{c+c1}{\PYGZsh{} Compute IIS for infeasible model.}
\PYG{n}{m}\PYG{o}{.}\PYG{n}{computeIIS}\PYG{p}{(}\PYG{p}{)}
\end{sphinxVerbatim}

\subsubsection{Model.feasRelax()}
\label{\detokenize{pyapiref:model-feasrelax}}\begin{quote}

\sphinxAtStartPar
\sphinxstylestrong{Synopsis}
\begin{quote}

\sphinxAtStartPar
\sphinxcode{\sphinxupquote{feasRelax(vars, lbpen, ubpen, constrs, rhspen, uppen=None)}}
\end{quote}

\sphinxAtStartPar
\sphinxstylestrong{Description}
\begin{quote}

\sphinxAtStartPar
Compute the feasibility relaxation of an infeasible model.
\end{quote}

\sphinxAtStartPar
\sphinxstylestrong{Arguments}
\begin{quote}

\sphinxAtStartPar
\sphinxcode{\sphinxupquote{vars}}
\begin{quote}

\sphinxAtStartPar
The variables to relax.

\sphinxAtStartPar
Possible values are {\hyperref[\detokenize{pyapiref:chappyapi-var}]{\sphinxcrossref{\DUrole{std,std-ref}{Var Class}}}} objects, {\hyperref[\detokenize{pyapiref:chappyapi-vararray}]{\sphinxcrossref{\DUrole{std,std-ref}{VarArray Class}}}} objects,
arrays, lists, dictionaries, or {\hyperref[\detokenize{pyapiref:chappyapi-util-tupledict}]{\sphinxcrossref{\DUrole{std,std-ref}{tupledict Class}}}} objects.

\sphinxAtStartPar
If \sphinxcode{\sphinxupquote{vars}} is a {\hyperref[\detokenize{pyapiref:chappyapi-var}]{\sphinxcrossref{\DUrole{std,std-ref}{Var Class}}}} object,
then \sphinxcode{\sphinxupquote{lbpen}} / \sphinxcode{\sphinxupquote{ubpen}} can be constant.

\sphinxAtStartPar
If \sphinxcode{\sphinxupquote{vars}} is a {\hyperref[\detokenize{pyapiref:chappyapi-vararray}]{\sphinxcrossref{\DUrole{std,std-ref}{VarArray Class}}}} object, an array or a list.
then \sphinxcode{\sphinxupquote{lbpen}} / \sphinxcode{\sphinxupquote{ubpen}} can correspondingly be array, list or constant.

\sphinxAtStartPar
If \sphinxcode{\sphinxupquote{vars}} is a {\hyperref[\detokenize{pyapiref:chappyapi-util-tupledict}]{\sphinxcrossref{\DUrole{std,std-ref}{tupledict Class}}}} object or a dictionary,
then \sphinxcode{\sphinxupquote{lbpen}} / \sphinxcode{\sphinxupquote{ubpen}} can be dictionary or constant.
\end{quote}

\sphinxAtStartPar
\sphinxcode{\sphinxupquote{lbpen}}
\begin{quote}

\sphinxAtStartPar
The penalty factor for the lower bound of the variables.

\sphinxAtStartPar
If \sphinxcode{\sphinxupquote{lbpen}} is \sphinxcode{\sphinxupquote{None}}, no relaxation of the variables’ lower bound is allowed.

\sphinxAtStartPar
If \sphinxcode{\sphinxupquote{lbpen}} is \sphinxcode{\sphinxupquote{COPT.INFINITY}}, the lower bound of the corresponding variable in \sphinxcode{\sphinxupquote{vars}} is not relaxed.
\end{quote}

\sphinxAtStartPar
\sphinxcode{\sphinxupquote{ubpen}}
\begin{quote}

\sphinxAtStartPar
The penalty factor for the upper bound of the variables.

\sphinxAtStartPar
If \sphinxcode{\sphinxupquote{ubpen}} is \sphinxcode{\sphinxupquote{None}}, no relaxation of the variables’ lower bound is allowed.

\sphinxAtStartPar
If \sphinxcode{\sphinxupquote{ubpen}} is \sphinxcode{\sphinxupquote{COPT.INFINITY}}, the upper bound of the corresponding variable in \sphinxcode{\sphinxupquote{vars}} is not relaxed.
\end{quote}

\sphinxAtStartPar
\sphinxcode{\sphinxupquote{constrs}}
\begin{quote}

\sphinxAtStartPar
The constraints to relax.

\sphinxAtStartPar
Possible values are {\hyperref[\detokenize{pyapiref:chappyapi-constraint}]{\sphinxcrossref{\DUrole{std,std-ref}{Constraint Class}}}} objects,
{\hyperref[\detokenize{pyapiref:chappyapi-constrarray}]{\sphinxcrossref{\DUrole{std,std-ref}{ConstrArray Class}}}} objects, arrays, lists, or dictionaries.

\sphinxAtStartPar
If \sphinxcode{\sphinxupquote{constrs}} is a {\hyperref[\detokenize{pyapiref:chappyapi-constraint}]{\sphinxcrossref{\DUrole{std,std-ref}{Constraint Class}}}} object,
then \sphinxcode{\sphinxupquote{rhspen}} / \sphinxcode{\sphinxupquote{uppen}} can be constant.

\sphinxAtStartPar
If \sphinxcode{\sphinxupquote{constrs}} is a {\hyperref[\detokenize{pyapiref:chappyapi-constrarray}]{\sphinxcrossref{\DUrole{std,std-ref}{ConstrArray Class}}}} object, an array, or a list,
then \sphinxcode{\sphinxupquote{rhspen}} / \sphinxcode{\sphinxupquote{uppen}} can correspondingly be array, list, or constant.

\sphinxAtStartPar
If \sphinxcode{\sphinxupquote{constrs}} is a dictionary,
then \sphinxcode{\sphinxupquote{rhspen}} / \sphinxcode{\sphinxupquote{uppen}} can be dictionary or constant.
\end{quote}

\sphinxAtStartPar
\sphinxcode{\sphinxupquote{rhspen}}
\begin{quote}

\sphinxAtStartPar
The penalty factor for the right\sphinxhyphen{}hand side of the constraints.

\sphinxAtStartPar
If \sphinxcode{\sphinxupquote{rhspen}} is \sphinxcode{\sphinxupquote{None}}, no relaxation of the constraint boundaries is applied.

\sphinxAtStartPar
If \sphinxcode{\sphinxupquote{rhspen}} is \sphinxcode{\sphinxupquote{COPT.INFINITY}}, the constraint’s right\sphinxhyphen{}hand side is not relaxed.
\end{quote}

\sphinxAtStartPar
\sphinxcode{\sphinxupquote{uppen}}
\begin{quote}

\sphinxAtStartPar
The penalty factor for the upper bound of bilateral constraints.

\sphinxAtStartPar
If \sphinxcode{\sphinxupquote{uppen}} is \sphinxcode{\sphinxupquote{None}}, the penalty factor is specified by \sphinxcode{\sphinxupquote{rhspen}} .

\sphinxAtStartPar
If \sphinxcode{\sphinxupquote{uppen}} is \sphinxcode{\sphinxupquote{COPT.INFINITY}}, the upper bound of the constraint is not relaxed.
\end{quote}
\end{quote}

\sphinxAtStartPar
\sphinxstylestrong{Example}
\end{quote}

\begin{sphinxVerbatim}[commandchars=\\\{\}]
\PYG{c+c1}{\PYGZsh{} Compute the feasibility relaxation for an infeasible model}
\PYG{n}{x} \PYG{o}{=} \PYG{n}{model}\PYG{o}{.}\PYG{n}{addVars}\PYG{p}{(}\PYG{l+m+mi}{3}\PYG{p}{,} \PYG{n}{ub}\PYG{o}{=}\PYG{l+m+mi}{1}\PYG{p}{)}
\PYG{n}{y} \PYG{o}{=} \PYG{n}{model}\PYG{o}{.}\PYG{n}{addVars}\PYG{p}{(}\PYG{l+m+mi}{3}\PYG{p}{,} \PYG{n}{ub}\PYG{o}{=}\PYG{l+m+mi}{2}\PYG{p}{)}
\PYG{n}{c} \PYG{o}{=} \PYG{n}{model}\PYG{o}{.}\PYG{n}{addConstrs}\PYG{p}{(}\PYG{n}{x}\PYG{p}{[}\PYG{n}{i}\PYG{p}{]} \PYG{o}{+} \PYG{n}{y}\PYG{p}{[}\PYG{n}{i}\PYG{p}{]} \PYG{o}{\PYGZlt{}}\PYG{o}{=} \PYG{l+m+mi}{1} \PYG{k}{for} \PYG{n}{i} \PYG{o+ow}{in} \PYG{n+nb}{range}\PYG{p}{(}\PYG{l+m+mi}{3}\PYG{p}{)}\PYG{p}{)}
\PYG{n}{model}\PYG{o}{.}\PYG{n}{feasRelax}\PYG{p}{(}\PYG{p}{[}\PYG{n}{x}\PYG{p}{[}\PYG{l+m+mi}{0}\PYG{p}{]}\PYG{p}{,} \PYG{n}{x}\PYG{p}{[}\PYG{l+m+mi}{1}\PYG{p}{]}\PYG{p}{]}\PYG{p}{,} \PYG{p}{[}\PYG{l+m+mf}{0.5}\PYG{p}{,} \PYG{l+m+mi}{1}\PYG{p}{]}\PYG{p}{,} \PYG{p}{[}\PYG{n}{COPT}\PYG{o}{.}\PYG{n}{INFINITY}\PYG{p}{,} \PYG{l+m+mi}{2}\PYG{p}{]}\PYG{p}{,} \PYG{n}{c}\PYG{p}{,} \PYG{l+m+mi}{1}\PYG{p}{)}
\end{sphinxVerbatim}

\subsubsection{Model.feasRelaxS()}
\label{\detokenize{pyapiref:model-feasrelaxs}}\begin{quote}

\sphinxAtStartPar
\sphinxstylestrong{Synopsis}
\begin{quote}

\sphinxAtStartPar
\sphinxcode{\sphinxupquote{feasRelaxS(vrelax, crelax)}}
\end{quote}

\sphinxAtStartPar
\sphinxstylestrong{Description}
\begin{quote}

\sphinxAtStartPar
Compute the feasibilty relaxation of an infeasible model in a simplified moode.
\end{quote}

\sphinxAtStartPar
\sphinxstylestrong{Arguments}
\begin{quote}

\sphinxAtStartPar
\sphinxcode{\sphinxupquote{vrelax}}
\begin{quote}

\sphinxAtStartPar
Whether to relax variables.
\end{quote}

\sphinxAtStartPar
\sphinxcode{\sphinxupquote{crelax}}
\begin{quote}

\sphinxAtStartPar
Whether to relax constraints.
\end{quote}
\end{quote}

\sphinxAtStartPar
\sphinxstylestrong{Example}
\end{quote}

\begin{sphinxVerbatim}[commandchars=\\\{\}]
\PYG{c+c1}{\PYGZsh{} Compute the feasibilty relaxation of model m}
\PYG{n}{m}\PYG{o}{.}\PYG{n}{feasRelaxS}\PYG{p}{(}\PYG{k+kc}{True}\PYG{p}{,} \PYG{k+kc}{True}\PYG{p}{)}
\end{sphinxVerbatim}

\subsubsection{Model.feasRelaxN()}
\label{\detokenize{pyapiref:model-feasrelaxn}}\begin{quote}

\sphinxAtStartPar
\sphinxstylestrong{Synopsis}
\begin{quote}

\sphinxAtStartPar
\sphinxcode{\sphinxupquote{feasRelaxN(vars, lbpen, ubpen, constrs, rhspen, uppen=None)}}
\end{quote}

\sphinxAtStartPar
\sphinxstylestrong{Description}
\begin{quote}

\sphinxAtStartPar
Computes the feasibility relaxation of an infeasible model formulated by
matrix modeling.

\sphinxAtStartPar
Either \sphinxcode{\sphinxupquote{vars}} or \sphinxcode{\sphinxupquote{constrs}} must be a matrix modeling object;
that is, at least one of the following conditions must hold:
\begin{itemize}
\item {} 
\sphinxAtStartPar
\sphinxcode{\sphinxupquote{vars}} is a {\hyperref[\detokenize{pyapiref:chappyapi-mvar}]{\sphinxcrossref{\DUrole{std,std-ref}{MVar Class}}}} object.

\item {} 
\sphinxAtStartPar
\sphinxcode{\sphinxupquote{constrs}} is a {\hyperref[\detokenize{pyapiref:chappyapi-mconstr}]{\sphinxcrossref{\DUrole{std,std-ref}{MConstr Class}}}} object.

\end{itemize}
\end{quote}

\sphinxAtStartPar
\sphinxstylestrong{Arguments}
\begin{quote}

\sphinxAtStartPar
\sphinxcode{\sphinxupquote{vars}}
\begin{quote}

\sphinxAtStartPar
The variables to relax.

\sphinxAtStartPar
Possible values are {\hyperref[\detokenize{pyapiref:chappyapi-mvar}]{\sphinxcrossref{\DUrole{std,std-ref}{MVar Class}}}} objects, {\hyperref[\detokenize{pyapiref:chappyapi-var}]{\sphinxcrossref{\DUrole{std,std-ref}{Var Class}}}} objects,
{\hyperref[\detokenize{pyapiref:chappyapi-vararray}]{\sphinxcrossref{\DUrole{std,std-ref}{VarArray Class}}}} objects, arrays, lists, dictionaries,
or {\hyperref[\detokenize{pyapiref:chappyapi-util-tupledict}]{\sphinxcrossref{\DUrole{std,std-ref}{tupledict Class}}}} objects.

\sphinxAtStartPar
If \sphinxcode{\sphinxupquote{vars}} is a {\hyperref[\detokenize{pyapiref:chappyapi-mvar}]{\sphinxcrossref{\DUrole{std,std-ref}{MVar Class}}}} object,
then \sphinxcode{\sphinxupquote{lbpen}} / \sphinxcode{\sphinxupquote{ubpen}} can be a {\hyperref[\detokenize{pyapiref:chappyapi-ndarray}]{\sphinxcrossref{\DUrole{std,std-ref}{NdArray Class}}}} object,
array, list, dictionary, or constant.

\sphinxAtStartPar
If \sphinxcode{\sphinxupquote{vars}} is a {\hyperref[\detokenize{pyapiref:chappyapi-var}]{\sphinxcrossref{\DUrole{std,std-ref}{Var Class}}}} object,
then \sphinxcode{\sphinxupquote{lbpen}} / \sphinxcode{\sphinxupquote{ubpen}} can be constant.

\sphinxAtStartPar
If \sphinxcode{\sphinxupquote{vars}} is a {\hyperref[\detokenize{pyapiref:chappyapi-vararray}]{\sphinxcrossref{\DUrole{std,std-ref}{VarArray Class}}}} object, an array, or a list,
then \sphinxcode{\sphinxupquote{lbpen}} / \sphinxcode{\sphinxupquote{ubpen}} can be array, list, or constant.

\sphinxAtStartPar
If \sphinxcode{\sphinxupquote{vars}} is a {\hyperref[\detokenize{pyapiref:chappyapi-util-tupledict}]{\sphinxcrossref{\DUrole{std,std-ref}{tupledict Class}}}} object or dictionary,
then \sphinxcode{\sphinxupquote{lbpen}} / \sphinxcode{\sphinxupquote{ubpen}} can be dictionary or constant.
\end{quote}

\sphinxAtStartPar
\sphinxcode{\sphinxupquote{lbpen}}
\begin{quote}

\sphinxAtStartPar
The penalty factor for the lower bound of the variables.

\sphinxAtStartPar
If \sphinxcode{\sphinxupquote{lbpen}} is \sphinxcode{\sphinxupquote{None}}, no relaxation of the variables’ lower bound is allowed.

\sphinxAtStartPar
If \sphinxcode{\sphinxupquote{lbpen}} is \sphinxcode{\sphinxupquote{COPT.INFINITY}}, the lower bound of the corresponding
variable in \sphinxcode{\sphinxupquote{vars}} is not relaxed.
\end{quote}

\sphinxAtStartPar
\sphinxcode{\sphinxupquote{ubpen}}
\begin{quote}

\sphinxAtStartPar
The penalty factor for the upper bound of the variables.

\sphinxAtStartPar
If \sphinxcode{\sphinxupquote{ubpen}} is \sphinxcode{\sphinxupquote{None}}, no relaxation of the variables’ lower bound is allowed.

\sphinxAtStartPar
If \sphinxcode{\sphinxupquote{ubpen}} is \sphinxcode{\sphinxupquote{COPT.INFINITY}}, the upper bound of the corresponding
variable in \sphinxcode{\sphinxupquote{vars}} is not relaxed.
\end{quote}

\sphinxAtStartPar
\sphinxcode{\sphinxupquote{constrs}}
\begin{quote}

\sphinxAtStartPar
The constraints to relax.

\sphinxAtStartPar
Possible values are {\hyperref[\detokenize{pyapiref:chappyapi-mconstr}]{\sphinxcrossref{\DUrole{std,std-ref}{MConstr Class}}}} objects,
{\hyperref[\detokenize{pyapiref:chappyapi-constraint}]{\sphinxcrossref{\DUrole{std,std-ref}{Constraint Class}}}} objects, {\hyperref[\detokenize{pyapiref:chappyapi-constrarray}]{\sphinxcrossref{\DUrole{std,std-ref}{ConstrArray Class}}}} objects,
arrays, lists, or dictionaries.

\sphinxAtStartPar
If \sphinxcode{\sphinxupquote{constrs}} is a {\hyperref[\detokenize{pyapiref:chappyapi-mconstr}]{\sphinxcrossref{\DUrole{std,std-ref}{MConstr Class}}}} object,
then \sphinxcode{\sphinxupquote{rhspen}} / \sphinxcode{\sphinxupquote{uppen}} can be a {\hyperref[\detokenize{pyapiref:chappyapi-ndarray}]{\sphinxcrossref{\DUrole{std,std-ref}{NdArray Class}}}} object,
array, list, dictionary, or constant.

\sphinxAtStartPar
If \sphinxcode{\sphinxupquote{constrs}} is a {\hyperref[\detokenize{pyapiref:chappyapi-constraint}]{\sphinxcrossref{\DUrole{std,std-ref}{Constraint Class}}}} object,
then \sphinxcode{\sphinxupquote{rhspen}} / \sphinxcode{\sphinxupquote{uppen}} can be constant.

\sphinxAtStartPar
If \sphinxcode{\sphinxupquote{constrs}} is a {\hyperref[\detokenize{pyapiref:chappyapi-constrarray}]{\sphinxcrossref{\DUrole{std,std-ref}{ConstrArray Class}}}} object, an array, or a list,
then \sphinxcode{\sphinxupquote{rhspen}} / \sphinxcode{\sphinxupquote{uppen}} can correspondingly be array, list, or constant.

\sphinxAtStartPar
If \sphinxcode{\sphinxupquote{constrs}} is a dictionary,
then \sphinxcode{\sphinxupquote{rhspen}} / \sphinxcode{\sphinxupquote{uppen}} can be dictionary or constant.
\end{quote}

\sphinxAtStartPar
\sphinxcode{\sphinxupquote{rhspen}}
\begin{quote}

\sphinxAtStartPar
The penalty factor for the right\sphinxhyphen{}hand side of the constraints.

\sphinxAtStartPar
If \sphinxcode{\sphinxupquote{rhspen}} is \sphinxcode{\sphinxupquote{None}}, no relaxation of the constraint boundaries is applied.

\sphinxAtStartPar
If \sphinxcode{\sphinxupquote{rhspen}} is \sphinxcode{\sphinxupquote{COPT.INFINITY}}, the constraint’s right\sphinxhyphen{}hand side is not relaxed.
\end{quote}

\sphinxAtStartPar
\sphinxcode{\sphinxupquote{uppen}}
\begin{quote}

\sphinxAtStartPar
The penalty factor for the upper bound of bilateral constraints.

\sphinxAtStartPar
If \sphinxcode{\sphinxupquote{uppen}} is \sphinxcode{\sphinxupquote{None}}, the penalty factor is specified by \sphinxcode{\sphinxupquote{rhspen}}.

\sphinxAtStartPar
If \sphinxcode{\sphinxupquote{uppen}} is \sphinxcode{\sphinxupquote{COPT.INFINITY}}, the upper bound of the constraint is not relaxed.
\end{quote}
\end{quote}

\sphinxAtStartPar
\sphinxstylestrong{Example}
\end{quote}

\begin{sphinxVerbatim}[commandchars=\\\{\}]
\PYG{c+c1}{\PYGZsh{} Compute the feasibility relaxation for an infeasible model formulated using matrix modeling}
\PYG{n}{mx} \PYG{o}{=} \PYG{n}{model}\PYG{o}{.}\PYG{n}{addMVar}\PYG{p}{(}\PYG{n}{shape}\PYG{o}{=}\PYG{l+m+mi}{3}\PYG{p}{,} \PYG{n}{nameprefix}\PYG{o}{=}\PYG{l+s+s2}{\PYGZdq{}}\PYG{l+s+s2}{x}\PYG{l+s+s2}{\PYGZdq{}}\PYG{p}{)}
\PYG{c+c1}{\PYGZsh{} Define constraints}
\PYG{n}{A} \PYG{o}{=} \PYG{n}{np}\PYG{o}{.}\PYG{n}{array}\PYG{p}{(}\PYG{p}{[}\PYG{p}{[}\PYG{l+m+mi}{1}\PYG{p}{,} \PYG{l+m+mi}{2}\PYG{p}{,} \PYG{l+m+mi}{3}\PYG{p}{]}\PYG{p}{,} \PYG{p}{[}\PYG{l+m+mi}{4}\PYG{p}{,} \PYG{l+m+mi}{5}\PYG{p}{,} \PYG{l+m+mi}{6}\PYG{p}{]}\PYG{p}{]}\PYG{p}{)}
\PYG{n}{b} \PYG{o}{=} \PYG{n}{np}\PYG{o}{.}\PYG{n}{array}\PYG{p}{(}\PYG{p}{[}\PYG{l+m+mi}{10}\PYG{p}{,} \PYG{l+m+mi}{20}\PYG{p}{]}\PYG{p}{)}
\PYG{n}{mc} \PYG{o}{=} \PYG{n}{model}\PYG{o}{.}\PYG{n}{addConstr}\PYG{p}{(}\PYG{n}{A} \PYG{o}{@} \PYG{n}{x} \PYG{o}{\PYGZlt{}}\PYG{o}{=} \PYG{n}{b}\PYG{p}{,} \PYG{n}{name}\PYG{o}{=}\PYG{l+s+s2}{\PYGZdq{}}\PYG{l+s+s2}{constrs}\PYG{l+s+s2}{\PYGZdq{}}\PYG{p}{)}
\PYG{c+c1}{\PYGZsh{} Define relaxation penalty factors}
\PYG{c+c1}{\PYGZsh{} The lower bound of the 2nd element of mx is not relaxed}
\PYG{n}{lbpen} \PYG{o}{=} \PYG{p}{[}\PYG{l+m+mf}{1.0}\PYG{p}{,} \PYG{n}{COPT}\PYG{o}{.}\PYG{n}{INFINITY}\PYG{p}{,} \PYG{l+m+mf}{0.5}\PYG{p}{]}
\PYG{c+c1}{\PYGZsh{} The upper bound of the 3rd element of mx is not relaxed}
\PYG{n}{ubpen} \PYG{o}{=} \PYG{n}{NdArray}\PYG{p}{(}\PYG{p}{[}\PYG{l+m+mf}{0.8}\PYG{p}{,} \PYG{l+m+mf}{1.5}\PYG{p}{,} \PYG{n}{COPT}\PYG{o}{.}\PYG{n}{INFINITY}\PYG{p}{]}\PYG{p}{)}
\PYG{c+c1}{\PYGZsh{} Relaxation penalty factor for mc}
\PYG{n}{rhspen} \PYG{o}{=} \PYG{l+m+mf}{2.0}
\PYG{c+c1}{\PYGZsh{} Compute feasibility relaxation for mx and mc in matrix modeling}
\PYG{n}{model}\PYG{o}{.}\PYG{n}{feasRelaxN}\PYG{p}{(}\PYG{n}{mx}\PYG{p}{,} \PYG{n}{lbpen}\PYG{p}{,} \PYG{n}{ubpen}\PYG{p}{,} \PYG{n}{mc}\PYG{p}{,} \PYG{n}{rhspen}\PYG{p}{)}
\end{sphinxVerbatim}

\subsubsection{Model.tune()}
\label{\detokenize{pyapiref:model-tune}}\begin{quote}

\sphinxAtStartPar
\sphinxstylestrong{Synopsis}
\begin{quote}

\sphinxAtStartPar
\sphinxcode{\sphinxupquote{tune()}}
\end{quote}

\sphinxAtStartPar
\sphinxstylestrong{Description}
\begin{quote}

\sphinxAtStartPar
Parameter tuning of the model.
\end{quote}

\sphinxAtStartPar
\sphinxstylestrong{Example}
\end{quote}

\subsubsection{Model.loadTuneParam()}
\label{\detokenize{pyapiref:model-loadtuneparam}}\begin{quote}

\sphinxAtStartPar
\sphinxstylestrong{Synopsis}
\begin{quote}

\sphinxAtStartPar
\sphinxcode{\sphinxupquote{loadTuneParam(idx)}}
\end{quote}

\sphinxAtStartPar
\sphinxstylestrong{Description}
\begin{quote}

\sphinxAtStartPar
Load the parameter tuning results of the specified number into the model.
\end{quote}

\sphinxAtStartPar
\sphinxstylestrong{Example}
\end{quote}

\subsubsection{Model.interrupt()}
\label{\detokenize{pyapiref:model-interrupt}}\begin{quote}

\sphinxAtStartPar
\sphinxstylestrong{Synopsis}
\begin{quote}

\sphinxAtStartPar
\sphinxcode{\sphinxupquote{interrupt()}}
\end{quote}

\sphinxAtStartPar
\sphinxstylestrong{Description}
\begin{quote}

\sphinxAtStartPar
Interrupt solving process of current problem.
\end{quote}

\sphinxAtStartPar
\sphinxstylestrong{Example}
\end{quote}

\begin{sphinxVerbatim}[commandchars=\\\{\}]
\PYG{c+c1}{\PYGZsh{} Interrupt the solving process.}
\PYG{n}{m}\PYG{o}{.}\PYG{n}{interrupt}\PYG{p}{(}\PYG{p}{)}
\end{sphinxVerbatim}

\subsubsection{Model.remove()}
\label{\detokenize{pyapiref:model-remove}}\begin{quote}

\sphinxAtStartPar
\sphinxstylestrong{Synopsis}
\begin{quote}

\sphinxAtStartPar
\sphinxcode{\sphinxupquote{remove(args)}}
\end{quote}

\sphinxAtStartPar
\sphinxstylestrong{Description}
\begin{quote}

\sphinxAtStartPar
Remove variables or constraints from a model.

\sphinxAtStartPar
To remove variable, then parameter \sphinxcode{\sphinxupquote{args}} can be {\hyperref[\detokenize{pyapiref:chappyapi-var}]{\sphinxcrossref{\DUrole{std,std-ref}{Var Class}}}} object, {\hyperref[\detokenize{pyapiref:chappyapi-vararray}]{\sphinxcrossref{\DUrole{std,std-ref}{VarArray Class}}}} object,
list, dictionary or {\hyperref[\detokenize{pyapiref:chappyapi-util-tupledict}]{\sphinxcrossref{\DUrole{std,std-ref}{tupledict Class}}}} object.

\sphinxAtStartPar
To remove linear constraint, then parameter \sphinxcode{\sphinxupquote{args}} can be {\hyperref[\detokenize{pyapiref:chappyapi-constraint}]{\sphinxcrossref{\DUrole{std,std-ref}{Constraint Class}}}} object,
{\hyperref[\detokenize{pyapiref:chappyapi-constrarray}]{\sphinxcrossref{\DUrole{std,std-ref}{ConstrArray Class}}}} object, list, dictionary or {\hyperref[\detokenize{pyapiref:chappyapi-util-tupledict}]{\sphinxcrossref{\DUrole{std,std-ref}{tupledict Class}}}} object.

\sphinxAtStartPar
To remove SOS constraint, then parameter \sphinxcode{\sphinxupquote{args}} can be {\hyperref[\detokenize{pyapiref:chappyapi-sos}]{\sphinxcrossref{\DUrole{std,std-ref}{SOS Class}}}} object, {\hyperref[\detokenize{pyapiref:chappyapi-sosarray}]{\sphinxcrossref{\DUrole{std,std-ref}{SOSArray Class}}}}
object, list, dictionary or {\hyperref[\detokenize{pyapiref:chappyapi-util-tupledict}]{\sphinxcrossref{\DUrole{std,std-ref}{tupledict Class}}}} object.

\sphinxAtStartPar
To remove Second\sphinxhyphen{}Order\sphinxhyphen{}Cone constraints, then parameter \sphinxcode{\sphinxupquote{args}} can be {\hyperref[\detokenize{pyapiref:chappyapi-cone}]{\sphinxcrossref{\DUrole{std,std-ref}{Cone Class}}}} object, {\hyperref[\detokenize{pyapiref:chappyapi-conearray}]{\sphinxcrossref{\DUrole{std,std-ref}{ConeArray Class}}}}
object, list, dictionary or {\hyperref[\detokenize{pyapiref:chappyapi-util-tupledict}]{\sphinxcrossref{\DUrole{std,std-ref}{tupledict Class}}}} object.

\sphinxAtStartPar
To remove exponential cone constraints, then parameter \sphinxcode{\sphinxupquote{args}} can be {\hyperref[\detokenize{pyapiref:chappyapi-expcone}]{\sphinxcrossref{\DUrole{std,std-ref}{ExpCone Class}}}} object, {\hyperref[\detokenize{pyapiref:chappyapi-expconearray}]{\sphinxcrossref{\DUrole{std,std-ref}{ExpConeArray Class}}}}
object, list, dictionary or {\hyperref[\detokenize{pyapiref:chappyapi-util-tupledict}]{\sphinxcrossref{\DUrole{std,std-ref}{tupledict Class}}}} object.

\sphinxAtStartPar
To remove affine cone constraints, then parameter \sphinxcode{\sphinxupquote{args}} can be {\hyperref[\detokenize{pyapiref:chappyapi-affinecone}]{\sphinxcrossref{\DUrole{std,std-ref}{AffineCone Class}}}} object, {\hyperref[\detokenize{pyapiref:chappyapi-affineconearray}]{\sphinxcrossref{\DUrole{std,std-ref}{AffineConeArray Class}}}}
object, list, dictionary or {\hyperref[\detokenize{pyapiref:chappyapi-util-tupledict}]{\sphinxcrossref{\DUrole{std,std-ref}{tupledict Class}}}} object.

\sphinxAtStartPar
To remove quadratic constraints, then parameter \sphinxcode{\sphinxupquote{args}} can be {\hyperref[\detokenize{pyapiref:chappyapi-qconstraint}]{\sphinxcrossref{\DUrole{std,std-ref}{QConstraint Class}}}} object, {\hyperref[\detokenize{pyapiref:chappyapi-qconstrarray}]{\sphinxcrossref{\DUrole{std,std-ref}{QConstrArray Class}}}}
object, list, dictionary or {\hyperref[\detokenize{pyapiref:chappyapi-util-tupledict}]{\sphinxcrossref{\DUrole{std,std-ref}{tupledict Class}}}} object.

\sphinxAtStartPar
To remove positive semi\sphinxhyphen{}definite constraints, then parameter \sphinxcode{\sphinxupquote{args}} can be {\hyperref[\detokenize{pyapiref:chappyapi-psdconstraint}]{\sphinxcrossref{\DUrole{std,std-ref}{PsdConstraint Class}}}} object, {\hyperref[\detokenize{pyapiref:chappyapi-psdconstrarray}]{\sphinxcrossref{\DUrole{std,std-ref}{PsdConstrArray Class}}}}
object, list, dictionary or {\hyperref[\detokenize{pyapiref:chappyapi-util-tupledict}]{\sphinxcrossref{\DUrole{std,std-ref}{tupledict Class}}}} object.

\sphinxAtStartPar
To remove indicator constraint, then parameter \sphinxcode{\sphinxupquote{args}} can be {\hyperref[\detokenize{pyapiref:chappyapi-genconstr}]{\sphinxcrossref{\DUrole{std,std-ref}{GenConstr Class}}}} object,
{\hyperref[\detokenize{pyapiref:chappyapi-genconstrarray}]{\sphinxcrossref{\DUrole{std,std-ref}{GenConstrArray Class}}}} object, list, dictionary or {\hyperref[\detokenize{pyapiref:chappyapi-util-tupledict}]{\sphinxcrossref{\DUrole{std,std-ref}{tupledict Class}}}} object.

\sphinxAtStartPar
To remove LMI constraint, then parameter \sphinxcode{\sphinxupquote{args}} can be {\hyperref[\detokenize{pyapiref:chappyapi-lmiconstraint}]{\sphinxcrossref{\DUrole{std,std-ref}{LmiConstraint Class}}}} object,
{\hyperref[\detokenize{pyapiref:chappyapi-lmiconstrarray}]{\sphinxcrossref{\DUrole{std,std-ref}{LmiConstrArray Class}}}} object, list, dictionary or {\hyperref[\detokenize{pyapiref:chappyapi-util-tupledict}]{\sphinxcrossref{\DUrole{std,std-ref}{tupledict Class}}}} object.

\sphinxAtStartPar
To remove matrix variables or matrix constriants, then parameter \sphinxcode{\sphinxupquote{args}} can be {\hyperref[\detokenize{pyapiref:chappyapi-mvar}]{\sphinxcrossref{\DUrole{std,std-ref}{MVar Class}}}} object, {\hyperref[\detokenize{pyapiref:chappyapi-mconstr}]{\sphinxcrossref{\DUrole{std,std-ref}{MConstr Class}}}} object, {\hyperref[\detokenize{pyapiref:chappyapi-mqconstr}]{\sphinxcrossref{\DUrole{std,std-ref}{MQConstr Class}}}} object or {\hyperref[\detokenize{pyapiref:chappyapi-mpsdconstr}]{\sphinxcrossref{\DUrole{std,std-ref}{MPsdConstr Class}}}} object.
\end{quote}

\sphinxAtStartPar
\sphinxstylestrong{Arguments}
\begin{quote}

\sphinxAtStartPar
\sphinxcode{\sphinxupquote{args}}
\begin{quote}

\sphinxAtStartPar
Variables or constraints to be removed.
\end{quote}
\end{quote}

\sphinxAtStartPar
\sphinxstylestrong{Example}
\end{quote}

\begin{sphinxVerbatim}[commandchars=\\\{\}]
\PYG{c+c1}{\PYGZsh{} Remove linear constraint conx}
\PYG{n}{m}\PYG{o}{.}\PYG{n}{remove}\PYG{p}{(}\PYG{n}{conx}\PYG{p}{)}
\PYG{c+c1}{\PYGZsh{} Remove variables x and y}
\PYG{n}{m}\PYG{o}{.}\PYG{n}{remove}\PYG{p}{(}\PYG{p}{[}\PYG{n}{x}\PYG{p}{,} \PYG{n}{y}\PYG{p}{]}\PYG{p}{)}
\end{sphinxVerbatim}

\subsubsection{Model.reset()}
\label{\detokenize{pyapiref:model-reset}}\begin{quote}

\sphinxAtStartPar
\sphinxstylestrong{Synopsis}
\begin{quote}

\sphinxAtStartPar
\sphinxcode{\sphinxupquote{reset()}}
\end{quote}

\sphinxAtStartPar
\sphinxstylestrong{Description}
\begin{quote}

\sphinxAtStartPar
Reset the result information of the model.
\end{quote}

\sphinxAtStartPar
\sphinxstylestrong{Example}
\end{quote}

\begin{sphinxVerbatim}[commandchars=\\\{\}]
\PYG{c+c1}{\PYGZsh{} Reset the result information in model.}
\PYG{n}{m}\PYG{o}{.}\PYG{n}{reset}\PYG{p}{(}\PYG{p}{)}
\end{sphinxVerbatim}

\subsubsection{Model.resetAll()}
\label{\detokenize{pyapiref:model-resetall}}\begin{quote}

\sphinxAtStartPar
\sphinxstylestrong{Synopsis}
\begin{quote}

\sphinxAtStartPar
\sphinxcode{\sphinxupquote{resetAll()}}
\end{quote}

\sphinxAtStartPar
\sphinxstylestrong{Description}
\begin{quote}

\sphinxAtStartPar
Reset the result and other additional information of the model, such as MIP start, IIS, etc.
By executing this function, the information that needs to be calculated of the model will be all cleared, only the original model itself will be kept.
\end{quote}

\sphinxAtStartPar
\sphinxstylestrong{Example}
\end{quote}

\begin{sphinxVerbatim}[commandchars=\\\{\}]
\PYG{c+c1}{\PYGZsh{} Reset the result and other additional information of the model}
\PYG{n}{m}\PYG{o}{.}\PYG{n}{resetAll}\PYG{p}{(}\PYG{p}{)}
\end{sphinxVerbatim}

\subsubsection{Model.clear()}
\label{\detokenize{pyapiref:model-clear}}\begin{quote}

\sphinxAtStartPar
\sphinxstylestrong{Synopsis}
\begin{quote}

\sphinxAtStartPar
\sphinxcode{\sphinxupquote{clear()}}
\end{quote}

\sphinxAtStartPar
\sphinxstylestrong{Description}
\begin{quote}

\sphinxAtStartPar
Clear all content of the whole model.
By executing this function, all content of the whole model will be cleared, including the added variables, objective, and constraints.
\end{quote}

\sphinxAtStartPar
\sphinxstylestrong{Example}
\end{quote}

\begin{sphinxVerbatim}[commandchars=\\\{\}]
\PYG{c+c1}{\PYGZsh{} Clear all content of the model.}
\PYG{n}{m}\PYG{o}{.}\PYG{n}{clear}\PYG{p}{(}\PYG{p}{)}
\end{sphinxVerbatim}

\subsubsection{Model.clone()}
\label{\detokenize{pyapiref:model-clone}}\begin{quote}

\sphinxAtStartPar
\sphinxstylestrong{Synopsis}
\begin{quote}

\sphinxAtStartPar
\sphinxcode{\sphinxupquote{clone()}}
\end{quote}

\sphinxAtStartPar
\sphinxstylestrong{Description}
\begin{quote}

\sphinxAtStartPar
Create a deep copy of an existing model. Return a {\hyperref[\detokenize{pyapiref:chappyapi-model}]{\sphinxcrossref{\DUrole{std,std-ref}{Model Class}}}} object.
\end{quote}

\sphinxAtStartPar
\sphinxstylestrong{Example}
\end{quote}

\begin{sphinxVerbatim}[commandchars=\\\{\}]
\PYG{c+c1}{\PYGZsh{} Create a deep copy of model}
\PYG{n}{mcopy} \PYG{o}{=} \PYG{n}{m}\PYG{o}{.}\PYG{n}{clone}\PYG{p}{(}\PYG{p}{)}
\end{sphinxVerbatim}

\subsubsection{Model.setCallback()}
\label{\detokenize{pyapiref:model-setcallback}}\label{\detokenize{pyapiref:chappyapi-model-setcallback}}\begin{quote}

\sphinxAtStartPar
\sphinxstylestrong{Synopsis}
\begin{quote}

\sphinxAtStartPar
\sphinxcode{\sphinxupquote{setCallback(cb, cbctx)}}
\end{quote}

\sphinxAtStartPar
\sphinxstylestrong{Description}
\begin{quote}

\sphinxAtStartPar
Set the callback function of the model.
\end{quote}

\sphinxAtStartPar
\sphinxstylestrong{Arguments}
\begin{quote}

\sphinxAtStartPar
\sphinxcode{\sphinxupquote{cb}}
\begin{quote}

\sphinxAtStartPar
Callback Class object.
\end{quote}

\sphinxAtStartPar
\sphinxcode{\sphinxupquote{cbctx}}
\begin{quote}

\sphinxAtStartPar
Callback context. Please refer to {\hyperref[\detokenize{constant:chapconst-cbc}]{\sphinxcrossref{\DUrole{std,std-ref}{Callback context}}}} .
\end{quote}
\end{quote}

\sphinxAtStartPar
\sphinxstylestrong{Example}
\end{quote}

\begin{sphinxVerbatim}[commandchars=\\\{\}]
\PYG{n}{cb} \PYG{o}{=} \PYG{n}{CoptCallback}\PYG{p}{(}\PYG{p}{)}
\PYG{n}{model}\PYG{o}{.}\PYG{n}{setCallback}\PYG{p}{(}\PYG{n}{cb}\PYG{p}{,} \PYG{n}{COPT}\PYG{o}{.}\PYG{n}{CBCONTEXT\PYGZus{}MIPSOL}\PYG{p}{)}
\end{sphinxVerbatim}

\subsection{Var Class}
\label{\detokenize{pyapiref:var-class}}\label{\detokenize{pyapiref:chappyapi-var}}
\sphinxAtStartPar
For easy access to information of variables, Var object provides methods such as \sphinxcode{\sphinxupquote{Var.LB}}.
The full list of information can be found in the {\hyperref[\detokenize{pyapiref:chappyapi-const-info}]{\sphinxcrossref{\DUrole{std,std-ref}{Information}}}} section.
For convenience, information can be accessed by names in original case or lowercase.

\sphinxAtStartPar
In addition, you can also access the value of the variable through \sphinxcode{\sphinxupquote{Var.x}}, the variable type through \sphinxcode{\sphinxupquote{Var.vtype}}, the name of the variable through \sphinxcode{\sphinxupquote{Var.name}},
the Reduced cost value of the variable in LP through \sphinxcode{\sphinxupquote{Var.rc}}, the basis status through \sphinxcode{\sphinxupquote{Var.basis}},
and the index of the variable in the coefficient matrix through \sphinxcode{\sphinxupquote{Var.index}} the Reduced cost value of the variable in LP through \sphinxcode{\sphinxupquote{Var.rc}}, the basis status through \sphinxcode{\sphinxupquote{Var.basis}}, and the index of the variable in the coefficient matrix through \sphinxcode{\sphinxupquote{Var.index}}.

\sphinxAtStartPar
For the model\sphinxhyphen{}related information of the variables, as well as the variable type and name, the user can set the corresponding information value in the form of \sphinxcode{\sphinxupquote{"Var.LB = 0.0"}}.

\sphinxAtStartPar
Var object contains related operations of COPT variables and provides the following methods:

\subsubsection{Var.getType()}
\label{\detokenize{pyapiref:var-gettype}}\begin{quote}

\sphinxAtStartPar
\sphinxstylestrong{Synopsis}
\begin{quote}

\sphinxAtStartPar
\sphinxcode{\sphinxupquote{getType()}}
\end{quote}

\sphinxAtStartPar
\sphinxstylestrong{Description}
\begin{quote}

\sphinxAtStartPar
Retrieve the type of variable.
\end{quote}

\sphinxAtStartPar
\sphinxstylestrong{Example}
\end{quote}

\begin{sphinxVerbatim}[commandchars=\\\{\}]
\PYG{c+c1}{\PYGZsh{} Retrieve the type of variable v}
\PYG{n}{vtype} \PYG{o}{=} \PYG{n}{v}\PYG{o}{.}\PYG{n}{getType}\PYG{p}{(}\PYG{p}{)}
\end{sphinxVerbatim}

\subsubsection{Var.getName()}
\label{\detokenize{pyapiref:var-getname}}\begin{quote}

\sphinxAtStartPar
\sphinxstylestrong{Synopsis}
\begin{quote}

\sphinxAtStartPar
\sphinxcode{\sphinxupquote{getName()}}
\end{quote}

\sphinxAtStartPar
\sphinxstylestrong{Description}
\begin{quote}

\sphinxAtStartPar
Retrieve the name of variable.
\end{quote}

\sphinxAtStartPar
\sphinxstylestrong{Example}
\end{quote}

\begin{sphinxVerbatim}[commandchars=\\\{\}]
\PYG{c+c1}{\PYGZsh{} Retrieve the name of variable v}
\PYG{n}{varname} \PYG{o}{=} \PYG{n}{v}\PYG{o}{.}\PYG{n}{getName}\PYG{p}{(}\PYG{p}{)}
\end{sphinxVerbatim}

\subsubsection{Var.getBasis()}
\label{\detokenize{pyapiref:var-getbasis}}\begin{quote}

\sphinxAtStartPar
\sphinxstylestrong{Synopsis}
\begin{quote}

\sphinxAtStartPar
\sphinxcode{\sphinxupquote{getBasis()}}
\end{quote}

\sphinxAtStartPar
\sphinxstylestrong{Description}
\begin{quote}

\sphinxAtStartPar
Retrieve the basis status of variable.
\end{quote}

\sphinxAtStartPar
\sphinxstylestrong{Example}
\end{quote}

\begin{sphinxVerbatim}[commandchars=\\\{\}]
\PYG{c+c1}{\PYGZsh{} Retrieve the basis status of variable v}
\PYG{n}{varbasis} \PYG{o}{=} \PYG{n}{v}\PYG{o}{.}\PYG{n}{getBasis}\PYG{p}{(}\PYG{p}{)}
\end{sphinxVerbatim}

\subsubsection{Var.getLowerIIS()}
\label{\detokenize{pyapiref:var-getloweriis}}\begin{quote}

\sphinxAtStartPar
\sphinxstylestrong{Synopsis}
\begin{quote}

\sphinxAtStartPar
\sphinxcode{\sphinxupquote{getLowerIIS()}}
\end{quote}

\sphinxAtStartPar
\sphinxstylestrong{Description}
\begin{quote}

\sphinxAtStartPar
Retrieve the IIS status of lower bound of variable.
\end{quote}

\sphinxAtStartPar
\sphinxstylestrong{Example}
\end{quote}

\begin{sphinxVerbatim}[commandchars=\\\{\}]
\PYG{c+c1}{\PYGZsh{} Retrieve the IIS status of lower bound of variable v}
\PYG{n}{lowerIIS} \PYG{o}{=} \PYG{n}{v}\PYG{o}{.}\PYG{n}{getLowerIIS}\PYG{p}{(}\PYG{p}{)}
\end{sphinxVerbatim}

\subsubsection{Var.getUpperIIS()}
\label{\detokenize{pyapiref:var-getupperiis}}\begin{quote}

\sphinxAtStartPar
\sphinxstylestrong{Synopsis}
\begin{quote}

\sphinxAtStartPar
\sphinxcode{\sphinxupquote{getUpperIIS()}}
\end{quote}

\sphinxAtStartPar
\sphinxstylestrong{Description}
\begin{quote}

\sphinxAtStartPar
Retrieve the IIS status of upper bound of variable.
\end{quote}

\sphinxAtStartPar
\sphinxstylestrong{Example}
\end{quote}

\begin{sphinxVerbatim}[commandchars=\\\{\}]
\PYG{c+c1}{\PYGZsh{} Retrieve the IIS status of upper bound of variable v}
\PYG{n}{upperIIS} \PYG{o}{=} \PYG{n}{v}\PYG{o}{.}\PYG{n}{getUpperIIS}\PYG{p}{(}\PYG{p}{)}
\end{sphinxVerbatim}

\subsubsection{Var.getIdx()}
\label{\detokenize{pyapiref:var-getidx}}\begin{quote}

\sphinxAtStartPar
\sphinxstylestrong{Synopsis}
\begin{quote}

\sphinxAtStartPar
\sphinxcode{\sphinxupquote{getIdx()}}
\end{quote}

\sphinxAtStartPar
\sphinxstylestrong{Description}
\begin{quote}

\sphinxAtStartPar
Retrieve the subscript of the variable in the coefficient matrix.
\end{quote}

\sphinxAtStartPar
\sphinxstylestrong{Example}
\end{quote}

\begin{sphinxVerbatim}[commandchars=\\\{\}]
\PYG{c+c1}{\PYGZsh{}  Retrieve the subscript of variable v}
\PYG{n}{vindex} \PYG{o}{=} \PYG{n}{v}\PYG{o}{.}\PYG{n}{getIdx}\PYG{p}{(}\PYG{p}{)}
\end{sphinxVerbatim}

\subsubsection{Var.setType()}
\label{\detokenize{pyapiref:var-settype}}\begin{quote}

\sphinxAtStartPar
\sphinxstylestrong{Synopsis}
\begin{quote}

\sphinxAtStartPar
\sphinxcode{\sphinxupquote{setType(newtype)}}
\end{quote}

\sphinxAtStartPar
\sphinxstylestrong{Description}
\begin{quote}

\sphinxAtStartPar
Set the type of variable.
\end{quote}

\sphinxAtStartPar
\sphinxstylestrong{Arguments}
\begin{quote}

\sphinxAtStartPar
\sphinxcode{\sphinxupquote{newtype}}
\begin{quote}

\sphinxAtStartPar
The type of variable to be set. Please refer to {\hyperref[\detokenize{constant:chapconst-vartype}]{\sphinxcrossref{\DUrole{std,std-ref}{Variable types}}}} section for possible values.
\end{quote}
\end{quote}

\sphinxAtStartPar
\sphinxstylestrong{Example}
\end{quote}

\begin{sphinxVerbatim}[commandchars=\\\{\}]
\PYG{c+c1}{\PYGZsh{} Set the type of variable v}
\PYG{n}{v}\PYG{o}{.}\PYG{n}{setType}\PYG{p}{(}\PYG{n}{COPT}\PYG{o}{.}\PYG{n}{BINARY}\PYG{p}{)}
\end{sphinxVerbatim}

\subsubsection{Var.setName()}
\label{\detokenize{pyapiref:var-setname}}\begin{quote}

\sphinxAtStartPar
\sphinxstylestrong{Synopsis}
\begin{quote}

\sphinxAtStartPar
\sphinxcode{\sphinxupquote{setName(newname)}}
\end{quote}

\sphinxAtStartPar
\sphinxstylestrong{Description}
\begin{quote}

\sphinxAtStartPar
Set the name of variable.
\end{quote}

\sphinxAtStartPar
\sphinxstylestrong{Arguments}
\begin{quote}

\sphinxAtStartPar
\sphinxcode{\sphinxupquote{newname}}
\begin{quote}

\sphinxAtStartPar
The name of variable to be set.
\end{quote}
\end{quote}

\sphinxAtStartPar
\sphinxstylestrong{Example}
\end{quote}

\begin{sphinxVerbatim}[commandchars=\\\{\}]
\PYG{c+c1}{\PYGZsh{} Set the name of variable v}
\PYG{n}{v}\PYG{o}{.}\PYG{n}{setName}\PYG{p}{(}\PYG{n}{COPT}\PYG{o}{.}\PYG{n}{BINARY}\PYG{p}{)}
\end{sphinxVerbatim}

\subsubsection{Var.getInfo()}
\label{\detokenize{pyapiref:var-getinfo}}\begin{quote}

\sphinxAtStartPar
\sphinxstylestrong{Synopsis}
\begin{quote}

\sphinxAtStartPar
\sphinxcode{\sphinxupquote{getInfo(infoname)}}
\end{quote}

\sphinxAtStartPar
\sphinxstylestrong{Description}
\begin{quote}

\sphinxAtStartPar
Retrieve specified information. Return a constant.
\end{quote}

\sphinxAtStartPar
\sphinxstylestrong{Arguments}
\begin{quote}

\sphinxAtStartPar
\sphinxcode{\sphinxupquote{infoname}}
\begin{quote}

\sphinxAtStartPar
The name of the information. Please refer to {\hyperref[\detokenize{pyapiref:chappyapi-const-info}]{\sphinxcrossref{\DUrole{std,std-ref}{Information}}}} for possible values.
\end{quote}
\end{quote}

\sphinxAtStartPar
\sphinxstylestrong{Example}
\end{quote}

\begin{sphinxVerbatim}[commandchars=\\\{\}]
\PYG{c+c1}{\PYGZsh{} Get lowerbound of variable x}
\PYG{n}{x}\PYG{o}{.}\PYG{n}{getInfo}\PYG{p}{(}\PYG{n}{COPT}\PYG{o}{.}\PYG{n}{Info}\PYG{o}{.}\PYG{n}{LB}\PYG{p}{)}
\end{sphinxVerbatim}

\subsubsection{Var.setInfo()}
\label{\detokenize{pyapiref:var-setinfo}}\begin{quote}

\sphinxAtStartPar
\sphinxstylestrong{Synopsis}
\begin{quote}

\sphinxAtStartPar
\sphinxcode{\sphinxupquote{setInfo(infoname, newval)}}
\end{quote}

\sphinxAtStartPar
\sphinxstylestrong{Description}
\begin{quote}

\sphinxAtStartPar
Set new information value for a variable.
\end{quote}

\sphinxAtStartPar
\sphinxstylestrong{Arguments}
\begin{quote}

\sphinxAtStartPar
\sphinxcode{\sphinxupquote{infoname}}
\begin{quote}

\sphinxAtStartPar
The name of the information to be set. Please refer to {\hyperref[\detokenize{pyapiref:chappyapi-const-info}]{\sphinxcrossref{\DUrole{std,std-ref}{Information}}}} section for possible values.
\end{quote}

\sphinxAtStartPar
\sphinxcode{\sphinxupquote{newval}}
\begin{quote}

\sphinxAtStartPar
New information value to set.
\end{quote}
\end{quote}

\sphinxAtStartPar
\sphinxstylestrong{Example}
\end{quote}

\begin{sphinxVerbatim}[commandchars=\\\{\}]
\PYG{c+c1}{\PYGZsh{} Set the lower bound of variable x}
\PYG{n}{x}\PYG{o}{.}\PYG{n}{setInfo}\PYG{p}{(}\PYG{n}{COPT}\PYG{o}{.}\PYG{n}{Info}\PYG{o}{.}\PYG{n}{LB}\PYG{p}{,} \PYG{l+m+mf}{1.0}\PYG{p}{)}
\end{sphinxVerbatim}

\subsubsection{Var.remove()}
\label{\detokenize{pyapiref:var-remove}}\begin{quote}

\sphinxAtStartPar
\sphinxstylestrong{Synopsis}
\begin{quote}

\sphinxAtStartPar
\sphinxcode{\sphinxupquote{remove()}}
\end{quote}

\sphinxAtStartPar
\sphinxstylestrong{Description}
\begin{quote}

\sphinxAtStartPar
Delete the variable from model.
\end{quote}

\sphinxAtStartPar
\sphinxstylestrong{Example}
\end{quote}

\begin{sphinxVerbatim}[commandchars=\\\{\}]
\PYG{c+c1}{\PYGZsh{} Delete variable \PYGZsq{}x\PYGZsq{}}
\PYG{n}{x}\PYG{o}{.}\PYG{n}{remove}\PYG{p}{(}\PYG{p}{)}
\end{sphinxVerbatim}

\subsection{VarArray Class}
\label{\detokenize{pyapiref:vararray-class}}\label{\detokenize{pyapiref:chappyapi-vararray}}
\sphinxAtStartPar
To facilitate users to operate on multiple {\hyperref[\detokenize{pyapiref:chappyapi-var}]{\sphinxcrossref{\DUrole{std,std-ref}{Var Class}}}} objects,
the Python interface of COPT provides VarArray object with the following methods:

\subsubsection{VarArray()}
\label{\detokenize{pyapiref:vararray}}\begin{quote}

\sphinxAtStartPar
\sphinxstylestrong{Synopsis}
\begin{quote}

\sphinxAtStartPar
\sphinxcode{\sphinxupquote{VarArray(vars=None)}}
\end{quote}

\sphinxAtStartPar
\sphinxstylestrong{Description}
\begin{quote}

\sphinxAtStartPar
Create a {\hyperref[\detokenize{pyapiref:chappyapi-vararray}]{\sphinxcrossref{\DUrole{std,std-ref}{VarArray Class}}}} object.

\sphinxAtStartPar
If parameter \sphinxcode{\sphinxupquote{vars}} is \sphinxcode{\sphinxupquote{None}}, then create an empty {\hyperref[\detokenize{pyapiref:chappyapi-vararray}]{\sphinxcrossref{\DUrole{std,std-ref}{VarArray Class}}}} object,
otherwise initialize the new created {\hyperref[\detokenize{pyapiref:chappyapi-vararray}]{\sphinxcrossref{\DUrole{std,std-ref}{VarArray Class}}}} object based on \sphinxcode{\sphinxupquote{vars}}.
\end{quote}

\sphinxAtStartPar
\sphinxstylestrong{Arguments}
\begin{quote}

\sphinxAtStartPar
\sphinxcode{\sphinxupquote{vars}}
\begin{quote}

\sphinxAtStartPar
Variables to be added. Optional, \sphinxcode{\sphinxupquote{None}} by default.
\sphinxcode{\sphinxupquote{vars}} can be {\hyperref[\detokenize{pyapiref:chappyapi-var}]{\sphinxcrossref{\DUrole{std,std-ref}{Var Class}}}} object, {\hyperref[\detokenize{pyapiref:chappyapi-vararray}]{\sphinxcrossref{\DUrole{std,std-ref}{VarArray Class}}}} object, list, dictionary or
{\hyperref[\detokenize{pyapiref:chappyapi-util-tupledict}]{\sphinxcrossref{\DUrole{std,std-ref}{tupledict Class}}}} object.
\end{quote}
\end{quote}

\sphinxAtStartPar
\sphinxstylestrong{Example}
\end{quote}

\begin{sphinxVerbatim}[commandchars=\\\{\}]
\PYG{c+c1}{\PYGZsh{} Create an empty VarArray object}
\PYG{n}{vararr} \PYG{o}{=} \PYG{n}{VarArray}\PYG{p}{(}\PYG{p}{)}
\PYG{c+c1}{\PYGZsh{} Create an empty VarArray object and initialize variables x, y.}
\PYG{n}{vararr} \PYG{o}{=} \PYG{n}{VarArray}\PYG{p}{(}\PYG{p}{[}\PYG{n}{x}\PYG{p}{,} \PYG{n}{y}\PYG{p}{]}\PYG{p}{)}
\end{sphinxVerbatim}

\subsubsection{VarArray.pushBack()}
\label{\detokenize{pyapiref:vararray-pushback}}\begin{quote}

\sphinxAtStartPar
\sphinxstylestrong{Synopsis}
\begin{quote}

\sphinxAtStartPar
\sphinxcode{\sphinxupquote{pushBack(vars)}}
\end{quote}

\sphinxAtStartPar
\sphinxstylestrong{Description}
\begin{quote}

\sphinxAtStartPar
Add single or multiple {\hyperref[\detokenize{pyapiref:chappyapi-var}]{\sphinxcrossref{\DUrole{std,std-ref}{Var Class}}}} objects.
\end{quote}

\sphinxAtStartPar
\sphinxstylestrong{Arguments}
\begin{quote}

\sphinxAtStartPar
\sphinxcode{\sphinxupquote{vars}}
\begin{quote}

\sphinxAtStartPar
Variables to be applied.
\sphinxcode{\sphinxupquote{vars}} can be {\hyperref[\detokenize{pyapiref:chappyapi-var}]{\sphinxcrossref{\DUrole{std,std-ref}{Var Class}}}} object, {\hyperref[\detokenize{pyapiref:chappyapi-vararray}]{\sphinxcrossref{\DUrole{std,std-ref}{VarArray Class}}}} object,
list, dictionary or {\hyperref[\detokenize{pyapiref:chappyapi-util-tupledict}]{\sphinxcrossref{\DUrole{std,std-ref}{tupledict Class}}}} object.
\end{quote}
\end{quote}

\sphinxAtStartPar
\sphinxstylestrong{Example}
\end{quote}

\begin{sphinxVerbatim}[commandchars=\\\{\}]
\PYG{c+c1}{\PYGZsh{} Add variable x to vararr}
\PYG{n}{vararr}\PYG{o}{.}\PYG{n}{pushBack}\PYG{p}{(}\PYG{n}{x}\PYG{p}{)}
\PYG{c+c1}{\PYGZsh{} Add variables x and y to vararr}
\PYG{n}{vararr}\PYG{o}{.}\PYG{n}{pushBack}\PYG{p}{(}\PYG{p}{[}\PYG{n}{x}\PYG{p}{,} \PYG{n}{y}\PYG{p}{]}\PYG{p}{)}
\end{sphinxVerbatim}

\subsubsection{VarArray.getVar()}
\label{\detokenize{pyapiref:vararray-getvar}}\begin{quote}

\sphinxAtStartPar
\sphinxstylestrong{Synopsis}
\begin{quote}

\sphinxAtStartPar
\sphinxcode{\sphinxupquote{getVar(idx)}}
\end{quote}

\sphinxAtStartPar
\sphinxstylestrong{Description}
\begin{quote}

\sphinxAtStartPar
Retrieve a variable from an index in a {\hyperref[\detokenize{pyapiref:chappyapi-vararray}]{\sphinxcrossref{\DUrole{std,std-ref}{VarArray Class}}}} object. Return a {\hyperref[\detokenize{pyapiref:chappyapi-var}]{\sphinxcrossref{\DUrole{std,std-ref}{Var Class}}}} object.
\end{quote}

\sphinxAtStartPar
\sphinxstylestrong{Arguments}
\begin{quote}

\sphinxAtStartPar
\sphinxcode{\sphinxupquote{idx}}
\begin{quote}

\sphinxAtStartPar
Subscript of the specified variable in {\hyperref[\detokenize{pyapiref:chappyapi-vararray}]{\sphinxcrossref{\DUrole{std,std-ref}{VarArray Class}}}} object, starting with 0.
\end{quote}
\end{quote}

\sphinxAtStartPar
\sphinxstylestrong{Example}
\end{quote}

\begin{sphinxVerbatim}[commandchars=\\\{\}]
\PYG{c+c1}{\PYGZsh{} Get the variable with subscript of 1 in vararr}
 \PYG{n}{vararr}\PYG{o}{.}\PYG{n}{getVar}\PYG{p}{(}\PYG{l+m+mi}{1}\PYG{p}{)}
\end{sphinxVerbatim}

\subsubsection{VarArray.getAll()}
\label{\detokenize{pyapiref:vararray-getall}}\begin{quote}

\sphinxAtStartPar
\sphinxstylestrong{Synopsis}
\begin{quote}

\sphinxAtStartPar
\sphinxcode{\sphinxupquote{getAll()}}
\end{quote}

\sphinxAtStartPar
\sphinxstylestrong{Description}
\begin{quote}

\sphinxAtStartPar
Retrieve all variables in {\hyperref[\detokenize{pyapiref:chappyapi-vararray}]{\sphinxcrossref{\DUrole{std,std-ref}{VarArray Class}}}} object. Returns a list object.
\end{quote}

\sphinxAtStartPar
\sphinxstylestrong{Example}
\end{quote}

\begin{sphinxVerbatim}[commandchars=\\\{\}]
\PYG{c+c1}{\PYGZsh{} Get all variables in \PYGZsq{}vararr\PYGZsq{}}
\PYG{n}{varall} \PYG{o}{=} \PYG{n}{vararr}\PYG{o}{.}\PYG{n}{getAll}\PYG{p}{(}\PYG{p}{)}
\end{sphinxVerbatim}

\subsubsection{VarArray.getSize()}
\label{\detokenize{pyapiref:vararray-getsize}}\begin{quote}

\sphinxAtStartPar
\sphinxstylestrong{Synopsis}
\begin{quote}

\sphinxAtStartPar
\sphinxcode{\sphinxupquote{getSize()}}
\end{quote}

\sphinxAtStartPar
\sphinxstylestrong{Description}
\begin{quote}

\sphinxAtStartPar
Retrieve the number of variables in {\hyperref[\detokenize{pyapiref:chappyapi-vararray}]{\sphinxcrossref{\DUrole{std,std-ref}{VarArray Class}}}} object.
\end{quote}

\sphinxAtStartPar
\sphinxstylestrong{Example}
\end{quote}

\begin{sphinxVerbatim}[commandchars=\\\{\}]
\PYG{c+c1}{\PYGZsh{} Retrieve the number of variables in vararr.}
\PYG{n}{arrsize} \PYG{o}{=} \PYG{n}{vararr}\PYG{o}{.}\PYG{n}{getSize}\PYG{p}{(}\PYG{p}{)}
\end{sphinxVerbatim}

\subsection{PsdVar Class}
\label{\detokenize{pyapiref:psdvar-class}}\label{\detokenize{pyapiref:chappyapi-psdvar}}
\sphinxAtStartPar
PsdVar object contains related operations of COPT positive semi\sphinxhyphen{}definite variables and
provides the following methods:

\subsubsection{PsdVar.getName()}
\label{\detokenize{pyapiref:psdvar-getname}}\begin{quote}

\sphinxAtStartPar
\sphinxstylestrong{Synopsis}
\begin{quote}

\sphinxAtStartPar
\sphinxcode{\sphinxupquote{getName()}}
\end{quote}

\sphinxAtStartPar
\sphinxstylestrong{Description}
\begin{quote}

\sphinxAtStartPar
Retrieve the name of positive semi\sphinxhyphen{}definite variable.
\end{quote}

\sphinxAtStartPar
\sphinxstylestrong{Example}
\end{quote}

\begin{sphinxVerbatim}[commandchars=\\\{\}]
\PYG{c+c1}{\PYGZsh{} Retrieve the name of variable v}
\PYG{n}{varname} \PYG{o}{=} \PYG{n}{v}\PYG{o}{.}\PYG{n}{getName}\PYG{p}{(}\PYG{p}{)}
\end{sphinxVerbatim}

\subsubsection{PsdVar.getIdx()}
\label{\detokenize{pyapiref:psdvar-getidx}}\begin{quote}

\sphinxAtStartPar
\sphinxstylestrong{Synopsis}
\begin{quote}

\sphinxAtStartPar
\sphinxcode{\sphinxupquote{getIdx()}}
\end{quote}

\sphinxAtStartPar
\sphinxstylestrong{Description}
\begin{quote}

\sphinxAtStartPar
Retrieve the subscript of the variable in the model.
\end{quote}

\sphinxAtStartPar
\sphinxstylestrong{Example}
\end{quote}

\begin{sphinxVerbatim}[commandchars=\\\{\}]
\PYG{c+c1}{\PYGZsh{} Retrieve the subscript of variable v}
\PYG{n}{vindex} \PYG{o}{=} \PYG{n}{v}\PYG{o}{.}\PYG{n}{getIdx}\PYG{p}{(}\PYG{p}{)}
\end{sphinxVerbatim}

\subsubsection{PsdVar.getDim()}
\label{\detokenize{pyapiref:psdvar-getdim}}\begin{quote}

\sphinxAtStartPar
\sphinxstylestrong{Synopsis}
\begin{quote}

\sphinxAtStartPar
\sphinxcode{\sphinxupquote{getDim()}}
\end{quote}

\sphinxAtStartPar
\sphinxstylestrong{Description}
\begin{quote}

\sphinxAtStartPar
Retrieve the dimension of positive semi\sphinxhyphen{}definite variable.
\end{quote}

\sphinxAtStartPar
\sphinxstylestrong{Example}
\end{quote}

\begin{sphinxVerbatim}[commandchars=\\\{\}]
\PYG{c+c1}{\PYGZsh{} Retrieve the dimension of variable \PYGZdq{}v\PYGZdq{}}
\PYG{n}{vdim} \PYG{o}{=} \PYG{n}{v}\PYG{o}{.}\PYG{n}{getDim}\PYG{p}{(}\PYG{p}{)}
\end{sphinxVerbatim}

\subsubsection{PsdVar.getLen()}
\label{\detokenize{pyapiref:psdvar-getlen}}\begin{quote}

\sphinxAtStartPar
\sphinxstylestrong{Synopsis}
\begin{quote}

\sphinxAtStartPar
\sphinxcode{\sphinxupquote{getLen()}}
\end{quote}

\sphinxAtStartPar
\sphinxstylestrong{Description}
\begin{quote}

\sphinxAtStartPar
Retrieve the length of the expanded positive semi\sphinxhyphen{}definite variable.
\end{quote}

\sphinxAtStartPar
\sphinxstylestrong{Example}
\end{quote}

\begin{sphinxVerbatim}[commandchars=\\\{\}]
\PYG{c+c1}{\PYGZsh{} Retrieve the length of the expanded positive semi\PYGZhy{}definite variable \PYGZdq{}v\PYGZdq{}}
\PYG{n}{vlen} \PYG{o}{=} \PYG{n}{v}\PYG{o}{.}\PYG{n}{getLen}\PYG{p}{(}\PYG{p}{)}
\end{sphinxVerbatim}

\subsubsection{PsdVar.setName()}
\label{\detokenize{pyapiref:psdvar-setname}}\begin{quote}

\sphinxAtStartPar
\sphinxstylestrong{Synopsis}
\begin{quote}

\sphinxAtStartPar
\sphinxcode{\sphinxupquote{setName(newname)}}
\end{quote}

\sphinxAtStartPar
\sphinxstylestrong{Description}
\begin{quote}

\sphinxAtStartPar
Set the name of positive semi\sphinxhyphen{}definite variable.
\end{quote}

\sphinxAtStartPar
\sphinxstylestrong{Arguments}
\begin{quote}

\sphinxAtStartPar
\sphinxcode{\sphinxupquote{newname}}
\begin{quote}

\sphinxAtStartPar
The name of positive semi\sphinxhyphen{}definite variable to be set.
\end{quote}
\end{quote}

\sphinxAtStartPar
\sphinxstylestrong{Example}
\end{quote}

\begin{sphinxVerbatim}[commandchars=\\\{\}]
\PYG{c+c1}{\PYGZsh{} Set the name of variable v}
\PYG{n}{v}\PYG{o}{.}\PYG{n}{setName}\PYG{p}{(}\PYG{l+s+s1}{\PYGZsq{}}\PYG{l+s+s1}{v}\PYG{l+s+s1}{\PYGZsq{}}\PYG{p}{)}
\end{sphinxVerbatim}

\subsubsection{PsdVar.getInfo()}
\label{\detokenize{pyapiref:psdvar-getinfo}}\begin{quote}

\sphinxAtStartPar
\sphinxstylestrong{Synopsis}
\begin{quote}

\sphinxAtStartPar
\sphinxcode{\sphinxupquote{getInfo(infoname)}}
\end{quote}

\sphinxAtStartPar
\sphinxstylestrong{Description}
\begin{quote}

\sphinxAtStartPar
Retrieve specified information of positive semi\sphinxhyphen{}definite variable. Return a list.
\end{quote}

\sphinxAtStartPar
\sphinxstylestrong{Arguments}
\begin{quote}

\sphinxAtStartPar
\sphinxcode{\sphinxupquote{infoname}}
\begin{quote}

\sphinxAtStartPar
The name of the information. Please refer to {\hyperref[\detokenize{pyapiref:chappyapi-const-info}]{\sphinxcrossref{\DUrole{std,std-ref}{Information}}}} for possible values.
\end{quote}
\end{quote}

\sphinxAtStartPar
\sphinxstylestrong{Example}
\end{quote}

\begin{sphinxVerbatim}[commandchars=\\\{\}]
\PYG{c+c1}{\PYGZsh{} Get solution values of positive semi\PYGZhy{}definite variable x}
\PYG{n}{sol} \PYG{o}{=} \PYG{n}{x}\PYG{o}{.}\PYG{n}{getInfo}\PYG{p}{(}\PYG{n}{COPT}\PYG{o}{.}\PYG{n}{Info}\PYG{o}{.}\PYG{n}{Value}\PYG{p}{)}
\end{sphinxVerbatim}

\subsubsection{PsdVar.remove()}
\label{\detokenize{pyapiref:psdvar-remove}}\begin{quote}

\sphinxAtStartPar
\sphinxstylestrong{Synopsis}
\begin{quote}

\sphinxAtStartPar
\sphinxcode{\sphinxupquote{remove()}}
\end{quote}

\sphinxAtStartPar
\sphinxstylestrong{Description}
\begin{quote}

\sphinxAtStartPar
Delete the positive semi\sphinxhyphen{}definite variable from model.
\end{quote}

\sphinxAtStartPar
\sphinxstylestrong{Example}
\end{quote}

\begin{sphinxVerbatim}[commandchars=\\\{\}]
\PYG{c+c1}{\PYGZsh{} Delete variable \PYGZsq{}x}
\PYG{n}{x}\PYG{o}{.}\PYG{n}{remove}\PYG{p}{(}\PYG{p}{)}
\end{sphinxVerbatim}

\subsubsection{PsdVar.diag()}
\label{\detokenize{pyapiref:psdvar-diag}}\begin{quote}

\sphinxAtStartPar
\sphinxstylestrong{Synopsis}
\begin{quote}

\sphinxAtStartPar
\sphinxcode{\sphinxupquote{diag(offset=0)}}
\end{quote}

\sphinxAtStartPar
\sphinxstylestrong{Description}
\begin{quote}

\sphinxAtStartPar
Retrieves the diagonal elements of the PSD variable.
\end{quote}

\sphinxAtStartPar
\sphinxstylestrong{Arguments}
\begin{quote}

\sphinxAtStartPar
\sphinxcode{\sphinxupquote{offset}}
\begin{quote}

\sphinxAtStartPar
The diagonal offset, default is 0.

\sphinxAtStartPar
If \sphinxcode{\sphinxupquote{offset}} \textgreater{} 0, it represents a downward diagonal offset.
If \sphinxcode{\sphinxupquote{offset}} \textless{} 0, it represents an upward diagonal offset.
\end{quote}
\end{quote}

\sphinxAtStartPar
\sphinxstylestrong{Return Value}
\begin{quote}

\sphinxAtStartPar
A \sphinxtitleref{PsdExpr} object.
\end{quote}

\sphinxAtStartPar
\sphinxstylestrong{Example}
\end{quote}

\begin{sphinxVerbatim}[commandchars=\\\{\}]
\PYG{c+c1}{\PYGZsh{} Retrieve the main diagonal elements of the PSD variable}
\PYG{n}{v}\PYG{o}{.}\PYG{n}{diag}\PYG{p}{(}\PYG{n}{offset}\PYG{o}{=}\PYG{l+m+mi}{0}\PYG{p}{)}
\end{sphinxVerbatim}

\subsubsection{PsdVar.pick()}
\label{\detokenize{pyapiref:psdvar-pick}}\begin{quote}

\sphinxAtStartPar
\sphinxstylestrong{Synopsis}
\begin{quote}

\sphinxAtStartPar
\sphinxcode{\sphinxupquote{pick(indexes)}}
\end{quote}

\sphinxAtStartPar
\sphinxstylestrong{Description}
\begin{quote}

\sphinxAtStartPar
Retrieves a PSD expression composed of the elements at the specified indices of the PSD variable.
\end{quote}

\sphinxAtStartPar
\sphinxstylestrong{Arguments}
\begin{quote}

\sphinxAtStartPar
\sphinxcode{\sphinxupquote{indexes}}
\begin{quote}

\sphinxAtStartPar
The array of specified indices.
\end{quote}
\end{quote}

\sphinxAtStartPar
\sphinxstylestrong{Return Value}
\begin{quote}

\sphinxAtStartPar
A \sphinxtitleref{PsdExpr} object.
\end{quote}

\sphinxAtStartPar
\sphinxstylestrong{Example}
\end{quote}

\begin{sphinxVerbatim}[commandchars=\\\{\}]
\PYG{c+c1}{\PYGZsh{} Retrieve a PSD expression composed of the element at index 0}
\PYG{n}{barX} \PYG{o}{=} \PYG{n}{model}\PYG{o}{.}\PYG{n}{addPsdVars}\PYG{p}{(}\PYG{l+m+mi}{3}\PYG{p}{,} \PYG{l+s+s2}{\PYGZdq{}}\PYG{l+s+s2}{BAR\PYGZus{}X}\PYG{l+s+s2}{\PYGZdq{}}\PYG{p}{)}
\PYG{n}{barX}\PYG{o}{.}\PYG{n}{pick}\PYG{p}{(}\PYG{p}{[}\PYG{l+m+mi}{0}\PYG{p}{]}\PYG{p}{)}
\end{sphinxVerbatim}

\subsubsection{PsdVar.sum()}
\label{\detokenize{pyapiref:psdvar-sum}}\begin{quote}

\sphinxAtStartPar
\sphinxstylestrong{Synopsis}
\begin{quote}

\sphinxAtStartPar
\sphinxcode{\sphinxupquote{sum()}}
\end{quote}

\sphinxAtStartPar
\sphinxstylestrong{Description}
\begin{quote}

\sphinxAtStartPar
Retrieves a PSD expression composed of the sum of all elements in the PSD variable.
\end{quote}

\sphinxAtStartPar
\sphinxstylestrong{Return Value}
\begin{quote}

\sphinxAtStartPar
A \sphinxtitleref{PsdExpr} object.
\end{quote}
\end{quote}

\subsubsection{PsdVar.toexpr()}
\label{\detokenize{pyapiref:psdvar-toexpr}}\begin{quote}

\sphinxAtStartPar
\sphinxstylestrong{Synopsis}
\begin{quote}

\sphinxAtStartPar
\sphinxcode{\sphinxupquote{toexpr()}}
\end{quote}

\sphinxAtStartPar
\sphinxstylestrong{Description}
\begin{quote}

\sphinxAtStartPar
Retrieves the PSD expression equivalent to the PSD variable.
\end{quote}

\sphinxAtStartPar
\sphinxstylestrong{Return Value}
\begin{quote}

\sphinxAtStartPar
A \sphinxtitleref{PsdExpr} object.
\end{quote}
\end{quote}

\subsubsection{PsdVar.shape}
\label{\detokenize{pyapiref:psdvar-shape}}\begin{quote}

\sphinxAtStartPar
\sphinxstylestrong{Synopsis}
\begin{quote}

\sphinxAtStartPar
\sphinxcode{\sphinxupquote{shape}}
\end{quote}

\sphinxAtStartPar
\sphinxstylestrong{Description}
\begin{quote}

\sphinxAtStartPar
Shape of the \sphinxcode{\sphinxupquote{PsdVar}} object.
\end{quote}

\sphinxAtStartPar
\sphinxstylestrong{Return value}
\begin{quote}

\sphinxAtStartPar
Integer tuple.
\end{quote}
\end{quote}

\subsubsection{PsdVar.size}
\label{\detokenize{pyapiref:psdvar-size}}\begin{quote}

\sphinxAtStartPar
\sphinxstylestrong{Synopsis}
\begin{quote}

\sphinxAtStartPar
\sphinxcode{\sphinxupquote{size}}
\end{quote}

\sphinxAtStartPar
\sphinxstylestrong{Description}
\begin{quote}

\sphinxAtStartPar
Size of the \sphinxcode{\sphinxupquote{PsdVar}} object.
\end{quote}

\sphinxAtStartPar
\sphinxstylestrong{Return value}
\begin{quote}

\sphinxAtStartPar
Integer tuple.
\end{quote}
\end{quote}

\subsubsection{PsdVar.dim}
\label{\detokenize{pyapiref:psdvar-dim}}\begin{quote}

\sphinxAtStartPar
\sphinxstylestrong{Synopsis}
\begin{quote}

\sphinxAtStartPar
\sphinxcode{\sphinxupquote{dim}}
\end{quote}

\sphinxAtStartPar
\sphinxstylestrong{Description}
\begin{quote}

\sphinxAtStartPar
Dimension of the \sphinxcode{\sphinxupquote{PsdVar}} object.
\end{quote}

\sphinxAtStartPar
\sphinxstylestrong{Return value}
\begin{quote}

\sphinxAtStartPar
Integer.
\end{quote}
\end{quote}

\subsubsection{PsdVar.len}
\label{\detokenize{pyapiref:psdvar-len}}\begin{quote}

\sphinxAtStartPar
\sphinxstylestrong{Synopsis}
\begin{quote}

\sphinxAtStartPar
\sphinxcode{\sphinxupquote{len}}
\end{quote}

\sphinxAtStartPar
\sphinxstylestrong{Description}
\begin{quote}

\sphinxAtStartPar
Flattened length of the \sphinxcode{\sphinxupquote{PsdVar}} object.
\end{quote}

\sphinxAtStartPar
\sphinxstylestrong{Return value}
\begin{quote}

\sphinxAtStartPar
Integer.
\end{quote}
\end{quote}

\subsection{PsdVarArray Class}
\label{\detokenize{pyapiref:psdvararray-class}}\label{\detokenize{pyapiref:chappyapi-psdvararray}}
\sphinxAtStartPar
To facilitate users to operate on multiple {\hyperref[\detokenize{pyapiref:chappyapi-psdvar}]{\sphinxcrossref{\DUrole{std,std-ref}{PsdVar Class}}}} objects,
the Python interface of COPT provides PsdVarArray object with the following methods:

\subsubsection{PsdVarArray()}
\label{\detokenize{pyapiref:psdvararray}}\begin{quote}

\sphinxAtStartPar
\sphinxstylestrong{Synopsis}
\begin{quote}

\sphinxAtStartPar
\sphinxcode{\sphinxupquote{PsdVarArray(vars=None)}}
\end{quote}

\sphinxAtStartPar
\sphinxstylestrong{Description}
\begin{quote}

\sphinxAtStartPar
Create a {\hyperref[\detokenize{pyapiref:chappyapi-psdvararray}]{\sphinxcrossref{\DUrole{std,std-ref}{PsdVarArray Class}}}} object.

\sphinxAtStartPar
If parameter \sphinxcode{\sphinxupquote{vars}} is \sphinxcode{\sphinxupquote{None}}, then create an empty {\hyperref[\detokenize{pyapiref:chappyapi-psdvararray}]{\sphinxcrossref{\DUrole{std,std-ref}{PsdVarArray Class}}}} object,
otherwise initialize the new created {\hyperref[\detokenize{pyapiref:chappyapi-psdvararray}]{\sphinxcrossref{\DUrole{std,std-ref}{PsdVarArray Class}}}} object based on \sphinxcode{\sphinxupquote{vars}}.
\end{quote}

\sphinxAtStartPar
\sphinxstylestrong{Arguments}
\begin{quote}

\sphinxAtStartPar
\sphinxcode{\sphinxupquote{vars}}
\begin{quote}

\sphinxAtStartPar
Positive semi\sphinxhyphen{}definite variables to be added. Optional, \sphinxcode{\sphinxupquote{None}} by default.
\sphinxcode{\sphinxupquote{vars}} can be {\hyperref[\detokenize{pyapiref:chappyapi-psdvar}]{\sphinxcrossref{\DUrole{std,std-ref}{PsdVar Class}}}} object, {\hyperref[\detokenize{pyapiref:chappyapi-psdvararray}]{\sphinxcrossref{\DUrole{std,std-ref}{PsdVarArray Class}}}} object, list, dictionary or
{\hyperref[\detokenize{pyapiref:chappyapi-util-tupledict}]{\sphinxcrossref{\DUrole{std,std-ref}{tupledict Class}}}} object.
\end{quote}
\end{quote}

\sphinxAtStartPar
\sphinxstylestrong{Example}
\end{quote}

\begin{sphinxVerbatim}[commandchars=\\\{\}]
\PYG{c+c1}{\PYGZsh{} Create an empty PsdVarArray object}
\PYG{n}{vararr} \PYG{o}{=} \PYG{n}{PsdVarArray}\PYG{p}{(}\PYG{p}{)}
\PYG{c+c1}{\PYGZsh{} Create a PsdVarArray object containing positive semi\PYGZhy{}definite variables x, y.}
\PYG{n}{vararr} \PYG{o}{=} \PYG{n}{PsdVarArray}\PYG{p}{(}\PYG{p}{[}\PYG{n}{x}\PYG{p}{,} \PYG{n}{y}\PYG{p}{]}\PYG{p}{)}
\end{sphinxVerbatim}

\subsubsection{PsdVarArray.pushBack()}
\label{\detokenize{pyapiref:psdvararray-pushback}}\begin{quote}

\sphinxAtStartPar
\sphinxstylestrong{Synopsis}
\begin{quote}

\sphinxAtStartPar
\sphinxcode{\sphinxupquote{pushBack(var)}}
\end{quote}

\sphinxAtStartPar
\sphinxstylestrong{Description}
\begin{quote}

\sphinxAtStartPar
Add single or multiple {\hyperref[\detokenize{pyapiref:chappyapi-psdvar}]{\sphinxcrossref{\DUrole{std,std-ref}{PsdVar Class}}}} objects.
\end{quote}

\sphinxAtStartPar
\sphinxstylestrong{Arguments}
\begin{quote}

\sphinxAtStartPar
\sphinxcode{\sphinxupquote{var}}
\begin{quote}

\sphinxAtStartPar
Postive semi\sphinxhyphen{}definite variables to be applied.
\sphinxcode{\sphinxupquote{vars}} can be {\hyperref[\detokenize{pyapiref:chappyapi-psdvar}]{\sphinxcrossref{\DUrole{std,std-ref}{PsdVar Class}}}} object, {\hyperref[\detokenize{pyapiref:chappyapi-psdvararray}]{\sphinxcrossref{\DUrole{std,std-ref}{PsdVarArray Class}}}} object,
list, dictionary or {\hyperref[\detokenize{pyapiref:chappyapi-util-tupledict}]{\sphinxcrossref{\DUrole{std,std-ref}{tupledict Class}}}} object.
\end{quote}
\end{quote}

\sphinxAtStartPar
\sphinxstylestrong{Example}
\end{quote}

\begin{sphinxVerbatim}[commandchars=\\\{\}]
\PYG{c+c1}{\PYGZsh{} Add variable x to vararr}
\PYG{n}{vararr}\PYG{o}{.}\PYG{n}{pushBack}\PYG{p}{(}\PYG{n}{x}\PYG{p}{)}
\PYG{c+c1}{\PYGZsh{} Add variables x and y to vararr}
\PYG{n}{vararr}\PYG{o}{.}\PYG{n}{pushBack}\PYG{p}{(}\PYG{p}{[}\PYG{n}{x}\PYG{p}{,} \PYG{n}{y}\PYG{p}{]}\PYG{p}{)}
\end{sphinxVerbatim}

\subsubsection{PsdVarArray.getPsdVar()}
\label{\detokenize{pyapiref:psdvararray-getpsdvar}}\begin{quote}

\sphinxAtStartPar
\sphinxstylestrong{Synopsis}
\begin{quote}

\sphinxAtStartPar
\sphinxcode{\sphinxupquote{getPsdVar(idx)}}
\end{quote}

\sphinxAtStartPar
\sphinxstylestrong{Description}
\begin{quote}

\sphinxAtStartPar
Retrieve a positive semi\sphinxhyphen{}definite variable from an index in a {\hyperref[\detokenize{pyapiref:chappyapi-psdvararray}]{\sphinxcrossref{\DUrole{std,std-ref}{PsdVarArray Class}}}} object.
Return a {\hyperref[\detokenize{pyapiref:chappyapi-psdvar}]{\sphinxcrossref{\DUrole{std,std-ref}{PsdVar Class}}}} object.
\end{quote}

\sphinxAtStartPar
\sphinxstylestrong{Arguments}
\begin{quote}

\sphinxAtStartPar
\sphinxcode{\sphinxupquote{idx}}
\begin{quote}

\sphinxAtStartPar
Subscript of the specified positive semi\sphinxhyphen{}definite variable in {\hyperref[\detokenize{pyapiref:chappyapi-psdvararray}]{\sphinxcrossref{\DUrole{std,std-ref}{PsdVarArray Class}}}} object, starting with 0.
\end{quote}
\end{quote}

\sphinxAtStartPar
\sphinxstylestrong{Example}
\end{quote}

\begin{sphinxVerbatim}[commandchars=\\\{\}]
\PYG{c+c1}{\PYGZsh{} Get the positive semi\PYGZhy{}definite variable with subscript of 1 in vararr}
\PYG{n}{var} \PYG{o}{=} \PYG{n}{vararr}\PYG{o}{.}\PYG{n}{getPsdVar}\PYG{p}{(}\PYG{l+m+mi}{1}\PYG{p}{)}
\end{sphinxVerbatim}

\subsubsection{PsdVarArray.getSize()}
\label{\detokenize{pyapiref:psdvararray-getsize}}\begin{quote}

\sphinxAtStartPar
\sphinxstylestrong{Synopsis}
\begin{quote}

\sphinxAtStartPar
\sphinxcode{\sphinxupquote{getSize()}}
\end{quote}

\sphinxAtStartPar
\sphinxstylestrong{Description}
\begin{quote}

\sphinxAtStartPar
Retrieve the number of positive semi\sphinxhyphen{}definite variables in {\hyperref[\detokenize{pyapiref:chappyapi-psdvararray}]{\sphinxcrossref{\DUrole{std,std-ref}{PsdVarArray Class}}}} object.
\end{quote}

\sphinxAtStartPar
\sphinxstylestrong{Example}
\end{quote}

\begin{sphinxVerbatim}[commandchars=\\\{\}]
\PYG{c+c1}{\PYGZsh{} Retrieve the number of variables in vararr.}
\PYG{n}{arrsize} \PYG{o}{=} \PYG{n}{vararr}\PYG{o}{.}\PYG{n}{getSize}\PYG{p}{(}\PYG{p}{)}
\end{sphinxVerbatim}

\subsection{SymMatrix Class}
\label{\detokenize{pyapiref:symmatrix-class}}\label{\detokenize{pyapiref:chappyapi-symmatrix}}
\sphinxAtStartPar
SymMatrix object contains related operations of COPT symmetric matrices and
provides the following methods:

\subsubsection{SymMatrix.getIdx()}
\label{\detokenize{pyapiref:symmatrix-getidx}}\begin{quote}

\sphinxAtStartPar
\sphinxstylestrong{Synopsis}
\begin{quote}

\sphinxAtStartPar
\sphinxcode{\sphinxupquote{getIdx()}}
\end{quote}

\sphinxAtStartPar
\sphinxstylestrong{Description}
\begin{quote}

\sphinxAtStartPar
Retrieve the subscript of the symmetric matrix in the model.
\end{quote}

\sphinxAtStartPar
\sphinxstylestrong{Example}
\end{quote}

\begin{sphinxVerbatim}[commandchars=\\\{\}]
\PYG{c+c1}{\PYGZsh{} Retrieve the subscript of symmetric matrix mat}
\PYG{n}{matidx} \PYG{o}{=} \PYG{n}{mat}\PYG{o}{.}\PYG{n}{getIdx}\PYG{p}{(}\PYG{p}{)}
\end{sphinxVerbatim}

\subsubsection{SymMatrix.getDim()}
\label{\detokenize{pyapiref:symmatrix-getdim}}\begin{quote}

\sphinxAtStartPar
\sphinxstylestrong{Synopsis}
\begin{quote}

\sphinxAtStartPar
\sphinxcode{\sphinxupquote{getDim()}}
\end{quote}

\sphinxAtStartPar
\sphinxstylestrong{Description}
\begin{quote}

\sphinxAtStartPar
Retrieve the dimension of symmetric matrix.
\end{quote}

\sphinxAtStartPar
\sphinxstylestrong{Example}
\end{quote}

\begin{sphinxVerbatim}[commandchars=\\\{\}]
\PYG{c+c1}{\PYGZsh{} Retrieve the dimension of symmetric matrix \PYGZdq{}mat\PYGZdq{}.}
\PYG{n}{matdim} \PYG{o}{=} \PYG{n}{mat}\PYG{o}{.}\PYG{n}{getDim}\PYG{p}{(}\PYG{p}{)}
\end{sphinxVerbatim}

\subsection{SymMatrixArray Class}
\label{\detokenize{pyapiref:symmatrixarray-class}}\label{\detokenize{pyapiref:chappyapi-symmatrixarray}}
\sphinxAtStartPar
To facilitate users to operate on multiple {\hyperref[\detokenize{pyapiref:chappyapi-symmatrix}]{\sphinxcrossref{\DUrole{std,std-ref}{SymMatrix Class}}}} objects,
the Python interface of COPT provides SymMatrixArray object with the following methods:

\subsubsection{SymMatrixArray()}
\label{\detokenize{pyapiref:symmatrixarray}}\begin{quote}

\sphinxAtStartPar
\sphinxstylestrong{Synopsis}
\begin{quote}

\sphinxAtStartPar
\sphinxcode{\sphinxupquote{SymMatrixArray(mats=None)}}
\end{quote}

\sphinxAtStartPar
\sphinxstylestrong{Description}
\begin{quote}

\sphinxAtStartPar
Create a {\hyperref[\detokenize{pyapiref:chappyapi-symmatrixarray}]{\sphinxcrossref{\DUrole{std,std-ref}{SymMatrixArray Class}}}} object.

\sphinxAtStartPar
If parameter \sphinxcode{\sphinxupquote{mats}} is \sphinxcode{\sphinxupquote{None}}, then create an empty {\hyperref[\detokenize{pyapiref:chappyapi-symmatrixarray}]{\sphinxcrossref{\DUrole{std,std-ref}{SymMatrixArray Class}}}} object,
otherwise initialize the new created {\hyperref[\detokenize{pyapiref:chappyapi-symmatrixarray}]{\sphinxcrossref{\DUrole{std,std-ref}{SymMatrixArray Class}}}} object based on \sphinxcode{\sphinxupquote{mats}}.
\end{quote}

\sphinxAtStartPar
\sphinxstylestrong{Arguments}
\begin{quote}

\sphinxAtStartPar
\sphinxcode{\sphinxupquote{mats}}
\begin{quote}

\sphinxAtStartPar
\sphinxcode{\sphinxupquote{mats}} can be {\hyperref[\detokenize{pyapiref:chappyapi-symmatrix}]{\sphinxcrossref{\DUrole{std,std-ref}{SymMatrix Class}}}} object, {\hyperref[\detokenize{pyapiref:chappyapi-symmatrixarray}]{\sphinxcrossref{\DUrole{std,std-ref}{SymMatrixArray Class}}}} object,
list, dictionary or {\hyperref[\detokenize{pyapiref:chappyapi-util-tupledict}]{\sphinxcrossref{\DUrole{std,std-ref}{tupledict Class}}}} object.
\end{quote}
\end{quote}

\sphinxAtStartPar
\sphinxstylestrong{Example}
\end{quote}

\begin{sphinxVerbatim}[commandchars=\\\{\}]
\PYG{c+c1}{\PYGZsh{} Create an empty SymMatrixArray object}
\PYG{n}{matarr} \PYG{o}{=} \PYG{n}{SymMatrixArray}\PYG{p}{(}\PYG{p}{)}
\PYG{c+c1}{\PYGZsh{} Create a SymMatrixArray object containing matx, maty.}
\PYG{n}{matarr} \PYG{o}{=} \PYG{n}{SymMatrixArray}\PYG{p}{(}\PYG{p}{[}\PYG{n}{matx}\PYG{p}{,} \PYG{n}{maty}\PYG{p}{]}\PYG{p}{)}
\end{sphinxVerbatim}

\subsubsection{SymMatrixArray.pushBack()}
\label{\detokenize{pyapiref:symmatrixarray-pushback}}\begin{quote}

\sphinxAtStartPar
\sphinxstylestrong{Synopsis}
\begin{quote}

\sphinxAtStartPar
\sphinxcode{\sphinxupquote{pushBack(mat)}}
\end{quote}

\sphinxAtStartPar
\sphinxstylestrong{Description}
\begin{quote}

\sphinxAtStartPar
Add single or multiple {\hyperref[\detokenize{pyapiref:chappyapi-symmatrix}]{\sphinxcrossref{\DUrole{std,std-ref}{SymMatrix Class}}}} objects.
\end{quote}

\sphinxAtStartPar
\sphinxstylestrong{Arguments}
\begin{quote}

\sphinxAtStartPar
\sphinxcode{\sphinxupquote{mat}}
\begin{quote}

\sphinxAtStartPar
Symmetric matrices to be applied.
\sphinxcode{\sphinxupquote{mat}} can be {\hyperref[\detokenize{pyapiref:chappyapi-symmatrix}]{\sphinxcrossref{\DUrole{std,std-ref}{SymMatrix Class}}}} object, {\hyperref[\detokenize{pyapiref:chappyapi-symmatrixarray}]{\sphinxcrossref{\DUrole{std,std-ref}{SymMatrixArray Class}}}} object,
list, dictionary or {\hyperref[\detokenize{pyapiref:chappyapi-util-tupledict}]{\sphinxcrossref{\DUrole{std,std-ref}{tupledict Class}}}} object.
\end{quote}
\end{quote}

\sphinxAtStartPar
\sphinxstylestrong{Example}
\end{quote}

\begin{sphinxVerbatim}[commandchars=\\\{\}]
\PYG{c+c1}{\PYGZsh{} Add symmetric matrix matx to matarr}
\PYG{n}{matarr}\PYG{o}{.}\PYG{n}{pushBack}\PYG{p}{(}\PYG{n}{matx}\PYG{p}{)}
\PYG{c+c1}{\PYGZsh{} Add symmetric matrices matx and maty to matarr}
\PYG{n}{matarr}\PYG{o}{.}\PYG{n}{pushBack}\PYG{p}{(}\PYG{p}{[}\PYG{n}{matx}\PYG{p}{,} \PYG{n}{maty}\PYG{p}{]}\PYG{p}{)}
\end{sphinxVerbatim}

\subsubsection{SymMatrixArray.getMatrix()}
\label{\detokenize{pyapiref:symmatrixarray-getmatrix}}\begin{quote}

\sphinxAtStartPar
\sphinxstylestrong{Synopsis}
\begin{quote}

\sphinxAtStartPar
\sphinxcode{\sphinxupquote{getMatrix(idx)}}
\end{quote}

\sphinxAtStartPar
\sphinxstylestrong{Description}
\begin{quote}

\sphinxAtStartPar
Retrieve a symmetric matrix from an index in a {\hyperref[\detokenize{pyapiref:chappyapi-symmatrixarray}]{\sphinxcrossref{\DUrole{std,std-ref}{SymMatrixArray Class}}}} object.
Return a {\hyperref[\detokenize{pyapiref:chappyapi-symmatrix}]{\sphinxcrossref{\DUrole{std,std-ref}{SymMatrix Class}}}} object.
\end{quote}

\sphinxAtStartPar
\sphinxstylestrong{Arguments}
\begin{quote}

\sphinxAtStartPar
\sphinxcode{\sphinxupquote{idx}}
\begin{quote}

\sphinxAtStartPar
Subscript of the specified symmetric matrix in {\hyperref[\detokenize{pyapiref:chappyapi-symmatrixarray}]{\sphinxcrossref{\DUrole{std,std-ref}{SymMatrixArray Class}}}} object, starting with 0.
\end{quote}
\end{quote}

\sphinxAtStartPar
\sphinxstylestrong{Example}
\end{quote}

\begin{sphinxVerbatim}[commandchars=\\\{\}]
\PYG{c+c1}{\PYGZsh{} Get the symmetric matrix with subscript of 1 in matarr}
\PYG{n}{mat} \PYG{o}{=} \PYG{n}{matarr}\PYG{o}{.}\PYG{n}{getMatrix}\PYG{p}{(}\PYG{l+m+mi}{1}\PYG{p}{)}
\end{sphinxVerbatim}

\subsubsection{SymMatrixArray.getSize()}
\label{\detokenize{pyapiref:symmatrixarray-getsize}}\begin{quote}

\sphinxAtStartPar
\sphinxstylestrong{Synopsis}
\begin{quote}

\sphinxAtStartPar
\sphinxcode{\sphinxupquote{getSize()}}
\end{quote}

\sphinxAtStartPar
\sphinxstylestrong{Description}
\begin{quote}

\sphinxAtStartPar
Retrieve the number of symmetric matrices in {\hyperref[\detokenize{pyapiref:chappyapi-symmatrixarray}]{\sphinxcrossref{\DUrole{std,std-ref}{SymMatrixArray Class}}}} object.
\end{quote}

\sphinxAtStartPar
\sphinxstylestrong{Example}
\end{quote}

\begin{sphinxVerbatim}[commandchars=\\\{\}]
\PYG{c+c1}{\PYGZsh{} Retrieve the number of symmetric matrices in matarr.}
\PYG{n}{arrsize} \PYG{o}{=} \PYG{n}{matarr}\PYG{o}{.}\PYG{n}{getSize}\PYG{p}{(}\PYG{p}{)}
\end{sphinxVerbatim}

\subsection{Constraint Class}
\label{\detokenize{pyapiref:constraint-class}}\label{\detokenize{pyapiref:chappyapi-constraint}}
\sphinxAtStartPar
For easy access to information of constraints, Constraint class provides methods such as \sphinxcode{\sphinxupquote{Constraint.LB}}.
The supported information can be found in {\hyperref[\detokenize{pyapiref:chappyapi-const-info}]{\sphinxcrossref{\DUrole{std,std-ref}{Information}}}} section.
For convenience, information can be queried by names in original case or lowercase.

\sphinxAtStartPar
In addition, you can access the name of the constraint through \sphinxcode{\sphinxupquote{Constraint.name}}, the dual value of the constraint in LP through \sphinxcode{\sphinxupquote{Constraint.pi}},
the basis status of the constraint through \sphinxcode{\sphinxupquote{Constraint.basis}}, and the index in the coefficient matrix through \sphinxcode{\sphinxupquote{Constraint.index}}.

\sphinxAtStartPar
For the model\sphinxhyphen{}related information and constraint name, the user can also set the corresponding information in the form of \sphinxcode{\sphinxupquote{"Constraint.lb = \sphinxhyphen{}100"}}.

\sphinxAtStartPar
Constraint object contains related operations of COPT constraints and provides the following methods:

\subsubsection{Constraint.getName()}
\label{\detokenize{pyapiref:constraint-getname}}\begin{quote}

\sphinxAtStartPar
\sphinxstylestrong{Synopsis}
\begin{quote}

\sphinxAtStartPar
\sphinxcode{\sphinxupquote{getName()}}
\end{quote}

\sphinxAtStartPar
\sphinxstylestrong{Description}
\begin{quote}

\sphinxAtStartPar
Retrieve the name of linear constraint.
\end{quote}

\sphinxAtStartPar
\sphinxstylestrong{Example}
\end{quote}

\begin{sphinxVerbatim}[commandchars=\\\{\}]
\PYG{c+c1}{\PYGZsh{} Retrieve the name of linear constraint \PYGZsq{}con\PYGZsq{}.}
\PYG{n}{conname} \PYG{o}{=} \PYG{n}{con}\PYG{o}{.}\PYG{n}{getName}\PYG{p}{(}\PYG{p}{)}
\end{sphinxVerbatim}

\subsubsection{Constraint.getBasis()}
\label{\detokenize{pyapiref:constraint-getbasis}}\begin{quote}

\sphinxAtStartPar
\sphinxstylestrong{Synopsis}
\begin{quote}

\sphinxAtStartPar
\sphinxcode{\sphinxupquote{getBasis()}}
\end{quote}

\sphinxAtStartPar
\sphinxstylestrong{Description}
\begin{quote}

\sphinxAtStartPar
Retrieve the basis status of linear constraint.
\end{quote}

\sphinxAtStartPar
\sphinxstylestrong{Example}
\end{quote}

\begin{sphinxVerbatim}[commandchars=\\\{\}]
\PYG{c+c1}{\PYGZsh{} Retrieve the basis status of linear constraint \PYGZsq{}con\PYGZsq{}.}
\PYG{n}{conbasis} \PYG{o}{=} \PYG{n}{con}\PYG{o}{.}\PYG{n}{getBasis}\PYG{p}{(}\PYG{p}{)}
\end{sphinxVerbatim}

\subsubsection{Constraint.getLowerIIS()}
\label{\detokenize{pyapiref:constraint-getloweriis}}\begin{quote}

\sphinxAtStartPar
\sphinxstylestrong{Synopsis}
\begin{quote}

\sphinxAtStartPar
\sphinxcode{\sphinxupquote{getLowerIIS()}}
\end{quote}

\sphinxAtStartPar
\sphinxstylestrong{Description}
\begin{quote}

\sphinxAtStartPar
Retrieve the IIS status of lower bound of linear constraint.
\end{quote}

\sphinxAtStartPar
\sphinxstylestrong{Example}
\end{quote}

\begin{sphinxVerbatim}[commandchars=\\\{\}]
\PYG{c+c1}{\PYGZsh{} Retrieve the IIS status of lower bound of linear constraint \PYGZsq{}con\PYGZsq{}.}
\PYG{n}{lowerIIS} \PYG{o}{=} \PYG{n}{con}\PYG{o}{.}\PYG{n}{getLowerIIS}\PYG{p}{(}\PYG{p}{)}
\end{sphinxVerbatim}

\subsubsection{Constraint.getUpperIIS()}
\label{\detokenize{pyapiref:constraint-getupperiis}}\begin{quote}

\sphinxAtStartPar
\sphinxstylestrong{Synopsis}
\begin{quote}

\sphinxAtStartPar
\sphinxcode{\sphinxupquote{getUpperIIS()}}
\end{quote}

\sphinxAtStartPar
\sphinxstylestrong{Description}
\begin{quote}

\sphinxAtStartPar
Retrieve the IIS status of upper bound of linear constraint.
\end{quote}

\sphinxAtStartPar
\sphinxstylestrong{Example}
\end{quote}

\begin{sphinxVerbatim}[commandchars=\\\{\}]
\PYG{c+c1}{\PYGZsh{} Retrieve the IIS status of upper bound of linear constraint \PYGZsq{}con\PYGZsq{}.}
\PYG{n}{upperIIS} \PYG{o}{=} \PYG{n}{con}\PYG{o}{.}\PYG{n}{getUpperIIS}\PYG{p}{(}\PYG{p}{)}
\end{sphinxVerbatim}

\subsubsection{Constraint.getIdx()}
\label{\detokenize{pyapiref:constraint-getidx}}\begin{quote}

\sphinxAtStartPar
\sphinxstylestrong{Synopsis}
\begin{quote}

\sphinxAtStartPar
\sphinxcode{\sphinxupquote{getIdx()}}
\end{quote}

\sphinxAtStartPar
\sphinxstylestrong{Description}
\begin{quote}

\sphinxAtStartPar
Retrieve the subscript of linear constraint in coefficient matrix.
\end{quote}

\sphinxAtStartPar
\sphinxstylestrong{Example}
\end{quote}

\begin{sphinxVerbatim}[commandchars=\\\{\}]
\PYG{c+c1}{\PYGZsh{} Retrieve the subscript of linear constraint con.}
\PYG{n}{conidx} \PYG{o}{=} \PYG{n}{con}\PYG{o}{.}\PYG{n}{getIdx}\PYG{p}{(}\PYG{p}{)}
\end{sphinxVerbatim}

\subsubsection{Constraint.setName()}
\label{\detokenize{pyapiref:constraint-setname}}\begin{quote}

\sphinxAtStartPar
\sphinxstylestrong{Synopsis}
\begin{quote}

\sphinxAtStartPar
\sphinxcode{\sphinxupquote{setName(newname)}}
\end{quote}

\sphinxAtStartPar
\sphinxstylestrong{Description}
\begin{quote}

\sphinxAtStartPar
Set the name of linear constraint.
\end{quote}

\sphinxAtStartPar
\sphinxstylestrong{Arguments}
\begin{quote}

\sphinxAtStartPar
\sphinxcode{\sphinxupquote{newname}}
\begin{quote}

\sphinxAtStartPar
The name of constraint to be set.
\end{quote}
\end{quote}

\sphinxAtStartPar
\sphinxstylestrong{Example}
\end{quote}

\begin{sphinxVerbatim}[commandchars=\\\{\}]
\PYG{c+c1}{\PYGZsh{} Set the name of linear constraint \PYGZsq{}con\PYGZsq{}.}
\PYG{n}{con}\PYG{o}{.}\PYG{n}{setName}\PYG{p}{(}\PYG{l+s+s1}{\PYGZsq{}}\PYG{l+s+s1}{con}\PYG{l+s+s1}{\PYGZsq{}}\PYG{p}{)}
\end{sphinxVerbatim}

\subsubsection{Constraint.getInfo()}
\label{\detokenize{pyapiref:constraint-getinfo}}\begin{quote}

\sphinxAtStartPar
\sphinxstylestrong{Synopsis}
\begin{quote}

\sphinxAtStartPar
\sphinxcode{\sphinxupquote{getInfo(infoname)}}
\end{quote}

\sphinxAtStartPar
\sphinxstylestrong{Description}
\begin{quote}

\sphinxAtStartPar
Retrieve specified information. Return a constant.
\end{quote}

\sphinxAtStartPar
\sphinxstylestrong{Arguments}
\begin{quote}

\sphinxAtStartPar
\sphinxcode{\sphinxupquote{infoname}}
\begin{quote}

\sphinxAtStartPar
Name of the information to be obtained. Please refer to {\hyperref[\detokenize{pyapiref:chappyapi-const-info}]{\sphinxcrossref{\DUrole{std,std-ref}{Information}}}} section for possible values.
\end{quote}
\end{quote}

\sphinxAtStartPar
\sphinxstylestrong{Example}
\end{quote}

\begin{sphinxVerbatim}[commandchars=\\\{\}]
\PYG{c+c1}{\PYGZsh{} Get the lower bound of linear constraint con}
\PYG{n}{conlb} \PYG{o}{=} \PYG{n}{con}\PYG{o}{.}\PYG{n}{getInfo}\PYG{p}{(}\PYG{n}{COPT}\PYG{o}{.}\PYG{n}{Info}\PYG{o}{.}\PYG{n}{LB}\PYG{p}{)}
\end{sphinxVerbatim}

\subsubsection{Constraint.setInfo()}
\label{\detokenize{pyapiref:constraint-setinfo}}\begin{quote}

\sphinxAtStartPar
\sphinxstylestrong{Synopsis}
\begin{quote}

\sphinxAtStartPar
\sphinxcode{\sphinxupquote{setInfo(infoname, newval)}}
\end{quote}

\sphinxAtStartPar
\sphinxstylestrong{Description}
\begin{quote}

\sphinxAtStartPar
Set new information value to the specified constraint.
\end{quote}

\sphinxAtStartPar
\sphinxstylestrong{Arguments}
\begin{quote}

\sphinxAtStartPar
\sphinxcode{\sphinxupquote{infoname}}
\begin{quote}

\sphinxAtStartPar
The name of the information to be set. Please refer to {\hyperref[\detokenize{pyapiref:chappyapi-const-info}]{\sphinxcrossref{\DUrole{std,std-ref}{Information}}}} section for possible values.
\end{quote}

\sphinxAtStartPar
\sphinxcode{\sphinxupquote{newval}}
\begin{quote}

\sphinxAtStartPar
New information value to be set.
\end{quote}
\end{quote}

\sphinxAtStartPar
\sphinxstylestrong{Example}
\end{quote}

\begin{sphinxVerbatim}[commandchars=\\\{\}]
\PYG{c+c1}{\PYGZsh{} Set the lower bound of linear constraint con}
\PYG{n}{con}\PYG{o}{.}\PYG{n}{setInfo}\PYG{p}{(}\PYG{n}{COPT}\PYG{o}{.}\PYG{n}{Info}\PYG{o}{.}\PYG{n}{LB}\PYG{p}{,} \PYG{l+m+mf}{1.0}\PYG{p}{)}
\end{sphinxVerbatim}

\subsubsection{Constraint.remove()}
\label{\detokenize{pyapiref:constraint-remove}}\begin{quote}

\sphinxAtStartPar
\sphinxstylestrong{Synopsis}
\begin{quote}

\sphinxAtStartPar
\sphinxcode{\sphinxupquote{remove()}}
\end{quote}

\sphinxAtStartPar
\sphinxstylestrong{Description}
\begin{quote}

\sphinxAtStartPar
Delete the linear constraint from model.
\end{quote}

\sphinxAtStartPar
\sphinxstylestrong{Example}
\end{quote}

\begin{sphinxVerbatim}[commandchars=\\\{\}]
\PYG{c+c1}{\PYGZsh{} Delete the linear constraint \PYGZsq{}conx\PYGZsq{}}
\PYG{n}{conx}\PYG{o}{.}\PYG{n}{remove}\PYG{p}{(}\PYG{p}{)}
\end{sphinxVerbatim}

\subsection{ConstrArray Class}
\label{\detokenize{pyapiref:constrarray-class}}\label{\detokenize{pyapiref:chappyapi-constrarray}}
\sphinxAtStartPar
To facilitate users to operate on multiple {\hyperref[\detokenize{pyapiref:chappyapi-constraint}]{\sphinxcrossref{\DUrole{std,std-ref}{Constraint Class}}}}  objects,
the Python interface of COPT provides ConstrArray class with the following methods:

\subsubsection{ConstrArray()}
\label{\detokenize{pyapiref:constrarray}}\begin{quote}

\sphinxAtStartPar
\sphinxstylestrong{Synopsis}
\begin{quote}

\sphinxAtStartPar
\sphinxcode{\sphinxupquote{ConstrArray(constrs=None)}}
\end{quote}

\sphinxAtStartPar
\sphinxstylestrong{Description}
\begin{quote}

\sphinxAtStartPar
Create a {\hyperref[\detokenize{pyapiref:chappyapi-constrarray}]{\sphinxcrossref{\DUrole{std,std-ref}{ConstrArray Class}}}} object.

\sphinxAtStartPar
If parameter \sphinxcode{\sphinxupquote{constrs}} is \sphinxcode{\sphinxupquote{None}}, the create an empty {\hyperref[\detokenize{pyapiref:chappyapi-constrarray}]{\sphinxcrossref{\DUrole{std,std-ref}{ConstrArray Class}}}} object,
otherwise initialize the newly created {\hyperref[\detokenize{pyapiref:chappyapi-constrarray}]{\sphinxcrossref{\DUrole{std,std-ref}{ConstrArray Class}}}} object with parameter \sphinxcode{\sphinxupquote{constrs}}
\end{quote}

\sphinxAtStartPar
\sphinxstylestrong{Arguments}
\begin{quote}

\sphinxAtStartPar
\sphinxcode{\sphinxupquote{constrs}}
\begin{quote}

\sphinxAtStartPar
Linear constraints to be added. \sphinxcode{\sphinxupquote{None}} by default.

\sphinxAtStartPar
\sphinxcode{\sphinxupquote{constrs}} can be {\hyperref[\detokenize{pyapiref:chappyapi-constraint}]{\sphinxcrossref{\DUrole{std,std-ref}{Constraint Class}}}} object,
{\hyperref[\detokenize{pyapiref:chappyapi-constrarray}]{\sphinxcrossref{\DUrole{std,std-ref}{ConstrArray Class}}}} object, list, dictionary or {\hyperref[\detokenize{pyapiref:chappyapi-util-tupledict}]{\sphinxcrossref{\DUrole{std,std-ref}{tupledict Class}}}} object.
\end{quote}
\end{quote}

\sphinxAtStartPar
\sphinxstylestrong{Example}
\end{quote}

\begin{sphinxVerbatim}[commandchars=\\\{\}]
\PYG{c+c1}{\PYGZsh{} Create an empty ConstrArray object}
\PYG{n}{conarr} \PYG{o}{=} \PYG{n}{ConstrArray}\PYG{p}{(}\PYG{p}{)}
\PYG{c+c1}{\PYGZsh{} Create an ConstrArray object initialized with linear constraint conx and cony}
\PYG{n}{conarr} \PYG{o}{=} \PYG{n}{ConstrArray}\PYG{p}{(}\PYG{p}{[}\PYG{n}{conx}\PYG{p}{,} \PYG{n}{cony}\PYG{p}{]}\PYG{p}{)}
\end{sphinxVerbatim}

\subsubsection{ConstrArray.pushBack()}
\label{\detokenize{pyapiref:constrarray-pushback}}\begin{quote}

\sphinxAtStartPar
\sphinxstylestrong{Synopsis}
\begin{quote}

\sphinxAtStartPar
\sphinxcode{\sphinxupquote{pushBack(constrs)}}
\end{quote}

\sphinxAtStartPar
\sphinxstylestrong{Description}
\begin{quote}

\sphinxAtStartPar
Add single or multiple {\hyperref[\detokenize{pyapiref:chappyapi-constraint}]{\sphinxcrossref{\DUrole{std,std-ref}{Constraint Class}}}} objects.
\end{quote}

\sphinxAtStartPar
\sphinxstylestrong{Arguments}
\begin{quote}

\sphinxAtStartPar
\sphinxcode{\sphinxupquote{constrs}}
\begin{quote}

\sphinxAtStartPar
Linear constraints to be applied.
\sphinxcode{\sphinxupquote{constrs}} can be {\hyperref[\detokenize{pyapiref:chappyapi-constraint}]{\sphinxcrossref{\DUrole{std,std-ref}{Constraint Class}}}} object, {\hyperref[\detokenize{pyapiref:chappyapi-constrarray}]{\sphinxcrossref{\DUrole{std,std-ref}{ConstrArray Class}}}} object,
list, dictionary or {\hyperref[\detokenize{pyapiref:chappyapi-util-tupledict}]{\sphinxcrossref{\DUrole{std,std-ref}{tupledict Class}}}} object.
\end{quote}
\end{quote}

\sphinxAtStartPar
\sphinxstylestrong{Example}
\end{quote}

\begin{sphinxVerbatim}[commandchars=\\\{\}]
\PYG{c+c1}{\PYGZsh{} Add linear constraint r to conarr}
\PYG{n}{conarr}\PYG{o}{.}\PYG{n}{pushBack}\PYG{p}{(}\PYG{n}{r}\PYG{p}{)}
\PYG{c+c1}{\PYGZsh{} Add linear constraint r0 and r1 to conarr}
\PYG{n}{conarr}\PYG{o}{.}\PYG{n}{pushBack}\PYG{p}{(}\PYG{p}{[}\PYG{n}{r0}\PYG{p}{,} \PYG{n}{r1}\PYG{p}{]}\PYG{p}{)}
\end{sphinxVerbatim}

\subsubsection{ConstrArray.getConstr()}
\label{\detokenize{pyapiref:constrarray-getconstr}}\begin{quote}

\sphinxAtStartPar
\sphinxstylestrong{Synopsis}
\begin{quote}

\sphinxAtStartPar
\sphinxcode{\sphinxupquote{getConstr(idx)}}
\end{quote}

\sphinxAtStartPar
\sphinxstylestrong{Description}
\begin{quote}

\sphinxAtStartPar
Retrieve the linear constraint according to its subscript in {\hyperref[\detokenize{pyapiref:chappyapi-constrarray}]{\sphinxcrossref{\DUrole{std,std-ref}{ConstrArray Class}}}} object.
Return a {\hyperref[\detokenize{pyapiref:chappyapi-constraint}]{\sphinxcrossref{\DUrole{std,std-ref}{Constraint Class}}}} object.
\end{quote}

\sphinxAtStartPar
\sphinxstylestrong{Arguments}
\begin{quote}

\sphinxAtStartPar
\sphinxcode{\sphinxupquote{idx}}
\begin{quote}

\sphinxAtStartPar
Subscript of the desired constraint in {\hyperref[\detokenize{pyapiref:chappyapi-constrarray}]{\sphinxcrossref{\DUrole{std,std-ref}{ConstrArray Class}}}} object, starting with 0.
\end{quote}
\end{quote}

\sphinxAtStartPar
\sphinxstylestrong{Example}
\end{quote}

\begin{sphinxVerbatim}[commandchars=\\\{\}]
\PYG{c+c1}{\PYGZsh{} Retrieve the linear constraint with subscript 1  in conarr}
\PYG{n}{conarr}\PYG{o}{.}\PYG{n}{getConstr}\PYG{p}{(}\PYG{l+m+mi}{1}\PYG{p}{)}
\end{sphinxVerbatim}

\subsubsection{ConstrArray.getAll()}
\label{\detokenize{pyapiref:constrarray-getall}}\begin{quote}

\sphinxAtStartPar
\sphinxstylestrong{Synopsis}
\begin{quote}

\sphinxAtStartPar
\sphinxcode{\sphinxupquote{getAll()}}
\end{quote}

\sphinxAtStartPar
\sphinxstylestrong{Description}
\begin{quote}

\sphinxAtStartPar
Retrieve all linear constraints in {\hyperref[\detokenize{pyapiref:chappyapi-constrarray}]{\sphinxcrossref{\DUrole{std,std-ref}{ConstrArray Class}}}} object. Returns a list object.
\end{quote}

\sphinxAtStartPar
\sphinxstylestrong{Example}
\end{quote}

\begin{sphinxVerbatim}[commandchars=\\\{\}]
\PYG{c+c1}{\PYGZsh{} Get all linear constraints in \PYGZsq{}conarr\PYGZsq{}}
\PYG{n}{cons} \PYG{o}{=} \PYG{n}{conarr}\PYG{o}{.}\PYG{n}{getAll}\PYG{p}{(}\PYG{p}{)}
\end{sphinxVerbatim}

\subsubsection{ConstrArray.getSize()}
\label{\detokenize{pyapiref:constrarray-getsize}}\begin{quote}

\sphinxAtStartPar
\sphinxstylestrong{Synopsis}
\begin{quote}

\sphinxAtStartPar
\sphinxcode{\sphinxupquote{getSize()}}
\end{quote}

\sphinxAtStartPar
\sphinxstylestrong{Description}
\begin{quote}

\sphinxAtStartPar
Get the number of elements in {\hyperref[\detokenize{pyapiref:chappyapi-constrarray}]{\sphinxcrossref{\DUrole{std,std-ref}{ConstrArray Class}}}} object.
\end{quote}

\sphinxAtStartPar
\sphinxstylestrong{Example}
\end{quote}

\begin{sphinxVerbatim}[commandchars=\\\{\}]
\PYG{c+c1}{\PYGZsh{} Get the number of linear constraints in conarr}
\PYG{n}{arrsize} \PYG{o}{=} \PYG{n}{conarr}\PYG{o}{.}\PYG{n}{getSize}\PYG{p}{(}\PYG{p}{)}
\end{sphinxVerbatim}

\subsection{ConstrBuilder Class}
\label{\detokenize{pyapiref:constrbuilder-class}}\label{\detokenize{pyapiref:chappyapi-constrbuilder}}
\sphinxAtStartPar
ConstrBuilder object contains operations related to temporary constraints when building constraints, and provides the following methods:

\subsubsection{ConstrBuilder()}
\label{\detokenize{pyapiref:constrbuilder}}\begin{quote}

\sphinxAtStartPar
\sphinxstylestrong{Synopsis}
\begin{quote}

\sphinxAtStartPar
\sphinxcode{\sphinxupquote{ConstrBuilder()}}
\end{quote}

\sphinxAtStartPar
\sphinxstylestrong{Description}
\begin{quote}

\sphinxAtStartPar
Create an empty {\hyperref[\detokenize{pyapiref:chappyapi-constrbuilder}]{\sphinxcrossref{\DUrole{std,std-ref}{ConstrBuilder Class}}}} object
\end{quote}

\sphinxAtStartPar
\sphinxstylestrong{Example}
\end{quote}

\begin{sphinxVerbatim}[commandchars=\\\{\}]
\PYG{c+c1}{\PYGZsh{} Create an empty linear constraint builder}
\PYG{n}{constrbuilder} \PYG{o}{=} \PYG{n}{ConstrBuilder}\PYG{p}{(}\PYG{p}{)}
\end{sphinxVerbatim}

\subsubsection{ConstrBuilder.setBuilder()}
\label{\detokenize{pyapiref:constrbuilder-setbuilder}}\begin{quote}

\sphinxAtStartPar
\sphinxstylestrong{Synopsis}
\begin{quote}

\sphinxAtStartPar
\sphinxcode{\sphinxupquote{setBuilder(expr, sense, rhs)}}
\end{quote}

\sphinxAtStartPar
\sphinxstylestrong{Description}
\begin{quote}

\sphinxAtStartPar
Set expression and constraint type for linear constraint builder.
\end{quote}

\sphinxAtStartPar
\sphinxstylestrong{Arguments}
\begin{quote}

\sphinxAtStartPar
\sphinxcode{\sphinxupquote{expr}}
\begin{quote}

\sphinxAtStartPar
The expression to be set, which can be {\hyperref[\detokenize{pyapiref:chappyapi-var}]{\sphinxcrossref{\DUrole{std,std-ref}{Var Class}}}} expression or
{\hyperref[\detokenize{pyapiref:chappyapi-linexpr}]{\sphinxcrossref{\DUrole{std,std-ref}{LinExpr Class}}}} expression.
\end{quote}

\sphinxAtStartPar
\sphinxcode{\sphinxupquote{sense}}
\begin{quote}

\sphinxAtStartPar
Sense of constraint. The full list of available tyeps can be found in {\hyperref[\detokenize{constant:chapconst-constrtype}]{\sphinxcrossref{\DUrole{std,std-ref}{Constraint type}}}} section.
\end{quote}

\sphinxAtStartPar
\sphinxcode{\sphinxupquote{rhs}}
\begin{quote}

\sphinxAtStartPar
Right hand side of constraint.
\end{quote}
\end{quote}

\sphinxAtStartPar
\sphinxstylestrong{Example}
\end{quote}

\begin{sphinxVerbatim}[commandchars=\\\{\}]
\PYG{c+c1}{\PYGZsh{} Set the expresson of linear constraint builder as: x+y==1}
\PYG{n}{constrbuilder}\PYG{o}{.}\PYG{n}{setBuilder}\PYG{p}{(}\PYG{n}{x} \PYG{o}{+} \PYG{n}{y}\PYG{p}{,} \PYG{n}{COPT}\PYG{o}{.}\PYG{n}{EQUAL}\PYG{p}{,} \PYG{l+m+mi}{1}\PYG{p}{)}
\end{sphinxVerbatim}

\subsubsection{ConstrBuilder.getExpr()}
\label{\detokenize{pyapiref:constrbuilder-getexpr}}\begin{quote}

\sphinxAtStartPar
\sphinxstylestrong{Synopsis}
\begin{quote}

\sphinxAtStartPar
\sphinxcode{\sphinxupquote{getExpr()}}
\end{quote}

\sphinxAtStartPar
\sphinxstylestrong{Description}
\begin{quote}

\sphinxAtStartPar
Retrieve the expression of a linear constraint builder object.
\end{quote}

\sphinxAtStartPar
\sphinxstylestrong{Example}
\end{quote}

\begin{sphinxVerbatim}[commandchars=\\\{\}]
\PYG{c+c1}{\PYGZsh{}  Retrieve the expression of a linear constraint builder}
\PYG{n}{linexpr} \PYG{o}{=} \PYG{n}{constrbuilder}\PYG{o}{.}\PYG{n}{getExpr}\PYG{p}{(}\PYG{p}{)}
\end{sphinxVerbatim}

\subsubsection{ConstrBuilder.getSense()}
\label{\detokenize{pyapiref:constrbuilder-getsense}}\begin{quote}

\sphinxAtStartPar
\sphinxstylestrong{Synopsis}
\begin{quote}

\sphinxAtStartPar
\sphinxcode{\sphinxupquote{getSense()}}
\end{quote}

\sphinxAtStartPar
\sphinxstylestrong{Description}
\begin{quote}

\sphinxAtStartPar
Retrieve the constraint sense of linear constraint builder object.
\end{quote}

\sphinxAtStartPar
\sphinxstylestrong{Example}
\end{quote}

\begin{sphinxVerbatim}[commandchars=\\\{\}]
\PYG{c+c1}{\PYGZsh{} Retrieve the constraint sense of linear constraint builder object.}
\PYG{n}{consense} \PYG{o}{=} \PYG{n}{constrbuilder}\PYG{o}{.}\PYG{n}{getSense}\PYG{p}{(}\PYG{p}{)}
\end{sphinxVerbatim}

\subsection{ConstrBuilderArray Class}
\label{\detokenize{pyapiref:constrbuilderarray-class}}\label{\detokenize{pyapiref:chappyapi-constrbuilderarray}}
\sphinxAtStartPar
To facilitate users to operate on multiple {\hyperref[\detokenize{pyapiref:chappyapi-constrbuilder}]{\sphinxcrossref{\DUrole{std,std-ref}{ConstrBuilder Class}}}} objects,
the Python interface of COPT provides ConstrArray object with the following methods:

\subsubsection{ConstrBuilderArray()}
\label{\detokenize{pyapiref:constrbuilderarray}}\begin{quote}

\sphinxAtStartPar
\sphinxstylestrong{Synopsis}
\begin{quote}

\sphinxAtStartPar
\sphinxcode{\sphinxupquote{ConstrBuilderArray(constrbuilders=None)}}
\end{quote}

\sphinxAtStartPar
\sphinxstylestrong{Description}
\begin{quote}

\sphinxAtStartPar
Create a {\hyperref[\detokenize{pyapiref:chappyapi-constrbuilderarray}]{\sphinxcrossref{\DUrole{std,std-ref}{ConstrBuilderArray Class}}}} object.

\sphinxAtStartPar
If parameter \sphinxcode{\sphinxupquote{constrbuilders}} is \sphinxcode{\sphinxupquote{None}}, then create an empty {\hyperref[\detokenize{pyapiref:chappyapi-constrbuilderarray}]{\sphinxcrossref{\DUrole{std,std-ref}{ConstrBuilderArray Class}}}} object,
otherwise initialize the newly created {\hyperref[\detokenize{pyapiref:chappyapi-constrbuilderarray}]{\sphinxcrossref{\DUrole{std,std-ref}{ConstrBuilderArray Class}}}} object by parameter \sphinxcode{\sphinxupquote{constrbuilders}}.
\end{quote}

\sphinxAtStartPar
\sphinxstylestrong{Arguments}
\begin{quote}
\begin{description}
\sphinxlineitem{\sphinxcode{\sphinxupquote{constrbuilders}}}
\sphinxAtStartPar
Linear constraint builder to be added. Optional, \sphinxcode{\sphinxupquote{None}} by default. It can be {\hyperref[\detokenize{pyapiref:chappyapi-constrbuilder}]{\sphinxcrossref{\DUrole{std,std-ref}{ConstrBuilder Class}}}} object,
{\hyperref[\detokenize{pyapiref:chappyapi-constrbuilderarray}]{\sphinxcrossref{\DUrole{std,std-ref}{ConstrBuilderArray Class}}}} object, list, dictionary or {\hyperref[\detokenize{pyapiref:chappyapi-util-tupledict}]{\sphinxcrossref{\DUrole{std,std-ref}{tupledict Class}}}} object.

\end{description}
\end{quote}

\sphinxAtStartPar
\sphinxstylestrong{Example}
\end{quote}

\begin{sphinxVerbatim}[commandchars=\\\{\}]
\PYG{c+c1}{\PYGZsh{} Create an empty ConstrBuilderArray object.}
\PYG{n}{conbuilderarr} \PYG{o}{=} \PYG{n}{ConstrBuilderArray}\PYG{p}{(}\PYG{p}{)}
\PYG{c+c1}{\PYGZsh{} Create a ConstrBuilderArray object and initialize it with builders: conbuilderx and conbuildery}
\PYG{n}{conbuilderarr} \PYG{o}{=} \PYG{n}{ConstrBuilderArray}\PYG{p}{(}\PYG{p}{[}\PYG{n}{conbuilderx}\PYG{p}{,} \PYG{n}{conbuildery}\PYG{p}{]}\PYG{p}{)}
\end{sphinxVerbatim}

\subsubsection{ConstrBuilderArray.pushBack()}
\label{\detokenize{pyapiref:constrbuilderarray-pushback}}\begin{quote}

\sphinxAtStartPar
\sphinxstylestrong{Synopsis}
\begin{quote}

\sphinxAtStartPar
\sphinxcode{\sphinxupquote{pushBack(constrbuilder)}}
\end{quote}

\sphinxAtStartPar
\sphinxstylestrong{Description}
\begin{quote}

\sphinxAtStartPar
Add single or multiple {\hyperref[\detokenize{pyapiref:chappyapi-constrbuilder}]{\sphinxcrossref{\DUrole{std,std-ref}{ConstrBuilder Class}}}} objects.
\end{quote}

\sphinxAtStartPar
\sphinxstylestrong{Arguments}
\begin{quote}

\sphinxAtStartPar
\sphinxcode{\sphinxupquote{constrbuilder}}
\begin{quote}

\sphinxAtStartPar
Builder of linear constraint to be added, which can be {\hyperref[\detokenize{pyapiref:chappyapi-constrbuilder}]{\sphinxcrossref{\DUrole{std,std-ref}{ConstrBuilder Class}}}} object,
{\hyperref[\detokenize{pyapiref:chappyapi-constrbuilderarray}]{\sphinxcrossref{\DUrole{std,std-ref}{ConstrBuilderArray Class}}}} object, list, dictionary or {\hyperref[\detokenize{pyapiref:chappyapi-util-tupledict}]{\sphinxcrossref{\DUrole{std,std-ref}{tupledict Class}}}} object.
\end{quote}
\end{quote}

\sphinxAtStartPar
\sphinxstylestrong{Example}
\end{quote}

\begin{sphinxVerbatim}[commandchars=\\\{\}]
\PYG{c+c1}{\PYGZsh{} Add linear constraint builder conbuilderx to conbuilderarr}
\PYG{n}{conbuilderarr}\PYG{o}{.}\PYG{n}{pushBack}\PYG{p}{(}\PYG{n}{conbuilderx}\PYG{p}{)}
\PYG{c+c1}{\PYGZsh{} Add linear constraint builders conbuilderx and conbuildery to conbuilderarr}
\PYG{n}{conbuilderarr}\PYG{o}{.}\PYG{n}{pushBack}\PYG{p}{(}\PYG{p}{[}\PYG{n}{conbuilderx}\PYG{p}{,} \PYG{n}{conbuildery}\PYG{p}{]}\PYG{p}{)}
\end{sphinxVerbatim}

\subsubsection{ConstrBuilderArray.getBuilder()}
\label{\detokenize{pyapiref:constrbuilderarray-getbuilder}}\begin{quote}

\sphinxAtStartPar
\sphinxstylestrong{Synopsis}
\begin{quote}

\sphinxAtStartPar
\sphinxcode{\sphinxupquote{getBuilder(idx)}}
\end{quote}

\sphinxAtStartPar
\sphinxstylestrong{Description}
\begin{quote}

\sphinxAtStartPar
Retrieve a temporary constraint from its index in {\hyperref[\detokenize{pyapiref:chappyapi-constrbuilderarray}]{\sphinxcrossref{\DUrole{std,std-ref}{ConstrBuilderArray Class}}}} object.
Return a {\hyperref[\detokenize{pyapiref:chappyapi-constrbuilder}]{\sphinxcrossref{\DUrole{std,std-ref}{ConstrBuilder Class}}}} object.

\sphinxAtStartPar
Retrieve the corresponding builder object according to the subscript of linear constraint builder in {\hyperref[\detokenize{pyapiref:chappyapi-constrbuilderarray}]{\sphinxcrossref{\DUrole{std,std-ref}{ConstrBuilderArray Class}}}} object.
\end{quote}

\sphinxAtStartPar
\sphinxstylestrong{Arguments}
\begin{quote}

\sphinxAtStartPar
\sphinxcode{\sphinxupquote{idx}}
\begin{quote}

\sphinxAtStartPar
Subscript of the linear constraint builder in the {\hyperref[\detokenize{pyapiref:chappyapi-constrbuilderarray}]{\sphinxcrossref{\DUrole{std,std-ref}{ConstrBuilderArray Class}}}} object, starting with 0.
\end{quote}
\end{quote}

\sphinxAtStartPar
\sphinxstylestrong{Example}
\end{quote}

\begin{sphinxVerbatim}[commandchars=\\\{\}]
\PYG{c+c1}{\PYGZsh{} Retrieve the builder with subscript 1 in conbuilderarr}
\PYG{n}{conbuilder} \PYG{o}{=} \PYG{n}{conbuilderarr}\PYG{o}{.}\PYG{n}{getBuilder}\PYG{p}{(}\PYG{l+m+mi}{1}\PYG{p}{)}
\end{sphinxVerbatim}

\subsubsection{ConstrBuilderArray.getSize()}
\label{\detokenize{pyapiref:constrbuilderarray-getsize}}\begin{quote}

\sphinxAtStartPar
\sphinxstylestrong{Synopsis}
\begin{quote}

\sphinxAtStartPar
\sphinxcode{\sphinxupquote{getSize()}}
\end{quote}

\sphinxAtStartPar
\sphinxstylestrong{Description}
\begin{quote}

\sphinxAtStartPar
Get the number of elements in {\hyperref[\detokenize{pyapiref:chappyapi-constrbuilderarray}]{\sphinxcrossref{\DUrole{std,std-ref}{ConstrBuilderArray Class}}}} object.
\end{quote}

\sphinxAtStartPar
\sphinxstylestrong{Example}
\end{quote}

\begin{sphinxVerbatim}[commandchars=\\\{\}]
\PYG{c+c1}{\PYGZsh{} Get the number of builders in conbuilderarr}
\PYG{n}{arrsize} \PYG{o}{=} \PYG{n}{conbuilderarr}\PYG{o}{.}\PYG{n}{getSize}\PYG{p}{(}\PYG{p}{)}
\end{sphinxVerbatim}

\subsection{QConstraint Class}
\label{\detokenize{pyapiref:qconstraint-class}}\label{\detokenize{pyapiref:chappyapi-qconstraint}}
\sphinxAtStartPar
For easy access to information of quadratic constraints, QConstraint class provides methods such as \sphinxcode{\sphinxupquote{QConstraint.index}}.
The supported information can be found in {\hyperref[\detokenize{pyapiref:chappyapi-const-info}]{\sphinxcrossref{\DUrole{std,std-ref}{Information}}}} section.
For convenience, information can be queried by names in original case or lowercase.

\sphinxAtStartPar
In addition, you can access the name of the quadratic constraint through \sphinxcode{\sphinxupquote{QConstraint.name}}, and the index in the model through \sphinxcode{\sphinxupquote{QConstraint.index}}.

\sphinxAtStartPar
For the model\sphinxhyphen{}related information and constraint name, the user can also set the corresponding information in the form of \sphinxcode{\sphinxupquote{"QConstraint.rhs = \sphinxhyphen{}100"}}.

\sphinxAtStartPar
QConstraint object contains related operations of COPT quadratic constraints and provides the following methods:

\subsubsection{QConstraint.getName()}
\label{\detokenize{pyapiref:qconstraint-getname}}\begin{quote}

\sphinxAtStartPar
\sphinxstylestrong{Synopsis}
\begin{quote}

\sphinxAtStartPar
\sphinxcode{\sphinxupquote{getName()}}
\end{quote}

\sphinxAtStartPar
\sphinxstylestrong{Description}
\begin{quote}

\sphinxAtStartPar
Retrieve the name of quadratic constraint.
\end{quote}

\sphinxAtStartPar
\sphinxstylestrong{Example}
\end{quote}

\begin{sphinxVerbatim}[commandchars=\\\{\}]
\PYG{c+c1}{\PYGZsh{} Retrieve the name of quadratic constraint \PYGZsq{}qcon\PYGZsq{}}
\PYG{n}{qconname} \PYG{o}{=} \PYG{n}{qcon}\PYG{o}{.}\PYG{n}{getName}\PYG{p}{(}\PYG{p}{)}
\end{sphinxVerbatim}

\subsubsection{QConstraint.getRhs()}
\label{\detokenize{pyapiref:qconstraint-getrhs}}\begin{quote}

\sphinxAtStartPar
\sphinxstylestrong{Synopsis}
\begin{quote}

\sphinxAtStartPar
\sphinxcode{\sphinxupquote{getRhs()}}
\end{quote}

\sphinxAtStartPar
\sphinxstylestrong{Description}
\begin{quote}

\sphinxAtStartPar
Retrieve the right hand side of quadratic constraint.
\end{quote}

\sphinxAtStartPar
\sphinxstylestrong{Example}
\end{quote}

\begin{sphinxVerbatim}[commandchars=\\\{\}]
\PYG{c+c1}{\PYGZsh{} Retrieve the RHS of quadratic constraint \PYGZsq{}qcon\PYGZsq{}}
\PYG{n}{qconrhs} \PYG{o}{=} \PYG{n}{qcon}\PYG{o}{.}\PYG{n}{getRhs}\PYG{p}{(}\PYG{p}{)}
\end{sphinxVerbatim}

\subsubsection{QConstraint.getSense()}
\label{\detokenize{pyapiref:qconstraint-getsense}}\begin{quote}

\sphinxAtStartPar
\sphinxstylestrong{Synopsis}
\begin{quote}

\sphinxAtStartPar
\sphinxcode{\sphinxupquote{getSense()}}
\end{quote}

\sphinxAtStartPar
\sphinxstylestrong{Description}
\begin{quote}

\sphinxAtStartPar
Retrieve the type of quadratic constraint.
\end{quote}

\sphinxAtStartPar
\sphinxstylestrong{Example}
\end{quote}

\begin{sphinxVerbatim}[commandchars=\\\{\}]
\PYG{c+c1}{\PYGZsh{} Retrieve the type of quadratic constraint \PYGZsq{}qcon\PYGZsq{}}
\PYG{n}{qconsense} \PYG{o}{=} \PYG{n}{qcon}\PYG{o}{.}\PYG{n}{getSense}\PYG{p}{(}\PYG{p}{)}
\end{sphinxVerbatim}

\subsubsection{QConstraint.getIdx()}
\label{\detokenize{pyapiref:qconstraint-getidx}}\begin{quote}

\sphinxAtStartPar
\sphinxstylestrong{Synopsis}
\begin{quote}

\sphinxAtStartPar
\sphinxcode{\sphinxupquote{getIdx()}}
\end{quote}

\sphinxAtStartPar
\sphinxstylestrong{Description}
\begin{quote}

\sphinxAtStartPar
Retrieve the subscript of quadratic constraint.
\end{quote}

\sphinxAtStartPar
\sphinxstylestrong{Example}
\end{quote}

\begin{sphinxVerbatim}[commandchars=\\\{\}]
\PYG{c+c1}{\PYGZsh{} Retrieve the subscript of quadratic constraint \PYGZsq{}qcon\PYGZsq{}}
\PYG{n}{qconidx} \PYG{o}{=} \PYG{n}{qcon}\PYG{o}{.}\PYG{n}{getIdx}\PYG{p}{(}\PYG{p}{)}
\end{sphinxVerbatim}

\subsubsection{QConstraint.setName()}
\label{\detokenize{pyapiref:qconstraint-setname}}\begin{quote}

\sphinxAtStartPar
\sphinxstylestrong{Synopsis}
\begin{quote}

\sphinxAtStartPar
\sphinxcode{\sphinxupquote{setName(newname)}}
\end{quote}

\sphinxAtStartPar
\sphinxstylestrong{Description}
\begin{quote}

\sphinxAtStartPar
Set the name of quadratic constraint.
\end{quote}

\sphinxAtStartPar
\sphinxstylestrong{Arguments}
\begin{quote}

\sphinxAtStartPar
\sphinxcode{\sphinxupquote{newname}}
\begin{quote}

\sphinxAtStartPar
The name of quadratic constraint to be set.
\end{quote}
\end{quote}

\sphinxAtStartPar
\sphinxstylestrong{Example}
\end{quote}

\begin{sphinxVerbatim}[commandchars=\\\{\}]
\PYG{c+c1}{\PYGZsh{} Set the name of quadratic constraint \PYGZsq{}qcon\PYGZsq{}.}
\PYG{n}{qcon}\PYG{o}{.}\PYG{n}{setName}\PYG{p}{(}\PYG{l+s+s1}{\PYGZsq{}}\PYG{l+s+s1}{qcon}\PYG{l+s+s1}{\PYGZsq{}}\PYG{p}{)}
\end{sphinxVerbatim}

\subsubsection{QConstraint.setRhs()}
\label{\detokenize{pyapiref:qconstraint-setrhs}}\begin{quote}

\sphinxAtStartPar
\sphinxstylestrong{Synopsis}
\begin{quote}

\sphinxAtStartPar
\sphinxcode{\sphinxupquote{setRhs(rhs)}}
\end{quote}

\sphinxAtStartPar
\sphinxstylestrong{Description}
\begin{quote}

\sphinxAtStartPar
Set the right hand side of quadratic constraint.
\end{quote}

\sphinxAtStartPar
\sphinxstylestrong{Arguments}
\begin{quote}

\sphinxAtStartPar
\sphinxcode{\sphinxupquote{rhs}}
\begin{quote}

\sphinxAtStartPar
The right hand side of quadratic constraint to be set.
\end{quote}
\end{quote}

\sphinxAtStartPar
\sphinxstylestrong{Example}
\end{quote}

\begin{sphinxVerbatim}[commandchars=\\\{\}]
\PYG{c+c1}{\PYGZsh{} Set the RHS of quadratic constraint \PYGZsq{}qcon\PYGZsq{} to 0.0}
\PYG{n}{qcon}\PYG{o}{.}\PYG{n}{setRhs}\PYG{p}{(}\PYG{l+m+mf}{0.0}\PYG{p}{)}
\end{sphinxVerbatim}

\subsubsection{QConstraint.setSense()}
\label{\detokenize{pyapiref:qconstraint-setsense}}\begin{quote}

\sphinxAtStartPar
\sphinxstylestrong{Synopsis}
\begin{quote}

\sphinxAtStartPar
\sphinxcode{\sphinxupquote{setSense(sense)}}
\end{quote}

\sphinxAtStartPar
\sphinxstylestrong{Description}
\begin{quote}

\sphinxAtStartPar
Set the sense of quadratic constraint.
\end{quote}

\sphinxAtStartPar
\sphinxstylestrong{Arguments}
\begin{quote}

\sphinxAtStartPar
\sphinxcode{\sphinxupquote{sense}}
\begin{quote}

\sphinxAtStartPar
The sense of quadratic constraint to be set.
\end{quote}
\end{quote}

\sphinxAtStartPar
\sphinxstylestrong{Example}
\end{quote}

\begin{sphinxVerbatim}[commandchars=\\\{\}]
\PYG{c+c1}{\PYGZsh{} Set the sense of quadratic constraint \PYGZsq{}qcon\PYGZsq{} to \PYGZlt{}=}
\PYG{n}{qcon}\PYG{o}{.}\PYG{n}{setSense}\PYG{p}{(}\PYG{n}{COPT}\PYG{o}{.}\PYG{n}{LESS\PYGZus{}EQUAL}\PYG{p}{)}
\end{sphinxVerbatim}

\subsubsection{QConstraint.getInfo()}
\label{\detokenize{pyapiref:qconstraint-getinfo}}\begin{quote}

\sphinxAtStartPar
\sphinxstylestrong{Synopsis}
\begin{quote}

\sphinxAtStartPar
\sphinxcode{\sphinxupquote{getInfo(infoname)}}
\end{quote}

\sphinxAtStartPar
\sphinxstylestrong{Description}
\begin{quote}

\sphinxAtStartPar
Retrieve specified information. Return a constant.
\end{quote}

\sphinxAtStartPar
\sphinxstylestrong{Arguments}
\begin{quote}

\sphinxAtStartPar
\sphinxcode{\sphinxupquote{infoname}}
\begin{quote}

\sphinxAtStartPar
Name of the information to be obtained.

\sphinxAtStartPar
Please refer to {\hyperref[\detokenize{pyapiref:chappyapi-const-info}]{\sphinxcrossref{\DUrole{std,std-ref}{Information}}}} section for possible values.
\end{quote}
\end{quote}

\sphinxAtStartPar
\sphinxstylestrong{Example}
\end{quote}

\begin{sphinxVerbatim}[commandchars=\\\{\}]
\PYG{c+c1}{\PYGZsh{} Get the row activity of quadratic constraint \PYGZsq{}qcon\PYGZsq{}}
\PYG{n}{qconlb} \PYG{o}{=} \PYG{n}{qcon}\PYG{o}{.}\PYG{n}{getInfo}\PYG{p}{(}\PYG{n}{COPT}\PYG{o}{.}\PYG{n}{Info}\PYG{o}{.}\PYG{n}{Slack}\PYG{p}{)}
\end{sphinxVerbatim}

\subsubsection{QConstraint.setInfo()}
\label{\detokenize{pyapiref:qconstraint-setinfo}}\begin{quote}

\sphinxAtStartPar
\sphinxstylestrong{Synopsis}
\begin{quote}

\sphinxAtStartPar
\sphinxcode{\sphinxupquote{setInfo(infoname, newval)}}
\end{quote}

\sphinxAtStartPar
\sphinxstylestrong{Description}
\begin{quote}

\sphinxAtStartPar
Set new information value to the specified quadratic constraint.
\end{quote}

\sphinxAtStartPar
\sphinxstylestrong{Arguments}
\begin{quote}

\sphinxAtStartPar
\sphinxcode{\sphinxupquote{infoname}}
\begin{quote}

\sphinxAtStartPar
The name of the information to be set. Please refer to {\hyperref[\detokenize{pyapiref:chappyapi-const-info}]{\sphinxcrossref{\DUrole{std,std-ref}{Information}}}} section for possible values.
\end{quote}

\sphinxAtStartPar
\sphinxcode{\sphinxupquote{newval}}
\begin{quote}

\sphinxAtStartPar
New information value to be set.
\end{quote}
\end{quote}

\sphinxAtStartPar
\sphinxstylestrong{Example}
\end{quote}

\begin{sphinxVerbatim}[commandchars=\\\{\}]
\PYG{c+c1}{\PYGZsh{} Set the lower bound of quadratic constraint \PYGZsq{}qcon\PYGZsq{}}
\PYG{n}{qcon}\PYG{o}{.}\PYG{n}{setInfo}\PYG{p}{(}\PYG{n}{COPT}\PYG{o}{.}\PYG{n}{Info}\PYG{o}{.}\PYG{n}{LB}\PYG{p}{,} \PYG{l+m+mf}{1.0}\PYG{p}{)}
\end{sphinxVerbatim}

\subsubsection{Constraint.remove()}
\label{\detokenize{pyapiref:id1}}\begin{quote}

\sphinxAtStartPar
\sphinxstylestrong{Synopsis}
\begin{quote}

\sphinxAtStartPar
\sphinxcode{\sphinxupquote{remove()}}
\end{quote}

\sphinxAtStartPar
\sphinxstylestrong{Description}
\begin{quote}

\sphinxAtStartPar
Delete the quadratic constraint from model.
\end{quote}

\sphinxAtStartPar
\sphinxstylestrong{Example}
\end{quote}

\begin{sphinxVerbatim}[commandchars=\\\{\}]
\PYG{c+c1}{\PYGZsh{} Delete the quadratic constraint \PYGZsq{}qconx\PYGZsq{}}
\PYG{n}{qconx}\PYG{o}{.}\PYG{n}{remove}\PYG{p}{(}\PYG{p}{)}
\end{sphinxVerbatim}

\subsection{QConstrArray Class}
\label{\detokenize{pyapiref:qconstrarray-class}}\label{\detokenize{pyapiref:chappyapi-qconstrarray}}
\sphinxAtStartPar
To facilitate users to operate on multiple {\hyperref[\detokenize{pyapiref:chappyapi-qconstraint}]{\sphinxcrossref{\DUrole{std,std-ref}{QConstraint Class}}}}  objects,
the Python interface of COPT provides QConstrArray class with the following methods:

\subsubsection{QConstrArray()}
\label{\detokenize{pyapiref:qconstrarray}}\begin{quote}

\sphinxAtStartPar
\sphinxstylestrong{Synopsis}
\begin{quote}

\sphinxAtStartPar
\sphinxcode{\sphinxupquote{QConstrArray(qconstrs=None)}}
\end{quote}

\sphinxAtStartPar
\sphinxstylestrong{Description}
\begin{quote}

\sphinxAtStartPar
Create a {\hyperref[\detokenize{pyapiref:chappyapi-qconstrarray}]{\sphinxcrossref{\DUrole{std,std-ref}{QConstrArray Class}}}} object.

\sphinxAtStartPar
If parameter \sphinxcode{\sphinxupquote{qconstrs}} is \sphinxcode{\sphinxupquote{None}}, the create an empty {\hyperref[\detokenize{pyapiref:chappyapi-qconstrarray}]{\sphinxcrossref{\DUrole{std,std-ref}{QConstrArray Class}}}} object,
otherwise initialize the newly created {\hyperref[\detokenize{pyapiref:chappyapi-qconstrarray}]{\sphinxcrossref{\DUrole{std,std-ref}{QConstrArray Class}}}} object with parameter \sphinxcode{\sphinxupquote{qconstrs}}
\end{quote}

\sphinxAtStartPar
\sphinxstylestrong{Arguments}
\begin{quote}

\sphinxAtStartPar
\sphinxcode{\sphinxupquote{qconstrs}}
\begin{quote}

\sphinxAtStartPar
Quadratic constraints to be added. \sphinxcode{\sphinxupquote{None}} by default.

\sphinxAtStartPar
\sphinxcode{\sphinxupquote{qconstrs}} can be {\hyperref[\detokenize{pyapiref:chappyapi-qconstraint}]{\sphinxcrossref{\DUrole{std,std-ref}{QConstraint Class}}}} object,
{\hyperref[\detokenize{pyapiref:chappyapi-qconstrarray}]{\sphinxcrossref{\DUrole{std,std-ref}{QConstrArray Class}}}} object, list, dictionary or {\hyperref[\detokenize{pyapiref:chappyapi-util-tupledict}]{\sphinxcrossref{\DUrole{std,std-ref}{tupledict Class}}}} object.
\end{quote}
\end{quote}

\sphinxAtStartPar
\sphinxstylestrong{Example}
\end{quote}

\begin{sphinxVerbatim}[commandchars=\\\{\}]
\PYG{c+c1}{\PYGZsh{} Create an empty QConstrArray object}
\PYG{n}{qconarr} \PYG{o}{=} \PYG{n}{QConstrArray}\PYG{p}{(}\PYG{p}{)}
\PYG{c+c1}{\PYGZsh{} Create an QConstrArray object initialized with quadratic constraint qconx and qcony}
\PYG{n}{qconarr} \PYG{o}{=} \PYG{n}{QConstrArray}\PYG{p}{(}\PYG{p}{[}\PYG{n}{qconx}\PYG{p}{,} \PYG{n}{qcony}\PYG{p}{]}\PYG{p}{)}
\end{sphinxVerbatim}

\subsubsection{QConstrArray.pushBack()}
\label{\detokenize{pyapiref:qconstrarray-pushback}}\begin{quote}

\sphinxAtStartPar
\sphinxstylestrong{Synopsis}
\begin{quote}

\sphinxAtStartPar
\sphinxcode{\sphinxupquote{pushBack(constr)}}
\end{quote}

\sphinxAtStartPar
\sphinxstylestrong{Description}
\begin{quote}

\sphinxAtStartPar
Add single or multiple {\hyperref[\detokenize{pyapiref:chappyapi-qconstraint}]{\sphinxcrossref{\DUrole{std,std-ref}{QConstraint Class}}}} object.
\end{quote}

\sphinxAtStartPar
\sphinxstylestrong{Arguments}
\begin{quote}

\sphinxAtStartPar
\sphinxcode{\sphinxupquote{constr}}
\begin{quote}

\sphinxAtStartPar
Quadratic constraints to be added. \sphinxcode{\sphinxupquote{None}} by default.

\sphinxAtStartPar
\sphinxcode{\sphinxupquote{qconstrs}} can be {\hyperref[\detokenize{pyapiref:chappyapi-qconstraint}]{\sphinxcrossref{\DUrole{std,std-ref}{QConstraint Class}}}} object,
{\hyperref[\detokenize{pyapiref:chappyapi-qconstrarray}]{\sphinxcrossref{\DUrole{std,std-ref}{QConstrArray Class}}}} object, list, dictionary or {\hyperref[\detokenize{pyapiref:chappyapi-util-tupledict}]{\sphinxcrossref{\DUrole{std,std-ref}{tupledict Class}}}} object.
\end{quote}
\end{quote}

\sphinxAtStartPar
\sphinxstylestrong{Example}
\end{quote}

\begin{sphinxVerbatim}[commandchars=\\\{\}]
\PYG{c+c1}{\PYGZsh{} Add quadratic constraint qr to conarr}
\PYG{n}{qconarr}\PYG{o}{.}\PYG{n}{pushBack}\PYG{p}{(}\PYG{n}{qr}\PYG{p}{)}
\PYG{c+c1}{\PYGZsh{} Add quadratic constraint qr0 and qr1 to qconarr}
\PYG{n}{qconarr}\PYG{o}{.}\PYG{n}{pushBack}\PYG{p}{(}\PYG{p}{[}\PYG{n}{qr0}\PYG{p}{,} \PYG{n}{qr1}\PYG{p}{]}\PYG{p}{)}
\end{sphinxVerbatim}

\subsubsection{QConstrArray.getQConstr()}
\label{\detokenize{pyapiref:qconstrarray-getqconstr}}\begin{quote}

\sphinxAtStartPar
\sphinxstylestrong{Synopsis}
\begin{quote}

\sphinxAtStartPar
\sphinxcode{\sphinxupquote{getQConstr(idx)}}
\end{quote}

\sphinxAtStartPar
\sphinxstylestrong{Description}
\begin{quote}

\sphinxAtStartPar
Retrieve the quadratic constraint according to its subscript in {\hyperref[\detokenize{pyapiref:chappyapi-qconstrarray}]{\sphinxcrossref{\DUrole{std,std-ref}{QConstrArray Class}}}} object.
Return a {\hyperref[\detokenize{pyapiref:chappyapi-qconstraint}]{\sphinxcrossref{\DUrole{std,std-ref}{QConstraint Class}}}} object.
\end{quote}

\sphinxAtStartPar
\sphinxstylestrong{Arguments}
\begin{quote}

\sphinxAtStartPar
\sphinxcode{\sphinxupquote{idx}}
\begin{quote}

\sphinxAtStartPar
Subscript of the desired quadratic constraint in {\hyperref[\detokenize{pyapiref:chappyapi-qconstrarray}]{\sphinxcrossref{\DUrole{std,std-ref}{QConstrArray Class}}}} object, starting with 0.
\end{quote}
\end{quote}

\sphinxAtStartPar
\sphinxstylestrong{Example}
\end{quote}

\begin{sphinxVerbatim}[commandchars=\\\{\}]
\PYG{c+c1}{\PYGZsh{} Retrieve the quadratic constraint with subscript 1  in qconarr}
\PYG{n}{qcon} \PYG{o}{=} \PYG{n}{qconarr}\PYG{o}{.}\PYG{n}{getQConstr}\PYG{p}{(}\PYG{l+m+mi}{1}\PYG{p}{)}
\end{sphinxVerbatim}

\subsubsection{QConstrArray.getSize()}
\label{\detokenize{pyapiref:qconstrarray-getsize}}\begin{quote}

\sphinxAtStartPar
\sphinxstylestrong{Synopsis}
\begin{quote}

\sphinxAtStartPar
\sphinxcode{\sphinxupquote{getSize()}}
\end{quote}

\sphinxAtStartPar
\sphinxstylestrong{Description}
\begin{quote}

\sphinxAtStartPar
Get the number of elements in {\hyperref[\detokenize{pyapiref:chappyapi-qconstrarray}]{\sphinxcrossref{\DUrole{std,std-ref}{QConstrArray Class}}}} object.
\end{quote}

\sphinxAtStartPar
\sphinxstylestrong{Example}
\end{quote}

\begin{sphinxVerbatim}[commandchars=\\\{\}]
\PYG{c+c1}{\PYGZsh{} Get the number of quadratic constraints in qconarr}
\PYG{n}{qarrsize} \PYG{o}{=} \PYG{n}{qconarr}\PYG{o}{.}\PYG{n}{getSize}\PYG{p}{(}\PYG{p}{)}
\end{sphinxVerbatim}

\subsection{QConstrBuilder Class}
\label{\detokenize{pyapiref:qconstrbuilder-class}}\label{\detokenize{pyapiref:chappyapi-qconstrbuilder}}
\sphinxAtStartPar
QConstrBuilder object contains operations related to temporary constraints when building quadratic constraints, and provides the following methods:

\subsubsection{QConstrBuilder()}
\label{\detokenize{pyapiref:qconstrbuilder}}\begin{quote}

\sphinxAtStartPar
\sphinxstylestrong{Synopsis}
\begin{quote}

\sphinxAtStartPar
\sphinxcode{\sphinxupquote{QConstrBuilder()}}
\end{quote}

\sphinxAtStartPar
\sphinxstylestrong{Description}
\begin{quote}

\sphinxAtStartPar
Create an empty {\hyperref[\detokenize{pyapiref:chappyapi-qconstrbuilder}]{\sphinxcrossref{\DUrole{std,std-ref}{QConstrBuilder Class}}}} object.
\end{quote}

\sphinxAtStartPar
\sphinxstylestrong{Example}
\end{quote}

\begin{sphinxVerbatim}[commandchars=\\\{\}]
\PYG{c+c1}{\PYGZsh{} Create an empty quadratic constraint builder}
\PYG{n}{qconstrbuilder} \PYG{o}{=} \PYG{n}{QConstrBuilder}\PYG{p}{(}\PYG{p}{)}
\end{sphinxVerbatim}

\subsubsection{QConstrBuilder.setBuilder()}
\label{\detokenize{pyapiref:qconstrbuilder-setbuilder}}\begin{quote}

\sphinxAtStartPar
\sphinxstylestrong{Synopsis}
\begin{quote}

\sphinxAtStartPar
\sphinxcode{\sphinxupquote{setBuilder(expr, sense, rhs)}}
\end{quote}

\sphinxAtStartPar
\sphinxstylestrong{Description}
\begin{quote}

\sphinxAtStartPar
Set expression, constraint type and RHS for quadratic constraint builder.
\end{quote}

\sphinxAtStartPar
\sphinxstylestrong{Arguments}
\begin{quote}

\sphinxAtStartPar
\sphinxcode{\sphinxupquote{expr}}
\begin{quote}

\sphinxAtStartPar
The expression to be set, which can be {\hyperref[\detokenize{pyapiref:chappyapi-var}]{\sphinxcrossref{\DUrole{std,std-ref}{Var Class}}}} object,
{\hyperref[\detokenize{pyapiref:chappyapi-linexpr}]{\sphinxcrossref{\DUrole{std,std-ref}{LinExpr Class}}}} object or {\hyperref[\detokenize{pyapiref:chappyapi-quadexpr}]{\sphinxcrossref{\DUrole{std,std-ref}{QuadExpr Class}}}} object.
\end{quote}

\sphinxAtStartPar
\sphinxcode{\sphinxupquote{sense}}
\begin{quote}

\sphinxAtStartPar
Sense of quadratic constraint. The full list of available tyeps can be found in {\hyperref[\detokenize{constant:chapconst-constrtype}]{\sphinxcrossref{\DUrole{std,std-ref}{Constraint type}}}} section.
\end{quote}

\sphinxAtStartPar
\sphinxcode{\sphinxupquote{rhs}}
\begin{quote}

\sphinxAtStartPar
Right hand side of quadratic constraint.
\end{quote}
\end{quote}

\sphinxAtStartPar
\sphinxstylestrong{Example}
\end{quote}

\begin{sphinxVerbatim}[commandchars=\\\{\}]
\PYG{c+c1}{\PYGZsh{} Set the expresson of quadratic constraint builder as: x+y, sense of constraint as equal and RHS as 1}
\PYG{n}{qconstrbuilder}\PYG{o}{.}\PYG{n}{setBuilder}\PYG{p}{(}\PYG{n}{x} \PYG{o}{+} \PYG{n}{y}\PYG{p}{,} \PYG{n}{COPT}\PYG{o}{.}\PYG{n}{LESS\PYGZus{}EQUAL}\PYG{p}{,} \PYG{l+m+mf}{1.0}\PYG{p}{)}
\end{sphinxVerbatim}

\subsubsection{QConstrBuilder.getQuadExpr()}
\label{\detokenize{pyapiref:qconstrbuilder-getquadexpr}}\begin{quote}

\sphinxAtStartPar
\sphinxstylestrong{Synopsis}
\begin{quote}

\sphinxAtStartPar
\sphinxcode{\sphinxupquote{getQuadExpr()}}
\end{quote}

\sphinxAtStartPar
\sphinxstylestrong{Description}
\begin{quote}

\sphinxAtStartPar
Retrieve the expression of a quadratic constraint builder object.
\end{quote}

\sphinxAtStartPar
\sphinxstylestrong{Example}
\end{quote}

\begin{sphinxVerbatim}[commandchars=\\\{\}]
\PYG{c+c1}{\PYGZsh{}  Retrieve the expression of a quadratic constraint builder}
\PYG{n}{quadexpr} \PYG{o}{=} \PYG{n}{constrbuilder}\PYG{o}{.}\PYG{n}{getQuadExpr}\PYG{p}{(}\PYG{p}{)}
\end{sphinxVerbatim}

\subsubsection{QConstrBuilder.getSense()}
\label{\detokenize{pyapiref:qconstrbuilder-getsense}}\begin{quote}

\sphinxAtStartPar
\sphinxstylestrong{Synopsis}
\begin{quote}

\sphinxAtStartPar
\sphinxcode{\sphinxupquote{getSense()}}
\end{quote}

\sphinxAtStartPar
\sphinxstylestrong{Description}
\begin{quote}

\sphinxAtStartPar
Retrieve the constraint sense of quadratic constraint builder object.
\end{quote}

\sphinxAtStartPar
\sphinxstylestrong{Example}
\end{quote}

\begin{sphinxVerbatim}[commandchars=\\\{\}]
\PYG{c+c1}{\PYGZsh{} Retrieve the constraint sense of quadratic constraint builder object.}
\PYG{n}{qconsense} \PYG{o}{=} \PYG{n}{qconstrbuilder}\PYG{o}{.}\PYG{n}{getSense}\PYG{p}{(}\PYG{p}{)}
\end{sphinxVerbatim}

\subsection{QConstrBuilderArray Class}
\label{\detokenize{pyapiref:qconstrbuilderarray-class}}\label{\detokenize{pyapiref:chappyapi-qconstrbuilderarray}}
\sphinxAtStartPar
To facilitate users to operate on multiple {\hyperref[\detokenize{pyapiref:chappyapi-qconstrbuilder}]{\sphinxcrossref{\DUrole{std,std-ref}{QConstrBuilder Class}}}} objects,
the Python interface of COPT provides QConstrArray object with the following methods:

\subsubsection{QConstrBuilderArray()}
\label{\detokenize{pyapiref:qconstrbuilderarray}}\begin{quote}

\sphinxAtStartPar
\sphinxstylestrong{Synopsis}
\begin{quote}

\sphinxAtStartPar
\sphinxcode{\sphinxupquote{QConstrBuilderArray(qconstrbuilders=None)}}
\end{quote}

\sphinxAtStartPar
\sphinxstylestrong{Description}
\begin{quote}

\sphinxAtStartPar
Create a {\hyperref[\detokenize{pyapiref:chappyapi-qconstrbuilderarray}]{\sphinxcrossref{\DUrole{std,std-ref}{QConstrBuilderArray Class}}}} object.

\sphinxAtStartPar
If parameter \sphinxcode{\sphinxupquote{qconstrbuilders}} is \sphinxcode{\sphinxupquote{None}}, then create an empty {\hyperref[\detokenize{pyapiref:chappyapi-qconstrbuilderarray}]{\sphinxcrossref{\DUrole{std,std-ref}{QConstrBuilderArray Class}}}} object,
otherwise initialize the newly created {\hyperref[\detokenize{pyapiref:chappyapi-qconstrbuilderarray}]{\sphinxcrossref{\DUrole{std,std-ref}{QConstrBuilderArray Class}}}} object by parameter \sphinxcode{\sphinxupquote{qconstrbuilders}}.
\end{quote}

\sphinxAtStartPar
\sphinxstylestrong{Arguments}
\begin{quote}

\sphinxAtStartPar
\sphinxcode{\sphinxupquote{qconstrbuilders}}
\begin{quote}

\sphinxAtStartPar
Quadratic constraint builder to be added. Optional, \sphinxcode{\sphinxupquote{None}} by default. It can be
{\hyperref[\detokenize{pyapiref:chappyapi-qconstrbuilder}]{\sphinxcrossref{\DUrole{std,std-ref}{QConstrBuilder Class}}}} object, {\hyperref[\detokenize{pyapiref:chappyapi-qconstrbuilderarray}]{\sphinxcrossref{\DUrole{std,std-ref}{QConstrBuilderArray Class}}}} object,
list, dictionary or {\hyperref[\detokenize{pyapiref:chappyapi-util-tupledict}]{\sphinxcrossref{\DUrole{std,std-ref}{tupledict Class}}}} object.
\end{quote}
\end{quote}

\sphinxAtStartPar
\sphinxstylestrong{Example}
\end{quote}

\begin{sphinxVerbatim}[commandchars=\\\{\}]
\PYG{c+c1}{\PYGZsh{} Create an empty QConstrBuilderArray object.}
\PYG{n}{qconbuilderarr} \PYG{o}{=} \PYG{n}{QConstrBuilderArray}\PYG{p}{(}\PYG{p}{)}
\PYG{c+c1}{\PYGZsh{} Create a QConstrBuilderArray object and initialize it with builders: qconbuilderx and qconbuildery}
\PYG{n}{qconbuilderarr} \PYG{o}{=} \PYG{n}{QConstrBuilderArray}\PYG{p}{(}\PYG{p}{[}\PYG{n}{qconbuilderx}\PYG{p}{,} \PYG{n}{qconbuildery}\PYG{p}{]}\PYG{p}{)}
\end{sphinxVerbatim}

\subsubsection{QConstrBuilderArray.pushBack()}
\label{\detokenize{pyapiref:qconstrbuilderarray-pushback}}\begin{quote}

\sphinxAtStartPar
\sphinxstylestrong{Synopsis}
\begin{quote}

\sphinxAtStartPar
\sphinxcode{\sphinxupquote{pushBack(qconstrbuilder)}}
\end{quote}

\sphinxAtStartPar
\sphinxstylestrong{Description}
\begin{quote}

\sphinxAtStartPar
Add single or multiple {\hyperref[\detokenize{pyapiref:chappyapi-qconstrbuilder}]{\sphinxcrossref{\DUrole{std,std-ref}{QConstrBuilder Class}}}} objects.
\end{quote}

\sphinxAtStartPar
\sphinxstylestrong{Arguments}
\begin{quote}

\sphinxAtStartPar
\sphinxcode{\sphinxupquote{qconstrbuilder}}
\begin{quote}

\sphinxAtStartPar
Builder of quadratic constraint to be added, which can be {\hyperref[\detokenize{pyapiref:chappyapi-qconstrbuilder}]{\sphinxcrossref{\DUrole{std,std-ref}{QConstrBuilder Class}}}} object,
{\hyperref[\detokenize{pyapiref:chappyapi-qconstrbuilderarray}]{\sphinxcrossref{\DUrole{std,std-ref}{QConstrBuilderArray Class}}}} object, list, dictionary or {\hyperref[\detokenize{pyapiref:chappyapi-util-tupledict}]{\sphinxcrossref{\DUrole{std,std-ref}{tupledict Class}}}} object.
\end{quote}
\end{quote}

\sphinxAtStartPar
\sphinxstylestrong{Example}
\end{quote}

\begin{sphinxVerbatim}[commandchars=\\\{\}]
\PYG{c+c1}{\PYGZsh{} Add quadratic constraint builder qconbuilderx to qconbuilderarr}
\PYG{n}{qconbuilderarr}\PYG{o}{.}\PYG{n}{pushBack}\PYG{p}{(}\PYG{n}{qconbuilderx}\PYG{p}{)}
\PYG{c+c1}{\PYGZsh{} Add quadratic constraint builders qconbuilderx and qconbuildery to qconbuilderarr}
\PYG{n}{qconbuilderarr}\PYG{o}{.}\PYG{n}{pushBack}\PYG{p}{(}\PYG{p}{[}\PYG{n}{qconbuilderx}\PYG{p}{,} \PYG{n}{qconbuildery}\PYG{p}{]}\PYG{p}{)}
\end{sphinxVerbatim}

\subsubsection{QConstrBuilderArray.getBuilder()}
\label{\detokenize{pyapiref:qconstrbuilderarray-getbuilder}}\begin{quote}

\sphinxAtStartPar
\sphinxstylestrong{Synopsis}
\begin{quote}

\sphinxAtStartPar
\sphinxcode{\sphinxupquote{getBuilder(idx)}}
\end{quote}

\sphinxAtStartPar
\sphinxstylestrong{Description}
\begin{quote}

\sphinxAtStartPar
Retrieve the corresponding builder object according to the subscript of quadratic constraint builder
in {\hyperref[\detokenize{pyapiref:chappyapi-qconstrbuilderarray}]{\sphinxcrossref{\DUrole{std,std-ref}{QConstrBuilderArray Class}}}} object.
\end{quote}

\sphinxAtStartPar
\sphinxstylestrong{Arguments}
\begin{quote}

\sphinxAtStartPar
\sphinxcode{\sphinxupquote{idx}}
\begin{quote}

\sphinxAtStartPar
Subscript of the quadratic constraint builder in the {\hyperref[\detokenize{pyapiref:chappyapi-qconstrbuilderarray}]{\sphinxcrossref{\DUrole{std,std-ref}{QConstrBuilderArray Class}}}} object, starting with 0.
\end{quote}
\end{quote}

\sphinxAtStartPar
\sphinxstylestrong{Example}
\end{quote}

\begin{sphinxVerbatim}[commandchars=\\\{\}]
\PYG{c+c1}{\PYGZsh{} Retrieve the builder with subscript 1 in qconbuilderarr}
\PYG{n}{qconbuilder} \PYG{o}{=} \PYG{n}{qconbuilderarr}\PYG{o}{.}\PYG{n}{getBuilder}\PYG{p}{(}\PYG{l+m+mi}{1}\PYG{p}{)}
\end{sphinxVerbatim}

\subsubsection{QConstrBuilderArray.getSize()}
\label{\detokenize{pyapiref:qconstrbuilderarray-getsize}}\begin{quote}

\sphinxAtStartPar
\sphinxstylestrong{Synopsis}
\begin{quote}

\sphinxAtStartPar
\sphinxcode{\sphinxupquote{getSize()}}
\end{quote}

\sphinxAtStartPar
\sphinxstylestrong{Description}
\begin{quote}

\sphinxAtStartPar
Get the number of elements in {\hyperref[\detokenize{pyapiref:chappyapi-qconstrbuilderarray}]{\sphinxcrossref{\DUrole{std,std-ref}{QConstrBuilderArray Class}}}} object.
\end{quote}

\sphinxAtStartPar
\sphinxstylestrong{Example}
\end{quote}

\begin{sphinxVerbatim}[commandchars=\\\{\}]
\PYG{c+c1}{\PYGZsh{} Get the number of builders in qconbuilderarr}
\PYG{n}{qarrsize} \PYG{o}{=} \PYG{n}{qconbuilderarr}\PYG{o}{.}\PYG{n}{getSize}\PYG{p}{(}\PYG{p}{)}
\end{sphinxVerbatim}

\subsection{PsdConstraint Class}
\label{\detokenize{pyapiref:psdconstraint-class}}\label{\detokenize{pyapiref:chappyapi-psdconstraint}}
\sphinxAtStartPar
PsdConstraint object contains related operations of COPT positive semi\sphinxhyphen{}definite constraints
and provides the following methods:

\subsubsection{PsdConstraint.getName()}
\label{\detokenize{pyapiref:psdconstraint-getname}}\begin{quote}

\sphinxAtStartPar
\sphinxstylestrong{Synopsis}
\begin{quote}

\sphinxAtStartPar
\sphinxcode{\sphinxupquote{getName()}}
\end{quote}

\sphinxAtStartPar
\sphinxstylestrong{Description}
\begin{quote}

\sphinxAtStartPar
Retrieve the name of positive semi\sphinxhyphen{}definite constraint.
\end{quote}

\sphinxAtStartPar
\sphinxstylestrong{Example}
\end{quote}

\begin{sphinxVerbatim}[commandchars=\\\{\}]
\PYG{c+c1}{\PYGZsh{} Retrieve the name of positive semi\PYGZhy{}definite constraint \PYGZsq{}con\PYGZsq{}.}
\PYG{n}{conname} \PYG{o}{=} \PYG{n}{con}\PYG{o}{.}\PYG{n}{getName}\PYG{p}{(}\PYG{p}{)}
\end{sphinxVerbatim}

\subsubsection{PsdConstraint.getIdx()}
\label{\detokenize{pyapiref:psdconstraint-getidx}}\begin{quote}

\sphinxAtStartPar
\sphinxstylestrong{Synopsis}
\begin{quote}

\sphinxAtStartPar
\sphinxcode{\sphinxupquote{getIdx()}}
\end{quote}

\sphinxAtStartPar
\sphinxstylestrong{Description}
\begin{quote}

\sphinxAtStartPar
Retrieve the subscript of positive semi\sphinxhyphen{}definite constraint in the model.
\end{quote}

\sphinxAtStartPar
\sphinxstylestrong{Example}
\end{quote}

\begin{sphinxVerbatim}[commandchars=\\\{\}]
\PYG{c+c1}{\PYGZsh{} Retrieve the subscript of positive semi\PYGZhy{}definite constraint con.}
\PYG{n}{conidx} \PYG{o}{=} \PYG{n}{con}\PYG{o}{.}\PYG{n}{getIdx}\PYG{p}{(}\PYG{p}{)}
\end{sphinxVerbatim}

\subsubsection{PsdConstraint.setName()}
\label{\detokenize{pyapiref:psdconstraint-setname}}\begin{quote}

\sphinxAtStartPar
\sphinxstylestrong{Synopsis}
\begin{quote}

\sphinxAtStartPar
\sphinxcode{\sphinxupquote{setName(newname)}}
\end{quote}

\sphinxAtStartPar
\sphinxstylestrong{Description}
\begin{quote}

\sphinxAtStartPar
Set the name of positive semi\sphinxhyphen{}definite constraint.
\end{quote}

\sphinxAtStartPar
\sphinxstylestrong{Arguments}
\begin{quote}

\sphinxAtStartPar
\sphinxcode{\sphinxupquote{newname}}
\begin{quote}

\sphinxAtStartPar
The name of positive semi\sphinxhyphen{}definite constraint to be set.
\end{quote}
\end{quote}

\sphinxAtStartPar
\sphinxstylestrong{Example}
\end{quote}

\begin{sphinxVerbatim}[commandchars=\\\{\}]
\PYG{c+c1}{\PYGZsh{} Set the name of positive semi\PYGZhy{}definite constraint \PYGZsq{}con\PYGZsq{}.}
\PYG{n}{con}\PYG{o}{.}\PYG{n}{setName}\PYG{p}{(}\PYG{l+s+s1}{\PYGZsq{}}\PYG{l+s+s1}{con}\PYG{l+s+s1}{\PYGZsq{}}\PYG{p}{)}
\end{sphinxVerbatim}

\subsubsection{PsdConstraint.getInfo()}
\label{\detokenize{pyapiref:psdconstraint-getinfo}}\begin{quote}

\sphinxAtStartPar
\sphinxstylestrong{Synopsis}
\begin{quote}

\sphinxAtStartPar
\sphinxcode{\sphinxupquote{getInfo(infoname)}}
\end{quote}

\sphinxAtStartPar
\sphinxstylestrong{Description}
\begin{quote}

\sphinxAtStartPar
Retrieve specified information. Return a constant.
\end{quote}

\sphinxAtStartPar
\sphinxstylestrong{Arguments}
\begin{quote}

\sphinxAtStartPar
\sphinxcode{\sphinxupquote{infoname}}
\begin{quote}

\sphinxAtStartPar
Name of the information to be obtained. Please refer to {\hyperref[\detokenize{pyapiref:chappyapi-const-info}]{\sphinxcrossref{\DUrole{std,std-ref}{Information}}}} section for possible values.
\end{quote}
\end{quote}

\sphinxAtStartPar
\sphinxstylestrong{Example}
\end{quote}

\begin{sphinxVerbatim}[commandchars=\\\{\}]
\PYG{c+c1}{\PYGZsh{} Get the lower bound of positive semi\PYGZhy{}definite constraint con}
\PYG{n}{conlb} \PYG{o}{=} \PYG{n}{con}\PYG{o}{.}\PYG{n}{getInfo}\PYG{p}{(}\PYG{n}{COPT}\PYG{o}{.}\PYG{n}{Info}\PYG{o}{.}\PYG{n}{LB}\PYG{p}{)}
\end{sphinxVerbatim}

\subsubsection{PsdConstraint.setInfo()}
\label{\detokenize{pyapiref:psdconstraint-setinfo}}\begin{quote}

\sphinxAtStartPar
\sphinxstylestrong{Synopsis}
\begin{quote}

\sphinxAtStartPar
\sphinxcode{\sphinxupquote{setInfo(infoname, newval)}}
\end{quote}

\sphinxAtStartPar
\sphinxstylestrong{Description}
\begin{quote}

\sphinxAtStartPar
Set new information value to the specified positive semi\sphinxhyphen{}definite constraint.
\end{quote}

\sphinxAtStartPar
\sphinxstylestrong{Arguments}
\begin{quote}

\sphinxAtStartPar
\sphinxcode{\sphinxupquote{infoname}}
\begin{quote}

\sphinxAtStartPar
The name of the information to be set. Please refer to {\hyperref[\detokenize{pyapiref:chappyapi-const-info}]{\sphinxcrossref{\DUrole{std,std-ref}{Information}}}} section for possible values.
\end{quote}

\sphinxAtStartPar
\sphinxcode{\sphinxupquote{newval}}
\begin{quote}

\sphinxAtStartPar
New information value to be set.
\end{quote}
\end{quote}

\sphinxAtStartPar
\sphinxstylestrong{Example}
\end{quote}

\begin{sphinxVerbatim}[commandchars=\\\{\}]
\PYG{c+c1}{\PYGZsh{} Set the lower bound of positive semi\PYGZhy{}definite constraint con}
\PYG{n}{con}\PYG{o}{.}\PYG{n}{setInfo}\PYG{p}{(}\PYG{n}{COPT}\PYG{o}{.}\PYG{n}{Info}\PYG{o}{.}\PYG{n}{LB}\PYG{p}{,} \PYG{l+m+mf}{1.0}\PYG{p}{)}
\end{sphinxVerbatim}

\subsubsection{PsdConstraint.remove()}
\label{\detokenize{pyapiref:psdconstraint-remove}}\begin{quote}

\sphinxAtStartPar
\sphinxstylestrong{Synopsis}
\begin{quote}

\sphinxAtStartPar
\sphinxcode{\sphinxupquote{remove()}}
\end{quote}

\sphinxAtStartPar
\sphinxstylestrong{Description}
\begin{quote}

\sphinxAtStartPar
Delete the positive semi\sphinxhyphen{}definite constraint from model.
\end{quote}

\sphinxAtStartPar
\sphinxstylestrong{Example}
\end{quote}

\begin{sphinxVerbatim}[commandchars=\\\{\}]
\PYG{c+c1}{\PYGZsh{} Delete the positive semi\PYGZhy{}definite constraint \PYGZsq{}conx\PYGZsq{}}
\PYG{n}{conx}\PYG{o}{.}\PYG{n}{remove}\PYG{p}{(}\PYG{p}{)}
\end{sphinxVerbatim}

\subsection{PsdConstrArray Class}
\label{\detokenize{pyapiref:psdconstrarray-class}}\label{\detokenize{pyapiref:chappyapi-psdconstrarray}}
\sphinxAtStartPar
To facilitate users to operate on multiple {\hyperref[\detokenize{pyapiref:chappyapi-psdconstraint}]{\sphinxcrossref{\DUrole{std,std-ref}{PsdConstraint Class}}}}  objects,
the Python interface of COPT provides PsdConstrArray class with the following methods:

\subsubsection{PsdConstrArray()}
\label{\detokenize{pyapiref:psdconstrarray}}\begin{quote}

\sphinxAtStartPar
\sphinxstylestrong{Synopsis}
\begin{quote}

\sphinxAtStartPar
\sphinxcode{\sphinxupquote{PsdConstrArray(constrs=None)}}
\end{quote}

\sphinxAtStartPar
\sphinxstylestrong{Description}
\begin{quote}

\sphinxAtStartPar
Create a {\hyperref[\detokenize{pyapiref:chappyapi-psdconstrarray}]{\sphinxcrossref{\DUrole{std,std-ref}{PsdConstrArray Class}}}} object.

\sphinxAtStartPar
If parameter \sphinxcode{\sphinxupquote{constrs}} is \sphinxcode{\sphinxupquote{None}}, the create an empty {\hyperref[\detokenize{pyapiref:chappyapi-psdconstrarray}]{\sphinxcrossref{\DUrole{std,std-ref}{PsdConstrArray Class}}}} object,
otherwise initialize the newly created {\hyperref[\detokenize{pyapiref:chappyapi-psdconstrarray}]{\sphinxcrossref{\DUrole{std,std-ref}{PsdConstrArray Class}}}} object with parameter \sphinxcode{\sphinxupquote{constrs}} .
\end{quote}

\sphinxAtStartPar
\sphinxstylestrong{Arguments}
\begin{quote}

\sphinxAtStartPar
\sphinxcode{\sphinxupquote{constrs}}
\begin{quote}

\sphinxAtStartPar
Positive semi\sphinxhyphen{}definite constraints to be added. \sphinxcode{\sphinxupquote{None}} by default.

\sphinxAtStartPar
\sphinxcode{\sphinxupquote{constrs}} can be {\hyperref[\detokenize{pyapiref:chappyapi-psdconstraint}]{\sphinxcrossref{\DUrole{std,std-ref}{PsdConstraint Class}}}} object,
{\hyperref[\detokenize{pyapiref:chappyapi-psdconstrarray}]{\sphinxcrossref{\DUrole{std,std-ref}{PsdConstrArray Class}}}} object, list, dictionary or {\hyperref[\detokenize{pyapiref:chappyapi-util-tupledict}]{\sphinxcrossref{\DUrole{std,std-ref}{tupledict Class}}}} object.
\end{quote}
\end{quote}

\sphinxAtStartPar
\sphinxstylestrong{Example}
\end{quote}

\begin{sphinxVerbatim}[commandchars=\\\{\}]
\PYG{c+c1}{\PYGZsh{} Create an empty PsdConstrArray object}
\PYG{n}{conarr} \PYG{o}{=} \PYG{n}{PsdConstrArray}\PYG{p}{(}\PYG{p}{)}
\PYG{c+c1}{\PYGZsh{} Create an PsdConstrArray object containing positive semi\PYGZhy{}definite constraint conx and cony}
\PYG{n}{conarr} \PYG{o}{=} \PYG{n}{PsdConstrArray}\PYG{p}{(}\PYG{p}{[}\PYG{n}{conx}\PYG{p}{,} \PYG{n}{cony}\PYG{p}{]}\PYG{p}{)}
\end{sphinxVerbatim}

\subsubsection{PsdConstrArray.pushBack()}
\label{\detokenize{pyapiref:psdconstrarray-pushback}}\begin{quote}

\sphinxAtStartPar
\sphinxstylestrong{Synopsis}
\begin{quote}

\sphinxAtStartPar
\sphinxcode{\sphinxupquote{pushBack(constr)}}
\end{quote}

\sphinxAtStartPar
\sphinxstylestrong{Description}
\begin{quote}

\sphinxAtStartPar
Add single or multiple {\hyperref[\detokenize{pyapiref:chappyapi-psdconstraint}]{\sphinxcrossref{\DUrole{std,std-ref}{PsdConstraint Class}}}} objects.
\end{quote}

\sphinxAtStartPar
\sphinxstylestrong{Arguments}
\begin{quote}

\sphinxAtStartPar
\sphinxcode{\sphinxupquote{constr}}
\begin{quote}

\sphinxAtStartPar
Positive semi\sphinxhyphen{}definite constraints to be applied.
\sphinxcode{\sphinxupquote{constrs}} can be {\hyperref[\detokenize{pyapiref:chappyapi-psdconstraint}]{\sphinxcrossref{\DUrole{std,std-ref}{PsdConstraint Class}}}} object, {\hyperref[\detokenize{pyapiref:chappyapi-psdconstrarray}]{\sphinxcrossref{\DUrole{std,std-ref}{PsdConstrArray Class}}}} object,
list, dictionary or {\hyperref[\detokenize{pyapiref:chappyapi-util-tupledict}]{\sphinxcrossref{\DUrole{std,std-ref}{tupledict Class}}}} object.
\end{quote}
\end{quote}

\sphinxAtStartPar
\sphinxstylestrong{Example}
\end{quote}

\begin{sphinxVerbatim}[commandchars=\\\{\}]
\PYG{c+c1}{\PYGZsh{} Add positive semi\PYGZhy{}definite constraint r to conarr}
\PYG{n}{conarr}\PYG{o}{.}\PYG{n}{pushBack}\PYG{p}{(}\PYG{n}{r}\PYG{p}{)}
\PYG{c+c1}{\PYGZsh{} Add positive semi\PYGZhy{}definite constraint r0 and r1 to conarr}
\PYG{n}{conarr}\PYG{o}{.}\PYG{n}{pushBack}\PYG{p}{(}\PYG{p}{[}\PYG{n}{r0}\PYG{p}{,} \PYG{n}{r1}\PYG{p}{]}\PYG{p}{)}
\end{sphinxVerbatim}

\subsubsection{PsdConstrArray.getPsdConstr()}
\label{\detokenize{pyapiref:psdconstrarray-getpsdconstr}}\begin{quote}

\sphinxAtStartPar
\sphinxstylestrong{Synopsis}
\begin{quote}

\sphinxAtStartPar
\sphinxcode{\sphinxupquote{getPsdConstr(idx)}}
\end{quote}

\sphinxAtStartPar
\sphinxstylestrong{Description}
\begin{quote}

\sphinxAtStartPar
Retrieve the positive semi\sphinxhyphen{}definite constraint according to its subscript in {\hyperref[\detokenize{pyapiref:chappyapi-psdconstrarray}]{\sphinxcrossref{\DUrole{std,std-ref}{PsdConstrArray Class}}}} object.
Return a {\hyperref[\detokenize{pyapiref:chappyapi-psdconstraint}]{\sphinxcrossref{\DUrole{std,std-ref}{PsdConstraint Class}}}} object.
\end{quote}

\sphinxAtStartPar
\sphinxstylestrong{Arguments}
\begin{quote}

\sphinxAtStartPar
\sphinxcode{\sphinxupquote{idx}}
\begin{quote}

\sphinxAtStartPar
Subscript of the desired positive semi\sphinxhyphen{}definite constraint in {\hyperref[\detokenize{pyapiref:chappyapi-psdconstrarray}]{\sphinxcrossref{\DUrole{std,std-ref}{PsdConstrArray Class}}}} object, starting with 0.
\end{quote}
\end{quote}

\sphinxAtStartPar
\sphinxstylestrong{Example}
\end{quote}

\begin{sphinxVerbatim}[commandchars=\\\{\}]
\PYG{c+c1}{\PYGZsh{} Retrieve the positive semi\PYGZhy{}definite constraint with subscript 1  in conarr}
\PYG{n}{con} \PYG{o}{=} \PYG{n}{conarr}\PYG{o}{.}\PYG{n}{getPsdConstr}\PYG{p}{(}\PYG{l+m+mi}{1}\PYG{p}{)}
\end{sphinxVerbatim}

\subsubsection{PsdConstrArray.getSize()}
\label{\detokenize{pyapiref:psdconstrarray-getsize}}\begin{quote}

\sphinxAtStartPar
\sphinxstylestrong{Synopsis}
\begin{quote}

\sphinxAtStartPar
\sphinxcode{\sphinxupquote{getSize()}}
\end{quote}

\sphinxAtStartPar
\sphinxstylestrong{Description}
\begin{quote}

\sphinxAtStartPar
Get the number of elements in {\hyperref[\detokenize{pyapiref:chappyapi-psdconstrarray}]{\sphinxcrossref{\DUrole{std,std-ref}{PsdConstrArray Class}}}} object.
\end{quote}

\sphinxAtStartPar
\sphinxstylestrong{Example}
\end{quote}

\begin{sphinxVerbatim}[commandchars=\\\{\}]
\PYG{c+c1}{\PYGZsh{} Get the number of positive semi\PYGZhy{}definite constraints in conarr}
\PYG{n}{arrsize} \PYG{o}{=} \PYG{n}{conarr}\PYG{o}{.}\PYG{n}{getSize}\PYG{p}{(}\PYG{p}{)}
\end{sphinxVerbatim}

\subsection{PsdConstrBuilder Class}
\label{\detokenize{pyapiref:psdconstrbuilder-class}}\label{\detokenize{pyapiref:chappyapi-psdconstrbuilder}}
\sphinxAtStartPar
PsdConstrBuilder object contains operations related to temporary constraints when building
positive semi\sphinxhyphen{}definite constraints, and provides the following methods:

\subsubsection{PsdConstrBuilder()}
\label{\detokenize{pyapiref:psdconstrbuilder}}\begin{quote}

\sphinxAtStartPar
\sphinxstylestrong{Synopsis}
\begin{quote}

\sphinxAtStartPar
\sphinxcode{\sphinxupquote{PsdConstrBuilder()}}
\end{quote}

\sphinxAtStartPar
\sphinxstylestrong{Description}
\begin{quote}

\sphinxAtStartPar
Create an empty {\hyperref[\detokenize{pyapiref:chappyapi-psdconstrbuilder}]{\sphinxcrossref{\DUrole{std,std-ref}{PsdConstrBuilder Class}}}} object
\end{quote}

\sphinxAtStartPar
\sphinxstylestrong{Example}
\end{quote}

\begin{sphinxVerbatim}[commandchars=\\\{\}]
\PYG{c+c1}{\PYGZsh{} Create an empty positive semi\PYGZhy{}definite constraint builder}
\PYG{n}{constrbuilder} \PYG{o}{=} \PYG{n}{PsdConstrBuilder}\PYG{p}{(}\PYG{p}{)}
\end{sphinxVerbatim}

\subsubsection{PsdConstrBuilder.setBuilder()}
\label{\detokenize{pyapiref:psdconstrbuilder-setbuilder}}\begin{quote}

\sphinxAtStartPar
\sphinxstylestrong{Synopsis}
\begin{quote}

\sphinxAtStartPar
\sphinxcode{\sphinxupquote{setBuilder(expr, sense, rhs)}}
\end{quote}

\sphinxAtStartPar
\sphinxstylestrong{Description}
\begin{quote}

\sphinxAtStartPar
Set expression, constraint type and right hand side for positive semi\sphinxhyphen{}definite constraint builder.
\end{quote}

\sphinxAtStartPar
\sphinxstylestrong{Arguments}
\begin{quote}

\sphinxAtStartPar
\sphinxcode{\sphinxupquote{expr}}
\begin{quote}

\sphinxAtStartPar
The expression to be set, which can be {\hyperref[\detokenize{pyapiref:chappyapi-psdvar}]{\sphinxcrossref{\DUrole{std,std-ref}{PsdVar Class}}}} expression or
{\hyperref[\detokenize{pyapiref:chappyapi-psdexpr}]{\sphinxcrossref{\DUrole{std,std-ref}{PsdExpr Class}}}} expression.
\end{quote}

\sphinxAtStartPar
\sphinxcode{\sphinxupquote{sense}}
\begin{quote}

\sphinxAtStartPar
Sense of constraint. The full list of available tyeps can be found in {\hyperref[\detokenize{constant:chapconst-constrtype}]{\sphinxcrossref{\DUrole{std,std-ref}{Constraint type}}}} section.
\end{quote}

\sphinxAtStartPar
\sphinxcode{\sphinxupquote{rhs}}
\begin{quote}

\sphinxAtStartPar
The right hand side of constraint.
\end{quote}
\end{quote}

\sphinxAtStartPar
\sphinxstylestrong{Example}
\end{quote}

\begin{sphinxVerbatim}[commandchars=\\\{\}]
\PYG{c+c1}{\PYGZsh{} Set the expresson of positive semi\PYGZhy{}definite constraint builder as: x + y == 1, and sense of constraint as equal}
\PYG{n}{constrbuilder}\PYG{o}{.}\PYG{n}{setBuilder}\PYG{p}{(}\PYG{n}{x} \PYG{o}{+} \PYG{n}{y}\PYG{p}{,} \PYG{n}{COPT}\PYG{o}{.}\PYG{n}{EQUAL}\PYG{p}{,} \PYG{l+m+mi}{1}\PYG{p}{)}
\end{sphinxVerbatim}

\subsubsection{PsdConstrBuilder.setRange()}
\label{\detokenize{pyapiref:psdconstrbuilder-setrange}}\begin{quote}

\sphinxAtStartPar
\sphinxstylestrong{Synopsis}
\begin{quote}

\sphinxAtStartPar
\sphinxcode{\sphinxupquote{setRange(expr, range)}}
\end{quote}

\sphinxAtStartPar
\sphinxstylestrong{Description}
\begin{quote}

\sphinxAtStartPar
Set a range positive semi\sphinxhyphen{}definite constraint builder where \sphinxcode{\sphinxupquote{expr}} is less than or equals to 0
and greater than or equals to \sphinxhyphen{} \sphinxcode{\sphinxupquote{range}}.
\end{quote}

\sphinxAtStartPar
\sphinxstylestrong{Arguments}
\begin{quote}

\sphinxAtStartPar
\sphinxcode{\sphinxupquote{expr}}
\begin{quote}

\sphinxAtStartPar
The expression to be set, which can be {\hyperref[\detokenize{pyapiref:chappyapi-psdvar}]{\sphinxcrossref{\DUrole{std,std-ref}{PsdVar Class}}}} expression or
{\hyperref[\detokenize{pyapiref:chappyapi-psdexpr}]{\sphinxcrossref{\DUrole{std,std-ref}{PsdExpr Class}}}} expression.
\end{quote}

\sphinxAtStartPar
\sphinxcode{\sphinxupquote{range}}
\begin{quote}

\sphinxAtStartPar
Range of constraint, nonnegative constant.
\end{quote}
\end{quote}

\sphinxAtStartPar
\sphinxstylestrong{Example}
\end{quote}

\begin{sphinxVerbatim}[commandchars=\\\{\}]
\PYG{c+c1}{\PYGZsh{} Set a range positive semi\PYGZhy{}definite constraint builder: \PYGZhy{}1 \PYGZlt{}= x + y \PYGZhy{} 1 \PYGZlt{}= 0}
\PYG{n}{constrbuilder}\PYG{o}{.}\PYG{n}{setRange}\PYG{p}{(}\PYG{n}{x} \PYG{o}{+} \PYG{n}{y} \PYG{o}{\PYGZhy{}} \PYG{l+m+mi}{1}\PYG{p}{,} \PYG{l+m+mi}{1}\PYG{p}{)}
\end{sphinxVerbatim}

\subsubsection{PsdConstrBuilder.getPsdExpr()}
\label{\detokenize{pyapiref:psdconstrbuilder-getpsdexpr}}\begin{quote}

\sphinxAtStartPar
\sphinxstylestrong{Synopsis}
\begin{quote}

\sphinxAtStartPar
\sphinxcode{\sphinxupquote{getPsdExpr()}}
\end{quote}

\sphinxAtStartPar
\sphinxstylestrong{Description}
\begin{quote}

\sphinxAtStartPar
Retrieve the expression of a positive semi\sphinxhyphen{}definite constraint builder object.
\end{quote}

\sphinxAtStartPar
\sphinxstylestrong{Example}
\end{quote}

\begin{sphinxVerbatim}[commandchars=\\\{\}]
\PYG{c+c1}{\PYGZsh{}  Retrieve the expression of a positive semi\PYGZhy{}definite constraint builder}
\PYG{n}{psdexpr} \PYG{o}{=} \PYG{n}{constrbuilder}\PYG{o}{.}\PYG{n}{getPsdExpr}\PYG{p}{(}\PYG{p}{)}
\end{sphinxVerbatim}

\subsubsection{PsdConstrBuilder.getSense()}
\label{\detokenize{pyapiref:psdconstrbuilder-getsense}}\begin{quote}

\sphinxAtStartPar
\sphinxstylestrong{Synopsis}
\begin{quote}

\sphinxAtStartPar
\sphinxcode{\sphinxupquote{getSense()}}
\end{quote}

\sphinxAtStartPar
\sphinxstylestrong{Description}
\begin{quote}

\sphinxAtStartPar
Retrieve the constraint sense of positive semi\sphinxhyphen{}definite constraint builder object.
\end{quote}

\sphinxAtStartPar
\sphinxstylestrong{Example}
\end{quote}

\begin{sphinxVerbatim}[commandchars=\\\{\}]
\PYG{c+c1}{\PYGZsh{} Retrieve the constraint sense of positive semi\PYGZhy{}definite constraint builder object.}
\PYG{n}{consense} \PYG{o}{=} \PYG{n}{constrbuilder}\PYG{o}{.}\PYG{n}{getSense}\PYG{p}{(}\PYG{p}{)}
\end{sphinxVerbatim}

\subsubsection{PsdConstrBuilder.getRange()}
\label{\detokenize{pyapiref:psdconstrbuilder-getrange}}\begin{quote}

\sphinxAtStartPar
\sphinxstylestrong{Synopsis}
\begin{quote}

\sphinxAtStartPar
\sphinxcode{\sphinxupquote{getRange()}}
\end{quote}

\sphinxAtStartPar
\sphinxstylestrong{Description}
\begin{quote}

\sphinxAtStartPar
Retrieve the range of positive semi\sphinxhyphen{}definite constraint builder object,
i.e. length from lower bound to upper bound of the constraint
\end{quote}

\sphinxAtStartPar
\sphinxstylestrong{Example}
\end{quote}

\begin{sphinxVerbatim}[commandchars=\\\{\}]
\PYG{c+c1}{\PYGZsh{} Retrieve the range of positive semi\PYGZhy{}definite constraint builder object}
\PYG{n}{rngval} \PYG{o}{=} \PYG{n}{constrbuilder}\PYG{o}{.}\PYG{n}{getRange}\PYG{p}{(}\PYG{p}{)}
\end{sphinxVerbatim}

\subsection{PsdConstrBuilderArray Class}
\label{\detokenize{pyapiref:psdconstrbuilderarray-class}}\label{\detokenize{pyapiref:chappyapi-psdconstrbuilderarray}}
\sphinxAtStartPar
To facilitate users to operate on multiple {\hyperref[\detokenize{pyapiref:chappyapi-psdconstrbuilder}]{\sphinxcrossref{\DUrole{std,std-ref}{PsdConstrBuilder Class}}}} objects,
the Python interface of COPT provides PsdConstrBuilderArray object with the following methods:

\subsubsection{PsdConstrBuilderArray()}
\label{\detokenize{pyapiref:psdconstrbuilderarray}}\begin{quote}

\sphinxAtStartPar
\sphinxstylestrong{Synopsis}
\begin{quote}

\sphinxAtStartPar
\sphinxcode{\sphinxupquote{PsdConstrBuilderArray(builders=None)}}
\end{quote}

\sphinxAtStartPar
\sphinxstylestrong{Description}
\begin{quote}

\sphinxAtStartPar
Create a {\hyperref[\detokenize{pyapiref:chappyapi-psdconstrbuilderarray}]{\sphinxcrossref{\DUrole{std,std-ref}{PsdConstrBuilderArray Class}}}} object.

\sphinxAtStartPar
If parameter \sphinxcode{\sphinxupquote{builders}} is \sphinxcode{\sphinxupquote{None}}, then create an empty {\hyperref[\detokenize{pyapiref:chappyapi-psdconstrbuilderarray}]{\sphinxcrossref{\DUrole{std,std-ref}{PsdConstrBuilderArray Class}}}} object,
otherwise initialize the newly created {\hyperref[\detokenize{pyapiref:chappyapi-psdconstrbuilderarray}]{\sphinxcrossref{\DUrole{std,std-ref}{PsdConstrBuilderArray Class}}}} object by parameter \sphinxcode{\sphinxupquote{builders}}.
\end{quote}

\sphinxAtStartPar
\sphinxstylestrong{Arguments}
\begin{quote}

\sphinxAtStartPar
\sphinxcode{\sphinxupquote{builders}}
\begin{quote}

\sphinxAtStartPar
Positive semi\sphinxhyphen{}definite constraint builder to be added. Optional, \sphinxcode{\sphinxupquote{None}} by default.
It can be {\hyperref[\detokenize{pyapiref:chappyapi-psdconstrbuilder}]{\sphinxcrossref{\DUrole{std,std-ref}{PsdConstrBuilder Class}}}} object, {\hyperref[\detokenize{pyapiref:chappyapi-psdconstrbuilderarray}]{\sphinxcrossref{\DUrole{std,std-ref}{PsdConstrBuilderArray Class}}}} object,
list, dictionary or {\hyperref[\detokenize{pyapiref:chappyapi-util-tupledict}]{\sphinxcrossref{\DUrole{std,std-ref}{tupledict Class}}}} object.
\end{quote}
\end{quote}

\sphinxAtStartPar
\sphinxstylestrong{Example}
\end{quote}

\begin{sphinxVerbatim}[commandchars=\\\{\}]
\PYG{c+c1}{\PYGZsh{} Create an empty PsdConstrBuilderArray object.}
\PYG{n}{conbuilderarr} \PYG{o}{=} \PYG{n}{PsdConstrBuilderArray}\PYG{p}{(}\PYG{p}{)}
\PYG{c+c1}{\PYGZsh{} Create a PsdConstrBuilderArray object containing builders: conbuilderx and conbuildery}
\PYG{n}{conbuilderarr} \PYG{o}{=} \PYG{n}{PsdConstrBuilderArray}\PYG{p}{(}\PYG{p}{[}\PYG{n}{conbuilderx}\PYG{p}{,} \PYG{n}{conbuildery}\PYG{p}{]}\PYG{p}{)}
\end{sphinxVerbatim}

\subsubsection{PsdConstrBuilderArray.pushBack()}
\label{\detokenize{pyapiref:psdconstrbuilderarray-pushback}}\begin{quote}

\sphinxAtStartPar
\sphinxstylestrong{Synopsis}
\begin{quote}

\sphinxAtStartPar
\sphinxcode{\sphinxupquote{pushBack(builder)}}
\end{quote}

\sphinxAtStartPar
\sphinxstylestrong{Description}
\begin{quote}

\sphinxAtStartPar
Add single or multiple {\hyperref[\detokenize{pyapiref:chappyapi-psdconstrbuilder}]{\sphinxcrossref{\DUrole{std,std-ref}{PsdConstrBuilder Class}}}} objects.
\end{quote}

\sphinxAtStartPar
\sphinxstylestrong{Arguments}
\begin{quote}

\sphinxAtStartPar
\sphinxcode{\sphinxupquote{builder}}
\begin{quote}

\sphinxAtStartPar
Builder of positive semi\sphinxhyphen{}definite constraint to be added, which can be {\hyperref[\detokenize{pyapiref:chappyapi-psdconstrbuilder}]{\sphinxcrossref{\DUrole{std,std-ref}{PsdConstrBuilder Class}}}} object,
{\hyperref[\detokenize{pyapiref:chappyapi-psdconstrbuilderarray}]{\sphinxcrossref{\DUrole{std,std-ref}{PsdConstrBuilderArray Class}}}} object, list, dictionary or {\hyperref[\detokenize{pyapiref:chappyapi-util-tupledict}]{\sphinxcrossref{\DUrole{std,std-ref}{tupledict Class}}}} object.
\end{quote}
\end{quote}

\sphinxAtStartPar
\sphinxstylestrong{Example}
\end{quote}

\begin{sphinxVerbatim}[commandchars=\\\{\}]
\PYG{c+c1}{\PYGZsh{} Add positive semi\PYGZhy{}definite constraint builder conbuilderx to conbuilderarr}
\PYG{n}{conbuilderarr}\PYG{o}{.}\PYG{n}{pushBack}\PYG{p}{(}\PYG{n}{conbuilderx}\PYG{p}{)}
\PYG{c+c1}{\PYGZsh{} Add positive semi\PYGZhy{}definite constraint builders conbuilderx and conbuildery to conbuilderarr}
\PYG{n}{conbuilderarr}\PYG{o}{.}\PYG{n}{pushBack}\PYG{p}{(}\PYG{p}{[}\PYG{n}{conbuilderx}\PYG{p}{,} \PYG{n}{conbuildery}\PYG{p}{]}\PYG{p}{)}
\end{sphinxVerbatim}

\subsubsection{PsdConstrBuilderArray.getBuilder()}
\label{\detokenize{pyapiref:psdconstrbuilderarray-getbuilder}}\begin{quote}

\sphinxAtStartPar
\sphinxstylestrong{Synopsis}
\begin{quote}

\sphinxAtStartPar
\sphinxcode{\sphinxupquote{getBuilder(idx)}}
\end{quote}

\sphinxAtStartPar
\sphinxstylestrong{Description}
\begin{quote}

\sphinxAtStartPar
Retrieve the corresponding builder object according to the subscript of positive semi\sphinxhyphen{}definite constraint builder
in {\hyperref[\detokenize{pyapiref:chappyapi-psdconstrbuilderarray}]{\sphinxcrossref{\DUrole{std,std-ref}{PsdConstrBuilderArray Class}}}} object.
\end{quote}

\sphinxAtStartPar
\sphinxstylestrong{Arguments}
\begin{quote}

\sphinxAtStartPar
\sphinxcode{\sphinxupquote{idx}}
\begin{quote}

\sphinxAtStartPar
Subscript of the positive semi\sphinxhyphen{}definite constraint builder in the {\hyperref[\detokenize{pyapiref:chappyapi-psdconstrbuilderarray}]{\sphinxcrossref{\DUrole{std,std-ref}{PsdConstrBuilderArray Class}}}} object, starting with 0.
\end{quote}
\end{quote}

\sphinxAtStartPar
\sphinxstylestrong{Example}
\end{quote}

\begin{sphinxVerbatim}[commandchars=\\\{\}]
\PYG{c+c1}{\PYGZsh{} Retrieve the builder with subscript 1 in conbuilderarr}
\PYG{n}{conbuilder} \PYG{o}{=} \PYG{n}{conbuilderarr}\PYG{o}{.}\PYG{n}{getBuilder}\PYG{p}{(}\PYG{l+m+mi}{1}\PYG{p}{)}
\end{sphinxVerbatim}

\subsubsection{PsdConstrBuilderArray.getSize()}
\label{\detokenize{pyapiref:psdconstrbuilderarray-getsize}}\begin{quote}

\sphinxAtStartPar
\sphinxstylestrong{Synopsis}
\begin{quote}

\sphinxAtStartPar
\sphinxcode{\sphinxupquote{getSize()}}
\end{quote}

\sphinxAtStartPar
\sphinxstylestrong{Description}
\begin{quote}

\sphinxAtStartPar
Get the number of elements in {\hyperref[\detokenize{pyapiref:chappyapi-psdconstrbuilderarray}]{\sphinxcrossref{\DUrole{std,std-ref}{PsdConstrBuilderArray Class}}}} object.
\end{quote}

\sphinxAtStartPar
\sphinxstylestrong{Example}
\end{quote}

\begin{sphinxVerbatim}[commandchars=\\\{\}]
\PYG{c+c1}{\PYGZsh{} Get the number of builders in conbuilderarr}
\PYG{n}{arrsize} \PYG{o}{=} \PYG{n}{conbuilderarr}\PYG{o}{.}\PYG{n}{getSize}\PYG{p}{(}\PYG{p}{)}
\end{sphinxVerbatim}

\subsection{LmiConstraint Class}
\label{\detokenize{pyapiref:lmiconstraint-class}}\label{\detokenize{pyapiref:chappyapi-lmiconstraint}}
\sphinxAtStartPar
LmiConstraint object contains related operations of COPT LMI (Linear Matrix Inequality) constraints
and provides the following methods:

\subsubsection{LmiConstraint.getName()}
\label{\detokenize{pyapiref:lmiconstraint-getname}}\begin{quote}

\sphinxAtStartPar
\sphinxstylestrong{Synopsis}
\begin{quote}

\sphinxAtStartPar
\sphinxcode{\sphinxupquote{getName()}}
\end{quote}

\sphinxAtStartPar
\sphinxstylestrong{Description}
\begin{quote}

\sphinxAtStartPar
Retrieve the name of the {\hyperref[\detokenize{pyapiref:chappyapi-lmiconstraint}]{\sphinxcrossref{\DUrole{std,std-ref}{LmiConstraint Class}}}} object.
\end{quote}

\sphinxAtStartPar
\sphinxstylestrong{Example}
\end{quote}

\begin{sphinxVerbatim}[commandchars=\\\{\}]
\PYG{c+c1}{\PYGZsh{} Get name of the LmiConstraint con}
\PYG{n}{conname} \PYG{o}{=} \PYG{n}{con}\PYG{o}{.}\PYG{n}{getName}\PYG{p}{(}\PYG{p}{)}
\end{sphinxVerbatim}

\subsubsection{LmiConstraint.getIdx()}
\label{\detokenize{pyapiref:lmiconstraint-getidx}}\begin{quote}

\sphinxAtStartPar
\sphinxstylestrong{Synopsis}
\begin{quote}

\sphinxAtStartPar
\sphinxcode{\sphinxupquote{getIdx()}}
\end{quote}

\sphinxAtStartPar
\sphinxstylestrong{Description}
\begin{quote}

\sphinxAtStartPar
Retrieve the subscript of the LMI constraint in the model.
\end{quote}

\sphinxAtStartPar
\sphinxstylestrong{Example}
\end{quote}

\begin{sphinxVerbatim}[commandchars=\\\{\}]
\PYG{c+c1}{\PYGZsh{} Retrieve the subscript of the LMI constraint}
\PYG{n}{conidx} \PYG{o}{=} \PYG{n}{con}\PYG{o}{.}\PYG{n}{getIdx}\PYG{p}{(}\PYG{p}{)}
\end{sphinxVerbatim}

\subsubsection{LmiConstraint.getDim()}
\label{\detokenize{pyapiref:lmiconstraint-getdim}}\begin{quote}

\sphinxAtStartPar
\sphinxstylestrong{Synopsis}
\begin{quote}

\sphinxAtStartPar
\sphinxcode{\sphinxupquote{getDim()}}
\end{quote}

\sphinxAtStartPar
\sphinxstylestrong{Description}
\begin{quote}

\sphinxAtStartPar
Retrieve the dimension of the LMI constraint.
\end{quote}

\sphinxAtStartPar
\sphinxstylestrong{Example}
\end{quote}

\begin{sphinxVerbatim}[commandchars=\\\{\}]
\PYG{c+c1}{\PYGZsh{} Retrieve the dimension of the LMI constraint}
\PYG{n}{conidx} \PYG{o}{=} \PYG{n}{con}\PYG{o}{.}\PYG{n}{getDim}\PYG{p}{(}\PYG{p}{)}
\end{sphinxVerbatim}

\subsubsection{LmiConstraint.getLen()}
\label{\detokenize{pyapiref:lmiconstraint-getlen}}\begin{quote}

\sphinxAtStartPar
\sphinxstylestrong{Synopsis}
\begin{quote}

\sphinxAtStartPar
\sphinxcode{\sphinxupquote{getLen()}}
\end{quote}

\sphinxAtStartPar
\sphinxstylestrong{Description}
\begin{quote}

\sphinxAtStartPar
Retrieve the flattened length of the LMI constraint.
\end{quote}

\sphinxAtStartPar
\sphinxstylestrong{Example}
\end{quote}

\begin{sphinxVerbatim}[commandchars=\\\{\}]
\PYG{c+c1}{\PYGZsh{} Retrieve the flattened length of the LMI constraint}
\PYG{n}{conidx} \PYG{o}{=} \PYG{n}{con}\PYG{o}{.}\PYG{n}{getDim}\PYG{p}{(}\PYG{p}{)}
\end{sphinxVerbatim}

\subsubsection{LmiConstraint.setName()}
\label{\detokenize{pyapiref:lmiconstraint-setname}}\begin{quote}

\sphinxAtStartPar
\sphinxstylestrong{Synopsis}
\begin{quote}

\sphinxAtStartPar
\sphinxcode{\sphinxupquote{setName(newname)}}
\end{quote}

\sphinxAtStartPar
\sphinxstylestrong{Description}
\begin{quote}

\sphinxAtStartPar
Set the name of the LMI constraint to \sphinxcode{\sphinxupquote{newname}} .
\end{quote}

\sphinxAtStartPar
\sphinxstylestrong{Arguments}
\begin{quote}

\sphinxAtStartPar
\sphinxcode{\sphinxupquote{newname}}
\begin{quote}

\sphinxAtStartPar
The name of the LMI constraint to be set.
\end{quote}
\end{quote}

\sphinxAtStartPar
\sphinxstylestrong{Example}
\end{quote}

\begin{sphinxVerbatim}[commandchars=\\\{\}]
\PYG{c+c1}{\PYGZsh{} Set the name of the LMI constraint}
\PYG{n}{con}\PYG{o}{.}\PYG{n}{setName}\PYG{p}{(}\PYG{l+s+s1}{\PYGZsq{}}\PYG{l+s+s1}{con}\PYG{l+s+s1}{\PYGZsq{}}\PYG{p}{)}
\end{sphinxVerbatim}

\subsubsection{LmiConstraint.setRhs()}
\label{\detokenize{pyapiref:lmiconstraint-setrhs}}\begin{quote}

\sphinxAtStartPar
\sphinxstylestrong{Synopsis}
\begin{quote}

\sphinxAtStartPar
\sphinxcode{\sphinxupquote{setRhs(mat)}}
\end{quote}

\sphinxAtStartPar
\sphinxstylestrong{Description}
\begin{quote}

\sphinxAtStartPar
Set the constant\sphinxhyphen{}term symmetric matrix of the LMI constraint.
\end{quote}

\sphinxAtStartPar
\sphinxstylestrong{Arguments}
\begin{quote}

\sphinxAtStartPar
\sphinxcode{\sphinxupquote{mat}}
\begin{quote}

\sphinxAtStartPar
The constant\sphinxhyphen{}term symmetric matrix of the LMI constraint to be set.
It should be {\hyperref[\detokenize{pyapiref:chappyapi-symmatrix}]{\sphinxcrossref{\DUrole{std,std-ref}{SymMatrix Class}}}} .
\end{quote}
\end{quote}

\sphinxAtStartPar
\sphinxstylestrong{Example}
\end{quote}

\begin{sphinxVerbatim}[commandchars=\\\{\}]
\PYG{c+c1}{\PYGZsh{} Set the constant\PYGZhy{}term symmetric of the LMI constraint}
\PYG{n}{D} \PYG{o}{=} \PYG{n}{m}\PYG{o}{.}\PYG{n}{addSparseMat}\PYG{p}{(}\PYG{l+m+mi}{2}\PYG{p}{,} \PYG{p}{[}\PYG{l+m+mi}{0}\PYG{p}{,} \PYG{l+m+mi}{1}\PYG{p}{]}\PYG{p}{,} \PYG{p}{[}\PYG{l+m+mi}{0}\PYG{p}{,} \PYG{l+m+mi}{1}\PYG{p}{]}\PYG{p}{,} \PYG{p}{[}\PYG{l+m+mf}{1.0}\PYG{p}{,} \PYG{l+m+mf}{1.0}\PYG{p}{]}\PYG{p}{)}
\PYG{n}{con}\PYG{o}{.}\PYG{n}{setRhs}\PYG{p}{(}\PYG{n}{D}\PYG{p}{)}
\end{sphinxVerbatim}

\subsubsection{LmiConstraint.getInfo()}
\label{\detokenize{pyapiref:lmiconstraint-getinfo}}\begin{quote}

\sphinxAtStartPar
\sphinxstylestrong{Synopsis}
\begin{quote}

\sphinxAtStartPar
\sphinxcode{\sphinxupquote{getInfo(infoname)}}
\end{quote}

\sphinxAtStartPar
\sphinxstylestrong{Description}
\begin{quote}

\sphinxAtStartPar
Retrieve the specified information with the name \sphinxcode{\sphinxupquote{infoname}} .
\end{quote}

\sphinxAtStartPar
\sphinxstylestrong{Arguments}
\begin{quote}

\sphinxAtStartPar
\sphinxcode{\sphinxupquote{infoname}}
\begin{quote}

\sphinxAtStartPar
Name of the information to be obtained. Please refer to {\hyperref[\detokenize{pyapiref:chappyapi-const-info}]{\sphinxcrossref{\DUrole{std,std-ref}{Information}}}} section for possible values.
\end{quote}
\end{quote}

\sphinxAtStartPar
\sphinxstylestrong{Example}
\end{quote}

\begin{sphinxVerbatim}[commandchars=\\\{\}]
\PYG{c+c1}{\PYGZsh{} Get slack of the LMI constraint con}
\PYG{n}{conlb} \PYG{o}{=} \PYG{n}{con}\PYG{o}{.}\PYG{n}{getInfo}\PYG{p}{(}\PYG{n}{COPT}\PYG{o}{.}\PYG{n}{Info}\PYG{o}{.}\PYG{n}{Slack}\PYG{p}{)}
\end{sphinxVerbatim}

\subsubsection{LmiConstraint.remove()}
\label{\detokenize{pyapiref:lmiconstraint-remove}}\begin{quote}

\sphinxAtStartPar
\sphinxstylestrong{Synopsis}
\begin{quote}

\sphinxAtStartPar
\sphinxcode{\sphinxupquote{remove()}}
\end{quote}

\sphinxAtStartPar
\sphinxstylestrong{Description}
\begin{quote}

\sphinxAtStartPar
Remove the current LMI constraint from the model.
\end{quote}

\sphinxAtStartPar
\sphinxstylestrong{Example}
\end{quote}

\begin{sphinxVerbatim}[commandchars=\\\{\}]
\PYG{c+c1}{\PYGZsh{} Remove the LMI constraint conx}
\PYG{n}{conx}\PYG{o}{.}\PYG{n}{remove}\PYG{p}{(}\PYG{p}{)}
\end{sphinxVerbatim}

\subsubsection{LmiConstraint.shape}
\label{\detokenize{pyapiref:lmiconstraint-shape}}\begin{quote}

\sphinxAtStartPar
\sphinxstylestrong{Synopsis}
\begin{quote}

\sphinxAtStartPar
\sphinxcode{\sphinxupquote{shape}}
\end{quote}

\sphinxAtStartPar
\sphinxstylestrong{Description}
\begin{quote}

\sphinxAtStartPar
Shape of the \sphinxcode{\sphinxupquote{LmiConstraint}} object.
\end{quote}

\sphinxAtStartPar
\sphinxstylestrong{Return value}
\begin{quote}

\sphinxAtStartPar
Integer tuple.
\end{quote}
\end{quote}

\subsubsection{LmiConstraint.size}
\label{\detokenize{pyapiref:lmiconstraint-size}}\begin{quote}

\sphinxAtStartPar
\sphinxstylestrong{Synopsis}
\begin{quote}

\sphinxAtStartPar
\sphinxcode{\sphinxupquote{size}}
\end{quote}

\sphinxAtStartPar
\sphinxstylestrong{Description}
\begin{quote}

\sphinxAtStartPar
Size of the \sphinxcode{\sphinxupquote{LmiConstraint}} object.
\end{quote}

\sphinxAtStartPar
\sphinxstylestrong{Return value}
\begin{quote}

\sphinxAtStartPar
Integer tuple.
\end{quote}
\end{quote}

\subsubsection{LmiConstraint.dim}
\label{\detokenize{pyapiref:lmiconstraint-dim}}\begin{quote}

\sphinxAtStartPar
\sphinxstylestrong{Synopsis}
\begin{quote}

\sphinxAtStartPar
\sphinxcode{\sphinxupquote{dim}}
\end{quote}

\sphinxAtStartPar
\sphinxstylestrong{Description}
\begin{quote}

\sphinxAtStartPar
Dimension of the \sphinxcode{\sphinxupquote{LmiConstraint}} object.
\end{quote}

\sphinxAtStartPar
\sphinxstylestrong{Return value}
\begin{quote}

\sphinxAtStartPar
Integer.
\end{quote}
\end{quote}

\subsubsection{LmiConstraint.len}
\label{\detokenize{pyapiref:lmiconstraint-len}}\begin{quote}

\sphinxAtStartPar
\sphinxstylestrong{Synopsis}
\begin{quote}

\sphinxAtStartPar
\sphinxcode{\sphinxupquote{len}}
\end{quote}

\sphinxAtStartPar
\sphinxstylestrong{Description}
\begin{quote}

\sphinxAtStartPar
Flattened length of the \sphinxcode{\sphinxupquote{LmiConstraint}} object.
\end{quote}

\sphinxAtStartPar
\sphinxstylestrong{Return value}
\begin{quote}

\sphinxAtStartPar
Integer.
\end{quote}
\end{quote}

\subsection{LmiConstrArray Class}
\label{\detokenize{pyapiref:lmiconstrarray-class}}\label{\detokenize{pyapiref:chappyapi-lmiconstrarray}}
\sphinxAtStartPar
To facilitate users to operate on multiple {\hyperref[\detokenize{pyapiref:chappyapi-lmiconstraint}]{\sphinxcrossref{\DUrole{std,std-ref}{LmiConstraint Class}}}}  objects,
the Python interface of COPT provides LmiConstrArray class with the following methods:

\subsubsection{LmiConstrArray()}
\label{\detokenize{pyapiref:lmiconstrarray}}\begin{quote}

\sphinxAtStartPar
\sphinxstylestrong{Synopsis}
\begin{quote}

\sphinxAtStartPar
\sphinxcode{\sphinxupquote{LmiConstrArray(constrs=None)}}
\end{quote}

\sphinxAtStartPar
\sphinxstylestrong{Description}
\begin{quote}

\sphinxAtStartPar
Create a {\hyperref[\detokenize{pyapiref:chappyapi-lmiconstrarray}]{\sphinxcrossref{\DUrole{std,std-ref}{LmiConstrArray Class}}}} object.

\sphinxAtStartPar
If parameter \sphinxcode{\sphinxupquote{constrs}} is \sphinxcode{\sphinxupquote{None}}, the create an empty {\hyperref[\detokenize{pyapiref:chappyapi-lmiconstrarray}]{\sphinxcrossref{\DUrole{std,std-ref}{LmiConstrArray Class}}}} object,
otherwise initialize the newly created {\hyperref[\detokenize{pyapiref:chappyapi-lmiconstrarray}]{\sphinxcrossref{\DUrole{std,std-ref}{LmiConstrArray Class}}}} object with parameter \sphinxcode{\sphinxupquote{constrs}} .
\end{quote}

\sphinxAtStartPar
\sphinxstylestrong{Arguments}
\begin{quote}

\sphinxAtStartPar
\sphinxcode{\sphinxupquote{constrs}}
\begin{quote}

\sphinxAtStartPar
Can be {\hyperref[\detokenize{pyapiref:chappyapi-lmiconstraint}]{\sphinxcrossref{\DUrole{std,std-ref}{LmiConstraint Class}}}} object,
{\hyperref[\detokenize{pyapiref:chappyapi-lmiconstrarray}]{\sphinxcrossref{\DUrole{std,std-ref}{LmiConstrArray Class}}}} object, list, dictionary or {\hyperref[\detokenize{pyapiref:chappyapi-util-tupledict}]{\sphinxcrossref{\DUrole{std,std-ref}{tupledict Class}}}} object.
\end{quote}
\end{quote}

\sphinxAtStartPar
\sphinxstylestrong{Example}
\end{quote}

\begin{sphinxVerbatim}[commandchars=\\\{\}]
\PYG{c+c1}{\PYGZsh{} Create an empty LmiConstrArray class object}
\PYG{n}{conarr} \PYG{o}{=} \PYG{n}{LmiConstrArray}\PYG{p}{(}\PYG{p}{)}
\PYG{c+c1}{\PYGZsh{} Create a LmiConstrArray class object and initialize it with the LMI constraints conx and cony}
\PYG{n}{conarr} \PYG{o}{=} \PYG{n}{LmiConstrArray}\PYG{p}{(}\PYG{p}{[}\PYG{n}{conx}\PYG{p}{,} \PYG{n}{cony}\PYG{p}{]}\PYG{p}{)}
\end{sphinxVerbatim}

\subsubsection{LmiConstrArray.pushBack()}
\label{\detokenize{pyapiref:lmiconstrarray-pushback}}\begin{quote}

\sphinxAtStartPar
\sphinxstylestrong{Synopsis}
\begin{quote}

\sphinxAtStartPar
\sphinxcode{\sphinxupquote{pushBack(constr)}}
\end{quote}

\sphinxAtStartPar
\sphinxstylestrong{Description}
\begin{quote}

\sphinxAtStartPar
Add single or multiple {\hyperref[\detokenize{pyapiref:chappyapi-lmiconstraint}]{\sphinxcrossref{\DUrole{std,std-ref}{LmiConstraint Class}}}} objects.
\end{quote}

\sphinxAtStartPar
\sphinxstylestrong{Arguments}
\begin{quote}

\sphinxAtStartPar
\sphinxcode{\sphinxupquote{constr}}
\begin{quote}

\sphinxAtStartPar
\sphinxcode{\sphinxupquote{constrs}} can be {\hyperref[\detokenize{pyapiref:chappyapi-lmiconstraint}]{\sphinxcrossref{\DUrole{std,std-ref}{LmiConstraint Class}}}} object, {\hyperref[\detokenize{pyapiref:chappyapi-lmiconstrarray}]{\sphinxcrossref{\DUrole{std,std-ref}{LmiConstrArray Class}}}} object,
list, dictionary or {\hyperref[\detokenize{pyapiref:chappyapi-util-tupledict}]{\sphinxcrossref{\DUrole{std,std-ref}{tupledict Class}}}} object.
\end{quote}
\end{quote}

\sphinxAtStartPar
\sphinxstylestrong{Example}
\end{quote}

\begin{sphinxVerbatim}[commandchars=\\\{\}]
\PYG{c+c1}{\PYGZsh{} Add LMI constraint r to conarr}
\PYG{n}{conarr}\PYG{o}{.}\PYG{n}{pushBack}\PYG{p}{(}\PYG{n}{r}\PYG{p}{)}
\PYG{c+c1}{\PYGZsh{} Add LMI constraint r0 and r1 to conarr}
\PYG{n}{conarr}\PYG{o}{.}\PYG{n}{pushBack}\PYG{p}{(}\PYG{p}{[}\PYG{n}{r0}\PYG{p}{,} \PYG{n}{r1}\PYG{p}{]}\PYG{p}{)}
\end{sphinxVerbatim}

\subsubsection{LmiConstrArray.getLmiConstr()}
\label{\detokenize{pyapiref:lmiconstrarray-getlmiconstr}}\begin{quote}

\sphinxAtStartPar
\sphinxstylestrong{Synopsis}
\begin{quote}

\sphinxAtStartPar
\sphinxcode{\sphinxupquote{getLmiConstr(idx)}}
\end{quote}

\sphinxAtStartPar
\sphinxstylestrong{Description}
\begin{quote}

\sphinxAtStartPar
Retrieve the LMI constraint according to its subscript in {\hyperref[\detokenize{pyapiref:chappyapi-lmiconstrarray}]{\sphinxcrossref{\DUrole{std,std-ref}{LmiConstrArray Class}}}} object.
Return a {\hyperref[\detokenize{pyapiref:chappyapi-lmiconstraint}]{\sphinxcrossref{\DUrole{std,std-ref}{LmiConstraint Class}}}} object.
\end{quote}

\sphinxAtStartPar
\sphinxstylestrong{Arguments}
\begin{quote}

\sphinxAtStartPar
\sphinxcode{\sphinxupquote{idx}}
\begin{quote}

\sphinxAtStartPar
Subscript of the desired LMI constraint in {\hyperref[\detokenize{pyapiref:chappyapi-lmiconstrarray}]{\sphinxcrossref{\DUrole{std,std-ref}{LmiConstrArray Class}}}} object, starting with 0.
\end{quote}
\end{quote}

\sphinxAtStartPar
\sphinxstylestrong{Example}
\end{quote}

\begin{sphinxVerbatim}[commandchars=\\\{\}]
\PYG{c+c1}{\PYGZsh{} Get the LMI constraint with subscript 1 in conarr}
\PYG{n}{con} \PYG{o}{=} \PYG{n}{conarr}\PYG{o}{.}\PYG{n}{getLmiConstr}\PYG{p}{(}\PYG{l+m+mi}{1}\PYG{p}{)}
\end{sphinxVerbatim}

\subsubsection{LmiConstrArray.getSize()}
\label{\detokenize{pyapiref:lmiconstrarray-getsize}}\begin{quote}

\sphinxAtStartPar
\sphinxstylestrong{Synopsis}
\begin{quote}

\sphinxAtStartPar
\sphinxcode{\sphinxupquote{getSize()}}
\end{quote}

\sphinxAtStartPar
\sphinxstylestrong{Description}
\begin{quote}

\sphinxAtStartPar
Get the number of elements in {\hyperref[\detokenize{pyapiref:chappyapi-lmiconstrarray}]{\sphinxcrossref{\DUrole{std,std-ref}{LmiConstrArray Class}}}} object.
\end{quote}

\sphinxAtStartPar
\sphinxstylestrong{Example}
\end{quote}

\begin{sphinxVerbatim}[commandchars=\\\{\}]
\PYG{c+c1}{\PYGZsh{} Get the number of LMI constraints in conarr}
\PYG{n}{arrsize} \PYG{o}{=} \PYG{n}{conarr}\PYG{o}{.}\PYG{n}{getSize}\PYG{p}{(}\PYG{p}{)}
\end{sphinxVerbatim}

\subsubsection{LmiConstrArray.reserve()}
\label{\detokenize{pyapiref:lmiconstrarray-reserve}}\begin{quote}

\sphinxAtStartPar
\sphinxstylestrong{Synopsis}
\begin{quote}

\sphinxAtStartPar
\sphinxcode{\sphinxupquote{reserve(n)}}
\end{quote}

\sphinxAtStartPar
\sphinxstylestrong{Description}
\begin{quote}

\sphinxAtStartPar
Reserve space for {\hyperref[\detokenize{pyapiref:chappyapi-lmiconstrarray}]{\sphinxcrossref{\DUrole{std,std-ref}{LmiConstrArray Class}}}} objects of size n.
\end{quote}

\sphinxAtStartPar
\sphinxstylestrong{Arguments}
\begin{quote}

\sphinxAtStartPar
\sphinxcode{\sphinxupquote{n}}
\begin{quote}

\sphinxAtStartPar
The number of elements in the object {\hyperref[\detokenize{pyapiref:chappyapi-lmiconstrarray}]{\sphinxcrossref{\DUrole{std,std-ref}{LmiConstrArray Class}}}} .
\end{quote}
\end{quote}
\end{quote}

\subsection{NlConstraint Class}
\label{\detokenize{pyapiref:nlconstraint-class}}\label{\detokenize{pyapiref:chappyapi-nlconstraint}}
\sphinxAtStartPar
The \sphinxcode{\sphinxupquote{NlConstraint}} class provides an interface for operations
on nonlinear constraints in the COPT. It offers the following member functions:

\subsubsection{NlConstraint.getName()}
\label{\detokenize{pyapiref:nlconstraint-getname}}\begin{quote}

\sphinxAtStartPar
\sphinxstylestrong{Synopsis}
\begin{quote}

\sphinxAtStartPar
\sphinxcode{\sphinxupquote{getName()}}
\end{quote}

\sphinxAtStartPar
\sphinxstylestrong{Description}
\begin{quote}

\sphinxAtStartPar
Retrieve the name of the nonlinear constraint.
\end{quote}

\sphinxAtStartPar
\sphinxstylestrong{Return Value}
\begin{quote}

\sphinxAtStartPar
A string.
\end{quote}
\end{quote}

\subsubsection{NlConstraint.getIdx()}
\label{\detokenize{pyapiref:nlconstraint-getidx}}\begin{quote}

\sphinxAtStartPar
\sphinxstylestrong{Synopsis}
\begin{quote}

\sphinxAtStartPar
\sphinxcode{\sphinxupquote{getIdx()}}
\end{quote}

\sphinxAtStartPar
\sphinxstylestrong{Description}
\begin{quote}

\sphinxAtStartPar
Retrieve the index of the nonlinear constraint in the model.
\end{quote}

\sphinxAtStartPar
\sphinxstylestrong{Return Value}
\begin{quote}

\sphinxAtStartPar
An integer.
\end{quote}
\end{quote}

\subsubsection{NlConstraint.setName()}
\label{\detokenize{pyapiref:nlconstraint-setname}}\begin{quote}

\sphinxAtStartPar
\sphinxstylestrong{Synopsis}
\begin{quote}

\sphinxAtStartPar
\sphinxcode{\sphinxupquote{setName(newname)}}
\end{quote}

\sphinxAtStartPar
\sphinxstylestrong{Description}
\begin{quote}

\sphinxAtStartPar
Set the name of the nonlinear constraint.
\end{quote}

\sphinxAtStartPar
\sphinxstylestrong{Arguments}
\begin{quote}

\sphinxAtStartPar
\sphinxcode{\sphinxupquote{newname}}
\begin{quote}

\sphinxAtStartPar
The new name of the constraint.
\end{quote}
\end{quote}
\end{quote}

\subsubsection{NlConstraint.getRhs()}
\label{\detokenize{pyapiref:nlconstraint-getrhs}}\begin{quote}

\sphinxAtStartPar
\sphinxstylestrong{Synopsis}
\begin{quote}

\sphinxAtStartPar
\sphinxcode{\sphinxupquote{getRhs()}}
\end{quote}

\sphinxAtStartPar
\sphinxstylestrong{Description}
\begin{quote}

\sphinxAtStartPar
Retrieve the right hand side of nonlinear constraint.
\end{quote}

\sphinxAtStartPar
\sphinxstylestrong{Return Value}
\begin{quote}

\sphinxAtStartPar
A double value.
\end{quote}
\end{quote}

\subsubsection{NlConstraint.getSense()}
\label{\detokenize{pyapiref:nlconstraint-getsense}}\begin{quote}

\sphinxAtStartPar
\sphinxstylestrong{Synopsis}
\begin{quote}

\sphinxAtStartPar
\sphinxcode{\sphinxupquote{getSense()}}
\end{quote}

\sphinxAtStartPar
\sphinxstylestrong{Description}
\begin{quote}

\sphinxAtStartPar
Retrieve the type of nonlinear constraint.
\end{quote}

\sphinxAtStartPar
\sphinxstylestrong{Return Value}
\begin{quote}

\sphinxAtStartPar
A string.
\end{quote}
\end{quote}

\subsubsection{NlConstraint.setRhs()}
\label{\detokenize{pyapiref:nlconstraint-setrhs}}\begin{quote}

\sphinxAtStartPar
\sphinxstylestrong{Synopsis}
\begin{quote}

\sphinxAtStartPar
\sphinxcode{\sphinxupquote{setRhs(rhs)}}
\end{quote}

\sphinxAtStartPar
\sphinxstylestrong{Description}
\begin{quote}

\sphinxAtStartPar
Set the right hand side value of the nonlinear constraint.
\end{quote}

\sphinxAtStartPar
\sphinxstylestrong{Arguments}
\begin{quote}

\sphinxAtStartPar
\sphinxcode{\sphinxupquote{rhs}}
\begin{quote}

\sphinxAtStartPar
The new right hand side value.
\end{quote}
\end{quote}
\end{quote}

\subsubsection{NlConstraint.setSense()}
\label{\detokenize{pyapiref:nlconstraint-setsense}}\begin{quote}

\sphinxAtStartPar
\sphinxstylestrong{Synopsis}
\begin{quote}

\sphinxAtStartPar
\sphinxcode{\sphinxupquote{setSense(sense)}}
\end{quote}

\sphinxAtStartPar
\sphinxstylestrong{Description}
\begin{quote}

\sphinxAtStartPar
Set the sense of nonlinear constraint.
\end{quote}

\sphinxAtStartPar
\sphinxstylestrong{Arguments}
\begin{quote}

\sphinxAtStartPar
\sphinxcode{\sphinxupquote{sense}}
\begin{quote}

\sphinxAtStartPar
The sense of nonlinear constraint to be set.
\end{quote}
\end{quote}
\end{quote}

\subsubsection{NlConstraint.getInfo()}
\label{\detokenize{pyapiref:nlconstraint-getinfo}}\begin{quote}

\sphinxAtStartPar
\sphinxstylestrong{Synopsis}
\begin{quote}

\sphinxAtStartPar
\sphinxcode{\sphinxupquote{getInfo(infoname)}}
\end{quote}

\sphinxAtStartPar
\sphinxstylestrong{Description}
\begin{quote}

\sphinxAtStartPar
Retrieve specified information of the nonlinear constraint.
\end{quote}

\sphinxAtStartPar
\sphinxstylestrong{Arguments}
\begin{quote}

\sphinxAtStartPar
\sphinxcode{\sphinxupquote{infoname}}
\begin{quote}

\sphinxAtStartPar
Name of the information to be retrieved.

\sphinxAtStartPar
See {\hyperref[\detokenize{pyapiref:chappyapi-const-info}]{\sphinxcrossref{\DUrole{std,std-ref}{Information}}}} for possible values.
\end{quote}
\end{quote}

\sphinxAtStartPar
\sphinxstylestrong{Return Value}
\begin{quote}

\sphinxAtStartPar
A double value.
\end{quote}
\end{quote}

\subsubsection{NlConstraint.setInfo()}
\label{\detokenize{pyapiref:nlconstraint-setinfo}}\begin{quote}

\sphinxAtStartPar
\sphinxstylestrong{Synopsis}
\begin{quote}

\sphinxAtStartPar
\sphinxcode{\sphinxupquote{setInfo(infoname, newval)}}
\end{quote}

\sphinxAtStartPar
\sphinxstylestrong{Description}
\begin{quote}

\sphinxAtStartPar
Set the specified information value of the nonlinear constraint.
\end{quote}

\sphinxAtStartPar
\sphinxstylestrong{Arguments}
\begin{quote}

\sphinxAtStartPar
\sphinxcode{\sphinxupquote{infoname}}
\begin{quote}

\sphinxAtStartPar
Name of the information to be set.

\sphinxAtStartPar
See {\hyperref[\detokenize{pyapiref:chappyapi-const-info}]{\sphinxcrossref{\DUrole{std,std-ref}{Information}}}} for possible values.
\end{quote}

\sphinxAtStartPar
\sphinxcode{\sphinxupquote{newval}}
\begin{quote}

\sphinxAtStartPar
The new information value.
\end{quote}
\end{quote}
\end{quote}

\subsubsection{NlConstraint.remove()}
\label{\detokenize{pyapiref:nlconstraint-remove}}\begin{quote}

\sphinxAtStartPar
\sphinxstylestrong{Synopsis}
\begin{quote}

\sphinxAtStartPar
\sphinxcode{\sphinxupquote{remove()}}
\end{quote}

\sphinxAtStartPar
\sphinxstylestrong{Description}
\begin{quote}

\sphinxAtStartPar
Remove the current nonlinear constraint from the model.
\end{quote}
\end{quote}

\subsection{NlConstrArray Class}
\label{\detokenize{pyapiref:nlconstrarray-class}}\label{\detokenize{pyapiref:chappyapi-nlconstrarray}}
\sphinxAtStartPar
The \sphinxcode{\sphinxupquote{NlConstrArray}} class is used to manage and operate on a
collection of nonlinear constraints. The following methods are provided:

\subsubsection{NlConstrArray()}
\label{\detokenize{pyapiref:nlconstrarray}}\begin{quote}

\sphinxAtStartPar
\sphinxstylestrong{Synopsis}
\begin{quote}

\sphinxAtStartPar
\sphinxcode{\sphinxupquote{NlConstrArray(constrs=None)}}
\end{quote}

\sphinxAtStartPar
\sphinxstylestrong{Description}
\begin{quote}

\sphinxAtStartPar
Construct a new nonlinear constraint array object.

\sphinxAtStartPar
The argument \sphinxcode{\sphinxupquote{constrs}} can be a single {\hyperref[\detokenize{pyapiref:chappyapi-nlconstraint}]{\sphinxcrossref{\DUrole{std,std-ref}{NlConstraint Class}}}} object,
an iterable of {\hyperref[\detokenize{pyapiref:chappyapi-nlconstraint}]{\sphinxcrossref{\DUrole{std,std-ref}{NlConstraint Class}}}} objects, or a mapping from arbitrary keys to constraint objects.
\end{quote}

\sphinxAtStartPar
\sphinxstylestrong{Arguments}
\begin{quote}

\sphinxAtStartPar
\sphinxcode{\sphinxupquote{constrs}}
\begin{quote}

\sphinxAtStartPar
Optional. Defaults to \sphinxcode{\sphinxupquote{None}}.
\end{quote}
\end{quote}
\end{quote}

\subsubsection{NlConstrArray.pushBack()}
\label{\detokenize{pyapiref:nlconstrarray-pushback}}\begin{quote}

\sphinxAtStartPar
\sphinxstylestrong{Synopsis}
\begin{quote}

\sphinxAtStartPar
\sphinxcode{\sphinxupquote{pushBack(constr)}}
\end{quote}

\sphinxAtStartPar
\sphinxstylestrong{Description}
\begin{quote}

\sphinxAtStartPar
Append one or more {\hyperref[\detokenize{pyapiref:chappyapi-nlconstraint}]{\sphinxcrossref{\DUrole{std,std-ref}{NlConstraint Class}}}} objects to the current array.
\end{quote}

\sphinxAtStartPar
\sphinxstylestrong{Arguments}
\begin{quote}

\sphinxAtStartPar
\sphinxcode{\sphinxupquote{qconstrs}}
\begin{quote}

\sphinxAtStartPar
Could be a single {\hyperref[\detokenize{pyapiref:chappyapi-nlconstraint}]{\sphinxcrossref{\DUrole{std,std-ref}{NlConstraint Class}}}} object,
an iterable of such objects, or a mapping from keys to constraints.
\end{quote}
\end{quote}
\end{quote}

\subsubsection{NlConstrArray.getNlConstr()}
\label{\detokenize{pyapiref:nlconstrarray-getnlconstr}}\begin{quote}

\sphinxAtStartPar
\sphinxstylestrong{Synopsis}
\begin{quote}

\sphinxAtStartPar
\sphinxcode{\sphinxupquote{getNlConstr(idx)}}
\end{quote}

\sphinxAtStartPar
\sphinxstylestrong{Description}
\begin{quote}

\sphinxAtStartPar
Retrieve the constraint object at the specified index.
\end{quote}

\sphinxAtStartPar
\sphinxstylestrong{Arguments}
\begin{quote}

\sphinxAtStartPar
\sphinxcode{\sphinxupquote{idx}}
\begin{quote}

\sphinxAtStartPar
The index of the constraint.
\end{quote}
\end{quote}

\sphinxAtStartPar
\sphinxstylestrong{Return Value}
\begin{quote}

\sphinxAtStartPar
The corresponding {\hyperref[\detokenize{pyapiref:chappyapi-nlconstraint}]{\sphinxcrossref{\DUrole{std,std-ref}{NlConstraint Class}}}} object.
\end{quote}
\end{quote}

\subsubsection{NlConstrArray.getSize()}
\label{\detokenize{pyapiref:nlconstrarray-getsize}}\begin{quote}

\sphinxAtStartPar
\sphinxstylestrong{Synopsis}
\begin{quote}

\sphinxAtStartPar
\sphinxcode{\sphinxupquote{getSize()}}
\end{quote}

\sphinxAtStartPar
\sphinxstylestrong{Description}
\begin{quote}

\sphinxAtStartPar
Obtain the number of elements in {\hyperref[\detokenize{pyapiref:chappyapi-nlconstrarray}]{\sphinxcrossref{\DUrole{std,std-ref}{NlConstrArray Class}}}} object.
\end{quote}

\sphinxAtStartPar
\sphinxstylestrong{Return Value}
\begin{quote}

\sphinxAtStartPar
An integer value.
\end{quote}
\end{quote}

\subsection{NlConstrBuilder Class}
\label{\detokenize{pyapiref:nlconstrbuilder-class}}\label{\detokenize{pyapiref:chappyapi-nlconstrbuilder}}
\sphinxAtStartPar
The \sphinxcode{\sphinxupquote{NlConstrBuilder}} class provides a builder interface for defining nonlinear
constraints in the COPT. It provides the following functions:

\subsubsection{NlConstrBuilder()}
\label{\detokenize{pyapiref:nlconstrbuilder}}\begin{quote}

\sphinxAtStartPar
\sphinxstylestrong{Synopsis}
\begin{quote}

\sphinxAtStartPar
\sphinxcode{\sphinxupquote{NlConstrBuilder()}}
\end{quote}

\sphinxAtStartPar
\sphinxstylestrong{Description}
\begin{quote}

\sphinxAtStartPar
Construct a new nonlinear constraint builder object.
\end{quote}
\end{quote}

\subsubsection{NlConstrBuilder.setBuilder()}
\label{\detokenize{pyapiref:nlconstrbuilder-setbuilder}}\begin{quote}

\sphinxAtStartPar
\sphinxstylestrong{Synopsis}
\begin{quote}

\sphinxAtStartPar
\sphinxcode{\sphinxupquote{setBuilder(expr, sense, rhs)}}
\end{quote}

\sphinxAtStartPar
\sphinxstylestrong{Description}
\begin{quote}

\sphinxAtStartPar
Define the nonlinear expression, constraint sense,
and right hand side for the builder.
\end{quote}

\sphinxAtStartPar
\sphinxstylestrong{Arguments}
\begin{quote}

\sphinxAtStartPar
\sphinxcode{\sphinxupquote{expr}}
\begin{quote}

\sphinxAtStartPar
The nonlinear expression.
\end{quote}

\sphinxAtStartPar
\sphinxcode{\sphinxupquote{sense}}
\begin{quote}

\sphinxAtStartPar
The constraint sense.
\end{quote}

\sphinxAtStartPar
\sphinxcode{\sphinxupquote{rhs}}
\begin{quote}

\sphinxAtStartPar
The right hand side.
\end{quote}
\end{quote}
\end{quote}

\subsubsection{NlConstrBuilder.getNlExpr()}
\label{\detokenize{pyapiref:nlconstrbuilder-getnlexpr}}\begin{quote}

\sphinxAtStartPar
\sphinxstylestrong{Synopsis}
\begin{quote}

\sphinxAtStartPar
\sphinxcode{\sphinxupquote{getNlExpr()}}
\end{quote}

\sphinxAtStartPar
\sphinxstylestrong{Description}
\begin{quote}

\sphinxAtStartPar
Retrieve the nonlinear expression in the builder.
\end{quote}

\sphinxAtStartPar
\sphinxstylestrong{Return Value}
\begin{quote}

\sphinxAtStartPar
A {\hyperref[\detokenize{pyapiref:chappyapi-nlexpr}]{\sphinxcrossref{\DUrole{std,std-ref}{NlExpr Class}}}} object.
\end{quote}
\end{quote}

\subsubsection{NlConstrBuilder.getSense()}
\label{\detokenize{pyapiref:nlconstrbuilder-getsense}}\begin{quote}

\sphinxAtStartPar
\sphinxstylestrong{Synopsis}
\begin{quote}

\sphinxAtStartPar
\sphinxcode{\sphinxupquote{getSense()}}
\end{quote}

\sphinxAtStartPar
\sphinxstylestrong{Description}
\begin{quote}

\sphinxAtStartPar
Retreive the constraint sense of the nonlinear constraint builder.
\end{quote}

\sphinxAtStartPar
\sphinxstylestrong{Return Value}
\begin{quote}

\sphinxAtStartPar
A string.
\end{quote}
\end{quote}

\subsection{NlConstrBuilderArray Class}
\label{\detokenize{pyapiref:nlconstrbuilderarray-class}}\label{\detokenize{pyapiref:chappyapi-nlconstrbuilderarray}}
\sphinxAtStartPar
The NlConstrBuilderArray class is used to
manage a collection of nonlinear constraint builders.

\sphinxAtStartPar
It provides the following member functions:

\subsubsection{NlConstrBuilderArray()}
\label{\detokenize{pyapiref:nlconstrbuilderarray}}\begin{quote}

\sphinxAtStartPar
\sphinxstylestrong{Synopsis}
\begin{quote}

\sphinxAtStartPar
\sphinxcode{\sphinxupquote{NlConstrBuilderArray(constrs=None)}}
\end{quote}

\sphinxAtStartPar
\sphinxstylestrong{Description}
\begin{quote}

\sphinxAtStartPar
Construct a new nonlinear constraint builder array.

\sphinxAtStartPar
The argument \sphinxcode{\sphinxupquote{constrs}} can be a single {\hyperref[\detokenize{pyapiref:chappyapi-nlconstrbuilder}]{\sphinxcrossref{\DUrole{std,std-ref}{NlConstrBuilder Class}}}} object,
an iterable of builder objects, or a mapping from arbitrary keys to builder objects.
\end{quote}

\sphinxAtStartPar
\sphinxstylestrong{Arguments}
\begin{quote}

\sphinxAtStartPar
\sphinxcode{\sphinxupquote{constrs}}
\begin{quote}

\sphinxAtStartPar
Optional. Defaults to \sphinxcode{\sphinxupquote{None}}.
\end{quote}
\end{quote}
\end{quote}

\subsubsection{NlConstrBuilderArray.pushBack()}
\label{\detokenize{pyapiref:nlconstrbuilderarray-pushback}}\begin{quote}

\sphinxAtStartPar
\sphinxstylestrong{Synopsis}
\begin{quote}

\sphinxAtStartPar
\sphinxcode{\sphinxupquote{pushBack(constrbuilder)}}
\end{quote}

\sphinxAtStartPar
\sphinxstylestrong{Description}
\begin{quote}

\sphinxAtStartPar
Append one or more {\hyperref[\detokenize{pyapiref:chappyapi-nlconstrbuilder}]{\sphinxcrossref{\DUrole{std,std-ref}{NlConstrBuilder Class}}}} objects to the array.
\end{quote}

\sphinxAtStartPar
\sphinxstylestrong{Arguments}
\begin{quote}

\sphinxAtStartPar
\sphinxcode{\sphinxupquote{constrbuilder}}
\begin{quote}

\sphinxAtStartPar
The builder(s) to append.

\sphinxAtStartPar
Can be a single {\hyperref[\detokenize{pyapiref:chappyapi-nlconstrbuilder}]{\sphinxcrossref{\DUrole{std,std-ref}{NlConstrBuilder Class}}}} object,
a {\hyperref[\detokenize{pyapiref:chappyapi-nlconstrbuilderarray}]{\sphinxcrossref{\DUrole{std,std-ref}{NlConstrBuilderArray Class}}}}, a list,
a dictionary, or a {\hyperref[\detokenize{pyapiref:chappyapi-util-tupledict}]{\sphinxcrossref{\DUrole{std,std-ref}{tupledict Class}}}}.
\end{quote}
\end{quote}
\end{quote}

\subsubsection{NlConstrBuilderArray.getBuilder()}
\label{\detokenize{pyapiref:nlconstrbuilderarray-getbuilder}}\begin{quote}

\sphinxAtStartPar
\sphinxstylestrong{Synopsis}
\begin{quote}

\sphinxAtStartPar
\sphinxcode{\sphinxupquote{getBuilder(idx)}}
\end{quote}

\sphinxAtStartPar
\sphinxstylestrong{Description}
\begin{quote}

\sphinxAtStartPar
Retrieve the builder object at the specified index.
\end{quote}

\sphinxAtStartPar
\sphinxstylestrong{Arguments}
\begin{quote}

\sphinxAtStartPar
\sphinxcode{\sphinxupquote{idx}}
\begin{quote}

\sphinxAtStartPar
The index of the builder.
\end{quote}
\end{quote}

\sphinxAtStartPar
\sphinxstylestrong{Return Value}
\begin{quote}

\sphinxAtStartPar
The corresponding {\hyperref[\detokenize{pyapiref:chappyapi-nlconstrbuilder}]{\sphinxcrossref{\DUrole{std,std-ref}{NlConstrBuilder Class}}}} object.
\end{quote}
\end{quote}

\subsubsection{NlConstrBuilderArray.getSize()}
\label{\detokenize{pyapiref:nlconstrbuilderarray-getsize}}\begin{quote}

\sphinxAtStartPar
\sphinxstylestrong{Synopsis}
\begin{quote}

\sphinxAtStartPar
\sphinxcode{\sphinxupquote{getSize()}}
\end{quote}

\sphinxAtStartPar
\sphinxstylestrong{Description}
\begin{quote}

\sphinxAtStartPar
Retrieve the number of elements
in the {\hyperref[\detokenize{pyapiref:chappyapi-nlconstrbuilderarray}]{\sphinxcrossref{\DUrole{std,std-ref}{NlConstrBuilderArray Class}}}} object.
\end{quote}

\sphinxAtStartPar
\sphinxstylestrong{Return Value}
\begin{quote}

\sphinxAtStartPar
An integer.
\end{quote}
\end{quote}

\subsection{SOS Class}
\label{\detokenize{pyapiref:sos-class}}\label{\detokenize{pyapiref:chappyapi-sos}}
\sphinxAtStartPar
SOS object contains related operations of COPT SOS constraints.
The following methods are provided:

\subsubsection{SOS.getIdx()}
\label{\detokenize{pyapiref:sos-getidx}}\begin{quote}

\sphinxAtStartPar
\sphinxstylestrong{Synopsis}
\begin{quote}

\sphinxAtStartPar
\sphinxcode{\sphinxupquote{getIdx()}}
\end{quote}

\sphinxAtStartPar
\sphinxstylestrong{Description}
\begin{quote}

\sphinxAtStartPar
Retrieve the subscript of SOS constraint in model.
\end{quote}

\sphinxAtStartPar
\sphinxstylestrong{Example}
\end{quote}

\begin{sphinxVerbatim}[commandchars=\\\{\}]
\PYG{c+c1}{\PYGZsh{} Retrieve the subscript of SOS constraint sosx.}
\PYG{n}{sosidx} \PYG{o}{=} \PYG{n}{sosx}\PYG{o}{.}\PYG{n}{getIdx}\PYG{p}{(}\PYG{p}{)}
\end{sphinxVerbatim}

\subsubsection{SOS.remove()}
\label{\detokenize{pyapiref:sos-remove}}\begin{quote}

\sphinxAtStartPar
\sphinxstylestrong{Synopsis}
\begin{quote}

\sphinxAtStartPar
\sphinxcode{\sphinxupquote{remove()}}
\end{quote}

\sphinxAtStartPar
\sphinxstylestrong{Description}
\begin{quote}

\sphinxAtStartPar
Delete the SOS constraint from model.
\end{quote}

\sphinxAtStartPar
\sphinxstylestrong{Example}
\end{quote}

\begin{sphinxVerbatim}[commandchars=\\\{\}]
\PYG{c+c1}{\PYGZsh{} Delete the SOS constraint \PYGZsq{}sosx\PYGZsq{}}
\PYG{n}{sosx}\PYG{o}{.}\PYG{n}{remove}\PYG{p}{(}\PYG{p}{)}
\end{sphinxVerbatim}

\subsection{SOSArray Class}
\label{\detokenize{pyapiref:sosarray-class}}\label{\detokenize{pyapiref:chappyapi-sosarray}}
\sphinxAtStartPar
To facilitate users to operate on a set of {\hyperref[\detokenize{pyapiref:chappyapi-sos}]{\sphinxcrossref{\DUrole{std,std-ref}{SOS Class}}}} objects, COPT designed SOSArray class in Python interface.
The following methods are provided:

\subsubsection{SOSArray()}
\label{\detokenize{pyapiref:sosarray}}\begin{quote}

\sphinxAtStartPar
\sphinxstylestrong{Synopsis}
\begin{quote}

\sphinxAtStartPar
\sphinxcode{\sphinxupquote{SOSArray(soss=None)}}
\end{quote}

\sphinxAtStartPar
\sphinxstylestrong{Description}
\begin{quote}

\sphinxAtStartPar
Create a {\hyperref[\detokenize{pyapiref:chappyapi-sosarray}]{\sphinxcrossref{\DUrole{std,std-ref}{SOSArray Class}}}} object.

\sphinxAtStartPar
If parameter \sphinxcode{\sphinxupquote{soss}} is \sphinxcode{\sphinxupquote{None}}, then build an empty {\hyperref[\detokenize{pyapiref:chappyapi-sosarray}]{\sphinxcrossref{\DUrole{std,std-ref}{SOSArray Class}}}} object,
otherwise initialize the newly created {\hyperref[\detokenize{pyapiref:chappyapi-sosarray}]{\sphinxcrossref{\DUrole{std,std-ref}{SOSArray Class}}}} object with \sphinxcode{\sphinxupquote{soss}}.
\end{quote}

\sphinxAtStartPar
\sphinxstylestrong{Arguments}
\begin{quote}

\sphinxAtStartPar
\sphinxcode{\sphinxupquote{soss}}
\begin{quote}

\sphinxAtStartPar
SOS constraint to be added. Optional, \sphinxcode{\sphinxupquote{None}} by default.
It can be {\hyperref[\detokenize{pyapiref:chappyapi-sos}]{\sphinxcrossref{\DUrole{std,std-ref}{SOS Class}}}} object, {\hyperref[\detokenize{pyapiref:chappyapi-sosarray}]{\sphinxcrossref{\DUrole{std,std-ref}{SOSArray Class}}}} object, list, dictionary or
{\hyperref[\detokenize{pyapiref:chappyapi-util-tupledict}]{\sphinxcrossref{\DUrole{std,std-ref}{tupledict Class}}}} object.
\end{quote}
\end{quote}

\sphinxAtStartPar
\sphinxstylestrong{Example}
\end{quote}

\begin{sphinxVerbatim}[commandchars=\\\{\}]
\PYG{c+c1}{\PYGZsh{} Create a new SOSArray object}
\PYG{n}{sosarr} \PYG{o}{=} \PYG{n}{SOSArray}\PYG{p}{(}\PYG{p}{)}
\PYG{c+c1}{\PYGZsh{} Create a SOSArray object, and initialize it with SOS constraints sosx and sosy.}
\PYG{n}{sosarr} \PYG{o}{=} \PYG{n}{SOSArray}\PYG{p}{(}\PYG{p}{[}\PYG{n}{sosx}\PYG{p}{,} \PYG{n}{sosy}\PYG{p}{]}\PYG{p}{)}
\end{sphinxVerbatim}

\subsubsection{SOSArray.pushBack()}
\label{\detokenize{pyapiref:sosarray-pushback}}\begin{quote}

\sphinxAtStartPar
\sphinxstylestrong{Synopsis}
\begin{quote}

\sphinxAtStartPar
\sphinxcode{\sphinxupquote{pushBack(sos)}}
\end{quote}

\sphinxAtStartPar
\sphinxstylestrong{Description}
\begin{quote}

\sphinxAtStartPar
Add one or multiple {\hyperref[\detokenize{pyapiref:chappyapi-sos}]{\sphinxcrossref{\DUrole{std,std-ref}{SOS Class}}}} objects.
\end{quote}

\sphinxAtStartPar
\sphinxstylestrong{Arguments}
\begin{quote}

\sphinxAtStartPar
\sphinxcode{\sphinxupquote{sos}}
\begin{quote}

\sphinxAtStartPar
SOS constraints to be added, which can be {\hyperref[\detokenize{pyapiref:chappyapi-sos}]{\sphinxcrossref{\DUrole{std,std-ref}{SOS Class}}}} object, {\hyperref[\detokenize{pyapiref:chappyapi-sosarray}]{\sphinxcrossref{\DUrole{std,std-ref}{SOSArray Class}}}} object, list, dictionary
or {\hyperref[\detokenize{pyapiref:chappyapi-util-tupledict}]{\sphinxcrossref{\DUrole{std,std-ref}{tupledict Class}}}} object.
\end{quote}
\end{quote}

\sphinxAtStartPar
\sphinxstylestrong{Example}
\end{quote}

\begin{sphinxVerbatim}[commandchars=\\\{\}]
\PYG{c+c1}{\PYGZsh{} Add SOS constraint sosx to sosarr}
\PYG{n}{sosarr}\PYG{o}{.}\PYG{n}{pushBack}\PYG{p}{(}\PYG{n}{sosx}\PYG{p}{)}
\PYG{c+c1}{\PYGZsh{} Add SOS constraints sosx and sosy to sosarr}
\PYG{n}{sosarr}\PYG{o}{.}\PYG{n}{pushBack}\PYG{p}{(}\PYG{p}{[}\PYG{n}{sosx}\PYG{p}{,} \PYG{n}{sosy}\PYG{p}{]}\PYG{p}{)}
\end{sphinxVerbatim}

\subsubsection{SOSArray.getSOS()}
\label{\detokenize{pyapiref:sosarray-getsos}}\begin{quote}

\sphinxAtStartPar
\sphinxstylestrong{Synopsis}
\begin{quote}

\sphinxAtStartPar
\sphinxcode{\sphinxupquote{getSOS(idx)}}
\end{quote}

\sphinxAtStartPar
\sphinxstylestrong{Description}
\begin{quote}

\sphinxAtStartPar
Retrieve the corresponding SOS constraint according to its subscript in {\hyperref[\detokenize{pyapiref:chappyapi-sosarray}]{\sphinxcrossref{\DUrole{std,std-ref}{SOSArray Class}}}} object and
return a {\hyperref[\detokenize{pyapiref:chappyapi-sos}]{\sphinxcrossref{\DUrole{std,std-ref}{SOS Class}}}} object.
\end{quote}

\sphinxAtStartPar
\sphinxstylestrong{Arguments}
\begin{quote}

\sphinxAtStartPar
\sphinxcode{\sphinxupquote{idx}}
\begin{quote}

\sphinxAtStartPar
Indice of the SOS constraint in {\hyperref[\detokenize{pyapiref:chappyapi-sosarray}]{\sphinxcrossref{\DUrole{std,std-ref}{SOSArray Class}}}} object, starting with 0.
\end{quote}
\end{quote}

\sphinxAtStartPar
\sphinxstylestrong{Example}
\end{quote}

\begin{sphinxVerbatim}[commandchars=\\\{\}]
\PYG{c+c1}{\PYGZsh{} Retrieve the SOS constraint with indice of 1 in sosarr}
\PYG{n}{sos} \PYG{o}{=} \PYG{n}{sosarr}\PYG{o}{.}\PYG{n}{getSOS}\PYG{p}{(}\PYG{l+m+mi}{1}\PYG{p}{)}
\end{sphinxVerbatim}

\subsubsection{SOSArray.getSize()}
\label{\detokenize{pyapiref:sosarray-getsize}}\begin{quote}

\sphinxAtStartPar
\sphinxstylestrong{Synopsis}
\begin{quote}

\sphinxAtStartPar
\sphinxcode{\sphinxupquote{getSize()}}
\end{quote}

\sphinxAtStartPar
\sphinxstylestrong{Description}
\begin{quote}

\sphinxAtStartPar
Retrieve the number of elements in {\hyperref[\detokenize{pyapiref:chappyapi-sosarray}]{\sphinxcrossref{\DUrole{std,std-ref}{SOSArray Class}}}} object.
\end{quote}

\sphinxAtStartPar
\sphinxstylestrong{Example}
\end{quote}

\begin{sphinxVerbatim}[commandchars=\\\{\}]
\PYG{c+c1}{\PYGZsh{} Retrieve the number of SOS constraints in sosarr.}
\PYG{n}{arrsize} \PYG{o}{=} \PYG{n}{sosarr}\PYG{o}{.}\PYG{n}{getSize}\PYG{p}{(}\PYG{p}{)}
\end{sphinxVerbatim}

\subsection{SOSBuilder Class}
\label{\detokenize{pyapiref:sosbuilder-class}}\label{\detokenize{pyapiref:chappyapi-sosbuilder}}
\sphinxAtStartPar
For easy access builders of SOS constraints, SOSBuilder class provides the following methods:

\subsubsection{SOSBuilder()}
\label{\detokenize{pyapiref:sosbuilder}}\begin{quote}

\sphinxAtStartPar
\sphinxstylestrong{Synopsis}
\begin{quote}

\sphinxAtStartPar
\sphinxcode{\sphinxupquote{SOSBuilder()}}
\end{quote}

\sphinxAtStartPar
\sphinxstylestrong{Description}
\begin{quote}

\sphinxAtStartPar
Create an empty {\hyperref[\detokenize{pyapiref:chappyapi-sosbuilder}]{\sphinxcrossref{\DUrole{std,std-ref}{SOSBuilder Class}}}} object.
\end{quote}

\sphinxAtStartPar
\sphinxstylestrong{Example}
\end{quote}

\begin{sphinxVerbatim}[commandchars=\\\{\}]
\PYG{c+c1}{\PYGZsh{} Create an empty SOSBuilder object.}
\PYG{n}{sosbuilder} \PYG{o}{=} \PYG{n}{SOSBuilder}\PYG{p}{(}\PYG{p}{)}
\end{sphinxVerbatim}

\subsubsection{SOSBuilder.setBuilder()}
\label{\detokenize{pyapiref:sosbuilder-setbuilder}}\begin{quote}

\sphinxAtStartPar
\sphinxstylestrong{Synopsis}
\begin{quote}

\sphinxAtStartPar
\sphinxcode{\sphinxupquote{setBuilder(sostype, vars, weights=None)}}
\end{quote}

\sphinxAtStartPar
\sphinxstylestrong{Description}
\begin{quote}

\sphinxAtStartPar
Set type, variable, weight of variable on {\hyperref[\detokenize{pyapiref:chappyapi-sosbuilder}]{\sphinxcrossref{\DUrole{std,std-ref}{SOSBuilder Class}}}} object.
\end{quote}

\sphinxAtStartPar
\sphinxstylestrong{Arguments}
\begin{quote}

\sphinxAtStartPar
\sphinxcode{\sphinxupquote{sostype}}
\begin{quote}

\sphinxAtStartPar
SOS constraint type. Full list of available types can be found in {\hyperref[\detokenize{constant:chapconst-sostype}]{\sphinxcrossref{\DUrole{std,std-ref}{SOS\sphinxhyphen{}constraint types}}}}.
\end{quote}

\sphinxAtStartPar
\sphinxcode{\sphinxupquote{vars}}
\begin{quote}

\sphinxAtStartPar
Variables of SOS constarint, which can be {\hyperref[\detokenize{pyapiref:chappyapi-vararray}]{\sphinxcrossref{\DUrole{std,std-ref}{VarArray Class}}}} object, list, dictionary or
{\hyperref[\detokenize{pyapiref:chappyapi-util-tupledict}]{\sphinxcrossref{\DUrole{std,std-ref}{tupledict Class}}}} object.
\end{quote}
\begin{description}
\sphinxlineitem{\sphinxcode{\sphinxupquote{weights}}}
\sphinxAtStartPar
Weights of variables in SOS constraint. Optional, \sphinxcode{\sphinxupquote{None}} by default.
Could be list, dictionary or {\hyperref[\detokenize{pyapiref:chappyapi-util-tupledict}]{\sphinxcrossref{\DUrole{std,std-ref}{tupledict Class}}}} object.

\end{description}
\end{quote}

\sphinxAtStartPar
\sphinxstylestrong{Example}
\end{quote}

\begin{sphinxVerbatim}[commandchars=\\\{\}]
\PYG{c+c1}{\PYGZsh{} Set the type of SOS constraint builder as SOS1, variables x and y, weights of variables as 1 and 2 respectively.}
\PYG{n}{sosbuilder}\PYG{o}{.}\PYG{n}{setBuilder}\PYG{p}{(}\PYG{n}{COPT}\PYG{o}{.}\PYG{n}{SOS\PYGZus{}TYPE1}\PYG{p}{,} \PYG{p}{[}\PYG{n}{x}\PYG{p}{,} \PYG{n}{y}\PYG{p}{]}\PYG{p}{,} \PYG{p}{[}\PYG{l+m+mi}{1}\PYG{p}{,} \PYG{l+m+mi}{2}\PYG{p}{]}\PYG{p}{)}
\end{sphinxVerbatim}

\subsubsection{SOSBuilder.getType()}
\label{\detokenize{pyapiref:sosbuilder-gettype}}\begin{quote}

\sphinxAtStartPar
\sphinxstylestrong{Synopsis}
\begin{quote}

\sphinxAtStartPar
\sphinxcode{\sphinxupquote{getType()}}
\end{quote}

\sphinxAtStartPar
\sphinxstylestrong{Description}
\begin{quote}

\sphinxAtStartPar
Retrieve the SOS constraint type of {\hyperref[\detokenize{pyapiref:chappyapi-sosbuilder}]{\sphinxcrossref{\DUrole{std,std-ref}{SOSBuilder Class}}}} object.
\end{quote}

\sphinxAtStartPar
\sphinxstylestrong{Example}
\end{quote}

\begin{sphinxVerbatim}[commandchars=\\\{\}]
\PYG{c+c1}{\PYGZsh{} Retrieve the type of SOS constraint builder sosx.}
\PYG{n}{sostype} \PYG{o}{=} \PYG{n}{sosbuilder}\PYG{o}{.}\PYG{n}{getType}\PYG{p}{(}\PYG{n}{sosx}\PYG{p}{)}
\end{sphinxVerbatim}

\subsubsection{SOSBuilder.getVar()}
\label{\detokenize{pyapiref:sosbuilder-getvar}}\begin{quote}

\sphinxAtStartPar
\sphinxstylestrong{Synopsis}
\begin{quote}

\sphinxAtStartPar
\sphinxcode{\sphinxupquote{getVar(idx)}}
\end{quote}

\sphinxAtStartPar
\sphinxstylestrong{Description}
\begin{quote}

\sphinxAtStartPar
Retrieve the corresponding variables according to its indice in {\hyperref[\detokenize{pyapiref:chappyapi-sosbuilder}]{\sphinxcrossref{\DUrole{std,std-ref}{SOSBuilder Class}}}} object, and return a
{\hyperref[\detokenize{pyapiref:chappyapi-var}]{\sphinxcrossref{\DUrole{std,std-ref}{Var Class}}}} object.
\end{quote}

\sphinxAtStartPar
\sphinxstylestrong{Arguments}
\begin{quote}

\sphinxAtStartPar
\sphinxcode{\sphinxupquote{idx}}
\begin{quote}

\sphinxAtStartPar
Indice of the variable in {\hyperref[\detokenize{pyapiref:chappyapi-sosbuilder}]{\sphinxcrossref{\DUrole{std,std-ref}{SOSBuilder Class}}}} object, starting with 0.
\end{quote}
\end{quote}

\sphinxAtStartPar
\sphinxstylestrong{Example}
\end{quote}

\begin{sphinxVerbatim}[commandchars=\\\{\}]
\PYG{c+c1}{\PYGZsh{} Retrieve the variable in SOS constraint builder sosx with indice of 1}
\PYG{n}{sosvar} \PYG{o}{=} \PYG{n}{sosx}\PYG{o}{.}\PYG{n}{getVar}\PYG{p}{(}\PYG{l+m+mi}{1}\PYG{p}{)}
\end{sphinxVerbatim}

\subsubsection{SOSBuilder.getVars()}
\label{\detokenize{pyapiref:sosbuilder-getvars}}\begin{quote}

\sphinxAtStartPar
\sphinxstylestrong{Synopsis}
\begin{quote}

\sphinxAtStartPar
\sphinxcode{\sphinxupquote{getVars()}}
\end{quote}

\sphinxAtStartPar
\sphinxstylestrong{Description}
\begin{quote}

\sphinxAtStartPar
Retrieve all variables in {\hyperref[\detokenize{pyapiref:chappyapi-sosbuilder}]{\sphinxcrossref{\DUrole{std,std-ref}{SOSBuilder Class}}}} objects, and return a {\hyperref[\detokenize{pyapiref:chappyapi-vararray}]{\sphinxcrossref{\DUrole{std,std-ref}{VarArray Class}}}} object.
\end{quote}

\sphinxAtStartPar
\sphinxstylestrong{Example}
\end{quote}

\begin{sphinxVerbatim}[commandchars=\\\{\}]
\PYG{c+c1}{\PYGZsh{} Retrieve all variables in SOS constraint builder sosx.}
\PYG{n}{sosvars} \PYG{o}{=} \PYG{n}{sosx}\PYG{o}{.}\PYG{n}{getVars}\PYG{p}{(}\PYG{p}{)}
\end{sphinxVerbatim}

\subsubsection{SOSBuilder.getWeight()}
\label{\detokenize{pyapiref:sosbuilder-getweight}}\begin{quote}

\sphinxAtStartPar
\sphinxstylestrong{Synopsis}
\begin{quote}

\sphinxAtStartPar
\sphinxcode{\sphinxupquote{getWeight(idx)}}
\end{quote}

\sphinxAtStartPar
\sphinxstylestrong{Description}
\begin{quote}

\sphinxAtStartPar
Retrieve the corresponding weight of variable according to its indice in {\hyperref[\detokenize{pyapiref:chappyapi-sosbuilder}]{\sphinxcrossref{\DUrole{std,std-ref}{SOSBuilder Class}}}} object.
\end{quote}

\sphinxAtStartPar
\sphinxstylestrong{Arguments}
\begin{quote}

\sphinxAtStartPar
\sphinxcode{\sphinxupquote{idx}}
\begin{quote}

\sphinxAtStartPar
Indice of the variable in {\hyperref[\detokenize{pyapiref:chappyapi-sosbuilder}]{\sphinxcrossref{\DUrole{std,std-ref}{SOSBuilder Class}}}} object, starting with 0.
\end{quote}
\end{quote}

\sphinxAtStartPar
\sphinxstylestrong{Example}
\end{quote}

\begin{sphinxVerbatim}[commandchars=\\\{\}]
\PYG{c+c1}{\PYGZsh{} Retrieve the corresponding weight of variable according to its indice in the SOS constraint builder sosx.}
\PYG{n}{sosweight} \PYG{o}{=} \PYG{n}{sosx}\PYG{o}{.}\PYG{n}{getWeight}\PYG{p}{(}\PYG{l+m+mi}{1}\PYG{p}{)}
\end{sphinxVerbatim}

\subsubsection{SOSBuilder.getWeights()}
\label{\detokenize{pyapiref:sosbuilder-getweights}}\begin{quote}

\sphinxAtStartPar
\sphinxstylestrong{Synopsis}
\begin{quote}

\sphinxAtStartPar
\sphinxcode{\sphinxupquote{getWeights()}}
\end{quote}

\sphinxAtStartPar
\sphinxstylestrong{Description}
\begin{quote}

\sphinxAtStartPar
Retrieve weights of all the variables in {\hyperref[\detokenize{pyapiref:chappyapi-sosbuilder}]{\sphinxcrossref{\DUrole{std,std-ref}{SOSBuilder Class}}}} object.
\end{quote}

\sphinxAtStartPar
\sphinxstylestrong{Example}
\end{quote}

\begin{sphinxVerbatim}[commandchars=\\\{\}]
\PYG{c+c1}{\PYGZsh{} Retrieve weights of all the variables in SOS constraint builder sosx.}
\PYG{n}{sosweights} \PYG{o}{=} \PYG{n}{sosx}\PYG{o}{.}\PYG{n}{getWeights}\PYG{p}{(}\PYG{p}{)}
\end{sphinxVerbatim}

\subsubsection{SOSBuilder.getSize()}
\label{\detokenize{pyapiref:sosbuilder-getsize}}\begin{quote}

\sphinxAtStartPar
\sphinxstylestrong{Synopsis}
\begin{quote}

\sphinxAtStartPar
\sphinxcode{\sphinxupquote{getSize()}}
\end{quote}

\sphinxAtStartPar
\sphinxstylestrong{Description}
\begin{quote}

\sphinxAtStartPar
Retrieve the number of elements in {\hyperref[\detokenize{pyapiref:chappyapi-sosbuilder}]{\sphinxcrossref{\DUrole{std,std-ref}{SOSBuilder Class}}}} object.
\end{quote}

\sphinxAtStartPar
\sphinxstylestrong{Example}
\end{quote}

\begin{sphinxVerbatim}[commandchars=\\\{\}]
\PYG{c+c1}{\PYGZsh{} Retrieve the number of elements in SOS constraint builder sosx.}
\PYG{n}{sossize} \PYG{o}{=} \PYG{n}{sosx}\PYG{o}{.}\PYG{n}{getSize}\PYG{p}{(}\PYG{p}{)}
\end{sphinxVerbatim}

\subsection{SOSBuilderArray Class}
\label{\detokenize{pyapiref:sosbuilderarray-class}}\label{\detokenize{pyapiref:chappyapi-sosbuilderarray}}
\sphinxAtStartPar
In order to facilitate users to operate on a set of {\hyperref[\detokenize{pyapiref:chappyapi-sosbuilder}]{\sphinxcrossref{\DUrole{std,std-ref}{SOSBuilder Class}}}} objects, COPT provides SOSBuilderArray class
in Python interface, providing the following methods:

\subsubsection{SOSBuilderArray()}
\label{\detokenize{pyapiref:sosbuilderarray}}\begin{quote}

\sphinxAtStartPar
\sphinxstylestrong{Synopsis}
\begin{quote}

\sphinxAtStartPar
\sphinxcode{\sphinxupquote{SOSBuilderArray(sosbuilders=None)}}
\end{quote}

\sphinxAtStartPar
\sphinxstylestrong{Description}
\begin{quote}

\sphinxAtStartPar
Create a {\hyperref[\detokenize{pyapiref:chappyapi-sosbuilderarray}]{\sphinxcrossref{\DUrole{std,std-ref}{SOSBuilderArray Class}}}} object.

\sphinxAtStartPar
If parameter \sphinxcode{\sphinxupquote{sosbuilders}} is \sphinxcode{\sphinxupquote{None}}, then create an empty {\hyperref[\detokenize{pyapiref:chappyapi-sosbuilderarray}]{\sphinxcrossref{\DUrole{std,std-ref}{SOSBuilderArray Class}}}} object,
otherwise initialize the newly created {\hyperref[\detokenize{pyapiref:chappyapi-sosbuilderarray}]{\sphinxcrossref{\DUrole{std,std-ref}{SOSBuilderArray Class}}}} object with parameter \sphinxcode{\sphinxupquote{sosbuilders}}.
\end{quote}

\sphinxAtStartPar
\sphinxstylestrong{Arguments}
\begin{quote}

\sphinxAtStartPar
\sphinxcode{\sphinxupquote{sosbuilders}}
\begin{quote}

\sphinxAtStartPar
SOS constraint builder to be added. Optional, \sphinxcode{\sphinxupquote{None}} by default.
Could be {\hyperref[\detokenize{pyapiref:chappyapi-sosbuilder}]{\sphinxcrossref{\DUrole{std,std-ref}{SOSBuilder Class}}}} object, {\hyperref[\detokenize{pyapiref:chappyapi-sosbuilderarray}]{\sphinxcrossref{\DUrole{std,std-ref}{SOSBuilderArray Class}}}} object, list, dictionary
or {\hyperref[\detokenize{pyapiref:chappyapi-util-tupledict}]{\sphinxcrossref{\DUrole{std,std-ref}{tupledict Class}}}} object.
\end{quote}
\end{quote}

\sphinxAtStartPar
\sphinxstylestrong{Example}
\end{quote}

\begin{sphinxVerbatim}[commandchars=\\\{\}]
\PYG{c+c1}{\PYGZsh{} Create an empty SOSBuilderArray object.}
\PYG{n}{sosbuilderarr} \PYG{o}{=} \PYG{n}{SOSBuilderArray}\PYG{p}{(}\PYG{p}{)}
\PYG{c+c1}{\PYGZsh{} Create a SOSBuilderArray object and initialize it with SOS constraint builder sosx and sosy}
\PYG{n}{sosbuilderarr} \PYG{o}{=} \PYG{n}{SOSBuilderArray}\PYG{p}{(}\PYG{p}{[}\PYG{n}{sosx}\PYG{p}{,} \PYG{n}{sosy}\PYG{p}{]}\PYG{p}{)}
\end{sphinxVerbatim}

\subsubsection{SOSBuilderArray.pushBack()}
\label{\detokenize{pyapiref:sosbuilderarray-pushback}}\begin{quote}

\sphinxAtStartPar
\sphinxstylestrong{Synopsis}
\begin{quote}

\sphinxAtStartPar
\sphinxcode{\sphinxupquote{pushBack(sosbuilder)}}
\end{quote}

\sphinxAtStartPar
\sphinxstylestrong{Description}
\begin{quote}

\sphinxAtStartPar
Add one or multiple {\hyperref[\detokenize{pyapiref:chappyapi-sosbuilder}]{\sphinxcrossref{\DUrole{std,std-ref}{SOSBuilder Class}}}} objects.
\end{quote}

\sphinxAtStartPar
\sphinxstylestrong{Arguments}
\begin{quote}

\sphinxAtStartPar
\sphinxcode{\sphinxupquote{sosbuilder}}
\begin{quote}

\sphinxAtStartPar
SOS constraint builderto be added. Cound be {\hyperref[\detokenize{pyapiref:chappyapi-sosbuilder}]{\sphinxcrossref{\DUrole{std,std-ref}{SOSBuilder Class}}}} object, {\hyperref[\detokenize{pyapiref:chappyapi-sosbuilderarray}]{\sphinxcrossref{\DUrole{std,std-ref}{SOSBuilderArray Class}}}} object,
list, dictionary or {\hyperref[\detokenize{pyapiref:chappyapi-util-tupledict}]{\sphinxcrossref{\DUrole{std,std-ref}{tupledict Class}}}} object.
\end{quote}
\end{quote}

\sphinxAtStartPar
\sphinxstylestrong{Example}
\end{quote}

\begin{sphinxVerbatim}[commandchars=\\\{\}]
\PYG{c+c1}{\PYGZsh{} Add SOS constraint builder sosx to sosbuilderarr}
\PYG{n}{sosbuilderarr}\PYG{o}{.}\PYG{n}{pushBack}\PYG{p}{(}\PYG{n}{sosx}\PYG{p}{)}
\end{sphinxVerbatim}

\subsubsection{SOSBuilderArray.getBuilder()}
\label{\detokenize{pyapiref:sosbuilderarray-getbuilder}}\begin{quote}

\sphinxAtStartPar
\sphinxstylestrong{Synopsis}
\begin{quote}

\sphinxAtStartPar
\sphinxcode{\sphinxupquote{getBuilder(idx)}}
\end{quote}

\sphinxAtStartPar
\sphinxstylestrong{Description}
\begin{quote}

\sphinxAtStartPar
Retrieve the corresponding builder according to the indice of SOS constraint builder in {\hyperref[\detokenize{pyapiref:chappyapi-sosbuilderarray}]{\sphinxcrossref{\DUrole{std,std-ref}{SOSBuilderArray Class}}}} object.
\end{quote}

\sphinxAtStartPar
\sphinxstylestrong{Arguments}
\begin{quote}

\sphinxAtStartPar
\sphinxcode{\sphinxupquote{idx}}
\begin{quote}

\sphinxAtStartPar
Indice of the SOS constraint builder in {\hyperref[\detokenize{pyapiref:chappyapi-sosbuilderarray}]{\sphinxcrossref{\DUrole{std,std-ref}{SOSBuilderArray Class}}}} object, starting with 0.
\end{quote}
\end{quote}

\sphinxAtStartPar
\sphinxstylestrong{Example}
\end{quote}

\begin{sphinxVerbatim}[commandchars=\\\{\}]
\PYG{c+c1}{\PYGZsh{} Retrieve the SOS constraint builder with indice of 1 in sosbuilderarr}
\PYG{n}{sosbuilder} \PYG{o}{=} \PYG{n}{sosbuilderarr}\PYG{o}{.}\PYG{n}{getBuilder}\PYG{p}{(}\PYG{l+m+mi}{1}\PYG{p}{)}
\end{sphinxVerbatim}

\subsubsection{SOSBuilderArray.getSize()}
\label{\detokenize{pyapiref:sosbuilderarray-getsize}}\begin{quote}

\sphinxAtStartPar
\sphinxstylestrong{Synopsis}
\begin{quote}

\sphinxAtStartPar
\sphinxcode{\sphinxupquote{getSize()}}
\end{quote}

\sphinxAtStartPar
\sphinxstylestrong{Description}
\begin{quote}

\sphinxAtStartPar
Retrieve the number of elements in {\hyperref[\detokenize{pyapiref:chappyapi-sosbuilderarray}]{\sphinxcrossref{\DUrole{std,std-ref}{SOSBuilderArray Class}}}} object.
\end{quote}

\sphinxAtStartPar
\sphinxstylestrong{Example}
\end{quote}

\begin{sphinxVerbatim}[commandchars=\\\{\}]
\PYG{c+c1}{\PYGZsh{} Retrieve the number of elements in sosbuilderarr}
\PYG{n}{sosbuildersize} \PYG{o}{=} \PYG{n}{sosbuilderarr}\PYG{o}{.}\PYG{n}{getSize}\PYG{p}{(}\PYG{p}{)}
\end{sphinxVerbatim}

\subsection{Cone Class}
\label{\detokenize{pyapiref:cone-class}}\label{\detokenize{pyapiref:chappyapi-cone}}
\sphinxAtStartPar
Cone object contains related operations of COPT Second\sphinxhyphen{}Order\sphinxhyphen{}Cone (SOC) constraints.
The following methods are provided:

\subsubsection{Cone.getIdx()}
\label{\detokenize{pyapiref:cone-getidx}}\begin{quote}

\sphinxAtStartPar
\sphinxstylestrong{Synopsis}
\begin{quote}

\sphinxAtStartPar
\sphinxcode{\sphinxupquote{getIdx()}}
\end{quote}

\sphinxAtStartPar
\sphinxstylestrong{Description}
\begin{quote}

\sphinxAtStartPar
Retrieve the subscript of SOC constraint in model.
\end{quote}

\sphinxAtStartPar
\sphinxstylestrong{Example}
\end{quote}

\begin{sphinxVerbatim}[commandchars=\\\{\}]
\PYG{c+c1}{\PYGZsh{} Retrieve the subscript of SOC constraint cone.}
\PYG{n}{coneidx} \PYG{o}{=} \PYG{n}{cone}\PYG{o}{.}\PYG{n}{getIdx}\PYG{p}{(}\PYG{p}{)}
\end{sphinxVerbatim}

\subsubsection{Cone.remove()}
\label{\detokenize{pyapiref:cone-remove}}\begin{quote}

\sphinxAtStartPar
\sphinxstylestrong{Synopsis}
\begin{quote}

\sphinxAtStartPar
\sphinxcode{\sphinxupquote{remove()}}
\end{quote}

\sphinxAtStartPar
\sphinxstylestrong{Description}
\begin{quote}

\sphinxAtStartPar
Delete the SOC constraint from model.
\end{quote}

\sphinxAtStartPar
\sphinxstylestrong{Example}
\end{quote}

\begin{sphinxVerbatim}[commandchars=\\\{\}]
\PYG{c+c1}{\PYGZsh{} Delete the SOC constraint \PYGZsq{}cone\PYGZsq{}}
\PYG{n}{cone}\PYG{o}{.}\PYG{n}{remove}\PYG{p}{(}\PYG{p}{)}
\end{sphinxVerbatim}

\subsection{ConeArray Class}
\label{\detokenize{pyapiref:conearray-class}}\label{\detokenize{pyapiref:chappyapi-conearray}}
\sphinxAtStartPar
To facilitate users to operate on a set of {\hyperref[\detokenize{pyapiref:chappyapi-cone}]{\sphinxcrossref{\DUrole{std,std-ref}{Cone Class}}}} objects, COPT designed ConeArray class
in Python interface. The following methods are provided:

\subsubsection{ConeArray()}
\label{\detokenize{pyapiref:conearray}}\begin{quote}

\sphinxAtStartPar
\sphinxstylestrong{Synopsis}
\begin{quote}

\sphinxAtStartPar
\sphinxcode{\sphinxupquote{ConeArray(cones=None)}}
\end{quote}

\sphinxAtStartPar
\sphinxstylestrong{Description}
\begin{quote}

\sphinxAtStartPar
Create a {\hyperref[\detokenize{pyapiref:chappyapi-conearray}]{\sphinxcrossref{\DUrole{std,std-ref}{ConeArray Class}}}} object.

\sphinxAtStartPar
If parameter \sphinxcode{\sphinxupquote{cones}} is \sphinxcode{\sphinxupquote{None}}, then build an empty {\hyperref[\detokenize{pyapiref:chappyapi-conearray}]{\sphinxcrossref{\DUrole{std,std-ref}{ConeArray Class}}}} object,
otherwise initialize the newly created {\hyperref[\detokenize{pyapiref:chappyapi-conearray}]{\sphinxcrossref{\DUrole{std,std-ref}{ConeArray Class}}}} object with \sphinxcode{\sphinxupquote{cones}}.
\end{quote}

\sphinxAtStartPar
\sphinxstylestrong{Arguments}
\begin{quote}

\sphinxAtStartPar
\sphinxcode{\sphinxupquote{cones}}
\begin{quote}

\sphinxAtStartPar
Second\sphinxhyphen{}Order\sphinxhyphen{}Cone constraint to be added. Optional, \sphinxcode{\sphinxupquote{None}} by default.
It can be {\hyperref[\detokenize{pyapiref:chappyapi-cone}]{\sphinxcrossref{\DUrole{std,std-ref}{Cone Class}}}} object, {\hyperref[\detokenize{pyapiref:chappyapi-conearray}]{\sphinxcrossref{\DUrole{std,std-ref}{ConeArray Class}}}} object, list, dictionary or
{\hyperref[\detokenize{pyapiref:chappyapi-util-tupledict}]{\sphinxcrossref{\DUrole{std,std-ref}{tupledict Class}}}} object.
\end{quote}
\end{quote}

\sphinxAtStartPar
\sphinxstylestrong{Example}
\end{quote}

\begin{sphinxVerbatim}[commandchars=\\\{\}]
\PYG{c+c1}{\PYGZsh{} Create a new ConeArray object}
\PYG{n}{conearr} \PYG{o}{=} \PYG{n}{ConeArray}\PYG{p}{(}\PYG{p}{)}
\PYG{c+c1}{\PYGZsh{} Create a ConeArray object, and initialize it with SOC constraints conex and coney.}
\PYG{n}{conearr} \PYG{o}{=} \PYG{n}{ConeArray}\PYG{p}{(}\PYG{p}{[}\PYG{n}{conex}\PYG{p}{,} \PYG{n}{coney}\PYG{p}{]}\PYG{p}{)}
\end{sphinxVerbatim}

\subsubsection{ConeArray.pushBack()}
\label{\detokenize{pyapiref:conearray-pushback}}\begin{quote}

\sphinxAtStartPar
\sphinxstylestrong{Synopsis}
\begin{quote}

\sphinxAtStartPar
\sphinxcode{\sphinxupquote{pushBack(cone)}}
\end{quote}

\sphinxAtStartPar
\sphinxstylestrong{Description}
\begin{quote}

\sphinxAtStartPar
Add one or multiple {\hyperref[\detokenize{pyapiref:chappyapi-cone}]{\sphinxcrossref{\DUrole{std,std-ref}{Cone Class}}}} objects.
\end{quote}

\sphinxAtStartPar
\sphinxstylestrong{Arguments}
\begin{quote}

\sphinxAtStartPar
\sphinxcode{\sphinxupquote{cone}}
\begin{quote}

\sphinxAtStartPar
Second\sphinxhyphen{}Order\sphinxhyphen{}Cone constraints to be added, which can be {\hyperref[\detokenize{pyapiref:chappyapi-cone}]{\sphinxcrossref{\DUrole{std,std-ref}{Cone Class}}}} object,
{\hyperref[\detokenize{pyapiref:chappyapi-conearray}]{\sphinxcrossref{\DUrole{std,std-ref}{ConeArray Class}}}} object, list, dictionary or {\hyperref[\detokenize{pyapiref:chappyapi-util-tupledict}]{\sphinxcrossref{\DUrole{std,std-ref}{tupledict Class}}}} object.
\end{quote}
\end{quote}

\sphinxAtStartPar
\sphinxstylestrong{Example}
\end{quote}

\begin{sphinxVerbatim}[commandchars=\\\{\}]
\PYG{c+c1}{\PYGZsh{} Add SOC constraint conex to conearr}
\PYG{n}{conearr}\PYG{o}{.}\PYG{n}{pushBack}\PYG{p}{(}\PYG{n}{conex}\PYG{p}{)}
\PYG{c+c1}{\PYGZsh{} Add SOC constraints conex and coney to conearr}
\PYG{n}{conearr}\PYG{o}{.}\PYG{n}{pushBack}\PYG{p}{(}\PYG{p}{[}\PYG{n}{conex}\PYG{p}{,} \PYG{n}{coney}\PYG{p}{]}\PYG{p}{)}
\end{sphinxVerbatim}

\subsubsection{ConeArray.getCone()}
\label{\detokenize{pyapiref:conearray-getcone}}\begin{quote}

\sphinxAtStartPar
\sphinxstylestrong{Synopsis}
\begin{quote}

\sphinxAtStartPar
\sphinxcode{\sphinxupquote{getCone(idx)}}
\end{quote}

\sphinxAtStartPar
\sphinxstylestrong{Description}
\begin{quote}

\sphinxAtStartPar
Retrieve the corresponding Second\sphinxhyphen{}Order\sphinxhyphen{}Cone (SOC) constraint according to its subscript in
{\hyperref[\detokenize{pyapiref:chappyapi-conearray}]{\sphinxcrossref{\DUrole{std,std-ref}{ConeArray Class}}}} object and return a {\hyperref[\detokenize{pyapiref:chappyapi-cone}]{\sphinxcrossref{\DUrole{std,std-ref}{Cone Class}}}} object.
\end{quote}

\sphinxAtStartPar
\sphinxstylestrong{Arguments}
\begin{quote}

\sphinxAtStartPar
\sphinxcode{\sphinxupquote{idx}}
\begin{quote}

\sphinxAtStartPar
Indice of the SOC constraint in {\hyperref[\detokenize{pyapiref:chappyapi-conearray}]{\sphinxcrossref{\DUrole{std,std-ref}{ConeArray Class}}}} object, starting with 0.
\end{quote}
\end{quote}

\sphinxAtStartPar
\sphinxstylestrong{Example}
\end{quote}

\begin{sphinxVerbatim}[commandchars=\\\{\}]
\PYG{c+c1}{\PYGZsh{} Retrieve the SOC constraint with indice of 1 in conearr}
\PYG{n}{cone} \PYG{o}{=} \PYG{n}{conearr}\PYG{o}{.}\PYG{n}{getCone}\PYG{p}{(}\PYG{l+m+mi}{1}\PYG{p}{)}
\end{sphinxVerbatim}

\subsubsection{ConeArray.getSize()}
\label{\detokenize{pyapiref:conearray-getsize}}\begin{quote}

\sphinxAtStartPar
\sphinxstylestrong{Synopsis}
\begin{quote}

\sphinxAtStartPar
\sphinxcode{\sphinxupquote{getSize()}}
\end{quote}

\sphinxAtStartPar
\sphinxstylestrong{Description}
\begin{quote}

\sphinxAtStartPar
Retrieve the number of elements in {\hyperref[\detokenize{pyapiref:chappyapi-conearray}]{\sphinxcrossref{\DUrole{std,std-ref}{ConeArray Class}}}} object.
\end{quote}

\sphinxAtStartPar
\sphinxstylestrong{Example}
\end{quote}

\begin{sphinxVerbatim}[commandchars=\\\{\}]
\PYG{c+c1}{\PYGZsh{} Retrieve the number of SOC constraints in conearr.}
\PYG{n}{arrsize} \PYG{o}{=} \PYG{n}{conearr}\PYG{o}{.}\PYG{n}{getSize}\PYG{p}{(}\PYG{p}{)}
\end{sphinxVerbatim}

\subsection{ConeBuilder Class}
\label{\detokenize{pyapiref:conebuilder-class}}\label{\detokenize{pyapiref:chappyapi-conebuilder}}
\sphinxAtStartPar
For easy access builders of Second\sphinxhyphen{}Order\sphinxhyphen{}Cone (SOC) constraints, ConeBuilder class provides the following methods:

\subsubsection{ConeBuilder()}
\label{\detokenize{pyapiref:conebuilder}}\begin{quote}

\sphinxAtStartPar
\sphinxstylestrong{Synopsis}
\begin{quote}

\sphinxAtStartPar
\sphinxcode{\sphinxupquote{ConeBuilder()}}
\end{quote}

\sphinxAtStartPar
\sphinxstylestrong{Description}
\begin{quote}

\sphinxAtStartPar
Create an empty {\hyperref[\detokenize{pyapiref:chappyapi-conebuilder}]{\sphinxcrossref{\DUrole{std,std-ref}{ConeBuilder Class}}}} object.
\end{quote}

\sphinxAtStartPar
\sphinxstylestrong{Example}
\end{quote}

\begin{sphinxVerbatim}[commandchars=\\\{\}]
\PYG{c+c1}{\PYGZsh{} Create an empty ConeBuilder object.}
\PYG{n}{conebuilder} \PYG{o}{=} \PYG{n}{ConeBuilder}\PYG{p}{(}\PYG{p}{)}
\end{sphinxVerbatim}

\subsubsection{ConeBuilder.setBuilder()}
\label{\detokenize{pyapiref:conebuilder-setbuilder}}\begin{quote}

\sphinxAtStartPar
\sphinxstylestrong{Synopsis}
\begin{quote}

\sphinxAtStartPar
\sphinxcode{\sphinxupquote{setBuilder(conetype, vars)}}
\end{quote}

\sphinxAtStartPar
\sphinxstylestrong{Description}
\begin{quote}

\sphinxAtStartPar
Set type, variables of {\hyperref[\detokenize{pyapiref:chappyapi-conebuilder}]{\sphinxcrossref{\DUrole{std,std-ref}{ConeBuilder Class}}}} object.
\end{quote}

\sphinxAtStartPar
\sphinxstylestrong{Arguments}
\begin{quote}

\sphinxAtStartPar
\sphinxcode{\sphinxupquote{conetype}}
\begin{quote}

\sphinxAtStartPar
Type of Second\sphinxhyphen{}Order\sphinxhyphen{}Cone (SOC) constraint. Full list of available types can be found in
{\hyperref[\detokenize{constant:chapconst-conetype}]{\sphinxcrossref{\DUrole{std,std-ref}{SOC\sphinxhyphen{}constraint types}}}}.
\end{quote}

\sphinxAtStartPar
\sphinxcode{\sphinxupquote{vars}}
\begin{quote}

\sphinxAtStartPar
Variables of SOC constarint, which can be {\hyperref[\detokenize{pyapiref:chappyapi-vararray}]{\sphinxcrossref{\DUrole{std,std-ref}{VarArray Class}}}} object, list, dictionary or
{\hyperref[\detokenize{pyapiref:chappyapi-util-tupledict}]{\sphinxcrossref{\DUrole{std,std-ref}{tupledict Class}}}} object.
\end{quote}
\end{quote}

\sphinxAtStartPar
\sphinxstylestrong{Example}
\end{quote}

\begin{sphinxVerbatim}[commandchars=\\\{\}]
\PYG{c+c1}{\PYGZsh{} Set type as regular, variables as [z, x, y] for SOC constraint builder}
\PYG{n}{conebuilder}\PYG{o}{.}\PYG{n}{setBuilder}\PYG{p}{(}\PYG{n}{COPT}\PYG{o}{.}\PYG{n}{CONE\PYGZus{}QUAD}\PYG{p}{,} \PYG{p}{[}\PYG{n}{z}\PYG{p}{,} \PYG{n}{x}\PYG{p}{,} \PYG{n}{y}\PYG{p}{]}\PYG{p}{)}
\end{sphinxVerbatim}

\subsubsection{ConeBuilder.getType()}
\label{\detokenize{pyapiref:conebuilder-gettype}}\begin{quote}

\sphinxAtStartPar
\sphinxstylestrong{Synopsis}
\begin{quote}

\sphinxAtStartPar
\sphinxcode{\sphinxupquote{getType()}}
\end{quote}

\sphinxAtStartPar
\sphinxstylestrong{Description}
\begin{quote}

\sphinxAtStartPar
Retrieve the Second\sphinxhyphen{}Order\sphinxhyphen{}Cone (SOC) constraint type of {\hyperref[\detokenize{pyapiref:chappyapi-conebuilder}]{\sphinxcrossref{\DUrole{std,std-ref}{ConeBuilder Class}}}} object.
\end{quote}

\sphinxAtStartPar
\sphinxstylestrong{Example}
\end{quote}

\begin{sphinxVerbatim}[commandchars=\\\{\}]
\PYG{c+c1}{\PYGZsh{} Retrieve the type of SOC constraint builder conex.}
\PYG{n}{conetype} \PYG{o}{=} \PYG{n}{conebuilder}\PYG{o}{.}\PYG{n}{getType}\PYG{p}{(}\PYG{n}{conex}\PYG{p}{)}
\end{sphinxVerbatim}

\subsubsection{ConeBuilder.getVar()}
\label{\detokenize{pyapiref:conebuilder-getvar}}\begin{quote}

\sphinxAtStartPar
\sphinxstylestrong{Synopsis}
\begin{quote}

\sphinxAtStartPar
\sphinxcode{\sphinxupquote{getVar(idx)}}
\end{quote}

\sphinxAtStartPar
\sphinxstylestrong{Description}
\begin{quote}

\sphinxAtStartPar
Retrieve the corresponding variables according to its indice in {\hyperref[\detokenize{pyapiref:chappyapi-conebuilder}]{\sphinxcrossref{\DUrole{std,std-ref}{ConeBuilder Class}}}} object,
and return a {\hyperref[\detokenize{pyapiref:chappyapi-var}]{\sphinxcrossref{\DUrole{std,std-ref}{Var Class}}}} object.
\end{quote}

\sphinxAtStartPar
\sphinxstylestrong{Arguments}
\begin{quote}

\sphinxAtStartPar
\sphinxcode{\sphinxupquote{idx}}
\begin{quote}

\sphinxAtStartPar
Indice of the variable in {\hyperref[\detokenize{pyapiref:chappyapi-conebuilder}]{\sphinxcrossref{\DUrole{std,std-ref}{ConeBuilder Class}}}} object, starting with 0.
\end{quote}
\end{quote}

\sphinxAtStartPar
\sphinxstylestrong{Example}
\end{quote}

\begin{sphinxVerbatim}[commandchars=\\\{\}]
\PYG{c+c1}{\PYGZsh{} Retrieve the variable in SOC constraint builder conex with indice of 1}
\PYG{n}{conevar} \PYG{o}{=} \PYG{n}{conex}\PYG{o}{.}\PYG{n}{getVar}\PYG{p}{(}\PYG{l+m+mi}{1}\PYG{p}{)}
\end{sphinxVerbatim}

\subsubsection{ConeBuilder.getVars()}
\label{\detokenize{pyapiref:conebuilder-getvars}}\begin{quote}

\sphinxAtStartPar
\sphinxstylestrong{Synopsis}
\begin{quote}

\sphinxAtStartPar
\sphinxcode{\sphinxupquote{getVars()}}
\end{quote}

\sphinxAtStartPar
\sphinxstylestrong{Description}
\begin{quote}

\sphinxAtStartPar
Retrieve all variables in {\hyperref[\detokenize{pyapiref:chappyapi-conebuilder}]{\sphinxcrossref{\DUrole{std,std-ref}{ConeBuilder Class}}}} objects, and return a {\hyperref[\detokenize{pyapiref:chappyapi-vararray}]{\sphinxcrossref{\DUrole{std,std-ref}{VarArray Class}}}} object.
\end{quote}

\sphinxAtStartPar
\sphinxstylestrong{Example}
\end{quote}

\begin{sphinxVerbatim}[commandchars=\\\{\}]
\PYG{c+c1}{\PYGZsh{} Retrieve all variables in SOC constraint builder conex.}
\PYG{n}{conevars} \PYG{o}{=} \PYG{n}{conex}\PYG{o}{.}\PYG{n}{getVars}\PYG{p}{(}\PYG{p}{)}
\end{sphinxVerbatim}

\subsubsection{ConeBuilder.getSize()}
\label{\detokenize{pyapiref:conebuilder-getsize}}\begin{quote}

\sphinxAtStartPar
\sphinxstylestrong{Synopsis}
\begin{quote}

\sphinxAtStartPar
\sphinxcode{\sphinxupquote{getSize()}}
\end{quote}

\sphinxAtStartPar
\sphinxstylestrong{Description}
\begin{quote}

\sphinxAtStartPar
Retrieve the number of elements in {\hyperref[\detokenize{pyapiref:chappyapi-conebuilder}]{\sphinxcrossref{\DUrole{std,std-ref}{ConeBuilder Class}}}} object.
\end{quote}

\sphinxAtStartPar
\sphinxstylestrong{Example}
\end{quote}

\begin{sphinxVerbatim}[commandchars=\\\{\}]
\PYG{c+c1}{\PYGZsh{} Retrieve the number of elements in SOC constraint builder conex.}
\PYG{n}{conesize} \PYG{o}{=} \PYG{n}{conex}\PYG{o}{.}\PYG{n}{getSize}\PYG{p}{(}\PYG{p}{)}
\end{sphinxVerbatim}

\subsection{ConeBuilderArray Class}
\label{\detokenize{pyapiref:conebuilderarray-class}}\label{\detokenize{pyapiref:chappyapi-conebuilderarray}}
\sphinxAtStartPar
In order to facilitate users to operate on a set of {\hyperref[\detokenize{pyapiref:chappyapi-conebuilder}]{\sphinxcrossref{\DUrole{std,std-ref}{ConeBuilder Class}}}} objects, COPT provides ConeBuilderArray class
in Python interface, providing the following methods:

\subsubsection{ConeBuilderArray()}
\label{\detokenize{pyapiref:conebuilderarray}}\begin{quote}

\sphinxAtStartPar
\sphinxstylestrong{Synopsis}
\begin{quote}

\sphinxAtStartPar
\sphinxcode{\sphinxupquote{ConeBuilderArray(conebuilders=None)}}
\end{quote}

\sphinxAtStartPar
\sphinxstylestrong{Description}
\begin{quote}

\sphinxAtStartPar
Create a {\hyperref[\detokenize{pyapiref:chappyapi-conebuilderarray}]{\sphinxcrossref{\DUrole{std,std-ref}{ConeBuilderArray Class}}}} object.

\sphinxAtStartPar
If parameter \sphinxcode{\sphinxupquote{conebuilders}} is \sphinxcode{\sphinxupquote{None}}, then create an empty {\hyperref[\detokenize{pyapiref:chappyapi-conebuilderarray}]{\sphinxcrossref{\DUrole{std,std-ref}{ConeBuilderArray Class}}}} object,
otherwise initialize the newly created {\hyperref[\detokenize{pyapiref:chappyapi-conebuilderarray}]{\sphinxcrossref{\DUrole{std,std-ref}{ConeBuilderArray Class}}}} object with parameter \sphinxcode{\sphinxupquote{conebuilders}}.
\end{quote}

\sphinxAtStartPar
\sphinxstylestrong{Arguments}
\begin{quote}

\sphinxAtStartPar
\sphinxcode{\sphinxupquote{conebuilders}}
\begin{quote}

\sphinxAtStartPar
SOC constraint builder to be added. Optional, \sphinxcode{\sphinxupquote{None}} by default.
Could be {\hyperref[\detokenize{pyapiref:chappyapi-conebuilder}]{\sphinxcrossref{\DUrole{std,std-ref}{ConeBuilder Class}}}} object, {\hyperref[\detokenize{pyapiref:chappyapi-conebuilderarray}]{\sphinxcrossref{\DUrole{std,std-ref}{ConeBuilderArray Class}}}} object, list, dictionary
or {\hyperref[\detokenize{pyapiref:chappyapi-util-tupledict}]{\sphinxcrossref{\DUrole{std,std-ref}{tupledict Class}}}} object.
\end{quote}
\end{quote}

\sphinxAtStartPar
\sphinxstylestrong{Example}
\end{quote}

\begin{sphinxVerbatim}[commandchars=\\\{\}]
\PYG{c+c1}{\PYGZsh{} Create an empty ConeBuilderArray object.}
\PYG{n}{conebuilderarr} \PYG{o}{=} \PYG{n}{ConeBuilderArray}\PYG{p}{(}\PYG{p}{)}
\PYG{c+c1}{\PYGZsh{} Create a ConeBuilderArray object and initialize it with SOC constraint builder conex and coney}
\PYG{n}{conebuilderarr} \PYG{o}{=} \PYG{n}{ConeBuilderArray}\PYG{p}{(}\PYG{p}{[}\PYG{n}{conex}\PYG{p}{,} \PYG{n}{coney}\PYG{p}{]}\PYG{p}{)}
\end{sphinxVerbatim}

\subsubsection{ConeBuilderArray.pushBack()}
\label{\detokenize{pyapiref:conebuilderarray-pushback}}\begin{quote}

\sphinxAtStartPar
\sphinxstylestrong{Synopsis}
\begin{quote}

\sphinxAtStartPar
\sphinxcode{\sphinxupquote{pushBack(conebuilder)}}
\end{quote}

\sphinxAtStartPar
\sphinxstylestrong{Description}
\begin{quote}

\sphinxAtStartPar
Add one or multiple {\hyperref[\detokenize{pyapiref:chappyapi-conebuilder}]{\sphinxcrossref{\DUrole{std,std-ref}{ConeBuilder Class}}}} objects.
\end{quote}

\sphinxAtStartPar
\sphinxstylestrong{Arguments}
\begin{quote}

\sphinxAtStartPar
\sphinxcode{\sphinxupquote{conebuilder}}
\begin{quote}

\sphinxAtStartPar
SOC constraint builder to be added. Could be {\hyperref[\detokenize{pyapiref:chappyapi-conebuilder}]{\sphinxcrossref{\DUrole{std,std-ref}{ConeBuilder Class}}}} object, {\hyperref[\detokenize{pyapiref:chappyapi-conebuilderarray}]{\sphinxcrossref{\DUrole{std,std-ref}{ConeBuilderArray Class}}}} object,
list, dictionary or {\hyperref[\detokenize{pyapiref:chappyapi-util-tupledict}]{\sphinxcrossref{\DUrole{std,std-ref}{tupledict Class}}}} object.
\end{quote}
\end{quote}

\sphinxAtStartPar
\sphinxstylestrong{Example}
\end{quote}

\begin{sphinxVerbatim}[commandchars=\\\{\}]
\PYG{c+c1}{\PYGZsh{} Add SOC constraint builder conex to conebuilderarr}
\PYG{n}{conebuilderarr}\PYG{o}{.}\PYG{n}{pushBack}\PYG{p}{(}\PYG{n}{conex}\PYG{p}{)}
\end{sphinxVerbatim}

\subsubsection{ConeBuilderArray.getBuilder()}
\label{\detokenize{pyapiref:conebuilderarray-getbuilder}}\begin{quote}

\sphinxAtStartPar
\sphinxstylestrong{Synopsis}
\begin{quote}

\sphinxAtStartPar
\sphinxcode{\sphinxupquote{getBuilder(idx)}}
\end{quote}

\sphinxAtStartPar
\sphinxstylestrong{Description}
\begin{quote}

\sphinxAtStartPar
Retrieve the corresponding builder according to the indice of SOC constraint builder in {\hyperref[\detokenize{pyapiref:chappyapi-conebuilderarray}]{\sphinxcrossref{\DUrole{std,std-ref}{ConeBuilderArray Class}}}} object.
\end{quote}

\sphinxAtStartPar
\sphinxstylestrong{Arguments}
\begin{quote}

\sphinxAtStartPar
\sphinxcode{\sphinxupquote{idx}}
\begin{quote}

\sphinxAtStartPar
Indice of the SOC constraint builder in {\hyperref[\detokenize{pyapiref:chappyapi-conebuilderarray}]{\sphinxcrossref{\DUrole{std,std-ref}{ConeBuilderArray Class}}}} object, starting with 0.
\end{quote}
\end{quote}

\sphinxAtStartPar
\sphinxstylestrong{Example}
\end{quote}

\begin{sphinxVerbatim}[commandchars=\\\{\}]
\PYG{c+c1}{\PYGZsh{} Retrieve the SOC constraint builder with indice of 1 in conebuilderarr}
\PYG{n}{ConeBuilder} \PYG{o}{=} \PYG{n}{conebuilderarr}\PYG{o}{.}\PYG{n}{getBuilder}\PYG{p}{(}\PYG{l+m+mi}{1}\PYG{p}{)}
\end{sphinxVerbatim}

\subsubsection{ConeBuilderArray.getSize()}
\label{\detokenize{pyapiref:conebuilderarray-getsize}}\begin{quote}

\sphinxAtStartPar
\sphinxstylestrong{Synopsis}
\begin{quote}

\sphinxAtStartPar
\sphinxcode{\sphinxupquote{getSize()}}
\end{quote}

\sphinxAtStartPar
\sphinxstylestrong{Description}
\begin{quote}

\sphinxAtStartPar
Retrieve the number of elements in {\hyperref[\detokenize{pyapiref:chappyapi-conebuilderarray}]{\sphinxcrossref{\DUrole{std,std-ref}{ConeBuilderArray Class}}}} object.
\end{quote}

\sphinxAtStartPar
\sphinxstylestrong{Example}
\end{quote}

\begin{sphinxVerbatim}[commandchars=\\\{\}]
\PYG{c+c1}{\PYGZsh{} Retrieve the number of elements in conebuilderarr}
\PYG{n}{conebuildersize} \PYG{o}{=} \PYG{n}{conebuilderarr}\PYG{o}{.}\PYG{n}{getSize}\PYG{p}{(}\PYG{p}{)}
\end{sphinxVerbatim}

\subsection{ExpCone Class}
\label{\detokenize{pyapiref:expcone-class}}\label{\detokenize{pyapiref:chappyapi-expcone}}
\sphinxAtStartPar
ExpCone object contains related operations of COPT exponential cone constraints.
The following methods are provided:

\subsubsection{ExpCone.getIdx()}
\label{\detokenize{pyapiref:expcone-getidx}}\begin{quote}

\sphinxAtStartPar
\sphinxstylestrong{Synopsis}
\begin{quote}

\sphinxAtStartPar
\sphinxcode{\sphinxupquote{getIdx()}}
\end{quote}

\sphinxAtStartPar
\sphinxstylestrong{Description}
\begin{quote}

\sphinxAtStartPar
Retrieve the subscript of exponential cone constraint in model.
\end{quote}

\sphinxAtStartPar
\sphinxstylestrong{Example}
\end{quote}

\begin{sphinxVerbatim}[commandchars=\\\{\}]
\PYG{c+c1}{\PYGZsh{} Retrieve the subscript of exponential cone constraint cone.}
\PYG{n}{coneidx} \PYG{o}{=} \PYG{n}{cone}\PYG{o}{.}\PYG{n}{getIdx}\PYG{p}{(}\PYG{p}{)}
\end{sphinxVerbatim}

\subsubsection{ExpCone.remove()}
\label{\detokenize{pyapiref:expcone-remove}}\begin{quote}

\sphinxAtStartPar
\sphinxstylestrong{Synopsis}
\begin{quote}

\sphinxAtStartPar
\sphinxcode{\sphinxupquote{remove()}}
\end{quote}

\sphinxAtStartPar
\sphinxstylestrong{Description}
\begin{quote}

\sphinxAtStartPar
Delete the exponential cone constraint from model.
\end{quote}

\sphinxAtStartPar
\sphinxstylestrong{Example}
\end{quote}

\begin{sphinxVerbatim}[commandchars=\\\{\}]
\PYG{c+c1}{\PYGZsh{} Delete the exponential cone constraint \PYGZsq{}cone\PYGZsq{}}
\PYG{n}{cone}\PYG{o}{.}\PYG{n}{remove}\PYG{p}{(}\PYG{p}{)}
\end{sphinxVerbatim}

\subsection{ExpConeArray Class}
\label{\detokenize{pyapiref:expconearray-class}}\label{\detokenize{pyapiref:chappyapi-expconearray}}
\sphinxAtStartPar
To facilitate users to operate on a set of {\hyperref[\detokenize{pyapiref:chappyapi-expcone}]{\sphinxcrossref{\DUrole{std,std-ref}{ExpCone Class}}}} objects, COPT designed ExpConeArray class
in Python interface. The following methods are provided:

\subsubsection{ExpConeArray()}
\label{\detokenize{pyapiref:expconearray}}\begin{quote}

\sphinxAtStartPar
\sphinxstylestrong{Synopsis}
\begin{quote}

\sphinxAtStartPar
\sphinxcode{\sphinxupquote{ExpConeArray(cones=None)}}
\end{quote}

\sphinxAtStartPar
\sphinxstylestrong{Description}
\begin{quote}

\sphinxAtStartPar
Create a {\hyperref[\detokenize{pyapiref:chappyapi-expconearray}]{\sphinxcrossref{\DUrole{std,std-ref}{ExpConeArray Class}}}} object.

\sphinxAtStartPar
If parameter \sphinxcode{\sphinxupquote{cones}} is \sphinxcode{\sphinxupquote{None}}, then build an empty {\hyperref[\detokenize{pyapiref:chappyapi-expconearray}]{\sphinxcrossref{\DUrole{std,std-ref}{ExpConeArray Class}}}} object,
otherwise initialize the newly created {\hyperref[\detokenize{pyapiref:chappyapi-expconearray}]{\sphinxcrossref{\DUrole{std,std-ref}{ExpConeArray Class}}}} object with \sphinxcode{\sphinxupquote{cones}}.
\end{quote}

\sphinxAtStartPar
\sphinxstylestrong{Arguments}
\begin{quote}

\sphinxAtStartPar
\sphinxcode{\sphinxupquote{cones}}
\begin{quote}

\sphinxAtStartPar
Exponential cone constraint to be added. Optional, \sphinxcode{\sphinxupquote{None}} by default.
It can be {\hyperref[\detokenize{pyapiref:chappyapi-expcone}]{\sphinxcrossref{\DUrole{std,std-ref}{ExpCone Class}}}} object, {\hyperref[\detokenize{pyapiref:chappyapi-expconearray}]{\sphinxcrossref{\DUrole{std,std-ref}{ExpConeArray Class}}}} object, list, dictionary or
{\hyperref[\detokenize{pyapiref:chappyapi-util-tupledict}]{\sphinxcrossref{\DUrole{std,std-ref}{tupledict Class}}}} object.
\end{quote}
\end{quote}

\sphinxAtStartPar
\sphinxstylestrong{Example}
\end{quote}

\begin{sphinxVerbatim}[commandchars=\\\{\}]
\PYG{c+c1}{\PYGZsh{} Create a new ExpConeArray object}
\PYG{n}{conearr} \PYG{o}{=} \PYG{n}{ExpConeArray}\PYG{p}{(}\PYG{p}{)}
\PYG{c+c1}{\PYGZsh{} Create a ExpConeArray object, and initialize it with exponential cone constraints conex and coney.}
\PYG{n}{conearr} \PYG{o}{=} \PYG{n}{ExpConeArray}\PYG{p}{(}\PYG{p}{[}\PYG{n}{conex}\PYG{p}{,} \PYG{n}{coney}\PYG{p}{]}\PYG{p}{)}
\end{sphinxVerbatim}

\subsubsection{ExpConeArray.pushBack()}
\label{\detokenize{pyapiref:expconearray-pushback}}\begin{quote}

\sphinxAtStartPar
\sphinxstylestrong{Synopsis}
\begin{quote}

\sphinxAtStartPar
\sphinxcode{\sphinxupquote{pushBack(cone)}}
\end{quote}

\sphinxAtStartPar
\sphinxstylestrong{Description}
\begin{quote}

\sphinxAtStartPar
Add one or multiple {\hyperref[\detokenize{pyapiref:chappyapi-expcone}]{\sphinxcrossref{\DUrole{std,std-ref}{ExpCone Class}}}} objects.
\end{quote}

\sphinxAtStartPar
\sphinxstylestrong{Arguments}
\begin{quote}

\sphinxAtStartPar
\sphinxcode{\sphinxupquote{cone}}
\begin{quote}

\sphinxAtStartPar
Exponential cone constraints to be added, which can be {\hyperref[\detokenize{pyapiref:chappyapi-expcone}]{\sphinxcrossref{\DUrole{std,std-ref}{ExpCone Class}}}} object,
{\hyperref[\detokenize{pyapiref:chappyapi-expconearray}]{\sphinxcrossref{\DUrole{std,std-ref}{ExpConeArray Class}}}} object, list, dictionary or {\hyperref[\detokenize{pyapiref:chappyapi-util-tupledict}]{\sphinxcrossref{\DUrole{std,std-ref}{tupledict Class}}}} object.
\end{quote}
\end{quote}

\sphinxAtStartPar
\sphinxstylestrong{Example}
\end{quote}

\begin{sphinxVerbatim}[commandchars=\\\{\}]
\PYG{c+c1}{\PYGZsh{} Add exponential cone constraint conex to conearr}
\PYG{n}{conearr}\PYG{o}{.}\PYG{n}{pushBack}\PYG{p}{(}\PYG{n}{conex}\PYG{p}{)}
\PYG{c+c1}{\PYGZsh{} Add exponential cone constraints conex and coney to conearr}
\PYG{n}{conearr}\PYG{o}{.}\PYG{n}{pushBack}\PYG{p}{(}\PYG{p}{[}\PYG{n}{conex}\PYG{p}{,} \PYG{n}{coney}\PYG{p}{]}\PYG{p}{)}
\end{sphinxVerbatim}

\subsubsection{ExpConeArray.getCone()}
\label{\detokenize{pyapiref:expconearray-getcone}}\begin{quote}

\sphinxAtStartPar
\sphinxstylestrong{Synopsis}
\begin{quote}

\sphinxAtStartPar
\sphinxcode{\sphinxupquote{getCone(idx)}}
\end{quote}

\sphinxAtStartPar
\sphinxstylestrong{Description}
\begin{quote}

\sphinxAtStartPar
Retrieve the corresponding exponential cone constraint according to its subscript in
{\hyperref[\detokenize{pyapiref:chappyapi-expconearray}]{\sphinxcrossref{\DUrole{std,std-ref}{ExpConeArray Class}}}} object and return a {\hyperref[\detokenize{pyapiref:chappyapi-expcone}]{\sphinxcrossref{\DUrole{std,std-ref}{ExpCone Class}}}} object.
\end{quote}

\sphinxAtStartPar
\sphinxstylestrong{Arguments}
\begin{quote}

\sphinxAtStartPar
\sphinxcode{\sphinxupquote{idx}}
\begin{quote}

\sphinxAtStartPar
Indice of the exponential cone constraint in {\hyperref[\detokenize{pyapiref:chappyapi-expconearray}]{\sphinxcrossref{\DUrole{std,std-ref}{ExpConeArray Class}}}} object, starting with 0.
\end{quote}
\end{quote}

\sphinxAtStartPar
\sphinxstylestrong{Example}
\end{quote}

\begin{sphinxVerbatim}[commandchars=\\\{\}]
\PYG{c+c1}{\PYGZsh{} Retrieve the exponential cone constraint with indice of 1 in conearr}
\PYG{n}{cone} \PYG{o}{=} \PYG{n}{conearr}\PYG{o}{.}\PYG{n}{getCone}\PYG{p}{(}\PYG{l+m+mi}{1}\PYG{p}{)}
\end{sphinxVerbatim}

\subsubsection{ExpConeArray.getSize()}
\label{\detokenize{pyapiref:expconearray-getsize}}\begin{quote}

\sphinxAtStartPar
\sphinxstylestrong{Synopsis}
\begin{quote}

\sphinxAtStartPar
\sphinxcode{\sphinxupquote{getSize()}}
\end{quote}

\sphinxAtStartPar
\sphinxstylestrong{Description}
\begin{quote}

\sphinxAtStartPar
Retrieve the number of elements in {\hyperref[\detokenize{pyapiref:chappyapi-expconearray}]{\sphinxcrossref{\DUrole{std,std-ref}{ExpConeArray Class}}}} object.
\end{quote}

\sphinxAtStartPar
\sphinxstylestrong{Example}
\end{quote}

\begin{sphinxVerbatim}[commandchars=\\\{\}]
\PYG{c+c1}{\PYGZsh{} Retrieve the number of exponential cone constraints in conearr.}
\PYG{n}{arrsize} \PYG{o}{=} \PYG{n}{conearr}\PYG{o}{.}\PYG{n}{getSize}\PYG{p}{(}\PYG{p}{)}
\end{sphinxVerbatim}

\subsection{ExpConeBuilder Class}
\label{\detokenize{pyapiref:expconebuilder-class}}\label{\detokenize{pyapiref:chappyapi-expconebuilder}}
\sphinxAtStartPar
For easy access builders of exponential cone constraints, ExpConeBuilder class provides the following methods:

\subsubsection{ExpConeBuilder()}
\label{\detokenize{pyapiref:expconebuilder}}\begin{quote}

\sphinxAtStartPar
\sphinxstylestrong{Synopsis}
\begin{quote}

\sphinxAtStartPar
\sphinxcode{\sphinxupquote{ExpConeBuilder()}}
\end{quote}

\sphinxAtStartPar
\sphinxstylestrong{Description}
\begin{quote}

\sphinxAtStartPar
Create an empty {\hyperref[\detokenize{pyapiref:chappyapi-expconebuilder}]{\sphinxcrossref{\DUrole{std,std-ref}{ExpConeBuilder Class}}}} object.
\end{quote}

\sphinxAtStartPar
\sphinxstylestrong{Example}
\end{quote}

\begin{sphinxVerbatim}[commandchars=\\\{\}]
\PYG{c+c1}{\PYGZsh{} Create an empty ExpConeBuilder object.}
\PYG{n}{ExpConeBuilder} \PYG{o}{=} \PYG{n}{ExpConeBuilder}\PYG{p}{(}\PYG{p}{)}
\end{sphinxVerbatim}

\subsubsection{ExpConeBuilder.setBuilder()}
\label{\detokenize{pyapiref:expconebuilder-setbuilder}}\begin{quote}

\sphinxAtStartPar
\sphinxstylestrong{Synopsis}
\begin{quote}

\sphinxAtStartPar
\sphinxcode{\sphinxupquote{setBuilder(conetype, vars)}}
\end{quote}

\sphinxAtStartPar
\sphinxstylestrong{Description}
\begin{quote}

\sphinxAtStartPar
Set type of variables of {\hyperref[\detokenize{pyapiref:chappyapi-expconebuilder}]{\sphinxcrossref{\DUrole{std,std-ref}{ExpConeBuilder Class}}}} object.
\end{quote}

\sphinxAtStartPar
\sphinxstylestrong{Arguments}
\begin{quote}

\sphinxAtStartPar
\sphinxcode{\sphinxupquote{conetype}}
\begin{quote}

\sphinxAtStartPar
Type of exponential cone constraint. Full list of available types can be found in
{\hyperref[\detokenize{constant:chapconst-expconetype}]{\sphinxcrossref{\DUrole{std,std-ref}{exponential cone types}}}}.
\end{quote}

\sphinxAtStartPar
\sphinxcode{\sphinxupquote{vars}}
\begin{quote}

\sphinxAtStartPar
Variables of exponential cone constarint, which can be {\hyperref[\detokenize{pyapiref:chappyapi-vararray}]{\sphinxcrossref{\DUrole{std,std-ref}{VarArray Class}}}} object, list, dictionary or
{\hyperref[\detokenize{pyapiref:chappyapi-util-tupledict}]{\sphinxcrossref{\DUrole{std,std-ref}{tupledict Class}}}} object.
\end{quote}
\end{quote}

\sphinxAtStartPar
\sphinxstylestrong{Example}
\end{quote}

\begin{sphinxVerbatim}[commandchars=\\\{\}]
\PYG{c+c1}{\PYGZsh{} Set type as regular, variables as [z, x, y] for exponential cone constraint builder}
\PYG{n}{ExpConeBuilder}\PYG{o}{.}\PYG{n}{setBuilder}\PYG{p}{(}\PYG{n}{COPT}\PYG{o}{.}\PYG{n}{EXPCONE\PYGZus{}PRIMAL}\PYG{p}{,} \PYG{p}{[}\PYG{n}{z}\PYG{p}{,} \PYG{n}{x}\PYG{p}{,} \PYG{n}{y}\PYG{p}{]}\PYG{p}{)}
\end{sphinxVerbatim}

\subsubsection{ExpConeBuilder.getType()}
\label{\detokenize{pyapiref:expconebuilder-gettype}}\begin{quote}

\sphinxAtStartPar
\sphinxstylestrong{Synopsis}
\begin{quote}

\sphinxAtStartPar
\sphinxcode{\sphinxupquote{getType()}}
\end{quote}

\sphinxAtStartPar
\sphinxstylestrong{Description}
\begin{quote}

\sphinxAtStartPar
Retrieve the exponential cone constraint type of {\hyperref[\detokenize{pyapiref:chappyapi-expconebuilder}]{\sphinxcrossref{\DUrole{std,std-ref}{ExpConeBuilder Class}}}} object.
\end{quote}

\sphinxAtStartPar
\sphinxstylestrong{Example}
\end{quote}

\begin{sphinxVerbatim}[commandchars=\\\{\}]
\PYG{c+c1}{\PYGZsh{} Retrieve the type of exponential cone constraint builder conex.}
\PYG{n}{conetype} \PYG{o}{=} \PYG{n}{ExpConeBuilder}\PYG{o}{.}\PYG{n}{getType}\PYG{p}{(}\PYG{n}{conex}\PYG{p}{)}
\end{sphinxVerbatim}

\subsubsection{ExpConeBuilder.getVar()}
\label{\detokenize{pyapiref:expconebuilder-getvar}}\begin{quote}

\sphinxAtStartPar
\sphinxstylestrong{Synopsis}
\begin{quote}

\sphinxAtStartPar
\sphinxcode{\sphinxupquote{getVar(idx)}}
\end{quote}

\sphinxAtStartPar
\sphinxstylestrong{Description}
\begin{quote}

\sphinxAtStartPar
Retrieve the corresponding variables according to its indice in {\hyperref[\detokenize{pyapiref:chappyapi-expconebuilder}]{\sphinxcrossref{\DUrole{std,std-ref}{ExpConeBuilder Class}}}} object,
and return a {\hyperref[\detokenize{pyapiref:chappyapi-var}]{\sphinxcrossref{\DUrole{std,std-ref}{Var Class}}}} object.
\end{quote}

\sphinxAtStartPar
\sphinxstylestrong{Arguments}
\begin{quote}

\sphinxAtStartPar
\sphinxcode{\sphinxupquote{idx}}
\begin{quote}

\sphinxAtStartPar
Indice of the variable in {\hyperref[\detokenize{pyapiref:chappyapi-expconebuilder}]{\sphinxcrossref{\DUrole{std,std-ref}{ExpConeBuilder Class}}}} object, starting with 0.
\end{quote}
\end{quote}

\sphinxAtStartPar
\sphinxstylestrong{Example}
\end{quote}

\begin{sphinxVerbatim}[commandchars=\\\{\}]
\PYG{c+c1}{\PYGZsh{} Retrieve the variable in exponential cone constraint builder conex with indice of 1}
\PYG{n}{conevar} \PYG{o}{=} \PYG{n}{conex}\PYG{o}{.}\PYG{n}{getVar}\PYG{p}{(}\PYG{l+m+mi}{1}\PYG{p}{)}
\end{sphinxVerbatim}

\subsubsection{ExpConeBuilder.getVars()}
\label{\detokenize{pyapiref:expconebuilder-getvars}}\begin{quote}

\sphinxAtStartPar
\sphinxstylestrong{Synopsis}
\begin{quote}

\sphinxAtStartPar
\sphinxcode{\sphinxupquote{getVars()}}
\end{quote}

\sphinxAtStartPar
\sphinxstylestrong{Description}
\begin{quote}

\sphinxAtStartPar
Retrieve all variables in {\hyperref[\detokenize{pyapiref:chappyapi-expconebuilder}]{\sphinxcrossref{\DUrole{std,std-ref}{ExpConeBuilder Class}}}} objects, and return a {\hyperref[\detokenize{pyapiref:chappyapi-vararray}]{\sphinxcrossref{\DUrole{std,std-ref}{VarArray Class}}}} object.
\end{quote}

\sphinxAtStartPar
\sphinxstylestrong{Example}
\end{quote}

\begin{sphinxVerbatim}[commandchars=\\\{\}]
\PYG{c+c1}{\PYGZsh{} Retrieve all variables in exponential cone constraint builder conex.}
\PYG{n}{conevars} \PYG{o}{=} \PYG{n}{conex}\PYG{o}{.}\PYG{n}{getVars}\PYG{p}{(}\PYG{p}{)}
\end{sphinxVerbatim}

\subsubsection{ExpConeBuilder.getSize()}
\label{\detokenize{pyapiref:expconebuilder-getsize}}\begin{quote}

\sphinxAtStartPar
\sphinxstylestrong{Synopsis}
\begin{quote}

\sphinxAtStartPar
\sphinxcode{\sphinxupquote{getSize()}}
\end{quote}

\sphinxAtStartPar
\sphinxstylestrong{Description}
\begin{quote}

\sphinxAtStartPar
Retrieve the number of elements in {\hyperref[\detokenize{pyapiref:chappyapi-expconebuilder}]{\sphinxcrossref{\DUrole{std,std-ref}{ExpConeBuilder Class}}}} object.
\end{quote}

\sphinxAtStartPar
\sphinxstylestrong{Example}
\end{quote}

\begin{sphinxVerbatim}[commandchars=\\\{\}]
\PYG{c+c1}{\PYGZsh{} Retrieve the number of elements in exponential cone constraint builder conex.}
\PYG{n}{conesize} \PYG{o}{=} \PYG{n}{conex}\PYG{o}{.}\PYG{n}{getSize}\PYG{p}{(}\PYG{p}{)}
\end{sphinxVerbatim}

\subsection{ExpConeBuilderArray Class}
\label{\detokenize{pyapiref:expconebuilderarray-class}}\label{\detokenize{pyapiref:chappyapi-expconebuilderarray}}
\sphinxAtStartPar
In order to facilitate users to operate on a set of {\hyperref[\detokenize{pyapiref:chappyapi-expconebuilder}]{\sphinxcrossref{\DUrole{std,std-ref}{ExpConeBuilder Class}}}} objects, COPT provides ExpConeBuilderArray class
in Python interface, providing the following methods:

\subsubsection{ExpConeBuilderArray()}
\label{\detokenize{pyapiref:expconebuilderarray}}\begin{quote}

\sphinxAtStartPar
\sphinxstylestrong{Synopsis}
\begin{quote}

\sphinxAtStartPar
\sphinxcode{\sphinxupquote{ExpConeBuilderArray(ExpConeBuilders=None)}}
\end{quote}

\sphinxAtStartPar
\sphinxstylestrong{Description}
\begin{quote}

\sphinxAtStartPar
Create a {\hyperref[\detokenize{pyapiref:chappyapi-expconebuilderarray}]{\sphinxcrossref{\DUrole{std,std-ref}{ExpConeBuilderArray Class}}}} object.

\sphinxAtStartPar
If parameter \sphinxcode{\sphinxupquote{ExpConeBuilders}} is \sphinxcode{\sphinxupquote{None}}, then create an empty {\hyperref[\detokenize{pyapiref:chappyapi-expconebuilderarray}]{\sphinxcrossref{\DUrole{std,std-ref}{ExpConeBuilderArray Class}}}} object,
otherwise initialize the newly created {\hyperref[\detokenize{pyapiref:chappyapi-expconebuilderarray}]{\sphinxcrossref{\DUrole{std,std-ref}{ExpConeBuilderArray Class}}}} object with parameter \sphinxcode{\sphinxupquote{ExpConeBuilders}}.
\end{quote}

\sphinxAtStartPar
\sphinxstylestrong{Arguments}
\begin{quote}

\sphinxAtStartPar
\sphinxcode{\sphinxupquote{ExpConeBuilders}}
\begin{quote}

\sphinxAtStartPar
Exponential cone constraint builder to be added. Optional, \sphinxcode{\sphinxupquote{None}} by default.
Could be {\hyperref[\detokenize{pyapiref:chappyapi-expconebuilder}]{\sphinxcrossref{\DUrole{std,std-ref}{ExpConeBuilder Class}}}} object, {\hyperref[\detokenize{pyapiref:chappyapi-expconebuilderarray}]{\sphinxcrossref{\DUrole{std,std-ref}{ExpConeBuilderArray Class}}}} object, list, dictionary
or {\hyperref[\detokenize{pyapiref:chappyapi-util-tupledict}]{\sphinxcrossref{\DUrole{std,std-ref}{tupledict Class}}}} object.
\end{quote}
\end{quote}

\sphinxAtStartPar
\sphinxstylestrong{Example}
\end{quote}

\begin{sphinxVerbatim}[commandchars=\\\{\}]
\PYG{c+c1}{\PYGZsh{} Create an empty ExpConeBuilderArray object.}
\PYG{n}{ExpConeBuilderarr} \PYG{o}{=} \PYG{n}{ExpConeBuilderArray}\PYG{p}{(}\PYG{p}{)}
\PYG{c+c1}{\PYGZsh{} Create a ExpConeBuilderArray object and initialize it with exponential cone constraint builder conex and coney}
\PYG{n}{ExpConeBuilderarr} \PYG{o}{=} \PYG{n}{ExpConeBuilderArray}\PYG{p}{(}\PYG{p}{[}\PYG{n}{conex}\PYG{p}{,} \PYG{n}{coney}\PYG{p}{]}\PYG{p}{)}
\end{sphinxVerbatim}

\subsubsection{ExpConeBuilderArray.pushBack()}
\label{\detokenize{pyapiref:expconebuilderarray-pushback}}\begin{quote}

\sphinxAtStartPar
\sphinxstylestrong{Synopsis}
\begin{quote}

\sphinxAtStartPar
\sphinxcode{\sphinxupquote{pushBack(conebuilder)}}
\end{quote}

\sphinxAtStartPar
\sphinxstylestrong{Description}
\begin{quote}

\sphinxAtStartPar
Add one or multiple {\hyperref[\detokenize{pyapiref:chappyapi-expconebuilder}]{\sphinxcrossref{\DUrole{std,std-ref}{ExpConeBuilder Class}}}} objects.
\end{quote}

\sphinxAtStartPar
\sphinxstylestrong{Arguments}
\begin{quote}

\sphinxAtStartPar
\sphinxcode{\sphinxupquote{conebuilder}}
\begin{quote}

\sphinxAtStartPar
Exponential cone constraint builder to be added. Could be {\hyperref[\detokenize{pyapiref:chappyapi-expconebuilder}]{\sphinxcrossref{\DUrole{std,std-ref}{ExpConeBuilder Class}}}} object, {\hyperref[\detokenize{pyapiref:chappyapi-expconebuilderarray}]{\sphinxcrossref{\DUrole{std,std-ref}{ExpConeBuilderArray Class}}}} object,
list, dictionary or {\hyperref[\detokenize{pyapiref:chappyapi-util-tupledict}]{\sphinxcrossref{\DUrole{std,std-ref}{tupledict Class}}}} object.
\end{quote}
\end{quote}

\sphinxAtStartPar
\sphinxstylestrong{Example}
\end{quote}

\begin{sphinxVerbatim}[commandchars=\\\{\}]
\PYG{c+c1}{\PYGZsh{} Add exponential cone constraint builder conex to ExpConeBuilderarr}
\PYG{n}{ExpConeBuilderarr}\PYG{o}{.}\PYG{n}{pushBack}\PYG{p}{(}\PYG{n}{conex}\PYG{p}{)}
\end{sphinxVerbatim}

\subsubsection{ExpConeBuilderArray.getBuilder()}
\label{\detokenize{pyapiref:expconebuilderarray-getbuilder}}\begin{quote}

\sphinxAtStartPar
\sphinxstylestrong{Synopsis}
\begin{quote}

\sphinxAtStartPar
\sphinxcode{\sphinxupquote{getBuilder(idx)}}
\end{quote}

\sphinxAtStartPar
\sphinxstylestrong{Description}
\begin{quote}

\sphinxAtStartPar
Retrieve the corresponding builder according to the indice of exponential cone constraint builder in {\hyperref[\detokenize{pyapiref:chappyapi-expconebuilderarray}]{\sphinxcrossref{\DUrole{std,std-ref}{ExpConeBuilderArray Class}}}} object.
\end{quote}

\sphinxAtStartPar
\sphinxstylestrong{Arguments}
\begin{quote}

\sphinxAtStartPar
\sphinxcode{\sphinxupquote{idx}}
\begin{quote}

\sphinxAtStartPar
Indice of the exponential cone constraint builder in {\hyperref[\detokenize{pyapiref:chappyapi-expconebuilderarray}]{\sphinxcrossref{\DUrole{std,std-ref}{ExpConeBuilderArray Class}}}} object, starting with 0.
\end{quote}
\end{quote}

\sphinxAtStartPar
\sphinxstylestrong{Example}
\end{quote}

\begin{sphinxVerbatim}[commandchars=\\\{\}]
\PYG{c+c1}{\PYGZsh{} Retrieve the exponential cone constraint builder with indice of 1 in ExpConeBuilderarr}
\PYG{n}{ExpExpConeBuilder} \PYG{o}{=} \PYG{n}{ExpConeBuilderarr}\PYG{o}{.}\PYG{n}{getBuilder}\PYG{p}{(}\PYG{l+m+mi}{1}\PYG{p}{)}
\end{sphinxVerbatim}

\subsubsection{ExpConeBuilderArray.getSize()}
\label{\detokenize{pyapiref:expconebuilderarray-getsize}}\begin{quote}

\sphinxAtStartPar
\sphinxstylestrong{Synopsis}
\begin{quote}

\sphinxAtStartPar
\sphinxcode{\sphinxupquote{getSize()}}
\end{quote}

\sphinxAtStartPar
\sphinxstylestrong{Description}
\begin{quote}

\sphinxAtStartPar
Retrieve the number of elements in {\hyperref[\detokenize{pyapiref:chappyapi-expconebuilderarray}]{\sphinxcrossref{\DUrole{std,std-ref}{ExpConeBuilderArray Class}}}} object.
\end{quote}

\sphinxAtStartPar
\sphinxstylestrong{Example}
\end{quote}

\begin{sphinxVerbatim}[commandchars=\\\{\}]
\PYG{c+c1}{\PYGZsh{} Retrieve the number of elements in ExpConeBuilderarr}
\PYG{n}{ExpConeBuildersize} \PYG{o}{=} \PYG{n}{ExpConeBuilderarr}\PYG{o}{.}\PYG{n}{getSize}\PYG{p}{(}\PYG{p}{)}
\end{sphinxVerbatim}

\subsection{AffineCone Class}
\label{\detokenize{pyapiref:affinecone-class}}\label{\detokenize{pyapiref:chappyapi-affinecone}}
\sphinxAtStartPar
The AffineCone class encapsulates operations related to affine cone in COPT.
The following methods are provided:

\subsubsection{AffineCone.getIdx()}
\label{\detokenize{pyapiref:affinecone-getidx}}\begin{quote}

\sphinxAtStartPar
\sphinxstylestrong{Synopsis}
\begin{quote}

\sphinxAtStartPar
\sphinxcode{\sphinxupquote{getIdx()}}
\end{quote}

\sphinxAtStartPar
\sphinxstylestrong{Description}
\begin{quote}

\sphinxAtStartPar
Retrieve the index of the affine cone in the model.
\end{quote}

\sphinxAtStartPar
\sphinxstylestrong{Example}
\end{quote}

\begin{sphinxVerbatim}[commandchars=\\\{\}]
\PYG{c+c1}{\PYGZsh{} Retrieve the index of the afcone}
\PYG{n}{afconeidx} \PYG{o}{=} \PYG{n}{afcone}\PYG{o}{.}\PYG{n}{getIdx}\PYG{p}{(}\PYG{p}{)}
\end{sphinxVerbatim}

\subsubsection{AffineCone.getName()}
\label{\detokenize{pyapiref:affinecone-getname}}\begin{quote}

\sphinxAtStartPar
\sphinxstylestrong{Synopsis}
\begin{quote}

\sphinxAtStartPar
\sphinxcode{\sphinxupquote{getName()}}
\end{quote}

\sphinxAtStartPar
\sphinxstylestrong{Description}
\begin{quote}

\sphinxAtStartPar
Retrieve the name of the affine cone.
\end{quote}

\sphinxAtStartPar
\sphinxstylestrong{Example}
\end{quote}

\begin{sphinxVerbatim}[commandchars=\\\{\}]
\PYG{c+c1}{\PYGZsh{} Retrieve the name of the afcone}
\PYG{n}{afconename} \PYG{o}{=} \PYG{n}{afcone}\PYG{o}{.}\PYG{n}{getName}\PYG{p}{(}\PYG{p}{)}
\end{sphinxVerbatim}

\subsubsection{AffineCone.setName()}
\label{\detokenize{pyapiref:affinecone-setname}}\begin{quote}

\sphinxAtStartPar
\sphinxstylestrong{Synopsis}
\begin{quote}

\sphinxAtStartPar
\sphinxcode{\sphinxupquote{setName(newname)}}
\end{quote}

\sphinxAtStartPar
\sphinxstylestrong{Description}
\begin{quote}

\sphinxAtStartPar
Set the name of the affine cone.
\end{quote}

\sphinxAtStartPar
\sphinxstylestrong{Arguments}
\begin{quote}

\sphinxAtStartPar
\sphinxcode{\sphinxupquote{newname}}
\begin{quote}

\sphinxAtStartPar
The name to set for the affine cone.
\end{quote}
\end{quote}

\sphinxAtStartPar
\sphinxstylestrong{Example}
\end{quote}

\begin{sphinxVerbatim}[commandchars=\\\{\}]
\PYG{c+c1}{\PYGZsh{} Set the name of the afcone}
\PYG{n}{afcone}\PYG{o}{.}\PYG{n}{setName}\PYG{p}{(}\PYG{l+s+s2}{\PYGZdq{}}\PYG{l+s+s2}{afcone}\PYG{l+s+s2}{\PYGZdq{}}\PYG{p}{)}
\end{sphinxVerbatim}

\subsubsection{AffineCone.remove()}
\label{\detokenize{pyapiref:affinecone-remove}}\begin{quote}

\sphinxAtStartPar
\sphinxstylestrong{Synopsis}
\begin{quote}

\sphinxAtStartPar
\sphinxcode{\sphinxupquote{remove()}}
\end{quote}

\sphinxAtStartPar
\sphinxstylestrong{Description}
\begin{quote}

\sphinxAtStartPar
Remove the current affine cone from the model.
\end{quote}

\sphinxAtStartPar
\sphinxstylestrong{Example}
\end{quote}

\begin{sphinxVerbatim}[commandchars=\\\{\}]
\PYG{c+c1}{\PYGZsh{} Remove the current afcone from the model}
\PYG{n}{afcone}\PYG{o}{.}\PYG{n}{remove}\PYG{p}{(}\PYG{p}{)}
\end{sphinxVerbatim}

\subsection{AffineConeArray Class}
\label{\detokenize{pyapiref:affineconearray-class}}\label{\detokenize{pyapiref:chappyapi-affineconearray}}
\sphinxAtStartPar
To facilitate user operations on a group of {\hyperref[\detokenize{pyapiref:chappyapi-affinecone}]{\sphinxcrossref{\DUrole{std,std-ref}{AffineCone Class}}}} objects,
COPT introduces the AffineConeArray class. The following methods are provided:

\subsubsection{AffineConeArray()}
\label{\detokenize{pyapiref:affineconearray}}\begin{quote}

\sphinxAtStartPar
\sphinxstylestrong{Synopsis}
\begin{quote}

\sphinxAtStartPar
\sphinxcode{\sphinxupquote{AffineConeArray(cones=None)}}
\end{quote}

\sphinxAtStartPar
\sphinxstylestrong{Description}
\begin{quote}

\sphinxAtStartPar
Creates an {\hyperref[\detokenize{pyapiref:chappyapi-affineconearray}]{\sphinxcrossref{\DUrole{std,std-ref}{AffineConeArray Class}}}} object.

\sphinxAtStartPar
If the argument \sphinxcode{\sphinxupquote{cones}} is \sphinxcode{\sphinxupquote{None}}, an empty {\hyperref[\detokenize{pyapiref:chappyapi-affineconearray}]{\sphinxcrossref{\DUrole{std,std-ref}{AffineConeArray Class}}}} object is created.
Otherwise, the new {\hyperref[\detokenize{pyapiref:chappyapi-affineconearray}]{\sphinxcrossref{\DUrole{std,std-ref}{AffineConeArray Class}}}} object is initialized with the argument \sphinxcode{\sphinxupquote{cones}}.
\end{quote}

\sphinxAtStartPar
\sphinxstylestrong{Arguments}
\begin{quote}

\sphinxAtStartPar
\sphinxcode{\sphinxupquote{cones}}
\begin{quote}

\sphinxAtStartPar
Affine cone constraints to be added.
This is an optional parameter, the default value is \sphinxcode{\sphinxupquote{None}}.

\sphinxAtStartPar
Acceptable values include {\hyperref[\detokenize{pyapiref:chappyapi-affinecone}]{\sphinxcrossref{\DUrole{std,std-ref}{AffineCone Class}}}} objects,
{\hyperref[\detokenize{pyapiref:chappyapi-affineconearray}]{\sphinxcrossref{\DUrole{std,std-ref}{AffineConeArray Class}}}} objects, lists, dictionaries, or {\hyperref[\detokenize{pyapiref:chappyapi-util-tupledict}]{\sphinxcrossref{\DUrole{std,std-ref}{tupledict Class}}}} objects.
\end{quote}
\end{quote}

\sphinxAtStartPar
\sphinxstylestrong{Example}
\end{quote}

\begin{sphinxVerbatim}[commandchars=\\\{\}]
\PYG{c+c1}{\PYGZsh{} Create an empty AffineConeArray object}
\PYG{n}{afconearr} \PYG{o}{=} \PYG{n}{AffineConeArray}\PYG{p}{(}\PYG{p}{)}
\PYG{c+c1}{\PYGZsh{} Create an AffineConeArray object initialized with afconex and afconey}
\PYG{n}{afconearr} \PYG{o}{=} \PYG{n}{AffineConeArray}\PYG{p}{(}\PYG{p}{[}\PYG{n}{afconex}\PYG{p}{,} \PYG{n}{afconey}\PYG{p}{]}\PYG{p}{)}
\end{sphinxVerbatim}

\subsubsection{AffineConeArray.pushBack()}
\label{\detokenize{pyapiref:affineconearray-pushback}}\begin{quote}

\sphinxAtStartPar
\sphinxstylestrong{Synopsis}
\begin{quote}

\sphinxAtStartPar
\sphinxcode{\sphinxupquote{pushBack(cone)}}
\end{quote}

\sphinxAtStartPar
\sphinxstylestrong{Description}
\begin{quote}

\sphinxAtStartPar
Adds one or more {\hyperref[\detokenize{pyapiref:chappyapi-affinecone}]{\sphinxcrossref{\DUrole{std,std-ref}{AffineCone Class}}}} objects.
\end{quote}

\sphinxAtStartPar
\sphinxstylestrong{Arguments}
\begin{quote}

\sphinxAtStartPar
\sphinxcode{\sphinxupquote{cone}}
\begin{quote}

\sphinxAtStartPar
Affine cone constraints to be added. Possible values include
{\hyperref[\detokenize{pyapiref:chappyapi-affinecone}]{\sphinxcrossref{\DUrole{std,std-ref}{AffineCone Class}}}} objects, {\hyperref[\detokenize{pyapiref:chappyapi-affineconearray}]{\sphinxcrossref{\DUrole{std,std-ref}{AffineConeArray Class}}}} objects,
lists, dictionaries, or {\hyperref[\detokenize{pyapiref:chappyapi-util-tupledict}]{\sphinxcrossref{\DUrole{std,std-ref}{tupledict Class}}}} objects.
\end{quote}
\end{quote}

\sphinxAtStartPar
\sphinxstylestrong{Example}
\end{quote}

\begin{sphinxVerbatim}[commandchars=\\\{\}]
\PYG{c+c1}{\PYGZsh{} Add afconex to afconearr}
\PYG{n}{afconearr}\PYG{o}{.}\PYG{n}{pushBack}\PYG{p}{(}\PYG{n}{afconex}\PYG{p}{)}
\PYG{c+c1}{\PYGZsh{} Add afconex and afconey to afconearr}
\PYG{n}{afconearr}\PYG{o}{.}\PYG{n}{pushBack}\PYG{p}{(}\PYG{p}{[}\PYG{n}{afconex}\PYG{p}{,} \PYG{n}{afconey}\PYG{p}{]}\PYG{p}{)}

\PYG{n}{dict\PYGZus{}cone} \PYG{o}{=} \PYG{p}{\PYGZob{}}\PYG{l+s+s1}{\PYGZsq{}}\PYG{l+s+s1}{AC}\PYG{l+s+s1}{\PYGZsq{}}\PYG{p}{:} \PYG{n}{afcone1}\PYG{p}{,} \PYG{l+s+s1}{\PYGZsq{}}\PYG{l+s+s1}{CA}\PYG{l+s+s1}{\PYGZsq{}}\PYG{p}{:} \PYG{n}{afcone2}\PYG{p}{\PYGZcb{}}
\PYG{n}{afconearr2} \PYG{o}{=} \PYG{n}{AffineConeArray}\PYG{p}{(}\PYG{p}{\PYGZob{}}\PYG{l+s+s1}{\PYGZsq{}}\PYG{l+s+s1}{AC}\PYG{l+s+s1}{\PYGZsq{}}\PYG{p}{:} \PYG{n}{afcone1}\PYG{p}{,} \PYG{l+s+s1}{\PYGZsq{}}\PYG{l+s+s1}{CA}\PYG{l+s+s1}{\PYGZsq{}}\PYG{p}{:} \PYG{n}{afcone2}\PYG{p}{\PYGZcb{}}\PYG{p}{)}
\end{sphinxVerbatim}

\subsubsection{AffineConeArray.getCone()}
\label{\detokenize{pyapiref:affineconearray-getcone}}\begin{quote}

\sphinxAtStartPar
\sphinxstylestrong{Synopsis}
\begin{quote}

\sphinxAtStartPar
\sphinxcode{\sphinxupquote{getCone(idx)}}
\end{quote}

\sphinxAtStartPar
\sphinxstylestrong{Description}
\begin{quote}

\sphinxAtStartPar
Retrieves the affine cone at the specified index in the {\hyperref[\detokenize{pyapiref:chappyapi-affineconearray}]{\sphinxcrossref{\DUrole{std,std-ref}{AffineConeArray Class}}}} object.
Returns an {\hyperref[\detokenize{pyapiref:chappyapi-affinecone}]{\sphinxcrossref{\DUrole{std,std-ref}{AffineCone Class}}}} object.
\end{quote}

\sphinxAtStartPar
\sphinxstylestrong{Arguments}
\begin{quote}

\sphinxAtStartPar
\sphinxcode{\sphinxupquote{idx}}
\begin{quote}

\sphinxAtStartPar
Index of the affine cone constraint in the {\hyperref[\detokenize{pyapiref:chappyapi-affineconearray}]{\sphinxcrossref{\DUrole{std,std-ref}{AffineConeArray Class}}}} object. Indexing starts at 0.
\end{quote}
\end{quote}

\sphinxAtStartPar
\sphinxstylestrong{Example}
\end{quote}

\begin{sphinxVerbatim}[commandchars=\\\{\}]
\PYG{c+c1}{\PYGZsh{} Retrieve the affine cone at index 1 in afconearr}
\PYG{n}{afcone} \PYG{o}{=} \PYG{n}{afconearr}\PYG{o}{.}\PYG{n}{getCone}\PYG{p}{(}\PYG{l+m+mi}{1}\PYG{p}{)}
\end{sphinxVerbatim}

\subsubsection{AffineConeArray.getSize()}
\label{\detokenize{pyapiref:affineconearray-getsize}}\begin{quote}

\sphinxAtStartPar
\sphinxstylestrong{Synopsis}
\begin{quote}

\sphinxAtStartPar
\sphinxcode{\sphinxupquote{getSize()}}
\end{quote}

\sphinxAtStartPar
\sphinxstylestrong{Description}
\begin{quote}

\sphinxAtStartPar
Retrieves the number of elements in the {\hyperref[\detokenize{pyapiref:chappyapi-affineconearray}]{\sphinxcrossref{\DUrole{std,std-ref}{AffineConeArray Class}}}} object.
\end{quote}

\sphinxAtStartPar
\sphinxstylestrong{Example}
\end{quote}

\begin{sphinxVerbatim}[commandchars=\\\{\}]
\PYG{c+c1}{\PYGZsh{} Get the number of affine cone in afconearr}
\PYG{n}{arrsize} \PYG{o}{=} \PYG{n}{afconearr}\PYG{o}{.}\PYG{n}{getSize}\PYG{p}{(}\PYG{p}{)}
\end{sphinxVerbatim}

\subsection{AffineConeBuilder Class}
\label{\detokenize{pyapiref:affineconebuilder-class}}\label{\detokenize{pyapiref:chappyapi-affineconebuilder}}
\sphinxAtStartPar
The \sphinxtitleref{AffineConeBuilder} class is an encapsulation of the builder for constructing
affine cone in COPT.  The following methods are provided:

\subsubsection{AffineConeBuilder()}
\label{\detokenize{pyapiref:affineconebuilder}}\begin{quote}

\sphinxAtStartPar
\sphinxstylestrong{Synopsis}
\begin{quote}

\sphinxAtStartPar
\sphinxcode{\sphinxupquote{AffineConeBuilder()}}
\end{quote}

\sphinxAtStartPar
\sphinxstylestrong{Description}
\begin{quote}

\sphinxAtStartPar
Creates an empty {\hyperref[\detokenize{pyapiref:chappyapi-affineconebuilder}]{\sphinxcrossref{\DUrole{std,std-ref}{AffineConeBuilder Class}}}} object.
\end{quote}

\sphinxAtStartPar
\sphinxstylestrong{Example}
\end{quote}

\begin{sphinxVerbatim}[commandchars=\\\{\}]
\PYG{c+c1}{\PYGZsh{} Create an empty AffineConeBuilder object}
\PYG{n}{afconebuilder} \PYG{o}{=} \PYG{n}{AffineConeBuilder}\PYG{p}{(}\PYG{p}{)}
\end{sphinxVerbatim}

\subsubsection{AffineConeBuilder.setBuilder()}
\label{\detokenize{pyapiref:affineconebuilder-setbuilder}}\begin{quote}

\sphinxAtStartPar
\sphinxstylestrong{Synopsis}
\begin{quote}

\sphinxAtStartPar
\sphinxcode{\sphinxupquote{setBuilder(conetype, exprs)}}
\end{quote}

\sphinxAtStartPar
\sphinxstylestrong{Description}
\begin{quote}

\sphinxAtStartPar
Sets up the affine cone builder according to \sphinxcode{\sphinxupquote{conetype}} and \sphinxcode{\sphinxupquote{exprs}} .
\end{quote}

\sphinxAtStartPar
\sphinxstylestrong{Arguments}
\begin{quote}

\sphinxAtStartPar
\sphinxcode{\sphinxupquote{conetype}}
\begin{quote}

\sphinxAtStartPar
The type of affine cone. For possible values, refer to
{\hyperref[\detokenize{constant:chapconst-conetype}]{\sphinxcrossref{\DUrole{std,std-ref}{Second\sphinxhyphen{}Order Cone Constraint Types}}}} and
{\hyperref[\detokenize{constant:chapconst-expconetype}]{\sphinxcrossref{\DUrole{std,std-ref}{Exponential Cone Constraint Types}}}}.
\end{quote}

\sphinxAtStartPar
\sphinxcode{\sphinxupquote{exprs}}
\begin{quote}

\sphinxAtStartPar
Mathematical expressions forming the affine cone. Possible values include
linear expressions or semidefinite expressions.
\end{quote}
\end{quote}

\sphinxAtStartPar
\sphinxstylestrong{Example}
\end{quote}

\begin{sphinxVerbatim}[commandchars=\\\{\}]
\PYG{c+c1}{\PYGZsh{} Set the affine cone builder\PYGZsq{}s type to standard second\PYGZhy{}order cone,}
\PYG{c+c1}{\PYGZsh{} composed of z, x+y, y}
\PYG{n}{afconebuilder} \PYG{o}{=} \PYG{n}{AffineConeBuilder}\PYG{p}{(}\PYG{p}{)}
\PYG{n}{afconebuilder}\PYG{o}{.}\PYG{n}{setBuilder}\PYG{p}{(}\PYG{n}{COPT}\PYG{o}{.}\PYG{n}{CONE\PYGZus{}QUAD}\PYG{p}{,} \PYG{p}{[}\PYG{n}{z}\PYG{p}{,} \PYG{n}{x}\PYG{o}{+}\PYG{n}{y}\PYG{p}{,} \PYG{n}{y}\PYG{p}{]}\PYG{p}{)}
\end{sphinxVerbatim}

\subsubsection{AffineConeBuilder.hasPsdTerm()}
\label{\detokenize{pyapiref:affineconebuilder-haspsdterm}}\begin{quote}

\sphinxAtStartPar
\sphinxstylestrong{Synopsis}
\begin{quote}

\sphinxAtStartPar
\sphinxcode{\sphinxupquote{hasPsdTerm()}}
\end{quote}

\sphinxAtStartPar
\sphinxstylestrong{Description}
\begin{quote}

\sphinxAtStartPar
Checks if the affine cone builder contains any semidefinite terms.
\end{quote}

\sphinxAtStartPar
\sphinxstylestrong{Return Value}
\begin{quote}

\sphinxAtStartPar
Bool. If \sphinxcode{\sphinxupquote{False}}, the affine cone contains no semidefinite terms, and linear expressions can be directly retrieved.
\end{quote}
\end{quote}

\subsubsection{AffineConeBuilder.getType()}
\label{\detokenize{pyapiref:affineconebuilder-gettype}}\begin{quote}

\sphinxAtStartPar
\sphinxstylestrong{Synopsis}
\begin{quote}

\sphinxAtStartPar
\sphinxcode{\sphinxupquote{getType()}}
\end{quote}

\sphinxAtStartPar
\sphinxstylestrong{Description}
\begin{quote}

\sphinxAtStartPar
Retrieves the affine cone type of the {\hyperref[\detokenize{pyapiref:chappyapi-affineconebuilder}]{\sphinxcrossref{\DUrole{std,std-ref}{AffineConeBuilder Class}}}} object.
\end{quote}

\sphinxAtStartPar
\sphinxstylestrong{Example}
\end{quote}

\begin{sphinxVerbatim}[commandchars=\\\{\}]
\PYG{c+c1}{\PYGZsh{} Retrieve the type of afconex}
\PYG{n}{afconetype} \PYG{o}{=} \PYG{n}{afconebuilder}\PYG{o}{.}\PYG{n}{getType}\PYG{p}{(}\PYG{n}{afconex}\PYG{p}{)}
\end{sphinxVerbatim}

\subsubsection{AffineConeBuilder.getExpr()}
\label{\detokenize{pyapiref:affineconebuilder-getexpr}}\begin{quote}

\sphinxAtStartPar
\sphinxstylestrong{Synopsis}
\begin{quote}

\sphinxAtStartPar
\sphinxcode{\sphinxupquote{getExpr(idx)}}
\end{quote}

\sphinxAtStartPar
\sphinxstylestrong{Description}
\begin{quote}

\sphinxAtStartPar
Retrieves the expression at the specified index in the affine cone builder.
\end{quote}

\sphinxAtStartPar
\sphinxstylestrong{Arguments}
\begin{quote}

\sphinxAtStartPar
\sphinxcode{\sphinxupquote{idx}}
\begin{quote}

\sphinxAtStartPar
The index of the expression in the affine cone builder, starting with \sphinxcode{\sphinxupquote{0}} .
\end{quote}
\end{quote}

\sphinxAtStartPar
\sphinxstylestrong{Return Value}
\begin{quote}

\sphinxAtStartPar
If semidefinite terms exist, returns a \sphinxtitleref{PsdExpr} object; otherwise, returns a \sphinxtitleref{LinExpr} object.
\end{quote}

\sphinxAtStartPar
\sphinxstylestrong{Example}
\end{quote}

\begin{sphinxVerbatim}[commandchars=\\\{\}]
\PYG{c+c1}{\PYGZsh{} Retrieve the expression at index 0 in the affine cone builder}
\PYG{n}{expr} \PYG{o}{=} \PYG{n}{afconstrbuilder}\PYG{o}{.}\PYG{n}{getExpr}\PYG{p}{(}\PYG{l+m+mi}{0}\PYG{p}{)}
\end{sphinxVerbatim}

\subsubsection{AffineConeBuilder.getPsdExpr()}
\label{\detokenize{pyapiref:affineconebuilder-getpsdexpr}}\begin{quote}

\sphinxAtStartPar
\sphinxstylestrong{Synopsis}
\begin{quote}

\sphinxAtStartPar
\sphinxcode{\sphinxupquote{getPsdExpr(idx)}}
\end{quote}

\sphinxAtStartPar
\sphinxstylestrong{Description}
\begin{quote}

\sphinxAtStartPar
Retrieves the semidefinite expression at the specified index in the affine cone builder.
\end{quote}

\sphinxAtStartPar
\sphinxstylestrong{Arguments}
\begin{quote}

\sphinxAtStartPar
\sphinxcode{\sphinxupquote{idx}}
\begin{quote}

\sphinxAtStartPar
The index of the semidefinite expression in the affine cone builder,
, starting with \sphinxcode{\sphinxupquote{0}} .
\end{quote}
\end{quote}

\sphinxAtStartPar
\sphinxstylestrong{Return Value}
\begin{quote}

\sphinxAtStartPar
Returns a \sphinxtitleref{PsdExpr} object.
\end{quote}

\sphinxAtStartPar
\sphinxstylestrong{Example}
\end{quote}

\begin{sphinxVerbatim}[commandchars=\\\{\}]
\PYG{c+c1}{\PYGZsh{} Retrieve the semidefinite expression at index 0 in the afconstrbuilder}
\PYG{n}{psdexpr} \PYG{o}{=} \PYG{n}{afconstrbuilder}\PYG{o}{.}\PYG{n}{getPsdExpr}\PYG{p}{(}\PYG{l+m+mi}{0}\PYG{p}{)}
\end{sphinxVerbatim}

\subsubsection{AffineConeBuilder.getExprs()}
\label{\detokenize{pyapiref:affineconebuilder-getexprs}}\begin{quote}

\sphinxAtStartPar
\sphinxstylestrong{Synopsis}
\begin{quote}

\sphinxAtStartPar
\sphinxcode{\sphinxupquote{getExprs(idx=None)}}
\end{quote}

\sphinxAtStartPar
\sphinxstylestrong{Description}
\begin{quote}

\sphinxAtStartPar
Retrieves a group of expressions from the affine cone builder at the specified indices.
\end{quote}

\sphinxAtStartPar
\sphinxstylestrong{Arguments}
\begin{quote}

\sphinxAtStartPar
\sphinxcode{\sphinxupquote{idx}}
\begin{quote}

\sphinxAtStartPar
The indices of the expressions in the affine cone builder.

\sphinxAtStartPar
If \sphinxcode{\sphinxupquote{idx}} is \sphinxcode{\sphinxupquote{None}}, all expressions in the affine cone are returned.
If \sphinxcode{\sphinxupquote{idx}} is a Python list, the expressions at the specified indices are returned.
\end{quote}
\end{quote}

\sphinxAtStartPar
\sphinxstylestrong{Return Value}
\begin{quote}

\sphinxAtStartPar
If semidefinite terms exist, returns \sphinxtitleref{PsdExpr} objects; otherwise, returns \sphinxtitleref{LinExpr} objects.
\end{quote}

\sphinxAtStartPar
\sphinxstylestrong{Example}
\end{quote}

\begin{sphinxVerbatim}[commandchars=\\\{\}]
\PYG{c+c1}{\PYGZsh{} Retrieve all expressions in the affine cone builder}
\PYG{n}{allexprs} \PYG{o}{=} \PYG{n}{afconstrbuilder}\PYG{o}{.}\PYG{n}{getExprs}\PYG{p}{(}\PYG{p}{)}
\end{sphinxVerbatim}

\subsubsection{AffineConeBuilder.getPsdExprs()}
\label{\detokenize{pyapiref:affineconebuilder-getpsdexprs}}\begin{quote}

\sphinxAtStartPar
\sphinxstylestrong{Synopsis}
\begin{quote}

\sphinxAtStartPar
\sphinxcode{\sphinxupquote{getPsdExprs(idx=None)}}
\end{quote}

\sphinxAtStartPar
\sphinxstylestrong{Description}
\begin{quote}

\sphinxAtStartPar
Retrieves a group of semidefinite expressions from the affine cone builder at the specified indices.
\end{quote}

\sphinxAtStartPar
\sphinxstylestrong{Arguments}
\begin{quote}

\sphinxAtStartPar
\sphinxcode{\sphinxupquote{idx}}
\begin{quote}

\sphinxAtStartPar
The indices of the semidefinite expressions in the affine cone builder.

\sphinxAtStartPar
If \sphinxcode{\sphinxupquote{idx}} is \sphinxcode{\sphinxupquote{None}}, all semidefinite expressions in the affine cone are returned.
If \sphinxcode{\sphinxupquote{idx}} is a Python list, the semidefinite expressions at the specified indices are returned.
\end{quote}
\end{quote}

\sphinxAtStartPar
\sphinxstylestrong{Return Value}
\begin{quote}

\sphinxAtStartPar
Returns \sphinxtitleref{PsdExpr} objects.
\end{quote}

\sphinxAtStartPar
\sphinxstylestrong{Example}
\end{quote}

\begin{sphinxVerbatim}[commandchars=\\\{\}]
\PYG{c+c1}{\PYGZsh{} Retrieve all semidefinite expressions in the afconstrbuilder}
\PYG{n}{allpsdexprs} \PYG{o}{=} \PYG{n}{afconstrbuilder}\PYG{o}{.}\PYG{n}{getPsdExprs}\PYG{p}{(}\PYG{p}{)}
\end{sphinxVerbatim}

\subsubsection{AffineConeBuilder.getSize()}
\label{\detokenize{pyapiref:affineconebuilder-getsize}}\begin{quote}

\sphinxAtStartPar
\sphinxstylestrong{Synopsis}
\begin{quote}

\sphinxAtStartPar
\sphinxcode{\sphinxupquote{getSize()}}
\end{quote}

\sphinxAtStartPar
\sphinxstylestrong{Description}
\begin{quote}

\sphinxAtStartPar
Retrieves the number of expressions in the affine cone builder.
\end{quote}

\sphinxAtStartPar
\sphinxstylestrong{Example}
\end{quote}

\begin{sphinxVerbatim}[commandchars=\\\{\}]
\PYG{c+c1}{\PYGZsh{} Retrieve the number of expressions in the afconebuilder}
\PYG{n}{afconesize} \PYG{o}{=} \PYG{n}{afconebuilder}\PYG{o}{.}\PYG{n}{getSize}\PYG{p}{(}\PYG{p}{)}
\end{sphinxVerbatim}

\subsection{AffineConeBuilderArray Class}
\label{\detokenize{pyapiref:affineconebuilderarray-class}}\label{\detokenize{pyapiref:chappyapi-affineconebuilderarray}}
\sphinxAtStartPar
To facilitate user operations on a group of {\hyperref[\detokenize{pyapiref:chappyapi-affineconebuilder}]{\sphinxcrossref{\DUrole{std,std-ref}{AffineConeBuilder Class}}}} objects,
COPT introduces the AffineConeBuilderArray class, providing the following methods:

\subsubsection{AffineConeBuilderArray()}
\label{\detokenize{pyapiref:affineconebuilderarray}}\begin{quote}

\sphinxAtStartPar
\sphinxstylestrong{Synopsis}
\begin{quote}

\sphinxAtStartPar
\sphinxcode{\sphinxupquote{AffineConeBuilderArray(conebuilders=None)}}
\end{quote}

\sphinxAtStartPar
\sphinxstylestrong{Description}
\begin{quote}

\sphinxAtStartPar
Creates an {\hyperref[\detokenize{pyapiref:chappyapi-affineconebuilderarray}]{\sphinxcrossref{\DUrole{std,std-ref}{AffineConeBuilderArray Class}}}} object.

\sphinxAtStartPar
If the argument \sphinxcode{\sphinxupquote{conebuilders}} is \sphinxcode{\sphinxupquote{None}}, an empty {\hyperref[\detokenize{pyapiref:chappyapi-affineconebuilderarray}]{\sphinxcrossref{\DUrole{std,std-ref}{AffineConeBuilderArray Class}}}} object is created.
Otherwise, the new {\hyperref[\detokenize{pyapiref:chappyapi-affineconebuilderarray}]{\sphinxcrossref{\DUrole{std,std-ref}{AffineConeBuilderArray Class}}}} object is initialized with the argument \sphinxcode{\sphinxupquote{conebuilders}}.
\end{quote}

\sphinxAtStartPar
\sphinxstylestrong{Arguments}
\begin{quote}

\sphinxAtStartPar
\sphinxcode{\sphinxupquote{conebuilders}}
\begin{quote}

\sphinxAtStartPar
Affine cone constraint builders to add. This is an optional argument, defaulting to \sphinxcode{\sphinxupquote{None}}.
Acceptable values include {\hyperref[\detokenize{pyapiref:chappyapi-affineconebuilder}]{\sphinxcrossref{\DUrole{std,std-ref}{AffineConeBuilder Class}}}} objects,
{\hyperref[\detokenize{pyapiref:chappyapi-affineconebuilderarray}]{\sphinxcrossref{\DUrole{std,std-ref}{AffineConeBuilderArray Class}}}} objects, lists, dictionaries, or {\hyperref[\detokenize{pyapiref:chappyapi-util-tupledict}]{\sphinxcrossref{\DUrole{std,std-ref}{tupledict Class}}}} objects.
\end{quote}
\end{quote}

\sphinxAtStartPar
\sphinxstylestrong{Example}
\end{quote}

\begin{sphinxVerbatim}[commandchars=\\\{\}]
\PYG{c+c1}{\PYGZsh{} Create an empty AffineConeBuilderArray object}
\PYG{n}{conebuilderarr} \PYG{o}{=} \PYG{n}{AffineConeBuilderArray}\PYG{p}{(}\PYG{p}{)}
\PYG{c+c1}{\PYGZsh{} Create an AffineConeBuilderArray object initialized with conex and coney}
\PYG{n}{conebuilderarr} \PYG{o}{=} \PYG{n}{AffineConeBuilderArray}\PYG{p}{(}\PYG{p}{[}\PYG{n}{conex}\PYG{p}{,} \PYG{n}{coney}\PYG{p}{]}\PYG{p}{)}
\end{sphinxVerbatim}

\subsubsection{AffineConeBuilderArray.pushBack()}
\label{\detokenize{pyapiref:affineconebuilderarray-pushback}}\begin{quote}

\sphinxAtStartPar
\sphinxstylestrong{Synopsis}
\begin{quote}

\sphinxAtStartPar
\sphinxcode{\sphinxupquote{pushBack(AffineConeBuilder)}}
\end{quote}

\sphinxAtStartPar
\sphinxstylestrong{Description}
\begin{quote}

\sphinxAtStartPar
Adds one or more {\hyperref[\detokenize{pyapiref:chappyapi-affineconebuilder}]{\sphinxcrossref{\DUrole{std,std-ref}{AffineConeBuilder Class}}}} objects.
\end{quote}

\sphinxAtStartPar
\sphinxstylestrong{Arguments}
\begin{quote}

\sphinxAtStartPar
\sphinxcode{\sphinxupquote{AffineConeBuilder}}
\begin{quote}

\sphinxAtStartPar
Affine cone builders to add. Acceptable values include
{\hyperref[\detokenize{pyapiref:chappyapi-affineconebuilder}]{\sphinxcrossref{\DUrole{std,std-ref}{AffineConeBuilder Class}}}} objects,
{\hyperref[\detokenize{pyapiref:chappyapi-affineconebuilderarray}]{\sphinxcrossref{\DUrole{std,std-ref}{AffineConeBuilderArray Class}}}} objects, lists, dictionaries, or
{\hyperref[\detokenize{pyapiref:chappyapi-util-tupledict}]{\sphinxcrossref{\DUrole{std,std-ref}{tupledict Class}}}} objects.
\end{quote}
\end{quote}

\sphinxAtStartPar
\sphinxstylestrong{Example}
\end{quote}

\begin{sphinxVerbatim}[commandchars=\\\{\}]
\PYG{c+c1}{\PYGZsh{} Add affine cone builder afconex to afconebuilderarr}
\PYG{n}{afconebuilderarr}\PYG{o}{.}\PYG{n}{pushBack}\PYG{p}{(}\PYG{n}{afconex}\PYG{p}{)}
\end{sphinxVerbatim}

\subsubsection{AffineConeBuilderArray.getBuilder()}
\label{\detokenize{pyapiref:affineconebuilderarray-getbuilder}}\begin{quote}

\sphinxAtStartPar
\sphinxstylestrong{Synopsis}
\begin{quote}

\sphinxAtStartPar
\sphinxcode{\sphinxupquote{getBuilder(idx)}}
\end{quote}

\sphinxAtStartPar
\sphinxstylestrong{Description}
\begin{quote}

\sphinxAtStartPar
Retrieves the affine cone builder at the specified index in the {\hyperref[\detokenize{pyapiref:chappyapi-affineconebuilderarray}]{\sphinxcrossref{\DUrole{std,std-ref}{AffineConeBuilderArray Class}}}} object.
\end{quote}

\sphinxAtStartPar
\sphinxstylestrong{Arguments}
\begin{quote}

\sphinxAtStartPar
\sphinxcode{\sphinxupquote{idx}}
\begin{quote}

\sphinxAtStartPar
The index of the affine cone builder in the {\hyperref[\detokenize{pyapiref:chappyapi-affineconebuilderarray}]{\sphinxcrossref{\DUrole{std,std-ref}{AffineConeBuilderArray Class}}}} object. Indexing starts at 0.
\end{quote}
\end{quote}

\sphinxAtStartPar
\sphinxstylestrong{Example}
\end{quote}

\begin{sphinxVerbatim}[commandchars=\\\{\}]
\PYG{c+c1}{\PYGZsh{} Retrieve the second\PYGZhy{}order affine cone builder at index 1 in afconebuilderarr}
\PYG{n}{affineConeBuilder} \PYG{o}{=} \PYG{n}{afconebuilderarr}\PYG{o}{.}\PYG{n}{getBuilder}\PYG{p}{(}\PYG{l+m+mi}{1}\PYG{p}{)}
\end{sphinxVerbatim}

\subsubsection{AffineConeBuilderArray.getSize()}
\label{\detokenize{pyapiref:affineconebuilderarray-getsize}}\begin{quote}

\sphinxAtStartPar
\sphinxstylestrong{Synopsis}
\begin{quote}

\sphinxAtStartPar
\sphinxcode{\sphinxupquote{getSize()}}
\end{quote}

\sphinxAtStartPar
\sphinxstylestrong{Description}
\begin{quote}

\sphinxAtStartPar
Retrieves the number of elements in the {\hyperref[\detokenize{pyapiref:chappyapi-affineconebuilderarray}]{\sphinxcrossref{\DUrole{std,std-ref}{AffineConeBuilderArray Class}}}} object.
\end{quote}

\sphinxAtStartPar
\sphinxstylestrong{Example}
\end{quote}

\begin{sphinxVerbatim}[commandchars=\\\{\}]
\PYG{c+c1}{\PYGZsh{} Get the number of elements in afconebuilderarr}
\PYG{n}{afconebuildersize} \PYG{o}{=} \PYG{n}{afconebuilderarr}\PYG{o}{.}\PYG{n}{getSize}\PYG{p}{(}\PYG{p}{)}
\end{sphinxVerbatim}

\subsection{GenConstr Class}
\label{\detokenize{pyapiref:genconstr-class}}\label{\detokenize{pyapiref:chappyapi-genconstr}}
\sphinxAtStartPar
For easy access to indicator constraints, COPT provides GenConstr class which containing the following methods:

\subsubsection{GenConstr.getName()}
\label{\detokenize{pyapiref:genconstr-getname}}\begin{quote}

\sphinxAtStartPar
\sphinxstylestrong{Synopsis}
\begin{quote}

\sphinxAtStartPar
\sphinxcode{\sphinxupquote{getName()}}
\end{quote}

\sphinxAtStartPar
\sphinxstylestrong{Description}
\begin{quote}

\sphinxAtStartPar
Retrieve the name of the indicator constraint in model.
\end{quote}

\sphinxAtStartPar
\sphinxstylestrong{Example}
\end{quote}

\begin{sphinxVerbatim}[commandchars=\\\{\}]
\PYG{c+c1}{\PYGZsh{} Retrieve the name of indicator constraint indicx}
\PYG{n}{indiname} \PYG{o}{=} \PYG{n}{indicx}\PYG{o}{.}\PYG{n}{getName}\PYG{p}{(}\PYG{p}{)}
\end{sphinxVerbatim}

\subsubsection{GenConstr.setName()}
\label{\detokenize{pyapiref:genconstr-setname}}\begin{quote}

\sphinxAtStartPar
\sphinxstylestrong{Synopsis}
\begin{quote}

\sphinxAtStartPar
\sphinxcode{\sphinxupquote{setName(newname)}}
\end{quote}

\sphinxAtStartPar
\sphinxstylestrong{Description}
\begin{quote}

\sphinxAtStartPar
Set the name of the indicator constraint in model with \sphinxcode{\sphinxupquote{newname}} .
\end{quote}

\sphinxAtStartPar
\sphinxstylestrong{Example}
\end{quote}

\begin{sphinxVerbatim}[commandchars=\\\{\}]
\PYG{c+c1}{\PYGZsh{} Set the name of the indicator constraint indicx with \PYGZdq{}if\PYGZdq{}}
\PYG{n}{indicx}\PYG{o}{.}\PYG{n}{setName}\PYG{p}{(}\PYG{l+s+s2}{\PYGZdq{}}\PYG{l+s+s2}{if}\PYG{l+s+s2}{\PYGZdq{}}\PYG{p}{)}
\end{sphinxVerbatim}

\subsubsection{GenConstr.getIdx()}
\label{\detokenize{pyapiref:genconstr-getidx}}\begin{quote}

\sphinxAtStartPar
\sphinxstylestrong{Synopsis}
\begin{quote}

\sphinxAtStartPar
\sphinxcode{\sphinxupquote{getIdx()}}
\end{quote}

\sphinxAtStartPar
\sphinxstylestrong{Description}
\begin{quote}

\sphinxAtStartPar
Retrieve the subscript of indicator constraint in model.
\end{quote}

\sphinxAtStartPar
\sphinxstylestrong{Example}
\end{quote}

\begin{sphinxVerbatim}[commandchars=\\\{\}]
\PYG{c+c1}{\PYGZsh{} Retrieve the indice of indicator constraint indicx}
\PYG{n}{indidx} \PYG{o}{=} \PYG{n}{indicx}\PYG{o}{.}\PYG{n}{getIdx}\PYG{p}{(}\PYG{p}{)}
\end{sphinxVerbatim}

\subsubsection{GenConstr.remove()}
\label{\detokenize{pyapiref:genconstr-remove}}\begin{quote}

\sphinxAtStartPar
\sphinxstylestrong{Synopsis}
\begin{quote}

\sphinxAtStartPar
\sphinxcode{\sphinxupquote{remove()}}
\end{quote}

\sphinxAtStartPar
\sphinxstylestrong{Description}
\begin{quote}

\sphinxAtStartPar
Delete the indicator constraint from model.
\end{quote}

\sphinxAtStartPar
\sphinxstylestrong{Example}
\end{quote}

\begin{sphinxVerbatim}[commandchars=\\\{\}]
\PYG{c+c1}{\PYGZsh{} Delete indicator constraint \PYGZsq{}indx\PYGZsq{}}
\PYG{n}{indx}\PYG{o}{.}\PYG{n}{remove}\PYG{p}{(}\PYG{p}{)}
\end{sphinxVerbatim}

\subsection{GenConstrArray Class}
\label{\detokenize{pyapiref:genconstrarray-class}}\label{\detokenize{pyapiref:chappyapi-genconstrarray}}
\sphinxAtStartPar
In order to facilitate users to operate on a set of {\hyperref[\detokenize{pyapiref:chappyapi-genconstr}]{\sphinxcrossref{\DUrole{std,std-ref}{GenConstr Class}}}} objects, COPT provides GenConstrArray
class in Python interface, providing the following methods:

\subsubsection{GenConstrArray()}
\label{\detokenize{pyapiref:genconstrarray}}\begin{quote}

\sphinxAtStartPar
\sphinxstylestrong{Synopsis}
\begin{quote}

\sphinxAtStartPar
\sphinxcode{\sphinxupquote{GenConstrArray(genconstrs=None)}}
\end{quote}

\sphinxAtStartPar
\sphinxstylestrong{Description}
\begin{quote}

\sphinxAtStartPar
Create a {\hyperref[\detokenize{pyapiref:chappyapi-genconstrarray}]{\sphinxcrossref{\DUrole{std,std-ref}{GenConstrArray Class}}}} object.

\sphinxAtStartPar
If parameter \sphinxcode{\sphinxupquote{genconstrs}} is \sphinxcode{\sphinxupquote{None}}, then create an empty {\hyperref[\detokenize{pyapiref:chappyapi-genconstrarray}]{\sphinxcrossref{\DUrole{std,std-ref}{GenConstrArray Class}}}} object,
otherwise initialize the newly created {\hyperref[\detokenize{pyapiref:chappyapi-genconstrarray}]{\sphinxcrossref{\DUrole{std,std-ref}{GenConstrArray Class}}}} object with parameter \sphinxcode{\sphinxupquote{genconstrs}}.
\end{quote}

\sphinxAtStartPar
\sphinxstylestrong{Arguments}
\begin{quote}

\sphinxAtStartPar
\sphinxcode{\sphinxupquote{genconstrs}}
\begin{quote}

\sphinxAtStartPar
Indicator constraint to be added. Optional, \sphinxcode{\sphinxupquote{None}} by default.
Could be {\hyperref[\detokenize{pyapiref:chappyapi-genconstr}]{\sphinxcrossref{\DUrole{std,std-ref}{GenConstr Class}}}} object, {\hyperref[\detokenize{pyapiref:chappyapi-genconstrarray}]{\sphinxcrossref{\DUrole{std,std-ref}{GenConstrArray Class}}}} object, list, dictionary
or {\hyperref[\detokenize{pyapiref:chappyapi-util-tupledict}]{\sphinxcrossref{\DUrole{std,std-ref}{tupledict Class}}}} object.
\end{quote}
\end{quote}

\sphinxAtStartPar
\sphinxstylestrong{Example}
\end{quote}

\begin{sphinxVerbatim}[commandchars=\\\{\}]
\PYG{c+c1}{\PYGZsh{} Create a new GenConstrArray object}
\PYG{n}{genconstrarr} \PYG{o}{=} \PYG{n}{GenConstrArray}\PYG{p}{(}\PYG{p}{)}
\PYG{c+c1}{\PYGZsh{} Create a GenConstrArray object and user indicator constraints genx and geny to initialize it.}
\PYG{n}{genconstrarr} \PYG{o}{=} \PYG{n}{GenConstrArray}\PYG{p}{(}\PYG{p}{[}\PYG{n}{genx}\PYG{p}{,} \PYG{n}{geny}\PYG{p}{]}\PYG{p}{)}
\end{sphinxVerbatim}

\subsubsection{GenConstrArray.pushBack()}
\label{\detokenize{pyapiref:genconstrarray-pushback}}\begin{quote}

\sphinxAtStartPar
\sphinxstylestrong{Synopsis}
\begin{quote}

\sphinxAtStartPar
\sphinxcode{\sphinxupquote{pushBack(genconstr)}}
\end{quote}

\sphinxAtStartPar
\sphinxstylestrong{Description}
\begin{quote}

\sphinxAtStartPar
Add one or multiple {\hyperref[\detokenize{pyapiref:chappyapi-genconstr}]{\sphinxcrossref{\DUrole{std,std-ref}{GenConstr Class}}}} objects.
\end{quote}

\sphinxAtStartPar
\sphinxstylestrong{Arguments}
\begin{quote}
\begin{description}
\sphinxlineitem{\sphinxcode{\sphinxupquote{constrs}}}
\sphinxAtStartPar
The indicator constraint to be added.
Cound be {\hyperref[\detokenize{pyapiref:chappyapi-genconstr}]{\sphinxcrossref{\DUrole{std,std-ref}{GenConstr Class}}}} object, {\hyperref[\detokenize{pyapiref:chappyapi-genconstrarray}]{\sphinxcrossref{\DUrole{std,std-ref}{GenConstrArray Class}}}} object, list, dictionary
or {\hyperref[\detokenize{pyapiref:chappyapi-util-tupledict}]{\sphinxcrossref{\DUrole{std,std-ref}{tupledict Class}}}} object.

\end{description}
\end{quote}

\sphinxAtStartPar
\sphinxstylestrong{Example}
\end{quote}

\begin{sphinxVerbatim}[commandchars=\\\{\}]
\PYG{c+c1}{\PYGZsh{} Aff indicator constraint genx to genconarr}
\PYG{n}{genconarr}\PYG{o}{.}\PYG{n}{pushBack}\PYG{p}{(}\PYG{n}{genx}\PYG{p}{)}
\PYG{c+c1}{\PYGZsh{} Add indicator constraint genx and geny to genconarr}
\PYG{n}{genconarr}\PYG{o}{.}\PYG{n}{pushBack}\PYG{p}{(}\PYG{p}{[}\PYG{n}{genx}\PYG{p}{,} \PYG{n}{geny}\PYG{p}{]}\PYG{p}{)}
\end{sphinxVerbatim}

\subsubsection{GenConstrArray.getGenConstr()}
\label{\detokenize{pyapiref:genconstrarray-getgenconstr}}\begin{quote}

\sphinxAtStartPar
\sphinxstylestrong{Synopsis}
\begin{quote}

\sphinxAtStartPar
\sphinxcode{\sphinxupquote{getGenConstr(idx)}}
\end{quote}

\sphinxAtStartPar
\sphinxstylestrong{Description}
\begin{quote}

\sphinxAtStartPar
Retrieve the corresponding indicator constraint according to the indice of indicator constraint in {\hyperref[\detokenize{pyapiref:chappyapi-genconstrarray}]{\sphinxcrossref{\DUrole{std,std-ref}{GenConstrArray Class}}}}
object, and return a {\hyperref[\detokenize{pyapiref:chappyapi-genconstr}]{\sphinxcrossref{\DUrole{std,std-ref}{GenConstr Class}}}} object.
\end{quote}

\sphinxAtStartPar
\sphinxstylestrong{Arguments}
\begin{quote}

\sphinxAtStartPar
\sphinxcode{\sphinxupquote{idx}}
\begin{quote}

\sphinxAtStartPar
Indice of the indicator constraint in {\hyperref[\detokenize{pyapiref:chappyapi-genconstrarray}]{\sphinxcrossref{\DUrole{std,std-ref}{GenConstrArray Class}}}}, starting with 0.
\end{quote}
\end{quote}

\sphinxAtStartPar
\sphinxstylestrong{Example}
\end{quote}

\begin{sphinxVerbatim}[commandchars=\\\{\}]
\PYG{c+c1}{\PYGZsh{} Retrieve the indicator constraint with indice of 1 in genconarr}
\PYG{n}{genconstr} \PYG{o}{=} \PYG{n}{genconarr}\PYG{o}{.}\PYG{n}{getGenConstr}\PYG{p}{(}\PYG{l+m+mi}{1}\PYG{p}{)}
\end{sphinxVerbatim}

\subsubsection{GenConstrArray.getSize()}
\label{\detokenize{pyapiref:genconstrarray-getsize}}\begin{quote}

\sphinxAtStartPar
\sphinxstylestrong{Synopsis}
\begin{quote}

\sphinxAtStartPar
\sphinxcode{\sphinxupquote{getSize()}}
\end{quote}

\sphinxAtStartPar
\sphinxstylestrong{Description}
\begin{quote}

\sphinxAtStartPar
Retrieve the number of elements in {\hyperref[\detokenize{pyapiref:chappyapi-genconstrarray}]{\sphinxcrossref{\DUrole{std,std-ref}{GenConstrArray Class}}}} object.
\end{quote}

\sphinxAtStartPar
\sphinxstylestrong{Example}
\end{quote}

\begin{sphinxVerbatim}[commandchars=\\\{\}]
\PYG{c+c1}{\PYGZsh{} Retrieve the number of elements in genconarr}
\PYG{n}{genconsize} \PYG{o}{=} \PYG{n}{genconarr}\PYG{o}{.}\PYG{n}{getSize}\PYG{p}{(}\PYG{p}{)}
\end{sphinxVerbatim}

\subsection{GenConstrBuilder Class}
\label{\detokenize{pyapiref:genconstrbuilder-class}}\label{\detokenize{pyapiref:chappyapi-genconstrbuilder}}
\sphinxAtStartPar
GenConstrBuilder object contains operations for building indicator constraints, and provides the following methods:

\subsubsection{GenConstrBuilder()}
\label{\detokenize{pyapiref:genconstrbuilder}}\begin{quote}

\sphinxAtStartPar
\sphinxstylestrong{Synopsis}
\begin{quote}

\sphinxAtStartPar
\sphinxcode{\sphinxupquote{GenConstrBuilder()}}
\end{quote}

\sphinxAtStartPar
\sphinxstylestrong{Description}
\begin{quote}

\sphinxAtStartPar
Create an empty {\hyperref[\detokenize{pyapiref:chappyapi-genconstrbuilder}]{\sphinxcrossref{\DUrole{std,std-ref}{GenConstrBuilder Class}}}} object.
\end{quote}

\sphinxAtStartPar
\sphinxstylestrong{Example}
\end{quote}

\begin{sphinxVerbatim}[commandchars=\\\{\}]
\PYG{c+c1}{\PYGZsh{} Create an empty GenConstrBuilder object}
\PYG{n}{genconbuilder} \PYG{o}{=} \PYG{n}{GenConstrBuilder}\PYG{p}{(}\PYG{p}{)}
\end{sphinxVerbatim}

\subsubsection{GenConstrBuilder.setBuilder()}
\label{\detokenize{pyapiref:genconstrbuilder-setbuilder}}\begin{quote}

\sphinxAtStartPar
\sphinxstylestrong{Synopsis}
\begin{quote}

\sphinxAtStartPar
\sphinxcode{\sphinxupquote{setBuilder(var, val, expr, sense, type=COPT.INDICATOR\_IF)}}
\end{quote}

\sphinxAtStartPar
\sphinxstylestrong{Description}
\begin{quote}

\sphinxAtStartPar
Set indicator variable, the value of indicator variable, the expression/sense of constraint, and the type of indicator
constraint of the {\hyperref[\detokenize{pyapiref:chappyapi-genconstrbuilder}]{\sphinxcrossref{\DUrole{std,std-ref}{GenConstrBuilder Class}}}} object.
\end{quote}

\sphinxAtStartPar
\sphinxstylestrong{Arguments}
\begin{quote}

\sphinxAtStartPar
\sphinxcode{\sphinxupquote{var}}
\begin{quote}

\sphinxAtStartPar
Indicator variable.
\end{quote}

\sphinxAtStartPar
\sphinxcode{\sphinxupquote{val}}
\begin{quote}

\sphinxAtStartPar
Value of an indicator variable.
\end{quote}

\sphinxAtStartPar
\sphinxcode{\sphinxupquote{expr}}
\begin{quote}

\sphinxAtStartPar
Expression of linear constraint, which can be {\hyperref[\detokenize{pyapiref:chappyapi-var}]{\sphinxcrossref{\DUrole{std,std-ref}{Var Class}}}} object or {\hyperref[\detokenize{pyapiref:chappyapi-linexpr}]{\sphinxcrossref{\DUrole{std,std-ref}{LinExpr Class}}}} object.
\end{quote}

\sphinxAtStartPar
\sphinxcode{\sphinxupquote{sense}}
\begin{quote}

\sphinxAtStartPar
Sense for the linear constraint. Please refer to {\hyperref[\detokenize{constant:chapconst-constrtype}]{\sphinxcrossref{\DUrole{std,std-ref}{Constraint type}}}} for possible values.
\end{quote}

\sphinxAtStartPar
\sphinxcode{\sphinxupquote{type}}
\begin{quote}

\sphinxAtStartPar
Type of the indicator constraint. The default value is \sphinxcode{\sphinxupquote{COPT.INDICATOR\_IF}} (If\sphinxhyphen{}Then).
Please refer to {\hyperref[\detokenize{constant:chapconst-indicatortype}]{\sphinxcrossref{\DUrole{std,std-ref}{Indicator constraint type}}}} for possible values.
\end{quote}
\end{quote}

\sphinxAtStartPar
\sphinxstylestrong{Example}
\end{quote}

\begin{sphinxVerbatim}[commandchars=\\\{\}]
\PYG{c+c1}{\PYGZsh{} Set indicator variable of indicator constraint builder to x. When x is true, the linear constraint x + y == 1 holds}
\PYG{n}{genconbuilder}\PYG{o}{.}\PYG{n}{setBuilder}\PYG{p}{(}\PYG{n}{x}\PYG{p}{,} \PYG{k+kc}{True}\PYG{p}{,} \PYG{n}{x} \PYG{o}{+} \PYG{n}{y} \PYG{o}{\PYGZhy{}} \PYG{l+m+mi}{1}\PYG{p}{,} \PYG{n}{COPT}\PYG{o}{.}\PYG{n}{EQUAL}\PYG{p}{)}
\end{sphinxVerbatim}

\subsubsection{GenConstrBuilder.getBinVar()}
\label{\detokenize{pyapiref:genconstrbuilder-getbinvar}}\begin{quote}

\sphinxAtStartPar
\sphinxstylestrong{Synopsis}
\begin{quote}

\sphinxAtStartPar
\sphinxcode{\sphinxupquote{getBinVar()}}
\end{quote}

\sphinxAtStartPar
\sphinxstylestrong{Description}
\begin{quote}

\sphinxAtStartPar
Retrieve the indicator variable of a {\hyperref[\detokenize{pyapiref:chappyapi-genconstrbuilder}]{\sphinxcrossref{\DUrole{std,std-ref}{GenConstrBuilder Class}}}} object.
\end{quote}

\sphinxAtStartPar
\sphinxstylestrong{Example}
\end{quote}

\begin{sphinxVerbatim}[commandchars=\\\{\}]
\PYG{c+c1}{\PYGZsh{} Retrieve the indicator variable of indicator constraint builder genbuilderx}
\PYG{n}{indvar} \PYG{o}{=} \PYG{n}{genbuilderx}\PYG{o}{.}\PYG{n}{getBinVar}\PYG{p}{(}\PYG{p}{)}
\end{sphinxVerbatim}

\subsubsection{GenConstrBuilder.getBinVal()}
\label{\detokenize{pyapiref:genconstrbuilder-getbinval}}\begin{quote}

\sphinxAtStartPar
\sphinxstylestrong{Synopsis}
\begin{quote}

\sphinxAtStartPar
\sphinxcode{\sphinxupquote{getBinVal()}}
\end{quote}

\sphinxAtStartPar
\sphinxstylestrong{Description}
\begin{quote}

\sphinxAtStartPar
Retrieve the value of indicator variable of a {\hyperref[\detokenize{pyapiref:chappyapi-genconstrbuilder}]{\sphinxcrossref{\DUrole{std,std-ref}{GenConstrBuilder Class}}}} object.
\end{quote}

\sphinxAtStartPar
\sphinxstylestrong{Example}
\end{quote}

\begin{sphinxVerbatim}[commandchars=\\\{\}]
\PYG{c+c1}{\PYGZsh{} Retrieve the value when the indicator variable of indicator constraint builder genbuilderx is valid}
\PYG{n}{indval} \PYG{o}{=} \PYG{n}{genbuilderx}\PYG{o}{.}\PYG{n}{getBinVal}\PYG{p}{(}\PYG{p}{)}
\end{sphinxVerbatim}

\subsubsection{GenConstrBuilder.getExpr()}
\label{\detokenize{pyapiref:genconstrbuilder-getexpr}}\begin{quote}

\sphinxAtStartPar
\sphinxstylestrong{Synopsis}
\begin{quote}

\sphinxAtStartPar
\sphinxcode{\sphinxupquote{getExpr()}}
\end{quote}

\sphinxAtStartPar
\sphinxstylestrong{Description}
\begin{quote}

\sphinxAtStartPar
Retrieve the linear expression of a {\hyperref[\detokenize{pyapiref:chappyapi-genconstrbuilder}]{\sphinxcrossref{\DUrole{std,std-ref}{GenConstrBuilder Class}}}} object.
\end{quote}

\sphinxAtStartPar
\sphinxstylestrong{Example}
\end{quote}

\begin{sphinxVerbatim}[commandchars=\\\{\}]
\PYG{c+c1}{\PYGZsh{} Retrieve the linear expression of indicator constraint builder genbuilderx}
\PYG{n}{linexpr} \PYG{o}{=} \PYG{n}{genbuilderx}\PYG{o}{.}\PYG{n}{getExpr}\PYG{p}{(}\PYG{p}{)}
\end{sphinxVerbatim}

\subsubsection{GenConstrBuilder.getSense()}
\label{\detokenize{pyapiref:genconstrbuilder-getsense}}\begin{quote}

\sphinxAtStartPar
\sphinxstylestrong{Synopsis}
\begin{quote}

\sphinxAtStartPar
\sphinxcode{\sphinxupquote{getSense()}}
\end{quote}

\sphinxAtStartPar
\sphinxstylestrong{Description}
\begin{quote}

\sphinxAtStartPar
Retrieve the sense for the linear constraint of a {\hyperref[\detokenize{pyapiref:chappyapi-genconstrbuilder}]{\sphinxcrossref{\DUrole{std,std-ref}{GenConstrBuilder Class}}}} object.
\end{quote}

\sphinxAtStartPar
\sphinxstylestrong{Example}
\end{quote}

\begin{sphinxVerbatim}[commandchars=\\\{\}]
\PYG{c+c1}{\PYGZsh{} Retrieve the sense for the linear constraint of indicator constraint builder genbuilderx}
\PYG{n}{linsense} \PYG{o}{=} \PYG{n}{genbuilderx}\PYG{o}{.}\PYG{n}{getSense}\PYG{p}{(}\PYG{p}{)}
\end{sphinxVerbatim}

\subsubsection{GenConstrBuilder.getIndType()}
\label{\detokenize{pyapiref:genconstrbuilder-getindtype}}\begin{quote}

\sphinxAtStartPar
\sphinxstylestrong{Synopsis}
\begin{quote}

\sphinxAtStartPar
\sphinxcode{\sphinxupquote{getIndType()}}
\end{quote}

\sphinxAtStartPar
\sphinxstylestrong{Description}
\begin{quote}

\sphinxAtStartPar
Retrieve the type for the indicator constraint of a {\hyperref[\detokenize{pyapiref:chappyapi-genconstrbuilder}]{\sphinxcrossref{\DUrole{std,std-ref}{GenConstrBuilder Class}}}} object.
\end{quote}

\sphinxAtStartPar
\sphinxstylestrong{Example}
\end{quote}

\begin{sphinxVerbatim}[commandchars=\\\{\}]
\PYG{c+c1}{\PYGZsh{} Retrieve the type for the indicator constraint of indicator constraint builder genbuilderx}
\PYG{n}{linsense} \PYG{o}{=} \PYG{n}{genbuilderx}\PYG{o}{.}\PYG{n}{getIndType}\PYG{p}{(}\PYG{p}{)}
\end{sphinxVerbatim}

\subsection{GenConstrBuilderArray Class}
\label{\detokenize{pyapiref:genconstrbuilderarray-class}}\label{\detokenize{pyapiref:chappyapi-genconstrbuilderarray}}
\sphinxAtStartPar
To facilitate users to operate on multiple {\hyperref[\detokenize{pyapiref:chappyapi-genconstrbuilder}]{\sphinxcrossref{\DUrole{std,std-ref}{GenConstrBuilder Class}}}} objects, the Python interface of COPT provides GenConstrBuilderArray object with the following methods:

\subsubsection{GenConstrBuilderArray()}
\label{\detokenize{pyapiref:genconstrbuilderarray}}\begin{quote}

\sphinxAtStartPar
\sphinxstylestrong{Synopsis}
\begin{quote}

\sphinxAtStartPar
\sphinxcode{\sphinxupquote{GenConstrBuilderArray(genconstrbuilders=None)}}
\end{quote}

\sphinxAtStartPar
\sphinxstylestrong{Description}
\begin{quote}

\sphinxAtStartPar
Create a {\hyperref[\detokenize{pyapiref:chappyapi-genconstrbuilderarray}]{\sphinxcrossref{\DUrole{std,std-ref}{GenConstrBuilderArray Class}}}} object.

\sphinxAtStartPar
If argument \sphinxcode{\sphinxupquote{genconstrbuilders}} is \sphinxcode{\sphinxupquote{None}}, then create an empty {\hyperref[\detokenize{pyapiref:chappyapi-genconstrbuilderarray}]{\sphinxcrossref{\DUrole{std,std-ref}{GenConstrBuilderArray Class}}}} object;
otherwise use the argument \sphinxcode{\sphinxupquote{genconstrbuilders}} to initialize the newly created {\hyperref[\detokenize{pyapiref:chappyapi-genconstrbuilderarray}]{\sphinxcrossref{\DUrole{std,std-ref}{GenConstrBuilderArray Class}}}} object.
\end{quote}

\sphinxAtStartPar
\sphinxstylestrong{Arguments}
\begin{quote}

\sphinxAtStartPar
\sphinxcode{\sphinxupquote{genconstrbuilders}}
\begin{quote}

\sphinxAtStartPar
Indicator constraint builder to add. Optional, \sphinxcode{\sphinxupquote{None}} by default. It can be {\hyperref[\detokenize{pyapiref:chappyapi-genconstrbuilder}]{\sphinxcrossref{\DUrole{std,std-ref}{GenConstrBuilder Class}}}} object,
{\hyperref[\detokenize{pyapiref:chappyapi-genconstrbuilderarray}]{\sphinxcrossref{\DUrole{std,std-ref}{GenConstrBuilderArray Class}}}} object, list, dict, or {\hyperref[\detokenize{pyapiref:chappyapi-util-tupledict}]{\sphinxcrossref{\DUrole{std,std-ref}{tupledict Class}}}} object.
\end{quote}
\end{quote}

\sphinxAtStartPar
\sphinxstylestrong{Example}
\end{quote}

\begin{sphinxVerbatim}[commandchars=\\\{\}]
\PYG{c+c1}{\PYGZsh{} Create an empty GenConstrBuilderArray object}
\PYG{n}{genbuilderarr} \PYG{o}{=} \PYG{n}{GenConstrBuilderArray}\PYG{p}{(}\PYG{p}{)}
\PYG{c+c1}{\PYGZsh{} Create a GenConstrBuilderArray object and use indicator constraint builder genbuilderx and genbuildery to initialize it.}
\PYG{n}{genbuilderarr} \PYG{o}{=} \PYG{n}{GenConstrBuilderArray}\PYG{p}{(}\PYG{p}{[}\PYG{n}{genbuilderx}\PYG{p}{,} \PYG{n}{genbuildery}\PYG{p}{]}\PYG{p}{)}
\end{sphinxVerbatim}

\subsubsection{GenConstrBuilderArray.pushBack()}
\label{\detokenize{pyapiref:genconstrbuilderarray-pushback}}\begin{quote}

\sphinxAtStartPar
\sphinxstylestrong{Synopsis}
\begin{quote}

\sphinxAtStartPar
\sphinxcode{\sphinxupquote{pushBack(genconstrbuilder)}}
\end{quote}

\sphinxAtStartPar
\sphinxstylestrong{Description}
\begin{quote}

\sphinxAtStartPar
Add single or multiple {\hyperref[\detokenize{pyapiref:chappyapi-genconstrbuilder}]{\sphinxcrossref{\DUrole{std,std-ref}{GenConstrBuilder Class}}}} objects.
\end{quote}

\sphinxAtStartPar
\sphinxstylestrong{Arguments}
\begin{quote}

\sphinxAtStartPar
\sphinxcode{\sphinxupquote{genconstrbuilder}}
\begin{quote}

\sphinxAtStartPar
Indicator constraint builders to add, which can be {\hyperref[\detokenize{pyapiref:chappyapi-genconstrbuilder}]{\sphinxcrossref{\DUrole{std,std-ref}{GenConstrBuilder Class}}}} object,
{\hyperref[\detokenize{pyapiref:chappyapi-genconstrbuilderarray}]{\sphinxcrossref{\DUrole{std,std-ref}{GenConstrBuilderArray Class}}}} object, list, dict, or {\hyperref[\detokenize{pyapiref:chappyapi-util-tupledict}]{\sphinxcrossref{\DUrole{std,std-ref}{tupledict Class}}}} object.
\end{quote}
\end{quote}

\sphinxAtStartPar
\sphinxstylestrong{Example}
\end{quote}

\begin{sphinxVerbatim}[commandchars=\\\{\}]
\PYG{c+c1}{\PYGZsh{} Add an indicator constraint builder to genbuilderarr}
\PYG{n}{genbuilderarr}\PYG{o}{.}\PYG{n}{pushBack}\PYG{p}{(}\PYG{n}{genbuilderx}\PYG{p}{)}
\PYG{c+c1}{\PYGZsh{} Add indicator constraint builders genbuilderx and genbuildery to genbuilderarr}
\PYG{n}{genbuilderarr}\PYG{o}{.}\PYG{n}{pushBack}\PYG{p}{(}\PYG{p}{[}\PYG{n}{genbuilderx}\PYG{p}{,} \PYG{n}{genbuildery}\PYG{p}{]}\PYG{p}{)}
\end{sphinxVerbatim}

\subsubsection{GenConstrBuilderArray.getBuilder()}
\label{\detokenize{pyapiref:genconstrbuilderarray-getbuilder}}\begin{quote}

\sphinxAtStartPar
\sphinxstylestrong{Synopsis}
\begin{quote}

\sphinxAtStartPar
\sphinxcode{\sphinxupquote{getBuilder(idx)}}
\end{quote}

\sphinxAtStartPar
\sphinxstylestrong{Description}
\begin{quote}

\sphinxAtStartPar
Retrieve the indicator constraint builder according to its index in the {\hyperref[\detokenize{pyapiref:chappyapi-genconstrbuilderarray}]{\sphinxcrossref{\DUrole{std,std-ref}{GenConstrBuilderArray Class}}}} object,
and return a {\hyperref[\detokenize{pyapiref:chappyapi-genconstrbuilder}]{\sphinxcrossref{\DUrole{std,std-ref}{GenConstrBuilder Class}}}} object.
\end{quote}

\sphinxAtStartPar
\sphinxstylestrong{Arguments}
\begin{quote}

\sphinxAtStartPar
\sphinxcode{\sphinxupquote{idx}}
\begin{quote}

\sphinxAtStartPar
Index of the indicator constraint builder in the {\hyperref[\detokenize{pyapiref:chappyapi-genconstrbuilderarray}]{\sphinxcrossref{\DUrole{std,std-ref}{GenConstrBuilderArray Class}}}} object, starting with 0.
\end{quote}
\end{quote}

\sphinxAtStartPar
\sphinxstylestrong{Example}
\end{quote}

\begin{sphinxVerbatim}[commandchars=\\\{\}]
\PYG{c+c1}{\PYGZsh{} Retrieve the indicator constraint builder whose index in genbuilderarr is 1}
\PYG{n}{genbuilder} \PYG{o}{=} \PYG{n}{genbuilderarr}\PYG{o}{.}\PYG{n}{getBuilder}\PYG{p}{(}\PYG{l+m+mi}{1}\PYG{p}{)}
\end{sphinxVerbatim}

\subsubsection{GenConstrBuilderArray.getSize()}
\label{\detokenize{pyapiref:genconstrbuilderarray-getsize}}\begin{quote}

\sphinxAtStartPar
\sphinxstylestrong{Synopsis}
\begin{quote}

\sphinxAtStartPar
\sphinxcode{\sphinxupquote{getSize()}}
\end{quote}

\sphinxAtStartPar
\sphinxstylestrong{Description}
\begin{quote}

\sphinxAtStartPar
Retrieve the number of elements in the {\hyperref[\detokenize{pyapiref:chappyapi-genconstrbuilderarray}]{\sphinxcrossref{\DUrole{std,std-ref}{GenConstrBuilderArray Class}}}} object.
\end{quote}

\sphinxAtStartPar
\sphinxstylestrong{Example}
\end{quote}

\begin{sphinxVerbatim}[commandchars=\\\{\}]
\PYG{c+c1}{\PYGZsh{} Retrieve the number of elements in genbuilderarr}
\PYG{n}{genbuildersize} \PYG{o}{=} \PYG{n}{genbuilderarr}\PYG{o}{.}\PYG{n}{getSize}\PYG{p}{(}\PYG{p}{)}
\end{sphinxVerbatim}

\subsection{Column Class}
\label{\detokenize{pyapiref:column-class}}\label{\detokenize{pyapiref:chappyapi-column}}
\sphinxAtStartPar
To facilitate users to model by column, the Python interface of COPT provides Column object with the following methods:

\subsubsection{Column()}
\label{\detokenize{pyapiref:column}}\begin{quote}

\sphinxAtStartPar
\sphinxstylestrong{Synopsis}
\begin{quote}

\sphinxAtStartPar
\sphinxcode{\sphinxupquote{Column(constrs=0.0, coeffs=None)}}
\end{quote}

\sphinxAtStartPar
\sphinxstylestrong{Description}
\begin{quote}

\sphinxAtStartPar
Create an {\hyperref[\detokenize{pyapiref:chappyapi-column}]{\sphinxcrossref{\DUrole{std,std-ref}{Column Class}}}} object.

\sphinxAtStartPar
If argument \sphinxcode{\sphinxupquote{constrs}} is \sphinxcode{\sphinxupquote{None}} and argument \sphinxcode{\sphinxupquote{coeffs}} is \sphinxcode{\sphinxupquote{None}}, then create an empty {\hyperref[\detokenize{pyapiref:chappyapi-column}]{\sphinxcrossref{\DUrole{std,std-ref}{Column Class}}}} object;
otherwise use the argument \sphinxcode{\sphinxupquote{constrs}} and \sphinxcode{\sphinxupquote{coeffs}} to initialize the newly created {\hyperref[\detokenize{pyapiref:chappyapi-column}]{\sphinxcrossref{\DUrole{std,std-ref}{Column Class}}}} object.
If argument \sphinxcode{\sphinxupquote{constrs}} is a {\hyperref[\detokenize{pyapiref:chappyapi-constraint}]{\sphinxcrossref{\DUrole{std,std-ref}{Constraint Class}}}} or {\hyperref[\detokenize{pyapiref:chappyapi-column}]{\sphinxcrossref{\DUrole{std,std-ref}{Column Class}}}} object, then argument \sphinxcode{\sphinxupquote{coeffs}} is a constant.
If argument \sphinxcode{\sphinxupquote{coeffs}} is \sphinxcode{\sphinxupquote{None}}, then it is considered to be constant 1.0;
If argument \sphinxcode{\sphinxupquote{constrs}} is a list and argument \sphinxcode{\sphinxupquote{coeffs}} is \sphinxcode{\sphinxupquote{None}}, then the elements of argument \sphinxcode{\sphinxupquote{constrs}} are constraint\sphinxhyphen{}coefficient pairs;
For other forms of arguments, call method \sphinxcode{\sphinxupquote{addTerms}} to initialize the newly created {\hyperref[\detokenize{pyapiref:chappyapi-column}]{\sphinxcrossref{\DUrole{std,std-ref}{Column Class}}}} object.
\end{quote}

\sphinxAtStartPar
\sphinxstylestrong{Arguments}
\begin{quote}

\sphinxAtStartPar
\sphinxcode{\sphinxupquote{constrs}}
\begin{quote}

\sphinxAtStartPar
Linear constraint.
\end{quote}

\sphinxAtStartPar
\sphinxcode{\sphinxupquote{coeffs}}
\begin{quote}

\sphinxAtStartPar
Coefficient for variables in the linear constraint.
\end{quote}
\end{quote}

\sphinxAtStartPar
\sphinxstylestrong{Example}
\end{quote}

\begin{sphinxVerbatim}[commandchars=\\\{\}]
\PYG{c+c1}{\PYGZsh{} Create an empty Column object}
\PYG{n}{col} \PYG{o}{=} \PYG{n}{Column}\PYG{p}{(}\PYG{p}{)}
\PYG{c+c1}{\PYGZsh{} Create a Column object and add two terms: coefficient is 2 in constraint conx and 3 in constraint cony}
\PYG{n}{col} \PYG{o}{=} \PYG{n}{Column}\PYG{p}{(}\PYG{p}{[}\PYG{p}{(}\PYG{n}{conx}\PYG{p}{,} \PYG{l+m+mi}{2}\PYG{p}{)}\PYG{p}{,} \PYG{p}{(}\PYG{n}{cony}\PYG{p}{,} \PYG{l+m+mi}{3}\PYG{p}{)}\PYG{p}{]}\PYG{p}{)}
\PYG{c+c1}{\PYGZsh{} Create a Column object and add two terms: coefficient is 1 in constraint conxx and 2 in constraint conyy}
\PYG{n}{col} \PYG{o}{=} \PYG{n}{Column}\PYG{p}{(}\PYG{p}{[}\PYG{n}{conxx}\PYG{p}{,} \PYG{n}{conyy}\PYG{p}{]}\PYG{p}{,} \PYG{p}{[}\PYG{l+m+mi}{1}\PYG{p}{,} \PYG{l+m+mi}{2}\PYG{p}{]}\PYG{p}{)}
\end{sphinxVerbatim}

\subsubsection{Column.getCoeff()}
\label{\detokenize{pyapiref:column-getcoeff}}\begin{quote}

\sphinxAtStartPar
\sphinxstylestrong{Synopsis}
\begin{quote}

\sphinxAtStartPar
\sphinxcode{\sphinxupquote{getCoeff(idx)}}
\end{quote}

\sphinxAtStartPar
\sphinxstylestrong{Description}
\begin{quote}

\sphinxAtStartPar
Retrieve the coefficient according to its index in the {\hyperref[\detokenize{pyapiref:chappyapi-column}]{\sphinxcrossref{\DUrole{std,std-ref}{Column Class}}}} object.
\end{quote}

\sphinxAtStartPar
\sphinxstylestrong{Arguments}
\begin{quote}

\sphinxAtStartPar
\sphinxcode{\sphinxupquote{idx}}
\begin{quote}

\sphinxAtStartPar
Index for the element, starting with 0.
\end{quote}
\end{quote}

\sphinxAtStartPar
\sphinxstylestrong{Example}
\end{quote}

\begin{sphinxVerbatim}[commandchars=\\\{\}]
\PYG{c+c1}{\PYGZsh{} Retrieve the coefficient whose index is 0 in col}
\PYG{n}{coeff} \PYG{o}{=} \PYG{n}{col}\PYG{o}{.}\PYG{n}{getCoeff}\PYG{p}{(}\PYG{l+m+mi}{0}\PYG{p}{)}
\end{sphinxVerbatim}

\subsubsection{Column.getConstr()}
\label{\detokenize{pyapiref:column-getconstr}}\begin{quote}

\sphinxAtStartPar
\sphinxstylestrong{Synopsis}
\begin{quote}

\sphinxAtStartPar
\sphinxcode{\sphinxupquote{getConstr(idx)}}
\end{quote}

\sphinxAtStartPar
\sphinxstylestrong{Description}
\begin{quote}

\sphinxAtStartPar
Retrieve the linear constraint according to its index in the {\hyperref[\detokenize{pyapiref:chappyapi-column}]{\sphinxcrossref{\DUrole{std,std-ref}{Column Class}}}} object.
\end{quote}

\sphinxAtStartPar
\sphinxstylestrong{Arguments}
\begin{quote}

\sphinxAtStartPar
\sphinxcode{\sphinxupquote{idx}}
\begin{quote}

\sphinxAtStartPar
Index for the element, starting with 0.
\end{quote}
\end{quote}

\sphinxAtStartPar
\sphinxstylestrong{Example}
\end{quote}

\begin{sphinxVerbatim}[commandchars=\\\{\}]
\PYG{c+c1}{\PYGZsh{} Retrieve the linear constraint whose index is 1 in col}
\PYG{n}{constr} \PYG{o}{=} \PYG{n}{col}\PYG{o}{.}\PYG{n}{getConstr}\PYG{p}{(}\PYG{l+m+mi}{1}\PYG{p}{)}
\end{sphinxVerbatim}

\subsubsection{Column.getSize()}
\label{\detokenize{pyapiref:column-getsize}}\begin{quote}

\sphinxAtStartPar
\sphinxstylestrong{Synopsis}
\begin{quote}

\sphinxAtStartPar
\sphinxcode{\sphinxupquote{getSize()}}
\end{quote}

\sphinxAtStartPar
\sphinxstylestrong{Description}
\begin{quote}

\sphinxAtStartPar
Retrieve the number of elements in the {\hyperref[\detokenize{pyapiref:chappyapi-column}]{\sphinxcrossref{\DUrole{std,std-ref}{Column Class}}}} object.
\end{quote}

\sphinxAtStartPar
\sphinxstylestrong{Example}
\end{quote}

\begin{sphinxVerbatim}[commandchars=\\\{\}]
\PYG{c+c1}{\PYGZsh{} Retrieve the number of elements in col}
\PYG{n}{colsize} \PYG{o}{=} \PYG{n}{col}\PYG{o}{.}\PYG{n}{getSize}\PYG{p}{(}\PYG{p}{)}
\end{sphinxVerbatim}

\subsubsection{Column.addTerm()}
\label{\detokenize{pyapiref:column-addterm}}\begin{quote}

\sphinxAtStartPar
\sphinxstylestrong{Synopsis}
\begin{quote}

\sphinxAtStartPar
\sphinxcode{\sphinxupquote{addTerm(constr, coeff=1.0)}}
\end{quote}

\sphinxAtStartPar
\sphinxstylestrong{Description}
\begin{quote}

\sphinxAtStartPar
Add a new term.
\end{quote}

\sphinxAtStartPar
\sphinxstylestrong{Arguments}
\begin{quote}

\sphinxAtStartPar
\sphinxcode{\sphinxupquote{constr}}
\begin{quote}

\sphinxAtStartPar
The linear constraint for the term to add.
\end{quote}

\sphinxAtStartPar
\sphinxcode{\sphinxupquote{coeff}}
\begin{quote}

\sphinxAtStartPar
The coefficient for the term to add. Optional, 1.0 by default.
\end{quote}
\end{quote}

\sphinxAtStartPar
\sphinxstylestrong{Example}
\end{quote}

\begin{sphinxVerbatim}[commandchars=\\\{\}]
\PYG{c+c1}{\PYGZsh{} Add an term to col, whose constraint is cony and coefficient is 2.0}
\PYG{n}{col}\PYG{o}{.}\PYG{n}{addTerm}\PYG{p}{(}\PYG{n}{cony}\PYG{p}{,} \PYG{l+m+mf}{2.0}\PYG{p}{)}
\PYG{c+c1}{\PYGZsh{} Add an term to col, whose constraint is conx and coefficient is 1.0}
\PYG{n}{col}\PYG{o}{.}\PYG{n}{addTerm}\PYG{p}{(}\PYG{n}{conx}\PYG{p}{)}
\end{sphinxVerbatim}

\subsubsection{Column.addTerms()}
\label{\detokenize{pyapiref:column-addterms}}\begin{quote}

\sphinxAtStartPar
\sphinxstylestrong{Synopsis}
\begin{quote}

\sphinxAtStartPar
\sphinxcode{\sphinxupquote{addTerms(constrs, coeffs)}}
\end{quote}

\sphinxAtStartPar
\sphinxstylestrong{Description}
\begin{quote}

\sphinxAtStartPar
Add single or multiple terms.

\sphinxAtStartPar
If argument \sphinxcode{\sphinxupquote{constrs}} is {\hyperref[\detokenize{pyapiref:chappyapi-constraint}]{\sphinxcrossref{\DUrole{std,std-ref}{Constraint Class}}}} object, then argument \sphinxcode{\sphinxupquote{coeffs}} is constant;
If argument \sphinxcode{\sphinxupquote{constrs}} is {\hyperref[\detokenize{pyapiref:chappyapi-constrarray}]{\sphinxcrossref{\DUrole{std,std-ref}{ConstrArray Class}}}} object or list, then argument \sphinxcode{\sphinxupquote{coeffs}} is constant or list;
If argument \sphinxcode{\sphinxupquote{constrs}} is dictionary or {\hyperref[\detokenize{pyapiref:chappyapi-util-tupledict}]{\sphinxcrossref{\DUrole{std,std-ref}{tupledict Class}}}} object, then argument \sphinxcode{\sphinxupquote{coeffs}} is constant, dict, or
{\hyperref[\detokenize{pyapiref:chappyapi-util-tupledict}]{\sphinxcrossref{\DUrole{std,std-ref}{tupledict Class}}}} object.
\end{quote}

\sphinxAtStartPar
\sphinxstylestrong{Arguments}
\begin{quote}

\sphinxAtStartPar
\sphinxcode{\sphinxupquote{constrs}}
\begin{quote}

\sphinxAtStartPar
The linear constraints for terms to add.
\end{quote}

\sphinxAtStartPar
\sphinxcode{\sphinxupquote{coeffs}}
\begin{quote}

\sphinxAtStartPar
The coefficients for terms to add.
\end{quote}
\end{quote}

\sphinxAtStartPar
\sphinxstylestrong{Example}
\end{quote}

\begin{sphinxVerbatim}[commandchars=\\\{\}]
\PYG{c+c1}{\PYGZsh{} Add two terms: constraint conx with coefficient 2.0, constraint cony with coefficient 3.0}
\PYG{n}{col}\PYG{o}{.}\PYG{n}{addTerms}\PYG{p}{(}\PYG{p}{[}\PYG{n}{conx}\PYG{p}{,} \PYG{n}{cony}\PYG{p}{]}\PYG{p}{,} \PYG{p}{[}\PYG{l+m+mf}{2.0}\PYG{p}{,} \PYG{l+m+mf}{3.0}\PYG{p}{]}\PYG{p}{)}
\end{sphinxVerbatim}

\subsubsection{Column.addColumn()}
\label{\detokenize{pyapiref:column-addcolumn}}\begin{quote}

\sphinxAtStartPar
\sphinxstylestrong{Synopsis}
\begin{quote}

\sphinxAtStartPar
\sphinxcode{\sphinxupquote{addColumn(col, mult=1.0)}}
\end{quote}

\sphinxAtStartPar
\sphinxstylestrong{Description}
\begin{quote}

\sphinxAtStartPar
Add a new column to current column.
\end{quote}

\sphinxAtStartPar
\sphinxstylestrong{Arguments}
\begin{quote}

\sphinxAtStartPar
\sphinxcode{\sphinxupquote{col}}
\begin{quote}

\sphinxAtStartPar
Column to add.
\end{quote}

\sphinxAtStartPar
\sphinxcode{\sphinxupquote{mult}}
\begin{quote}

\sphinxAtStartPar
Magnification coefficient for added column. Optional, 1.0 by default.
\end{quote}
\end{quote}

\sphinxAtStartPar
\sphinxstylestrong{Example}
\end{quote}

\begin{sphinxVerbatim}[commandchars=\\\{\}]
\PYG{c+c1}{\PYGZsh{} Add column coly to column colx. The magnification coefficient for coly is 2.0}
\PYG{n}{colx}\PYG{o}{.}\PYG{n}{addColumn}\PYG{p}{(}\PYG{n}{coly}\PYG{p}{,} \PYG{l+m+mf}{2.0}\PYG{p}{)}
\end{sphinxVerbatim}

\subsubsection{Column.clone()}
\label{\detokenize{pyapiref:column-clone}}\begin{quote}

\sphinxAtStartPar
\sphinxstylestrong{Synopsis}
\begin{quote}

\sphinxAtStartPar
\sphinxcode{\sphinxupquote{clone()}}
\end{quote}

\sphinxAtStartPar
\sphinxstylestrong{Description}
\begin{quote}

\sphinxAtStartPar
Create a deep copy of a column.
\end{quote}

\sphinxAtStartPar
\sphinxstylestrong{Example}
\end{quote}

\begin{sphinxVerbatim}[commandchars=\\\{\}]
\PYG{c+c1}{\PYGZsh{} Create a deep copy of column col}
\PYG{n}{colcopy} \PYG{o}{=} \PYG{n}{col}\PYG{o}{.}\PYG{n}{clone}\PYG{p}{(}\PYG{p}{)}
\end{sphinxVerbatim}

\subsubsection{Column.remove()}
\label{\detokenize{pyapiref:column-remove}}\begin{quote}

\sphinxAtStartPar
\sphinxstylestrong{Synopsis}
\begin{quote}

\sphinxAtStartPar
\sphinxcode{\sphinxupquote{remove(item)}}
\end{quote}

\sphinxAtStartPar
\sphinxstylestrong{Description}
\begin{quote}

\sphinxAtStartPar
Remove a term from a column.

\sphinxAtStartPar
If argument \sphinxcode{\sphinxupquote{item}} is a constant, then remove the term by its index;
otherwise argument \sphinxcode{\sphinxupquote{item}} is a {\hyperref[\detokenize{pyapiref:chappyapi-constraint}]{\sphinxcrossref{\DUrole{std,std-ref}{Constraint Class}}}} object.
\end{quote}

\sphinxAtStartPar
\sphinxstylestrong{Arguments}
\begin{quote}

\sphinxAtStartPar
\sphinxcode{\sphinxupquote{item}}
\begin{quote}

\sphinxAtStartPar
Constant index or the linear constraint for the term to be removed.
\end{quote}
\end{quote}

\sphinxAtStartPar
\sphinxstylestrong{Example}
\end{quote}

\begin{sphinxVerbatim}[commandchars=\\\{\}]
\PYG{c+c1}{\PYGZsh{} Remove the term whose index is 2 from column col}
\PYG{n}{col}\PYG{o}{.}\PYG{n}{remove}\PYG{p}{(}\PYG{l+m+mi}{2}\PYG{p}{)}
\PYG{c+c1}{\PYGZsh{} Remove the term of the linear constraint conx from col}
\PYG{n}{col}\PYG{o}{.}\PYG{n}{remove}\PYG{p}{(}\PYG{n}{conx}\PYG{p}{)}
\end{sphinxVerbatim}

\subsubsection{Column.clear()}
\label{\detokenize{pyapiref:column-clear}}\begin{quote}

\sphinxAtStartPar
\sphinxstylestrong{Synopsis}
\begin{quote}

\sphinxAtStartPar
\sphinxcode{\sphinxupquote{clear()}}
\end{quote}

\sphinxAtStartPar
\sphinxstylestrong{Description}
\begin{quote}

\sphinxAtStartPar
Remove all terms from a column.
\end{quote}

\sphinxAtStartPar
\sphinxstylestrong{Example}
\end{quote}

\begin{sphinxVerbatim}[commandchars=\\\{\}]
\PYG{c+c1}{\PYGZsh{} Remove all terms from column col}
\PYG{n}{col}\PYG{o}{.}\PYG{n}{clear}\PYG{p}{(}\PYG{p}{)}
\end{sphinxVerbatim}

\subsection{ColumnArray Class}
\label{\detokenize{pyapiref:columnarray-class}}\label{\detokenize{pyapiref:chappyapi-columnarray}}
\sphinxAtStartPar
To facilitate users to operate on multiple {\hyperref[\detokenize{pyapiref:chappyapi-column}]{\sphinxcrossref{\DUrole{std,std-ref}{Column Class}}}} objects, the Python interface of COPT provides ColumnArray object with the following methods:

\subsubsection{ColumnArray()}
\label{\detokenize{pyapiref:columnarray}}\begin{quote}

\sphinxAtStartPar
\sphinxstylestrong{Synopsis}
\begin{quote}

\sphinxAtStartPar
\sphinxcode{\sphinxupquote{ColumnArray(columns=None)}}
\end{quote}

\sphinxAtStartPar
\sphinxstylestrong{Description}
\begin{quote}

\sphinxAtStartPar
Create a {\hyperref[\detokenize{pyapiref:chappyapi-columnarray}]{\sphinxcrossref{\DUrole{std,std-ref}{ColumnArray Class}}}} object.

\sphinxAtStartPar
If argument \sphinxcode{\sphinxupquote{columns}} is \sphinxcode{\sphinxupquote{None}}, then create an empty {\hyperref[\detokenize{pyapiref:chappyapi-columnarray}]{\sphinxcrossref{\DUrole{std,std-ref}{ColumnArray Class}}}} object;
otherwise use argument \sphinxcode{\sphinxupquote{columns}} to initialize the newly created {\hyperref[\detokenize{pyapiref:chappyapi-columnarray}]{\sphinxcrossref{\DUrole{std,std-ref}{ColumnArray Class}}}} object.
\end{quote}

\sphinxAtStartPar
\sphinxstylestrong{Arguments}
\begin{quote}

\sphinxAtStartPar
\sphinxcode{\sphinxupquote{columns}}
\begin{quote}

\sphinxAtStartPar
Columns to add. Optional, \sphinxcode{\sphinxupquote{None}} by default. It can be {\hyperref[\detokenize{pyapiref:chappyapi-column}]{\sphinxcrossref{\DUrole{std,std-ref}{Column Class}}}} object,
{\hyperref[\detokenize{pyapiref:chappyapi-columnarray}]{\sphinxcrossref{\DUrole{std,std-ref}{ColumnArray Class}}}} object, list, dict, or {\hyperref[\detokenize{pyapiref:chappyapi-util-tupledict}]{\sphinxcrossref{\DUrole{std,std-ref}{tupledict Class}}}} object.
\end{quote}
\end{quote}

\sphinxAtStartPar
\sphinxstylestrong{Example}
\end{quote}

\begin{sphinxVerbatim}[commandchars=\\\{\}]
\PYG{c+c1}{\PYGZsh{} Create an empty ColumnArray object}
\PYG{n}{colarr} \PYG{o}{=} \PYG{n}{ColumnArray}\PYG{p}{(}\PYG{p}{)}
\PYG{c+c1}{\PYGZsh{} Create a ColumnArray object and use columns colx and coly to initialize it}
\PYG{n}{colarr} \PYG{o}{=} \PYG{n}{ColumnArray}\PYG{p}{(}\PYG{p}{[}\PYG{n}{colx}\PYG{p}{,} \PYG{n}{coly}\PYG{p}{]}\PYG{p}{)}
\end{sphinxVerbatim}

\subsubsection{ColumnArray.pushBack()}
\label{\detokenize{pyapiref:columnarray-pushback}}\begin{quote}

\sphinxAtStartPar
\sphinxstylestrong{Synopsis}
\begin{quote}

\sphinxAtStartPar
\sphinxcode{\sphinxupquote{pushBack(column)}}
\end{quote}

\sphinxAtStartPar
\sphinxstylestrong{Description}
\begin{quote}

\sphinxAtStartPar
Add single or multiple {\hyperref[\detokenize{pyapiref:chappyapi-column}]{\sphinxcrossref{\DUrole{std,std-ref}{Column Class}}}} objects.
\end{quote}

\sphinxAtStartPar
\sphinxstylestrong{Arguments}
\begin{quote}

\sphinxAtStartPar
\sphinxcode{\sphinxupquote{column}}
\begin{quote}

\sphinxAtStartPar
Columns to add, which can be {\hyperref[\detokenize{pyapiref:chappyapi-column}]{\sphinxcrossref{\DUrole{std,std-ref}{Column Class}}}} object, {\hyperref[\detokenize{pyapiref:chappyapi-columnarray}]{\sphinxcrossref{\DUrole{std,std-ref}{ColumnArray Class}}}} object,
list, dict, or {\hyperref[\detokenize{pyapiref:chappyapi-util-tupledict}]{\sphinxcrossref{\DUrole{std,std-ref}{tupledict Class}}}} object.
\end{quote}
\end{quote}

\sphinxAtStartPar
\sphinxstylestrong{Example}
\end{quote}

\begin{sphinxVerbatim}[commandchars=\\\{\}]
\PYG{c+c1}{\PYGZsh{} Add column colx to colarr}
\PYG{n}{colarr}\PYG{o}{.}\PYG{n}{pushBack}\PYG{p}{(}\PYG{n}{colx}\PYG{p}{)}
\PYG{c+c1}{\PYGZsh{} Add columns colx and coly to colarr}
\PYG{n}{colarr}\PYG{o}{.}\PYG{n}{pushBack}\PYG{p}{(}\PYG{p}{[}\PYG{n}{colx}\PYG{p}{,} \PYG{n}{coly}\PYG{p}{]}\PYG{p}{)}
\end{sphinxVerbatim}

\subsubsection{ColumnArray.getColumn()}
\label{\detokenize{pyapiref:columnarray-getcolumn}}\begin{quote}

\sphinxAtStartPar
\sphinxstylestrong{Synopsis}
\begin{quote}

\sphinxAtStartPar
\sphinxcode{\sphinxupquote{getColumn(idx)}}
\end{quote}

\sphinxAtStartPar
\sphinxstylestrong{Description}
\begin{quote}

\sphinxAtStartPar
Retrieve the column according to its index in a {\hyperref[\detokenize{pyapiref:chappyapi-columnarray}]{\sphinxcrossref{\DUrole{std,std-ref}{ColumnArray Class}}}} object.
Return a {\hyperref[\detokenize{pyapiref:chappyapi-column}]{\sphinxcrossref{\DUrole{std,std-ref}{Column Class}}}} object.
\end{quote}

\sphinxAtStartPar
\sphinxstylestrong{Arguments}
\begin{quote}

\sphinxAtStartPar
\sphinxcode{\sphinxupquote{idx}}
\begin{quote}

\sphinxAtStartPar
Index of the column in the {\hyperref[\detokenize{pyapiref:chappyapi-columnarray}]{\sphinxcrossref{\DUrole{std,std-ref}{ColumnArray Class}}}} object, starting with 0.
\end{quote}
\end{quote}

\sphinxAtStartPar
\sphinxstylestrong{Example}
\end{quote}

\begin{sphinxVerbatim}[commandchars=\\\{\}]
\PYG{c+c1}{\PYGZsh{} Retrieve the column whose index is 1 in colarr}
\PYG{n}{col} \PYG{o}{=} \PYG{n}{colarr}\PYG{o}{.}\PYG{n}{getColumn}\PYG{p}{(}\PYG{l+m+mi}{1}\PYG{p}{)}
\end{sphinxVerbatim}

\subsubsection{ColumnArray.getSize()}
\label{\detokenize{pyapiref:columnarray-getsize}}\begin{quote}

\sphinxAtStartPar
\sphinxstylestrong{Synopsis}
\begin{quote}

\sphinxAtStartPar
\sphinxcode{\sphinxupquote{getSize()}}
\end{quote}

\sphinxAtStartPar
\sphinxstylestrong{Description}
\begin{quote}

\sphinxAtStartPar
Retrieve the number of elements in a {\hyperref[\detokenize{pyapiref:chappyapi-columnarray}]{\sphinxcrossref{\DUrole{std,std-ref}{ColumnArray Class}}}} object.
\end{quote}

\sphinxAtStartPar
\sphinxstylestrong{Example}
\end{quote}

\begin{sphinxVerbatim}[commandchars=\\\{\}]
\PYG{c+c1}{\PYGZsh{} Retrieve the number of element in colarr}
\PYG{n}{colsize} \PYG{o}{=} \PYG{n}{colarr}\PYG{o}{.}\PYG{n}{getSize}\PYG{p}{(}\PYG{p}{)}
\end{sphinxVerbatim}

\subsubsection{ColumnArray.clear()}
\label{\detokenize{pyapiref:columnarray-clear}}\begin{quote}

\sphinxAtStartPar
\sphinxstylestrong{Synopsis}
\begin{quote}

\sphinxAtStartPar
\sphinxcode{\sphinxupquote{clear()}}
\end{quote}

\sphinxAtStartPar
\sphinxstylestrong{Description}
\begin{quote}

\sphinxAtStartPar
Remove all terms from a {\hyperref[\detokenize{pyapiref:chappyapi-columnarray}]{\sphinxcrossref{\DUrole{std,std-ref}{ColumnArray Class}}}} object.
\end{quote}

\sphinxAtStartPar
\sphinxstylestrong{Example}
\end{quote}

\begin{sphinxVerbatim}[commandchars=\\\{\}]
\PYG{c+c1}{\PYGZsh{} Remove all terms from colarr}
\PYG{n}{colarr}\PYG{o}{.}\PYG{n}{clear}\PYG{p}{(}\PYG{p}{)}
\end{sphinxVerbatim}

\subsection{MVar Class}
\label{\detokenize{pyapiref:mvar-class}}\label{\detokenize{pyapiref:chappyapi-mvar}}
\sphinxAtStartPar
The MVar class is used in COPT to build multi\sphinxhyphen{}dimensional variables and supports NumPy’s multi\sphinxhyphen{}dimensional array operations.
We recommend to generate it through the method \sphinxcode{\sphinxupquote{addVars}} or \sphinxcode{\sphinxupquote{addMVar}} of the model class, although it can also be converted and generated through the two built\sphinxhyphen{}in class methods \sphinxcode{\sphinxupquote{fromlist}} and \sphinxcode{\sphinxupquote{fromvar}} .
It supports a unified set of indexing and slicing rules.
For detailed usage, see {\hyperref[\detokenize{matrix:chapmatrixindex}]{\sphinxcrossref{\DUrole{std,std-ref}{Matrix Modeling: Indexing and Slicing}}}}.

\sphinxAtStartPar
The following member methods are provided:

\subsubsection{MVar.fromlist()}
\label{\detokenize{pyapiref:mvar-fromlist}}\begin{quote}

\sphinxAtStartPar
\sphinxstylestrong{Synopsis}
\begin{quote}

\sphinxAtStartPar
\sphinxcode{\sphinxupquote{fromlist(vars)}}
\end{quote}

\sphinxAtStartPar
\sphinxstylestrong{Description}
\begin{quote}

\sphinxAtStartPar
Generate a {\hyperref[\detokenize{pyapiref:chappyapi-mvar}]{\sphinxcrossref{\DUrole{std,std-ref}{MVar Class}}}} object from a set of {\hyperref[\detokenize{pyapiref:chappyapi-var}]{\sphinxcrossref{\DUrole{std,std-ref}{Var Class}}}} objects.
This is the class generation method and can be called directly without MVar object.
\end{quote}

\sphinxAtStartPar
\sphinxstylestrong{Arguments}
\begin{quote}

\sphinxAtStartPar
\sphinxcode{\sphinxupquote{vars}}
\begin{quote}

\sphinxAtStartPar
A set of Var objects, which can be a multi\sphinxhyphen{}dimensional list or ndarray.
\end{quote}
\end{quote}

\sphinxAtStartPar
\sphinxstylestrong{Return value}
\begin{quote}

\sphinxAtStartPar
new MVar object whose dimensions depend on the dimensions of the arguments vars.
\end{quote}

\sphinxAtStartPar
\sphinxstylestrong{Example}

\begin{sphinxVerbatim}[commandchars=\\\{\}]
\PYG{n+nb}{vars} \PYG{o}{=} \PYG{n}{model}\PYG{o}{.}\PYG{n}{addVars}\PYG{p}{(}\PYG{l+m+mi}{4}\PYG{p}{)}
\PYG{n}{mx\PYGZus{}1d} \PYG{o}{=} \PYG{n}{MVar}\PYG{o}{.}\PYG{n}{fromlist}\PYG{p}{(}\PYG{n+nb}{vars}\PYG{p}{)}
\PYG{n}{mx\PYGZus{}2d} \PYG{o}{=} \PYG{n}{MVar}\PYG{o}{.}\PYG{n}{fromlist}\PYG{p}{(}\PYG{p}{[}\PYG{n+nb}{vars}\PYG{p}{[}\PYG{l+m+mi}{0}\PYG{p}{]}\PYG{p}{,} \PYG{n+nb}{vars}\PYG{p}{[}\PYG{l+m+mi}{1}\PYG{p}{]}\PYG{p}{]}\PYG{p}{,} \PYG{p}{[}\PYG{n+nb}{vars}\PYG{p}{[}\PYG{l+m+mi}{2}\PYG{p}{]}\PYG{p}{,} \PYG{n+nb}{vars}\PYG{p}{[}\PYG{l+m+mi}{3}\PYG{p}{]}\PYG{p}{]}\PYG{p}{)}
\end{sphinxVerbatim}
\end{quote}

\subsubsection{MVar.fromvar()}
\label{\detokenize{pyapiref:mvar-fromvar}}\begin{quote}

\sphinxAtStartPar
\sphinxstylestrong{Synopsis}
\begin{quote}

\sphinxAtStartPar
\sphinxcode{\sphinxupquote{fromvar(var)}}
\end{quote}

\sphinxAtStartPar
\sphinxstylestrong{Description}
\begin{quote}

\sphinxAtStartPar
Generate a 0\sphinxhyphen{}dimensional {\hyperref[\detokenize{pyapiref:chappyapi-mvar}]{\sphinxcrossref{\DUrole{std,std-ref}{MVar Class}}}} object from a {\hyperref[\detokenize{pyapiref:chappyapi-var}]{\sphinxcrossref{\DUrole{std,std-ref}{Var Class}}}} object.
This is the class generation method and can be called directly without MVar object.
\end{quote}

\sphinxAtStartPar
\sphinxstylestrong{Arguments}
\begin{quote}

\sphinxAtStartPar
\sphinxcode{\sphinxupquote{var}}
\begin{quote}

\sphinxAtStartPar
A Var object.
\end{quote}
\end{quote}

\sphinxAtStartPar
\sphinxstylestrong{Return value}
\begin{quote}

\sphinxAtStartPar
The new 0\sphinxhyphen{}dimensional MVar object.
\end{quote}

\sphinxAtStartPar
\sphinxstylestrong{Example}

\begin{sphinxVerbatim}[commandchars=\\\{\}]
\PYG{n}{x} \PYG{o}{=} \PYG{n}{model}\PYG{o}{.}\PYG{n}{addVar}\PYG{p}{(}\PYG{p}{)}
\PYG{n}{mx\PYGZus{}0d} \PYG{o}{=} \PYG{n}{MVar}\PYG{o}{.}\PYG{n}{fromvar}\PYG{p}{(}\PYG{n}{x}\PYG{p}{)}
\end{sphinxVerbatim}
\end{quote}

\subsubsection{MVar.clone()}
\label{\detokenize{pyapiref:mvar-clone}}\begin{quote}

\sphinxAtStartPar
\sphinxstylestrong{Synopsis}
\begin{quote}

\sphinxAtStartPar
\sphinxcode{\sphinxupquote{clone()}}
\end{quote}

\sphinxAtStartPar
\sphinxstylestrong{Description}
\begin{quote}

\sphinxAtStartPar
Deep\sphinxhyphen{}copy a {\hyperref[\detokenize{pyapiref:chappyapi-mvar}]{\sphinxcrossref{\DUrole{std,std-ref}{MVar Class}}}} object.
\end{quote}

\sphinxAtStartPar
\sphinxstylestrong{Return value}
\begin{quote}

\sphinxAtStartPar
new MVar object
\end{quote}

\sphinxAtStartPar
\sphinxstylestrong{Example}

\begin{sphinxVerbatim}[commandchars=\\\{\}]
\PYG{c+c1}{\PYGZsh{} Create a 2\PYGZhy{}D variable and make a copy. Note that the actual variable is not incremented.}
\PYG{n}{mx} \PYG{o}{=} \PYG{n}{model}\PYG{o}{.}\PYG{n}{addMVar}\PYG{p}{(}\PYG{p}{(}\PYG{l+m+mi}{3}\PYG{p}{,} \PYG{l+m+mi}{2}\PYG{p}{)}\PYG{p}{,} \PYG{n}{nameprefix}\PYG{o}{=}\PYG{l+s+s2}{\PYGZdq{}}\PYG{l+s+s2}{mx}\PYG{l+s+s2}{\PYGZdq{}}\PYG{p}{)}
\PYG{n}{mx\PYGZus{}copy} \PYG{o}{=} \PYG{n}{mx}\PYG{o}{.}\PYG{n}{clone}\PYG{p}{(}\PYG{p}{)}
\end{sphinxVerbatim}
\end{quote}

\subsubsection{MVar.diagonal()}
\label{\detokenize{pyapiref:mvar-diagonal}}\begin{quote}

\sphinxAtStartPar
\sphinxstylestrong{Synopsis}
\begin{quote}

\sphinxAtStartPar
\sphinxcode{\sphinxupquote{diagonal(offset=0, axis1=0, axis2=1)}}
\end{quote}

\sphinxAtStartPar
\sphinxstylestrong{Description}
\begin{quote}

\sphinxAtStartPar
Generate a {\hyperref[\detokenize{pyapiref:chappyapi-mvar}]{\sphinxcrossref{\DUrole{std,std-ref}{MVar Class}}}} object whose elements are the elements on the diagonal of the original MVar object.
\end{quote}

\sphinxAtStartPar
\sphinxstylestrong{Arguments}
\begin{quote}

\sphinxAtStartPar
\sphinxcode{\sphinxupquote{offset}}
\begin{quote}

\sphinxAtStartPar
Optional parameter, indicating the offset of the diagonal, the default value is 0. If the value is greater than 0, it means the diagonal upward offset; if the value is less than 0, it means the diagonal downward offset.
\end{quote}

\sphinxAtStartPar
\sphinxcode{\sphinxupquote{axis1}}
\begin{quote}

\sphinxAtStartPar
Optional parameter, the axis to use as the first axis of the 2D sub MVar, from where the diagonal should start. The default first axis is 0.
\end{quote}

\sphinxAtStartPar
\sphinxcode{\sphinxupquote{axis2}}
\begin{quote}

\sphinxAtStartPar
Optional parameter, the axis to use as the second axis of the 2D sub MVar, from where the diagonal should start. The default second axis is 1.
\end{quote}
\end{quote}

\sphinxAtStartPar
\sphinxstylestrong{Return value}
\begin{quote}

\sphinxAtStartPar
new MVar object
\end{quote}

\sphinxAtStartPar
\sphinxstylestrong{Example}

\begin{sphinxVerbatim}[commandchars=\\\{\}]
\PYG{n}{mx} \PYG{o}{=} \PYG{n}{model}\PYG{o}{.}\PYG{n}{addMVar}\PYG{p}{(}\PYG{p}{(}\PYG{l+m+mi}{5}\PYG{p}{,} \PYG{l+m+mi}{5}\PYG{p}{)}\PYG{p}{,} \PYG{n}{nameprefix}\PYG{o}{=}\PYG{l+s+s2}{\PYGZdq{}}\PYG{l+s+s2}{mx}\PYG{l+s+s2}{\PYGZdq{}}\PYG{p}{)}
\PYG{n}{diag\PYGZus{}m0} \PYG{o}{=} \PYG{n}{mx}\PYG{o}{.}\PYG{n}{diagonal}\PYG{p}{(}\PYG{p}{)}
\PYG{n}{diag\PYGZus{}a1} \PYG{o}{=} \PYG{n}{mx}\PYG{o}{.}\PYG{n}{diagonal}\PYG{p}{(}\PYG{l+m+mi}{1}\PYG{p}{)}
\PYG{n}{diag\PYGZus{}b1} \PYG{o}{=} \PYG{n}{mx}\PYG{o}{.}\PYG{n}{diagonal}\PYG{p}{(}\PYG{o}{\PYGZhy{}}\PYG{l+m+mi}{1}\PYG{p}{)}
\end{sphinxVerbatim}
\end{quote}

\subsubsection{MVar.getInfo()}
\label{\detokenize{pyapiref:mvar-getinfo}}\begin{quote}

\sphinxAtStartPar
\sphinxstylestrong{Synopsis}
\begin{quote}

\sphinxAtStartPar
\sphinxcode{\sphinxupquote{getInfo(infoname)}}
\end{quote}

\sphinxAtStartPar
\sphinxstylestrong{Description}
\begin{quote}

\sphinxAtStartPar
Get the information value of each variable inside MVar.
\end{quote}

\sphinxAtStartPar
\sphinxstylestrong{Arguments}
\begin{quote}

\sphinxAtStartPar
\sphinxcode{\sphinxupquote{infoname}}
\begin{quote}

\sphinxAtStartPar
The name of the information being queried.
Please refer to {\hyperref[\detokenize{information:chapinfo}]{\sphinxcrossref{\DUrole{std,std-ref}{Information Section}}}} for possible values.
\end{quote}
\end{quote}

\sphinxAtStartPar
\sphinxstylestrong{Return value}
\begin{quote}

\sphinxAtStartPar
Returns a \sphinxcode{\sphinxupquote{NumPy}} ndarray with the same dimension as the \sphinxcode{\sphinxupquote{MVar}} object, whose elements are the information values of the corresponding variable.
\end{quote}

\sphinxAtStartPar
\sphinxstylestrong{Example}

\begin{sphinxVerbatim}[commandchars=\\\{\}]
\PYG{n}{mx} \PYG{o}{=} \PYG{n}{model}\PYG{o}{.}\PYG{n}{addMVar}\PYG{p}{(}\PYG{l+m+mi}{3}\PYG{p}{)}
\PYG{n+nb}{print}\PYG{p}{(}\PYG{n}{mx}\PYG{o}{.}\PYG{n}{getInfo}\PYG{p}{(}\PYG{l+s+s2}{\PYGZdq{}}\PYG{l+s+s2}{LB}\PYG{l+s+s2}{\PYGZdq{}}\PYG{p}{)}\PYG{p}{)}
\end{sphinxVerbatim}
\end{quote}

\subsubsection{MVar.item()}
\label{\detokenize{pyapiref:mvar-item}}\begin{quote}

\sphinxAtStartPar
\sphinxstylestrong{Synopsis}
\begin{quote}

\sphinxAtStartPar
\sphinxcode{\sphinxupquote{item()}}
\end{quote}

\sphinxAtStartPar
\sphinxstylestrong{Description}
\begin{quote}

\sphinxAtStartPar
Get the Var variable inside the 0\sphinxhyphen{}dimensional MVar. Raises a ValueError exception if the MVar object is not 0\sphinxhyphen{}dimensional.
\end{quote}

\sphinxAtStartPar
\sphinxstylestrong{Return value}
\begin{quote}

\sphinxAtStartPar
Returns the Var object.
\end{quote}

\sphinxAtStartPar
\sphinxstylestrong{Example}

\begin{sphinxVerbatim}[commandchars=\\\{\}]
\PYG{n}{mx} \PYG{o}{=} \PYG{n}{model}\PYG{o}{.}\PYG{n}{addMVar}\PYG{p}{(}\PYG{l+m+mi}{3}\PYG{p}{)}
\PYG{n}{var} \PYG{o}{=} \PYG{n}{mx}\PYG{p}{[}\PYG{l+m+mi}{0}\PYG{p}{]}\PYG{o}{.}\PYG{n}{item}\PYG{p}{(}\PYG{p}{)}
\end{sphinxVerbatim}
\end{quote}

\subsubsection{MVar.reshape()}
\label{\detokenize{pyapiref:mvar-reshape}}\begin{quote}

\sphinxAtStartPar
\sphinxstylestrong{Synopsis}
\begin{quote}

\sphinxAtStartPar
\sphinxcode{\sphinxupquote{reshape(shape, order=\textquotesingle{}C\textquotesingle{})}}
\end{quote}

\sphinxAtStartPar
\sphinxstylestrong{Description}
\begin{quote}

\sphinxAtStartPar
Returns a new MVar object whose elements remain unchanged but whose shape is transformed by the parameter shape.
\end{quote}

\sphinxAtStartPar
\sphinxstylestrong{Arguments}
\begin{quote}

\sphinxAtStartPar
\sphinxcode{\sphinxupquote{shape}}
\begin{quote}

\sphinxAtStartPar
The value is an integer, or a tuple of integers. which represents the shape of the new MVar object.
\end{quote}

\sphinxAtStartPar
\sphinxcode{\sphinxupquote{order}}
\begin{quote}

\sphinxAtStartPar
Optional parameter, the default is the character ‘C’, which means it is compatible with the C language, that is, it is stored in rows; it can also be set to the character ‘F’, that is, it is stored in columns, and it is compatible with the Fortune language.
\end{quote}
\end{quote}

\sphinxAtStartPar
\sphinxstylestrong{Return value}
\begin{quote}

\sphinxAtStartPar
Returns a new MVar object with the same elements as the original MVar object but with a different shape.
\end{quote}

\sphinxAtStartPar
\sphinxstylestrong{Example}

\begin{sphinxVerbatim}[commandchars=\\\{\}]
\PYG{n}{mx} \PYG{o}{=} \PYG{n}{model}\PYG{o}{.}\PYG{n}{addMVar}\PYG{p}{(}\PYG{l+m+mi}{6}\PYG{p}{)}
\PYG{n}{mx\PYGZus{}2x3} \PYG{o}{=} \PYG{n}{mx}\PYG{o}{.}\PYG{n}{reshape}\PYG{p}{(}\PYG{p}{(}\PYG{l+m+mi}{2}\PYG{p}{,} \PYG{l+m+mi}{3}\PYG{p}{)}\PYG{p}{)}
\end{sphinxVerbatim}
\end{quote}

\subsubsection{MVar.setInfo()}
\label{\detokenize{pyapiref:mvar-setinfo}}\begin{quote}

\sphinxAtStartPar
\sphinxstylestrong{Synopsis}
\begin{quote}

\sphinxAtStartPar
\sphinxcode{\sphinxupquote{setInfo(infoname, newval)}}
\end{quote}

\sphinxAtStartPar
\sphinxstylestrong{Description}
\begin{quote}

\sphinxAtStartPar
Sets the information value of each variable inside the MVar.
\end{quote}

\sphinxAtStartPar
\sphinxstylestrong{Arguments}
\begin{quote}

\sphinxAtStartPar
\sphinxcode{\sphinxupquote{infoname}}
\begin{quote}

\sphinxAtStartPar
The name of the information being queried.
Please refer to {\hyperref[\detokenize{information:chapinfo}]{\sphinxcrossref{\DUrole{std,std-ref}{Information Section}}}} for possible values.
\end{quote}

\sphinxAtStartPar
\sphinxcode{\sphinxupquote{newval}}
\begin{quote}

\sphinxAtStartPar
The information value to be set.
\end{quote}
\end{quote}

\sphinxAtStartPar
\sphinxstylestrong{Example}

\begin{sphinxVerbatim}[commandchars=\\\{\}]
\PYG{n}{mx} \PYG{o}{=} \PYG{n}{model}\PYG{o}{.}\PYG{n}{addMVar}\PYG{p}{(}\PYG{l+m+mi}{3}\PYG{p}{)}
\PYG{n}{mx}\PYG{o}{.}\PYG{n}{setInfo}\PYG{p}{(}\PYG{l+s+s2}{\PYGZdq{}}\PYG{l+s+s2}{ub}\PYG{l+s+s2}{\PYGZdq{}}\PYG{p}{,} \PYG{l+m+mf}{9.0}\PYG{p}{)}
\PYG{n}{mx}\PYG{o}{.}\PYG{n}{setInfo}\PYG{p}{(}\PYG{n}{COPT}\PYG{o}{.}\PYG{n}{Info}\PYG{o}{.}\PYG{n}{LB}\PYG{p}{,} \PYG{l+m+mf}{0.0}\PYG{p}{)}
\end{sphinxVerbatim}
\end{quote}

\subsubsection{MVar.sum()}
\label{\detokenize{pyapiref:mvar-sum}}\begin{quote}

\sphinxAtStartPar
\sphinxstylestrong{Synopsis}
\begin{quote}

\sphinxAtStartPar
\sphinxcode{\sphinxupquote{sum(axis=None)}}
\end{quote}

\sphinxAtStartPar
\sphinxstylestrong{Description}
\begin{quote}

\sphinxAtStartPar
Sum the variables in the MVar, returning a new {\hyperref[\detokenize{pyapiref:chappyapi-mlinexpr}]{\sphinxcrossref{\DUrole{std,std-ref}{MLinExpr Class}}}} object.
\end{quote}

\sphinxAtStartPar
\sphinxstylestrong{Arguments}
\begin{quote}

\sphinxAtStartPar
\sphinxcode{\sphinxupquote{axis}}
\begin{quote}

\sphinxAtStartPar
Optional integer parameter, the default value is None, that is, to sum up variables one by one. Otherwise, sum over the given axis.
\end{quote}
\end{quote}

\sphinxAtStartPar
\sphinxstylestrong{Return value}
\begin{quote}

\sphinxAtStartPar
Returns an MLinExpr object representing the sum of the corresponding variables.
\end{quote}

\sphinxAtStartPar
\sphinxstylestrong{Example}

\begin{sphinxVerbatim}[commandchars=\\\{\}]
\PYG{n}{mx} \PYG{o}{=} \PYG{n}{model}\PYG{o}{.}\PYG{n}{addMVar}\PYG{p}{(}\PYG{p}{(}\PYG{l+m+mi}{3}\PYG{p}{,} \PYG{l+m+mi}{5}\PYG{p}{)}\PYG{p}{)}
\PYG{n}{sum\PYGZus{}all} \PYG{o}{=} \PYG{n}{mx}\PYG{o}{.}\PYG{n}{sum}\PYG{p}{(}\PYG{p}{)} \PYG{c+c1}{\PYGZsh{}Return 0\PYGZhy{}dimensional MLinExpr object}
\PYG{n}{sum\PYGZus{}row} \PYG{o}{=} \PYG{n}{mx}\PYG{o}{.}\PYG{n}{sum}\PYG{p}{(}\PYG{n}{axis} \PYG{o}{=} \PYG{l+m+mi}{0}\PYG{p}{)} \PYG{c+c1}{\PYGZsh{}Return a 1\PYGZhy{}dimensional MLinExpr object with a shape of (5, )}
\end{sphinxVerbatim}
\end{quote}

\subsubsection{MVar.tolist()}
\label{\detokenize{pyapiref:mvar-tolist}}\begin{quote}

\sphinxAtStartPar
\sphinxstylestrong{Synopsis}
\begin{quote}

\sphinxAtStartPar
\sphinxcode{\sphinxupquote{tolist()}}
\end{quote}

\sphinxAtStartPar
\sphinxstylestrong{Description}
\begin{quote}

\sphinxAtStartPar
Convert an MVar object to a one\sphinxhyphen{}dimensional list of Var objects.
\end{quote}

\sphinxAtStartPar
\sphinxstylestrong{Return value}
\begin{quote}

\sphinxAtStartPar
Returns a one\sphinxhyphen{}dimensional list containing Var objects.
\end{quote}

\sphinxAtStartPar
\sphinxstylestrong{Example}

\begin{sphinxVerbatim}[commandchars=\\\{\}]
\PYG{n}{mx} \PYG{o}{=} \PYG{n}{model}\PYG{o}{.}\PYG{n}{addMVar}\PYG{p}{(}\PYG{p}{(}\PYG{l+m+mi}{3}\PYG{p}{,} \PYG{l+m+mi}{5}\PYG{p}{)}\PYG{p}{)}
\PYG{n+nb}{print}\PYG{p}{(}\PYG{n}{mx}\PYG{o}{.}\PYG{n}{tolist}\PYG{p}{(}\PYG{p}{)}\PYG{p}{)}
\end{sphinxVerbatim}
\end{quote}

\subsubsection{MVar.transpose()}
\label{\detokenize{pyapiref:mvar-transpose}}\begin{quote}

\sphinxAtStartPar
\sphinxstylestrong{Synopsis}
\begin{quote}

\sphinxAtStartPar
\sphinxcode{\sphinxupquote{transpose()}}
\end{quote}

\sphinxAtStartPar
\sphinxstylestrong{Description}
\begin{quote}

\sphinxAtStartPar
Generates a new MVar object that is the transpose of the original MVar object.
\end{quote}

\sphinxAtStartPar
\sphinxstylestrong{Return value}
\begin{quote}

\sphinxAtStartPar
Return the new MVar object.
\end{quote}

\sphinxAtStartPar
\sphinxstylestrong{Example}

\begin{sphinxVerbatim}[commandchars=\\\{\}]
\PYG{n}{mx} \PYG{o}{=} \PYG{n}{model}\PYG{o}{.}\PYG{n}{addMVar}\PYG{p}{(}\PYG{p}{(}\PYG{l+m+mi}{3}\PYG{p}{,} \PYG{l+m+mi}{5}\PYG{p}{)}\PYG{p}{)}
\PYG{n+nb}{print}\PYG{p}{(}\PYG{n}{mx}\PYG{o}{.}\PYG{n}{transpose}\PYG{p}{(}\PYG{p}{)}\PYG{o}{.}\PYG{n}{shape}\PYG{p}{)} \PYG{c+c1}{\PYGZsh{}its shape is (5, 3)}
\end{sphinxVerbatim}
\end{quote}

\subsubsection{MVar.ndim}
\label{\detokenize{pyapiref:mvar-ndim}}\begin{quote}

\sphinxAtStartPar
\sphinxstylestrong{Synopsis}
\begin{quote}

\sphinxAtStartPar
\sphinxcode{\sphinxupquote{ndim}}
\end{quote}

\sphinxAtStartPar
\sphinxstylestrong{Description}
\begin{quote}

\sphinxAtStartPar
Dimensions of the MVar object.
\end{quote}

\sphinxAtStartPar
\sphinxstylestrong{Return value}
\begin{quote}

\sphinxAtStartPar
Integer value.
\end{quote}

\sphinxAtStartPar
\sphinxstylestrong{Example}

\begin{sphinxVerbatim}[commandchars=\\\{\}]
\PYG{n}{mx} \PYG{o}{=} \PYG{n}{model}\PYG{o}{.}\PYG{n}{addMVar}\PYG{p}{(}\PYG{p}{(}\PYG{l+m+mi}{3}\PYG{p}{,} \PYG{l+m+mi}{5}\PYG{p}{)}\PYG{p}{)}
\PYG{n+nb}{print}\PYG{p}{(}\PYG{n}{mx}\PYG{o}{.}\PYG{n}{ndim}\PYG{p}{)} \PYG{c+c1}{\PYGZsh{}ndim = 2}
\end{sphinxVerbatim}
\end{quote}

\subsubsection{MVar.shape}
\label{\detokenize{pyapiref:mvar-shape}}\begin{quote}

\sphinxAtStartPar
\sphinxstylestrong{Synopsis}
\begin{quote}

\sphinxAtStartPar
\sphinxcode{\sphinxupquote{shape}}
\end{quote}

\sphinxAtStartPar
\sphinxstylestrong{Description}
\begin{quote}

\sphinxAtStartPar
The shape of the MVar object.
\end{quote}

\sphinxAtStartPar
\sphinxstylestrong{Return value}
\begin{quote}

\sphinxAtStartPar
Integer tuple.
\end{quote}

\sphinxAtStartPar
\sphinxstylestrong{Example}

\begin{sphinxVerbatim}[commandchars=\\\{\}]
\PYG{n}{mx} \PYG{o}{=} \PYG{n}{model}\PYG{o}{.}\PYG{n}{addMVar}\PYG{p}{(}\PYG{p}{(}\PYG{l+m+mi}{3}\PYG{p}{,}\PYG{p}{)}\PYG{p}{)}
\PYG{n+nb}{print}\PYG{p}{(}\PYG{n}{mx}\PYG{o}{.}\PYG{n}{shape}\PYG{p}{)} \PYG{c+c1}{\PYGZsh{} shape = (3, )}
\end{sphinxVerbatim}
\end{quote}

\subsubsection{MVar.size}
\label{\detokenize{pyapiref:mvar-size}}\begin{quote}

\sphinxAtStartPar
\sphinxstylestrong{Synopsis}
\begin{quote}

\sphinxAtStartPar
\sphinxcode{\sphinxupquote{size}}
\end{quote}

\sphinxAtStartPar
\sphinxstylestrong{Description}
\begin{quote}

\sphinxAtStartPar
The number of Var variables in the MVar object.
\end{quote}

\sphinxAtStartPar
\sphinxstylestrong{Return value}
\begin{quote}

\sphinxAtStartPar
Integer value.
\end{quote}

\sphinxAtStartPar
\sphinxstylestrong{Example}

\begin{sphinxVerbatim}[commandchars=\\\{\}]
\PYG{n}{mx} \PYG{o}{=} \PYG{n}{model}\PYG{o}{.}\PYG{n}{addMVar}\PYG{p}{(}\PYG{p}{(}\PYG{l+m+mi}{3}\PYG{p}{,} \PYG{l+m+mi}{4}\PYG{p}{)}\PYG{p}{)}
\PYG{n+nb}{print}\PYG{p}{(}\PYG{n}{mx}\PYG{o}{.}\PYG{n}{size}\PYG{p}{)} \PYG{c+c1}{\PYGZsh{} size = 12}
\end{sphinxVerbatim}
\end{quote}

\subsubsection{MVar.T}
\label{\detokenize{pyapiref:mvar-t}}\begin{quote}

\sphinxAtStartPar
\sphinxstylestrong{Synopsis}
\begin{quote}

\sphinxAtStartPar
\sphinxcode{\sphinxupquote{T}}
\end{quote}

\sphinxAtStartPar
\sphinxstylestrong{Description}
\begin{quote}

\sphinxAtStartPar
Transpose of the MVar object. Similar to the class method transpose().
\end{quote}

\sphinxAtStartPar
\sphinxstylestrong{Return value}
\begin{quote}

\sphinxAtStartPar
Returns the transposed MVar object.
\end{quote}

\sphinxAtStartPar
\sphinxstylestrong{Example}

\begin{sphinxVerbatim}[commandchars=\\\{\}]
\PYG{n}{mx} \PYG{o}{=} \PYG{n}{model}\PYG{o}{.}\PYG{n}{addMVar}\PYG{p}{(}\PYG{p}{(}\PYG{l+m+mi}{3}\PYG{p}{,} \PYG{l+m+mi}{4}\PYG{p}{)}\PYG{p}{)}
\PYG{n+nb}{print}\PYG{p}{(}\PYG{n}{mx}\PYG{o}{.}\PYG{n}{T}\PYG{o}{.}\PYG{n}{shape}\PYG{p}{)} \PYG{c+c1}{\PYGZsh{} shape = (4, 3)}
\end{sphinxVerbatim}
\end{quote}

\subsection{MConstr Class}
\label{\detokenize{pyapiref:mconstr-class}}\label{\detokenize{pyapiref:chappyapi-mconstr}}
\sphinxAtStartPar
The MConstr class holds multi\sphinxhyphen{}dimensional linear constraints in COPT and supports NumPy’s multi\sphinxhyphen{}dimensional array operations.
It is generated by the method \sphinxcode{\sphinxupquote{addConstrs}} or \sphinxcode{\sphinxupquote{addMConstr}} of the model class.
It supports a unified set of indexing and slicing rules.
For detailed usage, see {\hyperref[\detokenize{matrix:chapmatrixindex}]{\sphinxcrossref{\DUrole{std,std-ref}{Matrix Modeling: Indexing and Slicing}}}}.

\sphinxAtStartPar
The following member methods are provided:

\subsubsection{MConstr.getInfo()}
\label{\detokenize{pyapiref:mconstr-getinfo}}\begin{quote}

\sphinxAtStartPar
\sphinxstylestrong{Synopsis}
\begin{quote}

\sphinxAtStartPar
\sphinxcode{\sphinxupquote{getInfo(infoname)}}
\end{quote}

\sphinxAtStartPar
\sphinxstylestrong{Description}
\begin{quote}

\sphinxAtStartPar
Get the information value of each constraint within MConstr.
\end{quote}

\sphinxAtStartPar
\sphinxstylestrong{Arguments}
\begin{quote}

\sphinxAtStartPar
\sphinxcode{\sphinxupquote{infoname}}
\begin{quote}

\sphinxAtStartPar
The name of the information being queried.
Please refer to {\hyperref[\detokenize{information:chapinfo}]{\sphinxcrossref{\DUrole{std,std-ref}{Information Section}}}} for possible values.
\end{quote}
\end{quote}

\sphinxAtStartPar
\sphinxstylestrong{Return value}
\begin{quote}

\sphinxAtStartPar
Returns a NumPy ndarray with the same dimension as the MConstr object, whose elements are the attribute values of the corresponding constraint.
\end{quote}

\sphinxAtStartPar
\sphinxstylestrong{Example}

\begin{sphinxVerbatim}[commandchars=\\\{\}]
\PYG{n}{a} \PYG{o}{=} \PYG{n}{np}\PYG{o}{.}\PYG{n}{random}\PYG{o}{.}\PYG{n}{rand}\PYG{p}{(}\PYG{l+m+mi}{4}\PYG{p}{)}
\PYG{n}{mx} \PYG{o}{=} \PYG{n}{m}\PYG{o}{.}\PYG{n}{addMVar}\PYG{p}{(}\PYG{p}{(}\PYG{l+m+mi}{4}\PYG{p}{,} \PYG{l+m+mi}{3}\PYG{p}{)}\PYG{p}{,} \PYG{n}{nameprefix}\PYG{o}{=}\PYG{l+s+s2}{\PYGZdq{}}\PYG{l+s+s2}{mx}\PYG{l+s+s2}{\PYGZdq{}}\PYG{p}{)}
\PYG{n}{b} \PYG{o}{=} \PYG{n}{np}\PYG{o}{.}\PYG{n}{random}\PYG{o}{.}\PYG{n}{rand}\PYG{p}{(}\PYG{l+m+mi}{3}\PYG{p}{)}
\PYG{n}{mc} \PYG{o}{=} \PYG{n}{m}\PYG{o}{.}\PYG{n}{addConstrs}\PYG{p}{(}\PYG{n}{a} \PYG{o}{@} \PYG{n}{mx} \PYG{o}{\PYGZlt{}}\PYG{o}{=} \PYG{n}{b}\PYG{p}{)}
\PYG{n+nb}{print}\PYG{p}{(}\PYG{n}{mc}\PYG{o}{.}\PYG{n}{getInfo}\PYG{p}{(}\PYG{l+s+s2}{\PYGZdq{}}\PYG{l+s+s2}{pi}\PYG{l+s+s2}{\PYGZdq{}}\PYG{p}{)}\PYG{p}{)}
\end{sphinxVerbatim}
\end{quote}

\subsubsection{MConstr.item()}
\label{\detokenize{pyapiref:mconstr-item}}\begin{quote}

\sphinxAtStartPar
\sphinxstylestrong{Synopsis}
\begin{quote}

\sphinxAtStartPar
\sphinxcode{\sphinxupquote{item()}}
\end{quote}

\sphinxAtStartPar
\sphinxstylestrong{Description}
\begin{quote}

\sphinxAtStartPar
Get the constraint object in 0\sphinxhyphen{}dimensional MConstr. If the MConstr object is not 0\sphinxhyphen{}dimensional, a ValueError exception is raised.
\end{quote}

\sphinxAtStartPar
\sphinxstylestrong{Return value}
\begin{quote}

\sphinxAtStartPar
Returns the linear constraint object.
\end{quote}

\sphinxAtStartPar
\sphinxstylestrong{Example}

\begin{sphinxVerbatim}[commandchars=\\\{\}]
\PYG{n}{mc} \PYG{o}{=} \PYG{n}{m}\PYG{o}{.}\PYG{n}{addConstrs}\PYG{p}{(}\PYG{n}{a} \PYG{o}{@} \PYG{n}{mx} \PYG{o}{\PYGZlt{}}\PYG{o}{=} \PYG{n}{b}\PYG{p}{)}
\PYG{n+nb}{print}\PYG{p}{(}\PYG{n}{mc}\PYG{p}{[}\PYG{l+m+mi}{0}\PYG{p}{]}\PYG{o}{.}\PYG{n}{item}\PYG{p}{(}\PYG{p}{)}\PYG{p}{)}
\end{sphinxVerbatim}
\end{quote}

\subsubsection{MConstr.reshape()}
\label{\detokenize{pyapiref:mconstr-reshape}}\begin{quote}

\sphinxAtStartPar
\sphinxstylestrong{Synopsis}
\begin{quote}

\sphinxAtStartPar
\sphinxcode{\sphinxupquote{reshape(shape, order=\textquotesingle{}C\textquotesingle{})}}
\end{quote}

\sphinxAtStartPar
\sphinxstylestrong{Description}
\begin{quote}

\sphinxAtStartPar
Returns a new MConstr object whose elements remain unchanged but whose shape is transformed by the argument shape.
\end{quote}

\sphinxAtStartPar
\sphinxstylestrong{Arguments}
\begin{quote}

\sphinxAtStartPar
\sphinxcode{\sphinxupquote{shape}}
\begin{quote}

\sphinxAtStartPar
The value is an integer, or a tuple of integers. which represents the shape of the new MConstr object.
\end{quote}

\sphinxAtStartPar
\sphinxcode{\sphinxupquote{order}}
\begin{quote}

\sphinxAtStartPar
Optional parameter, the default is the character ‘C’, which means it is compatible with the C language,
that is, it is stored in rows; it can also be set to the character ‘F’, that is, it is stored in columns,
and it is compatible with the Fortune language.
\end{quote}
\end{quote}

\sphinxAtStartPar
\sphinxstylestrong{Return value}
\begin{quote}

\sphinxAtStartPar
Returns a new MConstr object with the same elements as the original MConstr object but with a different shape.
\end{quote}

\sphinxAtStartPar
\sphinxstylestrong{Example}

\begin{sphinxVerbatim}[commandchars=\\\{\}]
\PYG{n}{mc} \PYG{o}{=} \PYG{n}{m}\PYG{o}{.}\PYG{n}{addConstrs}\PYG{p}{(}\PYG{n}{a} \PYG{o}{@} \PYG{n}{mx} \PYG{o}{\PYGZlt{}}\PYG{o}{=} \PYG{n}{b}\PYG{p}{)}
\PYG{n}{mc\PYGZus{}2x2} \PYG{o}{=} \PYG{n}{mc}\PYG{o}{.}\PYG{n}{reshape}\PYG{p}{(}\PYG{p}{(}\PYG{l+m+mi}{2}\PYG{p}{,} \PYG{l+m+mi}{2}\PYG{p}{)}\PYG{p}{)}
\end{sphinxVerbatim}
\end{quote}

\subsubsection{MConstr.setInfo()}
\label{\detokenize{pyapiref:mconstr-setinfo}}\begin{quote}

\sphinxAtStartPar
\sphinxstylestrong{Synopsis}
\begin{quote}

\sphinxAtStartPar
\sphinxcode{\sphinxupquote{setInfo(infoname, newval)}}
\end{quote}

\sphinxAtStartPar
\sphinxstylestrong{Description}
\begin{quote}

\sphinxAtStartPar
Set information values for each constraint within MConstr.
\end{quote}

\sphinxAtStartPar
\sphinxstylestrong{Arguments}
\begin{quote}

\sphinxAtStartPar
\sphinxcode{\sphinxupquote{infoname}}
\begin{quote}

\sphinxAtStartPar
The name of the information to be set.
Please refer to {\hyperref[\detokenize{information:chapinfo}]{\sphinxcrossref{\DUrole{std,std-ref}{Information Section}}}} for possible values.
\end{quote}

\sphinxAtStartPar
\sphinxcode{\sphinxupquote{newval}}
\begin{quote}

\sphinxAtStartPar
The new value to be set.
\end{quote}
\end{quote}

\sphinxAtStartPar
\sphinxstylestrong{Example}

\begin{sphinxVerbatim}[commandchars=\\\{\}]
\PYG{n}{mc} \PYG{o}{=} \PYG{n}{model}\PYG{o}{.}\PYG{n}{addConstrs}\PYG{p}{(}\PYG{n}{a} \PYG{o}{@} \PYG{n}{mx} \PYG{o}{\PYGZlt{}}\PYG{o}{=} \PYG{n}{b}\PYG{p}{)}
\PYG{n}{mc}\PYG{o}{.}\PYG{n}{setInfo}\PYG{p}{(}\PYG{l+s+s2}{\PYGZdq{}}\PYG{l+s+s2}{obj}\PYG{l+s+s2}{\PYGZdq{}}\PYG{p}{,} \PYG{l+m+mf}{9.0}\PYG{p}{)}
\end{sphinxVerbatim}
\end{quote}

\subsubsection{MConstr.tolist()}
\label{\detokenize{pyapiref:mconstr-tolist}}\begin{quote}

\sphinxAtStartPar
\sphinxstylestrong{Synopsis}
\begin{quote}

\sphinxAtStartPar
\sphinxcode{\sphinxupquote{tolist()}}
\end{quote}

\sphinxAtStartPar
\sphinxstylestrong{Description}
\begin{quote}

\sphinxAtStartPar
Convert the MConstr object to a one\sphinxhyphen{}dimensional list of constraints.
\end{quote}

\sphinxAtStartPar
\sphinxstylestrong{Return value}
\begin{quote}

\sphinxAtStartPar
Returns a one\sphinxhyphen{}dimensional list containing Constraint objects.
\end{quote}

\sphinxAtStartPar
\sphinxstylestrong{Example}

\begin{sphinxVerbatim}[commandchars=\\\{\}]
\PYG{n}{mc} \PYG{o}{=} \PYG{n}{m}\PYG{o}{.}\PYG{n}{addConstrs}\PYG{p}{(}\PYG{n}{a} \PYG{o}{@} \PYG{n}{mx} \PYG{o}{\PYGZlt{}}\PYG{o}{=} \PYG{n}{b}\PYG{p}{)}
\PYG{n+nb}{print}\PYG{p}{(}\PYG{n}{mc}\PYG{o}{.}\PYG{n}{tolist}\PYG{p}{(}\PYG{p}{)}\PYG{p}{)}
\end{sphinxVerbatim}
\end{quote}

\subsubsection{MConstr.transpose()}
\label{\detokenize{pyapiref:mconstr-transpose}}\begin{quote}

\sphinxAtStartPar
\sphinxstylestrong{Synopsis}
\begin{quote}

\sphinxAtStartPar
\sphinxcode{\sphinxupquote{transpose()}}
\end{quote}

\sphinxAtStartPar
\sphinxstylestrong{Description}
\begin{quote}

\sphinxAtStartPar
Generates a new MConstr object that is the transpose of the original MConstr object.
\end{quote}

\sphinxAtStartPar
\sphinxstylestrong{Return value}
\begin{quote}

\sphinxAtStartPar
Returns the transposed MConstr object.
\end{quote}

\sphinxAtStartPar
\sphinxstylestrong{Example}

\begin{sphinxVerbatim}[commandchars=\\\{\}]
\PYG{n}{mc} \PYG{o}{=} \PYG{n}{m}\PYG{o}{.}\PYG{n}{addConstrs}\PYG{p}{(}\PYG{n}{a} \PYG{o}{@} \PYG{n}{mx} \PYG{o}{\PYGZlt{}}\PYG{o}{=} \PYG{n}{b}\PYG{p}{)}
\PYG{n+nb}{print}\PYG{p}{(}\PYG{n}{mc}\PYG{o}{.}\PYG{n}{transpose}\PYG{p}{(}\PYG{p}{)}\PYG{p}{)}
\end{sphinxVerbatim}
\end{quote}

\subsubsection{MConstr.ndim}
\label{\detokenize{pyapiref:mconstr-ndim}}\begin{quote}

\sphinxAtStartPar
\sphinxstylestrong{Synopsis}
\begin{quote}

\sphinxAtStartPar
\sphinxcode{\sphinxupquote{ndim}}
\end{quote}

\sphinxAtStartPar
\sphinxstylestrong{Description}
\begin{quote}

\sphinxAtStartPar
Dimensions of the MConstr object.
\end{quote}

\sphinxAtStartPar
\sphinxstylestrong{Return value}
\begin{quote}

\sphinxAtStartPar
Integer value.
\end{quote}

\sphinxAtStartPar
\sphinxstylestrong{Example}

\begin{sphinxVerbatim}[commandchars=\\\{\}]
\PYG{n}{mc} \PYG{o}{=} \PYG{n}{m}\PYG{o}{.}\PYG{n}{addConstrs}\PYG{p}{(}\PYG{n}{a} \PYG{o}{@} \PYG{n}{mx} \PYG{o}{\PYGZlt{}}\PYG{o}{=} \PYG{n}{b}\PYG{p}{)}
\PYG{n+nb}{print}\PYG{p}{(}\PYG{n}{mc}\PYG{o}{.}\PYG{n}{ndim}\PYG{p}{)}
\end{sphinxVerbatim}
\end{quote}

\subsubsection{MConstr.shape}
\label{\detokenize{pyapiref:mconstr-shape}}\begin{quote}

\sphinxAtStartPar
\sphinxstylestrong{Synopsis}
\begin{quote}

\sphinxAtStartPar
\sphinxcode{\sphinxupquote{shape}}
\end{quote}

\sphinxAtStartPar
\sphinxstylestrong{Description}
\begin{quote}

\sphinxAtStartPar
The shape of the MConstr object.
\end{quote}

\sphinxAtStartPar
\sphinxstylestrong{Return value}
\begin{quote}

\sphinxAtStartPar
Integer tuple.
\end{quote}

\sphinxAtStartPar
\sphinxstylestrong{Example}

\begin{sphinxVerbatim}[commandchars=\\\{\}]
\PYG{n}{mc} \PYG{o}{=} \PYG{n}{m}\PYG{o}{.}\PYG{n}{addConstrs}\PYG{p}{(}\PYG{n}{a} \PYG{o}{@} \PYG{n}{mx} \PYG{o}{\PYGZlt{}}\PYG{o}{=} \PYG{n}{b}\PYG{p}{)}
\PYG{n+nb}{print}\PYG{p}{(}\PYG{n}{mc}\PYG{o}{.}\PYG{n}{shape}\PYG{p}{)}
\end{sphinxVerbatim}
\end{quote}

\subsubsection{MConstr.size}
\label{\detokenize{pyapiref:mconstr-size}}\begin{quote}

\sphinxAtStartPar
\sphinxstylestrong{Synopsis}
\begin{quote}

\sphinxAtStartPar
\sphinxcode{\sphinxupquote{size}}
\end{quote}

\sphinxAtStartPar
\sphinxstylestrong{Description}
\begin{quote}

\sphinxAtStartPar
The number of constraints of the MConstr object.
\end{quote}

\sphinxAtStartPar
\sphinxstylestrong{Return value}
\begin{quote}

\sphinxAtStartPar
Integer value.
\end{quote}

\sphinxAtStartPar
\sphinxstylestrong{Example}

\begin{sphinxVerbatim}[commandchars=\\\{\}]
\PYG{n}{mc} \PYG{o}{=} \PYG{n}{m}\PYG{o}{.}\PYG{n}{addConstrs}\PYG{p}{(}\PYG{n}{a} \PYG{o}{@} \PYG{n}{mx} \PYG{o}{\PYGZlt{}}\PYG{o}{=} \PYG{n}{b}\PYG{p}{)}
\PYG{n+nb}{print}\PYG{p}{(}\PYG{n}{mc}\PYG{o}{.}\PYG{n}{size}\PYG{p}{)}
\end{sphinxVerbatim}
\end{quote}

\subsubsection{MConstr.T}
\label{\detokenize{pyapiref:mconstr-t}}\begin{quote}

\sphinxAtStartPar
\sphinxstylestrong{Synopsis}
\begin{quote}

\sphinxAtStartPar
\sphinxcode{\sphinxupquote{T}}
\end{quote}

\sphinxAtStartPar
\sphinxstylestrong{Description}
\begin{quote}

\sphinxAtStartPar
Transpose of the MConstr object. Similar to the class method transpose().
\end{quote}

\sphinxAtStartPar
\sphinxstylestrong{Return value}
\begin{quote}

\sphinxAtStartPar
Returns the transposed MConstr object.
\end{quote}

\sphinxAtStartPar
\sphinxstylestrong{Example}

\begin{sphinxVerbatim}[commandchars=\\\{\}]
\PYG{n}{A} \PYG{o}{=} \PYG{n}{np}\PYG{o}{.}\PYG{n}{ones}\PYG{p}{(}\PYG{p}{[}\PYG{l+m+mi}{2}\PYG{p}{,} \PYG{l+m+mi}{4}\PYG{p}{]}\PYG{p}{)}
\PYG{n}{mx} \PYG{o}{=} \PYG{n}{m}\PYG{o}{.}\PYG{n}{addMVar}\PYG{p}{(}\PYG{p}{(}\PYG{l+m+mi}{4}\PYG{p}{,} \PYG{l+m+mi}{3}\PYG{p}{)}\PYG{p}{,} \PYG{n}{nameprefix}\PYG{o}{=}\PYG{l+s+s2}{\PYGZdq{}}\PYG{l+s+s2}{mx}\PYG{l+s+s2}{\PYGZdq{}}\PYG{p}{)}
\PYG{n}{mc} \PYG{o}{=} \PYG{n}{m}\PYG{o}{.}\PYG{n}{addConstrs}\PYG{p}{(}\PYG{n}{A} \PYG{o}{@} \PYG{n}{X} \PYG{o}{==} \PYG{l+m+mf}{0.0}\PYG{p}{)}
\PYG{n+nb}{print}\PYG{p}{(}\PYG{n}{mc}\PYG{o}{.}\PYG{n}{T}\PYG{o}{.}\PYG{n}{shape}\PYG{p}{)} \PYG{c+c1}{\PYGZsh{} shape = (3, 2)}
\end{sphinxVerbatim}
\end{quote}

\subsection{MConstrBuilder Class}
\label{\detokenize{pyapiref:mconstrbuilder-class}}\label{\detokenize{pyapiref:chappyapi-mconstrbuilder}}
\sphinxAtStartPar
The MConstrBuilder class is used to build multi\sphinxhyphen{}dimensional linear constraints in COPT
and supports NumPy’s multi\sphinxhyphen{}dimensional array operations.
Users might create a MConstrBuilder object by its constructor with a list of
{\hyperref[\detokenize{pyapiref:chappyapi-constrbuilder}]{\sphinxcrossref{\DUrole{std,std-ref}{ConstrBuilder Class}}}} objects, or simply by overloaded comparison operator
of {\hyperref[\detokenize{pyapiref:chappyapi-mlinexpr}]{\sphinxcrossref{\DUrole{std,std-ref}{MLinExpr Class}}}}. The following member methods are provided:

\subsubsection{MConstrBuilder()}
\label{\detokenize{pyapiref:mconstrbuilder}}\begin{quote}

\sphinxAtStartPar
\sphinxstylestrong{Synopsis}
\begin{quote}

\sphinxAtStartPar
\sphinxcode{\sphinxupquote{MConstrBuilder(args, shape=None)}}
\end{quote}

\sphinxAtStartPar
\sphinxstylestrong{Description}
\begin{quote}

\sphinxAtStartPar
construtor of MConstrBuilder.
\end{quote}

\sphinxAtStartPar
\sphinxstylestrong{Arguments}
\begin{quote}

\sphinxAtStartPar
\sphinxcode{\sphinxupquote{args}}
\begin{quote}

\sphinxAtStartPar
one or a set of {\hyperref[\detokenize{pyapiref:chappyapi-constrbuilder}]{\sphinxcrossref{\DUrole{std,std-ref}{ConstrBuilder Class}}}} objects,
in form of Python list or NumPy ndarray.
\end{quote}

\sphinxAtStartPar
\sphinxcode{\sphinxupquote{shape}}
\begin{quote}

\sphinxAtStartPar
an integer, or tuple of integers, which is the shape of new MConstrBuilder object.
\end{quote}
\end{quote}

\sphinxAtStartPar
\sphinxstylestrong{Example}

\begin{sphinxVerbatim}[commandchars=\\\{\}]
\PYG{n+nb}{vars} \PYG{o}{=} \PYG{n}{m}\PYG{o}{.}\PYG{n}{addVars}\PYG{p}{(}\PYG{l+m+mi}{4}\PYG{p}{)}
\PYG{n}{builders} \PYG{o}{=} \PYG{p}{[}\PYG{n}{x} \PYG{o}{\PYGZlt{}}\PYG{o}{=} \PYG{l+m+mf}{1.0} \PYG{k}{for} \PYG{n}{x} \PYG{o+ow}{in} \PYG{n+nb}{vars}\PYG{p}{]}
\PYG{n}{mcb} \PYG{o}{=} \PYG{n}{MConstrBuilder}\PYG{p}{(}\PYG{n}{builders}\PYG{p}{,} \PYG{p}{(}\PYG{l+m+mi}{2}\PYG{p}{,} \PYG{l+m+mi}{2}\PYG{p}{)}\PYG{p}{)}

\PYG{c+c1}{\PYGZsh{} or by overloaded comparison operator of MVar}
\PYG{n}{mx} \PYG{o}{=} \PYG{n}{m}\PYG{o}{.}\PYG{n}{addMVar}\PYG{p}{(}\PYG{p}{(}\PYG{l+m+mi}{3}\PYG{p}{,} \PYG{l+m+mi}{2}\PYG{p}{)}\PYG{p}{)}
\PYG{n}{mcb} \PYG{o}{=} \PYG{n}{mx} \PYG{o}{\PYGZgt{}}\PYG{o}{=} \PYG{l+m+mf}{1.0}

\PYG{c+c1}{\PYGZsh{} or immediately passed to addConstrs()}
\PYG{n}{model}\PYG{o}{.}\PYG{n}{addConstrs}\PYG{p}{(}\PYG{n}{mx} \PYG{o}{\PYGZgt{}}\PYG{o}{=} \PYG{l+m+mf}{1.0}\PYG{p}{)}
\end{sphinxVerbatim}
\end{quote}

\subsection{MQConstr Class}
\label{\detokenize{pyapiref:mqconstr-class}}\label{\detokenize{pyapiref:chappyapi-mqconstr}}
\sphinxAtStartPar
The MQConstr class holds multi\sphinxhyphen{}dimensional quadratic constraints in COPT and
supports NumPy’s multi\sphinxhyphen{}dimensional array operations. In experimental version of
matrix modeling, it is generated by the method \sphinxcode{\sphinxupquote{addQConstr}} or
\sphinxcode{\sphinxupquote{addMQConstr}} of the model class.
It supports a unified set of indexing and slicing rules.
For detailed usage, see {\hyperref[\detokenize{matrix:chapmatrixindex}]{\sphinxcrossref{\DUrole{std,std-ref}{Matrix Modeling: Indexing and Slicing}}}}.

\sphinxAtStartPar
The following member methods are provided:

\subsubsection{MQConstr.getInfo()}
\label{\detokenize{pyapiref:mqconstr-getinfo}}\begin{quote}

\sphinxAtStartPar
\sphinxstylestrong{Synopsis}
\begin{quote}

\sphinxAtStartPar
\sphinxcode{\sphinxupquote{getInfo(infoname)}}
\end{quote}

\sphinxAtStartPar
\sphinxstylestrong{Description}
\begin{quote}

\sphinxAtStartPar
Get the information value of each quadratic constraint within MQConstr.
\end{quote}

\sphinxAtStartPar
\sphinxstylestrong{Arguments}
\begin{quote}

\sphinxAtStartPar
\sphinxcode{\sphinxupquote{infoname}}
\begin{quote}

\sphinxAtStartPar
The name of the information being queried.
Please refer to {\hyperref[\detokenize{information:chapinfo}]{\sphinxcrossref{\DUrole{std,std-ref}{Information Section}}}} for possible values.
\end{quote}
\end{quote}

\sphinxAtStartPar
\sphinxstylestrong{Return value}
\begin{quote}

\sphinxAtStartPar
Returns a built\sphinxhyphen{}in NdArray with the same shape as the MQConstr object,
whose elements are the information values of the corresponding
constraints.
\end{quote}

\sphinxAtStartPar
\sphinxstylestrong{Example}

\begin{sphinxVerbatim}[commandchars=\\\{\}]
\PYG{n}{mx} \PYG{o}{=} \PYG{n}{m}\PYG{o}{.}\PYG{n}{addMVar}\PYG{p}{(}\PYG{p}{(}\PYG{l+m+mi}{1}\PYG{p}{,} \PYG{l+m+mi}{3}\PYG{p}{)}\PYG{p}{,} \PYG{n}{nameprefix}\PYG{o}{=}\PYG{l+s+s2}{\PYGZdq{}}\PYG{l+s+s2}{mx}\PYG{l+s+s2}{\PYGZdq{}}\PYG{p}{)}
\PYG{n}{mc} \PYG{o}{=} \PYG{n}{m}\PYG{o}{.}\PYG{n}{addQConstr}\PYG{p}{(}\PYG{n}{mx} \PYG{o}{@} \PYG{n}{mx}\PYG{o}{.}\PYG{n}{T} \PYG{o}{\PYGZlt{}}\PYG{o}{=} \PYG{l+m+mf}{1.0}\PYG{p}{)}
\PYG{n+nb}{print}\PYG{p}{(}\PYG{n}{mc}\PYG{o}{.}\PYG{n}{getInfo}\PYG{p}{(}\PYG{l+s+s2}{\PYGZdq{}}\PYG{l+s+s2}{x}\PYG{l+s+s2}{\PYGZdq{}}\PYG{p}{)}\PYG{p}{)}
\end{sphinxVerbatim}
\end{quote}

\subsubsection{MQConstr.item()}
\label{\detokenize{pyapiref:mqconstr-item}}\begin{quote}

\sphinxAtStartPar
\sphinxstylestrong{Synopsis}
\begin{quote}

\sphinxAtStartPar
\sphinxcode{\sphinxupquote{item()}}
\end{quote}

\sphinxAtStartPar
\sphinxstylestrong{Description}
\begin{quote}

\sphinxAtStartPar
Get the quadratic constraint object in MQConstr object.  If the
MQConstr object has more than one item, an exception of ValueError is
raised.
\end{quote}

\sphinxAtStartPar
\sphinxstylestrong{Return value}
\begin{quote}

\sphinxAtStartPar
Returns the quadratic constraint object.
\end{quote}

\sphinxAtStartPar
\sphinxstylestrong{Example}

\begin{sphinxVerbatim}[commandchars=\\\{\}]
\PYG{n}{mc} \PYG{o}{=} \PYG{n}{m}\PYG{o}{.}\PYG{n}{addQConstr}\PYG{p}{(}\PYG{n}{mx} \PYG{o}{@} \PYG{n}{mx}\PYG{o}{.}\PYG{n}{T} \PYG{o}{\PYGZlt{}}\PYG{o}{=} \PYG{l+m+mf}{1.0}\PYG{p}{)}
\PYG{n+nb}{print}\PYG{p}{(}\PYG{n}{mc}\PYG{o}{.}\PYG{n}{item}\PYG{p}{(}\PYG{p}{)}\PYG{p}{)}
\end{sphinxVerbatim}
\end{quote}

\subsubsection{MQConstr.reshape()}
\label{\detokenize{pyapiref:mqconstr-reshape}}\begin{quote}

\sphinxAtStartPar
\sphinxstylestrong{Synopsis}
\begin{quote}

\sphinxAtStartPar
\sphinxcode{\sphinxupquote{reshape(shape, order=\textquotesingle{}C\textquotesingle{})}}
\end{quote}

\sphinxAtStartPar
\sphinxstylestrong{Description}
\begin{quote}

\sphinxAtStartPar
Returns a new MQConstr object whose elements remain unchanged but whose shape is transformed by the argument shape.
\end{quote}

\sphinxAtStartPar
\sphinxstylestrong{Arguments}
\begin{quote}

\sphinxAtStartPar
\sphinxcode{\sphinxupquote{shape}}
\begin{quote}

\sphinxAtStartPar
The value is an integer, or a tuple of integers. which represents the shape of the new MQConstr object.
\end{quote}

\sphinxAtStartPar
\sphinxcode{\sphinxupquote{order}}
\begin{quote}

\sphinxAtStartPar
Optional parameter, the default is the character ‘C’, which means it is compatible with the C language.
The other option ‘F’ is not implemented yet.
\end{quote}
\end{quote}

\sphinxAtStartPar
\sphinxstylestrong{Return value}
\begin{quote}

\sphinxAtStartPar
Returns a new MQConstr object with the same elements as the original MQConstr object but with a different shape.
\end{quote}

\sphinxAtStartPar
\sphinxstylestrong{Example}

\begin{sphinxVerbatim}[commandchars=\\\{\}]
\PYG{n}{mc} \PYG{o}{=} \PYG{n}{m}\PYG{o}{.}\PYG{n}{addQConstr}\PYG{p}{(}\PYG{n}{mx}\PYG{o}{.}\PYG{n}{T} \PYG{o}{@} \PYG{n}{mx} \PYG{o}{\PYGZlt{}}\PYG{o}{=} \PYG{l+m+mf}{1.0}\PYG{p}{)}
\PYG{n}{mc\PYGZus{}1x9} \PYG{o}{=} \PYG{n}{mc}\PYG{o}{.}\PYG{n}{reshape}\PYG{p}{(}\PYG{p}{(}\PYG{l+m+mi}{1}\PYG{p}{,} \PYG{l+m+mi}{9}\PYG{p}{)}\PYG{p}{)}
\end{sphinxVerbatim}
\end{quote}

\subsubsection{MQConstr.setInfo()}
\label{\detokenize{pyapiref:mqconstr-setinfo}}\begin{quote}

\sphinxAtStartPar
\sphinxstylestrong{Synopsis}
\begin{quote}

\sphinxAtStartPar
\sphinxcode{\sphinxupquote{setInfo(infoname, newval)}}
\end{quote}

\sphinxAtStartPar
\sphinxstylestrong{Description}
\begin{quote}

\sphinxAtStartPar
Set information values for each quadratic constraint within MQConstr.
\end{quote}

\sphinxAtStartPar
\sphinxstylestrong{Arguments}
\begin{quote}

\sphinxAtStartPar
\sphinxcode{\sphinxupquote{infoname}}
\begin{quote}

\sphinxAtStartPar
The name of the information to be set.
Please refer to {\hyperref[\detokenize{information:chapinfo}]{\sphinxcrossref{\DUrole{std,std-ref}{Information Section}}}} for possible values.
\end{quote}

\sphinxAtStartPar
\sphinxcode{\sphinxupquote{newval}}
\begin{quote}

\sphinxAtStartPar
The new value to be set.
\end{quote}
\end{quote}

\sphinxAtStartPar
\sphinxstylestrong{Example}

\begin{sphinxVerbatim}[commandchars=\\\{\}]
\PYG{n}{mc} \PYG{o}{=} \PYG{n}{model}\PYG{o}{.}\PYG{n}{addQConstr}\PYG{p}{(}\PYG{n}{mx} \PYG{o}{@} \PYG{n}{mx}\PYG{o}{.}\PYG{n}{T} \PYG{o}{\PYGZlt{}}\PYG{o}{=} \PYG{l+m+mf}{1.0}\PYG{p}{)}
\PYG{n}{mc}\PYG{o}{.}\PYG{n}{setInfo}\PYG{p}{(}\PYG{l+s+s2}{\PYGZdq{}}\PYG{l+s+s2}{LB}\PYG{l+s+s2}{\PYGZdq{}}\PYG{p}{,} \PYG{l+m+mf}{9.0}\PYG{p}{)}
\end{sphinxVerbatim}
\end{quote}

\subsubsection{MQConstr.tolist()}
\label{\detokenize{pyapiref:mqconstr-tolist}}\begin{quote}

\sphinxAtStartPar
\sphinxstylestrong{Synopsis}
\begin{quote}

\sphinxAtStartPar
\sphinxcode{\sphinxupquote{tolist()}}
\end{quote}

\sphinxAtStartPar
\sphinxstylestrong{Description}
\begin{quote}

\sphinxAtStartPar
Convert the MQConstr object to a one\sphinxhyphen{}dimensional list of quadratic constraints.
\end{quote}

\sphinxAtStartPar
\sphinxstylestrong{Return value}
\begin{quote}

\sphinxAtStartPar
Returns a one\sphinxhyphen{}dimensional list containing {\hyperref[\detokenize{pyapiref:chappyapi-qconstraint}]{\sphinxcrossref{\DUrole{std,std-ref}{QConstraint Class}}}} objects.
\end{quote}

\sphinxAtStartPar
\sphinxstylestrong{Example}

\begin{sphinxVerbatim}[commandchars=\\\{\}]
\PYG{n}{mc} \PYG{o}{=} \PYG{n}{m}\PYG{o}{.}\PYG{n}{addQConstr}\PYG{p}{(}\PYG{n}{mx}\PYG{o}{.}\PYG{n}{T} \PYG{o}{@} \PYG{n}{mx} \PYG{o}{\PYGZlt{}}\PYG{o}{=} \PYG{l+m+mf}{1.0}\PYG{p}{)}
\PYG{n+nb}{print}\PYG{p}{(}\PYG{n}{mc}\PYG{o}{.}\PYG{n}{tolist}\PYG{p}{(}\PYG{p}{)}\PYG{p}{)}
\end{sphinxVerbatim}
\end{quote}

\subsubsection{MQConstr.transpose()}
\label{\detokenize{pyapiref:mqconstr-transpose}}\begin{quote}

\sphinxAtStartPar
\sphinxstylestrong{Synopsis}
\begin{quote}

\sphinxAtStartPar
\sphinxcode{\sphinxupquote{transpose()}}
\end{quote}

\sphinxAtStartPar
\sphinxstylestrong{Description}
\begin{quote}

\sphinxAtStartPar
Generates a new MQConstr object that is the transpose of the original MQConstr object.
\end{quote}

\sphinxAtStartPar
\sphinxstylestrong{Return value}
\begin{quote}

\sphinxAtStartPar
Returns the transposed MQConstr object.
\end{quote}

\sphinxAtStartPar
\sphinxstylestrong{Example}

\begin{sphinxVerbatim}[commandchars=\\\{\}]
\PYG{n}{mc} \PYG{o}{=} \PYG{n}{m}\PYG{o}{.}\PYG{n}{addQConstr}\PYG{p}{(}\PYG{n}{mx}\PYG{o}{.}\PYG{n}{T} \PYG{o}{@} \PYG{n}{mx} \PYG{o}{\PYGZlt{}}\PYG{o}{=} \PYG{l+m+mf}{1.0}\PYG{p}{)}
\PYG{n+nb}{print}\PYG{p}{(}\PYG{n}{mc}\PYG{o}{.}\PYG{n}{transpose}\PYG{p}{(}\PYG{p}{)}\PYG{p}{)}
\end{sphinxVerbatim}
\end{quote}

\subsubsection{MQConstr.ndim}
\label{\detokenize{pyapiref:mqconstr-ndim}}\begin{quote}

\sphinxAtStartPar
\sphinxstylestrong{Synopsis}
\begin{quote}

\sphinxAtStartPar
\sphinxcode{\sphinxupquote{ndim}}
\end{quote}

\sphinxAtStartPar
\sphinxstylestrong{Description}
\begin{quote}

\sphinxAtStartPar
Dimensions of the MQConstr object.
\end{quote}

\sphinxAtStartPar
\sphinxstylestrong{Return value}
\begin{quote}

\sphinxAtStartPar
Integer value.
\end{quote}

\sphinxAtStartPar
\sphinxstylestrong{Example}

\begin{sphinxVerbatim}[commandchars=\\\{\}]
\PYG{n}{mc} \PYG{o}{=} \PYG{n}{m}\PYG{o}{.}\PYG{n}{addQConstr}\PYG{p}{(}\PYG{n}{mx}\PYG{o}{.}\PYG{n}{T} \PYG{o}{@} \PYG{n}{mx} \PYG{o}{\PYGZlt{}}\PYG{o}{=} \PYG{l+m+mf}{1.0}\PYG{p}{)}
\PYG{n+nb}{print}\PYG{p}{(}\PYG{n}{mc}\PYG{o}{.}\PYG{n}{ndim}\PYG{p}{)}
\end{sphinxVerbatim}
\end{quote}

\subsubsection{MQConstr.shape}
\label{\detokenize{pyapiref:mqconstr-shape}}\begin{quote}

\sphinxAtStartPar
\sphinxstylestrong{Synopsis}
\begin{quote}

\sphinxAtStartPar
\sphinxcode{\sphinxupquote{shape}}
\end{quote}

\sphinxAtStartPar
\sphinxstylestrong{Description}
\begin{quote}

\sphinxAtStartPar
The shape of the MQConstr object.
\end{quote}

\sphinxAtStartPar
\sphinxstylestrong{Return value}
\begin{quote}

\sphinxAtStartPar
Integer tuple.
\end{quote}

\sphinxAtStartPar
\sphinxstylestrong{Example}

\begin{sphinxVerbatim}[commandchars=\\\{\}]
\PYG{n}{mc} \PYG{o}{=} \PYG{n}{m}\PYG{o}{.}\PYG{n}{addQConstr}\PYG{p}{(}\PYG{n}{mx}\PYG{o}{.}\PYG{n}{T} \PYG{o}{@} \PYG{n}{mx} \PYG{o}{\PYGZlt{}}\PYG{o}{=} \PYG{l+m+mf}{1.0}\PYG{p}{)}
\PYG{n+nb}{print}\PYG{p}{(}\PYG{n}{mc}\PYG{o}{.}\PYG{n}{shape}\PYG{p}{)}
\end{sphinxVerbatim}
\end{quote}

\subsubsection{MQConstr.size}
\label{\detokenize{pyapiref:mqconstr-size}}\begin{quote}

\sphinxAtStartPar
\sphinxstylestrong{Synopsis}
\begin{quote}

\sphinxAtStartPar
\sphinxcode{\sphinxupquote{size}}
\end{quote}

\sphinxAtStartPar
\sphinxstylestrong{Description}
\begin{quote}

\sphinxAtStartPar
The number of constraints of the MQConstr object.
\end{quote}

\sphinxAtStartPar
\sphinxstylestrong{Return value}
\begin{quote}

\sphinxAtStartPar
Integer value.
\end{quote}

\sphinxAtStartPar
\sphinxstylestrong{Example}

\begin{sphinxVerbatim}[commandchars=\\\{\}]
\PYG{n}{mc} \PYG{o}{=} \PYG{n}{m}\PYG{o}{.}\PYG{n}{addQConstr}\PYG{p}{(}\PYG{n}{mx}\PYG{o}{.}\PYG{n}{T} \PYG{o}{@} \PYG{n}{mx} \PYG{o}{\PYGZlt{}}\PYG{o}{=} \PYG{l+m+mf}{1.0}\PYG{p}{)}
\PYG{n+nb}{print}\PYG{p}{(}\PYG{n}{mc}\PYG{o}{.}\PYG{n}{size}\PYG{p}{)}
\end{sphinxVerbatim}
\end{quote}

\subsubsection{MQConstr.T}
\label{\detokenize{pyapiref:mqconstr-t}}\begin{quote}

\sphinxAtStartPar
\sphinxstylestrong{Synopsis}
\begin{quote}

\sphinxAtStartPar
\sphinxcode{\sphinxupquote{T}}
\end{quote}

\sphinxAtStartPar
\sphinxstylestrong{Description}
\begin{quote}

\sphinxAtStartPar
Transpose of the MQConstr object. Similar to the class method \sphinxcode{\sphinxupquote{transpose()}} .
\end{quote}

\sphinxAtStartPar
\sphinxstylestrong{Return value}
\begin{quote}

\sphinxAtStartPar
Returns the transposed MQConstr object.
\end{quote}

\sphinxAtStartPar
\sphinxstylestrong{Example}

\begin{sphinxVerbatim}[commandchars=\\\{\}]
\PYG{n}{A} \PYG{o}{=} \PYG{n}{np}\PYG{o}{.}\PYG{n}{ones}\PYG{p}{(}\PYG{p}{[}\PYG{l+m+mi}{4}\PYG{p}{,} \PYG{l+m+mi}{3}\PYG{p}{]}\PYG{p}{)}
\PYG{n}{mx} \PYG{o}{=} \PYG{n}{model}\PYG{o}{.}\PYG{n}{addMVar}\PYG{p}{(}\PYG{p}{(}\PYG{l+m+mi}{3}\PYG{p}{,} \PYG{l+m+mi}{4}\PYG{p}{)}\PYG{p}{,} \PYG{n}{nameprefix}\PYG{o}{=}\PYG{l+s+s2}{\PYGZdq{}}\PYG{l+s+s2}{mx}\PYG{l+s+s2}{\PYGZdq{}}\PYG{p}{)}
\PYG{n}{mc} \PYG{o}{=} \PYG{n}{model}\PYG{o}{.}\PYG{n}{addQConstr}\PYG{p}{(}\PYG{n}{mx} \PYG{o}{@} \PYG{n}{A} \PYG{o}{@} \PYG{n}{mx} \PYG{o}{==} \PYG{l+m+mf}{0.0}\PYG{p}{)}
\PYG{n+nb}{print}\PYG{p}{(}\PYG{n}{mc}\PYG{o}{.}\PYG{n}{shape}\PYG{p}{)} \PYG{c+c1}{\PYGZsh{} shape = (3, 4)}
\PYG{n+nb}{print}\PYG{p}{(}\PYG{n}{mc}\PYG{o}{.}\PYG{n}{T}\PYG{o}{.}\PYG{n}{shape}\PYG{p}{)} \PYG{c+c1}{\PYGZsh{} shape = (4, 3)}
\end{sphinxVerbatim}
\end{quote}

\subsection{MQConstrBuilder Class}
\label{\detokenize{pyapiref:mqconstrbuilder-class}}\label{\detokenize{pyapiref:chappyapi-mqconstrbuilder}}
\sphinxAtStartPar
The MQConstrBuilder class is used to build multi\sphinxhyphen{}dimensional quadratic
constraints in COPT and supports NumPy’s multi\sphinxhyphen{}dimensional array operations.
Users might create a MQConstrBuilder object by its constructor with a list of
{\hyperref[\detokenize{pyapiref:chappyapi-qconstrbuilder}]{\sphinxcrossref{\DUrole{std,std-ref}{QConstrBuilder Class}}}} objects, or simply by overloaded comparison operator
of {\hyperref[\detokenize{pyapiref:chappyapi-mquadexpr}]{\sphinxcrossref{\DUrole{std,std-ref}{MQuadExpr Class}}}}. The following member methods are provided:

\subsubsection{MQConstrBuilder()}
\label{\detokenize{pyapiref:mqconstrbuilder}}\begin{quote}

\sphinxAtStartPar
\sphinxstylestrong{Synopsis}
\begin{quote}

\sphinxAtStartPar
\sphinxcode{\sphinxupquote{MQConstrBuilder(args, shape=None)}}
\end{quote}

\sphinxAtStartPar
\sphinxstylestrong{Description}
\begin{quote}

\sphinxAtStartPar
construtor of MQConstrBuilder.
\end{quote}

\sphinxAtStartPar
\sphinxstylestrong{Arguments}
\begin{quote}

\sphinxAtStartPar
\sphinxcode{\sphinxupquote{args}}
\begin{quote}

\sphinxAtStartPar
one or a set of {\hyperref[\detokenize{pyapiref:chappyapi-qconstrbuilder}]{\sphinxcrossref{\DUrole{std,std-ref}{QConstrBuilder Class}}}} objects,
in form of Python list or NumPy ndarray.
\end{quote}

\sphinxAtStartPar
\sphinxcode{\sphinxupquote{shape}}
\begin{quote}

\sphinxAtStartPar
an integer, or tuple of integers, which is the shape of new MQConstrBuilder object.
\end{quote}
\end{quote}

\sphinxAtStartPar
\sphinxstylestrong{Example}

\begin{sphinxVerbatim}[commandchars=\\\{\}]
\PYG{n}{x} \PYG{o}{=} \PYG{n}{model}\PYG{o}{.}\PYG{n}{addVar}\PYG{p}{(}\PYG{p}{)}
\PYG{n}{mqcb} \PYG{o}{=} \PYG{n}{MQConstrBuilder}\PYG{p}{(}\PYG{n}{x} \PYG{o}{*} \PYG{n}{x} \PYG{o}{\PYGZlt{}}\PYG{o}{=} \PYG{l+m+mf}{9.0}\PYG{p}{)}

\PYG{c+c1}{\PYGZsh{} or by overloaded comparison operator of MQuadExpr}
\PYG{n}{mx} \PYG{o}{=} \PYG{n}{model}\PYG{o}{.}\PYG{n}{addMVar}\PYG{p}{(}\PYG{l+m+mi}{3}\PYG{p}{,} \PYG{l+m+mi}{3}\PYG{p}{)}
\PYG{n}{mqcb} \PYG{o}{=} \PYG{n}{mx} \PYG{o}{@} \PYG{n}{mx} \PYG{o}{\PYGZgt{}}\PYG{o}{=} \PYG{l+m+mf}{1.0}

\PYG{c+c1}{\PYGZsh{} or immediately passed to addConstrs()}
\PYG{n}{ma} \PYG{o}{=} \PYG{n}{model}\PYG{o}{.}\PYG{n}{addMVar}\PYG{p}{(}\PYG{l+m+mi}{2}\PYG{p}{)}
\PYG{n}{A} \PYG{o}{=} \PYG{n}{np}\PYG{o}{.}\PYG{n}{full}\PYG{p}{(}\PYG{p}{(}\PYG{l+m+mi}{2}\PYG{p}{,}\PYG{l+m+mi}{3}\PYG{p}{)}\PYG{p}{,} \PYG{l+m+mi}{1}\PYG{p}{)}
\PYG{n}{mb} \PYG{o}{=} \PYG{n}{model}\PYG{o}{.}\PYG{n}{addMVar}\PYG{p}{(}\PYG{l+m+mi}{3}\PYG{p}{)}
\PYG{n}{model}\PYG{o}{.}\PYG{n}{addQConstr}\PYG{p}{(}\PYG{n}{ma} \PYG{o}{@} \PYG{n}{A} \PYG{o}{@} \PYG{n}{mb} \PYG{o}{\PYGZlt{}}\PYG{o}{=} \PYG{l+m+mf}{1.0}\PYG{p}{)}
\end{sphinxVerbatim}
\end{quote}

\subsection{MPsdConstr Class}
\label{\detokenize{pyapiref:mpsdconstr-class}}\label{\detokenize{pyapiref:chappyapi-mpsdconstr}}
\sphinxAtStartPar
The \sphinxtitleref{MPsdConstr} class in the COPT is used to construct multi\sphinxhyphen{}dimensional semidefinite constraints. It is created via the \sphinxtitleref{addConstr} or \sphinxtitleref{addConstrs} methods of the model class.
It supports a unified set of indexing and slicing rules.
For detailed usage, see {\hyperref[\detokenize{matrix:chapmatrixindex}]{\sphinxcrossref{\DUrole{std,std-ref}{Matrix Modeling: Indexing and Slicing}}}}.

\sphinxAtStartPar
The following member methods are provided:

\subsubsection{MPsdConstr.getInfo()}
\label{\detokenize{pyapiref:mpsdconstr-getinfo}}\begin{quote}

\sphinxAtStartPar
\sphinxstylestrong{Synopsis}
\begin{quote}

\sphinxAtStartPar
\sphinxcode{\sphinxupquote{getInfo(infoname)}}
\end{quote}

\sphinxAtStartPar
\sphinxstylestrong{Description}
\begin{quote}

\sphinxAtStartPar
Retrieves the information value for each semidefinite constraint in the \sphinxtitleref{MPsdConstr}.
\end{quote}

\sphinxAtStartPar
\sphinxstylestrong{Arguments}
\begin{quote}

\sphinxAtStartPar
\sphinxcode{\sphinxupquote{infoname}}
\begin{quote}

\sphinxAtStartPar
The name of the information to retrieve. Possible values are detailed in {\hyperref[\detokenize{information:chapinfo}]{\sphinxcrossref{\DUrole{std,std-ref}{Information Section}}}}.
\end{quote}
\end{quote}

\sphinxAtStartPar
\sphinxstylestrong{Return Value}
\begin{quote}

\sphinxAtStartPar
Returns the information value as a multi\sphinxhyphen{}dimensional array.
\end{quote}

\sphinxAtStartPar
\sphinxstylestrong{Example}

\begin{sphinxVerbatim}[commandchars=\\\{\}]
\PYG{n}{mpsdcon} \PYG{o}{=} \PYG{n}{model}\PYG{o}{.}\PYG{n}{addConstr}\PYG{p}{(}\PYG{n}{barX}\PYG{p}{[}\PYG{p}{:}\PYG{o}{\PYGZhy{}}\PYG{l+m+mi}{1}\PYG{p}{,} \PYG{p}{:}\PYG{o}{\PYGZhy{}}\PYG{l+m+mi}{1}\PYG{p}{]}\PYG{o}{.}\PYG{n}{sum}\PYG{p}{(}\PYG{p}{)} \PYG{o}{==} \PYG{l+m+mi}{1}\PYG{p}{)}
\PYG{n+nb}{print}\PYG{p}{(}\PYG{n}{mpsdcon}\PYG{o}{.}\PYG{n}{UB}\PYG{p}{)}
\PYG{n+nb}{print}\PYG{p}{(}\PYG{n}{mpsdcon}\PYG{o}{.}\PYG{n}{getInfo}\PYG{p}{(}\PYG{l+s+s2}{\PYGZdq{}}\PYG{l+s+s2}{UB}\PYG{l+s+s2}{\PYGZdq{}}\PYG{p}{)}\PYG{p}{)}
\end{sphinxVerbatim}
\end{quote}

\subsubsection{MPsdConstr.setInfo()}
\label{\detokenize{pyapiref:mpsdconstr-setinfo}}\begin{quote}

\sphinxAtStartPar
\sphinxstylestrong{Synopsis}
\begin{quote}

\sphinxAtStartPar
\sphinxcode{\sphinxupquote{setInfo(infoname, newval)}}
\end{quote}

\sphinxAtStartPar
\sphinxstylestrong{Description}
\begin{quote}

\sphinxAtStartPar
Sets the information value for each semidefinite constraint in the \sphinxtitleref{MPsdConstr}.
\end{quote}

\sphinxAtStartPar
\sphinxstylestrong{Arguments}
\begin{quote}

\sphinxAtStartPar
\sphinxcode{\sphinxupquote{infoname}}
\begin{quote}

\sphinxAtStartPar
The name of the information to set. Possible values are detailed in {\hyperref[\detokenize{information:chapinfo}]{\sphinxcrossref{\DUrole{std,std-ref}{Information Section}}}}.
\end{quote}

\sphinxAtStartPar
\sphinxcode{\sphinxupquote{newval}}
\begin{quote}

\sphinxAtStartPar
The new value to set for the information.
\end{quote}
\end{quote}

\sphinxAtStartPar
\sphinxstylestrong{Example}

\begin{sphinxVerbatim}[commandchars=\\\{\}]
\PYG{n}{barX} \PYG{o}{=} \PYG{n}{model}\PYG{o}{.}\PYG{n}{addPsdVars}\PYG{p}{(}\PYG{l+m+mi}{3}\PYG{p}{,} \PYG{l+s+s2}{\PYGZdq{}}\PYG{l+s+s2}{BAR\PYGZus{}X}\PYG{l+s+s2}{\PYGZdq{}}\PYG{p}{)}
\PYG{n}{mY} \PYG{o}{=} \PYG{n}{model}\PYG{o}{.}\PYG{n}{addMVar}\PYG{p}{(}\PYG{l+m+mi}{2}\PYG{p}{,} \PYG{n}{nameprefix}\PYG{o}{=}\PYG{l+s+s2}{\PYGZdq{}}\PYG{l+s+s2}{M\PYGZus{}X}\PYG{l+s+s2}{\PYGZdq{}}\PYG{p}{)}
\PYG{n}{mpsdCon} \PYG{o}{=} \PYG{n}{model}\PYG{o}{.}\PYG{n}{addConstrs}\PYG{p}{(}\PYG{n}{barX}\PYG{p}{[}\PYG{p}{:}\PYG{o}{\PYGZhy{}}\PYG{l+m+mi}{1}\PYG{p}{,} \PYG{p}{:}\PYG{o}{\PYGZhy{}}\PYG{l+m+mi}{1}\PYG{p}{]}\PYG{o}{.}\PYG{n}{diagonal}\PYG{p}{(}\PYG{p}{)} \PYG{o}{==} \PYG{n}{mY}\PYG{p}{)}
\PYG{n}{mpsdCon}\PYG{o}{.}\PYG{n}{UB} \PYG{o}{=} \PYG{l+m+mi}{10}
\PYG{n}{mpsdCon}\PYG{o}{.}\PYG{n}{setInfo}\PYG{p}{(}\PYG{l+s+s2}{\PYGZdq{}}\PYG{l+s+s2}{UB}\PYG{l+s+s2}{\PYGZdq{}}\PYG{p}{,} \PYG{l+m+mi}{10}\PYG{p}{)}
\end{sphinxVerbatim}
\end{quote}

\subsubsection{MPsdConstr.item()}
\label{\detokenize{pyapiref:mpsdconstr-item}}\begin{quote}

\sphinxAtStartPar
\sphinxstylestrong{Synopsis}
\begin{quote}

\sphinxAtStartPar
\sphinxcode{\sphinxupquote{item()}}
\end{quote}

\sphinxAtStartPar
\sphinxstylestrong{Description}
\begin{quote}

\sphinxAtStartPar
Retrieves the semidefinite constraint object in a 0\sphinxhyphen{}dimensional \sphinxtitleref{MPsdConstr}. If the \sphinxtitleref{MPsdConstr} object is not 0\sphinxhyphen{}dimensional, raises a \sphinxtitleref{ValueError}.
\end{quote}

\sphinxAtStartPar
\sphinxstylestrong{Return Value}
\begin{quote}

\sphinxAtStartPar
Returns a \sphinxtitleref{PsdConstraint} object.
\end{quote}

\sphinxAtStartPar
\sphinxstylestrong{Example}

\begin{sphinxVerbatim}[commandchars=\\\{\}]
\PYG{n}{barX} \PYG{o}{=} \PYG{n}{model}\PYG{o}{.}\PYG{n}{addPsdVars}\PYG{p}{(}\PYG{l+m+mi}{3}\PYG{p}{,} \PYG{l+s+s2}{\PYGZdq{}}\PYG{l+s+s2}{BAR\PYGZus{}X}\PYG{l+s+s2}{\PYGZdq{}}\PYG{p}{)}
\PYG{n}{mY} \PYG{o}{=} \PYG{n}{model}\PYG{o}{.}\PYG{n}{addMVar}\PYG{p}{(}\PYG{l+m+mi}{2}\PYG{p}{,} \PYG{n}{nameprefix}\PYG{o}{=}\PYG{l+s+s2}{\PYGZdq{}}\PYG{l+s+s2}{M\PYGZus{}X}\PYG{l+s+s2}{\PYGZdq{}}\PYG{p}{)}
\PYG{n}{mpsdCon} \PYG{o}{=} \PYG{n}{model}\PYG{o}{.}\PYG{n}{addConstrs}\PYG{p}{(}\PYG{n}{barX}\PYG{p}{[}\PYG{p}{:}\PYG{o}{\PYGZhy{}}\PYG{l+m+mi}{1}\PYG{p}{,} \PYG{p}{:}\PYG{o}{\PYGZhy{}}\PYG{l+m+mi}{1}\PYG{p}{]}\PYG{o}{.}\PYG{n}{diagonal}\PYG{p}{(}\PYG{p}{)} \PYG{o}{==} \PYG{n}{mY}\PYG{p}{)}
\PYG{n}{psdCon} \PYG{o}{=} \PYG{n}{mpsdCon}\PYG{p}{[}\PYG{l+m+mi}{1}\PYG{p}{]}\PYG{o}{.}\PYG{n}{item}\PYG{p}{(}\PYG{p}{)}
\end{sphinxVerbatim}
\end{quote}

\subsubsection{MPsdConstr.clone()}
\label{\detokenize{pyapiref:mpsdconstr-clone}}\begin{quote}

\sphinxAtStartPar
\sphinxstylestrong{Synopsis}
\begin{quote}

\sphinxAtStartPar
\sphinxcode{\sphinxupquote{clone()}}
\end{quote}

\sphinxAtStartPar
\sphinxstylestrong{Description}
\begin{quote}

\sphinxAtStartPar
Creates a deep copy of an {\hyperref[\detokenize{pyapiref:chappyapi-mpsdconstr}]{\sphinxcrossref{\DUrole{std,std-ref}{MPsdConstr Class}}}} object.
\end{quote}

\sphinxAtStartPar
\sphinxstylestrong{Return Value}
\begin{quote}

\sphinxAtStartPar
Returns a new \sphinxtitleref{MPsdConstr} object.
\end{quote}

\sphinxAtStartPar
\sphinxstylestrong{Example}

\begin{sphinxVerbatim}[commandchars=\\\{\}]
\PYG{n}{barX} \PYG{o}{=} \PYG{n}{model}\PYG{o}{.}\PYG{n}{addPsdVars}\PYG{p}{(}\PYG{l+m+mi}{3}\PYG{p}{,} \PYG{l+s+s2}{\PYGZdq{}}\PYG{l+s+s2}{BAR\PYGZus{}X}\PYG{l+s+s2}{\PYGZdq{}}\PYG{p}{)}
\PYG{n}{mY} \PYG{o}{=} \PYG{n}{model}\PYG{o}{.}\PYG{n}{addMVar}\PYG{p}{(}\PYG{l+m+mi}{2}\PYG{p}{,} \PYG{n}{nameprefix}\PYG{o}{=}\PYG{l+s+s2}{\PYGZdq{}}\PYG{l+s+s2}{M\PYGZus{}X}\PYG{l+s+s2}{\PYGZdq{}}\PYG{p}{)}
\PYG{n}{mpsdCon} \PYG{o}{=} \PYG{n}{model}\PYG{o}{.}\PYG{n}{addConstrs}\PYG{p}{(}\PYG{n}{barX}\PYG{p}{[}\PYG{p}{:}\PYG{o}{\PYGZhy{}}\PYG{l+m+mi}{1}\PYG{p}{,} \PYG{p}{:}\PYG{o}{\PYGZhy{}}\PYG{l+m+mi}{1}\PYG{p}{]}\PYG{o}{.}\PYG{n}{diagonal}\PYG{p}{(}\PYG{p}{)} \PYG{o}{==} \PYG{n}{mY}\PYG{p}{)}
\PYG{n}{mpsdCon\PYGZus{}copy} \PYG{o}{=} \PYG{n}{mpsdCon}\PYG{o}{.}\PYG{n}{clone}\PYG{p}{(}\PYG{p}{)}
\end{sphinxVerbatim}
\end{quote}

\subsection{MPsdConstrBuilder Class}
\label{\detokenize{pyapiref:mpsdconstrbuilder-class}}\label{\detokenize{pyapiref:chappyapi-mpsdconstrbuilder}}
\sphinxAtStartPar
The \sphinxtitleref{MPsdConstrBuilder} class in COPT is used as a generator for multi\sphinxhyphen{}dimensional semidefinite constraints.
It can be generated from a group of {\hyperref[\detokenize{pyapiref:chappyapi-psdconstrbuilder}]{\sphinxcrossref{\DUrole{std,std-ref}{PsdConstrBuilder Class}}}} objects
or by overloading the comparison operators of {\hyperref[\detokenize{pyapiref:chappyapi-mvar}]{\sphinxcrossref{\DUrole{std,std-ref}{MVar Class}}}} and {\hyperref[\detokenize{pyapiref:chappyapi-mpsdexpr}]{\sphinxcrossref{\DUrole{std,std-ref}{MPsdExpr Class}}}} objects.
The following methods are provided:

\subsubsection{MPsdConstrBuilder()}
\label{\detokenize{pyapiref:mpsdconstrbuilder}}\begin{quote}

\sphinxAtStartPar
\sphinxstylestrong{Synopsis}
\begin{quote}

\sphinxAtStartPar
\sphinxcode{\sphinxupquote{MPsdConstrBuilder()}}
\end{quote}

\sphinxAtStartPar
\sphinxstylestrong{Description}
\begin{quote}

\sphinxAtStartPar
Creates an empty {\hyperref[\detokenize{pyapiref:chappyapi-mpsdconstrbuilder}]{\sphinxcrossref{\DUrole{std,std-ref}{MPsdConstrBuilder Class}}}} object.
\end{quote}

\sphinxAtStartPar
\sphinxstylestrong{Example}
\end{quote}

\begin{sphinxVerbatim}[commandchars=\\\{\}]
\PYG{c+c1}{\PYGZsh{} Create an empty semidefinite constraint builder}
\PYG{n}{constrbuilder} \PYG{o}{=} \PYG{n}{MPsdConstrBuilder}\PYG{p}{(}\PYG{p}{)}
\end{sphinxVerbatim}

\subsubsection{MPsdConstrBuilder.setBuilder()}
\label{\detokenize{pyapiref:mpsdconstrbuilder-setbuilder}}\begin{quote}

\sphinxAtStartPar
\sphinxstylestrong{Synopsis}
\begin{quote}

\sphinxAtStartPar
\sphinxcode{\sphinxupquote{setBuilder(expr, sense, rhs)}}
\end{quote}

\sphinxAtStartPar
\sphinxstylestrong{Description}
\begin{quote}

\sphinxAtStartPar
Sets the expression and constraint type for the multi\sphinxhyphen{}dimensional semidefinite constraint builder.
\end{quote}

\sphinxAtStartPar
\sphinxstylestrong{Arguments}
\begin{quote}

\sphinxAtStartPar
\sphinxcode{\sphinxupquote{expr}}
\begin{quote}

\sphinxAtStartPar
The expression to set. Possible values include constants, {\hyperref[\detokenize{pyapiref:chappyapi-mvar}]{\sphinxcrossref{\DUrole{std,std-ref}{MVar Class}}}} objects,
{\hyperref[\detokenize{pyapiref:chappyapi-mlinexpr}]{\sphinxcrossref{\DUrole{std,std-ref}{MLinExpr Class}}}} objects, or {\hyperref[\detokenize{pyapiref:chappyapi-mpsdexpr}]{\sphinxcrossref{\DUrole{std,std-ref}{MPsdExpr Class}}}} objects.
\end{quote}

\sphinxAtStartPar
\sphinxcode{\sphinxupquote{sense}}
\begin{quote}

\sphinxAtStartPar
The constraint type. Possible values are detailed in {\hyperref[\detokenize{constant:chapconst-constrtype}]{\sphinxcrossref{\DUrole{std,std-ref}{Constraint Types}}}}.
\end{quote}

\sphinxAtStartPar
\sphinxcode{\sphinxupquote{rhs}}
\begin{quote}

\sphinxAtStartPar
The right\sphinxhyphen{}hand side of the constraint.
Possible values include constants, {\hyperref[\detokenize{pyapiref:chappyapi-mvar}]{\sphinxcrossref{\DUrole{std,std-ref}{MVar Class}}}} objects,
{\hyperref[\detokenize{pyapiref:chappyapi-mlinexpr}]{\sphinxcrossref{\DUrole{std,std-ref}{MLinExpr Class}}}} objects, or {\hyperref[\detokenize{pyapiref:chappyapi-mpsdexpr}]{\sphinxcrossref{\DUrole{std,std-ref}{MPsdExpr Class}}}} objects.
\end{quote}
\end{quote}

\sphinxAtStartPar
\sphinxstylestrong{Example}
\end{quote}

\begin{sphinxVerbatim}[commandchars=\\\{\}]
\PYG{c+c1}{\PYGZsh{} Set the expression for the multi\PYGZhy{}dimensional PSD constraint builder: x + y == 1}
\PYG{n}{constrbuilder}\PYG{o}{.}\PYG{n}{setBuilder}\PYG{p}{(}\PYG{n}{x} \PYG{o}{+} \PYG{n}{y}\PYG{p}{,} \PYG{n}{COPT}\PYG{o}{.}\PYG{n}{EQUAL}\PYG{p}{,} \PYG{l+m+mi}{1}\PYG{p}{)}
\end{sphinxVerbatim}

\subsubsection{MPsdConstrBuilder.setRange()}
\label{\detokenize{pyapiref:mpsdconstrbuilder-setrange}}\begin{quote}

\sphinxAtStartPar
\sphinxstylestrong{Synopsis}
\begin{quote}

\sphinxAtStartPar
\sphinxcode{\sphinxupquote{setRange(expr, range)}}
\end{quote}

\sphinxAtStartPar
\sphinxstylestrong{Description}
\begin{quote}

\sphinxAtStartPar
Sets the expression and the range for the multi\sphinxhyphen{}dimensional PSD constraint builder.
The format is \sphinxtitleref{expr} less than or equal to 0, and greater than or equal to \sphinxtitleref{\sphinxhyphen{}range}.
\end{quote}

\sphinxAtStartPar
\sphinxstylestrong{Arguments}
\begin{quote}

\sphinxAtStartPar
\sphinxcode{\sphinxupquote{expr}}
\begin{quote}

\sphinxAtStartPar
The expression to set. Possible values include {\hyperref[\detokenize{pyapiref:chappyapi-mpsdexpr}]{\sphinxcrossref{\DUrole{std,std-ref}{MPsdExpr Class}}}} objects.
\end{quote}

\sphinxAtStartPar
\sphinxcode{\sphinxupquote{range}}
\begin{quote}

\sphinxAtStartPar
The right\sphinxhyphen{}hand side of the range constraint. Must be non\sphinxhyphen{}negative.
\end{quote}
\end{quote}

\sphinxAtStartPar
\sphinxstylestrong{Example}
\end{quote}

\begin{sphinxVerbatim}[commandchars=\\\{\}]
\PYG{c+c1}{\PYGZsh{} Set the expression for the multi\PYGZhy{}dimensional PSD constraint builder: x + y \PYGZhy{} 1, with range 1}
\PYG{n}{constrbuilder}\PYG{o}{.}\PYG{n}{setRange}\PYG{p}{(}\PYG{n}{x} \PYG{o}{+} \PYG{n}{y} \PYG{o}{\PYGZhy{}} \PYG{l+m+mi}{1}\PYG{p}{,} \PYG{l+m+mi}{1}\PYG{p}{)}
\end{sphinxVerbatim}

\subsubsection{MPsdConstrBuilder.getPsdExpr()}
\label{\detokenize{pyapiref:mpsdconstrbuilder-getpsdexpr}}\begin{quote}

\sphinxAtStartPar
\sphinxstylestrong{Synopsis}
\begin{quote}

\sphinxAtStartPar
\sphinxcode{\sphinxupquote{getPsdExpr()}}
\end{quote}

\sphinxAtStartPar
\sphinxstylestrong{Description}
\begin{quote}

\sphinxAtStartPar
Retrieves the expression of the multi\sphinxhyphen{}dimensional PSD constraint builder.
\end{quote}

\sphinxAtStartPar
\sphinxstylestrong{Example}
\end{quote}

\begin{sphinxVerbatim}[commandchars=\\\{\}]
\PYG{c+c1}{\PYGZsh{} Retrieve the expression of the multi\PYGZhy{}dimensional PSD constraint builder}
\PYG{n}{psdexpr} \PYG{o}{=} \PYG{n}{constrbuilder}\PYG{o}{.}\PYG{n}{getPsdExpr}\PYG{p}{(}\PYG{p}{)}
\end{sphinxVerbatim}

\subsubsection{MPsdConstrBuilder.getSense()}
\label{\detokenize{pyapiref:mpsdconstrbuilder-getsense}}\begin{quote}

\sphinxAtStartPar
\sphinxstylestrong{Synopsis}
\begin{quote}

\sphinxAtStartPar
\sphinxcode{\sphinxupquote{getSense()}}
\end{quote}

\sphinxAtStartPar
\sphinxstylestrong{Description}
\begin{quote}

\sphinxAtStartPar
Retrieves the constraint type of the multi\sphinxhyphen{}dimensional PSD constraint builder.
\end{quote}

\sphinxAtStartPar
\sphinxstylestrong{Example}
\end{quote}

\begin{sphinxVerbatim}[commandchars=\\\{\}]
\PYG{c+c1}{\PYGZsh{} Retrieve the constraint type of the multi\PYGZhy{}dimensional PSD constraint builder}
\PYG{n}{consense} \PYG{o}{=} \PYG{n}{constrbuilder}\PYG{o}{.}\PYG{n}{getSense}\PYG{p}{(}\PYG{p}{)}
\end{sphinxVerbatim}

\subsubsection{MPsdConstrBuilder.getRange()}
\label{\detokenize{pyapiref:mpsdconstrbuilder-getrange}}\begin{quote}

\sphinxAtStartPar
\sphinxstylestrong{Synopsis}
\begin{quote}

\sphinxAtStartPar
\sphinxcode{\sphinxupquote{getRange()}}
\end{quote}

\sphinxAtStartPar
\sphinxstylestrong{Description}
\begin{quote}

\sphinxAtStartPar
Retrieves the range right\sphinxhyphen{}hand side of the multi\sphinxhyphen{}dimensional PSD constraint builder,
representing the range between the upper and lower bounds.
\end{quote}

\sphinxAtStartPar
\sphinxstylestrong{Example}
\end{quote}

\begin{sphinxVerbatim}[commandchars=\\\{\}]
\PYG{c+c1}{\PYGZsh{} Retrieve the range right\PYGZhy{}hand side of the multi\PYGZhy{}dimensional PSD constraint builder}
\PYG{n}{rngval} \PYG{o}{=} \PYG{n}{constrbuilder}\PYG{o}{.}\PYG{n}{getRange}\PYG{p}{(}\PYG{p}{)}
\end{sphinxVerbatim}

\subsection{MLinExpr Class}
\label{\detokenize{pyapiref:mlinexpr-class}}\label{\detokenize{pyapiref:chappyapi-mlinexpr}}
\sphinxAtStartPar
The MLinExpr class is used in COPT to build multi\sphinxhyphen{}dimensional linear expressions
and supports NumPy’s multi\sphinxhyphen{}dimensional array operations.
The MLinExpr object with an initial value of 0.0 can be generated by the class
generation method \sphinxcode{\sphinxupquote{zeros()}} , or by {\hyperref[\detokenize{pyapiref:chappyapi-mvar}]{\sphinxcrossref{\DUrole{std,std-ref}{MVar Class}}}} object to generate
a linear combination.
It supports a unified set of indexing and slicing rules.
For detailed usage, see {\hyperref[\detokenize{matrix:chapmatrixindex}]{\sphinxcrossref{\DUrole{std,std-ref}{Matrix Modeling: Indexing and Slicing}}}}.

\sphinxAtStartPar
The following member methods are provided:

\subsubsection{MLinExpr.zeros()}
\label{\detokenize{pyapiref:mlinexpr-zeros}}\begin{quote}

\sphinxAtStartPar
\sphinxstylestrong{Synopsis}
\begin{quote}

\sphinxAtStartPar
\sphinxcode{\sphinxupquote{zeros(shape)}}
\end{quote}

\sphinxAtStartPar
\sphinxstylestrong{Description}
\begin{quote}

\sphinxAtStartPar
This is the class generation method and can be called directly without the MLinExpr object.
\end{quote}

\sphinxAtStartPar
\sphinxstylestrong{Arguments}
\begin{quote}

\sphinxAtStartPar
\sphinxcode{\sphinxupquote{shape}}
\begin{quote}

\sphinxAtStartPar
The value is an integer, or a tuple of integers. which represents the shape of the new MLinExpr object.
\end{quote}
\end{quote}

\sphinxAtStartPar
\sphinxstylestrong{Return value}
\begin{quote}

\sphinxAtStartPar
new MLinExpr object.
\end{quote}

\sphinxAtStartPar
\sphinxstylestrong{Example}

\begin{sphinxVerbatim}[commandchars=\\\{\}]
\PYG{n}{mexpr} \PYG{o}{=} \PYG{n}{MLinExpr}\PYG{o}{.}\PYG{n}{zeros}\PYG{p}{(}\PYG{p}{(}\PYG{l+m+mi}{2}\PYG{p}{,}\PYG{l+m+mi}{3}\PYG{p}{)}\PYG{p}{)}
\PYG{n}{x} \PYG{o}{=} \PYG{n}{model}\PYG{o}{.}\PYG{n}{addVar}\PYG{p}{(}\PYG{p}{)}
\PYG{n}{mexpr} \PYG{o}{+}\PYG{o}{=} \PYG{n}{x}
\end{sphinxVerbatim}
\end{quote}

\subsubsection{MLinExpr.clear()}
\label{\detokenize{pyapiref:mlinexpr-clear}}\begin{quote}

\sphinxAtStartPar
\sphinxstylestrong{Synopsis}
\begin{quote}

\sphinxAtStartPar
\sphinxcode{\sphinxupquote{clear()}}
\end{quote}

\sphinxAtStartPar
\sphinxstylestrong{Description}
\begin{quote}

\sphinxAtStartPar
Resets each element of the {\hyperref[\detokenize{pyapiref:chappyapi-mlinexpr}]{\sphinxcrossref{\DUrole{std,std-ref}{MLinExpr Class}}}} object to 0.0.
\end{quote}

\sphinxAtStartPar
\sphinxstylestrong{Example}

\begin{sphinxVerbatim}[commandchars=\\\{\}]
\PYG{n}{mexpr} \PYG{o}{=} \PYG{l+m+mf}{2.0} \PYG{o}{*} \PYG{n}{model}\PYG{o}{.}\PYG{n}{addMVar}\PYG{p}{(}\PYG{l+m+mi}{3}\PYG{p}{)}
\PYG{n}{mexpr}\PYG{o}{.}\PYG{n}{clear}\PYG{p}{(}\PYG{p}{)}
\end{sphinxVerbatim}
\end{quote}

\subsubsection{MLinExpr.clone()}
\label{\detokenize{pyapiref:mlinexpr-clone}}\begin{quote}

\sphinxAtStartPar
\sphinxstylestrong{Synopsis}
\begin{quote}

\sphinxAtStartPar
\sphinxcode{\sphinxupquote{clone()}}
\end{quote}

\sphinxAtStartPar
\sphinxstylestrong{Description}
\begin{quote}

\sphinxAtStartPar
Deep\sphinxhyphen{}copy a {\hyperref[\detokenize{pyapiref:chappyapi-mlinexpr}]{\sphinxcrossref{\DUrole{std,std-ref}{MLinExpr Class}}}} object.
\end{quote}

\sphinxAtStartPar
\sphinxstylestrong{Return value}
\begin{quote}

\sphinxAtStartPar
new MLinExpr object
\end{quote}

\sphinxAtStartPar
\sphinxstylestrong{Example}

\begin{sphinxVerbatim}[commandchars=\\\{\}]
\PYG{n}{mexpr} \PYG{o}{=} \PYG{l+m+mf}{2.0} \PYG{o}{*} \PYG{n}{model}\PYG{o}{.}\PYG{n}{addMVar}\PYG{p}{(}\PYG{l+m+mi}{3}\PYG{p}{)}
\PYG{n}{mexpr\PYGZus{}copy} \PYG{o}{=} \PYG{n}{mexpr}\PYG{o}{.}\PYG{n}{clone}\PYG{p}{(}\PYG{p}{)}
\end{sphinxVerbatim}
\end{quote}

\subsubsection{MLinExpr.getValue()}
\label{\detokenize{pyapiref:mlinexpr-getvalue}}\begin{quote}

\sphinxAtStartPar
\sphinxstylestrong{Synopsis}
\begin{quote}

\sphinxAtStartPar
\sphinxcode{\sphinxupquote{getValue()}}
\end{quote}

\sphinxAtStartPar
\sphinxstylestrong{Description}
\begin{quote}

\sphinxAtStartPar
Get the evaluation of each linear expression within the {\hyperref[\detokenize{pyapiref:chappyapi-mlinexpr}]{\sphinxcrossref{\DUrole{std,std-ref}{MLinExpr Class}}}} object.
\end{quote}

\sphinxAtStartPar
\sphinxstylestrong{Return value}
\begin{quote}

\sphinxAtStartPar
Returns a NumPy ndarray of the same dimensions as the MLinExpr object whose elements are the evaluations of the corresponding expression.
\end{quote}

\sphinxAtStartPar
\sphinxstylestrong{Example}

\begin{sphinxVerbatim}[commandchars=\\\{\}]
\PYG{n}{a} \PYG{o}{=} \PYG{n}{np}\PYG{o}{.}\PYG{n}{random}\PYG{o}{.}\PYG{n}{rand}\PYG{p}{(}\PYG{l+m+mi}{4}\PYG{p}{)}
\PYG{n}{mx} \PYG{o}{=} \PYG{n}{m}\PYG{o}{.}\PYG{n}{addMVar}\PYG{p}{(}\PYG{p}{(}\PYG{l+m+mi}{4}\PYG{p}{,} \PYG{l+m+mi}{3}\PYG{p}{)}\PYG{p}{,} \PYG{n}{nameprefix}\PYG{o}{=}\PYG{l+s+s2}{\PYGZdq{}}\PYG{l+s+s2}{mx}\PYG{l+s+s2}{\PYGZdq{}}\PYG{p}{)}
\PYG{n}{mexpr} \PYG{o}{=} \PYG{n}{a} \PYG{o}{@} \PYG{n}{mx}
\PYG{n}{mc} \PYG{o}{=} \PYG{n}{m}\PYG{o}{.}\PYG{n}{addConstrs}\PYG{p}{(}\PYG{n}{mexpr} \PYG{o}{\PYGZlt{}}\PYG{o}{=} \PYG{l+m+mf}{1.0}\PYG{p}{)}
\PYG{n}{model}\PYG{o}{.}\PYG{n}{solve}\PYG{p}{(}\PYG{p}{)}
\PYG{n+nb}{print}\PYG{p}{(}\PYG{n}{mc}\PYG{o}{.}\PYG{n}{getValue}\PYG{p}{(}\PYG{p}{)}\PYG{p}{)}
\end{sphinxVerbatim}
\end{quote}

\subsubsection{MLinExpr.item()}
\label{\detokenize{pyapiref:mlinexpr-item}}\begin{quote}

\sphinxAtStartPar
\sphinxstylestrong{Synopsis}
\begin{quote}

\sphinxAtStartPar
\sphinxcode{\sphinxupquote{item()}}
\end{quote}

\sphinxAtStartPar
\sphinxstylestrong{Description}
\begin{quote}

\sphinxAtStartPar
Get the constraint object in 0\sphinxhyphen{}dimensional MLinExpr. Raises a ValueError exception if the MLinExpr object is not 0\sphinxhyphen{}dimensional.
\end{quote}

\sphinxAtStartPar
\sphinxstylestrong{Return value}
\begin{quote}

\sphinxAtStartPar
Returns the linear constraint object.
\end{quote}

\sphinxAtStartPar
\sphinxstylestrong{Example}

\begin{sphinxVerbatim}[commandchars=\\\{\}]
\PYG{n}{mexpr} \PYG{o}{=} \PYG{l+m+mf}{2.0} \PYG{o}{*} \PYG{n}{model}\PYG{o}{.}\PYG{n}{addMVar}\PYG{p}{(}\PYG{l+m+mi}{3}\PYG{p}{)}
\PYG{n+nb}{print}\PYG{p}{(}\PYG{n}{mexpr}\PYG{p}{[}\PYG{l+m+mi}{0}\PYG{p}{]}\PYG{o}{.}\PYG{n}{item}\PYG{p}{(}\PYG{p}{)}\PYG{p}{)}
\end{sphinxVerbatim}
\end{quote}

\subsubsection{MLinExpr.reshape()}
\label{\detokenize{pyapiref:mlinexpr-reshape}}\begin{quote}

\sphinxAtStartPar
\sphinxstylestrong{Synopsis}
\begin{quote}

\sphinxAtStartPar
\sphinxcode{\sphinxupquote{reshape(shape, order=\textquotesingle{}C\textquotesingle{})}}
\end{quote}

\sphinxAtStartPar
\sphinxstylestrong{Description}
\begin{quote}

\sphinxAtStartPar
Returns a new MLinExpr object whose elements remain unchanged but whose shape is transformed by the argument shape.
\end{quote}

\sphinxAtStartPar
\sphinxstylestrong{Arguments}
\begin{quote}

\sphinxAtStartPar
\sphinxcode{\sphinxupquote{shape}}
\begin{quote}

\sphinxAtStartPar
The value is an integer, or a tuple of integers. which represents the shape of the new MLinExpr object.
\end{quote}

\sphinxAtStartPar
\sphinxcode{\sphinxupquote{order}}
\begin{quote}

\sphinxAtStartPar
Optional parameter, the default is the character ‘C’, which means it is compatible with the C language, that is, it is stored in rows; it can also be set to the character ‘F’, that is, it is stored in columns, and it is compatible with the Fortune language.
\end{quote}
\end{quote}

\sphinxAtStartPar
\sphinxstylestrong{Return value}
\begin{quote}

\sphinxAtStartPar
Returns a new MLinExpr object with the same elements as the original MLinExpr object but with a different shape.
\end{quote}

\sphinxAtStartPar
\sphinxstylestrong{Example}

\begin{sphinxVerbatim}[commandchars=\\\{\}]
\PYG{n}{mc} \PYG{o}{=} \PYG{n}{m}\PYG{o}{.}\PYG{n}{addConstrs}\PYG{p}{(}\PYG{n}{a} \PYG{o}{@} \PYG{n}{mx} \PYG{o}{\PYGZlt{}}\PYG{o}{=} \PYG{n}{b}\PYG{p}{)}
\PYG{n}{mc\PYGZus{}2x2} \PYG{o}{=} \PYG{n}{mc}\PYG{o}{.}\PYG{n}{reshape}\PYG{p}{(}\PYG{p}{(}\PYG{l+m+mi}{2}\PYG{p}{,} \PYG{l+m+mi}{2}\PYG{p}{)}\PYG{p}{)}
\end{sphinxVerbatim}
\end{quote}

\subsubsection{MLinExpr.sum()}
\label{\detokenize{pyapiref:mlinexpr-sum}}\begin{quote}

\sphinxAtStartPar
\sphinxstylestrong{Synopsis}
\begin{quote}

\sphinxAtStartPar
\sphinxcode{\sphinxupquote{sum(axis=None)}}
\end{quote}

\sphinxAtStartPar
\sphinxstylestrong{Description}
\begin{quote}

\sphinxAtStartPar
Sum the variables in the MLinExpr object, returning a new {\hyperref[\detokenize{pyapiref:chappyapi-mlinexpr}]{\sphinxcrossref{\DUrole{std,std-ref}{MLinExpr Class}}}} object.
\end{quote}

\sphinxAtStartPar
\sphinxstylestrong{Arguments}
\begin{quote}

\sphinxAtStartPar
\sphinxcode{\sphinxupquote{axis}}
\begin{quote}

\sphinxAtStartPar
Optional integer parameter, the default value is None, that is, to sum up variables one by one. Otherwise, sum over the given axis.
\end{quote}
\end{quote}

\sphinxAtStartPar
\sphinxstylestrong{Return value}
\begin{quote}

\sphinxAtStartPar
Returns an MLinExpr object representing the sum of the corresponding linear expressions.
\end{quote}

\sphinxAtStartPar
\sphinxstylestrong{Example}

\begin{sphinxVerbatim}[commandchars=\\\{\}]
\PYG{n}{mexpr} \PYG{o}{=} \PYG{l+m+mf}{2.0} \PYG{o}{*} \PYG{n}{model}\PYG{o}{.}\PYG{n}{addMVar}\PYG{p}{(}\PYG{p}{(}\PYG{l+m+mi}{3}\PYG{p}{,} \PYG{l+m+mi}{5}\PYG{p}{)}\PYG{p}{)}
\PYG{n}{sum\PYGZus{}all} \PYG{o}{=} \PYG{n}{mexpr}\PYG{o}{.}\PYG{n}{sum}\PYG{p}{(}\PYG{p}{)} \PYG{c+c1}{\PYGZsh{}return 0\PYGZhy{}dimensional MLinExpr object}
\PYG{n}{sum\PYGZus{}row} \PYG{o}{=} \PYG{n}{mexpr}\PYG{o}{.}\PYG{n}{sum}\PYG{p}{(}\PYG{n}{axis} \PYG{o}{=} \PYG{l+m+mi}{0}\PYG{p}{)} \PYG{c+c1}{\PYGZsh{}Return a 1\PYGZhy{}dimensional MLinExpr object with a shape of (5, )}
\end{sphinxVerbatim}
\end{quote}

\subsubsection{MLinExpr.tolist()}
\label{\detokenize{pyapiref:mlinexpr-tolist}}\begin{quote}

\sphinxAtStartPar
\sphinxstylestrong{Synopsis}
\begin{quote}

\sphinxAtStartPar
\sphinxcode{\sphinxupquote{tolist()}}
\end{quote}

\sphinxAtStartPar
\sphinxstylestrong{Description}
\begin{quote}

\sphinxAtStartPar
Converts an MLinExpr object to a one\sphinxhyphen{}dimensional list whose elements are linear expressions.
\end{quote}

\sphinxAtStartPar
\sphinxstylestrong{Return value}
\begin{quote}

\sphinxAtStartPar
Return a 1D list containing {\hyperref[\detokenize{pyapiref:chappyapi-linexpr}]{\sphinxcrossref{\DUrole{std,std-ref}{LinExpr Class}}}}.
\end{quote}

\sphinxAtStartPar
\sphinxstylestrong{Example}

\begin{sphinxVerbatim}[commandchars=\\\{\}]
\PYG{n}{mexpr} \PYG{o}{=} \PYG{l+m+mf}{2.0} \PYG{o}{*} \PYG{n}{model}\PYG{o}{.}\PYG{n}{addMVar}\PYG{p}{(}\PYG{p}{(}\PYG{l+m+mi}{3}\PYG{p}{,} \PYG{l+m+mi}{5}\PYG{p}{)}\PYG{p}{)}
\PYG{n+nb}{print}\PYG{p}{(}\PYG{n}{mexpr}\PYG{o}{.}\PYG{n}{tolist}\PYG{p}{(}\PYG{p}{)}\PYG{p}{)}
\end{sphinxVerbatim}
\end{quote}

\subsubsection{MLinExpr.transpose()}
\label{\detokenize{pyapiref:mlinexpr-transpose}}\begin{quote}

\sphinxAtStartPar
\sphinxstylestrong{Synopsis}
\begin{quote}

\sphinxAtStartPar
\sphinxcode{\sphinxupquote{transpose()}}
\end{quote}

\sphinxAtStartPar
\sphinxstylestrong{Description}
\begin{quote}

\sphinxAtStartPar
Generates a new MLinExpr object that is the transpose of the original MLinExpr object.
\end{quote}

\sphinxAtStartPar
\sphinxstylestrong{Return value}
\begin{quote}

\sphinxAtStartPar
Returns the transposed MLinExpr object.
\end{quote}

\sphinxAtStartPar
\sphinxstylestrong{Example}

\begin{sphinxVerbatim}[commandchars=\\\{\}]
\PYG{n}{mexpr} \PYG{o}{=} \PYG{l+m+mf}{2.0} \PYG{o}{*} \PYG{n}{model}\PYG{o}{.}\PYG{n}{addMVar}\PYG{p}{(}\PYG{p}{(}\PYG{l+m+mi}{3}\PYG{p}{,} \PYG{l+m+mi}{5}\PYG{p}{)}\PYG{p}{)}
\PYG{n+nb}{print}\PYG{p}{(}\PYG{n}{mexpr}\PYG{o}{.}\PYG{n}{transpose}\PYG{p}{(}\PYG{p}{)}\PYG{p}{)}
\end{sphinxVerbatim}
\end{quote}

\subsubsection{MLinExpr.ndim}
\label{\detokenize{pyapiref:mlinexpr-ndim}}\begin{quote}

\sphinxAtStartPar
\sphinxstylestrong{Synopsis}
\begin{quote}

\sphinxAtStartPar
\sphinxcode{\sphinxupquote{ndim}}
\end{quote}

\sphinxAtStartPar
\sphinxstylestrong{Description}
\begin{quote}

\sphinxAtStartPar
{\hyperref[\detokenize{pyapiref:chappyapi-mlinexpr}]{\sphinxcrossref{\DUrole{std,std-ref}{MLinExpr Class}}}} Dimensions of the object.
\end{quote}

\sphinxAtStartPar
\sphinxstylestrong{Return value}
\begin{quote}

\sphinxAtStartPar
Integer value.
\end{quote}

\sphinxAtStartPar
\sphinxstylestrong{Example}

\begin{sphinxVerbatim}[commandchars=\\\{\}]
\PYG{n}{mexpr} \PYG{o}{=} \PYG{l+m+mf}{2.0} \PYG{o}{*} \PYG{n}{model}\PYG{o}{.}\PYG{n}{addMVar}\PYG{p}{(}\PYG{p}{(}\PYG{l+m+mi}{3}\PYG{p}{,} \PYG{l+m+mi}{5}\PYG{p}{)}\PYG{p}{)}
\PYG{n+nb}{print}\PYG{p}{(}\PYG{n}{mexpr}\PYG{o}{.}\PYG{n}{ndim}\PYG{p}{)}
\end{sphinxVerbatim}
\end{quote}

\subsubsection{MLinExpr.shape}
\label{\detokenize{pyapiref:mlinexpr-shape}}\begin{quote}

\sphinxAtStartPar
\sphinxstylestrong{Synopsis}
\begin{quote}

\sphinxAtStartPar
\sphinxcode{\sphinxupquote{shape}}
\end{quote}

\sphinxAtStartPar
\sphinxstylestrong{Description}
\begin{quote}

\sphinxAtStartPar
{\hyperref[\detokenize{pyapiref:chappyapi-mlinexpr}]{\sphinxcrossref{\DUrole{std,std-ref}{MLinExpr Class}}}} The shape of the object.
\end{quote}

\sphinxAtStartPar
\sphinxstylestrong{Return value}
\begin{quote}

\sphinxAtStartPar
Integer tuple.
\end{quote}

\sphinxAtStartPar
\sphinxstylestrong{Example}

\begin{sphinxVerbatim}[commandchars=\\\{\}]
\PYG{n}{mexpr} \PYG{o}{=} \PYG{l+m+mf}{2.0} \PYG{o}{*} \PYG{n}{model}\PYG{o}{.}\PYG{n}{addMVar}\PYG{p}{(}\PYG{p}{(}\PYG{l+m+mi}{3}\PYG{p}{,} \PYG{l+m+mi}{5}\PYG{p}{)}\PYG{p}{)}
\PYG{n+nb}{print}\PYG{p}{(}\PYG{n}{mexpr}\PYG{o}{.}\PYG{n}{shape}\PYG{p}{)}
\end{sphinxVerbatim}
\end{quote}

\subsubsection{MLinExpr.size}
\label{\detokenize{pyapiref:mlinexpr-size}}\begin{quote}

\sphinxAtStartPar
\sphinxstylestrong{Synopsis}
\begin{quote}

\sphinxAtStartPar
\sphinxcode{\sphinxupquote{size}}
\end{quote}

\sphinxAtStartPar
\sphinxstylestrong{Description}
\begin{quote}

\sphinxAtStartPar
The number of elements of {\hyperref[\detokenize{pyapiref:chappyapi-mlinexpr}]{\sphinxcrossref{\DUrole{std,std-ref}{MLinExpr Class}}}} object.
\end{quote}

\sphinxAtStartPar
\sphinxstylestrong{Return value}
\begin{quote}

\sphinxAtStartPar
Integer value.
\end{quote}

\sphinxAtStartPar
\sphinxstylestrong{Example}

\begin{sphinxVerbatim}[commandchars=\\\{\}]
\PYG{n}{mexpr} \PYG{o}{=} \PYG{l+m+mf}{2.0} \PYG{o}{*} \PYG{n}{model}\PYG{o}{.}\PYG{n}{addMVar}\PYG{p}{(}\PYG{p}{(}\PYG{l+m+mi}{3}\PYG{p}{,} \PYG{l+m+mi}{5}\PYG{p}{)}\PYG{p}{)}
\PYG{n+nb}{print}\PYG{p}{(}\PYG{n}{mexpr}\PYG{o}{.}\PYG{n}{size}\PYG{p}{)}
\end{sphinxVerbatim}
\end{quote}

\subsubsection{MLinExpr.T}
\label{\detokenize{pyapiref:mlinexpr-t}}\begin{quote}

\sphinxAtStartPar
\sphinxstylestrong{Synopsis}
\begin{quote}

\sphinxAtStartPar
\sphinxcode{\sphinxupquote{T}}
\end{quote}

\sphinxAtStartPar
\sphinxstylestrong{Description}
\begin{quote}

\sphinxAtStartPar
Transpose of {\hyperref[\detokenize{pyapiref:chappyapi-mlinexpr}]{\sphinxcrossref{\DUrole{std,std-ref}{MLinExpr Class}}}} object. Similar to the class method transpose().
\end{quote}

\sphinxAtStartPar
\sphinxstylestrong{Return value}
\begin{quote}

\sphinxAtStartPar
Returns the transposed MLinExpr object.
\end{quote}

\sphinxAtStartPar
\sphinxstylestrong{Example}

\begin{sphinxVerbatim}[commandchars=\\\{\}]
\PYG{n}{mexpr} \PYG{o}{=} \PYG{l+m+mf}{2.0} \PYG{o}{*} \PYG{n}{model}\PYG{o}{.}\PYG{n}{addMVar}\PYG{p}{(}\PYG{p}{(}\PYG{l+m+mi}{3}\PYG{p}{,} \PYG{l+m+mi}{5}\PYG{p}{)}\PYG{p}{)}
\PYG{n+nb}{print}\PYG{p}{(}\PYG{n}{mexpr}\PYG{o}{.}\PYG{n}{T}\PYG{o}{.}\PYG{n}{shape}\PYG{p}{)} \PYG{c+c1}{\PYGZsh{} The transposed shape is (5, 3)}
\end{sphinxVerbatim}
\end{quote}

\subsubsection{MLinExpr.\_\_eq\_\_()}
\label{\detokenize{pyapiref:mlinexpr-eq}}\begin{quote}

\sphinxAtStartPar
\sphinxstylestrong{Synopsis}
\begin{quote}

\sphinxAtStartPar
\sphinxcode{\sphinxupquote{\_\_eq\_\_()}}
\end{quote}

\sphinxAtStartPar
\sphinxstylestrong{Description}
\begin{quote}

\sphinxAtStartPar
Overload the == operator to build a {\hyperref[\detokenize{pyapiref:chappyapi-mconstrbuilder}]{\sphinxcrossref{\DUrole{std,std-ref}{MConstrBuilder Class}}}} object,
which can be passed as the first argument to \sphinxtitleref{Model.addConstrs}.
\end{quote}

\sphinxAtStartPar
\sphinxstylestrong{Return value}
\begin{quote}

\sphinxAtStartPar
a {\hyperref[\detokenize{pyapiref:chappyapi-mconstrbuilder}]{\sphinxcrossref{\DUrole{std,std-ref}{MConstrBuilder Class}}}} object.
\end{quote}

\sphinxAtStartPar
\sphinxstylestrong{Example}

\begin{sphinxVerbatim}[commandchars=\\\{\}]
\PYG{n}{model}\PYG{o}{.}\PYG{n}{addConstrs}\PYG{p}{(}\PYG{n}{A} \PYG{o}{@} \PYG{n}{x} \PYG{o}{==} \PYG{l+m+mf}{1.0}\PYG{p}{)}
\end{sphinxVerbatim}
\end{quote}

\subsubsection{MLinExpr.\_\_ge\_\_()}
\label{\detokenize{pyapiref:mlinexpr-ge}}\begin{quote}

\sphinxAtStartPar
\sphinxstylestrong{Synopsis}
\begin{quote}

\sphinxAtStartPar
\sphinxcode{\sphinxupquote{\_\_ge\_\_()}}
\end{quote}

\sphinxAtStartPar
\sphinxstylestrong{Description}
\begin{quote}

\sphinxAtStartPar
Overload the \textgreater{}= operator to build a {\hyperref[\detokenize{pyapiref:chappyapi-mconstrbuilder}]{\sphinxcrossref{\DUrole{std,std-ref}{MConstrBuilder Class}}}} object,
which can be passed as the first argument to \sphinxtitleref{Model.addConstrs}.
\end{quote}

\sphinxAtStartPar
\sphinxstylestrong{Return value}
\begin{quote}

\sphinxAtStartPar
a {\hyperref[\detokenize{pyapiref:chappyapi-mconstrbuilder}]{\sphinxcrossref{\DUrole{std,std-ref}{MConstrBuilder Class}}}} object.
\end{quote}

\sphinxAtStartPar
\sphinxstylestrong{Example}

\begin{sphinxVerbatim}[commandchars=\\\{\}]
\PYG{n}{model}\PYG{o}{.}\PYG{n}{addConstrs}\PYG{p}{(}\PYG{n}{A} \PYG{o}{@} \PYG{n}{x} \PYG{o}{\PYGZgt{}}\PYG{o}{=} \PYG{l+m+mf}{1.0}\PYG{p}{)}
\end{sphinxVerbatim}
\end{quote}

\subsubsection{MLinExpr.\_\_le\_\_()}
\label{\detokenize{pyapiref:mlinexpr-le}}\begin{quote}

\sphinxAtStartPar
\sphinxstylestrong{Synopsis}
\begin{quote}

\sphinxAtStartPar
\sphinxcode{\sphinxupquote{\_\_le\_\_()}}
\end{quote}

\sphinxAtStartPar
\sphinxstylestrong{Description}
\begin{quote}

\sphinxAtStartPar
Overload the \textless{}= operator to build a {\hyperref[\detokenize{pyapiref:chappyapi-mconstrbuilder}]{\sphinxcrossref{\DUrole{std,std-ref}{MConstrBuilder Class}}}} object,
which can be passed as the first argument to \sphinxtitleref{Model.addConstrs}.
\end{quote}

\sphinxAtStartPar
\sphinxstylestrong{Return value}
\begin{quote}

\sphinxAtStartPar
a {\hyperref[\detokenize{pyapiref:chappyapi-mconstrbuilder}]{\sphinxcrossref{\DUrole{std,std-ref}{MConstrBuilder Class}}}} object.
\end{quote}

\sphinxAtStartPar
\sphinxstylestrong{Example}

\begin{sphinxVerbatim}[commandchars=\\\{\}]
\PYG{n}{model}\PYG{o}{.}\PYG{n}{addConstrs}\PYG{p}{(}\PYG{n}{A} \PYG{o}{@} \PYG{n}{x} \PYG{o}{\PYGZlt{}}\PYG{o}{=} \PYG{l+m+mf}{1.0}\PYG{p}{)}
\end{sphinxVerbatim}
\end{quote}

\subsection{MQuadExpr Class}
\label{\detokenize{pyapiref:mquadexpr-class}}\label{\detokenize{pyapiref:chappyapi-mquadexpr}}
\sphinxAtStartPar
The MQuadExpr class is used in COPT to construct multi\sphinxhyphen{}dimensional quadratic expressions and supports NumPy’s multi\sphinxhyphen{}dimensional array operations.
The MQuadExpr object with an initial value of 0.0 can be generated by the class generation method \sphinxcode{\sphinxupquote{zeros()}} , or by pairing two {\hyperref[\detokenize{pyapiref:chappyapi-mvar}]{\sphinxcrossref{\DUrole{std,std-ref}{MVar Class}}}} objects are generated by (matrix) multiplication.
It supports a unified set of indexing and slicing rules.
For detailed usage, see {\hyperref[\detokenize{matrix:chapmatrixindex}]{\sphinxcrossref{\DUrole{std,std-ref}{Matrix Modeling: Indexing and Slicing}}}}.

\sphinxAtStartPar
The following member methods are provided:

\subsubsection{MQuadExpr.zeros()}
\label{\detokenize{pyapiref:mquadexpr-zeros}}\begin{quote}

\sphinxAtStartPar
\sphinxstylestrong{Synopsis}
\begin{quote}

\sphinxAtStartPar
\sphinxcode{\sphinxupquote{zeros(shape)}}
\end{quote}

\sphinxAtStartPar
\sphinxstylestrong{Description}
\begin{quote}

\sphinxAtStartPar
This is the class generation method and can be called directly without an MQuadExpr object.
\end{quote}

\sphinxAtStartPar
\sphinxstylestrong{Arguments}
\begin{quote}

\sphinxAtStartPar
\sphinxcode{\sphinxupquote{shape}}
\begin{quote}

\sphinxAtStartPar
The value is an integer, or a tuple of integers. which represents the shape of the new MQuadExpr object.
\end{quote}
\end{quote}

\sphinxAtStartPar
\sphinxstylestrong{Return value}
\begin{quote}

\sphinxAtStartPar
new MQuadExpr object.
\end{quote}

\sphinxAtStartPar
\sphinxstylestrong{Example}

\begin{sphinxVerbatim}[commandchars=\\\{\}]
\PYG{n}{mqx} \PYG{o}{=} \PYG{n}{MQuadExpr}\PYG{o}{.}\PYG{n}{zeros}\PYG{p}{(}\PYG{p}{(}\PYG{l+m+mi}{2}\PYG{p}{,}\PYG{l+m+mi}{3}\PYG{p}{)}\PYG{p}{)} \PYG{c+c1}{\PYGZsh{} shape = (2, 3)}
\PYG{n}{x} \PYG{o}{=} \PYG{n}{model}\PYG{o}{.}\PYG{n}{addVar}\PYG{p}{(}\PYG{p}{)}
\PYG{n}{mqx} \PYG{o}{+}\PYG{o}{=} \PYG{l+m+mf}{2.0} \PYG{o}{*} \PYG{n}{x} \PYG{o}{*} \PYG{n}{x}           \PYG{c+c1}{\PYGZsh{} broadcast scalar}
\PYG{n}{mqx} \PYG{o}{+}\PYG{o}{=} \PYG{n}{model}\PYG{o}{.}\PYG{n}{addMVar}\PYG{p}{(}\PYG{l+m+mi}{3}\PYG{p}{)}      \PYG{c+c1}{\PYGZsh{} broadcast MVar of shape (3,)}
\end{sphinxVerbatim}
\end{quote}

\subsubsection{MQuadExpr.clear()}
\label{\detokenize{pyapiref:mquadexpr-clear}}\begin{quote}

\sphinxAtStartPar
\sphinxstylestrong{Synopsis}
\begin{quote}

\sphinxAtStartPar
\sphinxcode{\sphinxupquote{clear()}}
\end{quote}

\sphinxAtStartPar
\sphinxstylestrong{Description}
\begin{quote}

\sphinxAtStartPar
Resets each element of the {\hyperref[\detokenize{pyapiref:chappyapi-mquadexpr}]{\sphinxcrossref{\DUrole{std,std-ref}{MQuadExpr Class}}}} object to 0.0.
\end{quote}

\sphinxAtStartPar
\sphinxstylestrong{Example}

\begin{sphinxVerbatim}[commandchars=\\\{\}]
\PYG{n}{ma} \PYG{o}{=} \PYG{n}{model}\PYG{o}{.}\PYG{n}{addMVar}\PYG{p}{(}\PYG{l+m+mi}{3}\PYG{p}{,} \PYG{n}{nameprefix}\PYG{o}{=}\PYG{l+s+s1}{\PYGZsq{}}\PYG{l+s+s1}{a}\PYG{l+s+s1}{\PYGZsq{}}\PYG{p}{)}
\PYG{n}{mb} \PYG{o}{=} \PYG{n}{model}\PYG{o}{.}\PYG{n}{addMVar}\PYG{p}{(}\PYG{l+m+mi}{3}\PYG{p}{,} \PYG{n}{nameprefix}\PYG{o}{=}\PYG{l+s+s1}{\PYGZsq{}}\PYG{l+s+s1}{b}\PYG{l+s+s1}{\PYGZsq{}}\PYG{p}{)}
\PYG{n}{mqx} \PYG{o}{=} \PYG{n}{ma} \PYG{o}{*} \PYG{n}{mb}    \PYG{c+c1}{\PYGZsh{} elementwise multiply, shape = (3,)}
\PYG{n}{mqx}\PYG{o}{.}\PYG{n}{clear}\PYG{p}{(}\PYG{p}{)}
\PYG{n+nb}{print}\PYG{p}{(}\PYG{n}{mqx}\PYG{p}{)}       \PYG{c+c1}{\PYGZsh{} result is [0.0, 0.0, 0.0]}
\end{sphinxVerbatim}
\end{quote}

\subsubsection{MQuadExpr.clone()}
\label{\detokenize{pyapiref:mquadexpr-clone}}\begin{quote}

\sphinxAtStartPar
\sphinxstylestrong{Synopsis}
\begin{quote}

\sphinxAtStartPar
\sphinxcode{\sphinxupquote{clone()}}
\end{quote}

\sphinxAtStartPar
\sphinxstylestrong{Description}
\begin{quote}

\sphinxAtStartPar
Deep\sphinxhyphen{}copy a {\hyperref[\detokenize{pyapiref:chappyapi-mquadexpr}]{\sphinxcrossref{\DUrole{std,std-ref}{MQuadExpr Class}}}} object.
\end{quote}

\sphinxAtStartPar
\sphinxstylestrong{Return value}
\begin{quote}

\sphinxAtStartPar
new MQuadExpr object
\end{quote}

\sphinxAtStartPar
\sphinxstylestrong{Example}

\begin{sphinxVerbatim}[commandchars=\\\{\}]
\PYG{n}{mx} \PYG{o}{=} \PYG{n}{model}\PYG{o}{.}\PYG{n}{addMVar}\PYG{p}{(}\PYG{p}{(}\PYG{l+m+mi}{3}\PYG{p}{,} \PYG{l+m+mi}{3}\PYG{p}{)}\PYG{p}{,} \PYG{n}{nameprefix}\PYG{o}{=}\PYG{l+s+s1}{\PYGZsq{}}\PYG{l+s+s1}{mx}\PYG{l+s+s1}{\PYGZsq{}}\PYG{p}{)}
\PYG{n}{mqx} \PYG{o}{=} \PYG{l+m+mf}{2.0} \PYG{o}{*} \PYG{n}{mx} \PYG{o}{@} \PYG{n}{mx}      \PYG{c+c1}{\PYGZsh{} matrix multiply, shape = (3, 3)}
\PYG{n}{mqx\PYGZus{}copy} \PYG{o}{=} \PYG{n}{mqx}\PYG{o}{.}\PYG{n}{clone}\PYG{p}{(}\PYG{p}{)}
\PYG{n}{mqx\PYGZus{}copy}\PYG{o}{.}\PYG{n}{clear}\PYG{p}{(}\PYG{p}{)}
\PYG{n+nb}{print}\PYG{p}{(}\PYG{n}{mqx}\PYG{p}{)}               \PYG{c+c1}{\PYGZsh{} mqx is untouched}
\end{sphinxVerbatim}
\end{quote}

\subsubsection{MQuadExpr.getValue()}
\label{\detokenize{pyapiref:mquadexpr-getvalue}}\begin{quote}

\sphinxAtStartPar
\sphinxstylestrong{Synopsis}
\begin{quote}

\sphinxAtStartPar
\sphinxcode{\sphinxupquote{getValue()}}
\end{quote}

\sphinxAtStartPar
\sphinxstylestrong{Description}
\begin{quote}

\sphinxAtStartPar
Get the evaluation of each linear expression within the {\hyperref[\detokenize{pyapiref:chappyapi-mquadexpr}]{\sphinxcrossref{\DUrole{std,std-ref}{MQuadExpr Class}}}} object.
\end{quote}

\sphinxAtStartPar
\sphinxstylestrong{Return value}
\begin{quote}

\sphinxAtStartPar
Returns a NumPy ndarray of the same dimensions as the MQuadExpr object whose elements are the evaluations of the corresponding expression.
\end{quote}

\sphinxAtStartPar
\sphinxstylestrong{Example}

\begin{sphinxVerbatim}[commandchars=\\\{\}]
\PYG{n}{A} \PYG{o}{=} \PYG{n}{np}\PYG{o}{.}\PYG{n}{eye}\PYG{p}{(}\PYG{l+m+mi}{3}\PYG{p}{)}
\PYG{n}{mx} \PYG{o}{=} \PYG{n}{m}\PYG{o}{.}\PYG{n}{addMVar}\PYG{p}{(}\PYG{l+m+mi}{3}\PYG{p}{,} \PYG{n}{nameprefix}\PYG{o}{=}\PYG{l+s+s2}{\PYGZdq{}}\PYG{l+s+s2}{mx}\PYG{l+s+s2}{\PYGZdq{}}\PYG{p}{)}
\PYG{n}{mqx} \PYG{o}{=} \PYG{n}{mx} \PYG{o}{@} \PYG{n}{A} \PYG{o}{@} \PYG{n}{mx}        \PYG{c+c1}{\PYGZsh{} 0\PYGZhy{}D MQuadExpr, shape = ()}
\PYG{n}{m}\PYG{o}{.}\PYG{n}{addQConstr}\PYG{p}{(}\PYG{n}{mqx} \PYG{o}{\PYGZlt{}}\PYG{o}{=} \PYG{l+m+mf}{9.0}\PYG{p}{)}
\PYG{n}{m}\PYG{o}{.}\PYG{n}{solve}\PYG{p}{(}\PYG{p}{)}
\PYG{n+nb}{print}\PYG{p}{(}\PYG{n}{mqx}\PYG{o}{.}\PYG{n}{getValue}\PYG{p}{(}\PYG{p}{)}\PYG{p}{)}
\end{sphinxVerbatim}
\end{quote}

\subsubsection{MQuadExpr.item()}
\label{\detokenize{pyapiref:mquadexpr-item}}\begin{quote}

\sphinxAtStartPar
\sphinxstylestrong{Synopsis}
\begin{quote}

\sphinxAtStartPar
\sphinxcode{\sphinxupquote{item()}}
\end{quote}

\sphinxAtStartPar
\sphinxstylestrong{Description}
\begin{quote}

\sphinxAtStartPar
Get the constraint object in the 0\sphinxhyphen{}dimensional MQuadExpr. Raises a ValueError exception if the MQuadExpr object is not 0\sphinxhyphen{}dimensional.
\end{quote}

\sphinxAtStartPar
\sphinxstylestrong{Return value}
\begin{quote}

\sphinxAtStartPar
Returns the linear constraint object.
\end{quote}

\sphinxAtStartPar
\sphinxstylestrong{Example}

\begin{sphinxVerbatim}[commandchars=\\\{\}]
\PYG{n}{x} \PYG{o}{=} \PYG{n}{m}\PYG{o}{.}\PYG{n}{addVar}\PYG{p}{(}\PYG{p}{)}
\PYG{n}{mqx} \PYG{o}{=} \PYG{n}{MQuadExpr}\PYG{o}{.}\PYG{n}{zeros}\PYG{p}{(}\PYG{l+m+mi}{3}\PYG{p}{)} \PYG{o}{+} \PYG{n}{x} \PYG{o}{*} \PYG{n}{x}
\PYG{n+nb}{print}\PYG{p}{(}\PYG{n}{mqx}\PYG{p}{[}\PYG{l+m+mi}{1}\PYG{p}{]}\PYG{o}{.}\PYG{n}{item}\PYG{p}{(}\PYG{p}{)}\PYG{p}{)} \PYG{c+c1}{\PYGZsh{} Return QuadExpr(x * x)}
\end{sphinxVerbatim}
\end{quote}

\subsubsection{MQuadExpr.reshape()}
\label{\detokenize{pyapiref:mquadexpr-reshape}}\begin{quote}

\sphinxAtStartPar
\sphinxstylestrong{Synopsis}
\begin{quote}

\sphinxAtStartPar
\sphinxcode{\sphinxupquote{reshape(shape, order=\textquotesingle{}C\textquotesingle{})}}
\end{quote}

\sphinxAtStartPar
\sphinxstylestrong{Description}
\begin{quote}

\sphinxAtStartPar
Returns a new MQuadExpr object whose elements are unchanged but whose shape is transformed by the argument shape.
\end{quote}

\sphinxAtStartPar
\sphinxstylestrong{Arguments}
\begin{quote}

\sphinxAtStartPar
\sphinxcode{\sphinxupquote{shape}}
\begin{quote}

\sphinxAtStartPar
The value is an integer, or a tuple of integers. which represents the shape of the new MQuadExpr object.
\end{quote}

\sphinxAtStartPar
\sphinxcode{\sphinxupquote{order}}
\begin{quote}

\sphinxAtStartPar
Optional parameter, the default is the character ‘C’, which means it is compatible with the C language, that is, it is stored in rows; it can also be set to the character ‘F’, that is, it is stored in columns, and it is compatible with the Fortune language.
\end{quote}
\end{quote}

\sphinxAtStartPar
\sphinxstylestrong{Return value}
\begin{quote}

\sphinxAtStartPar
Returns a new MQuadExpr object with the same elements as the original MQuadExpr object but with a different shape.
\end{quote}

\sphinxAtStartPar
\sphinxstylestrong{Example}

\begin{sphinxVerbatim}[commandchars=\\\{\}]
\PYG{n}{mqx} \PYG{o}{=} \PYG{n}{MQuadExpr}\PYG{o}{.}\PYG{n}{zeros}\PYG{p}{(}\PYG{l+m+mi}{6}\PYG{p}{)}
\PYG{n}{mqx\PYGZus{}2x3} \PYG{o}{=} \PYG{n}{mqx}\PYG{o}{.}\PYG{n}{reshape}\PYG{p}{(}\PYG{p}{(}\PYG{l+m+mi}{2}\PYG{p}{,} \PYG{l+m+mi}{3}\PYG{p}{)}\PYG{p}{)}
\end{sphinxVerbatim}
\end{quote}

\subsubsection{MQuadExpr.sum()}
\label{\detokenize{pyapiref:mquadexpr-sum}}\begin{quote}

\sphinxAtStartPar
\sphinxstylestrong{Synopsis}
\begin{quote}

\sphinxAtStartPar
\sphinxcode{\sphinxupquote{sum(axis=None)}}
\end{quote}

\sphinxAtStartPar
\sphinxstylestrong{Description}
\begin{quote}

\sphinxAtStartPar
Sums the variables in the MQuadExpr object, returning a new {\hyperref[\detokenize{pyapiref:chappyapi-mquadexpr}]{\sphinxcrossref{\DUrole{std,std-ref}{MQuadExpr Class}}}} object.
\end{quote}

\sphinxAtStartPar
\sphinxstylestrong{Arguments}
\begin{quote}

\sphinxAtStartPar
\sphinxcode{\sphinxupquote{axis}}
\begin{quote}

\sphinxAtStartPar
Optional integer parameter, the default value is None, that is, to sum up variables one by one. Otherwise, sum over the given axis.
\end{quote}
\end{quote}

\sphinxAtStartPar
\sphinxstylestrong{Return value}
\begin{quote}

\sphinxAtStartPar
Returns an MQuadExpr object representing the sum of the corresponding linear expressions.
\end{quote}

\sphinxAtStartPar
\sphinxstylestrong{Example}

\begin{sphinxVerbatim}[commandchars=\\\{\}]
\PYG{n}{ma} \PYG{o}{=} \PYG{n}{model}\PYG{o}{.}\PYG{n}{addMVar}\PYG{p}{(}\PYG{p}{(}\PYG{l+m+mi}{2}\PYG{p}{,} \PYG{l+m+mi}{3}\PYG{p}{)}\PYG{p}{,} \PYG{n}{nameprefix}\PYG{o}{=}\PYG{l+s+s1}{\PYGZsq{}}\PYG{l+s+s1}{ma}\PYG{l+s+s1}{\PYGZsq{}}\PYG{p}{)}
\PYG{n}{mb} \PYG{o}{=} \PYG{n}{model}\PYG{o}{.}\PYG{n}{addMVar}\PYG{p}{(}\PYG{p}{(}\PYG{l+m+mi}{3}\PYG{p}{,} \PYG{l+m+mi}{2}\PYG{p}{)}\PYG{p}{,} \PYG{n}{nameprefix}\PYG{o}{=}\PYG{l+s+s1}{\PYGZsq{}}\PYG{l+s+s1}{mb}\PYG{l+s+s1}{\PYGZsq{}}\PYG{p}{)}
\PYG{n}{mqx} \PYG{o}{=} \PYG{n}{ma} \PYG{o}{@} \PYG{n}{mb}
\PYG{n}{sum\PYGZus{}all} \PYG{o}{=} \PYG{n}{mqx}\PYG{o}{.}\PYG{n}{sum}\PYG{p}{(}\PYG{p}{)} \PYG{c+c1}{\PYGZsh{} Return 0\PYGZhy{}dimensional MQuadExpr object}
\PYG{n}{sum\PYGZus{}row} \PYG{o}{=} \PYG{n}{mqx}\PYG{o}{.}\PYG{n}{sum}\PYG{p}{(}\PYG{n}{axis} \PYG{o}{=} \PYG{l+m+mi}{0}\PYG{p}{)} \PYG{c+c1}{\PYGZsh{} Return a 1\PYGZhy{}dimensional MQuadExpr object with a shape of (2, )}
\end{sphinxVerbatim}
\end{quote}

\subsubsection{MQuadExpr.tolist()}
\label{\detokenize{pyapiref:mquadexpr-tolist}}\begin{quote}

\sphinxAtStartPar
\sphinxstylestrong{Synopsis}
\begin{quote}

\sphinxAtStartPar
\sphinxcode{\sphinxupquote{tolist()}}
\end{quote}

\sphinxAtStartPar
\sphinxstylestrong{Description}
\begin{quote}

\sphinxAtStartPar
Converts an MQuadExpr object to a one\sphinxhyphen{}dimensional list whose elements are linear expressions.
\end{quote}

\sphinxAtStartPar
\sphinxstylestrong{Return value}
\begin{quote}

\sphinxAtStartPar
Return a 1D list containing {\hyperref[\detokenize{pyapiref:chappyapi-linexpr}]{\sphinxcrossref{\DUrole{std,std-ref}{LinExpr Class}}}}.
\end{quote}

\sphinxAtStartPar
\sphinxstylestrong{Example}

\begin{sphinxVerbatim}[commandchars=\\\{\}]
\PYG{n+nb}{print}\PYG{p}{(}\PYG{n}{MQuadExpr}\PYG{o}{.}\PYG{n}{zeros}\PYG{p}{(}\PYG{p}{(}\PYG{l+m+mi}{2}\PYG{p}{,}\PYG{l+m+mi}{3}\PYG{p}{)}\PYG{p}{)}\PYG{o}{.}\PYG{n}{tolist}\PYG{p}{(}\PYG{p}{)}\PYG{p}{)} \PYG{c+c1}{\PYGZsh{} a list of length 6}
\end{sphinxVerbatim}
\end{quote}

\subsubsection{MQuadExpr.transpose()}
\label{\detokenize{pyapiref:mquadexpr-transpose}}\begin{quote}

\sphinxAtStartPar
\sphinxstylestrong{Synopsis}
\begin{quote}

\sphinxAtStartPar
\sphinxcode{\sphinxupquote{transpose()}}
\end{quote}

\sphinxAtStartPar
\sphinxstylestrong{Description}
\begin{quote}

\sphinxAtStartPar
Generates a new MQuadExpr object that is the transpose of the original MQuadExpr object.
\end{quote}

\sphinxAtStartPar
\sphinxstylestrong{Return value}
\begin{quote}

\sphinxAtStartPar
Returns the transposed MQuadExpr object.
\end{quote}

\sphinxAtStartPar
\sphinxstylestrong{Example}

\begin{sphinxVerbatim}[commandchars=\\\{\}]
\PYG{n}{mqx} \PYG{o}{=} \PYG{n}{MQuadExpr}\PYG{o}{.}\PYG{n}{zeros}\PYG{p}{(}\PYG{p}{(}\PYG{l+m+mi}{2}\PYG{p}{,}\PYG{l+m+mi}{3}\PYG{p}{)}\PYG{p}{)}
\PYG{n+nb}{print}\PYG{p}{(}\PYG{n}{mqx}\PYG{o}{.}\PYG{n}{transpose}\PYG{p}{(}\PYG{p}{)}\PYG{o}{.}\PYG{n}{shape}\PYG{p}{)} \PYG{c+c1}{\PYGZsh{} shape = (3, 2)}
\end{sphinxVerbatim}
\end{quote}

\subsubsection{MQuadExpr.ndim}
\label{\detokenize{pyapiref:mquadexpr-ndim}}\begin{quote}

\sphinxAtStartPar
\sphinxstylestrong{Synopsis}
\begin{quote}

\sphinxAtStartPar
\sphinxcode{\sphinxupquote{ndim}}
\end{quote}

\sphinxAtStartPar
\sphinxstylestrong{Description}
\begin{quote}

\sphinxAtStartPar
{\hyperref[\detokenize{pyapiref:chappyapi-mquadexpr}]{\sphinxcrossref{\DUrole{std,std-ref}{MQuadExpr Class}}}} Dimensions of the object.
\end{quote}

\sphinxAtStartPar
\sphinxstylestrong{Return value}
\begin{quote}

\sphinxAtStartPar
Integer value.
\end{quote}

\sphinxAtStartPar
\sphinxstylestrong{Example}

\begin{sphinxVerbatim}[commandchars=\\\{\}]
\PYG{n}{mqx} \PYG{o}{=} \PYG{n}{MQuadExpr}\PYG{o}{.}\PYG{n}{zeros}\PYG{p}{(}\PYG{p}{(}\PYG{l+m+mi}{2}\PYG{p}{,}\PYG{l+m+mi}{3}\PYG{p}{)}\PYG{p}{)}
\PYG{n+nb}{print}\PYG{p}{(}\PYG{n}{mqx}\PYG{o}{.}\PYG{n}{ndim}\PYG{p}{)}  \PYG{c+c1}{\PYGZsh{} ndim = 2}
\end{sphinxVerbatim}
\end{quote}

\subsubsection{MQuadExpr.shape}
\label{\detokenize{pyapiref:mquadexpr-shape}}\begin{quote}

\sphinxAtStartPar
\sphinxstylestrong{Synopsis}
\begin{quote}

\sphinxAtStartPar
\sphinxcode{\sphinxupquote{shape}}
\end{quote}

\sphinxAtStartPar
\sphinxstylestrong{Description}
\begin{quote}

\sphinxAtStartPar
{\hyperref[\detokenize{pyapiref:chappyapi-mquadexpr}]{\sphinxcrossref{\DUrole{std,std-ref}{MQuadExpr Class}}}} The shape of the object.
\end{quote}

\sphinxAtStartPar
\sphinxstylestrong{Return value}
\begin{quote}

\sphinxAtStartPar
Integer tuple.
\end{quote}

\sphinxAtStartPar
\sphinxstylestrong{Example}

\begin{sphinxVerbatim}[commandchars=\\\{\}]
\PYG{n+nb}{print}\PYG{p}{(}\PYG{n}{MQuadExpr}\PYG{o}{.}\PYG{n}{zeros}\PYG{p}{(}\PYG{p}{(}\PYG{l+m+mi}{2}\PYG{p}{,}\PYG{l+m+mi}{3}\PYG{p}{)}\PYG{p}{)}\PYG{o}{.}\PYG{n}{shape}\PYG{p}{)} \PYG{c+c1}{\PYGZsh{} shape = (2, 3)}
\end{sphinxVerbatim}
\end{quote}

\subsubsection{MQuadExpr.size}
\label{\detokenize{pyapiref:mquadexpr-size}}\begin{quote}

\sphinxAtStartPar
\sphinxstylestrong{Synopsis}
\begin{quote}

\sphinxAtStartPar
\sphinxcode{\sphinxupquote{size}}
\end{quote}

\sphinxAtStartPar
\sphinxstylestrong{Description}
\begin{quote}

\sphinxAtStartPar
The number of elements of {\hyperref[\detokenize{pyapiref:chappyapi-mquadexpr}]{\sphinxcrossref{\DUrole{std,std-ref}{MQuadExpr Class}}}} object.
\end{quote}

\sphinxAtStartPar
\sphinxstylestrong{Return value}
\begin{quote}

\sphinxAtStartPar
Integer value.
\end{quote}

\sphinxAtStartPar
\sphinxstylestrong{Example}

\begin{sphinxVerbatim}[commandchars=\\\{\}]
\PYG{n}{mqx} \PYG{o}{=} \PYG{n}{MQuadExpr}\PYG{o}{.}\PYG{n}{zeros}\PYG{p}{(}\PYG{p}{(}\PYG{l+m+mi}{2}\PYG{p}{,}\PYG{l+m+mi}{3}\PYG{p}{)}\PYG{p}{)}
\PYG{n+nb}{print}\PYG{p}{(}\PYG{n}{mqx}\PYG{o}{.}\PYG{n}{size}\PYG{p}{)}  \PYG{c+c1}{\PYGZsh{} size= 6}
\end{sphinxVerbatim}
\end{quote}

\subsubsection{MQuadExpr.T}
\label{\detokenize{pyapiref:mquadexpr-t}}\begin{quote}

\sphinxAtStartPar
\sphinxstylestrong{Synopsis}
\begin{quote}

\sphinxAtStartPar
\sphinxcode{\sphinxupquote{T}}
\end{quote}

\sphinxAtStartPar
\sphinxstylestrong{Description}
\begin{quote}

\sphinxAtStartPar
Transpose of {\hyperref[\detokenize{pyapiref:chappyapi-mquadexpr}]{\sphinxcrossref{\DUrole{std,std-ref}{MQuadExpr Class}}}} object. Similar to the class method transpose().
\end{quote}

\sphinxAtStartPar
\sphinxstylestrong{Return value}
\begin{quote}

\sphinxAtStartPar
Returns the transposed MQuadExpr object.
\end{quote}

\sphinxAtStartPar
\sphinxstylestrong{Example}

\begin{sphinxVerbatim}[commandchars=\\\{\}]
\PYG{n}{mqx} \PYG{o}{=} \PYG{n}{MQuadExpr}\PYG{o}{.}\PYG{n}{zeros}\PYG{p}{(}\PYG{p}{(}\PYG{l+m+mi}{2}\PYG{p}{,}\PYG{l+m+mi}{3}\PYG{p}{)}\PYG{p}{)}
\PYG{n+nb}{print}\PYG{p}{(}\PYG{n}{mqx}\PYG{o}{.}\PYG{n}{T}\PYG{o}{.}\PYG{n}{shape}\PYG{p}{)} \PYG{c+c1}{\PYGZsh{} shape = (3, 2)}
\end{sphinxVerbatim}
\end{quote}

\subsubsection{MQuadExpr.\_\_eq\_\_()}
\label{\detokenize{pyapiref:mquadexpr-eq}}\begin{quote}

\sphinxAtStartPar
\sphinxstylestrong{Synopsis}
\begin{quote}

\sphinxAtStartPar
\sphinxcode{\sphinxupquote{\_\_eq\_\_()}}
\end{quote}

\sphinxAtStartPar
\sphinxstylestrong{Description}
\begin{quote}

\sphinxAtStartPar
Overload the == operator to build a {\hyperref[\detokenize{pyapiref:chappyapi-mqconstrbuilder}]{\sphinxcrossref{\DUrole{std,std-ref}{MQConstrBuilder Class}}}} object,
which can be passed as the first argument to \sphinxtitleref{Model.addQConstr}.
\end{quote}

\sphinxAtStartPar
\sphinxstylestrong{Return value}
\begin{quote}

\sphinxAtStartPar
a {\hyperref[\detokenize{pyapiref:chappyapi-mqconstrbuilder}]{\sphinxcrossref{\DUrole{std,std-ref}{MQConstrBuilder Class}}}} object.
\end{quote}

\sphinxAtStartPar
\sphinxstylestrong{Example}

\begin{sphinxVerbatim}[commandchars=\\\{\}]
\PYG{n}{model}\PYG{o}{.}\PYG{n}{addQConstr}\PYG{p}{(}\PYG{n}{x} \PYG{o}{@} \PYG{n}{Q} \PYG{o}{@} \PYG{n}{y} \PYG{o}{==} \PYG{l+m+mf}{1.0}\PYG{p}{)}
\end{sphinxVerbatim}
\end{quote}

\subsubsection{MQuadExpr.\_\_ge\_\_()}
\label{\detokenize{pyapiref:mquadexpr-ge}}\begin{quote}

\sphinxAtStartPar
\sphinxstylestrong{Synopsis}
\begin{quote}

\sphinxAtStartPar
\sphinxcode{\sphinxupquote{\_\_ge\_\_()}}
\end{quote}

\sphinxAtStartPar
\sphinxstylestrong{Description}
\begin{quote}

\sphinxAtStartPar
Overload the \textgreater{}= operator to build a {\hyperref[\detokenize{pyapiref:chappyapi-mqconstrbuilder}]{\sphinxcrossref{\DUrole{std,std-ref}{MQConstrBuilder Class}}}} object,
which can be passed as the first argument to \sphinxtitleref{Model.addQConstr}.
\end{quote}

\sphinxAtStartPar
\sphinxstylestrong{Return value}
\begin{quote}

\sphinxAtStartPar
a {\hyperref[\detokenize{pyapiref:chappyapi-mqconstrbuilder}]{\sphinxcrossref{\DUrole{std,std-ref}{MQConstrBuilder Class}}}} object.
\end{quote}

\sphinxAtStartPar
\sphinxstylestrong{Example}

\begin{sphinxVerbatim}[commandchars=\\\{\}]
\PYG{n}{model}\PYG{o}{.}\PYG{n}{addQConstr}\PYG{p}{(}\PYG{n}{x} \PYG{o}{@} \PYG{n}{Q} \PYG{o}{@} \PYG{n}{y} \PYG{o}{\PYGZgt{}}\PYG{o}{=} \PYG{l+m+mf}{1.0}\PYG{p}{)}
\end{sphinxVerbatim}
\end{quote}

\subsubsection{MQuadExpr.\_\_le\_\_()}
\label{\detokenize{pyapiref:mquadexpr-le}}\begin{quote}

\sphinxAtStartPar
\sphinxstylestrong{Synopsis}
\begin{quote}

\sphinxAtStartPar
\sphinxcode{\sphinxupquote{\_\_le\_\_()}}
\end{quote}

\sphinxAtStartPar
\sphinxstylestrong{Description}
\begin{quote}

\sphinxAtStartPar
Overload the \textless{}= operator to build a {\hyperref[\detokenize{pyapiref:chappyapi-mqconstrbuilder}]{\sphinxcrossref{\DUrole{std,std-ref}{MQConstrBuilder Class}}}} object,
which can be passed as the first argument to \sphinxtitleref{Model.addQConstr}.
\end{quote}

\sphinxAtStartPar
\sphinxstylestrong{Return value}
\begin{quote}

\sphinxAtStartPar
a {\hyperref[\detokenize{pyapiref:chappyapi-mqconstrbuilder}]{\sphinxcrossref{\DUrole{std,std-ref}{MQConstrBuilder Class}}}} object.
\end{quote}

\sphinxAtStartPar
\sphinxstylestrong{Example}

\begin{sphinxVerbatim}[commandchars=\\\{\}]
\PYG{n}{model}\PYG{o}{.}\PYG{n}{addQConstr}\PYG{p}{(}\PYG{n}{x} \PYG{o}{@} \PYG{n}{Q} \PYG{o}{@} \PYG{n}{y} \PYG{o}{\PYGZlt{}}\PYG{o}{=} \PYG{l+m+mf}{1.0}\PYG{p}{)}
\end{sphinxVerbatim}
\end{quote}

\subsection{NdArray Class}
\label{\detokenize{pyapiref:ndarray-class}}\label{\detokenize{pyapiref:chappyapi-ndarray}}
\sphinxAtStartPar
The NdArray class is a built\sphinxhyphen{}in multi\sphinxhyphen{}dimensional array class in COPT.
It represents a table of elements of the same type, indexed by a tuple of integers.
It supports a unified set of indexing and slicing rules.
For detailed usage, see {\hyperref[\detokenize{matrix:chapmatrixindex}]{\sphinxcrossref{\DUrole{std,std-ref}{Matrix Modeling: Indexing and Slicing}}}}.

\sphinxAtStartPar
The following methods are provided:

\subsubsection{NdArray()}
\label{\detokenize{pyapiref:ndarray}}\begin{quote}

\sphinxAtStartPar
\sphinxstylestrong{Synopsis}
\begin{quote}

\sphinxAtStartPar
\sphinxcode{\sphinxupquote{NdArray(args=None, dtype=None, shape=None)}}
\end{quote}

\sphinxAtStartPar
\sphinxstylestrong{Description}
\begin{quote}

\sphinxAtStartPar
Create a {\hyperref[\detokenize{pyapiref:chappyapi-ndarray}]{\sphinxcrossref{\DUrole{std,std-ref}{NdArray Class}}}} object.
\end{quote}

\sphinxAtStartPar
\sphinxstylestrong{Return value}
\begin{quote}

\sphinxAtStartPar
Returns a NdArray object.
\end{quote}

\sphinxAtStartPar
\sphinxstylestrong{Example}

\begin{sphinxVerbatim}[commandchars=\\\{\}]
\PYG{c+c1}{\PYGZsh{} Create a NdArray object with a shape of 3x3 and initialize its elements to 0}
\PYG{n}{ndmat} \PYG{o}{=} \PYG{n}{NdArray}\PYG{p}{(}\PYG{n}{shape}\PYG{o}{=}\PYG{p}{(}\PYG{l+m+mi}{3}\PYG{p}{,} \PYG{l+m+mi}{3}\PYG{p}{)}\PYG{p}{)}
\end{sphinxVerbatim}
\end{quote}

\subsubsection{NdArray.item()}
\label{\detokenize{pyapiref:ndarray-item}}\begin{quote}

\sphinxAtStartPar
\sphinxstylestrong{Synopsis}
\begin{quote}

\sphinxAtStartPar
\sphinxcode{\sphinxupquote{item()}}
\end{quote}

\sphinxAtStartPar
\sphinxstylestrong{Description}
\begin{quote}

\sphinxAtStartPar
Gets the single element within a 0\sphinxhyphen{}dimensional NdArray object. If the NdArray object is not 0\sphinxhyphen{}dimensional, a \sphinxcode{\sphinxupquote{ValueError}} exception will be triggered.
\end{quote}

\sphinxAtStartPar
\sphinxstylestrong{Return value}
\begin{quote}

\sphinxAtStartPar
Returns the type of the elements in the 0\sphinxhyphen{}dimensional NdArray object. (For example: \sphinxcode{\sphinxupquote{"float"}} or \sphinxcode{\sphinxupquote{"int"}}, etc.)
\end{quote}

\sphinxAtStartPar
\sphinxstylestrong{Example}

\begin{sphinxVerbatim}[commandchars=\\\{\}]
\PYG{n}{ndmat} \PYG{o}{=} \PYG{n}{NdArray}\PYG{p}{(}\PYG{n}{args}\PYG{o}{=}\PYG{l+m+mf}{1.1}\PYG{p}{,} \PYG{n}{shape}\PYG{o}{=}\PYG{p}{(}\PYG{l+m+mi}{1}\PYG{p}{,}\PYG{p}{)}\PYG{p}{)}
\PYG{c+c1}{\PYGZsh{} Type of value is \PYGZdq{}float\PYGZdq{}}
\PYG{n}{value} \PYG{o}{=} \PYG{n}{ndmat}\PYG{o}{.}\PYG{n}{item}\PYG{p}{(}\PYG{p}{)}
\end{sphinxVerbatim}
\end{quote}

\subsubsection{NdArray.reshape()}
\label{\detokenize{pyapiref:ndarray-reshape}}\begin{quote}

\sphinxAtStartPar
\sphinxstylestrong{Synopsis}
\begin{quote}

\sphinxAtStartPar
\sphinxcode{\sphinxupquote{reshape(shape, order=\textquotesingle{}C\textquotesingle{})}}
\end{quote}

\sphinxAtStartPar
\sphinxstylestrong{Description}
\begin{quote}

\sphinxAtStartPar
Returns a new NdArray object whose elements remain unchanged where shape is transformed according to the shape in the arguments.
\end{quote}

\sphinxAtStartPar
\sphinxstylestrong{Arguments}
\begin{quote}

\sphinxAtStartPar
\sphinxcode{\sphinxupquote{shape}}
\begin{quote}

\sphinxAtStartPar
The value could be an integer or a tuple of integers, representing the shape of the new NdArray object.
\end{quote}

\sphinxAtStartPar
\sphinxcode{\sphinxupquote{order}}
\begin{quote}

\sphinxAtStartPar
Optional.
The default is the character ‘C’, which means it is compatible with C language (stored in rows). The current version does not yet support the character ‘F’ (stored in columns).
\end{quote}
\end{quote}

\sphinxAtStartPar
\sphinxstylestrong{Return value}
\begin{quote}

\sphinxAtStartPar
Returns a new NdArray object with the same elements as the original one but a different shape.
\end{quote}

\sphinxAtStartPar
\sphinxstylestrong{Example}

\begin{sphinxVerbatim}[commandchars=\\\{\}]
\PYG{n}{ndmat} \PYG{o}{=} \PYG{n}{NdArray}\PYG{p}{(}\PYG{n}{shape}\PYG{o}{=}\PYG{p}{(}\PYG{l+m+mi}{6}\PYG{p}{,}\PYG{p}{)}\PYG{p}{)}
\PYG{n}{ndmat\PYGZus{}2x3} \PYG{o}{=} \PYG{n}{ndmat}\PYG{o}{.}\PYG{n}{reshape}\PYG{p}{(}\PYG{p}{(}\PYG{l+m+mi}{2}\PYG{p}{,} \PYG{l+m+mi}{3}\PYG{p}{)}\PYG{p}{)}
\end{sphinxVerbatim}
\end{quote}

\subsubsection{NdArray.sum()}
\label{\detokenize{pyapiref:ndarray-sum}}\begin{quote}

\sphinxAtStartPar
\sphinxstylestrong{Synopsis}
\begin{quote}

\sphinxAtStartPar
\sphinxcode{\sphinxupquote{sum(axis=None)}}
\end{quote}

\sphinxAtStartPar
\sphinxstylestrong{Description}
\begin{quote}

\sphinxAtStartPar
Compute the sum of elements in NdArray along the given axis.

\sphinxAtStartPar
If \sphinxcode{\sphinxupquote{axis}} is not specified, all elements are summed;
otherwise, the summation is performed along the specified axis.
\end{quote}

\sphinxAtStartPar
\sphinxstylestrong{Arguments}
\begin{quote}

\sphinxAtStartPar
\sphinxcode{\sphinxupquote{axis}}
\begin{quote}

\sphinxAtStartPar
Optional integer argument specifying the axis.

\sphinxAtStartPar
The default value is \sphinxcode{\sphinxupquote{None}}.
\end{quote}
\end{quote}

\sphinxAtStartPar
\sphinxstylestrong{Return value}
\begin{quote}

\sphinxAtStartPar
If \sphinxcode{\sphinxupquote{axis}} is \sphinxcode{\sphinxupquote{None}}, returns a scalar result representing the sum of all elements.
The numeric type of the returned value is consistent with the element type of the NdArray.

\sphinxAtStartPar
If \sphinxcode{\sphinxupquote{axis}} is specified, returns a new NdArray object whose dimension
is reduced by one. The elements represent the sums computed along the given axis.
\end{quote}

\sphinxAtStartPar
\sphinxstylestrong{Example}

\begin{sphinxVerbatim}[commandchars=\\\{\}]
\PYG{n}{ndmat} \PYG{o}{=} \PYG{n}{NdArray}\PYG{p}{(}\PYG{n}{args}\PYG{o}{=}\PYG{l+m+mf}{1.1}\PYG{p}{,} \PYG{n}{shape}\PYG{o}{=}\PYG{p}{(}\PYG{l+m+mi}{2}\PYG{p}{,} \PYG{l+m+mi}{2}\PYG{p}{)}\PYG{p}{)}
\PYG{n}{sum\PYGZus{}all} \PYG{o}{=} \PYG{n}{ndmat}\PYG{o}{.}\PYG{n}{sum}\PYG{p}{(}\PYG{p}{)}        \PYG{c+c1}{\PYGZsh{} Returns a scalar result}
\PYG{n}{sum\PYGZus{}row} \PYG{o}{=} \PYG{n}{ndmat}\PYG{o}{.}\PYG{n}{sum}\PYG{p}{(}\PYG{n}{axis}\PYG{o}{=}\PYG{l+m+mi}{0}\PYG{p}{)}  \PYG{c+c1}{\PYGZsh{} Returns a 1\PYGZhy{}D NdArray}
\end{sphinxVerbatim}
\end{quote}

\subsubsection{NdArray.prod()}
\label{\detokenize{pyapiref:ndarray-prod}}\begin{quote}

\sphinxAtStartPar
\sphinxstylestrong{Synopsis}
\begin{quote}

\sphinxAtStartPar
\sphinxcode{\sphinxupquote{prod(axis=None)}}
\end{quote}

\sphinxAtStartPar
\sphinxstylestrong{Description}
\begin{quote}

\sphinxAtStartPar
Compute the product of elements in NdArray along the given axis.

\sphinxAtStartPar
If \sphinxcode{\sphinxupquote{axis}} is not specified, the product is taken over all elements;
otherwise, the product is computed along the specified axis.
\end{quote}

\sphinxAtStartPar
\sphinxstylestrong{Arguments}
\begin{quote}

\sphinxAtStartPar
\sphinxcode{\sphinxupquote{axis}}
\begin{quote}

\sphinxAtStartPar
Optional integer argument specifying the axis.

\sphinxAtStartPar
The default value is \sphinxcode{\sphinxupquote{None}}.
\end{quote}
\end{quote}

\sphinxAtStartPar
\sphinxstylestrong{Return value}
\begin{quote}

\sphinxAtStartPar
If \sphinxcode{\sphinxupquote{axis}} is \sphinxcode{\sphinxupquote{None}}, returns a scalar result representing the product of all elements.
The numeric type of the returned value is consistent with the element type of the NdArray.

\sphinxAtStartPar
If \sphinxcode{\sphinxupquote{axis}} is specified, returns a new NdArray object whose dimension
is reduced by one. The elements represent the products computed along the given axis.
\end{quote}

\sphinxAtStartPar
\sphinxstylestrong{Example}

\begin{sphinxVerbatim}[commandchars=\\\{\}]
\PYG{n}{ndmat} \PYG{o}{=} \PYG{n}{NdArray}\PYG{p}{(}\PYG{n}{args}\PYG{o}{=}\PYG{l+m+mf}{2.0}\PYG{p}{,} \PYG{n}{shape}\PYG{o}{=}\PYG{p}{(}\PYG{l+m+mi}{2}\PYG{p}{,} \PYG{l+m+mi}{2}\PYG{p}{)}\PYG{p}{)}
\PYG{n}{prod\PYGZus{}all} \PYG{o}{=} \PYG{n}{ndmat}\PYG{o}{.}\PYG{n}{prod}\PYG{p}{(}\PYG{p}{)}        \PYG{c+c1}{\PYGZsh{} Returns a scalar result}
\PYG{n}{prod\PYGZus{}row} \PYG{o}{=} \PYG{n}{ndmat}\PYG{o}{.}\PYG{n}{prod}\PYG{p}{(}\PYG{n}{axis}\PYG{o}{=}\PYG{l+m+mi}{0}\PYG{p}{)}  \PYG{c+c1}{\PYGZsh{} Returns a 1\PYGZhy{}D NdArray}
\end{sphinxVerbatim}
\end{quote}

\subsubsection{NdArray.dot()}
\label{\detokenize{pyapiref:ndarray-dot}}\begin{quote}

\sphinxAtStartPar
\sphinxstylestrong{Synopsis}
\begin{quote}

\sphinxAtStartPar
\sphinxcode{\sphinxupquote{dot(other)}}
\end{quote}

\sphinxAtStartPar
\sphinxstylestrong{Description}
\begin{quote}

\sphinxAtStartPar
Compute the dot product between the current NdArray and the given object.
\end{quote}

\sphinxAtStartPar
\sphinxstylestrong{Arguments}
\begin{quote}

\sphinxAtStartPar
\sphinxcode{\sphinxupquote{other}}
\begin{quote}

\sphinxAtStartPar
The object used for the dot product computation.

\sphinxAtStartPar
Acceptable values include an iterable object or an NdArray.
The length must match that of the current NdArray.
\end{quote}
\end{quote}

\sphinxAtStartPar
\sphinxstylestrong{Return value}
\begin{quote}

\sphinxAtStartPar
Returns a scalar result representing the dot product computation.
The numeric type of the returned value is consistent with the element type of the NdArray.
\end{quote}

\sphinxAtStartPar
\sphinxstylestrong{Example}

\begin{sphinxVerbatim}[commandchars=\\\{\}]
\PYG{n}{a} \PYG{o}{=} \PYG{n}{NdArray}\PYG{p}{(}\PYG{n}{args}\PYG{o}{=}\PYG{p}{[}\PYG{l+m+mf}{1.0}\PYG{p}{,} \PYG{l+m+mf}{2.0}\PYG{p}{,} \PYG{l+m+mf}{3.0}\PYG{p}{]}\PYG{p}{,} \PYG{n}{shape}\PYG{o}{=}\PYG{p}{(}\PYG{l+m+mi}{3}\PYG{p}{,} \PYG{p}{)}\PYG{p}{)}
\PYG{n}{b} \PYG{o}{=} \PYG{n}{NdArray}\PYG{p}{(}\PYG{n}{args}\PYG{o}{=}\PYG{p}{[}\PYG{l+m+mf}{4.0}\PYG{p}{,} \PYG{l+m+mf}{5.0}\PYG{p}{,} \PYG{l+m+mf}{6.0}\PYG{p}{]}\PYG{p}{,} \PYG{n}{shape}\PYG{o}{=}\PYG{p}{(}\PYG{l+m+mi}{3}\PYG{p}{,} \PYG{p}{)}\PYG{p}{)}

\PYG{n}{val} \PYG{o}{=} \PYG{n}{a}\PYG{o}{.}\PYG{n}{dot}\PYG{p}{(}\PYG{n}{b}\PYG{p}{)}  \PYG{c+c1}{\PYGZsh{} Returns a scalar result}
\end{sphinxVerbatim}
\end{quote}

\subsubsection{NdArray.tolist()}
\label{\detokenize{pyapiref:ndarray-tolist}}\begin{quote}

\sphinxAtStartPar
\sphinxstylestrong{Synopsis}
\begin{quote}

\sphinxAtStartPar
\sphinxcode{\sphinxupquote{tolist()}}
\end{quote}

\sphinxAtStartPar
\sphinxstylestrong{Description}
\begin{quote}

\sphinxAtStartPar
Convert a NdArray object into a list.
\end{quote}

\sphinxAtStartPar
\sphinxstylestrong{Return value}
\begin{quote}

\sphinxAtStartPar
Returns a list object.
\end{quote}

\sphinxAtStartPar
\sphinxstylestrong{Example}

\begin{sphinxVerbatim}[commandchars=\\\{\}]
\PYG{c+c1}{\PYGZsh{} Type of object mat\PYGZus{}tolist is \PYGZdq{}list\PYGZdq{}}
\PYG{n}{mat\PYGZus{}tolist} \PYG{o}{=} \PYG{n}{ndmat}\PYG{o}{.}\PYG{n}{tolist}\PYG{p}{(}\PYG{p}{)}
\end{sphinxVerbatim}
\end{quote}

\subsubsection{NdArray.tonumpy()}
\label{\detokenize{pyapiref:ndarray-tonumpy}}\begin{quote}

\sphinxAtStartPar
\sphinxstylestrong{Synopsis}
\begin{quote}

\sphinxAtStartPar
\sphinxcode{\sphinxupquote{tonumpy()}}
\end{quote}

\sphinxAtStartPar
\sphinxstylestrong{Description}
\begin{quote}

\sphinxAtStartPar
Convert a NdArray object to a NumPy ndarray object.
\end{quote}

\sphinxAtStartPar
\sphinxstylestrong{Return value}
\begin{quote}

\sphinxAtStartPar
Returns a numpy ndarray object.
\end{quote}

\sphinxAtStartPar
\sphinxstylestrong{Example}

\begin{sphinxVerbatim}[commandchars=\\\{\}]
\PYG{c+c1}{\PYGZsh{} Type of object mat\PYGZus{}tolist is \PYGZdq{}numpy.ndarray\PYGZdq{}}
\PYG{n}{mat\PYGZus{}tonumpy} \PYG{o}{=} \PYG{n}{ndmat}\PYG{o}{.}\PYG{n}{tonumpy}\PYG{p}{(}\PYG{p}{)}
\end{sphinxVerbatim}
\end{quote}

\subsubsection{NdArray.fill()}
\label{\detokenize{pyapiref:ndarray-fill}}\begin{quote}

\sphinxAtStartPar
\sphinxstylestrong{Synopsis}
\begin{quote}

\sphinxAtStartPar
\sphinxcode{\sphinxupquote{fill(value)}}
\end{quote}

\sphinxAtStartPar
\sphinxstylestrong{Description}
\begin{quote}

\sphinxAtStartPar
Fills each element in the NdArray object with the specified value.
\end{quote}

\sphinxAtStartPar
\sphinxstylestrong{Arguments}
\begin{quote}

\sphinxAtStartPar
\sphinxcode{\sphinxupquote{value}}
\begin{quote}

\sphinxAtStartPar
The new value of each element in the NdArray object.
\end{quote}
\end{quote}

\sphinxAtStartPar
\sphinxstylestrong{Return value}
\begin{quote}

\sphinxAtStartPar
Returns a NdArray object.
\end{quote}

\sphinxAtStartPar
\sphinxstylestrong{Example}

\begin{sphinxVerbatim}[commandchars=\\\{\}]
\PYG{n}{mat\PYGZus{}fillvalue} \PYG{o}{=} \PYG{n}{ndmat}\PYG{o}{.}\PYG{n}{fill}\PYG{p}{(}\PYG{l+m+mf}{100.0}\PYG{p}{)}
\end{sphinxVerbatim}
\end{quote}

\subsubsection{NdArray.expand()}
\label{\detokenize{pyapiref:ndarray-expand}}\begin{quote}

\sphinxAtStartPar
\sphinxstylestrong{Synopsis}
\begin{quote}

\sphinxAtStartPar
\sphinxcode{\sphinxupquote{expand(axis=0)}}
\end{quote}

\sphinxAtStartPar
\sphinxstylestrong{Description}
\begin{quote}

\sphinxAtStartPar
Expand the NdArray object into an N+1 dimensional shape on the axis.
\end{quote}

\sphinxAtStartPar
\sphinxstylestrong{Arguments}
\begin{quote}

\sphinxAtStartPar
\sphinxcode{\sphinxupquote{axis}}
\begin{quote}

\sphinxAtStartPar
The specified dimension, which defaults to 0 (the first dimension).
\end{quote}
\end{quote}

\sphinxAtStartPar
\sphinxstylestrong{Return value}
\begin{quote}

\sphinxAtStartPar
Returns a N+1 dimensional NdArray object.
\end{quote}

\sphinxAtStartPar
\sphinxstylestrong{Example}

\begin{sphinxVerbatim}[commandchars=\\\{\}]
\PYG{n}{mat\PYGZus{}1} \PYG{o}{=} \PYG{n}{ndmat}\PYG{o}{.}\PYG{n}{expand}\PYG{p}{(}\PYG{p}{)}
\end{sphinxVerbatim}
\end{quote}

\subsubsection{NdArray.squeeze()}
\label{\detokenize{pyapiref:ndarray-squeeze}}\begin{quote}

\sphinxAtStartPar
\sphinxstylestrong{Synopsis}
\begin{quote}

\sphinxAtStartPar
\sphinxcode{\sphinxupquote{squeeze(axis=0)}}
\end{quote}

\sphinxAtStartPar
\sphinxstylestrong{Description}
\begin{quote}

\sphinxAtStartPar
Reduce the NdArray object to an N\sphinxhyphen{}1 dimensional shape on the axis.
\end{quote}

\sphinxAtStartPar
\sphinxstylestrong{Arguments}
\begin{quote}

\sphinxAtStartPar
\sphinxcode{\sphinxupquote{axis}}
\begin{quote}

\sphinxAtStartPar
The specified dimension, default to 0 (the first dimension).
\end{quote}
\end{quote}

\sphinxAtStartPar
\sphinxstylestrong{Return value}
\begin{quote}

\sphinxAtStartPar
Returns a N\sphinxhyphen{}1 dimensional NdArray object.
\end{quote}

\sphinxAtStartPar
\sphinxstylestrong{Example}

\begin{sphinxVerbatim}[commandchars=\\\{\}]
\PYG{n}{mat\PYGZus{}1} \PYG{o}{=} \PYG{n}{ndmat}\PYG{o}{.}\PYG{n}{squeeze}\PYG{p}{(}\PYG{p}{)}
\end{sphinxVerbatim}
\end{quote}

\subsubsection{NdArray.flatten()}
\label{\detokenize{pyapiref:ndarray-flatten}}\begin{quote}

\sphinxAtStartPar
\sphinxstylestrong{Synopsis}
\begin{quote}

\sphinxAtStartPar
\sphinxcode{\sphinxupquote{flatten()}}
\end{quote}

\sphinxAtStartPar
\sphinxstylestrong{Description}
\begin{quote}

\sphinxAtStartPar
Expand a NdArray object into a one\sphinxhyphen{}dimensional shape.
\end{quote}

\sphinxAtStartPar
\sphinxstylestrong{Return value}
\begin{quote}

\sphinxAtStartPar
Returns a new one\sphinxhyphen{}dimensional NdArray object.
\end{quote}

\sphinxAtStartPar
\sphinxstylestrong{Example}

\begin{sphinxVerbatim}[commandchars=\\\{\}]
\PYG{n}{ndmat} \PYG{o}{=} \PYG{n}{NdArray}\PYG{p}{(}\PYG{n}{shape}\PYG{o}{=}\PYG{p}{(}\PYG{l+m+mi}{2}\PYG{p}{,} \PYG{l+m+mi}{2}\PYG{p}{)}\PYG{p}{)}
\PYG{c+c1}{\PYGZsh{} The shape of mat\PYGZus{}1 is (4,)}
\PYG{n}{mat\PYGZus{}1} \PYG{o}{=} \PYG{n}{ndmat}\PYG{o}{.}\PYG{n}{flatten}\PYG{p}{(}\PYG{p}{)}
\end{sphinxVerbatim}
\end{quote}

\subsubsection{NdArray.setItem()}
\label{\detokenize{pyapiref:ndarray-setitem}}\begin{quote}

\sphinxAtStartPar
\sphinxstylestrong{Synopsis}
\begin{quote}

\sphinxAtStartPar
\sphinxcode{\sphinxupquote{setItem(idx, value)}}
\end{quote}

\sphinxAtStartPar
\sphinxstylestrong{Description}
\begin{quote}

\sphinxAtStartPar
Sets the value of the element according to the given index in the NdArray object.
\end{quote}

\sphinxAtStartPar
\sphinxstylestrong{Arguments}
\begin{quote}

\sphinxAtStartPar
\sphinxcode{\sphinxupquote{idx}}
\begin{quote}

\sphinxAtStartPar
The specified one\sphinxhyphen{}dimensional index is the corresponding one after flattening NdArray into one dimension.
\end{quote}

\sphinxAtStartPar
\sphinxcode{\sphinxupquote{value}}
\begin{quote}

\sphinxAtStartPar
The new value of the specified element.
\end{quote}
\end{quote}

\sphinxAtStartPar
\sphinxstylestrong{Return value}
\begin{quote}

\sphinxAtStartPar
Returns a NdArray object in which the element value of \sphinxcode{\sphinxupquote{idx}} is set to \sphinxcode{\sphinxupquote{value}}.
\end{quote}

\sphinxAtStartPar
\sphinxstylestrong{Example}

\begin{sphinxVerbatim}[commandchars=\\\{\}]
\PYG{c+c1}{\PYGZsh{} Set the value of the element at index 0 to 100}
\PYG{n}{mat\PYGZus{}1} \PYG{o}{=} \PYG{n}{ndmat}\PYG{o}{.}\PYG{n}{setItem}\PYG{p}{(}\PYG{l+m+mi}{0}\PYG{p}{,} \PYG{l+m+mi}{100}\PYG{p}{)}
\end{sphinxVerbatim}
\end{quote}

\subsubsection{NdArray.transpose()}
\label{\detokenize{pyapiref:ndarray-transpose}}\begin{quote}

\sphinxAtStartPar
\sphinxstylestrong{Synopsis}
\begin{quote}

\sphinxAtStartPar
\sphinxcode{\sphinxupquote{transpose()}}
\end{quote}

\sphinxAtStartPar
\sphinxstylestrong{Description}
\begin{quote}

\sphinxAtStartPar
Generates a new NdArray object, which is the transpose of the original one.
\end{quote}

\sphinxAtStartPar
\sphinxstylestrong{Return value}
\begin{quote}

\sphinxAtStartPar
Returns a new NdArray object.
\end{quote}

\sphinxAtStartPar
\sphinxstylestrong{Example}

\begin{sphinxVerbatim}[commandchars=\\\{\}]
\PYG{n}{ndmat} \PYG{o}{=} \PYG{n}{NdArray}\PYG{p}{(}\PYG{n}{shape}\PYG{o}{=}\PYG{p}{(}\PYG{l+m+mi}{3}\PYG{p}{,} \PYG{l+m+mi}{5}\PYG{p}{)}\PYG{p}{)}
\PYG{n+nb}{print}\PYG{p}{(}\PYG{n}{ndmat}\PYG{o}{.}\PYG{n}{transpose}\PYG{p}{(}\PYG{p}{)}\PYG{o}{.}\PYG{n}{shape}\PYG{p}{)} \PYG{c+c1}{\PYGZsh{}shape=(5, 3)}
\end{sphinxVerbatim}
\end{quote}

\subsubsection{NdArray.diagonal()}
\label{\detokenize{pyapiref:ndarray-diagonal}}\begin{quote}

\sphinxAtStartPar
\sphinxstylestrong{Synopsis}
\begin{quote}

\sphinxAtStartPar
\sphinxcode{\sphinxupquote{diagonal(offset=0, axis1=0, axis2=1)}}
\end{quote}

\sphinxAtStartPar
\sphinxstylestrong{Description}
\begin{quote}
\begin{quote}

\sphinxAtStartPar
Generate a {\hyperref[\detokenize{pyapiref:chappyapi-ndarray}]{\sphinxcrossref{\DUrole{std,std-ref}{NdArray Class}}}} object, in which elements are the ones on the diagonal of the original NdArray object.
\end{quote}

\sphinxAtStartPar
\sphinxstylestrong{Arguments}
\begin{quote}

\sphinxAtStartPar
\sphinxcode{\sphinxupquote{offset}}
\begin{quote}

\sphinxAtStartPar
Optional argument, indicating the offset of the diagonal, the default value is 0.
If the value is greater than 0, it represents the upward offset of the diagonal. If the value is less than 0, it represents the downward offset of the diagonal.
\end{quote}

\sphinxAtStartPar
\sphinxcode{\sphinxupquote{axis1}}
\begin{quote}

\sphinxAtStartPar
Optional argument.
Axis to be used as the first axis of the 2\sphinxhyphen{}D sub\sphinxhyphen{}NdArrays from which the diagonals should be taken. Defaults to first axis (0).
\end{quote}

\sphinxAtStartPar
\sphinxcode{\sphinxupquote{axis2}}
\begin{quote}

\sphinxAtStartPar
Optional argument.
Axis to be used as the second axis of the 2\sphinxhyphen{}D sub\sphinxhyphen{}NdArrays from which the diagonals should be taken. Defaults to second axis (1).
\end{quote}
\end{quote}
\end{quote}

\sphinxAtStartPar
\sphinxstylestrong{Return value}
\begin{quote}

\sphinxAtStartPar
A new NdArray object.
\end{quote}

\sphinxAtStartPar
\sphinxstylestrong{Example}

\begin{sphinxVerbatim}[commandchars=\\\{\}]
\PYG{n}{ndmat} \PYG{o}{=} \PYG{n}{NdArray}\PYG{p}{(}\PYG{n}{shape}\PYG{o}{=}\PYG{p}{(}\PYG{l+m+mi}{3}\PYG{p}{,} \PYG{l+m+mi}{3}\PYG{p}{)}\PYG{p}{,}\PYG{n}{args}\PYG{o}{=}\PYG{p}{[}\PYG{p}{[}\PYG{l+m+mi}{1}\PYG{p}{,}\PYG{l+m+mi}{1}\PYG{p}{,}\PYG{l+m+mi}{1}\PYG{p}{]}\PYG{p}{,}\PYG{p}{[}\PYG{l+m+mi}{2}\PYG{p}{,}\PYG{l+m+mi}{2}\PYG{p}{,}\PYG{l+m+mi}{2}\PYG{p}{]}\PYG{p}{,}\PYG{p}{[}\PYG{l+m+mi}{3}\PYG{p}{,}\PYG{l+m+mi}{3}\PYG{p}{,}\PYG{l+m+mi}{3}\PYG{p}{]}\PYG{p}{]}\PYG{p}{)}

\PYG{n}{diag\PYGZus{}m0} \PYG{o}{=} \PYG{n}{ndmat}\PYG{o}{.}\PYG{n}{diagonal}\PYG{p}{(}\PYG{l+m+mi}{0}\PYG{p}{)}
\PYG{n}{diag\PYGZus{}a1} \PYG{o}{=} \PYG{n}{ndmat}\PYG{o}{.}\PYG{n}{diagonal}\PYG{p}{(}\PYG{l+m+mi}{1}\PYG{p}{)}
\PYG{n}{diag\PYGZus{}b1} \PYG{o}{=} \PYG{n}{ndmat}\PYG{o}{.}\PYG{n}{diagonal}\PYG{p}{(}\PYG{o}{\PYGZhy{}}\PYG{l+m+mi}{1}\PYG{p}{)}
\end{sphinxVerbatim}
\end{quote}

\subsubsection{NdArray.ndim}
\label{\detokenize{pyapiref:ndarray-ndim}}\begin{quote}

\sphinxAtStartPar
\sphinxstylestrong{Synopsis}
\begin{quote}

\sphinxAtStartPar
\sphinxcode{\sphinxupquote{ndim}}
\end{quote}

\sphinxAtStartPar
\sphinxstylestrong{Description}
\begin{quote}

\sphinxAtStartPar
The dimensions of the NdArray object.
\end{quote}

\sphinxAtStartPar
\sphinxstylestrong{Return value}
\begin{quote}

\sphinxAtStartPar
An integer Value
\end{quote}

\sphinxAtStartPar
\sphinxstylestrong{Example}

\begin{sphinxVerbatim}[commandchars=\\\{\}]
\PYG{n}{ndmat} \PYG{o}{=} \PYG{n}{NdArray}\PYG{p}{(}\PYG{p}{(}\PYG{l+m+mi}{3}\PYG{p}{,} \PYG{l+m+mi}{5}\PYG{p}{)}\PYG{p}{)}
\PYG{n+nb}{print}\PYG{p}{(}\PYG{n}{ndmat}\PYG{o}{.}\PYG{n}{ndim}\PYG{p}{)} \PYG{c+c1}{\PYGZsh{} ndim = 2}
\end{sphinxVerbatim}
\end{quote}

\subsubsection{NdArray.shape}
\label{\detokenize{pyapiref:ndarray-shape}}\begin{quote}

\sphinxAtStartPar
\sphinxstylestrong{Synopsis}
\begin{quote}

\sphinxAtStartPar
\sphinxcode{\sphinxupquote{shape}}
\end{quote}

\sphinxAtStartPar
\sphinxstylestrong{Description}
\begin{quote}

\sphinxAtStartPar
The shape of the NdArray object.
\end{quote}

\sphinxAtStartPar
\sphinxstylestrong{Return value}
\begin{quote}

\sphinxAtStartPar
An integer tuple.
\end{quote}

\sphinxAtStartPar
\sphinxstylestrong{Example}

\begin{sphinxVerbatim}[commandchars=\\\{\}]
\PYG{n}{ndmat} \PYG{o}{=} \PYG{n}{NdArray}\PYG{p}{(}\PYG{p}{(}\PYG{l+m+mi}{3}\PYG{p}{,} \PYG{l+m+mi}{5}\PYG{p}{)}\PYG{p}{)}
\PYG{n+nb}{print}\PYG{p}{(}\PYG{n}{ndmat}\PYG{o}{.}\PYG{n}{shape}\PYG{p}{)} \PYG{c+c1}{\PYGZsh{} shape = (3, 5)}
\end{sphinxVerbatim}
\end{quote}

\subsubsection{NdArray.pick()}
\label{\detokenize{pyapiref:ndarray-pick}}\begin{quote}

\sphinxAtStartPar
\sphinxstylestrong{Synopsis}
\begin{quote}

\sphinxAtStartPar
\sphinxcode{\sphinxupquote{pick(indexes)}}
\end{quote}

\sphinxAtStartPar
\sphinxstylestrong{Description}
\begin{quote}

\sphinxAtStartPar
Retrieves a new multi\sphinxhyphen{}dimensional array from the \sphinxtitleref{NdArray} object
according to the specified indices.
\end{quote}

\sphinxAtStartPar
\sphinxstylestrong{Arguments}
\begin{quote}

\sphinxAtStartPar
\sphinxcode{\sphinxupquote{indexes}}
\begin{quote}

\sphinxAtStartPar
The specified indices.
\end{quote}
\end{quote}

\sphinxAtStartPar
\sphinxstylestrong{Return Value}
\begin{quote}

\sphinxAtStartPar
Returns a new \sphinxtitleref{NdArray} object.
\end{quote}
\end{quote}

\subsubsection{NdArray.hstack()}
\label{\detokenize{pyapiref:ndarray-hstack}}\begin{quote}

\sphinxAtStartPar
\sphinxstylestrong{Synopsis}
\begin{quote}

\sphinxAtStartPar
\sphinxcode{\sphinxupquote{hstack(other)}}
\end{quote}

\sphinxAtStartPar
\sphinxstylestrong{Description}
\begin{quote}

\sphinxAtStartPar
Stacks another \sphinxtitleref{NdArray} object along the horizontal dimension (last dimension)
to form a new \sphinxtitleref{NdArray} object.
\end{quote}

\sphinxAtStartPar
\sphinxstylestrong{Arguments}
\begin{quote}

\sphinxAtStartPar
\sphinxcode{\sphinxupquote{other}}
\begin{quote}

\sphinxAtStartPar
Another \sphinxtitleref{NdArray} object.
\end{quote}
\end{quote}

\sphinxAtStartPar
\sphinxstylestrong{Return Value}
\begin{quote}

\sphinxAtStartPar
Returns a new \sphinxtitleref{NdArray} object.
\end{quote}
\end{quote}

\subsubsection{NdArray.vstack()}
\label{\detokenize{pyapiref:ndarray-vstack}}\begin{quote}

\sphinxAtStartPar
\sphinxstylestrong{Synopsis}
\begin{quote}

\sphinxAtStartPar
\sphinxcode{\sphinxupquote{vstack(other)}}
\end{quote}

\sphinxAtStartPar
\sphinxstylestrong{Description}
\begin{quote}

\sphinxAtStartPar
Stacks another \sphinxtitleref{NdArray} object along the vertical dimension to form a new \sphinxtitleref{NdArray} object.
\end{quote}

\sphinxAtStartPar
\sphinxstylestrong{Arguments}
\begin{quote}

\sphinxAtStartPar
\sphinxcode{\sphinxupquote{other}}
\begin{quote}

\sphinxAtStartPar
Another \sphinxtitleref{NdArray} object.
\end{quote}
\end{quote}

\sphinxAtStartPar
\sphinxstylestrong{Return Value}
\begin{quote}

\sphinxAtStartPar
Returns a new \sphinxtitleref{NdArray} object.
\end{quote}
\end{quote}

\subsubsection{NdArray.stack()}
\label{\detokenize{pyapiref:ndarray-stack}}\begin{quote}

\sphinxAtStartPar
\sphinxstylestrong{Synopsis}
\begin{quote}

\sphinxAtStartPar
\sphinxcode{\sphinxupquote{stack(other, axis)}}
\end{quote}

\sphinxAtStartPar
\sphinxstylestrong{Description}
\begin{quote}

\sphinxAtStartPar
Stacks another \sphinxtitleref{NdArray} object along the specified axis to form a new \sphinxtitleref{NdArray} object.
\end{quote}

\sphinxAtStartPar
\sphinxstylestrong{Arguments}
\begin{quote}

\sphinxAtStartPar
\sphinxcode{\sphinxupquote{other}}
\begin{quote}

\sphinxAtStartPar
Another \sphinxtitleref{NdArray} object.
\end{quote}

\sphinxAtStartPar
\sphinxcode{\sphinxupquote{axis}}
\begin{quote}

\sphinxAtStartPar
The specified axis index.
\end{quote}
\end{quote}

\sphinxAtStartPar
\sphinxstylestrong{Return Value}
\begin{quote}

\sphinxAtStartPar
Returns a new \sphinxtitleref{NdArray} object.
\end{quote}
\end{quote}

\subsubsection{NdArray.size}
\label{\detokenize{pyapiref:ndarray-size}}\begin{quote}

\sphinxAtStartPar
\sphinxstylestrong{Synopsis}
\begin{quote}

\sphinxAtStartPar
\sphinxcode{\sphinxupquote{size}}
\end{quote}

\sphinxAtStartPar
\sphinxstylestrong{Description}
\begin{quote}

\sphinxAtStartPar
The number of elements in the NdArray object.
\end{quote}

\sphinxAtStartPar
\sphinxstylestrong{Return value}
\begin{quote}

\sphinxAtStartPar
An integer value.
\end{quote}

\sphinxAtStartPar
\sphinxstylestrong{Example}

\begin{sphinxVerbatim}[commandchars=\\\{\}]
\PYG{n}{ndmat} \PYG{o}{=} \PYG{n}{NdArray}\PYG{p}{(}\PYG{p}{(}\PYG{l+m+mi}{3}\PYG{p}{,} \PYG{l+m+mi}{5}\PYG{p}{)}\PYG{p}{)}
\PYG{n+nb}{print}\PYG{p}{(}\PYG{n}{ndmat}\PYG{o}{.}\PYG{n}{size}\PYG{p}{)} \PYG{c+c1}{\PYGZsh{} size = 15}
\end{sphinxVerbatim}
\end{quote}

\subsubsection{NdArray.T}
\label{\detokenize{pyapiref:ndarray-t}}\begin{quote}

\sphinxAtStartPar
\sphinxstylestrong{Synopsis}
\begin{quote}

\sphinxAtStartPar
\sphinxcode{\sphinxupquote{T}}
\end{quote}

\sphinxAtStartPar
\sphinxstylestrong{Description}
\begin{quote}

\sphinxAtStartPar
Transpose of NdArray objects, similar to the method \sphinxcode{\sphinxupquote{NdArray.transpose()}} .
\end{quote}

\sphinxAtStartPar
\sphinxstylestrong{Return value}
\begin{quote}

\sphinxAtStartPar
Returns the transposed NdArray object.
\end{quote}

\sphinxAtStartPar
\sphinxstylestrong{Example}

\begin{sphinxVerbatim}[commandchars=\\\{\}]
\PYG{n}{ndmat} \PYG{o}{=} \PYG{n}{NdArray}\PYG{p}{(}\PYG{n}{shape}\PYG{o}{=}\PYG{p}{(}\PYG{l+m+mi}{3}\PYG{p}{,} \PYG{l+m+mi}{5}\PYG{p}{)}\PYG{p}{)}
\PYG{n+nb}{print}\PYG{p}{(}\PYG{n}{ndmat}\PYG{o}{.}\PYG{n}{T}\PYG{o}{.}\PYG{n}{shape}\PYG{p}{)} \PYG{c+c1}{\PYGZsh{} shape = (5, 3)}
\end{sphinxVerbatim}
\end{quote}

\subsection{ExprBuilder Class}
\label{\detokenize{pyapiref:exprbuilder-class}}\label{\detokenize{pyapiref:chappyapi-exprbuilder}}
\sphinxAtStartPar
ExprBuilder object contains operations related to building linear expressions, and provides the following methods:

\subsubsection{ExprBuilder()}
\label{\detokenize{pyapiref:exprbuilder}}\begin{quote}

\sphinxAtStartPar
\sphinxstylestrong{Synopsis}
\begin{quote}

\sphinxAtStartPar
\sphinxcode{\sphinxupquote{ExprBuilder(arg1=0.0, arg2=None)}}
\end{quote}

\sphinxAtStartPar
\sphinxstylestrong{Description}
\begin{quote}

\sphinxAtStartPar
Create a {\hyperref[\detokenize{pyapiref:chappyapi-exprbuilder}]{\sphinxcrossref{\DUrole{std,std-ref}{ExprBuilder Class}}}} object.

\sphinxAtStartPar
If argument \sphinxcode{\sphinxupquote{arg1}} is constant, argument \sphinxcode{\sphinxupquote{arg2}} is \sphinxcode{\sphinxupquote{None}}, then create a {\hyperref[\detokenize{pyapiref:chappyapi-exprbuilder}]{\sphinxcrossref{\DUrole{std,std-ref}{ExprBuilder Class}}}} object and initialize it using argument \sphinxcode{\sphinxupquote{arg1}}.
If argument \sphinxcode{\sphinxupquote{arg1}} is {\hyperref[\detokenize{pyapiref:chappyapi-var}]{\sphinxcrossref{\DUrole{std,std-ref}{Var Class}}}} or {\hyperref[\detokenize{pyapiref:chappyapi-exprbuilder}]{\sphinxcrossref{\DUrole{std,std-ref}{ExprBuilder Class}}}} object, and argument \sphinxcode{\sphinxupquote{arg2}} is constant or considered to be constant 1.0 when
argument \sphinxcode{\sphinxupquote{arg2}} is \sphinxcode{\sphinxupquote{None}}, then initialize the newly created {\hyperref[\detokenize{pyapiref:chappyapi-exprbuilder}]{\sphinxcrossref{\DUrole{std,std-ref}{ExprBuilder Class}}}} object using arguments \sphinxcode{\sphinxupquote{arg1}} and \sphinxcode{\sphinxupquote{arg2}}.
If argument \sphinxcode{\sphinxupquote{arg1}} and \sphinxcode{\sphinxupquote{arg2}} are list objects, then they are variables and coefficients used to initialize the newly created {\hyperref[\detokenize{pyapiref:chappyapi-exprbuilder}]{\sphinxcrossref{\DUrole{std,std-ref}{ExprBuilder Class}}}} object.
\end{quote}

\sphinxAtStartPar
\sphinxstylestrong{Arguments}
\begin{quote}

\sphinxAtStartPar
\sphinxcode{\sphinxupquote{arg1}}
\begin{quote}

\sphinxAtStartPar
Optional, 0.0 by default.
\end{quote}

\sphinxAtStartPar
\sphinxcode{\sphinxupquote{arg2}}
\begin{quote}

\sphinxAtStartPar
Optional, \sphinxcode{\sphinxupquote{None}} by default.
\end{quote}
\end{quote}

\sphinxAtStartPar
\sphinxstylestrong{Example}
\end{quote}

\begin{sphinxVerbatim}[commandchars=\\\{\}]
\PYG{c+c1}{\PYGZsh{} Create a new ExprBuilder object and initialize it to 0.0}
\PYG{n}{expr0} \PYG{o}{=} \PYG{n}{ExprBuilder}\PYG{p}{(}\PYG{p}{)}
\PYG{c+c1}{\PYGZsh{} Create a ExprBuilder object and initialize it to x + 2*y}
\PYG{n}{expr2} \PYG{o}{=} \PYG{n}{ExprBuilder}\PYG{p}{(}\PYG{p}{[}\PYG{n}{x}\PYG{p}{,} \PYG{n}{y}\PYG{p}{]}\PYG{p}{,} \PYG{p}{[}\PYG{l+m+mi}{1}\PYG{p}{,} \PYG{l+m+mi}{2}\PYG{p}{]}\PYG{p}{)}
\end{sphinxVerbatim}

\subsubsection{ExprBuilder.getSize()}
\label{\detokenize{pyapiref:exprbuilder-getsize}}\begin{quote}

\sphinxAtStartPar
\sphinxstylestrong{Synopsis}
\begin{quote}

\sphinxAtStartPar
\sphinxcode{\sphinxupquote{getSize()}}
\end{quote}

\sphinxAtStartPar
\sphinxstylestrong{Description}
\begin{quote}

\sphinxAtStartPar
Retrieve the number of terms in an expression builder.
\end{quote}

\sphinxAtStartPar
\sphinxstylestrong{Example}
\end{quote}

\begin{sphinxVerbatim}[commandchars=\\\{\}]
\PYG{c+c1}{\PYGZsh{} Retrieve the number of terms in expression builder \PYGZsq{}expr\PYGZsq{}}
\PYG{n}{exprsize} \PYG{o}{=} \PYG{n}{expr}\PYG{o}{.}\PYG{n}{getSize}\PYG{p}{(}\PYG{p}{)}
\end{sphinxVerbatim}

\subsubsection{ExprBuilder.getCoeff()}
\label{\detokenize{pyapiref:exprbuilder-getcoeff}}\begin{quote}

\sphinxAtStartPar
\sphinxstylestrong{Synopsis}
\begin{quote}

\sphinxAtStartPar
\sphinxcode{\sphinxupquote{getCoeff(idx)}}
\end{quote}

\sphinxAtStartPar
\sphinxstylestrong{Description}
\begin{quote}

\sphinxAtStartPar
Retrieve the coefficient of a variable by its index from an expression builder.
\end{quote}

\sphinxAtStartPar
\sphinxstylestrong{Arguments}
\begin{quote}

\sphinxAtStartPar
\sphinxcode{\sphinxupquote{idx}}
\begin{quote}

\sphinxAtStartPar
Index of the variable in the expression builder, starting with 0.
\end{quote}
\end{quote}

\sphinxAtStartPar
\sphinxstylestrong{Example}
\end{quote}

\begin{sphinxVerbatim}[commandchars=\\\{\}]
\PYG{c+c1}{\PYGZsh{} Retrieve the coefficient for the term at index 1 from expression builder \PYGZsq{}expr\PYGZsq{}}
\PYG{n}{coeff} \PYG{o}{=} \PYG{n}{expr}\PYG{o}{.}\PYG{n}{getCoeff}\PYG{p}{(}\PYG{l+m+mi}{1}\PYG{p}{)}
\end{sphinxVerbatim}

\subsubsection{ExprBuilder.getVar()}
\label{\detokenize{pyapiref:exprbuilder-getvar}}\begin{quote}

\sphinxAtStartPar
\sphinxstylestrong{Synopsis}
\begin{quote}

\sphinxAtStartPar
\sphinxcode{\sphinxupquote{getVar(idx)}}
\end{quote}

\sphinxAtStartPar
\sphinxstylestrong{Description}
\begin{quote}

\sphinxAtStartPar
Retrieve the variable by its index from an expression builder. Return a {\hyperref[\detokenize{pyapiref:chappyapi-var}]{\sphinxcrossref{\DUrole{std,std-ref}{Var Class}}}} object.
\end{quote}

\sphinxAtStartPar
\sphinxstylestrong{Arguments}
\begin{quote}

\sphinxAtStartPar
\sphinxcode{\sphinxupquote{idx}}
\begin{quote}

\sphinxAtStartPar
Index of the variable in the expression builder, starting with 0.
\end{quote}
\end{quote}

\sphinxAtStartPar
\sphinxstylestrong{Example}
\end{quote}

\begin{sphinxVerbatim}[commandchars=\\\{\}]
\PYG{c+c1}{\PYGZsh{} Retrieve the variable for the term at index 1 from expression builder \PYGZsq{}expr\PYGZsq{}}
\PYG{n}{x} \PYG{o}{=} \PYG{n}{expr}\PYG{o}{.}\PYG{n}{getVar}\PYG{p}{(}\PYG{l+m+mi}{1}\PYG{p}{)}
\end{sphinxVerbatim}

\subsubsection{ExprBuilder.getConstant()}
\label{\detokenize{pyapiref:exprbuilder-getconstant}}\begin{quote}

\sphinxAtStartPar
\sphinxstylestrong{Synopsis}
\begin{quote}

\sphinxAtStartPar
\sphinxcode{\sphinxupquote{getConstant()}}
\end{quote}

\sphinxAtStartPar
\sphinxstylestrong{Description}
\begin{quote}

\sphinxAtStartPar
Retrieve the constant term from an expression builder.
\end{quote}

\sphinxAtStartPar
\sphinxstylestrong{Example}
\end{quote}

\begin{sphinxVerbatim}[commandchars=\\\{\}]
\PYG{c+c1}{\PYGZsh{} Retrieve the constant term from linear expression builder \PYGZsq{}expr\PYGZsq{}}
\PYG{n}{constant} \PYG{o}{=} \PYG{n}{expr}\PYG{o}{.}\PYG{n}{getConstant}\PYG{p}{(}\PYG{p}{)}
\end{sphinxVerbatim}

\subsubsection{ExprBuilder.addTerm()}
\label{\detokenize{pyapiref:exprbuilder-addterm}}\begin{quote}

\sphinxAtStartPar
\sphinxstylestrong{Synopsis}
\begin{quote}

\sphinxAtStartPar
\sphinxcode{\sphinxupquote{addTerm(var, coeff=1.0)}}
\end{quote}

\sphinxAtStartPar
\sphinxstylestrong{Description}
\begin{quote}

\sphinxAtStartPar
Add a new term to current expression builder.
\end{quote}

\sphinxAtStartPar
\sphinxstylestrong{Arguments}
\begin{quote}

\sphinxAtStartPar
\sphinxcode{\sphinxupquote{var}}
\begin{quote}

\sphinxAtStartPar
Variable to add.
\end{quote}

\sphinxAtStartPar
\sphinxcode{\sphinxupquote{coeff}}
\begin{quote}

\sphinxAtStartPar
Magnification coefficient for added term. Optional, 1.0 by default.
\end{quote}
\end{quote}

\sphinxAtStartPar
\sphinxstylestrong{Example}
\end{quote}

\begin{sphinxVerbatim}[commandchars=\\\{\}]
\PYG{c+c1}{\PYGZsh{} Add term 2*x to linear expression builder \PYGZsq{}expr\PYGZsq{}}
\PYG{n}{expr}\PYG{o}{.}\PYG{n}{addTerm}\PYG{p}{(}\PYG{n}{x}\PYG{p}{,} \PYG{l+m+mf}{2.0}\PYG{p}{)}
\end{sphinxVerbatim}

\subsubsection{ExprBuilder.addExpr()}
\label{\detokenize{pyapiref:exprbuilder-addexpr}}\begin{quote}

\sphinxAtStartPar
\sphinxstylestrong{Synopsis}
\begin{quote}

\sphinxAtStartPar
\sphinxcode{\sphinxupquote{addExpr(expr, coeff=1.0)}}
\end{quote}

\sphinxAtStartPar
\sphinxstylestrong{Description}
\begin{quote}

\sphinxAtStartPar
Add new expression builder to the current one.
\end{quote}

\sphinxAtStartPar
\sphinxstylestrong{Arguments}
\begin{quote}

\sphinxAtStartPar
\sphinxcode{\sphinxupquote{expr}}
\begin{quote}

\sphinxAtStartPar
Expression builder to add.
\end{quote}

\sphinxAtStartPar
\sphinxcode{\sphinxupquote{coeff}}
\begin{quote}

\sphinxAtStartPar
Magnification coefficients for the added expression builder. Optional, 1.0 by default.
\end{quote}
\end{quote}

\sphinxAtStartPar
\sphinxstylestrong{Example}
\end{quote}

\begin{sphinxVerbatim}[commandchars=\\\{\}]
\PYG{c+c1}{\PYGZsh{} Add linear expression builder 2*x + 2*y to \PYGZsq{}expr\PYGZsq{}}
\PYG{n}{expr}\PYG{o}{.}\PYG{n}{addExpr}\PYG{p}{(}\PYG{n}{x} \PYG{o}{+} \PYG{n}{y}\PYG{p}{,} \PYG{l+m+mf}{2.0}\PYG{p}{)}
\end{sphinxVerbatim}

\subsubsection{ExprBuilder.clone()}
\label{\detokenize{pyapiref:exprbuilder-clone}}\begin{quote}

\sphinxAtStartPar
\sphinxstylestrong{Synopsis}
\begin{quote}

\sphinxAtStartPar
\sphinxcode{\sphinxupquote{clone()}}
\end{quote}

\sphinxAtStartPar
\sphinxstylestrong{Description}
\begin{quote}

\sphinxAtStartPar
Create a deep copy of the expression builder.
\end{quote}

\sphinxAtStartPar
\sphinxstylestrong{Example}
\end{quote}

\begin{sphinxVerbatim}[commandchars=\\\{\}]
\PYG{c+c1}{\PYGZsh{} Create a deep copy of expression builder \PYGZsq{}expr\PYGZsq{}}
\PYG{n}{exprcopy} \PYG{o}{=} \PYG{n}{expr}\PYG{o}{.}\PYG{n}{clone}\PYG{p}{(}\PYG{p}{)}
\end{sphinxVerbatim}

\subsubsection{ExprBuilder.getExpr()}
\label{\detokenize{pyapiref:exprbuilder-getexpr}}\begin{quote}

\sphinxAtStartPar
\sphinxstylestrong{Synopsis}
\begin{quote}

\sphinxAtStartPar
\sphinxcode{\sphinxupquote{getExpr()}}
\end{quote}

\sphinxAtStartPar
\sphinxstylestrong{Description}
\begin{quote}

\sphinxAtStartPar
Create a linear expression related to the expression builder. Returns a {\hyperref[\detokenize{pyapiref:chappyapi-linexpr}]{\sphinxcrossref{\DUrole{std,std-ref}{LinExpr Class}}}} object.
\end{quote}

\sphinxAtStartPar
\sphinxstylestrong{Example}
\end{quote}

\begin{sphinxVerbatim}[commandchars=\\\{\}]
\PYG{c+c1}{\PYGZsh{} Get the linear expression object related to expression builder \PYGZsq{}exprbuilder\PYGZsq{}}
\PYG{n}{expr} \PYG{o}{=} \PYG{n}{exprbuilder}\PYG{o}{.}\PYG{n}{getExpr}\PYG{p}{(}\PYG{p}{)}
\end{sphinxVerbatim}

\subsection{LinExpr Class}
\label{\detokenize{pyapiref:linexpr-class}}\label{\detokenize{pyapiref:chappyapi-linexpr}}
\sphinxAtStartPar
LinExpr object contains operations related to variables for building linear constraints, and provides the following methods:

\subsubsection{LinExpr()}
\label{\detokenize{pyapiref:linexpr}}\begin{quote}

\sphinxAtStartPar
\sphinxstylestrong{Synopsis}
\begin{quote}

\sphinxAtStartPar
\sphinxcode{\sphinxupquote{LinExpr(arg1=0.0, arg2=None)}}
\end{quote}

\sphinxAtStartPar
\sphinxstylestrong{Description}
\begin{quote}

\sphinxAtStartPar
Create a {\hyperref[\detokenize{pyapiref:chappyapi-linexpr}]{\sphinxcrossref{\DUrole{std,std-ref}{LinExpr Class}}}} object.

\sphinxAtStartPar
If argument \sphinxcode{\sphinxupquote{arg1}} is constant, argument \sphinxcode{\sphinxupquote{arg2}} is \sphinxcode{\sphinxupquote{None}}, then create a {\hyperref[\detokenize{pyapiref:chappyapi-linexpr}]{\sphinxcrossref{\DUrole{std,std-ref}{LinExpr Class}}}} object and initialize it using argument \sphinxcode{\sphinxupquote{arg1}}.
If argument \sphinxcode{\sphinxupquote{arg1}} is {\hyperref[\detokenize{pyapiref:chappyapi-var}]{\sphinxcrossref{\DUrole{std,std-ref}{Var Class}}}} or {\hyperref[\detokenize{pyapiref:chappyapi-linexpr}]{\sphinxcrossref{\DUrole{std,std-ref}{LinExpr Class}}}} object, and argument \sphinxcode{\sphinxupquote{arg2}} is constant or considered to be constant 1.0 when
argument \sphinxcode{\sphinxupquote{arg2}} is \sphinxcode{\sphinxupquote{None}}, then initialize the newly created {\hyperref[\detokenize{pyapiref:chappyapi-linexpr}]{\sphinxcrossref{\DUrole{std,std-ref}{LinExpr Class}}}} object using arguments \sphinxcode{\sphinxupquote{arg1}} and \sphinxcode{\sphinxupquote{arg2}}.
If argument \sphinxcode{\sphinxupquote{arg1}} is list object and argument \sphinxcode{\sphinxupquote{arg2}} is \sphinxcode{\sphinxupquote{None}}, then argument \sphinxcode{\sphinxupquote{arg1}} contains a list of variable\sphinxhyphen{}coefficient pairs and initialize the newly created {\hyperref[\detokenize{pyapiref:chappyapi-linexpr}]{\sphinxcrossref{\DUrole{std,std-ref}{LinExpr Class}}}} object using arguments \sphinxcode{\sphinxupquote{arg1}} and \sphinxcode{\sphinxupquote{arg2}}.
For other forms of arguments, call method \sphinxcode{\sphinxupquote{addTerms}} to initialize the newly created {\hyperref[\detokenize{pyapiref:chappyapi-linexpr}]{\sphinxcrossref{\DUrole{std,std-ref}{LinExpr Class}}}} object.
\end{quote}

\sphinxAtStartPar
\sphinxstylestrong{Arguments}
\begin{quote}

\sphinxAtStartPar
\sphinxcode{\sphinxupquote{arg1}}
\begin{quote}

\sphinxAtStartPar
Optional, 0.0 by default.
\end{quote}

\sphinxAtStartPar
\sphinxcode{\sphinxupquote{arg2}}
\begin{quote}

\sphinxAtStartPar
Optional, \sphinxcode{\sphinxupquote{None}} by default.
\end{quote}
\end{quote}

\sphinxAtStartPar
\sphinxstylestrong{Example}
\end{quote}

\begin{sphinxVerbatim}[commandchars=\\\{\}]
\PYG{c+c1}{\PYGZsh{} Create a new LinExpr object and initialize it to 0.0}
\PYG{n}{expr0} \PYG{o}{=} \PYG{n}{LinExpr}\PYG{p}{(}\PYG{p}{)}
\PYG{c+c1}{\PYGZsh{} Create a LinExpr object and initialize it to 2*x + 3*y}
\PYG{n}{expr1} \PYG{o}{=} \PYG{n}{LinExpr}\PYG{p}{(}\PYG{p}{[}\PYG{p}{(}\PYG{n}{x}\PYG{p}{,} \PYG{l+m+mi}{2}\PYG{p}{)}\PYG{p}{,} \PYG{p}{(}\PYG{n}{y}\PYG{p}{,} \PYG{l+m+mi}{3}\PYG{p}{)}\PYG{p}{]}\PYG{p}{)}
\PYG{c+c1}{\PYGZsh{} Create a LinExpr object and initialize it to x + 2*y}
\PYG{n}{expr2} \PYG{o}{=} \PYG{n}{LinExpr}\PYG{p}{(}\PYG{p}{[}\PYG{n}{x}\PYG{p}{,} \PYG{n}{y}\PYG{p}{]}\PYG{p}{,} \PYG{p}{[}\PYG{l+m+mi}{1}\PYG{p}{,} \PYG{l+m+mi}{2}\PYG{p}{]}\PYG{p}{)}
\end{sphinxVerbatim}

\subsubsection{LinExpr.setCoeff()}
\label{\detokenize{pyapiref:linexpr-setcoeff}}\begin{quote}

\sphinxAtStartPar
\sphinxstylestrong{Synopsis}
\begin{quote}

\sphinxAtStartPar
\sphinxcode{\sphinxupquote{setCoeff(idx, newval)}}
\end{quote}

\sphinxAtStartPar
\sphinxstylestrong{Description}
\begin{quote}

\sphinxAtStartPar
Set new coefficient value of a variable based on its index in an expression.
\end{quote}

\sphinxAtStartPar
\sphinxstylestrong{Arguments}
\begin{quote}

\sphinxAtStartPar
\sphinxcode{\sphinxupquote{idx}}
\begin{quote}

\sphinxAtStartPar
Index of the variable in the expression, starting with 0.
\end{quote}

\sphinxAtStartPar
\sphinxcode{\sphinxupquote{newval}}
\begin{quote}

\sphinxAtStartPar
New coefficient value of the variable.
\end{quote}
\end{quote}

\sphinxAtStartPar
\sphinxstylestrong{Example}
\end{quote}

\begin{sphinxVerbatim}[commandchars=\\\{\}]
\PYG{c+c1}{\PYGZsh{} Set the coefficient for the term at index 0 in expression expr to 1.0}
\PYG{n}{expr}\PYG{o}{.}\PYG{n}{setCoeff}\PYG{p}{(}\PYG{l+m+mi}{0}\PYG{p}{,} \PYG{l+m+mf}{1.0}\PYG{p}{)}
\end{sphinxVerbatim}

\subsubsection{LinExpr.getCoeff()}
\label{\detokenize{pyapiref:linexpr-getcoeff}}\begin{quote}

\sphinxAtStartPar
\sphinxstylestrong{Synopsis}
\begin{quote}

\sphinxAtStartPar
\sphinxcode{\sphinxupquote{getCoeff(idx)}}
\end{quote}

\sphinxAtStartPar
\sphinxstylestrong{Description}
\begin{quote}

\sphinxAtStartPar
Retrieve the coefficient of a variable by its index from an expression.
\end{quote}

\sphinxAtStartPar
\sphinxstylestrong{Arguments}
\begin{quote}

\sphinxAtStartPar
\sphinxcode{\sphinxupquote{idx}}
\begin{quote}

\sphinxAtStartPar
Index of the variable in the expression, starting with 0.
\end{quote}
\end{quote}

\sphinxAtStartPar
\sphinxstylestrong{Example}
\end{quote}

\begin{sphinxVerbatim}[commandchars=\\\{\}]
\PYG{c+c1}{\PYGZsh{} Retrieve the coefficient for the term at index 1 from expression expr}
\PYG{n}{coeff} \PYG{o}{=} \PYG{n}{expr}\PYG{o}{.}\PYG{n}{getCoeff}\PYG{p}{(}\PYG{l+m+mi}{1}\PYG{p}{)}
\end{sphinxVerbatim}

\subsubsection{LinExpr.getVar()}
\label{\detokenize{pyapiref:linexpr-getvar}}\begin{quote}

\sphinxAtStartPar
\sphinxstylestrong{Synopsis}
\begin{quote}

\sphinxAtStartPar
\sphinxcode{\sphinxupquote{getVar(idx)}}
\end{quote}

\sphinxAtStartPar
\sphinxstylestrong{Description}
\begin{quote}

\sphinxAtStartPar
Retrieve the variable by its index from an expression. Return a {\hyperref[\detokenize{pyapiref:chappyapi-var}]{\sphinxcrossref{\DUrole{std,std-ref}{Var Class}}}} object.
\end{quote}

\sphinxAtStartPar
\sphinxstylestrong{Arguments}
\begin{quote}

\sphinxAtStartPar
\sphinxcode{\sphinxupquote{idx}}
\begin{quote}

\sphinxAtStartPar
Index of the variable in the expression, starting with 0.
\end{quote}
\end{quote}

\sphinxAtStartPar
\sphinxstylestrong{Example}
\end{quote}

\begin{sphinxVerbatim}[commandchars=\\\{\}]
\PYG{c+c1}{\PYGZsh{} Retrieve the variable for the term at index 1 from expression expr}
\PYG{n}{x} \PYG{o}{=} \PYG{n}{expr}\PYG{o}{.}\PYG{n}{getVar}\PYG{p}{(}\PYG{l+m+mi}{1}\PYG{p}{)}
\end{sphinxVerbatim}

\subsubsection{LinExpr.getConstant()}
\label{\detokenize{pyapiref:linexpr-getconstant}}\begin{quote}

\sphinxAtStartPar
\sphinxstylestrong{Synopsis}
\begin{quote}

\sphinxAtStartPar
\sphinxcode{\sphinxupquote{getConstant()}}
\end{quote}

\sphinxAtStartPar
\sphinxstylestrong{Description}
\begin{quote}

\sphinxAtStartPar
Retrieve the constant term from an expression.
\end{quote}

\sphinxAtStartPar
\sphinxstylestrong{Example}
\end{quote}

\begin{sphinxVerbatim}[commandchars=\\\{\}]
\PYG{c+c1}{\PYGZsh{} Retrieve the constant term from linear expression expr}
\PYG{n}{constant} \PYG{o}{=} \PYG{n}{expr}\PYG{o}{.}\PYG{n}{getConstant}\PYG{p}{(}\PYG{p}{)}
\end{sphinxVerbatim}

\subsubsection{LinExpr.getValue()}
\label{\detokenize{pyapiref:linexpr-getvalue}}\begin{quote}

\sphinxAtStartPar
\sphinxstylestrong{Synopsis}
\begin{quote}

\sphinxAtStartPar
\sphinxcode{\sphinxupquote{getValue()}}
\end{quote}

\sphinxAtStartPar
\sphinxstylestrong{Description}
\begin{quote}

\sphinxAtStartPar
Retrieve the value of an expression computed using the current solution.
\end{quote}

\sphinxAtStartPar
\sphinxstylestrong{Example}
\end{quote}

\begin{sphinxVerbatim}[commandchars=\\\{\}]
\PYG{c+c1}{\PYGZsh{} Retrieve the value of expression expr for the current solution.}
\PYG{n}{val} \PYG{o}{=} \PYG{n}{expr}\PYG{o}{.}\PYG{n}{getValue}\PYG{p}{(}\PYG{p}{)}
\end{sphinxVerbatim}

\subsubsection{LinExpr.getSize()}
\label{\detokenize{pyapiref:linexpr-getsize}}\begin{quote}

\sphinxAtStartPar
\sphinxstylestrong{Synopsis}
\begin{quote}

\sphinxAtStartPar
\sphinxcode{\sphinxupquote{getSize()}}
\end{quote}

\sphinxAtStartPar
\sphinxstylestrong{Description}
\begin{quote}

\sphinxAtStartPar
Retrieve the number of terms in an expression.
\end{quote}

\sphinxAtStartPar
\sphinxstylestrong{Example}
\end{quote}

\begin{sphinxVerbatim}[commandchars=\\\{\}]
\PYG{c+c1}{\PYGZsh{} Retrieve the number of terms in expression expr}
\PYG{n}{exprsize} \PYG{o}{=} \PYG{n}{expr}\PYG{o}{.}\PYG{n}{getSize}\PYG{p}{(}\PYG{p}{)}
\end{sphinxVerbatim}

\subsubsection{LinExpr.setConstant()}
\label{\detokenize{pyapiref:linexpr-setconstant}}\begin{quote}

\sphinxAtStartPar
\sphinxstylestrong{Synopsis}
\begin{quote}

\sphinxAtStartPar
\sphinxcode{\sphinxupquote{setConstant(newval)}}
\end{quote}

\sphinxAtStartPar
\sphinxstylestrong{Description}
\begin{quote}

\sphinxAtStartPar
Set the constant term of linear expression.
\end{quote}

\sphinxAtStartPar
\sphinxstylestrong{Arguments}
\begin{quote}

\sphinxAtStartPar
\sphinxcode{\sphinxupquote{newval}}
\begin{quote}

\sphinxAtStartPar
Constant term to be set.
\end{quote}
\end{quote}

\sphinxAtStartPar
\sphinxstylestrong{Example}
\end{quote}

\begin{sphinxVerbatim}[commandchars=\\\{\}]
\PYG{c+c1}{\PYGZsh{} Set constant term of linear expression \PYGZsq{}expr\PYGZsq{} to 2.0}
\PYG{n}{expr}\PYG{o}{.}\PYG{n}{setConstant}\PYG{p}{(}\PYG{l+m+mf}{2.0}\PYG{p}{)}
\end{sphinxVerbatim}

\subsubsection{LinExpr.addConstant()}
\label{\detokenize{pyapiref:linexpr-addconstant}}\begin{quote}

\sphinxAtStartPar
\sphinxstylestrong{Synopsis}
\begin{quote}

\sphinxAtStartPar
\sphinxcode{\sphinxupquote{addConstant(newval)}}
\end{quote}

\sphinxAtStartPar
\sphinxstylestrong{Description}
\begin{quote}

\sphinxAtStartPar
Add a constant to an expression.
\end{quote}

\sphinxAtStartPar
\sphinxstylestrong{Arguments}
\begin{quote}

\sphinxAtStartPar
\sphinxcode{\sphinxupquote{newval}}
\begin{quote}

\sphinxAtStartPar
Constant to add.
\end{quote}
\end{quote}

\sphinxAtStartPar
\sphinxstylestrong{Example}
\end{quote}

\begin{sphinxVerbatim}[commandchars=\\\{\}]
\PYG{c+c1}{\PYGZsh{} Add constant 2.0 to linear expression \PYGZsq{}expr\PYGZsq{}}
\PYG{n}{expr}\PYG{o}{.}\PYG{n}{addConstant}\PYG{p}{(}\PYG{l+m+mf}{2.0}\PYG{p}{)}
\end{sphinxVerbatim}

\subsubsection{LinExpr.addTerm()}
\label{\detokenize{pyapiref:linexpr-addterm}}\begin{quote}

\sphinxAtStartPar
\sphinxstylestrong{Synopsis}
\begin{quote}

\sphinxAtStartPar
\sphinxcode{\sphinxupquote{addTerm(var, coeff=1.0)}}
\end{quote}

\sphinxAtStartPar
\sphinxstylestrong{Description}
\begin{quote}

\sphinxAtStartPar
Add a new term to current expression.
\end{quote}

\sphinxAtStartPar
\sphinxstylestrong{Arguments}
\begin{quote}

\sphinxAtStartPar
\sphinxcode{\sphinxupquote{var}}
\begin{quote}

\sphinxAtStartPar
Variable to add.
\end{quote}

\sphinxAtStartPar
\sphinxcode{\sphinxupquote{coeff}}
\begin{quote}

\sphinxAtStartPar
Magnification coefficient for added term. Optional, 1.0 by default.
\end{quote}
\end{quote}

\sphinxAtStartPar
\sphinxstylestrong{Example}
\end{quote}

\begin{sphinxVerbatim}[commandchars=\\\{\}]
\PYG{c+c1}{\PYGZsh{} Add term x to linear expression \PYGZsq{}expr\PYGZsq{}}
\PYG{n}{expr}\PYG{o}{.}\PYG{n}{addTerm}\PYG{p}{(}\PYG{n}{x}\PYG{p}{)}
\end{sphinxVerbatim}

\subsubsection{LinExpr.addTerms()}
\label{\detokenize{pyapiref:linexpr-addterms}}\begin{quote}

\sphinxAtStartPar
\sphinxstylestrong{Synopsis}
\begin{quote}

\sphinxAtStartPar
\sphinxcode{\sphinxupquote{addTerms(vars, coeffs)}}
\end{quote}

\sphinxAtStartPar
\sphinxstylestrong{Description}
\begin{quote}

\sphinxAtStartPar
Add a single term or multiple terms into an expression.

\sphinxAtStartPar
If argument \sphinxcode{\sphinxupquote{vars}} is {\hyperref[\detokenize{pyapiref:chappyapi-var}]{\sphinxcrossref{\DUrole{std,std-ref}{Var Class}}}} object, then argument \sphinxcode{\sphinxupquote{coeffs}} is constant;
If argument \sphinxcode{\sphinxupquote{vars}} is {\hyperref[\detokenize{pyapiref:chappyapi-vararray}]{\sphinxcrossref{\DUrole{std,std-ref}{VarArray Class}}}} object or list, then argument \sphinxcode{\sphinxupquote{coeffs}} is constant or list;
If argument \sphinxcode{\sphinxupquote{vars}} is dictionary or {\hyperref[\detokenize{pyapiref:chappyapi-util-tupledict}]{\sphinxcrossref{\DUrole{std,std-ref}{tupledict Class}}}} object, then argument \sphinxcode{\sphinxupquote{coeffs}} is constant, dict, or
{\hyperref[\detokenize{pyapiref:chappyapi-util-tupledict}]{\sphinxcrossref{\DUrole{std,std-ref}{tupledict Class}}}} object.
\end{quote}

\sphinxAtStartPar
\sphinxstylestrong{Arguments}
\begin{quote}

\sphinxAtStartPar
\sphinxcode{\sphinxupquote{vars}}
\begin{quote}

\sphinxAtStartPar
Variables to add.
\end{quote}

\sphinxAtStartPar
\sphinxcode{\sphinxupquote{coeffs}}
\begin{quote}

\sphinxAtStartPar
Coefficients for variables.
\end{quote}
\end{quote}

\sphinxAtStartPar
\sphinxstylestrong{Example}
\end{quote}

\begin{sphinxVerbatim}[commandchars=\\\{\}]
\PYG{c+c1}{\PYGZsh{} Add term 2*x + 2*y to linear expression \PYGZsq{}expr\PYGZsq{}}
\PYG{n}{expr}\PYG{o}{.}\PYG{n}{addTerms}\PYG{p}{(}\PYG{p}{[}\PYG{n}{x}\PYG{p}{,} \PYG{n}{y}\PYG{p}{]}\PYG{p}{,} \PYG{p}{[}\PYG{l+m+mf}{2.0}\PYG{p}{,} \PYG{l+m+mf}{3.0}\PYG{p}{]}\PYG{p}{)}
\end{sphinxVerbatim}

\subsubsection{LinExpr.addExpr()}
\label{\detokenize{pyapiref:linexpr-addexpr}}\begin{quote}

\sphinxAtStartPar
\sphinxstylestrong{Synopsis}
\begin{quote}

\sphinxAtStartPar
\sphinxcode{\sphinxupquote{addExpr(expr, coeff=1.0)}}
\end{quote}

\sphinxAtStartPar
\sphinxstylestrong{Description}
\begin{quote}

\sphinxAtStartPar
Add new expression to the current one.
\end{quote}

\sphinxAtStartPar
\sphinxstylestrong{Arguments}
\begin{quote}

\sphinxAtStartPar
\sphinxcode{\sphinxupquote{expr}}
\begin{quote}

\sphinxAtStartPar
Expression or expression builder to add.
\end{quote}

\sphinxAtStartPar
\sphinxcode{\sphinxupquote{coeff}}
\begin{quote}

\sphinxAtStartPar
Magnification coefficients for the added expression. Optional, 1.0 by default.
\end{quote}
\end{quote}

\sphinxAtStartPar
\sphinxstylestrong{Example}
\end{quote}

\begin{sphinxVerbatim}[commandchars=\\\{\}]
\PYG{c+c1}{\PYGZsh{} Add linear expression 2*x + 2*y to \PYGZsq{}expr\PYGZsq{}}
\PYG{n}{expr}\PYG{o}{.}\PYG{n}{addExpr}\PYG{p}{(}\PYG{n}{x} \PYG{o}{+} \PYG{n}{y}\PYG{p}{,} \PYG{l+m+mf}{2.0}\PYG{p}{)}
\end{sphinxVerbatim}

\subsubsection{LinExpr.clone()}
\label{\detokenize{pyapiref:linexpr-clone}}\begin{quote}

\sphinxAtStartPar
\sphinxstylestrong{Synopsis}
\begin{quote}

\sphinxAtStartPar
\sphinxcode{\sphinxupquote{clone()}}
\end{quote}

\sphinxAtStartPar
\sphinxstylestrong{Description}
\begin{quote}

\sphinxAtStartPar
Create a deep copy of the expression.
\end{quote}

\sphinxAtStartPar
\sphinxstylestrong{Example}
\end{quote}

\begin{sphinxVerbatim}[commandchars=\\\{\}]
\PYG{c+c1}{\PYGZsh{} Create a deep copy of expression expr}
\PYG{n}{exprcopy} \PYG{o}{=} \PYG{n}{expr}\PYG{o}{.}\PYG{n}{clone}\PYG{p}{(}\PYG{p}{)}
\end{sphinxVerbatim}

\subsubsection{LinExpr.reserve()}
\label{\detokenize{pyapiref:linexpr-reserve}}\begin{quote}

\sphinxAtStartPar
\sphinxstylestrong{Synopsis}
\begin{quote}

\sphinxAtStartPar
\sphinxcode{\sphinxupquote{reserve(n)}}
\end{quote}

\sphinxAtStartPar
\sphinxstylestrong{Description}
\begin{quote}

\sphinxAtStartPar
Pre\sphinxhyphen{}allocate space for linear expression object.
\end{quote}

\sphinxAtStartPar
\sphinxstylestrong{Arguments}
\begin{quote}

\sphinxAtStartPar
\sphinxcode{\sphinxupquote{n}}
\begin{quote}

\sphinxAtStartPar
Number of terms to be allocated.
\end{quote}
\end{quote}

\sphinxAtStartPar
\sphinxstylestrong{Example}
\end{quote}

\begin{sphinxVerbatim}[commandchars=\\\{\}]
\PYG{c+c1}{\PYGZsh{} Allocate 100 terms for linear expression \PYGZsq{}expr\PYGZsq{}}
\PYG{n}{expr}\PYG{o}{.}\PYG{n}{reserve}\PYG{p}{(}\PYG{l+m+mi}{100}\PYG{p}{)}
\end{sphinxVerbatim}

\subsubsection{LinExpr.remove()}
\label{\detokenize{pyapiref:linexpr-remove}}\begin{quote}

\sphinxAtStartPar
\sphinxstylestrong{Synopsis}
\begin{quote}

\sphinxAtStartPar
\sphinxcode{\sphinxupquote{remove(item)}}
\end{quote}

\sphinxAtStartPar
\sphinxstylestrong{Description}
\begin{quote}

\sphinxAtStartPar
Remove a term from a linear expression.

\sphinxAtStartPar
If argument \sphinxcode{\sphinxupquote{item}} is constant, then remove the term stored at index i of the expression; otherwise argument \sphinxcode{\sphinxupquote{item}} is
{\hyperref[\detokenize{pyapiref:chappyapi-var}]{\sphinxcrossref{\DUrole{std,std-ref}{Var Class}}}} object.
\end{quote}

\sphinxAtStartPar
\sphinxstylestrong{Arguments}
\begin{quote}

\sphinxAtStartPar
\sphinxcode{\sphinxupquote{item}}
\begin{quote}

\sphinxAtStartPar
Constant index or variable of the term to be removed.
\end{quote}
\end{quote}

\sphinxAtStartPar
\sphinxstylestrong{Example}
\end{quote}

\begin{sphinxVerbatim}[commandchars=\\\{\}]
\PYG{c+c1}{\PYGZsh{} Remove the term whose index is 2 from linear expression expr}
\PYG{n}{expr}\PYG{o}{.}\PYG{n}{remove}\PYG{p}{(}\PYG{l+m+mi}{2}\PYG{p}{)}
\PYG{c+c1}{\PYGZsh{} Remove the term whose variable is x from linear expression expr}
\PYG{n}{expr}\PYG{o}{.}\PYG{n}{remove}\PYG{p}{(}\PYG{n}{x}\PYG{p}{)}
\end{sphinxVerbatim}

\subsection{QuadExpr Class}
\label{\detokenize{pyapiref:quadexpr-class}}\label{\detokenize{pyapiref:chappyapi-quadexpr}}
\sphinxAtStartPar
QuadExpr object contains operations related to variables for building quadratic constraints, and provides the following methods:

\subsubsection{QuadExpr()}
\label{\detokenize{pyapiref:quadexpr}}\begin{quote}

\sphinxAtStartPar
\sphinxstylestrong{Synopsis}
\begin{quote}

\sphinxAtStartPar
\sphinxcode{\sphinxupquote{QuadExpr(expr=0.0)}}
\end{quote}

\sphinxAtStartPar
\sphinxstylestrong{Description}
\begin{quote}

\sphinxAtStartPar
Create a {\hyperref[\detokenize{pyapiref:chappyapi-quadexpr}]{\sphinxcrossref{\DUrole{std,std-ref}{QuadExpr Class}}}} object.

\sphinxAtStartPar
Argument \sphinxcode{\sphinxupquote{expr}} is constant, {\hyperref[\detokenize{pyapiref:chappyapi-var}]{\sphinxcrossref{\DUrole{std,std-ref}{Var Class}}}}, {\hyperref[\detokenize{pyapiref:chappyapi-linexpr}]{\sphinxcrossref{\DUrole{std,std-ref}{LinExpr Class}}}} object or {\hyperref[\detokenize{pyapiref:chappyapi-quadexpr}]{\sphinxcrossref{\DUrole{std,std-ref}{QuadExpr Class}}}} object.
\end{quote}

\sphinxAtStartPar
\sphinxstylestrong{Arguments}
\begin{quote}

\sphinxAtStartPar
\sphinxcode{\sphinxupquote{expr}}
\begin{quote}

\sphinxAtStartPar
Optional, 0.0 by default.
\end{quote}
\end{quote}

\sphinxAtStartPar
\sphinxstylestrong{Example}
\end{quote}

\begin{sphinxVerbatim}[commandchars=\\\{\}]
\PYG{c+c1}{\PYGZsh{} Create a new QuadExpr object and initialize it to 0.0}
\PYG{n}{quadexpr0} \PYG{o}{=} \PYG{n}{QuadExpr}\PYG{p}{(}\PYG{p}{)}
\PYG{c+c1}{\PYGZsh{} Create a QuadExpr object and initialize it to 2*x + 3*y}
\PYG{n}{quadexpr1} \PYG{o}{=} \PYG{n}{QuadExpr}\PYG{p}{(}\PYG{p}{[}\PYG{p}{(}\PYG{n}{x}\PYG{p}{,} \PYG{l+m+mi}{2}\PYG{p}{)}\PYG{p}{,} \PYG{p}{(}\PYG{n}{y}\PYG{p}{,} \PYG{l+m+mi}{3}\PYG{p}{)}\PYG{p}{]}\PYG{p}{)}
\PYG{c+c1}{\PYGZsh{} Create a QuadExpr object and initialize it to x*x + 2*x*y}
\PYG{n}{quadexpr2} \PYG{o}{=} \PYG{n}{QuadExpr}\PYG{p}{(}\PYG{n}{x}\PYG{o}{*}\PYG{n}{x} \PYG{o}{+} \PYG{l+m+mi}{2}\PYG{o}{*}\PYG{n}{x}\PYG{o}{*}\PYG{n}{y}\PYG{p}{)}
\end{sphinxVerbatim}

\subsubsection{QuadExpr.setCoeff()}
\label{\detokenize{pyapiref:quadexpr-setcoeff}}\begin{quote}

\sphinxAtStartPar
\sphinxstylestrong{Synopsis}
\begin{quote}

\sphinxAtStartPar
\sphinxcode{\sphinxupquote{setCoeff(idx, newval)}}
\end{quote}

\sphinxAtStartPar
\sphinxstylestrong{Description}
\begin{quote}

\sphinxAtStartPar
Set new coefficient value of a term based on its index in a quadratic expression.
\end{quote}

\sphinxAtStartPar
\sphinxstylestrong{Arguments}
\begin{quote}

\sphinxAtStartPar
\sphinxcode{\sphinxupquote{idx}}
\begin{quote}

\sphinxAtStartPar
Index of the term in the quadratic expression, starting with 0.
\end{quote}

\sphinxAtStartPar
\sphinxcode{\sphinxupquote{newval}}
\begin{quote}

\sphinxAtStartPar
New coefficient value of the term.
\end{quote}
\end{quote}

\sphinxAtStartPar
\sphinxstylestrong{Example}
\end{quote}

\begin{sphinxVerbatim}[commandchars=\\\{\}]
\PYG{c+c1}{\PYGZsh{} Set the coefficient for the term at index 0 in quadratic expression quadexpr to 1.0}
\PYG{n}{quadexpr}\PYG{o}{.}\PYG{n}{setCoeff}\PYG{p}{(}\PYG{l+m+mi}{0}\PYG{p}{,} \PYG{l+m+mf}{1.0}\PYG{p}{)}
\end{sphinxVerbatim}

\subsubsection{QuadExpr.getCoeff()}
\label{\detokenize{pyapiref:quadexpr-getcoeff}}\begin{quote}

\sphinxAtStartPar
\sphinxstylestrong{Synopsis}
\begin{quote}

\sphinxAtStartPar
\sphinxcode{\sphinxupquote{getCoeff(idx)}}
\end{quote}

\sphinxAtStartPar
\sphinxstylestrong{Description}
\begin{quote}

\sphinxAtStartPar
Retrieve the coefficient of a term by its index from a quadratic expression.
\end{quote}

\sphinxAtStartPar
\sphinxstylestrong{Arguments}
\begin{quote}

\sphinxAtStartPar
\sphinxcode{\sphinxupquote{idx}}
\begin{quote}

\sphinxAtStartPar
Index of the term in the quadratic expression, starting with 0.
\end{quote}
\end{quote}

\sphinxAtStartPar
\sphinxstylestrong{Example}
\end{quote}

\begin{sphinxVerbatim}[commandchars=\\\{\}]
\PYG{c+c1}{\PYGZsh{} Retrieve the coefficient for the term at index 1 from quadratic expression quadexpr}
\PYG{n}{coeff} \PYG{o}{=} \PYG{n}{quadexpr}\PYG{o}{.}\PYG{n}{getCoeff}\PYG{p}{(}\PYG{l+m+mi}{1}\PYG{p}{)}
\end{sphinxVerbatim}

\subsubsection{QuadExpr.getVar1()}
\label{\detokenize{pyapiref:quadexpr-getvar1}}\begin{quote}

\sphinxAtStartPar
\sphinxstylestrong{Synopsis}
\begin{quote}

\sphinxAtStartPar
\sphinxcode{\sphinxupquote{getVar1(idx)}}
\end{quote}

\sphinxAtStartPar
\sphinxstylestrong{Description}
\begin{quote}

\sphinxAtStartPar
Retrieve the first variable of a quadratic term by its index from an expression. Return a {\hyperref[\detokenize{pyapiref:chappyapi-var}]{\sphinxcrossref{\DUrole{std,std-ref}{Var Class}}}} object.
\end{quote}

\sphinxAtStartPar
\sphinxstylestrong{Arguments}
\begin{quote}

\sphinxAtStartPar
\sphinxcode{\sphinxupquote{idx}}
\begin{quote}

\sphinxAtStartPar
Index of the quadratic term in the expression, starting with 0.
\end{quote}
\end{quote}

\sphinxAtStartPar
\sphinxstylestrong{Example}
\end{quote}

\begin{sphinxVerbatim}[commandchars=\\\{\}]
\PYG{c+c1}{\PYGZsh{} Retrieve the first variable of a quadratic term at index 1 from quadratic expression quadexpr}
\PYG{n}{x} \PYG{o}{=} \PYG{n}{quadexpr}\PYG{o}{.}\PYG{n}{getVar1}\PYG{p}{(}\PYG{l+m+mi}{1}\PYG{p}{)}
\end{sphinxVerbatim}

\subsubsection{QuadExpr.getVar2()}
\label{\detokenize{pyapiref:quadexpr-getvar2}}\begin{quote}

\sphinxAtStartPar
\sphinxstylestrong{Synopsis}
\begin{quote}

\sphinxAtStartPar
\sphinxcode{\sphinxupquote{getVar2(idx)}}
\end{quote}

\sphinxAtStartPar
\sphinxstylestrong{Description}
\begin{quote}

\sphinxAtStartPar
Retrieve the second variable of a quadratic term by its index from an expression. Return a {\hyperref[\detokenize{pyapiref:chappyapi-var}]{\sphinxcrossref{\DUrole{std,std-ref}{Var Class}}}} object.
\end{quote}

\sphinxAtStartPar
\sphinxstylestrong{Arguments}
\begin{quote}

\sphinxAtStartPar
\sphinxcode{\sphinxupquote{idx}}
\begin{quote}

\sphinxAtStartPar
Index of the quadratic term in the expression, starting with 0.
\end{quote}
\end{quote}

\sphinxAtStartPar
\sphinxstylestrong{Example}
\end{quote}

\begin{sphinxVerbatim}[commandchars=\\\{\}]
\PYG{c+c1}{\PYGZsh{} Retrieve the first variable of a quadratic term at index 1 from quadratic expression quadexpr}
\PYG{n}{y} \PYG{o}{=} \PYG{n}{quadexpr}\PYG{o}{.}\PYG{n}{getVar2}\PYG{p}{(}\PYG{l+m+mi}{1}\PYG{p}{)}
\end{sphinxVerbatim}

\subsubsection{QuadExpr.getLinExpr()}
\label{\detokenize{pyapiref:quadexpr-getlinexpr}}\begin{quote}

\sphinxAtStartPar
\sphinxstylestrong{Synopsis}
\begin{quote}

\sphinxAtStartPar
\sphinxcode{\sphinxupquote{getLinExpr()}}
\end{quote}

\sphinxAtStartPar
\sphinxstylestrong{Description}
\begin{quote}

\sphinxAtStartPar
Retrieve the linear terms (if exist) fronm quadratic expression. Return a {\hyperref[\detokenize{pyapiref:chappyapi-linexpr}]{\sphinxcrossref{\DUrole{std,std-ref}{LinExpr Class}}}} object.
\end{quote}

\sphinxAtStartPar
\sphinxstylestrong{Example}
\end{quote}

\begin{sphinxVerbatim}[commandchars=\\\{\}]
\PYG{c+c1}{\PYGZsh{} Retrieve the linear terms from a quadratic expression quadexpr}
\PYG{n}{linexpr} \PYG{o}{=} \PYG{n}{quadexpr}\PYG{o}{.}\PYG{n}{getLinExpr}\PYG{p}{(}\PYG{p}{)}
\end{sphinxVerbatim}

\subsubsection{QuadExpr.getConstant()}
\label{\detokenize{pyapiref:quadexpr-getconstant}}\begin{quote}

\sphinxAtStartPar
\sphinxstylestrong{Synopsis}
\begin{quote}

\sphinxAtStartPar
\sphinxcode{\sphinxupquote{getConstant()}}
\end{quote}

\sphinxAtStartPar
\sphinxstylestrong{Description}
\begin{quote}

\sphinxAtStartPar
Retrieve the constant term from a quadratic expression.
\end{quote}

\sphinxAtStartPar
\sphinxstylestrong{Example}
\end{quote}

\begin{sphinxVerbatim}[commandchars=\\\{\}]
\PYG{c+c1}{\PYGZsh{} Retrieve the constant term from quadratic expression quadexpr}
\PYG{n}{constant} \PYG{o}{=} \PYG{n}{quadexpr}\PYG{o}{.}\PYG{n}{getConstant}\PYG{p}{(}\PYG{p}{)}
\end{sphinxVerbatim}

\subsubsection{QuadExpr.getValue()}
\label{\detokenize{pyapiref:quadexpr-getvalue}}\begin{quote}

\sphinxAtStartPar
\sphinxstylestrong{Synopsis}
\begin{quote}

\sphinxAtStartPar
\sphinxcode{\sphinxupquote{getValue()}}
\end{quote}

\sphinxAtStartPar
\sphinxstylestrong{Description}
\begin{quote}

\sphinxAtStartPar
Retrieve the value of a quadratic expression computed using the current solution.
\end{quote}

\sphinxAtStartPar
\sphinxstylestrong{Example}
\end{quote}

\begin{sphinxVerbatim}[commandchars=\\\{\}]
\PYG{c+c1}{\PYGZsh{} Retrieve the value of quadratic expression quadexpr for the current solution.}
\PYG{n}{val} \PYG{o}{=} \PYG{n}{quadexpr}\PYG{o}{.}\PYG{n}{getValue}\PYG{p}{(}\PYG{p}{)}
\end{sphinxVerbatim}

\subsubsection{QuadExpr.getSize()}
\label{\detokenize{pyapiref:quadexpr-getsize}}\begin{quote}

\sphinxAtStartPar
\sphinxstylestrong{Synopsis}
\begin{quote}

\sphinxAtStartPar
\sphinxcode{\sphinxupquote{getSize()}}
\end{quote}

\sphinxAtStartPar
\sphinxstylestrong{Description}
\begin{quote}

\sphinxAtStartPar
Retrieve the number of terms in a quadratic expression.
\end{quote}

\sphinxAtStartPar
\sphinxstylestrong{Example}
\end{quote}

\begin{sphinxVerbatim}[commandchars=\\\{\}]
\PYG{c+c1}{\PYGZsh{} Retrieve the number of terms in quadratic expression quadexpr}
\PYG{n}{exprsize} \PYG{o}{=} \PYG{n}{quadexpr}\PYG{o}{.}\PYG{n}{getSize}\PYG{p}{(}\PYG{p}{)}
\end{sphinxVerbatim}

\subsubsection{QuadExpr.setConstant()}
\label{\detokenize{pyapiref:quadexpr-setconstant}}\begin{quote}

\sphinxAtStartPar
\sphinxstylestrong{Synopsis}
\begin{quote}

\sphinxAtStartPar
\sphinxcode{\sphinxupquote{setConstant(newval)}}
\end{quote}

\sphinxAtStartPar
\sphinxstylestrong{Description}
\begin{quote}

\sphinxAtStartPar
Set the constant term of quadratic expression.
\end{quote}

\sphinxAtStartPar
\sphinxstylestrong{Arguments}
\begin{quote}

\sphinxAtStartPar
\sphinxcode{\sphinxupquote{newval}}
\begin{quote}

\sphinxAtStartPar
Constant to set.
\end{quote}
\end{quote}

\sphinxAtStartPar
\sphinxstylestrong{Example}
\end{quote}

\begin{sphinxVerbatim}[commandchars=\\\{\}]
\PYG{c+c1}{\PYGZsh{} Set constant term of quadratic expression \PYGZsq{}quadexpr\PYGZsq{} to 2.0}
\PYG{n}{quadexpr}\PYG{o}{.}\PYG{n}{setConstant}\PYG{p}{(}\PYG{l+m+mf}{2.0}\PYG{p}{)}
\end{sphinxVerbatim}

\subsubsection{QuadExpr.addConstant()}
\label{\detokenize{pyapiref:quadexpr-addconstant}}\begin{quote}

\sphinxAtStartPar
\sphinxstylestrong{Synopsis}
\begin{quote}

\sphinxAtStartPar
\sphinxcode{\sphinxupquote{addConstant(newval)}}
\end{quote}

\sphinxAtStartPar
\sphinxstylestrong{Description}
\begin{quote}

\sphinxAtStartPar
Add a constant to a quadratic expression.
\end{quote}

\sphinxAtStartPar
\sphinxstylestrong{Arguments}
\begin{quote}

\sphinxAtStartPar
\sphinxcode{\sphinxupquote{newval}}
\begin{quote}

\sphinxAtStartPar
Constant to add.
\end{quote}
\end{quote}

\sphinxAtStartPar
\sphinxstylestrong{Example}
\end{quote}

\begin{sphinxVerbatim}[commandchars=\\\{\}]
\PYG{c+c1}{\PYGZsh{} Add constant 2.0 to quadratic expression \PYGZsq{}quadexpr\PYGZsq{}}
\PYG{n}{quadexpr}\PYG{o}{.}\PYG{n}{addConstant}\PYG{p}{(}\PYG{l+m+mf}{2.0}\PYG{p}{)}
\end{sphinxVerbatim}

\subsubsection{QuadExpr.addTerm()}
\label{\detokenize{pyapiref:quadexpr-addterm}}\begin{quote}

\sphinxAtStartPar
\sphinxstylestrong{Synopsis}
\begin{quote}

\sphinxAtStartPar
\sphinxcode{\sphinxupquote{addTerm(coeff, var1, var2=None)}}
\end{quote}

\sphinxAtStartPar
\sphinxstylestrong{Description}
\begin{quote}

\sphinxAtStartPar
Add a new term to current quadratic expression.
\end{quote}

\sphinxAtStartPar
\sphinxstylestrong{Arguments}
\begin{quote}

\sphinxAtStartPar
\sphinxcode{\sphinxupquote{coeff}}
\begin{quote}

\sphinxAtStartPar
Magnification coefficient for added term. Optional, 1.0 by default.
\end{quote}

\sphinxAtStartPar
\sphinxcode{\sphinxupquote{var1}}
\begin{quote}

\sphinxAtStartPar
The first variable for added term.
\end{quote}

\sphinxAtStartPar
\sphinxcode{\sphinxupquote{var2}}
\begin{quote}

\sphinxAtStartPar
The second variable for added term, defaults to \sphinxcode{\sphinxupquote{None}}, i.e. add a linear term.
\end{quote}
\end{quote}

\sphinxAtStartPar
\sphinxstylestrong{Example}
\end{quote}

\begin{sphinxVerbatim}[commandchars=\\\{\}]
\PYG{c+c1}{\PYGZsh{} Add term x to quadratic expression \PYGZsq{}quadexpr\PYGZsq{}}
\PYG{n}{quadexpr}\PYG{o}{.}\PYG{n}{addTerm}\PYG{p}{(}\PYG{l+m+mf}{1.0}\PYG{p}{,} \PYG{n}{x}\PYG{p}{)}
\end{sphinxVerbatim}

\subsubsection{QuadExpr.addTerms()}
\label{\detokenize{pyapiref:quadexpr-addterms}}\begin{quote}

\sphinxAtStartPar
\sphinxstylestrong{Synopsis}
\begin{quote}

\sphinxAtStartPar
\sphinxcode{\sphinxupquote{addTerms(coeffs, vars1, vars2=None)}}
\end{quote}

\sphinxAtStartPar
\sphinxstylestrong{Description}
\begin{quote}

\sphinxAtStartPar
Add a single term or multiple terms into a quadratic expression.

\sphinxAtStartPar
If argument \sphinxcode{\sphinxupquote{vars}} is {\hyperref[\detokenize{pyapiref:chappyapi-var}]{\sphinxcrossref{\DUrole{std,std-ref}{Var Class}}}} object, then argument \sphinxcode{\sphinxupquote{vars2}} is
{\hyperref[\detokenize{pyapiref:chappyapi-var}]{\sphinxcrossref{\DUrole{std,std-ref}{Var Class}}}} object or \sphinxcode{\sphinxupquote{None}}, argument \sphinxcode{\sphinxupquote{coeffs}} is constant;
If argument \sphinxcode{\sphinxupquote{vars}} is {\hyperref[\detokenize{pyapiref:chappyapi-vararray}]{\sphinxcrossref{\DUrole{std,std-ref}{VarArray Class}}}} object or list, then argument
\sphinxcode{\sphinxupquote{vars2}} is {\hyperref[\detokenize{pyapiref:chappyapi-vararray}]{\sphinxcrossref{\DUrole{std,std-ref}{VarArray Class}}}} object, list or \sphinxcode{\sphinxupquote{None}}, argument \sphinxcode{\sphinxupquote{coeffs}}
is constant or list;
If argument \sphinxcode{\sphinxupquote{vars}} is dictionary or {\hyperref[\detokenize{pyapiref:chappyapi-util-tupledict}]{\sphinxcrossref{\DUrole{std,std-ref}{tupledict Class}}}} object, then
argument \sphinxcode{\sphinxupquote{vars2}} is dictionary, {\hyperref[\detokenize{pyapiref:chappyapi-util-tupledict}]{\sphinxcrossref{\DUrole{std,std-ref}{tupledict Class}}}} object or \sphinxcode{\sphinxupquote{None}},
argument \sphinxcode{\sphinxupquote{coeffs}} is constant, dictionary, or
{\hyperref[\detokenize{pyapiref:chappyapi-util-tupledict}]{\sphinxcrossref{\DUrole{std,std-ref}{tupledict Class}}}} object.
\end{quote}

\sphinxAtStartPar
\sphinxstylestrong{Arguments}
\begin{quote}

\sphinxAtStartPar
\sphinxcode{\sphinxupquote{coeffs}}
\begin{quote}

\sphinxAtStartPar
Coefficients for terms.
\end{quote}

\sphinxAtStartPar
\sphinxcode{\sphinxupquote{vars1}}
\begin{quote}

\sphinxAtStartPar
The first variable of each term.
\end{quote}

\sphinxAtStartPar
\sphinxcode{\sphinxupquote{vars2}}
\begin{quote}

\sphinxAtStartPar
The second variable of each term, defaults to \sphinxcode{\sphinxupquote{None}}, i.e. add a linear term.
\end{quote}
\end{quote}

\sphinxAtStartPar
\sphinxstylestrong{Example}
\end{quote}

\begin{sphinxVerbatim}[commandchars=\\\{\}]
\PYG{c+c1}{\PYGZsh{} Add term 2*x + 3y + 2*x*x + 3*x*y to quadratic expression \PYGZsq{}quadexpr\PYGZsq{}}
\PYG{c+c1}{\PYGZsh{} Note: Mixed format is supported by addTerms yet.}
\PYG{n}{quadexpr}\PYG{o}{.}\PYG{n}{addTerms}\PYG{p}{(}\PYG{p}{[}\PYG{l+m+mf}{2.0}\PYG{p}{,} \PYG{l+m+mf}{3.0}\PYG{p}{]}\PYG{p}{,} \PYG{p}{[}\PYG{n}{x}\PYG{p}{,} \PYG{n}{y}\PYG{p}{]}\PYG{p}{)}
\PYG{n}{quadexpr}\PYG{o}{.}\PYG{n}{addTerms}\PYG{p}{(}\PYG{p}{[}\PYG{l+m+mf}{2.0}\PYG{p}{,} \PYG{l+m+mf}{3.0}\PYG{p}{]}\PYG{p}{,} \PYG{p}{[}\PYG{n}{x}\PYG{p}{,} \PYG{n}{x}\PYG{p}{]}\PYG{p}{,} \PYG{p}{[}\PYG{n}{x}\PYG{p}{,} \PYG{n}{y}\PYG{p}{]}\PYG{p}{)}
\end{sphinxVerbatim}

\subsubsection{QuadExpr.addLinExpr()}
\label{\detokenize{pyapiref:quadexpr-addlinexpr}}\begin{quote}

\sphinxAtStartPar
\sphinxstylestrong{Synopsis}
\begin{quote}

\sphinxAtStartPar
\sphinxcode{\sphinxupquote{addLinExpr(expr, mult=1.0)}}
\end{quote}

\sphinxAtStartPar
\sphinxstylestrong{Description}
\begin{quote}

\sphinxAtStartPar
Add new linear expression to the current quadratic expression.
\end{quote}

\sphinxAtStartPar
\sphinxstylestrong{Arguments}
\begin{quote}

\sphinxAtStartPar
\sphinxcode{\sphinxupquote{expr}}
\begin{quote}

\sphinxAtStartPar
Linear expression or linear expression builder to add.
\end{quote}

\sphinxAtStartPar
\sphinxcode{\sphinxupquote{mult}}
\begin{quote}

\sphinxAtStartPar
Magnification coefficient for the added expression. Optional, 1.0 by default.
\end{quote}
\end{quote}

\sphinxAtStartPar
\sphinxstylestrong{Example}
\end{quote}

\begin{sphinxVerbatim}[commandchars=\\\{\}]
\PYG{c+c1}{\PYGZsh{} Add linear expression 2*x + 2*y to \PYGZsq{}quadexpr\PYGZsq{}}
\PYG{n}{quadexpr}\PYG{o}{.}\PYG{n}{addLinExpr}\PYG{p}{(}\PYG{n}{x} \PYG{o}{+} \PYG{n}{y}\PYG{p}{,} \PYG{l+m+mf}{2.0}\PYG{p}{)}
\end{sphinxVerbatim}

\subsubsection{QuadExpr.addQuadExpr()}
\label{\detokenize{pyapiref:quadexpr-addquadexpr}}\begin{quote}

\sphinxAtStartPar
\sphinxstylestrong{Synopsis}
\begin{quote}

\sphinxAtStartPar
\sphinxcode{\sphinxupquote{addQuadExpr(expr, mult=1.0)}}
\end{quote}

\sphinxAtStartPar
\sphinxstylestrong{Description}
\begin{quote}

\sphinxAtStartPar
Add new quadratic expression to the current one.
\end{quote}

\sphinxAtStartPar
\sphinxstylestrong{Arguments}
\begin{quote}

\sphinxAtStartPar
\sphinxcode{\sphinxupquote{expr}}
\begin{quote}

\sphinxAtStartPar
Expression or expression builder to add.
\end{quote}

\sphinxAtStartPar
\sphinxcode{\sphinxupquote{mult}}
\begin{quote}

\sphinxAtStartPar
Magnification coefficients for the added expression. Optional, 1.0 by default.
\end{quote}
\end{quote}

\sphinxAtStartPar
\sphinxstylestrong{Example}
\end{quote}

\begin{sphinxVerbatim}[commandchars=\\\{\}]
\PYG{c+c1}{\PYGZsh{} Add quadratic expression x*x + 2*y to \PYGZsq{}quadexpr\PYGZsq{}}
\PYG{n}{quadexpr}\PYG{o}{.}\PYG{n}{addQuadExpr}\PYG{p}{(}\PYG{n}{x}\PYG{o}{*}\PYG{n}{x} \PYG{o}{+} \PYG{l+m+mi}{2}\PYG{o}{*}\PYG{n}{y}\PYG{p}{,} \PYG{l+m+mf}{2.0}\PYG{p}{)}
\end{sphinxVerbatim}

\subsubsection{QuadExpr.clone()}
\label{\detokenize{pyapiref:quadexpr-clone}}\begin{quote}

\sphinxAtStartPar
\sphinxstylestrong{Synopsis}
\begin{quote}

\sphinxAtStartPar
\sphinxcode{\sphinxupquote{clone()}}
\end{quote}

\sphinxAtStartPar
\sphinxstylestrong{Description}
\begin{quote}

\sphinxAtStartPar
Create a deep copy of the expression.
\end{quote}

\sphinxAtStartPar
\sphinxstylestrong{Example}
\end{quote}

\begin{sphinxVerbatim}[commandchars=\\\{\}]
\PYG{c+c1}{\PYGZsh{} Create a deep copy of quadratic expression quadexpr}
\PYG{n}{exprcopy} \PYG{o}{=} \PYG{n}{quadexpr}\PYG{o}{.}\PYG{n}{clone}\PYG{p}{(}\PYG{p}{)}
\end{sphinxVerbatim}

\subsubsection{QuadExpr.reserve()}
\label{\detokenize{pyapiref:quadexpr-reserve}}\begin{quote}

\sphinxAtStartPar
\sphinxstylestrong{Synopsis}
\begin{quote}

\sphinxAtStartPar
\sphinxcode{\sphinxupquote{reserve(n)}}
\end{quote}

\sphinxAtStartPar
\sphinxstylestrong{Description}
\begin{quote}

\sphinxAtStartPar
Pre\sphinxhyphen{}allocate space for quadratic expression object.
\end{quote}

\sphinxAtStartPar
\sphinxstylestrong{Arguments}
\begin{quote}

\sphinxAtStartPar
\sphinxcode{\sphinxupquote{n}}
\begin{quote}

\sphinxAtStartPar
Number of terms to be allocated.
\end{quote}
\end{quote}

\sphinxAtStartPar
\sphinxstylestrong{Example}
\end{quote}

\begin{sphinxVerbatim}[commandchars=\\\{\}]
\PYG{c+c1}{\PYGZsh{} Allocate 100 terms for quadratic expression \PYGZsq{}expr\PYGZsq{}}
\PYG{n}{expr}\PYG{o}{.}\PYG{n}{reserve}\PYG{p}{(}\PYG{l+m+mi}{100}\PYG{p}{)}
\end{sphinxVerbatim}

\subsubsection{QuadExpr.remove()}
\label{\detokenize{pyapiref:quadexpr-remove}}\begin{quote}

\sphinxAtStartPar
\sphinxstylestrong{Synopsis}
\begin{quote}

\sphinxAtStartPar
\sphinxcode{\sphinxupquote{remove(item)}}
\end{quote}

\sphinxAtStartPar
\sphinxstylestrong{Description}
\begin{quote}

\sphinxAtStartPar
Remove a term from a quadratic expression.

\sphinxAtStartPar
If argument \sphinxcode{\sphinxupquote{item}} is constant, then remove the term stored at index i of the expression; otherwise argument \sphinxcode{\sphinxupquote{item}} is
{\hyperref[\detokenize{pyapiref:chappyapi-var}]{\sphinxcrossref{\DUrole{std,std-ref}{Var Class}}}} object.
\end{quote}

\sphinxAtStartPar
\sphinxstylestrong{Arguments}
\begin{quote}

\sphinxAtStartPar
\sphinxcode{\sphinxupquote{item}}
\begin{quote}

\sphinxAtStartPar
Constant index or variable of the term to be removed.
\end{quote}
\end{quote}

\sphinxAtStartPar
\sphinxstylestrong{Example}
\end{quote}

\begin{sphinxVerbatim}[commandchars=\\\{\}]
\PYG{c+c1}{\PYGZsh{} Remove the term whose index is 2 from quadratic expression quadexpr}
\PYG{n}{quadexpr}\PYG{o}{.}\PYG{n}{remove}\PYG{p}{(}\PYG{l+m+mi}{2}\PYG{p}{)}
\PYG{c+c1}{\PYGZsh{} Remove the terms one of which variable is x from quadratic expression quadexpr}
\PYG{n}{quadexpr}\PYG{o}{.}\PYG{n}{remove}\PYG{p}{(}\PYG{n}{x}\PYG{p}{)}
\end{sphinxVerbatim}

\subsection{PsdExpr Class}
\label{\detokenize{pyapiref:psdexpr-class}}\label{\detokenize{pyapiref:chappyapi-psdexpr}}
\sphinxAtStartPar
PsdExpr object contains operations related to variables for building positive semi\sphinxhyphen{}definite constraints,
and provides the following methods:

\subsubsection{PsdExpr()}
\label{\detokenize{pyapiref:psdexpr}}\begin{quote}

\sphinxAtStartPar
\sphinxstylestrong{Synopsis}
\begin{quote}

\sphinxAtStartPar
\sphinxcode{\sphinxupquote{PsdExpr(expr=0.0)}}
\end{quote}

\sphinxAtStartPar
\sphinxstylestrong{Description}
\begin{quote}

\sphinxAtStartPar
Create a {\hyperref[\detokenize{pyapiref:chappyapi-psdexpr}]{\sphinxcrossref{\DUrole{std,std-ref}{PsdExpr Class}}}} object.
\end{quote}

\sphinxAtStartPar
\sphinxstylestrong{Arguments}
\begin{quote}

\sphinxAtStartPar
\sphinxcode{\sphinxupquote{expr}}
\begin{quote}

\sphinxAtStartPar
Optional, 0.0 by default, which can be a constant, {\hyperref[\detokenize{pyapiref:chappyapi-var}]{\sphinxcrossref{\DUrole{std,std-ref}{Var Class}}}},
{\hyperref[\detokenize{pyapiref:chappyapi-linexpr}]{\sphinxcrossref{\DUrole{std,std-ref}{LinExpr Class}}}} object or {\hyperref[\detokenize{pyapiref:chappyapi-psdexpr}]{\sphinxcrossref{\DUrole{std,std-ref}{PsdExpr Class}}}} object.
\end{quote}
\end{quote}

\sphinxAtStartPar
\sphinxstylestrong{Example}
\end{quote}

\begin{sphinxVerbatim}[commandchars=\\\{\}]
\PYG{c+c1}{\PYGZsh{} Create a new PsdExpr object and initialize it to 0.0}
\PYG{n}{expr0} \PYG{o}{=} \PYG{n}{PsdExpr}\PYG{p}{(}\PYG{p}{)}
\PYG{c+c1}{\PYGZsh{} Create a PsdExpr object and initialize it to 2*x + 3*y}
\PYG{n}{expr1} \PYG{o}{=} \PYG{n}{PsdExpr}\PYG{p}{(}\PYG{l+m+mi}{2}\PYG{o}{*}\PYG{n}{x} \PYG{o}{+} \PYG{l+m+mi}{3}\PYG{o}{*}\PYG{n}{y}\PYG{p}{)}
\end{sphinxVerbatim}

\subsubsection{PsdExpr.setCoeff()}
\label{\detokenize{pyapiref:psdexpr-setcoeff}}\begin{quote}

\sphinxAtStartPar
\sphinxstylestrong{Synopsis}
\begin{quote}

\sphinxAtStartPar
\sphinxcode{\sphinxupquote{setCoeff(idx, mat)}}
\end{quote}

\sphinxAtStartPar
\sphinxstylestrong{Description}
\begin{quote}

\sphinxAtStartPar
Set the coefficient symmetric matrix corresponding to the specified index value \sphinxcode{\sphinxupquote{idx}} in the LMI expression.
\end{quote}

\sphinxAtStartPar
\sphinxstylestrong{Arguments}
\begin{quote}

\sphinxAtStartPar
\sphinxcode{\sphinxupquote{idx}}
\begin{quote}

\sphinxAtStartPar
Index of the positive semi\sphinxhyphen{}definite variable in the expression, starting with 0.
\end{quote}

\sphinxAtStartPar
\sphinxcode{\sphinxupquote{mat}}
\begin{quote}

\sphinxAtStartPar
New symmetric matrix coefficient of the positive semi\sphinxhyphen{}definite variable.
\end{quote}
\end{quote}

\sphinxAtStartPar
\sphinxstylestrong{Example}
\end{quote}

\begin{sphinxVerbatim}[commandchars=\\\{\}]
\PYG{c+c1}{\PYGZsh{} Set symmetric matrix for the positive semi\PYGZhy{}definite variable at index 0 in expression \PYGZdq{}expr\PYGZdq{} to mat}
\PYG{n}{expr}\PYG{o}{.}\PYG{n}{setCoeff}\PYG{p}{(}\PYG{l+m+mi}{0}\PYG{p}{,} \PYG{n}{mat}\PYG{p}{)}
\end{sphinxVerbatim}

\subsubsection{PsdExpr.getCoeff()}
\label{\detokenize{pyapiref:psdexpr-getcoeff}}\begin{quote}

\sphinxAtStartPar
\sphinxstylestrong{Synopsis}
\begin{quote}

\sphinxAtStartPar
\sphinxcode{\sphinxupquote{getCoeff(idx)}}
\end{quote}

\sphinxAtStartPar
\sphinxstylestrong{Description}
\begin{quote}

\sphinxAtStartPar
Retrieve the symmetric matrix coefficient of a positive semi\sphinxhyphen{}definite variable by its index from the expression.
\end{quote}

\sphinxAtStartPar
\sphinxstylestrong{Arguments}
\begin{quote}

\sphinxAtStartPar
\sphinxcode{\sphinxupquote{idx}}
\begin{quote}

\sphinxAtStartPar
Index of the positive semi\sphinxhyphen{}definite variable in the expression, starting with 0.
\end{quote}
\end{quote}

\sphinxAtStartPar
\sphinxstylestrong{Example}
\end{quote}

\begin{sphinxVerbatim}[commandchars=\\\{\}]
\PYG{c+c1}{\PYGZsh{} Retrieve the symmetric matrix coefficient for the positive semi\PYGZhy{}definite variable at index 1 from expression expr}
\PYG{n}{mat} \PYG{o}{=} \PYG{n}{expr}\PYG{o}{.}\PYG{n}{getCoeff}\PYG{p}{(}\PYG{l+m+mi}{1}\PYG{p}{)}
\end{sphinxVerbatim}

\subsubsection{PsdExpr.getPsdVar()}
\label{\detokenize{pyapiref:psdexpr-getpsdvar}}\begin{quote}

\sphinxAtStartPar
\sphinxstylestrong{Synopsis}
\begin{quote}

\sphinxAtStartPar
\sphinxcode{\sphinxupquote{getPsdVar(idx)}}
\end{quote}

\sphinxAtStartPar
\sphinxstylestrong{Description}
\begin{quote}

\sphinxAtStartPar
Retrieve a positive semi\sphinxhyphen{}definite variable by its index from the expression. Return a {\hyperref[\detokenize{pyapiref:chappyapi-psdvar}]{\sphinxcrossref{\DUrole{std,std-ref}{PsdVar Class}}}} object.
\end{quote}

\sphinxAtStartPar
\sphinxstylestrong{Arguments}
\begin{quote}

\sphinxAtStartPar
\sphinxcode{\sphinxupquote{idx}}
\begin{quote}

\sphinxAtStartPar
Index of the positive semi\sphinxhyphen{}definite variable in the expression, starting with 0.
\end{quote}
\end{quote}

\sphinxAtStartPar
\sphinxstylestrong{Example}
\end{quote}

\begin{sphinxVerbatim}[commandchars=\\\{\}]
\PYG{c+c1}{\PYGZsh{} Retrieve the positive semi\PYGZhy{}definite variable at index 1 from expression expr}
\PYG{n}{x} \PYG{o}{=} \PYG{n}{expr}\PYG{o}{.}\PYG{n}{getPsdVar}\PYG{p}{(}\PYG{l+m+mi}{1}\PYG{p}{)}
\end{sphinxVerbatim}

\subsubsection{PsdExpr.getLinExpr()}
\label{\detokenize{pyapiref:psdexpr-getlinexpr}}\begin{quote}

\sphinxAtStartPar
\sphinxstylestrong{Synopsis}
\begin{quote}

\sphinxAtStartPar
\sphinxcode{\sphinxupquote{getLinExpr()}}
\end{quote}

\sphinxAtStartPar
\sphinxstylestrong{Description}
\begin{quote}

\sphinxAtStartPar
Retrieve the linear terms (if exist) from positive semi\sphinxhyphen{}definite expression. Return a {\hyperref[\detokenize{pyapiref:chappyapi-linexpr}]{\sphinxcrossref{\DUrole{std,std-ref}{LinExpr Class}}}} object.
\end{quote}

\sphinxAtStartPar
\sphinxstylestrong{Example}
\end{quote}

\begin{sphinxVerbatim}[commandchars=\\\{\}]
\PYG{c+c1}{\PYGZsh{} Retrieve the linear terms from a positive semi\PYGZhy{}definite expression expr}
\PYG{n}{linexpr} \PYG{o}{=} \PYG{n}{expr}\PYG{o}{.}\PYG{n}{getLinExpr}\PYG{p}{(}\PYG{p}{)}
\end{sphinxVerbatim}

\subsubsection{PsdExpr.getConstant()}
\label{\detokenize{pyapiref:psdexpr-getconstant}}\begin{quote}

\sphinxAtStartPar
\sphinxstylestrong{Synopsis}
\begin{quote}

\sphinxAtStartPar
\sphinxcode{\sphinxupquote{getConstant()}}
\end{quote}

\sphinxAtStartPar
\sphinxstylestrong{Description}
\begin{quote}

\sphinxAtStartPar
Retrieve the constant term from a positive semi\sphinxhyphen{}definite expression.
\end{quote}

\sphinxAtStartPar
\sphinxstylestrong{Example}
\end{quote}

\begin{sphinxVerbatim}[commandchars=\\\{\}]
\PYG{c+c1}{\PYGZsh{} Retrieve the constant term from expression expr}
\PYG{n}{constant} \PYG{o}{=} \PYG{n}{expr}\PYG{o}{.}\PYG{n}{getConstant}\PYG{p}{(}\PYG{p}{)}
\end{sphinxVerbatim}

\subsubsection{PsdExpr.getValue()}
\label{\detokenize{pyapiref:psdexpr-getvalue}}\begin{quote}

\sphinxAtStartPar
\sphinxstylestrong{Synopsis}
\begin{quote}

\sphinxAtStartPar
\sphinxcode{\sphinxupquote{getValue()}}
\end{quote}

\sphinxAtStartPar
\sphinxstylestrong{Description}
\begin{quote}

\sphinxAtStartPar
Retrieve the value of a positive semi\sphinxhyphen{}definite expression computed using the current solution.
\end{quote}

\sphinxAtStartPar
\sphinxstylestrong{Example}
\end{quote}

\begin{sphinxVerbatim}[commandchars=\\\{\}]
\PYG{c+c1}{\PYGZsh{} Retrieve the value of positive semi\PYGZhy{}definite expression expr for the current solution.}
\PYG{n}{val} \PYG{o}{=} \PYG{n}{expr}\PYG{o}{.}\PYG{n}{getValue}\PYG{p}{(}\PYG{p}{)}
\end{sphinxVerbatim}

\subsubsection{PsdExpr.getSize()}
\label{\detokenize{pyapiref:psdexpr-getsize}}\begin{quote}

\sphinxAtStartPar
\sphinxstylestrong{Synopsis}
\begin{quote}

\sphinxAtStartPar
\sphinxcode{\sphinxupquote{getSize()}}
\end{quote}

\sphinxAtStartPar
\sphinxstylestrong{Description}
\begin{quote}

\sphinxAtStartPar
Retrieve the number of terms in a positive semi\sphinxhyphen{}definite expression.
\end{quote}

\sphinxAtStartPar
\sphinxstylestrong{Example}
\end{quote}

\begin{sphinxVerbatim}[commandchars=\\\{\}]
\PYG{c+c1}{\PYGZsh{} Retrieve the number of terms in expression expr}
\PYG{n}{exprsize} \PYG{o}{=} \PYG{n}{expr}\PYG{o}{.}\PYG{n}{getSize}\PYG{p}{(}\PYG{p}{)}
\end{sphinxVerbatim}

\subsubsection{PsdExpr.setConstant()}
\label{\detokenize{pyapiref:psdexpr-setconstant}}\begin{quote}

\sphinxAtStartPar
\sphinxstylestrong{Synopsis}
\begin{quote}

\sphinxAtStartPar
\sphinxcode{\sphinxupquote{setConstant(newval)}}
\end{quote}

\sphinxAtStartPar
\sphinxstylestrong{Description}
\begin{quote}

\sphinxAtStartPar
Set the constant term of positive semi\sphinxhyphen{}definite expression.
\end{quote}

\sphinxAtStartPar
\sphinxstylestrong{Arguments}
\begin{quote}

\sphinxAtStartPar
\sphinxcode{\sphinxupquote{newval}}
\begin{quote}

\sphinxAtStartPar
Constant to set.
\end{quote}
\end{quote}

\sphinxAtStartPar
\sphinxstylestrong{Example}
\end{quote}

\begin{sphinxVerbatim}[commandchars=\\\{\}]
\PYG{c+c1}{\PYGZsh{} Set constant term of expression \PYGZsq{}expr\PYGZsq{} to 2.0}
\PYG{n}{expr}\PYG{o}{.}\PYG{n}{setConstant}\PYG{p}{(}\PYG{l+m+mf}{2.0}\PYG{p}{)}
\end{sphinxVerbatim}

\subsubsection{PsdExpr.addConstant()}
\label{\detokenize{pyapiref:psdexpr-addconstant}}\begin{quote}

\sphinxAtStartPar
\sphinxstylestrong{Synopsis}
\begin{quote}

\sphinxAtStartPar
\sphinxcode{\sphinxupquote{addConstant(newval)}}
\end{quote}

\sphinxAtStartPar
\sphinxstylestrong{Description}
\begin{quote}

\sphinxAtStartPar
Add a constant to a positive semi\sphinxhyphen{}definite expression.
\end{quote}

\sphinxAtStartPar
\sphinxstylestrong{Arguments}
\begin{quote}

\sphinxAtStartPar
\sphinxcode{\sphinxupquote{newval}}
\begin{quote}

\sphinxAtStartPar
Constant to add.
\end{quote}
\end{quote}

\sphinxAtStartPar
\sphinxstylestrong{Example}
\end{quote}

\begin{sphinxVerbatim}[commandchars=\\\{\}]
\PYG{c+c1}{\PYGZsh{} Add constant 2.0 to expression \PYGZsq{}expr\PYGZsq{}}
\PYG{n}{expr}\PYG{o}{.}\PYG{n}{addConstant}\PYG{p}{(}\PYG{l+m+mf}{2.0}\PYG{p}{)}
\end{sphinxVerbatim}

\subsubsection{PsdExpr.addTerm()}
\label{\detokenize{pyapiref:psdexpr-addterm}}\begin{quote}

\sphinxAtStartPar
\sphinxstylestrong{Synopsis}
\begin{quote}

\sphinxAtStartPar
\sphinxcode{\sphinxupquote{addTerm(var, mat)}}
\end{quote}

\sphinxAtStartPar
\sphinxstylestrong{Description}
\begin{quote}

\sphinxAtStartPar
Add a new term to current positive semi\sphinxhyphen{}definite expression.
\end{quote}

\sphinxAtStartPar
\sphinxstylestrong{Arguments}
\begin{quote}

\sphinxAtStartPar
\sphinxcode{\sphinxupquote{var}}
\begin{quote}

\sphinxAtStartPar
The positive semi\sphinxhyphen{}definite variable to add.
\end{quote}

\sphinxAtStartPar
\sphinxcode{\sphinxupquote{mat}}
\begin{quote}

\sphinxAtStartPar
The symmetric matrix coefficient for the positive semi\sphinxhyphen{}definite variable.
\end{quote}
\end{quote}

\sphinxAtStartPar
\sphinxstylestrong{Example}
\end{quote}

\begin{sphinxVerbatim}[commandchars=\\\{\}]
\PYG{c+c1}{\PYGZsh{} Add positive semi\PYGZhy{}definite term C1 * X to expression \PYGZsq{}expr\PYGZsq{}}
\PYG{n}{expr}\PYG{o}{.}\PYG{n}{addTerm}\PYG{p}{(}\PYG{n}{X}\PYG{p}{,} \PYG{n}{C1}\PYG{p}{)}
\end{sphinxVerbatim}

\subsubsection{PsdExpr.addTerms()}
\label{\detokenize{pyapiref:psdexpr-addterms}}\begin{quote}

\sphinxAtStartPar
\sphinxstylestrong{Synopsis}
\begin{quote}

\sphinxAtStartPar
\sphinxcode{\sphinxupquote{addTerms(vars, mats)}}
\end{quote}

\sphinxAtStartPar
\sphinxstylestrong{Description}
\begin{quote}

\sphinxAtStartPar
Add a single term or multiple positive semi\sphinxhyphen{}definite terms into a positive semi\sphinxhyphen{}definite expression.

\sphinxAtStartPar
If argument \sphinxcode{\sphinxupquote{vars}} is {\hyperref[\detokenize{pyapiref:chappyapi-psdvar}]{\sphinxcrossref{\DUrole{std,std-ref}{PsdVar Class}}}} object,
then argument \sphinxcode{\sphinxupquote{mats}} is {\hyperref[\detokenize{pyapiref:chappyapi-symmatrix}]{\sphinxcrossref{\DUrole{std,std-ref}{SymMatrix Class}}}} object;
If argument \sphinxcode{\sphinxupquote{vars}} is {\hyperref[\detokenize{pyapiref:chappyapi-psdvararray}]{\sphinxcrossref{\DUrole{std,std-ref}{PsdVarArray Class}}}} object or list,
then argument \sphinxcode{\sphinxupquote{mats}} is {\hyperref[\detokenize{pyapiref:chappyapi-symmatrixarray}]{\sphinxcrossref{\DUrole{std,std-ref}{SymMatrixArray Class}}}} or list;
\end{quote}

\sphinxAtStartPar
\sphinxstylestrong{Arguments}
\begin{quote}

\sphinxAtStartPar
\sphinxcode{\sphinxupquote{vars}}
\begin{quote}

\sphinxAtStartPar
The positive semi\sphinxhyphen{}definite variables to add.
\end{quote}

\sphinxAtStartPar
\sphinxcode{\sphinxupquote{mats}}
\begin{quote}

\sphinxAtStartPar
The symmetric matrices of the positive semi\sphinxhyphen{}definite terms.
\end{quote}
\end{quote}

\sphinxAtStartPar
\sphinxstylestrong{Example}
\end{quote}

\begin{sphinxVerbatim}[commandchars=\\\{\}]
\PYG{c+c1}{\PYGZsh{} Add terms C1 * X1 + C2 * X2 to expression \PYGZsq{}expr\PYGZsq{}}
\PYG{n}{expr}\PYG{o}{.}\PYG{n}{addTerms}\PYG{p}{(}\PYG{p}{[}\PYG{n}{X1}\PYG{p}{,} \PYG{n}{X2}\PYG{p}{]}\PYG{p}{,} \PYG{p}{[}\PYG{n}{C1}\PYG{p}{,} \PYG{n}{C2}\PYG{p}{]}\PYG{p}{)}
\end{sphinxVerbatim}

\subsubsection{PsdExpr.addLinExpr()}
\label{\detokenize{pyapiref:psdexpr-addlinexpr}}\begin{quote}

\sphinxAtStartPar
\sphinxstylestrong{Synopsis}
\begin{quote}

\sphinxAtStartPar
\sphinxcode{\sphinxupquote{addLinExpr(expr, mult=1.0)}}
\end{quote}

\sphinxAtStartPar
\sphinxstylestrong{Description}
\begin{quote}

\sphinxAtStartPar
Add new linear expression to the current positive semi\sphinxhyphen{}definite expression.
\end{quote}

\sphinxAtStartPar
\sphinxstylestrong{Arguments}
\begin{quote}

\sphinxAtStartPar
\sphinxcode{\sphinxupquote{expr}}
\begin{quote}

\sphinxAtStartPar
Linear expression or linear expression builder to add.
\end{quote}

\sphinxAtStartPar
\sphinxcode{\sphinxupquote{mult}}
\begin{quote}

\sphinxAtStartPar
Magnification coefficient for the added expression. Optional, 1.0 by default.
\end{quote}
\end{quote}

\sphinxAtStartPar
\sphinxstylestrong{Example}
\end{quote}

\begin{sphinxVerbatim}[commandchars=\\\{\}]
\PYG{c+c1}{\PYGZsh{} Add linear expression 2*x + 2*y to \PYGZsq{}expr\PYGZsq{}}
\PYG{n}{expr}\PYG{o}{.}\PYG{n}{addLinExpr}\PYG{p}{(}\PYG{n}{x} \PYG{o}{+} \PYG{n}{y}\PYG{p}{,} \PYG{l+m+mf}{2.0}\PYG{p}{)}
\end{sphinxVerbatim}

\subsubsection{PsdExpr.addPsdExpr()}
\label{\detokenize{pyapiref:psdexpr-addpsdexpr}}\begin{quote}

\sphinxAtStartPar
\sphinxstylestrong{Synopsis}
\begin{quote}

\sphinxAtStartPar
\sphinxcode{\sphinxupquote{addPsdExpr(expr, mult=1.0)}}
\end{quote}

\sphinxAtStartPar
\sphinxstylestrong{Description}
\begin{quote}

\sphinxAtStartPar
Add new positive semi\sphinxhyphen{}definite expression to the current one.
\end{quote}

\sphinxAtStartPar
\sphinxstylestrong{Arguments}
\begin{quote}

\sphinxAtStartPar
\sphinxcode{\sphinxupquote{expr}}
\begin{quote}

\sphinxAtStartPar
Positive semi\sphinxhyphen{}definite expression or positive semi\sphinxhyphen{}definite expression builder to add.
\end{quote}

\sphinxAtStartPar
\sphinxcode{\sphinxupquote{mult}}
\begin{quote}

\sphinxAtStartPar
Magnification coefficient for the added positive semi\sphinxhyphen{}definite expression. Optional, 1.0 by default.
\end{quote}
\end{quote}

\sphinxAtStartPar
\sphinxstylestrong{Example}
\end{quote}

\begin{sphinxVerbatim}[commandchars=\\\{\}]
\PYG{c+c1}{\PYGZsh{} Add positive semi\PYGZhy{}definite expression C * X to \PYGZsq{}expr\PYGZsq{}}
\PYG{n}{expr}\PYG{o}{.}\PYG{n}{addPsdExpr}\PYG{p}{(}\PYG{n}{C}\PYG{o}{*}\PYG{n}{X}\PYG{p}{)}
\end{sphinxVerbatim}

\subsubsection{PsdExpr.addMExpr()}
\label{\detokenize{pyapiref:psdexpr-addmexpr}}\begin{quote}

\sphinxAtStartPar
\sphinxstylestrong{Synopsis}
\begin{quote}

\sphinxAtStartPar
\sphinxcode{\sphinxupquote{addMExpr(expr, mult=1.0)}}
\end{quote}

\sphinxAtStartPar
\sphinxstylestrong{Description}
\begin{quote}

\sphinxAtStartPar
Add a new multi\sphinxhyphen{}dimensional expression to the current PSD expression.
\end{quote}

\sphinxAtStartPar
\sphinxstylestrong{Arguments}
\begin{quote}

\sphinxAtStartPar
\sphinxcode{\sphinxupquote{expr}}
\begin{quote}

\sphinxAtStartPar
The multi\sphinxhyphen{}dimensional array expression or expression builder object to add.
\end{quote}

\sphinxAtStartPar
\sphinxcode{\sphinxupquote{mult}}
\begin{quote}

\sphinxAtStartPar
The scaling factor for the multi\sphinxhyphen{}dimensional array expression. Optional, defaulting to 1.0.
\end{quote}
\end{quote}

\sphinxAtStartPar
\sphinxstylestrong{Example}
\end{quote}

\begin{sphinxVerbatim}[commandchars=\\\{\}]
\PYG{c+c1}{\PYGZsh{} Add a multi\PYGZhy{}dimensional linear expression: 2.0 * A @ x to the PSD expression expr}
\PYG{n}{expr}\PYG{o}{.}\PYG{n}{addMExpr}\PYG{p}{(}\PYG{n}{A} \PYG{o}{@} \PYG{n}{x}\PYG{p}{,} \PYG{l+m+mf}{2.0}\PYG{p}{)}
\end{sphinxVerbatim}

\subsubsection{PsdExpr.clone()}
\label{\detokenize{pyapiref:psdexpr-clone}}\begin{quote}

\sphinxAtStartPar
\sphinxstylestrong{Synopsis}
\begin{quote}

\sphinxAtStartPar
\sphinxcode{\sphinxupquote{clone()}}
\end{quote}

\sphinxAtStartPar
\sphinxstylestrong{Description}
\begin{quote}

\sphinxAtStartPar
Create a deep copy of the expression.
\end{quote}

\sphinxAtStartPar
\sphinxstylestrong{Example}
\end{quote}

\begin{sphinxVerbatim}[commandchars=\\\{\}]
\PYG{c+c1}{\PYGZsh{} Create a deep copy of expression expr}
\PYG{n}{exprcopy} \PYG{o}{=} \PYG{n}{expr}\PYG{o}{.}\PYG{n}{clone}\PYG{p}{(}\PYG{p}{)}
\end{sphinxVerbatim}

\subsubsection{PsdExpr.reserve()}
\label{\detokenize{pyapiref:psdexpr-reserve}}\begin{quote}

\sphinxAtStartPar
\sphinxstylestrong{Synopsis}
\begin{quote}

\sphinxAtStartPar
\sphinxcode{\sphinxupquote{reserve(n)}}
\end{quote}

\sphinxAtStartPar
\sphinxstylestrong{Description}
\begin{quote}

\sphinxAtStartPar
Pre\sphinxhyphen{}allocate space for positive semi\sphinxhyphen{}definite expression object.
\end{quote}

\sphinxAtStartPar
\sphinxstylestrong{Arguments}
\begin{quote}

\sphinxAtStartPar
\sphinxcode{\sphinxupquote{n}}
\begin{quote}

\sphinxAtStartPar
Number of terms to be allocated.
\end{quote}
\end{quote}

\sphinxAtStartPar
\sphinxstylestrong{Example}
\end{quote}

\begin{sphinxVerbatim}[commandchars=\\\{\}]
\PYG{c+c1}{\PYGZsh{} Allocate 100 terms for positive semi\PYGZhy{}definite expression \PYGZsq{}expr\PYGZsq{}}
\PYG{n}{expr}\PYG{o}{.}\PYG{n}{reserve}\PYG{p}{(}\PYG{l+m+mi}{100}\PYG{p}{)}
\end{sphinxVerbatim}

\subsubsection{PsdExpr.remove()}
\label{\detokenize{pyapiref:psdexpr-remove}}\begin{quote}

\sphinxAtStartPar
\sphinxstylestrong{Synopsis}
\begin{quote}

\sphinxAtStartPar
\sphinxcode{\sphinxupquote{remove(item)}}
\end{quote}

\sphinxAtStartPar
\sphinxstylestrong{Description}
\begin{quote}

\sphinxAtStartPar
Remove a term from a positive semi\sphinxhyphen{}definite expression.

\sphinxAtStartPar
If argument \sphinxcode{\sphinxupquote{item}} is constant, then remove the term stored at index i of the expression; otherwise argument \sphinxcode{\sphinxupquote{item}} is {\hyperref[\detokenize{pyapiref:chappyapi-psdvar}]{\sphinxcrossref{\DUrole{std,std-ref}{PsdVar Class}}}} object.
\end{quote}

\sphinxAtStartPar
\sphinxstylestrong{Arguments}
\begin{quote}

\sphinxAtStartPar
\sphinxcode{\sphinxupquote{item}}
\begin{quote}

\sphinxAtStartPar
Constant index or {\hyperref[\detokenize{pyapiref:chappyapi-psdvar}]{\sphinxcrossref{\DUrole{std,std-ref}{PsdVar Class}}}} variable of the term to be removed.
\end{quote}
\end{quote}

\sphinxAtStartPar
\sphinxstylestrong{Example}
\end{quote}

\begin{sphinxVerbatim}[commandchars=\\\{\}]
\PYG{c+c1}{\PYGZsh{} Remove the term whose index is 2 from positive semi\PYGZhy{}definite expression expr}
\PYG{n}{expr}\PYG{o}{.}\PYG{n}{remove}\PYG{p}{(}\PYG{l+m+mi}{2}\PYG{p}{)}
\PYG{c+c1}{\PYGZsh{} Remove the terms one of which variable is x from positive semi\PYGZhy{}definite expression expr}
\PYG{n}{expr}\PYG{o}{.}\PYG{n}{remove}\PYG{p}{(}\PYG{n}{x}\PYG{p}{)}
\end{sphinxVerbatim}

\subsection{MPsdExpr Class}
\label{\detokenize{pyapiref:mpsdexpr-class}}\label{\detokenize{pyapiref:chappyapi-mpsdexpr}}
\sphinxAtStartPar
The \sphinxtitleref{MPsdExpr} class in COPT is used for operations that combine PSD variables into multi\sphinxhyphen{}dimensional PSD expressions. The following methods are provided:

\subsubsection{MPsdExpr.addTerm()}
\label{\detokenize{pyapiref:mpsdexpr-addterm}}\begin{quote}

\sphinxAtStartPar
\sphinxstylestrong{Synopsis}
\begin{quote}

\sphinxAtStartPar
\sphinxcode{\sphinxupquote{addTerm(var, mat)}}
\end{quote}

\sphinxAtStartPar
\sphinxstylestrong{Description}
\begin{quote}

\sphinxAtStartPar
Adds a new PSD term to the current multi\sphinxhyphen{}dimensional PSD expression.
\end{quote}

\sphinxAtStartPar
\sphinxstylestrong{Arguments}
\begin{quote}

\sphinxAtStartPar
\sphinxcode{\sphinxupquote{var}}
\begin{quote}

\sphinxAtStartPar
The semidefinite variable in the term to be added.
\end{quote}

\sphinxAtStartPar
\sphinxcode{\sphinxupquote{mat}}
\begin{quote}

\sphinxAtStartPar
The symmetric matrix in the term to be added.
\end{quote}
\end{quote}

\sphinxAtStartPar
\sphinxstylestrong{Example}
\end{quote}

\begin{sphinxVerbatim}[commandchars=\\\{\}]
\PYG{c+c1}{\PYGZsh{} Add the PSD term C1 * X to mexpr}
\PYG{n}{mexpr}\PYG{o}{.}\PYG{n}{addTerm}\PYG{p}{(}\PYG{n}{X}\PYG{p}{,} \PYG{n}{C1}\PYG{p}{)}
\end{sphinxVerbatim}

\subsubsection{MPsdExpr.addTerms()}
\label{\detokenize{pyapiref:mpsdexpr-addterms}}\begin{quote}

\sphinxAtStartPar
\sphinxstylestrong{Synopsis}
\begin{quote}

\sphinxAtStartPar
\sphinxcode{\sphinxupquote{addTerms(vars, coeffs)}}
\end{quote}

\sphinxAtStartPar
\sphinxstylestrong{Description}
\begin{quote}

\sphinxAtStartPar
Adds new terms to the multi\sphinxhyphen{}dimensional PSD expression object.
\end{quote}

\sphinxAtStartPar
\sphinxstylestrong{Arguments}
\begin{quote}

\sphinxAtStartPar
\sphinxcode{\sphinxupquote{vars}}
\begin{quote}

\sphinxAtStartPar
Multi\sphinxhyphen{}dimensional array variable objects. Possible value could be \sphinxtitleref{MVar}.
\end{quote}

\sphinxAtStartPar
\sphinxcode{\sphinxupquote{coeffs}}
\begin{quote}

\sphinxAtStartPar
Coefficient matrices for the terms. Possible values could be a floating number or \sphinxtitleref{NdArray}.
\end{quote}
\end{quote}

\sphinxAtStartPar
\sphinxstylestrong{Example}
\end{quote}

\begin{sphinxVerbatim}[commandchars=\\\{\}]
\PYG{c+c1}{\PYGZsh{} Add terms: mA @ mX to mpsdexpr}
\PYG{n}{mX} \PYG{o}{=} \PYG{n}{model}\PYG{o}{.}\PYG{n}{addMVar}\PYG{p}{(}\PYG{p}{(}\PYG{l+m+mi}{3}\PYG{p}{,}\PYG{l+m+mi}{3}\PYG{p}{)}\PYG{p}{,} \PYG{n}{nameprefix}\PYG{o}{=}\PYG{l+s+s2}{\PYGZdq{}}\PYG{l+s+s2}{M\PYGZus{}X}\PYG{l+s+s2}{\PYGZdq{}}\PYG{p}{)}
\PYG{n}{mA} \PYG{o}{=} \PYG{n}{cp}\PYG{o}{.}\PYG{n}{NdArray}\PYG{p}{(}\PYG{n}{np}\PYG{o}{.}\PYG{n}{ones}\PYG{p}{(}\PYG{n}{shape}\PYG{o}{=}\PYG{p}{(}\PYG{l+m+mi}{3}\PYG{p}{,}\PYG{l+m+mi}{3}\PYG{p}{)}\PYG{p}{)}\PYG{p}{)}
\PYG{n}{mpsdexpr}\PYG{o}{.}\PYG{n}{addTerms}\PYG{p}{(}\PYG{n}{mX}\PYG{p}{,} \PYG{n}{mA}\PYG{p}{)}
\end{sphinxVerbatim}

\subsubsection{MPsdExpr.addLinExpr()}
\label{\detokenize{pyapiref:mpsdexpr-addlinexpr}}\begin{quote}

\sphinxAtStartPar
\sphinxstylestrong{Synopsis}
\begin{quote}

\sphinxAtStartPar
\sphinxcode{\sphinxupquote{addLinExpr(expr, mult=1.0)}}
\end{quote}

\sphinxAtStartPar
\sphinxstylestrong{Description}
\begin{quote}

\sphinxAtStartPar
Adds a new linear expression to the current multi\sphinxhyphen{}dimensional PSD expression.
\end{quote}

\sphinxAtStartPar
\sphinxstylestrong{Arguments}
\begin{quote}

\sphinxAtStartPar
\sphinxcode{\sphinxupquote{expr}}
\begin{quote}

\sphinxAtStartPar
The linear expression or expression builder object to add.
Possible values could be \sphinxtitleref{LinExpr} or \sphinxtitleref{ExprBuilder}.
\end{quote}

\sphinxAtStartPar
\sphinxcode{\sphinxupquote{mult}}
\begin{quote}

\sphinxAtStartPar
The scaling factor for the linear expression.
Optional, defaulting to 1.0.
\end{quote}
\end{quote}

\sphinxAtStartPar
\sphinxstylestrong{Example}
\end{quote}

\begin{sphinxVerbatim}[commandchars=\\\{\}]
\PYG{c+c1}{\PYGZsh{} Add a linear expression: 2*x + 2*y to mexpr}
\PYG{n}{mexpr}\PYG{o}{.}\PYG{n}{addLinExpr}\PYG{p}{(}\PYG{n}{x} \PYG{o}{+} \PYG{n}{y}\PYG{p}{,} \PYG{l+m+mf}{2.0}\PYG{p}{)}
\end{sphinxVerbatim}

\subsubsection{MPsdExpr.addPsdExpr()}
\label{\detokenize{pyapiref:mpsdexpr-addpsdexpr}}\begin{quote}

\sphinxAtStartPar
\sphinxstylestrong{Synopsis}
\begin{quote}

\sphinxAtStartPar
\sphinxcode{\sphinxupquote{addPsdExpr(expr, mult=1.0)}}
\end{quote}

\sphinxAtStartPar
\sphinxstylestrong{Description}
\begin{quote}

\sphinxAtStartPar
Adds a new semidefinite expression to the current multi\sphinxhyphen{}dimensional PSD expression.
\end{quote}

\sphinxAtStartPar
\sphinxstylestrong{Arguments}
\begin{quote}

\sphinxAtStartPar
\sphinxcode{\sphinxupquote{expr}}
\begin{quote}

\sphinxAtStartPar
The semidefinite expression to add.
\end{quote}

\sphinxAtStartPar
\sphinxcode{\sphinxupquote{mult}}
\begin{quote}

\sphinxAtStartPar
The scaling factor for the PSD expression.
Optional, defaulting to 1.0.
\end{quote}
\end{quote}

\sphinxAtStartPar
\sphinxstylestrong{Example}
\end{quote}

\begin{sphinxVerbatim}[commandchars=\\\{\}]
\PYG{c+c1}{\PYGZsh{} Add a PSD expression: C * X to mexpr}
\PYG{n}{mexpr}\PYG{o}{.}\PYG{n}{addPsdExpr}\PYG{p}{(}\PYG{n}{C} \PYG{o}{*} \PYG{n}{X}\PYG{p}{)}
\end{sphinxVerbatim}

\subsubsection{MPsdExpr.addMExpr()}
\label{\detokenize{pyapiref:mpsdexpr-addmexpr}}\begin{quote}

\sphinxAtStartPar
\sphinxstylestrong{Synopsis}
\begin{quote}

\sphinxAtStartPar
\sphinxcode{\sphinxupquote{addMExpr(expr, mult=1.0)}}
\end{quote}

\sphinxAtStartPar
\sphinxstylestrong{Description}
\begin{quote}

\sphinxAtStartPar
Add a new multi\sphinxhyphen{}dimensional expression to the current multi\sphinxhyphen{}dimensional PSD expression.
\end{quote}

\sphinxAtStartPar
\sphinxstylestrong{Arguments}
\begin{quote}

\sphinxAtStartPar
\sphinxcode{\sphinxupquote{expr}}
\begin{quote}

\sphinxAtStartPar
The multi\sphinxhyphen{}dimensional expression or expression builder object to add.
\end{quote}

\sphinxAtStartPar
\sphinxcode{\sphinxupquote{mult}}
\begin{quote}

\sphinxAtStartPar
The scaling factor for the multi\sphinxhyphen{}dimensional expression.
Optional, defaulting to 1.0.
\end{quote}
\end{quote}

\sphinxAtStartPar
\sphinxstylestrong{Example}
\end{quote}

\begin{sphinxVerbatim}[commandchars=\\\{\}]
\PYG{c+c1}{\PYGZsh{} Add a multi\PYGZhy{}dimensional linear expression: 2.0 * A @ x to mexpr}
\PYG{n}{mexpr}\PYG{o}{.}\PYG{n}{addMExpr}\PYG{p}{(}\PYG{n}{A} \PYG{o}{@} \PYG{n}{x}\PYG{p}{,} \PYG{l+m+mf}{2.0}\PYG{p}{)}
\end{sphinxVerbatim}

\subsubsection{MPsdExpr.addMLinExpr()}
\label{\detokenize{pyapiref:mpsdexpr-addmlinexpr}}\begin{quote}

\sphinxAtStartPar
\sphinxstylestrong{Synopsis}
\begin{quote}

\sphinxAtStartPar
\sphinxcode{\sphinxupquote{addMLinExpr(exprs, mult=1.0)}}
\end{quote}

\sphinxAtStartPar
\sphinxstylestrong{Description}
\begin{quote}

\sphinxAtStartPar
Add corresponding linear expressions to each PSD expression in the \sphinxcode{\sphinxupquote{MPsdExpr}} object.
\end{quote}

\sphinxAtStartPar
\sphinxstylestrong{Arguments}
\begin{quote}

\sphinxAtStartPar
\sphinxcode{\sphinxupquote{exprs}}
\begin{quote}

\sphinxAtStartPar
The multi\sphinxhyphen{}dimensional linear expressions or expression builder objects to add.
\end{quote}

\sphinxAtStartPar
\sphinxcode{\sphinxupquote{mult}}
\begin{quote}

\sphinxAtStartPar
The same scaling factor for all linear expressions to add.
Optional, defaulting to 1.0.
\end{quote}
\end{quote}

\sphinxAtStartPar
\sphinxstylestrong{Example}
\end{quote}

\begin{sphinxVerbatim}[commandchars=\\\{\}]
\PYG{c+c1}{\PYGZsh{} Add a linear expression A @ x to each PSD expression in mexpr}
\PYG{n}{mexpr}\PYG{o}{.}\PYG{n}{addMLinExpr}\PYG{p}{(}\PYG{n}{A} \PYG{o}{@} \PYG{n}{x}\PYG{p}{)}
\end{sphinxVerbatim}

\subsubsection{MPsdExpr.addMPsdExpr()}
\label{\detokenize{pyapiref:mpsdexpr-addmpsdexpr}}\begin{quote}

\sphinxAtStartPar
\sphinxstylestrong{Synopsis}
\begin{quote}

\sphinxAtStartPar
\sphinxcode{\sphinxupquote{addMPsdExpr(exprs, mult=1.0)}}
\end{quote}

\sphinxAtStartPar
\sphinxstylestrong{Description}
\begin{quote}

\sphinxAtStartPar
Adds corresponding new PSD expressions to each PSD expression in the \sphinxcode{\sphinxupquote{MPsdExpr}} object.
\end{quote}

\sphinxAtStartPar
\sphinxstylestrong{Arguments}
\begin{quote}

\sphinxAtStartPar
\sphinxcode{\sphinxupquote{exprs}}
\begin{quote}

\sphinxAtStartPar
The new semidefinite expressions to add.
\end{quote}

\sphinxAtStartPar
\sphinxcode{\sphinxupquote{mult}}
\begin{quote}

\sphinxAtStartPar
The same scaling factor for all PSD expressions to add.
Optional, defaulting to 1.0.
\end{quote}
\end{quote}

\sphinxAtStartPar
\sphinxstylestrong{Example}
\end{quote}

\begin{sphinxVerbatim}[commandchars=\\\{\}]
\PYG{c+c1}{\PYGZsh{} Add a PSD expression: C * X to each PSD expression in mexpr}
\PYG{n}{mexpr}\PYG{o}{.}\PYG{n}{addMPsdExpr}\PYG{p}{(}\PYG{n}{C} \PYG{o}{*} \PYG{n}{X}\PYG{p}{)}
\end{sphinxVerbatim}

\subsubsection{MPsdExpr.item()}
\label{\detokenize{pyapiref:mpsdexpr-item}}\begin{quote}

\sphinxAtStartPar
\sphinxstylestrong{Synopsis}
\begin{quote}

\sphinxAtStartPar
\sphinxcode{\sphinxupquote{item()}}
\end{quote}

\sphinxAtStartPar
\sphinxstylestrong{Description}
\begin{quote}

\sphinxAtStartPar
Retrieves the \sphinxtitleref{PsdExpr} in a 0\sphinxhyphen{}dimensional semidefinite expression.
If the \sphinxtitleref{MPsdExpr} object is not 0\sphinxhyphen{}dimensional, raises a \sphinxtitleref{ValueError}.
\end{quote}

\sphinxAtStartPar
\sphinxstylestrong{Return Value}
\begin{quote}

\sphinxAtStartPar
Returns a \sphinxtitleref{PsdExpr} object.
\end{quote}

\sphinxAtStartPar
\sphinxstylestrong{Example}
\end{quote}

\begin{sphinxVerbatim}[commandchars=\\\{\}]
\PYG{n}{barX} \PYG{o}{=} \PYG{n}{model}\PYG{o}{.}\PYG{n}{addPsdVars}\PYG{p}{(}\PYG{l+m+mi}{3}\PYG{p}{,} \PYG{l+s+s2}{\PYGZdq{}}\PYG{l+s+s2}{BAR\PYGZus{}X}\PYG{l+s+s2}{\PYGZdq{}}\PYG{p}{)}
\PYG{n}{mpsdexpr} \PYG{o}{=} \PYG{n}{barX}\PYG{p}{[}\PYG{p}{:}\PYG{o}{\PYGZhy{}}\PYG{l+m+mi}{1}\PYG{p}{,} \PYG{p}{:}\PYG{o}{\PYGZhy{}}\PYG{l+m+mi}{1}\PYG{p}{]}
\PYG{n}{psdexpr1} \PYG{o}{=} \PYG{n}{mpsdexpr}\PYG{p}{[}\PYG{l+m+mi}{0}\PYG{p}{,}\PYG{l+m+mi}{0}\PYG{p}{]}\PYG{o}{.}\PYG{n}{item}\PYG{p}{(}\PYG{p}{)}
\PYG{n}{psdexpr2} \PYG{o}{=} \PYG{n}{mpsdexpr}\PYG{o}{.}\PYG{n}{sum}\PYG{p}{(}\PYG{p}{)}\PYG{o}{.}\PYG{n}{item}\PYG{p}{(}\PYG{p}{)}
\end{sphinxVerbatim}

\subsubsection{MPsdExpr.sum()}
\label{\detokenize{pyapiref:mpsdexpr-sum}}\begin{quote}

\sphinxAtStartPar
\sphinxstylestrong{Synopsis}
\begin{quote}

\sphinxAtStartPar
\sphinxcode{\sphinxupquote{sum(axis=None)}}
\end{quote}

\sphinxAtStartPar
\sphinxstylestrong{Description}
\begin{quote}

\sphinxAtStartPar
Computes the sum of the semidefinite terms along the specified axis
in the \sphinxtitleref{MPsdExpr} object.
\end{quote}

\sphinxAtStartPar
\sphinxstylestrong{Arguments}
\begin{quote}

\sphinxAtStartPar
\sphinxcode{\sphinxupquote{axis}}
\begin{quote}

\sphinxAtStartPar
Optional. Defaults to \sphinxtitleref{None}, which sums over all variables.
Otherwise, sums along the specified axis.
\end{quote}
\end{quote}

\sphinxAtStartPar
\sphinxstylestrong{Return Value}
\begin{quote}

\sphinxAtStartPar
Returns an {\hyperref[\detokenize{pyapiref:chappyapi-mpsdexpr}]{\sphinxcrossref{\DUrole{std,std-ref}{MPsdExpr Class}}}} object representing the sum of the corresponding
multi\sphinxhyphen{}dimensional PSD expressions.
\end{quote}
\end{quote}

\subsubsection{MPsdExpr.clear()}
\label{\detokenize{pyapiref:mpsdexpr-clear}}\begin{quote}

\sphinxAtStartPar
\sphinxstylestrong{Synopsis}
\begin{quote}

\sphinxAtStartPar
\sphinxcode{\sphinxupquote{clear()}}
\end{quote}

\sphinxAtStartPar
\sphinxstylestrong{Description}
\begin{quote}

\sphinxAtStartPar
Resets every element of the {\hyperref[\detokenize{pyapiref:chappyapi-mpsdexpr}]{\sphinxcrossref{\DUrole{std,std-ref}{MPsdExpr Class}}}} object to 0.0.
\end{quote}

\sphinxAtStartPar
\sphinxstylestrong{Example}
\end{quote}

\begin{sphinxVerbatim}[commandchars=\\\{\}]
\PYG{n}{barX} \PYG{o}{=} \PYG{n}{model}\PYG{o}{.}\PYG{n}{addPsdVars}\PYG{p}{(}\PYG{l+m+mi}{3}\PYG{p}{,} \PYG{l+s+s2}{\PYGZdq{}}\PYG{l+s+s2}{BAR\PYGZus{}X}\PYG{l+s+s2}{\PYGZdq{}}\PYG{p}{)}
\PYG{n}{mpsdexpr} \PYG{o}{=} \PYG{n}{barX}\PYG{p}{[}\PYG{p}{:}\PYG{o}{\PYGZhy{}}\PYG{l+m+mi}{1}\PYG{p}{,} \PYG{p}{:}\PYG{o}{\PYGZhy{}}\PYG{l+m+mi}{1}\PYG{p}{]}
\PYG{n}{mpsdexpr}\PYG{o}{.}\PYG{n}{clear}\PYG{p}{(}\PYG{p}{)}
\end{sphinxVerbatim}

\subsubsection{MPsdExpr.clone()}
\label{\detokenize{pyapiref:mpsdexpr-clone}}\begin{quote}

\sphinxAtStartPar
\sphinxstylestrong{Synopsis}
\begin{quote}

\sphinxAtStartPar
\sphinxcode{\sphinxupquote{clone()}}
\end{quote}

\sphinxAtStartPar
\sphinxstylestrong{Description}
\begin{quote}

\sphinxAtStartPar
Creates a deep copy of an {\hyperref[\detokenize{pyapiref:chappyapi-mpsdexpr}]{\sphinxcrossref{\DUrole{std,std-ref}{MPsdExpr Class}}}} object.
\end{quote}

\sphinxAtStartPar
\sphinxstylestrong{Return Value}
\begin{quote}

\sphinxAtStartPar
Returns a new \sphinxtitleref{MPsdExpr} object.
\end{quote}

\sphinxAtStartPar
\sphinxstylestrong{Example}
\end{quote}

\begin{sphinxVerbatim}[commandchars=\\\{\}]
\PYG{n}{barX} \PYG{o}{=} \PYG{n}{model}\PYG{o}{.}\PYG{n}{addPsdVars}\PYG{p}{(}\PYG{l+m+mi}{3}\PYG{p}{,} \PYG{l+s+s2}{\PYGZdq{}}\PYG{l+s+s2}{BAR\PYGZus{}X}\PYG{l+s+s2}{\PYGZdq{}}\PYG{p}{)}
\PYG{n}{mpsdexpr} \PYG{o}{=} \PYG{n}{barX}\PYG{p}{[}\PYG{p}{:}\PYG{o}{\PYGZhy{}}\PYG{l+m+mi}{1}\PYG{p}{,} \PYG{p}{:}\PYG{o}{\PYGZhy{}}\PYG{l+m+mi}{1}\PYG{p}{]}
\PYG{n}{newmpsdexpr} \PYG{o}{=} \PYG{n}{mpsdexpr}\PYG{o}{.}\PYG{n}{clone}\PYG{p}{(}\PYG{p}{)}
\end{sphinxVerbatim}

\subsection{LmiExpr Class}
\label{\detokenize{pyapiref:lmiexpr-class}}\label{\detokenize{pyapiref:chappyapi-lmiexpr}}
\sphinxAtStartPar
LmiExpr object contains operations related to variables for building LMI constraints,
and provides the following methods:

\subsubsection{LmiExpr()}
\label{\detokenize{pyapiref:lmiexpr}}\begin{quote}

\sphinxAtStartPar
\sphinxstylestrong{Synopsis}
\begin{quote}

\sphinxAtStartPar
\sphinxcode{\sphinxupquote{LmiExpr(arg1=None, arg2=None)}}
\end{quote}

\sphinxAtStartPar
\sphinxstylestrong{Description}
\begin{quote}

\sphinxAtStartPar
Create a {\hyperref[\detokenize{pyapiref:chappyapi-lmiexpr}]{\sphinxcrossref{\DUrole{std,std-ref}{LmiExpr Class}}}} object.
\end{quote}

\sphinxAtStartPar
\sphinxstylestrong{Arguments}
\begin{quote}

\sphinxAtStartPar
The default value of \sphinxcode{\sphinxupquote{arg1}} is \sphinxcode{\sphinxupquote{None}}, and the Possible values:
{\hyperref[\detokenize{pyapiref:chappyapi-var}]{\sphinxcrossref{\DUrole{std,std-ref}{Var Class}}}} object, or {\hyperref[\detokenize{pyapiref:chappyapi-symmatrix}]{\sphinxcrossref{\DUrole{std,std-ref}{SymMatrix Class}}}} object.

\sphinxAtStartPar
If the argument \sphinxcode{\sphinxupquote{arg1}} is a {\hyperref[\detokenize{pyapiref:chappyapi-var}]{\sphinxcrossref{\DUrole{std,std-ref}{Var Class}}}} object, then the argument \sphinxcode{\sphinxupquote{arg2}} is a {\hyperref[\detokenize{pyapiref:chappyapi-symmatrix}]{\sphinxcrossref{\DUrole{std,std-ref}{SymMatrix Class}}}} object.
\end{quote}
\end{quote}

\subsubsection{LmiExpr.setCoeff()}
\label{\detokenize{pyapiref:lmiexpr-setcoeff}}\begin{quote}

\sphinxAtStartPar
\sphinxstylestrong{Synopsis}
\begin{quote}

\sphinxAtStartPar
\sphinxcode{\sphinxupquote{setCoeff(idx, mat)}}
\end{quote}

\sphinxAtStartPar
\sphinxstylestrong{Description}
\begin{quote}

\sphinxAtStartPar
Set the coefficient matrix for the entry corresponding to the specified index \sphinxcode{\sphinxupquote{idx}} in the LMI expression.
\end{quote}

\sphinxAtStartPar
\sphinxstylestrong{Arguments}
\begin{quote}

\sphinxAtStartPar
\sphinxcode{\sphinxupquote{idx}}
\begin{quote}

\sphinxAtStartPar
The specified the index value. Starts with 0.
\end{quote}

\sphinxAtStartPar
\sphinxcode{\sphinxupquote{mat}}
\begin{quote}

\sphinxAtStartPar
The new coefficient symmetric matrix of the variable to be set, which must be a {\hyperref[\detokenize{pyapiref:chappyapi-symmatrix}]{\sphinxcrossref{\DUrole{std,std-ref}{SymMatrix Class}}}} class object.
\end{quote}
\end{quote}

\sphinxAtStartPar
\sphinxstylestrong{Example}
\end{quote}

\begin{sphinxVerbatim}[commandchars=\\\{\}]
\PYG{c+c1}{\PYGZsh{} Set the coefficient of the 0\PYGZhy{}th term of the LMI expression expr to the symmetric matrix mat}
\PYG{n}{expr}\PYG{o}{.}\PYG{n}{setCoeff}\PYG{p}{(}\PYG{l+m+mi}{0}\PYG{p}{,} \PYG{n}{mat}\PYG{p}{)}
\end{sphinxVerbatim}

\subsubsection{LmiExpr.getCoeff()}
\label{\detokenize{pyapiref:lmiexpr-getcoeff}}\begin{quote}

\sphinxAtStartPar
\sphinxstylestrong{Synopsis}
\begin{quote}

\sphinxAtStartPar
\sphinxcode{\sphinxupquote{getCoeff(idx)}}
\end{quote}

\sphinxAtStartPar
\sphinxstylestrong{Description}
\begin{quote}

\sphinxAtStartPar
Get the coefficient matrix for the entry corresponding to the specified index \sphinxcode{\sphinxupquote{idx}} in the LMI expression.
\end{quote}

\sphinxAtStartPar
\sphinxstylestrong{Arguments}
\begin{quote}

\sphinxAtStartPar
\sphinxcode{\sphinxupquote{idx}}
\begin{quote}

\sphinxAtStartPar
The specified the index value. Starts with 0.
\end{quote}
\end{quote}

\sphinxAtStartPar
\sphinxstylestrong{Example}
\end{quote}

\begin{sphinxVerbatim}[commandchars=\\\{\}]
\PYG{c+c1}{\PYGZsh{} Get the symmetric matrix coefficient of the 1st term of the LMI expression expr}
\PYG{n}{mat} \PYG{o}{=} \PYG{n}{expr}\PYG{o}{.}\PYG{n}{getCoeff}\PYG{p}{(}\PYG{l+m+mi}{1}\PYG{p}{)}
\end{sphinxVerbatim}

\subsubsection{LmiExpr.getVar()}
\label{\detokenize{pyapiref:lmiexpr-getvar}}\begin{quote}

\sphinxAtStartPar
\sphinxstylestrong{Synopsis}
\begin{quote}

\sphinxAtStartPar
\sphinxcode{\sphinxupquote{getVar(idx)}}
\end{quote}

\sphinxAtStartPar
\sphinxstylestrong{Description}
\begin{quote}

\sphinxAtStartPar
Get the variable in the entry corresponding to the specified index \sphinxcode{\sphinxupquote{idx}} in the LMI expression.
\end{quote}

\sphinxAtStartPar
\sphinxstylestrong{Arguments}
\begin{quote}

\sphinxAtStartPar
\sphinxcode{\sphinxupquote{idx}}
\begin{quote}

\sphinxAtStartPar
The specified the index value. Starts with 0.
\end{quote}
\end{quote}

\sphinxAtStartPar
\sphinxstylestrong{Example}
\end{quote}

\begin{sphinxVerbatim}[commandchars=\\\{\}]
\PYG{c+c1}{\PYGZsh{} Get the variable of the 1st item of the LMI expression expr}
\PYG{n}{mat} \PYG{o}{=} \PYG{n}{expr}\PYG{o}{.}\PYG{n}{getVar}\PYG{p}{(}\PYG{l+m+mi}{1}\PYG{p}{)}
\end{sphinxVerbatim}

\subsubsection{LmiExpr.getConstant()}
\label{\detokenize{pyapiref:lmiexpr-getconstant}}\begin{quote}

\sphinxAtStartPar
\sphinxstylestrong{Synopsis}
\begin{quote}

\sphinxAtStartPar
\sphinxcode{\sphinxupquote{getConstant()}}
\end{quote}

\sphinxAtStartPar
\sphinxstylestrong{Description}
\begin{quote}

\sphinxAtStartPar
Get the constant\sphinxhyphen{}term symmetric matrix in the LMI expression.
\end{quote}

\sphinxAtStartPar
\sphinxstylestrong{Example}
\end{quote}

\begin{sphinxVerbatim}[commandchars=\\\{\}]
\PYG{c+c1}{\PYGZsh{} Get the constant term of the LMI expression expr}
\PYG{n}{constant} \PYG{o}{=} \PYG{n}{expr}\PYG{o}{.}\PYG{n}{getConstant}\PYG{p}{(}\PYG{p}{)}
\end{sphinxVerbatim}

\subsubsection{LmiExpr.getSize()}
\label{\detokenize{pyapiref:lmiexpr-getsize}}\begin{quote}

\sphinxAtStartPar
\sphinxstylestrong{Synopsis}
\begin{quote}

\sphinxAtStartPar
\sphinxcode{\sphinxupquote{getSize()}}
\end{quote}

\sphinxAtStartPar
\sphinxstylestrong{Description}
\begin{quote}

\sphinxAtStartPar
Retrieve the number of terms in the LMI expression.
\end{quote}

\sphinxAtStartPar
\sphinxstylestrong{Example}
\end{quote}

\begin{sphinxVerbatim}[commandchars=\\\{\}]
\PYG{c+c1}{\PYGZsh{} Retrieve the number of terms in the LMI expression expr.}
\PYG{n}{val} \PYG{o}{=} \PYG{n}{expr}\PYG{o}{.}\PYG{n}{getSize}\PYG{p}{(}\PYG{p}{)}
\end{sphinxVerbatim}

\subsubsection{LmiExpr.setConstant()}
\label{\detokenize{pyapiref:lmiexpr-setconstant}}\begin{quote}

\sphinxAtStartPar
\sphinxstylestrong{Synopsis}
\begin{quote}

\sphinxAtStartPar
\sphinxcode{\sphinxupquote{setConstant(mat)}}
\end{quote}

\sphinxAtStartPar
\sphinxstylestrong{Description}
\begin{quote}

\sphinxAtStartPar
Set the constant\sphinxhyphen{}term symmetric matrix of the LMI expression.
\end{quote}

\sphinxAtStartPar
\sphinxstylestrong{Arguments}
\begin{quote}

\sphinxAtStartPar
\sphinxcode{\sphinxupquote{mat}}
\begin{quote}

\sphinxAtStartPar
The symmetric matrix corresponding to the constant item, must be a {\hyperref[\detokenize{pyapiref:chappyapi-symmatrix}]{\sphinxcrossref{\DUrole{std,std-ref}{SymMatrix Class}}}}  object.
\end{quote}
\end{quote}

\sphinxAtStartPar
\sphinxstylestrong{Example}
\end{quote}

\begin{sphinxVerbatim}[commandchars=\\\{\}]
\PYG{c+c1}{\PYGZsh{} Sets the constant term of the LMI expression expr to the symmetric matrix D1}
\PYG{n}{expr}\PYG{o}{.}\PYG{n}{setConstant}\PYG{p}{(}\PYG{n}{D1}\PYG{p}{)}
\end{sphinxVerbatim}

\subsubsection{LmiExpr.addConstant()}
\label{\detokenize{pyapiref:lmiexpr-addconstant}}\begin{quote}

\sphinxAtStartPar
\sphinxstylestrong{Synopsis}
\begin{quote}

\sphinxAtStartPar
\sphinxcode{\sphinxupquote{addConstant(mat)}}
\end{quote}

\sphinxAtStartPar
\sphinxstylestrong{Description}
\begin{quote}

\sphinxAtStartPar
Add the symmetric matrix to constant term of the LMI expression.
\end{quote}

\sphinxAtStartPar
\sphinxstylestrong{Arguments}
\begin{quote}

\sphinxAtStartPar
\sphinxcode{\sphinxupquote{mat}}
\begin{quote}

\sphinxAtStartPar
Matrix expression object added to constant term.
\end{quote}
\end{quote}

\sphinxAtStartPar
\sphinxstylestrong{Example}
\end{quote}

\begin{sphinxVerbatim}[commandchars=\\\{\}]
\PYG{c+c1}{\PYGZsh{} Add to the constant term of the LMI expression expr with symmetric matrix D2}
\PYG{n}{expr}\PYG{o}{.}\PYG{n}{addConstant}\PYG{p}{(}\PYG{n}{D2}\PYG{p}{)}
\end{sphinxVerbatim}

\subsubsection{LmiExpr.addTerm()}
\label{\detokenize{pyapiref:lmiexpr-addterm}}\begin{quote}

\sphinxAtStartPar
\sphinxstylestrong{Synopsis}
\begin{quote}

\sphinxAtStartPar
\sphinxcode{\sphinxupquote{addTerm(var, mat)}}
\end{quote}

\sphinxAtStartPar
\sphinxstylestrong{Description}
\begin{quote}

\sphinxAtStartPar
Add a new term to current LMI expression.
\end{quote}

\sphinxAtStartPar
\sphinxstylestrong{Arguments}
\begin{quote}

\sphinxAtStartPar
\sphinxcode{\sphinxupquote{var}}
\begin{quote}

\sphinxAtStartPar
The variable in the new item.
\end{quote}

\sphinxAtStartPar
\sphinxcode{\sphinxupquote{mat}}
\begin{quote}

\sphinxAtStartPar
The symmetric matrix as variable coefficients in the new term.
\end{quote}
\end{quote}

\sphinxAtStartPar
\sphinxstylestrong{Example}
\end{quote}

\begin{sphinxVerbatim}[commandchars=\\\{\}]
\PYG{c+c1}{\PYGZsh{} Add the term C1 * X to expression \PYGZsq{}expr\PYGZsq{}}
\PYG{n}{expr}\PYG{o}{.}\PYG{n}{addTerm}\PYG{p}{(}\PYG{n}{x}\PYG{p}{,} \PYG{n}{C1}\PYG{p}{)}
\end{sphinxVerbatim}

\subsubsection{PsdExpr.addTerms()}
\label{\detokenize{pyapiref:id2}}\begin{quote}

\sphinxAtStartPar
\sphinxstylestrong{Synopsis}
\begin{quote}

\sphinxAtStartPar
\sphinxcode{\sphinxupquote{addTerms(vars, mats)}}
\end{quote}

\sphinxAtStartPar
\sphinxstylestrong{Description}
\begin{quote}

\sphinxAtStartPar
Adds multiple new terms to the LMI expression.

\sphinxAtStartPar
If argument \sphinxcode{\sphinxupquote{vars}} is {\hyperref[\detokenize{pyapiref:chappyapi-psdvar}]{\sphinxcrossref{\DUrole{std,std-ref}{PsdVar Class}}}} object,
then argument \sphinxcode{\sphinxupquote{mats}} is {\hyperref[\detokenize{pyapiref:chappyapi-symmatrix}]{\sphinxcrossref{\DUrole{std,std-ref}{SymMatrix Class}}}} object;
If argument \sphinxcode{\sphinxupquote{vars}} is {\hyperref[\detokenize{pyapiref:chappyapi-psdvararray}]{\sphinxcrossref{\DUrole{std,std-ref}{PsdVarArray Class}}}} object or list,
then argument \sphinxcode{\sphinxupquote{mats}} is {\hyperref[\detokenize{pyapiref:chappyapi-symmatrixarray}]{\sphinxcrossref{\DUrole{std,std-ref}{SymMatrixArray Class}}}} or list;
\end{quote}

\sphinxAtStartPar
\sphinxstylestrong{Arguments}
\begin{quote}

\sphinxAtStartPar
\sphinxcode{\sphinxupquote{vars}}
\begin{quote}

\sphinxAtStartPar
Array of variables to add new items to.
\end{quote}

\sphinxAtStartPar
\sphinxcode{\sphinxupquote{mats}}
\begin{quote}

\sphinxAtStartPar
The symmetric array of matrices to add new items to.
\end{quote}
\end{quote}

\sphinxAtStartPar
\sphinxstylestrong{Example}
\end{quote}

\begin{sphinxVerbatim}[commandchars=\\\{\}]
\PYG{c+c1}{\PYGZsh{} Add terms x1 * C1 + x2 * C2 to expression \PYGZsq{}expr\PYGZsq{}}
\PYG{n}{expr}\PYG{o}{.}\PYG{n}{addTerms}\PYG{p}{(}\PYG{p}{[}\PYG{n}{x1}\PYG{p}{,} \PYG{n}{x2}\PYG{p}{]}\PYG{p}{,} \PYG{p}{[}\PYG{n}{C1}\PYG{p}{,} \PYG{n}{C2}\PYG{p}{]}\PYG{p}{)}
\end{sphinxVerbatim}

\subsubsection{LmiExpr.addLmiExpr()}
\label{\detokenize{pyapiref:lmiexpr-addlmiexpr}}\begin{quote}

\sphinxAtStartPar
\sphinxstylestrong{Synopsis}
\begin{quote}

\sphinxAtStartPar
\sphinxcode{\sphinxupquote{addLmiExpr(expr, mult=1.0)}}
\end{quote}

\sphinxAtStartPar
\sphinxstylestrong{Description}
\begin{quote}

\sphinxAtStartPar
Add a new LMI expression to the current LMI expression.
\end{quote}

\sphinxAtStartPar
\sphinxstylestrong{Arguments}
\begin{quote}

\sphinxAtStartPar
\sphinxcode{\sphinxupquote{expr}}
\begin{quote}

\sphinxAtStartPar
LMI expression to add.
\end{quote}

\sphinxAtStartPar
\sphinxcode{\sphinxupquote{mult}}
\begin{quote}

\sphinxAtStartPar
Magnification coefficient for the added LMI expression. Optional, 1.0 by default.
\end{quote}
\end{quote}

\sphinxAtStartPar
\sphinxstylestrong{Example}
\end{quote}

\begin{sphinxVerbatim}[commandchars=\\\{\}]
\PYG{c+c1}{\PYGZsh{} Add linear expression 2 * x * C to \PYGZsq{}expr\PYGZsq{}}
\PYG{n}{expr}\PYG{o}{.}\PYG{n}{addLinExpr}\PYG{p}{(}\PYG{n}{x} \PYG{o}{*} \PYG{n}{C}\PYG{p}{,} \PYG{l+m+mf}{2.0}\PYG{p}{)}
\end{sphinxVerbatim}

\subsubsection{LmiExpr.clone()}
\label{\detokenize{pyapiref:lmiexpr-clone}}\begin{quote}

\sphinxAtStartPar
\sphinxstylestrong{Synopsis}
\begin{quote}

\sphinxAtStartPar
\sphinxcode{\sphinxupquote{clone()}}
\end{quote}

\sphinxAtStartPar
\sphinxstylestrong{Description}
\begin{quote}

\sphinxAtStartPar
Create a deep copy of the expression.
\end{quote}

\sphinxAtStartPar
\sphinxstylestrong{Example}
\end{quote}

\begin{sphinxVerbatim}[commandchars=\\\{\}]
\PYG{c+c1}{\PYGZsh{} Create a deep copy of expression expr}
\PYG{n}{exprcopy} \PYG{o}{=} \PYG{n}{expr}\PYG{o}{.}\PYG{n}{clone}\PYG{p}{(}\PYG{p}{)}
\end{sphinxVerbatim}

\subsubsection{LmiExpr.reserve()}
\label{\detokenize{pyapiref:lmiexpr-reserve}}\begin{quote}

\sphinxAtStartPar
\sphinxstylestrong{Synopsis}
\begin{quote}

\sphinxAtStartPar
\sphinxcode{\sphinxupquote{reserve(n)}}
\end{quote}

\sphinxAtStartPar
\sphinxstylestrong{Description}
\begin{quote}

\sphinxAtStartPar
Pre\sphinxhyphen{}allocate space for LMI expression object.
\end{quote}

\sphinxAtStartPar
\sphinxstylestrong{Arguments}
\begin{quote}

\sphinxAtStartPar
\sphinxcode{\sphinxupquote{n}}
\begin{quote}

\sphinxAtStartPar
Number of terms to be allocated.
\end{quote}
\end{quote}

\sphinxAtStartPar
\sphinxstylestrong{Example}
\end{quote}

\begin{sphinxVerbatim}[commandchars=\\\{\}]
\PYG{c+c1}{\PYGZsh{} Allocate 100 terms for LMI expression \PYGZsq{}expr\PYGZsq{}}
\PYG{n}{expr}\PYG{o}{.}\PYG{n}{reserve}\PYG{p}{(}\PYG{l+m+mi}{100}\PYG{p}{)}
\end{sphinxVerbatim}

\subsubsection{LmiExpr.remove()}
\label{\detokenize{pyapiref:lmiexpr-remove}}\begin{quote}

\sphinxAtStartPar
\sphinxstylestrong{Synopsis}
\begin{quote}

\sphinxAtStartPar
\sphinxcode{\sphinxupquote{remove(item)}}
\end{quote}

\sphinxAtStartPar
\sphinxstylestrong{Description}
\begin{quote}

\sphinxAtStartPar
Remove a term from the LMI expression.

\sphinxAtStartPar
If argument \sphinxcode{\sphinxupquote{item}} is constant, then remove the term stored at index i of the expression; otherwise argument \sphinxcode{\sphinxupquote{item}} is {\hyperref[\detokenize{pyapiref:chappyapi-var}]{\sphinxcrossref{\DUrole{std,std-ref}{Var Class}}}} object.
\end{quote}

\sphinxAtStartPar
\sphinxstylestrong{Arguments}
\begin{quote}

\sphinxAtStartPar
\sphinxcode{\sphinxupquote{item}}
\begin{quote}

\sphinxAtStartPar
Constant index or {\hyperref[\detokenize{pyapiref:chappyapi-var}]{\sphinxcrossref{\DUrole{std,std-ref}{Var Class}}}} variable of the term to be removed.
\end{quote}
\end{quote}

\sphinxAtStartPar
\sphinxstylestrong{Example}
\end{quote}

\begin{sphinxVerbatim}[commandchars=\\\{\}]
\PYG{c+c1}{\PYGZsh{} Remove the term whose index is 2 from positive semi\PYGZhy{}definite expression expr}
\PYG{n}{expr}\PYG{o}{.}\PYG{n}{remove}\PYG{p}{(}\PYG{l+m+mi}{2}\PYG{p}{)}
\PYG{c+c1}{\PYGZsh{} Remove the terms one of which variable is x from LMI expression expr}
\PYG{n}{expr}\PYG{o}{.}\PYG{n}{remove}\PYG{p}{(}\PYG{n}{x}\PYG{p}{)}
\end{sphinxVerbatim}

\subsection{NlExpr Class}
\label{\detokenize{pyapiref:nlexpr-class}}\label{\detokenize{pyapiref:chappyapi-nlexpr}}
\sphinxAtStartPar
The \sphinxcode{\sphinxupquote{NlExpr}} class provides interfaces for constructing nonlinear expressions
in the COPT. It supports variable combination operations and offers the following member functions:

\subsubsection{NlExpr()}
\label{\detokenize{pyapiref:nlexpr}}\begin{quote}

\sphinxAtStartPar
\sphinxstylestrong{Synopsis}
\begin{quote}

\sphinxAtStartPar
\sphinxcode{\sphinxupquote{NlExpr(arg1=0.0, arg2=None)}}
\end{quote}

\sphinxAtStartPar
\sphinxstylestrong{Description}
\begin{quote}

\sphinxAtStartPar
Creates a {\hyperref[\detokenize{pyapiref:chappyapi-nlexpr}]{\sphinxcrossref{\DUrole{std,std-ref}{NlExpr Class}}}} object.

\sphinxAtStartPar
If \sphinxcode{\sphinxupquote{arg1}} is a float and \sphinxcode{\sphinxupquote{arg2}} is \sphinxcode{\sphinxupquote{None}}, a constant nonlinear
expression is created with value \sphinxcode{\sphinxupquote{arg1}}.

\sphinxAtStartPar
If \sphinxcode{\sphinxupquote{arg1}} is a {\hyperref[\detokenize{pyapiref:chappyapi-var}]{\sphinxcrossref{\DUrole{std,std-ref}{Var Class}}}}, {\hyperref[\detokenize{pyapiref:chappyapi-linexpr}]{\sphinxcrossref{\DUrole{std,std-ref}{LinExpr Class}}}},
{\hyperref[\detokenize{pyapiref:chappyapi-quadexpr}]{\sphinxcrossref{\DUrole{std,std-ref}{QuadExpr Class}}}}, or {\hyperref[\detokenize{pyapiref:chappyapi-exprbuilder}]{\sphinxcrossref{\DUrole{std,std-ref}{ExprBuilder Class}}}},
then \sphinxcode{\sphinxupquote{arg2}} can be a float or \sphinxcode{\sphinxupquote{None}}, indicating that the expression
is scaled by \sphinxcode{\sphinxupquote{arg2}} (default is 1.0).
\end{quote}

\sphinxAtStartPar
\sphinxstylestrong{Arguments}
\begin{quote}

\sphinxAtStartPar
\sphinxcode{\sphinxupquote{arg1}}
\begin{quote}

\sphinxAtStartPar
Optional argument. Default is 0.0.

\sphinxAtStartPar
Acceptable types include float, {\hyperref[\detokenize{pyapiref:chappyapi-var}]{\sphinxcrossref{\DUrole{std,std-ref}{Var Class}}}}, {\hyperref[\detokenize{pyapiref:chappyapi-linexpr}]{\sphinxcrossref{\DUrole{std,std-ref}{LinExpr Class}}}},
{\hyperref[\detokenize{pyapiref:chappyapi-quadexpr}]{\sphinxcrossref{\DUrole{std,std-ref}{QuadExpr Class}}}}, or {\hyperref[\detokenize{pyapiref:chappyapi-exprbuilder}]{\sphinxcrossref{\DUrole{std,std-ref}{ExprBuilder Class}}}}.
\end{quote}

\sphinxAtStartPar
\sphinxcode{\sphinxupquote{arg2}}
\begin{quote}

\sphinxAtStartPar
Coefficient of the expression.

\sphinxAtStartPar
Optional argument. Default is \sphinxcode{\sphinxupquote{None}}, which is interpreted as 1.0.
\end{quote}
\end{quote}

\sphinxAtStartPar
\sphinxstylestrong{Example}
\end{quote}

\begin{sphinxVerbatim}[commandchars=\\\{\}]
\PYG{c+c1}{\PYGZsh{} Variable term with coefficient: 2 * x}
\PYG{n}{expr1} \PYG{o}{=} \PYG{n}{NlExpr}\PYG{p}{(}\PYG{n}{x}\PYG{p}{,} \PYG{l+m+mi}{2}\PYG{p}{)}
\PYG{c+c1}{\PYGZsh{} Linear expression: 2 * x + 3 * y}
\PYG{n}{lin} \PYG{o}{=} \PYG{l+m+mi}{2} \PYG{o}{*} \PYG{n}{x} \PYG{o}{+} \PYG{l+m+mi}{3} \PYG{o}{*} \PYG{n}{y}
\PYG{n}{expr2} \PYG{o}{=} \PYG{n}{NlExpr}\PYG{p}{(}\PYG{n}{lin}\PYG{p}{)}
\end{sphinxVerbatim}

\subsubsection{NlExpr.getConstant()}
\label{\detokenize{pyapiref:nlexpr-getconstant}}\begin{quote}

\sphinxAtStartPar
\sphinxstylestrong{Synopsis}
\begin{quote}

\sphinxAtStartPar
\sphinxcode{\sphinxupquote{getConstant()}}
\end{quote}

\sphinxAtStartPar
\sphinxstylestrong{Description}
\begin{quote}

\sphinxAtStartPar
Returns the constant term in the nonlinear expression.
\end{quote}

\sphinxAtStartPar
\sphinxstylestrong{Example}
\end{quote}

\begin{sphinxVerbatim}[commandchars=\\\{\}]
\PYG{c+c1}{\PYGZsh{} Get the constant term of nexpr}
\PYG{n}{constant} \PYG{o}{=} \PYG{n}{nexpr}\PYG{o}{.}\PYG{n}{getConstant}\PYG{p}{(}\PYG{p}{)}
\end{sphinxVerbatim}

\subsubsection{NlExpr.getSize()}
\label{\detokenize{pyapiref:nlexpr-getsize}}\begin{quote}

\sphinxAtStartPar
\sphinxstylestrong{Synopsis}
\begin{quote}

\sphinxAtStartPar
\sphinxcode{\sphinxupquote{getSize()}}
\end{quote}

\sphinxAtStartPar
\sphinxstylestrong{Description}
\begin{quote}

\sphinxAtStartPar
Returns the number of terms in the nonlinear expression.
\end{quote}

\sphinxAtStartPar
\sphinxstylestrong{Example}
\end{quote}

\begin{sphinxVerbatim}[commandchars=\\\{\}]
\PYG{c+c1}{\PYGZsh{} Get the number of elements in nexpr}
\PYG{n}{exprsize} \PYG{o}{=} \PYG{n}{nexpr}\PYG{o}{.}\PYG{n}{getSize}\PYG{p}{(}\PYG{p}{)}
\end{sphinxVerbatim}

\subsubsection{NlExpr.getValue()}
\label{\detokenize{pyapiref:nlexpr-getvalue}}\begin{quote}

\sphinxAtStartPar
\sphinxstylestrong{Synopsis}
\begin{quote}

\sphinxAtStartPar
\sphinxcode{\sphinxupquote{getValue()}}
\end{quote}

\sphinxAtStartPar
\sphinxstylestrong{Description}
\begin{quote}

\sphinxAtStartPar
Returns the value of the nonlinear expression based on the current values of the variables.
\end{quote}

\sphinxAtStartPar
\sphinxstylestrong{Example}
\end{quote}

\begin{sphinxVerbatim}[commandchars=\\\{\}]
\PYG{c+c1}{\PYGZsh{} Get the current value of the nonlinear expression nlexpr}
\PYG{n}{val} \PYG{o}{=} \PYG{n}{nlexpr}\PYG{o}{.}\PYG{n}{getValue}\PYG{p}{(}\PYG{p}{)}
\end{sphinxVerbatim}

\subsubsection{NlExpr.setConstant()}
\label{\detokenize{pyapiref:nlexpr-setconstant}}\begin{quote}

\sphinxAtStartPar
\sphinxstylestrong{Synopsis}
\begin{quote}

\sphinxAtStartPar
\sphinxcode{\sphinxupquote{setConstant(newval)}}
\end{quote}

\sphinxAtStartPar
\sphinxstylestrong{Description}
\begin{quote}

\sphinxAtStartPar
Sets the constant term of the nonlinear expression.
\end{quote}

\sphinxAtStartPar
\sphinxstylestrong{Arguments}
\begin{quote}

\sphinxAtStartPar
\sphinxcode{\sphinxupquote{newval}}
\begin{quote}

\sphinxAtStartPar
The constant value to set.
\end{quote}
\end{quote}

\sphinxAtStartPar
\sphinxstylestrong{Example}
\end{quote}

\begin{sphinxVerbatim}[commandchars=\\\{\}]
\PYG{c+c1}{\PYGZsh{} Set constant of nexpr to 2.0}
\PYG{n}{nexpr}\PYG{o}{.}\PYG{n}{setConstant}\PYG{p}{(}\PYG{l+m+mf}{2.0}\PYG{p}{)}
\end{sphinxVerbatim}

\subsubsection{NlExpr.addConstant()}
\label{\detokenize{pyapiref:nlexpr-addconstant}}\begin{quote}

\sphinxAtStartPar
\sphinxstylestrong{Synopsis}
\begin{quote}

\sphinxAtStartPar
\sphinxcode{\sphinxupquote{addConstant(newval)}}
\end{quote}

\sphinxAtStartPar
\sphinxstylestrong{Description}
\begin{quote}

\sphinxAtStartPar
Adds a constant term to the nonlinear expression.
\end{quote}

\sphinxAtStartPar
\sphinxstylestrong{Arguments}
\begin{quote}

\sphinxAtStartPar
\sphinxcode{\sphinxupquote{newval}}
\begin{quote}

\sphinxAtStartPar
The constant value to add.
\end{quote}
\end{quote}

\sphinxAtStartPar
\sphinxstylestrong{Example}
\end{quote}

\begin{sphinxVerbatim}[commandchars=\\\{\}]
\PYG{c+c1}{\PYGZsh{} Add constant 2.0 to nexpr}
\PYG{n}{nexpr}\PYG{o}{.}\PYG{n}{addConstant}\PYG{p}{(}\PYG{l+m+mf}{2.0}\PYG{p}{)}
\end{sphinxVerbatim}

\subsubsection{NlExpr.addTerm()}
\label{\detokenize{pyapiref:nlexpr-addterm}}\begin{quote}

\sphinxAtStartPar
\sphinxstylestrong{Synopsis}
\begin{quote}

\sphinxAtStartPar
\sphinxcode{\sphinxupquote{addTerm(var, coeff=1.0)}}
\end{quote}

\sphinxAtStartPar
\sphinxstylestrong{Description}
\begin{quote}

\sphinxAtStartPar
Adds a single linear term to the nonlinear expression.
\end{quote}

\sphinxAtStartPar
\sphinxstylestrong{Arguments}
\begin{quote}

\sphinxAtStartPar
\sphinxcode{\sphinxupquote{var}}
\begin{quote}

\sphinxAtStartPar
Variable in the new term.
\end{quote}

\sphinxAtStartPar
\sphinxcode{\sphinxupquote{coeff}}
\begin{quote}

\sphinxAtStartPar
Coefficient of the new term. Optional, default is 1.0.
\end{quote}
\end{quote}

\sphinxAtStartPar
\sphinxstylestrong{Example}
\end{quote}

\begin{sphinxVerbatim}[commandchars=\\\{\}]
\PYG{c+c1}{\PYGZsh{} Add term 2*x to nexpr}
\PYG{n}{nexpr}\PYG{o}{.}\PYG{n}{addTerm}\PYG{p}{(}\PYG{n}{x}\PYG{p}{,} \PYG{l+m+mi}{2}\PYG{p}{)}
\end{sphinxVerbatim}

\subsubsection{NlExpr.addTerms()}
\label{\detokenize{pyapiref:nlexpr-addterms}}\begin{quote}

\sphinxAtStartPar
\sphinxstylestrong{Synopsis}
\begin{quote}

\sphinxAtStartPar
\sphinxcode{\sphinxupquote{addTerms(vars, coeffs)}}
\end{quote}

\sphinxAtStartPar
\sphinxstylestrong{Description}
\begin{quote}

\sphinxAtStartPar
Adds one or more linear terms to the nonlinear expression.

\sphinxAtStartPar
If \sphinxcode{\sphinxupquote{vars}} is a {\hyperref[\detokenize{pyapiref:chappyapi-var}]{\sphinxcrossref{\DUrole{std,std-ref}{Var Class}}}} object, \sphinxcode{\sphinxupquote{coeffs}} may be a single
constant.

\sphinxAtStartPar
If \sphinxcode{\sphinxupquote{vars}} is a {\hyperref[\detokenize{pyapiref:chappyapi-vararray}]{\sphinxcrossref{\DUrole{std,std-ref}{VarArray Class}}}} or an iterable of
{\hyperref[\detokenize{pyapiref:chappyapi-var}]{\sphinxcrossref{\DUrole{std,std-ref}{Var Class}}}} objects, \sphinxcode{\sphinxupquote{coeffs}} may be a single constant or a
list of constants corresponding to each variable.

\sphinxAtStartPar
If \sphinxcode{\sphinxupquote{vars}} is a mapping from arbitrary keys to variables,
\sphinxcode{\sphinxupquote{coeffs}} may be a single constant or a mapping with the same keys
specifying individual coefficients.
\end{quote}

\sphinxAtStartPar
\sphinxstylestrong{Arguments}
\begin{quote}

\sphinxAtStartPar
\sphinxcode{\sphinxupquote{vars}}
\begin{quote}

\sphinxAtStartPar
Variables to add.
\end{quote}

\sphinxAtStartPar
\sphinxcode{\sphinxupquote{coeffs}}
\begin{quote}

\sphinxAtStartPar
Coefficients for the terms.
\end{quote}
\end{quote}

\sphinxAtStartPar
\sphinxstylestrong{Example}
\end{quote}

\begin{sphinxVerbatim}[commandchars=\\\{\}]
\PYG{c+c1}{\PYGZsh{} Add terms 2*x + 3*y to nexpr}
\PYG{n}{nexpr}\PYG{o}{.}\PYG{n}{addTerms}\PYG{p}{(}\PYG{p}{[}\PYG{n}{x}\PYG{p}{,} \PYG{n}{y}\PYG{p}{]}\PYG{p}{,} \PYG{p}{[}\PYG{l+m+mf}{2.0}\PYG{p}{,} \PYG{l+m+mf}{3.0}\PYG{p}{]}\PYG{p}{)}
\end{sphinxVerbatim}

\subsubsection{NlExpr.addLinExpr()}
\label{\detokenize{pyapiref:nlexpr-addlinexpr}}\begin{quote}

\sphinxAtStartPar
\sphinxstylestrong{Synopsis}
\begin{quote}

\sphinxAtStartPar
\sphinxcode{\sphinxupquote{addLinExpr(expr, mult=1.0)}}
\end{quote}

\sphinxAtStartPar
\sphinxstylestrong{Description}
\begin{quote}

\sphinxAtStartPar
Adds a linear expression \sphinxcode{\sphinxupquote{expr}} to the current nonlinear expression,
scaling all terms by the given multiplier \sphinxcode{\sphinxupquote{mult}}.
\end{quote}

\sphinxAtStartPar
\sphinxstylestrong{Arguments}
\begin{quote}

\sphinxAtStartPar
\sphinxcode{\sphinxupquote{expr}}
\begin{quote}

\sphinxAtStartPar
Linear expression to be added.

\sphinxAtStartPar
Must be a {\hyperref[\detokenize{pyapiref:chappyapi-exprbuilder}]{\sphinxcrossref{\DUrole{std,std-ref}{ExprBuilder Class}}}} or {\hyperref[\detokenize{pyapiref:chappyapi-linexpr}]{\sphinxcrossref{\DUrole{std,std-ref}{LinExpr Class}}}} object.
\end{quote}

\sphinxAtStartPar
\sphinxcode{\sphinxupquote{mult}}
\begin{quote}

\sphinxAtStartPar
Scaling factor applied to all terms in \sphinxcode{\sphinxupquote{expr}}.

\sphinxAtStartPar
Optional argument. Default is 1.0.
\end{quote}
\end{quote}

\sphinxAtStartPar
\sphinxstylestrong{Example}
\end{quote}

\begin{sphinxVerbatim}[commandchars=\\\{\}]
\PYG{c+c1}{\PYGZsh{} Add a scaled linear expression to nexpr}
\PYG{n}{lexpr} \PYG{o}{=} \PYG{n}{LinExpr}\PYG{p}{(}\PYG{n}{x}\PYG{p}{,} \PYG{l+m+mf}{2.0}\PYG{p}{)} \PYG{o}{+} \PYG{n}{y} \PYG{o}{*} \PYG{l+m+mf}{3.0}
\PYG{n}{nexpr} \PYG{o}{=} \PYG{n}{NlExpr}\PYG{p}{(}\PYG{p}{)}
\PYG{n}{nexpr}\PYG{o}{.}\PYG{n}{addLinExpr}\PYG{p}{(}\PYG{n}{lexpr}\PYG{p}{)}
\PYG{n}{nexpr}\PYG{o}{.}\PYG{n}{addLinExpr}\PYG{p}{(}\PYG{n}{lexpr}\PYG{p}{,} \PYG{n}{mult}\PYG{o}{=}\PYG{l+m+mf}{0.5}\PYG{p}{)}
\end{sphinxVerbatim}

\subsubsection{NlExpr.addQuadExpr()}
\label{\detokenize{pyapiref:nlexpr-addquadexpr}}\begin{quote}

\sphinxAtStartPar
\sphinxstylestrong{Synopsis}
\begin{quote}

\sphinxAtStartPar
\sphinxcode{\sphinxupquote{addQuadExpr(expr, mult=1.0)}}
\end{quote}

\sphinxAtStartPar
\sphinxstylestrong{Description}
\begin{quote}

\sphinxAtStartPar
Adds a quadratic expression \sphinxcode{\sphinxupquote{expr}} to the current nonlinear expression,
with all terms scaled by \sphinxcode{\sphinxupquote{mult}}.
\end{quote}

\sphinxAtStartPar
\sphinxstylestrong{Arguments}
\begin{quote}

\sphinxAtStartPar
\sphinxcode{\sphinxupquote{expr}}
\begin{quote}

\sphinxAtStartPar
Quadratic expression to add.

\sphinxAtStartPar
Must be a {\hyperref[\detokenize{pyapiref:chappyapi-quadexpr}]{\sphinxcrossref{\DUrole{std,std-ref}{QuadExpr Class}}}} object.
\end{quote}

\sphinxAtStartPar
\sphinxcode{\sphinxupquote{mult}}
\begin{quote}

\sphinxAtStartPar
Scaling factor for all terms.

\sphinxAtStartPar
Optional. Default is 1.0.
\end{quote}
\end{quote}
\end{quote}

\subsubsection{NlExpr.addNlExpr()}
\label{\detokenize{pyapiref:nlexpr-addnlexpr}}\begin{quote}

\sphinxAtStartPar
\sphinxstylestrong{Synopsis}
\begin{quote}

\sphinxAtStartPar
\sphinxcode{\sphinxupquote{addNlExpr(expr, mult=1.0)}}
\end{quote}

\sphinxAtStartPar
\sphinxstylestrong{Description}
\begin{quote}

\sphinxAtStartPar
Adds another nonlinear expression \sphinxcode{\sphinxupquote{expr}} to the current expression,
with all terms scaled by \sphinxcode{\sphinxupquote{mult}}.
\end{quote}

\sphinxAtStartPar
\sphinxstylestrong{Arguments}
\begin{quote}

\sphinxAtStartPar
\sphinxcode{\sphinxupquote{expr}}
\begin{quote}

\sphinxAtStartPar
Nonlinear expression to add.

\sphinxAtStartPar
Must be a {\hyperref[\detokenize{pyapiref:chappyapi-nlexpr}]{\sphinxcrossref{\DUrole{std,std-ref}{NlExpr Class}}}} object.
\end{quote}

\sphinxAtStartPar
\sphinxcode{\sphinxupquote{mult}}
\begin{quote}

\sphinxAtStartPar
Scaling factor.

\sphinxAtStartPar
Optional. Default is 1.0.
\end{quote}
\end{quote}
\end{quote}

\subsubsection{NlExpr.clear()}
\label{\detokenize{pyapiref:nlexpr-clear}}\begin{quote}

\sphinxAtStartPar
\sphinxstylestrong{Synopsis}
\begin{quote}

\sphinxAtStartPar
\sphinxcode{\sphinxupquote{clear()}}
\end{quote}

\sphinxAtStartPar
\sphinxstylestrong{Description}
\begin{quote}

\sphinxAtStartPar
Clears all contents of the nonlinear expression.
\end{quote}
\end{quote}

\subsubsection{NlExpr.clone()}
\label{\detokenize{pyapiref:nlexpr-clone}}\begin{quote}

\sphinxAtStartPar
\sphinxstylestrong{Synopsis}
\begin{quote}

\sphinxAtStartPar
\sphinxcode{\sphinxupquote{clone()}}
\end{quote}

\sphinxAtStartPar
\sphinxstylestrong{Description}
\begin{quote}

\sphinxAtStartPar
Creates a deep copy of the current nonlinear expression.
\end{quote}

\sphinxAtStartPar
\sphinxstylestrong{Return Value}
\begin{quote}

\sphinxAtStartPar
A new {\hyperref[\detokenize{pyapiref:chappyapi-nlexpr}]{\sphinxcrossref{\DUrole{std,std-ref}{NlExpr Class}}}} object.
\end{quote}
\end{quote}

\subsubsection{NlExpr.negate()}
\label{\detokenize{pyapiref:nlexpr-negate}}\begin{quote}

\sphinxAtStartPar
\sphinxstylestrong{Synopsis}
\begin{quote}

\sphinxAtStartPar
\sphinxcode{\sphinxupquote{negate()}}
\end{quote}

\sphinxAtStartPar
\sphinxstylestrong{Description}
\begin{quote}

\sphinxAtStartPar
Negates the current nonlinear expression by flipping the sign of each term.
\end{quote}
\end{quote}

\subsubsection{NlExpr.reserve()}
\label{\detokenize{pyapiref:nlexpr-reserve}}\begin{quote}

\sphinxAtStartPar
\sphinxstylestrong{Synopsis}
\begin{quote}

\sphinxAtStartPar
\sphinxcode{\sphinxupquote{reserve(n)}}
\end{quote}

\sphinxAtStartPar
\sphinxstylestrong{Description}
\begin{quote}

\sphinxAtStartPar
Preallocates space for \sphinxcode{\sphinxupquote{n}} terms in the expression.
\end{quote}

\sphinxAtStartPar
\sphinxstylestrong{Arguments}
\begin{quote}

\sphinxAtStartPar
\sphinxcode{\sphinxupquote{n}}
\begin{quote}

\sphinxAtStartPar
Number of terms to reserve.
\end{quote}
\end{quote}
\end{quote}

\subsubsection{NlExpr.size()}
\label{\detokenize{pyapiref:nlexpr-size}}\begin{quote}

\sphinxAtStartPar
\sphinxstylestrong{Synopsis}
\begin{quote}

\sphinxAtStartPar
\sphinxcode{\sphinxupquote{size()}}
\end{quote}

\sphinxAtStartPar
\sphinxstylestrong{Description}
\begin{quote}

\sphinxAtStartPar
Returns the number of nonlinear terms in the expression.
\end{quote}

\sphinxAtStartPar
\sphinxstylestrong{Return Value}
\begin{quote}

\sphinxAtStartPar
An integer.
\end{quote}
\end{quote}

\subsection{nl Namespace}
\label{\detokenize{pyapiref:nl-namespace}}\label{\detokenize{pyapiref:chappyapi-nl}}

\subsubsection{nl.abs()}
\label{\detokenize{pyapiref:nl-abs}}\begin{quote}

\sphinxAtStartPar
\sphinxstylestrong{Synopsis}
\begin{quote}

\sphinxAtStartPar
\sphinxcode{\sphinxupquote{NlExpr abs(expr)}}
\end{quote}

\sphinxAtStartPar
\sphinxstylestrong{Description}
\begin{quote}

\sphinxAtStartPar
Calculate the absolute value of the input object.
\end{quote}

\sphinxAtStartPar
\sphinxstylestrong{Arguments}
\begin{quote}

\sphinxAtStartPar
\sphinxcode{\sphinxupquote{expr}}
\begin{quote}

\sphinxAtStartPar
The object to be calculated.

\sphinxAtStartPar
Possible values include: constant, {\hyperref[\detokenize{pyapiref:chappyapi-var}]{\sphinxcrossref{\DUrole{std,std-ref}{Var Class}}}} object,
{\hyperref[\detokenize{pyapiref:chappyapi-linexpr}]{\sphinxcrossref{\DUrole{std,std-ref}{LinExpr Class}}}}, {\hyperref[\detokenize{pyapiref:chappyapi-quadexpr}]{\sphinxcrossref{\DUrole{std,std-ref}{QuadExpr Class}}}},
{\hyperref[\detokenize{pyapiref:chappyapi-nlexpr}]{\sphinxcrossref{\DUrole{std,std-ref}{NlExpr Class}}}} object, or {\hyperref[\detokenize{pyapiref:chappyapi-exprbuilder}]{\sphinxcrossref{\DUrole{std,std-ref}{ExprBuilder Class}}}} object.
\end{quote}
\end{quote}

\sphinxAtStartPar
\sphinxstylestrong{Return Value}
\begin{quote}

\sphinxAtStartPar
Returns a {\hyperref[\detokenize{pyapiref:chappyapi-nlexpr}]{\sphinxcrossref{\DUrole{std,std-ref}{NlExpr Class}}}} object.
\end{quote}
\end{quote}

\subsubsection{nl.acos()}
\label{\detokenize{pyapiref:nl-acos}}\begin{quote}

\sphinxAtStartPar
\sphinxstylestrong{Synopsis}
\begin{quote}

\sphinxAtStartPar
\sphinxcode{\sphinxupquote{NlExpr acos(expr)}}
\end{quote}

\sphinxAtStartPar
\sphinxstylestrong{Description}
\begin{quote}

\sphinxAtStartPar
Calculate the arc cosine of the input object.
\end{quote}

\sphinxAtStartPar
\sphinxstylestrong{Arguments}
\begin{quote}

\sphinxAtStartPar
\sphinxcode{\sphinxupquote{expr}}
\begin{quote}

\sphinxAtStartPar
The object to be calculated.

\sphinxAtStartPar
Possible values include: constant, {\hyperref[\detokenize{pyapiref:chappyapi-var}]{\sphinxcrossref{\DUrole{std,std-ref}{Var Class}}}} object,
{\hyperref[\detokenize{pyapiref:chappyapi-linexpr}]{\sphinxcrossref{\DUrole{std,std-ref}{LinExpr Class}}}}, {\hyperref[\detokenize{pyapiref:chappyapi-quadexpr}]{\sphinxcrossref{\DUrole{std,std-ref}{QuadExpr Class}}}},
{\hyperref[\detokenize{pyapiref:chappyapi-nlexpr}]{\sphinxcrossref{\DUrole{std,std-ref}{NlExpr Class}}}} object, or {\hyperref[\detokenize{pyapiref:chappyapi-exprbuilder}]{\sphinxcrossref{\DUrole{std,std-ref}{ExprBuilder Class}}}} object.
\end{quote}
\end{quote}

\sphinxAtStartPar
\sphinxstylestrong{Return Value}
\begin{quote}

\sphinxAtStartPar
Returns a {\hyperref[\detokenize{pyapiref:chappyapi-nlexpr}]{\sphinxcrossref{\DUrole{std,std-ref}{NlExpr Class}}}} object.
\end{quote}
\end{quote}

\subsubsection{nl.acosh()}
\label{\detokenize{pyapiref:nl-acosh}}\begin{quote}

\sphinxAtStartPar
\sphinxstylestrong{Synopsis}
\begin{quote}

\sphinxAtStartPar
\sphinxcode{\sphinxupquote{NlExpr acosh(expr)}}
\end{quote}

\sphinxAtStartPar
\sphinxstylestrong{Description}
\begin{quote}

\sphinxAtStartPar
Calculate the inverse hyperbolic cosine of the input object.
\end{quote}

\sphinxAtStartPar
\sphinxstylestrong{Arguments}
\begin{quote}

\sphinxAtStartPar
\sphinxcode{\sphinxupquote{expr}}
\begin{quote}

\sphinxAtStartPar
The object to be calculated.

\sphinxAtStartPar
Possible values include: constant, {\hyperref[\detokenize{pyapiref:chappyapi-var}]{\sphinxcrossref{\DUrole{std,std-ref}{Var Class}}}} object,
{\hyperref[\detokenize{pyapiref:chappyapi-linexpr}]{\sphinxcrossref{\DUrole{std,std-ref}{LinExpr Class}}}}, {\hyperref[\detokenize{pyapiref:chappyapi-quadexpr}]{\sphinxcrossref{\DUrole{std,std-ref}{QuadExpr Class}}}},
{\hyperref[\detokenize{pyapiref:chappyapi-nlexpr}]{\sphinxcrossref{\DUrole{std,std-ref}{NlExpr Class}}}} object, or {\hyperref[\detokenize{pyapiref:chappyapi-exprbuilder}]{\sphinxcrossref{\DUrole{std,std-ref}{ExprBuilder Class}}}} object.
\end{quote}
\end{quote}

\sphinxAtStartPar
\sphinxstylestrong{Return Value}
\begin{quote}

\sphinxAtStartPar
Returns a {\hyperref[\detokenize{pyapiref:chappyapi-nlexpr}]{\sphinxcrossref{\DUrole{std,std-ref}{NlExpr Class}}}} object.
\end{quote}
\end{quote}

\subsubsection{nl.asin()}
\label{\detokenize{pyapiref:nl-asin}}\begin{quote}

\sphinxAtStartPar
\sphinxstylestrong{Synopsis}
\begin{quote}

\sphinxAtStartPar
\sphinxcode{\sphinxupquote{NlExpr asin(expr)}}
\end{quote}

\sphinxAtStartPar
\sphinxstylestrong{Description}
\begin{quote}

\sphinxAtStartPar
Calculate the arc sine of the input object.
\end{quote}

\sphinxAtStartPar
\sphinxstylestrong{Arguments}
\begin{quote}

\sphinxAtStartPar
\sphinxcode{\sphinxupquote{expr}}
\begin{quote}

\sphinxAtStartPar
The object to be calculated.

\sphinxAtStartPar
Possible values include: constant, {\hyperref[\detokenize{pyapiref:chappyapi-var}]{\sphinxcrossref{\DUrole{std,std-ref}{Var Class}}}} object,
{\hyperref[\detokenize{pyapiref:chappyapi-linexpr}]{\sphinxcrossref{\DUrole{std,std-ref}{LinExpr Class}}}}, {\hyperref[\detokenize{pyapiref:chappyapi-quadexpr}]{\sphinxcrossref{\DUrole{std,std-ref}{QuadExpr Class}}}},
{\hyperref[\detokenize{pyapiref:chappyapi-nlexpr}]{\sphinxcrossref{\DUrole{std,std-ref}{NlExpr Class}}}} object, or {\hyperref[\detokenize{pyapiref:chappyapi-exprbuilder}]{\sphinxcrossref{\DUrole{std,std-ref}{ExprBuilder Class}}}} object.
\end{quote}
\end{quote}

\sphinxAtStartPar
\sphinxstylestrong{Return Value}
\begin{quote}

\sphinxAtStartPar
Returns a {\hyperref[\detokenize{pyapiref:chappyapi-nlexpr}]{\sphinxcrossref{\DUrole{std,std-ref}{NlExpr Class}}}} object.
\end{quote}
\end{quote}

\subsubsection{nl.asinh()}
\label{\detokenize{pyapiref:nl-asinh}}\begin{quote}

\sphinxAtStartPar
\sphinxstylestrong{Synopsis}
\begin{quote}

\sphinxAtStartPar
\sphinxcode{\sphinxupquote{NlExpr asinh(expr)}}
\end{quote}

\sphinxAtStartPar
\sphinxstylestrong{Description}
\begin{quote}

\sphinxAtStartPar
Calculate the inverse hyperbolic sine of the input object.
\end{quote}

\sphinxAtStartPar
\sphinxstylestrong{Arguments}
\begin{quote}

\sphinxAtStartPar
\sphinxcode{\sphinxupquote{expr}}
\begin{quote}

\sphinxAtStartPar
The object to be calculated.

\sphinxAtStartPar
Possible values include: constant, {\hyperref[\detokenize{pyapiref:chappyapi-var}]{\sphinxcrossref{\DUrole{std,std-ref}{Var Class}}}} object,
{\hyperref[\detokenize{pyapiref:chappyapi-linexpr}]{\sphinxcrossref{\DUrole{std,std-ref}{LinExpr Class}}}}, {\hyperref[\detokenize{pyapiref:chappyapi-quadexpr}]{\sphinxcrossref{\DUrole{std,std-ref}{QuadExpr Class}}}},
{\hyperref[\detokenize{pyapiref:chappyapi-nlexpr}]{\sphinxcrossref{\DUrole{std,std-ref}{NlExpr Class}}}} object, or {\hyperref[\detokenize{pyapiref:chappyapi-exprbuilder}]{\sphinxcrossref{\DUrole{std,std-ref}{ExprBuilder Class}}}} object.
\end{quote}
\end{quote}

\sphinxAtStartPar
\sphinxstylestrong{Return Value}
\begin{quote}

\sphinxAtStartPar
Returns a {\hyperref[\detokenize{pyapiref:chappyapi-nlexpr}]{\sphinxcrossref{\DUrole{std,std-ref}{NlExpr Class}}}} object.
\end{quote}
\end{quote}

\subsubsection{nl.atan()}
\label{\detokenize{pyapiref:nl-atan}}\begin{quote}

\sphinxAtStartPar
\sphinxstylestrong{Synopsis}
\begin{quote}

\sphinxAtStartPar
\sphinxcode{\sphinxupquote{NlExpr atan(expr)}}
\end{quote}

\sphinxAtStartPar
\sphinxstylestrong{Description}
\begin{quote}

\sphinxAtStartPar
Calculate the arc tangent of the input object.
\end{quote}

\sphinxAtStartPar
\sphinxstylestrong{Arguments}
\begin{quote}

\sphinxAtStartPar
\sphinxcode{\sphinxupquote{expr}}
\begin{quote}

\sphinxAtStartPar
The object to be calculated.

\sphinxAtStartPar
Possible values include: constant, {\hyperref[\detokenize{pyapiref:chappyapi-var}]{\sphinxcrossref{\DUrole{std,std-ref}{Var Class}}}} object,
{\hyperref[\detokenize{pyapiref:chappyapi-linexpr}]{\sphinxcrossref{\DUrole{std,std-ref}{LinExpr Class}}}}, {\hyperref[\detokenize{pyapiref:chappyapi-quadexpr}]{\sphinxcrossref{\DUrole{std,std-ref}{QuadExpr Class}}}},
{\hyperref[\detokenize{pyapiref:chappyapi-nlexpr}]{\sphinxcrossref{\DUrole{std,std-ref}{NlExpr Class}}}} object, or {\hyperref[\detokenize{pyapiref:chappyapi-exprbuilder}]{\sphinxcrossref{\DUrole{std,std-ref}{ExprBuilder Class}}}} object.
\end{quote}
\end{quote}

\sphinxAtStartPar
\sphinxstylestrong{Return Value}
\begin{quote}

\sphinxAtStartPar
Returns a {\hyperref[\detokenize{pyapiref:chappyapi-nlexpr}]{\sphinxcrossref{\DUrole{std,std-ref}{NlExpr Class}}}} object.
\end{quote}
\end{quote}

\subsubsection{nl.atan2()}
\label{\detokenize{pyapiref:nl-atan2}}\begin{quote}

\sphinxAtStartPar
\sphinxstylestrong{Synopsis}
\begin{quote}

\sphinxAtStartPar
\sphinxcode{\sphinxupquote{NlExpr atan2(y, x)}}
\end{quote}

\sphinxAtStartPar
\sphinxstylestrong{Description}
\begin{quote}

\sphinxAtStartPar
Calculate the arc tangent of \sphinxcode{\sphinxupquote{y/x}}, considering the quadrant of the point (x, y).
\end{quote}

\sphinxAtStartPar
\sphinxstylestrong{Arguments}
\begin{quote}

\sphinxAtStartPar
\sphinxcode{\sphinxupquote{y}}
\begin{quote}

\sphinxAtStartPar
The numerator (y\sphinxhyphen{}coordinate) object.

\sphinxAtStartPar
Possible values include: constant, {\hyperref[\detokenize{pyapiref:chappyapi-var}]{\sphinxcrossref{\DUrole{std,std-ref}{Var Class}}}} object,
{\hyperref[\detokenize{pyapiref:chappyapi-linexpr}]{\sphinxcrossref{\DUrole{std,std-ref}{LinExpr Class}}}}, {\hyperref[\detokenize{pyapiref:chappyapi-quadexpr}]{\sphinxcrossref{\DUrole{std,std-ref}{QuadExpr Class}}}},
{\hyperref[\detokenize{pyapiref:chappyapi-nlexpr}]{\sphinxcrossref{\DUrole{std,std-ref}{NlExpr Class}}}} object, or {\hyperref[\detokenize{pyapiref:chappyapi-exprbuilder}]{\sphinxcrossref{\DUrole{std,std-ref}{ExprBuilder Class}}}} object.
\end{quote}

\sphinxAtStartPar
\sphinxcode{\sphinxupquote{x}}
\begin{quote}

\sphinxAtStartPar
The denominator (x\sphinxhyphen{}coordinate) object.

\sphinxAtStartPar
Possible values include: constant, {\hyperref[\detokenize{pyapiref:chappyapi-var}]{\sphinxcrossref{\DUrole{std,std-ref}{Var Class}}}} object,
{\hyperref[\detokenize{pyapiref:chappyapi-linexpr}]{\sphinxcrossref{\DUrole{std,std-ref}{LinExpr Class}}}}, {\hyperref[\detokenize{pyapiref:chappyapi-quadexpr}]{\sphinxcrossref{\DUrole{std,std-ref}{QuadExpr Class}}}},
{\hyperref[\detokenize{pyapiref:chappyapi-nlexpr}]{\sphinxcrossref{\DUrole{std,std-ref}{NlExpr Class}}}} object, or {\hyperref[\detokenize{pyapiref:chappyapi-exprbuilder}]{\sphinxcrossref{\DUrole{std,std-ref}{ExprBuilder Class}}}} object.
\end{quote}
\end{quote}

\sphinxAtStartPar
\sphinxstylestrong{Return Value}
\begin{quote}

\sphinxAtStartPar
Returns a {\hyperref[\detokenize{pyapiref:chappyapi-nlexpr}]{\sphinxcrossref{\DUrole{std,std-ref}{NlExpr Class}}}} object.
\end{quote}
\end{quote}

\subsubsection{nl.atanh()}
\label{\detokenize{pyapiref:nl-atanh}}\begin{quote}

\sphinxAtStartPar
\sphinxstylestrong{Synopsis}
\begin{quote}

\sphinxAtStartPar
\sphinxcode{\sphinxupquote{NlExpr atanh(expr)}}
\end{quote}

\sphinxAtStartPar
\sphinxstylestrong{Description}
\begin{quote}

\sphinxAtStartPar
Calculate the inverse hyperbolic tangent of the input object.
\end{quote}

\sphinxAtStartPar
\sphinxstylestrong{Arguments}
\begin{quote}

\sphinxAtStartPar
\sphinxcode{\sphinxupquote{expr}}
\begin{quote}

\sphinxAtStartPar
The object to be calculated.

\sphinxAtStartPar
Possible values include: constant, {\hyperref[\detokenize{pyapiref:chappyapi-var}]{\sphinxcrossref{\DUrole{std,std-ref}{Var Class}}}} object,
{\hyperref[\detokenize{pyapiref:chappyapi-linexpr}]{\sphinxcrossref{\DUrole{std,std-ref}{LinExpr Class}}}}, {\hyperref[\detokenize{pyapiref:chappyapi-quadexpr}]{\sphinxcrossref{\DUrole{std,std-ref}{QuadExpr Class}}}},
{\hyperref[\detokenize{pyapiref:chappyapi-nlexpr}]{\sphinxcrossref{\DUrole{std,std-ref}{NlExpr Class}}}} object, or {\hyperref[\detokenize{pyapiref:chappyapi-exprbuilder}]{\sphinxcrossref{\DUrole{std,std-ref}{ExprBuilder Class}}}} object.
\end{quote}
\end{quote}

\sphinxAtStartPar
\sphinxstylestrong{Return Value}
\begin{quote}

\sphinxAtStartPar
Returns a {\hyperref[\detokenize{pyapiref:chappyapi-nlexpr}]{\sphinxcrossref{\DUrole{std,std-ref}{NlExpr Class}}}} object.
\end{quote}
\end{quote}

\subsubsection{nl.neg()}
\label{\detokenize{pyapiref:nl-neg}}\begin{quote}

\sphinxAtStartPar
\sphinxstylestrong{Synopsis}
\begin{quote}

\sphinxAtStartPar
\sphinxcode{\sphinxupquote{NlExpr neg(expr)}}
\end{quote}

\sphinxAtStartPar
\sphinxstylestrong{Description}
\begin{quote}

\sphinxAtStartPar
Calculate the negation of the input object.
\end{quote}

\sphinxAtStartPar
\sphinxstylestrong{Arguments}
\begin{quote}

\sphinxAtStartPar
\sphinxcode{\sphinxupquote{expr}}
\begin{quote}

\sphinxAtStartPar
The object to be calculated.

\sphinxAtStartPar
Possible values include: constant, {\hyperref[\detokenize{pyapiref:chappyapi-var}]{\sphinxcrossref{\DUrole{std,std-ref}{Var Class}}}} object,
{\hyperref[\detokenize{pyapiref:chappyapi-linexpr}]{\sphinxcrossref{\DUrole{std,std-ref}{LinExpr Class}}}}, {\hyperref[\detokenize{pyapiref:chappyapi-quadexpr}]{\sphinxcrossref{\DUrole{std,std-ref}{QuadExpr Class}}}},
{\hyperref[\detokenize{pyapiref:chappyapi-nlexpr}]{\sphinxcrossref{\DUrole{std,std-ref}{NlExpr Class}}}} object, or {\hyperref[\detokenize{pyapiref:chappyapi-exprbuilder}]{\sphinxcrossref{\DUrole{std,std-ref}{ExprBuilder Class}}}} object.
\end{quote}
\end{quote}

\sphinxAtStartPar
\sphinxstylestrong{Return Value}
\begin{quote}

\sphinxAtStartPar
Returns a {\hyperref[\detokenize{pyapiref:chappyapi-nlexpr}]{\sphinxcrossref{\DUrole{std,std-ref}{NlExpr Class}}}} object.
\end{quote}
\end{quote}

\subsubsection{nl.pow()}
\label{\detokenize{pyapiref:nl-pow}}\begin{quote}

\sphinxAtStartPar
\sphinxstylestrong{Synopsis}
\begin{quote}

\sphinxAtStartPar
\sphinxcode{\sphinxupquote{NlExpr pow(base, expo)}}
\end{quote}

\sphinxAtStartPar
\sphinxstylestrong{Description}
\begin{quote}

\sphinxAtStartPar
Calculate the exponentiation of the base raised to the power of the exponent.
\end{quote}

\sphinxAtStartPar
\sphinxstylestrong{Arguments}
\begin{quote}

\sphinxAtStartPar
\sphinxcode{\sphinxupquote{base}}
\begin{quote}

\sphinxAtStartPar
The expression object to be used as the base.

\sphinxAtStartPar
Possible values include: constant, {\hyperref[\detokenize{pyapiref:chappyapi-var}]{\sphinxcrossref{\DUrole{std,std-ref}{Var Class}}}} object,
{\hyperref[\detokenize{pyapiref:chappyapi-linexpr}]{\sphinxcrossref{\DUrole{std,std-ref}{LinExpr Class}}}}, {\hyperref[\detokenize{pyapiref:chappyapi-quadexpr}]{\sphinxcrossref{\DUrole{std,std-ref}{QuadExpr Class}}}},
{\hyperref[\detokenize{pyapiref:chappyapi-nlexpr}]{\sphinxcrossref{\DUrole{std,std-ref}{NlExpr Class}}}} object, or {\hyperref[\detokenize{pyapiref:chappyapi-exprbuilder}]{\sphinxcrossref{\DUrole{std,std-ref}{ExprBuilder Class}}}} object.
\end{quote}

\sphinxAtStartPar
\sphinxcode{\sphinxupquote{expo}}
\begin{quote}

\sphinxAtStartPar
The expression object to be used as the exponent.

\sphinxAtStartPar
Possible values include: constant, {\hyperref[\detokenize{pyapiref:chappyapi-var}]{\sphinxcrossref{\DUrole{std,std-ref}{Var Class}}}} object,
{\hyperref[\detokenize{pyapiref:chappyapi-linexpr}]{\sphinxcrossref{\DUrole{std,std-ref}{LinExpr Class}}}}, {\hyperref[\detokenize{pyapiref:chappyapi-quadexpr}]{\sphinxcrossref{\DUrole{std,std-ref}{QuadExpr Class}}}},
{\hyperref[\detokenize{pyapiref:chappyapi-nlexpr}]{\sphinxcrossref{\DUrole{std,std-ref}{NlExpr Class}}}} object, or {\hyperref[\detokenize{pyapiref:chappyapi-exprbuilder}]{\sphinxcrossref{\DUrole{std,std-ref}{ExprBuilder Class}}}} object.
\end{quote}
\end{quote}

\sphinxAtStartPar
\sphinxstylestrong{Return Value}
\begin{quote}

\sphinxAtStartPar
Returns a {\hyperref[\detokenize{pyapiref:chappyapi-nlexpr}]{\sphinxcrossref{\DUrole{std,std-ref}{NlExpr Class}}}} object.
\end{quote}
\end{quote}

\subsubsection{nl.sin()}
\label{\detokenize{pyapiref:nl-sin}}\begin{quote}

\sphinxAtStartPar
\sphinxstylestrong{Synopsis}
\begin{quote}

\sphinxAtStartPar
\sphinxcode{\sphinxupquote{NlExpr sin(expr)}}
\end{quote}

\sphinxAtStartPar
\sphinxstylestrong{Description}
\begin{quote}

\sphinxAtStartPar
Calculate the sine of the input object.
\end{quote}

\sphinxAtStartPar
\sphinxstylestrong{Arguments}
\begin{quote}

\sphinxAtStartPar
\sphinxcode{\sphinxupquote{expr}}
\begin{quote}

\sphinxAtStartPar
The object to be calculated.

\sphinxAtStartPar
Possible values include: constant, {\hyperref[\detokenize{pyapiref:chappyapi-var}]{\sphinxcrossref{\DUrole{std,std-ref}{Var Class}}}} object,
{\hyperref[\detokenize{pyapiref:chappyapi-linexpr}]{\sphinxcrossref{\DUrole{std,std-ref}{LinExpr Class}}}}, {\hyperref[\detokenize{pyapiref:chappyapi-quadexpr}]{\sphinxcrossref{\DUrole{std,std-ref}{QuadExpr Class}}}},
{\hyperref[\detokenize{pyapiref:chappyapi-nlexpr}]{\sphinxcrossref{\DUrole{std,std-ref}{NlExpr Class}}}} object, or {\hyperref[\detokenize{pyapiref:chappyapi-exprbuilder}]{\sphinxcrossref{\DUrole{std,std-ref}{ExprBuilder Class}}}} object.
\end{quote}
\end{quote}

\sphinxAtStartPar
\sphinxstylestrong{Return Value}
\begin{quote}

\sphinxAtStartPar
Returns a {\hyperref[\detokenize{pyapiref:chappyapi-nlexpr}]{\sphinxcrossref{\DUrole{std,std-ref}{NlExpr Class}}}} object.
\end{quote}
\end{quote}

\subsubsection{nl.sinh()}
\label{\detokenize{pyapiref:nl-sinh}}\begin{quote}

\sphinxAtStartPar
\sphinxstylestrong{Synopsis}
\begin{quote}

\sphinxAtStartPar
\sphinxcode{\sphinxupquote{NlExpr sinh(expr)}}
\end{quote}

\sphinxAtStartPar
\sphinxstylestrong{Description}
\begin{quote}

\sphinxAtStartPar
Calculate the hyperbolic sine of the input object.
\end{quote}

\sphinxAtStartPar
\sphinxstylestrong{Arguments}
\begin{quote}

\sphinxAtStartPar
\sphinxcode{\sphinxupquote{expr}}
\begin{quote}

\sphinxAtStartPar
The object to be calculated.

\sphinxAtStartPar
Possible values include: constant, {\hyperref[\detokenize{pyapiref:chappyapi-var}]{\sphinxcrossref{\DUrole{std,std-ref}{Var Class}}}} object,
{\hyperref[\detokenize{pyapiref:chappyapi-linexpr}]{\sphinxcrossref{\DUrole{std,std-ref}{LinExpr Class}}}}, {\hyperref[\detokenize{pyapiref:chappyapi-quadexpr}]{\sphinxcrossref{\DUrole{std,std-ref}{QuadExpr Class}}}},
{\hyperref[\detokenize{pyapiref:chappyapi-nlexpr}]{\sphinxcrossref{\DUrole{std,std-ref}{NlExpr Class}}}} object, or {\hyperref[\detokenize{pyapiref:chappyapi-exprbuilder}]{\sphinxcrossref{\DUrole{std,std-ref}{ExprBuilder Class}}}} object.
\end{quote}
\end{quote}

\sphinxAtStartPar
\sphinxstylestrong{Return Value}
\begin{quote}

\sphinxAtStartPar
Returns a {\hyperref[\detokenize{pyapiref:chappyapi-nlexpr}]{\sphinxcrossref{\DUrole{std,std-ref}{NlExpr Class}}}} object.
\end{quote}
\end{quote}

\subsubsection{nl.sqrt()}
\label{\detokenize{pyapiref:nl-sqrt}}\begin{quote}

\sphinxAtStartPar
\sphinxstylestrong{Synopsis}
\begin{quote}

\sphinxAtStartPar
\sphinxcode{\sphinxupquote{NlExpr sqrt(expr)}}
\end{quote}

\sphinxAtStartPar
\sphinxstylestrong{Description}
\begin{quote}

\sphinxAtStartPar
Calculate the square root of the input object.
\end{quote}

\sphinxAtStartPar
\sphinxstylestrong{Arguments}
\begin{quote}

\sphinxAtStartPar
\sphinxcode{\sphinxupquote{expr}}
\begin{quote}

\sphinxAtStartPar
The object to be calculated.

\sphinxAtStartPar
Possible values include: constant, {\hyperref[\detokenize{pyapiref:chappyapi-var}]{\sphinxcrossref{\DUrole{std,std-ref}{Var Class}}}} object,
{\hyperref[\detokenize{pyapiref:chappyapi-linexpr}]{\sphinxcrossref{\DUrole{std,std-ref}{LinExpr Class}}}}, {\hyperref[\detokenize{pyapiref:chappyapi-quadexpr}]{\sphinxcrossref{\DUrole{std,std-ref}{QuadExpr Class}}}},
{\hyperref[\detokenize{pyapiref:chappyapi-nlexpr}]{\sphinxcrossref{\DUrole{std,std-ref}{NlExpr Class}}}} object, or {\hyperref[\detokenize{pyapiref:chappyapi-exprbuilder}]{\sphinxcrossref{\DUrole{std,std-ref}{ExprBuilder Class}}}} object.
\end{quote}
\end{quote}

\sphinxAtStartPar
\sphinxstylestrong{Return Value}
\begin{quote}

\sphinxAtStartPar
Returns a {\hyperref[\detokenize{pyapiref:chappyapi-nlexpr}]{\sphinxcrossref{\DUrole{std,std-ref}{NlExpr Class}}}} object.
\end{quote}
\end{quote}

\subsubsection{nl.sum()}
\label{\detokenize{pyapiref:nl-sum}}\begin{quote}

\sphinxAtStartPar
\sphinxstylestrong{Synopsis}
\begin{quote}

\sphinxAtStartPar
\sphinxcode{\sphinxupquote{NlExpr sum(op1, op2=None, op3=None, op4=None)}}
\end{quote}

\sphinxAtStartPar
\sphinxstylestrong{Description}
\begin{quote}

\sphinxAtStartPar
Calculate the sum of nonlinear expression objects.
\end{quote}

\sphinxAtStartPar
\sphinxstylestrong{Arguments}
\begin{quote}

\sphinxAtStartPar
\sphinxcode{\sphinxupquote{op1}}
\begin{quote}

\sphinxAtStartPar
Required argument.

\sphinxAtStartPar
The expression object(s) to be summed.

\sphinxAtStartPar
Possible values include a {\hyperref[\detokenize{pyapiref:chappyapi-nlexpr}]{\sphinxcrossref{\DUrole{std,std-ref}{NlExpr Class}}}} object or an iterable/mapping of {\hyperref[\detokenize{pyapiref:chappyapi-nlexpr}]{\sphinxcrossref{\DUrole{std,std-ref}{NlExpr Class}}}} objects.
\end{quote}

\sphinxAtStartPar
\sphinxcode{\sphinxupquote{op2}}, \sphinxcode{\sphinxupquote{op3}}, \sphinxcode{\sphinxupquote{op4}}
\begin{quote}

\sphinxAtStartPar
Optional arguments.

\sphinxAtStartPar
Additional expression objects to be summed.

\sphinxAtStartPar
Possible values include: constant, {\hyperref[\detokenize{pyapiref:chappyapi-var}]{\sphinxcrossref{\DUrole{std,std-ref}{Var Class}}}} object,
{\hyperref[\detokenize{pyapiref:chappyapi-linexpr}]{\sphinxcrossref{\DUrole{std,std-ref}{LinExpr Class}}}}, {\hyperref[\detokenize{pyapiref:chappyapi-quadexpr}]{\sphinxcrossref{\DUrole{std,std-ref}{QuadExpr Class}}}},
{\hyperref[\detokenize{pyapiref:chappyapi-nlexpr}]{\sphinxcrossref{\DUrole{std,std-ref}{NlExpr Class}}}} object, or {\hyperref[\detokenize{pyapiref:chappyapi-exprbuilder}]{\sphinxcrossref{\DUrole{std,std-ref}{ExprBuilder Class}}}} object.
\end{quote}
\end{quote}

\sphinxAtStartPar
\sphinxstylestrong{Return Value}
\begin{quote}

\sphinxAtStartPar
Returns a {\hyperref[\detokenize{pyapiref:chappyapi-nlexpr}]{\sphinxcrossref{\DUrole{std,std-ref}{NlExpr Class}}}} object.
\end{quote}
\end{quote}

\subsubsection{nl.tan()}
\label{\detokenize{pyapiref:nl-tan}}\begin{quote}

\sphinxAtStartPar
\sphinxstylestrong{Synopsis}
\begin{quote}

\sphinxAtStartPar
\sphinxcode{\sphinxupquote{NlExpr tan(expr)}}
\end{quote}

\sphinxAtStartPar
\sphinxstylestrong{Description}
\begin{quote}

\sphinxAtStartPar
Calculate the tangent of the input object.
\end{quote}

\sphinxAtStartPar
\sphinxstylestrong{Arguments}
\begin{quote}

\sphinxAtStartPar
\sphinxcode{\sphinxupquote{expr}}
\begin{quote}

\sphinxAtStartPar
The object to be calculated.

\sphinxAtStartPar
Possible values include: constant, {\hyperref[\detokenize{pyapiref:chappyapi-var}]{\sphinxcrossref{\DUrole{std,std-ref}{Var Class}}}} object,
{\hyperref[\detokenize{pyapiref:chappyapi-linexpr}]{\sphinxcrossref{\DUrole{std,std-ref}{LinExpr Class}}}}, {\hyperref[\detokenize{pyapiref:chappyapi-quadexpr}]{\sphinxcrossref{\DUrole{std,std-ref}{QuadExpr Class}}}},
{\hyperref[\detokenize{pyapiref:chappyapi-nlexpr}]{\sphinxcrossref{\DUrole{std,std-ref}{NlExpr Class}}}} object, or {\hyperref[\detokenize{pyapiref:chappyapi-exprbuilder}]{\sphinxcrossref{\DUrole{std,std-ref}{ExprBuilder Class}}}} object.
\end{quote}
\end{quote}

\sphinxAtStartPar
\sphinxstylestrong{Return Value}
\begin{quote}

\sphinxAtStartPar
Returns a {\hyperref[\detokenize{pyapiref:chappyapi-nlexpr}]{\sphinxcrossref{\DUrole{std,std-ref}{NlExpr Class}}}} object.
\end{quote}
\end{quote}

\subsubsection{nl.tanh()}
\label{\detokenize{pyapiref:nl-tanh}}\begin{quote}

\sphinxAtStartPar
\sphinxstylestrong{Synopsis}
\begin{quote}

\sphinxAtStartPar
\sphinxcode{\sphinxupquote{NlExpr tanh(expr)}}
\end{quote}

\sphinxAtStartPar
\sphinxstylestrong{Description}
\begin{quote}

\sphinxAtStartPar
Calculate the hyperbolic tangent of the input object.
\end{quote}

\sphinxAtStartPar
\sphinxstylestrong{Arguments}
\begin{quote}

\sphinxAtStartPar
\sphinxcode{\sphinxupquote{expr}}
\begin{quote}

\sphinxAtStartPar
The expression object to be calculated.

\sphinxAtStartPar
Possible values include: constant, {\hyperref[\detokenize{pyapiref:chappyapi-var}]{\sphinxcrossref{\DUrole{std,std-ref}{Var Class}}}} object,
{\hyperref[\detokenize{pyapiref:chappyapi-linexpr}]{\sphinxcrossref{\DUrole{std,std-ref}{LinExpr Class}}}}, {\hyperref[\detokenize{pyapiref:chappyapi-quadexpr}]{\sphinxcrossref{\DUrole{std,std-ref}{QuadExpr Class}}}},
{\hyperref[\detokenize{pyapiref:chappyapi-nlexpr}]{\sphinxcrossref{\DUrole{std,std-ref}{NlExpr Class}}}} object, or {\hyperref[\detokenize{pyapiref:chappyapi-exprbuilder}]{\sphinxcrossref{\DUrole{std,std-ref}{ExprBuilder Class}}}} object.
\end{quote}
\end{quote}

\sphinxAtStartPar
\sphinxstylestrong{Return Value}
\begin{quote}

\sphinxAtStartPar
Returns a {\hyperref[\detokenize{pyapiref:chappyapi-nlexpr}]{\sphinxcrossref{\DUrole{std,std-ref}{NlExpr Class}}}} object.
\end{quote}
\end{quote}

\subsection{CallbackBase Class}
\label{\detokenize{pyapiref:callbackbase-class}}\label{\detokenize{pyapiref:chappyapi-cbcbase}}
\sphinxAtStartPar
COPT CallbackBase class. This is an abstract class, the user needs to implement the function {\hyperref[\detokenize{pyapiref:chappyapi-cbcbase-callback}]{\sphinxcrossref{\DUrole{std,std-ref}{CallbackBase.callback()}}}} to create an instance.
The instance is passed in as the first argument of the method {\hyperref[\detokenize{pyapiref:chappyapi-model-setcallback}]{\sphinxcrossref{\DUrole{std,std-ref}{Model.setCallback()}}}}

\subsubsection{CallbackBase.where()}
\label{\detokenize{pyapiref:callbackbase-where}}\begin{quote}

\sphinxAtStartPar
\sphinxstylestrong{Synopsis}
\begin{quote}

\sphinxAtStartPar
\sphinxcode{\sphinxupquote{where()}}
\end{quote}

\sphinxAtStartPar
\sphinxstylestrong{Description}
\begin{quote}

\sphinxAtStartPar
Get context in callback.
\end{quote}

\sphinxAtStartPar
\sphinxstylestrong{Return value}
\begin{quote}

\sphinxAtStartPar
Returns an integer value.
\end{quote}
\end{quote}

\subsubsection{CallbackBase.callback()}
\label{\detokenize{pyapiref:callbackbase-callback}}\label{\detokenize{pyapiref:chappyapi-cbcbase-callback}}\begin{quote}

\sphinxAtStartPar
\sphinxstylestrong{Synopsis}
\begin{quote}

\sphinxAtStartPar
\sphinxcode{\sphinxupquote{callback()}}
\end{quote}

\sphinxAtStartPar
\sphinxstylestrong{Description}
\begin{quote}

\sphinxAtStartPar
The callback function is a pure virtual function which needs to be implemented by the user.
The user can describe the information that needs to be obtained or the operation that needs to be performed during the solution process.
\end{quote}

\sphinxAtStartPar
\sphinxstylestrong{Example}

\begin{sphinxVerbatim}[commandchars=\\\{\}]
\PYG{k}{class}\PYG{+w}{ }\PYG{n+nc}{CoptCallback}\PYG{p}{(}\PYG{n}{CallbackBase}\PYG{p}{)}\PYG{p}{:}
    \PYG{k}{def}\PYG{+w}{ }\PYG{n+nf+fm}{\PYGZus{}\PYGZus{}init\PYGZus{}\PYGZus{}}\PYG{p}{(}\PYG{n+nb+bp}{self}\PYG{p}{)}\PYG{p}{:}
        \PYG{n+nb}{super}\PYG{p}{(}\PYG{p}{)}\PYG{o}{.}\PYG{n+nf+fm}{\PYGZus{}\PYGZus{}init\PYGZus{}\PYGZus{}}\PYG{p}{(}\PYG{p}{)}
    \PYG{k}{def}\PYG{+w}{ }\PYG{n+nf}{callback}\PYG{p}{(}\PYG{n+nb+bp}{self}\PYG{p}{)}\PYG{p}{:}
        \PYG{c+c1}{\PYGZsh{} Get the objective value when finding a feasible MIP solution}
        \PYG{k}{if} \PYG{n+nb+bp}{self}\PYG{o}{.}\PYG{n}{where}\PYG{p}{(}\PYG{p}{)} \PYG{o}{==} \PYG{n}{COPT}\PYG{o}{.}\PYG{n}{CBCONTEXT\PYGZus{}MIPSOL}\PYG{p}{:}
            \PYG{n}{db} \PYG{o}{=} \PYG{n+nb+bp}{self}\PYG{o}{.}\PYG{n}{getInfo}\PYG{p}{(}\PYG{n}{COPT}\PYG{o}{.}\PYG{n}{CBInfo}\PYG{o}{.}\PYG{n}{MipCandObj}\PYG{p}{)}
\end{sphinxVerbatim}
\end{quote}

\subsubsection{CallbackBase.interrupt()}
\label{\detokenize{pyapiref:callbackbase-interrupt}}\begin{quote}

\sphinxAtStartPar
\sphinxstylestrong{Synopsis}
\begin{quote}

\sphinxAtStartPar
\sphinxcode{\sphinxupquote{interrupt()}}
\end{quote}

\sphinxAtStartPar
\sphinxstylestrong{Description}
\begin{quote}

\sphinxAtStartPar
Interrupt the callback process.
\end{quote}
\end{quote}

\subsubsection{CallbackBase.addUserCut()}
\label{\detokenize{pyapiref:callbackbase-addusercut}}\begin{quote}

\sphinxAtStartPar
\sphinxstylestrong{Synopsis}
\begin{quote}

\sphinxAtStartPar
\sphinxcode{\sphinxupquote{addUserCut(lhs, sense = None, rhs = None)}}
\end{quote}

\sphinxAtStartPar
\sphinxstylestrong{Description}
\begin{quote}

\sphinxAtStartPar
Add a user cut to the MIP model from within the callback function.
\end{quote}

\sphinxAtStartPar
\sphinxstylestrong{Arguments}
\begin{quote}

\sphinxAtStartPar
\sphinxcode{\sphinxupquote{lhs}}
\begin{quote}

\sphinxAtStartPar
Left\sphinxhyphen{}hand side expression for the new user cut. It can take the value of {\hyperref[\detokenize{pyapiref:chappyapi-var}]{\sphinxcrossref{\DUrole{std,std-ref}{Var Class}}}} object, {\hyperref[\detokenize{pyapiref:chappyapi-linexpr}]{\sphinxcrossref{\DUrole{std,std-ref}{LinExpr Class}}}} object, or {\hyperref[\detokenize{pyapiref:chappyapi-constrbuilder}]{\sphinxcrossref{\DUrole{std,std-ref}{ConstrBuilder Class}}}} .
\end{quote}

\sphinxAtStartPar
\sphinxcode{\sphinxupquote{sense}}
\begin{quote}

\sphinxAtStartPar
The sense of the new user cut. It supports for \sphinxcode{\sphinxupquote{LESS\_EQUAL}}, \sphinxcode{\sphinxupquote{GREATER\_EQUAL}}, \sphinxcode{\sphinxupquote{EQUAL}} and \sphinxcode{\sphinxupquote{FREE}} .

\sphinxAtStartPar
Optional. None by default.

\sphinxAtStartPar
The user cut added from within callback can only have a single comparison operator.
\end{quote}

\sphinxAtStartPar
\sphinxcode{\sphinxupquote{rhs}}
\begin{quote}

\sphinxAtStartPar
Right\sphinxhyphen{}hand side expression for the new user cut.

\sphinxAtStartPar
Optional. None by default.

\sphinxAtStartPar
It can be a constant, or {\hyperref[\detokenize{pyapiref:chappyapi-var}]{\sphinxcrossref{\DUrole{std,std-ref}{Var Class}}}} object, or {\hyperref[\detokenize{pyapiref:chappyapi-linexpr}]{\sphinxcrossref{\DUrole{std,std-ref}{LinExpr Class}}}} object.
\end{quote}
\end{quote}

\sphinxAtStartPar
\sphinxstylestrong{Example}

\begin{sphinxVerbatim}[commandchars=\\\{\}]
\PYG{n+nb+bp}{self}\PYG{o}{.}\PYG{n}{addUserCut}\PYG{p}{(}\PYG{n}{x}\PYG{o}{+}\PYG{n}{y} \PYG{o}{\PYGZlt{}}\PYG{o}{=} \PYG{l+m+mi}{1}\PYG{p}{)}
\end{sphinxVerbatim}
\end{quote}

\subsubsection{CallbackBase.addUserCuts()}
\label{\detokenize{pyapiref:callbackbase-addusercuts}}\begin{quote}

\sphinxAtStartPar
\sphinxstylestrong{Synopsis}
\begin{quote}

\sphinxAtStartPar
\sphinxcode{\sphinxupquote{addUserCuts(generator)}}
\end{quote}

\sphinxAtStartPar
\sphinxstylestrong{Description}
\begin{quote}

\sphinxAtStartPar
Add a set of user cuts to the MIP model from within the callback function.
\end{quote}

\sphinxAtStartPar
\sphinxstylestrong{Arguments}
\begin{quote}

\sphinxAtStartPar
\sphinxcode{\sphinxupquote{generator}}
\begin{quote}

\sphinxAtStartPar
Array of builders for user cuts. It can be {\hyperref[\detokenize{pyapiref:chappyapi-constrbuilderarray}]{\sphinxcrossref{\DUrole{std,std-ref}{ConstrBuilderArray Class}}}} object or {\hyperref[\detokenize{pyapiref:chappyapi-mconstrbuilder}]{\sphinxcrossref{\DUrole{std,std-ref}{MConstrBuilder Class}}}} object.
\end{quote}
\end{quote}

\sphinxAtStartPar
\sphinxstylestrong{Example}

\begin{sphinxVerbatim}[commandchars=\\\{\}]
\PYG{n+nb+bp}{self}\PYG{o}{.}\PYG{n}{addUserCuts}\PYG{p}{(}\PYG{n}{x}\PYG{p}{[}\PYG{n}{i}\PYG{p}{]}\PYG{o}{+}\PYG{n}{y}\PYG{p}{[}\PYG{n}{i}\PYG{p}{]} \PYG{o}{\PYGZlt{}}\PYG{o}{=} \PYG{l+m+mi}{1} \PYG{k}{for} \PYG{n}{i} \PYG{o+ow}{in} \PYG{n+nb}{range}\PYG{p}{(}\PYG{l+m+mi}{10}\PYG{p}{)}\PYG{p}{)}
\end{sphinxVerbatim}
\end{quote}

\subsubsection{CallbackBase.addLazyConstr()}
\label{\detokenize{pyapiref:callbackbase-addlazyconstr}}\begin{quote}

\sphinxAtStartPar
\sphinxstylestrong{Synopsis}
\begin{quote}

\sphinxAtStartPar
\sphinxcode{\sphinxupquote{addLazyConstr(lhs, sense = None, rhs = None)}}
\end{quote}

\sphinxAtStartPar
\sphinxstylestrong{Description}
\begin{quote}

\sphinxAtStartPar
Add a lazy constraint to the MIP model from within the callback function.
\end{quote}

\sphinxAtStartPar
\sphinxstylestrong{Arguments}
\begin{quote}

\sphinxAtStartPar
\sphinxcode{\sphinxupquote{lhs}}
\begin{quote}

\sphinxAtStartPar
Left\sphinxhyphen{}hand side expression for the new lazy constraint. It can take the value of {\hyperref[\detokenize{pyapiref:chappyapi-var}]{\sphinxcrossref{\DUrole{std,std-ref}{Var Class}}}} object, {\hyperref[\detokenize{pyapiref:chappyapi-linexpr}]{\sphinxcrossref{\DUrole{std,std-ref}{LinExpr Class}}}} object or {\hyperref[\detokenize{pyapiref:chappyapi-constrbuilder}]{\sphinxcrossref{\DUrole{std,std-ref}{ConstrBuilder Class}}}} object.
\end{quote}

\sphinxAtStartPar
\sphinxcode{\sphinxupquote{sense}}
\begin{quote}

\sphinxAtStartPar
The sense of the lazy constraint. It supports for \sphinxcode{\sphinxupquote{LESS\_EQUAL}}, \sphinxcode{\sphinxupquote{GREATER\_EQUAL}}, \sphinxcode{\sphinxupquote{EQUAL}} and \sphinxcode{\sphinxupquote{FREE}} .

\sphinxAtStartPar
Optional. None by default.

\sphinxAtStartPar
The lazy constraint added from within callback can only have a single comparison operator.
\end{quote}

\sphinxAtStartPar
\sphinxcode{\sphinxupquote{rhs}}
\begin{quote}

\sphinxAtStartPar
Right\sphinxhyphen{}hand side expression for the new lazy constraint.

\sphinxAtStartPar
Optional. None by default.

\sphinxAtStartPar
It can be a constant, or {\hyperref[\detokenize{pyapiref:chappyapi-var}]{\sphinxcrossref{\DUrole{std,std-ref}{Var Class}}}} object, or {\hyperref[\detokenize{pyapiref:chappyapi-linexpr}]{\sphinxcrossref{\DUrole{std,std-ref}{LinExpr Class}}}} object.
\end{quote}
\end{quote}

\sphinxAtStartPar
\sphinxstylestrong{Example}

\begin{sphinxVerbatim}[commandchars=\\\{\}]
\PYG{n+nb+bp}{self}\PYG{o}{.}\PYG{n}{addLazyConstr}\PYG{p}{(}\PYG{n}{x}\PYG{o}{+}\PYG{n}{y} \PYG{o}{\PYGZlt{}}\PYG{o}{=} \PYG{l+m+mi}{1}\PYG{p}{)}
\end{sphinxVerbatim}
\end{quote}

\subsubsection{CallbackBase.addLazyConstrs()}
\label{\detokenize{pyapiref:callbackbase-addlazyconstrs}}\begin{quote}

\sphinxAtStartPar
\sphinxstylestrong{Synopsis}
\begin{quote}

\sphinxAtStartPar
\sphinxcode{\sphinxupquote{addLazyConstrs(generator)}}
\end{quote}

\sphinxAtStartPar
\sphinxstylestrong{Description}
\begin{quote}

\sphinxAtStartPar
Add a set of lazy constraints to the MIP model from within the callback function.
\end{quote}

\sphinxAtStartPar
\sphinxstylestrong{Arguments}
\begin{quote}

\sphinxAtStartPar
\sphinxcode{\sphinxupquote{generator}}
\begin{quote}

\sphinxAtStartPar
Array of builders for lazy constraints. It can be {\hyperref[\detokenize{pyapiref:chappyapi-constrbuilderarray}]{\sphinxcrossref{\DUrole{std,std-ref}{ConstrBuilderArray Class}}}} object or {\hyperref[\detokenize{pyapiref:chappyapi-mconstrbuilder}]{\sphinxcrossref{\DUrole{std,std-ref}{MConstrBuilder Class}}}} object.
\end{quote}
\end{quote}

\sphinxAtStartPar
\sphinxstylestrong{Example}

\begin{sphinxVerbatim}[commandchars=\\\{\}]
\PYG{n+nb+bp}{self}\PYG{o}{.}\PYG{n}{addLazyConstrs}\PYG{p}{(}\PYG{n}{x}\PYG{p}{[}\PYG{n}{i}\PYG{p}{]} \PYG{o}{+} \PYG{n}{y}\PYG{p}{[}\PYG{n}{i}\PYG{p}{]} \PYG{o}{\PYGZlt{}}\PYG{o}{=} \PYG{l+m+mi}{1} \PYG{k}{for} \PYG{n}{i} \PYG{o+ow}{in} \PYG{n+nb}{range}\PYG{p}{(}\PYG{l+m+mi}{10}\PYG{p}{)}\PYG{p}{)}
\end{sphinxVerbatim}
\end{quote}

\subsubsection{CallbackBase.getInfo()}
\label{\detokenize{pyapiref:callbackbase-getinfo}}\begin{quote}

\sphinxAtStartPar
\sphinxstylestrong{Synopsis}
\begin{quote}

\sphinxAtStartPar
\sphinxcode{\sphinxupquote{getInfo(cbinfo)}}
\end{quote}

\sphinxAtStartPar
\sphinxstylestrong{Description}
\begin{quote}

\sphinxAtStartPar
Retrieve the value of the specified callback information.
\end{quote}

\sphinxAtStartPar
\sphinxstylestrong{Arguments}
\begin{quote}

\sphinxAtStartPar
\sphinxcode{\sphinxupquote{cbinfo}}
\begin{quote}

\sphinxAtStartPar
The name of the callback information. Please refer to {\hyperref[\detokenize{information:chapinfo-cbc}]{\sphinxcrossref{\DUrole{std,std-ref}{Callback Information}}}} for possible values.
\end{quote}
\end{quote}

\sphinxAtStartPar
\sphinxstylestrong{Return value}
\begin{quote}

\sphinxAtStartPar
Returns a constant(int\sphinxhyphen{}valued or double\sphinxhyphen{}valued).
\end{quote}

\sphinxAtStartPar
\sphinxstylestrong{Example}

\begin{sphinxVerbatim}[commandchars=\\\{\}]
\PYG{n}{db} \PYG{o}{=} \PYG{n+nb+bp}{self}\PYG{o}{.}\PYG{n}{getInfo}\PYG{p}{(}\PYG{n}{COPT}\PYG{o}{.}\PYG{n}{CBInfo}\PYG{o}{.}\PYG{n}{BestBnd}\PYG{p}{)}
\end{sphinxVerbatim}
\end{quote}

\subsubsection{CallbackBase.getRelaxSol()}
\label{\detokenize{pyapiref:callbackbase-getrelaxsol}}\begin{quote}

\sphinxAtStartPar
\sphinxstylestrong{Synopsis}
\begin{quote}

\sphinxAtStartPar
\sphinxcode{\sphinxupquote{getRelaxSol(vars)}}
\end{quote}

\sphinxAtStartPar
\sphinxstylestrong{Description}
\begin{quote}

\sphinxAtStartPar
Retrieve the LP\sphinxhyphen{}relaxation solution of the specified variables at the current node.

\sphinxAtStartPar
Note that this method can only be invoked if \sphinxcode{\sphinxupquote{CallbackBase.where() == COPT.CBCONTEXT\_MIPRELAX}}.
\end{quote}

\sphinxAtStartPar
\sphinxstylestrong{Arguments}
\begin{quote}

\sphinxAtStartPar
\sphinxcode{\sphinxupquote{vars}}
\begin{quote}

\sphinxAtStartPar
The variables to retrieve the LP\sphinxhyphen{}relaxation solution values.
\end{quote}
\end{quote}

\sphinxAtStartPar
\sphinxstylestrong{Return value}
\begin{quote}

\sphinxAtStartPar
When parameter \sphinxcode{\sphinxupquote{vars}} is {\hyperref[\detokenize{pyapiref:chappyapi-var}]{\sphinxcrossref{\DUrole{std,std-ref}{Var Class}}}} object, it returns a constant, which is the LP\sphinxhyphen{}relaxation solution value of the specified variable.

\sphinxAtStartPar
When parameter \sphinxcode{\sphinxupquote{vars}} is list or {\hyperref[\detokenize{pyapiref:chappyapi-vararray}]{\sphinxcrossref{\DUrole{std,std-ref}{VarArray Class}}}} object, it returns a list of constants, consisting of the solution of the specified variables.

\sphinxAtStartPar
When parameter \sphinxcode{\sphinxupquote{args}} is dictionary or {\hyperref[\detokenize{pyapiref:chappyapi-util-tupledict}]{\sphinxcrossref{\DUrole{std,std-ref}{tupledict Class}}}} object, it returns {\hyperref[\detokenize{pyapiref:chappyapi-util-tupledict}]{\sphinxcrossref{\DUrole{std,std-ref}{tupledict Class}}}} object(the indices of specified variables as key, the solutions of the specified variables as value).

\sphinxAtStartPar
When parameter \sphinxcode{\sphinxupquote{args}} is \sphinxcode{\sphinxupquote{None}} , it returns the LP\sphinxhyphen{}relaxation solution values of all variables.
\end{quote}

\sphinxAtStartPar
\sphinxstylestrong{Example}

\begin{sphinxVerbatim}[commandchars=\\\{\}]
\PYG{n}{vals} \PYG{o}{=} \PYG{n+nb+bp}{self}\PYG{o}{.}\PYG{n}{getRelaxSol}\PYG{p}{(}\PYG{n+nb}{vars}\PYG{p}{)}
\end{sphinxVerbatim}
\end{quote}

\subsubsection{CallbackBase.getIncumbent()}
\label{\detokenize{pyapiref:callbackbase-getincumbent}}\begin{quote}

\sphinxAtStartPar
\sphinxstylestrong{Synopsis}
\begin{quote}

\sphinxAtStartPar
\sphinxcode{\sphinxupquote{getIncumbent(vars)}}
\end{quote}

\sphinxAtStartPar
\sphinxstylestrong{Description}
\begin{quote}

\sphinxAtStartPar
Retrieve values from the best feasible solution of the specified variables.
\end{quote}

\sphinxAtStartPar
\sphinxstylestrong{Arguments}
\begin{quote}

\sphinxAtStartPar
\sphinxcode{\sphinxupquote{vars}}
\begin{quote}

\sphinxAtStartPar
The variables whose solution values to retrieve.
\end{quote}
\end{quote}

\sphinxAtStartPar
\sphinxstylestrong{Return value}
\begin{quote}

\sphinxAtStartPar
When parameter \sphinxcode{\sphinxupquote{vars}} is {\hyperref[\detokenize{pyapiref:chappyapi-var}]{\sphinxcrossref{\DUrole{std,std-ref}{Var Class}}}} object, it returns a constant, which is the solution value of the specified variable.

\sphinxAtStartPar
When parameter \sphinxcode{\sphinxupquote{vars}} is list or {\hyperref[\detokenize{pyapiref:chappyapi-vararray}]{\sphinxcrossref{\DUrole{std,std-ref}{VarArray Class}}}} object, it returns a list of constants, consisting of the solution of the specified variables.

\sphinxAtStartPar
When parameter \sphinxcode{\sphinxupquote{args}} is dictionary or {\hyperref[\detokenize{pyapiref:chappyapi-util-tupledict}]{\sphinxcrossref{\DUrole{std,std-ref}{tupledict Class}}}} object, it returns {\hyperref[\detokenize{pyapiref:chappyapi-util-tupledict}]{\sphinxcrossref{\DUrole{std,std-ref}{tupledict Class}}}} object(the indices of specified variables as key, the solutions of the specified variables as value).

\sphinxAtStartPar
When parameter \sphinxcode{\sphinxupquote{args}} is \sphinxcode{\sphinxupquote{None}} , it returns the incumbent values of all variables.
\end{quote}

\sphinxAtStartPar
\sphinxstylestrong{Example}

\begin{sphinxVerbatim}[commandchars=\\\{\}]
\PYG{n}{vals} \PYG{o}{=} \PYG{n+nb+bp}{self}\PYG{o}{.}\PYG{n}{getIncumbent}\PYG{p}{(}\PYG{n+nb}{vars}\PYG{p}{)}
\end{sphinxVerbatim}
\end{quote}

\subsubsection{CallbackBase.getSolution()}
\label{\detokenize{pyapiref:callbackbase-getsolution}}\begin{quote}

\sphinxAtStartPar
\sphinxstylestrong{Synopsis}
\begin{quote}

\sphinxAtStartPar
\sphinxcode{\sphinxupquote{getSolution(vars)}}
\end{quote}

\sphinxAtStartPar
\sphinxstylestrong{Description}
\begin{quote}

\sphinxAtStartPar
Retrieve values from the current solution of the specified variables.

\sphinxAtStartPar
Note that this method can only be invoked if \sphinxcode{\sphinxupquote{CallbackBase.where() == COPT.CBCONTEXT\_MIPSOL}}.
\end{quote}

\sphinxAtStartPar
\sphinxstylestrong{Arguments}
\begin{quote}

\sphinxAtStartPar
\sphinxcode{\sphinxupquote{vars}}
\begin{quote}

\sphinxAtStartPar
The variables whose solution values to retrieve.
\end{quote}
\end{quote}

\sphinxAtStartPar
\sphinxstylestrong{Return value}
\begin{quote}

\sphinxAtStartPar
When parameter \sphinxcode{\sphinxupquote{vars}} is {\hyperref[\detokenize{pyapiref:chappyapi-var}]{\sphinxcrossref{\DUrole{std,std-ref}{Var Class}}}} object, it returns a constant, which is the solution value of the specified variable.

\sphinxAtStartPar
When parameter \sphinxcode{\sphinxupquote{vars}} is list or {\hyperref[\detokenize{pyapiref:chappyapi-vararray}]{\sphinxcrossref{\DUrole{std,std-ref}{VarArray Class}}}} object, it returns a list of constants, consisting of the solution of the specified variables.

\sphinxAtStartPar
When parameter \sphinxcode{\sphinxupquote{args}} is dictionary or {\hyperref[\detokenize{pyapiref:chappyapi-util-tupledict}]{\sphinxcrossref{\DUrole{std,std-ref}{tupledict Class}}}} object, it returns {\hyperref[\detokenize{pyapiref:chappyapi-util-tupledict}]{\sphinxcrossref{\DUrole{std,std-ref}{tupledict Class}}}} object(the indices of specified variables as key, the solutions of the specified variables as value).

\sphinxAtStartPar
When parameter \sphinxcode{\sphinxupquote{args}} is \sphinxcode{\sphinxupquote{None}} , it returns the solution values of all variables.
\end{quote}

\sphinxAtStartPar
\sphinxstylestrong{Example}

\begin{sphinxVerbatim}[commandchars=\\\{\}]
\PYG{n}{vals} \PYG{o}{=} \PYG{n+nb+bp}{self}\PYG{o}{.}\PYG{n}{getSolution}\PYG{p}{(}\PYG{n+nb}{vars}\PYG{p}{)}
\end{sphinxVerbatim}
\end{quote}

\subsubsection{CallbackBase.setSolution()}
\label{\detokenize{pyapiref:callbackbase-setsolution}}\begin{quote}

\sphinxAtStartPar
\sphinxstylestrong{Synopsis}
\begin{quote}

\sphinxAtStartPar
\sphinxcode{\sphinxupquote{setSolution(vars, vals)}}
\end{quote}

\sphinxAtStartPar
\sphinxstylestrong{Description}
\begin{quote}

\sphinxAtStartPar
Set feasible solution values for the specified variables.

\sphinxAtStartPar
When parameter \sphinxcode{\sphinxupquote{vars}} is {\hyperref[\detokenize{pyapiref:chappyapi-var}]{\sphinxcrossref{\DUrole{std,std-ref}{Var Class}}}} object,
parameter \sphinxcode{\sphinxupquote{vals}} is constant;

\sphinxAtStartPar
When parameter \sphinxcode{\sphinxupquote{vars}} is dictionary or {\hyperref[\detokenize{pyapiref:chappyapi-util-tupledict}]{\sphinxcrossref{\DUrole{std,std-ref}{tupledict Class}}}} object,
parameter \sphinxcode{\sphinxupquote{vals}} can be constant, dictionary or {\hyperref[\detokenize{pyapiref:chappyapi-util-tupledict}]{\sphinxcrossref{\DUrole{std,std-ref}{tupledict Class}}}} object;

\sphinxAtStartPar
When parameter \sphinxcode{\sphinxupquote{vars}} is list, {\hyperref[\detokenize{pyapiref:chappyapi-vararray}]{\sphinxcrossref{\DUrole{std,std-ref}{VarArray Class}}}} object,
parameter \sphinxcode{\sphinxupquote{vals}} can be constant or list.
\end{quote}

\sphinxAtStartPar
\sphinxstylestrong{Arguments}
\begin{quote}

\sphinxAtStartPar
\sphinxcode{\sphinxupquote{vars}}
\begin{quote}

\sphinxAtStartPar
The variables to be set value.
\end{quote}

\sphinxAtStartPar
\sphinxcode{\sphinxupquote{vals}}
\begin{quote}

\sphinxAtStartPar
The values of the variables in the solution.
\end{quote}
\end{quote}

\sphinxAtStartPar
\sphinxstylestrong{Example}

\begin{sphinxVerbatim}[commandchars=\\\{\}]
\PYG{n+nb+bp}{self}\PYG{o}{.}\PYG{n}{setSolution}\PYG{p}{(}\PYG{n}{x}\PYG{p}{,} \PYG{l+m+mi}{1}\PYG{p}{)}
\end{sphinxVerbatim}
\end{quote}

\subsubsection{CallbackBase.loadSolution()}
\label{\detokenize{pyapiref:callbackbase-loadsolution}}\begin{quote}

\sphinxAtStartPar
\sphinxstylestrong{Synopsis}
\begin{quote}

\sphinxAtStartPar
\sphinxcode{\sphinxupquote{loadSolution()}}
\end{quote}

\sphinxAtStartPar
\sphinxstylestrong{Description}
\begin{quote}

\sphinxAtStartPar
Load the currently feasible solution into model.

\sphinxAtStartPar
Note that a complete solution is required here.
\end{quote}

\sphinxAtStartPar
\sphinxstylestrong{Example}

\begin{sphinxVerbatim}[commandchars=\\\{\}]
\PYG{n+nb+bp}{self}\PYG{o}{.}\PYG{n}{loadSolution}\PYG{p}{(}\PYG{p}{)}
\end{sphinxVerbatim}
\end{quote}

\subsection{NlpCallbackBase Class}
\label{\detokenize{pyapiref:nlpcallbackbase-class}}\label{\detokenize{pyapiref:chappyapi-nlpcbcbase}}
\sphinxAtStartPar
The NlpCallbackBase class provides an interface for defining evaluation callbacks
for nonlinear optimization models in COPT.

\sphinxAtStartPar
By inheriting this class and implementing the corresponding member methods,
users provide objective and constraint value evaluations, together with their
first\sphinxhyphen{}order and second\sphinxhyphen{}order derivative information to COPT during nonlinear optimization.

\subsubsection{NlpCallbackBase.EvalObj()}
\label{\detokenize{pyapiref:nlpcallbackbase-evalobj}}\begin{quote}

\sphinxAtStartPar
\sphinxstylestrong{Synopsis}
\begin{quote}

\sphinxAtStartPar
\sphinxcode{\sphinxupquote{EvalObj(xdata, outdata)}}
\end{quote}

\sphinxAtStartPar
\sphinxstylestrong{Description}
\begin{quote}

\sphinxAtStartPar
Evaluate the objective function value of the nonlinear model.
\end{quote}

\sphinxAtStartPar
\sphinxstylestrong{Arguments}
\begin{quote}

\sphinxAtStartPar
\sphinxcode{\sphinxupquote{xdata}}
\begin{quote}

\sphinxAtStartPar
A one\sphinxhyphen{}dimensional array containing the values of the decision variables.

\sphinxAtStartPar
Acceptable values: {\hyperref[\detokenize{pyapiref:chappyapi-ndarray}]{\sphinxcrossref{\DUrole{std,std-ref}{NdArray Class}}}}.
\end{quote}

\sphinxAtStartPar
\sphinxcode{\sphinxupquote{outdata}}
\begin{quote}

\sphinxAtStartPar
An output array used to store the objective function value
(this array contains exactly one element).

\sphinxAtStartPar
Acceptable values: {\hyperref[\detokenize{pyapiref:chappyapi-ndarray}]{\sphinxcrossref{\DUrole{std,std-ref}{NdArray Class}}}}.
\end{quote}
\end{quote}

\sphinxAtStartPar
\sphinxstylestrong{Example}

\begin{sphinxVerbatim}[commandchars=\\\{\}]
\PYG{k}{class}\PYG{+w}{ }\PYG{n+nc}{MyCallback}\PYG{p}{(}\PYG{n}{NlpCallbackBase}\PYG{p}{)}\PYG{p}{:}
    \PYG{k}{def}\PYG{+w}{ }\PYG{n+nf}{EvalObj}\PYG{p}{(}\PYG{n+nb+bp}{self}\PYG{p}{,} \PYG{n}{xdata}\PYG{p}{,} \PYG{n}{outdata}\PYG{p}{)}\PYG{p}{:}
        \PYG{n}{x} \PYG{o}{=} \PYG{n}{NdArray}\PYG{p}{(}\PYG{n}{xdata}\PYG{p}{)}
        \PYG{n}{outval} \PYG{o}{=} \PYG{n}{NdArray}\PYG{p}{(}\PYG{n}{outdata}\PYG{p}{)}
        \PYG{n}{outval}\PYG{p}{[}\PYG{l+m+mi}{0}\PYG{p}{]} \PYG{o}{=} \PYG{n}{x}\PYG{p}{[}\PYG{l+m+mi}{0}\PYG{p}{]} \PYG{o}{*} \PYG{n}{x}\PYG{p}{[}\PYG{l+m+mi}{0}\PYG{p}{]} \PYG{o}{+} \PYG{n}{x}\PYG{p}{[}\PYG{l+m+mi}{1}\PYG{p}{]} \PYG{o}{*} \PYG{n}{x}\PYG{p}{[}\PYG{l+m+mi}{1}\PYG{p}{]}
        \PYG{k}{return} \PYG{l+m+mi}{0}
\end{sphinxVerbatim}
\end{quote}

\subsubsection{NlpCallbackBase.EvalGrad()}
\label{\detokenize{pyapiref:nlpcallbackbase-evalgrad}}\begin{quote}

\sphinxAtStartPar
\sphinxstylestrong{Synopsis}
\begin{quote}

\sphinxAtStartPar
\sphinxcode{\sphinxupquote{EvalGrad(xdata, outdata)}}
\end{quote}

\sphinxAtStartPar
\sphinxstylestrong{Description}
\begin{quote}

\sphinxAtStartPar
Evaluate the gradient of the objective function of the nonlinear model.
\end{quote}

\sphinxAtStartPar
\sphinxstylestrong{Arguments}
\begin{quote}

\sphinxAtStartPar
\sphinxcode{\sphinxupquote{xdata}}
\begin{quote}

\sphinxAtStartPar
A one\sphinxhyphen{}dimensional array containing the current values of the decision variables.

\sphinxAtStartPar
Acceptable values: {\hyperref[\detokenize{pyapiref:chappyapi-ndarray}]{\sphinxcrossref{\DUrole{std,std-ref}{NdArray Class}}}}.
\end{quote}

\sphinxAtStartPar
\sphinxcode{\sphinxupquote{outdata}}
\begin{quote}

\sphinxAtStartPar
An output array used to store the objective function gradient.
The length of this array equals the number of variables, and
the element order corresponds to the variable order.

\sphinxAtStartPar
Acceptable values: {\hyperref[\detokenize{pyapiref:chappyapi-ndarray}]{\sphinxcrossref{\DUrole{std,std-ref}{NdArray Class}}}}.
\end{quote}
\end{quote}

\sphinxAtStartPar
\sphinxstylestrong{Return value}
\begin{quote}

\sphinxAtStartPar
Returns an integer status code. A value of \sphinxcode{\sphinxupquote{0}} indicates successful evaluation.
\end{quote}

\sphinxAtStartPar
\sphinxstylestrong{Example}

\begin{sphinxVerbatim}[commandchars=\\\{\}]
\PYG{k}{class}\PYG{+w}{ }\PYG{n+nc}{MyCallback}\PYG{p}{(}\PYG{n}{NlpCallbackBase}\PYG{p}{)}\PYG{p}{:}
    \PYG{k}{def}\PYG{+w}{ }\PYG{n+nf}{EvalGrad}\PYG{p}{(}\PYG{n+nb+bp}{self}\PYG{p}{,} \PYG{n}{xdata}\PYG{p}{,} \PYG{n}{outdata}\PYG{p}{)}\PYG{p}{:}
        \PYG{n}{x} \PYG{o}{=} \PYG{n}{NdArray}\PYG{p}{(}\PYG{n}{xdata}\PYG{p}{)}
        \PYG{n}{outval} \PYG{o}{=} \PYG{n}{NdArray}\PYG{p}{(}\PYG{n}{outdata}\PYG{p}{)}
        \PYG{n}{outval}\PYG{p}{[}\PYG{l+m+mi}{0}\PYG{p}{]} \PYG{o}{=} \PYG{l+m+mf}{2.0} \PYG{o}{*} \PYG{n}{x}\PYG{p}{[}\PYG{l+m+mi}{0}\PYG{p}{]}
        \PYG{n}{outval}\PYG{p}{[}\PYG{l+m+mi}{1}\PYG{p}{]} \PYG{o}{=} \PYG{l+m+mf}{2.0} \PYG{o}{*} \PYG{n}{x}\PYG{p}{[}\PYG{l+m+mi}{1}\PYG{p}{]}
        \PYG{k}{return} \PYG{l+m+mi}{0}
\end{sphinxVerbatim}
\end{quote}

\subsubsection{NlpCallbackBase.EvalCon()}
\label{\detokenize{pyapiref:nlpcallbackbase-evalcon}}\begin{quote}

\sphinxAtStartPar
\sphinxstylestrong{Synopsis}
\begin{quote}

\sphinxAtStartPar
\sphinxcode{\sphinxupquote{EvalCon(xdata, outdata)}}
\end{quote}

\sphinxAtStartPar
\sphinxstylestrong{Description}
\begin{quote}

\sphinxAtStartPar
Evaluate the constraint function values of the nonlinear model.
\end{quote}

\sphinxAtStartPar
\sphinxstylestrong{Arguments}
\begin{quote}

\sphinxAtStartPar
\sphinxcode{\sphinxupquote{xdata}}
\begin{quote}

\sphinxAtStartPar
A one\sphinxhyphen{}dimensional array containing the values of the decision variables.

\sphinxAtStartPar
Acceptable values: {\hyperref[\detokenize{pyapiref:chappyapi-ndarray}]{\sphinxcrossref{\DUrole{std,std-ref}{NdArray Class}}}}.
\end{quote}

\sphinxAtStartPar
\sphinxcode{\sphinxupquote{outdata}}
\begin{quote}

\sphinxAtStartPar
An output array used to store the constraint function values.
The length of this array equals the number of constraints, and
the element order corresponds to the constraint definition order in the model.

\sphinxAtStartPar
Acceptable values: {\hyperref[\detokenize{pyapiref:chappyapi-ndarray}]{\sphinxcrossref{\DUrole{std,std-ref}{NdArray Class}}}}.
\end{quote}
\end{quote}

\sphinxAtStartPar
\sphinxstylestrong{Example}

\begin{sphinxVerbatim}[commandchars=\\\{\}]
\PYG{k}{class}\PYG{+w}{ }\PYG{n+nc}{MyCallback}\PYG{p}{(}\PYG{n}{NlpCallbackBase}\PYG{p}{)}\PYG{p}{:}
    \PYG{k}{def}\PYG{+w}{ }\PYG{n+nf}{EvalCon}\PYG{p}{(}\PYG{n+nb+bp}{self}\PYG{p}{,} \PYG{n}{xdata}\PYG{p}{,} \PYG{n}{outdata}\PYG{p}{)}\PYG{p}{:}
        \PYG{n}{x} \PYG{o}{=} \PYG{n}{NdArray}\PYG{p}{(}\PYG{n}{xdata}\PYG{p}{)}
        \PYG{n}{outval} \PYG{o}{=} \PYG{n}{NdArray}\PYG{p}{(}\PYG{n}{outdata}\PYG{p}{)}
        \PYG{n}{outval}\PYG{p}{[}\PYG{l+m+mi}{0}\PYG{p}{]} \PYG{o}{=} \PYG{n}{x}\PYG{p}{[}\PYG{l+m+mi}{0}\PYG{p}{]} \PYG{o}{*} \PYG{n}{x}\PYG{p}{[}\PYG{l+m+mi}{0}\PYG{p}{]} \PYG{o}{+} \PYG{n}{x}\PYG{p}{[}\PYG{l+m+mi}{1}\PYG{p}{]} \PYG{o}{\PYGZhy{}} \PYG{l+m+mf}{1.0}
        \PYG{k}{return} \PYG{l+m+mi}{0}
\end{sphinxVerbatim}
\end{quote}

\subsubsection{NlpCallbackBase.EvalJac()}
\label{\detokenize{pyapiref:nlpcallbackbase-evaljac}}\begin{quote}

\sphinxAtStartPar
\sphinxstylestrong{Synopsis}
\begin{quote}

\sphinxAtStartPar
\sphinxcode{\sphinxupquote{EvalJac(xdata, outdata)}}
\end{quote}

\sphinxAtStartPar
\sphinxstylestrong{Description}
\begin{quote}

\sphinxAtStartPar
Evaluate the Jacobian matrix values of the nonlinear constraint functions.
\end{quote}

\sphinxAtStartPar
\sphinxstylestrong{Arguments}
\begin{quote}

\sphinxAtStartPar
\sphinxcode{\sphinxupquote{xdata}}
\begin{quote}

\sphinxAtStartPar
A one\sphinxhyphen{}dimensional array containing the values of the decision variables.

\sphinxAtStartPar
Acceptable values: {\hyperref[\detokenize{pyapiref:chappyapi-ndarray}]{\sphinxcrossref{\DUrole{std,std-ref}{NdArray Class}}}}.
\end{quote}

\sphinxAtStartPar
\sphinxcode{\sphinxupquote{outdata}}
\begin{quote}

\sphinxAtStartPar
An output array used to store the nonzero elements of the constraint Jacobian matrix.

\sphinxAtStartPar
Acceptable values: {\hyperref[\detokenize{pyapiref:chappyapi-ndarray}]{\sphinxcrossref{\DUrole{std,std-ref}{NdArray Class}}}}.

\sphinxAtStartPar
The length of this array equals the number of nonzero elements in the Jacobian matrix.
Each element corresponds to the position specified by \sphinxcode{\sphinxupquote{idxJacRow}} and
\sphinxcode{\sphinxupquote{idxJacCol}} provided in \sphinxcode{\sphinxupquote{Model.loadNlData}}.
\end{quote}
\end{quote}

\sphinxAtStartPar
\sphinxstylestrong{Example}

\begin{sphinxVerbatim}[commandchars=\\\{\}]
\PYG{k}{class}\PYG{+w}{ }\PYG{n+nc}{MyCallback}\PYG{p}{(}\PYG{n}{NlpCallbackBase}\PYG{p}{)}\PYG{p}{:}
    \PYG{k}{def}\PYG{+w}{ }\PYG{n+nf}{EvalJac}\PYG{p}{(}\PYG{n+nb+bp}{self}\PYG{p}{,} \PYG{n}{xdata}\PYG{p}{,} \PYG{n}{outdata}\PYG{p}{)}\PYG{p}{:}
        \PYG{n}{x} \PYG{o}{=} \PYG{n}{NdArray}\PYG{p}{(}\PYG{n}{xdata}\PYG{p}{)}
        \PYG{n}{outval} \PYG{o}{=} \PYG{n}{NdArray}\PYG{p}{(}\PYG{n}{outdata}\PYG{p}{)}
        \PYG{n}{outval}\PYG{p}{[}\PYG{l+m+mi}{0}\PYG{p}{]} \PYG{o}{=} \PYG{l+m+mf}{2.0} \PYG{o}{*} \PYG{n}{x}\PYG{p}{[}\PYG{l+m+mi}{0}\PYG{p}{]}
        \PYG{n}{outval}\PYG{p}{[}\PYG{l+m+mi}{1}\PYG{p}{]} \PYG{o}{=} \PYG{l+m+mf}{1.0}
        \PYG{k}{return} \PYG{l+m+mi}{0}
\end{sphinxVerbatim}
\end{quote}

\subsubsection{NlpCallbackBase.EvalHess()}
\label{\detokenize{pyapiref:nlpcallbackbase-evalhess}}\begin{quote}

\sphinxAtStartPar
\sphinxstylestrong{Synopsis}
\begin{quote}

\sphinxAtStartPar
\sphinxcode{\sphinxupquote{EvalHess(xdata, sigma, lambdata, outdata)}}
\end{quote}

\sphinxAtStartPar
\sphinxstylestrong{Description}
\begin{quote}

\sphinxAtStartPar
Evaluate the Hessian matrix values of the Lagrangian function
of the nonlinear model.
\end{quote}

\sphinxAtStartPar
\sphinxstylestrong{Arguments}
\begin{quote}

\sphinxAtStartPar
\sphinxcode{\sphinxupquote{xdata}}
\begin{quote}

\sphinxAtStartPar
A one\sphinxhyphen{}dimensional array containing the values of the decision variables.

\sphinxAtStartPar
Acceptable values: {\hyperref[\detokenize{pyapiref:chappyapi-ndarray}]{\sphinxcrossref{\DUrole{std,std-ref}{NdArray Class}}}}.
\end{quote}

\sphinxAtStartPar
\sphinxcode{\sphinxupquote{sigma}}
\begin{quote}

\sphinxAtStartPar
The weight coefficient of the objective in the Lagrangian function.
\end{quote}

\sphinxAtStartPar
\sphinxcode{\sphinxupquote{lambdata}}
\begin{quote}

\sphinxAtStartPar
A one\sphinxhyphen{}dimensional array containing the current Lagrange multiplier values
associated with the constraint functions.

\sphinxAtStartPar
Acceptable values: {\hyperref[\detokenize{pyapiref:chappyapi-ndarray}]{\sphinxcrossref{\DUrole{std,std-ref}{NdArray Class}}}}.
\end{quote}

\sphinxAtStartPar
\sphinxcode{\sphinxupquote{outdata}}
\begin{quote}

\sphinxAtStartPar
An output array used to store the nonzero elements of the Hessian matrix
of the Lagrangian function.

\sphinxAtStartPar
Acceptable values: {\hyperref[\detokenize{pyapiref:chappyapi-ndarray}]{\sphinxcrossref{\DUrole{std,std-ref}{NdArray Class}}}}.

\sphinxAtStartPar
The length of this array equals the number of nonzero elements in the Hessian matrix.
Each element corresponds to the position specified by \sphinxcode{\sphinxupquote{idxHessRow}} and
\sphinxcode{\sphinxupquote{idxHessCol}} provided in \sphinxcode{\sphinxupquote{Model.loadNlData}}.
\end{quote}
\end{quote}

\sphinxAtStartPar
\sphinxstylestrong{Example}

\begin{sphinxVerbatim}[commandchars=\\\{\}]
\PYG{k}{class}\PYG{+w}{ }\PYG{n+nc}{MyCallback}\PYG{p}{(}\PYG{n}{NlpCallbackBase}\PYG{p}{)}\PYG{p}{:}
    \PYG{k}{def}\PYG{+w}{ }\PYG{n+nf}{EvalHess}\PYG{p}{(}\PYG{n+nb+bp}{self}\PYG{p}{,} \PYG{n}{xdata}\PYG{p}{,} \PYG{n}{sigma}\PYG{p}{,} \PYG{n}{lambdata}\PYG{p}{,} \PYG{n}{outdata}\PYG{p}{)}\PYG{p}{:}
        \PYG{n}{x} \PYG{o}{=} \PYG{n}{NdArray}\PYG{p}{(}\PYG{n}{xdata}\PYG{p}{)}
        \PYG{n}{lam} \PYG{o}{=} \PYG{n}{NdArray}\PYG{p}{(}\PYG{n}{lambdata}\PYG{p}{)}
        \PYG{n}{outval} \PYG{o}{=} \PYG{n}{NdArray}\PYG{p}{(}\PYG{n}{outdata}\PYG{p}{)}
        \PYG{n}{outval}\PYG{p}{[}\PYG{l+m+mi}{0}\PYG{p}{]} \PYG{o}{=} \PYG{l+m+mf}{2.0} \PYG{o}{*} \PYG{n}{sigma} \PYG{o}{+} \PYG{l+m+mf}{2.0} \PYG{o}{*} \PYG{n}{lam}\PYG{p}{[}\PYG{l+m+mi}{0}\PYG{p}{]}
        \PYG{k}{return} \PYG{l+m+mi}{0}
\end{sphinxVerbatim}
\end{quote}

\subsection{GenConstrX Class}
\label{\detokenize{pyapiref:genconstrx-class}}
\sphinxAtStartPar
In the \sphinxcode{\sphinxupquote{Model}} class, the constraints added by \sphinxcode{\sphinxupquote{addGenConstrXXX}} (such as: \sphinxcode{\sphinxupquote{addGenConstrMax}}) will return a \sphinxcode{\sphinxupquote{GenConstrX}} object.

\subsubsection{GenConstrX.getAttr()}
\label{\detokenize{pyapiref:genconstrx-getattr}}\begin{quote}

\sphinxAtStartPar
\sphinxstylestrong{Synopsis}
\begin{quote}

\sphinxAtStartPar
\sphinxcode{\sphinxupquote{getAttr(attrname)}}
\end{quote}

\sphinxAtStartPar
\sphinxstylestrong{Description}
\begin{quote}

\sphinxAtStartPar
Get the attribute value of the \sphinxcode{\sphinxupquote{GenConstrX}} class object, support to get the type and name of the \sphinxcode{\sphinxupquote{GenConstrX}} class object.
\end{quote}

\sphinxAtStartPar
\sphinxstylestrong{Example}
\begin{quote}

\begin{sphinxVerbatim}[commandchars=\\\{\}]
\PYG{c+c1}{\PYGZsh{} Get the name of con\PYGZus{}max}
\PYG{n}{con\PYGZus{}max}\PYG{o}{.}\PYG{n}{getAttr}\PYG{p}{(}\PYG{l+s+s2}{\PYGZdq{}}\PYG{l+s+s2}{name}\PYG{l+s+s2}{\PYGZdq{}}\PYG{p}{)}
\PYG{c+c1}{\PYGZsh{} Get the type of con\PYGZus{}max}
\PYG{n}{con\PYGZus{}max}\PYG{o}{.}\PYG{n}{getAttr}\PYG{p}{(}\PYG{l+s+s2}{\PYGZdq{}}\PYG{l+s+s2}{type}\PYG{l+s+s2}{\PYGZdq{}}\PYG{p}{)}
\end{sphinxVerbatim}
\end{quote}
\end{quote}

\subsubsection{GenConstrX.setAttr()}
\label{\detokenize{pyapiref:genconstrx-setattr}}\begin{quote}

\sphinxAtStartPar
\sphinxstylestrong{Synopsis}
\begin{quote}

\sphinxAtStartPar
\sphinxcode{\sphinxupquote{setAttr(attrname)}}
\end{quote}

\sphinxAtStartPar
\sphinxstylestrong{Arguments}
\begin{quote}

\sphinxAtStartPar
Set the attribute value of the \sphinxcode{\sphinxupquote{GenConstrX}} class object, and support setting the name of the \sphinxcode{\sphinxupquote{GenConstrX}} class object.
\end{quote}

\sphinxAtStartPar
\sphinxstylestrong{Example}
\begin{quote}

\begin{sphinxVerbatim}[commandchars=\\\{\}]
\PYG{c+c1}{\PYGZsh{} Set the name of con\PYGZus{}max}
\PYG{n}{con\PYGZus{}max}\PYG{o}{.}\PYG{n}{setAttr}\PYG{p}{(}\PYG{l+s+s2}{\PYGZdq{}}\PYG{l+s+s2}{name}\PYG{l+s+s2}{\PYGZdq{}}\PYG{p}{)}
\end{sphinxVerbatim}
\end{quote}
\end{quote}

\subsection{CoptError Class}
\label{\detokenize{pyapiref:copterror-class}}\label{\detokenize{pyapiref:chappyapi-error}}
\sphinxAtStartPar
CoptError Class provides operations on error.
An exception of the CoptError is thrown when error occurs in a method call corresponding to the
underlying interface of solver.
The following attributes are provided to retrieve error information:
\begin{itemize}
\item {} 
\sphinxAtStartPar
CoptError.retcode
\begin{quote}

\sphinxAtStartPar
Error code.
\end{quote}

\item {} 
\sphinxAtStartPar
CoptError.message
\begin{quote}

\sphinxAtStartPar
Error message.
\end{quote}

\end{itemize}

\section{Helper Functions and Utilities}
\label{\detokenize{pyapiref:helper-functions-and-utilities}}\label{\detokenize{pyapiref:chappyapi-utility}}
\sphinxAtStartPar
Helper functions and utilities are encapsulated based on Python’s basic data types,
providing easy\sphinxhyphen{}to\sphinxhyphen{}use data types to facilitate the rapid construction of complex optimization models.
This section will explain its functions and usages.

\subsection{Helper Functions}
\label{\detokenize{pyapiref:helper-functions}}

\subsubsection{multidict()}
\label{\detokenize{pyapiref:multidict}}\begin{quote}

\sphinxAtStartPar
\sphinxstylestrong{Synopsis}
\begin{quote}

\sphinxAtStartPar
\sphinxcode{\sphinxupquote{multidict(data)}}
\end{quote}

\sphinxAtStartPar
\sphinxstylestrong{Description}
\begin{quote}

\sphinxAtStartPar
Split a single dictionary into keys and multiple dictionaries. Return keys and dictionaries.
\end{quote}

\sphinxAtStartPar
\sphinxstylestrong{Arguments}
\begin{quote}

\sphinxAtStartPar
\sphinxcode{\sphinxupquote{data}}
\begin{quote}

\sphinxAtStartPar
A Python dictionary to be applied. Each key should map to a list of \(n\) values.
\end{quote}
\end{quote}

\sphinxAtStartPar
\sphinxstylestrong{Example}
\end{quote}

\begin{sphinxVerbatim}[commandchars=\\\{\}]
\PYG{n}{keys}\PYG{p}{,} \PYG{n}{dict1}\PYG{p}{,} \PYG{n}{dict2} \PYG{o}{=} \PYG{n}{multidict}\PYG{p}{(}\PYG{p}{\PYGZob{}}
  \PYG{l+s+s2}{\PYGZdq{}}\PYG{l+s+s2}{hello}\PYG{l+s+s2}{\PYGZdq{}}\PYG{p}{:} \PYG{p}{[}\PYG{l+m+mi}{0}\PYG{p}{,} \PYG{l+m+mi}{1}\PYG{p}{]}\PYG{p}{,}
  \PYG{l+s+s2}{\PYGZdq{}}\PYG{l+s+s2}{world}\PYG{l+s+s2}{\PYGZdq{}}\PYG{p}{:} \PYG{p}{[}\PYG{l+m+mi}{2}\PYG{p}{,} \PYG{l+m+mi}{3}\PYG{p}{]}\PYG{p}{\PYGZcb{}}\PYG{p}{)}
\end{sphinxVerbatim}

\subsubsection{quicksum()}
\label{\detokenize{pyapiref:quicksum}}\begin{quote}

\sphinxAtStartPar
\sphinxstylestrong{Synopsis}
\begin{quote}

\sphinxAtStartPar
\sphinxcode{\sphinxupquote{quicksum(data)}}
\end{quote}

\sphinxAtStartPar
\sphinxstylestrong{Description}
\begin{quote}

\sphinxAtStartPar
Build expressions efficiently. Return a {\hyperref[\detokenize{pyapiref:chappyapi-linexpr}]{\sphinxcrossref{\DUrole{std,std-ref}{LinExpr Class}}}} object.
\end{quote}

\sphinxAtStartPar
\sphinxstylestrong{Arguments}
\begin{quote}

\sphinxAtStartPar
\sphinxcode{\sphinxupquote{data}}
\begin{quote}

\sphinxAtStartPar
Terms to add.
\end{quote}
\end{quote}

\sphinxAtStartPar
\sphinxstylestrong{Example}
\end{quote}

\begin{sphinxVerbatim}[commandchars=\\\{\}]
\PYG{n}{expr} \PYG{o}{=} \PYG{n}{quicksum}\PYG{p}{(}\PYG{n}{m}\PYG{o}{.}\PYG{n}{getVars}\PYG{p}{(}\PYG{p}{)}\PYG{p}{)}
\end{sphinxVerbatim}

\subsection{tuplelist Class}
\label{\detokenize{pyapiref:tuplelist-class}}\label{\detokenize{pyapiref:chappyapi-util-tuplelist}}
\sphinxAtStartPar
The tuplelist object is an encapsulation based on Python lists, and provides the following methods:

\subsubsection{tuplelist()}
\label{\detokenize{pyapiref:tuplelist}}\begin{quote}

\sphinxAtStartPar
\sphinxstylestrong{Synopsis}
\begin{quote}

\sphinxAtStartPar
\sphinxcode{\sphinxupquote{tuplelist(list)}}
\end{quote}

\sphinxAtStartPar
\sphinxstylestrong{Description}
\begin{quote}

\sphinxAtStartPar
Create and return a {\hyperref[\detokenize{pyapiref:chappyapi-util-tuplelist}]{\sphinxcrossref{\DUrole{std,std-ref}{tuplelist Class}}}} object.
\end{quote}

\sphinxAtStartPar
\sphinxstylestrong{Arguments}
\begin{quote}

\sphinxAtStartPar
\sphinxcode{\sphinxupquote{list}}
\begin{quote}

\sphinxAtStartPar
A Python list.
\end{quote}
\end{quote}

\sphinxAtStartPar
\sphinxstylestrong{Example}
\end{quote}

\begin{sphinxVerbatim}[commandchars=\\\{\}]
\PYG{n}{tl} \PYG{o}{=} \PYG{n}{tuplelist}\PYG{p}{(}\PYG{p}{[}\PYG{p}{(}\PYG{l+m+mi}{0}\PYG{p}{,} \PYG{l+m+mi}{1}\PYG{p}{)}\PYG{p}{,} \PYG{p}{(}\PYG{l+m+mi}{1}\PYG{p}{,} \PYG{l+m+mi}{2}\PYG{p}{)}\PYG{p}{]}\PYG{p}{)}
\PYG{n}{tl} \PYG{o}{=} \PYG{n}{tuplelist}\PYG{p}{(}\PYG{p}{[}\PYG{p}{(}\PYG{l+s+s1}{\PYGZsq{}}\PYG{l+s+s1}{a}\PYG{l+s+s1}{\PYGZsq{}}\PYG{p}{,} \PYG{l+s+s1}{\PYGZsq{}}\PYG{l+s+s1}{b}\PYG{l+s+s1}{\PYGZsq{}}\PYG{p}{)}\PYG{p}{,} \PYG{p}{(}\PYG{l+s+s1}{\PYGZsq{}}\PYG{l+s+s1}{b}\PYG{l+s+s1}{\PYGZsq{}}\PYG{p}{,} \PYG{l+s+s1}{\PYGZsq{}}\PYG{l+s+s1}{c}\PYG{l+s+s1}{\PYGZsq{}}\PYG{p}{)}\PYG{p}{]}\PYG{p}{)}
\end{sphinxVerbatim}

\subsubsection{tuplelist.add()}
\label{\detokenize{pyapiref:tuplelist-add}}\begin{quote}

\sphinxAtStartPar
\sphinxstylestrong{Synopsis}
\begin{quote}

\sphinxAtStartPar
\sphinxcode{\sphinxupquote{add(item)}}
\end{quote}

\sphinxAtStartPar
\sphinxstylestrong{Description}
\begin{quote}

\sphinxAtStartPar
Add an item to a {\hyperref[\detokenize{pyapiref:chappyapi-util-tuplelist}]{\sphinxcrossref{\DUrole{std,std-ref}{tuplelist Class}}}} object
\end{quote}

\sphinxAtStartPar
\sphinxstylestrong{Arguments}
\begin{quote}

\sphinxAtStartPar
\sphinxcode{\sphinxupquote{item}}
\begin{quote}

\sphinxAtStartPar
Item to add, which can be a Python tuple.
\end{quote}
\end{quote}

\sphinxAtStartPar
\sphinxstylestrong{Example}
\end{quote}

\begin{sphinxVerbatim}[commandchars=\\\{\}]
\PYG{n}{tl} \PYG{o}{=} \PYG{n}{tuplelist}\PYG{p}{(}\PYG{p}{[}\PYG{p}{(}\PYG{l+m+mi}{0}\PYG{p}{,} \PYG{l+m+mi}{1}\PYG{p}{)}\PYG{p}{,} \PYG{p}{(}\PYG{l+m+mi}{1}\PYG{p}{,} \PYG{l+m+mi}{2}\PYG{p}{)}\PYG{p}{]}\PYG{p}{)}
\PYG{n}{tl}\PYG{o}{.}\PYG{n}{add}\PYG{p}{(}\PYG{p}{(}\PYG{l+m+mi}{2}\PYG{p}{,} \PYG{l+m+mi}{3}\PYG{p}{)}\PYG{p}{)}
\end{sphinxVerbatim}

\subsubsection{tuplelist.select()}
\label{\detokenize{pyapiref:tuplelist-select}}\begin{quote}

\sphinxAtStartPar
\sphinxstylestrong{Synopsis}
\begin{quote}

\sphinxAtStartPar
\sphinxcode{\sphinxupquote{select(pattern)}}
\end{quote}

\sphinxAtStartPar
\sphinxstylestrong{Description}
\begin{quote}

\sphinxAtStartPar
Get all terms that match the specified pattern. Return a {\hyperref[\detokenize{pyapiref:chappyapi-util-tuplelist}]{\sphinxcrossref{\DUrole{std,std-ref}{tuplelist Class}}}} object.
\end{quote}

\sphinxAtStartPar
\sphinxstylestrong{Arguments}
\begin{quote}

\sphinxAtStartPar
\sphinxcode{\sphinxupquote{pattern}}
\begin{quote}

\sphinxAtStartPar
Specified pattern.
\end{quote}
\end{quote}

\sphinxAtStartPar
\sphinxstylestrong{Example}
\end{quote}

\begin{sphinxVerbatim}[commandchars=\\\{\}]
\PYG{n}{tl} \PYG{o}{=} \PYG{n}{tuplelist}\PYG{p}{(}\PYG{p}{[}\PYG{p}{(}\PYG{l+m+mi}{0}\PYG{p}{,} \PYG{l+m+mi}{1}\PYG{p}{)}\PYG{p}{,} \PYG{p}{(}\PYG{l+m+mi}{0}\PYG{p}{,} \PYG{l+m+mi}{2}\PYG{p}{)}\PYG{p}{,} \PYG{p}{(}\PYG{l+m+mi}{1}\PYG{p}{,} \PYG{l+m+mi}{2}\PYG{p}{)}\PYG{p}{]}\PYG{p}{)}
\PYG{n}{tl}\PYG{o}{.}\PYG{n}{select}\PYG{p}{(}\PYG{l+m+mi}{0}\PYG{p}{,} \PYG{l+s+s1}{\PYGZsq{}}\PYG{l+s+s1}{*}\PYG{l+s+s1}{\PYGZsq{}}\PYG{p}{)}
\end{sphinxVerbatim}

\subsubsection{repeat()}
\label{\detokenize{pyapiref:repeat}}\begin{quote}

\sphinxAtStartPar
\sphinxstylestrong{Synopsis}
\begin{quote}

\sphinxAtStartPar
\sphinxcode{\sphinxupquote{repeat(obj, repeats)}}
\end{quote}

\sphinxAtStartPar
\sphinxstylestrong{Description}
\begin{quote}

\sphinxAtStartPar
Repeats the specified object a given number of times.
\end{quote}

\sphinxAtStartPar
\sphinxstylestrong{Arguments}
\begin{quote}

\sphinxAtStartPar
\sphinxcode{\sphinxupquote{obj}}
\begin{quote}

\sphinxAtStartPar
A Python object.
\end{quote}

\sphinxAtStartPar
\sphinxcode{\sphinxupquote{repeats}}
\begin{quote}

\sphinxAtStartPar
The number of repetitions.
\end{quote}
\end{quote}
\end{quote}

\subsubsection{hstack()}
\label{\detokenize{pyapiref:hstack}}\begin{quote}

\sphinxAtStartPar
\sphinxstylestrong{Synopsis}
\begin{quote}

\sphinxAtStartPar
\sphinxcode{\sphinxupquote{hstack(left, right)}}
\end{quote}

\sphinxAtStartPar
\sphinxstylestrong{Description}
\begin{quote}

\sphinxAtStartPar
Horizontally stacks the input objects.
\end{quote}

\sphinxAtStartPar
\sphinxstylestrong{Arguments}
\begin{quote}

\sphinxAtStartPar
\sphinxcode{\sphinxupquote{left}}
\begin{quote}

\sphinxAtStartPar
A Python object.
\end{quote}

\sphinxAtStartPar
\sphinxcode{\sphinxupquote{right}}
\begin{quote}

\sphinxAtStartPar
A Python object.
\end{quote}
\end{quote}
\end{quote}

\subsubsection{vstack()}
\label{\detokenize{pyapiref:vstack}}\begin{quote}

\sphinxAtStartPar
\sphinxstylestrong{Synopsis}
\begin{quote}

\sphinxAtStartPar
\sphinxcode{\sphinxupquote{vstack(left, right)}}
\end{quote}

\sphinxAtStartPar
\sphinxstylestrong{Description}
\begin{quote}

\sphinxAtStartPar
Vertically stacks the input objects.
\end{quote}

\sphinxAtStartPar
\sphinxstylestrong{Arguments}
\begin{quote}

\sphinxAtStartPar
\sphinxcode{\sphinxupquote{left}}
\begin{quote}

\sphinxAtStartPar
A Python object.
\end{quote}

\sphinxAtStartPar
\sphinxcode{\sphinxupquote{right}}
\begin{quote}

\sphinxAtStartPar
A Python object.
\end{quote}
\end{quote}
\end{quote}

\subsubsection{stack()}
\label{\detokenize{pyapiref:stack}}\begin{quote}

\sphinxAtStartPar
\sphinxstylestrong{Synopsis}
\begin{quote}

\sphinxAtStartPar
\sphinxcode{\sphinxupquote{stack(left, right, axis)}}
\end{quote}

\sphinxAtStartPar
\sphinxstylestrong{Description}
\begin{quote}

\sphinxAtStartPar
Stacks the input objects along the specified axis.
\end{quote}

\sphinxAtStartPar
\sphinxstylestrong{Arguments}
\begin{quote}

\sphinxAtStartPar
\sphinxcode{\sphinxupquote{left}}
\begin{quote}

\sphinxAtStartPar
A Python object.
\end{quote}

\sphinxAtStartPar
\sphinxcode{\sphinxupquote{right}}
\begin{quote}

\sphinxAtStartPar
A Python object.
\end{quote}

\sphinxAtStartPar
\sphinxcode{\sphinxupquote{axis}}
\begin{quote}

\sphinxAtStartPar
The specified axis index.
\end{quote}
\end{quote}
\end{quote}

\subsection{tupledict Class}
\label{\detokenize{pyapiref:tupledict-class}}\label{\detokenize{pyapiref:chappyapi-util-tupledict}}
\sphinxAtStartPar
The tupledict class is an encapsulation based on Python dictionaries, and provides the following methods:

\subsubsection{tupledict()}
\label{\detokenize{pyapiref:tupledict}}\begin{quote}

\sphinxAtStartPar
\sphinxstylestrong{Synopsis}
\begin{quote}

\sphinxAtStartPar
\sphinxcode{\sphinxupquote{tupledict(args, kwargs)}}
\end{quote}

\sphinxAtStartPar
\sphinxstylestrong{Description}
\begin{quote}

\sphinxAtStartPar
Create and return a {\hyperref[\detokenize{pyapiref:chappyapi-util-tupledict}]{\sphinxcrossref{\DUrole{std,std-ref}{tupledict Class}}}} object.
\end{quote}

\sphinxAtStartPar
\sphinxstylestrong{Arguments}
\begin{quote}

\sphinxAtStartPar
\sphinxcode{\sphinxupquote{args}}
\begin{quote}

\sphinxAtStartPar
Positional arguments.
\end{quote}

\sphinxAtStartPar
\sphinxcode{\sphinxupquote{kwargs}}
\begin{quote}

\sphinxAtStartPar
Named arguments.
\end{quote}
\end{quote}

\sphinxAtStartPar
\sphinxstylestrong{Example}
\end{quote}

\begin{sphinxVerbatim}[commandchars=\\\{\}]
\PYG{n}{d} \PYG{o}{=} \PYG{n}{tupledict}\PYG{p}{(}\PYG{p}{[}\PYG{p}{(}\PYG{l+m+mi}{0}\PYG{p}{,} \PYG{l+s+s2}{\PYGZdq{}}\PYG{l+s+s2}{hello}\PYG{l+s+s2}{\PYGZdq{}}\PYG{p}{)}\PYG{p}{,} \PYG{p}{(}\PYG{l+m+mi}{1}\PYG{p}{,} \PYG{l+s+s2}{\PYGZdq{}}\PYG{l+s+s2}{world}\PYG{l+s+s2}{\PYGZdq{}}\PYG{p}{)}\PYG{p}{]}\PYG{p}{)}
\end{sphinxVerbatim}

\subsubsection{tupledict.select()}
\label{\detokenize{pyapiref:tupledict-select}}\begin{quote}

\sphinxAtStartPar
\sphinxstylestrong{Synopsis}
\begin{quote}

\sphinxAtStartPar
\sphinxcode{\sphinxupquote{select(pattern)}}
\end{quote}

\sphinxAtStartPar
\sphinxstylestrong{Description}
\begin{quote}

\sphinxAtStartPar
Get all terms that match the specified pattern. Return a {\hyperref[\detokenize{pyapiref:chappyapi-util-tupledict}]{\sphinxcrossref{\DUrole{std,std-ref}{tupledict Class}}}} object.
\end{quote}

\sphinxAtStartPar
\sphinxstylestrong{Arguments}
\begin{quote}

\sphinxAtStartPar
\sphinxcode{\sphinxupquote{pattern}}
\begin{quote}

\sphinxAtStartPar
Specified pattern.
\end{quote}
\end{quote}

\sphinxAtStartPar
\sphinxstylestrong{Example}
\end{quote}

\begin{sphinxVerbatim}[commandchars=\\\{\}]
\PYG{n}{d} \PYG{o}{=} \PYG{n}{tupledict}\PYG{p}{(}\PYG{p}{[}\PYG{p}{(}\PYG{l+m+mi}{0}\PYG{p}{,} \PYG{l+s+s2}{\PYGZdq{}}\PYG{l+s+s2}{hello}\PYG{l+s+s2}{\PYGZdq{}}\PYG{p}{)}\PYG{p}{,} \PYG{p}{(}\PYG{l+m+mi}{1}\PYG{p}{,} \PYG{l+s+s2}{\PYGZdq{}}\PYG{l+s+s2}{world}\PYG{l+s+s2}{\PYGZdq{}}\PYG{p}{)}\PYG{p}{]}\PYG{p}{)}
\PYG{n}{d}\PYG{o}{.}\PYG{n}{select}\PYG{p}{(}\PYG{p}{)}
\end{sphinxVerbatim}

\subsubsection{tupledict.sum()}
\label{\detokenize{pyapiref:tupledict-sum}}\begin{quote}

\sphinxAtStartPar
\sphinxstylestrong{Synopsis}
\begin{quote}

\sphinxAtStartPar
\sphinxcode{\sphinxupquote{sum(pattern)}}
\end{quote}

\sphinxAtStartPar
\sphinxstylestrong{Description}
\begin{quote}

\sphinxAtStartPar
Sum all terms that match the specified pattern. Return a {\hyperref[\detokenize{pyapiref:chappyapi-linexpr}]{\sphinxcrossref{\DUrole{std,std-ref}{LinExpr Class}}}} object.
\end{quote}

\sphinxAtStartPar
\sphinxstylestrong{Arguments}
\begin{quote}

\sphinxAtStartPar
\sphinxcode{\sphinxupquote{pattern}}
\begin{quote}

\sphinxAtStartPar
Specified pattern.
\end{quote}
\end{quote}

\sphinxAtStartPar
\sphinxstylestrong{Example}
\end{quote}

\begin{sphinxVerbatim}[commandchars=\\\{\}]
\PYG{n}{expr} \PYG{o}{=} \PYG{n}{x}\PYG{o}{.}\PYG{n}{sum}\PYG{p}{(}\PYG{p}{)}
\end{sphinxVerbatim}

\subsubsection{tupledict.prod()}
\label{\detokenize{pyapiref:tupledict-prod}}\begin{quote}

\sphinxAtStartPar
\sphinxstylestrong{Synopsis}
\begin{quote}

\sphinxAtStartPar
\sphinxcode{\sphinxupquote{prod(coeff, pattern)}}
\end{quote}

\sphinxAtStartPar
\sphinxstylestrong{Description}
\begin{quote}

\sphinxAtStartPar
Filter terms that match the specified pattern and multiply by coefficients. Return a {\hyperref[\detokenize{pyapiref:chappyapi-linexpr}]{\sphinxcrossref{\DUrole{std,std-ref}{LinExpr Class}}}} object.
\end{quote}

\sphinxAtStartPar
\sphinxstylestrong{Arguments}
\begin{quote}

\sphinxAtStartPar
\sphinxcode{\sphinxupquote{coeff}}
\begin{quote}

\sphinxAtStartPar
Coefficients, which can be a dict or a {\hyperref[\detokenize{pyapiref:chappyapi-util-tupledict}]{\sphinxcrossref{\DUrole{std,std-ref}{tupledict Class}}}} object.
\end{quote}

\sphinxAtStartPar
\sphinxcode{\sphinxupquote{pattern}}
\begin{quote}

\sphinxAtStartPar
Specified pattern.
\end{quote}
\end{quote}

\sphinxAtStartPar
\sphinxstylestrong{Example}
\end{quote}

\begin{sphinxVerbatim}[commandchars=\\\{\}]
\PYG{n}{coeff} \PYG{o}{=} \PYG{n+nb}{dict}\PYG{p}{(}\PYG{p}{[}\PYG{p}{(}\PYG{l+m+mi}{1}\PYG{p}{,} \PYG{l+m+mf}{0.1}\PYG{p}{)}\PYG{p}{,} \PYG{p}{(}\PYG{l+m+mi}{2}\PYG{p}{,} \PYG{l+m+mf}{0.2}\PYG{p}{)}\PYG{p}{]}\PYG{p}{)}
\PYG{n}{expr}  \PYG{o}{=} \PYG{n}{x}\PYG{o}{.}\PYG{n}{prod}\PYG{p}{(}\PYG{n}{coeff}\PYG{p}{)}
\end{sphinxVerbatim}

\subsection{ProbBuffer Class}
\label{\detokenize{pyapiref:probbuffer-class}}\label{\detokenize{pyapiref:chappyapi-util-probbuffer}}
\sphinxAtStartPar
The ProbBuffer is an encapsulation of buffer of string stream, and provides the following methods:

\subsubsection{ProbBuffer()}
\label{\detokenize{pyapiref:probbuffer}}\begin{quote}

\sphinxAtStartPar
\sphinxstylestrong{Synopsis}
\begin{quote}

\sphinxAtStartPar
\sphinxcode{\sphinxupquote{ProbBuffer(buff)}}
\end{quote}

\sphinxAtStartPar
\sphinxstylestrong{Description}
\begin{quote}

\sphinxAtStartPar
Create and return a {\hyperref[\detokenize{pyapiref:chappyapi-util-probbuffer}]{\sphinxcrossref{\DUrole{std,std-ref}{ProbBuffer Class}}}} object.
\end{quote}

\sphinxAtStartPar
\sphinxstylestrong{Arguments}
\begin{quote}

\sphinxAtStartPar
\sphinxcode{\sphinxupquote{buff}}
\begin{quote}

\sphinxAtStartPar
Size of buffer, defaults to \sphinxcode{\sphinxupquote{None}}, i.e. the buffer size is 0.
\end{quote}
\end{quote}

\sphinxAtStartPar
\sphinxstylestrong{Example}
\end{quote}

\begin{sphinxVerbatim}[commandchars=\\\{\}]
\PYG{c+c1}{\PYGZsh{} Create a buffer of size 100}
\PYG{n}{buff} \PYG{o}{=} \PYG{n}{ProbBuffer}\PYG{p}{(}\PYG{l+m+mi}{100}\PYG{p}{)}
\end{sphinxVerbatim}

\subsubsection{ProbBuffer.getData()}
\label{\detokenize{pyapiref:probbuffer-getdata}}\begin{quote}

\sphinxAtStartPar
\sphinxstylestrong{Synopsis}
\begin{quote}

\sphinxAtStartPar
\sphinxcode{\sphinxupquote{getData()}}
\end{quote}

\sphinxAtStartPar
\sphinxstylestrong{Description}
\begin{quote}

\sphinxAtStartPar
Get the contents of buffer.
\end{quote}

\sphinxAtStartPar
\sphinxstylestrong{Example}
\end{quote}

\begin{sphinxVerbatim}[commandchars=\\\{\}]
\PYG{c+c1}{\PYGZsh{} Print the contents in buffer}
\PYG{n+nb}{print}\PYG{p}{(}\PYG{n}{buff}\PYG{o}{.}\PYG{n}{getData}\PYG{p}{(}\PYG{p}{)}\PYG{p}{)}
\end{sphinxVerbatim}

\subsubsection{ProbBuffer.getSize()}
\label{\detokenize{pyapiref:probbuffer-getsize}}\begin{quote}

\sphinxAtStartPar
\sphinxstylestrong{Synopsis}
\begin{quote}

\sphinxAtStartPar
\sphinxcode{\sphinxupquote{getSize()}}
\end{quote}

\sphinxAtStartPar
\sphinxstylestrong{Description}
\begin{quote}

\sphinxAtStartPar
Get the size of the buffer.
\end{quote}

\sphinxAtStartPar
\sphinxstylestrong{Example}
\end{quote}

\begin{sphinxVerbatim}[commandchars=\\\{\}]
\PYG{c+c1}{\PYGZsh{} Get the size of the buffer}
\PYG{n+nb}{print}\PYG{p}{(}\PYG{n}{buff}\PYG{o}{.}\PYG{n}{getSize}\PYG{p}{(}\PYG{p}{)}\PYG{p}{)}
\end{sphinxVerbatim}

\subsubsection{ProbBuffer.resize()}
\label{\detokenize{pyapiref:probbuffer-resize}}\begin{quote}

\sphinxAtStartPar
\sphinxstylestrong{Synopsis}
\begin{quote}

\sphinxAtStartPar
\sphinxcode{\sphinxupquote{resize(sz)}}
\end{quote}

\sphinxAtStartPar
\sphinxstylestrong{Description}
\begin{quote}

\sphinxAtStartPar
Resize the size of the buffer.
\end{quote}

\sphinxAtStartPar
\sphinxstylestrong{Arguments}
\begin{quote}

\sphinxAtStartPar
\sphinxcode{\sphinxupquote{sz}}
\begin{quote}

\sphinxAtStartPar
New size of buffer.
\end{quote}
\end{quote}

\sphinxAtStartPar
\sphinxstylestrong{Example}
\end{quote}

\begin{sphinxVerbatim}[commandchars=\\\{\}]
\PYG{c+c1}{\PYGZsh{} Resize the size of buffer to 100}
\PYG{n}{buff}\PYG{o}{.}\PYG{n}{resize}\PYG{p}{(}\PYG{l+m+mi}{100}\PYG{p}{)}
\end{sphinxVerbatim}

\sphinxstepscope

\chapter{C++ API Reference}
\label{\detokenize{cppapiref:c-api-reference}}\label{\detokenize{cppapiref:chapcppapiref}}\label{\detokenize{cppapiref::doc}}
\sphinxAtStartPar
The \sphinxstylestrong{Cardinal Optimizer} provides C++ API library.
This chapter documents all COPT constants, including
parameters and attributes, and API functions for C++ applications.

\section{Constants}
\label{\detokenize{cppapiref:constants}}\label{\detokenize{cppapiref:chapcppapiref-const}}
\sphinxAtStartPar
All C++ constants are the same as C constants. Please refer to
{\hyperref[\detokenize{capiref:chapapi-const}]{\sphinxcrossref{\DUrole{std,std-ref}{C API Reference: Constants}}}} for more details.

\section{Attributes}
\label{\detokenize{cppapiref:attributes}}\label{\detokenize{cppapiref:chapcppapiref-attrs}}
\sphinxAtStartPar
All C++ attributes are the same as C attributes. Please refer to
{\hyperref[\detokenize{capiref:chapapi-attrs}]{\sphinxcrossref{\DUrole{std,std-ref}{C API Reference: Attributes}}}} for more details.

\sphinxAtStartPar
In the C++ API, user can get the attribute value by specifying the attribute
name. The provided functions are as follows, please refer to
{\hyperref[\detokenize{cppapiref:chapcppapiref-model}]{\sphinxcrossref{\DUrole{std,std-ref}{C++ API: Model Class}}}} for details.
\begin{itemize}
\item {} 
\sphinxAtStartPar
\sphinxcode{\sphinxupquote{Model::GetIntAttr()}} : Get value of a COPT integer attribute.

\item {} 
\sphinxAtStartPar
\sphinxcode{\sphinxupquote{Model::GetDblAttr()}} : Get value of a COPT double attribute.

\end{itemize}

\section{Information}
\label{\detokenize{cppapiref:information}}\label{\detokenize{cppapiref:chapcppapiref-info}}
\sphinxAtStartPar
All C++ information is the same as C information. Please refer to
{\hyperref[\detokenize{capiref:chapapi-info}]{\sphinxcrossref{\DUrole{std,std-ref}{C API Reference: Information}}}} .

\sphinxAtStartPar
In the C++ API, user can get or set the information value by specifying the information name.
The functions provided are as follows, please refer to {\hyperref[\detokenize{cppapiref:chapcppapiref-model}]{\sphinxcrossref{\DUrole{std,std-ref}{C++ Model class}}}} for details.
Information of variables or constraints can also be obtained/set through the \sphinxcode{\sphinxupquote{Get()/Set()}} function of themselves.
\begin{itemize}
\item {} 
\sphinxAtStartPar
\sphinxcode{\sphinxupquote{Model::Get()}} : Get the value of information related to the variables or constraints.

\item {} 
\sphinxAtStartPar
\sphinxcode{\sphinxupquote{Model::Set()}} : Set the value of information related to the variables or constraints.

\end{itemize}

\section{Parameters}
\label{\detokenize{cppapiref:parameters}}\label{\detokenize{cppapiref:chapcppapiref-params}}
\sphinxAtStartPar
All C++ parameters are the same as C parameters. Please refer to
{\hyperref[\detokenize{capiref:chapapi-param}]{\sphinxcrossref{\DUrole{std,std-ref}{C API Reference: Parameters}}}} for more details.

\sphinxAtStartPar
In the C++ API, user can get and set the parameter value by specifying the
parameter name. The provided functions are as follows, please refer to
{\hyperref[\detokenize{cppapiref:chapcppapiref-model}]{\sphinxcrossref{\DUrole{std,std-ref}{C++ Model class}}}} for details.
\begin{itemize}
\item {} 
\sphinxAtStartPar
Get detailed information of the specified parameter (current value/max/min):
\sphinxcode{\sphinxupquote{Model::GetParamInfo()}}

\item {} 
\sphinxAtStartPar
Get the current value of the specified integer/double parameter:
\sphinxcode{\sphinxupquote{Model::GetIntParam()}} / \sphinxcode{\sphinxupquote{Model::GetDblParam()}}

\item {} 
\sphinxAtStartPar
Set the specified integer/double parameter value:
\sphinxcode{\sphinxupquote{Model::SetIntParam()}} / \sphinxcode{\sphinxupquote{Model::SetDblParam()}}

\end{itemize}

\section{C++ Modeling Classes}
\label{\detokenize{cppapiref:c-modeling-classes}}\label{\detokenize{cppapiref:chapcppapiref-class}}
\sphinxAtStartPar
This chapter documents COPT C++ interface. Users may refer to
C++ classes described below for details of how to construct and
solve C++ models.

\subsection{Envr}
\label{\detokenize{cppapiref:envr}}\label{\detokenize{cppapiref:chapcppapiref-envr}}
\sphinxAtStartPar
Essentially, any C++ application using Cardinal Optimizer should start with a
COPT environment. COPT models are always associated with a COPT environment.
User must create an environment object before populating models.
User generally only need a single environment object in program.

\sphinxstepscope

\subsubsection{Envr::Envr()}
\label{\detokenize{cppapi/envr:envr-envr}}\label{\detokenize{cppapi/envr::doc}}\begin{quote}

\sphinxAtStartPar
Constructor of COPT Envr object.

\sphinxAtStartPar
\sphinxstylestrong{Synopsis}
\begin{quote}

\sphinxAtStartPar
\sphinxcode{\sphinxupquote{Envr()}}
\end{quote}
\end{quote}

\subsubsection{Envr::Envr()}
\label{\detokenize{cppapi/envr:id1}}\begin{quote}

\sphinxAtStartPar
Constructor of COPT Envr object, given a license folder.

\sphinxAtStartPar
\sphinxstylestrong{Synopsis}
\begin{quote}

\sphinxAtStartPar
\sphinxcode{\sphinxupquote{Envr(const char *szLicDir)}}
\end{quote}

\sphinxAtStartPar
\sphinxstylestrong{Arguments}
\begin{quote}

\sphinxAtStartPar
\sphinxcode{\sphinxupquote{szLicDir}}: directory having local license or client config file.
\end{quote}
\end{quote}

\subsubsection{Envr::Envr()}
\label{\detokenize{cppapi/envr:id2}}\begin{quote}

\sphinxAtStartPar
Constructor of COPT Envr object, given an Envr config object.

\sphinxAtStartPar
\sphinxstylestrong{Synopsis}
\begin{quote}

\sphinxAtStartPar
\sphinxcode{\sphinxupquote{Envr(const EnvrConfig \&config)}}
\end{quote}

\sphinxAtStartPar
\sphinxstylestrong{Arguments}
\begin{quote}

\sphinxAtStartPar
\sphinxcode{\sphinxupquote{config}}: Envr config object holding settings for remote connection.
\end{quote}
\end{quote}

\subsubsection{Envr::BindNumaCpu()}
\label{\detokenize{cppapi/envr:envr-bindnumacpu}}\begin{quote}

\sphinxAtStartPar
Bind the CPUs for the current process to a NUMA node.

\sphinxAtStartPar
\sphinxstylestrong{Synopsis}
\begin{quote}

\sphinxAtStartPar
\sphinxcode{\sphinxupquote{void BindNumaCpu(int numaNode)}}
\end{quote}

\sphinxAtStartPar
\sphinxstylestrong{Arguments}
\begin{quote}

\sphinxAtStartPar
\sphinxcode{\sphinxupquote{numaNode}}: ID of a NUMA node.
\end{quote}
\end{quote}

\subsubsection{Envr::BindNumaMem()}
\label{\detokenize{cppapi/envr:envr-bindnumamem}}\begin{quote}

\sphinxAtStartPar
Bind memory for the current process to a NUMA node (Linux only).

\sphinxAtStartPar
\sphinxstylestrong{Synopsis}
\begin{quote}

\sphinxAtStartPar
\sphinxcode{\sphinxupquote{void BindNumaMem(int numaNode)}}
\end{quote}

\sphinxAtStartPar
\sphinxstylestrong{Arguments}
\begin{quote}

\sphinxAtStartPar
\sphinxcode{\sphinxupquote{numaNode}}: the ID of a NUMA node.
\end{quote}
\end{quote}

\subsubsection{Envr::Close()}
\label{\detokenize{cppapi/envr:envr-close}}\begin{quote}

\sphinxAtStartPar
Close remote connection and token becomes invalid for all problems in current envr.

\sphinxAtStartPar
\sphinxstylestrong{Synopsis}
\begin{quote}

\sphinxAtStartPar
\sphinxcode{\sphinxupquote{void Close()}}
\end{quote}
\end{quote}

\subsubsection{Envr::CreateModel()}
\label{\detokenize{cppapi/envr:envr-createmodel}}\begin{quote}

\sphinxAtStartPar
Create a COPT model object.

\sphinxAtStartPar
\sphinxstylestrong{Synopsis}
\begin{quote}

\sphinxAtStartPar
\sphinxcode{\sphinxupquote{Model CreateModel(const char *szName)}}
\end{quote}

\sphinxAtStartPar
\sphinxstylestrong{Arguments}
\begin{quote}

\sphinxAtStartPar
\sphinxcode{\sphinxupquote{szName}}: customized model name.
\end{quote}

\sphinxAtStartPar
\sphinxstylestrong{Return}
\begin{quote}

\sphinxAtStartPar
a COPT model object.
\end{quote}
\end{quote}

\subsubsection{Envr::GetCpuAffinity()}
\label{\detokenize{cppapi/envr:envr-getcpuaffinity}}\begin{quote}

\sphinxAtStartPar
Get CPU affinity for the current process, which is saved in an integer array.

\sphinxAtStartPar
\sphinxstylestrong{Synopsis}
\begin{quote}

\sphinxAtStartPar
\sphinxcode{\sphinxupquote{int GetCpuAffinity(int *cpuList, int len)}}
\end{quote}

\sphinxAtStartPar
\sphinxstylestrong{Arguments}
\begin{quote}

\sphinxAtStartPar
\sphinxcode{\sphinxupquote{cpuList}}: a list of CPU IDs.

\sphinxAtStartPar
\sphinxcode{\sphinxupquote{len}}: length of the CPU list.
\end{quote}

\sphinxAtStartPar
\sphinxstylestrong{Return}
\begin{quote}

\sphinxAtStartPar
actual size of binding CPUs.
\end{quote}
\end{quote}

\subsubsection{Envr::GetNumaNodeCount()}
\label{\detokenize{cppapi/envr:envr-getnumanodecount}}\begin{quote}

\sphinxAtStartPar
Get count of NUMA nodes.

\sphinxAtStartPar
\sphinxstylestrong{Synopsis}
\begin{quote}

\sphinxAtStartPar
\sphinxcode{\sphinxupquote{int GetNumaNodeCount()}}
\end{quote}

\sphinxAtStartPar
\sphinxstylestrong{Return}
\begin{quote}

\sphinxAtStartPar
count of NUMA nodes.
\end{quote}
\end{quote}

\subsubsection{Envr::SetCpuAffinity()}
\label{\detokenize{cppapi/envr:envr-setcpuaffinity}}\begin{quote}

\sphinxAtStartPar
Set CPU affinity with given mask string.

\sphinxAtStartPar
\sphinxstylestrong{Synopsis}
\begin{quote}

\sphinxAtStartPar
\sphinxcode{\sphinxupquote{void SetCpuAffinity(const char *hexMask)}}
\end{quote}

\sphinxAtStartPar
\sphinxstylestrong{Arguments}
\begin{quote}

\sphinxAtStartPar
\sphinxcode{\sphinxupquote{hexMask}}: CPU mask string of hexadecimal characters.
\end{quote}
\end{quote}

\subsection{EnvrConfig}
\label{\detokenize{cppapiref:envrconfig}}
\sphinxAtStartPar
If user connects to COPT remote services, such as floating token server
or compute cluster, it is necessary to add config settings with EnvrConfig
object.

\sphinxstepscope

\subsubsection{EnvrConfig::EnvrConfig()}
\label{\detokenize{cppapi/envrconfig:envrconfig-envrconfig}}\label{\detokenize{cppapi/envrconfig::doc}}\begin{quote}

\sphinxAtStartPar
Constructor of COPT environment config object.

\sphinxAtStartPar
\sphinxstylestrong{Synopsis}
\begin{quote}

\sphinxAtStartPar
\sphinxcode{\sphinxupquote{EnvrConfig()}}
\end{quote}
\end{quote}

\subsubsection{EnvrConfig::Set()}
\label{\detokenize{cppapi/envrconfig:envrconfig-set}}\begin{quote}

\sphinxAtStartPar
Set config settings in terms of name\sphinxhyphen{}value pair.

\sphinxAtStartPar
\sphinxstylestrong{Synopsis}
\begin{quote}

\sphinxAtStartPar
\sphinxcode{\sphinxupquote{void Set(const char *szName, const char *szValue)}}
\end{quote}

\sphinxAtStartPar
\sphinxstylestrong{Arguments}
\begin{quote}

\sphinxAtStartPar
\sphinxcode{\sphinxupquote{szName}}: keyword of a config setting

\sphinxAtStartPar
\sphinxcode{\sphinxupquote{szValue}}: value of a config setting
\end{quote}
\end{quote}

\subsection{Model}
\label{\detokenize{cppapiref:model}}\label{\detokenize{cppapiref:chapcppapiref-model}}
\sphinxAtStartPar
In general, a COPT model consists of a set of variables, a (linear) objective
function on these variables, a set of constraints on there varaibles, etc. COPT
model class encapsulates all required methods for constructing a COPT model.

\sphinxstepscope

\subsubsection{Model::AddAffineCone()}
\label{\detokenize{cppapi/model:model-addaffinecone}}\label{\detokenize{cppapi/model::doc}}\begin{quote}

\sphinxAtStartPar
Add an affine cone constraint to model.

\sphinxAtStartPar
\sphinxstylestrong{Synopsis}
\begin{quote}

\sphinxAtStartPar
\sphinxcode{\sphinxupquote{AffineCone AddAffineCone(const AffineConeBuilder \&builder, const char *szName)}}
\end{quote}

\sphinxAtStartPar
\sphinxstylestrong{Arguments}
\begin{quote}

\sphinxAtStartPar
\sphinxcode{\sphinxupquote{builder}}: builder for new affine cone constraint.

\sphinxAtStartPar
\sphinxcode{\sphinxupquote{szName}}: optional, name of new affine cone constraint.
\end{quote}

\sphinxAtStartPar
\sphinxstylestrong{Return}
\begin{quote}

\sphinxAtStartPar
new affine cone constraint object.
\end{quote}
\end{quote}

\subsubsection{Model::AddAffineCone()}
\label{\detokenize{cppapi/model:id1}}\begin{quote}

\sphinxAtStartPar
Add an affine cone constraint to a model.

\sphinxAtStartPar
\sphinxstylestrong{Synopsis}
\begin{quote}

\sphinxAtStartPar
\sphinxcode{\sphinxupquote{AffineCone AddAffineCone(}}
\begin{quote}

\sphinxAtStartPar
\sphinxcode{\sphinxupquote{const MLinExpr\textless{}1\textgreater{} \&exprs,}}

\sphinxAtStartPar
\sphinxcode{\sphinxupquote{int type,}}

\sphinxAtStartPar
\sphinxcode{\sphinxupquote{const char *szName)}}
\end{quote}
\end{quote}

\sphinxAtStartPar
\sphinxstylestrong{Arguments}
\begin{quote}

\sphinxAtStartPar
\sphinxcode{\sphinxupquote{exprs}}: 1\sphinxhyphen{}dimensional array of linear expressions.

\sphinxAtStartPar
\sphinxcode{\sphinxupquote{type}}: type of an affine cone.

\sphinxAtStartPar
\sphinxcode{\sphinxupquote{szName}}: name of new affine cone constraint.
\end{quote}

\sphinxAtStartPar
\sphinxstylestrong{Return}
\begin{quote}

\sphinxAtStartPar
new affine cone constraint object.
\end{quote}
\end{quote}

\subsubsection{Model::AddAffineCone()}
\label{\detokenize{cppapi/model:id2}}\begin{quote}

\sphinxAtStartPar
Add an affine cone constraint to a model.

\sphinxAtStartPar
\sphinxstylestrong{Synopsis}
\begin{quote}

\sphinxAtStartPar
\sphinxcode{\sphinxupquote{AffineCone AddAffineCone(}}
\begin{quote}

\sphinxAtStartPar
\sphinxcode{\sphinxupquote{const MPsdExpr\textless{}1\textgreater{} \&exprs,}}

\sphinxAtStartPar
\sphinxcode{\sphinxupquote{int type,}}

\sphinxAtStartPar
\sphinxcode{\sphinxupquote{const char *szName)}}
\end{quote}
\end{quote}

\sphinxAtStartPar
\sphinxstylestrong{Arguments}
\begin{quote}

\sphinxAtStartPar
\sphinxcode{\sphinxupquote{exprs}}: 1\sphinxhyphen{}dimensional array of PSD expressions.

\sphinxAtStartPar
\sphinxcode{\sphinxupquote{type}}: type of an affine cone.

\sphinxAtStartPar
\sphinxcode{\sphinxupquote{szName}}: name of new affine cone constraint.
\end{quote}

\sphinxAtStartPar
\sphinxstylestrong{Return}
\begin{quote}

\sphinxAtStartPar
new affine cone constraint object.
\end{quote}
\end{quote}

\subsubsection{Model::AddAffineCones()}
\label{\detokenize{cppapi/model:model-addaffinecones}}\begin{quote}

\sphinxAtStartPar
Add a list of affine cone constraints to a model.

\sphinxAtStartPar
\sphinxstylestrong{Synopsis}
\begin{quote}

\sphinxAtStartPar
\sphinxcode{\sphinxupquote{AffineConeArray AddAffineCones(}}
\begin{quote}

\sphinxAtStartPar
\sphinxcode{\sphinxupquote{const MPsdExpr\textless{}2\textgreater{} \&exprs,}}

\sphinxAtStartPar
\sphinxcode{\sphinxupquote{int type,}}

\sphinxAtStartPar
\sphinxcode{\sphinxupquote{const char *szPrefix)}}
\end{quote}
\end{quote}

\sphinxAtStartPar
\sphinxstylestrong{Arguments}
\begin{quote}

\sphinxAtStartPar
\sphinxcode{\sphinxupquote{exprs}}: 2\sphinxhyphen{}dimensional array of PSD expressions.

\sphinxAtStartPar
\sphinxcode{\sphinxupquote{type}}: type of an affine cone.

\sphinxAtStartPar
\sphinxcode{\sphinxupquote{szPrefix}}: name prefix for new affine cone constraints.
\end{quote}

\sphinxAtStartPar
\sphinxstylestrong{Return}
\begin{quote}

\sphinxAtStartPar
array of new affine cone constraint objects.
\end{quote}
\end{quote}

\subsubsection{Model::AddAffineCones()}
\label{\detokenize{cppapi/model:id3}}\begin{quote}

\sphinxAtStartPar
Add a list of affine cone constraints to a model.

\sphinxAtStartPar
\sphinxstylestrong{Synopsis}
\begin{quote}

\sphinxAtStartPar
\sphinxcode{\sphinxupquote{AffineConeArray AddAffineCones(}}
\begin{quote}

\sphinxAtStartPar
\sphinxcode{\sphinxupquote{const MLinExpr\textless{}2\textgreater{} \&exprs,}}

\sphinxAtStartPar
\sphinxcode{\sphinxupquote{int type,}}

\sphinxAtStartPar
\sphinxcode{\sphinxupquote{const char *szPrefix)}}
\end{quote}
\end{quote}

\sphinxAtStartPar
\sphinxstylestrong{Arguments}
\begin{quote}

\sphinxAtStartPar
\sphinxcode{\sphinxupquote{exprs}}: 2\sphinxhyphen{}dimensional array of linear expressions.

\sphinxAtStartPar
\sphinxcode{\sphinxupquote{type}}: type of an affine cone.

\sphinxAtStartPar
\sphinxcode{\sphinxupquote{szPrefix}}: name prefix for new affine cone constraints.
\end{quote}

\sphinxAtStartPar
\sphinxstylestrong{Return}
\begin{quote}

\sphinxAtStartPar
array of new affine cone constraint objects.
\end{quote}
\end{quote}

\subsubsection{Model::AddCol()}
\label{\detokenize{cppapi/model:model-addcol}}\begin{quote}

\sphinxAtStartPar
With column data, add varaible to model via advanced interface.

\sphinxAtStartPar
\sphinxstylestrong{Synopsis}
\begin{quote}

\sphinxAtStartPar
\sphinxcode{\sphinxupquote{Var AddCol(}}
\begin{quote}

\sphinxAtStartPar
\sphinxcode{\sphinxupquote{double obj,}}

\sphinxAtStartPar
\sphinxcode{\sphinxupquote{int colCnt,}}

\sphinxAtStartPar
\sphinxcode{\sphinxupquote{const int *colIdx,}}

\sphinxAtStartPar
\sphinxcode{\sphinxupquote{const double *colElem,}}

\sphinxAtStartPar
\sphinxcode{\sphinxupquote{char vtype,}}

\sphinxAtStartPar
\sphinxcode{\sphinxupquote{double lb,}}

\sphinxAtStartPar
\sphinxcode{\sphinxupquote{double ub,}}

\sphinxAtStartPar
\sphinxcode{\sphinxupquote{const char *szName)}}
\end{quote}
\end{quote}

\sphinxAtStartPar
\sphinxstylestrong{Arguments}
\begin{quote}

\sphinxAtStartPar
\sphinxcode{\sphinxupquote{obj}}: coefficient of variable in objective function.

\sphinxAtStartPar
\sphinxcode{\sphinxupquote{colCnt}}: number of terms in the column data.

\sphinxAtStartPar
\sphinxcode{\sphinxupquote{colIdx}}: array of constraint indexes in the column data.

\sphinxAtStartPar
\sphinxcode{\sphinxupquote{colElem}}: array of coefficients in the column data.

\sphinxAtStartPar
\sphinxcode{\sphinxupquote{vtype}}: variable type.

\sphinxAtStartPar
\sphinxcode{\sphinxupquote{lb}}: lower bound of new variable.

\sphinxAtStartPar
\sphinxcode{\sphinxupquote{ub}}: upper bound of new variable.

\sphinxAtStartPar
\sphinxcode{\sphinxupquote{szName}}: name of new vaiable, with default value of empty.
\end{quote}

\sphinxAtStartPar
\sphinxstylestrong{Return}
\begin{quote}

\sphinxAtStartPar
new variable object.
\end{quote}
\end{quote}

\subsubsection{Model::AddCols()}
\label{\detokenize{cppapi/model:model-addcols}}\begin{quote}

\sphinxAtStartPar
Add variables to model via advanced interface.

\sphinxAtStartPar
\sphinxstylestrong{Synopsis}
\begin{quote}

\sphinxAtStartPar
\sphinxcode{\sphinxupquote{VarArray AddCols(}}
\begin{quote}

\sphinxAtStartPar
\sphinxcode{\sphinxupquote{int count,}}

\sphinxAtStartPar
\sphinxcode{\sphinxupquote{double *pobj,}}

\sphinxAtStartPar
\sphinxcode{\sphinxupquote{char *pvtype,}}

\sphinxAtStartPar
\sphinxcode{\sphinxupquote{double *plb,}}

\sphinxAtStartPar
\sphinxcode{\sphinxupquote{double *pub,}}

\sphinxAtStartPar
\sphinxcode{\sphinxupquote{char const* const *names)}}
\end{quote}
\end{quote}

\sphinxAtStartPar
\sphinxstylestrong{Arguments}
\begin{quote}

\sphinxAtStartPar
\sphinxcode{\sphinxupquote{count}}: number of varaibles added to model.

\sphinxAtStartPar
\sphinxcode{\sphinxupquote{pobj}}: coefficients of new variables in objective function.

\sphinxAtStartPar
\sphinxcode{\sphinxupquote{pvtype}}: variable types. If empty, all types are continuous.

\sphinxAtStartPar
\sphinxcode{\sphinxupquote{plb}}: lower bounds for new variables. If empty, lower bounds are set to zero.

\sphinxAtStartPar
\sphinxcode{\sphinxupquote{pub}}: upper bounds for new varaibles. If empty, upper bounds are set to infinity.

\sphinxAtStartPar
\sphinxcode{\sphinxupquote{names}}: array of names of new variables, with default value of empty.
\end{quote}

\sphinxAtStartPar
\sphinxstylestrong{Return}
\begin{quote}

\sphinxAtStartPar
new variable array objects.
\end{quote}
\end{quote}

\subsubsection{Model::AddCols()}
\label{\detokenize{cppapi/model:id4}}\begin{quote}

\sphinxAtStartPar
With coefficient matrix data in CSC format, add variables to model via advanced interface.

\sphinxAtStartPar
\sphinxstylestrong{Synopsis}
\begin{quote}

\sphinxAtStartPar
\sphinxcode{\sphinxupquote{VarArray AddCols(}}
\begin{quote}

\sphinxAtStartPar
\sphinxcode{\sphinxupquote{int count,}}

\sphinxAtStartPar
\sphinxcode{\sphinxupquote{double *pobj,}}

\sphinxAtStartPar
\sphinxcode{\sphinxupquote{int *colBeg,}}

\sphinxAtStartPar
\sphinxcode{\sphinxupquote{int *colCnt,}}

\sphinxAtStartPar
\sphinxcode{\sphinxupquote{int *colIdx,}}

\sphinxAtStartPar
\sphinxcode{\sphinxupquote{double *colElem,}}

\sphinxAtStartPar
\sphinxcode{\sphinxupquote{char *pvtype,}}

\sphinxAtStartPar
\sphinxcode{\sphinxupquote{double *plb,}}

\sphinxAtStartPar
\sphinxcode{\sphinxupquote{double *pub,}}

\sphinxAtStartPar
\sphinxcode{\sphinxupquote{char const* const *names)}}
\end{quote}
\end{quote}

\sphinxAtStartPar
\sphinxstylestrong{Arguments}
\begin{quote}

\sphinxAtStartPar
\sphinxcode{\sphinxupquote{count}}: number of varaibles added to model.

\sphinxAtStartPar
\sphinxcode{\sphinxupquote{pobj}}: coefficients of new variables in objective function.

\sphinxAtStartPar
\sphinxcode{\sphinxupquote{colBeg}}: indexes of begin elements in CSC format. If empty, variables are added without column data.

\sphinxAtStartPar
\sphinxcode{\sphinxupquote{colCnt}}: count of nonzero elements in each column. If empty, use colBeg to calcuate count.

\sphinxAtStartPar
\sphinxcode{\sphinxupquote{colIdx}}: constraint indexes of columns in CSC format.

\sphinxAtStartPar
\sphinxcode{\sphinxupquote{colElem}}: corresponding constraint coefficients of columns in CSC format.

\sphinxAtStartPar
\sphinxcode{\sphinxupquote{pvtype}}: variable types. If empty, all types are continuous.

\sphinxAtStartPar
\sphinxcode{\sphinxupquote{plb}}: lower bounds for new variables. If empty, lower bounds are set to zero.

\sphinxAtStartPar
\sphinxcode{\sphinxupquote{pub}}: upper bounds for new varaibles. If empty, upper bounds are set to infinity.

\sphinxAtStartPar
\sphinxcode{\sphinxupquote{names}}: array of names of new variables, with default value of empty.
\end{quote}

\sphinxAtStartPar
\sphinxstylestrong{Return}
\begin{quote}

\sphinxAtStartPar
new variable array objects.
\end{quote}
\end{quote}

\subsubsection{Model::AddCone()}
\label{\detokenize{cppapi/model:model-addcone}}\begin{quote}

\sphinxAtStartPar
Add a cone constraint to a model, given its dimension.

\sphinxAtStartPar
\sphinxstylestrong{Synopsis}
\begin{quote}

\sphinxAtStartPar
\sphinxcode{\sphinxupquote{Cone AddCone(}}
\begin{quote}

\sphinxAtStartPar
\sphinxcode{\sphinxupquote{int dim,}}

\sphinxAtStartPar
\sphinxcode{\sphinxupquote{int type,}}

\sphinxAtStartPar
\sphinxcode{\sphinxupquote{char *pvtype,}}

\sphinxAtStartPar
\sphinxcode{\sphinxupquote{const char *szPrefix)}}
\end{quote}
\end{quote}

\sphinxAtStartPar
\sphinxstylestrong{Arguments}
\begin{quote}

\sphinxAtStartPar
\sphinxcode{\sphinxupquote{dim}}: dimension of the cone constraint.

\sphinxAtStartPar
\sphinxcode{\sphinxupquote{type}}: type of the cone constraint.

\sphinxAtStartPar
\sphinxcode{\sphinxupquote{pvtype}}: types of variables in the cone.

\sphinxAtStartPar
\sphinxcode{\sphinxupquote{szPrefix}}: name prefix of variables in the cone.
\end{quote}

\sphinxAtStartPar
\sphinxstylestrong{Return}
\begin{quote}

\sphinxAtStartPar
object of new cone constraint.
\end{quote}
\end{quote}

\subsubsection{Model::AddCone()}
\label{\detokenize{cppapi/model:id5}}\begin{quote}

\sphinxAtStartPar
Add a cone constraint to model.

\sphinxAtStartPar
\sphinxstylestrong{Synopsis}
\begin{quote}

\sphinxAtStartPar
\sphinxcode{\sphinxupquote{Cone AddCone(const ConeBuilder \&builder)}}
\end{quote}

\sphinxAtStartPar
\sphinxstylestrong{Arguments}
\begin{quote}

\sphinxAtStartPar
\sphinxcode{\sphinxupquote{builder}}: builder for new cone constraint.
\end{quote}

\sphinxAtStartPar
\sphinxstylestrong{Return}
\begin{quote}

\sphinxAtStartPar
new cone constraint object.
\end{quote}
\end{quote}

\subsubsection{Model::AddCone()}
\label{\detokenize{cppapi/model:id6}}\begin{quote}

\sphinxAtStartPar
Add a cone constraint to model.

\sphinxAtStartPar
\sphinxstylestrong{Synopsis}
\begin{quote}

\sphinxAtStartPar
\sphinxcode{\sphinxupquote{Cone AddCone(const VarArray \&vars, int type)}}
\end{quote}

\sphinxAtStartPar
\sphinxstylestrong{Arguments}
\begin{quote}

\sphinxAtStartPar
\sphinxcode{\sphinxupquote{vars}}: variables that participate in the cone constraint.

\sphinxAtStartPar
\sphinxcode{\sphinxupquote{type}}: type of the cone constraint.
\end{quote}

\sphinxAtStartPar
\sphinxstylestrong{Return}
\begin{quote}

\sphinxAtStartPar
object of new cone constraint.
\end{quote}
\end{quote}

\subsubsection{Model::AddCone()}
\label{\detokenize{cppapi/model:id7}}\begin{quote}

\sphinxAtStartPar
Add a cone constraint to model.

\sphinxAtStartPar
\sphinxstylestrong{Synopsis}
\begin{quote}

\sphinxAtStartPar
\sphinxcode{\sphinxupquote{Cone AddCone(const MVar\textless{}1\textgreater{} \&vars, int type)}}
\end{quote}

\sphinxAtStartPar
\sphinxstylestrong{Arguments}
\begin{quote}

\sphinxAtStartPar
\sphinxcode{\sphinxupquote{vars}}: one\sphinxhyphen{}dimensional variables in the cone constraint.

\sphinxAtStartPar
\sphinxcode{\sphinxupquote{type}}: type of the cone constraint.
\end{quote}

\sphinxAtStartPar
\sphinxstylestrong{Return}
\begin{quote}

\sphinxAtStartPar
object of new cone constraint.
\end{quote}
\end{quote}

\subsubsection{Model::AddCones()}
\label{\detokenize{cppapi/model:model-addcones}}\begin{quote}

\sphinxAtStartPar
Add a list of cone constraints to a model.

\sphinxAtStartPar
\sphinxstylestrong{Synopsis}
\begin{quote}

\sphinxAtStartPar
\sphinxcode{\sphinxupquote{ConeArray AddCones(const MVar\textless{}2\textgreater{} \&vars, int type)}}
\end{quote}

\sphinxAtStartPar
\sphinxstylestrong{Arguments}
\begin{quote}

\sphinxAtStartPar
\sphinxcode{\sphinxupquote{vars}}: 2\sphinxhyphen{}dimensional array of variables.

\sphinxAtStartPar
\sphinxcode{\sphinxupquote{type}}: type of a cone.
\end{quote}

\sphinxAtStartPar
\sphinxstylestrong{Return}
\begin{quote}

\sphinxAtStartPar
array of new cone constraint objects.
\end{quote}
\end{quote}

\subsubsection{Model::AddConstr()}
\label{\detokenize{cppapi/model:model-addconstr}}\begin{quote}

\sphinxAtStartPar
Add a linear constraint to model.

\sphinxAtStartPar
\sphinxstylestrong{Synopsis}
\begin{quote}

\sphinxAtStartPar
\sphinxcode{\sphinxupquote{Constraint AddConstr(}}
\begin{quote}

\sphinxAtStartPar
\sphinxcode{\sphinxupquote{const Expr \&expr,}}

\sphinxAtStartPar
\sphinxcode{\sphinxupquote{char sense,}}

\sphinxAtStartPar
\sphinxcode{\sphinxupquote{double rhs,}}

\sphinxAtStartPar
\sphinxcode{\sphinxupquote{const char *szName)}}
\end{quote}
\end{quote}

\sphinxAtStartPar
\sphinxstylestrong{Arguments}
\begin{quote}

\sphinxAtStartPar
\sphinxcode{\sphinxupquote{expr}}: expression for the new contraint.

\sphinxAtStartPar
\sphinxcode{\sphinxupquote{sense}}: sense for new linear constraint, other than range sense.

\sphinxAtStartPar
\sphinxcode{\sphinxupquote{rhs}}: right hand side value for the new constraint.

\sphinxAtStartPar
\sphinxcode{\sphinxupquote{szName}}: optional, name of new constraint.
\end{quote}

\sphinxAtStartPar
\sphinxstylestrong{Return}
\begin{quote}

\sphinxAtStartPar
new constraint object.
\end{quote}
\end{quote}

\subsubsection{Model::AddConstr()}
\label{\detokenize{cppapi/model:id8}}\begin{quote}

\sphinxAtStartPar
Add a linear constraint to model.

\sphinxAtStartPar
\sphinxstylestrong{Synopsis}
\begin{quote}

\sphinxAtStartPar
\sphinxcode{\sphinxupquote{Constraint AddConstr(}}
\begin{quote}

\sphinxAtStartPar
\sphinxcode{\sphinxupquote{const Expr \&lhs,}}

\sphinxAtStartPar
\sphinxcode{\sphinxupquote{char sense,}}

\sphinxAtStartPar
\sphinxcode{\sphinxupquote{const Expr \&rhs,}}

\sphinxAtStartPar
\sphinxcode{\sphinxupquote{const char *szName)}}
\end{quote}
\end{quote}

\sphinxAtStartPar
\sphinxstylestrong{Arguments}
\begin{quote}

\sphinxAtStartPar
\sphinxcode{\sphinxupquote{lhs}}: left hand side expression for the new constraint.

\sphinxAtStartPar
\sphinxcode{\sphinxupquote{sense}}: sense for new linear constraint, other than range sense.

\sphinxAtStartPar
\sphinxcode{\sphinxupquote{rhs}}: right hand side expression for the new constraint.

\sphinxAtStartPar
\sphinxcode{\sphinxupquote{szName}}: optional, name of new constraint.
\end{quote}

\sphinxAtStartPar
\sphinxstylestrong{Return}
\begin{quote}

\sphinxAtStartPar
new constraint object.
\end{quote}
\end{quote}

\subsubsection{Model::AddConstr()}
\label{\detokenize{cppapi/model:id9}}\begin{quote}

\sphinxAtStartPar
Add a linear constraint to model.

\sphinxAtStartPar
\sphinxstylestrong{Synopsis}
\begin{quote}

\sphinxAtStartPar
\sphinxcode{\sphinxupquote{Constraint AddConstr(}}
\begin{quote}

\sphinxAtStartPar
\sphinxcode{\sphinxupquote{const Expr \&expr,}}

\sphinxAtStartPar
\sphinxcode{\sphinxupquote{double lb,}}

\sphinxAtStartPar
\sphinxcode{\sphinxupquote{double ub,}}

\sphinxAtStartPar
\sphinxcode{\sphinxupquote{const char *szName)}}
\end{quote}
\end{quote}

\sphinxAtStartPar
\sphinxstylestrong{Arguments}
\begin{quote}

\sphinxAtStartPar
\sphinxcode{\sphinxupquote{expr}}: expression for the new constraint.

\sphinxAtStartPar
\sphinxcode{\sphinxupquote{lb}}: lower bound for the new constraint.

\sphinxAtStartPar
\sphinxcode{\sphinxupquote{ub}}: upper bound for the new constraint

\sphinxAtStartPar
\sphinxcode{\sphinxupquote{szName}}: optional, name of new constraint.
\end{quote}

\sphinxAtStartPar
\sphinxstylestrong{Return}
\begin{quote}

\sphinxAtStartPar
new constraint object.
\end{quote}
\end{quote}

\subsubsection{Model::AddConstr()}
\label{\detokenize{cppapi/model:id10}}\begin{quote}

\sphinxAtStartPar
Add a linear constraint to a model.

\sphinxAtStartPar
\sphinxstylestrong{Synopsis}
\begin{quote}

\sphinxAtStartPar
\sphinxcode{\sphinxupquote{Constraint AddConstr(const ConstrBuilder \&builder, const char *szName)}}
\end{quote}

\sphinxAtStartPar
\sphinxstylestrong{Arguments}
\begin{quote}

\sphinxAtStartPar
\sphinxcode{\sphinxupquote{builder}}: builder for the new constraint.

\sphinxAtStartPar
\sphinxcode{\sphinxupquote{szName}}: optional, name of new constraint.
\end{quote}

\sphinxAtStartPar
\sphinxstylestrong{Return}
\begin{quote}

\sphinxAtStartPar
new constraint object.
\end{quote}
\end{quote}

\subsubsection{Model::AddConstrs()}
\label{\detokenize{cppapi/model:model-addconstrs}}\begin{quote}

\sphinxAtStartPar
Add linear constraints to model.

\sphinxAtStartPar
\sphinxstylestrong{Synopsis}
\begin{quote}

\sphinxAtStartPar
\sphinxcode{\sphinxupquote{ConstrArray AddConstrs(}}
\begin{quote}

\sphinxAtStartPar
\sphinxcode{\sphinxupquote{int count,}}

\sphinxAtStartPar
\sphinxcode{\sphinxupquote{char *pSense,}}

\sphinxAtStartPar
\sphinxcode{\sphinxupquote{double *pRhs,}}

\sphinxAtStartPar
\sphinxcode{\sphinxupquote{const char *szPrefix)}}
\end{quote}
\end{quote}

\sphinxAtStartPar
\sphinxstylestrong{Arguments}
\begin{quote}

\sphinxAtStartPar
\sphinxcode{\sphinxupquote{count}}: number of constraints added to model.

\sphinxAtStartPar
\sphinxcode{\sphinxupquote{pSense}}: sense array for new linear constraints, other than range sense.

\sphinxAtStartPar
\sphinxcode{\sphinxupquote{pRhs}}: right hand side values for new constraints.

\sphinxAtStartPar
\sphinxcode{\sphinxupquote{szPrefix}}: name prefix for new constraints.
\end{quote}

\sphinxAtStartPar
\sphinxstylestrong{Return}
\begin{quote}

\sphinxAtStartPar
array of new constraint objects.
\end{quote}
\end{quote}

\subsubsection{Model::AddConstrs()}
\label{\detokenize{cppapi/model:id11}}\begin{quote}

\sphinxAtStartPar
Add linear constraints to a model.

\sphinxAtStartPar
\sphinxstylestrong{Synopsis}
\begin{quote}

\sphinxAtStartPar
\sphinxcode{\sphinxupquote{ConstrArray AddConstrs(}}
\begin{quote}

\sphinxAtStartPar
\sphinxcode{\sphinxupquote{int count,}}

\sphinxAtStartPar
\sphinxcode{\sphinxupquote{double *pLower,}}

\sphinxAtStartPar
\sphinxcode{\sphinxupquote{double *pUpper,}}

\sphinxAtStartPar
\sphinxcode{\sphinxupquote{const char *szPrefix)}}
\end{quote}
\end{quote}

\sphinxAtStartPar
\sphinxstylestrong{Arguments}
\begin{quote}

\sphinxAtStartPar
\sphinxcode{\sphinxupquote{count}}: number of constraints added to the model.

\sphinxAtStartPar
\sphinxcode{\sphinxupquote{pLower}}: lower bounds of new constraints.

\sphinxAtStartPar
\sphinxcode{\sphinxupquote{pUpper}}: upper bounds of new constraints.

\sphinxAtStartPar
\sphinxcode{\sphinxupquote{szPrefix}}: name prefix for new constraints.
\end{quote}

\sphinxAtStartPar
\sphinxstylestrong{Return}
\begin{quote}

\sphinxAtStartPar
array of new constraint objects.
\end{quote}
\end{quote}

\subsubsection{Model::AddConstrs()}
\label{\detokenize{cppapi/model:id12}}\begin{quote}

\sphinxAtStartPar
Add linear constraints to a model.

\sphinxAtStartPar
\sphinxstylestrong{Synopsis}
\begin{quote}

\sphinxAtStartPar
\sphinxcode{\sphinxupquote{ConstrArray AddConstrs(}}
\begin{quote}

\sphinxAtStartPar
\sphinxcode{\sphinxupquote{int count,}}

\sphinxAtStartPar
\sphinxcode{\sphinxupquote{double *pLower,}}

\sphinxAtStartPar
\sphinxcode{\sphinxupquote{double *pUpper,}}

\sphinxAtStartPar
\sphinxcode{\sphinxupquote{const char *szNames,}}

\sphinxAtStartPar
\sphinxcode{\sphinxupquote{size\_t len)}}
\end{quote}
\end{quote}

\sphinxAtStartPar
\sphinxstylestrong{Arguments}
\begin{quote}

\sphinxAtStartPar
\sphinxcode{\sphinxupquote{count}}: number of constraints added to the model.

\sphinxAtStartPar
\sphinxcode{\sphinxupquote{pLower}}: lower bounds of new constraints.

\sphinxAtStartPar
\sphinxcode{\sphinxupquote{pUpper}}: upper bounds of new constraints.

\sphinxAtStartPar
\sphinxcode{\sphinxupquote{szNames}}: name buffer of new constraints.

\sphinxAtStartPar
\sphinxcode{\sphinxupquote{len}}: length of the name buffer.
\end{quote}

\sphinxAtStartPar
\sphinxstylestrong{Return}
\begin{quote}

\sphinxAtStartPar
array of new constraint objects.
\end{quote}
\end{quote}

\subsubsection{Model::AddConstrs()}
\label{\detokenize{cppapi/model:id13}}\begin{quote}

\sphinxAtStartPar
Add linear constraints to a model.

\sphinxAtStartPar
\sphinxstylestrong{Synopsis}
\begin{quote}

\sphinxAtStartPar
\sphinxcode{\sphinxupquote{ConstrArray AddConstrs(const ConstrBuilderArray \&builders, const char *szPrefix)}}
\end{quote}

\sphinxAtStartPar
\sphinxstylestrong{Arguments}
\begin{quote}

\sphinxAtStartPar
\sphinxcode{\sphinxupquote{builders}}: builders for new constraints.

\sphinxAtStartPar
\sphinxcode{\sphinxupquote{szPrefix}}: name prefix for new constraints.
\end{quote}

\sphinxAtStartPar
\sphinxstylestrong{Return}
\begin{quote}

\sphinxAtStartPar
array of new constraint objects.
\end{quote}
\end{quote}

\subsubsection{Model::AddConstrs()}
\label{\detokenize{cppapi/model:id14}}\begin{quote}

\sphinxAtStartPar
Add linear constraints to model.

\sphinxAtStartPar
\sphinxstylestrong{Synopsis}
\begin{quote}

\sphinxAtStartPar
\sphinxcode{\sphinxupquote{ConstrArray AddConstrs(}}
\begin{quote}

\sphinxAtStartPar
\sphinxcode{\sphinxupquote{const ConstrBuilderArray \&builders,}}

\sphinxAtStartPar
\sphinxcode{\sphinxupquote{const char *szNames,}}

\sphinxAtStartPar
\sphinxcode{\sphinxupquote{size\_t len)}}
\end{quote}
\end{quote}

\sphinxAtStartPar
\sphinxstylestrong{Arguments}
\begin{quote}

\sphinxAtStartPar
\sphinxcode{\sphinxupquote{builders}}: builders for new constraints.

\sphinxAtStartPar
\sphinxcode{\sphinxupquote{szNames}}: name buffer of new constraints.

\sphinxAtStartPar
\sphinxcode{\sphinxupquote{len}}: length of the name buffer.
\end{quote}

\sphinxAtStartPar
\sphinxstylestrong{Return}
\begin{quote}

\sphinxAtStartPar
array of new constraint objects.
\end{quote}
\end{quote}

\subsubsection{Model::AddDenseMat()}
\label{\detokenize{cppapi/model:model-adddensemat}}\begin{quote}

\sphinxAtStartPar
Add a dense symmetric matrix to a model.

\sphinxAtStartPar
\sphinxstylestrong{Synopsis}
\begin{quote}

\sphinxAtStartPar
\sphinxcode{\sphinxupquote{SymMatrix AddDenseMat(}}
\begin{quote}

\sphinxAtStartPar
\sphinxcode{\sphinxupquote{int dim,}}

\sphinxAtStartPar
\sphinxcode{\sphinxupquote{double *pVals,}}

\sphinxAtStartPar
\sphinxcode{\sphinxupquote{int len)}}
\end{quote}
\end{quote}

\sphinxAtStartPar
\sphinxstylestrong{Arguments}
\begin{quote}

\sphinxAtStartPar
\sphinxcode{\sphinxupquote{dim}}: dimension of the dense symmetric matrix.

\sphinxAtStartPar
\sphinxcode{\sphinxupquote{pVals}}: array of non\sphinxhyphen{}zero elements, filled in column\sphinxhyphen{}wise up to len or max length of symmetric matrix.

\sphinxAtStartPar
\sphinxcode{\sphinxupquote{len}}: length of value array.
\end{quote}

\sphinxAtStartPar
\sphinxstylestrong{Return}
\begin{quote}

\sphinxAtStartPar
new symmetric matrix object.
\end{quote}
\end{quote}

\subsubsection{Model::AddDenseMat()}
\label{\detokenize{cppapi/model:id15}}\begin{quote}

\sphinxAtStartPar
Add a dense symmetric matrix to a model.

\sphinxAtStartPar
\sphinxstylestrong{Synopsis}
\begin{quote}

\sphinxAtStartPar
\sphinxcode{\sphinxupquote{SymMatrix AddDenseMat(int dim, double val)}}
\end{quote}

\sphinxAtStartPar
\sphinxstylestrong{Arguments}
\begin{quote}

\sphinxAtStartPar
\sphinxcode{\sphinxupquote{dim}}: dimension of dense symmetric matrix.

\sphinxAtStartPar
\sphinxcode{\sphinxupquote{val}}: value to fill dense symmetric matrix.
\end{quote}

\sphinxAtStartPar
\sphinxstylestrong{Return}
\begin{quote}

\sphinxAtStartPar
new symmetric matrix object.
\end{quote}
\end{quote}

\subsubsection{Model::AddDiagMat()}
\label{\detokenize{cppapi/model:model-adddiagmat}}\begin{quote}

\sphinxAtStartPar
Add a diagonal matrix to a model.

\sphinxAtStartPar
\sphinxstylestrong{Synopsis}
\begin{quote}

\sphinxAtStartPar
\sphinxcode{\sphinxupquote{SymMatrix AddDiagMat(int dim, double val)}}
\end{quote}

\sphinxAtStartPar
\sphinxstylestrong{Arguments}
\begin{quote}

\sphinxAtStartPar
\sphinxcode{\sphinxupquote{dim}}: dimension of diagonal matrix.

\sphinxAtStartPar
\sphinxcode{\sphinxupquote{val}}: value to fill diagonal elements.
\end{quote}

\sphinxAtStartPar
\sphinxstylestrong{Return}
\begin{quote}

\sphinxAtStartPar
new diagonal matrix object.
\end{quote}
\end{quote}

\subsubsection{Model::AddDiagMat()}
\label{\detokenize{cppapi/model:id16}}\begin{quote}

\sphinxAtStartPar
Add a diagonal matrix to a model.

\sphinxAtStartPar
\sphinxstylestrong{Synopsis}
\begin{quote}

\sphinxAtStartPar
\sphinxcode{\sphinxupquote{SymMatrix AddDiagMat(}}
\begin{quote}

\sphinxAtStartPar
\sphinxcode{\sphinxupquote{int dim,}}

\sphinxAtStartPar
\sphinxcode{\sphinxupquote{double *pVals,}}

\sphinxAtStartPar
\sphinxcode{\sphinxupquote{int len)}}
\end{quote}
\end{quote}

\sphinxAtStartPar
\sphinxstylestrong{Arguments}
\begin{quote}

\sphinxAtStartPar
\sphinxcode{\sphinxupquote{dim}}: dimension of diagonal matrix.

\sphinxAtStartPar
\sphinxcode{\sphinxupquote{pVals}}: array of values of diagonal elements.

\sphinxAtStartPar
\sphinxcode{\sphinxupquote{len}}: length of value array.
\end{quote}

\sphinxAtStartPar
\sphinxstylestrong{Return}
\begin{quote}

\sphinxAtStartPar
new diagonal matrix object.
\end{quote}
\end{quote}

\subsubsection{Model::AddDiagMat()}
\label{\detokenize{cppapi/model:id17}}\begin{quote}

\sphinxAtStartPar
Add a diagonal matrix to a model.

\sphinxAtStartPar
\sphinxstylestrong{Synopsis}
\begin{quote}

\sphinxAtStartPar
\sphinxcode{\sphinxupquote{SymMatrix AddDiagMat(}}
\begin{quote}

\sphinxAtStartPar
\sphinxcode{\sphinxupquote{int dim,}}

\sphinxAtStartPar
\sphinxcode{\sphinxupquote{double val,}}

\sphinxAtStartPar
\sphinxcode{\sphinxupquote{int offset)}}
\end{quote}
\end{quote}

\sphinxAtStartPar
\sphinxstylestrong{Arguments}
\begin{quote}

\sphinxAtStartPar
\sphinxcode{\sphinxupquote{dim}}: dimension of diagonal matrix.

\sphinxAtStartPar
\sphinxcode{\sphinxupquote{val}}: value to fill diagonal elements.

\sphinxAtStartPar
\sphinxcode{\sphinxupquote{offset}}: shift distance against diagonal line.
\end{quote}

\sphinxAtStartPar
\sphinxstylestrong{Return}
\begin{quote}

\sphinxAtStartPar
new diagonal matrix object.
\end{quote}
\end{quote}

\subsubsection{Model::AddDiagMat()}
\label{\detokenize{cppapi/model:id18}}\begin{quote}

\sphinxAtStartPar
Add a diagonal matrix to a model.

\sphinxAtStartPar
\sphinxstylestrong{Synopsis}
\begin{quote}

\sphinxAtStartPar
\sphinxcode{\sphinxupquote{SymMatrix AddDiagMat(}}
\begin{quote}

\sphinxAtStartPar
\sphinxcode{\sphinxupquote{int dim,}}

\sphinxAtStartPar
\sphinxcode{\sphinxupquote{double *pVals,}}

\sphinxAtStartPar
\sphinxcode{\sphinxupquote{int len,}}

\sphinxAtStartPar
\sphinxcode{\sphinxupquote{int offset)}}
\end{quote}
\end{quote}

\sphinxAtStartPar
\sphinxstylestrong{Arguments}
\begin{quote}

\sphinxAtStartPar
\sphinxcode{\sphinxupquote{dim}}: dimension of diagonal matrix.

\sphinxAtStartPar
\sphinxcode{\sphinxupquote{pVals}}: array of values of diagonal elements.

\sphinxAtStartPar
\sphinxcode{\sphinxupquote{len}}: length of value array.

\sphinxAtStartPar
\sphinxcode{\sphinxupquote{offset}}: shift distance against diagonal line.
\end{quote}

\sphinxAtStartPar
\sphinxstylestrong{Return}
\begin{quote}

\sphinxAtStartPar
new diagonal matrix object.
\end{quote}
\end{quote}

\subsubsection{Model::AddExpCone()}
\label{\detokenize{cppapi/model:model-addexpcone}}\begin{quote}

\sphinxAtStartPar
Add an exponential cone constraint to a model.

\sphinxAtStartPar
\sphinxstylestrong{Synopsis}
\begin{quote}

\sphinxAtStartPar
\sphinxcode{\sphinxupquote{ExpCone AddExpCone(}}
\begin{quote}

\sphinxAtStartPar
\sphinxcode{\sphinxupquote{int type,}}

\sphinxAtStartPar
\sphinxcode{\sphinxupquote{char *pvtype,}}

\sphinxAtStartPar
\sphinxcode{\sphinxupquote{const char *szPrefix)}}
\end{quote}
\end{quote}

\sphinxAtStartPar
\sphinxstylestrong{Arguments}
\begin{quote}

\sphinxAtStartPar
\sphinxcode{\sphinxupquote{type}}: type of the exponential cone constraint.

\sphinxAtStartPar
\sphinxcode{\sphinxupquote{pvtype}}: types of variables in the exponential cone.

\sphinxAtStartPar
\sphinxcode{\sphinxupquote{szPrefix}}: name prefix of variables in the exponential cone.
\end{quote}

\sphinxAtStartPar
\sphinxstylestrong{Return}
\begin{quote}

\sphinxAtStartPar
object of new exponential cone constraint.
\end{quote}
\end{quote}

\subsubsection{Model::AddExpCone()}
\label{\detokenize{cppapi/model:id19}}\begin{quote}

\sphinxAtStartPar
Add an exponential cone constraint to model.

\sphinxAtStartPar
\sphinxstylestrong{Synopsis}
\begin{quote}

\sphinxAtStartPar
\sphinxcode{\sphinxupquote{ExpCone AddExpCone(const ExpConeBuilder \&builder)}}
\end{quote}

\sphinxAtStartPar
\sphinxstylestrong{Arguments}
\begin{quote}

\sphinxAtStartPar
\sphinxcode{\sphinxupquote{builder}}: builder for new exponential cone constraint.
\end{quote}

\sphinxAtStartPar
\sphinxstylestrong{Return}
\begin{quote}

\sphinxAtStartPar
new exponential cone constraint object.
\end{quote}
\end{quote}

\subsubsection{Model::AddExpCone()}
\label{\detokenize{cppapi/model:id20}}\begin{quote}

\sphinxAtStartPar
Add an exponential cone constraint to model.

\sphinxAtStartPar
\sphinxstylestrong{Synopsis}
\begin{quote}

\sphinxAtStartPar
\sphinxcode{\sphinxupquote{ExpCone AddExpCone(const VarArray \&vars, int type)}}
\end{quote}

\sphinxAtStartPar
\sphinxstylestrong{Arguments}
\begin{quote}

\sphinxAtStartPar
\sphinxcode{\sphinxupquote{vars}}: variables that participate in the exponential cone constraint.

\sphinxAtStartPar
\sphinxcode{\sphinxupquote{type}}: type of the exponential cone constraint.
\end{quote}

\sphinxAtStartPar
\sphinxstylestrong{Return}
\begin{quote}

\sphinxAtStartPar
object of new exponential cone constraint.
\end{quote}
\end{quote}

\subsubsection{Model::AddExpCone()}
\label{\detokenize{cppapi/model:id21}}\begin{quote}

\sphinxAtStartPar
Add an exponential cone constraint to model.

\sphinxAtStartPar
\sphinxstylestrong{Synopsis}
\begin{quote}

\sphinxAtStartPar
\sphinxcode{\sphinxupquote{ExpCone AddExpCone(const MVar\textless{}1\textgreater{} \&vars, int type)}}
\end{quote}

\sphinxAtStartPar
\sphinxstylestrong{Arguments}
\begin{quote}

\sphinxAtStartPar
\sphinxcode{\sphinxupquote{vars}}: 1\sphinxhyphen{}dimensional variables in the exponential cone constraint.

\sphinxAtStartPar
\sphinxcode{\sphinxupquote{type}}: type of the exponential cone constraint.
\end{quote}

\sphinxAtStartPar
\sphinxstylestrong{Return}
\begin{quote}

\sphinxAtStartPar
object of new exponential cone constraint.
\end{quote}
\end{quote}

\subsubsection{Model::AddExpCones()}
\label{\detokenize{cppapi/model:model-addexpcones}}\begin{quote}

\sphinxAtStartPar
Add a list of expoential cone constraints to a model.

\sphinxAtStartPar
\sphinxstylestrong{Synopsis}
\begin{quote}

\sphinxAtStartPar
\sphinxcode{\sphinxupquote{ExpConeArray AddExpCones(const MVar\textless{}2\textgreater{} \&vars, int type)}}
\end{quote}

\sphinxAtStartPar
\sphinxstylestrong{Arguments}
\begin{quote}

\sphinxAtStartPar
\sphinxcode{\sphinxupquote{vars}}: 2\sphinxhyphen{}dimensional array of variables.

\sphinxAtStartPar
\sphinxcode{\sphinxupquote{type}}: type of an expoential cone.
\end{quote}

\sphinxAtStartPar
\sphinxstylestrong{Return}
\begin{quote}

\sphinxAtStartPar
array of new expoential cone constraint objects.
\end{quote}
\end{quote}

\subsubsection{Model::AddEyeMat()}
\label{\detokenize{cppapi/model:model-addeyemat}}\begin{quote}

\sphinxAtStartPar
Add an identity matrix to a model.

\sphinxAtStartPar
\sphinxstylestrong{Synopsis}
\begin{quote}

\sphinxAtStartPar
\sphinxcode{\sphinxupquote{SymMatrix AddEyeMat(int dim)}}
\end{quote}

\sphinxAtStartPar
\sphinxstylestrong{Arguments}
\begin{quote}

\sphinxAtStartPar
\sphinxcode{\sphinxupquote{dim}}: dimension of identity matrix.
\end{quote}

\sphinxAtStartPar
\sphinxstylestrong{Return}
\begin{quote}

\sphinxAtStartPar
new identity matrix object.
\end{quote}
\end{quote}

\subsubsection{Model::AddGenConstrIndicator()}
\label{\detokenize{cppapi/model:model-addgenconstrindicator}}\begin{quote}

\sphinxAtStartPar
Add a general constraint of type indicator to model.

\sphinxAtStartPar
\sphinxstylestrong{Synopsis}
\begin{quote}

\sphinxAtStartPar
\sphinxcode{\sphinxupquote{GenConstr AddGenConstrIndicator(const GenConstrBuilder \&builder, const char *szName)}}
\end{quote}

\sphinxAtStartPar
\sphinxstylestrong{Arguments}
\begin{quote}

\sphinxAtStartPar
\sphinxcode{\sphinxupquote{builder}}: builder for the general constraint.

\sphinxAtStartPar
\sphinxcode{\sphinxupquote{szName}}: optional, name of new general constraint.
\end{quote}

\sphinxAtStartPar
\sphinxstylestrong{Return}
\begin{quote}

\sphinxAtStartPar
new general constraint object of type indicator.
\end{quote}
\end{quote}

\subsubsection{Model::AddGenConstrIndicator()}
\label{\detokenize{cppapi/model:id22}}\begin{quote}

\sphinxAtStartPar
Add a general constraint of type indicator to model.

\sphinxAtStartPar
\sphinxstylestrong{Synopsis}
\begin{quote}

\sphinxAtStartPar
\sphinxcode{\sphinxupquote{GenConstr AddGenConstrIndicator(}}
\begin{quote}

\sphinxAtStartPar
\sphinxcode{\sphinxupquote{Var binVar,}}

\sphinxAtStartPar
\sphinxcode{\sphinxupquote{int binVal,}}

\sphinxAtStartPar
\sphinxcode{\sphinxupquote{const ConstrBuilder \&builder,}}

\sphinxAtStartPar
\sphinxcode{\sphinxupquote{int type,}}

\sphinxAtStartPar
\sphinxcode{\sphinxupquote{const char *szName)}}
\end{quote}
\end{quote}

\sphinxAtStartPar
\sphinxstylestrong{Arguments}
\begin{quote}

\sphinxAtStartPar
\sphinxcode{\sphinxupquote{binVar}}: binary indicator variable.

\sphinxAtStartPar
\sphinxcode{\sphinxupquote{binVal}}: value for binary indicator variable that force a linear constraint to be satisfied(0 or 1).

\sphinxAtStartPar
\sphinxcode{\sphinxupquote{builder}}: builder for linear constraint.

\sphinxAtStartPar
\sphinxcode{\sphinxupquote{type}}: type of general constraint with default value of COPT\_INDICATOR\_IF.

\sphinxAtStartPar
\sphinxcode{\sphinxupquote{szName}}: optional, name of new general constraint.
\end{quote}

\sphinxAtStartPar
\sphinxstylestrong{Return}
\begin{quote}

\sphinxAtStartPar
new general constraint object of type indicator.
\end{quote}
\end{quote}

\subsubsection{Model::AddGenConstrIndicator()}
\label{\detokenize{cppapi/model:id23}}\begin{quote}

\sphinxAtStartPar
Add a general constraint of type indicator to model.

\sphinxAtStartPar
\sphinxstylestrong{Synopsis}
\begin{quote}

\sphinxAtStartPar
\sphinxcode{\sphinxupquote{GenConstr AddGenConstrIndicator(}}
\begin{quote}

\sphinxAtStartPar
\sphinxcode{\sphinxupquote{Var binVar,}}

\sphinxAtStartPar
\sphinxcode{\sphinxupquote{int binVal,}}

\sphinxAtStartPar
\sphinxcode{\sphinxupquote{const Expr \&expr,}}

\sphinxAtStartPar
\sphinxcode{\sphinxupquote{char sense,}}

\sphinxAtStartPar
\sphinxcode{\sphinxupquote{double rhs,}}

\sphinxAtStartPar
\sphinxcode{\sphinxupquote{int type,}}

\sphinxAtStartPar
\sphinxcode{\sphinxupquote{const char *szName)}}
\end{quote}
\end{quote}

\sphinxAtStartPar
\sphinxstylestrong{Arguments}
\begin{quote}

\sphinxAtStartPar
\sphinxcode{\sphinxupquote{binVar}}: binary indicator variable.

\sphinxAtStartPar
\sphinxcode{\sphinxupquote{binVal}}: value for binary indicator variable that force a linear constraint to be satisfied(0 or 1).

\sphinxAtStartPar
\sphinxcode{\sphinxupquote{expr}}: expression for new linear contraint.

\sphinxAtStartPar
\sphinxcode{\sphinxupquote{sense}}: sense for new linear constraint.

\sphinxAtStartPar
\sphinxcode{\sphinxupquote{rhs}}: right hand side value for new linear constraint.

\sphinxAtStartPar
\sphinxcode{\sphinxupquote{type}}: type of general constraint with default value of COPT\_INDICATOR\_IF.

\sphinxAtStartPar
\sphinxcode{\sphinxupquote{szName}}: optional, name of new general constraint.
\end{quote}

\sphinxAtStartPar
\sphinxstylestrong{Return}
\begin{quote}

\sphinxAtStartPar
new general constraint object of type indicator.
\end{quote}
\end{quote}

\subsubsection{Model::AddGenConstrIndicators()}
\label{\detokenize{cppapi/model:model-addgenconstrindicators}}\begin{quote}

\sphinxAtStartPar
Add general constraints of type indicator to model.

\sphinxAtStartPar
\sphinxstylestrong{Synopsis}
\begin{quote}

\sphinxAtStartPar
\sphinxcode{\sphinxupquote{GenConstrArray AddGenConstrIndicators(}}
\begin{quote}

\sphinxAtStartPar
\sphinxcode{\sphinxupquote{Var binVar,}}

\sphinxAtStartPar
\sphinxcode{\sphinxupquote{int binVal,}}

\sphinxAtStartPar
\sphinxcode{\sphinxupquote{const MConstrBuilder\textless{}N\textgreater{} \&builders,}}

\sphinxAtStartPar
\sphinxcode{\sphinxupquote{int type,}}

\sphinxAtStartPar
\sphinxcode{\sphinxupquote{const char *szPrefix)}}
\end{quote}
\end{quote}

\sphinxAtStartPar
\sphinxstylestrong{Arguments}
\begin{quote}

\sphinxAtStartPar
\sphinxcode{\sphinxupquote{binVar}}: binary indicator variable.

\sphinxAtStartPar
\sphinxcode{\sphinxupquote{binVal}}: value for binary indicator variable that force a linear constraint to be satisfied(0 or 1).

\sphinxAtStartPar
\sphinxcode{\sphinxupquote{builders}}: builder for linear constraints in matrix format.

\sphinxAtStartPar
\sphinxcode{\sphinxupquote{type}}: type of general constraint with default value of COPT\_INDICATOR\_IF.

\sphinxAtStartPar
\sphinxcode{\sphinxupquote{szPrefix}}: optional, name prefix of new general constraints.
\end{quote}

\sphinxAtStartPar
\sphinxstylestrong{Return}
\begin{quote}

\sphinxAtStartPar
new general constraint objects of type indicator.
\end{quote}
\end{quote}

\subsubsection{Model::AddGenConstrIndicators()}
\label{\detokenize{cppapi/model:id24}}\begin{quote}

\sphinxAtStartPar
Add a list of general constraints to model.

\sphinxAtStartPar
\sphinxstylestrong{Synopsis}
\begin{quote}

\sphinxAtStartPar
\sphinxcode{\sphinxupquote{GenConstrArray AddGenConstrIndicators(const GenConstrBuilderArray \&builders, const char *szPrefix)}}
\end{quote}

\sphinxAtStartPar
\sphinxstylestrong{Arguments}
\begin{quote}

\sphinxAtStartPar
\sphinxcode{\sphinxupquote{builders}}: builders for new general constraints.

\sphinxAtStartPar
\sphinxcode{\sphinxupquote{szPrefix}}: name prefix for new general constraints.
\end{quote}

\sphinxAtStartPar
\sphinxstylestrong{Return}
\begin{quote}

\sphinxAtStartPar
array of new general constraint objects.
\end{quote}
\end{quote}

\subsubsection{Model::AddGenConstrIndicators()}
\label{\detokenize{cppapi/model:id25}}\begin{quote}

\sphinxAtStartPar
Add a list of general constraints to model.

\sphinxAtStartPar
\sphinxstylestrong{Synopsis}
\begin{quote}

\sphinxAtStartPar
\sphinxcode{\sphinxupquote{GenConstrArray AddGenConstrIndicators(}}
\begin{quote}

\sphinxAtStartPar
\sphinxcode{\sphinxupquote{const GenConstrBuilderArray \&builders,}}

\sphinxAtStartPar
\sphinxcode{\sphinxupquote{const char *szNames,}}

\sphinxAtStartPar
\sphinxcode{\sphinxupquote{size\_t len)}}
\end{quote}
\end{quote}

\sphinxAtStartPar
\sphinxstylestrong{Arguments}
\begin{quote}

\sphinxAtStartPar
\sphinxcode{\sphinxupquote{builders}}: builders for new general constraints.

\sphinxAtStartPar
\sphinxcode{\sphinxupquote{szNames}}: name buffer of new general constraints.

\sphinxAtStartPar
\sphinxcode{\sphinxupquote{len}}: length of the name buffer.
\end{quote}

\sphinxAtStartPar
\sphinxstylestrong{Return}
\begin{quote}

\sphinxAtStartPar
array of new general constraint objects.
\end{quote}
\end{quote}

\subsubsection{Model::AddLazyConstr()}
\label{\detokenize{cppapi/model:model-addlazyconstr}}\begin{quote}

\sphinxAtStartPar
Add a lazy constraint to model.

\sphinxAtStartPar
\sphinxstylestrong{Synopsis}
\begin{quote}

\sphinxAtStartPar
\sphinxcode{\sphinxupquote{void AddLazyConstr(}}
\begin{quote}

\sphinxAtStartPar
\sphinxcode{\sphinxupquote{const Expr \&lhs,}}

\sphinxAtStartPar
\sphinxcode{\sphinxupquote{char sense,}}

\sphinxAtStartPar
\sphinxcode{\sphinxupquote{double rhs,}}

\sphinxAtStartPar
\sphinxcode{\sphinxupquote{const char *szName)}}
\end{quote}
\end{quote}

\sphinxAtStartPar
\sphinxstylestrong{Arguments}
\begin{quote}

\sphinxAtStartPar
\sphinxcode{\sphinxupquote{lhs}}: expression for lazy contraint.

\sphinxAtStartPar
\sphinxcode{\sphinxupquote{sense}}: sense for lazy constraint.

\sphinxAtStartPar
\sphinxcode{\sphinxupquote{rhs}}: right hand side value for lazy constraint.

\sphinxAtStartPar
\sphinxcode{\sphinxupquote{szName}}: optional, name of lazy constraint.
\end{quote}
\end{quote}

\subsubsection{Model::AddLazyConstr()}
\label{\detokenize{cppapi/model:id26}}\begin{quote}

\sphinxAtStartPar
Add a lazy constraint to model.

\sphinxAtStartPar
\sphinxstylestrong{Synopsis}
\begin{quote}

\sphinxAtStartPar
\sphinxcode{\sphinxupquote{void AddLazyConstr(}}
\begin{quote}

\sphinxAtStartPar
\sphinxcode{\sphinxupquote{const Expr \&lhs,}}

\sphinxAtStartPar
\sphinxcode{\sphinxupquote{char sense,}}

\sphinxAtStartPar
\sphinxcode{\sphinxupquote{const Expr \&rhs,}}

\sphinxAtStartPar
\sphinxcode{\sphinxupquote{const char *szName)}}
\end{quote}
\end{quote}

\sphinxAtStartPar
\sphinxstylestrong{Arguments}
\begin{quote}

\sphinxAtStartPar
\sphinxcode{\sphinxupquote{lhs}}: left hand side expression for lazy contraint.

\sphinxAtStartPar
\sphinxcode{\sphinxupquote{sense}}: sense for lazy constraint.

\sphinxAtStartPar
\sphinxcode{\sphinxupquote{rhs}}: right hand side expression for lazy contraint.

\sphinxAtStartPar
\sphinxcode{\sphinxupquote{szName}}: optional, name of lazy constraint.
\end{quote}
\end{quote}

\subsubsection{Model::AddLazyConstr()}
\label{\detokenize{cppapi/model:id27}}\begin{quote}

\sphinxAtStartPar
Add a lazy constraint to model.

\sphinxAtStartPar
\sphinxstylestrong{Synopsis}
\begin{quote}

\sphinxAtStartPar
\sphinxcode{\sphinxupquote{void AddLazyConstr(const ConstrBuilder \&builder, const char *szName)}}
\end{quote}

\sphinxAtStartPar
\sphinxstylestrong{Arguments}
\begin{quote}

\sphinxAtStartPar
\sphinxcode{\sphinxupquote{builder}}: builder for lazy contraint.

\sphinxAtStartPar
\sphinxcode{\sphinxupquote{szName}}: optional, name of lazy constraint.
\end{quote}
\end{quote}

\subsubsection{Model::AddLazyConstrs()}
\label{\detokenize{cppapi/model:model-addlazyconstrs}}\begin{quote}

\sphinxAtStartPar
Add lazy constraints to model.

\sphinxAtStartPar
\sphinxstylestrong{Synopsis}
\begin{quote}

\sphinxAtStartPar
\sphinxcode{\sphinxupquote{void AddLazyConstrs(const ConstrBuilderArray \&builders, const char *szPrefix)}}
\end{quote}

\sphinxAtStartPar
\sphinxstylestrong{Arguments}
\begin{quote}

\sphinxAtStartPar
\sphinxcode{\sphinxupquote{builders}}: array of builders for lazy contraints.

\sphinxAtStartPar
\sphinxcode{\sphinxupquote{szPrefix}}: name prefix of new lazy constraints.
\end{quote}
\end{quote}

\subsubsection{Model::AddLmiConstr()}
\label{\detokenize{cppapi/model:model-addlmiconstr}}\begin{quote}

\sphinxAtStartPar
Add an LMI constraint to model.

\sphinxAtStartPar
\sphinxstylestrong{Synopsis}
\begin{quote}

\sphinxAtStartPar
\sphinxcode{\sphinxupquote{LmiConstraint AddLmiConstr(const LmiExpr \&expr, const char *szName)}}
\end{quote}

\sphinxAtStartPar
\sphinxstylestrong{Arguments}
\begin{quote}

\sphinxAtStartPar
\sphinxcode{\sphinxupquote{expr}}: LMI expression for new LMI contraint.

\sphinxAtStartPar
\sphinxcode{\sphinxupquote{szName}}: optional, name of new LMI constraint.
\end{quote}

\sphinxAtStartPar
\sphinxstylestrong{Return}
\begin{quote}

\sphinxAtStartPar
new LMI constraint object.
\end{quote}
\end{quote}

\subsubsection{Model::AddMConstr()}
\label{\detokenize{cppapi/model:model-addmconstr}}\begin{quote}

\sphinxAtStartPar
Add a MConstr object in N\sphinxhyphen{}dimensions to model.

\sphinxAtStartPar
\sphinxstylestrong{Synopsis}
\begin{quote}

\sphinxAtStartPar
\sphinxcode{\sphinxupquote{template \textless{}int N\textgreater{} MConstr\textless{}N\textgreater{} AddMConstr(const MConstrBuilder\textless{}N\textgreater{} \&builder, const char *szPrefix)}}
\end{quote}

\sphinxAtStartPar
\sphinxstylestrong{Arguments}
\begin{quote}

\sphinxAtStartPar
\sphinxcode{\sphinxupquote{builder}}: builder for MConstr object.

\sphinxAtStartPar
\sphinxcode{\sphinxupquote{szPrefix}}: name prefix for constraints in MConstr object.
\end{quote}

\sphinxAtStartPar
\sphinxstylestrong{Return}
\begin{quote}

\sphinxAtStartPar
new MConstr object.
\end{quote}
\end{quote}

\subsubsection{Model::AddMConstr()}
\label{\detokenize{cppapi/model:id28}}\begin{quote}

\sphinxAtStartPar
Add a N\sphinxhyphen{}dimensional MConstr object to model.

\sphinxAtStartPar
\sphinxstylestrong{Synopsis}
\begin{quote}

\sphinxAtStartPar
\sphinxcode{\sphinxupquote{template \textless{}int N\textgreater{} MConstr\textless{}N\textgreater{} AddMConstr(}}
\begin{quote}

\sphinxAtStartPar
\sphinxcode{\sphinxupquote{const MLinExpr\textless{}N\textgreater{} \&exprs,}}

\sphinxAtStartPar
\sphinxcode{\sphinxupquote{char sense,}}

\sphinxAtStartPar
\sphinxcode{\sphinxupquote{double rhs,}}

\sphinxAtStartPar
\sphinxcode{\sphinxupquote{const char *szPrefix)}}
\end{quote}
\end{quote}

\sphinxAtStartPar
\sphinxstylestrong{Arguments}
\begin{quote}

\sphinxAtStartPar
\sphinxcode{\sphinxupquote{exprs}}: N\sphinxhyphen{}dimensional MLinExpr object.

\sphinxAtStartPar
\sphinxcode{\sphinxupquote{sense}}: sense for new linear constraints.

\sphinxAtStartPar
\sphinxcode{\sphinxupquote{rhs}}: double value at right side of the new linear constraints.

\sphinxAtStartPar
\sphinxcode{\sphinxupquote{szPrefix}}: name prefix for constraints in MConstr object.
\end{quote}

\sphinxAtStartPar
\sphinxstylestrong{Return}
\begin{quote}

\sphinxAtStartPar
new MConstr object.
\end{quote}
\end{quote}

\subsubsection{Model::AddMPsdConstr()}
\label{\detokenize{cppapi/model:model-addmpsdconstr}}\begin{quote}

\sphinxAtStartPar
Add a N\sphinxhyphen{}dimensional MPsdConstr object to model.

\sphinxAtStartPar
\sphinxstylestrong{Synopsis}
\begin{quote}

\sphinxAtStartPar
\sphinxcode{\sphinxupquote{template \textless{}int N\textgreater{} MPsdConstr\textless{}N\textgreater{} AddMPsdConstr(const MPsdConstrBuilder\textless{}N\textgreater{} \&builder, const char *szPrefix)}}
\end{quote}

\sphinxAtStartPar
\sphinxstylestrong{Arguments}
\begin{quote}

\sphinxAtStartPar
\sphinxcode{\sphinxupquote{builder}}: builder for MPsdConstr object.

\sphinxAtStartPar
\sphinxcode{\sphinxupquote{szPrefix}}: name prefix of PSD constraints in MPsdConstr object.
\end{quote}

\sphinxAtStartPar
\sphinxstylestrong{Return}
\begin{quote}

\sphinxAtStartPar
new MPsdConstr object.
\end{quote}
\end{quote}

\subsubsection{Model::AddMPsdConstr()}
\label{\detokenize{cppapi/model:id29}}\begin{quote}

\sphinxAtStartPar
Add a N\sphinxhyphen{}dimensional MPsdConstr object to model.

\sphinxAtStartPar
\sphinxstylestrong{Synopsis}
\begin{quote}

\sphinxAtStartPar
\sphinxcode{\sphinxupquote{template \textless{}int N\textgreater{} MPsdConstr\textless{}N\textgreater{} AddMPsdConstr(}}
\begin{quote}

\sphinxAtStartPar
\sphinxcode{\sphinxupquote{const MPsdExpr\textless{}N\textgreater{} \&exprs,}}

\sphinxAtStartPar
\sphinxcode{\sphinxupquote{char sense,}}

\sphinxAtStartPar
\sphinxcode{\sphinxupquote{double rhs,}}

\sphinxAtStartPar
\sphinxcode{\sphinxupquote{const char *szPrefix)}}
\end{quote}
\end{quote}

\sphinxAtStartPar
\sphinxstylestrong{Arguments}
\begin{quote}

\sphinxAtStartPar
\sphinxcode{\sphinxupquote{exprs}}: N\sphinxhyphen{}dimensional MPsdExpr object.

\sphinxAtStartPar
\sphinxcode{\sphinxupquote{sense}}: sense for new PSD constraints.

\sphinxAtStartPar
\sphinxcode{\sphinxupquote{rhs}}: double value at right side of the new PSD constraints.

\sphinxAtStartPar
\sphinxcode{\sphinxupquote{szPrefix}}: name prefix of PSD constraints in MPsdConstr object.
\end{quote}

\sphinxAtStartPar
\sphinxstylestrong{Return}
\begin{quote}

\sphinxAtStartPar
new MPsdConstr object.
\end{quote}
\end{quote}

\subsubsection{Model::AddMQConstr()}
\label{\detokenize{cppapi/model:model-addmqconstr}}\begin{quote}

\sphinxAtStartPar
Add a N\sphinxhyphen{}dimensional MQConstr object to model.

\sphinxAtStartPar
\sphinxstylestrong{Synopsis}
\begin{quote}

\sphinxAtStartPar
\sphinxcode{\sphinxupquote{template \textless{}int N\textgreater{} MQConstr\textless{}N\textgreater{} AddMQConstr(const MQConstrBuilder\textless{}N\textgreater{} \&builder, const char *szPrefix)}}
\end{quote}

\sphinxAtStartPar
\sphinxstylestrong{Arguments}
\begin{quote}

\sphinxAtStartPar
\sphinxcode{\sphinxupquote{builder}}: builder for MQConstr object.

\sphinxAtStartPar
\sphinxcode{\sphinxupquote{szPrefix}}: name prefix of quadratic constraints in MQConstr object.
\end{quote}

\sphinxAtStartPar
\sphinxstylestrong{Return}
\begin{quote}

\sphinxAtStartPar
new MQConstr object.
\end{quote}
\end{quote}

\subsubsection{Model::AddMQConstr()}
\label{\detokenize{cppapi/model:id30}}\begin{quote}

\sphinxAtStartPar
Add a N\sphinxhyphen{}dimensional MQConstr object to model.

\sphinxAtStartPar
\sphinxstylestrong{Synopsis}
\begin{quote}

\sphinxAtStartPar
\sphinxcode{\sphinxupquote{template \textless{}int N\textgreater{} MQConstr\textless{}N\textgreater{} AddMQConstr(}}
\begin{quote}

\sphinxAtStartPar
\sphinxcode{\sphinxupquote{const MQuadExpr\textless{}N\textgreater{} \&exprs,}}

\sphinxAtStartPar
\sphinxcode{\sphinxupquote{char sense,}}

\sphinxAtStartPar
\sphinxcode{\sphinxupquote{double rhs,}}

\sphinxAtStartPar
\sphinxcode{\sphinxupquote{const char *szPrefix)}}
\end{quote}
\end{quote}

\sphinxAtStartPar
\sphinxstylestrong{Arguments}
\begin{quote}

\sphinxAtStartPar
\sphinxcode{\sphinxupquote{exprs}}: N\sphinxhyphen{}dimensional MQuadExpr object.

\sphinxAtStartPar
\sphinxcode{\sphinxupquote{sense}}: sense for new quadratic constraints.

\sphinxAtStartPar
\sphinxcode{\sphinxupquote{rhs}}: double value at right side of the new quadratic constraints.

\sphinxAtStartPar
\sphinxcode{\sphinxupquote{szPrefix}}: name prefix of quadratic constraints in MQConstr object.
\end{quote}

\sphinxAtStartPar
\sphinxstylestrong{Return}
\begin{quote}

\sphinxAtStartPar
new MQConstr object.
\end{quote}
\end{quote}

\subsubsection{Model::AddMVar()}
\label{\detokenize{cppapi/model:model-addmvar}}\begin{quote}

\sphinxAtStartPar
Add a MVar object in N\sphinxhyphen{}dimensions to model.

\sphinxAtStartPar
\sphinxstylestrong{Synopsis}
\begin{quote}

\sphinxAtStartPar
\sphinxcode{\sphinxupquote{template \textless{}int N\textgreater{} MVar\textless{}N\textgreater{} AddMVar(}}
\begin{quote}

\sphinxAtStartPar
\sphinxcode{\sphinxupquote{const Shape\textless{}N\textgreater{} \&shp,}}

\sphinxAtStartPar
\sphinxcode{\sphinxupquote{char vtype,}}

\sphinxAtStartPar
\sphinxcode{\sphinxupquote{const char *szPrefix)}}
\end{quote}
\end{quote}

\sphinxAtStartPar
\sphinxstylestrong{Arguments}
\begin{quote}

\sphinxAtStartPar
\sphinxcode{\sphinxupquote{shp}}: shape of MVar object.

\sphinxAtStartPar
\sphinxcode{\sphinxupquote{vtype}}: type of variables in MVar object.

\sphinxAtStartPar
\sphinxcode{\sphinxupquote{szPrefix}}: name prefix of variables in MVar object.
\end{quote}

\sphinxAtStartPar
\sphinxstylestrong{Return}
\begin{quote}

\sphinxAtStartPar
new MVar object.
\end{quote}
\end{quote}

\subsubsection{Model::AddMVar()}
\label{\detokenize{cppapi/model:id31}}\begin{quote}

\sphinxAtStartPar
Add a MVar object in N\sphinxhyphen{}dimensions to model.

\sphinxAtStartPar
\sphinxstylestrong{Synopsis}
\begin{quote}

\sphinxAtStartPar
\sphinxcode{\sphinxupquote{MVar\textless{}N\textgreater{} AddMVar(}}
\begin{quote}

\sphinxAtStartPar
\sphinxcode{\sphinxupquote{const Shape\textless{}N\textgreater{} \&shp,}}

\sphinxAtStartPar
\sphinxcode{\sphinxupquote{double lb,}}

\sphinxAtStartPar
\sphinxcode{\sphinxupquote{double ub,}}

\sphinxAtStartPar
\sphinxcode{\sphinxupquote{double obj,}}

\sphinxAtStartPar
\sphinxcode{\sphinxupquote{char vtype,}}

\sphinxAtStartPar
\sphinxcode{\sphinxupquote{const char *szPrefix)}}
\end{quote}
\end{quote}

\sphinxAtStartPar
\sphinxstylestrong{Arguments}
\begin{quote}

\sphinxAtStartPar
\sphinxcode{\sphinxupquote{shp}}: shape of MVar object.

\sphinxAtStartPar
\sphinxcode{\sphinxupquote{lb}}: lower bound for variables in MVar object.

\sphinxAtStartPar
\sphinxcode{\sphinxupquote{ub}}: upper bound for variables in MVar object.

\sphinxAtStartPar
\sphinxcode{\sphinxupquote{obj}}: objective coefficient for variables in MVar object.

\sphinxAtStartPar
\sphinxcode{\sphinxupquote{vtype}}: type of variables in MVar object.

\sphinxAtStartPar
\sphinxcode{\sphinxupquote{szPrefix}}: name prefix of variables in MVar object.
\end{quote}

\sphinxAtStartPar
\sphinxstylestrong{Return}
\begin{quote}

\sphinxAtStartPar
new MVar object.
\end{quote}
\end{quote}

\subsubsection{Model::AddMVar()}
\label{\detokenize{cppapi/model:id32}}\begin{quote}

\sphinxAtStartPar
Add a MVar object in N\sphinxhyphen{}dimensions to model.

\sphinxAtStartPar
\sphinxstylestrong{Synopsis}
\begin{quote}

\sphinxAtStartPar
\sphinxcode{\sphinxupquote{MVar\textless{}N\textgreater{} AddMVar(}}
\begin{quote}

\sphinxAtStartPar
\sphinxcode{\sphinxupquote{const Shape\textless{}N\textgreater{} \&shp,}}

\sphinxAtStartPar
\sphinxcode{\sphinxupquote{double *plb,}}

\sphinxAtStartPar
\sphinxcode{\sphinxupquote{double *pub,}}

\sphinxAtStartPar
\sphinxcode{\sphinxupquote{double *pobj,}}

\sphinxAtStartPar
\sphinxcode{\sphinxupquote{char *pvtype,}}

\sphinxAtStartPar
\sphinxcode{\sphinxupquote{const char *szPrefix)}}
\end{quote}
\end{quote}

\sphinxAtStartPar
\sphinxstylestrong{Arguments}
\begin{quote}

\sphinxAtStartPar
\sphinxcode{\sphinxupquote{shp}}: shape of MVar object.

\sphinxAtStartPar
\sphinxcode{\sphinxupquote{plb}}: lower bounds for variables in MVar object. If NULL, lower bounds are 0.0.

\sphinxAtStartPar
\sphinxcode{\sphinxupquote{pub}}: upper bounds for variables in Mar object. If NULL, upper bounds are infinity or 1 for binary variables.

\sphinxAtStartPar
\sphinxcode{\sphinxupquote{pobj}}: objective coefficient for variables in MVar object. If NULL, objective coefficients are 0.0.

\sphinxAtStartPar
\sphinxcode{\sphinxupquote{pvtype}}: type of variables in MVar object. If NULL, variable types are continuous.

\sphinxAtStartPar
\sphinxcode{\sphinxupquote{szPrefix}}: name prefix of variables in MVar object.
\end{quote}

\sphinxAtStartPar
\sphinxstylestrong{Return}
\begin{quote}

\sphinxAtStartPar
new MVar object.
\end{quote}
\end{quote}

\subsubsection{Model::AddNlConstr()}
\label{\detokenize{cppapi/model:model-addnlconstr}}\begin{quote}

\sphinxAtStartPar
Add a nonlinear constraint to model.

\sphinxAtStartPar
\sphinxstylestrong{Synopsis}
\begin{quote}

\sphinxAtStartPar
\sphinxcode{\sphinxupquote{NlConstraint AddNlConstr(}}
\begin{quote}

\sphinxAtStartPar
\sphinxcode{\sphinxupquote{const NlExpr \&expr,}}

\sphinxAtStartPar
\sphinxcode{\sphinxupquote{char sense,}}

\sphinxAtStartPar
\sphinxcode{\sphinxupquote{double rhs,}}

\sphinxAtStartPar
\sphinxcode{\sphinxupquote{const char *szName)}}
\end{quote}
\end{quote}

\sphinxAtStartPar
\sphinxstylestrong{Arguments}
\begin{quote}

\sphinxAtStartPar
\sphinxcode{\sphinxupquote{expr}}: non\sphinxhyphen{}expression for the new contraint.

\sphinxAtStartPar
\sphinxcode{\sphinxupquote{sense}}: sense for new nonlinear constraint, other than range sense.

\sphinxAtStartPar
\sphinxcode{\sphinxupquote{rhs}}: right hand side value for the new constraint.

\sphinxAtStartPar
\sphinxcode{\sphinxupquote{szName}}: optional, name of new nonlinear constraint.
\end{quote}

\sphinxAtStartPar
\sphinxstylestrong{Return}
\begin{quote}

\sphinxAtStartPar
new nonlinear constraint object.
\end{quote}
\end{quote}

\subsubsection{Model::AddNlConstr()}
\label{\detokenize{cppapi/model:id33}}\begin{quote}

\sphinxAtStartPar
Add a nonlinear constraint to model.

\sphinxAtStartPar
\sphinxstylestrong{Synopsis}
\begin{quote}

\sphinxAtStartPar
\sphinxcode{\sphinxupquote{NlConstraint AddNlConstr(}}
\begin{quote}

\sphinxAtStartPar
\sphinxcode{\sphinxupquote{const NlExpr \&lhs,}}

\sphinxAtStartPar
\sphinxcode{\sphinxupquote{char sense,}}

\sphinxAtStartPar
\sphinxcode{\sphinxupquote{const NlExpr \&rhs,}}

\sphinxAtStartPar
\sphinxcode{\sphinxupquote{const char *szName)}}
\end{quote}
\end{quote}

\sphinxAtStartPar
\sphinxstylestrong{Arguments}
\begin{quote}

\sphinxAtStartPar
\sphinxcode{\sphinxupquote{lhs}}: left hand side nonlinear expression for the new constraint.

\sphinxAtStartPar
\sphinxcode{\sphinxupquote{sense}}: sense for new nonlinear constraint, other than range sense.

\sphinxAtStartPar
\sphinxcode{\sphinxupquote{rhs}}: right hand side nonlinear expression for the new constraint.

\sphinxAtStartPar
\sphinxcode{\sphinxupquote{szName}}: optional, name of new nonlinear constraint.
\end{quote}

\sphinxAtStartPar
\sphinxstylestrong{Return}
\begin{quote}

\sphinxAtStartPar
new nonlinear constraint object.
\end{quote}
\end{quote}

\subsubsection{Model::AddNlConstr()}
\label{\detokenize{cppapi/model:id34}}\begin{quote}

\sphinxAtStartPar
Add a nonlinear constraint to model.

\sphinxAtStartPar
\sphinxstylestrong{Synopsis}
\begin{quote}

\sphinxAtStartPar
\sphinxcode{\sphinxupquote{NlConstraint AddNlConstr(}}
\begin{quote}

\sphinxAtStartPar
\sphinxcode{\sphinxupquote{const NlExpr \&expr,}}

\sphinxAtStartPar
\sphinxcode{\sphinxupquote{double lb,}}

\sphinxAtStartPar
\sphinxcode{\sphinxupquote{double ub,}}

\sphinxAtStartPar
\sphinxcode{\sphinxupquote{const char *szName)}}
\end{quote}
\end{quote}

\sphinxAtStartPar
\sphinxstylestrong{Arguments}
\begin{quote}

\sphinxAtStartPar
\sphinxcode{\sphinxupquote{expr}}: nonlinear expression for the new constraint.

\sphinxAtStartPar
\sphinxcode{\sphinxupquote{lb}}: lower bound for the new nonlinear constraint.

\sphinxAtStartPar
\sphinxcode{\sphinxupquote{ub}}: upper bound for the new nonlinear constraint

\sphinxAtStartPar
\sphinxcode{\sphinxupquote{szName}}: optional, name of new constraint.
\end{quote}

\sphinxAtStartPar
\sphinxstylestrong{Return}
\begin{quote}

\sphinxAtStartPar
new nonlinear constraint object.
\end{quote}
\end{quote}

\subsubsection{Model::AddNlConstr()}
\label{\detokenize{cppapi/model:id35}}\begin{quote}

\sphinxAtStartPar
Add a nonlinear constraint to a model.

\sphinxAtStartPar
\sphinxstylestrong{Synopsis}
\begin{quote}

\sphinxAtStartPar
\sphinxcode{\sphinxupquote{NlConstraint AddNlConstr(const NlConstrBuilder \&builder, const char *szName)}}
\end{quote}

\sphinxAtStartPar
\sphinxstylestrong{Arguments}
\begin{quote}

\sphinxAtStartPar
\sphinxcode{\sphinxupquote{builder}}: builder for the new nonlinear constraint.

\sphinxAtStartPar
\sphinxcode{\sphinxupquote{szName}}: optional, name of new nonlinear constraint.
\end{quote}

\sphinxAtStartPar
\sphinxstylestrong{Return}
\begin{quote}

\sphinxAtStartPar
new nonlinear constraint object.
\end{quote}
\end{quote}

\subsubsection{Model::AddNlConstrs()}
\label{\detokenize{cppapi/model:model-addnlconstrs}}\begin{quote}

\sphinxAtStartPar
Add nonlinear constraints to a model.

\sphinxAtStartPar
\sphinxstylestrong{Synopsis}
\begin{quote}

\sphinxAtStartPar
\sphinxcode{\sphinxupquote{NlConstrArray AddNlConstrs(const NlConstrBuilderArray \&builders, const char *szPrefix)}}
\end{quote}

\sphinxAtStartPar
\sphinxstylestrong{Arguments}
\begin{quote}

\sphinxAtStartPar
\sphinxcode{\sphinxupquote{builders}}: builders for new nonlinear constraints.

\sphinxAtStartPar
\sphinxcode{\sphinxupquote{szPrefix}}: name prefix for new constraints.
\end{quote}

\sphinxAtStartPar
\sphinxstylestrong{Return}
\begin{quote}

\sphinxAtStartPar
array of new nonlinear constraint objects.
\end{quote}
\end{quote}

\subsubsection{Model::AddNlConstrs()}
\label{\detokenize{cppapi/model:id36}}\begin{quote}

\sphinxAtStartPar
Add nonlinear constraints to model.

\sphinxAtStartPar
\sphinxstylestrong{Synopsis}
\begin{quote}

\sphinxAtStartPar
\sphinxcode{\sphinxupquote{NlConstrArray AddNlConstrs(}}
\begin{quote}

\sphinxAtStartPar
\sphinxcode{\sphinxupquote{const NlConstrBuilderArray \&builders,}}

\sphinxAtStartPar
\sphinxcode{\sphinxupquote{const char *szNames,}}

\sphinxAtStartPar
\sphinxcode{\sphinxupquote{size\_t len)}}
\end{quote}
\end{quote}

\sphinxAtStartPar
\sphinxstylestrong{Arguments}
\begin{quote}

\sphinxAtStartPar
\sphinxcode{\sphinxupquote{builders}}: builders for new nonlinear constraints.

\sphinxAtStartPar
\sphinxcode{\sphinxupquote{szNames}}: name buffer of new constraints.

\sphinxAtStartPar
\sphinxcode{\sphinxupquote{len}}: length of the name buffer.
\end{quote}

\sphinxAtStartPar
\sphinxstylestrong{Return}
\begin{quote}

\sphinxAtStartPar
array of new nonlinear constraint objects.
\end{quote}
\end{quote}

\subsubsection{Model::AddOnesMat()}
\label{\detokenize{cppapi/model:model-addonesmat}}\begin{quote}

\sphinxAtStartPar
Add a dense symmetric matrix of value one to a model.

\sphinxAtStartPar
\sphinxstylestrong{Synopsis}
\begin{quote}

\sphinxAtStartPar
\sphinxcode{\sphinxupquote{SymMatrix AddOnesMat(int dim)}}
\end{quote}

\sphinxAtStartPar
\sphinxstylestrong{Arguments}
\begin{quote}

\sphinxAtStartPar
\sphinxcode{\sphinxupquote{dim}}: dimension of dense symmetric matrix.
\end{quote}

\sphinxAtStartPar
\sphinxstylestrong{Return}
\begin{quote}

\sphinxAtStartPar
new symmetric matrix object.
\end{quote}
\end{quote}

\subsubsection{Model::AddPsdConstr()}
\label{\detokenize{cppapi/model:model-addpsdconstr}}\begin{quote}

\sphinxAtStartPar
Add a PSD constraint to model.

\sphinxAtStartPar
\sphinxstylestrong{Synopsis}
\begin{quote}

\sphinxAtStartPar
\sphinxcode{\sphinxupquote{PsdConstraint AddPsdConstr(}}
\begin{quote}

\sphinxAtStartPar
\sphinxcode{\sphinxupquote{const PsdExpr \&expr,}}

\sphinxAtStartPar
\sphinxcode{\sphinxupquote{char sense,}}

\sphinxAtStartPar
\sphinxcode{\sphinxupquote{double rhs,}}

\sphinxAtStartPar
\sphinxcode{\sphinxupquote{const char *szName)}}
\end{quote}
\end{quote}

\sphinxAtStartPar
\sphinxstylestrong{Arguments}
\begin{quote}

\sphinxAtStartPar
\sphinxcode{\sphinxupquote{expr}}: PSD expression for new PSD contraint.

\sphinxAtStartPar
\sphinxcode{\sphinxupquote{sense}}: sense for new PSD constraint.

\sphinxAtStartPar
\sphinxcode{\sphinxupquote{rhs}}: double value at right side of the new PSD constraint.

\sphinxAtStartPar
\sphinxcode{\sphinxupquote{szName}}: optional, name of new PSD constraint.
\end{quote}

\sphinxAtStartPar
\sphinxstylestrong{Return}
\begin{quote}

\sphinxAtStartPar
new PSD constraint object.
\end{quote}
\end{quote}

\subsubsection{Model::AddPsdConstr()}
\label{\detokenize{cppapi/model:id37}}\begin{quote}

\sphinxAtStartPar
Add a PSD constraint to model.

\sphinxAtStartPar
\sphinxstylestrong{Synopsis}
\begin{quote}

\sphinxAtStartPar
\sphinxcode{\sphinxupquote{PsdConstraint AddPsdConstr(}}
\begin{quote}

\sphinxAtStartPar
\sphinxcode{\sphinxupquote{const PsdExpr \&expr,}}

\sphinxAtStartPar
\sphinxcode{\sphinxupquote{double lb,}}

\sphinxAtStartPar
\sphinxcode{\sphinxupquote{double ub,}}

\sphinxAtStartPar
\sphinxcode{\sphinxupquote{const char *szName)}}
\end{quote}
\end{quote}

\sphinxAtStartPar
\sphinxstylestrong{Arguments}
\begin{quote}

\sphinxAtStartPar
\sphinxcode{\sphinxupquote{expr}}: expression for new PSD constraint.

\sphinxAtStartPar
\sphinxcode{\sphinxupquote{lb}}: lower bound for ew PSD constraint.

\sphinxAtStartPar
\sphinxcode{\sphinxupquote{ub}}: upper bound for new PSD constraint

\sphinxAtStartPar
\sphinxcode{\sphinxupquote{szName}}: optional, name of new PSD constraint.
\end{quote}

\sphinxAtStartPar
\sphinxstylestrong{Return}
\begin{quote}

\sphinxAtStartPar
new PSD constraint object.
\end{quote}
\end{quote}

\subsubsection{Model::AddPsdConstr()}
\label{\detokenize{cppapi/model:id38}}\begin{quote}

\sphinxAtStartPar
Add a PSD constraint to model.

\sphinxAtStartPar
\sphinxstylestrong{Synopsis}
\begin{quote}

\sphinxAtStartPar
\sphinxcode{\sphinxupquote{PsdConstraint AddPsdConstr(}}
\begin{quote}

\sphinxAtStartPar
\sphinxcode{\sphinxupquote{const PsdExpr \&lhs,}}

\sphinxAtStartPar
\sphinxcode{\sphinxupquote{char sense,}}

\sphinxAtStartPar
\sphinxcode{\sphinxupquote{const PsdExpr \&rhs,}}

\sphinxAtStartPar
\sphinxcode{\sphinxupquote{const char *szName)}}
\end{quote}
\end{quote}

\sphinxAtStartPar
\sphinxstylestrong{Arguments}
\begin{quote}

\sphinxAtStartPar
\sphinxcode{\sphinxupquote{lhs}}: PSD expression at left side of new PSD constraint.

\sphinxAtStartPar
\sphinxcode{\sphinxupquote{sense}}: sense for new PSD constraint.

\sphinxAtStartPar
\sphinxcode{\sphinxupquote{rhs}}: PSD expression at right side of new PSD constraint.

\sphinxAtStartPar
\sphinxcode{\sphinxupquote{szName}}: optional, name of new PSD constraint.
\end{quote}

\sphinxAtStartPar
\sphinxstylestrong{Return}
\begin{quote}

\sphinxAtStartPar
new PSD constraint object.
\end{quote}
\end{quote}

\subsubsection{Model::AddPsdConstr()}
\label{\detokenize{cppapi/model:id39}}\begin{quote}

\sphinxAtStartPar
Add a PSD constraint to a model.

\sphinxAtStartPar
\sphinxstylestrong{Synopsis}
\begin{quote}

\sphinxAtStartPar
\sphinxcode{\sphinxupquote{PsdConstraint AddPsdConstr(const PsdConstrBuilder \&builder, const char *szName)}}
\end{quote}

\sphinxAtStartPar
\sphinxstylestrong{Arguments}
\begin{quote}

\sphinxAtStartPar
\sphinxcode{\sphinxupquote{builder}}: builder for new PSD constraint.

\sphinxAtStartPar
\sphinxcode{\sphinxupquote{szName}}: optional, name of new PSD constraint.
\end{quote}

\sphinxAtStartPar
\sphinxstylestrong{Return}
\begin{quote}

\sphinxAtStartPar
new PSD constraint object.
\end{quote}
\end{quote}

\subsubsection{Model::AddPsdVar()}
\label{\detokenize{cppapi/model:model-addpsdvar}}\begin{quote}

\sphinxAtStartPar
Add a new PSD variable to model.

\sphinxAtStartPar
\sphinxstylestrong{Synopsis}
\begin{quote}

\sphinxAtStartPar
\sphinxcode{\sphinxupquote{PsdVar AddPsdVar(int dim, const char *szName)}}
\end{quote}

\sphinxAtStartPar
\sphinxstylestrong{Arguments}
\begin{quote}

\sphinxAtStartPar
\sphinxcode{\sphinxupquote{dim}}: dimension of new PSD variable.

\sphinxAtStartPar
\sphinxcode{\sphinxupquote{szName}}: name of new PSD variable.
\end{quote}

\sphinxAtStartPar
\sphinxstylestrong{Return}
\begin{quote}

\sphinxAtStartPar
PSD variable object.
\end{quote}
\end{quote}

\subsubsection{Model::AddPsdVars()}
\label{\detokenize{cppapi/model:model-addpsdvars}}\begin{quote}

\sphinxAtStartPar
Add new PSD variables to model.

\sphinxAtStartPar
\sphinxstylestrong{Synopsis}
\begin{quote}

\sphinxAtStartPar
\sphinxcode{\sphinxupquote{PsdVarArray AddPsdVars(}}
\begin{quote}

\sphinxAtStartPar
\sphinxcode{\sphinxupquote{int count,}}

\sphinxAtStartPar
\sphinxcode{\sphinxupquote{int *pDim,}}

\sphinxAtStartPar
\sphinxcode{\sphinxupquote{const char *szPrefix)}}
\end{quote}
\end{quote}

\sphinxAtStartPar
\sphinxstylestrong{Arguments}
\begin{quote}

\sphinxAtStartPar
\sphinxcode{\sphinxupquote{count}}: number of new PSD variables.

\sphinxAtStartPar
\sphinxcode{\sphinxupquote{pDim}}: array of dimensions of new PSD variables.

\sphinxAtStartPar
\sphinxcode{\sphinxupquote{szPrefix}}: name prefix of new PSD variables.
\end{quote}

\sphinxAtStartPar
\sphinxstylestrong{Return}
\begin{quote}

\sphinxAtStartPar
array of PSD variable objects.
\end{quote}
\end{quote}

\subsubsection{Model::AddPsdVars()}
\label{\detokenize{cppapi/model:id40}}\begin{quote}

\sphinxAtStartPar
Add new PSD variables to model.

\sphinxAtStartPar
\sphinxstylestrong{Synopsis}
\begin{quote}

\sphinxAtStartPar
\sphinxcode{\sphinxupquote{PsdVarArray AddPsdVars(}}
\begin{quote}

\sphinxAtStartPar
\sphinxcode{\sphinxupquote{int count,}}

\sphinxAtStartPar
\sphinxcode{\sphinxupquote{int *pDim,}}

\sphinxAtStartPar
\sphinxcode{\sphinxupquote{const char *szNames,}}

\sphinxAtStartPar
\sphinxcode{\sphinxupquote{size\_t len)}}
\end{quote}
\end{quote}

\sphinxAtStartPar
\sphinxstylestrong{Arguments}
\begin{quote}

\sphinxAtStartPar
\sphinxcode{\sphinxupquote{count}}: number of new PSD variables.

\sphinxAtStartPar
\sphinxcode{\sphinxupquote{pDim}}: array of dimensions of new PSD variables.

\sphinxAtStartPar
\sphinxcode{\sphinxupquote{szNames}}: name buffer of new PSD variables.

\sphinxAtStartPar
\sphinxcode{\sphinxupquote{len}}: length of the name buffer.
\end{quote}

\sphinxAtStartPar
\sphinxstylestrong{Return}
\begin{quote}

\sphinxAtStartPar
array of PSD variable objects.
\end{quote}
\end{quote}

\subsubsection{Model::AddQConstr()}
\label{\detokenize{cppapi/model:model-addqconstr}}\begin{quote}

\sphinxAtStartPar
Add a quadratic constraint to model.

\sphinxAtStartPar
\sphinxstylestrong{Synopsis}
\begin{quote}

\sphinxAtStartPar
\sphinxcode{\sphinxupquote{QConstraint AddQConstr(}}
\begin{quote}

\sphinxAtStartPar
\sphinxcode{\sphinxupquote{const QuadExpr \&expr,}}

\sphinxAtStartPar
\sphinxcode{\sphinxupquote{char sense,}}

\sphinxAtStartPar
\sphinxcode{\sphinxupquote{double rhs,}}

\sphinxAtStartPar
\sphinxcode{\sphinxupquote{const char *szName)}}
\end{quote}
\end{quote}

\sphinxAtStartPar
\sphinxstylestrong{Arguments}
\begin{quote}

\sphinxAtStartPar
\sphinxcode{\sphinxupquote{expr}}: quadratic expression for the new contraint.

\sphinxAtStartPar
\sphinxcode{\sphinxupquote{sense}}: sense for new quadratic constraint.

\sphinxAtStartPar
\sphinxcode{\sphinxupquote{rhs}}: double value at right side of the new quadratic constraint.

\sphinxAtStartPar
\sphinxcode{\sphinxupquote{szName}}: optional, name of new quadratic constraint.
\end{quote}

\sphinxAtStartPar
\sphinxstylestrong{Return}
\begin{quote}

\sphinxAtStartPar
new quadratic constraint object.
\end{quote}
\end{quote}

\subsubsection{Model::AddQConstr()}
\label{\detokenize{cppapi/model:id41}}\begin{quote}

\sphinxAtStartPar
Add a quadratic constraint to model.

\sphinxAtStartPar
\sphinxstylestrong{Synopsis}
\begin{quote}

\sphinxAtStartPar
\sphinxcode{\sphinxupquote{QConstraint AddQConstr(}}
\begin{quote}

\sphinxAtStartPar
\sphinxcode{\sphinxupquote{const QuadExpr \&lhs,}}

\sphinxAtStartPar
\sphinxcode{\sphinxupquote{char sense,}}

\sphinxAtStartPar
\sphinxcode{\sphinxupquote{const QuadExpr \&rhs,}}

\sphinxAtStartPar
\sphinxcode{\sphinxupquote{const char *szName)}}
\end{quote}
\end{quote}

\sphinxAtStartPar
\sphinxstylestrong{Arguments}
\begin{quote}

\sphinxAtStartPar
\sphinxcode{\sphinxupquote{lhs}}: quadratic expression at left side of the new quadratic constraint.

\sphinxAtStartPar
\sphinxcode{\sphinxupquote{sense}}: sense for new quadratic constraint.

\sphinxAtStartPar
\sphinxcode{\sphinxupquote{rhs}}: quadratic expression at right side of the new quadratic constraint.

\sphinxAtStartPar
\sphinxcode{\sphinxupquote{szName}}: optional, name of new quadratic constraint.
\end{quote}

\sphinxAtStartPar
\sphinxstylestrong{Return}
\begin{quote}

\sphinxAtStartPar
new quadratic constraint object.
\end{quote}
\end{quote}

\subsubsection{Model::AddQConstr()}
\label{\detokenize{cppapi/model:id42}}\begin{quote}

\sphinxAtStartPar
Add a quadratic constraint to a model.

\sphinxAtStartPar
\sphinxstylestrong{Synopsis}
\begin{quote}

\sphinxAtStartPar
\sphinxcode{\sphinxupquote{QConstraint AddQConstr(const QConstrBuilder \&builder, const char *szName)}}
\end{quote}

\sphinxAtStartPar
\sphinxstylestrong{Arguments}
\begin{quote}

\sphinxAtStartPar
\sphinxcode{\sphinxupquote{builder}}: builder for the new quadratic constraint.

\sphinxAtStartPar
\sphinxcode{\sphinxupquote{szName}}: optional, name of new quadratic constraint.
\end{quote}

\sphinxAtStartPar
\sphinxstylestrong{Return}
\begin{quote}

\sphinxAtStartPar
new quadratic constraint object.
\end{quote}
\end{quote}

\subsubsection{Model::AddRow()}
\label{\detokenize{cppapi/model:model-addrow}}\begin{quote}

\sphinxAtStartPar
Add linear constraint to model via advanced interface.

\sphinxAtStartPar
\sphinxstylestrong{Synopsis}
\begin{quote}

\sphinxAtStartPar
\sphinxcode{\sphinxupquote{AddRow(}}
\begin{quote}

\sphinxAtStartPar
\sphinxcode{\sphinxupquote{int rowCnt,}}

\sphinxAtStartPar
\sphinxcode{\sphinxupquote{const int *rowIdx,}}

\sphinxAtStartPar
\sphinxcode{\sphinxupquote{const double *rowElem,}}

\sphinxAtStartPar
\sphinxcode{\sphinxupquote{char sense,}}

\sphinxAtStartPar
\sphinxcode{\sphinxupquote{double rhs,}}

\sphinxAtStartPar
\sphinxcode{\sphinxupquote{const char *szName)}}
\end{quote}
\end{quote}

\sphinxAtStartPar
\sphinxstylestrong{Arguments}
\begin{quote}

\sphinxAtStartPar
\sphinxcode{\sphinxupquote{rowCnt}}: number of terms in the linear constraint.

\sphinxAtStartPar
\sphinxcode{\sphinxupquote{rowIdx}}: array of variable indexes in the linear constraint.

\sphinxAtStartPar
\sphinxcode{\sphinxupquote{rowElem}}: array of coefficients in the linear constraint.

\sphinxAtStartPar
\sphinxcode{\sphinxupquote{sense}}: constraint sense, excluding range bound.

\sphinxAtStartPar
\sphinxcode{\sphinxupquote{rhs}}: right hand side value of constraint.

\sphinxAtStartPar
\sphinxcode{\sphinxupquote{szName}}: name of new constraint, with default value of empty.
\end{quote}
\end{quote}

\subsubsection{Model::AddRow()}
\label{\detokenize{cppapi/model:id43}}\begin{quote}

\sphinxAtStartPar
Add linear constraint to model via advanced interface.

\sphinxAtStartPar
\sphinxstylestrong{Synopsis}
\begin{quote}

\sphinxAtStartPar
\sphinxcode{\sphinxupquote{AddRow(}}
\begin{quote}

\sphinxAtStartPar
\sphinxcode{\sphinxupquote{int rowCnt,}}

\sphinxAtStartPar
\sphinxcode{\sphinxupquote{const int *rowIdx,}}

\sphinxAtStartPar
\sphinxcode{\sphinxupquote{const double *rowElem,}}

\sphinxAtStartPar
\sphinxcode{\sphinxupquote{double lb,}}

\sphinxAtStartPar
\sphinxcode{\sphinxupquote{double ub,}}

\sphinxAtStartPar
\sphinxcode{\sphinxupquote{const char *szName)}}
\end{quote}
\end{quote}

\sphinxAtStartPar
\sphinxstylestrong{Arguments}
\begin{quote}

\sphinxAtStartPar
\sphinxcode{\sphinxupquote{rowCnt}}: number of terms in the linear constraint.

\sphinxAtStartPar
\sphinxcode{\sphinxupquote{rowIdx}}: array of variable indexes in the linear constraint.

\sphinxAtStartPar
\sphinxcode{\sphinxupquote{rowElem}}: array of coefficients in the linear constraint.

\sphinxAtStartPar
\sphinxcode{\sphinxupquote{lb}}: lower bound of new constraint.

\sphinxAtStartPar
\sphinxcode{\sphinxupquote{ub}}: upper bound of new constraint.

\sphinxAtStartPar
\sphinxcode{\sphinxupquote{szName}}: name of new constraint, with default value of empty.
\end{quote}
\end{quote}

\subsubsection{Model::AddRows()}
\label{\detokenize{cppapi/model:model-addrows}}\begin{quote}

\sphinxAtStartPar
Add linear constraints to model via advanced interface.

\sphinxAtStartPar
\sphinxstylestrong{Synopsis}
\begin{quote}

\sphinxAtStartPar
\sphinxcode{\sphinxupquote{ConstrArray AddRows(}}
\begin{quote}

\sphinxAtStartPar
\sphinxcode{\sphinxupquote{int count,}}

\sphinxAtStartPar
\sphinxcode{\sphinxupquote{char *pSense,}}

\sphinxAtStartPar
\sphinxcode{\sphinxupquote{double *pRhs,}}

\sphinxAtStartPar
\sphinxcode{\sphinxupquote{char const* const *names)}}
\end{quote}
\end{quote}

\sphinxAtStartPar
\sphinxstylestrong{Arguments}
\begin{quote}

\sphinxAtStartPar
\sphinxcode{\sphinxupquote{count}}: number of constraints added to model.

\sphinxAtStartPar
\sphinxcode{\sphinxupquote{pSense}}: senses for new constraints, excluding range bound.

\sphinxAtStartPar
\sphinxcode{\sphinxupquote{pRhs}}: right hand side values for new constraints.

\sphinxAtStartPar
\sphinxcode{\sphinxupquote{names}}: array of names of new constraints, with default value of empty.
\end{quote}

\sphinxAtStartPar
\sphinxstylestrong{Return}
\begin{quote}

\sphinxAtStartPar
new constraint array objects.
\end{quote}
\end{quote}

\subsubsection{Model::AddRows()}
\label{\detokenize{cppapi/model:id44}}\begin{quote}

\sphinxAtStartPar
Add linear constraints to model via advanced interface.

\sphinxAtStartPar
\sphinxstylestrong{Synopsis}
\begin{quote}

\sphinxAtStartPar
\sphinxcode{\sphinxupquote{ConstrArray AddRows(}}
\begin{quote}

\sphinxAtStartPar
\sphinxcode{\sphinxupquote{int count,}}

\sphinxAtStartPar
\sphinxcode{\sphinxupquote{double *plb,}}

\sphinxAtStartPar
\sphinxcode{\sphinxupquote{double *pub,}}

\sphinxAtStartPar
\sphinxcode{\sphinxupquote{char const* const *names)}}
\end{quote}
\end{quote}

\sphinxAtStartPar
\sphinxstylestrong{Arguments}
\begin{quote}

\sphinxAtStartPar
\sphinxcode{\sphinxupquote{count}}: number of constraints added to model.

\sphinxAtStartPar
\sphinxcode{\sphinxupquote{plb}}: lower bounds for new constraints.

\sphinxAtStartPar
\sphinxcode{\sphinxupquote{pub}}: upper bounds for new constraints.

\sphinxAtStartPar
\sphinxcode{\sphinxupquote{names}}: array of names of new constraints, with default value of empty.
\end{quote}

\sphinxAtStartPar
\sphinxstylestrong{Return}
\begin{quote}

\sphinxAtStartPar
new constraint array objects.
\end{quote}
\end{quote}

\subsubsection{Model::AddRows()}
\label{\detokenize{cppapi/model:id45}}\begin{quote}

\sphinxAtStartPar
With coefficient matrix in CSR format, add linear constraints to model via advanced interface.

\sphinxAtStartPar
\sphinxstylestrong{Synopsis}
\begin{quote}

\sphinxAtStartPar
\sphinxcode{\sphinxupquote{ConstrArray AddRows(}}
\begin{quote}

\sphinxAtStartPar
\sphinxcode{\sphinxupquote{int count,}}

\sphinxAtStartPar
\sphinxcode{\sphinxupquote{int *rowBeg,}}

\sphinxAtStartPar
\sphinxcode{\sphinxupquote{int *rowCnt,}}

\sphinxAtStartPar
\sphinxcode{\sphinxupquote{int *rowIdx,}}

\sphinxAtStartPar
\sphinxcode{\sphinxupquote{double *rowElem,}}

\sphinxAtStartPar
\sphinxcode{\sphinxupquote{char *rowSense,}}

\sphinxAtStartPar
\sphinxcode{\sphinxupquote{double *rowBound,}}

\sphinxAtStartPar
\sphinxcode{\sphinxupquote{double *rowUpper,}}

\sphinxAtStartPar
\sphinxcode{\sphinxupquote{char const* const *names)}}
\end{quote}
\end{quote}

\sphinxAtStartPar
\sphinxstylestrong{Arguments}
\begin{quote}

\sphinxAtStartPar
\sphinxcode{\sphinxupquote{count}}: number of constraints added to model.

\sphinxAtStartPar
\sphinxcode{\sphinxupquote{rowBeg}}: indexes of begin elements in CSR format. If empty, constraints are added without coefficients.

\sphinxAtStartPar
\sphinxcode{\sphinxupquote{rowCnt}}: count of nonzero elements in each row. If empty, use rowBeg to calculate count.

\sphinxAtStartPar
\sphinxcode{\sphinxupquote{rowIdx}}: variable indexes of each row in CSR format.

\sphinxAtStartPar
\sphinxcode{\sphinxupquote{rowElem}}: corresponding coefficients of each row in CSR format.

\sphinxAtStartPar
\sphinxcode{\sphinxupquote{rowSense}}: senses for new constraints, including range bound. If empty, rowBound and rowUpper are used as lower and upper bounds respectively.

\sphinxAtStartPar
\sphinxcode{\sphinxupquote{rowBound}}: bounds for new constraints, used for sense of equal, less\_than and greater\_than.

\sphinxAtStartPar
\sphinxcode{\sphinxupquote{rowUpper}}: bounds for new constraints, used as range bound or upper bound.

\sphinxAtStartPar
\sphinxcode{\sphinxupquote{names}}: array of names of new constraints, with default value of empty.
\end{quote}

\sphinxAtStartPar
\sphinxstylestrong{Return}
\begin{quote}

\sphinxAtStartPar
new constraint array objects.
\end{quote}
\end{quote}

\subsubsection{Model::AddSos()}
\label{\detokenize{cppapi/model:model-addsos}}\begin{quote}

\sphinxAtStartPar
Add a SOS constraint to model.

\sphinxAtStartPar
\sphinxstylestrong{Synopsis}
\begin{quote}

\sphinxAtStartPar
\sphinxcode{\sphinxupquote{Sos AddSos(const SosBuilder \&builder)}}
\end{quote}

\sphinxAtStartPar
\sphinxstylestrong{Arguments}
\begin{quote}

\sphinxAtStartPar
\sphinxcode{\sphinxupquote{builder}}: builder for new SOS constraint.
\end{quote}

\sphinxAtStartPar
\sphinxstylestrong{Return}
\begin{quote}

\sphinxAtStartPar
new SOS constraint object.
\end{quote}
\end{quote}

\subsubsection{Model::AddSos()}
\label{\detokenize{cppapi/model:id46}}\begin{quote}

\sphinxAtStartPar
Add a SOS constraint to model.

\sphinxAtStartPar
\sphinxstylestrong{Synopsis}
\begin{quote}

\sphinxAtStartPar
\sphinxcode{\sphinxupquote{Sos AddSos(}}
\begin{quote}

\sphinxAtStartPar
\sphinxcode{\sphinxupquote{const VarArray \&vars,}}

\sphinxAtStartPar
\sphinxcode{\sphinxupquote{const double *pWeights,}}

\sphinxAtStartPar
\sphinxcode{\sphinxupquote{int type)}}
\end{quote}
\end{quote}

\sphinxAtStartPar
\sphinxstylestrong{Arguments}
\begin{quote}

\sphinxAtStartPar
\sphinxcode{\sphinxupquote{vars}}: variables that participate in the SOS constraint.

\sphinxAtStartPar
\sphinxcode{\sphinxupquote{pWeights}}: optional, weights for variables in the SOS constraint.

\sphinxAtStartPar
\sphinxcode{\sphinxupquote{type}}: type of SOS constraint.
\end{quote}

\sphinxAtStartPar
\sphinxstylestrong{Return}
\begin{quote}

\sphinxAtStartPar
new SOS constraint object.
\end{quote}
\end{quote}

\subsubsection{Model::AddSparseMat()}
\label{\detokenize{cppapi/model:model-addsparsemat}}\begin{quote}

\sphinxAtStartPar
Add a sparse symmetric matrix to a model.

\sphinxAtStartPar
\sphinxstylestrong{Synopsis}
\begin{quote}

\sphinxAtStartPar
\sphinxcode{\sphinxupquote{SymMatrix AddSparseMat(}}
\begin{quote}

\sphinxAtStartPar
\sphinxcode{\sphinxupquote{int dim,}}

\sphinxAtStartPar
\sphinxcode{\sphinxupquote{int nElems,}}

\sphinxAtStartPar
\sphinxcode{\sphinxupquote{int *pRows,}}

\sphinxAtStartPar
\sphinxcode{\sphinxupquote{int *pCols,}}

\sphinxAtStartPar
\sphinxcode{\sphinxupquote{double *pVals)}}
\end{quote}
\end{quote}

\sphinxAtStartPar
\sphinxstylestrong{Arguments}
\begin{quote}

\sphinxAtStartPar
\sphinxcode{\sphinxupquote{dim}}: dimension of the sparse symmetric matrix.

\sphinxAtStartPar
\sphinxcode{\sphinxupquote{nElems}}: number of non\sphinxhyphen{}zero elements in the sparse symmetric matrix.

\sphinxAtStartPar
\sphinxcode{\sphinxupquote{pRows}}: array of row indexes of non\sphinxhyphen{}zero elements.

\sphinxAtStartPar
\sphinxcode{\sphinxupquote{pCols}}: array of col indexes of non\sphinxhyphen{}zero elements.

\sphinxAtStartPar
\sphinxcode{\sphinxupquote{pVals}}: array of values of non\sphinxhyphen{}zero elements.
\end{quote}

\sphinxAtStartPar
\sphinxstylestrong{Return}
\begin{quote}

\sphinxAtStartPar
new symmetric matrix object.
\end{quote}
\end{quote}

\subsubsection{Model::AddSymMat()}
\label{\detokenize{cppapi/model:model-addsymmat}}\begin{quote}

\sphinxAtStartPar
Given a symmetric matrix expression, add results matrix to model.

\sphinxAtStartPar
\sphinxstylestrong{Synopsis}
\begin{quote}

\sphinxAtStartPar
\sphinxcode{\sphinxupquote{SymMatrix AddSymMat(const SymMatExpr \&expr)}}
\end{quote}

\sphinxAtStartPar
\sphinxstylestrong{Arguments}
\begin{quote}

\sphinxAtStartPar
\sphinxcode{\sphinxupquote{expr}}: symmetric matrix expression object.
\end{quote}

\sphinxAtStartPar
\sphinxstylestrong{Return}
\begin{quote}

\sphinxAtStartPar
results symmetric matrix object.
\end{quote}
\end{quote}

\subsubsection{Model::AddUserCut()}
\label{\detokenize{cppapi/model:model-addusercut}}\begin{quote}

\sphinxAtStartPar
Add a user cut to model.

\sphinxAtStartPar
\sphinxstylestrong{Synopsis}
\begin{quote}

\sphinxAtStartPar
\sphinxcode{\sphinxupquote{void AddUserCut(}}
\begin{quote}

\sphinxAtStartPar
\sphinxcode{\sphinxupquote{const Expr \&lhs,}}

\sphinxAtStartPar
\sphinxcode{\sphinxupquote{char sense,}}

\sphinxAtStartPar
\sphinxcode{\sphinxupquote{double rhs,}}

\sphinxAtStartPar
\sphinxcode{\sphinxupquote{const char *szName)}}
\end{quote}
\end{quote}

\sphinxAtStartPar
\sphinxstylestrong{Arguments}
\begin{quote}

\sphinxAtStartPar
\sphinxcode{\sphinxupquote{lhs}}: expression for user cut.

\sphinxAtStartPar
\sphinxcode{\sphinxupquote{sense}}: sense for user cut.

\sphinxAtStartPar
\sphinxcode{\sphinxupquote{rhs}}: right hand side value for user cut.

\sphinxAtStartPar
\sphinxcode{\sphinxupquote{szName}}: optional, name of user cut.
\end{quote}
\end{quote}

\subsubsection{Model::AddUserCut()}
\label{\detokenize{cppapi/model:id47}}\begin{quote}

\sphinxAtStartPar
Add a user cut to model.

\sphinxAtStartPar
\sphinxstylestrong{Synopsis}
\begin{quote}

\sphinxAtStartPar
\sphinxcode{\sphinxupquote{void AddUserCut(}}
\begin{quote}

\sphinxAtStartPar
\sphinxcode{\sphinxupquote{const Expr \&lhs,}}

\sphinxAtStartPar
\sphinxcode{\sphinxupquote{char sense,}}

\sphinxAtStartPar
\sphinxcode{\sphinxupquote{const Expr \&rhs,}}

\sphinxAtStartPar
\sphinxcode{\sphinxupquote{const char *szName)}}
\end{quote}
\end{quote}

\sphinxAtStartPar
\sphinxstylestrong{Arguments}
\begin{quote}

\sphinxAtStartPar
\sphinxcode{\sphinxupquote{lhs}}: left hand side expression for user cut.

\sphinxAtStartPar
\sphinxcode{\sphinxupquote{sense}}: sense for user cut.

\sphinxAtStartPar
\sphinxcode{\sphinxupquote{rhs}}: right hand side expression for user cut.

\sphinxAtStartPar
\sphinxcode{\sphinxupquote{szName}}: optional, name of user cut.
\end{quote}
\end{quote}

\subsubsection{Model::AddUserCut()}
\label{\detokenize{cppapi/model:id48}}\begin{quote}

\sphinxAtStartPar
Add a user cut to model.

\sphinxAtStartPar
\sphinxstylestrong{Synopsis}
\begin{quote}

\sphinxAtStartPar
\sphinxcode{\sphinxupquote{void AddUserCut(const ConstrBuilder \&builder, const char *szName)}}
\end{quote}

\sphinxAtStartPar
\sphinxstylestrong{Arguments}
\begin{quote}

\sphinxAtStartPar
\sphinxcode{\sphinxupquote{builder}}: builder for user cut.

\sphinxAtStartPar
\sphinxcode{\sphinxupquote{szName}}: optional, name of user cut.
\end{quote}
\end{quote}

\subsubsection{Model::AddUserCuts()}
\label{\detokenize{cppapi/model:model-addusercuts}}\begin{quote}

\sphinxAtStartPar
Add user cuts to model.

\sphinxAtStartPar
\sphinxstylestrong{Synopsis}
\begin{quote}

\sphinxAtStartPar
\sphinxcode{\sphinxupquote{void AddUserCuts(const ConstrBuilderArray \&builders, const char *szPrefix)}}
\end{quote}

\sphinxAtStartPar
\sphinxstylestrong{Arguments}
\begin{quote}

\sphinxAtStartPar
\sphinxcode{\sphinxupquote{builders}}: array of builders for user cuts.

\sphinxAtStartPar
\sphinxcode{\sphinxupquote{szPrefix}}: name prefix of new user cuts.
\end{quote}
\end{quote}

\subsubsection{Model::AddVar()}
\label{\detokenize{cppapi/model:model-addvar}}\begin{quote}

\sphinxAtStartPar
Add a variable to model.

\sphinxAtStartPar
\sphinxstylestrong{Synopsis}
\begin{quote}

\sphinxAtStartPar
\sphinxcode{\sphinxupquote{Var AddVar(}}
\begin{quote}

\sphinxAtStartPar
\sphinxcode{\sphinxupquote{double lb,}}

\sphinxAtStartPar
\sphinxcode{\sphinxupquote{double ub,}}

\sphinxAtStartPar
\sphinxcode{\sphinxupquote{double obj,}}

\sphinxAtStartPar
\sphinxcode{\sphinxupquote{char vtype,}}

\sphinxAtStartPar
\sphinxcode{\sphinxupquote{const char *szName)}}
\end{quote}
\end{quote}

\sphinxAtStartPar
\sphinxstylestrong{Arguments}
\begin{quote}

\sphinxAtStartPar
\sphinxcode{\sphinxupquote{lb}}: lower bound for new variable.

\sphinxAtStartPar
\sphinxcode{\sphinxupquote{ub}}: upper bound for new variable.

\sphinxAtStartPar
\sphinxcode{\sphinxupquote{obj}}: objective coefficient for new variable.

\sphinxAtStartPar
\sphinxcode{\sphinxupquote{vtype}}: variable type for new variable.

\sphinxAtStartPar
\sphinxcode{\sphinxupquote{szName}}: optional, name for new variable.
\end{quote}

\sphinxAtStartPar
\sphinxstylestrong{Return}
\begin{quote}

\sphinxAtStartPar
new variable object.
\end{quote}
\end{quote}

\subsubsection{Model::AddVar()}
\label{\detokenize{cppapi/model:id49}}\begin{quote}

\sphinxAtStartPar
Add a variable and the associated non\sphinxhyphen{}zero coefficients as column.

\sphinxAtStartPar
\sphinxstylestrong{Synopsis}
\begin{quote}

\sphinxAtStartPar
\sphinxcode{\sphinxupquote{Var AddVar(}}
\begin{quote}

\sphinxAtStartPar
\sphinxcode{\sphinxupquote{double lb,}}

\sphinxAtStartPar
\sphinxcode{\sphinxupquote{double ub,}}

\sphinxAtStartPar
\sphinxcode{\sphinxupquote{double obj,}}

\sphinxAtStartPar
\sphinxcode{\sphinxupquote{char vtype,}}

\sphinxAtStartPar
\sphinxcode{\sphinxupquote{const Column \&col,}}

\sphinxAtStartPar
\sphinxcode{\sphinxupquote{const char *szName)}}
\end{quote}
\end{quote}

\sphinxAtStartPar
\sphinxstylestrong{Arguments}
\begin{quote}

\sphinxAtStartPar
\sphinxcode{\sphinxupquote{lb}}: lower bound for new variable.

\sphinxAtStartPar
\sphinxcode{\sphinxupquote{ub}}: upper bound for new variable.

\sphinxAtStartPar
\sphinxcode{\sphinxupquote{obj}}: objective coefficient for new variable.

\sphinxAtStartPar
\sphinxcode{\sphinxupquote{vtype}}: variable type for new variable.

\sphinxAtStartPar
\sphinxcode{\sphinxupquote{col}}: column object for specifying a set of constraints to which the variable belongs.

\sphinxAtStartPar
\sphinxcode{\sphinxupquote{szName}}: optional, name for new variable.
\end{quote}

\sphinxAtStartPar
\sphinxstylestrong{Return}
\begin{quote}

\sphinxAtStartPar
new variable object.
\end{quote}
\end{quote}

\subsubsection{Model::AddVars()}
\label{\detokenize{cppapi/model:model-addvars}}\begin{quote}

\sphinxAtStartPar
Add new variables to model.

\sphinxAtStartPar
\sphinxstylestrong{Synopsis}
\begin{quote}

\sphinxAtStartPar
\sphinxcode{\sphinxupquote{VarArray AddVars(}}
\begin{quote}

\sphinxAtStartPar
\sphinxcode{\sphinxupquote{int count,}}

\sphinxAtStartPar
\sphinxcode{\sphinxupquote{char vtype,}}

\sphinxAtStartPar
\sphinxcode{\sphinxupquote{const char *szPrefix)}}
\end{quote}
\end{quote}

\sphinxAtStartPar
\sphinxstylestrong{Arguments}
\begin{quote}

\sphinxAtStartPar
\sphinxcode{\sphinxupquote{count}}: the number of variables to add.

\sphinxAtStartPar
\sphinxcode{\sphinxupquote{vtype}}: variable types for new variables.

\sphinxAtStartPar
\sphinxcode{\sphinxupquote{szPrefix}}: prefix part for names of new variables.
\end{quote}

\sphinxAtStartPar
\sphinxstylestrong{Return}
\begin{quote}

\sphinxAtStartPar
array of new variable objects.
\end{quote}
\end{quote}

\subsubsection{Model::AddVars()}
\label{\detokenize{cppapi/model:id50}}\begin{quote}

\sphinxAtStartPar
Add new variables to model.

\sphinxAtStartPar
\sphinxstylestrong{Synopsis}
\begin{quote}

\sphinxAtStartPar
\sphinxcode{\sphinxupquote{VarArray AddVars(}}
\begin{quote}

\sphinxAtStartPar
\sphinxcode{\sphinxupquote{int count,}}

\sphinxAtStartPar
\sphinxcode{\sphinxupquote{char vtype,}}

\sphinxAtStartPar
\sphinxcode{\sphinxupquote{const char *szNames,}}

\sphinxAtStartPar
\sphinxcode{\sphinxupquote{size\_t len)}}
\end{quote}
\end{quote}

\sphinxAtStartPar
\sphinxstylestrong{Arguments}
\begin{quote}

\sphinxAtStartPar
\sphinxcode{\sphinxupquote{count}}: the number of variables to add.

\sphinxAtStartPar
\sphinxcode{\sphinxupquote{vtype}}: variable types for new variables.

\sphinxAtStartPar
\sphinxcode{\sphinxupquote{szNames}}: name buffer for new variables.

\sphinxAtStartPar
\sphinxcode{\sphinxupquote{len}}: length of name buffer.
\end{quote}

\sphinxAtStartPar
\sphinxstylestrong{Return}
\begin{quote}

\sphinxAtStartPar
array of new variable objects.
\end{quote}
\end{quote}

\subsubsection{Model::AddVars()}
\label{\detokenize{cppapi/model:id51}}\begin{quote}

\sphinxAtStartPar
Add new variables to model.

\sphinxAtStartPar
\sphinxstylestrong{Synopsis}
\begin{quote}

\sphinxAtStartPar
\sphinxcode{\sphinxupquote{VarArray AddVars(}}
\begin{quote}

\sphinxAtStartPar
\sphinxcode{\sphinxupquote{int count,}}

\sphinxAtStartPar
\sphinxcode{\sphinxupquote{double lb,}}

\sphinxAtStartPar
\sphinxcode{\sphinxupquote{double ub,}}

\sphinxAtStartPar
\sphinxcode{\sphinxupquote{double obj,}}

\sphinxAtStartPar
\sphinxcode{\sphinxupquote{char vtype,}}

\sphinxAtStartPar
\sphinxcode{\sphinxupquote{const char *szPrefix)}}
\end{quote}
\end{quote}

\sphinxAtStartPar
\sphinxstylestrong{Arguments}
\begin{quote}

\sphinxAtStartPar
\sphinxcode{\sphinxupquote{count}}: the number of variables to add.

\sphinxAtStartPar
\sphinxcode{\sphinxupquote{lb}}: lower bound for new variable.

\sphinxAtStartPar
\sphinxcode{\sphinxupquote{ub}}: upper bound for new variable.

\sphinxAtStartPar
\sphinxcode{\sphinxupquote{obj}}: objective coefficient for new variable.

\sphinxAtStartPar
\sphinxcode{\sphinxupquote{vtype}}: variable types for new variables.

\sphinxAtStartPar
\sphinxcode{\sphinxupquote{szPrefix}}: prefix part for names of new variables.
\end{quote}

\sphinxAtStartPar
\sphinxstylestrong{Return}
\begin{quote}

\sphinxAtStartPar
array of new variable objects.
\end{quote}
\end{quote}

\subsubsection{Model::AddVars()}
\label{\detokenize{cppapi/model:id52}}\begin{quote}

\sphinxAtStartPar
Add new variables to model.

\sphinxAtStartPar
\sphinxstylestrong{Synopsis}
\begin{quote}

\sphinxAtStartPar
\sphinxcode{\sphinxupquote{VarArray AddVars(}}
\begin{quote}

\sphinxAtStartPar
\sphinxcode{\sphinxupquote{int count,}}

\sphinxAtStartPar
\sphinxcode{\sphinxupquote{double lb,}}

\sphinxAtStartPar
\sphinxcode{\sphinxupquote{double ub,}}

\sphinxAtStartPar
\sphinxcode{\sphinxupquote{double obj,}}

\sphinxAtStartPar
\sphinxcode{\sphinxupquote{char vtype,}}

\sphinxAtStartPar
\sphinxcode{\sphinxupquote{const char *szNames,}}

\sphinxAtStartPar
\sphinxcode{\sphinxupquote{size\_t len)}}
\end{quote}
\end{quote}

\sphinxAtStartPar
\sphinxstylestrong{Arguments}
\begin{quote}

\sphinxAtStartPar
\sphinxcode{\sphinxupquote{count}}: the number of variables to add.

\sphinxAtStartPar
\sphinxcode{\sphinxupquote{lb}}: lower bound for new variable.

\sphinxAtStartPar
\sphinxcode{\sphinxupquote{ub}}: upper bound for new variable.

\sphinxAtStartPar
\sphinxcode{\sphinxupquote{obj}}: objective coefficient for new variable.

\sphinxAtStartPar
\sphinxcode{\sphinxupquote{vtype}}: variable types for new variables.

\sphinxAtStartPar
\sphinxcode{\sphinxupquote{szNames}}: name buffer for new variables.

\sphinxAtStartPar
\sphinxcode{\sphinxupquote{len}}: length of name buffer.
\end{quote}

\sphinxAtStartPar
\sphinxstylestrong{Return}
\begin{quote}

\sphinxAtStartPar
array of new variable objects.
\end{quote}
\end{quote}

\subsubsection{Model::AddVars()}
\label{\detokenize{cppapi/model:id53}}\begin{quote}

\sphinxAtStartPar
Add new variables to model.

\sphinxAtStartPar
\sphinxstylestrong{Synopsis}
\begin{quote}

\sphinxAtStartPar
\sphinxcode{\sphinxupquote{VarArray AddVars(}}
\begin{quote}

\sphinxAtStartPar
\sphinxcode{\sphinxupquote{int count,}}

\sphinxAtStartPar
\sphinxcode{\sphinxupquote{double *plb,}}

\sphinxAtStartPar
\sphinxcode{\sphinxupquote{double *pub,}}

\sphinxAtStartPar
\sphinxcode{\sphinxupquote{double *pobj,}}

\sphinxAtStartPar
\sphinxcode{\sphinxupquote{char *pvtype,}}

\sphinxAtStartPar
\sphinxcode{\sphinxupquote{const char *szPrefix)}}
\end{quote}
\end{quote}

\sphinxAtStartPar
\sphinxstylestrong{Arguments}
\begin{quote}

\sphinxAtStartPar
\sphinxcode{\sphinxupquote{count}}: the number of variables to add.

\sphinxAtStartPar
\sphinxcode{\sphinxupquote{plb}}: lower bounds for new variables. if NULL, lower bounds are 0.0.

\sphinxAtStartPar
\sphinxcode{\sphinxupquote{pub}}: upper bounds for new variables. if NULL, upper bounds are infinity or 1 for binary variables.

\sphinxAtStartPar
\sphinxcode{\sphinxupquote{pobj}}: objective coefficients for new variables. if NULL, objective coefficients are 0.0.

\sphinxAtStartPar
\sphinxcode{\sphinxupquote{pvtype}}: variable types for new variables. if NULL, variable types are continuous.

\sphinxAtStartPar
\sphinxcode{\sphinxupquote{szPrefix}}: prefix part for names of new variables.
\end{quote}

\sphinxAtStartPar
\sphinxstylestrong{Return}
\begin{quote}

\sphinxAtStartPar
array of new variable objects.
\end{quote}
\end{quote}

\subsubsection{Model::AddVars()}
\label{\detokenize{cppapi/model:id54}}\begin{quote}

\sphinxAtStartPar
Add new decision variables to a model.

\sphinxAtStartPar
\sphinxstylestrong{Synopsis}
\begin{quote}

\sphinxAtStartPar
\sphinxcode{\sphinxupquote{VarArray AddVars(}}
\begin{quote}

\sphinxAtStartPar
\sphinxcode{\sphinxupquote{int count,}}

\sphinxAtStartPar
\sphinxcode{\sphinxupquote{double *plb,}}

\sphinxAtStartPar
\sphinxcode{\sphinxupquote{double *pub,}}

\sphinxAtStartPar
\sphinxcode{\sphinxupquote{double *pobj,}}

\sphinxAtStartPar
\sphinxcode{\sphinxupquote{char *pvtype,}}

\sphinxAtStartPar
\sphinxcode{\sphinxupquote{const char *szNames,}}

\sphinxAtStartPar
\sphinxcode{\sphinxupquote{size\_t len)}}
\end{quote}
\end{quote}

\sphinxAtStartPar
\sphinxstylestrong{Arguments}
\begin{quote}

\sphinxAtStartPar
\sphinxcode{\sphinxupquote{count}}: the number of variables to add.

\sphinxAtStartPar
\sphinxcode{\sphinxupquote{plb}}: lower bounds for new variables. if NULL, lower bounds are 0.0.

\sphinxAtStartPar
\sphinxcode{\sphinxupquote{pub}}: upper bounds for new variables. if NULL, upper bounds are infinity or 1 for binary variables.

\sphinxAtStartPar
\sphinxcode{\sphinxupquote{pobj}}: objective coefficients for new variables. if NULL, objective coefficients are 0.0.

\sphinxAtStartPar
\sphinxcode{\sphinxupquote{pvtype}}: variable types for new variables. if NULL, variable types are continuous.

\sphinxAtStartPar
\sphinxcode{\sphinxupquote{szNames}}: name buffer for new variables.

\sphinxAtStartPar
\sphinxcode{\sphinxupquote{len}}: length of name buffer.
\end{quote}

\sphinxAtStartPar
\sphinxstylestrong{Return}
\begin{quote}

\sphinxAtStartPar
array of new variable objects.
\end{quote}
\end{quote}

\subsubsection{Model::AddVars()}
\label{\detokenize{cppapi/model:id55}}\begin{quote}

\sphinxAtStartPar
Add new variables to model.

\sphinxAtStartPar
\sphinxstylestrong{Synopsis}
\begin{quote}

\sphinxAtStartPar
\sphinxcode{\sphinxupquote{VarArray AddVars(}}
\begin{quote}

\sphinxAtStartPar
\sphinxcode{\sphinxupquote{int count,}}

\sphinxAtStartPar
\sphinxcode{\sphinxupquote{double *plb,}}

\sphinxAtStartPar
\sphinxcode{\sphinxupquote{double *pub,}}

\sphinxAtStartPar
\sphinxcode{\sphinxupquote{double *pobj,}}

\sphinxAtStartPar
\sphinxcode{\sphinxupquote{char *pvtype,}}

\sphinxAtStartPar
\sphinxcode{\sphinxupquote{const ColumnArray \&cols,}}

\sphinxAtStartPar
\sphinxcode{\sphinxupquote{const char *szPrefix)}}
\end{quote}
\end{quote}

\sphinxAtStartPar
\sphinxstylestrong{Arguments}
\begin{quote}

\sphinxAtStartPar
\sphinxcode{\sphinxupquote{count}}: the number of variables to add.

\sphinxAtStartPar
\sphinxcode{\sphinxupquote{plb}}: lower bounds for new variables. if NULL, lower bounds are 0.0.

\sphinxAtStartPar
\sphinxcode{\sphinxupquote{pub}}: upper bounds for new variables. if NULL, upper bounds are infinity or 1 for binary variables.

\sphinxAtStartPar
\sphinxcode{\sphinxupquote{pobj}}: objective coefficients for new variables. if NULL, objective coefficients are 0.0.

\sphinxAtStartPar
\sphinxcode{\sphinxupquote{pvtype}}: variable types for new variables. if NULL, variable types are continuous.

\sphinxAtStartPar
\sphinxcode{\sphinxupquote{cols}}: column objects for specifying a set of constraints to which each new variable belongs.

\sphinxAtStartPar
\sphinxcode{\sphinxupquote{szPrefix}}: prefix part for names of new variables.
\end{quote}

\sphinxAtStartPar
\sphinxstylestrong{Return}
\begin{quote}

\sphinxAtStartPar
array of new variable objects.
\end{quote}
\end{quote}

\subsubsection{Model::Clear()}
\label{\detokenize{cppapi/model:model-clear}}\begin{quote}

\sphinxAtStartPar
Clear all settings including problem itself.

\sphinxAtStartPar
\sphinxstylestrong{Synopsis}
\begin{quote}

\sphinxAtStartPar
\sphinxcode{\sphinxupquote{void Clear()}}
\end{quote}
\end{quote}

\subsubsection{Model::Clone()}
\label{\detokenize{cppapi/model:model-clone}}\begin{quote}

\sphinxAtStartPar
Deep copy COPT model.

\sphinxAtStartPar
\sphinxstylestrong{Synopsis}
\begin{quote}

\sphinxAtStartPar
\sphinxcode{\sphinxupquote{Model Clone()}}
\end{quote}

\sphinxAtStartPar
\sphinxstylestrong{Return}
\begin{quote}

\sphinxAtStartPar
cloned model object.
\end{quote}
\end{quote}

\subsubsection{Model::ComputeIIS()}
\label{\detokenize{cppapi/model:model-computeiis}}\begin{quote}

\sphinxAtStartPar
Compute IIS for infeasible model.

\sphinxAtStartPar
\sphinxstylestrong{Synopsis}
\begin{quote}

\sphinxAtStartPar
\sphinxcode{\sphinxupquote{void ComputeIIS()}}
\end{quote}
\end{quote}

\subsubsection{Model::DelNlObj()}
\label{\detokenize{cppapi/model:model-delnlobj}}\begin{quote}

\sphinxAtStartPar
Delete nonlinear part of objective in model.

\sphinxAtStartPar
\sphinxstylestrong{Synopsis}
\begin{quote}

\sphinxAtStartPar
\sphinxcode{\sphinxupquote{void DelNlObj()}}
\end{quote}
\end{quote}

\subsubsection{Model::DelObjN()}
\label{\detokenize{cppapi/model:model-delobjn}}\begin{quote}

\sphinxAtStartPar
Delete linear part of specific multi\sphinxhyphen{}objective function in model.

\sphinxAtStartPar
\sphinxstylestrong{Synopsis}
\begin{quote}

\sphinxAtStartPar
\sphinxcode{\sphinxupquote{void DelObjN(int idx)}}
\end{quote}

\sphinxAtStartPar
\sphinxstylestrong{Arguments}
\begin{quote}

\sphinxAtStartPar
\sphinxcode{\sphinxupquote{idx}}: index of a multi\sphinxhyphen{}objective function.
\end{quote}
\end{quote}

\subsubsection{Model::DelPsdObj()}
\label{\detokenize{cppapi/model:model-delpsdobj}}\begin{quote}

\sphinxAtStartPar
Delete PSD part of objective in model.

\sphinxAtStartPar
\sphinxstylestrong{Synopsis}
\begin{quote}

\sphinxAtStartPar
\sphinxcode{\sphinxupquote{void DelPsdObj()}}
\end{quote}
\end{quote}

\subsubsection{Model::DelQuadObj()}
\label{\detokenize{cppapi/model:model-delquadobj}}\begin{quote}

\sphinxAtStartPar
Delete quadratic part of objective in model.

\sphinxAtStartPar
\sphinxstylestrong{Synopsis}
\begin{quote}

\sphinxAtStartPar
\sphinxcode{\sphinxupquote{void DelQuadObj()}}
\end{quote}
\end{quote}

\subsubsection{Model::FeasRelax()}
\label{\detokenize{cppapi/model:model-feasrelax}}\begin{quote}

\sphinxAtStartPar
Compute feasibility relaxation for infeasible model.

\sphinxAtStartPar
\sphinxstylestrong{Synopsis}
\begin{quote}

\sphinxAtStartPar
\sphinxcode{\sphinxupquote{void FeasRelax(}}
\begin{quote}

\sphinxAtStartPar
\sphinxcode{\sphinxupquote{const VarArray \&vars,}}

\sphinxAtStartPar
\sphinxcode{\sphinxupquote{double *pColLowPen,}}

\sphinxAtStartPar
\sphinxcode{\sphinxupquote{double *pColUppPen,}}

\sphinxAtStartPar
\sphinxcode{\sphinxupquote{const ConstrArray \&cons,}}

\sphinxAtStartPar
\sphinxcode{\sphinxupquote{double *pRowBndPen,}}

\sphinxAtStartPar
\sphinxcode{\sphinxupquote{double *pRowUppPen)}}
\end{quote}
\end{quote}

\sphinxAtStartPar
\sphinxstylestrong{Arguments}
\begin{quote}

\sphinxAtStartPar
\sphinxcode{\sphinxupquote{vars}}: an array of variables.

\sphinxAtStartPar
\sphinxcode{\sphinxupquote{pColLowPen}}: penalties for lower bounds of variables.

\sphinxAtStartPar
\sphinxcode{\sphinxupquote{pColUppPen}}: penalties for upper bounds of variables.

\sphinxAtStartPar
\sphinxcode{\sphinxupquote{cons}}: an array of constraints.

\sphinxAtStartPar
\sphinxcode{\sphinxupquote{pRowBndPen}}: penalties for right hand sides of constraints.

\sphinxAtStartPar
\sphinxcode{\sphinxupquote{pRowUppPen}}: penalties for upper bounds of range constraints.
\end{quote}
\end{quote}

\subsubsection{Model::FeasRelax()}
\label{\detokenize{cppapi/model:id56}}\begin{quote}

\sphinxAtStartPar
Compute feasibility relaxation for infeasible model.

\sphinxAtStartPar
\sphinxstylestrong{Synopsis}
\begin{quote}

\sphinxAtStartPar
\sphinxcode{\sphinxupquote{void FeasRelax(}}
\begin{quote}

\sphinxAtStartPar
\sphinxcode{\sphinxupquote{const MVar\textless{}1\textgreater{} \&vars,}}

\sphinxAtStartPar
\sphinxcode{\sphinxupquote{double *pColLowPen,}}

\sphinxAtStartPar
\sphinxcode{\sphinxupquote{double *pColUppPen,}}

\sphinxAtStartPar
\sphinxcode{\sphinxupquote{const MConstr\textless{}1\textgreater{} \&cons,}}

\sphinxAtStartPar
\sphinxcode{\sphinxupquote{double *pRowBndPen,}}

\sphinxAtStartPar
\sphinxcode{\sphinxupquote{double *pRowUppPen)}}
\end{quote}
\end{quote}

\sphinxAtStartPar
\sphinxstylestrong{Arguments}
\begin{quote}

\sphinxAtStartPar
\sphinxcode{\sphinxupquote{vars}}: one\sphinxhyphen{}dimensional variables.

\sphinxAtStartPar
\sphinxcode{\sphinxupquote{pColLowPen}}: penalties for lower bounds of variables.

\sphinxAtStartPar
\sphinxcode{\sphinxupquote{pColUppPen}}: penalties for upper bounds of variables.

\sphinxAtStartPar
\sphinxcode{\sphinxupquote{cons}}: one\sphinxhyphen{}dimensional constraints.

\sphinxAtStartPar
\sphinxcode{\sphinxupquote{pRowBndPen}}: penalties for right hand sides of constraints.

\sphinxAtStartPar
\sphinxcode{\sphinxupquote{pRowUppPen}}: penalties for upper bounds of range constraints.
\end{quote}
\end{quote}

\subsubsection{Model::FeasRelax()}
\label{\detokenize{cppapi/model:id57}}\begin{quote}

\sphinxAtStartPar
Compute feasibility relaxation for infeasible model.

\sphinxAtStartPar
\sphinxstylestrong{Synopsis}
\begin{quote}

\sphinxAtStartPar
\sphinxcode{\sphinxupquote{void FeasRelax(}}
\begin{quote}

\sphinxAtStartPar
\sphinxcode{\sphinxupquote{const MVar\textless{}N\textgreater{} \&vars,}}

\sphinxAtStartPar
\sphinxcode{\sphinxupquote{const NdArray\textless{}double, N\textgreater{} \&colLowPen,}}

\sphinxAtStartPar
\sphinxcode{\sphinxupquote{const NdArray\textless{}double, N\textgreater{} \&colUppPen,}}

\sphinxAtStartPar
\sphinxcode{\sphinxupquote{const MConstr\textless{}M\textgreater{} \&cons,}}

\sphinxAtStartPar
\sphinxcode{\sphinxupquote{const NdArray\textless{}double, M\textgreater{} \&rowBndPen,}}

\sphinxAtStartPar
\sphinxcode{\sphinxupquote{const NdArray\textless{}double, M\textgreater{} \&rowUppPen)}}
\end{quote}
\end{quote}

\sphinxAtStartPar
\sphinxstylestrong{Arguments}
\begin{quote}

\sphinxAtStartPar
\sphinxcode{\sphinxupquote{vars}}: N\sphinxhyphen{}dimensional variables.

\sphinxAtStartPar
\sphinxcode{\sphinxupquote{colLowPen}}: penalties for lower bounds of variables.

\sphinxAtStartPar
\sphinxcode{\sphinxupquote{colUppPen}}: penalties for upper bounds of variables.

\sphinxAtStartPar
\sphinxcode{\sphinxupquote{cons}}: M\sphinxhyphen{}dimensional constraints.

\sphinxAtStartPar
\sphinxcode{\sphinxupquote{rowBndPen}}: penalties for right hand sides of constraints.

\sphinxAtStartPar
\sphinxcode{\sphinxupquote{rowUppPen}}: penalties for upper bounds of range constraints.
\end{quote}
\end{quote}

\subsubsection{Model::FeasRelax()}
\label{\detokenize{cppapi/model:id58}}\begin{quote}

\sphinxAtStartPar
Compute feasibility relaxation for infeasible model.

\sphinxAtStartPar
\sphinxstylestrong{Synopsis}
\begin{quote}

\sphinxAtStartPar
\sphinxcode{\sphinxupquote{void FeasRelax(int ifRelaxVars, int ifRelaxCons)}}
\end{quote}

\sphinxAtStartPar
\sphinxstylestrong{Arguments}
\begin{quote}

\sphinxAtStartPar
\sphinxcode{\sphinxupquote{ifRelaxVars}}: whether to relax variables.

\sphinxAtStartPar
\sphinxcode{\sphinxupquote{ifRelaxCons}}: whether to relax constraints.
\end{quote}
\end{quote}

\subsubsection{Model::Get()}
\label{\detokenize{cppapi/model:model-get}}\begin{quote}

\sphinxAtStartPar
Query values of information associated with variables.

\sphinxAtStartPar
\sphinxstylestrong{Synopsis}
\begin{quote}

\sphinxAtStartPar
\sphinxcode{\sphinxupquote{int Get(}}
\begin{quote}

\sphinxAtStartPar
\sphinxcode{\sphinxupquote{const char *szName,}}

\sphinxAtStartPar
\sphinxcode{\sphinxupquote{const VarArray \&vars,}}

\sphinxAtStartPar
\sphinxcode{\sphinxupquote{double *pOut)}}
\end{quote}
\end{quote}

\sphinxAtStartPar
\sphinxstylestrong{Arguments}
\begin{quote}

\sphinxAtStartPar
\sphinxcode{\sphinxupquote{szName}}: name of information.

\sphinxAtStartPar
\sphinxcode{\sphinxupquote{vars}}: a list of desired variables.

\sphinxAtStartPar
\sphinxcode{\sphinxupquote{pOut}}: output array of information values.
\end{quote}

\sphinxAtStartPar
\sphinxstylestrong{Return}
\begin{quote}

\sphinxAtStartPar
the number of valid variables. If failed, return \sphinxhyphen{}1.
\end{quote}
\end{quote}

\subsubsection{Model::Get()}
\label{\detokenize{cppapi/model:id59}}\begin{quote}

\sphinxAtStartPar
Query values of information associated with constraints.

\sphinxAtStartPar
\sphinxstylestrong{Synopsis}
\begin{quote}

\sphinxAtStartPar
\sphinxcode{\sphinxupquote{int Get(}}
\begin{quote}

\sphinxAtStartPar
\sphinxcode{\sphinxupquote{const char *szName,}}

\sphinxAtStartPar
\sphinxcode{\sphinxupquote{const ConstrArray \&constrs,}}

\sphinxAtStartPar
\sphinxcode{\sphinxupquote{double *pOut)}}
\end{quote}
\end{quote}

\sphinxAtStartPar
\sphinxstylestrong{Arguments}
\begin{quote}

\sphinxAtStartPar
\sphinxcode{\sphinxupquote{szName}}: name of information.

\sphinxAtStartPar
\sphinxcode{\sphinxupquote{constrs}}: a list of desired constraints.

\sphinxAtStartPar
\sphinxcode{\sphinxupquote{pOut}}: output array of information values.
\end{quote}

\sphinxAtStartPar
\sphinxstylestrong{Return}
\begin{quote}

\sphinxAtStartPar
the number of valid constraints. If failed, return \sphinxhyphen{}1.
\end{quote}
\end{quote}

\subsubsection{Model::Get()}
\label{\detokenize{cppapi/model:id60}}\begin{quote}

\sphinxAtStartPar
Query values of information associated with nonlinear constraints.

\sphinxAtStartPar
\sphinxstylestrong{Synopsis}
\begin{quote}

\sphinxAtStartPar
\sphinxcode{\sphinxupquote{int Get(}}
\begin{quote}

\sphinxAtStartPar
\sphinxcode{\sphinxupquote{const char *szName,}}

\sphinxAtStartPar
\sphinxcode{\sphinxupquote{const NlConstrArray \&constrs,}}

\sphinxAtStartPar
\sphinxcode{\sphinxupquote{double *pOut)}}
\end{quote}
\end{quote}

\sphinxAtStartPar
\sphinxstylestrong{Arguments}
\begin{quote}

\sphinxAtStartPar
\sphinxcode{\sphinxupquote{szName}}: name of information.

\sphinxAtStartPar
\sphinxcode{\sphinxupquote{constrs}}: a list of desired nonlinear constraints.

\sphinxAtStartPar
\sphinxcode{\sphinxupquote{pOut}}: output array of information values.
\end{quote}

\sphinxAtStartPar
\sphinxstylestrong{Return}
\begin{quote}

\sphinxAtStartPar
the number of valid nonlinear constraints. If failed, return \sphinxhyphen{}1.
\end{quote}
\end{quote}

\subsubsection{Model::Get()}
\label{\detokenize{cppapi/model:id61}}\begin{quote}

\sphinxAtStartPar
Query values of information associated with quadratic constraints.

\sphinxAtStartPar
\sphinxstylestrong{Synopsis}
\begin{quote}

\sphinxAtStartPar
\sphinxcode{\sphinxupquote{int Get(}}
\begin{quote}

\sphinxAtStartPar
\sphinxcode{\sphinxupquote{const char *szName,}}

\sphinxAtStartPar
\sphinxcode{\sphinxupquote{const QConstrArray \&constrs,}}

\sphinxAtStartPar
\sphinxcode{\sphinxupquote{double *pOut)}}
\end{quote}
\end{quote}

\sphinxAtStartPar
\sphinxstylestrong{Arguments}
\begin{quote}

\sphinxAtStartPar
\sphinxcode{\sphinxupquote{szName}}: name of information.

\sphinxAtStartPar
\sphinxcode{\sphinxupquote{constrs}}: a list of desired quadratic constraints.

\sphinxAtStartPar
\sphinxcode{\sphinxupquote{pOut}}: output array of information values.
\end{quote}

\sphinxAtStartPar
\sphinxstylestrong{Return}
\begin{quote}

\sphinxAtStartPar
the number of valid quadratic constraints. If failed, return \sphinxhyphen{}1.
\end{quote}
\end{quote}

\subsubsection{Model::Get()}
\label{\detokenize{cppapi/model:id62}}\begin{quote}

\sphinxAtStartPar
Query values of information associated with PSD constraints.

\sphinxAtStartPar
\sphinxstylestrong{Synopsis}
\begin{quote}

\sphinxAtStartPar
\sphinxcode{\sphinxupquote{int Get(}}
\begin{quote}

\sphinxAtStartPar
\sphinxcode{\sphinxupquote{const char *szName,}}

\sphinxAtStartPar
\sphinxcode{\sphinxupquote{const PsdConstrArray \&constrs,}}

\sphinxAtStartPar
\sphinxcode{\sphinxupquote{double *pOut)}}
\end{quote}
\end{quote}

\sphinxAtStartPar
\sphinxstylestrong{Arguments}
\begin{quote}

\sphinxAtStartPar
\sphinxcode{\sphinxupquote{szName}}: name of information.

\sphinxAtStartPar
\sphinxcode{\sphinxupquote{constrs}}: a list of desired PSD constraints.

\sphinxAtStartPar
\sphinxcode{\sphinxupquote{pOut}}: output array of information values.
\end{quote}

\sphinxAtStartPar
\sphinxstylestrong{Return}
\begin{quote}

\sphinxAtStartPar
the number of valid PSD constraints. If failed, return \sphinxhyphen{}1.
\end{quote}
\end{quote}

\subsubsection{Model::GetAffineCone()}
\label{\detokenize{cppapi/model:model-getaffinecone}}\begin{quote}

\sphinxAtStartPar
Get an affine cone constraint of given index in model.

\sphinxAtStartPar
\sphinxstylestrong{Synopsis}
\begin{quote}

\sphinxAtStartPar
\sphinxcode{\sphinxupquote{AffineCone GetAffineCone(int idx)}}
\end{quote}

\sphinxAtStartPar
\sphinxstylestrong{Arguments}
\begin{quote}

\sphinxAtStartPar
\sphinxcode{\sphinxupquote{idx}}: index of the desired affine cone constraint.
\end{quote}

\sphinxAtStartPar
\sphinxstylestrong{Return}
\begin{quote}

\sphinxAtStartPar
the desired affine cone constraint object.
\end{quote}
\end{quote}

\subsubsection{Model::GetAffineConeBuilder()}
\label{\detokenize{cppapi/model:model-getaffineconebuilder}}\begin{quote}

\sphinxAtStartPar
Get builder of given affine cone constraint in model.

\sphinxAtStartPar
\sphinxstylestrong{Synopsis}
\begin{quote}

\sphinxAtStartPar
\sphinxcode{\sphinxupquote{AffineConeBuilder GetAffineConeBuilder(const AffineCone \&cone)}}
\end{quote}

\sphinxAtStartPar
\sphinxstylestrong{Arguments}
\begin{quote}

\sphinxAtStartPar
\sphinxcode{\sphinxupquote{cone}}: an affine cone constraint.
\end{quote}

\sphinxAtStartPar
\sphinxstylestrong{Return}
\begin{quote}

\sphinxAtStartPar
desired affine cone constraint builder.
\end{quote}
\end{quote}

\subsubsection{Model::GetAffineConeBuilders()}
\label{\detokenize{cppapi/model:model-getaffineconebuilders}}\begin{quote}

\sphinxAtStartPar
Get builders of desired affine cone constraints in model.

\sphinxAtStartPar
\sphinxstylestrong{Synopsis}
\begin{quote}

\sphinxAtStartPar
\sphinxcode{\sphinxupquote{AffineConeBuilderArray GetAffineConeBuilders(const AffineConeArray \&cones)}}
\end{quote}

\sphinxAtStartPar
\sphinxstylestrong{Arguments}
\begin{quote}

\sphinxAtStartPar
\sphinxcode{\sphinxupquote{cones}}: array of affine cone constraints.
\end{quote}

\sphinxAtStartPar
\sphinxstylestrong{Return}
\begin{quote}

\sphinxAtStartPar
array object of desired affine cone constraint builders.
\end{quote}
\end{quote}

\subsubsection{Model::GetAffineConeBuilders()}
\label{\detokenize{cppapi/model:id63}}\begin{quote}

\sphinxAtStartPar
Get builders of all affine cone constraints in model.

\sphinxAtStartPar
\sphinxstylestrong{Synopsis}
\begin{quote}

\sphinxAtStartPar
\sphinxcode{\sphinxupquote{AffineConeBuilderArray GetAffineConeBuilders()}}
\end{quote}

\sphinxAtStartPar
\sphinxstylestrong{Return}
\begin{quote}

\sphinxAtStartPar
array object of all affine cone constraint builders.
\end{quote}
\end{quote}

\subsubsection{Model::GetAffineConeByName()}
\label{\detokenize{cppapi/model:model-getaffineconebyname}}\begin{quote}

\sphinxAtStartPar
Get an affine cone constraint of given name in model.

\sphinxAtStartPar
\sphinxstylestrong{Synopsis}
\begin{quote}

\sphinxAtStartPar
\sphinxcode{\sphinxupquote{AffineCone GetAffineConeByName(const char *szName)}}
\end{quote}

\sphinxAtStartPar
\sphinxstylestrong{Arguments}
\begin{quote}

\sphinxAtStartPar
\sphinxcode{\sphinxupquote{szName}}: name of the affine cone constraint.
\end{quote}

\sphinxAtStartPar
\sphinxstylestrong{Return}
\begin{quote}

\sphinxAtStartPar
the desired affine cone constraint object.
\end{quote}
\end{quote}

\subsubsection{Model::GetAffineCones()}
\label{\detokenize{cppapi/model:model-getaffinecones}}\begin{quote}

\sphinxAtStartPar
Get all affine cone constraints in model.

\sphinxAtStartPar
\sphinxstylestrong{Synopsis}
\begin{quote}

\sphinxAtStartPar
\sphinxcode{\sphinxupquote{AffineConeArray GetAffineCones()}}
\end{quote}

\sphinxAtStartPar
\sphinxstylestrong{Return}
\begin{quote}

\sphinxAtStartPar
array object of affine cone constraints.
\end{quote}
\end{quote}

\subsubsection{Model::GetCoeff()}
\label{\detokenize{cppapi/model:model-getcoeff}}\begin{quote}

\sphinxAtStartPar
Get the coefficient of variable in a linear constraint.

\sphinxAtStartPar
\sphinxstylestrong{Synopsis}
\begin{quote}

\sphinxAtStartPar
\sphinxcode{\sphinxupquote{double GetCoeff(const Constraint \&constr, const Var \&var)}}
\end{quote}

\sphinxAtStartPar
\sphinxstylestrong{Arguments}
\begin{quote}

\sphinxAtStartPar
\sphinxcode{\sphinxupquote{constr}}: The requested constraint.

\sphinxAtStartPar
\sphinxcode{\sphinxupquote{var}}: The requested variable.
\end{quote}

\sphinxAtStartPar
\sphinxstylestrong{Return}
\begin{quote}

\sphinxAtStartPar
The requested coefficient.
\end{quote}
\end{quote}

\subsubsection{Model::GetCol()}
\label{\detokenize{cppapi/model:model-getcol}}\begin{quote}

\sphinxAtStartPar
Get a column object that have a list of constraints in which the variable participates.

\sphinxAtStartPar
\sphinxstylestrong{Synopsis}
\begin{quote}

\sphinxAtStartPar
\sphinxcode{\sphinxupquote{Column GetCol(const Var \&var)}}
\end{quote}

\sphinxAtStartPar
\sphinxstylestrong{Arguments}
\begin{quote}

\sphinxAtStartPar
\sphinxcode{\sphinxupquote{var}}: a variable object.
\end{quote}

\sphinxAtStartPar
\sphinxstylestrong{Return}
\begin{quote}

\sphinxAtStartPar
a column object associated with a variable.
\end{quote}
\end{quote}

\subsubsection{Model::GetColBasis()}
\label{\detokenize{cppapi/model:model-getcolbasis}}\begin{quote}

\sphinxAtStartPar
Get status of column basis.

\sphinxAtStartPar
\sphinxstylestrong{Synopsis}
\begin{quote}

\sphinxAtStartPar
\sphinxcode{\sphinxupquote{int GetColBasis(int *pBasis)}}
\end{quote}

\sphinxAtStartPar
\sphinxstylestrong{Arguments}
\begin{quote}

\sphinxAtStartPar
\sphinxcode{\sphinxupquote{pBasis}}: integer pointer to basis status.
\end{quote}

\sphinxAtStartPar
\sphinxstylestrong{Return}
\begin{quote}

\sphinxAtStartPar
number of columns.
\end{quote}
\end{quote}

\subsubsection{Model::GetCone()}
\label{\detokenize{cppapi/model:model-getcone}}\begin{quote}

\sphinxAtStartPar
Get a cone constraint of given index in model.

\sphinxAtStartPar
\sphinxstylestrong{Synopsis}
\begin{quote}

\sphinxAtStartPar
\sphinxcode{\sphinxupquote{Cone GetCone(int idx)}}
\end{quote}

\sphinxAtStartPar
\sphinxstylestrong{Arguments}
\begin{quote}

\sphinxAtStartPar
\sphinxcode{\sphinxupquote{idx}}: index of the desired cone constraint.
\end{quote}

\sphinxAtStartPar
\sphinxstylestrong{Return}
\begin{quote}

\sphinxAtStartPar
the desired cone constraint object.
\end{quote}
\end{quote}

\subsubsection{Model::GetConeBuilders()}
\label{\detokenize{cppapi/model:model-getconebuilders}}\begin{quote}

\sphinxAtStartPar
Get builders of all cone constraints in model.

\sphinxAtStartPar
\sphinxstylestrong{Synopsis}
\begin{quote}

\sphinxAtStartPar
\sphinxcode{\sphinxupquote{ConeBuilderArray GetConeBuilders()}}
\end{quote}

\sphinxAtStartPar
\sphinxstylestrong{Return}
\begin{quote}

\sphinxAtStartPar
array object of all cone constraint builders.
\end{quote}
\end{quote}

\subsubsection{Model::GetConeBuilders()}
\label{\detokenize{cppapi/model:id64}}\begin{quote}

\sphinxAtStartPar
Get builders of given cone constraints in model.

\sphinxAtStartPar
\sphinxstylestrong{Synopsis}
\begin{quote}

\sphinxAtStartPar
\sphinxcode{\sphinxupquote{ConeBuilderArray GetConeBuilders(const ConeArray \&cones)}}
\end{quote}

\sphinxAtStartPar
\sphinxstylestrong{Arguments}
\begin{quote}

\sphinxAtStartPar
\sphinxcode{\sphinxupquote{cones}}: array of cone constraints.
\end{quote}

\sphinxAtStartPar
\sphinxstylestrong{Return}
\begin{quote}

\sphinxAtStartPar
array object of desired cone constraint builders.
\end{quote}
\end{quote}

\subsubsection{Model::GetCones()}
\label{\detokenize{cppapi/model:model-getcones}}\begin{quote}

\sphinxAtStartPar
Get all cone constraints in model.

\sphinxAtStartPar
\sphinxstylestrong{Synopsis}
\begin{quote}

\sphinxAtStartPar
\sphinxcode{\sphinxupquote{ConeArray GetCones()}}
\end{quote}

\sphinxAtStartPar
\sphinxstylestrong{Return}
\begin{quote}

\sphinxAtStartPar
array object of cone constraints.
\end{quote}
\end{quote}

\subsubsection{Model::GetConstr()}
\label{\detokenize{cppapi/model:model-getconstr}}\begin{quote}

\sphinxAtStartPar
Get a constraint of given index in model.

\sphinxAtStartPar
\sphinxstylestrong{Synopsis}
\begin{quote}

\sphinxAtStartPar
\sphinxcode{\sphinxupquote{Constraint GetConstr(int idx)}}
\end{quote}

\sphinxAtStartPar
\sphinxstylestrong{Arguments}
\begin{quote}

\sphinxAtStartPar
\sphinxcode{\sphinxupquote{idx}}: index of the desired constraint.
\end{quote}

\sphinxAtStartPar
\sphinxstylestrong{Return}
\begin{quote}

\sphinxAtStartPar
the desired constraint object.
\end{quote}
\end{quote}

\subsubsection{Model::GetConstrBuilder()}
\label{\detokenize{cppapi/model:model-getconstrbuilder}}\begin{quote}

\sphinxAtStartPar
Get builder of a constraint in model, including variables and associated coefficients, sense and RHS.

\sphinxAtStartPar
\sphinxstylestrong{Synopsis}
\begin{quote}

\sphinxAtStartPar
\sphinxcode{\sphinxupquote{ConstrBuilder GetConstrBuilder(Constraint constr)}}
\end{quote}

\sphinxAtStartPar
\sphinxstylestrong{Arguments}
\begin{quote}

\sphinxAtStartPar
\sphinxcode{\sphinxupquote{constr}}: a constraint object.
\end{quote}

\sphinxAtStartPar
\sphinxstylestrong{Return}
\begin{quote}

\sphinxAtStartPar
constraint builder object.
\end{quote}
\end{quote}

\subsubsection{Model::GetConstrBuilders()}
\label{\detokenize{cppapi/model:model-getconstrbuilders}}\begin{quote}

\sphinxAtStartPar
Get builders of all constraints in model.

\sphinxAtStartPar
\sphinxstylestrong{Synopsis}
\begin{quote}

\sphinxAtStartPar
\sphinxcode{\sphinxupquote{ConstrBuilderArray GetConstrBuilders()}}
\end{quote}

\sphinxAtStartPar
\sphinxstylestrong{Return}
\begin{quote}

\sphinxAtStartPar
array object of constraint builders.
\end{quote}
\end{quote}

\subsubsection{Model::GetConstrByName()}
\label{\detokenize{cppapi/model:model-getconstrbyname}}\begin{quote}

\sphinxAtStartPar
Get a constraint of given name in model.

\sphinxAtStartPar
\sphinxstylestrong{Synopsis}
\begin{quote}

\sphinxAtStartPar
\sphinxcode{\sphinxupquote{Constraint GetConstrByName(const char *szName)}}
\end{quote}

\sphinxAtStartPar
\sphinxstylestrong{Arguments}
\begin{quote}

\sphinxAtStartPar
\sphinxcode{\sphinxupquote{szName}}: name of the desired constraint.
\end{quote}

\sphinxAtStartPar
\sphinxstylestrong{Return}
\begin{quote}

\sphinxAtStartPar
the desired constraint object.
\end{quote}
\end{quote}

\subsubsection{Model::GetConstrLowerIIS()}
\label{\detokenize{cppapi/model:model-getconstrloweriis}}\begin{quote}

\sphinxAtStartPar
Get IIS status of lower bounds of constraints.

\sphinxAtStartPar
\sphinxstylestrong{Synopsis}
\begin{quote}

\sphinxAtStartPar
\sphinxcode{\sphinxupquote{int GetConstrLowerIIS(const ConstrArray \&constrs, int *pLowerIIS)}}
\end{quote}

\sphinxAtStartPar
\sphinxstylestrong{Arguments}
\begin{quote}

\sphinxAtStartPar
\sphinxcode{\sphinxupquote{constrs}}: Array of constraints.

\sphinxAtStartPar
\sphinxcode{\sphinxupquote{pLowerIIS}}: IIS status of lower bounds of constraints.
\end{quote}

\sphinxAtStartPar
\sphinxstylestrong{Return}
\begin{quote}

\sphinxAtStartPar
Number of constraints.
\end{quote}
\end{quote}

\subsubsection{Model::GetConstrs()}
\label{\detokenize{cppapi/model:model-getconstrs}}\begin{quote}

\sphinxAtStartPar
Get all constraints in model.

\sphinxAtStartPar
\sphinxstylestrong{Synopsis}
\begin{quote}

\sphinxAtStartPar
\sphinxcode{\sphinxupquote{ConstrArray GetConstrs()}}
\end{quote}

\sphinxAtStartPar
\sphinxstylestrong{Return}
\begin{quote}

\sphinxAtStartPar
array object of constraints.
\end{quote}
\end{quote}

\subsubsection{Model::GetConstrUpperIIS()}
\label{\detokenize{cppapi/model:model-getconstrupperiis}}\begin{quote}

\sphinxAtStartPar
Get IIS status of upper bounds of constraints.

\sphinxAtStartPar
\sphinxstylestrong{Synopsis}
\begin{quote}

\sphinxAtStartPar
\sphinxcode{\sphinxupquote{int GetConstrUpperIIS(const ConstrArray \&constrs, int *pUpperIIS)}}
\end{quote}

\sphinxAtStartPar
\sphinxstylestrong{Arguments}
\begin{quote}

\sphinxAtStartPar
\sphinxcode{\sphinxupquote{constrs}}: Array of constraints.

\sphinxAtStartPar
\sphinxcode{\sphinxupquote{pUpperIIS}}: IIS status of upper bounds of constraints.
\end{quote}

\sphinxAtStartPar
\sphinxstylestrong{Return}
\begin{quote}

\sphinxAtStartPar
Number of constraints.
\end{quote}
\end{quote}

\subsubsection{Model::GetDblAttr()}
\label{\detokenize{cppapi/model:model-getdblattr}}\begin{quote}

\sphinxAtStartPar
Get value of a COPT double attribute.

\sphinxAtStartPar
\sphinxstylestrong{Synopsis}
\begin{quote}

\sphinxAtStartPar
\sphinxcode{\sphinxupquote{double GetDblAttr(const char *szAttr)}}
\end{quote}

\sphinxAtStartPar
\sphinxstylestrong{Arguments}
\begin{quote}

\sphinxAtStartPar
\sphinxcode{\sphinxupquote{szAttr}}: name of double attribute.
\end{quote}

\sphinxAtStartPar
\sphinxstylestrong{Return}
\begin{quote}

\sphinxAtStartPar
value of double attribute.
\end{quote}
\end{quote}

\subsubsection{Model::GetDblAttrN()}
\label{\detokenize{cppapi/model:model-getdblattrn}}\begin{quote}

\sphinxAtStartPar
Get value of a double attribute of a multi\sphinxhyphen{}objective function.

\sphinxAtStartPar
\sphinxstylestrong{Synopsis}
\begin{quote}

\sphinxAtStartPar
\sphinxcode{\sphinxupquote{double GetDblAttrN(int idx, const char *szAttr)}}
\end{quote}

\sphinxAtStartPar
\sphinxstylestrong{Arguments}
\begin{quote}

\sphinxAtStartPar
\sphinxcode{\sphinxupquote{idx}}: index of a multi\sphinxhyphen{}objective function.

\sphinxAtStartPar
\sphinxcode{\sphinxupquote{szAttr}}: name of double attribute.
\end{quote}

\sphinxAtStartPar
\sphinxstylestrong{Return}
\begin{quote}

\sphinxAtStartPar
value of double attribute.
\end{quote}
\end{quote}

\subsubsection{Model::GetDblParam()}
\label{\detokenize{cppapi/model:model-getdblparam}}\begin{quote}

\sphinxAtStartPar
Get value of a COPT double parameter.

\sphinxAtStartPar
\sphinxstylestrong{Synopsis}
\begin{quote}

\sphinxAtStartPar
\sphinxcode{\sphinxupquote{double GetDblParam(const char *szParam)}}
\end{quote}

\sphinxAtStartPar
\sphinxstylestrong{Arguments}
\begin{quote}

\sphinxAtStartPar
\sphinxcode{\sphinxupquote{szParam}}: name of double parameter.
\end{quote}

\sphinxAtStartPar
\sphinxstylestrong{Return}
\begin{quote}

\sphinxAtStartPar
value of double parameter.
\end{quote}
\end{quote}

\subsubsection{Model::GetDblParamN()}
\label{\detokenize{cppapi/model:model-getdblparamn}}\begin{quote}

\sphinxAtStartPar
Get value of a double parameter of a multi\sphinxhyphen{}objective function.

\sphinxAtStartPar
\sphinxstylestrong{Synopsis}
\begin{quote}

\sphinxAtStartPar
\sphinxcode{\sphinxupquote{double GetDblParamN(int idx, const char *szParam)}}
\end{quote}

\sphinxAtStartPar
\sphinxstylestrong{Arguments}
\begin{quote}

\sphinxAtStartPar
\sphinxcode{\sphinxupquote{idx}}: index of a multi\sphinxhyphen{}objective function.

\sphinxAtStartPar
\sphinxcode{\sphinxupquote{szParam}}: name of double parameter.
\end{quote}

\sphinxAtStartPar
\sphinxstylestrong{Return}
\begin{quote}

\sphinxAtStartPar
value of double parameter.
\end{quote}
\end{quote}

\subsubsection{Model::GetExpCone()}
\label{\detokenize{cppapi/model:model-getexpcone}}\begin{quote}

\sphinxAtStartPar
Get an exponential cone constraint of given index in model.

\sphinxAtStartPar
\sphinxstylestrong{Synopsis}
\begin{quote}

\sphinxAtStartPar
\sphinxcode{\sphinxupquote{ExpCone GetExpCone(int idx)}}
\end{quote}

\sphinxAtStartPar
\sphinxstylestrong{Arguments}
\begin{quote}

\sphinxAtStartPar
\sphinxcode{\sphinxupquote{idx}}: index of the desired exponential cone constraint.
\end{quote}

\sphinxAtStartPar
\sphinxstylestrong{Return}
\begin{quote}

\sphinxAtStartPar
the desired exponential cone constraint object.
\end{quote}
\end{quote}

\subsubsection{Model::GetExpConeBuilders()}
\label{\detokenize{cppapi/model:model-getexpconebuilders}}\begin{quote}

\sphinxAtStartPar
Get builders of all exponential cone constraints in model.

\sphinxAtStartPar
\sphinxstylestrong{Synopsis}
\begin{quote}

\sphinxAtStartPar
\sphinxcode{\sphinxupquote{ExpConeBuilderArray GetExpConeBuilders()}}
\end{quote}

\sphinxAtStartPar
\sphinxstylestrong{Return}
\begin{quote}

\sphinxAtStartPar
array object of all exponential cone constraint builders.
\end{quote}
\end{quote}

\subsubsection{Model::GetExpConeBuilders()}
\label{\detokenize{cppapi/model:id65}}\begin{quote}

\sphinxAtStartPar
Get builders of given exponential cone constraints in model.

\sphinxAtStartPar
\sphinxstylestrong{Synopsis}
\begin{quote}

\sphinxAtStartPar
\sphinxcode{\sphinxupquote{ExpConeBuilderArray GetExpConeBuilders(const ExpConeArray \&cones)}}
\end{quote}

\sphinxAtStartPar
\sphinxstylestrong{Arguments}
\begin{quote}

\sphinxAtStartPar
\sphinxcode{\sphinxupquote{cones}}: array of exponential cone constraints.
\end{quote}

\sphinxAtStartPar
\sphinxstylestrong{Return}
\begin{quote}

\sphinxAtStartPar
array object of desired exponential cone constraint builders.
\end{quote}
\end{quote}

\subsubsection{Model::GetExpCones()}
\label{\detokenize{cppapi/model:model-getexpcones}}\begin{quote}

\sphinxAtStartPar
Get all exponential cone constraints in model.

\sphinxAtStartPar
\sphinxstylestrong{Synopsis}
\begin{quote}

\sphinxAtStartPar
\sphinxcode{\sphinxupquote{ExpConeArray GetExpCones()}}
\end{quote}

\sphinxAtStartPar
\sphinxstylestrong{Return}
\begin{quote}

\sphinxAtStartPar
array object of exponential cone constraints.
\end{quote}
\end{quote}

\subsubsection{Model::GetGenConstr()}
\label{\detokenize{cppapi/model:model-getgenconstr}}\begin{quote}

\sphinxAtStartPar
Get a general constraint of given index in model.

\sphinxAtStartPar
\sphinxstylestrong{Synopsis}
\begin{quote}

\sphinxAtStartPar
\sphinxcode{\sphinxupquote{GenConstr GetGenConstr(int idx)}}
\end{quote}

\sphinxAtStartPar
\sphinxstylestrong{Arguments}
\begin{quote}

\sphinxAtStartPar
\sphinxcode{\sphinxupquote{idx}}: index of the desired general constraint.
\end{quote}

\sphinxAtStartPar
\sphinxstylestrong{Return}
\begin{quote}

\sphinxAtStartPar
the desired general constraint object.
\end{quote}
\end{quote}

\subsubsection{Model::GetGenConstrByName()}
\label{\detokenize{cppapi/model:model-getgenconstrbyname}}\begin{quote}

\sphinxAtStartPar
Get a general constraint of given name in model.

\sphinxAtStartPar
\sphinxstylestrong{Synopsis}
\begin{quote}

\sphinxAtStartPar
\sphinxcode{\sphinxupquote{GenConstr GetGenConstrByName(const char *szName)}}
\end{quote}

\sphinxAtStartPar
\sphinxstylestrong{Arguments}
\begin{quote}

\sphinxAtStartPar
\sphinxcode{\sphinxupquote{szName}}: name of the desired general constraint.
\end{quote}

\sphinxAtStartPar
\sphinxstylestrong{Return}
\begin{quote}

\sphinxAtStartPar
the desired general constraint object.
\end{quote}
\end{quote}

\subsubsection{Model::GetGenConstrIndicator()}
\label{\detokenize{cppapi/model:model-getgenconstrindicator}}\begin{quote}

\sphinxAtStartPar
Get builder of given general constraint of type indicator.

\sphinxAtStartPar
\sphinxstylestrong{Synopsis}
\begin{quote}

\sphinxAtStartPar
\sphinxcode{\sphinxupquote{GenConstrBuilder GetGenConstrIndicator(const GenConstr \&indicator)}}
\end{quote}

\sphinxAtStartPar
\sphinxstylestrong{Arguments}
\begin{quote}

\sphinxAtStartPar
\sphinxcode{\sphinxupquote{indicator}}: a general constraint of type indicator.
\end{quote}

\sphinxAtStartPar
\sphinxstylestrong{Return}
\begin{quote}

\sphinxAtStartPar
builder object of general constraint of type indicator.
\end{quote}
\end{quote}

\subsubsection{Model::GetGenConstrIndicators()}
\label{\detokenize{cppapi/model:model-getgenconstrindicators}}\begin{quote}

\sphinxAtStartPar
Get builders of all general constraints in model.

\sphinxAtStartPar
\sphinxstylestrong{Synopsis}
\begin{quote}

\sphinxAtStartPar
\sphinxcode{\sphinxupquote{GenConstrBuilderArray GetGenConstrIndicators()}}
\end{quote}

\sphinxAtStartPar
\sphinxstylestrong{Return}
\begin{quote}

\sphinxAtStartPar
array object of general constraint builders.
\end{quote}
\end{quote}

\subsubsection{Model::GetGenConstrs()}
\label{\detokenize{cppapi/model:model-getgenconstrs}}\begin{quote}

\sphinxAtStartPar
Get all general constraints in model.

\sphinxAtStartPar
\sphinxstylestrong{Synopsis}
\begin{quote}

\sphinxAtStartPar
\sphinxcode{\sphinxupquote{GenConstrArray GetGenConstrs()}}
\end{quote}

\sphinxAtStartPar
\sphinxstylestrong{Return}
\begin{quote}

\sphinxAtStartPar
array object of general constraints.
\end{quote}
\end{quote}

\subsubsection{Model::GetIndicatorIIS()}
\label{\detokenize{cppapi/model:model-getindicatoriis}}\begin{quote}

\sphinxAtStartPar
Get IIS status of indicator constraints.

\sphinxAtStartPar
\sphinxstylestrong{Synopsis}
\begin{quote}

\sphinxAtStartPar
\sphinxcode{\sphinxupquote{int GetIndicatorIIS(const GenConstrArray \&genconstrs, int *pIIS)}}
\end{quote}

\sphinxAtStartPar
\sphinxstylestrong{Arguments}
\begin{quote}

\sphinxAtStartPar
\sphinxcode{\sphinxupquote{genconstrs}}: Array of indicator constraints.

\sphinxAtStartPar
\sphinxcode{\sphinxupquote{pIIS}}: IIS status of indicator constraints.
\end{quote}

\sphinxAtStartPar
\sphinxstylestrong{Return}
\begin{quote}

\sphinxAtStartPar
Number of indicator constraints.
\end{quote}
\end{quote}

\subsubsection{Model::GetIntAttr()}
\label{\detokenize{cppapi/model:model-getintattr}}\begin{quote}

\sphinxAtStartPar
Get value of a COPT integer attribute.

\sphinxAtStartPar
\sphinxstylestrong{Synopsis}
\begin{quote}

\sphinxAtStartPar
\sphinxcode{\sphinxupquote{int GetIntAttr(const char *szAttr)}}
\end{quote}

\sphinxAtStartPar
\sphinxstylestrong{Arguments}
\begin{quote}

\sphinxAtStartPar
\sphinxcode{\sphinxupquote{szAttr}}: name of integer attribute.
\end{quote}

\sphinxAtStartPar
\sphinxstylestrong{Return}
\begin{quote}

\sphinxAtStartPar
value of integer attribute.
\end{quote}
\end{quote}

\subsubsection{Model::GetIntAttrN()}
\label{\detokenize{cppapi/model:model-getintattrn}}\begin{quote}

\sphinxAtStartPar
Get value of a integer attribute of a multi\sphinxhyphen{}objective function.

\sphinxAtStartPar
\sphinxstylestrong{Synopsis}
\begin{quote}

\sphinxAtStartPar
\sphinxcode{\sphinxupquote{int GetIntAttrN(int idx, const char *szAttr)}}
\end{quote}

\sphinxAtStartPar
\sphinxstylestrong{Arguments}
\begin{quote}

\sphinxAtStartPar
\sphinxcode{\sphinxupquote{idx}}: index of a multi\sphinxhyphen{}objective function.

\sphinxAtStartPar
\sphinxcode{\sphinxupquote{szAttr}}: name of integer attribute.
\end{quote}

\sphinxAtStartPar
\sphinxstylestrong{Return}
\begin{quote}

\sphinxAtStartPar
value of integer attribute.
\end{quote}
\end{quote}

\subsubsection{Model::GetIntParam()}
\label{\detokenize{cppapi/model:model-getintparam}}\begin{quote}

\sphinxAtStartPar
Get value of a COPT integer parameter.

\sphinxAtStartPar
\sphinxstylestrong{Synopsis}
\begin{quote}

\sphinxAtStartPar
\sphinxcode{\sphinxupquote{int GetIntParam(const char *szParam)}}
\end{quote}

\sphinxAtStartPar
\sphinxstylestrong{Arguments}
\begin{quote}

\sphinxAtStartPar
\sphinxcode{\sphinxupquote{szParam}}: name of integer parameter.
\end{quote}

\sphinxAtStartPar
\sphinxstylestrong{Return}
\begin{quote}

\sphinxAtStartPar
value of integer parameter.
\end{quote}
\end{quote}

\subsubsection{Model::GetIntParamN()}
\label{\detokenize{cppapi/model:model-getintparamn}}\begin{quote}

\sphinxAtStartPar
Get value of an integer parameter of a multi\sphinxhyphen{}objective function.

\sphinxAtStartPar
\sphinxstylestrong{Synopsis}
\begin{quote}

\sphinxAtStartPar
\sphinxcode{\sphinxupquote{int GetIntParamN(int idx, const char *szParam)}}
\end{quote}

\sphinxAtStartPar
\sphinxstylestrong{Arguments}
\begin{quote}

\sphinxAtStartPar
\sphinxcode{\sphinxupquote{idx}}: index of a multi\sphinxhyphen{}objective function.

\sphinxAtStartPar
\sphinxcode{\sphinxupquote{szParam}}: name of integer parameter.
\end{quote}

\sphinxAtStartPar
\sphinxstylestrong{Return}
\begin{quote}

\sphinxAtStartPar
value of integer parameter.
\end{quote}
\end{quote}

\subsubsection{Model::GetLmiCoeff()}
\label{\detokenize{cppapi/model:model-getlmicoeff}}\begin{quote}

\sphinxAtStartPar
Get the symmetric matrix of variable in LMI constraint.

\sphinxAtStartPar
\sphinxstylestrong{Synopsis}
\begin{quote}

\sphinxAtStartPar
\sphinxcode{\sphinxupquote{SymMatrix GetLmiCoeff(const LmiConstraint \&constr, const Var \&var)}}
\end{quote}

\sphinxAtStartPar
\sphinxstylestrong{Arguments}
\begin{quote}

\sphinxAtStartPar
\sphinxcode{\sphinxupquote{constr}}: The desired LMI constraint.

\sphinxAtStartPar
\sphinxcode{\sphinxupquote{var}}: The desired variable.
\end{quote}

\sphinxAtStartPar
\sphinxstylestrong{Return}
\begin{quote}

\sphinxAtStartPar
The associated coefficient matrix.
\end{quote}
\end{quote}

\subsubsection{Model::GetLmiConstr()}
\label{\detokenize{cppapi/model:model-getlmiconstr}}\begin{quote}

\sphinxAtStartPar
Get LMI constraint of given index in model.

\sphinxAtStartPar
\sphinxstylestrong{Synopsis}
\begin{quote}

\sphinxAtStartPar
\sphinxcode{\sphinxupquote{LmiConstraint GetLmiConstr(int idx)}}
\end{quote}

\sphinxAtStartPar
\sphinxstylestrong{Arguments}
\begin{quote}

\sphinxAtStartPar
\sphinxcode{\sphinxupquote{idx}}: index of desired LMI constraint.
\end{quote}

\sphinxAtStartPar
\sphinxstylestrong{Return}
\begin{quote}

\sphinxAtStartPar
LMI constraint object.
\end{quote}
\end{quote}

\subsubsection{Model::GetLmiConstrByName()}
\label{\detokenize{cppapi/model:model-getlmiconstrbyname}}\begin{quote}

\sphinxAtStartPar
Get LMI constraint of given name in model.

\sphinxAtStartPar
\sphinxstylestrong{Synopsis}
\begin{quote}

\sphinxAtStartPar
\sphinxcode{\sphinxupquote{LmiConstraint GetLmiConstrByName(const char *szName)}}
\end{quote}

\sphinxAtStartPar
\sphinxstylestrong{Arguments}
\begin{quote}

\sphinxAtStartPar
\sphinxcode{\sphinxupquote{szName}}: name of desired LMI constraint.
\end{quote}

\sphinxAtStartPar
\sphinxstylestrong{Return}
\begin{quote}

\sphinxAtStartPar
LMI constraint object.
\end{quote}
\end{quote}

\subsubsection{Model::GetLmiConstrs()}
\label{\detokenize{cppapi/model:model-getlmiconstrs}}\begin{quote}

\sphinxAtStartPar
Get all LMI constraints in model.

\sphinxAtStartPar
\sphinxstylestrong{Synopsis}
\begin{quote}

\sphinxAtStartPar
\sphinxcode{\sphinxupquote{LmiConstrArray GetLmiConstrs()}}
\end{quote}

\sphinxAtStartPar
\sphinxstylestrong{Return}
\begin{quote}

\sphinxAtStartPar
array object of LMI constraints.
\end{quote}
\end{quote}

\subsubsection{Model::GetLmiRhs()}
\label{\detokenize{cppapi/model:model-getlmirhs}}\begin{quote}

\sphinxAtStartPar
Get the symmetric matrix of constant of LMI constraint.

\sphinxAtStartPar
\sphinxstylestrong{Synopsis}
\begin{quote}

\sphinxAtStartPar
\sphinxcode{\sphinxupquote{SymMatrix GetLmiRhs(const LmiConstraint \&constr)}}
\end{quote}

\sphinxAtStartPar
\sphinxstylestrong{Arguments}
\begin{quote}

\sphinxAtStartPar
\sphinxcode{\sphinxupquote{constr}}: The desired LMI constraint.
\end{quote}

\sphinxAtStartPar
\sphinxstylestrong{Return}
\begin{quote}

\sphinxAtStartPar
matrix of constant term.
\end{quote}
\end{quote}

\subsubsection{Model::GetLmiRow()}
\label{\detokenize{cppapi/model:model-getlmirow}}\begin{quote}

\sphinxAtStartPar
Get variables and associated symmetric matrices that participate in a LMI constraint.

\sphinxAtStartPar
\sphinxstylestrong{Synopsis}
\begin{quote}

\sphinxAtStartPar
\sphinxcode{\sphinxupquote{LmiExpr GetLmiRow(const LmiConstraint \&constr)}}
\end{quote}

\sphinxAtStartPar
\sphinxstylestrong{Arguments}
\begin{quote}

\sphinxAtStartPar
\sphinxcode{\sphinxupquote{constr}}: given LMI constraint object.
\end{quote}

\sphinxAtStartPar
\sphinxstylestrong{Return}
\begin{quote}

\sphinxAtStartPar
pointer to LMI expression object of LMI constraint.
\end{quote}
\end{quote}

\subsubsection{Model::GetLpSolution()}
\label{\detokenize{cppapi/model:model-getlpsolution}}\begin{quote}

\sphinxAtStartPar
Get LP solution.

\sphinxAtStartPar
\sphinxstylestrong{Synopsis}
\begin{quote}

\sphinxAtStartPar
\sphinxcode{\sphinxupquote{void GetLpSolution(}}
\begin{quote}

\sphinxAtStartPar
\sphinxcode{\sphinxupquote{double *pValue,}}

\sphinxAtStartPar
\sphinxcode{\sphinxupquote{double *pSlack,}}

\sphinxAtStartPar
\sphinxcode{\sphinxupquote{double *pRowDual,}}

\sphinxAtStartPar
\sphinxcode{\sphinxupquote{double *pRedCost)}}
\end{quote}
\end{quote}

\sphinxAtStartPar
\sphinxstylestrong{Arguments}
\begin{quote}

\sphinxAtStartPar
\sphinxcode{\sphinxupquote{pValue}}: optional, double pointer to solution values.

\sphinxAtStartPar
\sphinxcode{\sphinxupquote{pSlack}}: optional, double poitner to slack values.

\sphinxAtStartPar
\sphinxcode{\sphinxupquote{pRowDual}}: optional, double pointer to dual values.

\sphinxAtStartPar
\sphinxcode{\sphinxupquote{pRedCost}}: optional, double pointer to reduced costs.
\end{quote}
\end{quote}

\subsubsection{Model::GetName()}
\label{\detokenize{cppapi/model:model-getname}}\begin{quote}

\sphinxAtStartPar
Get name of the model.

\sphinxAtStartPar
\sphinxstylestrong{Synopsis}
\begin{quote}

\sphinxAtStartPar
\sphinxcode{\sphinxupquote{const char *GetName()}}
\end{quote}

\sphinxAtStartPar
\sphinxstylestrong{Return}
\begin{quote}

\sphinxAtStartPar
model name.
\end{quote}
\end{quote}

\subsubsection{Model::GetNlConstr()}
\label{\detokenize{cppapi/model:model-getnlconstr}}\begin{quote}

\sphinxAtStartPar
Get a nonlinear constraint of given index in model.

\sphinxAtStartPar
\sphinxstylestrong{Synopsis}
\begin{quote}

\sphinxAtStartPar
\sphinxcode{\sphinxupquote{NlConstraint GetNlConstr(int idx)}}
\end{quote}

\sphinxAtStartPar
\sphinxstylestrong{Arguments}
\begin{quote}

\sphinxAtStartPar
\sphinxcode{\sphinxupquote{idx}}: index of the desired nonlinear constraint.
\end{quote}

\sphinxAtStartPar
\sphinxstylestrong{Return}
\begin{quote}

\sphinxAtStartPar
the desired nonlinear constraint object.
\end{quote}
\end{quote}

\subsubsection{Model::GetNlConstrBuilder()}
\label{\detokenize{cppapi/model:model-getnlconstrbuilder}}\begin{quote}

\sphinxAtStartPar
Get builder of a nonlinear constraint in model, including nonlinear expression, sense and RHS.

\sphinxAtStartPar
\sphinxstylestrong{Synopsis}
\begin{quote}

\sphinxAtStartPar
\sphinxcode{\sphinxupquote{NlConstrBuilder GetNlConstrBuilder(const NlConstraint \&constr)}}
\end{quote}

\sphinxAtStartPar
\sphinxstylestrong{Arguments}
\begin{quote}

\sphinxAtStartPar
\sphinxcode{\sphinxupquote{constr}}: a nonlinear constraint object.
\end{quote}

\sphinxAtStartPar
\sphinxstylestrong{Return}
\begin{quote}

\sphinxAtStartPar
nonlinear constraint builder object.
\end{quote}
\end{quote}

\subsubsection{Model::GetNlConstrBuilders()}
\label{\detokenize{cppapi/model:model-getnlconstrbuilders}}\begin{quote}

\sphinxAtStartPar
Get builders of all nonlinear constraints in model.

\sphinxAtStartPar
\sphinxstylestrong{Synopsis}
\begin{quote}

\sphinxAtStartPar
\sphinxcode{\sphinxupquote{NlConstrBuilderArray GetNlConstrBuilders()}}
\end{quote}

\sphinxAtStartPar
\sphinxstylestrong{Return}
\begin{quote}

\sphinxAtStartPar
array object of nonlinear constraint builders.
\end{quote}
\end{quote}

\subsubsection{Model::GetNlConstrByName()}
\label{\detokenize{cppapi/model:model-getnlconstrbyname}}\begin{quote}

\sphinxAtStartPar
Get a nonlinear constraint of given name in model.

\sphinxAtStartPar
\sphinxstylestrong{Synopsis}
\begin{quote}

\sphinxAtStartPar
\sphinxcode{\sphinxupquote{NlConstraint GetNlConstrByName(const char *szName)}}
\end{quote}

\sphinxAtStartPar
\sphinxstylestrong{Arguments}
\begin{quote}

\sphinxAtStartPar
\sphinxcode{\sphinxupquote{szName}}: name of the desired constraint.
\end{quote}

\sphinxAtStartPar
\sphinxstylestrong{Return}
\begin{quote}

\sphinxAtStartPar
the desired nonlinear constraint object.
\end{quote}
\end{quote}

\subsubsection{Model::GetNlConstrs()}
\label{\detokenize{cppapi/model:model-getnlconstrs}}\begin{quote}

\sphinxAtStartPar
Get all nonlinear constraints in model.

\sphinxAtStartPar
\sphinxstylestrong{Synopsis}
\begin{quote}

\sphinxAtStartPar
\sphinxcode{\sphinxupquote{NlConstrArray GetNlConstrs()}}
\end{quote}

\sphinxAtStartPar
\sphinxstylestrong{Return}
\begin{quote}

\sphinxAtStartPar
array object of nonlinear constraints.
\end{quote}
\end{quote}

\subsubsection{Model::GetNlObjective()}
\label{\detokenize{cppapi/model:model-getnlobjective}}\begin{quote}

\sphinxAtStartPar
Get nonlinear objective of model.

\sphinxAtStartPar
\sphinxstylestrong{Synopsis}
\begin{quote}

\sphinxAtStartPar
\sphinxcode{\sphinxupquote{NlExpr GetNlObjective()}}
\end{quote}

\sphinxAtStartPar
\sphinxstylestrong{Return}
\begin{quote}

\sphinxAtStartPar
a nonlinear expression object.
\end{quote}
\end{quote}

\subsubsection{Model::GetNlRow()}
\label{\detokenize{cppapi/model:model-getnlrow}}\begin{quote}

\sphinxAtStartPar
Get nonlinear expression of a nonlinear constraint.

\sphinxAtStartPar
\sphinxstylestrong{Synopsis}
\begin{quote}

\sphinxAtStartPar
\sphinxcode{\sphinxupquote{NlExpr GetNlRow(const NlConstraint \&constr)}}
\end{quote}

\sphinxAtStartPar
\sphinxstylestrong{Arguments}
\begin{quote}

\sphinxAtStartPar
\sphinxcode{\sphinxupquote{constr}}: a nonlinear constraint object.
\end{quote}

\sphinxAtStartPar
\sphinxstylestrong{Return}
\begin{quote}

\sphinxAtStartPar
output object of nonlinear expression.
\end{quote}
\end{quote}

\subsubsection{Model::GetObjective()}
\label{\detokenize{cppapi/model:model-getobjective}}\begin{quote}

\sphinxAtStartPar
Get linear expression of objective for model.

\sphinxAtStartPar
\sphinxstylestrong{Synopsis}
\begin{quote}

\sphinxAtStartPar
\sphinxcode{\sphinxupquote{Expr GetObjective()}}
\end{quote}

\sphinxAtStartPar
\sphinxstylestrong{Return}
\begin{quote}

\sphinxAtStartPar
a linear expression object.
\end{quote}
\end{quote}

\subsubsection{Model::GetObjectiveN()}
\label{\detokenize{cppapi/model:model-getobjectiven}}\begin{quote}

\sphinxAtStartPar
Get linear expression of a multi\sphinxhyphen{}objective function in model.

\sphinxAtStartPar
\sphinxstylestrong{Synopsis}
\begin{quote}

\sphinxAtStartPar
\sphinxcode{\sphinxupquote{Expr GetObjectiveN(int idx)}}
\end{quote}

\sphinxAtStartPar
\sphinxstylestrong{Arguments}
\begin{quote}

\sphinxAtStartPar
\sphinxcode{\sphinxupquote{idx}}: index of a multi\sphinxhyphen{}objective function.
\end{quote}

\sphinxAtStartPar
\sphinxstylestrong{Return}
\begin{quote}

\sphinxAtStartPar
a linear expression object.
\end{quote}
\end{quote}

\subsubsection{Model::GetObjParamN()}
\label{\detokenize{cppapi/model:model-getobjparamn}}\begin{quote}

\sphinxAtStartPar
Get value of objective parameter of a multi\sphinxhyphen{}objective function.

\sphinxAtStartPar
\sphinxstylestrong{Synopsis}
\begin{quote}

\sphinxAtStartPar
\sphinxcode{\sphinxupquote{double GetObjParamN(int idx, const char *szParam)}}
\end{quote}

\sphinxAtStartPar
\sphinxstylestrong{Arguments}
\begin{quote}

\sphinxAtStartPar
\sphinxcode{\sphinxupquote{idx}}: index of a multi\sphinxhyphen{}objective function.

\sphinxAtStartPar
\sphinxcode{\sphinxupquote{szParam}}: name of objective parameter, including priority, weight, abstol and reltol.
\end{quote}

\sphinxAtStartPar
\sphinxstylestrong{Return}
\begin{quote}

\sphinxAtStartPar
value of objective parameter.
\end{quote}
\end{quote}

\subsubsection{Model::GetParamAttrType()}
\label{\detokenize{cppapi/model:model-getparamattrtype}}\begin{quote}

\sphinxAtStartPar
Get type of a COPT parameter or attribute.

\sphinxAtStartPar
\sphinxstylestrong{Synopsis}
\begin{quote}

\sphinxAtStartPar
\sphinxcode{\sphinxupquote{int GetParamAttrType(const char *szName)}}
\end{quote}

\sphinxAtStartPar
\sphinxstylestrong{Arguments}
\begin{quote}

\sphinxAtStartPar
\sphinxcode{\sphinxupquote{szName}}: name of parameter or attribute.
\end{quote}

\sphinxAtStartPar
\sphinxstylestrong{Return}
\begin{quote}

\sphinxAtStartPar
type of parameter or attribute.
\end{quote}
\end{quote}

\subsubsection{Model::GetParamInfo()}
\label{\detokenize{cppapi/model:model-getparaminfo}}\begin{quote}

\sphinxAtStartPar
Get current, default, minimum, maximum of COPT integer parameter.

\sphinxAtStartPar
\sphinxstylestrong{Synopsis}
\begin{quote}

\sphinxAtStartPar
\sphinxcode{\sphinxupquote{void GetParamInfo(}}
\begin{quote}

\sphinxAtStartPar
\sphinxcode{\sphinxupquote{const char *szParam,}}

\sphinxAtStartPar
\sphinxcode{\sphinxupquote{int *pnCur,}}

\sphinxAtStartPar
\sphinxcode{\sphinxupquote{int *pnDef,}}

\sphinxAtStartPar
\sphinxcode{\sphinxupquote{int *pnMin,}}

\sphinxAtStartPar
\sphinxcode{\sphinxupquote{int *pnMax)}}
\end{quote}
\end{quote}

\sphinxAtStartPar
\sphinxstylestrong{Arguments}
\begin{quote}

\sphinxAtStartPar
\sphinxcode{\sphinxupquote{szParam}}: name of integer parameter.

\sphinxAtStartPar
\sphinxcode{\sphinxupquote{pnCur}}: out, current value of integer parameter.

\sphinxAtStartPar
\sphinxcode{\sphinxupquote{pnDef}}: out, default value of integer parameter.

\sphinxAtStartPar
\sphinxcode{\sphinxupquote{pnMin}}: out, minimum value of integer parameter.

\sphinxAtStartPar
\sphinxcode{\sphinxupquote{pnMax}}: out, maximum value of integer parameter.
\end{quote}
\end{quote}

\subsubsection{Model::GetParamInfo()}
\label{\detokenize{cppapi/model:id66}}\begin{quote}

\sphinxAtStartPar
Get current, default, minimum, maximum of COPT double parameter.

\sphinxAtStartPar
\sphinxstylestrong{Synopsis}
\begin{quote}

\sphinxAtStartPar
\sphinxcode{\sphinxupquote{void GetParamInfo(}}
\begin{quote}

\sphinxAtStartPar
\sphinxcode{\sphinxupquote{const char *szParam,}}

\sphinxAtStartPar
\sphinxcode{\sphinxupquote{double *pdCur,}}

\sphinxAtStartPar
\sphinxcode{\sphinxupquote{double *pdDef,}}

\sphinxAtStartPar
\sphinxcode{\sphinxupquote{double *pdMin,}}

\sphinxAtStartPar
\sphinxcode{\sphinxupquote{double *pdMax)}}
\end{quote}
\end{quote}

\sphinxAtStartPar
\sphinxstylestrong{Arguments}
\begin{quote}

\sphinxAtStartPar
\sphinxcode{\sphinxupquote{szParam}}: name of double parameter.

\sphinxAtStartPar
\sphinxcode{\sphinxupquote{pdCur}}: out, current value of double parameter.

\sphinxAtStartPar
\sphinxcode{\sphinxupquote{pdDef}}: out, default value of double parameter.

\sphinxAtStartPar
\sphinxcode{\sphinxupquote{pdMin}}: out, minimum value of double parameter.

\sphinxAtStartPar
\sphinxcode{\sphinxupquote{pdMax}}: out, maximum value of double parameter.
\end{quote}
\end{quote}

\subsubsection{Model::GetPoolObjVal()}
\label{\detokenize{cppapi/model:model-getpoolobjval}}\begin{quote}

\sphinxAtStartPar
Get the iSol\sphinxhyphen{}th objective value in solution pool.

\sphinxAtStartPar
\sphinxstylestrong{Synopsis}
\begin{quote}

\sphinxAtStartPar
\sphinxcode{\sphinxupquote{double GetPoolObjVal(int iSol)}}
\end{quote}

\sphinxAtStartPar
\sphinxstylestrong{Arguments}
\begin{quote}

\sphinxAtStartPar
\sphinxcode{\sphinxupquote{iSol}}: Index of solution.
\end{quote}

\sphinxAtStartPar
\sphinxstylestrong{Return}
\begin{quote}

\sphinxAtStartPar
The requested objective value.
\end{quote}
\end{quote}

\subsubsection{Model::GetPoolObjValN()}
\label{\detokenize{cppapi/model:model-getpoolobjvaln}}\begin{quote}

\sphinxAtStartPar
Get the objective value of required multi\sphinxhyphen{}objective function in solution pool.

\sphinxAtStartPar
\sphinxstylestrong{Synopsis}
\begin{quote}

\sphinxAtStartPar
\sphinxcode{\sphinxupquote{double GetPoolObjValN(int idx, int iSol)}}
\end{quote}

\sphinxAtStartPar
\sphinxstylestrong{Arguments}
\begin{quote}

\sphinxAtStartPar
\sphinxcode{\sphinxupquote{idx}}: index of a multi\sphinxhyphen{}objective function.

\sphinxAtStartPar
\sphinxcode{\sphinxupquote{iSol}}: index of solution.
\end{quote}

\sphinxAtStartPar
\sphinxstylestrong{Return}
\begin{quote}

\sphinxAtStartPar
value of required multi\sphinxhyphen{}objective function.
\end{quote}
\end{quote}

\subsubsection{Model::GetPoolSolution()}
\label{\detokenize{cppapi/model:model-getpoolsolution}}\begin{quote}

\sphinxAtStartPar
Get the iSol\sphinxhyphen{}th solution in solution pool.

\sphinxAtStartPar
\sphinxstylestrong{Synopsis}
\begin{quote}

\sphinxAtStartPar
\sphinxcode{\sphinxupquote{int GetPoolSolution(}}
\begin{quote}

\sphinxAtStartPar
\sphinxcode{\sphinxupquote{int iSol,}}

\sphinxAtStartPar
\sphinxcode{\sphinxupquote{const VarArray \&vars,}}

\sphinxAtStartPar
\sphinxcode{\sphinxupquote{double *pColVals)}}
\end{quote}
\end{quote}

\sphinxAtStartPar
\sphinxstylestrong{Arguments}
\begin{quote}

\sphinxAtStartPar
\sphinxcode{\sphinxupquote{iSol}}: Index of solution.

\sphinxAtStartPar
\sphinxcode{\sphinxupquote{vars}}: The requested variables.

\sphinxAtStartPar
\sphinxcode{\sphinxupquote{pColVals}}: Pointer to the requested solutions.
\end{quote}

\sphinxAtStartPar
\sphinxstylestrong{Return}
\begin{quote}

\sphinxAtStartPar
The length of requested solution array.
\end{quote}
\end{quote}

\subsubsection{Model::GetPsdCoeff()}
\label{\detokenize{cppapi/model:model-getpsdcoeff}}\begin{quote}

\sphinxAtStartPar
Get the symmetric matrix of PSD variable in a PSD constraint.

\sphinxAtStartPar
\sphinxstylestrong{Synopsis}
\begin{quote}

\sphinxAtStartPar
\sphinxcode{\sphinxupquote{SymMatrix GetPsdCoeff(const PsdConstraint \&constr, const PsdVar \&var)}}
\end{quote}

\sphinxAtStartPar
\sphinxstylestrong{Arguments}
\begin{quote}

\sphinxAtStartPar
\sphinxcode{\sphinxupquote{constr}}: The desired PSD constraint.

\sphinxAtStartPar
\sphinxcode{\sphinxupquote{var}}: The desired PSD variable.
\end{quote}

\sphinxAtStartPar
\sphinxstylestrong{Return}
\begin{quote}

\sphinxAtStartPar
The associated coefficient matrix.
\end{quote}
\end{quote}

\subsubsection{Model::GetPsdConstr()}
\label{\detokenize{cppapi/model:model-getpsdconstr}}\begin{quote}

\sphinxAtStartPar
Get PSD constraint of given index in model.

\sphinxAtStartPar
\sphinxstylestrong{Synopsis}
\begin{quote}

\sphinxAtStartPar
\sphinxcode{\sphinxupquote{PsdConstraint GetPsdConstr(int idx)}}
\end{quote}

\sphinxAtStartPar
\sphinxstylestrong{Arguments}
\begin{quote}

\sphinxAtStartPar
\sphinxcode{\sphinxupquote{idx}}: index of desired PSD constraint.
\end{quote}

\sphinxAtStartPar
\sphinxstylestrong{Return}
\begin{quote}

\sphinxAtStartPar
PSD constraint object.
\end{quote}
\end{quote}

\subsubsection{Model::GetPsdConstrBuilder()}
\label{\detokenize{cppapi/model:model-getpsdconstrbuilder}}\begin{quote}

\sphinxAtStartPar
Get builder of a PSD constraint in model, including PSD variables, sense and associated symmetric matrix.

\sphinxAtStartPar
\sphinxstylestrong{Synopsis}
\begin{quote}

\sphinxAtStartPar
\sphinxcode{\sphinxupquote{PsdConstrBuilder GetPsdConstrBuilder(const PsdConstraint \&constr)}}
\end{quote}

\sphinxAtStartPar
\sphinxstylestrong{Arguments}
\begin{quote}

\sphinxAtStartPar
\sphinxcode{\sphinxupquote{constr}}: PSD constraint object.
\end{quote}

\sphinxAtStartPar
\sphinxstylestrong{Return}
\begin{quote}

\sphinxAtStartPar
pointer to PSD constraint builder object.
\end{quote}
\end{quote}

\subsubsection{Model::GetPsdConstrBuilders()}
\label{\detokenize{cppapi/model:model-getpsdconstrbuilders}}\begin{quote}

\sphinxAtStartPar
Get builders of all PSD constraints in model.

\sphinxAtStartPar
\sphinxstylestrong{Synopsis}
\begin{quote}

\sphinxAtStartPar
\sphinxcode{\sphinxupquote{PsdConstrBuilderArray GetPsdConstrBuilders()}}
\end{quote}

\sphinxAtStartPar
\sphinxstylestrong{Return}
\begin{quote}

\sphinxAtStartPar
pointer to array object of PSD constraint builders.
\end{quote}
\end{quote}

\subsubsection{Model::GetPsdConstrByName()}
\label{\detokenize{cppapi/model:model-getpsdconstrbyname}}\begin{quote}

\sphinxAtStartPar
Get PSD constraint of given name in model.

\sphinxAtStartPar
\sphinxstylestrong{Synopsis}
\begin{quote}

\sphinxAtStartPar
\sphinxcode{\sphinxupquote{PsdConstraint GetPsdConstrByName(const char *szName)}}
\end{quote}

\sphinxAtStartPar
\sphinxstylestrong{Arguments}
\begin{quote}

\sphinxAtStartPar
\sphinxcode{\sphinxupquote{szName}}: name of desired PSD constraint.
\end{quote}

\sphinxAtStartPar
\sphinxstylestrong{Return}
\begin{quote}

\sphinxAtStartPar
PSD constraint object.
\end{quote}
\end{quote}

\subsubsection{Model::GetPsdConstrs()}
\label{\detokenize{cppapi/model:model-getpsdconstrs}}\begin{quote}

\sphinxAtStartPar
Get all PSD constraints in model.

\sphinxAtStartPar
\sphinxstylestrong{Synopsis}
\begin{quote}

\sphinxAtStartPar
\sphinxcode{\sphinxupquote{PsdConstrArray GetPsdConstrs()}}
\end{quote}

\sphinxAtStartPar
\sphinxstylestrong{Return}
\begin{quote}

\sphinxAtStartPar
pointer to array object of PSD constraints.
\end{quote}
\end{quote}

\subsubsection{Model::GetPsdObjective()}
\label{\detokenize{cppapi/model:model-getpsdobjective}}\begin{quote}

\sphinxAtStartPar
Get PSD objective of model.

\sphinxAtStartPar
\sphinxstylestrong{Synopsis}
\begin{quote}

\sphinxAtStartPar
\sphinxcode{\sphinxupquote{PsdExpr GetPsdObjective()}}
\end{quote}

\sphinxAtStartPar
\sphinxstylestrong{Return}
\begin{quote}

\sphinxAtStartPar
a PSD expression object.
\end{quote}
\end{quote}

\subsubsection{Model::GetPsdRow()}
\label{\detokenize{cppapi/model:model-getpsdrow}}\begin{quote}

\sphinxAtStartPar
Get PSD variables and associated symmetric matrices that participate in a PSD constraint.

\sphinxAtStartPar
\sphinxstylestrong{Synopsis}
\begin{quote}

\sphinxAtStartPar
\sphinxcode{\sphinxupquote{PsdExpr GetPsdRow(const PsdConstraint \&constr)}}
\end{quote}

\sphinxAtStartPar
\sphinxstylestrong{Arguments}
\begin{quote}

\sphinxAtStartPar
\sphinxcode{\sphinxupquote{constr}}: PSD constraint object.
\end{quote}

\sphinxAtStartPar
\sphinxstylestrong{Return}
\begin{quote}

\sphinxAtStartPar
pointer to PSD expression object of the PSD constraint.
\end{quote}
\end{quote}

\subsubsection{Model::GetPsdRow()}
\label{\detokenize{cppapi/model:id67}}\begin{quote}

\sphinxAtStartPar
Get PSD variables, associated symmetric matrix, LB/UB that participate in a PSD constraint.

\sphinxAtStartPar
\sphinxstylestrong{Synopsis}
\begin{quote}

\sphinxAtStartPar
\sphinxcode{\sphinxupquote{PsdExpr GetPsdRow(}}
\begin{quote}

\sphinxAtStartPar
\sphinxcode{\sphinxupquote{const PsdConstraint \&constr,}}

\sphinxAtStartPar
\sphinxcode{\sphinxupquote{double *pLower,}}

\sphinxAtStartPar
\sphinxcode{\sphinxupquote{double *pUpper)}}
\end{quote}
\end{quote}

\sphinxAtStartPar
\sphinxstylestrong{Arguments}
\begin{quote}

\sphinxAtStartPar
\sphinxcode{\sphinxupquote{constr}}: a PSD constraint object.

\sphinxAtStartPar
\sphinxcode{\sphinxupquote{pLower}}: pointer to double value of lower bound.

\sphinxAtStartPar
\sphinxcode{\sphinxupquote{pUpper}}: pointer to double value of upper bound.
\end{quote}

\sphinxAtStartPar
\sphinxstylestrong{Return}
\begin{quote}

\sphinxAtStartPar
pointer to PSD expression object of the PSD constraint.
\end{quote}
\end{quote}

\subsubsection{Model::GetPsdVar()}
\label{\detokenize{cppapi/model:model-getpsdvar}}\begin{quote}

\sphinxAtStartPar
Get a PSD variable of given index in model.

\sphinxAtStartPar
\sphinxstylestrong{Synopsis}
\begin{quote}

\sphinxAtStartPar
\sphinxcode{\sphinxupquote{PsdVar GetPsdVar(int idx)}}
\end{quote}

\sphinxAtStartPar
\sphinxstylestrong{Arguments}
\begin{quote}

\sphinxAtStartPar
\sphinxcode{\sphinxupquote{idx}}: index of the desired PSD variable.
\end{quote}

\sphinxAtStartPar
\sphinxstylestrong{Return}
\begin{quote}

\sphinxAtStartPar
the desired PSD variable object.
\end{quote}
\end{quote}

\subsubsection{Model::GetPsdVarByName()}
\label{\detokenize{cppapi/model:model-getpsdvarbyname}}\begin{quote}

\sphinxAtStartPar
Get a PSD variable of given name in model.

\sphinxAtStartPar
\sphinxstylestrong{Synopsis}
\begin{quote}

\sphinxAtStartPar
\sphinxcode{\sphinxupquote{PsdVar GetPsdVarByName(const char *szName)}}
\end{quote}

\sphinxAtStartPar
\sphinxstylestrong{Arguments}
\begin{quote}

\sphinxAtStartPar
\sphinxcode{\sphinxupquote{szName}}: name of the desired PSD variable.
\end{quote}

\sphinxAtStartPar
\sphinxstylestrong{Return}
\begin{quote}

\sphinxAtStartPar
the desired PSD variable object.
\end{quote}
\end{quote}

\subsubsection{Model::GetPsdVars()}
\label{\detokenize{cppapi/model:model-getpsdvars}}\begin{quote}

\sphinxAtStartPar
Get all PSD variables in model.

\sphinxAtStartPar
\sphinxstylestrong{Synopsis}
\begin{quote}

\sphinxAtStartPar
\sphinxcode{\sphinxupquote{PsdVarArray GetPsdVars()}}
\end{quote}

\sphinxAtStartPar
\sphinxstylestrong{Return}
\begin{quote}

\sphinxAtStartPar
array object of PSD variables.
\end{quote}
\end{quote}

\subsubsection{Model::GetQConstr()}
\label{\detokenize{cppapi/model:model-getqconstr}}\begin{quote}

\sphinxAtStartPar
Get a quadratic constraint of given index in model.

\sphinxAtStartPar
\sphinxstylestrong{Synopsis}
\begin{quote}

\sphinxAtStartPar
\sphinxcode{\sphinxupquote{QConstraint GetQConstr(int idx)}}
\end{quote}

\sphinxAtStartPar
\sphinxstylestrong{Arguments}
\begin{quote}

\sphinxAtStartPar
\sphinxcode{\sphinxupquote{idx}}: index of the desired quadratic constraint.
\end{quote}

\sphinxAtStartPar
\sphinxstylestrong{Return}
\begin{quote}

\sphinxAtStartPar
the desired quadratic constraint object.
\end{quote}
\end{quote}

\subsubsection{Model::GetQConstrBuilder()}
\label{\detokenize{cppapi/model:model-getqconstrbuilder}}\begin{quote}

\sphinxAtStartPar
Get builder of a constraint in model, including variables and associated coefficients, sense and RHS.

\sphinxAtStartPar
\sphinxstylestrong{Synopsis}
\begin{quote}

\sphinxAtStartPar
\sphinxcode{\sphinxupquote{QConstrBuilder GetQConstrBuilder(const QConstraint \&constr)}}
\end{quote}

\sphinxAtStartPar
\sphinxstylestrong{Arguments}
\begin{quote}

\sphinxAtStartPar
\sphinxcode{\sphinxupquote{constr}}: a constraint object.
\end{quote}

\sphinxAtStartPar
\sphinxstylestrong{Return}
\begin{quote}

\sphinxAtStartPar
constraint builder object.
\end{quote}
\end{quote}

\subsubsection{Model::GetQConstrBuilders()}
\label{\detokenize{cppapi/model:model-getqconstrbuilders}}\begin{quote}

\sphinxAtStartPar
Get builders of all constraints in model.

\sphinxAtStartPar
\sphinxstylestrong{Synopsis}
\begin{quote}

\sphinxAtStartPar
\sphinxcode{\sphinxupquote{QConstrBuilderArray GetQConstrBuilders()}}
\end{quote}

\sphinxAtStartPar
\sphinxstylestrong{Return}
\begin{quote}

\sphinxAtStartPar
array object of constraint builders.
\end{quote}
\end{quote}

\subsubsection{Model::GetQConstrByName()}
\label{\detokenize{cppapi/model:model-getqconstrbyname}}\begin{quote}

\sphinxAtStartPar
Get a quadratic constraint of given name in model.

\sphinxAtStartPar
\sphinxstylestrong{Synopsis}
\begin{quote}

\sphinxAtStartPar
\sphinxcode{\sphinxupquote{QConstraint GetQConstrByName(const char *szName)}}
\end{quote}

\sphinxAtStartPar
\sphinxstylestrong{Arguments}
\begin{quote}

\sphinxAtStartPar
\sphinxcode{\sphinxupquote{szName}}: name of the desired constraint.
\end{quote}

\sphinxAtStartPar
\sphinxstylestrong{Return}
\begin{quote}

\sphinxAtStartPar
the desired quadratic constraint object.
\end{quote}
\end{quote}

\subsubsection{Model::GetQConstrs()}
\label{\detokenize{cppapi/model:model-getqconstrs}}\begin{quote}

\sphinxAtStartPar
Get all quadratic constraints in model.

\sphinxAtStartPar
\sphinxstylestrong{Synopsis}
\begin{quote}

\sphinxAtStartPar
\sphinxcode{\sphinxupquote{QConstrArray GetQConstrs()}}
\end{quote}

\sphinxAtStartPar
\sphinxstylestrong{Return}
\begin{quote}

\sphinxAtStartPar
array object of quadratic constraints.
\end{quote}
\end{quote}

\subsubsection{Model::GetQuadObjective()}
\label{\detokenize{cppapi/model:model-getquadobjective}}\begin{quote}

\sphinxAtStartPar
Get quadratic objective of model.

\sphinxAtStartPar
\sphinxstylestrong{Synopsis}
\begin{quote}

\sphinxAtStartPar
\sphinxcode{\sphinxupquote{QuadExpr GetQuadObjective()}}
\end{quote}

\sphinxAtStartPar
\sphinxstylestrong{Return}
\begin{quote}

\sphinxAtStartPar
a quadratic expression object.
\end{quote}
\end{quote}

\subsubsection{Model::GetQuadRow()}
\label{\detokenize{cppapi/model:model-getquadrow}}\begin{quote}

\sphinxAtStartPar
Get two variables and associated coefficients that participate in a quadratic constraint.

\sphinxAtStartPar
\sphinxstylestrong{Synopsis}
\begin{quote}

\sphinxAtStartPar
\sphinxcode{\sphinxupquote{QuadExpr GetQuadRow(const QConstraint \&constr)}}
\end{quote}

\sphinxAtStartPar
\sphinxstylestrong{Arguments}
\begin{quote}

\sphinxAtStartPar
\sphinxcode{\sphinxupquote{constr}}: a quadratic constraint object.
\end{quote}

\sphinxAtStartPar
\sphinxstylestrong{Return}
\begin{quote}

\sphinxAtStartPar
quadratic expression object of the constraint.
\end{quote}
\end{quote}

\subsubsection{Model::GetQuadRow()}
\label{\detokenize{cppapi/model:id68}}\begin{quote}

\sphinxAtStartPar
Get two variables and associated coefficients that participate in a quadratic constraint.

\sphinxAtStartPar
\sphinxstylestrong{Synopsis}
\begin{quote}

\sphinxAtStartPar
\sphinxcode{\sphinxupquote{QuadExpr GetQuadRow(}}
\begin{quote}

\sphinxAtStartPar
\sphinxcode{\sphinxupquote{const QConstraint \&constr,}}

\sphinxAtStartPar
\sphinxcode{\sphinxupquote{char *pSense,}}

\sphinxAtStartPar
\sphinxcode{\sphinxupquote{double *pBound)}}
\end{quote}
\end{quote}

\sphinxAtStartPar
\sphinxstylestrong{Arguments}
\begin{quote}

\sphinxAtStartPar
\sphinxcode{\sphinxupquote{constr}}: a quadratic constraint object.

\sphinxAtStartPar
\sphinxcode{\sphinxupquote{pSense}}: sense of quadratic constraint.

\sphinxAtStartPar
\sphinxcode{\sphinxupquote{pBound}}: right hand side of quadratic constraint.
\end{quote}

\sphinxAtStartPar
\sphinxstylestrong{Return}
\begin{quote}

\sphinxAtStartPar
quadratic expression object of the constraint.
\end{quote}
\end{quote}

\subsubsection{Model::GetRow()}
\label{\detokenize{cppapi/model:model-getrow}}\begin{quote}

\sphinxAtStartPar
Get variables that participate in a constraint, and the associated coefficients.

\sphinxAtStartPar
\sphinxstylestrong{Synopsis}
\begin{quote}

\sphinxAtStartPar
\sphinxcode{\sphinxupquote{Expr GetRow(const Constraint \&constr)}}
\end{quote}

\sphinxAtStartPar
\sphinxstylestrong{Arguments}
\begin{quote}

\sphinxAtStartPar
\sphinxcode{\sphinxupquote{constr}}: a constraint object.
\end{quote}

\sphinxAtStartPar
\sphinxstylestrong{Return}
\begin{quote}

\sphinxAtStartPar
expression object of the constraint.
\end{quote}
\end{quote}

\subsubsection{Model::GetRowBasis()}
\label{\detokenize{cppapi/model:model-getrowbasis}}\begin{quote}

\sphinxAtStartPar
Get status of row basis.

\sphinxAtStartPar
\sphinxstylestrong{Synopsis}
\begin{quote}

\sphinxAtStartPar
\sphinxcode{\sphinxupquote{int GetRowBasis(int *pBasis)}}
\end{quote}

\sphinxAtStartPar
\sphinxstylestrong{Arguments}
\begin{quote}

\sphinxAtStartPar
\sphinxcode{\sphinxupquote{pBasis}}: integer pointer to basis status.
\end{quote}

\sphinxAtStartPar
\sphinxstylestrong{Return}
\begin{quote}

\sphinxAtStartPar
number of rows.
\end{quote}
\end{quote}

\subsubsection{Model::GetSolution()}
\label{\detokenize{cppapi/model:model-getsolution}}\begin{quote}

\sphinxAtStartPar
Get MIP solution.

\sphinxAtStartPar
\sphinxstylestrong{Synopsis}
\begin{quote}

\sphinxAtStartPar
\sphinxcode{\sphinxupquote{void GetSolution(double *pValue)}}
\end{quote}

\sphinxAtStartPar
\sphinxstylestrong{Arguments}
\begin{quote}

\sphinxAtStartPar
\sphinxcode{\sphinxupquote{pValue}}: double pointer to solution values.
\end{quote}
\end{quote}

\subsubsection{Model::GetSos()}
\label{\detokenize{cppapi/model:model-getsos}}\begin{quote}

\sphinxAtStartPar
Get a SOS constraint of given index in model.

\sphinxAtStartPar
\sphinxstylestrong{Synopsis}
\begin{quote}

\sphinxAtStartPar
\sphinxcode{\sphinxupquote{Sos GetSos(int idx)}}
\end{quote}

\sphinxAtStartPar
\sphinxstylestrong{Arguments}
\begin{quote}

\sphinxAtStartPar
\sphinxcode{\sphinxupquote{idx}}: index of the desired SOS constraint.
\end{quote}

\sphinxAtStartPar
\sphinxstylestrong{Return}
\begin{quote}

\sphinxAtStartPar
the desired SOS constraint object.
\end{quote}
\end{quote}

\subsubsection{Model::GetSosBuilders()}
\label{\detokenize{cppapi/model:model-getsosbuilders}}\begin{quote}

\sphinxAtStartPar
Get builders of all SOS constraints in model.

\sphinxAtStartPar
\sphinxstylestrong{Synopsis}
\begin{quote}

\sphinxAtStartPar
\sphinxcode{\sphinxupquote{SosBuilderArray GetSosBuilders()}}
\end{quote}

\sphinxAtStartPar
\sphinxstylestrong{Return}
\begin{quote}

\sphinxAtStartPar
array object of SOS constraint builders.
\end{quote}
\end{quote}

\subsubsection{Model::GetSosBuilders()}
\label{\detokenize{cppapi/model:id69}}\begin{quote}

\sphinxAtStartPar
Get builders of given SOS constraints in model.

\sphinxAtStartPar
\sphinxstylestrong{Synopsis}
\begin{quote}

\sphinxAtStartPar
\sphinxcode{\sphinxupquote{SosBuilderArray GetSosBuilders(const SosArray \&soss)}}
\end{quote}

\sphinxAtStartPar
\sphinxstylestrong{Arguments}
\begin{quote}

\sphinxAtStartPar
\sphinxcode{\sphinxupquote{soss}}: array of SOS constraints.
\end{quote}

\sphinxAtStartPar
\sphinxstylestrong{Return}
\begin{quote}

\sphinxAtStartPar
array object of desired SOS constraint builders.
\end{quote}
\end{quote}

\subsubsection{Model::GetSOSIIS()}
\label{\detokenize{cppapi/model:model-getsosiis}}\begin{quote}

\sphinxAtStartPar
Get IIS status of SOS constraints.

\sphinxAtStartPar
\sphinxstylestrong{Synopsis}
\begin{quote}

\sphinxAtStartPar
\sphinxcode{\sphinxupquote{int GetSOSIIS(const SosArray \&soss, int *pIIS)}}
\end{quote}

\sphinxAtStartPar
\sphinxstylestrong{Arguments}
\begin{quote}

\sphinxAtStartPar
\sphinxcode{\sphinxupquote{soss}}: Array of SOS constraints.

\sphinxAtStartPar
\sphinxcode{\sphinxupquote{pIIS}}: IIS status of SOS constraints.
\end{quote}

\sphinxAtStartPar
\sphinxstylestrong{Return}
\begin{quote}

\sphinxAtStartPar
Number of SOS constraints.
\end{quote}
\end{quote}

\subsubsection{Model::GetSoss()}
\label{\detokenize{cppapi/model:model-getsoss}}\begin{quote}

\sphinxAtStartPar
Get all SOS constraints in model.

\sphinxAtStartPar
\sphinxstylestrong{Synopsis}
\begin{quote}

\sphinxAtStartPar
\sphinxcode{\sphinxupquote{SosArray GetSoss()}}
\end{quote}

\sphinxAtStartPar
\sphinxstylestrong{Return}
\begin{quote}

\sphinxAtStartPar
array object of SOS constraints.
\end{quote}
\end{quote}

\subsubsection{Model::GetSymMat()}
\label{\detokenize{cppapi/model:model-getsymmat}}\begin{quote}

\sphinxAtStartPar
Get a symmetric matrix of given index in model.

\sphinxAtStartPar
\sphinxstylestrong{Synopsis}
\begin{quote}

\sphinxAtStartPar
\sphinxcode{\sphinxupquote{SymMatrix GetSymMat(int idx)}}
\end{quote}

\sphinxAtStartPar
\sphinxstylestrong{Arguments}
\begin{quote}

\sphinxAtStartPar
\sphinxcode{\sphinxupquote{idx}}: index of the desired symmetric matrix.
\end{quote}

\sphinxAtStartPar
\sphinxstylestrong{Return}
\begin{quote}

\sphinxAtStartPar
the desired symmetric matrix object.
\end{quote}
\end{quote}

\subsubsection{Model::GetVar()}
\label{\detokenize{cppapi/model:model-getvar}}\begin{quote}

\sphinxAtStartPar
Get a variable of given index in model.

\sphinxAtStartPar
\sphinxstylestrong{Synopsis}
\begin{quote}

\sphinxAtStartPar
\sphinxcode{\sphinxupquote{Var GetVar(int idx)}}
\end{quote}

\sphinxAtStartPar
\sphinxstylestrong{Arguments}
\begin{quote}

\sphinxAtStartPar
\sphinxcode{\sphinxupquote{idx}}: index of the desired variable.
\end{quote}

\sphinxAtStartPar
\sphinxstylestrong{Return}
\begin{quote}

\sphinxAtStartPar
the desired variable object.
\end{quote}
\end{quote}

\subsubsection{Model::GetVarByName()}
\label{\detokenize{cppapi/model:model-getvarbyname}}\begin{quote}

\sphinxAtStartPar
Get a variable of given name in model.

\sphinxAtStartPar
\sphinxstylestrong{Synopsis}
\begin{quote}

\sphinxAtStartPar
\sphinxcode{\sphinxupquote{Var GetVarByName(const char *szName)}}
\end{quote}

\sphinxAtStartPar
\sphinxstylestrong{Arguments}
\begin{quote}

\sphinxAtStartPar
\sphinxcode{\sphinxupquote{szName}}: name of the desired variable.
\end{quote}

\sphinxAtStartPar
\sphinxstylestrong{Return}
\begin{quote}

\sphinxAtStartPar
the desired variable object.
\end{quote}
\end{quote}

\subsubsection{Model::GetVarLowerIIS()}
\label{\detokenize{cppapi/model:model-getvarloweriis}}\begin{quote}

\sphinxAtStartPar
Get IIS status of lower bounds of variables.

\sphinxAtStartPar
\sphinxstylestrong{Synopsis}
\begin{quote}

\sphinxAtStartPar
\sphinxcode{\sphinxupquote{int GetVarLowerIIS(const VarArray \&vars, int *pLowerIIS)}}
\end{quote}

\sphinxAtStartPar
\sphinxstylestrong{Arguments}
\begin{quote}

\sphinxAtStartPar
\sphinxcode{\sphinxupquote{vars}}: Array of variables

\sphinxAtStartPar
\sphinxcode{\sphinxupquote{pLowerIIS}}: IIS status of lower bounds of variables.
\end{quote}

\sphinxAtStartPar
\sphinxstylestrong{Return}
\begin{quote}

\sphinxAtStartPar
Number of variables.
\end{quote}
\end{quote}

\subsubsection{Model::GetVars()}
\label{\detokenize{cppapi/model:model-getvars}}\begin{quote}

\sphinxAtStartPar
Get all variables in model.

\sphinxAtStartPar
\sphinxstylestrong{Synopsis}
\begin{quote}

\sphinxAtStartPar
\sphinxcode{\sphinxupquote{VarArray GetVars()}}
\end{quote}

\sphinxAtStartPar
\sphinxstylestrong{Return}
\begin{quote}

\sphinxAtStartPar
variable array object.
\end{quote}
\end{quote}

\subsubsection{Model::GetVarUpperIIS()}
\label{\detokenize{cppapi/model:model-getvarupperiis}}\begin{quote}

\sphinxAtStartPar
Get IIS status of upper bounds of variables.

\sphinxAtStartPar
\sphinxstylestrong{Synopsis}
\begin{quote}

\sphinxAtStartPar
\sphinxcode{\sphinxupquote{int GetVarUpperIIS(const VarArray \&vars, int *pUpperIIS)}}
\end{quote}

\sphinxAtStartPar
\sphinxstylestrong{Arguments}
\begin{quote}

\sphinxAtStartPar
\sphinxcode{\sphinxupquote{vars}}: Array of variables

\sphinxAtStartPar
\sphinxcode{\sphinxupquote{pUpperIIS}}: IIS status of upper bounds of variables.
\end{quote}

\sphinxAtStartPar
\sphinxstylestrong{Return}
\begin{quote}

\sphinxAtStartPar
Number of variables.
\end{quote}
\end{quote}

\subsubsection{Model::Interrupt()}
\label{\detokenize{cppapi/model:model-interrupt}}\begin{quote}

\sphinxAtStartPar
Interrupt optimization of current problem.

\sphinxAtStartPar
\sphinxstylestrong{Synopsis}
\begin{quote}

\sphinxAtStartPar
\sphinxcode{\sphinxupquote{void Interrupt()}}
\end{quote}
\end{quote}

\subsubsection{Model::LoadMatrix()}
\label{\detokenize{cppapi/model:model-loadmatrix}}\begin{quote}

\sphinxAtStartPar
Load matrix data to formulate problem.

\sphinxAtStartPar
\sphinxstylestrong{Synopsis}
\begin{quote}

\sphinxAtStartPar
\sphinxcode{\sphinxupquote{void LoadMatrix(}}
\begin{quote}

\sphinxAtStartPar
\sphinxcode{\sphinxupquote{int nCol,}}

\sphinxAtStartPar
\sphinxcode{\sphinxupquote{int nRow,}}

\sphinxAtStartPar
\sphinxcode{\sphinxupquote{double *pCost,}}

\sphinxAtStartPar
\sphinxcode{\sphinxupquote{int *pMatBeg,}}

\sphinxAtStartPar
\sphinxcode{\sphinxupquote{int *pMatCnt,}}

\sphinxAtStartPar
\sphinxcode{\sphinxupquote{int *pMatIdx,}}

\sphinxAtStartPar
\sphinxcode{\sphinxupquote{double *pMatElem,}}

\sphinxAtStartPar
\sphinxcode{\sphinxupquote{double *pLHS,}}

\sphinxAtStartPar
\sphinxcode{\sphinxupquote{double *pRHS,}}

\sphinxAtStartPar
\sphinxcode{\sphinxupquote{double *pLB,}}

\sphinxAtStartPar
\sphinxcode{\sphinxupquote{double *pUB,}}

\sphinxAtStartPar
\sphinxcode{\sphinxupquote{char *pType)}}
\end{quote}
\end{quote}

\sphinxAtStartPar
\sphinxstylestrong{Arguments}
\begin{quote}

\sphinxAtStartPar
\sphinxcode{\sphinxupquote{nCol}}: Number of columns.

\sphinxAtStartPar
\sphinxcode{\sphinxupquote{nRow}}: Number of rows.

\sphinxAtStartPar
\sphinxcode{\sphinxupquote{pCost}}: Objective cost.

\sphinxAtStartPar
\sphinxcode{\sphinxupquote{pMatBeg}}: Begin pointer.

\sphinxAtStartPar
\sphinxcode{\sphinxupquote{pMatCnt}}: Counter pointer.

\sphinxAtStartPar
\sphinxcode{\sphinxupquote{pMatIdx}}: Index pointer.

\sphinxAtStartPar
\sphinxcode{\sphinxupquote{pMatElem}}: Elements pointer.

\sphinxAtStartPar
\sphinxcode{\sphinxupquote{pLHS}}: Lower bounds of rows.

\sphinxAtStartPar
\sphinxcode{\sphinxupquote{pRHS}}: Upper bounds of rows.

\sphinxAtStartPar
\sphinxcode{\sphinxupquote{pLB}}: Lower bounds of columns.

\sphinxAtStartPar
\sphinxcode{\sphinxupquote{pUB}}: Upper bounds of columns.

\sphinxAtStartPar
\sphinxcode{\sphinxupquote{pType}}: Column types.
\end{quote}
\end{quote}

\subsubsection{Model::LoadMipStart()}
\label{\detokenize{cppapi/model:model-loadmipstart}}\begin{quote}

\sphinxAtStartPar
Load final initial values of variables to the problem.

\sphinxAtStartPar
\sphinxstylestrong{Synopsis}
\begin{quote}

\sphinxAtStartPar
\sphinxcode{\sphinxupquote{void LoadMipStart()}}
\end{quote}
\end{quote}

\subsubsection{Model::LoadNlData()}
\label{\detokenize{cppapi/model:model-loadnldata}}\begin{quote}

\sphinxAtStartPar
Load nonlinear data into problem with customized callback.

\sphinxAtStartPar
\sphinxstylestrong{Synopsis}
\begin{quote}

\sphinxAtStartPar
\sphinxcode{\sphinxupquote{void LoadNlData(}}
\begin{quote}

\sphinxAtStartPar
\sphinxcode{\sphinxupquote{int nCol,}}

\sphinxAtStartPar
\sphinxcode{\sphinxupquote{int nRow,}}

\sphinxAtStartPar
\sphinxcode{\sphinxupquote{int iObjSense,}}

\sphinxAtStartPar
\sphinxcode{\sphinxupquote{int nObjGrad,}}

\sphinxAtStartPar
\sphinxcode{\sphinxupquote{const int *objGradIndex,}}

\sphinxAtStartPar
\sphinxcode{\sphinxupquote{int nnzJac,}}

\sphinxAtStartPar
\sphinxcode{\sphinxupquote{const int *jacRowIndex,}}

\sphinxAtStartPar
\sphinxcode{\sphinxupquote{const int *jacColIndex,}}

\sphinxAtStartPar
\sphinxcode{\sphinxupquote{int nnzHess,}}

\sphinxAtStartPar
\sphinxcode{\sphinxupquote{const int *hessRowIndex,}}

\sphinxAtStartPar
\sphinxcode{\sphinxupquote{const int *hessColIndex,}}

\sphinxAtStartPar
\sphinxcode{\sphinxupquote{const double *colLower,}}

\sphinxAtStartPar
\sphinxcode{\sphinxupquote{const double *colUpper,}}

\sphinxAtStartPar
\sphinxcode{\sphinxupquote{const double *rowLower,}}

\sphinxAtStartPar
\sphinxcode{\sphinxupquote{const double *rowUpper,}}

\sphinxAtStartPar
\sphinxcode{\sphinxupquote{const double *initColVal,}}

\sphinxAtStartPar
\sphinxcode{\sphinxupquote{int evalType,}}

\sphinxAtStartPar
\sphinxcode{\sphinxupquote{INlpCallback *pcb)}}
\end{quote}
\end{quote}

\sphinxAtStartPar
\sphinxstylestrong{Arguments}
\begin{quote}

\sphinxAtStartPar
\sphinxcode{\sphinxupquote{nCol}}: number of columns.

\sphinxAtStartPar
\sphinxcode{\sphinxupquote{nRow}}: number of rows.

\sphinxAtStartPar
\sphinxcode{\sphinxupquote{iObjSense}}: sense of objective.

\sphinxAtStartPar
\sphinxcode{\sphinxupquote{nObjGrad}}: number of objective gradients.

\sphinxAtStartPar
\sphinxcode{\sphinxupquote{objGradIndex}}: a list of index for objective gradients.

\sphinxAtStartPar
\sphinxcode{\sphinxupquote{nnzJac}}: number of non\sphinxhyphen{}zeros in Jacobian matrix.

\sphinxAtStartPar
\sphinxcode{\sphinxupquote{jacRowIndex}}: a list of row index in Jacobian matrix.

\sphinxAtStartPar
\sphinxcode{\sphinxupquote{jacColIndex}}: a list of col index in Jacobian matrix.

\sphinxAtStartPar
\sphinxcode{\sphinxupquote{nnzHess}}: number of non\sphinxhyphen{}zeros in Hessian matrix.

\sphinxAtStartPar
\sphinxcode{\sphinxupquote{hessRowIndex}}: a list of row index in Hessian matrix.

\sphinxAtStartPar
\sphinxcode{\sphinxupquote{hessColIndex}}: a list of col index in Hessian matrix.

\sphinxAtStartPar
\sphinxcode{\sphinxupquote{colLower}}: lower bounds for columns.

\sphinxAtStartPar
\sphinxcode{\sphinxupquote{colUpper}}: upper bounds for columns.

\sphinxAtStartPar
\sphinxcode{\sphinxupquote{rowLower}}: lower bounds for rows.

\sphinxAtStartPar
\sphinxcode{\sphinxupquote{rowUpper}}: upper bounds for rows.

\sphinxAtStartPar
\sphinxcode{\sphinxupquote{initColVal}}: a list of initial column values.

\sphinxAtStartPar
\sphinxcode{\sphinxupquote{evalType}}: evaluation type of the nonlinear model.

\sphinxAtStartPar
\sphinxcode{\sphinxupquote{pcb}}: pointer to callback object for nonlinear model.
\end{quote}
\end{quote}

\subsubsection{Model::LoadTuneParam()}
\label{\detokenize{cppapi/model:model-loadtuneparam}}\begin{quote}

\sphinxAtStartPar
Load specified tuned parameters into model.

\sphinxAtStartPar
\sphinxstylestrong{Synopsis}
\begin{quote}

\sphinxAtStartPar
\sphinxcode{\sphinxupquote{void LoadTuneParam(int idx)}}
\end{quote}

\sphinxAtStartPar
\sphinxstylestrong{Arguments}
\begin{quote}

\sphinxAtStartPar
\sphinxcode{\sphinxupquote{idx}}: Index of tuned parameters.
\end{quote}
\end{quote}

\subsubsection{Model::Read()}
\label{\detokenize{cppapi/model:model-read}}\begin{quote}

\sphinxAtStartPar
Read problem, solution, basis, MIP start or COPT parameters from file.

\sphinxAtStartPar
\sphinxstylestrong{Synopsis}
\begin{quote}

\sphinxAtStartPar
\sphinxcode{\sphinxupquote{void Read(const char *szFileName)}}
\end{quote}

\sphinxAtStartPar
\sphinxstylestrong{Arguments}
\begin{quote}

\sphinxAtStartPar
\sphinxcode{\sphinxupquote{szFileName}}: an input file name.
\end{quote}
\end{quote}

\subsubsection{Model::ReadBasis()}
\label{\detokenize{cppapi/model:model-readbasis}}\begin{quote}

\sphinxAtStartPar
Read basis from file.

\sphinxAtStartPar
\sphinxstylestrong{Synopsis}
\begin{quote}

\sphinxAtStartPar
\sphinxcode{\sphinxupquote{void ReadBasis(const char *szFileName)}}
\end{quote}

\sphinxAtStartPar
\sphinxstylestrong{Arguments}
\begin{quote}

\sphinxAtStartPar
\sphinxcode{\sphinxupquote{szFileName}}: an input file name.
\end{quote}
\end{quote}

\subsubsection{Model::ReadBin()}
\label{\detokenize{cppapi/model:model-readbin}}\begin{quote}

\sphinxAtStartPar
Read problem in COPT binary format from file.

\sphinxAtStartPar
\sphinxstylestrong{Synopsis}
\begin{quote}

\sphinxAtStartPar
\sphinxcode{\sphinxupquote{void ReadBin(const char *szFileName)}}
\end{quote}

\sphinxAtStartPar
\sphinxstylestrong{Arguments}
\begin{quote}

\sphinxAtStartPar
\sphinxcode{\sphinxupquote{szFileName}}: an input file name.
\end{quote}
\end{quote}

\subsubsection{Model::ReadCbf()}
\label{\detokenize{cppapi/model:model-readcbf}}\begin{quote}

\sphinxAtStartPar
Read problem in CBF format from file.

\sphinxAtStartPar
\sphinxstylestrong{Synopsis}
\begin{quote}

\sphinxAtStartPar
\sphinxcode{\sphinxupquote{void ReadCbf(const char *szFileName)}}
\end{quote}

\sphinxAtStartPar
\sphinxstylestrong{Arguments}
\begin{quote}

\sphinxAtStartPar
\sphinxcode{\sphinxupquote{szFileName}}: an input file name.
\end{quote}
\end{quote}

\subsubsection{Model::ReadJsonSol()}
\label{\detokenize{cppapi/model:model-readjsonsol}}\begin{quote}

\sphinxAtStartPar
Read solution in format of JSON from file.

\sphinxAtStartPar
\sphinxstylestrong{Synopsis}
\begin{quote}

\sphinxAtStartPar
\sphinxcode{\sphinxupquote{void ReadJsonSol(const char *szFileName)}}
\end{quote}

\sphinxAtStartPar
\sphinxstylestrong{Arguments}
\begin{quote}

\sphinxAtStartPar
\sphinxcode{\sphinxupquote{szFileName}}: an input file name.
\end{quote}
\end{quote}

\subsubsection{Model::ReadLp()}
\label{\detokenize{cppapi/model:model-readlp}}\begin{quote}

\sphinxAtStartPar
Read problem in LP format from file.

\sphinxAtStartPar
\sphinxstylestrong{Synopsis}
\begin{quote}

\sphinxAtStartPar
\sphinxcode{\sphinxupquote{void ReadLp(const char *szFileName)}}
\end{quote}

\sphinxAtStartPar
\sphinxstylestrong{Arguments}
\begin{quote}

\sphinxAtStartPar
\sphinxcode{\sphinxupquote{szFileName}}: an input file name.
\end{quote}
\end{quote}

\subsubsection{Model::ReadMps()}
\label{\detokenize{cppapi/model:model-readmps}}\begin{quote}

\sphinxAtStartPar
Read problem in MPS format from file.

\sphinxAtStartPar
\sphinxstylestrong{Synopsis}
\begin{quote}

\sphinxAtStartPar
\sphinxcode{\sphinxupquote{void ReadMps(const char *szFileName)}}
\end{quote}

\sphinxAtStartPar
\sphinxstylestrong{Arguments}
\begin{quote}

\sphinxAtStartPar
\sphinxcode{\sphinxupquote{szFileName}}: an input file name.
\end{quote}
\end{quote}

\subsubsection{Model::ReadMst()}
\label{\detokenize{cppapi/model:model-readmst}}\begin{quote}

\sphinxAtStartPar
Read MIP start information from file.

\sphinxAtStartPar
\sphinxstylestrong{Synopsis}
\begin{quote}

\sphinxAtStartPar
\sphinxcode{\sphinxupquote{void ReadMst(const char *szFileName)}}
\end{quote}

\sphinxAtStartPar
\sphinxstylestrong{Arguments}
\begin{quote}

\sphinxAtStartPar
\sphinxcode{\sphinxupquote{szFileName}}: an input file name.
\end{quote}
\end{quote}

\subsubsection{Model::ReadOrd()}
\label{\detokenize{cppapi/model:model-readord}}\begin{quote}

\sphinxAtStartPar
Read branching order from file.

\sphinxAtStartPar
\sphinxstylestrong{Synopsis}
\begin{quote}

\sphinxAtStartPar
\sphinxcode{\sphinxupquote{void ReadOrd(const char *szFileName)}}
\end{quote}

\sphinxAtStartPar
\sphinxstylestrong{Arguments}
\begin{quote}

\sphinxAtStartPar
\sphinxcode{\sphinxupquote{szFileName}}: an input file name.
\end{quote}
\end{quote}

\subsubsection{Model::ReadParam()}
\label{\detokenize{cppapi/model:model-readparam}}\begin{quote}

\sphinxAtStartPar
Read COPT parameters from file.

\sphinxAtStartPar
\sphinxstylestrong{Synopsis}
\begin{quote}

\sphinxAtStartPar
\sphinxcode{\sphinxupquote{void ReadParam(const char *szFileName)}}
\end{quote}

\sphinxAtStartPar
\sphinxstylestrong{Arguments}
\begin{quote}

\sphinxAtStartPar
\sphinxcode{\sphinxupquote{szFileName}}: an input file name.
\end{quote}
\end{quote}

\subsubsection{Model::ReadSdpa()}
\label{\detokenize{cppapi/model:model-readsdpa}}\begin{quote}

\sphinxAtStartPar
Read problem in SDPA format from file.

\sphinxAtStartPar
\sphinxstylestrong{Synopsis}
\begin{quote}

\sphinxAtStartPar
\sphinxcode{\sphinxupquote{void ReadSdpa(const char *szFileName)}}
\end{quote}

\sphinxAtStartPar
\sphinxstylestrong{Arguments}
\begin{quote}

\sphinxAtStartPar
\sphinxcode{\sphinxupquote{szFileName}}: an input file name.
\end{quote}
\end{quote}

\subsubsection{Model::ReadSol()}
\label{\detokenize{cppapi/model:model-readsol}}\begin{quote}

\sphinxAtStartPar
Read solution from file.

\sphinxAtStartPar
\sphinxstylestrong{Synopsis}
\begin{quote}

\sphinxAtStartPar
\sphinxcode{\sphinxupquote{void ReadSol(const char *szFileName)}}
\end{quote}

\sphinxAtStartPar
\sphinxstylestrong{Arguments}
\begin{quote}

\sphinxAtStartPar
\sphinxcode{\sphinxupquote{szFileName}}: an input file name.
\end{quote}
\end{quote}

\subsubsection{Model::ReadTune()}
\label{\detokenize{cppapi/model:model-readtune}}\begin{quote}

\sphinxAtStartPar
Read tuning parameters from file.

\sphinxAtStartPar
\sphinxstylestrong{Synopsis}
\begin{quote}

\sphinxAtStartPar
\sphinxcode{\sphinxupquote{void ReadTune(const char *szFileName)}}
\end{quote}

\sphinxAtStartPar
\sphinxstylestrong{Arguments}
\begin{quote}

\sphinxAtStartPar
\sphinxcode{\sphinxupquote{szFileName}}: an input file name.
\end{quote}
\end{quote}

\subsubsection{Model::Remove()}
\label{\detokenize{cppapi/model:model-remove}}\begin{quote}

\sphinxAtStartPar
Remove a list of variables from model.

\sphinxAtStartPar
\sphinxstylestrong{Synopsis}
\begin{quote}

\sphinxAtStartPar
\sphinxcode{\sphinxupquote{void Remove(VarArray \&vars)}}
\end{quote}

\sphinxAtStartPar
\sphinxstylestrong{Arguments}
\begin{quote}

\sphinxAtStartPar
\sphinxcode{\sphinxupquote{vars}}: an array of variables.
\end{quote}
\end{quote}

\subsubsection{Model::Remove()}
\label{\detokenize{cppapi/model:id70}}\begin{quote}

\sphinxAtStartPar
Remove a list of constraints from model.

\sphinxAtStartPar
\sphinxstylestrong{Synopsis}
\begin{quote}

\sphinxAtStartPar
\sphinxcode{\sphinxupquote{void Remove(ConstrArray \&constrs)}}
\end{quote}

\sphinxAtStartPar
\sphinxstylestrong{Arguments}
\begin{quote}

\sphinxAtStartPar
\sphinxcode{\sphinxupquote{constrs}}: an array of constraints.
\end{quote}
\end{quote}

\subsubsection{Model::Remove()}
\label{\detokenize{cppapi/model:id71}}\begin{quote}

\sphinxAtStartPar
Remove a list of nonlinear constraints from model.

\sphinxAtStartPar
\sphinxstylestrong{Synopsis}
\begin{quote}

\sphinxAtStartPar
\sphinxcode{\sphinxupquote{void Remove(NlConstrArray \&constrs)}}
\end{quote}

\sphinxAtStartPar
\sphinxstylestrong{Arguments}
\begin{quote}

\sphinxAtStartPar
\sphinxcode{\sphinxupquote{constrs}}: an array of nonlinear constraints.
\end{quote}
\end{quote}

\subsubsection{Model::Remove()}
\label{\detokenize{cppapi/model:id72}}\begin{quote}

\sphinxAtStartPar
Remove a list of SOS constraints from model.

\sphinxAtStartPar
\sphinxstylestrong{Synopsis}
\begin{quote}

\sphinxAtStartPar
\sphinxcode{\sphinxupquote{void Remove(SosArray \&soss)}}
\end{quote}

\sphinxAtStartPar
\sphinxstylestrong{Arguments}
\begin{quote}

\sphinxAtStartPar
\sphinxcode{\sphinxupquote{soss}}: an array of SOS constraints.
\end{quote}
\end{quote}

\subsubsection{Model::Remove()}
\label{\detokenize{cppapi/model:id73}}\begin{quote}

\sphinxAtStartPar
Remove a list of gernal constraints from model.

\sphinxAtStartPar
\sphinxstylestrong{Synopsis}
\begin{quote}

\sphinxAtStartPar
\sphinxcode{\sphinxupquote{void Remove(GenConstrArray \&genConstrs)}}
\end{quote}

\sphinxAtStartPar
\sphinxstylestrong{Arguments}
\begin{quote}

\sphinxAtStartPar
\sphinxcode{\sphinxupquote{genConstrs}}: an array of general constraints.
\end{quote}
\end{quote}

\subsubsection{Model::Remove()}
\label{\detokenize{cppapi/model:id74}}\begin{quote}

\sphinxAtStartPar
Remove a list of cone constraints from model.

\sphinxAtStartPar
\sphinxstylestrong{Synopsis}
\begin{quote}

\sphinxAtStartPar
\sphinxcode{\sphinxupquote{void Remove(ConeArray \&cones)}}
\end{quote}

\sphinxAtStartPar
\sphinxstylestrong{Arguments}
\begin{quote}

\sphinxAtStartPar
\sphinxcode{\sphinxupquote{cones}}: an array of cone constraints.
\end{quote}
\end{quote}

\subsubsection{Model::Remove()}
\label{\detokenize{cppapi/model:id75}}\begin{quote}

\sphinxAtStartPar
Remove a list of exponential cone constraints from model.

\sphinxAtStartPar
\sphinxstylestrong{Synopsis}
\begin{quote}

\sphinxAtStartPar
\sphinxcode{\sphinxupquote{void Remove(ExpConeArray \&cones)}}
\end{quote}

\sphinxAtStartPar
\sphinxstylestrong{Arguments}
\begin{quote}

\sphinxAtStartPar
\sphinxcode{\sphinxupquote{cones}}: an array of exponential cone constraints.
\end{quote}
\end{quote}

\subsubsection{Model::Remove()}
\label{\detokenize{cppapi/model:id76}}\begin{quote}

\sphinxAtStartPar
Remove a list of affine cone constraints from model.

\sphinxAtStartPar
\sphinxstylestrong{Synopsis}
\begin{quote}

\sphinxAtStartPar
\sphinxcode{\sphinxupquote{void Remove(AffineConeArray \&cones)}}
\end{quote}

\sphinxAtStartPar
\sphinxstylestrong{Arguments}
\begin{quote}

\sphinxAtStartPar
\sphinxcode{\sphinxupquote{cones}}: an array of affine cone constraints.
\end{quote}
\end{quote}

\subsubsection{Model::Remove()}
\label{\detokenize{cppapi/model:id77}}\begin{quote}

\sphinxAtStartPar
Remove a list of quadratic constraints from model.

\sphinxAtStartPar
\sphinxstylestrong{Synopsis}
\begin{quote}

\sphinxAtStartPar
\sphinxcode{\sphinxupquote{void Remove(QConstrArray \&qconstrs)}}
\end{quote}

\sphinxAtStartPar
\sphinxstylestrong{Arguments}
\begin{quote}

\sphinxAtStartPar
\sphinxcode{\sphinxupquote{qconstrs}}: an array of quadratic constraints.
\end{quote}
\end{quote}

\subsubsection{Model::Remove()}
\label{\detokenize{cppapi/model:id78}}\begin{quote}

\sphinxAtStartPar
Remove a list of PSD variables from model.

\sphinxAtStartPar
\sphinxstylestrong{Synopsis}
\begin{quote}

\sphinxAtStartPar
\sphinxcode{\sphinxupquote{void Remove(PsdVarArray \&vars)}}
\end{quote}

\sphinxAtStartPar
\sphinxstylestrong{Arguments}
\begin{quote}

\sphinxAtStartPar
\sphinxcode{\sphinxupquote{vars}}: an array of PSD variables.
\end{quote}
\end{quote}

\subsubsection{Model::Remove()}
\label{\detokenize{cppapi/model:id79}}\begin{quote}

\sphinxAtStartPar
Remove a list of PSD constraints from model.

\sphinxAtStartPar
\sphinxstylestrong{Synopsis}
\begin{quote}

\sphinxAtStartPar
\sphinxcode{\sphinxupquote{void Remove(PsdConstrArray \&constrs)}}
\end{quote}

\sphinxAtStartPar
\sphinxstylestrong{Arguments}
\begin{quote}

\sphinxAtStartPar
\sphinxcode{\sphinxupquote{constrs}}: an array of PSD constraints.
\end{quote}
\end{quote}

\subsubsection{Model::Remove()}
\label{\detokenize{cppapi/model:id80}}\begin{quote}

\sphinxAtStartPar
Remove a list of LMI constraints from model.

\sphinxAtStartPar
\sphinxstylestrong{Synopsis}
\begin{quote}

\sphinxAtStartPar
\sphinxcode{\sphinxupquote{void Remove(LmiConstrArray \&constrs)}}
\end{quote}

\sphinxAtStartPar
\sphinxstylestrong{Arguments}
\begin{quote}

\sphinxAtStartPar
\sphinxcode{\sphinxupquote{constrs}}: an array of LMI constraints.
\end{quote}
\end{quote}

\subsubsection{Model::Reset()}
\label{\detokenize{cppapi/model:model-reset}}\begin{quote}

\sphinxAtStartPar
Reset solution of problem only.

\sphinxAtStartPar
\sphinxstylestrong{Synopsis}
\begin{quote}

\sphinxAtStartPar
\sphinxcode{\sphinxupquote{void Reset()}}
\end{quote}
\end{quote}

\subsubsection{Model::ResetAll()}
\label{\detokenize{cppapi/model:model-resetall}}\begin{quote}

\sphinxAtStartPar
Reset solution of problem, and additional information such as MIP start, etc.

\sphinxAtStartPar
\sphinxstylestrong{Synopsis}
\begin{quote}

\sphinxAtStartPar
\sphinxcode{\sphinxupquote{void ResetAll()}}
\end{quote}
\end{quote}

\subsubsection{Model::ResetObjParamN()}
\label{\detokenize{cppapi/model:model-resetobjparamn}}\begin{quote}

\sphinxAtStartPar
Reset objective parameters of a multi\sphinxhyphen{}objective function.

\sphinxAtStartPar
\sphinxstylestrong{Synopsis}
\begin{quote}

\sphinxAtStartPar
\sphinxcode{\sphinxupquote{void ResetObjParamN(int idx)}}
\end{quote}

\sphinxAtStartPar
\sphinxstylestrong{Arguments}
\begin{quote}

\sphinxAtStartPar
\sphinxcode{\sphinxupquote{idx}}: index of a multi\sphinxhyphen{}objective function.
\end{quote}
\end{quote}

\subsubsection{Model::ResetParam()}
\label{\detokenize{cppapi/model:model-resetparam}}\begin{quote}

\sphinxAtStartPar
Reset parameters to default settings.

\sphinxAtStartPar
\sphinxstylestrong{Synopsis}
\begin{quote}

\sphinxAtStartPar
\sphinxcode{\sphinxupquote{void ResetParam()}}
\end{quote}
\end{quote}

\subsubsection{Model::ResetParamN()}
\label{\detokenize{cppapi/model:model-resetparamn}}\begin{quote}

\sphinxAtStartPar
Reset double and integer parameters of a multi\sphinxhyphen{}objective function.

\sphinxAtStartPar
\sphinxstylestrong{Synopsis}
\begin{quote}

\sphinxAtStartPar
\sphinxcode{\sphinxupquote{void ResetParamN(int idx)}}
\end{quote}

\sphinxAtStartPar
\sphinxstylestrong{Arguments}
\begin{quote}

\sphinxAtStartPar
\sphinxcode{\sphinxupquote{idx}}: index of a multi\sphinxhyphen{}objective function.
\end{quote}
\end{quote}

\subsubsection{Model::Set()}
\label{\detokenize{cppapi/model:model-set}}\begin{quote}

\sphinxAtStartPar
Set values of information associated with variables.

\sphinxAtStartPar
\sphinxstylestrong{Synopsis}
\begin{quote}

\sphinxAtStartPar
\sphinxcode{\sphinxupquote{void Set(}}
\begin{quote}

\sphinxAtStartPar
\sphinxcode{\sphinxupquote{const char *szName,}}

\sphinxAtStartPar
\sphinxcode{\sphinxupquote{const VarArray \&vars,}}

\sphinxAtStartPar
\sphinxcode{\sphinxupquote{double *pVals,}}

\sphinxAtStartPar
\sphinxcode{\sphinxupquote{int len)}}
\end{quote}
\end{quote}

\sphinxAtStartPar
\sphinxstylestrong{Arguments}
\begin{quote}

\sphinxAtStartPar
\sphinxcode{\sphinxupquote{szName}}: name of double information.

\sphinxAtStartPar
\sphinxcode{\sphinxupquote{vars}}: a list of desired variables.

\sphinxAtStartPar
\sphinxcode{\sphinxupquote{pVals}}: array of information values.

\sphinxAtStartPar
\sphinxcode{\sphinxupquote{len}}: length of value array.
\end{quote}
\end{quote}

\subsubsection{Model::Set()}
\label{\detokenize{cppapi/model:id81}}\begin{quote}

\sphinxAtStartPar
Set values of information associated with constraints.

\sphinxAtStartPar
\sphinxstylestrong{Synopsis}
\begin{quote}

\sphinxAtStartPar
\sphinxcode{\sphinxupquote{void Set(}}
\begin{quote}

\sphinxAtStartPar
\sphinxcode{\sphinxupquote{const char *szName,}}

\sphinxAtStartPar
\sphinxcode{\sphinxupquote{const ConstrArray \&constrs,}}

\sphinxAtStartPar
\sphinxcode{\sphinxupquote{double *pVals,}}

\sphinxAtStartPar
\sphinxcode{\sphinxupquote{int len)}}
\end{quote}
\end{quote}

\sphinxAtStartPar
\sphinxstylestrong{Arguments}
\begin{quote}

\sphinxAtStartPar
\sphinxcode{\sphinxupquote{szName}}: name of double information.

\sphinxAtStartPar
\sphinxcode{\sphinxupquote{constrs}}: a list of desired constraints.

\sphinxAtStartPar
\sphinxcode{\sphinxupquote{pVals}}: array of information values.

\sphinxAtStartPar
\sphinxcode{\sphinxupquote{len}}: length of value array.
\end{quote}
\end{quote}

\subsubsection{Model::Set()}
\label{\detokenize{cppapi/model:id82}}\begin{quote}

\sphinxAtStartPar
Set values of information associated with nonlinear constraints.

\sphinxAtStartPar
\sphinxstylestrong{Synopsis}
\begin{quote}

\sphinxAtStartPar
\sphinxcode{\sphinxupquote{void Set(}}
\begin{quote}

\sphinxAtStartPar
\sphinxcode{\sphinxupquote{const char *szName,}}

\sphinxAtStartPar
\sphinxcode{\sphinxupquote{const NlConstrArray \&constrs,}}

\sphinxAtStartPar
\sphinxcode{\sphinxupquote{double *pVals,}}

\sphinxAtStartPar
\sphinxcode{\sphinxupquote{int len)}}
\end{quote}
\end{quote}

\sphinxAtStartPar
\sphinxstylestrong{Arguments}
\begin{quote}

\sphinxAtStartPar
\sphinxcode{\sphinxupquote{szName}}: name of double information.

\sphinxAtStartPar
\sphinxcode{\sphinxupquote{constrs}}: a list of desired nonlinear constraints.

\sphinxAtStartPar
\sphinxcode{\sphinxupquote{pVals}}: array of values of information.

\sphinxAtStartPar
\sphinxcode{\sphinxupquote{len}}: length of value array.
\end{quote}
\end{quote}

\subsubsection{Model::Set()}
\label{\detokenize{cppapi/model:id83}}\begin{quote}

\sphinxAtStartPar
Set values of information associated with PSD constraints.

\sphinxAtStartPar
\sphinxstylestrong{Synopsis}
\begin{quote}

\sphinxAtStartPar
\sphinxcode{\sphinxupquote{void Set(}}
\begin{quote}

\sphinxAtStartPar
\sphinxcode{\sphinxupquote{const char *szName,}}

\sphinxAtStartPar
\sphinxcode{\sphinxupquote{const PsdConstrArray \&constrs,}}

\sphinxAtStartPar
\sphinxcode{\sphinxupquote{double *pVals,}}

\sphinxAtStartPar
\sphinxcode{\sphinxupquote{int len)}}
\end{quote}
\end{quote}

\sphinxAtStartPar
\sphinxstylestrong{Arguments}
\begin{quote}

\sphinxAtStartPar
\sphinxcode{\sphinxupquote{szName}}: name of double information.

\sphinxAtStartPar
\sphinxcode{\sphinxupquote{constrs}}: a list of desired PSD constraints.

\sphinxAtStartPar
\sphinxcode{\sphinxupquote{pVals}}: array of values of information.

\sphinxAtStartPar
\sphinxcode{\sphinxupquote{len}}: length of value array.
\end{quote}
\end{quote}

\subsubsection{Model::SetBasis()}
\label{\detokenize{cppapi/model:model-setbasis}}\begin{quote}

\sphinxAtStartPar
Set column and row basis status to model.

\sphinxAtStartPar
\sphinxstylestrong{Synopsis}
\begin{quote}

\sphinxAtStartPar
\sphinxcode{\sphinxupquote{void SetBasis(int *pColBasis, int *pRowBasis)}}
\end{quote}

\sphinxAtStartPar
\sphinxstylestrong{Arguments}
\begin{quote}

\sphinxAtStartPar
\sphinxcode{\sphinxupquote{pColBasis}}: pointer to status of column basis.

\sphinxAtStartPar
\sphinxcode{\sphinxupquote{pRowBasis}}: pointer to status of row basis.
\end{quote}
\end{quote}

\subsubsection{Model::SetCallback()}
\label{\detokenize{cppapi/model:model-setcallback}}\begin{quote}

\sphinxAtStartPar
Set user callback to COPT model.

\sphinxAtStartPar
\sphinxstylestrong{Synopsis}
\begin{quote}

\sphinxAtStartPar
\sphinxcode{\sphinxupquote{void SetCallback(ICallback *pcb, int cbctx)}}
\end{quote}

\sphinxAtStartPar
\sphinxstylestrong{Arguments}
\begin{quote}

\sphinxAtStartPar
\sphinxcode{\sphinxupquote{pcb}}: pointer to user callback object.

\sphinxAtStartPar
\sphinxcode{\sphinxupquote{cbctx}}: COPT callback context, which is declared in public C header.
\end{quote}
\end{quote}

\subsubsection{Model::SetCoeff()}
\label{\detokenize{cppapi/model:model-setcoeff}}\begin{quote}

\sphinxAtStartPar
Set the coefficient of a variable in a linear constraint.

\sphinxAtStartPar
\sphinxstylestrong{Synopsis}
\begin{quote}

\sphinxAtStartPar
\sphinxcode{\sphinxupquote{void SetCoeff(}}
\begin{quote}

\sphinxAtStartPar
\sphinxcode{\sphinxupquote{const Constraint \&constr,}}

\sphinxAtStartPar
\sphinxcode{\sphinxupquote{const Var \&var,}}

\sphinxAtStartPar
\sphinxcode{\sphinxupquote{double newVal)}}
\end{quote}
\end{quote}

\sphinxAtStartPar
\sphinxstylestrong{Arguments}
\begin{quote}

\sphinxAtStartPar
\sphinxcode{\sphinxupquote{constr}}: The requested constraint.

\sphinxAtStartPar
\sphinxcode{\sphinxupquote{var}}: The requested variable.

\sphinxAtStartPar
\sphinxcode{\sphinxupquote{newVal}}: New coefficient.
\end{quote}
\end{quote}

\subsubsection{Model::SetCoeffs()}
\label{\detokenize{cppapi/model:model-setcoeffs}}\begin{quote}

\sphinxAtStartPar
Set a list of coefficients in the model.

\sphinxAtStartPar
\sphinxstylestrong{Synopsis}
\begin{quote}

\sphinxAtStartPar
\sphinxcode{\sphinxupquote{void SetCoeffs(}}
\begin{quote}

\sphinxAtStartPar
\sphinxcode{\sphinxupquote{const ConstrArray \&constrs,}}

\sphinxAtStartPar
\sphinxcode{\sphinxupquote{const VarArray \&vars,}}

\sphinxAtStartPar
\sphinxcode{\sphinxupquote{double *vals,}}

\sphinxAtStartPar
\sphinxcode{\sphinxupquote{int len)}}
\end{quote}
\end{quote}

\sphinxAtStartPar
\sphinxstylestrong{Arguments}
\begin{quote}

\sphinxAtStartPar
\sphinxcode{\sphinxupquote{constrs}}: A list of constraints for coefficients to be set.

\sphinxAtStartPar
\sphinxcode{\sphinxupquote{vars}}: A list of vars for coefficients to be set.

\sphinxAtStartPar
\sphinxcode{\sphinxupquote{vals}}: New values for coefficients.

\sphinxAtStartPar
\sphinxcode{\sphinxupquote{len}}: Length of values.
\end{quote}
\end{quote}

\subsubsection{Model::SetCost()}
\label{\detokenize{cppapi/model:model-setcost}}\begin{quote}

\sphinxAtStartPar
Set objective for model

\sphinxAtStartPar
\sphinxstylestrong{Synopsis}
\begin{quote}

\sphinxAtStartPar
\sphinxcode{\sphinxupquote{void SetCost(}}
\begin{quote}

\sphinxAtStartPar
\sphinxcode{\sphinxupquote{int num,}}

\sphinxAtStartPar
\sphinxcode{\sphinxupquote{const int *list,}}

\sphinxAtStartPar
\sphinxcode{\sphinxupquote{const double *obj)}}
\end{quote}
\end{quote}

\sphinxAtStartPar
\sphinxstylestrong{Arguments}
\begin{quote}

\sphinxAtStartPar
\sphinxcode{\sphinxupquote{num}}: number of terms in objective.

\sphinxAtStartPar
\sphinxcode{\sphinxupquote{list}}: variable indexes in objective.

\sphinxAtStartPar
\sphinxcode{\sphinxupquote{obj}}: corresponding coefficients in objective.
\end{quote}
\end{quote}

\subsubsection{Model::SetDblParam()}
\label{\detokenize{cppapi/model:model-setdblparam}}\begin{quote}

\sphinxAtStartPar
Set value of a COPT double parameter.

\sphinxAtStartPar
\sphinxstylestrong{Synopsis}
\begin{quote}

\sphinxAtStartPar
\sphinxcode{\sphinxupquote{void SetDblParam(const char *szParam, double dVal)}}
\end{quote}

\sphinxAtStartPar
\sphinxstylestrong{Arguments}
\begin{quote}

\sphinxAtStartPar
\sphinxcode{\sphinxupquote{szParam}}: name of integer parameter.

\sphinxAtStartPar
\sphinxcode{\sphinxupquote{dVal}}: double value.
\end{quote}
\end{quote}

\subsubsection{Model::SetDblParamN()}
\label{\detokenize{cppapi/model:model-setdblparamn}}\begin{quote}

\sphinxAtStartPar
Set value of a double parameter of a multi\sphinxhyphen{}objective function.

\sphinxAtStartPar
\sphinxstylestrong{Synopsis}
\begin{quote}

\sphinxAtStartPar
\sphinxcode{\sphinxupquote{void SetDblParamN(}}
\begin{quote}

\sphinxAtStartPar
\sphinxcode{\sphinxupquote{int idx,}}

\sphinxAtStartPar
\sphinxcode{\sphinxupquote{const char *szParam,}}

\sphinxAtStartPar
\sphinxcode{\sphinxupquote{double val)}}
\end{quote}
\end{quote}

\sphinxAtStartPar
\sphinxstylestrong{Arguments}
\begin{quote}

\sphinxAtStartPar
\sphinxcode{\sphinxupquote{idx}}: index of a multi\sphinxhyphen{}objective function.

\sphinxAtStartPar
\sphinxcode{\sphinxupquote{szParam}}: name of double parameter.

\sphinxAtStartPar
\sphinxcode{\sphinxupquote{val}}: new value of double parameter.
\end{quote}
\end{quote}

\subsubsection{Model::SetIntParam()}
\label{\detokenize{cppapi/model:model-setintparam}}\begin{quote}

\sphinxAtStartPar
Set value of a COPT integer parameter.

\sphinxAtStartPar
\sphinxstylestrong{Synopsis}
\begin{quote}

\sphinxAtStartPar
\sphinxcode{\sphinxupquote{void SetIntParam(const char *szParam, int nVal)}}
\end{quote}

\sphinxAtStartPar
\sphinxstylestrong{Arguments}
\begin{quote}

\sphinxAtStartPar
\sphinxcode{\sphinxupquote{szParam}}: name of integer parameter.

\sphinxAtStartPar
\sphinxcode{\sphinxupquote{nVal}}: integer value.
\end{quote}
\end{quote}

\subsubsection{Model::SetIntParamN()}
\label{\detokenize{cppapi/model:model-setintparamn}}\begin{quote}

\sphinxAtStartPar
Set value of an integer parameter of a multi\sphinxhyphen{}objective function.

\sphinxAtStartPar
\sphinxstylestrong{Synopsis}
\begin{quote}

\sphinxAtStartPar
\sphinxcode{\sphinxupquote{void SetIntParamN(}}
\begin{quote}

\sphinxAtStartPar
\sphinxcode{\sphinxupquote{int idx,}}

\sphinxAtStartPar
\sphinxcode{\sphinxupquote{const char *szParam,}}

\sphinxAtStartPar
\sphinxcode{\sphinxupquote{int val)}}
\end{quote}
\end{quote}

\sphinxAtStartPar
\sphinxstylestrong{Arguments}
\begin{quote}

\sphinxAtStartPar
\sphinxcode{\sphinxupquote{idx}}: index of a multi\sphinxhyphen{}objective function.

\sphinxAtStartPar
\sphinxcode{\sphinxupquote{szParam}}: name of integer parameter.

\sphinxAtStartPar
\sphinxcode{\sphinxupquote{val}}: new value of integer parameter.
\end{quote}
\end{quote}

\subsubsection{Model::SetLmiCoeff()}
\label{\detokenize{cppapi/model:model-setlmicoeff}}\begin{quote}

\sphinxAtStartPar
Set the coefficient matrix of a variable in LMI constraint.

\sphinxAtStartPar
\sphinxstylestrong{Synopsis}
\begin{quote}

\sphinxAtStartPar
\sphinxcode{\sphinxupquote{void SetLmiCoeff(}}
\begin{quote}

\sphinxAtStartPar
\sphinxcode{\sphinxupquote{const LmiConstraint \&constr,}}

\sphinxAtStartPar
\sphinxcode{\sphinxupquote{const Var \&var,}}

\sphinxAtStartPar
\sphinxcode{\sphinxupquote{const SymMatrix \&mat)}}
\end{quote}
\end{quote}

\sphinxAtStartPar
\sphinxstylestrong{Arguments}
\begin{quote}

\sphinxAtStartPar
\sphinxcode{\sphinxupquote{constr}}: The desired LMI constraint.

\sphinxAtStartPar
\sphinxcode{\sphinxupquote{var}}: The desired variable.

\sphinxAtStartPar
\sphinxcode{\sphinxupquote{mat}}: new coefficient matrix.
\end{quote}
\end{quote}

\subsubsection{Model::SetLmiRhs()}
\label{\detokenize{cppapi/model:model-setlmirhs}}\begin{quote}

\sphinxAtStartPar
Set constant matrix of LMI constraint.

\sphinxAtStartPar
\sphinxstylestrong{Synopsis}
\begin{quote}

\sphinxAtStartPar
\sphinxcode{\sphinxupquote{void SetLmiRhs(const LmiConstraint \&constr, const SymMatrix \&mat)}}
\end{quote}

\sphinxAtStartPar
\sphinxstylestrong{Arguments}
\begin{quote}

\sphinxAtStartPar
\sphinxcode{\sphinxupquote{constr}}: The desired LMI constraint.

\sphinxAtStartPar
\sphinxcode{\sphinxupquote{mat}}: new constant matrix.
\end{quote}
\end{quote}

\subsubsection{Model::SetLpSolution()}
\label{\detokenize{cppapi/model:model-setlpsolution}}\begin{quote}

\sphinxAtStartPar
Set LP solution.

\sphinxAtStartPar
\sphinxstylestrong{Synopsis}
\begin{quote}

\sphinxAtStartPar
\sphinxcode{\sphinxupquote{void SetLpSolution(}}
\begin{quote}

\sphinxAtStartPar
\sphinxcode{\sphinxupquote{double *pValue,}}

\sphinxAtStartPar
\sphinxcode{\sphinxupquote{double *pSlack,}}

\sphinxAtStartPar
\sphinxcode{\sphinxupquote{double *pRowDual,}}

\sphinxAtStartPar
\sphinxcode{\sphinxupquote{double *pRedCost)}}
\end{quote}
\end{quote}

\sphinxAtStartPar
\sphinxstylestrong{Arguments}
\begin{quote}

\sphinxAtStartPar
\sphinxcode{\sphinxupquote{pValue}}: double pointer to solution values.

\sphinxAtStartPar
\sphinxcode{\sphinxupquote{pSlack}}: double poitner to slack values.

\sphinxAtStartPar
\sphinxcode{\sphinxupquote{pRowDual}}: double pointer to dual values.

\sphinxAtStartPar
\sphinxcode{\sphinxupquote{pRedCost}}: double pointer to reduced costs.
\end{quote}
\end{quote}

\subsubsection{Model::SetMipStart()}
\label{\detokenize{cppapi/model:model-setmipstart}}\begin{quote}

\sphinxAtStartPar
Set initial values for variables of given number, starting from the first one.

\sphinxAtStartPar
\sphinxstylestrong{Synopsis}
\begin{quote}

\sphinxAtStartPar
\sphinxcode{\sphinxupquote{void SetMipStart(int count, double *pVals)}}
\end{quote}

\sphinxAtStartPar
\sphinxstylestrong{Arguments}
\begin{quote}

\sphinxAtStartPar
\sphinxcode{\sphinxupquote{count}}: the number of variables to set.

\sphinxAtStartPar
\sphinxcode{\sphinxupquote{pVals}}: pointer to initial values of variables.
\end{quote}
\end{quote}

\subsubsection{Model::SetMipStart()}
\label{\detokenize{cppapi/model:id84}}\begin{quote}

\sphinxAtStartPar
Set initial value for the specified variable.

\sphinxAtStartPar
\sphinxstylestrong{Synopsis}
\begin{quote}

\sphinxAtStartPar
\sphinxcode{\sphinxupquote{void SetMipStart(const Var \&var, double val)}}
\end{quote}

\sphinxAtStartPar
\sphinxstylestrong{Arguments}
\begin{quote}

\sphinxAtStartPar
\sphinxcode{\sphinxupquote{var}}: an interested variable.

\sphinxAtStartPar
\sphinxcode{\sphinxupquote{val}}: initial value of the variable.
\end{quote}
\end{quote}

\subsubsection{Model::SetMipStart()}
\label{\detokenize{cppapi/model:id85}}\begin{quote}

\sphinxAtStartPar
Set initial values for an array of variables.

\sphinxAtStartPar
\sphinxstylestrong{Synopsis}
\begin{quote}

\sphinxAtStartPar
\sphinxcode{\sphinxupquote{void SetMipStart(const VarArray \&vars, double *pVals)}}
\end{quote}

\sphinxAtStartPar
\sphinxstylestrong{Arguments}
\begin{quote}

\sphinxAtStartPar
\sphinxcode{\sphinxupquote{vars}}: a list of interested variables.

\sphinxAtStartPar
\sphinxcode{\sphinxupquote{pVals}}: pointer to initial values of variables.
\end{quote}
\end{quote}

\subsubsection{Model::SetNames()}
\label{\detokenize{cppapi/model:model-setnames}}\begin{quote}

\sphinxAtStartPar
Set names for given variables in model.

\sphinxAtStartPar
\sphinxstylestrong{Synopsis}
\begin{quote}

\sphinxAtStartPar
\sphinxcode{\sphinxupquote{void SetNames(}}
\begin{quote}

\sphinxAtStartPar
\sphinxcode{\sphinxupquote{const VarArray \&vars,}}

\sphinxAtStartPar
\sphinxcode{\sphinxupquote{const char *szNames,}}

\sphinxAtStartPar
\sphinxcode{\sphinxupquote{size\_t len)}}
\end{quote}
\end{quote}

\sphinxAtStartPar
\sphinxstylestrong{Arguments}
\begin{quote}

\sphinxAtStartPar
\sphinxcode{\sphinxupquote{vars}}: array object of variables.

\sphinxAtStartPar
\sphinxcode{\sphinxupquote{szNames}}: name buffer for variables.

\sphinxAtStartPar
\sphinxcode{\sphinxupquote{len}}: length of name buffer.
\end{quote}
\end{quote}

\subsubsection{Model::SetNames()}
\label{\detokenize{cppapi/model:id86}}\begin{quote}

\sphinxAtStartPar
Set names for given constraints in model.

\sphinxAtStartPar
\sphinxstylestrong{Synopsis}
\begin{quote}

\sphinxAtStartPar
\sphinxcode{\sphinxupquote{void SetNames(}}
\begin{quote}

\sphinxAtStartPar
\sphinxcode{\sphinxupquote{const ConstrArray \&cons,}}

\sphinxAtStartPar
\sphinxcode{\sphinxupquote{const char *szNames,}}

\sphinxAtStartPar
\sphinxcode{\sphinxupquote{size\_t len)}}
\end{quote}
\end{quote}

\sphinxAtStartPar
\sphinxstylestrong{Arguments}
\begin{quote}

\sphinxAtStartPar
\sphinxcode{\sphinxupquote{cons}}: array object of constraints.

\sphinxAtStartPar
\sphinxcode{\sphinxupquote{szNames}}: name buffer for constraints.

\sphinxAtStartPar
\sphinxcode{\sphinxupquote{len}}: length of name buffer.
\end{quote}
\end{quote}

\subsubsection{Model::SetNames()}
\label{\detokenize{cppapi/model:id87}}\begin{quote}

\sphinxAtStartPar
Set names for given gernal constraints in model.

\sphinxAtStartPar
\sphinxstylestrong{Synopsis}
\begin{quote}

\sphinxAtStartPar
\sphinxcode{\sphinxupquote{void SetNames(}}
\begin{quote}

\sphinxAtStartPar
\sphinxcode{\sphinxupquote{const GenConstrArray \&genConstrs,}}

\sphinxAtStartPar
\sphinxcode{\sphinxupquote{const char *szNames,}}

\sphinxAtStartPar
\sphinxcode{\sphinxupquote{size\_t len)}}
\end{quote}
\end{quote}

\sphinxAtStartPar
\sphinxstylestrong{Arguments}
\begin{quote}

\sphinxAtStartPar
\sphinxcode{\sphinxupquote{genConstrs}}: array object of genernal constrains.

\sphinxAtStartPar
\sphinxcode{\sphinxupquote{szNames}}: name buffer for general constraints.

\sphinxAtStartPar
\sphinxcode{\sphinxupquote{len}}: length of name buffer.
\end{quote}
\end{quote}

\subsubsection{Model::SetNames()}
\label{\detokenize{cppapi/model:id88}}\begin{quote}

\sphinxAtStartPar
Set names for given affine cone constraints in model.

\sphinxAtStartPar
\sphinxstylestrong{Synopsis}
\begin{quote}

\sphinxAtStartPar
\sphinxcode{\sphinxupquote{void SetNames(}}
\begin{quote}

\sphinxAtStartPar
\sphinxcode{\sphinxupquote{const AffineConeArray \&cons,}}

\sphinxAtStartPar
\sphinxcode{\sphinxupquote{const char *szNames,}}

\sphinxAtStartPar
\sphinxcode{\sphinxupquote{size\_t len)}}
\end{quote}
\end{quote}

\sphinxAtStartPar
\sphinxstylestrong{Arguments}
\begin{quote}

\sphinxAtStartPar
\sphinxcode{\sphinxupquote{cons}}: array object of affine cone constraints.

\sphinxAtStartPar
\sphinxcode{\sphinxupquote{szNames}}: name buffer for affine cone constraints.

\sphinxAtStartPar
\sphinxcode{\sphinxupquote{len}}: length of name buffer.
\end{quote}
\end{quote}

\subsubsection{Model::SetNames()}
\label{\detokenize{cppapi/model:id89}}\begin{quote}

\sphinxAtStartPar
Set names for given nonlinear constraints in model.

\sphinxAtStartPar
\sphinxstylestrong{Synopsis}
\begin{quote}

\sphinxAtStartPar
\sphinxcode{\sphinxupquote{void SetNames(}}
\begin{quote}

\sphinxAtStartPar
\sphinxcode{\sphinxupquote{const NlConstrArray \&cons,}}

\sphinxAtStartPar
\sphinxcode{\sphinxupquote{const char *szNames,}}

\sphinxAtStartPar
\sphinxcode{\sphinxupquote{size\_t len)}}
\end{quote}
\end{quote}

\sphinxAtStartPar
\sphinxstylestrong{Arguments}
\begin{quote}

\sphinxAtStartPar
\sphinxcode{\sphinxupquote{cons}}: array object of nonlinear constraints.

\sphinxAtStartPar
\sphinxcode{\sphinxupquote{szNames}}: name buffer for nonlinear constraints.

\sphinxAtStartPar
\sphinxcode{\sphinxupquote{len}}: length of name buffer.
\end{quote}
\end{quote}

\subsubsection{Model::SetNames()}
\label{\detokenize{cppapi/model:id90}}\begin{quote}

\sphinxAtStartPar
Set names for given quadratic constraints in model.

\sphinxAtStartPar
\sphinxstylestrong{Synopsis}
\begin{quote}

\sphinxAtStartPar
\sphinxcode{\sphinxupquote{void SetNames(}}
\begin{quote}

\sphinxAtStartPar
\sphinxcode{\sphinxupquote{const QConstrArray \&cons,}}

\sphinxAtStartPar
\sphinxcode{\sphinxupquote{const char *szNames,}}

\sphinxAtStartPar
\sphinxcode{\sphinxupquote{size\_t len)}}
\end{quote}
\end{quote}

\sphinxAtStartPar
\sphinxstylestrong{Arguments}
\begin{quote}

\sphinxAtStartPar
\sphinxcode{\sphinxupquote{cons}}: array object of quadratic constraints.

\sphinxAtStartPar
\sphinxcode{\sphinxupquote{szNames}}: name buffer for quadratic constraints.

\sphinxAtStartPar
\sphinxcode{\sphinxupquote{len}}: length of name buffer.
\end{quote}
\end{quote}

\subsubsection{Model::SetNames()}
\label{\detokenize{cppapi/model:id91}}\begin{quote}

\sphinxAtStartPar
Set names for given PSD variables in model.

\sphinxAtStartPar
\sphinxstylestrong{Synopsis}
\begin{quote}

\sphinxAtStartPar
\sphinxcode{\sphinxupquote{void SetNames(}}
\begin{quote}

\sphinxAtStartPar
\sphinxcode{\sphinxupquote{const PsdVarArray \&vars,}}

\sphinxAtStartPar
\sphinxcode{\sphinxupquote{const char *szNames,}}

\sphinxAtStartPar
\sphinxcode{\sphinxupquote{size\_t len)}}
\end{quote}
\end{quote}

\sphinxAtStartPar
\sphinxstylestrong{Arguments}
\begin{quote}

\sphinxAtStartPar
\sphinxcode{\sphinxupquote{vars}}: array object of PSD variables.

\sphinxAtStartPar
\sphinxcode{\sphinxupquote{szNames}}: name buffer for PSD variables.

\sphinxAtStartPar
\sphinxcode{\sphinxupquote{len}}: length of name buffer.
\end{quote}
\end{quote}

\subsubsection{Model::SetNames()}
\label{\detokenize{cppapi/model:id92}}\begin{quote}

\sphinxAtStartPar
Set names for given PSD constraints in model.

\sphinxAtStartPar
\sphinxstylestrong{Synopsis}
\begin{quote}

\sphinxAtStartPar
\sphinxcode{\sphinxupquote{void SetNames(}}
\begin{quote}

\sphinxAtStartPar
\sphinxcode{\sphinxupquote{const PsdConstrArray \&cons,}}

\sphinxAtStartPar
\sphinxcode{\sphinxupquote{const char *szNames,}}

\sphinxAtStartPar
\sphinxcode{\sphinxupquote{size\_t len)}}
\end{quote}
\end{quote}

\sphinxAtStartPar
\sphinxstylestrong{Arguments}
\begin{quote}

\sphinxAtStartPar
\sphinxcode{\sphinxupquote{cons}}: array object of PSD constraints.

\sphinxAtStartPar
\sphinxcode{\sphinxupquote{szNames}}: name buffer for PSD constraints.

\sphinxAtStartPar
\sphinxcode{\sphinxupquote{len}}: length of name buffer.
\end{quote}
\end{quote}

\subsubsection{Model::SetNames()}
\label{\detokenize{cppapi/model:id93}}\begin{quote}

\sphinxAtStartPar
Set names for given LMI constraints in model.

\sphinxAtStartPar
\sphinxstylestrong{Synopsis}
\begin{quote}

\sphinxAtStartPar
\sphinxcode{\sphinxupquote{void SetNames(}}
\begin{quote}

\sphinxAtStartPar
\sphinxcode{\sphinxupquote{const LmiConstrArray \&cons,}}

\sphinxAtStartPar
\sphinxcode{\sphinxupquote{const char *szNames,}}

\sphinxAtStartPar
\sphinxcode{\sphinxupquote{size\_t len)}}
\end{quote}
\end{quote}

\sphinxAtStartPar
\sphinxstylestrong{Arguments}
\begin{quote}

\sphinxAtStartPar
\sphinxcode{\sphinxupquote{cons}}: array object of LMI constraints.

\sphinxAtStartPar
\sphinxcode{\sphinxupquote{szNames}}: name buffer for LMI constraints.

\sphinxAtStartPar
\sphinxcode{\sphinxupquote{len}}: length of name buffer.
\end{quote}
\end{quote}

\subsubsection{Model::SetNlObjective()}
\label{\detokenize{cppapi/model:model-setnlobjective}}\begin{quote}

\sphinxAtStartPar
Set nonlinear objective for model.

\sphinxAtStartPar
\sphinxstylestrong{Synopsis}
\begin{quote}

\sphinxAtStartPar
\sphinxcode{\sphinxupquote{void SetNlObjective(const NlExpr \&expr, int sense)}}
\end{quote}

\sphinxAtStartPar
\sphinxstylestrong{Arguments}
\begin{quote}

\sphinxAtStartPar
\sphinxcode{\sphinxupquote{expr}}: nonlinear expression of the objective.

\sphinxAtStartPar
\sphinxcode{\sphinxupquote{sense}}: optimization sense. optional, default value 0 does not change COPT sense.
\end{quote}
\end{quote}

\subsubsection{Model::SetNlPrimalStart()}
\label{\detokenize{cppapi/model:model-setnlprimalstart}}\begin{quote}

\sphinxAtStartPar
Given count, set initial values for variables of NLP from beginning.

\sphinxAtStartPar
\sphinxstylestrong{Synopsis}
\begin{quote}

\sphinxAtStartPar
\sphinxcode{\sphinxupquote{void SetNlPrimalStart(int count, double *pVals)}}
\end{quote}

\sphinxAtStartPar
\sphinxstylestrong{Arguments}
\begin{quote}

\sphinxAtStartPar
\sphinxcode{\sphinxupquote{count}}: the number of variables to set.

\sphinxAtStartPar
\sphinxcode{\sphinxupquote{pVals}}: pointer to initial values of variables.
\end{quote}
\end{quote}

\subsubsection{Model::SetNlPrimalStart()}
\label{\detokenize{cppapi/model:id94}}\begin{quote}

\sphinxAtStartPar
Set initial value for the specified variable of NLP.

\sphinxAtStartPar
\sphinxstylestrong{Synopsis}
\begin{quote}

\sphinxAtStartPar
\sphinxcode{\sphinxupquote{void SetNlPrimalStart(const Var \&var, double val)}}
\end{quote}

\sphinxAtStartPar
\sphinxstylestrong{Arguments}
\begin{quote}

\sphinxAtStartPar
\sphinxcode{\sphinxupquote{var}}: an interested variable.

\sphinxAtStartPar
\sphinxcode{\sphinxupquote{val}}: initial value of the variable.
\end{quote}
\end{quote}

\subsubsection{Model::SetNlPrimalStart()}
\label{\detokenize{cppapi/model:id95}}\begin{quote}

\sphinxAtStartPar
Set initial values for an array of variables of NLP.

\sphinxAtStartPar
\sphinxstylestrong{Synopsis}
\begin{quote}

\sphinxAtStartPar
\sphinxcode{\sphinxupquote{void SetNlPrimalStart(const VarArray \&vars, double *pVals)}}
\end{quote}

\sphinxAtStartPar
\sphinxstylestrong{Arguments}
\begin{quote}

\sphinxAtStartPar
\sphinxcode{\sphinxupquote{vars}}: a list of interested variables.

\sphinxAtStartPar
\sphinxcode{\sphinxupquote{pVals}}: pointer to initial values of variables.
\end{quote}
\end{quote}

\subsubsection{Model::SetObjConst()}
\label{\detokenize{cppapi/model:model-setobjconst}}\begin{quote}

\sphinxAtStartPar
Set objective constant.

\sphinxAtStartPar
\sphinxstylestrong{Synopsis}
\begin{quote}

\sphinxAtStartPar
\sphinxcode{\sphinxupquote{void SetObjConst(double constant)}}
\end{quote}

\sphinxAtStartPar
\sphinxstylestrong{Arguments}
\begin{quote}

\sphinxAtStartPar
\sphinxcode{\sphinxupquote{constant}}: constant value to set.
\end{quote}
\end{quote}

\subsubsection{Model::SetObjective()}
\label{\detokenize{cppapi/model:model-setobjective}}\begin{quote}

\sphinxAtStartPar
Set objective for model.

\sphinxAtStartPar
\sphinxstylestrong{Synopsis}
\begin{quote}

\sphinxAtStartPar
\sphinxcode{\sphinxupquote{void SetObjective(const Expr \&expr, int sense)}}
\end{quote}

\sphinxAtStartPar
\sphinxstylestrong{Arguments}
\begin{quote}

\sphinxAtStartPar
\sphinxcode{\sphinxupquote{expr}}: expression of the objective.

\sphinxAtStartPar
\sphinxcode{\sphinxupquote{sense}}: optimization sense. optional, default value 0 does not change COPT sense.
\end{quote}
\end{quote}

\subsubsection{Model::SetObjectiveN()}
\label{\detokenize{cppapi/model:model-setobjectiven}}\begin{quote}

\sphinxAtStartPar
Set a multi\sphinxhyphen{}objective function in model.

\sphinxAtStartPar
\sphinxstylestrong{Synopsis}
\begin{quote}

\sphinxAtStartPar
\sphinxcode{\sphinxupquote{void SetObjectiveN(}}
\begin{quote}

\sphinxAtStartPar
\sphinxcode{\sphinxupquote{int idx,}}

\sphinxAtStartPar
\sphinxcode{\sphinxupquote{const Expr \&expr,}}

\sphinxAtStartPar
\sphinxcode{\sphinxupquote{int sense,}}

\sphinxAtStartPar
\sphinxcode{\sphinxupquote{double priority,}}

\sphinxAtStartPar
\sphinxcode{\sphinxupquote{double weight,}}

\sphinxAtStartPar
\sphinxcode{\sphinxupquote{double abstol,}}

\sphinxAtStartPar
\sphinxcode{\sphinxupquote{double reltol)}}
\end{quote}
\end{quote}

\sphinxAtStartPar
\sphinxstylestrong{Arguments}
\begin{quote}

\sphinxAtStartPar
\sphinxcode{\sphinxupquote{idx}}: index of a multi\sphinxhyphen{}objective function.

\sphinxAtStartPar
\sphinxcode{\sphinxupquote{expr}}: linear expression of the multi\sphinxhyphen{}objective function.

\sphinxAtStartPar
\sphinxcode{\sphinxupquote{sense}}: optimization sense. optional, default value 0 does not change COPT sense.

\sphinxAtStartPar
\sphinxcode{\sphinxupquote{priority}}: an optional objective parameter. Default value is 0.0.

\sphinxAtStartPar
\sphinxcode{\sphinxupquote{weight}}: an optional objective parameter. Default value is 1.0.

\sphinxAtStartPar
\sphinxcode{\sphinxupquote{abstol}}: absolute tolerance is an optional objective parameter. Default value is 1e\sphinxhyphen{}6.

\sphinxAtStartPar
\sphinxcode{\sphinxupquote{reltol}}: relative tolerance is an optional objective parameter. Default value is 0.0.
\end{quote}
\end{quote}

\subsubsection{Model::SetObjParamN()}
\label{\detokenize{cppapi/model:model-setobjparamn}}\begin{quote}

\sphinxAtStartPar
Set value of objective parameter of a multi\sphinxhyphen{}objective function.

\sphinxAtStartPar
\sphinxstylestrong{Synopsis}
\begin{quote}

\sphinxAtStartPar
\sphinxcode{\sphinxupquote{void SetObjParamN(}}
\begin{quote}

\sphinxAtStartPar
\sphinxcode{\sphinxupquote{int idx,}}

\sphinxAtStartPar
\sphinxcode{\sphinxupquote{const char *szParam,}}

\sphinxAtStartPar
\sphinxcode{\sphinxupquote{double val)}}
\end{quote}
\end{quote}

\sphinxAtStartPar
\sphinxstylestrong{Arguments}
\begin{quote}

\sphinxAtStartPar
\sphinxcode{\sphinxupquote{idx}}: index of a multi\sphinxhyphen{}objective function.

\sphinxAtStartPar
\sphinxcode{\sphinxupquote{szParam}}: name of objective parameter, including priority, weight, abstol and reltol.

\sphinxAtStartPar
\sphinxcode{\sphinxupquote{val}}: new value of objective parameter.
\end{quote}
\end{quote}

\subsubsection{Model::SetObjSense()}
\label{\detokenize{cppapi/model:model-setobjsense}}\begin{quote}

\sphinxAtStartPar
Set objective sense for model.

\sphinxAtStartPar
\sphinxstylestrong{Synopsis}
\begin{quote}

\sphinxAtStartPar
\sphinxcode{\sphinxupquote{void SetObjSense(int sense)}}
\end{quote}

\sphinxAtStartPar
\sphinxstylestrong{Arguments}
\begin{quote}

\sphinxAtStartPar
\sphinxcode{\sphinxupquote{sense}}: the objective sense.
\end{quote}
\end{quote}

\subsubsection{Model::SetPsdCoeff()}
\label{\detokenize{cppapi/model:model-setpsdcoeff}}\begin{quote}

\sphinxAtStartPar
Set the coefficient matrix of a PSD variable in a PSD constraint.

\sphinxAtStartPar
\sphinxstylestrong{Synopsis}
\begin{quote}

\sphinxAtStartPar
\sphinxcode{\sphinxupquote{void SetPsdCoeff(}}
\begin{quote}

\sphinxAtStartPar
\sphinxcode{\sphinxupquote{const PsdConstraint \&constr,}}

\sphinxAtStartPar
\sphinxcode{\sphinxupquote{const PsdVar \&var,}}

\sphinxAtStartPar
\sphinxcode{\sphinxupquote{const SymMatrix \&mat)}}
\end{quote}
\end{quote}

\sphinxAtStartPar
\sphinxstylestrong{Arguments}
\begin{quote}

\sphinxAtStartPar
\sphinxcode{\sphinxupquote{constr}}: The desired PSD constraint.

\sphinxAtStartPar
\sphinxcode{\sphinxupquote{var}}: The desired PSD variable.

\sphinxAtStartPar
\sphinxcode{\sphinxupquote{mat}}: new coefficient matrix.
\end{quote}
\end{quote}

\subsubsection{Model::SetPsdObjective()}
\label{\detokenize{cppapi/model:model-setpsdobjective}}\begin{quote}

\sphinxAtStartPar
Set PSD objective for model.

\sphinxAtStartPar
\sphinxstylestrong{Synopsis}
\begin{quote}

\sphinxAtStartPar
\sphinxcode{\sphinxupquote{void SetPsdObjective(const PsdExpr \&expr, int sense)}}
\end{quote}

\sphinxAtStartPar
\sphinxstylestrong{Arguments}
\begin{quote}

\sphinxAtStartPar
\sphinxcode{\sphinxupquote{expr}}: PSD expression of the objective.

\sphinxAtStartPar
\sphinxcode{\sphinxupquote{sense}}: optimization sense. optional, default value 0 does not change COPT sense.
\end{quote}
\end{quote}

\subsubsection{Model::SetQuadCost()}
\label{\detokenize{cppapi/model:model-setquadcost}}\begin{quote}

\sphinxAtStartPar
Set quadratic objective of the model.

\sphinxAtStartPar
\sphinxstylestrong{Synopsis}
\begin{quote}

\sphinxAtStartPar
\sphinxcode{\sphinxupquote{void SetQuadCost(}}
\begin{quote}

\sphinxAtStartPar
\sphinxcode{\sphinxupquote{int nQElem,}}

\sphinxAtStartPar
\sphinxcode{\sphinxupquote{const int *pMatRow,}}

\sphinxAtStartPar
\sphinxcode{\sphinxupquote{const int *pMatCol,}}

\sphinxAtStartPar
\sphinxcode{\sphinxupquote{const double *pMatElem)}}
\end{quote}
\end{quote}

\sphinxAtStartPar
\sphinxstylestrong{Arguments}
\begin{quote}

\sphinxAtStartPar
\sphinxcode{\sphinxupquote{nQElem}}: number of Q objective elements.

\sphinxAtStartPar
\sphinxcode{\sphinxupquote{pMatRow}}: row indices of Q objective elements.

\sphinxAtStartPar
\sphinxcode{\sphinxupquote{pMatCol}}: column indices of Q objective elements.

\sphinxAtStartPar
\sphinxcode{\sphinxupquote{pMatElem}}: nonzero elements of Q objective.
\end{quote}
\end{quote}

\subsubsection{Model::SetQuadObjective()}
\label{\detokenize{cppapi/model:model-setquadobjective}}\begin{quote}

\sphinxAtStartPar
Set quadratic objective for model.

\sphinxAtStartPar
\sphinxstylestrong{Synopsis}
\begin{quote}

\sphinxAtStartPar
\sphinxcode{\sphinxupquote{void SetQuadObjective(const QuadExpr \&expr, int sense)}}
\end{quote}

\sphinxAtStartPar
\sphinxstylestrong{Arguments}
\begin{quote}

\sphinxAtStartPar
\sphinxcode{\sphinxupquote{expr}}: quadratic expression of the objective.

\sphinxAtStartPar
\sphinxcode{\sphinxupquote{sense}}: optimization sense. optional, default value 0 does not change COPT sense.
\end{quote}
\end{quote}

\subsubsection{Model::SetSlackBasis()}
\label{\detokenize{cppapi/model:model-setslackbasis}}\begin{quote}

\sphinxAtStartPar
Set slack basis to model.

\sphinxAtStartPar
\sphinxstylestrong{Synopsis}
\begin{quote}

\sphinxAtStartPar
\sphinxcode{\sphinxupquote{void SetSlackBasis()}}
\end{quote}
\end{quote}

\subsubsection{Model::SetSolverLogCallback()}
\label{\detokenize{cppapi/model:model-setsolverlogcallback}}\begin{quote}

\sphinxAtStartPar
Set log callback for COPT.

\sphinxAtStartPar
\sphinxstylestrong{Synopsis}
\begin{quote}

\sphinxAtStartPar
\sphinxcode{\sphinxupquote{void SetSolverLogCallback(ILogCallback *pcb)}}
\end{quote}

\sphinxAtStartPar
\sphinxstylestrong{Arguments}
\begin{quote}

\sphinxAtStartPar
\sphinxcode{\sphinxupquote{pcb}}: pointer to ILogCallback object.
\end{quote}
\end{quote}

\subsubsection{Model::SetSolverLogFile()}
\label{\detokenize{cppapi/model:model-setsolverlogfile}}\begin{quote}

\sphinxAtStartPar
Set log file for COPT.

\sphinxAtStartPar
\sphinxstylestrong{Synopsis}
\begin{quote}

\sphinxAtStartPar
\sphinxcode{\sphinxupquote{void SetSolverLogFile(const char *szLogFile)}}
\end{quote}

\sphinxAtStartPar
\sphinxstylestrong{Arguments}
\begin{quote}

\sphinxAtStartPar
\sphinxcode{\sphinxupquote{szLogFile}}: log file name.
\end{quote}
\end{quote}

\subsubsection{Model::Solve()}
\label{\detokenize{cppapi/model:model-solve}}\begin{quote}

\sphinxAtStartPar
Solve the model as MIP.

\sphinxAtStartPar
\sphinxstylestrong{Synopsis}
\begin{quote}

\sphinxAtStartPar
\sphinxcode{\sphinxupquote{void Solve()}}
\end{quote}
\end{quote}

\subsubsection{Model::SolveLp()}
\label{\detokenize{cppapi/model:model-solvelp}}\begin{quote}

\sphinxAtStartPar
Solve the model as LP.

\sphinxAtStartPar
\sphinxstylestrong{Synopsis}
\begin{quote}

\sphinxAtStartPar
\sphinxcode{\sphinxupquote{void SolveLp()}}
\end{quote}
\end{quote}

\subsubsection{Model::Tune()}
\label{\detokenize{cppapi/model:model-tune}}\begin{quote}

\sphinxAtStartPar
Tune model.

\sphinxAtStartPar
\sphinxstylestrong{Synopsis}
\begin{quote}

\sphinxAtStartPar
\sphinxcode{\sphinxupquote{void Tune()}}
\end{quote}
\end{quote}

\subsubsection{Model::Write()}
\label{\detokenize{cppapi/model:model-write}}\begin{quote}

\sphinxAtStartPar
Output problem, solution, basis, MIP start or modified COPT parameters to file.

\sphinxAtStartPar
\sphinxstylestrong{Synopsis}
\begin{quote}

\sphinxAtStartPar
\sphinxcode{\sphinxupquote{void Write(const char *szFileName)}}
\end{quote}

\sphinxAtStartPar
\sphinxstylestrong{Arguments}
\begin{quote}

\sphinxAtStartPar
\sphinxcode{\sphinxupquote{szFileName}}: an output file name.
\end{quote}
\end{quote}

\subsubsection{Model::WriteBasis()}
\label{\detokenize{cppapi/model:model-writebasis}}\begin{quote}

\sphinxAtStartPar
Output optimal basis to a file of type ‘.bas’.

\sphinxAtStartPar
\sphinxstylestrong{Synopsis}
\begin{quote}

\sphinxAtStartPar
\sphinxcode{\sphinxupquote{void WriteBasis(const char *szFileName)}}
\end{quote}

\sphinxAtStartPar
\sphinxstylestrong{Arguments}
\begin{quote}

\sphinxAtStartPar
\sphinxcode{\sphinxupquote{szFileName}}: an output file name.
\end{quote}
\end{quote}

\subsubsection{Model::WriteBin()}
\label{\detokenize{cppapi/model:model-writebin}}\begin{quote}

\sphinxAtStartPar
Output problem to a file as COPT binary format.

\sphinxAtStartPar
\sphinxstylestrong{Synopsis}
\begin{quote}

\sphinxAtStartPar
\sphinxcode{\sphinxupquote{void WriteBin(const char *szFileName)}}
\end{quote}

\sphinxAtStartPar
\sphinxstylestrong{Arguments}
\begin{quote}

\sphinxAtStartPar
\sphinxcode{\sphinxupquote{szFileName}}: an output file name.
\end{quote}
\end{quote}

\subsubsection{Model::WriteIIS()}
\label{\detokenize{cppapi/model:model-writeiis}}\begin{quote}

\sphinxAtStartPar
Output IIS to file.

\sphinxAtStartPar
\sphinxstylestrong{Synopsis}
\begin{quote}

\sphinxAtStartPar
\sphinxcode{\sphinxupquote{void WriteIIS(const char *szFileName)}}
\end{quote}

\sphinxAtStartPar
\sphinxstylestrong{Arguments}
\begin{quote}

\sphinxAtStartPar
\sphinxcode{\sphinxupquote{szFileName}}: Output file name.
\end{quote}
\end{quote}

\subsubsection{Model::WriteJsonSol()}
\label{\detokenize{cppapi/model:model-writejsonsol}}\begin{quote}

\sphinxAtStartPar
Output solution to a file of type ‘.json’.

\sphinxAtStartPar
\sphinxstylestrong{Synopsis}
\begin{quote}

\sphinxAtStartPar
\sphinxcode{\sphinxupquote{void WriteJsonSol(const char *szFileName)}}
\end{quote}

\sphinxAtStartPar
\sphinxstylestrong{Arguments}
\begin{quote}

\sphinxAtStartPar
\sphinxcode{\sphinxupquote{szFileName}}: an output file name.
\end{quote}
\end{quote}

\subsubsection{Model::WriteLp()}
\label{\detokenize{cppapi/model:model-writelp}}\begin{quote}

\sphinxAtStartPar
Output problem to a file as LP format.

\sphinxAtStartPar
\sphinxstylestrong{Synopsis}
\begin{quote}

\sphinxAtStartPar
\sphinxcode{\sphinxupquote{void WriteLp(const char *szFileName)}}
\end{quote}

\sphinxAtStartPar
\sphinxstylestrong{Arguments}
\begin{quote}

\sphinxAtStartPar
\sphinxcode{\sphinxupquote{szFileName}}: an output file name.
\end{quote}
\end{quote}

\subsubsection{Model::WriteMps()}
\label{\detokenize{cppapi/model:model-writemps}}\begin{quote}

\sphinxAtStartPar
Output problem to a file as MPS format.

\sphinxAtStartPar
\sphinxstylestrong{Synopsis}
\begin{quote}

\sphinxAtStartPar
\sphinxcode{\sphinxupquote{void WriteMps(const char *szFileName)}}
\end{quote}

\sphinxAtStartPar
\sphinxstylestrong{Arguments}
\begin{quote}

\sphinxAtStartPar
\sphinxcode{\sphinxupquote{szFileName}}: an output file name.
\end{quote}
\end{quote}

\subsubsection{Model::WriteMpsStr()}
\label{\detokenize{cppapi/model:model-writempsstr}}\begin{quote}

\sphinxAtStartPar
Output problem to a buffer as MPS format.

\sphinxAtStartPar
\sphinxstylestrong{Synopsis}
\begin{quote}

\sphinxAtStartPar
\sphinxcode{\sphinxupquote{ProbBuffer WriteMpsStr()}}
\end{quote}

\sphinxAtStartPar
\sphinxstylestrong{Return}
\begin{quote}

\sphinxAtStartPar
output problem buffer.
\end{quote}
\end{quote}

\subsubsection{Model::WriteMst()}
\label{\detokenize{cppapi/model:model-writemst}}\begin{quote}

\sphinxAtStartPar
Output MIP start information to a file of type ‘.mst’.

\sphinxAtStartPar
\sphinxstylestrong{Synopsis}
\begin{quote}

\sphinxAtStartPar
\sphinxcode{\sphinxupquote{void WriteMst(const char *szFileName)}}
\end{quote}

\sphinxAtStartPar
\sphinxstylestrong{Arguments}
\begin{quote}

\sphinxAtStartPar
\sphinxcode{\sphinxupquote{szFileName}}: an output file name.
\end{quote}
\end{quote}

\subsubsection{Model::WriteNL()}
\label{\detokenize{cppapi/model:model-writenl}}\begin{quote}

\sphinxAtStartPar
Output problem to a file as NL format.

\sphinxAtStartPar
\sphinxstylestrong{Synopsis}
\begin{quote}

\sphinxAtStartPar
\sphinxcode{\sphinxupquote{void WriteNL(const char *szFileName)}}
\end{quote}

\sphinxAtStartPar
\sphinxstylestrong{Arguments}
\begin{quote}

\sphinxAtStartPar
\sphinxcode{\sphinxupquote{szFileName}}: an output file name.
\end{quote}
\end{quote}

\subsubsection{Model::WriteOrd()}
\label{\detokenize{cppapi/model:model-writeord}}\begin{quote}

\sphinxAtStartPar
Output branching order to file.

\sphinxAtStartPar
\sphinxstylestrong{Synopsis}
\begin{quote}

\sphinxAtStartPar
\sphinxcode{\sphinxupquote{void WriteOrd(const char *szFileName)}}
\end{quote}

\sphinxAtStartPar
\sphinxstylestrong{Arguments}
\begin{quote}

\sphinxAtStartPar
\sphinxcode{\sphinxupquote{szFileName}}: Output file name.
\end{quote}
\end{quote}

\subsubsection{Model::WriteParam()}
\label{\detokenize{cppapi/model:model-writeparam}}\begin{quote}

\sphinxAtStartPar
Output modified COPT parameters to a file of type ‘.par’.

\sphinxAtStartPar
\sphinxstylestrong{Synopsis}
\begin{quote}

\sphinxAtStartPar
\sphinxcode{\sphinxupquote{void WriteParam(const char *szFileName)}}
\end{quote}

\sphinxAtStartPar
\sphinxstylestrong{Arguments}
\begin{quote}

\sphinxAtStartPar
\sphinxcode{\sphinxupquote{szFileName}}: an output file name.
\end{quote}
\end{quote}

\subsubsection{Model::WritePoolSol()}
\label{\detokenize{cppapi/model:model-writepoolsol}}\begin{quote}

\sphinxAtStartPar
Output selected pool solution to a file of type ‘.sol’.

\sphinxAtStartPar
\sphinxstylestrong{Synopsis}
\begin{quote}

\sphinxAtStartPar
\sphinxcode{\sphinxupquote{void WritePoolSol(int iSol, const char *szFileName)}}
\end{quote}

\sphinxAtStartPar
\sphinxstylestrong{Arguments}
\begin{quote}

\sphinxAtStartPar
\sphinxcode{\sphinxupquote{iSol}}: index of pool solution.

\sphinxAtStartPar
\sphinxcode{\sphinxupquote{szFileName}}: an output file name.
\end{quote}
\end{quote}

\subsubsection{Model::WriteRelax()}
\label{\detokenize{cppapi/model:model-writerelax}}\begin{quote}

\sphinxAtStartPar
Output feasibility relaxation problem to file.

\sphinxAtStartPar
\sphinxstylestrong{Synopsis}
\begin{quote}

\sphinxAtStartPar
\sphinxcode{\sphinxupquote{void WriteRelax(const char *szFileName)}}
\end{quote}

\sphinxAtStartPar
\sphinxstylestrong{Arguments}
\begin{quote}

\sphinxAtStartPar
\sphinxcode{\sphinxupquote{szFileName}}: Output file name.
\end{quote}
\end{quote}

\subsubsection{Model::WriteSol()}
\label{\detokenize{cppapi/model:model-writesol}}\begin{quote}

\sphinxAtStartPar
Output solution to a file of type ‘.sol’.

\sphinxAtStartPar
\sphinxstylestrong{Synopsis}
\begin{quote}

\sphinxAtStartPar
\sphinxcode{\sphinxupquote{void WriteSol(const char *szFileName)}}
\end{quote}

\sphinxAtStartPar
\sphinxstylestrong{Arguments}
\begin{quote}

\sphinxAtStartPar
\sphinxcode{\sphinxupquote{szFileName}}: an output file name.
\end{quote}
\end{quote}

\subsubsection{Model::WriteTuneParam()}
\label{\detokenize{cppapi/model:model-writetuneparam}}\begin{quote}

\sphinxAtStartPar
Output specified tuned parameters to a file of type ‘.par’.

\sphinxAtStartPar
\sphinxstylestrong{Synopsis}
\begin{quote}

\sphinxAtStartPar
\sphinxcode{\sphinxupquote{void WriteTuneParam(int idx, const char *szFileName)}}
\end{quote}

\sphinxAtStartPar
\sphinxstylestrong{Arguments}
\begin{quote}

\sphinxAtStartPar
\sphinxcode{\sphinxupquote{idx}}: Index of tuned parameters.

\sphinxAtStartPar
\sphinxcode{\sphinxupquote{szFileName}}: Output file name.
\end{quote}
\end{quote}

\subsection{Var}
\label{\detokenize{cppapiref:var}}\label{\detokenize{cppapiref:chapcppapiref-var}}
\sphinxAtStartPar
COPT variable object. Variables are always associated with a particular model.
User creates a variable object by adding a variable to a model, rather than
by using constructor of Var class.

\sphinxstepscope

\subsubsection{Var::Get()}
\label{\detokenize{cppapi/var:var-get}}\label{\detokenize{cppapi/var::doc}}\begin{quote}

\sphinxAtStartPar
Get information of the variable. Support informations of “Value”, “RedCost”, “PrimalRay”, “LB”, “UB”, “Obj” and “BranchFactor”.

\sphinxAtStartPar
\sphinxstylestrong{Synopsis}
\begin{quote}

\sphinxAtStartPar
\sphinxcode{\sphinxupquote{double Get(const char *szInfo)}}
\end{quote}

\sphinxAtStartPar
\sphinxstylestrong{Arguments}
\begin{quote}

\sphinxAtStartPar
\sphinxcode{\sphinxupquote{szInfo}}: name of information.
\end{quote}

\sphinxAtStartPar
\sphinxstylestrong{Return}
\begin{quote}

\sphinxAtStartPar
value of information
\end{quote}
\end{quote}

\subsubsection{Var::GetBasis()}
\label{\detokenize{cppapi/var:var-getbasis}}\begin{quote}

\sphinxAtStartPar
Get basis status of this variable.

\sphinxAtStartPar
\sphinxstylestrong{Synopsis}
\begin{quote}

\sphinxAtStartPar
\sphinxcode{\sphinxupquote{int GetBasis()}}
\end{quote}

\sphinxAtStartPar
\sphinxstylestrong{Return}
\begin{quote}

\sphinxAtStartPar
basis status.
\end{quote}
\end{quote}

\subsubsection{Var::GetIdx()}
\label{\detokenize{cppapi/var:var-getidx}}\begin{quote}

\sphinxAtStartPar
Get index of the variable.

\sphinxAtStartPar
\sphinxstylestrong{Synopsis}
\begin{quote}

\sphinxAtStartPar
\sphinxcode{\sphinxupquote{inline int GetIdx()}}
\end{quote}

\sphinxAtStartPar
\sphinxstylestrong{Return}
\begin{quote}

\sphinxAtStartPar
variable index.
\end{quote}
\end{quote}

\subsubsection{Var::GetLowerIIS()}
\label{\detokenize{cppapi/var:var-getloweriis}}\begin{quote}

\sphinxAtStartPar
Get IIS status for lower bound of the variable.

\sphinxAtStartPar
\sphinxstylestrong{Synopsis}
\begin{quote}

\sphinxAtStartPar
\sphinxcode{\sphinxupquote{int GetLowerIIS()}}
\end{quote}

\sphinxAtStartPar
\sphinxstylestrong{Return}
\begin{quote}

\sphinxAtStartPar
IIS status.
\end{quote}
\end{quote}

\subsubsection{Var::GetName()}
\label{\detokenize{cppapi/var:var-getname}}\begin{quote}

\sphinxAtStartPar
Get name of the variable.

\sphinxAtStartPar
\sphinxstylestrong{Synopsis}
\begin{quote}

\sphinxAtStartPar
\sphinxcode{\sphinxupquote{const char *GetName()}}
\end{quote}

\sphinxAtStartPar
\sphinxstylestrong{Return}
\begin{quote}

\sphinxAtStartPar
variable name.
\end{quote}
\end{quote}

\subsubsection{Var::GetType()}
\label{\detokenize{cppapi/var:var-gettype}}\begin{quote}

\sphinxAtStartPar
Get type of the variable.

\sphinxAtStartPar
\sphinxstylestrong{Synopsis}
\begin{quote}

\sphinxAtStartPar
\sphinxcode{\sphinxupquote{char GetType()}}
\end{quote}

\sphinxAtStartPar
\sphinxstylestrong{Return}
\begin{quote}

\sphinxAtStartPar
variable type.
\end{quote}
\end{quote}

\subsubsection{Var::GetUpperIIS()}
\label{\detokenize{cppapi/var:var-getupperiis}}\begin{quote}

\sphinxAtStartPar
Get IIS status for upper bound of the variable.

\sphinxAtStartPar
\sphinxstylestrong{Synopsis}
\begin{quote}

\sphinxAtStartPar
\sphinxcode{\sphinxupquote{int GetUpperIIS()}}
\end{quote}

\sphinxAtStartPar
\sphinxstylestrong{Return}
\begin{quote}

\sphinxAtStartPar
IIS status.
\end{quote}
\end{quote}

\subsubsection{Var::Remove()}
\label{\detokenize{cppapi/var:var-remove}}\begin{quote}

\sphinxAtStartPar
Remove variable from model.

\sphinxAtStartPar
\sphinxstylestrong{Synopsis}
\begin{quote}

\sphinxAtStartPar
\sphinxcode{\sphinxupquote{void Remove()}}
\end{quote}
\end{quote}

\subsubsection{Var::Set()}
\label{\detokenize{cppapi/var:var-set}}\begin{quote}

\sphinxAtStartPar
Set information value of the variable. Support informations of “LB”, “UB”, “Obj” and “BranchFactor”.

\sphinxAtStartPar
\sphinxstylestrong{Synopsis}
\begin{quote}

\sphinxAtStartPar
\sphinxcode{\sphinxupquote{void Set(const char *szInfo, double value)}}
\end{quote}

\sphinxAtStartPar
\sphinxstylestrong{Arguments}
\begin{quote}

\sphinxAtStartPar
\sphinxcode{\sphinxupquote{szInfo}}: name of information.

\sphinxAtStartPar
\sphinxcode{\sphinxupquote{value}}: new information value.
\end{quote}
\end{quote}

\subsubsection{Var::SetName()}
\label{\detokenize{cppapi/var:var-setname}}\begin{quote}

\sphinxAtStartPar
Set name of the variable.

\sphinxAtStartPar
\sphinxstylestrong{Synopsis}
\begin{quote}

\sphinxAtStartPar
\sphinxcode{\sphinxupquote{void SetName(const char *szName)}}
\end{quote}

\sphinxAtStartPar
\sphinxstylestrong{Arguments}
\begin{quote}

\sphinxAtStartPar
\sphinxcode{\sphinxupquote{szName}}: variable name.
\end{quote}
\end{quote}

\subsubsection{Var::SetType()}
\label{\detokenize{cppapi/var:var-settype}}\begin{quote}

\sphinxAtStartPar
Set type of the variable.

\sphinxAtStartPar
\sphinxstylestrong{Synopsis}
\begin{quote}

\sphinxAtStartPar
\sphinxcode{\sphinxupquote{void SetType(char type)}}
\end{quote}

\sphinxAtStartPar
\sphinxstylestrong{Arguments}
\begin{quote}

\sphinxAtStartPar
\sphinxcode{\sphinxupquote{type}}: variable type.
\end{quote}
\end{quote}

\subsection{VarArray}
\label{\detokenize{cppapiref:vararray}}\label{\detokenize{cppapiref:chapcppapiref-vararray}}
\sphinxAtStartPar
COPT variable array object. To store and access a set of C++
{\hyperref[\detokenize{cppapiref:chapcppapiref-var}]{\sphinxcrossref{\DUrole{std,std-ref}{Var}}}} objects, Cardinal Optimizer provides C++ VarArray
class, which defines the following methods.

\sphinxstepscope

\subsubsection{VarArray::GetVar()}
\label{\detokenize{cppapi/vararray:vararray-getvar}}\label{\detokenize{cppapi/vararray::doc}}\begin{quote}

\sphinxAtStartPar
Get i\sphinxhyphen{}th variable object.

\sphinxAtStartPar
\sphinxstylestrong{Synopsis}
\begin{quote}

\sphinxAtStartPar
\sphinxcode{\sphinxupquote{Var \&GetVar(int i)}}
\end{quote}

\sphinxAtStartPar
\sphinxstylestrong{Arguments}
\begin{quote}

\sphinxAtStartPar
\sphinxcode{\sphinxupquote{i}}: index of the variable.
\end{quote}

\sphinxAtStartPar
\sphinxstylestrong{Return}
\begin{quote}

\sphinxAtStartPar
variable object with index i.
\end{quote}
\end{quote}

\subsubsection{VarArray::PushBack()}
\label{\detokenize{cppapi/vararray:vararray-pushback}}\begin{quote}

\sphinxAtStartPar
Add a variable object to variable array.

\sphinxAtStartPar
\sphinxstylestrong{Synopsis}
\begin{quote}

\sphinxAtStartPar
\sphinxcode{\sphinxupquote{void PushBack(const Var \&var)}}
\end{quote}

\sphinxAtStartPar
\sphinxstylestrong{Arguments}
\begin{quote}

\sphinxAtStartPar
\sphinxcode{\sphinxupquote{var}}: a variable object.
\end{quote}
\end{quote}

\subsubsection{VarArray::Reserve()}
\label{\detokenize{cppapi/vararray:vararray-reserve}}\begin{quote}

\sphinxAtStartPar
Reserve capacity to contain at least n items.

\sphinxAtStartPar
\sphinxstylestrong{Synopsis}
\begin{quote}

\sphinxAtStartPar
\sphinxcode{\sphinxupquote{void Reserve(int n)}}
\end{quote}

\sphinxAtStartPar
\sphinxstylestrong{Arguments}
\begin{quote}

\sphinxAtStartPar
\sphinxcode{\sphinxupquote{n}}: minimum capacity for variable object.
\end{quote}
\end{quote}

\subsubsection{VarArray::Size()}
\label{\detokenize{cppapi/vararray:vararray-size}}\begin{quote}

\sphinxAtStartPar
Get the number of variable objects.

\sphinxAtStartPar
\sphinxstylestrong{Synopsis}
\begin{quote}

\sphinxAtStartPar
\sphinxcode{\sphinxupquote{int Size()}}
\end{quote}

\sphinxAtStartPar
\sphinxstylestrong{Return}
\begin{quote}

\sphinxAtStartPar
number of variable objects.
\end{quote}
\end{quote}

\subsection{Expr}
\label{\detokenize{cppapiref:expr}}\label{\detokenize{cppapiref:chapcppapiref-expr}}
\sphinxAtStartPar
COPT linear expression object. A linear expression consists of a constant
term, a list of terms of variables and associated coefficients. Linear
expressions are used to build constraints.

\sphinxstepscope

\subsubsection{Expr::Expr()}
\label{\detokenize{cppapi/expr:expr-expr}}\label{\detokenize{cppapi/expr::doc}}\begin{quote}

\sphinxAtStartPar
Constructor of a constant linear expression.

\sphinxAtStartPar
\sphinxstylestrong{Synopsis}
\begin{quote}

\sphinxAtStartPar
\sphinxcode{\sphinxupquote{Expr(double constant)}}
\end{quote}

\sphinxAtStartPar
\sphinxstylestrong{Arguments}
\begin{quote}

\sphinxAtStartPar
\sphinxcode{\sphinxupquote{constant}}: constant value in expression object.
\end{quote}
\end{quote}

\subsubsection{Expr::Expr()}
\label{\detokenize{cppapi/expr:id1}}\begin{quote}

\sphinxAtStartPar
Constructor of a linear expression with one term.

\sphinxAtStartPar
\sphinxstylestrong{Synopsis}
\begin{quote}

\sphinxAtStartPar
\sphinxcode{\sphinxupquote{Expr(const Var \&var, double coeff)}}
\end{quote}

\sphinxAtStartPar
\sphinxstylestrong{Arguments}
\begin{quote}

\sphinxAtStartPar
\sphinxcode{\sphinxupquote{var}}: variable for the added term.

\sphinxAtStartPar
\sphinxcode{\sphinxupquote{coeff}}: coefficent for the added term.
\end{quote}
\end{quote}

\subsubsection{Expr::AddConstant()}
\label{\detokenize{cppapi/expr:expr-addconstant}}\begin{quote}

\sphinxAtStartPar
Add constant for the expression.

\sphinxAtStartPar
\sphinxstylestrong{Synopsis}
\begin{quote}

\sphinxAtStartPar
\sphinxcode{\sphinxupquote{void AddConstant(double constant)}}
\end{quote}

\sphinxAtStartPar
\sphinxstylestrong{Arguments}
\begin{quote}

\sphinxAtStartPar
\sphinxcode{\sphinxupquote{constant}}: the value of the constant.
\end{quote}
\end{quote}

\subsubsection{Expr::AddExpr()}
\label{\detokenize{cppapi/expr:expr-addexpr}}\begin{quote}

\sphinxAtStartPar
Add an expression to self.

\sphinxAtStartPar
\sphinxstylestrong{Synopsis}
\begin{quote}

\sphinxAtStartPar
\sphinxcode{\sphinxupquote{void AddExpr(const Expr \&expr, double mult)}}
\end{quote}

\sphinxAtStartPar
\sphinxstylestrong{Arguments}
\begin{quote}

\sphinxAtStartPar
\sphinxcode{\sphinxupquote{expr}}: expression to be added.

\sphinxAtStartPar
\sphinxcode{\sphinxupquote{mult}}: optional, constant multiplier, default value is 1.0.
\end{quote}
\end{quote}

\subsubsection{Expr::AddTerm()}
\label{\detokenize{cppapi/expr:expr-addterm}}\begin{quote}

\sphinxAtStartPar
Add a term to expression object.

\sphinxAtStartPar
\sphinxstylestrong{Synopsis}
\begin{quote}

\sphinxAtStartPar
\sphinxcode{\sphinxupquote{void AddTerm(const Var \&var, double coeff)}}
\end{quote}

\sphinxAtStartPar
\sphinxstylestrong{Arguments}
\begin{quote}

\sphinxAtStartPar
\sphinxcode{\sphinxupquote{var}}: a variable for new term.

\sphinxAtStartPar
\sphinxcode{\sphinxupquote{coeff}}: coefficient for new term.
\end{quote}
\end{quote}

\subsubsection{Expr::AddTerms()}
\label{\detokenize{cppapi/expr:expr-addterms}}\begin{quote}

\sphinxAtStartPar
Add terms to expression object.

\sphinxAtStartPar
\sphinxstylestrong{Synopsis}
\begin{quote}

\sphinxAtStartPar
\sphinxcode{\sphinxupquote{int AddTerms(}}
\begin{quote}

\sphinxAtStartPar
\sphinxcode{\sphinxupquote{const VarArray \&vars,}}

\sphinxAtStartPar
\sphinxcode{\sphinxupquote{double *pCoeff,}}

\sphinxAtStartPar
\sphinxcode{\sphinxupquote{int len)}}
\end{quote}
\end{quote}

\sphinxAtStartPar
\sphinxstylestrong{Arguments}
\begin{quote}

\sphinxAtStartPar
\sphinxcode{\sphinxupquote{vars}}: variables for added terms.

\sphinxAtStartPar
\sphinxcode{\sphinxupquote{pCoeff}}: coefficient array for added terms.

\sphinxAtStartPar
\sphinxcode{\sphinxupquote{len}}: length of coefficient array.
\end{quote}

\sphinxAtStartPar
\sphinxstylestrong{Return}
\begin{quote}

\sphinxAtStartPar
number of added terms.
\end{quote}
\end{quote}

\subsubsection{Expr::Clear()}
\label{\detokenize{cppapi/expr:expr-clear}}\begin{quote}

\sphinxAtStartPar
Clear linear expression object.

\sphinxAtStartPar
\sphinxstylestrong{Synopsis}
\begin{quote}

\sphinxAtStartPar
\sphinxcode{\sphinxupquote{void Clear()}}
\end{quote}
\end{quote}

\subsubsection{Expr::Clone()}
\label{\detokenize{cppapi/expr:expr-clone}}\begin{quote}

\sphinxAtStartPar
Deep copy linear expression object.

\sphinxAtStartPar
\sphinxstylestrong{Synopsis}
\begin{quote}

\sphinxAtStartPar
\sphinxcode{\sphinxupquote{Expr Clone()}}
\end{quote}

\sphinxAtStartPar
\sphinxstylestrong{Return}
\begin{quote}

\sphinxAtStartPar
cloned expression object.
\end{quote}
\end{quote}

\subsubsection{Expr::Evaluate()}
\label{\detokenize{cppapi/expr:expr-evaluate}}\begin{quote}

\sphinxAtStartPar
Evaluate linear expression after solving.

\sphinxAtStartPar
\sphinxstylestrong{Synopsis}
\begin{quote}

\sphinxAtStartPar
\sphinxcode{\sphinxupquote{double Evaluate()}}
\end{quote}

\sphinxAtStartPar
\sphinxstylestrong{Return}
\begin{quote}

\sphinxAtStartPar
value of linear expression.
\end{quote}
\end{quote}

\subsubsection{Expr::GetCoeff()}
\label{\detokenize{cppapi/expr:expr-getcoeff}}\begin{quote}

\sphinxAtStartPar
Get coefficient from the i\sphinxhyphen{}th term in expression.

\sphinxAtStartPar
\sphinxstylestrong{Synopsis}
\begin{quote}

\sphinxAtStartPar
\sphinxcode{\sphinxupquote{double GetCoeff(int i)}}
\end{quote}

\sphinxAtStartPar
\sphinxstylestrong{Arguments}
\begin{quote}

\sphinxAtStartPar
\sphinxcode{\sphinxupquote{i}}: index of the term.
\end{quote}

\sphinxAtStartPar
\sphinxstylestrong{Return}
\begin{quote}

\sphinxAtStartPar
coefficient of the i\sphinxhyphen{}th term in expression object.
\end{quote}
\end{quote}

\subsubsection{Expr::GetConstant()}
\label{\detokenize{cppapi/expr:expr-getconstant}}\begin{quote}

\sphinxAtStartPar
Get constant in expression.

\sphinxAtStartPar
\sphinxstylestrong{Synopsis}
\begin{quote}

\sphinxAtStartPar
\sphinxcode{\sphinxupquote{double GetConstant()}}
\end{quote}

\sphinxAtStartPar
\sphinxstylestrong{Return}
\begin{quote}

\sphinxAtStartPar
constant in expression.
\end{quote}
\end{quote}

\subsubsection{Expr::GetVar()}
\label{\detokenize{cppapi/expr:expr-getvar}}\begin{quote}

\sphinxAtStartPar
Get variable from the i\sphinxhyphen{}th term in expression.

\sphinxAtStartPar
\sphinxstylestrong{Synopsis}
\begin{quote}

\sphinxAtStartPar
\sphinxcode{\sphinxupquote{Var \&GetVar(int i)}}
\end{quote}

\sphinxAtStartPar
\sphinxstylestrong{Arguments}
\begin{quote}

\sphinxAtStartPar
\sphinxcode{\sphinxupquote{i}}: index of the term.
\end{quote}

\sphinxAtStartPar
\sphinxstylestrong{Return}
\begin{quote}

\sphinxAtStartPar
variable of the i\sphinxhyphen{}th term in expression object.
\end{quote}
\end{quote}

\subsubsection{Expr::operator*=()}
\label{\detokenize{cppapi/expr:expr-operator}}\begin{quote}

\sphinxAtStartPar
Multiply a constant to self.

\sphinxAtStartPar
\sphinxstylestrong{Synopsis}
\begin{quote}

\sphinxAtStartPar
\sphinxcode{\sphinxupquote{void operator*=(double c)}}
\end{quote}

\sphinxAtStartPar
\sphinxstylestrong{Arguments}
\begin{quote}

\sphinxAtStartPar
\sphinxcode{\sphinxupquote{c}}: constant multiplier.
\end{quote}
\end{quote}

\subsubsection{Expr::operator*()}
\label{\detokenize{cppapi/expr:id2}}\begin{quote}

\sphinxAtStartPar
Multiply constant and return new expression.

\sphinxAtStartPar
\sphinxstylestrong{Synopsis}
\begin{quote}

\sphinxAtStartPar
\sphinxcode{\sphinxupquote{Expr operator*(double c)}}
\end{quote}

\sphinxAtStartPar
\sphinxstylestrong{Arguments}
\begin{quote}

\sphinxAtStartPar
\sphinxcode{\sphinxupquote{c}}: constant multiplier.
\end{quote}

\sphinxAtStartPar
\sphinxstylestrong{Return}
\begin{quote}

\sphinxAtStartPar
result expression.
\end{quote}
\end{quote}

\subsubsection{Expr::operator*()}
\label{\detokenize{cppapi/expr:id3}}\begin{quote}

\sphinxAtStartPar
Multiply a variable and return new quadratic expression object.

\sphinxAtStartPar
\sphinxstylestrong{Synopsis}
\begin{quote}

\sphinxAtStartPar
\sphinxcode{\sphinxupquote{QuadExpr operator*(const Var \&var)}}
\end{quote}

\sphinxAtStartPar
\sphinxstylestrong{Arguments}
\begin{quote}

\sphinxAtStartPar
\sphinxcode{\sphinxupquote{var}}: variable object.
\end{quote}

\sphinxAtStartPar
\sphinxstylestrong{Return}
\begin{quote}

\sphinxAtStartPar
result quadratic expression.
\end{quote}
\end{quote}

\subsubsection{Expr::operator*()}
\label{\detokenize{cppapi/expr:id4}}\begin{quote}

\sphinxAtStartPar
Multiply a linear expression and return new quadratic expression object.

\sphinxAtStartPar
\sphinxstylestrong{Synopsis}
\begin{quote}

\sphinxAtStartPar
\sphinxcode{\sphinxupquote{QuadExpr operator*(const Expr \&other)}}
\end{quote}

\sphinxAtStartPar
\sphinxstylestrong{Arguments}
\begin{quote}

\sphinxAtStartPar
\sphinxcode{\sphinxupquote{other}}: linear expression object.
\end{quote}

\sphinxAtStartPar
\sphinxstylestrong{Return}
\begin{quote}

\sphinxAtStartPar
result quadratic expression.
\end{quote}
\end{quote}

\subsubsection{Expr::operator/()}
\label{\detokenize{cppapi/expr:id5}}\begin{quote}

\sphinxAtStartPar
Devided by a constant and return new expression.

\sphinxAtStartPar
\sphinxstylestrong{Synopsis}
\begin{quote}

\sphinxAtStartPar
\sphinxcode{\sphinxupquote{Expr operator/(double c)}}
\end{quote}

\sphinxAtStartPar
\sphinxstylestrong{Arguments}
\begin{quote}

\sphinxAtStartPar
\sphinxcode{\sphinxupquote{c}}: constant divisor.
\end{quote}

\sphinxAtStartPar
\sphinxstylestrong{Return}
\begin{quote}

\sphinxAtStartPar
result expression.
\end{quote}
\end{quote}

\subsubsection{Expr::operator/()}
\label{\detokenize{cppapi/expr:id6}}\begin{quote}

\sphinxAtStartPar
Devided by a variable and return new nonlinear expression.

\sphinxAtStartPar
\sphinxstylestrong{Synopsis}
\begin{quote}

\sphinxAtStartPar
\sphinxcode{\sphinxupquote{NlExpr operator/(const Var \&var)}}
\end{quote}

\sphinxAtStartPar
\sphinxstylestrong{Arguments}
\begin{quote}

\sphinxAtStartPar
\sphinxcode{\sphinxupquote{var}}: a variable as divisor.
\end{quote}

\sphinxAtStartPar
\sphinxstylestrong{Return}
\begin{quote}

\sphinxAtStartPar
result nonlinear expression.
\end{quote}
\end{quote}

\subsubsection{Expr::operator/()}
\label{\detokenize{cppapi/expr:id7}}\begin{quote}

\sphinxAtStartPar
Devided by a linear expression and return new nonlinear expression.

\sphinxAtStartPar
\sphinxstylestrong{Synopsis}
\begin{quote}

\sphinxAtStartPar
\sphinxcode{\sphinxupquote{NlExpr operator/(const Expr \&other)}}
\end{quote}

\sphinxAtStartPar
\sphinxstylestrong{Arguments}
\begin{quote}

\sphinxAtStartPar
\sphinxcode{\sphinxupquote{other}}: a linear expression as divisor.
\end{quote}

\sphinxAtStartPar
\sphinxstylestrong{Return}
\begin{quote}

\sphinxAtStartPar
result nonlinear expression.
\end{quote}
\end{quote}

\subsubsection{Expr::operator+=()}
\label{\detokenize{cppapi/expr:id8}}\begin{quote}

\sphinxAtStartPar
Add an expression to self.

\sphinxAtStartPar
\sphinxstylestrong{Synopsis}
\begin{quote}

\sphinxAtStartPar
\sphinxcode{\sphinxupquote{void operator+=(const Expr \&expr)}}
\end{quote}

\sphinxAtStartPar
\sphinxstylestrong{Arguments}
\begin{quote}

\sphinxAtStartPar
\sphinxcode{\sphinxupquote{expr}}: expression to be added.
\end{quote}
\end{quote}

\subsubsection{Expr::operator+()}
\label{\detokenize{cppapi/expr:id9}}\begin{quote}

\sphinxAtStartPar
Add expression and return new expression.

\sphinxAtStartPar
\sphinxstylestrong{Synopsis}
\begin{quote}

\sphinxAtStartPar
\sphinxcode{\sphinxupquote{Expr operator+(const Expr \&other)}}
\end{quote}

\sphinxAtStartPar
\sphinxstylestrong{Arguments}
\begin{quote}

\sphinxAtStartPar
\sphinxcode{\sphinxupquote{other}}: other expression to add.
\end{quote}

\sphinxAtStartPar
\sphinxstylestrong{Return}
\begin{quote}

\sphinxAtStartPar
result expression.
\end{quote}
\end{quote}

\subsubsection{Expr::operator\sphinxhyphen{}=()}
\label{\detokenize{cppapi/expr:id10}}\begin{quote}

\sphinxAtStartPar
Substract an expression from self.

\sphinxAtStartPar
\sphinxstylestrong{Synopsis}
\begin{quote}

\sphinxAtStartPar
\sphinxcode{\sphinxupquote{void operator\sphinxhyphen{}=(const Expr \&expr)}}
\end{quote}

\sphinxAtStartPar
\sphinxstylestrong{Arguments}
\begin{quote}

\sphinxAtStartPar
\sphinxcode{\sphinxupquote{expr}}: expression to be substracted.
\end{quote}
\end{quote}

\subsubsection{Expr::operator\sphinxhyphen{}()}
\label{\detokenize{cppapi/expr:id11}}\begin{quote}

\sphinxAtStartPar
Substract expression and return new expression.

\sphinxAtStartPar
\sphinxstylestrong{Synopsis}
\begin{quote}

\sphinxAtStartPar
\sphinxcode{\sphinxupquote{Expr operator\sphinxhyphen{}(const Expr \&other)}}
\end{quote}

\sphinxAtStartPar
\sphinxstylestrong{Arguments}
\begin{quote}

\sphinxAtStartPar
\sphinxcode{\sphinxupquote{other}}: other expression to substract.
\end{quote}

\sphinxAtStartPar
\sphinxstylestrong{Return}
\begin{quote}

\sphinxAtStartPar
result expression.
\end{quote}
\end{quote}

\subsubsection{Expr::Remove()}
\label{\detokenize{cppapi/expr:expr-remove}}\begin{quote}

\sphinxAtStartPar
Remove i\sphinxhyphen{}th term from expression object.

\sphinxAtStartPar
\sphinxstylestrong{Synopsis}
\begin{quote}

\sphinxAtStartPar
\sphinxcode{\sphinxupquote{void Remove(int i)}}
\end{quote}

\sphinxAtStartPar
\sphinxstylestrong{Arguments}
\begin{quote}

\sphinxAtStartPar
\sphinxcode{\sphinxupquote{i}}: index of the term to be removed.
\end{quote}
\end{quote}

\subsubsection{Expr::Remove()}
\label{\detokenize{cppapi/expr:id12}}\begin{quote}

\sphinxAtStartPar
Remove the term associated with variable from expression.

\sphinxAtStartPar
\sphinxstylestrong{Synopsis}
\begin{quote}

\sphinxAtStartPar
\sphinxcode{\sphinxupquote{void Remove(const Var \&var)}}
\end{quote}

\sphinxAtStartPar
\sphinxstylestrong{Arguments}
\begin{quote}

\sphinxAtStartPar
\sphinxcode{\sphinxupquote{var}}: a variable whose term should be removed.
\end{quote}
\end{quote}

\subsubsection{Expr::Reserve()}
\label{\detokenize{cppapi/expr:expr-reserve}}\begin{quote}

\sphinxAtStartPar
Reserve capacity to contain at least n items.

\sphinxAtStartPar
\sphinxstylestrong{Synopsis}
\begin{quote}

\sphinxAtStartPar
\sphinxcode{\sphinxupquote{void Reserve(size\_t n)}}
\end{quote}

\sphinxAtStartPar
\sphinxstylestrong{Arguments}
\begin{quote}

\sphinxAtStartPar
\sphinxcode{\sphinxupquote{n}}: minimum capacity for linear expression object.
\end{quote}
\end{quote}

\subsubsection{Expr::SetCoeff()}
\label{\detokenize{cppapi/expr:expr-setcoeff}}\begin{quote}

\sphinxAtStartPar
Set coefficient for the i\sphinxhyphen{}th term in expression.

\sphinxAtStartPar
\sphinxstylestrong{Synopsis}
\begin{quote}

\sphinxAtStartPar
\sphinxcode{\sphinxupquote{void SetCoeff(int i, double val)}}
\end{quote}

\sphinxAtStartPar
\sphinxstylestrong{Arguments}
\begin{quote}

\sphinxAtStartPar
\sphinxcode{\sphinxupquote{i}}: index of the term.

\sphinxAtStartPar
\sphinxcode{\sphinxupquote{val}}: coefficient of the term.
\end{quote}
\end{quote}

\subsubsection{Expr::SetConstant()}
\label{\detokenize{cppapi/expr:expr-setconstant}}\begin{quote}

\sphinxAtStartPar
Set constant for the expression.

\sphinxAtStartPar
\sphinxstylestrong{Synopsis}
\begin{quote}

\sphinxAtStartPar
\sphinxcode{\sphinxupquote{void SetConstant(double constant)}}
\end{quote}

\sphinxAtStartPar
\sphinxstylestrong{Arguments}
\begin{quote}

\sphinxAtStartPar
\sphinxcode{\sphinxupquote{constant}}: the value of the constant.
\end{quote}
\end{quote}

\subsubsection{Expr::Size()}
\label{\detokenize{cppapi/expr:expr-size}}\begin{quote}

\sphinxAtStartPar
Get number of terms in expression.

\sphinxAtStartPar
\sphinxstylestrong{Synopsis}
\begin{quote}

\sphinxAtStartPar
\sphinxcode{\sphinxupquote{size\_t Size()}}
\end{quote}

\sphinxAtStartPar
\sphinxstylestrong{Return}
\begin{quote}

\sphinxAtStartPar
number of terms.
\end{quote}
\end{quote}

\subsection{Constraint}
\label{\detokenize{cppapiref:constraint}}\label{\detokenize{cppapiref:chapcppapiref-constraint}}
\sphinxAtStartPar
COPT constraint object. Constraints are always associated with a particular
model.  User creates a constraint object by adding a constraint to a model,
rather than by using constructor of Constraint class.

\sphinxstepscope

\subsubsection{Constraint::Get()}
\label{\detokenize{cppapi/constraint:constraint-get}}\label{\detokenize{cppapi/constraint::doc}}\begin{quote}

\sphinxAtStartPar
Get information value of the constraint. Support informations of “Dual”, “Slack”, “DualFarkas”, “LB”, “UB”.

\sphinxAtStartPar
\sphinxstylestrong{Synopsis}
\begin{quote}

\sphinxAtStartPar
\sphinxcode{\sphinxupquote{double Get(const char *szInfo)}}
\end{quote}

\sphinxAtStartPar
\sphinxstylestrong{Arguments}
\begin{quote}

\sphinxAtStartPar
\sphinxcode{\sphinxupquote{szInfo}}: name of the information being queried.
\end{quote}

\sphinxAtStartPar
\sphinxstylestrong{Return}
\begin{quote}

\sphinxAtStartPar
value of information.
\end{quote}
\end{quote}

\subsubsection{Constraint::GetBasis()}
\label{\detokenize{cppapi/constraint:constraint-getbasis}}\begin{quote}

\sphinxAtStartPar
Get basis status of this constraint.

\sphinxAtStartPar
\sphinxstylestrong{Synopsis}
\begin{quote}

\sphinxAtStartPar
\sphinxcode{\sphinxupquote{int GetBasis()}}
\end{quote}

\sphinxAtStartPar
\sphinxstylestrong{Return}
\begin{quote}

\sphinxAtStartPar
basis status.
\end{quote}
\end{quote}

\subsubsection{Constraint::GetIdx()}
\label{\detokenize{cppapi/constraint:constraint-getidx}}\begin{quote}

\sphinxAtStartPar
Get index of the constraint.

\sphinxAtStartPar
\sphinxstylestrong{Synopsis}
\begin{quote}

\sphinxAtStartPar
\sphinxcode{\sphinxupquote{int GetIdx()}}
\end{quote}

\sphinxAtStartPar
\sphinxstylestrong{Return}
\begin{quote}

\sphinxAtStartPar
the index of the constraint.
\end{quote}
\end{quote}

\subsubsection{Constraint::GetLowerIIS()}
\label{\detokenize{cppapi/constraint:constraint-getloweriis}}\begin{quote}

\sphinxAtStartPar
Get IIS status for lower bound of the constraint.

\sphinxAtStartPar
\sphinxstylestrong{Synopsis}
\begin{quote}

\sphinxAtStartPar
\sphinxcode{\sphinxupquote{int GetLowerIIS()}}
\end{quote}

\sphinxAtStartPar
\sphinxstylestrong{Return}
\begin{quote}

\sphinxAtStartPar
IIS status.
\end{quote}
\end{quote}

\subsubsection{Constraint::GetName()}
\label{\detokenize{cppapi/constraint:constraint-getname}}\begin{quote}

\sphinxAtStartPar
Get name of the constraint.

\sphinxAtStartPar
\sphinxstylestrong{Synopsis}
\begin{quote}

\sphinxAtStartPar
\sphinxcode{\sphinxupquote{const char *GetName()}}
\end{quote}

\sphinxAtStartPar
\sphinxstylestrong{Return}
\begin{quote}

\sphinxAtStartPar
the name of the constraint.
\end{quote}
\end{quote}

\subsubsection{Constraint::GetUpperIIS()}
\label{\detokenize{cppapi/constraint:constraint-getupperiis}}\begin{quote}

\sphinxAtStartPar
Get IIS status for upper bound of the constraint.

\sphinxAtStartPar
\sphinxstylestrong{Synopsis}
\begin{quote}

\sphinxAtStartPar
\sphinxcode{\sphinxupquote{int GetUpperIIS()}}
\end{quote}

\sphinxAtStartPar
\sphinxstylestrong{Return}
\begin{quote}

\sphinxAtStartPar
IIS status.
\end{quote}
\end{quote}

\subsubsection{Constraint::Remove()}
\label{\detokenize{cppapi/constraint:constraint-remove}}\begin{quote}

\sphinxAtStartPar
Remove this constraint from model.

\sphinxAtStartPar
\sphinxstylestrong{Synopsis}
\begin{quote}

\sphinxAtStartPar
\sphinxcode{\sphinxupquote{void Remove()}}
\end{quote}
\end{quote}

\subsubsection{Constraint::Set()}
\label{\detokenize{cppapi/constraint:constraint-set}}\begin{quote}

\sphinxAtStartPar
Set information value of the constraint. Support informations of “LB” and “UB”.

\sphinxAtStartPar
\sphinxstylestrong{Synopsis}
\begin{quote}

\sphinxAtStartPar
\sphinxcode{\sphinxupquote{void Set(const char *szInfo, double value)}}
\end{quote}

\sphinxAtStartPar
\sphinxstylestrong{Arguments}
\begin{quote}

\sphinxAtStartPar
\sphinxcode{\sphinxupquote{szInfo}}: name of the information.

\sphinxAtStartPar
\sphinxcode{\sphinxupquote{value}}: new information value.
\end{quote}
\end{quote}

\subsubsection{Constraint::SetName()}
\label{\detokenize{cppapi/constraint:constraint-setname}}\begin{quote}

\sphinxAtStartPar
Set name for the constraint.

\sphinxAtStartPar
\sphinxstylestrong{Synopsis}
\begin{quote}

\sphinxAtStartPar
\sphinxcode{\sphinxupquote{void SetName(const char *szName)}}
\end{quote}

\sphinxAtStartPar
\sphinxstylestrong{Arguments}
\begin{quote}

\sphinxAtStartPar
\sphinxcode{\sphinxupquote{szName}}: the name to set.
\end{quote}
\end{quote}

\subsection{ConstrArray}
\label{\detokenize{cppapiref:constrarray}}\label{\detokenize{cppapiref:chapcppapiref-constrarray}}
\sphinxAtStartPar
COPT constraint array object. To store and access a set of C++
{\hyperref[\detokenize{cppapiref:chapcppapiref-constraint}]{\sphinxcrossref{\DUrole{std,std-ref}{Constraint}}}} objects, Cardinal Optimizer provides C++
ConstrArray class, which defines the following methods.

\sphinxstepscope

\subsubsection{ConstrArray::GetConstr()}
\label{\detokenize{cppapi/constrarray:constrarray-getconstr}}\label{\detokenize{cppapi/constrarray::doc}}\begin{quote}

\sphinxAtStartPar
Get i\sphinxhyphen{}th constraint object.

\sphinxAtStartPar
\sphinxstylestrong{Synopsis}
\begin{quote}

\sphinxAtStartPar
\sphinxcode{\sphinxupquote{Constraint \&GetConstr(int i)}}
\end{quote}

\sphinxAtStartPar
\sphinxstylestrong{Arguments}
\begin{quote}

\sphinxAtStartPar
\sphinxcode{\sphinxupquote{i}}: index of the constraint.
\end{quote}

\sphinxAtStartPar
\sphinxstylestrong{Return}
\begin{quote}

\sphinxAtStartPar
constraint object with index i.
\end{quote}
\end{quote}

\subsubsection{ConstrArray::PushBack()}
\label{\detokenize{cppapi/constrarray:constrarray-pushback}}\begin{quote}

\sphinxAtStartPar
Add a constraint object to constraint array.

\sphinxAtStartPar
\sphinxstylestrong{Synopsis}
\begin{quote}

\sphinxAtStartPar
\sphinxcode{\sphinxupquote{void PushBack(const Constraint \&constr)}}
\end{quote}

\sphinxAtStartPar
\sphinxstylestrong{Arguments}
\begin{quote}

\sphinxAtStartPar
\sphinxcode{\sphinxupquote{constr}}: a constraint object.
\end{quote}
\end{quote}

\subsubsection{ConstrArray::Reserve()}
\label{\detokenize{cppapi/constrarray:constrarray-reserve}}\begin{quote}

\sphinxAtStartPar
Reserve capacity to contain at least n items.

\sphinxAtStartPar
\sphinxstylestrong{Synopsis}
\begin{quote}

\sphinxAtStartPar
\sphinxcode{\sphinxupquote{void Reserve(int n)}}
\end{quote}

\sphinxAtStartPar
\sphinxstylestrong{Arguments}
\begin{quote}

\sphinxAtStartPar
\sphinxcode{\sphinxupquote{n}}: minimum capacity for Constraint object.
\end{quote}
\end{quote}

\subsubsection{ConstrArray::Size()}
\label{\detokenize{cppapi/constrarray:constrarray-size}}\begin{quote}

\sphinxAtStartPar
Get the number of constraint objects.

\sphinxAtStartPar
\sphinxstylestrong{Synopsis}
\begin{quote}

\sphinxAtStartPar
\sphinxcode{\sphinxupquote{int Size()}}
\end{quote}

\sphinxAtStartPar
\sphinxstylestrong{Return}
\begin{quote}

\sphinxAtStartPar
number of constraint objects.
\end{quote}
\end{quote}

\subsection{ConstrBuilder}
\label{\detokenize{cppapiref:constrbuilder}}\label{\detokenize{cppapiref:chapcppapiref-constrbuilder}}
\sphinxAtStartPar
COPT constraint builder object. To help building a constraint, given a linear
expression, constraint sense and right\sphinxhyphen{}hand side value, Cardinal Optimizer
provides C++ ConstrBuilder class, which defines the following methods.

\sphinxstepscope

\subsubsection{ConstrBuilder::GetExpr()}
\label{\detokenize{cppapi/constrbuilder:constrbuilder-getexpr}}\label{\detokenize{cppapi/constrbuilder::doc}}\begin{quote}

\sphinxAtStartPar
Get expression associated with constraint.

\sphinxAtStartPar
\sphinxstylestrong{Synopsis}
\begin{quote}

\sphinxAtStartPar
\sphinxcode{\sphinxupquote{const Expr \&GetExpr()}}
\end{quote}

\sphinxAtStartPar
\sphinxstylestrong{Return}
\begin{quote}

\sphinxAtStartPar
expression object.
\end{quote}
\end{quote}

\subsubsection{ConstrBuilder::GetRange()}
\label{\detokenize{cppapi/constrbuilder:constrbuilder-getrange}}\begin{quote}

\sphinxAtStartPar
Get range from lower bound to upper bound of range constraint.

\sphinxAtStartPar
\sphinxstylestrong{Synopsis}
\begin{quote}

\sphinxAtStartPar
\sphinxcode{\sphinxupquote{double GetRange()}}
\end{quote}

\sphinxAtStartPar
\sphinxstylestrong{Return}
\begin{quote}

\sphinxAtStartPar
length from lower bound to upper bound of the constraint.
\end{quote}
\end{quote}

\subsubsection{ConstrBuilder::GetSense()}
\label{\detokenize{cppapi/constrbuilder:constrbuilder-getsense}}\begin{quote}

\sphinxAtStartPar
Get sense associated with constraint.

\sphinxAtStartPar
\sphinxstylestrong{Synopsis}
\begin{quote}

\sphinxAtStartPar
\sphinxcode{\sphinxupquote{char GetSense()}}
\end{quote}

\sphinxAtStartPar
\sphinxstylestrong{Return}
\begin{quote}

\sphinxAtStartPar
constraint sense.
\end{quote}
\end{quote}

\subsubsection{ConstrBuilder::Set()}
\label{\detokenize{cppapi/constrbuilder:constrbuilder-set}}\begin{quote}

\sphinxAtStartPar
Set detail of a constraint to its builder object.

\sphinxAtStartPar
\sphinxstylestrong{Synopsis}
\begin{quote}

\sphinxAtStartPar
\sphinxcode{\sphinxupquote{void Set(}}
\begin{quote}

\sphinxAtStartPar
\sphinxcode{\sphinxupquote{const Expr \&expr,}}

\sphinxAtStartPar
\sphinxcode{\sphinxupquote{char sense,}}

\sphinxAtStartPar
\sphinxcode{\sphinxupquote{double rhs)}}
\end{quote}
\end{quote}

\sphinxAtStartPar
\sphinxstylestrong{Arguments}
\begin{quote}

\sphinxAtStartPar
\sphinxcode{\sphinxupquote{expr}}: expression object at one side of the constraint

\sphinxAtStartPar
\sphinxcode{\sphinxupquote{sense}}: constraint sense other than COPT\_RANGE.

\sphinxAtStartPar
\sphinxcode{\sphinxupquote{rhs}}: constant of right side of the constraint.
\end{quote}
\end{quote}

\subsubsection{ConstrBuilder::SetRange()}
\label{\detokenize{cppapi/constrbuilder:constrbuilder-setrange}}\begin{quote}

\sphinxAtStartPar
Set a range constraint to its builder.

\sphinxAtStartPar
\sphinxstylestrong{Synopsis}
\begin{quote}

\sphinxAtStartPar
\sphinxcode{\sphinxupquote{void SetRange(const Expr \&expr, double range)}}
\end{quote}

\sphinxAtStartPar
\sphinxstylestrong{Arguments}
\begin{quote}

\sphinxAtStartPar
\sphinxcode{\sphinxupquote{expr}}: expression object, whose constant is negative upper bound.

\sphinxAtStartPar
\sphinxcode{\sphinxupquote{range}}: length from lower bound to upper bound of the constraint. Must greater than 0.
\end{quote}
\end{quote}

\subsection{ConstrBuilderArray}
\label{\detokenize{cppapiref:constrbuilderarray}}\label{\detokenize{cppapiref:chapcppapiref-constrbuilderarray}}
\sphinxAtStartPar
COPT constraint builder array object. To store and access a set of C++
{\hyperref[\detokenize{cppapiref:chapcppapiref-constrbuilder}]{\sphinxcrossref{\DUrole{std,std-ref}{ConstrBuilder}}}} objects, Cardinal Optimizer provides
C++ ConstrBuilderArray class, which defines the following methods.

\sphinxstepscope

\subsubsection{ConstrBuilderArray::GetBuilder()}
\label{\detokenize{cppapi/constrbuilderarray:constrbuilderarray-getbuilder}}\label{\detokenize{cppapi/constrbuilderarray::doc}}\begin{quote}

\sphinxAtStartPar
Get i\sphinxhyphen{}th constraint builder object.

\sphinxAtStartPar
\sphinxstylestrong{Synopsis}
\begin{quote}

\sphinxAtStartPar
\sphinxcode{\sphinxupquote{ConstrBuilder \&GetBuilder(int i)}}
\end{quote}

\sphinxAtStartPar
\sphinxstylestrong{Arguments}
\begin{quote}

\sphinxAtStartPar
\sphinxcode{\sphinxupquote{i}}: index of the constraint builder.
\end{quote}

\sphinxAtStartPar
\sphinxstylestrong{Return}
\begin{quote}

\sphinxAtStartPar
constraint builder object with index i.
\end{quote}
\end{quote}

\subsubsection{ConstrBuilderArray::PushBack()}
\label{\detokenize{cppapi/constrbuilderarray:constrbuilderarray-pushback}}\begin{quote}

\sphinxAtStartPar
Add a constraint builder object to constraint builder array.

\sphinxAtStartPar
\sphinxstylestrong{Synopsis}
\begin{quote}

\sphinxAtStartPar
\sphinxcode{\sphinxupquote{void PushBack(const ConstrBuilder \&builder)}}
\end{quote}

\sphinxAtStartPar
\sphinxstylestrong{Arguments}
\begin{quote}

\sphinxAtStartPar
\sphinxcode{\sphinxupquote{builder}}: a constraint builder object.
\end{quote}
\end{quote}

\subsubsection{ConstrBuilderArray::Reserve()}
\label{\detokenize{cppapi/constrbuilderarray:constrbuilderarray-reserve}}\begin{quote}

\sphinxAtStartPar
Reserve capacity to contain at least n items.

\sphinxAtStartPar
\sphinxstylestrong{Synopsis}
\begin{quote}

\sphinxAtStartPar
\sphinxcode{\sphinxupquote{void Reserve(int n)}}
\end{quote}

\sphinxAtStartPar
\sphinxstylestrong{Arguments}
\begin{quote}

\sphinxAtStartPar
\sphinxcode{\sphinxupquote{n}}: minimum capacity for constraint builder object.
\end{quote}
\end{quote}

\subsubsection{ConstrBuilderArray::Size()}
\label{\detokenize{cppapi/constrbuilderarray:constrbuilderarray-size}}\begin{quote}

\sphinxAtStartPar
Get the number of constraint builder objects.

\sphinxAtStartPar
\sphinxstylestrong{Synopsis}
\begin{quote}

\sphinxAtStartPar
\sphinxcode{\sphinxupquote{int Size()}}
\end{quote}

\sphinxAtStartPar
\sphinxstylestrong{Return}
\begin{quote}

\sphinxAtStartPar
number of constraint builder objects.
\end{quote}
\end{quote}

\subsection{Column}
\label{\detokenize{cppapiref:column}}\label{\detokenize{cppapiref:chapcppapiref-column}}
\sphinxAtStartPar
COPT column object. A column consists of a list of constraints and associated
coefficients.  Columns are used to represent the set of constraints in which a
variable participates, and the asssociated coefficents.

\sphinxstepscope

\subsubsection{Column::Column()}
\label{\detokenize{cppapi/column:column-column}}\label{\detokenize{cppapi/column::doc}}\begin{quote}

\sphinxAtStartPar
Constructor of column.

\sphinxAtStartPar
\sphinxstylestrong{Synopsis}
\begin{quote}

\sphinxAtStartPar
\sphinxcode{\sphinxupquote{Column()}}
\end{quote}
\end{quote}

\subsubsection{Column::AddColumn()}
\label{\detokenize{cppapi/column:column-addcolumn}}\begin{quote}

\sphinxAtStartPar
Add a column to self.

\sphinxAtStartPar
\sphinxstylestrong{Synopsis}
\begin{quote}

\sphinxAtStartPar
\sphinxcode{\sphinxupquote{void AddColumn(const Column \&col, double mult)}}
\end{quote}

\sphinxAtStartPar
\sphinxstylestrong{Arguments}
\begin{quote}

\sphinxAtStartPar
\sphinxcode{\sphinxupquote{col}}: column object to be added.

\sphinxAtStartPar
\sphinxcode{\sphinxupquote{mult}}: multiply constant.
\end{quote}
\end{quote}

\subsubsection{Column::AddTerm()}
\label{\detokenize{cppapi/column:column-addterm}}\begin{quote}

\sphinxAtStartPar
Add a term to column object.

\sphinxAtStartPar
\sphinxstylestrong{Synopsis}
\begin{quote}

\sphinxAtStartPar
\sphinxcode{\sphinxupquote{void AddTerm(const Constraint \&constr, double coeff)}}
\end{quote}

\sphinxAtStartPar
\sphinxstylestrong{Arguments}
\begin{quote}

\sphinxAtStartPar
\sphinxcode{\sphinxupquote{constr}}: a constraint for new term.

\sphinxAtStartPar
\sphinxcode{\sphinxupquote{coeff}}: coefficient for new term.
\end{quote}
\end{quote}

\subsubsection{Column::AddTerms()}
\label{\detokenize{cppapi/column:column-addterms}}\begin{quote}

\sphinxAtStartPar
Add terms to column object.

\sphinxAtStartPar
\sphinxstylestrong{Synopsis}
\begin{quote}

\sphinxAtStartPar
\sphinxcode{\sphinxupquote{int AddTerms(}}
\begin{quote}

\sphinxAtStartPar
\sphinxcode{\sphinxupquote{const ConstrArray \&constrs,}}

\sphinxAtStartPar
\sphinxcode{\sphinxupquote{double *pCoeff,}}

\sphinxAtStartPar
\sphinxcode{\sphinxupquote{int len)}}
\end{quote}
\end{quote}

\sphinxAtStartPar
\sphinxstylestrong{Arguments}
\begin{quote}

\sphinxAtStartPar
\sphinxcode{\sphinxupquote{constrs}}: constraints for added terms.

\sphinxAtStartPar
\sphinxcode{\sphinxupquote{pCoeff}}: coefficients for added terms.

\sphinxAtStartPar
\sphinxcode{\sphinxupquote{len}}: number of terms to be added.
\end{quote}

\sphinxAtStartPar
\sphinxstylestrong{Return}
\begin{quote}

\sphinxAtStartPar
number of added terms.
\end{quote}
\end{quote}

\subsubsection{Column::Clear()}
\label{\detokenize{cppapi/column:column-clear}}\begin{quote}

\sphinxAtStartPar
Clear all terms.

\sphinxAtStartPar
\sphinxstylestrong{Synopsis}
\begin{quote}

\sphinxAtStartPar
\sphinxcode{\sphinxupquote{void Clear()}}
\end{quote}
\end{quote}

\subsubsection{Column::Clone()}
\label{\detokenize{cppapi/column:column-clone}}\begin{quote}

\sphinxAtStartPar
Deep copy column object.

\sphinxAtStartPar
\sphinxstylestrong{Synopsis}
\begin{quote}

\sphinxAtStartPar
\sphinxcode{\sphinxupquote{Column Clone()}}
\end{quote}

\sphinxAtStartPar
\sphinxstylestrong{Return}
\begin{quote}

\sphinxAtStartPar
cloned column object.
\end{quote}
\end{quote}

\subsubsection{Column::GetCoeff()}
\label{\detokenize{cppapi/column:column-getcoeff}}\begin{quote}

\sphinxAtStartPar
Get coefficient from the i\sphinxhyphen{}th term in column object.

\sphinxAtStartPar
\sphinxstylestrong{Synopsis}
\begin{quote}

\sphinxAtStartPar
\sphinxcode{\sphinxupquote{double GetCoeff(int i)}}
\end{quote}

\sphinxAtStartPar
\sphinxstylestrong{Arguments}
\begin{quote}

\sphinxAtStartPar
\sphinxcode{\sphinxupquote{i}}: index of the term.
\end{quote}

\sphinxAtStartPar
\sphinxstylestrong{Return}
\begin{quote}

\sphinxAtStartPar
coefficient of the i\sphinxhyphen{}th term in column object.
\end{quote}
\end{quote}

\subsubsection{Column::GetConstr()}
\label{\detokenize{cppapi/column:column-getconstr}}\begin{quote}

\sphinxAtStartPar
Get constraint from the i\sphinxhyphen{}th term in column object.

\sphinxAtStartPar
\sphinxstylestrong{Synopsis}
\begin{quote}

\sphinxAtStartPar
\sphinxcode{\sphinxupquote{Constraint GetConstr(int i)}}
\end{quote}

\sphinxAtStartPar
\sphinxstylestrong{Arguments}
\begin{quote}

\sphinxAtStartPar
\sphinxcode{\sphinxupquote{i}}: index of the term.
\end{quote}

\sphinxAtStartPar
\sphinxstylestrong{Return}
\begin{quote}

\sphinxAtStartPar
constraint of the i\sphinxhyphen{}th term in column object.
\end{quote}
\end{quote}

\subsubsection{Column::Remove()}
\label{\detokenize{cppapi/column:column-remove}}\begin{quote}

\sphinxAtStartPar
Remove i\sphinxhyphen{}th term from column object.

\sphinxAtStartPar
\sphinxstylestrong{Synopsis}
\begin{quote}

\sphinxAtStartPar
\sphinxcode{\sphinxupquote{void Remove(int i)}}
\end{quote}

\sphinxAtStartPar
\sphinxstylestrong{Arguments}
\begin{quote}

\sphinxAtStartPar
\sphinxcode{\sphinxupquote{i}}: index of the term to be removed.
\end{quote}
\end{quote}

\subsubsection{Column::Remove()}
\label{\detokenize{cppapi/column:id1}}\begin{quote}

\sphinxAtStartPar
Remove the term associated with constraint from column object.

\sphinxAtStartPar
\sphinxstylestrong{Synopsis}
\begin{quote}

\sphinxAtStartPar
\sphinxcode{\sphinxupquote{bool Remove(const Constraint \&constr)}}
\end{quote}

\sphinxAtStartPar
\sphinxstylestrong{Arguments}
\begin{quote}

\sphinxAtStartPar
\sphinxcode{\sphinxupquote{constr}}: a constraint whose term should be removed.
\end{quote}

\sphinxAtStartPar
\sphinxstylestrong{Return}
\begin{quote}

\sphinxAtStartPar
true if constraint exits in column object.
\end{quote}
\end{quote}

\subsubsection{Column::Reserve()}
\label{\detokenize{cppapi/column:column-reserve}}\begin{quote}

\sphinxAtStartPar
Reserve capacity to contain at least n items.

\sphinxAtStartPar
\sphinxstylestrong{Synopsis}
\begin{quote}

\sphinxAtStartPar
\sphinxcode{\sphinxupquote{void Reserve(int n)}}
\end{quote}

\sphinxAtStartPar
\sphinxstylestrong{Arguments}
\begin{quote}

\sphinxAtStartPar
\sphinxcode{\sphinxupquote{n}}: minimum capacity for Column object.
\end{quote}
\end{quote}

\subsubsection{Column::Size()}
\label{\detokenize{cppapi/column:column-size}}\begin{quote}

\sphinxAtStartPar
Get number of terms in column object.

\sphinxAtStartPar
\sphinxstylestrong{Synopsis}
\begin{quote}

\sphinxAtStartPar
\sphinxcode{\sphinxupquote{int Size()}}
\end{quote}

\sphinxAtStartPar
\sphinxstylestrong{Return}
\begin{quote}

\sphinxAtStartPar
number of terms.
\end{quote}
\end{quote}

\subsection{ColumnArray}
\label{\detokenize{cppapiref:columnarray}}\label{\detokenize{cppapiref:chapcppapiref-columnarray}}
\sphinxAtStartPar
COPT column array object. To store and access a set of C++
{\hyperref[\detokenize{cppapiref:chapcppapiref-column}]{\sphinxcrossref{\DUrole{std,std-ref}{Column}}}} objects, Cardinal Optimizer provides C++
ColumnArray class, which defines the following methods.

\sphinxstepscope

\subsubsection{ColumnArray::Clear()}
\label{\detokenize{cppapi/columnarray:columnarray-clear}}\label{\detokenize{cppapi/columnarray::doc}}\begin{quote}

\sphinxAtStartPar
Clear all column objects.

\sphinxAtStartPar
\sphinxstylestrong{Synopsis}
\begin{quote}

\sphinxAtStartPar
\sphinxcode{\sphinxupquote{void Clear()}}
\end{quote}
\end{quote}

\subsubsection{ColumnArray::GetColumn()}
\label{\detokenize{cppapi/columnarray:columnarray-getcolumn}}\begin{quote}

\sphinxAtStartPar
Get i\sphinxhyphen{}th column object.

\sphinxAtStartPar
\sphinxstylestrong{Synopsis}
\begin{quote}

\sphinxAtStartPar
\sphinxcode{\sphinxupquote{Column \&GetColumn(int i)}}
\end{quote}

\sphinxAtStartPar
\sphinxstylestrong{Arguments}
\begin{quote}

\sphinxAtStartPar
\sphinxcode{\sphinxupquote{i}}: index of the column.
\end{quote}

\sphinxAtStartPar
\sphinxstylestrong{Return}
\begin{quote}

\sphinxAtStartPar
column object with index i.
\end{quote}
\end{quote}

\subsubsection{ColumnArray::PushBack()}
\label{\detokenize{cppapi/columnarray:columnarray-pushback}}\begin{quote}

\sphinxAtStartPar
Add a column object to column array.

\sphinxAtStartPar
\sphinxstylestrong{Synopsis}
\begin{quote}

\sphinxAtStartPar
\sphinxcode{\sphinxupquote{void PushBack(const Column \&col)}}
\end{quote}

\sphinxAtStartPar
\sphinxstylestrong{Arguments}
\begin{quote}

\sphinxAtStartPar
\sphinxcode{\sphinxupquote{col}}: a column object.
\end{quote}
\end{quote}

\subsubsection{ColumnArray::Reserve()}
\label{\detokenize{cppapi/columnarray:columnarray-reserve}}\begin{quote}

\sphinxAtStartPar
Reserve capacity to contain at least n items.

\sphinxAtStartPar
\sphinxstylestrong{Synopsis}
\begin{quote}

\sphinxAtStartPar
\sphinxcode{\sphinxupquote{void Reserve(int n)}}
\end{quote}

\sphinxAtStartPar
\sphinxstylestrong{Arguments}
\begin{quote}

\sphinxAtStartPar
\sphinxcode{\sphinxupquote{n}}: minimum capacity for linear expression object.
\end{quote}
\end{quote}

\subsubsection{ColumnArray::Size()}
\label{\detokenize{cppapi/columnarray:columnarray-size}}\begin{quote}

\sphinxAtStartPar
Get the number of column objects.

\sphinxAtStartPar
\sphinxstylestrong{Synopsis}
\begin{quote}

\sphinxAtStartPar
\sphinxcode{\sphinxupquote{int Size()}}
\end{quote}

\sphinxAtStartPar
\sphinxstylestrong{Return}
\begin{quote}

\sphinxAtStartPar
number of column objects.
\end{quote}
\end{quote}

\subsection{Sos}
\label{\detokenize{cppapiref:sos}}\label{\detokenize{cppapiref:chapcppapiref-sos}}
\sphinxAtStartPar
COPT SOS constraint object. SOS constraints are always associated with a
particular model.  User creates an SOS constraint object by adding an SOS
constraint to a model, rather than by using constructor of Sos class.

\sphinxAtStartPar
An SOS constraint can be type 1 or 2 (\sphinxcode{\sphinxupquote{COPT\_SOS\_TYPE1}} or
\sphinxcode{\sphinxupquote{COPT\_SOS\_TYPE2}}).

\sphinxstepscope

\subsubsection{Sos::GetIdx()}
\label{\detokenize{cppapi/sos:sos-getidx}}\label{\detokenize{cppapi/sos::doc}}\begin{quote}

\sphinxAtStartPar
Get the index of SOS constraint.

\sphinxAtStartPar
\sphinxstylestrong{Synopsis}
\begin{quote}

\sphinxAtStartPar
\sphinxcode{\sphinxupquote{int GetIdx()}}
\end{quote}

\sphinxAtStartPar
\sphinxstylestrong{Return}
\begin{quote}

\sphinxAtStartPar
index of SOS constraint.
\end{quote}
\end{quote}

\subsubsection{Sos::GetIIS()}
\label{\detokenize{cppapi/sos:sos-getiis}}\begin{quote}

\sphinxAtStartPar
Get IIS status of the SOS constraint.

\sphinxAtStartPar
\sphinxstylestrong{Synopsis}
\begin{quote}

\sphinxAtStartPar
\sphinxcode{\sphinxupquote{int GetIIS()}}
\end{quote}

\sphinxAtStartPar
\sphinxstylestrong{Return}
\begin{quote}

\sphinxAtStartPar
IIS status.
\end{quote}
\end{quote}

\subsubsection{Sos::Remove()}
\label{\detokenize{cppapi/sos:sos-remove}}\begin{quote}

\sphinxAtStartPar
Remove the SOS constraint from model.

\sphinxAtStartPar
\sphinxstylestrong{Synopsis}
\begin{quote}

\sphinxAtStartPar
\sphinxcode{\sphinxupquote{void Remove()}}
\end{quote}
\end{quote}

\subsection{SosArray}
\label{\detokenize{cppapiref:sosarray}}\label{\detokenize{cppapiref:chapcppapiref-sosarray}}
\sphinxAtStartPar
COPT SOS constraint array object. To store and access a set of C++
{\hyperref[\detokenize{cppapiref:chapcppapiref-sos}]{\sphinxcrossref{\DUrole{std,std-ref}{Sos}}}} objects, Cardinal Optimizer provides C++ SosArray
class, which defines the following methods.

\sphinxstepscope

\subsubsection{SosArray::GetSos()}
\label{\detokenize{cppapi/sosarray:sosarray-getsos}}\label{\detokenize{cppapi/sosarray::doc}}\begin{quote}

\sphinxAtStartPar
Get i\sphinxhyphen{}th SOS constraint object.

\sphinxAtStartPar
\sphinxstylestrong{Synopsis}
\begin{quote}

\sphinxAtStartPar
\sphinxcode{\sphinxupquote{Sos \&GetSos(int i)}}
\end{quote}

\sphinxAtStartPar
\sphinxstylestrong{Arguments}
\begin{quote}

\sphinxAtStartPar
\sphinxcode{\sphinxupquote{i}}: index of the SOS constraint.
\end{quote}

\sphinxAtStartPar
\sphinxstylestrong{Return}
\begin{quote}

\sphinxAtStartPar
SOS constraint object with index i.
\end{quote}
\end{quote}

\subsubsection{SosArray::PushBack()}
\label{\detokenize{cppapi/sosarray:sosarray-pushback}}\begin{quote}

\sphinxAtStartPar
Add a SOS constraint object to SOS constraint array.

\sphinxAtStartPar
\sphinxstylestrong{Synopsis}
\begin{quote}

\sphinxAtStartPar
\sphinxcode{\sphinxupquote{void PushBack(const Sos \&sos)}}
\end{quote}

\sphinxAtStartPar
\sphinxstylestrong{Arguments}
\begin{quote}

\sphinxAtStartPar
\sphinxcode{\sphinxupquote{sos}}: a SOS constraint object.
\end{quote}
\end{quote}

\subsubsection{SosArray::Size()}
\label{\detokenize{cppapi/sosarray:sosarray-size}}\begin{quote}

\sphinxAtStartPar
Get the number of SOS constraint objects.

\sphinxAtStartPar
\sphinxstylestrong{Synopsis}
\begin{quote}

\sphinxAtStartPar
\sphinxcode{\sphinxupquote{int Size()}}
\end{quote}

\sphinxAtStartPar
\sphinxstylestrong{Return}
\begin{quote}

\sphinxAtStartPar
number of SOS constraint objects.
\end{quote}
\end{quote}

\subsection{SosBuilder}
\label{\detokenize{cppapiref:sosbuilder}}\label{\detokenize{cppapiref:chapcppapiref-sosbuilder}}
\sphinxAtStartPar
COPT SOS constraint builder object. To help building an SOS constraint,
given the SOS type, a set of variables and associated weights, Cardinal
Optimizer provides C++ SosBuilder class, which defines the following methods.

\sphinxstepscope

\subsubsection{SosBuilder::GetSize()}
\label{\detokenize{cppapi/sosbuilder:sosbuilder-getsize}}\label{\detokenize{cppapi/sosbuilder::doc}}\begin{quote}

\sphinxAtStartPar
Get number of terms in SOS constraint.

\sphinxAtStartPar
\sphinxstylestrong{Synopsis}
\begin{quote}

\sphinxAtStartPar
\sphinxcode{\sphinxupquote{int GetSize()}}
\end{quote}

\sphinxAtStartPar
\sphinxstylestrong{Return}
\begin{quote}

\sphinxAtStartPar
number of terms.
\end{quote}
\end{quote}

\subsubsection{SosBuilder::GetType()}
\label{\detokenize{cppapi/sosbuilder:sosbuilder-gettype}}\begin{quote}

\sphinxAtStartPar
Get type of SOS constraint.

\sphinxAtStartPar
\sphinxstylestrong{Synopsis}
\begin{quote}

\sphinxAtStartPar
\sphinxcode{\sphinxupquote{int GetType()}}
\end{quote}

\sphinxAtStartPar
\sphinxstylestrong{Return}
\begin{quote}

\sphinxAtStartPar
type of SOS constraint.
\end{quote}
\end{quote}

\subsubsection{SosBuilder::GetVar()}
\label{\detokenize{cppapi/sosbuilder:sosbuilder-getvar}}\begin{quote}

\sphinxAtStartPar
Get variable from the i\sphinxhyphen{}th term in SOS constraint.

\sphinxAtStartPar
\sphinxstylestrong{Synopsis}
\begin{quote}

\sphinxAtStartPar
\sphinxcode{\sphinxupquote{Var GetVar(int i)}}
\end{quote}

\sphinxAtStartPar
\sphinxstylestrong{Arguments}
\begin{quote}

\sphinxAtStartPar
\sphinxcode{\sphinxupquote{i}}: index of the term.
\end{quote}

\sphinxAtStartPar
\sphinxstylestrong{Return}
\begin{quote}

\sphinxAtStartPar
variable of the i\sphinxhyphen{}th term in SOS constraint.
\end{quote}
\end{quote}

\subsubsection{SosBuilder::GetVars()}
\label{\detokenize{cppapi/sosbuilder:sosbuilder-getvars}}\begin{quote}

\sphinxAtStartPar
Get variables of all terms in SOS constraint.

\sphinxAtStartPar
\sphinxstylestrong{Synopsis}
\begin{quote}

\sphinxAtStartPar
\sphinxcode{\sphinxupquote{VarArray GetVars()}}
\end{quote}

\sphinxAtStartPar
\sphinxstylestrong{Return}
\begin{quote}

\sphinxAtStartPar
variable array object.
\end{quote}
\end{quote}

\subsubsection{SosBuilder::GetWeight()}
\label{\detokenize{cppapi/sosbuilder:sosbuilder-getweight}}\begin{quote}

\sphinxAtStartPar
Get weight from the i\sphinxhyphen{}th term in SOS constraint.

\sphinxAtStartPar
\sphinxstylestrong{Synopsis}
\begin{quote}

\sphinxAtStartPar
\sphinxcode{\sphinxupquote{double GetWeight(int i)}}
\end{quote}

\sphinxAtStartPar
\sphinxstylestrong{Arguments}
\begin{quote}

\sphinxAtStartPar
\sphinxcode{\sphinxupquote{i}}: index of the term.
\end{quote}

\sphinxAtStartPar
\sphinxstylestrong{Return}
\begin{quote}

\sphinxAtStartPar
weight of the i\sphinxhyphen{}th term in SOS constraint.
\end{quote}
\end{quote}

\subsubsection{SosBuilder::GetWeights()}
\label{\detokenize{cppapi/sosbuilder:sosbuilder-getweights}}\begin{quote}

\sphinxAtStartPar
Get weights of all terms in SOS constraint.

\sphinxAtStartPar
\sphinxstylestrong{Synopsis}
\begin{quote}

\sphinxAtStartPar
\sphinxcode{\sphinxupquote{double GetWeights()}}
\end{quote}

\sphinxAtStartPar
\sphinxstylestrong{Return}
\begin{quote}

\sphinxAtStartPar
pointer to array of weights.
\end{quote}
\end{quote}

\subsubsection{SosBuilder::Set()}
\label{\detokenize{cppapi/sosbuilder:sosbuilder-set}}\begin{quote}

\sphinxAtStartPar
Set variables and weights of SOS constraint.

\sphinxAtStartPar
\sphinxstylestrong{Synopsis}
\begin{quote}

\sphinxAtStartPar
\sphinxcode{\sphinxupquote{void Set(}}
\begin{quote}

\sphinxAtStartPar
\sphinxcode{\sphinxupquote{const VarArray \&vars,}}

\sphinxAtStartPar
\sphinxcode{\sphinxupquote{const double *pWeights,}}

\sphinxAtStartPar
\sphinxcode{\sphinxupquote{int type)}}
\end{quote}
\end{quote}

\sphinxAtStartPar
\sphinxstylestrong{Arguments}
\begin{quote}

\sphinxAtStartPar
\sphinxcode{\sphinxupquote{vars}}: variable array object.

\sphinxAtStartPar
\sphinxcode{\sphinxupquote{pWeights}}: pointer to array of weights.

\sphinxAtStartPar
\sphinxcode{\sphinxupquote{type}}: type of SOS constraint.
\end{quote}
\end{quote}

\subsection{SosBuilderArray}
\label{\detokenize{cppapiref:sosbuilderarray}}\label{\detokenize{cppapiref:chapcppapiref-sosbuilderarray}}
\sphinxAtStartPar
COPT SOS constraint builder array object. To store and access a set of C++
{\hyperref[\detokenize{cppapiref:chapcppapiref-sosbuilder}]{\sphinxcrossref{\DUrole{std,std-ref}{SosBuilder}}}} objects, Cardinal Optimizer provides C++
SosBuilderArray class, which defines the following methods.

\sphinxstepscope

\subsubsection{SosBuilderArray::GetBuilder()}
\label{\detokenize{cppapi/sosbuilderarray:sosbuilderarray-getbuilder}}\label{\detokenize{cppapi/sosbuilderarray::doc}}\begin{quote}

\sphinxAtStartPar
Get i\sphinxhyphen{}th SOS constraint builder object.

\sphinxAtStartPar
\sphinxstylestrong{Synopsis}
\begin{quote}

\sphinxAtStartPar
\sphinxcode{\sphinxupquote{SosBuilder \&GetBuilder(int i)}}
\end{quote}

\sphinxAtStartPar
\sphinxstylestrong{Arguments}
\begin{quote}

\sphinxAtStartPar
\sphinxcode{\sphinxupquote{i}}: index of the SOS constraint builder.
\end{quote}

\sphinxAtStartPar
\sphinxstylestrong{Return}
\begin{quote}

\sphinxAtStartPar
SOS constraint builder object with index i.
\end{quote}
\end{quote}

\subsubsection{SosBuilderArray::PushBack()}
\label{\detokenize{cppapi/sosbuilderarray:sosbuilderarray-pushback}}\begin{quote}

\sphinxAtStartPar
Add a SOS constraint builder object to SOS constraint builder array.

\sphinxAtStartPar
\sphinxstylestrong{Synopsis}
\begin{quote}

\sphinxAtStartPar
\sphinxcode{\sphinxupquote{void PushBack(const SosBuilder \&builder)}}
\end{quote}

\sphinxAtStartPar
\sphinxstylestrong{Arguments}
\begin{quote}

\sphinxAtStartPar
\sphinxcode{\sphinxupquote{builder}}: a SOS constraint builder object.
\end{quote}
\end{quote}

\subsubsection{SosBuilderArray::Size()}
\label{\detokenize{cppapi/sosbuilderarray:sosbuilderarray-size}}\begin{quote}

\sphinxAtStartPar
Get the number of SOS constraint builder objects.

\sphinxAtStartPar
\sphinxstylestrong{Synopsis}
\begin{quote}

\sphinxAtStartPar
\sphinxcode{\sphinxupquote{int Size()}}
\end{quote}

\sphinxAtStartPar
\sphinxstylestrong{Return}
\begin{quote}

\sphinxAtStartPar
number of SOS constraint builder objects.
\end{quote}
\end{quote}

\subsection{GenConstr}
\label{\detokenize{cppapiref:genconstr}}\label{\detokenize{cppapiref:chapcppapiref-genconstr}}
\sphinxAtStartPar
COPT general constraint object. General constraints are always associated
with a particular model.  User creates a general constraint object by adding
a general constraint to a model, rather than by using constructor of GenConstr
class.

\sphinxstepscope

\subsubsection{GenConstr::GetIdx()}
\label{\detokenize{cppapi/genconstr:genconstr-getidx}}\label{\detokenize{cppapi/genconstr::doc}}\begin{quote}

\sphinxAtStartPar
Get the index of the general constraint.

\sphinxAtStartPar
\sphinxstylestrong{Synopsis}
\begin{quote}

\sphinxAtStartPar
\sphinxcode{\sphinxupquote{int GetIdx()}}
\end{quote}

\sphinxAtStartPar
\sphinxstylestrong{Return}
\begin{quote}

\sphinxAtStartPar
index of the general constraint.
\end{quote}
\end{quote}

\subsubsection{GenConstr::GetIIS()}
\label{\detokenize{cppapi/genconstr:genconstr-getiis}}\begin{quote}

\sphinxAtStartPar
Get IIS status of the general constraint.

\sphinxAtStartPar
\sphinxstylestrong{Synopsis}
\begin{quote}

\sphinxAtStartPar
\sphinxcode{\sphinxupquote{int GetIIS()}}
\end{quote}

\sphinxAtStartPar
\sphinxstylestrong{Return}
\begin{quote}

\sphinxAtStartPar
IIS status.
\end{quote}
\end{quote}

\subsubsection{GenConstr::GetName()}
\label{\detokenize{cppapi/genconstr:genconstr-getname}}\begin{quote}

\sphinxAtStartPar
Get name of gerneral constraint.

\sphinxAtStartPar
\sphinxstylestrong{Synopsis}
\begin{quote}

\sphinxAtStartPar
\sphinxcode{\sphinxupquote{const char *GetName()}}
\end{quote}

\sphinxAtStartPar
\sphinxstylestrong{Return}
\begin{quote}

\sphinxAtStartPar
the name of general constraint.
\end{quote}
\end{quote}

\subsubsection{GenConstr::Remove()}
\label{\detokenize{cppapi/genconstr:genconstr-remove}}\begin{quote}

\sphinxAtStartPar
Remove the general constraint from model.

\sphinxAtStartPar
\sphinxstylestrong{Synopsis}
\begin{quote}

\sphinxAtStartPar
\sphinxcode{\sphinxupquote{void Remove()}}
\end{quote}
\end{quote}

\subsubsection{GenConstr::SetName()}
\label{\detokenize{cppapi/genconstr:genconstr-setname}}\begin{quote}

\sphinxAtStartPar
Set name of general constraint.

\sphinxAtStartPar
\sphinxstylestrong{Synopsis}
\begin{quote}

\sphinxAtStartPar
\sphinxcode{\sphinxupquote{void SetName(const char *szName)}}
\end{quote}

\sphinxAtStartPar
\sphinxstylestrong{Arguments}
\begin{quote}

\sphinxAtStartPar
\sphinxcode{\sphinxupquote{szName}}: new name to set.
\end{quote}
\end{quote}

\subsection{GenConstrArray}
\label{\detokenize{cppapiref:genconstrarray}}\label{\detokenize{cppapiref:chapcppapiref-genconstrarray}}
\sphinxAtStartPar
COPT general constraint array object. To store and access a set of C++
{\hyperref[\detokenize{cppapiref:chapcppapiref-genconstr}]{\sphinxcrossref{\DUrole{std,std-ref}{GenConstr}}}} objects, Cardinal Optimizer provides C++
GenConstrArray class, which defines the following methods.

\sphinxstepscope

\subsubsection{GenConstrArray::GetGenConstr()}
\label{\detokenize{cppapi/genconstrarray:genconstrarray-getgenconstr}}\label{\detokenize{cppapi/genconstrarray::doc}}\begin{quote}

\sphinxAtStartPar
Get i\sphinxhyphen{}th general constraint object.

\sphinxAtStartPar
\sphinxstylestrong{Synopsis}
\begin{quote}

\sphinxAtStartPar
\sphinxcode{\sphinxupquote{GenConstr \&GetGenConstr(int i)}}
\end{quote}

\sphinxAtStartPar
\sphinxstylestrong{Arguments}
\begin{quote}

\sphinxAtStartPar
\sphinxcode{\sphinxupquote{i}}: index of the general constraint.
\end{quote}

\sphinxAtStartPar
\sphinxstylestrong{Return}
\begin{quote}

\sphinxAtStartPar
general constraint object with index i.
\end{quote}
\end{quote}

\subsubsection{GenConstrArray::PushBack()}
\label{\detokenize{cppapi/genconstrarray:genconstrarray-pushback}}\begin{quote}

\sphinxAtStartPar
Add a general constraint object to general constraint array.

\sphinxAtStartPar
\sphinxstylestrong{Synopsis}
\begin{quote}

\sphinxAtStartPar
\sphinxcode{\sphinxupquote{void PushBack(const GenConstr \&constr)}}
\end{quote}

\sphinxAtStartPar
\sphinxstylestrong{Arguments}
\begin{quote}

\sphinxAtStartPar
\sphinxcode{\sphinxupquote{constr}}: a general constraint object.
\end{quote}
\end{quote}

\subsubsection{GenConstrArray::Reserve()}
\label{\detokenize{cppapi/genconstrarray:genconstrarray-reserve}}\begin{quote}

\sphinxAtStartPar
Reserve capacity to contain at least n items.

\sphinxAtStartPar
\sphinxstylestrong{Synopsis}
\begin{quote}

\sphinxAtStartPar
\sphinxcode{\sphinxupquote{void Reserve(int n)}}
\end{quote}

\sphinxAtStartPar
\sphinxstylestrong{Arguments}
\begin{quote}

\sphinxAtStartPar
\sphinxcode{\sphinxupquote{n}}: minimum capacity for general constraint objects.
\end{quote}
\end{quote}

\subsubsection{GenConstrArray::Size()}
\label{\detokenize{cppapi/genconstrarray:genconstrarray-size}}\begin{quote}

\sphinxAtStartPar
Get the number of general constraint objects.

\sphinxAtStartPar
\sphinxstylestrong{Synopsis}
\begin{quote}

\sphinxAtStartPar
\sphinxcode{\sphinxupquote{int Size()}}
\end{quote}

\sphinxAtStartPar
\sphinxstylestrong{Return}
\begin{quote}

\sphinxAtStartPar
number of general constraint objects.
\end{quote}
\end{quote}

\subsection{GenConstrBuilder}
\label{\detokenize{cppapiref:genconstrbuilder}}\label{\detokenize{cppapiref:chapcppapiref-genconstrbuilder}}
\sphinxAtStartPar
COPT general constraint builder object. To help building a general
constraint, given a binary variable and associated value, a linear
expression and constraint sense, Cardinal Optimizer provides C++
GenConstrBuilder class, which defines the following methods.

\sphinxstepscope

\subsubsection{GenConstrBuilder::GetBinVal()}
\label{\detokenize{cppapi/genconstrbuilder:genconstrbuilder-getbinval}}\label{\detokenize{cppapi/genconstrbuilder::doc}}\begin{quote}

\sphinxAtStartPar
Get binary value associated with general constraint.

\sphinxAtStartPar
\sphinxstylestrong{Synopsis}
\begin{quote}

\sphinxAtStartPar
\sphinxcode{\sphinxupquote{int GetBinVal()}}
\end{quote}

\sphinxAtStartPar
\sphinxstylestrong{Return}
\begin{quote}

\sphinxAtStartPar
binary value.
\end{quote}
\end{quote}

\subsubsection{GenConstrBuilder::GetBinVar()}
\label{\detokenize{cppapi/genconstrbuilder:genconstrbuilder-getbinvar}}\begin{quote}

\sphinxAtStartPar
Get binary variable associated with general constraint.

\sphinxAtStartPar
\sphinxstylestrong{Synopsis}
\begin{quote}

\sphinxAtStartPar
\sphinxcode{\sphinxupquote{Var GetBinVar()}}
\end{quote}

\sphinxAtStartPar
\sphinxstylestrong{Return}
\begin{quote}

\sphinxAtStartPar
binary vaiable object.
\end{quote}
\end{quote}

\subsubsection{GenConstrBuilder::GetExpr()}
\label{\detokenize{cppapi/genconstrbuilder:genconstrbuilder-getexpr}}\begin{quote}

\sphinxAtStartPar
Get expression associated with general constraint.

\sphinxAtStartPar
\sphinxstylestrong{Synopsis}
\begin{quote}

\sphinxAtStartPar
\sphinxcode{\sphinxupquote{const Expr \&GetExpr()}}
\end{quote}

\sphinxAtStartPar
\sphinxstylestrong{Return}
\begin{quote}

\sphinxAtStartPar
expression object.
\end{quote}
\end{quote}

\subsubsection{GenConstrBuilder::GetIndType()}
\label{\detokenize{cppapi/genconstrbuilder:genconstrbuilder-getindtype}}\begin{quote}

\sphinxAtStartPar
Get type of general constraint.

\sphinxAtStartPar
\sphinxstylestrong{Synopsis}
\begin{quote}

\sphinxAtStartPar
\sphinxcode{\sphinxupquote{int GetIndType()}}
\end{quote}

\sphinxAtStartPar
\sphinxstylestrong{Return}
\begin{quote}

\sphinxAtStartPar
type of GenConstr, COPT\_INDICATOR\_IF, COPT\_INDICATOR\_ONLYIF and COPT\_INDICATOR\_IFANDONLYIF.
\end{quote}
\end{quote}

\subsubsection{GenConstrBuilder::GetSense()}
\label{\detokenize{cppapi/genconstrbuilder:genconstrbuilder-getsense}}\begin{quote}

\sphinxAtStartPar
Get sense associated with general constraint.

\sphinxAtStartPar
\sphinxstylestrong{Synopsis}
\begin{quote}

\sphinxAtStartPar
\sphinxcode{\sphinxupquote{char GetSense()}}
\end{quote}

\sphinxAtStartPar
\sphinxstylestrong{Return}
\begin{quote}

\sphinxAtStartPar
constraint sense.
\end{quote}
\end{quote}

\subsubsection{GenConstrBuilder::Set()}
\label{\detokenize{cppapi/genconstrbuilder:genconstrbuilder-set}}\begin{quote}

\sphinxAtStartPar
Set binary variable, binary value, expression, sense and type of general constraint.

\sphinxAtStartPar
\sphinxstylestrong{Synopsis}
\begin{quote}

\sphinxAtStartPar
\sphinxcode{\sphinxupquote{void Set(}}
\begin{quote}

\sphinxAtStartPar
\sphinxcode{\sphinxupquote{Var bvar,}}

\sphinxAtStartPar
\sphinxcode{\sphinxupquote{int bval,}}

\sphinxAtStartPar
\sphinxcode{\sphinxupquote{const Expr \&expr,}}

\sphinxAtStartPar
\sphinxcode{\sphinxupquote{char sense,}}

\sphinxAtStartPar
\sphinxcode{\sphinxupquote{int type)}}
\end{quote}
\end{quote}

\sphinxAtStartPar
\sphinxstylestrong{Arguments}
\begin{quote}

\sphinxAtStartPar
\sphinxcode{\sphinxupquote{bvar}}: binary variable.

\sphinxAtStartPar
\sphinxcode{\sphinxupquote{bval}}: binary value.

\sphinxAtStartPar
\sphinxcode{\sphinxupquote{expr}}: expression object.

\sphinxAtStartPar
\sphinxcode{\sphinxupquote{sense}}: general constraint sense.

\sphinxAtStartPar
\sphinxcode{\sphinxupquote{type}}: type of general constraint with default value of COPT\_INDICATOR\_IF.
\end{quote}
\end{quote}

\subsection{GenConstrBuilderArray}
\label{\detokenize{cppapiref:genconstrbuilderarray}}\label{\detokenize{cppapiref:chapcppapiref-genconstrbuilderarray}}
\sphinxAtStartPar
COPT general constraint builder array object. To store and access a set of C++
{\hyperref[\detokenize{cppapiref:chapcppapiref-genconstrbuilder}]{\sphinxcrossref{\DUrole{std,std-ref}{GenConstrBuilder}}}} objects, Cardinal Optimizer provides C++
GenConstrBuilderArray class, which defines the following methods.

\sphinxstepscope

\subsubsection{GenConstrBuilderArray::GetBuilder()}
\label{\detokenize{cppapi/genconstrbuilderarray:genconstrbuilderarray-getbuilder}}\label{\detokenize{cppapi/genconstrbuilderarray::doc}}\begin{quote}

\sphinxAtStartPar
Get i\sphinxhyphen{}th general constraint builder object.

\sphinxAtStartPar
\sphinxstylestrong{Synopsis}
\begin{quote}

\sphinxAtStartPar
\sphinxcode{\sphinxupquote{GenConstrBuilder \&GetBuilder(int i)}}
\end{quote}

\sphinxAtStartPar
\sphinxstylestrong{Arguments}
\begin{quote}

\sphinxAtStartPar
\sphinxcode{\sphinxupquote{i}}: index of the general constraint builder.
\end{quote}

\sphinxAtStartPar
\sphinxstylestrong{Return}
\begin{quote}

\sphinxAtStartPar
general constraint builder object with index i.
\end{quote}
\end{quote}

\subsubsection{GenConstrBuilderArray::PushBack()}
\label{\detokenize{cppapi/genconstrbuilderarray:genconstrbuilderarray-pushback}}\begin{quote}

\sphinxAtStartPar
Add a general constraint builder object to general constraint builder array.

\sphinxAtStartPar
\sphinxstylestrong{Synopsis}
\begin{quote}

\sphinxAtStartPar
\sphinxcode{\sphinxupquote{void PushBack(const GenConstrBuilder \&builder)}}
\end{quote}

\sphinxAtStartPar
\sphinxstylestrong{Arguments}
\begin{quote}

\sphinxAtStartPar
\sphinxcode{\sphinxupquote{builder}}: a general constraint builder object.
\end{quote}
\end{quote}

\subsubsection{GenConstrBuilderArray::Size()}
\label{\detokenize{cppapi/genconstrbuilderarray:genconstrbuilderarray-size}}\begin{quote}

\sphinxAtStartPar
Get the number of general constraint builder objects.

\sphinxAtStartPar
\sphinxstylestrong{Synopsis}
\begin{quote}

\sphinxAtStartPar
\sphinxcode{\sphinxupquote{int Size()}}
\end{quote}

\sphinxAtStartPar
\sphinxstylestrong{Return}
\begin{quote}

\sphinxAtStartPar
number of general constraint builder objects.
\end{quote}
\end{quote}

\subsection{Cone}
\label{\detokenize{cppapiref:cone}}\label{\detokenize{cppapiref:chapcppapiref-cone}}
\sphinxAtStartPar
COPT cone constraint object. Cone constraints are always associated
with a particular model.  User creates a cone constraint object by adding
a cone constraint to a model, rather than by using constructor of Cone class.

\sphinxAtStartPar
A cone constraint can be regular or rotated
(\sphinxcode{\sphinxupquote{COPT\_CONE\_QUAD}} or \sphinxcode{\sphinxupquote{COPT\_CONE\_RQUAD}}).

\sphinxstepscope

\subsubsection{Cone::GetIdx()}
\label{\detokenize{cppapi/cone:cone-getidx}}\label{\detokenize{cppapi/cone::doc}}\begin{quote}

\sphinxAtStartPar
Get the index of a cone constraint.

\sphinxAtStartPar
\sphinxstylestrong{Synopsis}
\begin{quote}

\sphinxAtStartPar
\sphinxcode{\sphinxupquote{int GetIdx()}}
\end{quote}

\sphinxAtStartPar
\sphinxstylestrong{Return}
\begin{quote}

\sphinxAtStartPar
index of a cone constraint.
\end{quote}
\end{quote}

\subsubsection{Cone::Remove()}
\label{\detokenize{cppapi/cone:cone-remove}}\begin{quote}

\sphinxAtStartPar
Remove the cone constraint from model.

\sphinxAtStartPar
\sphinxstylestrong{Synopsis}
\begin{quote}

\sphinxAtStartPar
\sphinxcode{\sphinxupquote{void Remove()}}
\end{quote}
\end{quote}

\subsection{ConeArray}
\label{\detokenize{cppapiref:conearray}}\label{\detokenize{cppapiref:chapcppapiref-conearray}}
\sphinxAtStartPar
COPT cone constraint array object. To store and access a set of
C++ {\hyperref[\detokenize{cppapiref:chapcppapiref-cone}]{\sphinxcrossref{\DUrole{std,std-ref}{Cone}}}} objects, Cardinal Optimizer provides
C++ ConeArray class, which defines the following methods.

\sphinxstepscope

\subsubsection{ConeArray::GetCone()}
\label{\detokenize{cppapi/conearray:conearray-getcone}}\label{\detokenize{cppapi/conearray::doc}}\begin{quote}

\sphinxAtStartPar
Get i\sphinxhyphen{}th cone constraint object.

\sphinxAtStartPar
\sphinxstylestrong{Synopsis}
\begin{quote}

\sphinxAtStartPar
\sphinxcode{\sphinxupquote{Cone \&GetCone(int i)}}
\end{quote}

\sphinxAtStartPar
\sphinxstylestrong{Arguments}
\begin{quote}

\sphinxAtStartPar
\sphinxcode{\sphinxupquote{i}}: index of the cone constraint.
\end{quote}

\sphinxAtStartPar
\sphinxstylestrong{Return}
\begin{quote}

\sphinxAtStartPar
cone constraint object with index i.
\end{quote}
\end{quote}

\subsubsection{ConeArray::PushBack()}
\label{\detokenize{cppapi/conearray:conearray-pushback}}\begin{quote}

\sphinxAtStartPar
Add a cone constraint object to cone constraint array.

\sphinxAtStartPar
\sphinxstylestrong{Synopsis}
\begin{quote}

\sphinxAtStartPar
\sphinxcode{\sphinxupquote{void PushBack(const Cone \&cone)}}
\end{quote}

\sphinxAtStartPar
\sphinxstylestrong{Arguments}
\begin{quote}

\sphinxAtStartPar
\sphinxcode{\sphinxupquote{cone}}: a cone constraint object.
\end{quote}
\end{quote}

\subsubsection{ConeArray::Size()}
\label{\detokenize{cppapi/conearray:conearray-size}}\begin{quote}

\sphinxAtStartPar
Get the number of cone constraint objects.

\sphinxAtStartPar
\sphinxstylestrong{Synopsis}
\begin{quote}

\sphinxAtStartPar
\sphinxcode{\sphinxupquote{int Size()}}
\end{quote}

\sphinxAtStartPar
\sphinxstylestrong{Return}
\begin{quote}

\sphinxAtStartPar
number of cone constraint objects.
\end{quote}
\end{quote}

\subsection{ConeBuilder}
\label{\detokenize{cppapiref:conebuilder}}\label{\detokenize{cppapiref:chapcppapiref-conebuilder}}
\sphinxAtStartPar
COPT cone constraint builder object. To help building a cone constraint,
given the cone type and a set of variables, Cardinal Optimizer provides
C++ ConeBuilder class, which defines the following methods.

\sphinxstepscope

\subsubsection{ConeBuilder::GetSize()}
\label{\detokenize{cppapi/conebuilder:conebuilder-getsize}}\label{\detokenize{cppapi/conebuilder::doc}}\begin{quote}

\sphinxAtStartPar
Get number of vars in a cone constraint.

\sphinxAtStartPar
\sphinxstylestrong{Synopsis}
\begin{quote}

\sphinxAtStartPar
\sphinxcode{\sphinxupquote{int GetSize()}}
\end{quote}

\sphinxAtStartPar
\sphinxstylestrong{Return}
\begin{quote}

\sphinxAtStartPar
number of vars.
\end{quote}
\end{quote}

\subsubsection{ConeBuilder::GetType()}
\label{\detokenize{cppapi/conebuilder:conebuilder-gettype}}\begin{quote}

\sphinxAtStartPar
Get type of a cone constraint.

\sphinxAtStartPar
\sphinxstylestrong{Synopsis}
\begin{quote}

\sphinxAtStartPar
\sphinxcode{\sphinxupquote{int GetType()}}
\end{quote}

\sphinxAtStartPar
\sphinxstylestrong{Return}
\begin{quote}

\sphinxAtStartPar
type of a cone constraint.
\end{quote}
\end{quote}

\subsubsection{ConeBuilder::GetVar()}
\label{\detokenize{cppapi/conebuilder:conebuilder-getvar}}\begin{quote}

\sphinxAtStartPar
Get the i\sphinxhyphen{}th variable in a cone constraint.

\sphinxAtStartPar
\sphinxstylestrong{Synopsis}
\begin{quote}

\sphinxAtStartPar
\sphinxcode{\sphinxupquote{Var GetVar(int i)}}
\end{quote}

\sphinxAtStartPar
\sphinxstylestrong{Arguments}
\begin{quote}

\sphinxAtStartPar
\sphinxcode{\sphinxupquote{i}}: index of vars in a cone constraint.
\end{quote}

\sphinxAtStartPar
\sphinxstylestrong{Return}
\begin{quote}

\sphinxAtStartPar
the i\sphinxhyphen{}th variable in a cone constraint.
\end{quote}
\end{quote}

\subsubsection{ConeBuilder::GetVars()}
\label{\detokenize{cppapi/conebuilder:conebuilder-getvars}}\begin{quote}

\sphinxAtStartPar
Get all variables in a cone constraint.

\sphinxAtStartPar
\sphinxstylestrong{Synopsis}
\begin{quote}

\sphinxAtStartPar
\sphinxcode{\sphinxupquote{VarArray GetVars()}}
\end{quote}

\sphinxAtStartPar
\sphinxstylestrong{Return}
\begin{quote}

\sphinxAtStartPar
variable array object.
\end{quote}
\end{quote}

\subsubsection{ConeBuilder::Set()}
\label{\detokenize{cppapi/conebuilder:conebuilder-set}}\begin{quote}

\sphinxAtStartPar
Set variables of a cone constraint.

\sphinxAtStartPar
\sphinxstylestrong{Synopsis}
\begin{quote}

\sphinxAtStartPar
\sphinxcode{\sphinxupquote{void Set(const VarArray \&vars, int type)}}
\end{quote}

\sphinxAtStartPar
\sphinxstylestrong{Arguments}
\begin{quote}

\sphinxAtStartPar
\sphinxcode{\sphinxupquote{vars}}: variable array object.

\sphinxAtStartPar
\sphinxcode{\sphinxupquote{type}}: type of cone constraint.
\end{quote}
\end{quote}

\subsection{ConeBuilderArray}
\label{\detokenize{cppapiref:conebuilderarray}}\label{\detokenize{cppapiref:chapcppapiref-conebuilderarray}}
\sphinxAtStartPar
COPT cone constraint builder array object. To store and access a
set of C++ {\hyperref[\detokenize{cppapiref:chapcppapiref-conebuilder}]{\sphinxcrossref{\DUrole{std,std-ref}{ConeBuilder}}}} objects, Cardinal Optimizer
provides C++ ConeBuilderArray class, which defines the following methods.

\sphinxstepscope

\subsubsection{ConeBuilderArray::GetBuilder()}
\label{\detokenize{cppapi/conebuilderarray:conebuilderarray-getbuilder}}\label{\detokenize{cppapi/conebuilderarray::doc}}\begin{quote}

\sphinxAtStartPar
Get i\sphinxhyphen{}th cone constraint builder object.

\sphinxAtStartPar
\sphinxstylestrong{Synopsis}
\begin{quote}

\sphinxAtStartPar
\sphinxcode{\sphinxupquote{ConeBuilder \&GetBuilder(int i)}}
\end{quote}

\sphinxAtStartPar
\sphinxstylestrong{Arguments}
\begin{quote}

\sphinxAtStartPar
\sphinxcode{\sphinxupquote{i}}: index of the cone constraint builder.
\end{quote}

\sphinxAtStartPar
\sphinxstylestrong{Return}
\begin{quote}

\sphinxAtStartPar
cone constraint builder object with index i.
\end{quote}
\end{quote}

\subsubsection{ConeBuilderArray::PushBack()}
\label{\detokenize{cppapi/conebuilderarray:conebuilderarray-pushback}}\begin{quote}

\sphinxAtStartPar
Add a cone constraint builder object to cone constraint builder array.

\sphinxAtStartPar
\sphinxstylestrong{Synopsis}
\begin{quote}

\sphinxAtStartPar
\sphinxcode{\sphinxupquote{void PushBack(const ConeBuilder \&builder)}}
\end{quote}

\sphinxAtStartPar
\sphinxstylestrong{Arguments}
\begin{quote}

\sphinxAtStartPar
\sphinxcode{\sphinxupquote{builder}}: a cone constraint builder object.
\end{quote}
\end{quote}

\subsubsection{ConeBuilderArray::Size()}
\label{\detokenize{cppapi/conebuilderarray:conebuilderarray-size}}\begin{quote}

\sphinxAtStartPar
Get the number of cone constraint builder objects.

\sphinxAtStartPar
\sphinxstylestrong{Synopsis}
\begin{quote}

\sphinxAtStartPar
\sphinxcode{\sphinxupquote{int Size()}}
\end{quote}

\sphinxAtStartPar
\sphinxstylestrong{Return}
\begin{quote}

\sphinxAtStartPar
number of cone constraint builder objects.
\end{quote}
\end{quote}

\subsection{ExpCone}
\label{\detokenize{cppapiref:expcone}}\label{\detokenize{cppapiref:chapcppapiref-expcone}}
\sphinxAtStartPar
COPT exponential cone constraint object. ExpCone constraints are always associated
with a particular model.  User creates an exponential cone constraint object by adding
an expcone constraint to a model, rather than by using constructor of ExpCone class.

\sphinxstepscope

\subsubsection{ExpCone::GetIdx()}
\label{\detokenize{cppapi/expcone:expcone-getidx}}\label{\detokenize{cppapi/expcone::doc}}\begin{quote}

\sphinxAtStartPar
Get the index of an exponential cone constraint.

\sphinxAtStartPar
\sphinxstylestrong{Synopsis}
\begin{quote}

\sphinxAtStartPar
\sphinxcode{\sphinxupquote{int GetIdx()}}
\end{quote}

\sphinxAtStartPar
\sphinxstylestrong{Return}
\begin{quote}

\sphinxAtStartPar
index of an exponential cone constraint.
\end{quote}
\end{quote}

\subsubsection{ExpCone::Remove()}
\label{\detokenize{cppapi/expcone:expcone-remove}}\begin{quote}

\sphinxAtStartPar
Remove the exponential cone constraint from model.

\sphinxAtStartPar
\sphinxstylestrong{Synopsis}
\begin{quote}

\sphinxAtStartPar
\sphinxcode{\sphinxupquote{void Remove()}}
\end{quote}
\end{quote}

\subsection{ExpConeArray}
\label{\detokenize{cppapiref:expconearray}}\label{\detokenize{cppapiref:chapcppapiref-expconearray}}
\sphinxAtStartPar
COPT exponential cone constraint array object. To store and access a set of
C++ {\hyperref[\detokenize{cppapiref:chapcppapiref-expcone}]{\sphinxcrossref{\DUrole{std,std-ref}{ExpCone}}}} objects, Cardinal Optimizer provides
C++ ExpConeArray class, which defines the following methods.

\sphinxstepscope

\subsubsection{ExpConeArray::GetCone()}
\label{\detokenize{cppapi/expconearray:expconearray-getcone}}\label{\detokenize{cppapi/expconearray::doc}}\begin{quote}

\sphinxAtStartPar
Get i\sphinxhyphen{}th exponential cone constraint object.

\sphinxAtStartPar
\sphinxstylestrong{Synopsis}
\begin{quote}

\sphinxAtStartPar
\sphinxcode{\sphinxupquote{ExpCone \&GetCone(int i)}}
\end{quote}

\sphinxAtStartPar
\sphinxstylestrong{Arguments}
\begin{quote}

\sphinxAtStartPar
\sphinxcode{\sphinxupquote{i}}: index of the exponential cone constraint.
\end{quote}

\sphinxAtStartPar
\sphinxstylestrong{Return}
\begin{quote}

\sphinxAtStartPar
exponential cone constraint object with index i.
\end{quote}
\end{quote}

\subsubsection{ExpConeArray::PushBack()}
\label{\detokenize{cppapi/expconearray:expconearray-pushback}}\begin{quote}

\sphinxAtStartPar
Add an exponential cone constraint object to exponential cone constraint array.

\sphinxAtStartPar
\sphinxstylestrong{Synopsis}
\begin{quote}

\sphinxAtStartPar
\sphinxcode{\sphinxupquote{void PushBack(const ExpCone \&cone)}}
\end{quote}

\sphinxAtStartPar
\sphinxstylestrong{Arguments}
\begin{quote}

\sphinxAtStartPar
\sphinxcode{\sphinxupquote{cone}}: an exponential constraint object.
\end{quote}
\end{quote}

\subsubsection{ExpConeArray::Size()}
\label{\detokenize{cppapi/expconearray:expconearray-size}}\begin{quote}

\sphinxAtStartPar
Get the number of exponential cone constraint objects.

\sphinxAtStartPar
\sphinxstylestrong{Synopsis}
\begin{quote}

\sphinxAtStartPar
\sphinxcode{\sphinxupquote{int Size()}}
\end{quote}

\sphinxAtStartPar
\sphinxstylestrong{Return}
\begin{quote}

\sphinxAtStartPar
number of exponential cone constraint objects.
\end{quote}
\end{quote}

\subsection{ExpConeBuilder}
\label{\detokenize{cppapiref:expconebuilder}}\label{\detokenize{cppapiref:chapcppapiref-expconebuilder}}
\sphinxAtStartPar
COPT exponential cone constraint builder object. To help building a exponential cone constraint,
given the expcone type and a set of variables, Cardinal Optimizer provides
C++ ExpConeBuilder class, which defines the following methods.

\sphinxstepscope

\subsubsection{ExpConeBuilder::GetSize()}
\label{\detokenize{cppapi/expconebuilder:expconebuilder-getsize}}\label{\detokenize{cppapi/expconebuilder::doc}}\begin{quote}

\sphinxAtStartPar
Get number of vars in an exponential cone constraint.

\sphinxAtStartPar
\sphinxstylestrong{Synopsis}
\begin{quote}

\sphinxAtStartPar
\sphinxcode{\sphinxupquote{int GetSize()}}
\end{quote}

\sphinxAtStartPar
\sphinxstylestrong{Return}
\begin{quote}

\sphinxAtStartPar
number of vars.
\end{quote}
\end{quote}

\subsubsection{ExpConeBuilder::GetType()}
\label{\detokenize{cppapi/expconebuilder:expconebuilder-gettype}}\begin{quote}

\sphinxAtStartPar
Get type of an exponential cone constraint.

\sphinxAtStartPar
\sphinxstylestrong{Synopsis}
\begin{quote}

\sphinxAtStartPar
\sphinxcode{\sphinxupquote{int GetType()}}
\end{quote}

\sphinxAtStartPar
\sphinxstylestrong{Return}
\begin{quote}

\sphinxAtStartPar
type of an exponential cone constraint.
\end{quote}
\end{quote}

\subsubsection{ExpConeBuilder::GetVar()}
\label{\detokenize{cppapi/expconebuilder:expconebuilder-getvar}}\begin{quote}

\sphinxAtStartPar
Get the i\sphinxhyphen{}th variable in an exponential cone constraint.

\sphinxAtStartPar
\sphinxstylestrong{Synopsis}
\begin{quote}

\sphinxAtStartPar
\sphinxcode{\sphinxupquote{Var GetVar(int i)}}
\end{quote}

\sphinxAtStartPar
\sphinxstylestrong{Arguments}
\begin{quote}

\sphinxAtStartPar
\sphinxcode{\sphinxupquote{i}}: index of vars in an exponential cone constraint.
\end{quote}

\sphinxAtStartPar
\sphinxstylestrong{Return}
\begin{quote}

\sphinxAtStartPar
the i\sphinxhyphen{}th variable in an exponential cone constraint.
\end{quote}
\end{quote}

\subsubsection{ExpConeBuilder::GetVars()}
\label{\detokenize{cppapi/expconebuilder:expconebuilder-getvars}}\begin{quote}

\sphinxAtStartPar
Get all variables in an exponential cone constraint.

\sphinxAtStartPar
\sphinxstylestrong{Synopsis}
\begin{quote}

\sphinxAtStartPar
\sphinxcode{\sphinxupquote{VarArray GetVars()}}
\end{quote}

\sphinxAtStartPar
\sphinxstylestrong{Return}
\begin{quote}

\sphinxAtStartPar
variable array object.
\end{quote}
\end{quote}

\subsubsection{ExpConeBuilder::Set()}
\label{\detokenize{cppapi/expconebuilder:expconebuilder-set}}\begin{quote}

\sphinxAtStartPar
Set variables of an exponential cone constraint.

\sphinxAtStartPar
\sphinxstylestrong{Synopsis}
\begin{quote}

\sphinxAtStartPar
\sphinxcode{\sphinxupquote{void Set(const VarArray \&vars, int type)}}
\end{quote}

\sphinxAtStartPar
\sphinxstylestrong{Arguments}
\begin{quote}

\sphinxAtStartPar
\sphinxcode{\sphinxupquote{vars}}: variable array object.

\sphinxAtStartPar
\sphinxcode{\sphinxupquote{type}}: type of exponential cone constraint.
\end{quote}
\end{quote}

\subsection{ExpConeBuilderArray}
\label{\detokenize{cppapiref:expconebuilderarray}}\label{\detokenize{cppapiref:chapcppapiref-expconebuilderarray}}
\sphinxAtStartPar
COPT exponential cone constraint builder array object. To store and access a
set of C++ {\hyperref[\detokenize{cppapiref:chapcppapiref-expconebuilder}]{\sphinxcrossref{\DUrole{std,std-ref}{ExpConeBuilder}}}} objects, Cardinal Optimizer
provides C++ ExpConeBuilderArray class, which defines the following methods.

\sphinxstepscope

\subsubsection{ExpConeBuilderArray::GetBuilder()}
\label{\detokenize{cppapi/expconebuilderarray:expconebuilderarray-getbuilder}}\label{\detokenize{cppapi/expconebuilderarray::doc}}\begin{quote}

\sphinxAtStartPar
Get i\sphinxhyphen{}th exponential cone constraint builder object.

\sphinxAtStartPar
\sphinxstylestrong{Synopsis}
\begin{quote}

\sphinxAtStartPar
\sphinxcode{\sphinxupquote{ExpConeBuilder \&GetBuilder(int i)}}
\end{quote}

\sphinxAtStartPar
\sphinxstylestrong{Arguments}
\begin{quote}

\sphinxAtStartPar
\sphinxcode{\sphinxupquote{i}}: index of the exponential cone constraint builder.
\end{quote}

\sphinxAtStartPar
\sphinxstylestrong{Return}
\begin{quote}

\sphinxAtStartPar
exponential cone constraint builder object with index i.
\end{quote}
\end{quote}

\subsubsection{ExpConeBuilderArray::PushBack()}
\label{\detokenize{cppapi/expconebuilderarray:expconebuilderarray-pushback}}\begin{quote}

\sphinxAtStartPar
Add an exponential cone constraint builder object to exponential cone constraint builder array.

\sphinxAtStartPar
\sphinxstylestrong{Synopsis}
\begin{quote}

\sphinxAtStartPar
\sphinxcode{\sphinxupquote{void PushBack(const ExpConeBuilder \&builder)}}
\end{quote}

\sphinxAtStartPar
\sphinxstylestrong{Arguments}
\begin{quote}

\sphinxAtStartPar
\sphinxcode{\sphinxupquote{builder}}: an exponential cone constraint builder object.
\end{quote}
\end{quote}

\subsubsection{ExpConeBuilderArray::Size()}
\label{\detokenize{cppapi/expconebuilderarray:expconebuilderarray-size}}\begin{quote}

\sphinxAtStartPar
Get the number of exponential cone constraint builder objects.

\sphinxAtStartPar
\sphinxstylestrong{Synopsis}
\begin{quote}

\sphinxAtStartPar
\sphinxcode{\sphinxupquote{int Size()}}
\end{quote}

\sphinxAtStartPar
\sphinxstylestrong{Return}
\begin{quote}

\sphinxAtStartPar
number of exponential cone constraint builder objects.
\end{quote}
\end{quote}

\subsection{AffineCone}
\label{\detokenize{cppapiref:affinecone}}\label{\detokenize{cppapiref:chapcppapiref-affinecone}}
\sphinxAtStartPar
The \sphinxtitleref{AffineCone} class in COPT encapsulates operations related to affine cone constraints.
The following methods are provided:

\sphinxstepscope

\subsubsection{AffineCone::GetIdx()}
\label{\detokenize{cppapi/affinecone:affinecone-getidx}}\label{\detokenize{cppapi/affinecone::doc}}\begin{quote}

\sphinxAtStartPar
Get the index of an affine cone constraint.

\sphinxAtStartPar
\sphinxstylestrong{Synopsis}
\begin{quote}

\sphinxAtStartPar
\sphinxcode{\sphinxupquote{int GetIdx()}}
\end{quote}

\sphinxAtStartPar
\sphinxstylestrong{Return}
\begin{quote}

\sphinxAtStartPar
index of an affine cone constraint.
\end{quote}
\end{quote}

\subsubsection{AffineCone::GetName()}
\label{\detokenize{cppapi/affinecone:affinecone-getname}}\begin{quote}

\sphinxAtStartPar
Get name of affine cone constraint.

\sphinxAtStartPar
\sphinxstylestrong{Synopsis}
\begin{quote}

\sphinxAtStartPar
\sphinxcode{\sphinxupquote{const char *GetName()}}
\end{quote}

\sphinxAtStartPar
\sphinxstylestrong{Return}
\begin{quote}

\sphinxAtStartPar
name of affine cone constraint.
\end{quote}
\end{quote}

\subsubsection{AffineCone::Remove()}
\label{\detokenize{cppapi/affinecone:affinecone-remove}}\begin{quote}

\sphinxAtStartPar
Remove the affine cone constraint from model.

\sphinxAtStartPar
\sphinxstylestrong{Synopsis}
\begin{quote}

\sphinxAtStartPar
\sphinxcode{\sphinxupquote{void Remove()}}
\end{quote}
\end{quote}

\subsubsection{AffineCone::SetName()}
\label{\detokenize{cppapi/affinecone:affinecone-setname}}\begin{quote}

\sphinxAtStartPar
Set name of affine cone constraint.

\sphinxAtStartPar
\sphinxstylestrong{Synopsis}
\begin{quote}

\sphinxAtStartPar
\sphinxcode{\sphinxupquote{void SetName(const char *szName)}}
\end{quote}

\sphinxAtStartPar
\sphinxstylestrong{Arguments}
\begin{quote}

\sphinxAtStartPar
\sphinxcode{\sphinxupquote{szName}}: name of affine cone constraint.
\end{quote}
\end{quote}

\subsection{AffineConeArray}
\label{\detokenize{cppapiref:affineconearray}}\label{\detokenize{cppapiref:chapcppapiref-affineconearray}}
\sphinxAtStartPar
To facilitate user operations on a group of C++ {\hyperref[\detokenize{cppapiref:chapcppapiref-affinecone}]{\sphinxcrossref{\DUrole{std,std-ref}{AffineCone}}}} objects,
the C++ interface of COPT introduces the \sphinxtitleref{AffineConeArray} class.
The following methods are provided:

\sphinxstepscope

\subsubsection{AffineConeArray::GetCone()}
\label{\detokenize{cppapi/affineconearray:affineconearray-getcone}}\label{\detokenize{cppapi/affineconearray::doc}}\begin{quote}

\sphinxAtStartPar
Get i\sphinxhyphen{}th affine cone constraint object.

\sphinxAtStartPar
\sphinxstylestrong{Synopsis}
\begin{quote}

\sphinxAtStartPar
\sphinxcode{\sphinxupquote{AffineCone \&GetCone(int i)}}
\end{quote}

\sphinxAtStartPar
\sphinxstylestrong{Arguments}
\begin{quote}

\sphinxAtStartPar
\sphinxcode{\sphinxupquote{i}}: index of the affine cone constraint.
\end{quote}

\sphinxAtStartPar
\sphinxstylestrong{Return}
\begin{quote}

\sphinxAtStartPar
affine cone constraint object with index i.
\end{quote}
\end{quote}

\subsubsection{AffineConeArray::PushBack()}
\label{\detokenize{cppapi/affineconearray:affineconearray-pushback}}\begin{quote}

\sphinxAtStartPar
Add an affine cone constraint object to affine cone constraint array.

\sphinxAtStartPar
\sphinxstylestrong{Synopsis}
\begin{quote}

\sphinxAtStartPar
\sphinxcode{\sphinxupquote{void PushBack(const AffineCone \&cone)}}
\end{quote}

\sphinxAtStartPar
\sphinxstylestrong{Arguments}
\begin{quote}

\sphinxAtStartPar
\sphinxcode{\sphinxupquote{cone}}: an affine constraint object.
\end{quote}
\end{quote}

\subsubsection{AffineConeArray::Reserve()}
\label{\detokenize{cppapi/affineconearray:affineconearray-reserve}}\begin{quote}

\sphinxAtStartPar
Reserve capacity to contain at least n items.

\sphinxAtStartPar
\sphinxstylestrong{Synopsis}
\begin{quote}

\sphinxAtStartPar
\sphinxcode{\sphinxupquote{void Reserve(int n)}}
\end{quote}

\sphinxAtStartPar
\sphinxstylestrong{Arguments}
\begin{quote}

\sphinxAtStartPar
\sphinxcode{\sphinxupquote{n}}: minimum capacity for affine cone constraint objects.
\end{quote}
\end{quote}

\subsubsection{AffineConeArray::Size()}
\label{\detokenize{cppapi/affineconearray:affineconearray-size}}\begin{quote}

\sphinxAtStartPar
Get the number of affine cone constraint objects.

\sphinxAtStartPar
\sphinxstylestrong{Synopsis}
\begin{quote}

\sphinxAtStartPar
\sphinxcode{\sphinxupquote{int Size()}}
\end{quote}

\sphinxAtStartPar
\sphinxstylestrong{Return}
\begin{quote}

\sphinxAtStartPar
number of affine cone constraint objects.
\end{quote}
\end{quote}

\subsection{AffineConeBuilder}
\label{\detokenize{cppapiref:affineconebuilder}}\label{\detokenize{cppapiref:chapcppapiref-affineconebuilder}}
\sphinxAtStartPar
The \sphinxtitleref{AffineConeBuilder} class in COPT encapsulates the builder for constructing affine cone constraints.
The following methods are provided:

\sphinxstepscope

\subsubsection{AffineConeBuilder::GetExpr()}
\label{\detokenize{cppapi/affineconebuilder:affineconebuilder-getexpr}}\label{\detokenize{cppapi/affineconebuilder::doc}}\begin{quote}

\sphinxAtStartPar
Get the i\sphinxhyphen{}th linear expression in an affine cone constraint.

\sphinxAtStartPar
\sphinxstylestrong{Synopsis}
\begin{quote}

\sphinxAtStartPar
\sphinxcode{\sphinxupquote{const Expr \&GetExpr(int i)}}
\end{quote}

\sphinxAtStartPar
\sphinxstylestrong{Arguments}
\begin{quote}

\sphinxAtStartPar
\sphinxcode{\sphinxupquote{i}}: index of linear expression in an affine cone constraint.
\end{quote}

\sphinxAtStartPar
\sphinxstylestrong{Return}
\begin{quote}

\sphinxAtStartPar
the i\sphinxhyphen{}th linear expression in an affine cone constraint.
\end{quote}
\end{quote}

\subsubsection{AffineConeBuilder::GetExprs()}
\label{\detokenize{cppapi/affineconebuilder:affineconebuilder-getexprs}}\begin{quote}

\sphinxAtStartPar
Get all linear expressions in an affine cone constraint.

\sphinxAtStartPar
\sphinxstylestrong{Synopsis}
\begin{quote}

\sphinxAtStartPar
\sphinxcode{\sphinxupquote{MLinExpr\textless{}1\textgreater{} GetExprs()}}
\end{quote}

\sphinxAtStartPar
\sphinxstylestrong{Return}
\begin{quote}

\sphinxAtStartPar
MLinExpr object.
\end{quote}
\end{quote}

\subsubsection{AffineConeBuilder::GetPsdExpr()}
\label{\detokenize{cppapi/affineconebuilder:affineconebuilder-getpsdexpr}}\begin{quote}

\sphinxAtStartPar
Get the i\sphinxhyphen{}th PSD expression in an affine cone constraint.

\sphinxAtStartPar
\sphinxstylestrong{Synopsis}
\begin{quote}

\sphinxAtStartPar
\sphinxcode{\sphinxupquote{const PsdExpr \&GetPsdExpr(int i)}}
\end{quote}

\sphinxAtStartPar
\sphinxstylestrong{Arguments}
\begin{quote}

\sphinxAtStartPar
\sphinxcode{\sphinxupquote{i}}: index of PSD expression in an affine cone constraint.
\end{quote}

\sphinxAtStartPar
\sphinxstylestrong{Return}
\begin{quote}

\sphinxAtStartPar
the i\sphinxhyphen{}th PSD expression in an affine cone constraint.
\end{quote}
\end{quote}

\subsubsection{AffineConeBuilder::GetPsdExprs()}
\label{\detokenize{cppapi/affineconebuilder:affineconebuilder-getpsdexprs}}\begin{quote}

\sphinxAtStartPar
Get all PSD expressions in an affine cone constraint.

\sphinxAtStartPar
\sphinxstylestrong{Synopsis}
\begin{quote}

\sphinxAtStartPar
\sphinxcode{\sphinxupquote{const MPsdExpr\textless{}1\textgreater{} \&GetPsdExprs()}}
\end{quote}

\sphinxAtStartPar
\sphinxstylestrong{Return}
\begin{quote}

\sphinxAtStartPar
MPsdExpr object.
\end{quote}
\end{quote}

\subsubsection{AffineConeBuilder::GetSize()}
\label{\detokenize{cppapi/affineconebuilder:affineconebuilder-getsize}}\begin{quote}

\sphinxAtStartPar
Get number of expressions in an affine cone constraint.

\sphinxAtStartPar
\sphinxstylestrong{Synopsis}
\begin{quote}

\sphinxAtStartPar
\sphinxcode{\sphinxupquote{int GetSize()}}
\end{quote}

\sphinxAtStartPar
\sphinxstylestrong{Return}
\begin{quote}

\sphinxAtStartPar
number of expressions.
\end{quote}
\end{quote}

\subsubsection{AffineConeBuilder::GetType()}
\label{\detokenize{cppapi/affineconebuilder:affineconebuilder-gettype}}\begin{quote}

\sphinxAtStartPar
Get type of an affine cone constraint.

\sphinxAtStartPar
\sphinxstylestrong{Synopsis}
\begin{quote}

\sphinxAtStartPar
\sphinxcode{\sphinxupquote{int GetType()}}
\end{quote}

\sphinxAtStartPar
\sphinxstylestrong{Return}
\begin{quote}

\sphinxAtStartPar
type of an affine cone constraint.
\end{quote}
\end{quote}

\subsubsection{AffineConeBuilder::HasPsdTerm()}
\label{\detokenize{cppapi/affineconebuilder:affineconebuilder-haspsdterm}}\begin{quote}

\sphinxAtStartPar
Check whether affine cone has PSD terms.

\sphinxAtStartPar
\sphinxstylestrong{Synopsis}
\begin{quote}

\sphinxAtStartPar
\sphinxcode{\sphinxupquote{bool HasPsdTerm()}}
\end{quote}

\sphinxAtStartPar
\sphinxstylestrong{Return}
\begin{quote}

\sphinxAtStartPar
flag to indicate whether affine cone has PSD terms.
\end{quote}
\end{quote}

\subsubsection{AffineConeBuilder::Set()}
\label{\detokenize{cppapi/affineconebuilder:affineconebuilder-set}}\begin{quote}

\sphinxAtStartPar
Set linear expressions of an affine cone constraint.

\sphinxAtStartPar
\sphinxstylestrong{Synopsis}
\begin{quote}

\sphinxAtStartPar
\sphinxcode{\sphinxupquote{void Set(const MLinExpr\textless{}1\textgreater{} \&exprs, int type)}}
\end{quote}

\sphinxAtStartPar
\sphinxstylestrong{Arguments}
\begin{quote}

\sphinxAtStartPar
\sphinxcode{\sphinxupquote{exprs}}: 1\sphinxhyphen{}dimensional MLinExpr object.

\sphinxAtStartPar
\sphinxcode{\sphinxupquote{type}}: type of affine cone constraint.
\end{quote}
\end{quote}

\subsubsection{AffineConeBuilder::Set()}
\label{\detokenize{cppapi/affineconebuilder:id1}}\begin{quote}

\sphinxAtStartPar
Set PSD expressions of an affine cone constraint.

\sphinxAtStartPar
\sphinxstylestrong{Synopsis}
\begin{quote}

\sphinxAtStartPar
\sphinxcode{\sphinxupquote{void Set(const MPsdExpr\textless{}1\textgreater{} \&exprs, int type)}}
\end{quote}

\sphinxAtStartPar
\sphinxstylestrong{Arguments}
\begin{quote}

\sphinxAtStartPar
\sphinxcode{\sphinxupquote{exprs}}: 1\sphinxhyphen{}dimensional MPsdExpr object.

\sphinxAtStartPar
\sphinxcode{\sphinxupquote{type}}: type of affine cone constraint.
\end{quote}
\end{quote}

\subsection{AffineConeBuilderArray}
\label{\detokenize{cppapiref:affineconebuilderarray}}\label{\detokenize{cppapiref:chapcppapiref-affineconebuilderarray}}
\sphinxAtStartPar
To facilitate user operations on a group of C++ {\hyperref[\detokenize{cppapiref:chapcppapiref-affineconebuilder}]{\sphinxcrossref{\DUrole{std,std-ref}{AffineConeBuilder}}}} objects,
the C++ interface of COPT introduces the \sphinxtitleref{AffineConeBuilderArray} class.
The following methods are provided:

\sphinxstepscope

\subsubsection{AffineConeBuilderArray::GetBuilder()}
\label{\detokenize{cppapi/affineconebuilderarray:affineconebuilderarray-getbuilder}}\label{\detokenize{cppapi/affineconebuilderarray::doc}}\begin{quote}

\sphinxAtStartPar
Get i\sphinxhyphen{}th affine cone constraint builder object.

\sphinxAtStartPar
\sphinxstylestrong{Synopsis}
\begin{quote}

\sphinxAtStartPar
\sphinxcode{\sphinxupquote{AffineConeBuilder \&GetBuilder(int i)}}
\end{quote}

\sphinxAtStartPar
\sphinxstylestrong{Arguments}
\begin{quote}

\sphinxAtStartPar
\sphinxcode{\sphinxupquote{i}}: index of the affine cone constraint builder.
\end{quote}

\sphinxAtStartPar
\sphinxstylestrong{Return}
\begin{quote}

\sphinxAtStartPar
affine cone constraint builder object with index i.
\end{quote}
\end{quote}

\subsubsection{AffineConeBuilderArray::PushBack()}
\label{\detokenize{cppapi/affineconebuilderarray:affineconebuilderarray-pushback}}\begin{quote}

\sphinxAtStartPar
Add an affine cone constraint builder object to affine cone constraint builder array.

\sphinxAtStartPar
\sphinxstylestrong{Synopsis}
\begin{quote}

\sphinxAtStartPar
\sphinxcode{\sphinxupquote{void PushBack(const AffineConeBuilder \&builder)}}
\end{quote}

\sphinxAtStartPar
\sphinxstylestrong{Arguments}
\begin{quote}

\sphinxAtStartPar
\sphinxcode{\sphinxupquote{builder}}: an affine cone constraint builder object.
\end{quote}
\end{quote}

\subsubsection{AffineConeBuilderArray::Size()}
\label{\detokenize{cppapi/affineconebuilderarray:affineconebuilderarray-size}}\begin{quote}

\sphinxAtStartPar
Get the number of affine cone constraint builder objects.

\sphinxAtStartPar
\sphinxstylestrong{Synopsis}
\begin{quote}

\sphinxAtStartPar
\sphinxcode{\sphinxupquote{int Size()}}
\end{quote}

\sphinxAtStartPar
\sphinxstylestrong{Return}
\begin{quote}

\sphinxAtStartPar
number of affine cone constraint builder objects.
\end{quote}
\end{quote}

\subsection{QuadExpr}
\label{\detokenize{cppapiref:quadexpr}}\label{\detokenize{cppapiref:chapcppapiref-quadexpr}}
\sphinxAtStartPar
COPT quadratic expression object. A quadratic expression consists of
a linear expression, a list of variable pairs and associated coefficients
of quadratic terms. Quadratic expressions are used to build quadratic
constraints and objectives.

\sphinxstepscope

\subsubsection{QuadExpr::QuadExpr()}
\label{\detokenize{cppapi/quadexpr:quadexpr-quadexpr}}\label{\detokenize{cppapi/quadexpr::doc}}\begin{quote}

\sphinxAtStartPar
Constructor of a quadratic expression with a constant.

\sphinxAtStartPar
\sphinxstylestrong{Synopsis}
\begin{quote}

\sphinxAtStartPar
\sphinxcode{\sphinxupquote{QuadExpr(double constant)}}
\end{quote}

\sphinxAtStartPar
\sphinxstylestrong{Arguments}
\begin{quote}

\sphinxAtStartPar
\sphinxcode{\sphinxupquote{constant}}: constant value in quadratic expression object.
\end{quote}
\end{quote}

\subsubsection{QuadExpr::QuadExpr()}
\label{\detokenize{cppapi/quadexpr:id1}}\begin{quote}

\sphinxAtStartPar
Constructor of a quadratic expression with one term.

\sphinxAtStartPar
\sphinxstylestrong{Synopsis}
\begin{quote}

\sphinxAtStartPar
\sphinxcode{\sphinxupquote{QuadExpr(const Var \&var, double coeff)}}
\end{quote}

\sphinxAtStartPar
\sphinxstylestrong{Arguments}
\begin{quote}

\sphinxAtStartPar
\sphinxcode{\sphinxupquote{var}}: variable for the added term.

\sphinxAtStartPar
\sphinxcode{\sphinxupquote{coeff}}: coefficent for the added term.
\end{quote}
\end{quote}

\subsubsection{QuadExpr::QuadExpr()}
\label{\detokenize{cppapi/quadexpr:id2}}\begin{quote}

\sphinxAtStartPar
Constructor of a quadratic expression with a linear expression.

\sphinxAtStartPar
\sphinxstylestrong{Synopsis}
\begin{quote}

\sphinxAtStartPar
\sphinxcode{\sphinxupquote{QuadExpr(const Expr \&expr)}}
\end{quote}

\sphinxAtStartPar
\sphinxstylestrong{Arguments}
\begin{quote}

\sphinxAtStartPar
\sphinxcode{\sphinxupquote{expr}}: input linear expression.
\end{quote}
\end{quote}

\subsubsection{QuadExpr::QuadExpr()}
\label{\detokenize{cppapi/quadexpr:id3}}\begin{quote}

\sphinxAtStartPar
Constructor of a quadratic expression with two linear expression.

\sphinxAtStartPar
\sphinxstylestrong{Synopsis}
\begin{quote}

\sphinxAtStartPar
\sphinxcode{\sphinxupquote{QuadExpr(const Expr \&expr, const Var \&var)}}
\end{quote}

\sphinxAtStartPar
\sphinxstylestrong{Arguments}
\begin{quote}

\sphinxAtStartPar
\sphinxcode{\sphinxupquote{expr}}: one linear expression.

\sphinxAtStartPar
\sphinxcode{\sphinxupquote{var}}: another variable.
\end{quote}
\end{quote}

\subsubsection{QuadExpr::QuadExpr()}
\label{\detokenize{cppapi/quadexpr:id4}}\begin{quote}

\sphinxAtStartPar
Constructor of a quadratic expression with two linear expression.

\sphinxAtStartPar
\sphinxstylestrong{Synopsis}
\begin{quote}

\sphinxAtStartPar
\sphinxcode{\sphinxupquote{QuadExpr(const Expr \&left, const Expr \&right)}}
\end{quote}

\sphinxAtStartPar
\sphinxstylestrong{Arguments}
\begin{quote}

\sphinxAtStartPar
\sphinxcode{\sphinxupquote{left}}: one linear expression.

\sphinxAtStartPar
\sphinxcode{\sphinxupquote{right}}: another linear expression.
\end{quote}
\end{quote}

\subsubsection{QuadExpr::AddConstant()}
\label{\detokenize{cppapi/quadexpr:quadexpr-addconstant}}\begin{quote}

\sphinxAtStartPar
Add constant for the expression.

\sphinxAtStartPar
\sphinxstylestrong{Synopsis}
\begin{quote}

\sphinxAtStartPar
\sphinxcode{\sphinxupquote{void AddConstant(double constant)}}
\end{quote}

\sphinxAtStartPar
\sphinxstylestrong{Arguments}
\begin{quote}

\sphinxAtStartPar
\sphinxcode{\sphinxupquote{constant}}: the value of the constant.
\end{quote}
\end{quote}

\subsubsection{QuadExpr::AddLinExpr()}
\label{\detokenize{cppapi/quadexpr:quadexpr-addlinexpr}}\begin{quote}

\sphinxAtStartPar
Add a linear expression to self.

\sphinxAtStartPar
\sphinxstylestrong{Synopsis}
\begin{quote}

\sphinxAtStartPar
\sphinxcode{\sphinxupquote{void AddLinExpr(const Expr \&expr, double mult)}}
\end{quote}

\sphinxAtStartPar
\sphinxstylestrong{Arguments}
\begin{quote}

\sphinxAtStartPar
\sphinxcode{\sphinxupquote{expr}}: linear expression to be added.

\sphinxAtStartPar
\sphinxcode{\sphinxupquote{mult}}: optional, constant multiplier, default value is 1.0.
\end{quote}
\end{quote}

\subsubsection{QuadExpr::AddQuadExpr()}
\label{\detokenize{cppapi/quadexpr:quadexpr-addquadexpr}}\begin{quote}

\sphinxAtStartPar
Add a quadratic expression to self.

\sphinxAtStartPar
\sphinxstylestrong{Synopsis}
\begin{quote}

\sphinxAtStartPar
\sphinxcode{\sphinxupquote{void AddQuadExpr(const QuadExpr \&expr, double mult)}}
\end{quote}

\sphinxAtStartPar
\sphinxstylestrong{Arguments}
\begin{quote}

\sphinxAtStartPar
\sphinxcode{\sphinxupquote{expr}}: quadratic expression to be added.

\sphinxAtStartPar
\sphinxcode{\sphinxupquote{mult}}: optional, constant multiplier, default value is 1.0.
\end{quote}
\end{quote}

\subsubsection{QuadExpr::AddTerm()}
\label{\detokenize{cppapi/quadexpr:quadexpr-addterm}}\begin{quote}

\sphinxAtStartPar
Add a linear term to expression object.

\sphinxAtStartPar
\sphinxstylestrong{Synopsis}
\begin{quote}

\sphinxAtStartPar
\sphinxcode{\sphinxupquote{void AddTerm(const Var \&var, double coeff)}}
\end{quote}

\sphinxAtStartPar
\sphinxstylestrong{Arguments}
\begin{quote}

\sphinxAtStartPar
\sphinxcode{\sphinxupquote{var}}: variable of new linear term.

\sphinxAtStartPar
\sphinxcode{\sphinxupquote{coeff}}: coefficient of new linear term.
\end{quote}
\end{quote}

\subsubsection{QuadExpr::AddTerm()}
\label{\detokenize{cppapi/quadexpr:id5}}\begin{quote}

\sphinxAtStartPar
Add a quadratic term to expression object.

\sphinxAtStartPar
\sphinxstylestrong{Synopsis}
\begin{quote}

\sphinxAtStartPar
\sphinxcode{\sphinxupquote{void AddTerm(}}
\begin{quote}

\sphinxAtStartPar
\sphinxcode{\sphinxupquote{const Var \&var1,}}

\sphinxAtStartPar
\sphinxcode{\sphinxupquote{const Var \&var2,}}

\sphinxAtStartPar
\sphinxcode{\sphinxupquote{double coeff)}}
\end{quote}
\end{quote}

\sphinxAtStartPar
\sphinxstylestrong{Arguments}
\begin{quote}

\sphinxAtStartPar
\sphinxcode{\sphinxupquote{var1}}: first variable of new quadratic term.

\sphinxAtStartPar
\sphinxcode{\sphinxupquote{var2}}: second variable of new quadratic term.

\sphinxAtStartPar
\sphinxcode{\sphinxupquote{coeff}}: coefficient of new quadratic term.
\end{quote}
\end{quote}

\subsubsection{QuadExpr::AddTerms()}
\label{\detokenize{cppapi/quadexpr:quadexpr-addterms}}\begin{quote}

\sphinxAtStartPar
Add linear terms to expression object.

\sphinxAtStartPar
\sphinxstylestrong{Synopsis}
\begin{quote}

\sphinxAtStartPar
\sphinxcode{\sphinxupquote{int AddTerms(}}
\begin{quote}

\sphinxAtStartPar
\sphinxcode{\sphinxupquote{const VarArray \&vars,}}

\sphinxAtStartPar
\sphinxcode{\sphinxupquote{double *pCoeff,}}

\sphinxAtStartPar
\sphinxcode{\sphinxupquote{int len)}}
\end{quote}
\end{quote}

\sphinxAtStartPar
\sphinxstylestrong{Arguments}
\begin{quote}

\sphinxAtStartPar
\sphinxcode{\sphinxupquote{vars}}: variables for added linear terms.

\sphinxAtStartPar
\sphinxcode{\sphinxupquote{pCoeff}}: coefficient array for added linear terms.

\sphinxAtStartPar
\sphinxcode{\sphinxupquote{len}}: length of coefficient array.
\end{quote}

\sphinxAtStartPar
\sphinxstylestrong{Return}
\begin{quote}

\sphinxAtStartPar
number of added linear terms.
\end{quote}
\end{quote}

\subsubsection{QuadExpr::AddTerms()}
\label{\detokenize{cppapi/quadexpr:id6}}\begin{quote}

\sphinxAtStartPar
Add quadratic terms to expression object.

\sphinxAtStartPar
\sphinxstylestrong{Synopsis}
\begin{quote}

\sphinxAtStartPar
\sphinxcode{\sphinxupquote{int AddTerms(}}
\begin{quote}

\sphinxAtStartPar
\sphinxcode{\sphinxupquote{const VarArray \&vars1,}}

\sphinxAtStartPar
\sphinxcode{\sphinxupquote{const VarArray \&vars2,}}

\sphinxAtStartPar
\sphinxcode{\sphinxupquote{double *pCoeff,}}

\sphinxAtStartPar
\sphinxcode{\sphinxupquote{int len)}}
\end{quote}
\end{quote}

\sphinxAtStartPar
\sphinxstylestrong{Arguments}
\begin{quote}

\sphinxAtStartPar
\sphinxcode{\sphinxupquote{vars1}}: first set of variables for added quadratic terms.

\sphinxAtStartPar
\sphinxcode{\sphinxupquote{vars2}}: second set of variables for added quadratic terms.

\sphinxAtStartPar
\sphinxcode{\sphinxupquote{pCoeff}}: coefficient array for added quadratic terms.

\sphinxAtStartPar
\sphinxcode{\sphinxupquote{len}}: length of coefficient array.
\end{quote}

\sphinxAtStartPar
\sphinxstylestrong{Return}
\begin{quote}

\sphinxAtStartPar
number of added quadratic terms.
\end{quote}
\end{quote}

\subsubsection{QuadExpr::Clear()}
\label{\detokenize{cppapi/quadexpr:quadexpr-clear}}\begin{quote}

\sphinxAtStartPar
Clear quadratic expression object.

\sphinxAtStartPar
\sphinxstylestrong{Synopsis}
\begin{quote}

\sphinxAtStartPar
\sphinxcode{\sphinxupquote{void Clear()}}
\end{quote}
\end{quote}

\subsubsection{QuadExpr::Clone()}
\label{\detokenize{cppapi/quadexpr:quadexpr-clone}}\begin{quote}

\sphinxAtStartPar
Deep copy quadratic expression object.

\sphinxAtStartPar
\sphinxstylestrong{Synopsis}
\begin{quote}

\sphinxAtStartPar
\sphinxcode{\sphinxupquote{QuadExpr Clone()}}
\end{quote}

\sphinxAtStartPar
\sphinxstylestrong{Return}
\begin{quote}

\sphinxAtStartPar
cloned quadratic expression object.
\end{quote}
\end{quote}

\subsubsection{QuadExpr::Evaluate()}
\label{\detokenize{cppapi/quadexpr:quadexpr-evaluate}}\begin{quote}

\sphinxAtStartPar
Evaluate quadratic expression after solving

\sphinxAtStartPar
\sphinxstylestrong{Synopsis}
\begin{quote}

\sphinxAtStartPar
\sphinxcode{\sphinxupquote{double Evaluate()}}
\end{quote}

\sphinxAtStartPar
\sphinxstylestrong{Return}
\begin{quote}

\sphinxAtStartPar
value of quadratic expression
\end{quote}
\end{quote}

\subsubsection{QuadExpr::GetCoeff()}
\label{\detokenize{cppapi/quadexpr:quadexpr-getcoeff}}\begin{quote}

\sphinxAtStartPar
Get coefficient from the i\sphinxhyphen{}th term in quadratic expression.

\sphinxAtStartPar
\sphinxstylestrong{Synopsis}
\begin{quote}

\sphinxAtStartPar
\sphinxcode{\sphinxupquote{double GetCoeff(int i)}}
\end{quote}

\sphinxAtStartPar
\sphinxstylestrong{Arguments}
\begin{quote}

\sphinxAtStartPar
\sphinxcode{\sphinxupquote{i}}: index of the quadratic term.
\end{quote}

\sphinxAtStartPar
\sphinxstylestrong{Return}
\begin{quote}

\sphinxAtStartPar
coefficient of the i\sphinxhyphen{}th quadratic term in quadratic expression object.
\end{quote}
\end{quote}

\subsubsection{QuadExpr::GetConstant()}
\label{\detokenize{cppapi/quadexpr:quadexpr-getconstant}}\begin{quote}

\sphinxAtStartPar
Get constant in quadratic expression.

\sphinxAtStartPar
\sphinxstylestrong{Synopsis}
\begin{quote}

\sphinxAtStartPar
\sphinxcode{\sphinxupquote{double GetConstant()}}
\end{quote}

\sphinxAtStartPar
\sphinxstylestrong{Return}
\begin{quote}

\sphinxAtStartPar
constant in quadratic expression.
\end{quote}
\end{quote}

\subsubsection{QuadExpr::GetLinExpr()}
\label{\detokenize{cppapi/quadexpr:quadexpr-getlinexpr}}\begin{quote}

\sphinxAtStartPar
Get linear expression in quadratic expression.

\sphinxAtStartPar
\sphinxstylestrong{Synopsis}
\begin{quote}

\sphinxAtStartPar
\sphinxcode{\sphinxupquote{Expr \&GetLinExpr()}}
\end{quote}

\sphinxAtStartPar
\sphinxstylestrong{Return}
\begin{quote}

\sphinxAtStartPar
linear expression object.
\end{quote}
\end{quote}

\subsubsection{QuadExpr::GetVar1()}
\label{\detokenize{cppapi/quadexpr:quadexpr-getvar1}}\begin{quote}

\sphinxAtStartPar
Get the first variable from the i\sphinxhyphen{}th term in quadratic expression.

\sphinxAtStartPar
\sphinxstylestrong{Synopsis}
\begin{quote}

\sphinxAtStartPar
\sphinxcode{\sphinxupquote{Var \&GetVar1(int i)}}
\end{quote}

\sphinxAtStartPar
\sphinxstylestrong{Arguments}
\begin{quote}

\sphinxAtStartPar
\sphinxcode{\sphinxupquote{i}}: index of the term.
\end{quote}

\sphinxAtStartPar
\sphinxstylestrong{Return}
\begin{quote}

\sphinxAtStartPar
the first variable of the i\sphinxhyphen{}th term in quadratic expression object.
\end{quote}
\end{quote}

\subsubsection{QuadExpr::GetVar2()}
\label{\detokenize{cppapi/quadexpr:quadexpr-getvar2}}\begin{quote}

\sphinxAtStartPar
Get the second variable from the i\sphinxhyphen{}th term in quadratic expression.

\sphinxAtStartPar
\sphinxstylestrong{Synopsis}
\begin{quote}

\sphinxAtStartPar
\sphinxcode{\sphinxupquote{Var \&GetVar2(int i)}}
\end{quote}

\sphinxAtStartPar
\sphinxstylestrong{Arguments}
\begin{quote}

\sphinxAtStartPar
\sphinxcode{\sphinxupquote{i}}: index of the term.
\end{quote}

\sphinxAtStartPar
\sphinxstylestrong{Return}
\begin{quote}

\sphinxAtStartPar
the second variable of the i\sphinxhyphen{}th term in quadratic expression object.
\end{quote}
\end{quote}

\subsubsection{QuadExpr::operator*=()}
\label{\detokenize{cppapi/quadexpr:quadexpr-operator}}\begin{quote}

\sphinxAtStartPar
Multiply a constant to self.

\sphinxAtStartPar
\sphinxstylestrong{Synopsis}
\begin{quote}

\sphinxAtStartPar
\sphinxcode{\sphinxupquote{void operator*=(double c)}}
\end{quote}

\sphinxAtStartPar
\sphinxstylestrong{Arguments}
\begin{quote}

\sphinxAtStartPar
\sphinxcode{\sphinxupquote{c}}: constant multiplier.
\end{quote}
\end{quote}

\subsubsection{QuadExpr::operator*()}
\label{\detokenize{cppapi/quadexpr:id7}}\begin{quote}

\sphinxAtStartPar
Multiply constant and return new expression.

\sphinxAtStartPar
\sphinxstylestrong{Synopsis}
\begin{quote}

\sphinxAtStartPar
\sphinxcode{\sphinxupquote{QuadExpr operator*(double c)}}
\end{quote}

\sphinxAtStartPar
\sphinxstylestrong{Arguments}
\begin{quote}

\sphinxAtStartPar
\sphinxcode{\sphinxupquote{c}}: constant multiplier.
\end{quote}

\sphinxAtStartPar
\sphinxstylestrong{Return}
\begin{quote}

\sphinxAtStartPar
result expression.
\end{quote}
\end{quote}

\subsubsection{QuadExpr::operator*()}
\label{\detokenize{cppapi/quadexpr:id8}}\begin{quote}

\sphinxAtStartPar
Multiply a variable and return new nonlinear expression.

\sphinxAtStartPar
\sphinxstylestrong{Synopsis}
\begin{quote}

\sphinxAtStartPar
\sphinxcode{\sphinxupquote{NlExpr operator*(const Var \&var)}}
\end{quote}

\sphinxAtStartPar
\sphinxstylestrong{Arguments}
\begin{quote}

\sphinxAtStartPar
\sphinxcode{\sphinxupquote{var}}: a variable as multiplier.
\end{quote}

\sphinxAtStartPar
\sphinxstylestrong{Return}
\begin{quote}

\sphinxAtStartPar
result nonlinear expression.
\end{quote}
\end{quote}

\subsubsection{QuadExpr::operator*()}
\label{\detokenize{cppapi/quadexpr:id9}}\begin{quote}

\sphinxAtStartPar
Multiply a linear expression and return new nonlinear expression.

\sphinxAtStartPar
\sphinxstylestrong{Synopsis}
\begin{quote}

\sphinxAtStartPar
\sphinxcode{\sphinxupquote{NlExpr operator*(const Expr \&expr)}}
\end{quote}

\sphinxAtStartPar
\sphinxstylestrong{Arguments}
\begin{quote}

\sphinxAtStartPar
\sphinxcode{\sphinxupquote{expr}}: a linear expression as multiplier.
\end{quote}

\sphinxAtStartPar
\sphinxstylestrong{Return}
\begin{quote}

\sphinxAtStartPar
result nonlinear expression.
\end{quote}
\end{quote}

\subsubsection{QuadExpr::operator*()}
\label{\detokenize{cppapi/quadexpr:id10}}\begin{quote}

\sphinxAtStartPar
Multiply a quadratic expression and return new nonlinear expression.

\sphinxAtStartPar
\sphinxstylestrong{Synopsis}
\begin{quote}

\sphinxAtStartPar
\sphinxcode{\sphinxupquote{NlExpr operator*(const QuadExpr \&expr)}}
\end{quote}

\sphinxAtStartPar
\sphinxstylestrong{Arguments}
\begin{quote}

\sphinxAtStartPar
\sphinxcode{\sphinxupquote{expr}}: a quadratic expression as multiplier.
\end{quote}

\sphinxAtStartPar
\sphinxstylestrong{Return}
\begin{quote}

\sphinxAtStartPar
result nonlinear expression.
\end{quote}
\end{quote}

\subsubsection{QuadExpr::operator/()}
\label{\detokenize{cppapi/quadexpr:id11}}\begin{quote}

\sphinxAtStartPar
Devided by a constant and return new quadratic expression.

\sphinxAtStartPar
\sphinxstylestrong{Synopsis}
\begin{quote}

\sphinxAtStartPar
\sphinxcode{\sphinxupquote{QuadExpr operator/(double c)}}
\end{quote}

\sphinxAtStartPar
\sphinxstylestrong{Arguments}
\begin{quote}

\sphinxAtStartPar
\sphinxcode{\sphinxupquote{c}}: constant divisor.
\end{quote}

\sphinxAtStartPar
\sphinxstylestrong{Return}
\begin{quote}

\sphinxAtStartPar
result expression.
\end{quote}
\end{quote}

\subsubsection{QuadExpr::operator/()}
\label{\detokenize{cppapi/quadexpr:id12}}\begin{quote}

\sphinxAtStartPar
Devided by a variable and return new nonlinear expression.

\sphinxAtStartPar
\sphinxstylestrong{Synopsis}
\begin{quote}

\sphinxAtStartPar
\sphinxcode{\sphinxupquote{NlExpr operator/(const Var \&var)}}
\end{quote}

\sphinxAtStartPar
\sphinxstylestrong{Arguments}
\begin{quote}

\sphinxAtStartPar
\sphinxcode{\sphinxupquote{var}}: a variable as divisor.
\end{quote}

\sphinxAtStartPar
\sphinxstylestrong{Return}
\begin{quote}

\sphinxAtStartPar
result nonlinear expression.
\end{quote}
\end{quote}

\subsubsection{QuadExpr::operator/()}
\label{\detokenize{cppapi/quadexpr:id13}}\begin{quote}

\sphinxAtStartPar
Devided by a linear expression and return new nonlinear expression.

\sphinxAtStartPar
\sphinxstylestrong{Synopsis}
\begin{quote}

\sphinxAtStartPar
\sphinxcode{\sphinxupquote{NlExpr operator/(const Expr \&other)}}
\end{quote}

\sphinxAtStartPar
\sphinxstylestrong{Arguments}
\begin{quote}

\sphinxAtStartPar
\sphinxcode{\sphinxupquote{other}}: a linear expression as divisor.
\end{quote}

\sphinxAtStartPar
\sphinxstylestrong{Return}
\begin{quote}

\sphinxAtStartPar
result nonlinear expression.
\end{quote}
\end{quote}

\subsubsection{QuadExpr::operator/()}
\label{\detokenize{cppapi/quadexpr:id14}}\begin{quote}

\sphinxAtStartPar
Devided by a quadratic expression and return new nonlinear expression.

\sphinxAtStartPar
\sphinxstylestrong{Synopsis}
\begin{quote}

\sphinxAtStartPar
\sphinxcode{\sphinxupquote{NlExpr operator/(const QuadExpr \&other)}}
\end{quote}

\sphinxAtStartPar
\sphinxstylestrong{Arguments}
\begin{quote}

\sphinxAtStartPar
\sphinxcode{\sphinxupquote{other}}: a quadratic expression as divisor.
\end{quote}

\sphinxAtStartPar
\sphinxstylestrong{Return}
\begin{quote}

\sphinxAtStartPar
result nonlinear expression.
\end{quote}
\end{quote}

\subsubsection{QuadExpr::operator+=()}
\label{\detokenize{cppapi/quadexpr:id15}}\begin{quote}

\sphinxAtStartPar
Add an expression to self.

\sphinxAtStartPar
\sphinxstylestrong{Synopsis}
\begin{quote}

\sphinxAtStartPar
\sphinxcode{\sphinxupquote{void operator+=(const QuadExpr \&expr)}}
\end{quote}

\sphinxAtStartPar
\sphinxstylestrong{Arguments}
\begin{quote}

\sphinxAtStartPar
\sphinxcode{\sphinxupquote{expr}}: expression to be added.
\end{quote}
\end{quote}

\subsubsection{QuadExpr::operator+()}
\label{\detokenize{cppapi/quadexpr:id16}}\begin{quote}

\sphinxAtStartPar
Add expression and return new expression.

\sphinxAtStartPar
\sphinxstylestrong{Synopsis}
\begin{quote}

\sphinxAtStartPar
\sphinxcode{\sphinxupquote{QuadExpr operator+(const QuadExpr \&other)}}
\end{quote}

\sphinxAtStartPar
\sphinxstylestrong{Arguments}
\begin{quote}

\sphinxAtStartPar
\sphinxcode{\sphinxupquote{other}}: other expression to add.
\end{quote}

\sphinxAtStartPar
\sphinxstylestrong{Return}
\begin{quote}

\sphinxAtStartPar
result expression.
\end{quote}
\end{quote}

\subsubsection{QuadExpr::operator\sphinxhyphen{}=()}
\label{\detokenize{cppapi/quadexpr:id17}}\begin{quote}

\sphinxAtStartPar
Substract an expression from self.

\sphinxAtStartPar
\sphinxstylestrong{Synopsis}
\begin{quote}

\sphinxAtStartPar
\sphinxcode{\sphinxupquote{void operator\sphinxhyphen{}=(const QuadExpr \&expr)}}
\end{quote}

\sphinxAtStartPar
\sphinxstylestrong{Arguments}
\begin{quote}

\sphinxAtStartPar
\sphinxcode{\sphinxupquote{expr}}: expression to be substracted.
\end{quote}
\end{quote}

\subsubsection{QuadExpr::operator\sphinxhyphen{}()}
\label{\detokenize{cppapi/quadexpr:id18}}\begin{quote}

\sphinxAtStartPar
Substract expression and return new expression.

\sphinxAtStartPar
\sphinxstylestrong{Synopsis}
\begin{quote}

\sphinxAtStartPar
\sphinxcode{\sphinxupquote{QuadExpr operator\sphinxhyphen{}(const QuadExpr \&other)}}
\end{quote}

\sphinxAtStartPar
\sphinxstylestrong{Arguments}
\begin{quote}

\sphinxAtStartPar
\sphinxcode{\sphinxupquote{other}}: other expression to substract.
\end{quote}

\sphinxAtStartPar
\sphinxstylestrong{Return}
\begin{quote}

\sphinxAtStartPar
result expression.
\end{quote}
\end{quote}

\subsubsection{QuadExpr::Remove()}
\label{\detokenize{cppapi/quadexpr:quadexpr-remove}}\begin{quote}

\sphinxAtStartPar
Remove i\sphinxhyphen{}th term from expression object.

\sphinxAtStartPar
\sphinxstylestrong{Synopsis}
\begin{quote}

\sphinxAtStartPar
\sphinxcode{\sphinxupquote{void Remove(int i)}}
\end{quote}

\sphinxAtStartPar
\sphinxstylestrong{Arguments}
\begin{quote}

\sphinxAtStartPar
\sphinxcode{\sphinxupquote{i}}: index of the term to be removed.
\end{quote}
\end{quote}

\subsubsection{QuadExpr::Remove()}
\label{\detokenize{cppapi/quadexpr:id19}}\begin{quote}

\sphinxAtStartPar
Remove the term associated with variable from expression.

\sphinxAtStartPar
\sphinxstylestrong{Synopsis}
\begin{quote}

\sphinxAtStartPar
\sphinxcode{\sphinxupquote{void Remove(const Var \&var)}}
\end{quote}

\sphinxAtStartPar
\sphinxstylestrong{Arguments}
\begin{quote}

\sphinxAtStartPar
\sphinxcode{\sphinxupquote{var}}: a variable whose term should be removed.
\end{quote}
\end{quote}

\subsubsection{QuadExpr::Reserve()}
\label{\detokenize{cppapi/quadexpr:quadexpr-reserve}}\begin{quote}

\sphinxAtStartPar
Reserve capacity to contain at least n items.

\sphinxAtStartPar
\sphinxstylestrong{Synopsis}
\begin{quote}

\sphinxAtStartPar
\sphinxcode{\sphinxupquote{void Reserve(size\_t n)}}
\end{quote}

\sphinxAtStartPar
\sphinxstylestrong{Arguments}
\begin{quote}

\sphinxAtStartPar
\sphinxcode{\sphinxupquote{n}}: minimum capacity for quadratic expression object.
\end{quote}
\end{quote}

\subsubsection{QuadExpr::SetCoeff()}
\label{\detokenize{cppapi/quadexpr:quadexpr-setcoeff}}\begin{quote}

\sphinxAtStartPar
Set coefficient of the i\sphinxhyphen{}th term in quadratic expression.

\sphinxAtStartPar
\sphinxstylestrong{Synopsis}
\begin{quote}

\sphinxAtStartPar
\sphinxcode{\sphinxupquote{void SetCoeff(int i, double val)}}
\end{quote}

\sphinxAtStartPar
\sphinxstylestrong{Arguments}
\begin{quote}

\sphinxAtStartPar
\sphinxcode{\sphinxupquote{i}}: index of the quadratic term.

\sphinxAtStartPar
\sphinxcode{\sphinxupquote{val}}: coefficient of the term.
\end{quote}
\end{quote}

\subsubsection{QuadExpr::SetConstant()}
\label{\detokenize{cppapi/quadexpr:quadexpr-setconstant}}\begin{quote}

\sphinxAtStartPar
Set constant for the expression.

\sphinxAtStartPar
\sphinxstylestrong{Synopsis}
\begin{quote}

\sphinxAtStartPar
\sphinxcode{\sphinxupquote{void SetConstant(double constant)}}
\end{quote}

\sphinxAtStartPar
\sphinxstylestrong{Arguments}
\begin{quote}

\sphinxAtStartPar
\sphinxcode{\sphinxupquote{constant}}: the value of the constant.
\end{quote}
\end{quote}

\subsubsection{QuadExpr::Size()}
\label{\detokenize{cppapi/quadexpr:quadexpr-size}}\begin{quote}

\sphinxAtStartPar
Get number of terms in expression.

\sphinxAtStartPar
\sphinxstylestrong{Synopsis}
\begin{quote}

\sphinxAtStartPar
\sphinxcode{\sphinxupquote{size\_t Size()}}
\end{quote}

\sphinxAtStartPar
\sphinxstylestrong{Return}
\begin{quote}

\sphinxAtStartPar
number of terms.
\end{quote}
\end{quote}

\subsection{QConstraint}
\label{\detokenize{cppapiref:qconstraint}}\label{\detokenize{cppapiref:chapcppapiref-qconstr}}
\sphinxAtStartPar
COPT quadratic constraint object. Quadratic constraints are always
associated with a particular model. User creates a quadratic constraint
object by adding a quadratic constraint to a model, rather than by
using constructor of QConstraint class.

\sphinxstepscope

\subsubsection{QConstraint::Get()}
\label{\detokenize{cppapi/qconstraint:qconstraint-get}}\label{\detokenize{cppapi/qconstraint::doc}}\begin{quote}

\sphinxAtStartPar
Get information value of the quadratic constraint. Support related quadratic informations.

\sphinxAtStartPar
\sphinxstylestrong{Synopsis}
\begin{quote}

\sphinxAtStartPar
\sphinxcode{\sphinxupquote{double Get(const char *szInfo)}}
\end{quote}

\sphinxAtStartPar
\sphinxstylestrong{Arguments}
\begin{quote}

\sphinxAtStartPar
\sphinxcode{\sphinxupquote{szInfo}}: name of the information being queried.
\end{quote}

\sphinxAtStartPar
\sphinxstylestrong{Return}
\begin{quote}

\sphinxAtStartPar
value of information.
\end{quote}
\end{quote}

\subsubsection{QConstraint::GetIdx()}
\label{\detokenize{cppapi/qconstraint:qconstraint-getidx}}\begin{quote}

\sphinxAtStartPar
Get index of the quadratic constraint.

\sphinxAtStartPar
\sphinxstylestrong{Synopsis}
\begin{quote}

\sphinxAtStartPar
\sphinxcode{\sphinxupquote{int GetIdx()}}
\end{quote}

\sphinxAtStartPar
\sphinxstylestrong{Return}
\begin{quote}

\sphinxAtStartPar
the index of the quadratic constraint.
\end{quote}
\end{quote}

\subsubsection{QConstraint::GetName()}
\label{\detokenize{cppapi/qconstraint:qconstraint-getname}}\begin{quote}

\sphinxAtStartPar
Get name of the quadratic constraint.

\sphinxAtStartPar
\sphinxstylestrong{Synopsis}
\begin{quote}

\sphinxAtStartPar
\sphinxcode{\sphinxupquote{const char *GetName()}}
\end{quote}

\sphinxAtStartPar
\sphinxstylestrong{Return}
\begin{quote}

\sphinxAtStartPar
the name of the quadratic constraint.
\end{quote}
\end{quote}

\subsubsection{QConstraint::GetRhs()}
\label{\detokenize{cppapi/qconstraint:qconstraint-getrhs}}\begin{quote}

\sphinxAtStartPar
Get rhs of quadratic constraint.

\sphinxAtStartPar
\sphinxstylestrong{Synopsis}
\begin{quote}

\sphinxAtStartPar
\sphinxcode{\sphinxupquote{double GetRhs()}}
\end{quote}

\sphinxAtStartPar
\sphinxstylestrong{Return}
\begin{quote}

\sphinxAtStartPar
rhs of quadratic constraint.
\end{quote}
\end{quote}

\subsubsection{QConstraint::GetSense()}
\label{\detokenize{cppapi/qconstraint:qconstraint-getsense}}\begin{quote}

\sphinxAtStartPar
Get sense of quadratic constraint.

\sphinxAtStartPar
\sphinxstylestrong{Synopsis}
\begin{quote}

\sphinxAtStartPar
\sphinxcode{\sphinxupquote{char GetSense()}}
\end{quote}

\sphinxAtStartPar
\sphinxstylestrong{Return}
\begin{quote}

\sphinxAtStartPar
sense of quadratic constraint.
\end{quote}
\end{quote}

\subsubsection{QConstraint::Remove()}
\label{\detokenize{cppapi/qconstraint:qconstraint-remove}}\begin{quote}

\sphinxAtStartPar
Remove this quadratic constraint from model.

\sphinxAtStartPar
\sphinxstylestrong{Synopsis}
\begin{quote}

\sphinxAtStartPar
\sphinxcode{\sphinxupquote{void Remove()}}
\end{quote}
\end{quote}

\subsubsection{QConstraint::Set()}
\label{\detokenize{cppapi/qconstraint:qconstraint-set}}\begin{quote}

\sphinxAtStartPar
Set information value of the quadratic constraint. Support related quadratic informations.

\sphinxAtStartPar
\sphinxstylestrong{Synopsis}
\begin{quote}

\sphinxAtStartPar
\sphinxcode{\sphinxupquote{void Set(const char *szInfo, double value)}}
\end{quote}

\sphinxAtStartPar
\sphinxstylestrong{Arguments}
\begin{quote}

\sphinxAtStartPar
\sphinxcode{\sphinxupquote{szInfo}}: name of the information.

\sphinxAtStartPar
\sphinxcode{\sphinxupquote{value}}: new information value.
\end{quote}
\end{quote}

\subsubsection{QConstraint::SetName()}
\label{\detokenize{cppapi/qconstraint:qconstraint-setname}}\begin{quote}

\sphinxAtStartPar
Set name of a quadratic constraint.

\sphinxAtStartPar
\sphinxstylestrong{Synopsis}
\begin{quote}

\sphinxAtStartPar
\sphinxcode{\sphinxupquote{void SetName(const char *szName)}}
\end{quote}

\sphinxAtStartPar
\sphinxstylestrong{Arguments}
\begin{quote}

\sphinxAtStartPar
\sphinxcode{\sphinxupquote{szName}}: the name to set.
\end{quote}
\end{quote}

\subsubsection{QConstraint::SetRhs()}
\label{\detokenize{cppapi/qconstraint:qconstraint-setrhs}}\begin{quote}

\sphinxAtStartPar
Set rhs of quadratic constraint.

\sphinxAtStartPar
\sphinxstylestrong{Synopsis}
\begin{quote}

\sphinxAtStartPar
\sphinxcode{\sphinxupquote{void SetRhs(double rhs)}}
\end{quote}

\sphinxAtStartPar
\sphinxstylestrong{Arguments}
\begin{quote}

\sphinxAtStartPar
\sphinxcode{\sphinxupquote{rhs}}: rhs of quadratic constraint.
\end{quote}
\end{quote}

\subsubsection{QConstraint::SetSense()}
\label{\detokenize{cppapi/qconstraint:qconstraint-setsense}}\begin{quote}

\sphinxAtStartPar
Set sense of quadratic constraint.

\sphinxAtStartPar
\sphinxstylestrong{Synopsis}
\begin{quote}

\sphinxAtStartPar
\sphinxcode{\sphinxupquote{void SetSense(char sense)}}
\end{quote}

\sphinxAtStartPar
\sphinxstylestrong{Arguments}
\begin{quote}

\sphinxAtStartPar
\sphinxcode{\sphinxupquote{sense}}: sense of quadratic constraint.
\end{quote}
\end{quote}

\subsection{QConstrArray}
\label{\detokenize{cppapiref:qconstrarray}}\label{\detokenize{cppapiref:chapcppapiref-qconstrarray}}
\sphinxAtStartPar
COPT quadratic constraint array object. To store and access a
set of C++ {\hyperref[\detokenize{cppapiref:chapcppapiref-qconstr}]{\sphinxcrossref{\DUrole{std,std-ref}{QConstraint}}}} objects, Cardinal Optimizer
provides C++ QConstrArray class, which defines the following methods.

\sphinxstepscope

\subsubsection{QConstrArray::GetQConstr()}
\label{\detokenize{cppapi/qconstrarray:qconstrarray-getqconstr}}\label{\detokenize{cppapi/qconstrarray::doc}}\begin{quote}

\sphinxAtStartPar
Get i\sphinxhyphen{}th quadratic constraint object.

\sphinxAtStartPar
\sphinxstylestrong{Synopsis}
\begin{quote}

\sphinxAtStartPar
\sphinxcode{\sphinxupquote{QConstraint \&GetQConstr(int idx)}}
\end{quote}

\sphinxAtStartPar
\sphinxstylestrong{Arguments}
\begin{quote}

\sphinxAtStartPar
\sphinxcode{\sphinxupquote{idx}}: index of the quadratic constraint.
\end{quote}

\sphinxAtStartPar
\sphinxstylestrong{Return}
\begin{quote}

\sphinxAtStartPar
quadratic constraint object with index i.
\end{quote}
\end{quote}

\subsubsection{QConstrArray::PushBack()}
\label{\detokenize{cppapi/qconstrarray:qconstrarray-pushback}}\begin{quote}

\sphinxAtStartPar
Add a quadratic constraint object to constraint array.

\sphinxAtStartPar
\sphinxstylestrong{Synopsis}
\begin{quote}

\sphinxAtStartPar
\sphinxcode{\sphinxupquote{void PushBack(const QConstraint \&constr)}}
\end{quote}

\sphinxAtStartPar
\sphinxstylestrong{Arguments}
\begin{quote}

\sphinxAtStartPar
\sphinxcode{\sphinxupquote{constr}}: a quadratic constraint object.
\end{quote}
\end{quote}

\subsubsection{QConstrArray::Reserve()}
\label{\detokenize{cppapi/qconstrarray:qconstrarray-reserve}}\begin{quote}

\sphinxAtStartPar
Reserve capacity to contain at least n items.

\sphinxAtStartPar
\sphinxstylestrong{Synopsis}
\begin{quote}

\sphinxAtStartPar
\sphinxcode{\sphinxupquote{void Reserve(int n)}}
\end{quote}

\sphinxAtStartPar
\sphinxstylestrong{Arguments}
\begin{quote}

\sphinxAtStartPar
\sphinxcode{\sphinxupquote{n}}: minimum capacity for quadratic constraint objects.
\end{quote}
\end{quote}

\subsubsection{QConstrArray::Size()}
\label{\detokenize{cppapi/qconstrarray:qconstrarray-size}}\begin{quote}

\sphinxAtStartPar
Get the number of quadratic constraint objects.

\sphinxAtStartPar
\sphinxstylestrong{Synopsis}
\begin{quote}

\sphinxAtStartPar
\sphinxcode{\sphinxupquote{int Size()}}
\end{quote}

\sphinxAtStartPar
\sphinxstylestrong{Return}
\begin{quote}

\sphinxAtStartPar
number of quadratic constraint objects.
\end{quote}
\end{quote}

\subsection{QConstrBuilder}
\label{\detokenize{cppapiref:qconstrbuilder}}\label{\detokenize{cppapiref:chapcppapiref-qconstrbuilder}}
\sphinxAtStartPar
COPT quadratic constraint builder object. To help building a quadratic
constraint, given a quadratic expression, constraint sense and right\sphinxhyphen{}hand
side value, Cardinal Optimizer provides C++ QConstrBuilder
class, which defines the following methods.

\sphinxstepscope

\subsubsection{QConstrBuilder::GetQuadExpr()}
\label{\detokenize{cppapi/qconstrbuilder:qconstrbuilder-getquadexpr}}\label{\detokenize{cppapi/qconstrbuilder::doc}}\begin{quote}

\sphinxAtStartPar
Get expression associated with quadratic constraint.

\sphinxAtStartPar
\sphinxstylestrong{Synopsis}
\begin{quote}

\sphinxAtStartPar
\sphinxcode{\sphinxupquote{const QuadExpr \&GetQuadExpr()}}
\end{quote}

\sphinxAtStartPar
\sphinxstylestrong{Return}
\begin{quote}

\sphinxAtStartPar
quadratic expression object.
\end{quote}
\end{quote}

\subsubsection{QConstrBuilder::GetSense()}
\label{\detokenize{cppapi/qconstrbuilder:qconstrbuilder-getsense}}\begin{quote}

\sphinxAtStartPar
Get sense associated with quadratic constraint.

\sphinxAtStartPar
\sphinxstylestrong{Synopsis}
\begin{quote}

\sphinxAtStartPar
\sphinxcode{\sphinxupquote{char GetSense()}}
\end{quote}

\sphinxAtStartPar
\sphinxstylestrong{Return}
\begin{quote}

\sphinxAtStartPar
quadratic constraint sense.
\end{quote}
\end{quote}

\subsubsection{QConstrBuilder::Set()}
\label{\detokenize{cppapi/qconstrbuilder:qconstrbuilder-set}}\begin{quote}

\sphinxAtStartPar
Set detail of a quadratic constraint to its builder object.

\sphinxAtStartPar
\sphinxstylestrong{Synopsis}
\begin{quote}

\sphinxAtStartPar
\sphinxcode{\sphinxupquote{void Set(}}
\begin{quote}

\sphinxAtStartPar
\sphinxcode{\sphinxupquote{const QuadExpr \&expr,}}

\sphinxAtStartPar
\sphinxcode{\sphinxupquote{char sense,}}

\sphinxAtStartPar
\sphinxcode{\sphinxupquote{double rhs)}}
\end{quote}
\end{quote}

\sphinxAtStartPar
\sphinxstylestrong{Arguments}
\begin{quote}

\sphinxAtStartPar
\sphinxcode{\sphinxupquote{expr}}: expression object at one side of the quadratic constraint.

\sphinxAtStartPar
\sphinxcode{\sphinxupquote{sense}}: quadratic constraint sense.

\sphinxAtStartPar
\sphinxcode{\sphinxupquote{rhs}}: constant of right side of quadratic constraint.
\end{quote}
\end{quote}

\subsection{QConstrBuilderArray}
\label{\detokenize{cppapiref:qconstrbuilderarray}}\label{\detokenize{cppapiref:chapcppapiref-qconstrbuilderarray}}
\sphinxAtStartPar
COPT quadratic constraint builder array object. To store and access a set
of C++ {\hyperref[\detokenize{cppapiref:chapcppapiref-qconstrbuilder}]{\sphinxcrossref{\DUrole{std,std-ref}{QConstrBuilder}}}} objects, Cardinal Optimizer provides
C++ QConstrBuilderArray class, which defines the following methods.

\sphinxstepscope

\subsubsection{QConstrBuilderArray::GetBuilder()}
\label{\detokenize{cppapi/qconstrbuilderarray:qconstrbuilderarray-getbuilder}}\label{\detokenize{cppapi/qconstrbuilderarray::doc}}\begin{quote}

\sphinxAtStartPar
Get i\sphinxhyphen{}th quadratic constraint builder object.

\sphinxAtStartPar
\sphinxstylestrong{Synopsis}
\begin{quote}

\sphinxAtStartPar
\sphinxcode{\sphinxupquote{QConstrBuilder \&GetBuilder(int idx)}}
\end{quote}

\sphinxAtStartPar
\sphinxstylestrong{Arguments}
\begin{quote}

\sphinxAtStartPar
\sphinxcode{\sphinxupquote{idx}}: index of the quadratic constraint builder.
\end{quote}

\sphinxAtStartPar
\sphinxstylestrong{Return}
\begin{quote}

\sphinxAtStartPar
quadratic constraint builder object with index i.
\end{quote}
\end{quote}

\subsubsection{QConstrBuilderArray::PushBack()}
\label{\detokenize{cppapi/qconstrbuilderarray:qconstrbuilderarray-pushback}}\begin{quote}

\sphinxAtStartPar
Add a quadratic constraint builder object to quadratic constraint builder array.

\sphinxAtStartPar
\sphinxstylestrong{Synopsis}
\begin{quote}

\sphinxAtStartPar
\sphinxcode{\sphinxupquote{void PushBack(const QConstrBuilder \&builder)}}
\end{quote}

\sphinxAtStartPar
\sphinxstylestrong{Arguments}
\begin{quote}

\sphinxAtStartPar
\sphinxcode{\sphinxupquote{builder}}: a quadratic constraint builder object.
\end{quote}
\end{quote}

\subsubsection{QConstrBuilderArray::Reserve()}
\label{\detokenize{cppapi/qconstrbuilderarray:qconstrbuilderarray-reserve}}\begin{quote}

\sphinxAtStartPar
Reserve capacity to contain at least n items.

\sphinxAtStartPar
\sphinxstylestrong{Synopsis}
\begin{quote}

\sphinxAtStartPar
\sphinxcode{\sphinxupquote{void Reserve(int n)}}
\end{quote}

\sphinxAtStartPar
\sphinxstylestrong{Arguments}
\begin{quote}

\sphinxAtStartPar
\sphinxcode{\sphinxupquote{n}}: minimum capacity for quadratic constraint builder object.
\end{quote}
\end{quote}

\subsubsection{QConstrBuilderArray::Size()}
\label{\detokenize{cppapi/qconstrbuilderarray:qconstrbuilderarray-size}}\begin{quote}

\sphinxAtStartPar
Get the number of quadratic constraint builder objects.

\sphinxAtStartPar
\sphinxstylestrong{Synopsis}
\begin{quote}

\sphinxAtStartPar
\sphinxcode{\sphinxupquote{int Size()}}
\end{quote}

\sphinxAtStartPar
\sphinxstylestrong{Return}
\begin{quote}

\sphinxAtStartPar
number of quadratic constraint builder objects.
\end{quote}
\end{quote}

\subsection{PsdVar}
\label{\detokenize{cppapiref:psdvar}}\label{\detokenize{cppapiref:chapcppapiref-psdvar}}
\sphinxAtStartPar
COPT PSD variable object. PSD variables are always associated with a
particular model.  User creates a PSD variable object by adding a PSD
variable to model, rather than by constructor of PsdVar class.

\sphinxstepscope

\subsubsection{PsdVar::Get()}
\label{\detokenize{cppapi/psdvar:psdvar-get}}\label{\detokenize{cppapi/psdvar::doc}}\begin{quote}

\sphinxAtStartPar
Get information values of PSD variable.

\sphinxAtStartPar
\sphinxstylestrong{Synopsis}
\begin{quote}

\sphinxAtStartPar
\sphinxcode{\sphinxupquote{double Get(const char *szInfo, int sz)}}
\end{quote}

\sphinxAtStartPar
\sphinxstylestrong{Arguments}
\begin{quote}

\sphinxAtStartPar
\sphinxcode{\sphinxupquote{szInfo}}: name of informatioin.

\sphinxAtStartPar
\sphinxcode{\sphinxupquote{sz}}: length of the output array.
\end{quote}

\sphinxAtStartPar
\sphinxstylestrong{Return}
\begin{quote}

\sphinxAtStartPar
output array of information values.
\end{quote}
\end{quote}

\subsubsection{PsdVar::Get()}
\label{\detokenize{cppapi/psdvar:id1}}\begin{quote}

\sphinxAtStartPar
Get information values of PSD variable.

\sphinxAtStartPar
\sphinxstylestrong{Synopsis}
\begin{quote}

\sphinxAtStartPar
\sphinxcode{\sphinxupquote{NdArray\textless{}double, 2\textgreater{} Get(const char *szInfo)}}
\end{quote}

\sphinxAtStartPar
\sphinxstylestrong{Arguments}
\begin{quote}

\sphinxAtStartPar
\sphinxcode{\sphinxupquote{szInfo}}: name of informatioin.
\end{quote}

\sphinxAtStartPar
\sphinxstylestrong{Return}
\begin{quote}

\sphinxAtStartPar
2\sphinxhyphen{}dimensional NdArray with related information values.
\end{quote}
\end{quote}

\subsubsection{PsdVar::GetDim()}
\label{\detokenize{cppapi/psdvar:psdvar-getdim}}\begin{quote}

\sphinxAtStartPar
Get dimension of PSD variable.

\sphinxAtStartPar
\sphinxstylestrong{Synopsis}
\begin{quote}

\sphinxAtStartPar
\sphinxcode{\sphinxupquote{int GetDim()}}
\end{quote}

\sphinxAtStartPar
\sphinxstylestrong{Return}
\begin{quote}

\sphinxAtStartPar
dimension of PSD variable.
\end{quote}
\end{quote}

\subsubsection{PsdVar::GetIdx()}
\label{\detokenize{cppapi/psdvar:psdvar-getidx}}\begin{quote}

\sphinxAtStartPar
Get index of PSD variable.

\sphinxAtStartPar
\sphinxstylestrong{Synopsis}
\begin{quote}

\sphinxAtStartPar
\sphinxcode{\sphinxupquote{int GetIdx()}}
\end{quote}

\sphinxAtStartPar
\sphinxstylestrong{Return}
\begin{quote}

\sphinxAtStartPar
index of PSD variable.
\end{quote}
\end{quote}

\subsubsection{PsdVar::GetLen()}
\label{\detokenize{cppapi/psdvar:psdvar-getlen}}\begin{quote}

\sphinxAtStartPar
Get length of PSD variable.

\sphinxAtStartPar
\sphinxstylestrong{Synopsis}
\begin{quote}

\sphinxAtStartPar
\sphinxcode{\sphinxupquote{int GetLen()}}
\end{quote}

\sphinxAtStartPar
\sphinxstylestrong{Return}
\begin{quote}

\sphinxAtStartPar
length of PSD variable.
\end{quote}
\end{quote}

\subsubsection{PsdVar::GetName()}
\label{\detokenize{cppapi/psdvar:psdvar-getname}}\begin{quote}

\sphinxAtStartPar
Get name of PSD variable.

\sphinxAtStartPar
\sphinxstylestrong{Synopsis}
\begin{quote}

\sphinxAtStartPar
\sphinxcode{\sphinxupquote{const char *GetName()}}
\end{quote}

\sphinxAtStartPar
\sphinxstylestrong{Return}
\begin{quote}

\sphinxAtStartPar
name of PSD variable.
\end{quote}
\end{quote}

\subsubsection{PsdVar::GetSize()}
\label{\detokenize{cppapi/psdvar:psdvar-getsize}}\begin{quote}

\sphinxAtStartPar
Get size or length of PSD variable.

\sphinxAtStartPar
\sphinxstylestrong{Synopsis}
\begin{quote}

\sphinxAtStartPar
\sphinxcode{\sphinxupquote{int GetSize()}}
\end{quote}

\sphinxAtStartPar
\sphinxstylestrong{Return}
\begin{quote}

\sphinxAtStartPar
Size of PSD variable.
\end{quote}
\end{quote}

\subsubsection{PsdVar::Remove()}
\label{\detokenize{cppapi/psdvar:psdvar-remove}}\begin{quote}

\sphinxAtStartPar
Remove PSD variable from model.

\sphinxAtStartPar
\sphinxstylestrong{Synopsis}
\begin{quote}

\sphinxAtStartPar
\sphinxcode{\sphinxupquote{void Remove()}}
\end{quote}
\end{quote}

\subsubsection{PsdVar::SetName()}
\label{\detokenize{cppapi/psdvar:psdvar-setname}}\begin{quote}

\sphinxAtStartPar
Set name of PSD variable.

\sphinxAtStartPar
\sphinxstylestrong{Synopsis}
\begin{quote}

\sphinxAtStartPar
\sphinxcode{\sphinxupquote{void SetName(const char *szName)}}
\end{quote}

\sphinxAtStartPar
\sphinxstylestrong{Arguments}
\begin{quote}

\sphinxAtStartPar
\sphinxcode{\sphinxupquote{szName}}: name of PSD variable.
\end{quote}
\end{quote}

\subsection{PsdVarArray}
\label{\detokenize{cppapiref:psdvararray}}\label{\detokenize{cppapiref:chapcppapiref-psdvararray}}
\sphinxAtStartPar
COPT PSD variable array object. To store and access a set of
{\hyperref[\detokenize{cppapiref:chapcppapiref-psdvar}]{\sphinxcrossref{\DUrole{std,std-ref}{PsdVar}}}} objects, Cardinal Optimizer provides PsdVarArray
class, which defines the following methods.

\sphinxstepscope

\subsubsection{PsdVarArray::GetPsdVar()}
\label{\detokenize{cppapi/psdvararray:psdvararray-getpsdvar}}\label{\detokenize{cppapi/psdvararray::doc}}\begin{quote}

\sphinxAtStartPar
Get idx\sphinxhyphen{}th PSD variable object.

\sphinxAtStartPar
\sphinxstylestrong{Synopsis}
\begin{quote}

\sphinxAtStartPar
\sphinxcode{\sphinxupquote{PsdVar \&GetPsdVar(int idx)}}
\end{quote}

\sphinxAtStartPar
\sphinxstylestrong{Arguments}
\begin{quote}

\sphinxAtStartPar
\sphinxcode{\sphinxupquote{idx}}: index of the PSD variable.
\end{quote}

\sphinxAtStartPar
\sphinxstylestrong{Return}
\begin{quote}

\sphinxAtStartPar
PSD variable object with index idx.
\end{quote}
\end{quote}

\subsubsection{PsdVarArray::PushBack()}
\label{\detokenize{cppapi/psdvararray:psdvararray-pushback}}\begin{quote}

\sphinxAtStartPar
Add a PSD variable object to PSD variable array.

\sphinxAtStartPar
\sphinxstylestrong{Synopsis}
\begin{quote}

\sphinxAtStartPar
\sphinxcode{\sphinxupquote{void PushBack(const PsdVar \&var)}}
\end{quote}

\sphinxAtStartPar
\sphinxstylestrong{Arguments}
\begin{quote}

\sphinxAtStartPar
\sphinxcode{\sphinxupquote{var}}: a PSD variable object.
\end{quote}
\end{quote}

\subsubsection{PsdVarArray::Reserve()}
\label{\detokenize{cppapi/psdvararray:psdvararray-reserve}}\begin{quote}

\sphinxAtStartPar
Reserve capacity to contain at least n items.

\sphinxAtStartPar
\sphinxstylestrong{Synopsis}
\begin{quote}

\sphinxAtStartPar
\sphinxcode{\sphinxupquote{void Reserve(int n)}}
\end{quote}

\sphinxAtStartPar
\sphinxstylestrong{Arguments}
\begin{quote}

\sphinxAtStartPar
\sphinxcode{\sphinxupquote{n}}: minimum capacity for PSD variable object.
\end{quote}
\end{quote}

\subsubsection{PsdVarArray::Size()}
\label{\detokenize{cppapi/psdvararray:psdvararray-size}}\begin{quote}

\sphinxAtStartPar
Get the number of PSD variable objects.

\sphinxAtStartPar
\sphinxstylestrong{Synopsis}
\begin{quote}

\sphinxAtStartPar
\sphinxcode{\sphinxupquote{int Size()}}
\end{quote}

\sphinxAtStartPar
\sphinxstylestrong{Return}
\begin{quote}

\sphinxAtStartPar
number of PSD variable objects.
\end{quote}
\end{quote}

\subsection{PsdExpr}
\label{\detokenize{cppapiref:psdexpr}}\label{\detokenize{cppapiref:chapcppapiref-psdexpr}}
\sphinxAtStartPar
COPT PSD expression object. A PSD expression consists of a linear expression,
a list of PSD variables and associated coefficient matrices of PSD terms.
PSD expressions are used to build PSD constraints and objectives.

\sphinxstepscope

\subsubsection{PsdExpr::PsdExpr()}
\label{\detokenize{cppapi/psdexpr:psdexpr-psdexpr}}\label{\detokenize{cppapi/psdexpr::doc}}\begin{quote}

\sphinxAtStartPar
Constructor of a PSD expression with a constant.

\sphinxAtStartPar
\sphinxstylestrong{Synopsis}
\begin{quote}

\sphinxAtStartPar
\sphinxcode{\sphinxupquote{PsdExpr(double constant)}}
\end{quote}

\sphinxAtStartPar
\sphinxstylestrong{Arguments}
\begin{quote}

\sphinxAtStartPar
\sphinxcode{\sphinxupquote{constant}}: constant value in PSD expression object.
\end{quote}
\end{quote}

\subsubsection{PsdExpr::PsdExpr()}
\label{\detokenize{cppapi/psdexpr:id1}}\begin{quote}

\sphinxAtStartPar
Constructor of a PSD expression with one term.

\sphinxAtStartPar
\sphinxstylestrong{Synopsis}
\begin{quote}

\sphinxAtStartPar
\sphinxcode{\sphinxupquote{PsdExpr(const Var \&var, double coeff)}}
\end{quote}

\sphinxAtStartPar
\sphinxstylestrong{Arguments}
\begin{quote}

\sphinxAtStartPar
\sphinxcode{\sphinxupquote{var}}: variable for the added term.

\sphinxAtStartPar
\sphinxcode{\sphinxupquote{coeff}}: coefficent for the added term.
\end{quote}
\end{quote}

\subsubsection{PsdExpr::PsdExpr()}
\label{\detokenize{cppapi/psdexpr:id2}}\begin{quote}

\sphinxAtStartPar
Constructor of a PSD expression with a linear expression.

\sphinxAtStartPar
\sphinxstylestrong{Synopsis}
\begin{quote}

\sphinxAtStartPar
\sphinxcode{\sphinxupquote{PsdExpr(const Expr \&expr)}}
\end{quote}

\sphinxAtStartPar
\sphinxstylestrong{Arguments}
\begin{quote}

\sphinxAtStartPar
\sphinxcode{\sphinxupquote{expr}}: input linear expression.
\end{quote}
\end{quote}

\subsubsection{PsdExpr::PsdExpr()}
\label{\detokenize{cppapi/psdexpr:id3}}\begin{quote}

\sphinxAtStartPar
Constructor of a PSD expression with a MExpression object.

\sphinxAtStartPar
\sphinxstylestrong{Synopsis}
\begin{quote}

\sphinxAtStartPar
\sphinxcode{\sphinxupquote{PsdExpr(const MExpression \&expr)}}
\end{quote}

\sphinxAtStartPar
\sphinxstylestrong{Arguments}
\begin{quote}

\sphinxAtStartPar
\sphinxcode{\sphinxupquote{expr}}: input MExpression object.
\end{quote}
\end{quote}

\subsubsection{PsdExpr::PsdExpr()}
\label{\detokenize{cppapi/psdexpr:id4}}\begin{quote}

\sphinxAtStartPar
Constructor of a PSD expression with one term.

\sphinxAtStartPar
\sphinxstylestrong{Synopsis}
\begin{quote}

\sphinxAtStartPar
\sphinxcode{\sphinxupquote{PsdExpr(const PsdVar \&var, const SymMatrix \&mat)}}
\end{quote}

\sphinxAtStartPar
\sphinxstylestrong{Arguments}
\begin{quote}

\sphinxAtStartPar
\sphinxcode{\sphinxupquote{var}}: PSD variable for the added term.

\sphinxAtStartPar
\sphinxcode{\sphinxupquote{mat}}: coefficient matrix for the added term.
\end{quote}
\end{quote}

\subsubsection{PsdExpr::PsdExpr()}
\label{\detokenize{cppapi/psdexpr:id5}}\begin{quote}

\sphinxAtStartPar
Constructor of a PSD expression with one term.

\sphinxAtStartPar
\sphinxstylestrong{Synopsis}
\begin{quote}

\sphinxAtStartPar
\sphinxcode{\sphinxupquote{PsdExpr(const PsdVar \&var, const SymMatExpr \&expr)}}
\end{quote}

\sphinxAtStartPar
\sphinxstylestrong{Arguments}
\begin{quote}

\sphinxAtStartPar
\sphinxcode{\sphinxupquote{var}}: PSD variable for the added term.

\sphinxAtStartPar
\sphinxcode{\sphinxupquote{expr}}: coefficient expression of symmetric matrices for the added term.
\end{quote}
\end{quote}

\subsubsection{PsdExpr::AddConstant()}
\label{\detokenize{cppapi/psdexpr:psdexpr-addconstant}}\begin{quote}

\sphinxAtStartPar
Add constant to the PSD expression.

\sphinxAtStartPar
\sphinxstylestrong{Synopsis}
\begin{quote}

\sphinxAtStartPar
\sphinxcode{\sphinxupquote{void AddConstant(double constant)}}
\end{quote}

\sphinxAtStartPar
\sphinxstylestrong{Arguments}
\begin{quote}

\sphinxAtStartPar
\sphinxcode{\sphinxupquote{constant}}: value to be added.
\end{quote}
\end{quote}

\subsubsection{PsdExpr::AddLinExpr()}
\label{\detokenize{cppapi/psdexpr:psdexpr-addlinexpr}}\begin{quote}

\sphinxAtStartPar
Add a linear expression to self.

\sphinxAtStartPar
\sphinxstylestrong{Synopsis}
\begin{quote}

\sphinxAtStartPar
\sphinxcode{\sphinxupquote{void AddLinExpr(const Expr \&expr, double mult)}}
\end{quote}

\sphinxAtStartPar
\sphinxstylestrong{Arguments}
\begin{quote}

\sphinxAtStartPar
\sphinxcode{\sphinxupquote{expr}}: linear expression to be added.

\sphinxAtStartPar
\sphinxcode{\sphinxupquote{mult}}: optional, constant multiplier, default value is 1.0.
\end{quote}
\end{quote}

\subsubsection{PsdExpr::AddMExpr()}
\label{\detokenize{cppapi/psdexpr:psdexpr-addmexpr}}\begin{quote}

\sphinxAtStartPar
Add MExpression to PsdExpr object.

\sphinxAtStartPar
\sphinxstylestrong{Synopsis}
\begin{quote}

\sphinxAtStartPar
\sphinxcode{\sphinxupquote{void AddMExpr(const MExpression \&expr, double mult)}}
\end{quote}

\sphinxAtStartPar
\sphinxstylestrong{Arguments}
\begin{quote}

\sphinxAtStartPar
\sphinxcode{\sphinxupquote{expr}}: MExpression object.

\sphinxAtStartPar
\sphinxcode{\sphinxupquote{mult}}: the multiplier of MExpression, default value is 1.0.
\end{quote}
\end{quote}

\subsubsection{PsdExpr::AddPsdExpr()}
\label{\detokenize{cppapi/psdexpr:psdexpr-addpsdexpr}}\begin{quote}

\sphinxAtStartPar
Add a PSD expression to self.

\sphinxAtStartPar
\sphinxstylestrong{Synopsis}
\begin{quote}

\sphinxAtStartPar
\sphinxcode{\sphinxupquote{void AddPsdExpr(const PsdExpr \&expr, double mult)}}
\end{quote}

\sphinxAtStartPar
\sphinxstylestrong{Arguments}
\begin{quote}

\sphinxAtStartPar
\sphinxcode{\sphinxupquote{expr}}: PSD expression to be added.

\sphinxAtStartPar
\sphinxcode{\sphinxupquote{mult}}: optional, constant multiplier, default value is 1.0.
\end{quote}
\end{quote}

\subsubsection{PsdExpr::AddTerm()}
\label{\detokenize{cppapi/psdexpr:psdexpr-addterm}}\begin{quote}

\sphinxAtStartPar
Add a linear term to PSD expression object.

\sphinxAtStartPar
\sphinxstylestrong{Synopsis}
\begin{quote}

\sphinxAtStartPar
\sphinxcode{\sphinxupquote{void AddTerm(const Var \&var, double coeff)}}
\end{quote}

\sphinxAtStartPar
\sphinxstylestrong{Arguments}
\begin{quote}

\sphinxAtStartPar
\sphinxcode{\sphinxupquote{var}}: variable of new linear term.

\sphinxAtStartPar
\sphinxcode{\sphinxupquote{coeff}}: coefficient of new linear term.
\end{quote}
\end{quote}

\subsubsection{PsdExpr::AddTerm()}
\label{\detokenize{cppapi/psdexpr:id6}}\begin{quote}

\sphinxAtStartPar
Add a PSD term to PSD expression object.

\sphinxAtStartPar
\sphinxstylestrong{Synopsis}
\begin{quote}

\sphinxAtStartPar
\sphinxcode{\sphinxupquote{void AddTerm(const PsdVar \&var, const SymMatrix \&mat)}}
\end{quote}

\sphinxAtStartPar
\sphinxstylestrong{Arguments}
\begin{quote}

\sphinxAtStartPar
\sphinxcode{\sphinxupquote{var}}: PSD variable of new PSD term.

\sphinxAtStartPar
\sphinxcode{\sphinxupquote{mat}}: coefficient matrix of new PSD term.
\end{quote}
\end{quote}

\subsubsection{PsdExpr::AddTerm()}
\label{\detokenize{cppapi/psdexpr:id7}}\begin{quote}

\sphinxAtStartPar
Add a PSD term to PSD expression object.

\sphinxAtStartPar
\sphinxstylestrong{Synopsis}
\begin{quote}

\sphinxAtStartPar
\sphinxcode{\sphinxupquote{void AddTerm(const PsdVar \&var, const SymMatExpr \&expr)}}
\end{quote}

\sphinxAtStartPar
\sphinxstylestrong{Arguments}
\begin{quote}

\sphinxAtStartPar
\sphinxcode{\sphinxupquote{var}}: PSD variable of new PSD term.

\sphinxAtStartPar
\sphinxcode{\sphinxupquote{expr}}: coefficient expression of symmetric matrices of new PSD term.
\end{quote}
\end{quote}

\subsubsection{PsdExpr::AddTerms()}
\label{\detokenize{cppapi/psdexpr:psdexpr-addterms}}\begin{quote}

\sphinxAtStartPar
Add linear terms to PSD expression object.

\sphinxAtStartPar
\sphinxstylestrong{Synopsis}
\begin{quote}

\sphinxAtStartPar
\sphinxcode{\sphinxupquote{int AddTerms(}}
\begin{quote}

\sphinxAtStartPar
\sphinxcode{\sphinxupquote{const VarArray \&vars,}}

\sphinxAtStartPar
\sphinxcode{\sphinxupquote{double *pCoeff,}}

\sphinxAtStartPar
\sphinxcode{\sphinxupquote{int len)}}
\end{quote}
\end{quote}

\sphinxAtStartPar
\sphinxstylestrong{Arguments}
\begin{quote}

\sphinxAtStartPar
\sphinxcode{\sphinxupquote{vars}}: variables for added linear terms.

\sphinxAtStartPar
\sphinxcode{\sphinxupquote{pCoeff}}: coefficient array for added linear terms.

\sphinxAtStartPar
\sphinxcode{\sphinxupquote{len}}: length of coefficient array.
\end{quote}

\sphinxAtStartPar
\sphinxstylestrong{Return}
\begin{quote}

\sphinxAtStartPar
number of added linear terms.
\end{quote}
\end{quote}

\subsubsection{PsdExpr::AddTerms()}
\label{\detokenize{cppapi/psdexpr:id8}}\begin{quote}

\sphinxAtStartPar
Add PSD terms to PSD expression object.

\sphinxAtStartPar
\sphinxstylestrong{Synopsis}
\begin{quote}

\sphinxAtStartPar
\sphinxcode{\sphinxupquote{int AddTerms(const PsdVarArray \&vars, const SymMatrixArray \&mats)}}
\end{quote}

\sphinxAtStartPar
\sphinxstylestrong{Arguments}
\begin{quote}

\sphinxAtStartPar
\sphinxcode{\sphinxupquote{vars}}: PSD variables for added PSD terms.

\sphinxAtStartPar
\sphinxcode{\sphinxupquote{mats}}: coefficient matrixes for added PSD terms.
\end{quote}

\sphinxAtStartPar
\sphinxstylestrong{Return}
\begin{quote}

\sphinxAtStartPar
number of added PSD terms.
\end{quote}
\end{quote}

\subsubsection{PsdExpr::Clone()}
\label{\detokenize{cppapi/psdexpr:psdexpr-clone}}\begin{quote}

\sphinxAtStartPar
Deep copy PSD expression object.

\sphinxAtStartPar
\sphinxstylestrong{Synopsis}
\begin{quote}

\sphinxAtStartPar
\sphinxcode{\sphinxupquote{PsdExpr Clone()}}
\end{quote}

\sphinxAtStartPar
\sphinxstylestrong{Return}
\begin{quote}

\sphinxAtStartPar
cloned PSD expression object.
\end{quote}
\end{quote}

\subsubsection{PsdExpr::Evaluate()}
\label{\detokenize{cppapi/psdexpr:psdexpr-evaluate}}\begin{quote}

\sphinxAtStartPar
evaluate PSD expression after solving

\sphinxAtStartPar
\sphinxstylestrong{Synopsis}
\begin{quote}

\sphinxAtStartPar
\sphinxcode{\sphinxupquote{double Evaluate()}}
\end{quote}

\sphinxAtStartPar
\sphinxstylestrong{Return}
\begin{quote}

\sphinxAtStartPar
value of PSD expression
\end{quote}
\end{quote}

\subsubsection{PsdExpr::GetCoeff()}
\label{\detokenize{cppapi/psdexpr:psdexpr-getcoeff}}\begin{quote}

\sphinxAtStartPar
Get coefficient from the i\sphinxhyphen{}th term in PSD expression.

\sphinxAtStartPar
\sphinxstylestrong{Synopsis}
\begin{quote}

\sphinxAtStartPar
\sphinxcode{\sphinxupquote{SymMatExpr \&GetCoeff(int i)}}
\end{quote}

\sphinxAtStartPar
\sphinxstylestrong{Arguments}
\begin{quote}

\sphinxAtStartPar
\sphinxcode{\sphinxupquote{i}}: index of the PSD term.
\end{quote}

\sphinxAtStartPar
\sphinxstylestrong{Return}
\begin{quote}

\sphinxAtStartPar
coefficient of the i\sphinxhyphen{}th PSD term in PSD expression object.
\end{quote}
\end{quote}

\subsubsection{PsdExpr::GetConstant()}
\label{\detokenize{cppapi/psdexpr:psdexpr-getconstant}}\begin{quote}

\sphinxAtStartPar
Get constant in PSD expression.

\sphinxAtStartPar
\sphinxstylestrong{Synopsis}
\begin{quote}

\sphinxAtStartPar
\sphinxcode{\sphinxupquote{double GetConstant()}}
\end{quote}

\sphinxAtStartPar
\sphinxstylestrong{Return}
\begin{quote}

\sphinxAtStartPar
constant in PSD expression.
\end{quote}
\end{quote}

\subsubsection{PsdExpr::GetLinExpr()}
\label{\detokenize{cppapi/psdexpr:psdexpr-getlinexpr}}\begin{quote}

\sphinxAtStartPar
Get linear expression in PSD expression.

\sphinxAtStartPar
\sphinxstylestrong{Synopsis}
\begin{quote}

\sphinxAtStartPar
\sphinxcode{\sphinxupquote{Expr \&GetLinExpr()}}
\end{quote}

\sphinxAtStartPar
\sphinxstylestrong{Return}
\begin{quote}

\sphinxAtStartPar
linear expression object.
\end{quote}
\end{quote}

\subsubsection{PsdExpr::GetPsdVar()}
\label{\detokenize{cppapi/psdexpr:psdexpr-getpsdvar}}\begin{quote}

\sphinxAtStartPar
Get the PSD variable from the i\sphinxhyphen{}th term in PSD expression.

\sphinxAtStartPar
\sphinxstylestrong{Synopsis}
\begin{quote}

\sphinxAtStartPar
\sphinxcode{\sphinxupquote{PsdVar \&GetPsdVar(int i)}}
\end{quote}

\sphinxAtStartPar
\sphinxstylestrong{Arguments}
\begin{quote}

\sphinxAtStartPar
\sphinxcode{\sphinxupquote{i}}: index of the term.
\end{quote}

\sphinxAtStartPar
\sphinxstylestrong{Return}
\begin{quote}

\sphinxAtStartPar
the first variable of the i\sphinxhyphen{}th term in PSD expression object.
\end{quote}
\end{quote}

\subsubsection{PsdExpr::operator*=()}
\label{\detokenize{cppapi/psdexpr:psdexpr-operator}}\begin{quote}

\sphinxAtStartPar
Multiply a constant to self.

\sphinxAtStartPar
\sphinxstylestrong{Synopsis}
\begin{quote}

\sphinxAtStartPar
\sphinxcode{\sphinxupquote{void operator*=(double c)}}
\end{quote}

\sphinxAtStartPar
\sphinxstylestrong{Arguments}
\begin{quote}

\sphinxAtStartPar
\sphinxcode{\sphinxupquote{c}}: constant multiplier.
\end{quote}
\end{quote}

\subsubsection{PsdExpr::operator*()}
\label{\detokenize{cppapi/psdexpr:id9}}\begin{quote}

\sphinxAtStartPar
Multiply constant and return new expression.

\sphinxAtStartPar
\sphinxstylestrong{Synopsis}
\begin{quote}

\sphinxAtStartPar
\sphinxcode{\sphinxupquote{PsdExpr operator*(double c)}}
\end{quote}

\sphinxAtStartPar
\sphinxstylestrong{Arguments}
\begin{quote}

\sphinxAtStartPar
\sphinxcode{\sphinxupquote{c}}: constant multiplier.
\end{quote}

\sphinxAtStartPar
\sphinxstylestrong{Return}
\begin{quote}

\sphinxAtStartPar
result expression.
\end{quote}
\end{quote}

\subsubsection{PsdExpr::operator+=()}
\label{\detokenize{cppapi/psdexpr:id10}}\begin{quote}

\sphinxAtStartPar
Add an expression to self.

\sphinxAtStartPar
\sphinxstylestrong{Synopsis}
\begin{quote}

\sphinxAtStartPar
\sphinxcode{\sphinxupquote{void operator+=(const PsdExpr \&expr)}}
\end{quote}

\sphinxAtStartPar
\sphinxstylestrong{Arguments}
\begin{quote}

\sphinxAtStartPar
\sphinxcode{\sphinxupquote{expr}}: expression to be added.
\end{quote}
\end{quote}

\subsubsection{PsdExpr::operator+()}
\label{\detokenize{cppapi/psdexpr:id11}}\begin{quote}

\sphinxAtStartPar
Add expression and return new expression.

\sphinxAtStartPar
\sphinxstylestrong{Synopsis}
\begin{quote}

\sphinxAtStartPar
\sphinxcode{\sphinxupquote{PsdExpr operator+(const PsdExpr \&other)}}
\end{quote}

\sphinxAtStartPar
\sphinxstylestrong{Arguments}
\begin{quote}

\sphinxAtStartPar
\sphinxcode{\sphinxupquote{other}}: other expression to add.
\end{quote}

\sphinxAtStartPar
\sphinxstylestrong{Return}
\begin{quote}

\sphinxAtStartPar
result expression.
\end{quote}
\end{quote}

\subsubsection{PsdExpr::operator\sphinxhyphen{}=()}
\label{\detokenize{cppapi/psdexpr:id12}}\begin{quote}

\sphinxAtStartPar
Substract an expression from self.

\sphinxAtStartPar
\sphinxstylestrong{Synopsis}
\begin{quote}

\sphinxAtStartPar
\sphinxcode{\sphinxupquote{void operator\sphinxhyphen{}=(const PsdExpr \&expr)}}
\end{quote}

\sphinxAtStartPar
\sphinxstylestrong{Arguments}
\begin{quote}

\sphinxAtStartPar
\sphinxcode{\sphinxupquote{expr}}: expression to be substracted.
\end{quote}
\end{quote}

\subsubsection{PsdExpr::operator\sphinxhyphen{}()}
\label{\detokenize{cppapi/psdexpr:id13}}\begin{quote}

\sphinxAtStartPar
Substract expression and return new expression.

\sphinxAtStartPar
\sphinxstylestrong{Synopsis}
\begin{quote}

\sphinxAtStartPar
\sphinxcode{\sphinxupquote{PsdExpr operator\sphinxhyphen{}(const PsdExpr \&other)}}
\end{quote}

\sphinxAtStartPar
\sphinxstylestrong{Arguments}
\begin{quote}

\sphinxAtStartPar
\sphinxcode{\sphinxupquote{other}}: other expression to substract.
\end{quote}

\sphinxAtStartPar
\sphinxstylestrong{Return}
\begin{quote}

\sphinxAtStartPar
result expression.
\end{quote}
\end{quote}

\subsubsection{PsdExpr::Remove()}
\label{\detokenize{cppapi/psdexpr:psdexpr-remove}}\begin{quote}

\sphinxAtStartPar
Remove i\sphinxhyphen{}th term from PSD expression object.

\sphinxAtStartPar
\sphinxstylestrong{Synopsis}
\begin{quote}

\sphinxAtStartPar
\sphinxcode{\sphinxupquote{void Remove(int idx)}}
\end{quote}

\sphinxAtStartPar
\sphinxstylestrong{Arguments}
\begin{quote}

\sphinxAtStartPar
\sphinxcode{\sphinxupquote{idx}}: index of the term to be removed.
\end{quote}
\end{quote}

\subsubsection{PsdExpr::Remove()}
\label{\detokenize{cppapi/psdexpr:id14}}\begin{quote}

\sphinxAtStartPar
Remove the term associated with variable from PSD expression.

\sphinxAtStartPar
\sphinxstylestrong{Synopsis}
\begin{quote}

\sphinxAtStartPar
\sphinxcode{\sphinxupquote{void Remove(const Var \&var)}}
\end{quote}

\sphinxAtStartPar
\sphinxstylestrong{Arguments}
\begin{quote}

\sphinxAtStartPar
\sphinxcode{\sphinxupquote{var}}: a variable whose term should be removed.
\end{quote}
\end{quote}

\subsubsection{PsdExpr::Remove()}
\label{\detokenize{cppapi/psdexpr:id15}}\begin{quote}

\sphinxAtStartPar
Remove the term associated with PSD variable from PSD expression.

\sphinxAtStartPar
\sphinxstylestrong{Synopsis}
\begin{quote}

\sphinxAtStartPar
\sphinxcode{\sphinxupquote{void Remove(const PsdVar \&var)}}
\end{quote}

\sphinxAtStartPar
\sphinxstylestrong{Arguments}
\begin{quote}

\sphinxAtStartPar
\sphinxcode{\sphinxupquote{var}}: a PSD variable whose term should be removed.
\end{quote}
\end{quote}

\subsubsection{PsdExpr::Reserve()}
\label{\detokenize{cppapi/psdexpr:psdexpr-reserve}}\begin{quote}

\sphinxAtStartPar
Reserve capacity to contain at least n items.

\sphinxAtStartPar
\sphinxstylestrong{Synopsis}
\begin{quote}

\sphinxAtStartPar
\sphinxcode{\sphinxupquote{void Reserve(size\_t n)}}
\end{quote}

\sphinxAtStartPar
\sphinxstylestrong{Arguments}
\begin{quote}

\sphinxAtStartPar
\sphinxcode{\sphinxupquote{n}}: minimum capacity for PSD expression object.
\end{quote}
\end{quote}

\subsubsection{PsdExpr::SetCoeff()}
\label{\detokenize{cppapi/psdexpr:psdexpr-setcoeff}}\begin{quote}

\sphinxAtStartPar
Set coefficient matrix of the i\sphinxhyphen{}th term in PSD expression.

\sphinxAtStartPar
\sphinxstylestrong{Synopsis}
\begin{quote}

\sphinxAtStartPar
\sphinxcode{\sphinxupquote{void SetCoeff(int i, const SymMatrix \&mat)}}
\end{quote}

\sphinxAtStartPar
\sphinxstylestrong{Arguments}
\begin{quote}

\sphinxAtStartPar
\sphinxcode{\sphinxupquote{i}}: index of the PSD term.

\sphinxAtStartPar
\sphinxcode{\sphinxupquote{mat}}: coefficient matrix of the term.
\end{quote}
\end{quote}

\subsubsection{PsdExpr::SetConstant()}
\label{\detokenize{cppapi/psdexpr:psdexpr-setconstant}}\begin{quote}

\sphinxAtStartPar
Set constant for the PSD expression.

\sphinxAtStartPar
\sphinxstylestrong{Synopsis}
\begin{quote}

\sphinxAtStartPar
\sphinxcode{\sphinxupquote{void SetConstant(double constant)}}
\end{quote}

\sphinxAtStartPar
\sphinxstylestrong{Arguments}
\begin{quote}

\sphinxAtStartPar
\sphinxcode{\sphinxupquote{constant}}: the value of the constant.
\end{quote}
\end{quote}

\subsubsection{PsdExpr::Size()}
\label{\detokenize{cppapi/psdexpr:psdexpr-size}}\begin{quote}

\sphinxAtStartPar
Get number of PSD terms in expression.

\sphinxAtStartPar
\sphinxstylestrong{Synopsis}
\begin{quote}

\sphinxAtStartPar
\sphinxcode{\sphinxupquote{size\_t Size()}}
\end{quote}

\sphinxAtStartPar
\sphinxstylestrong{Return}
\begin{quote}

\sphinxAtStartPar
number of PSD terms.
\end{quote}
\end{quote}

\subsection{PsdConstraint}
\label{\detokenize{cppapiref:psdconstraint}}\label{\detokenize{cppapiref:chapcppapiref-psdconstraint}}
\sphinxAtStartPar
COPT PSD constraint object. PSD constraints are always associated with a
particular model.  User creates a PSD constraint object by adding a PSD
constraint to model, rather than by constructor of PsdConstraint class.

\sphinxstepscope

\subsubsection{PsdConstraint::Get()}
\label{\detokenize{cppapi/psdconstraint:psdconstraint-get}}\label{\detokenize{cppapi/psdconstraint::doc}}\begin{quote}

\sphinxAtStartPar
Get information value of the PSD constraint. Support related PSd informations.

\sphinxAtStartPar
\sphinxstylestrong{Synopsis}
\begin{quote}

\sphinxAtStartPar
\sphinxcode{\sphinxupquote{double Get(const char *szInfo)}}
\end{quote}

\sphinxAtStartPar
\sphinxstylestrong{Arguments}
\begin{quote}

\sphinxAtStartPar
\sphinxcode{\sphinxupquote{szInfo}}: name of queried information.
\end{quote}

\sphinxAtStartPar
\sphinxstylestrong{Return}
\begin{quote}

\sphinxAtStartPar
value of information.
\end{quote}
\end{quote}

\subsubsection{PsdConstraint::GetIdx()}
\label{\detokenize{cppapi/psdconstraint:psdconstraint-getidx}}\begin{quote}

\sphinxAtStartPar
Get index of the PSD constraint.

\sphinxAtStartPar
\sphinxstylestrong{Synopsis}
\begin{quote}

\sphinxAtStartPar
\sphinxcode{\sphinxupquote{int GetIdx()}}
\end{quote}

\sphinxAtStartPar
\sphinxstylestrong{Return}
\begin{quote}

\sphinxAtStartPar
the index of the PSD constraint.
\end{quote}
\end{quote}

\subsubsection{PsdConstraint::GetName()}
\label{\detokenize{cppapi/psdconstraint:psdconstraint-getname}}\begin{quote}

\sphinxAtStartPar
Get name of the PSD constraint.

\sphinxAtStartPar
\sphinxstylestrong{Synopsis}
\begin{quote}

\sphinxAtStartPar
\sphinxcode{\sphinxupquote{const char *GetName()}}
\end{quote}

\sphinxAtStartPar
\sphinxstylestrong{Return}
\begin{quote}

\sphinxAtStartPar
the name of the PSD constraint.
\end{quote}
\end{quote}

\subsubsection{PsdConstraint::Remove()}
\label{\detokenize{cppapi/psdconstraint:psdconstraint-remove}}\begin{quote}

\sphinxAtStartPar
Remove this PSD constraint from model.

\sphinxAtStartPar
\sphinxstylestrong{Synopsis}
\begin{quote}

\sphinxAtStartPar
\sphinxcode{\sphinxupquote{void Remove()}}
\end{quote}
\end{quote}

\subsubsection{PsdConstraint::Set()}
\label{\detokenize{cppapi/psdconstraint:psdconstraint-set}}\begin{quote}

\sphinxAtStartPar
Set information value of the PSD constraint. Support related PSD informations.

\sphinxAtStartPar
\sphinxstylestrong{Synopsis}
\begin{quote}

\sphinxAtStartPar
\sphinxcode{\sphinxupquote{void Set(const char *szInfo, double value)}}
\end{quote}

\sphinxAtStartPar
\sphinxstylestrong{Arguments}
\begin{quote}

\sphinxAtStartPar
\sphinxcode{\sphinxupquote{szInfo}}: name of queried information.

\sphinxAtStartPar
\sphinxcode{\sphinxupquote{value}}: new information value.
\end{quote}
\end{quote}

\subsubsection{PsdConstraint::SetName()}
\label{\detokenize{cppapi/psdconstraint:psdconstraint-setname}}\begin{quote}

\sphinxAtStartPar
Set name of a PSD constraint.

\sphinxAtStartPar
\sphinxstylestrong{Synopsis}
\begin{quote}

\sphinxAtStartPar
\sphinxcode{\sphinxupquote{void SetName(const char *szName)}}
\end{quote}

\sphinxAtStartPar
\sphinxstylestrong{Arguments}
\begin{quote}

\sphinxAtStartPar
\sphinxcode{\sphinxupquote{szName}}: the name to set.
\end{quote}
\end{quote}

\subsection{PsdConstrArray}
\label{\detokenize{cppapiref:psdconstrarray}}\label{\detokenize{cppapiref:chapcppapiref-psdconstrarray}}
\sphinxAtStartPar
COPT PSD constraint array object. To store and access a set of
{\hyperref[\detokenize{cppapiref:chapcppapiref-psdconstraint}]{\sphinxcrossref{\DUrole{std,std-ref}{PsdConstraint}}}} objects, Cardinal Optimizer provides
PsdConstrArray class, which defines the following methods.

\sphinxstepscope

\subsubsection{PsdConstrArray::GetPsdConstr()}
\label{\detokenize{cppapi/psdconstrarray:psdconstrarray-getpsdconstr}}\label{\detokenize{cppapi/psdconstrarray::doc}}\begin{quote}

\sphinxAtStartPar
Get idx\sphinxhyphen{}th PSD constraint object.

\sphinxAtStartPar
\sphinxstylestrong{Synopsis}
\begin{quote}

\sphinxAtStartPar
\sphinxcode{\sphinxupquote{PsdConstraint \&GetPsdConstr(int idx)}}
\end{quote}

\sphinxAtStartPar
\sphinxstylestrong{Arguments}
\begin{quote}

\sphinxAtStartPar
\sphinxcode{\sphinxupquote{idx}}: index of the PSD constraint.
\end{quote}

\sphinxAtStartPar
\sphinxstylestrong{Return}
\begin{quote}

\sphinxAtStartPar
PSD constraint object with index idx.
\end{quote}
\end{quote}

\subsubsection{PsdConstrArray::PushBack()}
\label{\detokenize{cppapi/psdconstrarray:psdconstrarray-pushback}}\begin{quote}

\sphinxAtStartPar
Add a PSD constraint object to PSD constraint array.

\sphinxAtStartPar
\sphinxstylestrong{Synopsis}
\begin{quote}

\sphinxAtStartPar
\sphinxcode{\sphinxupquote{void PushBack(const PsdConstraint \&constr)}}
\end{quote}

\sphinxAtStartPar
\sphinxstylestrong{Arguments}
\begin{quote}

\sphinxAtStartPar
\sphinxcode{\sphinxupquote{constr}}: a PSD constraint object.
\end{quote}
\end{quote}

\subsubsection{PsdConstrArray::Reserve()}
\label{\detokenize{cppapi/psdconstrarray:psdconstrarray-reserve}}\begin{quote}

\sphinxAtStartPar
Reserve capacity to contain at least n items.

\sphinxAtStartPar
\sphinxstylestrong{Synopsis}
\begin{quote}

\sphinxAtStartPar
\sphinxcode{\sphinxupquote{void Reserve(int n)}}
\end{quote}

\sphinxAtStartPar
\sphinxstylestrong{Arguments}
\begin{quote}

\sphinxAtStartPar
\sphinxcode{\sphinxupquote{n}}: minimum capacity for PSD constraint objects.
\end{quote}
\end{quote}

\subsubsection{PsdConstrArray::Size()}
\label{\detokenize{cppapi/psdconstrarray:psdconstrarray-size}}\begin{quote}

\sphinxAtStartPar
Get the number of PSD constraint objects.

\sphinxAtStartPar
\sphinxstylestrong{Synopsis}
\begin{quote}

\sphinxAtStartPar
\sphinxcode{\sphinxupquote{int Size()}}
\end{quote}

\sphinxAtStartPar
\sphinxstylestrong{Return}
\begin{quote}

\sphinxAtStartPar
number of PSD constraint objects.
\end{quote}
\end{quote}

\subsection{PsdConstrBuilder}
\label{\detokenize{cppapiref:psdconstrbuilder}}\label{\detokenize{cppapiref:chapcppapiref-psdconstrbuilder}}
\sphinxAtStartPar
COPT PSD constraint builder object. To help building a PSD constraint,
given a PSD expression, constraint sense and right\sphinxhyphen{}hand side value, Cardinal
Optimizer provides PsdConstrBuilder class, which defines the following methods.

\sphinxstepscope

\subsubsection{PsdConstrBuilder::GetPsdExpr()}
\label{\detokenize{cppapi/psdconstrbuilder:psdconstrbuilder-getpsdexpr}}\label{\detokenize{cppapi/psdconstrbuilder::doc}}\begin{quote}

\sphinxAtStartPar
Get expression associated with PSD constraint.

\sphinxAtStartPar
\sphinxstylestrong{Synopsis}
\begin{quote}

\sphinxAtStartPar
\sphinxcode{\sphinxupquote{const PsdExpr \&GetPsdExpr()}}
\end{quote}

\sphinxAtStartPar
\sphinxstylestrong{Return}
\begin{quote}

\sphinxAtStartPar
PSD expression object.
\end{quote}
\end{quote}

\subsubsection{PsdConstrBuilder::GetRange()}
\label{\detokenize{cppapi/psdconstrbuilder:psdconstrbuilder-getrange}}\begin{quote}

\sphinxAtStartPar
Get range from lower bound to upper bound of range constraint.

\sphinxAtStartPar
\sphinxstylestrong{Synopsis}
\begin{quote}

\sphinxAtStartPar
\sphinxcode{\sphinxupquote{double GetRange()}}
\end{quote}

\sphinxAtStartPar
\sphinxstylestrong{Return}
\begin{quote}

\sphinxAtStartPar
length from lower bound to upper bound of the constraint.
\end{quote}
\end{quote}

\subsubsection{PsdConstrBuilder::GetSense()}
\label{\detokenize{cppapi/psdconstrbuilder:psdconstrbuilder-getsense}}\begin{quote}

\sphinxAtStartPar
Get sense associated with PSD constraint.

\sphinxAtStartPar
\sphinxstylestrong{Synopsis}
\begin{quote}

\sphinxAtStartPar
\sphinxcode{\sphinxupquote{char GetSense()}}
\end{quote}

\sphinxAtStartPar
\sphinxstylestrong{Return}
\begin{quote}

\sphinxAtStartPar
PSD constraint sense.
\end{quote}
\end{quote}

\subsubsection{PsdConstrBuilder::Set()}
\label{\detokenize{cppapi/psdconstrbuilder:psdconstrbuilder-set}}\begin{quote}

\sphinxAtStartPar
Set detail of a PSD constraint to its builder object.

\sphinxAtStartPar
\sphinxstylestrong{Synopsis}
\begin{quote}

\sphinxAtStartPar
\sphinxcode{\sphinxupquote{void Set(}}
\begin{quote}

\sphinxAtStartPar
\sphinxcode{\sphinxupquote{const PsdExpr \&expr,}}

\sphinxAtStartPar
\sphinxcode{\sphinxupquote{char sense,}}

\sphinxAtStartPar
\sphinxcode{\sphinxupquote{double rhs)}}
\end{quote}
\end{quote}

\sphinxAtStartPar
\sphinxstylestrong{Arguments}
\begin{quote}

\sphinxAtStartPar
\sphinxcode{\sphinxupquote{expr}}: expression object at one side of the PSD constraint.

\sphinxAtStartPar
\sphinxcode{\sphinxupquote{sense}}: PSD constraint sense, other than COPT\_RANGE.

\sphinxAtStartPar
\sphinxcode{\sphinxupquote{rhs}}: constant at right side of the PSD constraint.
\end{quote}
\end{quote}

\subsubsection{PsdConstrBuilder::SetRange()}
\label{\detokenize{cppapi/psdconstrbuilder:psdconstrbuilder-setrange}}\begin{quote}

\sphinxAtStartPar
Set a range constraint to its builder.

\sphinxAtStartPar
\sphinxstylestrong{Synopsis}
\begin{quote}

\sphinxAtStartPar
\sphinxcode{\sphinxupquote{void SetRange(const PsdExpr \&expr, double range)}}
\end{quote}

\sphinxAtStartPar
\sphinxstylestrong{Arguments}
\begin{quote}

\sphinxAtStartPar
\sphinxcode{\sphinxupquote{expr}}: PSD expression object, whose constant is negative upper bound.

\sphinxAtStartPar
\sphinxcode{\sphinxupquote{range}}: length from lower bound to upper bound of the constraint. Must greater than 0.
\end{quote}
\end{quote}

\subsection{PsdConstrBuilderArray}
\label{\detokenize{cppapiref:psdconstrbuilderarray}}\label{\detokenize{cppapiref:chapcppapiref-psdconstrbuilderarray}}
\sphinxAtStartPar
COPT PSD constraint builder array object. To store and access a set of
{\hyperref[\detokenize{cppapiref:chapcppapiref-psdconstrbuilder}]{\sphinxcrossref{\DUrole{std,std-ref}{PsdConstrBuilder}}}} objects, Cardinal Optimizer provides
PsdConstrBuilderArray class, which defines the following methods.

\sphinxstepscope

\subsubsection{PsdConstrBuilderArray::GetBuilder()}
\label{\detokenize{cppapi/psdconstrbuilderarray:psdconstrbuilderarray-getbuilder}}\label{\detokenize{cppapi/psdconstrbuilderarray::doc}}\begin{quote}

\sphinxAtStartPar
Get idx\sphinxhyphen{}th PSD constraint builder object.

\sphinxAtStartPar
\sphinxstylestrong{Synopsis}
\begin{quote}

\sphinxAtStartPar
\sphinxcode{\sphinxupquote{PsdConstrBuilder \&GetBuilder(int idx)}}
\end{quote}

\sphinxAtStartPar
\sphinxstylestrong{Arguments}
\begin{quote}

\sphinxAtStartPar
\sphinxcode{\sphinxupquote{idx}}: index of the PSD constraint builder.
\end{quote}

\sphinxAtStartPar
\sphinxstylestrong{Return}
\begin{quote}

\sphinxAtStartPar
PSD constraint builder object with index idx.
\end{quote}
\end{quote}

\subsubsection{PsdConstrBuilderArray::PushBack()}
\label{\detokenize{cppapi/psdconstrbuilderarray:psdconstrbuilderarray-pushback}}\begin{quote}

\sphinxAtStartPar
Add a PSD constraint builder object to PSD constraint builder array.

\sphinxAtStartPar
\sphinxstylestrong{Synopsis}
\begin{quote}

\sphinxAtStartPar
\sphinxcode{\sphinxupquote{void PushBack(const PsdConstrBuilder \&builder)}}
\end{quote}

\sphinxAtStartPar
\sphinxstylestrong{Arguments}
\begin{quote}

\sphinxAtStartPar
\sphinxcode{\sphinxupquote{builder}}: a PSD constraint builder object.
\end{quote}
\end{quote}

\subsubsection{PsdConstrBuilderArray::Reserve()}
\label{\detokenize{cppapi/psdconstrbuilderarray:psdconstrbuilderarray-reserve}}\begin{quote}

\sphinxAtStartPar
Reserve capacity to contain at least n items.

\sphinxAtStartPar
\sphinxstylestrong{Synopsis}
\begin{quote}

\sphinxAtStartPar
\sphinxcode{\sphinxupquote{void Reserve(int n)}}
\end{quote}

\sphinxAtStartPar
\sphinxstylestrong{Arguments}
\begin{quote}

\sphinxAtStartPar
\sphinxcode{\sphinxupquote{n}}: minimum capacity for PSD constraint builder object.
\end{quote}
\end{quote}

\subsubsection{PsdConstrBuilderArray::Size()}
\label{\detokenize{cppapi/psdconstrbuilderarray:psdconstrbuilderarray-size}}\begin{quote}

\sphinxAtStartPar
Get the number of PSD constraint builder objects.

\sphinxAtStartPar
\sphinxstylestrong{Synopsis}
\begin{quote}

\sphinxAtStartPar
\sphinxcode{\sphinxupquote{int Size()}}
\end{quote}

\sphinxAtStartPar
\sphinxstylestrong{Return}
\begin{quote}

\sphinxAtStartPar
number of PSD constraint builder objects.
\end{quote}
\end{quote}

\subsection{LmiConstraint}
\label{\detokenize{cppapiref:lmiconstraint}}\label{\detokenize{cppapiref:chapcppapiref-lmiconstraint}}
\sphinxAtStartPar
COPT LMI constraint object. LMI constraints are always associated with a
particular model.  User creates a LMI constraint object by adding a LMI
constraint to model, rather than by constructor of LmiConstraint class.

\sphinxstepscope

\subsubsection{LmiConstraint::Get()}
\label{\detokenize{cppapi/lmiconstraint:lmiconstraint-get}}\label{\detokenize{cppapi/lmiconstraint::doc}}\begin{quote}

\sphinxAtStartPar
Get information values of LMI expression.

\sphinxAtStartPar
\sphinxstylestrong{Synopsis}
\begin{quote}

\sphinxAtStartPar
\sphinxcode{\sphinxupquote{double Get(const char *szInfo, int len)}}
\end{quote}

\sphinxAtStartPar
\sphinxstylestrong{Arguments}
\begin{quote}

\sphinxAtStartPar
\sphinxcode{\sphinxupquote{szInfo}}: name of queried information.

\sphinxAtStartPar
\sphinxcode{\sphinxupquote{len}}: length of output array.
\end{quote}

\sphinxAtStartPar
\sphinxstylestrong{Return}
\begin{quote}

\sphinxAtStartPar
output list of information values.
\end{quote}
\end{quote}

\subsubsection{LmiConstraint::GetDim()}
\label{\detokenize{cppapi/lmiconstraint:lmiconstraint-getdim}}\begin{quote}

\sphinxAtStartPar
Get dimension of LMI constraint.

\sphinxAtStartPar
\sphinxstylestrong{Synopsis}
\begin{quote}

\sphinxAtStartPar
\sphinxcode{\sphinxupquote{int GetDim()}}
\end{quote}

\sphinxAtStartPar
\sphinxstylestrong{Return}
\begin{quote}

\sphinxAtStartPar
dimension of LMI constraint.
\end{quote}
\end{quote}

\subsubsection{LmiConstraint::GetIdx()}
\label{\detokenize{cppapi/lmiconstraint:lmiconstraint-getidx}}\begin{quote}

\sphinxAtStartPar
Get index of LMI constraint.

\sphinxAtStartPar
\sphinxstylestrong{Synopsis}
\begin{quote}

\sphinxAtStartPar
\sphinxcode{\sphinxupquote{int GetIdx()}}
\end{quote}

\sphinxAtStartPar
\sphinxstylestrong{Return}
\begin{quote}

\sphinxAtStartPar
index of LMI constraint.
\end{quote}
\end{quote}

\subsubsection{LmiConstraint::GetLen()}
\label{\detokenize{cppapi/lmiconstraint:lmiconstraint-getlen}}\begin{quote}

\sphinxAtStartPar
Get length of LMI constraint.

\sphinxAtStartPar
\sphinxstylestrong{Synopsis}
\begin{quote}

\sphinxAtStartPar
\sphinxcode{\sphinxupquote{int GetLen()}}
\end{quote}

\sphinxAtStartPar
\sphinxstylestrong{Return}
\begin{quote}

\sphinxAtStartPar
length of LMI constraint.
\end{quote}
\end{quote}

\subsubsection{LmiConstraint::GetName()}
\label{\detokenize{cppapi/lmiconstraint:lmiconstraint-getname}}\begin{quote}

\sphinxAtStartPar
Get name of LMI constraint.

\sphinxAtStartPar
\sphinxstylestrong{Synopsis}
\begin{quote}

\sphinxAtStartPar
\sphinxcode{\sphinxupquote{const char *GetName()}}
\end{quote}

\sphinxAtStartPar
\sphinxstylestrong{Return}
\begin{quote}

\sphinxAtStartPar
the name of LMI constraint.
\end{quote}
\end{quote}

\subsubsection{LmiConstraint::Remove()}
\label{\detokenize{cppapi/lmiconstraint:lmiconstraint-remove}}\begin{quote}

\sphinxAtStartPar
Remove this LMI constraint from model.

\sphinxAtStartPar
\sphinxstylestrong{Synopsis}
\begin{quote}

\sphinxAtStartPar
\sphinxcode{\sphinxupquote{void Remove()}}
\end{quote}
\end{quote}

\subsubsection{LmiConstraint::SetName()}
\label{\detokenize{cppapi/lmiconstraint:lmiconstraint-setname}}\begin{quote}

\sphinxAtStartPar
Set name of LMI constraint.

\sphinxAtStartPar
\sphinxstylestrong{Synopsis}
\begin{quote}

\sphinxAtStartPar
\sphinxcode{\sphinxupquote{void SetName(const char *szName)}}
\end{quote}

\sphinxAtStartPar
\sphinxstylestrong{Arguments}
\begin{quote}

\sphinxAtStartPar
\sphinxcode{\sphinxupquote{szName}}: new name to set.
\end{quote}
\end{quote}

\subsubsection{LmiConstraint::SetRhs()}
\label{\detokenize{cppapi/lmiconstraint:lmiconstraint-setrhs}}\begin{quote}

\sphinxAtStartPar
Set constant term of LMI constraint.

\sphinxAtStartPar
\sphinxstylestrong{Synopsis}
\begin{quote}

\sphinxAtStartPar
\sphinxcode{\sphinxupquote{void SetRhs(const SymMatrix \&mat)}}
\end{quote}

\sphinxAtStartPar
\sphinxstylestrong{Arguments}
\begin{quote}

\sphinxAtStartPar
\sphinxcode{\sphinxupquote{mat}}: new symmetric matrix for constant term.
\end{quote}
\end{quote}

\subsection{LmiConstrArray}
\label{\detokenize{cppapiref:lmiconstrarray}}\label{\detokenize{cppapiref:chapcppapiref-lmiconstrarray}}
\sphinxAtStartPar
COPT LMI constraint array object. To store and access a set of
{\hyperref[\detokenize{cppapiref:chapcppapiref-lmiconstraint}]{\sphinxcrossref{\DUrole{std,std-ref}{LmiConstraint}}}} objects, Cardinal Optimizer provides
LmiConstrArray class, which defines the following methods.

\sphinxstepscope

\subsubsection{LmiConstrArray::GetLmiConstr()}
\label{\detokenize{cppapi/lmiconstrarray:lmiconstrarray-getlmiconstr}}\label{\detokenize{cppapi/lmiconstrarray::doc}}\begin{quote}

\sphinxAtStartPar
Get idx\sphinxhyphen{}th LMI constraint object.

\sphinxAtStartPar
\sphinxstylestrong{Synopsis}
\begin{quote}

\sphinxAtStartPar
\sphinxcode{\sphinxupquote{LmiConstraint \&GetLmiConstr(int idx)}}
\end{quote}

\sphinxAtStartPar
\sphinxstylestrong{Arguments}
\begin{quote}

\sphinxAtStartPar
\sphinxcode{\sphinxupquote{idx}}: index of the LMI constraint.
\end{quote}

\sphinxAtStartPar
\sphinxstylestrong{Return}
\begin{quote}

\sphinxAtStartPar
LMI constraint object with index idx.
\end{quote}
\end{quote}

\subsubsection{LmiConstrArray::PushBack()}
\label{\detokenize{cppapi/lmiconstrarray:lmiconstrarray-pushback}}\begin{quote}

\sphinxAtStartPar
Add a LMI constraint to LMI constraint array.

\sphinxAtStartPar
\sphinxstylestrong{Synopsis}
\begin{quote}

\sphinxAtStartPar
\sphinxcode{\sphinxupquote{void PushBack(const LmiConstraint \&constr)}}
\end{quote}

\sphinxAtStartPar
\sphinxstylestrong{Arguments}
\begin{quote}

\sphinxAtStartPar
\sphinxcode{\sphinxupquote{constr}}: LMI constraint object.
\end{quote}
\end{quote}

\subsubsection{LmiConstrArray::Reserve()}
\label{\detokenize{cppapi/lmiconstrarray:lmiconstrarray-reserve}}\begin{quote}

\sphinxAtStartPar
Reserve capacity to contain at least n items.

\sphinxAtStartPar
\sphinxstylestrong{Synopsis}
\begin{quote}

\sphinxAtStartPar
\sphinxcode{\sphinxupquote{void Reserve(int n)}}
\end{quote}

\sphinxAtStartPar
\sphinxstylestrong{Arguments}
\begin{quote}

\sphinxAtStartPar
\sphinxcode{\sphinxupquote{n}}: capacity number of LMI constraint objects.
\end{quote}
\end{quote}

\subsubsection{LmiConstrArray::Size()}
\label{\detokenize{cppapi/lmiconstrarray:lmiconstrarray-size}}\begin{quote}

\sphinxAtStartPar
Get the number of LMI constraint objects.

\sphinxAtStartPar
\sphinxstylestrong{Synopsis}
\begin{quote}

\sphinxAtStartPar
\sphinxcode{\sphinxupquote{int Size()}}
\end{quote}

\sphinxAtStartPar
\sphinxstylestrong{Return}
\begin{quote}

\sphinxAtStartPar
number of LMI constraint objects.
\end{quote}
\end{quote}

\subsection{LmiExpr}
\label{\detokenize{cppapiref:lmiexpr}}\label{\detokenize{cppapiref:chapcppapiref-lmiexpr}}
\sphinxAtStartPar
COPT LMI expression object. A LMI expression consists of
a list of variables, associated coefficient matrices of LMI term, and constant matrices.
LMI expressions are used to build LMI constraints.

\sphinxstepscope

\subsubsection{LmiExpr::LmiExpr()}
\label{\detokenize{cppapi/lmiexpr:lmiexpr-lmiexpr}}\label{\detokenize{cppapi/lmiexpr::doc}}\begin{quote}

\sphinxAtStartPar
Default constructor of a LMI expression.

\sphinxAtStartPar
\sphinxstylestrong{Synopsis}
\begin{quote}

\sphinxAtStartPar
\sphinxcode{\sphinxupquote{LmiExpr()}}
\end{quote}
\end{quote}

\subsubsection{LmiExpr::LmiExpr()}
\label{\detokenize{cppapi/lmiexpr:id1}}\begin{quote}

\sphinxAtStartPar
Constructor of LMI expression with given symmetric matrix.

\sphinxAtStartPar
\sphinxstylestrong{Synopsis}
\begin{quote}

\sphinxAtStartPar
\sphinxcode{\sphinxupquote{LmiExpr(const SymMatrix \&mat)}}
\end{quote}

\sphinxAtStartPar
\sphinxstylestrong{Arguments}
\begin{quote}

\sphinxAtStartPar
\sphinxcode{\sphinxupquote{mat}}: symmetric matrix as constant term of LMI expression.
\end{quote}
\end{quote}

\subsubsection{LmiExpr::LmiExpr()}
\label{\detokenize{cppapi/lmiexpr:id2}}\begin{quote}

\sphinxAtStartPar
Constructor of LMI expression with given matrix expression.

\sphinxAtStartPar
\sphinxstylestrong{Synopsis}
\begin{quote}

\sphinxAtStartPar
\sphinxcode{\sphinxupquote{LmiExpr(const SymMatExpr \&expr)}}
\end{quote}

\sphinxAtStartPar
\sphinxstylestrong{Arguments}
\begin{quote}

\sphinxAtStartPar
\sphinxcode{\sphinxupquote{expr}}: matrix expression as constant term of LMI expression.
\end{quote}
\end{quote}

\subsubsection{LmiExpr::LmiExpr()}
\label{\detokenize{cppapi/lmiexpr:id3}}\begin{quote}

\sphinxAtStartPar
Constructor of LMI expression with one term.

\sphinxAtStartPar
\sphinxstylestrong{Synopsis}
\begin{quote}

\sphinxAtStartPar
\sphinxcode{\sphinxupquote{LmiExpr(const Var \&var, const SymMatrix \&mat)}}
\end{quote}

\sphinxAtStartPar
\sphinxstylestrong{Arguments}
\begin{quote}

\sphinxAtStartPar
\sphinxcode{\sphinxupquote{var}}: variable of the added term.

\sphinxAtStartPar
\sphinxcode{\sphinxupquote{mat}}: coefficient matrix of the added term.
\end{quote}
\end{quote}

\subsubsection{LmiExpr::LmiExpr()}
\label{\detokenize{cppapi/lmiexpr:id4}}\begin{quote}

\sphinxAtStartPar
Constructor of LMI expression with one term.

\sphinxAtStartPar
\sphinxstylestrong{Synopsis}
\begin{quote}

\sphinxAtStartPar
\sphinxcode{\sphinxupquote{LmiExpr(const Var \&var, const SymMatExpr \&expr)}}
\end{quote}

\sphinxAtStartPar
\sphinxstylestrong{Arguments}
\begin{quote}

\sphinxAtStartPar
\sphinxcode{\sphinxupquote{var}}: variable of the added term.

\sphinxAtStartPar
\sphinxcode{\sphinxupquote{expr}}: coefficient expression of symmetric matrices for the added term.
\end{quote}
\end{quote}

\subsubsection{LmiExpr::AddConstant()}
\label{\detokenize{cppapi/lmiexpr:lmiexpr-addconstant}}\begin{quote}

\sphinxAtStartPar
Add to constant term of LMI expression.

\sphinxAtStartPar
\sphinxstylestrong{Synopsis}
\begin{quote}

\sphinxAtStartPar
\sphinxcode{\sphinxupquote{void AddConstant(const SymMatExpr \&expr)}}
\end{quote}

\sphinxAtStartPar
\sphinxstylestrong{Arguments}
\begin{quote}

\sphinxAtStartPar
\sphinxcode{\sphinxupquote{expr}}: matrix expression object added to constant term.
\end{quote}
\end{quote}

\subsubsection{LmiExpr::AddLmiExpr()}
\label{\detokenize{cppapi/lmiexpr:lmiexpr-addlmiexpr}}\begin{quote}

\sphinxAtStartPar
Add a LMI expression to self.

\sphinxAtStartPar
\sphinxstylestrong{Synopsis}
\begin{quote}

\sphinxAtStartPar
\sphinxcode{\sphinxupquote{void AddLmiExpr(const LmiExpr \&expr, double mult)}}
\end{quote}

\sphinxAtStartPar
\sphinxstylestrong{Arguments}
\begin{quote}

\sphinxAtStartPar
\sphinxcode{\sphinxupquote{expr}}: LMI expression to be added.

\sphinxAtStartPar
\sphinxcode{\sphinxupquote{mult}}: optional, constant multiplier, default value is 1.0.
\end{quote}
\end{quote}

\subsubsection{LmiExpr::AddTerm()}
\label{\detokenize{cppapi/lmiexpr:lmiexpr-addterm}}\begin{quote}

\sphinxAtStartPar
Add a term to LMI expression.

\sphinxAtStartPar
\sphinxstylestrong{Synopsis}
\begin{quote}

\sphinxAtStartPar
\sphinxcode{\sphinxupquote{void AddTerm(const Var \&var, const SymMatrix \&mat)}}
\end{quote}

\sphinxAtStartPar
\sphinxstylestrong{Arguments}
\begin{quote}

\sphinxAtStartPar
\sphinxcode{\sphinxupquote{var}}: variable of new LMI term.

\sphinxAtStartPar
\sphinxcode{\sphinxupquote{mat}}: coefficient matrix object of new LMI term.
\end{quote}
\end{quote}

\subsubsection{LmiExpr::AddTerm()}
\label{\detokenize{cppapi/lmiexpr:id5}}\begin{quote}

\sphinxAtStartPar
Add a term to LMI expression.

\sphinxAtStartPar
\sphinxstylestrong{Synopsis}
\begin{quote}

\sphinxAtStartPar
\sphinxcode{\sphinxupquote{void AddTerm(const Var \&var, const SymMatExpr \&expr)}}
\end{quote}

\sphinxAtStartPar
\sphinxstylestrong{Arguments}
\begin{quote}

\sphinxAtStartPar
\sphinxcode{\sphinxupquote{var}}: variable of new LMI term.

\sphinxAtStartPar
\sphinxcode{\sphinxupquote{expr}}: coefficient expression object of symmetric matrices of new LMI term.
\end{quote}
\end{quote}

\subsubsection{LmiExpr::AddTerms()}
\label{\detokenize{cppapi/lmiexpr:lmiexpr-addterms}}\begin{quote}

\sphinxAtStartPar
Add LMI terms to LMI expression.

\sphinxAtStartPar
\sphinxstylestrong{Synopsis}
\begin{quote}

\sphinxAtStartPar
\sphinxcode{\sphinxupquote{int AddTerms(const VarArray \&vars, const SymMatrixArray \&mats)}}
\end{quote}

\sphinxAtStartPar
\sphinxstylestrong{Arguments}
\begin{quote}

\sphinxAtStartPar
\sphinxcode{\sphinxupquote{vars}}: variables of added LMI terms.

\sphinxAtStartPar
\sphinxcode{\sphinxupquote{mats}}: coefficient matrix objects for added LMI terms.
\end{quote}

\sphinxAtStartPar
\sphinxstylestrong{Return}
\begin{quote}

\sphinxAtStartPar
number of added LMI terms.
\end{quote}
\end{quote}

\subsubsection{LmiExpr::Clone()}
\label{\detokenize{cppapi/lmiexpr:lmiexpr-clone}}\begin{quote}

\sphinxAtStartPar
Deep copy LMI expression.

\sphinxAtStartPar
\sphinxstylestrong{Synopsis}
\begin{quote}

\sphinxAtStartPar
\sphinxcode{\sphinxupquote{LmiExpr Clone()}}
\end{quote}

\sphinxAtStartPar
\sphinxstylestrong{Return}
\begin{quote}

\sphinxAtStartPar
cloned LMI expression object.
\end{quote}
\end{quote}

\subsubsection{LmiExpr::GetCoeff()}
\label{\detokenize{cppapi/lmiexpr:lmiexpr-getcoeff}}\begin{quote}

\sphinxAtStartPar
Get coefficient from the i\sphinxhyphen{}th term in LMI expression.

\sphinxAtStartPar
\sphinxstylestrong{Synopsis}
\begin{quote}

\sphinxAtStartPar
\sphinxcode{\sphinxupquote{SymMatExpr \&GetCoeff(int i)}}
\end{quote}

\sphinxAtStartPar
\sphinxstylestrong{Arguments}
\begin{quote}

\sphinxAtStartPar
\sphinxcode{\sphinxupquote{i}}: index of the LMI term.
\end{quote}

\sphinxAtStartPar
\sphinxstylestrong{Return}
\begin{quote}

\sphinxAtStartPar
coefficient matrix expression object of the i\sphinxhyphen{}th LMI term in LMI expression.
\end{quote}
\end{quote}

\subsubsection{LmiExpr::GetConstant()}
\label{\detokenize{cppapi/lmiexpr:lmiexpr-getconstant}}\begin{quote}

\sphinxAtStartPar
Get constant term in LMI expression.

\sphinxAtStartPar
\sphinxstylestrong{Synopsis}
\begin{quote}

\sphinxAtStartPar
\sphinxcode{\sphinxupquote{SymMatExpr \&GetConstant()}}
\end{quote}

\sphinxAtStartPar
\sphinxstylestrong{Return}
\begin{quote}

\sphinxAtStartPar
matrix expression object in LMI expression.
\end{quote}
\end{quote}

\subsubsection{LmiExpr::GetVar()}
\label{\detokenize{cppapi/lmiexpr:lmiexpr-getvar}}\begin{quote}

\sphinxAtStartPar
Get variable from the i\sphinxhyphen{}th term in LMI expression.

\sphinxAtStartPar
\sphinxstylestrong{Synopsis}
\begin{quote}

\sphinxAtStartPar
\sphinxcode{\sphinxupquote{Var \&GetVar(int i)}}
\end{quote}

\sphinxAtStartPar
\sphinxstylestrong{Arguments}
\begin{quote}

\sphinxAtStartPar
\sphinxcode{\sphinxupquote{i}}: index of the term.
\end{quote}

\sphinxAtStartPar
\sphinxstylestrong{Return}
\begin{quote}

\sphinxAtStartPar
variable of the i\sphinxhyphen{}th term in LMI expression object.
\end{quote}
\end{quote}

\subsubsection{LmiExpr::operator*=()}
\label{\detokenize{cppapi/lmiexpr:lmiexpr-operator}}\begin{quote}

\sphinxAtStartPar
Multiply a double constant to self.

\sphinxAtStartPar
\sphinxstylestrong{Synopsis}
\begin{quote}

\sphinxAtStartPar
\sphinxcode{\sphinxupquote{void operator*=(double c)}}
\end{quote}

\sphinxAtStartPar
\sphinxstylestrong{Arguments}
\begin{quote}

\sphinxAtStartPar
\sphinxcode{\sphinxupquote{c}}: constant multiplier.
\end{quote}
\end{quote}

\subsubsection{LmiExpr::operator*()}
\label{\detokenize{cppapi/lmiexpr:id6}}\begin{quote}

\sphinxAtStartPar
Multiply double constant and return new expression.

\sphinxAtStartPar
\sphinxstylestrong{Synopsis}
\begin{quote}

\sphinxAtStartPar
\sphinxcode{\sphinxupquote{LmiExpr operator*(double c)}}
\end{quote}

\sphinxAtStartPar
\sphinxstylestrong{Arguments}
\begin{quote}

\sphinxAtStartPar
\sphinxcode{\sphinxupquote{c}}: constant multiplier.
\end{quote}

\sphinxAtStartPar
\sphinxstylestrong{Return}
\begin{quote}

\sphinxAtStartPar
result expression.
\end{quote}
\end{quote}

\subsubsection{LmiExpr::operator+=()}
\label{\detokenize{cppapi/lmiexpr:id7}}\begin{quote}

\sphinxAtStartPar
Add a symmetric matrix or LMI expression to self.

\sphinxAtStartPar
\sphinxstylestrong{Synopsis}
\begin{quote}

\sphinxAtStartPar
\sphinxcode{\sphinxupquote{void operator+=(const LmiExpr \&expr)}}
\end{quote}

\sphinxAtStartPar
\sphinxstylestrong{Arguments}
\begin{quote}

\sphinxAtStartPar
\sphinxcode{\sphinxupquote{expr}}: symmetric matrix or LMI expression to be added.
\end{quote}
\end{quote}

\subsubsection{LmiExpr::operator+()}
\label{\detokenize{cppapi/lmiexpr:id8}}\begin{quote}

\sphinxAtStartPar
Add a symmetric matrix or LMI expression and return new LMI expression.

\sphinxAtStartPar
\sphinxstylestrong{Synopsis}
\begin{quote}

\sphinxAtStartPar
\sphinxcode{\sphinxupquote{LmiExpr operator+(const LmiExpr \&other)}}
\end{quote}

\sphinxAtStartPar
\sphinxstylestrong{Arguments}
\begin{quote}

\sphinxAtStartPar
\sphinxcode{\sphinxupquote{other}}: other symmetric matrix or LMI expression to add.
\end{quote}

\sphinxAtStartPar
\sphinxstylestrong{Return}
\begin{quote}

\sphinxAtStartPar
result expression.
\end{quote}
\end{quote}

\subsubsection{LmiExpr::operator\sphinxhyphen{}=()}
\label{\detokenize{cppapi/lmiexpr:id9}}\begin{quote}

\sphinxAtStartPar
Substract a symmetric matrix or LMI expression from self.

\sphinxAtStartPar
\sphinxstylestrong{Synopsis}
\begin{quote}

\sphinxAtStartPar
\sphinxcode{\sphinxupquote{void operator\sphinxhyphen{}=(const LmiExpr \&expr)}}
\end{quote}

\sphinxAtStartPar
\sphinxstylestrong{Arguments}
\begin{quote}

\sphinxAtStartPar
\sphinxcode{\sphinxupquote{expr}}: symmetric matrix or LMI expression to be substracted.
\end{quote}
\end{quote}

\subsubsection{LmiExpr::operator\sphinxhyphen{}()}
\label{\detokenize{cppapi/lmiexpr:id10}}\begin{quote}

\sphinxAtStartPar
Substract a symmetric matrix or LMI expression and return new expression.

\sphinxAtStartPar
\sphinxstylestrong{Synopsis}
\begin{quote}

\sphinxAtStartPar
\sphinxcode{\sphinxupquote{LmiExpr operator\sphinxhyphen{}(const LmiExpr \&other)}}
\end{quote}

\sphinxAtStartPar
\sphinxstylestrong{Arguments}
\begin{quote}

\sphinxAtStartPar
\sphinxcode{\sphinxupquote{other}}: other symmetric matrix or LMI expression to substract.
\end{quote}

\sphinxAtStartPar
\sphinxstylestrong{Return}
\begin{quote}

\sphinxAtStartPar
result expression.
\end{quote}
\end{quote}

\subsubsection{LmiExpr::Remove()}
\label{\detokenize{cppapi/lmiexpr:lmiexpr-remove}}\begin{quote}

\sphinxAtStartPar
Remove i\sphinxhyphen{}th term from LMI expression.

\sphinxAtStartPar
\sphinxstylestrong{Synopsis}
\begin{quote}

\sphinxAtStartPar
\sphinxcode{\sphinxupquote{void Remove(int idx)}}
\end{quote}

\sphinxAtStartPar
\sphinxstylestrong{Arguments}
\begin{quote}

\sphinxAtStartPar
\sphinxcode{\sphinxupquote{idx}}: index of the term to be removed.
\end{quote}
\end{quote}

\subsubsection{LmiExpr::Remove()}
\label{\detokenize{cppapi/lmiexpr:id11}}\begin{quote}

\sphinxAtStartPar
Remove the term associated with given variable from LMI expression.

\sphinxAtStartPar
\sphinxstylestrong{Synopsis}
\begin{quote}

\sphinxAtStartPar
\sphinxcode{\sphinxupquote{void Remove(const Var \&var)}}
\end{quote}

\sphinxAtStartPar
\sphinxstylestrong{Arguments}
\begin{quote}

\sphinxAtStartPar
\sphinxcode{\sphinxupquote{var}}: a variable whose term should be removed.
\end{quote}
\end{quote}

\subsubsection{LmiExpr::Reserve()}
\label{\detokenize{cppapi/lmiexpr:lmiexpr-reserve}}\begin{quote}

\sphinxAtStartPar
Reserve capacity to contain at least n items.

\sphinxAtStartPar
\sphinxstylestrong{Synopsis}
\begin{quote}

\sphinxAtStartPar
\sphinxcode{\sphinxupquote{void Reserve(size\_t n)}}
\end{quote}

\sphinxAtStartPar
\sphinxstylestrong{Arguments}
\begin{quote}

\sphinxAtStartPar
\sphinxcode{\sphinxupquote{n}}: capacity number of LMI expression.
\end{quote}
\end{quote}

\subsubsection{LmiExpr::SetCoeff()}
\label{\detokenize{cppapi/lmiexpr:lmiexpr-setcoeff}}\begin{quote}

\sphinxAtStartPar
Set coefficient matrix of the i\sphinxhyphen{}th term in LMI expression.

\sphinxAtStartPar
\sphinxstylestrong{Synopsis}
\begin{quote}

\sphinxAtStartPar
\sphinxcode{\sphinxupquote{void SetCoeff(int i, const SymMatrix \&mat)}}
\end{quote}

\sphinxAtStartPar
\sphinxstylestrong{Arguments}
\begin{quote}

\sphinxAtStartPar
\sphinxcode{\sphinxupquote{i}}: index of the LMI term.

\sphinxAtStartPar
\sphinxcode{\sphinxupquote{mat}}: new coefficient matrix object.
\end{quote}
\end{quote}

\subsubsection{LmiExpr::SetConstant()}
\label{\detokenize{cppapi/lmiexpr:lmiexpr-setconstant}}\begin{quote}

\sphinxAtStartPar
Set constant term of the LMI expression.

\sphinxAtStartPar
\sphinxstylestrong{Synopsis}
\begin{quote}

\sphinxAtStartPar
\sphinxcode{\sphinxupquote{void SetConstant(const SymMatrix \&mat)}}
\end{quote}

\sphinxAtStartPar
\sphinxstylestrong{Arguments}
\begin{quote}

\sphinxAtStartPar
\sphinxcode{\sphinxupquote{mat}}: new matrix object.
\end{quote}
\end{quote}

\subsubsection{LmiExpr::Size()}
\label{\detokenize{cppapi/lmiexpr:lmiexpr-size}}\begin{quote}

\sphinxAtStartPar
Get number of terms in LMI expression.

\sphinxAtStartPar
\sphinxstylestrong{Synopsis}
\begin{quote}

\sphinxAtStartPar
\sphinxcode{\sphinxupquote{size\_t Size()}}
\end{quote}

\sphinxAtStartPar
\sphinxstylestrong{Return}
\begin{quote}

\sphinxAtStartPar
number of LMI terms.
\end{quote}
\end{quote}

\subsection{SymMatrix}
\label{\detokenize{cppapiref:symmatrix}}\label{\detokenize{cppapiref:chapcppapiref-symmatrix}}
\sphinxAtStartPar
COPT symmetric matrix object. Symmetric matrices are always associated with a
particular model. User creates a symmetric matrix object by adding a symmetric
matrix to model, rather than by constructor of SymMatrix class.

\sphinxAtStartPar
Symmetric matrices are used as coefficient matrices of PSD terms in
PSD expressions, PSD constraints or PSD objectives.

\sphinxstepscope

\subsubsection{SymMatrix::GetDim()}
\label{\detokenize{cppapi/symmatrix:symmatrix-getdim}}\label{\detokenize{cppapi/symmatrix::doc}}\begin{quote}

\sphinxAtStartPar
Get the dimension of a symmetric matrix.

\sphinxAtStartPar
\sphinxstylestrong{Synopsis}
\begin{quote}

\sphinxAtStartPar
\sphinxcode{\sphinxupquote{int GetDim()}}
\end{quote}

\sphinxAtStartPar
\sphinxstylestrong{Return}
\begin{quote}

\sphinxAtStartPar
dimension of a symmetric matrix.
\end{quote}
\end{quote}

\subsubsection{SymMatrix::GetIdx()}
\label{\detokenize{cppapi/symmatrix:symmatrix-getidx}}\begin{quote}

\sphinxAtStartPar
Get the index of a symmetric matrix.

\sphinxAtStartPar
\sphinxstylestrong{Synopsis}
\begin{quote}

\sphinxAtStartPar
\sphinxcode{\sphinxupquote{int GetIdx()}}
\end{quote}

\sphinxAtStartPar
\sphinxstylestrong{Return}
\begin{quote}

\sphinxAtStartPar
index of a symmetric matrix.
\end{quote}
\end{quote}

\subsection{SymMatrixArray}
\label{\detokenize{cppapiref:symmatrixarray}}\label{\detokenize{cppapiref:chapcppapiref-symmatrixarray}}
\sphinxAtStartPar
COPT symmetric matrix object. To store and access a set of
{\hyperref[\detokenize{cppapiref:chapcppapiref-symmatrix}]{\sphinxcrossref{\DUrole{std,std-ref}{SymMatrix}}}} objects, Cardinal Optimizer provides
SymMatrixArray class, which defines the following methods.

\sphinxstepscope

\subsubsection{SymMatrixArray::GetMatrix()}
\label{\detokenize{cppapi/symmatrixarray:symmatrixarray-getmatrix}}\label{\detokenize{cppapi/symmatrixarray::doc}}\begin{quote}

\sphinxAtStartPar
Get i\sphinxhyphen{}th SymMatrix object.

\sphinxAtStartPar
\sphinxstylestrong{Synopsis}
\begin{quote}

\sphinxAtStartPar
\sphinxcode{\sphinxupquote{SymMatrix \&GetMatrix(int idx)}}
\end{quote}

\sphinxAtStartPar
\sphinxstylestrong{Arguments}
\begin{quote}

\sphinxAtStartPar
\sphinxcode{\sphinxupquote{idx}}: index of the SymMatrix object.
\end{quote}

\sphinxAtStartPar
\sphinxstylestrong{Return}
\begin{quote}

\sphinxAtStartPar
SymMatrix object with index idx.
\end{quote}
\end{quote}

\subsubsection{SymMatrixArray::PushBack()}
\label{\detokenize{cppapi/symmatrixarray:symmatrixarray-pushback}}\begin{quote}

\sphinxAtStartPar
Add a SymMatrix object to SymMatrix array.

\sphinxAtStartPar
\sphinxstylestrong{Synopsis}
\begin{quote}

\sphinxAtStartPar
\sphinxcode{\sphinxupquote{void PushBack(const SymMatrix \&mat)}}
\end{quote}

\sphinxAtStartPar
\sphinxstylestrong{Arguments}
\begin{quote}

\sphinxAtStartPar
\sphinxcode{\sphinxupquote{mat}}: a SymMatrix object.
\end{quote}
\end{quote}

\subsubsection{SymMatrixArray::Reserve()}
\label{\detokenize{cppapi/symmatrixarray:symmatrixarray-reserve}}\begin{quote}

\sphinxAtStartPar
Reserve capacity to contain at least n items.

\sphinxAtStartPar
\sphinxstylestrong{Synopsis}
\begin{quote}

\sphinxAtStartPar
\sphinxcode{\sphinxupquote{void Reserve(int n)}}
\end{quote}

\sphinxAtStartPar
\sphinxstylestrong{Arguments}
\begin{quote}

\sphinxAtStartPar
\sphinxcode{\sphinxupquote{n}}: minimum capacity for symmetric matrix object.
\end{quote}
\end{quote}

\subsubsection{SymMatrixArray::Size()}
\label{\detokenize{cppapi/symmatrixarray:symmatrixarray-size}}\begin{quote}

\sphinxAtStartPar
Get the number of SymMatrix objects.

\sphinxAtStartPar
\sphinxstylestrong{Synopsis}
\begin{quote}

\sphinxAtStartPar
\sphinxcode{\sphinxupquote{int Size()}}
\end{quote}

\sphinxAtStartPar
\sphinxstylestrong{Return}
\begin{quote}

\sphinxAtStartPar
number of SymMatrix objects.
\end{quote}
\end{quote}

\subsection{SymMatExpr}
\label{\detokenize{cppapiref:symmatexpr}}\label{\detokenize{cppapiref:chapcppapiref-symmatexpr}}
\sphinxAtStartPar
COPT symmetric matrix expression object. A symmetric matrix expression is a
linear combination of symmetric matrices, which is still a symmetric matrix.
However, by doing so, we are able to delay computing the final matrix
until setting PSD constraints or PSD objective.

\sphinxstepscope

\subsubsection{SymMatExpr::SymMatExpr()}
\label{\detokenize{cppapi/symmatexpr:symmatexpr-symmatexpr}}\label{\detokenize{cppapi/symmatexpr::doc}}\begin{quote}

\sphinxAtStartPar
Constructor of a symmetric matrix expression.

\sphinxAtStartPar
\sphinxstylestrong{Synopsis}
\begin{quote}

\sphinxAtStartPar
\sphinxcode{\sphinxupquote{SymMatExpr()}}
\end{quote}
\end{quote}

\subsubsection{SymMatExpr::SymMatExpr()}
\label{\detokenize{cppapi/symmatexpr:id1}}\begin{quote}

\sphinxAtStartPar
Constructor of a symmetric matrix expression with one term.

\sphinxAtStartPar
\sphinxstylestrong{Synopsis}
\begin{quote}

\sphinxAtStartPar
\sphinxcode{\sphinxupquote{SymMatExpr(const SymMatrix \&mat, double coeff)}}
\end{quote}

\sphinxAtStartPar
\sphinxstylestrong{Arguments}
\begin{quote}

\sphinxAtStartPar
\sphinxcode{\sphinxupquote{mat}}: symmetric matrix of the added term.

\sphinxAtStartPar
\sphinxcode{\sphinxupquote{coeff}}: optional, coefficent for the added term. Its default value is 1.0.
\end{quote}
\end{quote}

\subsubsection{SymMatExpr::AddSymMatExpr()}
\label{\detokenize{cppapi/symmatexpr:symmatexpr-addsymmatexpr}}\begin{quote}

\sphinxAtStartPar
Add a symmetric matrix expression to self.

\sphinxAtStartPar
\sphinxstylestrong{Synopsis}
\begin{quote}

\sphinxAtStartPar
\sphinxcode{\sphinxupquote{void AddSymMatExpr(const SymMatExpr \&expr, double mult)}}
\end{quote}

\sphinxAtStartPar
\sphinxstylestrong{Arguments}
\begin{quote}

\sphinxAtStartPar
\sphinxcode{\sphinxupquote{expr}}: symmetric matrix expression to be added.

\sphinxAtStartPar
\sphinxcode{\sphinxupquote{mult}}: optional, constant multiplier, default value is 1.0.
\end{quote}
\end{quote}

\subsubsection{SymMatExpr::AddTerm()}
\label{\detokenize{cppapi/symmatexpr:symmatexpr-addterm}}\begin{quote}

\sphinxAtStartPar
Add a term to symmetric matrix expression object.

\sphinxAtStartPar
\sphinxstylestrong{Synopsis}
\begin{quote}

\sphinxAtStartPar
\sphinxcode{\sphinxupquote{bool AddTerm(const SymMatrix \&mat, double coeff)}}
\end{quote}

\sphinxAtStartPar
\sphinxstylestrong{Arguments}
\begin{quote}

\sphinxAtStartPar
\sphinxcode{\sphinxupquote{mat}}: symmetric matrix of the new term.

\sphinxAtStartPar
\sphinxcode{\sphinxupquote{coeff}}: coefficient of the new term.
\end{quote}

\sphinxAtStartPar
\sphinxstylestrong{Return}
\begin{quote}

\sphinxAtStartPar
True if the term is added successfully.
\end{quote}
\end{quote}

\subsubsection{SymMatExpr::AddTerms()}
\label{\detokenize{cppapi/symmatexpr:symmatexpr-addterms}}\begin{quote}

\sphinxAtStartPar
Add multiple terms to expression object.

\sphinxAtStartPar
\sphinxstylestrong{Synopsis}
\begin{quote}

\sphinxAtStartPar
\sphinxcode{\sphinxupquote{int AddTerms(}}
\begin{quote}

\sphinxAtStartPar
\sphinxcode{\sphinxupquote{const SymMatrixArray \&mats,}}

\sphinxAtStartPar
\sphinxcode{\sphinxupquote{double *pCoeff,}}

\sphinxAtStartPar
\sphinxcode{\sphinxupquote{int len)}}
\end{quote}
\end{quote}

\sphinxAtStartPar
\sphinxstylestrong{Arguments}
\begin{quote}

\sphinxAtStartPar
\sphinxcode{\sphinxupquote{mats}}: symmetric matrix array object for added terms.

\sphinxAtStartPar
\sphinxcode{\sphinxupquote{pCoeff}}: coefficient array for added terms.

\sphinxAtStartPar
\sphinxcode{\sphinxupquote{len}}: length of coefficient array.
\end{quote}

\sphinxAtStartPar
\sphinxstylestrong{Return}
\begin{quote}

\sphinxAtStartPar
Number of added terms. If negative, fail to add one of terms.
\end{quote}
\end{quote}

\subsubsection{SymMatExpr::Clone()}
\label{\detokenize{cppapi/symmatexpr:symmatexpr-clone}}\begin{quote}

\sphinxAtStartPar
Deep copy symmetric matrix expression object.

\sphinxAtStartPar
\sphinxstylestrong{Synopsis}
\begin{quote}

\sphinxAtStartPar
\sphinxcode{\sphinxupquote{SymMatExpr Clone()}}
\end{quote}

\sphinxAtStartPar
\sphinxstylestrong{Return}
\begin{quote}

\sphinxAtStartPar
cloned expression object.
\end{quote}
\end{quote}

\subsubsection{SymMatExpr::GetCoeff()}
\label{\detokenize{cppapi/symmatexpr:symmatexpr-getcoeff}}\begin{quote}

\sphinxAtStartPar
Get coefficient of the i\sphinxhyphen{}th term in expression object.

\sphinxAtStartPar
\sphinxstylestrong{Synopsis}
\begin{quote}

\sphinxAtStartPar
\sphinxcode{\sphinxupquote{double GetCoeff(int i)}}
\end{quote}

\sphinxAtStartPar
\sphinxstylestrong{Arguments}
\begin{quote}

\sphinxAtStartPar
\sphinxcode{\sphinxupquote{i}}: index of the term.
\end{quote}

\sphinxAtStartPar
\sphinxstylestrong{Return}
\begin{quote}

\sphinxAtStartPar
coefficient of the i\sphinxhyphen{}th term.
\end{quote}
\end{quote}

\subsubsection{SymMatExpr::GetDim()}
\label{\detokenize{cppapi/symmatexpr:symmatexpr-getdim}}\begin{quote}

\sphinxAtStartPar
Get dimension of symmetric matrix in expression.

\sphinxAtStartPar
\sphinxstylestrong{Synopsis}
\begin{quote}

\sphinxAtStartPar
\sphinxcode{\sphinxupquote{int GetDim()}}
\end{quote}

\sphinxAtStartPar
\sphinxstylestrong{Return}
\begin{quote}

\sphinxAtStartPar
dimension of symmetric matrix.
\end{quote}
\end{quote}

\subsubsection{SymMatExpr::GetSymMat()}
\label{\detokenize{cppapi/symmatexpr:symmatexpr-getsymmat}}\begin{quote}

\sphinxAtStartPar
Get symmetric matrix of the i\sphinxhyphen{}th term in expression object.

\sphinxAtStartPar
\sphinxstylestrong{Synopsis}
\begin{quote}

\sphinxAtStartPar
\sphinxcode{\sphinxupquote{SymMatrix \&GetSymMat(int i)}}
\end{quote}

\sphinxAtStartPar
\sphinxstylestrong{Arguments}
\begin{quote}

\sphinxAtStartPar
\sphinxcode{\sphinxupquote{i}}: index of the term.
\end{quote}

\sphinxAtStartPar
\sphinxstylestrong{Return}
\begin{quote}

\sphinxAtStartPar
the symmetric matrix of the i\sphinxhyphen{}th term.
\end{quote}
\end{quote}

\subsubsection{SymMatExpr::operator*=()}
\label{\detokenize{cppapi/symmatexpr:symmatexpr-operator}}\begin{quote}

\sphinxAtStartPar
Multiply a constant to self.

\sphinxAtStartPar
\sphinxstylestrong{Synopsis}
\begin{quote}

\sphinxAtStartPar
\sphinxcode{\sphinxupquote{void operator*=(double c)}}
\end{quote}

\sphinxAtStartPar
\sphinxstylestrong{Arguments}
\begin{quote}

\sphinxAtStartPar
\sphinxcode{\sphinxupquote{c}}: constant multiplier.
\end{quote}
\end{quote}

\subsubsection{SymMatExpr::operator*()}
\label{\detokenize{cppapi/symmatexpr:id2}}\begin{quote}

\sphinxAtStartPar
Multiply constant and return new expression.

\sphinxAtStartPar
\sphinxstylestrong{Synopsis}
\begin{quote}

\sphinxAtStartPar
\sphinxcode{\sphinxupquote{SymMatExpr operator*(double c)}}
\end{quote}

\sphinxAtStartPar
\sphinxstylestrong{Arguments}
\begin{quote}

\sphinxAtStartPar
\sphinxcode{\sphinxupquote{c}}: constant multiplier.
\end{quote}

\sphinxAtStartPar
\sphinxstylestrong{Return}
\begin{quote}

\sphinxAtStartPar
result expression.
\end{quote}
\end{quote}

\subsubsection{SymMatExpr::operator+=()}
\label{\detokenize{cppapi/symmatexpr:id3}}\begin{quote}

\sphinxAtStartPar
Add a symmetric matrix expression to self.

\sphinxAtStartPar
\sphinxstylestrong{Synopsis}
\begin{quote}

\sphinxAtStartPar
\sphinxcode{\sphinxupquote{void operator+=(const SymMatExpr \&expr)}}
\end{quote}

\sphinxAtStartPar
\sphinxstylestrong{Arguments}
\begin{quote}

\sphinxAtStartPar
\sphinxcode{\sphinxupquote{expr}}: symmetric matrix expression to be added.
\end{quote}
\end{quote}

\subsubsection{SymMatExpr::operator+()}
\label{\detokenize{cppapi/symmatexpr:id4}}\begin{quote}

\sphinxAtStartPar
Add expression and return new expression.

\sphinxAtStartPar
\sphinxstylestrong{Synopsis}
\begin{quote}

\sphinxAtStartPar
\sphinxcode{\sphinxupquote{SymMatExpr operator+(const SymMatExpr \&other)}}
\end{quote}

\sphinxAtStartPar
\sphinxstylestrong{Arguments}
\begin{quote}

\sphinxAtStartPar
\sphinxcode{\sphinxupquote{other}}: other expression to add.
\end{quote}

\sphinxAtStartPar
\sphinxstylestrong{Return}
\begin{quote}

\sphinxAtStartPar
result expression.
\end{quote}
\end{quote}

\subsubsection{SymMatExpr::operator\sphinxhyphen{}=()}
\label{\detokenize{cppapi/symmatexpr:id5}}\begin{quote}

\sphinxAtStartPar
Substract a symmetric matrix expression from self.

\sphinxAtStartPar
\sphinxstylestrong{Synopsis}
\begin{quote}

\sphinxAtStartPar
\sphinxcode{\sphinxupquote{void operator\sphinxhyphen{}=(const SymMatExpr \&expr)}}
\end{quote}

\sphinxAtStartPar
\sphinxstylestrong{Arguments}
\begin{quote}

\sphinxAtStartPar
\sphinxcode{\sphinxupquote{expr}}: symmetric matrix to be substracted.
\end{quote}
\end{quote}

\subsubsection{SymMatExpr::operator\sphinxhyphen{}()}
\label{\detokenize{cppapi/symmatexpr:id6}}\begin{quote}

\sphinxAtStartPar
Substract expression and return new expression.

\sphinxAtStartPar
\sphinxstylestrong{Synopsis}
\begin{quote}

\sphinxAtStartPar
\sphinxcode{\sphinxupquote{SymMatExpr operator\sphinxhyphen{}(const SymMatExpr \&other)}}
\end{quote}

\sphinxAtStartPar
\sphinxstylestrong{Arguments}
\begin{quote}

\sphinxAtStartPar
\sphinxcode{\sphinxupquote{other}}: other expression to substract.
\end{quote}

\sphinxAtStartPar
\sphinxstylestrong{Return}
\begin{quote}

\sphinxAtStartPar
result expression.
\end{quote}
\end{quote}

\subsubsection{SymMatExpr::Remove()}
\label{\detokenize{cppapi/symmatexpr:symmatexpr-remove}}\begin{quote}

\sphinxAtStartPar
Remove i\sphinxhyphen{}th term from expression object.

\sphinxAtStartPar
\sphinxstylestrong{Synopsis}
\begin{quote}

\sphinxAtStartPar
\sphinxcode{\sphinxupquote{void Remove(int idx)}}
\end{quote}

\sphinxAtStartPar
\sphinxstylestrong{Arguments}
\begin{quote}

\sphinxAtStartPar
\sphinxcode{\sphinxupquote{idx}}: index of the term to be removed.
\end{quote}
\end{quote}

\subsubsection{SymMatExpr::Remove()}
\label{\detokenize{cppapi/symmatexpr:id7}}\begin{quote}

\sphinxAtStartPar
Remove the term associated with the symmetric matrix.

\sphinxAtStartPar
\sphinxstylestrong{Synopsis}
\begin{quote}

\sphinxAtStartPar
\sphinxcode{\sphinxupquote{void Remove(const SymMatrix \&mat)}}
\end{quote}

\sphinxAtStartPar
\sphinxstylestrong{Arguments}
\begin{quote}

\sphinxAtStartPar
\sphinxcode{\sphinxupquote{mat}}: a symmetric matrix whose term should be removed.
\end{quote}
\end{quote}

\subsubsection{SymMatExpr::Reserve()}
\label{\detokenize{cppapi/symmatexpr:symmatexpr-reserve}}\begin{quote}

\sphinxAtStartPar
Reserve capacity to contain at least n items.

\sphinxAtStartPar
\sphinxstylestrong{Synopsis}
\begin{quote}

\sphinxAtStartPar
\sphinxcode{\sphinxupquote{void Reserve(size\_t n)}}
\end{quote}

\sphinxAtStartPar
\sphinxstylestrong{Arguments}
\begin{quote}

\sphinxAtStartPar
\sphinxcode{\sphinxupquote{n}}: minimum capacity for expression object.
\end{quote}
\end{quote}

\subsubsection{SymMatExpr::SetCoeff()}
\label{\detokenize{cppapi/symmatexpr:symmatexpr-setcoeff}}\begin{quote}

\sphinxAtStartPar
Set coefficient for the i\sphinxhyphen{}th term in expression object.

\sphinxAtStartPar
\sphinxstylestrong{Synopsis}
\begin{quote}

\sphinxAtStartPar
\sphinxcode{\sphinxupquote{void SetCoeff(int i, double val)}}
\end{quote}

\sphinxAtStartPar
\sphinxstylestrong{Arguments}
\begin{quote}

\sphinxAtStartPar
\sphinxcode{\sphinxupquote{i}}: index of the term.

\sphinxAtStartPar
\sphinxcode{\sphinxupquote{val}}: coefficient of the term.
\end{quote}
\end{quote}

\subsubsection{SymMatExpr::Size()}
\label{\detokenize{cppapi/symmatexpr:symmatexpr-size}}\begin{quote}

\sphinxAtStartPar
Get number of terms in expression.

\sphinxAtStartPar
\sphinxstylestrong{Synopsis}
\begin{quote}

\sphinxAtStartPar
\sphinxcode{\sphinxupquote{size\_t Size()}}
\end{quote}

\sphinxAtStartPar
\sphinxstylestrong{Return}
\begin{quote}

\sphinxAtStartPar
number of terms.
\end{quote}
\end{quote}

\subsection{MVar}
\label{\detokenize{cppapiref:mvar}}\label{\detokenize{cppapiref:chapcppapiref-mvar}}
\sphinxAtStartPar
COPT multi\sphinxhyphen{}dimensional variable object. It is used to construct
multi\sphinxhyphen{}dimensional variables and support operations with the built\sphinxhyphen{}in
multi\sphinxhyphen{}dimensional array {\hyperref[\detokenize{cppapiref:chapcppapiref-ndarray}]{\sphinxcrossref{\DUrole{std,std-ref}{NdArray}}}} in COPT. It can be
created by calling the method \sphinxcode{\sphinxupquote{addMVar}} of {\hyperref[\detokenize{cppapiref:chapcppapiref-model}]{\sphinxcrossref{\DUrole{std,std-ref}{Model}}}}.
The following methods are provided:

\sphinxstepscope

\subsubsection{MVar::MVar()}
\label{\detokenize{cppapi/mvar:mvar-mvar}}\label{\detokenize{cppapi/mvar::doc}}\begin{quote}

\sphinxAtStartPar
Construct a MVar object with the given shape, filling with the given variable.

\sphinxAtStartPar
\sphinxstylestrong{Synopsis}
\begin{quote}

\sphinxAtStartPar
\sphinxcode{\sphinxupquote{MVar(const Shape\textless{}N\textgreater{} \&shp, const Var \&var)}}
\end{quote}

\sphinxAtStartPar
\sphinxstylestrong{Arguments}
\begin{quote}

\sphinxAtStartPar
\sphinxcode{\sphinxupquote{shp}}: shape of MVar.

\sphinxAtStartPar
\sphinxcode{\sphinxupquote{var}}: Variable object.
\end{quote}
\end{quote}

\subsubsection{MVar::MVar()}
\label{\detokenize{cppapi/mvar:id1}}\begin{quote}

\sphinxAtStartPar
Construct a MVar object with the given shape, filling with an array of variables.

\sphinxAtStartPar
\sphinxstylestrong{Synopsis}
\begin{quote}

\sphinxAtStartPar
\sphinxcode{\sphinxupquote{MVar(const Shape\textless{}N\textgreater{} \&shp, const VarArray \&vars)}}
\end{quote}

\sphinxAtStartPar
\sphinxstylestrong{Arguments}
\begin{quote}

\sphinxAtStartPar
\sphinxcode{\sphinxupquote{shp}}: shape of MVar.

\sphinxAtStartPar
\sphinxcode{\sphinxupquote{vars}}: an array of variables.
\end{quote}
\end{quote}

\subsubsection{MVar::Clone()}
\label{\detokenize{cppapi/mvar:mvar-clone}}\begin{quote}

\sphinxAtStartPar
Clone MVar object.

\sphinxAtStartPar
\sphinxstylestrong{Synopsis}
\begin{quote}

\sphinxAtStartPar
\sphinxcode{\sphinxupquote{MVar Clone()}}
\end{quote}

\sphinxAtStartPar
\sphinxstylestrong{Return}
\begin{quote}

\sphinxAtStartPar
new MVar object.
\end{quote}
\end{quote}

\subsubsection{MVar::Diagonal()}
\label{\detokenize{cppapi/mvar:mvar-diagonal}}\begin{quote}

\sphinxAtStartPar
Get diagonals of MVar object.

\sphinxAtStartPar
\sphinxstylestrong{Synopsis}
\begin{quote}

\sphinxAtStartPar
\sphinxcode{\sphinxupquote{MVar\textless{}N \sphinxhyphen{} 1\textgreater{} Diagonal(}}
\begin{quote}

\sphinxAtStartPar
\sphinxcode{\sphinxupquote{int offset,}}

\sphinxAtStartPar
\sphinxcode{\sphinxupquote{int axis1,}}

\sphinxAtStartPar
\sphinxcode{\sphinxupquote{int axis2)}}
\end{quote}
\end{quote}

\sphinxAtStartPar
\sphinxstylestrong{Arguments}
\begin{quote}

\sphinxAtStartPar
\sphinxcode{\sphinxupquote{offset}}: offset of the diagonal from the main diagonal. Can be positive or negative.

\sphinxAtStartPar
\sphinxcode{\sphinxupquote{axis1}}: 1st axis of MVar.

\sphinxAtStartPar
\sphinxcode{\sphinxupquote{axis2}}: 2nd axis of MVar.
\end{quote}

\sphinxAtStartPar
\sphinxstylestrong{Return}
\begin{quote}

\sphinxAtStartPar
(N\sphinxhyphen{}1)\sphinxhyphen{}dimensional diagonals.
\end{quote}
\end{quote}

\subsubsection{MVar::Expand()}
\label{\detokenize{cppapi/mvar:mvar-expand}}\begin{quote}

\sphinxAtStartPar
Expand shape of MVar object.

\sphinxAtStartPar
\sphinxstylestrong{Synopsis}
\begin{quote}

\sphinxAtStartPar
\sphinxcode{\sphinxupquote{MVar\textless{}N + 1\textgreater{} Expand(int axis)}}
\end{quote}

\sphinxAtStartPar
\sphinxstylestrong{Arguments}
\begin{quote}

\sphinxAtStartPar
\sphinxcode{\sphinxupquote{axis}}: axis of MVar.
\end{quote}

\sphinxAtStartPar
\sphinxstylestrong{Return}
\begin{quote}

\sphinxAtStartPar
MVar object of (N+1)\sphinxhyphen{}dimensional shape.
\end{quote}
\end{quote}

\subsubsection{MVar::Flatten()}
\label{\detokenize{cppapi/mvar:mvar-flatten}}\begin{quote}

\sphinxAtStartPar
Flatten a MVar object to a 1\sphinxhyphen{}dimensional shape.

\sphinxAtStartPar
\sphinxstylestrong{Synopsis}
\begin{quote}

\sphinxAtStartPar
\sphinxcode{\sphinxupquote{MVar\textless{}1\textgreater{} Flatten()}}
\end{quote}

\sphinxAtStartPar
\sphinxstylestrong{Return}
\begin{quote}

\sphinxAtStartPar
a MVar object collapsed into one dimension.
\end{quote}
\end{quote}

\subsubsection{MVar::Get()}
\label{\detokenize{cppapi/mvar:mvar-get}}\begin{quote}

\sphinxAtStartPar
Get values of information associated with variables in MVar object.

\sphinxAtStartPar
\sphinxstylestrong{Synopsis}
\begin{quote}

\sphinxAtStartPar
\sphinxcode{\sphinxupquote{NdArray\textless{}double, N\textgreater{} Get(const char *szInfo)}}
\end{quote}

\sphinxAtStartPar
\sphinxstylestrong{Arguments}
\begin{quote}

\sphinxAtStartPar
\sphinxcode{\sphinxupquote{szInfo}}: name of information.
\end{quote}

\sphinxAtStartPar
\sphinxstylestrong{Return}
\begin{quote}

\sphinxAtStartPar
multi\sphinxhyphen{}dimensional array of information of variables.
\end{quote}
\end{quote}

\subsubsection{MVar::GetBasis()}
\label{\detokenize{cppapi/mvar:mvar-getbasis}}\begin{quote}

\sphinxAtStartPar
Get basis of variables in MVar object.

\sphinxAtStartPar
\sphinxstylestrong{Synopsis}
\begin{quote}

\sphinxAtStartPar
\sphinxcode{\sphinxupquote{NdArray\textless{}int, N\textgreater{} GetBasis()}}
\end{quote}

\sphinxAtStartPar
\sphinxstylestrong{Return}
\begin{quote}

\sphinxAtStartPar
multi\sphinxhyphen{}dimensional array of basis of variables.
\end{quote}
\end{quote}

\subsubsection{MVar::GetDim()}
\label{\detokenize{cppapi/mvar:mvar-getdim}}\begin{quote}

\sphinxAtStartPar
Get i\sphinxhyphen{}th dimension of MVar object.

\sphinxAtStartPar
\sphinxstylestrong{Synopsis}
\begin{quote}

\sphinxAtStartPar
\sphinxcode{\sphinxupquote{size\_t GetDim(int i)}}
\end{quote}

\sphinxAtStartPar
\sphinxstylestrong{Arguments}
\begin{quote}

\sphinxAtStartPar
\sphinxcode{\sphinxupquote{i}}: index of dimension
\end{quote}

\sphinxAtStartPar
\sphinxstylestrong{Return}
\begin{quote}

\sphinxAtStartPar
i\sphinxhyphen{}th dimension.
\end{quote}
\end{quote}

\subsubsection{MVar::GetIdx()}
\label{\detokenize{cppapi/mvar:mvar-getidx}}\begin{quote}

\sphinxAtStartPar
Get indexes of variables in MVar object.

\sphinxAtStartPar
\sphinxstylestrong{Synopsis}
\begin{quote}

\sphinxAtStartPar
\sphinxcode{\sphinxupquote{NdArray\textless{}int, N\textgreater{} GetIdx()}}
\end{quote}

\sphinxAtStartPar
\sphinxstylestrong{Return}
\begin{quote}

\sphinxAtStartPar
multi\sphinxhyphen{}dimensional array of indexes of variables.
\end{quote}
\end{quote}

\subsubsection{MVar::GetLowerIIS()}
\label{\detokenize{cppapi/mvar:mvar-getloweriis}}\begin{quote}

\sphinxAtStartPar
Get IIS status of lower bound of variables in MVar object.

\sphinxAtStartPar
\sphinxstylestrong{Synopsis}
\begin{quote}

\sphinxAtStartPar
\sphinxcode{\sphinxupquote{NdArray\textless{}int, N\textgreater{} GetLowerIIS()}}
\end{quote}

\sphinxAtStartPar
\sphinxstylestrong{Return}
\begin{quote}

\sphinxAtStartPar
multi\sphinxhyphen{}dimensional array of IIS status of lower bounds of variables.
\end{quote}
\end{quote}

\subsubsection{MVar::GetND()}
\label{\detokenize{cppapi/mvar:mvar-getnd}}\begin{quote}

\sphinxAtStartPar
Get number of dimensions of MVar object.

\sphinxAtStartPar
\sphinxstylestrong{Synopsis}
\begin{quote}

\sphinxAtStartPar
\sphinxcode{\sphinxupquote{int GetND()}}
\end{quote}

\sphinxAtStartPar
\sphinxstylestrong{Return}
\begin{quote}

\sphinxAtStartPar
number of dimensions.
\end{quote}
\end{quote}

\subsubsection{MVar::GetShape()}
\label{\detokenize{cppapi/mvar:mvar-getshape}}\begin{quote}

\sphinxAtStartPar
Get shape of MVar object.

\sphinxAtStartPar
\sphinxstylestrong{Synopsis}
\begin{quote}

\sphinxAtStartPar
\sphinxcode{\sphinxupquote{Shape\textless{}N\textgreater{} GetShape()}}
\end{quote}

\sphinxAtStartPar
\sphinxstylestrong{Return}
\begin{quote}

\sphinxAtStartPar
shape object.
\end{quote}
\end{quote}

\subsubsection{MVar::GetSize()}
\label{\detokenize{cppapi/mvar:mvar-getsize}}\begin{quote}

\sphinxAtStartPar
Get size of MVar object.

\sphinxAtStartPar
\sphinxstylestrong{Synopsis}
\begin{quote}

\sphinxAtStartPar
\sphinxcode{\sphinxupquote{size\_t GetSize()}}
\end{quote}

\sphinxAtStartPar
\sphinxstylestrong{Return}
\begin{quote}

\sphinxAtStartPar
number of vars.
\end{quote}
\end{quote}

\subsubsection{MVar::GetType()}
\label{\detokenize{cppapi/mvar:mvar-gettype}}\begin{quote}

\sphinxAtStartPar
Get types of variables in MVar object.

\sphinxAtStartPar
\sphinxstylestrong{Synopsis}
\begin{quote}

\sphinxAtStartPar
\sphinxcode{\sphinxupquote{NdArray\textless{}char, N\textgreater{} GetType()}}
\end{quote}

\sphinxAtStartPar
\sphinxstylestrong{Return}
\begin{quote}

\sphinxAtStartPar
multi\sphinxhyphen{}dimensional array of types of variables.
\end{quote}
\end{quote}

\subsubsection{MVar::GetUpperIIS()}
\label{\detokenize{cppapi/mvar:mvar-getupperiis}}\begin{quote}

\sphinxAtStartPar
Get IIS status of upper bound of variables in MVar object.

\sphinxAtStartPar
\sphinxstylestrong{Synopsis}
\begin{quote}

\sphinxAtStartPar
\sphinxcode{\sphinxupquote{NdArray\textless{}int, N\textgreater{} GetUpperIIS()}}
\end{quote}

\sphinxAtStartPar
\sphinxstylestrong{Return}
\begin{quote}

\sphinxAtStartPar
multi\sphinxhyphen{}dimensional array of IIS status of upper bounds of variables.
\end{quote}
\end{quote}

\subsubsection{MVar::Item()}
\label{\detokenize{cppapi/mvar:mvar-item}}\begin{quote}

\sphinxAtStartPar
Get variable of given index from MVar object.

\sphinxAtStartPar
\sphinxstylestrong{Synopsis}
\begin{quote}

\sphinxAtStartPar
\sphinxcode{\sphinxupquote{Var \&Item(size\_t idx)}}
\end{quote}

\sphinxAtStartPar
\sphinxstylestrong{Arguments}
\begin{quote}

\sphinxAtStartPar
\sphinxcode{\sphinxupquote{idx}}: index of var.
\end{quote}

\sphinxAtStartPar
\sphinxstylestrong{Return}
\begin{quote}

\sphinxAtStartPar
Var object.
\end{quote}
\end{quote}

\subsubsection{MVar::operator{[}{]}()}
\label{\detokenize{cppapi/mvar:mvar-operator}}\begin{quote}

\sphinxAtStartPar
Get variable of given index from MVar object.

\sphinxAtStartPar
\sphinxstylestrong{Synopsis}
\begin{quote}

\sphinxAtStartPar
\sphinxcode{\sphinxupquote{Var \&operator{[}{]}(size\_t idx)}}
\end{quote}

\sphinxAtStartPar
\sphinxstylestrong{Arguments}
\begin{quote}

\sphinxAtStartPar
\sphinxcode{\sphinxupquote{idx}}: index of var.
\end{quote}

\sphinxAtStartPar
\sphinxstylestrong{Return}
\begin{quote}

\sphinxAtStartPar
Var object.
\end{quote}
\end{quote}

\subsubsection{MVar::operator{[}{]}()}
\label{\detokenize{cppapi/mvar:id2}}\begin{quote}

\sphinxAtStartPar
Get sub\sphinxhyphen{}arrays of MVar object, given view object.

\sphinxAtStartPar
\sphinxstylestrong{Synopsis}
\begin{quote}

\sphinxAtStartPar
\sphinxcode{\sphinxupquote{MVar operator{[}{]}(const View \&view)}}
\end{quote}

\sphinxAtStartPar
\sphinxstylestrong{Arguments}
\begin{quote}

\sphinxAtStartPar
\sphinxcode{\sphinxupquote{view}}: view of multi\sphinxhyphen{}dimensional array.
\end{quote}

\sphinxAtStartPar
\sphinxstylestrong{Return}
\begin{quote}

\sphinxAtStartPar
sub\sphinxhyphen{}arrays of MVar object.
\end{quote}
\end{quote}

\subsubsection{MVar::Pick()}
\label{\detokenize{cppapi/mvar:mvar-pick}}\begin{quote}

\sphinxAtStartPar
Given a list of indexes, get variables from MVar object.

\sphinxAtStartPar
\sphinxstylestrong{Synopsis}
\begin{quote}

\sphinxAtStartPar
\sphinxcode{\sphinxupquote{MVar\textless{}1\textgreater{} Pick(const NdArray\textless{}int, 1\textgreater{} \&indexes)}}
\end{quote}

\sphinxAtStartPar
\sphinxstylestrong{Arguments}
\begin{quote}

\sphinxAtStartPar
\sphinxcode{\sphinxupquote{indexes}}: indexes of elements.
\end{quote}

\sphinxAtStartPar
\sphinxstylestrong{Return}
\begin{quote}

\sphinxAtStartPar
one\sphinxhyphen{}dimensional array of desired variables.
\end{quote}
\end{quote}

\subsubsection{MVar::Pick()}
\label{\detokenize{cppapi/mvar:id3}}\begin{quote}

\sphinxAtStartPar
Given a list of indexes, get variables from MVar object.

\sphinxAtStartPar
\sphinxstylestrong{Synopsis}
\begin{quote}

\sphinxAtStartPar
\sphinxcode{\sphinxupquote{MVar\textless{}1\textgreater{} Pick(const NdArray\textless{}int, 2\textgreater{} \&idxrows)}}
\end{quote}

\sphinxAtStartPar
\sphinxstylestrong{Arguments}
\begin{quote}

\sphinxAtStartPar
\sphinxcode{\sphinxupquote{idxrows}}: indexes in format of 2\sphinxhyphen{}dimensional array, where each row is position of element.
\end{quote}

\sphinxAtStartPar
\sphinxstylestrong{Return}
\begin{quote}

\sphinxAtStartPar
one\sphinxhyphen{}dimensional array of desired variables.
\end{quote}
\end{quote}

\subsubsection{MVar::Repeat()}
\label{\detokenize{cppapi/mvar:mvar-repeat}}\begin{quote}

\sphinxAtStartPar
Repeat each element of MVar along given axis.

\sphinxAtStartPar
\sphinxstylestrong{Synopsis}
\begin{quote}

\sphinxAtStartPar
\sphinxcode{\sphinxupquote{MVar\textless{}N\textgreater{} Repeat(size\_t repeats, int axis)}}
\end{quote}

\sphinxAtStartPar
\sphinxstylestrong{Arguments}
\begin{quote}

\sphinxAtStartPar
\sphinxcode{\sphinxupquote{repeats}}: number of repetitions for each element.

\sphinxAtStartPar
\sphinxcode{\sphinxupquote{axis}}: axis of MVar.
\end{quote}

\sphinxAtStartPar
\sphinxstylestrong{Return}
\begin{quote}

\sphinxAtStartPar
new MVar object.
\end{quote}
\end{quote}

\subsubsection{MVar::RepeatBlock()}
\label{\detokenize{cppapi/mvar:mvar-repeatblock}}\begin{quote}

\sphinxAtStartPar
Repeat an MVar a number of times along given axis.

\sphinxAtStartPar
\sphinxstylestrong{Synopsis}
\begin{quote}

\sphinxAtStartPar
\sphinxcode{\sphinxupquote{MVar\textless{}N\textgreater{} RepeatBlock(size\_t repeats, int axis)}}
\end{quote}

\sphinxAtStartPar
\sphinxstylestrong{Arguments}
\begin{quote}

\sphinxAtStartPar
\sphinxcode{\sphinxupquote{repeats}}: number of repetitions.

\sphinxAtStartPar
\sphinxcode{\sphinxupquote{axis}}: axis of MVar.
\end{quote}

\sphinxAtStartPar
\sphinxstylestrong{Return}
\begin{quote}

\sphinxAtStartPar
new MVar object.
\end{quote}
\end{quote}

\subsubsection{MVar::Represent()}
\label{\detokenize{cppapi/mvar:mvar-represent}}\begin{quote}

\sphinxAtStartPar
String representation of MVar object.

\sphinxAtStartPar
\sphinxstylestrong{Synopsis}
\begin{quote}

\sphinxAtStartPar
\sphinxcode{\sphinxupquote{std::string Represent(size\_t maxlen)}}
\end{quote}

\sphinxAtStartPar
\sphinxstylestrong{Arguments}
\begin{quote}

\sphinxAtStartPar
\sphinxcode{\sphinxupquote{maxlen}}: max length of representation.
\end{quote}

\sphinxAtStartPar
\sphinxstylestrong{Return}
\begin{quote}

\sphinxAtStartPar
string object.
\end{quote}
\end{quote}

\subsubsection{MVar::Reshape()}
\label{\detokenize{cppapi/mvar:mvar-reshape}}\begin{quote}

\sphinxAtStartPar
Reshape MVar object to new shape.

\sphinxAtStartPar
\sphinxstylestrong{Synopsis}
\begin{quote}

\sphinxAtStartPar
\sphinxcode{\sphinxupquote{template \textless{}int M\textgreater{} MVar\textless{}M\textgreater{} Reshape(const Shape\textless{}M\textgreater{} \&shape)}}
\end{quote}

\sphinxAtStartPar
\sphinxstylestrong{Arguments}
\begin{quote}

\sphinxAtStartPar
\sphinxcode{\sphinxupquote{shape}}: new shape of M\sphinxhyphen{}dimensions.
\end{quote}

\sphinxAtStartPar
\sphinxstylestrong{Return}
\begin{quote}

\sphinxAtStartPar
M\sphinxhyphen{}dimensional MVar object.
\end{quote}
\end{quote}

\subsubsection{MVar::Set()}
\label{\detokenize{cppapi/mvar:mvar-set}}\begin{quote}

\sphinxAtStartPar
Set values of information associated with variables in MVar object.

\sphinxAtStartPar
\sphinxstylestrong{Synopsis}
\begin{quote}

\sphinxAtStartPar
\sphinxcode{\sphinxupquote{void Set(const char *szInfo, double val)}}
\end{quote}

\sphinxAtStartPar
\sphinxstylestrong{Arguments}
\begin{quote}

\sphinxAtStartPar
\sphinxcode{\sphinxupquote{szInfo}}: name of information.

\sphinxAtStartPar
\sphinxcode{\sphinxupquote{val}}: value of information.
\end{quote}
\end{quote}

\subsubsection{MVar::Set()}
\label{\detokenize{cppapi/mvar:id4}}\begin{quote}

\sphinxAtStartPar
Set values of information associated with variables in MVar object.

\sphinxAtStartPar
\sphinxstylestrong{Synopsis}
\begin{quote}

\sphinxAtStartPar
\sphinxcode{\sphinxupquote{void Set(const char *szInfo, const NdArray\textless{}double, N\textgreater{} \&vals)}}
\end{quote}

\sphinxAtStartPar
\sphinxstylestrong{Arguments}
\begin{quote}

\sphinxAtStartPar
\sphinxcode{\sphinxupquote{szInfo}}: name of information.

\sphinxAtStartPar
\sphinxcode{\sphinxupquote{vals}}: multi\sphinxhyphen{}dimensional array of values of information.
\end{quote}
\end{quote}

\subsubsection{MVar::SetItem()}
\label{\detokenize{cppapi/mvar:mvar-setitem}}\begin{quote}

\sphinxAtStartPar
Set variable of given index to MVar object.

\sphinxAtStartPar
\sphinxstylestrong{Synopsis}
\begin{quote}

\sphinxAtStartPar
\sphinxcode{\sphinxupquote{void SetItem(size\_t idx, const Var \&var)}}
\end{quote}

\sphinxAtStartPar
\sphinxstylestrong{Arguments}
\begin{quote}

\sphinxAtStartPar
\sphinxcode{\sphinxupquote{idx}}: index of element.

\sphinxAtStartPar
\sphinxcode{\sphinxupquote{var}}: Var object.
\end{quote}
\end{quote}

\subsubsection{MVar::Squeeze()}
\label{\detokenize{cppapi/mvar:mvar-squeeze}}\begin{quote}

\sphinxAtStartPar
Remove axis of length 1 from shape of MVar object.

\sphinxAtStartPar
\sphinxstylestrong{Synopsis}
\begin{quote}

\sphinxAtStartPar
\sphinxcode{\sphinxupquote{MVar\textless{}N \sphinxhyphen{} 1\textgreater{} Squeeze(int axis)}}
\end{quote}

\sphinxAtStartPar
\sphinxstylestrong{Arguments}
\begin{quote}

\sphinxAtStartPar
\sphinxcode{\sphinxupquote{axis}}: axis of MVar, where the length is 1.
\end{quote}

\sphinxAtStartPar
\sphinxstylestrong{Return}
\begin{quote}

\sphinxAtStartPar
MVar object of (N\sphinxhyphen{}1)\sphinxhyphen{}dimensional shape.
\end{quote}
\end{quote}

\subsubsection{MVar::Stack()}
\label{\detokenize{cppapi/mvar:mvar-stack}}\begin{quote}

\sphinxAtStartPar
Stack with other MVar object along given axis.

\sphinxAtStartPar
\sphinxstylestrong{Synopsis}
\begin{quote}

\sphinxAtStartPar
\sphinxcode{\sphinxupquote{MVar\textless{}N\textgreater{} Stack(const MVar\textless{}N\textgreater{} \&other, int axis)}}
\end{quote}

\sphinxAtStartPar
\sphinxstylestrong{Arguments}
\begin{quote}

\sphinxAtStartPar
\sphinxcode{\sphinxupquote{other}}: a MVar object.

\sphinxAtStartPar
\sphinxcode{\sphinxupquote{axis}}: an axis of MVar.
\end{quote}

\sphinxAtStartPar
\sphinxstylestrong{Return}
\begin{quote}

\sphinxAtStartPar
the result MVar object.
\end{quote}
\end{quote}

\subsubsection{MVar::Sum()}
\label{\detokenize{cppapi/mvar:mvar-sum}}\begin{quote}

\sphinxAtStartPar
Sum of all variables in MVar object.

\sphinxAtStartPar
\sphinxstylestrong{Synopsis}
\begin{quote}

\sphinxAtStartPar
\sphinxcode{\sphinxupquote{MLinExpr\textless{}0\textgreater{} Sum()}}
\end{quote}

\sphinxAtStartPar
\sphinxstylestrong{Return}
\begin{quote}

\sphinxAtStartPar
sum in zero dimension.
\end{quote}
\end{quote}

\subsubsection{MVar::Sum()}
\label{\detokenize{cppapi/mvar:id5}}\begin{quote}

\sphinxAtStartPar
Sum of variables at given axis of MVar object.

\sphinxAtStartPar
\sphinxstylestrong{Synopsis}
\begin{quote}

\sphinxAtStartPar
\sphinxcode{\sphinxupquote{MLinExpr\textless{}N \sphinxhyphen{} 1\textgreater{} Sum(int axis)}}
\end{quote}

\sphinxAtStartPar
\sphinxstylestrong{Arguments}
\begin{quote}

\sphinxAtStartPar
\sphinxcode{\sphinxupquote{axis}}: axis of MVar.
\end{quote}

\sphinxAtStartPar
\sphinxstylestrong{Return}
\begin{quote}

\sphinxAtStartPar
MLinExpr object in (N\sphinxhyphen{}1)\sphinxhyphen{}dimension.
\end{quote}
\end{quote}

\subsubsection{MVar::Transpose()}
\label{\detokenize{cppapi/mvar:mvar-transpose}}\begin{quote}

\sphinxAtStartPar
Perform matrix transpose of MVar object.

\sphinxAtStartPar
\sphinxstylestrong{Synopsis}
\begin{quote}

\sphinxAtStartPar
\sphinxcode{\sphinxupquote{MVar\textless{}N\textgreater{} Transpose()}}
\end{quote}

\sphinxAtStartPar
\sphinxstylestrong{Return}
\begin{quote}

\sphinxAtStartPar
transposed MVar object.
\end{quote}
\end{quote}

\subsection{MConstr}
\label{\detokenize{cppapiref:mconstr}}\label{\detokenize{cppapiref:chapcppapiref-mconstr}}
\sphinxAtStartPar
COPT multi\sphinxhyphen{}dimensional linear constraint object. It can be created by calling
the method \sphinxcode{\sphinxupquote{addMConstr}} of {\hyperref[\detokenize{cppapiref:chapcppapiref-model}]{\sphinxcrossref{\DUrole{std,std-ref}{Model}}}}. The following methods
are provided:

\sphinxstepscope

\subsubsection{MConstr::MConstr()}
\label{\detokenize{cppapi/mconstr:mconstr-mconstr}}\label{\detokenize{cppapi/mconstr::doc}}\begin{quote}

\sphinxAtStartPar
Construct a MConstr object with the given shape, filling with the given constraint.

\sphinxAtStartPar
\sphinxstylestrong{Synopsis}
\begin{quote}

\sphinxAtStartPar
\sphinxcode{\sphinxupquote{MConstr(const Shape\textless{}N\textgreater{} \&shp, const Constraint \&con)}}
\end{quote}

\sphinxAtStartPar
\sphinxstylestrong{Arguments}
\begin{quote}

\sphinxAtStartPar
\sphinxcode{\sphinxupquote{shp}}: shape of MConstr.

\sphinxAtStartPar
\sphinxcode{\sphinxupquote{con}}: Constraint object.
\end{quote}
\end{quote}

\subsubsection{MConstr::MConstr()}
\label{\detokenize{cppapi/mconstr:id1}}\begin{quote}

\sphinxAtStartPar
Construct a MConstr object with the given shape, filling with an array of constraints.

\sphinxAtStartPar
\sphinxstylestrong{Synopsis}
\begin{quote}

\sphinxAtStartPar
\sphinxcode{\sphinxupquote{MConstr(const Shape\textless{}N\textgreater{} \&shp, const ConstrArray \&cons)}}
\end{quote}

\sphinxAtStartPar
\sphinxstylestrong{Arguments}
\begin{quote}

\sphinxAtStartPar
\sphinxcode{\sphinxupquote{shp}}: shape of MConstr.

\sphinxAtStartPar
\sphinxcode{\sphinxupquote{cons}}: an array of constraints.
\end{quote}
\end{quote}

\subsubsection{MConstr::Clone()}
\label{\detokenize{cppapi/mconstr:mconstr-clone}}\begin{quote}

\sphinxAtStartPar
Clone MConstr object.

\sphinxAtStartPar
\sphinxstylestrong{Synopsis}
\begin{quote}

\sphinxAtStartPar
\sphinxcode{\sphinxupquote{MConstr Clone()}}
\end{quote}

\sphinxAtStartPar
\sphinxstylestrong{Return}
\begin{quote}

\sphinxAtStartPar
new MConstr object.
\end{quote}
\end{quote}

\subsubsection{MConstr::Diagonal()}
\label{\detokenize{cppapi/mconstr:mconstr-diagonal}}\begin{quote}

\sphinxAtStartPar
Get diagonals of MConstr object.

\sphinxAtStartPar
\sphinxstylestrong{Synopsis}
\begin{quote}

\sphinxAtStartPar
\sphinxcode{\sphinxupquote{MConstr\textless{}N \sphinxhyphen{} 1\textgreater{} Diagonal(}}
\begin{quote}

\sphinxAtStartPar
\sphinxcode{\sphinxupquote{int offset,}}

\sphinxAtStartPar
\sphinxcode{\sphinxupquote{int axis1,}}

\sphinxAtStartPar
\sphinxcode{\sphinxupquote{int axis2)}}
\end{quote}
\end{quote}

\sphinxAtStartPar
\sphinxstylestrong{Arguments}
\begin{quote}

\sphinxAtStartPar
\sphinxcode{\sphinxupquote{offset}}: offset of the diagonal from the main diagonal. Can be positive or negative.

\sphinxAtStartPar
\sphinxcode{\sphinxupquote{axis1}}: 1st axis of MConstr.

\sphinxAtStartPar
\sphinxcode{\sphinxupquote{axis2}}: 2nd axis of MConstr.
\end{quote}

\sphinxAtStartPar
\sphinxstylestrong{Return}
\begin{quote}

\sphinxAtStartPar
(N\sphinxhyphen{}1)\sphinxhyphen{}dimensional diagonals.
\end{quote}
\end{quote}

\subsubsection{MConstr::Expand()}
\label{\detokenize{cppapi/mconstr:mconstr-expand}}\begin{quote}

\sphinxAtStartPar
Expand shape of MConstr object.

\sphinxAtStartPar
\sphinxstylestrong{Synopsis}
\begin{quote}

\sphinxAtStartPar
\sphinxcode{\sphinxupquote{MConstr\textless{}N + 1\textgreater{} Expand(int axis)}}
\end{quote}

\sphinxAtStartPar
\sphinxstylestrong{Arguments}
\begin{quote}

\sphinxAtStartPar
\sphinxcode{\sphinxupquote{axis}}: axis of MConstr.
\end{quote}

\sphinxAtStartPar
\sphinxstylestrong{Return}
\begin{quote}

\sphinxAtStartPar
MConstr object of (N+1)\sphinxhyphen{}dimensional shape.
\end{quote}
\end{quote}

\subsubsection{MConstr::Flatten()}
\label{\detokenize{cppapi/mconstr:mconstr-flatten}}\begin{quote}

\sphinxAtStartPar
Flatten a MConstr object to a 1\sphinxhyphen{}dimensional shape.

\sphinxAtStartPar
\sphinxstylestrong{Synopsis}
\begin{quote}

\sphinxAtStartPar
\sphinxcode{\sphinxupquote{MConstr\textless{}1\textgreater{} Flatten()}}
\end{quote}

\sphinxAtStartPar
\sphinxstylestrong{Return}
\begin{quote}

\sphinxAtStartPar
a MConstr object collapsed into one dimension.
\end{quote}
\end{quote}

\subsubsection{MConstr::Get()}
\label{\detokenize{cppapi/mconstr:mconstr-get}}\begin{quote}

\sphinxAtStartPar
Get values of information associated with constraints in MConstr object.

\sphinxAtStartPar
\sphinxstylestrong{Synopsis}
\begin{quote}

\sphinxAtStartPar
\sphinxcode{\sphinxupquote{NdArray\textless{}double, N\textgreater{} Get(const char *szInfo)}}
\end{quote}

\sphinxAtStartPar
\sphinxstylestrong{Arguments}
\begin{quote}

\sphinxAtStartPar
\sphinxcode{\sphinxupquote{szInfo}}: name of information.
\end{quote}

\sphinxAtStartPar
\sphinxstylestrong{Return}
\begin{quote}

\sphinxAtStartPar
multi\sphinxhyphen{}dimensional array of information of constraints.
\end{quote}
\end{quote}

\subsubsection{MConstr::GetBasis()}
\label{\detokenize{cppapi/mconstr:mconstr-getbasis}}\begin{quote}

\sphinxAtStartPar
Get basis of constraints in MConstr object.

\sphinxAtStartPar
\sphinxstylestrong{Synopsis}
\begin{quote}

\sphinxAtStartPar
\sphinxcode{\sphinxupquote{NdArray\textless{}int, N\textgreater{} GetBasis()}}
\end{quote}

\sphinxAtStartPar
\sphinxstylestrong{Return}
\begin{quote}

\sphinxAtStartPar
multi\sphinxhyphen{}dimensional array of basis of constraints.
\end{quote}
\end{quote}

\subsubsection{MConstr::GetDim()}
\label{\detokenize{cppapi/mconstr:mconstr-getdim}}\begin{quote}

\sphinxAtStartPar
Get i\sphinxhyphen{}th dimension of MConstr object.

\sphinxAtStartPar
\sphinxstylestrong{Synopsis}
\begin{quote}

\sphinxAtStartPar
\sphinxcode{\sphinxupquote{size\_t GetDim(int i)}}
\end{quote}

\sphinxAtStartPar
\sphinxstylestrong{Arguments}
\begin{quote}

\sphinxAtStartPar
\sphinxcode{\sphinxupquote{i}}: index of dimension
\end{quote}

\sphinxAtStartPar
\sphinxstylestrong{Return}
\begin{quote}

\sphinxAtStartPar
i\sphinxhyphen{}th dimension.
\end{quote}
\end{quote}

\subsubsection{MConstr::GetIdx()}
\label{\detokenize{cppapi/mconstr:mconstr-getidx}}\begin{quote}

\sphinxAtStartPar
Get index of constraints in MConstr object.

\sphinxAtStartPar
\sphinxstylestrong{Synopsis}
\begin{quote}

\sphinxAtStartPar
\sphinxcode{\sphinxupquote{NdArray\textless{}int, N\textgreater{} GetIdx()}}
\end{quote}

\sphinxAtStartPar
\sphinxstylestrong{Return}
\begin{quote}

\sphinxAtStartPar
multi\sphinxhyphen{}dimensional array of indexes of constraints.
\end{quote}
\end{quote}

\subsubsection{MConstr::GetLowerIIS()}
\label{\detokenize{cppapi/mconstr:mconstr-getloweriis}}\begin{quote}

\sphinxAtStartPar
Get IIS status of lower bound of constraints in MConstr object.

\sphinxAtStartPar
\sphinxstylestrong{Synopsis}
\begin{quote}

\sphinxAtStartPar
\sphinxcode{\sphinxupquote{NdArray\textless{}int, N\textgreater{} GetLowerIIS()}}
\end{quote}

\sphinxAtStartPar
\sphinxstylestrong{Return}
\begin{quote}

\sphinxAtStartPar
multi\sphinxhyphen{}dimensional array of IIS status of lower bounds of constraints.
\end{quote}
\end{quote}

\subsubsection{MConstr::GetND()}
\label{\detokenize{cppapi/mconstr:mconstr-getnd}}\begin{quote}

\sphinxAtStartPar
Get number of dimensions of MConstr object.

\sphinxAtStartPar
\sphinxstylestrong{Synopsis}
\begin{quote}

\sphinxAtStartPar
\sphinxcode{\sphinxupquote{int GetND()}}
\end{quote}

\sphinxAtStartPar
\sphinxstylestrong{Return}
\begin{quote}

\sphinxAtStartPar
number of dimensions.
\end{quote}
\end{quote}

\subsubsection{MConstr::GetShape()}
\label{\detokenize{cppapi/mconstr:mconstr-getshape}}\begin{quote}

\sphinxAtStartPar
Get shape of MConstr object.

\sphinxAtStartPar
\sphinxstylestrong{Synopsis}
\begin{quote}

\sphinxAtStartPar
\sphinxcode{\sphinxupquote{Shape\textless{}N\textgreater{} GetShape()}}
\end{quote}

\sphinxAtStartPar
\sphinxstylestrong{Return}
\begin{quote}

\sphinxAtStartPar
shape object.
\end{quote}
\end{quote}

\subsubsection{MConstr::GetSize()}
\label{\detokenize{cppapi/mconstr:mconstr-getsize}}\begin{quote}

\sphinxAtStartPar
Get size of MConstr object.

\sphinxAtStartPar
\sphinxstylestrong{Synopsis}
\begin{quote}

\sphinxAtStartPar
\sphinxcode{\sphinxupquote{size\_t GetSize()}}
\end{quote}

\sphinxAtStartPar
\sphinxstylestrong{Return}
\begin{quote}

\sphinxAtStartPar
number of constraints.
\end{quote}
\end{quote}

\subsubsection{MConstr::GetUpperIIS()}
\label{\detokenize{cppapi/mconstr:mconstr-getupperiis}}\begin{quote}

\sphinxAtStartPar
Get IIS status of upper bound of constraints in MConstr object.

\sphinxAtStartPar
\sphinxstylestrong{Synopsis}
\begin{quote}

\sphinxAtStartPar
\sphinxcode{\sphinxupquote{NdArray\textless{}int, N\textgreater{} GetUpperIIS()}}
\end{quote}

\sphinxAtStartPar
\sphinxstylestrong{Return}
\begin{quote}

\sphinxAtStartPar
multi\sphinxhyphen{}dimensional array of IIS status of upper bounds of constraints.
\end{quote}
\end{quote}

\subsubsection{MConstr::Item()}
\label{\detokenize{cppapi/mconstr:mconstr-item}}\begin{quote}

\sphinxAtStartPar
Get constraint of given index from MConstr object.

\sphinxAtStartPar
\sphinxstylestrong{Synopsis}
\begin{quote}

\sphinxAtStartPar
\sphinxcode{\sphinxupquote{Constraint \&Item(size\_t idx)}}
\end{quote}

\sphinxAtStartPar
\sphinxstylestrong{Arguments}
\begin{quote}

\sphinxAtStartPar
\sphinxcode{\sphinxupquote{idx}}: index of constraint.
\end{quote}

\sphinxAtStartPar
\sphinxstylestrong{Return}
\begin{quote}

\sphinxAtStartPar
Constraint object.
\end{quote}
\end{quote}

\subsubsection{MConstr::operator{[}{]}()}
\label{\detokenize{cppapi/mconstr:mconstr-operator}}\begin{quote}

\sphinxAtStartPar
Get constraint of given index from MConstr object.

\sphinxAtStartPar
\sphinxstylestrong{Synopsis}
\begin{quote}

\sphinxAtStartPar
\sphinxcode{\sphinxupquote{Constraint \&operator{[}{]}(size\_t idx)}}
\end{quote}

\sphinxAtStartPar
\sphinxstylestrong{Arguments}
\begin{quote}

\sphinxAtStartPar
\sphinxcode{\sphinxupquote{idx}}: index of constraint.
\end{quote}

\sphinxAtStartPar
\sphinxstylestrong{Return}
\begin{quote}

\sphinxAtStartPar
Constraint object.
\end{quote}
\end{quote}

\subsubsection{MConstr::operator{[}{]}()}
\label{\detokenize{cppapi/mconstr:id2}}\begin{quote}

\sphinxAtStartPar
Get constraints of given view from MConstr object.

\sphinxAtStartPar
\sphinxstylestrong{Synopsis}
\begin{quote}

\sphinxAtStartPar
\sphinxcode{\sphinxupquote{MConstr operator{[}{]}(const View \&view)}}
\end{quote}

\sphinxAtStartPar
\sphinxstylestrong{Arguments}
\begin{quote}

\sphinxAtStartPar
\sphinxcode{\sphinxupquote{view}}: view of multi\sphinxhyphen{}dimensional array.
\end{quote}

\sphinxAtStartPar
\sphinxstylestrong{Return}
\begin{quote}

\sphinxAtStartPar
new MConstr object.
\end{quote}
\end{quote}

\subsubsection{MConstr::Pick()}
\label{\detokenize{cppapi/mconstr:mconstr-pick}}\begin{quote}

\sphinxAtStartPar
Given a list of indexes, get constraints from MConstr object.

\sphinxAtStartPar
\sphinxstylestrong{Synopsis}
\begin{quote}

\sphinxAtStartPar
\sphinxcode{\sphinxupquote{MConstr\textless{}1\textgreater{} Pick(const NdArray\textless{}int, 1\textgreater{} \&indexes)}}
\end{quote}

\sphinxAtStartPar
\sphinxstylestrong{Arguments}
\begin{quote}

\sphinxAtStartPar
\sphinxcode{\sphinxupquote{indexes}}: indexes of elements.
\end{quote}

\sphinxAtStartPar
\sphinxstylestrong{Return}
\begin{quote}

\sphinxAtStartPar
one\sphinxhyphen{}dimensional array of desired constraints.
\end{quote}
\end{quote}

\subsubsection{MConstr::Pick()}
\label{\detokenize{cppapi/mconstr:id3}}\begin{quote}

\sphinxAtStartPar
Given a list of indexes, get constraints from MConstr object.

\sphinxAtStartPar
\sphinxstylestrong{Synopsis}
\begin{quote}

\sphinxAtStartPar
\sphinxcode{\sphinxupquote{MConstr\textless{}1\textgreater{} Pick(const NdArray\textless{}int, 2\textgreater{} \&idxrows)}}
\end{quote}

\sphinxAtStartPar
\sphinxstylestrong{Arguments}
\begin{quote}

\sphinxAtStartPar
\sphinxcode{\sphinxupquote{idxrows}}: indexes in format of 2\sphinxhyphen{}dimensional array, where each row is position of element.
\end{quote}

\sphinxAtStartPar
\sphinxstylestrong{Return}
\begin{quote}

\sphinxAtStartPar
one\sphinxhyphen{}dimensional array of desired constraints.
\end{quote}
\end{quote}

\subsubsection{MConstr::Represent()}
\label{\detokenize{cppapi/mconstr:mconstr-represent}}\begin{quote}

\sphinxAtStartPar
String representation of MConstr object.

\sphinxAtStartPar
\sphinxstylestrong{Synopsis}
\begin{quote}

\sphinxAtStartPar
\sphinxcode{\sphinxupquote{std::string Represent(size\_t maxlen)}}
\end{quote}

\sphinxAtStartPar
\sphinxstylestrong{Arguments}
\begin{quote}

\sphinxAtStartPar
\sphinxcode{\sphinxupquote{maxlen}}: max length of representation.
\end{quote}

\sphinxAtStartPar
\sphinxstylestrong{Return}
\begin{quote}

\sphinxAtStartPar
string object.
\end{quote}
\end{quote}

\subsubsection{MConstr::Reshape()}
\label{\detokenize{cppapi/mconstr:mconstr-reshape}}\begin{quote}

\sphinxAtStartPar
Reshape MConstr object to new shape.

\sphinxAtStartPar
\sphinxstylestrong{Synopsis}
\begin{quote}

\sphinxAtStartPar
\sphinxcode{\sphinxupquote{template \textless{}int M\textgreater{} MConstr\textless{}M\textgreater{} Reshape(const Shape\textless{}M\textgreater{} \&shape)}}
\end{quote}

\sphinxAtStartPar
\sphinxstylestrong{Arguments}
\begin{quote}

\sphinxAtStartPar
\sphinxcode{\sphinxupquote{shape}}: new shape of M\sphinxhyphen{}dimensions.
\end{quote}

\sphinxAtStartPar
\sphinxstylestrong{Return}
\begin{quote}

\sphinxAtStartPar
M\sphinxhyphen{}dimensional MConstr object.
\end{quote}
\end{quote}

\subsubsection{MConstr::Set()}
\label{\detokenize{cppapi/mconstr:mconstr-set}}\begin{quote}

\sphinxAtStartPar
Set values of information associated with constraints in MConstr object.

\sphinxAtStartPar
\sphinxstylestrong{Synopsis}
\begin{quote}

\sphinxAtStartPar
\sphinxcode{\sphinxupquote{void Set(const char *szInfo, double val)}}
\end{quote}

\sphinxAtStartPar
\sphinxstylestrong{Arguments}
\begin{quote}

\sphinxAtStartPar
\sphinxcode{\sphinxupquote{szInfo}}: name of information.

\sphinxAtStartPar
\sphinxcode{\sphinxupquote{val}}: value of information.
\end{quote}
\end{quote}

\subsubsection{MConstr::Set()}
\label{\detokenize{cppapi/mconstr:id4}}\begin{quote}

\sphinxAtStartPar
Set values of information associated with constraints in MConstr object.

\sphinxAtStartPar
\sphinxstylestrong{Synopsis}
\begin{quote}

\sphinxAtStartPar
\sphinxcode{\sphinxupquote{void Set(const char *szInfo, const NdArray\textless{}double, N\textgreater{} \&vals)}}
\end{quote}

\sphinxAtStartPar
\sphinxstylestrong{Arguments}
\begin{quote}

\sphinxAtStartPar
\sphinxcode{\sphinxupquote{szInfo}}: name of information.

\sphinxAtStartPar
\sphinxcode{\sphinxupquote{vals}}: multi\sphinxhyphen{}dimensional array of values of information.
\end{quote}
\end{quote}

\subsubsection{MConstr::SetItem()}
\label{\detokenize{cppapi/mconstr:mconstr-setitem}}\begin{quote}

\sphinxAtStartPar
Set constraint of given index to MConstr object.

\sphinxAtStartPar
\sphinxstylestrong{Synopsis}
\begin{quote}

\sphinxAtStartPar
\sphinxcode{\sphinxupquote{void SetItem(size\_t idx, const Constraint \&con)}}
\end{quote}

\sphinxAtStartPar
\sphinxstylestrong{Arguments}
\begin{quote}

\sphinxAtStartPar
\sphinxcode{\sphinxupquote{idx}}: index of element.

\sphinxAtStartPar
\sphinxcode{\sphinxupquote{con}}: constraint object.
\end{quote}
\end{quote}

\subsubsection{MConstr::Squeeze()}
\label{\detokenize{cppapi/mconstr:mconstr-squeeze}}\begin{quote}

\sphinxAtStartPar
Remove axis of length 1 from shape of MConstr object.

\sphinxAtStartPar
\sphinxstylestrong{Synopsis}
\begin{quote}

\sphinxAtStartPar
\sphinxcode{\sphinxupquote{MConstr\textless{}N \sphinxhyphen{} 1\textgreater{} Squeeze(int axis)}}
\end{quote}

\sphinxAtStartPar
\sphinxstylestrong{Arguments}
\begin{quote}

\sphinxAtStartPar
\sphinxcode{\sphinxupquote{axis}}: axis of MConstr, where the length is 1.
\end{quote}

\sphinxAtStartPar
\sphinxstylestrong{Return}
\begin{quote}

\sphinxAtStartPar
MConstr object of (N\sphinxhyphen{}1)\sphinxhyphen{}dimensional shape.
\end{quote}
\end{quote}

\subsubsection{MConstr::Stack()}
\label{\detokenize{cppapi/mconstr:mconstr-stack}}\begin{quote}

\sphinxAtStartPar
Stack with other MConstr object along given axis.

\sphinxAtStartPar
\sphinxstylestrong{Synopsis}
\begin{quote}

\sphinxAtStartPar
\sphinxcode{\sphinxupquote{MConstr\textless{}N\textgreater{} Stack(const MConstr\textless{}N\textgreater{} \&other, int axis)}}
\end{quote}

\sphinxAtStartPar
\sphinxstylestrong{Arguments}
\begin{quote}

\sphinxAtStartPar
\sphinxcode{\sphinxupquote{other}}: a MConstr object.

\sphinxAtStartPar
\sphinxcode{\sphinxupquote{axis}}: an axis of MConstr.
\end{quote}

\sphinxAtStartPar
\sphinxstylestrong{Return}
\begin{quote}

\sphinxAtStartPar
the result MConstr object.
\end{quote}
\end{quote}

\subsubsection{MConstr::Transpose()}
\label{\detokenize{cppapi/mconstr:mconstr-transpose}}\begin{quote}

\sphinxAtStartPar
Perform matrix transpose of MConstr object.

\sphinxAtStartPar
\sphinxstylestrong{Synopsis}
\begin{quote}

\sphinxAtStartPar
\sphinxcode{\sphinxupquote{MConstr\textless{}N\textgreater{} Transpose()}}
\end{quote}

\sphinxAtStartPar
\sphinxstylestrong{Return}
\begin{quote}

\sphinxAtStartPar
transposed MConstr object.
\end{quote}
\end{quote}

\subsection{MConstrBuilder}
\label{\detokenize{cppapiref:mconstrbuilder}}\label{\detokenize{cppapiref:chapcppapiref-mconstrbuilder}}
\sphinxAtStartPar
COPT builder object of multi\sphinxhyphen{}dimensional linear constraints. It is used to
generate multi\sphinxhyphen{}dimensional linear constraints and support operations with the
built\sphinxhyphen{}in multi\sphinxhyphen{}dimensional array {\hyperref[\detokenize{cppapiref:chapcppapiref-ndarray}]{\sphinxcrossref{\DUrole{std,std-ref}{NdArray}}}} in COPT.
It is recommended to create MConstrBuilder object by comparing two objects, one of
which should be {\hyperref[\detokenize{cppapiref:chapcppapiref-mvar}]{\sphinxcrossref{\DUrole{std,std-ref}{MVar}}}} object or {\hyperref[\detokenize{cppapiref:chapcppapiref-mlinexpr}]{\sphinxcrossref{\DUrole{std,std-ref}{MLinExpr}}}} object,
by comparison operators. The following methods are provided:

\sphinxstepscope

\subsubsection{MConstrBuilder::MConstrBuilder()}
\label{\detokenize{cppapi/mconstrbuilder:mconstrbuilder-mconstrbuilder}}\label{\detokenize{cppapi/mconstrbuilder::doc}}\begin{quote}

\sphinxAtStartPar
Construct a MConstrBuilder object with the given shape.

\sphinxAtStartPar
\sphinxstylestrong{Synopsis}
\begin{quote}

\sphinxAtStartPar
\sphinxcode{\sphinxupquote{MConstrBuilder(const Shape\textless{}N\textgreater{} \&shp)}}
\end{quote}

\sphinxAtStartPar
\sphinxstylestrong{Arguments}
\begin{quote}

\sphinxAtStartPar
\sphinxcode{\sphinxupquote{shp}}: shape of MConstrBuilder.
\end{quote}
\end{quote}

\subsubsection{MConstrBuilder::Flatten()}
\label{\detokenize{cppapi/mconstrbuilder:mconstrbuilder-flatten}}\begin{quote}

\sphinxAtStartPar
Flatten a MConstrBuilder object to a 1\sphinxhyphen{}dimensional shape.

\sphinxAtStartPar
\sphinxstylestrong{Synopsis}
\begin{quote}

\sphinxAtStartPar
\sphinxcode{\sphinxupquote{MConstrBuilder\textless{}1\textgreater{} Flatten()}}
\end{quote}

\sphinxAtStartPar
\sphinxstylestrong{Return}
\begin{quote}

\sphinxAtStartPar
a MConstrBuilder object collapsed into one dimension.
\end{quote}
\end{quote}

\subsubsection{MConstrBuilder::GetExpr()}
\label{\detokenize{cppapi/mconstrbuilder:mconstrbuilder-getexpr}}\begin{quote}

\sphinxAtStartPar
Get N\sphinxhyphen{}dimensional linear expressions associated with N\sphinxhyphen{}dimensional constraints.

\sphinxAtStartPar
\sphinxstylestrong{Synopsis}
\begin{quote}

\sphinxAtStartPar
\sphinxcode{\sphinxupquote{const MLinExpr\textless{}N\textgreater{} \&GetExpr()}}
\end{quote}

\sphinxAtStartPar
\sphinxstylestrong{Return}
\begin{quote}

\sphinxAtStartPar
MLinExpr object.
\end{quote}
\end{quote}

\subsubsection{MConstrBuilder::GetND()}
\label{\detokenize{cppapi/mconstrbuilder:mconstrbuilder-getnd}}\begin{quote}

\sphinxAtStartPar
Get number of dimensions of MConstrBuilder object.

\sphinxAtStartPar
\sphinxstylestrong{Synopsis}
\begin{quote}

\sphinxAtStartPar
\sphinxcode{\sphinxupquote{int GetND()}}
\end{quote}

\sphinxAtStartPar
\sphinxstylestrong{Return}
\begin{quote}

\sphinxAtStartPar
number of dimensions.
\end{quote}
\end{quote}

\subsubsection{MConstrBuilder::GetRange()}
\label{\detokenize{cppapi/mconstrbuilder:mconstrbuilder-getrange}}\begin{quote}

\sphinxAtStartPar
Get range from lower bound to upper bound of N\sphinxhyphen{}dimensional range constraints.

\sphinxAtStartPar
\sphinxstylestrong{Synopsis}
\begin{quote}

\sphinxAtStartPar
\sphinxcode{\sphinxupquote{double GetRange()}}
\end{quote}

\sphinxAtStartPar
\sphinxstylestrong{Return}
\begin{quote}

\sphinxAtStartPar
length from lower bound to upper bound of range constraints.
\end{quote}
\end{quote}

\subsubsection{MConstrBuilder::GetSense()}
\label{\detokenize{cppapi/mconstrbuilder:mconstrbuilder-getsense}}\begin{quote}

\sphinxAtStartPar
Get sense associated with N\sphinxhyphen{}dimensional constraints.

\sphinxAtStartPar
\sphinxstylestrong{Synopsis}
\begin{quote}

\sphinxAtStartPar
\sphinxcode{\sphinxupquote{char GetSense()}}
\end{quote}

\sphinxAtStartPar
\sphinxstylestrong{Return}
\begin{quote}

\sphinxAtStartPar
constraint sense.
\end{quote}
\end{quote}

\subsubsection{MConstrBuilder::Set()}
\label{\detokenize{cppapi/mconstrbuilder:mconstrbuilder-set}}\begin{quote}

\sphinxAtStartPar
Set N\sphinxhyphen{}dimensional constraints to its builder object.

\sphinxAtStartPar
\sphinxstylestrong{Synopsis}
\begin{quote}

\sphinxAtStartPar
\sphinxcode{\sphinxupquote{void Set(}}
\begin{quote}

\sphinxAtStartPar
\sphinxcode{\sphinxupquote{const MLinExpr\textless{}N\textgreater{} \&expr,}}

\sphinxAtStartPar
\sphinxcode{\sphinxupquote{char sense,}}

\sphinxAtStartPar
\sphinxcode{\sphinxupquote{double rhs)}}
\end{quote}
\end{quote}

\sphinxAtStartPar
\sphinxstylestrong{Arguments}
\begin{quote}

\sphinxAtStartPar
\sphinxcode{\sphinxupquote{expr}}: MLinExpr object

\sphinxAtStartPar
\sphinxcode{\sphinxupquote{sense}}: constraint sense other than COPT\_RANGE.

\sphinxAtStartPar
\sphinxcode{\sphinxupquote{rhs}}: constant of right side of constraints.
\end{quote}
\end{quote}

\subsubsection{MConstrBuilder::Set()}
\label{\detokenize{cppapi/mconstrbuilder:id1}}\begin{quote}

\sphinxAtStartPar
Set N\sphinxhyphen{}dimensional constraints to its builder object.

\sphinxAtStartPar
\sphinxstylestrong{Synopsis}
\begin{quote}

\sphinxAtStartPar
\sphinxcode{\sphinxupquote{template \textless{}class T\textgreater{} void Set(}}
\begin{quote}

\sphinxAtStartPar
\sphinxcode{\sphinxupquote{const MLinExpr\textless{}N\textgreater{} \&expr,}}

\sphinxAtStartPar
\sphinxcode{\sphinxupquote{char sense,}}

\sphinxAtStartPar
\sphinxcode{\sphinxupquote{const NdArray\textless{}T, N\textgreater{} \&rhs)}}
\end{quote}
\end{quote}

\sphinxAtStartPar
\sphinxstylestrong{Arguments}
\begin{quote}

\sphinxAtStartPar
\sphinxcode{\sphinxupquote{expr}}: MLinExpr object

\sphinxAtStartPar
\sphinxcode{\sphinxupquote{sense}}: constraint sense other than COPT\_RANGE.

\sphinxAtStartPar
\sphinxcode{\sphinxupquote{rhs}}: N\sphinxhyphen{}dimensional constants at right side of constraints.
\end{quote}
\end{quote}

\subsubsection{MConstrBuilder::Set()}
\label{\detokenize{cppapi/mconstrbuilder:id2}}\begin{quote}

\sphinxAtStartPar
Set N\sphinxhyphen{}dimensional constraints to its builder object.

\sphinxAtStartPar
\sphinxstylestrong{Synopsis}
\begin{quote}

\sphinxAtStartPar
\sphinxcode{\sphinxupquote{template \textless{}int M\textgreater{} void Set(}}
\begin{quote}

\sphinxAtStartPar
\sphinxcode{\sphinxupquote{const MLinExpr\textless{}N\textgreater{} \&expr,}}

\sphinxAtStartPar
\sphinxcode{\sphinxupquote{char sense,}}

\sphinxAtStartPar
\sphinxcode{\sphinxupquote{const MVar\textless{}M\textgreater{} \&rhs)}}
\end{quote}
\end{quote}

\sphinxAtStartPar
\sphinxstylestrong{Arguments}
\begin{quote}

\sphinxAtStartPar
\sphinxcode{\sphinxupquote{expr}}: MLinExpr object

\sphinxAtStartPar
\sphinxcode{\sphinxupquote{sense}}: constraint sense other than COPT\_RANGE.

\sphinxAtStartPar
\sphinxcode{\sphinxupquote{rhs}}: MVar object at right side of constraints.
\end{quote}
\end{quote}

\subsubsection{MConstrBuilder::Set()}
\label{\detokenize{cppapi/mconstrbuilder:id3}}\begin{quote}

\sphinxAtStartPar
Set N\sphinxhyphen{}dimensional constraints to its builder object.

\sphinxAtStartPar
\sphinxstylestrong{Synopsis}
\begin{quote}

\sphinxAtStartPar
\sphinxcode{\sphinxupquote{template \textless{}int M\textgreater{} void Set(}}
\begin{quote}

\sphinxAtStartPar
\sphinxcode{\sphinxupquote{const MLinExpr\textless{}N\textgreater{} \&expr,}}

\sphinxAtStartPar
\sphinxcode{\sphinxupquote{char sense,}}

\sphinxAtStartPar
\sphinxcode{\sphinxupquote{const MLinExpr\textless{}M\textgreater{} \&rhs)}}
\end{quote}
\end{quote}

\sphinxAtStartPar
\sphinxstylestrong{Arguments}
\begin{quote}

\sphinxAtStartPar
\sphinxcode{\sphinxupquote{expr}}: MLinExpr object

\sphinxAtStartPar
\sphinxcode{\sphinxupquote{sense}}: constraint sense other than COPT\_RANGE.

\sphinxAtStartPar
\sphinxcode{\sphinxupquote{rhs}}: MLinExpr object at right side of constraints.
\end{quote}
\end{quote}

\subsubsection{MConstrBuilder::SetRange()}
\label{\detokenize{cppapi/mconstrbuilder:mconstrbuilder-setrange}}\begin{quote}

\sphinxAtStartPar
Set N\sphinxhyphen{}dimensional range constraints to its builder object.

\sphinxAtStartPar
\sphinxstylestrong{Synopsis}
\begin{quote}

\sphinxAtStartPar
\sphinxcode{\sphinxupquote{void SetRange(const MLinExpr\textless{}N\textgreater{} \&expr, double range)}}
\end{quote}

\sphinxAtStartPar
\sphinxstylestrong{Arguments}
\begin{quote}

\sphinxAtStartPar
\sphinxcode{\sphinxupquote{expr}}: MLinExpr object.

\sphinxAtStartPar
\sphinxcode{\sphinxupquote{range}}: length from lower bound to upper bound of the constraint. Must greater than 0.
\end{quote}
\end{quote}

\subsection{MExpression}
\label{\detokenize{cppapiref:mexpression}}\label{\detokenize{cppapiref:chapcppapiref-mexpression}}
\sphinxAtStartPar
The MExpression class is a generalized version of {\hyperref[\detokenize{cppapiref:chapcppapiref-expr}]{\sphinxcrossref{\DUrole{std,std-ref}{Expr}}}}.
It represents a linear expression and supports most of methods in Expr class.
In addition, it supports linear combination of multi\sphinxhyphen{}dimensional objects,
such as {\hyperref[\detokenize{cppapiref:chapcppapiref-mvar}]{\sphinxcrossref{\DUrole{std,std-ref}{MVar}}}} object and {\hyperref[\detokenize{cppapiref:chapcppapiref-ndarray}]{\sphinxcrossref{\DUrole{std,std-ref}{NdArray}}}} object.
The following methods are provided:

\sphinxstepscope

\subsubsection{MExpression::MExpression()}
\label{\detokenize{cppapi/mexpression:mexpression-mexpression}}\label{\detokenize{cppapi/mexpression::doc}}\begin{quote}

\sphinxAtStartPar
Construct a MExpression object with the given constant.

\sphinxAtStartPar
\sphinxstylestrong{Synopsis}
\begin{quote}

\sphinxAtStartPar
\sphinxcode{\sphinxupquote{MExpression(double constant)}}
\end{quote}

\sphinxAtStartPar
\sphinxstylestrong{Arguments}
\begin{quote}

\sphinxAtStartPar
\sphinxcode{\sphinxupquote{constant}}: constant number.
\end{quote}
\end{quote}

\subsubsection{MExpression::MExpression()}
\label{\detokenize{cppapi/mexpression:id1}}\begin{quote}

\sphinxAtStartPar
Construct a MExpression object with the given linear expression.

\sphinxAtStartPar
\sphinxstylestrong{Synopsis}
\begin{quote}

\sphinxAtStartPar
\sphinxcode{\sphinxupquote{MExpression(const Expr \&expr)}}
\end{quote}

\sphinxAtStartPar
\sphinxstylestrong{Arguments}
\begin{quote}

\sphinxAtStartPar
\sphinxcode{\sphinxupquote{expr}}: a linear expression.
\end{quote}
\end{quote}

\subsubsection{MExpression::MExpression()}
\label{\detokenize{cppapi/mexpression:id2}}\begin{quote}

\sphinxAtStartPar
Construct a MExpression object with the given variable.

\sphinxAtStartPar
\sphinxstylestrong{Synopsis}
\begin{quote}

\sphinxAtStartPar
\sphinxcode{\sphinxupquote{MExpression(const Var \&var)}}
\end{quote}

\sphinxAtStartPar
\sphinxstylestrong{Arguments}
\begin{quote}

\sphinxAtStartPar
\sphinxcode{\sphinxupquote{var}}: variable object.
\end{quote}
\end{quote}

\subsubsection{MExpression::AddConstant()}
\label{\detokenize{cppapi/mexpression:mexpression-addconstant}}\begin{quote}

\sphinxAtStartPar
Add constant for the expression.

\sphinxAtStartPar
\sphinxstylestrong{Synopsis}
\begin{quote}

\sphinxAtStartPar
\sphinxcode{\sphinxupquote{void AddConstant(double constant)}}
\end{quote}

\sphinxAtStartPar
\sphinxstylestrong{Arguments}
\begin{quote}

\sphinxAtStartPar
\sphinxcode{\sphinxupquote{constant}}: the value of the constant.
\end{quote}
\end{quote}

\subsubsection{MExpression::AddExpr()}
\label{\detokenize{cppapi/mexpression:mexpression-addexpr}}\begin{quote}

\sphinxAtStartPar
Add a linear expression to MExpression object.

\sphinxAtStartPar
\sphinxstylestrong{Synopsis}
\begin{quote}

\sphinxAtStartPar
\sphinxcode{\sphinxupquote{void AddExpr(const Expr \&expr, double mult)}}
\end{quote}

\sphinxAtStartPar
\sphinxstylestrong{Arguments}
\begin{quote}

\sphinxAtStartPar
\sphinxcode{\sphinxupquote{expr}}: linear expression object.

\sphinxAtStartPar
\sphinxcode{\sphinxupquote{mult}}: the multiplier of linear expression, default value is 1.0.
\end{quote}
\end{quote}

\subsubsection{MExpression::AddMExpr()}
\label{\detokenize{cppapi/mexpression:mexpression-addmexpr}}\begin{quote}

\sphinxAtStartPar
Add MExpression to MExpression object.

\sphinxAtStartPar
\sphinxstylestrong{Synopsis}
\begin{quote}

\sphinxAtStartPar
\sphinxcode{\sphinxupquote{void AddMExpr(const MExpression \&expr, double mult)}}
\end{quote}

\sphinxAtStartPar
\sphinxstylestrong{Arguments}
\begin{quote}

\sphinxAtStartPar
\sphinxcode{\sphinxupquote{expr}}: MExpression object.

\sphinxAtStartPar
\sphinxcode{\sphinxupquote{mult}}: the multiplier of MExpression, default value is 1.0.
\end{quote}
\end{quote}

\subsubsection{MExpression::AddTerm()}
\label{\detokenize{cppapi/mexpression:mexpression-addterm}}\begin{quote}

\sphinxAtStartPar
Add a linear term to MExpression object.

\sphinxAtStartPar
\sphinxstylestrong{Synopsis}
\begin{quote}

\sphinxAtStartPar
\sphinxcode{\sphinxupquote{void AddTerm(const Var \&var, double coeff)}}
\end{quote}

\sphinxAtStartPar
\sphinxstylestrong{Arguments}
\begin{quote}

\sphinxAtStartPar
\sphinxcode{\sphinxupquote{var}}: variable of new term.

\sphinxAtStartPar
\sphinxcode{\sphinxupquote{coeff}}: coefficient of new term.
\end{quote}
\end{quote}

\subsubsection{MExpression::Clone()}
\label{\detokenize{cppapi/mexpression:mexpression-clone}}\begin{quote}

\sphinxAtStartPar
Clone MExpression object.

\sphinxAtStartPar
\sphinxstylestrong{Synopsis}
\begin{quote}

\sphinxAtStartPar
\sphinxcode{\sphinxupquote{MExpression Clone()}}
\end{quote}

\sphinxAtStartPar
\sphinxstylestrong{Return}
\begin{quote}

\sphinxAtStartPar
new MExpression object.
\end{quote}
\end{quote}

\subsubsection{MExpression::Evaluate()}
\label{\detokenize{cppapi/mexpression:mexpression-evaluate}}\begin{quote}

\sphinxAtStartPar
Evaluate MExpression after solving.

\sphinxAtStartPar
\sphinxstylestrong{Synopsis}
\begin{quote}

\sphinxAtStartPar
\sphinxcode{\sphinxupquote{double Evaluate()}}
\end{quote}

\sphinxAtStartPar
\sphinxstylestrong{Return}
\begin{quote}

\sphinxAtStartPar
value of MExpression object.
\end{quote}
\end{quote}

\subsubsection{MExpression::GetConstant()}
\label{\detokenize{cppapi/mexpression:mexpression-getconstant}}\begin{quote}

\sphinxAtStartPar
Get constant in expression.

\sphinxAtStartPar
\sphinxstylestrong{Synopsis}
\begin{quote}

\sphinxAtStartPar
\sphinxcode{\sphinxupquote{double GetConstant()}}
\end{quote}

\sphinxAtStartPar
\sphinxstylestrong{Return}
\begin{quote}

\sphinxAtStartPar
constant in expression.
\end{quote}
\end{quote}

\subsubsection{MExpression::Represent()}
\label{\detokenize{cppapi/mexpression:mexpression-represent}}\begin{quote}

\sphinxAtStartPar
String representation of MExpression object.

\sphinxAtStartPar
\sphinxstylestrong{Synopsis}
\begin{quote}

\sphinxAtStartPar
\sphinxcode{\sphinxupquote{std::string Represent(size\_t maxlen)}}
\end{quote}

\sphinxAtStartPar
\sphinxstylestrong{Arguments}
\begin{quote}

\sphinxAtStartPar
\sphinxcode{\sphinxupquote{maxlen}}: max length of representation.
\end{quote}

\sphinxAtStartPar
\sphinxstylestrong{Return}
\begin{quote}

\sphinxAtStartPar
string object.
\end{quote}
\end{quote}

\subsubsection{MExpression::SetConstant()}
\label{\detokenize{cppapi/mexpression:mexpression-setconstant}}\begin{quote}

\sphinxAtStartPar
Set constant for the expression.

\sphinxAtStartPar
\sphinxstylestrong{Synopsis}
\begin{quote}

\sphinxAtStartPar
\sphinxcode{\sphinxupquote{void SetConstant(double constant)}}
\end{quote}

\sphinxAtStartPar
\sphinxstylestrong{Arguments}
\begin{quote}

\sphinxAtStartPar
\sphinxcode{\sphinxupquote{constant}}: the value of the constant.
\end{quote}
\end{quote}

\subsection{MLinExpr}
\label{\detokenize{cppapiref:mlinexpr}}\label{\detokenize{cppapiref:chapcppapiref-mlinexpr}}
\sphinxAtStartPar
COPT multi\sphinxhyphen{}dimensional linear expression object. It is used to construct
multi\sphinxhyphen{}dimensional linear expressions and perform operations with the built\sphinxhyphen{}in
multi\sphinxhyphen{}dimensional array {\hyperref[\detokenize{cppapiref:chapcppapiref-ndarray}]{\sphinxcrossref{\DUrole{std,std-ref}{NdArray}}}} in COPT. Its elements
are {\hyperref[\detokenize{cppapiref:chapcppapiref-mexpression}]{\sphinxcrossref{\DUrole{std,std-ref}{MExpression}}}} objects.  It can be created by linear
combination of {\hyperref[\detokenize{cppapiref:chapcppapiref-mvar}]{\sphinxcrossref{\DUrole{std,std-ref}{MVar}}}} objects. The following methods are
provided:

\sphinxstepscope

\subsubsection{MLinExpr::MLinExpr()}
\label{\detokenize{cppapi/mlinexpr:mlinexpr-mlinexpr}}\label{\detokenize{cppapi/mlinexpr::doc}}\begin{quote}

\sphinxAtStartPar
Construct a MLinExpr object with the given shape and a constant.

\sphinxAtStartPar
\sphinxstylestrong{Synopsis}
\begin{quote}

\sphinxAtStartPar
\sphinxcode{\sphinxupquote{MLinExpr(const Shape\textless{}N\textgreater{} \&shp, double constant)}}
\end{quote}

\sphinxAtStartPar
\sphinxstylestrong{Arguments}
\begin{quote}

\sphinxAtStartPar
\sphinxcode{\sphinxupquote{shp}}: shape of MLinExpr.

\sphinxAtStartPar
\sphinxcode{\sphinxupquote{constant}}: constant number.
\end{quote}
\end{quote}

\subsubsection{MLinExpr::MLinExpr()}
\label{\detokenize{cppapi/mlinexpr:id1}}\begin{quote}

\sphinxAtStartPar
Construct a MLinExpr object with the given shape and a linear expression.

\sphinxAtStartPar
\sphinxstylestrong{Synopsis}
\begin{quote}

\sphinxAtStartPar
\sphinxcode{\sphinxupquote{MLinExpr(const Shape\textless{}N\textgreater{} \&shp, const Expr \&expr)}}
\end{quote}

\sphinxAtStartPar
\sphinxstylestrong{Arguments}
\begin{quote}

\sphinxAtStartPar
\sphinxcode{\sphinxupquote{shp}}: shape of MLinExpr.

\sphinxAtStartPar
\sphinxcode{\sphinxupquote{expr}}: a linear expression.
\end{quote}
\end{quote}

\subsubsection{MLinExpr::MLinExpr()}
\label{\detokenize{cppapi/mlinexpr:id2}}\begin{quote}

\sphinxAtStartPar
Construct a MLinExpr object with the given shape and a MExpression object.

\sphinxAtStartPar
\sphinxstylestrong{Synopsis}
\begin{quote}

\sphinxAtStartPar
\sphinxcode{\sphinxupquote{MLinExpr(const Shape\textless{}N\textgreater{} \&shp, const MExpression \&expr)}}
\end{quote}

\sphinxAtStartPar
\sphinxstylestrong{Arguments}
\begin{quote}

\sphinxAtStartPar
\sphinxcode{\sphinxupquote{shp}}: shape of MLinExpr.

\sphinxAtStartPar
\sphinxcode{\sphinxupquote{expr}}: a MExpression object.
\end{quote}
\end{quote}

\subsubsection{MLinExpr::AddConstant()}
\label{\detokenize{cppapi/mlinexpr:mlinexpr-addconstant}}\begin{quote}

\sphinxAtStartPar
Add constant to each expression in MLinExpr object.

\sphinxAtStartPar
\sphinxstylestrong{Synopsis}
\begin{quote}

\sphinxAtStartPar
\sphinxcode{\sphinxupquote{void AddConstant(double constant)}}
\end{quote}

\sphinxAtStartPar
\sphinxstylestrong{Arguments}
\begin{quote}

\sphinxAtStartPar
\sphinxcode{\sphinxupquote{constant}}: the value of the constant.
\end{quote}
\end{quote}

\subsubsection{MLinExpr::AddConstant()}
\label{\detokenize{cppapi/mlinexpr:id3}}\begin{quote}

\sphinxAtStartPar
Add constants to each expression in MLinExpr object.

\sphinxAtStartPar
\sphinxstylestrong{Synopsis}
\begin{quote}

\sphinxAtStartPar
\sphinxcode{\sphinxupquote{template \textless{}class T\textgreater{} void AddConstant(const NdArray\textless{}T, N\textgreater{} \&constants)}}
\end{quote}

\sphinxAtStartPar
\sphinxstylestrong{Arguments}
\begin{quote}

\sphinxAtStartPar
\sphinxcode{\sphinxupquote{constants}}: N\sphinxhyphen{}dimension NdArray object.
\end{quote}
\end{quote}

\subsubsection{MLinExpr::AddExpr()}
\label{\detokenize{cppapi/mlinexpr:mlinexpr-addexpr}}\begin{quote}

\sphinxAtStartPar
Add a linear expression to each expression in MLinExpr object.

\sphinxAtStartPar
\sphinxstylestrong{Synopsis}
\begin{quote}

\sphinxAtStartPar
\sphinxcode{\sphinxupquote{void AddExpr(const Expr \&expr, double mult)}}
\end{quote}

\sphinxAtStartPar
\sphinxstylestrong{Arguments}
\begin{quote}

\sphinxAtStartPar
\sphinxcode{\sphinxupquote{expr}}: linear expression object.

\sphinxAtStartPar
\sphinxcode{\sphinxupquote{mult}}: the multiplier of linear expression, default value is 1.0.
\end{quote}
\end{quote}

\subsubsection{MLinExpr::AddMExpr()}
\label{\detokenize{cppapi/mlinexpr:mlinexpr-addmexpr}}\begin{quote}

\sphinxAtStartPar
Add MExpression to each expression in MLinExpr object.

\sphinxAtStartPar
\sphinxstylestrong{Synopsis}
\begin{quote}

\sphinxAtStartPar
\sphinxcode{\sphinxupquote{void AddMExpr(const MExpression \&expr, double mult)}}
\end{quote}

\sphinxAtStartPar
\sphinxstylestrong{Arguments}
\begin{quote}

\sphinxAtStartPar
\sphinxcode{\sphinxupquote{expr}}: MExpression object.

\sphinxAtStartPar
\sphinxcode{\sphinxupquote{mult}}: the multiplier of MExpression, default value is 1.0.
\end{quote}
\end{quote}

\subsubsection{MLinExpr::AddMLinExpr()}
\label{\detokenize{cppapi/mlinexpr:mlinexpr-addmlinexpr}}\begin{quote}

\sphinxAtStartPar
Add linear expressions to MLinExpr object.

\sphinxAtStartPar
\sphinxstylestrong{Synopsis}
\begin{quote}

\sphinxAtStartPar
\sphinxcode{\sphinxupquote{void AddMLinExpr(const MLinExpr\textless{}N\textgreater{} \&exprs, double mult)}}
\end{quote}

\sphinxAtStartPar
\sphinxstylestrong{Arguments}
\begin{quote}

\sphinxAtStartPar
\sphinxcode{\sphinxupquote{exprs}}: N\sphinxhyphen{}dimension MLinExpr object.

\sphinxAtStartPar
\sphinxcode{\sphinxupquote{mult}}: the same multiplier for added linear expressions, default value is 1.0.
\end{quote}
\end{quote}

\subsubsection{MLinExpr::AddTerm()}
\label{\detokenize{cppapi/mlinexpr:mlinexpr-addterm}}\begin{quote}

\sphinxAtStartPar
Add a linear term to MLinExpr object.

\sphinxAtStartPar
\sphinxstylestrong{Synopsis}
\begin{quote}

\sphinxAtStartPar
\sphinxcode{\sphinxupquote{void AddTerm(const Var \&var, double coeff)}}
\end{quote}

\sphinxAtStartPar
\sphinxstylestrong{Arguments}
\begin{quote}

\sphinxAtStartPar
\sphinxcode{\sphinxupquote{var}}: variable of new term.

\sphinxAtStartPar
\sphinxcode{\sphinxupquote{coeff}}: coefficient of new term.
\end{quote}
\end{quote}

\subsubsection{MLinExpr::AddTerms()}
\label{\detokenize{cppapi/mlinexpr:mlinexpr-addterms}}\begin{quote}

\sphinxAtStartPar
Add terms to expressions in MLinExpr object.

\sphinxAtStartPar
\sphinxstylestrong{Synopsis}
\begin{quote}

\sphinxAtStartPar
\sphinxcode{\sphinxupquote{void AddTerms(const MVar\textless{}N\textgreater{} \&vars, double mult)}}
\end{quote}

\sphinxAtStartPar
\sphinxstylestrong{Arguments}
\begin{quote}

\sphinxAtStartPar
\sphinxcode{\sphinxupquote{vars}}: N\sphinxhyphen{}dimension MVar object for added terms.

\sphinxAtStartPar
\sphinxcode{\sphinxupquote{mult}}: the same coefficient for added terms, default value 1.0.
\end{quote}
\end{quote}

\subsubsection{MLinExpr::AddTerms()}
\label{\detokenize{cppapi/mlinexpr:id4}}\begin{quote}

\sphinxAtStartPar
Add terms to expressions in MLinExpr object.

\sphinxAtStartPar
\sphinxstylestrong{Synopsis}
\begin{quote}

\sphinxAtStartPar
\sphinxcode{\sphinxupquote{void AddTerms(const MVar\textless{}N\textgreater{} \&vars, const NdArray\textless{}double, N\textgreater{} \&coeffs)}}
\end{quote}

\sphinxAtStartPar
\sphinxstylestrong{Arguments}
\begin{quote}

\sphinxAtStartPar
\sphinxcode{\sphinxupquote{vars}}: N\sphinxhyphen{}dimension MVar object for added terms.

\sphinxAtStartPar
\sphinxcode{\sphinxupquote{coeffs}}: N\sphinxhyphen{}dimension NdArray object of coefficients for added terms.
\end{quote}
\end{quote}

\subsubsection{MLinExpr::Clear()}
\label{\detokenize{cppapi/mlinexpr:mlinexpr-clear}}\begin{quote}

\sphinxAtStartPar
Clear MLinExpr object.

\sphinxAtStartPar
\sphinxstylestrong{Synopsis}
\begin{quote}

\sphinxAtStartPar
\sphinxcode{\sphinxupquote{void Clear()}}
\end{quote}
\end{quote}

\subsubsection{MLinExpr::Clone()}
\label{\detokenize{cppapi/mlinexpr:mlinexpr-clone}}\begin{quote}

\sphinxAtStartPar
Clone MLinExpr object.

\sphinxAtStartPar
\sphinxstylestrong{Synopsis}
\begin{quote}

\sphinxAtStartPar
\sphinxcode{\sphinxupquote{MLinExpr Clone()}}
\end{quote}

\sphinxAtStartPar
\sphinxstylestrong{Return}
\begin{quote}

\sphinxAtStartPar
new MLinExpr object.
\end{quote}
\end{quote}

\subsubsection{MLinExpr::Diagonal()}
\label{\detokenize{cppapi/mlinexpr:mlinexpr-diagonal}}\begin{quote}

\sphinxAtStartPar
Get diagonals of MLinExpr object.

\sphinxAtStartPar
\sphinxstylestrong{Synopsis}
\begin{quote}

\sphinxAtStartPar
\sphinxcode{\sphinxupquote{MLinExpr\textless{}N \sphinxhyphen{} 1\textgreater{} Diagonal(}}
\begin{quote}

\sphinxAtStartPar
\sphinxcode{\sphinxupquote{int offset,}}

\sphinxAtStartPar
\sphinxcode{\sphinxupquote{int axis1,}}

\sphinxAtStartPar
\sphinxcode{\sphinxupquote{int axis2)}}
\end{quote}
\end{quote}

\sphinxAtStartPar
\sphinxstylestrong{Arguments}
\begin{quote}

\sphinxAtStartPar
\sphinxcode{\sphinxupquote{offset}}: offset of the diagonal from the main diagonal. Can be positive or negative.

\sphinxAtStartPar
\sphinxcode{\sphinxupquote{axis1}}: 1st axis of MLinExpr.

\sphinxAtStartPar
\sphinxcode{\sphinxupquote{axis2}}: 2nd axis of MLinExpr.
\end{quote}

\sphinxAtStartPar
\sphinxstylestrong{Return}
\begin{quote}

\sphinxAtStartPar
(N\sphinxhyphen{}1)\sphinxhyphen{}dimensional diagonals.
\end{quote}
\end{quote}

\subsubsection{MLinExpr::Evaluate()}
\label{\detokenize{cppapi/mlinexpr:mlinexpr-evaluate}}\begin{quote}

\sphinxAtStartPar
Evaluate MLinExpr object after solving.

\sphinxAtStartPar
\sphinxstylestrong{Synopsis}
\begin{quote}

\sphinxAtStartPar
\sphinxcode{\sphinxupquote{NdArray\textless{}double, N\textgreater{} Evaluate()}}
\end{quote}

\sphinxAtStartPar
\sphinxstylestrong{Return}
\begin{quote}

\sphinxAtStartPar
NdArray object storing value of each linear expression.
\end{quote}
\end{quote}

\subsubsection{MLinExpr::Expand()}
\label{\detokenize{cppapi/mlinexpr:mlinexpr-expand}}\begin{quote}

\sphinxAtStartPar
Expand shape of MLinExpr object.

\sphinxAtStartPar
\sphinxstylestrong{Synopsis}
\begin{quote}

\sphinxAtStartPar
\sphinxcode{\sphinxupquote{MLinExpr\textless{}N + 1\textgreater{} Expand(int axis)}}
\end{quote}

\sphinxAtStartPar
\sphinxstylestrong{Arguments}
\begin{quote}

\sphinxAtStartPar
\sphinxcode{\sphinxupquote{axis}}: axis of MLinExpr.
\end{quote}

\sphinxAtStartPar
\sphinxstylestrong{Return}
\begin{quote}

\sphinxAtStartPar
MLinExpr object of (N+1)\sphinxhyphen{}dimensional shape.
\end{quote}
\end{quote}

\subsubsection{MLinExpr::Flatten()}
\label{\detokenize{cppapi/mlinexpr:mlinexpr-flatten}}\begin{quote}

\sphinxAtStartPar
Flatten a MLinExpr object to a 1\sphinxhyphen{}dimensional shape.

\sphinxAtStartPar
\sphinxstylestrong{Synopsis}
\begin{quote}

\sphinxAtStartPar
\sphinxcode{\sphinxupquote{MLinExpr\textless{}1\textgreater{} Flatten()}}
\end{quote}

\sphinxAtStartPar
\sphinxstylestrong{Return}
\begin{quote}

\sphinxAtStartPar
a MLinExpr object collapsed into one dimension.
\end{quote}
\end{quote}

\subsubsection{MLinExpr::GetDim()}
\label{\detokenize{cppapi/mlinexpr:mlinexpr-getdim}}\begin{quote}

\sphinxAtStartPar
Get i\sphinxhyphen{}th dimension of MLinExpr object.

\sphinxAtStartPar
\sphinxstylestrong{Synopsis}
\begin{quote}

\sphinxAtStartPar
\sphinxcode{\sphinxupquote{size\_t GetDim(int i)}}
\end{quote}

\sphinxAtStartPar
\sphinxstylestrong{Arguments}
\begin{quote}

\sphinxAtStartPar
\sphinxcode{\sphinxupquote{i}}: index of dimension
\end{quote}

\sphinxAtStartPar
\sphinxstylestrong{Return}
\begin{quote}

\sphinxAtStartPar
i\sphinxhyphen{}th dimension.
\end{quote}
\end{quote}

\subsubsection{MLinExpr::GetND()}
\label{\detokenize{cppapi/mlinexpr:mlinexpr-getnd}}\begin{quote}

\sphinxAtStartPar
Get number of dimensions of MLinExpr object.

\sphinxAtStartPar
\sphinxstylestrong{Synopsis}
\begin{quote}

\sphinxAtStartPar
\sphinxcode{\sphinxupquote{int GetND()}}
\end{quote}

\sphinxAtStartPar
\sphinxstylestrong{Return}
\begin{quote}

\sphinxAtStartPar
number of dimensions.
\end{quote}
\end{quote}

\subsubsection{MLinExpr::GetShape()}
\label{\detokenize{cppapi/mlinexpr:mlinexpr-getshape}}\begin{quote}

\sphinxAtStartPar
Get shape of MLinExpr object.

\sphinxAtStartPar
\sphinxstylestrong{Synopsis}
\begin{quote}

\sphinxAtStartPar
\sphinxcode{\sphinxupquote{Shape\textless{}N\textgreater{} GetShape()}}
\end{quote}

\sphinxAtStartPar
\sphinxstylestrong{Return}
\begin{quote}

\sphinxAtStartPar
shape object.
\end{quote}
\end{quote}

\subsubsection{MLinExpr::GetSize()}
\label{\detokenize{cppapi/mlinexpr:mlinexpr-getsize}}\begin{quote}

\sphinxAtStartPar
Get size of MLinExpr object.

\sphinxAtStartPar
\sphinxstylestrong{Synopsis}
\begin{quote}

\sphinxAtStartPar
\sphinxcode{\sphinxupquote{size\_t GetSize()}}
\end{quote}

\sphinxAtStartPar
\sphinxstylestrong{Return}
\begin{quote}

\sphinxAtStartPar
number of MExpressions.
\end{quote}
\end{quote}

\subsubsection{MLinExpr::Item()}
\label{\detokenize{cppapi/mlinexpr:mlinexpr-item}}\begin{quote}

\sphinxAtStartPar
Get expression of given index from MLinExpr object.

\sphinxAtStartPar
\sphinxstylestrong{Synopsis}
\begin{quote}

\sphinxAtStartPar
\sphinxcode{\sphinxupquote{MExpression \&Item(size\_t idx)}}
\end{quote}

\sphinxAtStartPar
\sphinxstylestrong{Arguments}
\begin{quote}

\sphinxAtStartPar
\sphinxcode{\sphinxupquote{idx}}: index of expression.
\end{quote}

\sphinxAtStartPar
\sphinxstylestrong{Return}
\begin{quote}

\sphinxAtStartPar
expression object.
\end{quote}
\end{quote}

\subsubsection{MLinExpr::Item()}
\label{\detokenize{cppapi/mlinexpr:id5}}\begin{quote}

\sphinxAtStartPar
Get sub\sphinxhyphen{}arrays of MLinExpr object, given view object.

\sphinxAtStartPar
\sphinxstylestrong{Synopsis}
\begin{quote}

\sphinxAtStartPar
\sphinxcode{\sphinxupquote{MLinExpr Item(const View \&view)}}
\end{quote}

\sphinxAtStartPar
\sphinxstylestrong{Arguments}
\begin{quote}

\sphinxAtStartPar
\sphinxcode{\sphinxupquote{view}}: view of multi\sphinxhyphen{}dimensional array.
\end{quote}

\sphinxAtStartPar
\sphinxstylestrong{Return}
\begin{quote}

\sphinxAtStartPar
sub\sphinxhyphen{}arrays of MLinExpr object.
\end{quote}
\end{quote}

\subsubsection{MLinExpr::operator{[}{]}()}
\label{\detokenize{cppapi/mlinexpr:mlinexpr-operator}}\begin{quote}

\sphinxAtStartPar
Get expression of given index from MLinExpr object.

\sphinxAtStartPar
\sphinxstylestrong{Synopsis}
\begin{quote}

\sphinxAtStartPar
\sphinxcode{\sphinxupquote{MExpression \&operator{[}{]}(size\_t i)}}
\end{quote}

\sphinxAtStartPar
\sphinxstylestrong{Arguments}
\begin{quote}

\sphinxAtStartPar
\sphinxcode{\sphinxupquote{i}}: index of expression.
\end{quote}

\sphinxAtStartPar
\sphinxstylestrong{Return}
\begin{quote}

\sphinxAtStartPar
expression object.
\end{quote}
\end{quote}

\subsubsection{MLinExpr::operator{[}{]}()}
\label{\detokenize{cppapi/mlinexpr:id6}}\begin{quote}

\sphinxAtStartPar
Get constraints of given view from MLinExpr object.

\sphinxAtStartPar
\sphinxstylestrong{Synopsis}
\begin{quote}

\sphinxAtStartPar
\sphinxcode{\sphinxupquote{MLinExpr operator{[}{]}(const View \&view)}}
\end{quote}

\sphinxAtStartPar
\sphinxstylestrong{Arguments}
\begin{quote}

\sphinxAtStartPar
\sphinxcode{\sphinxupquote{view}}: view of multi\sphinxhyphen{}dimensional array.
\end{quote}

\sphinxAtStartPar
\sphinxstylestrong{Return}
\begin{quote}

\sphinxAtStartPar
new MLinExpr object.
\end{quote}
\end{quote}

\subsubsection{MLinExpr::Pick()}
\label{\detokenize{cppapi/mlinexpr:mlinexpr-pick}}\begin{quote}

\sphinxAtStartPar
Given a list of indexes, get linear expressions from MLinExpr object.

\sphinxAtStartPar
\sphinxstylestrong{Synopsis}
\begin{quote}

\sphinxAtStartPar
\sphinxcode{\sphinxupquote{MLinExpr\textless{}1\textgreater{} Pick(const NdArray\textless{}int, 1\textgreater{} \&indexes)}}
\end{quote}

\sphinxAtStartPar
\sphinxstylestrong{Arguments}
\begin{quote}

\sphinxAtStartPar
\sphinxcode{\sphinxupquote{indexes}}: indexes of elements.
\end{quote}

\sphinxAtStartPar
\sphinxstylestrong{Return}
\begin{quote}

\sphinxAtStartPar
one\sphinxhyphen{}dimensional array of desired linear expressions.
\end{quote}
\end{quote}

\subsubsection{MLinExpr::Pick()}
\label{\detokenize{cppapi/mlinexpr:id7}}\begin{quote}

\sphinxAtStartPar
Given a list of indexes, get linear expressions from MLinExpr object.

\sphinxAtStartPar
\sphinxstylestrong{Synopsis}
\begin{quote}

\sphinxAtStartPar
\sphinxcode{\sphinxupquote{MLinExpr\textless{}1\textgreater{} Pick(const NdArray\textless{}int, 2\textgreater{} \&idxrows)}}
\end{quote}

\sphinxAtStartPar
\sphinxstylestrong{Arguments}
\begin{quote}

\sphinxAtStartPar
\sphinxcode{\sphinxupquote{idxrows}}: indexes in format of 2\sphinxhyphen{}dimensional array, where each row is position of element.
\end{quote}

\sphinxAtStartPar
\sphinxstylestrong{Return}
\begin{quote}

\sphinxAtStartPar
one\sphinxhyphen{}dimensional array of desired linear expressions.
\end{quote}
\end{quote}

\subsubsection{MLinExpr::Repeat()}
\label{\detokenize{cppapi/mlinexpr:mlinexpr-repeat}}\begin{quote}

\sphinxAtStartPar
Repeat each element of MLinExpr along given axis.

\sphinxAtStartPar
\sphinxstylestrong{Synopsis}
\begin{quote}

\sphinxAtStartPar
\sphinxcode{\sphinxupquote{MLinExpr\textless{}N\textgreater{} Repeat(size\_t repeats, int axis)}}
\end{quote}

\sphinxAtStartPar
\sphinxstylestrong{Arguments}
\begin{quote}

\sphinxAtStartPar
\sphinxcode{\sphinxupquote{repeats}}: number of repetitions for each element.

\sphinxAtStartPar
\sphinxcode{\sphinxupquote{axis}}: axis of MLinExpr.
\end{quote}

\sphinxAtStartPar
\sphinxstylestrong{Return}
\begin{quote}

\sphinxAtStartPar
new MLinExpr object.
\end{quote}
\end{quote}

\subsubsection{MLinExpr::RepeatBlock()}
\label{\detokenize{cppapi/mlinexpr:mlinexpr-repeatblock}}\begin{quote}

\sphinxAtStartPar
Repeat an MLinExpr a number of times along given axis.

\sphinxAtStartPar
\sphinxstylestrong{Synopsis}
\begin{quote}

\sphinxAtStartPar
\sphinxcode{\sphinxupquote{MLinExpr\textless{}N\textgreater{} RepeatBlock(size\_t repeats, int axis)}}
\end{quote}

\sphinxAtStartPar
\sphinxstylestrong{Arguments}
\begin{quote}

\sphinxAtStartPar
\sphinxcode{\sphinxupquote{repeats}}: number of repetitions.

\sphinxAtStartPar
\sphinxcode{\sphinxupquote{axis}}: axis of MLinExpr.
\end{quote}

\sphinxAtStartPar
\sphinxstylestrong{Return}
\begin{quote}

\sphinxAtStartPar
new MLinExpr object.
\end{quote}
\end{quote}

\subsubsection{MLinExpr::Represent()}
\label{\detokenize{cppapi/mlinexpr:mlinexpr-represent}}\begin{quote}

\sphinxAtStartPar
String representation of MLinExpr object.

\sphinxAtStartPar
\sphinxstylestrong{Synopsis}
\begin{quote}

\sphinxAtStartPar
\sphinxcode{\sphinxupquote{std::string Represent(size\_t maxlen)}}
\end{quote}

\sphinxAtStartPar
\sphinxstylestrong{Arguments}
\begin{quote}

\sphinxAtStartPar
\sphinxcode{\sphinxupquote{maxlen}}: max length of representation.
\end{quote}

\sphinxAtStartPar
\sphinxstylestrong{Return}
\begin{quote}

\sphinxAtStartPar
string object.
\end{quote}
\end{quote}

\subsubsection{MLinExpr::Reshape()}
\label{\detokenize{cppapi/mlinexpr:mlinexpr-reshape}}\begin{quote}

\sphinxAtStartPar
Reshape MLinExpr object to new shape.

\sphinxAtStartPar
\sphinxstylestrong{Synopsis}
\begin{quote}

\sphinxAtStartPar
\sphinxcode{\sphinxupquote{template \textless{}int M\textgreater{} MLinExpr\textless{}M\textgreater{} Reshape(const Shape\textless{}M\textgreater{} \&shape)}}
\end{quote}

\sphinxAtStartPar
\sphinxstylestrong{Arguments}
\begin{quote}

\sphinxAtStartPar
\sphinxcode{\sphinxupquote{shape}}: new shape of M\sphinxhyphen{}dimensions.
\end{quote}

\sphinxAtStartPar
\sphinxstylestrong{Return}
\begin{quote}

\sphinxAtStartPar
M\sphinxhyphen{}dimensional MLinExpr object.
\end{quote}
\end{quote}

\subsubsection{MLinExpr::SetItem()}
\label{\detokenize{cppapi/mlinexpr:mlinexpr-setitem}}\begin{quote}

\sphinxAtStartPar
Set expression of given index to MLinExpr object.

\sphinxAtStartPar
\sphinxstylestrong{Synopsis}
\begin{quote}

\sphinxAtStartPar
\sphinxcode{\sphinxupquote{void SetItem(size\_t idx, const MExpression \&expr)}}
\end{quote}

\sphinxAtStartPar
\sphinxstylestrong{Arguments}
\begin{quote}

\sphinxAtStartPar
\sphinxcode{\sphinxupquote{idx}}: index of element.

\sphinxAtStartPar
\sphinxcode{\sphinxupquote{expr}}: MExpression object.
\end{quote}
\end{quote}

\subsubsection{MLinExpr::Squeeze()}
\label{\detokenize{cppapi/mlinexpr:mlinexpr-squeeze}}\begin{quote}

\sphinxAtStartPar
Remove axis of length 1 from shape of MLinExpr object.

\sphinxAtStartPar
\sphinxstylestrong{Synopsis}
\begin{quote}

\sphinxAtStartPar
\sphinxcode{\sphinxupquote{MLinExpr\textless{}N \sphinxhyphen{} 1\textgreater{} Squeeze(int axis)}}
\end{quote}

\sphinxAtStartPar
\sphinxstylestrong{Arguments}
\begin{quote}

\sphinxAtStartPar
\sphinxcode{\sphinxupquote{axis}}: axis of MLinExpr, where the length is 1.
\end{quote}

\sphinxAtStartPar
\sphinxstylestrong{Return}
\begin{quote}

\sphinxAtStartPar
MLinExpr object of (N\sphinxhyphen{}1)\sphinxhyphen{}dimensional shape.
\end{quote}
\end{quote}

\subsubsection{MLinExpr::Stack()}
\label{\detokenize{cppapi/mlinexpr:mlinexpr-stack}}\begin{quote}

\sphinxAtStartPar
Stack with other MLinExpr object along given axis.

\sphinxAtStartPar
\sphinxstylestrong{Synopsis}
\begin{quote}

\sphinxAtStartPar
\sphinxcode{\sphinxupquote{MLinExpr\textless{}N\textgreater{} Stack(const MLinExpr\textless{}N\textgreater{} \&other, int axis)}}
\end{quote}

\sphinxAtStartPar
\sphinxstylestrong{Arguments}
\begin{quote}

\sphinxAtStartPar
\sphinxcode{\sphinxupquote{other}}: a MLinExpr object.

\sphinxAtStartPar
\sphinxcode{\sphinxupquote{axis}}: an axis of MLinExpr.
\end{quote}

\sphinxAtStartPar
\sphinxstylestrong{Return}
\begin{quote}

\sphinxAtStartPar
the result MLinExpr object.
\end{quote}
\end{quote}

\subsubsection{MLinExpr::Stack()}
\label{\detokenize{cppapi/mlinexpr:id8}}\begin{quote}

\sphinxAtStartPar
Stack with other MVar object along given axis.

\sphinxAtStartPar
\sphinxstylestrong{Synopsis}
\begin{quote}

\sphinxAtStartPar
\sphinxcode{\sphinxupquote{MLinExpr\textless{}N\textgreater{} Stack(const MVar\textless{}N\textgreater{} \&other, int axis)}}
\end{quote}

\sphinxAtStartPar
\sphinxstylestrong{Arguments}
\begin{quote}

\sphinxAtStartPar
\sphinxcode{\sphinxupquote{other}}: a MVar object.

\sphinxAtStartPar
\sphinxcode{\sphinxupquote{axis}}: an axis of MLinExpr.
\end{quote}

\sphinxAtStartPar
\sphinxstylestrong{Return}
\begin{quote}

\sphinxAtStartPar
the result MLinExpr object.
\end{quote}
\end{quote}

\subsubsection{MLinExpr::Stack()}
\label{\detokenize{cppapi/mlinexpr:id9}}\begin{quote}

\sphinxAtStartPar
Stack with other NdArray object along given axis.

\sphinxAtStartPar
\sphinxstylestrong{Synopsis}
\begin{quote}

\sphinxAtStartPar
\sphinxcode{\sphinxupquote{template \textless{}class T\textgreater{} MLinExpr\textless{}N\textgreater{} Stack(const NdArray\textless{}T, N\textgreater{} \&other, int axis)}}
\end{quote}

\sphinxAtStartPar
\sphinxstylestrong{Arguments}
\begin{quote}

\sphinxAtStartPar
\sphinxcode{\sphinxupquote{other}}: a NdArray object.

\sphinxAtStartPar
\sphinxcode{\sphinxupquote{axis}}: an axis of MLinExpr.
\end{quote}

\sphinxAtStartPar
\sphinxstylestrong{Return}
\begin{quote}

\sphinxAtStartPar
the result MLinExpr object.
\end{quote}
\end{quote}

\subsubsection{MLinExpr::SubConstant()}
\label{\detokenize{cppapi/mlinexpr:mlinexpr-subconstant}}\begin{quote}

\sphinxAtStartPar
Substract constants from each expression in MLinExpr object.

\sphinxAtStartPar
\sphinxstylestrong{Synopsis}
\begin{quote}

\sphinxAtStartPar
\sphinxcode{\sphinxupquote{template \textless{}class T\textgreater{} void SubConstant(const NdArray\textless{}T, N\textgreater{} \&constants)}}
\end{quote}

\sphinxAtStartPar
\sphinxstylestrong{Arguments}
\begin{quote}

\sphinxAtStartPar
\sphinxcode{\sphinxupquote{constants}}: N\sphinxhyphen{}dimension NdArray object.
\end{quote}
\end{quote}

\subsubsection{MLinExpr::Sum()}
\label{\detokenize{cppapi/mlinexpr:mlinexpr-sum}}\begin{quote}

\sphinxAtStartPar
Sum of all expressions in MLinExpr object.

\sphinxAtStartPar
\sphinxstylestrong{Synopsis}
\begin{quote}

\sphinxAtStartPar
\sphinxcode{\sphinxupquote{MLinExpr\textless{}0\textgreater{} Sum()}}
\end{quote}

\sphinxAtStartPar
\sphinxstylestrong{Return}
\begin{quote}

\sphinxAtStartPar
sum in zero dimension.
\end{quote}
\end{quote}

\subsubsection{MLinExpr::Sum()}
\label{\detokenize{cppapi/mlinexpr:id10}}\begin{quote}

\sphinxAtStartPar
Sum of expressions at given axis of MLinExpr object.

\sphinxAtStartPar
\sphinxstylestrong{Synopsis}
\begin{quote}

\sphinxAtStartPar
\sphinxcode{\sphinxupquote{MLinExpr\textless{}N \sphinxhyphen{} 1\textgreater{} Sum(int axis)}}
\end{quote}

\sphinxAtStartPar
\sphinxstylestrong{Arguments}
\begin{quote}

\sphinxAtStartPar
\sphinxcode{\sphinxupquote{axis}}: axis of MLinExpr
\end{quote}

\sphinxAtStartPar
\sphinxstylestrong{Return}
\begin{quote}

\sphinxAtStartPar
MLinExpr object in (N\sphinxhyphen{}1)\sphinxhyphen{}dimension.
\end{quote}
\end{quote}

\subsubsection{MLinExpr::Transpose()}
\label{\detokenize{cppapi/mlinexpr:mlinexpr-transpose}}\begin{quote}

\sphinxAtStartPar
Perform matrix transpose of MLinExpr object.

\sphinxAtStartPar
\sphinxstylestrong{Synopsis}
\begin{quote}

\sphinxAtStartPar
\sphinxcode{\sphinxupquote{MLinExpr\textless{}N\textgreater{} Transpose()}}
\end{quote}

\sphinxAtStartPar
\sphinxstylestrong{Return}
\begin{quote}

\sphinxAtStartPar
transposed MLinExpr object.
\end{quote}
\end{quote}

\subsection{MPsdConstr Class}
\label{\detokenize{cppapiref:mpsdconstr-class}}\label{\detokenize{cppapiref:chapcppapiref-mpsdconstr}}
\sphinxAtStartPar
The \sphinxtitleref{MPsdConstr} class in COPT represents multi\sphinxhyphen{}dimensional semidefinite constraints.
It is generated through the methods \sphinxcode{\sphinxupquote{addConstrs}} or \sphinxcode{\sphinxupquote{addConstr}} of
{\hyperref[\detokenize{cppapiref:chapcppapiref-model}]{\sphinxcrossref{\DUrole{std,std-ref}{Model}}}}. The following member methods are provided:

\sphinxstepscope

\subsubsection{MPsdConstr::MPsdConstr()}
\label{\detokenize{cppapi/mpsdconstr:mpsdconstr-mpsdconstr}}\label{\detokenize{cppapi/mpsdconstr::doc}}\begin{quote}

\sphinxAtStartPar
Construct a MPsdConstr object with the given shape, filling with the given PSD constraint.

\sphinxAtStartPar
\sphinxstylestrong{Synopsis}
\begin{quote}

\sphinxAtStartPar
\sphinxcode{\sphinxupquote{MPsdConstr(const Shape\textless{}N\textgreater{} \&shp, const PsdConstraint \&con)}}
\end{quote}

\sphinxAtStartPar
\sphinxstylestrong{Arguments}
\begin{quote}

\sphinxAtStartPar
\sphinxcode{\sphinxupquote{shp}}: shape of MPsdConstr.

\sphinxAtStartPar
\sphinxcode{\sphinxupquote{con}}: PSD constraint object.
\end{quote}
\end{quote}

\subsubsection{MPsdConstr::MPsdConstr()}
\label{\detokenize{cppapi/mpsdconstr:id1}}\begin{quote}

\sphinxAtStartPar
Construct a MPsdConstr object with the given shape, filling with an array of PSD constraints.

\sphinxAtStartPar
\sphinxstylestrong{Synopsis}
\begin{quote}

\sphinxAtStartPar
\sphinxcode{\sphinxupquote{MPsdConstr(const Shape\textless{}N\textgreater{} \&shp, const PsdConstrArray \&cons)}}
\end{quote}

\sphinxAtStartPar
\sphinxstylestrong{Arguments}
\begin{quote}

\sphinxAtStartPar
\sphinxcode{\sphinxupquote{shp}}: shape of MPsdConstr.

\sphinxAtStartPar
\sphinxcode{\sphinxupquote{cons}}: an array of PSD constraints.
\end{quote}
\end{quote}

\subsubsection{MPsdConstr::Clone()}
\label{\detokenize{cppapi/mpsdconstr:mpsdconstr-clone}}\begin{quote}

\sphinxAtStartPar
Clone MPsdConstr object.

\sphinxAtStartPar
\sphinxstylestrong{Synopsis}
\begin{quote}

\sphinxAtStartPar
\sphinxcode{\sphinxupquote{MPsdConstr Clone()}}
\end{quote}

\sphinxAtStartPar
\sphinxstylestrong{Return}
\begin{quote}

\sphinxAtStartPar
new MPsdConstr object.
\end{quote}
\end{quote}

\subsubsection{MPsdConstr::Diagonal()}
\label{\detokenize{cppapi/mpsdconstr:mpsdconstr-diagonal}}\begin{quote}

\sphinxAtStartPar
Get diagonals of MPsdConstr object.

\sphinxAtStartPar
\sphinxstylestrong{Synopsis}
\begin{quote}

\sphinxAtStartPar
\sphinxcode{\sphinxupquote{MPsdConstr\textless{}N \sphinxhyphen{} 1\textgreater{} Diagonal(}}
\begin{quote}

\sphinxAtStartPar
\sphinxcode{\sphinxupquote{int offset,}}

\sphinxAtStartPar
\sphinxcode{\sphinxupquote{int axis1,}}

\sphinxAtStartPar
\sphinxcode{\sphinxupquote{int axis2)}}
\end{quote}
\end{quote}

\sphinxAtStartPar
\sphinxstylestrong{Arguments}
\begin{quote}

\sphinxAtStartPar
\sphinxcode{\sphinxupquote{offset}}: offset of the diagonal from the main diagonal. Can be positive or negative.

\sphinxAtStartPar
\sphinxcode{\sphinxupquote{axis1}}: 1st axis of MPsdConstr.

\sphinxAtStartPar
\sphinxcode{\sphinxupquote{axis2}}: 2nd axis of MPsdConstr.
\end{quote}

\sphinxAtStartPar
\sphinxstylestrong{Return}
\begin{quote}

\sphinxAtStartPar
(N\sphinxhyphen{}1)\sphinxhyphen{}dimensional diagonals.
\end{quote}
\end{quote}

\subsubsection{MPsdConstr::Expand()}
\label{\detokenize{cppapi/mpsdconstr:mpsdconstr-expand}}\begin{quote}

\sphinxAtStartPar
Expand shape of MPsdConstr object.

\sphinxAtStartPar
\sphinxstylestrong{Synopsis}
\begin{quote}

\sphinxAtStartPar
\sphinxcode{\sphinxupquote{MPsdConstr\textless{}N + 1\textgreater{} Expand(int axis)}}
\end{quote}

\sphinxAtStartPar
\sphinxstylestrong{Arguments}
\begin{quote}

\sphinxAtStartPar
\sphinxcode{\sphinxupquote{axis}}: axis of MPsdConstr.
\end{quote}

\sphinxAtStartPar
\sphinxstylestrong{Return}
\begin{quote}

\sphinxAtStartPar
MPsdConstr object of (N+1)\sphinxhyphen{}dimensional shape.
\end{quote}
\end{quote}

\subsubsection{MPsdConstr::Flatten()}
\label{\detokenize{cppapi/mpsdconstr:mpsdconstr-flatten}}\begin{quote}

\sphinxAtStartPar
Flatten a MPsdConstr object to a 1\sphinxhyphen{}dimensional shape.

\sphinxAtStartPar
\sphinxstylestrong{Synopsis}
\begin{quote}

\sphinxAtStartPar
\sphinxcode{\sphinxupquote{MPsdConstr\textless{}1\textgreater{} Flatten()}}
\end{quote}

\sphinxAtStartPar
\sphinxstylestrong{Return}
\begin{quote}

\sphinxAtStartPar
a MPsdConstr object collapsed into one dimension.
\end{quote}
\end{quote}

\subsubsection{MPsdConstr::Get()}
\label{\detokenize{cppapi/mpsdconstr:mpsdconstr-get}}\begin{quote}

\sphinxAtStartPar
Get values of information associated with PSD constraints in MPsdConstr object.

\sphinxAtStartPar
\sphinxstylestrong{Synopsis}
\begin{quote}

\sphinxAtStartPar
\sphinxcode{\sphinxupquote{NdArray\textless{}double, N\textgreater{} Get(const char *szInfo)}}
\end{quote}

\sphinxAtStartPar
\sphinxstylestrong{Arguments}
\begin{quote}

\sphinxAtStartPar
\sphinxcode{\sphinxupquote{szInfo}}: name of information.
\end{quote}

\sphinxAtStartPar
\sphinxstylestrong{Return}
\begin{quote}

\sphinxAtStartPar
multi\sphinxhyphen{}dimensional array of information of PSD constraints.
\end{quote}
\end{quote}

\subsubsection{MPsdConstr::GetDim()}
\label{\detokenize{cppapi/mpsdconstr:mpsdconstr-getdim}}\begin{quote}

\sphinxAtStartPar
Get i\sphinxhyphen{}th dimension of MPsdConstr object.

\sphinxAtStartPar
\sphinxstylestrong{Synopsis}
\begin{quote}

\sphinxAtStartPar
\sphinxcode{\sphinxupquote{size\_t GetDim(int i)}}
\end{quote}

\sphinxAtStartPar
\sphinxstylestrong{Arguments}
\begin{quote}

\sphinxAtStartPar
\sphinxcode{\sphinxupquote{i}}: index of dimension
\end{quote}

\sphinxAtStartPar
\sphinxstylestrong{Return}
\begin{quote}

\sphinxAtStartPar
i\sphinxhyphen{}th dimension.
\end{quote}
\end{quote}

\subsubsection{MPsdConstr::GetIdx()}
\label{\detokenize{cppapi/mpsdconstr:mpsdconstr-getidx}}\begin{quote}

\sphinxAtStartPar
Get index of PSD constraints in MPsdConstr object.

\sphinxAtStartPar
\sphinxstylestrong{Synopsis}
\begin{quote}

\sphinxAtStartPar
\sphinxcode{\sphinxupquote{NdArray\textless{}int, N\textgreater{} GetIdx()}}
\end{quote}

\sphinxAtStartPar
\sphinxstylestrong{Return}
\begin{quote}

\sphinxAtStartPar
multi\sphinxhyphen{}dimensional array of indexes of PSD constraints.
\end{quote}
\end{quote}

\subsubsection{MPsdConstr::GetND()}
\label{\detokenize{cppapi/mpsdconstr:mpsdconstr-getnd}}\begin{quote}

\sphinxAtStartPar
Get number of dimensions of MPsdConstr object.

\sphinxAtStartPar
\sphinxstylestrong{Synopsis}
\begin{quote}

\sphinxAtStartPar
\sphinxcode{\sphinxupquote{int GetND()}}
\end{quote}

\sphinxAtStartPar
\sphinxstylestrong{Return}
\begin{quote}

\sphinxAtStartPar
number of dimensions.
\end{quote}
\end{quote}

\subsubsection{MPsdConstr::GetShape()}
\label{\detokenize{cppapi/mpsdconstr:mpsdconstr-getshape}}\begin{quote}

\sphinxAtStartPar
Get shape of MPsdConstr object.

\sphinxAtStartPar
\sphinxstylestrong{Synopsis}
\begin{quote}

\sphinxAtStartPar
\sphinxcode{\sphinxupquote{Shape\textless{}N\textgreater{} GetShape()}}
\end{quote}

\sphinxAtStartPar
\sphinxstylestrong{Return}
\begin{quote}

\sphinxAtStartPar
shape object.
\end{quote}
\end{quote}

\subsubsection{MPsdConstr::GetSize()}
\label{\detokenize{cppapi/mpsdconstr:mpsdconstr-getsize}}\begin{quote}

\sphinxAtStartPar
Get size of MPsdConstr object.

\sphinxAtStartPar
\sphinxstylestrong{Synopsis}
\begin{quote}

\sphinxAtStartPar
\sphinxcode{\sphinxupquote{size\_t GetSize()}}
\end{quote}

\sphinxAtStartPar
\sphinxstylestrong{Return}
\begin{quote}

\sphinxAtStartPar
number of QConstraints
\end{quote}
\end{quote}

\subsubsection{MPsdConstr::Item()}
\label{\detokenize{cppapi/mpsdconstr:mpsdconstr-item}}\begin{quote}

\sphinxAtStartPar
Get PSD constraint of given index from MPsdConstr object.

\sphinxAtStartPar
\sphinxstylestrong{Synopsis}
\begin{quote}

\sphinxAtStartPar
\sphinxcode{\sphinxupquote{PsdConstraint \&Item(size\_t idx)}}
\end{quote}

\sphinxAtStartPar
\sphinxstylestrong{Arguments}
\begin{quote}

\sphinxAtStartPar
\sphinxcode{\sphinxupquote{idx}}: index of PSD constraint.
\end{quote}

\sphinxAtStartPar
\sphinxstylestrong{Return}
\begin{quote}

\sphinxAtStartPar
PsdConstraint object.
\end{quote}
\end{quote}

\subsubsection{MPsdConstr::operator{[}{]}()}
\label{\detokenize{cppapi/mpsdconstr:mpsdconstr-operator}}\begin{quote}

\sphinxAtStartPar
Get PSD constraint of given index from MPsdConstr object.

\sphinxAtStartPar
\sphinxstylestrong{Synopsis}
\begin{quote}

\sphinxAtStartPar
\sphinxcode{\sphinxupquote{PsdConstraint \&operator{[}{]}(size\_t idx)}}
\end{quote}

\sphinxAtStartPar
\sphinxstylestrong{Arguments}
\begin{quote}

\sphinxAtStartPar
\sphinxcode{\sphinxupquote{idx}}: index of PSD constraint.
\end{quote}

\sphinxAtStartPar
\sphinxstylestrong{Return}
\begin{quote}

\sphinxAtStartPar
PsdConstraint object.
\end{quote}
\end{quote}

\subsubsection{MPsdConstr::operator{[}{]}()}
\label{\detokenize{cppapi/mpsdconstr:id2}}\begin{quote}

\sphinxAtStartPar
Get constraints of given view from MPsdConstr object.

\sphinxAtStartPar
\sphinxstylestrong{Synopsis}
\begin{quote}

\sphinxAtStartPar
\sphinxcode{\sphinxupquote{MPsdConstr operator{[}{]}(const View \&view)}}
\end{quote}

\sphinxAtStartPar
\sphinxstylestrong{Arguments}
\begin{quote}

\sphinxAtStartPar
\sphinxcode{\sphinxupquote{view}}: view of multi\sphinxhyphen{}dimensional array.
\end{quote}

\sphinxAtStartPar
\sphinxstylestrong{Return}
\begin{quote}

\sphinxAtStartPar
new MPsdConstr object.
\end{quote}
\end{quote}

\subsubsection{MPsdConstr::Pick()}
\label{\detokenize{cppapi/mpsdconstr:mpsdconstr-pick}}\begin{quote}

\sphinxAtStartPar
Given a list of indexes, get PSD constraints from MPsdConstr object.

\sphinxAtStartPar
\sphinxstylestrong{Synopsis}
\begin{quote}

\sphinxAtStartPar
\sphinxcode{\sphinxupquote{MPsdConstr\textless{}1\textgreater{} Pick(const NdArray\textless{}int, 1\textgreater{} \&indexes)}}
\end{quote}

\sphinxAtStartPar
\sphinxstylestrong{Arguments}
\begin{quote}

\sphinxAtStartPar
\sphinxcode{\sphinxupquote{indexes}}: indexes of elements.
\end{quote}

\sphinxAtStartPar
\sphinxstylestrong{Return}
\begin{quote}

\sphinxAtStartPar
one\sphinxhyphen{}dimensional array of desired PSD constraints.
\end{quote}
\end{quote}

\subsubsection{MPsdConstr::Pick()}
\label{\detokenize{cppapi/mpsdconstr:id3}}\begin{quote}

\sphinxAtStartPar
Given a list of indexes, get PSD constraints from MPsdConstr object.

\sphinxAtStartPar
\sphinxstylestrong{Synopsis}
\begin{quote}

\sphinxAtStartPar
\sphinxcode{\sphinxupquote{MPsdConstr\textless{}1\textgreater{} Pick(const NdArray\textless{}int, 2\textgreater{} \&idxrows)}}
\end{quote}

\sphinxAtStartPar
\sphinxstylestrong{Arguments}
\begin{quote}

\sphinxAtStartPar
\sphinxcode{\sphinxupquote{idxrows}}: indexes in format of 2\sphinxhyphen{}dimensional array, where each row is position of element.
\end{quote}

\sphinxAtStartPar
\sphinxstylestrong{Return}
\begin{quote}

\sphinxAtStartPar
one\sphinxhyphen{}dimensional array of desired PSD constraints.
\end{quote}
\end{quote}

\subsubsection{MPsdConstr::Represent()}
\label{\detokenize{cppapi/mpsdconstr:mpsdconstr-represent}}\begin{quote}

\sphinxAtStartPar
String representation of MPsdConstr object.

\sphinxAtStartPar
\sphinxstylestrong{Synopsis}
\begin{quote}

\sphinxAtStartPar
\sphinxcode{\sphinxupquote{std::string Represent(size\_t maxlen)}}
\end{quote}

\sphinxAtStartPar
\sphinxstylestrong{Arguments}
\begin{quote}

\sphinxAtStartPar
\sphinxcode{\sphinxupquote{maxlen}}: max length of representation.
\end{quote}

\sphinxAtStartPar
\sphinxstylestrong{Return}
\begin{quote}

\sphinxAtStartPar
string object.
\end{quote}
\end{quote}

\subsubsection{MPsdConstr::Reshape()}
\label{\detokenize{cppapi/mpsdconstr:mpsdconstr-reshape}}\begin{quote}

\sphinxAtStartPar
Reshape MPsdConstr object to new shape.

\sphinxAtStartPar
\sphinxstylestrong{Synopsis}
\begin{quote}

\sphinxAtStartPar
\sphinxcode{\sphinxupquote{template \textless{}int M\textgreater{} MPsdConstr\textless{}M\textgreater{} Reshape(const Shape\textless{}M\textgreater{} \&shape)}}
\end{quote}

\sphinxAtStartPar
\sphinxstylestrong{Arguments}
\begin{quote}

\sphinxAtStartPar
\sphinxcode{\sphinxupquote{shape}}: new shape of M\sphinxhyphen{}dimensions.
\end{quote}

\sphinxAtStartPar
\sphinxstylestrong{Return}
\begin{quote}

\sphinxAtStartPar
M\sphinxhyphen{}dimensional MPsdConstr object.
\end{quote}
\end{quote}

\subsubsection{MPsdConstr::Set()}
\label{\detokenize{cppapi/mpsdconstr:mpsdconstr-set}}\begin{quote}

\sphinxAtStartPar
Set values of information associated with PSD constraints in MPsdConstr object.

\sphinxAtStartPar
\sphinxstylestrong{Synopsis}
\begin{quote}

\sphinxAtStartPar
\sphinxcode{\sphinxupquote{void Set(const char *szInfo, const NdArray\textless{}double, N\textgreater{} \&vals)}}
\end{quote}

\sphinxAtStartPar
\sphinxstylestrong{Arguments}
\begin{quote}

\sphinxAtStartPar
\sphinxcode{\sphinxupquote{szInfo}}: name of information.

\sphinxAtStartPar
\sphinxcode{\sphinxupquote{vals}}: multi\sphinxhyphen{}dimensional array of values of information.
\end{quote}
\end{quote}

\subsubsection{MPsdConstr::Set()}
\label{\detokenize{cppapi/mpsdconstr:id4}}\begin{quote}

\sphinxAtStartPar
Set values of information associated with PSD constraints in MPsdConstr object.

\sphinxAtStartPar
\sphinxstylestrong{Synopsis}
\begin{quote}

\sphinxAtStartPar
\sphinxcode{\sphinxupquote{void Set(const char *szInfo, double val)}}
\end{quote}

\sphinxAtStartPar
\sphinxstylestrong{Arguments}
\begin{quote}

\sphinxAtStartPar
\sphinxcode{\sphinxupquote{szInfo}}: name of information.

\sphinxAtStartPar
\sphinxcode{\sphinxupquote{val}}: value of information.
\end{quote}
\end{quote}

\subsubsection{MPsdConstr::SetItem()}
\label{\detokenize{cppapi/mpsdconstr:mpsdconstr-setitem}}\begin{quote}

\sphinxAtStartPar
Set PSD constraint of given index to MPsdConstr object.

\sphinxAtStartPar
\sphinxstylestrong{Synopsis}
\begin{quote}

\sphinxAtStartPar
\sphinxcode{\sphinxupquote{void SetItem(size\_t idx, const PsdConstraint \&con)}}
\end{quote}

\sphinxAtStartPar
\sphinxstylestrong{Arguments}
\begin{quote}

\sphinxAtStartPar
\sphinxcode{\sphinxupquote{idx}}: index of element.

\sphinxAtStartPar
\sphinxcode{\sphinxupquote{con}}: PSD constraint object.
\end{quote}
\end{quote}

\subsubsection{MPsdConstr::Squeeze()}
\label{\detokenize{cppapi/mpsdconstr:mpsdconstr-squeeze}}\begin{quote}

\sphinxAtStartPar
Remove axis of length 1 from shape of MPsdConstr object.

\sphinxAtStartPar
\sphinxstylestrong{Synopsis}
\begin{quote}

\sphinxAtStartPar
\sphinxcode{\sphinxupquote{MPsdConstr\textless{}N \sphinxhyphen{} 1\textgreater{} Squeeze(int axis)}}
\end{quote}

\sphinxAtStartPar
\sphinxstylestrong{Arguments}
\begin{quote}

\sphinxAtStartPar
\sphinxcode{\sphinxupquote{axis}}: axis of MPsdConstr, where the length is 1.
\end{quote}

\sphinxAtStartPar
\sphinxstylestrong{Return}
\begin{quote}

\sphinxAtStartPar
MPsdConstr object of (N\sphinxhyphen{}1)\sphinxhyphen{}dimensional shape.
\end{quote}
\end{quote}

\subsubsection{MPsdConstr::Stack()}
\label{\detokenize{cppapi/mpsdconstr:mpsdconstr-stack}}\begin{quote}

\sphinxAtStartPar
Stack with other MPsdConstr object along given axis.

\sphinxAtStartPar
\sphinxstylestrong{Synopsis}
\begin{quote}

\sphinxAtStartPar
\sphinxcode{\sphinxupquote{MPsdConstr\textless{}N\textgreater{} Stack(const MPsdConstr\textless{}N\textgreater{} \&other, int axis)}}
\end{quote}

\sphinxAtStartPar
\sphinxstylestrong{Arguments}
\begin{quote}

\sphinxAtStartPar
\sphinxcode{\sphinxupquote{other}}: a MPsdConstr object.

\sphinxAtStartPar
\sphinxcode{\sphinxupquote{axis}}: an axis of MPsdConstr.
\end{quote}

\sphinxAtStartPar
\sphinxstylestrong{Return}
\begin{quote}

\sphinxAtStartPar
the result MPsdConstr object.
\end{quote}
\end{quote}

\subsubsection{MPsdConstr::Transpose()}
\label{\detokenize{cppapi/mpsdconstr:mpsdconstr-transpose}}\begin{quote}

\sphinxAtStartPar
Perform matrix transpose of MPsdConstr object.

\sphinxAtStartPar
\sphinxstylestrong{Synopsis}
\begin{quote}

\sphinxAtStartPar
\sphinxcode{\sphinxupquote{MPsdConstr\textless{}N\textgreater{} Transpose()}}
\end{quote}

\sphinxAtStartPar
\sphinxstylestrong{Return}
\begin{quote}

\sphinxAtStartPar
transposed MPsdConstr object.
\end{quote}
\end{quote}

\subsection{MPsdConstrBuilder}
\label{\detokenize{cppapiref:mpsdconstrbuilder}}\label{\detokenize{cppapiref:chapcppapiref-mpsdconstrbuilder}}
\sphinxAtStartPar
The \sphinxtitleref{MPsdConstrBuilder} class in COPT serves as a builder for multi\sphinxhyphen{}dimensional
semidefinite constraints. It is used to generate multi\sphinxhyphen{}dimensional semidefinite
constraints and supports operations with the built\sphinxhyphen{}in multi\sphinxhyphen{}dimensional array
{\hyperref[\detokenize{cppapiref:chapcppapiref-ndarray}]{\sphinxcrossref{\DUrole{std,std-ref}{NdArray}}}}.
An \sphinxtitleref{MPsdConstrBuilder} object can be created through comparison operations between
two objects, one of which can be an {\hyperref[\detokenize{cppapiref:chapcppapiref-mpsdexpr}]{\sphinxcrossref{\DUrole{std,std-ref}{MPsdExpr Class}}}} object.
The following member methods are provided:

\sphinxstepscope

\subsubsection{MPsdConstrBuilder::MPsdConstrBuilder()}
\label{\detokenize{cppapi/mpsdconstrbuilder:mpsdconstrbuilder-mpsdconstrbuilder}}\label{\detokenize{cppapi/mpsdconstrbuilder::doc}}\begin{quote}

\sphinxAtStartPar
Construct a MPsdConstrBuilder object with the given shape.

\sphinxAtStartPar
\sphinxstylestrong{Synopsis}
\begin{quote}

\sphinxAtStartPar
\sphinxcode{\sphinxupquote{MPsdConstrBuilder(const Shape\textless{}N\textgreater{} \&shp)}}
\end{quote}

\sphinxAtStartPar
\sphinxstylestrong{Arguments}
\begin{quote}

\sphinxAtStartPar
\sphinxcode{\sphinxupquote{shp}}: shape of MPsdConstrBuilder.
\end{quote}
\end{quote}

\subsubsection{MPsdConstrBuilder::Flatten()}
\label{\detokenize{cppapi/mpsdconstrbuilder:mpsdconstrbuilder-flatten}}\begin{quote}

\sphinxAtStartPar
Flatten a MPsdConstrBuilder object to a 1\sphinxhyphen{}dimensional shape.

\sphinxAtStartPar
\sphinxstylestrong{Synopsis}
\begin{quote}

\sphinxAtStartPar
\sphinxcode{\sphinxupquote{MPsdConstrBuilder\textless{}1\textgreater{} Flatten()}}
\end{quote}

\sphinxAtStartPar
\sphinxstylestrong{Return}
\begin{quote}

\sphinxAtStartPar
a MPsdConstrBuilder object collapsed into one dimension.
\end{quote}
\end{quote}

\subsubsection{MPsdConstrBuilder::GetND()}
\label{\detokenize{cppapi/mpsdconstrbuilder:mpsdconstrbuilder-getnd}}\begin{quote}

\sphinxAtStartPar
Get number of dimensions of MPsdConstrBuilder object.

\sphinxAtStartPar
\sphinxstylestrong{Synopsis}
\begin{quote}

\sphinxAtStartPar
\sphinxcode{\sphinxupquote{int GetND()}}
\end{quote}

\sphinxAtStartPar
\sphinxstylestrong{Return}
\begin{quote}

\sphinxAtStartPar
number of dimensions.
\end{quote}
\end{quote}

\subsubsection{MPsdConstrBuilder::GetPsdExpr()}
\label{\detokenize{cppapi/mpsdconstrbuilder:mpsdconstrbuilder-getpsdexpr}}\begin{quote}

\sphinxAtStartPar
Get N\sphinxhyphen{}dimensional PSD expressions associated with N\sphinxhyphen{}dimensional PSD constraints.

\sphinxAtStartPar
\sphinxstylestrong{Synopsis}
\begin{quote}

\sphinxAtStartPar
\sphinxcode{\sphinxupquote{const MPsdExpr\textless{}N\textgreater{} \&GetPsdExpr()}}
\end{quote}

\sphinxAtStartPar
\sphinxstylestrong{Return}
\begin{quote}

\sphinxAtStartPar
MPsdExpr object.
\end{quote}
\end{quote}

\subsubsection{MPsdConstrBuilder::GetRange()}
\label{\detokenize{cppapi/mpsdconstrbuilder:mpsdconstrbuilder-getrange}}\begin{quote}

\sphinxAtStartPar
Get range from lower bound to upper bound of N\sphinxhyphen{}dimensional range PSD constraints.

\sphinxAtStartPar
\sphinxstylestrong{Synopsis}
\begin{quote}

\sphinxAtStartPar
\sphinxcode{\sphinxupquote{double GetRange()}}
\end{quote}

\sphinxAtStartPar
\sphinxstylestrong{Return}
\begin{quote}

\sphinxAtStartPar
length from lower bound to upper bound of range PSD constraints.
\end{quote}
\end{quote}

\subsubsection{MPsdConstrBuilder::GetSense()}
\label{\detokenize{cppapi/mpsdconstrbuilder:mpsdconstrbuilder-getsense}}\begin{quote}

\sphinxAtStartPar
Get sense associated with N\sphinxhyphen{}dimensional PSD constraints.

\sphinxAtStartPar
\sphinxstylestrong{Synopsis}
\begin{quote}

\sphinxAtStartPar
\sphinxcode{\sphinxupquote{char GetSense()}}
\end{quote}

\sphinxAtStartPar
\sphinxstylestrong{Return}
\begin{quote}

\sphinxAtStartPar
PSD constraint sense.
\end{quote}
\end{quote}

\subsubsection{MPsdConstrBuilder::Set()}
\label{\detokenize{cppapi/mpsdconstrbuilder:mpsdconstrbuilder-set}}\begin{quote}

\sphinxAtStartPar
Set N\sphinxhyphen{}dimensional PSD constraints to its builder object.

\sphinxAtStartPar
\sphinxstylestrong{Synopsis}
\begin{quote}

\sphinxAtStartPar
\sphinxcode{\sphinxupquote{template \textless{}int M\textgreater{} void Set(}}
\begin{quote}

\sphinxAtStartPar
\sphinxcode{\sphinxupquote{const MPsdExpr\textless{}N\textgreater{} \&exprs,}}

\sphinxAtStartPar
\sphinxcode{\sphinxupquote{char sense,}}

\sphinxAtStartPar
\sphinxcode{\sphinxupquote{const MPsdExpr\textless{}M\textgreater{} \&rhs)}}
\end{quote}
\end{quote}

\sphinxAtStartPar
\sphinxstylestrong{Arguments}
\begin{quote}

\sphinxAtStartPar
\sphinxcode{\sphinxupquote{exprs}}: MPsdExpr object

\sphinxAtStartPar
\sphinxcode{\sphinxupquote{sense}}: PSD constraint sense other than COPT\_RANGE.

\sphinxAtStartPar
\sphinxcode{\sphinxupquote{rhs}}: MPsdExpr object at right side of PSD constraints.
\end{quote}
\end{quote}

\subsubsection{MPsdConstrBuilder::Set()}
\label{\detokenize{cppapi/mpsdconstrbuilder:id1}}\begin{quote}

\sphinxAtStartPar
Set N\sphinxhyphen{}dimensional PSD constraints to its builder object.

\sphinxAtStartPar
\sphinxstylestrong{Synopsis}
\begin{quote}

\sphinxAtStartPar
\sphinxcode{\sphinxupquote{void Set(}}
\begin{quote}

\sphinxAtStartPar
\sphinxcode{\sphinxupquote{const MPsdExpr\textless{}N\textgreater{} \&exprs,}}

\sphinxAtStartPar
\sphinxcode{\sphinxupquote{char sense,}}

\sphinxAtStartPar
\sphinxcode{\sphinxupquote{double rhs)}}
\end{quote}
\end{quote}

\sphinxAtStartPar
\sphinxstylestrong{Arguments}
\begin{quote}

\sphinxAtStartPar
\sphinxcode{\sphinxupquote{exprs}}: MPsdExpr object

\sphinxAtStartPar
\sphinxcode{\sphinxupquote{sense}}: PSD constraint sense other than COPT\_RANGE.

\sphinxAtStartPar
\sphinxcode{\sphinxupquote{rhs}}: constant of right side of PSD constraints.
\end{quote}
\end{quote}

\subsubsection{MPsdConstrBuilder::Set()}
\label{\detokenize{cppapi/mpsdconstrbuilder:id2}}\begin{quote}

\sphinxAtStartPar
Set N\sphinxhyphen{}dimensional PSD constraints to its builder object.

\sphinxAtStartPar
\sphinxstylestrong{Synopsis}
\begin{quote}

\sphinxAtStartPar
\sphinxcode{\sphinxupquote{template \textless{}class T\textgreater{} void Set(}}
\begin{quote}

\sphinxAtStartPar
\sphinxcode{\sphinxupquote{const MPsdExpr\textless{}N\textgreater{} \&exprs,}}

\sphinxAtStartPar
\sphinxcode{\sphinxupquote{char sense,}}

\sphinxAtStartPar
\sphinxcode{\sphinxupquote{const NdArray\textless{}T, N\textgreater{} \&rhs)}}
\end{quote}
\end{quote}

\sphinxAtStartPar
\sphinxstylestrong{Arguments}
\begin{quote}

\sphinxAtStartPar
\sphinxcode{\sphinxupquote{exprs}}: MPsdExpr object

\sphinxAtStartPar
\sphinxcode{\sphinxupquote{sense}}: PSD constraint sense other than COPT\_RANGE.

\sphinxAtStartPar
\sphinxcode{\sphinxupquote{rhs}}: N\sphinxhyphen{}dimensional constants at right side of PSD constraints.
\end{quote}
\end{quote}

\subsubsection{MPsdConstrBuilder::Set()}
\label{\detokenize{cppapi/mpsdconstrbuilder:id3}}\begin{quote}

\sphinxAtStartPar
Set N\sphinxhyphen{}dimensional PSD constraints to its builder object.

\sphinxAtStartPar
\sphinxstylestrong{Synopsis}
\begin{quote}

\sphinxAtStartPar
\sphinxcode{\sphinxupquote{template \textless{}int M\textgreater{} void Set(}}
\begin{quote}

\sphinxAtStartPar
\sphinxcode{\sphinxupquote{const MPsdExpr\textless{}N\textgreater{} \&exprs,}}

\sphinxAtStartPar
\sphinxcode{\sphinxupquote{char sense,}}

\sphinxAtStartPar
\sphinxcode{\sphinxupquote{const MVar\textless{}M\textgreater{} \&rhs)}}
\end{quote}
\end{quote}

\sphinxAtStartPar
\sphinxstylestrong{Arguments}
\begin{quote}

\sphinxAtStartPar
\sphinxcode{\sphinxupquote{exprs}}: MPsdExpr object

\sphinxAtStartPar
\sphinxcode{\sphinxupquote{sense}}: PSD constraint sense other than COPT\_RANGE.

\sphinxAtStartPar
\sphinxcode{\sphinxupquote{rhs}}: MVar object at right side of PSD constraints.
\end{quote}
\end{quote}

\subsubsection{MPsdConstrBuilder::Set()}
\label{\detokenize{cppapi/mpsdconstrbuilder:id4}}\begin{quote}

\sphinxAtStartPar
Set N\sphinxhyphen{}dimensional PSD constraints to its builder object.

\sphinxAtStartPar
\sphinxstylestrong{Synopsis}
\begin{quote}

\sphinxAtStartPar
\sphinxcode{\sphinxupquote{template \textless{}int M\textgreater{} void Set(}}
\begin{quote}

\sphinxAtStartPar
\sphinxcode{\sphinxupquote{const MPsdExpr\textless{}N\textgreater{} \&exprs,}}

\sphinxAtStartPar
\sphinxcode{\sphinxupquote{char sense,}}

\sphinxAtStartPar
\sphinxcode{\sphinxupquote{const MLinExpr\textless{}M\textgreater{} \&rhs)}}
\end{quote}
\end{quote}

\sphinxAtStartPar
\sphinxstylestrong{Arguments}
\begin{quote}

\sphinxAtStartPar
\sphinxcode{\sphinxupquote{exprs}}: MPsdExpr object

\sphinxAtStartPar
\sphinxcode{\sphinxupquote{sense}}: PSD constraint sense other than COPT\_RANGE.

\sphinxAtStartPar
\sphinxcode{\sphinxupquote{rhs}}: MLinExpr object at right side of PSD constraints.
\end{quote}
\end{quote}

\subsubsection{MPsdConstrBuilder::SetRange()}
\label{\detokenize{cppapi/mpsdconstrbuilder:mpsdconstrbuilder-setrange}}\begin{quote}

\sphinxAtStartPar
Set N\sphinxhyphen{}dimensional range PSD constraints to its builder object.

\sphinxAtStartPar
\sphinxstylestrong{Synopsis}
\begin{quote}

\sphinxAtStartPar
\sphinxcode{\sphinxupquote{void SetRange(const MPsdExpr\textless{}N\textgreater{} \&exprs, double range)}}
\end{quote}

\sphinxAtStartPar
\sphinxstylestrong{Arguments}
\begin{quote}

\sphinxAtStartPar
\sphinxcode{\sphinxupquote{exprs}}: MPsdExpr object.

\sphinxAtStartPar
\sphinxcode{\sphinxupquote{range}}: length from lower bound to upper bound of PSD constraint. Must greater than 0.
\end{quote}
\end{quote}

\subsection{MPsdExpr Class}
\label{\detokenize{cppapiref:mpsdexpr-class}}\label{\detokenize{cppapiref:chapcppapiref-mpsdexpr}}
\sphinxAtStartPar
The \sphinxtitleref{MPsdExpr} class in COPT represents multi\sphinxhyphen{}dimensional semidefinite expressions.
It is used to construct multi\sphinxhyphen{}dimensional semidefinite expressions and perform
operations with the built\sphinxhyphen{}in multi\sphinxhyphen{}dimensional array {\hyperref[\detokenize{cppapiref:chapcppapiref-ndarray}]{\sphinxcrossref{\DUrole{std,std-ref}{NdArray}}}}
in COPT.

\sphinxAtStartPar
The elements of \sphinxtitleref{MPsdExpr} are either {\hyperref[\detokenize{cppapiref:chapcppapiref-psdexpr}]{\sphinxcrossref{\DUrole{std,std-ref}{PsdExpr}}}} objects or
their multi\sphinxhyphen{}dimensional linear combinations.

\sphinxAtStartPar
The following member methods are provided:

\sphinxstepscope

\subsubsection{MPsdExpr::MPsdExpr()}
\label{\detokenize{cppapi/mpsdexpr:mpsdexpr-mpsdexpr}}\label{\detokenize{cppapi/mpsdexpr::doc}}\begin{quote}

\sphinxAtStartPar
Construct a MPsdExpr object with the given shape and a constant.

\sphinxAtStartPar
\sphinxstylestrong{Synopsis}
\begin{quote}

\sphinxAtStartPar
\sphinxcode{\sphinxupquote{MPsdExpr(const Shape\textless{}N\textgreater{} \&shp, double constant)}}
\end{quote}

\sphinxAtStartPar
\sphinxstylestrong{Arguments}
\begin{quote}

\sphinxAtStartPar
\sphinxcode{\sphinxupquote{shp}}: shape of MPsdExpr.

\sphinxAtStartPar
\sphinxcode{\sphinxupquote{constant}}: constant number, default vlaue is 0.0.
\end{quote}
\end{quote}

\subsubsection{MPsdExpr::MPsdExpr()}
\label{\detokenize{cppapi/mpsdexpr:id1}}\begin{quote}

\sphinxAtStartPar
Construct a MPsdExpr object with the given shape and a PSD expression.

\sphinxAtStartPar
\sphinxstylestrong{Synopsis}
\begin{quote}

\sphinxAtStartPar
\sphinxcode{\sphinxupquote{MPsdExpr(const Shape\textless{}N\textgreater{} \&shp, const PsdExpr \&expr)}}
\end{quote}

\sphinxAtStartPar
\sphinxstylestrong{Arguments}
\begin{quote}

\sphinxAtStartPar
\sphinxcode{\sphinxupquote{shp}}: shape of MPsdExpr.

\sphinxAtStartPar
\sphinxcode{\sphinxupquote{expr}}: a PSD expression.
\end{quote}
\end{quote}

\subsubsection{MPsdExpr::AddConstant()}
\label{\detokenize{cppapi/mpsdexpr:mpsdexpr-addconstant}}\begin{quote}

\sphinxAtStartPar
Add constant to each quadratic expression in MPsdExpr object.

\sphinxAtStartPar
\sphinxstylestrong{Synopsis}
\begin{quote}

\sphinxAtStartPar
\sphinxcode{\sphinxupquote{void AddConstant(double constant)}}
\end{quote}

\sphinxAtStartPar
\sphinxstylestrong{Arguments}
\begin{quote}

\sphinxAtStartPar
\sphinxcode{\sphinxupquote{constant}}: the value of the constant.
\end{quote}
\end{quote}

\subsubsection{MPsdExpr::AddConstant()}
\label{\detokenize{cppapi/mpsdexpr:id2}}\begin{quote}

\sphinxAtStartPar
Add constants to each PSD expression in MPsdExpr object.

\sphinxAtStartPar
\sphinxstylestrong{Synopsis}
\begin{quote}

\sphinxAtStartPar
\sphinxcode{\sphinxupquote{template \textless{}class T\textgreater{} void AddConstant(const NdArray\textless{}T, N\textgreater{} \&constants)}}
\end{quote}

\sphinxAtStartPar
\sphinxstylestrong{Arguments}
\begin{quote}

\sphinxAtStartPar
\sphinxcode{\sphinxupquote{constants}}: N\sphinxhyphen{}dimension NdArray object.
\end{quote}
\end{quote}

\subsubsection{MPsdExpr::AddLinExpr()}
\label{\detokenize{cppapi/mpsdexpr:mpsdexpr-addlinexpr}}\begin{quote}

\sphinxAtStartPar
Add a linear expression to each PsdExpr in MPsdExpr object.

\sphinxAtStartPar
\sphinxstylestrong{Synopsis}
\begin{quote}

\sphinxAtStartPar
\sphinxcode{\sphinxupquote{void AddLinExpr(const Expr \&expr, double mult)}}
\end{quote}

\sphinxAtStartPar
\sphinxstylestrong{Arguments}
\begin{quote}

\sphinxAtStartPar
\sphinxcode{\sphinxupquote{expr}}: linear expression object.

\sphinxAtStartPar
\sphinxcode{\sphinxupquote{mult}}: the multiplier of linear expression, default value is 1.0.
\end{quote}
\end{quote}

\subsubsection{MPsdExpr::AddMExpr()}
\label{\detokenize{cppapi/mpsdexpr:mpsdexpr-addmexpr}}\begin{quote}

\sphinxAtStartPar
Add MExpression to each PSD expression in MPsdExpr object.

\sphinxAtStartPar
\sphinxstylestrong{Synopsis}
\begin{quote}

\sphinxAtStartPar
\sphinxcode{\sphinxupquote{void AddMExpr(const MExpression \&expr, double mult)}}
\end{quote}

\sphinxAtStartPar
\sphinxstylestrong{Arguments}
\begin{quote}

\sphinxAtStartPar
\sphinxcode{\sphinxupquote{expr}}: MExpression object.

\sphinxAtStartPar
\sphinxcode{\sphinxupquote{mult}}: the multiplier of MExpression, default value is 1.0.
\end{quote}
\end{quote}

\subsubsection{MPsdExpr::AddMLinExpr()}
\label{\detokenize{cppapi/mpsdexpr:mpsdexpr-addmlinexpr}}\begin{quote}

\sphinxAtStartPar
Add linear expressions to MPsdExpr object.

\sphinxAtStartPar
\sphinxstylestrong{Synopsis}
\begin{quote}

\sphinxAtStartPar
\sphinxcode{\sphinxupquote{void AddMLinExpr(const MLinExpr\textless{}N\textgreater{} \&exprs, double mult)}}
\end{quote}

\sphinxAtStartPar
\sphinxstylestrong{Arguments}
\begin{quote}

\sphinxAtStartPar
\sphinxcode{\sphinxupquote{exprs}}: N\sphinxhyphen{}dimension MLinExpr object.

\sphinxAtStartPar
\sphinxcode{\sphinxupquote{mult}}: the same multiplier for added linear expressions, default value is 1.0.
\end{quote}
\end{quote}

\subsubsection{MPsdExpr::AddMPsdExpr()}
\label{\detokenize{cppapi/mpsdexpr:mpsdexpr-addmpsdexpr}}\begin{quote}

\sphinxAtStartPar
Add PSD expressions to MPsdExpr object.

\sphinxAtStartPar
\sphinxstylestrong{Synopsis}
\begin{quote}

\sphinxAtStartPar
\sphinxcode{\sphinxupquote{void AddMPsdExpr(const MPsdExpr\textless{}N\textgreater{} \&exprs, double mult)}}
\end{quote}

\sphinxAtStartPar
\sphinxstylestrong{Arguments}
\begin{quote}

\sphinxAtStartPar
\sphinxcode{\sphinxupquote{exprs}}: N\sphinxhyphen{}dimension MPsdExpr object.

\sphinxAtStartPar
\sphinxcode{\sphinxupquote{mult}}: the same multiplier for added PSD expressions, default value is 1.0.
\end{quote}
\end{quote}

\subsubsection{MPsdExpr::AddPsdExpr()}
\label{\detokenize{cppapi/mpsdexpr:mpsdexpr-addpsdexpr}}\begin{quote}

\sphinxAtStartPar
Add a PSD expression to each PSD expression in MPsdExpr object.

\sphinxAtStartPar
\sphinxstylestrong{Synopsis}
\begin{quote}

\sphinxAtStartPar
\sphinxcode{\sphinxupquote{void AddPsdExpr(const PsdExpr \&expr, double mult)}}
\end{quote}

\sphinxAtStartPar
\sphinxstylestrong{Arguments}
\begin{quote}

\sphinxAtStartPar
\sphinxcode{\sphinxupquote{expr}}: PSD expression object.

\sphinxAtStartPar
\sphinxcode{\sphinxupquote{mult}}: the multiplier of PSD expression, default value is 1.0.
\end{quote}
\end{quote}

\subsubsection{MPsdExpr::AddTerm()}
\label{\detokenize{cppapi/mpsdexpr:mpsdexpr-addterm}}\begin{quote}

\sphinxAtStartPar
Add a linear term to MPsdExpr object.

\sphinxAtStartPar
\sphinxstylestrong{Synopsis}
\begin{quote}

\sphinxAtStartPar
\sphinxcode{\sphinxupquote{void AddTerm(const Var \&var, double coeff)}}
\end{quote}

\sphinxAtStartPar
\sphinxstylestrong{Arguments}
\begin{quote}

\sphinxAtStartPar
\sphinxcode{\sphinxupquote{var}}: variable of new term.

\sphinxAtStartPar
\sphinxcode{\sphinxupquote{coeff}}: coefficient of new term.
\end{quote}
\end{quote}

\subsubsection{MPsdExpr::AddTerm()}
\label{\detokenize{cppapi/mpsdexpr:id3}}\begin{quote}

\sphinxAtStartPar
Add a PSD term to MPsdExpr object.

\sphinxAtStartPar
\sphinxstylestrong{Synopsis}
\begin{quote}

\sphinxAtStartPar
\sphinxcode{\sphinxupquote{void AddTerm(const PsdVar \&var, const SymMatrix \&mat)}}
\end{quote}

\sphinxAtStartPar
\sphinxstylestrong{Arguments}
\begin{quote}

\sphinxAtStartPar
\sphinxcode{\sphinxupquote{var}}: PSD variable of new PSD term.

\sphinxAtStartPar
\sphinxcode{\sphinxupquote{mat}}: coefficient matrix of new PSD term.
\end{quote}
\end{quote}

\subsubsection{MPsdExpr::AddTerm()}
\label{\detokenize{cppapi/mpsdexpr:id4}}\begin{quote}

\sphinxAtStartPar
Add a PSD term to MPsdExpr object.

\sphinxAtStartPar
\sphinxstylestrong{Synopsis}
\begin{quote}

\sphinxAtStartPar
\sphinxcode{\sphinxupquote{void AddTerm(const PsdVar \&var, const SymMatExpr \&expr)}}
\end{quote}

\sphinxAtStartPar
\sphinxstylestrong{Arguments}
\begin{quote}

\sphinxAtStartPar
\sphinxcode{\sphinxupquote{var}}: PSD variable of new PSD term.

\sphinxAtStartPar
\sphinxcode{\sphinxupquote{expr}}: coefficient expression of symmetric matrices of new PSD term.
\end{quote}
\end{quote}

\subsubsection{MPsdExpr::AddTerms()}
\label{\detokenize{cppapi/mpsdexpr:mpsdexpr-addterms}}\begin{quote}

\sphinxAtStartPar
Add terms to PSD expressions in MPsdExpr object.

\sphinxAtStartPar
\sphinxstylestrong{Synopsis}
\begin{quote}

\sphinxAtStartPar
\sphinxcode{\sphinxupquote{void AddTerms(const MVar\textless{}N\textgreater{} \&vars, double mult)}}
\end{quote}

\sphinxAtStartPar
\sphinxstylestrong{Arguments}
\begin{quote}

\sphinxAtStartPar
\sphinxcode{\sphinxupquote{vars}}: N\sphinxhyphen{}dimension MVar object for added terms.

\sphinxAtStartPar
\sphinxcode{\sphinxupquote{mult}}: the same coefficient for added terms, default value 1.0.
\end{quote}
\end{quote}

\subsubsection{MPsdExpr::AddTerms()}
\label{\detokenize{cppapi/mpsdexpr:id5}}\begin{quote}

\sphinxAtStartPar
Add terms to PSD expressions in MPsdExpr object.

\sphinxAtStartPar
\sphinxstylestrong{Synopsis}
\begin{quote}

\sphinxAtStartPar
\sphinxcode{\sphinxupquote{void AddTerms(const MVar\textless{}N\textgreater{} \&vars, const NdArray\textless{}double, N\textgreater{} \&coeffs)}}
\end{quote}

\sphinxAtStartPar
\sphinxstylestrong{Arguments}
\begin{quote}

\sphinxAtStartPar
\sphinxcode{\sphinxupquote{vars}}: N\sphinxhyphen{}dimension MVar object for added terms.

\sphinxAtStartPar
\sphinxcode{\sphinxupquote{coeffs}}: N\sphinxhyphen{}dimension NdArray object of coefficients for added terms.
\end{quote}
\end{quote}

\subsubsection{MPsdExpr::Clear()}
\label{\detokenize{cppapi/mpsdexpr:mpsdexpr-clear}}\begin{quote}

\sphinxAtStartPar
Clear MPsdExpr object.

\sphinxAtStartPar
\sphinxstylestrong{Synopsis}
\begin{quote}

\sphinxAtStartPar
\sphinxcode{\sphinxupquote{void Clear()}}
\end{quote}
\end{quote}

\subsubsection{MPsdExpr::Clone()}
\label{\detokenize{cppapi/mpsdexpr:mpsdexpr-clone}}\begin{quote}

\sphinxAtStartPar
Clone MPsdExpr object.

\sphinxAtStartPar
\sphinxstylestrong{Synopsis}
\begin{quote}

\sphinxAtStartPar
\sphinxcode{\sphinxupquote{MPsdExpr Clone()}}
\end{quote}

\sphinxAtStartPar
\sphinxstylestrong{Return}
\begin{quote}

\sphinxAtStartPar
new MPsdExpr object.
\end{quote}
\end{quote}

\subsubsection{MPsdExpr::Diagonal()}
\label{\detokenize{cppapi/mpsdexpr:mpsdexpr-diagonal}}\begin{quote}

\sphinxAtStartPar
Get diagonals of MPsdExpr object.

\sphinxAtStartPar
\sphinxstylestrong{Synopsis}
\begin{quote}

\sphinxAtStartPar
\sphinxcode{\sphinxupquote{MPsdExpr\textless{}N \sphinxhyphen{} 1\textgreater{} Diagonal(}}
\begin{quote}

\sphinxAtStartPar
\sphinxcode{\sphinxupquote{int offset,}}

\sphinxAtStartPar
\sphinxcode{\sphinxupquote{int axis1,}}

\sphinxAtStartPar
\sphinxcode{\sphinxupquote{int axis2)}}
\end{quote}
\end{quote}

\sphinxAtStartPar
\sphinxstylestrong{Arguments}
\begin{quote}

\sphinxAtStartPar
\sphinxcode{\sphinxupquote{offset}}: offset of the diagonal from the main diagonal. Can be positive or negative.

\sphinxAtStartPar
\sphinxcode{\sphinxupquote{axis1}}: 1st axis of MPsdExpr.

\sphinxAtStartPar
\sphinxcode{\sphinxupquote{axis2}}: 2nd axis of MPsdExpr.
\end{quote}

\sphinxAtStartPar
\sphinxstylestrong{Return}
\begin{quote}

\sphinxAtStartPar
(N\sphinxhyphen{}1)\sphinxhyphen{}dimensional diagonals.
\end{quote}
\end{quote}

\subsubsection{MPsdExpr::Evaluate()}
\label{\detokenize{cppapi/mpsdexpr:mpsdexpr-evaluate}}\begin{quote}

\sphinxAtStartPar
Evaluate MPsdExpr object after solving.

\sphinxAtStartPar
\sphinxstylestrong{Synopsis}
\begin{quote}

\sphinxAtStartPar
\sphinxcode{\sphinxupquote{NdArray\textless{}double, N\textgreater{} Evaluate()}}
\end{quote}

\sphinxAtStartPar
\sphinxstylestrong{Return}
\begin{quote}

\sphinxAtStartPar
NdArray object storing value of each PSD expression.
\end{quote}
\end{quote}

\subsubsection{MPsdExpr::Expand()}
\label{\detokenize{cppapi/mpsdexpr:mpsdexpr-expand}}\begin{quote}

\sphinxAtStartPar
Expand shape of MPsdExpr object.

\sphinxAtStartPar
\sphinxstylestrong{Synopsis}
\begin{quote}

\sphinxAtStartPar
\sphinxcode{\sphinxupquote{MPsdExpr\textless{}N + 1\textgreater{} Expand(int axis)}}
\end{quote}

\sphinxAtStartPar
\sphinxstylestrong{Arguments}
\begin{quote}

\sphinxAtStartPar
\sphinxcode{\sphinxupquote{axis}}: axis of MPsdExpr.
\end{quote}

\sphinxAtStartPar
\sphinxstylestrong{Return}
\begin{quote}

\sphinxAtStartPar
MPsdExpr object of (N+1)\sphinxhyphen{}dimensional shape.
\end{quote}
\end{quote}

\subsubsection{MPsdExpr::Flatten()}
\label{\detokenize{cppapi/mpsdexpr:mpsdexpr-flatten}}\begin{quote}

\sphinxAtStartPar
Flatten a MPsdExpr object to a 1\sphinxhyphen{}dimensional shape.

\sphinxAtStartPar
\sphinxstylestrong{Synopsis}
\begin{quote}

\sphinxAtStartPar
\sphinxcode{\sphinxupquote{MPsdExpr\textless{}1\textgreater{} Flatten()}}
\end{quote}

\sphinxAtStartPar
\sphinxstylestrong{Return}
\begin{quote}

\sphinxAtStartPar
a MPsdExpr object collapsed into one dimension.
\end{quote}
\end{quote}

\subsubsection{MPsdExpr::GetDim()}
\label{\detokenize{cppapi/mpsdexpr:mpsdexpr-getdim}}\begin{quote}

\sphinxAtStartPar
Get i\sphinxhyphen{}th dimension of MPsdExpr object.

\sphinxAtStartPar
\sphinxstylestrong{Synopsis}
\begin{quote}

\sphinxAtStartPar
\sphinxcode{\sphinxupquote{size\_t GetDim(int i)}}
\end{quote}

\sphinxAtStartPar
\sphinxstylestrong{Arguments}
\begin{quote}

\sphinxAtStartPar
\sphinxcode{\sphinxupquote{i}}: index of dimension
\end{quote}

\sphinxAtStartPar
\sphinxstylestrong{Return}
\begin{quote}

\sphinxAtStartPar
i\sphinxhyphen{}th dimension.
\end{quote}
\end{quote}

\subsubsection{MPsdExpr::GetND()}
\label{\detokenize{cppapi/mpsdexpr:mpsdexpr-getnd}}\begin{quote}

\sphinxAtStartPar
Get number of dimensions of MPsdExpr object.

\sphinxAtStartPar
\sphinxstylestrong{Synopsis}
\begin{quote}

\sphinxAtStartPar
\sphinxcode{\sphinxupquote{int GetND()}}
\end{quote}

\sphinxAtStartPar
\sphinxstylestrong{Return}
\begin{quote}

\sphinxAtStartPar
number of dimensions.
\end{quote}
\end{quote}

\subsubsection{MPsdExpr::GetShape()}
\label{\detokenize{cppapi/mpsdexpr:mpsdexpr-getshape}}\begin{quote}

\sphinxAtStartPar
Get shape of MPsdExpr object.

\sphinxAtStartPar
\sphinxstylestrong{Synopsis}
\begin{quote}

\sphinxAtStartPar
\sphinxcode{\sphinxupquote{Shape\textless{}N\textgreater{} GetShape()}}
\end{quote}

\sphinxAtStartPar
\sphinxstylestrong{Return}
\begin{quote}

\sphinxAtStartPar
shape object.
\end{quote}
\end{quote}

\subsubsection{MPsdExpr::GetSize()}
\label{\detokenize{cppapi/mpsdexpr:mpsdexpr-getsize}}\begin{quote}

\sphinxAtStartPar
Get size of MPsdExpr object.

\sphinxAtStartPar
\sphinxstylestrong{Synopsis}
\begin{quote}

\sphinxAtStartPar
\sphinxcode{\sphinxupquote{size\_t GetSize()}}
\end{quote}

\sphinxAtStartPar
\sphinxstylestrong{Return}
\begin{quote}

\sphinxAtStartPar
number of PSD expressions.
\end{quote}
\end{quote}

\subsubsection{MPsdExpr::Item()}
\label{\detokenize{cppapi/mpsdexpr:mpsdexpr-item}}\begin{quote}

\sphinxAtStartPar
Get PSD expression of given index from MPsdExpr object.

\sphinxAtStartPar
\sphinxstylestrong{Synopsis}
\begin{quote}

\sphinxAtStartPar
\sphinxcode{\sphinxupquote{PsdExpr \&Item(size\_t idx)}}
\end{quote}

\sphinxAtStartPar
\sphinxstylestrong{Arguments}
\begin{quote}

\sphinxAtStartPar
\sphinxcode{\sphinxupquote{idx}}: index of PSD expression.
\end{quote}

\sphinxAtStartPar
\sphinxstylestrong{Return}
\begin{quote}

\sphinxAtStartPar
PSD expression object.
\end{quote}
\end{quote}

\subsubsection{MPsdExpr::Item()}
\label{\detokenize{cppapi/mpsdexpr:id6}}\begin{quote}

\sphinxAtStartPar
Get sub\sphinxhyphen{}arrays of MPsdExpr object, given view object.

\sphinxAtStartPar
\sphinxstylestrong{Synopsis}
\begin{quote}

\sphinxAtStartPar
\sphinxcode{\sphinxupquote{MPsdExpr Item(const View \&view)}}
\end{quote}

\sphinxAtStartPar
\sphinxstylestrong{Arguments}
\begin{quote}

\sphinxAtStartPar
\sphinxcode{\sphinxupquote{view}}: view of multi\sphinxhyphen{}dimensional array.
\end{quote}

\sphinxAtStartPar
\sphinxstylestrong{Return}
\begin{quote}

\sphinxAtStartPar
sub\sphinxhyphen{}arrays of MPsdExpr object.
\end{quote}
\end{quote}

\subsubsection{MPsdExpr::operator{[}{]}()}
\label{\detokenize{cppapi/mpsdexpr:mpsdexpr-operator}}\begin{quote}

\sphinxAtStartPar
Get PSD expression of given index from MPsdExpr object.

\sphinxAtStartPar
\sphinxstylestrong{Synopsis}
\begin{quote}

\sphinxAtStartPar
\sphinxcode{\sphinxupquote{PsdExpr \&operator{[}{]}(size\_t i)}}
\end{quote}

\sphinxAtStartPar
\sphinxstylestrong{Arguments}
\begin{quote}

\sphinxAtStartPar
\sphinxcode{\sphinxupquote{i}}: index of PSD expression.
\end{quote}

\sphinxAtStartPar
\sphinxstylestrong{Return}
\begin{quote}

\sphinxAtStartPar
PSD expression object.
\end{quote}
\end{quote}

\subsubsection{MPsdExpr::operator{[}{]}()}
\label{\detokenize{cppapi/mpsdexpr:id7}}\begin{quote}

\sphinxAtStartPar
Get PSD expressions of given view from MPsdExpr object.

\sphinxAtStartPar
\sphinxstylestrong{Synopsis}
\begin{quote}

\sphinxAtStartPar
\sphinxcode{\sphinxupquote{MPsdExpr operator{[}{]}(const View \&view)}}
\end{quote}

\sphinxAtStartPar
\sphinxstylestrong{Arguments}
\begin{quote}

\sphinxAtStartPar
\sphinxcode{\sphinxupquote{view}}: view of multi\sphinxhyphen{}dimensional array.
\end{quote}

\sphinxAtStartPar
\sphinxstylestrong{Return}
\begin{quote}

\sphinxAtStartPar
new MPsdExpr object.
\end{quote}
\end{quote}

\subsubsection{MPsdExpr::Pick()}
\label{\detokenize{cppapi/mpsdexpr:mpsdexpr-pick}}\begin{quote}

\sphinxAtStartPar
Given a list of indexes, get PSD expressions from MPsdExpr object.

\sphinxAtStartPar
\sphinxstylestrong{Synopsis}
\begin{quote}

\sphinxAtStartPar
\sphinxcode{\sphinxupquote{MPsdExpr\textless{}1\textgreater{} Pick(const NdArray\textless{}int, 1\textgreater{} \&indexes)}}
\end{quote}

\sphinxAtStartPar
\sphinxstylestrong{Arguments}
\begin{quote}

\sphinxAtStartPar
\sphinxcode{\sphinxupquote{indexes}}: indexes of elements.
\end{quote}

\sphinxAtStartPar
\sphinxstylestrong{Return}
\begin{quote}

\sphinxAtStartPar
one\sphinxhyphen{}dimensional array of desired PSD expressions.
\end{quote}
\end{quote}

\subsubsection{MPsdExpr::Pick()}
\label{\detokenize{cppapi/mpsdexpr:id8}}\begin{quote}

\sphinxAtStartPar
Given a list of indexes, get PSD expressions from MPsdExpr object.

\sphinxAtStartPar
\sphinxstylestrong{Synopsis}
\begin{quote}

\sphinxAtStartPar
\sphinxcode{\sphinxupquote{MPsdExpr\textless{}1\textgreater{} Pick(const NdArray\textless{}int, 2\textgreater{} \&idxrows)}}
\end{quote}

\sphinxAtStartPar
\sphinxstylestrong{Arguments}
\begin{quote}

\sphinxAtStartPar
\sphinxcode{\sphinxupquote{idxrows}}: indexes in format of 2\sphinxhyphen{}dimensional array, where each row is position of element.
\end{quote}

\sphinxAtStartPar
\sphinxstylestrong{Return}
\begin{quote}

\sphinxAtStartPar
one\sphinxhyphen{}dimensional array of desired PSD expressions.
\end{quote}
\end{quote}

\subsubsection{MPsdExpr::Repeat()}
\label{\detokenize{cppapi/mpsdexpr:mpsdexpr-repeat}}\begin{quote}

\sphinxAtStartPar
Repeat each element of MPsdExpr along given axis.

\sphinxAtStartPar
\sphinxstylestrong{Synopsis}
\begin{quote}

\sphinxAtStartPar
\sphinxcode{\sphinxupquote{MPsdExpr\textless{}N\textgreater{} Repeat(size\_t repeats, int axis)}}
\end{quote}

\sphinxAtStartPar
\sphinxstylestrong{Arguments}
\begin{quote}

\sphinxAtStartPar
\sphinxcode{\sphinxupquote{repeats}}: number of repetitions for each element.

\sphinxAtStartPar
\sphinxcode{\sphinxupquote{axis}}: axis of MPsdExpr.
\end{quote}

\sphinxAtStartPar
\sphinxstylestrong{Return}
\begin{quote}

\sphinxAtStartPar
new MPsdExpr object.
\end{quote}
\end{quote}

\subsubsection{MPsdExpr::RepeatBlock()}
\label{\detokenize{cppapi/mpsdexpr:mpsdexpr-repeatblock}}\begin{quote}

\sphinxAtStartPar
Repeat an MPsdExpr a number of times along given axis.

\sphinxAtStartPar
\sphinxstylestrong{Synopsis}
\begin{quote}

\sphinxAtStartPar
\sphinxcode{\sphinxupquote{MPsdExpr\textless{}N\textgreater{} RepeatBlock(size\_t repeats, int axis)}}
\end{quote}

\sphinxAtStartPar
\sphinxstylestrong{Arguments}
\begin{quote}

\sphinxAtStartPar
\sphinxcode{\sphinxupquote{repeats}}: number of repetitions.

\sphinxAtStartPar
\sphinxcode{\sphinxupquote{axis}}: axis of MPsdExpr.
\end{quote}

\sphinxAtStartPar
\sphinxstylestrong{Return}
\begin{quote}

\sphinxAtStartPar
new MPsdExpr object.
\end{quote}
\end{quote}

\subsubsection{MPsdExpr::Represent()}
\label{\detokenize{cppapi/mpsdexpr:mpsdexpr-represent}}\begin{quote}

\sphinxAtStartPar
String representation of MPsdExpr object.

\sphinxAtStartPar
\sphinxstylestrong{Synopsis}
\begin{quote}

\sphinxAtStartPar
\sphinxcode{\sphinxupquote{std::string Represent(size\_t maxlen)}}
\end{quote}

\sphinxAtStartPar
\sphinxstylestrong{Arguments}
\begin{quote}

\sphinxAtStartPar
\sphinxcode{\sphinxupquote{maxlen}}: max length of representation.
\end{quote}

\sphinxAtStartPar
\sphinxstylestrong{Return}
\begin{quote}

\sphinxAtStartPar
string object.
\end{quote}
\end{quote}

\subsubsection{MPsdExpr::Reshape()}
\label{\detokenize{cppapi/mpsdexpr:mpsdexpr-reshape}}\begin{quote}

\sphinxAtStartPar
Reshape MPsdExpr object to new shape.

\sphinxAtStartPar
\sphinxstylestrong{Synopsis}
\begin{quote}

\sphinxAtStartPar
\sphinxcode{\sphinxupquote{template \textless{}int M\textgreater{} MPsdExpr\textless{}M\textgreater{} Reshape(const Shape\textless{}M\textgreater{} \&shape)}}
\end{quote}

\sphinxAtStartPar
\sphinxstylestrong{Arguments}
\begin{quote}

\sphinxAtStartPar
\sphinxcode{\sphinxupquote{shape}}: new shape of M\sphinxhyphen{}dimensions.
\end{quote}

\sphinxAtStartPar
\sphinxstylestrong{Return}
\begin{quote}

\sphinxAtStartPar
M\sphinxhyphen{}dimensional MPsdExpr object.
\end{quote}
\end{quote}

\subsubsection{MPsdExpr::SetItem()}
\label{\detokenize{cppapi/mpsdexpr:mpsdexpr-setitem}}\begin{quote}

\sphinxAtStartPar
Set PSD expression of given index to MPsdExpr object.

\sphinxAtStartPar
\sphinxstylestrong{Synopsis}
\begin{quote}

\sphinxAtStartPar
\sphinxcode{\sphinxupquote{void SetItem(size\_t idx, const PsdExpr \&expr)}}
\end{quote}

\sphinxAtStartPar
\sphinxstylestrong{Arguments}
\begin{quote}

\sphinxAtStartPar
\sphinxcode{\sphinxupquote{idx}}: index of element.

\sphinxAtStartPar
\sphinxcode{\sphinxupquote{expr}}: PSD expression object.
\end{quote}
\end{quote}

\subsubsection{MPsdExpr::Squeeze()}
\label{\detokenize{cppapi/mpsdexpr:mpsdexpr-squeeze}}\begin{quote}

\sphinxAtStartPar
Remove axis of length 1 from shape of MPsdExpr object.

\sphinxAtStartPar
\sphinxstylestrong{Synopsis}
\begin{quote}

\sphinxAtStartPar
\sphinxcode{\sphinxupquote{MPsdExpr\textless{}N \sphinxhyphen{} 1\textgreater{} Squeeze(int axis)}}
\end{quote}

\sphinxAtStartPar
\sphinxstylestrong{Arguments}
\begin{quote}

\sphinxAtStartPar
\sphinxcode{\sphinxupquote{axis}}: axis of MPsdExpr, where the length is 1.
\end{quote}

\sphinxAtStartPar
\sphinxstylestrong{Return}
\begin{quote}

\sphinxAtStartPar
MPsdExpr object of (N\sphinxhyphen{}1)\sphinxhyphen{}dimensional shape.
\end{quote}
\end{quote}

\subsubsection{MPsdExpr::Stack()}
\label{\detokenize{cppapi/mpsdexpr:mpsdexpr-stack}}\begin{quote}

\sphinxAtStartPar
Stack with other MPsdExpr object along given axis.

\sphinxAtStartPar
\sphinxstylestrong{Synopsis}
\begin{quote}

\sphinxAtStartPar
\sphinxcode{\sphinxupquote{MPsdExpr\textless{}N\textgreater{} Stack(const MPsdExpr\textless{}N\textgreater{} \&other, int axis)}}
\end{quote}

\sphinxAtStartPar
\sphinxstylestrong{Arguments}
\begin{quote}

\sphinxAtStartPar
\sphinxcode{\sphinxupquote{other}}: a MPsdExpr object.

\sphinxAtStartPar
\sphinxcode{\sphinxupquote{axis}}: an axis of MPsdExpr.
\end{quote}

\sphinxAtStartPar
\sphinxstylestrong{Return}
\begin{quote}

\sphinxAtStartPar
the result MPsdExpr object.
\end{quote}
\end{quote}

\subsubsection{MPsdExpr::Stack()}
\label{\detokenize{cppapi/mpsdexpr:id9}}\begin{quote}

\sphinxAtStartPar
Stack with other MLinExpr object along given axis.

\sphinxAtStartPar
\sphinxstylestrong{Synopsis}
\begin{quote}

\sphinxAtStartPar
\sphinxcode{\sphinxupquote{MPsdExpr\textless{}N\textgreater{} Stack(const MLinExpr\textless{}N\textgreater{} \&other, int axis)}}
\end{quote}

\sphinxAtStartPar
\sphinxstylestrong{Arguments}
\begin{quote}

\sphinxAtStartPar
\sphinxcode{\sphinxupquote{other}}: a MLinExpr object.

\sphinxAtStartPar
\sphinxcode{\sphinxupquote{axis}}: an axis of MPsdExpr.
\end{quote}

\sphinxAtStartPar
\sphinxstylestrong{Return}
\begin{quote}

\sphinxAtStartPar
the result MPsdExpr object.
\end{quote}
\end{quote}

\subsubsection{MPsdExpr::Stack()}
\label{\detokenize{cppapi/mpsdexpr:id10}}\begin{quote}

\sphinxAtStartPar
Stack with other MVar object along given axis.

\sphinxAtStartPar
\sphinxstylestrong{Synopsis}
\begin{quote}

\sphinxAtStartPar
\sphinxcode{\sphinxupquote{MPsdExpr\textless{}N\textgreater{} Stack(const MVar\textless{}N\textgreater{} \&other, int axis)}}
\end{quote}

\sphinxAtStartPar
\sphinxstylestrong{Arguments}
\begin{quote}

\sphinxAtStartPar
\sphinxcode{\sphinxupquote{other}}: a MVar object.

\sphinxAtStartPar
\sphinxcode{\sphinxupquote{axis}}: an axis of MPsdExpr.
\end{quote}

\sphinxAtStartPar
\sphinxstylestrong{Return}
\begin{quote}

\sphinxAtStartPar
the result MPsdExpr object.
\end{quote}
\end{quote}

\subsubsection{MPsdExpr::Stack()}
\label{\detokenize{cppapi/mpsdexpr:id11}}\begin{quote}

\sphinxAtStartPar
Stack with other NdArray object along given axis.

\sphinxAtStartPar
\sphinxstylestrong{Synopsis}
\begin{quote}

\sphinxAtStartPar
\sphinxcode{\sphinxupquote{template \textless{}class T\textgreater{} MPsdExpr\textless{}N\textgreater{} Stack(const NdArray\textless{}T, N\textgreater{} \&other, int axis)}}
\end{quote}

\sphinxAtStartPar
\sphinxstylestrong{Arguments}
\begin{quote}

\sphinxAtStartPar
\sphinxcode{\sphinxupquote{other}}: a NdArray object.

\sphinxAtStartPar
\sphinxcode{\sphinxupquote{axis}}: an axis of MPsdExpr.
\end{quote}

\sphinxAtStartPar
\sphinxstylestrong{Return}
\begin{quote}

\sphinxAtStartPar
the result MPsdExpr object.
\end{quote}
\end{quote}

\subsubsection{MPsdExpr::SubConstant()}
\label{\detokenize{cppapi/mpsdexpr:mpsdexpr-subconstant}}\begin{quote}

\sphinxAtStartPar
Substract constants from each PSD expression in MPsdExpr object.

\sphinxAtStartPar
\sphinxstylestrong{Synopsis}
\begin{quote}

\sphinxAtStartPar
\sphinxcode{\sphinxupquote{template \textless{}class T\textgreater{} void SubConstant(const NdArray\textless{}T, N\textgreater{} \&constants)}}
\end{quote}

\sphinxAtStartPar
\sphinxstylestrong{Arguments}
\begin{quote}

\sphinxAtStartPar
\sphinxcode{\sphinxupquote{constants}}: N\sphinxhyphen{}dimension NdArray object.
\end{quote}
\end{quote}

\subsubsection{MPsdExpr::Sum()}
\label{\detokenize{cppapi/mpsdexpr:mpsdexpr-sum}}\begin{quote}

\sphinxAtStartPar
Sum of all quadratic expressions in MPsdExpr object.

\sphinxAtStartPar
\sphinxstylestrong{Synopsis}
\begin{quote}

\sphinxAtStartPar
\sphinxcode{\sphinxupquote{MPsdExpr\textless{}0\textgreater{} Sum()}}
\end{quote}

\sphinxAtStartPar
\sphinxstylestrong{Return}
\begin{quote}

\sphinxAtStartPar
sum in zero dimension.
\end{quote}
\end{quote}

\subsubsection{MPsdExpr::Sum()}
\label{\detokenize{cppapi/mpsdexpr:id12}}\begin{quote}

\sphinxAtStartPar
Sum of quadratic expressions at given axis of MPsdExpr object.

\sphinxAtStartPar
\sphinxstylestrong{Synopsis}
\begin{quote}

\sphinxAtStartPar
\sphinxcode{\sphinxupquote{MPsdExpr\textless{}N \sphinxhyphen{} 1\textgreater{} Sum(int axis)}}
\end{quote}

\sphinxAtStartPar
\sphinxstylestrong{Arguments}
\begin{quote}

\sphinxAtStartPar
\sphinxcode{\sphinxupquote{axis}}: axis of MPsdExpr.
\end{quote}

\sphinxAtStartPar
\sphinxstylestrong{Return}
\begin{quote}

\sphinxAtStartPar
MPsdExpr object in (N\sphinxhyphen{}1)\sphinxhyphen{}dimension.
\end{quote}
\end{quote}

\subsubsection{MPsdExpr::Transpose()}
\label{\detokenize{cppapi/mpsdexpr:mpsdexpr-transpose}}\begin{quote}

\sphinxAtStartPar
Perform matrix transpose of MPsdExpr object.

\sphinxAtStartPar
\sphinxstylestrong{Synopsis}
\begin{quote}

\sphinxAtStartPar
\sphinxcode{\sphinxupquote{MPsdExpr\textless{}N\textgreater{} Transpose()}}
\end{quote}

\sphinxAtStartPar
\sphinxstylestrong{Return}
\begin{quote}

\sphinxAtStartPar
transposed MPsdExpr object.
\end{quote}
\end{quote}

\subsection{MQConstr}
\label{\detokenize{cppapiref:mqconstr}}\label{\detokenize{cppapiref:chapcppapiref-mqconstr}}
\sphinxAtStartPar
The MQConstr class is a COPT multi\sphinxhyphen{}dimensional quadratic constraint object.
It can be created by calling the method \sphinxcode{\sphinxupquote{addMQConstr}} of {\hyperref[\detokenize{cppapiref:chapcppapiref-model}]{\sphinxcrossref{\DUrole{std,std-ref}{Model}}}}.
The following methods are provided:

\sphinxstepscope

\subsubsection{MQConstr::MQConstr()}
\label{\detokenize{cppapi/mqconstr:mqconstr-mqconstr}}\label{\detokenize{cppapi/mqconstr::doc}}\begin{quote}

\sphinxAtStartPar
Construct a MQConstr object with the given shape, filling with the given quadratic constraint.

\sphinxAtStartPar
\sphinxstylestrong{Synopsis}
\begin{quote}

\sphinxAtStartPar
\sphinxcode{\sphinxupquote{MQConstr(const Shape\textless{}N\textgreater{} \&shp, const QConstraint \&con)}}
\end{quote}

\sphinxAtStartPar
\sphinxstylestrong{Arguments}
\begin{quote}

\sphinxAtStartPar
\sphinxcode{\sphinxupquote{shp}}: shape of MQConstr.

\sphinxAtStartPar
\sphinxcode{\sphinxupquote{con}}: Quadratic constraint object.
\end{quote}
\end{quote}

\subsubsection{MQConstr::MQConstr()}
\label{\detokenize{cppapi/mqconstr:id1}}\begin{quote}

\sphinxAtStartPar
Construct a MQConstr object with the given shape, filling with an array of quadratic constraints.

\sphinxAtStartPar
\sphinxstylestrong{Synopsis}
\begin{quote}

\sphinxAtStartPar
\sphinxcode{\sphinxupquote{MQConstr(const Shape\textless{}N\textgreater{} \&shp, const QConstrArray \&cons)}}
\end{quote}

\sphinxAtStartPar
\sphinxstylestrong{Arguments}
\begin{quote}

\sphinxAtStartPar
\sphinxcode{\sphinxupquote{shp}}: shape of MQConstr.

\sphinxAtStartPar
\sphinxcode{\sphinxupquote{cons}}: an array of quadratic constraints.
\end{quote}
\end{quote}

\subsubsection{MQConstr::Clone()}
\label{\detokenize{cppapi/mqconstr:mqconstr-clone}}\begin{quote}

\sphinxAtStartPar
Clone MQConstr object.

\sphinxAtStartPar
\sphinxstylestrong{Synopsis}
\begin{quote}

\sphinxAtStartPar
\sphinxcode{\sphinxupquote{MQConstr Clone()}}
\end{quote}

\sphinxAtStartPar
\sphinxstylestrong{Return}
\begin{quote}

\sphinxAtStartPar
new MQConstr object.
\end{quote}
\end{quote}

\subsubsection{MQConstr::Diagonal()}
\label{\detokenize{cppapi/mqconstr:mqconstr-diagonal}}\begin{quote}

\sphinxAtStartPar
Get diagonals of MQConstr object.

\sphinxAtStartPar
\sphinxstylestrong{Synopsis}
\begin{quote}

\sphinxAtStartPar
\sphinxcode{\sphinxupquote{MQConstr\textless{}N \sphinxhyphen{} 1\textgreater{} Diagonal(}}
\begin{quote}

\sphinxAtStartPar
\sphinxcode{\sphinxupquote{int offset,}}

\sphinxAtStartPar
\sphinxcode{\sphinxupquote{int axis1,}}

\sphinxAtStartPar
\sphinxcode{\sphinxupquote{int axis2)}}
\end{quote}
\end{quote}

\sphinxAtStartPar
\sphinxstylestrong{Arguments}
\begin{quote}

\sphinxAtStartPar
\sphinxcode{\sphinxupquote{offset}}: offset of the diagonal from the main diagonal. Can be positive or negative.

\sphinxAtStartPar
\sphinxcode{\sphinxupquote{axis1}}: 1st axis of MQConstr.

\sphinxAtStartPar
\sphinxcode{\sphinxupquote{axis2}}: 2nd axis of MQConstr.
\end{quote}

\sphinxAtStartPar
\sphinxstylestrong{Return}
\begin{quote}

\sphinxAtStartPar
(N\sphinxhyphen{}1)\sphinxhyphen{}dimensional diagonals.
\end{quote}
\end{quote}

\subsubsection{MQConstr::Expand()}
\label{\detokenize{cppapi/mqconstr:mqconstr-expand}}\begin{quote}

\sphinxAtStartPar
Expand shape of MQConstr object.

\sphinxAtStartPar
\sphinxstylestrong{Synopsis}
\begin{quote}

\sphinxAtStartPar
\sphinxcode{\sphinxupquote{MQConstr\textless{}N + 1\textgreater{} Expand(int axis)}}
\end{quote}

\sphinxAtStartPar
\sphinxstylestrong{Arguments}
\begin{quote}

\sphinxAtStartPar
\sphinxcode{\sphinxupquote{axis}}: axis of MQConstr.
\end{quote}

\sphinxAtStartPar
\sphinxstylestrong{Return}
\begin{quote}

\sphinxAtStartPar
MQConstr object of (N+1)\sphinxhyphen{}dimensional shape.
\end{quote}
\end{quote}

\subsubsection{MQConstr::Flatten()}
\label{\detokenize{cppapi/mqconstr:mqconstr-flatten}}\begin{quote}

\sphinxAtStartPar
Flatten a MQConstr object to a 1\sphinxhyphen{}dimensional shape.

\sphinxAtStartPar
\sphinxstylestrong{Synopsis}
\begin{quote}

\sphinxAtStartPar
\sphinxcode{\sphinxupquote{MQConstr\textless{}1\textgreater{} Flatten()}}
\end{quote}

\sphinxAtStartPar
\sphinxstylestrong{Return}
\begin{quote}

\sphinxAtStartPar
a MQConstr object collapsed into one dimension.
\end{quote}
\end{quote}

\subsubsection{MQConstr::Get()}
\label{\detokenize{cppapi/mqconstr:mqconstr-get}}\begin{quote}

\sphinxAtStartPar
Get values of information associated with quadratic constraints in MQConstr object.

\sphinxAtStartPar
\sphinxstylestrong{Synopsis}
\begin{quote}

\sphinxAtStartPar
\sphinxcode{\sphinxupquote{NdArray\textless{}double, N\textgreater{} Get(const char *szInfo)}}
\end{quote}

\sphinxAtStartPar
\sphinxstylestrong{Arguments}
\begin{quote}

\sphinxAtStartPar
\sphinxcode{\sphinxupquote{szInfo}}: name of information.
\end{quote}

\sphinxAtStartPar
\sphinxstylestrong{Return}
\begin{quote}

\sphinxAtStartPar
multi\sphinxhyphen{}dimensional array of information of quadratic constraints.
\end{quote}
\end{quote}

\subsubsection{MQConstr::GetDim()}
\label{\detokenize{cppapi/mqconstr:mqconstr-getdim}}\begin{quote}

\sphinxAtStartPar
Get i\sphinxhyphen{}th dimension of MQConstr object.

\sphinxAtStartPar
\sphinxstylestrong{Synopsis}
\begin{quote}

\sphinxAtStartPar
\sphinxcode{\sphinxupquote{size\_t GetDim(int i)}}
\end{quote}

\sphinxAtStartPar
\sphinxstylestrong{Arguments}
\begin{quote}

\sphinxAtStartPar
\sphinxcode{\sphinxupquote{i}}: index of dimension
\end{quote}

\sphinxAtStartPar
\sphinxstylestrong{Return}
\begin{quote}

\sphinxAtStartPar
i\sphinxhyphen{}th dimension.
\end{quote}
\end{quote}

\subsubsection{MQConstr::GetIdx()}
\label{\detokenize{cppapi/mqconstr:mqconstr-getidx}}\begin{quote}

\sphinxAtStartPar
Get index of quadratic constraints in MQConstr object.

\sphinxAtStartPar
\sphinxstylestrong{Synopsis}
\begin{quote}

\sphinxAtStartPar
\sphinxcode{\sphinxupquote{NdArray\textless{}int, N\textgreater{} GetIdx()}}
\end{quote}

\sphinxAtStartPar
\sphinxstylestrong{Return}
\begin{quote}

\sphinxAtStartPar
multi\sphinxhyphen{}dimensional array of indexes of quadratic constraints.
\end{quote}
\end{quote}

\subsubsection{MQConstr::GetND()}
\label{\detokenize{cppapi/mqconstr:mqconstr-getnd}}\begin{quote}

\sphinxAtStartPar
Get number of dimensions of MQConstr object.

\sphinxAtStartPar
\sphinxstylestrong{Synopsis}
\begin{quote}

\sphinxAtStartPar
\sphinxcode{\sphinxupquote{int GetND()}}
\end{quote}

\sphinxAtStartPar
\sphinxstylestrong{Return}
\begin{quote}

\sphinxAtStartPar
number of dimensions.
\end{quote}
\end{quote}

\subsubsection{MQConstr::GetRhs()}
\label{\detokenize{cppapi/mqconstr:mqconstr-getrhs}}\begin{quote}

\sphinxAtStartPar
Get RHS of quadratic constraints in MQConstr object.

\sphinxAtStartPar
\sphinxstylestrong{Synopsis}
\begin{quote}

\sphinxAtStartPar
\sphinxcode{\sphinxupquote{NdArray\textless{}double, N\textgreater{} GetRhs()}}
\end{quote}

\sphinxAtStartPar
\sphinxstylestrong{Return}
\begin{quote}

\sphinxAtStartPar
multi\sphinxhyphen{}dimensional array of RHS of quadratic constraints.
\end{quote}
\end{quote}

\subsubsection{MQConstr::GetSense()}
\label{\detokenize{cppapi/mqconstr:mqconstr-getsense}}\begin{quote}

\sphinxAtStartPar
Get senses of quadratic constraints in MQConstr object.

\sphinxAtStartPar
\sphinxstylestrong{Synopsis}
\begin{quote}

\sphinxAtStartPar
\sphinxcode{\sphinxupquote{NdArray\textless{}char, N\textgreater{} GetSense()}}
\end{quote}

\sphinxAtStartPar
\sphinxstylestrong{Return}
\begin{quote}

\sphinxAtStartPar
multi\sphinxhyphen{}dimensional array of senses of quadratic constraints.
\end{quote}
\end{quote}

\subsubsection{MQConstr::GetShape()}
\label{\detokenize{cppapi/mqconstr:mqconstr-getshape}}\begin{quote}

\sphinxAtStartPar
Get shape of MQConstr object.

\sphinxAtStartPar
\sphinxstylestrong{Synopsis}
\begin{quote}

\sphinxAtStartPar
\sphinxcode{\sphinxupquote{Shape\textless{}N\textgreater{} GetShape()}}
\end{quote}

\sphinxAtStartPar
\sphinxstylestrong{Return}
\begin{quote}

\sphinxAtStartPar
shape object.
\end{quote}
\end{quote}

\subsubsection{MQConstr::GetSize()}
\label{\detokenize{cppapi/mqconstr:mqconstr-getsize}}\begin{quote}

\sphinxAtStartPar
Get size of MQConstr object.

\sphinxAtStartPar
\sphinxstylestrong{Synopsis}
\begin{quote}

\sphinxAtStartPar
\sphinxcode{\sphinxupquote{size\_t GetSize()}}
\end{quote}

\sphinxAtStartPar
\sphinxstylestrong{Return}
\begin{quote}

\sphinxAtStartPar
number of QConstraints
\end{quote}
\end{quote}

\subsubsection{MQConstr::Item()}
\label{\detokenize{cppapi/mqconstr:mqconstr-item}}\begin{quote}

\sphinxAtStartPar
Get quadratic constraint of given index from MQConstr object.

\sphinxAtStartPar
\sphinxstylestrong{Synopsis}
\begin{quote}

\sphinxAtStartPar
\sphinxcode{\sphinxupquote{QConstraint \&Item(size\_t idx)}}
\end{quote}

\sphinxAtStartPar
\sphinxstylestrong{Arguments}
\begin{quote}

\sphinxAtStartPar
\sphinxcode{\sphinxupquote{idx}}: index of quadratic constraint.
\end{quote}

\sphinxAtStartPar
\sphinxstylestrong{Return}
\begin{quote}

\sphinxAtStartPar
QConstraint object.
\end{quote}
\end{quote}

\subsubsection{MQConstr::operator{[}{]}()}
\label{\detokenize{cppapi/mqconstr:mqconstr-operator}}\begin{quote}

\sphinxAtStartPar
Get quadratic constraint of given index from MQConstr object.

\sphinxAtStartPar
\sphinxstylestrong{Synopsis}
\begin{quote}

\sphinxAtStartPar
\sphinxcode{\sphinxupquote{QConstraint \&operator{[}{]}(size\_t idx)}}
\end{quote}

\sphinxAtStartPar
\sphinxstylestrong{Arguments}
\begin{quote}

\sphinxAtStartPar
\sphinxcode{\sphinxupquote{idx}}: index of quadratic constraint.
\end{quote}

\sphinxAtStartPar
\sphinxstylestrong{Return}
\begin{quote}

\sphinxAtStartPar
QConstraint object.
\end{quote}
\end{quote}

\subsubsection{MQConstr::operator{[}{]}()}
\label{\detokenize{cppapi/mqconstr:id2}}\begin{quote}

\sphinxAtStartPar
Get constraints of given view from MQConstr object.

\sphinxAtStartPar
\sphinxstylestrong{Synopsis}
\begin{quote}

\sphinxAtStartPar
\sphinxcode{\sphinxupquote{MQConstr operator{[}{]}(const View \&view)}}
\end{quote}

\sphinxAtStartPar
\sphinxstylestrong{Arguments}
\begin{quote}

\sphinxAtStartPar
\sphinxcode{\sphinxupquote{view}}: view of multi\sphinxhyphen{}dimensional array.
\end{quote}

\sphinxAtStartPar
\sphinxstylestrong{Return}
\begin{quote}

\sphinxAtStartPar
new MQConstr object.
\end{quote}
\end{quote}

\subsubsection{MQConstr::Pick()}
\label{\detokenize{cppapi/mqconstr:mqconstr-pick}}\begin{quote}

\sphinxAtStartPar
Given a list of indexes, get quadratic constraints from MQConstr object.

\sphinxAtStartPar
\sphinxstylestrong{Synopsis}
\begin{quote}

\sphinxAtStartPar
\sphinxcode{\sphinxupquote{MQConstr\textless{}1\textgreater{} Pick(const NdArray\textless{}int, 1\textgreater{} \&indexes)}}
\end{quote}

\sphinxAtStartPar
\sphinxstylestrong{Arguments}
\begin{quote}

\sphinxAtStartPar
\sphinxcode{\sphinxupquote{indexes}}: indexes of elements.
\end{quote}

\sphinxAtStartPar
\sphinxstylestrong{Return}
\begin{quote}

\sphinxAtStartPar
one\sphinxhyphen{}dimensional array of desired quadratic constraints.
\end{quote}
\end{quote}

\subsubsection{MQConstr::Pick()}
\label{\detokenize{cppapi/mqconstr:id3}}\begin{quote}

\sphinxAtStartPar
Given a list of indexes, get quadratic constraints from MQConstr object.

\sphinxAtStartPar
\sphinxstylestrong{Synopsis}
\begin{quote}

\sphinxAtStartPar
\sphinxcode{\sphinxupquote{MQConstr\textless{}1\textgreater{} Pick(const NdArray\textless{}int, 2\textgreater{} \&idxrows)}}
\end{quote}

\sphinxAtStartPar
\sphinxstylestrong{Arguments}
\begin{quote}

\sphinxAtStartPar
\sphinxcode{\sphinxupquote{idxrows}}: indexes in format of 2\sphinxhyphen{}dimensional array, where each row is position of element.
\end{quote}

\sphinxAtStartPar
\sphinxstylestrong{Return}
\begin{quote}

\sphinxAtStartPar
one\sphinxhyphen{}dimensional array of desired quadratic constraints.
\end{quote}
\end{quote}

\subsubsection{MQConstr::Represent()}
\label{\detokenize{cppapi/mqconstr:mqconstr-represent}}\begin{quote}

\sphinxAtStartPar
String representation of MQConstr object.

\sphinxAtStartPar
\sphinxstylestrong{Synopsis}
\begin{quote}

\sphinxAtStartPar
\sphinxcode{\sphinxupquote{std::string Represent(size\_t maxlen)}}
\end{quote}

\sphinxAtStartPar
\sphinxstylestrong{Arguments}
\begin{quote}

\sphinxAtStartPar
\sphinxcode{\sphinxupquote{maxlen}}: max length of representation.
\end{quote}

\sphinxAtStartPar
\sphinxstylestrong{Return}
\begin{quote}

\sphinxAtStartPar
string object.
\end{quote}
\end{quote}

\subsubsection{MQConstr::Reshape()}
\label{\detokenize{cppapi/mqconstr:mqconstr-reshape}}\begin{quote}

\sphinxAtStartPar
Reshape MQConstr object to new shape.

\sphinxAtStartPar
\sphinxstylestrong{Synopsis}
\begin{quote}

\sphinxAtStartPar
\sphinxcode{\sphinxupquote{template \textless{}int M\textgreater{} MQConstr\textless{}M\textgreater{} Reshape(const Shape\textless{}M\textgreater{} \&shape)}}
\end{quote}

\sphinxAtStartPar
\sphinxstylestrong{Arguments}
\begin{quote}

\sphinxAtStartPar
\sphinxcode{\sphinxupquote{shape}}: new shape of M\sphinxhyphen{}dimensions.
\end{quote}

\sphinxAtStartPar
\sphinxstylestrong{Return}
\begin{quote}

\sphinxAtStartPar
M\sphinxhyphen{}dimensional MQConstr object.
\end{quote}
\end{quote}

\subsubsection{MQConstr::Set()}
\label{\detokenize{cppapi/mqconstr:mqconstr-set}}\begin{quote}

\sphinxAtStartPar
Set values of information associated with quadratic constraints in MQConstr object.

\sphinxAtStartPar
\sphinxstylestrong{Synopsis}
\begin{quote}

\sphinxAtStartPar
\sphinxcode{\sphinxupquote{void Set(const char *szInfo, double val)}}
\end{quote}

\sphinxAtStartPar
\sphinxstylestrong{Arguments}
\begin{quote}

\sphinxAtStartPar
\sphinxcode{\sphinxupquote{szInfo}}: name of information.

\sphinxAtStartPar
\sphinxcode{\sphinxupquote{val}}: value of information.
\end{quote}
\end{quote}

\subsubsection{MQConstr::Set()}
\label{\detokenize{cppapi/mqconstr:id4}}\begin{quote}

\sphinxAtStartPar
Set values of information associated with quadratic constraints in MQConstr object.

\sphinxAtStartPar
\sphinxstylestrong{Synopsis}
\begin{quote}

\sphinxAtStartPar
\sphinxcode{\sphinxupquote{void Set(const char *szInfo, const NdArray\textless{}double, N\textgreater{} \&vals)}}
\end{quote}

\sphinxAtStartPar
\sphinxstylestrong{Arguments}
\begin{quote}

\sphinxAtStartPar
\sphinxcode{\sphinxupquote{szInfo}}: name of information.

\sphinxAtStartPar
\sphinxcode{\sphinxupquote{vals}}: multi\sphinxhyphen{}dimensional array of values of information.
\end{quote}
\end{quote}

\subsubsection{MQConstr::SetItem()}
\label{\detokenize{cppapi/mqconstr:mqconstr-setitem}}\begin{quote}

\sphinxAtStartPar
Set quadratic constraint of given index to MQConstr object.

\sphinxAtStartPar
\sphinxstylestrong{Synopsis}
\begin{quote}

\sphinxAtStartPar
\sphinxcode{\sphinxupquote{void SetItem(size\_t idx, const QConstraint \&con)}}
\end{quote}

\sphinxAtStartPar
\sphinxstylestrong{Arguments}
\begin{quote}

\sphinxAtStartPar
\sphinxcode{\sphinxupquote{idx}}: index of element.

\sphinxAtStartPar
\sphinxcode{\sphinxupquote{con}}: quadratic constraint object.
\end{quote}
\end{quote}

\subsubsection{MQConstr::Squeeze()}
\label{\detokenize{cppapi/mqconstr:mqconstr-squeeze}}\begin{quote}

\sphinxAtStartPar
Remove axis of length 1 from shape of MQConstr object.

\sphinxAtStartPar
\sphinxstylestrong{Synopsis}
\begin{quote}

\sphinxAtStartPar
\sphinxcode{\sphinxupquote{MQConstr\textless{}N \sphinxhyphen{} 1\textgreater{} Squeeze(int axis)}}
\end{quote}

\sphinxAtStartPar
\sphinxstylestrong{Arguments}
\begin{quote}

\sphinxAtStartPar
\sphinxcode{\sphinxupquote{axis}}: axis of MQConstr, where the length is 1.
\end{quote}

\sphinxAtStartPar
\sphinxstylestrong{Return}
\begin{quote}

\sphinxAtStartPar
MQConstr object of (N\sphinxhyphen{}1)\sphinxhyphen{}dimensional shape.
\end{quote}
\end{quote}

\subsubsection{MQConstr::Stack()}
\label{\detokenize{cppapi/mqconstr:mqconstr-stack}}\begin{quote}

\sphinxAtStartPar
Stack with other MQConstr object along given axis.

\sphinxAtStartPar
\sphinxstylestrong{Synopsis}
\begin{quote}

\sphinxAtStartPar
\sphinxcode{\sphinxupquote{MQConstr\textless{}N\textgreater{} Stack(const MQConstr\textless{}N\textgreater{} \&other, int axis)}}
\end{quote}

\sphinxAtStartPar
\sphinxstylestrong{Arguments}
\begin{quote}

\sphinxAtStartPar
\sphinxcode{\sphinxupquote{other}}: a MQConstr object.

\sphinxAtStartPar
\sphinxcode{\sphinxupquote{axis}}: an axis of MQConstr.
\end{quote}

\sphinxAtStartPar
\sphinxstylestrong{Return}
\begin{quote}

\sphinxAtStartPar
the result MQConstr object.
\end{quote}
\end{quote}

\subsubsection{MQConstr::Transpose()}
\label{\detokenize{cppapi/mqconstr:mqconstr-transpose}}\begin{quote}

\sphinxAtStartPar
Perform matrix transpose of MQConstr object.

\sphinxAtStartPar
\sphinxstylestrong{Synopsis}
\begin{quote}

\sphinxAtStartPar
\sphinxcode{\sphinxupquote{MQConstr\textless{}N\textgreater{} Transpose()}}
\end{quote}

\sphinxAtStartPar
\sphinxstylestrong{Return}
\begin{quote}

\sphinxAtStartPar
transposed MQConstr object.
\end{quote}
\end{quote}

\subsection{MQConstrBuilder}
\label{\detokenize{cppapiref:mqconstrbuilder}}\label{\detokenize{cppapiref:chapcppapiref-mqconstrbuilder}}
\sphinxAtStartPar
The MQConstrBuilder class is a COPT builder object of multi\sphinxhyphen{}dimensional
quadratic constraints.  It is used to generate multi\sphinxhyphen{}dimensional quadratic
constraints and supports operations with the built\sphinxhyphen{}in multi\sphinxhyphen{}dimensional array
{\hyperref[\detokenize{cppapiref:chapcppapiref-ndarray}]{\sphinxcrossref{\DUrole{std,std-ref}{NdArray}}}} in COPT. It can be created by comparing two
objects, one of which should be {\hyperref[\detokenize{cppapiref:chapcppapiref-mquadexpr}]{\sphinxcrossref{\DUrole{std,std-ref}{MQuadExpr}}}} object,
by comparison operators. The following methods are provided:

\sphinxstepscope

\subsubsection{MQConstrBuilder::MQConstrBuilder()}
\label{\detokenize{cppapi/mqconstrbuilder:mqconstrbuilder-mqconstrbuilder}}\label{\detokenize{cppapi/mqconstrbuilder::doc}}\begin{quote}

\sphinxAtStartPar
Construct a MQConstrBuilder object with the given shape.

\sphinxAtStartPar
\sphinxstylestrong{Synopsis}
\begin{quote}

\sphinxAtStartPar
\sphinxcode{\sphinxupquote{MQConstrBuilder(const Shape\textless{}N\textgreater{} \&shp)}}
\end{quote}

\sphinxAtStartPar
\sphinxstylestrong{Arguments}
\begin{quote}

\sphinxAtStartPar
\sphinxcode{\sphinxupquote{shp}}: shape of MQConstrBuilder.
\end{quote}
\end{quote}

\subsubsection{MQConstrBuilder::Flatten()}
\label{\detokenize{cppapi/mqconstrbuilder:mqconstrbuilder-flatten}}\begin{quote}

\sphinxAtStartPar
Flatten a MQConstrBuilder object to a 1\sphinxhyphen{}dimensional shape.

\sphinxAtStartPar
\sphinxstylestrong{Synopsis}
\begin{quote}

\sphinxAtStartPar
\sphinxcode{\sphinxupquote{MQConstrBuilder\textless{}1\textgreater{} Flatten()}}
\end{quote}

\sphinxAtStartPar
\sphinxstylestrong{Return}
\begin{quote}

\sphinxAtStartPar
a MQConstrBuilder object collapsed into one dimension.
\end{quote}
\end{quote}

\subsubsection{MQConstrBuilder::GetND()}
\label{\detokenize{cppapi/mqconstrbuilder:mqconstrbuilder-getnd}}\begin{quote}

\sphinxAtStartPar
Get number of dimensions of MQConstrBuilder object.

\sphinxAtStartPar
\sphinxstylestrong{Synopsis}
\begin{quote}

\sphinxAtStartPar
\sphinxcode{\sphinxupquote{int GetND()}}
\end{quote}

\sphinxAtStartPar
\sphinxstylestrong{Return}
\begin{quote}

\sphinxAtStartPar
number of dimensions.
\end{quote}
\end{quote}

\subsubsection{MQConstrBuilder::GetQuadExpr()}
\label{\detokenize{cppapi/mqconstrbuilder:mqconstrbuilder-getquadexpr}}\begin{quote}

\sphinxAtStartPar
Get N\sphinxhyphen{}dimensional quadratic expressions associated with N\sphinxhyphen{}dimensional quadratic constraints.

\sphinxAtStartPar
\sphinxstylestrong{Synopsis}
\begin{quote}

\sphinxAtStartPar
\sphinxcode{\sphinxupquote{const MQuadExpr\textless{}N\textgreater{} \&GetQuadExpr()}}
\end{quote}

\sphinxAtStartPar
\sphinxstylestrong{Return}
\begin{quote}

\sphinxAtStartPar
MQuadExpr object.
\end{quote}
\end{quote}

\subsubsection{MQConstrBuilder::GetSense()}
\label{\detokenize{cppapi/mqconstrbuilder:mqconstrbuilder-getsense}}\begin{quote}

\sphinxAtStartPar
Get sense associated with N\sphinxhyphen{}dimensional quadratic constraints.

\sphinxAtStartPar
\sphinxstylestrong{Synopsis}
\begin{quote}

\sphinxAtStartPar
\sphinxcode{\sphinxupquote{char GetSense()}}
\end{quote}

\sphinxAtStartPar
\sphinxstylestrong{Return}
\begin{quote}

\sphinxAtStartPar
quadratic constraint sense.
\end{quote}
\end{quote}

\subsubsection{MQConstrBuilder::Set()}
\label{\detokenize{cppapi/mqconstrbuilder:mqconstrbuilder-set}}\begin{quote}

\sphinxAtStartPar
Set N\sphinxhyphen{}dimensional quadratic constraints to its builder object.

\sphinxAtStartPar
\sphinxstylestrong{Synopsis}
\begin{quote}

\sphinxAtStartPar
\sphinxcode{\sphinxupquote{void Set(}}
\begin{quote}

\sphinxAtStartPar
\sphinxcode{\sphinxupquote{const MQuadExpr\textless{}N\textgreater{} \&expr,}}

\sphinxAtStartPar
\sphinxcode{\sphinxupquote{char sense,}}

\sphinxAtStartPar
\sphinxcode{\sphinxupquote{double rhs)}}
\end{quote}
\end{quote}

\sphinxAtStartPar
\sphinxstylestrong{Arguments}
\begin{quote}

\sphinxAtStartPar
\sphinxcode{\sphinxupquote{expr}}: MQuadExpr object

\sphinxAtStartPar
\sphinxcode{\sphinxupquote{sense}}: constraint sense other than COPT\_RANGE.

\sphinxAtStartPar
\sphinxcode{\sphinxupquote{rhs}}: constant of right side of quadratic constraints.
\end{quote}
\end{quote}

\subsubsection{MQConstrBuilder::Set()}
\label{\detokenize{cppapi/mqconstrbuilder:id1}}\begin{quote}

\sphinxAtStartPar
Set N\sphinxhyphen{}dimensional quadratic constraints to its builder object.

\sphinxAtStartPar
\sphinxstylestrong{Synopsis}
\begin{quote}

\sphinxAtStartPar
\sphinxcode{\sphinxupquote{template \textless{}class T\textgreater{} void Set(}}
\begin{quote}

\sphinxAtStartPar
\sphinxcode{\sphinxupquote{const MQuadExpr\textless{}N\textgreater{} \&expr,}}

\sphinxAtStartPar
\sphinxcode{\sphinxupquote{char sense,}}

\sphinxAtStartPar
\sphinxcode{\sphinxupquote{const NdArray\textless{}T, N\textgreater{} \&rhs)}}
\end{quote}
\end{quote}

\sphinxAtStartPar
\sphinxstylestrong{Arguments}
\begin{quote}

\sphinxAtStartPar
\sphinxcode{\sphinxupquote{expr}}: MQuadExpr object

\sphinxAtStartPar
\sphinxcode{\sphinxupquote{sense}}: constraint sense other than COPT\_RANGE.

\sphinxAtStartPar
\sphinxcode{\sphinxupquote{rhs}}: N\sphinxhyphen{}dimensional constants at right side of quadratic constraints.
\end{quote}
\end{quote}

\subsubsection{MQConstrBuilder::Set()}
\label{\detokenize{cppapi/mqconstrbuilder:id2}}\begin{quote}

\sphinxAtStartPar
Set N\sphinxhyphen{}dimensional quadratic constraints to its builder object.

\sphinxAtStartPar
\sphinxstylestrong{Synopsis}
\begin{quote}

\sphinxAtStartPar
\sphinxcode{\sphinxupquote{template \textless{}int M\textgreater{} void Set(}}
\begin{quote}

\sphinxAtStartPar
\sphinxcode{\sphinxupquote{const MQuadExpr\textless{}N\textgreater{} \&expr,}}

\sphinxAtStartPar
\sphinxcode{\sphinxupquote{char sense,}}

\sphinxAtStartPar
\sphinxcode{\sphinxupquote{const MVar\textless{}M\textgreater{} \&rhs)}}
\end{quote}
\end{quote}

\sphinxAtStartPar
\sphinxstylestrong{Arguments}
\begin{quote}

\sphinxAtStartPar
\sphinxcode{\sphinxupquote{expr}}: MQuadExpr object

\sphinxAtStartPar
\sphinxcode{\sphinxupquote{sense}}: constraint sense other than COPT\_RANGE.

\sphinxAtStartPar
\sphinxcode{\sphinxupquote{rhs}}: MVar object at right side of quadratic constraints.
\end{quote}
\end{quote}

\subsubsection{MQConstrBuilder::Set()}
\label{\detokenize{cppapi/mqconstrbuilder:id3}}\begin{quote}

\sphinxAtStartPar
Set N\sphinxhyphen{}dimensional quadratic constraints to its builder object.

\sphinxAtStartPar
\sphinxstylestrong{Synopsis}
\begin{quote}

\sphinxAtStartPar
\sphinxcode{\sphinxupquote{template \textless{}int M\textgreater{} void Set(}}
\begin{quote}

\sphinxAtStartPar
\sphinxcode{\sphinxupquote{const MQuadExpr\textless{}N\textgreater{} \&expr,}}

\sphinxAtStartPar
\sphinxcode{\sphinxupquote{char sense,}}

\sphinxAtStartPar
\sphinxcode{\sphinxupquote{const MLinExpr\textless{}M\textgreater{} \&rhs)}}
\end{quote}
\end{quote}

\sphinxAtStartPar
\sphinxstylestrong{Arguments}
\begin{quote}

\sphinxAtStartPar
\sphinxcode{\sphinxupquote{expr}}: MQuadExpr object

\sphinxAtStartPar
\sphinxcode{\sphinxupquote{sense}}: constraint sense other than COPT\_RANGE.

\sphinxAtStartPar
\sphinxcode{\sphinxupquote{rhs}}: MLinExpr object at right side of quadratic constraints.
\end{quote}
\end{quote}

\subsubsection{MQConstrBuilder::Set()}
\label{\detokenize{cppapi/mqconstrbuilder:id4}}\begin{quote}

\sphinxAtStartPar
Set N\sphinxhyphen{}dimensional quadratic constraints to its builder object.

\sphinxAtStartPar
\sphinxstylestrong{Synopsis}
\begin{quote}

\sphinxAtStartPar
\sphinxcode{\sphinxupquote{template \textless{}int M\textgreater{} void Set(}}
\begin{quote}

\sphinxAtStartPar
\sphinxcode{\sphinxupquote{const MQuadExpr\textless{}N\textgreater{} \&expr,}}

\sphinxAtStartPar
\sphinxcode{\sphinxupquote{char sense,}}

\sphinxAtStartPar
\sphinxcode{\sphinxupquote{const MQuadExpr\textless{}M\textgreater{} \&rhs)}}
\end{quote}
\end{quote}

\sphinxAtStartPar
\sphinxstylestrong{Arguments}
\begin{quote}

\sphinxAtStartPar
\sphinxcode{\sphinxupquote{expr}}: MQuadExpr object

\sphinxAtStartPar
\sphinxcode{\sphinxupquote{sense}}: constraint sense other than COPT\_RANGE.

\sphinxAtStartPar
\sphinxcode{\sphinxupquote{rhs}}: MQuadExpr object at right side of quadratic constraints.
\end{quote}
\end{quote}

\subsection{MQExpression}
\label{\detokenize{cppapiref:mqexpression}}\label{\detokenize{cppapiref:chapcppapiref-mqexpression}}
\sphinxAtStartPar
The MQExpression class is a generalized version of {\hyperref[\detokenize{cppapiref:chapcppapiref-quadexpr}]{\sphinxcrossref{\DUrole{std,std-ref}{QuadExpr}}}}.
It represents a quadratic expression and supports most of methods in QuadExpr class.
In addition, it supports quadratic combination of multi\sphinxhyphen{}dimensional objects,
such as {\hyperref[\detokenize{cppapiref:chapcppapiref-mvar}]{\sphinxcrossref{\DUrole{std,std-ref}{MVar}}}} object and {\hyperref[\detokenize{cppapiref:chapcppapiref-mlinexpr}]{\sphinxcrossref{\DUrole{std,std-ref}{MLinExpr}}}} object.
The following methods are provided:

\sphinxstepscope

\subsubsection{MQExpression::MQExpression()}
\label{\detokenize{cppapi/mqexpression:mqexpression-mqexpression}}\label{\detokenize{cppapi/mqexpression::doc}}\begin{quote}

\sphinxAtStartPar
Construct a MQExpression object with the given constant.

\sphinxAtStartPar
\sphinxstylestrong{Synopsis}
\begin{quote}

\sphinxAtStartPar
\sphinxcode{\sphinxupquote{MQExpression(double constant)}}
\end{quote}

\sphinxAtStartPar
\sphinxstylestrong{Arguments}
\begin{quote}

\sphinxAtStartPar
\sphinxcode{\sphinxupquote{constant}}: constant number.
\end{quote}
\end{quote}

\subsubsection{MQExpression::MQExpression()}
\label{\detokenize{cppapi/mqexpression:id1}}\begin{quote}

\sphinxAtStartPar
Construct a MQExpression object with the given quadratic expression.

\sphinxAtStartPar
\sphinxstylestrong{Synopsis}
\begin{quote}

\sphinxAtStartPar
\sphinxcode{\sphinxupquote{MQExpression(const QuadExpr \&expr)}}
\end{quote}

\sphinxAtStartPar
\sphinxstylestrong{Arguments}
\begin{quote}

\sphinxAtStartPar
\sphinxcode{\sphinxupquote{expr}}: a quadratic expression.
\end{quote}
\end{quote}

\subsubsection{MQExpression::AddConstant()}
\label{\detokenize{cppapi/mqexpression:mqexpression-addconstant}}\begin{quote}

\sphinxAtStartPar
Add constant for the expression.

\sphinxAtStartPar
\sphinxstylestrong{Synopsis}
\begin{quote}

\sphinxAtStartPar
\sphinxcode{\sphinxupquote{void AddConstant(double constant)}}
\end{quote}

\sphinxAtStartPar
\sphinxstylestrong{Arguments}
\begin{quote}

\sphinxAtStartPar
\sphinxcode{\sphinxupquote{constant}}: the value of the constant.
\end{quote}
\end{quote}

\subsubsection{MQExpression::AddExpr()}
\label{\detokenize{cppapi/mqexpression:mqexpression-addexpr}}\begin{quote}

\sphinxAtStartPar
Add a linear expression to MQExpression object.

\sphinxAtStartPar
\sphinxstylestrong{Synopsis}
\begin{quote}

\sphinxAtStartPar
\sphinxcode{\sphinxupquote{void AddExpr(const Expr \&expr, double mult)}}
\end{quote}

\sphinxAtStartPar
\sphinxstylestrong{Arguments}
\begin{quote}

\sphinxAtStartPar
\sphinxcode{\sphinxupquote{expr}}: linear expression object.

\sphinxAtStartPar
\sphinxcode{\sphinxupquote{mult}}: the multiplier of linear expression, default value is 1.0.
\end{quote}
\end{quote}

\subsubsection{MQExpression::AddMExpr()}
\label{\detokenize{cppapi/mqexpression:mqexpression-addmexpr}}\begin{quote}

\sphinxAtStartPar
Add MExpression to MQExpression object.

\sphinxAtStartPar
\sphinxstylestrong{Synopsis}
\begin{quote}

\sphinxAtStartPar
\sphinxcode{\sphinxupquote{void AddMExpr(const MExpression \&expr, double mult)}}
\end{quote}

\sphinxAtStartPar
\sphinxstylestrong{Arguments}
\begin{quote}

\sphinxAtStartPar
\sphinxcode{\sphinxupquote{expr}}: MExpression object.

\sphinxAtStartPar
\sphinxcode{\sphinxupquote{mult}}: the multiplier of MExpression, default value is 1.0.
\end{quote}
\end{quote}

\subsubsection{MQExpression::AddMQExpr()}
\label{\detokenize{cppapi/mqexpression:mqexpression-addmqexpr}}\begin{quote}

\sphinxAtStartPar
Add MQExpression to MQExpression object.

\sphinxAtStartPar
\sphinxstylestrong{Synopsis}
\begin{quote}

\sphinxAtStartPar
\sphinxcode{\sphinxupquote{void AddMQExpr(const MQExpression \&expr, double mult)}}
\end{quote}

\sphinxAtStartPar
\sphinxstylestrong{Arguments}
\begin{quote}

\sphinxAtStartPar
\sphinxcode{\sphinxupquote{expr}}: MQExpression object.

\sphinxAtStartPar
\sphinxcode{\sphinxupquote{mult}}: the multiplier of MQExpression, default value is 1.0.
\end{quote}
\end{quote}

\subsubsection{MQExpression::AddQuadExpr()}
\label{\detokenize{cppapi/mqexpression:mqexpression-addquadexpr}}\begin{quote}

\sphinxAtStartPar
Add a quadratic expression to MQExpression object.

\sphinxAtStartPar
\sphinxstylestrong{Synopsis}
\begin{quote}

\sphinxAtStartPar
\sphinxcode{\sphinxupquote{void AddQuadExpr(const MExpression \&left, const MExpression \&right)}}
\end{quote}

\sphinxAtStartPar
\sphinxstylestrong{Arguments}
\begin{quote}

\sphinxAtStartPar
\sphinxcode{\sphinxupquote{left}}: left MExpression object.

\sphinxAtStartPar
\sphinxcode{\sphinxupquote{right}}: right MExpression object.
\end{quote}
\end{quote}

\subsubsection{MQExpression::AddQuadExpr()}
\label{\detokenize{cppapi/mqexpression:id2}}\begin{quote}

\sphinxAtStartPar
Add a quadratic expression to MQExpression object.

\sphinxAtStartPar
\sphinxstylestrong{Synopsis}
\begin{quote}

\sphinxAtStartPar
\sphinxcode{\sphinxupquote{void AddQuadExpr(const QuadExpr \&expr, double mult)}}
\end{quote}

\sphinxAtStartPar
\sphinxstylestrong{Arguments}
\begin{quote}

\sphinxAtStartPar
\sphinxcode{\sphinxupquote{expr}}: quadratic expression object.

\sphinxAtStartPar
\sphinxcode{\sphinxupquote{mult}}: the multiplier of quadratic expression, default value is 1.0.
\end{quote}
\end{quote}

\subsubsection{MQExpression::AddQuadExpr()}
\label{\detokenize{cppapi/mqexpression:id3}}\begin{quote}

\sphinxAtStartPar
Add a quadratic expression to MQExpression object.

\sphinxAtStartPar
\sphinxstylestrong{Synopsis}
\begin{quote}

\sphinxAtStartPar
\sphinxcode{\sphinxupquote{void AddQuadExpr(const MExpression \&expr, const Var \&var)}}
\end{quote}

\sphinxAtStartPar
\sphinxstylestrong{Arguments}
\begin{quote}

\sphinxAtStartPar
\sphinxcode{\sphinxupquote{expr}}: MExpression object.

\sphinxAtStartPar
\sphinxcode{\sphinxupquote{var}}: Var object.
\end{quote}
\end{quote}

\subsubsection{MQExpression::AddQuadExpr()}
\label{\detokenize{cppapi/mqexpression:id4}}\begin{quote}

\sphinxAtStartPar
Add a quadratic expression to MQExpression object.

\sphinxAtStartPar
\sphinxstylestrong{Synopsis}
\begin{quote}

\sphinxAtStartPar
\sphinxcode{\sphinxupquote{void AddQuadExpr(const MExpression \&left, const Expr \&right)}}
\end{quote}

\sphinxAtStartPar
\sphinxstylestrong{Arguments}
\begin{quote}

\sphinxAtStartPar
\sphinxcode{\sphinxupquote{left}}: MExpression object.

\sphinxAtStartPar
\sphinxcode{\sphinxupquote{right}}: Expr object.
\end{quote}
\end{quote}

\subsubsection{MQExpression::AddTerm()}
\label{\detokenize{cppapi/mqexpression:mqexpression-addterm}}\begin{quote}

\sphinxAtStartPar
Add a linear term to MQExpression object.

\sphinxAtStartPar
\sphinxstylestrong{Synopsis}
\begin{quote}

\sphinxAtStartPar
\sphinxcode{\sphinxupquote{void AddTerm(const Var \&var, double coeff)}}
\end{quote}

\sphinxAtStartPar
\sphinxstylestrong{Arguments}
\begin{quote}

\sphinxAtStartPar
\sphinxcode{\sphinxupquote{var}}: variable of new term.

\sphinxAtStartPar
\sphinxcode{\sphinxupquote{coeff}}: coefficient of new term.
\end{quote}
\end{quote}

\subsubsection{MQExpression::AddTerm()}
\label{\detokenize{cppapi/mqexpression:id5}}\begin{quote}

\sphinxAtStartPar
Add a quadratic term to MQExpression object.

\sphinxAtStartPar
\sphinxstylestrong{Synopsis}
\begin{quote}

\sphinxAtStartPar
\sphinxcode{\sphinxupquote{void AddTerm(}}
\begin{quote}

\sphinxAtStartPar
\sphinxcode{\sphinxupquote{const Var \&var1,}}

\sphinxAtStartPar
\sphinxcode{\sphinxupquote{const Var \&var2,}}

\sphinxAtStartPar
\sphinxcode{\sphinxupquote{double coeff)}}
\end{quote}
\end{quote}

\sphinxAtStartPar
\sphinxstylestrong{Arguments}
\begin{quote}

\sphinxAtStartPar
\sphinxcode{\sphinxupquote{var1}}: first variable of new quadratic term.

\sphinxAtStartPar
\sphinxcode{\sphinxupquote{var2}}: second variable of new quadratic term.

\sphinxAtStartPar
\sphinxcode{\sphinxupquote{coeff}}: coefficient of new quadratic term.
\end{quote}
\end{quote}

\subsubsection{MQExpression::Clone()}
\label{\detokenize{cppapi/mqexpression:mqexpression-clone}}\begin{quote}

\sphinxAtStartPar
Clone MQExpression object.

\sphinxAtStartPar
\sphinxstylestrong{Synopsis}
\begin{quote}

\sphinxAtStartPar
\sphinxcode{\sphinxupquote{MQExpression Clone()}}
\end{quote}

\sphinxAtStartPar
\sphinxstylestrong{Return}
\begin{quote}

\sphinxAtStartPar
new MQExpression object.
\end{quote}
\end{quote}

\subsubsection{MQExpression::Evaluate()}
\label{\detokenize{cppapi/mqexpression:mqexpression-evaluate}}\begin{quote}

\sphinxAtStartPar
Evaluate MQExpression after solving.

\sphinxAtStartPar
\sphinxstylestrong{Synopsis}
\begin{quote}

\sphinxAtStartPar
\sphinxcode{\sphinxupquote{double Evaluate()}}
\end{quote}

\sphinxAtStartPar
\sphinxstylestrong{Return}
\begin{quote}

\sphinxAtStartPar
value of MQExpression object.
\end{quote}
\end{quote}

\subsubsection{MQExpression::GetConstant()}
\label{\detokenize{cppapi/mqexpression:mqexpression-getconstant}}\begin{quote}

\sphinxAtStartPar
Get constant in expression.

\sphinxAtStartPar
\sphinxstylestrong{Synopsis}
\begin{quote}

\sphinxAtStartPar
\sphinxcode{\sphinxupquote{double GetConstant()}}
\end{quote}

\sphinxAtStartPar
\sphinxstylestrong{Return}
\begin{quote}

\sphinxAtStartPar
constant in expression.
\end{quote}
\end{quote}

\subsubsection{MQExpression::Represent()}
\label{\detokenize{cppapi/mqexpression:mqexpression-represent}}\begin{quote}

\sphinxAtStartPar
String representation of MQExpression object.

\sphinxAtStartPar
\sphinxstylestrong{Synopsis}
\begin{quote}

\sphinxAtStartPar
\sphinxcode{\sphinxupquote{std::string Represent(size\_t maxlen)}}
\end{quote}

\sphinxAtStartPar
\sphinxstylestrong{Arguments}
\begin{quote}

\sphinxAtStartPar
\sphinxcode{\sphinxupquote{maxlen}}: max length of representation.
\end{quote}

\sphinxAtStartPar
\sphinxstylestrong{Return}
\begin{quote}

\sphinxAtStartPar
string object.
\end{quote}
\end{quote}

\subsection{MQuadExpr}
\label{\detokenize{cppapiref:mquadexpr}}\label{\detokenize{cppapiref:chapcppapiref-mquadexpr}}
\sphinxAtStartPar
COPT multi\sphinxhyphen{}dimensional quadratic expression object.It is used to construct
multi\sphinxhyphen{}dimensional quadratic expressions and perform operations with the
multi\sphinxhyphen{}dimensional array built in COPT {\hyperref[\detokenize{cppapiref:chapcppapiref-ndarray}]{\sphinxcrossref{\DUrole{std,std-ref}{NdArray}}}} .
Its elements are {\hyperref[\detokenize{cppapiref:chapcppapiref-mqexpression}]{\sphinxcrossref{\DUrole{std,std-ref}{MQExpression}}}} objects.
It can be created by quadratic combination of {\hyperref[\detokenize{cppapiref:chapcppapiref-mvar}]{\sphinxcrossref{\DUrole{std,std-ref}{MVar}}}} objects.
The following methods are provided:

\sphinxstepscope

\subsubsection{MQuadExpr::MQuadExpr()}
\label{\detokenize{cppapi/mquadexpr:mquadexpr-mquadexpr}}\label{\detokenize{cppapi/mquadexpr::doc}}\begin{quote}

\sphinxAtStartPar
Construct a MQuadExpr object with the given shape and a constant.

\sphinxAtStartPar
\sphinxstylestrong{Synopsis}
\begin{quote}

\sphinxAtStartPar
\sphinxcode{\sphinxupquote{MQuadExpr(const Shape\textless{}N\textgreater{} \&shp, double constant)}}
\end{quote}

\sphinxAtStartPar
\sphinxstylestrong{Arguments}
\begin{quote}

\sphinxAtStartPar
\sphinxcode{\sphinxupquote{shp}}: shape of MQuadExpr.

\sphinxAtStartPar
\sphinxcode{\sphinxupquote{constant}}: constant number.
\end{quote}
\end{quote}

\subsubsection{MQuadExpr::MQuadExpr()}
\label{\detokenize{cppapi/mquadexpr:id1}}\begin{quote}

\sphinxAtStartPar
Construct a MQuadExpr object with the given shape and a quadratic expression.

\sphinxAtStartPar
\sphinxstylestrong{Synopsis}
\begin{quote}

\sphinxAtStartPar
\sphinxcode{\sphinxupquote{MQuadExpr(const Shape\textless{}N\textgreater{} \&shp, const QuadExpr \&expr)}}
\end{quote}

\sphinxAtStartPar
\sphinxstylestrong{Arguments}
\begin{quote}

\sphinxAtStartPar
\sphinxcode{\sphinxupquote{shp}}: shape of MQuadExpr.

\sphinxAtStartPar
\sphinxcode{\sphinxupquote{expr}}: a quadratic expression.
\end{quote}
\end{quote}

\subsubsection{MQuadExpr::MQuadExpr()}
\label{\detokenize{cppapi/mquadexpr:id2}}\begin{quote}

\sphinxAtStartPar
Construct a MQuadExpr object with the given shape and a MQExpression object.

\sphinxAtStartPar
\sphinxstylestrong{Synopsis}
\begin{quote}

\sphinxAtStartPar
\sphinxcode{\sphinxupquote{MQuadExpr(const Shape\textless{}N\textgreater{} \&shp, const MQExpression \&expr)}}
\end{quote}

\sphinxAtStartPar
\sphinxstylestrong{Arguments}
\begin{quote}

\sphinxAtStartPar
\sphinxcode{\sphinxupquote{shp}}: shape of MQuadExpr.

\sphinxAtStartPar
\sphinxcode{\sphinxupquote{expr}}: a MQExpression object.
\end{quote}
\end{quote}

\subsubsection{MQuadExpr::AddConstant()}
\label{\detokenize{cppapi/mquadexpr:mquadexpr-addconstant}}\begin{quote}

\sphinxAtStartPar
Add constant to each quadratic expression in MQuadExpr object.

\sphinxAtStartPar
\sphinxstylestrong{Synopsis}
\begin{quote}

\sphinxAtStartPar
\sphinxcode{\sphinxupquote{void AddConstant(double constant)}}
\end{quote}

\sphinxAtStartPar
\sphinxstylestrong{Arguments}
\begin{quote}

\sphinxAtStartPar
\sphinxcode{\sphinxupquote{constant}}: the value of the constant.
\end{quote}
\end{quote}

\subsubsection{MQuadExpr::AddConstant()}
\label{\detokenize{cppapi/mquadexpr:id3}}\begin{quote}

\sphinxAtStartPar
Add constants to each quadratic expression in MQuadExpr object.

\sphinxAtStartPar
\sphinxstylestrong{Synopsis}
\begin{quote}

\sphinxAtStartPar
\sphinxcode{\sphinxupquote{template \textless{}class T\textgreater{} void AddConstant(const NdArray\textless{}T, N\textgreater{} \&constants)}}
\end{quote}

\sphinxAtStartPar
\sphinxstylestrong{Arguments}
\begin{quote}

\sphinxAtStartPar
\sphinxcode{\sphinxupquote{constants}}: N\sphinxhyphen{}dimension NdArray object.
\end{quote}
\end{quote}

\subsubsection{MQuadExpr::AddExpr()}
\label{\detokenize{cppapi/mquadexpr:mquadexpr-addexpr}}\begin{quote}

\sphinxAtStartPar
Add a linear expression to each quadratic expression in MQuadExpr object.

\sphinxAtStartPar
\sphinxstylestrong{Synopsis}
\begin{quote}

\sphinxAtStartPar
\sphinxcode{\sphinxupquote{void AddExpr(const Expr \&expr, double mult)}}
\end{quote}

\sphinxAtStartPar
\sphinxstylestrong{Arguments}
\begin{quote}

\sphinxAtStartPar
\sphinxcode{\sphinxupquote{expr}}: linear expression object.

\sphinxAtStartPar
\sphinxcode{\sphinxupquote{mult}}: the multiplier of linear expression, default value is 1.0.
\end{quote}
\end{quote}

\subsubsection{MQuadExpr::AddMExpr()}
\label{\detokenize{cppapi/mquadexpr:mquadexpr-addmexpr}}\begin{quote}

\sphinxAtStartPar
Add MExpression to each quadratic expression in MQuadExpr object.

\sphinxAtStartPar
\sphinxstylestrong{Synopsis}
\begin{quote}

\sphinxAtStartPar
\sphinxcode{\sphinxupquote{void AddMExpr(const MExpression \&expr, double mult)}}
\end{quote}

\sphinxAtStartPar
\sphinxstylestrong{Arguments}
\begin{quote}

\sphinxAtStartPar
\sphinxcode{\sphinxupquote{expr}}: MExpression object.

\sphinxAtStartPar
\sphinxcode{\sphinxupquote{mult}}: the multiplier of MExpression, default value is 1.0.
\end{quote}
\end{quote}

\subsubsection{MQuadExpr::AddMLinExpr()}
\label{\detokenize{cppapi/mquadexpr:mquadexpr-addmlinexpr}}\begin{quote}

\sphinxAtStartPar
Add linear expressions to MQuadExpr object.

\sphinxAtStartPar
\sphinxstylestrong{Synopsis}
\begin{quote}

\sphinxAtStartPar
\sphinxcode{\sphinxupquote{void AddMLinExpr(const MLinExpr\textless{}N\textgreater{} \&exprs, double mult)}}
\end{quote}

\sphinxAtStartPar
\sphinxstylestrong{Arguments}
\begin{quote}

\sphinxAtStartPar
\sphinxcode{\sphinxupquote{exprs}}: N\sphinxhyphen{}dimension MLinExpr object.

\sphinxAtStartPar
\sphinxcode{\sphinxupquote{mult}}: the same multiplier for added linear expressions, default value is 1.0.
\end{quote}
\end{quote}

\subsubsection{MQuadExpr::AddMQExpr()}
\label{\detokenize{cppapi/mquadexpr:mquadexpr-addmqexpr}}\begin{quote}

\sphinxAtStartPar
Add MQExpression to each quadratic expression in MQuadExpr object.

\sphinxAtStartPar
\sphinxstylestrong{Synopsis}
\begin{quote}

\sphinxAtStartPar
\sphinxcode{\sphinxupquote{void AddMQExpr(const MQExpression \&expr, double mult)}}
\end{quote}

\sphinxAtStartPar
\sphinxstylestrong{Arguments}
\begin{quote}

\sphinxAtStartPar
\sphinxcode{\sphinxupquote{expr}}: MQExpression object.

\sphinxAtStartPar
\sphinxcode{\sphinxupquote{mult}}: the multiplier of MQExpression, default value is 1.0.
\end{quote}
\end{quote}

\subsubsection{MQuadExpr::AddMQuadExpr()}
\label{\detokenize{cppapi/mquadexpr:mquadexpr-addmquadexpr}}\begin{quote}

\sphinxAtStartPar
Add quadratic expressions to MQuadExpr object.

\sphinxAtStartPar
\sphinxstylestrong{Synopsis}
\begin{quote}

\sphinxAtStartPar
\sphinxcode{\sphinxupquote{void AddMQuadExpr(const MQuadExpr\textless{}N\textgreater{} \&exprs, double mult)}}
\end{quote}

\sphinxAtStartPar
\sphinxstylestrong{Arguments}
\begin{quote}

\sphinxAtStartPar
\sphinxcode{\sphinxupquote{exprs}}: N\sphinxhyphen{}dimension MQuadExpr object.

\sphinxAtStartPar
\sphinxcode{\sphinxupquote{mult}}: the same multiplier for added quadratic expressions, default value is 1.0.
\end{quote}
\end{quote}

\subsubsection{MQuadExpr::AddQuadExpr()}
\label{\detokenize{cppapi/mquadexpr:mquadexpr-addquadexpr}}\begin{quote}

\sphinxAtStartPar
Add a quadratic expression to each quadratic expression in MQuadExpr object.

\sphinxAtStartPar
\sphinxstylestrong{Synopsis}
\begin{quote}

\sphinxAtStartPar
\sphinxcode{\sphinxupquote{void AddQuadExpr(const QuadExpr \&expr, double mult)}}
\end{quote}

\sphinxAtStartPar
\sphinxstylestrong{Arguments}
\begin{quote}

\sphinxAtStartPar
\sphinxcode{\sphinxupquote{expr}}: quadratic expression object.

\sphinxAtStartPar
\sphinxcode{\sphinxupquote{mult}}: the multiplier of quadratic expression, default value is 1.0.
\end{quote}
\end{quote}

\subsubsection{MQuadExpr::AddTerm()}
\label{\detokenize{cppapi/mquadexpr:mquadexpr-addterm}}\begin{quote}

\sphinxAtStartPar
Add a linear term to MQuadExpr object.

\sphinxAtStartPar
\sphinxstylestrong{Synopsis}
\begin{quote}

\sphinxAtStartPar
\sphinxcode{\sphinxupquote{void AddTerm(const Var \&var, double coeff)}}
\end{quote}

\sphinxAtStartPar
\sphinxstylestrong{Arguments}
\begin{quote}

\sphinxAtStartPar
\sphinxcode{\sphinxupquote{var}}: variable of new term.

\sphinxAtStartPar
\sphinxcode{\sphinxupquote{coeff}}: coefficient of new term.
\end{quote}
\end{quote}

\subsubsection{MQuadExpr::AddTerm()}
\label{\detokenize{cppapi/mquadexpr:id4}}\begin{quote}

\sphinxAtStartPar
Add a quadratic term to MQuadExpr object.

\sphinxAtStartPar
\sphinxstylestrong{Synopsis}
\begin{quote}

\sphinxAtStartPar
\sphinxcode{\sphinxupquote{void AddTerm(}}
\begin{quote}

\sphinxAtStartPar
\sphinxcode{\sphinxupquote{const Var \&var1,}}

\sphinxAtStartPar
\sphinxcode{\sphinxupquote{const Var \&var2,}}

\sphinxAtStartPar
\sphinxcode{\sphinxupquote{double coeff)}}
\end{quote}
\end{quote}

\sphinxAtStartPar
\sphinxstylestrong{Arguments}
\begin{quote}

\sphinxAtStartPar
\sphinxcode{\sphinxupquote{var1}}: first variable of new quadratic term.

\sphinxAtStartPar
\sphinxcode{\sphinxupquote{var2}}: second variable of new quadratic term.

\sphinxAtStartPar
\sphinxcode{\sphinxupquote{coeff}}: coefficient of new quadratic term.
\end{quote}
\end{quote}

\subsubsection{MQuadExpr::AddTerms()}
\label{\detokenize{cppapi/mquadexpr:mquadexpr-addterms}}\begin{quote}

\sphinxAtStartPar
Add terms to quadratic expressions in MQuadExpr object.

\sphinxAtStartPar
\sphinxstylestrong{Synopsis}
\begin{quote}

\sphinxAtStartPar
\sphinxcode{\sphinxupquote{void AddTerms(const MVar\textless{}N\textgreater{} \&vars, double mult)}}
\end{quote}

\sphinxAtStartPar
\sphinxstylestrong{Arguments}
\begin{quote}

\sphinxAtStartPar
\sphinxcode{\sphinxupquote{vars}}: N\sphinxhyphen{}dimension MVar object for added terms.

\sphinxAtStartPar
\sphinxcode{\sphinxupquote{mult}}: the same coefficient for added terms, default value 1.0.
\end{quote}
\end{quote}

\subsubsection{MQuadExpr::AddTerms()}
\label{\detokenize{cppapi/mquadexpr:id5}}\begin{quote}

\sphinxAtStartPar
Add terms to quadratic expressions in MQuadExpr object.

\sphinxAtStartPar
\sphinxstylestrong{Synopsis}
\begin{quote}

\sphinxAtStartPar
\sphinxcode{\sphinxupquote{void AddTerms(const MVar\textless{}N\textgreater{} \&vars, const NdArray\textless{}double, N\textgreater{} \&coeffs)}}
\end{quote}

\sphinxAtStartPar
\sphinxstylestrong{Arguments}
\begin{quote}

\sphinxAtStartPar
\sphinxcode{\sphinxupquote{vars}}: N\sphinxhyphen{}dimension MVar object for added terms.

\sphinxAtStartPar
\sphinxcode{\sphinxupquote{coeffs}}: N\sphinxhyphen{}dimension NdArray object of coefficients for added terms.
\end{quote}
\end{quote}

\subsubsection{MQuadExpr::Clear()}
\label{\detokenize{cppapi/mquadexpr:mquadexpr-clear}}\begin{quote}

\sphinxAtStartPar
Clear MQuadExpr object.

\sphinxAtStartPar
\sphinxstylestrong{Synopsis}
\begin{quote}

\sphinxAtStartPar
\sphinxcode{\sphinxupquote{void Clear()}}
\end{quote}
\end{quote}

\subsubsection{MQuadExpr::Clone()}
\label{\detokenize{cppapi/mquadexpr:mquadexpr-clone}}\begin{quote}

\sphinxAtStartPar
Clone MQuadExpr object.

\sphinxAtStartPar
\sphinxstylestrong{Synopsis}
\begin{quote}

\sphinxAtStartPar
\sphinxcode{\sphinxupquote{MQuadExpr Clone()}}
\end{quote}

\sphinxAtStartPar
\sphinxstylestrong{Return}
\begin{quote}

\sphinxAtStartPar
new MQuadExpr object.
\end{quote}
\end{quote}

\subsubsection{MQuadExpr::Diagonal()}
\label{\detokenize{cppapi/mquadexpr:mquadexpr-diagonal}}\begin{quote}

\sphinxAtStartPar
Get diagonals of MQuadExpr object.

\sphinxAtStartPar
\sphinxstylestrong{Synopsis}
\begin{quote}

\sphinxAtStartPar
\sphinxcode{\sphinxupquote{MQuadExpr\textless{}N \sphinxhyphen{} 1\textgreater{} Diagonal(}}
\begin{quote}

\sphinxAtStartPar
\sphinxcode{\sphinxupquote{int offset,}}

\sphinxAtStartPar
\sphinxcode{\sphinxupquote{int axis1,}}

\sphinxAtStartPar
\sphinxcode{\sphinxupquote{int axis2)}}
\end{quote}
\end{quote}

\sphinxAtStartPar
\sphinxstylestrong{Arguments}
\begin{quote}

\sphinxAtStartPar
\sphinxcode{\sphinxupquote{offset}}: offset of the diagonal from the main diagonal. Can be positive or negative.

\sphinxAtStartPar
\sphinxcode{\sphinxupquote{axis1}}: 1st axis of MQuadExpr.

\sphinxAtStartPar
\sphinxcode{\sphinxupquote{axis2}}: 2nd axis of MQuadExpr.
\end{quote}

\sphinxAtStartPar
\sphinxstylestrong{Return}
\begin{quote}

\sphinxAtStartPar
(N\sphinxhyphen{}1)\sphinxhyphen{}dimensional diagonals.
\end{quote}
\end{quote}

\subsubsection{MQuadExpr::Evaluate()}
\label{\detokenize{cppapi/mquadexpr:mquadexpr-evaluate}}\begin{quote}

\sphinxAtStartPar
Evaluate MQuadExpr object after solving.

\sphinxAtStartPar
\sphinxstylestrong{Synopsis}
\begin{quote}

\sphinxAtStartPar
\sphinxcode{\sphinxupquote{NdArray\textless{}double, N\textgreater{} Evaluate()}}
\end{quote}

\sphinxAtStartPar
\sphinxstylestrong{Return}
\begin{quote}

\sphinxAtStartPar
NdArray object storing value of each quadratic expression.
\end{quote}
\end{quote}

\subsubsection{MQuadExpr::Expand()}
\label{\detokenize{cppapi/mquadexpr:mquadexpr-expand}}\begin{quote}

\sphinxAtStartPar
Expand shape of MQuadExpr object.

\sphinxAtStartPar
\sphinxstylestrong{Synopsis}
\begin{quote}

\sphinxAtStartPar
\sphinxcode{\sphinxupquote{MQuadExpr\textless{}N + 1\textgreater{} Expand(int axis)}}
\end{quote}

\sphinxAtStartPar
\sphinxstylestrong{Arguments}
\begin{quote}

\sphinxAtStartPar
\sphinxcode{\sphinxupquote{axis}}: axis of MQuadExpr.
\end{quote}

\sphinxAtStartPar
\sphinxstylestrong{Return}
\begin{quote}

\sphinxAtStartPar
MQuadExpr object of (N+1)\sphinxhyphen{}dimensional shape.
\end{quote}
\end{quote}

\subsubsection{MQuadExpr::Flatten()}
\label{\detokenize{cppapi/mquadexpr:mquadexpr-flatten}}\begin{quote}

\sphinxAtStartPar
Flatten a MQuadExpr object to a 1\sphinxhyphen{}dimensional shape.

\sphinxAtStartPar
\sphinxstylestrong{Synopsis}
\begin{quote}

\sphinxAtStartPar
\sphinxcode{\sphinxupquote{MQuadExpr\textless{}1\textgreater{} Flatten()}}
\end{quote}

\sphinxAtStartPar
\sphinxstylestrong{Return}
\begin{quote}

\sphinxAtStartPar
a MQuadExpr object collapsed into one dimension.
\end{quote}
\end{quote}

\subsubsection{MQuadExpr::GetDim()}
\label{\detokenize{cppapi/mquadexpr:mquadexpr-getdim}}\begin{quote}

\sphinxAtStartPar
Get i\sphinxhyphen{}th dimension of MQuadExpr object.

\sphinxAtStartPar
\sphinxstylestrong{Synopsis}
\begin{quote}

\sphinxAtStartPar
\sphinxcode{\sphinxupquote{size\_t GetDim(int i)}}
\end{quote}

\sphinxAtStartPar
\sphinxstylestrong{Arguments}
\begin{quote}

\sphinxAtStartPar
\sphinxcode{\sphinxupquote{i}}: index of dimension
\end{quote}

\sphinxAtStartPar
\sphinxstylestrong{Return}
\begin{quote}

\sphinxAtStartPar
i\sphinxhyphen{}th dimension.
\end{quote}
\end{quote}

\subsubsection{MQuadExpr::GetND()}
\label{\detokenize{cppapi/mquadexpr:mquadexpr-getnd}}\begin{quote}

\sphinxAtStartPar
Get number of dimensions of MQuadExpr object.

\sphinxAtStartPar
\sphinxstylestrong{Synopsis}
\begin{quote}

\sphinxAtStartPar
\sphinxcode{\sphinxupquote{int GetND()}}
\end{quote}

\sphinxAtStartPar
\sphinxstylestrong{Return}
\begin{quote}

\sphinxAtStartPar
number of dimensions.
\end{quote}
\end{quote}

\subsubsection{MQuadExpr::GetShape()}
\label{\detokenize{cppapi/mquadexpr:mquadexpr-getshape}}\begin{quote}

\sphinxAtStartPar
Get shape of MQuadExpr object.

\sphinxAtStartPar
\sphinxstylestrong{Synopsis}
\begin{quote}

\sphinxAtStartPar
\sphinxcode{\sphinxupquote{Shape\textless{}N\textgreater{} GetShape()}}
\end{quote}

\sphinxAtStartPar
\sphinxstylestrong{Return}
\begin{quote}

\sphinxAtStartPar
shape object.
\end{quote}
\end{quote}

\subsubsection{MQuadExpr::GetSize()}
\label{\detokenize{cppapi/mquadexpr:mquadexpr-getsize}}\begin{quote}

\sphinxAtStartPar
Get size of MQuadExpr object.

\sphinxAtStartPar
\sphinxstylestrong{Synopsis}
\begin{quote}

\sphinxAtStartPar
\sphinxcode{\sphinxupquote{size\_t GetSize()}}
\end{quote}

\sphinxAtStartPar
\sphinxstylestrong{Return}
\begin{quote}

\sphinxAtStartPar
number of MQExpressions.
\end{quote}
\end{quote}

\subsubsection{MQuadExpr::Item()}
\label{\detokenize{cppapi/mquadexpr:mquadexpr-item}}\begin{quote}

\sphinxAtStartPar
Get quadratic expression of given index from MQuadExpr object.

\sphinxAtStartPar
\sphinxstylestrong{Synopsis}
\begin{quote}

\sphinxAtStartPar
\sphinxcode{\sphinxupquote{MQExpression \&Item(size\_t idx)}}
\end{quote}

\sphinxAtStartPar
\sphinxstylestrong{Arguments}
\begin{quote}

\sphinxAtStartPar
\sphinxcode{\sphinxupquote{idx}}: index of quadratic expression.
\end{quote}

\sphinxAtStartPar
\sphinxstylestrong{Return}
\begin{quote}

\sphinxAtStartPar
quadratic expression object.
\end{quote}
\end{quote}

\subsubsection{MQuadExpr::Item()}
\label{\detokenize{cppapi/mquadexpr:id6}}\begin{quote}

\sphinxAtStartPar
Get sub\sphinxhyphen{}arrays of MQuadExpr object, given view object.

\sphinxAtStartPar
\sphinxstylestrong{Synopsis}
\begin{quote}

\sphinxAtStartPar
\sphinxcode{\sphinxupquote{MQuadExpr Item(const View \&view)}}
\end{quote}

\sphinxAtStartPar
\sphinxstylestrong{Arguments}
\begin{quote}

\sphinxAtStartPar
\sphinxcode{\sphinxupquote{view}}: view of multi\sphinxhyphen{}dimensional array.
\end{quote}

\sphinxAtStartPar
\sphinxstylestrong{Return}
\begin{quote}

\sphinxAtStartPar
sub\sphinxhyphen{}arrays of MQuadExpr object.
\end{quote}
\end{quote}

\subsubsection{MQuadExpr::operator{[}{]}()}
\label{\detokenize{cppapi/mquadexpr:mquadexpr-operator}}\begin{quote}

\sphinxAtStartPar
Get quadratic expression of given index from MQuadExpr object.

\sphinxAtStartPar
\sphinxstylestrong{Synopsis}
\begin{quote}

\sphinxAtStartPar
\sphinxcode{\sphinxupquote{MQExpression \&operator{[}{]}(size\_t i)}}
\end{quote}

\sphinxAtStartPar
\sphinxstylestrong{Arguments}
\begin{quote}

\sphinxAtStartPar
\sphinxcode{\sphinxupquote{i}}: index of quadratic expression.
\end{quote}

\sphinxAtStartPar
\sphinxstylestrong{Return}
\begin{quote}

\sphinxAtStartPar
quadratic expression object.
\end{quote}
\end{quote}

\subsubsection{MQuadExpr::operator{[}{]}()}
\label{\detokenize{cppapi/mquadexpr:id7}}\begin{quote}

\sphinxAtStartPar
Get constraints of given view from MQuadExpr object.

\sphinxAtStartPar
\sphinxstylestrong{Synopsis}
\begin{quote}

\sphinxAtStartPar
\sphinxcode{\sphinxupquote{MQuadExpr operator{[}{]}(const View \&view)}}
\end{quote}

\sphinxAtStartPar
\sphinxstylestrong{Arguments}
\begin{quote}

\sphinxAtStartPar
\sphinxcode{\sphinxupquote{view}}: view of multi\sphinxhyphen{}dimensional array.
\end{quote}

\sphinxAtStartPar
\sphinxstylestrong{Return}
\begin{quote}

\sphinxAtStartPar
new MQuadExpr object.
\end{quote}
\end{quote}

\subsubsection{MQuadExpr::Pick()}
\label{\detokenize{cppapi/mquadexpr:mquadexpr-pick}}\begin{quote}

\sphinxAtStartPar
Given a list of indexes, get quadratic expressions from MQuadExpr object.

\sphinxAtStartPar
\sphinxstylestrong{Synopsis}
\begin{quote}

\sphinxAtStartPar
\sphinxcode{\sphinxupquote{MQuadExpr\textless{}1\textgreater{} Pick(const NdArray\textless{}int, 1\textgreater{} \&indexes)}}
\end{quote}

\sphinxAtStartPar
\sphinxstylestrong{Arguments}
\begin{quote}

\sphinxAtStartPar
\sphinxcode{\sphinxupquote{indexes}}: indexes of elements.
\end{quote}

\sphinxAtStartPar
\sphinxstylestrong{Return}
\begin{quote}

\sphinxAtStartPar
one\sphinxhyphen{}dimensional array of desired quadratic expressions.
\end{quote}
\end{quote}

\subsubsection{MQuadExpr::Pick()}
\label{\detokenize{cppapi/mquadexpr:id8}}\begin{quote}

\sphinxAtStartPar
Given a list of indexes, get quadratic expressions from MQuadExpr object.

\sphinxAtStartPar
\sphinxstylestrong{Synopsis}
\begin{quote}

\sphinxAtStartPar
\sphinxcode{\sphinxupquote{MQuadExpr\textless{}1\textgreater{} Pick(const NdArray\textless{}int, 2\textgreater{} \&idxrows)}}
\end{quote}

\sphinxAtStartPar
\sphinxstylestrong{Arguments}
\begin{quote}

\sphinxAtStartPar
\sphinxcode{\sphinxupquote{idxrows}}: indexes in format of 2\sphinxhyphen{}dimensional array, where each row is position of element.
\end{quote}

\sphinxAtStartPar
\sphinxstylestrong{Return}
\begin{quote}

\sphinxAtStartPar
one\sphinxhyphen{}dimensional array of desired quadratic expressions.
\end{quote}
\end{quote}

\subsubsection{MQuadExpr::Repeat()}
\label{\detokenize{cppapi/mquadexpr:mquadexpr-repeat}}\begin{quote}

\sphinxAtStartPar
Repeat each element of MQuadExpr along given axis.

\sphinxAtStartPar
\sphinxstylestrong{Synopsis}
\begin{quote}

\sphinxAtStartPar
\sphinxcode{\sphinxupquote{MQuadExpr\textless{}N\textgreater{} Repeat(size\_t repeats, int axis)}}
\end{quote}

\sphinxAtStartPar
\sphinxstylestrong{Arguments}
\begin{quote}

\sphinxAtStartPar
\sphinxcode{\sphinxupquote{repeats}}: number of repetitions for each element.

\sphinxAtStartPar
\sphinxcode{\sphinxupquote{axis}}: axis of MQuadExpr.
\end{quote}

\sphinxAtStartPar
\sphinxstylestrong{Return}
\begin{quote}

\sphinxAtStartPar
new MQuadExpr object.
\end{quote}
\end{quote}

\subsubsection{MQuadExpr::RepeatBlock()}
\label{\detokenize{cppapi/mquadexpr:mquadexpr-repeatblock}}\begin{quote}

\sphinxAtStartPar
Repeat an MQuadExpr a number of times along given axis.

\sphinxAtStartPar
\sphinxstylestrong{Synopsis}
\begin{quote}

\sphinxAtStartPar
\sphinxcode{\sphinxupquote{MQuadExpr\textless{}N\textgreater{} RepeatBlock(size\_t repeats, int axis)}}
\end{quote}

\sphinxAtStartPar
\sphinxstylestrong{Arguments}
\begin{quote}

\sphinxAtStartPar
\sphinxcode{\sphinxupquote{repeats}}: number of repetitions.

\sphinxAtStartPar
\sphinxcode{\sphinxupquote{axis}}: axis of MQuadExpr.
\end{quote}

\sphinxAtStartPar
\sphinxstylestrong{Return}
\begin{quote}

\sphinxAtStartPar
new MQuadExpr object.
\end{quote}
\end{quote}

\subsubsection{MQuadExpr::Represent()}
\label{\detokenize{cppapi/mquadexpr:mquadexpr-represent}}\begin{quote}

\sphinxAtStartPar
String representation of MQuadExpr object.

\sphinxAtStartPar
\sphinxstylestrong{Synopsis}
\begin{quote}

\sphinxAtStartPar
\sphinxcode{\sphinxupquote{std::string Represent(size\_t maxlen)}}
\end{quote}

\sphinxAtStartPar
\sphinxstylestrong{Arguments}
\begin{quote}

\sphinxAtStartPar
\sphinxcode{\sphinxupquote{maxlen}}: max length of representation.
\end{quote}

\sphinxAtStartPar
\sphinxstylestrong{Return}
\begin{quote}

\sphinxAtStartPar
string object.
\end{quote}
\end{quote}

\subsubsection{MQuadExpr::Reshape()}
\label{\detokenize{cppapi/mquadexpr:mquadexpr-reshape}}\begin{quote}

\sphinxAtStartPar
Reshape MQuadExpr object to new shape.

\sphinxAtStartPar
\sphinxstylestrong{Synopsis}
\begin{quote}

\sphinxAtStartPar
\sphinxcode{\sphinxupquote{template \textless{}int M\textgreater{} MQuadExpr\textless{}M\textgreater{} Reshape(const Shape\textless{}M\textgreater{} \&shape)}}
\end{quote}

\sphinxAtStartPar
\sphinxstylestrong{Arguments}
\begin{quote}

\sphinxAtStartPar
\sphinxcode{\sphinxupquote{shape}}: new shape of M\sphinxhyphen{}dimensions.
\end{quote}

\sphinxAtStartPar
\sphinxstylestrong{Return}
\begin{quote}

\sphinxAtStartPar
M\sphinxhyphen{}dimensional MQuadExpr object.
\end{quote}
\end{quote}

\subsubsection{MQuadExpr::SetItem()}
\label{\detokenize{cppapi/mquadexpr:mquadexpr-setitem}}\begin{quote}

\sphinxAtStartPar
Set quadratic expression of given index to MQuadExpr object.

\sphinxAtStartPar
\sphinxstylestrong{Synopsis}
\begin{quote}

\sphinxAtStartPar
\sphinxcode{\sphinxupquote{void SetItem(size\_t idx, const MQExpression \&expr)}}
\end{quote}

\sphinxAtStartPar
\sphinxstylestrong{Arguments}
\begin{quote}

\sphinxAtStartPar
\sphinxcode{\sphinxupquote{idx}}: index of element.

\sphinxAtStartPar
\sphinxcode{\sphinxupquote{expr}}: MQExpression object.
\end{quote}
\end{quote}

\subsubsection{MQuadExpr::Squeeze()}
\label{\detokenize{cppapi/mquadexpr:mquadexpr-squeeze}}\begin{quote}

\sphinxAtStartPar
Remove axis of length 1 from shape of MQuadExpr object.

\sphinxAtStartPar
\sphinxstylestrong{Synopsis}
\begin{quote}

\sphinxAtStartPar
\sphinxcode{\sphinxupquote{MQuadExpr\textless{}N \sphinxhyphen{} 1\textgreater{} Squeeze(int axis)}}
\end{quote}

\sphinxAtStartPar
\sphinxstylestrong{Arguments}
\begin{quote}

\sphinxAtStartPar
\sphinxcode{\sphinxupquote{axis}}: axis of MQuadExpr, where the length is 1.
\end{quote}

\sphinxAtStartPar
\sphinxstylestrong{Return}
\begin{quote}

\sphinxAtStartPar
MQuadExpr object of (N\sphinxhyphen{}1)\sphinxhyphen{}dimensional shape.
\end{quote}
\end{quote}

\subsubsection{MQuadExpr::Stack()}
\label{\detokenize{cppapi/mquadexpr:mquadexpr-stack}}\begin{quote}

\sphinxAtStartPar
Stack with other MQuadExpr object along given axis.

\sphinxAtStartPar
\sphinxstylestrong{Synopsis}
\begin{quote}

\sphinxAtStartPar
\sphinxcode{\sphinxupquote{MQuadExpr\textless{}N\textgreater{} Stack(const MQuadExpr\textless{}N\textgreater{} \&other, int axis)}}
\end{quote}

\sphinxAtStartPar
\sphinxstylestrong{Arguments}
\begin{quote}

\sphinxAtStartPar
\sphinxcode{\sphinxupquote{other}}: a MQuadExpr object.

\sphinxAtStartPar
\sphinxcode{\sphinxupquote{axis}}: an axis of MQuadExpr.
\end{quote}

\sphinxAtStartPar
\sphinxstylestrong{Return}
\begin{quote}

\sphinxAtStartPar
the result MQuadExpr object.
\end{quote}
\end{quote}

\subsubsection{MQuadExpr::Stack()}
\label{\detokenize{cppapi/mquadexpr:id9}}\begin{quote}

\sphinxAtStartPar
Stack with other MLinExpr object along given axis.

\sphinxAtStartPar
\sphinxstylestrong{Synopsis}
\begin{quote}

\sphinxAtStartPar
\sphinxcode{\sphinxupquote{MQuadExpr\textless{}N\textgreater{} Stack(const MLinExpr\textless{}N\textgreater{} \&other, int axis)}}
\end{quote}

\sphinxAtStartPar
\sphinxstylestrong{Arguments}
\begin{quote}

\sphinxAtStartPar
\sphinxcode{\sphinxupquote{other}}: a MLinExpr object.

\sphinxAtStartPar
\sphinxcode{\sphinxupquote{axis}}: an axis of MQuadExpr.
\end{quote}

\sphinxAtStartPar
\sphinxstylestrong{Return}
\begin{quote}

\sphinxAtStartPar
the result MQuadExpr object.
\end{quote}
\end{quote}

\subsubsection{MQuadExpr::Stack()}
\label{\detokenize{cppapi/mquadexpr:id10}}\begin{quote}

\sphinxAtStartPar
Stack with other MVar object along given axis.

\sphinxAtStartPar
\sphinxstylestrong{Synopsis}
\begin{quote}

\sphinxAtStartPar
\sphinxcode{\sphinxupquote{MQuadExpr\textless{}N\textgreater{} Stack(const MVar\textless{}N\textgreater{} \&other, int axis)}}
\end{quote}

\sphinxAtStartPar
\sphinxstylestrong{Arguments}
\begin{quote}

\sphinxAtStartPar
\sphinxcode{\sphinxupquote{other}}: a MVar object.

\sphinxAtStartPar
\sphinxcode{\sphinxupquote{axis}}: an axis of MQuadExpr.
\end{quote}

\sphinxAtStartPar
\sphinxstylestrong{Return}
\begin{quote}

\sphinxAtStartPar
the result MQuadExpr object.
\end{quote}
\end{quote}

\subsubsection{MQuadExpr::Stack()}
\label{\detokenize{cppapi/mquadexpr:id11}}\begin{quote}

\sphinxAtStartPar
Stack with other NdArray object along given axis.

\sphinxAtStartPar
\sphinxstylestrong{Synopsis}
\begin{quote}

\sphinxAtStartPar
\sphinxcode{\sphinxupquote{template \textless{}class T\textgreater{} MQuadExpr\textless{}N\textgreater{} Stack(const NdArray\textless{}T, N\textgreater{} \&other, int axis)}}
\end{quote}

\sphinxAtStartPar
\sphinxstylestrong{Arguments}
\begin{quote}

\sphinxAtStartPar
\sphinxcode{\sphinxupquote{other}}: a NdArray object.

\sphinxAtStartPar
\sphinxcode{\sphinxupquote{axis}}: an axis of MQuadExpr.
\end{quote}

\sphinxAtStartPar
\sphinxstylestrong{Return}
\begin{quote}

\sphinxAtStartPar
the result MQuadExpr object.
\end{quote}
\end{quote}

\subsubsection{MQuadExpr::SubConstant()}
\label{\detokenize{cppapi/mquadexpr:mquadexpr-subconstant}}\begin{quote}

\sphinxAtStartPar
Substract constants from each quadratic expression in MQuadExpr object.

\sphinxAtStartPar
\sphinxstylestrong{Synopsis}
\begin{quote}

\sphinxAtStartPar
\sphinxcode{\sphinxupquote{template \textless{}class T\textgreater{} void SubConstant(const NdArray\textless{}T, N\textgreater{} \&constants)}}
\end{quote}

\sphinxAtStartPar
\sphinxstylestrong{Arguments}
\begin{quote}

\sphinxAtStartPar
\sphinxcode{\sphinxupquote{constants}}: N\sphinxhyphen{}dimension NdArray object.
\end{quote}
\end{quote}

\subsubsection{MQuadExpr::Sum()}
\label{\detokenize{cppapi/mquadexpr:mquadexpr-sum}}\begin{quote}

\sphinxAtStartPar
Sum of all quadratic expressions in MQuadExpr object.

\sphinxAtStartPar
\sphinxstylestrong{Synopsis}
\begin{quote}

\sphinxAtStartPar
\sphinxcode{\sphinxupquote{MQuadExpr\textless{}0\textgreater{} Sum()}}
\end{quote}

\sphinxAtStartPar
\sphinxstylestrong{Return}
\begin{quote}

\sphinxAtStartPar
sum in zero dimension.
\end{quote}
\end{quote}

\subsubsection{MQuadExpr::Sum()}
\label{\detokenize{cppapi/mquadexpr:id12}}\begin{quote}

\sphinxAtStartPar
Sum of quadratic expressions at given axis of MQuadExpr object.

\sphinxAtStartPar
\sphinxstylestrong{Synopsis}
\begin{quote}

\sphinxAtStartPar
\sphinxcode{\sphinxupquote{MQuadExpr\textless{}N \sphinxhyphen{} 1\textgreater{} Sum(int axis)}}
\end{quote}

\sphinxAtStartPar
\sphinxstylestrong{Arguments}
\begin{quote}

\sphinxAtStartPar
\sphinxcode{\sphinxupquote{axis}}: axis of MQuadExpr.
\end{quote}

\sphinxAtStartPar
\sphinxstylestrong{Return}
\begin{quote}

\sphinxAtStartPar
MQuadExpr object in (N\sphinxhyphen{}1)\sphinxhyphen{}dimension.
\end{quote}
\end{quote}

\subsubsection{MQuadExpr::Transpose()}
\label{\detokenize{cppapi/mquadexpr:mquadexpr-transpose}}\begin{quote}

\sphinxAtStartPar
Perform matrix transpose of MQuadExpr object.

\sphinxAtStartPar
\sphinxstylestrong{Synopsis}
\begin{quote}

\sphinxAtStartPar
\sphinxcode{\sphinxupquote{MQuadExpr\textless{}N\textgreater{} Transpose()}}
\end{quote}

\sphinxAtStartPar
\sphinxstylestrong{Return}
\begin{quote}

\sphinxAtStartPar
transposed MQuadExpr object.
\end{quote}
\end{quote}

\subsection{NlExpr Class}
\label{\detokenize{cppapiref:nlexpr-class}}\label{\detokenize{cppapiref:chapcppapiref-nlexpr}}
\sphinxAtStartPar
The \sphinxcode{\sphinxupquote{NlExpr}} class represents nonlinear expressions in COPT.
The following methods are provided:

\sphinxstepscope

\subsubsection{NlExpr::NlExpr()}
\label{\detokenize{cppapi/nlexpr:nlexpr-nlexpr}}\label{\detokenize{cppapi/nlexpr::doc}}\begin{quote}

\sphinxAtStartPar
Constructor of a nonlinear expression with a constant.

\sphinxAtStartPar
\sphinxstylestrong{Synopsis}
\begin{quote}

\sphinxAtStartPar
\sphinxcode{\sphinxupquote{NlExpr(double constant)}}
\end{quote}

\sphinxAtStartPar
\sphinxstylestrong{Arguments}
\begin{quote}

\sphinxAtStartPar
\sphinxcode{\sphinxupquote{constant}}: constant value in nonlinear expression object.
\end{quote}
\end{quote}

\subsubsection{NlExpr::NlExpr()}
\label{\detokenize{cppapi/nlexpr:id1}}\begin{quote}

\sphinxAtStartPar
Constructor of a nonlinear expression with one linear term.

\sphinxAtStartPar
\sphinxstylestrong{Synopsis}
\begin{quote}

\sphinxAtStartPar
\sphinxcode{\sphinxupquote{NlExpr(const Var \&var, double coeff)}}
\end{quote}

\sphinxAtStartPar
\sphinxstylestrong{Arguments}
\begin{quote}

\sphinxAtStartPar
\sphinxcode{\sphinxupquote{var}}: variable for the added term.

\sphinxAtStartPar
\sphinxcode{\sphinxupquote{coeff}}: coefficent for the added term.
\end{quote}
\end{quote}

\subsubsection{NlExpr::NlExpr()}
\label{\detokenize{cppapi/nlexpr:id2}}\begin{quote}

\sphinxAtStartPar
Constructor of a nonlinear expression with a linear expression.

\sphinxAtStartPar
\sphinxstylestrong{Synopsis}
\begin{quote}

\sphinxAtStartPar
\sphinxcode{\sphinxupquote{NlExpr(const Expr \&expr)}}
\end{quote}

\sphinxAtStartPar
\sphinxstylestrong{Arguments}
\begin{quote}

\sphinxAtStartPar
\sphinxcode{\sphinxupquote{expr}}: linear expression.
\end{quote}
\end{quote}

\subsubsection{NlExpr::NlExpr()}
\label{\detokenize{cppapi/nlexpr:id3}}\begin{quote}

\sphinxAtStartPar
Constructor of a nonlinear expression with a quadratic expression.

\sphinxAtStartPar
\sphinxstylestrong{Synopsis}
\begin{quote}

\sphinxAtStartPar
\sphinxcode{\sphinxupquote{NlExpr(const QuadExpr \&expr)}}
\end{quote}

\sphinxAtStartPar
\sphinxstylestrong{Arguments}
\begin{quote}

\sphinxAtStartPar
\sphinxcode{\sphinxupquote{expr}}: quadratic expression.
\end{quote}
\end{quote}

\subsubsection{NlExpr::AddConstant()}
\label{\detokenize{cppapi/nlexpr:nlexpr-addconstant}}\begin{quote}

\sphinxAtStartPar
Add constant to the nonlinear expression.

\sphinxAtStartPar
\sphinxstylestrong{Synopsis}
\begin{quote}

\sphinxAtStartPar
\sphinxcode{\sphinxupquote{void AddConstant(double constant)}}
\end{quote}

\sphinxAtStartPar
\sphinxstylestrong{Arguments}
\begin{quote}

\sphinxAtStartPar
\sphinxcode{\sphinxupquote{constant}}: value to be added.
\end{quote}
\end{quote}

\subsubsection{NlExpr::AddLinExpr()}
\label{\detokenize{cppapi/nlexpr:nlexpr-addlinexpr}}\begin{quote}

\sphinxAtStartPar
Add a linear expression to self.

\sphinxAtStartPar
\sphinxstylestrong{Synopsis}
\begin{quote}

\sphinxAtStartPar
\sphinxcode{\sphinxupquote{void AddLinExpr(const Expr \&expr, double mult)}}
\end{quote}

\sphinxAtStartPar
\sphinxstylestrong{Arguments}
\begin{quote}

\sphinxAtStartPar
\sphinxcode{\sphinxupquote{expr}}: linear expression to be added.

\sphinxAtStartPar
\sphinxcode{\sphinxupquote{mult}}: optional, constant multiplier, default value is 1.0.
\end{quote}
\end{quote}

\subsubsection{NlExpr::AddNlExpr()}
\label{\detokenize{cppapi/nlexpr:nlexpr-addnlexpr}}\begin{quote}

\sphinxAtStartPar
Add a nonlinear expression to self.

\sphinxAtStartPar
\sphinxstylestrong{Synopsis}
\begin{quote}

\sphinxAtStartPar
\sphinxcode{\sphinxupquote{void AddNlExpr(const NlExpr \&expr, double mult)}}
\end{quote}

\sphinxAtStartPar
\sphinxstylestrong{Arguments}
\begin{quote}

\sphinxAtStartPar
\sphinxcode{\sphinxupquote{expr}}: nonlinear expression to be added.

\sphinxAtStartPar
\sphinxcode{\sphinxupquote{mult}}: optional, constant multiplier, default value is 1.0.
\end{quote}
\end{quote}

\subsubsection{NlExpr::AddQuadExpr()}
\label{\detokenize{cppapi/nlexpr:nlexpr-addquadexpr}}\begin{quote}

\sphinxAtStartPar
Add a quadratic expression to self.

\sphinxAtStartPar
\sphinxstylestrong{Synopsis}
\begin{quote}

\sphinxAtStartPar
\sphinxcode{\sphinxupquote{void AddQuadExpr(const QuadExpr \&expr, double mult)}}
\end{quote}

\sphinxAtStartPar
\sphinxstylestrong{Arguments}
\begin{quote}

\sphinxAtStartPar
\sphinxcode{\sphinxupquote{expr}}: quadratic expression to be added.

\sphinxAtStartPar
\sphinxcode{\sphinxupquote{mult}}: optional, constant multiplier, default value is 1.0.
\end{quote}
\end{quote}

\subsubsection{NlExpr::AddTerm()}
\label{\detokenize{cppapi/nlexpr:nlexpr-addterm}}\begin{quote}

\sphinxAtStartPar
Add a linear term to nonlinear expression object.

\sphinxAtStartPar
\sphinxstylestrong{Synopsis}
\begin{quote}

\sphinxAtStartPar
\sphinxcode{\sphinxupquote{void AddTerm(const Var \&var, double coeff)}}
\end{quote}

\sphinxAtStartPar
\sphinxstylestrong{Arguments}
\begin{quote}

\sphinxAtStartPar
\sphinxcode{\sphinxupquote{var}}: variable of new linear term.

\sphinxAtStartPar
\sphinxcode{\sphinxupquote{coeff}}: coefficient of new linear term.
\end{quote}
\end{quote}

\subsubsection{NlExpr::AddTerms()}
\label{\detokenize{cppapi/nlexpr:nlexpr-addterms}}\begin{quote}

\sphinxAtStartPar
Add linear terms to nonlinear expression object.

\sphinxAtStartPar
\sphinxstylestrong{Synopsis}
\begin{quote}

\sphinxAtStartPar
\sphinxcode{\sphinxupquote{int AddTerms(}}
\begin{quote}

\sphinxAtStartPar
\sphinxcode{\sphinxupquote{const VarArray \&vars,}}

\sphinxAtStartPar
\sphinxcode{\sphinxupquote{double *pCoeff,}}

\sphinxAtStartPar
\sphinxcode{\sphinxupquote{int len)}}
\end{quote}
\end{quote}

\sphinxAtStartPar
\sphinxstylestrong{Arguments}
\begin{quote}

\sphinxAtStartPar
\sphinxcode{\sphinxupquote{vars}}: variables for added linear terms.

\sphinxAtStartPar
\sphinxcode{\sphinxupquote{pCoeff}}: coefficient array for added linear terms.

\sphinxAtStartPar
\sphinxcode{\sphinxupquote{len}}: length of coefficient array.
\end{quote}

\sphinxAtStartPar
\sphinxstylestrong{Return}
\begin{quote}

\sphinxAtStartPar
number of added linear terms.
\end{quote}
\end{quote}

\subsubsection{NlExpr::Clear()}
\label{\detokenize{cppapi/nlexpr:nlexpr-clear}}\begin{quote}

\sphinxAtStartPar
Clear nonlinear expression object.

\sphinxAtStartPar
\sphinxstylestrong{Synopsis}
\begin{quote}

\sphinxAtStartPar
\sphinxcode{\sphinxupquote{void Clear()}}
\end{quote}
\end{quote}

\subsubsection{NlExpr::Clone()}
\label{\detokenize{cppapi/nlexpr:nlexpr-clone}}\begin{quote}

\sphinxAtStartPar
Deep copy nonlinear expression object.

\sphinxAtStartPar
\sphinxstylestrong{Synopsis}
\begin{quote}

\sphinxAtStartPar
\sphinxcode{\sphinxupquote{NlExpr Clone()}}
\end{quote}

\sphinxAtStartPar
\sphinxstylestrong{Return}
\begin{quote}

\sphinxAtStartPar
cloned nonlinear expression object.
\end{quote}
\end{quote}

\subsubsection{NlExpr::Evaluate()}
\label{\detokenize{cppapi/nlexpr:nlexpr-evaluate}}\begin{quote}

\sphinxAtStartPar
Evaluate nonlinear expression after solving.

\sphinxAtStartPar
\sphinxstylestrong{Synopsis}
\begin{quote}

\sphinxAtStartPar
\sphinxcode{\sphinxupquote{double Evaluate()}}
\end{quote}

\sphinxAtStartPar
\sphinxstylestrong{Return}
\begin{quote}

\sphinxAtStartPar
value of nonlinear expression.
\end{quote}
\end{quote}

\subsubsection{NlExpr::GetConstant()}
\label{\detokenize{cppapi/nlexpr:nlexpr-getconstant}}\begin{quote}

\sphinxAtStartPar
Get constant in nonlinear expression.

\sphinxAtStartPar
\sphinxstylestrong{Synopsis}
\begin{quote}

\sphinxAtStartPar
\sphinxcode{\sphinxupquote{double GetConstant()}}
\end{quote}

\sphinxAtStartPar
\sphinxstylestrong{Return}
\begin{quote}

\sphinxAtStartPar
constant in nonlinear expression.
\end{quote}
\end{quote}

\subsubsection{NlExpr::GetLinExpr()}
\label{\detokenize{cppapi/nlexpr:nlexpr-getlinexpr}}\begin{quote}

\sphinxAtStartPar
Get linear expression in nonlinear expression.

\sphinxAtStartPar
\sphinxstylestrong{Synopsis}
\begin{quote}

\sphinxAtStartPar
\sphinxcode{\sphinxupquote{Expr \&GetLinExpr()}}
\end{quote}

\sphinxAtStartPar
\sphinxstylestrong{Return}
\begin{quote}

\sphinxAtStartPar
linear expression object.
\end{quote}
\end{quote}

\subsubsection{NlExpr::Negate()}
\label{\detokenize{cppapi/nlexpr:nlexpr-negate}}\begin{quote}

\sphinxAtStartPar
Negate self.

\sphinxAtStartPar
\sphinxstylestrong{Synopsis}
\begin{quote}

\sphinxAtStartPar
\sphinxcode{\sphinxupquote{void Negate()}}
\end{quote}
\end{quote}

\subsubsection{NlExpr::operator*=()}
\label{\detokenize{cppapi/nlexpr:nlexpr-operator}}\begin{quote}

\sphinxAtStartPar
Multiply a nonlinear expression to self.

\sphinxAtStartPar
\sphinxstylestrong{Synopsis}
\begin{quote}

\sphinxAtStartPar
\sphinxcode{\sphinxupquote{void operator*=(const NlExpr \&expr)}}
\end{quote}

\sphinxAtStartPar
\sphinxstylestrong{Arguments}
\begin{quote}

\sphinxAtStartPar
\sphinxcode{\sphinxupquote{expr}}: nonlinear expression for multiplication, including double, Var, Expr, QuadExpr and NlExpr.
\end{quote}
\end{quote}

\subsubsection{NlExpr::operator*()}
\label{\detokenize{cppapi/nlexpr:id4}}\begin{quote}

\sphinxAtStartPar
Multiply expression and return new nonlinear expression.

\sphinxAtStartPar
\sphinxstylestrong{Synopsis}
\begin{quote}

\sphinxAtStartPar
\sphinxcode{\sphinxupquote{NlExpr operator*(const NlExpr \&other)}}
\end{quote}

\sphinxAtStartPar
\sphinxstylestrong{Arguments}
\begin{quote}

\sphinxAtStartPar
\sphinxcode{\sphinxupquote{other}}: operand of multiplication, including double, Var, Expr, QuadExpr and NlExpr.
\end{quote}

\sphinxAtStartPar
\sphinxstylestrong{Return}
\begin{quote}

\sphinxAtStartPar
result expression.
\end{quote}
\end{quote}

\subsubsection{NlExpr::operator/=()}
\label{\detokenize{cppapi/nlexpr:id5}}\begin{quote}

\sphinxAtStartPar
Divide a nonlinear expression by self.

\sphinxAtStartPar
\sphinxstylestrong{Synopsis}
\begin{quote}

\sphinxAtStartPar
\sphinxcode{\sphinxupquote{void operator/=(const NlExpr \&expr)}}
\end{quote}

\sphinxAtStartPar
\sphinxstylestrong{Arguments}
\begin{quote}

\sphinxAtStartPar
\sphinxcode{\sphinxupquote{expr}}: nonlinear expression for multiplication, including double, Var, Expr, QuadExpr and NlExpr.
\end{quote}
\end{quote}

\subsubsection{NlExpr::operator/()}
\label{\detokenize{cppapi/nlexpr:id6}}\begin{quote}

\sphinxAtStartPar
Divide an expression and return new nonlinear expression.

\sphinxAtStartPar
\sphinxstylestrong{Synopsis}
\begin{quote}

\sphinxAtStartPar
\sphinxcode{\sphinxupquote{NlExpr operator/(const NlExpr \&other)}}
\end{quote}

\sphinxAtStartPar
\sphinxstylestrong{Arguments}
\begin{quote}

\sphinxAtStartPar
\sphinxcode{\sphinxupquote{other}}: operand of division, including double, Var, Expr, QuadExpr and NlExpr.
\end{quote}

\sphinxAtStartPar
\sphinxstylestrong{Return}
\begin{quote}

\sphinxAtStartPar
result expression.
\end{quote}
\end{quote}

\subsubsection{NlExpr::operator+=()}
\label{\detokenize{cppapi/nlexpr:id7}}\begin{quote}

\sphinxAtStartPar
Add an expression to self.

\sphinxAtStartPar
\sphinxstylestrong{Synopsis}
\begin{quote}

\sphinxAtStartPar
\sphinxcode{\sphinxupquote{void operator+=(const NlExpr \&expr)}}
\end{quote}

\sphinxAtStartPar
\sphinxstylestrong{Arguments}
\begin{quote}

\sphinxAtStartPar
\sphinxcode{\sphinxupquote{expr}}: nonlinear expression for addition, including double, Var, Expr, QuadExpr and NlExpr.
\end{quote}
\end{quote}

\subsubsection{NlExpr::operator+()}
\label{\detokenize{cppapi/nlexpr:id8}}\begin{quote}

\sphinxAtStartPar
Add expression and return new expression.

\sphinxAtStartPar
\sphinxstylestrong{Synopsis}
\begin{quote}

\sphinxAtStartPar
\sphinxcode{\sphinxupquote{NlExpr operator+(const NlExpr \&other)}}
\end{quote}

\sphinxAtStartPar
\sphinxstylestrong{Arguments}
\begin{quote}

\sphinxAtStartPar
\sphinxcode{\sphinxupquote{other}}: operand of addition, including double, Var, Expr, QuadExpr and NlExpr.
\end{quote}

\sphinxAtStartPar
\sphinxstylestrong{Return}
\begin{quote}

\sphinxAtStartPar
result expression.
\end{quote}
\end{quote}

\subsubsection{NlExpr::operator\sphinxhyphen{}=()}
\label{\detokenize{cppapi/nlexpr:id9}}\begin{quote}

\sphinxAtStartPar
Substract an expression from self.

\sphinxAtStartPar
\sphinxstylestrong{Synopsis}
\begin{quote}

\sphinxAtStartPar
\sphinxcode{\sphinxupquote{void operator\sphinxhyphen{}=(const NlExpr \&expr)}}
\end{quote}

\sphinxAtStartPar
\sphinxstylestrong{Arguments}
\begin{quote}

\sphinxAtStartPar
\sphinxcode{\sphinxupquote{expr}}: nonlinear expression for substraction, including double, Var, Expr, QuadExpr and NlExpr.
\end{quote}
\end{quote}

\subsubsection{NlExpr::operator\sphinxhyphen{}()}
\label{\detokenize{cppapi/nlexpr:id10}}\begin{quote}

\sphinxAtStartPar
Substract expression and return new expression.

\sphinxAtStartPar
\sphinxstylestrong{Synopsis}
\begin{quote}

\sphinxAtStartPar
\sphinxcode{\sphinxupquote{NlExpr operator\sphinxhyphen{}(const NlExpr \&other)}}
\end{quote}

\sphinxAtStartPar
\sphinxstylestrong{Arguments}
\begin{quote}

\sphinxAtStartPar
\sphinxcode{\sphinxupquote{other}}: operand of substraction, including double, Var, Expr, QuadExpr and NlExpr.
\end{quote}

\sphinxAtStartPar
\sphinxstylestrong{Return}
\begin{quote}

\sphinxAtStartPar
result expression.
\end{quote}
\end{quote}

\subsubsection{NlExpr::Reserve()}
\label{\detokenize{cppapi/nlexpr:nlexpr-reserve}}\begin{quote}

\sphinxAtStartPar
Reserve capacity to contain at least n items.

\sphinxAtStartPar
\sphinxstylestrong{Synopsis}
\begin{quote}

\sphinxAtStartPar
\sphinxcode{\sphinxupquote{void Reserve(size\_t n)}}
\end{quote}

\sphinxAtStartPar
\sphinxstylestrong{Arguments}
\begin{quote}

\sphinxAtStartPar
\sphinxcode{\sphinxupquote{n}}: minimum capacity for nonlinear expression object.
\end{quote}
\end{quote}

\subsubsection{NlExpr::SetConstant()}
\label{\detokenize{cppapi/nlexpr:nlexpr-setconstant}}\begin{quote}

\sphinxAtStartPar
Set constant for the nonlinear expression.

\sphinxAtStartPar
\sphinxstylestrong{Synopsis}
\begin{quote}

\sphinxAtStartPar
\sphinxcode{\sphinxupquote{void SetConstant(double constant)}}
\end{quote}

\sphinxAtStartPar
\sphinxstylestrong{Arguments}
\begin{quote}

\sphinxAtStartPar
\sphinxcode{\sphinxupquote{constant}}: the value of the constant.
\end{quote}
\end{quote}

\subsubsection{NlExpr::Size()}
\label{\detokenize{cppapi/nlexpr:nlexpr-size}}\begin{quote}

\sphinxAtStartPar
Get size of tokens in nonlinear expression.

\sphinxAtStartPar
\sphinxstylestrong{Synopsis}
\begin{quote}

\sphinxAtStartPar
\sphinxcode{\sphinxupquote{size\_t Size()}}
\end{quote}

\sphinxAtStartPar
\sphinxstylestrong{Return}
\begin{quote}

\sphinxAtStartPar
size of none\sphinxhyphen{}linear tokens.
\end{quote}
\end{quote}

\subsection{NlConstraint Class}
\label{\detokenize{cppapiref:nlconstraint-class}}\label{\detokenize{cppapiref:chapcppapiref-nlconstraint}}
\sphinxAtStartPar
The \sphinxcode{\sphinxupquote{NlConstraint}} class provides an interface for operations on nonlinear constraints in COPT.
The following methods are provided:

\sphinxstepscope

\subsubsection{NlConstraint::Get()}
\label{\detokenize{cppapi/nlconstraint:nlconstraint-get}}\label{\detokenize{cppapi/nlconstraint::doc}}\begin{quote}

\sphinxAtStartPar
Get information value of the constraint. Support informations of “LB”, “UB”, “Slack”.

\sphinxAtStartPar
\sphinxstylestrong{Synopsis}
\begin{quote}

\sphinxAtStartPar
\sphinxcode{\sphinxupquote{double Get(const char *szInfo)}}
\end{quote}

\sphinxAtStartPar
\sphinxstylestrong{Arguments}
\begin{quote}

\sphinxAtStartPar
\sphinxcode{\sphinxupquote{szInfo}}: name of the information being queried.
\end{quote}

\sphinxAtStartPar
\sphinxstylestrong{Return}
\begin{quote}

\sphinxAtStartPar
value of information.
\end{quote}
\end{quote}

\subsubsection{NlConstraint::GetIdx()}
\label{\detokenize{cppapi/nlconstraint:nlconstraint-getidx}}\begin{quote}

\sphinxAtStartPar
Get index of nonlinear constraint.

\sphinxAtStartPar
\sphinxstylestrong{Synopsis}
\begin{quote}

\sphinxAtStartPar
\sphinxcode{\sphinxupquote{int GetIdx()}}
\end{quote}

\sphinxAtStartPar
\sphinxstylestrong{Return}
\begin{quote}

\sphinxAtStartPar
the index of nonlinear constraint.
\end{quote}
\end{quote}

\subsubsection{NlConstraint::GetName()}
\label{\detokenize{cppapi/nlconstraint:nlconstraint-getname}}\begin{quote}

\sphinxAtStartPar
Get name of nonlinear constraint.

\sphinxAtStartPar
\sphinxstylestrong{Synopsis}
\begin{quote}

\sphinxAtStartPar
\sphinxcode{\sphinxupquote{const char *GetName()}}
\end{quote}

\sphinxAtStartPar
\sphinxstylestrong{Return}
\begin{quote}

\sphinxAtStartPar
the name of nonlinear constraint.
\end{quote}
\end{quote}

\subsubsection{NlConstraint::Remove()}
\label{\detokenize{cppapi/nlconstraint:nlconstraint-remove}}\begin{quote}

\sphinxAtStartPar
Remove this nonlinear constraint from model.

\sphinxAtStartPar
\sphinxstylestrong{Synopsis}
\begin{quote}

\sphinxAtStartPar
\sphinxcode{\sphinxupquote{void Remove()}}
\end{quote}
\end{quote}

\subsubsection{NlConstraint::Set()}
\label{\detokenize{cppapi/nlconstraint:nlconstraint-set}}\begin{quote}

\sphinxAtStartPar
Set information value of nonlinear constraint. Support informations of “LB” and “UB”.

\sphinxAtStartPar
\sphinxstylestrong{Synopsis}
\begin{quote}

\sphinxAtStartPar
\sphinxcode{\sphinxupquote{void Set(const char *szInfo, double value)}}
\end{quote}

\sphinxAtStartPar
\sphinxstylestrong{Arguments}
\begin{quote}

\sphinxAtStartPar
\sphinxcode{\sphinxupquote{szInfo}}: name of the information.

\sphinxAtStartPar
\sphinxcode{\sphinxupquote{value}}: new information value.
\end{quote}
\end{quote}

\subsubsection{NlConstraint::SetName()}
\label{\detokenize{cppapi/nlconstraint:nlconstraint-setname}}\begin{quote}

\sphinxAtStartPar
Set name for nonlinear constraint.

\sphinxAtStartPar
\sphinxstylestrong{Synopsis}
\begin{quote}

\sphinxAtStartPar
\sphinxcode{\sphinxupquote{void SetName(const char *szName)}}
\end{quote}

\sphinxAtStartPar
\sphinxstylestrong{Arguments}
\begin{quote}

\sphinxAtStartPar
\sphinxcode{\sphinxupquote{szName}}: the name to set.
\end{quote}
\end{quote}

\subsection{NlConstrArray Class}
\label{\detokenize{cppapiref:nlconstrarray-class}}\label{\detokenize{cppapiref:chapcppapiref-nlconstrarray}}
\sphinxAtStartPar
To facilitate operations on a group of C++ {\hyperref[\detokenize{cppapiref:chapcppapiref-nlconstraint}]{\sphinxcrossref{\DUrole{std,std-ref}{NlConstraint Class}}}} objects,
the \sphinxcode{\sphinxupquote{NlConstrArray}} class is provided in the COPT C++ API.
The following methods are provided:

\sphinxstepscope

\subsubsection{NlConstrArray::GetNlConstr()}
\label{\detokenize{cppapi/nlconstrarray:nlconstrarray-getnlconstr}}\label{\detokenize{cppapi/nlconstrarray::doc}}\begin{quote}

\sphinxAtStartPar
Get idx\sphinxhyphen{}th nonlinear constraint object.

\sphinxAtStartPar
\sphinxstylestrong{Synopsis}
\begin{quote}

\sphinxAtStartPar
\sphinxcode{\sphinxupquote{NlConstraint \&GetNlConstr(int idx)}}
\end{quote}

\sphinxAtStartPar
\sphinxstylestrong{Arguments}
\begin{quote}

\sphinxAtStartPar
\sphinxcode{\sphinxupquote{idx}}: index of the nonlinear constraint.
\end{quote}

\sphinxAtStartPar
\sphinxstylestrong{Return}
\begin{quote}

\sphinxAtStartPar
nonlinear constraint object with index value.
\end{quote}
\end{quote}

\subsubsection{NlConstrArray::PushBack()}
\label{\detokenize{cppapi/nlconstrarray:nlconstrarray-pushback}}\begin{quote}

\sphinxAtStartPar
Add a nonlinear constraint to nonlinear constraint array.

\sphinxAtStartPar
\sphinxstylestrong{Synopsis}
\begin{quote}

\sphinxAtStartPar
\sphinxcode{\sphinxupquote{void PushBack(const NlConstraint \&constr)}}
\end{quote}

\sphinxAtStartPar
\sphinxstylestrong{Arguments}
\begin{quote}

\sphinxAtStartPar
\sphinxcode{\sphinxupquote{constr}}: nonlinear constraint object.
\end{quote}
\end{quote}

\subsubsection{NlConstrArray::Reserve()}
\label{\detokenize{cppapi/nlconstrarray:nlconstrarray-reserve}}\begin{quote}

\sphinxAtStartPar
Reserve capacity to contain at least n items.

\sphinxAtStartPar
\sphinxstylestrong{Synopsis}
\begin{quote}

\sphinxAtStartPar
\sphinxcode{\sphinxupquote{void Reserve(int n)}}
\end{quote}

\sphinxAtStartPar
\sphinxstylestrong{Arguments}
\begin{quote}

\sphinxAtStartPar
\sphinxcode{\sphinxupquote{n}}: capacity of nonlinear constraint objects.
\end{quote}
\end{quote}

\subsubsection{NlConstrArray::Size()}
\label{\detokenize{cppapi/nlconstrarray:nlconstrarray-size}}\begin{quote}

\sphinxAtStartPar
Get the number of nonlinear constraint objects.

\sphinxAtStartPar
\sphinxstylestrong{Synopsis}
\begin{quote}

\sphinxAtStartPar
\sphinxcode{\sphinxupquote{int Size()}}
\end{quote}

\sphinxAtStartPar
\sphinxstylestrong{Return}
\begin{quote}

\sphinxAtStartPar
number of nonlinear constraint objects.
\end{quote}
\end{quote}

\subsection{NlConstrBuilder Class}
\label{\detokenize{cppapiref:nlconstrbuilder-class}}\label{\detokenize{cppapiref:chapcppapiref-nlconstrbuilder}}
\sphinxAtStartPar
COPT NlConstraint builder object. To help building a nonlinear constraint, given a nonlinear
expression, constraint sense and right\sphinxhyphen{}hand side value, Cardinal Optimizer
provides C++ NlConstrBuilder class, which defines the following methods.

\sphinxstepscope

\subsubsection{NlConstrBuilder::GetNlExpr()}
\label{\detokenize{cppapi/nlconstrbuilder:nlconstrbuilder-getnlexpr}}\label{\detokenize{cppapi/nlconstrbuilder::doc}}\begin{quote}

\sphinxAtStartPar
Get nonlinear expression associated with constraint.

\sphinxAtStartPar
\sphinxstylestrong{Synopsis}
\begin{quote}

\sphinxAtStartPar
\sphinxcode{\sphinxupquote{const NlExpr \&GetNlExpr()}}
\end{quote}

\sphinxAtStartPar
\sphinxstylestrong{Return}
\begin{quote}

\sphinxAtStartPar
nonlinear expression object.
\end{quote}
\end{quote}

\subsubsection{NlConstrBuilder::GetRange()}
\label{\detokenize{cppapi/nlconstrbuilder:nlconstrbuilder-getrange}}\begin{quote}

\sphinxAtStartPar
Get range from lower bound to upper bound of range constraint.

\sphinxAtStartPar
\sphinxstylestrong{Synopsis}
\begin{quote}

\sphinxAtStartPar
\sphinxcode{\sphinxupquote{double GetRange()}}
\end{quote}

\sphinxAtStartPar
\sphinxstylestrong{Return}
\begin{quote}

\sphinxAtStartPar
length from lower bound to upper bound of nonlinear constraint.
\end{quote}
\end{quote}

\subsubsection{NlConstrBuilder::GetSense()}
\label{\detokenize{cppapi/nlconstrbuilder:nlconstrbuilder-getsense}}\begin{quote}

\sphinxAtStartPar
Get sense associated with nonlinear constraint.

\sphinxAtStartPar
\sphinxstylestrong{Synopsis}
\begin{quote}

\sphinxAtStartPar
\sphinxcode{\sphinxupquote{char GetSense()}}
\end{quote}

\sphinxAtStartPar
\sphinxstylestrong{Return}
\begin{quote}

\sphinxAtStartPar
nonlinear constraint sense.
\end{quote}
\end{quote}

\subsubsection{NlConstrBuilder::Set()}
\label{\detokenize{cppapi/nlconstrbuilder:nlconstrbuilder-set}}\begin{quote}

\sphinxAtStartPar
Set detail of a nonlinear constraint to its builder object.

\sphinxAtStartPar
\sphinxstylestrong{Synopsis}
\begin{quote}

\sphinxAtStartPar
\sphinxcode{\sphinxupquote{void Set(}}
\begin{quote}

\sphinxAtStartPar
\sphinxcode{\sphinxupquote{const NlExpr \&expr,}}

\sphinxAtStartPar
\sphinxcode{\sphinxupquote{char sense,}}

\sphinxAtStartPar
\sphinxcode{\sphinxupquote{double rhs)}}
\end{quote}
\end{quote}

\sphinxAtStartPar
\sphinxstylestrong{Arguments}
\begin{quote}

\sphinxAtStartPar
\sphinxcode{\sphinxupquote{expr}}: nonlinear expression object at one side of nonlinear constraint

\sphinxAtStartPar
\sphinxcode{\sphinxupquote{sense}}: constraint sense other than COPT\_RANGE.

\sphinxAtStartPar
\sphinxcode{\sphinxupquote{rhs}}: constant of right side of nonlinear constraint.
\end{quote}
\end{quote}

\subsubsection{NlConstrBuilder::SetRange()}
\label{\detokenize{cppapi/nlconstrbuilder:nlconstrbuilder-setrange}}\begin{quote}

\sphinxAtStartPar
Set a range constraint to nonlinear constraint builder.

\sphinxAtStartPar
\sphinxstylestrong{Synopsis}
\begin{quote}

\sphinxAtStartPar
\sphinxcode{\sphinxupquote{void SetRange(const NlExpr \&expr, double range)}}
\end{quote}

\sphinxAtStartPar
\sphinxstylestrong{Arguments}
\begin{quote}

\sphinxAtStartPar
\sphinxcode{\sphinxupquote{expr}}: nonlinear expression object, whose constant is negative upper bound.

\sphinxAtStartPar
\sphinxcode{\sphinxupquote{range}}: length from lower bound to upper bound of nonlinear constraint. Must greater than 0.
\end{quote}
\end{quote}

\subsection{NlConstrBuilderArray Class}
\label{\detokenize{cppapiref:nlconstrbuilderarray-class}}\label{\detokenize{cppapiref:chapcppapiref-nlconstrbuilderarray}}
\sphinxAtStartPar
To facilitate operations on a group of C++ {\hyperref[\detokenize{cppapiref:chapcppapiref-nlconstrbuilder}]{\sphinxcrossref{\DUrole{std,std-ref}{NlConstrBuilder Class}}}} objects,
the \sphinxcode{\sphinxupquote{NlConstrBuilderArray}} class is provided in the COPT C++ API.
The following methods are provided:

\sphinxstepscope

\subsubsection{NlConstrBuilderArray::GetBuilder()}
\label{\detokenize{cppapi/nlconstrbuilderarray:nlconstrbuilderarray-getbuilder}}\label{\detokenize{cppapi/nlconstrbuilderarray::doc}}\begin{quote}

\sphinxAtStartPar
Get idx\sphinxhyphen{}th nonlinear constraint builder object.

\sphinxAtStartPar
\sphinxstylestrong{Synopsis}
\begin{quote}

\sphinxAtStartPar
\sphinxcode{\sphinxupquote{NlConstrBuilder \&GetBuilder(int idx)}}
\end{quote}

\sphinxAtStartPar
\sphinxstylestrong{Arguments}
\begin{quote}

\sphinxAtStartPar
\sphinxcode{\sphinxupquote{idx}}: index of the nonlinear constraint builder.
\end{quote}

\sphinxAtStartPar
\sphinxstylestrong{Return}
\begin{quote}

\sphinxAtStartPar
nonlinear constraint builder object with index idx.
\end{quote}
\end{quote}

\subsubsection{NlConstrBuilderArray::PushBack()}
\label{\detokenize{cppapi/nlconstrbuilderarray:nlconstrbuilderarray-pushback}}\begin{quote}

\sphinxAtStartPar
Add a nonlinear constraint builder object to nonlinear constraint builder array.

\sphinxAtStartPar
\sphinxstylestrong{Synopsis}
\begin{quote}

\sphinxAtStartPar
\sphinxcode{\sphinxupquote{void PushBack(const NlConstrBuilder \&builder)}}
\end{quote}

\sphinxAtStartPar
\sphinxstylestrong{Arguments}
\begin{quote}

\sphinxAtStartPar
\sphinxcode{\sphinxupquote{builder}}: a nonlinear constraint builder object.
\end{quote}
\end{quote}

\subsubsection{NlConstrBuilderArray::Reserve()}
\label{\detokenize{cppapi/nlconstrbuilderarray:nlconstrbuilderarray-reserve}}\begin{quote}

\sphinxAtStartPar
Reserve capacity to contain at least n items.

\sphinxAtStartPar
\sphinxstylestrong{Synopsis}
\begin{quote}

\sphinxAtStartPar
\sphinxcode{\sphinxupquote{void Reserve(int n)}}
\end{quote}

\sphinxAtStartPar
\sphinxstylestrong{Arguments}
\begin{quote}

\sphinxAtStartPar
\sphinxcode{\sphinxupquote{n}}: minimum capacity for nonlinear constraint builder object.
\end{quote}
\end{quote}

\subsubsection{NlConstrBuilderArray::Size()}
\label{\detokenize{cppapi/nlconstrbuilderarray:nlconstrbuilderarray-size}}\begin{quote}

\sphinxAtStartPar
Get the number of nonlinear constraint builder objects.

\sphinxAtStartPar
\sphinxstylestrong{Synopsis}
\begin{quote}

\sphinxAtStartPar
\sphinxcode{\sphinxupquote{int Size()}}
\end{quote}

\sphinxAtStartPar
\sphinxstylestrong{Return}
\begin{quote}

\sphinxAtStartPar
number of nonlinear constraint builder objects.
\end{quote}
\end{quote}

\subsection{NL Namespace}
\label{\detokenize{cppapiref:nl-namespace}}\label{\detokenize{cppapiref:chapcppapiref-nl}}
\sphinxAtStartPar
The \sphinxcode{\sphinxupquote{NL}} namespace provides common nonlinear functions for constructing nonlinear expressions.
The following methods are provided:

\sphinxstepscope

\subsubsection{NL::Abs()}
\label{\detokenize{cppapi/nlp:nl-abs}}\label{\detokenize{cppapi/nlp::doc}}\begin{quote}

\sphinxAtStartPar
Calculate absolute value of a nonlinear expression.

\sphinxAtStartPar
\sphinxstylestrong{Synopsis}
\begin{quote}

\sphinxAtStartPar
\sphinxcode{\sphinxupquote{NlExpr Abs(const NlExpr \&expr)}}
\end{quote}

\sphinxAtStartPar
\sphinxstylestrong{Arguments}
\begin{quote}

\sphinxAtStartPar
\sphinxcode{\sphinxupquote{expr}}: a nonlinear expression.
\end{quote}

\sphinxAtStartPar
\sphinxstylestrong{Return}
\begin{quote}

\sphinxAtStartPar
result as a nonlinear expression.
\end{quote}
\end{quote}

\subsubsection{NL::ACos()}
\label{\detokenize{cppapi/nlp:nl-acos}}\begin{quote}

\sphinxAtStartPar
Calculate arccosine of a nonlinear expression.

\sphinxAtStartPar
\sphinxstylestrong{Synopsis}
\begin{quote}

\sphinxAtStartPar
\sphinxcode{\sphinxupquote{NlExpr ACos(const NlExpr \&expr)}}
\end{quote}

\sphinxAtStartPar
\sphinxstylestrong{Arguments}
\begin{quote}

\sphinxAtStartPar
\sphinxcode{\sphinxupquote{expr}}: a nonlinear expression.
\end{quote}

\sphinxAtStartPar
\sphinxstylestrong{Return}
\begin{quote}

\sphinxAtStartPar
result as a nonlinear expression.
\end{quote}
\end{quote}

\subsubsection{NL::ACosH()}
\label{\detokenize{cppapi/nlp:nl-acosh}}\begin{quote}

\sphinxAtStartPar
Calculate inverse hyperbolic cosine of a nonlinear expression.

\sphinxAtStartPar
\sphinxstylestrong{Synopsis}
\begin{quote}

\sphinxAtStartPar
\sphinxcode{\sphinxupquote{NlExpr ACosH(const NlExpr \&expr)}}
\end{quote}

\sphinxAtStartPar
\sphinxstylestrong{Arguments}
\begin{quote}

\sphinxAtStartPar
\sphinxcode{\sphinxupquote{expr}}: a nonlinear expression.
\end{quote}

\sphinxAtStartPar
\sphinxstylestrong{Return}
\begin{quote}

\sphinxAtStartPar
result as a nonlinear expression.
\end{quote}
\end{quote}

\subsubsection{NL::ASin()}
\label{\detokenize{cppapi/nlp:nl-asin}}\begin{quote}

\sphinxAtStartPar
Calculate arcsine of a nonlinear expression.

\sphinxAtStartPar
\sphinxstylestrong{Synopsis}
\begin{quote}

\sphinxAtStartPar
\sphinxcode{\sphinxupquote{NlExpr ASin(const NlExpr \&expr)}}
\end{quote}

\sphinxAtStartPar
\sphinxstylestrong{Arguments}
\begin{quote}

\sphinxAtStartPar
\sphinxcode{\sphinxupquote{expr}}: a nonlinear expression.
\end{quote}

\sphinxAtStartPar
\sphinxstylestrong{Return}
\begin{quote}

\sphinxAtStartPar
result as a nonlinear expression.
\end{quote}
\end{quote}

\subsubsection{NL::ASinH()}
\label{\detokenize{cppapi/nlp:nl-asinh}}\begin{quote}

\sphinxAtStartPar
Calculate inverse hyperbolic sine of a nonlinear expression.

\sphinxAtStartPar
\sphinxstylestrong{Synopsis}
\begin{quote}

\sphinxAtStartPar
\sphinxcode{\sphinxupquote{NlExpr ASinH(const NlExpr \&expr)}}
\end{quote}

\sphinxAtStartPar
\sphinxstylestrong{Arguments}
\begin{quote}

\sphinxAtStartPar
\sphinxcode{\sphinxupquote{expr}}: a nonlinear expression.
\end{quote}

\sphinxAtStartPar
\sphinxstylestrong{Return}
\begin{quote}

\sphinxAtStartPar
result as a nonlinear expression.
\end{quote}
\end{quote}

\subsubsection{NL::ATan()}
\label{\detokenize{cppapi/nlp:nl-atan}}\begin{quote}

\sphinxAtStartPar
Calculate arctangent of a nonlinear expression.

\sphinxAtStartPar
\sphinxstylestrong{Synopsis}
\begin{quote}

\sphinxAtStartPar
\sphinxcode{\sphinxupquote{NlExpr ATan(const NlExpr \&expr)}}
\end{quote}

\sphinxAtStartPar
\sphinxstylestrong{Arguments}
\begin{quote}

\sphinxAtStartPar
\sphinxcode{\sphinxupquote{expr}}: a nonlinear expression.
\end{quote}

\sphinxAtStartPar
\sphinxstylestrong{Return}
\begin{quote}

\sphinxAtStartPar
result as a nonlinear expression.
\end{quote}
\end{quote}

\subsubsection{NL::ATan2()}
\label{\detokenize{cppapi/nlp:nl-atan2}}\begin{quote}

\sphinxAtStartPar
Calculate two\sphinxhyphen{}argument arctangent of a nonlinear expression.

\sphinxAtStartPar
\sphinxstylestrong{Synopsis}
\begin{quote}

\sphinxAtStartPar
\sphinxcode{\sphinxupquote{NlExpr ATan2(const NlExpr \&y, const NlExpr \&x)}}
\end{quote}

\sphinxAtStartPar
\sphinxstylestrong{Arguments}
\begin{quote}

\sphinxAtStartPar
\sphinxcode{\sphinxupquote{y}}: y coordinate as a nonlinear expression.

\sphinxAtStartPar
\sphinxcode{\sphinxupquote{x}}: x coordinate as a nonlinear expression.
\end{quote}

\sphinxAtStartPar
\sphinxstylestrong{Return}
\begin{quote}

\sphinxAtStartPar
result as a nonlinear expression.
\end{quote}
\end{quote}

\subsubsection{NL::ATanH()}
\label{\detokenize{cppapi/nlp:nl-atanh}}\begin{quote}

\sphinxAtStartPar
Calculate inverse hyperbolic tangent of a nonlinear expression.

\sphinxAtStartPar
\sphinxstylestrong{Synopsis}
\begin{quote}

\sphinxAtStartPar
\sphinxcode{\sphinxupquote{NlExpr ATanH(const NlExpr \&expr)}}
\end{quote}

\sphinxAtStartPar
\sphinxstylestrong{Arguments}
\begin{quote}

\sphinxAtStartPar
\sphinxcode{\sphinxupquote{expr}}: a nonlinear expression.
\end{quote}

\sphinxAtStartPar
\sphinxstylestrong{Return}
\begin{quote}

\sphinxAtStartPar
result as a nonlinear expression.
\end{quote}
\end{quote}

\subsubsection{NL::Ceil()}
\label{\detokenize{cppapi/nlp:nl-ceil}}\begin{quote}

\sphinxAtStartPar
Calculate ceiling value of a nonlinear expression.

\sphinxAtStartPar
\sphinxstylestrong{Synopsis}
\begin{quote}

\sphinxAtStartPar
\sphinxcode{\sphinxupquote{NlExpr Ceil(const NlExpr \&expr)}}
\end{quote}

\sphinxAtStartPar
\sphinxstylestrong{Arguments}
\begin{quote}

\sphinxAtStartPar
\sphinxcode{\sphinxupquote{expr}}: a nonlinear expression.
\end{quote}

\sphinxAtStartPar
\sphinxstylestrong{Return}
\begin{quote}

\sphinxAtStartPar
result as a nonlinear expression.
\end{quote}
\end{quote}

\subsubsection{NL::Cos()}
\label{\detokenize{cppapi/nlp:nl-cos}}\begin{quote}

\sphinxAtStartPar
Calculate cosine of a nonlinear expression.

\sphinxAtStartPar
\sphinxstylestrong{Synopsis}
\begin{quote}

\sphinxAtStartPar
\sphinxcode{\sphinxupquote{NlExpr Cos(const NlExpr \&expr)}}
\end{quote}

\sphinxAtStartPar
\sphinxstylestrong{Arguments}
\begin{quote}

\sphinxAtStartPar
\sphinxcode{\sphinxupquote{expr}}: a nonlinear expression.
\end{quote}

\sphinxAtStartPar
\sphinxstylestrong{Return}
\begin{quote}

\sphinxAtStartPar
result as a nonlinear expression.
\end{quote}
\end{quote}

\subsubsection{NL::CosH()}
\label{\detokenize{cppapi/nlp:nl-cosh}}\begin{quote}

\sphinxAtStartPar
Calculate hyperbolic cosine of a nonlinear expression.

\sphinxAtStartPar
\sphinxstylestrong{Synopsis}
\begin{quote}

\sphinxAtStartPar
\sphinxcode{\sphinxupquote{NlExpr CosH(const NlExpr \&expr)}}
\end{quote}

\sphinxAtStartPar
\sphinxstylestrong{Arguments}
\begin{quote}

\sphinxAtStartPar
\sphinxcode{\sphinxupquote{expr}}: a nonlinear expression.
\end{quote}

\sphinxAtStartPar
\sphinxstylestrong{Return}
\begin{quote}

\sphinxAtStartPar
result as a nonlinear expression.
\end{quote}
\end{quote}

\subsubsection{NL::Exp()}
\label{\detokenize{cppapi/nlp:nl-exp}}\begin{quote}

\sphinxAtStartPar
Calculate exponential function of a nonlinear expression.

\sphinxAtStartPar
\sphinxstylestrong{Synopsis}
\begin{quote}

\sphinxAtStartPar
\sphinxcode{\sphinxupquote{NlExpr Exp(const NlExpr \&expo)}}
\end{quote}

\sphinxAtStartPar
\sphinxstylestrong{Arguments}
\begin{quote}

\sphinxAtStartPar
\sphinxcode{\sphinxupquote{expo}}: exponent as a nonlinear expression.
\end{quote}

\sphinxAtStartPar
\sphinxstylestrong{Return}
\begin{quote}

\sphinxAtStartPar
result as a nonlinear expression.
\end{quote}
\end{quote}

\subsubsection{NL::Floor()}
\label{\detokenize{cppapi/nlp:nl-floor}}\begin{quote}

\sphinxAtStartPar
Calculate floor value of a nonlinear expression.

\sphinxAtStartPar
\sphinxstylestrong{Synopsis}
\begin{quote}

\sphinxAtStartPar
\sphinxcode{\sphinxupquote{NlExpr Floor(const NlExpr \&expr)}}
\end{quote}

\sphinxAtStartPar
\sphinxstylestrong{Arguments}
\begin{quote}

\sphinxAtStartPar
\sphinxcode{\sphinxupquote{expr}}: a nonlinear expression.
\end{quote}

\sphinxAtStartPar
\sphinxstylestrong{Return}
\begin{quote}

\sphinxAtStartPar
result as a nonlinear expression.
\end{quote}
\end{quote}

\subsubsection{NL::Log10()}
\label{\detokenize{cppapi/nlp:nl-log10}}\begin{quote}

\sphinxAtStartPar
Calculate logarithmic function of a nonlinear expression with base 10.

\sphinxAtStartPar
\sphinxstylestrong{Synopsis}
\begin{quote}

\sphinxAtStartPar
\sphinxcode{\sphinxupquote{NlExpr Log10(const NlExpr \&expr)}}
\end{quote}

\sphinxAtStartPar
\sphinxstylestrong{Arguments}
\begin{quote}

\sphinxAtStartPar
\sphinxcode{\sphinxupquote{expr}}: a nonlinear expression.
\end{quote}

\sphinxAtStartPar
\sphinxstylestrong{Return}
\begin{quote}

\sphinxAtStartPar
result as a nonlinear expression.
\end{quote}
\end{quote}

\subsubsection{NL::Log()}
\label{\detokenize{cppapi/nlp:nl-log}}\begin{quote}

\sphinxAtStartPar
Calculate nature logarithmic function of a nonlinear expression.

\sphinxAtStartPar
\sphinxstylestrong{Synopsis}
\begin{quote}

\sphinxAtStartPar
\sphinxcode{\sphinxupquote{NlExpr Log(const NlExpr \&expr)}}
\end{quote}

\sphinxAtStartPar
\sphinxstylestrong{Arguments}
\begin{quote}

\sphinxAtStartPar
\sphinxcode{\sphinxupquote{expr}}: a nonlinear expression.
\end{quote}

\sphinxAtStartPar
\sphinxstylestrong{Return}
\begin{quote}

\sphinxAtStartPar
result as a nonlinear expression.
\end{quote}
\end{quote}

\subsubsection{NL::Neg()}
\label{\detokenize{cppapi/nlp:nl-neg}}\begin{quote}

\sphinxAtStartPar
Calculate negative value of a nonlinear expression.

\sphinxAtStartPar
\sphinxstylestrong{Synopsis}
\begin{quote}

\sphinxAtStartPar
\sphinxcode{\sphinxupquote{NlExpr Neg(const NlExpr \&expr)}}
\end{quote}

\sphinxAtStartPar
\sphinxstylestrong{Arguments}
\begin{quote}

\sphinxAtStartPar
\sphinxcode{\sphinxupquote{expr}}: a nonlinear expression.
\end{quote}

\sphinxAtStartPar
\sphinxstylestrong{Return}
\begin{quote}

\sphinxAtStartPar
result as a nonlinear expression.
\end{quote}
\end{quote}

\subsubsection{NL::Pow()}
\label{\detokenize{cppapi/nlp:nl-pow}}\begin{quote}

\sphinxAtStartPar
Calculate power function of a nonlinear expression.

\sphinxAtStartPar
\sphinxstylestrong{Synopsis}
\begin{quote}

\sphinxAtStartPar
\sphinxcode{\sphinxupquote{NlExpr Pow(const NlExpr \&base, const NlExpr \&expo)}}
\end{quote}

\sphinxAtStartPar
\sphinxstylestrong{Arguments}
\begin{quote}

\sphinxAtStartPar
\sphinxcode{\sphinxupquote{base}}: base as a nonlinear expression.

\sphinxAtStartPar
\sphinxcode{\sphinxupquote{expo}}: exponent as a nonlinear expression.
\end{quote}

\sphinxAtStartPar
\sphinxstylestrong{Return}
\begin{quote}

\sphinxAtStartPar
result as a nonlinear expression.
\end{quote}
\end{quote}

\subsubsection{NL::Sin()}
\label{\detokenize{cppapi/nlp:nl-sin}}\begin{quote}

\sphinxAtStartPar
Calculate sine of a nonlinear expression.

\sphinxAtStartPar
\sphinxstylestrong{Synopsis}
\begin{quote}

\sphinxAtStartPar
\sphinxcode{\sphinxupquote{NlExpr Sin(const NlExpr \&expr)}}
\end{quote}

\sphinxAtStartPar
\sphinxstylestrong{Arguments}
\begin{quote}

\sphinxAtStartPar
\sphinxcode{\sphinxupquote{expr}}: a nonlinear expression.
\end{quote}

\sphinxAtStartPar
\sphinxstylestrong{Return}
\begin{quote}

\sphinxAtStartPar
result as a nonlinear expression.
\end{quote}
\end{quote}

\subsubsection{NL::SinH()}
\label{\detokenize{cppapi/nlp:nl-sinh}}\begin{quote}

\sphinxAtStartPar
Calculate hyperbolic sine of a nonlinear expression.

\sphinxAtStartPar
\sphinxstylestrong{Synopsis}
\begin{quote}

\sphinxAtStartPar
\sphinxcode{\sphinxupquote{NlExpr SinH(const NlExpr \&expr)}}
\end{quote}

\sphinxAtStartPar
\sphinxstylestrong{Arguments}
\begin{quote}

\sphinxAtStartPar
\sphinxcode{\sphinxupquote{expr}}: a nonlinear expression.
\end{quote}

\sphinxAtStartPar
\sphinxstylestrong{Return}
\begin{quote}

\sphinxAtStartPar
result as a nonlinear expression.
\end{quote}
\end{quote}

\subsubsection{NL::Sqrt()}
\label{\detokenize{cppapi/nlp:nl-sqrt}}\begin{quote}

\sphinxAtStartPar
Calculate square root of a nonlinear expression.

\sphinxAtStartPar
\sphinxstylestrong{Synopsis}
\begin{quote}

\sphinxAtStartPar
\sphinxcode{\sphinxupquote{NlExpr Sqrt(const NlExpr \&expr)}}
\end{quote}

\sphinxAtStartPar
\sphinxstylestrong{Arguments}
\begin{quote}

\sphinxAtStartPar
\sphinxcode{\sphinxupquote{expr}}: a nonlinear expression.
\end{quote}

\sphinxAtStartPar
\sphinxstylestrong{Return}
\begin{quote}

\sphinxAtStartPar
result as a nonlinear expression.
\end{quote}
\end{quote}

\subsubsection{NL::Sum()}
\label{\detokenize{cppapi/nlp:nl-sum}}\begin{quote}

\sphinxAtStartPar
Sum of nonlinear expressions.

\sphinxAtStartPar
\sphinxstylestrong{Synopsis}
\begin{quote}

\sphinxAtStartPar
\sphinxcode{\sphinxupquote{NlExpr Sum(const std::vector\textless{}NlExpr*\textgreater{} \&exprs)}}
\end{quote}

\sphinxAtStartPar
\sphinxstylestrong{Arguments}
\begin{quote}

\sphinxAtStartPar
\sphinxcode{\sphinxupquote{exprs}}: vector of nonlinear expressions.
\end{quote}

\sphinxAtStartPar
\sphinxstylestrong{Return}
\begin{quote}

\sphinxAtStartPar
result as a nonlinear expression.
\end{quote}
\end{quote}

\subsubsection{NL::Sum()}
\label{\detokenize{cppapi/nlp:id1}}\begin{quote}

\sphinxAtStartPar
Sum of nonlinear expressions.

\sphinxAtStartPar
\sphinxstylestrong{Synopsis}
\begin{quote}

\sphinxAtStartPar
\sphinxcode{\sphinxupquote{NlExpr Sum(}}
\begin{quote}

\sphinxAtStartPar
\sphinxcode{\sphinxupquote{const NlExpr \&op1,}}

\sphinxAtStartPar
\sphinxcode{\sphinxupquote{const NlExpr \&op2,}}

\sphinxAtStartPar
\sphinxcode{\sphinxupquote{const NlExpr \&op3)}}
\end{quote}
\end{quote}

\sphinxAtStartPar
\sphinxstylestrong{Arguments}
\begin{quote}

\sphinxAtStartPar
\sphinxcode{\sphinxupquote{op1}}: first nonlinear expression.

\sphinxAtStartPar
\sphinxcode{\sphinxupquote{op2}}: second nonlinear expression.

\sphinxAtStartPar
\sphinxcode{\sphinxupquote{op3}}: third nonlinear expression.
\end{quote}

\sphinxAtStartPar
\sphinxstylestrong{Return}
\begin{quote}

\sphinxAtStartPar
result as a nonlinear expression.
\end{quote}
\end{quote}

\subsubsection{NL::Sum()}
\label{\detokenize{cppapi/nlp:id2}}\begin{quote}

\sphinxAtStartPar
Sum of nonlinear expressions.

\sphinxAtStartPar
\sphinxstylestrong{Synopsis}
\begin{quote}

\sphinxAtStartPar
\sphinxcode{\sphinxupquote{NlExpr Sum(}}
\begin{quote}

\sphinxAtStartPar
\sphinxcode{\sphinxupquote{const NlExpr \&op1,}}

\sphinxAtStartPar
\sphinxcode{\sphinxupquote{const NlExpr \&op2,}}

\sphinxAtStartPar
\sphinxcode{\sphinxupquote{const NlExpr \&op3,}}

\sphinxAtStartPar
\sphinxcode{\sphinxupquote{const NlExpr \&op4)}}
\end{quote}
\end{quote}

\sphinxAtStartPar
\sphinxstylestrong{Arguments}
\begin{quote}

\sphinxAtStartPar
\sphinxcode{\sphinxupquote{op1}}: first nonlinear expression.

\sphinxAtStartPar
\sphinxcode{\sphinxupquote{op2}}: second nonlinear expression.

\sphinxAtStartPar
\sphinxcode{\sphinxupquote{op3}}: third nonlinear expression.

\sphinxAtStartPar
\sphinxcode{\sphinxupquote{op4}}: fourth nonlinear expression.
\end{quote}

\sphinxAtStartPar
\sphinxstylestrong{Return}
\begin{quote}

\sphinxAtStartPar
result as a nonlinear expression.
\end{quote}
\end{quote}

\subsubsection{NL::Tan()}
\label{\detokenize{cppapi/nlp:nl-tan}}\begin{quote}

\sphinxAtStartPar
Calculate tangent of a nonlinear expression.

\sphinxAtStartPar
\sphinxstylestrong{Synopsis}
\begin{quote}

\sphinxAtStartPar
\sphinxcode{\sphinxupquote{NlExpr Tan(const NlExpr \&expr)}}
\end{quote}

\sphinxAtStartPar
\sphinxstylestrong{Arguments}
\begin{quote}

\sphinxAtStartPar
\sphinxcode{\sphinxupquote{expr}}: a nonlinear expression.
\end{quote}

\sphinxAtStartPar
\sphinxstylestrong{Return}
\begin{quote}

\sphinxAtStartPar
result as a nonlinear expression.
\end{quote}
\end{quote}

\subsubsection{NL::TanH()}
\label{\detokenize{cppapi/nlp:nl-tanh}}\begin{quote}

\sphinxAtStartPar
Calculat hyperbolic tangent of a nonlinear expression.

\sphinxAtStartPar
\sphinxstylestrong{Synopsis}
\begin{quote}

\sphinxAtStartPar
\sphinxcode{\sphinxupquote{NlExpr TanH(const NlExpr \&expr)}}
\end{quote}

\sphinxAtStartPar
\sphinxstylestrong{Arguments}
\begin{quote}

\sphinxAtStartPar
\sphinxcode{\sphinxupquote{expr}}: a nonlinear expression.
\end{quote}

\sphinxAtStartPar
\sphinxstylestrong{Return}
\begin{quote}

\sphinxAtStartPar
result as a nonlinear expression.
\end{quote}
\end{quote}

\subsection{NdArray}
\label{\detokenize{cppapiref:ndarray}}\label{\detokenize{cppapiref:chapcppapiref-ndarray}}
\sphinxAtStartPar
The NdArray class is a built\sphinxhyphen{}in multi\sphinxhyphen{}dimensional array in COPT. It represents a
table of elements of the same type, indexed by a tuple of integers. The following
methods are provided:

\sphinxstepscope

\subsubsection{NdArray::NdArray()}
\label{\detokenize{cppapi/ndarray:ndarray-ndarray}}\label{\detokenize{cppapi/ndarray::doc}}\begin{quote}

\sphinxAtStartPar
Construct an NdArray object with the given shape, filling with the given element.

\sphinxAtStartPar
\sphinxstylestrong{Synopsis}
\begin{quote}

\sphinxAtStartPar
\sphinxcode{\sphinxupquote{NdArray(const Shape\textless{}N\textgreater{} \&shape, const T \&val)}}
\end{quote}

\sphinxAtStartPar
\sphinxstylestrong{Arguments}
\begin{quote}

\sphinxAtStartPar
\sphinxcode{\sphinxupquote{shape}}: shape of NdArray.

\sphinxAtStartPar
\sphinxcode{\sphinxupquote{val}}: value of element.
\end{quote}
\end{quote}

\subsubsection{NdArray::NdArray()}
\label{\detokenize{cppapi/ndarray:id1}}\begin{quote}

\sphinxAtStartPar
Construct an NdArray object with the given shape, filling with an array of type T.

\sphinxAtStartPar
\sphinxstylestrong{Synopsis}
\begin{quote}

\sphinxAtStartPar
\sphinxcode{\sphinxupquote{NdArray(}}
\begin{quote}

\sphinxAtStartPar
\sphinxcode{\sphinxupquote{const Shape\textless{}N\textgreater{} \&shape,}}

\sphinxAtStartPar
\sphinxcode{\sphinxupquote{const T *data,}}

\sphinxAtStartPar
\sphinxcode{\sphinxupquote{size\_t sz)}}
\end{quote}
\end{quote}

\sphinxAtStartPar
\sphinxstylestrong{Arguments}
\begin{quote}

\sphinxAtStartPar
\sphinxcode{\sphinxupquote{shape}}: shape of NdArray.

\sphinxAtStartPar
\sphinxcode{\sphinxupquote{data}}: an array of elements.

\sphinxAtStartPar
\sphinxcode{\sphinxupquote{sz}}: size of elements.
\end{quote}
\end{quote}

\subsubsection{NdArray::NdArray()}
\label{\detokenize{cppapi/ndarray:id2}}\begin{quote}

\sphinxAtStartPar
Construct an NdArray object with the given shape and the filling function.

\sphinxAtStartPar
\sphinxstylestrong{Synopsis}
\begin{quote}

\sphinxAtStartPar
\sphinxcode{\sphinxupquote{NdArray(const Shape\textless{}N\textgreater{} \&shape)}}
\end{quote}

\sphinxAtStartPar
\sphinxstylestrong{Arguments}
\begin{quote}

\sphinxAtStartPar
\sphinxcode{\sphinxupquote{shape}}: shape of NdArray.
\end{quote}
\end{quote}

\subsubsection{NdArray::Clone()}
\label{\detokenize{cppapi/ndarray:ndarray-clone}}\begin{quote}

\sphinxAtStartPar
Clone NdArray object.

\sphinxAtStartPar
\sphinxstylestrong{Synopsis}
\begin{quote}

\sphinxAtStartPar
\sphinxcode{\sphinxupquote{NdArray\textless{}T, N\textgreater{} Clone()}}
\end{quote}

\sphinxAtStartPar
\sphinxstylestrong{Return}
\begin{quote}

\sphinxAtStartPar
new NdArray object.
\end{quote}
\end{quote}

\subsubsection{NdArray::Diagonal()}
\label{\detokenize{cppapi/ndarray:ndarray-diagonal}}\begin{quote}

\sphinxAtStartPar
Get diagonals of NdArray object.

\sphinxAtStartPar
\sphinxstylestrong{Synopsis}
\begin{quote}

\sphinxAtStartPar
\sphinxcode{\sphinxupquote{NdArray\textless{}T, N \sphinxhyphen{} 1\textgreater{} Diagonal(}}
\begin{quote}

\sphinxAtStartPar
\sphinxcode{\sphinxupquote{int offset,}}

\sphinxAtStartPar
\sphinxcode{\sphinxupquote{int axis1,}}

\sphinxAtStartPar
\sphinxcode{\sphinxupquote{int axis2)}}
\end{quote}
\end{quote}

\sphinxAtStartPar
\sphinxstylestrong{Arguments}
\begin{quote}

\sphinxAtStartPar
\sphinxcode{\sphinxupquote{offset}}: offset of the diagonal from the main diagonal. Can be positive or negative.

\sphinxAtStartPar
\sphinxcode{\sphinxupquote{axis1}}: 1st axis of NdArray.

\sphinxAtStartPar
\sphinxcode{\sphinxupquote{axis2}}: 2nd axis of NdArray.
\end{quote}

\sphinxAtStartPar
\sphinxstylestrong{Return}
\begin{quote}

\sphinxAtStartPar
(N\sphinxhyphen{}1)\sphinxhyphen{}dimensional diagonals.
\end{quote}
\end{quote}

\subsubsection{NdArray::Dot()}
\label{\detokenize{cppapi/ndarray:ndarray-dot}}\begin{quote}

\sphinxAtStartPar
Dot product with another NdArray object of type double.

\sphinxAtStartPar
\sphinxstylestrong{Synopsis}
\begin{quote}

\sphinxAtStartPar
\sphinxcode{\sphinxupquote{double Dot(const NdArray\textless{}double, N\textgreater{} \&other)}}
\end{quote}

\sphinxAtStartPar
\sphinxstylestrong{Arguments}
\begin{quote}

\sphinxAtStartPar
\sphinxcode{\sphinxupquote{other}}: another NdArray object of type double.
\end{quote}

\sphinxAtStartPar
\sphinxstylestrong{Return}
\begin{quote}

\sphinxAtStartPar
dot product value.
\end{quote}
\end{quote}

\subsubsection{NdArray::Dot()}
\label{\detokenize{cppapi/ndarray:id3}}\begin{quote}

\sphinxAtStartPar
Dot product with another NdArray object of type int.

\sphinxAtStartPar
\sphinxstylestrong{Synopsis}
\begin{quote}

\sphinxAtStartPar
\sphinxcode{\sphinxupquote{T Dot(const NdArray\textless{}int, N\textgreater{} \&other)}}
\end{quote}

\sphinxAtStartPar
\sphinxstylestrong{Arguments}
\begin{quote}

\sphinxAtStartPar
\sphinxcode{\sphinxupquote{other}}: another NdArray object of type int.
\end{quote}

\sphinxAtStartPar
\sphinxstylestrong{Return}
\begin{quote}

\sphinxAtStartPar
dot product value.
\end{quote}
\end{quote}

\subsubsection{NdArray::Expand()}
\label{\detokenize{cppapi/ndarray:ndarray-expand}}\begin{quote}

\sphinxAtStartPar
Expand shape of NdArray object.

\sphinxAtStartPar
\sphinxstylestrong{Synopsis}
\begin{quote}

\sphinxAtStartPar
\sphinxcode{\sphinxupquote{NdArray\textless{}T, N + 1\textgreater{} Expand(int axis)}}
\end{quote}

\sphinxAtStartPar
\sphinxstylestrong{Arguments}
\begin{quote}

\sphinxAtStartPar
\sphinxcode{\sphinxupquote{axis}}: axis of NdArray.
\end{quote}

\sphinxAtStartPar
\sphinxstylestrong{Return}
\begin{quote}

\sphinxAtStartPar
NdArray object in (N+1)\sphinxhyphen{}dimensions.
\end{quote}
\end{quote}

\subsubsection{NdArray::Fill()}
\label{\detokenize{cppapi/ndarray:ndarray-fill}}\begin{quote}

\sphinxAtStartPar
Fill NdArray object with given value.

\sphinxAtStartPar
\sphinxstylestrong{Synopsis}
\begin{quote}

\sphinxAtStartPar
\sphinxcode{\sphinxupquote{void Fill(const T \&val)}}
\end{quote}

\sphinxAtStartPar
\sphinxstylestrong{Arguments}
\begin{quote}

\sphinxAtStartPar
\sphinxcode{\sphinxupquote{val}}: new value.
\end{quote}
\end{quote}

\subsubsection{NdArray::Flatten()}
\label{\detokenize{cppapi/ndarray:ndarray-flatten}}\begin{quote}

\sphinxAtStartPar
Flatten an NdArray object to a 1\sphinxhyphen{}dimensional shape.

\sphinxAtStartPar
\sphinxstylestrong{Synopsis}
\begin{quote}

\sphinxAtStartPar
\sphinxcode{\sphinxupquote{NdArray\textless{}T, 1\textgreater{} Flatten()}}
\end{quote}

\sphinxAtStartPar
\sphinxstylestrong{Return}
\begin{quote}

\sphinxAtStartPar
An NdArray object collapsed into one dimension.
\end{quote}
\end{quote}

\subsubsection{NdArray::GetDim()}
\label{\detokenize{cppapi/ndarray:ndarray-getdim}}\begin{quote}

\sphinxAtStartPar
Get i\sphinxhyphen{}th dimension of NdArray object.

\sphinxAtStartPar
\sphinxstylestrong{Synopsis}
\begin{quote}

\sphinxAtStartPar
\sphinxcode{\sphinxupquote{size\_t GetDim(int i)}}
\end{quote}

\sphinxAtStartPar
\sphinxstylestrong{Arguments}
\begin{quote}

\sphinxAtStartPar
\sphinxcode{\sphinxupquote{i}}: index of dimension
\end{quote}

\sphinxAtStartPar
\sphinxstylestrong{Return}
\begin{quote}

\sphinxAtStartPar
i\sphinxhyphen{}th dimension.
\end{quote}
\end{quote}

\subsubsection{NdArray::GetND()}
\label{\detokenize{cppapi/ndarray:ndarray-getnd}}\begin{quote}

\sphinxAtStartPar
Get number of dimensions of NdArray object.

\sphinxAtStartPar
\sphinxstylestrong{Synopsis}
\begin{quote}

\sphinxAtStartPar
\sphinxcode{\sphinxupquote{int GetND()}}
\end{quote}

\sphinxAtStartPar
\sphinxstylestrong{Return}
\begin{quote}

\sphinxAtStartPar
number of dimensions.
\end{quote}
\end{quote}

\subsubsection{NdArray::GetShape()}
\label{\detokenize{cppapi/ndarray:ndarray-getshape}}\begin{quote}

\sphinxAtStartPar
Get shape of NdArray object.

\sphinxAtStartPar
\sphinxstylestrong{Synopsis}
\begin{quote}

\sphinxAtStartPar
\sphinxcode{\sphinxupquote{const Shape\textless{}N\textgreater{} \&GetShape()}}
\end{quote}

\sphinxAtStartPar
\sphinxstylestrong{Return}
\begin{quote}

\sphinxAtStartPar
shape object.
\end{quote}
\end{quote}

\subsubsection{NdArray::GetSize()}
\label{\detokenize{cppapi/ndarray:ndarray-getsize}}\begin{quote}

\sphinxAtStartPar
Get size of NdArray object.

\sphinxAtStartPar
\sphinxstylestrong{Synopsis}
\begin{quote}

\sphinxAtStartPar
\sphinxcode{\sphinxupquote{size\_t GetSize()}}
\end{quote}

\sphinxAtStartPar
\sphinxstylestrong{Return}
\begin{quote}

\sphinxAtStartPar
number of elements.
\end{quote}
\end{quote}

\subsubsection{NdArray::Item()}
\label{\detokenize{cppapi/ndarray:ndarray-item}}\begin{quote}

\sphinxAtStartPar
Get element of given index from NdArray object.

\sphinxAtStartPar
\sphinxstylestrong{Synopsis}
\begin{quote}

\sphinxAtStartPar
\sphinxcode{\sphinxupquote{T \&Item(size\_t idx)}}
\end{quote}

\sphinxAtStartPar
\sphinxstylestrong{Arguments}
\begin{quote}

\sphinxAtStartPar
\sphinxcode{\sphinxupquote{idx}}: index of element.
\end{quote}

\sphinxAtStartPar
\sphinxstylestrong{Return}
\begin{quote}

\sphinxAtStartPar
value of element.
\end{quote}
\end{quote}

\subsubsection{NdArray::Item()}
\label{\detokenize{cppapi/ndarray:id4}}\begin{quote}

\sphinxAtStartPar
Get sub\sphinxhyphen{}array of NdArray object, given View object.

\sphinxAtStartPar
\sphinxstylestrong{Synopsis}
\begin{quote}

\sphinxAtStartPar
\sphinxcode{\sphinxupquote{NdArray\textless{}T, N\textgreater{} Item(const View \&view)}}
\end{quote}

\sphinxAtStartPar
\sphinxstylestrong{Arguments}
\begin{quote}

\sphinxAtStartPar
\sphinxcode{\sphinxupquote{view}}: View object.
\end{quote}

\sphinxAtStartPar
\sphinxstylestrong{Return}
\begin{quote}

\sphinxAtStartPar
sub NdArray without copying underlying data.
\end{quote}
\end{quote}

\subsubsection{NdArray::operator{[}{]}()}
\label{\detokenize{cppapi/ndarray:ndarray-operator}}\begin{quote}

\sphinxAtStartPar
Get element of given index from NdArray object.

\sphinxAtStartPar
\sphinxstylestrong{Synopsis}
\begin{quote}

\sphinxAtStartPar
\sphinxcode{\sphinxupquote{T \&operator{[}{]}(size\_t idx)}}
\end{quote}

\sphinxAtStartPar
\sphinxstylestrong{Arguments}
\begin{quote}

\sphinxAtStartPar
\sphinxcode{\sphinxupquote{idx}}: index of element.
\end{quote}

\sphinxAtStartPar
\sphinxstylestrong{Return}
\begin{quote}

\sphinxAtStartPar
value of element.
\end{quote}
\end{quote}

\subsubsection{NdArray::operator{[}{]}()}
\label{\detokenize{cppapi/ndarray:id5}}\begin{quote}

\sphinxAtStartPar
Get constraints of given view from NdArray object.

\sphinxAtStartPar
\sphinxstylestrong{Synopsis}
\begin{quote}

\sphinxAtStartPar
\sphinxcode{\sphinxupquote{NdArray\textless{}T, N\textgreater{} operator{[}{]}(const View \&view)}}
\end{quote}

\sphinxAtStartPar
\sphinxstylestrong{Arguments}
\begin{quote}

\sphinxAtStartPar
\sphinxcode{\sphinxupquote{view}}: view of multi\sphinxhyphen{}dimensional array.
\end{quote}

\sphinxAtStartPar
\sphinxstylestrong{Return}
\begin{quote}

\sphinxAtStartPar
new NdArray object.
\end{quote}
\end{quote}

\subsubsection{NdArray::Pick()}
\label{\detokenize{cppapi/ndarray:ndarray-pick}}\begin{quote}

\sphinxAtStartPar
Given a list of indexes, get elements from NdArray object.

\sphinxAtStartPar
\sphinxstylestrong{Synopsis}
\begin{quote}

\sphinxAtStartPar
\sphinxcode{\sphinxupquote{NdArray\textless{}T, 1\textgreater{} Pick(const NdArray\textless{}int, 1\textgreater{} \&indexes)}}
\end{quote}

\sphinxAtStartPar
\sphinxstylestrong{Arguments}
\begin{quote}

\sphinxAtStartPar
\sphinxcode{\sphinxupquote{indexes}}: indexes of elements.
\end{quote}

\sphinxAtStartPar
\sphinxstylestrong{Return}
\begin{quote}

\sphinxAtStartPar
one\sphinxhyphen{}dimensional array of desired elements.
\end{quote}
\end{quote}

\subsubsection{NdArray::Pick()}
\label{\detokenize{cppapi/ndarray:id6}}\begin{quote}

\sphinxAtStartPar
Given a list of indexes, get elements from NdArray object.

\sphinxAtStartPar
\sphinxstylestrong{Synopsis}
\begin{quote}

\sphinxAtStartPar
\sphinxcode{\sphinxupquote{NdArray\textless{}T, 1\textgreater{} Pick(const NdArray\textless{}int, 2\textgreater{} \&idxrows)}}
\end{quote}

\sphinxAtStartPar
\sphinxstylestrong{Arguments}
\begin{quote}

\sphinxAtStartPar
\sphinxcode{\sphinxupquote{idxrows}}: indexes in format of 2\sphinxhyphen{}dimensional array, where each row is position of element.
\end{quote}

\sphinxAtStartPar
\sphinxstylestrong{Return}
\begin{quote}

\sphinxAtStartPar
one\sphinxhyphen{}dimensional array of desired elements.
\end{quote}
\end{quote}

\subsubsection{NdArray::Prod()}
\label{\detokenize{cppapi/ndarray:ndarray-prod}}\begin{quote}

\sphinxAtStartPar
Product of all elements in NdArray object.

\sphinxAtStartPar
\sphinxstylestrong{Synopsis}
\begin{quote}

\sphinxAtStartPar
\sphinxcode{\sphinxupquote{T Prod()}}
\end{quote}

\sphinxAtStartPar
\sphinxstylestrong{Return}
\begin{quote}

\sphinxAtStartPar
product value.
\end{quote}
\end{quote}

\subsubsection{NdArray::Prod()}
\label{\detokenize{cppapi/ndarray:id7}}\begin{quote}

\sphinxAtStartPar
Product of elements at given axis of NdArray object.

\sphinxAtStartPar
\sphinxstylestrong{Synopsis}
\begin{quote}

\sphinxAtStartPar
\sphinxcode{\sphinxupquote{NdArray\textless{}T, N \sphinxhyphen{} 1\textgreater{} Prod(int axis)}}
\end{quote}

\sphinxAtStartPar
\sphinxstylestrong{Arguments}
\begin{quote}

\sphinxAtStartPar
\sphinxcode{\sphinxupquote{axis}}: axis of NdArray.
\end{quote}

\sphinxAtStartPar
\sphinxstylestrong{Return}
\begin{quote}

\sphinxAtStartPar
(N\sphinxhyphen{}1)\sphinxhyphen{}dimensional NdArray object.
\end{quote}
\end{quote}

\subsubsection{NdArray::Repeat()}
\label{\detokenize{cppapi/ndarray:ndarray-repeat}}\begin{quote}

\sphinxAtStartPar
Repeat each element of an array along given axis.

\sphinxAtStartPar
\sphinxstylestrong{Synopsis}
\begin{quote}

\sphinxAtStartPar
\sphinxcode{\sphinxupquote{NdArray\textless{}T, N\textgreater{} Repeat(size\_t repeats, int axis)}}
\end{quote}

\sphinxAtStartPar
\sphinxstylestrong{Arguments}
\begin{quote}

\sphinxAtStartPar
\sphinxcode{\sphinxupquote{repeats}}: number of repetitions for each element.

\sphinxAtStartPar
\sphinxcode{\sphinxupquote{axis}}: axis of NdArray.
\end{quote}

\sphinxAtStartPar
\sphinxstylestrong{Return}
\begin{quote}

\sphinxAtStartPar
new NdArray object.
\end{quote}
\end{quote}

\subsubsection{NdArray::RepeatBlock()}
\label{\detokenize{cppapi/ndarray:ndarray-repeatblock}}\begin{quote}

\sphinxAtStartPar
Repeat an array a number of times along given axis.

\sphinxAtStartPar
\sphinxstylestrong{Synopsis}
\begin{quote}

\sphinxAtStartPar
\sphinxcode{\sphinxupquote{NdArray\textless{}T, N\textgreater{} RepeatBlock(size\_t repeats, int axis)}}
\end{quote}

\sphinxAtStartPar
\sphinxstylestrong{Arguments}
\begin{quote}

\sphinxAtStartPar
\sphinxcode{\sphinxupquote{repeats}}: number of repetitions.

\sphinxAtStartPar
\sphinxcode{\sphinxupquote{axis}}: axis of NdArray.
\end{quote}

\sphinxAtStartPar
\sphinxstylestrong{Return}
\begin{quote}

\sphinxAtStartPar
new NdArray object.
\end{quote}
\end{quote}

\subsubsection{NdArray::Represent()}
\label{\detokenize{cppapi/ndarray:ndarray-represent}}\begin{quote}

\sphinxAtStartPar
String representation of NdArray object.

\sphinxAtStartPar
\sphinxstylestrong{Synopsis}
\begin{quote}

\sphinxAtStartPar
\sphinxcode{\sphinxupquote{std::string Represent(size\_t maxlen)}}
\end{quote}

\sphinxAtStartPar
\sphinxstylestrong{Arguments}
\begin{quote}

\sphinxAtStartPar
\sphinxcode{\sphinxupquote{maxlen}}: max length of representation string.
\end{quote}

\sphinxAtStartPar
\sphinxstylestrong{Return}
\begin{quote}

\sphinxAtStartPar
representation string object.
\end{quote}
\end{quote}

\subsubsection{NdArray::Reshape()}
\label{\detokenize{cppapi/ndarray:ndarray-reshape}}\begin{quote}

\sphinxAtStartPar
Reshape NdArray object to new shape.

\sphinxAtStartPar
\sphinxstylestrong{Synopsis}
\begin{quote}

\sphinxAtStartPar
\sphinxcode{\sphinxupquote{template \textless{}int M\textgreater{} NdArray\textless{}T, M\textgreater{} Reshape(const Shape\textless{}M\textgreater{} \&shape)}}
\end{quote}

\sphinxAtStartPar
\sphinxstylestrong{Arguments}
\begin{quote}

\sphinxAtStartPar
\sphinxcode{\sphinxupquote{shape}}: new shape of M\sphinxhyphen{}dimensions.
\end{quote}

\sphinxAtStartPar
\sphinxstylestrong{Return}
\begin{quote}

\sphinxAtStartPar
M\sphinxhyphen{}dimensional NdArray object.
\end{quote}
\end{quote}

\subsubsection{NdArray::SetItem()}
\label{\detokenize{cppapi/ndarray:ndarray-setitem}}\begin{quote}

\sphinxAtStartPar
Set element of given index to NdArray object.

\sphinxAtStartPar
\sphinxstylestrong{Synopsis}
\begin{quote}

\sphinxAtStartPar
\sphinxcode{\sphinxupquote{void SetItem(size\_t idx, const T \&val)}}
\end{quote}

\sphinxAtStartPar
\sphinxstylestrong{Arguments}
\begin{quote}

\sphinxAtStartPar
\sphinxcode{\sphinxupquote{idx}}: index of element.

\sphinxAtStartPar
\sphinxcode{\sphinxupquote{val}}: value of element.
\end{quote}
\end{quote}

\subsubsection{NdArray::Squeeze()}
\label{\detokenize{cppapi/ndarray:ndarray-squeeze}}\begin{quote}

\sphinxAtStartPar
Remove axis of length 1 from shape of NdArray object.

\sphinxAtStartPar
\sphinxstylestrong{Synopsis}
\begin{quote}

\sphinxAtStartPar
\sphinxcode{\sphinxupquote{NdArray\textless{}T, N \sphinxhyphen{} 1\textgreater{} Squeeze(int axis)}}
\end{quote}

\sphinxAtStartPar
\sphinxstylestrong{Arguments}
\begin{quote}

\sphinxAtStartPar
\sphinxcode{\sphinxupquote{axis}}: axis of NdArray, where the length is 1.
\end{quote}

\sphinxAtStartPar
\sphinxstylestrong{Return}
\begin{quote}

\sphinxAtStartPar
(N\sphinxhyphen{}1)\sphinxhyphen{}dimensional NdArray object.
\end{quote}
\end{quote}

\subsubsection{NdArray::Sum()}
\label{\detokenize{cppapi/ndarray:ndarray-sum}}\begin{quote}

\sphinxAtStartPar
Sum of all elements in NdArray object.

\sphinxAtStartPar
\sphinxstylestrong{Synopsis}
\begin{quote}

\sphinxAtStartPar
\sphinxcode{\sphinxupquote{T Sum()}}
\end{quote}

\sphinxAtStartPar
\sphinxstylestrong{Return}
\begin{quote}

\sphinxAtStartPar
sum value.
\end{quote}
\end{quote}

\subsubsection{NdArray::Sum()}
\label{\detokenize{cppapi/ndarray:id8}}\begin{quote}

\sphinxAtStartPar
Sum of elements at given axis of NdArray object.

\sphinxAtStartPar
\sphinxstylestrong{Synopsis}
\begin{quote}

\sphinxAtStartPar
\sphinxcode{\sphinxupquote{NdArray\textless{}T, N \sphinxhyphen{} 1\textgreater{} Sum(int axis)}}
\end{quote}

\sphinxAtStartPar
\sphinxstylestrong{Arguments}
\begin{quote}

\sphinxAtStartPar
\sphinxcode{\sphinxupquote{axis}}: axis of NdArray.
\end{quote}

\sphinxAtStartPar
\sphinxstylestrong{Return}
\begin{quote}

\sphinxAtStartPar
(N\sphinxhyphen{}1)\sphinxhyphen{}dimensional NdArray object.
\end{quote}
\end{quote}

\subsubsection{NdArray::Transpose()}
\label{\detokenize{cppapi/ndarray:ndarray-transpose}}\begin{quote}

\sphinxAtStartPar
Perform matrix transpose of NdArray object.

\sphinxAtStartPar
\sphinxstylestrong{Synopsis}
\begin{quote}

\sphinxAtStartPar
\sphinxcode{\sphinxupquote{NdArray\textless{}T, N\textgreater{} Transpose()}}
\end{quote}

\sphinxAtStartPar
\sphinxstylestrong{Return}
\begin{quote}

\sphinxAtStartPar
transposed NdArray object.
\end{quote}
\end{quote}

\subsection{Shape}
\label{\detokenize{cppapiref:shape}}\label{\detokenize{cppapiref:chapcppapiref-shape}}
\sphinxAtStartPar
The Shape class encapsulates a tuple of integers, indicating the size of array
along each dimension. It refers to dimensions of built\sphinxhyphen{}in {\hyperref[\detokenize{cppapiref:chapcppapiref-ndarray}]{\sphinxcrossref{\DUrole{std,std-ref}{NdArray}}}}
in COPT. The following methods are provided:

\sphinxstepscope

\subsubsection{Shape::Expand()}
\label{\detokenize{cppapi/ndim:shape-expand}}\label{\detokenize{cppapi/ndim::doc}}\begin{quote}

\sphinxAtStartPar
Expand shape of Shape object.

\sphinxAtStartPar
\sphinxstylestrong{Synopsis}
\begin{quote}

\sphinxAtStartPar
\sphinxcode{\sphinxupquote{Shape\textless{}N + 1\textgreater{} Expand(int axis)}}
\end{quote}

\sphinxAtStartPar
\sphinxstylestrong{Arguments}
\begin{quote}

\sphinxAtStartPar
\sphinxcode{\sphinxupquote{axis}}: given axis.
\end{quote}

\sphinxAtStartPar
\sphinxstylestrong{Return}
\begin{quote}

\sphinxAtStartPar
Shape object in (N+1)\sphinxhyphen{}dimensions.
\end{quote}
\end{quote}

\subsubsection{Shape::GetDim()}
\label{\detokenize{cppapi/ndim:shape-getdim}}\begin{quote}

\sphinxAtStartPar
Get i\sphinxhyphen{}th dimension in Shape object.

\sphinxAtStartPar
\sphinxstylestrong{Synopsis}
\begin{quote}

\sphinxAtStartPar
\sphinxcode{\sphinxupquote{size\_t GetDim(int i)}}
\end{quote}

\sphinxAtStartPar
\sphinxstylestrong{Arguments}
\begin{quote}

\sphinxAtStartPar
\sphinxcode{\sphinxupquote{i}}: index of dimensions.
\end{quote}

\sphinxAtStartPar
\sphinxstylestrong{Return}
\begin{quote}

\sphinxAtStartPar
the i\sphinxhyphen{}th dimension.
\end{quote}
\end{quote}

\subsubsection{Shape::GetND()}
\label{\detokenize{cppapi/ndim:shape-getnd}}\begin{quote}

\sphinxAtStartPar
Get number of dimensions in Shape object.

\sphinxAtStartPar
\sphinxstylestrong{Synopsis}
\begin{quote}

\sphinxAtStartPar
\sphinxcode{\sphinxupquote{int GetND()}}
\end{quote}

\sphinxAtStartPar
\sphinxstylestrong{Return}
\begin{quote}

\sphinxAtStartPar
number of dimensions.
\end{quote}
\end{quote}

\subsubsection{Shape::GetSize()}
\label{\detokenize{cppapi/ndim:shape-getsize}}\begin{quote}

\sphinxAtStartPar
Get size of Shape object.

\sphinxAtStartPar
\sphinxstylestrong{Synopsis}
\begin{quote}

\sphinxAtStartPar
\sphinxcode{\sphinxupquote{size\_t GetSize()}}
\end{quote}

\sphinxAtStartPar
\sphinxstylestrong{Return}
\begin{quote}

\sphinxAtStartPar
size of shape.
\end{quote}
\end{quote}

\subsubsection{Shape::GetStart()}
\label{\detokenize{cppapi/ndim:shape-getstart}}\begin{quote}

\sphinxAtStartPar
Get the i\sphinxhyphen{}th start postion in Shape object.

\sphinxAtStartPar
\sphinxstylestrong{Synopsis}
\begin{quote}

\sphinxAtStartPar
\sphinxcode{\sphinxupquote{int64\_t GetStart(int i)}}
\end{quote}

\sphinxAtStartPar
\sphinxstylestrong{Arguments}
\begin{quote}

\sphinxAtStartPar
\sphinxcode{\sphinxupquote{i}}: index of dimensions.
\end{quote}

\sphinxAtStartPar
\sphinxstylestrong{Return}
\begin{quote}

\sphinxAtStartPar
start position in i\sphinxhyphen{}th dimension.
\end{quote}
\end{quote}

\subsubsection{Shape::GetStride()}
\label{\detokenize{cppapi/ndim:shape-getstride}}\begin{quote}

\sphinxAtStartPar
Get i\sphinxhyphen{}th stride in Shape object.

\sphinxAtStartPar
\sphinxstylestrong{Synopsis}
\begin{quote}

\sphinxAtStartPar
\sphinxcode{\sphinxupquote{int64\_t GetStride(int i)}}
\end{quote}

\sphinxAtStartPar
\sphinxstylestrong{Arguments}
\begin{quote}

\sphinxAtStartPar
\sphinxcode{\sphinxupquote{i}}: index of dimensions.
\end{quote}

\sphinxAtStartPar
\sphinxstylestrong{Return}
\begin{quote}

\sphinxAtStartPar
stride in i\sphinxhyphen{}th dimension.
\end{quote}
\end{quote}

\subsubsection{Shape::operator!=()}
\label{\detokenize{cppapi/ndim:shape-operator}}\begin{quote}

\sphinxAtStartPar
Use operator ‘!=’ to compare with other Shape object.

\sphinxAtStartPar
\sphinxstylestrong{Synopsis}
\begin{quote}

\sphinxAtStartPar
\sphinxcode{\sphinxupquote{bool operator!=(const Shape\textless{}N\textgreater{} \&other)}}
\end{quote}

\sphinxAtStartPar
\sphinxstylestrong{Arguments}
\begin{quote}

\sphinxAtStartPar
\sphinxcode{\sphinxupquote{other}}: other Shape object.
\end{quote}

\sphinxAtStartPar
\sphinxstylestrong{Return}
\begin{quote}

\sphinxAtStartPar
True if shape is not the same.
\end{quote}
\end{quote}

\subsubsection{Shape::operator==()}
\label{\detokenize{cppapi/ndim:id1}}\begin{quote}

\sphinxAtStartPar
Use operator ‘==’ to compare with other Shape object.

\sphinxAtStartPar
\sphinxstylestrong{Synopsis}
\begin{quote}

\sphinxAtStartPar
\sphinxcode{\sphinxupquote{bool operator==(const Shape\textless{}N\textgreater{} \&other)}}
\end{quote}

\sphinxAtStartPar
\sphinxstylestrong{Arguments}
\begin{quote}

\sphinxAtStartPar
\sphinxcode{\sphinxupquote{other}}: other Shape object.
\end{quote}

\sphinxAtStartPar
\sphinxstylestrong{Return}
\begin{quote}

\sphinxAtStartPar
True if shape is the same.
\end{quote}
\end{quote}

\subsubsection{Shape::Rebuild()}
\label{\detokenize{cppapi/ndim:shape-rebuild}}\begin{quote}

\sphinxAtStartPar
Rebuild Shape object, that is, keep dimensions while reset strides and starts.

\sphinxAtStartPar
\sphinxstylestrong{Synopsis}
\begin{quote}

\sphinxAtStartPar
\sphinxcode{\sphinxupquote{Shape\textless{}N\textgreater{} Rebuild()}}
\end{quote}

\sphinxAtStartPar
\sphinxstylestrong{Return}
\begin{quote}

\sphinxAtStartPar
new Shape object.
\end{quote}
\end{quote}

\subsubsection{Shape::Represent()}
\label{\detokenize{cppapi/ndim:shape-represent}}\begin{quote}

\sphinxAtStartPar
String representation of Shape object.

\sphinxAtStartPar
\sphinxstylestrong{Synopsis}
\begin{quote}

\sphinxAtStartPar
\sphinxcode{\sphinxupquote{std::string Represent(int type)}}
\end{quote}

\sphinxAtStartPar
\sphinxstylestrong{Arguments}
\begin{quote}

\sphinxAtStartPar
\sphinxcode{\sphinxupquote{type}}: 0: dimensions; 1: strides; 2: starts.
\end{quote}

\sphinxAtStartPar
\sphinxstylestrong{Return}
\begin{quote}

\sphinxAtStartPar
string object.
\end{quote}
\end{quote}

\subsubsection{Shape::Squeeze()}
\label{\detokenize{cppapi/ndim:shape-squeeze}}\begin{quote}

\sphinxAtStartPar
Remove axis of length 1 from Shape object.

\sphinxAtStartPar
\sphinxstylestrong{Synopsis}
\begin{quote}

\sphinxAtStartPar
\sphinxcode{\sphinxupquote{Shape\textless{}N \sphinxhyphen{} 1\textgreater{} Squeeze(int axis)}}
\end{quote}

\sphinxAtStartPar
\sphinxstylestrong{Arguments}
\begin{quote}

\sphinxAtStartPar
\sphinxcode{\sphinxupquote{axis}}: given axis, where the length is 1.
\end{quote}

\sphinxAtStartPar
\sphinxstylestrong{Return}
\begin{quote}

\sphinxAtStartPar
Shape object in (N\sphinxhyphen{}1)\sphinxhyphen{}dimensions.
\end{quote}
\end{quote}

\subsection{View}
\label{\detokenize{cppapiref:view}}\label{\detokenize{cppapiref:chapcppapiref-view}}
\sphinxAtStartPar
The View class is used to perform slicing operations on multi\sphinxhyphen{}dimensional arrays.
The following methods are provided:

\sphinxstepscope

\subsubsection{View::View()}
\label{\detokenize{cppapi/view:view-view}}\label{\detokenize{cppapi/view::doc}}\begin{quote}

\sphinxAtStartPar
Create an empty View object. Call make\_view() for your convenience.

\sphinxAtStartPar
\sphinxstylestrong{Synopsis}
\begin{quote}

\sphinxAtStartPar
\sphinxcode{\sphinxupquote{View()}}
\end{quote}
\end{quote}

\subsubsection{View::AddFull()}
\label{\detokenize{cppapi/view:view-addfull}}\begin{quote}

\sphinxAtStartPar
Create full view object at current dimension.

\sphinxAtStartPar
\sphinxstylestrong{Synopsis}
\begin{quote}

\sphinxAtStartPar
\sphinxcode{\sphinxupquote{View \&AddFull()}}
\end{quote}

\sphinxAtStartPar
\sphinxstylestrong{Return}
\begin{quote}

\sphinxAtStartPar
View object.
\end{quote}
\end{quote}

\subsubsection{View::AddScalar()}
\label{\detokenize{cppapi/view:view-addscalar}}\begin{quote}

\sphinxAtStartPar
Create view object of given index at current dimension.

\sphinxAtStartPar
\sphinxstylestrong{Synopsis}
\begin{quote}

\sphinxAtStartPar
\sphinxcode{\sphinxupquote{View \&AddScalar(int64\_t n)}}
\end{quote}

\sphinxAtStartPar
\sphinxstylestrong{Arguments}
\begin{quote}

\sphinxAtStartPar
\sphinxcode{\sphinxupquote{n}}: given index.
\end{quote}

\sphinxAtStartPar
\sphinxstylestrong{Return}
\begin{quote}

\sphinxAtStartPar
View object.
\end{quote}
\end{quote}

\subsubsection{View::AddSlice()}
\label{\detokenize{cppapi/view:view-addslice}}\begin{quote}

\sphinxAtStartPar
Create view object of slice at current dimension.

\sphinxAtStartPar
\sphinxstylestrong{Synopsis}
\begin{quote}

\sphinxAtStartPar
\sphinxcode{\sphinxupquote{View \&AddSlice(int64\_t start)}}
\end{quote}

\sphinxAtStartPar
\sphinxstylestrong{Arguments}
\begin{quote}

\sphinxAtStartPar
\sphinxcode{\sphinxupquote{start}}: start index, inclusive.
\end{quote}

\sphinxAtStartPar
\sphinxstylestrong{Return}
\begin{quote}

\sphinxAtStartPar
View object.
\end{quote}
\end{quote}

\subsubsection{View::AddSlice()}
\label{\detokenize{cppapi/view:id1}}\begin{quote}

\sphinxAtStartPar
Create view object of slice at current dimension.

\sphinxAtStartPar
\sphinxstylestrong{Synopsis}
\begin{quote}

\sphinxAtStartPar
\sphinxcode{\sphinxupquote{View \&AddSlice(int64\_t start, int64\_t stop)}}
\end{quote}

\sphinxAtStartPar
\sphinxstylestrong{Arguments}
\begin{quote}

\sphinxAtStartPar
\sphinxcode{\sphinxupquote{start}}: start index, inclusive.

\sphinxAtStartPar
\sphinxcode{\sphinxupquote{stop}}: stop index, exclusive.
\end{quote}

\sphinxAtStartPar
\sphinxstylestrong{Return}
\begin{quote}

\sphinxAtStartPar
View object.
\end{quote}
\end{quote}

\subsubsection{View::AddSlice()}
\label{\detokenize{cppapi/view:id2}}\begin{quote}

\sphinxAtStartPar
Create view object of slice at current dimension.

\sphinxAtStartPar
\sphinxstylestrong{Synopsis}
\begin{quote}

\sphinxAtStartPar
\sphinxcode{\sphinxupquote{View \&AddSlice(}}
\begin{quote}

\sphinxAtStartPar
\sphinxcode{\sphinxupquote{int64\_t start,}}

\sphinxAtStartPar
\sphinxcode{\sphinxupquote{int64\_t stop,}}

\sphinxAtStartPar
\sphinxcode{\sphinxupquote{int64\_t step)}}
\end{quote}
\end{quote}

\sphinxAtStartPar
\sphinxstylestrong{Arguments}
\begin{quote}

\sphinxAtStartPar
\sphinxcode{\sphinxupquote{start}}: start index, inclusive.

\sphinxAtStartPar
\sphinxcode{\sphinxupquote{stop}}: stop index, exclusive.

\sphinxAtStartPar
\sphinxcode{\sphinxupquote{step}}: step size between start and stop index. It can be negative.
\end{quote}

\sphinxAtStartPar
\sphinxstylestrong{Return}
\begin{quote}

\sphinxAtStartPar
View object.
\end{quote}
\end{quote}

\subsubsection{View::operator()()}
\label{\detokenize{cppapi/view:view-operator}}\begin{quote}

\sphinxAtStartPar
Create full view object at current dimension.

\sphinxAtStartPar
\sphinxstylestrong{Synopsis}
\begin{quote}

\sphinxAtStartPar
\sphinxcode{\sphinxupquote{View \&operator()()}}
\end{quote}

\sphinxAtStartPar
\sphinxstylestrong{Return}
\begin{quote}

\sphinxAtStartPar
View object.
\end{quote}
\end{quote}

\subsubsection{View::operator()()}
\label{\detokenize{cppapi/view:id3}}\begin{quote}

\sphinxAtStartPar
Create view object of slice at current dimension.

\sphinxAtStartPar
\sphinxstylestrong{Synopsis}
\begin{quote}

\sphinxAtStartPar
\sphinxcode{\sphinxupquote{View \&operator()(int64\_t start, int64\_t stop)}}
\end{quote}

\sphinxAtStartPar
\sphinxstylestrong{Arguments}
\begin{quote}

\sphinxAtStartPar
\sphinxcode{\sphinxupquote{start}}: start index, inclusive.

\sphinxAtStartPar
\sphinxcode{\sphinxupquote{stop}}: stop index, exclusive.
\end{quote}

\sphinxAtStartPar
\sphinxstylestrong{Return}
\begin{quote}

\sphinxAtStartPar
View object.
\end{quote}
\end{quote}

\subsubsection{View::operator()()}
\label{\detokenize{cppapi/view:id4}}\begin{quote}

\sphinxAtStartPar
Create view object of slice at current dimension.

\sphinxAtStartPar
\sphinxstylestrong{Synopsis}
\begin{quote}

\sphinxAtStartPar
\sphinxcode{\sphinxupquote{View \&operator()(}}
\begin{quote}

\sphinxAtStartPar
\sphinxcode{\sphinxupquote{int64\_t start,}}

\sphinxAtStartPar
\sphinxcode{\sphinxupquote{int64\_t stop,}}

\sphinxAtStartPar
\sphinxcode{\sphinxupquote{int64\_t step)}}
\end{quote}
\end{quote}

\sphinxAtStartPar
\sphinxstylestrong{Arguments}
\begin{quote}

\sphinxAtStartPar
\sphinxcode{\sphinxupquote{start}}: start index, inclusive.

\sphinxAtStartPar
\sphinxcode{\sphinxupquote{stop}}: stop index, exclusive.

\sphinxAtStartPar
\sphinxcode{\sphinxupquote{step}}: step size between start and stop index. It can be negative.
\end{quote}

\sphinxAtStartPar
\sphinxstylestrong{Return}
\begin{quote}

\sphinxAtStartPar
View object.
\end{quote}
\end{quote}

\subsubsection{View::operator()()}
\label{\detokenize{cppapi/view:id5}}\begin{quote}

\sphinxAtStartPar
Create view object of given index at current dimension.

\sphinxAtStartPar
\sphinxstylestrong{Synopsis}
\begin{quote}

\sphinxAtStartPar
\sphinxcode{\sphinxupquote{View \&operator()(int64\_t n)}}
\end{quote}

\sphinxAtStartPar
\sphinxstylestrong{Arguments}
\begin{quote}

\sphinxAtStartPar
\sphinxcode{\sphinxupquote{n}}: given index.
\end{quote}

\sphinxAtStartPar
\sphinxstylestrong{Return}
\begin{quote}

\sphinxAtStartPar
View object.
\end{quote}
\end{quote}

\subsection{CallbackBase}
\label{\detokenize{cppapiref:callbackbase}}\label{\detokenize{cppapiref:chapcppapiref-callbackbase}}
\sphinxAtStartPar
COPT Callback abstract base object. Users must implment its virtual method
\sphinxcode{\sphinxupquote{virtual void CallbackBase::callback()}} to instantiate an instance, which
pass to \sphinxcode{\sphinxupquote{Model::SetCallback(ICallback* pcb, int cbctx)}} as the first
parameter. Subclass of CallbackBase inherits the following methods:

\sphinxstepscope

\subsubsection{CallbackBase::AddLazyConstr()}
\label{\detokenize{cppapi/callbackbase:callbackbase-addlazyconstr}}\label{\detokenize{cppapi/callbackbase::doc}}\begin{quote}

\sphinxAtStartPar
Add a lazy constraint to model.

\sphinxAtStartPar
\sphinxstylestrong{Synopsis}
\begin{quote}

\sphinxAtStartPar
\sphinxcode{\sphinxupquote{void AddLazyConstr(}}
\begin{quote}

\sphinxAtStartPar
\sphinxcode{\sphinxupquote{const Expr \&lhs,}}

\sphinxAtStartPar
\sphinxcode{\sphinxupquote{char sense,}}

\sphinxAtStartPar
\sphinxcode{\sphinxupquote{double rhs)}}
\end{quote}
\end{quote}

\sphinxAtStartPar
\sphinxstylestrong{Arguments}
\begin{quote}

\sphinxAtStartPar
\sphinxcode{\sphinxupquote{lhs}}: expression for lazy contraint.

\sphinxAtStartPar
\sphinxcode{\sphinxupquote{sense}}: sense for lazy constraint.

\sphinxAtStartPar
\sphinxcode{\sphinxupquote{rhs}}: right hand side value for lazy constraint.
\end{quote}
\end{quote}

\subsubsection{CallbackBase::AddLazyConstr()}
\label{\detokenize{cppapi/callbackbase:id1}}\begin{quote}

\sphinxAtStartPar
Add a lazy constraint to model.

\sphinxAtStartPar
\sphinxstylestrong{Synopsis}
\begin{quote}

\sphinxAtStartPar
\sphinxcode{\sphinxupquote{void AddLazyConstr(}}
\begin{quote}

\sphinxAtStartPar
\sphinxcode{\sphinxupquote{const Expr \&lhs,}}

\sphinxAtStartPar
\sphinxcode{\sphinxupquote{char sense,}}

\sphinxAtStartPar
\sphinxcode{\sphinxupquote{const Expr \&rhs)}}
\end{quote}
\end{quote}

\sphinxAtStartPar
\sphinxstylestrong{Arguments}
\begin{quote}

\sphinxAtStartPar
\sphinxcode{\sphinxupquote{lhs}}: left hand side expression for lazy contraint.

\sphinxAtStartPar
\sphinxcode{\sphinxupquote{sense}}: sense for lazy constraint.

\sphinxAtStartPar
\sphinxcode{\sphinxupquote{rhs}}: right hand side expression for lazy contraint.
\end{quote}
\end{quote}

\subsubsection{CallbackBase::AddLazyConstr()}
\label{\detokenize{cppapi/callbackbase:id2}}\begin{quote}

\sphinxAtStartPar
Add a lazy constraint to model.

\sphinxAtStartPar
\sphinxstylestrong{Synopsis}
\begin{quote}

\sphinxAtStartPar
\sphinxcode{\sphinxupquote{void AddLazyConstr(const ConstrBuilder \&builder)}}
\end{quote}

\sphinxAtStartPar
\sphinxstylestrong{Arguments}
\begin{quote}

\sphinxAtStartPar
\sphinxcode{\sphinxupquote{builder}}: builder for lazy contraint.
\end{quote}
\end{quote}

\subsubsection{CallbackBase::AddLazyConstrs()}
\label{\detokenize{cppapi/callbackbase:callbackbase-addlazyconstrs}}\begin{quote}

\sphinxAtStartPar
Add lazy constraints to model.

\sphinxAtStartPar
\sphinxstylestrong{Synopsis}
\begin{quote}

\sphinxAtStartPar
\sphinxcode{\sphinxupquote{void AddLazyConstrs(const ConstrBuilderArray \&builders)}}
\end{quote}

\sphinxAtStartPar
\sphinxstylestrong{Arguments}
\begin{quote}

\sphinxAtStartPar
\sphinxcode{\sphinxupquote{builders}}: array of builders for lazy contraints.
\end{quote}
\end{quote}

\subsubsection{CallbackBase::AddUserCut()}
\label{\detokenize{cppapi/callbackbase:callbackbase-addusercut}}\begin{quote}

\sphinxAtStartPar
Add a user cut to model.

\sphinxAtStartPar
\sphinxstylestrong{Synopsis}
\begin{quote}

\sphinxAtStartPar
\sphinxcode{\sphinxupquote{void AddUserCut(}}
\begin{quote}

\sphinxAtStartPar
\sphinxcode{\sphinxupquote{const Expr \&lhs,}}

\sphinxAtStartPar
\sphinxcode{\sphinxupquote{char sense,}}

\sphinxAtStartPar
\sphinxcode{\sphinxupquote{double rhs)}}
\end{quote}
\end{quote}

\sphinxAtStartPar
\sphinxstylestrong{Arguments}
\begin{quote}

\sphinxAtStartPar
\sphinxcode{\sphinxupquote{lhs}}: expression for user cut.

\sphinxAtStartPar
\sphinxcode{\sphinxupquote{sense}}: sense for user cut.

\sphinxAtStartPar
\sphinxcode{\sphinxupquote{rhs}}: right hand side value for user cut.
\end{quote}
\end{quote}

\subsubsection{CallbackBase::AddUserCut()}
\label{\detokenize{cppapi/callbackbase:id3}}\begin{quote}

\sphinxAtStartPar
Add a user cut to model.

\sphinxAtStartPar
\sphinxstylestrong{Synopsis}
\begin{quote}

\sphinxAtStartPar
\sphinxcode{\sphinxupquote{void AddUserCut(}}
\begin{quote}

\sphinxAtStartPar
\sphinxcode{\sphinxupquote{const Expr \&lhs,}}

\sphinxAtStartPar
\sphinxcode{\sphinxupquote{char sense,}}

\sphinxAtStartPar
\sphinxcode{\sphinxupquote{const Expr \&rhs)}}
\end{quote}
\end{quote}

\sphinxAtStartPar
\sphinxstylestrong{Arguments}
\begin{quote}

\sphinxAtStartPar
\sphinxcode{\sphinxupquote{lhs}}: left hand side expression for user cut.

\sphinxAtStartPar
\sphinxcode{\sphinxupquote{sense}}: sense for user cut.

\sphinxAtStartPar
\sphinxcode{\sphinxupquote{rhs}}: right hand side expression for user cut.
\end{quote}
\end{quote}

\subsubsection{CallbackBase::AddUserCut()}
\label{\detokenize{cppapi/callbackbase:id4}}\begin{quote}

\sphinxAtStartPar
Add a user cut to model.

\sphinxAtStartPar
\sphinxstylestrong{Synopsis}
\begin{quote}

\sphinxAtStartPar
\sphinxcode{\sphinxupquote{void AddUserCut(const ConstrBuilder \&builder)}}
\end{quote}

\sphinxAtStartPar
\sphinxstylestrong{Arguments}
\begin{quote}

\sphinxAtStartPar
\sphinxcode{\sphinxupquote{builder}}: builder for user cut.
\end{quote}
\end{quote}

\subsubsection{CallbackBase::AddUserCuts()}
\label{\detokenize{cppapi/callbackbase:callbackbase-addusercuts}}\begin{quote}

\sphinxAtStartPar
Add user cuts to model.

\sphinxAtStartPar
\sphinxstylestrong{Synopsis}
\begin{quote}

\sphinxAtStartPar
\sphinxcode{\sphinxupquote{void AddUserCuts(const ConstrBuilderArray \&builders)}}
\end{quote}

\sphinxAtStartPar
\sphinxstylestrong{Arguments}
\begin{quote}

\sphinxAtStartPar
\sphinxcode{\sphinxupquote{builders}}: array of builders for user cuts.
\end{quote}
\end{quote}

\subsubsection{CallbackBase::GetDblInfo()}
\label{\detokenize{cppapi/callbackbase:callbackbase-getdblinfo}}\begin{quote}

\sphinxAtStartPar
Get double value of given information name in callback.

\sphinxAtStartPar
\sphinxstylestrong{Synopsis}
\begin{quote}

\sphinxAtStartPar
\sphinxcode{\sphinxupquote{double GetDblInfo(const char *cbinfo)}}
\end{quote}

\sphinxAtStartPar
\sphinxstylestrong{Arguments}
\begin{quote}

\sphinxAtStartPar
\sphinxcode{\sphinxupquote{cbinfo}}: name of callback info.
\end{quote}

\sphinxAtStartPar
\sphinxstylestrong{Return}
\begin{quote}

\sphinxAtStartPar
value of desired information.
\end{quote}
\end{quote}

\subsubsection{CallbackBase::GetIncumbent()}
\label{\detokenize{cppapi/callbackbase:callbackbase-getincumbent}}\begin{quote}

\sphinxAtStartPar
Get best feasible solution of given variable in callback.

\sphinxAtStartPar
\sphinxstylestrong{Synopsis}
\begin{quote}

\sphinxAtStartPar
\sphinxcode{\sphinxupquote{double GetIncumbent(Var \&var)}}
\end{quote}

\sphinxAtStartPar
\sphinxstylestrong{Arguments}
\begin{quote}

\sphinxAtStartPar
\sphinxcode{\sphinxupquote{var}}: given variable.
\end{quote}

\sphinxAtStartPar
\sphinxstylestrong{Return}
\begin{quote}

\sphinxAtStartPar
best feasible solution of given variable.
\end{quote}
\end{quote}

\subsubsection{CallbackBase::GetIncumbent()}
\label{\detokenize{cppapi/callbackbase:id5}}\begin{quote}

\sphinxAtStartPar
Get best feasible solution of variables in callback.

\sphinxAtStartPar
\sphinxstylestrong{Synopsis}
\begin{quote}

\sphinxAtStartPar
\sphinxcode{\sphinxupquote{int GetIncumbent(VarArray \&vars, double *pOut)}}
\end{quote}

\sphinxAtStartPar
\sphinxstylestrong{Arguments}
\begin{quote}

\sphinxAtStartPar
\sphinxcode{\sphinxupquote{vars}}: an array of variables.

\sphinxAtStartPar
\sphinxcode{\sphinxupquote{pOut}}: best feasible solution of desired variables.
\end{quote}

\sphinxAtStartPar
\sphinxstylestrong{Return}
\begin{quote}

\sphinxAtStartPar
the number of valid variables. If failed, return \sphinxhyphen{}1.
\end{quote}
\end{quote}

\subsubsection{CallbackBase::GetIncumbent()}
\label{\detokenize{cppapi/callbackbase:id6}}\begin{quote}

\sphinxAtStartPar
Get best feasible solution of all variables in callback.

\sphinxAtStartPar
\sphinxstylestrong{Synopsis}
\begin{quote}

\sphinxAtStartPar
\sphinxcode{\sphinxupquote{int GetIncumbent(double *pOut, int len)}}
\end{quote}

\sphinxAtStartPar
\sphinxstylestrong{Arguments}
\begin{quote}

\sphinxAtStartPar
\sphinxcode{\sphinxupquote{pOut}}: optional, output best feasible solution of all variables.

\sphinxAtStartPar
\sphinxcode{\sphinxupquote{len}}: the length of output array. The solution is written up to number of len.
\end{quote}

\sphinxAtStartPar
\sphinxstylestrong{Return}
\begin{quote}

\sphinxAtStartPar
number of all variables. Return \sphinxhyphen{}1 if error occurs.
\end{quote}
\end{quote}

\subsubsection{CallbackBase::GetIntInfo()}
\label{\detokenize{cppapi/callbackbase:callbackbase-getintinfo}}\begin{quote}

\sphinxAtStartPar
Get integer value of given information name in callback.

\sphinxAtStartPar
\sphinxstylestrong{Synopsis}
\begin{quote}

\sphinxAtStartPar
\sphinxcode{\sphinxupquote{int GetIntInfo(const char *cbinfo)}}
\end{quote}

\sphinxAtStartPar
\sphinxstylestrong{Arguments}
\begin{quote}

\sphinxAtStartPar
\sphinxcode{\sphinxupquote{cbinfo}}: name of callback info.
\end{quote}

\sphinxAtStartPar
\sphinxstylestrong{Return}
\begin{quote}

\sphinxAtStartPar
value of desired information.
\end{quote}
\end{quote}

\subsubsection{CallbackBase::GetRelaxSol()}
\label{\detokenize{cppapi/callbackbase:callbackbase-getrelaxsol}}\begin{quote}

\sphinxAtStartPar
Get LP\sphinxhyphen{}relaxation solution of given variable in callback.

\sphinxAtStartPar
\sphinxstylestrong{Synopsis}
\begin{quote}

\sphinxAtStartPar
\sphinxcode{\sphinxupquote{double GetRelaxSol(Var \&var)}}
\end{quote}

\sphinxAtStartPar
\sphinxstylestrong{Arguments}
\begin{quote}

\sphinxAtStartPar
\sphinxcode{\sphinxupquote{var}}: given variable.
\end{quote}

\sphinxAtStartPar
\sphinxstylestrong{Return}
\begin{quote}

\sphinxAtStartPar
LP\sphinxhyphen{}relaxation solution of given variable.
\end{quote}
\end{quote}

\subsubsection{CallbackBase::GetRelaxSol()}
\label{\detokenize{cppapi/callbackbase:id7}}\begin{quote}

\sphinxAtStartPar
Get LP\sphinxhyphen{}relaxation solution of variables in callback.

\sphinxAtStartPar
\sphinxstylestrong{Synopsis}
\begin{quote}

\sphinxAtStartPar
\sphinxcode{\sphinxupquote{int GetRelaxSol(VarArray \&vars, double *pOut)}}
\end{quote}

\sphinxAtStartPar
\sphinxstylestrong{Arguments}
\begin{quote}

\sphinxAtStartPar
\sphinxcode{\sphinxupquote{vars}}: an array of variables.

\sphinxAtStartPar
\sphinxcode{\sphinxupquote{pOut}}: LP\sphinxhyphen{}relaxation solution of desired variables.
\end{quote}

\sphinxAtStartPar
\sphinxstylestrong{Return}
\begin{quote}

\sphinxAtStartPar
the number of valid variables. If failed, return \sphinxhyphen{}1.
\end{quote}
\end{quote}

\subsubsection{CallbackBase::GetRelaxSol()}
\label{\detokenize{cppapi/callbackbase:id8}}\begin{quote}

\sphinxAtStartPar
Get LP\sphinxhyphen{}relaxation solution of all variables in callback.

\sphinxAtStartPar
\sphinxstylestrong{Synopsis}
\begin{quote}

\sphinxAtStartPar
\sphinxcode{\sphinxupquote{int GetRelaxSol(double *pOut, int len)}}
\end{quote}

\sphinxAtStartPar
\sphinxstylestrong{Arguments}
\begin{quote}

\sphinxAtStartPar
\sphinxcode{\sphinxupquote{pOut}}: optional, output LP\sphinxhyphen{}relaxation solution of all variables.

\sphinxAtStartPar
\sphinxcode{\sphinxupquote{len}}: the length of output array. The solution is written up to number of len.
\end{quote}

\sphinxAtStartPar
\sphinxstylestrong{Return}
\begin{quote}

\sphinxAtStartPar
number of all variables. Return \sphinxhyphen{}1 if error occurs.
\end{quote}
\end{quote}

\subsubsection{CallbackBase::GetSolution()}
\label{\detokenize{cppapi/callbackbase:callbackbase-getsolution}}\begin{quote}

\sphinxAtStartPar
Get solution of given variable in callback.

\sphinxAtStartPar
\sphinxstylestrong{Synopsis}
\begin{quote}

\sphinxAtStartPar
\sphinxcode{\sphinxupquote{double GetSolution(Var \&var)}}
\end{quote}

\sphinxAtStartPar
\sphinxstylestrong{Arguments}
\begin{quote}

\sphinxAtStartPar
\sphinxcode{\sphinxupquote{var}}: given variable.
\end{quote}

\sphinxAtStartPar
\sphinxstylestrong{Return}
\begin{quote}

\sphinxAtStartPar
solution of desired variable.
\end{quote}
\end{quote}

\subsubsection{CallbackBase::GetSolution()}
\label{\detokenize{cppapi/callbackbase:id9}}\begin{quote}

\sphinxAtStartPar
Get solution of variables in callback.

\sphinxAtStartPar
\sphinxstylestrong{Synopsis}
\begin{quote}

\sphinxAtStartPar
\sphinxcode{\sphinxupquote{int GetSolution(VarArray \&vars, double *pOut)}}
\end{quote}

\sphinxAtStartPar
\sphinxstylestrong{Arguments}
\begin{quote}

\sphinxAtStartPar
\sphinxcode{\sphinxupquote{vars}}: an array of variables.

\sphinxAtStartPar
\sphinxcode{\sphinxupquote{pOut}}: solution of desired variables.
\end{quote}

\sphinxAtStartPar
\sphinxstylestrong{Return}
\begin{quote}

\sphinxAtStartPar
the number of valid variables. If failed, return \sphinxhyphen{}1.
\end{quote}
\end{quote}

\subsubsection{CallbackBase::GetSolution()}
\label{\detokenize{cppapi/callbackbase:id10}}\begin{quote}

\sphinxAtStartPar
Get solution of all variables in callback.

\sphinxAtStartPar
\sphinxstylestrong{Synopsis}
\begin{quote}

\sphinxAtStartPar
\sphinxcode{\sphinxupquote{int GetSolution(double *pOut, int len)}}
\end{quote}

\sphinxAtStartPar
\sphinxstylestrong{Arguments}
\begin{quote}

\sphinxAtStartPar
\sphinxcode{\sphinxupquote{pOut}}: optional, output solution of all variables.

\sphinxAtStartPar
\sphinxcode{\sphinxupquote{len}}: the length of output array. The solution is written up to number of len.
\end{quote}

\sphinxAtStartPar
\sphinxstylestrong{Return}
\begin{quote}

\sphinxAtStartPar
number of all variables. Return \sphinxhyphen{}1 if error occurs.
\end{quote}
\end{quote}

\subsubsection{CallbackBase::Interrupt()}
\label{\detokenize{cppapi/callbackbase:callbackbase-interrupt}}\begin{quote}

\sphinxAtStartPar
Interrupt solving problems in callback

\sphinxAtStartPar
\sphinxstylestrong{Synopsis}
\begin{quote}

\sphinxAtStartPar
\sphinxcode{\sphinxupquote{void Interrupt()}}
\end{quote}
\end{quote}

\subsubsection{CallbackBase::LoadSolution()}
\label{\detokenize{cppapi/callbackbase:callbackbase-loadsolution}}\begin{quote}

\sphinxAtStartPar
Load customized solution to model.

\sphinxAtStartPar
\sphinxstylestrong{Synopsis}
\begin{quote}

\sphinxAtStartPar
\sphinxcode{\sphinxupquote{double LoadSolution()}}
\end{quote}

\sphinxAtStartPar
\sphinxstylestrong{Return}
\begin{quote}

\sphinxAtStartPar
objective value of given solution.
\end{quote}
\end{quote}

\subsubsection{CallbackBase::SetSolution()}
\label{\detokenize{cppapi/callbackbase:callbackbase-setsolution}}\begin{quote}

\sphinxAtStartPar
Set solution of a given variable in callback.

\sphinxAtStartPar
\sphinxstylestrong{Synopsis}
\begin{quote}

\sphinxAtStartPar
\sphinxcode{\sphinxupquote{void SetSolution(Var \&var, double val)}}
\end{quote}

\sphinxAtStartPar
\sphinxstylestrong{Arguments}
\begin{quote}

\sphinxAtStartPar
\sphinxcode{\sphinxupquote{var}}: a variable object.

\sphinxAtStartPar
\sphinxcode{\sphinxupquote{val}}: double value.
\end{quote}
\end{quote}

\subsubsection{CallbackBase::SetSolution()}
\label{\detokenize{cppapi/callbackbase:id11}}\begin{quote}

\sphinxAtStartPar
Set solution of variables in callback.

\sphinxAtStartPar
\sphinxstylestrong{Synopsis}
\begin{quote}

\sphinxAtStartPar
\sphinxcode{\sphinxupquote{void SetSolution(}}
\begin{quote}

\sphinxAtStartPar
\sphinxcode{\sphinxupquote{VarArray \&vars,}}

\sphinxAtStartPar
\sphinxcode{\sphinxupquote{const double *vals,}}

\sphinxAtStartPar
\sphinxcode{\sphinxupquote{int len)}}
\end{quote}
\end{quote}

\sphinxAtStartPar
\sphinxstylestrong{Arguments}
\begin{quote}

\sphinxAtStartPar
\sphinxcode{\sphinxupquote{vars}}: an array of variable objects.

\sphinxAtStartPar
\sphinxcode{\sphinxupquote{vals}}: an array of double values.

\sphinxAtStartPar
\sphinxcode{\sphinxupquote{len}}: length of array of double values.
\end{quote}
\end{quote}

\subsubsection{CallbackBase::Where()}
\label{\detokenize{cppapi/callbackbase:callbackbase-where}}\begin{quote}

\sphinxAtStartPar
Get context in callback.

\sphinxAtStartPar
\sphinxstylestrong{Synopsis}
\begin{quote}

\sphinxAtStartPar
\sphinxcode{\sphinxupquote{int Where()}}
\end{quote}

\sphinxAtStartPar
\sphinxstylestrong{Return}
\begin{quote}

\sphinxAtStartPar
integer value of context.
\end{quote}
\end{quote}

\subsection{NlpCallbackBase}
\label{\detokenize{cppapiref:nlpcallbackbase}}\label{\detokenize{cppapiref:chapcppapiref-nlpcallbackbase}}
\sphinxAtStartPar
The NlpCallbackBase class provides an interface for defining evaluation callbacks
for nonlinear optimization models in COPT.

\sphinxAtStartPar
By inheriting this class and implementing the corresponding member methods,
users provide objective and constraint value evaluations, together with their
first\sphinxhyphen{}order and second\sphinxhyphen{}order derivative information to COPT during nonlinear optimization.

\sphinxstepscope

\subsubsection{NlpCallbackBase::EvalObj(INdArray\textless{}double, 1\textgreater{}* xdata, INdArray\textless{}double, 1\textgreater{}* outdata)}
\label{\detokenize{cppapi/nlpcallbackbase:nlpcallbackbase-evalobj-indarray-double-1-xdata-indarray-double-1-outdata}}\label{\detokenize{cppapi/nlpcallbackbase::doc}}\begin{quote}

\sphinxAtStartPar
Evaluate objective of nonlinear model.

\sphinxAtStartPar
\sphinxstylestrong{Synopsis}
\begin{quote}

\sphinxAtStartPar
\sphinxcode{\sphinxupquote{void EvalObj(}}
\begin{quote}

\sphinxAtStartPar
\sphinxcode{\sphinxupquote{INdArray\textless{}double, 1\textgreater{} *xdata,}}

\sphinxAtStartPar
\sphinxcode{\sphinxupquote{INdArray\textless{}double, 1\textgreater{} *outdata)}}
\end{quote}
\end{quote}

\sphinxAtStartPar
\sphinxstylestrong{Arguments}
\begin{quote}

\sphinxAtStartPar
\sphinxcode{\sphinxupquote{xdata}}: pointer to NdArray of variable values.

\sphinxAtStartPar
\sphinxcode{\sphinxupquote{outdata}}: pointer to output NdArray of objective value.
\end{quote}
\end{quote}

\subsubsection{NlpCallbackBase::EvalGrad(INdArray\textless{}double, 1\textgreater{}* xdata, INdArray\textless{}double, 1\textgreater{}* outdata)}
\label{\detokenize{cppapi/nlpcallbackbase:nlpcallbackbase-evalgrad-indarray-double-1-xdata-indarray-double-1-outdata}}\begin{quote}

\sphinxAtStartPar
Evaluate objective gradient of nonlinear model.

\sphinxAtStartPar
\sphinxstylestrong{Synopsis}
\begin{quote}

\sphinxAtStartPar
\sphinxcode{\sphinxupquote{void EvalGrad(}}
\begin{quote}

\sphinxAtStartPar
\sphinxcode{\sphinxupquote{INdArray\textless{}double, 1\textgreater{} *xdata,}}

\sphinxAtStartPar
\sphinxcode{\sphinxupquote{INdArray\textless{}double, 1\textgreater{} *outdata)}}
\end{quote}
\end{quote}

\sphinxAtStartPar
\sphinxstylestrong{Arguments}
\begin{quote}

\sphinxAtStartPar
\sphinxcode{\sphinxupquote{xdata}}: pointer to NdArray of variable values.

\sphinxAtStartPar
\sphinxcode{\sphinxupquote{outdata}}: pointer to output NdArray for objective gradient.
\end{quote}
\end{quote}

\subsubsection{NlpCallbackBase::EvalCon(INdArray\textless{}double, 1\textgreater{}* xdata, INdArray\textless{}double, 1\textgreater{}* outdata)}
\label{\detokenize{cppapi/nlpcallbackbase:nlpcallbackbase-evalcon-indarray-double-1-xdata-indarray-double-1-outdata}}\begin{quote}

\sphinxAtStartPar
Evaluate constraint values of nonlinear model.

\sphinxAtStartPar
\sphinxstylestrong{Synopsis}
\begin{quote}

\sphinxAtStartPar
\sphinxcode{\sphinxupquote{void EvalCon(}}
\begin{quote}

\sphinxAtStartPar
\sphinxcode{\sphinxupquote{INdArray\textless{}double, 1\textgreater{} *xdata,}}

\sphinxAtStartPar
\sphinxcode{\sphinxupquote{INdArray\textless{}double, 1\textgreater{} *outdata)}}
\end{quote}
\end{quote}

\sphinxAtStartPar
\sphinxstylestrong{Arguments}
\begin{quote}

\sphinxAtStartPar
\sphinxcode{\sphinxupquote{xdata}}: pointer to NdArray of variable values.

\sphinxAtStartPar
\sphinxcode{\sphinxupquote{outdata}}: pointer to output NdArray for constraint values.
\end{quote}
\end{quote}

\subsubsection{NlpCallbackBase::EvalJac(INdArray\textless{}double, 1\textgreater{}* xdata, INdArray\textless{}double, 1\textgreater{}* outdata)}
\label{\detokenize{cppapi/nlpcallbackbase:nlpcallbackbase-evaljac-indarray-double-1-xdata-indarray-double-1-outdata}}\begin{quote}

\sphinxAtStartPar
Evaluate Jacobian values of the nonlinear constraint functions.

\sphinxAtStartPar
\sphinxstylestrong{Synopsis}
\begin{quote}

\sphinxAtStartPar
\sphinxcode{\sphinxupquote{void EvalJac(}}
\begin{quote}

\sphinxAtStartPar
\sphinxcode{\sphinxupquote{INdArray\textless{}double, 1\textgreater{} *xdata,}}

\sphinxAtStartPar
\sphinxcode{\sphinxupquote{INdArray\textless{}double, 1\textgreater{} *outdata)}}
\end{quote}
\end{quote}

\sphinxAtStartPar
\sphinxstylestrong{Arguments}
\begin{quote}

\sphinxAtStartPar
\sphinxcode{\sphinxupquote{xdata}}: pointer to NdArray of variable values.

\sphinxAtStartPar
\sphinxcode{\sphinxupquote{outdata}}: pointer to output NdArray for nonzero Jacobian entries.
\end{quote}
\end{quote}

\subsubsection{NlpCallbackBase::EvalHess(INdArray\textless{}double, 1\textgreater{}* xdata, double sigma, INdArray\textless{}double, 1\textgreater{}* lambdata, INdArray\textless{}double, 1\textgreater{}* outdata)}
\label{\detokenize{cppapi/nlpcallbackbase:nlpcallbackbase-evalhess-indarray-double-1-xdata-double-sigma-indarray-double-1-lambdata-indarray-double-1-outdata}}\begin{quote}

\sphinxAtStartPar
Evaluate Hessian values of the Lagrangian function of the nonlinear model.

\sphinxAtStartPar
\sphinxstylestrong{Synopsis}
\begin{quote}

\sphinxAtStartPar
\sphinxcode{\sphinxupquote{void EvalHess(}}
\begin{quote}

\sphinxAtStartPar
\sphinxcode{\sphinxupquote{INdArray\textless{}double, 1\textgreater{} *xdata,}}

\sphinxAtStartPar
\sphinxcode{\sphinxupquote{double sigma,}}

\sphinxAtStartPar
\sphinxcode{\sphinxupquote{INdArray\textless{}double, 1\textgreater{} *lambdata,}}

\sphinxAtStartPar
\sphinxcode{\sphinxupquote{INdArray\textless{}double, 1\textgreater{} *outdata)}}
\end{quote}
\end{quote}

\sphinxAtStartPar
\sphinxstylestrong{Arguments}
\begin{quote}

\sphinxAtStartPar
\sphinxcode{\sphinxupquote{xdata}}: pointer to NdArray of variable values.

\sphinxAtStartPar
\sphinxcode{\sphinxupquote{sigma}}: weight on the objective in the Lagrangian function.

\sphinxAtStartPar
\sphinxcode{\sphinxupquote{lambdata}}: pointer to NdArray of Lagrange multiplier associated with the constraint functions.

\sphinxAtStartPar
\sphinxcode{\sphinxupquote{outdata}}: pointer to output NdArray for nonzero Hessian entries.
\end{quote}
\end{quote}

\subsection{ProbBuffer}
\label{\detokenize{cppapiref:probbuffer}}\label{\detokenize{cppapiref:chapcppapiref-probbuffer}}
\sphinxAtStartPar
Buffer object for COPT problem. ProbBuffer object holds the (MPS) problem
in string format.

\sphinxstepscope

\subsubsection{ProbBuffer::ProbBuffer()}
\label{\detokenize{cppapi/probbuffer:probbuffer-probbuffer}}\label{\detokenize{cppapi/probbuffer::doc}}\begin{quote}

\sphinxAtStartPar
Constructor of ProbBuffer object.

\sphinxAtStartPar
\sphinxstylestrong{Synopsis}
\begin{quote}

\sphinxAtStartPar
\sphinxcode{\sphinxupquote{ProbBuffer(int sz)}}
\end{quote}

\sphinxAtStartPar
\sphinxstylestrong{Arguments}
\begin{quote}

\sphinxAtStartPar
\sphinxcode{\sphinxupquote{sz}}: initial size of the problem buffer.
\end{quote}
\end{quote}

\subsubsection{ProbBuffer::GetData()}
\label{\detokenize{cppapi/probbuffer:probbuffer-getdata}}\begin{quote}

\sphinxAtStartPar
Get string of problem in problem buffer.

\sphinxAtStartPar
\sphinxstylestrong{Synopsis}
\begin{quote}

\sphinxAtStartPar
\sphinxcode{\sphinxupquote{char *GetData()}}
\end{quote}

\sphinxAtStartPar
\sphinxstylestrong{Return}
\begin{quote}

\sphinxAtStartPar
string of problem in problem buffer.
\end{quote}
\end{quote}

\subsubsection{ProbBuffer::Resize()}
\label{\detokenize{cppapi/probbuffer:probbuffer-resize}}\begin{quote}

\sphinxAtStartPar
Resize buffer to given size, and zero\sphinxhyphen{}ended

\sphinxAtStartPar
\sphinxstylestrong{Synopsis}
\begin{quote}

\sphinxAtStartPar
\sphinxcode{\sphinxupquote{void Resize(int sz)}}
\end{quote}

\sphinxAtStartPar
\sphinxstylestrong{Arguments}
\begin{quote}

\sphinxAtStartPar
\sphinxcode{\sphinxupquote{sz}}: new buffer size.
\end{quote}
\end{quote}

\subsubsection{ProbBuffer::Size()}
\label{\detokenize{cppapi/probbuffer:probbuffer-size}}\begin{quote}

\sphinxAtStartPar
Get the size of problem buffer.

\sphinxAtStartPar
\sphinxstylestrong{Synopsis}
\begin{quote}

\sphinxAtStartPar
\sphinxcode{\sphinxupquote{int Size()}}
\end{quote}

\sphinxAtStartPar
\sphinxstylestrong{Return}
\begin{quote}

\sphinxAtStartPar
size of problem buffer.
\end{quote}
\end{quote}

\sphinxstepscope

\chapter{C\# API Reference}
\label{\detokenize{csharpapiref:c-api-reference}}\label{\detokenize{csharpapiref:chapcsharpapiref}}\label{\detokenize{csharpapiref::doc}}
\sphinxAtStartPar
The \sphinxstylestrong{Cardinal Optimizer} provides C\# API library.
This chapter documents all COPT C\# constants and API functions for C\# applications.

\section{Constants}
\label{\detokenize{csharpapiref:constants}}\label{\detokenize{csharpapiref:chapcsharpapiref-const}}
\sphinxAtStartPar
There are four types of constants defined in \sphinxstylestrong{Cardinal Optimizer}.
They are general constants, information constants, attributes and parameters.

\subsection{General Constants}
\label{\detokenize{csharpapiref:general-constants}}\label{\detokenize{csharpapiref:chapcsharpapiref-general}}
\sphinxAtStartPar
For the contents of C\# general constants, see
{\hyperref[\detokenize{constant:chapconst}]{\sphinxcrossref{\DUrole{std,std-ref}{General Constants}}}}.

\sphinxAtStartPar
General constants are defined in \sphinxcode{\sphinxupquote{Consts}} class. User may refer general
constants with namespace, that is, \sphinxcode{\sphinxupquote{Copt.Consts.XXXX}}.

\subsection{Attributes}
\label{\detokenize{csharpapiref:attributes}}\label{\detokenize{csharpapiref:chapcsharpapiref-attrs}}
\sphinxAtStartPar
For the contents of C\# attribute constants, see {\hyperref[\detokenize{attribute:chapattrs}]{\sphinxcrossref{\DUrole{std,std-ref}{Attributes}}}}.

\sphinxAtStartPar
All COPT C\# attributes are defined in \sphinxcode{\sphinxupquote{DblAttr}} and \sphinxcode{\sphinxupquote{IntAttr}} classes.
User may refer double attributes by \sphinxcode{\sphinxupquote{Copt.DblAttr.XXXX}}, and integer
attributes by \sphinxcode{\sphinxupquote{Copt.IntAttr.XXXX}}.

\sphinxAtStartPar
In the C\# API, user can get the attribute value by specifying the attribute
name.  The two functions of obtaining attribute values are as follows, please
refer to {\hyperref[\detokenize{csharpapiref:chapcsharpapiref-model}]{\sphinxcrossref{\DUrole{std,std-ref}{C\# API: Model Class}}}} for details.
\begin{itemize}
\item {} 
\sphinxAtStartPar
\sphinxcode{\sphinxupquote{Model.GetIntAttr()}}: Get value of a COPT integer attribute.

\item {} 
\sphinxAtStartPar
\sphinxcode{\sphinxupquote{Model.GetDblAttr()}}: Get value of a COPT double attribute.

\end{itemize}

\subsection{Information}
\label{\detokenize{csharpapiref:information}}\label{\detokenize{csharpapiref:chapcsharpapiref-info}}
\sphinxAtStartPar
For the content of C\# information constants, see {\hyperref[\detokenize{information:chapinfo}]{\sphinxcrossref{\DUrole{std,std-ref}{Information}}}}.

\sphinxAtStartPar
In the C\# API, information constants are defined in the \sphinxcode{\sphinxupquote{DblInfo}} class.
Users can access information constants through the prefix \sphinxcode{\sphinxupquote{Copt}} in the namespace (usually can be omitted) \sphinxcode{\sphinxupquote{Copt.DblInfo.}}

\sphinxAtStartPar
For instance, \sphinxcode{\sphinxupquote{Copt.DblInfo.Obj}} is the coefficients of variables in the objective function.

\subsection{Callback Information}
\label{\detokenize{csharpapiref:callback-information}}
\sphinxAtStartPar
For the content of C\# API callback information class constants, see {\hyperref[\detokenize{information:chapinfo-cbc}]{\sphinxcrossref{\DUrole{std,std-ref}{Callback Information}}}}.

\sphinxAtStartPar
In the C\# API, callback\sphinxhyphen{}related information constants are defined in the \sphinxcode{\sphinxupquote{CbInfo}} class.
Users can access information constants through the prefix \sphinxcode{\sphinxupquote{Copt}} in the namespace (usually can be omitted) \sphinxcode{\sphinxupquote{Copt.CbInfo.}}

\sphinxAtStartPar
For instance, \sphinxcode{\sphinxupquote{Copt.CbInfo.BestObj}} is the current best objective.

\subsection{Parameters}
\label{\detokenize{csharpapiref:parameters}}\label{\detokenize{csharpapiref:chapcsharpapiref-params}}
\sphinxAtStartPar
For the contents of C\# parameters constants, see
{\hyperref[\detokenize{parameter:chapparams}]{\sphinxcrossref{\DUrole{std,std-ref}{Parameters}}}}.

\sphinxAtStartPar
All COPT C\# parameters are defined in DblParam and IntParam classes.
User may refer double parameters by \sphinxcode{\sphinxupquote{Copt.DblParam.XXXX}}, and integer
parameters by \sphinxcode{\sphinxupquote{Copt.IntParam.XXXX}}.

\sphinxAtStartPar
In the C\# API, user can get and set the parameter value by specifying the
parameter name. The provided functions are as follows, please refer to
{\hyperref[\detokenize{csharpapiref:chapcsharpapiref-model}]{\sphinxcrossref{\DUrole{std,std-ref}{C\# API: Model Class}}}} for details.
\begin{itemize}
\item {} 
\sphinxAtStartPar
Get detailed information of the specified parameter (current value/max/min): \sphinxcode{\sphinxupquote{Model.GetParamInfo()}}

\item {} 
\sphinxAtStartPar
Get the current value of the specified integer/double parameter: \sphinxcode{\sphinxupquote{Model.GetIntParam()}} / \sphinxcode{\sphinxupquote{Model.GetDblParam()}}

\item {} 
\sphinxAtStartPar
Set the specified integer/double parameter value: \sphinxcode{\sphinxupquote{Model.SetIntParam()}} / \sphinxcode{\sphinxupquote{Model.SetDblParam()}}

\end{itemize}

\section{C\# Modeling Classes}
\label{\detokenize{csharpapiref:c-modeling-classes}}
\sphinxAtStartPar
This chapter documents COPT C\# interface. Users may refer to
C\# classes described below for details of how to construct and
solve C\# models.

\subsection{Envr}
\label{\detokenize{csharpapiref:envr}}\label{\detokenize{csharpapiref:chapcsharpapiref-envr}}
\sphinxAtStartPar
Essentially, any C\# application using Cardinal Optimizer should start with a
COPT environment. COPT models are always associated with a COPT environment.
User must create an environment object before populating models.
User generally only need a single environment object in program.

\sphinxstepscope

\subsubsection{Envr.Envr()}
\label{\detokenize{csapi/envr:envr-envr}}\label{\detokenize{csapi/envr::doc}}\begin{quote}

\sphinxAtStartPar
Constructor of COPT Envr object.

\sphinxAtStartPar
\sphinxstylestrong{Synopsis}
\begin{quote}

\sphinxAtStartPar
\sphinxcode{\sphinxupquote{Envr()}}
\end{quote}
\end{quote}

\subsubsection{Envr.Envr()}
\label{\detokenize{csapi/envr:id1}}\begin{quote}

\sphinxAtStartPar
Constructor of COPT Envr object, given a license folder.

\sphinxAtStartPar
\sphinxstylestrong{Synopsis}
\begin{quote}

\sphinxAtStartPar
\sphinxcode{\sphinxupquote{Envr(string licDir)}}
\end{quote}

\sphinxAtStartPar
\sphinxstylestrong{Arguments}
\begin{quote}

\sphinxAtStartPar
\sphinxcode{\sphinxupquote{licDir}}: directory having local license or client config file.
\end{quote}
\end{quote}

\subsubsection{Envr.Envr()}
\label{\detokenize{csapi/envr:id2}}\begin{quote}

\sphinxAtStartPar
Constructor of COPT Envr object, given an Envr config object.

\sphinxAtStartPar
\sphinxstylestrong{Synopsis}
\begin{quote}

\sphinxAtStartPar
\sphinxcode{\sphinxupquote{Envr(EnvrConfig config)}}
\end{quote}

\sphinxAtStartPar
\sphinxstylestrong{Arguments}
\begin{quote}

\sphinxAtStartPar
\sphinxcode{\sphinxupquote{config}}: Envr config object holding settings for remote connection.
\end{quote}
\end{quote}

\subsubsection{Envr.BindNumaCpu()}
\label{\detokenize{csapi/envr:envr-bindnumacpu}}\begin{quote}

\sphinxAtStartPar
Bind the CPUs for the current process to a NUMA node.

\sphinxAtStartPar
\sphinxstylestrong{Synopsis}
\begin{quote}

\sphinxAtStartPar
\sphinxcode{\sphinxupquote{void BindNumaCpu(int numaNode)}}
\end{quote}

\sphinxAtStartPar
\sphinxstylestrong{Arguments}
\begin{quote}

\sphinxAtStartPar
\sphinxcode{\sphinxupquote{numaNode}}: ID of a NUMA node.
\end{quote}
\end{quote}

\subsubsection{Envr.BindNumaMem()}
\label{\detokenize{csapi/envr:envr-bindnumamem}}\begin{quote}

\sphinxAtStartPar
Bind memory for the current process to a NUMA node (Linux only).

\sphinxAtStartPar
\sphinxstylestrong{Synopsis}
\begin{quote}

\sphinxAtStartPar
\sphinxcode{\sphinxupquote{void BindNumaMem(int numaNode)}}
\end{quote}

\sphinxAtStartPar
\sphinxstylestrong{Arguments}
\begin{quote}

\sphinxAtStartPar
\sphinxcode{\sphinxupquote{numaNode}}: the ID of a NUMA node.
\end{quote}
\end{quote}

\subsubsection{Envr.Close()}
\label{\detokenize{csapi/envr:envr-close}}\begin{quote}

\sphinxAtStartPar
close remote connection and token becomes invalid for all problems in current envr.

\sphinxAtStartPar
\sphinxstylestrong{Synopsis}
\begin{quote}

\sphinxAtStartPar
\sphinxcode{\sphinxupquote{void Close()}}
\end{quote}
\end{quote}

\subsubsection{Envr.CreateModel()}
\label{\detokenize{csapi/envr:envr-createmodel}}\begin{quote}

\sphinxAtStartPar
Create a model object.

\sphinxAtStartPar
\sphinxstylestrong{Synopsis}
\begin{quote}

\sphinxAtStartPar
\sphinxcode{\sphinxupquote{Model CreateModel(string name)}}
\end{quote}

\sphinxAtStartPar
\sphinxstylestrong{Arguments}
\begin{quote}

\sphinxAtStartPar
\sphinxcode{\sphinxupquote{name}}: customized model name.
\end{quote}

\sphinxAtStartPar
\sphinxstylestrong{Return}
\begin{quote}

\sphinxAtStartPar
a model object.
\end{quote}
\end{quote}

\subsubsection{Envr.GetCpuAffinity()}
\label{\detokenize{csapi/envr:envr-getcpuaffinity}}\begin{quote}

\sphinxAtStartPar
Get CPU affinity for the current process, which is saved in an integer array.

\sphinxAtStartPar
\sphinxstylestrong{Synopsis}
\begin{quote}

\sphinxAtStartPar
\sphinxcode{\sphinxupquote{int{[}{]} GetCpuAffinity()}}
\end{quote}

\sphinxAtStartPar
\sphinxstylestrong{Return}
\begin{quote}

\sphinxAtStartPar
an integer array of CPU IDs.
\end{quote}
\end{quote}

\subsubsection{Envr.GetNumaNodeCount()}
\label{\detokenize{csapi/envr:envr-getnumanodecount}}\begin{quote}

\sphinxAtStartPar
Get count of NUMA nodes.

\sphinxAtStartPar
\sphinxstylestrong{Synopsis}
\begin{quote}

\sphinxAtStartPar
\sphinxcode{\sphinxupquote{int GetNumaNodeCount()}}
\end{quote}

\sphinxAtStartPar
\sphinxstylestrong{Return}
\begin{quote}

\sphinxAtStartPar
count of NUMA nodes.
\end{quote}
\end{quote}

\subsubsection{Envr.SetCpuAffinity()}
\label{\detokenize{csapi/envr:envr-setcpuaffinity}}\begin{quote}

\sphinxAtStartPar
Set CPU affinity with given mask string.

\sphinxAtStartPar
\sphinxstylestrong{Synopsis}
\begin{quote}

\sphinxAtStartPar
\sphinxcode{\sphinxupquote{void SetCpuAffinity(string hexMask)}}
\end{quote}

\sphinxAtStartPar
\sphinxstylestrong{Arguments}
\begin{quote}

\sphinxAtStartPar
\sphinxcode{\sphinxupquote{hexMask}}: CPU mask string of hexadecimal characters.
\end{quote}
\end{quote}

\subsection{EnvrConfig}
\label{\detokenize{csharpapiref:envrconfig}}
\sphinxAtStartPar
If user connects to COPT remote services, such as floating token server or compute cluster,
it is necessary to add config settings with EnvrConfig object.

\sphinxstepscope

\subsubsection{EnvrConfig.EnvrConfig()}
\label{\detokenize{csapi/envrconfig:envrconfig-envrconfig}}\label{\detokenize{csapi/envrconfig::doc}}\begin{quote}

\sphinxAtStartPar
Constructor of envr config object.

\sphinxAtStartPar
\sphinxstylestrong{Synopsis}
\begin{quote}

\sphinxAtStartPar
\sphinxcode{\sphinxupquote{EnvrConfig()}}
\end{quote}
\end{quote}

\subsubsection{EnvrConfig.Set()}
\label{\detokenize{csapi/envrconfig:envrconfig-set}}\begin{quote}

\sphinxAtStartPar
Set config settings in terms of name\sphinxhyphen{}value pair.

\sphinxAtStartPar
\sphinxstylestrong{Synopsis}
\begin{quote}

\sphinxAtStartPar
\sphinxcode{\sphinxupquote{void Set(string name, string value)}}
\end{quote}

\sphinxAtStartPar
\sphinxstylestrong{Arguments}
\begin{quote}

\sphinxAtStartPar
\sphinxcode{\sphinxupquote{name}}: keyword of a config setting.

\sphinxAtStartPar
\sphinxcode{\sphinxupquote{value}}: value of a config setting.
\end{quote}
\end{quote}

\subsection{Model}
\label{\detokenize{csharpapiref:model}}\label{\detokenize{csharpapiref:chapcsharpapiref-model}}
\sphinxAtStartPar
In general, a COPT model consists of a set of variables, a (linear) objective
function on these variables, a set of constraints on there varaibles, etc.
COPT model class encapsulates all required methods for constructing a COPT model.

\sphinxstepscope

\subsubsection{Model.Model()}
\label{\detokenize{csapi/model:model-model}}\label{\detokenize{csapi/model::doc}}\begin{quote}

\sphinxAtStartPar
Constructor of model.

\sphinxAtStartPar
\sphinxstylestrong{Synopsis}
\begin{quote}

\sphinxAtStartPar
\sphinxcode{\sphinxupquote{Model(Envr env, string name)}}
\end{quote}

\sphinxAtStartPar
\sphinxstylestrong{Arguments}
\begin{quote}

\sphinxAtStartPar
\sphinxcode{\sphinxupquote{env}}: associated environment object.

\sphinxAtStartPar
\sphinxcode{\sphinxupquote{name}}: string of model name.
\end{quote}
\end{quote}

\subsubsection{Model.AddAffineCone()}
\label{\detokenize{csapi/model:model-addaffinecone}}\begin{quote}

\sphinxAtStartPar
Add an affine cone constraint to model.

\sphinxAtStartPar
\sphinxstylestrong{Synopsis}
\begin{quote}

\sphinxAtStartPar
\sphinxcode{\sphinxupquote{AffineCone AddAffineCone(AffineConeBuilder builder, string name)}}
\end{quote}

\sphinxAtStartPar
\sphinxstylestrong{Arguments}
\begin{quote}

\sphinxAtStartPar
\sphinxcode{\sphinxupquote{builder}}: builder for new affine cone constraint.

\sphinxAtStartPar
\sphinxcode{\sphinxupquote{name}}: optional, name of new affine cone constraint.
\end{quote}

\sphinxAtStartPar
\sphinxstylestrong{Return}
\begin{quote}

\sphinxAtStartPar
new affine cone constraint object.
\end{quote}
\end{quote}

\subsubsection{Model.AddAffineCone()}
\label{\detokenize{csapi/model:id1}}\begin{quote}

\sphinxAtStartPar
Add an affine cone constraint to a model.

\sphinxAtStartPar
\sphinxstylestrong{Synopsis}
\begin{quote}

\sphinxAtStartPar
\sphinxcode{\sphinxupquote{AffineCone AddAffineCone(}}
\begin{quote}

\sphinxAtStartPar
\sphinxcode{\sphinxupquote{MLinExpr expr,}}

\sphinxAtStartPar
\sphinxcode{\sphinxupquote{int type,}}

\sphinxAtStartPar
\sphinxcode{\sphinxupquote{string name)}}
\end{quote}
\end{quote}

\sphinxAtStartPar
\sphinxstylestrong{Arguments}
\begin{quote}

\sphinxAtStartPar
\sphinxcode{\sphinxupquote{expr}}: 1\sphinxhyphen{}dimensional array of linear expressions.

\sphinxAtStartPar
\sphinxcode{\sphinxupquote{type}}: type of an affine cone.

\sphinxAtStartPar
\sphinxcode{\sphinxupquote{name}}: name of new affine cone constraint.
\end{quote}

\sphinxAtStartPar
\sphinxstylestrong{Return}
\begin{quote}

\sphinxAtStartPar
new affine cone constraint object.
\end{quote}
\end{quote}

\subsubsection{Model.AddAffineCone()}
\label{\detokenize{csapi/model:id2}}\begin{quote}

\sphinxAtStartPar
Add an affine cone constraint to a model.

\sphinxAtStartPar
\sphinxstylestrong{Synopsis}
\begin{quote}

\sphinxAtStartPar
\sphinxcode{\sphinxupquote{AffineCone AddAffineCone(}}
\begin{quote}

\sphinxAtStartPar
\sphinxcode{\sphinxupquote{MVar vars,}}

\sphinxAtStartPar
\sphinxcode{\sphinxupquote{int type,}}

\sphinxAtStartPar
\sphinxcode{\sphinxupquote{string name)}}
\end{quote}
\end{quote}

\sphinxAtStartPar
\sphinxstylestrong{Arguments}
\begin{quote}

\sphinxAtStartPar
\sphinxcode{\sphinxupquote{vars}}: 1\sphinxhyphen{}dimensional array of variables.

\sphinxAtStartPar
\sphinxcode{\sphinxupquote{type}}: type of an affine cone.

\sphinxAtStartPar
\sphinxcode{\sphinxupquote{name}}: name of new affine cone constraint.
\end{quote}

\sphinxAtStartPar
\sphinxstylestrong{Return}
\begin{quote}

\sphinxAtStartPar
new affine cone constraint object.
\end{quote}
\end{quote}

\subsubsection{Model.AddAffineCone()}
\label{\detokenize{csapi/model:id3}}\begin{quote}

\sphinxAtStartPar
Add an affine cone constraint to model.

\sphinxAtStartPar
\sphinxstylestrong{Synopsis}
\begin{quote}

\sphinxAtStartPar
\sphinxcode{\sphinxupquote{AffineCone AddAffineCone(}}
\begin{quote}

\sphinxAtStartPar
\sphinxcode{\sphinxupquote{Expr{[}{]} exprs,}}

\sphinxAtStartPar
\sphinxcode{\sphinxupquote{int type,}}

\sphinxAtStartPar
\sphinxcode{\sphinxupquote{string name)}}
\end{quote}
\end{quote}

\sphinxAtStartPar
\sphinxstylestrong{Arguments}
\begin{quote}

\sphinxAtStartPar
\sphinxcode{\sphinxupquote{exprs}}: linear expressions that participate in the affine cone constraint.

\sphinxAtStartPar
\sphinxcode{\sphinxupquote{type}}: type of an affine cone constraint.

\sphinxAtStartPar
\sphinxcode{\sphinxupquote{name}}: name of new affine cone constraint.
\end{quote}

\sphinxAtStartPar
\sphinxstylestrong{Return}
\begin{quote}

\sphinxAtStartPar
new affine cone constraint object.
\end{quote}
\end{quote}

\subsubsection{Model.AddAffineCone()}
\label{\detokenize{csapi/model:id4}}\begin{quote}

\sphinxAtStartPar
Add an affine cone constraint to a model.

\sphinxAtStartPar
\sphinxstylestrong{Synopsis}
\begin{quote}

\sphinxAtStartPar
\sphinxcode{\sphinxupquote{AffineCone AddAffineCone(}}
\begin{quote}

\sphinxAtStartPar
\sphinxcode{\sphinxupquote{Var{[}{]} vars,}}

\sphinxAtStartPar
\sphinxcode{\sphinxupquote{int type,}}

\sphinxAtStartPar
\sphinxcode{\sphinxupquote{string name)}}
\end{quote}
\end{quote}

\sphinxAtStartPar
\sphinxstylestrong{Arguments}
\begin{quote}

\sphinxAtStartPar
\sphinxcode{\sphinxupquote{vars}}: variables that participate in the affine cone constraint.

\sphinxAtStartPar
\sphinxcode{\sphinxupquote{type}}: type of an affine cone.

\sphinxAtStartPar
\sphinxcode{\sphinxupquote{name}}: name of new affine cone constraint.
\end{quote}

\sphinxAtStartPar
\sphinxstylestrong{Return}
\begin{quote}

\sphinxAtStartPar
new affine cone constraint object.
\end{quote}
\end{quote}

\subsubsection{Model.AddAffineCone()}
\label{\detokenize{csapi/model:id5}}\begin{quote}

\sphinxAtStartPar
Add an affine cone constraint to a model.

\sphinxAtStartPar
\sphinxstylestrong{Synopsis}
\begin{quote}

\sphinxAtStartPar
\sphinxcode{\sphinxupquote{AffineCone AddAffineCone(}}
\begin{quote}

\sphinxAtStartPar
\sphinxcode{\sphinxupquote{MPsdExpr expr,}}

\sphinxAtStartPar
\sphinxcode{\sphinxupquote{int type,}}

\sphinxAtStartPar
\sphinxcode{\sphinxupquote{string name)}}
\end{quote}
\end{quote}

\sphinxAtStartPar
\sphinxstylestrong{Arguments}
\begin{quote}

\sphinxAtStartPar
\sphinxcode{\sphinxupquote{expr}}: 1\sphinxhyphen{}dimensional array of PSD expressions.

\sphinxAtStartPar
\sphinxcode{\sphinxupquote{type}}: type of an affine cone.

\sphinxAtStartPar
\sphinxcode{\sphinxupquote{name}}: name of new affine cone constraint.
\end{quote}

\sphinxAtStartPar
\sphinxstylestrong{Return}
\begin{quote}

\sphinxAtStartPar
new affine cone constraint object.
\end{quote}
\end{quote}

\subsubsection{Model.AddAffineCone()}
\label{\detokenize{csapi/model:id6}}\begin{quote}

\sphinxAtStartPar
Add an affine cone constraint to model.

\sphinxAtStartPar
\sphinxstylestrong{Synopsis}
\begin{quote}

\sphinxAtStartPar
\sphinxcode{\sphinxupquote{AffineCone AddAffineCone(}}
\begin{quote}

\sphinxAtStartPar
\sphinxcode{\sphinxupquote{PsdExpr{[}{]} exprs,}}

\sphinxAtStartPar
\sphinxcode{\sphinxupquote{int type,}}

\sphinxAtStartPar
\sphinxcode{\sphinxupquote{string name)}}
\end{quote}
\end{quote}

\sphinxAtStartPar
\sphinxstylestrong{Arguments}
\begin{quote}

\sphinxAtStartPar
\sphinxcode{\sphinxupquote{exprs}}: PSD expressions that participate in the affine cone constraint.

\sphinxAtStartPar
\sphinxcode{\sphinxupquote{type}}: type of an affine cone constraint.

\sphinxAtStartPar
\sphinxcode{\sphinxupquote{name}}: name of new affine cone constraint.
\end{quote}

\sphinxAtStartPar
\sphinxstylestrong{Return}
\begin{quote}

\sphinxAtStartPar
new affine cone constraint object.
\end{quote}
\end{quote}

\subsubsection{Model.AddAffineCones()}
\label{\detokenize{csapi/model:model-addaffinecones}}\begin{quote}

\sphinxAtStartPar
Add a list of affine cone constraints to a model.

\sphinxAtStartPar
\sphinxstylestrong{Synopsis}
\begin{quote}

\sphinxAtStartPar
\sphinxcode{\sphinxupquote{AffineConeArray AddAffineCones(}}
\begin{quote}

\sphinxAtStartPar
\sphinxcode{\sphinxupquote{MLinExpr exprs,}}

\sphinxAtStartPar
\sphinxcode{\sphinxupquote{int type,}}

\sphinxAtStartPar
\sphinxcode{\sphinxupquote{string prefix)}}
\end{quote}
\end{quote}

\sphinxAtStartPar
\sphinxstylestrong{Arguments}
\begin{quote}

\sphinxAtStartPar
\sphinxcode{\sphinxupquote{exprs}}: 2\sphinxhyphen{}dimensional array of linear expressions.

\sphinxAtStartPar
\sphinxcode{\sphinxupquote{type}}: type of an affine cone.

\sphinxAtStartPar
\sphinxcode{\sphinxupquote{prefix}}: name prefix for new affine cone constraints.
\end{quote}

\sphinxAtStartPar
\sphinxstylestrong{Return}
\begin{quote}

\sphinxAtStartPar
array of new affine cone constraint objects.
\end{quote}
\end{quote}

\subsubsection{Model.AddAffineCones()}
\label{\detokenize{csapi/model:id7}}\begin{quote}

\sphinxAtStartPar
Add a list of affine cone constraints to a model.

\sphinxAtStartPar
\sphinxstylestrong{Synopsis}
\begin{quote}

\sphinxAtStartPar
\sphinxcode{\sphinxupquote{AffineConeArray AddAffineCones(}}
\begin{quote}

\sphinxAtStartPar
\sphinxcode{\sphinxupquote{MVar vars,}}

\sphinxAtStartPar
\sphinxcode{\sphinxupquote{int type,}}

\sphinxAtStartPar
\sphinxcode{\sphinxupquote{string prefix)}}
\end{quote}
\end{quote}

\sphinxAtStartPar
\sphinxstylestrong{Arguments}
\begin{quote}

\sphinxAtStartPar
\sphinxcode{\sphinxupquote{vars}}: 2\sphinxhyphen{}dimensional array of varissblespp

\sphinxAtStartPar
\sphinxcode{\sphinxupquote{type}}: type of affine cones.

\sphinxAtStartPar
\sphinxcode{\sphinxupquote{prefix}}: name prefix for new affine cone constraints.
\end{quote}

\sphinxAtStartPar
\sphinxstylestrong{Return}
\begin{quote}

\sphinxAtStartPar
array of new affine cone constraint objects.
\end{quote}
\end{quote}

\subsubsection{Model.AddAffineCones()}
\label{\detokenize{csapi/model:id8}}\begin{quote}

\sphinxAtStartPar
Add a list of affine cone constraints to a model.

\sphinxAtStartPar
\sphinxstylestrong{Synopsis}
\begin{quote}

\sphinxAtStartPar
\sphinxcode{\sphinxupquote{AffineConeArray AddAffineCones(}}
\begin{quote}

\sphinxAtStartPar
\sphinxcode{\sphinxupquote{MPsdExpr exprs,}}

\sphinxAtStartPar
\sphinxcode{\sphinxupquote{int type,}}

\sphinxAtStartPar
\sphinxcode{\sphinxupquote{string prefix)}}
\end{quote}
\end{quote}

\sphinxAtStartPar
\sphinxstylestrong{Arguments}
\begin{quote}

\sphinxAtStartPar
\sphinxcode{\sphinxupquote{exprs}}: 2\sphinxhyphen{}dimensional array of PSD expressions.

\sphinxAtStartPar
\sphinxcode{\sphinxupquote{type}}: type of an affine cone.

\sphinxAtStartPar
\sphinxcode{\sphinxupquote{prefix}}: name prefix for new affine cone constraints.
\end{quote}

\sphinxAtStartPar
\sphinxstylestrong{Return}
\begin{quote}

\sphinxAtStartPar
array of new affine cone constraint objects.
\end{quote}
\end{quote}

\subsubsection{Model.AddCone()}
\label{\detokenize{csapi/model:model-addcone}}\begin{quote}

\sphinxAtStartPar
Add a cone constraint to model.

\sphinxAtStartPar
\sphinxstylestrong{Synopsis}
\begin{quote}

\sphinxAtStartPar
\sphinxcode{\sphinxupquote{Cone AddCone(}}
\begin{quote}

\sphinxAtStartPar
\sphinxcode{\sphinxupquote{int dim,}}

\sphinxAtStartPar
\sphinxcode{\sphinxupquote{int type,}}

\sphinxAtStartPar
\sphinxcode{\sphinxupquote{char{[}{]} pvtype,}}

\sphinxAtStartPar
\sphinxcode{\sphinxupquote{string prefix)}}
\end{quote}
\end{quote}

\sphinxAtStartPar
\sphinxstylestrong{Arguments}
\begin{quote}

\sphinxAtStartPar
\sphinxcode{\sphinxupquote{dim}}: dimension of the cone constraint.

\sphinxAtStartPar
\sphinxcode{\sphinxupquote{type}}: type of the cone constraint.

\sphinxAtStartPar
\sphinxcode{\sphinxupquote{pvtype}}: type of variables in the cone.

\sphinxAtStartPar
\sphinxcode{\sphinxupquote{prefix}}: optional, name prefix of variables in the cone, default value is “ConeV”.
\end{quote}

\sphinxAtStartPar
\sphinxstylestrong{Return}
\begin{quote}

\sphinxAtStartPar
new cone constraint object.
\end{quote}
\end{quote}

\subsubsection{Model.AddCone()}
\label{\detokenize{csapi/model:id9}}\begin{quote}

\sphinxAtStartPar
Add a cone constraint to model.

\sphinxAtStartPar
\sphinxstylestrong{Synopsis}
\begin{quote}

\sphinxAtStartPar
\sphinxcode{\sphinxupquote{Cone AddCone(ConeBuilder builder)}}
\end{quote}

\sphinxAtStartPar
\sphinxstylestrong{Arguments}
\begin{quote}

\sphinxAtStartPar
\sphinxcode{\sphinxupquote{builder}}: builder for new cone constraint.
\end{quote}

\sphinxAtStartPar
\sphinxstylestrong{Return}
\begin{quote}

\sphinxAtStartPar
new cone constraint object.
\end{quote}
\end{quote}

\subsubsection{Model.AddCone()}
\label{\detokenize{csapi/model:id10}}\begin{quote}

\sphinxAtStartPar
Add a cone constraint to model.

\sphinxAtStartPar
\sphinxstylestrong{Synopsis}
\begin{quote}

\sphinxAtStartPar
\sphinxcode{\sphinxupquote{Cone AddCone(Var{[}{]} vars, int type)}}
\end{quote}

\sphinxAtStartPar
\sphinxstylestrong{Arguments}
\begin{quote}

\sphinxAtStartPar
\sphinxcode{\sphinxupquote{vars}}: variables that participate in the cone constraint.

\sphinxAtStartPar
\sphinxcode{\sphinxupquote{type}}: type of the cone constraint.
\end{quote}

\sphinxAtStartPar
\sphinxstylestrong{Return}
\begin{quote}

\sphinxAtStartPar
new cone constraint object.
\end{quote}
\end{quote}

\subsubsection{Model.AddCone()}
\label{\detokenize{csapi/model:id11}}\begin{quote}

\sphinxAtStartPar
Add a cone constraint to model.

\sphinxAtStartPar
\sphinxstylestrong{Synopsis}
\begin{quote}

\sphinxAtStartPar
\sphinxcode{\sphinxupquote{Cone AddCone(VarArray vars, int type)}}
\end{quote}

\sphinxAtStartPar
\sphinxstylestrong{Arguments}
\begin{quote}

\sphinxAtStartPar
\sphinxcode{\sphinxupquote{vars}}: variables that participate in the cone constraint.

\sphinxAtStartPar
\sphinxcode{\sphinxupquote{type}}: type of a cone constraint.
\end{quote}

\sphinxAtStartPar
\sphinxstylestrong{Return}
\begin{quote}

\sphinxAtStartPar
new cone constraint object.
\end{quote}
\end{quote}

\subsubsection{Model.AddConstr()}
\label{\detokenize{csapi/model:model-addconstr}}\begin{quote}

\sphinxAtStartPar
Add a linear constraint to model.

\sphinxAtStartPar
\sphinxstylestrong{Synopsis}
\begin{quote}

\sphinxAtStartPar
\sphinxcode{\sphinxupquote{Constraint AddConstr(}}
\begin{quote}

\sphinxAtStartPar
\sphinxcode{\sphinxupquote{Expr expr,}}

\sphinxAtStartPar
\sphinxcode{\sphinxupquote{char sense,}}

\sphinxAtStartPar
\sphinxcode{\sphinxupquote{double rhs,}}

\sphinxAtStartPar
\sphinxcode{\sphinxupquote{string name)}}
\end{quote}
\end{quote}

\sphinxAtStartPar
\sphinxstylestrong{Arguments}
\begin{quote}

\sphinxAtStartPar
\sphinxcode{\sphinxupquote{expr}}: expression for the new contraint.

\sphinxAtStartPar
\sphinxcode{\sphinxupquote{sense}}: sense for new linear constraint, other than range sense.

\sphinxAtStartPar
\sphinxcode{\sphinxupquote{rhs}}: right hand side value for the new constraint.

\sphinxAtStartPar
\sphinxcode{\sphinxupquote{name}}: optional, name of new constraint.
\end{quote}

\sphinxAtStartPar
\sphinxstylestrong{Return}
\begin{quote}

\sphinxAtStartPar
new constraint object.
\end{quote}
\end{quote}

\subsubsection{Model.AddConstr()}
\label{\detokenize{csapi/model:id12}}\begin{quote}

\sphinxAtStartPar
Add a linear constraint to model.

\sphinxAtStartPar
\sphinxstylestrong{Synopsis}
\begin{quote}

\sphinxAtStartPar
\sphinxcode{\sphinxupquote{Constraint AddConstr(}}
\begin{quote}

\sphinxAtStartPar
\sphinxcode{\sphinxupquote{Expr expr,}}

\sphinxAtStartPar
\sphinxcode{\sphinxupquote{char sense,}}

\sphinxAtStartPar
\sphinxcode{\sphinxupquote{Var var,}}

\sphinxAtStartPar
\sphinxcode{\sphinxupquote{string name)}}
\end{quote}
\end{quote}

\sphinxAtStartPar
\sphinxstylestrong{Arguments}
\begin{quote}

\sphinxAtStartPar
\sphinxcode{\sphinxupquote{expr}}: expression for the new contraint.

\sphinxAtStartPar
\sphinxcode{\sphinxupquote{sense}}: sense for new linear constraint, other than range sense.

\sphinxAtStartPar
\sphinxcode{\sphinxupquote{var}}: variable as right hand side for the new constraint.

\sphinxAtStartPar
\sphinxcode{\sphinxupquote{name}}: optional, name of new constraint.
\end{quote}

\sphinxAtStartPar
\sphinxstylestrong{Return}
\begin{quote}

\sphinxAtStartPar
new constraint object.
\end{quote}
\end{quote}

\subsubsection{Model.AddConstr()}
\label{\detokenize{csapi/model:id13}}\begin{quote}

\sphinxAtStartPar
Add a linear constraint to model.

\sphinxAtStartPar
\sphinxstylestrong{Synopsis}
\begin{quote}

\sphinxAtStartPar
\sphinxcode{\sphinxupquote{Constraint AddConstr(}}
\begin{quote}

\sphinxAtStartPar
\sphinxcode{\sphinxupquote{Expr lhs,}}

\sphinxAtStartPar
\sphinxcode{\sphinxupquote{char sense,}}

\sphinxAtStartPar
\sphinxcode{\sphinxupquote{Expr rhs,}}

\sphinxAtStartPar
\sphinxcode{\sphinxupquote{string name)}}
\end{quote}
\end{quote}

\sphinxAtStartPar
\sphinxstylestrong{Arguments}
\begin{quote}

\sphinxAtStartPar
\sphinxcode{\sphinxupquote{lhs}}: left hand side expression for the new constraint.

\sphinxAtStartPar
\sphinxcode{\sphinxupquote{sense}}: sense for new linear constraint, other than range sense.

\sphinxAtStartPar
\sphinxcode{\sphinxupquote{rhs}}: right hand side expression for the new constraint.

\sphinxAtStartPar
\sphinxcode{\sphinxupquote{name}}: optional, name of new constraint.
\end{quote}

\sphinxAtStartPar
\sphinxstylestrong{Return}
\begin{quote}

\sphinxAtStartPar
new constraint object.
\end{quote}
\end{quote}

\subsubsection{Model.AddConstr()}
\label{\detokenize{csapi/model:id14}}\begin{quote}

\sphinxAtStartPar
Add a linear constraint to model.

\sphinxAtStartPar
\sphinxstylestrong{Synopsis}
\begin{quote}

\sphinxAtStartPar
\sphinxcode{\sphinxupquote{Constraint AddConstr(}}
\begin{quote}

\sphinxAtStartPar
\sphinxcode{\sphinxupquote{Expr expr,}}

\sphinxAtStartPar
\sphinxcode{\sphinxupquote{double lb,}}

\sphinxAtStartPar
\sphinxcode{\sphinxupquote{double ub,}}

\sphinxAtStartPar
\sphinxcode{\sphinxupquote{string name)}}
\end{quote}
\end{quote}

\sphinxAtStartPar
\sphinxstylestrong{Arguments}
\begin{quote}

\sphinxAtStartPar
\sphinxcode{\sphinxupquote{expr}}: expression for the new constraint.

\sphinxAtStartPar
\sphinxcode{\sphinxupquote{lb}}: lower bound for the new constraint.

\sphinxAtStartPar
\sphinxcode{\sphinxupquote{ub}}: upper bound for the new constraint

\sphinxAtStartPar
\sphinxcode{\sphinxupquote{name}}: optional, name of new constraint.
\end{quote}

\sphinxAtStartPar
\sphinxstylestrong{Return}
\begin{quote}

\sphinxAtStartPar
new constraint object.
\end{quote}
\end{quote}

\subsubsection{Model.AddConstr()}
\label{\detokenize{csapi/model:id15}}\begin{quote}

\sphinxAtStartPar
Add a linear constraint to a model.

\sphinxAtStartPar
\sphinxstylestrong{Synopsis}
\begin{quote}

\sphinxAtStartPar
\sphinxcode{\sphinxupquote{Constraint AddConstr(ConstrBuilder builder, string name)}}
\end{quote}

\sphinxAtStartPar
\sphinxstylestrong{Arguments}
\begin{quote}

\sphinxAtStartPar
\sphinxcode{\sphinxupquote{builder}}: builder for the new constraint.

\sphinxAtStartPar
\sphinxcode{\sphinxupquote{name}}: optional, name of new constraint.
\end{quote}

\sphinxAtStartPar
\sphinxstylestrong{Return}
\begin{quote}

\sphinxAtStartPar
new constraint object.
\end{quote}
\end{quote}

\subsubsection{Model.AddConstrs()}
\label{\detokenize{csapi/model:model-addconstrs}}\begin{quote}

\sphinxAtStartPar
Add linear constraints to model.

\sphinxAtStartPar
\sphinxstylestrong{Synopsis}
\begin{quote}

\sphinxAtStartPar
\sphinxcode{\sphinxupquote{ConstrArray AddConstrs(}}
\begin{quote}

\sphinxAtStartPar
\sphinxcode{\sphinxupquote{int count,}}

\sphinxAtStartPar
\sphinxcode{\sphinxupquote{char{[}{]} senses,}}

\sphinxAtStartPar
\sphinxcode{\sphinxupquote{double{[}{]} rhss,}}

\sphinxAtStartPar
\sphinxcode{\sphinxupquote{string prefix)}}
\end{quote}
\end{quote}

\sphinxAtStartPar
\sphinxstylestrong{Arguments}
\begin{quote}

\sphinxAtStartPar
\sphinxcode{\sphinxupquote{count}}: number of constraints added to model.

\sphinxAtStartPar
\sphinxcode{\sphinxupquote{senses}}: sense array for new linear constraints, other than range sense.

\sphinxAtStartPar
\sphinxcode{\sphinxupquote{rhss}}: right hand side values for new variables.

\sphinxAtStartPar
\sphinxcode{\sphinxupquote{prefix}}: optional, name prefix for new constraints, default value is ‘R’.
\end{quote}

\sphinxAtStartPar
\sphinxstylestrong{Return}
\begin{quote}

\sphinxAtStartPar
array of new constraint objects.
\end{quote}
\end{quote}

\subsubsection{Model.AddConstrs()}
\label{\detokenize{csapi/model:id16}}\begin{quote}

\sphinxAtStartPar
Add linear constraints to a model.

\sphinxAtStartPar
\sphinxstylestrong{Synopsis}
\begin{quote}

\sphinxAtStartPar
\sphinxcode{\sphinxupquote{ConstrArray AddConstrs(}}
\begin{quote}

\sphinxAtStartPar
\sphinxcode{\sphinxupquote{int count,}}

\sphinxAtStartPar
\sphinxcode{\sphinxupquote{double{[}{]} lbs,}}

\sphinxAtStartPar
\sphinxcode{\sphinxupquote{double{[}{]} ubs,}}

\sphinxAtStartPar
\sphinxcode{\sphinxupquote{string prefix)}}
\end{quote}
\end{quote}

\sphinxAtStartPar
\sphinxstylestrong{Arguments}
\begin{quote}

\sphinxAtStartPar
\sphinxcode{\sphinxupquote{count}}: number of constraints added to the model.

\sphinxAtStartPar
\sphinxcode{\sphinxupquote{lbs}}: lower bounds of new constraints.

\sphinxAtStartPar
\sphinxcode{\sphinxupquote{ubs}}: upper bounds of new constraints.

\sphinxAtStartPar
\sphinxcode{\sphinxupquote{prefix}}: optional, name prefix for new constraints, default value is ‘R’.
\end{quote}

\sphinxAtStartPar
\sphinxstylestrong{Return}
\begin{quote}

\sphinxAtStartPar
array of new constraint objects.
\end{quote}
\end{quote}

\subsubsection{Model.AddConstrs()}
\label{\detokenize{csapi/model:id17}}\begin{quote}

\sphinxAtStartPar
Add linear constraints to a model.

\sphinxAtStartPar
\sphinxstylestrong{Synopsis}
\begin{quote}

\sphinxAtStartPar
\sphinxcode{\sphinxupquote{ConstrArray AddConstrs(ConstrBuilderArray builders, string prefix)}}
\end{quote}

\sphinxAtStartPar
\sphinxstylestrong{Arguments}
\begin{quote}

\sphinxAtStartPar
\sphinxcode{\sphinxupquote{builders}}: builders for new constraints.

\sphinxAtStartPar
\sphinxcode{\sphinxupquote{prefix}}: optional, name prefix for new constraints, default value is ‘R’.
\end{quote}

\sphinxAtStartPar
\sphinxstylestrong{Return}
\begin{quote}

\sphinxAtStartPar
array of new constraint objects.
\end{quote}
\end{quote}

\subsubsection{Model.AddDenseMat()}
\label{\detokenize{csapi/model:model-adddensemat}}\begin{quote}

\sphinxAtStartPar
Add a dense symmetric matrix to a model.

\sphinxAtStartPar
\sphinxstylestrong{Synopsis}
\begin{quote}

\sphinxAtStartPar
\sphinxcode{\sphinxupquote{SymMatrix AddDenseMat(int dim, double{[}{]} vals)}}
\end{quote}

\sphinxAtStartPar
\sphinxstylestrong{Arguments}
\begin{quote}

\sphinxAtStartPar
\sphinxcode{\sphinxupquote{dim}}: dimension of the dense symmetric matrix.

\sphinxAtStartPar
\sphinxcode{\sphinxupquote{vals}}: array of non\sphinxhyphen{}zero elements, filled in column\sphinxhyphen{}wise up to len or max length of symmetric matrix.
\end{quote}

\sphinxAtStartPar
\sphinxstylestrong{Return}
\begin{quote}

\sphinxAtStartPar
new symmetric matrix object.
\end{quote}
\end{quote}

\subsubsection{Model.AddDenseMat()}
\label{\detokenize{csapi/model:id18}}\begin{quote}

\sphinxAtStartPar
Add a dense symmetric matrix to a model.

\sphinxAtStartPar
\sphinxstylestrong{Synopsis}
\begin{quote}

\sphinxAtStartPar
\sphinxcode{\sphinxupquote{SymMatrix AddDenseMat(int dim, double val)}}
\end{quote}

\sphinxAtStartPar
\sphinxstylestrong{Arguments}
\begin{quote}

\sphinxAtStartPar
\sphinxcode{\sphinxupquote{dim}}: dimension of dense symmetric matrix.

\sphinxAtStartPar
\sphinxcode{\sphinxupquote{val}}: value to fill dense symmetric matrix.
\end{quote}

\sphinxAtStartPar
\sphinxstylestrong{Return}
\begin{quote}

\sphinxAtStartPar
new symmetric matrix object.
\end{quote}
\end{quote}

\subsubsection{Model.AddDiagMat()}
\label{\detokenize{csapi/model:model-adddiagmat}}\begin{quote}

\sphinxAtStartPar
Add a diagonal matrix to a model.

\sphinxAtStartPar
\sphinxstylestrong{Synopsis}
\begin{quote}

\sphinxAtStartPar
\sphinxcode{\sphinxupquote{SymMatrix AddDiagMat(int dim, double val)}}
\end{quote}

\sphinxAtStartPar
\sphinxstylestrong{Arguments}
\begin{quote}

\sphinxAtStartPar
\sphinxcode{\sphinxupquote{dim}}: dimension of diagonal matrix.

\sphinxAtStartPar
\sphinxcode{\sphinxupquote{val}}: value to fill diagonal elements.
\end{quote}

\sphinxAtStartPar
\sphinxstylestrong{Return}
\begin{quote}

\sphinxAtStartPar
new diagonal matrix object.
\end{quote}
\end{quote}

\subsubsection{Model.AddDiagMat()}
\label{\detokenize{csapi/model:id19}}\begin{quote}

\sphinxAtStartPar
Add a diagonal matrix to a model.

\sphinxAtStartPar
\sphinxstylestrong{Synopsis}
\begin{quote}

\sphinxAtStartPar
\sphinxcode{\sphinxupquote{SymMatrix AddDiagMat(int dim, double{[}{]} vals)}}
\end{quote}

\sphinxAtStartPar
\sphinxstylestrong{Arguments}
\begin{quote}

\sphinxAtStartPar
\sphinxcode{\sphinxupquote{dim}}: dimension of diagonal matrix.

\sphinxAtStartPar
\sphinxcode{\sphinxupquote{vals}}: array of values of diagonal elements.
\end{quote}

\sphinxAtStartPar
\sphinxstylestrong{Return}
\begin{quote}

\sphinxAtStartPar
new diagonal matrix object.
\end{quote}
\end{quote}

\subsubsection{Model.AddDiagMat()}
\label{\detokenize{csapi/model:id20}}\begin{quote}

\sphinxAtStartPar
Add a diagonal matrix to a model.

\sphinxAtStartPar
\sphinxstylestrong{Synopsis}
\begin{quote}

\sphinxAtStartPar
\sphinxcode{\sphinxupquote{SymMatrix AddDiagMat(}}
\begin{quote}

\sphinxAtStartPar
\sphinxcode{\sphinxupquote{int dim,}}

\sphinxAtStartPar
\sphinxcode{\sphinxupquote{double val,}}

\sphinxAtStartPar
\sphinxcode{\sphinxupquote{int offset)}}
\end{quote}
\end{quote}

\sphinxAtStartPar
\sphinxstylestrong{Arguments}
\begin{quote}

\sphinxAtStartPar
\sphinxcode{\sphinxupquote{dim}}: dimension of diagonal matrix.

\sphinxAtStartPar
\sphinxcode{\sphinxupquote{val}}: value to fill diagonal elements.

\sphinxAtStartPar
\sphinxcode{\sphinxupquote{offset}}: shift distance against diagonal line.
\end{quote}

\sphinxAtStartPar
\sphinxstylestrong{Return}
\begin{quote}

\sphinxAtStartPar
new diagonal matrix object.
\end{quote}
\end{quote}

\subsubsection{Model.AddDiagMat()}
\label{\detokenize{csapi/model:id21}}\begin{quote}

\sphinxAtStartPar
Add a diagonal matrix to a model.

\sphinxAtStartPar
\sphinxstylestrong{Synopsis}
\begin{quote}

\sphinxAtStartPar
\sphinxcode{\sphinxupquote{SymMatrix AddDiagMat(}}
\begin{quote}

\sphinxAtStartPar
\sphinxcode{\sphinxupquote{int dim,}}

\sphinxAtStartPar
\sphinxcode{\sphinxupquote{double{[}{]} vals,}}

\sphinxAtStartPar
\sphinxcode{\sphinxupquote{int offset)}}
\end{quote}
\end{quote}

\sphinxAtStartPar
\sphinxstylestrong{Arguments}
\begin{quote}

\sphinxAtStartPar
\sphinxcode{\sphinxupquote{dim}}: dimension of diagonal matrix.

\sphinxAtStartPar
\sphinxcode{\sphinxupquote{vals}}: array of values of diagonal elements.

\sphinxAtStartPar
\sphinxcode{\sphinxupquote{offset}}: shift distance against diagonal line.
\end{quote}

\sphinxAtStartPar
\sphinxstylestrong{Return}
\begin{quote}

\sphinxAtStartPar
new diagonal matrix object.
\end{quote}
\end{quote}

\subsubsection{Model.AddExpCone()}
\label{\detokenize{csapi/model:model-addexpcone}}\begin{quote}

\sphinxAtStartPar
Add an exponential cone constraint to model.

\sphinxAtStartPar
\sphinxstylestrong{Synopsis}
\begin{quote}

\sphinxAtStartPar
\sphinxcode{\sphinxupquote{ExpCone AddExpCone(}}
\begin{quote}

\sphinxAtStartPar
\sphinxcode{\sphinxupquote{int type,}}

\sphinxAtStartPar
\sphinxcode{\sphinxupquote{char{[}{]} pvtype,}}

\sphinxAtStartPar
\sphinxcode{\sphinxupquote{string prefix)}}
\end{quote}
\end{quote}

\sphinxAtStartPar
\sphinxstylestrong{Arguments}
\begin{quote}

\sphinxAtStartPar
\sphinxcode{\sphinxupquote{type}}: type of the exponential cone constraint.

\sphinxAtStartPar
\sphinxcode{\sphinxupquote{pvtype}}: type of variables in the exponential cone.

\sphinxAtStartPar
\sphinxcode{\sphinxupquote{prefix}}: optional, name prefix of variables in the exponential cone, default value is “ExpConeV”.
\end{quote}

\sphinxAtStartPar
\sphinxstylestrong{Return}
\begin{quote}

\sphinxAtStartPar
new exponential cone constraint object.
\end{quote}
\end{quote}

\subsubsection{Model.AddExpCone()}
\label{\detokenize{csapi/model:id22}}\begin{quote}

\sphinxAtStartPar
Add an exponential cone constraint to model.

\sphinxAtStartPar
\sphinxstylestrong{Synopsis}
\begin{quote}

\sphinxAtStartPar
\sphinxcode{\sphinxupquote{ExpCone AddExpCone(ExpConeBuilder builder)}}
\end{quote}

\sphinxAtStartPar
\sphinxstylestrong{Arguments}
\begin{quote}

\sphinxAtStartPar
\sphinxcode{\sphinxupquote{builder}}: builder for new exponential cone constraint.
\end{quote}

\sphinxAtStartPar
\sphinxstylestrong{Return}
\begin{quote}

\sphinxAtStartPar
new exponential cone constraint object.
\end{quote}
\end{quote}

\subsubsection{Model.AddExpCone()}
\label{\detokenize{csapi/model:id23}}\begin{quote}

\sphinxAtStartPar
Add an exponential cone constraint to model.

\sphinxAtStartPar
\sphinxstylestrong{Synopsis}
\begin{quote}

\sphinxAtStartPar
\sphinxcode{\sphinxupquote{ExpCone AddExpCone(Var{[}{]} vars, int type)}}
\end{quote}

\sphinxAtStartPar
\sphinxstylestrong{Arguments}
\begin{quote}

\sphinxAtStartPar
\sphinxcode{\sphinxupquote{vars}}: variables that participate in the exponential cone constraint.

\sphinxAtStartPar
\sphinxcode{\sphinxupquote{type}}: type of the exponential cone constraint.
\end{quote}

\sphinxAtStartPar
\sphinxstylestrong{Return}
\begin{quote}

\sphinxAtStartPar
new exponential cone constraint object.
\end{quote}
\end{quote}

\subsubsection{Model.AddExpCone()}
\label{\detokenize{csapi/model:id24}}\begin{quote}

\sphinxAtStartPar
Add an exponential cone constraint to model.

\sphinxAtStartPar
\sphinxstylestrong{Synopsis}
\begin{quote}

\sphinxAtStartPar
\sphinxcode{\sphinxupquote{ExpCone AddExpCone(VarArray vars, int type)}}
\end{quote}

\sphinxAtStartPar
\sphinxstylestrong{Arguments}
\begin{quote}

\sphinxAtStartPar
\sphinxcode{\sphinxupquote{vars}}: variables that participate in the exponential cone constraint.

\sphinxAtStartPar
\sphinxcode{\sphinxupquote{type}}: type of an exponential cone constraint.
\end{quote}

\sphinxAtStartPar
\sphinxstylestrong{Return}
\begin{quote}

\sphinxAtStartPar
new exponential cone constraint object.
\end{quote}
\end{quote}

\subsubsection{Model.AddEyeMat()}
\label{\detokenize{csapi/model:model-addeyemat}}\begin{quote}

\sphinxAtStartPar
Add an identity matrix to a model.

\sphinxAtStartPar
\sphinxstylestrong{Synopsis}
\begin{quote}

\sphinxAtStartPar
\sphinxcode{\sphinxupquote{SymMatrix AddEyeMat(int dim)}}
\end{quote}

\sphinxAtStartPar
\sphinxstylestrong{Arguments}
\begin{quote}

\sphinxAtStartPar
\sphinxcode{\sphinxupquote{dim}}: dimension of identity matrix.
\end{quote}

\sphinxAtStartPar
\sphinxstylestrong{Return}
\begin{quote}

\sphinxAtStartPar
new identity matrix object.
\end{quote}
\end{quote}

\subsubsection{Model.AddGenConstrIndicator()}
\label{\detokenize{csapi/model:model-addgenconstrindicator}}\begin{quote}

\sphinxAtStartPar
Add a general constraint of type indicator to model.

\sphinxAtStartPar
\sphinxstylestrong{Synopsis}
\begin{quote}

\sphinxAtStartPar
\sphinxcode{\sphinxupquote{GenConstr AddGenConstrIndicator(GenConstrBuilder builder, string name)}}
\end{quote}

\sphinxAtStartPar
\sphinxstylestrong{Arguments}
\begin{quote}

\sphinxAtStartPar
\sphinxcode{\sphinxupquote{builder}}: builder for the general constraint.

\sphinxAtStartPar
\sphinxcode{\sphinxupquote{name}}: optional, name of new general constraint.
\end{quote}

\sphinxAtStartPar
\sphinxstylestrong{Return}
\begin{quote}

\sphinxAtStartPar
new general constraint object of type indicator.
\end{quote}
\end{quote}

\subsubsection{Model.AddGenConstrIndicator()}
\label{\detokenize{csapi/model:id25}}\begin{quote}

\sphinxAtStartPar
Add a general constraint of type indicator to model.

\sphinxAtStartPar
\sphinxstylestrong{Synopsis}
\begin{quote}

\sphinxAtStartPar
\sphinxcode{\sphinxupquote{GenConstr AddGenConstrIndicator(}}
\begin{quote}

\sphinxAtStartPar
\sphinxcode{\sphinxupquote{Var binvar,}}

\sphinxAtStartPar
\sphinxcode{\sphinxupquote{int binval,}}

\sphinxAtStartPar
\sphinxcode{\sphinxupquote{ConstrBuilder builder,}}

\sphinxAtStartPar
\sphinxcode{\sphinxupquote{int type,}}

\sphinxAtStartPar
\sphinxcode{\sphinxupquote{string name)}}
\end{quote}
\end{quote}

\sphinxAtStartPar
\sphinxstylestrong{Arguments}
\begin{quote}

\sphinxAtStartPar
\sphinxcode{\sphinxupquote{binvar}}: binary indicator variable.

\sphinxAtStartPar
\sphinxcode{\sphinxupquote{binval}}: value for binary indicator variable that force a linear constraint to be satisfied(0 or 1).

\sphinxAtStartPar
\sphinxcode{\sphinxupquote{builder}}: builder for linear constraint.

\sphinxAtStartPar
\sphinxcode{\sphinxupquote{type}}: type of general constraint.

\sphinxAtStartPar
\sphinxcode{\sphinxupquote{name}}: optional, name of new general constraint.
\end{quote}

\sphinxAtStartPar
\sphinxstylestrong{Return}
\begin{quote}

\sphinxAtStartPar
new general constraint object of type indicator.
\end{quote}
\end{quote}

\subsubsection{Model.AddGenConstrIndicator()}
\label{\detokenize{csapi/model:id26}}\begin{quote}

\sphinxAtStartPar
Add a general constraint of type indicator to model.

\sphinxAtStartPar
\sphinxstylestrong{Synopsis}
\begin{quote}

\sphinxAtStartPar
\sphinxcode{\sphinxupquote{GenConstr AddGenConstrIndicator(}}
\begin{quote}

\sphinxAtStartPar
\sphinxcode{\sphinxupquote{Var binvar,}}

\sphinxAtStartPar
\sphinxcode{\sphinxupquote{int binval,}}

\sphinxAtStartPar
\sphinxcode{\sphinxupquote{Expr expr,}}

\sphinxAtStartPar
\sphinxcode{\sphinxupquote{char sense,}}

\sphinxAtStartPar
\sphinxcode{\sphinxupquote{double rhs,}}

\sphinxAtStartPar
\sphinxcode{\sphinxupquote{int type,}}

\sphinxAtStartPar
\sphinxcode{\sphinxupquote{string name)}}
\end{quote}
\end{quote}

\sphinxAtStartPar
\sphinxstylestrong{Arguments}
\begin{quote}

\sphinxAtStartPar
\sphinxcode{\sphinxupquote{binvar}}: binary indicator variable.

\sphinxAtStartPar
\sphinxcode{\sphinxupquote{binval}}: value for binary indicator variable that force a linear constraint to be satisfied(0 or 1).

\sphinxAtStartPar
\sphinxcode{\sphinxupquote{expr}}: expression for new linear contraint.

\sphinxAtStartPar
\sphinxcode{\sphinxupquote{sense}}: sense for new linear constraint.

\sphinxAtStartPar
\sphinxcode{\sphinxupquote{rhs}}: right hand side value for new linear constraint.

\sphinxAtStartPar
\sphinxcode{\sphinxupquote{type}}: type of general constraint.

\sphinxAtStartPar
\sphinxcode{\sphinxupquote{name}}: optional, name of new general constraint.
\end{quote}

\sphinxAtStartPar
\sphinxstylestrong{Return}
\begin{quote}

\sphinxAtStartPar
new general constraint object of type indicator.
\end{quote}
\end{quote}

\subsubsection{Model.AddGenConstrIndicators()}
\label{\detokenize{csapi/model:model-addgenconstrindicators}}\begin{quote}

\sphinxAtStartPar
Add general constraints to a model.

\sphinxAtStartPar
\sphinxstylestrong{Synopsis}
\begin{quote}

\sphinxAtStartPar
\sphinxcode{\sphinxupquote{GenConstrArray AddGenConstrIndicators(GenConstrBuilderArray builders, string prefix)}}
\end{quote}

\sphinxAtStartPar
\sphinxstylestrong{Arguments}
\begin{quote}

\sphinxAtStartPar
\sphinxcode{\sphinxupquote{builders}}: builders for new general constraints.

\sphinxAtStartPar
\sphinxcode{\sphinxupquote{prefix}}: optional, name prefix for new general constraints.
\end{quote}

\sphinxAtStartPar
\sphinxstylestrong{Return}
\begin{quote}

\sphinxAtStartPar
array of new general constraint objects.
\end{quote}
\end{quote}

\subsubsection{Model.AddLazyConstr()}
\label{\detokenize{csapi/model:model-addlazyconstr}}\begin{quote}

\sphinxAtStartPar
Add a lazy constraint to model.

\sphinxAtStartPar
\sphinxstylestrong{Synopsis}
\begin{quote}

\sphinxAtStartPar
\sphinxcode{\sphinxupquote{void AddLazyConstr(}}
\begin{quote}

\sphinxAtStartPar
\sphinxcode{\sphinxupquote{Expr lhs,}}

\sphinxAtStartPar
\sphinxcode{\sphinxupquote{char sense,}}

\sphinxAtStartPar
\sphinxcode{\sphinxupquote{double rhs,}}

\sphinxAtStartPar
\sphinxcode{\sphinxupquote{string name)}}
\end{quote}
\end{quote}

\sphinxAtStartPar
\sphinxstylestrong{Arguments}
\begin{quote}

\sphinxAtStartPar
\sphinxcode{\sphinxupquote{lhs}}: expression for lazy contraint.

\sphinxAtStartPar
\sphinxcode{\sphinxupquote{sense}}: sense for lazy constraint.

\sphinxAtStartPar
\sphinxcode{\sphinxupquote{rhs}}: right hand side value for lazy constraint.

\sphinxAtStartPar
\sphinxcode{\sphinxupquote{name}}: optional, name of lazy constraint.
\end{quote}
\end{quote}

\subsubsection{Model.AddLazyConstr()}
\label{\detokenize{csapi/model:id27}}\begin{quote}

\sphinxAtStartPar
Add a lazy constraint to model.

\sphinxAtStartPar
\sphinxstylestrong{Synopsis}
\begin{quote}

\sphinxAtStartPar
\sphinxcode{\sphinxupquote{void AddLazyConstr(}}
\begin{quote}

\sphinxAtStartPar
\sphinxcode{\sphinxupquote{Expr lhs,}}

\sphinxAtStartPar
\sphinxcode{\sphinxupquote{char sense,}}

\sphinxAtStartPar
\sphinxcode{\sphinxupquote{Expr rhs,}}

\sphinxAtStartPar
\sphinxcode{\sphinxupquote{string name)}}
\end{quote}
\end{quote}

\sphinxAtStartPar
\sphinxstylestrong{Arguments}
\begin{quote}

\sphinxAtStartPar
\sphinxcode{\sphinxupquote{lhs}}: left hand side expression for lazy contraint.

\sphinxAtStartPar
\sphinxcode{\sphinxupquote{sense}}: sense for lazy constraint.

\sphinxAtStartPar
\sphinxcode{\sphinxupquote{rhs}}: right hand side expression for lazy contraint.

\sphinxAtStartPar
\sphinxcode{\sphinxupquote{name}}: optional, name of lazy constraint.
\end{quote}
\end{quote}

\subsubsection{Model.AddLazyConstr()}
\label{\detokenize{csapi/model:id28}}\begin{quote}

\sphinxAtStartPar
Add a lazy constraint to model.

\sphinxAtStartPar
\sphinxstylestrong{Synopsis}
\begin{quote}

\sphinxAtStartPar
\sphinxcode{\sphinxupquote{void AddLazyConstr(ConstrBuilder builder, string name)}}
\end{quote}

\sphinxAtStartPar
\sphinxstylestrong{Arguments}
\begin{quote}

\sphinxAtStartPar
\sphinxcode{\sphinxupquote{builder}}: builder for lazy contraint.

\sphinxAtStartPar
\sphinxcode{\sphinxupquote{name}}: optional, name of lazy constraint.
\end{quote}
\end{quote}

\subsubsection{Model.AddLazyConstrs()}
\label{\detokenize{csapi/model:model-addlazyconstrs}}\begin{quote}

\sphinxAtStartPar
Add lazy constraints to model.

\sphinxAtStartPar
\sphinxstylestrong{Synopsis}
\begin{quote}

\sphinxAtStartPar
\sphinxcode{\sphinxupquote{void AddLazyConstrs(ConstrBuilderArray builders, string prefix)}}
\end{quote}

\sphinxAtStartPar
\sphinxstylestrong{Arguments}
\begin{quote}

\sphinxAtStartPar
\sphinxcode{\sphinxupquote{builders}}: array of builders for lazy contraints.

\sphinxAtStartPar
\sphinxcode{\sphinxupquote{prefix}}: name prefix of new lazy constraints.
\end{quote}
\end{quote}

\subsubsection{Model.AddLmiConstr()}
\label{\detokenize{csapi/model:model-addlmiconstr}}\begin{quote}

\sphinxAtStartPar
Add an LMI constraint to model.

\sphinxAtStartPar
\sphinxstylestrong{Synopsis}
\begin{quote}

\sphinxAtStartPar
\sphinxcode{\sphinxupquote{LmiConstraint AddLmiConstr(LmiExpr expr, string name)}}
\end{quote}

\sphinxAtStartPar
\sphinxstylestrong{Arguments}
\begin{quote}

\sphinxAtStartPar
\sphinxcode{\sphinxupquote{expr}}: LMI expression for new LMI contraint.

\sphinxAtStartPar
\sphinxcode{\sphinxupquote{name}}: optional, name of new LMI constraint.
\end{quote}

\sphinxAtStartPar
\sphinxstylestrong{Return}
\begin{quote}

\sphinxAtStartPar
new LMI constraint object.
\end{quote}
\end{quote}

\subsubsection{Model.AddMConstr()}
\label{\detokenize{csapi/model:model-addmconstr}}\begin{quote}

\sphinxAtStartPar
Add a MConstr object in N\sphinxhyphen{}dimensions to model.

\sphinxAtStartPar
\sphinxstylestrong{Synopsis}
\begin{quote}

\sphinxAtStartPar
\sphinxcode{\sphinxupquote{MConstr AddMConstr(MConstrBuilder builder, string prefix)}}
\end{quote}

\sphinxAtStartPar
\sphinxstylestrong{Arguments}
\begin{quote}

\sphinxAtStartPar
\sphinxcode{\sphinxupquote{builder}}: builder for MConstr object.

\sphinxAtStartPar
\sphinxcode{\sphinxupquote{prefix}}: name prefix for constraints in MConstr object.
\end{quote}

\sphinxAtStartPar
\sphinxstylestrong{Return}
\begin{quote}

\sphinxAtStartPar
new MConstr object.
\end{quote}
\end{quote}

\subsubsection{Model.AddMConstr()}
\label{\detokenize{csapi/model:id29}}\begin{quote}

\sphinxAtStartPar
Add a N\sphinxhyphen{}dimensional MConstr object to model.

\sphinxAtStartPar
\sphinxstylestrong{Synopsis}
\begin{quote}

\sphinxAtStartPar
\sphinxcode{\sphinxupquote{MConstr AddMConstr(}}
\begin{quote}

\sphinxAtStartPar
\sphinxcode{\sphinxupquote{MLinExpr exprs,}}

\sphinxAtStartPar
\sphinxcode{\sphinxupquote{char sense,}}

\sphinxAtStartPar
\sphinxcode{\sphinxupquote{double rhs,}}

\sphinxAtStartPar
\sphinxcode{\sphinxupquote{string prefix)}}
\end{quote}
\end{quote}

\sphinxAtStartPar
\sphinxstylestrong{Arguments}
\begin{quote}

\sphinxAtStartPar
\sphinxcode{\sphinxupquote{exprs}}: N\sphinxhyphen{}dimensional MLinExpr object.

\sphinxAtStartPar
\sphinxcode{\sphinxupquote{sense}}: sense for new linear constraints.

\sphinxAtStartPar
\sphinxcode{\sphinxupquote{rhs}}: double value at right side of the new linear constraints.

\sphinxAtStartPar
\sphinxcode{\sphinxupquote{prefix}}: name prefix for constraints in MConstr object.
\end{quote}

\sphinxAtStartPar
\sphinxstylestrong{Return}
\begin{quote}

\sphinxAtStartPar
new MConstr object.
\end{quote}
\end{quote}

\subsubsection{Model.AddMPsdConstr()}
\label{\detokenize{csapi/model:model-addmpsdconstr}}\begin{quote}

\sphinxAtStartPar
Add a N\sphinxhyphen{}dimensional MPsdConstr object to model.

\sphinxAtStartPar
\sphinxstylestrong{Synopsis}
\begin{quote}

\sphinxAtStartPar
\sphinxcode{\sphinxupquote{MPsdConstr AddMPsdConstr(}}
\begin{quote}

\sphinxAtStartPar
\sphinxcode{\sphinxupquote{MPsdExpr exprs,}}

\sphinxAtStartPar
\sphinxcode{\sphinxupquote{char sense,}}

\sphinxAtStartPar
\sphinxcode{\sphinxupquote{double rhs,}}

\sphinxAtStartPar
\sphinxcode{\sphinxupquote{string prefix)}}
\end{quote}
\end{quote}

\sphinxAtStartPar
\sphinxstylestrong{Arguments}
\begin{quote}

\sphinxAtStartPar
\sphinxcode{\sphinxupquote{exprs}}: N\sphinxhyphen{}dimensional MPsdExpr object.

\sphinxAtStartPar
\sphinxcode{\sphinxupquote{sense}}: sense for new PSD constraints.

\sphinxAtStartPar
\sphinxcode{\sphinxupquote{rhs}}: double value at right side of the new PSD constraints.

\sphinxAtStartPar
\sphinxcode{\sphinxupquote{prefix}}: name prefix of PSD constraints in MPsdConstr object.
\end{quote}

\sphinxAtStartPar
\sphinxstylestrong{Return}
\begin{quote}

\sphinxAtStartPar
new MPsdConstr object.
\end{quote}
\end{quote}

\subsubsection{Model.AddMPsdConstr()}
\label{\detokenize{csapi/model:id30}}\begin{quote}

\sphinxAtStartPar
Add a N\sphinxhyphen{}dimensional MPsdConstr object to model.

\sphinxAtStartPar
\sphinxstylestrong{Synopsis}
\begin{quote}

\sphinxAtStartPar
\sphinxcode{\sphinxupquote{MPsdConstr AddMPsdConstr(MPsdConstrBuilder builder, string prefix)}}
\end{quote}

\sphinxAtStartPar
\sphinxstylestrong{Arguments}
\begin{quote}

\sphinxAtStartPar
\sphinxcode{\sphinxupquote{builder}}: builder for MPsdConstr object.

\sphinxAtStartPar
\sphinxcode{\sphinxupquote{prefix}}: name prefix of PSD constraints in MPsdConstr object.
\end{quote}

\sphinxAtStartPar
\sphinxstylestrong{Return}
\begin{quote}

\sphinxAtStartPar
new MPsdConstr object.
\end{quote}
\end{quote}

\subsubsection{Model.AddMQConstr()}
\label{\detokenize{csapi/model:model-addmqconstr}}\begin{quote}

\sphinxAtStartPar
Add a N\sphinxhyphen{}dimensional MQConstr object to model.

\sphinxAtStartPar
\sphinxstylestrong{Synopsis}
\begin{quote}

\sphinxAtStartPar
\sphinxcode{\sphinxupquote{MQConstr AddMQConstr(MQConstrBuilder builder, string prefix)}}
\end{quote}

\sphinxAtStartPar
\sphinxstylestrong{Arguments}
\begin{quote}

\sphinxAtStartPar
\sphinxcode{\sphinxupquote{builder}}: builder for MQConstr object.

\sphinxAtStartPar
\sphinxcode{\sphinxupquote{prefix}}: name prefix of quadratic constraints in MQConstr object.
\end{quote}

\sphinxAtStartPar
\sphinxstylestrong{Return}
\begin{quote}

\sphinxAtStartPar
new MQConstr object.
\end{quote}
\end{quote}

\subsubsection{Model.AddMQConstr()}
\label{\detokenize{csapi/model:id31}}\begin{quote}

\sphinxAtStartPar
Add a N\sphinxhyphen{}dimensional MQConstr object to model.

\sphinxAtStartPar
\sphinxstylestrong{Synopsis}
\begin{quote}

\sphinxAtStartPar
\sphinxcode{\sphinxupquote{MQConstr AddMQConstr(}}
\begin{quote}

\sphinxAtStartPar
\sphinxcode{\sphinxupquote{MQuadExpr exprs,}}

\sphinxAtStartPar
\sphinxcode{\sphinxupquote{char sense,}}

\sphinxAtStartPar
\sphinxcode{\sphinxupquote{double rhs,}}

\sphinxAtStartPar
\sphinxcode{\sphinxupquote{string prefix)}}
\end{quote}
\end{quote}

\sphinxAtStartPar
\sphinxstylestrong{Arguments}
\begin{quote}

\sphinxAtStartPar
\sphinxcode{\sphinxupquote{exprs}}: N\sphinxhyphen{}dimensional MQuadExpr object.

\sphinxAtStartPar
\sphinxcode{\sphinxupquote{sense}}: sense for new quadratic constraints.

\sphinxAtStartPar
\sphinxcode{\sphinxupquote{rhs}}: double value at right side of the new quadratic constraints.

\sphinxAtStartPar
\sphinxcode{\sphinxupquote{prefix}}: name prefix of quadratic constraints in MQConstr object.
\end{quote}

\sphinxAtStartPar
\sphinxstylestrong{Return}
\begin{quote}

\sphinxAtStartPar
new MQConstr object.
\end{quote}
\end{quote}

\subsubsection{Model.AddMVar()}
\label{\detokenize{csapi/model:model-addmvar}}\begin{quote}

\sphinxAtStartPar
Add a MVar object in N\sphinxhyphen{}dimensions to model.

\sphinxAtStartPar
\sphinxstylestrong{Synopsis}
\begin{quote}

\sphinxAtStartPar
\sphinxcode{\sphinxupquote{MVar AddMVar(}}
\begin{quote}

\sphinxAtStartPar
\sphinxcode{\sphinxupquote{Shape shp,}}

\sphinxAtStartPar
\sphinxcode{\sphinxupquote{char vtype,}}

\sphinxAtStartPar
\sphinxcode{\sphinxupquote{string prefix)}}
\end{quote}
\end{quote}

\sphinxAtStartPar
\sphinxstylestrong{Arguments}
\begin{quote}

\sphinxAtStartPar
\sphinxcode{\sphinxupquote{shp}}: shape of MVar object.

\sphinxAtStartPar
\sphinxcode{\sphinxupquote{vtype}}: type of variables in MVar object.

\sphinxAtStartPar
\sphinxcode{\sphinxupquote{prefix}}: name prefix of variables in MVar object.
\end{quote}

\sphinxAtStartPar
\sphinxstylestrong{Return}
\begin{quote}

\sphinxAtStartPar
new MVar object.
\end{quote}
\end{quote}

\subsubsection{Model.AddMVar()}
\label{\detokenize{csapi/model:id32}}\begin{quote}

\sphinxAtStartPar
Add a MVar object in N\sphinxhyphen{}dimensions to model.

\sphinxAtStartPar
\sphinxstylestrong{Synopsis}
\begin{quote}

\sphinxAtStartPar
\sphinxcode{\sphinxupquote{MVar AddMVar(}}
\begin{quote}

\sphinxAtStartPar
\sphinxcode{\sphinxupquote{Shape shp,}}

\sphinxAtStartPar
\sphinxcode{\sphinxupquote{double lb,}}

\sphinxAtStartPar
\sphinxcode{\sphinxupquote{double ub,}}

\sphinxAtStartPar
\sphinxcode{\sphinxupquote{double obj,}}

\sphinxAtStartPar
\sphinxcode{\sphinxupquote{char vtype,}}

\sphinxAtStartPar
\sphinxcode{\sphinxupquote{string prefix)}}
\end{quote}
\end{quote}

\sphinxAtStartPar
\sphinxstylestrong{Arguments}
\begin{quote}

\sphinxAtStartPar
\sphinxcode{\sphinxupquote{shp}}: shape of MVar object.

\sphinxAtStartPar
\sphinxcode{\sphinxupquote{lb}}: lower bound for variables in MVar object.

\sphinxAtStartPar
\sphinxcode{\sphinxupquote{ub}}: upper bound for variables in MVar object.

\sphinxAtStartPar
\sphinxcode{\sphinxupquote{obj}}: objective coefficient for variables in MVar object.

\sphinxAtStartPar
\sphinxcode{\sphinxupquote{vtype}}: type of variables in MVar object.

\sphinxAtStartPar
\sphinxcode{\sphinxupquote{prefix}}: name prefix of variables in MVar object.
\end{quote}

\sphinxAtStartPar
\sphinxstylestrong{Return}
\begin{quote}

\sphinxAtStartPar
new MVar object.
\end{quote}
\end{quote}

\subsubsection{Model.AddMVar()}
\label{\detokenize{csapi/model:id33}}\begin{quote}

\sphinxAtStartPar
Add a MVar object in N\sphinxhyphen{}dimensions to model.

\sphinxAtStartPar
\sphinxstylestrong{Synopsis}
\begin{quote}

\sphinxAtStartPar
\sphinxcode{\sphinxupquote{MVar AddMVar(}}
\begin{quote}

\sphinxAtStartPar
\sphinxcode{\sphinxupquote{Shape shp,}}

\sphinxAtStartPar
\sphinxcode{\sphinxupquote{double{[}{]} lbs,}}

\sphinxAtStartPar
\sphinxcode{\sphinxupquote{double{[}{]} ubs,}}

\sphinxAtStartPar
\sphinxcode{\sphinxupquote{double{[}{]} objs,}}

\sphinxAtStartPar
\sphinxcode{\sphinxupquote{char{[}{]} types,}}

\sphinxAtStartPar
\sphinxcode{\sphinxupquote{string prefix)}}
\end{quote}
\end{quote}

\sphinxAtStartPar
\sphinxstylestrong{Arguments}
\begin{quote}

\sphinxAtStartPar
\sphinxcode{\sphinxupquote{shp}}: shape of MVar object.

\sphinxAtStartPar
\sphinxcode{\sphinxupquote{lbs}}: lower bounds for new variables. if NULL, lower bounds are 0.0.

\sphinxAtStartPar
\sphinxcode{\sphinxupquote{ubs}}: upper bounds for new variables. if NULL, upper bounds are infinity or 1 for binary variables.

\sphinxAtStartPar
\sphinxcode{\sphinxupquote{objs}}: objective coefficients for new variables. if NULL, objective coefficients are 0.0.

\sphinxAtStartPar
\sphinxcode{\sphinxupquote{types}}: variable types for new variables. if NULL, variable types are continuous.

\sphinxAtStartPar
\sphinxcode{\sphinxupquote{prefix}}: name prefix of variables in MVar object.
\end{quote}

\sphinxAtStartPar
\sphinxstylestrong{Return}
\begin{quote}

\sphinxAtStartPar
new MVar object.
\end{quote}
\end{quote}

\subsubsection{Model.AddNlConstr()}
\label{\detokenize{csapi/model:model-addnlconstr}}\begin{quote}

\sphinxAtStartPar
Add a nonlinear constraint to model.

\sphinxAtStartPar
\sphinxstylestrong{Synopsis}
\begin{quote}

\sphinxAtStartPar
\sphinxcode{\sphinxupquote{NlConstraint AddNlConstr(}}
\begin{quote}

\sphinxAtStartPar
\sphinxcode{\sphinxupquote{NlExpr expr,}}

\sphinxAtStartPar
\sphinxcode{\sphinxupquote{char sense,}}

\sphinxAtStartPar
\sphinxcode{\sphinxupquote{double rhs,}}

\sphinxAtStartPar
\sphinxcode{\sphinxupquote{string name)}}
\end{quote}
\end{quote}

\sphinxAtStartPar
\sphinxstylestrong{Arguments}
\begin{quote}

\sphinxAtStartPar
\sphinxcode{\sphinxupquote{expr}}: non\sphinxhyphen{}expression for the new contraint.

\sphinxAtStartPar
\sphinxcode{\sphinxupquote{sense}}: sense for new nonlinear constraint, other than range sense.

\sphinxAtStartPar
\sphinxcode{\sphinxupquote{rhs}}: right hand side value for the new constraint.

\sphinxAtStartPar
\sphinxcode{\sphinxupquote{name}}: optional, name of new nonlinear constraint.
\end{quote}

\sphinxAtStartPar
\sphinxstylestrong{Return}
\begin{quote}

\sphinxAtStartPar
new nonlinear constraint object.
\end{quote}
\end{quote}

\subsubsection{Model.AddNlConstr()}
\label{\detokenize{csapi/model:id34}}\begin{quote}

\sphinxAtStartPar
Add a nonlinear constraint to model.

\sphinxAtStartPar
\sphinxstylestrong{Synopsis}
\begin{quote}

\sphinxAtStartPar
\sphinxcode{\sphinxupquote{NlConstraint AddNlConstr(}}
\begin{quote}

\sphinxAtStartPar
\sphinxcode{\sphinxupquote{NlExpr lhs,}}

\sphinxAtStartPar
\sphinxcode{\sphinxupquote{char sense,}}

\sphinxAtStartPar
\sphinxcode{\sphinxupquote{NlExpr rhs,}}

\sphinxAtStartPar
\sphinxcode{\sphinxupquote{string name)}}
\end{quote}
\end{quote}

\sphinxAtStartPar
\sphinxstylestrong{Arguments}
\begin{quote}

\sphinxAtStartPar
\sphinxcode{\sphinxupquote{lhs}}: left hand side nonlinear expression for the new constraint.

\sphinxAtStartPar
\sphinxcode{\sphinxupquote{sense}}: sense for new nonlinear constraint, other than range sense.

\sphinxAtStartPar
\sphinxcode{\sphinxupquote{rhs}}: right hand side nonlinear expression for the new constraint.

\sphinxAtStartPar
\sphinxcode{\sphinxupquote{name}}: optional, name of new nonlinear constraint.
\end{quote}

\sphinxAtStartPar
\sphinxstylestrong{Return}
\begin{quote}

\sphinxAtStartPar
new nonlinear constraint object.
\end{quote}
\end{quote}

\subsubsection{Model.AddNlConstr()}
\label{\detokenize{csapi/model:id35}}\begin{quote}

\sphinxAtStartPar
Add a nonlinear constraint to model.

\sphinxAtStartPar
\sphinxstylestrong{Synopsis}
\begin{quote}

\sphinxAtStartPar
\sphinxcode{\sphinxupquote{NlConstraint AddNlConstr(}}
\begin{quote}

\sphinxAtStartPar
\sphinxcode{\sphinxupquote{NlExpr expr,}}

\sphinxAtStartPar
\sphinxcode{\sphinxupquote{double lb,}}

\sphinxAtStartPar
\sphinxcode{\sphinxupquote{double ub,}}

\sphinxAtStartPar
\sphinxcode{\sphinxupquote{string name)}}
\end{quote}
\end{quote}

\sphinxAtStartPar
\sphinxstylestrong{Arguments}
\begin{quote}

\sphinxAtStartPar
\sphinxcode{\sphinxupquote{expr}}: nonlinear expression for the new constraint.

\sphinxAtStartPar
\sphinxcode{\sphinxupquote{lb}}: lower bound for the new nonlinear constraint.

\sphinxAtStartPar
\sphinxcode{\sphinxupquote{ub}}: upper bound for the new nonlinear constraint

\sphinxAtStartPar
\sphinxcode{\sphinxupquote{name}}: optional, name of new constraint.
\end{quote}

\sphinxAtStartPar
\sphinxstylestrong{Return}
\begin{quote}

\sphinxAtStartPar
new nonlinear constraint object.
\end{quote}
\end{quote}

\subsubsection{Model.AddNlConstr()}
\label{\detokenize{csapi/model:id36}}\begin{quote}

\sphinxAtStartPar
Add a nonlinear constraint to a model.

\sphinxAtStartPar
\sphinxstylestrong{Synopsis}
\begin{quote}

\sphinxAtStartPar
\sphinxcode{\sphinxupquote{NlConstraint AddNlConstr(NlConstrBuilder builder, string name)}}
\end{quote}

\sphinxAtStartPar
\sphinxstylestrong{Arguments}
\begin{quote}

\sphinxAtStartPar
\sphinxcode{\sphinxupquote{builder}}: builder for the new nonlinear constraint.

\sphinxAtStartPar
\sphinxcode{\sphinxupquote{name}}: optional, name of new nonlinear constraint.
\end{quote}

\sphinxAtStartPar
\sphinxstylestrong{Return}
\begin{quote}

\sphinxAtStartPar
new nonlinear constraint object.
\end{quote}
\end{quote}

\subsubsection{Model.AddNlConstrs()}
\label{\detokenize{csapi/model:model-addnlconstrs}}\begin{quote}

\sphinxAtStartPar
Add nonlinear constraints to a model.

\sphinxAtStartPar
\sphinxstylestrong{Synopsis}
\begin{quote}

\sphinxAtStartPar
\sphinxcode{\sphinxupquote{NlConstrArray AddNlConstrs(NlConstrBuilderArray builders, string prefix)}}
\end{quote}

\sphinxAtStartPar
\sphinxstylestrong{Arguments}
\begin{quote}

\sphinxAtStartPar
\sphinxcode{\sphinxupquote{builders}}: builders for new nonlinear constraints.

\sphinxAtStartPar
\sphinxcode{\sphinxupquote{prefix}}: name prefix for new constraints.
\end{quote}

\sphinxAtStartPar
\sphinxstylestrong{Return}
\begin{quote}

\sphinxAtStartPar
array of new nonlinear constraint objects.
\end{quote}
\end{quote}

\subsubsection{Model.AddOnesMat()}
\label{\detokenize{csapi/model:model-addonesmat}}\begin{quote}

\sphinxAtStartPar
Add a dense symmetric matrix of value one to a model.

\sphinxAtStartPar
\sphinxstylestrong{Synopsis}
\begin{quote}

\sphinxAtStartPar
\sphinxcode{\sphinxupquote{SymMatrix AddOnesMat(int dim)}}
\end{quote}

\sphinxAtStartPar
\sphinxstylestrong{Arguments}
\begin{quote}

\sphinxAtStartPar
\sphinxcode{\sphinxupquote{dim}}: dimension of dense symmetric matrix.
\end{quote}

\sphinxAtStartPar
\sphinxstylestrong{Return}
\begin{quote}

\sphinxAtStartPar
new symmetric matrix object.
\end{quote}
\end{quote}

\subsubsection{Model.AddPsdConstr()}
\label{\detokenize{csapi/model:model-addpsdconstr}}\begin{quote}

\sphinxAtStartPar
Add a PSD constraint to model.

\sphinxAtStartPar
\sphinxstylestrong{Synopsis}
\begin{quote}

\sphinxAtStartPar
\sphinxcode{\sphinxupquote{PsdConstraint AddPsdConstr(}}
\begin{quote}

\sphinxAtStartPar
\sphinxcode{\sphinxupquote{PsdExpr expr,}}

\sphinxAtStartPar
\sphinxcode{\sphinxupquote{char sense,}}

\sphinxAtStartPar
\sphinxcode{\sphinxupquote{double rhs,}}

\sphinxAtStartPar
\sphinxcode{\sphinxupquote{string name)}}
\end{quote}
\end{quote}

\sphinxAtStartPar
\sphinxstylestrong{Arguments}
\begin{quote}

\sphinxAtStartPar
\sphinxcode{\sphinxupquote{expr}}: PSD expression for new PSD contraint.

\sphinxAtStartPar
\sphinxcode{\sphinxupquote{sense}}: sense for new PSD constraint.

\sphinxAtStartPar
\sphinxcode{\sphinxupquote{rhs}}: double value at right side of the new PSD constraint.

\sphinxAtStartPar
\sphinxcode{\sphinxupquote{name}}: optional, name of new PSD constraint.
\end{quote}

\sphinxAtStartPar
\sphinxstylestrong{Return}
\begin{quote}

\sphinxAtStartPar
new PSD constraint object.
\end{quote}
\end{quote}

\subsubsection{Model.AddPsdConstr()}
\label{\detokenize{csapi/model:id37}}\begin{quote}

\sphinxAtStartPar
Add a PSD constraint to model.

\sphinxAtStartPar
\sphinxstylestrong{Synopsis}
\begin{quote}

\sphinxAtStartPar
\sphinxcode{\sphinxupquote{PsdConstraint AddPsdConstr(}}
\begin{quote}

\sphinxAtStartPar
\sphinxcode{\sphinxupquote{PsdExpr expr,}}

\sphinxAtStartPar
\sphinxcode{\sphinxupquote{double lb,}}

\sphinxAtStartPar
\sphinxcode{\sphinxupquote{double ub,}}

\sphinxAtStartPar
\sphinxcode{\sphinxupquote{string name)}}
\end{quote}
\end{quote}

\sphinxAtStartPar
\sphinxstylestrong{Arguments}
\begin{quote}

\sphinxAtStartPar
\sphinxcode{\sphinxupquote{expr}}: expression for new PSD constraint.

\sphinxAtStartPar
\sphinxcode{\sphinxupquote{lb}}: lower bound for ew PSD constraint.

\sphinxAtStartPar
\sphinxcode{\sphinxupquote{ub}}: upper bound for new PSD constraint

\sphinxAtStartPar
\sphinxcode{\sphinxupquote{name}}: optional, name of new PSD constraint.
\end{quote}

\sphinxAtStartPar
\sphinxstylestrong{Return}
\begin{quote}

\sphinxAtStartPar
new PSD constraint object.
\end{quote}
\end{quote}

\subsubsection{Model.AddPsdConstr()}
\label{\detokenize{csapi/model:id38}}\begin{quote}

\sphinxAtStartPar
Add a PSD constraint to model.

\sphinxAtStartPar
\sphinxstylestrong{Synopsis}
\begin{quote}

\sphinxAtStartPar
\sphinxcode{\sphinxupquote{PsdConstraint AddPsdConstr(}}
\begin{quote}

\sphinxAtStartPar
\sphinxcode{\sphinxupquote{PsdExpr lhs,}}

\sphinxAtStartPar
\sphinxcode{\sphinxupquote{char sense,}}

\sphinxAtStartPar
\sphinxcode{\sphinxupquote{PsdExpr rhs,}}

\sphinxAtStartPar
\sphinxcode{\sphinxupquote{string name)}}
\end{quote}
\end{quote}

\sphinxAtStartPar
\sphinxstylestrong{Arguments}
\begin{quote}

\sphinxAtStartPar
\sphinxcode{\sphinxupquote{lhs}}: PSD expression at left side of new PSD constraint.

\sphinxAtStartPar
\sphinxcode{\sphinxupquote{sense}}: sense for new PSD constraint.

\sphinxAtStartPar
\sphinxcode{\sphinxupquote{rhs}}: PSD expression at right side of new PSD constraint.

\sphinxAtStartPar
\sphinxcode{\sphinxupquote{name}}: optional, name of new PSD constraint.
\end{quote}

\sphinxAtStartPar
\sphinxstylestrong{Return}
\begin{quote}

\sphinxAtStartPar
new PSD constraint object.
\end{quote}
\end{quote}

\subsubsection{Model.AddPsdConstr()}
\label{\detokenize{csapi/model:id39}}\begin{quote}

\sphinxAtStartPar
Add a PSD constraint to a model.

\sphinxAtStartPar
\sphinxstylestrong{Synopsis}
\begin{quote}

\sphinxAtStartPar
\sphinxcode{\sphinxupquote{PsdConstraint AddPsdConstr(PsdConstrBuilder builder, string name)}}
\end{quote}

\sphinxAtStartPar
\sphinxstylestrong{Arguments}
\begin{quote}

\sphinxAtStartPar
\sphinxcode{\sphinxupquote{builder}}: builder for new PSD constraint.

\sphinxAtStartPar
\sphinxcode{\sphinxupquote{name}}: optional, name of new PSD constraint.
\end{quote}

\sphinxAtStartPar
\sphinxstylestrong{Return}
\begin{quote}

\sphinxAtStartPar
new PSD constraint object.
\end{quote}
\end{quote}

\subsubsection{Model.AddPsdVar()}
\label{\detokenize{csapi/model:model-addpsdvar}}\begin{quote}

\sphinxAtStartPar
Add a new PSD variable to model.

\sphinxAtStartPar
\sphinxstylestrong{Synopsis}
\begin{quote}

\sphinxAtStartPar
\sphinxcode{\sphinxupquote{PsdVar AddPsdVar(int dim, string name)}}
\end{quote}

\sphinxAtStartPar
\sphinxstylestrong{Arguments}
\begin{quote}

\sphinxAtStartPar
\sphinxcode{\sphinxupquote{dim}}: dimension of new PSD variable.

\sphinxAtStartPar
\sphinxcode{\sphinxupquote{name}}: name of new PSD variable.
\end{quote}

\sphinxAtStartPar
\sphinxstylestrong{Return}
\begin{quote}

\sphinxAtStartPar
PSD variable object.
\end{quote}
\end{quote}

\subsubsection{Model.AddPsdVars()}
\label{\detokenize{csapi/model:model-addpsdvars}}\begin{quote}

\sphinxAtStartPar
Add new PSD variables to model.

\sphinxAtStartPar
\sphinxstylestrong{Synopsis}
\begin{quote}

\sphinxAtStartPar
\sphinxcode{\sphinxupquote{PsdVarArray AddPsdVars(}}
\begin{quote}

\sphinxAtStartPar
\sphinxcode{\sphinxupquote{int count,}}

\sphinxAtStartPar
\sphinxcode{\sphinxupquote{int{[}{]} dims,}}

\sphinxAtStartPar
\sphinxcode{\sphinxupquote{string prefix)}}
\end{quote}
\end{quote}

\sphinxAtStartPar
\sphinxstylestrong{Arguments}
\begin{quote}

\sphinxAtStartPar
\sphinxcode{\sphinxupquote{count}}: number of new PSD variables.

\sphinxAtStartPar
\sphinxcode{\sphinxupquote{dims}}: array of dimensions of new PSD variables.

\sphinxAtStartPar
\sphinxcode{\sphinxupquote{prefix}}: name prefix of new PSD variables, default prefix is PSD\_V.
\end{quote}

\sphinxAtStartPar
\sphinxstylestrong{Return}
\begin{quote}

\sphinxAtStartPar
array of PSD variable objects.
\end{quote}
\end{quote}

\subsubsection{Model.AddQConstr()}
\label{\detokenize{csapi/model:model-addqconstr}}\begin{quote}

\sphinxAtStartPar
Add a quadratic constraint to model.

\sphinxAtStartPar
\sphinxstylestrong{Synopsis}
\begin{quote}

\sphinxAtStartPar
\sphinxcode{\sphinxupquote{QConstraint AddQConstr(}}
\begin{quote}

\sphinxAtStartPar
\sphinxcode{\sphinxupquote{QuadExpr expr,}}

\sphinxAtStartPar
\sphinxcode{\sphinxupquote{char sense,}}

\sphinxAtStartPar
\sphinxcode{\sphinxupquote{double rhs,}}

\sphinxAtStartPar
\sphinxcode{\sphinxupquote{string name)}}
\end{quote}
\end{quote}

\sphinxAtStartPar
\sphinxstylestrong{Arguments}
\begin{quote}

\sphinxAtStartPar
\sphinxcode{\sphinxupquote{expr}}: quadratic expression for the new contraint.

\sphinxAtStartPar
\sphinxcode{\sphinxupquote{sense}}: sense for new quadratic constraint.

\sphinxAtStartPar
\sphinxcode{\sphinxupquote{rhs}}: double value at right side of the new quadratic constraint.

\sphinxAtStartPar
\sphinxcode{\sphinxupquote{name}}: optional, name of new quadratic constraint.
\end{quote}

\sphinxAtStartPar
\sphinxstylestrong{Return}
\begin{quote}

\sphinxAtStartPar
new quadratic constraint object.
\end{quote}
\end{quote}

\subsubsection{Model.AddQConstr()}
\label{\detokenize{csapi/model:id40}}\begin{quote}

\sphinxAtStartPar
Add a quadratic constraint to model.

\sphinxAtStartPar
\sphinxstylestrong{Synopsis}
\begin{quote}

\sphinxAtStartPar
\sphinxcode{\sphinxupquote{QConstraint AddQConstr(}}
\begin{quote}

\sphinxAtStartPar
\sphinxcode{\sphinxupquote{QuadExpr lhs,}}

\sphinxAtStartPar
\sphinxcode{\sphinxupquote{char sense,}}

\sphinxAtStartPar
\sphinxcode{\sphinxupquote{QuadExpr rhs,}}

\sphinxAtStartPar
\sphinxcode{\sphinxupquote{string name)}}
\end{quote}
\end{quote}

\sphinxAtStartPar
\sphinxstylestrong{Arguments}
\begin{quote}

\sphinxAtStartPar
\sphinxcode{\sphinxupquote{lhs}}: quadratic expression at left side of the new quadratic constraint.

\sphinxAtStartPar
\sphinxcode{\sphinxupquote{sense}}: sense for new quadratic constraint.

\sphinxAtStartPar
\sphinxcode{\sphinxupquote{rhs}}: quadratic expression at right side of the new quadratic constraint.

\sphinxAtStartPar
\sphinxcode{\sphinxupquote{name}}: optional, name of new quadratic constraint.
\end{quote}

\sphinxAtStartPar
\sphinxstylestrong{Return}
\begin{quote}

\sphinxAtStartPar
new quadratic constraint object.
\end{quote}
\end{quote}

\subsubsection{Model.AddQConstr()}
\label{\detokenize{csapi/model:id41}}\begin{quote}

\sphinxAtStartPar
Add a quadratic constraint to a model.

\sphinxAtStartPar
\sphinxstylestrong{Synopsis}
\begin{quote}

\sphinxAtStartPar
\sphinxcode{\sphinxupquote{QConstraint AddQConstr(QConstrBuilder builder, string name)}}
\end{quote}

\sphinxAtStartPar
\sphinxstylestrong{Arguments}
\begin{quote}

\sphinxAtStartPar
\sphinxcode{\sphinxupquote{builder}}: builder for the new quadratic constraint.

\sphinxAtStartPar
\sphinxcode{\sphinxupquote{name}}: optional, name of new quadratic constraint.
\end{quote}

\sphinxAtStartPar
\sphinxstylestrong{Return}
\begin{quote}

\sphinxAtStartPar
new quadratic constraint object.
\end{quote}
\end{quote}

\subsubsection{Model.AddSos()}
\label{\detokenize{csapi/model:model-addsos}}\begin{quote}

\sphinxAtStartPar
Add a SOS constraint to model.

\sphinxAtStartPar
\sphinxstylestrong{Synopsis}
\begin{quote}

\sphinxAtStartPar
\sphinxcode{\sphinxupquote{Sos AddSos(SosBuilder builder)}}
\end{quote}

\sphinxAtStartPar
\sphinxstylestrong{Arguments}
\begin{quote}

\sphinxAtStartPar
\sphinxcode{\sphinxupquote{builder}}: builder for new SOS constraint.
\end{quote}

\sphinxAtStartPar
\sphinxstylestrong{Return}
\begin{quote}

\sphinxAtStartPar
new SOS constraint object.
\end{quote}
\end{quote}

\subsubsection{Model.AddSos()}
\label{\detokenize{csapi/model:id42}}\begin{quote}

\sphinxAtStartPar
Add a SOS constraint to model.

\sphinxAtStartPar
\sphinxstylestrong{Synopsis}
\begin{quote}

\sphinxAtStartPar
\sphinxcode{\sphinxupquote{Sos AddSos(}}
\begin{quote}

\sphinxAtStartPar
\sphinxcode{\sphinxupquote{Var{[}{]} vars,}}

\sphinxAtStartPar
\sphinxcode{\sphinxupquote{double{[}{]} weights,}}

\sphinxAtStartPar
\sphinxcode{\sphinxupquote{int type)}}
\end{quote}
\end{quote}

\sphinxAtStartPar
\sphinxstylestrong{Arguments}
\begin{quote}

\sphinxAtStartPar
\sphinxcode{\sphinxupquote{vars}}: variables that participate in the SOS constraint.

\sphinxAtStartPar
\sphinxcode{\sphinxupquote{weights}}: weights for variables in the SOS constraint.

\sphinxAtStartPar
\sphinxcode{\sphinxupquote{type}}: type of SOS constraint.
\end{quote}

\sphinxAtStartPar
\sphinxstylestrong{Return}
\begin{quote}

\sphinxAtStartPar
new SOS constraint object.
\end{quote}
\end{quote}

\subsubsection{Model.AddSos()}
\label{\detokenize{csapi/model:id43}}\begin{quote}

\sphinxAtStartPar
Add a SOS constraint to model.

\sphinxAtStartPar
\sphinxstylestrong{Synopsis}
\begin{quote}

\sphinxAtStartPar
\sphinxcode{\sphinxupquote{Sos AddSos(}}
\begin{quote}

\sphinxAtStartPar
\sphinxcode{\sphinxupquote{VarArray vars,}}

\sphinxAtStartPar
\sphinxcode{\sphinxupquote{double{[}{]} weights,}}

\sphinxAtStartPar
\sphinxcode{\sphinxupquote{int type)}}
\end{quote}
\end{quote}

\sphinxAtStartPar
\sphinxstylestrong{Arguments}
\begin{quote}

\sphinxAtStartPar
\sphinxcode{\sphinxupquote{vars}}: variables that participate in the SOS constraint.

\sphinxAtStartPar
\sphinxcode{\sphinxupquote{weights}}: weights for variables in the SOS constraint.

\sphinxAtStartPar
\sphinxcode{\sphinxupquote{type}}: type of SOS constraint.
\end{quote}

\sphinxAtStartPar
\sphinxstylestrong{Return}
\begin{quote}

\sphinxAtStartPar
new SOS constraint object.
\end{quote}
\end{quote}

\subsubsection{Model.AddSparseMat()}
\label{\detokenize{csapi/model:model-addsparsemat}}\begin{quote}

\sphinxAtStartPar
Add a sparse symmetric matrix to a model.

\sphinxAtStartPar
\sphinxstylestrong{Synopsis}
\begin{quote}

\sphinxAtStartPar
\sphinxcode{\sphinxupquote{SymMatrix AddSparseMat(}}
\begin{quote}

\sphinxAtStartPar
\sphinxcode{\sphinxupquote{int dim,}}

\sphinxAtStartPar
\sphinxcode{\sphinxupquote{int nElems,}}

\sphinxAtStartPar
\sphinxcode{\sphinxupquote{int{[}{]} rows,}}

\sphinxAtStartPar
\sphinxcode{\sphinxupquote{int{[}{]} cols,}}

\sphinxAtStartPar
\sphinxcode{\sphinxupquote{double{[}{]} vals)}}
\end{quote}
\end{quote}

\sphinxAtStartPar
\sphinxstylestrong{Arguments}
\begin{quote}

\sphinxAtStartPar
\sphinxcode{\sphinxupquote{dim}}: dimension of the sparse symmetric matrix.

\sphinxAtStartPar
\sphinxcode{\sphinxupquote{nElems}}: number of non\sphinxhyphen{}zero elements in the sparse symmetric matrix.

\sphinxAtStartPar
\sphinxcode{\sphinxupquote{rows}}: array of row indexes of non\sphinxhyphen{}zero elements.

\sphinxAtStartPar
\sphinxcode{\sphinxupquote{cols}}: array of col indexes of non\sphinxhyphen{}zero elements.

\sphinxAtStartPar
\sphinxcode{\sphinxupquote{vals}}: array of values of non\sphinxhyphen{}zero elements.
\end{quote}

\sphinxAtStartPar
\sphinxstylestrong{Return}
\begin{quote}

\sphinxAtStartPar
new symmetric matrix object.
\end{quote}
\end{quote}

\subsubsection{Model.AddSymMat()}
\label{\detokenize{csapi/model:model-addsymmat}}\begin{quote}

\sphinxAtStartPar
Given a symmetric matrix expression, add results matrix to model.

\sphinxAtStartPar
\sphinxstylestrong{Synopsis}
\begin{quote}

\sphinxAtStartPar
\sphinxcode{\sphinxupquote{SymMatrix AddSymMat(SymMatExpr expr)}}
\end{quote}

\sphinxAtStartPar
\sphinxstylestrong{Arguments}
\begin{quote}

\sphinxAtStartPar
\sphinxcode{\sphinxupquote{expr}}: symmetric matrix expression object.
\end{quote}

\sphinxAtStartPar
\sphinxstylestrong{Return}
\begin{quote}

\sphinxAtStartPar
results symmetric matrix object.
\end{quote}
\end{quote}

\subsubsection{Model.AddUserCut()}
\label{\detokenize{csapi/model:model-addusercut}}\begin{quote}

\sphinxAtStartPar
Add a user cut to model.

\sphinxAtStartPar
\sphinxstylestrong{Synopsis}
\begin{quote}

\sphinxAtStartPar
\sphinxcode{\sphinxupquote{void AddUserCut(}}
\begin{quote}

\sphinxAtStartPar
\sphinxcode{\sphinxupquote{Expr lhs,}}

\sphinxAtStartPar
\sphinxcode{\sphinxupquote{char sense,}}

\sphinxAtStartPar
\sphinxcode{\sphinxupquote{double rhs,}}

\sphinxAtStartPar
\sphinxcode{\sphinxupquote{string name)}}
\end{quote}
\end{quote}

\sphinxAtStartPar
\sphinxstylestrong{Arguments}
\begin{quote}

\sphinxAtStartPar
\sphinxcode{\sphinxupquote{lhs}}: expression for user cut.

\sphinxAtStartPar
\sphinxcode{\sphinxupquote{sense}}: sense for user cut.

\sphinxAtStartPar
\sphinxcode{\sphinxupquote{rhs}}: right hand side value for user cut.

\sphinxAtStartPar
\sphinxcode{\sphinxupquote{name}}: optional, name of user cut.
\end{quote}
\end{quote}

\subsubsection{Model.AddUserCut()}
\label{\detokenize{csapi/model:id44}}\begin{quote}

\sphinxAtStartPar
Add a user cut to model.

\sphinxAtStartPar
\sphinxstylestrong{Synopsis}
\begin{quote}

\sphinxAtStartPar
\sphinxcode{\sphinxupquote{void AddUserCut(}}
\begin{quote}

\sphinxAtStartPar
\sphinxcode{\sphinxupquote{Expr lhs,}}

\sphinxAtStartPar
\sphinxcode{\sphinxupquote{char sense,}}

\sphinxAtStartPar
\sphinxcode{\sphinxupquote{Expr rhs,}}

\sphinxAtStartPar
\sphinxcode{\sphinxupquote{string name)}}
\end{quote}
\end{quote}

\sphinxAtStartPar
\sphinxstylestrong{Arguments}
\begin{quote}

\sphinxAtStartPar
\sphinxcode{\sphinxupquote{lhs}}: left hand side expression for user cut.

\sphinxAtStartPar
\sphinxcode{\sphinxupquote{sense}}: sense for user cut.

\sphinxAtStartPar
\sphinxcode{\sphinxupquote{rhs}}: right hand side expression for user cut.

\sphinxAtStartPar
\sphinxcode{\sphinxupquote{name}}: optional, name of user cut.
\end{quote}
\end{quote}

\subsubsection{Model.AddUserCut()}
\label{\detokenize{csapi/model:id45}}\begin{quote}

\sphinxAtStartPar
Add a user cut to model.

\sphinxAtStartPar
\sphinxstylestrong{Synopsis}
\begin{quote}

\sphinxAtStartPar
\sphinxcode{\sphinxupquote{void AddUserCut(ConstrBuilder builder, string name)}}
\end{quote}

\sphinxAtStartPar
\sphinxstylestrong{Arguments}
\begin{quote}

\sphinxAtStartPar
\sphinxcode{\sphinxupquote{builder}}: builder for user cut.

\sphinxAtStartPar
\sphinxcode{\sphinxupquote{name}}: optional, name of user cut.
\end{quote}
\end{quote}

\subsubsection{Model.AddUserCuts()}
\label{\detokenize{csapi/model:model-addusercuts}}\begin{quote}

\sphinxAtStartPar
Add user cuts to model.

\sphinxAtStartPar
\sphinxstylestrong{Synopsis}
\begin{quote}

\sphinxAtStartPar
\sphinxcode{\sphinxupquote{void AddUserCuts(ConstrBuilderArray builders, string prefix)}}
\end{quote}

\sphinxAtStartPar
\sphinxstylestrong{Arguments}
\begin{quote}

\sphinxAtStartPar
\sphinxcode{\sphinxupquote{builders}}: array of builders for user cuts.

\sphinxAtStartPar
\sphinxcode{\sphinxupquote{prefix}}: name prefix of new user cuts.
\end{quote}
\end{quote}

\subsubsection{Model.AddVar()}
\label{\detokenize{csapi/model:model-addvar}}\begin{quote}

\sphinxAtStartPar
Add a variable and the associated non\sphinxhyphen{}zero coefficients as column.

\sphinxAtStartPar
\sphinxstylestrong{Synopsis}
\begin{quote}

\sphinxAtStartPar
\sphinxcode{\sphinxupquote{Var AddVar(}}
\begin{quote}

\sphinxAtStartPar
\sphinxcode{\sphinxupquote{double lb,}}

\sphinxAtStartPar
\sphinxcode{\sphinxupquote{double ub,}}

\sphinxAtStartPar
\sphinxcode{\sphinxupquote{double obj,}}

\sphinxAtStartPar
\sphinxcode{\sphinxupquote{char vtype,}}

\sphinxAtStartPar
\sphinxcode{\sphinxupquote{string name)}}
\end{quote}
\end{quote}

\sphinxAtStartPar
\sphinxstylestrong{Arguments}
\begin{quote}

\sphinxAtStartPar
\sphinxcode{\sphinxupquote{lb}}: lower bound for new variable.

\sphinxAtStartPar
\sphinxcode{\sphinxupquote{ub}}: upper bound for new variable.

\sphinxAtStartPar
\sphinxcode{\sphinxupquote{obj}}: objective coefficient for new variable.

\sphinxAtStartPar
\sphinxcode{\sphinxupquote{vtype}}: variable type for new variable.

\sphinxAtStartPar
\sphinxcode{\sphinxupquote{name}}: optional, name for new variable.
\end{quote}

\sphinxAtStartPar
\sphinxstylestrong{Return}
\begin{quote}

\sphinxAtStartPar
new variable object.
\end{quote}
\end{quote}

\subsubsection{Model.AddVar()}
\label{\detokenize{csapi/model:id46}}\begin{quote}

\sphinxAtStartPar
Add a variable and the associated non\sphinxhyphen{}zero coefficients as column.

\sphinxAtStartPar
\sphinxstylestrong{Synopsis}
\begin{quote}

\sphinxAtStartPar
\sphinxcode{\sphinxupquote{Var AddVar(}}
\begin{quote}

\sphinxAtStartPar
\sphinxcode{\sphinxupquote{double lb,}}

\sphinxAtStartPar
\sphinxcode{\sphinxupquote{double ub,}}

\sphinxAtStartPar
\sphinxcode{\sphinxupquote{double obj,}}

\sphinxAtStartPar
\sphinxcode{\sphinxupquote{char vtype,}}

\sphinxAtStartPar
\sphinxcode{\sphinxupquote{Column col,}}

\sphinxAtStartPar
\sphinxcode{\sphinxupquote{string name)}}
\end{quote}
\end{quote}

\sphinxAtStartPar
\sphinxstylestrong{Arguments}
\begin{quote}

\sphinxAtStartPar
\sphinxcode{\sphinxupquote{lb}}: lower bound for new variable.

\sphinxAtStartPar
\sphinxcode{\sphinxupquote{ub}}: upper bound for new variable.

\sphinxAtStartPar
\sphinxcode{\sphinxupquote{obj}}: objective coefficient for new variable.

\sphinxAtStartPar
\sphinxcode{\sphinxupquote{vtype}}: variable type for new variable.

\sphinxAtStartPar
\sphinxcode{\sphinxupquote{col}}: column object for specifying a set of constraints to which the variable belongs.

\sphinxAtStartPar
\sphinxcode{\sphinxupquote{name}}: optional, name for new variable.
\end{quote}

\sphinxAtStartPar
\sphinxstylestrong{Return}
\begin{quote}

\sphinxAtStartPar
new variable object.
\end{quote}
\end{quote}

\subsubsection{Model.AddVars()}
\label{\detokenize{csapi/model:model-addvars}}\begin{quote}

\sphinxAtStartPar
Add new variables to model.

\sphinxAtStartPar
\sphinxstylestrong{Synopsis}
\begin{quote}

\sphinxAtStartPar
\sphinxcode{\sphinxupquote{VarArray AddVars(}}
\begin{quote}

\sphinxAtStartPar
\sphinxcode{\sphinxupquote{int count,}}

\sphinxAtStartPar
\sphinxcode{\sphinxupquote{char vtype,}}

\sphinxAtStartPar
\sphinxcode{\sphinxupquote{string prefix)}}
\end{quote}
\end{quote}

\sphinxAtStartPar
\sphinxstylestrong{Arguments}
\begin{quote}

\sphinxAtStartPar
\sphinxcode{\sphinxupquote{count}}: the number of variables to add.

\sphinxAtStartPar
\sphinxcode{\sphinxupquote{vtype}}: variable types for new variables.

\sphinxAtStartPar
\sphinxcode{\sphinxupquote{prefix}}: optional, prefix part for names of new variables, default value is ‘C’.
\end{quote}

\sphinxAtStartPar
\sphinxstylestrong{Return}
\begin{quote}

\sphinxAtStartPar
array of new variable objects.
\end{quote}
\end{quote}

\subsubsection{Model.AddVars()}
\label{\detokenize{csapi/model:id47}}\begin{quote}

\sphinxAtStartPar
Add new variables to model.

\sphinxAtStartPar
\sphinxstylestrong{Synopsis}
\begin{quote}

\sphinxAtStartPar
\sphinxcode{\sphinxupquote{VarArray AddVars(}}
\begin{quote}

\sphinxAtStartPar
\sphinxcode{\sphinxupquote{int count,}}

\sphinxAtStartPar
\sphinxcode{\sphinxupquote{double lb,}}

\sphinxAtStartPar
\sphinxcode{\sphinxupquote{double ub,}}

\sphinxAtStartPar
\sphinxcode{\sphinxupquote{double obj,}}

\sphinxAtStartPar
\sphinxcode{\sphinxupquote{char vtype,}}

\sphinxAtStartPar
\sphinxcode{\sphinxupquote{string prefix)}}
\end{quote}
\end{quote}

\sphinxAtStartPar
\sphinxstylestrong{Arguments}
\begin{quote}

\sphinxAtStartPar
\sphinxcode{\sphinxupquote{count}}: the number of variables to add.

\sphinxAtStartPar
\sphinxcode{\sphinxupquote{lb}}: lower bound for new variables.

\sphinxAtStartPar
\sphinxcode{\sphinxupquote{ub}}: upper bound for new variables.

\sphinxAtStartPar
\sphinxcode{\sphinxupquote{obj}}: objective coefficient for new variables.

\sphinxAtStartPar
\sphinxcode{\sphinxupquote{vtype}}: variable type for new variables.

\sphinxAtStartPar
\sphinxcode{\sphinxupquote{prefix}}: optional, prefix part for names of new variables, default value is ‘C’.
\end{quote}

\sphinxAtStartPar
\sphinxstylestrong{Return}
\begin{quote}

\sphinxAtStartPar
array of new variable objects.
\end{quote}
\end{quote}

\subsubsection{Model.AddVars()}
\label{\detokenize{csapi/model:id48}}\begin{quote}

\sphinxAtStartPar
Add new variables to model.

\sphinxAtStartPar
\sphinxstylestrong{Synopsis}
\begin{quote}

\sphinxAtStartPar
\sphinxcode{\sphinxupquote{VarArray AddVars(}}
\begin{quote}

\sphinxAtStartPar
\sphinxcode{\sphinxupquote{int count,}}

\sphinxAtStartPar
\sphinxcode{\sphinxupquote{double{[}{]} lbs,}}

\sphinxAtStartPar
\sphinxcode{\sphinxupquote{double{[}{]} ubs,}}

\sphinxAtStartPar
\sphinxcode{\sphinxupquote{double{[}{]} objs,}}

\sphinxAtStartPar
\sphinxcode{\sphinxupquote{char{[}{]} types,}}

\sphinxAtStartPar
\sphinxcode{\sphinxupquote{string prefix)}}
\end{quote}
\end{quote}

\sphinxAtStartPar
\sphinxstylestrong{Arguments}
\begin{quote}

\sphinxAtStartPar
\sphinxcode{\sphinxupquote{count}}: the number of variables to add.

\sphinxAtStartPar
\sphinxcode{\sphinxupquote{lbs}}: lower bounds for new variables. if NULL, lower bounds are 0.0.

\sphinxAtStartPar
\sphinxcode{\sphinxupquote{ubs}}: upper bounds for new variables. if NULL, upper bounds are infinity or 1 for binary variables.

\sphinxAtStartPar
\sphinxcode{\sphinxupquote{objs}}: objective coefficients for new variables. if NULL, objective coefficients are 0.0.

\sphinxAtStartPar
\sphinxcode{\sphinxupquote{types}}: variable types for new variables. if NULL, variable types are continuous.

\sphinxAtStartPar
\sphinxcode{\sphinxupquote{prefix}}: optional, prefix part for names of new variables, default value is ‘C’.
\end{quote}

\sphinxAtStartPar
\sphinxstylestrong{Return}
\begin{quote}

\sphinxAtStartPar
array of new variable objects.
\end{quote}
\end{quote}

\subsubsection{Model.AddVars()}
\label{\detokenize{csapi/model:id49}}\begin{quote}

\sphinxAtStartPar
Add new variables to model.

\sphinxAtStartPar
\sphinxstylestrong{Synopsis}
\begin{quote}

\sphinxAtStartPar
\sphinxcode{\sphinxupquote{VarArray AddVars(}}
\begin{quote}

\sphinxAtStartPar
\sphinxcode{\sphinxupquote{double{[}{]} lbs,}}

\sphinxAtStartPar
\sphinxcode{\sphinxupquote{double{[}{]} ubs,}}

\sphinxAtStartPar
\sphinxcode{\sphinxupquote{double{[}{]} objs,}}

\sphinxAtStartPar
\sphinxcode{\sphinxupquote{char{[}{]} types,}}

\sphinxAtStartPar
\sphinxcode{\sphinxupquote{Column{[}{]} cols,}}

\sphinxAtStartPar
\sphinxcode{\sphinxupquote{string prefix)}}
\end{quote}
\end{quote}

\sphinxAtStartPar
\sphinxstylestrong{Arguments}
\begin{quote}

\sphinxAtStartPar
\sphinxcode{\sphinxupquote{lbs}}: lower bounds for new variables. if NULL, lower bounds are 0.0.

\sphinxAtStartPar
\sphinxcode{\sphinxupquote{ubs}}: upper bounds for new variables. if NULL, upper bounds are infinity or 1 for binary variables.

\sphinxAtStartPar
\sphinxcode{\sphinxupquote{objs}}: objective coefficients for new variables. if NULL, objective coefficients are 0.0.

\sphinxAtStartPar
\sphinxcode{\sphinxupquote{types}}: variable types for new variables. if NULL, variable types are continuous.

\sphinxAtStartPar
\sphinxcode{\sphinxupquote{cols}}: column objects for specifying a set of constraints to which each new variable belongs.

\sphinxAtStartPar
\sphinxcode{\sphinxupquote{prefix}}: optional, prefix part for names of new variables, default value is ‘C’.
\end{quote}

\sphinxAtStartPar
\sphinxstylestrong{Return}
\begin{quote}

\sphinxAtStartPar
array of new variable objects.
\end{quote}
\end{quote}

\subsubsection{Model.AddVars()}
\label{\detokenize{csapi/model:id50}}\begin{quote}

\sphinxAtStartPar
Add new variables to model.

\sphinxAtStartPar
\sphinxstylestrong{Synopsis}
\begin{quote}

\sphinxAtStartPar
\sphinxcode{\sphinxupquote{VarArray AddVars(}}
\begin{quote}

\sphinxAtStartPar
\sphinxcode{\sphinxupquote{double{[}{]} lbs,}}

\sphinxAtStartPar
\sphinxcode{\sphinxupquote{double{[}{]} ubs,}}

\sphinxAtStartPar
\sphinxcode{\sphinxupquote{double{[}{]} objs,}}

\sphinxAtStartPar
\sphinxcode{\sphinxupquote{char{[}{]} types,}}

\sphinxAtStartPar
\sphinxcode{\sphinxupquote{ColumnArray cols,}}

\sphinxAtStartPar
\sphinxcode{\sphinxupquote{string prefix)}}
\end{quote}
\end{quote}

\sphinxAtStartPar
\sphinxstylestrong{Arguments}
\begin{quote}

\sphinxAtStartPar
\sphinxcode{\sphinxupquote{lbs}}: lower bounds for new variables. if NULL, lower bounds are 0.0.

\sphinxAtStartPar
\sphinxcode{\sphinxupquote{ubs}}: upper bounds for new variables. if NULL, upper bounds are infinity or 1 for binary variables.

\sphinxAtStartPar
\sphinxcode{\sphinxupquote{objs}}: objective coefficients for new variables. if NULL, objective coefficients are 0.0.

\sphinxAtStartPar
\sphinxcode{\sphinxupquote{types}}: variable types for new variables. if NULL, variable types are continuous.

\sphinxAtStartPar
\sphinxcode{\sphinxupquote{cols}}: columnarray for specifying a set of constraints to which each new variable belongs.

\sphinxAtStartPar
\sphinxcode{\sphinxupquote{prefix}}: optional, prefix part for names of new variables, default value is ‘C’.
\end{quote}

\sphinxAtStartPar
\sphinxstylestrong{Return}
\begin{quote}

\sphinxAtStartPar
array of new variable objects.
\end{quote}
\end{quote}

\subsubsection{Model.Clear()}
\label{\detokenize{csapi/model:model-clear}}\begin{quote}

\sphinxAtStartPar
Clear all settings including problem itself.

\sphinxAtStartPar
\sphinxstylestrong{Synopsis}
\begin{quote}

\sphinxAtStartPar
\sphinxcode{\sphinxupquote{void Clear()}}
\end{quote}
\end{quote}

\subsubsection{Model.Clone()}
\label{\detokenize{csapi/model:model-clone}}\begin{quote}

\sphinxAtStartPar
Deep copy COPT model.

\sphinxAtStartPar
\sphinxstylestrong{Synopsis}
\begin{quote}

\sphinxAtStartPar
\sphinxcode{\sphinxupquote{Model Clone()}}
\end{quote}

\sphinxAtStartPar
\sphinxstylestrong{Return}
\begin{quote}

\sphinxAtStartPar
cloned model object.
\end{quote}
\end{quote}

\subsubsection{Model.ComputeIIS()}
\label{\detokenize{csapi/model:model-computeiis}}\begin{quote}

\sphinxAtStartPar
Compute IIS for model

\sphinxAtStartPar
\sphinxstylestrong{Synopsis}
\begin{quote}

\sphinxAtStartPar
\sphinxcode{\sphinxupquote{void ComputeIIS()}}
\end{quote}
\end{quote}

\subsubsection{Model.DelNlObj()}
\label{\detokenize{csapi/model:model-delnlobj}}\begin{quote}

\sphinxAtStartPar
delete nonlinear part of objective in model.

\sphinxAtStartPar
\sphinxstylestrong{Synopsis}
\begin{quote}

\sphinxAtStartPar
\sphinxcode{\sphinxupquote{void DelNlObj()}}
\end{quote}
\end{quote}

\subsubsection{Model.DelObjN()}
\label{\detokenize{csapi/model:model-delobjn}}\begin{quote}

\sphinxAtStartPar
Delete linear part of specific multi\sphinxhyphen{}objective function in model.

\sphinxAtStartPar
\sphinxstylestrong{Synopsis}
\begin{quote}

\sphinxAtStartPar
\sphinxcode{\sphinxupquote{void DelObjN(int idx)}}
\end{quote}

\sphinxAtStartPar
\sphinxstylestrong{Arguments}
\begin{quote}

\sphinxAtStartPar
\sphinxcode{\sphinxupquote{idx}}: index of a multi\sphinxhyphen{}objective function.
\end{quote}
\end{quote}

\subsubsection{Model.DelPsdObj()}
\label{\detokenize{csapi/model:model-delpsdobj}}\begin{quote}

\sphinxAtStartPar
delete PSD part of objective in model.

\sphinxAtStartPar
\sphinxstylestrong{Synopsis}
\begin{quote}

\sphinxAtStartPar
\sphinxcode{\sphinxupquote{void DelPsdObj()}}
\end{quote}
\end{quote}

\subsubsection{Model.DelQuadObj()}
\label{\detokenize{csapi/model:model-delquadobj}}\begin{quote}

\sphinxAtStartPar
delete quadratic part of objective in model.

\sphinxAtStartPar
\sphinxstylestrong{Synopsis}
\begin{quote}

\sphinxAtStartPar
\sphinxcode{\sphinxupquote{void DelQuadObj()}}
\end{quote}
\end{quote}

\subsubsection{Model.FeasRelax()}
\label{\detokenize{csapi/model:model-feasrelax}}\begin{quote}

\sphinxAtStartPar
Compute feasibility relaxation for infeasible model.

\sphinxAtStartPar
\sphinxstylestrong{Synopsis}
\begin{quote}

\sphinxAtStartPar
\sphinxcode{\sphinxupquote{void FeasRelax(}}
\begin{quote}

\sphinxAtStartPar
\sphinxcode{\sphinxupquote{VarArray vars,}}

\sphinxAtStartPar
\sphinxcode{\sphinxupquote{double{[}{]} colLowPen,}}

\sphinxAtStartPar
\sphinxcode{\sphinxupquote{double{[}{]} colUppPen,}}

\sphinxAtStartPar
\sphinxcode{\sphinxupquote{ConstrArray cons,}}

\sphinxAtStartPar
\sphinxcode{\sphinxupquote{double{[}{]} rowBndPen,}}

\sphinxAtStartPar
\sphinxcode{\sphinxupquote{double{[}{]} rowUppPen)}}
\end{quote}
\end{quote}

\sphinxAtStartPar
\sphinxstylestrong{Arguments}
\begin{quote}

\sphinxAtStartPar
\sphinxcode{\sphinxupquote{vars}}: an array of variables.

\sphinxAtStartPar
\sphinxcode{\sphinxupquote{colLowPen}}: penalties for lower bounds of variables.

\sphinxAtStartPar
\sphinxcode{\sphinxupquote{colUppPen}}: penalties for upper bounds of variables.

\sphinxAtStartPar
\sphinxcode{\sphinxupquote{cons}}: an array of constraints.

\sphinxAtStartPar
\sphinxcode{\sphinxupquote{rowBndPen}}: penalties for right hand sides of constraints.

\sphinxAtStartPar
\sphinxcode{\sphinxupquote{rowUppPen}}: penalties for upper bounds of range constraints.
\end{quote}
\end{quote}

\subsubsection{Model.FeasRelax()}
\label{\detokenize{csapi/model:id51}}\begin{quote}

\sphinxAtStartPar
Compute feasibility relaxation for infeasible model.

\sphinxAtStartPar
\sphinxstylestrong{Synopsis}
\begin{quote}

\sphinxAtStartPar
\sphinxcode{\sphinxupquote{void FeasRelax(int ifRelaxVars, int ifRelaxCons)}}
\end{quote}

\sphinxAtStartPar
\sphinxstylestrong{Arguments}
\begin{quote}

\sphinxAtStartPar
\sphinxcode{\sphinxupquote{ifRelaxVars}}: whether to relax variables.

\sphinxAtStartPar
\sphinxcode{\sphinxupquote{ifRelaxCons}}: whether to relax constraints.
\end{quote}
\end{quote}

\subsubsection{Model.Get()}
\label{\detokenize{csapi/model:model-get}}\begin{quote}

\sphinxAtStartPar
Query values of information associated with variables.

\sphinxAtStartPar
\sphinxstylestrong{Synopsis}
\begin{quote}

\sphinxAtStartPar
\sphinxcode{\sphinxupquote{double{[}{]} Get(string name, Var{[}{]} vars)}}
\end{quote}

\sphinxAtStartPar
\sphinxstylestrong{Arguments}
\begin{quote}

\sphinxAtStartPar
\sphinxcode{\sphinxupquote{name}}: name of information.

\sphinxAtStartPar
\sphinxcode{\sphinxupquote{vars}}: a list of interested variables.
\end{quote}

\sphinxAtStartPar
\sphinxstylestrong{Return}
\begin{quote}

\sphinxAtStartPar
values of information.
\end{quote}
\end{quote}

\subsubsection{Model.Get()}
\label{\detokenize{csapi/model:id52}}\begin{quote}

\sphinxAtStartPar
Query values of information associated with variables.

\sphinxAtStartPar
\sphinxstylestrong{Synopsis}
\begin{quote}

\sphinxAtStartPar
\sphinxcode{\sphinxupquote{double{[}{]} Get(string name, VarArray vars)}}
\end{quote}

\sphinxAtStartPar
\sphinxstylestrong{Arguments}
\begin{quote}

\sphinxAtStartPar
\sphinxcode{\sphinxupquote{name}}: name of information.

\sphinxAtStartPar
\sphinxcode{\sphinxupquote{vars}}: array of interested variables.
\end{quote}

\sphinxAtStartPar
\sphinxstylestrong{Return}
\begin{quote}

\sphinxAtStartPar
values of information.
\end{quote}
\end{quote}

\subsubsection{Model.Get()}
\label{\detokenize{csapi/model:id53}}\begin{quote}

\sphinxAtStartPar
Query values of information associated with constraints.

\sphinxAtStartPar
\sphinxstylestrong{Synopsis}
\begin{quote}

\sphinxAtStartPar
\sphinxcode{\sphinxupquote{double{[}{]} Get(string name, Constraint{[}{]} constrs)}}
\end{quote}

\sphinxAtStartPar
\sphinxstylestrong{Arguments}
\begin{quote}

\sphinxAtStartPar
\sphinxcode{\sphinxupquote{name}}: name of information.

\sphinxAtStartPar
\sphinxcode{\sphinxupquote{constrs}}: a list of interested constraints.
\end{quote}

\sphinxAtStartPar
\sphinxstylestrong{Return}
\begin{quote}

\sphinxAtStartPar
values of information.
\end{quote}
\end{quote}

\subsubsection{Model.Get()}
\label{\detokenize{csapi/model:id54}}\begin{quote}

\sphinxAtStartPar
Query values of information associated with constraints.

\sphinxAtStartPar
\sphinxstylestrong{Synopsis}
\begin{quote}

\sphinxAtStartPar
\sphinxcode{\sphinxupquote{double{[}{]} Get(string name, ConstrArray constrs)}}
\end{quote}

\sphinxAtStartPar
\sphinxstylestrong{Arguments}
\begin{quote}

\sphinxAtStartPar
\sphinxcode{\sphinxupquote{name}}: name of information.

\sphinxAtStartPar
\sphinxcode{\sphinxupquote{constrs}}: array of interested constraints.
\end{quote}

\sphinxAtStartPar
\sphinxstylestrong{Return}
\begin{quote}

\sphinxAtStartPar
values of information.
\end{quote}
\end{quote}

\subsubsection{Model.Get()}
\label{\detokenize{csapi/model:id55}}\begin{quote}

\sphinxAtStartPar
Query values of information associated with nonlinear constraints.

\sphinxAtStartPar
\sphinxstylestrong{Synopsis}
\begin{quote}

\sphinxAtStartPar
\sphinxcode{\sphinxupquote{double{[}{]} Get(string name, NlConstraint{[}{]} constrs)}}
\end{quote}

\sphinxAtStartPar
\sphinxstylestrong{Arguments}
\begin{quote}

\sphinxAtStartPar
\sphinxcode{\sphinxupquote{name}}: name of information.

\sphinxAtStartPar
\sphinxcode{\sphinxupquote{constrs}}: array of desired nonlinear constraints.
\end{quote}

\sphinxAtStartPar
\sphinxstylestrong{Return}
\begin{quote}

\sphinxAtStartPar
output array of information values.
\end{quote}
\end{quote}

\subsubsection{Model.Get()}
\label{\detokenize{csapi/model:id56}}\begin{quote}

\sphinxAtStartPar
Query values of information associated with nonlinear constraints.

\sphinxAtStartPar
\sphinxstylestrong{Synopsis}
\begin{quote}

\sphinxAtStartPar
\sphinxcode{\sphinxupquote{double{[}{]} Get(string name, NlConstrArray constrs)}}
\end{quote}

\sphinxAtStartPar
\sphinxstylestrong{Arguments}
\begin{quote}

\sphinxAtStartPar
\sphinxcode{\sphinxupquote{name}}: name of information.

\sphinxAtStartPar
\sphinxcode{\sphinxupquote{constrs}}: a list of desired nonlinear constraints.
\end{quote}

\sphinxAtStartPar
\sphinxstylestrong{Return}
\begin{quote}

\sphinxAtStartPar
output array of information values.
\end{quote}
\end{quote}

\subsubsection{Model.Get()}
\label{\detokenize{csapi/model:id57}}\begin{quote}

\sphinxAtStartPar
Query values of information associated with quadratic constraints.

\sphinxAtStartPar
\sphinxstylestrong{Synopsis}
\begin{quote}

\sphinxAtStartPar
\sphinxcode{\sphinxupquote{double{[}{]} Get(string name, QConstraint{[}{]} constrs)}}
\end{quote}

\sphinxAtStartPar
\sphinxstylestrong{Arguments}
\begin{quote}

\sphinxAtStartPar
\sphinxcode{\sphinxupquote{name}}: name of information.

\sphinxAtStartPar
\sphinxcode{\sphinxupquote{constrs}}: a list of interested quadratic constraints.
\end{quote}

\sphinxAtStartPar
\sphinxstylestrong{Return}
\begin{quote}

\sphinxAtStartPar
values of information.
\end{quote}
\end{quote}

\subsubsection{Model.Get()}
\label{\detokenize{csapi/model:id58}}\begin{quote}

\sphinxAtStartPar
Query values of information associated with quadratic constraints.

\sphinxAtStartPar
\sphinxstylestrong{Synopsis}
\begin{quote}

\sphinxAtStartPar
\sphinxcode{\sphinxupquote{double{[}{]} Get(string name, QConstrArray constrs)}}
\end{quote}

\sphinxAtStartPar
\sphinxstylestrong{Arguments}
\begin{quote}

\sphinxAtStartPar
\sphinxcode{\sphinxupquote{name}}: name of information.

\sphinxAtStartPar
\sphinxcode{\sphinxupquote{constrs}}: array of interested quadratic constraints.
\end{quote}

\sphinxAtStartPar
\sphinxstylestrong{Return}
\begin{quote}

\sphinxAtStartPar
values of information.
\end{quote}
\end{quote}

\subsubsection{Model.Get()}
\label{\detokenize{csapi/model:id59}}\begin{quote}

\sphinxAtStartPar
Query values of information associated with PSD constraints.

\sphinxAtStartPar
\sphinxstylestrong{Synopsis}
\begin{quote}

\sphinxAtStartPar
\sphinxcode{\sphinxupquote{double{[}{]} Get(string name, PsdConstraint{[}{]} constrs)}}
\end{quote}

\sphinxAtStartPar
\sphinxstylestrong{Arguments}
\begin{quote}

\sphinxAtStartPar
\sphinxcode{\sphinxupquote{name}}: name of information.

\sphinxAtStartPar
\sphinxcode{\sphinxupquote{constrs}}: a list of desired PSD constraints.
\end{quote}

\sphinxAtStartPar
\sphinxstylestrong{Return}
\begin{quote}

\sphinxAtStartPar
output array of information values.
\end{quote}
\end{quote}

\subsubsection{Model.Get()}
\label{\detokenize{csapi/model:id60}}\begin{quote}

\sphinxAtStartPar
Query values of information associated with PSD constraints.

\sphinxAtStartPar
\sphinxstylestrong{Synopsis}
\begin{quote}

\sphinxAtStartPar
\sphinxcode{\sphinxupquote{double{[}{]} Get(string name, PsdConstrArray constrs)}}
\end{quote}

\sphinxAtStartPar
\sphinxstylestrong{Arguments}
\begin{quote}

\sphinxAtStartPar
\sphinxcode{\sphinxupquote{name}}: name of information.

\sphinxAtStartPar
\sphinxcode{\sphinxupquote{constrs}}: a list of desired PSD constraints.
\end{quote}

\sphinxAtStartPar
\sphinxstylestrong{Return}
\begin{quote}

\sphinxAtStartPar
output array of information values.
\end{quote}
\end{quote}

\subsubsection{Model.GetAffineCone()}
\label{\detokenize{csapi/model:model-getaffinecone}}\begin{quote}

\sphinxAtStartPar
Get an affine cone constraint of given index in model.

\sphinxAtStartPar
\sphinxstylestrong{Synopsis}
\begin{quote}

\sphinxAtStartPar
\sphinxcode{\sphinxupquote{AffineCone GetAffineCone(int idx)}}
\end{quote}

\sphinxAtStartPar
\sphinxstylestrong{Arguments}
\begin{quote}

\sphinxAtStartPar
\sphinxcode{\sphinxupquote{idx}}: index of the desired affine cone constraint.
\end{quote}

\sphinxAtStartPar
\sphinxstylestrong{Return}
\begin{quote}

\sphinxAtStartPar
the desired affine cone constraint object.
\end{quote}
\end{quote}

\subsubsection{Model.GetAffineConeBuilder()}
\label{\detokenize{csapi/model:model-getaffineconebuilder}}\begin{quote}

\sphinxAtStartPar
Get builder of given affine cone constraint in model.

\sphinxAtStartPar
\sphinxstylestrong{Synopsis}
\begin{quote}

\sphinxAtStartPar
\sphinxcode{\sphinxupquote{AffineConeBuilder GetAffineConeBuilder(AffineCone cone)}}
\end{quote}

\sphinxAtStartPar
\sphinxstylestrong{Arguments}
\begin{quote}

\sphinxAtStartPar
\sphinxcode{\sphinxupquote{cone}}: affine cone constraint.
\end{quote}

\sphinxAtStartPar
\sphinxstylestrong{Return}
\begin{quote}

\sphinxAtStartPar
desired affine cone constraint builder.
\end{quote}
\end{quote}

\subsubsection{Model.GetAffineConeBuilders()}
\label{\detokenize{csapi/model:model-getaffineconebuilders}}\begin{quote}

\sphinxAtStartPar
Get builders of all affine cone constraints in model.

\sphinxAtStartPar
\sphinxstylestrong{Synopsis}
\begin{quote}

\sphinxAtStartPar
\sphinxcode{\sphinxupquote{AffineConeBuilderArray GetAffineConeBuilders()}}
\end{quote}

\sphinxAtStartPar
\sphinxstylestrong{Return}
\begin{quote}

\sphinxAtStartPar
array object of affine cone constraint builders.
\end{quote}
\end{quote}

\subsubsection{Model.GetAffineConeBuilders()}
\label{\detokenize{csapi/model:id61}}\begin{quote}

\sphinxAtStartPar
Get builders of given affine cone constraints in model.

\sphinxAtStartPar
\sphinxstylestrong{Synopsis}
\begin{quote}

\sphinxAtStartPar
\sphinxcode{\sphinxupquote{AffineConeBuilderArray GetAffineConeBuilders(AffineCone{[}{]} cones)}}
\end{quote}

\sphinxAtStartPar
\sphinxstylestrong{Arguments}
\begin{quote}

\sphinxAtStartPar
\sphinxcode{\sphinxupquote{cones}}: array of affine cone constraints.
\end{quote}

\sphinxAtStartPar
\sphinxstylestrong{Return}
\begin{quote}

\sphinxAtStartPar
array object of desired affine cone constraint builders.
\end{quote}
\end{quote}

\subsubsection{Model.GetAffineConeBuilders()}
\label{\detokenize{csapi/model:id62}}\begin{quote}

\sphinxAtStartPar
Get builders of given affine cone constraints in model.

\sphinxAtStartPar
\sphinxstylestrong{Synopsis}
\begin{quote}

\sphinxAtStartPar
\sphinxcode{\sphinxupquote{AffineConeBuilderArray GetAffineConeBuilders(AffineConeArray cones)}}
\end{quote}

\sphinxAtStartPar
\sphinxstylestrong{Arguments}
\begin{quote}

\sphinxAtStartPar
\sphinxcode{\sphinxupquote{cones}}: array of affine cone constraints.
\end{quote}

\sphinxAtStartPar
\sphinxstylestrong{Return}
\begin{quote}

\sphinxAtStartPar
array object of desired affine cone constraint builders.
\end{quote}
\end{quote}

\subsubsection{Model.GetAffineConeByName()}
\label{\detokenize{csapi/model:model-getaffineconebyname}}\begin{quote}

\sphinxAtStartPar
Get an affine cone constraint of given name in model.

\sphinxAtStartPar
\sphinxstylestrong{Synopsis}
\begin{quote}

\sphinxAtStartPar
\sphinxcode{\sphinxupquote{AffineCone GetAffineConeByName(string name)}}
\end{quote}

\sphinxAtStartPar
\sphinxstylestrong{Arguments}
\begin{quote}

\sphinxAtStartPar
\sphinxcode{\sphinxupquote{name}}: name of the desired affine cone constraint.
\end{quote}

\sphinxAtStartPar
\sphinxstylestrong{Return}
\begin{quote}

\sphinxAtStartPar
the desired affine cone constraint object.
\end{quote}
\end{quote}

\subsubsection{Model.GetAffineCones()}
\label{\detokenize{csapi/model:model-getaffinecones}}\begin{quote}

\sphinxAtStartPar
Get all affine cone constraints in model.

\sphinxAtStartPar
\sphinxstylestrong{Synopsis}
\begin{quote}

\sphinxAtStartPar
\sphinxcode{\sphinxupquote{AffineConeArray GetAffineCones()}}
\end{quote}

\sphinxAtStartPar
\sphinxstylestrong{Return}
\begin{quote}

\sphinxAtStartPar
array object of affine cone constraints.
\end{quote}
\end{quote}

\subsubsection{Model.GetCoeff()}
\label{\detokenize{csapi/model:model-getcoeff}}\begin{quote}

\sphinxAtStartPar
Get the coefficient of variable in linear constraint.

\sphinxAtStartPar
\sphinxstylestrong{Synopsis}
\begin{quote}

\sphinxAtStartPar
\sphinxcode{\sphinxupquote{double GetCoeff(Constraint constr, Var var)}}
\end{quote}

\sphinxAtStartPar
\sphinxstylestrong{Arguments}
\begin{quote}

\sphinxAtStartPar
\sphinxcode{\sphinxupquote{constr}}: The requested constraint.

\sphinxAtStartPar
\sphinxcode{\sphinxupquote{var}}: The requested variable.
\end{quote}

\sphinxAtStartPar
\sphinxstylestrong{Return}
\begin{quote}

\sphinxAtStartPar
The requested coefficient.
\end{quote}
\end{quote}

\subsubsection{Model.GetCol()}
\label{\detokenize{csapi/model:model-getcol}}\begin{quote}

\sphinxAtStartPar
Get a column object that have a list of constraints in which the variable participates.

\sphinxAtStartPar
\sphinxstylestrong{Synopsis}
\begin{quote}

\sphinxAtStartPar
\sphinxcode{\sphinxupquote{Column GetCol(Var var)}}
\end{quote}

\sphinxAtStartPar
\sphinxstylestrong{Arguments}
\begin{quote}

\sphinxAtStartPar
\sphinxcode{\sphinxupquote{var}}: a variable object.
\end{quote}

\sphinxAtStartPar
\sphinxstylestrong{Return}
\begin{quote}

\sphinxAtStartPar
a column object associated with a variable.
\end{quote}
\end{quote}

\subsubsection{Model.GetColBasis()}
\label{\detokenize{csapi/model:model-getcolbasis}}\begin{quote}

\sphinxAtStartPar
Get status of column basis.

\sphinxAtStartPar
\sphinxstylestrong{Synopsis}
\begin{quote}

\sphinxAtStartPar
\sphinxcode{\sphinxupquote{int{[}{]} GetColBasis()}}
\end{quote}

\sphinxAtStartPar
\sphinxstylestrong{Return}
\begin{quote}

\sphinxAtStartPar
basis status.
\end{quote}
\end{quote}

\subsubsection{Model.GetCone()}
\label{\detokenize{csapi/model:model-getcone}}\begin{quote}

\sphinxAtStartPar
Get a cone constraint of given index in model.

\sphinxAtStartPar
\sphinxstylestrong{Synopsis}
\begin{quote}

\sphinxAtStartPar
\sphinxcode{\sphinxupquote{Cone GetCone(int idx)}}
\end{quote}

\sphinxAtStartPar
\sphinxstylestrong{Arguments}
\begin{quote}

\sphinxAtStartPar
\sphinxcode{\sphinxupquote{idx}}: index of the desired cone constraint.
\end{quote}

\sphinxAtStartPar
\sphinxstylestrong{Return}
\begin{quote}

\sphinxAtStartPar
the desired cone constraint object.
\end{quote}
\end{quote}

\subsubsection{Model.GetConeBuilders()}
\label{\detokenize{csapi/model:model-getconebuilders}}\begin{quote}

\sphinxAtStartPar
Get builders of all cone constraints in model.

\sphinxAtStartPar
\sphinxstylestrong{Synopsis}
\begin{quote}

\sphinxAtStartPar
\sphinxcode{\sphinxupquote{ConeBuilderArray GetConeBuilders()}}
\end{quote}

\sphinxAtStartPar
\sphinxstylestrong{Return}
\begin{quote}

\sphinxAtStartPar
array object of cone constraint builders.
\end{quote}
\end{quote}

\subsubsection{Model.GetConeBuilders()}
\label{\detokenize{csapi/model:id63}}\begin{quote}

\sphinxAtStartPar
Get builders of given cone constraints in model.

\sphinxAtStartPar
\sphinxstylestrong{Synopsis}
\begin{quote}

\sphinxAtStartPar
\sphinxcode{\sphinxupquote{ConeBuilderArray GetConeBuilders(Cone{[}{]} cones)}}
\end{quote}

\sphinxAtStartPar
\sphinxstylestrong{Arguments}
\begin{quote}

\sphinxAtStartPar
\sphinxcode{\sphinxupquote{cones}}: array of cone constraints.
\end{quote}

\sphinxAtStartPar
\sphinxstylestrong{Return}
\begin{quote}

\sphinxAtStartPar
array object of desired cone constraint builders.
\end{quote}
\end{quote}

\subsubsection{Model.GetConeBuilders()}
\label{\detokenize{csapi/model:id64}}\begin{quote}

\sphinxAtStartPar
Get builders of given cone constraints in model.

\sphinxAtStartPar
\sphinxstylestrong{Synopsis}
\begin{quote}

\sphinxAtStartPar
\sphinxcode{\sphinxupquote{ConeBuilderArray GetConeBuilders(ConeArray cones)}}
\end{quote}

\sphinxAtStartPar
\sphinxstylestrong{Arguments}
\begin{quote}

\sphinxAtStartPar
\sphinxcode{\sphinxupquote{cones}}: array of cone constraints.
\end{quote}

\sphinxAtStartPar
\sphinxstylestrong{Return}
\begin{quote}

\sphinxAtStartPar
array object of desired cone constraint builders.
\end{quote}
\end{quote}

\subsubsection{Model.GetCones()}
\label{\detokenize{csapi/model:model-getcones}}\begin{quote}

\sphinxAtStartPar
Get all cone constraints in model.

\sphinxAtStartPar
\sphinxstylestrong{Synopsis}
\begin{quote}

\sphinxAtStartPar
\sphinxcode{\sphinxupquote{ConeArray GetCones()}}
\end{quote}

\sphinxAtStartPar
\sphinxstylestrong{Return}
\begin{quote}

\sphinxAtStartPar
array object of cone constraints.
\end{quote}
\end{quote}

\subsubsection{Model.GetConstr()}
\label{\detokenize{csapi/model:model-getconstr}}\begin{quote}

\sphinxAtStartPar
Get a constraint of given index in model.

\sphinxAtStartPar
\sphinxstylestrong{Synopsis}
\begin{quote}

\sphinxAtStartPar
\sphinxcode{\sphinxupquote{Constraint GetConstr(int idx)}}
\end{quote}

\sphinxAtStartPar
\sphinxstylestrong{Arguments}
\begin{quote}

\sphinxAtStartPar
\sphinxcode{\sphinxupquote{idx}}: index of the desired constraint.
\end{quote}

\sphinxAtStartPar
\sphinxstylestrong{Return}
\begin{quote}

\sphinxAtStartPar
the desired constraint object.
\end{quote}
\end{quote}

\subsubsection{Model.GetConstrBuilder()}
\label{\detokenize{csapi/model:model-getconstrbuilder}}\begin{quote}

\sphinxAtStartPar
Get builder of a constraint in model, including variables and associated coefficients, sense and RHS.

\sphinxAtStartPar
\sphinxstylestrong{Synopsis}
\begin{quote}

\sphinxAtStartPar
\sphinxcode{\sphinxupquote{ConstrBuilder GetConstrBuilder(Constraint constr)}}
\end{quote}

\sphinxAtStartPar
\sphinxstylestrong{Arguments}
\begin{quote}

\sphinxAtStartPar
\sphinxcode{\sphinxupquote{constr}}: a constraint object.
\end{quote}

\sphinxAtStartPar
\sphinxstylestrong{Return}
\begin{quote}

\sphinxAtStartPar
constraint builder object.
\end{quote}
\end{quote}

\subsubsection{Model.GetConstrBuilders()}
\label{\detokenize{csapi/model:model-getconstrbuilders}}\begin{quote}

\sphinxAtStartPar
Get builders of all constraints in model.

\sphinxAtStartPar
\sphinxstylestrong{Synopsis}
\begin{quote}

\sphinxAtStartPar
\sphinxcode{\sphinxupquote{ConstrBuilderArray GetConstrBuilders()}}
\end{quote}

\sphinxAtStartPar
\sphinxstylestrong{Return}
\begin{quote}

\sphinxAtStartPar
array object of constraint builders.
\end{quote}
\end{quote}

\subsubsection{Model.GetConstrByName()}
\label{\detokenize{csapi/model:model-getconstrbyname}}\begin{quote}

\sphinxAtStartPar
Get a constraint of given name in model.

\sphinxAtStartPar
\sphinxstylestrong{Synopsis}
\begin{quote}

\sphinxAtStartPar
\sphinxcode{\sphinxupquote{Constraint GetConstrByName(string name)}}
\end{quote}

\sphinxAtStartPar
\sphinxstylestrong{Arguments}
\begin{quote}

\sphinxAtStartPar
\sphinxcode{\sphinxupquote{name}}: name of the desired constraint.
\end{quote}

\sphinxAtStartPar
\sphinxstylestrong{Return}
\begin{quote}

\sphinxAtStartPar
the desired constraint object.
\end{quote}
\end{quote}

\subsubsection{Model.GetConstrLowerIIS()}
\label{\detokenize{csapi/model:model-getconstrloweriis}}\begin{quote}

\sphinxAtStartPar
Get IIS status of lower bounds of constraints.

\sphinxAtStartPar
\sphinxstylestrong{Synopsis}
\begin{quote}

\sphinxAtStartPar
\sphinxcode{\sphinxupquote{int{[}{]} GetConstrLowerIIS(ConstrArray constrs)}}
\end{quote}

\sphinxAtStartPar
\sphinxstylestrong{Arguments}
\begin{quote}

\sphinxAtStartPar
\sphinxcode{\sphinxupquote{constrs}}: Array of constraints.
\end{quote}

\sphinxAtStartPar
\sphinxstylestrong{Return}
\begin{quote}

\sphinxAtStartPar
IIS status of lower bounds of constraints.
\end{quote}
\end{quote}

\subsubsection{Model.GetConstrLowerIIS()}
\label{\detokenize{csapi/model:id65}}\begin{quote}

\sphinxAtStartPar
Get IIS status of lower bounds of constraints.

\sphinxAtStartPar
\sphinxstylestrong{Synopsis}
\begin{quote}

\sphinxAtStartPar
\sphinxcode{\sphinxupquote{int{[}{]} GetConstrLowerIIS(Constraint{[}{]} constrs)}}
\end{quote}

\sphinxAtStartPar
\sphinxstylestrong{Arguments}
\begin{quote}

\sphinxAtStartPar
\sphinxcode{\sphinxupquote{constrs}}: Array of constraints.
\end{quote}

\sphinxAtStartPar
\sphinxstylestrong{Return}
\begin{quote}

\sphinxAtStartPar
IIS status of lower bounds of constraints.
\end{quote}
\end{quote}

\subsubsection{Model.GetConstrs()}
\label{\detokenize{csapi/model:model-getconstrs}}\begin{quote}

\sphinxAtStartPar
Get all constraints in model.

\sphinxAtStartPar
\sphinxstylestrong{Synopsis}
\begin{quote}

\sphinxAtStartPar
\sphinxcode{\sphinxupquote{ConstrArray GetConstrs()}}
\end{quote}

\sphinxAtStartPar
\sphinxstylestrong{Return}
\begin{quote}

\sphinxAtStartPar
array object of constraints.
\end{quote}
\end{quote}

\subsubsection{Model.GetConstrUpperIIS()}
\label{\detokenize{csapi/model:model-getconstrupperiis}}\begin{quote}

\sphinxAtStartPar
Get IIS status of upper bounds of constraints.

\sphinxAtStartPar
\sphinxstylestrong{Synopsis}
\begin{quote}

\sphinxAtStartPar
\sphinxcode{\sphinxupquote{int{[}{]} GetConstrUpperIIS(ConstrArray constrs)}}
\end{quote}

\sphinxAtStartPar
\sphinxstylestrong{Arguments}
\begin{quote}

\sphinxAtStartPar
\sphinxcode{\sphinxupquote{constrs}}: Array of constraints.
\end{quote}

\sphinxAtStartPar
\sphinxstylestrong{Return}
\begin{quote}

\sphinxAtStartPar
IIS status of upper bounds of constraints.
\end{quote}
\end{quote}

\subsubsection{Model.GetConstrUpperIIS()}
\label{\detokenize{csapi/model:id66}}\begin{quote}

\sphinxAtStartPar
Get IIS status of upper bounds of constraints.

\sphinxAtStartPar
\sphinxstylestrong{Synopsis}
\begin{quote}

\sphinxAtStartPar
\sphinxcode{\sphinxupquote{int{[}{]} GetConstrUpperIIS(Constraint{[}{]} constrs)}}
\end{quote}

\sphinxAtStartPar
\sphinxstylestrong{Arguments}
\begin{quote}

\sphinxAtStartPar
\sphinxcode{\sphinxupquote{constrs}}: Array of constraints.
\end{quote}

\sphinxAtStartPar
\sphinxstylestrong{Return}
\begin{quote}

\sphinxAtStartPar
IIS status of upper bounds of constraints.
\end{quote}
\end{quote}

\subsubsection{Model.GetDblAttr()}
\label{\detokenize{csapi/model:model-getdblattr}}\begin{quote}

\sphinxAtStartPar
Get value of a COPT double attribute.

\sphinxAtStartPar
\sphinxstylestrong{Synopsis}
\begin{quote}

\sphinxAtStartPar
\sphinxcode{\sphinxupquote{double GetDblAttr(string attr)}}
\end{quote}

\sphinxAtStartPar
\sphinxstylestrong{Arguments}
\begin{quote}

\sphinxAtStartPar
\sphinxcode{\sphinxupquote{attr}}: name of double attribute.
\end{quote}

\sphinxAtStartPar
\sphinxstylestrong{Return}
\begin{quote}

\sphinxAtStartPar
value of double attribute.
\end{quote}
\end{quote}

\subsubsection{Model.GetDblAttrN()}
\label{\detokenize{csapi/model:model-getdblattrn}}\begin{quote}

\sphinxAtStartPar
Get value of a double attribute of a multi\sphinxhyphen{}objective function.

\sphinxAtStartPar
\sphinxstylestrong{Synopsis}
\begin{quote}

\sphinxAtStartPar
\sphinxcode{\sphinxupquote{double GetDblAttrN(int idx, string attr)}}
\end{quote}

\sphinxAtStartPar
\sphinxstylestrong{Arguments}
\begin{quote}

\sphinxAtStartPar
\sphinxcode{\sphinxupquote{idx}}: index of a multi\sphinxhyphen{}objective function.

\sphinxAtStartPar
\sphinxcode{\sphinxupquote{attr}}: name of double attribute.
\end{quote}

\sphinxAtStartPar
\sphinxstylestrong{Return}
\begin{quote}

\sphinxAtStartPar
value of double attribute.
\end{quote}
\end{quote}

\subsubsection{Model.GetDblParam()}
\label{\detokenize{csapi/model:model-getdblparam}}\begin{quote}

\sphinxAtStartPar
Get value of a COPT double parameter.

\sphinxAtStartPar
\sphinxstylestrong{Synopsis}
\begin{quote}

\sphinxAtStartPar
\sphinxcode{\sphinxupquote{double GetDblParam(string param)}}
\end{quote}

\sphinxAtStartPar
\sphinxstylestrong{Arguments}
\begin{quote}

\sphinxAtStartPar
\sphinxcode{\sphinxupquote{param}}: name of integer parameter.
\end{quote}

\sphinxAtStartPar
\sphinxstylestrong{Return}
\begin{quote}

\sphinxAtStartPar
value of double parameter.
\end{quote}
\end{quote}

\subsubsection{Model.GetDblParamN()}
\label{\detokenize{csapi/model:model-getdblparamn}}\begin{quote}

\sphinxAtStartPar
Get value of a double parameter of a multi\sphinxhyphen{}objective function.

\sphinxAtStartPar
\sphinxstylestrong{Synopsis}
\begin{quote}

\sphinxAtStartPar
\sphinxcode{\sphinxupquote{double GetDblParamN(int idx, string param)}}
\end{quote}

\sphinxAtStartPar
\sphinxstylestrong{Arguments}
\begin{quote}

\sphinxAtStartPar
\sphinxcode{\sphinxupquote{idx}}: index of a multi\sphinxhyphen{}objective function.

\sphinxAtStartPar
\sphinxcode{\sphinxupquote{param}}: name of double parameter.
\end{quote}

\sphinxAtStartPar
\sphinxstylestrong{Return}
\begin{quote}

\sphinxAtStartPar
value of double parameter.
\end{quote}
\end{quote}

\subsubsection{Model.GetExpCone()}
\label{\detokenize{csapi/model:model-getexpcone}}\begin{quote}

\sphinxAtStartPar
Get an exponential cone constraint of given index in model.

\sphinxAtStartPar
\sphinxstylestrong{Synopsis}
\begin{quote}

\sphinxAtStartPar
\sphinxcode{\sphinxupquote{ExpCone GetExpCone(int idx)}}
\end{quote}

\sphinxAtStartPar
\sphinxstylestrong{Arguments}
\begin{quote}

\sphinxAtStartPar
\sphinxcode{\sphinxupquote{idx}}: index of the desired exponential cone constraint.
\end{quote}

\sphinxAtStartPar
\sphinxstylestrong{Return}
\begin{quote}

\sphinxAtStartPar
the desired exponential cone constraint object.
\end{quote}
\end{quote}

\subsubsection{Model.GetExpConeBuilders()}
\label{\detokenize{csapi/model:model-getexpconebuilders}}\begin{quote}

\sphinxAtStartPar
Get builders of all exponential cone constraints in model.

\sphinxAtStartPar
\sphinxstylestrong{Synopsis}
\begin{quote}

\sphinxAtStartPar
\sphinxcode{\sphinxupquote{ExpConeBuilderArray GetExpConeBuilders()}}
\end{quote}

\sphinxAtStartPar
\sphinxstylestrong{Return}
\begin{quote}

\sphinxAtStartPar
array object of exponential cone constraint builders.
\end{quote}
\end{quote}

\subsubsection{Model.GetExpConeBuilders()}
\label{\detokenize{csapi/model:id67}}\begin{quote}

\sphinxAtStartPar
Get builders of given exponential cone constraints in model.

\sphinxAtStartPar
\sphinxstylestrong{Synopsis}
\begin{quote}

\sphinxAtStartPar
\sphinxcode{\sphinxupquote{ExpConeBuilderArray GetExpConeBuilders(ExpCone{[}{]} cones)}}
\end{quote}

\sphinxAtStartPar
\sphinxstylestrong{Arguments}
\begin{quote}

\sphinxAtStartPar
\sphinxcode{\sphinxupquote{cones}}: array of exponential cone constraints.
\end{quote}

\sphinxAtStartPar
\sphinxstylestrong{Return}
\begin{quote}

\sphinxAtStartPar
array object of desired exponential cone constraint builders.
\end{quote}
\end{quote}

\subsubsection{Model.GetExpConeBuilders()}
\label{\detokenize{csapi/model:id68}}\begin{quote}

\sphinxAtStartPar
Get builders of given exponential cone constraints in model.

\sphinxAtStartPar
\sphinxstylestrong{Synopsis}
\begin{quote}

\sphinxAtStartPar
\sphinxcode{\sphinxupquote{ExpConeBuilderArray GetExpConeBuilders(ExpConeArray cones)}}
\end{quote}

\sphinxAtStartPar
\sphinxstylestrong{Arguments}
\begin{quote}

\sphinxAtStartPar
\sphinxcode{\sphinxupquote{cones}}: array of exponential cone constraints.
\end{quote}

\sphinxAtStartPar
\sphinxstylestrong{Return}
\begin{quote}

\sphinxAtStartPar
array object of desired exponential cone constraint builders.
\end{quote}
\end{quote}

\subsubsection{Model.GetExpCones()}
\label{\detokenize{csapi/model:model-getexpcones}}\begin{quote}

\sphinxAtStartPar
Get all exponential cone constraints in model.

\sphinxAtStartPar
\sphinxstylestrong{Synopsis}
\begin{quote}

\sphinxAtStartPar
\sphinxcode{\sphinxupquote{ExpConeArray GetExpCones()}}
\end{quote}

\sphinxAtStartPar
\sphinxstylestrong{Return}
\begin{quote}

\sphinxAtStartPar
array object of exponential cone constraints.
\end{quote}
\end{quote}

\subsubsection{Model.GetGenConstr()}
\label{\detokenize{csapi/model:model-getgenconstr}}\begin{quote}

\sphinxAtStartPar
Get a general constraint of given index in model.

\sphinxAtStartPar
\sphinxstylestrong{Synopsis}
\begin{quote}

\sphinxAtStartPar
\sphinxcode{\sphinxupquote{GenConstr GetGenConstr(int idx)}}
\end{quote}

\sphinxAtStartPar
\sphinxstylestrong{Arguments}
\begin{quote}

\sphinxAtStartPar
\sphinxcode{\sphinxupquote{idx}}: index of the desired general constraint.
\end{quote}

\sphinxAtStartPar
\sphinxstylestrong{Return}
\begin{quote}

\sphinxAtStartPar
the desired general constraint object.
\end{quote}
\end{quote}

\subsubsection{Model.GetGenConstrByName()}
\label{\detokenize{csapi/model:model-getgenconstrbyname}}\begin{quote}

\sphinxAtStartPar
Get a general constraint of given name in model.

\sphinxAtStartPar
\sphinxstylestrong{Synopsis}
\begin{quote}

\sphinxAtStartPar
\sphinxcode{\sphinxupquote{GenConstr GetGenConstrByName(string name)}}
\end{quote}

\sphinxAtStartPar
\sphinxstylestrong{Arguments}
\begin{quote}

\sphinxAtStartPar
\sphinxcode{\sphinxupquote{name}}: name of the desired general constraint.
\end{quote}

\sphinxAtStartPar
\sphinxstylestrong{Return}
\begin{quote}

\sphinxAtStartPar
the desired general constraint object.
\end{quote}
\end{quote}

\subsubsection{Model.GetGenConstrIndicator()}
\label{\detokenize{csapi/model:model-getgenconstrindicator}}\begin{quote}

\sphinxAtStartPar
Get builder of given general constraint of type indicator.

\sphinxAtStartPar
\sphinxstylestrong{Synopsis}
\begin{quote}

\sphinxAtStartPar
\sphinxcode{\sphinxupquote{GenConstrBuilder GetGenConstrIndicator(GenConstr indicator)}}
\end{quote}

\sphinxAtStartPar
\sphinxstylestrong{Arguments}
\begin{quote}

\sphinxAtStartPar
\sphinxcode{\sphinxupquote{indicator}}: a general constraint of type indicator.
\end{quote}

\sphinxAtStartPar
\sphinxstylestrong{Return}
\begin{quote}

\sphinxAtStartPar
builder object of general constraint of type indicator.
\end{quote}
\end{quote}

\subsubsection{Model.GetGenConstrIndicators()}
\label{\detokenize{csapi/model:model-getgenconstrindicators}}\begin{quote}

\sphinxAtStartPar
Get builders of all general constraints in model.

\sphinxAtStartPar
\sphinxstylestrong{Synopsis}
\begin{quote}

\sphinxAtStartPar
\sphinxcode{\sphinxupquote{GenConstrBuilderArray GetGenConstrIndicators()}}
\end{quote}

\sphinxAtStartPar
\sphinxstylestrong{Return}
\begin{quote}

\sphinxAtStartPar
array object of general constraint builders.
\end{quote}
\end{quote}

\subsubsection{Model.GetGenConstrs()}
\label{\detokenize{csapi/model:model-getgenconstrs}}\begin{quote}

\sphinxAtStartPar
Get all general constraints in model.

\sphinxAtStartPar
\sphinxstylestrong{Synopsis}
\begin{quote}

\sphinxAtStartPar
\sphinxcode{\sphinxupquote{GenConstrArray GetGenConstrs()}}
\end{quote}

\sphinxAtStartPar
\sphinxstylestrong{Return}
\begin{quote}

\sphinxAtStartPar
array object of general constraints.
\end{quote}
\end{quote}

\subsubsection{Model.GetIndicatorIIS()}
\label{\detokenize{csapi/model:model-getindicatoriis}}\begin{quote}

\sphinxAtStartPar
Get IIS status of indicator constraints.

\sphinxAtStartPar
\sphinxstylestrong{Synopsis}
\begin{quote}

\sphinxAtStartPar
\sphinxcode{\sphinxupquote{int{[}{]} GetIndicatorIIS(GenConstrArray genconstrs)}}
\end{quote}

\sphinxAtStartPar
\sphinxstylestrong{Arguments}
\begin{quote}

\sphinxAtStartPar
\sphinxcode{\sphinxupquote{genconstrs}}: Array of indicator constraints.
\end{quote}

\sphinxAtStartPar
\sphinxstylestrong{Return}
\begin{quote}

\sphinxAtStartPar
IIS status of indicator constraints.
\end{quote}
\end{quote}

\subsubsection{Model.GetIndicatorIIS()}
\label{\detokenize{csapi/model:id69}}\begin{quote}

\sphinxAtStartPar
Get IIS status of indicator constraints.

\sphinxAtStartPar
\sphinxstylestrong{Synopsis}
\begin{quote}

\sphinxAtStartPar
\sphinxcode{\sphinxupquote{int{[}{]} GetIndicatorIIS(GenConstr{[}{]} genconstrs)}}
\end{quote}

\sphinxAtStartPar
\sphinxstylestrong{Arguments}
\begin{quote}

\sphinxAtStartPar
\sphinxcode{\sphinxupquote{genconstrs}}: Array of indicator constraints.
\end{quote}

\sphinxAtStartPar
\sphinxstylestrong{Return}
\begin{quote}

\sphinxAtStartPar
IIS status of indicator constraints.
\end{quote}
\end{quote}

\subsubsection{Model.GetIntAttr()}
\label{\detokenize{csapi/model:model-getintattr}}\begin{quote}

\sphinxAtStartPar
Get value of a COPT integer attribute.

\sphinxAtStartPar
\sphinxstylestrong{Synopsis}
\begin{quote}

\sphinxAtStartPar
\sphinxcode{\sphinxupquote{int GetIntAttr(string attr)}}
\end{quote}

\sphinxAtStartPar
\sphinxstylestrong{Arguments}
\begin{quote}

\sphinxAtStartPar
\sphinxcode{\sphinxupquote{attr}}: name of integer attribute.
\end{quote}

\sphinxAtStartPar
\sphinxstylestrong{Return}
\begin{quote}

\sphinxAtStartPar
value of integer attribute.
\end{quote}
\end{quote}

\subsubsection{Model.GetIntAttrN()}
\label{\detokenize{csapi/model:model-getintattrn}}\begin{quote}

\sphinxAtStartPar
Get value of a integer attribute of a multi\sphinxhyphen{}objective function.

\sphinxAtStartPar
\sphinxstylestrong{Synopsis}
\begin{quote}

\sphinxAtStartPar
\sphinxcode{\sphinxupquote{int GetIntAttrN(int idx, string attr)}}
\end{quote}

\sphinxAtStartPar
\sphinxstylestrong{Arguments}
\begin{quote}

\sphinxAtStartPar
\sphinxcode{\sphinxupquote{idx}}: index of a multi\sphinxhyphen{}objective function.

\sphinxAtStartPar
\sphinxcode{\sphinxupquote{attr}}: name of integer attribute.
\end{quote}

\sphinxAtStartPar
\sphinxstylestrong{Return}
\begin{quote}

\sphinxAtStartPar
value of integer attribute.
\end{quote}
\end{quote}

\subsubsection{Model.GetIntParam()}
\label{\detokenize{csapi/model:model-getintparam}}\begin{quote}

\sphinxAtStartPar
Get value of a COPT integer parameter.

\sphinxAtStartPar
\sphinxstylestrong{Synopsis}
\begin{quote}

\sphinxAtStartPar
\sphinxcode{\sphinxupquote{int GetIntParam(string param)}}
\end{quote}

\sphinxAtStartPar
\sphinxstylestrong{Arguments}
\begin{quote}

\sphinxAtStartPar
\sphinxcode{\sphinxupquote{param}}: name of integer parameter.
\end{quote}

\sphinxAtStartPar
\sphinxstylestrong{Return}
\begin{quote}

\sphinxAtStartPar
value of integer parameter.
\end{quote}
\end{quote}

\subsubsection{Model.GetIntParamN()}
\label{\detokenize{csapi/model:model-getintparamn}}\begin{quote}

\sphinxAtStartPar
Get value of an integer parameter of a multi\sphinxhyphen{}objective function.

\sphinxAtStartPar
\sphinxstylestrong{Synopsis}
\begin{quote}

\sphinxAtStartPar
\sphinxcode{\sphinxupquote{int GetIntParamN(int idx, string param)}}
\end{quote}

\sphinxAtStartPar
\sphinxstylestrong{Arguments}
\begin{quote}

\sphinxAtStartPar
\sphinxcode{\sphinxupquote{idx}}: index of a multi\sphinxhyphen{}objective function.

\sphinxAtStartPar
\sphinxcode{\sphinxupquote{param}}: name of integer parameter.
\end{quote}

\sphinxAtStartPar
\sphinxstylestrong{Return}
\begin{quote}

\sphinxAtStartPar
value of integer parameter.
\end{quote}
\end{quote}

\subsubsection{Model.GetLmiCoeff()}
\label{\detokenize{csapi/model:model-getlmicoeff}}\begin{quote}

\sphinxAtStartPar
Get the symmetric matrix of variable in LMI constraint.

\sphinxAtStartPar
\sphinxstylestrong{Synopsis}
\begin{quote}

\sphinxAtStartPar
\sphinxcode{\sphinxupquote{SymMatrix GetLmiCoeff(LmiConstraint constr, Var var)}}
\end{quote}

\sphinxAtStartPar
\sphinxstylestrong{Arguments}
\begin{quote}

\sphinxAtStartPar
\sphinxcode{\sphinxupquote{constr}}: The desired LMI constraint.

\sphinxAtStartPar
\sphinxcode{\sphinxupquote{var}}: The desired variable.
\end{quote}

\sphinxAtStartPar
\sphinxstylestrong{Return}
\begin{quote}

\sphinxAtStartPar
The associated coefficient matrix.
\end{quote}
\end{quote}

\subsubsection{Model.GetLmiConstr()}
\label{\detokenize{csapi/model:model-getlmiconstr}}\begin{quote}

\sphinxAtStartPar
Get LMI constraint of given index in model.

\sphinxAtStartPar
\sphinxstylestrong{Synopsis}
\begin{quote}

\sphinxAtStartPar
\sphinxcode{\sphinxupquote{LmiConstraint GetLmiConstr(int idx)}}
\end{quote}

\sphinxAtStartPar
\sphinxstylestrong{Arguments}
\begin{quote}

\sphinxAtStartPar
\sphinxcode{\sphinxupquote{idx}}: index of desired LMI constraint.
\end{quote}

\sphinxAtStartPar
\sphinxstylestrong{Return}
\begin{quote}

\sphinxAtStartPar
LMI constraint object.
\end{quote}
\end{quote}

\subsubsection{Model.GetLmiConstrByName()}
\label{\detokenize{csapi/model:model-getlmiconstrbyname}}\begin{quote}

\sphinxAtStartPar
Get LMI constraint of given name in model.

\sphinxAtStartPar
\sphinxstylestrong{Synopsis}
\begin{quote}

\sphinxAtStartPar
\sphinxcode{\sphinxupquote{LmiConstraint GetLmiConstrByName(string name)}}
\end{quote}

\sphinxAtStartPar
\sphinxstylestrong{Arguments}
\begin{quote}

\sphinxAtStartPar
\sphinxcode{\sphinxupquote{name}}: name of desired LMI constraint.
\end{quote}

\sphinxAtStartPar
\sphinxstylestrong{Return}
\begin{quote}

\sphinxAtStartPar
LMI constraint object.
\end{quote}
\end{quote}

\subsubsection{Model.GetLmiConstrs()}
\label{\detokenize{csapi/model:model-getlmiconstrs}}\begin{quote}

\sphinxAtStartPar
Get all LMI constraints in model.

\sphinxAtStartPar
\sphinxstylestrong{Synopsis}
\begin{quote}

\sphinxAtStartPar
\sphinxcode{\sphinxupquote{LmiConstrArray GetLmiConstrs()}}
\end{quote}

\sphinxAtStartPar
\sphinxstylestrong{Return}
\begin{quote}

\sphinxAtStartPar
array object of LMI constraints.
\end{quote}
\end{quote}

\subsubsection{Model.GetLmiRhs()}
\label{\detokenize{csapi/model:model-getlmirhs}}\begin{quote}

\sphinxAtStartPar
Get the symmetric matrix of constant of LMI constraint.

\sphinxAtStartPar
\sphinxstylestrong{Synopsis}
\begin{quote}

\sphinxAtStartPar
\sphinxcode{\sphinxupquote{SymMatrix GetLmiRhs(LmiConstraint constr)}}
\end{quote}

\sphinxAtStartPar
\sphinxstylestrong{Arguments}
\begin{quote}

\sphinxAtStartPar
\sphinxcode{\sphinxupquote{constr}}: The desired LMI constraint.
\end{quote}

\sphinxAtStartPar
\sphinxstylestrong{Return}
\begin{quote}

\sphinxAtStartPar
matrix of constant term.
\end{quote}
\end{quote}

\subsubsection{Model.GetLmiRow()}
\label{\detokenize{csapi/model:model-getlmirow}}\begin{quote}

\sphinxAtStartPar
Get variables and associated symmetric matrices that participate in a LMI constraint.

\sphinxAtStartPar
\sphinxstylestrong{Synopsis}
\begin{quote}

\sphinxAtStartPar
\sphinxcode{\sphinxupquote{LmiExpr GetLmiRow(LmiConstraint constr)}}
\end{quote}

\sphinxAtStartPar
\sphinxstylestrong{Arguments}
\begin{quote}

\sphinxAtStartPar
\sphinxcode{\sphinxupquote{constr}}: given LMI constraint object.
\end{quote}

\sphinxAtStartPar
\sphinxstylestrong{Return}
\begin{quote}

\sphinxAtStartPar
LMI expression object of LMI constraint.
\end{quote}
\end{quote}

\subsubsection{Model.GetLpSolution()}
\label{\detokenize{csapi/model:model-getlpsolution}}\begin{quote}

\sphinxAtStartPar
Get LP solution.

\sphinxAtStartPar
\sphinxstylestrong{Synopsis}
\begin{quote}

\sphinxAtStartPar
\sphinxcode{\sphinxupquote{void GetLpSolution(}}
\begin{quote}

\sphinxAtStartPar
\sphinxcode{\sphinxupquote{out double{[}{]} value,}}

\sphinxAtStartPar
\sphinxcode{\sphinxupquote{out double{[}{]} slack,}}

\sphinxAtStartPar
\sphinxcode{\sphinxupquote{out double{[}{]} rowDual,}}

\sphinxAtStartPar
\sphinxcode{\sphinxupquote{out double{[}{]} redCost)}}
\end{quote}
\end{quote}

\sphinxAtStartPar
\sphinxstylestrong{Arguments}
\begin{quote}

\sphinxAtStartPar
\sphinxcode{\sphinxupquote{value}}: out, solution values.

\sphinxAtStartPar
\sphinxcode{\sphinxupquote{slack}}: out, slack values.

\sphinxAtStartPar
\sphinxcode{\sphinxupquote{rowDual}}: out, dual values.

\sphinxAtStartPar
\sphinxcode{\sphinxupquote{redCost}}: out, reduced costs.
\end{quote}
\end{quote}

\subsubsection{Model.GetNlConstr()}
\label{\detokenize{csapi/model:model-getnlconstr}}\begin{quote}

\sphinxAtStartPar
Get a nonlinear constraint of given index in model.

\sphinxAtStartPar
\sphinxstylestrong{Synopsis}
\begin{quote}

\sphinxAtStartPar
\sphinxcode{\sphinxupquote{NlConstraint GetNlConstr(int idx)}}
\end{quote}

\sphinxAtStartPar
\sphinxstylestrong{Arguments}
\begin{quote}

\sphinxAtStartPar
\sphinxcode{\sphinxupquote{idx}}: index of the desired nonlinear constraint.
\end{quote}

\sphinxAtStartPar
\sphinxstylestrong{Return}
\begin{quote}

\sphinxAtStartPar
the desired nonlinear constraint object.
\end{quote}
\end{quote}

\subsubsection{Model.GetNlConstrBuilder()}
\label{\detokenize{csapi/model:model-getnlconstrbuilder}}\begin{quote}

\sphinxAtStartPar
Get builder of a nonlinear constraint in model, including nonlinear expression, sense and RHS.

\sphinxAtStartPar
\sphinxstylestrong{Synopsis}
\begin{quote}

\sphinxAtStartPar
\sphinxcode{\sphinxupquote{NlConstrBuilder GetNlConstrBuilder(NlConstraint constr)}}
\end{quote}

\sphinxAtStartPar
\sphinxstylestrong{Arguments}
\begin{quote}

\sphinxAtStartPar
\sphinxcode{\sphinxupquote{constr}}: a nonlinear constraint object.
\end{quote}

\sphinxAtStartPar
\sphinxstylestrong{Return}
\begin{quote}

\sphinxAtStartPar
nonlinear constraint builder object.
\end{quote}
\end{quote}

\subsubsection{Model.GetNlConstrBuilders()}
\label{\detokenize{csapi/model:model-getnlconstrbuilders}}\begin{quote}

\sphinxAtStartPar
Get builders of all nonlinear constraints in model.

\sphinxAtStartPar
\sphinxstylestrong{Synopsis}
\begin{quote}

\sphinxAtStartPar
\sphinxcode{\sphinxupquote{NlConstrBuilderArray GetNlConstrBuilders()}}
\end{quote}

\sphinxAtStartPar
\sphinxstylestrong{Return}
\begin{quote}

\sphinxAtStartPar
array object of nonlinear constraint builders.
\end{quote}
\end{quote}

\subsubsection{Model.GetNlConstrByName()}
\label{\detokenize{csapi/model:model-getnlconstrbyname}}\begin{quote}

\sphinxAtStartPar
Get a nonlinear constraint of given name in model.

\sphinxAtStartPar
\sphinxstylestrong{Synopsis}
\begin{quote}

\sphinxAtStartPar
\sphinxcode{\sphinxupquote{NlConstraint GetNlConstrByName(string name)}}
\end{quote}

\sphinxAtStartPar
\sphinxstylestrong{Arguments}
\begin{quote}

\sphinxAtStartPar
\sphinxcode{\sphinxupquote{name}}: name of the desired constraint.
\end{quote}

\sphinxAtStartPar
\sphinxstylestrong{Return}
\begin{quote}

\sphinxAtStartPar
the desired nonlinear constraint object.
\end{quote}
\end{quote}

\subsubsection{Model.GetNlConstrs()}
\label{\detokenize{csapi/model:model-getnlconstrs}}\begin{quote}

\sphinxAtStartPar
Get all nonlinear constraints in model.

\sphinxAtStartPar
\sphinxstylestrong{Synopsis}
\begin{quote}

\sphinxAtStartPar
\sphinxcode{\sphinxupquote{NlConstrArray GetNlConstrs()}}
\end{quote}

\sphinxAtStartPar
\sphinxstylestrong{Return}
\begin{quote}

\sphinxAtStartPar
array object of nonlinear constraints.
\end{quote}
\end{quote}

\subsubsection{Model.GetNlObjective()}
\label{\detokenize{csapi/model:model-getnlobjective}}\begin{quote}

\sphinxAtStartPar
Get nonlinear objective of model.

\sphinxAtStartPar
\sphinxstylestrong{Synopsis}
\begin{quote}

\sphinxAtStartPar
\sphinxcode{\sphinxupquote{NlExpr GetNlObjective()}}
\end{quote}

\sphinxAtStartPar
\sphinxstylestrong{Return}
\begin{quote}

\sphinxAtStartPar
a nonlinear expression object.
\end{quote}
\end{quote}

\subsubsection{Model.GetNlRow()}
\label{\detokenize{csapi/model:model-getnlrow}}\begin{quote}

\sphinxAtStartPar
Get nonlinear expression of a nonlinear constraint.

\sphinxAtStartPar
\sphinxstylestrong{Synopsis}
\begin{quote}

\sphinxAtStartPar
\sphinxcode{\sphinxupquote{NlExpr GetNlRow(NlConstraint constr)}}
\end{quote}

\sphinxAtStartPar
\sphinxstylestrong{Arguments}
\begin{quote}

\sphinxAtStartPar
\sphinxcode{\sphinxupquote{constr}}: a nonlinear constraint object.
\end{quote}

\sphinxAtStartPar
\sphinxstylestrong{Return}
\begin{quote}

\sphinxAtStartPar
output object of nonlinear expression.
\end{quote}
\end{quote}

\subsubsection{Model.GetObjective()}
\label{\detokenize{csapi/model:model-getobjective}}\begin{quote}

\sphinxAtStartPar
Get linear expression of objective for model.

\sphinxAtStartPar
\sphinxstylestrong{Synopsis}
\begin{quote}

\sphinxAtStartPar
\sphinxcode{\sphinxupquote{Expr GetObjective()}}
\end{quote}

\sphinxAtStartPar
\sphinxstylestrong{Return}
\begin{quote}

\sphinxAtStartPar
an linear expression object.
\end{quote}
\end{quote}

\subsubsection{Model.GetObjectiveN()}
\label{\detokenize{csapi/model:model-getobjectiven}}\begin{quote}

\sphinxAtStartPar
Get linear expression of a multi\sphinxhyphen{}objective function in model.

\sphinxAtStartPar
\sphinxstylestrong{Synopsis}
\begin{quote}

\sphinxAtStartPar
\sphinxcode{\sphinxupquote{Expr GetObjectiveN(int idx)}}
\end{quote}

\sphinxAtStartPar
\sphinxstylestrong{Arguments}
\begin{quote}

\sphinxAtStartPar
\sphinxcode{\sphinxupquote{idx}}: index of a multi\sphinxhyphen{}objective function.
\end{quote}

\sphinxAtStartPar
\sphinxstylestrong{Return}
\begin{quote}

\sphinxAtStartPar
a linear expression object.
\end{quote}
\end{quote}

\subsubsection{Model.GetObjParamN()}
\label{\detokenize{csapi/model:model-getobjparamn}}\begin{quote}

\sphinxAtStartPar
Get value of objective parameter of a multi\sphinxhyphen{}objective function.

\sphinxAtStartPar
\sphinxstylestrong{Synopsis}
\begin{quote}

\sphinxAtStartPar
\sphinxcode{\sphinxupquote{double GetObjParamN(int idx, string param)}}
\end{quote}

\sphinxAtStartPar
\sphinxstylestrong{Arguments}
\begin{quote}

\sphinxAtStartPar
\sphinxcode{\sphinxupquote{idx}}: index of a multi\sphinxhyphen{}objective function.

\sphinxAtStartPar
\sphinxcode{\sphinxupquote{param}}: name of objective parameter, including priority, weight, abstol and reltol.
\end{quote}

\sphinxAtStartPar
\sphinxstylestrong{Return}
\begin{quote}

\sphinxAtStartPar
value of objective parameter.
\end{quote}
\end{quote}

\subsubsection{Model.GetParamInfo()}
\label{\detokenize{csapi/model:model-getparaminfo}}\begin{quote}

\sphinxAtStartPar
Get current, default, minimum, maximum of COPT integer parameter.

\sphinxAtStartPar
\sphinxstylestrong{Synopsis}
\begin{quote}

\sphinxAtStartPar
\sphinxcode{\sphinxupquote{void GetParamInfo(}}
\begin{quote}

\sphinxAtStartPar
\sphinxcode{\sphinxupquote{string name,}}

\sphinxAtStartPar
\sphinxcode{\sphinxupquote{out int cur,}}

\sphinxAtStartPar
\sphinxcode{\sphinxupquote{out int def,}}

\sphinxAtStartPar
\sphinxcode{\sphinxupquote{out int min,}}

\sphinxAtStartPar
\sphinxcode{\sphinxupquote{out int max)}}
\end{quote}
\end{quote}

\sphinxAtStartPar
\sphinxstylestrong{Arguments}
\begin{quote}

\sphinxAtStartPar
\sphinxcode{\sphinxupquote{name}}: name of integer parameter.

\sphinxAtStartPar
\sphinxcode{\sphinxupquote{cur}}: out, current value of integer parameter.

\sphinxAtStartPar
\sphinxcode{\sphinxupquote{def}}: out, default value of integer parameter.

\sphinxAtStartPar
\sphinxcode{\sphinxupquote{min}}: out, minimum value of integer parameter.

\sphinxAtStartPar
\sphinxcode{\sphinxupquote{max}}: out, maximum value of integer parameter.
\end{quote}
\end{quote}

\subsubsection{Model.GetParamInfo()}
\label{\detokenize{csapi/model:id70}}\begin{quote}

\sphinxAtStartPar
Get current, default, minimum, maximum of COPT double parameter.

\sphinxAtStartPar
\sphinxstylestrong{Synopsis}
\begin{quote}

\sphinxAtStartPar
\sphinxcode{\sphinxupquote{void GetParamInfo(}}
\begin{quote}

\sphinxAtStartPar
\sphinxcode{\sphinxupquote{string name,}}

\sphinxAtStartPar
\sphinxcode{\sphinxupquote{out double cur,}}

\sphinxAtStartPar
\sphinxcode{\sphinxupquote{out double def,}}

\sphinxAtStartPar
\sphinxcode{\sphinxupquote{out double min,}}

\sphinxAtStartPar
\sphinxcode{\sphinxupquote{out double max)}}
\end{quote}
\end{quote}

\sphinxAtStartPar
\sphinxstylestrong{Arguments}
\begin{quote}

\sphinxAtStartPar
\sphinxcode{\sphinxupquote{name}}: name of integer parameter.

\sphinxAtStartPar
\sphinxcode{\sphinxupquote{cur}}: out, current value of double parameter.

\sphinxAtStartPar
\sphinxcode{\sphinxupquote{def}}: out, default value of double parameter.

\sphinxAtStartPar
\sphinxcode{\sphinxupquote{min}}: out, minimum value of double parameter.

\sphinxAtStartPar
\sphinxcode{\sphinxupquote{max}}: out, maximum value of double parameter.
\end{quote}
\end{quote}

\subsubsection{Model.GetPoolObjVal()}
\label{\detokenize{csapi/model:model-getpoolobjval}}\begin{quote}

\sphinxAtStartPar
Get the idx\sphinxhyphen{}th objective value in solution pool.

\sphinxAtStartPar
\sphinxstylestrong{Synopsis}
\begin{quote}

\sphinxAtStartPar
\sphinxcode{\sphinxupquote{double GetPoolObjVal(int idx)}}
\end{quote}

\sphinxAtStartPar
\sphinxstylestrong{Arguments}
\begin{quote}

\sphinxAtStartPar
\sphinxcode{\sphinxupquote{idx}}: Index of solution.
\end{quote}

\sphinxAtStartPar
\sphinxstylestrong{Return}
\begin{quote}

\sphinxAtStartPar
The requested objective value.
\end{quote}
\end{quote}

\subsubsection{Model.GetPoolObjValN()}
\label{\detokenize{csapi/model:model-getpoolobjvaln}}\begin{quote}

\sphinxAtStartPar
Get the objective value of required multi\sphinxhyphen{}objective function in solution pool.

\sphinxAtStartPar
\sphinxstylestrong{Synopsis}
\begin{quote}

\sphinxAtStartPar
\sphinxcode{\sphinxupquote{double GetPoolObjValN(int idx, int iSol)}}
\end{quote}

\sphinxAtStartPar
\sphinxstylestrong{Arguments}
\begin{quote}

\sphinxAtStartPar
\sphinxcode{\sphinxupquote{idx}}: index of a multi\sphinxhyphen{}objective function.

\sphinxAtStartPar
\sphinxcode{\sphinxupquote{iSol}}: index of solution.
\end{quote}

\sphinxAtStartPar
\sphinxstylestrong{Return}
\begin{quote}

\sphinxAtStartPar
value of required multi\sphinxhyphen{}objective function.
\end{quote}
\end{quote}

\subsubsection{Model.GetPoolSolution()}
\label{\detokenize{csapi/model:model-getpoolsolution}}\begin{quote}

\sphinxAtStartPar
Get the idx\sphinxhyphen{}th solution in solution pool.

\sphinxAtStartPar
\sphinxstylestrong{Synopsis}
\begin{quote}

\sphinxAtStartPar
\sphinxcode{\sphinxupquote{double{[}{]} GetPoolSolution(int idx, VarArray vars)}}
\end{quote}

\sphinxAtStartPar
\sphinxstylestrong{Arguments}
\begin{quote}

\sphinxAtStartPar
\sphinxcode{\sphinxupquote{idx}}: Index of solution.

\sphinxAtStartPar
\sphinxcode{\sphinxupquote{vars}}: The requested variables.
\end{quote}

\sphinxAtStartPar
\sphinxstylestrong{Return}
\begin{quote}

\sphinxAtStartPar
The requested solution.
\end{quote}
\end{quote}

\subsubsection{Model.GetPoolSolution()}
\label{\detokenize{csapi/model:id71}}\begin{quote}

\sphinxAtStartPar
Get the idx\sphinxhyphen{}th solution in solution pool.

\sphinxAtStartPar
\sphinxstylestrong{Synopsis}
\begin{quote}

\sphinxAtStartPar
\sphinxcode{\sphinxupquote{double{[}{]} GetPoolSolution(int idx, Var{[}{]} vars)}}
\end{quote}

\sphinxAtStartPar
\sphinxstylestrong{Arguments}
\begin{quote}

\sphinxAtStartPar
\sphinxcode{\sphinxupquote{idx}}: Index of solution.

\sphinxAtStartPar
\sphinxcode{\sphinxupquote{vars}}: The requested variables.
\end{quote}

\sphinxAtStartPar
\sphinxstylestrong{Return}
\begin{quote}

\sphinxAtStartPar
The requested solution.
\end{quote}
\end{quote}

\subsubsection{Model.GetPsdCoeff()}
\label{\detokenize{csapi/model:model-getpsdcoeff}}\begin{quote}

\sphinxAtStartPar
Get the symmetric matrix of PSD variable in a PSD constraint.

\sphinxAtStartPar
\sphinxstylestrong{Synopsis}
\begin{quote}

\sphinxAtStartPar
\sphinxcode{\sphinxupquote{SymMatrix GetPsdCoeff(PsdConstraint constr, PsdVar var)}}
\end{quote}

\sphinxAtStartPar
\sphinxstylestrong{Arguments}
\begin{quote}

\sphinxAtStartPar
\sphinxcode{\sphinxupquote{constr}}: The desired PSD constraint.

\sphinxAtStartPar
\sphinxcode{\sphinxupquote{var}}: The desired PSD variable.
\end{quote}

\sphinxAtStartPar
\sphinxstylestrong{Return}
\begin{quote}

\sphinxAtStartPar
The associated coefficient matrix.
\end{quote}
\end{quote}

\subsubsection{Model.GetPsdConstr()}
\label{\detokenize{csapi/model:model-getpsdconstr}}\begin{quote}

\sphinxAtStartPar
Get a PSD constraint of given index in model.

\sphinxAtStartPar
\sphinxstylestrong{Synopsis}
\begin{quote}

\sphinxAtStartPar
\sphinxcode{\sphinxupquote{PsdConstraint GetPsdConstr(int idx)}}
\end{quote}

\sphinxAtStartPar
\sphinxstylestrong{Arguments}
\begin{quote}

\sphinxAtStartPar
\sphinxcode{\sphinxupquote{idx}}: index of desired PSD constraint.
\end{quote}

\sphinxAtStartPar
\sphinxstylestrong{Return}
\begin{quote}

\sphinxAtStartPar
PSD constraint object.
\end{quote}
\end{quote}

\subsubsection{Model.GetPsdConstrBuilder()}
\label{\detokenize{csapi/model:model-getpsdconstrbuilder}}\begin{quote}

\sphinxAtStartPar
Get builder of a PSD constraint in model, including PSD variables, sense and associated symmtric matrices.

\sphinxAtStartPar
\sphinxstylestrong{Synopsis}
\begin{quote}

\sphinxAtStartPar
\sphinxcode{\sphinxupquote{PsdConstrBuilder GetPsdConstrBuilder(PsdConstraint constr)}}
\end{quote}

\sphinxAtStartPar
\sphinxstylestrong{Arguments}
\begin{quote}

\sphinxAtStartPar
\sphinxcode{\sphinxupquote{constr}}: PSD constraint object.
\end{quote}

\sphinxAtStartPar
\sphinxstylestrong{Return}
\begin{quote}

\sphinxAtStartPar
PSD constraint builder object.
\end{quote}
\end{quote}

\subsubsection{Model.GetPsdConstrBuilders()}
\label{\detokenize{csapi/model:model-getpsdconstrbuilders}}\begin{quote}

\sphinxAtStartPar
Get builders of all PSD constraints in model.

\sphinxAtStartPar
\sphinxstylestrong{Synopsis}
\begin{quote}

\sphinxAtStartPar
\sphinxcode{\sphinxupquote{PsdConstrBuilderArray GetPsdConstrBuilders()}}
\end{quote}

\sphinxAtStartPar
\sphinxstylestrong{Return}
\begin{quote}

\sphinxAtStartPar
array object of PSD constraint builders.
\end{quote}
\end{quote}

\subsubsection{Model.GetPsdConstrByName()}
\label{\detokenize{csapi/model:model-getpsdconstrbyname}}\begin{quote}

\sphinxAtStartPar
Get a PSD constraint of given name in model.

\sphinxAtStartPar
\sphinxstylestrong{Synopsis}
\begin{quote}

\sphinxAtStartPar
\sphinxcode{\sphinxupquote{PsdConstraint GetPsdConstrByName(string name)}}
\end{quote}

\sphinxAtStartPar
\sphinxstylestrong{Arguments}
\begin{quote}

\sphinxAtStartPar
\sphinxcode{\sphinxupquote{name}}: name of desired PSD constraint.
\end{quote}

\sphinxAtStartPar
\sphinxstylestrong{Return}
\begin{quote}

\sphinxAtStartPar
PSD constraint object.
\end{quote}
\end{quote}

\subsubsection{Model.GetPsdConstrs()}
\label{\detokenize{csapi/model:model-getpsdconstrs}}\begin{quote}

\sphinxAtStartPar
Get all PSD constraints in model.

\sphinxAtStartPar
\sphinxstylestrong{Synopsis}
\begin{quote}

\sphinxAtStartPar
\sphinxcode{\sphinxupquote{PsdConstrArray GetPsdConstrs()}}
\end{quote}

\sphinxAtStartPar
\sphinxstylestrong{Return}
\begin{quote}

\sphinxAtStartPar
array object of PSD constraints.
\end{quote}
\end{quote}

\subsubsection{Model.GetPsdObjective()}
\label{\detokenize{csapi/model:model-getpsdobjective}}\begin{quote}

\sphinxAtStartPar
Get PSD objective of model.

\sphinxAtStartPar
\sphinxstylestrong{Synopsis}
\begin{quote}

\sphinxAtStartPar
\sphinxcode{\sphinxupquote{PsdExpr GetPsdObjective()}}
\end{quote}

\sphinxAtStartPar
\sphinxstylestrong{Return}
\begin{quote}

\sphinxAtStartPar
a PSD expression object.
\end{quote}
\end{quote}

\subsubsection{Model.GetPsdRow()}
\label{\detokenize{csapi/model:model-getpsdrow}}\begin{quote}

\sphinxAtStartPar
Get PSD variables and associated symmetric matrices that participate in a PSD constraint.

\sphinxAtStartPar
\sphinxstylestrong{Synopsis}
\begin{quote}

\sphinxAtStartPar
\sphinxcode{\sphinxupquote{PsdExpr GetPsdRow(PsdConstraint constr)}}
\end{quote}

\sphinxAtStartPar
\sphinxstylestrong{Arguments}
\begin{quote}

\sphinxAtStartPar
\sphinxcode{\sphinxupquote{constr}}: PSD constraint object.
\end{quote}

\sphinxAtStartPar
\sphinxstylestrong{Return}
\begin{quote}

\sphinxAtStartPar
PSD expression object of the PSD constraint.
\end{quote}
\end{quote}

\subsubsection{Model.GetPsdSolution()}
\label{\detokenize{csapi/model:model-getpsdsolution}}\begin{quote}

\sphinxAtStartPar
Get PSD solution.

\sphinxAtStartPar
\sphinxstylestrong{Synopsis}
\begin{quote}

\sphinxAtStartPar
\sphinxcode{\sphinxupquote{void GetPsdSolution(}}
\begin{quote}

\sphinxAtStartPar
\sphinxcode{\sphinxupquote{out double{[}{]} psdValue,}}

\sphinxAtStartPar
\sphinxcode{\sphinxupquote{out double{[}{]} psdSlack,}}

\sphinxAtStartPar
\sphinxcode{\sphinxupquote{out double{[}{]} psdRowDual,}}

\sphinxAtStartPar
\sphinxcode{\sphinxupquote{out double{[}{]} psdRedCost)}}
\end{quote}
\end{quote}

\sphinxAtStartPar
\sphinxstylestrong{Arguments}
\begin{quote}

\sphinxAtStartPar
\sphinxcode{\sphinxupquote{psdValue}}: out, solution of PSD variables.

\sphinxAtStartPar
\sphinxcode{\sphinxupquote{psdSlack}}: out, slack of PSD constraints.

\sphinxAtStartPar
\sphinxcode{\sphinxupquote{psdRowDual}}: out, dual of PSD constraints.

\sphinxAtStartPar
\sphinxcode{\sphinxupquote{psdRedCost}}: out, reduced costs of PSD variables.
\end{quote}
\end{quote}

\subsubsection{Model.GetPsdVar()}
\label{\detokenize{csapi/model:model-getpsdvar}}\begin{quote}

\sphinxAtStartPar
Get a PSD variable of given index in model.

\sphinxAtStartPar
\sphinxstylestrong{Synopsis}
\begin{quote}

\sphinxAtStartPar
\sphinxcode{\sphinxupquote{PsdVar GetPsdVar(int idx)}}
\end{quote}

\sphinxAtStartPar
\sphinxstylestrong{Arguments}
\begin{quote}

\sphinxAtStartPar
\sphinxcode{\sphinxupquote{idx}}: index of the desired PSD variable.
\end{quote}

\sphinxAtStartPar
\sphinxstylestrong{Return}
\begin{quote}

\sphinxAtStartPar
the desired PSD variable object.
\end{quote}
\end{quote}

\subsubsection{Model.GetPsdVarByName()}
\label{\detokenize{csapi/model:model-getpsdvarbyname}}\begin{quote}

\sphinxAtStartPar
Get a PSD variable of given name in model.

\sphinxAtStartPar
\sphinxstylestrong{Synopsis}
\begin{quote}

\sphinxAtStartPar
\sphinxcode{\sphinxupquote{PsdVar GetPsdVarByName(string name)}}
\end{quote}

\sphinxAtStartPar
\sphinxstylestrong{Arguments}
\begin{quote}

\sphinxAtStartPar
\sphinxcode{\sphinxupquote{name}}: name of the desired PSD variable.
\end{quote}

\sphinxAtStartPar
\sphinxstylestrong{Return}
\begin{quote}

\sphinxAtStartPar
the desired PSD variable object.
\end{quote}
\end{quote}

\subsubsection{Model.GetPsdVars()}
\label{\detokenize{csapi/model:model-getpsdvars}}\begin{quote}

\sphinxAtStartPar
Get all PSD variables in model.

\sphinxAtStartPar
\sphinxstylestrong{Synopsis}
\begin{quote}

\sphinxAtStartPar
\sphinxcode{\sphinxupquote{PsdVarArray GetPsdVars()}}
\end{quote}

\sphinxAtStartPar
\sphinxstylestrong{Return}
\begin{quote}

\sphinxAtStartPar
array object of PSD variables.
\end{quote}
\end{quote}

\subsubsection{Model.GetQConstr()}
\label{\detokenize{csapi/model:model-getqconstr}}\begin{quote}

\sphinxAtStartPar
Get a quadratic constraint of given index in model.

\sphinxAtStartPar
\sphinxstylestrong{Synopsis}
\begin{quote}

\sphinxAtStartPar
\sphinxcode{\sphinxupquote{QConstraint GetQConstr(int idx)}}
\end{quote}

\sphinxAtStartPar
\sphinxstylestrong{Arguments}
\begin{quote}

\sphinxAtStartPar
\sphinxcode{\sphinxupquote{idx}}: index of the desired quadratic constraint.
\end{quote}

\sphinxAtStartPar
\sphinxstylestrong{Return}
\begin{quote}

\sphinxAtStartPar
the desired quadratic constraint object.
\end{quote}
\end{quote}

\subsubsection{Model.GetQConstrBuilder()}
\label{\detokenize{csapi/model:model-getqconstrbuilder}}\begin{quote}

\sphinxAtStartPar
Get builder of a constraint in model, including variables and associated coefficients, sense and RHS.

\sphinxAtStartPar
\sphinxstylestrong{Synopsis}
\begin{quote}

\sphinxAtStartPar
\sphinxcode{\sphinxupquote{QConstrBuilder GetQConstrBuilder(QConstraint constr)}}
\end{quote}

\sphinxAtStartPar
\sphinxstylestrong{Arguments}
\begin{quote}

\sphinxAtStartPar
\sphinxcode{\sphinxupquote{constr}}: a constraint object.
\end{quote}

\sphinxAtStartPar
\sphinxstylestrong{Return}
\begin{quote}

\sphinxAtStartPar
constraint builder object.
\end{quote}
\end{quote}

\subsubsection{Model.GetQConstrBuilders()}
\label{\detokenize{csapi/model:model-getqconstrbuilders}}\begin{quote}

\sphinxAtStartPar
Get builders of all constraints in model.

\sphinxAtStartPar
\sphinxstylestrong{Synopsis}
\begin{quote}

\sphinxAtStartPar
\sphinxcode{\sphinxupquote{QConstrBuilderArray GetQConstrBuilders()}}
\end{quote}

\sphinxAtStartPar
\sphinxstylestrong{Return}
\begin{quote}

\sphinxAtStartPar
array object of constraint builders.
\end{quote}
\end{quote}

\subsubsection{Model.GetQConstrByName()}
\label{\detokenize{csapi/model:model-getqconstrbyname}}\begin{quote}

\sphinxAtStartPar
Get a quadratic constraint of given name in model.

\sphinxAtStartPar
\sphinxstylestrong{Synopsis}
\begin{quote}

\sphinxAtStartPar
\sphinxcode{\sphinxupquote{QConstraint GetQConstrByName(string name)}}
\end{quote}

\sphinxAtStartPar
\sphinxstylestrong{Arguments}
\begin{quote}

\sphinxAtStartPar
\sphinxcode{\sphinxupquote{name}}: name of the desired constraint.
\end{quote}

\sphinxAtStartPar
\sphinxstylestrong{Return}
\begin{quote}

\sphinxAtStartPar
the desired quadratic constraint object.
\end{quote}
\end{quote}

\subsubsection{Model.GetQConstrs()}
\label{\detokenize{csapi/model:model-getqconstrs}}\begin{quote}

\sphinxAtStartPar
Get all quadratic constraints in model.

\sphinxAtStartPar
\sphinxstylestrong{Synopsis}
\begin{quote}

\sphinxAtStartPar
\sphinxcode{\sphinxupquote{QConstrArray GetQConstrs()}}
\end{quote}

\sphinxAtStartPar
\sphinxstylestrong{Return}
\begin{quote}

\sphinxAtStartPar
array object of quadratic constraints.
\end{quote}
\end{quote}

\subsubsection{Model.GetQuadObjective()}
\label{\detokenize{csapi/model:model-getquadobjective}}\begin{quote}

\sphinxAtStartPar
Get quadratic objective of model.

\sphinxAtStartPar
\sphinxstylestrong{Synopsis}
\begin{quote}

\sphinxAtStartPar
\sphinxcode{\sphinxupquote{QuadExpr GetQuadObjective()}}
\end{quote}

\sphinxAtStartPar
\sphinxstylestrong{Return}
\begin{quote}

\sphinxAtStartPar
a quadratic expression object.
\end{quote}
\end{quote}

\subsubsection{Model.GetQuadRow()}
\label{\detokenize{csapi/model:model-getquadrow}}\begin{quote}

\sphinxAtStartPar
Get two variables and associated coefficients that participate in a quadratic constraint.

\sphinxAtStartPar
\sphinxstylestrong{Synopsis}
\begin{quote}

\sphinxAtStartPar
\sphinxcode{\sphinxupquote{QuadExpr GetQuadRow(QConstraint constr)}}
\end{quote}

\sphinxAtStartPar
\sphinxstylestrong{Arguments}
\begin{quote}

\sphinxAtStartPar
\sphinxcode{\sphinxupquote{constr}}: a quadratic constraint object.
\end{quote}

\sphinxAtStartPar
\sphinxstylestrong{Return}
\begin{quote}

\sphinxAtStartPar
quadratic expression object of the constraint.
\end{quote}
\end{quote}

\subsubsection{Model.GetRow()}
\label{\detokenize{csapi/model:model-getrow}}\begin{quote}

\sphinxAtStartPar
Get variables that participate in a constraint, and the associated coefficients.

\sphinxAtStartPar
\sphinxstylestrong{Synopsis}
\begin{quote}

\sphinxAtStartPar
\sphinxcode{\sphinxupquote{Expr GetRow(Constraint constr)}}
\end{quote}

\sphinxAtStartPar
\sphinxstylestrong{Arguments}
\begin{quote}

\sphinxAtStartPar
\sphinxcode{\sphinxupquote{constr}}: a constraint object.
\end{quote}

\sphinxAtStartPar
\sphinxstylestrong{Return}
\begin{quote}

\sphinxAtStartPar
expression object of the constraint.
\end{quote}
\end{quote}

\subsubsection{Model.GetRowBasis()}
\label{\detokenize{csapi/model:model-getrowbasis}}\begin{quote}

\sphinxAtStartPar
Get status of row basis.

\sphinxAtStartPar
\sphinxstylestrong{Synopsis}
\begin{quote}

\sphinxAtStartPar
\sphinxcode{\sphinxupquote{int{[}{]} GetRowBasis()}}
\end{quote}

\sphinxAtStartPar
\sphinxstylestrong{Return}
\begin{quote}

\sphinxAtStartPar
basis status.
\end{quote}
\end{quote}

\subsubsection{Model.GetSolution()}
\label{\detokenize{csapi/model:model-getsolution}}\begin{quote}

\sphinxAtStartPar
Get MIP solution.

\sphinxAtStartPar
\sphinxstylestrong{Synopsis}
\begin{quote}

\sphinxAtStartPar
\sphinxcode{\sphinxupquote{double{[}{]} GetSolution()}}
\end{quote}

\sphinxAtStartPar
\sphinxstylestrong{Return}
\begin{quote}

\sphinxAtStartPar
solution values.
\end{quote}
\end{quote}

\subsubsection{Model.GetSos()}
\label{\detokenize{csapi/model:model-getsos}}\begin{quote}

\sphinxAtStartPar
Get a SOS constraint of given index in model.

\sphinxAtStartPar
\sphinxstylestrong{Synopsis}
\begin{quote}

\sphinxAtStartPar
\sphinxcode{\sphinxupquote{Sos GetSos(int idx)}}
\end{quote}

\sphinxAtStartPar
\sphinxstylestrong{Arguments}
\begin{quote}

\sphinxAtStartPar
\sphinxcode{\sphinxupquote{idx}}: index of the desired SOS constraint.
\end{quote}

\sphinxAtStartPar
\sphinxstylestrong{Return}
\begin{quote}

\sphinxAtStartPar
the desired SOS constraint object.
\end{quote}
\end{quote}

\subsubsection{Model.GetSosBuilders()}
\label{\detokenize{csapi/model:model-getsosbuilders}}\begin{quote}

\sphinxAtStartPar
Get builders of all SOS constraints in model.

\sphinxAtStartPar
\sphinxstylestrong{Synopsis}
\begin{quote}

\sphinxAtStartPar
\sphinxcode{\sphinxupquote{SosBuilderArray GetSosBuilders()}}
\end{quote}

\sphinxAtStartPar
\sphinxstylestrong{Return}
\begin{quote}

\sphinxAtStartPar
array object of SOS constraint builders.
\end{quote}
\end{quote}

\subsubsection{Model.GetSosBuilders()}
\label{\detokenize{csapi/model:id72}}\begin{quote}

\sphinxAtStartPar
Get builders of given SOS constraints in model.

\sphinxAtStartPar
\sphinxstylestrong{Synopsis}
\begin{quote}

\sphinxAtStartPar
\sphinxcode{\sphinxupquote{SosBuilderArray GetSosBuilders(Sos{[}{]} soss)}}
\end{quote}

\sphinxAtStartPar
\sphinxstylestrong{Arguments}
\begin{quote}

\sphinxAtStartPar
\sphinxcode{\sphinxupquote{soss}}: array of SOS constraints.
\end{quote}

\sphinxAtStartPar
\sphinxstylestrong{Return}
\begin{quote}

\sphinxAtStartPar
array object of desired SOS constraint builders.
\end{quote}
\end{quote}

\subsubsection{Model.GetSosBuilders()}
\label{\detokenize{csapi/model:id73}}\begin{quote}

\sphinxAtStartPar
Get builders of given SOS constraints in model.

\sphinxAtStartPar
\sphinxstylestrong{Synopsis}
\begin{quote}

\sphinxAtStartPar
\sphinxcode{\sphinxupquote{SosBuilderArray GetSosBuilders(SosArray soss)}}
\end{quote}

\sphinxAtStartPar
\sphinxstylestrong{Arguments}
\begin{quote}

\sphinxAtStartPar
\sphinxcode{\sphinxupquote{soss}}: array of SOS constraints.
\end{quote}

\sphinxAtStartPar
\sphinxstylestrong{Return}
\begin{quote}

\sphinxAtStartPar
array object of desired SOS constraint builders.
\end{quote}
\end{quote}

\subsubsection{Model.GetSOSIIS()}
\label{\detokenize{csapi/model:model-getsosiis}}\begin{quote}

\sphinxAtStartPar
Get IIS status of SOS constraints.

\sphinxAtStartPar
\sphinxstylestrong{Synopsis}
\begin{quote}

\sphinxAtStartPar
\sphinxcode{\sphinxupquote{int{[}{]} GetSOSIIS(SosArray soss)}}
\end{quote}

\sphinxAtStartPar
\sphinxstylestrong{Arguments}
\begin{quote}

\sphinxAtStartPar
\sphinxcode{\sphinxupquote{soss}}: Array of SOS constraints.
\end{quote}

\sphinxAtStartPar
\sphinxstylestrong{Return}
\begin{quote}

\sphinxAtStartPar
IIS status of SOS constraints.
\end{quote}
\end{quote}

\subsubsection{Model.GetSOSIIS()}
\label{\detokenize{csapi/model:id74}}\begin{quote}

\sphinxAtStartPar
Get IIS status of SOS constraints.

\sphinxAtStartPar
\sphinxstylestrong{Synopsis}
\begin{quote}

\sphinxAtStartPar
\sphinxcode{\sphinxupquote{int{[}{]} GetSOSIIS(Sos{[}{]} soss)}}
\end{quote}

\sphinxAtStartPar
\sphinxstylestrong{Arguments}
\begin{quote}

\sphinxAtStartPar
\sphinxcode{\sphinxupquote{soss}}: Array of SOS constraints.
\end{quote}

\sphinxAtStartPar
\sphinxstylestrong{Return}
\begin{quote}

\sphinxAtStartPar
IIS status of SOS constraints.
\end{quote}
\end{quote}

\subsubsection{Model.GetSoss()}
\label{\detokenize{csapi/model:model-getsoss}}\begin{quote}

\sphinxAtStartPar
Get all SOS constraints in model.

\sphinxAtStartPar
\sphinxstylestrong{Synopsis}
\begin{quote}

\sphinxAtStartPar
\sphinxcode{\sphinxupquote{SosArray GetSoss()}}
\end{quote}

\sphinxAtStartPar
\sphinxstylestrong{Return}
\begin{quote}

\sphinxAtStartPar
array object of SOS constraints.
\end{quote}
\end{quote}

\subsubsection{Model.GetSymMat()}
\label{\detokenize{csapi/model:model-getsymmat}}\begin{quote}

\sphinxAtStartPar
Get a symmetric matrix of given index in model.

\sphinxAtStartPar
\sphinxstylestrong{Synopsis}
\begin{quote}

\sphinxAtStartPar
\sphinxcode{\sphinxupquote{SymMatrix GetSymMat(int idx)}}
\end{quote}

\sphinxAtStartPar
\sphinxstylestrong{Arguments}
\begin{quote}

\sphinxAtStartPar
\sphinxcode{\sphinxupquote{idx}}: index of the desired symmetric matrix.
\end{quote}

\sphinxAtStartPar
\sphinxstylestrong{Return}
\begin{quote}

\sphinxAtStartPar
the desired symmetric matrix object.
\end{quote}
\end{quote}

\subsubsection{Model.GetVar()}
\label{\detokenize{csapi/model:model-getvar}}\begin{quote}

\sphinxAtStartPar
Get a variable of given index in model.

\sphinxAtStartPar
\sphinxstylestrong{Synopsis}
\begin{quote}

\sphinxAtStartPar
\sphinxcode{\sphinxupquote{Var GetVar(int idx)}}
\end{quote}

\sphinxAtStartPar
\sphinxstylestrong{Arguments}
\begin{quote}

\sphinxAtStartPar
\sphinxcode{\sphinxupquote{idx}}: index of the desired variable.
\end{quote}

\sphinxAtStartPar
\sphinxstylestrong{Return}
\begin{quote}

\sphinxAtStartPar
the desired variable object.
\end{quote}
\end{quote}

\subsubsection{Model.GetVarByName()}
\label{\detokenize{csapi/model:model-getvarbyname}}\begin{quote}

\sphinxAtStartPar
Get a variable of given name in model.

\sphinxAtStartPar
\sphinxstylestrong{Synopsis}
\begin{quote}

\sphinxAtStartPar
\sphinxcode{\sphinxupquote{Var GetVarByName(string name)}}
\end{quote}

\sphinxAtStartPar
\sphinxstylestrong{Arguments}
\begin{quote}

\sphinxAtStartPar
\sphinxcode{\sphinxupquote{name}}: name of the desired variable.
\end{quote}

\sphinxAtStartPar
\sphinxstylestrong{Return}
\begin{quote}

\sphinxAtStartPar
the desired variable object.
\end{quote}
\end{quote}

\subsubsection{Model.GetVarLowerIIS()}
\label{\detokenize{csapi/model:model-getvarloweriis}}\begin{quote}

\sphinxAtStartPar
Get IIS status of lower bounds of variables.

\sphinxAtStartPar
\sphinxstylestrong{Synopsis}
\begin{quote}

\sphinxAtStartPar
\sphinxcode{\sphinxupquote{int{[}{]} GetVarLowerIIS(VarArray vars)}}
\end{quote}

\sphinxAtStartPar
\sphinxstylestrong{Arguments}
\begin{quote}

\sphinxAtStartPar
\sphinxcode{\sphinxupquote{vars}}: Array of variables.
\end{quote}

\sphinxAtStartPar
\sphinxstylestrong{Return}
\begin{quote}

\sphinxAtStartPar
IIS status of lower bounds of variables.
\end{quote}
\end{quote}

\subsubsection{Model.GetVarLowerIIS()}
\label{\detokenize{csapi/model:id75}}\begin{quote}

\sphinxAtStartPar
Get IIS status of lower bounds of variables.

\sphinxAtStartPar
\sphinxstylestrong{Synopsis}
\begin{quote}

\sphinxAtStartPar
\sphinxcode{\sphinxupquote{int{[}{]} GetVarLowerIIS(Var{[}{]} vars)}}
\end{quote}

\sphinxAtStartPar
\sphinxstylestrong{Arguments}
\begin{quote}

\sphinxAtStartPar
\sphinxcode{\sphinxupquote{vars}}: Array of variables.
\end{quote}

\sphinxAtStartPar
\sphinxstylestrong{Return}
\begin{quote}

\sphinxAtStartPar
IIS status of lower bounds of variables.
\end{quote}
\end{quote}

\subsubsection{Model.GetVars()}
\label{\detokenize{csapi/model:model-getvars}}\begin{quote}

\sphinxAtStartPar
Get all variables in model.

\sphinxAtStartPar
\sphinxstylestrong{Synopsis}
\begin{quote}

\sphinxAtStartPar
\sphinxcode{\sphinxupquote{VarArray GetVars()}}
\end{quote}

\sphinxAtStartPar
\sphinxstylestrong{Return}
\begin{quote}

\sphinxAtStartPar
variable array object.
\end{quote}
\end{quote}

\subsubsection{Model.GetVarUpperIIS()}
\label{\detokenize{csapi/model:model-getvarupperiis}}\begin{quote}

\sphinxAtStartPar
Get IIS status of upper bounds of variables.

\sphinxAtStartPar
\sphinxstylestrong{Synopsis}
\begin{quote}

\sphinxAtStartPar
\sphinxcode{\sphinxupquote{int{[}{]} GetVarUpperIIS(VarArray vars)}}
\end{quote}

\sphinxAtStartPar
\sphinxstylestrong{Arguments}
\begin{quote}

\sphinxAtStartPar
\sphinxcode{\sphinxupquote{vars}}: Array of variables.
\end{quote}

\sphinxAtStartPar
\sphinxstylestrong{Return}
\begin{quote}

\sphinxAtStartPar
IIS status of upper bounds of variables.
\end{quote}
\end{quote}

\subsubsection{Model.GetVarUpperIIS()}
\label{\detokenize{csapi/model:id76}}\begin{quote}

\sphinxAtStartPar
Get IIS status of upper bounds of variables.

\sphinxAtStartPar
\sphinxstylestrong{Synopsis}
\begin{quote}

\sphinxAtStartPar
\sphinxcode{\sphinxupquote{int{[}{]} GetVarUpperIIS(Var{[}{]} vars)}}
\end{quote}

\sphinxAtStartPar
\sphinxstylestrong{Arguments}
\begin{quote}

\sphinxAtStartPar
\sphinxcode{\sphinxupquote{vars}}: Array of variables.
\end{quote}

\sphinxAtStartPar
\sphinxstylestrong{Return}
\begin{quote}

\sphinxAtStartPar
IIS status of upper bounds of variables.
\end{quote}
\end{quote}

\subsubsection{Model.Interrupt()}
\label{\detokenize{csapi/model:model-interrupt}}\begin{quote}

\sphinxAtStartPar
Interrupt optimization of current problem.

\sphinxAtStartPar
\sphinxstylestrong{Synopsis}
\begin{quote}

\sphinxAtStartPar
\sphinxcode{\sphinxupquote{void Interrupt()}}
\end{quote}
\end{quote}

\subsubsection{Model.LoadMipStart()}
\label{\detokenize{csapi/model:model-loadmipstart}}\begin{quote}

\sphinxAtStartPar
Load final initial values of variables to the problem.

\sphinxAtStartPar
\sphinxstylestrong{Synopsis}
\begin{quote}

\sphinxAtStartPar
\sphinxcode{\sphinxupquote{void LoadMipStart()}}
\end{quote}
\end{quote}

\subsubsection{Model.LoadTuneParam()}
\label{\detokenize{csapi/model:model-loadtuneparam}}\begin{quote}

\sphinxAtStartPar
Load specified tuned parameters into model.

\sphinxAtStartPar
\sphinxstylestrong{Synopsis}
\begin{quote}

\sphinxAtStartPar
\sphinxcode{\sphinxupquote{void LoadTuneParam(int idx)}}
\end{quote}

\sphinxAtStartPar
\sphinxstylestrong{Arguments}
\begin{quote}

\sphinxAtStartPar
\sphinxcode{\sphinxupquote{idx}}: Index of tuned parameters.
\end{quote}
\end{quote}

\subsubsection{Model.Read()}
\label{\detokenize{csapi/model:model-read}}\begin{quote}

\sphinxAtStartPar
Read problem, solution, basis, MIP start or COPT parameters from file.

\sphinxAtStartPar
\sphinxstylestrong{Synopsis}
\begin{quote}

\sphinxAtStartPar
\sphinxcode{\sphinxupquote{void Read(string filename)}}
\end{quote}

\sphinxAtStartPar
\sphinxstylestrong{Arguments}
\begin{quote}

\sphinxAtStartPar
\sphinxcode{\sphinxupquote{filename}}: an input file name.
\end{quote}
\end{quote}

\subsubsection{Model.ReadBasis()}
\label{\detokenize{csapi/model:model-readbasis}}\begin{quote}

\sphinxAtStartPar
Read basis from file.

\sphinxAtStartPar
\sphinxstylestrong{Synopsis}
\begin{quote}

\sphinxAtStartPar
\sphinxcode{\sphinxupquote{void ReadBasis(string filename)}}
\end{quote}

\sphinxAtStartPar
\sphinxstylestrong{Arguments}
\begin{quote}

\sphinxAtStartPar
\sphinxcode{\sphinxupquote{filename}}: an input file name
\end{quote}
\end{quote}

\subsubsection{Model.ReadBin()}
\label{\detokenize{csapi/model:model-readbin}}\begin{quote}

\sphinxAtStartPar
Read problem in COPT binary format from file.

\sphinxAtStartPar
\sphinxstylestrong{Synopsis}
\begin{quote}

\sphinxAtStartPar
\sphinxcode{\sphinxupquote{void ReadBin(string filename)}}
\end{quote}

\sphinxAtStartPar
\sphinxstylestrong{Arguments}
\begin{quote}

\sphinxAtStartPar
\sphinxcode{\sphinxupquote{filename}}: an input file name.
\end{quote}
\end{quote}

\subsubsection{Model.ReadCbf()}
\label{\detokenize{csapi/model:model-readcbf}}\begin{quote}

\sphinxAtStartPar
Read problem in CBF format from file.

\sphinxAtStartPar
\sphinxstylestrong{Synopsis}
\begin{quote}

\sphinxAtStartPar
\sphinxcode{\sphinxupquote{void ReadCbf(string filename)}}
\end{quote}

\sphinxAtStartPar
\sphinxstylestrong{Arguments}
\begin{quote}

\sphinxAtStartPar
\sphinxcode{\sphinxupquote{filename}}: an input file name.
\end{quote}
\end{quote}

\subsubsection{Model.ReadJsonSol()}
\label{\detokenize{csapi/model:model-readjsonsol}}\begin{quote}

\sphinxAtStartPar
Read solution in format of JSON from file.

\sphinxAtStartPar
\sphinxstylestrong{Synopsis}
\begin{quote}

\sphinxAtStartPar
\sphinxcode{\sphinxupquote{void ReadJsonSol(string filename)}}
\end{quote}

\sphinxAtStartPar
\sphinxstylestrong{Arguments}
\begin{quote}

\sphinxAtStartPar
\sphinxcode{\sphinxupquote{filename}}: an input file name.
\end{quote}
\end{quote}

\subsubsection{Model.ReadLp()}
\label{\detokenize{csapi/model:model-readlp}}\begin{quote}

\sphinxAtStartPar
Read problem in LP format from file.

\sphinxAtStartPar
\sphinxstylestrong{Synopsis}
\begin{quote}

\sphinxAtStartPar
\sphinxcode{\sphinxupquote{void ReadLp(string filename)}}
\end{quote}

\sphinxAtStartPar
\sphinxstylestrong{Arguments}
\begin{quote}

\sphinxAtStartPar
\sphinxcode{\sphinxupquote{filename}}: an input file name.
\end{quote}
\end{quote}

\subsubsection{Model.ReadMps()}
\label{\detokenize{csapi/model:model-readmps}}\begin{quote}

\sphinxAtStartPar
Read problem in MPS format from file.

\sphinxAtStartPar
\sphinxstylestrong{Synopsis}
\begin{quote}

\sphinxAtStartPar
\sphinxcode{\sphinxupquote{void ReadMps(string filename)}}
\end{quote}

\sphinxAtStartPar
\sphinxstylestrong{Arguments}
\begin{quote}

\sphinxAtStartPar
\sphinxcode{\sphinxupquote{filename}}: an input file name.
\end{quote}
\end{quote}

\subsubsection{Model.ReadMst()}
\label{\detokenize{csapi/model:model-readmst}}\begin{quote}

\sphinxAtStartPar
Read MIP start information from file.

\sphinxAtStartPar
\sphinxstylestrong{Synopsis}
\begin{quote}

\sphinxAtStartPar
\sphinxcode{\sphinxupquote{void ReadMst(string filename)}}
\end{quote}

\sphinxAtStartPar
\sphinxstylestrong{Arguments}
\begin{quote}

\sphinxAtStartPar
\sphinxcode{\sphinxupquote{filename}}: an input file name.
\end{quote}
\end{quote}

\subsubsection{Model.ReadOrd()}
\label{\detokenize{csapi/model:model-readord}}\begin{quote}

\sphinxAtStartPar
Read branching order from file.

\sphinxAtStartPar
\sphinxstylestrong{Synopsis}
\begin{quote}

\sphinxAtStartPar
\sphinxcode{\sphinxupquote{void ReadOrd(string filename)}}
\end{quote}

\sphinxAtStartPar
\sphinxstylestrong{Arguments}
\begin{quote}

\sphinxAtStartPar
\sphinxcode{\sphinxupquote{filename}}: an input file name.
\end{quote}
\end{quote}

\subsubsection{Model.ReadParam()}
\label{\detokenize{csapi/model:model-readparam}}\begin{quote}

\sphinxAtStartPar
Read COPT parameters from file.

\sphinxAtStartPar
\sphinxstylestrong{Synopsis}
\begin{quote}

\sphinxAtStartPar
\sphinxcode{\sphinxupquote{void ReadParam(string filename)}}
\end{quote}

\sphinxAtStartPar
\sphinxstylestrong{Arguments}
\begin{quote}

\sphinxAtStartPar
\sphinxcode{\sphinxupquote{filename}}: an input file name.
\end{quote}
\end{quote}

\subsubsection{Model.ReadSdpa()}
\label{\detokenize{csapi/model:model-readsdpa}}\begin{quote}

\sphinxAtStartPar
Read problem in SDPA format from file.

\sphinxAtStartPar
\sphinxstylestrong{Synopsis}
\begin{quote}

\sphinxAtStartPar
\sphinxcode{\sphinxupquote{void ReadSdpa(string filename)}}
\end{quote}

\sphinxAtStartPar
\sphinxstylestrong{Arguments}
\begin{quote}

\sphinxAtStartPar
\sphinxcode{\sphinxupquote{filename}}: an input file name.
\end{quote}
\end{quote}

\subsubsection{Model.ReadSol()}
\label{\detokenize{csapi/model:model-readsol}}\begin{quote}

\sphinxAtStartPar
Read solution from file.

\sphinxAtStartPar
\sphinxstylestrong{Synopsis}
\begin{quote}

\sphinxAtStartPar
\sphinxcode{\sphinxupquote{void ReadSol(string filename)}}
\end{quote}

\sphinxAtStartPar
\sphinxstylestrong{Arguments}
\begin{quote}

\sphinxAtStartPar
\sphinxcode{\sphinxupquote{filename}}: an input file name.
\end{quote}
\end{quote}

\subsubsection{Model.ReadTune()}
\label{\detokenize{csapi/model:model-readtune}}\begin{quote}

\sphinxAtStartPar
Read tuning parameters from file.

\sphinxAtStartPar
\sphinxstylestrong{Synopsis}
\begin{quote}

\sphinxAtStartPar
\sphinxcode{\sphinxupquote{void ReadTune(string filename)}}
\end{quote}

\sphinxAtStartPar
\sphinxstylestrong{Arguments}
\begin{quote}

\sphinxAtStartPar
\sphinxcode{\sphinxupquote{filename}}: an input file name.
\end{quote}
\end{quote}

\subsubsection{Model.Remove()}
\label{\detokenize{csapi/model:model-remove}}\begin{quote}

\sphinxAtStartPar
Remove an array of variables from model.

\sphinxAtStartPar
\sphinxstylestrong{Synopsis}
\begin{quote}

\sphinxAtStartPar
\sphinxcode{\sphinxupquote{void Remove(Var{[}{]} vars)}}
\end{quote}

\sphinxAtStartPar
\sphinxstylestrong{Arguments}
\begin{quote}

\sphinxAtStartPar
\sphinxcode{\sphinxupquote{vars}}: a list of variables.
\end{quote}
\end{quote}

\subsubsection{Model.Remove()}
\label{\detokenize{csapi/model:id77}}\begin{quote}

\sphinxAtStartPar
Remove an array of variables from model.

\sphinxAtStartPar
\sphinxstylestrong{Synopsis}
\begin{quote}

\sphinxAtStartPar
\sphinxcode{\sphinxupquote{void Remove(VarArray vars)}}
\end{quote}

\sphinxAtStartPar
\sphinxstylestrong{Arguments}
\begin{quote}

\sphinxAtStartPar
\sphinxcode{\sphinxupquote{vars}}: array of variables.
\end{quote}
\end{quote}

\subsubsection{Model.Remove()}
\label{\detokenize{csapi/model:id78}}\begin{quote}

\sphinxAtStartPar
Remove a list of constraints from model.

\sphinxAtStartPar
\sphinxstylestrong{Synopsis}
\begin{quote}

\sphinxAtStartPar
\sphinxcode{\sphinxupquote{void Remove(Constraint{[}{]} constrs)}}
\end{quote}

\sphinxAtStartPar
\sphinxstylestrong{Arguments}
\begin{quote}

\sphinxAtStartPar
\sphinxcode{\sphinxupquote{constrs}}: a list of constraints.
\end{quote}
\end{quote}

\subsubsection{Model.Remove()}
\label{\detokenize{csapi/model:id79}}\begin{quote}

\sphinxAtStartPar
Remove a list of constraints from model.

\sphinxAtStartPar
\sphinxstylestrong{Synopsis}
\begin{quote}

\sphinxAtStartPar
\sphinxcode{\sphinxupquote{void Remove(ConstrArray constrs)}}
\end{quote}

\sphinxAtStartPar
\sphinxstylestrong{Arguments}
\begin{quote}

\sphinxAtStartPar
\sphinxcode{\sphinxupquote{constrs}}: an array of constraints.
\end{quote}
\end{quote}

\subsubsection{Model.Remove()}
\label{\detokenize{csapi/model:id80}}\begin{quote}

\sphinxAtStartPar
Remove an array of nonlinear constraints from model.

\sphinxAtStartPar
\sphinxstylestrong{Synopsis}
\begin{quote}

\sphinxAtStartPar
\sphinxcode{\sphinxupquote{void Remove(NlConstraint{[}{]} constrs)}}
\end{quote}

\sphinxAtStartPar
\sphinxstylestrong{Arguments}
\begin{quote}

\sphinxAtStartPar
\sphinxcode{\sphinxupquote{constrs}}: array of nonlinear constraints.
\end{quote}
\end{quote}

\subsubsection{Model.Remove()}
\label{\detokenize{csapi/model:id81}}\begin{quote}

\sphinxAtStartPar
Remove a list of nonlinear constraints from model.

\sphinxAtStartPar
\sphinxstylestrong{Synopsis}
\begin{quote}

\sphinxAtStartPar
\sphinxcode{\sphinxupquote{void Remove(NlConstrArray constrs)}}
\end{quote}

\sphinxAtStartPar
\sphinxstylestrong{Arguments}
\begin{quote}

\sphinxAtStartPar
\sphinxcode{\sphinxupquote{constrs}}: array object of nonlinear constraints.
\end{quote}
\end{quote}

\subsubsection{Model.Remove()}
\label{\detokenize{csapi/model:id82}}\begin{quote}

\sphinxAtStartPar
Remove a list of SOS constraints from model.

\sphinxAtStartPar
\sphinxstylestrong{Synopsis}
\begin{quote}

\sphinxAtStartPar
\sphinxcode{\sphinxupquote{void Remove(Sos{[}{]} soss)}}
\end{quote}

\sphinxAtStartPar
\sphinxstylestrong{Arguments}
\begin{quote}

\sphinxAtStartPar
\sphinxcode{\sphinxupquote{soss}}: a list of SOS constraints.
\end{quote}
\end{quote}

\subsubsection{Model.Remove()}
\label{\detokenize{csapi/model:id83}}\begin{quote}

\sphinxAtStartPar
Remove a list of SOS constraints from model.

\sphinxAtStartPar
\sphinxstylestrong{Synopsis}
\begin{quote}

\sphinxAtStartPar
\sphinxcode{\sphinxupquote{void Remove(SosArray soss)}}
\end{quote}

\sphinxAtStartPar
\sphinxstylestrong{Arguments}
\begin{quote}

\sphinxAtStartPar
\sphinxcode{\sphinxupquote{soss}}: an array of SOS constraints.
\end{quote}
\end{quote}

\subsubsection{Model.Remove()}
\label{\detokenize{csapi/model:id84}}\begin{quote}

\sphinxAtStartPar
Remove a list of cone constraints from model.

\sphinxAtStartPar
\sphinxstylestrong{Synopsis}
\begin{quote}

\sphinxAtStartPar
\sphinxcode{\sphinxupquote{void Remove(Cone{[}{]} cones)}}
\end{quote}

\sphinxAtStartPar
\sphinxstylestrong{Arguments}
\begin{quote}

\sphinxAtStartPar
\sphinxcode{\sphinxupquote{cones}}: a list of cone constraints.
\end{quote}
\end{quote}

\subsubsection{Model.Remove()}
\label{\detokenize{csapi/model:id85}}\begin{quote}

\sphinxAtStartPar
Remove a list of cone constraints from model.

\sphinxAtStartPar
\sphinxstylestrong{Synopsis}
\begin{quote}

\sphinxAtStartPar
\sphinxcode{\sphinxupquote{void Remove(ConeArray cones)}}
\end{quote}

\sphinxAtStartPar
\sphinxstylestrong{Arguments}
\begin{quote}

\sphinxAtStartPar
\sphinxcode{\sphinxupquote{cones}}: an array of cone constraints.
\end{quote}
\end{quote}

\subsubsection{Model.Remove()}
\label{\detokenize{csapi/model:id86}}\begin{quote}

\sphinxAtStartPar
Remove a list of exponential cone constraints from model.

\sphinxAtStartPar
\sphinxstylestrong{Synopsis}
\begin{quote}

\sphinxAtStartPar
\sphinxcode{\sphinxupquote{void Remove(ExpCone{[}{]} cones)}}
\end{quote}

\sphinxAtStartPar
\sphinxstylestrong{Arguments}
\begin{quote}

\sphinxAtStartPar
\sphinxcode{\sphinxupquote{cones}}: a list of exponential cone constraints.
\end{quote}
\end{quote}

\subsubsection{Model.Remove()}
\label{\detokenize{csapi/model:id87}}\begin{quote}

\sphinxAtStartPar
Remove a list of exponential cone constraints from model.

\sphinxAtStartPar
\sphinxstylestrong{Synopsis}
\begin{quote}

\sphinxAtStartPar
\sphinxcode{\sphinxupquote{void Remove(ExpConeArray cones)}}
\end{quote}

\sphinxAtStartPar
\sphinxstylestrong{Arguments}
\begin{quote}

\sphinxAtStartPar
\sphinxcode{\sphinxupquote{cones}}: an array of exponential cone constraints.
\end{quote}
\end{quote}

\subsubsection{Model.Remove()}
\label{\detokenize{csapi/model:id88}}\begin{quote}

\sphinxAtStartPar
Remove a list of affine cone constraints from model.

\sphinxAtStartPar
\sphinxstylestrong{Synopsis}
\begin{quote}

\sphinxAtStartPar
\sphinxcode{\sphinxupquote{void Remove(AffineCone{[}{]} cones)}}
\end{quote}

\sphinxAtStartPar
\sphinxstylestrong{Arguments}
\begin{quote}

\sphinxAtStartPar
\sphinxcode{\sphinxupquote{cones}}: a list of affine cone constraints.
\end{quote}
\end{quote}

\subsubsection{Model.Remove()}
\label{\detokenize{csapi/model:id89}}\begin{quote}

\sphinxAtStartPar
Remove an array of affine cone constraints from model.

\sphinxAtStartPar
\sphinxstylestrong{Synopsis}
\begin{quote}

\sphinxAtStartPar
\sphinxcode{\sphinxupquote{void Remove(AffineConeArray cones)}}
\end{quote}

\sphinxAtStartPar
\sphinxstylestrong{Arguments}
\begin{quote}

\sphinxAtStartPar
\sphinxcode{\sphinxupquote{cones}}: an array of affine cone constraints.
\end{quote}
\end{quote}

\subsubsection{Model.Remove()}
\label{\detokenize{csapi/model:id90}}\begin{quote}

\sphinxAtStartPar
Remove a list of gernal constraints from model.

\sphinxAtStartPar
\sphinxstylestrong{Synopsis}
\begin{quote}

\sphinxAtStartPar
\sphinxcode{\sphinxupquote{void Remove(GenConstr{[}{]} genConstrs)}}
\end{quote}

\sphinxAtStartPar
\sphinxstylestrong{Arguments}
\begin{quote}

\sphinxAtStartPar
\sphinxcode{\sphinxupquote{genConstrs}}: a list of general constraints.
\end{quote}
\end{quote}

\subsubsection{Model.Remove()}
\label{\detokenize{csapi/model:id91}}\begin{quote}

\sphinxAtStartPar
Remove a list of gernal constraints from model.

\sphinxAtStartPar
\sphinxstylestrong{Synopsis}
\begin{quote}

\sphinxAtStartPar
\sphinxcode{\sphinxupquote{void Remove(GenConstrArray genConstrs)}}
\end{quote}

\sphinxAtStartPar
\sphinxstylestrong{Arguments}
\begin{quote}

\sphinxAtStartPar
\sphinxcode{\sphinxupquote{genConstrs}}: an array of general constraints.
\end{quote}
\end{quote}

\subsubsection{Model.Remove()}
\label{\detokenize{csapi/model:id92}}\begin{quote}

\sphinxAtStartPar
Remove a list of quadratic constraints from model.

\sphinxAtStartPar
\sphinxstylestrong{Synopsis}
\begin{quote}

\sphinxAtStartPar
\sphinxcode{\sphinxupquote{void Remove(QConstraint{[}{]} qconstrs)}}
\end{quote}

\sphinxAtStartPar
\sphinxstylestrong{Arguments}
\begin{quote}

\sphinxAtStartPar
\sphinxcode{\sphinxupquote{qconstrs}}: an array of quadratic constraints.
\end{quote}
\end{quote}

\subsubsection{Model.Remove()}
\label{\detokenize{csapi/model:id93}}\begin{quote}

\sphinxAtStartPar
Remove a list of quadratic constraints from model.

\sphinxAtStartPar
\sphinxstylestrong{Synopsis}
\begin{quote}

\sphinxAtStartPar
\sphinxcode{\sphinxupquote{void Remove(QConstrArray qconstrs)}}
\end{quote}

\sphinxAtStartPar
\sphinxstylestrong{Arguments}
\begin{quote}

\sphinxAtStartPar
\sphinxcode{\sphinxupquote{qconstrs}}: an array of quadratic constraints.
\end{quote}
\end{quote}

\subsubsection{Model.Remove()}
\label{\detokenize{csapi/model:id94}}\begin{quote}

\sphinxAtStartPar
Remove a list of PSD variables from model.

\sphinxAtStartPar
\sphinxstylestrong{Synopsis}
\begin{quote}

\sphinxAtStartPar
\sphinxcode{\sphinxupquote{void Remove(PsdVar{[}{]} vars)}}
\end{quote}

\sphinxAtStartPar
\sphinxstylestrong{Arguments}
\begin{quote}

\sphinxAtStartPar
\sphinxcode{\sphinxupquote{vars}}: an array of PSD variables.
\end{quote}
\end{quote}

\subsubsection{Model.Remove()}
\label{\detokenize{csapi/model:id95}}\begin{quote}

\sphinxAtStartPar
Remove a list of PSD variables from model.

\sphinxAtStartPar
\sphinxstylestrong{Synopsis}
\begin{quote}

\sphinxAtStartPar
\sphinxcode{\sphinxupquote{void Remove(PsdVarArray vars)}}
\end{quote}

\sphinxAtStartPar
\sphinxstylestrong{Arguments}
\begin{quote}

\sphinxAtStartPar
\sphinxcode{\sphinxupquote{vars}}: an array of PSD variables.
\end{quote}
\end{quote}

\subsubsection{Model.Remove()}
\label{\detokenize{csapi/model:id96}}\begin{quote}

\sphinxAtStartPar
Remove a list of PSD constraints from model.

\sphinxAtStartPar
\sphinxstylestrong{Synopsis}
\begin{quote}

\sphinxAtStartPar
\sphinxcode{\sphinxupquote{void Remove(PsdConstraint{[}{]} constrs)}}
\end{quote}

\sphinxAtStartPar
\sphinxstylestrong{Arguments}
\begin{quote}

\sphinxAtStartPar
\sphinxcode{\sphinxupquote{constrs}}: an array of PSD constraints.
\end{quote}
\end{quote}

\subsubsection{Model.Remove()}
\label{\detokenize{csapi/model:id97}}\begin{quote}

\sphinxAtStartPar
Remove a list of PSD constraints from model.

\sphinxAtStartPar
\sphinxstylestrong{Synopsis}
\begin{quote}

\sphinxAtStartPar
\sphinxcode{\sphinxupquote{void Remove(PsdConstrArray constrs)}}
\end{quote}

\sphinxAtStartPar
\sphinxstylestrong{Arguments}
\begin{quote}

\sphinxAtStartPar
\sphinxcode{\sphinxupquote{constrs}}: an array of PSD constraints.
\end{quote}
\end{quote}

\subsubsection{Model.Remove()}
\label{\detokenize{csapi/model:id98}}\begin{quote}

\sphinxAtStartPar
Remove a list of LMI constraints from model.

\sphinxAtStartPar
\sphinxstylestrong{Synopsis}
\begin{quote}

\sphinxAtStartPar
\sphinxcode{\sphinxupquote{void Remove(LmiConstrArray constrs)}}
\end{quote}

\sphinxAtStartPar
\sphinxstylestrong{Arguments}
\begin{quote}

\sphinxAtStartPar
\sphinxcode{\sphinxupquote{constrs}}: an array of LMI constraints.
\end{quote}
\end{quote}

\subsubsection{Model.Remove()}
\label{\detokenize{csapi/model:id99}}\begin{quote}

\sphinxAtStartPar
Remove a list of LMI constraints from model.

\sphinxAtStartPar
\sphinxstylestrong{Synopsis}
\begin{quote}

\sphinxAtStartPar
\sphinxcode{\sphinxupquote{void Remove(LmiConstraint{[}{]} constrs)}}
\end{quote}

\sphinxAtStartPar
\sphinxstylestrong{Arguments}
\begin{quote}

\sphinxAtStartPar
\sphinxcode{\sphinxupquote{constrs}}: an array of LMI constraints.
\end{quote}
\end{quote}

\subsubsection{Model.Reset()}
\label{\detokenize{csapi/model:model-reset}}\begin{quote}

\sphinxAtStartPar
Reset solution of problem only.

\sphinxAtStartPar
\sphinxstylestrong{Synopsis}
\begin{quote}

\sphinxAtStartPar
\sphinxcode{\sphinxupquote{void Reset()}}
\end{quote}
\end{quote}

\subsubsection{Model.ResetAll()}
\label{\detokenize{csapi/model:model-resetall}}\begin{quote}

\sphinxAtStartPar
Reset solution in problem, and additional information such as MIP start, etc.

\sphinxAtStartPar
\sphinxstylestrong{Synopsis}
\begin{quote}

\sphinxAtStartPar
\sphinxcode{\sphinxupquote{void ResetAll()}}
\end{quote}
\end{quote}

\subsubsection{Model.ResetObjParamN()}
\label{\detokenize{csapi/model:model-resetobjparamn}}\begin{quote}

\sphinxAtStartPar
Reset objective parameters of a multi\sphinxhyphen{}objective function.

\sphinxAtStartPar
\sphinxstylestrong{Synopsis}
\begin{quote}

\sphinxAtStartPar
\sphinxcode{\sphinxupquote{void ResetObjParamN(int idx)}}
\end{quote}

\sphinxAtStartPar
\sphinxstylestrong{Arguments}
\begin{quote}

\sphinxAtStartPar
\sphinxcode{\sphinxupquote{idx}}: index of a multi\sphinxhyphen{}objective function.
\end{quote}
\end{quote}

\subsubsection{Model.ResetParam()}
\label{\detokenize{csapi/model:model-resetparam}}\begin{quote}

\sphinxAtStartPar
Reset parameters to default settings.

\sphinxAtStartPar
\sphinxstylestrong{Synopsis}
\begin{quote}

\sphinxAtStartPar
\sphinxcode{\sphinxupquote{void ResetParam()}}
\end{quote}
\end{quote}

\subsubsection{Model.ResetParamN()}
\label{\detokenize{csapi/model:model-resetparamn}}\begin{quote}

\sphinxAtStartPar
Reset double and integer parameters of a multi\sphinxhyphen{}objective function.

\sphinxAtStartPar
\sphinxstylestrong{Synopsis}
\begin{quote}

\sphinxAtStartPar
\sphinxcode{\sphinxupquote{void ResetParamN(int idx)}}
\end{quote}

\sphinxAtStartPar
\sphinxstylestrong{Arguments}
\begin{quote}

\sphinxAtStartPar
\sphinxcode{\sphinxupquote{idx}}: index of a multi\sphinxhyphen{}objective function.
\end{quote}
\end{quote}

\subsubsection{Model.Set()}
\label{\detokenize{csapi/model:model-set}}\begin{quote}

\sphinxAtStartPar
Set values of information associated with variables.

\sphinxAtStartPar
\sphinxstylestrong{Synopsis}
\begin{quote}

\sphinxAtStartPar
\sphinxcode{\sphinxupquote{void Set(}}
\begin{quote}

\sphinxAtStartPar
\sphinxcode{\sphinxupquote{string name,}}

\sphinxAtStartPar
\sphinxcode{\sphinxupquote{Var{[}{]} vars,}}

\sphinxAtStartPar
\sphinxcode{\sphinxupquote{double{[}{]} vals)}}
\end{quote}
\end{quote}

\sphinxAtStartPar
\sphinxstylestrong{Arguments}
\begin{quote}

\sphinxAtStartPar
\sphinxcode{\sphinxupquote{name}}: name of information.

\sphinxAtStartPar
\sphinxcode{\sphinxupquote{vars}}: a list of interested variables.

\sphinxAtStartPar
\sphinxcode{\sphinxupquote{vals}}: values of information.
\end{quote}
\end{quote}

\subsubsection{Model.Set()}
\label{\detokenize{csapi/model:id100}}\begin{quote}

\sphinxAtStartPar
Set values of information associated with variables.

\sphinxAtStartPar
\sphinxstylestrong{Synopsis}
\begin{quote}

\sphinxAtStartPar
\sphinxcode{\sphinxupquote{void Set(}}
\begin{quote}

\sphinxAtStartPar
\sphinxcode{\sphinxupquote{string name,}}

\sphinxAtStartPar
\sphinxcode{\sphinxupquote{VarArray vars,}}

\sphinxAtStartPar
\sphinxcode{\sphinxupquote{double{[}{]} vals)}}
\end{quote}
\end{quote}

\sphinxAtStartPar
\sphinxstylestrong{Arguments}
\begin{quote}

\sphinxAtStartPar
\sphinxcode{\sphinxupquote{name}}: name of information.

\sphinxAtStartPar
\sphinxcode{\sphinxupquote{vars}}: array of interested variables.

\sphinxAtStartPar
\sphinxcode{\sphinxupquote{vals}}: values of information.
\end{quote}
\end{quote}

\subsubsection{Model.Set()}
\label{\detokenize{csapi/model:id101}}\begin{quote}

\sphinxAtStartPar
Set values of information associated with constraints.

\sphinxAtStartPar
\sphinxstylestrong{Synopsis}
\begin{quote}

\sphinxAtStartPar
\sphinxcode{\sphinxupquote{void Set(}}
\begin{quote}

\sphinxAtStartPar
\sphinxcode{\sphinxupquote{string name,}}

\sphinxAtStartPar
\sphinxcode{\sphinxupquote{Constraint{[}{]} constrs,}}

\sphinxAtStartPar
\sphinxcode{\sphinxupquote{double{[}{]} vals)}}
\end{quote}
\end{quote}

\sphinxAtStartPar
\sphinxstylestrong{Arguments}
\begin{quote}

\sphinxAtStartPar
\sphinxcode{\sphinxupquote{name}}: name of information.

\sphinxAtStartPar
\sphinxcode{\sphinxupquote{constrs}}: a list of interested constraints.

\sphinxAtStartPar
\sphinxcode{\sphinxupquote{vals}}: values of information.
\end{quote}
\end{quote}

\subsubsection{Model.Set()}
\label{\detokenize{csapi/model:id102}}\begin{quote}

\sphinxAtStartPar
Set values of information associated with constraints.

\sphinxAtStartPar
\sphinxstylestrong{Synopsis}
\begin{quote}

\sphinxAtStartPar
\sphinxcode{\sphinxupquote{void Set(}}
\begin{quote}

\sphinxAtStartPar
\sphinxcode{\sphinxupquote{string name,}}

\sphinxAtStartPar
\sphinxcode{\sphinxupquote{ConstrArray constrs,}}

\sphinxAtStartPar
\sphinxcode{\sphinxupquote{double{[}{]} vals)}}
\end{quote}
\end{quote}

\sphinxAtStartPar
\sphinxstylestrong{Arguments}
\begin{quote}

\sphinxAtStartPar
\sphinxcode{\sphinxupquote{name}}: name of information.

\sphinxAtStartPar
\sphinxcode{\sphinxupquote{constrs}}: array of interested constraints.

\sphinxAtStartPar
\sphinxcode{\sphinxupquote{vals}}: values of information.
\end{quote}
\end{quote}

\subsubsection{Model.Set()}
\label{\detokenize{csapi/model:id103}}\begin{quote}

\sphinxAtStartPar
Set values of information associated with nonlinear constraints.

\sphinxAtStartPar
\sphinxstylestrong{Synopsis}
\begin{quote}

\sphinxAtStartPar
\sphinxcode{\sphinxupquote{void Set(}}
\begin{quote}

\sphinxAtStartPar
\sphinxcode{\sphinxupquote{string name,}}

\sphinxAtStartPar
\sphinxcode{\sphinxupquote{NlConstraint{[}{]} constrs,}}

\sphinxAtStartPar
\sphinxcode{\sphinxupquote{double{[}{]} vals)}}
\end{quote}
\end{quote}

\sphinxAtStartPar
\sphinxstylestrong{Arguments}
\begin{quote}

\sphinxAtStartPar
\sphinxcode{\sphinxupquote{name}}: name of double information.

\sphinxAtStartPar
\sphinxcode{\sphinxupquote{constrs}}: array of desired nonlinear constraints.

\sphinxAtStartPar
\sphinxcode{\sphinxupquote{vals}}: array of values of information.
\end{quote}
\end{quote}

\subsubsection{Model.Set()}
\label{\detokenize{csapi/model:id104}}\begin{quote}

\sphinxAtStartPar
Set values of information associated with nonlinear constraints.

\sphinxAtStartPar
\sphinxstylestrong{Synopsis}
\begin{quote}

\sphinxAtStartPar
\sphinxcode{\sphinxupquote{void Set(}}
\begin{quote}

\sphinxAtStartPar
\sphinxcode{\sphinxupquote{string name,}}

\sphinxAtStartPar
\sphinxcode{\sphinxupquote{NlConstrArray constrs,}}

\sphinxAtStartPar
\sphinxcode{\sphinxupquote{double{[}{]} vals)}}
\end{quote}
\end{quote}

\sphinxAtStartPar
\sphinxstylestrong{Arguments}
\begin{quote}

\sphinxAtStartPar
\sphinxcode{\sphinxupquote{name}}: name of double information.

\sphinxAtStartPar
\sphinxcode{\sphinxupquote{constrs}}: a list of desired nonlinear constraints.

\sphinxAtStartPar
\sphinxcode{\sphinxupquote{vals}}: array of values of information.
\end{quote}
\end{quote}

\subsubsection{Model.Set()}
\label{\detokenize{csapi/model:id105}}\begin{quote}

\sphinxAtStartPar
Set values of information associated with PSD constraints.

\sphinxAtStartPar
\sphinxstylestrong{Synopsis}
\begin{quote}

\sphinxAtStartPar
\sphinxcode{\sphinxupquote{void Set(}}
\begin{quote}

\sphinxAtStartPar
\sphinxcode{\sphinxupquote{string name,}}

\sphinxAtStartPar
\sphinxcode{\sphinxupquote{PsdConstraint{[}{]} constrs,}}

\sphinxAtStartPar
\sphinxcode{\sphinxupquote{double{[}{]} vals)}}
\end{quote}
\end{quote}

\sphinxAtStartPar
\sphinxstylestrong{Arguments}
\begin{quote}

\sphinxAtStartPar
\sphinxcode{\sphinxupquote{name}}: name of information.

\sphinxAtStartPar
\sphinxcode{\sphinxupquote{constrs}}: a list of desired PSD constraints.

\sphinxAtStartPar
\sphinxcode{\sphinxupquote{vals}}: array of values of information.
\end{quote}
\end{quote}

\subsubsection{Model.Set()}
\label{\detokenize{csapi/model:id106}}\begin{quote}

\sphinxAtStartPar
Set values of information associated with PSD constraints.

\sphinxAtStartPar
\sphinxstylestrong{Synopsis}
\begin{quote}

\sphinxAtStartPar
\sphinxcode{\sphinxupquote{void Set(}}
\begin{quote}

\sphinxAtStartPar
\sphinxcode{\sphinxupquote{string name,}}

\sphinxAtStartPar
\sphinxcode{\sphinxupquote{PsdConstrArray constrs,}}

\sphinxAtStartPar
\sphinxcode{\sphinxupquote{double{[}{]} vals)}}
\end{quote}
\end{quote}

\sphinxAtStartPar
\sphinxstylestrong{Arguments}
\begin{quote}

\sphinxAtStartPar
\sphinxcode{\sphinxupquote{name}}: name of information.

\sphinxAtStartPar
\sphinxcode{\sphinxupquote{constrs}}: a list of desired PSD constraints.

\sphinxAtStartPar
\sphinxcode{\sphinxupquote{vals}}: array of values of information.
\end{quote}
\end{quote}

\subsubsection{Model.SetBasis()}
\label{\detokenize{csapi/model:model-setbasis}}\begin{quote}

\sphinxAtStartPar
Set column and row basis status to model.

\sphinxAtStartPar
\sphinxstylestrong{Synopsis}
\begin{quote}

\sphinxAtStartPar
\sphinxcode{\sphinxupquote{void SetBasis(int{[}{]} colbasis, int{[}{]} rowbasis)}}
\end{quote}

\sphinxAtStartPar
\sphinxstylestrong{Arguments}
\begin{quote}

\sphinxAtStartPar
\sphinxcode{\sphinxupquote{colbasis}}: status of column basis.

\sphinxAtStartPar
\sphinxcode{\sphinxupquote{rowbasis}}: status of row basis.
\end{quote}
\end{quote}

\subsubsection{Model.SetCallback()}
\label{\detokenize{csapi/model:model-setcallback}}\begin{quote}

\sphinxAtStartPar
Set user callback to COPT model.

\sphinxAtStartPar
\sphinxstylestrong{Synopsis}
\begin{quote}

\sphinxAtStartPar
\sphinxcode{\sphinxupquote{void SetCallback(CallbackBase cb, int cbctx)}}
\end{quote}

\sphinxAtStartPar
\sphinxstylestrong{Arguments}
\begin{quote}

\sphinxAtStartPar
\sphinxcode{\sphinxupquote{cb}}: user callback instance, inheriting from CallbackBase class.

\sphinxAtStartPar
\sphinxcode{\sphinxupquote{cbctx}}: COPT callback context.
\end{quote}
\end{quote}

\subsubsection{Model.SetCoeff()}
\label{\detokenize{csapi/model:model-setcoeff}}\begin{quote}

\sphinxAtStartPar
Set the coefficient of a variable in a linear constraint.

\sphinxAtStartPar
\sphinxstylestrong{Synopsis}
\begin{quote}

\sphinxAtStartPar
\sphinxcode{\sphinxupquote{void SetCoeff(}}
\begin{quote}

\sphinxAtStartPar
\sphinxcode{\sphinxupquote{Constraint constr,}}

\sphinxAtStartPar
\sphinxcode{\sphinxupquote{Var var,}}

\sphinxAtStartPar
\sphinxcode{\sphinxupquote{double newVal)}}
\end{quote}
\end{quote}

\sphinxAtStartPar
\sphinxstylestrong{Arguments}
\begin{quote}

\sphinxAtStartPar
\sphinxcode{\sphinxupquote{constr}}: The requested constraint.

\sphinxAtStartPar
\sphinxcode{\sphinxupquote{var}}: The requested variable.

\sphinxAtStartPar
\sphinxcode{\sphinxupquote{newVal}}: New coefficient.
\end{quote}
\end{quote}

\subsubsection{Model.SetCoeffs()}
\label{\detokenize{csapi/model:model-setcoeffs}}\begin{quote}

\sphinxAtStartPar
Set a list of coefficients in the model.

\sphinxAtStartPar
\sphinxstylestrong{Synopsis}
\begin{quote}

\sphinxAtStartPar
\sphinxcode{\sphinxupquote{void SetCoeffs(}}
\begin{quote}

\sphinxAtStartPar
\sphinxcode{\sphinxupquote{Constraint{[}{]} constrs,}}

\sphinxAtStartPar
\sphinxcode{\sphinxupquote{Var{[}{]} vars,}}

\sphinxAtStartPar
\sphinxcode{\sphinxupquote{double{[}{]} vals)}}
\end{quote}
\end{quote}

\sphinxAtStartPar
\sphinxstylestrong{Arguments}
\begin{quote}

\sphinxAtStartPar
\sphinxcode{\sphinxupquote{constrs}}: Array of constraints for coefficients to be set.

\sphinxAtStartPar
\sphinxcode{\sphinxupquote{vars}}: Array of vars for coefficients to be set.

\sphinxAtStartPar
\sphinxcode{\sphinxupquote{vals}}: New values for coefficients.
\end{quote}
\end{quote}

\subsubsection{Model.SetCoeffs()}
\label{\detokenize{csapi/model:id107}}\begin{quote}

\sphinxAtStartPar
Set a list of coefficients in the model.

\sphinxAtStartPar
\sphinxstylestrong{Synopsis}
\begin{quote}

\sphinxAtStartPar
\sphinxcode{\sphinxupquote{void SetCoeffs(}}
\begin{quote}

\sphinxAtStartPar
\sphinxcode{\sphinxupquote{ConstrArray constrs,}}

\sphinxAtStartPar
\sphinxcode{\sphinxupquote{VarArray vars,}}

\sphinxAtStartPar
\sphinxcode{\sphinxupquote{double{[}{]} vals)}}
\end{quote}
\end{quote}

\sphinxAtStartPar
\sphinxstylestrong{Arguments}
\begin{quote}

\sphinxAtStartPar
\sphinxcode{\sphinxupquote{constrs}}: A list of constraints for coefficients to be set.

\sphinxAtStartPar
\sphinxcode{\sphinxupquote{vars}}: A list of vars for coefficients to be set.

\sphinxAtStartPar
\sphinxcode{\sphinxupquote{vals}}: New values for coefficients.
\end{quote}
\end{quote}

\subsubsection{Model.SetDblParam()}
\label{\detokenize{csapi/model:model-setdblparam}}\begin{quote}

\sphinxAtStartPar
Set value of a COPT double parameter.

\sphinxAtStartPar
\sphinxstylestrong{Synopsis}
\begin{quote}

\sphinxAtStartPar
\sphinxcode{\sphinxupquote{void SetDblParam(string param, double val)}}
\end{quote}

\sphinxAtStartPar
\sphinxstylestrong{Arguments}
\begin{quote}

\sphinxAtStartPar
\sphinxcode{\sphinxupquote{param}}: name of integer parameter.

\sphinxAtStartPar
\sphinxcode{\sphinxupquote{val}}: double value.
\end{quote}
\end{quote}

\subsubsection{Model.SetDblParamN()}
\label{\detokenize{csapi/model:model-setdblparamn}}\begin{quote}

\sphinxAtStartPar
Set value of a double parameter of a multi\sphinxhyphen{}objective function.

\sphinxAtStartPar
\sphinxstylestrong{Synopsis}
\begin{quote}

\sphinxAtStartPar
\sphinxcode{\sphinxupquote{void SetDblParamN(}}
\begin{quote}

\sphinxAtStartPar
\sphinxcode{\sphinxupquote{int idx,}}

\sphinxAtStartPar
\sphinxcode{\sphinxupquote{string param,}}

\sphinxAtStartPar
\sphinxcode{\sphinxupquote{double val)}}
\end{quote}
\end{quote}

\sphinxAtStartPar
\sphinxstylestrong{Arguments}
\begin{quote}

\sphinxAtStartPar
\sphinxcode{\sphinxupquote{idx}}: index of a multi\sphinxhyphen{}objective function.

\sphinxAtStartPar
\sphinxcode{\sphinxupquote{param}}: name of double parameter.

\sphinxAtStartPar
\sphinxcode{\sphinxupquote{val}}: new value of double parameter.
\end{quote}
\end{quote}

\subsubsection{Model.SetIntParam()}
\label{\detokenize{csapi/model:model-setintparam}}\begin{quote}

\sphinxAtStartPar
Set value of a COPT integer parameter.

\sphinxAtStartPar
\sphinxstylestrong{Synopsis}
\begin{quote}

\sphinxAtStartPar
\sphinxcode{\sphinxupquote{void SetIntParam(string param, int val)}}
\end{quote}

\sphinxAtStartPar
\sphinxstylestrong{Arguments}
\begin{quote}

\sphinxAtStartPar
\sphinxcode{\sphinxupquote{param}}: name of integer parameter.

\sphinxAtStartPar
\sphinxcode{\sphinxupquote{val}}: integer value.
\end{quote}
\end{quote}

\subsubsection{Model.SetIntParamN()}
\label{\detokenize{csapi/model:model-setintparamn}}\begin{quote}

\sphinxAtStartPar
Set value of an integer parameter of a multi\sphinxhyphen{}objective function.

\sphinxAtStartPar
\sphinxstylestrong{Synopsis}
\begin{quote}

\sphinxAtStartPar
\sphinxcode{\sphinxupquote{void SetIntParamN(}}
\begin{quote}

\sphinxAtStartPar
\sphinxcode{\sphinxupquote{int idx,}}

\sphinxAtStartPar
\sphinxcode{\sphinxupquote{string param,}}

\sphinxAtStartPar
\sphinxcode{\sphinxupquote{int val)}}
\end{quote}
\end{quote}

\sphinxAtStartPar
\sphinxstylestrong{Arguments}
\begin{quote}

\sphinxAtStartPar
\sphinxcode{\sphinxupquote{idx}}: index of a multi\sphinxhyphen{}objective function.

\sphinxAtStartPar
\sphinxcode{\sphinxupquote{param}}: name of integer parameter.

\sphinxAtStartPar
\sphinxcode{\sphinxupquote{val}}: new value of integer parameter.
\end{quote}
\end{quote}

\subsubsection{Model.SetLmiCoeff()}
\label{\detokenize{csapi/model:model-setlmicoeff}}\begin{quote}

\sphinxAtStartPar
Set the coefficient matrix of a variable in LMI constraint.

\sphinxAtStartPar
\sphinxstylestrong{Synopsis}
\begin{quote}

\sphinxAtStartPar
\sphinxcode{\sphinxupquote{void SetLmiCoeff(}}
\begin{quote}

\sphinxAtStartPar
\sphinxcode{\sphinxupquote{LmiConstraint constr,}}

\sphinxAtStartPar
\sphinxcode{\sphinxupquote{Var var,}}

\sphinxAtStartPar
\sphinxcode{\sphinxupquote{SymMatrix mat)}}
\end{quote}
\end{quote}

\sphinxAtStartPar
\sphinxstylestrong{Arguments}
\begin{quote}

\sphinxAtStartPar
\sphinxcode{\sphinxupquote{constr}}: The desired LMI constraint.

\sphinxAtStartPar
\sphinxcode{\sphinxupquote{var}}: The desired variable.

\sphinxAtStartPar
\sphinxcode{\sphinxupquote{mat}}: new coefficient matrix.
\end{quote}
\end{quote}

\subsubsection{Model.SetLmiRhs()}
\label{\detokenize{csapi/model:model-setlmirhs}}\begin{quote}

\sphinxAtStartPar
Set constant matrix of LMI constraint.

\sphinxAtStartPar
\sphinxstylestrong{Synopsis}
\begin{quote}

\sphinxAtStartPar
\sphinxcode{\sphinxupquote{void SetLmiRhs(LmiConstraint constr, SymMatrix mat)}}
\end{quote}

\sphinxAtStartPar
\sphinxstylestrong{Arguments}
\begin{quote}

\sphinxAtStartPar
\sphinxcode{\sphinxupquote{constr}}: The desired LMI constraint.

\sphinxAtStartPar
\sphinxcode{\sphinxupquote{mat}}: new constant matrix.
\end{quote}
\end{quote}

\subsubsection{Model.SetLpSolution()}
\label{\detokenize{csapi/model:model-setlpsolution}}\begin{quote}

\sphinxAtStartPar
Set LP solution.

\sphinxAtStartPar
\sphinxstylestrong{Synopsis}
\begin{quote}

\sphinxAtStartPar
\sphinxcode{\sphinxupquote{void SetLpSolution(}}
\begin{quote}

\sphinxAtStartPar
\sphinxcode{\sphinxupquote{double{[}{]} value,}}

\sphinxAtStartPar
\sphinxcode{\sphinxupquote{double{[}{]} slack,}}

\sphinxAtStartPar
\sphinxcode{\sphinxupquote{double{[}{]} rowDual,}}

\sphinxAtStartPar
\sphinxcode{\sphinxupquote{double{[}{]} redCost)}}
\end{quote}
\end{quote}

\sphinxAtStartPar
\sphinxstylestrong{Arguments}
\begin{quote}

\sphinxAtStartPar
\sphinxcode{\sphinxupquote{value}}: solution values.

\sphinxAtStartPar
\sphinxcode{\sphinxupquote{slack}}: slack values.

\sphinxAtStartPar
\sphinxcode{\sphinxupquote{rowDual}}: dual values.

\sphinxAtStartPar
\sphinxcode{\sphinxupquote{redCost}}: reduced costs.
\end{quote}
\end{quote}

\subsubsection{Model.SetMipStart()}
\label{\detokenize{csapi/model:model-setmipstart}}\begin{quote}

\sphinxAtStartPar
Set initial values for variables of given number, starting from the first one.

\sphinxAtStartPar
\sphinxstylestrong{Synopsis}
\begin{quote}

\sphinxAtStartPar
\sphinxcode{\sphinxupquote{void SetMipStart(int count, double{[}{]} vals)}}
\end{quote}

\sphinxAtStartPar
\sphinxstylestrong{Arguments}
\begin{quote}

\sphinxAtStartPar
\sphinxcode{\sphinxupquote{count}}: the number of variables to set.

\sphinxAtStartPar
\sphinxcode{\sphinxupquote{vals}}: values of variables.
\end{quote}
\end{quote}

\subsubsection{Model.SetMipStart()}
\label{\detokenize{csapi/model:id108}}\begin{quote}

\sphinxAtStartPar
Set initial value for the specified variable.

\sphinxAtStartPar
\sphinxstylestrong{Synopsis}
\begin{quote}

\sphinxAtStartPar
\sphinxcode{\sphinxupquote{void SetMipStart(Var var, double val)}}
\end{quote}

\sphinxAtStartPar
\sphinxstylestrong{Arguments}
\begin{quote}

\sphinxAtStartPar
\sphinxcode{\sphinxupquote{var}}: an interested variable.

\sphinxAtStartPar
\sphinxcode{\sphinxupquote{val}}: initial value of the variable.
\end{quote}
\end{quote}

\subsubsection{Model.SetMipStart()}
\label{\detokenize{csapi/model:id109}}\begin{quote}

\sphinxAtStartPar
Set initial value for the specified variable.

\sphinxAtStartPar
\sphinxstylestrong{Synopsis}
\begin{quote}

\sphinxAtStartPar
\sphinxcode{\sphinxupquote{void SetMipStart(Var{[}{]} vars, double{[}{]} vals)}}
\end{quote}

\sphinxAtStartPar
\sphinxstylestrong{Arguments}
\begin{quote}

\sphinxAtStartPar
\sphinxcode{\sphinxupquote{vars}}: a list of interested variables.

\sphinxAtStartPar
\sphinxcode{\sphinxupquote{vals}}: initial values of the variables.
\end{quote}
\end{quote}

\subsubsection{Model.SetMipStart()}
\label{\detokenize{csapi/model:id110}}\begin{quote}

\sphinxAtStartPar
Set initial value for the specified variable.

\sphinxAtStartPar
\sphinxstylestrong{Synopsis}
\begin{quote}

\sphinxAtStartPar
\sphinxcode{\sphinxupquote{void SetMipStart(VarArray vars, double{[}{]} vals)}}
\end{quote}

\sphinxAtStartPar
\sphinxstylestrong{Arguments}
\begin{quote}

\sphinxAtStartPar
\sphinxcode{\sphinxupquote{vars}}: a list of interested variables.

\sphinxAtStartPar
\sphinxcode{\sphinxupquote{vals}}: initial values of the variables.
\end{quote}
\end{quote}

\subsubsection{Model.SetNames()}
\label{\detokenize{csapi/model:model-setnames}}\begin{quote}

\sphinxAtStartPar
Set names for given variables in model.

\sphinxAtStartPar
\sphinxstylestrong{Synopsis}
\begin{quote}

\sphinxAtStartPar
\sphinxcode{\sphinxupquote{void SetNames(Var{[}{]} vars, string{[}{]} names)}}
\end{quote}

\sphinxAtStartPar
\sphinxstylestrong{Arguments}
\begin{quote}

\sphinxAtStartPar
\sphinxcode{\sphinxupquote{vars}}: array of variables.

\sphinxAtStartPar
\sphinxcode{\sphinxupquote{names}}: string array of names for variables.
\end{quote}
\end{quote}

\subsubsection{Model.SetNames()}
\label{\detokenize{csapi/model:id111}}\begin{quote}

\sphinxAtStartPar
Set names for given variables in model.

\sphinxAtStartPar
\sphinxstylestrong{Synopsis}
\begin{quote}

\sphinxAtStartPar
\sphinxcode{\sphinxupquote{void SetNames(VarArray vars, string{[}{]} names)}}
\end{quote}

\sphinxAtStartPar
\sphinxstylestrong{Arguments}
\begin{quote}

\sphinxAtStartPar
\sphinxcode{\sphinxupquote{vars}}: a list of variables.

\sphinxAtStartPar
\sphinxcode{\sphinxupquote{names}}: string array of names for variables.
\end{quote}
\end{quote}

\subsubsection{Model.SetNames()}
\label{\detokenize{csapi/model:id112}}\begin{quote}

\sphinxAtStartPar
Set names for given constraints in model.

\sphinxAtStartPar
\sphinxstylestrong{Synopsis}
\begin{quote}

\sphinxAtStartPar
\sphinxcode{\sphinxupquote{void SetNames(Constraint{[}{]} cons, string{[}{]} names)}}
\end{quote}

\sphinxAtStartPar
\sphinxstylestrong{Arguments}
\begin{quote}

\sphinxAtStartPar
\sphinxcode{\sphinxupquote{cons}}: array of constraints.

\sphinxAtStartPar
\sphinxcode{\sphinxupquote{names}}: string array of names for constraints.
\end{quote}
\end{quote}

\subsubsection{Model.SetNames()}
\label{\detokenize{csapi/model:id113}}\begin{quote}

\sphinxAtStartPar
Set names for given constraints in model.

\sphinxAtStartPar
\sphinxstylestrong{Synopsis}
\begin{quote}

\sphinxAtStartPar
\sphinxcode{\sphinxupquote{void SetNames(ConstrArray cons, string{[}{]} names)}}
\end{quote}

\sphinxAtStartPar
\sphinxstylestrong{Arguments}
\begin{quote}

\sphinxAtStartPar
\sphinxcode{\sphinxupquote{cons}}: a list of constraints.

\sphinxAtStartPar
\sphinxcode{\sphinxupquote{names}}: string array of names for constraints.
\end{quote}
\end{quote}

\subsubsection{Model.SetNames()}
\label{\detokenize{csapi/model:id114}}\begin{quote}

\sphinxAtStartPar
Set names for given general constraints in model.

\sphinxAtStartPar
\sphinxstylestrong{Synopsis}
\begin{quote}

\sphinxAtStartPar
\sphinxcode{\sphinxupquote{void SetNames(GenConstr{[}{]} genConstrs, string{[}{]} names)}}
\end{quote}

\sphinxAtStartPar
\sphinxstylestrong{Arguments}
\begin{quote}

\sphinxAtStartPar
\sphinxcode{\sphinxupquote{genConstrs}}: array of general constraints.

\sphinxAtStartPar
\sphinxcode{\sphinxupquote{names}}: string array of names for general constraints.
\end{quote}
\end{quote}

\subsubsection{Model.SetNames()}
\label{\detokenize{csapi/model:id115}}\begin{quote}

\sphinxAtStartPar
Set names for given general constraints in model.

\sphinxAtStartPar
\sphinxstylestrong{Synopsis}
\begin{quote}

\sphinxAtStartPar
\sphinxcode{\sphinxupquote{void SetNames(GenConstrArray genConstrs, string{[}{]} names)}}
\end{quote}

\sphinxAtStartPar
\sphinxstylestrong{Arguments}
\begin{quote}

\sphinxAtStartPar
\sphinxcode{\sphinxupquote{genConstrs}}: a list of general constraints.

\sphinxAtStartPar
\sphinxcode{\sphinxupquote{names}}: string array of names for general constraints.
\end{quote}
\end{quote}

\subsubsection{Model.SetNames()}
\label{\detokenize{csapi/model:id116}}\begin{quote}

\sphinxAtStartPar
Set names for given nonlinear constraints in model.

\sphinxAtStartPar
\sphinxstylestrong{Synopsis}
\begin{quote}

\sphinxAtStartPar
\sphinxcode{\sphinxupquote{void SetNames(NlConstraint{[}{]} cons, string{[}{]} names)}}
\end{quote}

\sphinxAtStartPar
\sphinxstylestrong{Arguments}
\begin{quote}

\sphinxAtStartPar
\sphinxcode{\sphinxupquote{cons}}: array of nonlinear constraints.

\sphinxAtStartPar
\sphinxcode{\sphinxupquote{names}}: string array of names for nonlinear constraints.
\end{quote}
\end{quote}

\subsubsection{Model.SetNames()}
\label{\detokenize{csapi/model:id117}}\begin{quote}

\sphinxAtStartPar
Set names for given nonlinear constraints in model.

\sphinxAtStartPar
\sphinxstylestrong{Synopsis}
\begin{quote}

\sphinxAtStartPar
\sphinxcode{\sphinxupquote{void SetNames(NlConstrArray cons, string{[}{]} names)}}
\end{quote}

\sphinxAtStartPar
\sphinxstylestrong{Arguments}
\begin{quote}

\sphinxAtStartPar
\sphinxcode{\sphinxupquote{cons}}: array object of nonlinear constraints.

\sphinxAtStartPar
\sphinxcode{\sphinxupquote{names}}: string array of names for nonlinear constraints.
\end{quote}
\end{quote}

\subsubsection{Model.SetNames()}
\label{\detokenize{csapi/model:id118}}\begin{quote}

\sphinxAtStartPar
Set names for given quadratic constraints in model.

\sphinxAtStartPar
\sphinxstylestrong{Synopsis}
\begin{quote}

\sphinxAtStartPar
\sphinxcode{\sphinxupquote{void SetNames(QConstraint{[}{]} cons, string{[}{]} names)}}
\end{quote}

\sphinxAtStartPar
\sphinxstylestrong{Arguments}
\begin{quote}

\sphinxAtStartPar
\sphinxcode{\sphinxupquote{cons}}: array of quadratic constraints.

\sphinxAtStartPar
\sphinxcode{\sphinxupquote{names}}: string array of names for quadratic constraints.
\end{quote}
\end{quote}

\subsubsection{Model.SetNames()}
\label{\detokenize{csapi/model:id119}}\begin{quote}

\sphinxAtStartPar
Set names for given quadratic constraints in model.

\sphinxAtStartPar
\sphinxstylestrong{Synopsis}
\begin{quote}

\sphinxAtStartPar
\sphinxcode{\sphinxupquote{void SetNames(QConstrArray cons, string{[}{]} names)}}
\end{quote}

\sphinxAtStartPar
\sphinxstylestrong{Arguments}
\begin{quote}

\sphinxAtStartPar
\sphinxcode{\sphinxupquote{cons}}: a list of quadratic constraints.

\sphinxAtStartPar
\sphinxcode{\sphinxupquote{names}}: string array of names for quadratic constraints.
\end{quote}
\end{quote}

\subsubsection{Model.SetNames()}
\label{\detokenize{csapi/model:id120}}\begin{quote}

\sphinxAtStartPar
Set names for given PSD variables in model.

\sphinxAtStartPar
\sphinxstylestrong{Synopsis}
\begin{quote}

\sphinxAtStartPar
\sphinxcode{\sphinxupquote{void SetNames(PsdVar{[}{]} vars, string{[}{]} names)}}
\end{quote}

\sphinxAtStartPar
\sphinxstylestrong{Arguments}
\begin{quote}

\sphinxAtStartPar
\sphinxcode{\sphinxupquote{vars}}: array of PSD variables.

\sphinxAtStartPar
\sphinxcode{\sphinxupquote{names}}: string array of names for PSD variables.
\end{quote}
\end{quote}

\subsubsection{Model.SetNames()}
\label{\detokenize{csapi/model:id121}}\begin{quote}

\sphinxAtStartPar
Set names for given PSD variables in model.

\sphinxAtStartPar
\sphinxstylestrong{Synopsis}
\begin{quote}

\sphinxAtStartPar
\sphinxcode{\sphinxupquote{void SetNames(PsdVarArray vars, string{[}{]} names)}}
\end{quote}

\sphinxAtStartPar
\sphinxstylestrong{Arguments}
\begin{quote}

\sphinxAtStartPar
\sphinxcode{\sphinxupquote{vars}}: a list of PSD variables.

\sphinxAtStartPar
\sphinxcode{\sphinxupquote{names}}: string array of names for PSD variables.
\end{quote}
\end{quote}

\subsubsection{Model.SetNames()}
\label{\detokenize{csapi/model:id122}}\begin{quote}

\sphinxAtStartPar
Set names for given PSD constraints in model.

\sphinxAtStartPar
\sphinxstylestrong{Synopsis}
\begin{quote}

\sphinxAtStartPar
\sphinxcode{\sphinxupquote{void SetNames(PsdConstraint{[}{]} cons, string{[}{]} names)}}
\end{quote}

\sphinxAtStartPar
\sphinxstylestrong{Arguments}
\begin{quote}

\sphinxAtStartPar
\sphinxcode{\sphinxupquote{cons}}: array of PSD constraints.

\sphinxAtStartPar
\sphinxcode{\sphinxupquote{names}}: string array of names for PSD constraints.
\end{quote}
\end{quote}

\subsubsection{Model.SetNames()}
\label{\detokenize{csapi/model:id123}}\begin{quote}

\sphinxAtStartPar
Set names for given PSD constraints in model.

\sphinxAtStartPar
\sphinxstylestrong{Synopsis}
\begin{quote}

\sphinxAtStartPar
\sphinxcode{\sphinxupquote{void SetNames(PsdConstrArray cons, string{[}{]} names)}}
\end{quote}

\sphinxAtStartPar
\sphinxstylestrong{Arguments}
\begin{quote}

\sphinxAtStartPar
\sphinxcode{\sphinxupquote{cons}}: a list of PSD constraints.

\sphinxAtStartPar
\sphinxcode{\sphinxupquote{names}}: string array of names for PSD constraints.
\end{quote}
\end{quote}

\subsubsection{Model.SetNames()}
\label{\detokenize{csapi/model:id124}}\begin{quote}

\sphinxAtStartPar
Set names for given LMI constraints in model.

\sphinxAtStartPar
\sphinxstylestrong{Synopsis}
\begin{quote}

\sphinxAtStartPar
\sphinxcode{\sphinxupquote{void SetNames(LmiConstraint{[}{]} cons, string{[}{]} names)}}
\end{quote}

\sphinxAtStartPar
\sphinxstylestrong{Arguments}
\begin{quote}

\sphinxAtStartPar
\sphinxcode{\sphinxupquote{cons}}: array of LMI constraints.

\sphinxAtStartPar
\sphinxcode{\sphinxupquote{names}}: string array of names for LMI constraints.
\end{quote}
\end{quote}

\subsubsection{Model.SetNames()}
\label{\detokenize{csapi/model:id125}}\begin{quote}

\sphinxAtStartPar
Set names for given LMI constraints in model.

\sphinxAtStartPar
\sphinxstylestrong{Synopsis}
\begin{quote}

\sphinxAtStartPar
\sphinxcode{\sphinxupquote{void SetNames(LmiConstrArray cons, string{[}{]} names)}}
\end{quote}

\sphinxAtStartPar
\sphinxstylestrong{Arguments}
\begin{quote}

\sphinxAtStartPar
\sphinxcode{\sphinxupquote{cons}}: a list of LMI constraints.

\sphinxAtStartPar
\sphinxcode{\sphinxupquote{names}}: string array of names for LMI constraints.
\end{quote}
\end{quote}

\subsubsection{Model.SetNames()}
\label{\detokenize{csapi/model:id126}}\begin{quote}

\sphinxAtStartPar
Set names for given affine cone constraints in model.

\sphinxAtStartPar
\sphinxstylestrong{Synopsis}
\begin{quote}

\sphinxAtStartPar
\sphinxcode{\sphinxupquote{void SetNames(AffineConeArray cones, string{[}{]} names)}}
\end{quote}

\sphinxAtStartPar
\sphinxstylestrong{Arguments}
\begin{quote}

\sphinxAtStartPar
\sphinxcode{\sphinxupquote{cones}}: an array of affine cone constraints.

\sphinxAtStartPar
\sphinxcode{\sphinxupquote{names}}: string array of names for affine cone constraints.
\end{quote}
\end{quote}

\subsubsection{Model.SetNlObjective()}
\label{\detokenize{csapi/model:model-setnlobjective}}\begin{quote}

\sphinxAtStartPar
Set nonlinear objective for model.

\sphinxAtStartPar
\sphinxstylestrong{Synopsis}
\begin{quote}

\sphinxAtStartPar
\sphinxcode{\sphinxupquote{void SetNlObjective(NlExpr expr, int sense)}}
\end{quote}

\sphinxAtStartPar
\sphinxstylestrong{Arguments}
\begin{quote}

\sphinxAtStartPar
\sphinxcode{\sphinxupquote{expr}}: nonlinear expression of the objective.

\sphinxAtStartPar
\sphinxcode{\sphinxupquote{sense}}: optimization sense. optional, default value 0 does not change COPT sense.
\end{quote}
\end{quote}

\subsubsection{Model.SetNlPrimalStart()}
\label{\detokenize{csapi/model:model-setnlprimalstart}}\begin{quote}

\sphinxAtStartPar
Given count, set initial values for variables of NLP from beginning.

\sphinxAtStartPar
\sphinxstylestrong{Synopsis}
\begin{quote}

\sphinxAtStartPar
\sphinxcode{\sphinxupquote{void SetNlPrimalStart(int count, double{[}{]} vals)}}
\end{quote}

\sphinxAtStartPar
\sphinxstylestrong{Arguments}
\begin{quote}

\sphinxAtStartPar
\sphinxcode{\sphinxupquote{count}}: the number of variables to set.

\sphinxAtStartPar
\sphinxcode{\sphinxupquote{vals}}: initial values of variables.
\end{quote}
\end{quote}

\subsubsection{Model.SetNlPrimalStart()}
\label{\detokenize{csapi/model:id127}}\begin{quote}

\sphinxAtStartPar
Set initial value for the specified variable of NLP.

\sphinxAtStartPar
\sphinxstylestrong{Synopsis}
\begin{quote}

\sphinxAtStartPar
\sphinxcode{\sphinxupquote{void SetNlPrimalStart(Var var, double val)}}
\end{quote}

\sphinxAtStartPar
\sphinxstylestrong{Arguments}
\begin{quote}

\sphinxAtStartPar
\sphinxcode{\sphinxupquote{var}}: an interested variable.

\sphinxAtStartPar
\sphinxcode{\sphinxupquote{val}}: initial value of the variable.
\end{quote}
\end{quote}

\subsubsection{Model.SetNlPrimalStart()}
\label{\detokenize{csapi/model:id128}}\begin{quote}

\sphinxAtStartPar
Set initial values for an array of variables of NLP.

\sphinxAtStartPar
\sphinxstylestrong{Synopsis}
\begin{quote}

\sphinxAtStartPar
\sphinxcode{\sphinxupquote{void SetNlPrimalStart(Var{[}{]} vars, double{[}{]} vals)}}
\end{quote}

\sphinxAtStartPar
\sphinxstylestrong{Arguments}
\begin{quote}

\sphinxAtStartPar
\sphinxcode{\sphinxupquote{vars}}: array of interested variables.

\sphinxAtStartPar
\sphinxcode{\sphinxupquote{vals}}: initial values of variables.
\end{quote}
\end{quote}

\subsubsection{Model.SetNlPrimalStart()}
\label{\detokenize{csapi/model:id129}}\begin{quote}

\sphinxAtStartPar
Set initial values for variable array of NLP.

\sphinxAtStartPar
\sphinxstylestrong{Synopsis}
\begin{quote}

\sphinxAtStartPar
\sphinxcode{\sphinxupquote{void SetNlPrimalStart(VarArray vars, double{[}{]} vals)}}
\end{quote}

\sphinxAtStartPar
\sphinxstylestrong{Arguments}
\begin{quote}

\sphinxAtStartPar
\sphinxcode{\sphinxupquote{vars}}: a list of interested variables.

\sphinxAtStartPar
\sphinxcode{\sphinxupquote{vals}}: initial values of variables.
\end{quote}
\end{quote}

\subsubsection{Model.SetObjConst()}
\label{\detokenize{csapi/model:model-setobjconst}}\begin{quote}

\sphinxAtStartPar
Set objective constant.

\sphinxAtStartPar
\sphinxstylestrong{Synopsis}
\begin{quote}

\sphinxAtStartPar
\sphinxcode{\sphinxupquote{void SetObjConst(double constant)}}
\end{quote}

\sphinxAtStartPar
\sphinxstylestrong{Arguments}
\begin{quote}

\sphinxAtStartPar
\sphinxcode{\sphinxupquote{constant}}: constant value to set.
\end{quote}
\end{quote}

\subsubsection{Model.SetObjective()}
\label{\detokenize{csapi/model:model-setobjective}}\begin{quote}

\sphinxAtStartPar
Set objective for model.

\sphinxAtStartPar
\sphinxstylestrong{Synopsis}
\begin{quote}

\sphinxAtStartPar
\sphinxcode{\sphinxupquote{void SetObjective(MExpression expr, int sense)}}
\end{quote}

\sphinxAtStartPar
\sphinxstylestrong{Arguments}
\begin{quote}

\sphinxAtStartPar
\sphinxcode{\sphinxupquote{expr}}: expression of the objective.

\sphinxAtStartPar
\sphinxcode{\sphinxupquote{sense}}: optimization sense. optional, default value 0 does not change COPT sense
\end{quote}
\end{quote}

\subsubsection{Model.SetObjectiveN()}
\label{\detokenize{csapi/model:model-setobjectiven}}\begin{quote}

\sphinxAtStartPar
Set a multi\sphinxhyphen{}objective function in model.

\sphinxAtStartPar
\sphinxstylestrong{Synopsis}
\begin{quote}

\sphinxAtStartPar
\sphinxcode{\sphinxupquote{void SetObjectiveN(}}
\begin{quote}

\sphinxAtStartPar
\sphinxcode{\sphinxupquote{int idx,}}

\sphinxAtStartPar
\sphinxcode{\sphinxupquote{MExpression expr,}}

\sphinxAtStartPar
\sphinxcode{\sphinxupquote{int sense,}}

\sphinxAtStartPar
\sphinxcode{\sphinxupquote{double priority,}}

\sphinxAtStartPar
\sphinxcode{\sphinxupquote{double weight,}}

\sphinxAtStartPar
\sphinxcode{\sphinxupquote{double abstol,}}

\sphinxAtStartPar
\sphinxcode{\sphinxupquote{double reltol)}}
\end{quote}
\end{quote}

\sphinxAtStartPar
\sphinxstylestrong{Arguments}
\begin{quote}

\sphinxAtStartPar
\sphinxcode{\sphinxupquote{idx}}: index of a multi\sphinxhyphen{}objective function.

\sphinxAtStartPar
\sphinxcode{\sphinxupquote{expr}}: linear expression of the multi\sphinxhyphen{}objective function.

\sphinxAtStartPar
\sphinxcode{\sphinxupquote{sense}}: optimization sense. optional, default value 0 does not change COPT sense.

\sphinxAtStartPar
\sphinxcode{\sphinxupquote{priority}}: an optional objective parameter. Default value is 0.0.

\sphinxAtStartPar
\sphinxcode{\sphinxupquote{weight}}: an optional objective parameter. Default value is 1.0.

\sphinxAtStartPar
\sphinxcode{\sphinxupquote{abstol}}: absolute tolerance is an optional objective parameter. Default value is 1e\sphinxhyphen{}6.

\sphinxAtStartPar
\sphinxcode{\sphinxupquote{reltol}}: relative tolerance is an optional objective parameter. Default value is 0.0.
\end{quote}
\end{quote}

\subsubsection{Model.SetObjParamN()}
\label{\detokenize{csapi/model:model-setobjparamn}}\begin{quote}

\sphinxAtStartPar
Set value of objective parameter of a multi\sphinxhyphen{}objective function.

\sphinxAtStartPar
\sphinxstylestrong{Synopsis}
\begin{quote}

\sphinxAtStartPar
\sphinxcode{\sphinxupquote{void SetObjParamN(}}
\begin{quote}

\sphinxAtStartPar
\sphinxcode{\sphinxupquote{int idx,}}

\sphinxAtStartPar
\sphinxcode{\sphinxupquote{string param,}}

\sphinxAtStartPar
\sphinxcode{\sphinxupquote{double val)}}
\end{quote}
\end{quote}

\sphinxAtStartPar
\sphinxstylestrong{Arguments}
\begin{quote}

\sphinxAtStartPar
\sphinxcode{\sphinxupquote{idx}}: index of a multi\sphinxhyphen{}objective function.

\sphinxAtStartPar
\sphinxcode{\sphinxupquote{param}}: name of objective parameter, including priority, weight, abstol and reltol.

\sphinxAtStartPar
\sphinxcode{\sphinxupquote{val}}: new value of objective parameter.
\end{quote}
\end{quote}

\subsubsection{Model.SetObjSense()}
\label{\detokenize{csapi/model:model-setobjsense}}\begin{quote}

\sphinxAtStartPar
Set objective sense for model.

\sphinxAtStartPar
\sphinxstylestrong{Synopsis}
\begin{quote}

\sphinxAtStartPar
\sphinxcode{\sphinxupquote{void SetObjSense(int sense)}}
\end{quote}

\sphinxAtStartPar
\sphinxstylestrong{Arguments}
\begin{quote}

\sphinxAtStartPar
\sphinxcode{\sphinxupquote{sense}}: the objective sense.
\end{quote}
\end{quote}

\subsubsection{Model.SetPsdCoeff()}
\label{\detokenize{csapi/model:model-setpsdcoeff}}\begin{quote}

\sphinxAtStartPar
Set the coefficient matrix of a PSD variable in a PSD constraint.

\sphinxAtStartPar
\sphinxstylestrong{Synopsis}
\begin{quote}

\sphinxAtStartPar
\sphinxcode{\sphinxupquote{void SetPsdCoeff(}}
\begin{quote}

\sphinxAtStartPar
\sphinxcode{\sphinxupquote{PsdConstraint constr,}}

\sphinxAtStartPar
\sphinxcode{\sphinxupquote{PsdVar var,}}

\sphinxAtStartPar
\sphinxcode{\sphinxupquote{SymMatrix mat)}}
\end{quote}
\end{quote}

\sphinxAtStartPar
\sphinxstylestrong{Arguments}
\begin{quote}

\sphinxAtStartPar
\sphinxcode{\sphinxupquote{constr}}: The desired PSD constraint.

\sphinxAtStartPar
\sphinxcode{\sphinxupquote{var}}: The desired PSD variable.

\sphinxAtStartPar
\sphinxcode{\sphinxupquote{mat}}: new coefficient matrix.
\end{quote}
\end{quote}

\subsubsection{Model.SetPsdObjective()}
\label{\detokenize{csapi/model:model-setpsdobjective}}\begin{quote}

\sphinxAtStartPar
Set PSD objective for model.

\sphinxAtStartPar
\sphinxstylestrong{Synopsis}
\begin{quote}

\sphinxAtStartPar
\sphinxcode{\sphinxupquote{void SetPsdObjective(PsdExpr expr, int sense)}}
\end{quote}

\sphinxAtStartPar
\sphinxstylestrong{Arguments}
\begin{quote}

\sphinxAtStartPar
\sphinxcode{\sphinxupquote{expr}}: PSD expression of the objective.

\sphinxAtStartPar
\sphinxcode{\sphinxupquote{sense}}: optimization sense. optional, default value 0 does not change COPT sense.
\end{quote}
\end{quote}

\subsubsection{Model.SetQuadObjective()}
\label{\detokenize{csapi/model:model-setquadobjective}}\begin{quote}

\sphinxAtStartPar
Set quadratic objective for model.

\sphinxAtStartPar
\sphinxstylestrong{Synopsis}
\begin{quote}

\sphinxAtStartPar
\sphinxcode{\sphinxupquote{void SetQuadObjective(QuadExpr expr, int sense)}}
\end{quote}

\sphinxAtStartPar
\sphinxstylestrong{Arguments}
\begin{quote}

\sphinxAtStartPar
\sphinxcode{\sphinxupquote{expr}}: quadratic expression of the objective.

\sphinxAtStartPar
\sphinxcode{\sphinxupquote{sense}}: optimization sense. optional, default value 0 does not change COPT sense.
\end{quote}
\end{quote}

\subsubsection{Model.SetSlackBasis()}
\label{\detokenize{csapi/model:model-setslackbasis}}\begin{quote}

\sphinxAtStartPar
Set slack basis to model.

\sphinxAtStartPar
\sphinxstylestrong{Synopsis}
\begin{quote}

\sphinxAtStartPar
\sphinxcode{\sphinxupquote{void SetSlackBasis()}}
\end{quote}
\end{quote}

\subsubsection{Model.SetSolverLogFile()}
\label{\detokenize{csapi/model:model-setsolverlogfile}}\begin{quote}

\sphinxAtStartPar
Set log file for COPT.

\sphinxAtStartPar
\sphinxstylestrong{Synopsis}
\begin{quote}

\sphinxAtStartPar
\sphinxcode{\sphinxupquote{void SetSolverLogFile(string filename)}}
\end{quote}

\sphinxAtStartPar
\sphinxstylestrong{Arguments}
\begin{quote}

\sphinxAtStartPar
\sphinxcode{\sphinxupquote{filename}}: log file name.
\end{quote}
\end{quote}

\subsubsection{Model.Solve()}
\label{\detokenize{csapi/model:model-solve}}\begin{quote}

\sphinxAtStartPar
Solve the model as MIP.

\sphinxAtStartPar
\sphinxstylestrong{Synopsis}
\begin{quote}

\sphinxAtStartPar
\sphinxcode{\sphinxupquote{void Solve()}}
\end{quote}
\end{quote}

\subsubsection{Model.SolveLp()}
\label{\detokenize{csapi/model:model-solvelp}}\begin{quote}

\sphinxAtStartPar
Solve the model as LP.

\sphinxAtStartPar
\sphinxstylestrong{Synopsis}
\begin{quote}

\sphinxAtStartPar
\sphinxcode{\sphinxupquote{void SolveLp()}}
\end{quote}
\end{quote}

\subsubsection{Model.Tune()}
\label{\detokenize{csapi/model:model-tune}}\begin{quote}

\sphinxAtStartPar
Tune model.

\sphinxAtStartPar
\sphinxstylestrong{Synopsis}
\begin{quote}

\sphinxAtStartPar
\sphinxcode{\sphinxupquote{void Tune()}}
\end{quote}
\end{quote}

\subsubsection{Model.Write()}
\label{\detokenize{csapi/model:model-write}}\begin{quote}

\sphinxAtStartPar
Output problem, solution, basis, MIP start or modified COPT parameters to file.

\sphinxAtStartPar
\sphinxstylestrong{Synopsis}
\begin{quote}

\sphinxAtStartPar
\sphinxcode{\sphinxupquote{void Write(string filename)}}
\end{quote}

\sphinxAtStartPar
\sphinxstylestrong{Arguments}
\begin{quote}

\sphinxAtStartPar
\sphinxcode{\sphinxupquote{filename}}: an output file name.
\end{quote}
\end{quote}

\subsubsection{Model.WriteBasis()}
\label{\detokenize{csapi/model:model-writebasis}}\begin{quote}

\sphinxAtStartPar
Output optimal basis to a file of type ‘.bas’.

\sphinxAtStartPar
\sphinxstylestrong{Synopsis}
\begin{quote}

\sphinxAtStartPar
\sphinxcode{\sphinxupquote{void WriteBasis(string filename)}}
\end{quote}

\sphinxAtStartPar
\sphinxstylestrong{Arguments}
\begin{quote}

\sphinxAtStartPar
\sphinxcode{\sphinxupquote{filename}}: an output file name.
\end{quote}
\end{quote}

\subsubsection{Model.WriteBin()}
\label{\detokenize{csapi/model:model-writebin}}\begin{quote}

\sphinxAtStartPar
Output problem to a file as COPT binary format.

\sphinxAtStartPar
\sphinxstylestrong{Synopsis}
\begin{quote}

\sphinxAtStartPar
\sphinxcode{\sphinxupquote{void WriteBin(string filename)}}
\end{quote}

\sphinxAtStartPar
\sphinxstylestrong{Arguments}
\begin{quote}

\sphinxAtStartPar
\sphinxcode{\sphinxupquote{filename}}: an output file name.
\end{quote}
\end{quote}

\subsubsection{Model.WriteIIS()}
\label{\detokenize{csapi/model:model-writeiis}}\begin{quote}

\sphinxAtStartPar
Output IIS to file.

\sphinxAtStartPar
\sphinxstylestrong{Synopsis}
\begin{quote}

\sphinxAtStartPar
\sphinxcode{\sphinxupquote{void WriteIIS(string filename)}}
\end{quote}

\sphinxAtStartPar
\sphinxstylestrong{Arguments}
\begin{quote}

\sphinxAtStartPar
\sphinxcode{\sphinxupquote{filename}}: Output file name.
\end{quote}
\end{quote}

\subsubsection{Model.WriteJsonSol()}
\label{\detokenize{csapi/model:model-writejsonsol}}\begin{quote}

\sphinxAtStartPar
Output solution to a file of type ‘.json’.

\sphinxAtStartPar
\sphinxstylestrong{Synopsis}
\begin{quote}

\sphinxAtStartPar
\sphinxcode{\sphinxupquote{void WriteJsonSol(string filename)}}
\end{quote}

\sphinxAtStartPar
\sphinxstylestrong{Arguments}
\begin{quote}

\sphinxAtStartPar
\sphinxcode{\sphinxupquote{filename}}: an output file name.
\end{quote}
\end{quote}

\subsubsection{Model.WriteLp()}
\label{\detokenize{csapi/model:model-writelp}}\begin{quote}

\sphinxAtStartPar
Output problem to a file as LP format.

\sphinxAtStartPar
\sphinxstylestrong{Synopsis}
\begin{quote}

\sphinxAtStartPar
\sphinxcode{\sphinxupquote{void WriteLp(string filename)}}
\end{quote}

\sphinxAtStartPar
\sphinxstylestrong{Arguments}
\begin{quote}

\sphinxAtStartPar
\sphinxcode{\sphinxupquote{filename}}: an output file name.
\end{quote}
\end{quote}

\subsubsection{Model.WriteMps()}
\label{\detokenize{csapi/model:model-writemps}}\begin{quote}

\sphinxAtStartPar
Output problem to a file as MPS format.

\sphinxAtStartPar
\sphinxstylestrong{Synopsis}
\begin{quote}

\sphinxAtStartPar
\sphinxcode{\sphinxupquote{void WriteMps(string filename)}}
\end{quote}

\sphinxAtStartPar
\sphinxstylestrong{Arguments}
\begin{quote}

\sphinxAtStartPar
\sphinxcode{\sphinxupquote{filename}}: an output file name.
\end{quote}
\end{quote}

\subsubsection{Model.WriteMpsStr()}
\label{\detokenize{csapi/model:model-writempsstr}}\begin{quote}

\sphinxAtStartPar
Output MPS problem to problem buffer.

\sphinxAtStartPar
\sphinxstylestrong{Synopsis}
\begin{quote}

\sphinxAtStartPar
\sphinxcode{\sphinxupquote{ProbBuffer WriteMpsStr()}}
\end{quote}

\sphinxAtStartPar
\sphinxstylestrong{Return}
\begin{quote}

\sphinxAtStartPar
problem buffer for string of MPS problem.
\end{quote}
\end{quote}

\subsubsection{Model.WriteMst()}
\label{\detokenize{csapi/model:model-writemst}}\begin{quote}

\sphinxAtStartPar
Output MIP start information to a file of type ‘.mst’.

\sphinxAtStartPar
\sphinxstylestrong{Synopsis}
\begin{quote}

\sphinxAtStartPar
\sphinxcode{\sphinxupquote{void WriteMst(string filename)}}
\end{quote}

\sphinxAtStartPar
\sphinxstylestrong{Arguments}
\begin{quote}

\sphinxAtStartPar
\sphinxcode{\sphinxupquote{filename}}: an output file name.
\end{quote}
\end{quote}

\subsubsection{Model.WriteNL()}
\label{\detokenize{csapi/model:model-writenl}}\begin{quote}

\sphinxAtStartPar
Output problem to a file as NL format.

\sphinxAtStartPar
\sphinxstylestrong{Synopsis}
\begin{quote}

\sphinxAtStartPar
\sphinxcode{\sphinxupquote{void WriteNL(string filename)}}
\end{quote}

\sphinxAtStartPar
\sphinxstylestrong{Arguments}
\begin{quote}

\sphinxAtStartPar
\sphinxcode{\sphinxupquote{filename}}: an output file name.
\end{quote}
\end{quote}

\subsubsection{Model.WriteOrd()}
\label{\detokenize{csapi/model:model-writeord}}\begin{quote}

\sphinxAtStartPar
Output branching order to file.

\sphinxAtStartPar
\sphinxstylestrong{Synopsis}
\begin{quote}

\sphinxAtStartPar
\sphinxcode{\sphinxupquote{void WriteOrd(string filename)}}
\end{quote}

\sphinxAtStartPar
\sphinxstylestrong{Arguments}
\begin{quote}

\sphinxAtStartPar
\sphinxcode{\sphinxupquote{filename}}: Output file name.
\end{quote}
\end{quote}

\subsubsection{Model.WriteParam()}
\label{\detokenize{csapi/model:model-writeparam}}\begin{quote}

\sphinxAtStartPar
Output modified COPT parameters to a file of type ‘.par’.

\sphinxAtStartPar
\sphinxstylestrong{Synopsis}
\begin{quote}

\sphinxAtStartPar
\sphinxcode{\sphinxupquote{void WriteParam(string filename)}}
\end{quote}

\sphinxAtStartPar
\sphinxstylestrong{Arguments}
\begin{quote}

\sphinxAtStartPar
\sphinxcode{\sphinxupquote{filename}}: an output file name.
\end{quote}
\end{quote}

\subsubsection{Model.WritePoolSol()}
\label{\detokenize{csapi/model:model-writepoolsol}}\begin{quote}

\sphinxAtStartPar
Output selected pool solution to a file of type ‘.sol’.

\sphinxAtStartPar
\sphinxstylestrong{Synopsis}
\begin{quote}

\sphinxAtStartPar
\sphinxcode{\sphinxupquote{void WritePoolSol(int idx, string filename)}}
\end{quote}

\sphinxAtStartPar
\sphinxstylestrong{Arguments}
\begin{quote}

\sphinxAtStartPar
\sphinxcode{\sphinxupquote{idx}}: index of pool solution.

\sphinxAtStartPar
\sphinxcode{\sphinxupquote{filename}}: an output file name.
\end{quote}
\end{quote}

\subsubsection{Model.WriteRelax()}
\label{\detokenize{csapi/model:model-writerelax}}\begin{quote}

\sphinxAtStartPar
Output feasibility relaxation problem to file.

\sphinxAtStartPar
\sphinxstylestrong{Synopsis}
\begin{quote}

\sphinxAtStartPar
\sphinxcode{\sphinxupquote{void WriteRelax(string filename)}}
\end{quote}

\sphinxAtStartPar
\sphinxstylestrong{Arguments}
\begin{quote}

\sphinxAtStartPar
\sphinxcode{\sphinxupquote{filename}}: Output file name.
\end{quote}
\end{quote}

\subsubsection{Model.WriteSol()}
\label{\detokenize{csapi/model:model-writesol}}\begin{quote}

\sphinxAtStartPar
Output solution to a file of type ‘.sol’.

\sphinxAtStartPar
\sphinxstylestrong{Synopsis}
\begin{quote}

\sphinxAtStartPar
\sphinxcode{\sphinxupquote{void WriteSol(string filename)}}
\end{quote}

\sphinxAtStartPar
\sphinxstylestrong{Arguments}
\begin{quote}

\sphinxAtStartPar
\sphinxcode{\sphinxupquote{filename}}: an output file name.
\end{quote}
\end{quote}

\subsubsection{Model.WriteTuneParam()}
\label{\detokenize{csapi/model:model-writetuneparam}}\begin{quote}

\sphinxAtStartPar
Output specified tuned parameters to a file of type ‘.par’.

\sphinxAtStartPar
\sphinxstylestrong{Synopsis}
\begin{quote}

\sphinxAtStartPar
\sphinxcode{\sphinxupquote{void WriteTuneParam(int idx, string filename)}}
\end{quote}

\sphinxAtStartPar
\sphinxstylestrong{Arguments}
\begin{quote}

\sphinxAtStartPar
\sphinxcode{\sphinxupquote{idx}}: Index of tuned parameters.

\sphinxAtStartPar
\sphinxcode{\sphinxupquote{filename}}: Output file name.
\end{quote}
\end{quote}

\subsection{Var}
\label{\detokenize{csharpapiref:var}}\label{\detokenize{csharpapiref:chapcsharpapiref-var}}
\sphinxAtStartPar
COPT variable object. Variables are always associated with a particular model.
User creates a variable object by adding a variable to a model, rather than
by using constructor of Var class.

\sphinxstepscope

\subsubsection{Var.Get()}
\label{\detokenize{csapi/var:var-get}}\label{\detokenize{csapi/var::doc}}\begin{quote}

\sphinxAtStartPar
Get information value of the variable. Support informations of “Value”, “RedCost”, “PrimalRay”, “LB”, “UB”, “Obj” and “BranchFactor”.

\sphinxAtStartPar
\sphinxstylestrong{Synopsis}
\begin{quote}

\sphinxAtStartPar
\sphinxcode{\sphinxupquote{double Get(string info)}}
\end{quote}

\sphinxAtStartPar
\sphinxstylestrong{Arguments}
\begin{quote}

\sphinxAtStartPar
\sphinxcode{\sphinxupquote{info}}: information name.
\end{quote}

\sphinxAtStartPar
\sphinxstylestrong{Return}
\begin{quote}

\sphinxAtStartPar
information value.
\end{quote}
\end{quote}

\subsubsection{Var.GetBasis()}
\label{\detokenize{csapi/var:var-getbasis}}\begin{quote}

\sphinxAtStartPar
Get basis status of the variable.

\sphinxAtStartPar
\sphinxstylestrong{Synopsis}
\begin{quote}

\sphinxAtStartPar
\sphinxcode{\sphinxupquote{int GetBasis()}}
\end{quote}

\sphinxAtStartPar
\sphinxstylestrong{Return}
\begin{quote}

\sphinxAtStartPar
Basis status.
\end{quote}
\end{quote}

\subsubsection{Var.GetIdx()}
\label{\detokenize{csapi/var:var-getidx}}\begin{quote}

\sphinxAtStartPar
Get index of the variable.

\sphinxAtStartPar
\sphinxstylestrong{Synopsis}
\begin{quote}

\sphinxAtStartPar
\sphinxcode{\sphinxupquote{int GetIdx()}}
\end{quote}

\sphinxAtStartPar
\sphinxstylestrong{Return}
\begin{quote}

\sphinxAtStartPar
variable index.
\end{quote}
\end{quote}

\subsubsection{Var.GetLowerIIS()}
\label{\detokenize{csapi/var:var-getloweriis}}\begin{quote}

\sphinxAtStartPar
Get IIS status for lower bound of the variable.

\sphinxAtStartPar
\sphinxstylestrong{Synopsis}
\begin{quote}

\sphinxAtStartPar
\sphinxcode{\sphinxupquote{int GetLowerIIS()}}
\end{quote}

\sphinxAtStartPar
\sphinxstylestrong{Return}
\begin{quote}

\sphinxAtStartPar
IIS status.
\end{quote}
\end{quote}

\subsubsection{Var.GetName()}
\label{\detokenize{csapi/var:var-getname}}\begin{quote}

\sphinxAtStartPar
Get name of the variable.

\sphinxAtStartPar
\sphinxstylestrong{Synopsis}
\begin{quote}

\sphinxAtStartPar
\sphinxcode{\sphinxupquote{string GetName()}}
\end{quote}

\sphinxAtStartPar
\sphinxstylestrong{Return}
\begin{quote}

\sphinxAtStartPar
variable name.
\end{quote}
\end{quote}

\subsubsection{Var.GetType()}
\label{\detokenize{csapi/var:var-gettype}}\begin{quote}

\sphinxAtStartPar
Get type of the variable.

\sphinxAtStartPar
\sphinxstylestrong{Synopsis}
\begin{quote}

\sphinxAtStartPar
\sphinxcode{\sphinxupquote{char GetType()}}
\end{quote}

\sphinxAtStartPar
\sphinxstylestrong{Return}
\begin{quote}

\sphinxAtStartPar
variable type.
\end{quote}
\end{quote}

\subsubsection{Var.GetUpperIIS()}
\label{\detokenize{csapi/var:var-getupperiis}}\begin{quote}

\sphinxAtStartPar
Get IIS status for upper bound of the variable.

\sphinxAtStartPar
\sphinxstylestrong{Synopsis}
\begin{quote}

\sphinxAtStartPar
\sphinxcode{\sphinxupquote{int GetUpperIIS()}}
\end{quote}

\sphinxAtStartPar
\sphinxstylestrong{Return}
\begin{quote}

\sphinxAtStartPar
IIS status.
\end{quote}
\end{quote}

\subsubsection{Var.Remove()}
\label{\detokenize{csapi/var:var-remove}}\begin{quote}

\sphinxAtStartPar
Remove variable from model.

\sphinxAtStartPar
\sphinxstylestrong{Synopsis}
\begin{quote}

\sphinxAtStartPar
\sphinxcode{\sphinxupquote{void Remove()}}
\end{quote}
\end{quote}

\subsubsection{Var.Set()}
\label{\detokenize{csapi/var:var-set}}\begin{quote}

\sphinxAtStartPar
Set information value of the variable. Support informations of “LB”, “UB”, “Obj” and “BranchFactor”.

\sphinxAtStartPar
\sphinxstylestrong{Synopsis}
\begin{quote}

\sphinxAtStartPar
\sphinxcode{\sphinxupquote{void Set(string info, double val)}}
\end{quote}

\sphinxAtStartPar
\sphinxstylestrong{Arguments}
\begin{quote}

\sphinxAtStartPar
\sphinxcode{\sphinxupquote{info}}: information name.

\sphinxAtStartPar
\sphinxcode{\sphinxupquote{val}}: new information value.
\end{quote}
\end{quote}

\subsubsection{Var.SetName()}
\label{\detokenize{csapi/var:var-setname}}\begin{quote}

\sphinxAtStartPar
Set name of the variable.

\sphinxAtStartPar
\sphinxstylestrong{Synopsis}
\begin{quote}

\sphinxAtStartPar
\sphinxcode{\sphinxupquote{void SetName(string name)}}
\end{quote}

\sphinxAtStartPar
\sphinxstylestrong{Arguments}
\begin{quote}

\sphinxAtStartPar
\sphinxcode{\sphinxupquote{name}}: variable name.
\end{quote}
\end{quote}

\subsubsection{Var.SetType()}
\label{\detokenize{csapi/var:var-settype}}\begin{quote}

\sphinxAtStartPar
Set type of the variable.

\sphinxAtStartPar
\sphinxstylestrong{Synopsis}
\begin{quote}

\sphinxAtStartPar
\sphinxcode{\sphinxupquote{void SetType(char vtype)}}
\end{quote}

\sphinxAtStartPar
\sphinxstylestrong{Arguments}
\begin{quote}

\sphinxAtStartPar
\sphinxcode{\sphinxupquote{vtype}}: variable type.
\end{quote}
\end{quote}

\subsection{VarArray}
\label{\detokenize{csharpapiref:vararray}}\label{\detokenize{csharpapiref:chapcsharpapiref-vararray}}
\sphinxAtStartPar
COPT variable array object. To store and access a set of C\# {\hyperref[\detokenize{csharpapiref:chapcsharpapiref-var}]{\sphinxcrossref{\DUrole{std,std-ref}{Var}}}} objects,
Cardinal Optimizer provides C\# VarArray class, which defines the following methods.

\sphinxstepscope

\subsubsection{VarArray.VarArray()}
\label{\detokenize{csapi/vararray:vararray-vararray}}\label{\detokenize{csapi/vararray::doc}}\begin{quote}

\sphinxAtStartPar
Constructor of vararray.

\sphinxAtStartPar
\sphinxstylestrong{Synopsis}
\begin{quote}

\sphinxAtStartPar
\sphinxcode{\sphinxupquote{VarArray()}}
\end{quote}
\end{quote}

\subsubsection{VarArray.GetVar()}
\label{\detokenize{csapi/vararray:vararray-getvar}}\begin{quote}

\sphinxAtStartPar
Get idx\sphinxhyphen{}th variable object.

\sphinxAtStartPar
\sphinxstylestrong{Synopsis}
\begin{quote}

\sphinxAtStartPar
\sphinxcode{\sphinxupquote{Var GetVar(int idx)}}
\end{quote}

\sphinxAtStartPar
\sphinxstylestrong{Arguments}
\begin{quote}

\sphinxAtStartPar
\sphinxcode{\sphinxupquote{idx}}: index of the variable.
\end{quote}

\sphinxAtStartPar
\sphinxstylestrong{Return}
\begin{quote}

\sphinxAtStartPar
variable object with index idx.
\end{quote}
\end{quote}

\subsubsection{VarArray.PushBack()}
\label{\detokenize{csapi/vararray:vararray-pushback}}\begin{quote}

\sphinxAtStartPar
Add a variable object to variable array.

\sphinxAtStartPar
\sphinxstylestrong{Synopsis}
\begin{quote}

\sphinxAtStartPar
\sphinxcode{\sphinxupquote{void PushBack(Var var)}}
\end{quote}

\sphinxAtStartPar
\sphinxstylestrong{Arguments}
\begin{quote}

\sphinxAtStartPar
\sphinxcode{\sphinxupquote{var}}: a variable object.
\end{quote}
\end{quote}

\subsubsection{VarArray.Size()}
\label{\detokenize{csapi/vararray:vararray-size}}\begin{quote}

\sphinxAtStartPar
Get the number of variable objects.

\sphinxAtStartPar
\sphinxstylestrong{Synopsis}
\begin{quote}

\sphinxAtStartPar
\sphinxcode{\sphinxupquote{int Size()}}
\end{quote}

\sphinxAtStartPar
\sphinxstylestrong{Return}
\begin{quote}

\sphinxAtStartPar
number of variable objects.
\end{quote}
\end{quote}

\subsection{Expr}
\label{\detokenize{csharpapiref:expr}}\label{\detokenize{csharpapiref:chapcsharpapiref-expr}}
\sphinxAtStartPar
COPT linear expression object. A linear expression consists of a constant term, a list of
terms of variables and associated coefficients. Linear expressions are used to build
constraints.

\sphinxstepscope

\subsubsection{Expr.Expr()}
\label{\detokenize{csapi/expr:expr-expr}}\label{\detokenize{csapi/expr::doc}}\begin{quote}

\sphinxAtStartPar
Constructor of a constant linear expression with default constant value 0.

\sphinxAtStartPar
\sphinxstylestrong{Synopsis}
\begin{quote}

\sphinxAtStartPar
\sphinxcode{\sphinxupquote{Expr(double constant)}}
\end{quote}

\sphinxAtStartPar
\sphinxstylestrong{Arguments}
\begin{quote}

\sphinxAtStartPar
\sphinxcode{\sphinxupquote{constant}}: optional, constant value in expression object.
\end{quote}
\end{quote}

\subsubsection{Expr.Expr()}
\label{\detokenize{csapi/expr:id1}}\begin{quote}

\sphinxAtStartPar
Constructor of a linear expression with one term.

\sphinxAtStartPar
\sphinxstylestrong{Synopsis}
\begin{quote}

\sphinxAtStartPar
\sphinxcode{\sphinxupquote{Expr(Var var, double coeff)}}
\end{quote}

\sphinxAtStartPar
\sphinxstylestrong{Arguments}
\begin{quote}

\sphinxAtStartPar
\sphinxcode{\sphinxupquote{var}}: variable for the added term.

\sphinxAtStartPar
\sphinxcode{\sphinxupquote{coeff}}: coefficent for the added term with default value 1.0.
\end{quote}
\end{quote}

\subsubsection{Expr.AddConstant()}
\label{\detokenize{csapi/expr:expr-addconstant}}\begin{quote}

\sphinxAtStartPar
Add extra constant to the expression.

\sphinxAtStartPar
\sphinxstylestrong{Synopsis}
\begin{quote}

\sphinxAtStartPar
\sphinxcode{\sphinxupquote{void AddConstant(double constant)}}
\end{quote}

\sphinxAtStartPar
\sphinxstylestrong{Arguments}
\begin{quote}

\sphinxAtStartPar
\sphinxcode{\sphinxupquote{constant}}: delta value to be added to expression constant.
\end{quote}
\end{quote}

\subsubsection{Expr.AddExpr()}
\label{\detokenize{csapi/expr:expr-addexpr}}\begin{quote}

\sphinxAtStartPar
Add an expression to self.

\sphinxAtStartPar
\sphinxstylestrong{Synopsis}
\begin{quote}

\sphinxAtStartPar
\sphinxcode{\sphinxupquote{void AddExpr(Expr expr, double mult)}}
\end{quote}

\sphinxAtStartPar
\sphinxstylestrong{Arguments}
\begin{quote}

\sphinxAtStartPar
\sphinxcode{\sphinxupquote{expr}}: expression to be added.

\sphinxAtStartPar
\sphinxcode{\sphinxupquote{mult}}: multiply constant.
\end{quote}
\end{quote}

\subsubsection{Expr.AddTerm()}
\label{\detokenize{csapi/expr:expr-addterm}}\begin{quote}

\sphinxAtStartPar
Add a term to expression object.

\sphinxAtStartPar
\sphinxstylestrong{Synopsis}
\begin{quote}

\sphinxAtStartPar
\sphinxcode{\sphinxupquote{void AddTerm(Var var, double coeff)}}
\end{quote}

\sphinxAtStartPar
\sphinxstylestrong{Arguments}
\begin{quote}

\sphinxAtStartPar
\sphinxcode{\sphinxupquote{var}}: a variable for new term.

\sphinxAtStartPar
\sphinxcode{\sphinxupquote{coeff}}: coefficient for new term.
\end{quote}
\end{quote}

\subsubsection{Expr.AddTerms()}
\label{\detokenize{csapi/expr:expr-addterms}}\begin{quote}

\sphinxAtStartPar
Add terms to expression object.

\sphinxAtStartPar
\sphinxstylestrong{Synopsis}
\begin{quote}

\sphinxAtStartPar
\sphinxcode{\sphinxupquote{void AddTerms(Var{[}{]} vars, double coeff)}}
\end{quote}

\sphinxAtStartPar
\sphinxstylestrong{Arguments}
\begin{quote}

\sphinxAtStartPar
\sphinxcode{\sphinxupquote{vars}}: variables for added terms.

\sphinxAtStartPar
\sphinxcode{\sphinxupquote{coeff}}: coefficient array for added terms with default value 1.0.
\end{quote}
\end{quote}

\subsubsection{Expr.AddTerms()}
\label{\detokenize{csapi/expr:id2}}\begin{quote}

\sphinxAtStartPar
Add terms to expression object.

\sphinxAtStartPar
\sphinxstylestrong{Synopsis}
\begin{quote}

\sphinxAtStartPar
\sphinxcode{\sphinxupquote{void AddTerms(Var{[}{]} vars, double{[}{]} coeffs)}}
\end{quote}

\sphinxAtStartPar
\sphinxstylestrong{Arguments}
\begin{quote}

\sphinxAtStartPar
\sphinxcode{\sphinxupquote{vars}}: variables for added terms.

\sphinxAtStartPar
\sphinxcode{\sphinxupquote{coeffs}}: coefficients array for added terms.
\end{quote}
\end{quote}

\subsubsection{Expr.AddTerms()}
\label{\detokenize{csapi/expr:id3}}\begin{quote}

\sphinxAtStartPar
Add terms to expression object.

\sphinxAtStartPar
\sphinxstylestrong{Synopsis}
\begin{quote}

\sphinxAtStartPar
\sphinxcode{\sphinxupquote{void AddTerms(VarArray vars, double coeff)}}
\end{quote}

\sphinxAtStartPar
\sphinxstylestrong{Arguments}
\begin{quote}

\sphinxAtStartPar
\sphinxcode{\sphinxupquote{vars}}: variables for added terms.

\sphinxAtStartPar
\sphinxcode{\sphinxupquote{coeff}}: coefficient array for added terms with default value 1.0.
\end{quote}
\end{quote}

\subsubsection{Expr.AddTerms()}
\label{\detokenize{csapi/expr:id4}}\begin{quote}

\sphinxAtStartPar
Add terms to expression object.

\sphinxAtStartPar
\sphinxstylestrong{Synopsis}
\begin{quote}

\sphinxAtStartPar
\sphinxcode{\sphinxupquote{void AddTerms(VarArray vars, double{[}{]} coeffs)}}
\end{quote}

\sphinxAtStartPar
\sphinxstylestrong{Arguments}
\begin{quote}

\sphinxAtStartPar
\sphinxcode{\sphinxupquote{vars}}: variables for added terms.

\sphinxAtStartPar
\sphinxcode{\sphinxupquote{coeffs}}: coefficients array for added terms.
\end{quote}
\end{quote}

\subsubsection{Expr.Clone()}
\label{\detokenize{csapi/expr:expr-clone}}\begin{quote}

\sphinxAtStartPar
Deep copy linear expression object.

\sphinxAtStartPar
\sphinxstylestrong{Synopsis}
\begin{quote}

\sphinxAtStartPar
\sphinxcode{\sphinxupquote{Expr Clone()}}
\end{quote}

\sphinxAtStartPar
\sphinxstylestrong{Return}
\begin{quote}

\sphinxAtStartPar
cloned linear expression object.
\end{quote}
\end{quote}

\subsubsection{Expr.Divide()}
\label{\detokenize{csapi/expr:expr-divide}}\begin{quote}

\sphinxAtStartPar
Divide itself by double constant.

\sphinxAtStartPar
\sphinxstylestrong{Synopsis}
\begin{quote}

\sphinxAtStartPar
\sphinxcode{\sphinxupquote{void Divide(double c)}}
\end{quote}

\sphinxAtStartPar
\sphinxstylestrong{Arguments}
\begin{quote}

\sphinxAtStartPar
\sphinxcode{\sphinxupquote{c}}: constant operand.
\end{quote}
\end{quote}

\subsubsection{Expr.Evaluate()}
\label{\detokenize{csapi/expr:expr-evaluate}}\begin{quote}

\sphinxAtStartPar
Evaluate linear expression after solving.

\sphinxAtStartPar
\sphinxstylestrong{Synopsis}
\begin{quote}

\sphinxAtStartPar
\sphinxcode{\sphinxupquote{double Evaluate()}}
\end{quote}

\sphinxAtStartPar
\sphinxstylestrong{Return}
\begin{quote}

\sphinxAtStartPar
value of linear expression.
\end{quote}
\end{quote}

\subsubsection{Expr.GetCoeff()}
\label{\detokenize{csapi/expr:expr-getcoeff}}\begin{quote}

\sphinxAtStartPar
Get coefficient from the i\sphinxhyphen{}th term in expression.

\sphinxAtStartPar
\sphinxstylestrong{Synopsis}
\begin{quote}

\sphinxAtStartPar
\sphinxcode{\sphinxupquote{double GetCoeff(int i)}}
\end{quote}

\sphinxAtStartPar
\sphinxstylestrong{Arguments}
\begin{quote}

\sphinxAtStartPar
\sphinxcode{\sphinxupquote{i}}: index of the term.
\end{quote}

\sphinxAtStartPar
\sphinxstylestrong{Return}
\begin{quote}

\sphinxAtStartPar
coefficient of the i\sphinxhyphen{}th term in expression object.
\end{quote}
\end{quote}

\subsubsection{Expr.GetConstant()}
\label{\detokenize{csapi/expr:expr-getconstant}}\begin{quote}

\sphinxAtStartPar
Get constant in expression.

\sphinxAtStartPar
\sphinxstylestrong{Synopsis}
\begin{quote}

\sphinxAtStartPar
\sphinxcode{\sphinxupquote{double GetConstant()}}
\end{quote}

\sphinxAtStartPar
\sphinxstylestrong{Return}
\begin{quote}

\sphinxAtStartPar
constant in expression.
\end{quote}
\end{quote}

\subsubsection{Expr.GetVar()}
\label{\detokenize{csapi/expr:expr-getvar}}\begin{quote}

\sphinxAtStartPar
Get variable from the i\sphinxhyphen{}th term in expression.

\sphinxAtStartPar
\sphinxstylestrong{Synopsis}
\begin{quote}

\sphinxAtStartPar
\sphinxcode{\sphinxupquote{Var GetVar(int i)}}
\end{quote}

\sphinxAtStartPar
\sphinxstylestrong{Arguments}
\begin{quote}

\sphinxAtStartPar
\sphinxcode{\sphinxupquote{i}}: index of the term.
\end{quote}

\sphinxAtStartPar
\sphinxstylestrong{Return}
\begin{quote}

\sphinxAtStartPar
variable of the i\sphinxhyphen{}th term in expression object.
\end{quote}
\end{quote}

\subsubsection{Expr.Multiply()}
\label{\detokenize{csapi/expr:expr-multiply}}\begin{quote}

\sphinxAtStartPar
Multiply itself by a constant.

\sphinxAtStartPar
\sphinxstylestrong{Synopsis}
\begin{quote}

\sphinxAtStartPar
\sphinxcode{\sphinxupquote{void Multiply(double c)}}
\end{quote}

\sphinxAtStartPar
\sphinxstylestrong{Arguments}
\begin{quote}

\sphinxAtStartPar
\sphinxcode{\sphinxupquote{c}}: constant operand.
\end{quote}
\end{quote}

\subsubsection{Expr.Remove()}
\label{\detokenize{csapi/expr:expr-remove}}\begin{quote}

\sphinxAtStartPar
Remove idx\sphinxhyphen{}th term from expression object.

\sphinxAtStartPar
\sphinxstylestrong{Synopsis}
\begin{quote}

\sphinxAtStartPar
\sphinxcode{\sphinxupquote{void Remove(int idx)}}
\end{quote}

\sphinxAtStartPar
\sphinxstylestrong{Arguments}
\begin{quote}

\sphinxAtStartPar
\sphinxcode{\sphinxupquote{idx}}: index of the term to be removed.
\end{quote}
\end{quote}

\subsubsection{Expr.Remove()}
\label{\detokenize{csapi/expr:id5}}\begin{quote}

\sphinxAtStartPar
Remove the term associated with variable from expression.

\sphinxAtStartPar
\sphinxstylestrong{Synopsis}
\begin{quote}

\sphinxAtStartPar
\sphinxcode{\sphinxupquote{void Remove(Var var)}}
\end{quote}

\sphinxAtStartPar
\sphinxstylestrong{Arguments}
\begin{quote}

\sphinxAtStartPar
\sphinxcode{\sphinxupquote{var}}: a variable whose term should be removed.
\end{quote}
\end{quote}

\subsubsection{Expr.SetCoeff()}
\label{\detokenize{csapi/expr:expr-setcoeff}}\begin{quote}

\sphinxAtStartPar
Set coefficient for the i\sphinxhyphen{}th term in expression.

\sphinxAtStartPar
\sphinxstylestrong{Synopsis}
\begin{quote}

\sphinxAtStartPar
\sphinxcode{\sphinxupquote{void SetCoeff(int i, double val)}}
\end{quote}

\sphinxAtStartPar
\sphinxstylestrong{Arguments}
\begin{quote}

\sphinxAtStartPar
\sphinxcode{\sphinxupquote{i}}: index of the term.

\sphinxAtStartPar
\sphinxcode{\sphinxupquote{val}}: coefficient of the term.
\end{quote}
\end{quote}

\subsubsection{Expr.SetConstant()}
\label{\detokenize{csapi/expr:expr-setconstant}}\begin{quote}

\sphinxAtStartPar
Set constant for the expression.

\sphinxAtStartPar
\sphinxstylestrong{Synopsis}
\begin{quote}

\sphinxAtStartPar
\sphinxcode{\sphinxupquote{void SetConstant(double constant)}}
\end{quote}

\sphinxAtStartPar
\sphinxstylestrong{Arguments}
\begin{quote}

\sphinxAtStartPar
\sphinxcode{\sphinxupquote{constant}}: the value of the constant.
\end{quote}
\end{quote}

\subsubsection{Expr.Size()}
\label{\detokenize{csapi/expr:expr-size}}\begin{quote}

\sphinxAtStartPar
Get number of terms in expression.

\sphinxAtStartPar
\sphinxstylestrong{Synopsis}
\begin{quote}

\sphinxAtStartPar
\sphinxcode{\sphinxupquote{long Size()}}
\end{quote}

\sphinxAtStartPar
\sphinxstylestrong{Return}
\begin{quote}

\sphinxAtStartPar
number of terms.
\end{quote}
\end{quote}

\subsection{Constraint}
\label{\detokenize{csharpapiref:constraint}}\label{\detokenize{csharpapiref:chapcsharpapiref-constraint}}
\sphinxAtStartPar
COPT constraint object. Constraints are always associated with a particular model.
User creates a constraint object by adding a constraint to a model,
rather than by using constructor of Constraint class.

\sphinxstepscope

\subsubsection{Constraint.Get()}
\label{\detokenize{csapi/constraint:constraint-get}}\label{\detokenize{csapi/constraint::doc}}\begin{quote}

\sphinxAtStartPar
Get information value of the constraint. Support informations of “Dual”, “Slack”, “LB”, “UB”.

\sphinxAtStartPar
\sphinxstylestrong{Synopsis}
\begin{quote}

\sphinxAtStartPar
\sphinxcode{\sphinxupquote{double Get(string info)}}
\end{quote}

\sphinxAtStartPar
\sphinxstylestrong{Arguments}
\begin{quote}

\sphinxAtStartPar
\sphinxcode{\sphinxupquote{info}}: name of information being queried.
\end{quote}

\sphinxAtStartPar
\sphinxstylestrong{Return}
\begin{quote}

\sphinxAtStartPar
information value.
\end{quote}
\end{quote}

\subsubsection{Constraint.GetBasis()}
\label{\detokenize{csapi/constraint:constraint-getbasis}}\begin{quote}

\sphinxAtStartPar
Get basis status of this constraint.

\sphinxAtStartPar
\sphinxstylestrong{Synopsis}
\begin{quote}

\sphinxAtStartPar
\sphinxcode{\sphinxupquote{int GetBasis()}}
\end{quote}

\sphinxAtStartPar
\sphinxstylestrong{Return}
\begin{quote}

\sphinxAtStartPar
basis status.
\end{quote}
\end{quote}

\subsubsection{Constraint.GetIdx()}
\label{\detokenize{csapi/constraint:constraint-getidx}}\begin{quote}

\sphinxAtStartPar
Get index of the constraint.

\sphinxAtStartPar
\sphinxstylestrong{Synopsis}
\begin{quote}

\sphinxAtStartPar
\sphinxcode{\sphinxupquote{int GetIdx()}}
\end{quote}

\sphinxAtStartPar
\sphinxstylestrong{Return}
\begin{quote}

\sphinxAtStartPar
the index of the constraint.
\end{quote}
\end{quote}

\subsubsection{Constraint.GetLowerIIS()}
\label{\detokenize{csapi/constraint:constraint-getloweriis}}\begin{quote}

\sphinxAtStartPar
Get IIS status for lower bound of the constraint.

\sphinxAtStartPar
\sphinxstylestrong{Synopsis}
\begin{quote}

\sphinxAtStartPar
\sphinxcode{\sphinxupquote{int GetLowerIIS()}}
\end{quote}

\sphinxAtStartPar
\sphinxstylestrong{Return}
\begin{quote}

\sphinxAtStartPar
IIS status.
\end{quote}
\end{quote}

\subsubsection{Constraint.GetName()}
\label{\detokenize{csapi/constraint:constraint-getname}}\begin{quote}

\sphinxAtStartPar
Get name of the constraint.

\sphinxAtStartPar
\sphinxstylestrong{Synopsis}
\begin{quote}

\sphinxAtStartPar
\sphinxcode{\sphinxupquote{string GetName()}}
\end{quote}

\sphinxAtStartPar
\sphinxstylestrong{Return}
\begin{quote}

\sphinxAtStartPar
the name of the constraint.
\end{quote}
\end{quote}

\subsubsection{Constraint.GetUpperIIS()}
\label{\detokenize{csapi/constraint:constraint-getupperiis}}\begin{quote}

\sphinxAtStartPar
Get IIS status for upper bound of the constraint.

\sphinxAtStartPar
\sphinxstylestrong{Synopsis}
\begin{quote}

\sphinxAtStartPar
\sphinxcode{\sphinxupquote{int GetUpperIIS()}}
\end{quote}

\sphinxAtStartPar
\sphinxstylestrong{Return}
\begin{quote}

\sphinxAtStartPar
IIS status.
\end{quote}
\end{quote}

\subsubsection{Constraint.Remove()}
\label{\detokenize{csapi/constraint:constraint-remove}}\begin{quote}

\sphinxAtStartPar
Remove this constraint from model.

\sphinxAtStartPar
\sphinxstylestrong{Synopsis}
\begin{quote}

\sphinxAtStartPar
\sphinxcode{\sphinxupquote{void Remove()}}
\end{quote}
\end{quote}

\subsubsection{Constraint.Set()}
\label{\detokenize{csapi/constraint:constraint-set}}\begin{quote}

\sphinxAtStartPar
Set information value of the constraint. Support informations of “LB” and “UB”.

\sphinxAtStartPar
\sphinxstylestrong{Synopsis}
\begin{quote}

\sphinxAtStartPar
\sphinxcode{\sphinxupquote{void Set(string info, double val)}}
\end{quote}

\sphinxAtStartPar
\sphinxstylestrong{Arguments}
\begin{quote}

\sphinxAtStartPar
\sphinxcode{\sphinxupquote{info}}: name of information.

\sphinxAtStartPar
\sphinxcode{\sphinxupquote{val}}: new information value.
\end{quote}
\end{quote}

\subsubsection{Constraint.SetName()}
\label{\detokenize{csapi/constraint:constraint-setname}}\begin{quote}

\sphinxAtStartPar
Set name for the constraint.

\sphinxAtStartPar
\sphinxstylestrong{Synopsis}
\begin{quote}

\sphinxAtStartPar
\sphinxcode{\sphinxupquote{void SetName(string name)}}
\end{quote}

\sphinxAtStartPar
\sphinxstylestrong{Arguments}
\begin{quote}

\sphinxAtStartPar
\sphinxcode{\sphinxupquote{name}}: the name to set.
\end{quote}
\end{quote}

\subsection{ConstrArray}
\label{\detokenize{csharpapiref:constrarray}}\label{\detokenize{csharpapiref:chapcsharpapiref-constrarray}}
\sphinxAtStartPar
COPT constraint array object. To store and access a set of C\# {\hyperref[\detokenize{csharpapiref:chapcsharpapiref-constraint}]{\sphinxcrossref{\DUrole{std,std-ref}{Constraint}}}}
objects, Cardinal Optimizer provides C\# ConstrArray class, which defines the following methods.

\sphinxstepscope

\subsubsection{ConstrArray.ConstrArray()}
\label{\detokenize{csapi/constrarray:constrarray-constrarray}}\label{\detokenize{csapi/constrarray::doc}}\begin{quote}

\sphinxAtStartPar
Constructor of constrarray object.

\sphinxAtStartPar
\sphinxstylestrong{Synopsis}
\begin{quote}

\sphinxAtStartPar
\sphinxcode{\sphinxupquote{ConstrArray()}}
\end{quote}
\end{quote}

\subsubsection{ConstrArray.GetConstr()}
\label{\detokenize{csapi/constrarray:constrarray-getconstr}}\begin{quote}

\sphinxAtStartPar
Get idx\sphinxhyphen{}th constraint object.

\sphinxAtStartPar
\sphinxstylestrong{Synopsis}
\begin{quote}

\sphinxAtStartPar
\sphinxcode{\sphinxupquote{Constraint GetConstr(int idx)}}
\end{quote}

\sphinxAtStartPar
\sphinxstylestrong{Arguments}
\begin{quote}

\sphinxAtStartPar
\sphinxcode{\sphinxupquote{idx}}: index of the constraint.
\end{quote}

\sphinxAtStartPar
\sphinxstylestrong{Return}
\begin{quote}

\sphinxAtStartPar
constraint object with index idx.
\end{quote}
\end{quote}

\subsubsection{ConstrArray.PushBack()}
\label{\detokenize{csapi/constrarray:constrarray-pushback}}\begin{quote}

\sphinxAtStartPar
Add a constraint object to constraint array.

\sphinxAtStartPar
\sphinxstylestrong{Synopsis}
\begin{quote}

\sphinxAtStartPar
\sphinxcode{\sphinxupquote{void PushBack(Constraint constr)}}
\end{quote}

\sphinxAtStartPar
\sphinxstylestrong{Arguments}
\begin{quote}

\sphinxAtStartPar
\sphinxcode{\sphinxupquote{constr}}: a constraint object.
\end{quote}
\end{quote}

\subsubsection{ConstrArray.Size()}
\label{\detokenize{csapi/constrarray:constrarray-size}}\begin{quote}

\sphinxAtStartPar
Get the number of constraint objects.

\sphinxAtStartPar
\sphinxstylestrong{Synopsis}
\begin{quote}

\sphinxAtStartPar
\sphinxcode{\sphinxupquote{int Size()}}
\end{quote}

\sphinxAtStartPar
\sphinxstylestrong{Return}
\begin{quote}

\sphinxAtStartPar
number of constraint objects.
\end{quote}
\end{quote}

\subsection{ConstrBuilder}
\label{\detokenize{csharpapiref:constrbuilder}}\label{\detokenize{csharpapiref:chapcsharpapiref-constrbuilder}}
\sphinxAtStartPar
COPT constraint builder object. To help building a constraint, given a linear expression,
constraint sense and right\sphinxhyphen{}hand side value, Cardinal Optimizer provides C\# ConstrBuilder
class, which defines the following methods.

\sphinxstepscope

\subsubsection{ConstrBuilder.ConstrBuilder()}
\label{\detokenize{csapi/constrbuilder:constrbuilder-constrbuilder}}\label{\detokenize{csapi/constrbuilder::doc}}\begin{quote}

\sphinxAtStartPar
Constructor of constrbuilder object.

\sphinxAtStartPar
\sphinxstylestrong{Synopsis}
\begin{quote}

\sphinxAtStartPar
\sphinxcode{\sphinxupquote{ConstrBuilder()}}
\end{quote}
\end{quote}

\subsubsection{ConstrBuilder.GetExpr()}
\label{\detokenize{csapi/constrbuilder:constrbuilder-getexpr}}\begin{quote}

\sphinxAtStartPar
Get expression associated with constraint.

\sphinxAtStartPar
\sphinxstylestrong{Synopsis}
\begin{quote}

\sphinxAtStartPar
\sphinxcode{\sphinxupquote{Expr GetExpr()}}
\end{quote}

\sphinxAtStartPar
\sphinxstylestrong{Return}
\begin{quote}

\sphinxAtStartPar
expression object.
\end{quote}
\end{quote}

\subsubsection{ConstrBuilder.GetRange()}
\label{\detokenize{csapi/constrbuilder:constrbuilder-getrange}}\begin{quote}

\sphinxAtStartPar
Get range from lower bound to upper bound of range constraint.

\sphinxAtStartPar
\sphinxstylestrong{Synopsis}
\begin{quote}

\sphinxAtStartPar
\sphinxcode{\sphinxupquote{double GetRange()}}
\end{quote}

\sphinxAtStartPar
\sphinxstylestrong{Return}
\begin{quote}

\sphinxAtStartPar
length from lower bound to upper bound of the constraint.
\end{quote}
\end{quote}

\subsubsection{ConstrBuilder.GetSense()}
\label{\detokenize{csapi/constrbuilder:constrbuilder-getsense}}\begin{quote}

\sphinxAtStartPar
Get sense associated with constraint.

\sphinxAtStartPar
\sphinxstylestrong{Synopsis}
\begin{quote}

\sphinxAtStartPar
\sphinxcode{\sphinxupquote{char GetSense()}}
\end{quote}

\sphinxAtStartPar
\sphinxstylestrong{Return}
\begin{quote}

\sphinxAtStartPar
constraint sense.
\end{quote}
\end{quote}

\subsubsection{ConstrBuilder.Set()}
\label{\detokenize{csapi/constrbuilder:constrbuilder-set}}\begin{quote}

\sphinxAtStartPar
Set detail of a constraint to its builder object.

\sphinxAtStartPar
\sphinxstylestrong{Synopsis}
\begin{quote}

\sphinxAtStartPar
\sphinxcode{\sphinxupquote{void Set(}}
\begin{quote}

\sphinxAtStartPar
\sphinxcode{\sphinxupquote{Expr expr,}}

\sphinxAtStartPar
\sphinxcode{\sphinxupquote{char sense,}}

\sphinxAtStartPar
\sphinxcode{\sphinxupquote{double rhs)}}
\end{quote}
\end{quote}

\sphinxAtStartPar
\sphinxstylestrong{Arguments}
\begin{quote}

\sphinxAtStartPar
\sphinxcode{\sphinxupquote{expr}}: expression object at one side of the constraint

\sphinxAtStartPar
\sphinxcode{\sphinxupquote{sense}}: constraint sense other than COPT\_RANGE.

\sphinxAtStartPar
\sphinxcode{\sphinxupquote{rhs}}: constant of right side of the constraint.
\end{quote}
\end{quote}

\subsubsection{ConstrBuilder.SetRange()}
\label{\detokenize{csapi/constrbuilder:constrbuilder-setrange}}\begin{quote}

\sphinxAtStartPar
Set a range constraint to its builder.

\sphinxAtStartPar
\sphinxstylestrong{Synopsis}
\begin{quote}

\sphinxAtStartPar
\sphinxcode{\sphinxupquote{void SetRange(Expr expr, double range)}}
\end{quote}

\sphinxAtStartPar
\sphinxstylestrong{Arguments}
\begin{quote}

\sphinxAtStartPar
\sphinxcode{\sphinxupquote{expr}}: expression object, whose constant is negative upper bound.

\sphinxAtStartPar
\sphinxcode{\sphinxupquote{range}}: length from lower bound to upper bound of the constraint. Must greater than 0.
\end{quote}
\end{quote}

\subsection{ConstrBuilderArray}
\label{\detokenize{csharpapiref:constrbuilderarray}}\label{\detokenize{csharpapiref:chapcsharpapiref-constrbuilderarray}}
\sphinxAtStartPar
COPT constraint builder array object. To store and access a set of C\# {\hyperref[\detokenize{csharpapiref:chapcsharpapiref-constrbuilder}]{\sphinxcrossref{\DUrole{std,std-ref}{ConstrBuilder}}}}
objects, Cardinal Optimizer provides C\# ConstrBuilderArray class, which defines the following methods.

\sphinxstepscope

\subsubsection{ConstrBuilderArray.ConstrBuilderArray()}
\label{\detokenize{csapi/constrbuilderarray:constrbuilderarray-constrbuilderarray}}\label{\detokenize{csapi/constrbuilderarray::doc}}\begin{quote}

\sphinxAtStartPar
Constructor of constrbuilderarray object.

\sphinxAtStartPar
\sphinxstylestrong{Synopsis}
\begin{quote}

\sphinxAtStartPar
\sphinxcode{\sphinxupquote{ConstrBuilderArray()}}
\end{quote}
\end{quote}

\subsubsection{ConstrBuilderArray.GetBuilder()}
\label{\detokenize{csapi/constrbuilderarray:constrbuilderarray-getbuilder}}\begin{quote}

\sphinxAtStartPar
Get idx\sphinxhyphen{}th constraint builder object.

\sphinxAtStartPar
\sphinxstylestrong{Synopsis}
\begin{quote}

\sphinxAtStartPar
\sphinxcode{\sphinxupquote{ConstrBuilder GetBuilder(int idx)}}
\end{quote}

\sphinxAtStartPar
\sphinxstylestrong{Arguments}
\begin{quote}

\sphinxAtStartPar
\sphinxcode{\sphinxupquote{idx}}: index of the constraint builder.
\end{quote}

\sphinxAtStartPar
\sphinxstylestrong{Return}
\begin{quote}

\sphinxAtStartPar
constraint builder object with index idx.
\end{quote}
\end{quote}

\subsubsection{ConstrBuilderArray.PushBack()}
\label{\detokenize{csapi/constrbuilderarray:constrbuilderarray-pushback}}\begin{quote}

\sphinxAtStartPar
Add a constraint builder object to constraint builder array.

\sphinxAtStartPar
\sphinxstylestrong{Synopsis}
\begin{quote}

\sphinxAtStartPar
\sphinxcode{\sphinxupquote{void PushBack(ConstrBuilder builder)}}
\end{quote}

\sphinxAtStartPar
\sphinxstylestrong{Arguments}
\begin{quote}

\sphinxAtStartPar
\sphinxcode{\sphinxupquote{builder}}: a constraint builder object.
\end{quote}
\end{quote}

\subsubsection{ConstrBuilderArray.Size()}
\label{\detokenize{csapi/constrbuilderarray:constrbuilderarray-size}}\begin{quote}

\sphinxAtStartPar
Get the number of constraint builder objects.

\sphinxAtStartPar
\sphinxstylestrong{Synopsis}
\begin{quote}

\sphinxAtStartPar
\sphinxcode{\sphinxupquote{int Size()}}
\end{quote}

\sphinxAtStartPar
\sphinxstylestrong{Return}
\begin{quote}

\sphinxAtStartPar
number of constraint builder objects.
\end{quote}
\end{quote}

\subsection{Column}
\label{\detokenize{csharpapiref:column}}\label{\detokenize{csharpapiref:chapcsharpapiref-column}}
\sphinxAtStartPar
COPT column object. A column consists of a list of constraints and associated coefficients.
Columns are used to represent the set of constraints in which a variable participates,
and the asssociated coefficents.

\sphinxstepscope

\subsubsection{Column.Column()}
\label{\detokenize{csapi/column:column-column}}\label{\detokenize{csapi/column::doc}}\begin{quote}

\sphinxAtStartPar
Constructor of column.

\sphinxAtStartPar
\sphinxstylestrong{Synopsis}
\begin{quote}

\sphinxAtStartPar
\sphinxcode{\sphinxupquote{Column()}}
\end{quote}
\end{quote}

\subsubsection{Column.AddColumn()}
\label{\detokenize{csapi/column:column-addcolumn}}\begin{quote}

\sphinxAtStartPar
Add a column to self.

\sphinxAtStartPar
\sphinxstylestrong{Synopsis}
\begin{quote}

\sphinxAtStartPar
\sphinxcode{\sphinxupquote{void AddColumn(Column col, double mult)}}
\end{quote}

\sphinxAtStartPar
\sphinxstylestrong{Arguments}
\begin{quote}

\sphinxAtStartPar
\sphinxcode{\sphinxupquote{col}}: column object to be added.

\sphinxAtStartPar
\sphinxcode{\sphinxupquote{mult}}: multiply constant.
\end{quote}
\end{quote}

\subsubsection{Column.AddTerm()}
\label{\detokenize{csapi/column:column-addterm}}\begin{quote}

\sphinxAtStartPar
Add a term to column object.

\sphinxAtStartPar
\sphinxstylestrong{Synopsis}
\begin{quote}

\sphinxAtStartPar
\sphinxcode{\sphinxupquote{void AddTerm(Constraint constr, double coeff)}}
\end{quote}

\sphinxAtStartPar
\sphinxstylestrong{Arguments}
\begin{quote}

\sphinxAtStartPar
\sphinxcode{\sphinxupquote{constr}}: a constraint for new term.

\sphinxAtStartPar
\sphinxcode{\sphinxupquote{coeff}}: coefficient for new term.
\end{quote}
\end{quote}

\subsubsection{Column.AddTerms()}
\label{\detokenize{csapi/column:column-addterms}}\begin{quote}

\sphinxAtStartPar
Add terms to column object.

\sphinxAtStartPar
\sphinxstylestrong{Synopsis}
\begin{quote}

\sphinxAtStartPar
\sphinxcode{\sphinxupquote{void AddTerms(Constraint{[}{]} constrs, double coeff)}}
\end{quote}

\sphinxAtStartPar
\sphinxstylestrong{Arguments}
\begin{quote}

\sphinxAtStartPar
\sphinxcode{\sphinxupquote{constrs}}: constraints for added terms.

\sphinxAtStartPar
\sphinxcode{\sphinxupquote{coeff}}: coefficient for added terms,default value is 1.
\end{quote}
\end{quote}

\subsubsection{Column.AddTerms()}
\label{\detokenize{csapi/column:id1}}\begin{quote}

\sphinxAtStartPar
Add terms to column object.

\sphinxAtStartPar
\sphinxstylestrong{Synopsis}
\begin{quote}

\sphinxAtStartPar
\sphinxcode{\sphinxupquote{void AddTerms(Constraint{[}{]} constrs, double{[}{]} coeffs)}}
\end{quote}

\sphinxAtStartPar
\sphinxstylestrong{Arguments}
\begin{quote}

\sphinxAtStartPar
\sphinxcode{\sphinxupquote{constrs}}: constraints for added terms.

\sphinxAtStartPar
\sphinxcode{\sphinxupquote{coeffs}}: coefficients for added terms.
\end{quote}
\end{quote}

\subsubsection{Column.AddTerms()}
\label{\detokenize{csapi/column:id2}}\begin{quote}

\sphinxAtStartPar
Add terms to column object.

\sphinxAtStartPar
\sphinxstylestrong{Synopsis}
\begin{quote}

\sphinxAtStartPar
\sphinxcode{\sphinxupquote{void AddTerms(ConstrArray constrs, double coeff)}}
\end{quote}

\sphinxAtStartPar
\sphinxstylestrong{Arguments}
\begin{quote}

\sphinxAtStartPar
\sphinxcode{\sphinxupquote{constrs}}: constraints for added terms.

\sphinxAtStartPar
\sphinxcode{\sphinxupquote{coeff}}: coefficient for added terms,default value is 1.
\end{quote}
\end{quote}

\subsubsection{Column.AddTerms()}
\label{\detokenize{csapi/column:id3}}\begin{quote}

\sphinxAtStartPar
Add terms to column object.

\sphinxAtStartPar
\sphinxstylestrong{Synopsis}
\begin{quote}

\sphinxAtStartPar
\sphinxcode{\sphinxupquote{void AddTerms(ConstrArray constrs, double{[}{]} coeffs)}}
\end{quote}

\sphinxAtStartPar
\sphinxstylestrong{Arguments}
\begin{quote}

\sphinxAtStartPar
\sphinxcode{\sphinxupquote{constrs}}: constraints for added terms.

\sphinxAtStartPar
\sphinxcode{\sphinxupquote{coeffs}}: coefficients for added terms.
\end{quote}
\end{quote}

\subsubsection{Column.Clear()}
\label{\detokenize{csapi/column:column-clear}}\begin{quote}

\sphinxAtStartPar
Clear all terms.

\sphinxAtStartPar
\sphinxstylestrong{Synopsis}
\begin{quote}

\sphinxAtStartPar
\sphinxcode{\sphinxupquote{void Clear()}}
\end{quote}
\end{quote}

\subsubsection{Column.Clone()}
\label{\detokenize{csapi/column:column-clone}}\begin{quote}

\sphinxAtStartPar
Deep copy column object.

\sphinxAtStartPar
\sphinxstylestrong{Synopsis}
\begin{quote}

\sphinxAtStartPar
\sphinxcode{\sphinxupquote{Column Clone()}}
\end{quote}

\sphinxAtStartPar
\sphinxstylestrong{Return}
\begin{quote}

\sphinxAtStartPar
cloned column object.
\end{quote}
\end{quote}

\subsubsection{Column.GetCoeff()}
\label{\detokenize{csapi/column:column-getcoeff}}\begin{quote}

\sphinxAtStartPar
Get coefficient from the i\sphinxhyphen{}th term in column object.

\sphinxAtStartPar
\sphinxstylestrong{Synopsis}
\begin{quote}

\sphinxAtStartPar
\sphinxcode{\sphinxupquote{double GetCoeff(int i)}}
\end{quote}

\sphinxAtStartPar
\sphinxstylestrong{Arguments}
\begin{quote}

\sphinxAtStartPar
\sphinxcode{\sphinxupquote{i}}: index of the term.
\end{quote}

\sphinxAtStartPar
\sphinxstylestrong{Return}
\begin{quote}

\sphinxAtStartPar
coefficient of the i\sphinxhyphen{}th term in column object.
\end{quote}
\end{quote}

\subsubsection{Column.GetConstr()}
\label{\detokenize{csapi/column:column-getconstr}}\begin{quote}

\sphinxAtStartPar
Get constraint from the i\sphinxhyphen{}th term in column object.

\sphinxAtStartPar
\sphinxstylestrong{Synopsis}
\begin{quote}

\sphinxAtStartPar
\sphinxcode{\sphinxupquote{Constraint GetConstr(int i)}}
\end{quote}

\sphinxAtStartPar
\sphinxstylestrong{Arguments}
\begin{quote}

\sphinxAtStartPar
\sphinxcode{\sphinxupquote{i}}: index of the term.
\end{quote}

\sphinxAtStartPar
\sphinxstylestrong{Return}
\begin{quote}

\sphinxAtStartPar
constraint of the i\sphinxhyphen{}th term in column object.
\end{quote}
\end{quote}

\subsubsection{Column.Remove()}
\label{\detokenize{csapi/column:column-remove}}\begin{quote}

\sphinxAtStartPar
Remove idx\sphinxhyphen{}th term from column object.

\sphinxAtStartPar
\sphinxstylestrong{Synopsis}
\begin{quote}

\sphinxAtStartPar
\sphinxcode{\sphinxupquote{void Remove(int idx)}}
\end{quote}

\sphinxAtStartPar
\sphinxstylestrong{Arguments}
\begin{quote}

\sphinxAtStartPar
\sphinxcode{\sphinxupquote{idx}}: index of the term to be removed.
\end{quote}
\end{quote}

\subsubsection{Column.Remove()}
\label{\detokenize{csapi/column:id4}}\begin{quote}

\sphinxAtStartPar
Remove the term associated with constraint from column object.

\sphinxAtStartPar
\sphinxstylestrong{Synopsis}
\begin{quote}

\sphinxAtStartPar
\sphinxcode{\sphinxupquote{void Remove(Constraint constr)}}
\end{quote}

\sphinxAtStartPar
\sphinxstylestrong{Arguments}
\begin{quote}

\sphinxAtStartPar
\sphinxcode{\sphinxupquote{constr}}: a constraint whose term should be removed.
\end{quote}
\end{quote}

\subsubsection{Column.Size()}
\label{\detokenize{csapi/column:column-size}}\begin{quote}

\sphinxAtStartPar
Get number of terms in column object.

\sphinxAtStartPar
\sphinxstylestrong{Synopsis}
\begin{quote}

\sphinxAtStartPar
\sphinxcode{\sphinxupquote{int Size()}}
\end{quote}

\sphinxAtStartPar
\sphinxstylestrong{Return}
\begin{quote}

\sphinxAtStartPar
number of terms.
\end{quote}
\end{quote}

\subsection{ColumnArray}
\label{\detokenize{csharpapiref:columnarray}}\label{\detokenize{csharpapiref:chapcsharpapiref-columnarray}}
\sphinxAtStartPar
COPT column array object. To store and access a set of C\#
{\hyperref[\detokenize{csharpapiref:chapcsharpapiref-column}]{\sphinxcrossref{\DUrole{std,std-ref}{Column}}}} objects, Cardinal Optimizer provides C\#
ColumnArray class, which defines the following methods.

\sphinxstepscope

\subsubsection{ColumnArray.ColumnArray()}
\label{\detokenize{csapi/columnarray:columnarray-columnarray}}\label{\detokenize{csapi/columnarray::doc}}\begin{quote}

\sphinxAtStartPar
Constructor of columnarray object.

\sphinxAtStartPar
\sphinxstylestrong{Synopsis}
\begin{quote}

\sphinxAtStartPar
\sphinxcode{\sphinxupquote{ColumnArray()}}
\end{quote}
\end{quote}

\subsubsection{ColumnArray.Clear()}
\label{\detokenize{csapi/columnarray:columnarray-clear}}\begin{quote}

\sphinxAtStartPar
Clear all column objects.

\sphinxAtStartPar
\sphinxstylestrong{Synopsis}
\begin{quote}

\sphinxAtStartPar
\sphinxcode{\sphinxupquote{void Clear()}}
\end{quote}
\end{quote}

\subsubsection{ColumnArray.GetColumn()}
\label{\detokenize{csapi/columnarray:columnarray-getcolumn}}\begin{quote}

\sphinxAtStartPar
Get idx\sphinxhyphen{}th column object.

\sphinxAtStartPar
\sphinxstylestrong{Synopsis}
\begin{quote}

\sphinxAtStartPar
\sphinxcode{\sphinxupquote{Column GetColumn(int idx)}}
\end{quote}

\sphinxAtStartPar
\sphinxstylestrong{Arguments}
\begin{quote}

\sphinxAtStartPar
\sphinxcode{\sphinxupquote{idx}}: index of the column.
\end{quote}

\sphinxAtStartPar
\sphinxstylestrong{Return}
\begin{quote}

\sphinxAtStartPar
column object with index idx.
\end{quote}
\end{quote}

\subsubsection{ColumnArray.PushBack()}
\label{\detokenize{csapi/columnarray:columnarray-pushback}}\begin{quote}

\sphinxAtStartPar
Add a column object to column array.

\sphinxAtStartPar
\sphinxstylestrong{Synopsis}
\begin{quote}

\sphinxAtStartPar
\sphinxcode{\sphinxupquote{void PushBack(Column col)}}
\end{quote}

\sphinxAtStartPar
\sphinxstylestrong{Arguments}
\begin{quote}

\sphinxAtStartPar
\sphinxcode{\sphinxupquote{col}}: a column object.
\end{quote}
\end{quote}

\subsubsection{ColumnArray.Size()}
\label{\detokenize{csapi/columnarray:columnarray-size}}\begin{quote}

\sphinxAtStartPar
Get the number of column objects.

\sphinxAtStartPar
\sphinxstylestrong{Synopsis}
\begin{quote}

\sphinxAtStartPar
\sphinxcode{\sphinxupquote{int Size()}}
\end{quote}

\sphinxAtStartPar
\sphinxstylestrong{Return}
\begin{quote}

\sphinxAtStartPar
number of column objects.
\end{quote}
\end{quote}

\subsection{Sos}
\label{\detokenize{csharpapiref:sos}}\label{\detokenize{csharpapiref:chapcsharpapiref-sos}}
\sphinxAtStartPar
COPT SOS constraint object. SOS constraints are always associated with a particular model.
User creates an SOS constraint object by adding an SOS constraint to a model,
rather than by using constructor of Sos class.

\sphinxAtStartPar
An SOS constraint can be type 1 or 2 (\sphinxcode{\sphinxupquote{COPT\_SOS\_TYPE1}} or \sphinxcode{\sphinxupquote{COPT\_SOS\_TYPE2}}).

\sphinxstepscope

\subsubsection{Sos.GetIdx()}
\label{\detokenize{csapi/sos:sos-getidx}}\label{\detokenize{csapi/sos::doc}}\begin{quote}

\sphinxAtStartPar
Get the index of SOS constraint.

\sphinxAtStartPar
\sphinxstylestrong{Synopsis}
\begin{quote}

\sphinxAtStartPar
\sphinxcode{\sphinxupquote{int GetIdx()}}
\end{quote}

\sphinxAtStartPar
\sphinxstylestrong{Return}
\begin{quote}

\sphinxAtStartPar
index of SOS constraint.
\end{quote}
\end{quote}

\subsubsection{Sos.GetIIS()}
\label{\detokenize{csapi/sos:sos-getiis}}\begin{quote}

\sphinxAtStartPar
Get IIS status of the SOS constraint.

\sphinxAtStartPar
\sphinxstylestrong{Synopsis}
\begin{quote}

\sphinxAtStartPar
\sphinxcode{\sphinxupquote{int GetIIS()}}
\end{quote}

\sphinxAtStartPar
\sphinxstylestrong{Return}
\begin{quote}

\sphinxAtStartPar
IIS status.
\end{quote}
\end{quote}

\subsubsection{Sos.Remove()}
\label{\detokenize{csapi/sos:sos-remove}}\begin{quote}

\sphinxAtStartPar
Remove the SOS constraint from model.

\sphinxAtStartPar
\sphinxstylestrong{Synopsis}
\begin{quote}

\sphinxAtStartPar
\sphinxcode{\sphinxupquote{void Remove()}}
\end{quote}
\end{quote}

\subsection{SosArray}
\label{\detokenize{csharpapiref:sosarray}}\label{\detokenize{csharpapiref:chapcsharpapiref-sosarray}}
\sphinxAtStartPar
COPT SOS constraint array object. To store and access a set of C\# {\hyperref[\detokenize{csharpapiref:chapcsharpapiref-sos}]{\sphinxcrossref{\DUrole{std,std-ref}{Sos}}}}
objects, Cardinal Optimizer provides C\# SosArray class, which defines the following methods.

\sphinxstepscope

\subsubsection{SosArray.SosArray()}
\label{\detokenize{csapi/sosarray:sosarray-sosarray}}\label{\detokenize{csapi/sosarray::doc}}\begin{quote}

\sphinxAtStartPar
Constructor of sosarray object.

\sphinxAtStartPar
\sphinxstylestrong{Synopsis}
\begin{quote}

\sphinxAtStartPar
\sphinxcode{\sphinxupquote{SosArray()}}
\end{quote}
\end{quote}

\subsubsection{SosArray.GetSos()}
\label{\detokenize{csapi/sosarray:sosarray-getsos}}\begin{quote}

\sphinxAtStartPar
Get idx\sphinxhyphen{}th SOS object.

\sphinxAtStartPar
\sphinxstylestrong{Synopsis}
\begin{quote}

\sphinxAtStartPar
\sphinxcode{\sphinxupquote{Sos GetSos(int idx)}}
\end{quote}

\sphinxAtStartPar
\sphinxstylestrong{Arguments}
\begin{quote}

\sphinxAtStartPar
\sphinxcode{\sphinxupquote{idx}}: index of SOS.
\end{quote}

\sphinxAtStartPar
\sphinxstylestrong{Return}
\begin{quote}

\sphinxAtStartPar
SOS object with index idx.
\end{quote}
\end{quote}

\subsubsection{SosArray.PushBack()}
\label{\detokenize{csapi/sosarray:sosarray-pushback}}\begin{quote}

\sphinxAtStartPar
Add a SOS constraint object to SOS constraint array.

\sphinxAtStartPar
\sphinxstylestrong{Synopsis}
\begin{quote}

\sphinxAtStartPar
\sphinxcode{\sphinxupquote{void PushBack(Sos sos)}}
\end{quote}

\sphinxAtStartPar
\sphinxstylestrong{Arguments}
\begin{quote}

\sphinxAtStartPar
\sphinxcode{\sphinxupquote{sos}}: a SOS constraint object.
\end{quote}
\end{quote}

\subsubsection{SosArray.Size()}
\label{\detokenize{csapi/sosarray:sosarray-size}}\begin{quote}

\sphinxAtStartPar
Get the number of SOS constraint objects.

\sphinxAtStartPar
\sphinxstylestrong{Synopsis}
\begin{quote}

\sphinxAtStartPar
\sphinxcode{\sphinxupquote{int Size()}}
\end{quote}

\sphinxAtStartPar
\sphinxstylestrong{Return}
\begin{quote}

\sphinxAtStartPar
number of SOS constraint objects.
\end{quote}
\end{quote}

\subsection{SosBuilder}
\label{\detokenize{csharpapiref:sosbuilder}}\label{\detokenize{csharpapiref:chapcsharpapiref-sosbuilder}}
\sphinxAtStartPar
COPT SOS constraint builder object. To help building an SOS constraint, given the SOS type,
a set of variables and associated weights, Cardinal Optimizer provides C\# SosBuilder
class, which defines the following methods.

\sphinxstepscope

\subsubsection{SosBuilder.SosBuilder()}
\label{\detokenize{csapi/sosbuilder:sosbuilder-sosbuilder}}\label{\detokenize{csapi/sosbuilder::doc}}\begin{quote}

\sphinxAtStartPar
Constructor of sosbuilder object.

\sphinxAtStartPar
\sphinxstylestrong{Synopsis}
\begin{quote}

\sphinxAtStartPar
\sphinxcode{\sphinxupquote{SosBuilder()}}
\end{quote}
\end{quote}

\subsubsection{SosBuilder.GetSize()}
\label{\detokenize{csapi/sosbuilder:sosbuilder-getsize}}\begin{quote}

\sphinxAtStartPar
Get number of terms in SOS constraint.

\sphinxAtStartPar
\sphinxstylestrong{Synopsis}
\begin{quote}

\sphinxAtStartPar
\sphinxcode{\sphinxupquote{int GetSize()}}
\end{quote}

\sphinxAtStartPar
\sphinxstylestrong{Return}
\begin{quote}

\sphinxAtStartPar
number of terms.
\end{quote}
\end{quote}

\subsubsection{SosBuilder.GetType()}
\label{\detokenize{csapi/sosbuilder:sosbuilder-gettype}}\begin{quote}

\sphinxAtStartPar
Get type of SOS constraint.

\sphinxAtStartPar
\sphinxstylestrong{Synopsis}
\begin{quote}

\sphinxAtStartPar
\sphinxcode{\sphinxupquote{int GetType()}}
\end{quote}

\sphinxAtStartPar
\sphinxstylestrong{Return}
\begin{quote}

\sphinxAtStartPar
type of SOS constraint.
\end{quote}
\end{quote}

\subsubsection{SosBuilder.GetVar()}
\label{\detokenize{csapi/sosbuilder:sosbuilder-getvar}}\begin{quote}

\sphinxAtStartPar
Get variable from the idx\sphinxhyphen{}th term in SOS constraint.

\sphinxAtStartPar
\sphinxstylestrong{Synopsis}
\begin{quote}

\sphinxAtStartPar
\sphinxcode{\sphinxupquote{Var GetVar(int idx)}}
\end{quote}

\sphinxAtStartPar
\sphinxstylestrong{Arguments}
\begin{quote}

\sphinxAtStartPar
\sphinxcode{\sphinxupquote{idx}}: index of the term.
\end{quote}

\sphinxAtStartPar
\sphinxstylestrong{Return}
\begin{quote}

\sphinxAtStartPar
variable of the idx\sphinxhyphen{}th term in SOS constraint.
\end{quote}
\end{quote}

\subsubsection{SosBuilder.GetVars()}
\label{\detokenize{csapi/sosbuilder:sosbuilder-getvars}}\begin{quote}

\sphinxAtStartPar
Get all variables in a SOS constraint.

\sphinxAtStartPar
\sphinxstylestrong{Synopsis}
\begin{quote}

\sphinxAtStartPar
\sphinxcode{\sphinxupquote{VarArray GetVars()}}
\end{quote}

\sphinxAtStartPar
\sphinxstylestrong{Return}
\begin{quote}

\sphinxAtStartPar
variables in a SOS constraint.
\end{quote}
\end{quote}

\subsubsection{SosBuilder.GetWeight()}
\label{\detokenize{csapi/sosbuilder:sosbuilder-getweight}}\begin{quote}

\sphinxAtStartPar
Get weight from the idx\sphinxhyphen{}th term in SOS constraint.

\sphinxAtStartPar
\sphinxstylestrong{Synopsis}
\begin{quote}

\sphinxAtStartPar
\sphinxcode{\sphinxupquote{double GetWeight(int idx)}}
\end{quote}

\sphinxAtStartPar
\sphinxstylestrong{Arguments}
\begin{quote}

\sphinxAtStartPar
\sphinxcode{\sphinxupquote{idx}}: index of the term.
\end{quote}

\sphinxAtStartPar
\sphinxstylestrong{Return}
\begin{quote}

\sphinxAtStartPar
weight of the idx\sphinxhyphen{}th term in SOS constraint.
\end{quote}
\end{quote}

\subsubsection{SosBuilder.GetWeights()}
\label{\detokenize{csapi/sosbuilder:sosbuilder-getweights}}\begin{quote}

\sphinxAtStartPar
Get weights of all terms in SOS constraint.

\sphinxAtStartPar
\sphinxstylestrong{Synopsis}
\begin{quote}

\sphinxAtStartPar
\sphinxcode{\sphinxupquote{double{[}{]} GetWeights()}}
\end{quote}

\sphinxAtStartPar
\sphinxstylestrong{Return}
\begin{quote}

\sphinxAtStartPar
array of weights.
\end{quote}
\end{quote}

\subsubsection{SosBuilder.Set()}
\label{\detokenize{csapi/sosbuilder:sosbuilder-set}}\begin{quote}

\sphinxAtStartPar
Set variables and weights of SOS constraint.

\sphinxAtStartPar
\sphinxstylestrong{Synopsis}
\begin{quote}

\sphinxAtStartPar
\sphinxcode{\sphinxupquote{void Set(}}
\begin{quote}

\sphinxAtStartPar
\sphinxcode{\sphinxupquote{VarArray vars,}}

\sphinxAtStartPar
\sphinxcode{\sphinxupquote{double{[}{]} weights,}}

\sphinxAtStartPar
\sphinxcode{\sphinxupquote{int type)}}
\end{quote}
\end{quote}

\sphinxAtStartPar
\sphinxstylestrong{Arguments}
\begin{quote}

\sphinxAtStartPar
\sphinxcode{\sphinxupquote{vars}}: variable array object.

\sphinxAtStartPar
\sphinxcode{\sphinxupquote{weights}}: pointer to array of weights.

\sphinxAtStartPar
\sphinxcode{\sphinxupquote{type}}: type of SOS constraint.
\end{quote}
\end{quote}

\subsection{SosBuilderArray}
\label{\detokenize{csharpapiref:sosbuilderarray}}\label{\detokenize{csharpapiref:chapcsharpapiref-sosbuilderarray}}
\sphinxAtStartPar
COPT SOS constraint builder array object. To store and access a set of C\# {\hyperref[\detokenize{csharpapiref:chapcsharpapiref-sosbuilder}]{\sphinxcrossref{\DUrole{std,std-ref}{SosBuilder}}}}
objects, Cardinal Optimizer provides C\# SosBuilderArray class, which defines the following methods.

\sphinxstepscope

\subsubsection{SosBuilderArray.SosBuilderArray()}
\label{\detokenize{csapi/sosbuilderarray:sosbuilderarray-sosbuilderarray}}\label{\detokenize{csapi/sosbuilderarray::doc}}\begin{quote}

\sphinxAtStartPar
Constructor of sosbuilderarray object.

\sphinxAtStartPar
\sphinxstylestrong{Synopsis}
\begin{quote}

\sphinxAtStartPar
\sphinxcode{\sphinxupquote{SosBuilderArray()}}
\end{quote}
\end{quote}

\subsubsection{SosBuilderArray.GetBuilder()}
\label{\detokenize{csapi/sosbuilderarray:sosbuilderarray-getbuilder}}\begin{quote}

\sphinxAtStartPar
Get idx\sphinxhyphen{}th SOS constraint builder object.

\sphinxAtStartPar
\sphinxstylestrong{Synopsis}
\begin{quote}

\sphinxAtStartPar
\sphinxcode{\sphinxupquote{SosBuilder GetBuilder(int idx)}}
\end{quote}

\sphinxAtStartPar
\sphinxstylestrong{Arguments}
\begin{quote}

\sphinxAtStartPar
\sphinxcode{\sphinxupquote{idx}}: index of the SOS constraint builder.
\end{quote}

\sphinxAtStartPar
\sphinxstylestrong{Return}
\begin{quote}

\sphinxAtStartPar
SOS constraint builder object with index idx.
\end{quote}
\end{quote}

\subsubsection{SosBuilderArray.PushBack()}
\label{\detokenize{csapi/sosbuilderarray:sosbuilderarray-pushback}}\begin{quote}

\sphinxAtStartPar
Add a SOS constraint builder object to SOS constraint builder array.

\sphinxAtStartPar
\sphinxstylestrong{Synopsis}
\begin{quote}

\sphinxAtStartPar
\sphinxcode{\sphinxupquote{void PushBack(SosBuilder builder)}}
\end{quote}

\sphinxAtStartPar
\sphinxstylestrong{Arguments}
\begin{quote}

\sphinxAtStartPar
\sphinxcode{\sphinxupquote{builder}}: a SOS constraint builder object.
\end{quote}
\end{quote}

\subsubsection{SosBuilderArray.Size()}
\label{\detokenize{csapi/sosbuilderarray:sosbuilderarray-size}}\begin{quote}

\sphinxAtStartPar
Get the number of SOS constraint builder objects.

\sphinxAtStartPar
\sphinxstylestrong{Synopsis}
\begin{quote}

\sphinxAtStartPar
\sphinxcode{\sphinxupquote{int Size()}}
\end{quote}

\sphinxAtStartPar
\sphinxstylestrong{Return}
\begin{quote}

\sphinxAtStartPar
number of SOS constraint builder objects.
\end{quote}
\end{quote}

\subsection{GenConstr}
\label{\detokenize{csharpapiref:genconstr}}\label{\detokenize{csharpapiref:chapcsharpapiref-genconstr}}
\sphinxAtStartPar
COPT general constraint object. General constraints are always associated with a particular model.
User creates a general constraint object by adding a general constraint to a model,
rather than by using constructor of GenConstr class.

\sphinxstepscope

\subsubsection{GenConstr.GetIdx()}
\label{\detokenize{csapi/genconstr:genconstr-getidx}}\label{\detokenize{csapi/genconstr::doc}}\begin{quote}

\sphinxAtStartPar
Get the index of the general constraint.

\sphinxAtStartPar
\sphinxstylestrong{Synopsis}
\begin{quote}

\sphinxAtStartPar
\sphinxcode{\sphinxupquote{int GetIdx()}}
\end{quote}

\sphinxAtStartPar
\sphinxstylestrong{Return}
\begin{quote}

\sphinxAtStartPar
index of the general constraint.
\end{quote}
\end{quote}

\subsubsection{GenConstr.GetIIS()}
\label{\detokenize{csapi/genconstr:genconstr-getiis}}\begin{quote}

\sphinxAtStartPar
Get IIS status of the general constraint.

\sphinxAtStartPar
\sphinxstylestrong{Synopsis}
\begin{quote}

\sphinxAtStartPar
\sphinxcode{\sphinxupquote{int GetIIS()}}
\end{quote}

\sphinxAtStartPar
\sphinxstylestrong{Return}
\begin{quote}

\sphinxAtStartPar
IIS status.
\end{quote}
\end{quote}

\subsubsection{GenConstr.GetName()}
\label{\detokenize{csapi/genconstr:genconstr-getname}}\begin{quote}

\sphinxAtStartPar
Get name of general constraint.

\sphinxAtStartPar
\sphinxstylestrong{Synopsis}
\begin{quote}

\sphinxAtStartPar
\sphinxcode{\sphinxupquote{string GetName()}}
\end{quote}

\sphinxAtStartPar
\sphinxstylestrong{Return}
\begin{quote}

\sphinxAtStartPar
the name of general constraint.
\end{quote}
\end{quote}

\subsubsection{GenConstr.Remove()}
\label{\detokenize{csapi/genconstr:genconstr-remove}}\begin{quote}

\sphinxAtStartPar
Remove the general constraint from model.

\sphinxAtStartPar
\sphinxstylestrong{Synopsis}
\begin{quote}

\sphinxAtStartPar
\sphinxcode{\sphinxupquote{void Remove()}}
\end{quote}
\end{quote}

\subsubsection{GenConstr.SetName()}
\label{\detokenize{csapi/genconstr:genconstr-setname}}\begin{quote}

\sphinxAtStartPar
Set name for general constraint.

\sphinxAtStartPar
\sphinxstylestrong{Synopsis}
\begin{quote}

\sphinxAtStartPar
\sphinxcode{\sphinxupquote{void SetName(string name)}}
\end{quote}

\sphinxAtStartPar
\sphinxstylestrong{Arguments}
\begin{quote}

\sphinxAtStartPar
\sphinxcode{\sphinxupquote{name}}: the name to set.
\end{quote}
\end{quote}

\subsection{GenConstrArray}
\label{\detokenize{csharpapiref:genconstrarray}}\label{\detokenize{csharpapiref:chapcsharpapiref-genconstrarray}}
\sphinxAtStartPar
COPT general constraint array object. To store and access a set of C\# {\hyperref[\detokenize{csharpapiref:chapcsharpapiref-genconstr}]{\sphinxcrossref{\DUrole{std,std-ref}{GenConstr}}}}
objects, Cardinal Optimizer provides C\# GenConstrArray class, which defines the following methods.

\sphinxstepscope

\subsubsection{GenConstrArray.GenConstrArray()}
\label{\detokenize{csapi/genconstrarray:genconstrarray-genconstrarray}}\label{\detokenize{csapi/genconstrarray::doc}}\begin{quote}

\sphinxAtStartPar
Constructor of genconstrarray.

\sphinxAtStartPar
\sphinxstylestrong{Synopsis}
\begin{quote}

\sphinxAtStartPar
\sphinxcode{\sphinxupquote{GenConstrArray()}}
\end{quote}
\end{quote}

\subsubsection{GenConstrArray.GetGenConstr()}
\label{\detokenize{csapi/genconstrarray:genconstrarray-getgenconstr}}\begin{quote}

\sphinxAtStartPar
Get idx\sphinxhyphen{}th general constraint object.

\sphinxAtStartPar
\sphinxstylestrong{Synopsis}
\begin{quote}

\sphinxAtStartPar
\sphinxcode{\sphinxupquote{GenConstr GetGenConstr(int idx)}}
\end{quote}

\sphinxAtStartPar
\sphinxstylestrong{Arguments}
\begin{quote}

\sphinxAtStartPar
\sphinxcode{\sphinxupquote{idx}}: index of the general constraint.
\end{quote}

\sphinxAtStartPar
\sphinxstylestrong{Return}
\begin{quote}

\sphinxAtStartPar
general constraint object with index idx.
\end{quote}
\end{quote}

\subsubsection{GenConstrArray.PushBack()}
\label{\detokenize{csapi/genconstrarray:genconstrarray-pushback}}\begin{quote}

\sphinxAtStartPar
Add a general constraint object to general constraint array.

\sphinxAtStartPar
\sphinxstylestrong{Synopsis}
\begin{quote}

\sphinxAtStartPar
\sphinxcode{\sphinxupquote{void PushBack(GenConstr genconstr)}}
\end{quote}

\sphinxAtStartPar
\sphinxstylestrong{Arguments}
\begin{quote}

\sphinxAtStartPar
\sphinxcode{\sphinxupquote{genconstr}}: a general constraint object.
\end{quote}
\end{quote}

\subsubsection{GenConstrArray.Reserve()}
\label{\detokenize{csapi/genconstrarray:genconstrarray-reserve}}\begin{quote}

\sphinxAtStartPar
Reserve capacity to contain at least n items.

\sphinxAtStartPar
\sphinxstylestrong{Synopsis}
\begin{quote}

\sphinxAtStartPar
\sphinxcode{\sphinxupquote{void Reserve(int n)}}
\end{quote}

\sphinxAtStartPar
\sphinxstylestrong{Arguments}
\begin{quote}

\sphinxAtStartPar
\sphinxcode{\sphinxupquote{n}}: capacity number of general constraint object.
\end{quote}
\end{quote}

\subsubsection{GenConstrArray.Size()}
\label{\detokenize{csapi/genconstrarray:genconstrarray-size}}\begin{quote}

\sphinxAtStartPar
Get the number of general constraint objects.

\sphinxAtStartPar
\sphinxstylestrong{Synopsis}
\begin{quote}

\sphinxAtStartPar
\sphinxcode{\sphinxupquote{int Size()}}
\end{quote}

\sphinxAtStartPar
\sphinxstylestrong{Return}
\begin{quote}

\sphinxAtStartPar
number of general constraint objects.
\end{quote}
\end{quote}

\subsection{GenConstrBuilder}
\label{\detokenize{csharpapiref:genconstrbuilder}}\label{\detokenize{csharpapiref:chapcsharpapiref-genconstrbuilder}}
\sphinxAtStartPar
COPT general constraint builder object. To help building a general constraint, given a binary variable
and associated value, a linear expression and constraint sense, Cardinal Optimizer provides C\#
GenConstrBuilder class, which defines the following methods.

\sphinxstepscope

\subsubsection{GenConstrBuilder.GenConstrBuilder()}
\label{\detokenize{csapi/genconstrbuilder:genconstrbuilder-genconstrbuilder}}\label{\detokenize{csapi/genconstrbuilder::doc}}\begin{quote}

\sphinxAtStartPar
Constructor of genconstrbuilder.

\sphinxAtStartPar
\sphinxstylestrong{Synopsis}
\begin{quote}

\sphinxAtStartPar
\sphinxcode{\sphinxupquote{GenConstrBuilder()}}
\end{quote}
\end{quote}

\subsubsection{GenConstrBuilder.GetBinVal()}
\label{\detokenize{csapi/genconstrbuilder:genconstrbuilder-getbinval}}\begin{quote}

\sphinxAtStartPar
Get binary value associated with general constraint.

\sphinxAtStartPar
\sphinxstylestrong{Synopsis}
\begin{quote}

\sphinxAtStartPar
\sphinxcode{\sphinxupquote{int GetBinVal()}}
\end{quote}

\sphinxAtStartPar
\sphinxstylestrong{Return}
\begin{quote}

\sphinxAtStartPar
binary value.
\end{quote}
\end{quote}

\subsubsection{GenConstrBuilder.GetBinVar()}
\label{\detokenize{csapi/genconstrbuilder:genconstrbuilder-getbinvar}}\begin{quote}

\sphinxAtStartPar
Get binary variable associated with general constraint.

\sphinxAtStartPar
\sphinxstylestrong{Synopsis}
\begin{quote}

\sphinxAtStartPar
\sphinxcode{\sphinxupquote{Var GetBinVar()}}
\end{quote}

\sphinxAtStartPar
\sphinxstylestrong{Return}
\begin{quote}

\sphinxAtStartPar
binary vaiable object.
\end{quote}
\end{quote}

\subsubsection{GenConstrBuilder.GetExpr()}
\label{\detokenize{csapi/genconstrbuilder:genconstrbuilder-getexpr}}\begin{quote}

\sphinxAtStartPar
Get expression associated with general constraint.

\sphinxAtStartPar
\sphinxstylestrong{Synopsis}
\begin{quote}

\sphinxAtStartPar
\sphinxcode{\sphinxupquote{Expr GetExpr()}}
\end{quote}

\sphinxAtStartPar
\sphinxstylestrong{Return}
\begin{quote}

\sphinxAtStartPar
expression object.
\end{quote}
\end{quote}

\subsubsection{GenConstrBuilder.GetIndType()}
\label{\detokenize{csapi/genconstrbuilder:genconstrbuilder-getindtype}}\begin{quote}

\sphinxAtStartPar
Get type of general constraint.

\sphinxAtStartPar
\sphinxstylestrong{Synopsis}
\begin{quote}

\sphinxAtStartPar
\sphinxcode{\sphinxupquote{int GetIndType()}}
\end{quote}

\sphinxAtStartPar
\sphinxstylestrong{Return}
\begin{quote}

\sphinxAtStartPar
type of general constraint.
\end{quote}
\end{quote}

\subsubsection{GenConstrBuilder.GetSense()}
\label{\detokenize{csapi/genconstrbuilder:genconstrbuilder-getsense}}\begin{quote}

\sphinxAtStartPar
Get sense associated with general constraint.

\sphinxAtStartPar
\sphinxstylestrong{Synopsis}
\begin{quote}

\sphinxAtStartPar
\sphinxcode{\sphinxupquote{char GetSense()}}
\end{quote}

\sphinxAtStartPar
\sphinxstylestrong{Return}
\begin{quote}

\sphinxAtStartPar
constraint sense.
\end{quote}
\end{quote}

\subsubsection{GenConstrBuilder.Set()}
\label{\detokenize{csapi/genconstrbuilder:genconstrbuilder-set}}\begin{quote}

\sphinxAtStartPar
Set binary variable, binary value, expression and sense of general constraint.

\sphinxAtStartPar
\sphinxstylestrong{Synopsis}
\begin{quote}

\sphinxAtStartPar
\sphinxcode{\sphinxupquote{void Set(}}
\begin{quote}

\sphinxAtStartPar
\sphinxcode{\sphinxupquote{Var binvar,}}

\sphinxAtStartPar
\sphinxcode{\sphinxupquote{int binval,}}

\sphinxAtStartPar
\sphinxcode{\sphinxupquote{Expr expr,}}

\sphinxAtStartPar
\sphinxcode{\sphinxupquote{char sense,}}

\sphinxAtStartPar
\sphinxcode{\sphinxupquote{int type)}}
\end{quote}
\end{quote}

\sphinxAtStartPar
\sphinxstylestrong{Arguments}
\begin{quote}

\sphinxAtStartPar
\sphinxcode{\sphinxupquote{binvar}}: binary variable.

\sphinxAtStartPar
\sphinxcode{\sphinxupquote{binval}}: binary value.

\sphinxAtStartPar
\sphinxcode{\sphinxupquote{expr}}: expression object.

\sphinxAtStartPar
\sphinxcode{\sphinxupquote{sense}}: general constraint sense.

\sphinxAtStartPar
\sphinxcode{\sphinxupquote{type}}: type of general constraint.
\end{quote}
\end{quote}

\subsection{GenConstrBuilderArray}
\label{\detokenize{csharpapiref:genconstrbuilderarray}}\label{\detokenize{csharpapiref:chapcsharpapiref-genconstrbuilderarray}}
\sphinxAtStartPar
COPT general constraint builder array object. To store and access a set of C\#
{\hyperref[\detokenize{csharpapiref:chapcsharpapiref-genconstrbuilder}]{\sphinxcrossref{\DUrole{std,std-ref}{GenConstrBuilder}}}} objects, Cardinal Optimizer provides C\#
GenConstrBuilderArray class, which defines the following methods.

\sphinxstepscope

\subsubsection{GenConstrBuilderArray.GenConstrBuilderArray()}
\label{\detokenize{csapi/genconstrbuilderarray:genconstrbuilderarray-genconstrbuilderarray}}\label{\detokenize{csapi/genconstrbuilderarray::doc}}\begin{quote}

\sphinxAtStartPar
Constructor of genconstrbuilderarray.

\sphinxAtStartPar
\sphinxstylestrong{Synopsis}
\begin{quote}

\sphinxAtStartPar
\sphinxcode{\sphinxupquote{GenConstrBuilderArray()}}
\end{quote}
\end{quote}

\subsubsection{GenConstrBuilderArray.GetBuilder()}
\label{\detokenize{csapi/genconstrbuilderarray:genconstrbuilderarray-getbuilder}}\begin{quote}

\sphinxAtStartPar
Get idx\sphinxhyphen{}th general constraint builder object.

\sphinxAtStartPar
\sphinxstylestrong{Synopsis}
\begin{quote}

\sphinxAtStartPar
\sphinxcode{\sphinxupquote{GenConstrBuilder GetBuilder(int idx)}}
\end{quote}

\sphinxAtStartPar
\sphinxstylestrong{Arguments}
\begin{quote}

\sphinxAtStartPar
\sphinxcode{\sphinxupquote{idx}}: index of the general constraint builder.
\end{quote}

\sphinxAtStartPar
\sphinxstylestrong{Return}
\begin{quote}

\sphinxAtStartPar
general constraint builder object with index idx.
\end{quote}
\end{quote}

\subsubsection{GenConstrBuilderArray.PushBack()}
\label{\detokenize{csapi/genconstrbuilderarray:genconstrbuilderarray-pushback}}\begin{quote}

\sphinxAtStartPar
Add a general constraint builder object to general constraint builder array.

\sphinxAtStartPar
\sphinxstylestrong{Synopsis}
\begin{quote}

\sphinxAtStartPar
\sphinxcode{\sphinxupquote{void PushBack(GenConstrBuilder builder)}}
\end{quote}

\sphinxAtStartPar
\sphinxstylestrong{Arguments}
\begin{quote}

\sphinxAtStartPar
\sphinxcode{\sphinxupquote{builder}}: a general constraint builder object.
\end{quote}
\end{quote}

\subsubsection{GenConstrBuilderArray.Size()}
\label{\detokenize{csapi/genconstrbuilderarray:genconstrbuilderarray-size}}\begin{quote}

\sphinxAtStartPar
Get the number of general constraint builder objects.

\sphinxAtStartPar
\sphinxstylestrong{Synopsis}
\begin{quote}

\sphinxAtStartPar
\sphinxcode{\sphinxupquote{int Size()}}
\end{quote}

\sphinxAtStartPar
\sphinxstylestrong{Return}
\begin{quote}

\sphinxAtStartPar
number of general constraint builder objects.
\end{quote}
\end{quote}

\subsection{Cone}
\label{\detokenize{csharpapiref:cone}}\label{\detokenize{csharpapiref:chapcsharpapiref-cone}}
\sphinxAtStartPar
COPT cone constraint object. Cone constraints are always associated with a particular model.
User creates a cone constraint object by adding a cone constraint to a model,
rather than by using constructor of Cone class.

\sphinxAtStartPar
A cone constraint can be regular or rotated (\sphinxcode{\sphinxupquote{COPT\_CONE\_QUAD}} or \sphinxcode{\sphinxupquote{COPT\_CONE\_RQUAD}}).

\sphinxstepscope

\subsubsection{Cone.GetIdx()}
\label{\detokenize{csapi/cone:cone-getidx}}\label{\detokenize{csapi/cone::doc}}\begin{quote}

\sphinxAtStartPar
Get the index of a cone constraint.

\sphinxAtStartPar
\sphinxstylestrong{Synopsis}
\begin{quote}

\sphinxAtStartPar
\sphinxcode{\sphinxupquote{int GetIdx()}}
\end{quote}

\sphinxAtStartPar
\sphinxstylestrong{Return}
\begin{quote}

\sphinxAtStartPar
index of the cone constraint.
\end{quote}
\end{quote}

\subsubsection{Cone.Remove()}
\label{\detokenize{csapi/cone:cone-remove}}\begin{quote}

\sphinxAtStartPar
Remove the cone constraint from model.

\sphinxAtStartPar
\sphinxstylestrong{Synopsis}
\begin{quote}

\sphinxAtStartPar
\sphinxcode{\sphinxupquote{void Remove()}}
\end{quote}
\end{quote}

\subsection{ConeArray}
\label{\detokenize{csharpapiref:conearray}}\label{\detokenize{csharpapiref:chapcsharpapiref-conearray}}
\sphinxAtStartPar
COPT cone constraint array object. To store and access a set of C\# {\hyperref[\detokenize{csharpapiref:chapcsharpapiref-cone}]{\sphinxcrossref{\DUrole{std,std-ref}{Cone}}}}
objects, Cardinal Optimizer provides C\# ConeArray class, which defines the following methods.

\sphinxstepscope

\subsubsection{ConeArray.ConeArray()}
\label{\detokenize{csapi/conearray:conearray-conearray}}\label{\detokenize{csapi/conearray::doc}}\begin{quote}

\sphinxAtStartPar
Constructor of ConeArray object.

\sphinxAtStartPar
\sphinxstylestrong{Synopsis}
\begin{quote}

\sphinxAtStartPar
\sphinxcode{\sphinxupquote{ConeArray()}}
\end{quote}
\end{quote}

\subsubsection{ConeArray.GetCone()}
\label{\detokenize{csapi/conearray:conearray-getcone}}\begin{quote}

\sphinxAtStartPar
Get idx\sphinxhyphen{}th cone object.

\sphinxAtStartPar
\sphinxstylestrong{Synopsis}
\begin{quote}

\sphinxAtStartPar
\sphinxcode{\sphinxupquote{Cone GetCone(int idx)}}
\end{quote}

\sphinxAtStartPar
\sphinxstylestrong{Arguments}
\begin{quote}

\sphinxAtStartPar
\sphinxcode{\sphinxupquote{idx}}: index of cone.
\end{quote}

\sphinxAtStartPar
\sphinxstylestrong{Return}
\begin{quote}

\sphinxAtStartPar
cone object with index idx.
\end{quote}
\end{quote}

\subsubsection{ConeArray.PushBack()}
\label{\detokenize{csapi/conearray:conearray-pushback}}\begin{quote}

\sphinxAtStartPar
Add a cone constraint object to cone constraint array.

\sphinxAtStartPar
\sphinxstylestrong{Synopsis}
\begin{quote}

\sphinxAtStartPar
\sphinxcode{\sphinxupquote{void PushBack(Cone cone)}}
\end{quote}

\sphinxAtStartPar
\sphinxstylestrong{Arguments}
\begin{quote}

\sphinxAtStartPar
\sphinxcode{\sphinxupquote{cone}}: a cone constraint object.
\end{quote}
\end{quote}

\subsubsection{ConeArray.Size()}
\label{\detokenize{csapi/conearray:conearray-size}}\begin{quote}

\sphinxAtStartPar
Get the number of cone constraint objects.

\sphinxAtStartPar
\sphinxstylestrong{Synopsis}
\begin{quote}

\sphinxAtStartPar
\sphinxcode{\sphinxupquote{int Size()}}
\end{quote}

\sphinxAtStartPar
\sphinxstylestrong{Return}
\begin{quote}

\sphinxAtStartPar
number of cone constraint objects.
\end{quote}
\end{quote}

\subsection{ConeBuilder}
\label{\detokenize{csharpapiref:conebuilder}}\label{\detokenize{csharpapiref:chapcsharpapiref-conebuilder}}
\sphinxAtStartPar
COPT cone constraint builder object. To help building a cone constraint, given the cone type and
a set of variables, Cardinal Optimizer provides C\# ConeBuilder class, which defines the following methods.

\sphinxstepscope

\subsubsection{ConeBuilder.ConeBuilder()}
\label{\detokenize{csapi/conebuilder:conebuilder-conebuilder}}\label{\detokenize{csapi/conebuilder::doc}}\begin{quote}

\sphinxAtStartPar
Constructor of ConeBuilder object.

\sphinxAtStartPar
\sphinxstylestrong{Synopsis}
\begin{quote}

\sphinxAtStartPar
\sphinxcode{\sphinxupquote{ConeBuilder()}}
\end{quote}
\end{quote}

\subsubsection{ConeBuilder.GetSize()}
\label{\detokenize{csapi/conebuilder:conebuilder-getsize}}\begin{quote}

\sphinxAtStartPar
Get number of variables in a cone constraint.

\sphinxAtStartPar
\sphinxstylestrong{Synopsis}
\begin{quote}

\sphinxAtStartPar
\sphinxcode{\sphinxupquote{int GetSize()}}
\end{quote}

\sphinxAtStartPar
\sphinxstylestrong{Return}
\begin{quote}

\sphinxAtStartPar
number of variables.
\end{quote}
\end{quote}

\subsubsection{ConeBuilder.GetType()}
\label{\detokenize{csapi/conebuilder:conebuilder-gettype}}\begin{quote}

\sphinxAtStartPar
Get type of a cone constraint.

\sphinxAtStartPar
\sphinxstylestrong{Synopsis}
\begin{quote}

\sphinxAtStartPar
\sphinxcode{\sphinxupquote{int GetType()}}
\end{quote}

\sphinxAtStartPar
\sphinxstylestrong{Return}
\begin{quote}

\sphinxAtStartPar
type of the cone constraint.
\end{quote}
\end{quote}

\subsubsection{ConeBuilder.GetVar()}
\label{\detokenize{csapi/conebuilder:conebuilder-getvar}}\begin{quote}

\sphinxAtStartPar
Get i\sphinxhyphen{}th variable in a cone constraint.

\sphinxAtStartPar
\sphinxstylestrong{Synopsis}
\begin{quote}

\sphinxAtStartPar
\sphinxcode{\sphinxupquote{Var GetVar(int idx)}}
\end{quote}

\sphinxAtStartPar
\sphinxstylestrong{Arguments}
\begin{quote}

\sphinxAtStartPar
\sphinxcode{\sphinxupquote{idx}}: index of vars.
\end{quote}

\sphinxAtStartPar
\sphinxstylestrong{Return}
\begin{quote}

\sphinxAtStartPar
the i\sphinxhyphen{}th variable in a cone constraint.
\end{quote}
\end{quote}

\subsubsection{ConeBuilder.GetVars()}
\label{\detokenize{csapi/conebuilder:conebuilder-getvars}}\begin{quote}

\sphinxAtStartPar
Get all variables in a cone constraint.

\sphinxAtStartPar
\sphinxstylestrong{Synopsis}
\begin{quote}

\sphinxAtStartPar
\sphinxcode{\sphinxupquote{VarArray GetVars()}}
\end{quote}

\sphinxAtStartPar
\sphinxstylestrong{Return}
\begin{quote}

\sphinxAtStartPar
variables in a cone constraint.
\end{quote}
\end{quote}

\subsubsection{ConeBuilder.Set()}
\label{\detokenize{csapi/conebuilder:conebuilder-set}}\begin{quote}

\sphinxAtStartPar
Set variables and type of a cone constraint.

\sphinxAtStartPar
\sphinxstylestrong{Synopsis}
\begin{quote}

\sphinxAtStartPar
\sphinxcode{\sphinxupquote{void Set(VarArray vars, int type)}}
\end{quote}

\sphinxAtStartPar
\sphinxstylestrong{Arguments}
\begin{quote}

\sphinxAtStartPar
\sphinxcode{\sphinxupquote{vars}}: variable array object.

\sphinxAtStartPar
\sphinxcode{\sphinxupquote{type}}: type of a cone constraint.
\end{quote}
\end{quote}

\subsection{ConeBuilderArray}
\label{\detokenize{csharpapiref:conebuilderarray}}\label{\detokenize{csharpapiref:chapcsharpapiref-conebuilderarray}}
\sphinxAtStartPar
COPT cone constraint builder array object. To store and access a set of C\# {\hyperref[\detokenize{csharpapiref:chapcsharpapiref-conebuilder}]{\sphinxcrossref{\DUrole{std,std-ref}{ConeBuilder}}}}
objects, Cardinal Optimizer provides C\# ConeBuilderArray class, which defines the following methods.

\sphinxstepscope

\subsubsection{ConeBuilderArray.ConeBuilderArray()}
\label{\detokenize{csapi/conebuilderarray:conebuilderarray-conebuilderarray}}\label{\detokenize{csapi/conebuilderarray::doc}}\begin{quote}

\sphinxAtStartPar
Constructor of ConeBuilderArray object.

\sphinxAtStartPar
\sphinxstylestrong{Synopsis}
\begin{quote}

\sphinxAtStartPar
\sphinxcode{\sphinxupquote{ConeBuilderArray()}}
\end{quote}
\end{quote}

\subsubsection{ConeBuilderArray.GetBuilder()}
\label{\detokenize{csapi/conebuilderarray:conebuilderarray-getbuilder}}\begin{quote}

\sphinxAtStartPar
Get idx\sphinxhyphen{}th cone constraint builder object.

\sphinxAtStartPar
\sphinxstylestrong{Synopsis}
\begin{quote}

\sphinxAtStartPar
\sphinxcode{\sphinxupquote{ConeBuilder GetBuilder(int idx)}}
\end{quote}

\sphinxAtStartPar
\sphinxstylestrong{Arguments}
\begin{quote}

\sphinxAtStartPar
\sphinxcode{\sphinxupquote{idx}}: index of the cone constraint builder.
\end{quote}

\sphinxAtStartPar
\sphinxstylestrong{Return}
\begin{quote}

\sphinxAtStartPar
cone constraint builder object with index idx.
\end{quote}
\end{quote}

\subsubsection{ConeBuilderArray.PushBack()}
\label{\detokenize{csapi/conebuilderarray:conebuilderarray-pushback}}\begin{quote}

\sphinxAtStartPar
Add a cone constraint builder object to cone constraint builder array.

\sphinxAtStartPar
\sphinxstylestrong{Synopsis}
\begin{quote}

\sphinxAtStartPar
\sphinxcode{\sphinxupquote{void PushBack(ConeBuilder builder)}}
\end{quote}

\sphinxAtStartPar
\sphinxstylestrong{Arguments}
\begin{quote}

\sphinxAtStartPar
\sphinxcode{\sphinxupquote{builder}}: a cone constraint builder object.
\end{quote}
\end{quote}

\subsubsection{ConeBuilderArray.Size()}
\label{\detokenize{csapi/conebuilderarray:conebuilderarray-size}}\begin{quote}

\sphinxAtStartPar
Get the number of cone constraint builder objects.

\sphinxAtStartPar
\sphinxstylestrong{Synopsis}
\begin{quote}

\sphinxAtStartPar
\sphinxcode{\sphinxupquote{int Size()}}
\end{quote}

\sphinxAtStartPar
\sphinxstylestrong{Return}
\begin{quote}

\sphinxAtStartPar
number of cone constraint builder objects.
\end{quote}
\end{quote}

\subsection{ExpCone}
\label{\detokenize{csharpapiref:expcone}}\label{\detokenize{csharpapiref:chapcsharpapiref-expcone}}
\sphinxAtStartPar
COPT exponential cone constraint object. Exponential Cone constraints are always associated with a particular model.
User creates an exponential cone constraint object by adding an exponential cone constraint to a model,
rather than by using constructor of ExpCone class.

\sphinxstepscope

\subsubsection{ExpCone.GetIdx()}
\label{\detokenize{csapi/expcone:expcone-getidx}}\label{\detokenize{csapi/expcone::doc}}\begin{quote}

\sphinxAtStartPar
Get the index of an exponential cone constraint.

\sphinxAtStartPar
\sphinxstylestrong{Synopsis}
\begin{quote}

\sphinxAtStartPar
\sphinxcode{\sphinxupquote{int GetIdx()}}
\end{quote}

\sphinxAtStartPar
\sphinxstylestrong{Return}
\begin{quote}

\sphinxAtStartPar
index of the exponential cone constraint.
\end{quote}
\end{quote}

\subsubsection{ExpCone.Remove()}
\label{\detokenize{csapi/expcone:expcone-remove}}\begin{quote}

\sphinxAtStartPar
Remove the exponential cone constraint from model.

\sphinxAtStartPar
\sphinxstylestrong{Synopsis}
\begin{quote}

\sphinxAtStartPar
\sphinxcode{\sphinxupquote{void Remove()}}
\end{quote}
\end{quote}

\subsection{ExpConeArray}
\label{\detokenize{csharpapiref:expconearray}}\label{\detokenize{csharpapiref:chapcsharpapiref-expconearray}}
\sphinxAtStartPar
COPT exponential cone constraint array object. To store and access a set of C\# {\hyperref[\detokenize{csharpapiref:chapcsharpapiref-expcone}]{\sphinxcrossref{\DUrole{std,std-ref}{ExpCone}}}}
objects, Cardinal Optimizer provides C\# ExpConeArray class, which defines the following methods.

\sphinxstepscope

\subsubsection{ExpConeArray.ExpConeArray()}
\label{\detokenize{csapi/expconearray:expconearray-expconearray}}\label{\detokenize{csapi/expconearray::doc}}\begin{quote}

\sphinxAtStartPar
Constructor of ExpConeArray object.

\sphinxAtStartPar
\sphinxstylestrong{Synopsis}
\begin{quote}

\sphinxAtStartPar
\sphinxcode{\sphinxupquote{ExpConeArray()}}
\end{quote}
\end{quote}

\subsubsection{ExpConeArray.GetCone()}
\label{\detokenize{csapi/expconearray:expconearray-getcone}}\begin{quote}

\sphinxAtStartPar
Get idx\sphinxhyphen{}th exponential cone object.

\sphinxAtStartPar
\sphinxstylestrong{Synopsis}
\begin{quote}

\sphinxAtStartPar
\sphinxcode{\sphinxupquote{ExpCone GetCone(int idx)}}
\end{quote}

\sphinxAtStartPar
\sphinxstylestrong{Arguments}
\begin{quote}

\sphinxAtStartPar
\sphinxcode{\sphinxupquote{idx}}: index of exponential cone.
\end{quote}

\sphinxAtStartPar
\sphinxstylestrong{Return}
\begin{quote}

\sphinxAtStartPar
exponential cone object with index idx.
\end{quote}
\end{quote}

\subsubsection{ExpConeArray.PushBack()}
\label{\detokenize{csapi/expconearray:expconearray-pushback}}\begin{quote}

\sphinxAtStartPar
Add an exponential cone constraint object to exponential cone constraint array.

\sphinxAtStartPar
\sphinxstylestrong{Synopsis}
\begin{quote}

\sphinxAtStartPar
\sphinxcode{\sphinxupquote{void PushBack(ExpCone cone)}}
\end{quote}

\sphinxAtStartPar
\sphinxstylestrong{Arguments}
\begin{quote}

\sphinxAtStartPar
\sphinxcode{\sphinxupquote{cone}}: an exponential cone constraint object.
\end{quote}
\end{quote}

\subsubsection{ExpConeArray.Size()}
\label{\detokenize{csapi/expconearray:expconearray-size}}\begin{quote}

\sphinxAtStartPar
Get the number of exponential cone constraint objects.

\sphinxAtStartPar
\sphinxstylestrong{Synopsis}
\begin{quote}

\sphinxAtStartPar
\sphinxcode{\sphinxupquote{int Size()}}
\end{quote}

\sphinxAtStartPar
\sphinxstylestrong{Return}
\begin{quote}

\sphinxAtStartPar
number of exponential cone constraint objects.
\end{quote}
\end{quote}

\subsection{ExpConeBuilder}
\label{\detokenize{csharpapiref:expconebuilder}}\label{\detokenize{csharpapiref:chapcsharpapiref-expconebuilder}}
\sphinxAtStartPar
COPT exponential cone constraint builder object. To help building an exponential cone constraint, given the exponential
cone type and a set of variables, Cardinal Optimizer provides C\# ExpConeBuilder class, which defines the following methods.

\sphinxstepscope

\subsubsection{ExpConeBuilder.ExpConeBuilder()}
\label{\detokenize{csapi/expconebuilder:expconebuilder-expconebuilder}}\label{\detokenize{csapi/expconebuilder::doc}}\begin{quote}

\sphinxAtStartPar
Constructor of ExpConeBuilder object.

\sphinxAtStartPar
\sphinxstylestrong{Synopsis}
\begin{quote}

\sphinxAtStartPar
\sphinxcode{\sphinxupquote{ExpConeBuilder()}}
\end{quote}
\end{quote}

\subsubsection{ExpConeBuilder.GetSize()}
\label{\detokenize{csapi/expconebuilder:expconebuilder-getsize}}\begin{quote}

\sphinxAtStartPar
Get number of variables in an exponential cone constraint.

\sphinxAtStartPar
\sphinxstylestrong{Synopsis}
\begin{quote}

\sphinxAtStartPar
\sphinxcode{\sphinxupquote{int GetSize()}}
\end{quote}

\sphinxAtStartPar
\sphinxstylestrong{Return}
\begin{quote}

\sphinxAtStartPar
number of variables.
\end{quote}
\end{quote}

\subsubsection{ExpConeBuilder.GetType()}
\label{\detokenize{csapi/expconebuilder:expconebuilder-gettype}}\begin{quote}

\sphinxAtStartPar
Get type of an exponential cone constraint.

\sphinxAtStartPar
\sphinxstylestrong{Synopsis}
\begin{quote}

\sphinxAtStartPar
\sphinxcode{\sphinxupquote{int GetType()}}
\end{quote}

\sphinxAtStartPar
\sphinxstylestrong{Return}
\begin{quote}

\sphinxAtStartPar
type of the exponential cone constraint.
\end{quote}
\end{quote}

\subsubsection{ExpConeBuilder.GetVar()}
\label{\detokenize{csapi/expconebuilder:expconebuilder-getvar}}\begin{quote}

\sphinxAtStartPar
Get i\sphinxhyphen{}th variable in an exponential cone constraint.

\sphinxAtStartPar
\sphinxstylestrong{Synopsis}
\begin{quote}

\sphinxAtStartPar
\sphinxcode{\sphinxupquote{Var GetVar(int idx)}}
\end{quote}

\sphinxAtStartPar
\sphinxstylestrong{Arguments}
\begin{quote}

\sphinxAtStartPar
\sphinxcode{\sphinxupquote{idx}}: index of vars.
\end{quote}

\sphinxAtStartPar
\sphinxstylestrong{Return}
\begin{quote}

\sphinxAtStartPar
the i\sphinxhyphen{}th variable in an exponential cone constraint.
\end{quote}
\end{quote}

\subsubsection{ExpConeBuilder.GetVars()}
\label{\detokenize{csapi/expconebuilder:expconebuilder-getvars}}\begin{quote}

\sphinxAtStartPar
Get all variables in an exponential cone constraint.

\sphinxAtStartPar
\sphinxstylestrong{Synopsis}
\begin{quote}

\sphinxAtStartPar
\sphinxcode{\sphinxupquote{VarArray GetVars()}}
\end{quote}

\sphinxAtStartPar
\sphinxstylestrong{Return}
\begin{quote}

\sphinxAtStartPar
variables in an exponential cone constraint.
\end{quote}
\end{quote}

\subsubsection{ExpConeBuilder.Set()}
\label{\detokenize{csapi/expconebuilder:expconebuilder-set}}\begin{quote}

\sphinxAtStartPar
Set variables and type of an exponential cone constraint.

\sphinxAtStartPar
\sphinxstylestrong{Synopsis}
\begin{quote}

\sphinxAtStartPar
\sphinxcode{\sphinxupquote{void Set(VarArray vars, int type)}}
\end{quote}

\sphinxAtStartPar
\sphinxstylestrong{Arguments}
\begin{quote}

\sphinxAtStartPar
\sphinxcode{\sphinxupquote{vars}}: variable array object.

\sphinxAtStartPar
\sphinxcode{\sphinxupquote{type}}: type of an exponential cone constraint.
\end{quote}
\end{quote}

\subsection{ExpConeBuilderArray}
\label{\detokenize{csharpapiref:expconebuilderarray}}\label{\detokenize{csharpapiref:chapcsharpapiref-expconebuilderarray}}
\sphinxAtStartPar
COPT exponential cone constraint builder array object. To store and access a set of C\# {\hyperref[\detokenize{csharpapiref:chapcsharpapiref-expconebuilder}]{\sphinxcrossref{\DUrole{std,std-ref}{ExpConeBuilder}}}}
objects, Cardinal Optimizer provides C\# ExpConeBuilderArray class, which defines the following methods.

\sphinxstepscope

\subsubsection{ExpConeBuilderArray.ExpConeBuilderArray()}
\label{\detokenize{csapi/expconebuilderarray:expconebuilderarray-expconebuilderarray}}\label{\detokenize{csapi/expconebuilderarray::doc}}\begin{quote}

\sphinxAtStartPar
Constructor of ExpConeBuilderArray object.

\sphinxAtStartPar
\sphinxstylestrong{Synopsis}
\begin{quote}

\sphinxAtStartPar
\sphinxcode{\sphinxupquote{ExpConeBuilderArray()}}
\end{quote}
\end{quote}

\subsubsection{ExpConeBuilderArray.GetBuilder()}
\label{\detokenize{csapi/expconebuilderarray:expconebuilderarray-getbuilder}}\begin{quote}

\sphinxAtStartPar
Get idx\sphinxhyphen{}th exponential cone constraint builder object.

\sphinxAtStartPar
\sphinxstylestrong{Synopsis}
\begin{quote}

\sphinxAtStartPar
\sphinxcode{\sphinxupquote{ExpConeBuilder GetBuilder(int idx)}}
\end{quote}

\sphinxAtStartPar
\sphinxstylestrong{Arguments}
\begin{quote}

\sphinxAtStartPar
\sphinxcode{\sphinxupquote{idx}}: index of the exponential cone constraint builder.
\end{quote}

\sphinxAtStartPar
\sphinxstylestrong{Return}
\begin{quote}

\sphinxAtStartPar
exponential cone constraint builder object with index idx.
\end{quote}
\end{quote}

\subsubsection{ExpConeBuilderArray.PushBack()}
\label{\detokenize{csapi/expconebuilderarray:expconebuilderarray-pushback}}\begin{quote}

\sphinxAtStartPar
Add an exponential cone constraint builder object to exponential cone constraint builder array.

\sphinxAtStartPar
\sphinxstylestrong{Synopsis}
\begin{quote}

\sphinxAtStartPar
\sphinxcode{\sphinxupquote{void PushBack(ExpConeBuilder builder)}}
\end{quote}

\sphinxAtStartPar
\sphinxstylestrong{Arguments}
\begin{quote}

\sphinxAtStartPar
\sphinxcode{\sphinxupquote{builder}}: an exponential cone constraint builder object.
\end{quote}
\end{quote}

\subsubsection{ExpConeBuilderArray.Size()}
\label{\detokenize{csapi/expconebuilderarray:expconebuilderarray-size}}\begin{quote}

\sphinxAtStartPar
Get the number of exponential cone constraint builder objects.

\sphinxAtStartPar
\sphinxstylestrong{Synopsis}
\begin{quote}

\sphinxAtStartPar
\sphinxcode{\sphinxupquote{int Size()}}
\end{quote}

\sphinxAtStartPar
\sphinxstylestrong{Return}
\begin{quote}

\sphinxAtStartPar
number of exponential cone constraint builder objects.
\end{quote}
\end{quote}

\subsection{AffineCone}
\label{\detokenize{csharpapiref:affinecone}}\label{\detokenize{csharpapiref:chapcsharpapiref-affinecone}}
\sphinxAtStartPar
The \sphinxtitleref{AffineCone} class in COPT encapsulates operations related to affine cones.
The following methods are provided:

\sphinxstepscope

\subsubsection{AffineCone.GetIdx()}
\label{\detokenize{csapi/affinecone:affinecone-getidx}}\label{\detokenize{csapi/affinecone::doc}}\begin{quote}

\sphinxAtStartPar
Get the index of an affine cone constraint.

\sphinxAtStartPar
\sphinxstylestrong{Synopsis}
\begin{quote}

\sphinxAtStartPar
\sphinxcode{\sphinxupquote{int GetIdx()}}
\end{quote}

\sphinxAtStartPar
\sphinxstylestrong{Return}
\begin{quote}

\sphinxAtStartPar
index of the affine cone constraint.
\end{quote}
\end{quote}

\subsubsection{AffineCone.GetName()}
\label{\detokenize{csapi/affinecone:affinecone-getname}}\begin{quote}

\sphinxAtStartPar
Get name of the affine cone.

\sphinxAtStartPar
\sphinxstylestrong{Synopsis}
\begin{quote}

\sphinxAtStartPar
\sphinxcode{\sphinxupquote{string GetName()}}
\end{quote}

\sphinxAtStartPar
\sphinxstylestrong{Return}
\begin{quote}

\sphinxAtStartPar
affine cone name.
\end{quote}
\end{quote}

\subsubsection{AffineCone.Remove()}
\label{\detokenize{csapi/affinecone:affinecone-remove}}\begin{quote}

\sphinxAtStartPar
Remove the affine cone constraint from model.

\sphinxAtStartPar
\sphinxstylestrong{Synopsis}
\begin{quote}

\sphinxAtStartPar
\sphinxcode{\sphinxupquote{void Remove()}}
\end{quote}
\end{quote}

\subsubsection{AffineCone.SetName()}
\label{\detokenize{csapi/affinecone:affinecone-setname}}\begin{quote}

\sphinxAtStartPar
Set name of the affine cone.

\sphinxAtStartPar
\sphinxstylestrong{Synopsis}
\begin{quote}

\sphinxAtStartPar
\sphinxcode{\sphinxupquote{void SetName(string name)}}
\end{quote}

\sphinxAtStartPar
\sphinxstylestrong{Arguments}
\begin{quote}

\sphinxAtStartPar
\sphinxcode{\sphinxupquote{name}}: affinecone name.
\end{quote}
\end{quote}

\subsection{AffineConeArray}
\label{\detokenize{csharpapiref:affineconearray}}\label{\detokenize{csharpapiref:chapcsharpapiref-affineconearray}}
\sphinxAtStartPar
To facilitate operations on a group of C\# {\hyperref[\detokenize{csharpapiref:chapcsharpapiref-affinecone}]{\sphinxcrossref{\DUrole{std,std-ref}{AffineCone}}}} objects, the C\# interface of COPT introduces the \sphinxtitleref{AffineConeArray} class. The following methods are provided:

\sphinxstepscope

\subsubsection{AffineConeArray.AffineConeArray()}
\label{\detokenize{csapi/affineconearray:affineconearray-affineconearray}}\label{\detokenize{csapi/affineconearray::doc}}\begin{quote}

\sphinxAtStartPar
Constructor of AffineConeArray object.

\sphinxAtStartPar
\sphinxstylestrong{Synopsis}
\begin{quote}

\sphinxAtStartPar
\sphinxcode{\sphinxupquote{AffineConeArray()}}
\end{quote}
\end{quote}

\subsubsection{AffineConeArray.GetCone()}
\label{\detokenize{csapi/affineconearray:affineconearray-getcone}}\begin{quote}

\sphinxAtStartPar
Get idx\sphinxhyphen{}th affine cone object.

\sphinxAtStartPar
\sphinxstylestrong{Synopsis}
\begin{quote}

\sphinxAtStartPar
\sphinxcode{\sphinxupquote{AffineCone GetCone(int idx)}}
\end{quote}

\sphinxAtStartPar
\sphinxstylestrong{Arguments}
\begin{quote}

\sphinxAtStartPar
\sphinxcode{\sphinxupquote{idx}}: index of affine cone.
\end{quote}

\sphinxAtStartPar
\sphinxstylestrong{Return}
\begin{quote}

\sphinxAtStartPar
affine cone object with index idx.
\end{quote}
\end{quote}

\subsubsection{AffineConeArray.PushBack()}
\label{\detokenize{csapi/affineconearray:affineconearray-pushback}}\begin{quote}

\sphinxAtStartPar
Add an affine cone constraint object to affine cone constraint array.

\sphinxAtStartPar
\sphinxstylestrong{Synopsis}
\begin{quote}

\sphinxAtStartPar
\sphinxcode{\sphinxupquote{void PushBack(AffineCone cone)}}
\end{quote}

\sphinxAtStartPar
\sphinxstylestrong{Arguments}
\begin{quote}

\sphinxAtStartPar
\sphinxcode{\sphinxupquote{cone}}: an affine cone constraint object.
\end{quote}
\end{quote}

\subsubsection{AffineConeArray.Size()}
\label{\detokenize{csapi/affineconearray:affineconearray-size}}\begin{quote}

\sphinxAtStartPar
Get the number of affine cone constraint objects.

\sphinxAtStartPar
\sphinxstylestrong{Synopsis}
\begin{quote}

\sphinxAtStartPar
\sphinxcode{\sphinxupquote{int Size()}}
\end{quote}

\sphinxAtStartPar
\sphinxstylestrong{Return}
\begin{quote}

\sphinxAtStartPar
number of affine cone constraint objects.
\end{quote}
\end{quote}

\subsection{AffineConeBuilder}
\label{\detokenize{csharpapiref:affineconebuilder}}\label{\detokenize{csharpapiref:chapcsharpapiref-affineconebuilder}}
\sphinxAtStartPar
The \sphinxtitleref{AffineConeBuilder} class in COPT encapsulates the builder for constructing affine cones.
The following methods are provided:

\sphinxstepscope

\subsubsection{AffineConeBuilder.AffineConeBuilder()}
\label{\detokenize{csapi/affineconebuilder:affineconebuilder-affineconebuilder}}\label{\detokenize{csapi/affineconebuilder::doc}}\begin{quote}

\sphinxAtStartPar
Constructor of AffineConeBuilder object.

\sphinxAtStartPar
\sphinxstylestrong{Synopsis}
\begin{quote}

\sphinxAtStartPar
\sphinxcode{\sphinxupquote{AffineConeBuilder()}}
\end{quote}
\end{quote}

\subsubsection{AffineConeBuilder.GetExpr()}
\label{\detokenize{csapi/affineconebuilder:affineconebuilder-getexpr}}\begin{quote}

\sphinxAtStartPar
Get i\sphinxhyphen{}th linear expression in an affine cone constraint.

\sphinxAtStartPar
\sphinxstylestrong{Synopsis}
\begin{quote}

\sphinxAtStartPar
\sphinxcode{\sphinxupquote{Expr GetExpr(int idx)}}
\end{quote}

\sphinxAtStartPar
\sphinxstylestrong{Arguments}
\begin{quote}

\sphinxAtStartPar
\sphinxcode{\sphinxupquote{idx}}: index of linear expression.
\end{quote}

\sphinxAtStartPar
\sphinxstylestrong{Return}
\begin{quote}

\sphinxAtStartPar
the i\sphinxhyphen{}th linear expression in an affine cone constraint.
\end{quote}
\end{quote}

\subsubsection{AffineConeBuilder.GetExprs()}
\label{\detokenize{csapi/affineconebuilder:affineconebuilder-getexprs}}\begin{quote}

\sphinxAtStartPar
Get all linear expressions in an affine cone constraint.

\sphinxAtStartPar
\sphinxstylestrong{Synopsis}
\begin{quote}

\sphinxAtStartPar
\sphinxcode{\sphinxupquote{Expr{[}{]} GetExprs()}}
\end{quote}

\sphinxAtStartPar
\sphinxstylestrong{Return}
\begin{quote}

\sphinxAtStartPar
array of linear expressions.
\end{quote}
\end{quote}

\subsubsection{AffineConeBuilder.GetPsdExpr()}
\label{\detokenize{csapi/affineconebuilder:affineconebuilder-getpsdexpr}}\begin{quote}

\sphinxAtStartPar
Get idx\sphinxhyphen{}th PSD expression in an affine cone constraint.

\sphinxAtStartPar
\sphinxstylestrong{Synopsis}
\begin{quote}

\sphinxAtStartPar
\sphinxcode{\sphinxupquote{PsdExpr GetPsdExpr(int idx)}}
\end{quote}

\sphinxAtStartPar
\sphinxstylestrong{Arguments}
\begin{quote}

\sphinxAtStartPar
\sphinxcode{\sphinxupquote{idx}}: index of PSD expression.
\end{quote}

\sphinxAtStartPar
\sphinxstylestrong{Return}
\begin{quote}

\sphinxAtStartPar
the idx\sphinxhyphen{}th PSD expression in an affine cone constraint.
\end{quote}
\end{quote}

\subsubsection{AffineConeBuilder.GetPsdExprs()}
\label{\detokenize{csapi/affineconebuilder:affineconebuilder-getpsdexprs}}\begin{quote}

\sphinxAtStartPar
Get all PSD expressions in an affine cone constraint.

\sphinxAtStartPar
\sphinxstylestrong{Synopsis}
\begin{quote}

\sphinxAtStartPar
\sphinxcode{\sphinxupquote{PsdExpr{[}{]} GetPsdExprs()}}
\end{quote}

\sphinxAtStartPar
\sphinxstylestrong{Return}
\begin{quote}

\sphinxAtStartPar
array of PSD expressions.
\end{quote}
\end{quote}

\subsubsection{AffineConeBuilder.GetSize()}
\label{\detokenize{csapi/affineconebuilder:affineconebuilder-getsize}}\begin{quote}

\sphinxAtStartPar
Get number of variables in an affine cone constraint.

\sphinxAtStartPar
\sphinxstylestrong{Synopsis}
\begin{quote}

\sphinxAtStartPar
\sphinxcode{\sphinxupquote{int GetSize()}}
\end{quote}

\sphinxAtStartPar
\sphinxstylestrong{Return}
\begin{quote}

\sphinxAtStartPar
number of variables.
\end{quote}
\end{quote}

\subsubsection{AffineConeBuilder.GetType()}
\label{\detokenize{csapi/affineconebuilder:affineconebuilder-gettype}}\begin{quote}

\sphinxAtStartPar
Get type of an affine cone constraint.

\sphinxAtStartPar
\sphinxstylestrong{Synopsis}
\begin{quote}

\sphinxAtStartPar
\sphinxcode{\sphinxupquote{int GetType()}}
\end{quote}

\sphinxAtStartPar
\sphinxstylestrong{Return}
\begin{quote}

\sphinxAtStartPar
type of the affine cone constraint.
\end{quote}
\end{quote}

\subsubsection{AffineConeBuilder.HasPsdTerm()}
\label{\detokenize{csapi/affineconebuilder:affineconebuilder-haspsdterm}}\begin{quote}

\sphinxAtStartPar
Check whether affine cone has PSD terms.

\sphinxAtStartPar
\sphinxstylestrong{Synopsis}
\begin{quote}

\sphinxAtStartPar
\sphinxcode{\sphinxupquote{bool HasPsdTerm()}}
\end{quote}

\sphinxAtStartPar
\sphinxstylestrong{Return}
\begin{quote}

\sphinxAtStartPar
flag to indicate whether affine cone has PSD terms.
\end{quote}
\end{quote}

\subsubsection{AffineConeBuilder.Set()}
\label{\detokenize{csapi/affineconebuilder:affineconebuilder-set}}\begin{quote}

\sphinxAtStartPar
Set linear expressions and type of an affine cone constraint.

\sphinxAtStartPar
\sphinxstylestrong{Synopsis}
\begin{quote}

\sphinxAtStartPar
\sphinxcode{\sphinxupquote{void Set(Expr{[}{]} exprs, int type)}}
\end{quote}

\sphinxAtStartPar
\sphinxstylestrong{Arguments}
\begin{quote}

\sphinxAtStartPar
\sphinxcode{\sphinxupquote{exprs}}: array of linear expressions.

\sphinxAtStartPar
\sphinxcode{\sphinxupquote{type}}: type of an affine cone constraint.
\end{quote}
\end{quote}

\subsubsection{AffineConeBuilder.Set()}
\label{\detokenize{csapi/affineconebuilder:id1}}\begin{quote}

\sphinxAtStartPar
Set PSD expressions and type of an affine cone constraint.

\sphinxAtStartPar
\sphinxstylestrong{Synopsis}
\begin{quote}

\sphinxAtStartPar
\sphinxcode{\sphinxupquote{void Set(PsdExpr{[}{]} exprs, int type)}}
\end{quote}

\sphinxAtStartPar
\sphinxstylestrong{Arguments}
\begin{quote}

\sphinxAtStartPar
\sphinxcode{\sphinxupquote{exprs}}: array of PSD expressions.

\sphinxAtStartPar
\sphinxcode{\sphinxupquote{type}}: type of an affine cone constraint.
\end{quote}
\end{quote}

\subsection{AffineConeBuilderArray}
\label{\detokenize{csharpapiref:affineconebuilderarray}}\label{\detokenize{csharpapiref:chapcsharpapiref-affineconebuilderarray}}
\sphinxAtStartPar
To facilitate operations on a group of C\# {\hyperref[\detokenize{csharpapiref:chapcsharpapiref-affineconebuilder}]{\sphinxcrossref{\DUrole{std,std-ref}{AffineConeBuilder}}}} objects,
the C\# interface of COPT introduces the \sphinxtitleref{AffineConeBuilderArray} class.
The following methods are provided:

\sphinxstepscope

\subsubsection{AffineConeBuilderArray.AffineConeBuilderArray()}
\label{\detokenize{csapi/affineconebuilderarray:affineconebuilderarray-affineconebuilderarray}}\label{\detokenize{csapi/affineconebuilderarray::doc}}\begin{quote}

\sphinxAtStartPar
Constructor of AffineConeBuilderArray object.

\sphinxAtStartPar
\sphinxstylestrong{Synopsis}
\begin{quote}

\sphinxAtStartPar
\sphinxcode{\sphinxupquote{AffineConeBuilderArray()}}
\end{quote}
\end{quote}

\subsubsection{AffineConeBuilderArray.GetBuilder()}
\label{\detokenize{csapi/affineconebuilderarray:affineconebuilderarray-getbuilder}}\begin{quote}

\sphinxAtStartPar
Get idx\sphinxhyphen{}th affine cone constraint builder object.

\sphinxAtStartPar
\sphinxstylestrong{Synopsis}
\begin{quote}

\sphinxAtStartPar
\sphinxcode{\sphinxupquote{AffineConeBuilder GetBuilder(int idx)}}
\end{quote}

\sphinxAtStartPar
\sphinxstylestrong{Arguments}
\begin{quote}

\sphinxAtStartPar
\sphinxcode{\sphinxupquote{idx}}: index of the affine cone constraint builder.
\end{quote}

\sphinxAtStartPar
\sphinxstylestrong{Return}
\begin{quote}

\sphinxAtStartPar
affine cone constraint builder object with index idx.
\end{quote}
\end{quote}

\subsubsection{AffineConeBuilderArray.PushBack()}
\label{\detokenize{csapi/affineconebuilderarray:affineconebuilderarray-pushback}}\begin{quote}

\sphinxAtStartPar
Add an affine cone constraint builder object to affine cone constraint builder array.

\sphinxAtStartPar
\sphinxstylestrong{Synopsis}
\begin{quote}

\sphinxAtStartPar
\sphinxcode{\sphinxupquote{void PushBack(AffineConeBuilder builder)}}
\end{quote}

\sphinxAtStartPar
\sphinxstylestrong{Arguments}
\begin{quote}

\sphinxAtStartPar
\sphinxcode{\sphinxupquote{builder}}: an affine cone constraint builder object.
\end{quote}
\end{quote}

\subsubsection{AffineConeBuilderArray.Size()}
\label{\detokenize{csapi/affineconebuilderarray:affineconebuilderarray-size}}\begin{quote}

\sphinxAtStartPar
Get the number of affine cone constraint builder objects.

\sphinxAtStartPar
\sphinxstylestrong{Synopsis}
\begin{quote}

\sphinxAtStartPar
\sphinxcode{\sphinxupquote{int Size()}}
\end{quote}

\sphinxAtStartPar
\sphinxstylestrong{Return}
\begin{quote}

\sphinxAtStartPar
number of affine cone constraint builder objects.
\end{quote}
\end{quote}

\subsection{QuadExpr}
\label{\detokenize{csharpapiref:quadexpr}}\label{\detokenize{csharpapiref:chapcsharpapiref-quadexpr}}
\sphinxAtStartPar
COPT quadratic expression object. A quadratic expression consists of a linear expression,
a list of variable pairs and associated coefficients of quadratic terms. Quadratic expressions
are used to build quadratic constraints and objectives.

\sphinxstepscope

\subsubsection{QuadExpr.QuadExpr()}
\label{\detokenize{csapi/quadexpr:quadexpr-quadexpr}}\label{\detokenize{csapi/quadexpr::doc}}\begin{quote}

\sphinxAtStartPar
Constructor of a quadratic expression with default constant value 0.

\sphinxAtStartPar
\sphinxstylestrong{Synopsis}
\begin{quote}

\sphinxAtStartPar
\sphinxcode{\sphinxupquote{QuadExpr(double constant)}}
\end{quote}

\sphinxAtStartPar
\sphinxstylestrong{Arguments}
\begin{quote}

\sphinxAtStartPar
\sphinxcode{\sphinxupquote{constant}}: optional, constant value in quadratic expression object.
\end{quote}
\end{quote}

\subsubsection{QuadExpr.QuadExpr()}
\label{\detokenize{csapi/quadexpr:id1}}\begin{quote}

\sphinxAtStartPar
Constructor of a quadratic expression with one linear term.

\sphinxAtStartPar
\sphinxstylestrong{Synopsis}
\begin{quote}

\sphinxAtStartPar
\sphinxcode{\sphinxupquote{QuadExpr(Var var, double coeff)}}
\end{quote}

\sphinxAtStartPar
\sphinxstylestrong{Arguments}
\begin{quote}

\sphinxAtStartPar
\sphinxcode{\sphinxupquote{var}}: variable of the added linear term.

\sphinxAtStartPar
\sphinxcode{\sphinxupquote{coeff}}: coefficent for the added linear term with default value 1.0.
\end{quote}
\end{quote}

\subsubsection{QuadExpr.QuadExpr()}
\label{\detokenize{csapi/quadexpr:id2}}\begin{quote}

\sphinxAtStartPar
Constructor of a quadratic expression with a linear expression.

\sphinxAtStartPar
\sphinxstylestrong{Synopsis}
\begin{quote}

\sphinxAtStartPar
\sphinxcode{\sphinxupquote{QuadExpr(Expr expr)}}
\end{quote}

\sphinxAtStartPar
\sphinxstylestrong{Arguments}
\begin{quote}

\sphinxAtStartPar
\sphinxcode{\sphinxupquote{expr}}: linear expression added to the quadratic expression.
\end{quote}
\end{quote}

\subsubsection{QuadExpr.QuadExpr()}
\label{\detokenize{csapi/quadexpr:id3}}\begin{quote}

\sphinxAtStartPar
Constructor of a quadratic expression with two linear expression.

\sphinxAtStartPar
\sphinxstylestrong{Synopsis}
\begin{quote}

\sphinxAtStartPar
\sphinxcode{\sphinxupquote{QuadExpr(Expr expr, Var var)}}
\end{quote}

\sphinxAtStartPar
\sphinxstylestrong{Arguments}
\begin{quote}

\sphinxAtStartPar
\sphinxcode{\sphinxupquote{expr}}: one linear expression.

\sphinxAtStartPar
\sphinxcode{\sphinxupquote{var}}: another variable.
\end{quote}
\end{quote}

\subsubsection{QuadExpr.QuadExpr()}
\label{\detokenize{csapi/quadexpr:id4}}\begin{quote}

\sphinxAtStartPar
Constructor of a quadratic expression with two linear expression.

\sphinxAtStartPar
\sphinxstylestrong{Synopsis}
\begin{quote}

\sphinxAtStartPar
\sphinxcode{\sphinxupquote{QuadExpr(Expr left, Expr right)}}
\end{quote}

\sphinxAtStartPar
\sphinxstylestrong{Arguments}
\begin{quote}

\sphinxAtStartPar
\sphinxcode{\sphinxupquote{left}}: one linear expression.

\sphinxAtStartPar
\sphinxcode{\sphinxupquote{right}}: another linear expression.
\end{quote}
\end{quote}

\subsubsection{QuadExpr.AddConstant()}
\label{\detokenize{csapi/quadexpr:quadexpr-addconstant}}\begin{quote}

\sphinxAtStartPar
Add a constant to the quadratic expression.

\sphinxAtStartPar
\sphinxstylestrong{Synopsis}
\begin{quote}

\sphinxAtStartPar
\sphinxcode{\sphinxupquote{void AddConstant(double constant)}}
\end{quote}

\sphinxAtStartPar
\sphinxstylestrong{Arguments}
\begin{quote}

\sphinxAtStartPar
\sphinxcode{\sphinxupquote{constant}}: value to be added.
\end{quote}
\end{quote}

\subsubsection{QuadExpr.AddLinExpr()}
\label{\detokenize{csapi/quadexpr:quadexpr-addlinexpr}}\begin{quote}

\sphinxAtStartPar
Add a linear expression to self.

\sphinxAtStartPar
\sphinxstylestrong{Synopsis}
\begin{quote}

\sphinxAtStartPar
\sphinxcode{\sphinxupquote{void AddLinExpr(Expr expr)}}
\end{quote}

\sphinxAtStartPar
\sphinxstylestrong{Arguments}
\begin{quote}

\sphinxAtStartPar
\sphinxcode{\sphinxupquote{expr}}: linear expression to be added.
\end{quote}
\end{quote}

\subsubsection{QuadExpr.AddLinExpr()}
\label{\detokenize{csapi/quadexpr:id5}}\begin{quote}

\sphinxAtStartPar
Add a linear expression to self.

\sphinxAtStartPar
\sphinxstylestrong{Synopsis}
\begin{quote}

\sphinxAtStartPar
\sphinxcode{\sphinxupquote{void AddLinExpr(Expr expr, double mult)}}
\end{quote}

\sphinxAtStartPar
\sphinxstylestrong{Arguments}
\begin{quote}

\sphinxAtStartPar
\sphinxcode{\sphinxupquote{expr}}: linear expression to be added.

\sphinxAtStartPar
\sphinxcode{\sphinxupquote{mult}}: multiplier constant.
\end{quote}
\end{quote}

\subsubsection{QuadExpr.AddQuadExpr()}
\label{\detokenize{csapi/quadexpr:quadexpr-addquadexpr}}\begin{quote}

\sphinxAtStartPar
Add a quadratic expression to self.

\sphinxAtStartPar
\sphinxstylestrong{Synopsis}
\begin{quote}

\sphinxAtStartPar
\sphinxcode{\sphinxupquote{void AddQuadExpr(QuadExpr expr)}}
\end{quote}

\sphinxAtStartPar
\sphinxstylestrong{Arguments}
\begin{quote}

\sphinxAtStartPar
\sphinxcode{\sphinxupquote{expr}}: quadratic expression to be added.
\end{quote}
\end{quote}

\subsubsection{QuadExpr.AddQuadExpr()}
\label{\detokenize{csapi/quadexpr:id6}}\begin{quote}

\sphinxAtStartPar
Add a quadratic expression to self.

\sphinxAtStartPar
\sphinxstylestrong{Synopsis}
\begin{quote}

\sphinxAtStartPar
\sphinxcode{\sphinxupquote{void AddQuadExpr(QuadExpr expr, double mult)}}
\end{quote}

\sphinxAtStartPar
\sphinxstylestrong{Arguments}
\begin{quote}

\sphinxAtStartPar
\sphinxcode{\sphinxupquote{expr}}: quadratic expression to be added.

\sphinxAtStartPar
\sphinxcode{\sphinxupquote{mult}}: multiplier constant.
\end{quote}
\end{quote}

\subsubsection{QuadExpr.AddTerm()}
\label{\detokenize{csapi/quadexpr:quadexpr-addterm}}\begin{quote}

\sphinxAtStartPar
Add a term to quadratic expression object.

\sphinxAtStartPar
\sphinxstylestrong{Synopsis}
\begin{quote}

\sphinxAtStartPar
\sphinxcode{\sphinxupquote{void AddTerm(Var var, double coeff)}}
\end{quote}

\sphinxAtStartPar
\sphinxstylestrong{Arguments}
\begin{quote}

\sphinxAtStartPar
\sphinxcode{\sphinxupquote{var}}: a variable of new term.

\sphinxAtStartPar
\sphinxcode{\sphinxupquote{coeff}}: coefficient of new term.
\end{quote}
\end{quote}

\subsubsection{QuadExpr.AddTerm()}
\label{\detokenize{csapi/quadexpr:id7}}\begin{quote}

\sphinxAtStartPar
Add a quadratic term to expression object.

\sphinxAtStartPar
\sphinxstylestrong{Synopsis}
\begin{quote}

\sphinxAtStartPar
\sphinxcode{\sphinxupquote{void AddTerm(}}
\begin{quote}

\sphinxAtStartPar
\sphinxcode{\sphinxupquote{Var var1,}}

\sphinxAtStartPar
\sphinxcode{\sphinxupquote{Var var2,}}

\sphinxAtStartPar
\sphinxcode{\sphinxupquote{double coeff)}}
\end{quote}
\end{quote}

\sphinxAtStartPar
\sphinxstylestrong{Arguments}
\begin{quote}

\sphinxAtStartPar
\sphinxcode{\sphinxupquote{var1}}: first variable of new quadratic term.

\sphinxAtStartPar
\sphinxcode{\sphinxupquote{var2}}: second variable of new quadratic term.

\sphinxAtStartPar
\sphinxcode{\sphinxupquote{coeff}}: coefficient of new quadratic term.
\end{quote}
\end{quote}

\subsubsection{QuadExpr.AddTerms()}
\label{\detokenize{csapi/quadexpr:quadexpr-addterms}}\begin{quote}

\sphinxAtStartPar
Add linear terms to quadratic expression object.

\sphinxAtStartPar
\sphinxstylestrong{Synopsis}
\begin{quote}

\sphinxAtStartPar
\sphinxcode{\sphinxupquote{void AddTerms(Var{[}{]} vars, double coeff)}}
\end{quote}

\sphinxAtStartPar
\sphinxstylestrong{Arguments}
\begin{quote}

\sphinxAtStartPar
\sphinxcode{\sphinxupquote{vars}}: variables of added linear terms.

\sphinxAtStartPar
\sphinxcode{\sphinxupquote{coeff}}: one coefficient for added linear terms.
\end{quote}
\end{quote}

\subsubsection{QuadExpr.AddTerms()}
\label{\detokenize{csapi/quadexpr:id8}}\begin{quote}

\sphinxAtStartPar
Add linear terms to quadratic expression object.

\sphinxAtStartPar
\sphinxstylestrong{Synopsis}
\begin{quote}

\sphinxAtStartPar
\sphinxcode{\sphinxupquote{void AddTerms(Var{[}{]} vars, double{[}{]} coeffs)}}
\end{quote}

\sphinxAtStartPar
\sphinxstylestrong{Arguments}
\begin{quote}

\sphinxAtStartPar
\sphinxcode{\sphinxupquote{vars}}: variables of added linear terms.

\sphinxAtStartPar
\sphinxcode{\sphinxupquote{coeffs}}: coefficients of added linear terms.
\end{quote}
\end{quote}

\subsubsection{QuadExpr.AddTerms()}
\label{\detokenize{csapi/quadexpr:id9}}\begin{quote}

\sphinxAtStartPar
Add linear terms to quadratic expression object.

\sphinxAtStartPar
\sphinxstylestrong{Synopsis}
\begin{quote}

\sphinxAtStartPar
\sphinxcode{\sphinxupquote{void AddTerms(VarArray vars, double coeff)}}
\end{quote}

\sphinxAtStartPar
\sphinxstylestrong{Arguments}
\begin{quote}

\sphinxAtStartPar
\sphinxcode{\sphinxupquote{vars}}: variables of added linear terms.

\sphinxAtStartPar
\sphinxcode{\sphinxupquote{coeff}}: one coefficient for added linear terms.
\end{quote}
\end{quote}

\subsubsection{QuadExpr.AddTerms()}
\label{\detokenize{csapi/quadexpr:id10}}\begin{quote}

\sphinxAtStartPar
Add linear terms to quadratic expression object.

\sphinxAtStartPar
\sphinxstylestrong{Synopsis}
\begin{quote}

\sphinxAtStartPar
\sphinxcode{\sphinxupquote{void AddTerms(VarArray vars, double{[}{]} coeffs)}}
\end{quote}

\sphinxAtStartPar
\sphinxstylestrong{Arguments}
\begin{quote}

\sphinxAtStartPar
\sphinxcode{\sphinxupquote{vars}}: variables of added terms.

\sphinxAtStartPar
\sphinxcode{\sphinxupquote{coeffs}}: coefficients of added terms.
\end{quote}
\end{quote}

\subsubsection{QuadExpr.AddTerms()}
\label{\detokenize{csapi/quadexpr:id11}}\begin{quote}

\sphinxAtStartPar
Add quadratic terms to expression object.

\sphinxAtStartPar
\sphinxstylestrong{Synopsis}
\begin{quote}

\sphinxAtStartPar
\sphinxcode{\sphinxupquote{void AddTerms(}}
\begin{quote}

\sphinxAtStartPar
\sphinxcode{\sphinxupquote{VarArray vars1,}}

\sphinxAtStartPar
\sphinxcode{\sphinxupquote{VarArray vars2,}}

\sphinxAtStartPar
\sphinxcode{\sphinxupquote{double{[}{]} coeffs)}}
\end{quote}
\end{quote}

\sphinxAtStartPar
\sphinxstylestrong{Arguments}
\begin{quote}

\sphinxAtStartPar
\sphinxcode{\sphinxupquote{vars1}}: first set of variables for added quadratic terms.

\sphinxAtStartPar
\sphinxcode{\sphinxupquote{vars2}}: second set of variables for added quadratic terms.

\sphinxAtStartPar
\sphinxcode{\sphinxupquote{coeffs}}: coefficient array for added quadratic terms.
\end{quote}
\end{quote}

\subsubsection{QuadExpr.AddTerms()}
\label{\detokenize{csapi/quadexpr:id12}}\begin{quote}

\sphinxAtStartPar
Add quadratic terms to expression object.

\sphinxAtStartPar
\sphinxstylestrong{Synopsis}
\begin{quote}

\sphinxAtStartPar
\sphinxcode{\sphinxupquote{void AddTerms(}}
\begin{quote}

\sphinxAtStartPar
\sphinxcode{\sphinxupquote{Var{[}{]} vars1,}}

\sphinxAtStartPar
\sphinxcode{\sphinxupquote{Var{[}{]} vars2,}}

\sphinxAtStartPar
\sphinxcode{\sphinxupquote{double{[}{]} coeffs)}}
\end{quote}
\end{quote}

\sphinxAtStartPar
\sphinxstylestrong{Arguments}
\begin{quote}

\sphinxAtStartPar
\sphinxcode{\sphinxupquote{vars1}}: first set of variables for added quadratic terms.

\sphinxAtStartPar
\sphinxcode{\sphinxupquote{vars2}}: second set of variables for added quadratic terms.

\sphinxAtStartPar
\sphinxcode{\sphinxupquote{coeffs}}: coefficient array for added quadratic terms.
\end{quote}
\end{quote}

\subsubsection{QuadExpr.Clone()}
\label{\detokenize{csapi/quadexpr:quadexpr-clone}}\begin{quote}

\sphinxAtStartPar
Deep copy quadratic expression object.

\sphinxAtStartPar
\sphinxstylestrong{Synopsis}
\begin{quote}

\sphinxAtStartPar
\sphinxcode{\sphinxupquote{QuadExpr Clone()}}
\end{quote}

\sphinxAtStartPar
\sphinxstylestrong{Return}
\begin{quote}

\sphinxAtStartPar
cloned quadratic expression object.
\end{quote}
\end{quote}

\subsubsection{QuadExpr.Divide()}
\label{\detokenize{csapi/quadexpr:quadexpr-divide}}\begin{quote}

\sphinxAtStartPar
Divide itself by a constant.

\sphinxAtStartPar
\sphinxstylestrong{Synopsis}
\begin{quote}

\sphinxAtStartPar
\sphinxcode{\sphinxupquote{void Divide(double c)}}
\end{quote}

\sphinxAtStartPar
\sphinxstylestrong{Arguments}
\begin{quote}

\sphinxAtStartPar
\sphinxcode{\sphinxupquote{c}}: constant operand.
\end{quote}
\end{quote}

\subsubsection{QuadExpr.Evaluate()}
\label{\detokenize{csapi/quadexpr:quadexpr-evaluate}}\begin{quote}

\sphinxAtStartPar
evaluate quadratic expression after solving.

\sphinxAtStartPar
\sphinxstylestrong{Synopsis}
\begin{quote}

\sphinxAtStartPar
\sphinxcode{\sphinxupquote{double Evaluate()}}
\end{quote}

\sphinxAtStartPar
\sphinxstylestrong{Return}
\begin{quote}

\sphinxAtStartPar
value of quadratic expression.
\end{quote}
\end{quote}

\subsubsection{QuadExpr.GetCoeff()}
\label{\detokenize{csapi/quadexpr:quadexpr-getcoeff}}\begin{quote}

\sphinxAtStartPar
Get coefficient from the i\sphinxhyphen{}th term in quadratic expression.

\sphinxAtStartPar
\sphinxstylestrong{Synopsis}
\begin{quote}

\sphinxAtStartPar
\sphinxcode{\sphinxupquote{double GetCoeff(int i)}}
\end{quote}

\sphinxAtStartPar
\sphinxstylestrong{Arguments}
\begin{quote}

\sphinxAtStartPar
\sphinxcode{\sphinxupquote{i}}: index of the term.
\end{quote}

\sphinxAtStartPar
\sphinxstylestrong{Return}
\begin{quote}

\sphinxAtStartPar
coefficient of the i\sphinxhyphen{}th term in quadratic expression object.
\end{quote}
\end{quote}

\subsubsection{QuadExpr.GetConstant()}
\label{\detokenize{csapi/quadexpr:quadexpr-getconstant}}\begin{quote}

\sphinxAtStartPar
Get constant in quadratic expression.

\sphinxAtStartPar
\sphinxstylestrong{Synopsis}
\begin{quote}

\sphinxAtStartPar
\sphinxcode{\sphinxupquote{double GetConstant()}}
\end{quote}

\sphinxAtStartPar
\sphinxstylestrong{Return}
\begin{quote}

\sphinxAtStartPar
constant in quadratic expression.
\end{quote}
\end{quote}

\subsubsection{QuadExpr.GetLinExpr()}
\label{\detokenize{csapi/quadexpr:quadexpr-getlinexpr}}\begin{quote}

\sphinxAtStartPar
Get linear expression in quadratic expression.

\sphinxAtStartPar
\sphinxstylestrong{Synopsis}
\begin{quote}

\sphinxAtStartPar
\sphinxcode{\sphinxupquote{Expr GetLinExpr()}}
\end{quote}

\sphinxAtStartPar
\sphinxstylestrong{Return}
\begin{quote}

\sphinxAtStartPar
linear expression object.
\end{quote}
\end{quote}

\subsubsection{QuadExpr.GetVar1()}
\label{\detokenize{csapi/quadexpr:quadexpr-getvar1}}\begin{quote}

\sphinxAtStartPar
Get first variable from the i\sphinxhyphen{}th term in quadratic expression.

\sphinxAtStartPar
\sphinxstylestrong{Synopsis}
\begin{quote}

\sphinxAtStartPar
\sphinxcode{\sphinxupquote{Var GetVar1(int i)}}
\end{quote}

\sphinxAtStartPar
\sphinxstylestrong{Arguments}
\begin{quote}

\sphinxAtStartPar
\sphinxcode{\sphinxupquote{i}}: index of the term.
\end{quote}

\sphinxAtStartPar
\sphinxstylestrong{Return}
\begin{quote}

\sphinxAtStartPar
first variable of the i\sphinxhyphen{}th term in quadratic expression object.
\end{quote}
\end{quote}

\subsubsection{QuadExpr.GetVar2()}
\label{\detokenize{csapi/quadexpr:quadexpr-getvar2}}\begin{quote}

\sphinxAtStartPar
Get second variable from the i\sphinxhyphen{}th term in quadratic expression.

\sphinxAtStartPar
\sphinxstylestrong{Synopsis}
\begin{quote}

\sphinxAtStartPar
\sphinxcode{\sphinxupquote{Var GetVar2(int i)}}
\end{quote}

\sphinxAtStartPar
\sphinxstylestrong{Arguments}
\begin{quote}

\sphinxAtStartPar
\sphinxcode{\sphinxupquote{i}}: index of the term.
\end{quote}

\sphinxAtStartPar
\sphinxstylestrong{Return}
\begin{quote}

\sphinxAtStartPar
second variable of the i\sphinxhyphen{}th term in quadratic expression object.
\end{quote}
\end{quote}

\subsubsection{QuadExpr.Multiply()}
\label{\detokenize{csapi/quadexpr:quadexpr-multiply}}\begin{quote}

\sphinxAtStartPar
Multiply itself by a constant.

\sphinxAtStartPar
\sphinxstylestrong{Synopsis}
\begin{quote}

\sphinxAtStartPar
\sphinxcode{\sphinxupquote{void Multiply(double c)}}
\end{quote}

\sphinxAtStartPar
\sphinxstylestrong{Arguments}
\begin{quote}

\sphinxAtStartPar
\sphinxcode{\sphinxupquote{c}}: constant operand.
\end{quote}
\end{quote}

\subsubsection{QuadExpr.Remove()}
\label{\detokenize{csapi/quadexpr:quadexpr-remove}}\begin{quote}

\sphinxAtStartPar
Remove idx\sphinxhyphen{}th term from quadratic expression object.

\sphinxAtStartPar
\sphinxstylestrong{Synopsis}
\begin{quote}

\sphinxAtStartPar
\sphinxcode{\sphinxupquote{void Remove(int idx)}}
\end{quote}

\sphinxAtStartPar
\sphinxstylestrong{Arguments}
\begin{quote}

\sphinxAtStartPar
\sphinxcode{\sphinxupquote{idx}}: index of the term to be removed.
\end{quote}
\end{quote}

\subsubsection{QuadExpr.Remove()}
\label{\detokenize{csapi/quadexpr:id13}}\begin{quote}

\sphinxAtStartPar
Remove the term associated with variable from quadratic expression.

\sphinxAtStartPar
\sphinxstylestrong{Synopsis}
\begin{quote}

\sphinxAtStartPar
\sphinxcode{\sphinxupquote{void Remove(Var var)}}
\end{quote}

\sphinxAtStartPar
\sphinxstylestrong{Arguments}
\begin{quote}

\sphinxAtStartPar
\sphinxcode{\sphinxupquote{var}}: a variable whose term should be removed.
\end{quote}
\end{quote}

\subsubsection{QuadExpr.SetCoeff()}
\label{\detokenize{csapi/quadexpr:quadexpr-setcoeff}}\begin{quote}

\sphinxAtStartPar
Set coefficient of the i\sphinxhyphen{}th term in quadratic expression.

\sphinxAtStartPar
\sphinxstylestrong{Synopsis}
\begin{quote}

\sphinxAtStartPar
\sphinxcode{\sphinxupquote{void SetCoeff(int i, double val)}}
\end{quote}

\sphinxAtStartPar
\sphinxstylestrong{Arguments}
\begin{quote}

\sphinxAtStartPar
\sphinxcode{\sphinxupquote{i}}: index of the quadratic term.

\sphinxAtStartPar
\sphinxcode{\sphinxupquote{val}}: coefficient of the term.
\end{quote}
\end{quote}

\subsubsection{QuadExpr.SetConstant()}
\label{\detokenize{csapi/quadexpr:quadexpr-setconstant}}\begin{quote}

\sphinxAtStartPar
Set constant for the quadratic expression.

\sphinxAtStartPar
\sphinxstylestrong{Synopsis}
\begin{quote}

\sphinxAtStartPar
\sphinxcode{\sphinxupquote{void SetConstant(double constant)}}
\end{quote}

\sphinxAtStartPar
\sphinxstylestrong{Arguments}
\begin{quote}

\sphinxAtStartPar
\sphinxcode{\sphinxupquote{constant}}: the value of the constant.
\end{quote}
\end{quote}

\subsubsection{QuadExpr.Size()}
\label{\detokenize{csapi/quadexpr:quadexpr-size}}\begin{quote}

\sphinxAtStartPar
Get number of terms in quadratic expression.

\sphinxAtStartPar
\sphinxstylestrong{Synopsis}
\begin{quote}

\sphinxAtStartPar
\sphinxcode{\sphinxupquote{long Size()}}
\end{quote}

\sphinxAtStartPar
\sphinxstylestrong{Return}
\begin{quote}

\sphinxAtStartPar
number of quadratic terms.
\end{quote}
\end{quote}

\subsection{QConstraint}
\label{\detokenize{csharpapiref:qconstraint}}\label{\detokenize{csharpapiref:chapcsharpapiref-qconstr}}
\sphinxAtStartPar
COPT quadratic constraint object. Quadratic constraints are always associated with a particular model.
User creates a quadratic constraint object by adding a quadratic constraint to a model,
rather than by using constructor of QConstraint class.

\sphinxstepscope

\subsubsection{QConstraint.Get()}
\label{\detokenize{csapi/qconstraint:qconstraint-get}}\label{\detokenize{csapi/qconstraint::doc}}\begin{quote}

\sphinxAtStartPar
Get information value of the quadratic constraint.

\sphinxAtStartPar
\sphinxstylestrong{Synopsis}
\begin{quote}

\sphinxAtStartPar
\sphinxcode{\sphinxupquote{double Get(string info)}}
\end{quote}

\sphinxAtStartPar
\sphinxstylestrong{Arguments}
\begin{quote}

\sphinxAtStartPar
\sphinxcode{\sphinxupquote{info}}: name of the information being queried.
\end{quote}

\sphinxAtStartPar
\sphinxstylestrong{Return}
\begin{quote}

\sphinxAtStartPar
information value.
\end{quote}
\end{quote}

\subsubsection{QConstraint.GetIdx()}
\label{\detokenize{csapi/qconstraint:qconstraint-getidx}}\begin{quote}

\sphinxAtStartPar
Get index of the quadratic constraint.

\sphinxAtStartPar
\sphinxstylestrong{Synopsis}
\begin{quote}

\sphinxAtStartPar
\sphinxcode{\sphinxupquote{int GetIdx()}}
\end{quote}

\sphinxAtStartPar
\sphinxstylestrong{Return}
\begin{quote}

\sphinxAtStartPar
the index of the quadratic constraint.
\end{quote}
\end{quote}

\subsubsection{QConstraint.GetName()}
\label{\detokenize{csapi/qconstraint:qconstraint-getname}}\begin{quote}

\sphinxAtStartPar
Get name of the quadratic constraint.

\sphinxAtStartPar
\sphinxstylestrong{Synopsis}
\begin{quote}

\sphinxAtStartPar
\sphinxcode{\sphinxupquote{string GetName()}}
\end{quote}

\sphinxAtStartPar
\sphinxstylestrong{Return}
\begin{quote}

\sphinxAtStartPar
the name of the quadratic constraint.
\end{quote}
\end{quote}

\subsubsection{QConstraint.GetRhs()}
\label{\detokenize{csapi/qconstraint:qconstraint-getrhs}}\begin{quote}

\sphinxAtStartPar
Get rhs of quadratic constraint.

\sphinxAtStartPar
\sphinxstylestrong{Synopsis}
\begin{quote}

\sphinxAtStartPar
\sphinxcode{\sphinxupquote{double GetRhs()}}
\end{quote}

\sphinxAtStartPar
\sphinxstylestrong{Return}
\begin{quote}

\sphinxAtStartPar
rhs of quadratic constraint.
\end{quote}
\end{quote}

\subsubsection{QConstraint.GetSense()}
\label{\detokenize{csapi/qconstraint:qconstraint-getsense}}\begin{quote}

\sphinxAtStartPar
Get rhs of quadratic constraint.

\sphinxAtStartPar
\sphinxstylestrong{Synopsis}
\begin{quote}

\sphinxAtStartPar
\sphinxcode{\sphinxupquote{char GetSense()}}
\end{quote}

\sphinxAtStartPar
\sphinxstylestrong{Return}
\begin{quote}

\sphinxAtStartPar
rhs of quadratic constraint.
\end{quote}
\end{quote}

\subsubsection{QConstraint.Remove()}
\label{\detokenize{csapi/qconstraint:qconstraint-remove}}\begin{quote}

\sphinxAtStartPar
Remove this constraint from model.

\sphinxAtStartPar
\sphinxstylestrong{Synopsis}
\begin{quote}

\sphinxAtStartPar
\sphinxcode{\sphinxupquote{void Remove()}}
\end{quote}
\end{quote}

\subsubsection{QConstraint.Set()}
\label{\detokenize{csapi/qconstraint:qconstraint-set}}\begin{quote}

\sphinxAtStartPar
Set information value of the quadratic constraint.

\sphinxAtStartPar
\sphinxstylestrong{Synopsis}
\begin{quote}

\sphinxAtStartPar
\sphinxcode{\sphinxupquote{void Set(string attr, double val)}}
\end{quote}

\sphinxAtStartPar
\sphinxstylestrong{Arguments}
\begin{quote}

\sphinxAtStartPar
\sphinxcode{\sphinxupquote{attr}}: name of the information.

\sphinxAtStartPar
\sphinxcode{\sphinxupquote{val}}: new information value.
\end{quote}
\end{quote}

\subsubsection{QConstraint.SetName()}
\label{\detokenize{csapi/qconstraint:qconstraint-setname}}\begin{quote}

\sphinxAtStartPar
Set name of quadratic constraint.

\sphinxAtStartPar
\sphinxstylestrong{Synopsis}
\begin{quote}

\sphinxAtStartPar
\sphinxcode{\sphinxupquote{void SetName(string name)}}
\end{quote}

\sphinxAtStartPar
\sphinxstylestrong{Arguments}
\begin{quote}

\sphinxAtStartPar
\sphinxcode{\sphinxupquote{name}}: the name to set.
\end{quote}
\end{quote}

\subsubsection{QConstraint.SetRhs()}
\label{\detokenize{csapi/qconstraint:qconstraint-setrhs}}\begin{quote}

\sphinxAtStartPar
Set rhs of quadratic constraint.

\sphinxAtStartPar
\sphinxstylestrong{Synopsis}
\begin{quote}

\sphinxAtStartPar
\sphinxcode{\sphinxupquote{void SetRhs(double rhs)}}
\end{quote}

\sphinxAtStartPar
\sphinxstylestrong{Arguments}
\begin{quote}

\sphinxAtStartPar
\sphinxcode{\sphinxupquote{rhs}}: rhs of quadratic constraint.
\end{quote}
\end{quote}

\subsubsection{QConstraint.SetSense()}
\label{\detokenize{csapi/qconstraint:qconstraint-setsense}}\begin{quote}

\sphinxAtStartPar
Set sense of quadratic constraint.

\sphinxAtStartPar
\sphinxstylestrong{Synopsis}
\begin{quote}

\sphinxAtStartPar
\sphinxcode{\sphinxupquote{void SetSense(char sense)}}
\end{quote}

\sphinxAtStartPar
\sphinxstylestrong{Arguments}
\begin{quote}

\sphinxAtStartPar
\sphinxcode{\sphinxupquote{sense}}: sense of quadratic constraint.
\end{quote}
\end{quote}

\subsection{QConstrArray}
\label{\detokenize{csharpapiref:qconstrarray}}\label{\detokenize{csharpapiref:chapcsharpapiref-qconstrarray}}
\sphinxAtStartPar
COPT quadratic constraint array object. To store and access a set of C\# {\hyperref[\detokenize{csharpapiref:chapcsharpapiref-qconstr}]{\sphinxcrossref{\DUrole{std,std-ref}{QConstraint}}}}
objects, Cardinal Optimizer provides C\# QConstrArray class, which defines the following methods.

\sphinxstepscope

\subsubsection{QConstrArray.QConstrArray()}
\label{\detokenize{csapi/qconstrarray:qconstrarray-qconstrarray}}\label{\detokenize{csapi/qconstrarray::doc}}\begin{quote}

\sphinxAtStartPar
Constructor of qconstrarray object.

\sphinxAtStartPar
\sphinxstylestrong{Synopsis}
\begin{quote}

\sphinxAtStartPar
\sphinxcode{\sphinxupquote{QConstrArray()}}
\end{quote}
\end{quote}

\subsubsection{QConstrArray.GetQConstr()}
\label{\detokenize{csapi/qconstrarray:qconstrarray-getqconstr}}\begin{quote}

\sphinxAtStartPar
Get idx\sphinxhyphen{}th quadratic constraint object.

\sphinxAtStartPar
\sphinxstylestrong{Synopsis}
\begin{quote}

\sphinxAtStartPar
\sphinxcode{\sphinxupquote{QConstraint GetQConstr(int idx)}}
\end{quote}

\sphinxAtStartPar
\sphinxstylestrong{Arguments}
\begin{quote}

\sphinxAtStartPar
\sphinxcode{\sphinxupquote{idx}}: index of the quadratic constraint.
\end{quote}

\sphinxAtStartPar
\sphinxstylestrong{Return}
\begin{quote}

\sphinxAtStartPar
constraint object with index ‘idx’.
\end{quote}
\end{quote}

\subsubsection{QConstrArray.PushBack()}
\label{\detokenize{csapi/qconstrarray:qconstrarray-pushback}}\begin{quote}

\sphinxAtStartPar
Add a quadratic constraint object to array.

\sphinxAtStartPar
\sphinxstylestrong{Synopsis}
\begin{quote}

\sphinxAtStartPar
\sphinxcode{\sphinxupquote{void PushBack(QConstraint constr)}}
\end{quote}

\sphinxAtStartPar
\sphinxstylestrong{Arguments}
\begin{quote}

\sphinxAtStartPar
\sphinxcode{\sphinxupquote{constr}}: a quadratic constraint object.
\end{quote}
\end{quote}

\subsubsection{QConstrArray.Size()}
\label{\detokenize{csapi/qconstrarray:qconstrarray-size}}\begin{quote}

\sphinxAtStartPar
Get the number of quadratic constraint objects.

\sphinxAtStartPar
\sphinxstylestrong{Synopsis}
\begin{quote}

\sphinxAtStartPar
\sphinxcode{\sphinxupquote{int Size()}}
\end{quote}

\sphinxAtStartPar
\sphinxstylestrong{Return}
\begin{quote}

\sphinxAtStartPar
number of quadratic constraint objects.
\end{quote}
\end{quote}

\subsection{QConstrBuilder}
\label{\detokenize{csharpapiref:qconstrbuilder}}\label{\detokenize{csharpapiref:chapcsharpapiref-qconstrbuilder}}
\sphinxAtStartPar
COPT quadratic constraint builder object. To help building a quadratic constraint, given a quadratic
expression, constraint sense and right\sphinxhyphen{}hand side value, Cardinal Optimizer provides C\# ConeBuilder
class, which defines the following methods.

\sphinxstepscope

\subsubsection{QConstrBuilder.QConstrBuilder()}
\label{\detokenize{csapi/qconstrbuilder:qconstrbuilder-qconstrbuilder}}\label{\detokenize{csapi/qconstrbuilder::doc}}\begin{quote}

\sphinxAtStartPar
Constructor of QConstrBuilder object.

\sphinxAtStartPar
\sphinxstylestrong{Synopsis}
\begin{quote}

\sphinxAtStartPar
\sphinxcode{\sphinxupquote{QConstrBuilder()}}
\end{quote}
\end{quote}

\subsubsection{QConstrBuilder.GetQuadExpr()}
\label{\detokenize{csapi/qconstrbuilder:qconstrbuilder-getquadexpr}}\begin{quote}

\sphinxAtStartPar
Get quadratic expression associated with constraint.

\sphinxAtStartPar
\sphinxstylestrong{Synopsis}
\begin{quote}

\sphinxAtStartPar
\sphinxcode{\sphinxupquote{QuadExpr GetQuadExpr()}}
\end{quote}

\sphinxAtStartPar
\sphinxstylestrong{Return}
\begin{quote}

\sphinxAtStartPar
quadratic expression object.
\end{quote}
\end{quote}

\subsubsection{QConstrBuilder.GetSense()}
\label{\detokenize{csapi/qconstrbuilder:qconstrbuilder-getsense}}\begin{quote}

\sphinxAtStartPar
Get sense associated with quadratic constraint.

\sphinxAtStartPar
\sphinxstylestrong{Synopsis}
\begin{quote}

\sphinxAtStartPar
\sphinxcode{\sphinxupquote{char GetSense()}}
\end{quote}

\sphinxAtStartPar
\sphinxstylestrong{Return}
\begin{quote}

\sphinxAtStartPar
quadratic constraint sense.
\end{quote}
\end{quote}

\subsubsection{QConstrBuilder.Set()}
\label{\detokenize{csapi/qconstrbuilder:qconstrbuilder-set}}\begin{quote}

\sphinxAtStartPar
Set detail of a quadratic constraint to its builder object.

\sphinxAtStartPar
\sphinxstylestrong{Synopsis}
\begin{quote}

\sphinxAtStartPar
\sphinxcode{\sphinxupquote{void Set(}}
\begin{quote}

\sphinxAtStartPar
\sphinxcode{\sphinxupquote{QuadExpr expr,}}

\sphinxAtStartPar
\sphinxcode{\sphinxupquote{char sense,}}

\sphinxAtStartPar
\sphinxcode{\sphinxupquote{double rhs)}}
\end{quote}
\end{quote}

\sphinxAtStartPar
\sphinxstylestrong{Arguments}
\begin{quote}

\sphinxAtStartPar
\sphinxcode{\sphinxupquote{expr}}: expression object at one side of the quadratic constraint.

\sphinxAtStartPar
\sphinxcode{\sphinxupquote{sense}}: quadratic constraint sense.

\sphinxAtStartPar
\sphinxcode{\sphinxupquote{rhs}}: constant of right side of quadratic constraint.
\end{quote}
\end{quote}

\subsection{QConstrBuilderArray}
\label{\detokenize{csharpapiref:qconstrbuilderarray}}\label{\detokenize{csharpapiref:chapcsharpapiref-qconstrbuilderarray}}
\sphinxAtStartPar
COPT quadratic constraint builder array object. To store and access a set of C\# {\hyperref[\detokenize{csharpapiref:chapcsharpapiref-qconstrbuilder}]{\sphinxcrossref{\DUrole{std,std-ref}{QConstrBuilder}}}}
objects, Cardinal Optimizer provides C\# QConstrBuilderArray class, which defines the following methods.

\sphinxstepscope

\subsubsection{QConstrBuilderArray.QConstrBuilderArray()}
\label{\detokenize{csapi/qconstrbuilderarray:qconstrbuilderarray-qconstrbuilderarray}}\label{\detokenize{csapi/qconstrbuilderarray::doc}}\begin{quote}

\sphinxAtStartPar
QConstructor of constrbuilderarray object.

\sphinxAtStartPar
\sphinxstylestrong{Synopsis}
\begin{quote}

\sphinxAtStartPar
\sphinxcode{\sphinxupquote{QConstrBuilderArray()}}
\end{quote}
\end{quote}

\subsubsection{QConstrBuilderArray.GetBuilder()}
\label{\detokenize{csapi/qconstrbuilderarray:qconstrbuilderarray-getbuilder}}\begin{quote}

\sphinxAtStartPar
Get idx\sphinxhyphen{}th quadratic constraint builder object.

\sphinxAtStartPar
\sphinxstylestrong{Synopsis}
\begin{quote}

\sphinxAtStartPar
\sphinxcode{\sphinxupquote{QConstrBuilder GetBuilder(int idx)}}
\end{quote}

\sphinxAtStartPar
\sphinxstylestrong{Arguments}
\begin{quote}

\sphinxAtStartPar
\sphinxcode{\sphinxupquote{idx}}: index of the quadratic constraint builder.
\end{quote}

\sphinxAtStartPar
\sphinxstylestrong{Return}
\begin{quote}

\sphinxAtStartPar
constraint builder object with index ‘idx’.
\end{quote}
\end{quote}

\subsubsection{QConstrBuilderArray.PushBack()}
\label{\detokenize{csapi/qconstrbuilderarray:qconstrbuilderarray-pushback}}\begin{quote}

\sphinxAtStartPar
Add a quadratic constraint builder object to constraint builder array.

\sphinxAtStartPar
\sphinxstylestrong{Synopsis}
\begin{quote}

\sphinxAtStartPar
\sphinxcode{\sphinxupquote{void PushBack(QConstrBuilder builder)}}
\end{quote}

\sphinxAtStartPar
\sphinxstylestrong{Arguments}
\begin{quote}

\sphinxAtStartPar
\sphinxcode{\sphinxupquote{builder}}: a quadratic constraint builder object.
\end{quote}
\end{quote}

\subsubsection{QConstrBuilderArray.Size()}
\label{\detokenize{csapi/qconstrbuilderarray:qconstrbuilderarray-size}}\begin{quote}

\sphinxAtStartPar
Get the number of quadratic constraint builder objects.

\sphinxAtStartPar
\sphinxstylestrong{Synopsis}
\begin{quote}

\sphinxAtStartPar
\sphinxcode{\sphinxupquote{int Size()}}
\end{quote}

\sphinxAtStartPar
\sphinxstylestrong{Return}
\begin{quote}

\sphinxAtStartPar
number of quadratic constraint builder objects.
\end{quote}
\end{quote}

\subsection{PsdVar}
\label{\detokenize{csharpapiref:psdvar}}\label{\detokenize{csharpapiref:chapcsharpapiref-psdvar}}
\sphinxAtStartPar
COPT PSD variable object. PSD variables are always associated with a particular model.
User creates a PSD variable object by adding a PSD variable to model, rather than
by constructor of PsdVar class.

\sphinxstepscope

\subsubsection{PsdVar.Diagonal()}
\label{\detokenize{csapi/psdvar:psdvar-diagonal}}\label{\detokenize{csapi/psdvar::doc}}\begin{quote}

\sphinxAtStartPar
Get diagonals of PSD variable.

\sphinxAtStartPar
\sphinxstylestrong{Synopsis}
\begin{quote}

\sphinxAtStartPar
\sphinxcode{\sphinxupquote{MPsdExpr Diagonal(int offset)}}
\end{quote}

\sphinxAtStartPar
\sphinxstylestrong{Arguments}
\begin{quote}

\sphinxAtStartPar
\sphinxcode{\sphinxupquote{offset}}: offset of the diagonal from the main diagonal. Can be positive or negative.
\end{quote}

\sphinxAtStartPar
\sphinxstylestrong{Return}
\begin{quote}

\sphinxAtStartPar
one\sphinxhyphen{}dimensional MPsdExpr object of diagonals.
\end{quote}
\end{quote}

\subsubsection{PsdVar.Get()}
\label{\detokenize{csapi/psdvar:psdvar-get}}\begin{quote}

\sphinxAtStartPar
Get information values of PSD variable.

\sphinxAtStartPar
\sphinxstylestrong{Synopsis}
\begin{quote}

\sphinxAtStartPar
\sphinxcode{\sphinxupquote{double{[}{]} Get(string info)}}
\end{quote}

\sphinxAtStartPar
\sphinxstylestrong{Arguments}
\begin{quote}

\sphinxAtStartPar
\sphinxcode{\sphinxupquote{info}}: name of information.
\end{quote}

\sphinxAtStartPar
\sphinxstylestrong{Return}
\begin{quote}

\sphinxAtStartPar
array of information values.
\end{quote}
\end{quote}

\subsubsection{PsdVar.GetDim()}
\label{\detokenize{csapi/psdvar:psdvar-getdim}}\begin{quote}

\sphinxAtStartPar
Get dimension of PSD variable.

\sphinxAtStartPar
\sphinxstylestrong{Synopsis}
\begin{quote}

\sphinxAtStartPar
\sphinxcode{\sphinxupquote{int GetDim()}}
\end{quote}

\sphinxAtStartPar
\sphinxstylestrong{Return}
\begin{quote}

\sphinxAtStartPar
dimension of PSD variable.
\end{quote}
\end{quote}

\subsubsection{PsdVar.GetIdx()}
\label{\detokenize{csapi/psdvar:psdvar-getidx}}\begin{quote}

\sphinxAtStartPar
Get index of PSD variable.

\sphinxAtStartPar
\sphinxstylestrong{Synopsis}
\begin{quote}

\sphinxAtStartPar
\sphinxcode{\sphinxupquote{int GetIdx()}}
\end{quote}

\sphinxAtStartPar
\sphinxstylestrong{Return}
\begin{quote}

\sphinxAtStartPar
index of PSD variable.
\end{quote}
\end{quote}

\subsubsection{PsdVar.GetInfo()}
\label{\detokenize{csapi/psdvar:psdvar-getinfo}}\begin{quote}

\sphinxAtStartPar
Get information values of PSD variable.

\sphinxAtStartPar
\sphinxstylestrong{Synopsis}
\begin{quote}

\sphinxAtStartPar
\sphinxcode{\sphinxupquote{NdArray\textless{}double\textgreater{} GetInfo(string name)}}
\end{quote}

\sphinxAtStartPar
\sphinxstylestrong{Arguments}
\begin{quote}

\sphinxAtStartPar
\sphinxcode{\sphinxupquote{name}}: name of informatioin.
\end{quote}

\sphinxAtStartPar
\sphinxstylestrong{Return}
\begin{quote}

\sphinxAtStartPar
2\sphinxhyphen{}dimensional NdArray with related information values.
\end{quote}
\end{quote}

\subsubsection{PsdVar.GetItem()}
\label{\detokenize{csapi/psdvar:psdvar-getitem}}\begin{quote}

\sphinxAtStartPar
Get corresponding MPsdExpr of given index from PSD variable.

\sphinxAtStartPar
\sphinxstylestrong{Synopsis}
\begin{quote}

\sphinxAtStartPar
\sphinxcode{\sphinxupquote{PsdExpr GetItem(int i, int j)}}
\end{quote}

\sphinxAtStartPar
\sphinxstylestrong{Arguments}
\begin{quote}

\sphinxAtStartPar
\sphinxcode{\sphinxupquote{i}}: row index.

\sphinxAtStartPar
\sphinxcode{\sphinxupquote{j}}: col index.
\end{quote}

\sphinxAtStartPar
\sphinxstylestrong{Return}
\begin{quote}

\sphinxAtStartPar
new PsdExpr object.
\end{quote}
\end{quote}

\subsubsection{PsdVar.GetLen()}
\label{\detokenize{csapi/psdvar:psdvar-getlen}}\begin{quote}

\sphinxAtStartPar
Get length of PSD variable.

\sphinxAtStartPar
\sphinxstylestrong{Synopsis}
\begin{quote}

\sphinxAtStartPar
\sphinxcode{\sphinxupquote{int GetLen()}}
\end{quote}

\sphinxAtStartPar
\sphinxstylestrong{Return}
\begin{quote}

\sphinxAtStartPar
length of PSD variable.
\end{quote}
\end{quote}

\subsubsection{PsdVar.GetName()}
\label{\detokenize{csapi/psdvar:psdvar-getname}}\begin{quote}

\sphinxAtStartPar
Get name of PSD variable.

\sphinxAtStartPar
\sphinxstylestrong{Synopsis}
\begin{quote}

\sphinxAtStartPar
\sphinxcode{\sphinxupquote{string GetName()}}
\end{quote}

\sphinxAtStartPar
\sphinxstylestrong{Return}
\begin{quote}

\sphinxAtStartPar
name of PSD variable.
\end{quote}
\end{quote}

\subsubsection{PsdVar.GetSize()}
\label{\detokenize{csapi/psdvar:psdvar-getsize}}\begin{quote}

\sphinxAtStartPar
Get size or length of PSD variable.

\sphinxAtStartPar
\sphinxstylestrong{Synopsis}
\begin{quote}

\sphinxAtStartPar
\sphinxcode{\sphinxupquote{int GetSize()}}
\end{quote}

\sphinxAtStartPar
\sphinxstylestrong{Return}
\begin{quote}

\sphinxAtStartPar
Size of PSD variable.
\end{quote}
\end{quote}

\subsubsection{PsdVar.Pick()}
\label{\detokenize{csapi/psdvar:psdvar-pick}}\begin{quote}

\sphinxAtStartPar
Given a list of indexes, pick elements of PSD variable in format of MPsdExpr.

\sphinxAtStartPar
\sphinxstylestrong{Synopsis}
\begin{quote}

\sphinxAtStartPar
\sphinxcode{\sphinxupquote{MPsdExpr Pick(NdArray\textless{}int\textgreater{} indexes)}}
\end{quote}

\sphinxAtStartPar
\sphinxstylestrong{Arguments}
\begin{quote}

\sphinxAtStartPar
\sphinxcode{\sphinxupquote{indexes}}: one or two dimensional indexes of elements. If two dimensional, each row is position of an element.
\end{quote}

\sphinxAtStartPar
\sphinxstylestrong{Return}
\begin{quote}

\sphinxAtStartPar
one\sphinxhyphen{}dimensional MPsdExpr object.
\end{quote}
\end{quote}

\subsubsection{PsdVar.Remove()}
\label{\detokenize{csapi/psdvar:psdvar-remove}}\begin{quote}

\sphinxAtStartPar
Remove PSD variable from model.

\sphinxAtStartPar
\sphinxstylestrong{Synopsis}
\begin{quote}

\sphinxAtStartPar
\sphinxcode{\sphinxupquote{void Remove()}}
\end{quote}
\end{quote}

\subsubsection{PsdVar.Slice()}
\label{\detokenize{csapi/psdvar:psdvar-slice}}\begin{quote}

\sphinxAtStartPar
Get corresponding MPsdExpr of given view from PSD variable.

\sphinxAtStartPar
\sphinxstylestrong{Synopsis}
\begin{quote}

\sphinxAtStartPar
\sphinxcode{\sphinxupquote{MPsdExpr Slice(View view)}}
\end{quote}

\sphinxAtStartPar
\sphinxstylestrong{Arguments}
\begin{quote}

\sphinxAtStartPar
\sphinxcode{\sphinxupquote{view}}: view of PSD variable.
\end{quote}

\sphinxAtStartPar
\sphinxstylestrong{Return}
\begin{quote}

\sphinxAtStartPar
new MPsdExpr object.
\end{quote}
\end{quote}

\subsubsection{PsdVar.Sum()}
\label{\detokenize{csapi/psdvar:psdvar-sum}}\begin{quote}

\sphinxAtStartPar
Sum of elements in PSD variable.

\sphinxAtStartPar
\sphinxstylestrong{Synopsis}
\begin{quote}

\sphinxAtStartPar
\sphinxcode{\sphinxupquote{PsdExpr Sum()}}
\end{quote}

\sphinxAtStartPar
\sphinxstylestrong{Return}
\begin{quote}

\sphinxAtStartPar
PsdExpr object.
\end{quote}
\end{quote}

\subsubsection{PsdVar.ToExpr()}
\label{\detokenize{csapi/psdvar:psdvar-toexpr}}\begin{quote}

\sphinxAtStartPar
convert PSD variable to 2\sphinxhyphen{}dimensional MPsdExpr.

\sphinxAtStartPar
\sphinxstylestrong{Synopsis}
\begin{quote}

\sphinxAtStartPar
\sphinxcode{\sphinxupquote{MPsdExpr ToExpr()}}
\end{quote}

\sphinxAtStartPar
\sphinxstylestrong{Return}
\begin{quote}

\sphinxAtStartPar
2\sphinxhyphen{}dimensional MPsdExpr object.
\end{quote}
\end{quote}

\subsection{PsdVarArray}
\label{\detokenize{csharpapiref:psdvararray}}\label{\detokenize{csharpapiref:chapcsharpapiref-psdvararray}}
\sphinxAtStartPar
COPT PSD variable array object. To store and access a set of {\hyperref[\detokenize{csharpapiref:chapcsharpapiref-psdvar}]{\sphinxcrossref{\DUrole{std,std-ref}{PsdVar}}}} objects,
Cardinal Optimizer provides PsdVarArray class, which defines the following methods.

\sphinxstepscope

\subsubsection{PsdVarArray.PsdVarArray()}
\label{\detokenize{csapi/psdvararray:psdvararray-psdvararray}}\label{\detokenize{csapi/psdvararray::doc}}\begin{quote}

\sphinxAtStartPar
Constructor of PsdVarArray.

\sphinxAtStartPar
\sphinxstylestrong{Synopsis}
\begin{quote}

\sphinxAtStartPar
\sphinxcode{\sphinxupquote{PsdVarArray()}}
\end{quote}
\end{quote}

\subsubsection{PsdVarArray.GetPsdVar()}
\label{\detokenize{csapi/psdvararray:psdvararray-getpsdvar}}\begin{quote}

\sphinxAtStartPar
Get idx\sphinxhyphen{}th PSD variable object.

\sphinxAtStartPar
\sphinxstylestrong{Synopsis}
\begin{quote}

\sphinxAtStartPar
\sphinxcode{\sphinxupquote{PsdVar GetPsdVar(int idx)}}
\end{quote}

\sphinxAtStartPar
\sphinxstylestrong{Arguments}
\begin{quote}

\sphinxAtStartPar
\sphinxcode{\sphinxupquote{idx}}: index of the PSD variable.
\end{quote}

\sphinxAtStartPar
\sphinxstylestrong{Return}
\begin{quote}

\sphinxAtStartPar
PSD variable object with index idx.
\end{quote}
\end{quote}

\subsubsection{PsdVarArray.PushBack()}
\label{\detokenize{csapi/psdvararray:psdvararray-pushback}}\begin{quote}

\sphinxAtStartPar
Add a PSD variable object to PSD variable array.

\sphinxAtStartPar
\sphinxstylestrong{Synopsis}
\begin{quote}

\sphinxAtStartPar
\sphinxcode{\sphinxupquote{void PushBack(PsdVar var)}}
\end{quote}

\sphinxAtStartPar
\sphinxstylestrong{Arguments}
\begin{quote}

\sphinxAtStartPar
\sphinxcode{\sphinxupquote{var}}: a PSD variable object.
\end{quote}
\end{quote}

\subsubsection{PsdVarArray.Reserve()}
\label{\detokenize{csapi/psdvararray:psdvararray-reserve}}\begin{quote}

\sphinxAtStartPar
Reserve capacity to contain at least n items.

\sphinxAtStartPar
\sphinxstylestrong{Synopsis}
\begin{quote}

\sphinxAtStartPar
\sphinxcode{\sphinxupquote{void Reserve(int n)}}
\end{quote}

\sphinxAtStartPar
\sphinxstylestrong{Arguments}
\begin{quote}

\sphinxAtStartPar
\sphinxcode{\sphinxupquote{n}}: minimum capacity for PSD variable object.
\end{quote}
\end{quote}

\subsubsection{PsdVarArray.Size()}
\label{\detokenize{csapi/psdvararray:psdvararray-size}}\begin{quote}

\sphinxAtStartPar
Get the number of PSD variable objects.

\sphinxAtStartPar
\sphinxstylestrong{Synopsis}
\begin{quote}

\sphinxAtStartPar
\sphinxcode{\sphinxupquote{int Size()}}
\end{quote}

\sphinxAtStartPar
\sphinxstylestrong{Return}
\begin{quote}

\sphinxAtStartPar
number of PSD variable objects.
\end{quote}
\end{quote}

\subsection{PsdExpr}
\label{\detokenize{csharpapiref:psdexpr}}\label{\detokenize{csharpapiref:chapcsharpapiref-psdexpr}}
\sphinxAtStartPar
COPT PSD expression object. A PSD expression consists of a linear expression,
a list of PSD variables and associated coefficient matrices of PSD terms. PSD expressions
are used to build PSD constraints and objectives.

\sphinxstepscope

\subsubsection{PsdExpr.PsdExpr()}
\label{\detokenize{csapi/psdexpr:psdexpr-psdexpr}}\label{\detokenize{csapi/psdexpr::doc}}\begin{quote}

\sphinxAtStartPar
Constructor of a PSD expression with default constant value 0.

\sphinxAtStartPar
\sphinxstylestrong{Synopsis}
\begin{quote}

\sphinxAtStartPar
\sphinxcode{\sphinxupquote{PsdExpr(double constant)}}
\end{quote}

\sphinxAtStartPar
\sphinxstylestrong{Arguments}
\begin{quote}

\sphinxAtStartPar
\sphinxcode{\sphinxupquote{constant}}: optional, constant value in PSD expression object.
\end{quote}
\end{quote}

\subsubsection{PsdExpr.PsdExpr()}
\label{\detokenize{csapi/psdexpr:id1}}\begin{quote}

\sphinxAtStartPar
Constructor of a PSD expression with one term.

\sphinxAtStartPar
\sphinxstylestrong{Synopsis}
\begin{quote}

\sphinxAtStartPar
\sphinxcode{\sphinxupquote{PsdExpr(Var var, double coeff)}}
\end{quote}

\sphinxAtStartPar
\sphinxstylestrong{Arguments}
\begin{quote}

\sphinxAtStartPar
\sphinxcode{\sphinxupquote{var}}: variable for the added term.

\sphinxAtStartPar
\sphinxcode{\sphinxupquote{coeff}}: coefficent for the added term.
\end{quote}
\end{quote}

\subsubsection{PsdExpr.PsdExpr()}
\label{\detokenize{csapi/psdexpr:id2}}\begin{quote}

\sphinxAtStartPar
Constructor of a PSD expression with a linear expression.

\sphinxAtStartPar
\sphinxstylestrong{Synopsis}
\begin{quote}

\sphinxAtStartPar
\sphinxcode{\sphinxupquote{PsdExpr(Expr expr)}}
\end{quote}

\sphinxAtStartPar
\sphinxstylestrong{Arguments}
\begin{quote}

\sphinxAtStartPar
\sphinxcode{\sphinxupquote{expr}}: input linear expression.
\end{quote}
\end{quote}

\subsubsection{PsdExpr.PsdExpr()}
\label{\detokenize{csapi/psdexpr:id3}}\begin{quote}

\sphinxAtStartPar
Constructor of a PSD expression with one term.

\sphinxAtStartPar
\sphinxstylestrong{Synopsis}
\begin{quote}

\sphinxAtStartPar
\sphinxcode{\sphinxupquote{PsdExpr(PsdVar var, SymMatrix mat)}}
\end{quote}

\sphinxAtStartPar
\sphinxstylestrong{Arguments}
\begin{quote}

\sphinxAtStartPar
\sphinxcode{\sphinxupquote{var}}: PSD variable for the added term.

\sphinxAtStartPar
\sphinxcode{\sphinxupquote{mat}}: coefficient matrix for the added term.
\end{quote}
\end{quote}

\subsubsection{PsdExpr.PsdExpr()}
\label{\detokenize{csapi/psdexpr:id4}}\begin{quote}

\sphinxAtStartPar
Constructor of a PSD expression with one term.

\sphinxAtStartPar
\sphinxstylestrong{Synopsis}
\begin{quote}

\sphinxAtStartPar
\sphinxcode{\sphinxupquote{PsdExpr(PsdVar var, SymMatExpr expr)}}
\end{quote}

\sphinxAtStartPar
\sphinxstylestrong{Arguments}
\begin{quote}

\sphinxAtStartPar
\sphinxcode{\sphinxupquote{var}}: PSD variable for the added term.

\sphinxAtStartPar
\sphinxcode{\sphinxupquote{expr}}: coefficient expression of symmetric matrices of new PSD term.
\end{quote}
\end{quote}

\subsubsection{PsdExpr.AddConstant()}
\label{\detokenize{csapi/psdexpr:psdexpr-addconstant}}\begin{quote}

\sphinxAtStartPar
Add constant to the PSD expression.

\sphinxAtStartPar
\sphinxstylestrong{Synopsis}
\begin{quote}

\sphinxAtStartPar
\sphinxcode{\sphinxupquote{void AddConstant(double constant)}}
\end{quote}

\sphinxAtStartPar
\sphinxstylestrong{Arguments}
\begin{quote}

\sphinxAtStartPar
\sphinxcode{\sphinxupquote{constant}}: value to be added.
\end{quote}
\end{quote}

\subsubsection{PsdExpr.AddLinExpr()}
\label{\detokenize{csapi/psdexpr:psdexpr-addlinexpr}}\begin{quote}

\sphinxAtStartPar
Add a linear expression to PSD expression object.

\sphinxAtStartPar
\sphinxstylestrong{Synopsis}
\begin{quote}

\sphinxAtStartPar
\sphinxcode{\sphinxupquote{void AddLinExpr(Expr expr)}}
\end{quote}

\sphinxAtStartPar
\sphinxstylestrong{Arguments}
\begin{quote}

\sphinxAtStartPar
\sphinxcode{\sphinxupquote{expr}}: linear expression to be added.
\end{quote}
\end{quote}

\subsubsection{PsdExpr.AddLinExpr()}
\label{\detokenize{csapi/psdexpr:id5}}\begin{quote}

\sphinxAtStartPar
Add a linear expression to PSD expression object.

\sphinxAtStartPar
\sphinxstylestrong{Synopsis}
\begin{quote}

\sphinxAtStartPar
\sphinxcode{\sphinxupquote{void AddLinExpr(Expr expr, double mult)}}
\end{quote}

\sphinxAtStartPar
\sphinxstylestrong{Arguments}
\begin{quote}

\sphinxAtStartPar
\sphinxcode{\sphinxupquote{expr}}: linear expression to be added.

\sphinxAtStartPar
\sphinxcode{\sphinxupquote{mult}}: multiplier constant.
\end{quote}
\end{quote}

\subsubsection{PsdExpr.AddPsdExpr()}
\label{\detokenize{csapi/psdexpr:psdexpr-addpsdexpr}}\begin{quote}

\sphinxAtStartPar
Add a PSD expression to self.

\sphinxAtStartPar
\sphinxstylestrong{Synopsis}
\begin{quote}

\sphinxAtStartPar
\sphinxcode{\sphinxupquote{void AddPsdExpr(PsdExpr expr)}}
\end{quote}

\sphinxAtStartPar
\sphinxstylestrong{Arguments}
\begin{quote}

\sphinxAtStartPar
\sphinxcode{\sphinxupquote{expr}}: PSD expression to be added.
\end{quote}
\end{quote}

\subsubsection{PsdExpr.AddPsdExpr()}
\label{\detokenize{csapi/psdexpr:id6}}\begin{quote}

\sphinxAtStartPar
Add a PSD expression to self.

\sphinxAtStartPar
\sphinxstylestrong{Synopsis}
\begin{quote}

\sphinxAtStartPar
\sphinxcode{\sphinxupquote{void AddPsdExpr(PsdExpr expr, double mult)}}
\end{quote}

\sphinxAtStartPar
\sphinxstylestrong{Arguments}
\begin{quote}

\sphinxAtStartPar
\sphinxcode{\sphinxupquote{expr}}: PSD expression to be added.

\sphinxAtStartPar
\sphinxcode{\sphinxupquote{mult}}: multiplier constant.
\end{quote}
\end{quote}

\subsubsection{PsdExpr.AddTerm()}
\label{\detokenize{csapi/psdexpr:psdexpr-addterm}}\begin{quote}

\sphinxAtStartPar
Add a linear term to PSD expression object.

\sphinxAtStartPar
\sphinxstylestrong{Synopsis}
\begin{quote}

\sphinxAtStartPar
\sphinxcode{\sphinxupquote{void AddTerm(Var var, double coeff)}}
\end{quote}

\sphinxAtStartPar
\sphinxstylestrong{Arguments}
\begin{quote}

\sphinxAtStartPar
\sphinxcode{\sphinxupquote{var}}: variable of new linear term.

\sphinxAtStartPar
\sphinxcode{\sphinxupquote{coeff}}: coefficient of new linear term.
\end{quote}
\end{quote}

\subsubsection{PsdExpr.AddTerm()}
\label{\detokenize{csapi/psdexpr:id7}}\begin{quote}

\sphinxAtStartPar
Add a PSD term to PSD expression object.

\sphinxAtStartPar
\sphinxstylestrong{Synopsis}
\begin{quote}

\sphinxAtStartPar
\sphinxcode{\sphinxupquote{void AddTerm(PsdVar var, SymMatrix mat)}}
\end{quote}

\sphinxAtStartPar
\sphinxstylestrong{Arguments}
\begin{quote}

\sphinxAtStartPar
\sphinxcode{\sphinxupquote{var}}: PSD variable of new PSD term.

\sphinxAtStartPar
\sphinxcode{\sphinxupquote{mat}}: coefficient matrix of new PSD term.
\end{quote}
\end{quote}

\subsubsection{PsdExpr.AddTerm()}
\label{\detokenize{csapi/psdexpr:id8}}\begin{quote}

\sphinxAtStartPar
Add a PSD term to PSD expression object.

\sphinxAtStartPar
\sphinxstylestrong{Synopsis}
\begin{quote}

\sphinxAtStartPar
\sphinxcode{\sphinxupquote{void AddTerm(PsdVar var, SymMatExpr expr)}}
\end{quote}

\sphinxAtStartPar
\sphinxstylestrong{Arguments}
\begin{quote}

\sphinxAtStartPar
\sphinxcode{\sphinxupquote{var}}: PSD variable of new PSD term.

\sphinxAtStartPar
\sphinxcode{\sphinxupquote{expr}}: coefficient expression of symmetric matrices of new PSD term.
\end{quote}
\end{quote}

\subsubsection{PsdExpr.AddTerms()}
\label{\detokenize{csapi/psdexpr:psdexpr-addterms}}\begin{quote}

\sphinxAtStartPar
Add linear terms to PSD expression object.

\sphinxAtStartPar
\sphinxstylestrong{Synopsis}
\begin{quote}

\sphinxAtStartPar
\sphinxcode{\sphinxupquote{void AddTerms(Var{[}{]} vars, double coeff)}}
\end{quote}

\sphinxAtStartPar
\sphinxstylestrong{Arguments}
\begin{quote}

\sphinxAtStartPar
\sphinxcode{\sphinxupquote{vars}}: variables of added linear terms.

\sphinxAtStartPar
\sphinxcode{\sphinxupquote{coeff}}: one coefficient for added linear terms.
\end{quote}
\end{quote}

\subsubsection{PsdExpr.AddTerms()}
\label{\detokenize{csapi/psdexpr:id9}}\begin{quote}

\sphinxAtStartPar
Add linear terms to PSD expression object.

\sphinxAtStartPar
\sphinxstylestrong{Synopsis}
\begin{quote}

\sphinxAtStartPar
\sphinxcode{\sphinxupquote{void AddTerms(Var{[}{]} vars, double{[}{]} coeffs)}}
\end{quote}

\sphinxAtStartPar
\sphinxstylestrong{Arguments}
\begin{quote}

\sphinxAtStartPar
\sphinxcode{\sphinxupquote{vars}}: variables for added linear terms.

\sphinxAtStartPar
\sphinxcode{\sphinxupquote{coeffs}}: coefficient array for added linear terms.
\end{quote}
\end{quote}

\subsubsection{PsdExpr.AddTerms()}
\label{\detokenize{csapi/psdexpr:id10}}\begin{quote}

\sphinxAtStartPar
Add linear terms to PSD expression object.

\sphinxAtStartPar
\sphinxstylestrong{Synopsis}
\begin{quote}

\sphinxAtStartPar
\sphinxcode{\sphinxupquote{void AddTerms(VarArray vars, double coeff)}}
\end{quote}

\sphinxAtStartPar
\sphinxstylestrong{Arguments}
\begin{quote}

\sphinxAtStartPar
\sphinxcode{\sphinxupquote{vars}}: variables of added linear terms.

\sphinxAtStartPar
\sphinxcode{\sphinxupquote{coeff}}: one coefficient for added linear terms.
\end{quote}
\end{quote}

\subsubsection{PsdExpr.AddTerms()}
\label{\detokenize{csapi/psdexpr:id11}}\begin{quote}

\sphinxAtStartPar
Add linear terms to PSD expression object.

\sphinxAtStartPar
\sphinxstylestrong{Synopsis}
\begin{quote}

\sphinxAtStartPar
\sphinxcode{\sphinxupquote{void AddTerms(VarArray vars, double{[}{]} coeffs)}}
\end{quote}

\sphinxAtStartPar
\sphinxstylestrong{Arguments}
\begin{quote}

\sphinxAtStartPar
\sphinxcode{\sphinxupquote{vars}}: variables of added terms.

\sphinxAtStartPar
\sphinxcode{\sphinxupquote{coeffs}}: coefficients of added terms.
\end{quote}
\end{quote}

\subsubsection{PsdExpr.AddTerms()}
\label{\detokenize{csapi/psdexpr:id12}}\begin{quote}

\sphinxAtStartPar
Add PSD terms to PSD expression object.

\sphinxAtStartPar
\sphinxstylestrong{Synopsis}
\begin{quote}

\sphinxAtStartPar
\sphinxcode{\sphinxupquote{void AddTerms(PsdVarArray vars, SymMatrixArray mats)}}
\end{quote}

\sphinxAtStartPar
\sphinxstylestrong{Arguments}
\begin{quote}

\sphinxAtStartPar
\sphinxcode{\sphinxupquote{vars}}: PSD variables for added PSD terms.

\sphinxAtStartPar
\sphinxcode{\sphinxupquote{mats}}: coefficient matrixes for added PSD terms.
\end{quote}
\end{quote}

\subsubsection{PsdExpr.AddTerms()}
\label{\detokenize{csapi/psdexpr:id13}}\begin{quote}

\sphinxAtStartPar
Add PSD terms to PSD expression object.

\sphinxAtStartPar
\sphinxstylestrong{Synopsis}
\begin{quote}

\sphinxAtStartPar
\sphinxcode{\sphinxupquote{void AddTerms(PsdVar{[}{]} vars, SymMatrix{[}{]} mats)}}
\end{quote}

\sphinxAtStartPar
\sphinxstylestrong{Arguments}
\begin{quote}

\sphinxAtStartPar
\sphinxcode{\sphinxupquote{vars}}: PSD variables for added PSD terms.

\sphinxAtStartPar
\sphinxcode{\sphinxupquote{mats}}: coefficient matrixes for added PSD terms.
\end{quote}
\end{quote}

\subsubsection{PsdExpr.Clone()}
\label{\detokenize{csapi/psdexpr:psdexpr-clone}}\begin{quote}

\sphinxAtStartPar
Deep copy PSD expression object.

\sphinxAtStartPar
\sphinxstylestrong{Synopsis}
\begin{quote}

\sphinxAtStartPar
\sphinxcode{\sphinxupquote{PsdExpr Clone()}}
\end{quote}

\sphinxAtStartPar
\sphinxstylestrong{Return}
\begin{quote}

\sphinxAtStartPar
cloned PSD expression object.
\end{quote}
\end{quote}

\subsubsection{PsdExpr.Evaluate()}
\label{\detokenize{csapi/psdexpr:psdexpr-evaluate}}\begin{quote}

\sphinxAtStartPar
Evaluate PSD expression after solving.

\sphinxAtStartPar
\sphinxstylestrong{Synopsis}
\begin{quote}

\sphinxAtStartPar
\sphinxcode{\sphinxupquote{double Evaluate()}}
\end{quote}

\sphinxAtStartPar
\sphinxstylestrong{Return}
\begin{quote}

\sphinxAtStartPar
value of PSD expression.
\end{quote}
\end{quote}

\subsubsection{PsdExpr.GetCoeff()}
\label{\detokenize{csapi/psdexpr:psdexpr-getcoeff}}\begin{quote}

\sphinxAtStartPar
Get coefficient from the i\sphinxhyphen{}th term in PSD expression.

\sphinxAtStartPar
\sphinxstylestrong{Synopsis}
\begin{quote}

\sphinxAtStartPar
\sphinxcode{\sphinxupquote{SymMatExpr GetCoeff(int i)}}
\end{quote}

\sphinxAtStartPar
\sphinxstylestrong{Arguments}
\begin{quote}

\sphinxAtStartPar
\sphinxcode{\sphinxupquote{i}}: index of the PSD term.
\end{quote}

\sphinxAtStartPar
\sphinxstylestrong{Return}
\begin{quote}

\sphinxAtStartPar
coefficient expression of the i\sphinxhyphen{}th PSD term.
\end{quote}
\end{quote}

\subsubsection{PsdExpr.GetConstant()}
\label{\detokenize{csapi/psdexpr:psdexpr-getconstant}}\begin{quote}

\sphinxAtStartPar
Get constant in PSD expression.

\sphinxAtStartPar
\sphinxstylestrong{Synopsis}
\begin{quote}

\sphinxAtStartPar
\sphinxcode{\sphinxupquote{double GetConstant()}}
\end{quote}

\sphinxAtStartPar
\sphinxstylestrong{Return}
\begin{quote}

\sphinxAtStartPar
constant in PSD expression.
\end{quote}
\end{quote}

\subsubsection{PsdExpr.GetLinExpr()}
\label{\detokenize{csapi/psdexpr:psdexpr-getlinexpr}}\begin{quote}

\sphinxAtStartPar
Get linear expression in PSD expression.

\sphinxAtStartPar
\sphinxstylestrong{Synopsis}
\begin{quote}

\sphinxAtStartPar
\sphinxcode{\sphinxupquote{Expr GetLinExpr()}}
\end{quote}

\sphinxAtStartPar
\sphinxstylestrong{Return}
\begin{quote}

\sphinxAtStartPar
linear expression object.
\end{quote}
\end{quote}

\subsubsection{PsdExpr.GetPsdVar()}
\label{\detokenize{csapi/psdexpr:psdexpr-getpsdvar}}\begin{quote}

\sphinxAtStartPar
Get the PSD variable from the i\sphinxhyphen{}th term in PSD expression.

\sphinxAtStartPar
\sphinxstylestrong{Synopsis}
\begin{quote}

\sphinxAtStartPar
\sphinxcode{\sphinxupquote{PsdVar GetPsdVar(int i)}}
\end{quote}

\sphinxAtStartPar
\sphinxstylestrong{Arguments}
\begin{quote}

\sphinxAtStartPar
\sphinxcode{\sphinxupquote{i}}: index of the term.
\end{quote}

\sphinxAtStartPar
\sphinxstylestrong{Return}
\begin{quote}

\sphinxAtStartPar
the first variable of the i\sphinxhyphen{}th term in PSD expression object.
\end{quote}
\end{quote}

\subsubsection{PsdExpr.Multiply()}
\label{\detokenize{csapi/psdexpr:psdexpr-multiply}}\begin{quote}

\sphinxAtStartPar
Multiply itself by a constant.

\sphinxAtStartPar
\sphinxstylestrong{Synopsis}
\begin{quote}

\sphinxAtStartPar
\sphinxcode{\sphinxupquote{void Multiply(double c)}}
\end{quote}

\sphinxAtStartPar
\sphinxstylestrong{Arguments}
\begin{quote}

\sphinxAtStartPar
\sphinxcode{\sphinxupquote{c}}: constant operand.
\end{quote}
\end{quote}

\subsubsection{PsdExpr.Remove()}
\label{\detokenize{csapi/psdexpr:psdexpr-remove}}\begin{quote}

\sphinxAtStartPar
Remove i\sphinxhyphen{}th term from PSD expression object.

\sphinxAtStartPar
\sphinxstylestrong{Synopsis}
\begin{quote}

\sphinxAtStartPar
\sphinxcode{\sphinxupquote{void Remove(int idx)}}
\end{quote}

\sphinxAtStartPar
\sphinxstylestrong{Arguments}
\begin{quote}

\sphinxAtStartPar
\sphinxcode{\sphinxupquote{idx}}: index of the term to be removed.
\end{quote}
\end{quote}

\subsubsection{PsdExpr.Remove()}
\label{\detokenize{csapi/psdexpr:id14}}\begin{quote}

\sphinxAtStartPar
Remove the term associated with variable from PSD expression.

\sphinxAtStartPar
\sphinxstylestrong{Synopsis}
\begin{quote}

\sphinxAtStartPar
\sphinxcode{\sphinxupquote{void Remove(Var var)}}
\end{quote}

\sphinxAtStartPar
\sphinxstylestrong{Arguments}
\begin{quote}

\sphinxAtStartPar
\sphinxcode{\sphinxupquote{var}}: a variable whose term should be removed.
\end{quote}
\end{quote}

\subsubsection{PsdExpr.Remove()}
\label{\detokenize{csapi/psdexpr:id15}}\begin{quote}

\sphinxAtStartPar
Remove the term associated with PSD variable from PSD expression.

\sphinxAtStartPar
\sphinxstylestrong{Synopsis}
\begin{quote}

\sphinxAtStartPar
\sphinxcode{\sphinxupquote{void Remove(PsdVar var)}}
\end{quote}

\sphinxAtStartPar
\sphinxstylestrong{Arguments}
\begin{quote}

\sphinxAtStartPar
\sphinxcode{\sphinxupquote{var}}: a PSD variable whose term should be removed.
\end{quote}
\end{quote}

\subsubsection{PsdExpr.SetCoeff()}
\label{\detokenize{csapi/psdexpr:psdexpr-setcoeff}}\begin{quote}

\sphinxAtStartPar
Set coefficient matrix of the i\sphinxhyphen{}th term in PSD expression.

\sphinxAtStartPar
\sphinxstylestrong{Synopsis}
\begin{quote}

\sphinxAtStartPar
\sphinxcode{\sphinxupquote{void SetCoeff(int i, SymMatrix mat)}}
\end{quote}

\sphinxAtStartPar
\sphinxstylestrong{Arguments}
\begin{quote}

\sphinxAtStartPar
\sphinxcode{\sphinxupquote{i}}: index of the PSD term.

\sphinxAtStartPar
\sphinxcode{\sphinxupquote{mat}}: coefficient matrix of the term.
\end{quote}
\end{quote}

\subsubsection{PsdExpr.SetConstant()}
\label{\detokenize{csapi/psdexpr:psdexpr-setconstant}}\begin{quote}

\sphinxAtStartPar
Set constant for the PSD expression.

\sphinxAtStartPar
\sphinxstylestrong{Synopsis}
\begin{quote}

\sphinxAtStartPar
\sphinxcode{\sphinxupquote{void SetConstant(double constant)}}
\end{quote}

\sphinxAtStartPar
\sphinxstylestrong{Arguments}
\begin{quote}

\sphinxAtStartPar
\sphinxcode{\sphinxupquote{constant}}: the value of the constant.
\end{quote}
\end{quote}

\subsubsection{PsdExpr.Size()}
\label{\detokenize{csapi/psdexpr:psdexpr-size}}\begin{quote}

\sphinxAtStartPar
Get number of PSD terms in expression.

\sphinxAtStartPar
\sphinxstylestrong{Synopsis}
\begin{quote}

\sphinxAtStartPar
\sphinxcode{\sphinxupquote{long Size()}}
\end{quote}

\sphinxAtStartPar
\sphinxstylestrong{Return}
\begin{quote}

\sphinxAtStartPar
number of PSD terms.
\end{quote}
\end{quote}

\subsection{PsdConstraint}
\label{\detokenize{csharpapiref:psdconstraint}}\label{\detokenize{csharpapiref:chapcsharpapiref-psdconstraint}}
\sphinxAtStartPar
COPT PSD constraint object. PSD constraints are always associated with a particular model.
User creates a PSD constraint object by adding a PSD constraint to model,
rather than by constructor of PsdConstraint class.

\sphinxstepscope

\subsubsection{PsdConstraint.Get()}
\label{\detokenize{csapi/psdconstraint:psdconstraint-get}}\label{\detokenize{csapi/psdconstraint::doc}}\begin{quote}

\sphinxAtStartPar
Get information value of the PSD constraint. Support related PSD informations.

\sphinxAtStartPar
\sphinxstylestrong{Synopsis}
\begin{quote}

\sphinxAtStartPar
\sphinxcode{\sphinxupquote{double Get(string info)}}
\end{quote}

\sphinxAtStartPar
\sphinxstylestrong{Arguments}
\begin{quote}

\sphinxAtStartPar
\sphinxcode{\sphinxupquote{info}}: name of queried information.
\end{quote}

\sphinxAtStartPar
\sphinxstylestrong{Return}
\begin{quote}

\sphinxAtStartPar
information value.
\end{quote}
\end{quote}

\subsubsection{PsdConstraint.GetIdx()}
\label{\detokenize{csapi/psdconstraint:psdconstraint-getidx}}\begin{quote}

\sphinxAtStartPar
Get index of the PSD constraint.

\sphinxAtStartPar
\sphinxstylestrong{Synopsis}
\begin{quote}

\sphinxAtStartPar
\sphinxcode{\sphinxupquote{int GetIdx()}}
\end{quote}

\sphinxAtStartPar
\sphinxstylestrong{Return}
\begin{quote}

\sphinxAtStartPar
the index of the PSD constraint.
\end{quote}
\end{quote}

\subsubsection{PsdConstraint.GetName()}
\label{\detokenize{csapi/psdconstraint:psdconstraint-getname}}\begin{quote}

\sphinxAtStartPar
Get name of the PSD constraint.

\sphinxAtStartPar
\sphinxstylestrong{Synopsis}
\begin{quote}

\sphinxAtStartPar
\sphinxcode{\sphinxupquote{string GetName()}}
\end{quote}

\sphinxAtStartPar
\sphinxstylestrong{Return}
\begin{quote}

\sphinxAtStartPar
the name of the PSD constraint.
\end{quote}
\end{quote}

\subsubsection{PsdConstraint.Remove()}
\label{\detokenize{csapi/psdconstraint:psdconstraint-remove}}\begin{quote}

\sphinxAtStartPar
Remove this PSD constraint from model.

\sphinxAtStartPar
\sphinxstylestrong{Synopsis}
\begin{quote}

\sphinxAtStartPar
\sphinxcode{\sphinxupquote{void Remove()}}
\end{quote}
\end{quote}

\subsubsection{PsdConstraint.Set()}
\label{\detokenize{csapi/psdconstraint:psdconstraint-set}}\begin{quote}

\sphinxAtStartPar
Set information value of the PSD constraint. Support related PSD informations.

\sphinxAtStartPar
\sphinxstylestrong{Synopsis}
\begin{quote}

\sphinxAtStartPar
\sphinxcode{\sphinxupquote{void Set(string info, double value)}}
\end{quote}

\sphinxAtStartPar
\sphinxstylestrong{Arguments}
\begin{quote}

\sphinxAtStartPar
\sphinxcode{\sphinxupquote{info}}: name of queried information.

\sphinxAtStartPar
\sphinxcode{\sphinxupquote{value}}: new information value.
\end{quote}
\end{quote}

\subsubsection{PsdConstraint.SetName()}
\label{\detokenize{csapi/psdconstraint:psdconstraint-setname}}\begin{quote}

\sphinxAtStartPar
Set name of a PSD constraint.

\sphinxAtStartPar
\sphinxstylestrong{Synopsis}
\begin{quote}

\sphinxAtStartPar
\sphinxcode{\sphinxupquote{void SetName(string name)}}
\end{quote}

\sphinxAtStartPar
\sphinxstylestrong{Arguments}
\begin{quote}

\sphinxAtStartPar
\sphinxcode{\sphinxupquote{name}}: the name to set.
\end{quote}
\end{quote}

\subsection{PsdConstrArray}
\label{\detokenize{csharpapiref:psdconstrarray}}\label{\detokenize{csharpapiref:chapcsharpapiref-psdconstrarray}}
\sphinxAtStartPar
COPT PSD constraint array object. To store and access a set of {\hyperref[\detokenize{csharpapiref:chapcsharpapiref-psdconstraint}]{\sphinxcrossref{\DUrole{std,std-ref}{PsdConstraint}}}}
objects, Cardinal Optimizer provides PsdConstrArray class, which defines the following methods.

\sphinxstepscope

\subsubsection{PsdConstrArray.PsdConstrArray()}
\label{\detokenize{csapi/psdconstrarray:psdconstrarray-psdconstrarray}}\label{\detokenize{csapi/psdconstrarray::doc}}\begin{quote}

\sphinxAtStartPar
Constructor of PsdConstrArray object.

\sphinxAtStartPar
\sphinxstylestrong{Synopsis}
\begin{quote}

\sphinxAtStartPar
\sphinxcode{\sphinxupquote{PsdConstrArray()}}
\end{quote}
\end{quote}

\subsubsection{PsdConstrArray.GetPsdConstr()}
\label{\detokenize{csapi/psdconstrarray:psdconstrarray-getpsdconstr}}\begin{quote}

\sphinxAtStartPar
Get idx\sphinxhyphen{}th PSD constraint object.

\sphinxAtStartPar
\sphinxstylestrong{Synopsis}
\begin{quote}

\sphinxAtStartPar
\sphinxcode{\sphinxupquote{PsdConstraint GetPsdConstr(int idx)}}
\end{quote}

\sphinxAtStartPar
\sphinxstylestrong{Arguments}
\begin{quote}

\sphinxAtStartPar
\sphinxcode{\sphinxupquote{idx}}: index of the PSD constraint.
\end{quote}

\sphinxAtStartPar
\sphinxstylestrong{Return}
\begin{quote}

\sphinxAtStartPar
PSD constraint object with index idx.
\end{quote}
\end{quote}

\subsubsection{PsdConstrArray.PushBack()}
\label{\detokenize{csapi/psdconstrarray:psdconstrarray-pushback}}\begin{quote}

\sphinxAtStartPar
Add a PSD constraint object to PSD constraint array.

\sphinxAtStartPar
\sphinxstylestrong{Synopsis}
\begin{quote}

\sphinxAtStartPar
\sphinxcode{\sphinxupquote{void PushBack(PsdConstraint constr)}}
\end{quote}

\sphinxAtStartPar
\sphinxstylestrong{Arguments}
\begin{quote}

\sphinxAtStartPar
\sphinxcode{\sphinxupquote{constr}}: a PSD constraint object.
\end{quote}
\end{quote}

\subsubsection{PsdConstrArray.Reserve()}
\label{\detokenize{csapi/psdconstrarray:psdconstrarray-reserve}}\begin{quote}

\sphinxAtStartPar
Reserve capacity to contain at least n items.

\sphinxAtStartPar
\sphinxstylestrong{Synopsis}
\begin{quote}

\sphinxAtStartPar
\sphinxcode{\sphinxupquote{void Reserve(int n)}}
\end{quote}

\sphinxAtStartPar
\sphinxstylestrong{Arguments}
\begin{quote}

\sphinxAtStartPar
\sphinxcode{\sphinxupquote{n}}: minimum capacity for PSD constraint objects.
\end{quote}
\end{quote}

\subsubsection{PsdConstrArray.Size()}
\label{\detokenize{csapi/psdconstrarray:psdconstrarray-size}}\begin{quote}

\sphinxAtStartPar
Get the number of PSD constraint objects.

\sphinxAtStartPar
\sphinxstylestrong{Synopsis}
\begin{quote}

\sphinxAtStartPar
\sphinxcode{\sphinxupquote{int Size()}}
\end{quote}

\sphinxAtStartPar
\sphinxstylestrong{Return}
\begin{quote}

\sphinxAtStartPar
number of PSD constraint objects.
\end{quote}
\end{quote}

\subsection{PsdConstrBuilder}
\label{\detokenize{csharpapiref:psdconstrbuilder}}\label{\detokenize{csharpapiref:chapcsharpapiref-psdconstrbuilder}}
\sphinxAtStartPar
COPT PSD constraint builder object. To help building a PSD constraint, given a PSD
expression, constraint sense and right\sphinxhyphen{}hand side value, Cardinal Optimizer provides PsdConstrBuilder
class, which defines the following methods.

\sphinxstepscope

\subsubsection{PsdConstrBuilder.PsdConstrBuilder()}
\label{\detokenize{csapi/psdconstrbuilder:psdconstrbuilder-psdconstrbuilder}}\label{\detokenize{csapi/psdconstrbuilder::doc}}\begin{quote}

\sphinxAtStartPar
Constructor of PsdConstrBuilder object.

\sphinxAtStartPar
\sphinxstylestrong{Synopsis}
\begin{quote}

\sphinxAtStartPar
\sphinxcode{\sphinxupquote{PsdConstrBuilder()}}
\end{quote}
\end{quote}

\subsubsection{PsdConstrBuilder.GetPsdExpr()}
\label{\detokenize{csapi/psdconstrbuilder:psdconstrbuilder-getpsdexpr}}\begin{quote}

\sphinxAtStartPar
Get expression associated with PSD constraint.

\sphinxAtStartPar
\sphinxstylestrong{Synopsis}
\begin{quote}

\sphinxAtStartPar
\sphinxcode{\sphinxupquote{PsdExpr GetPsdExpr()}}
\end{quote}

\sphinxAtStartPar
\sphinxstylestrong{Return}
\begin{quote}

\sphinxAtStartPar
PSD expression object.
\end{quote}
\end{quote}

\subsubsection{PsdConstrBuilder.GetRange()}
\label{\detokenize{csapi/psdconstrbuilder:psdconstrbuilder-getrange}}\begin{quote}

\sphinxAtStartPar
Get range from lower bound to upper bound of range constraint.

\sphinxAtStartPar
\sphinxstylestrong{Synopsis}
\begin{quote}

\sphinxAtStartPar
\sphinxcode{\sphinxupquote{double GetRange()}}
\end{quote}

\sphinxAtStartPar
\sphinxstylestrong{Return}
\begin{quote}

\sphinxAtStartPar
length from lower bound to upper bound of the constraint.
\end{quote}
\end{quote}

\subsubsection{PsdConstrBuilder.GetSense()}
\label{\detokenize{csapi/psdconstrbuilder:psdconstrbuilder-getsense}}\begin{quote}

\sphinxAtStartPar
Get sense associated with PSD constraint.

\sphinxAtStartPar
\sphinxstylestrong{Synopsis}
\begin{quote}

\sphinxAtStartPar
\sphinxcode{\sphinxupquote{char GetSense()}}
\end{quote}

\sphinxAtStartPar
\sphinxstylestrong{Return}
\begin{quote}

\sphinxAtStartPar
PSD constraint sense.
\end{quote}
\end{quote}

\subsubsection{PsdConstrBuilder.Set()}
\label{\detokenize{csapi/psdconstrbuilder:psdconstrbuilder-set}}\begin{quote}

\sphinxAtStartPar
Set detail of a PSD constraint to its builder object.

\sphinxAtStartPar
\sphinxstylestrong{Synopsis}
\begin{quote}

\sphinxAtStartPar
\sphinxcode{\sphinxupquote{void Set(}}
\begin{quote}

\sphinxAtStartPar
\sphinxcode{\sphinxupquote{PsdExpr expr,}}

\sphinxAtStartPar
\sphinxcode{\sphinxupquote{char sense,}}

\sphinxAtStartPar
\sphinxcode{\sphinxupquote{double rhs)}}
\end{quote}
\end{quote}

\sphinxAtStartPar
\sphinxstylestrong{Arguments}
\begin{quote}

\sphinxAtStartPar
\sphinxcode{\sphinxupquote{expr}}: expression object at one side of the PSD constraint.

\sphinxAtStartPar
\sphinxcode{\sphinxupquote{sense}}: PSD constraint sense, other than COPT\_RANGE.

\sphinxAtStartPar
\sphinxcode{\sphinxupquote{rhs}}: constant at right side of the PSD constraint.
\end{quote}
\end{quote}

\subsubsection{PsdConstrBuilder.SetRange()}
\label{\detokenize{csapi/psdconstrbuilder:psdconstrbuilder-setrange}}\begin{quote}

\sphinxAtStartPar
Set a range constraint to its builder.

\sphinxAtStartPar
\sphinxstylestrong{Synopsis}
\begin{quote}

\sphinxAtStartPar
\sphinxcode{\sphinxupquote{void SetRange(PsdExpr expr, double range)}}
\end{quote}

\sphinxAtStartPar
\sphinxstylestrong{Arguments}
\begin{quote}

\sphinxAtStartPar
\sphinxcode{\sphinxupquote{expr}}: PSD expression object, whose constant is negative upper bound.

\sphinxAtStartPar
\sphinxcode{\sphinxupquote{range}}: length from lower bound to upper bound of the constraint. Must greater than 0.
\end{quote}
\end{quote}

\subsection{PsdConstrBuilderArray}
\label{\detokenize{csharpapiref:psdconstrbuilderarray}}\label{\detokenize{csharpapiref:chapcsharpapiref-psdconstrbuilderarray}}
\sphinxAtStartPar
COPT PSD constraint builder array object. To store and access a set of {\hyperref[\detokenize{csharpapiref:chapcsharpapiref-psdconstrbuilder}]{\sphinxcrossref{\DUrole{std,std-ref}{PsdConstrBuilder}}}}
objects, Cardinal Optimizer provides PsdConstrBuilderArray class, which defines the following methods.

\sphinxstepscope

\subsubsection{PsdConstrBuilderArray.PsdConstrBuilderArray()}
\label{\detokenize{csapi/psdconstrbuilderarray:psdconstrbuilderarray-psdconstrbuilderarray}}\label{\detokenize{csapi/psdconstrbuilderarray::doc}}\begin{quote}

\sphinxAtStartPar
Constructor of PsdConstrBuilderArray object.

\sphinxAtStartPar
\sphinxstylestrong{Synopsis}
\begin{quote}

\sphinxAtStartPar
\sphinxcode{\sphinxupquote{PsdConstrBuilderArray()}}
\end{quote}
\end{quote}

\subsubsection{PsdConstrBuilderArray.GetBuilder()}
\label{\detokenize{csapi/psdconstrbuilderarray:psdconstrbuilderarray-getbuilder}}\begin{quote}

\sphinxAtStartPar
Get idx\sphinxhyphen{}th PSD constraint builder object.

\sphinxAtStartPar
\sphinxstylestrong{Synopsis}
\begin{quote}

\sphinxAtStartPar
\sphinxcode{\sphinxupquote{PsdConstrBuilder GetBuilder(int idx)}}
\end{quote}

\sphinxAtStartPar
\sphinxstylestrong{Arguments}
\begin{quote}

\sphinxAtStartPar
\sphinxcode{\sphinxupquote{idx}}: index of the PSD constraint builder.
\end{quote}

\sphinxAtStartPar
\sphinxstylestrong{Return}
\begin{quote}

\sphinxAtStartPar
PSD constraint builder object with index idx.
\end{quote}
\end{quote}

\subsubsection{PsdConstrBuilderArray.PushBack()}
\label{\detokenize{csapi/psdconstrbuilderarray:psdconstrbuilderarray-pushback}}\begin{quote}

\sphinxAtStartPar
Add a PSD constraint builder object to PSD constraint builder array.

\sphinxAtStartPar
\sphinxstylestrong{Synopsis}
\begin{quote}

\sphinxAtStartPar
\sphinxcode{\sphinxupquote{void PushBack(PsdConstrBuilder builder)}}
\end{quote}

\sphinxAtStartPar
\sphinxstylestrong{Arguments}
\begin{quote}

\sphinxAtStartPar
\sphinxcode{\sphinxupquote{builder}}: a PSD constraint builder object.
\end{quote}
\end{quote}

\subsubsection{PsdConstrBuilderArray.Reserve()}
\label{\detokenize{csapi/psdconstrbuilderarray:psdconstrbuilderarray-reserve}}\begin{quote}

\sphinxAtStartPar
Reserve capacity to contain at least n items.

\sphinxAtStartPar
\sphinxstylestrong{Synopsis}
\begin{quote}

\sphinxAtStartPar
\sphinxcode{\sphinxupquote{void Reserve(int n)}}
\end{quote}

\sphinxAtStartPar
\sphinxstylestrong{Arguments}
\begin{quote}

\sphinxAtStartPar
\sphinxcode{\sphinxupquote{n}}: minimum capacity for PSD constraint builder object.
\end{quote}
\end{quote}

\subsubsection{PsdConstrBuilderArray.Size()}
\label{\detokenize{csapi/psdconstrbuilderarray:psdconstrbuilderarray-size}}\begin{quote}

\sphinxAtStartPar
Get the number of PSD constraint builder objects.

\sphinxAtStartPar
\sphinxstylestrong{Synopsis}
\begin{quote}

\sphinxAtStartPar
\sphinxcode{\sphinxupquote{int Size()}}
\end{quote}

\sphinxAtStartPar
\sphinxstylestrong{Return}
\begin{quote}

\sphinxAtStartPar
number of PSD constraint builder objects.
\end{quote}
\end{quote}

\subsection{LmiConstraint}
\label{\detokenize{csharpapiref:lmiconstraint}}\label{\detokenize{csharpapiref:chapcsharpapiref-lmiconstraint}}
\sphinxAtStartPar
COPT LMI constraint object. LMI constraints are always associated with a
particular model.  User creates a LMI constraint object by adding a LMI
constraint to model, rather than by constructor of LmiConstraint class.

\sphinxstepscope

\subsubsection{LmiConstraint.Get()}
\label{\detokenize{csapi/lmiconstraint:lmiconstraint-get}}\label{\detokenize{csapi/lmiconstraint::doc}}\begin{quote}

\sphinxAtStartPar
Get information values of LMI constraint.

\sphinxAtStartPar
\sphinxstylestrong{Synopsis}
\begin{quote}

\sphinxAtStartPar
\sphinxcode{\sphinxupquote{double{[}{]} Get(string info)}}
\end{quote}

\sphinxAtStartPar
\sphinxstylestrong{Arguments}
\begin{quote}

\sphinxAtStartPar
\sphinxcode{\sphinxupquote{info}}: name of information.
\end{quote}

\sphinxAtStartPar
\sphinxstylestrong{Return}
\begin{quote}

\sphinxAtStartPar
output array of information values.
\end{quote}
\end{quote}

\subsubsection{LmiConstraint.GetDim()}
\label{\detokenize{csapi/lmiconstraint:lmiconstraint-getdim}}\begin{quote}

\sphinxAtStartPar
Get dimension of LMI constraint.

\sphinxAtStartPar
\sphinxstylestrong{Synopsis}
\begin{quote}

\sphinxAtStartPar
\sphinxcode{\sphinxupquote{int GetDim()}}
\end{quote}

\sphinxAtStartPar
\sphinxstylestrong{Return}
\begin{quote}

\sphinxAtStartPar
dimension of LMI constraint.
\end{quote}
\end{quote}

\subsubsection{LmiConstraint.GetIdx()}
\label{\detokenize{csapi/lmiconstraint:lmiconstraint-getidx}}\begin{quote}

\sphinxAtStartPar
Get index of LMI constraint.

\sphinxAtStartPar
\sphinxstylestrong{Synopsis}
\begin{quote}

\sphinxAtStartPar
\sphinxcode{\sphinxupquote{int GetIdx()}}
\end{quote}

\sphinxAtStartPar
\sphinxstylestrong{Return}
\begin{quote}

\sphinxAtStartPar
index of LMI constraint.
\end{quote}
\end{quote}

\subsubsection{LmiConstraint.GetLen()}
\label{\detokenize{csapi/lmiconstraint:lmiconstraint-getlen}}\begin{quote}

\sphinxAtStartPar
Get length of LMI constraint.

\sphinxAtStartPar
\sphinxstylestrong{Synopsis}
\begin{quote}

\sphinxAtStartPar
\sphinxcode{\sphinxupquote{int GetLen()}}
\end{quote}

\sphinxAtStartPar
\sphinxstylestrong{Return}
\begin{quote}

\sphinxAtStartPar
length of LMI constraint.
\end{quote}
\end{quote}

\subsubsection{LmiConstraint.GetName()}
\label{\detokenize{csapi/lmiconstraint:lmiconstraint-getname}}\begin{quote}

\sphinxAtStartPar
Get name of LMI constraint.

\sphinxAtStartPar
\sphinxstylestrong{Synopsis}
\begin{quote}

\sphinxAtStartPar
\sphinxcode{\sphinxupquote{string GetName()}}
\end{quote}

\sphinxAtStartPar
\sphinxstylestrong{Return}
\begin{quote}

\sphinxAtStartPar
name of LMI constraint.
\end{quote}
\end{quote}

\subsubsection{LmiConstraint.Remove()}
\label{\detokenize{csapi/lmiconstraint:lmiconstraint-remove}}\begin{quote}

\sphinxAtStartPar
Remove this LMI constraint from model.

\sphinxAtStartPar
\sphinxstylestrong{Synopsis}
\begin{quote}

\sphinxAtStartPar
\sphinxcode{\sphinxupquote{void Remove()}}
\end{quote}
\end{quote}

\subsubsection{LmiConstraint.SetRhs()}
\label{\detokenize{csapi/lmiconstraint:lmiconstraint-setrhs}}\begin{quote}

\sphinxAtStartPar
Set constant term of LMI constraint.

\sphinxAtStartPar
\sphinxstylestrong{Synopsis}
\begin{quote}

\sphinxAtStartPar
\sphinxcode{\sphinxupquote{void SetRhs(SymMatrix mat)}}
\end{quote}

\sphinxAtStartPar
\sphinxstylestrong{Arguments}
\begin{quote}

\sphinxAtStartPar
\sphinxcode{\sphinxupquote{mat}}: new symmetric matrix for constant term.
\end{quote}
\end{quote}

\subsection{LmiConstrArray}
\label{\detokenize{csharpapiref:lmiconstrarray}}\label{\detokenize{csharpapiref:chapcsharpapiref-lmiconstrarray}}
\sphinxAtStartPar
COPT LMI constraint array object. To store and access a set of
{\hyperref[\detokenize{csharpapiref:chapcsharpapiref-lmiconstraint}]{\sphinxcrossref{\DUrole{std,std-ref}{LmiConstraint}}}} objects, Cardinal Optimizer provides
LmiConstrArray class, which defines the following methods.

\sphinxstepscope

\subsubsection{LmiConstrArray.LmiConstrArray()}
\label{\detokenize{csapi/lmiconstrarray:lmiconstrarray-lmiconstrarray}}\label{\detokenize{csapi/lmiconstrarray::doc}}\begin{quote}

\sphinxAtStartPar
Constructor of LmiConstrArray.

\sphinxAtStartPar
\sphinxstylestrong{Synopsis}
\begin{quote}

\sphinxAtStartPar
\sphinxcode{\sphinxupquote{LmiConstrArray()}}
\end{quote}
\end{quote}

\subsubsection{LmiConstrArray.GetLmiConstr()}
\label{\detokenize{csapi/lmiconstrarray:lmiconstrarray-getlmiconstr}}\begin{quote}

\sphinxAtStartPar
Get idx\sphinxhyphen{}th LMI constraint object.

\sphinxAtStartPar
\sphinxstylestrong{Synopsis}
\begin{quote}

\sphinxAtStartPar
\sphinxcode{\sphinxupquote{LmiConstraint GetLmiConstr(int idx)}}
\end{quote}

\sphinxAtStartPar
\sphinxstylestrong{Arguments}
\begin{quote}

\sphinxAtStartPar
\sphinxcode{\sphinxupquote{idx}}: index of LMI constraint.
\end{quote}

\sphinxAtStartPar
\sphinxstylestrong{Return}
\begin{quote}

\sphinxAtStartPar
LMI constraint object with index idx.
\end{quote}
\end{quote}

\subsubsection{LmiConstrArray.PushBack()}
\label{\detokenize{csapi/lmiconstrarray:lmiconstrarray-pushback}}\begin{quote}

\sphinxAtStartPar
Add an LMI constraint to LMI constraint array.

\sphinxAtStartPar
\sphinxstylestrong{Synopsis}
\begin{quote}

\sphinxAtStartPar
\sphinxcode{\sphinxupquote{void PushBack(LmiConstraint constr)}}
\end{quote}

\sphinxAtStartPar
\sphinxstylestrong{Arguments}
\begin{quote}

\sphinxAtStartPar
\sphinxcode{\sphinxupquote{constr}}: LMI constraint object.
\end{quote}
\end{quote}

\subsubsection{LmiConstrArray.Reserve()}
\label{\detokenize{csapi/lmiconstrarray:lmiconstrarray-reserve}}\begin{quote}

\sphinxAtStartPar
Reserve capacity to contain at least n items.

\sphinxAtStartPar
\sphinxstylestrong{Synopsis}
\begin{quote}

\sphinxAtStartPar
\sphinxcode{\sphinxupquote{void Reserve(int n)}}
\end{quote}

\sphinxAtStartPar
\sphinxstylestrong{Arguments}
\begin{quote}

\sphinxAtStartPar
\sphinxcode{\sphinxupquote{n}}: capacity number of LMI constraint object.
\end{quote}
\end{quote}

\subsubsection{LmiConstrArray.Size()}
\label{\detokenize{csapi/lmiconstrarray:lmiconstrarray-size}}\begin{quote}

\sphinxAtStartPar
Get the number of LMI constraint objects.

\sphinxAtStartPar
\sphinxstylestrong{Synopsis}
\begin{quote}

\sphinxAtStartPar
\sphinxcode{\sphinxupquote{int Size()}}
\end{quote}

\sphinxAtStartPar
\sphinxstylestrong{Return}
\begin{quote}

\sphinxAtStartPar
number of LMI constraint objects.
\end{quote}
\end{quote}

\subsection{LmiExpr}
\label{\detokenize{csharpapiref:lmiexpr}}\label{\detokenize{csharpapiref:chapcsharpapiref-lmiexpr}}
\sphinxAtStartPar
COPT LMI expression object. A LMI expression consists of
a list of variables, associated coefficient matrices of LMI term, and constant matrices.
LMI expressions are used to build LMI constraints.

\sphinxstepscope

\subsubsection{LmiExpr.LmiExpr()}
\label{\detokenize{csapi/lmiexpr:lmiexpr-lmiexpr}}\label{\detokenize{csapi/lmiexpr::doc}}\begin{quote}

\sphinxAtStartPar
Default constructor of LMI expression

\sphinxAtStartPar
\sphinxstylestrong{Synopsis}
\begin{quote}

\sphinxAtStartPar
\sphinxcode{\sphinxupquote{LmiExpr()}}
\end{quote}
\end{quote}

\subsubsection{LmiExpr.LmiExpr()}
\label{\detokenize{csapi/lmiexpr:id1}}\begin{quote}

\sphinxAtStartPar
Constructor of LMI expression with constant term.

\sphinxAtStartPar
\sphinxstylestrong{Synopsis}
\begin{quote}

\sphinxAtStartPar
\sphinxcode{\sphinxupquote{LmiExpr(SymMatrix mat)}}
\end{quote}

\sphinxAtStartPar
\sphinxstylestrong{Arguments}
\begin{quote}

\sphinxAtStartPar
\sphinxcode{\sphinxupquote{mat}}: symmetric matrix object.
\end{quote}
\end{quote}

\subsubsection{LmiExpr.LmiExpr()}
\label{\detokenize{csapi/lmiexpr:id2}}\begin{quote}

\sphinxAtStartPar
Constructor of LMI expression with matrix expression.

\sphinxAtStartPar
\sphinxstylestrong{Synopsis}
\begin{quote}

\sphinxAtStartPar
\sphinxcode{\sphinxupquote{LmiExpr(SymMatExpr expr)}}
\end{quote}

\sphinxAtStartPar
\sphinxstylestrong{Arguments}
\begin{quote}

\sphinxAtStartPar
\sphinxcode{\sphinxupquote{expr}}: symmetric matrix expression.
\end{quote}
\end{quote}

\subsubsection{LmiExpr.LmiExpr()}
\label{\detokenize{csapi/lmiexpr:id3}}\begin{quote}

\sphinxAtStartPar
Constructor of LMI expression with one term.

\sphinxAtStartPar
\sphinxstylestrong{Synopsis}
\begin{quote}

\sphinxAtStartPar
\sphinxcode{\sphinxupquote{LmiExpr(Var var, SymMatrix mat)}}
\end{quote}

\sphinxAtStartPar
\sphinxstylestrong{Arguments}
\begin{quote}

\sphinxAtStartPar
\sphinxcode{\sphinxupquote{var}}: variable of the added term.

\sphinxAtStartPar
\sphinxcode{\sphinxupquote{mat}}: coefficient matrix of the added term.
\end{quote}
\end{quote}

\subsubsection{LmiExpr.LmiExpr()}
\label{\detokenize{csapi/lmiexpr:id4}}\begin{quote}

\sphinxAtStartPar
Constructor of LMI expression with one term.

\sphinxAtStartPar
\sphinxstylestrong{Synopsis}
\begin{quote}

\sphinxAtStartPar
\sphinxcode{\sphinxupquote{LmiExpr(Var var, SymMatExpr expr)}}
\end{quote}

\sphinxAtStartPar
\sphinxstylestrong{Arguments}
\begin{quote}

\sphinxAtStartPar
\sphinxcode{\sphinxupquote{var}}: variable of the added term.

\sphinxAtStartPar
\sphinxcode{\sphinxupquote{expr}}: coefficient expression of symmetric matrices of new LMI term.
\end{quote}
\end{quote}

\subsubsection{LmiExpr.AddConstant()}
\label{\detokenize{csapi/lmiexpr:lmiexpr-addconstant}}\begin{quote}

\sphinxAtStartPar
Add to constant term of the LMI expression.

\sphinxAtStartPar
\sphinxstylestrong{Synopsis}
\begin{quote}

\sphinxAtStartPar
\sphinxcode{\sphinxupquote{void AddConstant(SymMatExpr expr)}}
\end{quote}

\sphinxAtStartPar
\sphinxstylestrong{Arguments}
\begin{quote}

\sphinxAtStartPar
\sphinxcode{\sphinxupquote{expr}}: matrix expression added to the constant term.
\end{quote}
\end{quote}

\subsubsection{LmiExpr.AddLmiExpr()}
\label{\detokenize{csapi/lmiexpr:lmiexpr-addlmiexpr}}\begin{quote}

\sphinxAtStartPar
Add an LMI expression to self.

\sphinxAtStartPar
\sphinxstylestrong{Synopsis}
\begin{quote}

\sphinxAtStartPar
\sphinxcode{\sphinxupquote{void AddLmiExpr(LmiExpr expr, double mult)}}
\end{quote}

\sphinxAtStartPar
\sphinxstylestrong{Arguments}
\begin{quote}

\sphinxAtStartPar
\sphinxcode{\sphinxupquote{expr}}: LMI expression to be added.

\sphinxAtStartPar
\sphinxcode{\sphinxupquote{mult}}: optional, constant multiplier, default value is 1.0.
\end{quote}
\end{quote}

\subsubsection{LmiExpr.AddTerm()}
\label{\detokenize{csapi/lmiexpr:lmiexpr-addterm}}\begin{quote}

\sphinxAtStartPar
Add a term to LMI expression object.

\sphinxAtStartPar
\sphinxstylestrong{Synopsis}
\begin{quote}

\sphinxAtStartPar
\sphinxcode{\sphinxupquote{void AddTerm(Var var, SymMatrix mat)}}
\end{quote}

\sphinxAtStartPar
\sphinxstylestrong{Arguments}
\begin{quote}

\sphinxAtStartPar
\sphinxcode{\sphinxupquote{var}}: variable of new LMI term.

\sphinxAtStartPar
\sphinxcode{\sphinxupquote{mat}}: coefficient matrix of new LMI term.
\end{quote}
\end{quote}

\subsubsection{LmiExpr.AddTerm()}
\label{\detokenize{csapi/lmiexpr:id5}}\begin{quote}

\sphinxAtStartPar
Add a term to LMI expression object.

\sphinxAtStartPar
\sphinxstylestrong{Synopsis}
\begin{quote}

\sphinxAtStartPar
\sphinxcode{\sphinxupquote{void AddTerm(Var var, SymMatExpr expr)}}
\end{quote}

\sphinxAtStartPar
\sphinxstylestrong{Arguments}
\begin{quote}

\sphinxAtStartPar
\sphinxcode{\sphinxupquote{var}}: variable of new LMI term.

\sphinxAtStartPar
\sphinxcode{\sphinxupquote{expr}}: coefficient expression of symmetric matrices of new LMI term.
\end{quote}
\end{quote}

\subsubsection{LmiExpr.AddTerms()}
\label{\detokenize{csapi/lmiexpr:lmiexpr-addterms}}\begin{quote}

\sphinxAtStartPar
Add LMI terms to LMI expression object.

\sphinxAtStartPar
\sphinxstylestrong{Synopsis}
\begin{quote}

\sphinxAtStartPar
\sphinxcode{\sphinxupquote{void AddTerms(VarArray vars, SymMatrixArray mats)}}
\end{quote}

\sphinxAtStartPar
\sphinxstylestrong{Arguments}
\begin{quote}

\sphinxAtStartPar
\sphinxcode{\sphinxupquote{vars}}: variables for added LMI terms.

\sphinxAtStartPar
\sphinxcode{\sphinxupquote{mats}}: coefficient matrices for added LMI terms.
\end{quote}
\end{quote}

\subsubsection{LmiExpr.AddTerms()}
\label{\detokenize{csapi/lmiexpr:id6}}\begin{quote}

\sphinxAtStartPar
Add LMI terms to LMI expression object.

\sphinxAtStartPar
\sphinxstylestrong{Synopsis}
\begin{quote}

\sphinxAtStartPar
\sphinxcode{\sphinxupquote{void AddTerms(Var{[}{]} vars, SymMatrix{[}{]} mats)}}
\end{quote}

\sphinxAtStartPar
\sphinxstylestrong{Arguments}
\begin{quote}

\sphinxAtStartPar
\sphinxcode{\sphinxupquote{vars}}: variables for added LMI terms.

\sphinxAtStartPar
\sphinxcode{\sphinxupquote{mats}}: coefficient matrices for added LMI terms.
\end{quote}
\end{quote}

\subsubsection{LmiExpr.Clone()}
\label{\detokenize{csapi/lmiexpr:lmiexpr-clone}}\begin{quote}

\sphinxAtStartPar
Deep copy LMI expression.

\sphinxAtStartPar
\sphinxstylestrong{Synopsis}
\begin{quote}

\sphinxAtStartPar
\sphinxcode{\sphinxupquote{LmiExpr Clone()}}
\end{quote}

\sphinxAtStartPar
\sphinxstylestrong{Return}
\begin{quote}

\sphinxAtStartPar
cloned LMI expression object.
\end{quote}
\end{quote}

\subsubsection{LmiExpr.GetCoeff()}
\label{\detokenize{csapi/lmiexpr:lmiexpr-getcoeff}}\begin{quote}

\sphinxAtStartPar
Get coefficient from the i\sphinxhyphen{}th term in LMI expression.

\sphinxAtStartPar
\sphinxstylestrong{Synopsis}
\begin{quote}

\sphinxAtStartPar
\sphinxcode{\sphinxupquote{SymMatExpr GetCoeff(int i)}}
\end{quote}

\sphinxAtStartPar
\sphinxstylestrong{Arguments}
\begin{quote}

\sphinxAtStartPar
\sphinxcode{\sphinxupquote{i}}: index of the LMI term.
\end{quote}

\sphinxAtStartPar
\sphinxstylestrong{Return}
\begin{quote}

\sphinxAtStartPar
coefficient expression of the i\sphinxhyphen{}th LMI term.
\end{quote}
\end{quote}

\subsubsection{LmiExpr.GetConstant()}
\label{\detokenize{csapi/lmiexpr:lmiexpr-getconstant}}\begin{quote}

\sphinxAtStartPar
Get constant term in LMI expression.

\sphinxAtStartPar
\sphinxstylestrong{Synopsis}
\begin{quote}

\sphinxAtStartPar
\sphinxcode{\sphinxupquote{SymMatExpr GetConstant()}}
\end{quote}

\sphinxAtStartPar
\sphinxstylestrong{Return}
\begin{quote}

\sphinxAtStartPar
symmetric matrix expression object.
\end{quote}
\end{quote}

\subsubsection{LmiExpr.GetVar()}
\label{\detokenize{csapi/lmiexpr:lmiexpr-getvar}}\begin{quote}

\sphinxAtStartPar
Get variable from the i\sphinxhyphen{}th term in LMI expression.

\sphinxAtStartPar
\sphinxstylestrong{Synopsis}
\begin{quote}

\sphinxAtStartPar
\sphinxcode{\sphinxupquote{Var GetVar(int i)}}
\end{quote}

\sphinxAtStartPar
\sphinxstylestrong{Arguments}
\begin{quote}

\sphinxAtStartPar
\sphinxcode{\sphinxupquote{i}}: index of the term.
\end{quote}

\sphinxAtStartPar
\sphinxstylestrong{Return}
\begin{quote}

\sphinxAtStartPar
variable of the i\sphinxhyphen{}th term in LMI expression object.
\end{quote}
\end{quote}

\subsubsection{LmiExpr.Multiply()}
\label{\detokenize{csapi/lmiexpr:lmiexpr-multiply}}\begin{quote}

\sphinxAtStartPar
Multiply itself by a constant.

\sphinxAtStartPar
\sphinxstylestrong{Synopsis}
\begin{quote}

\sphinxAtStartPar
\sphinxcode{\sphinxupquote{void Multiply(double c)}}
\end{quote}

\sphinxAtStartPar
\sphinxstylestrong{Arguments}
\begin{quote}

\sphinxAtStartPar
\sphinxcode{\sphinxupquote{c}}: constant operand.
\end{quote}
\end{quote}

\subsubsection{LmiExpr.Remove()}
\label{\detokenize{csapi/lmiexpr:lmiexpr-remove}}\begin{quote}

\sphinxAtStartPar
Remove i\sphinxhyphen{}th term from LMI expression object.

\sphinxAtStartPar
\sphinxstylestrong{Synopsis}
\begin{quote}

\sphinxAtStartPar
\sphinxcode{\sphinxupquote{void Remove(int idx)}}
\end{quote}

\sphinxAtStartPar
\sphinxstylestrong{Arguments}
\begin{quote}

\sphinxAtStartPar
\sphinxcode{\sphinxupquote{idx}}: index of the term to be removed.
\end{quote}
\end{quote}

\subsubsection{LmiExpr.Remove()}
\label{\detokenize{csapi/lmiexpr:id7}}\begin{quote}

\sphinxAtStartPar
Remove the term associated with variable from LMI expression.

\sphinxAtStartPar
\sphinxstylestrong{Synopsis}
\begin{quote}

\sphinxAtStartPar
\sphinxcode{\sphinxupquote{void Remove(Var var)}}
\end{quote}

\sphinxAtStartPar
\sphinxstylestrong{Arguments}
\begin{quote}

\sphinxAtStartPar
\sphinxcode{\sphinxupquote{var}}: a variable whose term should be removed.
\end{quote}
\end{quote}

\subsubsection{LmiExpr.SetCoeff()}
\label{\detokenize{csapi/lmiexpr:lmiexpr-setcoeff}}\begin{quote}

\sphinxAtStartPar
Set coefficient matrix of the i\sphinxhyphen{}th term in LMI expression.

\sphinxAtStartPar
\sphinxstylestrong{Synopsis}
\begin{quote}

\sphinxAtStartPar
\sphinxcode{\sphinxupquote{void SetCoeff(int i, SymMatrix mat)}}
\end{quote}

\sphinxAtStartPar
\sphinxstylestrong{Arguments}
\begin{quote}

\sphinxAtStartPar
\sphinxcode{\sphinxupquote{i}}: index of the LMI term.

\sphinxAtStartPar
\sphinxcode{\sphinxupquote{mat}}: coefficient matrix of the term.
\end{quote}
\end{quote}

\subsubsection{LmiExpr.SetConstant()}
\label{\detokenize{csapi/lmiexpr:lmiexpr-setconstant}}\begin{quote}

\sphinxAtStartPar
Set constant term of the LMI expression.

\sphinxAtStartPar
\sphinxstylestrong{Synopsis}
\begin{quote}

\sphinxAtStartPar
\sphinxcode{\sphinxupquote{void SetConstant(SymMatrix mat)}}
\end{quote}

\sphinxAtStartPar
\sphinxstylestrong{Arguments}
\begin{quote}

\sphinxAtStartPar
\sphinxcode{\sphinxupquote{mat}}: symmetric matrix of the constant term.
\end{quote}
\end{quote}

\subsubsection{LmiExpr.Size()}
\label{\detokenize{csapi/lmiexpr:lmiexpr-size}}\begin{quote}

\sphinxAtStartPar
Get number of LMI terms in expression.

\sphinxAtStartPar
\sphinxstylestrong{Synopsis}
\begin{quote}

\sphinxAtStartPar
\sphinxcode{\sphinxupquote{long Size()}}
\end{quote}

\sphinxAtStartPar
\sphinxstylestrong{Return}
\begin{quote}

\sphinxAtStartPar
number of LMI terms.
\end{quote}
\end{quote}

\subsection{SymMatrix}
\label{\detokenize{csharpapiref:symmatrix}}\label{\detokenize{csharpapiref:chapcsharpapiref-symmatrix}}
\sphinxAtStartPar
COPT symmetric matrix object. Symmetric matrices are always associated with a particular model.
User creates a symmetric matrix object by adding a symmetric matrix to model,
rather than by constructor of SymMatrix class.

\sphinxAtStartPar
Symmetric matrices are used as coefficient matrices of PSD terms in PSD expressions, PSD constraints or PSD objectives.

\sphinxstepscope

\subsubsection{SymMatrix.GetDim()}
\label{\detokenize{csapi/symmatrix:symmatrix-getdim}}\label{\detokenize{csapi/symmatrix::doc}}\begin{quote}

\sphinxAtStartPar
Get the dimension of a symmetric matrix.

\sphinxAtStartPar
\sphinxstylestrong{Synopsis}
\begin{quote}

\sphinxAtStartPar
\sphinxcode{\sphinxupquote{int GetDim()}}
\end{quote}

\sphinxAtStartPar
\sphinxstylestrong{Return}
\begin{quote}

\sphinxAtStartPar
dimension of a symmetric matrix.
\end{quote}
\end{quote}

\subsubsection{SymMatrix.GetIdx()}
\label{\detokenize{csapi/symmatrix:symmatrix-getidx}}\begin{quote}

\sphinxAtStartPar
Get the index of a symmetric matrix.

\sphinxAtStartPar
\sphinxstylestrong{Synopsis}
\begin{quote}

\sphinxAtStartPar
\sphinxcode{\sphinxupquote{int GetIdx()}}
\end{quote}

\sphinxAtStartPar
\sphinxstylestrong{Return}
\begin{quote}

\sphinxAtStartPar
index of a symmetric matrix.
\end{quote}
\end{quote}

\subsection{SymMatrixArray}
\label{\detokenize{csharpapiref:symmatrixarray}}\label{\detokenize{csharpapiref:chapcsharpapiref-symmatrixarray}}
\sphinxAtStartPar
COPT symmetric matrix object. To store and access a set of {\hyperref[\detokenize{csharpapiref:chapcsharpapiref-symmatrix}]{\sphinxcrossref{\DUrole{std,std-ref}{SymMatrix}}}}
objects, Cardinal Optimizer provides SymMatrixArray class, which defines the following methods.

\sphinxstepscope

\subsubsection{SymMatrixArray.SymMatrixArray()}
\label{\detokenize{csapi/symmatrixarray:symmatrixarray-symmatrixarray}}\label{\detokenize{csapi/symmatrixarray::doc}}\begin{quote}

\sphinxAtStartPar
Constructor of SymMatrixAarray.

\sphinxAtStartPar
\sphinxstylestrong{Synopsis}
\begin{quote}

\sphinxAtStartPar
\sphinxcode{\sphinxupquote{SymMatrixArray()}}
\end{quote}
\end{quote}

\subsubsection{SymMatrixArray.GetMatrix()}
\label{\detokenize{csapi/symmatrixarray:symmatrixarray-getmatrix}}\begin{quote}

\sphinxAtStartPar
Get i\sphinxhyphen{}th SymMatrix object.

\sphinxAtStartPar
\sphinxstylestrong{Synopsis}
\begin{quote}

\sphinxAtStartPar
\sphinxcode{\sphinxupquote{SymMatrix GetMatrix(int idx)}}
\end{quote}

\sphinxAtStartPar
\sphinxstylestrong{Arguments}
\begin{quote}

\sphinxAtStartPar
\sphinxcode{\sphinxupquote{idx}}: index of the SymMatrix object.
\end{quote}

\sphinxAtStartPar
\sphinxstylestrong{Return}
\begin{quote}

\sphinxAtStartPar
SymMatrix object with index idx.
\end{quote}
\end{quote}

\subsubsection{SymMatrixArray.PushBack()}
\label{\detokenize{csapi/symmatrixarray:symmatrixarray-pushback}}\begin{quote}

\sphinxAtStartPar
Add a SymMatrix object to SymMatrix array.

\sphinxAtStartPar
\sphinxstylestrong{Synopsis}
\begin{quote}

\sphinxAtStartPar
\sphinxcode{\sphinxupquote{void PushBack(SymMatrix mat)}}
\end{quote}

\sphinxAtStartPar
\sphinxstylestrong{Arguments}
\begin{quote}

\sphinxAtStartPar
\sphinxcode{\sphinxupquote{mat}}: a SymMatrix object.
\end{quote}
\end{quote}

\subsubsection{SymMatrixArray.Reserve()}
\label{\detokenize{csapi/symmatrixarray:symmatrixarray-reserve}}\begin{quote}

\sphinxAtStartPar
Reserve capacity to contain at least n items.

\sphinxAtStartPar
\sphinxstylestrong{Synopsis}
\begin{quote}

\sphinxAtStartPar
\sphinxcode{\sphinxupquote{void Reserve(int n)}}
\end{quote}

\sphinxAtStartPar
\sphinxstylestrong{Arguments}
\begin{quote}

\sphinxAtStartPar
\sphinxcode{\sphinxupquote{n}}: minimum capacity for symmetric matrix object.
\end{quote}
\end{quote}

\subsubsection{SymMatrixArray.Size()}
\label{\detokenize{csapi/symmatrixarray:symmatrixarray-size}}\begin{quote}

\sphinxAtStartPar
Get the number of SymMatrix objects.

\sphinxAtStartPar
\sphinxstylestrong{Synopsis}
\begin{quote}

\sphinxAtStartPar
\sphinxcode{\sphinxupquote{int Size()}}
\end{quote}

\sphinxAtStartPar
\sphinxstylestrong{Return}
\begin{quote}

\sphinxAtStartPar
number of SymMatrix objects.
\end{quote}
\end{quote}

\subsection{SymMatExpr}
\label{\detokenize{csharpapiref:symmatexpr}}\label{\detokenize{csharpapiref:chapcsharpapiref-symmatexpr}}
\sphinxAtStartPar
COPT symmetric matrix expression object. A symmetric matrix expression is a
linear combination of symmetric matrices, which is still a symmetric matrix.
However, by doing so, we are able to delay computing the final matrix
until setting PSD constraints or PSD objective.

\sphinxstepscope

\subsubsection{SymMatExpr.SymMatExpr()}
\label{\detokenize{csapi/symmatexpr:symmatexpr-symmatexpr}}\label{\detokenize{csapi/symmatexpr::doc}}\begin{quote}

\sphinxAtStartPar
Constructor of a symmetric matrix expression.

\sphinxAtStartPar
\sphinxstylestrong{Synopsis}
\begin{quote}

\sphinxAtStartPar
\sphinxcode{\sphinxupquote{SymMatExpr()}}
\end{quote}
\end{quote}

\subsubsection{SymMatExpr.SymMatExpr()}
\label{\detokenize{csapi/symmatexpr:id1}}\begin{quote}

\sphinxAtStartPar
Constructor of a symmetric matrix expression with one term.

\sphinxAtStartPar
\sphinxstylestrong{Synopsis}
\begin{quote}

\sphinxAtStartPar
\sphinxcode{\sphinxupquote{SymMatExpr(SymMatrix mat, double coeff)}}
\end{quote}

\sphinxAtStartPar
\sphinxstylestrong{Arguments}
\begin{quote}

\sphinxAtStartPar
\sphinxcode{\sphinxupquote{mat}}: symmetric matrix of the added term.

\sphinxAtStartPar
\sphinxcode{\sphinxupquote{coeff}}: optional, coefficent for the added term. Its default value is 1.0.
\end{quote}
\end{quote}

\subsubsection{SymMatExpr.AddSymMatExpr()}
\label{\detokenize{csapi/symmatexpr:symmatexpr-addsymmatexpr}}\begin{quote}

\sphinxAtStartPar
Add a symmetric matrix expression to self.

\sphinxAtStartPar
\sphinxstylestrong{Synopsis}
\begin{quote}

\sphinxAtStartPar
\sphinxcode{\sphinxupquote{void AddSymMatExpr(SymMatExpr expr, double mult)}}
\end{quote}

\sphinxAtStartPar
\sphinxstylestrong{Arguments}
\begin{quote}

\sphinxAtStartPar
\sphinxcode{\sphinxupquote{expr}}: symmetric matrix expression to be added.

\sphinxAtStartPar
\sphinxcode{\sphinxupquote{mult}}: optional, constant multiplier, default value is 1.0.
\end{quote}
\end{quote}

\subsubsection{SymMatExpr.AddTerm()}
\label{\detokenize{csapi/symmatexpr:symmatexpr-addterm}}\begin{quote}

\sphinxAtStartPar
Add a term to symmetric matrix expression object.

\sphinxAtStartPar
\sphinxstylestrong{Synopsis}
\begin{quote}

\sphinxAtStartPar
\sphinxcode{\sphinxupquote{bool AddTerm(SymMatrix mat, double coeff)}}
\end{quote}

\sphinxAtStartPar
\sphinxstylestrong{Arguments}
\begin{quote}

\sphinxAtStartPar
\sphinxcode{\sphinxupquote{mat}}: symmetric matrix of the new term.

\sphinxAtStartPar
\sphinxcode{\sphinxupquote{coeff}}: coefficient of the new term.
\end{quote}

\sphinxAtStartPar
\sphinxstylestrong{Return}
\begin{quote}

\sphinxAtStartPar
True if the term is added successfully.
\end{quote}
\end{quote}

\subsubsection{SymMatExpr.AddTerms()}
\label{\detokenize{csapi/symmatexpr:symmatexpr-addterms}}\begin{quote}

\sphinxAtStartPar
Add multiple terms to expression object.

\sphinxAtStartPar
\sphinxstylestrong{Synopsis}
\begin{quote}

\sphinxAtStartPar
\sphinxcode{\sphinxupquote{int AddTerms(SymMatrix{[}{]} mats, double coeff)}}
\end{quote}

\sphinxAtStartPar
\sphinxstylestrong{Arguments}
\begin{quote}

\sphinxAtStartPar
\sphinxcode{\sphinxupquote{mats}}: symmetric matrix array object for added terms.

\sphinxAtStartPar
\sphinxcode{\sphinxupquote{coeff}}: optional, common coefficient for added terms, default is 1.0.
\end{quote}

\sphinxAtStartPar
\sphinxstylestrong{Return}
\begin{quote}

\sphinxAtStartPar
Number of added terms. If negative, fail to add one of terms.
\end{quote}
\end{quote}

\subsubsection{SymMatExpr.AddTerms()}
\label{\detokenize{csapi/symmatexpr:id2}}\begin{quote}

\sphinxAtStartPar
Add multiple terms to expression object.

\sphinxAtStartPar
\sphinxstylestrong{Synopsis}
\begin{quote}

\sphinxAtStartPar
\sphinxcode{\sphinxupquote{int AddTerms(SymMatrixArray mats, double{[}{]} coeffs)}}
\end{quote}

\sphinxAtStartPar
\sphinxstylestrong{Arguments}
\begin{quote}

\sphinxAtStartPar
\sphinxcode{\sphinxupquote{mats}}: symmetric matrix array object for added terms.

\sphinxAtStartPar
\sphinxcode{\sphinxupquote{coeffs}}: coefficient array for added terms.
\end{quote}

\sphinxAtStartPar
\sphinxstylestrong{Return}
\begin{quote}

\sphinxAtStartPar
Number of added terms. If negative, fail to add one of terms.
\end{quote}
\end{quote}

\subsubsection{SymMatExpr.AddTerms()}
\label{\detokenize{csapi/symmatexpr:id3}}\begin{quote}

\sphinxAtStartPar
Add multiple terms to expression object.

\sphinxAtStartPar
\sphinxstylestrong{Synopsis}
\begin{quote}

\sphinxAtStartPar
\sphinxcode{\sphinxupquote{int AddTerms(SymMatrix{[}{]} mats, double{[}{]} coeffs)}}
\end{quote}

\sphinxAtStartPar
\sphinxstylestrong{Arguments}
\begin{quote}

\sphinxAtStartPar
\sphinxcode{\sphinxupquote{mats}}: symmetric matrix array object for added terms.

\sphinxAtStartPar
\sphinxcode{\sphinxupquote{coeffs}}: coefficient array for added terms.
\end{quote}

\sphinxAtStartPar
\sphinxstylestrong{Return}
\begin{quote}

\sphinxAtStartPar
Number of added terms. If negative, fail to add one of terms.
\end{quote}
\end{quote}

\subsubsection{SymMatExpr.Clone()}
\label{\detokenize{csapi/symmatexpr:symmatexpr-clone}}\begin{quote}

\sphinxAtStartPar
Deep copy symmetric matrix expression object.

\sphinxAtStartPar
\sphinxstylestrong{Synopsis}
\begin{quote}

\sphinxAtStartPar
\sphinxcode{\sphinxupquote{SymMatExpr Clone()}}
\end{quote}

\sphinxAtStartPar
\sphinxstylestrong{Return}
\begin{quote}

\sphinxAtStartPar
cloned expression object.
\end{quote}
\end{quote}

\subsubsection{SymMatExpr.GetCoeff()}
\label{\detokenize{csapi/symmatexpr:symmatexpr-getcoeff}}\begin{quote}

\sphinxAtStartPar
Get coefficient of the i\sphinxhyphen{}th term in expression object.

\sphinxAtStartPar
\sphinxstylestrong{Synopsis}
\begin{quote}

\sphinxAtStartPar
\sphinxcode{\sphinxupquote{double GetCoeff(int i)}}
\end{quote}

\sphinxAtStartPar
\sphinxstylestrong{Arguments}
\begin{quote}

\sphinxAtStartPar
\sphinxcode{\sphinxupquote{i}}: index of the term.
\end{quote}

\sphinxAtStartPar
\sphinxstylestrong{Return}
\begin{quote}

\sphinxAtStartPar
coefficient of the i\sphinxhyphen{}th term.
\end{quote}
\end{quote}

\subsubsection{SymMatExpr.GetDim()}
\label{\detokenize{csapi/symmatexpr:symmatexpr-getdim}}\begin{quote}

\sphinxAtStartPar
Get dimension of symmetric matrix in expression.

\sphinxAtStartPar
\sphinxstylestrong{Synopsis}
\begin{quote}

\sphinxAtStartPar
\sphinxcode{\sphinxupquote{int GetDim()}}
\end{quote}

\sphinxAtStartPar
\sphinxstylestrong{Return}
\begin{quote}

\sphinxAtStartPar
dimension of symmetric matrix.
\end{quote}
\end{quote}

\subsubsection{SymMatExpr.GetSymMat()}
\label{\detokenize{csapi/symmatexpr:symmatexpr-getsymmat}}\begin{quote}

\sphinxAtStartPar
Get symmetric matrix of the i\sphinxhyphen{}th term in expression object.

\sphinxAtStartPar
\sphinxstylestrong{Synopsis}
\begin{quote}

\sphinxAtStartPar
\sphinxcode{\sphinxupquote{SymMatrix GetSymMat(int i)}}
\end{quote}

\sphinxAtStartPar
\sphinxstylestrong{Arguments}
\begin{quote}

\sphinxAtStartPar
\sphinxcode{\sphinxupquote{i}}: index of the term.
\end{quote}

\sphinxAtStartPar
\sphinxstylestrong{Return}
\begin{quote}

\sphinxAtStartPar
the symmetric matrix of the i\sphinxhyphen{}th term.
\end{quote}
\end{quote}

\subsubsection{SymMatExpr.Multiply()}
\label{\detokenize{csapi/symmatexpr:symmatexpr-multiply}}\begin{quote}

\sphinxAtStartPar
Multiply itself by a constant.

\sphinxAtStartPar
\sphinxstylestrong{Synopsis}
\begin{quote}

\sphinxAtStartPar
\sphinxcode{\sphinxupquote{void Multiply(double c)}}
\end{quote}

\sphinxAtStartPar
\sphinxstylestrong{Arguments}
\begin{quote}

\sphinxAtStartPar
\sphinxcode{\sphinxupquote{c}}: constant operand.
\end{quote}
\end{quote}

\subsubsection{SymMatExpr.Remove()}
\label{\detokenize{csapi/symmatexpr:symmatexpr-remove}}\begin{quote}

\sphinxAtStartPar
Remove i\sphinxhyphen{}th term from expression object.

\sphinxAtStartPar
\sphinxstylestrong{Synopsis}
\begin{quote}

\sphinxAtStartPar
\sphinxcode{\sphinxupquote{void Remove(int idx)}}
\end{quote}

\sphinxAtStartPar
\sphinxstylestrong{Arguments}
\begin{quote}

\sphinxAtStartPar
\sphinxcode{\sphinxupquote{idx}}: index of the term to be removed.
\end{quote}
\end{quote}

\subsubsection{SymMatExpr.Remove()}
\label{\detokenize{csapi/symmatexpr:id4}}\begin{quote}

\sphinxAtStartPar
Remove the term associated with the symmetric matrix.

\sphinxAtStartPar
\sphinxstylestrong{Synopsis}
\begin{quote}

\sphinxAtStartPar
\sphinxcode{\sphinxupquote{void Remove(SymMatrix mat)}}
\end{quote}

\sphinxAtStartPar
\sphinxstylestrong{Arguments}
\begin{quote}

\sphinxAtStartPar
\sphinxcode{\sphinxupquote{mat}}: a symmetric matrix whose term should be removed.
\end{quote}
\end{quote}

\subsubsection{SymMatExpr.Reserve()}
\label{\detokenize{csapi/symmatexpr:symmatexpr-reserve}}\begin{quote}

\sphinxAtStartPar
Reserve capacity to contain at least n items.

\sphinxAtStartPar
\sphinxstylestrong{Synopsis}
\begin{quote}

\sphinxAtStartPar
\sphinxcode{\sphinxupquote{void Reserve(int n)}}
\end{quote}

\sphinxAtStartPar
\sphinxstylestrong{Arguments}
\begin{quote}

\sphinxAtStartPar
\sphinxcode{\sphinxupquote{n}}: minimum capacity for expression object.
\end{quote}
\end{quote}

\subsubsection{SymMatExpr.SetCoeff()}
\label{\detokenize{csapi/symmatexpr:symmatexpr-setcoeff}}\begin{quote}

\sphinxAtStartPar
Set coefficient for the i\sphinxhyphen{}th term in expression object.

\sphinxAtStartPar
\sphinxstylestrong{Synopsis}
\begin{quote}

\sphinxAtStartPar
\sphinxcode{\sphinxupquote{void SetCoeff(int i, double val)}}
\end{quote}

\sphinxAtStartPar
\sphinxstylestrong{Arguments}
\begin{quote}

\sphinxAtStartPar
\sphinxcode{\sphinxupquote{i}}: index of the term.

\sphinxAtStartPar
\sphinxcode{\sphinxupquote{val}}: coefficient of the term.
\end{quote}
\end{quote}

\subsubsection{SymMatExpr.Size()}
\label{\detokenize{csapi/symmatexpr:symmatexpr-size}}\begin{quote}

\sphinxAtStartPar
Get number of terms in expression.

\sphinxAtStartPar
\sphinxstylestrong{Synopsis}
\begin{quote}

\sphinxAtStartPar
\sphinxcode{\sphinxupquote{long Size()}}
\end{quote}

\sphinxAtStartPar
\sphinxstylestrong{Return}
\begin{quote}

\sphinxAtStartPar
number of terms.
\end{quote}
\end{quote}

\subsection{MVar}
\label{\detokenize{csharpapiref:mvar}}\label{\detokenize{csharpapiref:chapcsharpapiref-mvar}}
\sphinxAtStartPar
COPT multi\sphinxhyphen{}dimensional variable object. It is used to construct
multi\sphinxhyphen{}dimensional variables and support operations with the built\sphinxhyphen{}in
multi\sphinxhyphen{}dimensional array {\hyperref[\detokenize{csharpapiref:chapcsharpapiref-ndarray}]{\sphinxcrossref{\DUrole{std,std-ref}{NdArray}}}} in COPT. It can be
created by calling the method \sphinxcode{\sphinxupquote{addMVar}} of {\hyperref[\detokenize{csharpapiref:chapcsharpapiref-model}]{\sphinxcrossref{\DUrole{std,std-ref}{Model}}}}.
The following methods are provided:

\sphinxstepscope

\subsubsection{MVar.Clone()}
\label{\detokenize{csapi/mvar:mvar-clone}}\label{\detokenize{csapi/mvar::doc}}\begin{quote}

\sphinxAtStartPar
Clone MVar object.

\sphinxAtStartPar
\sphinxstylestrong{Synopsis}
\begin{quote}

\sphinxAtStartPar
\sphinxcode{\sphinxupquote{MVar Clone()}}
\end{quote}

\sphinxAtStartPar
\sphinxstylestrong{Return}
\begin{quote}

\sphinxAtStartPar
new MVar object.
\end{quote}
\end{quote}

\subsubsection{MVar.Diagonal()}
\label{\detokenize{csapi/mvar:mvar-diagonal}}\begin{quote}

\sphinxAtStartPar
Get diagonals of MVar object.

\sphinxAtStartPar
\sphinxstylestrong{Synopsis}
\begin{quote}

\sphinxAtStartPar
\sphinxcode{\sphinxupquote{MVar Diagonal(}}
\begin{quote}

\sphinxAtStartPar
\sphinxcode{\sphinxupquote{int offset,}}

\sphinxAtStartPar
\sphinxcode{\sphinxupquote{int axis1,}}

\sphinxAtStartPar
\sphinxcode{\sphinxupquote{int axis2)}}
\end{quote}
\end{quote}

\sphinxAtStartPar
\sphinxstylestrong{Arguments}
\begin{quote}

\sphinxAtStartPar
\sphinxcode{\sphinxupquote{offset}}: offset of the diagonal from the main diagonal. Can be positive or negative.

\sphinxAtStartPar
\sphinxcode{\sphinxupquote{axis1}}: 1st axis of MVar.

\sphinxAtStartPar
\sphinxcode{\sphinxupquote{axis2}}: 2nd axis of MVar.
\end{quote}

\sphinxAtStartPar
\sphinxstylestrong{Return}
\begin{quote}

\sphinxAtStartPar
(N\sphinxhyphen{}1)\sphinxhyphen{}dimensional diagonals.
\end{quote}
\end{quote}

\subsubsection{MVar.Expand()}
\label{\detokenize{csapi/mvar:mvar-expand}}\begin{quote}

\sphinxAtStartPar
Expand shape of MVar object.

\sphinxAtStartPar
\sphinxstylestrong{Synopsis}
\begin{quote}

\sphinxAtStartPar
\sphinxcode{\sphinxupquote{MVar Expand(int axis)}}
\end{quote}

\sphinxAtStartPar
\sphinxstylestrong{Arguments}
\begin{quote}

\sphinxAtStartPar
\sphinxcode{\sphinxupquote{axis}}: axis of MVar.
\end{quote}

\sphinxAtStartPar
\sphinxstylestrong{Return}
\begin{quote}

\sphinxAtStartPar
MVar object of (N+1)\sphinxhyphen{}dimensional shape.
\end{quote}
\end{quote}

\subsubsection{MVar.Flatten()}
\label{\detokenize{csapi/mvar:mvar-flatten}}\begin{quote}

\sphinxAtStartPar
Flatten a MVar object to a 1\sphinxhyphen{}dimensional shape.

\sphinxAtStartPar
\sphinxstylestrong{Synopsis}
\begin{quote}

\sphinxAtStartPar
\sphinxcode{\sphinxupquote{MVar Flatten()}}
\end{quote}

\sphinxAtStartPar
\sphinxstylestrong{Return}
\begin{quote}

\sphinxAtStartPar
a MVar object collapsed into one dimension.
\end{quote}
\end{quote}

\subsubsection{MVar.Get()}
\label{\detokenize{csapi/mvar:mvar-get}}\begin{quote}

\sphinxAtStartPar
Get values of information associated with variables in MVar object.

\sphinxAtStartPar
\sphinxstylestrong{Synopsis}
\begin{quote}

\sphinxAtStartPar
\sphinxcode{\sphinxupquote{NdArray\textless{}double\textgreater{} Get(string info)}}
\end{quote}

\sphinxAtStartPar
\sphinxstylestrong{Arguments}
\begin{quote}

\sphinxAtStartPar
\sphinxcode{\sphinxupquote{info}}: name of information.
\end{quote}

\sphinxAtStartPar
\sphinxstylestrong{Return}
\begin{quote}

\sphinxAtStartPar
multi\sphinxhyphen{}dimensional array of information of variables.
\end{quote}
\end{quote}

\subsubsection{MVar.GetBasis()}
\label{\detokenize{csapi/mvar:mvar-getbasis}}\begin{quote}

\sphinxAtStartPar
Get basis of variables in MVar object.

\sphinxAtStartPar
\sphinxstylestrong{Synopsis}
\begin{quote}

\sphinxAtStartPar
\sphinxcode{\sphinxupquote{NdArray\textless{}int\textgreater{} GetBasis()}}
\end{quote}

\sphinxAtStartPar
\sphinxstylestrong{Return}
\begin{quote}

\sphinxAtStartPar
multi\sphinxhyphen{}dimensional array of basis of variables.
\end{quote}
\end{quote}

\subsubsection{MVar.GetDim()}
\label{\detokenize{csapi/mvar:mvar-getdim}}\begin{quote}

\sphinxAtStartPar
Get i\sphinxhyphen{}th dimension of MVar object.

\sphinxAtStartPar
\sphinxstylestrong{Synopsis}
\begin{quote}

\sphinxAtStartPar
\sphinxcode{\sphinxupquote{long GetDim(int i)}}
\end{quote}

\sphinxAtStartPar
\sphinxstylestrong{Arguments}
\begin{quote}

\sphinxAtStartPar
\sphinxcode{\sphinxupquote{i}}: index of dimension
\end{quote}

\sphinxAtStartPar
\sphinxstylestrong{Return}
\begin{quote}

\sphinxAtStartPar
i\sphinxhyphen{}th dimension.
\end{quote}
\end{quote}

\subsubsection{MVar.GetIdx()}
\label{\detokenize{csapi/mvar:mvar-getidx}}\begin{quote}

\sphinxAtStartPar
Get indexes of variables in MVar object.

\sphinxAtStartPar
\sphinxstylestrong{Synopsis}
\begin{quote}

\sphinxAtStartPar
\sphinxcode{\sphinxupquote{NdArray\textless{}int\textgreater{} GetIdx()}}
\end{quote}

\sphinxAtStartPar
\sphinxstylestrong{Return}
\begin{quote}

\sphinxAtStartPar
multi\sphinxhyphen{}dimensional array of indexes of variables.
\end{quote}
\end{quote}

\subsubsection{MVar.GetItem()}
\label{\detokenize{csapi/mvar:mvar-getitem}}\begin{quote}

\sphinxAtStartPar
Get variable of given index from MVar object.

\sphinxAtStartPar
\sphinxstylestrong{Synopsis}
\begin{quote}

\sphinxAtStartPar
\sphinxcode{\sphinxupquote{Var GetItem(long idx)}}
\end{quote}

\sphinxAtStartPar
\sphinxstylestrong{Arguments}
\begin{quote}

\sphinxAtStartPar
\sphinxcode{\sphinxupquote{idx}}: index of variable.
\end{quote}

\sphinxAtStartPar
\sphinxstylestrong{Return}
\begin{quote}

\sphinxAtStartPar
Var object.
\end{quote}
\end{quote}

\subsubsection{MVar.GetItem()}
\label{\detokenize{csapi/mvar:id1}}\begin{quote}

\sphinxAtStartPar
Get sub\sphinxhyphen{}arrays of MVar object, given view object.

\sphinxAtStartPar
\sphinxstylestrong{Synopsis}
\begin{quote}

\sphinxAtStartPar
\sphinxcode{\sphinxupquote{MVar GetItem(View view)}}
\end{quote}

\sphinxAtStartPar
\sphinxstylestrong{Arguments}
\begin{quote}

\sphinxAtStartPar
\sphinxcode{\sphinxupquote{view}}: view of multi\sphinxhyphen{}dimensional array.
\end{quote}

\sphinxAtStartPar
\sphinxstylestrong{Return}
\begin{quote}

\sphinxAtStartPar
sub\sphinxhyphen{}arrays of MVar object.
\end{quote}
\end{quote}

\subsubsection{MVar.GetLowerIIS()}
\label{\detokenize{csapi/mvar:mvar-getloweriis}}\begin{quote}

\sphinxAtStartPar
Get IIS status of lower bound of variables in MVar object.

\sphinxAtStartPar
\sphinxstylestrong{Synopsis}
\begin{quote}

\sphinxAtStartPar
\sphinxcode{\sphinxupquote{NdArray\textless{}int\textgreater{} GetLowerIIS()}}
\end{quote}

\sphinxAtStartPar
\sphinxstylestrong{Return}
\begin{quote}

\sphinxAtStartPar
multi\sphinxhyphen{}dimensional array of IIS status of lower bounds of variables.
\end{quote}
\end{quote}

\subsubsection{MVar.GetND()}
\label{\detokenize{csapi/mvar:mvar-getnd}}\begin{quote}

\sphinxAtStartPar
Get number of dimensions of MVar object.

\sphinxAtStartPar
\sphinxstylestrong{Synopsis}
\begin{quote}

\sphinxAtStartPar
\sphinxcode{\sphinxupquote{int GetND()}}
\end{quote}

\sphinxAtStartPar
\sphinxstylestrong{Return}
\begin{quote}

\sphinxAtStartPar
number of dimensions.
\end{quote}
\end{quote}

\subsubsection{MVar.GetShape()}
\label{\detokenize{csapi/mvar:mvar-getshape}}\begin{quote}

\sphinxAtStartPar
Get shape of MVar object.

\sphinxAtStartPar
\sphinxstylestrong{Synopsis}
\begin{quote}

\sphinxAtStartPar
\sphinxcode{\sphinxupquote{Shape GetShape()}}
\end{quote}

\sphinxAtStartPar
\sphinxstylestrong{Return}
\begin{quote}

\sphinxAtStartPar
shape object.
\end{quote}
\end{quote}

\subsubsection{MVar.GetSize()}
\label{\detokenize{csapi/mvar:mvar-getsize}}\begin{quote}

\sphinxAtStartPar
Get size of MVar object.

\sphinxAtStartPar
\sphinxstylestrong{Synopsis}
\begin{quote}

\sphinxAtStartPar
\sphinxcode{\sphinxupquote{long GetSize()}}
\end{quote}

\sphinxAtStartPar
\sphinxstylestrong{Return}
\begin{quote}

\sphinxAtStartPar
number of vars.
\end{quote}
\end{quote}

\subsubsection{MVar.GetType()}
\label{\detokenize{csapi/mvar:mvar-gettype}}\begin{quote}

\sphinxAtStartPar
Get types of variables in MVar object.

\sphinxAtStartPar
\sphinxstylestrong{Synopsis}
\begin{quote}

\sphinxAtStartPar
\sphinxcode{\sphinxupquote{NdArray\textless{}char\textgreater{} GetType()}}
\end{quote}

\sphinxAtStartPar
\sphinxstylestrong{Return}
\begin{quote}

\sphinxAtStartPar
multi\sphinxhyphen{}dimensional array of types of variables.
\end{quote}
\end{quote}

\subsubsection{MVar.GetUpperIIS()}
\label{\detokenize{csapi/mvar:mvar-getupperiis}}\begin{quote}

\sphinxAtStartPar
Get IIS status of upper bound of variables in MVar object.

\sphinxAtStartPar
\sphinxstylestrong{Synopsis}
\begin{quote}

\sphinxAtStartPar
\sphinxcode{\sphinxupquote{NdArray\textless{}int\textgreater{} GetUpperIIS()}}
\end{quote}

\sphinxAtStartPar
\sphinxstylestrong{Return}
\begin{quote}

\sphinxAtStartPar
multi\sphinxhyphen{}dimensional array of IIS status of upper bounds of variables.
\end{quote}
\end{quote}

\subsubsection{MVar.HStack()}
\label{\detokenize{csapi/mvar:mvar-hstack}}\begin{quote}

\sphinxAtStartPar
Stack with other MVar object along horizontal axis.

\sphinxAtStartPar
\sphinxstylestrong{Synopsis}
\begin{quote}

\sphinxAtStartPar
\sphinxcode{\sphinxupquote{MVar HStack(MVar other)}}
\end{quote}

\sphinxAtStartPar
\sphinxstylestrong{Arguments}
\begin{quote}

\sphinxAtStartPar
\sphinxcode{\sphinxupquote{other}}: a MVar object.
\end{quote}

\sphinxAtStartPar
\sphinxstylestrong{Return}
\begin{quote}

\sphinxAtStartPar
the result MVar object.
\end{quote}
\end{quote}

\subsubsection{MVar.Pick()}
\label{\detokenize{csapi/mvar:mvar-pick}}\begin{quote}

\sphinxAtStartPar
Given a list of indexes, get variables from MVar object.

\sphinxAtStartPar
\sphinxstylestrong{Synopsis}
\begin{quote}

\sphinxAtStartPar
\sphinxcode{\sphinxupquote{MVar Pick(NdArray\textless{}int\textgreater{} indexes)}}
\end{quote}

\sphinxAtStartPar
\sphinxstylestrong{Arguments}
\begin{quote}

\sphinxAtStartPar
\sphinxcode{\sphinxupquote{indexes}}: one or two dimensional indexes of elements. If two dimensional, each row is position of an element.
\end{quote}

\sphinxAtStartPar
\sphinxstylestrong{Return}
\begin{quote}

\sphinxAtStartPar
one\sphinxhyphen{}dimensional array of desired variables.
\end{quote}
\end{quote}

\subsubsection{MVar.Repeat()}
\label{\detokenize{csapi/mvar:mvar-repeat}}\begin{quote}

\sphinxAtStartPar
Repeat each element of MVar along given axis.

\sphinxAtStartPar
\sphinxstylestrong{Synopsis}
\begin{quote}

\sphinxAtStartPar
\sphinxcode{\sphinxupquote{MVar Repeat(long repeats, int axis)}}
\end{quote}

\sphinxAtStartPar
\sphinxstylestrong{Arguments}
\begin{quote}

\sphinxAtStartPar
\sphinxcode{\sphinxupquote{repeats}}: number of repetitions for each element.

\sphinxAtStartPar
\sphinxcode{\sphinxupquote{axis}}: axis of MVar.
\end{quote}

\sphinxAtStartPar
\sphinxstylestrong{Return}
\begin{quote}

\sphinxAtStartPar
new MVar object.
\end{quote}
\end{quote}

\subsubsection{MVar.RepeatBlock()}
\label{\detokenize{csapi/mvar:mvar-repeatblock}}\begin{quote}

\sphinxAtStartPar
Repeat an MVar a number of times along given axis.

\sphinxAtStartPar
\sphinxstylestrong{Synopsis}
\begin{quote}

\sphinxAtStartPar
\sphinxcode{\sphinxupquote{MVar RepeatBlock(long repeats, int axis)}}
\end{quote}

\sphinxAtStartPar
\sphinxstylestrong{Arguments}
\begin{quote}

\sphinxAtStartPar
\sphinxcode{\sphinxupquote{repeats}}: number of repetitions.

\sphinxAtStartPar
\sphinxcode{\sphinxupquote{axis}}: axis of MVar.
\end{quote}

\sphinxAtStartPar
\sphinxstylestrong{Return}
\begin{quote}

\sphinxAtStartPar
new MVar object.
\end{quote}
\end{quote}

\subsubsection{MVar.Represent()}
\label{\detokenize{csapi/mvar:mvar-represent}}\begin{quote}

\sphinxAtStartPar
String representation of MVar object.

\sphinxAtStartPar
\sphinxstylestrong{Synopsis}
\begin{quote}

\sphinxAtStartPar
\sphinxcode{\sphinxupquote{string Represent(int maxlen)}}
\end{quote}

\sphinxAtStartPar
\sphinxstylestrong{Arguments}
\begin{quote}

\sphinxAtStartPar
\sphinxcode{\sphinxupquote{maxlen}}: maximum buffer length for representations string.
\end{quote}

\sphinxAtStartPar
\sphinxstylestrong{Return}
\begin{quote}

\sphinxAtStartPar
string object.
\end{quote}
\end{quote}

\subsubsection{MVar.Reshape()}
\label{\detokenize{csapi/mvar:mvar-reshape}}\begin{quote}

\sphinxAtStartPar
Reshape MVar object to new shape.

\sphinxAtStartPar
\sphinxstylestrong{Synopsis}
\begin{quote}

\sphinxAtStartPar
\sphinxcode{\sphinxupquote{MVar Reshape(Shape shp)}}
\end{quote}

\sphinxAtStartPar
\sphinxstylestrong{Arguments}
\begin{quote}

\sphinxAtStartPar
\sphinxcode{\sphinxupquote{shp}}: new shape of M\sphinxhyphen{}dimensions.
\end{quote}

\sphinxAtStartPar
\sphinxstylestrong{Return}
\begin{quote}

\sphinxAtStartPar
M\sphinxhyphen{}dimensional MVar object.
\end{quote}
\end{quote}

\subsubsection{MVar.Set()}
\label{\detokenize{csapi/mvar:mvar-set}}\begin{quote}

\sphinxAtStartPar
Set values of information associated with variables in MVar object.

\sphinxAtStartPar
\sphinxstylestrong{Synopsis}
\begin{quote}

\sphinxAtStartPar
\sphinxcode{\sphinxupquote{void Set(string info, double val)}}
\end{quote}

\sphinxAtStartPar
\sphinxstylestrong{Arguments}
\begin{quote}

\sphinxAtStartPar
\sphinxcode{\sphinxupquote{info}}: name of information.

\sphinxAtStartPar
\sphinxcode{\sphinxupquote{val}}: value of information.
\end{quote}
\end{quote}

\subsubsection{MVar.Set()}
\label{\detokenize{csapi/mvar:id2}}\begin{quote}

\sphinxAtStartPar
Set values of information associated with variables in MVar object.

\sphinxAtStartPar
\sphinxstylestrong{Synopsis}
\begin{quote}

\sphinxAtStartPar
\sphinxcode{\sphinxupquote{void Set(string info, NdArray\textless{}double\textgreater{} vals)}}
\end{quote}

\sphinxAtStartPar
\sphinxstylestrong{Arguments}
\begin{quote}

\sphinxAtStartPar
\sphinxcode{\sphinxupquote{info}}: name of information.

\sphinxAtStartPar
\sphinxcode{\sphinxupquote{vals}}: multi\sphinxhyphen{}dimensional array of values of information.
\end{quote}
\end{quote}

\subsubsection{MVar.SetItem()}
\label{\detokenize{csapi/mvar:mvar-setitem}}\begin{quote}

\sphinxAtStartPar
Set variable of given index to MVar object.

\sphinxAtStartPar
\sphinxstylestrong{Synopsis}
\begin{quote}

\sphinxAtStartPar
\sphinxcode{\sphinxupquote{void SetItem(long idx, Var var)}}
\end{quote}

\sphinxAtStartPar
\sphinxstylestrong{Arguments}
\begin{quote}

\sphinxAtStartPar
\sphinxcode{\sphinxupquote{idx}}: index of element.

\sphinxAtStartPar
\sphinxcode{\sphinxupquote{var}}: Var object.
\end{quote}
\end{quote}

\subsubsection{MVar.Squeeze()}
\label{\detokenize{csapi/mvar:mvar-squeeze}}\begin{quote}

\sphinxAtStartPar
Remove axis of length 1 from shape of MVar object.

\sphinxAtStartPar
\sphinxstylestrong{Synopsis}
\begin{quote}

\sphinxAtStartPar
\sphinxcode{\sphinxupquote{MVar Squeeze(int axis)}}
\end{quote}

\sphinxAtStartPar
\sphinxstylestrong{Arguments}
\begin{quote}

\sphinxAtStartPar
\sphinxcode{\sphinxupquote{axis}}: axis of MVar, where the length is 1.
\end{quote}

\sphinxAtStartPar
\sphinxstylestrong{Return}
\begin{quote}

\sphinxAtStartPar
MVar object of (N\sphinxhyphen{}1)\sphinxhyphen{}dimensional shape.
\end{quote}
\end{quote}

\subsubsection{MVar.Stack()}
\label{\detokenize{csapi/mvar:mvar-stack}}\begin{quote}

\sphinxAtStartPar
Stack with other MVar object along given axis.

\sphinxAtStartPar
\sphinxstylestrong{Synopsis}
\begin{quote}

\sphinxAtStartPar
\sphinxcode{\sphinxupquote{MVar Stack(MVar other, int axis)}}
\end{quote}

\sphinxAtStartPar
\sphinxstylestrong{Arguments}
\begin{quote}

\sphinxAtStartPar
\sphinxcode{\sphinxupquote{other}}: a MVar object.

\sphinxAtStartPar
\sphinxcode{\sphinxupquote{axis}}: an axis of MVar.
\end{quote}

\sphinxAtStartPar
\sphinxstylestrong{Return}
\begin{quote}

\sphinxAtStartPar
the result MVar object.
\end{quote}
\end{quote}

\subsubsection{MVar.Sum()}
\label{\detokenize{csapi/mvar:mvar-sum}}\begin{quote}

\sphinxAtStartPar
Sum of all variables in MVar object.

\sphinxAtStartPar
\sphinxstylestrong{Synopsis}
\begin{quote}

\sphinxAtStartPar
\sphinxcode{\sphinxupquote{MLinExpr Sum()}}
\end{quote}

\sphinxAtStartPar
\sphinxstylestrong{Return}
\begin{quote}

\sphinxAtStartPar
sum in zero dimension.
\end{quote}
\end{quote}

\subsubsection{MVar.Sum()}
\label{\detokenize{csapi/mvar:id3}}\begin{quote}

\sphinxAtStartPar
Sum of variables at given axis of MVar object.

\sphinxAtStartPar
\sphinxstylestrong{Synopsis}
\begin{quote}

\sphinxAtStartPar
\sphinxcode{\sphinxupquote{MLinExpr Sum(int axis)}}
\end{quote}

\sphinxAtStartPar
\sphinxstylestrong{Arguments}
\begin{quote}

\sphinxAtStartPar
\sphinxcode{\sphinxupquote{axis}}: axis of MVar.
\end{quote}

\sphinxAtStartPar
\sphinxstylestrong{Return}
\begin{quote}

\sphinxAtStartPar
MLinExpr object in (N\sphinxhyphen{}1)\sphinxhyphen{}dimension.
\end{quote}
\end{quote}

\subsubsection{MVar.Transpose()}
\label{\detokenize{csapi/mvar:mvar-transpose}}\begin{quote}

\sphinxAtStartPar
Perform matrix transpose of MVar object.

\sphinxAtStartPar
\sphinxstylestrong{Synopsis}
\begin{quote}

\sphinxAtStartPar
\sphinxcode{\sphinxupquote{MVar Transpose()}}
\end{quote}

\sphinxAtStartPar
\sphinxstylestrong{Return}
\begin{quote}

\sphinxAtStartPar
transposed MVar object.
\end{quote}
\end{quote}

\subsubsection{MVar.VStack()}
\label{\detokenize{csapi/mvar:mvar-vstack}}\begin{quote}

\sphinxAtStartPar
Stack with other MVar object along vertical axis.

\sphinxAtStartPar
\sphinxstylestrong{Synopsis}
\begin{quote}

\sphinxAtStartPar
\sphinxcode{\sphinxupquote{MVar VStack(MVar other)}}
\end{quote}

\sphinxAtStartPar
\sphinxstylestrong{Arguments}
\begin{quote}

\sphinxAtStartPar
\sphinxcode{\sphinxupquote{other}}: a MVar object.
\end{quote}

\sphinxAtStartPar
\sphinxstylestrong{Return}
\begin{quote}

\sphinxAtStartPar
the result MVar object.
\end{quote}
\end{quote}

\subsection{MConstr}
\label{\detokenize{csharpapiref:mconstr}}\label{\detokenize{csharpapiref:chapcsharpapiref-mconstr}}
\sphinxAtStartPar
COPT multi\sphinxhyphen{}dimensional linear constraint object. It can be created by calling
the method \sphinxcode{\sphinxupquote{addMConstr}} of {\hyperref[\detokenize{csharpapiref:chapcsharpapiref-model}]{\sphinxcrossref{\DUrole{std,std-ref}{Model}}}}. The following methods
are provided:

\sphinxstepscope

\subsubsection{MConstr.Clone()}
\label{\detokenize{csapi/mconstr:mconstr-clone}}\label{\detokenize{csapi/mconstr::doc}}\begin{quote}

\sphinxAtStartPar
Clone MConstr object.

\sphinxAtStartPar
\sphinxstylestrong{Synopsis}
\begin{quote}

\sphinxAtStartPar
\sphinxcode{\sphinxupquote{MConstr Clone()}}
\end{quote}

\sphinxAtStartPar
\sphinxstylestrong{Return}
\begin{quote}

\sphinxAtStartPar
new MConstr object.
\end{quote}
\end{quote}

\subsubsection{MConstr.Diagonal()}
\label{\detokenize{csapi/mconstr:mconstr-diagonal}}\begin{quote}

\sphinxAtStartPar
Get diagonals of MConstr object.

\sphinxAtStartPar
\sphinxstylestrong{Synopsis}
\begin{quote}

\sphinxAtStartPar
\sphinxcode{\sphinxupquote{MConstr Diagonal(}}
\begin{quote}

\sphinxAtStartPar
\sphinxcode{\sphinxupquote{int offset,}}

\sphinxAtStartPar
\sphinxcode{\sphinxupquote{int axis1,}}

\sphinxAtStartPar
\sphinxcode{\sphinxupquote{int axis2)}}
\end{quote}
\end{quote}

\sphinxAtStartPar
\sphinxstylestrong{Arguments}
\begin{quote}

\sphinxAtStartPar
\sphinxcode{\sphinxupquote{offset}}: offset of the diagonal from the main diagonal. Can be positive or negative.

\sphinxAtStartPar
\sphinxcode{\sphinxupquote{axis1}}: 1st axis of MConstr.

\sphinxAtStartPar
\sphinxcode{\sphinxupquote{axis2}}: 2nd axis of MConstr.
\end{quote}

\sphinxAtStartPar
\sphinxstylestrong{Return}
\begin{quote}

\sphinxAtStartPar
(N\sphinxhyphen{}1)\sphinxhyphen{}dimensional diagonals.
\end{quote}
\end{quote}

\subsubsection{MConstr.Expand()}
\label{\detokenize{csapi/mconstr:mconstr-expand}}\begin{quote}

\sphinxAtStartPar
Expand shape of MConstr object.

\sphinxAtStartPar
\sphinxstylestrong{Synopsis}
\begin{quote}

\sphinxAtStartPar
\sphinxcode{\sphinxupquote{MConstr Expand(int axis)}}
\end{quote}

\sphinxAtStartPar
\sphinxstylestrong{Arguments}
\begin{quote}

\sphinxAtStartPar
\sphinxcode{\sphinxupquote{axis}}: axis of MConstr.
\end{quote}

\sphinxAtStartPar
\sphinxstylestrong{Return}
\begin{quote}

\sphinxAtStartPar
MConstr object of (N+1)\sphinxhyphen{}dimensional shape.
\end{quote}
\end{quote}

\subsubsection{MConstr.Flatten()}
\label{\detokenize{csapi/mconstr:mconstr-flatten}}\begin{quote}

\sphinxAtStartPar
Flatten a MConstr object to a 1\sphinxhyphen{}dimensional shape.

\sphinxAtStartPar
\sphinxstylestrong{Synopsis}
\begin{quote}

\sphinxAtStartPar
\sphinxcode{\sphinxupquote{MConstr Flatten()}}
\end{quote}

\sphinxAtStartPar
\sphinxstylestrong{Return}
\begin{quote}

\sphinxAtStartPar
a MConstr object collapsed into one dimension.
\end{quote}
\end{quote}

\subsubsection{MConstr.Get()}
\label{\detokenize{csapi/mconstr:mconstr-get}}\begin{quote}

\sphinxAtStartPar
Get values of information associated with constraints in MConstr object.

\sphinxAtStartPar
\sphinxstylestrong{Synopsis}
\begin{quote}

\sphinxAtStartPar
\sphinxcode{\sphinxupquote{NdArray\textless{}double\textgreater{} Get(string info)}}
\end{quote}

\sphinxAtStartPar
\sphinxstylestrong{Arguments}
\begin{quote}

\sphinxAtStartPar
\sphinxcode{\sphinxupquote{info}}: name of information.
\end{quote}

\sphinxAtStartPar
\sphinxstylestrong{Return}
\begin{quote}

\sphinxAtStartPar
multi\sphinxhyphen{}dimensional array of information of constraints.
\end{quote}
\end{quote}

\subsubsection{MConstr.GetDim()}
\label{\detokenize{csapi/mconstr:mconstr-getdim}}\begin{quote}

\sphinxAtStartPar
Get i\sphinxhyphen{}th dimension of MConstr object.

\sphinxAtStartPar
\sphinxstylestrong{Synopsis}
\begin{quote}

\sphinxAtStartPar
\sphinxcode{\sphinxupquote{long GetDim(int i)}}
\end{quote}

\sphinxAtStartPar
\sphinxstylestrong{Arguments}
\begin{quote}

\sphinxAtStartPar
\sphinxcode{\sphinxupquote{i}}: index of dimension
\end{quote}

\sphinxAtStartPar
\sphinxstylestrong{Return}
\begin{quote}

\sphinxAtStartPar
i\sphinxhyphen{}th dimension.
\end{quote}
\end{quote}

\subsubsection{MConstr.GetIdx()}
\label{\detokenize{csapi/mconstr:mconstr-getidx}}\begin{quote}

\sphinxAtStartPar
Get index of constraints in MConstr object.

\sphinxAtStartPar
\sphinxstylestrong{Synopsis}
\begin{quote}

\sphinxAtStartPar
\sphinxcode{\sphinxupquote{NdArray\textless{}int\textgreater{} GetIdx()}}
\end{quote}

\sphinxAtStartPar
\sphinxstylestrong{Return}
\begin{quote}

\sphinxAtStartPar
multi\sphinxhyphen{}dimensional array of indexes of constraints.
\end{quote}
\end{quote}

\subsubsection{MConstr.GetItem()}
\label{\detokenize{csapi/mconstr:mconstr-getitem}}\begin{quote}

\sphinxAtStartPar
Get constraint of given index from MConstr object.

\sphinxAtStartPar
\sphinxstylestrong{Synopsis}
\begin{quote}

\sphinxAtStartPar
\sphinxcode{\sphinxupquote{Constraint GetItem(long idx)}}
\end{quote}

\sphinxAtStartPar
\sphinxstylestrong{Arguments}
\begin{quote}

\sphinxAtStartPar
\sphinxcode{\sphinxupquote{idx}}: index of constraint.
\end{quote}

\sphinxAtStartPar
\sphinxstylestrong{Return}
\begin{quote}

\sphinxAtStartPar
Constraint object.
\end{quote}
\end{quote}

\subsubsection{MConstr.GetItem()}
\label{\detokenize{csapi/mconstr:id1}}\begin{quote}

\sphinxAtStartPar
Get sub\sphinxhyphen{}arrays of MConstr object, given view object.

\sphinxAtStartPar
\sphinxstylestrong{Synopsis}
\begin{quote}

\sphinxAtStartPar
\sphinxcode{\sphinxupquote{MConstr GetItem(View view)}}
\end{quote}

\sphinxAtStartPar
\sphinxstylestrong{Arguments}
\begin{quote}

\sphinxAtStartPar
\sphinxcode{\sphinxupquote{view}}: view of multi\sphinxhyphen{}dimensional array.
\end{quote}

\sphinxAtStartPar
\sphinxstylestrong{Return}
\begin{quote}

\sphinxAtStartPar
sub\sphinxhyphen{}arrays of MConstr object.
\end{quote}
\end{quote}

\subsubsection{MConstr.GetLowerIIS()}
\label{\detokenize{csapi/mconstr:mconstr-getloweriis}}\begin{quote}

\sphinxAtStartPar
Get IIS status of lower bound of constraints in MConstr object.

\sphinxAtStartPar
\sphinxstylestrong{Synopsis}
\begin{quote}

\sphinxAtStartPar
\sphinxcode{\sphinxupquote{NdArray\textless{}int\textgreater{} GetLowerIIS()}}
\end{quote}

\sphinxAtStartPar
\sphinxstylestrong{Return}
\begin{quote}

\sphinxAtStartPar
multi\sphinxhyphen{}dimensional array of IIS status of lower bounds of constraints.
\end{quote}
\end{quote}

\subsubsection{MConstr.GetND()}
\label{\detokenize{csapi/mconstr:mconstr-getnd}}\begin{quote}

\sphinxAtStartPar
Get number of dimensions of MConstr object.

\sphinxAtStartPar
\sphinxstylestrong{Synopsis}
\begin{quote}

\sphinxAtStartPar
\sphinxcode{\sphinxupquote{int GetND()}}
\end{quote}

\sphinxAtStartPar
\sphinxstylestrong{Return}
\begin{quote}

\sphinxAtStartPar
number of dimensions.
\end{quote}
\end{quote}

\subsubsection{MConstr.GetShape()}
\label{\detokenize{csapi/mconstr:mconstr-getshape}}\begin{quote}

\sphinxAtStartPar
Get shape of MConstr object.

\sphinxAtStartPar
\sphinxstylestrong{Synopsis}
\begin{quote}

\sphinxAtStartPar
\sphinxcode{\sphinxupquote{Shape GetShape()}}
\end{quote}

\sphinxAtStartPar
\sphinxstylestrong{Return}
\begin{quote}

\sphinxAtStartPar
shape object.
\end{quote}
\end{quote}

\subsubsection{MConstr.GetSize()}
\label{\detokenize{csapi/mconstr:mconstr-getsize}}\begin{quote}

\sphinxAtStartPar
Get size of MConstr object.

\sphinxAtStartPar
\sphinxstylestrong{Synopsis}
\begin{quote}

\sphinxAtStartPar
\sphinxcode{\sphinxupquote{long GetSize()}}
\end{quote}

\sphinxAtStartPar
\sphinxstylestrong{Return}
\begin{quote}

\sphinxAtStartPar
number of vars.
\end{quote}
\end{quote}

\subsubsection{MConstr.GetUpperIIS()}
\label{\detokenize{csapi/mconstr:mconstr-getupperiis}}\begin{quote}

\sphinxAtStartPar
Get IIS status of upper bound of constraints in MConstr object.

\sphinxAtStartPar
\sphinxstylestrong{Synopsis}
\begin{quote}

\sphinxAtStartPar
\sphinxcode{\sphinxupquote{NdArray\textless{}int\textgreater{} GetUpperIIS()}}
\end{quote}

\sphinxAtStartPar
\sphinxstylestrong{Return}
\begin{quote}

\sphinxAtStartPar
multi\sphinxhyphen{}dimensional array of IIS status of upper bounds of constraints.
\end{quote}
\end{quote}

\subsubsection{MConstr.HStack()}
\label{\detokenize{csapi/mconstr:mconstr-hstack}}\begin{quote}

\sphinxAtStartPar
Stack with other MConstr object along horizontal axis.

\sphinxAtStartPar
\sphinxstylestrong{Synopsis}
\begin{quote}

\sphinxAtStartPar
\sphinxcode{\sphinxupquote{MConstr HStack(MConstr other)}}
\end{quote}

\sphinxAtStartPar
\sphinxstylestrong{Arguments}
\begin{quote}

\sphinxAtStartPar
\sphinxcode{\sphinxupquote{other}}: a MConstr object.
\end{quote}

\sphinxAtStartPar
\sphinxstylestrong{Return}
\begin{quote}

\sphinxAtStartPar
the result MConstr object.
\end{quote}
\end{quote}

\subsubsection{MConstr.Pick()}
\label{\detokenize{csapi/mconstr:mconstr-pick}}\begin{quote}

\sphinxAtStartPar
Given a list of indexes, get constraints from MConstr object.

\sphinxAtStartPar
\sphinxstylestrong{Synopsis}
\begin{quote}

\sphinxAtStartPar
\sphinxcode{\sphinxupquote{MConstr Pick(NdArray\textless{}int\textgreater{} indexes)}}
\end{quote}

\sphinxAtStartPar
\sphinxstylestrong{Arguments}
\begin{quote}

\sphinxAtStartPar
\sphinxcode{\sphinxupquote{indexes}}: one or two dimensional indexes of elements. If two dimensional, each row is position of an element.
\end{quote}

\sphinxAtStartPar
\sphinxstylestrong{Return}
\begin{quote}

\sphinxAtStartPar
one\sphinxhyphen{}dimensional array of desired constraints.
\end{quote}
\end{quote}

\subsubsection{MConstr.Represent()}
\label{\detokenize{csapi/mconstr:mconstr-represent}}\begin{quote}

\sphinxAtStartPar
String representation of MConstr object.

\sphinxAtStartPar
\sphinxstylestrong{Synopsis}
\begin{quote}

\sphinxAtStartPar
\sphinxcode{\sphinxupquote{string Represent(int maxlen)}}
\end{quote}

\sphinxAtStartPar
\sphinxstylestrong{Arguments}
\begin{quote}

\sphinxAtStartPar
\sphinxcode{\sphinxupquote{maxlen}}: maximum buffer length for representations string.
\end{quote}

\sphinxAtStartPar
\sphinxstylestrong{Return}
\begin{quote}

\sphinxAtStartPar
string object.
\end{quote}
\end{quote}

\subsubsection{MConstr.Reshape()}
\label{\detokenize{csapi/mconstr:mconstr-reshape}}\begin{quote}

\sphinxAtStartPar
Reshape MConstr object to new shape.

\sphinxAtStartPar
\sphinxstylestrong{Synopsis}
\begin{quote}

\sphinxAtStartPar
\sphinxcode{\sphinxupquote{MConstr Reshape(Shape shp)}}
\end{quote}

\sphinxAtStartPar
\sphinxstylestrong{Arguments}
\begin{quote}

\sphinxAtStartPar
\sphinxcode{\sphinxupquote{shp}}: new shape of M\sphinxhyphen{}dimensions.
\end{quote}

\sphinxAtStartPar
\sphinxstylestrong{Return}
\begin{quote}

\sphinxAtStartPar
M\sphinxhyphen{}dimensional MConstr object.
\end{quote}
\end{quote}

\subsubsection{MConstr.Set()}
\label{\detokenize{csapi/mconstr:mconstr-set}}\begin{quote}

\sphinxAtStartPar
Set values of information associated with constraints in MConstr object.

\sphinxAtStartPar
\sphinxstylestrong{Synopsis}
\begin{quote}

\sphinxAtStartPar
\sphinxcode{\sphinxupquote{void Set(string info, NdArray\textless{}double\textgreater{} vals)}}
\end{quote}

\sphinxAtStartPar
\sphinxstylestrong{Arguments}
\begin{quote}

\sphinxAtStartPar
\sphinxcode{\sphinxupquote{info}}: name of information.

\sphinxAtStartPar
\sphinxcode{\sphinxupquote{vals}}: multi\sphinxhyphen{}dimensional array of values of information.
\end{quote}
\end{quote}

\subsubsection{MConstr.Set()}
\label{\detokenize{csapi/mconstr:id2}}\begin{quote}

\sphinxAtStartPar
Set values of information associated with constraints in MConstr object.

\sphinxAtStartPar
\sphinxstylestrong{Synopsis}
\begin{quote}

\sphinxAtStartPar
\sphinxcode{\sphinxupquote{void Set(string info, double val)}}
\end{quote}

\sphinxAtStartPar
\sphinxstylestrong{Arguments}
\begin{quote}

\sphinxAtStartPar
\sphinxcode{\sphinxupquote{info}}: name of information.

\sphinxAtStartPar
\sphinxcode{\sphinxupquote{val}}: value of information.
\end{quote}
\end{quote}

\subsubsection{MConstr.SetItem()}
\label{\detokenize{csapi/mconstr:mconstr-setitem}}\begin{quote}

\sphinxAtStartPar
Set constraint of given index to MConstr object.

\sphinxAtStartPar
\sphinxstylestrong{Synopsis}
\begin{quote}

\sphinxAtStartPar
\sphinxcode{\sphinxupquote{void SetItem(long idx, Constraint constr)}}
\end{quote}

\sphinxAtStartPar
\sphinxstylestrong{Arguments}
\begin{quote}

\sphinxAtStartPar
\sphinxcode{\sphinxupquote{idx}}: index of element.

\sphinxAtStartPar
\sphinxcode{\sphinxupquote{constr}}: Constraint object.
\end{quote}
\end{quote}

\subsubsection{MConstr.Squeeze()}
\label{\detokenize{csapi/mconstr:mconstr-squeeze}}\begin{quote}

\sphinxAtStartPar
Remove axis of length 1 from shape of MConstr object.

\sphinxAtStartPar
\sphinxstylestrong{Synopsis}
\begin{quote}

\sphinxAtStartPar
\sphinxcode{\sphinxupquote{MConstr Squeeze(int axis)}}
\end{quote}

\sphinxAtStartPar
\sphinxstylestrong{Arguments}
\begin{quote}

\sphinxAtStartPar
\sphinxcode{\sphinxupquote{axis}}: axis of MConstr, where the length is 1.
\end{quote}

\sphinxAtStartPar
\sphinxstylestrong{Return}
\begin{quote}

\sphinxAtStartPar
MConstr object of (N\sphinxhyphen{}1)\sphinxhyphen{}dimensional shape.
\end{quote}
\end{quote}

\subsubsection{MConstr.Stack()}
\label{\detokenize{csapi/mconstr:mconstr-stack}}\begin{quote}

\sphinxAtStartPar
Stack with other MConstr object along given axis.

\sphinxAtStartPar
\sphinxstylestrong{Synopsis}
\begin{quote}

\sphinxAtStartPar
\sphinxcode{\sphinxupquote{MConstr Stack(MConstr other, int axis)}}
\end{quote}

\sphinxAtStartPar
\sphinxstylestrong{Arguments}
\begin{quote}

\sphinxAtStartPar
\sphinxcode{\sphinxupquote{other}}: a MConstr object.

\sphinxAtStartPar
\sphinxcode{\sphinxupquote{axis}}: an axis of MConstr.
\end{quote}

\sphinxAtStartPar
\sphinxstylestrong{Return}
\begin{quote}

\sphinxAtStartPar
the result MConstr object.
\end{quote}
\end{quote}

\subsubsection{MConstr.Transpose()}
\label{\detokenize{csapi/mconstr:mconstr-transpose}}\begin{quote}

\sphinxAtStartPar
Perform matrix transpose of MConstr object.

\sphinxAtStartPar
\sphinxstylestrong{Synopsis}
\begin{quote}

\sphinxAtStartPar
\sphinxcode{\sphinxupquote{MConstr Transpose()}}
\end{quote}

\sphinxAtStartPar
\sphinxstylestrong{Return}
\begin{quote}

\sphinxAtStartPar
transposed MConstr object.
\end{quote}
\end{quote}

\subsubsection{MConstr.VStack()}
\label{\detokenize{csapi/mconstr:mconstr-vstack}}\begin{quote}

\sphinxAtStartPar
Stack with other MConstr object along vertical axis.

\sphinxAtStartPar
\sphinxstylestrong{Synopsis}
\begin{quote}

\sphinxAtStartPar
\sphinxcode{\sphinxupquote{MConstr VStack(MConstr other)}}
\end{quote}

\sphinxAtStartPar
\sphinxstylestrong{Arguments}
\begin{quote}

\sphinxAtStartPar
\sphinxcode{\sphinxupquote{other}}: a MConstr object.
\end{quote}

\sphinxAtStartPar
\sphinxstylestrong{Return}
\begin{quote}

\sphinxAtStartPar
the result MConstr object.
\end{quote}
\end{quote}

\subsection{MConstrBuilder}
\label{\detokenize{csharpapiref:mconstrbuilder}}\label{\detokenize{csharpapiref:chapcsharpapiref-mconstrbuilder}}
\sphinxAtStartPar
COPT builder object of multi\sphinxhyphen{}dimensional linear constraints. It is used to
generate multi\sphinxhyphen{}dimensional linear constraints and support operations with the
built\sphinxhyphen{}in multi\sphinxhyphen{}dimensional array {\hyperref[\detokenize{csharpapiref:chapcsharpapiref-ndarray}]{\sphinxcrossref{\DUrole{std,std-ref}{NdArray}}}} in COPT.
It is recommended to create MConstrBuilder object by comparing two objects, one of
which should be {\hyperref[\detokenize{csharpapiref:chapcsharpapiref-mvar}]{\sphinxcrossref{\DUrole{std,std-ref}{MVar}}}} object or
{\hyperref[\detokenize{csharpapiref:chapcsharpapiref-mlinexpr}]{\sphinxcrossref{\DUrole{std,std-ref}{MLinExpr}}}} object, by comparison operators.
The following methods are provided:

\sphinxstepscope

\subsubsection{MConstrBuilder.MConstrBuilder()}
\label{\detokenize{csapi/mconstrbuilder:mconstrbuilder-mconstrbuilder}}\label{\detokenize{csapi/mconstrbuilder::doc}}\begin{quote}

\sphinxAtStartPar
Construct a MConstrBuilder object with the given shape.

\sphinxAtStartPar
\sphinxstylestrong{Synopsis}
\begin{quote}

\sphinxAtStartPar
\sphinxcode{\sphinxupquote{MConstrBuilder(Shape shp)}}
\end{quote}

\sphinxAtStartPar
\sphinxstylestrong{Arguments}
\begin{quote}

\sphinxAtStartPar
\sphinxcode{\sphinxupquote{shp}}: shape of MConstrBuilder.
\end{quote}
\end{quote}

\subsubsection{MConstrBuilder.Flatten()}
\label{\detokenize{csapi/mconstrbuilder:mconstrbuilder-flatten}}\begin{quote}

\sphinxAtStartPar
Flatten a MConstrBuilder object to a 1\sphinxhyphen{}dimensional shape.

\sphinxAtStartPar
\sphinxstylestrong{Synopsis}
\begin{quote}

\sphinxAtStartPar
\sphinxcode{\sphinxupquote{MConstrBuilder Flatten()}}
\end{quote}

\sphinxAtStartPar
\sphinxstylestrong{Return}
\begin{quote}

\sphinxAtStartPar
a MConstrBuilder object collapsed into one dimension.
\end{quote}
\end{quote}

\subsubsection{MConstrBuilder.GetExpr()}
\label{\detokenize{csapi/mconstrbuilder:mconstrbuilder-getexpr}}\begin{quote}

\sphinxAtStartPar
Get N\sphinxhyphen{}dimensional linear expressions associated with N\sphinxhyphen{}dimensional constraints.

\sphinxAtStartPar
\sphinxstylestrong{Synopsis}
\begin{quote}

\sphinxAtStartPar
\sphinxcode{\sphinxupquote{MLinExpr GetExpr()}}
\end{quote}

\sphinxAtStartPar
\sphinxstylestrong{Return}
\begin{quote}

\sphinxAtStartPar
MLinExpr object.
\end{quote}
\end{quote}

\subsubsection{MConstrBuilder.GetND()}
\label{\detokenize{csapi/mconstrbuilder:mconstrbuilder-getnd}}\begin{quote}

\sphinxAtStartPar
Get number of dimensions of MConstrBuilder object.

\sphinxAtStartPar
\sphinxstylestrong{Synopsis}
\begin{quote}

\sphinxAtStartPar
\sphinxcode{\sphinxupquote{int GetND()}}
\end{quote}

\sphinxAtStartPar
\sphinxstylestrong{Return}
\begin{quote}

\sphinxAtStartPar
number of dimensions.
\end{quote}
\end{quote}

\subsubsection{MConstrBuilder.GetRange()}
\label{\detokenize{csapi/mconstrbuilder:mconstrbuilder-getrange}}\begin{quote}

\sphinxAtStartPar
Get range from lower bound to upper bound of N\sphinxhyphen{}dimensional range constraints.

\sphinxAtStartPar
\sphinxstylestrong{Synopsis}
\begin{quote}

\sphinxAtStartPar
\sphinxcode{\sphinxupquote{double GetRange()}}
\end{quote}

\sphinxAtStartPar
\sphinxstylestrong{Return}
\begin{quote}

\sphinxAtStartPar
length from lower bound to upper bound of range constraints.
\end{quote}
\end{quote}

\subsubsection{MConstrBuilder.GetSense()}
\label{\detokenize{csapi/mconstrbuilder:mconstrbuilder-getsense}}\begin{quote}

\sphinxAtStartPar
Get sense associated with N\sphinxhyphen{}dimensional constraints.

\sphinxAtStartPar
\sphinxstylestrong{Synopsis}
\begin{quote}

\sphinxAtStartPar
\sphinxcode{\sphinxupquote{char GetSense()}}
\end{quote}

\sphinxAtStartPar
\sphinxstylestrong{Return}
\begin{quote}

\sphinxAtStartPar
constraint sense.
\end{quote}
\end{quote}

\subsubsection{MConstrBuilder.Set\textless{}T\textgreater{}()}
\label{\detokenize{csapi/mconstrbuilder:mconstrbuilder-set-t}}\begin{quote}

\sphinxAtStartPar
Set N\sphinxhyphen{}dimensional constraints to its builder object.

\sphinxAtStartPar
\sphinxstylestrong{Synopsis}
\begin{quote}

\sphinxAtStartPar
\sphinxcode{\sphinxupquote{void Set\textless{}T\textgreater{}(}}
\begin{quote}

\sphinxAtStartPar
\sphinxcode{\sphinxupquote{MLinExpr expr,}}

\sphinxAtStartPar
\sphinxcode{\sphinxupquote{char sense,}}

\sphinxAtStartPar
\sphinxcode{\sphinxupquote{NdArray\textless{}T\textgreater{} rhs)}}
\end{quote}
\end{quote}

\sphinxAtStartPar
\sphinxstylestrong{Arguments}
\begin{quote}

\sphinxAtStartPar
\sphinxcode{\sphinxupquote{expr}}: MLinExpr object

\sphinxAtStartPar
\sphinxcode{\sphinxupquote{sense}}: constraint sense other than COPT\_RANGE.

\sphinxAtStartPar
\sphinxcode{\sphinxupquote{rhs}}: N\sphinxhyphen{}dimensional constants at right side of constraints.
\end{quote}
\end{quote}

\subsubsection{MConstrBuilder.Set()}
\label{\detokenize{csapi/mconstrbuilder:mconstrbuilder-set}}\begin{quote}

\sphinxAtStartPar
Set N\sphinxhyphen{}dimensional constraints to its builder object.

\sphinxAtStartPar
\sphinxstylestrong{Synopsis}
\begin{quote}

\sphinxAtStartPar
\sphinxcode{\sphinxupquote{void Set(}}
\begin{quote}

\sphinxAtStartPar
\sphinxcode{\sphinxupquote{MLinExpr expr,}}

\sphinxAtStartPar
\sphinxcode{\sphinxupquote{char sense,}}

\sphinxAtStartPar
\sphinxcode{\sphinxupquote{double rhs)}}
\end{quote}
\end{quote}

\sphinxAtStartPar
\sphinxstylestrong{Arguments}
\begin{quote}

\sphinxAtStartPar
\sphinxcode{\sphinxupquote{expr}}: MLinExpr object

\sphinxAtStartPar
\sphinxcode{\sphinxupquote{sense}}: constraint sense other than COPT\_RANGE.

\sphinxAtStartPar
\sphinxcode{\sphinxupquote{rhs}}: constant of right side of constraints.
\end{quote}
\end{quote}

\subsubsection{MConstrBuilder.Set()}
\label{\detokenize{csapi/mconstrbuilder:id1}}\begin{quote}

\sphinxAtStartPar
Set N\sphinxhyphen{}dimensional constraints to its builder object.

\sphinxAtStartPar
\sphinxstylestrong{Synopsis}
\begin{quote}

\sphinxAtStartPar
\sphinxcode{\sphinxupquote{void Set(}}
\begin{quote}

\sphinxAtStartPar
\sphinxcode{\sphinxupquote{MLinExpr expr,}}

\sphinxAtStartPar
\sphinxcode{\sphinxupquote{char sense,}}

\sphinxAtStartPar
\sphinxcode{\sphinxupquote{MVar rhs)}}
\end{quote}
\end{quote}

\sphinxAtStartPar
\sphinxstylestrong{Arguments}
\begin{quote}

\sphinxAtStartPar
\sphinxcode{\sphinxupquote{expr}}: MLinExpr object

\sphinxAtStartPar
\sphinxcode{\sphinxupquote{sense}}: constraint sense other than COPT\_RANGE.

\sphinxAtStartPar
\sphinxcode{\sphinxupquote{rhs}}: MVar object at right side of constraints.
\end{quote}
\end{quote}

\subsubsection{MConstrBuilder.Set()}
\label{\detokenize{csapi/mconstrbuilder:id2}}\begin{quote}

\sphinxAtStartPar
Set N\sphinxhyphen{}dimensional constraints to its builder object.

\sphinxAtStartPar
\sphinxstylestrong{Synopsis}
\begin{quote}

\sphinxAtStartPar
\sphinxcode{\sphinxupquote{void Set(}}
\begin{quote}

\sphinxAtStartPar
\sphinxcode{\sphinxupquote{MLinExpr expr,}}

\sphinxAtStartPar
\sphinxcode{\sphinxupquote{char sense,}}

\sphinxAtStartPar
\sphinxcode{\sphinxupquote{MLinExpr rhs)}}
\end{quote}
\end{quote}

\sphinxAtStartPar
\sphinxstylestrong{Arguments}
\begin{quote}

\sphinxAtStartPar
\sphinxcode{\sphinxupquote{expr}}: MLinExpr object

\sphinxAtStartPar
\sphinxcode{\sphinxupquote{sense}}: constraint sense other than COPT\_RANGE.

\sphinxAtStartPar
\sphinxcode{\sphinxupquote{rhs}}: MLinExpr object at right side of constraints.
\end{quote}
\end{quote}

\subsubsection{MConstrBuilder.SetRange()}
\label{\detokenize{csapi/mconstrbuilder:mconstrbuilder-setrange}}\begin{quote}

\sphinxAtStartPar
Set N\sphinxhyphen{}dimensional range constraints to its builder object.

\sphinxAtStartPar
\sphinxstylestrong{Synopsis}
\begin{quote}

\sphinxAtStartPar
\sphinxcode{\sphinxupquote{void SetRange(MLinExpr expr, double range)}}
\end{quote}

\sphinxAtStartPar
\sphinxstylestrong{Arguments}
\begin{quote}

\sphinxAtStartPar
\sphinxcode{\sphinxupquote{expr}}: MLinExpr object.

\sphinxAtStartPar
\sphinxcode{\sphinxupquote{range}}: length from lower bound to upper bound of the constraint. Must greater than 0.
\end{quote}
\end{quote}

\subsection{MExpression}
\label{\detokenize{csharpapiref:mexpression}}\label{\detokenize{csharpapiref:chapcsharpapiref-mexpression}}
\sphinxAtStartPar
The MExpression class is a generalized version of {\hyperref[\detokenize{csharpapiref:chapcsharpapiref-expr}]{\sphinxcrossref{\DUrole{std,std-ref}{Expr}}}}.
It represents a linear expression and supports most of methods in Expr class.
In addition, it supports linear combination of multi\sphinxhyphen{}dimensional objects,
such as {\hyperref[\detokenize{csharpapiref:chapcsharpapiref-mvar}]{\sphinxcrossref{\DUrole{std,std-ref}{MVar}}}} object and {\hyperref[\detokenize{csharpapiref:chapcsharpapiref-ndarray}]{\sphinxcrossref{\DUrole{std,std-ref}{NdArray}}}} object.
The following methods are provided:

\sphinxstepscope

\subsubsection{MExpression.MExpression()}
\label{\detokenize{csapi/mexpression:mexpression-mexpression}}\label{\detokenize{csapi/mexpression::doc}}\begin{quote}

\sphinxAtStartPar
Construct a MExpression object with the given constant.

\sphinxAtStartPar
\sphinxstylestrong{Synopsis}
\begin{quote}

\sphinxAtStartPar
\sphinxcode{\sphinxupquote{MExpression(double constant)}}
\end{quote}

\sphinxAtStartPar
\sphinxstylestrong{Arguments}
\begin{quote}

\sphinxAtStartPar
\sphinxcode{\sphinxupquote{constant}}: constant number.
\end{quote}
\end{quote}

\subsubsection{MExpression.MExpression()}
\label{\detokenize{csapi/mexpression:id1}}\begin{quote}

\sphinxAtStartPar
Construct a MExpression object with the given variable.

\sphinxAtStartPar
\sphinxstylestrong{Synopsis}
\begin{quote}

\sphinxAtStartPar
\sphinxcode{\sphinxupquote{MExpression(Var var)}}
\end{quote}

\sphinxAtStartPar
\sphinxstylestrong{Arguments}
\begin{quote}

\sphinxAtStartPar
\sphinxcode{\sphinxupquote{var}}: variable object.
\end{quote}
\end{quote}

\subsubsection{MExpression.MExpression()}
\label{\detokenize{csapi/mexpression:id2}}\begin{quote}

\sphinxAtStartPar
Construct a MExpression object with the given linear expression.

\sphinxAtStartPar
\sphinxstylestrong{Synopsis}
\begin{quote}

\sphinxAtStartPar
\sphinxcode{\sphinxupquote{MExpression(Expr expr)}}
\end{quote}

\sphinxAtStartPar
\sphinxstylestrong{Arguments}
\begin{quote}

\sphinxAtStartPar
\sphinxcode{\sphinxupquote{expr}}: a linear expression.
\end{quote}
\end{quote}

\subsubsection{MExpression.AddConstant()}
\label{\detokenize{csapi/mexpression:mexpression-addconstant}}\begin{quote}

\sphinxAtStartPar
Add constant for the expression.

\sphinxAtStartPar
\sphinxstylestrong{Synopsis}
\begin{quote}

\sphinxAtStartPar
\sphinxcode{\sphinxupquote{void AddConstant(double constant)}}
\end{quote}

\sphinxAtStartPar
\sphinxstylestrong{Arguments}
\begin{quote}

\sphinxAtStartPar
\sphinxcode{\sphinxupquote{constant}}: the value of the constant.
\end{quote}
\end{quote}

\subsubsection{MExpression.AddExpr()}
\label{\detokenize{csapi/mexpression:mexpression-addexpr}}\begin{quote}

\sphinxAtStartPar
Add a linear expression to MExpression object.

\sphinxAtStartPar
\sphinxstylestrong{Synopsis}
\begin{quote}

\sphinxAtStartPar
\sphinxcode{\sphinxupquote{void AddExpr(Expr expr, double mult)}}
\end{quote}

\sphinxAtStartPar
\sphinxstylestrong{Arguments}
\begin{quote}

\sphinxAtStartPar
\sphinxcode{\sphinxupquote{expr}}: linear expression object.

\sphinxAtStartPar
\sphinxcode{\sphinxupquote{mult}}: the multiplier of linear expression, default value is 1.0.
\end{quote}
\end{quote}

\subsubsection{MExpression.AddMExpr()}
\label{\detokenize{csapi/mexpression:mexpression-addmexpr}}\begin{quote}

\sphinxAtStartPar
Add MExpression to MExpression object.

\sphinxAtStartPar
\sphinxstylestrong{Synopsis}
\begin{quote}

\sphinxAtStartPar
\sphinxcode{\sphinxupquote{void AddMExpr(MExpression expr, double mult)}}
\end{quote}

\sphinxAtStartPar
\sphinxstylestrong{Arguments}
\begin{quote}

\sphinxAtStartPar
\sphinxcode{\sphinxupquote{expr}}: MExpression object.

\sphinxAtStartPar
\sphinxcode{\sphinxupquote{mult}}: the multiplier of MExpression, default value is 1.0.
\end{quote}
\end{quote}

\subsubsection{MExpression.AddTerm()}
\label{\detokenize{csapi/mexpression:mexpression-addterm}}\begin{quote}

\sphinxAtStartPar
Add a linear term to MExpression object.

\sphinxAtStartPar
\sphinxstylestrong{Synopsis}
\begin{quote}

\sphinxAtStartPar
\sphinxcode{\sphinxupquote{void AddTerm(Var var, double coeff)}}
\end{quote}

\sphinxAtStartPar
\sphinxstylestrong{Arguments}
\begin{quote}

\sphinxAtStartPar
\sphinxcode{\sphinxupquote{var}}: variable of new term.

\sphinxAtStartPar
\sphinxcode{\sphinxupquote{coeff}}: coefficient of new term.
\end{quote}
\end{quote}

\subsubsection{MExpression.Clone()}
\label{\detokenize{csapi/mexpression:mexpression-clone}}\begin{quote}

\sphinxAtStartPar
Clone MExpression object.

\sphinxAtStartPar
\sphinxstylestrong{Synopsis}
\begin{quote}

\sphinxAtStartPar
\sphinxcode{\sphinxupquote{MExpression Clone()}}
\end{quote}

\sphinxAtStartPar
\sphinxstylestrong{Return}
\begin{quote}

\sphinxAtStartPar
new MExpression object.
\end{quote}
\end{quote}

\subsubsection{MExpression.Evaluate()}
\label{\detokenize{csapi/mexpression:mexpression-evaluate}}\begin{quote}

\sphinxAtStartPar
evaluate MExpression object after solving.

\sphinxAtStartPar
\sphinxstylestrong{Synopsis}
\begin{quote}

\sphinxAtStartPar
\sphinxcode{\sphinxupquote{double Evaluate()}}
\end{quote}

\sphinxAtStartPar
\sphinxstylestrong{Return}
\begin{quote}

\sphinxAtStartPar
value of MExpression object.
\end{quote}
\end{quote}

\subsubsection{MExpression.GetConstant()}
\label{\detokenize{csapi/mexpression:mexpression-getconstant}}\begin{quote}

\sphinxAtStartPar
Get constant in expression.

\sphinxAtStartPar
\sphinxstylestrong{Synopsis}
\begin{quote}

\sphinxAtStartPar
\sphinxcode{\sphinxupquote{double GetConstant()}}
\end{quote}

\sphinxAtStartPar
\sphinxstylestrong{Return}
\begin{quote}

\sphinxAtStartPar
constant in expression.
\end{quote}
\end{quote}

\subsubsection{MExpression.Represent()}
\label{\detokenize{csapi/mexpression:mexpression-represent}}\begin{quote}

\sphinxAtStartPar
String representation of MExpression object.

\sphinxAtStartPar
\sphinxstylestrong{Synopsis}
\begin{quote}

\sphinxAtStartPar
\sphinxcode{\sphinxupquote{string Represent(int maxlen)}}
\end{quote}

\sphinxAtStartPar
\sphinxstylestrong{Arguments}
\begin{quote}

\sphinxAtStartPar
\sphinxcode{\sphinxupquote{maxlen}}: maximum buffer length for representations string.
\end{quote}

\sphinxAtStartPar
\sphinxstylestrong{Return}
\begin{quote}

\sphinxAtStartPar
string object.
\end{quote}
\end{quote}

\subsubsection{MExpression.SetConstant()}
\label{\detokenize{csapi/mexpression:mexpression-setconstant}}\begin{quote}

\sphinxAtStartPar
Set constant for the expression.

\sphinxAtStartPar
\sphinxstylestrong{Synopsis}
\begin{quote}

\sphinxAtStartPar
\sphinxcode{\sphinxupquote{void SetConstant(double constant)}}
\end{quote}

\sphinxAtStartPar
\sphinxstylestrong{Arguments}
\begin{quote}

\sphinxAtStartPar
\sphinxcode{\sphinxupquote{constant}}: the value of the constant.
\end{quote}
\end{quote}

\subsection{MLinExpr}
\label{\detokenize{csharpapiref:mlinexpr}}\label{\detokenize{csharpapiref:chapcsharpapiref-mlinexpr}}
\sphinxAtStartPar
COPT multi\sphinxhyphen{}dimensional linear expression object. It is used to construct
multi\sphinxhyphen{}dimensional linear expressions and perform operations with the built\sphinxhyphen{}in
multi\sphinxhyphen{}dimensional array {\hyperref[\detokenize{csharpapiref:chapcsharpapiref-ndarray}]{\sphinxcrossref{\DUrole{std,std-ref}{NdArray}}}} in COPT. Its elements
are {\hyperref[\detokenize{csharpapiref:chapcsharpapiref-mexpression}]{\sphinxcrossref{\DUrole{std,std-ref}{MExpression}}}} objects.  It can be created by linear
combination of {\hyperref[\detokenize{csharpapiref:chapcsharpapiref-mvar}]{\sphinxcrossref{\DUrole{std,std-ref}{MVar}}}} objects. The following methods are
provided:

\sphinxstepscope

\subsubsection{MLinExpr.AddConstant()}
\label{\detokenize{csapi/mlinexpr:mlinexpr-addconstant}}\label{\detokenize{csapi/mlinexpr::doc}}\begin{quote}

\sphinxAtStartPar
Add constant to each expression in MLinExpr object.

\sphinxAtStartPar
\sphinxstylestrong{Synopsis}
\begin{quote}

\sphinxAtStartPar
\sphinxcode{\sphinxupquote{void AddConstant(double constant)}}
\end{quote}

\sphinxAtStartPar
\sphinxstylestrong{Arguments}
\begin{quote}

\sphinxAtStartPar
\sphinxcode{\sphinxupquote{constant}}: the value of the constant.
\end{quote}
\end{quote}

\subsubsection{MLinExpr.AddConstant()}
\label{\detokenize{csapi/mlinexpr:id1}}\begin{quote}

\sphinxAtStartPar
Add constants to each expression in MLinExpr object.

\sphinxAtStartPar
\sphinxstylestrong{Synopsis}
\begin{quote}

\sphinxAtStartPar
\sphinxcode{\sphinxupquote{void AddConstant(NdArray\textless{}double\textgreater{} constants)}}
\end{quote}

\sphinxAtStartPar
\sphinxstylestrong{Arguments}
\begin{quote}

\sphinxAtStartPar
\sphinxcode{\sphinxupquote{constants}}: N\sphinxhyphen{}dimension NdArray object.
\end{quote}
\end{quote}

\subsubsection{MLinExpr.AddExpr()}
\label{\detokenize{csapi/mlinexpr:mlinexpr-addexpr}}\begin{quote}

\sphinxAtStartPar
Add a linear expression to each expression in MLinExpr object.

\sphinxAtStartPar
\sphinxstylestrong{Synopsis}
\begin{quote}

\sphinxAtStartPar
\sphinxcode{\sphinxupquote{void AddExpr(Expr expr, double mult)}}
\end{quote}

\sphinxAtStartPar
\sphinxstylestrong{Arguments}
\begin{quote}

\sphinxAtStartPar
\sphinxcode{\sphinxupquote{expr}}: linear expression object.

\sphinxAtStartPar
\sphinxcode{\sphinxupquote{mult}}: the multiplier of linear expression, default value is 1.0.
\end{quote}
\end{quote}

\subsubsection{MLinExpr.AddMExpr()}
\label{\detokenize{csapi/mlinexpr:mlinexpr-addmexpr}}\begin{quote}

\sphinxAtStartPar
Add MExpression to each expression in MLinExpr object.

\sphinxAtStartPar
\sphinxstylestrong{Synopsis}
\begin{quote}

\sphinxAtStartPar
\sphinxcode{\sphinxupquote{void AddMExpr(MExpression expr, double mult)}}
\end{quote}

\sphinxAtStartPar
\sphinxstylestrong{Arguments}
\begin{quote}

\sphinxAtStartPar
\sphinxcode{\sphinxupquote{expr}}: MExpression object.

\sphinxAtStartPar
\sphinxcode{\sphinxupquote{mult}}: the multiplier of MExpression, default value is 1.0.
\end{quote}
\end{quote}

\subsubsection{MLinExpr.AddMLinExpr()}
\label{\detokenize{csapi/mlinexpr:mlinexpr-addmlinexpr}}\begin{quote}

\sphinxAtStartPar
Add linear expressions to MLinExpr object.

\sphinxAtStartPar
\sphinxstylestrong{Synopsis}
\begin{quote}

\sphinxAtStartPar
\sphinxcode{\sphinxupquote{void AddMLinExpr(MLinExpr exprs, double mult)}}
\end{quote}

\sphinxAtStartPar
\sphinxstylestrong{Arguments}
\begin{quote}

\sphinxAtStartPar
\sphinxcode{\sphinxupquote{exprs}}: N\sphinxhyphen{}dimension MLinExpr object.

\sphinxAtStartPar
\sphinxcode{\sphinxupquote{mult}}: the same multiplier for added linear expressions, default value is 1.0.
\end{quote}
\end{quote}

\subsubsection{MLinExpr.AddTerm()}
\label{\detokenize{csapi/mlinexpr:mlinexpr-addterm}}\begin{quote}

\sphinxAtStartPar
Add a linear term to MLinExpr object.

\sphinxAtStartPar
\sphinxstylestrong{Synopsis}
\begin{quote}

\sphinxAtStartPar
\sphinxcode{\sphinxupquote{void AddTerm(Var var, double coeff)}}
\end{quote}

\sphinxAtStartPar
\sphinxstylestrong{Arguments}
\begin{quote}

\sphinxAtStartPar
\sphinxcode{\sphinxupquote{var}}: variable of new term.

\sphinxAtStartPar
\sphinxcode{\sphinxupquote{coeff}}: coefficient of new term.
\end{quote}
\end{quote}

\subsubsection{MLinExpr.AddTerms()}
\label{\detokenize{csapi/mlinexpr:mlinexpr-addterms}}\begin{quote}

\sphinxAtStartPar
Add terms to expressions in MLinExpr object.

\sphinxAtStartPar
\sphinxstylestrong{Synopsis}
\begin{quote}

\sphinxAtStartPar
\sphinxcode{\sphinxupquote{void AddTerms(MVar vars, NdArray\textless{}double\textgreater{} coeffs)}}
\end{quote}

\sphinxAtStartPar
\sphinxstylestrong{Arguments}
\begin{quote}

\sphinxAtStartPar
\sphinxcode{\sphinxupquote{vars}}: N\sphinxhyphen{}dimension MVar object for added terms.

\sphinxAtStartPar
\sphinxcode{\sphinxupquote{coeffs}}: N\sphinxhyphen{}dimension NdArray object of coefficients for added terms.
\end{quote}
\end{quote}

\subsubsection{MLinExpr.AddTerms()}
\label{\detokenize{csapi/mlinexpr:id2}}\begin{quote}

\sphinxAtStartPar
Add terms to expressions in MLinExpr object.

\sphinxAtStartPar
\sphinxstylestrong{Synopsis}
\begin{quote}

\sphinxAtStartPar
\sphinxcode{\sphinxupquote{void AddTerms(MVar vars, double mult)}}
\end{quote}

\sphinxAtStartPar
\sphinxstylestrong{Arguments}
\begin{quote}

\sphinxAtStartPar
\sphinxcode{\sphinxupquote{vars}}: N\sphinxhyphen{}dimension MVar object for added terms.

\sphinxAtStartPar
\sphinxcode{\sphinxupquote{mult}}: the same coefficient for added terms, default value 1.0.
\end{quote}
\end{quote}

\subsubsection{MLinExpr.Clear()}
\label{\detokenize{csapi/mlinexpr:mlinexpr-clear}}\begin{quote}

\sphinxAtStartPar
Clear MLinExpr object.

\sphinxAtStartPar
\sphinxstylestrong{Synopsis}
\begin{quote}

\sphinxAtStartPar
\sphinxcode{\sphinxupquote{void Clear()}}
\end{quote}
\end{quote}

\subsubsection{MLinExpr.Clone()}
\label{\detokenize{csapi/mlinexpr:mlinexpr-clone}}\begin{quote}

\sphinxAtStartPar
Clone MLinExpr object.

\sphinxAtStartPar
\sphinxstylestrong{Synopsis}
\begin{quote}

\sphinxAtStartPar
\sphinxcode{\sphinxupquote{MLinExpr Clone()}}
\end{quote}

\sphinxAtStartPar
\sphinxstylestrong{Return}
\begin{quote}

\sphinxAtStartPar
new MLinExpr object.
\end{quote}
\end{quote}

\subsubsection{MLinExpr.Diagonal()}
\label{\detokenize{csapi/mlinexpr:mlinexpr-diagonal}}\begin{quote}

\sphinxAtStartPar
Get diagonals of MLinExpr object.

\sphinxAtStartPar
\sphinxstylestrong{Synopsis}
\begin{quote}

\sphinxAtStartPar
\sphinxcode{\sphinxupquote{MLinExpr Diagonal(}}
\begin{quote}

\sphinxAtStartPar
\sphinxcode{\sphinxupquote{int offset,}}

\sphinxAtStartPar
\sphinxcode{\sphinxupquote{int axis1,}}

\sphinxAtStartPar
\sphinxcode{\sphinxupquote{int axis2)}}
\end{quote}
\end{quote}

\sphinxAtStartPar
\sphinxstylestrong{Arguments}
\begin{quote}

\sphinxAtStartPar
\sphinxcode{\sphinxupquote{offset}}: offset of the diagonal from the main diagonal. Can be positive or negative.

\sphinxAtStartPar
\sphinxcode{\sphinxupquote{axis1}}: 1st axis of MLinExpr.

\sphinxAtStartPar
\sphinxcode{\sphinxupquote{axis2}}: 2nd axis of MLinExpr.
\end{quote}

\sphinxAtStartPar
\sphinxstylestrong{Return}
\begin{quote}

\sphinxAtStartPar
(N\sphinxhyphen{}1)\sphinxhyphen{}dimensional diagonals.
\end{quote}
\end{quote}

\subsubsection{MLinExpr.Evaluate()}
\label{\detokenize{csapi/mlinexpr:mlinexpr-evaluate}}\begin{quote}

\sphinxAtStartPar
Evaluate MLinExpr object after solving.

\sphinxAtStartPar
\sphinxstylestrong{Synopsis}
\begin{quote}

\sphinxAtStartPar
\sphinxcode{\sphinxupquote{double Evaluate()}}
\end{quote}

\sphinxAtStartPar
\sphinxstylestrong{Return}
\begin{quote}

\sphinxAtStartPar
NdArray object storing value of each linear expression.
\end{quote}
\end{quote}

\subsubsection{MLinExpr.Expand()}
\label{\detokenize{csapi/mlinexpr:mlinexpr-expand}}\begin{quote}

\sphinxAtStartPar
Expand shape of MLinExpr object.

\sphinxAtStartPar
\sphinxstylestrong{Synopsis}
\begin{quote}

\sphinxAtStartPar
\sphinxcode{\sphinxupquote{MLinExpr Expand(int axis)}}
\end{quote}

\sphinxAtStartPar
\sphinxstylestrong{Arguments}
\begin{quote}

\sphinxAtStartPar
\sphinxcode{\sphinxupquote{axis}}: axis of MLinExpr.
\end{quote}

\sphinxAtStartPar
\sphinxstylestrong{Return}
\begin{quote}

\sphinxAtStartPar
MLinExpr object of (N+1)\sphinxhyphen{}dimensional shape.
\end{quote}
\end{quote}

\subsubsection{MLinExpr.Flatten()}
\label{\detokenize{csapi/mlinexpr:mlinexpr-flatten}}\begin{quote}

\sphinxAtStartPar
Flatten a MLinExpr object to a 1\sphinxhyphen{}dimensional shape.

\sphinxAtStartPar
\sphinxstylestrong{Synopsis}
\begin{quote}

\sphinxAtStartPar
\sphinxcode{\sphinxupquote{MLinExpr Flatten()}}
\end{quote}

\sphinxAtStartPar
\sphinxstylestrong{Return}
\begin{quote}

\sphinxAtStartPar
a MLinExpr object collapsed into one dimension.
\end{quote}
\end{quote}

\subsubsection{MLinExpr.GetDim()}
\label{\detokenize{csapi/mlinexpr:mlinexpr-getdim}}\begin{quote}

\sphinxAtStartPar
Get i\sphinxhyphen{}th dimension of MLinExpr object.

\sphinxAtStartPar
\sphinxstylestrong{Synopsis}
\begin{quote}

\sphinxAtStartPar
\sphinxcode{\sphinxupquote{long GetDim(int i)}}
\end{quote}

\sphinxAtStartPar
\sphinxstylestrong{Arguments}
\begin{quote}

\sphinxAtStartPar
\sphinxcode{\sphinxupquote{i}}: index of dimension
\end{quote}

\sphinxAtStartPar
\sphinxstylestrong{Return}
\begin{quote}

\sphinxAtStartPar
i\sphinxhyphen{}th dimension.
\end{quote}
\end{quote}

\subsubsection{MLinExpr.GetItem()}
\label{\detokenize{csapi/mlinexpr:mlinexpr-getitem}}\begin{quote}

\sphinxAtStartPar
Get expression of given index from MLinExpr object.

\sphinxAtStartPar
\sphinxstylestrong{Synopsis}
\begin{quote}

\sphinxAtStartPar
\sphinxcode{\sphinxupquote{MExpression GetItem(long idx)}}
\end{quote}

\sphinxAtStartPar
\sphinxstylestrong{Arguments}
\begin{quote}

\sphinxAtStartPar
\sphinxcode{\sphinxupquote{idx}}: index of expression.
\end{quote}

\sphinxAtStartPar
\sphinxstylestrong{Return}
\begin{quote}

\sphinxAtStartPar
MExpression object.
\end{quote}
\end{quote}

\subsubsection{MLinExpr.GetItem()}
\label{\detokenize{csapi/mlinexpr:id3}}\begin{quote}

\sphinxAtStartPar
Get sub\sphinxhyphen{}arrays of MLinExpr object, given view object.

\sphinxAtStartPar
\sphinxstylestrong{Synopsis}
\begin{quote}

\sphinxAtStartPar
\sphinxcode{\sphinxupquote{MLinExpr GetItem(View view)}}
\end{quote}

\sphinxAtStartPar
\sphinxstylestrong{Arguments}
\begin{quote}

\sphinxAtStartPar
\sphinxcode{\sphinxupquote{view}}: view of multi\sphinxhyphen{}dimensional array.
\end{quote}

\sphinxAtStartPar
\sphinxstylestrong{Return}
\begin{quote}

\sphinxAtStartPar
sub\sphinxhyphen{}arrays of MLinExpr object.
\end{quote}
\end{quote}

\subsubsection{MLinExpr.GetND()}
\label{\detokenize{csapi/mlinexpr:mlinexpr-getnd}}\begin{quote}

\sphinxAtStartPar
Get number of dimensions of MLinExpr object.

\sphinxAtStartPar
\sphinxstylestrong{Synopsis}
\begin{quote}

\sphinxAtStartPar
\sphinxcode{\sphinxupquote{int GetND()}}
\end{quote}

\sphinxAtStartPar
\sphinxstylestrong{Return}
\begin{quote}

\sphinxAtStartPar
number of dimensions.
\end{quote}
\end{quote}

\subsubsection{MLinExpr.GetShape()}
\label{\detokenize{csapi/mlinexpr:mlinexpr-getshape}}\begin{quote}

\sphinxAtStartPar
Get shape of MLinExpr object.

\sphinxAtStartPar
\sphinxstylestrong{Synopsis}
\begin{quote}

\sphinxAtStartPar
\sphinxcode{\sphinxupquote{Shape GetShape()}}
\end{quote}

\sphinxAtStartPar
\sphinxstylestrong{Return}
\begin{quote}

\sphinxAtStartPar
shape object.
\end{quote}
\end{quote}

\subsubsection{MLinExpr.GetSize()}
\label{\detokenize{csapi/mlinexpr:mlinexpr-getsize}}\begin{quote}

\sphinxAtStartPar
Get size of MLinExpr object.

\sphinxAtStartPar
\sphinxstylestrong{Synopsis}
\begin{quote}

\sphinxAtStartPar
\sphinxcode{\sphinxupquote{long GetSize()}}
\end{quote}

\sphinxAtStartPar
\sphinxstylestrong{Return}
\begin{quote}

\sphinxAtStartPar
number of linear expressions.
\end{quote}
\end{quote}

\subsubsection{MLinExpr.HStack\textless{}T\textgreater{}()}
\label{\detokenize{csapi/mlinexpr:mlinexpr-hstack-t}}\begin{quote}

\sphinxAtStartPar
Stack with other NdArray object along horizontal axis.

\sphinxAtStartPar
\sphinxstylestrong{Synopsis}
\begin{quote}

\sphinxAtStartPar
\sphinxcode{\sphinxupquote{MLinExpr HStack\textless{}T\textgreater{}(NdArray\textless{}T\textgreater{} other)}}
\end{quote}

\sphinxAtStartPar
\sphinxstylestrong{Arguments}
\begin{quote}

\sphinxAtStartPar
\sphinxcode{\sphinxupquote{other}}: a NdArray object.
\end{quote}

\sphinxAtStartPar
\sphinxstylestrong{Return}
\begin{quote}

\sphinxAtStartPar
the result MLinExpr object.
\end{quote}
\end{quote}

\subsubsection{MLinExpr.HStack()}
\label{\detokenize{csapi/mlinexpr:mlinexpr-hstack}}\begin{quote}

\sphinxAtStartPar
Stack with other MLinExpr object along horizontal axis.

\sphinxAtStartPar
\sphinxstylestrong{Synopsis}
\begin{quote}

\sphinxAtStartPar
\sphinxcode{\sphinxupquote{MLinExpr HStack(MLinExpr other)}}
\end{quote}

\sphinxAtStartPar
\sphinxstylestrong{Arguments}
\begin{quote}

\sphinxAtStartPar
\sphinxcode{\sphinxupquote{other}}: a MLinExpr object.
\end{quote}

\sphinxAtStartPar
\sphinxstylestrong{Return}
\begin{quote}

\sphinxAtStartPar
the result MLinExpr object.
\end{quote}
\end{quote}

\subsubsection{MLinExpr.HStack()}
\label{\detokenize{csapi/mlinexpr:id4}}\begin{quote}

\sphinxAtStartPar
Stack with other MVar object along horizontal axis.

\sphinxAtStartPar
\sphinxstylestrong{Synopsis}
\begin{quote}

\sphinxAtStartPar
\sphinxcode{\sphinxupquote{MLinExpr HStack(MVar other)}}
\end{quote}

\sphinxAtStartPar
\sphinxstylestrong{Arguments}
\begin{quote}

\sphinxAtStartPar
\sphinxcode{\sphinxupquote{other}}: a MVar object.
\end{quote}

\sphinxAtStartPar
\sphinxstylestrong{Return}
\begin{quote}

\sphinxAtStartPar
the result MLinExpr object.
\end{quote}
\end{quote}

\subsubsection{MLinExpr.Pick()}
\label{\detokenize{csapi/mlinexpr:mlinexpr-pick}}\begin{quote}

\sphinxAtStartPar
Given a list of indexes, get linear expressions from MLinExpr object.

\sphinxAtStartPar
\sphinxstylestrong{Synopsis}
\begin{quote}

\sphinxAtStartPar
\sphinxcode{\sphinxupquote{MLinExpr Pick(NdArray\textless{}int\textgreater{} indexes)}}
\end{quote}

\sphinxAtStartPar
\sphinxstylestrong{Arguments}
\begin{quote}

\sphinxAtStartPar
\sphinxcode{\sphinxupquote{indexes}}: one or two dimensional indexes of elements. If two dimensional, each row is position of an element.
\end{quote}

\sphinxAtStartPar
\sphinxstylestrong{Return}
\begin{quote}

\sphinxAtStartPar
one\sphinxhyphen{}dimensional array of desired linear expressions.
\end{quote}
\end{quote}

\subsubsection{MLinExpr.Repeat()}
\label{\detokenize{csapi/mlinexpr:mlinexpr-repeat}}\begin{quote}

\sphinxAtStartPar
Repeat each element of MLinExpr along given axis.

\sphinxAtStartPar
\sphinxstylestrong{Synopsis}
\begin{quote}

\sphinxAtStartPar
\sphinxcode{\sphinxupquote{MLinExpr Repeat(long repeats, int axis)}}
\end{quote}

\sphinxAtStartPar
\sphinxstylestrong{Arguments}
\begin{quote}

\sphinxAtStartPar
\sphinxcode{\sphinxupquote{repeats}}: number of repetitions for each element.

\sphinxAtStartPar
\sphinxcode{\sphinxupquote{axis}}: axis of MLinExpr.
\end{quote}

\sphinxAtStartPar
\sphinxstylestrong{Return}
\begin{quote}

\sphinxAtStartPar
new MLinExpr object.
\end{quote}
\end{quote}

\subsubsection{MLinExpr.RepeatBlock()}
\label{\detokenize{csapi/mlinexpr:mlinexpr-repeatblock}}\begin{quote}

\sphinxAtStartPar
Repeat an MLinExpr a number of times along given axis.

\sphinxAtStartPar
\sphinxstylestrong{Synopsis}
\begin{quote}

\sphinxAtStartPar
\sphinxcode{\sphinxupquote{MLinExpr RepeatBlock(long repeats, int axis)}}
\end{quote}

\sphinxAtStartPar
\sphinxstylestrong{Arguments}
\begin{quote}

\sphinxAtStartPar
\sphinxcode{\sphinxupquote{repeats}}: number of repetitions.

\sphinxAtStartPar
\sphinxcode{\sphinxupquote{axis}}: axis of MLinExpr.
\end{quote}

\sphinxAtStartPar
\sphinxstylestrong{Return}
\begin{quote}

\sphinxAtStartPar
new MLinExpr object.
\end{quote}
\end{quote}

\subsubsection{MLinExpr.Represent()}
\label{\detokenize{csapi/mlinexpr:mlinexpr-represent}}\begin{quote}

\sphinxAtStartPar
String representation of MLinExpr object.

\sphinxAtStartPar
\sphinxstylestrong{Synopsis}
\begin{quote}

\sphinxAtStartPar
\sphinxcode{\sphinxupquote{string Represent(int maxlen)}}
\end{quote}

\sphinxAtStartPar
\sphinxstylestrong{Arguments}
\begin{quote}

\sphinxAtStartPar
\sphinxcode{\sphinxupquote{maxlen}}: maximum buffer length for representations string.
\end{quote}

\sphinxAtStartPar
\sphinxstylestrong{Return}
\begin{quote}

\sphinxAtStartPar
string object.
\end{quote}
\end{quote}

\subsubsection{MLinExpr.Reshape()}
\label{\detokenize{csapi/mlinexpr:mlinexpr-reshape}}\begin{quote}

\sphinxAtStartPar
Reshape MLinExpr object to new shape.

\sphinxAtStartPar
\sphinxstylestrong{Synopsis}
\begin{quote}

\sphinxAtStartPar
\sphinxcode{\sphinxupquote{MLinExpr Reshape(Shape shp)}}
\end{quote}

\sphinxAtStartPar
\sphinxstylestrong{Arguments}
\begin{quote}

\sphinxAtStartPar
\sphinxcode{\sphinxupquote{shp}}: new shape of M\sphinxhyphen{}dimensions.
\end{quote}

\sphinxAtStartPar
\sphinxstylestrong{Return}
\begin{quote}

\sphinxAtStartPar
M\sphinxhyphen{}dimensional MLinExpr object.
\end{quote}
\end{quote}

\subsubsection{MLinExpr.SetItem()}
\label{\detokenize{csapi/mlinexpr:mlinexpr-setitem}}\begin{quote}

\sphinxAtStartPar
Set expression of given index to MLinExpr object.

\sphinxAtStartPar
\sphinxstylestrong{Synopsis}
\begin{quote}

\sphinxAtStartPar
\sphinxcode{\sphinxupquote{void SetItem(long idx, MExpression expr)}}
\end{quote}

\sphinxAtStartPar
\sphinxstylestrong{Arguments}
\begin{quote}

\sphinxAtStartPar
\sphinxcode{\sphinxupquote{idx}}: index of element.

\sphinxAtStartPar
\sphinxcode{\sphinxupquote{expr}}: MExpression object.
\end{quote}
\end{quote}

\subsubsection{MLinExpr.Squeeze()}
\label{\detokenize{csapi/mlinexpr:mlinexpr-squeeze}}\begin{quote}

\sphinxAtStartPar
Remove axis of length 1 from shape of MLinExpr object.

\sphinxAtStartPar
\sphinxstylestrong{Synopsis}
\begin{quote}

\sphinxAtStartPar
\sphinxcode{\sphinxupquote{MLinExpr Squeeze(int axis)}}
\end{quote}

\sphinxAtStartPar
\sphinxstylestrong{Arguments}
\begin{quote}

\sphinxAtStartPar
\sphinxcode{\sphinxupquote{axis}}: axis of MLinExpr, where the length is 1.
\end{quote}

\sphinxAtStartPar
\sphinxstylestrong{Return}
\begin{quote}

\sphinxAtStartPar
MLinExpr object of (N\sphinxhyphen{}1)\sphinxhyphen{}dimensional shape.
\end{quote}
\end{quote}

\subsubsection{MLinExpr.Stack\textless{}T\textgreater{}()}
\label{\detokenize{csapi/mlinexpr:mlinexpr-stack-t}}\begin{quote}

\sphinxAtStartPar
Stack with other NdArray object along given axis.

\sphinxAtStartPar
\sphinxstylestrong{Synopsis}
\begin{quote}

\sphinxAtStartPar
\sphinxcode{\sphinxupquote{MLinExpr Stack\textless{}T\textgreater{}(NdArray\textless{}T\textgreater{} other, int axis)}}
\end{quote}

\sphinxAtStartPar
\sphinxstylestrong{Arguments}
\begin{quote}

\sphinxAtStartPar
\sphinxcode{\sphinxupquote{other}}: a NdArray object.

\sphinxAtStartPar
\sphinxcode{\sphinxupquote{axis}}: an axis of MLinExpr.
\end{quote}

\sphinxAtStartPar
\sphinxstylestrong{Return}
\begin{quote}

\sphinxAtStartPar
the result MLinExpr object.
\end{quote}
\end{quote}

\subsubsection{MLinExpr.Stack()}
\label{\detokenize{csapi/mlinexpr:mlinexpr-stack}}\begin{quote}

\sphinxAtStartPar
Stack with other MLinExpr object along given axis.

\sphinxAtStartPar
\sphinxstylestrong{Synopsis}
\begin{quote}

\sphinxAtStartPar
\sphinxcode{\sphinxupquote{MLinExpr Stack(MLinExpr other, int axis)}}
\end{quote}

\sphinxAtStartPar
\sphinxstylestrong{Arguments}
\begin{quote}

\sphinxAtStartPar
\sphinxcode{\sphinxupquote{other}}: a MLinExpr object.

\sphinxAtStartPar
\sphinxcode{\sphinxupquote{axis}}: an axis of MLinExpr.
\end{quote}

\sphinxAtStartPar
\sphinxstylestrong{Return}
\begin{quote}

\sphinxAtStartPar
the result MLinExpr object.
\end{quote}
\end{quote}

\subsubsection{MLinExpr.Stack()}
\label{\detokenize{csapi/mlinexpr:id5}}\begin{quote}

\sphinxAtStartPar
Stack with other MLinExpr object along given axis.

\sphinxAtStartPar
\sphinxstylestrong{Synopsis}
\begin{quote}

\sphinxAtStartPar
\sphinxcode{\sphinxupquote{MLinExpr Stack(MVar other, int axis)}}
\end{quote}

\sphinxAtStartPar
\sphinxstylestrong{Arguments}
\begin{quote}

\sphinxAtStartPar
\sphinxcode{\sphinxupquote{other}}: a MVar object.

\sphinxAtStartPar
\sphinxcode{\sphinxupquote{axis}}: an axis of MLinExpr.
\end{quote}

\sphinxAtStartPar
\sphinxstylestrong{Return}
\begin{quote}

\sphinxAtStartPar
the result MLinExpr object.
\end{quote}
\end{quote}

\subsubsection{MLinExpr.SubConstant()}
\label{\detokenize{csapi/mlinexpr:mlinexpr-subconstant}}\begin{quote}

\sphinxAtStartPar
Substract constants from each expression in MLinExpr object.

\sphinxAtStartPar
\sphinxstylestrong{Synopsis}
\begin{quote}

\sphinxAtStartPar
\sphinxcode{\sphinxupquote{void SubConstant(NdArray\textless{}double\textgreater{} constants)}}
\end{quote}

\sphinxAtStartPar
\sphinxstylestrong{Arguments}
\begin{quote}

\sphinxAtStartPar
\sphinxcode{\sphinxupquote{constants}}: N\sphinxhyphen{}dimension NdArray object.
\end{quote}
\end{quote}

\subsubsection{MLinExpr.Sum()}
\label{\detokenize{csapi/mlinexpr:mlinexpr-sum}}\begin{quote}

\sphinxAtStartPar
Sum of all expressions in MLinExpr object.

\sphinxAtStartPar
\sphinxstylestrong{Synopsis}
\begin{quote}

\sphinxAtStartPar
\sphinxcode{\sphinxupquote{MLinExpr Sum()}}
\end{quote}

\sphinxAtStartPar
\sphinxstylestrong{Return}
\begin{quote}

\sphinxAtStartPar
sum in zero dimension.
\end{quote}
\end{quote}

\subsubsection{MLinExpr.Sum()}
\label{\detokenize{csapi/mlinexpr:id6}}\begin{quote}

\sphinxAtStartPar
Sum of variables at given axis of MLinExpr object.

\sphinxAtStartPar
\sphinxstylestrong{Synopsis}
\begin{quote}

\sphinxAtStartPar
\sphinxcode{\sphinxupquote{MLinExpr Sum(int axis)}}
\end{quote}

\sphinxAtStartPar
\sphinxstylestrong{Arguments}
\begin{quote}

\sphinxAtStartPar
\sphinxcode{\sphinxupquote{axis}}: axis of MLinExpr.
\end{quote}

\sphinxAtStartPar
\sphinxstylestrong{Return}
\begin{quote}

\sphinxAtStartPar
MLinExpr object in (N\sphinxhyphen{}1)\sphinxhyphen{}dimension.
\end{quote}
\end{quote}

\subsubsection{MLinExpr.Transpose()}
\label{\detokenize{csapi/mlinexpr:mlinexpr-transpose}}\begin{quote}

\sphinxAtStartPar
Perform matrix transpose of MLinExpr object.

\sphinxAtStartPar
\sphinxstylestrong{Synopsis}
\begin{quote}

\sphinxAtStartPar
\sphinxcode{\sphinxupquote{MLinExpr Transpose()}}
\end{quote}

\sphinxAtStartPar
\sphinxstylestrong{Return}
\begin{quote}

\sphinxAtStartPar
transposed MLinExpr object.
\end{quote}
\end{quote}

\subsubsection{MLinExpr.VStack\textless{}T\textgreater{}()}
\label{\detokenize{csapi/mlinexpr:mlinexpr-vstack-t}}\begin{quote}

\sphinxAtStartPar
Stack with other NdArray object along vertical axis.

\sphinxAtStartPar
\sphinxstylestrong{Synopsis}
\begin{quote}

\sphinxAtStartPar
\sphinxcode{\sphinxupquote{MLinExpr VStack\textless{}T\textgreater{}(NdArray\textless{}T\textgreater{} other)}}
\end{quote}

\sphinxAtStartPar
\sphinxstylestrong{Arguments}
\begin{quote}

\sphinxAtStartPar
\sphinxcode{\sphinxupquote{other}}: a NdArray object.
\end{quote}

\sphinxAtStartPar
\sphinxstylestrong{Return}
\begin{quote}

\sphinxAtStartPar
the result MLinExpr object.
\end{quote}
\end{quote}

\subsubsection{MLinExpr.VStack()}
\label{\detokenize{csapi/mlinexpr:mlinexpr-vstack}}\begin{quote}

\sphinxAtStartPar
Stack with other MLinExpr object along vertical axis.

\sphinxAtStartPar
\sphinxstylestrong{Synopsis}
\begin{quote}

\sphinxAtStartPar
\sphinxcode{\sphinxupquote{MLinExpr VStack(MLinExpr other)}}
\end{quote}

\sphinxAtStartPar
\sphinxstylestrong{Arguments}
\begin{quote}

\sphinxAtStartPar
\sphinxcode{\sphinxupquote{other}}: a MLinExpr object.
\end{quote}

\sphinxAtStartPar
\sphinxstylestrong{Return}
\begin{quote}

\sphinxAtStartPar
the result MLinExpr object.
\end{quote}
\end{quote}

\subsubsection{MLinExpr.VStack()}
\label{\detokenize{csapi/mlinexpr:id7}}\begin{quote}

\sphinxAtStartPar
Stack with other MVar object along vertical axis.

\sphinxAtStartPar
\sphinxstylestrong{Synopsis}
\begin{quote}

\sphinxAtStartPar
\sphinxcode{\sphinxupquote{MLinExpr VStack(MVar other)}}
\end{quote}

\sphinxAtStartPar
\sphinxstylestrong{Arguments}
\begin{quote}

\sphinxAtStartPar
\sphinxcode{\sphinxupquote{other}}: a MVar object.
\end{quote}

\sphinxAtStartPar
\sphinxstylestrong{Return}
\begin{quote}

\sphinxAtStartPar
the result MLinExpr object.
\end{quote}
\end{quote}

\subsection{MPsdConstr}
\label{\detokenize{csharpapiref:mpsdconstr}}\label{\detokenize{csharpapiref:chapcsharpapiref-mpsdconstr}}
\sphinxAtStartPar
The \sphinxtitleref{MPsdConstr} class in COPT represents multi\sphinxhyphen{}dimensional semidefinite constraints.
It is generated through the methods \sphinxcode{\sphinxupquote{addConstrs}} or \sphinxcode{\sphinxupquote{addConstr}} of {\hyperref[\detokenize{csharpapiref:chapcsharpapiref-model}]{\sphinxcrossref{\DUrole{std,std-ref}{Model}}}}.

\sphinxAtStartPar
The following member methods are provided:

\sphinxstepscope

\subsubsection{MPsdConstr.Clone()}
\label{\detokenize{csapi/mpsdconstr:mpsdconstr-clone}}\label{\detokenize{csapi/mpsdconstr::doc}}\begin{quote}

\sphinxAtStartPar
Clone MPsdConstr object.

\sphinxAtStartPar
\sphinxstylestrong{Synopsis}
\begin{quote}

\sphinxAtStartPar
\sphinxcode{\sphinxupquote{MPsdConstr Clone()}}
\end{quote}

\sphinxAtStartPar
\sphinxstylestrong{Return}
\begin{quote}

\sphinxAtStartPar
new MPsdConstr object.
\end{quote}
\end{quote}

\subsubsection{MPsdConstr.Diagonal()}
\label{\detokenize{csapi/mpsdconstr:mpsdconstr-diagonal}}\begin{quote}

\sphinxAtStartPar
Get diagonals of MPsdConstr object.

\sphinxAtStartPar
\sphinxstylestrong{Synopsis}
\begin{quote}

\sphinxAtStartPar
\sphinxcode{\sphinxupquote{MPsdConstr Diagonal(}}
\begin{quote}

\sphinxAtStartPar
\sphinxcode{\sphinxupquote{int offset,}}

\sphinxAtStartPar
\sphinxcode{\sphinxupquote{int axis1,}}

\sphinxAtStartPar
\sphinxcode{\sphinxupquote{int axis2)}}
\end{quote}
\end{quote}

\sphinxAtStartPar
\sphinxstylestrong{Arguments}
\begin{quote}

\sphinxAtStartPar
\sphinxcode{\sphinxupquote{offset}}: offset of the diagonal from the main diagonal. Can be positive or negative.

\sphinxAtStartPar
\sphinxcode{\sphinxupquote{axis1}}: 1st axis of MPsdConstr.

\sphinxAtStartPar
\sphinxcode{\sphinxupquote{axis2}}: 2nd axis of MPsdConstr.
\end{quote}

\sphinxAtStartPar
\sphinxstylestrong{Return}
\begin{quote}

\sphinxAtStartPar
(N\sphinxhyphen{}1)\sphinxhyphen{}dimensional diagonals.
\end{quote}
\end{quote}

\subsubsection{MPsdConstr.Expand()}
\label{\detokenize{csapi/mpsdconstr:mpsdconstr-expand}}\begin{quote}

\sphinxAtStartPar
Expand shape of MPsdConstr object.

\sphinxAtStartPar
\sphinxstylestrong{Synopsis}
\begin{quote}

\sphinxAtStartPar
\sphinxcode{\sphinxupquote{MPsdConstr Expand(int axis)}}
\end{quote}

\sphinxAtStartPar
\sphinxstylestrong{Arguments}
\begin{quote}

\sphinxAtStartPar
\sphinxcode{\sphinxupquote{axis}}: axis of MPsdConstr.
\end{quote}

\sphinxAtStartPar
\sphinxstylestrong{Return}
\begin{quote}

\sphinxAtStartPar
MPsdConstr object of (N+1)\sphinxhyphen{}dimensional shape.
\end{quote}
\end{quote}

\subsubsection{MPsdConstr.Flatten()}
\label{\detokenize{csapi/mpsdconstr:mpsdconstr-flatten}}\begin{quote}

\sphinxAtStartPar
Flatten a MPsdConstr object to a 1\sphinxhyphen{}dimensional shape.

\sphinxAtStartPar
\sphinxstylestrong{Synopsis}
\begin{quote}

\sphinxAtStartPar
\sphinxcode{\sphinxupquote{MPsdConstr Flatten()}}
\end{quote}

\sphinxAtStartPar
\sphinxstylestrong{Return}
\begin{quote}

\sphinxAtStartPar
a MPsdConstr object collapsed into one dimension.
\end{quote}
\end{quote}

\subsubsection{MPsdConstr.Get()}
\label{\detokenize{csapi/mpsdconstr:mpsdconstr-get}}\begin{quote}

\sphinxAtStartPar
Get values of information associated with PSD constraints in MPsdConstr object.

\sphinxAtStartPar
\sphinxstylestrong{Synopsis}
\begin{quote}

\sphinxAtStartPar
\sphinxcode{\sphinxupquote{NdArray\textless{}double\textgreater{} Get(string info)}}
\end{quote}

\sphinxAtStartPar
\sphinxstylestrong{Arguments}
\begin{quote}

\sphinxAtStartPar
\sphinxcode{\sphinxupquote{info}}: name of information.
\end{quote}

\sphinxAtStartPar
\sphinxstylestrong{Return}
\begin{quote}

\sphinxAtStartPar
multi\sphinxhyphen{}dimensional array of information of PSD constraints.
\end{quote}
\end{quote}

\subsubsection{MPsdConstr.GetDim()}
\label{\detokenize{csapi/mpsdconstr:mpsdconstr-getdim}}\begin{quote}

\sphinxAtStartPar
Get i\sphinxhyphen{}th dimension of MPsdConstr object.

\sphinxAtStartPar
\sphinxstylestrong{Synopsis}
\begin{quote}

\sphinxAtStartPar
\sphinxcode{\sphinxupquote{long GetDim(int i)}}
\end{quote}

\sphinxAtStartPar
\sphinxstylestrong{Arguments}
\begin{quote}

\sphinxAtStartPar
\sphinxcode{\sphinxupquote{i}}: index of dimension
\end{quote}

\sphinxAtStartPar
\sphinxstylestrong{Return}
\begin{quote}

\sphinxAtStartPar
i\sphinxhyphen{}th dimension.
\end{quote}
\end{quote}

\subsubsection{MPsdConstr.GetIdx()}
\label{\detokenize{csapi/mpsdconstr:mpsdconstr-getidx}}\begin{quote}

\sphinxAtStartPar
Get index of PSD constraints in MPsdConstr object.

\sphinxAtStartPar
\sphinxstylestrong{Synopsis}
\begin{quote}

\sphinxAtStartPar
\sphinxcode{\sphinxupquote{NdArray\textless{}int\textgreater{} GetIdx()}}
\end{quote}

\sphinxAtStartPar
\sphinxstylestrong{Return}
\begin{quote}

\sphinxAtStartPar
multi\sphinxhyphen{}dimensional array of indexes of PSD constraints.
\end{quote}
\end{quote}

\subsubsection{MPsdConstr.GetItem()}
\label{\detokenize{csapi/mpsdconstr:mpsdconstr-getitem}}\begin{quote}

\sphinxAtStartPar
Get PSD constraint of given index from MPsdConstr object.

\sphinxAtStartPar
\sphinxstylestrong{Synopsis}
\begin{quote}

\sphinxAtStartPar
\sphinxcode{\sphinxupquote{PsdConstraint GetItem(long idx)}}
\end{quote}

\sphinxAtStartPar
\sphinxstylestrong{Arguments}
\begin{quote}

\sphinxAtStartPar
\sphinxcode{\sphinxupquote{idx}}: index of var.
\end{quote}

\sphinxAtStartPar
\sphinxstylestrong{Return}
\begin{quote}

\sphinxAtStartPar
PsdConstraint object.
\end{quote}
\end{quote}

\subsubsection{MPsdConstr.GetItem()}
\label{\detokenize{csapi/mpsdconstr:id1}}\begin{quote}

\sphinxAtStartPar
Get sub\sphinxhyphen{}arrays of MPsdConstr object, given view object.

\sphinxAtStartPar
\sphinxstylestrong{Synopsis}
\begin{quote}

\sphinxAtStartPar
\sphinxcode{\sphinxupquote{MPsdConstr GetItem(View view)}}
\end{quote}

\sphinxAtStartPar
\sphinxstylestrong{Arguments}
\begin{quote}

\sphinxAtStartPar
\sphinxcode{\sphinxupquote{view}}: view of multi\sphinxhyphen{}dimensional array.
\end{quote}

\sphinxAtStartPar
\sphinxstylestrong{Return}
\begin{quote}

\sphinxAtStartPar
sub\sphinxhyphen{}arrays of MPsdConstr object.
\end{quote}
\end{quote}

\subsubsection{MPsdConstr.GetND()}
\label{\detokenize{csapi/mpsdconstr:mpsdconstr-getnd}}\begin{quote}

\sphinxAtStartPar
Get number of dimensions of MPsdConstr object.

\sphinxAtStartPar
\sphinxstylestrong{Synopsis}
\begin{quote}

\sphinxAtStartPar
\sphinxcode{\sphinxupquote{int GetND()}}
\end{quote}

\sphinxAtStartPar
\sphinxstylestrong{Return}
\begin{quote}

\sphinxAtStartPar
number of dimensions.
\end{quote}
\end{quote}

\subsubsection{MPsdConstr.GetShape()}
\label{\detokenize{csapi/mpsdconstr:mpsdconstr-getshape}}\begin{quote}

\sphinxAtStartPar
Get shape of MPsdConstr object.

\sphinxAtStartPar
\sphinxstylestrong{Synopsis}
\begin{quote}

\sphinxAtStartPar
\sphinxcode{\sphinxupquote{Shape GetShape()}}
\end{quote}

\sphinxAtStartPar
\sphinxstylestrong{Return}
\begin{quote}

\sphinxAtStartPar
shape object.
\end{quote}
\end{quote}

\subsubsection{MPsdConstr.GetSize()}
\label{\detokenize{csapi/mpsdconstr:mpsdconstr-getsize}}\begin{quote}

\sphinxAtStartPar
Get size of MPsdConstr object.

\sphinxAtStartPar
\sphinxstylestrong{Synopsis}
\begin{quote}

\sphinxAtStartPar
\sphinxcode{\sphinxupquote{long GetSize()}}
\end{quote}

\sphinxAtStartPar
\sphinxstylestrong{Return}
\begin{quote}

\sphinxAtStartPar
number of vars.
\end{quote}
\end{quote}

\subsubsection{MPsdConstr.HStack()}
\label{\detokenize{csapi/mpsdconstr:mpsdconstr-hstack}}\begin{quote}

\sphinxAtStartPar
Stack with other MPsdConstr object along horizontal axis.

\sphinxAtStartPar
\sphinxstylestrong{Synopsis}
\begin{quote}

\sphinxAtStartPar
\sphinxcode{\sphinxupquote{MPsdConstr HStack(MPsdConstr other)}}
\end{quote}

\sphinxAtStartPar
\sphinxstylestrong{Arguments}
\begin{quote}

\sphinxAtStartPar
\sphinxcode{\sphinxupquote{other}}: a MPsdConstr object.
\end{quote}

\sphinxAtStartPar
\sphinxstylestrong{Return}
\begin{quote}

\sphinxAtStartPar
the result MPsdConstr object.
\end{quote}
\end{quote}

\subsubsection{MPsdConstr.Pick()}
\label{\detokenize{csapi/mpsdconstr:mpsdconstr-pick}}\begin{quote}

\sphinxAtStartPar
Given a list of indexes, get PSD constraints from MPsdConstr object.

\sphinxAtStartPar
\sphinxstylestrong{Synopsis}
\begin{quote}

\sphinxAtStartPar
\sphinxcode{\sphinxupquote{MPsdConstr Pick(NdArray\textless{}int\textgreater{} indexes)}}
\end{quote}

\sphinxAtStartPar
\sphinxstylestrong{Arguments}
\begin{quote}

\sphinxAtStartPar
\sphinxcode{\sphinxupquote{indexes}}: one or two dimensional indexes of elements. If two dimensional, each row is position of an element.
\end{quote}

\sphinxAtStartPar
\sphinxstylestrong{Return}
\begin{quote}

\sphinxAtStartPar
one\sphinxhyphen{}dimensional array of desired PSD constraints.
\end{quote}
\end{quote}

\subsubsection{MPsdConstr.Represent()}
\label{\detokenize{csapi/mpsdconstr:mpsdconstr-represent}}\begin{quote}

\sphinxAtStartPar
String representation of MPsdConstr object.

\sphinxAtStartPar
\sphinxstylestrong{Synopsis}
\begin{quote}

\sphinxAtStartPar
\sphinxcode{\sphinxupquote{string Represent(int maxlen)}}
\end{quote}

\sphinxAtStartPar
\sphinxstylestrong{Arguments}
\begin{quote}

\sphinxAtStartPar
\sphinxcode{\sphinxupquote{maxlen}}: maximum buffer length for representations string.
\end{quote}

\sphinxAtStartPar
\sphinxstylestrong{Return}
\begin{quote}

\sphinxAtStartPar
string object.
\end{quote}
\end{quote}

\subsubsection{MPsdConstr.Reshape()}
\label{\detokenize{csapi/mpsdconstr:mpsdconstr-reshape}}\begin{quote}

\sphinxAtStartPar
Reshape MPsdConstr object to new shape.

\sphinxAtStartPar
\sphinxstylestrong{Synopsis}
\begin{quote}

\sphinxAtStartPar
\sphinxcode{\sphinxupquote{MPsdConstr Reshape(Shape shp)}}
\end{quote}

\sphinxAtStartPar
\sphinxstylestrong{Arguments}
\begin{quote}

\sphinxAtStartPar
\sphinxcode{\sphinxupquote{shp}}: new shape of M\sphinxhyphen{}dimensions.
\end{quote}

\sphinxAtStartPar
\sphinxstylestrong{Return}
\begin{quote}

\sphinxAtStartPar
M\sphinxhyphen{}dimensional MPsdConstr object.
\end{quote}
\end{quote}

\subsubsection{MPsdConstr.Set()}
\label{\detokenize{csapi/mpsdconstr:mpsdconstr-set}}\begin{quote}

\sphinxAtStartPar
Set values of information associated with PSD constraints in MPsdConstr object.

\sphinxAtStartPar
\sphinxstylestrong{Synopsis}
\begin{quote}

\sphinxAtStartPar
\sphinxcode{\sphinxupquote{void Set(string info, double val)}}
\end{quote}

\sphinxAtStartPar
\sphinxstylestrong{Arguments}
\begin{quote}

\sphinxAtStartPar
\sphinxcode{\sphinxupquote{info}}: name of information.

\sphinxAtStartPar
\sphinxcode{\sphinxupquote{val}}: value of information.
\end{quote}
\end{quote}

\subsubsection{MPsdConstr.Set()}
\label{\detokenize{csapi/mpsdconstr:id2}}\begin{quote}

\sphinxAtStartPar
Set values of information associated with PSD constraints in MPsdConstr object.

\sphinxAtStartPar
\sphinxstylestrong{Synopsis}
\begin{quote}

\sphinxAtStartPar
\sphinxcode{\sphinxupquote{void Set(string info, NdArray\textless{}double\textgreater{} vals)}}
\end{quote}

\sphinxAtStartPar
\sphinxstylestrong{Arguments}
\begin{quote}

\sphinxAtStartPar
\sphinxcode{\sphinxupquote{info}}: name of information.

\sphinxAtStartPar
\sphinxcode{\sphinxupquote{vals}}: multi\sphinxhyphen{}dimensional array of values of information.
\end{quote}
\end{quote}

\subsubsection{MPsdConstr.SetItem()}
\label{\detokenize{csapi/mpsdconstr:mpsdconstr-setitem}}\begin{quote}

\sphinxAtStartPar
Set PSD constraint of given index to MPsdConstr object.

\sphinxAtStartPar
\sphinxstylestrong{Synopsis}
\begin{quote}

\sphinxAtStartPar
\sphinxcode{\sphinxupquote{void SetItem(long idx, PsdConstraint constr)}}
\end{quote}

\sphinxAtStartPar
\sphinxstylestrong{Arguments}
\begin{quote}

\sphinxAtStartPar
\sphinxcode{\sphinxupquote{idx}}: index of element.

\sphinxAtStartPar
\sphinxcode{\sphinxupquote{constr}}: PsdConstraint object.
\end{quote}
\end{quote}

\subsubsection{MPsdConstr.Squeeze()}
\label{\detokenize{csapi/mpsdconstr:mpsdconstr-squeeze}}\begin{quote}

\sphinxAtStartPar
Remove axis of length 1 from shape of MPsdConstr object.

\sphinxAtStartPar
\sphinxstylestrong{Synopsis}
\begin{quote}

\sphinxAtStartPar
\sphinxcode{\sphinxupquote{MPsdConstr Squeeze(int axis)}}
\end{quote}

\sphinxAtStartPar
\sphinxstylestrong{Arguments}
\begin{quote}

\sphinxAtStartPar
\sphinxcode{\sphinxupquote{axis}}: axis of MPsdConstr, where the length is 1.
\end{quote}

\sphinxAtStartPar
\sphinxstylestrong{Return}
\begin{quote}

\sphinxAtStartPar
MPsdConstr object of (N\sphinxhyphen{}1)\sphinxhyphen{}dimensional shape.
\end{quote}
\end{quote}

\subsubsection{MPsdConstr.Stack()}
\label{\detokenize{csapi/mpsdconstr:mpsdconstr-stack}}\begin{quote}

\sphinxAtStartPar
Stack with other MPsdConstr object along given axis.

\sphinxAtStartPar
\sphinxstylestrong{Synopsis}
\begin{quote}

\sphinxAtStartPar
\sphinxcode{\sphinxupquote{MPsdConstr Stack(MPsdConstr other, int axis)}}
\end{quote}

\sphinxAtStartPar
\sphinxstylestrong{Arguments}
\begin{quote}

\sphinxAtStartPar
\sphinxcode{\sphinxupquote{other}}: a MPsdConstr object.

\sphinxAtStartPar
\sphinxcode{\sphinxupquote{axis}}: an axis of MPsdConstr.
\end{quote}

\sphinxAtStartPar
\sphinxstylestrong{Return}
\begin{quote}

\sphinxAtStartPar
the result MPsdConstr object.
\end{quote}
\end{quote}

\subsubsection{MPsdConstr.Transpose()}
\label{\detokenize{csapi/mpsdconstr:mpsdconstr-transpose}}\begin{quote}

\sphinxAtStartPar
Perform matrix transpose of MPsdConstr object.

\sphinxAtStartPar
\sphinxstylestrong{Synopsis}
\begin{quote}

\sphinxAtStartPar
\sphinxcode{\sphinxupquote{MPsdConstr Transpose()}}
\end{quote}

\sphinxAtStartPar
\sphinxstylestrong{Return}
\begin{quote}

\sphinxAtStartPar
transposed MPsdConstr object.
\end{quote}
\end{quote}

\subsubsection{MPsdConstr.VStack()}
\label{\detokenize{csapi/mpsdconstr:mpsdconstr-vstack}}\begin{quote}

\sphinxAtStartPar
Stack with other MPsdConstr object along vertical axis.

\sphinxAtStartPar
\sphinxstylestrong{Synopsis}
\begin{quote}

\sphinxAtStartPar
\sphinxcode{\sphinxupquote{MPsdConstr VStack(MPsdConstr other)}}
\end{quote}

\sphinxAtStartPar
\sphinxstylestrong{Arguments}
\begin{quote}

\sphinxAtStartPar
\sphinxcode{\sphinxupquote{other}}: a MPsdConstr object.
\end{quote}

\sphinxAtStartPar
\sphinxstylestrong{Return}
\begin{quote}

\sphinxAtStartPar
the result MPsdConstr object.
\end{quote}
\end{quote}

\subsection{MPsdConstrBuilder}
\label{\detokenize{csharpapiref:mpsdconstrbuilder}}\label{\detokenize{csharpapiref:chapcsharpapiref-mpsdconstrbuilder}}
\sphinxAtStartPar
The \sphinxtitleref{MPsdConstrBuilder} class in COPT serves as a builder for multi\sphinxhyphen{}dimensional
semidefinite constraints. It is used to generate multi\sphinxhyphen{}dimensional semidefinite
constraints and supports operations with the built\sphinxhyphen{}in multi\sphinxhyphen{}dimensional array
{\hyperref[\detokenize{csharpapiref:chapcsharpapiref-ndarray}]{\sphinxcrossref{\DUrole{std,std-ref}{NdArray}}}}.
An \sphinxtitleref{MPsdConstrBuilder} object can be created through comparison operations
between two objects, one of which can be an {\hyperref[\detokenize{csharpapiref:chapcsharpapiref-mpsdexpr}]{\sphinxcrossref{\DUrole{std,std-ref}{MPsdExpr}}}} object.
The following member methods are provided:

\sphinxstepscope

\subsubsection{MPsdConstrBuilder.MPsdConstrBuilder()}
\label{\detokenize{csapi/mpsdconstrbuilder:mpsdconstrbuilder-mpsdconstrbuilder}}\label{\detokenize{csapi/mpsdconstrbuilder::doc}}\begin{quote}

\sphinxAtStartPar
Construct a MPsdConstrBuilder object with the given shape.

\sphinxAtStartPar
\sphinxstylestrong{Synopsis}
\begin{quote}

\sphinxAtStartPar
\sphinxcode{\sphinxupquote{MPsdConstrBuilder(Shape shp)}}
\end{quote}

\sphinxAtStartPar
\sphinxstylestrong{Arguments}
\begin{quote}

\sphinxAtStartPar
\sphinxcode{\sphinxupquote{shp}}: shape of MPsdConstrBuilder.
\end{quote}
\end{quote}

\subsubsection{MPsdConstrBuilder.Flatten()}
\label{\detokenize{csapi/mpsdconstrbuilder:mpsdconstrbuilder-flatten}}\begin{quote}

\sphinxAtStartPar
Flatten a MPsdConstrBuilder object to a 1\sphinxhyphen{}dimensional shape.

\sphinxAtStartPar
\sphinxstylestrong{Synopsis}
\begin{quote}

\sphinxAtStartPar
\sphinxcode{\sphinxupquote{MPsdConstrBuilder Flatten()}}
\end{quote}

\sphinxAtStartPar
\sphinxstylestrong{Return}
\begin{quote}

\sphinxAtStartPar
a MPsdConstrBuilder object collapsed into one dimension.
\end{quote}
\end{quote}

\subsubsection{MPsdConstrBuilder.GetND()}
\label{\detokenize{csapi/mpsdconstrbuilder:mpsdconstrbuilder-getnd}}\begin{quote}

\sphinxAtStartPar
Get number of dimensions of MPsdConstrBuilder object.

\sphinxAtStartPar
\sphinxstylestrong{Synopsis}
\begin{quote}

\sphinxAtStartPar
\sphinxcode{\sphinxupquote{int GetND()}}
\end{quote}

\sphinxAtStartPar
\sphinxstylestrong{Return}
\begin{quote}

\sphinxAtStartPar
number of dimensions.
\end{quote}
\end{quote}

\subsubsection{MPsdConstrBuilder.GetPsdExpr()}
\label{\detokenize{csapi/mpsdconstrbuilder:mpsdconstrbuilder-getpsdexpr}}\begin{quote}

\sphinxAtStartPar
Get N\sphinxhyphen{}dimensional PSD expressions associated with N\sphinxhyphen{}dimensional constraints.

\sphinxAtStartPar
\sphinxstylestrong{Synopsis}
\begin{quote}

\sphinxAtStartPar
\sphinxcode{\sphinxupquote{MPsdExpr GetPsdExpr()}}
\end{quote}

\sphinxAtStartPar
\sphinxstylestrong{Return}
\begin{quote}

\sphinxAtStartPar
MPsdExpr object.
\end{quote}
\end{quote}

\subsubsection{MPsdConstrBuilder.GetRange()}
\label{\detokenize{csapi/mpsdconstrbuilder:mpsdconstrbuilder-getrange}}\begin{quote}

\sphinxAtStartPar
Get range from lower bound to upper bound of N\sphinxhyphen{}dimensional range constraints.

\sphinxAtStartPar
\sphinxstylestrong{Synopsis}
\begin{quote}

\sphinxAtStartPar
\sphinxcode{\sphinxupquote{double GetRange()}}
\end{quote}

\sphinxAtStartPar
\sphinxstylestrong{Return}
\begin{quote}

\sphinxAtStartPar
length from lower bound to upper bound of range constraints.
\end{quote}
\end{quote}

\subsubsection{MPsdConstrBuilder.GetSense()}
\label{\detokenize{csapi/mpsdconstrbuilder:mpsdconstrbuilder-getsense}}\begin{quote}

\sphinxAtStartPar
Get sense associated with N\sphinxhyphen{}dimensional PSD constraints.

\sphinxAtStartPar
\sphinxstylestrong{Synopsis}
\begin{quote}

\sphinxAtStartPar
\sphinxcode{\sphinxupquote{char GetSense()}}
\end{quote}

\sphinxAtStartPar
\sphinxstylestrong{Return}
\begin{quote}

\sphinxAtStartPar
PSD constraint sense.
\end{quote}
\end{quote}

\subsubsection{MPsdConstrBuilder.Set()}
\label{\detokenize{csapi/mpsdconstrbuilder:mpsdconstrbuilder-set}}\begin{quote}

\sphinxAtStartPar
Set N\sphinxhyphen{}dimensional PSD constraints to its builder object.

\sphinxAtStartPar
\sphinxstylestrong{Synopsis}
\begin{quote}

\sphinxAtStartPar
\sphinxcode{\sphinxupquote{void Set(}}
\begin{quote}

\sphinxAtStartPar
\sphinxcode{\sphinxupquote{MPsdExpr expr,}}

\sphinxAtStartPar
\sphinxcode{\sphinxupquote{char sense,}}

\sphinxAtStartPar
\sphinxcode{\sphinxupquote{NdArray\textless{}double\textgreater{} rhs)}}
\end{quote}
\end{quote}

\sphinxAtStartPar
\sphinxstylestrong{Arguments}
\begin{quote}

\sphinxAtStartPar
\sphinxcode{\sphinxupquote{expr}}: MPsdExpr object

\sphinxAtStartPar
\sphinxcode{\sphinxupquote{sense}}: constraint sense other than COPT\_RANGE.

\sphinxAtStartPar
\sphinxcode{\sphinxupquote{rhs}}: N\sphinxhyphen{}dimensional constants at right side of constraints.
\end{quote}
\end{quote}

\subsubsection{MPsdConstrBuilder.Set()}
\label{\detokenize{csapi/mpsdconstrbuilder:id1}}\begin{quote}

\sphinxAtStartPar
Set N\sphinxhyphen{}dimensional PSD constraints to its builder object.

\sphinxAtStartPar
\sphinxstylestrong{Synopsis}
\begin{quote}

\sphinxAtStartPar
\sphinxcode{\sphinxupquote{void Set(}}
\begin{quote}

\sphinxAtStartPar
\sphinxcode{\sphinxupquote{MPsdExpr expr,}}

\sphinxAtStartPar
\sphinxcode{\sphinxupquote{char sense,}}

\sphinxAtStartPar
\sphinxcode{\sphinxupquote{double rhs)}}
\end{quote}
\end{quote}

\sphinxAtStartPar
\sphinxstylestrong{Arguments}
\begin{quote}

\sphinxAtStartPar
\sphinxcode{\sphinxupquote{expr}}: MPsdExpr object

\sphinxAtStartPar
\sphinxcode{\sphinxupquote{sense}}: constraint sense other than COPT\_RANGE.

\sphinxAtStartPar
\sphinxcode{\sphinxupquote{rhs}}: constant of right side of constraints.
\end{quote}
\end{quote}

\subsubsection{MPsdConstrBuilder.Set()}
\label{\detokenize{csapi/mpsdconstrbuilder:id2}}\begin{quote}

\sphinxAtStartPar
Set N\sphinxhyphen{}dimensional PSD constraints to its builder object.

\sphinxAtStartPar
\sphinxstylestrong{Synopsis}
\begin{quote}

\sphinxAtStartPar
\sphinxcode{\sphinxupquote{void Set(}}
\begin{quote}

\sphinxAtStartPar
\sphinxcode{\sphinxupquote{MPsdExpr expr,}}

\sphinxAtStartPar
\sphinxcode{\sphinxupquote{char sense,}}

\sphinxAtStartPar
\sphinxcode{\sphinxupquote{MVar rhs)}}
\end{quote}
\end{quote}

\sphinxAtStartPar
\sphinxstylestrong{Arguments}
\begin{quote}

\sphinxAtStartPar
\sphinxcode{\sphinxupquote{expr}}: MPsdExpr object

\sphinxAtStartPar
\sphinxcode{\sphinxupquote{sense}}: constraint sense other than COPT\_RANGE.

\sphinxAtStartPar
\sphinxcode{\sphinxupquote{rhs}}: MVar object at right side of constraints.
\end{quote}
\end{quote}

\subsubsection{MPsdConstrBuilder.Set()}
\label{\detokenize{csapi/mpsdconstrbuilder:id3}}\begin{quote}

\sphinxAtStartPar
Set N\sphinxhyphen{}dimensional PSD constraints to its builder object.

\sphinxAtStartPar
\sphinxstylestrong{Synopsis}
\begin{quote}

\sphinxAtStartPar
\sphinxcode{\sphinxupquote{void Set(}}
\begin{quote}

\sphinxAtStartPar
\sphinxcode{\sphinxupquote{MPsdExpr expr,}}

\sphinxAtStartPar
\sphinxcode{\sphinxupquote{char sense,}}

\sphinxAtStartPar
\sphinxcode{\sphinxupquote{MLinExpr rhs)}}
\end{quote}
\end{quote}

\sphinxAtStartPar
\sphinxstylestrong{Arguments}
\begin{quote}

\sphinxAtStartPar
\sphinxcode{\sphinxupquote{expr}}: MPsdExpr object

\sphinxAtStartPar
\sphinxcode{\sphinxupquote{sense}}: constraint sense other than COPT\_RANGE.

\sphinxAtStartPar
\sphinxcode{\sphinxupquote{rhs}}: MLinExpr object at right side of constraints.
\end{quote}
\end{quote}

\subsubsection{MPsdConstrBuilder.Set()}
\label{\detokenize{csapi/mpsdconstrbuilder:id4}}\begin{quote}

\sphinxAtStartPar
Set N\sphinxhyphen{}dimensional PSD constraints to its builder object.

\sphinxAtStartPar
\sphinxstylestrong{Synopsis}
\begin{quote}

\sphinxAtStartPar
\sphinxcode{\sphinxupquote{void Set(}}
\begin{quote}

\sphinxAtStartPar
\sphinxcode{\sphinxupquote{MPsdExpr expr,}}

\sphinxAtStartPar
\sphinxcode{\sphinxupquote{char sense,}}

\sphinxAtStartPar
\sphinxcode{\sphinxupquote{MPsdExpr rhs)}}
\end{quote}
\end{quote}

\sphinxAtStartPar
\sphinxstylestrong{Arguments}
\begin{quote}

\sphinxAtStartPar
\sphinxcode{\sphinxupquote{expr}}: MPsdExpr object

\sphinxAtStartPar
\sphinxcode{\sphinxupquote{sense}}: PSD constraint sense other than COPT\_RANGE.

\sphinxAtStartPar
\sphinxcode{\sphinxupquote{rhs}}: MPsdExpr object at right side of PSD constraints.
\end{quote}
\end{quote}

\subsubsection{MPsdConstrBuilder.SetRange()}
\label{\detokenize{csapi/mpsdconstrbuilder:mpsdconstrbuilder-setrange}}\begin{quote}

\sphinxAtStartPar
Set N\sphinxhyphen{}dimensional range PSD constraints to its builder object.

\sphinxAtStartPar
\sphinxstylestrong{Synopsis}
\begin{quote}

\sphinxAtStartPar
\sphinxcode{\sphinxupquote{void SetRange(MPsdExpr expr, double range)}}
\end{quote}

\sphinxAtStartPar
\sphinxstylestrong{Arguments}
\begin{quote}

\sphinxAtStartPar
\sphinxcode{\sphinxupquote{expr}}: MPsdExpr object.

\sphinxAtStartPar
\sphinxcode{\sphinxupquote{range}}: length from lower bound to upper bound of PSD constraint. Must greater than 0.
\end{quote}
\end{quote}

\subsection{MPsdExpr}
\label{\detokenize{csharpapiref:mpsdexpr}}\label{\detokenize{csharpapiref:chapcsharpapiref-mpsdexpr}}
\sphinxAtStartPar
The \sphinxtitleref{MPsdExpr} class in COPT represents multi\sphinxhyphen{}dimensional semidefinite expressions.
It is used to construct multi\sphinxhyphen{}dimensional semidefinite expressions and perform
operations with the built\sphinxhyphen{}in multi\sphinxhyphen{}dimensional array
{\hyperref[\detokenize{csharpapiref:chapcsharpapiref-ndarray}]{\sphinxcrossref{\DUrole{std,std-ref}{NdArray}}}} in COPT.
The elements of {\hyperref[\detokenize{csharpapiref:chapcsharpapiref-mpsdexpr}]{\sphinxcrossref{\DUrole{std,std-ref}{MPsdExpr}}}} are either {\hyperref[\detokenize{csharpapiref:chapcsharpapiref-psdexpr}]{\sphinxcrossref{\DUrole{std,std-ref}{PsdExpr}}}} objects
or their multi\sphinxhyphen{}dimensional linear combinations.
The following member methods are provided:

\sphinxstepscope

\subsubsection{MPsdExpr.AddConstant()}
\label{\detokenize{csapi/mpsdexpr:mpsdexpr-addconstant}}\label{\detokenize{csapi/mpsdexpr::doc}}\begin{quote}

\sphinxAtStartPar
Add constant to each quadratic expression in MPsdExpr object.

\sphinxAtStartPar
\sphinxstylestrong{Synopsis}
\begin{quote}

\sphinxAtStartPar
\sphinxcode{\sphinxupquote{void AddConstant(double constant)}}
\end{quote}

\sphinxAtStartPar
\sphinxstylestrong{Arguments}
\begin{quote}

\sphinxAtStartPar
\sphinxcode{\sphinxupquote{constant}}: the value of the constant.
\end{quote}
\end{quote}

\subsubsection{MPsdExpr.AddConstant()}
\label{\detokenize{csapi/mpsdexpr:id1}}\begin{quote}

\sphinxAtStartPar
Add constants to each PSD expression in MPsdExpr object.

\sphinxAtStartPar
\sphinxstylestrong{Synopsis}
\begin{quote}

\sphinxAtStartPar
\sphinxcode{\sphinxupquote{void AddConstant(NdArray\textless{}double\textgreater{} constants)}}
\end{quote}

\sphinxAtStartPar
\sphinxstylestrong{Arguments}
\begin{quote}

\sphinxAtStartPar
\sphinxcode{\sphinxupquote{constants}}: N\sphinxhyphen{}dimension NdArray object.
\end{quote}
\end{quote}

\subsubsection{MPsdExpr.AddLinExpr()}
\label{\detokenize{csapi/mpsdexpr:mpsdexpr-addlinexpr}}\begin{quote}

\sphinxAtStartPar
Add a linear expression to each PsdExpr in MPsdExpr object.

\sphinxAtStartPar
\sphinxstylestrong{Synopsis}
\begin{quote}

\sphinxAtStartPar
\sphinxcode{\sphinxupquote{void AddLinExpr(Expr expr, double mult)}}
\end{quote}

\sphinxAtStartPar
\sphinxstylestrong{Arguments}
\begin{quote}

\sphinxAtStartPar
\sphinxcode{\sphinxupquote{expr}}: linear expression object.

\sphinxAtStartPar
\sphinxcode{\sphinxupquote{mult}}: the multiplier of linear expression, default value is 1.0.
\end{quote}
\end{quote}

\subsubsection{MPsdExpr.AddMExpr()}
\label{\detokenize{csapi/mpsdexpr:mpsdexpr-addmexpr}}\begin{quote}

\sphinxAtStartPar
Add MExpression to each PSD expression in MPsdExpr object.

\sphinxAtStartPar
\sphinxstylestrong{Synopsis}
\begin{quote}

\sphinxAtStartPar
\sphinxcode{\sphinxupquote{void AddMExpr(MExpression expr, double mult)}}
\end{quote}

\sphinxAtStartPar
\sphinxstylestrong{Arguments}
\begin{quote}

\sphinxAtStartPar
\sphinxcode{\sphinxupquote{expr}}: MExpression object.

\sphinxAtStartPar
\sphinxcode{\sphinxupquote{mult}}: the multiplier of MExpression, default value is 1.0.
\end{quote}
\end{quote}

\subsubsection{MPsdExpr.AddMLinExpr()}
\label{\detokenize{csapi/mpsdexpr:mpsdexpr-addmlinexpr}}\begin{quote}

\sphinxAtStartPar
Add linear expressions to MPsdExpr object.

\sphinxAtStartPar
\sphinxstylestrong{Synopsis}
\begin{quote}

\sphinxAtStartPar
\sphinxcode{\sphinxupquote{void AddMLinExpr(MLinExpr exprs, double mult)}}
\end{quote}

\sphinxAtStartPar
\sphinxstylestrong{Arguments}
\begin{quote}

\sphinxAtStartPar
\sphinxcode{\sphinxupquote{exprs}}: N\sphinxhyphen{}dimension MLinExpr object.

\sphinxAtStartPar
\sphinxcode{\sphinxupquote{mult}}: the same multiplier for added linear expressions, default value is 1.0.
\end{quote}
\end{quote}

\subsubsection{MPsdExpr.AddMPsdExpr()}
\label{\detokenize{csapi/mpsdexpr:mpsdexpr-addmpsdexpr}}\begin{quote}

\sphinxAtStartPar
Add PSD expressions to MPsdExpr object.

\sphinxAtStartPar
\sphinxstylestrong{Synopsis}
\begin{quote}

\sphinxAtStartPar
\sphinxcode{\sphinxupquote{void AddMPsdExpr(MPsdExpr exprs, double mult)}}
\end{quote}

\sphinxAtStartPar
\sphinxstylestrong{Arguments}
\begin{quote}

\sphinxAtStartPar
\sphinxcode{\sphinxupquote{exprs}}: N\sphinxhyphen{}dimension MPsdExpr object.

\sphinxAtStartPar
\sphinxcode{\sphinxupquote{mult}}: the same multiplier for added PSD expressions, default value is 1.0.
\end{quote}
\end{quote}

\subsubsection{MPsdExpr.AddPsdExpr()}
\label{\detokenize{csapi/mpsdexpr:mpsdexpr-addpsdexpr}}\begin{quote}

\sphinxAtStartPar
Add a PSD expression to each PSD expression in MPsdExpr object.

\sphinxAtStartPar
\sphinxstylestrong{Synopsis}
\begin{quote}

\sphinxAtStartPar
\sphinxcode{\sphinxupquote{void AddPsdExpr(PsdExpr expr, double mult)}}
\end{quote}

\sphinxAtStartPar
\sphinxstylestrong{Arguments}
\begin{quote}

\sphinxAtStartPar
\sphinxcode{\sphinxupquote{expr}}: PSD expression object.

\sphinxAtStartPar
\sphinxcode{\sphinxupquote{mult}}: the multiplier of PSD expression, default value is 1.0.
\end{quote}
\end{quote}

\subsubsection{MPsdExpr.AddTerm()}
\label{\detokenize{csapi/mpsdexpr:mpsdexpr-addterm}}\begin{quote}

\sphinxAtStartPar
Add a PSD term to MPsdExpr object.

\sphinxAtStartPar
\sphinxstylestrong{Synopsis}
\begin{quote}

\sphinxAtStartPar
\sphinxcode{\sphinxupquote{void AddTerm(PsdVar var, SymMatExpr expr)}}
\end{quote}

\sphinxAtStartPar
\sphinxstylestrong{Arguments}
\begin{quote}

\sphinxAtStartPar
\sphinxcode{\sphinxupquote{var}}: PSD variable of new PSD term.

\sphinxAtStartPar
\sphinxcode{\sphinxupquote{expr}}: coefficient expression of symmetric matrices of new PSD term.
\end{quote}
\end{quote}

\subsubsection{MPsdExpr.AddTerm()}
\label{\detokenize{csapi/mpsdexpr:id2}}\begin{quote}

\sphinxAtStartPar
Add a linear term to MPsdExpr object.

\sphinxAtStartPar
\sphinxstylestrong{Synopsis}
\begin{quote}

\sphinxAtStartPar
\sphinxcode{\sphinxupquote{void AddTerm(Var var, double coeff)}}
\end{quote}

\sphinxAtStartPar
\sphinxstylestrong{Arguments}
\begin{quote}

\sphinxAtStartPar
\sphinxcode{\sphinxupquote{var}}: variable of new term.

\sphinxAtStartPar
\sphinxcode{\sphinxupquote{coeff}}: coefficient of new term.
\end{quote}
\end{quote}

\subsubsection{MPsdExpr.AddTerm()}
\label{\detokenize{csapi/mpsdexpr:id3}}\begin{quote}

\sphinxAtStartPar
Add a PSD term to MPsdExpr object.

\sphinxAtStartPar
\sphinxstylestrong{Synopsis}
\begin{quote}

\sphinxAtStartPar
\sphinxcode{\sphinxupquote{void AddTerm(PsdVar var, SymMatrix mat)}}
\end{quote}

\sphinxAtStartPar
\sphinxstylestrong{Arguments}
\begin{quote}

\sphinxAtStartPar
\sphinxcode{\sphinxupquote{var}}: PSD variable of new PSD term.

\sphinxAtStartPar
\sphinxcode{\sphinxupquote{mat}}: coefficient matrix of new PSD term.
\end{quote}
\end{quote}

\subsubsection{MPsdExpr.AddTerms()}
\label{\detokenize{csapi/mpsdexpr:mpsdexpr-addterms}}\begin{quote}

\sphinxAtStartPar
Add terms to PSD expressions in MPsdExpr object.

\sphinxAtStartPar
\sphinxstylestrong{Synopsis}
\begin{quote}

\sphinxAtStartPar
\sphinxcode{\sphinxupquote{void AddTerms(MVar vars, double mult)}}
\end{quote}

\sphinxAtStartPar
\sphinxstylestrong{Arguments}
\begin{quote}

\sphinxAtStartPar
\sphinxcode{\sphinxupquote{vars}}: N\sphinxhyphen{}dimension MVar object for added terms.

\sphinxAtStartPar
\sphinxcode{\sphinxupquote{mult}}: the same coefficient for added terms, default value 1.0.
\end{quote}
\end{quote}

\subsubsection{MPsdExpr.AddTerms()}
\label{\detokenize{csapi/mpsdexpr:id4}}\begin{quote}

\sphinxAtStartPar
Add terms to PSD expressions in MPsdExpr object.

\sphinxAtStartPar
\sphinxstylestrong{Synopsis}
\begin{quote}

\sphinxAtStartPar
\sphinxcode{\sphinxupquote{void AddTerms(MVar vars, NdArray\textless{}double\textgreater{} coeffs)}}
\end{quote}

\sphinxAtStartPar
\sphinxstylestrong{Arguments}
\begin{quote}

\sphinxAtStartPar
\sphinxcode{\sphinxupquote{vars}}: N\sphinxhyphen{}dimension MVar object for added terms.

\sphinxAtStartPar
\sphinxcode{\sphinxupquote{coeffs}}: N\sphinxhyphen{}dimension NdArray object of coefficients for added terms.
\end{quote}
\end{quote}

\subsubsection{MPsdExpr.Clear()}
\label{\detokenize{csapi/mpsdexpr:mpsdexpr-clear}}\begin{quote}

\sphinxAtStartPar
Clear MPsdExpr object.

\sphinxAtStartPar
\sphinxstylestrong{Synopsis}
\begin{quote}

\sphinxAtStartPar
\sphinxcode{\sphinxupquote{void Clear()}}
\end{quote}
\end{quote}

\subsubsection{MPsdExpr.Clone()}
\label{\detokenize{csapi/mpsdexpr:mpsdexpr-clone}}\begin{quote}

\sphinxAtStartPar
Clone MPsdExpr object.

\sphinxAtStartPar
\sphinxstylestrong{Synopsis}
\begin{quote}

\sphinxAtStartPar
\sphinxcode{\sphinxupquote{MPsdExpr Clone()}}
\end{quote}

\sphinxAtStartPar
\sphinxstylestrong{Return}
\begin{quote}

\sphinxAtStartPar
new MPsdExpr object.
\end{quote}
\end{quote}

\subsubsection{MPsdExpr.Diagonal()}
\label{\detokenize{csapi/mpsdexpr:mpsdexpr-diagonal}}\begin{quote}

\sphinxAtStartPar
Get diagonals of MPsdExpr object.

\sphinxAtStartPar
\sphinxstylestrong{Synopsis}
\begin{quote}

\sphinxAtStartPar
\sphinxcode{\sphinxupquote{MPsdExpr Diagonal(}}
\begin{quote}

\sphinxAtStartPar
\sphinxcode{\sphinxupquote{int offset,}}

\sphinxAtStartPar
\sphinxcode{\sphinxupquote{int axis1,}}

\sphinxAtStartPar
\sphinxcode{\sphinxupquote{int axis2)}}
\end{quote}
\end{quote}

\sphinxAtStartPar
\sphinxstylestrong{Arguments}
\begin{quote}

\sphinxAtStartPar
\sphinxcode{\sphinxupquote{offset}}: offset of the diagonal from the main diagonal. Can be positive or negative.

\sphinxAtStartPar
\sphinxcode{\sphinxupquote{axis1}}: 1st axis of MPsdExpr.

\sphinxAtStartPar
\sphinxcode{\sphinxupquote{axis2}}: 2nd axis of MPsdExpr.
\end{quote}

\sphinxAtStartPar
\sphinxstylestrong{Return}
\begin{quote}

\sphinxAtStartPar
(N\sphinxhyphen{}1)\sphinxhyphen{}dimensional diagonals.
\end{quote}
\end{quote}

\subsubsection{MPsdExpr.Evaluate()}
\label{\detokenize{csapi/mpsdexpr:mpsdexpr-evaluate}}\begin{quote}

\sphinxAtStartPar
Evaluate MPsdExpr object after solving.

\sphinxAtStartPar
\sphinxstylestrong{Synopsis}
\begin{quote}

\sphinxAtStartPar
\sphinxcode{\sphinxupquote{double Evaluate()}}
\end{quote}

\sphinxAtStartPar
\sphinxstylestrong{Return}
\begin{quote}

\sphinxAtStartPar
NdArray object storing value of each PSD expression.
\end{quote}
\end{quote}

\subsubsection{MPsdExpr.Expand()}
\label{\detokenize{csapi/mpsdexpr:mpsdexpr-expand}}\begin{quote}

\sphinxAtStartPar
Expand shape of MPsdExpr object.

\sphinxAtStartPar
\sphinxstylestrong{Synopsis}
\begin{quote}

\sphinxAtStartPar
\sphinxcode{\sphinxupquote{MPsdExpr Expand(int axis)}}
\end{quote}

\sphinxAtStartPar
\sphinxstylestrong{Arguments}
\begin{quote}

\sphinxAtStartPar
\sphinxcode{\sphinxupquote{axis}}: axis of MPsdExpr.
\end{quote}

\sphinxAtStartPar
\sphinxstylestrong{Return}
\begin{quote}

\sphinxAtStartPar
MPsdExpr object of (N+1)\sphinxhyphen{}dimensional shape.
\end{quote}
\end{quote}

\subsubsection{MPsdExpr.Flatten()}
\label{\detokenize{csapi/mpsdexpr:mpsdexpr-flatten}}\begin{quote}

\sphinxAtStartPar
Flatten a MPsdExpr object to a 1\sphinxhyphen{}dimensional shape.

\sphinxAtStartPar
\sphinxstylestrong{Synopsis}
\begin{quote}

\sphinxAtStartPar
\sphinxcode{\sphinxupquote{MPsdExpr Flatten()}}
\end{quote}

\sphinxAtStartPar
\sphinxstylestrong{Return}
\begin{quote}

\sphinxAtStartPar
a MPsdExpr object collapsed into one dimension.
\end{quote}
\end{quote}

\subsubsection{MPsdExpr.GetDim()}
\label{\detokenize{csapi/mpsdexpr:mpsdexpr-getdim}}\begin{quote}

\sphinxAtStartPar
Get i\sphinxhyphen{}th dimension of MPsdExpr object.

\sphinxAtStartPar
\sphinxstylestrong{Synopsis}
\begin{quote}

\sphinxAtStartPar
\sphinxcode{\sphinxupquote{long GetDim(int i)}}
\end{quote}

\sphinxAtStartPar
\sphinxstylestrong{Arguments}
\begin{quote}

\sphinxAtStartPar
\sphinxcode{\sphinxupquote{i}}: index of dimension
\end{quote}

\sphinxAtStartPar
\sphinxstylestrong{Return}
\begin{quote}

\sphinxAtStartPar
i\sphinxhyphen{}th dimension.
\end{quote}
\end{quote}

\subsubsection{MPsdExpr.GetItem()}
\label{\detokenize{csapi/mpsdexpr:mpsdexpr-getitem}}\begin{quote}

\sphinxAtStartPar
Get PSD expression of given index from MPsdExpr object.

\sphinxAtStartPar
\sphinxstylestrong{Synopsis}
\begin{quote}

\sphinxAtStartPar
\sphinxcode{\sphinxupquote{PsdExpr GetItem(long idx)}}
\end{quote}

\sphinxAtStartPar
\sphinxstylestrong{Arguments}
\begin{quote}

\sphinxAtStartPar
\sphinxcode{\sphinxupquote{idx}}: index of PSD expression.
\end{quote}

\sphinxAtStartPar
\sphinxstylestrong{Return}
\begin{quote}

\sphinxAtStartPar
PSD expression object.
\end{quote}
\end{quote}

\subsubsection{MPsdExpr.GetItem()}
\label{\detokenize{csapi/mpsdexpr:id5}}\begin{quote}

\sphinxAtStartPar
Get sub\sphinxhyphen{}arrays of MPsdExpr object, given view object.

\sphinxAtStartPar
\sphinxstylestrong{Synopsis}
\begin{quote}

\sphinxAtStartPar
\sphinxcode{\sphinxupquote{MPsdExpr GetItem(View view)}}
\end{quote}

\sphinxAtStartPar
\sphinxstylestrong{Arguments}
\begin{quote}

\sphinxAtStartPar
\sphinxcode{\sphinxupquote{view}}: view of multi\sphinxhyphen{}dimensional array.
\end{quote}

\sphinxAtStartPar
\sphinxstylestrong{Return}
\begin{quote}

\sphinxAtStartPar
sub\sphinxhyphen{}arrays of MPsdExpr object.
\end{quote}
\end{quote}

\subsubsection{MPsdExpr.GetND()}
\label{\detokenize{csapi/mpsdexpr:mpsdexpr-getnd}}\begin{quote}

\sphinxAtStartPar
Get number of dimensions of MPsdExpr object.

\sphinxAtStartPar
\sphinxstylestrong{Synopsis}
\begin{quote}

\sphinxAtStartPar
\sphinxcode{\sphinxupquote{int GetND()}}
\end{quote}

\sphinxAtStartPar
\sphinxstylestrong{Return}
\begin{quote}

\sphinxAtStartPar
number of dimensions.
\end{quote}
\end{quote}

\subsubsection{MPsdExpr.GetShape()}
\label{\detokenize{csapi/mpsdexpr:mpsdexpr-getshape}}\begin{quote}

\sphinxAtStartPar
Get shape of MPsdExpr object.

\sphinxAtStartPar
\sphinxstylestrong{Synopsis}
\begin{quote}

\sphinxAtStartPar
\sphinxcode{\sphinxupquote{Shape GetShape()}}
\end{quote}

\sphinxAtStartPar
\sphinxstylestrong{Return}
\begin{quote}

\sphinxAtStartPar
shape object.
\end{quote}
\end{quote}

\subsubsection{MPsdExpr.GetSize()}
\label{\detokenize{csapi/mpsdexpr:mpsdexpr-getsize}}\begin{quote}

\sphinxAtStartPar
Get size of MPsdExpr object.

\sphinxAtStartPar
\sphinxstylestrong{Synopsis}
\begin{quote}

\sphinxAtStartPar
\sphinxcode{\sphinxupquote{long GetSize()}}
\end{quote}

\sphinxAtStartPar
\sphinxstylestrong{Return}
\begin{quote}

\sphinxAtStartPar
number of linear expressions.
\end{quote}
\end{quote}

\subsubsection{MPsdExpr.HStack\textless{}T\textgreater{}()}
\label{\detokenize{csapi/mpsdexpr:mpsdexpr-hstack-t}}\begin{quote}

\sphinxAtStartPar
Stack with other NdArray object along horizontal axis.

\sphinxAtStartPar
\sphinxstylestrong{Synopsis}
\begin{quote}

\sphinxAtStartPar
\sphinxcode{\sphinxupquote{MPsdExpr HStack\textless{}T\textgreater{}(NdArray\textless{}T\textgreater{} other)}}
\end{quote}

\sphinxAtStartPar
\sphinxstylestrong{Arguments}
\begin{quote}

\sphinxAtStartPar
\sphinxcode{\sphinxupquote{other}}: a NdArray object.
\end{quote}

\sphinxAtStartPar
\sphinxstylestrong{Return}
\begin{quote}

\sphinxAtStartPar
the result MPsdExpr object.
\end{quote}
\end{quote}

\subsubsection{MPsdExpr.HStack()}
\label{\detokenize{csapi/mpsdexpr:mpsdexpr-hstack}}\begin{quote}

\sphinxAtStartPar
Stack with other MPsdExpr object along horizontal axis.

\sphinxAtStartPar
\sphinxstylestrong{Synopsis}
\begin{quote}

\sphinxAtStartPar
\sphinxcode{\sphinxupquote{MPsdExpr HStack(MPsdExpr other)}}
\end{quote}

\sphinxAtStartPar
\sphinxstylestrong{Arguments}
\begin{quote}

\sphinxAtStartPar
\sphinxcode{\sphinxupquote{other}}: a MPsdExpr object.
\end{quote}

\sphinxAtStartPar
\sphinxstylestrong{Return}
\begin{quote}

\sphinxAtStartPar
the result MPsdExpr object.
\end{quote}
\end{quote}

\subsubsection{MPsdExpr.HStack()}
\label{\detokenize{csapi/mpsdexpr:id6}}\begin{quote}

\sphinxAtStartPar
Stack with other MLinExpr object along horizontal axis.

\sphinxAtStartPar
\sphinxstylestrong{Synopsis}
\begin{quote}

\sphinxAtStartPar
\sphinxcode{\sphinxupquote{MPsdExpr HStack(MLinExpr other)}}
\end{quote}

\sphinxAtStartPar
\sphinxstylestrong{Arguments}
\begin{quote}

\sphinxAtStartPar
\sphinxcode{\sphinxupquote{other}}: a MLinExpr object.
\end{quote}

\sphinxAtStartPar
\sphinxstylestrong{Return}
\begin{quote}

\sphinxAtStartPar
the result MPsdExpr object.
\end{quote}
\end{quote}

\subsubsection{MPsdExpr.HStack()}
\label{\detokenize{csapi/mpsdexpr:id7}}\begin{quote}

\sphinxAtStartPar
Stack with other MVar object along horizontal axis.

\sphinxAtStartPar
\sphinxstylestrong{Synopsis}
\begin{quote}

\sphinxAtStartPar
\sphinxcode{\sphinxupquote{MPsdExpr HStack(MVar other)}}
\end{quote}

\sphinxAtStartPar
\sphinxstylestrong{Arguments}
\begin{quote}

\sphinxAtStartPar
\sphinxcode{\sphinxupquote{other}}: a MVar object.
\end{quote}

\sphinxAtStartPar
\sphinxstylestrong{Return}
\begin{quote}

\sphinxAtStartPar
the result MPsdExpr object.
\end{quote}
\end{quote}

\subsubsection{MPsdExpr.Pick()}
\label{\detokenize{csapi/mpsdexpr:mpsdexpr-pick}}\begin{quote}

\sphinxAtStartPar
Given a list of indexes, get PSD expressions from MPsdExpr object.

\sphinxAtStartPar
\sphinxstylestrong{Synopsis}
\begin{quote}

\sphinxAtStartPar
\sphinxcode{\sphinxupquote{MPsdExpr Pick(NdArray\textless{}int\textgreater{} indexes)}}
\end{quote}

\sphinxAtStartPar
\sphinxstylestrong{Arguments}
\begin{quote}

\sphinxAtStartPar
\sphinxcode{\sphinxupquote{indexes}}: one or two dimensional indexes of elements. If two dimensional, each row is position of an element.
\end{quote}

\sphinxAtStartPar
\sphinxstylestrong{Return}
\begin{quote}

\sphinxAtStartPar
one\sphinxhyphen{}dimensional array of desired PSD expressions.
\end{quote}
\end{quote}

\subsubsection{MPsdExpr.Repeat()}
\label{\detokenize{csapi/mpsdexpr:mpsdexpr-repeat}}\begin{quote}

\sphinxAtStartPar
Repeat each element of MPsdExpr along given axis.

\sphinxAtStartPar
\sphinxstylestrong{Synopsis}
\begin{quote}

\sphinxAtStartPar
\sphinxcode{\sphinxupquote{MPsdExpr Repeat(long repeats, int axis)}}
\end{quote}

\sphinxAtStartPar
\sphinxstylestrong{Arguments}
\begin{quote}

\sphinxAtStartPar
\sphinxcode{\sphinxupquote{repeats}}: number of repetitions for each element.

\sphinxAtStartPar
\sphinxcode{\sphinxupquote{axis}}: axis of MPsdExpr.
\end{quote}

\sphinxAtStartPar
\sphinxstylestrong{Return}
\begin{quote}

\sphinxAtStartPar
new MPsdExpr object.
\end{quote}
\end{quote}

\subsubsection{MPsdExpr.RepeatBlock()}
\label{\detokenize{csapi/mpsdexpr:mpsdexpr-repeatblock}}\begin{quote}

\sphinxAtStartPar
Repeat an MPsdExpr a number of times along given axis.

\sphinxAtStartPar
\sphinxstylestrong{Synopsis}
\begin{quote}

\sphinxAtStartPar
\sphinxcode{\sphinxupquote{MPsdExpr RepeatBlock(long repeats, int axis)}}
\end{quote}

\sphinxAtStartPar
\sphinxstylestrong{Arguments}
\begin{quote}

\sphinxAtStartPar
\sphinxcode{\sphinxupquote{repeats}}: number of repetitions.

\sphinxAtStartPar
\sphinxcode{\sphinxupquote{axis}}: axis of MPsdExpr.
\end{quote}

\sphinxAtStartPar
\sphinxstylestrong{Return}
\begin{quote}

\sphinxAtStartPar
new MPsdExpr object.
\end{quote}
\end{quote}

\subsubsection{MPsdExpr.Represent()}
\label{\detokenize{csapi/mpsdexpr:mpsdexpr-represent}}\begin{quote}

\sphinxAtStartPar
String representation of MPsdExpr object.

\sphinxAtStartPar
\sphinxstylestrong{Synopsis}
\begin{quote}

\sphinxAtStartPar
\sphinxcode{\sphinxupquote{string Represent(int maxlen)}}
\end{quote}

\sphinxAtStartPar
\sphinxstylestrong{Arguments}
\begin{quote}

\sphinxAtStartPar
\sphinxcode{\sphinxupquote{maxlen}}: maximum buffer length for representations string.
\end{quote}

\sphinxAtStartPar
\sphinxstylestrong{Return}
\begin{quote}

\sphinxAtStartPar
string object.
\end{quote}
\end{quote}

\subsubsection{MPsdExpr.Reshape()}
\label{\detokenize{csapi/mpsdexpr:mpsdexpr-reshape}}\begin{quote}

\sphinxAtStartPar
Reshape MPsdExpr object to new shape.

\sphinxAtStartPar
\sphinxstylestrong{Synopsis}
\begin{quote}

\sphinxAtStartPar
\sphinxcode{\sphinxupquote{MPsdExpr Reshape(Shape shp)}}
\end{quote}

\sphinxAtStartPar
\sphinxstylestrong{Arguments}
\begin{quote}

\sphinxAtStartPar
\sphinxcode{\sphinxupquote{shp}}: new shape of M\sphinxhyphen{}dimensions.
\end{quote}

\sphinxAtStartPar
\sphinxstylestrong{Return}
\begin{quote}

\sphinxAtStartPar
M\sphinxhyphen{}dimensional MPsdExpr object.
\end{quote}
\end{quote}

\subsubsection{MPsdExpr.SetItem()}
\label{\detokenize{csapi/mpsdexpr:mpsdexpr-setitem}}\begin{quote}

\sphinxAtStartPar
Set expression of given index to MPsdExpr object.

\sphinxAtStartPar
\sphinxstylestrong{Synopsis}
\begin{quote}

\sphinxAtStartPar
\sphinxcode{\sphinxupquote{void SetItem(long idx, MExpression expr)}}
\end{quote}

\sphinxAtStartPar
\sphinxstylestrong{Arguments}
\begin{quote}

\sphinxAtStartPar
\sphinxcode{\sphinxupquote{idx}}: index of element.

\sphinxAtStartPar
\sphinxcode{\sphinxupquote{expr}}: MExpression object.
\end{quote}
\end{quote}

\subsubsection{MPsdExpr.SetItem()}
\label{\detokenize{csapi/mpsdexpr:id8}}\begin{quote}

\sphinxAtStartPar
Set PSD expression of given index to MPsdExpr object.

\sphinxAtStartPar
\sphinxstylestrong{Synopsis}
\begin{quote}

\sphinxAtStartPar
\sphinxcode{\sphinxupquote{void SetItem(long idx, PsdExpr expr)}}
\end{quote}

\sphinxAtStartPar
\sphinxstylestrong{Arguments}
\begin{quote}

\sphinxAtStartPar
\sphinxcode{\sphinxupquote{idx}}: index of element.

\sphinxAtStartPar
\sphinxcode{\sphinxupquote{expr}}: PSD expression object.
\end{quote}
\end{quote}

\subsubsection{MPsdExpr.Squeeze()}
\label{\detokenize{csapi/mpsdexpr:mpsdexpr-squeeze}}\begin{quote}

\sphinxAtStartPar
Remove axis of length 1 from shape of MPsdExpr object.

\sphinxAtStartPar
\sphinxstylestrong{Synopsis}
\begin{quote}

\sphinxAtStartPar
\sphinxcode{\sphinxupquote{MPsdExpr Squeeze(int axis)}}
\end{quote}

\sphinxAtStartPar
\sphinxstylestrong{Arguments}
\begin{quote}

\sphinxAtStartPar
\sphinxcode{\sphinxupquote{axis}}: axis of MPsdExpr, where the length is 1.
\end{quote}

\sphinxAtStartPar
\sphinxstylestrong{Return}
\begin{quote}

\sphinxAtStartPar
MPsdExpr object of (N\sphinxhyphen{}1)\sphinxhyphen{}dimensional shape.
\end{quote}
\end{quote}

\subsubsection{MPsdExpr.Stack\textless{}T\textgreater{}()}
\label{\detokenize{csapi/mpsdexpr:mpsdexpr-stack-t}}\begin{quote}

\sphinxAtStartPar
Stack with other NdArray object along given axis.

\sphinxAtStartPar
\sphinxstylestrong{Synopsis}
\begin{quote}

\sphinxAtStartPar
\sphinxcode{\sphinxupquote{MPsdExpr Stack\textless{}T\textgreater{}(NdArray\textless{}T\textgreater{} other, int axis)}}
\end{quote}

\sphinxAtStartPar
\sphinxstylestrong{Arguments}
\begin{quote}

\sphinxAtStartPar
\sphinxcode{\sphinxupquote{other}}: a NdArray object.

\sphinxAtStartPar
\sphinxcode{\sphinxupquote{axis}}: an axis of MPsdExpr.
\end{quote}

\sphinxAtStartPar
\sphinxstylestrong{Return}
\begin{quote}

\sphinxAtStartPar
the result MPsdExpr object.
\end{quote}
\end{quote}

\subsubsection{MPsdExpr.Stack()}
\label{\detokenize{csapi/mpsdexpr:mpsdexpr-stack}}\begin{quote}

\sphinxAtStartPar
Stack with other MPsdExpr object along given axis.

\sphinxAtStartPar
\sphinxstylestrong{Synopsis}
\begin{quote}

\sphinxAtStartPar
\sphinxcode{\sphinxupquote{MPsdExpr Stack(MPsdExpr other, int axis)}}
\end{quote}

\sphinxAtStartPar
\sphinxstylestrong{Arguments}
\begin{quote}

\sphinxAtStartPar
\sphinxcode{\sphinxupquote{other}}: a MPsdExpr object.

\sphinxAtStartPar
\sphinxcode{\sphinxupquote{axis}}: an axis of MPsdExpr.
\end{quote}

\sphinxAtStartPar
\sphinxstylestrong{Return}
\begin{quote}

\sphinxAtStartPar
the result MPsdExpr object.
\end{quote}
\end{quote}

\subsubsection{MPsdExpr.Stack()}
\label{\detokenize{csapi/mpsdexpr:id9}}\begin{quote}

\sphinxAtStartPar
Stack with other MLinExpr object along given axis.

\sphinxAtStartPar
\sphinxstylestrong{Synopsis}
\begin{quote}

\sphinxAtStartPar
\sphinxcode{\sphinxupquote{MPsdExpr Stack(MLinExpr other, int axis)}}
\end{quote}

\sphinxAtStartPar
\sphinxstylestrong{Arguments}
\begin{quote}

\sphinxAtStartPar
\sphinxcode{\sphinxupquote{other}}: a MLinExpr object.

\sphinxAtStartPar
\sphinxcode{\sphinxupquote{axis}}: an axis of MPsdExpr.
\end{quote}

\sphinxAtStartPar
\sphinxstylestrong{Return}
\begin{quote}

\sphinxAtStartPar
the result MPsdExpr object.
\end{quote}
\end{quote}

\subsubsection{MPsdExpr.Stack()}
\label{\detokenize{csapi/mpsdexpr:id10}}\begin{quote}

\sphinxAtStartPar
Stack with other MPsdExpr object along given axis.

\sphinxAtStartPar
\sphinxstylestrong{Synopsis}
\begin{quote}

\sphinxAtStartPar
\sphinxcode{\sphinxupquote{MPsdExpr Stack(MVar other, int axis)}}
\end{quote}

\sphinxAtStartPar
\sphinxstylestrong{Arguments}
\begin{quote}

\sphinxAtStartPar
\sphinxcode{\sphinxupquote{other}}: a MVar object.

\sphinxAtStartPar
\sphinxcode{\sphinxupquote{axis}}: an axis of MPsdExpr.
\end{quote}

\sphinxAtStartPar
\sphinxstylestrong{Return}
\begin{quote}

\sphinxAtStartPar
the result MPsdExpr object.
\end{quote}
\end{quote}

\subsubsection{MPsdExpr.SubConstant()}
\label{\detokenize{csapi/mpsdexpr:mpsdexpr-subconstant}}\begin{quote}

\sphinxAtStartPar
Substract constants from each PSD expression in MPsdExpr object.

\sphinxAtStartPar
\sphinxstylestrong{Synopsis}
\begin{quote}

\sphinxAtStartPar
\sphinxcode{\sphinxupquote{void SubConstant(NdArray\textless{}double\textgreater{} constants)}}
\end{quote}

\sphinxAtStartPar
\sphinxstylestrong{Arguments}
\begin{quote}

\sphinxAtStartPar
\sphinxcode{\sphinxupquote{constants}}: N\sphinxhyphen{}dimension NdArray object.
\end{quote}
\end{quote}

\subsubsection{MPsdExpr.Sum()}
\label{\detokenize{csapi/mpsdexpr:mpsdexpr-sum}}\begin{quote}

\sphinxAtStartPar
Sum of all expressions in MPsdExpr object.

\sphinxAtStartPar
\sphinxstylestrong{Synopsis}
\begin{quote}

\sphinxAtStartPar
\sphinxcode{\sphinxupquote{MPsdExpr Sum()}}
\end{quote}

\sphinxAtStartPar
\sphinxstylestrong{Return}
\begin{quote}

\sphinxAtStartPar
sum in zero dimension.
\end{quote}
\end{quote}

\subsubsection{MPsdExpr.Sum()}
\label{\detokenize{csapi/mpsdexpr:id11}}\begin{quote}

\sphinxAtStartPar
Sum of variables at given axis of MPsdExpr object.

\sphinxAtStartPar
\sphinxstylestrong{Synopsis}
\begin{quote}

\sphinxAtStartPar
\sphinxcode{\sphinxupquote{MPsdExpr Sum(int axis)}}
\end{quote}

\sphinxAtStartPar
\sphinxstylestrong{Arguments}
\begin{quote}

\sphinxAtStartPar
\sphinxcode{\sphinxupquote{axis}}: axis of MPsdExpr.
\end{quote}

\sphinxAtStartPar
\sphinxstylestrong{Return}
\begin{quote}

\sphinxAtStartPar
MPsdExpr object in (N\sphinxhyphen{}1)\sphinxhyphen{}dimension.
\end{quote}
\end{quote}

\subsubsection{MPsdExpr.Transpose()}
\label{\detokenize{csapi/mpsdexpr:mpsdexpr-transpose}}\begin{quote}

\sphinxAtStartPar
Perform matrix transpose of MPsdExpr object.

\sphinxAtStartPar
\sphinxstylestrong{Synopsis}
\begin{quote}

\sphinxAtStartPar
\sphinxcode{\sphinxupquote{MPsdExpr Transpose()}}
\end{quote}

\sphinxAtStartPar
\sphinxstylestrong{Return}
\begin{quote}

\sphinxAtStartPar
transposed MPsdExpr object.
\end{quote}
\end{quote}

\subsubsection{MPsdExpr.VStack\textless{}T\textgreater{}()}
\label{\detokenize{csapi/mpsdexpr:mpsdexpr-vstack-t}}\begin{quote}

\sphinxAtStartPar
Stack with other NdArray object along vertical axis.

\sphinxAtStartPar
\sphinxstylestrong{Synopsis}
\begin{quote}

\sphinxAtStartPar
\sphinxcode{\sphinxupquote{MPsdExpr VStack\textless{}T\textgreater{}(NdArray\textless{}T\textgreater{} other)}}
\end{quote}

\sphinxAtStartPar
\sphinxstylestrong{Arguments}
\begin{quote}

\sphinxAtStartPar
\sphinxcode{\sphinxupquote{other}}: a NdArray object.
\end{quote}

\sphinxAtStartPar
\sphinxstylestrong{Return}
\begin{quote}

\sphinxAtStartPar
the result MPsdExpr object.
\end{quote}
\end{quote}

\subsubsection{MPsdExpr.VStack()}
\label{\detokenize{csapi/mpsdexpr:mpsdexpr-vstack}}\begin{quote}

\sphinxAtStartPar
Stack with other MPsdExpr object along vertical axis.

\sphinxAtStartPar
\sphinxstylestrong{Synopsis}
\begin{quote}

\sphinxAtStartPar
\sphinxcode{\sphinxupquote{MPsdExpr VStack(MPsdExpr other)}}
\end{quote}

\sphinxAtStartPar
\sphinxstylestrong{Arguments}
\begin{quote}

\sphinxAtStartPar
\sphinxcode{\sphinxupquote{other}}: a MPsdExpr object.
\end{quote}

\sphinxAtStartPar
\sphinxstylestrong{Return}
\begin{quote}

\sphinxAtStartPar
the result MPsdExpr object.
\end{quote}
\end{quote}

\subsubsection{MPsdExpr.VStack()}
\label{\detokenize{csapi/mpsdexpr:id12}}\begin{quote}

\sphinxAtStartPar
Stack with other MLinExpr object along vertical axis.

\sphinxAtStartPar
\sphinxstylestrong{Synopsis}
\begin{quote}

\sphinxAtStartPar
\sphinxcode{\sphinxupquote{MPsdExpr VStack(MLinExpr other)}}
\end{quote}

\sphinxAtStartPar
\sphinxstylestrong{Arguments}
\begin{quote}

\sphinxAtStartPar
\sphinxcode{\sphinxupquote{other}}: a MLinExpr object.
\end{quote}

\sphinxAtStartPar
\sphinxstylestrong{Return}
\begin{quote}

\sphinxAtStartPar
the result MPsdExpr object.
\end{quote}
\end{quote}

\subsubsection{MPsdExpr.VStack()}
\label{\detokenize{csapi/mpsdexpr:id13}}\begin{quote}

\sphinxAtStartPar
Stack with other MVar object along vertical axis.

\sphinxAtStartPar
\sphinxstylestrong{Synopsis}
\begin{quote}

\sphinxAtStartPar
\sphinxcode{\sphinxupquote{MPsdExpr VStack(MVar other)}}
\end{quote}

\sphinxAtStartPar
\sphinxstylestrong{Arguments}
\begin{quote}

\sphinxAtStartPar
\sphinxcode{\sphinxupquote{other}}: a MVar object.
\end{quote}

\sphinxAtStartPar
\sphinxstylestrong{Return}
\begin{quote}

\sphinxAtStartPar
the result MPsdExpr object.
\end{quote}
\end{quote}

\subsection{MQConstr}
\label{\detokenize{csharpapiref:mqconstr}}\label{\detokenize{csharpapiref:chapcsharpapiref-mqconstr}}
\sphinxAtStartPar
The MQConstr class is a COPT multi\sphinxhyphen{}dimensional quadratic constraint object.
It can be created by calling the method \sphinxcode{\sphinxupquote{addMQConstr}} of {\hyperref[\detokenize{csharpapiref:chapcsharpapiref-model}]{\sphinxcrossref{\DUrole{std,std-ref}{Model}}}}.
The following methods are provided:

\sphinxstepscope

\subsubsection{MQConstr.Clone()}
\label{\detokenize{csapi/mqconstr:mqconstr-clone}}\label{\detokenize{csapi/mqconstr::doc}}\begin{quote}

\sphinxAtStartPar
Clone MQConstr object.

\sphinxAtStartPar
\sphinxstylestrong{Synopsis}
\begin{quote}

\sphinxAtStartPar
\sphinxcode{\sphinxupquote{MQConstr Clone()}}
\end{quote}

\sphinxAtStartPar
\sphinxstylestrong{Return}
\begin{quote}

\sphinxAtStartPar
new MQConstr object.
\end{quote}
\end{quote}

\subsubsection{MQConstr.Diagonal()}
\label{\detokenize{csapi/mqconstr:mqconstr-diagonal}}\begin{quote}

\sphinxAtStartPar
Get diagonals of MQConstr object.

\sphinxAtStartPar
\sphinxstylestrong{Synopsis}
\begin{quote}

\sphinxAtStartPar
\sphinxcode{\sphinxupquote{MQConstr Diagonal(}}
\begin{quote}

\sphinxAtStartPar
\sphinxcode{\sphinxupquote{int offset,}}

\sphinxAtStartPar
\sphinxcode{\sphinxupquote{int axis1,}}

\sphinxAtStartPar
\sphinxcode{\sphinxupquote{int axis2)}}
\end{quote}
\end{quote}

\sphinxAtStartPar
\sphinxstylestrong{Arguments}
\begin{quote}

\sphinxAtStartPar
\sphinxcode{\sphinxupquote{offset}}: offset of the diagonal from the main diagonal. Can be positive or negative.

\sphinxAtStartPar
\sphinxcode{\sphinxupquote{axis1}}: 1st axis of MQConstr.

\sphinxAtStartPar
\sphinxcode{\sphinxupquote{axis2}}: 2nd axis of MQConstr.
\end{quote}

\sphinxAtStartPar
\sphinxstylestrong{Return}
\begin{quote}

\sphinxAtStartPar
(N\sphinxhyphen{}1)\sphinxhyphen{}dimensional diagonals.
\end{quote}
\end{quote}

\subsubsection{MQConstr.Expand()}
\label{\detokenize{csapi/mqconstr:mqconstr-expand}}\begin{quote}

\sphinxAtStartPar
Expand shape of MQConstr object.

\sphinxAtStartPar
\sphinxstylestrong{Synopsis}
\begin{quote}

\sphinxAtStartPar
\sphinxcode{\sphinxupquote{MQConstr Expand(int axis)}}
\end{quote}

\sphinxAtStartPar
\sphinxstylestrong{Arguments}
\begin{quote}

\sphinxAtStartPar
\sphinxcode{\sphinxupquote{axis}}: axis of MQConstr.
\end{quote}

\sphinxAtStartPar
\sphinxstylestrong{Return}
\begin{quote}

\sphinxAtStartPar
MQConstr object of (N+1)\sphinxhyphen{}dimensional shape.
\end{quote}
\end{quote}

\subsubsection{MQConstr.Flatten()}
\label{\detokenize{csapi/mqconstr:mqconstr-flatten}}\begin{quote}

\sphinxAtStartPar
Flatten a MQConstr object to a 1\sphinxhyphen{}dimensional shape.

\sphinxAtStartPar
\sphinxstylestrong{Synopsis}
\begin{quote}

\sphinxAtStartPar
\sphinxcode{\sphinxupquote{MQConstr Flatten()}}
\end{quote}

\sphinxAtStartPar
\sphinxstylestrong{Return}
\begin{quote}

\sphinxAtStartPar
a MQConstr object collapsed into one dimension.
\end{quote}
\end{quote}

\subsubsection{MQConstr.Get()}
\label{\detokenize{csapi/mqconstr:mqconstr-get}}\begin{quote}

\sphinxAtStartPar
Get values of information associated with quadratic constraints in MQConstr object.

\sphinxAtStartPar
\sphinxstylestrong{Synopsis}
\begin{quote}

\sphinxAtStartPar
\sphinxcode{\sphinxupquote{NdArray\textless{}double\textgreater{} Get(string info)}}
\end{quote}

\sphinxAtStartPar
\sphinxstylestrong{Arguments}
\begin{quote}

\sphinxAtStartPar
\sphinxcode{\sphinxupquote{info}}: name of information.
\end{quote}

\sphinxAtStartPar
\sphinxstylestrong{Return}
\begin{quote}

\sphinxAtStartPar
multi\sphinxhyphen{}dimensional array of information of quadratic constraints.
\end{quote}
\end{quote}

\subsubsection{MQConstr.GetDim()}
\label{\detokenize{csapi/mqconstr:mqconstr-getdim}}\begin{quote}

\sphinxAtStartPar
Get i\sphinxhyphen{}th dimension of MQConstr object.

\sphinxAtStartPar
\sphinxstylestrong{Synopsis}
\begin{quote}

\sphinxAtStartPar
\sphinxcode{\sphinxupquote{long GetDim(int i)}}
\end{quote}

\sphinxAtStartPar
\sphinxstylestrong{Arguments}
\begin{quote}

\sphinxAtStartPar
\sphinxcode{\sphinxupquote{i}}: index of dimension
\end{quote}

\sphinxAtStartPar
\sphinxstylestrong{Return}
\begin{quote}

\sphinxAtStartPar
i\sphinxhyphen{}th dimension.
\end{quote}
\end{quote}

\subsubsection{MQConstr.GetIdx()}
\label{\detokenize{csapi/mqconstr:mqconstr-getidx}}\begin{quote}

\sphinxAtStartPar
Get index of quadratic constraints in MQConstr object.

\sphinxAtStartPar
\sphinxstylestrong{Synopsis}
\begin{quote}

\sphinxAtStartPar
\sphinxcode{\sphinxupquote{NdArray\textless{}int\textgreater{} GetIdx()}}
\end{quote}

\sphinxAtStartPar
\sphinxstylestrong{Return}
\begin{quote}

\sphinxAtStartPar
multi\sphinxhyphen{}dimensional array of indexes of quadratic constraints.
\end{quote}
\end{quote}

\subsubsection{MQConstr.GetItem()}
\label{\detokenize{csapi/mqconstr:mqconstr-getitem}}\begin{quote}

\sphinxAtStartPar
Get quadratic constraint of given index from MQConstr object.

\sphinxAtStartPar
\sphinxstylestrong{Synopsis}
\begin{quote}

\sphinxAtStartPar
\sphinxcode{\sphinxupquote{QConstraint GetItem(long idx)}}
\end{quote}

\sphinxAtStartPar
\sphinxstylestrong{Arguments}
\begin{quote}

\sphinxAtStartPar
\sphinxcode{\sphinxupquote{idx}}: index of quadratic constraint.
\end{quote}

\sphinxAtStartPar
\sphinxstylestrong{Return}
\begin{quote}

\sphinxAtStartPar
QConstraint object.
\end{quote}
\end{quote}

\subsubsection{MQConstr.GetItem()}
\label{\detokenize{csapi/mqconstr:id1}}\begin{quote}

\sphinxAtStartPar
Get sub\sphinxhyphen{}arrays of MQConstr object, given view object.

\sphinxAtStartPar
\sphinxstylestrong{Synopsis}
\begin{quote}

\sphinxAtStartPar
\sphinxcode{\sphinxupquote{MQConstr GetItem(View view)}}
\end{quote}

\sphinxAtStartPar
\sphinxstylestrong{Arguments}
\begin{quote}

\sphinxAtStartPar
\sphinxcode{\sphinxupquote{view}}: view of multi\sphinxhyphen{}dimensional array.
\end{quote}

\sphinxAtStartPar
\sphinxstylestrong{Return}
\begin{quote}

\sphinxAtStartPar
sub\sphinxhyphen{}arrays of MQConstr object.
\end{quote}
\end{quote}

\subsubsection{MQConstr.GetND()}
\label{\detokenize{csapi/mqconstr:mqconstr-getnd}}\begin{quote}

\sphinxAtStartPar
Get number of dimensions of MQConstr object.

\sphinxAtStartPar
\sphinxstylestrong{Synopsis}
\begin{quote}

\sphinxAtStartPar
\sphinxcode{\sphinxupquote{int GetND()}}
\end{quote}

\sphinxAtStartPar
\sphinxstylestrong{Return}
\begin{quote}

\sphinxAtStartPar
number of dimensions.
\end{quote}
\end{quote}

\subsubsection{MQConstr.GetRhs()}
\label{\detokenize{csapi/mqconstr:mqconstr-getrhs}}\begin{quote}

\sphinxAtStartPar
Get RHS of quadratic constraints in MQConstr object.

\sphinxAtStartPar
\sphinxstylestrong{Synopsis}
\begin{quote}

\sphinxAtStartPar
\sphinxcode{\sphinxupquote{NdArray\textless{}double\textgreater{} GetRhs()}}
\end{quote}

\sphinxAtStartPar
\sphinxstylestrong{Return}
\begin{quote}

\sphinxAtStartPar
multi\sphinxhyphen{}dimensional array of RHS of quadratic constraints.
\end{quote}
\end{quote}

\subsubsection{MQConstr.GetSense()}
\label{\detokenize{csapi/mqconstr:mqconstr-getsense}}\begin{quote}

\sphinxAtStartPar
Get senses of quadratic constraints in MQConstr object.

\sphinxAtStartPar
\sphinxstylestrong{Synopsis}
\begin{quote}

\sphinxAtStartPar
\sphinxcode{\sphinxupquote{NdArray\textless{}char\textgreater{} GetSense()}}
\end{quote}

\sphinxAtStartPar
\sphinxstylestrong{Return}
\begin{quote}

\sphinxAtStartPar
multi\sphinxhyphen{}dimensional array of senses of quadratic constraints.
\end{quote}
\end{quote}

\subsubsection{MQConstr.GetShape()}
\label{\detokenize{csapi/mqconstr:mqconstr-getshape}}\begin{quote}

\sphinxAtStartPar
Get shape of MQConstr object.

\sphinxAtStartPar
\sphinxstylestrong{Synopsis}
\begin{quote}

\sphinxAtStartPar
\sphinxcode{\sphinxupquote{Shape GetShape()}}
\end{quote}

\sphinxAtStartPar
\sphinxstylestrong{Return}
\begin{quote}

\sphinxAtStartPar
shape object.
\end{quote}
\end{quote}

\subsubsection{MQConstr.GetSize()}
\label{\detokenize{csapi/mqconstr:mqconstr-getsize}}\begin{quote}

\sphinxAtStartPar
Get size of MQConstr object.

\sphinxAtStartPar
\sphinxstylestrong{Synopsis}
\begin{quote}

\sphinxAtStartPar
\sphinxcode{\sphinxupquote{long GetSize()}}
\end{quote}

\sphinxAtStartPar
\sphinxstylestrong{Return}
\begin{quote}

\sphinxAtStartPar
number of vars.
\end{quote}
\end{quote}

\subsubsection{MQConstr.HStack()}
\label{\detokenize{csapi/mqconstr:mqconstr-hstack}}\begin{quote}

\sphinxAtStartPar
Stack with other MQConstr object along horizontal axis.

\sphinxAtStartPar
\sphinxstylestrong{Synopsis}
\begin{quote}

\sphinxAtStartPar
\sphinxcode{\sphinxupquote{MQConstr HStack(MQConstr other)}}
\end{quote}

\sphinxAtStartPar
\sphinxstylestrong{Arguments}
\begin{quote}

\sphinxAtStartPar
\sphinxcode{\sphinxupquote{other}}: a MQConstr object.
\end{quote}

\sphinxAtStartPar
\sphinxstylestrong{Return}
\begin{quote}

\sphinxAtStartPar
the result MQConstr object.
\end{quote}
\end{quote}

\subsubsection{MQConstr.Pick()}
\label{\detokenize{csapi/mqconstr:mqconstr-pick}}\begin{quote}

\sphinxAtStartPar
Given a list of indexes, get quadratic constraints from MQConstr object.

\sphinxAtStartPar
\sphinxstylestrong{Synopsis}
\begin{quote}

\sphinxAtStartPar
\sphinxcode{\sphinxupquote{MQConstr Pick(NdArray\textless{}int\textgreater{} indexes)}}
\end{quote}

\sphinxAtStartPar
\sphinxstylestrong{Arguments}
\begin{quote}

\sphinxAtStartPar
\sphinxcode{\sphinxupquote{indexes}}: one or two dimensional indexes of elements. If two dimensional, each row is position of an element.
\end{quote}

\sphinxAtStartPar
\sphinxstylestrong{Return}
\begin{quote}

\sphinxAtStartPar
one\sphinxhyphen{}dimensional array of desired quadratic constraints.
\end{quote}
\end{quote}

\subsubsection{MQConstr.Represent()}
\label{\detokenize{csapi/mqconstr:mqconstr-represent}}\begin{quote}

\sphinxAtStartPar
String representation of MQConstr object.

\sphinxAtStartPar
\sphinxstylestrong{Synopsis}
\begin{quote}

\sphinxAtStartPar
\sphinxcode{\sphinxupquote{string Represent(uint maxlen)}}
\end{quote}

\sphinxAtStartPar
\sphinxstylestrong{Arguments}
\begin{quote}

\sphinxAtStartPar
\sphinxcode{\sphinxupquote{maxlen}}: maximum buffer length for representations string.
\end{quote}

\sphinxAtStartPar
\sphinxstylestrong{Return}
\begin{quote}

\sphinxAtStartPar
string object.
\end{quote}
\end{quote}

\subsubsection{MQConstr.Reshape()}
\label{\detokenize{csapi/mqconstr:mqconstr-reshape}}\begin{quote}

\sphinxAtStartPar
Reshape MQConstr object to new shape.

\sphinxAtStartPar
\sphinxstylestrong{Synopsis}
\begin{quote}

\sphinxAtStartPar
\sphinxcode{\sphinxupquote{MQConstr Reshape(Shape shp)}}
\end{quote}

\sphinxAtStartPar
\sphinxstylestrong{Arguments}
\begin{quote}

\sphinxAtStartPar
\sphinxcode{\sphinxupquote{shp}}: new shape of M\sphinxhyphen{}dimensions.
\end{quote}

\sphinxAtStartPar
\sphinxstylestrong{Return}
\begin{quote}

\sphinxAtStartPar
M\sphinxhyphen{}dimensional MQConstr object.
\end{quote}
\end{quote}

\subsubsection{MQConstr.Set()}
\label{\detokenize{csapi/mqconstr:mqconstr-set}}\begin{quote}

\sphinxAtStartPar
Set values of information associated with quadratic constraints in MQConstr object.

\sphinxAtStartPar
\sphinxstylestrong{Synopsis}
\begin{quote}

\sphinxAtStartPar
\sphinxcode{\sphinxupquote{void Set(string info, NdArray\textless{}double\textgreater{} vals)}}
\end{quote}

\sphinxAtStartPar
\sphinxstylestrong{Arguments}
\begin{quote}

\sphinxAtStartPar
\sphinxcode{\sphinxupquote{info}}: name of information.

\sphinxAtStartPar
\sphinxcode{\sphinxupquote{vals}}: multi\sphinxhyphen{}dimensional array of values of information.
\end{quote}
\end{quote}

\subsubsection{MQConstr.Set()}
\label{\detokenize{csapi/mqconstr:id2}}\begin{quote}

\sphinxAtStartPar
Set values of information associated with quadratic constraints in MQConstr object.

\sphinxAtStartPar
\sphinxstylestrong{Synopsis}
\begin{quote}

\sphinxAtStartPar
\sphinxcode{\sphinxupquote{void Set(string info, double val)}}
\end{quote}

\sphinxAtStartPar
\sphinxstylestrong{Arguments}
\begin{quote}

\sphinxAtStartPar
\sphinxcode{\sphinxupquote{info}}: name of information.

\sphinxAtStartPar
\sphinxcode{\sphinxupquote{val}}: value of information.
\end{quote}
\end{quote}

\subsubsection{MQConstr.SetItem()}
\label{\detokenize{csapi/mqconstr:mqconstr-setitem}}\begin{quote}

\sphinxAtStartPar
Set quadratic constraint of given index to MQConstr object.

\sphinxAtStartPar
\sphinxstylestrong{Synopsis}
\begin{quote}

\sphinxAtStartPar
\sphinxcode{\sphinxupquote{void SetItem(long idx, QConstraint constr)}}
\end{quote}

\sphinxAtStartPar
\sphinxstylestrong{Arguments}
\begin{quote}

\sphinxAtStartPar
\sphinxcode{\sphinxupquote{idx}}: index of element.

\sphinxAtStartPar
\sphinxcode{\sphinxupquote{constr}}: quadratic constraint object.
\end{quote}
\end{quote}

\subsubsection{MQConstr.Squeeze()}
\label{\detokenize{csapi/mqconstr:mqconstr-squeeze}}\begin{quote}

\sphinxAtStartPar
Remove axis of length 1 from shape of MQConstr object.

\sphinxAtStartPar
\sphinxstylestrong{Synopsis}
\begin{quote}

\sphinxAtStartPar
\sphinxcode{\sphinxupquote{MQConstr Squeeze(int axis)}}
\end{quote}

\sphinxAtStartPar
\sphinxstylestrong{Arguments}
\begin{quote}

\sphinxAtStartPar
\sphinxcode{\sphinxupquote{axis}}: axis of MQConstr, where the length is 1.
\end{quote}

\sphinxAtStartPar
\sphinxstylestrong{Return}
\begin{quote}

\sphinxAtStartPar
MQConstr object of (N\sphinxhyphen{}1)\sphinxhyphen{}dimensional shape.
\end{quote}
\end{quote}

\subsubsection{MQConstr.Stack()}
\label{\detokenize{csapi/mqconstr:mqconstr-stack}}\begin{quote}

\sphinxAtStartPar
Stack with other MQConstr object along given axis.

\sphinxAtStartPar
\sphinxstylestrong{Synopsis}
\begin{quote}

\sphinxAtStartPar
\sphinxcode{\sphinxupquote{MQConstr Stack(MQConstr other, int axis)}}
\end{quote}

\sphinxAtStartPar
\sphinxstylestrong{Arguments}
\begin{quote}

\sphinxAtStartPar
\sphinxcode{\sphinxupquote{other}}: a MQConstr object.

\sphinxAtStartPar
\sphinxcode{\sphinxupquote{axis}}: an axis of MQConstr.
\end{quote}

\sphinxAtStartPar
\sphinxstylestrong{Return}
\begin{quote}

\sphinxAtStartPar
the result MQConstr object.
\end{quote}
\end{quote}

\subsubsection{MQConstr.Transpose()}
\label{\detokenize{csapi/mqconstr:mqconstr-transpose}}\begin{quote}

\sphinxAtStartPar
Perform matrix transpose of MQConstr object.

\sphinxAtStartPar
\sphinxstylestrong{Synopsis}
\begin{quote}

\sphinxAtStartPar
\sphinxcode{\sphinxupquote{MQConstr Transpose()}}
\end{quote}

\sphinxAtStartPar
\sphinxstylestrong{Return}
\begin{quote}

\sphinxAtStartPar
transposed MQConstr object.
\end{quote}
\end{quote}

\subsubsection{MQConstr.VStack()}
\label{\detokenize{csapi/mqconstr:mqconstr-vstack}}\begin{quote}

\sphinxAtStartPar
Stack with other MQConstr object along vertical axis.

\sphinxAtStartPar
\sphinxstylestrong{Synopsis}
\begin{quote}

\sphinxAtStartPar
\sphinxcode{\sphinxupquote{MQConstr VStack(MQConstr other)}}
\end{quote}

\sphinxAtStartPar
\sphinxstylestrong{Arguments}
\begin{quote}

\sphinxAtStartPar
\sphinxcode{\sphinxupquote{other}}: a MQConstr object.
\end{quote}

\sphinxAtStartPar
\sphinxstylestrong{Return}
\begin{quote}

\sphinxAtStartPar
the result MQConstr object.
\end{quote}
\end{quote}

\subsection{MQConstrBuilder}
\label{\detokenize{csharpapiref:mqconstrbuilder}}\label{\detokenize{csharpapiref:chapcsharpapiref-mqconstrbuilder}}
\sphinxAtStartPar
The MQConstrBuilder class is a COPT builder object of multi\sphinxhyphen{}dimensional
quadratic constraints.  It is used to generate multi\sphinxhyphen{}dimensional quadratic
constraints and supports operations with the built\sphinxhyphen{}in multi\sphinxhyphen{}dimensional array
{\hyperref[\detokenize{csharpapiref:chapcsharpapiref-ndarray}]{\sphinxcrossref{\DUrole{std,std-ref}{NdArray}}}} in COPT. It can be created by comparing two
objects, one of which should be {\hyperref[\detokenize{csharpapiref:chapcsharpapiref-mquadexpr}]{\sphinxcrossref{\DUrole{std,std-ref}{MQuadExpr}}}} object,
by comparison operators. The following methods are provided:

\sphinxstepscope

\subsubsection{MQConstrBuilder.MQConstrBuilder()}
\label{\detokenize{csapi/mqconstrbuilder:mqconstrbuilder-mqconstrbuilder}}\label{\detokenize{csapi/mqconstrbuilder::doc}}\begin{quote}

\sphinxAtStartPar
Construct a MQConstrBuilder object with the given shape.

\sphinxAtStartPar
\sphinxstylestrong{Synopsis}
\begin{quote}

\sphinxAtStartPar
\sphinxcode{\sphinxupquote{MQConstrBuilder(Shape shp)}}
\end{quote}

\sphinxAtStartPar
\sphinxstylestrong{Arguments}
\begin{quote}

\sphinxAtStartPar
\sphinxcode{\sphinxupquote{shp}}: shape of MQConstrBuilder.
\end{quote}
\end{quote}

\subsubsection{MQConstrBuilder.Flatten()}
\label{\detokenize{csapi/mqconstrbuilder:mqconstrbuilder-flatten}}\begin{quote}

\sphinxAtStartPar
Flatten a MQConstrBuilder object to a 1\sphinxhyphen{}dimensional shape.

\sphinxAtStartPar
\sphinxstylestrong{Synopsis}
\begin{quote}

\sphinxAtStartPar
\sphinxcode{\sphinxupquote{MQConstrBuilder Flatten()}}
\end{quote}

\sphinxAtStartPar
\sphinxstylestrong{Return}
\begin{quote}

\sphinxAtStartPar
a MQConstrBuilder object collapsed into one dimension.
\end{quote}
\end{quote}

\subsubsection{MQConstrBuilder.GetND()}
\label{\detokenize{csapi/mqconstrbuilder:mqconstrbuilder-getnd}}\begin{quote}

\sphinxAtStartPar
Get number of dimensions of MQConstrBuilder object.

\sphinxAtStartPar
\sphinxstylestrong{Synopsis}
\begin{quote}

\sphinxAtStartPar
\sphinxcode{\sphinxupquote{int GetND()}}
\end{quote}

\sphinxAtStartPar
\sphinxstylestrong{Return}
\begin{quote}

\sphinxAtStartPar
number of dimensions.
\end{quote}
\end{quote}

\subsubsection{MQConstrBuilder.GetQuadExpr()}
\label{\detokenize{csapi/mqconstrbuilder:mqconstrbuilder-getquadexpr}}\begin{quote}

\sphinxAtStartPar
Get N\sphinxhyphen{}dimensional quadratic expressions associated with N\sphinxhyphen{}dimensional quadratic constraints.

\sphinxAtStartPar
\sphinxstylestrong{Synopsis}
\begin{quote}

\sphinxAtStartPar
\sphinxcode{\sphinxupquote{MQuadExpr GetQuadExpr()}}
\end{quote}

\sphinxAtStartPar
\sphinxstylestrong{Return}
\begin{quote}

\sphinxAtStartPar
MQuadExpr object.
\end{quote}
\end{quote}

\subsubsection{MQConstrBuilder.GetSense()}
\label{\detokenize{csapi/mqconstrbuilder:mqconstrbuilder-getsense}}\begin{quote}

\sphinxAtStartPar
Get sense associated with N\sphinxhyphen{}dimensional constraints.

\sphinxAtStartPar
\sphinxstylestrong{Synopsis}
\begin{quote}

\sphinxAtStartPar
\sphinxcode{\sphinxupquote{char GetSense()}}
\end{quote}

\sphinxAtStartPar
\sphinxstylestrong{Return}
\begin{quote}

\sphinxAtStartPar
constraint sense.
\end{quote}
\end{quote}

\subsubsection{MQConstrBuilder.Set()}
\label{\detokenize{csapi/mqconstrbuilder:mqconstrbuilder-set}}\begin{quote}

\sphinxAtStartPar
Set N\sphinxhyphen{}dimensional quadratic constraints to its builder object.

\sphinxAtStartPar
\sphinxstylestrong{Synopsis}
\begin{quote}

\sphinxAtStartPar
\sphinxcode{\sphinxupquote{void Set(}}
\begin{quote}

\sphinxAtStartPar
\sphinxcode{\sphinxupquote{MQuadExpr expr,}}

\sphinxAtStartPar
\sphinxcode{\sphinxupquote{char sense,}}

\sphinxAtStartPar
\sphinxcode{\sphinxupquote{double rhs)}}
\end{quote}
\end{quote}

\sphinxAtStartPar
\sphinxstylestrong{Arguments}
\begin{quote}

\sphinxAtStartPar
\sphinxcode{\sphinxupquote{expr}}: MQuadExpr object

\sphinxAtStartPar
\sphinxcode{\sphinxupquote{sense}}: constraint sense other than COPT\_RANGE.

\sphinxAtStartPar
\sphinxcode{\sphinxupquote{rhs}}: constant of right side of constraints.
\end{quote}
\end{quote}

\subsubsection{MQConstrBuilder.Set()}
\label{\detokenize{csapi/mqconstrbuilder:id1}}\begin{quote}

\sphinxAtStartPar
Set N\sphinxhyphen{}dimensional quadratic constraints to its builder object.

\sphinxAtStartPar
\sphinxstylestrong{Synopsis}
\begin{quote}

\sphinxAtStartPar
\sphinxcode{\sphinxupquote{void Set(}}
\begin{quote}

\sphinxAtStartPar
\sphinxcode{\sphinxupquote{MQuadExpr expr,}}

\sphinxAtStartPar
\sphinxcode{\sphinxupquote{char sense,}}

\sphinxAtStartPar
\sphinxcode{\sphinxupquote{MVar rhs)}}
\end{quote}
\end{quote}

\sphinxAtStartPar
\sphinxstylestrong{Arguments}
\begin{quote}

\sphinxAtStartPar
\sphinxcode{\sphinxupquote{expr}}: MQuadExpr object

\sphinxAtStartPar
\sphinxcode{\sphinxupquote{sense}}: constraint sense other than COPT\_RANGE.

\sphinxAtStartPar
\sphinxcode{\sphinxupquote{rhs}}: MVar object at right side of quadratic constraints.
\end{quote}
\end{quote}

\subsubsection{MQConstrBuilder.Set()}
\label{\detokenize{csapi/mqconstrbuilder:id2}}\begin{quote}

\sphinxAtStartPar
Set N\sphinxhyphen{}dimensional quadratic constraints to its builder object.

\sphinxAtStartPar
\sphinxstylestrong{Synopsis}
\begin{quote}

\sphinxAtStartPar
\sphinxcode{\sphinxupquote{void Set(}}
\begin{quote}

\sphinxAtStartPar
\sphinxcode{\sphinxupquote{MQuadExpr expr,}}

\sphinxAtStartPar
\sphinxcode{\sphinxupquote{char sense,}}

\sphinxAtStartPar
\sphinxcode{\sphinxupquote{MLinExpr rhs)}}
\end{quote}
\end{quote}

\sphinxAtStartPar
\sphinxstylestrong{Arguments}
\begin{quote}

\sphinxAtStartPar
\sphinxcode{\sphinxupquote{expr}}: MQuadExpr object

\sphinxAtStartPar
\sphinxcode{\sphinxupquote{sense}}: constraint sense other than COPT\_RANGE.

\sphinxAtStartPar
\sphinxcode{\sphinxupquote{rhs}}: MLinExpr object at right side of quadratic constraints.
\end{quote}
\end{quote}

\subsubsection{MQConstrBuilder.Set()}
\label{\detokenize{csapi/mqconstrbuilder:id3}}\begin{quote}

\sphinxAtStartPar
Set N\sphinxhyphen{}dimensional quadratic constraints to its builder object.

\sphinxAtStartPar
\sphinxstylestrong{Synopsis}
\begin{quote}

\sphinxAtStartPar
\sphinxcode{\sphinxupquote{void Set(}}
\begin{quote}

\sphinxAtStartPar
\sphinxcode{\sphinxupquote{MQuadExpr expr,}}

\sphinxAtStartPar
\sphinxcode{\sphinxupquote{char sense,}}

\sphinxAtStartPar
\sphinxcode{\sphinxupquote{MQuadExpr rhs)}}
\end{quote}
\end{quote}

\sphinxAtStartPar
\sphinxstylestrong{Arguments}
\begin{quote}

\sphinxAtStartPar
\sphinxcode{\sphinxupquote{expr}}: MQuadExpr object

\sphinxAtStartPar
\sphinxcode{\sphinxupquote{sense}}: constraint sense other than COPT\_RANGE.

\sphinxAtStartPar
\sphinxcode{\sphinxupquote{rhs}}: MQuadExpr object at right side of quadratic constraints.
\end{quote}
\end{quote}

\subsubsection{MQConstrBuilder.Set()}
\label{\detokenize{csapi/mqconstrbuilder:id4}}\begin{quote}

\sphinxAtStartPar
Set N\sphinxhyphen{}dimensional quadratic constraints to its builder object.

\sphinxAtStartPar
\sphinxstylestrong{Synopsis}
\begin{quote}

\sphinxAtStartPar
\sphinxcode{\sphinxupquote{void Set(}}
\begin{quote}

\sphinxAtStartPar
\sphinxcode{\sphinxupquote{MQuadExpr expr,}}

\sphinxAtStartPar
\sphinxcode{\sphinxupquote{char sense,}}

\sphinxAtStartPar
\sphinxcode{\sphinxupquote{NdArray\textless{}double\textgreater{} rhs)}}
\end{quote}
\end{quote}

\sphinxAtStartPar
\sphinxstylestrong{Arguments}
\begin{quote}

\sphinxAtStartPar
\sphinxcode{\sphinxupquote{expr}}: MQuadExpr object

\sphinxAtStartPar
\sphinxcode{\sphinxupquote{sense}}: constraint sense other than COPT\_RANGE.

\sphinxAtStartPar
\sphinxcode{\sphinxupquote{rhs}}: N\sphinxhyphen{}dimensional constants at right side of quadratic constraints.
\end{quote}
\end{quote}

\subsection{MQExpression}
\label{\detokenize{csharpapiref:mqexpression}}\label{\detokenize{csharpapiref:chapcsharpapiref-mqexpression}}
\sphinxAtStartPar
The MQExpression class is a generalized version of {\hyperref[\detokenize{csharpapiref:chapcsharpapiref-quadexpr}]{\sphinxcrossref{\DUrole{std,std-ref}{QuadExpr}}}}.
It represents a quadratic expression and supports most of methods in QuadExpr class.
In addition, it supports quadratic combination of multi\sphinxhyphen{}dimensional objects,
such as {\hyperref[\detokenize{csharpapiref:chapcsharpapiref-mvar}]{\sphinxcrossref{\DUrole{std,std-ref}{MVar}}}} object and {\hyperref[\detokenize{csharpapiref:chapcsharpapiref-mlinexpr}]{\sphinxcrossref{\DUrole{std,std-ref}{MLinExpr}}}} object.
The following methods are provided:

\sphinxstepscope

\subsubsection{MQExpression.MQExpression()}
\label{\detokenize{csapi/mqexpression:mqexpression-mqexpression}}\label{\detokenize{csapi/mqexpression::doc}}\begin{quote}

\sphinxAtStartPar
Construct a MQExpression object with the given constant.

\sphinxAtStartPar
\sphinxstylestrong{Synopsis}
\begin{quote}

\sphinxAtStartPar
\sphinxcode{\sphinxupquote{MQExpression(double constant)}}
\end{quote}

\sphinxAtStartPar
\sphinxstylestrong{Arguments}
\begin{quote}

\sphinxAtStartPar
\sphinxcode{\sphinxupquote{constant}}: constant number.
\end{quote}
\end{quote}

\subsubsection{MQExpression.MQExpression()}
\label{\detokenize{csapi/mqexpression:id1}}\begin{quote}

\sphinxAtStartPar
Construct a MQExpression object with the given quadratic expression.

\sphinxAtStartPar
\sphinxstylestrong{Synopsis}
\begin{quote}

\sphinxAtStartPar
\sphinxcode{\sphinxupquote{MQExpression(QuadExpr expr)}}
\end{quote}

\sphinxAtStartPar
\sphinxstylestrong{Arguments}
\begin{quote}

\sphinxAtStartPar
\sphinxcode{\sphinxupquote{expr}}: a quadratic expression.
\end{quote}
\end{quote}

\subsubsection{MQExpression.AddConstant()}
\label{\detokenize{csapi/mqexpression:mqexpression-addconstant}}\begin{quote}

\sphinxAtStartPar
Add constant for the expression.

\sphinxAtStartPar
\sphinxstylestrong{Synopsis}
\begin{quote}

\sphinxAtStartPar
\sphinxcode{\sphinxupquote{void AddConstant(double constant)}}
\end{quote}

\sphinxAtStartPar
\sphinxstylestrong{Arguments}
\begin{quote}

\sphinxAtStartPar
\sphinxcode{\sphinxupquote{constant}}: the value of the constant.
\end{quote}
\end{quote}

\subsubsection{MQExpression.AddExpr()}
\label{\detokenize{csapi/mqexpression:mqexpression-addexpr}}\begin{quote}

\sphinxAtStartPar
Add a linear expression to MQExpression object.

\sphinxAtStartPar
\sphinxstylestrong{Synopsis}
\begin{quote}

\sphinxAtStartPar
\sphinxcode{\sphinxupquote{void AddExpr(Expr expr, double mult)}}
\end{quote}

\sphinxAtStartPar
\sphinxstylestrong{Arguments}
\begin{quote}

\sphinxAtStartPar
\sphinxcode{\sphinxupquote{expr}}: linear expression object.

\sphinxAtStartPar
\sphinxcode{\sphinxupquote{mult}}: the multiplier of linear expression, default value is 1.0.
\end{quote}
\end{quote}

\subsubsection{MQExpression.AddMExpr()}
\label{\detokenize{csapi/mqexpression:mqexpression-addmexpr}}\begin{quote}

\sphinxAtStartPar
Add MExpression to MQExpression object.

\sphinxAtStartPar
\sphinxstylestrong{Synopsis}
\begin{quote}

\sphinxAtStartPar
\sphinxcode{\sphinxupquote{void AddMExpr(MExpression expr, double mult)}}
\end{quote}

\sphinxAtStartPar
\sphinxstylestrong{Arguments}
\begin{quote}

\sphinxAtStartPar
\sphinxcode{\sphinxupquote{expr}}: MExpression object.

\sphinxAtStartPar
\sphinxcode{\sphinxupquote{mult}}: the multiplier of MExpression, default value is 1.0.
\end{quote}
\end{quote}

\subsubsection{MQExpression.AddMQExpr()}
\label{\detokenize{csapi/mqexpression:mqexpression-addmqexpr}}\begin{quote}

\sphinxAtStartPar
Add MQExpression to MQExpression object.

\sphinxAtStartPar
\sphinxstylestrong{Synopsis}
\begin{quote}

\sphinxAtStartPar
\sphinxcode{\sphinxupquote{void AddMQExpr(MQExpression expr, double mult)}}
\end{quote}

\sphinxAtStartPar
\sphinxstylestrong{Arguments}
\begin{quote}

\sphinxAtStartPar
\sphinxcode{\sphinxupquote{expr}}: MQExpression object.

\sphinxAtStartPar
\sphinxcode{\sphinxupquote{mult}}: the multiplier of MQExpression, default value is 1.0.
\end{quote}
\end{quote}

\subsubsection{MQExpression.AddQuadExpr()}
\label{\detokenize{csapi/mqexpression:mqexpression-addquadexpr}}\begin{quote}

\sphinxAtStartPar
Add a quadratic expression to MQExpression object.

\sphinxAtStartPar
\sphinxstylestrong{Synopsis}
\begin{quote}

\sphinxAtStartPar
\sphinxcode{\sphinxupquote{void AddQuadExpr(MExpression left, Expr right)}}
\end{quote}

\sphinxAtStartPar
\sphinxstylestrong{Arguments}
\begin{quote}

\sphinxAtStartPar
\sphinxcode{\sphinxupquote{left}}: MExpression object.

\sphinxAtStartPar
\sphinxcode{\sphinxupquote{right}}: Expr object.
\end{quote}
\end{quote}

\subsubsection{MQExpression.AddQuadExpr()}
\label{\detokenize{csapi/mqexpression:id2}}\begin{quote}

\sphinxAtStartPar
Add a quadratic expression to MQExpression object.

\sphinxAtStartPar
\sphinxstylestrong{Synopsis}
\begin{quote}

\sphinxAtStartPar
\sphinxcode{\sphinxupquote{void AddQuadExpr(MExpression left, MExpression right)}}
\end{quote}

\sphinxAtStartPar
\sphinxstylestrong{Arguments}
\begin{quote}

\sphinxAtStartPar
\sphinxcode{\sphinxupquote{left}}: left MExpression object.

\sphinxAtStartPar
\sphinxcode{\sphinxupquote{right}}: right MExpression object.
\end{quote}
\end{quote}

\subsubsection{MQExpression.AddQuadExpr()}
\label{\detokenize{csapi/mqexpression:id3}}\begin{quote}

\sphinxAtStartPar
Add a quadratic expression to MQExpression object.

\sphinxAtStartPar
\sphinxstylestrong{Synopsis}
\begin{quote}

\sphinxAtStartPar
\sphinxcode{\sphinxupquote{void AddQuadExpr(QuadExpr expr, double mult)}}
\end{quote}

\sphinxAtStartPar
\sphinxstylestrong{Arguments}
\begin{quote}

\sphinxAtStartPar
\sphinxcode{\sphinxupquote{expr}}: quadratic expression object.

\sphinxAtStartPar
\sphinxcode{\sphinxupquote{mult}}: the multiplier of quadratic expression, default value is 1.0.
\end{quote}
\end{quote}

\subsubsection{MQExpression.AddQuadExpr()}
\label{\detokenize{csapi/mqexpression:id4}}\begin{quote}

\sphinxAtStartPar
Add a quadratic expression to MQExpression object.

\sphinxAtStartPar
\sphinxstylestrong{Synopsis}
\begin{quote}

\sphinxAtStartPar
\sphinxcode{\sphinxupquote{void AddQuadExpr(MExpression expr, Var var)}}
\end{quote}

\sphinxAtStartPar
\sphinxstylestrong{Arguments}
\begin{quote}

\sphinxAtStartPar
\sphinxcode{\sphinxupquote{expr}}: MExpression object.

\sphinxAtStartPar
\sphinxcode{\sphinxupquote{var}}: Var object.
\end{quote}
\end{quote}

\subsubsection{MQExpression.AddTerm()}
\label{\detokenize{csapi/mqexpression:mqexpression-addterm}}\begin{quote}

\sphinxAtStartPar
Add a linear term to MQExpression object.

\sphinxAtStartPar
\sphinxstylestrong{Synopsis}
\begin{quote}

\sphinxAtStartPar
\sphinxcode{\sphinxupquote{void AddTerm(Var var, double coeff)}}
\end{quote}

\sphinxAtStartPar
\sphinxstylestrong{Arguments}
\begin{quote}

\sphinxAtStartPar
\sphinxcode{\sphinxupquote{var}}: variable of new term.

\sphinxAtStartPar
\sphinxcode{\sphinxupquote{coeff}}: coefficient of new term.
\end{quote}
\end{quote}

\subsubsection{MQExpression.AddTerm()}
\label{\detokenize{csapi/mqexpression:id5}}\begin{quote}

\sphinxAtStartPar
Add a quadratic term to MQExpression object.

\sphinxAtStartPar
\sphinxstylestrong{Synopsis}
\begin{quote}

\sphinxAtStartPar
\sphinxcode{\sphinxupquote{void AddTerm(}}
\begin{quote}

\sphinxAtStartPar
\sphinxcode{\sphinxupquote{Var var1,}}

\sphinxAtStartPar
\sphinxcode{\sphinxupquote{Var var2,}}

\sphinxAtStartPar
\sphinxcode{\sphinxupquote{double coeff)}}
\end{quote}
\end{quote}

\sphinxAtStartPar
\sphinxstylestrong{Arguments}
\begin{quote}

\sphinxAtStartPar
\sphinxcode{\sphinxupquote{var1}}: first variable of new quadratic term.

\sphinxAtStartPar
\sphinxcode{\sphinxupquote{var2}}: second variable of new quadratic term.

\sphinxAtStartPar
\sphinxcode{\sphinxupquote{coeff}}: coefficient of new quadratic term.
\end{quote}
\end{quote}

\subsubsection{MQExpression.Clone()}
\label{\detokenize{csapi/mqexpression:mqexpression-clone}}\begin{quote}

\sphinxAtStartPar
Clone MQExpression object.

\sphinxAtStartPar
\sphinxstylestrong{Synopsis}
\begin{quote}

\sphinxAtStartPar
\sphinxcode{\sphinxupquote{MQExpression Clone()}}
\end{quote}

\sphinxAtStartPar
\sphinxstylestrong{Return}
\begin{quote}

\sphinxAtStartPar
new MQExpression object.
\end{quote}
\end{quote}

\subsubsection{MQExpression.Evaluate()}
\label{\detokenize{csapi/mqexpression:mqexpression-evaluate}}\begin{quote}

\sphinxAtStartPar
evaluate MQExpression object after solving.

\sphinxAtStartPar
\sphinxstylestrong{Synopsis}
\begin{quote}

\sphinxAtStartPar
\sphinxcode{\sphinxupquote{double Evaluate()}}
\end{quote}

\sphinxAtStartPar
\sphinxstylestrong{Return}
\begin{quote}

\sphinxAtStartPar
value of MQExpression object.
\end{quote}
\end{quote}

\subsubsection{MQExpression.GetConstant()}
\label{\detokenize{csapi/mqexpression:mqexpression-getconstant}}\begin{quote}

\sphinxAtStartPar
Get constant in expression.

\sphinxAtStartPar
\sphinxstylestrong{Synopsis}
\begin{quote}

\sphinxAtStartPar
\sphinxcode{\sphinxupquote{double GetConstant()}}
\end{quote}

\sphinxAtStartPar
\sphinxstylestrong{Return}
\begin{quote}

\sphinxAtStartPar
constant in expression.
\end{quote}
\end{quote}

\subsubsection{MQExpression.Represent()}
\label{\detokenize{csapi/mqexpression:mqexpression-represent}}\begin{quote}

\sphinxAtStartPar
String representation of MQExpression object.

\sphinxAtStartPar
\sphinxstylestrong{Synopsis}
\begin{quote}

\sphinxAtStartPar
\sphinxcode{\sphinxupquote{string Represent(uint maxlen)}}
\end{quote}

\sphinxAtStartPar
\sphinxstylestrong{Arguments}
\begin{quote}

\sphinxAtStartPar
\sphinxcode{\sphinxupquote{maxlen}}: maximum buffer length for representations string.
\end{quote}

\sphinxAtStartPar
\sphinxstylestrong{Return}
\begin{quote}

\sphinxAtStartPar
string object.
\end{quote}
\end{quote}

\subsection{MQuadExpr}
\label{\detokenize{csharpapiref:mquadexpr}}\label{\detokenize{csharpapiref:chapcsharpapiref-mquadexpr}}
\sphinxAtStartPar
COPT multi\sphinxhyphen{}dimensional quadratic expression object.It is used to construct
multi\sphinxhyphen{}dimensional quadratic expressions and perform operations with the
multi\sphinxhyphen{}dimensional array built in COPT {\hyperref[\detokenize{csharpapiref:chapcsharpapiref-ndarray}]{\sphinxcrossref{\DUrole{std,std-ref}{NdArray}}}} .
Its elements are {\hyperref[\detokenize{csharpapiref:chapcsharpapiref-mqexpression}]{\sphinxcrossref{\DUrole{std,std-ref}{MQExpression}}}} objects.
It can be created by quadratic combination of {\hyperref[\detokenize{csharpapiref:chapcsharpapiref-mvar}]{\sphinxcrossref{\DUrole{std,std-ref}{MVar}}}} objects.
The following methods are provided:

\sphinxstepscope

\subsubsection{MQuadExpr.AddConstant()}
\label{\detokenize{csapi/mquadexpr:mquadexpr-addconstant}}\label{\detokenize{csapi/mquadexpr::doc}}\begin{quote}

\sphinxAtStartPar
Add constant to each quadratic expression in MQuadExpr object.

\sphinxAtStartPar
\sphinxstylestrong{Synopsis}
\begin{quote}

\sphinxAtStartPar
\sphinxcode{\sphinxupquote{void AddConstant(double constant)}}
\end{quote}

\sphinxAtStartPar
\sphinxstylestrong{Arguments}
\begin{quote}

\sphinxAtStartPar
\sphinxcode{\sphinxupquote{constant}}: the value of the constant.
\end{quote}
\end{quote}

\subsubsection{MQuadExpr.AddConstant()}
\label{\detokenize{csapi/mquadexpr:id1}}\begin{quote}

\sphinxAtStartPar
Add constants to each quadratic expression in MQuadExpr object.

\sphinxAtStartPar
\sphinxstylestrong{Synopsis}
\begin{quote}

\sphinxAtStartPar
\sphinxcode{\sphinxupquote{void AddConstant(NdArray\textless{}double\textgreater{} constants)}}
\end{quote}

\sphinxAtStartPar
\sphinxstylestrong{Arguments}
\begin{quote}

\sphinxAtStartPar
\sphinxcode{\sphinxupquote{constants}}: N\sphinxhyphen{}dimension NdArray object.
\end{quote}
\end{quote}

\subsubsection{MQuadExpr.AddExpr()}
\label{\detokenize{csapi/mquadexpr:mquadexpr-addexpr}}\begin{quote}

\sphinxAtStartPar
Add a linear expression to each quadratic expression in MQuadExpr object.

\sphinxAtStartPar
\sphinxstylestrong{Synopsis}
\begin{quote}

\sphinxAtStartPar
\sphinxcode{\sphinxupquote{void AddExpr(Expr expr, double mult)}}
\end{quote}

\sphinxAtStartPar
\sphinxstylestrong{Arguments}
\begin{quote}

\sphinxAtStartPar
\sphinxcode{\sphinxupquote{expr}}: linear expression object.

\sphinxAtStartPar
\sphinxcode{\sphinxupquote{mult}}: the multiplier of linear expression, default value is 1.0.
\end{quote}
\end{quote}

\subsubsection{MQuadExpr.AddMExpr()}
\label{\detokenize{csapi/mquadexpr:mquadexpr-addmexpr}}\begin{quote}

\sphinxAtStartPar
Add MExpression to each quadratic expression in MQuadExpr object.

\sphinxAtStartPar
\sphinxstylestrong{Synopsis}
\begin{quote}

\sphinxAtStartPar
\sphinxcode{\sphinxupquote{void AddMExpr(MExpression expr, double mult)}}
\end{quote}

\sphinxAtStartPar
\sphinxstylestrong{Arguments}
\begin{quote}

\sphinxAtStartPar
\sphinxcode{\sphinxupquote{expr}}: MExpression object.

\sphinxAtStartPar
\sphinxcode{\sphinxupquote{mult}}: the multiplier of MExpression, default value is 1.0.
\end{quote}
\end{quote}

\subsubsection{MQuadExpr.AddMLinExpr()}
\label{\detokenize{csapi/mquadexpr:mquadexpr-addmlinexpr}}\begin{quote}

\sphinxAtStartPar
Add linear expressions to MQuadExpr object.

\sphinxAtStartPar
\sphinxstylestrong{Synopsis}
\begin{quote}

\sphinxAtStartPar
\sphinxcode{\sphinxupquote{void AddMLinExpr(MLinExpr exprs, double mult)}}
\end{quote}

\sphinxAtStartPar
\sphinxstylestrong{Arguments}
\begin{quote}

\sphinxAtStartPar
\sphinxcode{\sphinxupquote{exprs}}: N\sphinxhyphen{}dimension MLinExpr object.

\sphinxAtStartPar
\sphinxcode{\sphinxupquote{mult}}: the same multiplier for added linear expressions, default value is 1.0.
\end{quote}
\end{quote}

\subsubsection{MQuadExpr.AddMQExpr()}
\label{\detokenize{csapi/mquadexpr:mquadexpr-addmqexpr}}\begin{quote}

\sphinxAtStartPar
Add MQExpression to each quadratic expression in MQuadExpr object.

\sphinxAtStartPar
\sphinxstylestrong{Synopsis}
\begin{quote}

\sphinxAtStartPar
\sphinxcode{\sphinxupquote{void AddMQExpr(MQExpression expr, double mult)}}
\end{quote}

\sphinxAtStartPar
\sphinxstylestrong{Arguments}
\begin{quote}

\sphinxAtStartPar
\sphinxcode{\sphinxupquote{expr}}: MQExpression object.

\sphinxAtStartPar
\sphinxcode{\sphinxupquote{mult}}: the multiplier of MQExpression, default value is 1.0.
\end{quote}
\end{quote}

\subsubsection{MQuadExpr.AddMQuadExpr()}
\label{\detokenize{csapi/mquadexpr:mquadexpr-addmquadexpr}}\begin{quote}

\sphinxAtStartPar
Add quadratic expressions to MQuadExpr object.

\sphinxAtStartPar
\sphinxstylestrong{Synopsis}
\begin{quote}

\sphinxAtStartPar
\sphinxcode{\sphinxupquote{void AddMQuadExpr(MQuadExpr exprs, double mult)}}
\end{quote}

\sphinxAtStartPar
\sphinxstylestrong{Arguments}
\begin{quote}

\sphinxAtStartPar
\sphinxcode{\sphinxupquote{exprs}}: N\sphinxhyphen{}dimension MQuadExpr object.

\sphinxAtStartPar
\sphinxcode{\sphinxupquote{mult}}: the same multiplier for added quadratic expressions, default value is 1.0.
\end{quote}
\end{quote}

\subsubsection{MQuadExpr.AddQuadExpr()}
\label{\detokenize{csapi/mquadexpr:mquadexpr-addquadexpr}}\begin{quote}

\sphinxAtStartPar
Add a quadratic expression to each quadratic expression in MQuadExpr object.

\sphinxAtStartPar
\sphinxstylestrong{Synopsis}
\begin{quote}

\sphinxAtStartPar
\sphinxcode{\sphinxupquote{void AddQuadExpr(QuadExpr expr, double mult)}}
\end{quote}

\sphinxAtStartPar
\sphinxstylestrong{Arguments}
\begin{quote}

\sphinxAtStartPar
\sphinxcode{\sphinxupquote{expr}}: quadratic expression object.

\sphinxAtStartPar
\sphinxcode{\sphinxupquote{mult}}: the multiplier of quadratic expression, default value is 1.0.
\end{quote}
\end{quote}

\subsubsection{MQuadExpr.AddTerm()}
\label{\detokenize{csapi/mquadexpr:mquadexpr-addterm}}\begin{quote}

\sphinxAtStartPar
Add a linear term to MQuadExpr object.

\sphinxAtStartPar
\sphinxstylestrong{Synopsis}
\begin{quote}

\sphinxAtStartPar
\sphinxcode{\sphinxupquote{void AddTerm(Var var, double coeff)}}
\end{quote}

\sphinxAtStartPar
\sphinxstylestrong{Arguments}
\begin{quote}

\sphinxAtStartPar
\sphinxcode{\sphinxupquote{var}}: variable of new term.

\sphinxAtStartPar
\sphinxcode{\sphinxupquote{coeff}}: coefficient of new term.
\end{quote}
\end{quote}

\subsubsection{MQuadExpr.AddTerm()}
\label{\detokenize{csapi/mquadexpr:id2}}\begin{quote}

\sphinxAtStartPar
Add a quadratic term to MQuadExpr object.

\sphinxAtStartPar
\sphinxstylestrong{Synopsis}
\begin{quote}

\sphinxAtStartPar
\sphinxcode{\sphinxupquote{void AddTerm(}}
\begin{quote}

\sphinxAtStartPar
\sphinxcode{\sphinxupquote{Var var1,}}

\sphinxAtStartPar
\sphinxcode{\sphinxupquote{Var var2,}}

\sphinxAtStartPar
\sphinxcode{\sphinxupquote{double coeff)}}
\end{quote}
\end{quote}

\sphinxAtStartPar
\sphinxstylestrong{Arguments}
\begin{quote}

\sphinxAtStartPar
\sphinxcode{\sphinxupquote{var1}}: first variable of new quadratic term.

\sphinxAtStartPar
\sphinxcode{\sphinxupquote{var2}}: second variable of new quadratic term.

\sphinxAtStartPar
\sphinxcode{\sphinxupquote{coeff}}: coefficient of new quadratic term.
\end{quote}
\end{quote}

\subsubsection{MQuadExpr.AddTerms()}
\label{\detokenize{csapi/mquadexpr:mquadexpr-addterms}}\begin{quote}

\sphinxAtStartPar
Add terms to quadratic expressions in MQuadExpr object.

\sphinxAtStartPar
\sphinxstylestrong{Synopsis}
\begin{quote}

\sphinxAtStartPar
\sphinxcode{\sphinxupquote{void AddTerms(MVar vars, double mult)}}
\end{quote}

\sphinxAtStartPar
\sphinxstylestrong{Arguments}
\begin{quote}

\sphinxAtStartPar
\sphinxcode{\sphinxupquote{vars}}: N\sphinxhyphen{}dimension MVar object for added terms.

\sphinxAtStartPar
\sphinxcode{\sphinxupquote{mult}}: the same coefficient for added terms, default value 1.0.
\end{quote}
\end{quote}

\subsubsection{MQuadExpr.AddTerms()}
\label{\detokenize{csapi/mquadexpr:id3}}\begin{quote}

\sphinxAtStartPar
Add terms to quadratic expressions in MQuadExpr object.

\sphinxAtStartPar
\sphinxstylestrong{Synopsis}
\begin{quote}

\sphinxAtStartPar
\sphinxcode{\sphinxupquote{void AddTerms(MVar vars, NdArray\textless{}double\textgreater{} coeffs)}}
\end{quote}

\sphinxAtStartPar
\sphinxstylestrong{Arguments}
\begin{quote}

\sphinxAtStartPar
\sphinxcode{\sphinxupquote{vars}}: N\sphinxhyphen{}dimension MVar object for added terms.

\sphinxAtStartPar
\sphinxcode{\sphinxupquote{coeffs}}: N\sphinxhyphen{}dimension NdArray object of coefficients for added terms.
\end{quote}
\end{quote}

\subsubsection{MQuadExpr.Clear()}
\label{\detokenize{csapi/mquadexpr:mquadexpr-clear}}\begin{quote}

\sphinxAtStartPar
Clear MQuadExpr object.

\sphinxAtStartPar
\sphinxstylestrong{Synopsis}
\begin{quote}

\sphinxAtStartPar
\sphinxcode{\sphinxupquote{void Clear()}}
\end{quote}
\end{quote}

\subsubsection{MQuadExpr.Clone()}
\label{\detokenize{csapi/mquadexpr:mquadexpr-clone}}\begin{quote}

\sphinxAtStartPar
Clone MQuadExpr object.

\sphinxAtStartPar
\sphinxstylestrong{Synopsis}
\begin{quote}

\sphinxAtStartPar
\sphinxcode{\sphinxupquote{MQuadExpr Clone()}}
\end{quote}

\sphinxAtStartPar
\sphinxstylestrong{Return}
\begin{quote}

\sphinxAtStartPar
new MQuadExpr object.
\end{quote}
\end{quote}

\subsubsection{MQuadExpr.Diagonal()}
\label{\detokenize{csapi/mquadexpr:mquadexpr-diagonal}}\begin{quote}

\sphinxAtStartPar
Get diagonals of MQuadExpr object.

\sphinxAtStartPar
\sphinxstylestrong{Synopsis}
\begin{quote}

\sphinxAtStartPar
\sphinxcode{\sphinxupquote{MQuadExpr Diagonal(}}
\begin{quote}

\sphinxAtStartPar
\sphinxcode{\sphinxupquote{int offset,}}

\sphinxAtStartPar
\sphinxcode{\sphinxupquote{int axis1,}}

\sphinxAtStartPar
\sphinxcode{\sphinxupquote{int axis2)}}
\end{quote}
\end{quote}

\sphinxAtStartPar
\sphinxstylestrong{Arguments}
\begin{quote}

\sphinxAtStartPar
\sphinxcode{\sphinxupquote{offset}}: offset of the diagonal from the main diagonal. Can be positive or negative.

\sphinxAtStartPar
\sphinxcode{\sphinxupquote{axis1}}: 1st axis of MQuadExpr.

\sphinxAtStartPar
\sphinxcode{\sphinxupquote{axis2}}: 2nd axis of MQuadExpr.
\end{quote}

\sphinxAtStartPar
\sphinxstylestrong{Return}
\begin{quote}

\sphinxAtStartPar
(N\sphinxhyphen{}1)\sphinxhyphen{}dimensional diagonals.
\end{quote}
\end{quote}

\subsubsection{MQuadExpr.Evaluate()}
\label{\detokenize{csapi/mquadexpr:mquadexpr-evaluate}}\begin{quote}

\sphinxAtStartPar
Evaluate MQuadExpr object after solving.

\sphinxAtStartPar
\sphinxstylestrong{Synopsis}
\begin{quote}

\sphinxAtStartPar
\sphinxcode{\sphinxupquote{double Evaluate()}}
\end{quote}

\sphinxAtStartPar
\sphinxstylestrong{Return}
\begin{quote}

\sphinxAtStartPar
NdArray object storing value of each quadratic expression.
\end{quote}
\end{quote}

\subsubsection{MQuadExpr.Expand()}
\label{\detokenize{csapi/mquadexpr:mquadexpr-expand}}\begin{quote}

\sphinxAtStartPar
Expand shape of MQuadExpr object.

\sphinxAtStartPar
\sphinxstylestrong{Synopsis}
\begin{quote}

\sphinxAtStartPar
\sphinxcode{\sphinxupquote{MQuadExpr Expand(int axis)}}
\end{quote}

\sphinxAtStartPar
\sphinxstylestrong{Arguments}
\begin{quote}

\sphinxAtStartPar
\sphinxcode{\sphinxupquote{axis}}: axis of MQuadExpr.
\end{quote}

\sphinxAtStartPar
\sphinxstylestrong{Return}
\begin{quote}

\sphinxAtStartPar
MQuadExpr object of (N+1)\sphinxhyphen{}dimensional shape.
\end{quote}
\end{quote}

\subsubsection{MQuadExpr.Flatten()}
\label{\detokenize{csapi/mquadexpr:mquadexpr-flatten}}\begin{quote}

\sphinxAtStartPar
Flatten a MQuadExpr object to a 1\sphinxhyphen{}dimensional shape.

\sphinxAtStartPar
\sphinxstylestrong{Synopsis}
\begin{quote}

\sphinxAtStartPar
\sphinxcode{\sphinxupquote{MQuadExpr Flatten()}}
\end{quote}

\sphinxAtStartPar
\sphinxstylestrong{Return}
\begin{quote}

\sphinxAtStartPar
a MQuadExpr object collapsed into one dimension.
\end{quote}
\end{quote}

\subsubsection{MQuadExpr.GetDim()}
\label{\detokenize{csapi/mquadexpr:mquadexpr-getdim}}\begin{quote}

\sphinxAtStartPar
Get i\sphinxhyphen{}th dimension of MQuadExpr object.

\sphinxAtStartPar
\sphinxstylestrong{Synopsis}
\begin{quote}

\sphinxAtStartPar
\sphinxcode{\sphinxupquote{uint GetDim(int i)}}
\end{quote}

\sphinxAtStartPar
\sphinxstylestrong{Arguments}
\begin{quote}

\sphinxAtStartPar
\sphinxcode{\sphinxupquote{i}}: index of dimension
\end{quote}

\sphinxAtStartPar
\sphinxstylestrong{Return}
\begin{quote}

\sphinxAtStartPar
i\sphinxhyphen{}th dimension.
\end{quote}
\end{quote}

\subsubsection{MQuadExpr.GetItem()}
\label{\detokenize{csapi/mquadexpr:mquadexpr-getitem}}\begin{quote}

\sphinxAtStartPar
Get quadratic expression of given index from MQuadExpr object.

\sphinxAtStartPar
\sphinxstylestrong{Synopsis}
\begin{quote}

\sphinxAtStartPar
\sphinxcode{\sphinxupquote{MQExpression GetItem(long idx)}}
\end{quote}

\sphinxAtStartPar
\sphinxstylestrong{Arguments}
\begin{quote}

\sphinxAtStartPar
\sphinxcode{\sphinxupquote{idx}}: index of quadratic expression.
\end{quote}

\sphinxAtStartPar
\sphinxstylestrong{Return}
\begin{quote}

\sphinxAtStartPar
MQExpression object.
\end{quote}
\end{quote}

\subsubsection{MQuadExpr.GetItem()}
\label{\detokenize{csapi/mquadexpr:id4}}\begin{quote}

\sphinxAtStartPar
Get sub\sphinxhyphen{}arrays of MQuadExpr object, given view object.

\sphinxAtStartPar
\sphinxstylestrong{Synopsis}
\begin{quote}

\sphinxAtStartPar
\sphinxcode{\sphinxupquote{MQuadExpr GetItem(View view)}}
\end{quote}

\sphinxAtStartPar
\sphinxstylestrong{Arguments}
\begin{quote}

\sphinxAtStartPar
\sphinxcode{\sphinxupquote{view}}: view of multi\sphinxhyphen{}dimensional array.
\end{quote}

\sphinxAtStartPar
\sphinxstylestrong{Return}
\begin{quote}

\sphinxAtStartPar
sub\sphinxhyphen{}arrays of MQuadExpr object.
\end{quote}
\end{quote}

\subsubsection{MQuadExpr.GetND()}
\label{\detokenize{csapi/mquadexpr:mquadexpr-getnd}}\begin{quote}

\sphinxAtStartPar
Get number of dimensions of MQuadExpr object.

\sphinxAtStartPar
\sphinxstylestrong{Synopsis}
\begin{quote}

\sphinxAtStartPar
\sphinxcode{\sphinxupquote{int GetND()}}
\end{quote}

\sphinxAtStartPar
\sphinxstylestrong{Return}
\begin{quote}

\sphinxAtStartPar
number of dimensions.
\end{quote}
\end{quote}

\subsubsection{MQuadExpr.GetShape()}
\label{\detokenize{csapi/mquadexpr:mquadexpr-getshape}}\begin{quote}

\sphinxAtStartPar
Get shape of MQuadExpr object.

\sphinxAtStartPar
\sphinxstylestrong{Synopsis}
\begin{quote}

\sphinxAtStartPar
\sphinxcode{\sphinxupquote{Shape GetShape()}}
\end{quote}

\sphinxAtStartPar
\sphinxstylestrong{Return}
\begin{quote}

\sphinxAtStartPar
shape object.
\end{quote}
\end{quote}

\subsubsection{MQuadExpr.GetSize()}
\label{\detokenize{csapi/mquadexpr:mquadexpr-getsize}}\begin{quote}

\sphinxAtStartPar
Get size of MQuadExpr object.

\sphinxAtStartPar
\sphinxstylestrong{Synopsis}
\begin{quote}

\sphinxAtStartPar
\sphinxcode{\sphinxupquote{uint GetSize()}}
\end{quote}

\sphinxAtStartPar
\sphinxstylestrong{Return}
\begin{quote}

\sphinxAtStartPar
number of linear expressions.
\end{quote}
\end{quote}

\subsubsection{MQuadExpr.HStack\textless{}T\textgreater{}()}
\label{\detokenize{csapi/mquadexpr:mquadexpr-hstack-t}}\begin{quote}

\sphinxAtStartPar
Stack with other NdArray object along horizontal axis.

\sphinxAtStartPar
\sphinxstylestrong{Synopsis}
\begin{quote}

\sphinxAtStartPar
\sphinxcode{\sphinxupquote{MQuadExpr HStack\textless{}T\textgreater{}(NdArray\textless{}T\textgreater{} other)}}
\end{quote}

\sphinxAtStartPar
\sphinxstylestrong{Arguments}
\begin{quote}

\sphinxAtStartPar
\sphinxcode{\sphinxupquote{other}}: a NdArray object.
\end{quote}

\sphinxAtStartPar
\sphinxstylestrong{Return}
\begin{quote}

\sphinxAtStartPar
the result MQuadExpr object.
\end{quote}
\end{quote}

\subsubsection{MQuadExpr.HStack()}
\label{\detokenize{csapi/mquadexpr:mquadexpr-hstack}}\begin{quote}

\sphinxAtStartPar
Stack with other MQuadExpr object along horizontal axis.

\sphinxAtStartPar
\sphinxstylestrong{Synopsis}
\begin{quote}

\sphinxAtStartPar
\sphinxcode{\sphinxupquote{MQuadExpr HStack(MQuadExpr other)}}
\end{quote}

\sphinxAtStartPar
\sphinxstylestrong{Arguments}
\begin{quote}

\sphinxAtStartPar
\sphinxcode{\sphinxupquote{other}}: a MQuadExpr object.
\end{quote}

\sphinxAtStartPar
\sphinxstylestrong{Return}
\begin{quote}

\sphinxAtStartPar
the result MQuadExpr object.
\end{quote}
\end{quote}

\subsubsection{MQuadExpr.HStack()}
\label{\detokenize{csapi/mquadexpr:id5}}\begin{quote}

\sphinxAtStartPar
Stack with other MLinExpr object along horizontal axis.

\sphinxAtStartPar
\sphinxstylestrong{Synopsis}
\begin{quote}

\sphinxAtStartPar
\sphinxcode{\sphinxupquote{MQuadExpr HStack(MLinExpr other)}}
\end{quote}

\sphinxAtStartPar
\sphinxstylestrong{Arguments}
\begin{quote}

\sphinxAtStartPar
\sphinxcode{\sphinxupquote{other}}: a MLinExpr object.
\end{quote}

\sphinxAtStartPar
\sphinxstylestrong{Return}
\begin{quote}

\sphinxAtStartPar
the result MQuadExpr object.
\end{quote}
\end{quote}

\subsubsection{MQuadExpr.HStack()}
\label{\detokenize{csapi/mquadexpr:id6}}\begin{quote}

\sphinxAtStartPar
Stack with other MVar object along horizontal axis.

\sphinxAtStartPar
\sphinxstylestrong{Synopsis}
\begin{quote}

\sphinxAtStartPar
\sphinxcode{\sphinxupquote{MQuadExpr HStack(MVar other)}}
\end{quote}

\sphinxAtStartPar
\sphinxstylestrong{Arguments}
\begin{quote}

\sphinxAtStartPar
\sphinxcode{\sphinxupquote{other}}: a MVar object.
\end{quote}

\sphinxAtStartPar
\sphinxstylestrong{Return}
\begin{quote}

\sphinxAtStartPar
the result MQuadExpr object.
\end{quote}
\end{quote}

\subsubsection{MQuadExpr.Pick()}
\label{\detokenize{csapi/mquadexpr:mquadexpr-pick}}\begin{quote}

\sphinxAtStartPar
Given a list of indexes, get quadratic expressions from MQuadExpr object.

\sphinxAtStartPar
\sphinxstylestrong{Synopsis}
\begin{quote}

\sphinxAtStartPar
\sphinxcode{\sphinxupquote{MQuadExpr Pick(NdArray\textless{}int\textgreater{} indexes)}}
\end{quote}

\sphinxAtStartPar
\sphinxstylestrong{Arguments}
\begin{quote}

\sphinxAtStartPar
\sphinxcode{\sphinxupquote{indexes}}: one or two dimensional indexes of elements. If two dimensional, each row is position of an element.
\end{quote}

\sphinxAtStartPar
\sphinxstylestrong{Return}
\begin{quote}

\sphinxAtStartPar
one\sphinxhyphen{}dimensional array of desired quadratic expressions.
\end{quote}
\end{quote}

\subsubsection{MQuadExpr.Repeat()}
\label{\detokenize{csapi/mquadexpr:mquadexpr-repeat}}\begin{quote}

\sphinxAtStartPar
Repeat each element of MQuadExpr along given axis.

\sphinxAtStartPar
\sphinxstylestrong{Synopsis}
\begin{quote}

\sphinxAtStartPar
\sphinxcode{\sphinxupquote{MQuadExpr Repeat(long repeats, int axis)}}
\end{quote}

\sphinxAtStartPar
\sphinxstylestrong{Arguments}
\begin{quote}

\sphinxAtStartPar
\sphinxcode{\sphinxupquote{repeats}}: number of repetitions for each element.

\sphinxAtStartPar
\sphinxcode{\sphinxupquote{axis}}: axis of MQuadExpr.
\end{quote}

\sphinxAtStartPar
\sphinxstylestrong{Return}
\begin{quote}

\sphinxAtStartPar
new MQuadExpr object.
\end{quote}
\end{quote}

\subsubsection{MQuadExpr.RepeatBlock()}
\label{\detokenize{csapi/mquadexpr:mquadexpr-repeatblock}}\begin{quote}

\sphinxAtStartPar
Repeat an MQuadExpr a number of times along given axis.

\sphinxAtStartPar
\sphinxstylestrong{Synopsis}
\begin{quote}

\sphinxAtStartPar
\sphinxcode{\sphinxupquote{MQuadExpr RepeatBlock(long repeats, int axis)}}
\end{quote}

\sphinxAtStartPar
\sphinxstylestrong{Arguments}
\begin{quote}

\sphinxAtStartPar
\sphinxcode{\sphinxupquote{repeats}}: number of repetitions.

\sphinxAtStartPar
\sphinxcode{\sphinxupquote{axis}}: axis of MQuadExpr.
\end{quote}

\sphinxAtStartPar
\sphinxstylestrong{Return}
\begin{quote}

\sphinxAtStartPar
new MQuadExpr object.
\end{quote}
\end{quote}

\subsubsection{MQuadExpr.Represent()}
\label{\detokenize{csapi/mquadexpr:mquadexpr-represent}}\begin{quote}

\sphinxAtStartPar
String representation of MQuadExpr object.

\sphinxAtStartPar
\sphinxstylestrong{Synopsis}
\begin{quote}

\sphinxAtStartPar
\sphinxcode{\sphinxupquote{string Represent(int maxlen)}}
\end{quote}

\sphinxAtStartPar
\sphinxstylestrong{Arguments}
\begin{quote}

\sphinxAtStartPar
\sphinxcode{\sphinxupquote{maxlen}}: maximum buffer length for representations string.
\end{quote}

\sphinxAtStartPar
\sphinxstylestrong{Return}
\begin{quote}

\sphinxAtStartPar
string object.
\end{quote}
\end{quote}

\subsubsection{MQuadExpr.Reshape()}
\label{\detokenize{csapi/mquadexpr:mquadexpr-reshape}}\begin{quote}

\sphinxAtStartPar
Reshape MQuadExpr object to new shape.

\sphinxAtStartPar
\sphinxstylestrong{Synopsis}
\begin{quote}

\sphinxAtStartPar
\sphinxcode{\sphinxupquote{MQuadExpr Reshape(Shape shp)}}
\end{quote}

\sphinxAtStartPar
\sphinxstylestrong{Arguments}
\begin{quote}

\sphinxAtStartPar
\sphinxcode{\sphinxupquote{shp}}: new shape of M\sphinxhyphen{}dimensions.
\end{quote}

\sphinxAtStartPar
\sphinxstylestrong{Return}
\begin{quote}

\sphinxAtStartPar
M\sphinxhyphen{}dimensional MQuadExpr object.
\end{quote}
\end{quote}

\subsubsection{MQuadExpr.SetItem()}
\label{\detokenize{csapi/mquadexpr:mquadexpr-setitem}}\begin{quote}

\sphinxAtStartPar
Set quadratic expression of given index to MQuadExpr object.

\sphinxAtStartPar
\sphinxstylestrong{Synopsis}
\begin{quote}

\sphinxAtStartPar
\sphinxcode{\sphinxupquote{void SetItem(long idx, MQExpression expr)}}
\end{quote}

\sphinxAtStartPar
\sphinxstylestrong{Arguments}
\begin{quote}

\sphinxAtStartPar
\sphinxcode{\sphinxupquote{idx}}: index of element.

\sphinxAtStartPar
\sphinxcode{\sphinxupquote{expr}}: MQExpression object.
\end{quote}
\end{quote}

\subsubsection{MQuadExpr.Squeeze()}
\label{\detokenize{csapi/mquadexpr:mquadexpr-squeeze}}\begin{quote}

\sphinxAtStartPar
Remove axis of length 1 from shape of MQuadExpr object.

\sphinxAtStartPar
\sphinxstylestrong{Synopsis}
\begin{quote}

\sphinxAtStartPar
\sphinxcode{\sphinxupquote{MQuadExpr Squeeze(int axis)}}
\end{quote}

\sphinxAtStartPar
\sphinxstylestrong{Arguments}
\begin{quote}

\sphinxAtStartPar
\sphinxcode{\sphinxupquote{axis}}: axis of MQuadExpr, where the length is 1.
\end{quote}

\sphinxAtStartPar
\sphinxstylestrong{Return}
\begin{quote}

\sphinxAtStartPar
MQuadExpr object of (N\sphinxhyphen{}1)\sphinxhyphen{}dimensional shape.
\end{quote}
\end{quote}

\subsubsection{MQuadExpr.Stack\textless{}T\textgreater{}()}
\label{\detokenize{csapi/mquadexpr:mquadexpr-stack-t}}\begin{quote}

\sphinxAtStartPar
Stack with other NdArray object along given axis.

\sphinxAtStartPar
\sphinxstylestrong{Synopsis}
\begin{quote}

\sphinxAtStartPar
\sphinxcode{\sphinxupquote{MQuadExpr Stack\textless{}T\textgreater{}(NdArray\textless{}T\textgreater{} other, int axis)}}
\end{quote}

\sphinxAtStartPar
\sphinxstylestrong{Arguments}
\begin{quote}

\sphinxAtStartPar
\sphinxcode{\sphinxupquote{other}}: a NdArray object.

\sphinxAtStartPar
\sphinxcode{\sphinxupquote{axis}}: an axis of MQuadExpr.
\end{quote}

\sphinxAtStartPar
\sphinxstylestrong{Return}
\begin{quote}

\sphinxAtStartPar
the result MQuadExpr object.
\end{quote}
\end{quote}

\subsubsection{MQuadExpr.Stack()}
\label{\detokenize{csapi/mquadexpr:mquadexpr-stack}}\begin{quote}

\sphinxAtStartPar
Stack with other MQuadExpr object along given axis.

\sphinxAtStartPar
\sphinxstylestrong{Synopsis}
\begin{quote}

\sphinxAtStartPar
\sphinxcode{\sphinxupquote{MQuadExpr Stack(MQuadExpr other, int axis)}}
\end{quote}

\sphinxAtStartPar
\sphinxstylestrong{Arguments}
\begin{quote}

\sphinxAtStartPar
\sphinxcode{\sphinxupquote{other}}: a MQuadExpr object.

\sphinxAtStartPar
\sphinxcode{\sphinxupquote{axis}}: an axis of MQuadExpr.
\end{quote}

\sphinxAtStartPar
\sphinxstylestrong{Return}
\begin{quote}

\sphinxAtStartPar
the result MQuadExpr object.
\end{quote}
\end{quote}

\subsubsection{MQuadExpr.Stack()}
\label{\detokenize{csapi/mquadexpr:id7}}\begin{quote}

\sphinxAtStartPar
Stack with other MLinExpr object along given axis.

\sphinxAtStartPar
\sphinxstylestrong{Synopsis}
\begin{quote}

\sphinxAtStartPar
\sphinxcode{\sphinxupquote{MQuadExpr Stack(MLinExpr other, int axis)}}
\end{quote}

\sphinxAtStartPar
\sphinxstylestrong{Arguments}
\begin{quote}

\sphinxAtStartPar
\sphinxcode{\sphinxupquote{other}}: a MLinExpr object.

\sphinxAtStartPar
\sphinxcode{\sphinxupquote{axis}}: an axis of MQuadExpr.
\end{quote}

\sphinxAtStartPar
\sphinxstylestrong{Return}
\begin{quote}

\sphinxAtStartPar
the result MQuadExpr object.
\end{quote}
\end{quote}

\subsubsection{MQuadExpr.Stack()}
\label{\detokenize{csapi/mquadexpr:id8}}\begin{quote}

\sphinxAtStartPar
Stack with other MQuadExpr object along given axis.

\sphinxAtStartPar
\sphinxstylestrong{Synopsis}
\begin{quote}

\sphinxAtStartPar
\sphinxcode{\sphinxupquote{MQuadExpr Stack(MVar other, int axis)}}
\end{quote}

\sphinxAtStartPar
\sphinxstylestrong{Arguments}
\begin{quote}

\sphinxAtStartPar
\sphinxcode{\sphinxupquote{other}}: a MVar object.

\sphinxAtStartPar
\sphinxcode{\sphinxupquote{axis}}: an axis of MQuadExpr.
\end{quote}

\sphinxAtStartPar
\sphinxstylestrong{Return}
\begin{quote}

\sphinxAtStartPar
the result MQuadExpr object.
\end{quote}
\end{quote}

\subsubsection{MQuadExpr.SubConstant()}
\label{\detokenize{csapi/mquadexpr:mquadexpr-subconstant}}\begin{quote}

\sphinxAtStartPar
Substract constants from each quadratic expression in MQuadExpr object.

\sphinxAtStartPar
\sphinxstylestrong{Synopsis}
\begin{quote}

\sphinxAtStartPar
\sphinxcode{\sphinxupquote{void SubConstant(NdArray\textless{}double\textgreater{} constants)}}
\end{quote}

\sphinxAtStartPar
\sphinxstylestrong{Arguments}
\begin{quote}

\sphinxAtStartPar
\sphinxcode{\sphinxupquote{constants}}: N\sphinxhyphen{}dimension NdArray object.
\end{quote}
\end{quote}

\subsubsection{MQuadExpr.Sum()}
\label{\detokenize{csapi/mquadexpr:mquadexpr-sum}}\begin{quote}

\sphinxAtStartPar
Sum of all expressions in MQuadExpr object.

\sphinxAtStartPar
\sphinxstylestrong{Synopsis}
\begin{quote}

\sphinxAtStartPar
\sphinxcode{\sphinxupquote{MQuadExpr Sum()}}
\end{quote}

\sphinxAtStartPar
\sphinxstylestrong{Return}
\begin{quote}

\sphinxAtStartPar
sum in zero dimension.
\end{quote}
\end{quote}

\subsubsection{MQuadExpr.Sum()}
\label{\detokenize{csapi/mquadexpr:id9}}\begin{quote}

\sphinxAtStartPar
Sum of variables at given axis of MQuadExpr object.

\sphinxAtStartPar
\sphinxstylestrong{Synopsis}
\begin{quote}

\sphinxAtStartPar
\sphinxcode{\sphinxupquote{MQuadExpr Sum(int axis)}}
\end{quote}

\sphinxAtStartPar
\sphinxstylestrong{Arguments}
\begin{quote}

\sphinxAtStartPar
\sphinxcode{\sphinxupquote{axis}}: axis of MQuadExpr.
\end{quote}

\sphinxAtStartPar
\sphinxstylestrong{Return}
\begin{quote}

\sphinxAtStartPar
MQuadExpr object in (N\sphinxhyphen{}1)\sphinxhyphen{}dimension.
\end{quote}
\end{quote}

\subsubsection{MQuadExpr.Transpose()}
\label{\detokenize{csapi/mquadexpr:mquadexpr-transpose}}\begin{quote}

\sphinxAtStartPar
Perform matrix transpose of MQuadExpr object.

\sphinxAtStartPar
\sphinxstylestrong{Synopsis}
\begin{quote}

\sphinxAtStartPar
\sphinxcode{\sphinxupquote{MQuadExpr Transpose()}}
\end{quote}

\sphinxAtStartPar
\sphinxstylestrong{Return}
\begin{quote}

\sphinxAtStartPar
transposed MQuadExpr object.
\end{quote}
\end{quote}

\subsubsection{MQuadExpr.VStack\textless{}T\textgreater{}()}
\label{\detokenize{csapi/mquadexpr:mquadexpr-vstack-t}}\begin{quote}

\sphinxAtStartPar
Stack with other NdArray object along vertical axis.

\sphinxAtStartPar
\sphinxstylestrong{Synopsis}
\begin{quote}

\sphinxAtStartPar
\sphinxcode{\sphinxupquote{MQuadExpr VStack\textless{}T\textgreater{}(NdArray\textless{}T\textgreater{} other)}}
\end{quote}

\sphinxAtStartPar
\sphinxstylestrong{Arguments}
\begin{quote}

\sphinxAtStartPar
\sphinxcode{\sphinxupquote{other}}: a NdArray object.
\end{quote}

\sphinxAtStartPar
\sphinxstylestrong{Return}
\begin{quote}

\sphinxAtStartPar
the result MQuadExpr object.
\end{quote}
\end{quote}

\subsubsection{MQuadExpr.VStack()}
\label{\detokenize{csapi/mquadexpr:mquadexpr-vstack}}\begin{quote}

\sphinxAtStartPar
Stack with other MQuadExpr object along vertical axis.

\sphinxAtStartPar
\sphinxstylestrong{Synopsis}
\begin{quote}

\sphinxAtStartPar
\sphinxcode{\sphinxupquote{MQuadExpr VStack(MQuadExpr other)}}
\end{quote}

\sphinxAtStartPar
\sphinxstylestrong{Arguments}
\begin{quote}

\sphinxAtStartPar
\sphinxcode{\sphinxupquote{other}}: a MQuadExpr object.
\end{quote}

\sphinxAtStartPar
\sphinxstylestrong{Return}
\begin{quote}

\sphinxAtStartPar
the result MQuadExpr object.
\end{quote}
\end{quote}

\subsubsection{MQuadExpr.VStack()}
\label{\detokenize{csapi/mquadexpr:id10}}\begin{quote}

\sphinxAtStartPar
Stack with other MLinExpr object along vertical axis.

\sphinxAtStartPar
\sphinxstylestrong{Synopsis}
\begin{quote}

\sphinxAtStartPar
\sphinxcode{\sphinxupquote{MQuadExpr VStack(MLinExpr other)}}
\end{quote}

\sphinxAtStartPar
\sphinxstylestrong{Arguments}
\begin{quote}

\sphinxAtStartPar
\sphinxcode{\sphinxupquote{other}}: a MLinExpr object.
\end{quote}

\sphinxAtStartPar
\sphinxstylestrong{Return}
\begin{quote}

\sphinxAtStartPar
the result MQuadExpr object.
\end{quote}
\end{quote}

\subsubsection{MQuadExpr.VStack()}
\label{\detokenize{csapi/mquadexpr:id11}}\begin{quote}

\sphinxAtStartPar
Stack with other MVar object along vertical axis.

\sphinxAtStartPar
\sphinxstylestrong{Synopsis}
\begin{quote}

\sphinxAtStartPar
\sphinxcode{\sphinxupquote{MQuadExpr VStack(MVar other)}}
\end{quote}

\sphinxAtStartPar
\sphinxstylestrong{Arguments}
\begin{quote}

\sphinxAtStartPar
\sphinxcode{\sphinxupquote{other}}: a MVar object.
\end{quote}

\sphinxAtStartPar
\sphinxstylestrong{Return}
\begin{quote}

\sphinxAtStartPar
the result MQuadExpr object.
\end{quote}
\end{quote}

\subsection{NdArray}
\label{\detokenize{csharpapiref:ndarray}}\label{\detokenize{csharpapiref:chapcsharpapiref-ndarray}}
\sphinxAtStartPar
The NdArray class is a built\sphinxhyphen{}in multi\sphinxhyphen{}dimensional array in COPT. It represents a
table of elements of the same type, indexed by a tuple of integers. The following
methods are provided:

\sphinxstepscope

\subsubsection{NdArray\textless{}T\textgreater{}.Diagonal()}
\label{\detokenize{csapi/ndarray:ndarray-t-diagonal}}\label{\detokenize{csapi/ndarray::doc}}\begin{quote}

\sphinxAtStartPar
Get data type of NdArray object. Get diagonals of NdArray object.

\sphinxAtStartPar
\sphinxstylestrong{Synopsis}
\begin{quote}

\sphinxAtStartPar
\sphinxcode{\sphinxupquote{NdArray\textless{}T\textgreater{} Diagonal(}}
\begin{quote}

\sphinxAtStartPar
\sphinxcode{\sphinxupquote{int offset,}}

\sphinxAtStartPar
\sphinxcode{\sphinxupquote{int axis1,}}

\sphinxAtStartPar
\sphinxcode{\sphinxupquote{int axis2)}}
\end{quote}
\end{quote}

\sphinxAtStartPar
\sphinxstylestrong{Arguments}
\begin{quote}

\sphinxAtStartPar
\sphinxcode{\sphinxupquote{offset}}: offset of the diagonal from the main diagonal. Can be positive or negative.

\sphinxAtStartPar
\sphinxcode{\sphinxupquote{axis1}}: 1st axis of NdArray.

\sphinxAtStartPar
\sphinxcode{\sphinxupquote{axis2}}: 2nd axis of NdArray.
\end{quote}

\sphinxAtStartPar
\sphinxstylestrong{Return}
\begin{quote}

\sphinxAtStartPar
data type of elements. (N\sphinxhyphen{}1)\sphinxhyphen{}dimensional diagonals.
\end{quote}
\end{quote}

\subsubsection{NdArray\textless{}T\textgreater{}.Dot()}
\label{\detokenize{csapi/ndarray:ndarray-t-dot}}\begin{quote}

\sphinxAtStartPar
Dot product with another NdArray object of type double.

\sphinxAtStartPar
\sphinxstylestrong{Synopsis}
\begin{quote}

\sphinxAtStartPar
\sphinxcode{\sphinxupquote{double Dot(NdArray\textless{}double\textgreater{} other)}}
\end{quote}

\sphinxAtStartPar
\sphinxstylestrong{Arguments}
\begin{quote}

\sphinxAtStartPar
\sphinxcode{\sphinxupquote{other}}: another NdArray object of type double.
\end{quote}

\sphinxAtStartPar
\sphinxstylestrong{Return}
\begin{quote}

\sphinxAtStartPar
dot product value.
\end{quote}
\end{quote}

\subsubsection{NdArray\textless{}T\textgreater{}.Dot()}
\label{\detokenize{csapi/ndarray:id1}}\begin{quote}

\sphinxAtStartPar
Dot product with another NdArray object of type int.

\sphinxAtStartPar
\sphinxstylestrong{Synopsis}
\begin{quote}

\sphinxAtStartPar
\sphinxcode{\sphinxupquote{T Dot(NdArray\textless{}int\textgreater{} other)}}
\end{quote}

\sphinxAtStartPar
\sphinxstylestrong{Arguments}
\begin{quote}

\sphinxAtStartPar
\sphinxcode{\sphinxupquote{other}}: another NdArray object of type int.
\end{quote}

\sphinxAtStartPar
\sphinxstylestrong{Return}
\begin{quote}

\sphinxAtStartPar
dot product value.
\end{quote}
\end{quote}

\subsubsection{NdArray\textless{}T\textgreater{}.Expand()}
\label{\detokenize{csapi/ndarray:ndarray-t-expand}}\begin{quote}

\sphinxAtStartPar
Expand shape of NdArray object.

\sphinxAtStartPar
\sphinxstylestrong{Synopsis}
\begin{quote}

\sphinxAtStartPar
\sphinxcode{\sphinxupquote{NdArray\textless{}T\textgreater{} Expand(int axis)}}
\end{quote}

\sphinxAtStartPar
\sphinxstylestrong{Arguments}
\begin{quote}

\sphinxAtStartPar
\sphinxcode{\sphinxupquote{axis}}: axis of NdArray.
\end{quote}

\sphinxAtStartPar
\sphinxstylestrong{Return}
\begin{quote}

\sphinxAtStartPar
NdArray object in (N+1)\sphinxhyphen{}dimensions.
\end{quote}
\end{quote}

\subsubsection{NdArray\textless{}T\textgreater{}.Fill()}
\label{\detokenize{csapi/ndarray:ndarray-t-fill}}\begin{quote}

\sphinxAtStartPar
Fill NdArray object with given value.

\sphinxAtStartPar
\sphinxstylestrong{Synopsis}
\begin{quote}

\sphinxAtStartPar
\sphinxcode{\sphinxupquote{void Fill(T val)}}
\end{quote}

\sphinxAtStartPar
\sphinxstylestrong{Arguments}
\begin{quote}

\sphinxAtStartPar
\sphinxcode{\sphinxupquote{val}}: new value.
\end{quote}
\end{quote}

\subsubsection{NdArray\textless{}T\textgreater{}.Flatten()}
\label{\detokenize{csapi/ndarray:ndarray-t-flatten}}\begin{quote}

\sphinxAtStartPar
Flatten an NdArray object to a 1\sphinxhyphen{}dimensional shape.

\sphinxAtStartPar
\sphinxstylestrong{Synopsis}
\begin{quote}

\sphinxAtStartPar
\sphinxcode{\sphinxupquote{NdArray\textless{}T\textgreater{} Flatten()}}
\end{quote}

\sphinxAtStartPar
\sphinxstylestrong{Return}
\begin{quote}

\sphinxAtStartPar
An NdArray object collapsed into one dimension.
\end{quote}
\end{quote}

\subsubsection{NdArray\textless{}T\textgreater{}.GetDim()}
\label{\detokenize{csapi/ndarray:ndarray-t-getdim}}\begin{quote}

\sphinxAtStartPar
Get i\sphinxhyphen{}th dimension in NdArray object.

\sphinxAtStartPar
\sphinxstylestrong{Synopsis}
\begin{quote}

\sphinxAtStartPar
\sphinxcode{\sphinxupquote{long GetDim(int i)}}
\end{quote}

\sphinxAtStartPar
\sphinxstylestrong{Arguments}
\begin{quote}

\sphinxAtStartPar
\sphinxcode{\sphinxupquote{i}}: index of dimensions.
\end{quote}

\sphinxAtStartPar
\sphinxstylestrong{Return}
\begin{quote}

\sphinxAtStartPar
the i\sphinxhyphen{}th dimension.
\end{quote}
\end{quote}

\subsubsection{NdArray\textless{}T\textgreater{}.GetItem()}
\label{\detokenize{csapi/ndarray:ndarray-t-getitem}}\begin{quote}

\sphinxAtStartPar
Get element of given index from NdArray object.

\sphinxAtStartPar
\sphinxstylestrong{Synopsis}
\begin{quote}

\sphinxAtStartPar
\sphinxcode{\sphinxupquote{T GetItem(long idx)}}
\end{quote}

\sphinxAtStartPar
\sphinxstylestrong{Arguments}
\begin{quote}

\sphinxAtStartPar
\sphinxcode{\sphinxupquote{idx}}: index of element.
\end{quote}

\sphinxAtStartPar
\sphinxstylestrong{Return}
\begin{quote}

\sphinxAtStartPar
value of element.
\end{quote}
\end{quote}

\subsubsection{NdArray\textless{}T\textgreater{}.GetItem()}
\label{\detokenize{csapi/ndarray:id2}}\begin{quote}

\sphinxAtStartPar
Get sub\sphinxhyphen{}array of NdArray object, given View object.

\sphinxAtStartPar
\sphinxstylestrong{Synopsis}
\begin{quote}

\sphinxAtStartPar
\sphinxcode{\sphinxupquote{NdArray\textless{}T\textgreater{} GetItem(View view)}}
\end{quote}

\sphinxAtStartPar
\sphinxstylestrong{Arguments}
\begin{quote}

\sphinxAtStartPar
\sphinxcode{\sphinxupquote{view}}: View object.
\end{quote}

\sphinxAtStartPar
\sphinxstylestrong{Return}
\begin{quote}

\sphinxAtStartPar
sub NdArray without copying underlying data.
\end{quote}
\end{quote}

\subsubsection{NdArray\textless{}T\textgreater{}.GetND()}
\label{\detokenize{csapi/ndarray:ndarray-t-getnd}}\begin{quote}

\sphinxAtStartPar
Get number of dimensions in NdArray object.

\sphinxAtStartPar
\sphinxstylestrong{Synopsis}
\begin{quote}

\sphinxAtStartPar
\sphinxcode{\sphinxupquote{int GetND()}}
\end{quote}

\sphinxAtStartPar
\sphinxstylestrong{Return}
\begin{quote}

\sphinxAtStartPar
number of dimensions.
\end{quote}
\end{quote}

\subsubsection{NdArray\textless{}T\textgreater{}.GetShape()}
\label{\detokenize{csapi/ndarray:ndarray-t-getshape}}\begin{quote}

\sphinxAtStartPar
Get shape of NdArray object.

\sphinxAtStartPar
\sphinxstylestrong{Synopsis}
\begin{quote}

\sphinxAtStartPar
\sphinxcode{\sphinxupquote{Shape GetShape()}}
\end{quote}

\sphinxAtStartPar
\sphinxstylestrong{Return}
\begin{quote}

\sphinxAtStartPar
shape object.
\end{quote}
\end{quote}

\subsubsection{NdArray\textless{}T\textgreater{}.GetSize()}
\label{\detokenize{csapi/ndarray:ndarray-t-getsize}}\begin{quote}

\sphinxAtStartPar
Get size of NdArray object.

\sphinxAtStartPar
\sphinxstylestrong{Synopsis}
\begin{quote}

\sphinxAtStartPar
\sphinxcode{\sphinxupquote{long GetSize()}}
\end{quote}

\sphinxAtStartPar
\sphinxstylestrong{Return}
\begin{quote}

\sphinxAtStartPar
size of NdArray.
\end{quote}
\end{quote}

\subsubsection{NdArray\textless{}T\textgreater{}.NdArray()}
\label{\detokenize{csapi/ndarray:ndarray-t-ndarray}}\begin{quote}

\sphinxAtStartPar
Constructor of NdArray object.

\sphinxAtStartPar
\sphinxstylestrong{Synopsis}
\begin{quote}

\sphinxAtStartPar
\sphinxcode{\sphinxupquote{NdArray(Shape shp, T val)}}
\end{quote}

\sphinxAtStartPar
\sphinxstylestrong{Arguments}
\begin{quote}

\sphinxAtStartPar
\sphinxcode{\sphinxupquote{shp}}: shape of NdArray object.

\sphinxAtStartPar
\sphinxcode{\sphinxupquote{val}}: initial value of elements, including int, long, float and double.
\end{quote}
\end{quote}

\subsubsection{NdArray\textless{}T\textgreater{}.Pick()}
\label{\detokenize{csapi/ndarray:ndarray-t-pick}}\begin{quote}

\sphinxAtStartPar
Given a list of indexes, get elements from NdArray object.

\sphinxAtStartPar
\sphinxstylestrong{Synopsis}
\begin{quote}

\sphinxAtStartPar
\sphinxcode{\sphinxupquote{NdArray\textless{}T\textgreater{} Pick(NdArray\textless{}int\textgreater{} indexes)}}
\end{quote}

\sphinxAtStartPar
\sphinxstylestrong{Arguments}
\begin{quote}

\sphinxAtStartPar
\sphinxcode{\sphinxupquote{indexes}}: one or two dimensional indexes of elements. if two dimensional, each row is position of an element.
\end{quote}

\sphinxAtStartPar
\sphinxstylestrong{Return}
\begin{quote}

\sphinxAtStartPar
one\sphinxhyphen{}dimensional array of desired elements.
\end{quote}
\end{quote}

\subsubsection{NdArray\textless{}T\textgreater{}.Prod()}
\label{\detokenize{csapi/ndarray:ndarray-t-prod}}\begin{quote}

\sphinxAtStartPar
Product of all elements in NdArray object.

\sphinxAtStartPar
\sphinxstylestrong{Synopsis}
\begin{quote}

\sphinxAtStartPar
\sphinxcode{\sphinxupquote{T Prod()}}
\end{quote}

\sphinxAtStartPar
\sphinxstylestrong{Return}
\begin{quote}

\sphinxAtStartPar
product value.
\end{quote}
\end{quote}

\subsubsection{NdArray\textless{}T\textgreater{}.Prod()}
\label{\detokenize{csapi/ndarray:id3}}\begin{quote}

\sphinxAtStartPar
Prodcut of elements at given axis of NdArray object.

\sphinxAtStartPar
\sphinxstylestrong{Synopsis}
\begin{quote}

\sphinxAtStartPar
\sphinxcode{\sphinxupquote{NdArray\textless{}T, N \sphinxhyphen{} 1\textgreater{} Prod(int axis)}}
\end{quote}

\sphinxAtStartPar
\sphinxstylestrong{Arguments}
\begin{quote}

\sphinxAtStartPar
\sphinxcode{\sphinxupquote{axis}}: axis of NdArray.
\end{quote}

\sphinxAtStartPar
\sphinxstylestrong{Return}
\begin{quote}

\sphinxAtStartPar
(N\sphinxhyphen{}1)\sphinxhyphen{}dimensional NdArray object.
\end{quote}
\end{quote}

\subsubsection{NdArray\textless{}T\textgreater{}.Repeat()}
\label{\detokenize{csapi/ndarray:ndarray-t-repeat}}\begin{quote}

\sphinxAtStartPar
Repeat each element of an array along given axis.

\sphinxAtStartPar
\sphinxstylestrong{Synopsis}
\begin{quote}

\sphinxAtStartPar
\sphinxcode{\sphinxupquote{NdArray\textless{}T\textgreater{} Repeat(long repeats, int axis)}}
\end{quote}

\sphinxAtStartPar
\sphinxstylestrong{Arguments}
\begin{quote}

\sphinxAtStartPar
\sphinxcode{\sphinxupquote{repeats}}: number of repetitions for each element.

\sphinxAtStartPar
\sphinxcode{\sphinxupquote{axis}}: axis of NdArray.
\end{quote}

\sphinxAtStartPar
\sphinxstylestrong{Return}
\begin{quote}

\sphinxAtStartPar
new NdArray object.
\end{quote}
\end{quote}

\subsubsection{NdArray\textless{}T\textgreater{}.RepeatBlock()}
\label{\detokenize{csapi/ndarray:ndarray-t-repeatblock}}\begin{quote}

\sphinxAtStartPar
Repeat an array a number of times along given axis.

\sphinxAtStartPar
\sphinxstylestrong{Synopsis}
\begin{quote}

\sphinxAtStartPar
\sphinxcode{\sphinxupquote{NdArray\textless{}T\textgreater{} RepeatBlock(long repeats, int axis)}}
\end{quote}

\sphinxAtStartPar
\sphinxstylestrong{Arguments}
\begin{quote}

\sphinxAtStartPar
\sphinxcode{\sphinxupquote{repeats}}: number of repetitions.

\sphinxAtStartPar
\sphinxcode{\sphinxupquote{axis}}: axis of NdArray.
\end{quote}

\sphinxAtStartPar
\sphinxstylestrong{Return}
\begin{quote}

\sphinxAtStartPar
new NdArray object.
\end{quote}
\end{quote}

\subsubsection{NdArray\textless{}T\textgreater{}.Represent()}
\label{\detokenize{csapi/ndarray:ndarray-t-represent}}\begin{quote}

\sphinxAtStartPar
String representation of NdArray object.

\sphinxAtStartPar
\sphinxstylestrong{Synopsis}
\begin{quote}

\sphinxAtStartPar
\sphinxcode{\sphinxupquote{string Represent(int maxlen)}}
\end{quote}

\sphinxAtStartPar
\sphinxstylestrong{Arguments}
\begin{quote}

\sphinxAtStartPar
\sphinxcode{\sphinxupquote{maxlen}}: maximum buffer length for representations string.
\end{quote}

\sphinxAtStartPar
\sphinxstylestrong{Return}
\begin{quote}

\sphinxAtStartPar
representation string object.
\end{quote}
\end{quote}

\subsubsection{NdArray\textless{}T\textgreater{}.Reshape()}
\label{\detokenize{csapi/ndarray:ndarray-t-reshape}}\begin{quote}

\sphinxAtStartPar
Reshape NdArray object to new shape.

\sphinxAtStartPar
\sphinxstylestrong{Synopsis}
\begin{quote}

\sphinxAtStartPar
\sphinxcode{\sphinxupquote{NdArray\textless{}T\textgreater{} Reshape(Shape shp)}}
\end{quote}

\sphinxAtStartPar
\sphinxstylestrong{Arguments}
\begin{quote}

\sphinxAtStartPar
\sphinxcode{\sphinxupquote{shp}}: new shape of M\sphinxhyphen{}dimensions.
\end{quote}

\sphinxAtStartPar
\sphinxstylestrong{Return}
\begin{quote}

\sphinxAtStartPar
M\sphinxhyphen{}dimensional NdArray object.
\end{quote}
\end{quote}

\subsubsection{NdArray\textless{}T\textgreater{}.SetItem()}
\label{\detokenize{csapi/ndarray:ndarray-t-setitem}}\begin{quote}

\sphinxAtStartPar
Set element of given index to NdArray object.

\sphinxAtStartPar
\sphinxstylestrong{Synopsis}
\begin{quote}

\sphinxAtStartPar
\sphinxcode{\sphinxupquote{void SetItem(long idx, T val)}}
\end{quote}

\sphinxAtStartPar
\sphinxstylestrong{Arguments}
\begin{quote}

\sphinxAtStartPar
\sphinxcode{\sphinxupquote{idx}}: index of element.

\sphinxAtStartPar
\sphinxcode{\sphinxupquote{val}}: value of element.
\end{quote}
\end{quote}

\subsubsection{NdArray\textless{}T\textgreater{}.Squeeze()}
\label{\detokenize{csapi/ndarray:ndarray-t-squeeze}}\begin{quote}

\sphinxAtStartPar
Remove axis of length 1 from shape of NdArray object.

\sphinxAtStartPar
\sphinxstylestrong{Synopsis}
\begin{quote}

\sphinxAtStartPar
\sphinxcode{\sphinxupquote{NdArray\textless{}T\textgreater{} Squeeze(int axis)}}
\end{quote}

\sphinxAtStartPar
\sphinxstylestrong{Arguments}
\begin{quote}

\sphinxAtStartPar
\sphinxcode{\sphinxupquote{axis}}: axis of NdArray, where the length is 1.
\end{quote}

\sphinxAtStartPar
\sphinxstylestrong{Return}
\begin{quote}

\sphinxAtStartPar
(N\sphinxhyphen{}1)\sphinxhyphen{}dimensional NdArray object.
\end{quote}
\end{quote}

\subsubsection{NdArray\textless{}T\textgreater{}.Sum()}
\label{\detokenize{csapi/ndarray:ndarray-t-sum}}\begin{quote}

\sphinxAtStartPar
Sum of all elements in NdArray object.

\sphinxAtStartPar
\sphinxstylestrong{Synopsis}
\begin{quote}

\sphinxAtStartPar
\sphinxcode{\sphinxupquote{T Sum()}}
\end{quote}

\sphinxAtStartPar
\sphinxstylestrong{Return}
\begin{quote}

\sphinxAtStartPar
sum value.
\end{quote}
\end{quote}

\subsubsection{NdArray\textless{}T\textgreater{}.Sum()}
\label{\detokenize{csapi/ndarray:id4}}\begin{quote}

\sphinxAtStartPar
Sum of elements at given axis of NdArray object.

\sphinxAtStartPar
\sphinxstylestrong{Synopsis}
\begin{quote}

\sphinxAtStartPar
\sphinxcode{\sphinxupquote{NdArray\textless{}T, N \sphinxhyphen{} 1\textgreater{} Sum(int axis)}}
\end{quote}

\sphinxAtStartPar
\sphinxstylestrong{Arguments}
\begin{quote}

\sphinxAtStartPar
\sphinxcode{\sphinxupquote{axis}}: axis of NdArray.
\end{quote}

\sphinxAtStartPar
\sphinxstylestrong{Return}
\begin{quote}

\sphinxAtStartPar
(N\sphinxhyphen{}1)\sphinxhyphen{}dimensional NdArray object.
\end{quote}
\end{quote}

\subsubsection{NdArray\textless{}T\textgreater{}.Transpose()}
\label{\detokenize{csapi/ndarray:ndarray-t-transpose}}\begin{quote}

\sphinxAtStartPar
Perform matrix transpose of NdArray object.

\sphinxAtStartPar
\sphinxstylestrong{Synopsis}
\begin{quote}

\sphinxAtStartPar
\sphinxcode{\sphinxupquote{NdArray\textless{}T\textgreater{} Transpose()}}
\end{quote}

\sphinxAtStartPar
\sphinxstylestrong{Return}
\begin{quote}

\sphinxAtStartPar
transposed NdArray object.
\end{quote}
\end{quote}

\subsection{Shape}
\label{\detokenize{csharpapiref:shape}}\label{\detokenize{csharpapiref:chapcsharpapiref-shape}}
\sphinxAtStartPar
The Shape class encapsulates a tuple of integers, indicating the size of array
along each dimension. It refers to dimensions of built\sphinxhyphen{}in {\hyperref[\detokenize{csharpapiref:chapcsharpapiref-ndarray}]{\sphinxcrossref{\DUrole{std,std-ref}{NdArray}}}}
in COPT. The following methods are provided:

\sphinxstepscope

\subsubsection{Shape.Shape()}
\label{\detokenize{csapi/shape:shape-shape}}\label{\detokenize{csapi/shape::doc}}\begin{quote}

\sphinxAtStartPar
Constructor of Shape object.

\sphinxAtStartPar
\sphinxstylestrong{Synopsis}
\begin{quote}

\sphinxAtStartPar
\sphinxcode{\sphinxupquote{Shape()}}
\end{quote}
\end{quote}

\subsubsection{Shape.Expand()}
\label{\detokenize{csapi/shape:shape-expand}}\begin{quote}

\sphinxAtStartPar
Expand shape of Shape object.

\sphinxAtStartPar
\sphinxstylestrong{Synopsis}
\begin{quote}

\sphinxAtStartPar
\sphinxcode{\sphinxupquote{Shape Expand(int axis)}}
\end{quote}

\sphinxAtStartPar
\sphinxstylestrong{Arguments}
\begin{quote}

\sphinxAtStartPar
\sphinxcode{\sphinxupquote{axis}}: given axis.
\end{quote}

\sphinxAtStartPar
\sphinxstylestrong{Return}
\begin{quote}

\sphinxAtStartPar
Shape object in (N+1)\sphinxhyphen{}dimensions.
\end{quote}
\end{quote}

\subsubsection{Shape.GetDim()}
\label{\detokenize{csapi/shape:shape-getdim}}\begin{quote}

\sphinxAtStartPar
Get i\sphinxhyphen{}th dimension in Shape object.

\sphinxAtStartPar
\sphinxstylestrong{Synopsis}
\begin{quote}

\sphinxAtStartPar
\sphinxcode{\sphinxupquote{long GetDim(int i)}}
\end{quote}

\sphinxAtStartPar
\sphinxstylestrong{Arguments}
\begin{quote}

\sphinxAtStartPar
\sphinxcode{\sphinxupquote{i}}: index of dimensions.
\end{quote}

\sphinxAtStartPar
\sphinxstylestrong{Return}
\begin{quote}

\sphinxAtStartPar
the i\sphinxhyphen{}th dimension.
\end{quote}
\end{quote}

\subsubsection{Shape.GetND()}
\label{\detokenize{csapi/shape:shape-getnd}}\begin{quote}

\sphinxAtStartPar
Get number of dimensions in Shape object.

\sphinxAtStartPar
\sphinxstylestrong{Synopsis}
\begin{quote}

\sphinxAtStartPar
\sphinxcode{\sphinxupquote{int GetND()}}
\end{quote}

\sphinxAtStartPar
\sphinxstylestrong{Return}
\begin{quote}

\sphinxAtStartPar
number of dimensions.
\end{quote}
\end{quote}

\subsubsection{Shape.GetSize()}
\label{\detokenize{csapi/shape:shape-getsize}}\begin{quote}

\sphinxAtStartPar
Get size of Shape object.

\sphinxAtStartPar
\sphinxstylestrong{Synopsis}
\begin{quote}

\sphinxAtStartPar
\sphinxcode{\sphinxupquote{long GetSize()}}
\end{quote}

\sphinxAtStartPar
\sphinxstylestrong{Return}
\begin{quote}

\sphinxAtStartPar
size of shape.
\end{quote}
\end{quote}

\subsubsection{Shape.GetStart()}
\label{\detokenize{csapi/shape:shape-getstart}}\begin{quote}

\sphinxAtStartPar
Get the i\sphinxhyphen{}th start postion in Shape object.

\sphinxAtStartPar
\sphinxstylestrong{Synopsis}
\begin{quote}

\sphinxAtStartPar
\sphinxcode{\sphinxupquote{int GetStart(int i)}}
\end{quote}

\sphinxAtStartPar
\sphinxstylestrong{Arguments}
\begin{quote}

\sphinxAtStartPar
\sphinxcode{\sphinxupquote{i}}: index of dimensions.
\end{quote}

\sphinxAtStartPar
\sphinxstylestrong{Return}
\begin{quote}

\sphinxAtStartPar
start position in i\sphinxhyphen{}th dimension.
\end{quote}
\end{quote}

\subsubsection{Shape.GetStride()}
\label{\detokenize{csapi/shape:shape-getstride}}\begin{quote}

\sphinxAtStartPar
Get i\sphinxhyphen{}th stride in Shape object.

\sphinxAtStartPar
\sphinxstylestrong{Synopsis}
\begin{quote}

\sphinxAtStartPar
\sphinxcode{\sphinxupquote{int GetStride(int i)}}
\end{quote}

\sphinxAtStartPar
\sphinxstylestrong{Arguments}
\begin{quote}

\sphinxAtStartPar
\sphinxcode{\sphinxupquote{i}}: index of dimensions.
\end{quote}

\sphinxAtStartPar
\sphinxstylestrong{Return}
\begin{quote}

\sphinxAtStartPar
stride in i\sphinxhyphen{}th dimension.
\end{quote}
\end{quote}

\subsubsection{Shape.Rebuild()}
\label{\detokenize{csapi/shape:shape-rebuild}}\begin{quote}

\sphinxAtStartPar
Rebuild Shape object, that is, keep dimensions while reset strides and starts.

\sphinxAtStartPar
\sphinxstylestrong{Synopsis}
\begin{quote}

\sphinxAtStartPar
\sphinxcode{\sphinxupquote{Shape Rebuild()}}
\end{quote}

\sphinxAtStartPar
\sphinxstylestrong{Return}
\begin{quote}

\sphinxAtStartPar
new Shape object.
\end{quote}
\end{quote}

\subsubsection{Shape.Represent()}
\label{\detokenize{csapi/shape:shape-represent}}\begin{quote}

\sphinxAtStartPar
String representation of Shape object.

\sphinxAtStartPar
\sphinxstylestrong{Synopsis}
\begin{quote}

\sphinxAtStartPar
\sphinxcode{\sphinxupquote{string Represent(int type)}}
\end{quote}

\sphinxAtStartPar
\sphinxstylestrong{Arguments}
\begin{quote}

\sphinxAtStartPar
\sphinxcode{\sphinxupquote{type}}: 0: dimensions; 1: strides; 2: starts.
\end{quote}

\sphinxAtStartPar
\sphinxstylestrong{Return}
\begin{quote}

\sphinxAtStartPar
string object.
\end{quote}
\end{quote}

\subsubsection{Shape.Squeeze()}
\label{\detokenize{csapi/shape:shape-squeeze}}\begin{quote}

\sphinxAtStartPar
Remove axis of length 1 from Shape object.

\sphinxAtStartPar
\sphinxstylestrong{Synopsis}
\begin{quote}

\sphinxAtStartPar
\sphinxcode{\sphinxupquote{Shape Squeeze(int axis)}}
\end{quote}

\sphinxAtStartPar
\sphinxstylestrong{Arguments}
\begin{quote}

\sphinxAtStartPar
\sphinxcode{\sphinxupquote{axis}}: given axis, where the length is 1.
\end{quote}

\sphinxAtStartPar
\sphinxstylestrong{Return}
\begin{quote}

\sphinxAtStartPar
Shape object in (N\sphinxhyphen{}1)\sphinxhyphen{}dimensions.
\end{quote}
\end{quote}

\subsection{View}
\label{\detokenize{csharpapiref:view}}\label{\detokenize{csharpapiref:chapcsharpapiref-view}}
\sphinxAtStartPar
The View class is used to perform slicing operations on multi\sphinxhyphen{}dimensional arrays.
The following methods are provided:

\sphinxstepscope

\subsubsection{View.View()}
\label{\detokenize{csapi/view:view-view}}\label{\detokenize{csapi/view::doc}}\begin{quote}

\sphinxAtStartPar
Constructor of View object.

\sphinxAtStartPar
\sphinxstylestrong{Synopsis}
\begin{quote}

\sphinxAtStartPar
\sphinxcode{\sphinxupquote{View()}}
\end{quote}
\end{quote}

\subsubsection{View.AddFull()}
\label{\detokenize{csapi/view:view-addfull}}\begin{quote}

\sphinxAtStartPar
Create full view object at current dimension.

\sphinxAtStartPar
\sphinxstylestrong{Synopsis}
\begin{quote}

\sphinxAtStartPar
\sphinxcode{\sphinxupquote{View AddFull()}}
\end{quote}

\sphinxAtStartPar
\sphinxstylestrong{Return}
\begin{quote}

\sphinxAtStartPar
View object.
\end{quote}
\end{quote}

\subsubsection{View.AddScalar()}
\label{\detokenize{csapi/view:view-addscalar}}\begin{quote}

\sphinxAtStartPar
Create View object of given index at current dimension.

\sphinxAtStartPar
\sphinxstylestrong{Synopsis}
\begin{quote}

\sphinxAtStartPar
\sphinxcode{\sphinxupquote{View AddScalar(long n)}}
\end{quote}

\sphinxAtStartPar
\sphinxstylestrong{Arguments}
\begin{quote}

\sphinxAtStartPar
\sphinxcode{\sphinxupquote{n}}: given index.
\end{quote}

\sphinxAtStartPar
\sphinxstylestrong{Return}
\begin{quote}

\sphinxAtStartPar
View object.
\end{quote}
\end{quote}

\subsubsection{View.AddSlice()}
\label{\detokenize{csapi/view:view-addslice}}\begin{quote}

\sphinxAtStartPar
Create view object of slice at current dimension.

\sphinxAtStartPar
\sphinxstylestrong{Synopsis}
\begin{quote}

\sphinxAtStartPar
\sphinxcode{\sphinxupquote{View AddSlice(long start)}}
\end{quote}

\sphinxAtStartPar
\sphinxstylestrong{Arguments}
\begin{quote}

\sphinxAtStartPar
\sphinxcode{\sphinxupquote{start}}: start index, inclusive.
\end{quote}

\sphinxAtStartPar
\sphinxstylestrong{Return}
\begin{quote}

\sphinxAtStartPar
View object.
\end{quote}
\end{quote}

\subsubsection{View.AddSlice()}
\label{\detokenize{csapi/view:id1}}\begin{quote}

\sphinxAtStartPar
Create view object of slice at current dimension.

\sphinxAtStartPar
\sphinxstylestrong{Synopsis}
\begin{quote}

\sphinxAtStartPar
\sphinxcode{\sphinxupquote{View AddSlice(long start, long stop)}}
\end{quote}

\sphinxAtStartPar
\sphinxstylestrong{Arguments}
\begin{quote}

\sphinxAtStartPar
\sphinxcode{\sphinxupquote{start}}: start index, inclusive.

\sphinxAtStartPar
\sphinxcode{\sphinxupquote{stop}}: stop index, exclusive.
\end{quote}

\sphinxAtStartPar
\sphinxstylestrong{Return}
\begin{quote}

\sphinxAtStartPar
View object.
\end{quote}
\end{quote}

\subsubsection{View.AddSlice()}
\label{\detokenize{csapi/view:id2}}\begin{quote}

\sphinxAtStartPar
Create view object of slice at current dimension.

\sphinxAtStartPar
\sphinxstylestrong{Synopsis}
\begin{quote}

\sphinxAtStartPar
\sphinxcode{\sphinxupquote{View AddSlice(}}
\begin{quote}

\sphinxAtStartPar
\sphinxcode{\sphinxupquote{long start,}}

\sphinxAtStartPar
\sphinxcode{\sphinxupquote{long stop,}}

\sphinxAtStartPar
\sphinxcode{\sphinxupquote{long step,}}

\sphinxAtStartPar
\sphinxcode{\sphinxupquote{int flag)}}
\end{quote}
\end{quote}

\sphinxAtStartPar
\sphinxstylestrong{Arguments}
\begin{quote}

\sphinxAtStartPar
\sphinxcode{\sphinxupquote{start}}: start index, inclusive.

\sphinxAtStartPar
\sphinxcode{\sphinxupquote{stop}}: stop index, exclusive.

\sphinxAtStartPar
\sphinxcode{\sphinxupquote{step}}: step size between start and stop index. It can be negative.

\sphinxAtStartPar
\sphinxcode{\sphinxupquote{flag}}: optional, flag for slicing type. Default is 0.
\end{quote}

\sphinxAtStartPar
\sphinxstylestrong{Return}
\begin{quote}

\sphinxAtStartPar
View object.
\end{quote}
\end{quote}

\subsection{NlExpr Class}
\label{\detokenize{csharpapiref:nlexpr-class}}\label{\detokenize{csharpapiref:chapcsharpapiref-nlexpr}}
\sphinxAtStartPar
COPT nonlinear expression object.
The \sphinxcode{\sphinxupquote{NlExpr}} class represents nonlinear expressions in COPT. The nonlinear expressions
are used to build nonlinear constraints. The following methods are provided:

\sphinxstepscope

\subsubsection{NlExpr.NlExpr()}
\label{\detokenize{csapi/nlexpr:nlexpr-nlexpr}}\label{\detokenize{csapi/nlexpr::doc}}\begin{quote}

\sphinxAtStartPar
Constructor of a nonlinear expression with a constant.

\sphinxAtStartPar
\sphinxstylestrong{Synopsis}
\begin{quote}

\sphinxAtStartPar
\sphinxcode{\sphinxupquote{NlExpr(double constant)}}
\end{quote}

\sphinxAtStartPar
\sphinxstylestrong{Arguments}
\begin{quote}

\sphinxAtStartPar
\sphinxcode{\sphinxupquote{constant}}: constant value in nonlinear expression object.
\end{quote}
\end{quote}

\subsubsection{NlExpr.NlExpr()}
\label{\detokenize{csapi/nlexpr:id1}}\begin{quote}

\sphinxAtStartPar
Constructor of a nonlinear expression with one linear term.

\sphinxAtStartPar
\sphinxstylestrong{Synopsis}
\begin{quote}

\sphinxAtStartPar
\sphinxcode{\sphinxupquote{NlExpr(Var var, double coeff)}}
\end{quote}

\sphinxAtStartPar
\sphinxstylestrong{Arguments}
\begin{quote}

\sphinxAtStartPar
\sphinxcode{\sphinxupquote{var}}: variable for the added term.

\sphinxAtStartPar
\sphinxcode{\sphinxupquote{coeff}}: optional, coefficent for the added term.
\end{quote}
\end{quote}

\subsubsection{NlExpr.NlExpr()}
\label{\detokenize{csapi/nlexpr:id2}}\begin{quote}

\sphinxAtStartPar
Constructor of a nonlinear expression with a linear expression.

\sphinxAtStartPar
\sphinxstylestrong{Synopsis}
\begin{quote}

\sphinxAtStartPar
\sphinxcode{\sphinxupquote{NlExpr(Expr expr)}}
\end{quote}

\sphinxAtStartPar
\sphinxstylestrong{Arguments}
\begin{quote}

\sphinxAtStartPar
\sphinxcode{\sphinxupquote{expr}}: linear expression.
\end{quote}
\end{quote}

\subsubsection{NlExpr.NlExpr()}
\label{\detokenize{csapi/nlexpr:id3}}\begin{quote}

\sphinxAtStartPar
Constructor of a nonlinear expression with a quadratic expression.

\sphinxAtStartPar
\sphinxstylestrong{Synopsis}
\begin{quote}

\sphinxAtStartPar
\sphinxcode{\sphinxupquote{NlExpr(QuadExpr expr)}}
\end{quote}

\sphinxAtStartPar
\sphinxstylestrong{Arguments}
\begin{quote}

\sphinxAtStartPar
\sphinxcode{\sphinxupquote{expr}}: quadratic expression.
\end{quote}
\end{quote}

\subsubsection{NlExpr.AddConstant()}
\label{\detokenize{csapi/nlexpr:nlexpr-addconstant}}\begin{quote}

\sphinxAtStartPar
Add constant to the nonlinear expression.

\sphinxAtStartPar
\sphinxstylestrong{Synopsis}
\begin{quote}

\sphinxAtStartPar
\sphinxcode{\sphinxupquote{void AddConstant(double constant)}}
\end{quote}

\sphinxAtStartPar
\sphinxstylestrong{Arguments}
\begin{quote}

\sphinxAtStartPar
\sphinxcode{\sphinxupquote{constant}}: value to be added.
\end{quote}
\end{quote}

\subsubsection{NlExpr.AddLinExpr()}
\label{\detokenize{csapi/nlexpr:nlexpr-addlinexpr}}\begin{quote}

\sphinxAtStartPar
Add a linear expression to self.

\sphinxAtStartPar
\sphinxstylestrong{Synopsis}
\begin{quote}

\sphinxAtStartPar
\sphinxcode{\sphinxupquote{void AddLinExpr(Expr expr, double mult)}}
\end{quote}

\sphinxAtStartPar
\sphinxstylestrong{Arguments}
\begin{quote}

\sphinxAtStartPar
\sphinxcode{\sphinxupquote{expr}}: linear expression to be added.

\sphinxAtStartPar
\sphinxcode{\sphinxupquote{mult}}: optional, constant multiplier, default value is 1.0.
\end{quote}
\end{quote}

\subsubsection{NlExpr.AddNlExpr()}
\label{\detokenize{csapi/nlexpr:nlexpr-addnlexpr}}\begin{quote}

\sphinxAtStartPar
Add a nonlinear expression to self.

\sphinxAtStartPar
\sphinxstylestrong{Synopsis}
\begin{quote}

\sphinxAtStartPar
\sphinxcode{\sphinxupquote{void AddNlExpr(NlExpr expr, double mult)}}
\end{quote}

\sphinxAtStartPar
\sphinxstylestrong{Arguments}
\begin{quote}

\sphinxAtStartPar
\sphinxcode{\sphinxupquote{expr}}: nonlinear expression to be added.

\sphinxAtStartPar
\sphinxcode{\sphinxupquote{mult}}: optional, constant multiplier, default value is 1.0.
\end{quote}
\end{quote}

\subsubsection{NlExpr.AddQuadExpr()}
\label{\detokenize{csapi/nlexpr:nlexpr-addquadexpr}}\begin{quote}

\sphinxAtStartPar
Add a quadratic expression to self.

\sphinxAtStartPar
\sphinxstylestrong{Synopsis}
\begin{quote}

\sphinxAtStartPar
\sphinxcode{\sphinxupquote{void AddQuadExpr(QuadExpr expr, double mult)}}
\end{quote}

\sphinxAtStartPar
\sphinxstylestrong{Arguments}
\begin{quote}

\sphinxAtStartPar
\sphinxcode{\sphinxupquote{expr}}: quadratic expression to be added.

\sphinxAtStartPar
\sphinxcode{\sphinxupquote{mult}}: optional, constant multiplier, default value is 1.0.
\end{quote}
\end{quote}

\subsubsection{NlExpr.AddTerm()}
\label{\detokenize{csapi/nlexpr:nlexpr-addterm}}\begin{quote}

\sphinxAtStartPar
Add a linear term to nonlinear expression object.

\sphinxAtStartPar
\sphinxstylestrong{Synopsis}
\begin{quote}

\sphinxAtStartPar
\sphinxcode{\sphinxupquote{void AddTerm(Var var, double coeff)}}
\end{quote}

\sphinxAtStartPar
\sphinxstylestrong{Arguments}
\begin{quote}

\sphinxAtStartPar
\sphinxcode{\sphinxupquote{var}}: variable of new linear term.

\sphinxAtStartPar
\sphinxcode{\sphinxupquote{coeff}}: optional, coefficient of new linear term.
\end{quote}
\end{quote}

\subsubsection{NlExpr.AddTerms()}
\label{\detokenize{csapi/nlexpr:nlexpr-addterms}}\begin{quote}

\sphinxAtStartPar
Add linear terms to nonlinear expression object.

\sphinxAtStartPar
\sphinxstylestrong{Synopsis}
\begin{quote}

\sphinxAtStartPar
\sphinxcode{\sphinxupquote{void AddTerms(VarArray vars, double{[}{]} coeffs)}}
\end{quote}

\sphinxAtStartPar
\sphinxstylestrong{Arguments}
\begin{quote}

\sphinxAtStartPar
\sphinxcode{\sphinxupquote{vars}}: variables for added linear terms.

\sphinxAtStartPar
\sphinxcode{\sphinxupquote{coeffs}}: coefficient array for added linear terms.
\end{quote}
\end{quote}

\subsubsection{NlExpr.AddTerms()}
\label{\detokenize{csapi/nlexpr:id4}}\begin{quote}

\sphinxAtStartPar
Add linear terms to nonlinear expression object.

\sphinxAtStartPar
\sphinxstylestrong{Synopsis}
\begin{quote}

\sphinxAtStartPar
\sphinxcode{\sphinxupquote{void AddTerms(Var{[}{]} vars, double{[}{]} coeffs)}}
\end{quote}

\sphinxAtStartPar
\sphinxstylestrong{Arguments}
\begin{quote}

\sphinxAtStartPar
\sphinxcode{\sphinxupquote{vars}}: variable array for added linear terms.

\sphinxAtStartPar
\sphinxcode{\sphinxupquote{coeffs}}: coefficient array for added linear terms.
\end{quote}
\end{quote}

\subsubsection{NlExpr.Clear()}
\label{\detokenize{csapi/nlexpr:nlexpr-clear}}\begin{quote}

\sphinxAtStartPar
Clear nonlinear expression object.

\sphinxAtStartPar
\sphinxstylestrong{Synopsis}
\begin{quote}

\sphinxAtStartPar
\sphinxcode{\sphinxupquote{void Clear()}}
\end{quote}
\end{quote}

\subsubsection{NlExpr.Clone()}
\label{\detokenize{csapi/nlexpr:nlexpr-clone}}\begin{quote}

\sphinxAtStartPar
Deep copy nonlinear expression object.

\sphinxAtStartPar
\sphinxstylestrong{Synopsis}
\begin{quote}

\sphinxAtStartPar
\sphinxcode{\sphinxupquote{NlExpr Clone()}}
\end{quote}

\sphinxAtStartPar
\sphinxstylestrong{Return}
\begin{quote}

\sphinxAtStartPar
cloned nonlinear expression object.
\end{quote}
\end{quote}

\subsubsection{NlExpr.Divide()}
\label{\detokenize{csapi/nlexpr:nlexpr-divide}}\begin{quote}

\sphinxAtStartPar
Divide itself by an expression.

\sphinxAtStartPar
\sphinxstylestrong{Synopsis}
\begin{quote}

\sphinxAtStartPar
\sphinxcode{\sphinxupquote{void Divide(NlExpr expr)}}
\end{quote}

\sphinxAtStartPar
\sphinxstylestrong{Arguments}
\begin{quote}

\sphinxAtStartPar
\sphinxcode{\sphinxupquote{expr}}: expression operand, including NlExpr, QuadExpr, Expr, Var and constant.
\end{quote}
\end{quote}

\subsubsection{NlExpr.Evaluate()}
\label{\detokenize{csapi/nlexpr:nlexpr-evaluate}}\begin{quote}

\sphinxAtStartPar
Evaluate nonlinear expression after solving.

\sphinxAtStartPar
\sphinxstylestrong{Synopsis}
\begin{quote}

\sphinxAtStartPar
\sphinxcode{\sphinxupquote{double Evaluate()}}
\end{quote}

\sphinxAtStartPar
\sphinxstylestrong{Return}
\begin{quote}

\sphinxAtStartPar
value of nonlinear expression.
\end{quote}
\end{quote}

\subsubsection{NlExpr.GetConstant()}
\label{\detokenize{csapi/nlexpr:nlexpr-getconstant}}\begin{quote}

\sphinxAtStartPar
Get constant in nonlinear expression.

\sphinxAtStartPar
\sphinxstylestrong{Synopsis}
\begin{quote}

\sphinxAtStartPar
\sphinxcode{\sphinxupquote{double GetConstant()}}
\end{quote}

\sphinxAtStartPar
\sphinxstylestrong{Return}
\begin{quote}

\sphinxAtStartPar
constant in nonlinear expression.
\end{quote}
\end{quote}

\subsubsection{NlExpr.GetLinExpr()}
\label{\detokenize{csapi/nlexpr:nlexpr-getlinexpr}}\begin{quote}

\sphinxAtStartPar
Get linear expression of nonlinear expression.

\sphinxAtStartPar
\sphinxstylestrong{Synopsis}
\begin{quote}

\sphinxAtStartPar
\sphinxcode{\sphinxupquote{Expr GetLinExpr()}}
\end{quote}

\sphinxAtStartPar
\sphinxstylestrong{Return}
\begin{quote}

\sphinxAtStartPar
linear expression object.
\end{quote}
\end{quote}

\subsubsection{NlExpr.Multiply()}
\label{\detokenize{csapi/nlexpr:nlexpr-multiply}}\begin{quote}

\sphinxAtStartPar
Multiply itself by an expression.

\sphinxAtStartPar
\sphinxstylestrong{Synopsis}
\begin{quote}

\sphinxAtStartPar
\sphinxcode{\sphinxupquote{void Multiply(NlExpr expr)}}
\end{quote}

\sphinxAtStartPar
\sphinxstylestrong{Arguments}
\begin{quote}

\sphinxAtStartPar
\sphinxcode{\sphinxupquote{expr}}: expression operand, including NlExpr, QuadExpr, Expr, Var and constant.
\end{quote}
\end{quote}

\subsubsection{NlExpr.Negate()}
\label{\detokenize{csapi/nlexpr:nlexpr-negate}}\begin{quote}

\sphinxAtStartPar
Negate itself.

\sphinxAtStartPar
\sphinxstylestrong{Synopsis}
\begin{quote}

\sphinxAtStartPar
\sphinxcode{\sphinxupquote{void Negate()}}
\end{quote}
\end{quote}

\subsubsection{NlExpr.Reserve()}
\label{\detokenize{csapi/nlexpr:nlexpr-reserve}}\begin{quote}

\sphinxAtStartPar
Reserve capacity to contain at least n items.

\sphinxAtStartPar
\sphinxstylestrong{Synopsis}
\begin{quote}

\sphinxAtStartPar
\sphinxcode{\sphinxupquote{void Reserve(int n)}}
\end{quote}

\sphinxAtStartPar
\sphinxstylestrong{Arguments}
\begin{quote}

\sphinxAtStartPar
\sphinxcode{\sphinxupquote{n}}: capacity of nonlinear constraint objects.
\end{quote}
\end{quote}

\subsubsection{NlExpr.SetConstant()}
\label{\detokenize{csapi/nlexpr:nlexpr-setconstant}}\begin{quote}

\sphinxAtStartPar
Set constant for the nonlinear expression.

\sphinxAtStartPar
\sphinxstylestrong{Synopsis}
\begin{quote}

\sphinxAtStartPar
\sphinxcode{\sphinxupquote{void SetConstant(double constant)}}
\end{quote}

\sphinxAtStartPar
\sphinxstylestrong{Arguments}
\begin{quote}

\sphinxAtStartPar
\sphinxcode{\sphinxupquote{constant}}: the value of the constant.
\end{quote}
\end{quote}

\subsubsection{NlExpr.Size()}
\label{\detokenize{csapi/nlexpr:nlexpr-size}}\begin{quote}

\sphinxAtStartPar
Get size of tokens in nonlinear expression.

\sphinxAtStartPar
\sphinxstylestrong{Synopsis}
\begin{quote}

\sphinxAtStartPar
\sphinxcode{\sphinxupquote{long Size()}}
\end{quote}

\sphinxAtStartPar
\sphinxstylestrong{Return}
\begin{quote}

\sphinxAtStartPar
size of none\sphinxhyphen{}linear tokens.
\end{quote}
\end{quote}

\subsection{NlConstraint Class}
\label{\detokenize{csharpapiref:nlconstraint-class}}\label{\detokenize{csharpapiref:chapcsharpapiref-nlconstraint}}
\sphinxAtStartPar
COPT nonlinear constraint object. The \sphinxcode{\sphinxupquote{NlConstraint}} object is always associated with a
particular model. User creates a \sphinxcode{\sphinxupquote{NlConstraint}} object by adding a nonlinear
constraint to model, rather than by constructor of \sphinxcode{\sphinxupquote{NlConstraint}} class.

\sphinxstepscope

\subsubsection{NlConstraint.Get()}
\label{\detokenize{csapi/nlconstraint:nlconstraint-get}}\label{\detokenize{csapi/nlconstraint::doc}}\begin{quote}

\sphinxAtStartPar
Get information value of the nonlinear constraint. Support informations of “LB”, “UB”, “Slack”.

\sphinxAtStartPar
\sphinxstylestrong{Synopsis}
\begin{quote}

\sphinxAtStartPar
\sphinxcode{\sphinxupquote{double Get(string info)}}
\end{quote}

\sphinxAtStartPar
\sphinxstylestrong{Arguments}
\begin{quote}

\sphinxAtStartPar
\sphinxcode{\sphinxupquote{info}}: name of the information being queried.
\end{quote}

\sphinxAtStartPar
\sphinxstylestrong{Return}
\begin{quote}

\sphinxAtStartPar
value of information.
\end{quote}
\end{quote}

\subsubsection{NlConstraint.GetIdx()}
\label{\detokenize{csapi/nlconstraint:nlconstraint-getidx}}\begin{quote}

\sphinxAtStartPar
Get index of nonlinear constraint.

\sphinxAtStartPar
\sphinxstylestrong{Synopsis}
\begin{quote}

\sphinxAtStartPar
\sphinxcode{\sphinxupquote{int GetIdx()}}
\end{quote}

\sphinxAtStartPar
\sphinxstylestrong{Return}
\begin{quote}

\sphinxAtStartPar
the index of nonlinear constraint.
\end{quote}
\end{quote}

\subsubsection{NlConstraint.GetName()}
\label{\detokenize{csapi/nlconstraint:nlconstraint-getname}}\begin{quote}

\sphinxAtStartPar
Get name of nonlinear constraint.

\sphinxAtStartPar
\sphinxstylestrong{Synopsis}
\begin{quote}

\sphinxAtStartPar
\sphinxcode{\sphinxupquote{string GetName()}}
\end{quote}

\sphinxAtStartPar
\sphinxstylestrong{Return}
\begin{quote}

\sphinxAtStartPar
the name of nonlinear constraint.
\end{quote}
\end{quote}

\subsubsection{NlConstraint.Remove()}
\label{\detokenize{csapi/nlconstraint:nlconstraint-remove}}\begin{quote}

\sphinxAtStartPar
Remove this nonlinear constraint from model.

\sphinxAtStartPar
\sphinxstylestrong{Synopsis}
\begin{quote}

\sphinxAtStartPar
\sphinxcode{\sphinxupquote{void Remove()}}
\end{quote}
\end{quote}

\subsubsection{NlConstraint.Set()}
\label{\detokenize{csapi/nlconstraint:nlconstraint-set}}\begin{quote}

\sphinxAtStartPar
Set information value of nonlinear constraint. Support informations of “LB” and “UB”.

\sphinxAtStartPar
\sphinxstylestrong{Synopsis}
\begin{quote}

\sphinxAtStartPar
\sphinxcode{\sphinxupquote{void Set(string info, double val)}}
\end{quote}

\sphinxAtStartPar
\sphinxstylestrong{Arguments}
\begin{quote}

\sphinxAtStartPar
\sphinxcode{\sphinxupquote{info}}: name of the information.

\sphinxAtStartPar
\sphinxcode{\sphinxupquote{val}}: new information value.
\end{quote}
\end{quote}

\subsubsection{NlConstraint.SetName()}
\label{\detokenize{csapi/nlconstraint:nlconstraint-setname}}\begin{quote}

\sphinxAtStartPar
Set name for nonlinear constraint.

\sphinxAtStartPar
\sphinxstylestrong{Synopsis}
\begin{quote}

\sphinxAtStartPar
\sphinxcode{\sphinxupquote{void SetName(string name)}}
\end{quote}

\sphinxAtStartPar
\sphinxstylestrong{Arguments}
\begin{quote}

\sphinxAtStartPar
\sphinxcode{\sphinxupquote{name}}: the name to set.
\end{quote}
\end{quote}

\subsection{NlConstrArray Class}
\label{\detokenize{csharpapiref:nlconstrarray-class}}\label{\detokenize{csharpapiref:chapcsharpapiref-nlconstrarray}}
\sphinxAtStartPar
COPT nonlinear constraint array object. To store and access a set of
{\hyperref[\detokenize{csharpapiref:chapcsharpapiref-nlconstraint}]{\sphinxcrossref{\DUrole{std,std-ref}{NlConstraint Class}}}} objects, Cardinal Optimizer provides
NlConstrArray class, which defines the following methods.

\sphinxstepscope

\subsubsection{NlConstrArray.NlConstrArray()}
\label{\detokenize{csapi/nlconstrarray:nlconstrarray-nlconstrarray}}\label{\detokenize{csapi/nlconstrarray::doc}}\begin{quote}

\sphinxAtStartPar
Constructor of NlConstrArray object.

\sphinxAtStartPar
\sphinxstylestrong{Synopsis}
\begin{quote}

\sphinxAtStartPar
\sphinxcode{\sphinxupquote{NlConstrArray()}}
\end{quote}
\end{quote}

\subsubsection{NlConstrArray.GetNlConstr()}
\label{\detokenize{csapi/nlconstrarray:nlconstrarray-getnlconstr}}\begin{quote}

\sphinxAtStartPar
Get idx\sphinxhyphen{}th nonlinear constraint object.

\sphinxAtStartPar
\sphinxstylestrong{Synopsis}
\begin{quote}

\sphinxAtStartPar
\sphinxcode{\sphinxupquote{NlConstraint GetNlConstr(int idx)}}
\end{quote}

\sphinxAtStartPar
\sphinxstylestrong{Arguments}
\begin{quote}

\sphinxAtStartPar
\sphinxcode{\sphinxupquote{idx}}: index of the nonlinear constraint.
\end{quote}

\sphinxAtStartPar
\sphinxstylestrong{Return}
\begin{quote}

\sphinxAtStartPar
nonlinear constraint object with index value.
\end{quote}
\end{quote}

\subsubsection{NlConstrArray.PushBack()}
\label{\detokenize{csapi/nlconstrarray:nlconstrarray-pushback}}\begin{quote}

\sphinxAtStartPar
Add a nonlinear constraint to nonlinear constraint array.

\sphinxAtStartPar
\sphinxstylestrong{Synopsis}
\begin{quote}

\sphinxAtStartPar
\sphinxcode{\sphinxupquote{void PushBack(NlConstraint constr)}}
\end{quote}

\sphinxAtStartPar
\sphinxstylestrong{Arguments}
\begin{quote}

\sphinxAtStartPar
\sphinxcode{\sphinxupquote{constr}}: nonlinear constraint object.
\end{quote}
\end{quote}

\subsubsection{NlConstrArray.Reserve()}
\label{\detokenize{csapi/nlconstrarray:nlconstrarray-reserve}}\begin{quote}

\sphinxAtStartPar
Reserve capacity to contain at least n items.

\sphinxAtStartPar
\sphinxstylestrong{Synopsis}
\begin{quote}

\sphinxAtStartPar
\sphinxcode{\sphinxupquote{void Reserve(int n)}}
\end{quote}

\sphinxAtStartPar
\sphinxstylestrong{Arguments}
\begin{quote}

\sphinxAtStartPar
\sphinxcode{\sphinxupquote{n}}: capacity of nonlinear constraint objects.
\end{quote}
\end{quote}

\subsubsection{NlConstrArray.Size()}
\label{\detokenize{csapi/nlconstrarray:nlconstrarray-size}}\begin{quote}

\sphinxAtStartPar
Get the number of nonlinear constraint objects.

\sphinxAtStartPar
\sphinxstylestrong{Synopsis}
\begin{quote}

\sphinxAtStartPar
\sphinxcode{\sphinxupquote{int Size()}}
\end{quote}

\sphinxAtStartPar
\sphinxstylestrong{Return}
\begin{quote}

\sphinxAtStartPar
number of nonlinear constraint objects.
\end{quote}
\end{quote}

\subsection{NlConstrBuilder Class}
\label{\detokenize{csharpapiref:nlconstrbuilder-class}}\label{\detokenize{csharpapiref:chapcsharpapiref-nlconstrbuilder}}
\sphinxAtStartPar
COPT nonlinear constraint builder object. To help building a nonlinear constraint, given a
nonlinear expression, constraint sense and right\sphinxhyphen{}hand side value, Cardinal Optimizer provides
NlConstrBuilder class, which defines the following methods.

\sphinxstepscope

\subsubsection{NlConstrBuilder.NlConstrBuilder()}
\label{\detokenize{csapi/nlconstrbuilder:nlconstrbuilder-nlconstrbuilder}}\label{\detokenize{csapi/nlconstrbuilder::doc}}\begin{quote}

\sphinxAtStartPar
Constructor of NlConstrBuilder object.

\sphinxAtStartPar
\sphinxstylestrong{Synopsis}
\begin{quote}

\sphinxAtStartPar
\sphinxcode{\sphinxupquote{NlConstrBuilder()}}
\end{quote}
\end{quote}

\subsubsection{NlConstrBuilder.GetNlExpr()}
\label{\detokenize{csapi/nlconstrbuilder:nlconstrbuilder-getnlexpr}}\begin{quote}

\sphinxAtStartPar
Get nonlinear expression associated with constraint.

\sphinxAtStartPar
\sphinxstylestrong{Synopsis}
\begin{quote}

\sphinxAtStartPar
\sphinxcode{\sphinxupquote{NlExpr GetNlExpr()}}
\end{quote}

\sphinxAtStartPar
\sphinxstylestrong{Return}
\begin{quote}

\sphinxAtStartPar
nonlinear expression object.
\end{quote}
\end{quote}

\subsubsection{NlConstrBuilder.GetRange()}
\label{\detokenize{csapi/nlconstrbuilder:nlconstrbuilder-getrange}}\begin{quote}

\sphinxAtStartPar
Get range from lower bound to upper bound of range constraint.

\sphinxAtStartPar
\sphinxstylestrong{Synopsis}
\begin{quote}

\sphinxAtStartPar
\sphinxcode{\sphinxupquote{double GetRange()}}
\end{quote}

\sphinxAtStartPar
\sphinxstylestrong{Return}
\begin{quote}

\sphinxAtStartPar
length from lower bound to upper bound of nonlinear constraint.
\end{quote}
\end{quote}

\subsubsection{NlConstrBuilder.GetSense()}
\label{\detokenize{csapi/nlconstrbuilder:nlconstrbuilder-getsense}}\begin{quote}

\sphinxAtStartPar
Get sense associated with nonlinear constraint.

\sphinxAtStartPar
\sphinxstylestrong{Synopsis}
\begin{quote}

\sphinxAtStartPar
\sphinxcode{\sphinxupquote{char GetSense()}}
\end{quote}

\sphinxAtStartPar
\sphinxstylestrong{Return}
\begin{quote}

\sphinxAtStartPar
nonlinear constraint sense.
\end{quote}
\end{quote}

\subsubsection{NlConstrBuilder.Set()}
\label{\detokenize{csapi/nlconstrbuilder:nlconstrbuilder-set}}\begin{quote}

\sphinxAtStartPar
Set detail of a nonlinear constraint to its builder object.

\sphinxAtStartPar
\sphinxstylestrong{Synopsis}
\begin{quote}

\sphinxAtStartPar
\sphinxcode{\sphinxupquote{void Set(}}
\begin{quote}

\sphinxAtStartPar
\sphinxcode{\sphinxupquote{NlExpr expr,}}

\sphinxAtStartPar
\sphinxcode{\sphinxupquote{char sense,}}

\sphinxAtStartPar
\sphinxcode{\sphinxupquote{double rhs)}}
\end{quote}
\end{quote}

\sphinxAtStartPar
\sphinxstylestrong{Arguments}
\begin{quote}

\sphinxAtStartPar
\sphinxcode{\sphinxupquote{expr}}: nonlinear expression object at one side of nonlinear constraint

\sphinxAtStartPar
\sphinxcode{\sphinxupquote{sense}}: constraint sense other than COPT\_RANGE.

\sphinxAtStartPar
\sphinxcode{\sphinxupquote{rhs}}: constant of right side of nonlinear constraint.
\end{quote}
\end{quote}

\subsubsection{NlConstrBuilder.SetRange()}
\label{\detokenize{csapi/nlconstrbuilder:nlconstrbuilder-setrange}}\begin{quote}

\sphinxAtStartPar
Set a range constraint to nonlinear constraint builder.

\sphinxAtStartPar
\sphinxstylestrong{Synopsis}
\begin{quote}

\sphinxAtStartPar
\sphinxcode{\sphinxupquote{void SetRange(NlExpr expr, double range)}}
\end{quote}

\sphinxAtStartPar
\sphinxstylestrong{Arguments}
\begin{quote}

\sphinxAtStartPar
\sphinxcode{\sphinxupquote{expr}}: nonlinear expression object, whose constant is negative upper bound.

\sphinxAtStartPar
\sphinxcode{\sphinxupquote{range}}: length from lower bound to upper bound of nonlinear constraint. Must greater than 0.
\end{quote}
\end{quote}

\subsection{NlConstrBuilderArray Class}
\label{\detokenize{csharpapiref:nlconstrbuilderarray-class}}\label{\detokenize{csharpapiref:chapcsharpapiref-nlconstrbuilderarray}}
\sphinxAtStartPar
COPT nonlinear constraint builder array object. To store and access a set of
{\hyperref[\detokenize{csharpapiref:chapcsharpapiref-nlconstrbuilder}]{\sphinxcrossref{\DUrole{std,std-ref}{NlConstrBuilder Class}}}} objects, Cardinal Optimizer provides
\sphinxcode{\sphinxupquote{NlConstrBuilderArray}} class, which defines the following methods.

\sphinxstepscope

\subsubsection{NlConstrBuilderArray.NlConstrBuilderArray()}
\label{\detokenize{csapi/nlconstrbuilderarray:nlconstrbuilderarray-nlconstrbuilderarray}}\label{\detokenize{csapi/nlconstrbuilderarray::doc}}\begin{quote}

\sphinxAtStartPar
Constructor of NlConstrBuilderArray object.

\sphinxAtStartPar
\sphinxstylestrong{Synopsis}
\begin{quote}

\sphinxAtStartPar
\sphinxcode{\sphinxupquote{NlConstrBuilderArray()}}
\end{quote}
\end{quote}

\subsubsection{NlConstrBuilderArray.GetBuilder()}
\label{\detokenize{csapi/nlconstrbuilderarray:nlconstrbuilderarray-getbuilder}}\begin{quote}

\sphinxAtStartPar
Get idx\sphinxhyphen{}th nonlinear constraint builder object.

\sphinxAtStartPar
\sphinxstylestrong{Synopsis}
\begin{quote}

\sphinxAtStartPar
\sphinxcode{\sphinxupquote{NlConstrBuilder GetBuilder(int idx)}}
\end{quote}

\sphinxAtStartPar
\sphinxstylestrong{Arguments}
\begin{quote}

\sphinxAtStartPar
\sphinxcode{\sphinxupquote{idx}}: index of the nonlinear constraint builder.
\end{quote}

\sphinxAtStartPar
\sphinxstylestrong{Return}
\begin{quote}

\sphinxAtStartPar
nonlinear constraint builder object with index idx.
\end{quote}
\end{quote}

\subsubsection{NlConstrBuilderArray.PushBack()}
\label{\detokenize{csapi/nlconstrbuilderarray:nlconstrbuilderarray-pushback}}\begin{quote}

\sphinxAtStartPar
Add a nonlinear constraint builder object to nonlinear constraint builder array.

\sphinxAtStartPar
\sphinxstylestrong{Synopsis}
\begin{quote}

\sphinxAtStartPar
\sphinxcode{\sphinxupquote{void PushBack(NlConstrBuilder builder)}}
\end{quote}

\sphinxAtStartPar
\sphinxstylestrong{Arguments}
\begin{quote}

\sphinxAtStartPar
\sphinxcode{\sphinxupquote{builder}}: a nonlinear constraint builder object.
\end{quote}
\end{quote}

\subsubsection{NlConstrBuilderArray.Reserve()}
\label{\detokenize{csapi/nlconstrbuilderarray:nlconstrbuilderarray-reserve}}\begin{quote}

\sphinxAtStartPar
Reserve capacity to contain at least n items.

\sphinxAtStartPar
\sphinxstylestrong{Synopsis}
\begin{quote}

\sphinxAtStartPar
\sphinxcode{\sphinxupquote{void Reserve(int n)}}
\end{quote}

\sphinxAtStartPar
\sphinxstylestrong{Arguments}
\begin{quote}

\sphinxAtStartPar
\sphinxcode{\sphinxupquote{n}}: capacity of nonlinear constraint objects.
\end{quote}
\end{quote}

\subsubsection{NlConstrBuilderArray.Size()}
\label{\detokenize{csapi/nlconstrbuilderarray:nlconstrbuilderarray-size}}\begin{quote}

\sphinxAtStartPar
Get the number of nonlinear constraint builder objects.

\sphinxAtStartPar
\sphinxstylestrong{Synopsis}
\begin{quote}

\sphinxAtStartPar
\sphinxcode{\sphinxupquote{int Size()}}
\end{quote}

\sphinxAtStartPar
\sphinxstylestrong{Return}
\begin{quote}

\sphinxAtStartPar
number of nonlinear constraint builder objects.
\end{quote}
\end{quote}

\subsection{NL Namespace}
\label{\detokenize{csharpapiref:nl-namespace}}\label{\detokenize{csharpapiref:chapcsharpapiref-nl}}
\sphinxAtStartPar
Common nonlinear functions in the \sphinxcode{\sphinxupquote{NL}} namespace are provided for
constructing nonlinear expressions.  The following methods are provided:

\sphinxstepscope

\subsubsection{NL.Abs()}
\label{\detokenize{csapi/nl:nl-abs}}\label{\detokenize{csapi/nl::doc}}\begin{quote}

\sphinxAtStartPar
Calculate absolute value of a nonlinear expression.

\sphinxAtStartPar
\sphinxstylestrong{Synopsis}
\begin{quote}

\sphinxAtStartPar
\sphinxcode{\sphinxupquote{static NlExpr Abs(NlExpr expr)}}
\end{quote}

\sphinxAtStartPar
\sphinxstylestrong{Arguments}
\begin{quote}

\sphinxAtStartPar
\sphinxcode{\sphinxupquote{expr}}: a nonlinear expression.
\end{quote}

\sphinxAtStartPar
\sphinxstylestrong{Return}
\begin{quote}

\sphinxAtStartPar
result as a nonlinear expression.
\end{quote}
\end{quote}

\subsubsection{NL.ACos()}
\label{\detokenize{csapi/nl:nl-acos}}\begin{quote}

\sphinxAtStartPar
Calculate arccosine of a nonlinear expression.

\sphinxAtStartPar
\sphinxstylestrong{Synopsis}
\begin{quote}

\sphinxAtStartPar
\sphinxcode{\sphinxupquote{static NlExpr ACos(NlExpr expr)}}
\end{quote}

\sphinxAtStartPar
\sphinxstylestrong{Arguments}
\begin{quote}

\sphinxAtStartPar
\sphinxcode{\sphinxupquote{expr}}: a nonlinear expression.
\end{quote}

\sphinxAtStartPar
\sphinxstylestrong{Return}
\begin{quote}

\sphinxAtStartPar
result as a nonlinear expression.
\end{quote}
\end{quote}

\subsubsection{NL.ACosH()}
\label{\detokenize{csapi/nl:nl-acosh}}\begin{quote}

\sphinxAtStartPar
Calculate inverse hyperbolic cosine of a nonlinear expression.

\sphinxAtStartPar
\sphinxstylestrong{Synopsis}
\begin{quote}

\sphinxAtStartPar
\sphinxcode{\sphinxupquote{static NlExpr ACosH(NlExpr expr)}}
\end{quote}

\sphinxAtStartPar
\sphinxstylestrong{Arguments}
\begin{quote}

\sphinxAtStartPar
\sphinxcode{\sphinxupquote{expr}}: a nonlinear expression.
\end{quote}

\sphinxAtStartPar
\sphinxstylestrong{Return}
\begin{quote}

\sphinxAtStartPar
result as a nonlinear expression.
\end{quote}
\end{quote}

\subsubsection{NL.ASin()}
\label{\detokenize{csapi/nl:nl-asin}}\begin{quote}

\sphinxAtStartPar
Calculate arcsine of a nonlinear expression.

\sphinxAtStartPar
\sphinxstylestrong{Synopsis}
\begin{quote}

\sphinxAtStartPar
\sphinxcode{\sphinxupquote{static NlExpr ASin(NlExpr expr)}}
\end{quote}

\sphinxAtStartPar
\sphinxstylestrong{Arguments}
\begin{quote}

\sphinxAtStartPar
\sphinxcode{\sphinxupquote{expr}}: a nonlinear expression.
\end{quote}

\sphinxAtStartPar
\sphinxstylestrong{Return}
\begin{quote}

\sphinxAtStartPar
result as a nonlinear expression.
\end{quote}
\end{quote}

\subsubsection{NL.ASinH()}
\label{\detokenize{csapi/nl:nl-asinh}}\begin{quote}

\sphinxAtStartPar
Calculate inverse hyperbolic sine of a nonlinear expression.

\sphinxAtStartPar
\sphinxstylestrong{Synopsis}
\begin{quote}

\sphinxAtStartPar
\sphinxcode{\sphinxupquote{static NlExpr ASinH(NlExpr expr)}}
\end{quote}

\sphinxAtStartPar
\sphinxstylestrong{Arguments}
\begin{quote}

\sphinxAtStartPar
\sphinxcode{\sphinxupquote{expr}}: a nonlinear expression.
\end{quote}

\sphinxAtStartPar
\sphinxstylestrong{Return}
\begin{quote}

\sphinxAtStartPar
result as a nonlinear expression.
\end{quote}
\end{quote}

\subsubsection{NL.ATan()}
\label{\detokenize{csapi/nl:nl-atan}}\begin{quote}

\sphinxAtStartPar
Calculate arctangent of a nonlinear expression.

\sphinxAtStartPar
\sphinxstylestrong{Synopsis}
\begin{quote}

\sphinxAtStartPar
\sphinxcode{\sphinxupquote{static NlExpr ATan(NlExpr expr)}}
\end{quote}

\sphinxAtStartPar
\sphinxstylestrong{Arguments}
\begin{quote}

\sphinxAtStartPar
\sphinxcode{\sphinxupquote{expr}}: a nonlinear expression.
\end{quote}

\sphinxAtStartPar
\sphinxstylestrong{Return}
\begin{quote}

\sphinxAtStartPar
result as a nonlinear expression.
\end{quote}
\end{quote}

\subsubsection{NL.ATan2()}
\label{\detokenize{csapi/nl:nl-atan2}}\begin{quote}

\sphinxAtStartPar
Calculate two\sphinxhyphen{}argument arctangent of a nonlinear expression.

\sphinxAtStartPar
\sphinxstylestrong{Synopsis}
\begin{quote}

\sphinxAtStartPar
\sphinxcode{\sphinxupquote{static NlExpr ATan2(NlExpr y, NlExpr x)}}
\end{quote}

\sphinxAtStartPar
\sphinxstylestrong{Arguments}
\begin{quote}

\sphinxAtStartPar
\sphinxcode{\sphinxupquote{y}}: y coordinate as a nonlinear expression.

\sphinxAtStartPar
\sphinxcode{\sphinxupquote{x}}: x coordinate as a nonlinear expression.
\end{quote}

\sphinxAtStartPar
\sphinxstylestrong{Return}
\begin{quote}

\sphinxAtStartPar
result as a nonlinear expression.
\end{quote}
\end{quote}

\subsubsection{NL.ATanH()}
\label{\detokenize{csapi/nl:nl-atanh}}\begin{quote}

\sphinxAtStartPar
Calculate inverse hyperbolic tangent of a nonlinear expression.

\sphinxAtStartPar
\sphinxstylestrong{Synopsis}
\begin{quote}

\sphinxAtStartPar
\sphinxcode{\sphinxupquote{static NlExpr ATanH(NlExpr expr)}}
\end{quote}

\sphinxAtStartPar
\sphinxstylestrong{Arguments}
\begin{quote}

\sphinxAtStartPar
\sphinxcode{\sphinxupquote{expr}}: a nonlinear expression.
\end{quote}

\sphinxAtStartPar
\sphinxstylestrong{Return}
\begin{quote}

\sphinxAtStartPar
result as a nonlinear expression.
\end{quote}
\end{quote}

\subsubsection{NL.Ceil()}
\label{\detokenize{csapi/nl:nl-ceil}}\begin{quote}

\sphinxAtStartPar
Calculate ceiling value of a nonlinear expression.

\sphinxAtStartPar
\sphinxstylestrong{Synopsis}
\begin{quote}

\sphinxAtStartPar
\sphinxcode{\sphinxupquote{static NlExpr Ceil(NlExpr expr)}}
\end{quote}

\sphinxAtStartPar
\sphinxstylestrong{Arguments}
\begin{quote}

\sphinxAtStartPar
\sphinxcode{\sphinxupquote{expr}}: a nonlinear expression.
\end{quote}

\sphinxAtStartPar
\sphinxstylestrong{Return}
\begin{quote}

\sphinxAtStartPar
result as a nonlinear expression.
\end{quote}
\end{quote}

\subsubsection{NL.Cos()}
\label{\detokenize{csapi/nl:nl-cos}}\begin{quote}

\sphinxAtStartPar
Calculate cosine of a nonlinear expression.

\sphinxAtStartPar
\sphinxstylestrong{Synopsis}
\begin{quote}

\sphinxAtStartPar
\sphinxcode{\sphinxupquote{static NlExpr Cos(NlExpr expr)}}
\end{quote}

\sphinxAtStartPar
\sphinxstylestrong{Arguments}
\begin{quote}

\sphinxAtStartPar
\sphinxcode{\sphinxupquote{expr}}: a nonlinear expression.
\end{quote}

\sphinxAtStartPar
\sphinxstylestrong{Return}
\begin{quote}

\sphinxAtStartPar
result as a nonlinear expression.
\end{quote}
\end{quote}

\subsubsection{NL.CosH()}
\label{\detokenize{csapi/nl:nl-cosh}}\begin{quote}

\sphinxAtStartPar
Calculate hyperbolic cosine of a nonlinear expression.

\sphinxAtStartPar
\sphinxstylestrong{Synopsis}
\begin{quote}

\sphinxAtStartPar
\sphinxcode{\sphinxupquote{static NlExpr CosH(NlExpr expr)}}
\end{quote}

\sphinxAtStartPar
\sphinxstylestrong{Arguments}
\begin{quote}

\sphinxAtStartPar
\sphinxcode{\sphinxupquote{expr}}: a nonlinear expression.
\end{quote}

\sphinxAtStartPar
\sphinxstylestrong{Return}
\begin{quote}

\sphinxAtStartPar
result as a nonlinear expression.
\end{quote}
\end{quote}

\subsubsection{NL.Exp()}
\label{\detokenize{csapi/nl:nl-exp}}\begin{quote}

\sphinxAtStartPar
Calculate exponential function of a nonlinear expression.

\sphinxAtStartPar
\sphinxstylestrong{Synopsis}
\begin{quote}

\sphinxAtStartPar
\sphinxcode{\sphinxupquote{static NlExpr Exp(NlExpr expo)}}
\end{quote}

\sphinxAtStartPar
\sphinxstylestrong{Arguments}
\begin{quote}

\sphinxAtStartPar
\sphinxcode{\sphinxupquote{expo}}: exponent as a nonlinear expression.
\end{quote}

\sphinxAtStartPar
\sphinxstylestrong{Return}
\begin{quote}

\sphinxAtStartPar
result as a nonlinear expression.
\end{quote}
\end{quote}

\subsubsection{NL.Floor()}
\label{\detokenize{csapi/nl:nl-floor}}\begin{quote}

\sphinxAtStartPar
Calculate floor value of a nonlinear expression.

\sphinxAtStartPar
\sphinxstylestrong{Synopsis}
\begin{quote}

\sphinxAtStartPar
\sphinxcode{\sphinxupquote{static NlExpr Floor(NlExpr expr)}}
\end{quote}

\sphinxAtStartPar
\sphinxstylestrong{Arguments}
\begin{quote}

\sphinxAtStartPar
\sphinxcode{\sphinxupquote{expr}}: a nonlinear expression.
\end{quote}

\sphinxAtStartPar
\sphinxstylestrong{Return}
\begin{quote}

\sphinxAtStartPar
result as a nonlinear expression.
\end{quote}
\end{quote}

\subsubsection{NL.Log10()}
\label{\detokenize{csapi/nl:nl-log10}}\begin{quote}

\sphinxAtStartPar
Calculate logarithmic function of a nonlinear expression with base 10.

\sphinxAtStartPar
\sphinxstylestrong{Synopsis}
\begin{quote}

\sphinxAtStartPar
\sphinxcode{\sphinxupquote{static NlExpr Log10(NlExpr expr)}}
\end{quote}

\sphinxAtStartPar
\sphinxstylestrong{Arguments}
\begin{quote}

\sphinxAtStartPar
\sphinxcode{\sphinxupquote{expr}}: a nonlinear expression.
\end{quote}

\sphinxAtStartPar
\sphinxstylestrong{Return}
\begin{quote}

\sphinxAtStartPar
result as a nonlinear expression.
\end{quote}
\end{quote}

\subsubsection{NL.Log()}
\label{\detokenize{csapi/nl:nl-log}}\begin{quote}

\sphinxAtStartPar
Calculate nature logarithmic function of a nonlinear expression.

\sphinxAtStartPar
\sphinxstylestrong{Synopsis}
\begin{quote}

\sphinxAtStartPar
\sphinxcode{\sphinxupquote{static NlExpr Log(NlExpr expr)}}
\end{quote}

\sphinxAtStartPar
\sphinxstylestrong{Arguments}
\begin{quote}

\sphinxAtStartPar
\sphinxcode{\sphinxupquote{expr}}: a nonlinear expression.
\end{quote}

\sphinxAtStartPar
\sphinxstylestrong{Return}
\begin{quote}

\sphinxAtStartPar
result as a nonlinear expression.
\end{quote}
\end{quote}

\subsubsection{NL.Neg()}
\label{\detokenize{csapi/nl:nl-neg}}\begin{quote}

\sphinxAtStartPar
Calculate negative value of a nonlinear expression.

\sphinxAtStartPar
\sphinxstylestrong{Synopsis}
\begin{quote}

\sphinxAtStartPar
\sphinxcode{\sphinxupquote{static NlExpr Neg(NlExpr expr)}}
\end{quote}

\sphinxAtStartPar
\sphinxstylestrong{Arguments}
\begin{quote}

\sphinxAtStartPar
\sphinxcode{\sphinxupquote{expr}}: a nonlinear expression.
\end{quote}

\sphinxAtStartPar
\sphinxstylestrong{Return}
\begin{quote}

\sphinxAtStartPar
result as a nonlinear expression.
\end{quote}
\end{quote}

\subsubsection{NL.Pow()}
\label{\detokenize{csapi/nl:nl-pow}}\begin{quote}

\sphinxAtStartPar
Calculate power function of a nonlinear expression.

\sphinxAtStartPar
\sphinxstylestrong{Synopsis}
\begin{quote}

\sphinxAtStartPar
\sphinxcode{\sphinxupquote{static NlExpr Pow(NlExpr bas, NlExpr expo)}}
\end{quote}

\sphinxAtStartPar
\sphinxstylestrong{Arguments}
\begin{quote}

\sphinxAtStartPar
\sphinxcode{\sphinxupquote{bas}}: base as a nonlinear expression.

\sphinxAtStartPar
\sphinxcode{\sphinxupquote{expo}}: exponent as a nonlinear expression.
\end{quote}

\sphinxAtStartPar
\sphinxstylestrong{Return}
\begin{quote}

\sphinxAtStartPar
result as a nonlinear expression.
\end{quote}
\end{quote}

\subsubsection{NL.Sin()}
\label{\detokenize{csapi/nl:nl-sin}}\begin{quote}

\sphinxAtStartPar
Calculate sine of a nonlinear expression.

\sphinxAtStartPar
\sphinxstylestrong{Synopsis}
\begin{quote}

\sphinxAtStartPar
\sphinxcode{\sphinxupquote{static NlExpr Sin(NlExpr expr)}}
\end{quote}

\sphinxAtStartPar
\sphinxstylestrong{Arguments}
\begin{quote}

\sphinxAtStartPar
\sphinxcode{\sphinxupquote{expr}}: a nonlinear expression.
\end{quote}

\sphinxAtStartPar
\sphinxstylestrong{Return}
\begin{quote}

\sphinxAtStartPar
result as a nonlinear expression.
\end{quote}
\end{quote}

\subsubsection{NL.SinH()}
\label{\detokenize{csapi/nl:nl-sinh}}\begin{quote}

\sphinxAtStartPar
Calculate hyperbolic sine of a nonlinear expression.

\sphinxAtStartPar
\sphinxstylestrong{Synopsis}
\begin{quote}

\sphinxAtStartPar
\sphinxcode{\sphinxupquote{static NlExpr SinH(NlExpr expr)}}
\end{quote}

\sphinxAtStartPar
\sphinxstylestrong{Arguments}
\begin{quote}

\sphinxAtStartPar
\sphinxcode{\sphinxupquote{expr}}: a nonlinear expression.
\end{quote}

\sphinxAtStartPar
\sphinxstylestrong{Return}
\begin{quote}

\sphinxAtStartPar
result as a nonlinear expression.
\end{quote}
\end{quote}

\subsubsection{NL.Sqrt()}
\label{\detokenize{csapi/nl:nl-sqrt}}\begin{quote}

\sphinxAtStartPar
Calculate square root of a nonlinear expression.

\sphinxAtStartPar
\sphinxstylestrong{Synopsis}
\begin{quote}

\sphinxAtStartPar
\sphinxcode{\sphinxupquote{static NlExpr Sqrt(NlExpr expr)}}
\end{quote}

\sphinxAtStartPar
\sphinxstylestrong{Arguments}
\begin{quote}

\sphinxAtStartPar
\sphinxcode{\sphinxupquote{expr}}: a nonlinear expression.
\end{quote}

\sphinxAtStartPar
\sphinxstylestrong{Return}
\begin{quote}

\sphinxAtStartPar
result as a nonlinear expression.
\end{quote}
\end{quote}

\subsubsection{NL.Sum()}
\label{\detokenize{csapi/nl:nl-sum}}\begin{quote}

\sphinxAtStartPar
Sum of nonlinear expressions.

\sphinxAtStartPar
\sphinxstylestrong{Synopsis}
\begin{quote}

\sphinxAtStartPar
\sphinxcode{\sphinxupquote{static NlExpr Sum(}}
\begin{quote}

\sphinxAtStartPar
\sphinxcode{\sphinxupquote{NlExpr op1,}}

\sphinxAtStartPar
\sphinxcode{\sphinxupquote{NlExpr op2,}}

\sphinxAtStartPar
\sphinxcode{\sphinxupquote{NlExpr op3,}}

\sphinxAtStartPar
\sphinxcode{\sphinxupquote{NlExpr op4)}}
\end{quote}
\end{quote}

\sphinxAtStartPar
\sphinxstylestrong{Arguments}
\begin{quote}

\sphinxAtStartPar
\sphinxcode{\sphinxupquote{op1}}: first nonlinear expression.

\sphinxAtStartPar
\sphinxcode{\sphinxupquote{op2}}: second nonlinear expression.

\sphinxAtStartPar
\sphinxcode{\sphinxupquote{op3}}: third nonlinear expression.

\sphinxAtStartPar
\sphinxcode{\sphinxupquote{op4}}: fourth nonlinear expression.
\end{quote}

\sphinxAtStartPar
\sphinxstylestrong{Return}
\begin{quote}

\sphinxAtStartPar
result as a nonlinear expression.
\end{quote}
\end{quote}

\subsubsection{NL.Sum()}
\label{\detokenize{csapi/nl:id1}}\begin{quote}

\sphinxAtStartPar
Sum of nonlinear expressions.

\sphinxAtStartPar
\sphinxstylestrong{Synopsis}
\begin{quote}

\sphinxAtStartPar
\sphinxcode{\sphinxupquote{static NlExpr Sum(}}
\begin{quote}

\sphinxAtStartPar
\sphinxcode{\sphinxupquote{NlExpr op1,}}

\sphinxAtStartPar
\sphinxcode{\sphinxupquote{NlExpr op2,}}

\sphinxAtStartPar
\sphinxcode{\sphinxupquote{NlExpr op3)}}
\end{quote}
\end{quote}

\sphinxAtStartPar
\sphinxstylestrong{Arguments}
\begin{quote}

\sphinxAtStartPar
\sphinxcode{\sphinxupquote{op1}}: first nonlinear expression.

\sphinxAtStartPar
\sphinxcode{\sphinxupquote{op2}}: second nonlinear expression.

\sphinxAtStartPar
\sphinxcode{\sphinxupquote{op3}}: third nonlinear expression.
\end{quote}

\sphinxAtStartPar
\sphinxstylestrong{Return}
\begin{quote}

\sphinxAtStartPar
result as a nonlinear expression.
\end{quote}
\end{quote}

\subsubsection{NL.Sum()}
\label{\detokenize{csapi/nl:id2}}\begin{quote}

\sphinxAtStartPar
Sum of nonlinear expressions.

\sphinxAtStartPar
\sphinxstylestrong{Synopsis}
\begin{quote}

\sphinxAtStartPar
\sphinxcode{\sphinxupquote{static NlExpr Sum(NlExpr{[}{]} exprs)}}
\end{quote}

\sphinxAtStartPar
\sphinxstylestrong{Arguments}
\begin{quote}

\sphinxAtStartPar
\sphinxcode{\sphinxupquote{exprs}}: array of nonlinear expressions.
\end{quote}

\sphinxAtStartPar
\sphinxstylestrong{Return}
\begin{quote}

\sphinxAtStartPar
result as a nonlinear expression.
\end{quote}
\end{quote}

\subsubsection{NL.Tan()}
\label{\detokenize{csapi/nl:nl-tan}}\begin{quote}

\sphinxAtStartPar
Calculate tangent of a nonlinear expression.

\sphinxAtStartPar
\sphinxstylestrong{Synopsis}
\begin{quote}

\sphinxAtStartPar
\sphinxcode{\sphinxupquote{static NlExpr Tan(NlExpr expr)}}
\end{quote}

\sphinxAtStartPar
\sphinxstylestrong{Arguments}
\begin{quote}

\sphinxAtStartPar
\sphinxcode{\sphinxupquote{expr}}: a nonlinear expression.
\end{quote}

\sphinxAtStartPar
\sphinxstylestrong{Return}
\begin{quote}

\sphinxAtStartPar
result as a nonlinear expression.
\end{quote}
\end{quote}

\subsubsection{NL.TanH()}
\label{\detokenize{csapi/nl:nl-tanh}}\begin{quote}

\sphinxAtStartPar
Calculat hyperbolic tangent of a nonlinear expression.

\sphinxAtStartPar
\sphinxstylestrong{Synopsis}
\begin{quote}

\sphinxAtStartPar
\sphinxcode{\sphinxupquote{static NlExpr TanH(NlExpr expr)}}
\end{quote}

\sphinxAtStartPar
\sphinxstylestrong{Arguments}
\begin{quote}

\sphinxAtStartPar
\sphinxcode{\sphinxupquote{expr}}: a nonlinear expression.
\end{quote}

\sphinxAtStartPar
\sphinxstylestrong{Return}
\begin{quote}

\sphinxAtStartPar
result as a nonlinear expression.
\end{quote}
\end{quote}

\subsection{CallbackBase}
\label{\detokenize{csharpapiref:callbackbase}}\label{\detokenize{csharpapiref:chapcsharpapiref-callback}}
\sphinxAtStartPar
COPT Callback abstract base object. Users must implment its virtual method
\sphinxcode{\sphinxupquote{virtual void CallbackBase::callback()}} to instantiate an instance, which
pass to \sphinxcode{\sphinxupquote{Model::SetCallback(CallbackBase cb, int cbctx)}} as the first
parameter. Subclass of CallbackBase inherits the following member methods:

\sphinxstepscope

\subsubsection{CallbackBase.CallbackBase()}
\label{\detokenize{csapi/callbackbase:callbackbase-callbackbase}}\label{\detokenize{csapi/callbackbase::doc}}\begin{quote}

\sphinxAtStartPar
Constructor of CallbackBase, implementing ICallback interface.

\sphinxAtStartPar
\sphinxstylestrong{Synopsis}
\begin{quote}

\sphinxAtStartPar
\sphinxcode{\sphinxupquote{CallbackBase()}}
\end{quote}
\end{quote}

\subsubsection{CallbackBase.AddLazyConstr()}
\label{\detokenize{csapi/callbackbase:callbackbase-addlazyconstr}}\begin{quote}

\sphinxAtStartPar
Add a lazy constraint to model.

\sphinxAtStartPar
\sphinxstylestrong{Synopsis}
\begin{quote}

\sphinxAtStartPar
\sphinxcode{\sphinxupquote{void AddLazyConstr(}}
\begin{quote}

\sphinxAtStartPar
\sphinxcode{\sphinxupquote{Expr lhs,}}

\sphinxAtStartPar
\sphinxcode{\sphinxupquote{char sense,}}

\sphinxAtStartPar
\sphinxcode{\sphinxupquote{double rhs)}}
\end{quote}
\end{quote}

\sphinxAtStartPar
\sphinxstylestrong{Arguments}
\begin{quote}

\sphinxAtStartPar
\sphinxcode{\sphinxupquote{lhs}}: expression for lazy contraint.

\sphinxAtStartPar
\sphinxcode{\sphinxupquote{sense}}: sense for lazy constraint.

\sphinxAtStartPar
\sphinxcode{\sphinxupquote{rhs}}: right hand side value for lazy constraint.
\end{quote}
\end{quote}

\subsubsection{CallbackBase.AddLazyConstr()}
\label{\detokenize{csapi/callbackbase:id1}}\begin{quote}

\sphinxAtStartPar
Add a lazy constraint to model.

\sphinxAtStartPar
\sphinxstylestrong{Synopsis}
\begin{quote}

\sphinxAtStartPar
\sphinxcode{\sphinxupquote{void AddLazyConstr(}}
\begin{quote}

\sphinxAtStartPar
\sphinxcode{\sphinxupquote{Expr lhs,}}

\sphinxAtStartPar
\sphinxcode{\sphinxupquote{char sense,}}

\sphinxAtStartPar
\sphinxcode{\sphinxupquote{Expr rhs)}}
\end{quote}
\end{quote}

\sphinxAtStartPar
\sphinxstylestrong{Arguments}
\begin{quote}

\sphinxAtStartPar
\sphinxcode{\sphinxupquote{lhs}}: left hand side expression for lazy contraint.

\sphinxAtStartPar
\sphinxcode{\sphinxupquote{sense}}: sense for lazy constraint.

\sphinxAtStartPar
\sphinxcode{\sphinxupquote{rhs}}: right hand side expression for lazy contraint.
\end{quote}
\end{quote}

\subsubsection{CallbackBase.AddLazyConstr()}
\label{\detokenize{csapi/callbackbase:id2}}\begin{quote}

\sphinxAtStartPar
Add a lazy constraint to model.

\sphinxAtStartPar
\sphinxstylestrong{Synopsis}
\begin{quote}

\sphinxAtStartPar
\sphinxcode{\sphinxupquote{void AddLazyConstr(ConstrBuilder builder)}}
\end{quote}

\sphinxAtStartPar
\sphinxstylestrong{Arguments}
\begin{quote}

\sphinxAtStartPar
\sphinxcode{\sphinxupquote{builder}}: builder for lazy contraint.
\end{quote}
\end{quote}

\subsubsection{CallbackBase.AddLazyConstrs()}
\label{\detokenize{csapi/callbackbase:callbackbase-addlazyconstrs}}\begin{quote}

\sphinxAtStartPar
Add lazy constraints to model.

\sphinxAtStartPar
\sphinxstylestrong{Synopsis}
\begin{quote}

\sphinxAtStartPar
\sphinxcode{\sphinxupquote{void AddLazyConstrs(ConstrBuilderArray builders)}}
\end{quote}

\sphinxAtStartPar
\sphinxstylestrong{Arguments}
\begin{quote}

\sphinxAtStartPar
\sphinxcode{\sphinxupquote{builders}}: array of builders for lazy contraints.
\end{quote}
\end{quote}

\subsubsection{CallbackBase.AddUserCut()}
\label{\detokenize{csapi/callbackbase:callbackbase-addusercut}}\begin{quote}

\sphinxAtStartPar
Add a user cut to model.

\sphinxAtStartPar
\sphinxstylestrong{Synopsis}
\begin{quote}

\sphinxAtStartPar
\sphinxcode{\sphinxupquote{void AddUserCut(}}
\begin{quote}

\sphinxAtStartPar
\sphinxcode{\sphinxupquote{Expr lhs,}}

\sphinxAtStartPar
\sphinxcode{\sphinxupquote{char sense,}}

\sphinxAtStartPar
\sphinxcode{\sphinxupquote{double rhs)}}
\end{quote}
\end{quote}

\sphinxAtStartPar
\sphinxstylestrong{Arguments}
\begin{quote}

\sphinxAtStartPar
\sphinxcode{\sphinxupquote{lhs}}: expression for user cut.

\sphinxAtStartPar
\sphinxcode{\sphinxupquote{sense}}: sense for user cut.

\sphinxAtStartPar
\sphinxcode{\sphinxupquote{rhs}}: right hand side value for user cut.
\end{quote}
\end{quote}

\subsubsection{CallbackBase.AddUserCut()}
\label{\detokenize{csapi/callbackbase:id3}}\begin{quote}

\sphinxAtStartPar
Add a user cut to model.

\sphinxAtStartPar
\sphinxstylestrong{Synopsis}
\begin{quote}

\sphinxAtStartPar
\sphinxcode{\sphinxupquote{void AddUserCut(}}
\begin{quote}

\sphinxAtStartPar
\sphinxcode{\sphinxupquote{Expr lhs,}}

\sphinxAtStartPar
\sphinxcode{\sphinxupquote{char sense,}}

\sphinxAtStartPar
\sphinxcode{\sphinxupquote{Expr rhs)}}
\end{quote}
\end{quote}

\sphinxAtStartPar
\sphinxstylestrong{Arguments}
\begin{quote}

\sphinxAtStartPar
\sphinxcode{\sphinxupquote{lhs}}: left hand side expression for user cut.

\sphinxAtStartPar
\sphinxcode{\sphinxupquote{sense}}: sense for user cut.

\sphinxAtStartPar
\sphinxcode{\sphinxupquote{rhs}}: right hand side expression for user cut.
\end{quote}
\end{quote}

\subsubsection{CallbackBase.AddUserCut()}
\label{\detokenize{csapi/callbackbase:id4}}\begin{quote}

\sphinxAtStartPar
Add a user cut to model.

\sphinxAtStartPar
\sphinxstylestrong{Synopsis}
\begin{quote}

\sphinxAtStartPar
\sphinxcode{\sphinxupquote{void AddUserCut(ConstrBuilder builder)}}
\end{quote}

\sphinxAtStartPar
\sphinxstylestrong{Arguments}
\begin{quote}

\sphinxAtStartPar
\sphinxcode{\sphinxupquote{builder}}: builder for user cut.
\end{quote}
\end{quote}

\subsubsection{CallbackBase.AddUserCuts()}
\label{\detokenize{csapi/callbackbase:callbackbase-addusercuts}}\begin{quote}

\sphinxAtStartPar
Add user cuts to model.

\sphinxAtStartPar
\sphinxstylestrong{Synopsis}
\begin{quote}

\sphinxAtStartPar
\sphinxcode{\sphinxupquote{void AddUserCuts(ConstrBuilderArray builders)}}
\end{quote}

\sphinxAtStartPar
\sphinxstylestrong{Arguments}
\begin{quote}

\sphinxAtStartPar
\sphinxcode{\sphinxupquote{builders}}: array of builders for user cuts.
\end{quote}
\end{quote}

\subsubsection{CallbackBase.callback()}
\label{\detokenize{csapi/callbackbase:callbackbase-callback}}\begin{quote}

\sphinxAtStartPar
Pure virtual function defined in ICallback interface. User must implement it.

\sphinxAtStartPar
\sphinxstylestrong{Synopsis}
\begin{quote}

\sphinxAtStartPar
\sphinxcode{\sphinxupquote{void callback()}}
\end{quote}
\end{quote}

\subsubsection{CallbackBase.GetDblInfo()}
\label{\detokenize{csapi/callbackbase:callbackbase-getdblinfo}}\begin{quote}

\sphinxAtStartPar
Get double value of given information name in callback.

\sphinxAtStartPar
\sphinxstylestrong{Synopsis}
\begin{quote}

\sphinxAtStartPar
\sphinxcode{\sphinxupquote{double GetDblInfo(string cbinfo)}}
\end{quote}

\sphinxAtStartPar
\sphinxstylestrong{Arguments}
\begin{quote}

\sphinxAtStartPar
\sphinxcode{\sphinxupquote{cbinfo}}: name of callback info.
\end{quote}

\sphinxAtStartPar
\sphinxstylestrong{Return}
\begin{quote}

\sphinxAtStartPar
value of desired information.
\end{quote}
\end{quote}

\subsubsection{CallbackBase.GetIncumbent()}
\label{\detokenize{csapi/callbackbase:callbackbase-getincumbent}}\begin{quote}

\sphinxAtStartPar
Get best feasible solution of given variable in callback.

\sphinxAtStartPar
\sphinxstylestrong{Synopsis}
\begin{quote}

\sphinxAtStartPar
\sphinxcode{\sphinxupquote{double GetIncumbent(Var var)}}
\end{quote}

\sphinxAtStartPar
\sphinxstylestrong{Arguments}
\begin{quote}

\sphinxAtStartPar
\sphinxcode{\sphinxupquote{var}}: given variable.
\end{quote}

\sphinxAtStartPar
\sphinxstylestrong{Return}
\begin{quote}

\sphinxAtStartPar
best feasible solution of given variable.
\end{quote}
\end{quote}

\subsubsection{CallbackBase.GetIncumbent()}
\label{\detokenize{csapi/callbackbase:id5}}\begin{quote}

\sphinxAtStartPar
Get best feasible solution of variables in callback.

\sphinxAtStartPar
\sphinxstylestrong{Synopsis}
\begin{quote}

\sphinxAtStartPar
\sphinxcode{\sphinxupquote{double{[}{]} GetIncumbent(VarArray vars)}}
\end{quote}

\sphinxAtStartPar
\sphinxstylestrong{Arguments}
\begin{quote}

\sphinxAtStartPar
\sphinxcode{\sphinxupquote{vars}}: an array of variables.
\end{quote}

\sphinxAtStartPar
\sphinxstylestrong{Return}
\begin{quote}

\sphinxAtStartPar
best feasible solution of desired variables.
\end{quote}
\end{quote}

\subsubsection{CallbackBase.GetIncumbent()}
\label{\detokenize{csapi/callbackbase:id6}}\begin{quote}

\sphinxAtStartPar
Get best feasible solution of variables in callback.

\sphinxAtStartPar
\sphinxstylestrong{Synopsis}
\begin{quote}

\sphinxAtStartPar
\sphinxcode{\sphinxupquote{double{[}{]} GetIncumbent(Var{[}{]} vars)}}
\end{quote}

\sphinxAtStartPar
\sphinxstylestrong{Arguments}
\begin{quote}

\sphinxAtStartPar
\sphinxcode{\sphinxupquote{vars}}: an array of variables.
\end{quote}

\sphinxAtStartPar
\sphinxstylestrong{Return}
\begin{quote}

\sphinxAtStartPar
best feasible solution of desired variables.
\end{quote}
\end{quote}

\subsubsection{CallbackBase.GetIncumbent()}
\label{\detokenize{csapi/callbackbase:id7}}\begin{quote}

\sphinxAtStartPar
Get best feasible solution of all variables in callback.

\sphinxAtStartPar
\sphinxstylestrong{Synopsis}
\begin{quote}

\sphinxAtStartPar
\sphinxcode{\sphinxupquote{double{[}{]} GetIncumbent()}}
\end{quote}

\sphinxAtStartPar
\sphinxstylestrong{Return}
\begin{quote}

\sphinxAtStartPar
best feasible solution of all variables.
\end{quote}
\end{quote}

\subsubsection{CallbackBase.GetIntInfo()}
\label{\detokenize{csapi/callbackbase:callbackbase-getintinfo}}\begin{quote}

\sphinxAtStartPar
Get integer value of given information name in callback.

\sphinxAtStartPar
\sphinxstylestrong{Synopsis}
\begin{quote}

\sphinxAtStartPar
\sphinxcode{\sphinxupquote{int GetIntInfo(string cbinfo)}}
\end{quote}

\sphinxAtStartPar
\sphinxstylestrong{Arguments}
\begin{quote}

\sphinxAtStartPar
\sphinxcode{\sphinxupquote{cbinfo}}: name of callback info.
\end{quote}

\sphinxAtStartPar
\sphinxstylestrong{Return}
\begin{quote}

\sphinxAtStartPar
value of desired information.
\end{quote}
\end{quote}

\subsubsection{CallbackBase.GetRelaxSol()}
\label{\detokenize{csapi/callbackbase:callbackbase-getrelaxsol}}\begin{quote}

\sphinxAtStartPar
Get LP\sphinxhyphen{}relaxation solution of given variable in callback.

\sphinxAtStartPar
\sphinxstylestrong{Synopsis}
\begin{quote}

\sphinxAtStartPar
\sphinxcode{\sphinxupquote{double GetRelaxSol(Var var)}}
\end{quote}

\sphinxAtStartPar
\sphinxstylestrong{Arguments}
\begin{quote}

\sphinxAtStartPar
\sphinxcode{\sphinxupquote{var}}: given variable.
\end{quote}

\sphinxAtStartPar
\sphinxstylestrong{Return}
\begin{quote}

\sphinxAtStartPar
LP\sphinxhyphen{}relaxation solution of given variable.
\end{quote}
\end{quote}

\subsubsection{CallbackBase.GetRelaxSol()}
\label{\detokenize{csapi/callbackbase:id8}}\begin{quote}

\sphinxAtStartPar
Get LP\sphinxhyphen{}relaxation solution of variables in callback.

\sphinxAtStartPar
\sphinxstylestrong{Synopsis}
\begin{quote}

\sphinxAtStartPar
\sphinxcode{\sphinxupquote{double{[}{]} GetRelaxSol(VarArray vars)}}
\end{quote}

\sphinxAtStartPar
\sphinxstylestrong{Arguments}
\begin{quote}

\sphinxAtStartPar
\sphinxcode{\sphinxupquote{vars}}: an array of variables.
\end{quote}

\sphinxAtStartPar
\sphinxstylestrong{Return}
\begin{quote}

\sphinxAtStartPar
LP\sphinxhyphen{}relaxation solution of variables.
\end{quote}
\end{quote}

\subsubsection{CallbackBase.GetRelaxSol()}
\label{\detokenize{csapi/callbackbase:id9}}\begin{quote}

\sphinxAtStartPar
Get LP\sphinxhyphen{}relaxation solution of variables in callback.

\sphinxAtStartPar
\sphinxstylestrong{Synopsis}
\begin{quote}

\sphinxAtStartPar
\sphinxcode{\sphinxupquote{double{[}{]} GetRelaxSol(Var{[}{]} vars)}}
\end{quote}

\sphinxAtStartPar
\sphinxstylestrong{Arguments}
\begin{quote}

\sphinxAtStartPar
\sphinxcode{\sphinxupquote{vars}}: an array of variables.
\end{quote}

\sphinxAtStartPar
\sphinxstylestrong{Return}
\begin{quote}

\sphinxAtStartPar
LP\sphinxhyphen{}relaxation solution of variables.
\end{quote}
\end{quote}

\subsubsection{CallbackBase.GetRelaxSol()}
\label{\detokenize{csapi/callbackbase:id10}}\begin{quote}

\sphinxAtStartPar
Get LP\sphinxhyphen{}relaxation solution of all variables in callback.

\sphinxAtStartPar
\sphinxstylestrong{Synopsis}
\begin{quote}

\sphinxAtStartPar
\sphinxcode{\sphinxupquote{double{[}{]} GetRelaxSol()}}
\end{quote}

\sphinxAtStartPar
\sphinxstylestrong{Return}
\begin{quote}

\sphinxAtStartPar
LP\sphinxhyphen{}relaxation solution of all variables.
\end{quote}
\end{quote}

\subsubsection{CallbackBase.GetSolution()}
\label{\detokenize{csapi/callbackbase:callbackbase-getsolution}}\begin{quote}

\sphinxAtStartPar
Get solution of given variable in callback.

\sphinxAtStartPar
\sphinxstylestrong{Synopsis}
\begin{quote}

\sphinxAtStartPar
\sphinxcode{\sphinxupquote{double GetSolution(Var var)}}
\end{quote}

\sphinxAtStartPar
\sphinxstylestrong{Arguments}
\begin{quote}

\sphinxAtStartPar
\sphinxcode{\sphinxupquote{var}}: given variable.
\end{quote}

\sphinxAtStartPar
\sphinxstylestrong{Return}
\begin{quote}

\sphinxAtStartPar
solution of given variable.
\end{quote}
\end{quote}

\subsubsection{CallbackBase.GetSolution()}
\label{\detokenize{csapi/callbackbase:id11}}\begin{quote}

\sphinxAtStartPar
Get solution of variables in callback.

\sphinxAtStartPar
\sphinxstylestrong{Synopsis}
\begin{quote}

\sphinxAtStartPar
\sphinxcode{\sphinxupquote{double{[}{]} GetSolution(VarArray vars)}}
\end{quote}

\sphinxAtStartPar
\sphinxstylestrong{Arguments}
\begin{quote}

\sphinxAtStartPar
\sphinxcode{\sphinxupquote{vars}}: an array of variables.
\end{quote}

\sphinxAtStartPar
\sphinxstylestrong{Return}
\begin{quote}

\sphinxAtStartPar
solution of variables.
\end{quote}
\end{quote}

\subsubsection{CallbackBase.GetSolution()}
\label{\detokenize{csapi/callbackbase:id12}}\begin{quote}

\sphinxAtStartPar
Get solution of variables in callback.

\sphinxAtStartPar
\sphinxstylestrong{Synopsis}
\begin{quote}

\sphinxAtStartPar
\sphinxcode{\sphinxupquote{double{[}{]} GetSolution(Var{[}{]} vars)}}
\end{quote}

\sphinxAtStartPar
\sphinxstylestrong{Arguments}
\begin{quote}

\sphinxAtStartPar
\sphinxcode{\sphinxupquote{vars}}: an array of variables.
\end{quote}

\sphinxAtStartPar
\sphinxstylestrong{Return}
\begin{quote}

\sphinxAtStartPar
solution of variables.
\end{quote}
\end{quote}

\subsubsection{CallbackBase.GetSolution()}
\label{\detokenize{csapi/callbackbase:id13}}\begin{quote}

\sphinxAtStartPar
Get solution of all variables in callback.

\sphinxAtStartPar
\sphinxstylestrong{Synopsis}
\begin{quote}

\sphinxAtStartPar
\sphinxcode{\sphinxupquote{double{[}{]} GetSolution()}}
\end{quote}

\sphinxAtStartPar
\sphinxstylestrong{Return}
\begin{quote}

\sphinxAtStartPar
solution of all variables.
\end{quote}
\end{quote}

\subsubsection{CallbackBase.Interrupt()}
\label{\detokenize{csapi/callbackbase:callbackbase-interrupt}}\begin{quote}

\sphinxAtStartPar
Interrupt solving problems in callback

\sphinxAtStartPar
\sphinxstylestrong{Synopsis}
\begin{quote}

\sphinxAtStartPar
\sphinxcode{\sphinxupquote{void Interrupt()}}
\end{quote}
\end{quote}

\subsubsection{CallbackBase.LoadSolution()}
\label{\detokenize{csapi/callbackbase:callbackbase-loadsolution}}\begin{quote}

\sphinxAtStartPar
Load customized solution to model.

\sphinxAtStartPar
\sphinxstylestrong{Synopsis}
\begin{quote}

\sphinxAtStartPar
\sphinxcode{\sphinxupquote{double LoadSolution()}}
\end{quote}

\sphinxAtStartPar
\sphinxstylestrong{Return}
\begin{quote}

\sphinxAtStartPar
objective value of given solution.
\end{quote}
\end{quote}

\subsubsection{CallbackBase.SetSolution()}
\label{\detokenize{csapi/callbackbase:callbackbase-setsolution}}\begin{quote}

\sphinxAtStartPar
Set solution of a given variable in callback.

\sphinxAtStartPar
\sphinxstylestrong{Synopsis}
\begin{quote}

\sphinxAtStartPar
\sphinxcode{\sphinxupquote{void SetSolution(Var var, double val)}}
\end{quote}

\sphinxAtStartPar
\sphinxstylestrong{Arguments}
\begin{quote}

\sphinxAtStartPar
\sphinxcode{\sphinxupquote{var}}: a variable object.

\sphinxAtStartPar
\sphinxcode{\sphinxupquote{val}}: double value.
\end{quote}
\end{quote}

\subsubsection{CallbackBase.SetSolution()}
\label{\detokenize{csapi/callbackbase:id14}}\begin{quote}

\sphinxAtStartPar
Set solution of variables in callback.

\sphinxAtStartPar
\sphinxstylestrong{Synopsis}
\begin{quote}

\sphinxAtStartPar
\sphinxcode{\sphinxupquote{void SetSolution(VarArray vars, double{[}{]} vals)}}
\end{quote}

\sphinxAtStartPar
\sphinxstylestrong{Arguments}
\begin{quote}

\sphinxAtStartPar
\sphinxcode{\sphinxupquote{vars}}: an array of variable objects.

\sphinxAtStartPar
\sphinxcode{\sphinxupquote{vals}}: an array of double values.
\end{quote}
\end{quote}

\subsubsection{CallbackBase.SetSolution()}
\label{\detokenize{csapi/callbackbase:id15}}\begin{quote}

\sphinxAtStartPar
Set solution of variables in callback.

\sphinxAtStartPar
\sphinxstylestrong{Synopsis}
\begin{quote}

\sphinxAtStartPar
\sphinxcode{\sphinxupquote{void SetSolution(Var{[}{]} vars, double{[}{]} vals)}}
\end{quote}

\sphinxAtStartPar
\sphinxstylestrong{Arguments}
\begin{quote}

\sphinxAtStartPar
\sphinxcode{\sphinxupquote{vars}}: an array of variable objects.

\sphinxAtStartPar
\sphinxcode{\sphinxupquote{vals}}: an array of double values.
\end{quote}
\end{quote}

\subsubsection{CallbackBase.Where()}
\label{\detokenize{csapi/callbackbase:callbackbase-where}}\begin{quote}

\sphinxAtStartPar
Get context in callback.

\sphinxAtStartPar
\sphinxstylestrong{Synopsis}
\begin{quote}

\sphinxAtStartPar
\sphinxcode{\sphinxupquote{int Where()}}
\end{quote}

\sphinxAtStartPar
\sphinxstylestrong{Return}
\begin{quote}

\sphinxAtStartPar
integer value of context.
\end{quote}
\end{quote}

\subsection{ProbBuffer}
\label{\detokenize{csharpapiref:probbuffer}}\label{\detokenize{csharpapiref:chapcsharpapiref-probbuffer}}
\sphinxAtStartPar
Buffer object for COPT problem. ProbBuffer object holds the (MPS) problem in string format.

\sphinxstepscope

\subsubsection{ProbBuffer.ProbBuffer()}
\label{\detokenize{csapi/probbuffer:probbuffer-probbuffer}}\label{\detokenize{csapi/probbuffer::doc}}\begin{quote}

\sphinxAtStartPar
Constructor of ProbBuffer object.

\sphinxAtStartPar
\sphinxstylestrong{Synopsis}
\begin{quote}

\sphinxAtStartPar
\sphinxcode{\sphinxupquote{ProbBuffer(int sz)}}
\end{quote}

\sphinxAtStartPar
\sphinxstylestrong{Arguments}
\begin{quote}

\sphinxAtStartPar
\sphinxcode{\sphinxupquote{sz}}: initial size of the problem buffer.
\end{quote}
\end{quote}

\subsubsection{ProbBuffer.GetData()}
\label{\detokenize{csapi/probbuffer:probbuffer-getdata}}\begin{quote}

\sphinxAtStartPar
Get string of problem in problem buffer.

\sphinxAtStartPar
\sphinxstylestrong{Synopsis}
\begin{quote}

\sphinxAtStartPar
\sphinxcode{\sphinxupquote{string GetData()}}
\end{quote}

\sphinxAtStartPar
\sphinxstylestrong{Return}
\begin{quote}

\sphinxAtStartPar
string of problem in problem buffer.
\end{quote}
\end{quote}

\subsubsection{ProbBuffer.Resize()}
\label{\detokenize{csapi/probbuffer:probbuffer-resize}}\begin{quote}

\sphinxAtStartPar
Resize buffer to given size, and zero\sphinxhyphen{}ended

\sphinxAtStartPar
\sphinxstylestrong{Synopsis}
\begin{quote}

\sphinxAtStartPar
\sphinxcode{\sphinxupquote{void Resize(int sz)}}
\end{quote}

\sphinxAtStartPar
\sphinxstylestrong{Arguments}
\begin{quote}

\sphinxAtStartPar
\sphinxcode{\sphinxupquote{sz}}: new buffer size.
\end{quote}
\end{quote}

\subsubsection{ProbBuffer.Size()}
\label{\detokenize{csapi/probbuffer:probbuffer-size}}\begin{quote}

\sphinxAtStartPar
Get the size of problem buffer.

\sphinxAtStartPar
\sphinxstylestrong{Synopsis}
\begin{quote}

\sphinxAtStartPar
\sphinxcode{\sphinxupquote{int Size()}}
\end{quote}

\sphinxAtStartPar
\sphinxstylestrong{Return}
\begin{quote}

\sphinxAtStartPar
size of problem buffer.
\end{quote}
\end{quote}

\subsection{CoptException}
\label{\detokenize{csharpapiref:coptexception}}
\sphinxAtStartPar
Copt exception object.

\sphinxstepscope

\subsubsection{CoptException.CoptException()}
\label{\detokenize{csapi/coptexception:coptexception-coptexception}}\label{\detokenize{csapi/coptexception::doc}}\begin{quote}

\sphinxAtStartPar
Constructor of COPT Exception class.

\sphinxAtStartPar
\sphinxstylestrong{Synopsis}
\begin{quote}

\sphinxAtStartPar
\sphinxcode{\sphinxupquote{CoptException(int code, string msg)}}
\end{quote}

\sphinxAtStartPar
\sphinxstylestrong{Arguments}
\begin{quote}

\sphinxAtStartPar
\sphinxcode{\sphinxupquote{code}}: error code for exception.

\sphinxAtStartPar
\sphinxcode{\sphinxupquote{msg}}: error message for exception.
\end{quote}
\end{quote}

\subsubsection{CoptException.GetCode()}
\label{\detokenize{csapi/coptexception:coptexception-getcode}}\begin{quote}

\sphinxAtStartPar
Get the error code associated with the exception.

\sphinxAtStartPar
\sphinxstylestrong{Synopsis}
\begin{quote}

\sphinxAtStartPar
\sphinxcode{\sphinxupquote{int GetCode()}}
\end{quote}

\sphinxAtStartPar
\sphinxstylestrong{Return}
\begin{quote}

\sphinxAtStartPar
the error code.
\end{quote}
\end{quote}

\sphinxstepscope

\chapter{Java API Reference}
\label{\detokenize{javaapiref:java-api-reference}}\label{\detokenize{javaapiref:chapjavaapiref}}\label{\detokenize{javaapiref::doc}}
\sphinxAtStartPar
The \sphinxstylestrong{Cardinal Optimizer} provides a Java API library.
This chapter documents all COPT Java constants and API functions for Java applications.

\section{Constants}
\label{\detokenize{javaapiref:constants}}\label{\detokenize{javaapiref:chapjavaapiref-const}}
\sphinxAtStartPar
There are four types of constants defined in \sphinxstylestrong{Cardinal Optimizer}.
They are general constants, information constants, attributes constants and parameters constants.

\subsection{General Constants}
\label{\detokenize{javaapiref:general-constants}}\label{\detokenize{javaapiref:chapjavaapiref-const-general}}
\sphinxAtStartPar
For the contents of Java general constants, see {\hyperref[\detokenize{constant:chapconst}]{\sphinxcrossref{\DUrole{std,std-ref}{General Constants}}}}.

\sphinxAtStartPar
General constants are defined in \sphinxcode{\sphinxupquote{Consts}} class. User may refer general constants
with namespace, that is, \sphinxcode{\sphinxupquote{copt.Consts.XXXX}}.

\subsection{Attributes}
\label{\detokenize{javaapiref:attributes}}\label{\detokenize{javaapiref:chapjavaapiref-const-attrs}}
\sphinxAtStartPar
For the contents of Java attribute constants, see {\hyperref[\detokenize{attribute:chapattrs}]{\sphinxcrossref{\DUrole{std,std-ref}{Attributes}}}}.

\sphinxAtStartPar
All COPT Java attributes are defined in \sphinxcode{\sphinxupquote{DblAttr}} and \sphinxcode{\sphinxupquote{IntAttr}} classes.
User may refer double attributes by \sphinxcode{\sphinxupquote{copt.DblAttr.XXXX}}, and integer
attributes by \sphinxcode{\sphinxupquote{copt.IntAttr.XXXX}}.

\sphinxAtStartPar
In the Java API, user can get the attribute value by specifying the attribute name. The two functions of obtaining attribute values are as follows, please refer to {\hyperref[\detokenize{javaapiref:chapjavaapiref-model}]{\sphinxcrossref{\DUrole{std,std-ref}{Java API: Model Class}}}} for details.
\begin{itemize}
\item {} 
\sphinxAtStartPar
\sphinxcode{\sphinxupquote{Model.getIntAttr()}}: Get value of a COPT integer attribute.

\item {} 
\sphinxAtStartPar
\sphinxcode{\sphinxupquote{Model.getDblAttr()}}: Get value of a COPT double attribute.

\end{itemize}

\subsection{Information}
\label{\detokenize{javaapiref:information}}\label{\detokenize{javaapiref:chapjavaapiref-const-info}}
\sphinxAtStartPar
For the content of Java information constants, see {\hyperref[\detokenize{information:chapinfo}]{\sphinxcrossref{\DUrole{std,std-ref}{Information}}}}.

\sphinxAtStartPar
In the Java API, information constants are defined in the \sphinxcode{\sphinxupquote{DblInfo}} class.
Users can access information constants through the prefix \sphinxcode{\sphinxupquote{copt}} in the namespace (usually can be omitted) \sphinxcode{\sphinxupquote{copt.DblInfo.}}

\sphinxAtStartPar
For instance, \sphinxcode{\sphinxupquote{copt.DblInfo.Obj}} is the coefficients of variables in the objective function.

\subsection{Callback Information}
\label{\detokenize{javaapiref:callback-information}}
\sphinxAtStartPar
For the content of Java API callback information class constants, see {\hyperref[\detokenize{information:chapinfo-cbc}]{\sphinxcrossref{\DUrole{std,std-ref}{Callback Information}}}}.

\sphinxAtStartPar
In the Java API, callback\sphinxhyphen{}related information constants are defined in the \sphinxcode{\sphinxupquote{CbInfo}} class.
Users can access information constants through the prefix \sphinxcode{\sphinxupquote{copt}} in the namespace (usually can be omitted) \sphinxcode{\sphinxupquote{copt.CbInfo.}}

\sphinxAtStartPar
For instance, \sphinxcode{\sphinxupquote{copt.CbInfo.BestObj}} is the current best objective.

\subsection{Parameters}
\label{\detokenize{javaapiref:parameters}}\label{\detokenize{javaapiref:chapjavaapiref-params}}
\sphinxAtStartPar
For the contents of Java parameters constants, see {\hyperref[\detokenize{parameter:chapparams}]{\sphinxcrossref{\DUrole{std,std-ref}{Parameters}}}}.

\sphinxAtStartPar
All COPT Java parameters are defined in DblParam and IntParam classes. User may refer double parameters by \sphinxcode{\sphinxupquote{copt.DblParam.XXXX}}, and integer parameters by \sphinxcode{\sphinxupquote{copt.IntParam.XXXX}}.

\sphinxAtStartPar
In the Java API, user can get and set the parameter value by specifying the parameter name.
The provided functions are as follows, please refer to {\hyperref[\detokenize{javaapiref:chapjavaapiref-model}]{\sphinxcrossref{\DUrole{std,std-ref}{Java API: Model Class}}}} for details.
\begin{itemize}
\item {} 
\sphinxAtStartPar
Get detailed information of the specified parameter (current value/max/min): \sphinxcode{\sphinxupquote{Model.getParamInfo()}}

\item {} 
\sphinxAtStartPar
Get the current value of the specified integer/double parameter: \sphinxcode{\sphinxupquote{Model.getIntParam()}} / \sphinxcode{\sphinxupquote{Model.getDblParam()}}

\item {} 
\sphinxAtStartPar
Set the specified integer/double parameter value: \sphinxcode{\sphinxupquote{Model.setIntParam()}} / \sphinxcode{\sphinxupquote{Model.setDblParam()}}

\end{itemize}

\section{Java Modeling Classes}
\label{\detokenize{javaapiref:java-modeling-classes}}\label{\detokenize{javaapiref:chapjavaapiref-class}}
\sphinxAtStartPar
This chapter documents COPT Java interface. Users may refer to
Java classes described below for details of how to construct and
solve Java models.

\subsection{Envr}
\label{\detokenize{javaapiref:envr}}\label{\detokenize{javaapiref:chapjavaapiref-envr}}
\sphinxAtStartPar
Essentially, any Java application using Cardinal Optimizer should start with a
COPT environment. COPT models are always associated with a COPT environment.
User must create an environment object before populating models.
User generally only need a single environment object in program.

\sphinxstepscope

\subsubsection{Envr.Envr()}
\label{\detokenize{javaapi/Envr:envr-envr}}\label{\detokenize{javaapi/Envr::doc}}\begin{quote}

\sphinxAtStartPar
Constructor of COPT Envr object.

\sphinxAtStartPar
\sphinxstylestrong{Synopsis}
\begin{quote}

\sphinxAtStartPar
\sphinxcode{\sphinxupquote{Envr()}}
\end{quote}
\end{quote}

\subsubsection{Envr.Envr()}
\label{\detokenize{javaapi/Envr:id1}}\begin{quote}

\sphinxAtStartPar
Constructor of COPT Envr object, given a license folder.

\sphinxAtStartPar
\sphinxstylestrong{Synopsis}
\begin{quote}

\sphinxAtStartPar
\sphinxcode{\sphinxupquote{Envr(String licDir)}}
\end{quote}

\sphinxAtStartPar
\sphinxstylestrong{Arguments}
\begin{quote}

\sphinxAtStartPar
\sphinxcode{\sphinxupquote{licDir}}: directory having local license or client config file.
\end{quote}
\end{quote}

\subsubsection{Envr.Envr()}
\label{\detokenize{javaapi/Envr:id2}}\begin{quote}

\sphinxAtStartPar
Constructor of COPT Envr object, given an Envr config object.

\sphinxAtStartPar
\sphinxstylestrong{Synopsis}
\begin{quote}

\sphinxAtStartPar
\sphinxcode{\sphinxupquote{Envr(EnvrConfig config)}}
\end{quote}

\sphinxAtStartPar
\sphinxstylestrong{Arguments}
\begin{quote}

\sphinxAtStartPar
\sphinxcode{\sphinxupquote{config}}: Envr config object holding settings for remote connection.
\end{quote}
\end{quote}

\subsubsection{Envr.bindNumaCpu()}
\label{\detokenize{javaapi/Envr:envr-bindnumacpu}}\begin{quote}

\sphinxAtStartPar
Bind the CPUs for the current process to a NUMA node.

\sphinxAtStartPar
\sphinxstylestrong{Synopsis}
\begin{quote}

\sphinxAtStartPar
\sphinxcode{\sphinxupquote{void bindNumaCpu(int numaNode)}}
\end{quote}

\sphinxAtStartPar
\sphinxstylestrong{Arguments}
\begin{quote}

\sphinxAtStartPar
\sphinxcode{\sphinxupquote{numaNode}}: ID of a NUMA node.
\end{quote}
\end{quote}

\subsubsection{Envr.bindNumaMem()}
\label{\detokenize{javaapi/Envr:envr-bindnumamem}}\begin{quote}

\sphinxAtStartPar
Bind memory for the current process to a NUMA node (Linux only).

\sphinxAtStartPar
\sphinxstylestrong{Synopsis}
\begin{quote}

\sphinxAtStartPar
\sphinxcode{\sphinxupquote{void bindNumaMem(int numaNode)}}
\end{quote}

\sphinxAtStartPar
\sphinxstylestrong{Arguments}
\begin{quote}

\sphinxAtStartPar
\sphinxcode{\sphinxupquote{numaNode}}: the ID of a NUMA node.
\end{quote}
\end{quote}

\subsubsection{Envr.close()}
\label{\detokenize{javaapi/Envr:envr-close}}\begin{quote}

\sphinxAtStartPar
close remote connection and token becomes invalid for all problems in current envr.

\sphinxAtStartPar
\sphinxstylestrong{Synopsis}
\begin{quote}

\sphinxAtStartPar
\sphinxcode{\sphinxupquote{void close()}}
\end{quote}
\end{quote}

\subsubsection{Envr.createModel()}
\label{\detokenize{javaapi/Envr:envr-createmodel}}\begin{quote}

\sphinxAtStartPar
Create a model object.

\sphinxAtStartPar
\sphinxstylestrong{Synopsis}
\begin{quote}

\sphinxAtStartPar
\sphinxcode{\sphinxupquote{Model createModel(String name)}}
\end{quote}

\sphinxAtStartPar
\sphinxstylestrong{Arguments}
\begin{quote}

\sphinxAtStartPar
\sphinxcode{\sphinxupquote{name}}: customized model name.
\end{quote}

\sphinxAtStartPar
\sphinxstylestrong{Return}
\begin{quote}

\sphinxAtStartPar
a model object.
\end{quote}
\end{quote}

\subsubsection{Envr.getCpuAffinity()}
\label{\detokenize{javaapi/Envr:envr-getcpuaffinity}}\begin{quote}

\sphinxAtStartPar
Get CPU affinity for the current process, which is saved in an integer array.

\sphinxAtStartPar
\sphinxstylestrong{Synopsis}
\begin{quote}

\sphinxAtStartPar
\sphinxcode{\sphinxupquote{int{[}{]} getCpuAffinity()}}
\end{quote}

\sphinxAtStartPar
\sphinxstylestrong{Return}
\begin{quote}

\sphinxAtStartPar
an integer array of CPU IDs.
\end{quote}
\end{quote}

\subsubsection{Envr.getNumaNodeCount()}
\label{\detokenize{javaapi/Envr:envr-getnumanodecount}}\begin{quote}

\sphinxAtStartPar
Get count of NUMA nodes.

\sphinxAtStartPar
\sphinxstylestrong{Synopsis}
\begin{quote}

\sphinxAtStartPar
\sphinxcode{\sphinxupquote{int getNumaNodeCount()}}
\end{quote}

\sphinxAtStartPar
\sphinxstylestrong{Return}
\begin{quote}

\sphinxAtStartPar
count of NUMA nodes.
\end{quote}
\end{quote}

\subsubsection{Envr.setCpuAffinity()}
\label{\detokenize{javaapi/Envr:envr-setcpuaffinity}}\begin{quote}

\sphinxAtStartPar
Set CPU affinity with given mask string.

\sphinxAtStartPar
\sphinxstylestrong{Synopsis}
\begin{quote}

\sphinxAtStartPar
\sphinxcode{\sphinxupquote{void setCpuAffinity(String hexMask)}}
\end{quote}

\sphinxAtStartPar
\sphinxstylestrong{Arguments}
\begin{quote}

\sphinxAtStartPar
\sphinxcode{\sphinxupquote{hexMask}}: CPU mask string of hexadecimal characters.
\end{quote}
\end{quote}

\subsection{EnvrConfig}
\label{\detokenize{javaapiref:envrconfig}}
\sphinxAtStartPar
If user connects to COPT remote services, such as floating token server or compute cluster,
it is necessary to add config settings with EnvrConfig object.

\sphinxstepscope

\subsubsection{EnvrConfig.EnvrConfig()}
\label{\detokenize{javaapi/EnvrConfig:envrconfig-envrconfig}}\label{\detokenize{javaapi/EnvrConfig::doc}}\begin{quote}

\sphinxAtStartPar
Constructor of envr config object.

\sphinxAtStartPar
\sphinxstylestrong{Synopsis}
\begin{quote}

\sphinxAtStartPar
\sphinxcode{\sphinxupquote{EnvrConfig()}}
\end{quote}
\end{quote}

\subsubsection{EnvrConfig.set()}
\label{\detokenize{javaapi/EnvrConfig:envrconfig-set}}\begin{quote}

\sphinxAtStartPar
Set config settings in terms of name\sphinxhyphen{}value pair.

\sphinxAtStartPar
\sphinxstylestrong{Synopsis}
\begin{quote}

\sphinxAtStartPar
\sphinxcode{\sphinxupquote{void set(String name, String value)}}
\end{quote}

\sphinxAtStartPar
\sphinxstylestrong{Arguments}
\begin{quote}

\sphinxAtStartPar
\sphinxcode{\sphinxupquote{name}}: keyword of a config setting.

\sphinxAtStartPar
\sphinxcode{\sphinxupquote{value}}: value of a config setting.
\end{quote}
\end{quote}

\subsection{Model}
\label{\detokenize{javaapiref:model}}\label{\detokenize{javaapiref:chapjavaapiref-model}}
\sphinxAtStartPar
In general, a COPT model consists of a set of variables, a (linear) objective
function on these variables, a set of constraints on there varaibles, etc.
COPT model class encapsulates all required methods for constructing a COPT model.

\sphinxstepscope

\subsubsection{Model.Model()}
\label{\detokenize{javaapi/Model:model-model}}\label{\detokenize{javaapi/Model::doc}}\begin{quote}

\sphinxAtStartPar
Constructor of model.

\sphinxAtStartPar
\sphinxstylestrong{Synopsis}
\begin{quote}

\sphinxAtStartPar
\sphinxcode{\sphinxupquote{Model(Envr env, String name)}}
\end{quote}

\sphinxAtStartPar
\sphinxstylestrong{Arguments}
\begin{quote}

\sphinxAtStartPar
\sphinxcode{\sphinxupquote{env}}: associated environment object.

\sphinxAtStartPar
\sphinxcode{\sphinxupquote{name}}: string of model name.
\end{quote}
\end{quote}

\subsubsection{Model.addAffineCone()}
\label{\detokenize{javaapi/Model:model-addaffinecone}}\begin{quote}

\sphinxAtStartPar
Add an affine cone constraint to model.

\sphinxAtStartPar
\sphinxstylestrong{Synopsis}
\begin{quote}

\sphinxAtStartPar
\sphinxcode{\sphinxupquote{AffineCone addAffineCone(AffineConeBuilder builder, String name)}}
\end{quote}

\sphinxAtStartPar
\sphinxstylestrong{Arguments}
\begin{quote}

\sphinxAtStartPar
\sphinxcode{\sphinxupquote{builder}}: builder for new affine cone constraint.

\sphinxAtStartPar
\sphinxcode{\sphinxupquote{name}}: name of new affine cone.
\end{quote}

\sphinxAtStartPar
\sphinxstylestrong{Return}
\begin{quote}

\sphinxAtStartPar
new affine cone constraint object.
\end{quote}
\end{quote}

\subsubsection{Model.addAffineCone()}
\label{\detokenize{javaapi/Model:id1}}\begin{quote}

\sphinxAtStartPar
Add an affine cone constraint to model.

\sphinxAtStartPar
\sphinxstylestrong{Synopsis}
\begin{quote}

\sphinxAtStartPar
\sphinxcode{\sphinxupquote{AffineCone addAffineCone(}}
\begin{quote}

\sphinxAtStartPar
\sphinxcode{\sphinxupquote{Expr{[}{]} exprs,}}

\sphinxAtStartPar
\sphinxcode{\sphinxupquote{int type,}}

\sphinxAtStartPar
\sphinxcode{\sphinxupquote{String name)}}
\end{quote}
\end{quote}

\sphinxAtStartPar
\sphinxstylestrong{Arguments}
\begin{quote}

\sphinxAtStartPar
\sphinxcode{\sphinxupquote{exprs}}: linear expressioins that participate in the affine cone constraint.

\sphinxAtStartPar
\sphinxcode{\sphinxupquote{type}}: type of an affine cone constraint.

\sphinxAtStartPar
\sphinxcode{\sphinxupquote{name}}: name of new affine cone.
\end{quote}

\sphinxAtStartPar
\sphinxstylestrong{Return}
\begin{quote}

\sphinxAtStartPar
new affine cone constraint object.
\end{quote}
\end{quote}

\subsubsection{Model.addAffineCone()}
\label{\detokenize{javaapi/Model:id2}}\begin{quote}

\sphinxAtStartPar
Add an affine cone constraint to model.

\sphinxAtStartPar
\sphinxstylestrong{Synopsis}
\begin{quote}

\sphinxAtStartPar
\sphinxcode{\sphinxupquote{AffineCone addAffineCone(}}
\begin{quote}

\sphinxAtStartPar
\sphinxcode{\sphinxupquote{PsdExpr{[}{]} exprs,}}

\sphinxAtStartPar
\sphinxcode{\sphinxupquote{int type,}}

\sphinxAtStartPar
\sphinxcode{\sphinxupquote{String name)}}
\end{quote}
\end{quote}

\sphinxAtStartPar
\sphinxstylestrong{Arguments}
\begin{quote}

\sphinxAtStartPar
\sphinxcode{\sphinxupquote{exprs}}: PSD expressioins that participate in the affine cone constraint.

\sphinxAtStartPar
\sphinxcode{\sphinxupquote{type}}: type of an affine cone constraint.

\sphinxAtStartPar
\sphinxcode{\sphinxupquote{name}}: name of new affine cone.
\end{quote}

\sphinxAtStartPar
\sphinxstylestrong{Return}
\begin{quote}

\sphinxAtStartPar
new affine cone constraint object.
\end{quote}
\end{quote}

\subsubsection{Model.addCone()}
\label{\detokenize{javaapi/Model:model-addcone}}\begin{quote}

\sphinxAtStartPar
Add a cone constraint to model.

\sphinxAtStartPar
\sphinxstylestrong{Synopsis}
\begin{quote}

\sphinxAtStartPar
\sphinxcode{\sphinxupquote{Cone addCone(}}
\begin{quote}

\sphinxAtStartPar
\sphinxcode{\sphinxupquote{int dim,}}

\sphinxAtStartPar
\sphinxcode{\sphinxupquote{int type,}}

\sphinxAtStartPar
\sphinxcode{\sphinxupquote{char{[}{]} pvtype,}}

\sphinxAtStartPar
\sphinxcode{\sphinxupquote{String prefix)}}
\end{quote}
\end{quote}

\sphinxAtStartPar
\sphinxstylestrong{Arguments}
\begin{quote}

\sphinxAtStartPar
\sphinxcode{\sphinxupquote{dim}}: dimension of the cone constraint.

\sphinxAtStartPar
\sphinxcode{\sphinxupquote{type}}: type of a cone constraint.

\sphinxAtStartPar
\sphinxcode{\sphinxupquote{pvtype}}: types of variables in the cone.

\sphinxAtStartPar
\sphinxcode{\sphinxupquote{prefix}}: name prefix of variables in the cone.
\end{quote}

\sphinxAtStartPar
\sphinxstylestrong{Return}
\begin{quote}

\sphinxAtStartPar
new cone constraint object.
\end{quote}
\end{quote}

\subsubsection{Model.addCone()}
\label{\detokenize{javaapi/Model:id3}}\begin{quote}

\sphinxAtStartPar
Add a cone constraint to model.

\sphinxAtStartPar
\sphinxstylestrong{Synopsis}
\begin{quote}

\sphinxAtStartPar
\sphinxcode{\sphinxupquote{Cone addCone(ConeBuilder builder)}}
\end{quote}

\sphinxAtStartPar
\sphinxstylestrong{Arguments}
\begin{quote}

\sphinxAtStartPar
\sphinxcode{\sphinxupquote{builder}}: builder for new cone constraint.
\end{quote}

\sphinxAtStartPar
\sphinxstylestrong{Return}
\begin{quote}

\sphinxAtStartPar
new cone constraint object.
\end{quote}
\end{quote}

\subsubsection{Model.addCone()}
\label{\detokenize{javaapi/Model:id4}}\begin{quote}

\sphinxAtStartPar
Add a cone constraint to model.

\sphinxAtStartPar
\sphinxstylestrong{Synopsis}
\begin{quote}

\sphinxAtStartPar
\sphinxcode{\sphinxupquote{Cone addCone(Var{[}{]} vars, int type)}}
\end{quote}

\sphinxAtStartPar
\sphinxstylestrong{Arguments}
\begin{quote}

\sphinxAtStartPar
\sphinxcode{\sphinxupquote{vars}}: variables that participate in the cone constraint.

\sphinxAtStartPar
\sphinxcode{\sphinxupquote{type}}: type of a cone constraint.
\end{quote}

\sphinxAtStartPar
\sphinxstylestrong{Return}
\begin{quote}

\sphinxAtStartPar
new cone constraint object.
\end{quote}
\end{quote}

\subsubsection{Model.addCone()}
\label{\detokenize{javaapi/Model:id5}}\begin{quote}

\sphinxAtStartPar
Add a cone constraint to model.

\sphinxAtStartPar
\sphinxstylestrong{Synopsis}
\begin{quote}

\sphinxAtStartPar
\sphinxcode{\sphinxupquote{Cone addCone(VarArray vars, int type)}}
\end{quote}

\sphinxAtStartPar
\sphinxstylestrong{Arguments}
\begin{quote}

\sphinxAtStartPar
\sphinxcode{\sphinxupquote{vars}}: variables that participate in the cone constraint.

\sphinxAtStartPar
\sphinxcode{\sphinxupquote{type}}: type of a cone constraint.
\end{quote}

\sphinxAtStartPar
\sphinxstylestrong{Return}
\begin{quote}

\sphinxAtStartPar
new cone constraint object.
\end{quote}
\end{quote}

\subsubsection{Model.addConstr()}
\label{\detokenize{javaapi/Model:model-addconstr}}\begin{quote}

\sphinxAtStartPar
Add a linear constraint to model.

\sphinxAtStartPar
\sphinxstylestrong{Synopsis}
\begin{quote}

\sphinxAtStartPar
\sphinxcode{\sphinxupquote{Constraint addConstr(}}
\begin{quote}

\sphinxAtStartPar
\sphinxcode{\sphinxupquote{Expr expr,}}

\sphinxAtStartPar
\sphinxcode{\sphinxupquote{char sense,}}

\sphinxAtStartPar
\sphinxcode{\sphinxupquote{double rhs,}}

\sphinxAtStartPar
\sphinxcode{\sphinxupquote{String name)}}
\end{quote}
\end{quote}

\sphinxAtStartPar
\sphinxstylestrong{Arguments}
\begin{quote}

\sphinxAtStartPar
\sphinxcode{\sphinxupquote{expr}}: expression for the new contraint.

\sphinxAtStartPar
\sphinxcode{\sphinxupquote{sense}}: sense for new linear constraint, other than range sense.

\sphinxAtStartPar
\sphinxcode{\sphinxupquote{rhs}}: right hand side value for the new constraint.

\sphinxAtStartPar
\sphinxcode{\sphinxupquote{name}}: name of new constraint.
\end{quote}

\sphinxAtStartPar
\sphinxstylestrong{Return}
\begin{quote}

\sphinxAtStartPar
new constraint object.
\end{quote}
\end{quote}

\subsubsection{Model.addConstr()}
\label{\detokenize{javaapi/Model:id6}}\begin{quote}

\sphinxAtStartPar
Add a linear constraint to model.

\sphinxAtStartPar
\sphinxstylestrong{Synopsis}
\begin{quote}

\sphinxAtStartPar
\sphinxcode{\sphinxupquote{Constraint addConstr(}}
\begin{quote}

\sphinxAtStartPar
\sphinxcode{\sphinxupquote{Expr expr,}}

\sphinxAtStartPar
\sphinxcode{\sphinxupquote{char sense,}}

\sphinxAtStartPar
\sphinxcode{\sphinxupquote{Var var,}}

\sphinxAtStartPar
\sphinxcode{\sphinxupquote{String name)}}
\end{quote}
\end{quote}

\sphinxAtStartPar
\sphinxstylestrong{Arguments}
\begin{quote}

\sphinxAtStartPar
\sphinxcode{\sphinxupquote{expr}}: expression for the new contraint.

\sphinxAtStartPar
\sphinxcode{\sphinxupquote{sense}}: sense for new linear constraint, other than range sense.

\sphinxAtStartPar
\sphinxcode{\sphinxupquote{var}}: varible for the new constraint.

\sphinxAtStartPar
\sphinxcode{\sphinxupquote{name}}: name of new constraint.
\end{quote}

\sphinxAtStartPar
\sphinxstylestrong{Return}
\begin{quote}

\sphinxAtStartPar
new constraint object.
\end{quote}
\end{quote}

\subsubsection{Model.addConstr()}
\label{\detokenize{javaapi/Model:id7}}\begin{quote}

\sphinxAtStartPar
Add a linear constraint to model.

\sphinxAtStartPar
\sphinxstylestrong{Synopsis}
\begin{quote}

\sphinxAtStartPar
\sphinxcode{\sphinxupquote{Constraint addConstr(}}
\begin{quote}

\sphinxAtStartPar
\sphinxcode{\sphinxupquote{Expr lhs,}}

\sphinxAtStartPar
\sphinxcode{\sphinxupquote{char sense,}}

\sphinxAtStartPar
\sphinxcode{\sphinxupquote{Expr rhs,}}

\sphinxAtStartPar
\sphinxcode{\sphinxupquote{String name)}}
\end{quote}
\end{quote}

\sphinxAtStartPar
\sphinxstylestrong{Arguments}
\begin{quote}

\sphinxAtStartPar
\sphinxcode{\sphinxupquote{lhs}}: left hand side expression for the new constraint.

\sphinxAtStartPar
\sphinxcode{\sphinxupquote{sense}}: sense for new linear constraint, other than range sense.

\sphinxAtStartPar
\sphinxcode{\sphinxupquote{rhs}}: right hand side expression for the new constraint.

\sphinxAtStartPar
\sphinxcode{\sphinxupquote{name}}: name of new constraint.
\end{quote}

\sphinxAtStartPar
\sphinxstylestrong{Return}
\begin{quote}

\sphinxAtStartPar
new constraint object.
\end{quote}
\end{quote}

\subsubsection{Model.addConstr()}
\label{\detokenize{javaapi/Model:id8}}\begin{quote}

\sphinxAtStartPar
Add a linear constraint to model.

\sphinxAtStartPar
\sphinxstylestrong{Synopsis}
\begin{quote}

\sphinxAtStartPar
\sphinxcode{\sphinxupquote{Constraint addConstr(}}
\begin{quote}

\sphinxAtStartPar
\sphinxcode{\sphinxupquote{Expr expr,}}

\sphinxAtStartPar
\sphinxcode{\sphinxupquote{double lb,}}

\sphinxAtStartPar
\sphinxcode{\sphinxupquote{double ub,}}

\sphinxAtStartPar
\sphinxcode{\sphinxupquote{String name)}}
\end{quote}
\end{quote}

\sphinxAtStartPar
\sphinxstylestrong{Arguments}
\begin{quote}

\sphinxAtStartPar
\sphinxcode{\sphinxupquote{expr}}: expression for the new constraint.

\sphinxAtStartPar
\sphinxcode{\sphinxupquote{lb}}: lower bound for the new constraint.

\sphinxAtStartPar
\sphinxcode{\sphinxupquote{ub}}: upper bound for the new constraint

\sphinxAtStartPar
\sphinxcode{\sphinxupquote{name}}: name of new constraint.
\end{quote}

\sphinxAtStartPar
\sphinxstylestrong{Return}
\begin{quote}

\sphinxAtStartPar
new constraint object.
\end{quote}
\end{quote}

\subsubsection{Model.addConstr()}
\label{\detokenize{javaapi/Model:id9}}\begin{quote}

\sphinxAtStartPar
Add a linear constraint to a model.

\sphinxAtStartPar
\sphinxstylestrong{Synopsis}
\begin{quote}

\sphinxAtStartPar
\sphinxcode{\sphinxupquote{Constraint addConstr(ConstrBuilder builder, String name)}}
\end{quote}

\sphinxAtStartPar
\sphinxstylestrong{Arguments}
\begin{quote}

\sphinxAtStartPar
\sphinxcode{\sphinxupquote{builder}}: builder for the new constraint.

\sphinxAtStartPar
\sphinxcode{\sphinxupquote{name}}: name of new constraint.
\end{quote}

\sphinxAtStartPar
\sphinxstylestrong{Return}
\begin{quote}

\sphinxAtStartPar
new constraint object.
\end{quote}
\end{quote}

\subsubsection{Model.addConstrs()}
\label{\detokenize{javaapi/Model:model-addconstrs}}\begin{quote}

\sphinxAtStartPar
Add linear constraints to model.

\sphinxAtStartPar
\sphinxstylestrong{Synopsis}
\begin{quote}

\sphinxAtStartPar
\sphinxcode{\sphinxupquote{ConstrArray addConstrs(}}
\begin{quote}

\sphinxAtStartPar
\sphinxcode{\sphinxupquote{int count,}}

\sphinxAtStartPar
\sphinxcode{\sphinxupquote{char{[}{]} senses,}}

\sphinxAtStartPar
\sphinxcode{\sphinxupquote{double{[}{]} rhss,}}

\sphinxAtStartPar
\sphinxcode{\sphinxupquote{String prefix)}}
\end{quote}
\end{quote}

\sphinxAtStartPar
\sphinxstylestrong{Arguments}
\begin{quote}

\sphinxAtStartPar
\sphinxcode{\sphinxupquote{count}}: number of constraints added to model.

\sphinxAtStartPar
\sphinxcode{\sphinxupquote{senses}}: sense array for new linear constraints, other than range sense.

\sphinxAtStartPar
\sphinxcode{\sphinxupquote{rhss}}: right hand side values for new variables.

\sphinxAtStartPar
\sphinxcode{\sphinxupquote{prefix}}: name prefix for new constraints.
\end{quote}

\sphinxAtStartPar
\sphinxstylestrong{Return}
\begin{quote}

\sphinxAtStartPar
array of new constraint objects.
\end{quote}
\end{quote}

\subsubsection{Model.addConstrs()}
\label{\detokenize{javaapi/Model:id10}}\begin{quote}

\sphinxAtStartPar
Add linear constraints to a model.

\sphinxAtStartPar
\sphinxstylestrong{Synopsis}
\begin{quote}

\sphinxAtStartPar
\sphinxcode{\sphinxupquote{ConstrArray addConstrs(}}
\begin{quote}

\sphinxAtStartPar
\sphinxcode{\sphinxupquote{int count,}}

\sphinxAtStartPar
\sphinxcode{\sphinxupquote{double{[}{]} lbs,}}

\sphinxAtStartPar
\sphinxcode{\sphinxupquote{double{[}{]} ubs,}}

\sphinxAtStartPar
\sphinxcode{\sphinxupquote{String prefix)}}
\end{quote}
\end{quote}

\sphinxAtStartPar
\sphinxstylestrong{Arguments}
\begin{quote}

\sphinxAtStartPar
\sphinxcode{\sphinxupquote{count}}: number of constraints added to the model.

\sphinxAtStartPar
\sphinxcode{\sphinxupquote{lbs}}: lower bounds of new constraints.

\sphinxAtStartPar
\sphinxcode{\sphinxupquote{ubs}}: upper bounds of new constraints.

\sphinxAtStartPar
\sphinxcode{\sphinxupquote{prefix}}: name prefix for new constraints.
\end{quote}

\sphinxAtStartPar
\sphinxstylestrong{Return}
\begin{quote}

\sphinxAtStartPar
array of new constraint objects.
\end{quote}
\end{quote}

\subsubsection{Model.addConstrs()}
\label{\detokenize{javaapi/Model:id11}}\begin{quote}

\sphinxAtStartPar
Add linear constraints to a model.

\sphinxAtStartPar
\sphinxstylestrong{Synopsis}
\begin{quote}

\sphinxAtStartPar
\sphinxcode{\sphinxupquote{ConstrArray addConstrs(ConstrBuilderArray builders, String prefix)}}
\end{quote}

\sphinxAtStartPar
\sphinxstylestrong{Arguments}
\begin{quote}

\sphinxAtStartPar
\sphinxcode{\sphinxupquote{builders}}: builders for new constraints.

\sphinxAtStartPar
\sphinxcode{\sphinxupquote{prefix}}: name prefix for new constraints.
\end{quote}

\sphinxAtStartPar
\sphinxstylestrong{Return}
\begin{quote}

\sphinxAtStartPar
array of new constraint objects.
\end{quote}
\end{quote}

\subsubsection{Model.addDenseMat()}
\label{\detokenize{javaapi/Model:model-adddensemat}}\begin{quote}

\sphinxAtStartPar
Add a dense symmetric matrix to a model.

\sphinxAtStartPar
\sphinxstylestrong{Synopsis}
\begin{quote}

\sphinxAtStartPar
\sphinxcode{\sphinxupquote{SymMatrix addDenseMat(int dim, double{[}{]} vals)}}
\end{quote}

\sphinxAtStartPar
\sphinxstylestrong{Arguments}
\begin{quote}

\sphinxAtStartPar
\sphinxcode{\sphinxupquote{dim}}: dimension of the dense symmetric matrix.

\sphinxAtStartPar
\sphinxcode{\sphinxupquote{vals}}: array of non\sphinxhyphen{}zero elements, filled in column\sphinxhyphen{}wise up to len or max length of symmetric matrix.
\end{quote}

\sphinxAtStartPar
\sphinxstylestrong{Return}
\begin{quote}

\sphinxAtStartPar
new symmetric matrix object.
\end{quote}
\end{quote}

\subsubsection{Model.addDenseMat()}
\label{\detokenize{javaapi/Model:id12}}\begin{quote}

\sphinxAtStartPar
Add a dense symmetric matrix to a model.

\sphinxAtStartPar
\sphinxstylestrong{Synopsis}
\begin{quote}

\sphinxAtStartPar
\sphinxcode{\sphinxupquote{SymMatrix addDenseMat(int dim, double val)}}
\end{quote}

\sphinxAtStartPar
\sphinxstylestrong{Arguments}
\begin{quote}

\sphinxAtStartPar
\sphinxcode{\sphinxupquote{dim}}: dimension of dense symmetric matrix.

\sphinxAtStartPar
\sphinxcode{\sphinxupquote{val}}: value to fill dense symmetric matrix.
\end{quote}

\sphinxAtStartPar
\sphinxstylestrong{Return}
\begin{quote}

\sphinxAtStartPar
new symmetric matrix object.
\end{quote}
\end{quote}

\subsubsection{Model.addDiagMat()}
\label{\detokenize{javaapi/Model:model-adddiagmat}}\begin{quote}

\sphinxAtStartPar
Add a diagonal matrix to a model.

\sphinxAtStartPar
\sphinxstylestrong{Synopsis}
\begin{quote}

\sphinxAtStartPar
\sphinxcode{\sphinxupquote{SymMatrix addDiagMat(int dim, double{[}{]} vals)}}
\end{quote}

\sphinxAtStartPar
\sphinxstylestrong{Arguments}
\begin{quote}

\sphinxAtStartPar
\sphinxcode{\sphinxupquote{dim}}: dimension of diagonal matrix.

\sphinxAtStartPar
\sphinxcode{\sphinxupquote{vals}}: array of values of diagonal elements.
\end{quote}

\sphinxAtStartPar
\sphinxstylestrong{Return}
\begin{quote}

\sphinxAtStartPar
new diagonal matrix object.
\end{quote}
\end{quote}

\subsubsection{Model.addDiagMat()}
\label{\detokenize{javaapi/Model:id13}}\begin{quote}

\sphinxAtStartPar
Add a diagonal matrix to a model.

\sphinxAtStartPar
\sphinxstylestrong{Synopsis}
\begin{quote}

\sphinxAtStartPar
\sphinxcode{\sphinxupquote{SymMatrix addDiagMat(}}
\begin{quote}

\sphinxAtStartPar
\sphinxcode{\sphinxupquote{int dim,}}

\sphinxAtStartPar
\sphinxcode{\sphinxupquote{double val,}}

\sphinxAtStartPar
\sphinxcode{\sphinxupquote{int offset)}}
\end{quote}
\end{quote}

\sphinxAtStartPar
\sphinxstylestrong{Arguments}
\begin{quote}

\sphinxAtStartPar
\sphinxcode{\sphinxupquote{dim}}: dimension of diagonal matrix.

\sphinxAtStartPar
\sphinxcode{\sphinxupquote{val}}: value to fill diagonal elements.

\sphinxAtStartPar
\sphinxcode{\sphinxupquote{offset}}: shift distance against diagonal line.
\end{quote}

\sphinxAtStartPar
\sphinxstylestrong{Return}
\begin{quote}

\sphinxAtStartPar
new diagonal matrix object.
\end{quote}
\end{quote}

\subsubsection{Model.addDiagMat()}
\label{\detokenize{javaapi/Model:id14}}\begin{quote}

\sphinxAtStartPar
Add a diagonal matrix to a model.

\sphinxAtStartPar
\sphinxstylestrong{Synopsis}
\begin{quote}

\sphinxAtStartPar
\sphinxcode{\sphinxupquote{SymMatrix addDiagMat(}}
\begin{quote}

\sphinxAtStartPar
\sphinxcode{\sphinxupquote{int dim,}}

\sphinxAtStartPar
\sphinxcode{\sphinxupquote{double{[}{]} vals,}}

\sphinxAtStartPar
\sphinxcode{\sphinxupquote{int offset)}}
\end{quote}
\end{quote}

\sphinxAtStartPar
\sphinxstylestrong{Arguments}
\begin{quote}

\sphinxAtStartPar
\sphinxcode{\sphinxupquote{dim}}: dimension of diagonal matrix.

\sphinxAtStartPar
\sphinxcode{\sphinxupquote{vals}}: array of values of diagonal elements.

\sphinxAtStartPar
\sphinxcode{\sphinxupquote{offset}}: shift distance against diagonal line.
\end{quote}

\sphinxAtStartPar
\sphinxstylestrong{Return}
\begin{quote}

\sphinxAtStartPar
new diagonal matrix object.
\end{quote}
\end{quote}

\subsubsection{Model.addDiagMat()}
\label{\detokenize{javaapi/Model:id15}}\begin{quote}

\sphinxAtStartPar
Add a diagonal matrix to a model.

\sphinxAtStartPar
\sphinxstylestrong{Synopsis}
\begin{quote}

\sphinxAtStartPar
\sphinxcode{\sphinxupquote{SymMatrix addDiagMat(int dim, double val)}}
\end{quote}

\sphinxAtStartPar
\sphinxstylestrong{Arguments}
\begin{quote}

\sphinxAtStartPar
\sphinxcode{\sphinxupquote{dim}}: dimension of diagonal matrix.

\sphinxAtStartPar
\sphinxcode{\sphinxupquote{val}}: value to fill diagonal elements.
\end{quote}

\sphinxAtStartPar
\sphinxstylestrong{Return}
\begin{quote}

\sphinxAtStartPar
new diagonal matrix object.
\end{quote}
\end{quote}

\subsubsection{Model.addExpCone()}
\label{\detokenize{javaapi/Model:model-addexpcone}}\begin{quote}

\sphinxAtStartPar
Add an exponential cone constraint to model.

\sphinxAtStartPar
\sphinxstylestrong{Synopsis}
\begin{quote}

\sphinxAtStartPar
\sphinxcode{\sphinxupquote{ExpCone addExpCone(}}
\begin{quote}

\sphinxAtStartPar
\sphinxcode{\sphinxupquote{int type,}}

\sphinxAtStartPar
\sphinxcode{\sphinxupquote{char{[}{]} pvtype,}}

\sphinxAtStartPar
\sphinxcode{\sphinxupquote{String prefix)}}
\end{quote}
\end{quote}

\sphinxAtStartPar
\sphinxstylestrong{Arguments}
\begin{quote}

\sphinxAtStartPar
\sphinxcode{\sphinxupquote{type}}: type of an exponential cone constraint.

\sphinxAtStartPar
\sphinxcode{\sphinxupquote{pvtype}}: types of variables in the exponential cone.

\sphinxAtStartPar
\sphinxcode{\sphinxupquote{prefix}}: name prefix of variables in the exponential cone.
\end{quote}

\sphinxAtStartPar
\sphinxstylestrong{Return}
\begin{quote}

\sphinxAtStartPar
new exponential cone constraint object.
\end{quote}
\end{quote}

\subsubsection{Model.addExpCone()}
\label{\detokenize{javaapi/Model:id16}}\begin{quote}

\sphinxAtStartPar
Add an exponential cone constraint to model.

\sphinxAtStartPar
\sphinxstylestrong{Synopsis}
\begin{quote}

\sphinxAtStartPar
\sphinxcode{\sphinxupquote{ExpCone addExpCone(ExpConeBuilder builder)}}
\end{quote}

\sphinxAtStartPar
\sphinxstylestrong{Arguments}
\begin{quote}

\sphinxAtStartPar
\sphinxcode{\sphinxupquote{builder}}: builder for new exponential cone constraint.
\end{quote}

\sphinxAtStartPar
\sphinxstylestrong{Return}
\begin{quote}

\sphinxAtStartPar
new exponential cone constraint object.
\end{quote}
\end{quote}

\subsubsection{Model.addExpCone()}
\label{\detokenize{javaapi/Model:id17}}\begin{quote}

\sphinxAtStartPar
Add an exponential cone constraint to model.

\sphinxAtStartPar
\sphinxstylestrong{Synopsis}
\begin{quote}

\sphinxAtStartPar
\sphinxcode{\sphinxupquote{ExpCone addExpCone(Var{[}{]} vars, int type)}}
\end{quote}

\sphinxAtStartPar
\sphinxstylestrong{Arguments}
\begin{quote}

\sphinxAtStartPar
\sphinxcode{\sphinxupquote{vars}}: variables that participate in the exponential cone constraint.

\sphinxAtStartPar
\sphinxcode{\sphinxupquote{type}}: type of an exponential cone constraint.
\end{quote}

\sphinxAtStartPar
\sphinxstylestrong{Return}
\begin{quote}

\sphinxAtStartPar
new exponential cone constraint object.
\end{quote}
\end{quote}

\subsubsection{Model.addExpCone()}
\label{\detokenize{javaapi/Model:id18}}\begin{quote}

\sphinxAtStartPar
Add an exponential cone constraint to model.

\sphinxAtStartPar
\sphinxstylestrong{Synopsis}
\begin{quote}

\sphinxAtStartPar
\sphinxcode{\sphinxupquote{ExpCone addExpCone(VarArray vars, int type)}}
\end{quote}

\sphinxAtStartPar
\sphinxstylestrong{Arguments}
\begin{quote}

\sphinxAtStartPar
\sphinxcode{\sphinxupquote{vars}}: variables that participate in the exponential cone constraint.

\sphinxAtStartPar
\sphinxcode{\sphinxupquote{type}}: type of an exponential cone constraint.
\end{quote}

\sphinxAtStartPar
\sphinxstylestrong{Return}
\begin{quote}

\sphinxAtStartPar
new exponential cone constraint object.
\end{quote}
\end{quote}

\subsubsection{Model.addEyeMat()}
\label{\detokenize{javaapi/Model:model-addeyemat}}\begin{quote}

\sphinxAtStartPar
Add an identity matrix to a model.

\sphinxAtStartPar
\sphinxstylestrong{Synopsis}
\begin{quote}

\sphinxAtStartPar
\sphinxcode{\sphinxupquote{SymMatrix addEyeMat(int dim)}}
\end{quote}

\sphinxAtStartPar
\sphinxstylestrong{Arguments}
\begin{quote}

\sphinxAtStartPar
\sphinxcode{\sphinxupquote{dim}}: dimension of identity matrix.
\end{quote}

\sphinxAtStartPar
\sphinxstylestrong{Return}
\begin{quote}

\sphinxAtStartPar
new identity matrix object.
\end{quote}
\end{quote}

\subsubsection{Model.addGenConstrIndicator()}
\label{\detokenize{javaapi/Model:model-addgenconstrindicator}}\begin{quote}

\sphinxAtStartPar
Add a general constraint of type indicator to model.

\sphinxAtStartPar
\sphinxstylestrong{Synopsis}
\begin{quote}

\sphinxAtStartPar
\sphinxcode{\sphinxupquote{GenConstr addGenConstrIndicator(GenConstrBuilder builder, String name)}}
\end{quote}

\sphinxAtStartPar
\sphinxstylestrong{Arguments}
\begin{quote}

\sphinxAtStartPar
\sphinxcode{\sphinxupquote{builder}}: builder for the general constraint.

\sphinxAtStartPar
\sphinxcode{\sphinxupquote{name}}: name of new general constraint.
\end{quote}

\sphinxAtStartPar
\sphinxstylestrong{Return}
\begin{quote}

\sphinxAtStartPar
new general constraint object of type indicator.
\end{quote}
\end{quote}

\subsubsection{Model.addGenConstrIndicator()}
\label{\detokenize{javaapi/Model:id19}}\begin{quote}

\sphinxAtStartPar
Add a general constraint of type indicator to model.

\sphinxAtStartPar
\sphinxstylestrong{Synopsis}
\begin{quote}

\sphinxAtStartPar
\sphinxcode{\sphinxupquote{GenConstr addGenConstrIndicator(}}
\begin{quote}

\sphinxAtStartPar
\sphinxcode{\sphinxupquote{Var binvar,}}

\sphinxAtStartPar
\sphinxcode{\sphinxupquote{int binval,}}

\sphinxAtStartPar
\sphinxcode{\sphinxupquote{ConstrBuilder builder,}}

\sphinxAtStartPar
\sphinxcode{\sphinxupquote{int type,}}

\sphinxAtStartPar
\sphinxcode{\sphinxupquote{String name)}}
\end{quote}
\end{quote}

\sphinxAtStartPar
\sphinxstylestrong{Arguments}
\begin{quote}

\sphinxAtStartPar
\sphinxcode{\sphinxupquote{binvar}}: binary indicator variable.

\sphinxAtStartPar
\sphinxcode{\sphinxupquote{binval}}: value for binary indicator variable that force a linear constraint to be satisfied(0 or 1).

\sphinxAtStartPar
\sphinxcode{\sphinxupquote{builder}}: builder for linear constraint.

\sphinxAtStartPar
\sphinxcode{\sphinxupquote{type}}: type of general constraint.

\sphinxAtStartPar
\sphinxcode{\sphinxupquote{name}}: name of new general constraint.
\end{quote}

\sphinxAtStartPar
\sphinxstylestrong{Return}
\begin{quote}

\sphinxAtStartPar
new general constraint object of type indicator.
\end{quote}
\end{quote}

\subsubsection{Model.addGenConstrIndicator()}
\label{\detokenize{javaapi/Model:id20}}\begin{quote}

\sphinxAtStartPar
Add a general constraint of type indicator to model.

\sphinxAtStartPar
\sphinxstylestrong{Synopsis}
\begin{quote}

\sphinxAtStartPar
\sphinxcode{\sphinxupquote{GenConstr addGenConstrIndicator(}}
\begin{quote}

\sphinxAtStartPar
\sphinxcode{\sphinxupquote{Var binvar,}}

\sphinxAtStartPar
\sphinxcode{\sphinxupquote{int binval,}}

\sphinxAtStartPar
\sphinxcode{\sphinxupquote{Expr expr,}}

\sphinxAtStartPar
\sphinxcode{\sphinxupquote{char sense,}}

\sphinxAtStartPar
\sphinxcode{\sphinxupquote{double rhs,}}

\sphinxAtStartPar
\sphinxcode{\sphinxupquote{int type,}}

\sphinxAtStartPar
\sphinxcode{\sphinxupquote{String name)}}
\end{quote}
\end{quote}

\sphinxAtStartPar
\sphinxstylestrong{Arguments}
\begin{quote}

\sphinxAtStartPar
\sphinxcode{\sphinxupquote{binvar}}: binary indicator variable.

\sphinxAtStartPar
\sphinxcode{\sphinxupquote{binval}}: value for binary indicator variable that force a linear constraint to be satisfied(0 or 1).

\sphinxAtStartPar
\sphinxcode{\sphinxupquote{expr}}: expression for new linear contraint.

\sphinxAtStartPar
\sphinxcode{\sphinxupquote{sense}}: sense for new linear constraint.

\sphinxAtStartPar
\sphinxcode{\sphinxupquote{rhs}}: right hand side value for new linear constraint.

\sphinxAtStartPar
\sphinxcode{\sphinxupquote{type}}: type of general constraint.

\sphinxAtStartPar
\sphinxcode{\sphinxupquote{name}}: name of new general constraint.
\end{quote}

\sphinxAtStartPar
\sphinxstylestrong{Return}
\begin{quote}

\sphinxAtStartPar
new general constraint object of type indicator.
\end{quote}
\end{quote}

\subsubsection{Model.addGenConstrIndicators()}
\label{\detokenize{javaapi/Model:model-addgenconstrindicators}}\begin{quote}

\sphinxAtStartPar
Add general constraints to a model.

\sphinxAtStartPar
\sphinxstylestrong{Synopsis}
\begin{quote}

\sphinxAtStartPar
\sphinxcode{\sphinxupquote{GenConstrArray addGenConstrIndicators(GenConstrBuilderArray builders, String prefix)}}
\end{quote}

\sphinxAtStartPar
\sphinxstylestrong{Arguments}
\begin{quote}

\sphinxAtStartPar
\sphinxcode{\sphinxupquote{builders}}: builders for new general constraints.

\sphinxAtStartPar
\sphinxcode{\sphinxupquote{prefix}}: name prefix for new general constraints.
\end{quote}

\sphinxAtStartPar
\sphinxstylestrong{Return}
\begin{quote}

\sphinxAtStartPar
array of new general constraint objects.
\end{quote}
\end{quote}

\subsubsection{Model.addLazyConstr()}
\label{\detokenize{javaapi/Model:model-addlazyconstr}}\begin{quote}

\sphinxAtStartPar
Add a lazy constraint to model.

\sphinxAtStartPar
\sphinxstylestrong{Synopsis}
\begin{quote}

\sphinxAtStartPar
\sphinxcode{\sphinxupquote{void addLazyConstr(}}
\begin{quote}

\sphinxAtStartPar
\sphinxcode{\sphinxupquote{Expr lhs,}}

\sphinxAtStartPar
\sphinxcode{\sphinxupquote{char sense,}}

\sphinxAtStartPar
\sphinxcode{\sphinxupquote{double rhs,}}

\sphinxAtStartPar
\sphinxcode{\sphinxupquote{String name)}}
\end{quote}
\end{quote}

\sphinxAtStartPar
\sphinxstylestrong{Arguments}
\begin{quote}

\sphinxAtStartPar
\sphinxcode{\sphinxupquote{lhs}}: expression for lazy contraint.

\sphinxAtStartPar
\sphinxcode{\sphinxupquote{sense}}: sense for lazy constraint.

\sphinxAtStartPar
\sphinxcode{\sphinxupquote{rhs}}: right hand side value for lazy constraint.

\sphinxAtStartPar
\sphinxcode{\sphinxupquote{name}}: name of lazy constraint.
\end{quote}
\end{quote}

\subsubsection{Model.addLazyConstr()}
\label{\detokenize{javaapi/Model:id21}}\begin{quote}

\sphinxAtStartPar
Add a lazy constraint to model.

\sphinxAtStartPar
\sphinxstylestrong{Synopsis}
\begin{quote}

\sphinxAtStartPar
\sphinxcode{\sphinxupquote{void addLazyConstr(}}
\begin{quote}

\sphinxAtStartPar
\sphinxcode{\sphinxupquote{Expr lhs,}}

\sphinxAtStartPar
\sphinxcode{\sphinxupquote{char sense,}}

\sphinxAtStartPar
\sphinxcode{\sphinxupquote{Expr rhs,}}

\sphinxAtStartPar
\sphinxcode{\sphinxupquote{String name)}}
\end{quote}
\end{quote}

\sphinxAtStartPar
\sphinxstylestrong{Arguments}
\begin{quote}

\sphinxAtStartPar
\sphinxcode{\sphinxupquote{lhs}}: left hand side expression for lazy contraint.

\sphinxAtStartPar
\sphinxcode{\sphinxupquote{sense}}: sense for lazy constraint.

\sphinxAtStartPar
\sphinxcode{\sphinxupquote{rhs}}: right hand side expression for lazy contraint.

\sphinxAtStartPar
\sphinxcode{\sphinxupquote{name}}: name of lazy constraint.
\end{quote}
\end{quote}

\subsubsection{Model.addLazyConstr()}
\label{\detokenize{javaapi/Model:id22}}\begin{quote}

\sphinxAtStartPar
Add a lazy constraint to model.

\sphinxAtStartPar
\sphinxstylestrong{Synopsis}
\begin{quote}

\sphinxAtStartPar
\sphinxcode{\sphinxupquote{void addLazyConstr(ConstrBuilder builder, String name)}}
\end{quote}

\sphinxAtStartPar
\sphinxstylestrong{Arguments}
\begin{quote}

\sphinxAtStartPar
\sphinxcode{\sphinxupquote{builder}}: builder for lazy contraint.

\sphinxAtStartPar
\sphinxcode{\sphinxupquote{name}}: name of lazy constraint.
\end{quote}
\end{quote}

\subsubsection{Model.addLazyConstrs()}
\label{\detokenize{javaapi/Model:model-addlazyconstrs}}\begin{quote}

\sphinxAtStartPar
Add lazy constraints to model.

\sphinxAtStartPar
\sphinxstylestrong{Synopsis}
\begin{quote}

\sphinxAtStartPar
\sphinxcode{\sphinxupquote{void addLazyConstrs(ConstrBuilderArray builders, String prefix)}}
\end{quote}

\sphinxAtStartPar
\sphinxstylestrong{Arguments}
\begin{quote}

\sphinxAtStartPar
\sphinxcode{\sphinxupquote{builders}}: array of builders for lazy contraints.

\sphinxAtStartPar
\sphinxcode{\sphinxupquote{prefix}}: name prefix of new lazy constraints.
\end{quote}
\end{quote}

\subsubsection{Model.addLmiConstr()}
\label{\detokenize{javaapi/Model:model-addlmiconstr}}\begin{quote}

\sphinxAtStartPar
Add an LMI constraint to model.

\sphinxAtStartPar
\sphinxstylestrong{Synopsis}
\begin{quote}

\sphinxAtStartPar
\sphinxcode{\sphinxupquote{LmiConstraint addLmiConstr(LmiExpr expr, String name)}}
\end{quote}

\sphinxAtStartPar
\sphinxstylestrong{Arguments}
\begin{quote}

\sphinxAtStartPar
\sphinxcode{\sphinxupquote{expr}}: LMI expression for new LMI contraint.

\sphinxAtStartPar
\sphinxcode{\sphinxupquote{name}}: name of new LMI constraint.
\end{quote}

\sphinxAtStartPar
\sphinxstylestrong{Return}
\begin{quote}

\sphinxAtStartPar
new LMI constraint object.
\end{quote}
\end{quote}

\subsubsection{Model.addNlConstr()}
\label{\detokenize{javaapi/Model:model-addnlconstr}}\begin{quote}

\sphinxAtStartPar
Add a nonlinear constraint to model.

\sphinxAtStartPar
\sphinxstylestrong{Synopsis}
\begin{quote}

\sphinxAtStartPar
\sphinxcode{\sphinxupquote{NlConstraint addNlConstr(}}
\begin{quote}

\sphinxAtStartPar
\sphinxcode{\sphinxupquote{NlExpr expr,}}

\sphinxAtStartPar
\sphinxcode{\sphinxupquote{char sense,}}

\sphinxAtStartPar
\sphinxcode{\sphinxupquote{double rhs,}}

\sphinxAtStartPar
\sphinxcode{\sphinxupquote{String name)}}
\end{quote}
\end{quote}

\sphinxAtStartPar
\sphinxstylestrong{Arguments}
\begin{quote}

\sphinxAtStartPar
\sphinxcode{\sphinxupquote{expr}}: non\sphinxhyphen{}expression for the new contraint.

\sphinxAtStartPar
\sphinxcode{\sphinxupquote{sense}}: sense for new nonlinear constraint, other than range sense.

\sphinxAtStartPar
\sphinxcode{\sphinxupquote{rhs}}: right hand side value for the new constraint.

\sphinxAtStartPar
\sphinxcode{\sphinxupquote{name}}: name of new nonlinear constraint.
\end{quote}

\sphinxAtStartPar
\sphinxstylestrong{Return}
\begin{quote}

\sphinxAtStartPar
new nonlinear constraint object.
\end{quote}
\end{quote}

\subsubsection{Model.addNlConstr()}
\label{\detokenize{javaapi/Model:id23}}\begin{quote}

\sphinxAtStartPar
Add a nonlinear constraint to model.

\sphinxAtStartPar
\sphinxstylestrong{Synopsis}
\begin{quote}

\sphinxAtStartPar
\sphinxcode{\sphinxupquote{NlConstraint addNlConstr(}}
\begin{quote}

\sphinxAtStartPar
\sphinxcode{\sphinxupquote{NlExpr lhs,}}

\sphinxAtStartPar
\sphinxcode{\sphinxupquote{char sense,}}

\sphinxAtStartPar
\sphinxcode{\sphinxupquote{NlExpr rhs,}}

\sphinxAtStartPar
\sphinxcode{\sphinxupquote{String name)}}
\end{quote}
\end{quote}

\sphinxAtStartPar
\sphinxstylestrong{Arguments}
\begin{quote}

\sphinxAtStartPar
\sphinxcode{\sphinxupquote{lhs}}: left hand side nonlinear expression for the new constraint.

\sphinxAtStartPar
\sphinxcode{\sphinxupquote{sense}}: sense for new nonlinear constraint, other than range sense.

\sphinxAtStartPar
\sphinxcode{\sphinxupquote{rhs}}: right hand side nonlinear expression for the new constraint.

\sphinxAtStartPar
\sphinxcode{\sphinxupquote{name}}: name of new nonlinear constraint.
\end{quote}

\sphinxAtStartPar
\sphinxstylestrong{Return}
\begin{quote}

\sphinxAtStartPar
new nonlinear constraint object.
\end{quote}
\end{quote}

\subsubsection{Model.addNlConstr()}
\label{\detokenize{javaapi/Model:id24}}\begin{quote}

\sphinxAtStartPar
Add a nonlinear constraint to model.

\sphinxAtStartPar
\sphinxstylestrong{Synopsis}
\begin{quote}

\sphinxAtStartPar
\sphinxcode{\sphinxupquote{NlConstraint addNlConstr(}}
\begin{quote}

\sphinxAtStartPar
\sphinxcode{\sphinxupquote{NlExpr expr,}}

\sphinxAtStartPar
\sphinxcode{\sphinxupquote{double lb,}}

\sphinxAtStartPar
\sphinxcode{\sphinxupquote{double ub,}}

\sphinxAtStartPar
\sphinxcode{\sphinxupquote{String name)}}
\end{quote}
\end{quote}

\sphinxAtStartPar
\sphinxstylestrong{Arguments}
\begin{quote}

\sphinxAtStartPar
\sphinxcode{\sphinxupquote{expr}}: nonlinear expression for the new constraint.

\sphinxAtStartPar
\sphinxcode{\sphinxupquote{lb}}: lower bound for the new nonlinear constraint.

\sphinxAtStartPar
\sphinxcode{\sphinxupquote{ub}}: upper bound for the new nonlinear constraint

\sphinxAtStartPar
\sphinxcode{\sphinxupquote{name}}: name of new constraint.
\end{quote}

\sphinxAtStartPar
\sphinxstylestrong{Return}
\begin{quote}

\sphinxAtStartPar
new nonlinear constraint object.
\end{quote}
\end{quote}

\subsubsection{Model.addNlConstr()}
\label{\detokenize{javaapi/Model:id25}}\begin{quote}

\sphinxAtStartPar
Add a nonlinear constraint to a model.

\sphinxAtStartPar
\sphinxstylestrong{Synopsis}
\begin{quote}

\sphinxAtStartPar
\sphinxcode{\sphinxupquote{NlConstraint addNlConstr(NlConstrBuilder builder, String name)}}
\end{quote}

\sphinxAtStartPar
\sphinxstylestrong{Arguments}
\begin{quote}

\sphinxAtStartPar
\sphinxcode{\sphinxupquote{builder}}: builder for the new nonlinear constraint.

\sphinxAtStartPar
\sphinxcode{\sphinxupquote{name}}: name of new nonlinear constraint.
\end{quote}

\sphinxAtStartPar
\sphinxstylestrong{Return}
\begin{quote}

\sphinxAtStartPar
new nonlinear constraint object.
\end{quote}
\end{quote}

\subsubsection{Model.addNlConstrs()}
\label{\detokenize{javaapi/Model:model-addnlconstrs}}\begin{quote}

\sphinxAtStartPar
Add nonlinear constraints to a model.

\sphinxAtStartPar
\sphinxstylestrong{Synopsis}
\begin{quote}

\sphinxAtStartPar
\sphinxcode{\sphinxupquote{NlConstrArray addNlConstrs(NlConstrBuilderArray builders, String prefix)}}
\end{quote}

\sphinxAtStartPar
\sphinxstylestrong{Arguments}
\begin{quote}

\sphinxAtStartPar
\sphinxcode{\sphinxupquote{builders}}: builders for new nonlinear constraints.

\sphinxAtStartPar
\sphinxcode{\sphinxupquote{prefix}}: name prefix for new constraints.
\end{quote}

\sphinxAtStartPar
\sphinxstylestrong{Return}
\begin{quote}

\sphinxAtStartPar
array of new nonlinear constraint objects.
\end{quote}
\end{quote}

\subsubsection{Model.addOnesMat()}
\label{\detokenize{javaapi/Model:model-addonesmat}}\begin{quote}

\sphinxAtStartPar
Add a dense symmetric matrix of value one to a model.

\sphinxAtStartPar
\sphinxstylestrong{Synopsis}
\begin{quote}

\sphinxAtStartPar
\sphinxcode{\sphinxupquote{SymMatrix addOnesMat(int dim)}}
\end{quote}

\sphinxAtStartPar
\sphinxstylestrong{Arguments}
\begin{quote}

\sphinxAtStartPar
\sphinxcode{\sphinxupquote{dim}}: dimension of dense symmetric matrix.
\end{quote}

\sphinxAtStartPar
\sphinxstylestrong{Return}
\begin{quote}

\sphinxAtStartPar
new symmetric matrix object.
\end{quote}
\end{quote}

\subsubsection{Model.addPsdConstr()}
\label{\detokenize{javaapi/Model:model-addpsdconstr}}\begin{quote}

\sphinxAtStartPar
Add a PSD constraint to model.

\sphinxAtStartPar
\sphinxstylestrong{Synopsis}
\begin{quote}

\sphinxAtStartPar
\sphinxcode{\sphinxupquote{PsdConstraint addPsdConstr(}}
\begin{quote}

\sphinxAtStartPar
\sphinxcode{\sphinxupquote{PsdExpr expr,}}

\sphinxAtStartPar
\sphinxcode{\sphinxupquote{char sense,}}

\sphinxAtStartPar
\sphinxcode{\sphinxupquote{double rhs,}}

\sphinxAtStartPar
\sphinxcode{\sphinxupquote{String name)}}
\end{quote}
\end{quote}

\sphinxAtStartPar
\sphinxstylestrong{Arguments}
\begin{quote}

\sphinxAtStartPar
\sphinxcode{\sphinxupquote{expr}}: PSD expression for new PSD contraint.

\sphinxAtStartPar
\sphinxcode{\sphinxupquote{sense}}: sense for new PSD constraint.

\sphinxAtStartPar
\sphinxcode{\sphinxupquote{rhs}}: double value at right side of the new PSD constraint.

\sphinxAtStartPar
\sphinxcode{\sphinxupquote{name}}: name of new PSD constraint.
\end{quote}

\sphinxAtStartPar
\sphinxstylestrong{Return}
\begin{quote}

\sphinxAtStartPar
new PSD constraint object.
\end{quote}
\end{quote}

\subsubsection{Model.addPsdConstr()}
\label{\detokenize{javaapi/Model:id26}}\begin{quote}

\sphinxAtStartPar
Add a PSD constraint to model.

\sphinxAtStartPar
\sphinxstylestrong{Synopsis}
\begin{quote}

\sphinxAtStartPar
\sphinxcode{\sphinxupquote{PsdConstraint addPsdConstr(}}
\begin{quote}

\sphinxAtStartPar
\sphinxcode{\sphinxupquote{PsdExpr expr,}}

\sphinxAtStartPar
\sphinxcode{\sphinxupquote{double lb,}}

\sphinxAtStartPar
\sphinxcode{\sphinxupquote{double ub,}}

\sphinxAtStartPar
\sphinxcode{\sphinxupquote{String name)}}
\end{quote}
\end{quote}

\sphinxAtStartPar
\sphinxstylestrong{Arguments}
\begin{quote}

\sphinxAtStartPar
\sphinxcode{\sphinxupquote{expr}}: expression for new PSD constraint.

\sphinxAtStartPar
\sphinxcode{\sphinxupquote{lb}}: lower bound for ew PSD constraint.

\sphinxAtStartPar
\sphinxcode{\sphinxupquote{ub}}: upper bound for new PSD constraint

\sphinxAtStartPar
\sphinxcode{\sphinxupquote{name}}: name of new PSD constraint.
\end{quote}

\sphinxAtStartPar
\sphinxstylestrong{Return}
\begin{quote}

\sphinxAtStartPar
new PSD constraint object.
\end{quote}
\end{quote}

\subsubsection{Model.addPsdConstr()}
\label{\detokenize{javaapi/Model:id27}}\begin{quote}

\sphinxAtStartPar
Add a PSD constraint to model.

\sphinxAtStartPar
\sphinxstylestrong{Synopsis}
\begin{quote}

\sphinxAtStartPar
\sphinxcode{\sphinxupquote{PsdConstraint addPsdConstr(}}
\begin{quote}

\sphinxAtStartPar
\sphinxcode{\sphinxupquote{PsdExpr lhs,}}

\sphinxAtStartPar
\sphinxcode{\sphinxupquote{char sense,}}

\sphinxAtStartPar
\sphinxcode{\sphinxupquote{PsdExpr rhs,}}

\sphinxAtStartPar
\sphinxcode{\sphinxupquote{String name)}}
\end{quote}
\end{quote}

\sphinxAtStartPar
\sphinxstylestrong{Arguments}
\begin{quote}

\sphinxAtStartPar
\sphinxcode{\sphinxupquote{lhs}}: PSD expression at left side of new PSD constraint.

\sphinxAtStartPar
\sphinxcode{\sphinxupquote{sense}}: sense for new PSD constraint.

\sphinxAtStartPar
\sphinxcode{\sphinxupquote{rhs}}: PSD expression at right side of new PSD constraint.

\sphinxAtStartPar
\sphinxcode{\sphinxupquote{name}}: name of new PSD constraint.
\end{quote}

\sphinxAtStartPar
\sphinxstylestrong{Return}
\begin{quote}

\sphinxAtStartPar
new PSD constraint object.
\end{quote}
\end{quote}

\subsubsection{Model.addPsdConstr()}
\label{\detokenize{javaapi/Model:id28}}\begin{quote}

\sphinxAtStartPar
Add a PSD constraint to a model.

\sphinxAtStartPar
\sphinxstylestrong{Synopsis}
\begin{quote}

\sphinxAtStartPar
\sphinxcode{\sphinxupquote{PsdConstraint addPsdConstr(PsdConstrBuilder builder, String name)}}
\end{quote}

\sphinxAtStartPar
\sphinxstylestrong{Arguments}
\begin{quote}

\sphinxAtStartPar
\sphinxcode{\sphinxupquote{builder}}: builder for new PSD constraint.

\sphinxAtStartPar
\sphinxcode{\sphinxupquote{name}}: name of new PSD constraint.
\end{quote}

\sphinxAtStartPar
\sphinxstylestrong{Return}
\begin{quote}

\sphinxAtStartPar
new PSD constraint object.
\end{quote}
\end{quote}

\subsubsection{Model.addPsdVar()}
\label{\detokenize{javaapi/Model:model-addpsdvar}}\begin{quote}

\sphinxAtStartPar
Add a new PSD variable to model.

\sphinxAtStartPar
\sphinxstylestrong{Synopsis}
\begin{quote}

\sphinxAtStartPar
\sphinxcode{\sphinxupquote{PsdVar addPsdVar(int dim, String name)}}
\end{quote}

\sphinxAtStartPar
\sphinxstylestrong{Arguments}
\begin{quote}

\sphinxAtStartPar
\sphinxcode{\sphinxupquote{dim}}: dimension of new PSD variable.

\sphinxAtStartPar
\sphinxcode{\sphinxupquote{name}}: name of new PSD variable.
\end{quote}

\sphinxAtStartPar
\sphinxstylestrong{Return}
\begin{quote}

\sphinxAtStartPar
PSD variable object.
\end{quote}
\end{quote}

\subsubsection{Model.addPsdVars()}
\label{\detokenize{javaapi/Model:model-addpsdvars}}\begin{quote}

\sphinxAtStartPar
Add new PSD variables to model.

\sphinxAtStartPar
\sphinxstylestrong{Synopsis}
\begin{quote}

\sphinxAtStartPar
\sphinxcode{\sphinxupquote{PsdVarArray addPsdVars(}}
\begin{quote}

\sphinxAtStartPar
\sphinxcode{\sphinxupquote{int count,}}

\sphinxAtStartPar
\sphinxcode{\sphinxupquote{int{[}{]} dims,}}

\sphinxAtStartPar
\sphinxcode{\sphinxupquote{String prefix)}}
\end{quote}
\end{quote}

\sphinxAtStartPar
\sphinxstylestrong{Arguments}
\begin{quote}

\sphinxAtStartPar
\sphinxcode{\sphinxupquote{count}}: number of new PSD variables.

\sphinxAtStartPar
\sphinxcode{\sphinxupquote{dims}}: array of dimensions of new PSD variables.

\sphinxAtStartPar
\sphinxcode{\sphinxupquote{prefix}}: name prefix of new PSD variables.
\end{quote}

\sphinxAtStartPar
\sphinxstylestrong{Return}
\begin{quote}

\sphinxAtStartPar
array of PSD variable objects.
\end{quote}
\end{quote}

\subsubsection{Model.addQConstr()}
\label{\detokenize{javaapi/Model:model-addqconstr}}\begin{quote}

\sphinxAtStartPar
Add a quadratic constraint to model.

\sphinxAtStartPar
\sphinxstylestrong{Synopsis}
\begin{quote}

\sphinxAtStartPar
\sphinxcode{\sphinxupquote{QConstraint addQConstr(}}
\begin{quote}

\sphinxAtStartPar
\sphinxcode{\sphinxupquote{QuadExpr expr,}}

\sphinxAtStartPar
\sphinxcode{\sphinxupquote{char sense,}}

\sphinxAtStartPar
\sphinxcode{\sphinxupquote{double rhs,}}

\sphinxAtStartPar
\sphinxcode{\sphinxupquote{String name)}}
\end{quote}
\end{quote}

\sphinxAtStartPar
\sphinxstylestrong{Arguments}
\begin{quote}

\sphinxAtStartPar
\sphinxcode{\sphinxupquote{expr}}: quadratic expression for the new contraint.

\sphinxAtStartPar
\sphinxcode{\sphinxupquote{sense}}: sense for new quadratic constraint.

\sphinxAtStartPar
\sphinxcode{\sphinxupquote{rhs}}: double value at right side of the new quadratic constraint.

\sphinxAtStartPar
\sphinxcode{\sphinxupquote{name}}: name of new quadratic constraint.
\end{quote}

\sphinxAtStartPar
\sphinxstylestrong{Return}
\begin{quote}

\sphinxAtStartPar
new quadratic constraint object.
\end{quote}
\end{quote}

\subsubsection{Model.addQConstr()}
\label{\detokenize{javaapi/Model:id29}}\begin{quote}

\sphinxAtStartPar
Add a quadratic constraint to model.

\sphinxAtStartPar
\sphinxstylestrong{Synopsis}
\begin{quote}

\sphinxAtStartPar
\sphinxcode{\sphinxupquote{QConstraint addQConstr(}}
\begin{quote}

\sphinxAtStartPar
\sphinxcode{\sphinxupquote{QuadExpr lhs,}}

\sphinxAtStartPar
\sphinxcode{\sphinxupquote{char sense,}}

\sphinxAtStartPar
\sphinxcode{\sphinxupquote{QuadExpr rhs,}}

\sphinxAtStartPar
\sphinxcode{\sphinxupquote{String name)}}
\end{quote}
\end{quote}

\sphinxAtStartPar
\sphinxstylestrong{Arguments}
\begin{quote}

\sphinxAtStartPar
\sphinxcode{\sphinxupquote{lhs}}: quadratic expression at left side of new quadratic constraint.

\sphinxAtStartPar
\sphinxcode{\sphinxupquote{sense}}: sense for new quadratic constraint.

\sphinxAtStartPar
\sphinxcode{\sphinxupquote{rhs}}: quadratic expression at right side of new quadratic constraint.

\sphinxAtStartPar
\sphinxcode{\sphinxupquote{name}}: name of new quadratic constraint.
\end{quote}

\sphinxAtStartPar
\sphinxstylestrong{Return}
\begin{quote}

\sphinxAtStartPar
new quadratic constraint object.
\end{quote}
\end{quote}

\subsubsection{Model.addQConstr()}
\label{\detokenize{javaapi/Model:id30}}\begin{quote}

\sphinxAtStartPar
Add a quadratic constraint to a model.

\sphinxAtStartPar
\sphinxstylestrong{Synopsis}
\begin{quote}

\sphinxAtStartPar
\sphinxcode{\sphinxupquote{QConstraint addQConstr(QConstrBuilder builder, String name)}}
\end{quote}

\sphinxAtStartPar
\sphinxstylestrong{Arguments}
\begin{quote}

\sphinxAtStartPar
\sphinxcode{\sphinxupquote{builder}}: builder for the new quadratic constraint.

\sphinxAtStartPar
\sphinxcode{\sphinxupquote{name}}: name of new quadratic constraint.
\end{quote}

\sphinxAtStartPar
\sphinxstylestrong{Return}
\begin{quote}

\sphinxAtStartPar
new quadratic constraint object.
\end{quote}
\end{quote}

\subsubsection{Model.addSos()}
\label{\detokenize{javaapi/Model:model-addsos}}\begin{quote}

\sphinxAtStartPar
Add a SOS constraint to model.

\sphinxAtStartPar
\sphinxstylestrong{Synopsis}
\begin{quote}

\sphinxAtStartPar
\sphinxcode{\sphinxupquote{Sos addSos(SosBuilder builder)}}
\end{quote}

\sphinxAtStartPar
\sphinxstylestrong{Arguments}
\begin{quote}

\sphinxAtStartPar
\sphinxcode{\sphinxupquote{builder}}: builder for new SOS constraint.
\end{quote}

\sphinxAtStartPar
\sphinxstylestrong{Return}
\begin{quote}

\sphinxAtStartPar
new SOS constraint object.
\end{quote}
\end{quote}

\subsubsection{Model.addSos()}
\label{\detokenize{javaapi/Model:id31}}\begin{quote}

\sphinxAtStartPar
Add a SOS constraint to model.

\sphinxAtStartPar
\sphinxstylestrong{Synopsis}
\begin{quote}

\sphinxAtStartPar
\sphinxcode{\sphinxupquote{Sos addSos(}}
\begin{quote}

\sphinxAtStartPar
\sphinxcode{\sphinxupquote{Var{[}{]} vars,}}

\sphinxAtStartPar
\sphinxcode{\sphinxupquote{double{[}{]} weights,}}

\sphinxAtStartPar
\sphinxcode{\sphinxupquote{int type)}}
\end{quote}
\end{quote}

\sphinxAtStartPar
\sphinxstylestrong{Arguments}
\begin{quote}

\sphinxAtStartPar
\sphinxcode{\sphinxupquote{vars}}: variables that participate in the SOS constraint.

\sphinxAtStartPar
\sphinxcode{\sphinxupquote{weights}}: weights for variables in the SOS constraint.

\sphinxAtStartPar
\sphinxcode{\sphinxupquote{type}}: type of SOS constraint.
\end{quote}

\sphinxAtStartPar
\sphinxstylestrong{Return}
\begin{quote}

\sphinxAtStartPar
new SOS constraint object.
\end{quote}
\end{quote}

\subsubsection{Model.addSos()}
\label{\detokenize{javaapi/Model:id32}}\begin{quote}

\sphinxAtStartPar
Add a SOS constraint to model.

\sphinxAtStartPar
\sphinxstylestrong{Synopsis}
\begin{quote}

\sphinxAtStartPar
\sphinxcode{\sphinxupquote{Sos addSos(}}
\begin{quote}

\sphinxAtStartPar
\sphinxcode{\sphinxupquote{VarArray vars,}}

\sphinxAtStartPar
\sphinxcode{\sphinxupquote{double{[}{]} weights,}}

\sphinxAtStartPar
\sphinxcode{\sphinxupquote{int type)}}
\end{quote}
\end{quote}

\sphinxAtStartPar
\sphinxstylestrong{Arguments}
\begin{quote}

\sphinxAtStartPar
\sphinxcode{\sphinxupquote{vars}}: variables that participate in the SOS constraint.

\sphinxAtStartPar
\sphinxcode{\sphinxupquote{weights}}: weights for variables in the SOS constraint.

\sphinxAtStartPar
\sphinxcode{\sphinxupquote{type}}: type of SOS constraint.
\end{quote}

\sphinxAtStartPar
\sphinxstylestrong{Return}
\begin{quote}

\sphinxAtStartPar
new SOS constraint object.
\end{quote}
\end{quote}

\subsubsection{Model.addSparseMat()}
\label{\detokenize{javaapi/Model:model-addsparsemat}}\begin{quote}

\sphinxAtStartPar
Add a sparse symmetric matrix to a model.

\sphinxAtStartPar
\sphinxstylestrong{Synopsis}
\begin{quote}

\sphinxAtStartPar
\sphinxcode{\sphinxupquote{SymMatrix addSparseMat(}}
\begin{quote}

\sphinxAtStartPar
\sphinxcode{\sphinxupquote{int dim,}}

\sphinxAtStartPar
\sphinxcode{\sphinxupquote{int nElems,}}

\sphinxAtStartPar
\sphinxcode{\sphinxupquote{int{[}{]} rows,}}

\sphinxAtStartPar
\sphinxcode{\sphinxupquote{int{[}{]} cols,}}

\sphinxAtStartPar
\sphinxcode{\sphinxupquote{double{[}{]} vals)}}
\end{quote}
\end{quote}

\sphinxAtStartPar
\sphinxstylestrong{Arguments}
\begin{quote}

\sphinxAtStartPar
\sphinxcode{\sphinxupquote{dim}}: dimension of the sparse symmetric matrix.

\sphinxAtStartPar
\sphinxcode{\sphinxupquote{nElems}}: number of non\sphinxhyphen{}zero elements in the sparse symmetric matrix.

\sphinxAtStartPar
\sphinxcode{\sphinxupquote{rows}}: array of row indexes of non\sphinxhyphen{}zero elements.

\sphinxAtStartPar
\sphinxcode{\sphinxupquote{cols}}: array of col indexes of non\sphinxhyphen{}zero elements.

\sphinxAtStartPar
\sphinxcode{\sphinxupquote{vals}}: array of values of non\sphinxhyphen{}zero elements.
\end{quote}

\sphinxAtStartPar
\sphinxstylestrong{Return}
\begin{quote}

\sphinxAtStartPar
new symmetric matrix object.
\end{quote}
\end{quote}

\subsubsection{Model.addSymMat()}
\label{\detokenize{javaapi/Model:model-addsymmat}}\begin{quote}

\sphinxAtStartPar
Given a symmetric matrix expression, add results matrix to model.

\sphinxAtStartPar
\sphinxstylestrong{Synopsis}
\begin{quote}

\sphinxAtStartPar
\sphinxcode{\sphinxupquote{SymMatrix addSymMat(SymMatExpr expr)}}
\end{quote}

\sphinxAtStartPar
\sphinxstylestrong{Arguments}
\begin{quote}

\sphinxAtStartPar
\sphinxcode{\sphinxupquote{expr}}: symmetric matrix expression object.
\end{quote}

\sphinxAtStartPar
\sphinxstylestrong{Return}
\begin{quote}

\sphinxAtStartPar
results symmetric matrix object.
\end{quote}
\end{quote}

\subsubsection{Model.addUserCut()}
\label{\detokenize{javaapi/Model:model-addusercut}}\begin{quote}

\sphinxAtStartPar
Add a user cut to model.

\sphinxAtStartPar
\sphinxstylestrong{Synopsis}
\begin{quote}

\sphinxAtStartPar
\sphinxcode{\sphinxupquote{void addUserCut(}}
\begin{quote}

\sphinxAtStartPar
\sphinxcode{\sphinxupquote{Expr lhs,}}

\sphinxAtStartPar
\sphinxcode{\sphinxupquote{char sense,}}

\sphinxAtStartPar
\sphinxcode{\sphinxupquote{double rhs,}}

\sphinxAtStartPar
\sphinxcode{\sphinxupquote{String name)}}
\end{quote}
\end{quote}

\sphinxAtStartPar
\sphinxstylestrong{Arguments}
\begin{quote}

\sphinxAtStartPar
\sphinxcode{\sphinxupquote{lhs}}: expression for user cut.

\sphinxAtStartPar
\sphinxcode{\sphinxupquote{sense}}: sense for user cut.

\sphinxAtStartPar
\sphinxcode{\sphinxupquote{rhs}}: right hand side value for user cut.

\sphinxAtStartPar
\sphinxcode{\sphinxupquote{name}}: name of user cut.
\end{quote}
\end{quote}

\subsubsection{Model.addUserCut()}
\label{\detokenize{javaapi/Model:id33}}\begin{quote}

\sphinxAtStartPar
Add a user cut to model.

\sphinxAtStartPar
\sphinxstylestrong{Synopsis}
\begin{quote}

\sphinxAtStartPar
\sphinxcode{\sphinxupquote{void addUserCut(}}
\begin{quote}

\sphinxAtStartPar
\sphinxcode{\sphinxupquote{Expr lhs,}}

\sphinxAtStartPar
\sphinxcode{\sphinxupquote{char sense,}}

\sphinxAtStartPar
\sphinxcode{\sphinxupquote{Expr rhs,}}

\sphinxAtStartPar
\sphinxcode{\sphinxupquote{String name)}}
\end{quote}
\end{quote}

\sphinxAtStartPar
\sphinxstylestrong{Arguments}
\begin{quote}

\sphinxAtStartPar
\sphinxcode{\sphinxupquote{lhs}}: left hand side expression for user cut.

\sphinxAtStartPar
\sphinxcode{\sphinxupquote{sense}}: sense for user cut.

\sphinxAtStartPar
\sphinxcode{\sphinxupquote{rhs}}: right hand side expression for user cut.

\sphinxAtStartPar
\sphinxcode{\sphinxupquote{name}}: name of user cut.
\end{quote}
\end{quote}

\subsubsection{Model.addUserCut()}
\label{\detokenize{javaapi/Model:id34}}\begin{quote}

\sphinxAtStartPar
Add a user cut to model.

\sphinxAtStartPar
\sphinxstylestrong{Synopsis}
\begin{quote}

\sphinxAtStartPar
\sphinxcode{\sphinxupquote{void addUserCut(ConstrBuilder builder, String name)}}
\end{quote}

\sphinxAtStartPar
\sphinxstylestrong{Arguments}
\begin{quote}

\sphinxAtStartPar
\sphinxcode{\sphinxupquote{builder}}: builder for user cut.

\sphinxAtStartPar
\sphinxcode{\sphinxupquote{name}}: name of user cut.
\end{quote}
\end{quote}

\subsubsection{Model.addUserCuts()}
\label{\detokenize{javaapi/Model:model-addusercuts}}\begin{quote}

\sphinxAtStartPar
Add user cuts to model.

\sphinxAtStartPar
\sphinxstylestrong{Synopsis}
\begin{quote}

\sphinxAtStartPar
\sphinxcode{\sphinxupquote{void addUserCuts(ConstrBuilderArray builders, String prefix)}}
\end{quote}

\sphinxAtStartPar
\sphinxstylestrong{Arguments}
\begin{quote}

\sphinxAtStartPar
\sphinxcode{\sphinxupquote{builders}}: array of builders for user cuts.

\sphinxAtStartPar
\sphinxcode{\sphinxupquote{prefix}}: name prefix of new user cuts.
\end{quote}
\end{quote}

\subsubsection{Model.addVar()}
\label{\detokenize{javaapi/Model:model-addvar}}\begin{quote}

\sphinxAtStartPar
Add a variable and the associated non\sphinxhyphen{}zero coefficients as column.

\sphinxAtStartPar
\sphinxstylestrong{Synopsis}
\begin{quote}

\sphinxAtStartPar
\sphinxcode{\sphinxupquote{Var addVar(}}
\begin{quote}

\sphinxAtStartPar
\sphinxcode{\sphinxupquote{double lb,}}

\sphinxAtStartPar
\sphinxcode{\sphinxupquote{double ub,}}

\sphinxAtStartPar
\sphinxcode{\sphinxupquote{double obj,}}

\sphinxAtStartPar
\sphinxcode{\sphinxupquote{char vtype,}}

\sphinxAtStartPar
\sphinxcode{\sphinxupquote{String name)}}
\end{quote}
\end{quote}

\sphinxAtStartPar
\sphinxstylestrong{Arguments}
\begin{quote}

\sphinxAtStartPar
\sphinxcode{\sphinxupquote{lb}}: lower bound for new variable.

\sphinxAtStartPar
\sphinxcode{\sphinxupquote{ub}}: upper bound for new variable.

\sphinxAtStartPar
\sphinxcode{\sphinxupquote{obj}}: objective coefficient for new variable.

\sphinxAtStartPar
\sphinxcode{\sphinxupquote{vtype}}: variable type for new variable.

\sphinxAtStartPar
\sphinxcode{\sphinxupquote{name}}: name for new variable.
\end{quote}

\sphinxAtStartPar
\sphinxstylestrong{Return}
\begin{quote}

\sphinxAtStartPar
new variable object.
\end{quote}
\end{quote}

\subsubsection{Model.addVar()}
\label{\detokenize{javaapi/Model:id35}}\begin{quote}

\sphinxAtStartPar
Add a variable and the associated non\sphinxhyphen{}zero coefficients as column.

\sphinxAtStartPar
\sphinxstylestrong{Synopsis}
\begin{quote}

\sphinxAtStartPar
\sphinxcode{\sphinxupquote{Var addVar(}}
\begin{quote}

\sphinxAtStartPar
\sphinxcode{\sphinxupquote{double lb,}}

\sphinxAtStartPar
\sphinxcode{\sphinxupquote{double ub,}}

\sphinxAtStartPar
\sphinxcode{\sphinxupquote{double obj,}}

\sphinxAtStartPar
\sphinxcode{\sphinxupquote{char vtype,}}

\sphinxAtStartPar
\sphinxcode{\sphinxupquote{Column col,}}

\sphinxAtStartPar
\sphinxcode{\sphinxupquote{String name)}}
\end{quote}
\end{quote}

\sphinxAtStartPar
\sphinxstylestrong{Arguments}
\begin{quote}

\sphinxAtStartPar
\sphinxcode{\sphinxupquote{lb}}: lower bound for new variable.

\sphinxAtStartPar
\sphinxcode{\sphinxupquote{ub}}: upper bound for new variable.

\sphinxAtStartPar
\sphinxcode{\sphinxupquote{obj}}: objective coefficient for new variable.

\sphinxAtStartPar
\sphinxcode{\sphinxupquote{vtype}}: variable type for new variable.

\sphinxAtStartPar
\sphinxcode{\sphinxupquote{col}}: column object for specifying a set of constraints to which the variable belongs.

\sphinxAtStartPar
\sphinxcode{\sphinxupquote{name}}: name for new variable.
\end{quote}

\sphinxAtStartPar
\sphinxstylestrong{Return}
\begin{quote}

\sphinxAtStartPar
new variable object.
\end{quote}
\end{quote}

\subsubsection{Model.addVars()}
\label{\detokenize{javaapi/Model:model-addvars}}\begin{quote}

\sphinxAtStartPar
Add new variables to model.

\sphinxAtStartPar
\sphinxstylestrong{Synopsis}
\begin{quote}

\sphinxAtStartPar
\sphinxcode{\sphinxupquote{VarArray addVars(}}
\begin{quote}

\sphinxAtStartPar
\sphinxcode{\sphinxupquote{int count,}}

\sphinxAtStartPar
\sphinxcode{\sphinxupquote{char vtype,}}

\sphinxAtStartPar
\sphinxcode{\sphinxupquote{String prefix)}}
\end{quote}
\end{quote}

\sphinxAtStartPar
\sphinxstylestrong{Arguments}
\begin{quote}

\sphinxAtStartPar
\sphinxcode{\sphinxupquote{count}}: the number of variables to add.

\sphinxAtStartPar
\sphinxcode{\sphinxupquote{vtype}}: variable types for new variables.

\sphinxAtStartPar
\sphinxcode{\sphinxupquote{prefix}}: prefix part for names of new variables.
\end{quote}

\sphinxAtStartPar
\sphinxstylestrong{Return}
\begin{quote}

\sphinxAtStartPar
array of new variable objects.
\end{quote}
\end{quote}

\subsubsection{Model.addVars()}
\label{\detokenize{javaapi/Model:id36}}\begin{quote}

\sphinxAtStartPar
Add new variables to model.

\sphinxAtStartPar
\sphinxstylestrong{Synopsis}
\begin{quote}

\sphinxAtStartPar
\sphinxcode{\sphinxupquote{VarArray addVars(}}
\begin{quote}

\sphinxAtStartPar
\sphinxcode{\sphinxupquote{int count,}}

\sphinxAtStartPar
\sphinxcode{\sphinxupquote{double lb,}}

\sphinxAtStartPar
\sphinxcode{\sphinxupquote{double ub,}}

\sphinxAtStartPar
\sphinxcode{\sphinxupquote{double obj,}}

\sphinxAtStartPar
\sphinxcode{\sphinxupquote{char vtype,}}

\sphinxAtStartPar
\sphinxcode{\sphinxupquote{String prefix)}}
\end{quote}
\end{quote}

\sphinxAtStartPar
\sphinxstylestrong{Arguments}
\begin{quote}

\sphinxAtStartPar
\sphinxcode{\sphinxupquote{count}}: the number of variables to add.

\sphinxAtStartPar
\sphinxcode{\sphinxupquote{lb}}: lower bound for new variables.

\sphinxAtStartPar
\sphinxcode{\sphinxupquote{ub}}: upper bound for new variables.

\sphinxAtStartPar
\sphinxcode{\sphinxupquote{obj}}: objective coefficient for new variables.

\sphinxAtStartPar
\sphinxcode{\sphinxupquote{vtype}}: variable type for new variables.

\sphinxAtStartPar
\sphinxcode{\sphinxupquote{prefix}}: prefix part for names of new variables.
\end{quote}

\sphinxAtStartPar
\sphinxstylestrong{Return}
\begin{quote}

\sphinxAtStartPar
array of new variable objects.
\end{quote}
\end{quote}

\subsubsection{Model.addVars()}
\label{\detokenize{javaapi/Model:id37}}\begin{quote}

\sphinxAtStartPar
Add new variables to model.

\sphinxAtStartPar
\sphinxstylestrong{Synopsis}
\begin{quote}

\sphinxAtStartPar
\sphinxcode{\sphinxupquote{VarArray addVars(}}
\begin{quote}

\sphinxAtStartPar
\sphinxcode{\sphinxupquote{int count,}}

\sphinxAtStartPar
\sphinxcode{\sphinxupquote{double{[}{]} lbs,}}

\sphinxAtStartPar
\sphinxcode{\sphinxupquote{double{[}{]} ubs,}}

\sphinxAtStartPar
\sphinxcode{\sphinxupquote{double{[}{]} objs,}}

\sphinxAtStartPar
\sphinxcode{\sphinxupquote{char{[}{]} types,}}

\sphinxAtStartPar
\sphinxcode{\sphinxupquote{String prefix)}}
\end{quote}
\end{quote}

\sphinxAtStartPar
\sphinxstylestrong{Arguments}
\begin{quote}

\sphinxAtStartPar
\sphinxcode{\sphinxupquote{count}}: the number of variables to add.

\sphinxAtStartPar
\sphinxcode{\sphinxupquote{lbs}}: lower bounds for new variables. if NULL, lower bounds are 0.0.

\sphinxAtStartPar
\sphinxcode{\sphinxupquote{ubs}}: upper bounds for new variables. if NULL, upper bounds are infinity or 1 for binary variables.

\sphinxAtStartPar
\sphinxcode{\sphinxupquote{objs}}: objective coefficients for new variables. if NULL, objective coefficients are 0.0.

\sphinxAtStartPar
\sphinxcode{\sphinxupquote{types}}: variable types for new variables. if NULL, variable types are continuous.

\sphinxAtStartPar
\sphinxcode{\sphinxupquote{prefix}}: prefix part for names of new variables.
\end{quote}

\sphinxAtStartPar
\sphinxstylestrong{Return}
\begin{quote}

\sphinxAtStartPar
array of new variable objects.
\end{quote}
\end{quote}

\subsubsection{Model.addVars()}
\label{\detokenize{javaapi/Model:id38}}\begin{quote}

\sphinxAtStartPar
Add new variables to model.

\sphinxAtStartPar
\sphinxstylestrong{Synopsis}
\begin{quote}

\sphinxAtStartPar
\sphinxcode{\sphinxupquote{VarArray addVars(}}
\begin{quote}

\sphinxAtStartPar
\sphinxcode{\sphinxupquote{double{[}{]} lbs,}}

\sphinxAtStartPar
\sphinxcode{\sphinxupquote{double{[}{]} ubs,}}

\sphinxAtStartPar
\sphinxcode{\sphinxupquote{double{[}{]} objs,}}

\sphinxAtStartPar
\sphinxcode{\sphinxupquote{char{[}{]} types,}}

\sphinxAtStartPar
\sphinxcode{\sphinxupquote{Column{[}{]} cols,}}

\sphinxAtStartPar
\sphinxcode{\sphinxupquote{String prefix)}}
\end{quote}
\end{quote}

\sphinxAtStartPar
\sphinxstylestrong{Arguments}
\begin{quote}

\sphinxAtStartPar
\sphinxcode{\sphinxupquote{lbs}}: lower bounds for new variables. if NULL, lower bounds are 0.0.

\sphinxAtStartPar
\sphinxcode{\sphinxupquote{ubs}}: upper bounds for new variables. if NULL, upper bounds are infinity or 1 for binary variables.

\sphinxAtStartPar
\sphinxcode{\sphinxupquote{objs}}: objective coefficients for new variables. if NULL, objective coefficients are 0.0.

\sphinxAtStartPar
\sphinxcode{\sphinxupquote{types}}: variable types for new variables. if NULL, variable types are continuous.

\sphinxAtStartPar
\sphinxcode{\sphinxupquote{cols}}: column objects for specifying a set of constraints to which each new variable belongs.

\sphinxAtStartPar
\sphinxcode{\sphinxupquote{prefix}}: prefix part for names of new variables.
\end{quote}

\sphinxAtStartPar
\sphinxstylestrong{Return}
\begin{quote}

\sphinxAtStartPar
array of new variable objects.
\end{quote}
\end{quote}

\subsubsection{Model.addVars()}
\label{\detokenize{javaapi/Model:id39}}\begin{quote}

\sphinxAtStartPar
Add new variables to model.

\sphinxAtStartPar
\sphinxstylestrong{Synopsis}
\begin{quote}

\sphinxAtStartPar
\sphinxcode{\sphinxupquote{VarArray addVars(}}
\begin{quote}

\sphinxAtStartPar
\sphinxcode{\sphinxupquote{double{[}{]} lbs,}}

\sphinxAtStartPar
\sphinxcode{\sphinxupquote{double{[}{]} ubs,}}

\sphinxAtStartPar
\sphinxcode{\sphinxupquote{double{[}{]} objs,}}

\sphinxAtStartPar
\sphinxcode{\sphinxupquote{char{[}{]} types,}}

\sphinxAtStartPar
\sphinxcode{\sphinxupquote{ColumnArray cols,}}

\sphinxAtStartPar
\sphinxcode{\sphinxupquote{String prefix)}}
\end{quote}
\end{quote}

\sphinxAtStartPar
\sphinxstylestrong{Arguments}
\begin{quote}

\sphinxAtStartPar
\sphinxcode{\sphinxupquote{lbs}}: lower bounds for new variables. if NULL, lower bounds are 0.0.

\sphinxAtStartPar
\sphinxcode{\sphinxupquote{ubs}}: upper bounds for new variables. if NULL, upper bounds are infinity or 1 for binary variables.

\sphinxAtStartPar
\sphinxcode{\sphinxupquote{objs}}: objective coefficients for new variables. if NULL, objective coefficients are 0.0.

\sphinxAtStartPar
\sphinxcode{\sphinxupquote{types}}: variable types for new variables. if NULL, variable types are continuous.

\sphinxAtStartPar
\sphinxcode{\sphinxupquote{cols}}: columnarray for specifying a set of constraints to which each new variable belongs.

\sphinxAtStartPar
\sphinxcode{\sphinxupquote{prefix}}: prefix part for names of new variables.
\end{quote}

\sphinxAtStartPar
\sphinxstylestrong{Return}
\begin{quote}

\sphinxAtStartPar
array of new variable objects.
\end{quote}
\end{quote}

\subsubsection{Model.clear()}
\label{\detokenize{javaapi/Model:model-clear}}\begin{quote}

\sphinxAtStartPar
Clear all settings including problem itself.

\sphinxAtStartPar
\sphinxstylestrong{Synopsis}
\begin{quote}

\sphinxAtStartPar
\sphinxcode{\sphinxupquote{void clear()}}
\end{quote}
\end{quote}

\subsubsection{Model.clone()}
\label{\detokenize{javaapi/Model:model-clone}}\begin{quote}

\sphinxAtStartPar
Deep copy COPT model.

\sphinxAtStartPar
\sphinxstylestrong{Synopsis}
\begin{quote}

\sphinxAtStartPar
\sphinxcode{\sphinxupquote{Model clone()}}
\end{quote}

\sphinxAtStartPar
\sphinxstylestrong{Return}
\begin{quote}

\sphinxAtStartPar
cloned model object.
\end{quote}
\end{quote}

\subsubsection{Model.computeIIS()}
\label{\detokenize{javaapi/Model:model-computeiis}}\begin{quote}

\sphinxAtStartPar
Compute IIS for infeasible model.

\sphinxAtStartPar
\sphinxstylestrong{Synopsis}
\begin{quote}

\sphinxAtStartPar
\sphinxcode{\sphinxupquote{void computeIIS()}}
\end{quote}
\end{quote}

\subsubsection{Model.delNlObj()}
\label{\detokenize{javaapi/Model:model-delnlobj}}\begin{quote}

\sphinxAtStartPar
delete nonlinear part of objective in model.

\sphinxAtStartPar
\sphinxstylestrong{Synopsis}
\begin{quote}

\sphinxAtStartPar
\sphinxcode{\sphinxupquote{void delNlObj()}}
\end{quote}
\end{quote}

\subsubsection{Model.delObjN()}
\label{\detokenize{javaapi/Model:model-delobjn}}\begin{quote}

\sphinxAtStartPar
Delete linear part of specific multi\sphinxhyphen{}objective function in model.

\sphinxAtStartPar
\sphinxstylestrong{Synopsis}
\begin{quote}

\sphinxAtStartPar
\sphinxcode{\sphinxupquote{void delObjN(int idx)}}
\end{quote}

\sphinxAtStartPar
\sphinxstylestrong{Arguments}
\begin{quote}

\sphinxAtStartPar
\sphinxcode{\sphinxupquote{idx}}: index of a multi\sphinxhyphen{}objective function.
\end{quote}
\end{quote}

\subsubsection{Model.delPsdObj()}
\label{\detokenize{javaapi/Model:model-delpsdobj}}\begin{quote}

\sphinxAtStartPar
delete PSD part of objective in model.

\sphinxAtStartPar
\sphinxstylestrong{Synopsis}
\begin{quote}

\sphinxAtStartPar
\sphinxcode{\sphinxupquote{void delPsdObj()}}
\end{quote}
\end{quote}

\subsubsection{Model.delQuadObj()}
\label{\detokenize{javaapi/Model:model-delquadobj}}\begin{quote}

\sphinxAtStartPar
delete quadratic part of objective in model.

\sphinxAtStartPar
\sphinxstylestrong{Synopsis}
\begin{quote}

\sphinxAtStartPar
\sphinxcode{\sphinxupquote{void delQuadObj()}}
\end{quote}
\end{quote}

\subsubsection{Model.feasRelax()}
\label{\detokenize{javaapi/Model:model-feasrelax}}\begin{quote}

\sphinxAtStartPar
Compute feasibility relaxation for infeasible model.

\sphinxAtStartPar
\sphinxstylestrong{Synopsis}
\begin{quote}

\sphinxAtStartPar
\sphinxcode{\sphinxupquote{void feasRelax(}}
\begin{quote}

\sphinxAtStartPar
\sphinxcode{\sphinxupquote{VarArray vars,}}

\sphinxAtStartPar
\sphinxcode{\sphinxupquote{double{[}{]} colLowPen,}}

\sphinxAtStartPar
\sphinxcode{\sphinxupquote{double{[}{]} colUppPen,}}

\sphinxAtStartPar
\sphinxcode{\sphinxupquote{ConstrArray cons,}}

\sphinxAtStartPar
\sphinxcode{\sphinxupquote{double{[}{]} rowBndPen,}}

\sphinxAtStartPar
\sphinxcode{\sphinxupquote{double{[}{]} rowUppPen)}}
\end{quote}
\end{quote}

\sphinxAtStartPar
\sphinxstylestrong{Arguments}
\begin{quote}

\sphinxAtStartPar
\sphinxcode{\sphinxupquote{vars}}: an array of variables.

\sphinxAtStartPar
\sphinxcode{\sphinxupquote{colLowPen}}: penalties for lower bounds of variables.

\sphinxAtStartPar
\sphinxcode{\sphinxupquote{colUppPen}}: penalties for upper bounds of variables.

\sphinxAtStartPar
\sphinxcode{\sphinxupquote{cons}}: an array of constraints.

\sphinxAtStartPar
\sphinxcode{\sphinxupquote{rowBndPen}}: penalties for right hand sides of constraints.

\sphinxAtStartPar
\sphinxcode{\sphinxupquote{rowUppPen}}: penalties for upper bounds of range constraints.
\end{quote}
\end{quote}

\subsubsection{Model.feasRelax()}
\label{\detokenize{javaapi/Model:id40}}\begin{quote}

\sphinxAtStartPar
Compute feasibility relaxation for infeasible model.

\sphinxAtStartPar
\sphinxstylestrong{Synopsis}
\begin{quote}

\sphinxAtStartPar
\sphinxcode{\sphinxupquote{void feasRelax(int ifRelaxVars, int ifRelaxCons)}}
\end{quote}

\sphinxAtStartPar
\sphinxstylestrong{Arguments}
\begin{quote}

\sphinxAtStartPar
\sphinxcode{\sphinxupquote{ifRelaxVars}}: whether to relax variables.

\sphinxAtStartPar
\sphinxcode{\sphinxupquote{ifRelaxCons}}: whether to relax constraints.
\end{quote}
\end{quote}

\subsubsection{Model.get()}
\label{\detokenize{javaapi/Model:model-get}}\begin{quote}

\sphinxAtStartPar
Query values of information associated with variables.

\sphinxAtStartPar
\sphinxstylestrong{Synopsis}
\begin{quote}

\sphinxAtStartPar
\sphinxcode{\sphinxupquote{double{[}{]} get(String name, Var{[}{]} vars)}}
\end{quote}

\sphinxAtStartPar
\sphinxstylestrong{Arguments}
\begin{quote}

\sphinxAtStartPar
\sphinxcode{\sphinxupquote{name}}: name of information.

\sphinxAtStartPar
\sphinxcode{\sphinxupquote{vars}}: a list of interested variables.
\end{quote}

\sphinxAtStartPar
\sphinxstylestrong{Return}
\begin{quote}

\sphinxAtStartPar
values of information.
\end{quote}
\end{quote}

\subsubsection{Model.get()}
\label{\detokenize{javaapi/Model:id41}}\begin{quote}

\sphinxAtStartPar
Query values of information associated with variables.

\sphinxAtStartPar
\sphinxstylestrong{Synopsis}
\begin{quote}

\sphinxAtStartPar
\sphinxcode{\sphinxupquote{double{[}{]} get(String name, VarArray vars)}}
\end{quote}

\sphinxAtStartPar
\sphinxstylestrong{Arguments}
\begin{quote}

\sphinxAtStartPar
\sphinxcode{\sphinxupquote{name}}: name of information.

\sphinxAtStartPar
\sphinxcode{\sphinxupquote{vars}}: array of interested variables.
\end{quote}

\sphinxAtStartPar
\sphinxstylestrong{Return}
\begin{quote}

\sphinxAtStartPar
values of information.
\end{quote}
\end{quote}

\subsubsection{Model.get()}
\label{\detokenize{javaapi/Model:id42}}\begin{quote}

\sphinxAtStartPar
Query values of information associated with constraints.

\sphinxAtStartPar
\sphinxstylestrong{Synopsis}
\begin{quote}

\sphinxAtStartPar
\sphinxcode{\sphinxupquote{double{[}{]} get(String name, Constraint{[}{]} constrs)}}
\end{quote}

\sphinxAtStartPar
\sphinxstylestrong{Arguments}
\begin{quote}

\sphinxAtStartPar
\sphinxcode{\sphinxupquote{name}}: name of information.

\sphinxAtStartPar
\sphinxcode{\sphinxupquote{constrs}}: a list of interested constraints.
\end{quote}

\sphinxAtStartPar
\sphinxstylestrong{Return}
\begin{quote}

\sphinxAtStartPar
values of information.
\end{quote}
\end{quote}

\subsubsection{Model.get()}
\label{\detokenize{javaapi/Model:id43}}\begin{quote}

\sphinxAtStartPar
Query values of information associated with constraints.

\sphinxAtStartPar
\sphinxstylestrong{Synopsis}
\begin{quote}

\sphinxAtStartPar
\sphinxcode{\sphinxupquote{double{[}{]} get(String name, ConstrArray constrs)}}
\end{quote}

\sphinxAtStartPar
\sphinxstylestrong{Arguments}
\begin{quote}

\sphinxAtStartPar
\sphinxcode{\sphinxupquote{name}}: name of information.

\sphinxAtStartPar
\sphinxcode{\sphinxupquote{constrs}}: array of interested constraints.
\end{quote}

\sphinxAtStartPar
\sphinxstylestrong{Return}
\begin{quote}

\sphinxAtStartPar
values of information.
\end{quote}
\end{quote}

\subsubsection{Model.get()}
\label{\detokenize{javaapi/Model:id44}}\begin{quote}

\sphinxAtStartPar
Query values of information associated with nonlinear constraints.

\sphinxAtStartPar
\sphinxstylestrong{Synopsis}
\begin{quote}

\sphinxAtStartPar
\sphinxcode{\sphinxupquote{double{[}{]} get(String name, NlConstraint{[}{]} constrs)}}
\end{quote}

\sphinxAtStartPar
\sphinxstylestrong{Arguments}
\begin{quote}

\sphinxAtStartPar
\sphinxcode{\sphinxupquote{name}}: name of information.

\sphinxAtStartPar
\sphinxcode{\sphinxupquote{constrs}}: array of desired nonlinear constraints.
\end{quote}

\sphinxAtStartPar
\sphinxstylestrong{Return}
\begin{quote}

\sphinxAtStartPar
output array of information values.
\end{quote}
\end{quote}

\subsubsection{Model.get()}
\label{\detokenize{javaapi/Model:id45}}\begin{quote}

\sphinxAtStartPar
Query values of information associated with nonlinear constraints.

\sphinxAtStartPar
\sphinxstylestrong{Synopsis}
\begin{quote}

\sphinxAtStartPar
\sphinxcode{\sphinxupquote{double{[}{]} get(String name, NlConstrArray constrs)}}
\end{quote}

\sphinxAtStartPar
\sphinxstylestrong{Arguments}
\begin{quote}

\sphinxAtStartPar
\sphinxcode{\sphinxupquote{name}}: name of information.

\sphinxAtStartPar
\sphinxcode{\sphinxupquote{constrs}}: an array object of desired nonlinear constraints.
\end{quote}

\sphinxAtStartPar
\sphinxstylestrong{Return}
\begin{quote}

\sphinxAtStartPar
output array of information values.
\end{quote}
\end{quote}

\subsubsection{Model.get()}
\label{\detokenize{javaapi/Model:id46}}\begin{quote}

\sphinxAtStartPar
Query values of information associated with quadratic constraints.

\sphinxAtStartPar
\sphinxstylestrong{Synopsis}
\begin{quote}

\sphinxAtStartPar
\sphinxcode{\sphinxupquote{double{[}{]} get(String name, QConstraint{[}{]} constrs)}}
\end{quote}

\sphinxAtStartPar
\sphinxstylestrong{Arguments}
\begin{quote}

\sphinxAtStartPar
\sphinxcode{\sphinxupquote{name}}: name of information.

\sphinxAtStartPar
\sphinxcode{\sphinxupquote{constrs}}: a list of interested quadratic constraints.
\end{quote}

\sphinxAtStartPar
\sphinxstylestrong{Return}
\begin{quote}

\sphinxAtStartPar
values of information.
\end{quote}
\end{quote}

\subsubsection{Model.get()}
\label{\detokenize{javaapi/Model:id47}}\begin{quote}

\sphinxAtStartPar
Query values of information associated with quadratic constraints.

\sphinxAtStartPar
\sphinxstylestrong{Synopsis}
\begin{quote}

\sphinxAtStartPar
\sphinxcode{\sphinxupquote{double{[}{]} get(String name, QConstrArray constrs)}}
\end{quote}

\sphinxAtStartPar
\sphinxstylestrong{Arguments}
\begin{quote}

\sphinxAtStartPar
\sphinxcode{\sphinxupquote{name}}: name of information.

\sphinxAtStartPar
\sphinxcode{\sphinxupquote{constrs}}: array of interested quadratic constraints.
\end{quote}

\sphinxAtStartPar
\sphinxstylestrong{Return}
\begin{quote}

\sphinxAtStartPar
values of information.
\end{quote}
\end{quote}

\subsubsection{Model.get()}
\label{\detokenize{javaapi/Model:id48}}\begin{quote}

\sphinxAtStartPar
Query values of information associated with PSD constraints.

\sphinxAtStartPar
\sphinxstylestrong{Synopsis}
\begin{quote}

\sphinxAtStartPar
\sphinxcode{\sphinxupquote{double{[}{]} get(String name, PsdConstraint{[}{]} constrs)}}
\end{quote}

\sphinxAtStartPar
\sphinxstylestrong{Arguments}
\begin{quote}

\sphinxAtStartPar
\sphinxcode{\sphinxupquote{name}}: name of information.

\sphinxAtStartPar
\sphinxcode{\sphinxupquote{constrs}}: a list of desired PSD constraints.
\end{quote}

\sphinxAtStartPar
\sphinxstylestrong{Return}
\begin{quote}

\sphinxAtStartPar
output array of information values.
\end{quote}
\end{quote}

\subsubsection{Model.get()}
\label{\detokenize{javaapi/Model:id49}}\begin{quote}

\sphinxAtStartPar
Query values of information associated with PSD constraints.

\sphinxAtStartPar
\sphinxstylestrong{Synopsis}
\begin{quote}

\sphinxAtStartPar
\sphinxcode{\sphinxupquote{double{[}{]} get(String name, PsdConstrArray constrs)}}
\end{quote}

\sphinxAtStartPar
\sphinxstylestrong{Arguments}
\begin{quote}

\sphinxAtStartPar
\sphinxcode{\sphinxupquote{name}}: name of information.

\sphinxAtStartPar
\sphinxcode{\sphinxupquote{constrs}}: a list of desired PSD constraints.
\end{quote}

\sphinxAtStartPar
\sphinxstylestrong{Return}
\begin{quote}

\sphinxAtStartPar
output array of information values.
\end{quote}
\end{quote}

\subsubsection{Model.getAffineCone()}
\label{\detokenize{javaapi/Model:model-getaffinecone}}\begin{quote}

\sphinxAtStartPar
Get an affine cone constraint of given index in model.

\sphinxAtStartPar
\sphinxstylestrong{Synopsis}
\begin{quote}

\sphinxAtStartPar
\sphinxcode{\sphinxupquote{AffineCone getAffineCone(int idx)}}
\end{quote}

\sphinxAtStartPar
\sphinxstylestrong{Arguments}
\begin{quote}

\sphinxAtStartPar
\sphinxcode{\sphinxupquote{idx}}: index of the desired affine cone constraint.
\end{quote}

\sphinxAtStartPar
\sphinxstylestrong{Return}
\begin{quote}

\sphinxAtStartPar
the desired affine cone constraint object.
\end{quote}
\end{quote}

\subsubsection{Model.GetAffineConeBuilder()}
\label{\detokenize{javaapi/Model:model-getaffineconebuilder}}\begin{quote}

\sphinxAtStartPar
Get builder of given affine cone constraint in model.

\sphinxAtStartPar
\sphinxstylestrong{Synopsis}
\begin{quote}

\sphinxAtStartPar
\sphinxcode{\sphinxupquote{AffineConeBuilder GetAffineConeBuilder(AffineCone cone)}}
\end{quote}

\sphinxAtStartPar
\sphinxstylestrong{Arguments}
\begin{quote}

\sphinxAtStartPar
\sphinxcode{\sphinxupquote{cone}}: affine cone constraint.
\end{quote}

\sphinxAtStartPar
\sphinxstylestrong{Return}
\begin{quote}

\sphinxAtStartPar
desired affine cone constraint builder.
\end{quote}
\end{quote}

\subsubsection{Model.GetAffineConeBuilders()}
\label{\detokenize{javaapi/Model:model-getaffineconebuilders}}\begin{quote}

\sphinxAtStartPar
Get builders of all affine cone constraints in model.

\sphinxAtStartPar
\sphinxstylestrong{Synopsis}
\begin{quote}

\sphinxAtStartPar
\sphinxcode{\sphinxupquote{AffineConeBuilderArray GetAffineConeBuilders()}}
\end{quote}

\sphinxAtStartPar
\sphinxstylestrong{Return}
\begin{quote}

\sphinxAtStartPar
array object of affine cone constraint builders.
\end{quote}
\end{quote}

\subsubsection{Model.getAffineConeBuilders()}
\label{\detokenize{javaapi/Model:id50}}\begin{quote}

\sphinxAtStartPar
Get builders of given affine cone constraints in model.

\sphinxAtStartPar
\sphinxstylestrong{Synopsis}
\begin{quote}

\sphinxAtStartPar
\sphinxcode{\sphinxupquote{AffineConeBuilderArray getAffineConeBuilders(AffineCone{[}{]} cones)}}
\end{quote}

\sphinxAtStartPar
\sphinxstylestrong{Arguments}
\begin{quote}

\sphinxAtStartPar
\sphinxcode{\sphinxupquote{cones}}: array of affine cone constraints.
\end{quote}

\sphinxAtStartPar
\sphinxstylestrong{Return}
\begin{quote}

\sphinxAtStartPar
array object of desired affine cone constraint builders.
\end{quote}
\end{quote}

\subsubsection{Model.getAffineConeBuilders()}
\label{\detokenize{javaapi/Model:id51}}\begin{quote}

\sphinxAtStartPar
Get builders of given affine cone constraints in model.

\sphinxAtStartPar
\sphinxstylestrong{Synopsis}
\begin{quote}

\sphinxAtStartPar
\sphinxcode{\sphinxupquote{AffineConeBuilderArray getAffineConeBuilders(AffineConeArray cones)}}
\end{quote}

\sphinxAtStartPar
\sphinxstylestrong{Arguments}
\begin{quote}

\sphinxAtStartPar
\sphinxcode{\sphinxupquote{cones}}: array of affine cone constraints.
\end{quote}

\sphinxAtStartPar
\sphinxstylestrong{Return}
\begin{quote}

\sphinxAtStartPar
array object of desired affine cone constraint builders.
\end{quote}
\end{quote}

\subsubsection{Model.getAffineConeByName()}
\label{\detokenize{javaapi/Model:model-getaffineconebyname}}\begin{quote}

\sphinxAtStartPar
Get an affine cone constraint of given name in model.

\sphinxAtStartPar
\sphinxstylestrong{Synopsis}
\begin{quote}

\sphinxAtStartPar
\sphinxcode{\sphinxupquote{AffineCone getAffineConeByName(String name)}}
\end{quote}

\sphinxAtStartPar
\sphinxstylestrong{Arguments}
\begin{quote}

\sphinxAtStartPar
\sphinxcode{\sphinxupquote{name}}: name of the desired affine cone constraint.
\end{quote}

\sphinxAtStartPar
\sphinxstylestrong{Return}
\begin{quote}

\sphinxAtStartPar
the desired affine cone constraint object.
\end{quote}
\end{quote}

\subsubsection{Model.getAffineCones()}
\label{\detokenize{javaapi/Model:model-getaffinecones}}\begin{quote}

\sphinxAtStartPar
Get all affine cone constraints in model.

\sphinxAtStartPar
\sphinxstylestrong{Synopsis}
\begin{quote}

\sphinxAtStartPar
\sphinxcode{\sphinxupquote{AffineConeArray getAffineCones()}}
\end{quote}

\sphinxAtStartPar
\sphinxstylestrong{Return}
\begin{quote}

\sphinxAtStartPar
array object of affine cone constraints.
\end{quote}
\end{quote}

\subsubsection{Model.getCoeff()}
\label{\detokenize{javaapi/Model:model-getcoeff}}\begin{quote}

\sphinxAtStartPar
Get the coefficient of variable in linear constraint.

\sphinxAtStartPar
\sphinxstylestrong{Synopsis}
\begin{quote}

\sphinxAtStartPar
\sphinxcode{\sphinxupquote{double getCoeff(Constraint constr, Var var)}}
\end{quote}

\sphinxAtStartPar
\sphinxstylestrong{Arguments}
\begin{quote}

\sphinxAtStartPar
\sphinxcode{\sphinxupquote{constr}}: The requested constraint.

\sphinxAtStartPar
\sphinxcode{\sphinxupquote{var}}: The requested variable.
\end{quote}

\sphinxAtStartPar
\sphinxstylestrong{Return}
\begin{quote}

\sphinxAtStartPar
The requested coefficient.
\end{quote}
\end{quote}

\subsubsection{Model.getCol()}
\label{\detokenize{javaapi/Model:model-getcol}}\begin{quote}

\sphinxAtStartPar
Get a column object that have a list of constraints in which the variable participates.

\sphinxAtStartPar
\sphinxstylestrong{Synopsis}
\begin{quote}

\sphinxAtStartPar
\sphinxcode{\sphinxupquote{Column getCol(Var var)}}
\end{quote}

\sphinxAtStartPar
\sphinxstylestrong{Arguments}
\begin{quote}

\sphinxAtStartPar
\sphinxcode{\sphinxupquote{var}}: a variable object.
\end{quote}

\sphinxAtStartPar
\sphinxstylestrong{Return}
\begin{quote}

\sphinxAtStartPar
a column object associated with a variable.
\end{quote}
\end{quote}

\subsubsection{Model.getColBasis()}
\label{\detokenize{javaapi/Model:model-getcolbasis}}\begin{quote}

\sphinxAtStartPar
Get status of column basis.

\sphinxAtStartPar
\sphinxstylestrong{Synopsis}
\begin{quote}

\sphinxAtStartPar
\sphinxcode{\sphinxupquote{int{[}{]} getColBasis()}}
\end{quote}

\sphinxAtStartPar
\sphinxstylestrong{Return}
\begin{quote}

\sphinxAtStartPar
basis status.
\end{quote}
\end{quote}

\subsubsection{Model.getCone()}
\label{\detokenize{javaapi/Model:model-getcone}}\begin{quote}

\sphinxAtStartPar
Get a cone constraint of given index in model.

\sphinxAtStartPar
\sphinxstylestrong{Synopsis}
\begin{quote}

\sphinxAtStartPar
\sphinxcode{\sphinxupquote{Cone getCone(int idx)}}
\end{quote}

\sphinxAtStartPar
\sphinxstylestrong{Arguments}
\begin{quote}

\sphinxAtStartPar
\sphinxcode{\sphinxupquote{idx}}: index of the desired cone constraint.
\end{quote}

\sphinxAtStartPar
\sphinxstylestrong{Return}
\begin{quote}

\sphinxAtStartPar
the desired cone constraint object.
\end{quote}
\end{quote}

\subsubsection{Model.getConeBuilders()}
\label{\detokenize{javaapi/Model:model-getconebuilders}}\begin{quote}

\sphinxAtStartPar
Get builders of all cone constraints in model.

\sphinxAtStartPar
\sphinxstylestrong{Synopsis}
\begin{quote}

\sphinxAtStartPar
\sphinxcode{\sphinxupquote{ConeBuilderArray getConeBuilders()}}
\end{quote}

\sphinxAtStartPar
\sphinxstylestrong{Return}
\begin{quote}

\sphinxAtStartPar
array object of cone constraint builders.
\end{quote}
\end{quote}

\subsubsection{Model.getConeBuilders()}
\label{\detokenize{javaapi/Model:id52}}\begin{quote}

\sphinxAtStartPar
Get builders of given cone constraints in model.

\sphinxAtStartPar
\sphinxstylestrong{Synopsis}
\begin{quote}

\sphinxAtStartPar
\sphinxcode{\sphinxupquote{ConeBuilderArray getConeBuilders(Cone{[}{]} cones)}}
\end{quote}

\sphinxAtStartPar
\sphinxstylestrong{Arguments}
\begin{quote}

\sphinxAtStartPar
\sphinxcode{\sphinxupquote{cones}}: array of cone constraints.
\end{quote}

\sphinxAtStartPar
\sphinxstylestrong{Return}
\begin{quote}

\sphinxAtStartPar
array object of desired cone constraint builders.
\end{quote}
\end{quote}

\subsubsection{Model.getConeBuilders()}
\label{\detokenize{javaapi/Model:id53}}\begin{quote}

\sphinxAtStartPar
Get builders of given cone constraints in model.

\sphinxAtStartPar
\sphinxstylestrong{Synopsis}
\begin{quote}

\sphinxAtStartPar
\sphinxcode{\sphinxupquote{ConeBuilderArray getConeBuilders(ConeArray cones)}}
\end{quote}

\sphinxAtStartPar
\sphinxstylestrong{Arguments}
\begin{quote}

\sphinxAtStartPar
\sphinxcode{\sphinxupquote{cones}}: array of cone constraints.
\end{quote}

\sphinxAtStartPar
\sphinxstylestrong{Return}
\begin{quote}

\sphinxAtStartPar
array object of desired cone constraint builders.
\end{quote}
\end{quote}

\subsubsection{Model.getCones()}
\label{\detokenize{javaapi/Model:model-getcones}}\begin{quote}

\sphinxAtStartPar
Get all cone constraints in model.

\sphinxAtStartPar
\sphinxstylestrong{Synopsis}
\begin{quote}

\sphinxAtStartPar
\sphinxcode{\sphinxupquote{ConeArray getCones()}}
\end{quote}

\sphinxAtStartPar
\sphinxstylestrong{Return}
\begin{quote}

\sphinxAtStartPar
array object of cone constraints.
\end{quote}
\end{quote}

\subsubsection{Model.getConstr()}
\label{\detokenize{javaapi/Model:model-getconstr}}\begin{quote}

\sphinxAtStartPar
Get a constraint of given index in model.

\sphinxAtStartPar
\sphinxstylestrong{Synopsis}
\begin{quote}

\sphinxAtStartPar
\sphinxcode{\sphinxupquote{Constraint getConstr(int idx)}}
\end{quote}

\sphinxAtStartPar
\sphinxstylestrong{Arguments}
\begin{quote}

\sphinxAtStartPar
\sphinxcode{\sphinxupquote{idx}}: index of the desired constraint.
\end{quote}

\sphinxAtStartPar
\sphinxstylestrong{Return}
\begin{quote}

\sphinxAtStartPar
the desired constraint object.
\end{quote}
\end{quote}

\subsubsection{Model.getConstrBuilder()}
\label{\detokenize{javaapi/Model:model-getconstrbuilder}}\begin{quote}

\sphinxAtStartPar
Get builder of a constraint in model, including variables and associated coefficients, sense and RHS.

\sphinxAtStartPar
\sphinxstylestrong{Synopsis}
\begin{quote}

\sphinxAtStartPar
\sphinxcode{\sphinxupquote{ConstrBuilder getConstrBuilder(Constraint constr)}}
\end{quote}

\sphinxAtStartPar
\sphinxstylestrong{Arguments}
\begin{quote}

\sphinxAtStartPar
\sphinxcode{\sphinxupquote{constr}}: a constraint object.
\end{quote}

\sphinxAtStartPar
\sphinxstylestrong{Return}
\begin{quote}

\sphinxAtStartPar
constraint builder object.
\end{quote}
\end{quote}

\subsubsection{Model.getConstrBuilders()}
\label{\detokenize{javaapi/Model:model-getconstrbuilders}}\begin{quote}

\sphinxAtStartPar
Get builders of all constraints in model.

\sphinxAtStartPar
\sphinxstylestrong{Synopsis}
\begin{quote}

\sphinxAtStartPar
\sphinxcode{\sphinxupquote{ConstrBuilderArray getConstrBuilders()}}
\end{quote}

\sphinxAtStartPar
\sphinxstylestrong{Return}
\begin{quote}

\sphinxAtStartPar
array object of constraint builders.
\end{quote}
\end{quote}

\subsubsection{Model.getConstrByName()}
\label{\detokenize{javaapi/Model:model-getconstrbyname}}\begin{quote}

\sphinxAtStartPar
Get a constraint of given name in model.

\sphinxAtStartPar
\sphinxstylestrong{Synopsis}
\begin{quote}

\sphinxAtStartPar
\sphinxcode{\sphinxupquote{Constraint getConstrByName(String name)}}
\end{quote}

\sphinxAtStartPar
\sphinxstylestrong{Arguments}
\begin{quote}

\sphinxAtStartPar
\sphinxcode{\sphinxupquote{name}}: name of the desired constraint.
\end{quote}

\sphinxAtStartPar
\sphinxstylestrong{Return}
\begin{quote}

\sphinxAtStartPar
the desired constraint object.
\end{quote}
\end{quote}

\subsubsection{Model.getConstrLowerIIS()}
\label{\detokenize{javaapi/Model:model-getconstrloweriis}}\begin{quote}

\sphinxAtStartPar
Get IIS status of lower bounds of constraints.

\sphinxAtStartPar
\sphinxstylestrong{Synopsis}
\begin{quote}

\sphinxAtStartPar
\sphinxcode{\sphinxupquote{int{[}{]} getConstrLowerIIS(ConstrArray constrs)}}
\end{quote}

\sphinxAtStartPar
\sphinxstylestrong{Arguments}
\begin{quote}

\sphinxAtStartPar
\sphinxcode{\sphinxupquote{constrs}}: Array of constraints.
\end{quote}

\sphinxAtStartPar
\sphinxstylestrong{Return}
\begin{quote}

\sphinxAtStartPar
IIS status of lower bounds of constraints.
\end{quote}
\end{quote}

\subsubsection{Model.getConstrLowerIIS()}
\label{\detokenize{javaapi/Model:id54}}\begin{quote}

\sphinxAtStartPar
Get IIS status of lower bounds of constraints.

\sphinxAtStartPar
\sphinxstylestrong{Synopsis}
\begin{quote}

\sphinxAtStartPar
\sphinxcode{\sphinxupquote{int{[}{]} getConstrLowerIIS(Constraint{[}{]} constrs)}}
\end{quote}

\sphinxAtStartPar
\sphinxstylestrong{Arguments}
\begin{quote}

\sphinxAtStartPar
\sphinxcode{\sphinxupquote{constrs}}: Array of constraints.
\end{quote}

\sphinxAtStartPar
\sphinxstylestrong{Return}
\begin{quote}

\sphinxAtStartPar
IIS status of lower bounds of constraints.
\end{quote}
\end{quote}

\subsubsection{Model.getConstrs()}
\label{\detokenize{javaapi/Model:model-getconstrs}}\begin{quote}

\sphinxAtStartPar
Get all constraints in model.

\sphinxAtStartPar
\sphinxstylestrong{Synopsis}
\begin{quote}

\sphinxAtStartPar
\sphinxcode{\sphinxupquote{ConstrArray getConstrs()}}
\end{quote}

\sphinxAtStartPar
\sphinxstylestrong{Return}
\begin{quote}

\sphinxAtStartPar
array object of constraints.
\end{quote}
\end{quote}

\subsubsection{Model.getConstrUpperIIS()}
\label{\detokenize{javaapi/Model:model-getconstrupperiis}}\begin{quote}

\sphinxAtStartPar
Get IIS status of upper bounds of constraints.

\sphinxAtStartPar
\sphinxstylestrong{Synopsis}
\begin{quote}

\sphinxAtStartPar
\sphinxcode{\sphinxupquote{int{[}{]} getConstrUpperIIS(ConstrArray constrs)}}
\end{quote}

\sphinxAtStartPar
\sphinxstylestrong{Arguments}
\begin{quote}

\sphinxAtStartPar
\sphinxcode{\sphinxupquote{constrs}}: Array of constraints.
\end{quote}

\sphinxAtStartPar
\sphinxstylestrong{Return}
\begin{quote}

\sphinxAtStartPar
IIS status of upper bounds of constraints.
\end{quote}
\end{quote}

\subsubsection{Model.getConstrUpperIIS()}
\label{\detokenize{javaapi/Model:id55}}\begin{quote}

\sphinxAtStartPar
Get IIS status of upper bounds of constraints.

\sphinxAtStartPar
\sphinxstylestrong{Synopsis}
\begin{quote}

\sphinxAtStartPar
\sphinxcode{\sphinxupquote{int{[}{]} getConstrUpperIIS(Constraint{[}{]} constrs)}}
\end{quote}

\sphinxAtStartPar
\sphinxstylestrong{Arguments}
\begin{quote}

\sphinxAtStartPar
\sphinxcode{\sphinxupquote{constrs}}: Array of constraints.
\end{quote}

\sphinxAtStartPar
\sphinxstylestrong{Return}
\begin{quote}

\sphinxAtStartPar
IIS status of upper bounds of constraints.
\end{quote}
\end{quote}

\subsubsection{Model.getDblAttr()}
\label{\detokenize{javaapi/Model:model-getdblattr}}\begin{quote}

\sphinxAtStartPar
Get value of a COPT double attribute.

\sphinxAtStartPar
\sphinxstylestrong{Synopsis}
\begin{quote}

\sphinxAtStartPar
\sphinxcode{\sphinxupquote{double getDblAttr(String attr)}}
\end{quote}

\sphinxAtStartPar
\sphinxstylestrong{Arguments}
\begin{quote}

\sphinxAtStartPar
\sphinxcode{\sphinxupquote{attr}}: name of double attribute.
\end{quote}

\sphinxAtStartPar
\sphinxstylestrong{Return}
\begin{quote}

\sphinxAtStartPar
value of double attribute.
\end{quote}
\end{quote}

\subsubsection{Model.getDblAttrN()}
\label{\detokenize{javaapi/Model:model-getdblattrn}}\begin{quote}

\sphinxAtStartPar
Get value of a double attribute of a multi\sphinxhyphen{}objective function.

\sphinxAtStartPar
\sphinxstylestrong{Synopsis}
\begin{quote}

\sphinxAtStartPar
\sphinxcode{\sphinxupquote{double getDblAttrN(int idx, String attr)}}
\end{quote}

\sphinxAtStartPar
\sphinxstylestrong{Arguments}
\begin{quote}

\sphinxAtStartPar
\sphinxcode{\sphinxupquote{idx}}: index of a multi\sphinxhyphen{}objective function.

\sphinxAtStartPar
\sphinxcode{\sphinxupquote{attr}}: name of double attribute.
\end{quote}

\sphinxAtStartPar
\sphinxstylestrong{Return}
\begin{quote}

\sphinxAtStartPar
value of double attribute.
\end{quote}
\end{quote}

\subsubsection{Model.getDblParam()}
\label{\detokenize{javaapi/Model:model-getdblparam}}\begin{quote}

\sphinxAtStartPar
Get value of a COPT double parameter.

\sphinxAtStartPar
\sphinxstylestrong{Synopsis}
\begin{quote}

\sphinxAtStartPar
\sphinxcode{\sphinxupquote{double getDblParam(String param)}}
\end{quote}

\sphinxAtStartPar
\sphinxstylestrong{Arguments}
\begin{quote}

\sphinxAtStartPar
\sphinxcode{\sphinxupquote{param}}: name of double parameter.
\end{quote}

\sphinxAtStartPar
\sphinxstylestrong{Return}
\begin{quote}

\sphinxAtStartPar
value of double parameter.
\end{quote}
\end{quote}

\subsubsection{Model.getDblParamInfo()}
\label{\detokenize{javaapi/Model:model-getdblparaminfo}}\begin{quote}

\sphinxAtStartPar
Get current, default, minimum, maximum of COPT double parameter.

\sphinxAtStartPar
\sphinxstylestrong{Synopsis}
\begin{quote}

\sphinxAtStartPar
\sphinxcode{\sphinxupquote{double{[}{]} getDblParamInfo(String name)}}
\end{quote}

\sphinxAtStartPar
\sphinxstylestrong{Arguments}
\begin{quote}

\sphinxAtStartPar
\sphinxcode{\sphinxupquote{name}}: name of integer parameter.
\end{quote}

\sphinxAtStartPar
\sphinxstylestrong{Return}
\begin{quote}

\sphinxAtStartPar
current, default, minimum, maximum of COPT double parameter.
\end{quote}
\end{quote}

\subsubsection{Model.getDblParamN()}
\label{\detokenize{javaapi/Model:model-getdblparamn}}\begin{quote}

\sphinxAtStartPar
Get value of a double parameter of a multi\sphinxhyphen{}objective function.

\sphinxAtStartPar
\sphinxstylestrong{Synopsis}
\begin{quote}

\sphinxAtStartPar
\sphinxcode{\sphinxupquote{double getDblParamN(int idx, String param)}}
\end{quote}

\sphinxAtStartPar
\sphinxstylestrong{Arguments}
\begin{quote}

\sphinxAtStartPar
\sphinxcode{\sphinxupquote{idx}}: index of a multi\sphinxhyphen{}objective function.

\sphinxAtStartPar
\sphinxcode{\sphinxupquote{param}}: name of double parameter.
\end{quote}

\sphinxAtStartPar
\sphinxstylestrong{Return}
\begin{quote}

\sphinxAtStartPar
value of double parameter.
\end{quote}
\end{quote}

\subsubsection{Model.getExpCone()}
\label{\detokenize{javaapi/Model:model-getexpcone}}\begin{quote}

\sphinxAtStartPar
Get an exponential cone constraint of given index in model.

\sphinxAtStartPar
\sphinxstylestrong{Synopsis}
\begin{quote}

\sphinxAtStartPar
\sphinxcode{\sphinxupquote{ExpCone getExpCone(int idx)}}
\end{quote}

\sphinxAtStartPar
\sphinxstylestrong{Arguments}
\begin{quote}

\sphinxAtStartPar
\sphinxcode{\sphinxupquote{idx}}: index of the desired exponential cone constraint.
\end{quote}

\sphinxAtStartPar
\sphinxstylestrong{Return}
\begin{quote}

\sphinxAtStartPar
the desired exponential cone constraint object.
\end{quote}
\end{quote}

\subsubsection{Model.getExpConeBuilders()}
\label{\detokenize{javaapi/Model:model-getexpconebuilders}}\begin{quote}

\sphinxAtStartPar
Get builders of all exponential cone constraints in model.

\sphinxAtStartPar
\sphinxstylestrong{Synopsis}
\begin{quote}

\sphinxAtStartPar
\sphinxcode{\sphinxupquote{ExpConeBuilderArray getExpConeBuilders()}}
\end{quote}

\sphinxAtStartPar
\sphinxstylestrong{Return}
\begin{quote}

\sphinxAtStartPar
array object of exponential cone constraint builders.
\end{quote}
\end{quote}

\subsubsection{Model.getExpConeBuilders()}
\label{\detokenize{javaapi/Model:id56}}\begin{quote}

\sphinxAtStartPar
Get builders of given exponential cone constraints in model.

\sphinxAtStartPar
\sphinxstylestrong{Synopsis}
\begin{quote}

\sphinxAtStartPar
\sphinxcode{\sphinxupquote{ExpConeBuilderArray getExpConeBuilders(ExpCone{[}{]} cones)}}
\end{quote}

\sphinxAtStartPar
\sphinxstylestrong{Arguments}
\begin{quote}

\sphinxAtStartPar
\sphinxcode{\sphinxupquote{cones}}: array of exponential cone constraints.
\end{quote}

\sphinxAtStartPar
\sphinxstylestrong{Return}
\begin{quote}

\sphinxAtStartPar
array object of desired exponential cone constraint builders.
\end{quote}
\end{quote}

\subsubsection{Model.getExpConeBuilders()}
\label{\detokenize{javaapi/Model:id57}}\begin{quote}

\sphinxAtStartPar
Get builders of given exponential cone constraints in model.

\sphinxAtStartPar
\sphinxstylestrong{Synopsis}
\begin{quote}

\sphinxAtStartPar
\sphinxcode{\sphinxupquote{ExpConeBuilderArray getExpConeBuilders(ExpConeArray cones)}}
\end{quote}

\sphinxAtStartPar
\sphinxstylestrong{Arguments}
\begin{quote}

\sphinxAtStartPar
\sphinxcode{\sphinxupquote{cones}}: array of exponential cone constraints.
\end{quote}

\sphinxAtStartPar
\sphinxstylestrong{Return}
\begin{quote}

\sphinxAtStartPar
array object of desired exponential cone constraint builders.
\end{quote}
\end{quote}

\subsubsection{Model.getExpCones()}
\label{\detokenize{javaapi/Model:model-getexpcones}}\begin{quote}

\sphinxAtStartPar
Get all exponential cone constraints in model.

\sphinxAtStartPar
\sphinxstylestrong{Synopsis}
\begin{quote}

\sphinxAtStartPar
\sphinxcode{\sphinxupquote{ExpConeArray getExpCones()}}
\end{quote}

\sphinxAtStartPar
\sphinxstylestrong{Return}
\begin{quote}

\sphinxAtStartPar
array object of exponential cone constraints.
\end{quote}
\end{quote}

\subsubsection{Model.getGenConstr()}
\label{\detokenize{javaapi/Model:model-getgenconstr}}\begin{quote}

\sphinxAtStartPar
Get a general constraint of given index in model.

\sphinxAtStartPar
\sphinxstylestrong{Synopsis}
\begin{quote}

\sphinxAtStartPar
\sphinxcode{\sphinxupquote{GenConstr getGenConstr(int idx)}}
\end{quote}

\sphinxAtStartPar
\sphinxstylestrong{Arguments}
\begin{quote}

\sphinxAtStartPar
\sphinxcode{\sphinxupquote{idx}}: index of the desired general constraint.
\end{quote}

\sphinxAtStartPar
\sphinxstylestrong{Return}
\begin{quote}

\sphinxAtStartPar
the desired general constraint object.
\end{quote}
\end{quote}

\subsubsection{Model.getGenConstrByName()}
\label{\detokenize{javaapi/Model:model-getgenconstrbyname}}\begin{quote}

\sphinxAtStartPar
Get a general constraint of given name in model.

\sphinxAtStartPar
\sphinxstylestrong{Synopsis}
\begin{quote}

\sphinxAtStartPar
\sphinxcode{\sphinxupquote{GenConstr getGenConstrByName(String name)}}
\end{quote}

\sphinxAtStartPar
\sphinxstylestrong{Arguments}
\begin{quote}

\sphinxAtStartPar
\sphinxcode{\sphinxupquote{name}}: name of the desired general constraint.
\end{quote}

\sphinxAtStartPar
\sphinxstylestrong{Return}
\begin{quote}

\sphinxAtStartPar
the desired general constraint object.
\end{quote}
\end{quote}

\subsubsection{Model.getGenConstrIndicator()}
\label{\detokenize{javaapi/Model:model-getgenconstrindicator}}\begin{quote}

\sphinxAtStartPar
Get builder of given general constraint of type indicator.

\sphinxAtStartPar
\sphinxstylestrong{Synopsis}
\begin{quote}

\sphinxAtStartPar
\sphinxcode{\sphinxupquote{GenConstrBuilder getGenConstrIndicator(GenConstr indicator)}}
\end{quote}

\sphinxAtStartPar
\sphinxstylestrong{Arguments}
\begin{quote}

\sphinxAtStartPar
\sphinxcode{\sphinxupquote{indicator}}: a general constraint of type indicator.
\end{quote}

\sphinxAtStartPar
\sphinxstylestrong{Return}
\begin{quote}

\sphinxAtStartPar
builder object of general constraint of type indicator.
\end{quote}
\end{quote}

\subsubsection{Model.getGenConstrIndicators()}
\label{\detokenize{javaapi/Model:model-getgenconstrindicators}}\begin{quote}

\sphinxAtStartPar
Get builders of all general constraints in model.

\sphinxAtStartPar
\sphinxstylestrong{Synopsis}
\begin{quote}

\sphinxAtStartPar
\sphinxcode{\sphinxupquote{GenConstrBuilderArray getGenConstrIndicators()}}
\end{quote}

\sphinxAtStartPar
\sphinxstylestrong{Return}
\begin{quote}

\sphinxAtStartPar
array object of general constraint builders.
\end{quote}
\end{quote}

\subsubsection{Model.getGenConstrs()}
\label{\detokenize{javaapi/Model:model-getgenconstrs}}\begin{quote}

\sphinxAtStartPar
Get all general constraints in model.

\sphinxAtStartPar
\sphinxstylestrong{Synopsis}
\begin{quote}

\sphinxAtStartPar
\sphinxcode{\sphinxupquote{GenConstrArray getGenConstrs()}}
\end{quote}

\sphinxAtStartPar
\sphinxstylestrong{Return}
\begin{quote}

\sphinxAtStartPar
array object of general constraints.
\end{quote}
\end{quote}

\subsubsection{Model.getIndicatorIIS()}
\label{\detokenize{javaapi/Model:model-getindicatoriis}}\begin{quote}

\sphinxAtStartPar
Get IIS status of indicator constraints.

\sphinxAtStartPar
\sphinxstylestrong{Synopsis}
\begin{quote}

\sphinxAtStartPar
\sphinxcode{\sphinxupquote{int{[}{]} getIndicatorIIS(GenConstrArray genconstrs)}}
\end{quote}

\sphinxAtStartPar
\sphinxstylestrong{Arguments}
\begin{quote}

\sphinxAtStartPar
\sphinxcode{\sphinxupquote{genconstrs}}: Array of indicator constraints.
\end{quote}

\sphinxAtStartPar
\sphinxstylestrong{Return}
\begin{quote}

\sphinxAtStartPar
IIS status of indicator constraints.
\end{quote}
\end{quote}

\subsubsection{Model.getIndicatorIIS()}
\label{\detokenize{javaapi/Model:id58}}\begin{quote}

\sphinxAtStartPar
Get IIS status of indicator constraints.

\sphinxAtStartPar
\sphinxstylestrong{Synopsis}
\begin{quote}

\sphinxAtStartPar
\sphinxcode{\sphinxupquote{int{[}{]} getIndicatorIIS(GenConstr{[}{]} genconstrs)}}
\end{quote}

\sphinxAtStartPar
\sphinxstylestrong{Arguments}
\begin{quote}

\sphinxAtStartPar
\sphinxcode{\sphinxupquote{genconstrs}}: Array of indicator constraints.
\end{quote}

\sphinxAtStartPar
\sphinxstylestrong{Return}
\begin{quote}

\sphinxAtStartPar
IIS status of indicator constraints.
\end{quote}
\end{quote}

\subsubsection{Model.getIntAttr()}
\label{\detokenize{javaapi/Model:model-getintattr}}\begin{quote}

\sphinxAtStartPar
Get value of a COPT integer attribute

\sphinxAtStartPar
\sphinxstylestrong{Synopsis}
\begin{quote}

\sphinxAtStartPar
\sphinxcode{\sphinxupquote{int getIntAttr(String attr)}}
\end{quote}

\sphinxAtStartPar
\sphinxstylestrong{Arguments}
\begin{quote}

\sphinxAtStartPar
\sphinxcode{\sphinxupquote{attr}}: name of integer attribute.
\end{quote}

\sphinxAtStartPar
\sphinxstylestrong{Return}
\begin{quote}

\sphinxAtStartPar
value of integer attribute.
\end{quote}
\end{quote}

\subsubsection{Model.getIntAttrN()}
\label{\detokenize{javaapi/Model:model-getintattrn}}\begin{quote}

\sphinxAtStartPar
Get value of a integer attribute of a multi\sphinxhyphen{}objective function.

\sphinxAtStartPar
\sphinxstylestrong{Synopsis}
\begin{quote}

\sphinxAtStartPar
\sphinxcode{\sphinxupquote{int getIntAttrN(int idx, String attr)}}
\end{quote}

\sphinxAtStartPar
\sphinxstylestrong{Arguments}
\begin{quote}

\sphinxAtStartPar
\sphinxcode{\sphinxupquote{idx}}: index of a multi\sphinxhyphen{}objective function.

\sphinxAtStartPar
\sphinxcode{\sphinxupquote{attr}}: name of integer attribute.
\end{quote}

\sphinxAtStartPar
\sphinxstylestrong{Return}
\begin{quote}

\sphinxAtStartPar
value of integer attribute.
\end{quote}
\end{quote}

\subsubsection{Model.getIntParam()}
\label{\detokenize{javaapi/Model:model-getintparam}}\begin{quote}

\sphinxAtStartPar
Get value of a COPT integer parameter.

\sphinxAtStartPar
\sphinxstylestrong{Synopsis}
\begin{quote}

\sphinxAtStartPar
\sphinxcode{\sphinxupquote{int getIntParam(String param)}}
\end{quote}

\sphinxAtStartPar
\sphinxstylestrong{Arguments}
\begin{quote}

\sphinxAtStartPar
\sphinxcode{\sphinxupquote{param}}: name of integer parameter.
\end{quote}

\sphinxAtStartPar
\sphinxstylestrong{Return}
\begin{quote}

\sphinxAtStartPar
value of integer parameter.
\end{quote}
\end{quote}

\subsubsection{Model.getIntParamInfo()}
\label{\detokenize{javaapi/Model:model-getintparaminfo}}\begin{quote}

\sphinxAtStartPar
Get current, default, minimum, maximum of COPT integer parameter.

\sphinxAtStartPar
\sphinxstylestrong{Synopsis}
\begin{quote}

\sphinxAtStartPar
\sphinxcode{\sphinxupquote{int{[}{]} getIntParamInfo(String name)}}
\end{quote}

\sphinxAtStartPar
\sphinxstylestrong{Arguments}
\begin{quote}

\sphinxAtStartPar
\sphinxcode{\sphinxupquote{name}}: name of integer parameter.
\end{quote}

\sphinxAtStartPar
\sphinxstylestrong{Return}
\begin{quote}

\sphinxAtStartPar
current, default, minimum, maximum of COPT integer parameter.
\end{quote}
\end{quote}

\subsubsection{Model.getIntParamN()}
\label{\detokenize{javaapi/Model:model-getintparamn}}\begin{quote}

\sphinxAtStartPar
Get value of an integer parameter of a multi\sphinxhyphen{}objective function.

\sphinxAtStartPar
\sphinxstylestrong{Synopsis}
\begin{quote}

\sphinxAtStartPar
\sphinxcode{\sphinxupquote{int getIntParamN(int idx, String param)}}
\end{quote}

\sphinxAtStartPar
\sphinxstylestrong{Arguments}
\begin{quote}

\sphinxAtStartPar
\sphinxcode{\sphinxupquote{idx}}: index of a multi\sphinxhyphen{}objective function.

\sphinxAtStartPar
\sphinxcode{\sphinxupquote{param}}: name of integer parameter.
\end{quote}

\sphinxAtStartPar
\sphinxstylestrong{Return}
\begin{quote}

\sphinxAtStartPar
value of integer parameter.
\end{quote}
\end{quote}

\subsubsection{Model.getLmiCoeff()}
\label{\detokenize{javaapi/Model:model-getlmicoeff}}\begin{quote}

\sphinxAtStartPar
Get the symmetric matrix of variable in LMI constraint.

\sphinxAtStartPar
\sphinxstylestrong{Synopsis}
\begin{quote}

\sphinxAtStartPar
\sphinxcode{\sphinxupquote{SymMatrix getLmiCoeff(LmiConstraint constr, Var var)}}
\end{quote}

\sphinxAtStartPar
\sphinxstylestrong{Arguments}
\begin{quote}

\sphinxAtStartPar
\sphinxcode{\sphinxupquote{constr}}: The desired LMI constraint.

\sphinxAtStartPar
\sphinxcode{\sphinxupquote{var}}: The desired variable.
\end{quote}

\sphinxAtStartPar
\sphinxstylestrong{Return}
\begin{quote}

\sphinxAtStartPar
The associated coefficient matrix.
\end{quote}
\end{quote}

\subsubsection{Model.getLmiConstr()}
\label{\detokenize{javaapi/Model:model-getlmiconstr}}\begin{quote}

\sphinxAtStartPar
Get LMI constraint of given index in model.

\sphinxAtStartPar
\sphinxstylestrong{Synopsis}
\begin{quote}

\sphinxAtStartPar
\sphinxcode{\sphinxupquote{LmiConstraint getLmiConstr(int idx)}}
\end{quote}

\sphinxAtStartPar
\sphinxstylestrong{Arguments}
\begin{quote}

\sphinxAtStartPar
\sphinxcode{\sphinxupquote{idx}}: index of desired LMI constraint.
\end{quote}

\sphinxAtStartPar
\sphinxstylestrong{Return}
\begin{quote}

\sphinxAtStartPar
LMI constraint object.
\end{quote}
\end{quote}

\subsubsection{Model.getLmiConstrByName()}
\label{\detokenize{javaapi/Model:model-getlmiconstrbyname}}\begin{quote}

\sphinxAtStartPar
Get LMI constraint of given name in model.

\sphinxAtStartPar
\sphinxstylestrong{Synopsis}
\begin{quote}

\sphinxAtStartPar
\sphinxcode{\sphinxupquote{LmiConstraint getLmiConstrByName(String name)}}
\end{quote}

\sphinxAtStartPar
\sphinxstylestrong{Arguments}
\begin{quote}

\sphinxAtStartPar
\sphinxcode{\sphinxupquote{name}}: name of desired LMI constraint.
\end{quote}

\sphinxAtStartPar
\sphinxstylestrong{Return}
\begin{quote}

\sphinxAtStartPar
LMI constraint object.
\end{quote}
\end{quote}

\subsubsection{Model.getLmiConstrs()}
\label{\detokenize{javaapi/Model:model-getlmiconstrs}}\begin{quote}

\sphinxAtStartPar
Get all LMI constraints in model.

\sphinxAtStartPar
\sphinxstylestrong{Synopsis}
\begin{quote}

\sphinxAtStartPar
\sphinxcode{\sphinxupquote{LmiConstrArray getLmiConstrs()}}
\end{quote}

\sphinxAtStartPar
\sphinxstylestrong{Return}
\begin{quote}

\sphinxAtStartPar
array object of LMI constraints.
\end{quote}
\end{quote}

\subsubsection{Model.getLmiRhs()}
\label{\detokenize{javaapi/Model:model-getlmirhs}}\begin{quote}

\sphinxAtStartPar
Get the symmetric matrix of constant of LMI constraint.

\sphinxAtStartPar
\sphinxstylestrong{Synopsis}
\begin{quote}

\sphinxAtStartPar
\sphinxcode{\sphinxupquote{SymMatrix getLmiRhs(LmiConstraint constr)}}
\end{quote}

\sphinxAtStartPar
\sphinxstylestrong{Arguments}
\begin{quote}

\sphinxAtStartPar
\sphinxcode{\sphinxupquote{constr}}: The desired LMI constraint.
\end{quote}

\sphinxAtStartPar
\sphinxstylestrong{Return}
\begin{quote}

\sphinxAtStartPar
matrix of constant term.
\end{quote}
\end{quote}

\subsubsection{Model.getLmiRow()}
\label{\detokenize{javaapi/Model:model-getlmirow}}\begin{quote}

\sphinxAtStartPar
Get variables and associated symmetric matrices that participate in a LMI constraint.

\sphinxAtStartPar
\sphinxstylestrong{Synopsis}
\begin{quote}

\sphinxAtStartPar
\sphinxcode{\sphinxupquote{LmiExpr getLmiRow(LmiConstraint constr)}}
\end{quote}

\sphinxAtStartPar
\sphinxstylestrong{Arguments}
\begin{quote}

\sphinxAtStartPar
\sphinxcode{\sphinxupquote{constr}}: given LMI constraint object.
\end{quote}

\sphinxAtStartPar
\sphinxstylestrong{Return}
\begin{quote}

\sphinxAtStartPar
LMI expression object of the LMI constraint.
\end{quote}
\end{quote}

\subsubsection{Model.getLmiSolution()}
\label{\detokenize{javaapi/Model:model-getlmisolution}}\begin{quote}

\sphinxAtStartPar
Get LMI solution.

\sphinxAtStartPar
\sphinxstylestrong{Synopsis}
\begin{quote}

\sphinxAtStartPar
\sphinxcode{\sphinxupquote{Object{[}{]} getLmiSolution()}}
\end{quote}

\sphinxAtStartPar
\sphinxstylestrong{Return}
\begin{quote}

\sphinxAtStartPar
slack and dual values.
\end{quote}
\end{quote}

\subsubsection{Model.getLpSolution()}
\label{\detokenize{javaapi/Model:model-getlpsolution}}\begin{quote}

\sphinxAtStartPar
Get LP solution.

\sphinxAtStartPar
\sphinxstylestrong{Synopsis}
\begin{quote}

\sphinxAtStartPar
\sphinxcode{\sphinxupquote{Object{[}{]} getLpSolution()}}
\end{quote}

\sphinxAtStartPar
\sphinxstylestrong{Return}
\begin{quote}

\sphinxAtStartPar
solution, slack, dual and reduced values.
\end{quote}
\end{quote}

\subsubsection{Model.getNlConstr()}
\label{\detokenize{javaapi/Model:model-getnlconstr}}\begin{quote}

\sphinxAtStartPar
Get a nonlinear constraint of given index in model.

\sphinxAtStartPar
\sphinxstylestrong{Synopsis}
\begin{quote}

\sphinxAtStartPar
\sphinxcode{\sphinxupquote{NlConstraint getNlConstr(int idx)}}
\end{quote}

\sphinxAtStartPar
\sphinxstylestrong{Arguments}
\begin{quote}

\sphinxAtStartPar
\sphinxcode{\sphinxupquote{idx}}: index of the desired nonlinear constraint.
\end{quote}

\sphinxAtStartPar
\sphinxstylestrong{Return}
\begin{quote}

\sphinxAtStartPar
the desired nonlinear constraint object.
\end{quote}
\end{quote}

\subsubsection{Model.getNlConstrBuilder()}
\label{\detokenize{javaapi/Model:model-getnlconstrbuilder}}\begin{quote}

\sphinxAtStartPar
Get builder of a nonlinear constraint in model, including nonlinear expression, sense and RHS.

\sphinxAtStartPar
\sphinxstylestrong{Synopsis}
\begin{quote}

\sphinxAtStartPar
\sphinxcode{\sphinxupquote{NlConstrBuilder getNlConstrBuilder(NlConstraint constr)}}
\end{quote}

\sphinxAtStartPar
\sphinxstylestrong{Arguments}
\begin{quote}

\sphinxAtStartPar
\sphinxcode{\sphinxupquote{constr}}: a nonlinear constraint object.
\end{quote}

\sphinxAtStartPar
\sphinxstylestrong{Return}
\begin{quote}

\sphinxAtStartPar
nonlinear constraint builder object.
\end{quote}
\end{quote}

\subsubsection{Model.getNlConstrBuilders()}
\label{\detokenize{javaapi/Model:model-getnlconstrbuilders}}\begin{quote}

\sphinxAtStartPar
Get builders of all nonlinear constraints in model.

\sphinxAtStartPar
\sphinxstylestrong{Synopsis}
\begin{quote}

\sphinxAtStartPar
\sphinxcode{\sphinxupquote{NlConstrBuilderArray getNlConstrBuilders()}}
\end{quote}

\sphinxAtStartPar
\sphinxstylestrong{Return}
\begin{quote}

\sphinxAtStartPar
array object of nonlinear constraint builders.
\end{quote}
\end{quote}

\subsubsection{Model.getNlConstrByName()}
\label{\detokenize{javaapi/Model:model-getnlconstrbyname}}\begin{quote}

\sphinxAtStartPar
Get a nonlinear constraint of given name in model.

\sphinxAtStartPar
\sphinxstylestrong{Synopsis}
\begin{quote}

\sphinxAtStartPar
\sphinxcode{\sphinxupquote{NlConstraint getNlConstrByName(String name)}}
\end{quote}

\sphinxAtStartPar
\sphinxstylestrong{Arguments}
\begin{quote}

\sphinxAtStartPar
\sphinxcode{\sphinxupquote{name}}: name of the desired constraint.
\end{quote}

\sphinxAtStartPar
\sphinxstylestrong{Return}
\begin{quote}

\sphinxAtStartPar
the desired nonlinear constraint object.
\end{quote}
\end{quote}

\subsubsection{Model.getNlConstrs()}
\label{\detokenize{javaapi/Model:model-getnlconstrs}}\begin{quote}

\sphinxAtStartPar
Get all nonlinear constraints in model.

\sphinxAtStartPar
\sphinxstylestrong{Synopsis}
\begin{quote}

\sphinxAtStartPar
\sphinxcode{\sphinxupquote{NlConstrArray getNlConstrs()}}
\end{quote}

\sphinxAtStartPar
\sphinxstylestrong{Return}
\begin{quote}

\sphinxAtStartPar
array object of nonlinear constraints.
\end{quote}
\end{quote}

\subsubsection{Model.getNlObjective()}
\label{\detokenize{javaapi/Model:model-getnlobjective}}\begin{quote}

\sphinxAtStartPar
Get nonlinear objective of model.

\sphinxAtStartPar
\sphinxstylestrong{Synopsis}
\begin{quote}

\sphinxAtStartPar
\sphinxcode{\sphinxupquote{NlExpr getNlObjective()}}
\end{quote}

\sphinxAtStartPar
\sphinxstylestrong{Return}
\begin{quote}

\sphinxAtStartPar
a nonlinear expression object.
\end{quote}
\end{quote}

\subsubsection{Model.getNlRow()}
\label{\detokenize{javaapi/Model:model-getnlrow}}\begin{quote}

\sphinxAtStartPar
Get nonlinear expression of a nonlinear constraint.

\sphinxAtStartPar
\sphinxstylestrong{Synopsis}
\begin{quote}

\sphinxAtStartPar
\sphinxcode{\sphinxupquote{NlExpr getNlRow(NlConstraint constr)}}
\end{quote}

\sphinxAtStartPar
\sphinxstylestrong{Arguments}
\begin{quote}

\sphinxAtStartPar
\sphinxcode{\sphinxupquote{constr}}: a nonlinear constraint object.
\end{quote}

\sphinxAtStartPar
\sphinxstylestrong{Return}
\begin{quote}

\sphinxAtStartPar
output object of nonlinear expression.
\end{quote}
\end{quote}

\subsubsection{Model.getObjective()}
\label{\detokenize{javaapi/Model:model-getobjective}}\begin{quote}

\sphinxAtStartPar
Get linear expression of objective for model.

\sphinxAtStartPar
\sphinxstylestrong{Synopsis}
\begin{quote}

\sphinxAtStartPar
\sphinxcode{\sphinxupquote{Expr getObjective()}}
\end{quote}

\sphinxAtStartPar
\sphinxstylestrong{Return}
\begin{quote}

\sphinxAtStartPar
an linear expression object.
\end{quote}
\end{quote}

\subsubsection{Model.getObjectiveN()}
\label{\detokenize{javaapi/Model:model-getobjectiven}}\begin{quote}

\sphinxAtStartPar
Get linear expression of a multi\sphinxhyphen{}objective function in model.

\sphinxAtStartPar
\sphinxstylestrong{Synopsis}
\begin{quote}

\sphinxAtStartPar
\sphinxcode{\sphinxupquote{Expr getObjectiveN(int idx)}}
\end{quote}

\sphinxAtStartPar
\sphinxstylestrong{Arguments}
\begin{quote}

\sphinxAtStartPar
\sphinxcode{\sphinxupquote{idx}}: index of a multi\sphinxhyphen{}objective function.
\end{quote}

\sphinxAtStartPar
\sphinxstylestrong{Return}
\begin{quote}

\sphinxAtStartPar
a linear expression object.
\end{quote}
\end{quote}

\subsubsection{Model.getObjParamN()}
\label{\detokenize{javaapi/Model:model-getobjparamn}}\begin{quote}

\sphinxAtStartPar
Get value of objective parameter of a multi\sphinxhyphen{}objective function.

\sphinxAtStartPar
\sphinxstylestrong{Synopsis}
\begin{quote}

\sphinxAtStartPar
\sphinxcode{\sphinxupquote{double getObjParamN(int idx, String param)}}
\end{quote}

\sphinxAtStartPar
\sphinxstylestrong{Arguments}
\begin{quote}

\sphinxAtStartPar
\sphinxcode{\sphinxupquote{idx}}: index of a multi\sphinxhyphen{}objective function.

\sphinxAtStartPar
\sphinxcode{\sphinxupquote{param}}: name of objective parameter, including priority, weight, abstol and reltol.
\end{quote}

\sphinxAtStartPar
\sphinxstylestrong{Return}
\begin{quote}

\sphinxAtStartPar
value of objective parameter.
\end{quote}
\end{quote}

\subsubsection{Model.getPoolObjVal()}
\label{\detokenize{javaapi/Model:model-getpoolobjval}}\begin{quote}

\sphinxAtStartPar
Get the idx\sphinxhyphen{}th objective value in solution pool.

\sphinxAtStartPar
\sphinxstylestrong{Synopsis}
\begin{quote}

\sphinxAtStartPar
\sphinxcode{\sphinxupquote{double getPoolObjVal(int idx)}}
\end{quote}

\sphinxAtStartPar
\sphinxstylestrong{Arguments}
\begin{quote}

\sphinxAtStartPar
\sphinxcode{\sphinxupquote{idx}}: Index of solution.
\end{quote}

\sphinxAtStartPar
\sphinxstylestrong{Return}
\begin{quote}

\sphinxAtStartPar
The requested objective value.
\end{quote}
\end{quote}

\subsubsection{Model.getPoolObjValN()}
\label{\detokenize{javaapi/Model:model-getpoolobjvaln}}\begin{quote}

\sphinxAtStartPar
Get the objective value of required multi\sphinxhyphen{}objective function in solution pool.

\sphinxAtStartPar
\sphinxstylestrong{Synopsis}
\begin{quote}

\sphinxAtStartPar
\sphinxcode{\sphinxupquote{double getPoolObjValN(int idx, int iSol)}}
\end{quote}

\sphinxAtStartPar
\sphinxstylestrong{Arguments}
\begin{quote}

\sphinxAtStartPar
\sphinxcode{\sphinxupquote{idx}}: index of a multi\sphinxhyphen{}objective function.

\sphinxAtStartPar
\sphinxcode{\sphinxupquote{iSol}}: index of solution.
\end{quote}

\sphinxAtStartPar
\sphinxstylestrong{Return}
\begin{quote}

\sphinxAtStartPar
value of required multi\sphinxhyphen{}objective function.
\end{quote}
\end{quote}

\subsubsection{Model.getPoolSolution()}
\label{\detokenize{javaapi/Model:model-getpoolsolution}}\begin{quote}

\sphinxAtStartPar
Get the idx\sphinxhyphen{}th solution in solution pool.

\sphinxAtStartPar
\sphinxstylestrong{Synopsis}
\begin{quote}

\sphinxAtStartPar
\sphinxcode{\sphinxupquote{double{[}{]} getPoolSolution(int idx, VarArray vars)}}
\end{quote}

\sphinxAtStartPar
\sphinxstylestrong{Arguments}
\begin{quote}

\sphinxAtStartPar
\sphinxcode{\sphinxupquote{idx}}: Index of solution.

\sphinxAtStartPar
\sphinxcode{\sphinxupquote{vars}}: The requested variables.
\end{quote}

\sphinxAtStartPar
\sphinxstylestrong{Return}
\begin{quote}

\sphinxAtStartPar
The requested solution.
\end{quote}
\end{quote}

\subsubsection{Model.getPoolSolution()}
\label{\detokenize{javaapi/Model:id59}}\begin{quote}

\sphinxAtStartPar
Get the idx\sphinxhyphen{}th solution in solution pool.

\sphinxAtStartPar
\sphinxstylestrong{Synopsis}
\begin{quote}

\sphinxAtStartPar
\sphinxcode{\sphinxupquote{double{[}{]} getPoolSolution(int idx, Var{[}{]} vars)}}
\end{quote}

\sphinxAtStartPar
\sphinxstylestrong{Arguments}
\begin{quote}

\sphinxAtStartPar
\sphinxcode{\sphinxupquote{idx}}: Index of solution.

\sphinxAtStartPar
\sphinxcode{\sphinxupquote{vars}}: The requested variables.
\end{quote}

\sphinxAtStartPar
\sphinxstylestrong{Return}
\begin{quote}

\sphinxAtStartPar
The requested solution.
\end{quote}
\end{quote}

\subsubsection{Model.getPsdCoeff()}
\label{\detokenize{javaapi/Model:model-getpsdcoeff}}\begin{quote}

\sphinxAtStartPar
Get the symmetric matrix of PSD variable in a PSD constraint.

\sphinxAtStartPar
\sphinxstylestrong{Synopsis}
\begin{quote}

\sphinxAtStartPar
\sphinxcode{\sphinxupquote{SymMatrix getPsdCoeff(PsdConstraint constr, PsdVar var)}}
\end{quote}

\sphinxAtStartPar
\sphinxstylestrong{Arguments}
\begin{quote}

\sphinxAtStartPar
\sphinxcode{\sphinxupquote{constr}}: The desired PSD constraint.

\sphinxAtStartPar
\sphinxcode{\sphinxupquote{var}}: The desired PSD variable.
\end{quote}

\sphinxAtStartPar
\sphinxstylestrong{Return}
\begin{quote}

\sphinxAtStartPar
The associated coefficient matrix.
\end{quote}
\end{quote}

\subsubsection{Model.getPsdConstr()}
\label{\detokenize{javaapi/Model:model-getpsdconstr}}\begin{quote}

\sphinxAtStartPar
Get PSD constraint of given index in model.

\sphinxAtStartPar
\sphinxstylestrong{Synopsis}
\begin{quote}

\sphinxAtStartPar
\sphinxcode{\sphinxupquote{PsdConstraint getPsdConstr(int idx)}}
\end{quote}

\sphinxAtStartPar
\sphinxstylestrong{Arguments}
\begin{quote}

\sphinxAtStartPar
\sphinxcode{\sphinxupquote{idx}}: index of desired PSD constraint.
\end{quote}

\sphinxAtStartPar
\sphinxstylestrong{Return}
\begin{quote}

\sphinxAtStartPar
PSD constraint object.
\end{quote}
\end{quote}

\subsubsection{Model.getPsdConstrBuilder()}
\label{\detokenize{javaapi/Model:model-getpsdconstrbuilder}}\begin{quote}

\sphinxAtStartPar
Get builder of a PSD constraint in model, including PSD variables, sense and associated symmtric matrix.

\sphinxAtStartPar
\sphinxstylestrong{Synopsis}
\begin{quote}

\sphinxAtStartPar
\sphinxcode{\sphinxupquote{PsdConstrBuilder getPsdConstrBuilder(PsdConstraint constr)}}
\end{quote}

\sphinxAtStartPar
\sphinxstylestrong{Arguments}
\begin{quote}

\sphinxAtStartPar
\sphinxcode{\sphinxupquote{constr}}: a PSD constraint object.
\end{quote}

\sphinxAtStartPar
\sphinxstylestrong{Return}
\begin{quote}

\sphinxAtStartPar
PSD constraint builder object.
\end{quote}
\end{quote}

\subsubsection{Model.getPsdConstrBuilders()}
\label{\detokenize{javaapi/Model:model-getpsdconstrbuilders}}\begin{quote}

\sphinxAtStartPar
Get builders of all PSD constraints in model.

\sphinxAtStartPar
\sphinxstylestrong{Synopsis}
\begin{quote}

\sphinxAtStartPar
\sphinxcode{\sphinxupquote{PsdConstrBuilderArray getPsdConstrBuilders()}}
\end{quote}

\sphinxAtStartPar
\sphinxstylestrong{Return}
\begin{quote}

\sphinxAtStartPar
array object of PSD constraint builders.
\end{quote}
\end{quote}

\subsubsection{Model.getPsdConstrByName()}
\label{\detokenize{javaapi/Model:model-getpsdconstrbyname}}\begin{quote}

\sphinxAtStartPar
Get PSD constraint of given name in model.

\sphinxAtStartPar
\sphinxstylestrong{Synopsis}
\begin{quote}

\sphinxAtStartPar
\sphinxcode{\sphinxupquote{PsdConstraint getPsdConstrByName(String name)}}
\end{quote}

\sphinxAtStartPar
\sphinxstylestrong{Arguments}
\begin{quote}

\sphinxAtStartPar
\sphinxcode{\sphinxupquote{name}}: name of desired PSD constraint.
\end{quote}

\sphinxAtStartPar
\sphinxstylestrong{Return}
\begin{quote}

\sphinxAtStartPar
PSD constraint object.
\end{quote}
\end{quote}

\subsubsection{Model.getPsdConstrs()}
\label{\detokenize{javaapi/Model:model-getpsdconstrs}}\begin{quote}

\sphinxAtStartPar
Get all PSD constraints in model.

\sphinxAtStartPar
\sphinxstylestrong{Synopsis}
\begin{quote}

\sphinxAtStartPar
\sphinxcode{\sphinxupquote{PsdConstrArray getPsdConstrs()}}
\end{quote}

\sphinxAtStartPar
\sphinxstylestrong{Return}
\begin{quote}

\sphinxAtStartPar
array object of PSD constraints.
\end{quote}
\end{quote}

\subsubsection{Model.getPsdObjective()}
\label{\detokenize{javaapi/Model:model-getpsdobjective}}\begin{quote}

\sphinxAtStartPar
Get PSD objective of model.

\sphinxAtStartPar
\sphinxstylestrong{Synopsis}
\begin{quote}

\sphinxAtStartPar
\sphinxcode{\sphinxupquote{PsdExpr getPsdObjective()}}
\end{quote}

\sphinxAtStartPar
\sphinxstylestrong{Return}
\begin{quote}

\sphinxAtStartPar
a PSD expression object.
\end{quote}
\end{quote}

\subsubsection{Model.getPsdRow()}
\label{\detokenize{javaapi/Model:model-getpsdrow}}\begin{quote}

\sphinxAtStartPar
Get PSD variables and associated symmetric matrices that participate in a PSD constraint.

\sphinxAtStartPar
\sphinxstylestrong{Synopsis}
\begin{quote}

\sphinxAtStartPar
\sphinxcode{\sphinxupquote{PsdExpr getPsdRow(PsdConstraint constr)}}
\end{quote}

\sphinxAtStartPar
\sphinxstylestrong{Arguments}
\begin{quote}

\sphinxAtStartPar
\sphinxcode{\sphinxupquote{constr}}: a PSD constraint object.
\end{quote}

\sphinxAtStartPar
\sphinxstylestrong{Return}
\begin{quote}

\sphinxAtStartPar
PSD expression object of the PSD constraint.
\end{quote}
\end{quote}

\subsubsection{Model.getPsdSolution()}
\label{\detokenize{javaapi/Model:model-getpsdsolution}}\begin{quote}

\sphinxAtStartPar
Get PSD solution.

\sphinxAtStartPar
\sphinxstylestrong{Synopsis}
\begin{quote}

\sphinxAtStartPar
\sphinxcode{\sphinxupquote{Object{[}{]} getPsdSolution()}}
\end{quote}

\sphinxAtStartPar
\sphinxstylestrong{Return}
\begin{quote}

\sphinxAtStartPar
solution, slack, dual and reduced values.
\end{quote}
\end{quote}

\subsubsection{Model.getPsdVar()}
\label{\detokenize{javaapi/Model:model-getpsdvar}}\begin{quote}

\sphinxAtStartPar
Get a PSD variable of given index in model.

\sphinxAtStartPar
\sphinxstylestrong{Synopsis}
\begin{quote}

\sphinxAtStartPar
\sphinxcode{\sphinxupquote{PsdVar getPsdVar(int idx)}}
\end{quote}

\sphinxAtStartPar
\sphinxstylestrong{Arguments}
\begin{quote}

\sphinxAtStartPar
\sphinxcode{\sphinxupquote{idx}}: index of the desired PSD variable.
\end{quote}

\sphinxAtStartPar
\sphinxstylestrong{Return}
\begin{quote}

\sphinxAtStartPar
the desired PSD variable object.
\end{quote}
\end{quote}

\subsubsection{Model.getPsdVarByName()}
\label{\detokenize{javaapi/Model:model-getpsdvarbyname}}\begin{quote}

\sphinxAtStartPar
Get a PSD variable of given name in model.

\sphinxAtStartPar
\sphinxstylestrong{Synopsis}
\begin{quote}

\sphinxAtStartPar
\sphinxcode{\sphinxupquote{PsdVar getPsdVarByName(String name)}}
\end{quote}

\sphinxAtStartPar
\sphinxstylestrong{Arguments}
\begin{quote}

\sphinxAtStartPar
\sphinxcode{\sphinxupquote{name}}: name of the desired PSD variable.
\end{quote}

\sphinxAtStartPar
\sphinxstylestrong{Return}
\begin{quote}

\sphinxAtStartPar
the desired PSD variable object.
\end{quote}
\end{quote}

\subsubsection{Model.getPsdVars()}
\label{\detokenize{javaapi/Model:model-getpsdvars}}\begin{quote}

\sphinxAtStartPar
Get all PSD variables in model.

\sphinxAtStartPar
\sphinxstylestrong{Synopsis}
\begin{quote}

\sphinxAtStartPar
\sphinxcode{\sphinxupquote{PsdVarArray getPsdVars()}}
\end{quote}

\sphinxAtStartPar
\sphinxstylestrong{Return}
\begin{quote}

\sphinxAtStartPar
array object of PSD variables.
\end{quote}
\end{quote}

\subsubsection{Model.getQConstr()}
\label{\detokenize{javaapi/Model:model-getqconstr}}\begin{quote}

\sphinxAtStartPar
Get a quadratic constraint of given index in model.

\sphinxAtStartPar
\sphinxstylestrong{Synopsis}
\begin{quote}

\sphinxAtStartPar
\sphinxcode{\sphinxupquote{QConstraint getQConstr(int idx)}}
\end{quote}

\sphinxAtStartPar
\sphinxstylestrong{Arguments}
\begin{quote}

\sphinxAtStartPar
\sphinxcode{\sphinxupquote{idx}}: index of the desired quadratic constraint.
\end{quote}

\sphinxAtStartPar
\sphinxstylestrong{Return}
\begin{quote}

\sphinxAtStartPar
the desired quadratic constraint object.
\end{quote}
\end{quote}

\subsubsection{Model.getQConstrBuilder()}
\label{\detokenize{javaapi/Model:model-getqconstrbuilder}}\begin{quote}

\sphinxAtStartPar
Get builder of a quadratic constraint in model, including variables and associated coefficients, sense and RHS.

\sphinxAtStartPar
\sphinxstylestrong{Synopsis}
\begin{quote}

\sphinxAtStartPar
\sphinxcode{\sphinxupquote{QConstrBuilder getQConstrBuilder(QConstraint constr)}}
\end{quote}

\sphinxAtStartPar
\sphinxstylestrong{Arguments}
\begin{quote}

\sphinxAtStartPar
\sphinxcode{\sphinxupquote{constr}}: a constraint object.
\end{quote}

\sphinxAtStartPar
\sphinxstylestrong{Return}
\begin{quote}

\sphinxAtStartPar
constraint builder object.
\end{quote}
\end{quote}

\subsubsection{Model.getQConstrBuilders()}
\label{\detokenize{javaapi/Model:model-getqconstrbuilders}}\begin{quote}

\sphinxAtStartPar
Get builders of all constraints in model.

\sphinxAtStartPar
\sphinxstylestrong{Synopsis}
\begin{quote}

\sphinxAtStartPar
\sphinxcode{\sphinxupquote{QConstrBuilderArray getQConstrBuilders()}}
\end{quote}

\sphinxAtStartPar
\sphinxstylestrong{Return}
\begin{quote}

\sphinxAtStartPar
array object of constraint builders.
\end{quote}
\end{quote}

\subsubsection{Model.getQConstrByName()}
\label{\detokenize{javaapi/Model:model-getqconstrbyname}}\begin{quote}

\sphinxAtStartPar
Get a quadratic constraint of given name in model.

\sphinxAtStartPar
\sphinxstylestrong{Synopsis}
\begin{quote}

\sphinxAtStartPar
\sphinxcode{\sphinxupquote{QConstraint getQConstrByName(String name)}}
\end{quote}

\sphinxAtStartPar
\sphinxstylestrong{Arguments}
\begin{quote}

\sphinxAtStartPar
\sphinxcode{\sphinxupquote{name}}: name of the desired constraint.
\end{quote}

\sphinxAtStartPar
\sphinxstylestrong{Return}
\begin{quote}

\sphinxAtStartPar
the desired quadratic constraint object.
\end{quote}
\end{quote}

\subsubsection{Model.getQConstrs()}
\label{\detokenize{javaapi/Model:model-getqconstrs}}\begin{quote}

\sphinxAtStartPar
Get all quadratic constraints in model.

\sphinxAtStartPar
\sphinxstylestrong{Synopsis}
\begin{quote}

\sphinxAtStartPar
\sphinxcode{\sphinxupquote{QConstrArray getQConstrs()}}
\end{quote}

\sphinxAtStartPar
\sphinxstylestrong{Return}
\begin{quote}

\sphinxAtStartPar
array object of quadratic constraints.
\end{quote}
\end{quote}

\subsubsection{Model.getQuadObjective()}
\label{\detokenize{javaapi/Model:model-getquadobjective}}\begin{quote}

\sphinxAtStartPar
Get quadratic objective of model.

\sphinxAtStartPar
\sphinxstylestrong{Synopsis}
\begin{quote}

\sphinxAtStartPar
\sphinxcode{\sphinxupquote{QuadExpr getQuadObjective()}}
\end{quote}

\sphinxAtStartPar
\sphinxstylestrong{Return}
\begin{quote}

\sphinxAtStartPar
a quadratic expression object.
\end{quote}
\end{quote}

\subsubsection{Model.getQuadRow()}
\label{\detokenize{javaapi/Model:model-getquadrow}}\begin{quote}

\sphinxAtStartPar
Get quadratic expression that participate in quadratic constraint.

\sphinxAtStartPar
\sphinxstylestrong{Synopsis}
\begin{quote}

\sphinxAtStartPar
\sphinxcode{\sphinxupquote{QuadExpr getQuadRow(QConstraint constr)}}
\end{quote}

\sphinxAtStartPar
\sphinxstylestrong{Arguments}
\begin{quote}

\sphinxAtStartPar
\sphinxcode{\sphinxupquote{constr}}: a quadratic constraint object.
\end{quote}

\sphinxAtStartPar
\sphinxstylestrong{Return}
\begin{quote}

\sphinxAtStartPar
quadratic expression object of the constraint.
\end{quote}
\end{quote}

\subsubsection{Model.getRow()}
\label{\detokenize{javaapi/Model:model-getrow}}\begin{quote}

\sphinxAtStartPar
Get variables that participate in a constraint, and the associated coefficients.

\sphinxAtStartPar
\sphinxstylestrong{Synopsis}
\begin{quote}

\sphinxAtStartPar
\sphinxcode{\sphinxupquote{Expr getRow(Constraint constr)}}
\end{quote}

\sphinxAtStartPar
\sphinxstylestrong{Arguments}
\begin{quote}

\sphinxAtStartPar
\sphinxcode{\sphinxupquote{constr}}: a constraint object.
\end{quote}

\sphinxAtStartPar
\sphinxstylestrong{Return}
\begin{quote}

\sphinxAtStartPar
expression object of the constraint.
\end{quote}
\end{quote}

\subsubsection{Model.getRowBasis()}
\label{\detokenize{javaapi/Model:model-getrowbasis}}\begin{quote}

\sphinxAtStartPar
Get status of row basis.

\sphinxAtStartPar
\sphinxstylestrong{Synopsis}
\begin{quote}

\sphinxAtStartPar
\sphinxcode{\sphinxupquote{int{[}{]} getRowBasis()}}
\end{quote}

\sphinxAtStartPar
\sphinxstylestrong{Return}
\begin{quote}

\sphinxAtStartPar
basis status.
\end{quote}
\end{quote}

\subsubsection{Model.getSolution()}
\label{\detokenize{javaapi/Model:model-getsolution}}\begin{quote}

\sphinxAtStartPar
Get MIP solution.

\sphinxAtStartPar
\sphinxstylestrong{Synopsis}
\begin{quote}

\sphinxAtStartPar
\sphinxcode{\sphinxupquote{double{[}{]} getSolution()}}
\end{quote}

\sphinxAtStartPar
\sphinxstylestrong{Return}
\begin{quote}

\sphinxAtStartPar
solution values.
\end{quote}
\end{quote}

\subsubsection{Model.getSos()}
\label{\detokenize{javaapi/Model:model-getsos}}\begin{quote}

\sphinxAtStartPar
Get a SOS constraint of given index in model.

\sphinxAtStartPar
\sphinxstylestrong{Synopsis}
\begin{quote}

\sphinxAtStartPar
\sphinxcode{\sphinxupquote{Sos getSos(int idx)}}
\end{quote}

\sphinxAtStartPar
\sphinxstylestrong{Arguments}
\begin{quote}

\sphinxAtStartPar
\sphinxcode{\sphinxupquote{idx}}: index of the desired SOS constraint.
\end{quote}

\sphinxAtStartPar
\sphinxstylestrong{Return}
\begin{quote}

\sphinxAtStartPar
the desired SOS constraint object.
\end{quote}
\end{quote}

\subsubsection{Model.getSosBuilders()}
\label{\detokenize{javaapi/Model:model-getsosbuilders}}\begin{quote}

\sphinxAtStartPar
Get builders of all SOS constraints in model.

\sphinxAtStartPar
\sphinxstylestrong{Synopsis}
\begin{quote}

\sphinxAtStartPar
\sphinxcode{\sphinxupquote{SosBuilderArray getSosBuilders()}}
\end{quote}

\sphinxAtStartPar
\sphinxstylestrong{Return}
\begin{quote}

\sphinxAtStartPar
array object of SOS constraint builders.
\end{quote}
\end{quote}

\subsubsection{Model.getSosBuilders()}
\label{\detokenize{javaapi/Model:id60}}\begin{quote}

\sphinxAtStartPar
Get builders of given SOS constraints in model.

\sphinxAtStartPar
\sphinxstylestrong{Synopsis}
\begin{quote}

\sphinxAtStartPar
\sphinxcode{\sphinxupquote{SosBuilderArray getSosBuilders(Sos{[}{]} soss)}}
\end{quote}

\sphinxAtStartPar
\sphinxstylestrong{Arguments}
\begin{quote}

\sphinxAtStartPar
\sphinxcode{\sphinxupquote{soss}}: array of SOS constraints.
\end{quote}

\sphinxAtStartPar
\sphinxstylestrong{Return}
\begin{quote}

\sphinxAtStartPar
array object of desired SOS constraint builders.
\end{quote}
\end{quote}

\subsubsection{Model.getSosBuilders()}
\label{\detokenize{javaapi/Model:id61}}\begin{quote}

\sphinxAtStartPar
Get builders of given SOS constraints in model.

\sphinxAtStartPar
\sphinxstylestrong{Synopsis}
\begin{quote}

\sphinxAtStartPar
\sphinxcode{\sphinxupquote{SosBuilderArray getSosBuilders(SosArray soss)}}
\end{quote}

\sphinxAtStartPar
\sphinxstylestrong{Arguments}
\begin{quote}

\sphinxAtStartPar
\sphinxcode{\sphinxupquote{soss}}: array of SOS constraints.
\end{quote}

\sphinxAtStartPar
\sphinxstylestrong{Return}
\begin{quote}

\sphinxAtStartPar
array object of desired SOS constraint builders.
\end{quote}
\end{quote}

\subsubsection{Model.getSOSIIS()}
\label{\detokenize{javaapi/Model:model-getsosiis}}\begin{quote}

\sphinxAtStartPar
Get IIS status of SOS constraints.

\sphinxAtStartPar
\sphinxstylestrong{Synopsis}
\begin{quote}

\sphinxAtStartPar
\sphinxcode{\sphinxupquote{int{[}{]} getSOSIIS(SosArray soss)}}
\end{quote}

\sphinxAtStartPar
\sphinxstylestrong{Arguments}
\begin{quote}

\sphinxAtStartPar
\sphinxcode{\sphinxupquote{soss}}: Array of SOS constraints.
\end{quote}

\sphinxAtStartPar
\sphinxstylestrong{Return}
\begin{quote}

\sphinxAtStartPar
IIS status of SOS constraints.
\end{quote}
\end{quote}

\subsubsection{Model.getSOSIIS()}
\label{\detokenize{javaapi/Model:id62}}\begin{quote}

\sphinxAtStartPar
Get IIS status of SOS constraints.

\sphinxAtStartPar
\sphinxstylestrong{Synopsis}
\begin{quote}

\sphinxAtStartPar
\sphinxcode{\sphinxupquote{int{[}{]} getSOSIIS(Sos{[}{]} soss)}}
\end{quote}

\sphinxAtStartPar
\sphinxstylestrong{Arguments}
\begin{quote}

\sphinxAtStartPar
\sphinxcode{\sphinxupquote{soss}}: Array of SOS constraints.
\end{quote}

\sphinxAtStartPar
\sphinxstylestrong{Return}
\begin{quote}

\sphinxAtStartPar
IIS status of SOS constraints.
\end{quote}
\end{quote}

\subsubsection{Model.getSoss()}
\label{\detokenize{javaapi/Model:model-getsoss}}\begin{quote}

\sphinxAtStartPar
Get all SOS constraints in model.

\sphinxAtStartPar
\sphinxstylestrong{Synopsis}
\begin{quote}

\sphinxAtStartPar
\sphinxcode{\sphinxupquote{SosArray getSoss()}}
\end{quote}

\sphinxAtStartPar
\sphinxstylestrong{Return}
\begin{quote}

\sphinxAtStartPar
array object of SOS constraints.
\end{quote}
\end{quote}

\subsubsection{Model.getSymMat()}
\label{\detokenize{javaapi/Model:model-getsymmat}}\begin{quote}

\sphinxAtStartPar
Get a symmetric matrix of given index in model.

\sphinxAtStartPar
\sphinxstylestrong{Synopsis}
\begin{quote}

\sphinxAtStartPar
\sphinxcode{\sphinxupquote{SymMatrix getSymMat(int idx)}}
\end{quote}

\sphinxAtStartPar
\sphinxstylestrong{Arguments}
\begin{quote}

\sphinxAtStartPar
\sphinxcode{\sphinxupquote{idx}}: index of the desired symmetric matrix.
\end{quote}

\sphinxAtStartPar
\sphinxstylestrong{Return}
\begin{quote}

\sphinxAtStartPar
the desired symmetric matrix object.
\end{quote}
\end{quote}

\subsubsection{Model.getVar()}
\label{\detokenize{javaapi/Model:model-getvar}}\begin{quote}

\sphinxAtStartPar
Get a variable of given index in model.

\sphinxAtStartPar
\sphinxstylestrong{Synopsis}
\begin{quote}

\sphinxAtStartPar
\sphinxcode{\sphinxupquote{Var getVar(int idx)}}
\end{quote}

\sphinxAtStartPar
\sphinxstylestrong{Arguments}
\begin{quote}

\sphinxAtStartPar
\sphinxcode{\sphinxupquote{idx}}: index of the desired variable.
\end{quote}

\sphinxAtStartPar
\sphinxstylestrong{Return}
\begin{quote}

\sphinxAtStartPar
the desired variable object.
\end{quote}
\end{quote}

\subsubsection{Model.getVarByName()}
\label{\detokenize{javaapi/Model:model-getvarbyname}}\begin{quote}

\sphinxAtStartPar
Get a variable of given name in model.

\sphinxAtStartPar
\sphinxstylestrong{Synopsis}
\begin{quote}

\sphinxAtStartPar
\sphinxcode{\sphinxupquote{Var getVarByName(String name)}}
\end{quote}

\sphinxAtStartPar
\sphinxstylestrong{Arguments}
\begin{quote}

\sphinxAtStartPar
\sphinxcode{\sphinxupquote{name}}: name of the desired variable.
\end{quote}

\sphinxAtStartPar
\sphinxstylestrong{Return}
\begin{quote}

\sphinxAtStartPar
the desired variable object.
\end{quote}
\end{quote}

\subsubsection{Model.getVarLowerIIS()}
\label{\detokenize{javaapi/Model:model-getvarloweriis}}\begin{quote}

\sphinxAtStartPar
Get IIS status of lower bounds of variables.

\sphinxAtStartPar
\sphinxstylestrong{Synopsis}
\begin{quote}

\sphinxAtStartPar
\sphinxcode{\sphinxupquote{int{[}{]} getVarLowerIIS(VarArray vars)}}
\end{quote}

\sphinxAtStartPar
\sphinxstylestrong{Arguments}
\begin{quote}

\sphinxAtStartPar
\sphinxcode{\sphinxupquote{vars}}: Array of variables.
\end{quote}

\sphinxAtStartPar
\sphinxstylestrong{Return}
\begin{quote}

\sphinxAtStartPar
IIS status of lower bounds of variables.
\end{quote}
\end{quote}

\subsubsection{Model.getVarLowerIIS()}
\label{\detokenize{javaapi/Model:id63}}\begin{quote}

\sphinxAtStartPar
Get IIS status of lower bounds of variables.

\sphinxAtStartPar
\sphinxstylestrong{Synopsis}
\begin{quote}

\sphinxAtStartPar
\sphinxcode{\sphinxupquote{int{[}{]} getVarLowerIIS(Var{[}{]} vars)}}
\end{quote}

\sphinxAtStartPar
\sphinxstylestrong{Arguments}
\begin{quote}

\sphinxAtStartPar
\sphinxcode{\sphinxupquote{vars}}: Array of variables.
\end{quote}

\sphinxAtStartPar
\sphinxstylestrong{Return}
\begin{quote}

\sphinxAtStartPar
IIS status of lower bounds of variables.
\end{quote}
\end{quote}

\subsubsection{Model.getVars()}
\label{\detokenize{javaapi/Model:model-getvars}}\begin{quote}

\sphinxAtStartPar
Get all variables in model.

\sphinxAtStartPar
\sphinxstylestrong{Synopsis}
\begin{quote}

\sphinxAtStartPar
\sphinxcode{\sphinxupquote{VarArray getVars()}}
\end{quote}

\sphinxAtStartPar
\sphinxstylestrong{Return}
\begin{quote}

\sphinxAtStartPar
variable array object.
\end{quote}
\end{quote}

\subsubsection{Model.getVarUpperIIS()}
\label{\detokenize{javaapi/Model:model-getvarupperiis}}\begin{quote}

\sphinxAtStartPar
Get IIS status of upper bounds of variables.

\sphinxAtStartPar
\sphinxstylestrong{Synopsis}
\begin{quote}

\sphinxAtStartPar
\sphinxcode{\sphinxupquote{int{[}{]} getVarUpperIIS(VarArray vars)}}
\end{quote}

\sphinxAtStartPar
\sphinxstylestrong{Arguments}
\begin{quote}

\sphinxAtStartPar
\sphinxcode{\sphinxupquote{vars}}: Array of variables.
\end{quote}

\sphinxAtStartPar
\sphinxstylestrong{Return}
\begin{quote}

\sphinxAtStartPar
IIS status of upper bounds of variables.
\end{quote}
\end{quote}

\subsubsection{Model.getVarUpperIIS()}
\label{\detokenize{javaapi/Model:id64}}\begin{quote}

\sphinxAtStartPar
Get IIS status of upper bounds of variables.

\sphinxAtStartPar
\sphinxstylestrong{Synopsis}
\begin{quote}

\sphinxAtStartPar
\sphinxcode{\sphinxupquote{int{[}{]} getVarUpperIIS(Var{[}{]} vars)}}
\end{quote}

\sphinxAtStartPar
\sphinxstylestrong{Arguments}
\begin{quote}

\sphinxAtStartPar
\sphinxcode{\sphinxupquote{vars}}: Array of variables.
\end{quote}

\sphinxAtStartPar
\sphinxstylestrong{Return}
\begin{quote}

\sphinxAtStartPar
IIS status of upper bounds of variables.
\end{quote}
\end{quote}

\subsubsection{Model.interrupt()}
\label{\detokenize{javaapi/Model:model-interrupt}}\begin{quote}

\sphinxAtStartPar
Interrupt optimization of current problem.

\sphinxAtStartPar
\sphinxstylestrong{Synopsis}
\begin{quote}

\sphinxAtStartPar
\sphinxcode{\sphinxupquote{void interrupt()}}
\end{quote}
\end{quote}

\subsubsection{Model.loadMipStart()}
\label{\detokenize{javaapi/Model:model-loadmipstart}}\begin{quote}

\sphinxAtStartPar
Load final initial values of variables to the problem.

\sphinxAtStartPar
\sphinxstylestrong{Synopsis}
\begin{quote}

\sphinxAtStartPar
\sphinxcode{\sphinxupquote{void loadMipStart()}}
\end{quote}
\end{quote}

\subsubsection{Model.loadTuneParam()}
\label{\detokenize{javaapi/Model:model-loadtuneparam}}\begin{quote}

\sphinxAtStartPar
Load specified tuned parameters into model.

\sphinxAtStartPar
\sphinxstylestrong{Synopsis}
\begin{quote}

\sphinxAtStartPar
\sphinxcode{\sphinxupquote{void loadTuneParam(int idx)}}
\end{quote}

\sphinxAtStartPar
\sphinxstylestrong{Arguments}
\begin{quote}

\sphinxAtStartPar
\sphinxcode{\sphinxupquote{idx}}: Index of tuned parameters.
\end{quote}
\end{quote}

\subsubsection{Model.read()}
\label{\detokenize{javaapi/Model:model-read}}\begin{quote}

\sphinxAtStartPar
Read problem, solution, basis, MIP start or COPT parameters from file.

\sphinxAtStartPar
\sphinxstylestrong{Synopsis}
\begin{quote}

\sphinxAtStartPar
\sphinxcode{\sphinxupquote{void read(String filename)}}
\end{quote}

\sphinxAtStartPar
\sphinxstylestrong{Arguments}
\begin{quote}

\sphinxAtStartPar
\sphinxcode{\sphinxupquote{filename}}: an input file name.
\end{quote}
\end{quote}

\subsubsection{Model.readBasis()}
\label{\detokenize{javaapi/Model:model-readbasis}}\begin{quote}

\sphinxAtStartPar
Read basis from file.

\sphinxAtStartPar
\sphinxstylestrong{Synopsis}
\begin{quote}

\sphinxAtStartPar
\sphinxcode{\sphinxupquote{void readBasis(String filename)}}
\end{quote}

\sphinxAtStartPar
\sphinxstylestrong{Arguments}
\begin{quote}

\sphinxAtStartPar
\sphinxcode{\sphinxupquote{filename}}: an input file name.
\end{quote}
\end{quote}

\subsubsection{Model.readBin()}
\label{\detokenize{javaapi/Model:model-readbin}}\begin{quote}

\sphinxAtStartPar
Read problem in COPT binary format from file.

\sphinxAtStartPar
\sphinxstylestrong{Synopsis}
\begin{quote}

\sphinxAtStartPar
\sphinxcode{\sphinxupquote{void readBin(String filename)}}
\end{quote}

\sphinxAtStartPar
\sphinxstylestrong{Arguments}
\begin{quote}

\sphinxAtStartPar
\sphinxcode{\sphinxupquote{filename}}: an input file name.
\end{quote}
\end{quote}

\subsubsection{Model.readCbf()}
\label{\detokenize{javaapi/Model:model-readcbf}}\begin{quote}

\sphinxAtStartPar
Read problem in CBF format from file.

\sphinxAtStartPar
\sphinxstylestrong{Synopsis}
\begin{quote}

\sphinxAtStartPar
\sphinxcode{\sphinxupquote{void readCbf(String filename)}}
\end{quote}

\sphinxAtStartPar
\sphinxstylestrong{Arguments}
\begin{quote}

\sphinxAtStartPar
\sphinxcode{\sphinxupquote{filename}}: an input file name.
\end{quote}
\end{quote}

\subsubsection{Model.readJsonSol()}
\label{\detokenize{javaapi/Model:model-readjsonsol}}\begin{quote}

\sphinxAtStartPar
Read solution in format of JSON from file.

\sphinxAtStartPar
\sphinxstylestrong{Synopsis}
\begin{quote}

\sphinxAtStartPar
\sphinxcode{\sphinxupquote{void readJsonSol(String filename)}}
\end{quote}

\sphinxAtStartPar
\sphinxstylestrong{Arguments}
\begin{quote}

\sphinxAtStartPar
\sphinxcode{\sphinxupquote{filename}}: an input file name.
\end{quote}
\end{quote}

\subsubsection{Model.readLp()}
\label{\detokenize{javaapi/Model:model-readlp}}\begin{quote}

\sphinxAtStartPar
Read problem in LP format from file.

\sphinxAtStartPar
\sphinxstylestrong{Synopsis}
\begin{quote}

\sphinxAtStartPar
\sphinxcode{\sphinxupquote{void readLp(String filename)}}
\end{quote}

\sphinxAtStartPar
\sphinxstylestrong{Arguments}
\begin{quote}

\sphinxAtStartPar
\sphinxcode{\sphinxupquote{filename}}: an input file name.
\end{quote}
\end{quote}

\subsubsection{Model.readMps()}
\label{\detokenize{javaapi/Model:model-readmps}}\begin{quote}

\sphinxAtStartPar
Read problem in MPS format from file.

\sphinxAtStartPar
\sphinxstylestrong{Synopsis}
\begin{quote}

\sphinxAtStartPar
\sphinxcode{\sphinxupquote{void readMps(String filename)}}
\end{quote}

\sphinxAtStartPar
\sphinxstylestrong{Arguments}
\begin{quote}

\sphinxAtStartPar
\sphinxcode{\sphinxupquote{filename}}: an input file name.
\end{quote}
\end{quote}

\subsubsection{Model.readMst()}
\label{\detokenize{javaapi/Model:model-readmst}}\begin{quote}

\sphinxAtStartPar
Read MIP start information from file.

\sphinxAtStartPar
\sphinxstylestrong{Synopsis}
\begin{quote}

\sphinxAtStartPar
\sphinxcode{\sphinxupquote{void readMst(String filename)}}
\end{quote}

\sphinxAtStartPar
\sphinxstylestrong{Arguments}
\begin{quote}

\sphinxAtStartPar
\sphinxcode{\sphinxupquote{filename}}: an input file name.
\end{quote}
\end{quote}

\subsubsection{Model.readOrd()}
\label{\detokenize{javaapi/Model:model-readord}}\begin{quote}

\sphinxAtStartPar
Read branching order from file.

\sphinxAtStartPar
\sphinxstylestrong{Synopsis}
\begin{quote}

\sphinxAtStartPar
\sphinxcode{\sphinxupquote{void readOrd(String filename)}}
\end{quote}

\sphinxAtStartPar
\sphinxstylestrong{Arguments}
\begin{quote}

\sphinxAtStartPar
\sphinxcode{\sphinxupquote{filename}}: an input file name.
\end{quote}
\end{quote}

\subsubsection{Model.readParam()}
\label{\detokenize{javaapi/Model:model-readparam}}\begin{quote}

\sphinxAtStartPar
Read COPT parameters from file.

\sphinxAtStartPar
\sphinxstylestrong{Synopsis}
\begin{quote}

\sphinxAtStartPar
\sphinxcode{\sphinxupquote{void readParam(String filename)}}
\end{quote}

\sphinxAtStartPar
\sphinxstylestrong{Arguments}
\begin{quote}

\sphinxAtStartPar
\sphinxcode{\sphinxupquote{filename}}: an input file name.
\end{quote}
\end{quote}

\subsubsection{Model.readSdpa()}
\label{\detokenize{javaapi/Model:model-readsdpa}}\begin{quote}

\sphinxAtStartPar
Read problem in SDPA format from file.

\sphinxAtStartPar
\sphinxstylestrong{Synopsis}
\begin{quote}

\sphinxAtStartPar
\sphinxcode{\sphinxupquote{void readSdpa(String filename)}}
\end{quote}

\sphinxAtStartPar
\sphinxstylestrong{Arguments}
\begin{quote}

\sphinxAtStartPar
\sphinxcode{\sphinxupquote{filename}}: an input file name.
\end{quote}
\end{quote}

\subsubsection{Model.readSol()}
\label{\detokenize{javaapi/Model:model-readsol}}\begin{quote}

\sphinxAtStartPar
Read solution from file.

\sphinxAtStartPar
\sphinxstylestrong{Synopsis}
\begin{quote}

\sphinxAtStartPar
\sphinxcode{\sphinxupquote{void readSol(String filename)}}
\end{quote}

\sphinxAtStartPar
\sphinxstylestrong{Arguments}
\begin{quote}

\sphinxAtStartPar
\sphinxcode{\sphinxupquote{filename}}: an input file name.
\end{quote}
\end{quote}

\subsubsection{Model.readTune()}
\label{\detokenize{javaapi/Model:model-readtune}}\begin{quote}

\sphinxAtStartPar
Read tuning parameters from file.

\sphinxAtStartPar
\sphinxstylestrong{Synopsis}
\begin{quote}

\sphinxAtStartPar
\sphinxcode{\sphinxupquote{void readTune(String filename)}}
\end{quote}

\sphinxAtStartPar
\sphinxstylestrong{Arguments}
\begin{quote}

\sphinxAtStartPar
\sphinxcode{\sphinxupquote{filename}}: an input file name.
\end{quote}
\end{quote}

\subsubsection{Model.remove()}
\label{\detokenize{javaapi/Model:model-remove}}\begin{quote}

\sphinxAtStartPar
Remove an array of variables from model.

\sphinxAtStartPar
\sphinxstylestrong{Synopsis}
\begin{quote}

\sphinxAtStartPar
\sphinxcode{\sphinxupquote{void remove(Var{[}{]} vars)}}
\end{quote}

\sphinxAtStartPar
\sphinxstylestrong{Arguments}
\begin{quote}

\sphinxAtStartPar
\sphinxcode{\sphinxupquote{vars}}: a list of variables.
\end{quote}
\end{quote}

\subsubsection{Model.remove()}
\label{\detokenize{javaapi/Model:id65}}\begin{quote}

\sphinxAtStartPar
Remove array of variables from model.

\sphinxAtStartPar
\sphinxstylestrong{Synopsis}
\begin{quote}

\sphinxAtStartPar
\sphinxcode{\sphinxupquote{void remove(VarArray vars)}}
\end{quote}

\sphinxAtStartPar
\sphinxstylestrong{Arguments}
\begin{quote}

\sphinxAtStartPar
\sphinxcode{\sphinxupquote{vars}}: an array of variables.
\end{quote}
\end{quote}

\subsubsection{Model.remove()}
\label{\detokenize{javaapi/Model:id66}}\begin{quote}

\sphinxAtStartPar
Remove a list of constraints from model.

\sphinxAtStartPar
\sphinxstylestrong{Synopsis}
\begin{quote}

\sphinxAtStartPar
\sphinxcode{\sphinxupquote{void remove(Constraint{[}{]} constrs)}}
\end{quote}

\sphinxAtStartPar
\sphinxstylestrong{Arguments}
\begin{quote}

\sphinxAtStartPar
\sphinxcode{\sphinxupquote{constrs}}: a list of constraints.
\end{quote}
\end{quote}

\subsubsection{Model.remove()}
\label{\detokenize{javaapi/Model:id67}}\begin{quote}

\sphinxAtStartPar
Remove a list of constraints from model.

\sphinxAtStartPar
\sphinxstylestrong{Synopsis}
\begin{quote}

\sphinxAtStartPar
\sphinxcode{\sphinxupquote{void remove(ConstrArray constrs)}}
\end{quote}

\sphinxAtStartPar
\sphinxstylestrong{Arguments}
\begin{quote}

\sphinxAtStartPar
\sphinxcode{\sphinxupquote{constrs}}: an array of constraints.
\end{quote}
\end{quote}

\subsubsection{Model.remove()}
\label{\detokenize{javaapi/Model:id68}}\begin{quote}

\sphinxAtStartPar
Remove an array of nonlinear constraints from model.

\sphinxAtStartPar
\sphinxstylestrong{Synopsis}
\begin{quote}

\sphinxAtStartPar
\sphinxcode{\sphinxupquote{void remove(NlConstraint{[}{]} constrs)}}
\end{quote}

\sphinxAtStartPar
\sphinxstylestrong{Arguments}
\begin{quote}

\sphinxAtStartPar
\sphinxcode{\sphinxupquote{constrs}}: array of nonlinear constraints.
\end{quote}
\end{quote}

\subsubsection{Model.remove()}
\label{\detokenize{javaapi/Model:id69}}\begin{quote}

\sphinxAtStartPar
Remove a list of nonlinear constraints from model.

\sphinxAtStartPar
\sphinxstylestrong{Synopsis}
\begin{quote}

\sphinxAtStartPar
\sphinxcode{\sphinxupquote{void remove(NlConstrArray constrs)}}
\end{quote}

\sphinxAtStartPar
\sphinxstylestrong{Arguments}
\begin{quote}

\sphinxAtStartPar
\sphinxcode{\sphinxupquote{constrs}}: array object of nonlinear constraints.
\end{quote}
\end{quote}

\subsubsection{Model.remove()}
\label{\detokenize{javaapi/Model:id70}}\begin{quote}

\sphinxAtStartPar
Remove a list of SOS constraints from model.

\sphinxAtStartPar
\sphinxstylestrong{Synopsis}
\begin{quote}

\sphinxAtStartPar
\sphinxcode{\sphinxupquote{void remove(Sos{[}{]} soss)}}
\end{quote}

\sphinxAtStartPar
\sphinxstylestrong{Arguments}
\begin{quote}

\sphinxAtStartPar
\sphinxcode{\sphinxupquote{soss}}: a list of SOS constraints.
\end{quote}
\end{quote}

\subsubsection{Model.remove()}
\label{\detokenize{javaapi/Model:id71}}\begin{quote}

\sphinxAtStartPar
Remove a list of SOS constraints from model.

\sphinxAtStartPar
\sphinxstylestrong{Synopsis}
\begin{quote}

\sphinxAtStartPar
\sphinxcode{\sphinxupquote{void remove(SosArray soss)}}
\end{quote}

\sphinxAtStartPar
\sphinxstylestrong{Arguments}
\begin{quote}

\sphinxAtStartPar
\sphinxcode{\sphinxupquote{soss}}: an array of SOS constraints.
\end{quote}
\end{quote}

\subsubsection{Model.remove()}
\label{\detokenize{javaapi/Model:id72}}\begin{quote}

\sphinxAtStartPar
Remove a list of cone constraints from model.

\sphinxAtStartPar
\sphinxstylestrong{Synopsis}
\begin{quote}

\sphinxAtStartPar
\sphinxcode{\sphinxupquote{void remove(Cone{[}{]} cones)}}
\end{quote}

\sphinxAtStartPar
\sphinxstylestrong{Arguments}
\begin{quote}

\sphinxAtStartPar
\sphinxcode{\sphinxupquote{cones}}: a list of cone constraints.
\end{quote}
\end{quote}

\subsubsection{Model.remove()}
\label{\detokenize{javaapi/Model:id73}}\begin{quote}

\sphinxAtStartPar
Remove a list of cone constraints from model.

\sphinxAtStartPar
\sphinxstylestrong{Synopsis}
\begin{quote}

\sphinxAtStartPar
\sphinxcode{\sphinxupquote{void remove(ConeArray cones)}}
\end{quote}

\sphinxAtStartPar
\sphinxstylestrong{Arguments}
\begin{quote}

\sphinxAtStartPar
\sphinxcode{\sphinxupquote{cones}}: an array of cone constraints.
\end{quote}
\end{quote}

\subsubsection{Model.remove()}
\label{\detokenize{javaapi/Model:id74}}\begin{quote}

\sphinxAtStartPar
Remove a list of exponential cone constraints from model.

\sphinxAtStartPar
\sphinxstylestrong{Synopsis}
\begin{quote}

\sphinxAtStartPar
\sphinxcode{\sphinxupquote{void remove(ExpCone{[}{]} cones)}}
\end{quote}

\sphinxAtStartPar
\sphinxstylestrong{Arguments}
\begin{quote}

\sphinxAtStartPar
\sphinxcode{\sphinxupquote{cones}}: a list of exponential cone constraints.
\end{quote}
\end{quote}

\subsubsection{Model.remove()}
\label{\detokenize{javaapi/Model:id75}}\begin{quote}

\sphinxAtStartPar
Remove a list of exponential cone constraints from model.

\sphinxAtStartPar
\sphinxstylestrong{Synopsis}
\begin{quote}

\sphinxAtStartPar
\sphinxcode{\sphinxupquote{void remove(ExpConeArray cones)}}
\end{quote}

\sphinxAtStartPar
\sphinxstylestrong{Arguments}
\begin{quote}

\sphinxAtStartPar
\sphinxcode{\sphinxupquote{cones}}: an array of exponential cone constraints.
\end{quote}
\end{quote}

\subsubsection{Model.remove()}
\label{\detokenize{javaapi/Model:id76}}\begin{quote}

\sphinxAtStartPar
Remove a list of affine cone constraints from model.

\sphinxAtStartPar
\sphinxstylestrong{Synopsis}
\begin{quote}

\sphinxAtStartPar
\sphinxcode{\sphinxupquote{void remove(AffineCone{[}{]} cones)}}
\end{quote}

\sphinxAtStartPar
\sphinxstylestrong{Arguments}
\begin{quote}

\sphinxAtStartPar
\sphinxcode{\sphinxupquote{cones}}: a list of affine cone constraints.
\end{quote}
\end{quote}

\subsubsection{Model.remove()}
\label{\detokenize{javaapi/Model:id77}}\begin{quote}

\sphinxAtStartPar
Remove an array of affine cone constraints from model.

\sphinxAtStartPar
\sphinxstylestrong{Synopsis}
\begin{quote}

\sphinxAtStartPar
\sphinxcode{\sphinxupquote{void remove(AffineConeArray cones)}}
\end{quote}

\sphinxAtStartPar
\sphinxstylestrong{Arguments}
\begin{quote}

\sphinxAtStartPar
\sphinxcode{\sphinxupquote{cones}}: an array of affine cone constraints.
\end{quote}
\end{quote}

\subsubsection{Model.remove()}
\label{\detokenize{javaapi/Model:id78}}\begin{quote}

\sphinxAtStartPar
Remove a list of gernal constraints from model.

\sphinxAtStartPar
\sphinxstylestrong{Synopsis}
\begin{quote}

\sphinxAtStartPar
\sphinxcode{\sphinxupquote{void remove(GenConstr{[}{]} genConstrs)}}
\end{quote}

\sphinxAtStartPar
\sphinxstylestrong{Arguments}
\begin{quote}

\sphinxAtStartPar
\sphinxcode{\sphinxupquote{genConstrs}}: a list of general constraints.
\end{quote}
\end{quote}

\subsubsection{Model.remove()}
\label{\detokenize{javaapi/Model:id79}}\begin{quote}

\sphinxAtStartPar
Remove a list of gernal constraints from model.

\sphinxAtStartPar
\sphinxstylestrong{Synopsis}
\begin{quote}

\sphinxAtStartPar
\sphinxcode{\sphinxupquote{void remove(GenConstrArray genConstrs)}}
\end{quote}

\sphinxAtStartPar
\sphinxstylestrong{Arguments}
\begin{quote}

\sphinxAtStartPar
\sphinxcode{\sphinxupquote{genConstrs}}: an array of general constraints.
\end{quote}
\end{quote}

\subsubsection{Model.remove()}
\label{\detokenize{javaapi/Model:id80}}\begin{quote}

\sphinxAtStartPar
Remove a list of quadratic constraints from model.

\sphinxAtStartPar
\sphinxstylestrong{Synopsis}
\begin{quote}

\sphinxAtStartPar
\sphinxcode{\sphinxupquote{void remove(QConstraint{[}{]} qconstrs)}}
\end{quote}

\sphinxAtStartPar
\sphinxstylestrong{Arguments}
\begin{quote}

\sphinxAtStartPar
\sphinxcode{\sphinxupquote{qconstrs}}: an array of quadratic constraints.
\end{quote}
\end{quote}

\subsubsection{Model.remove()}
\label{\detokenize{javaapi/Model:id81}}\begin{quote}

\sphinxAtStartPar
Remove a list of quadratic constraints from model.

\sphinxAtStartPar
\sphinxstylestrong{Synopsis}
\begin{quote}

\sphinxAtStartPar
\sphinxcode{\sphinxupquote{void remove(QConstrArray qconstrs)}}
\end{quote}

\sphinxAtStartPar
\sphinxstylestrong{Arguments}
\begin{quote}

\sphinxAtStartPar
\sphinxcode{\sphinxupquote{qconstrs}}: an array of quadratic constraints.
\end{quote}
\end{quote}

\subsubsection{Model.remove()}
\label{\detokenize{javaapi/Model:id82}}\begin{quote}

\sphinxAtStartPar
Remove a list of PSD variables from model.

\sphinxAtStartPar
\sphinxstylestrong{Synopsis}
\begin{quote}

\sphinxAtStartPar
\sphinxcode{\sphinxupquote{void remove(PsdVar{[}{]} vars)}}
\end{quote}

\sphinxAtStartPar
\sphinxstylestrong{Arguments}
\begin{quote}

\sphinxAtStartPar
\sphinxcode{\sphinxupquote{vars}}: an array of PSD variables.
\end{quote}
\end{quote}

\subsubsection{Model.remove()}
\label{\detokenize{javaapi/Model:id83}}\begin{quote}

\sphinxAtStartPar
Remove a list of PSD variables from model.

\sphinxAtStartPar
\sphinxstylestrong{Synopsis}
\begin{quote}

\sphinxAtStartPar
\sphinxcode{\sphinxupquote{void remove(PsdVarArray vars)}}
\end{quote}

\sphinxAtStartPar
\sphinxstylestrong{Arguments}
\begin{quote}

\sphinxAtStartPar
\sphinxcode{\sphinxupquote{vars}}: an array of PSD variables.
\end{quote}
\end{quote}

\subsubsection{Model.remove()}
\label{\detokenize{javaapi/Model:id84}}\begin{quote}

\sphinxAtStartPar
Remove a list of PSD constraints from model.

\sphinxAtStartPar
\sphinxstylestrong{Synopsis}
\begin{quote}

\sphinxAtStartPar
\sphinxcode{\sphinxupquote{void remove(PsdConstraint{[}{]} constrs)}}
\end{quote}

\sphinxAtStartPar
\sphinxstylestrong{Arguments}
\begin{quote}

\sphinxAtStartPar
\sphinxcode{\sphinxupquote{constrs}}: an array of PSD constraints.
\end{quote}
\end{quote}

\subsubsection{Model.remove()}
\label{\detokenize{javaapi/Model:id85}}\begin{quote}

\sphinxAtStartPar
Remove a list of PSD constraints from model.

\sphinxAtStartPar
\sphinxstylestrong{Synopsis}
\begin{quote}

\sphinxAtStartPar
\sphinxcode{\sphinxupquote{void remove(PsdConstrArray constrs)}}
\end{quote}

\sphinxAtStartPar
\sphinxstylestrong{Arguments}
\begin{quote}

\sphinxAtStartPar
\sphinxcode{\sphinxupquote{constrs}}: an array of PSD constraints.
\end{quote}
\end{quote}

\subsubsection{Model.remove()}
\label{\detokenize{javaapi/Model:id86}}\begin{quote}

\sphinxAtStartPar
Remove a list of LMI constraints from model.

\sphinxAtStartPar
\sphinxstylestrong{Synopsis}
\begin{quote}

\sphinxAtStartPar
\sphinxcode{\sphinxupquote{void remove(LmiConstrArray constrs)}}
\end{quote}

\sphinxAtStartPar
\sphinxstylestrong{Arguments}
\begin{quote}

\sphinxAtStartPar
\sphinxcode{\sphinxupquote{constrs}}: an array of LMI constraints.
\end{quote}
\end{quote}

\subsubsection{Model.remove()}
\label{\detokenize{javaapi/Model:id87}}\begin{quote}

\sphinxAtStartPar
Remove a list of LMI constraints from model.

\sphinxAtStartPar
\sphinxstylestrong{Synopsis}
\begin{quote}

\sphinxAtStartPar
\sphinxcode{\sphinxupquote{void remove(LmiConstraint{[}{]} constrs)}}
\end{quote}

\sphinxAtStartPar
\sphinxstylestrong{Arguments}
\begin{quote}

\sphinxAtStartPar
\sphinxcode{\sphinxupquote{constrs}}: an array of LMI constraints.
\end{quote}
\end{quote}

\subsubsection{Model.reset()}
\label{\detokenize{javaapi/Model:model-reset}}\begin{quote}

\sphinxAtStartPar
Reset solution of problem only.

\sphinxAtStartPar
\sphinxstylestrong{Synopsis}
\begin{quote}

\sphinxAtStartPar
\sphinxcode{\sphinxupquote{void reset()}}
\end{quote}
\end{quote}

\subsubsection{Model.resetAll()}
\label{\detokenize{javaapi/Model:model-resetall}}\begin{quote}

\sphinxAtStartPar
Reset solution of problem, and additional information such as MIP start, etc.

\sphinxAtStartPar
\sphinxstylestrong{Synopsis}
\begin{quote}

\sphinxAtStartPar
\sphinxcode{\sphinxupquote{void resetAll()}}
\end{quote}
\end{quote}

\subsubsection{Model.resetObjParamN()}
\label{\detokenize{javaapi/Model:model-resetobjparamn}}\begin{quote}

\sphinxAtStartPar
Reset objective parameters of a multi\sphinxhyphen{}objective function.

\sphinxAtStartPar
\sphinxstylestrong{Synopsis}
\begin{quote}

\sphinxAtStartPar
\sphinxcode{\sphinxupquote{void resetObjParamN(int idx)}}
\end{quote}

\sphinxAtStartPar
\sphinxstylestrong{Arguments}
\begin{quote}

\sphinxAtStartPar
\sphinxcode{\sphinxupquote{idx}}: index of a multi\sphinxhyphen{}objective function.
\end{quote}
\end{quote}

\subsubsection{Model.resetParam()}
\label{\detokenize{javaapi/Model:model-resetparam}}\begin{quote}

\sphinxAtStartPar
Reset parameters to default settings.

\sphinxAtStartPar
\sphinxstylestrong{Synopsis}
\begin{quote}

\sphinxAtStartPar
\sphinxcode{\sphinxupquote{void resetParam()}}
\end{quote}
\end{quote}

\subsubsection{Model.resetParamN()}
\label{\detokenize{javaapi/Model:model-resetparamn}}\begin{quote}

\sphinxAtStartPar
Reset double and integer parameters of a multi\sphinxhyphen{}objective function.

\sphinxAtStartPar
\sphinxstylestrong{Synopsis}
\begin{quote}

\sphinxAtStartPar
\sphinxcode{\sphinxupquote{void resetParamN(int idx)}}
\end{quote}

\sphinxAtStartPar
\sphinxstylestrong{Arguments}
\begin{quote}

\sphinxAtStartPar
\sphinxcode{\sphinxupquote{idx}}: index of a multi\sphinxhyphen{}objective function.
\end{quote}
\end{quote}

\subsubsection{Model.set()}
\label{\detokenize{javaapi/Model:model-set}}\begin{quote}

\sphinxAtStartPar
Set values of information associated with variables.

\sphinxAtStartPar
\sphinxstylestrong{Synopsis}
\begin{quote}

\sphinxAtStartPar
\sphinxcode{\sphinxupquote{void set(}}
\begin{quote}

\sphinxAtStartPar
\sphinxcode{\sphinxupquote{String name,}}

\sphinxAtStartPar
\sphinxcode{\sphinxupquote{Var{[}{]} vars,}}

\sphinxAtStartPar
\sphinxcode{\sphinxupquote{double{[}{]} vals)}}
\end{quote}
\end{quote}

\sphinxAtStartPar
\sphinxstylestrong{Arguments}
\begin{quote}

\sphinxAtStartPar
\sphinxcode{\sphinxupquote{name}}: name of information.

\sphinxAtStartPar
\sphinxcode{\sphinxupquote{vars}}: a list of interested variables.

\sphinxAtStartPar
\sphinxcode{\sphinxupquote{vals}}: values of information.
\end{quote}
\end{quote}

\subsubsection{Model.set()}
\label{\detokenize{javaapi/Model:id88}}\begin{quote}

\sphinxAtStartPar
Set values of information associated with variables.

\sphinxAtStartPar
\sphinxstylestrong{Synopsis}
\begin{quote}

\sphinxAtStartPar
\sphinxcode{\sphinxupquote{void set(}}
\begin{quote}

\sphinxAtStartPar
\sphinxcode{\sphinxupquote{String name,}}

\sphinxAtStartPar
\sphinxcode{\sphinxupquote{VarArray vars,}}

\sphinxAtStartPar
\sphinxcode{\sphinxupquote{double{[}{]} vals)}}
\end{quote}
\end{quote}

\sphinxAtStartPar
\sphinxstylestrong{Arguments}
\begin{quote}

\sphinxAtStartPar
\sphinxcode{\sphinxupquote{name}}: name of information.

\sphinxAtStartPar
\sphinxcode{\sphinxupquote{vars}}: array of interested variables.

\sphinxAtStartPar
\sphinxcode{\sphinxupquote{vals}}: values of information.
\end{quote}
\end{quote}

\subsubsection{Model.set()}
\label{\detokenize{javaapi/Model:id89}}\begin{quote}

\sphinxAtStartPar
Set values of information associated with constraints.

\sphinxAtStartPar
\sphinxstylestrong{Synopsis}
\begin{quote}

\sphinxAtStartPar
\sphinxcode{\sphinxupquote{void set(}}
\begin{quote}

\sphinxAtStartPar
\sphinxcode{\sphinxupquote{String name,}}

\sphinxAtStartPar
\sphinxcode{\sphinxupquote{Constraint{[}{]} constrs,}}

\sphinxAtStartPar
\sphinxcode{\sphinxupquote{double{[}{]} vals)}}
\end{quote}
\end{quote}

\sphinxAtStartPar
\sphinxstylestrong{Arguments}
\begin{quote}

\sphinxAtStartPar
\sphinxcode{\sphinxupquote{name}}: name of information.

\sphinxAtStartPar
\sphinxcode{\sphinxupquote{constrs}}: a list of interested constraints.

\sphinxAtStartPar
\sphinxcode{\sphinxupquote{vals}}: values of information.
\end{quote}
\end{quote}

\subsubsection{Model.set()}
\label{\detokenize{javaapi/Model:id90}}\begin{quote}

\sphinxAtStartPar
Set values of information associated with constraints.

\sphinxAtStartPar
\sphinxstylestrong{Synopsis}
\begin{quote}

\sphinxAtStartPar
\sphinxcode{\sphinxupquote{void set(}}
\begin{quote}

\sphinxAtStartPar
\sphinxcode{\sphinxupquote{String name,}}

\sphinxAtStartPar
\sphinxcode{\sphinxupquote{ConstrArray constrs,}}

\sphinxAtStartPar
\sphinxcode{\sphinxupquote{double{[}{]} vals)}}
\end{quote}
\end{quote}

\sphinxAtStartPar
\sphinxstylestrong{Arguments}
\begin{quote}

\sphinxAtStartPar
\sphinxcode{\sphinxupquote{name}}: name of information.

\sphinxAtStartPar
\sphinxcode{\sphinxupquote{constrs}}: array of interested constraints.

\sphinxAtStartPar
\sphinxcode{\sphinxupquote{vals}}: values of information.
\end{quote}
\end{quote}

\subsubsection{Model.set()}
\label{\detokenize{javaapi/Model:id91}}\begin{quote}

\sphinxAtStartPar
Set values of information associated with nonlinear constraints.

\sphinxAtStartPar
\sphinxstylestrong{Synopsis}
\begin{quote}

\sphinxAtStartPar
\sphinxcode{\sphinxupquote{void set(}}
\begin{quote}

\sphinxAtStartPar
\sphinxcode{\sphinxupquote{String name,}}

\sphinxAtStartPar
\sphinxcode{\sphinxupquote{NlConstraint{[}{]} constrs,}}

\sphinxAtStartPar
\sphinxcode{\sphinxupquote{double{[}{]} vals)}}
\end{quote}
\end{quote}

\sphinxAtStartPar
\sphinxstylestrong{Arguments}
\begin{quote}

\sphinxAtStartPar
\sphinxcode{\sphinxupquote{name}}: name of double information.

\sphinxAtStartPar
\sphinxcode{\sphinxupquote{constrs}}: array of desired nonlinear constraints.

\sphinxAtStartPar
\sphinxcode{\sphinxupquote{vals}}: array of values of information.
\end{quote}
\end{quote}

\subsubsection{Model.set()}
\label{\detokenize{javaapi/Model:id92}}\begin{quote}

\sphinxAtStartPar
Set values of information associated with nonlinear constraints.

\sphinxAtStartPar
\sphinxstylestrong{Synopsis}
\begin{quote}

\sphinxAtStartPar
\sphinxcode{\sphinxupquote{void set(}}
\begin{quote}

\sphinxAtStartPar
\sphinxcode{\sphinxupquote{String name,}}

\sphinxAtStartPar
\sphinxcode{\sphinxupquote{NlConstrArray constrs,}}

\sphinxAtStartPar
\sphinxcode{\sphinxupquote{double{[}{]} vals)}}
\end{quote}
\end{quote}

\sphinxAtStartPar
\sphinxstylestrong{Arguments}
\begin{quote}

\sphinxAtStartPar
\sphinxcode{\sphinxupquote{name}}: name of double information.

\sphinxAtStartPar
\sphinxcode{\sphinxupquote{constrs}}: an array object of desired nonlinear constraints.

\sphinxAtStartPar
\sphinxcode{\sphinxupquote{vals}}: array of values of information.
\end{quote}
\end{quote}

\subsubsection{Model.set()}
\label{\detokenize{javaapi/Model:id93}}\begin{quote}

\sphinxAtStartPar
Set values of information associated with PSD constraints.

\sphinxAtStartPar
\sphinxstylestrong{Synopsis}
\begin{quote}

\sphinxAtStartPar
\sphinxcode{\sphinxupquote{void set(}}
\begin{quote}

\sphinxAtStartPar
\sphinxcode{\sphinxupquote{String name,}}

\sphinxAtStartPar
\sphinxcode{\sphinxupquote{PsdConstraint{[}{]} constrs,}}

\sphinxAtStartPar
\sphinxcode{\sphinxupquote{double{[}{]} vals)}}
\end{quote}
\end{quote}

\sphinxAtStartPar
\sphinxstylestrong{Arguments}
\begin{quote}

\sphinxAtStartPar
\sphinxcode{\sphinxupquote{name}}: name of information.

\sphinxAtStartPar
\sphinxcode{\sphinxupquote{constrs}}: a list of desired PSD constraints.

\sphinxAtStartPar
\sphinxcode{\sphinxupquote{vals}}: array of values of information.
\end{quote}
\end{quote}

\subsubsection{Model.set()}
\label{\detokenize{javaapi/Model:id94}}\begin{quote}

\sphinxAtStartPar
Set values of information associated with PSD constraints.

\sphinxAtStartPar
\sphinxstylestrong{Synopsis}
\begin{quote}

\sphinxAtStartPar
\sphinxcode{\sphinxupquote{void set(}}
\begin{quote}

\sphinxAtStartPar
\sphinxcode{\sphinxupquote{String name,}}

\sphinxAtStartPar
\sphinxcode{\sphinxupquote{PsdConstrArray constrs,}}

\sphinxAtStartPar
\sphinxcode{\sphinxupquote{double{[}{]} vals)}}
\end{quote}
\end{quote}

\sphinxAtStartPar
\sphinxstylestrong{Arguments}
\begin{quote}

\sphinxAtStartPar
\sphinxcode{\sphinxupquote{name}}: name of information.

\sphinxAtStartPar
\sphinxcode{\sphinxupquote{constrs}}: a list of desired PSD constraints.

\sphinxAtStartPar
\sphinxcode{\sphinxupquote{vals}}: array of values of information.
\end{quote}
\end{quote}

\subsubsection{Model.setBasis()}
\label{\detokenize{javaapi/Model:model-setbasis}}\begin{quote}

\sphinxAtStartPar
Set column and row basis status to model.

\sphinxAtStartPar
\sphinxstylestrong{Synopsis}
\begin{quote}

\sphinxAtStartPar
\sphinxcode{\sphinxupquote{void setBasis(int{[}{]} colbasis, int{[}{]} rowbasis)}}
\end{quote}

\sphinxAtStartPar
\sphinxstylestrong{Arguments}
\begin{quote}

\sphinxAtStartPar
\sphinxcode{\sphinxupquote{colbasis}}: status of column basis.

\sphinxAtStartPar
\sphinxcode{\sphinxupquote{rowbasis}}: status of row basis.
\end{quote}
\end{quote}

\subsubsection{Model.setCallback()}
\label{\detokenize{javaapi/Model:model-setcallback}}\begin{quote}

\sphinxAtStartPar
Set user callback to COPT model.

\sphinxAtStartPar
\sphinxstylestrong{Synopsis}
\begin{quote}

\sphinxAtStartPar
\sphinxcode{\sphinxupquote{void setCallback(CallbackBase cb, int cbctx)}}
\end{quote}

\sphinxAtStartPar
\sphinxstylestrong{Arguments}
\begin{quote}

\sphinxAtStartPar
\sphinxcode{\sphinxupquote{cb}}: user callback instance, inheriting from CallbackBase class.

\sphinxAtStartPar
\sphinxcode{\sphinxupquote{cbctx}}: COPT callback context.
\end{quote}
\end{quote}

\subsubsection{Model.setCoeff()}
\label{\detokenize{javaapi/Model:model-setcoeff}}\begin{quote}

\sphinxAtStartPar
Set the coefficient of a variable in a linear constraint.

\sphinxAtStartPar
\sphinxstylestrong{Synopsis}
\begin{quote}

\sphinxAtStartPar
\sphinxcode{\sphinxupquote{void setCoeff(}}
\begin{quote}

\sphinxAtStartPar
\sphinxcode{\sphinxupquote{Constraint constr,}}

\sphinxAtStartPar
\sphinxcode{\sphinxupquote{Var var,}}

\sphinxAtStartPar
\sphinxcode{\sphinxupquote{double newVal)}}
\end{quote}
\end{quote}

\sphinxAtStartPar
\sphinxstylestrong{Arguments}
\begin{quote}

\sphinxAtStartPar
\sphinxcode{\sphinxupquote{constr}}: The requested constraint.

\sphinxAtStartPar
\sphinxcode{\sphinxupquote{var}}: The requested variable.

\sphinxAtStartPar
\sphinxcode{\sphinxupquote{newVal}}: New coefficient.
\end{quote}
\end{quote}

\subsubsection{Model.setCoeffs()}
\label{\detokenize{javaapi/Model:model-setcoeffs}}\begin{quote}

\sphinxAtStartPar
Set a list of coefficients in the model.

\sphinxAtStartPar
\sphinxstylestrong{Synopsis}
\begin{quote}

\sphinxAtStartPar
\sphinxcode{\sphinxupquote{void setCoeffs(}}
\begin{quote}

\sphinxAtStartPar
\sphinxcode{\sphinxupquote{Constraint{[}{]} constrs,}}

\sphinxAtStartPar
\sphinxcode{\sphinxupquote{Var{[}{]} vars,}}

\sphinxAtStartPar
\sphinxcode{\sphinxupquote{double{[}{]} vals)}}
\end{quote}
\end{quote}

\sphinxAtStartPar
\sphinxstylestrong{Arguments}
\begin{quote}

\sphinxAtStartPar
\sphinxcode{\sphinxupquote{constrs}}: Array of constraints for coefficients to be set.

\sphinxAtStartPar
\sphinxcode{\sphinxupquote{vars}}: Array of vars for coefficients to be set.

\sphinxAtStartPar
\sphinxcode{\sphinxupquote{vals}}: New values for coefficients.
\end{quote}
\end{quote}

\subsubsection{Model.setCoeffs()}
\label{\detokenize{javaapi/Model:id95}}\begin{quote}

\sphinxAtStartPar
Set a list of coefficients in the model.

\sphinxAtStartPar
\sphinxstylestrong{Synopsis}
\begin{quote}

\sphinxAtStartPar
\sphinxcode{\sphinxupquote{void setCoeffs(}}
\begin{quote}

\sphinxAtStartPar
\sphinxcode{\sphinxupquote{ConstrArray constrs,}}

\sphinxAtStartPar
\sphinxcode{\sphinxupquote{VarArray vars,}}

\sphinxAtStartPar
\sphinxcode{\sphinxupquote{double{[}{]} vals)}}
\end{quote}
\end{quote}

\sphinxAtStartPar
\sphinxstylestrong{Arguments}
\begin{quote}

\sphinxAtStartPar
\sphinxcode{\sphinxupquote{constrs}}: A list of constraints for coefficients to be set.

\sphinxAtStartPar
\sphinxcode{\sphinxupquote{vars}}: A list of vars for coefficients to be set.

\sphinxAtStartPar
\sphinxcode{\sphinxupquote{vals}}: New values for coefficients.
\end{quote}
\end{quote}

\subsubsection{Model.setDblParam()}
\label{\detokenize{javaapi/Model:model-setdblparam}}\begin{quote}

\sphinxAtStartPar
Set value of a COPT double parameter.

\sphinxAtStartPar
\sphinxstylestrong{Synopsis}
\begin{quote}

\sphinxAtStartPar
\sphinxcode{\sphinxupquote{void setDblParam(String param, double val)}}
\end{quote}

\sphinxAtStartPar
\sphinxstylestrong{Arguments}
\begin{quote}

\sphinxAtStartPar
\sphinxcode{\sphinxupquote{param}}: name of double parameter.

\sphinxAtStartPar
\sphinxcode{\sphinxupquote{val}}: double value.
\end{quote}
\end{quote}

\subsubsection{Model.setDblParamN()}
\label{\detokenize{javaapi/Model:model-setdblparamn}}\begin{quote}

\sphinxAtStartPar
Set value of a double parameter of a multi\sphinxhyphen{}objective function.

\sphinxAtStartPar
\sphinxstylestrong{Synopsis}
\begin{quote}

\sphinxAtStartPar
\sphinxcode{\sphinxupquote{void setDblParamN(}}
\begin{quote}

\sphinxAtStartPar
\sphinxcode{\sphinxupquote{int idx,}}

\sphinxAtStartPar
\sphinxcode{\sphinxupquote{String param,}}

\sphinxAtStartPar
\sphinxcode{\sphinxupquote{double val)}}
\end{quote}
\end{quote}

\sphinxAtStartPar
\sphinxstylestrong{Arguments}
\begin{quote}

\sphinxAtStartPar
\sphinxcode{\sphinxupquote{idx}}: index of a multi\sphinxhyphen{}objective function.

\sphinxAtStartPar
\sphinxcode{\sphinxupquote{param}}: name of double parameter.

\sphinxAtStartPar
\sphinxcode{\sphinxupquote{val}}: new value of double parameter.
\end{quote}
\end{quote}

\subsubsection{Model.setIntParam()}
\label{\detokenize{javaapi/Model:model-setintparam}}\begin{quote}

\sphinxAtStartPar
Set value of a COPT integer parameter.

\sphinxAtStartPar
\sphinxstylestrong{Synopsis}
\begin{quote}

\sphinxAtStartPar
\sphinxcode{\sphinxupquote{void setIntParam(String param, int val)}}
\end{quote}

\sphinxAtStartPar
\sphinxstylestrong{Arguments}
\begin{quote}

\sphinxAtStartPar
\sphinxcode{\sphinxupquote{param}}: name of integer parameter.

\sphinxAtStartPar
\sphinxcode{\sphinxupquote{val}}: integer value.
\end{quote}
\end{quote}

\subsubsection{Model.setIntParamN()}
\label{\detokenize{javaapi/Model:model-setintparamn}}\begin{quote}

\sphinxAtStartPar
Set value of an integer parameter of a multi\sphinxhyphen{}objective function.

\sphinxAtStartPar
\sphinxstylestrong{Synopsis}
\begin{quote}

\sphinxAtStartPar
\sphinxcode{\sphinxupquote{void setIntParamN(}}
\begin{quote}

\sphinxAtStartPar
\sphinxcode{\sphinxupquote{int idx,}}

\sphinxAtStartPar
\sphinxcode{\sphinxupquote{String param,}}

\sphinxAtStartPar
\sphinxcode{\sphinxupquote{int val)}}
\end{quote}
\end{quote}

\sphinxAtStartPar
\sphinxstylestrong{Arguments}
\begin{quote}

\sphinxAtStartPar
\sphinxcode{\sphinxupquote{idx}}: index of a multi\sphinxhyphen{}objective function.

\sphinxAtStartPar
\sphinxcode{\sphinxupquote{param}}: name of integer parameter.

\sphinxAtStartPar
\sphinxcode{\sphinxupquote{val}}: new value of integer parameter.
\end{quote}
\end{quote}

\subsubsection{Model.setLmiCoeff()}
\label{\detokenize{javaapi/Model:model-setlmicoeff}}\begin{quote}

\sphinxAtStartPar
Set the coefficient matrix of a variable in LMI constraint.

\sphinxAtStartPar
\sphinxstylestrong{Synopsis}
\begin{quote}

\sphinxAtStartPar
\sphinxcode{\sphinxupquote{void setLmiCoeff(}}
\begin{quote}

\sphinxAtStartPar
\sphinxcode{\sphinxupquote{LmiConstraint constr,}}

\sphinxAtStartPar
\sphinxcode{\sphinxupquote{Var var,}}

\sphinxAtStartPar
\sphinxcode{\sphinxupquote{SymMatrix mat)}}
\end{quote}
\end{quote}

\sphinxAtStartPar
\sphinxstylestrong{Arguments}
\begin{quote}

\sphinxAtStartPar
\sphinxcode{\sphinxupquote{constr}}: The desired LMI constraint.

\sphinxAtStartPar
\sphinxcode{\sphinxupquote{var}}: The desired variable.

\sphinxAtStartPar
\sphinxcode{\sphinxupquote{mat}}: new coefficient matrix.
\end{quote}
\end{quote}

\subsubsection{Model.setLmiRhs()}
\label{\detokenize{javaapi/Model:model-setlmirhs}}\begin{quote}

\sphinxAtStartPar
Set constant matrix of LMI constraint.

\sphinxAtStartPar
\sphinxstylestrong{Synopsis}
\begin{quote}

\sphinxAtStartPar
\sphinxcode{\sphinxupquote{void setLmiRhs(LmiConstraint constr, SymMatrix mat)}}
\end{quote}

\sphinxAtStartPar
\sphinxstylestrong{Arguments}
\begin{quote}

\sphinxAtStartPar
\sphinxcode{\sphinxupquote{constr}}: The desired LMI constraint.

\sphinxAtStartPar
\sphinxcode{\sphinxupquote{mat}}: new constant matrix.
\end{quote}
\end{quote}

\subsubsection{Model.setLpSolution()}
\label{\detokenize{javaapi/Model:model-setlpsolution}}\begin{quote}

\sphinxAtStartPar
Set LP solution.

\sphinxAtStartPar
\sphinxstylestrong{Synopsis}
\begin{quote}

\sphinxAtStartPar
\sphinxcode{\sphinxupquote{void setLpSolution(}}
\begin{quote}

\sphinxAtStartPar
\sphinxcode{\sphinxupquote{double{[}{]} value,}}

\sphinxAtStartPar
\sphinxcode{\sphinxupquote{double{[}{]} slack,}}

\sphinxAtStartPar
\sphinxcode{\sphinxupquote{double{[}{]} rowDual,}}

\sphinxAtStartPar
\sphinxcode{\sphinxupquote{double{[}{]} redCost)}}
\end{quote}
\end{quote}

\sphinxAtStartPar
\sphinxstylestrong{Arguments}
\begin{quote}

\sphinxAtStartPar
\sphinxcode{\sphinxupquote{value}}: solution of variables.

\sphinxAtStartPar
\sphinxcode{\sphinxupquote{slack}}: slack of constraints.

\sphinxAtStartPar
\sphinxcode{\sphinxupquote{rowDual}}: dual value of constraints.

\sphinxAtStartPar
\sphinxcode{\sphinxupquote{redCost}}: dual value of variables.
\end{quote}
\end{quote}

\subsubsection{Model.setMipStart()}
\label{\detokenize{javaapi/Model:model-setmipstart}}\begin{quote}

\sphinxAtStartPar
Set initial values for variables of given number, starting from the first one.

\sphinxAtStartPar
\sphinxstylestrong{Synopsis}
\begin{quote}

\sphinxAtStartPar
\sphinxcode{\sphinxupquote{void setMipStart(int count, double{[}{]} vals)}}
\end{quote}

\sphinxAtStartPar
\sphinxstylestrong{Arguments}
\begin{quote}

\sphinxAtStartPar
\sphinxcode{\sphinxupquote{count}}: the number of variables to set.

\sphinxAtStartPar
\sphinxcode{\sphinxupquote{vals}}: values of variables.
\end{quote}
\end{quote}

\subsubsection{Model.setMipStart()}
\label{\detokenize{javaapi/Model:id96}}\begin{quote}

\sphinxAtStartPar
Set initial value for the specified variable.

\sphinxAtStartPar
\sphinxstylestrong{Synopsis}
\begin{quote}

\sphinxAtStartPar
\sphinxcode{\sphinxupquote{void setMipStart(Var var, double val)}}
\end{quote}

\sphinxAtStartPar
\sphinxstylestrong{Arguments}
\begin{quote}

\sphinxAtStartPar
\sphinxcode{\sphinxupquote{var}}: an interested variable.

\sphinxAtStartPar
\sphinxcode{\sphinxupquote{val}}: initial value of the variable.
\end{quote}
\end{quote}

\subsubsection{Model.setMipStart()}
\label{\detokenize{javaapi/Model:id97}}\begin{quote}

\sphinxAtStartPar
Set initial value for the specified variable.

\sphinxAtStartPar
\sphinxstylestrong{Synopsis}
\begin{quote}

\sphinxAtStartPar
\sphinxcode{\sphinxupquote{void setMipStart(Var{[}{]} vars, double{[}{]} vals)}}
\end{quote}

\sphinxAtStartPar
\sphinxstylestrong{Arguments}
\begin{quote}

\sphinxAtStartPar
\sphinxcode{\sphinxupquote{vars}}: a list of interested variables.

\sphinxAtStartPar
\sphinxcode{\sphinxupquote{vals}}: initial values of the variables.
\end{quote}
\end{quote}

\subsubsection{Model.setMipStart()}
\label{\detokenize{javaapi/Model:id98}}\begin{quote}

\sphinxAtStartPar
Set initial values for an array of variables.

\sphinxAtStartPar
\sphinxstylestrong{Synopsis}
\begin{quote}

\sphinxAtStartPar
\sphinxcode{\sphinxupquote{void setMipStart(VarArray vars, double{[}{]} vals)}}
\end{quote}

\sphinxAtStartPar
\sphinxstylestrong{Arguments}
\begin{quote}

\sphinxAtStartPar
\sphinxcode{\sphinxupquote{vars}}: a list of interested variables.

\sphinxAtStartPar
\sphinxcode{\sphinxupquote{vals}}: initial values of variables.
\end{quote}
\end{quote}

\subsubsection{Model.setNames()}
\label{\detokenize{javaapi/Model:model-setnames}}\begin{quote}

\sphinxAtStartPar
Set names for given variables in model.

\sphinxAtStartPar
\sphinxstylestrong{Synopsis}
\begin{quote}

\sphinxAtStartPar
\sphinxcode{\sphinxupquote{void setNames(Var{[}{]} vars, String{[}{]} names)}}
\end{quote}

\sphinxAtStartPar
\sphinxstylestrong{Arguments}
\begin{quote}

\sphinxAtStartPar
\sphinxcode{\sphinxupquote{vars}}: Array of variables.

\sphinxAtStartPar
\sphinxcode{\sphinxupquote{names}}: String array of names for variables.
\end{quote}
\end{quote}

\subsubsection{Model.setNames()}
\label{\detokenize{javaapi/Model:id99}}\begin{quote}

\sphinxAtStartPar
Set names for given variables in model.

\sphinxAtStartPar
\sphinxstylestrong{Synopsis}
\begin{quote}

\sphinxAtStartPar
\sphinxcode{\sphinxupquote{void setNames(VarArray vars, String{[}{]} names)}}
\end{quote}

\sphinxAtStartPar
\sphinxstylestrong{Arguments}
\begin{quote}

\sphinxAtStartPar
\sphinxcode{\sphinxupquote{vars}}: A list of variables.

\sphinxAtStartPar
\sphinxcode{\sphinxupquote{names}}: string array of names for variables.
\end{quote}
\end{quote}

\subsubsection{Model.setNames()}
\label{\detokenize{javaapi/Model:id100}}\begin{quote}

\sphinxAtStartPar
Set names for given constraints in model.

\sphinxAtStartPar
\sphinxstylestrong{Synopsis}
\begin{quote}

\sphinxAtStartPar
\sphinxcode{\sphinxupquote{void setNames(Constraint{[}{]} cons, String{[}{]} names)}}
\end{quote}

\sphinxAtStartPar
\sphinxstylestrong{Arguments}
\begin{quote}

\sphinxAtStartPar
\sphinxcode{\sphinxupquote{cons}}: array of constraints.

\sphinxAtStartPar
\sphinxcode{\sphinxupquote{names}}: string array of names for constraints.
\end{quote}
\end{quote}

\subsubsection{Model.setNames()}
\label{\detokenize{javaapi/Model:id101}}\begin{quote}

\sphinxAtStartPar
Set names for given constraints in model.

\sphinxAtStartPar
\sphinxstylestrong{Synopsis}
\begin{quote}

\sphinxAtStartPar
\sphinxcode{\sphinxupquote{void setNames(ConstrArray cons, String{[}{]} names)}}
\end{quote}

\sphinxAtStartPar
\sphinxstylestrong{Arguments}
\begin{quote}

\sphinxAtStartPar
\sphinxcode{\sphinxupquote{cons}}: a list of constraints.

\sphinxAtStartPar
\sphinxcode{\sphinxupquote{names}}: string array of names for constraints.
\end{quote}
\end{quote}

\subsubsection{Model.setNames()}
\label{\detokenize{javaapi/Model:id102}}\begin{quote}

\sphinxAtStartPar
Set names for given general constraints in model.

\sphinxAtStartPar
\sphinxstylestrong{Synopsis}
\begin{quote}

\sphinxAtStartPar
\sphinxcode{\sphinxupquote{void setNames(GenConstr{[}{]} genConstrs, String{[}{]} names)}}
\end{quote}

\sphinxAtStartPar
\sphinxstylestrong{Arguments}
\begin{quote}

\sphinxAtStartPar
\sphinxcode{\sphinxupquote{genConstrs}}: array of general constraints.

\sphinxAtStartPar
\sphinxcode{\sphinxupquote{names}}: string array of names for general constraints.
\end{quote}
\end{quote}

\subsubsection{Model.setNames()}
\label{\detokenize{javaapi/Model:id103}}\begin{quote}

\sphinxAtStartPar
Set names for given general constraints in model.

\sphinxAtStartPar
\sphinxstylestrong{Synopsis}
\begin{quote}

\sphinxAtStartPar
\sphinxcode{\sphinxupquote{void setNames(GenConstrArray genConstrs, String{[}{]} names)}}
\end{quote}

\sphinxAtStartPar
\sphinxstylestrong{Arguments}
\begin{quote}

\sphinxAtStartPar
\sphinxcode{\sphinxupquote{genConstrs}}: a list of general constraints.

\sphinxAtStartPar
\sphinxcode{\sphinxupquote{names}}: string array of names for general constraints.
\end{quote}
\end{quote}

\subsubsection{Model.setNames()}
\label{\detokenize{javaapi/Model:id104}}\begin{quote}

\sphinxAtStartPar
Set names for given nonlinear constraints in model.

\sphinxAtStartPar
\sphinxstylestrong{Synopsis}
\begin{quote}

\sphinxAtStartPar
\sphinxcode{\sphinxupquote{void setNames(NlConstraint{[}{]} cons, String{[}{]} names)}}
\end{quote}

\sphinxAtStartPar
\sphinxstylestrong{Arguments}
\begin{quote}

\sphinxAtStartPar
\sphinxcode{\sphinxupquote{cons}}: array of nonlinear constraints.

\sphinxAtStartPar
\sphinxcode{\sphinxupquote{names}}: string array of names for nonlinear constraints.
\end{quote}
\end{quote}

\subsubsection{Model.setNames()}
\label{\detokenize{javaapi/Model:id105}}\begin{quote}

\sphinxAtStartPar
Set names for given nonlinear constraints in model.

\sphinxAtStartPar
\sphinxstylestrong{Synopsis}
\begin{quote}

\sphinxAtStartPar
\sphinxcode{\sphinxupquote{void setNames(NlConstrArray cons, String{[}{]} names)}}
\end{quote}

\sphinxAtStartPar
\sphinxstylestrong{Arguments}
\begin{quote}

\sphinxAtStartPar
\sphinxcode{\sphinxupquote{cons}}: array object of nonlinear constraints.

\sphinxAtStartPar
\sphinxcode{\sphinxupquote{names}}: string array of names for nonlinear constraints.
\end{quote}
\end{quote}

\subsubsection{Model.setNames()}
\label{\detokenize{javaapi/Model:id106}}\begin{quote}

\sphinxAtStartPar
Set names for given quadratic constraints in model.

\sphinxAtStartPar
\sphinxstylestrong{Synopsis}
\begin{quote}

\sphinxAtStartPar
\sphinxcode{\sphinxupquote{void setNames(QConstraint{[}{]} cons, String{[}{]} names)}}
\end{quote}

\sphinxAtStartPar
\sphinxstylestrong{Arguments}
\begin{quote}

\sphinxAtStartPar
\sphinxcode{\sphinxupquote{cons}}: array of quadratic constraints.

\sphinxAtStartPar
\sphinxcode{\sphinxupquote{names}}: string array of names for quadratic constraints.
\end{quote}
\end{quote}

\subsubsection{Model.setNames()}
\label{\detokenize{javaapi/Model:id107}}\begin{quote}

\sphinxAtStartPar
Set names for given quadratic constraints in model.

\sphinxAtStartPar
\sphinxstylestrong{Synopsis}
\begin{quote}

\sphinxAtStartPar
\sphinxcode{\sphinxupquote{void setNames(QConstrArray cons, String{[}{]} names)}}
\end{quote}

\sphinxAtStartPar
\sphinxstylestrong{Arguments}
\begin{quote}

\sphinxAtStartPar
\sphinxcode{\sphinxupquote{cons}}: a list of quadratic constraints.

\sphinxAtStartPar
\sphinxcode{\sphinxupquote{names}}: string array of names for quadratic constraints.
\end{quote}
\end{quote}

\subsubsection{Model.setNames()}
\label{\detokenize{javaapi/Model:id108}}\begin{quote}

\sphinxAtStartPar
Set names for given PSD variables in model.

\sphinxAtStartPar
\sphinxstylestrong{Synopsis}
\begin{quote}

\sphinxAtStartPar
\sphinxcode{\sphinxupquote{void setNames(PsdVar{[}{]} vars, String{[}{]} names)}}
\end{quote}

\sphinxAtStartPar
\sphinxstylestrong{Arguments}
\begin{quote}

\sphinxAtStartPar
\sphinxcode{\sphinxupquote{vars}}: array of PSD variables.

\sphinxAtStartPar
\sphinxcode{\sphinxupquote{names}}: string array of names for PSD variables.
\end{quote}
\end{quote}

\subsubsection{Model.setNames()}
\label{\detokenize{javaapi/Model:id109}}\begin{quote}

\sphinxAtStartPar
Set names for given PSD variables in model.

\sphinxAtStartPar
\sphinxstylestrong{Synopsis}
\begin{quote}

\sphinxAtStartPar
\sphinxcode{\sphinxupquote{void setNames(PsdVarArray vars, String{[}{]} names)}}
\end{quote}

\sphinxAtStartPar
\sphinxstylestrong{Arguments}
\begin{quote}

\sphinxAtStartPar
\sphinxcode{\sphinxupquote{vars}}: a list of PSD variables.

\sphinxAtStartPar
\sphinxcode{\sphinxupquote{names}}: string array of names for PSD variables.
\end{quote}
\end{quote}

\subsubsection{Model.setNames()}
\label{\detokenize{javaapi/Model:id110}}\begin{quote}

\sphinxAtStartPar
Set names for given PSD constraints in model.

\sphinxAtStartPar
\sphinxstylestrong{Synopsis}
\begin{quote}

\sphinxAtStartPar
\sphinxcode{\sphinxupquote{void setNames(PsdConstraint{[}{]} cons, String{[}{]} names)}}
\end{quote}

\sphinxAtStartPar
\sphinxstylestrong{Arguments}
\begin{quote}

\sphinxAtStartPar
\sphinxcode{\sphinxupquote{cons}}: array of PSD constraints.

\sphinxAtStartPar
\sphinxcode{\sphinxupquote{names}}: string array of names for PSD constraints.
\end{quote}
\end{quote}

\subsubsection{Model.setNames()}
\label{\detokenize{javaapi/Model:id111}}\begin{quote}

\sphinxAtStartPar
Set names for given PSD constraints in model.

\sphinxAtStartPar
\sphinxstylestrong{Synopsis}
\begin{quote}

\sphinxAtStartPar
\sphinxcode{\sphinxupquote{void setNames(PsdConstrArray cons, String{[}{]} names)}}
\end{quote}

\sphinxAtStartPar
\sphinxstylestrong{Arguments}
\begin{quote}

\sphinxAtStartPar
\sphinxcode{\sphinxupquote{cons}}: a list of PSD constraints.

\sphinxAtStartPar
\sphinxcode{\sphinxupquote{names}}: string array of names for PSD constraints.
\end{quote}
\end{quote}

\subsubsection{Model.setNames()}
\label{\detokenize{javaapi/Model:id112}}\begin{quote}

\sphinxAtStartPar
Set names for given LMI constraints in model.

\sphinxAtStartPar
\sphinxstylestrong{Synopsis}
\begin{quote}

\sphinxAtStartPar
\sphinxcode{\sphinxupquote{void setNames(LmiConstraint{[}{]} cons, String{[}{]} names)}}
\end{quote}

\sphinxAtStartPar
\sphinxstylestrong{Arguments}
\begin{quote}

\sphinxAtStartPar
\sphinxcode{\sphinxupquote{cons}}: array of LMI constraints.

\sphinxAtStartPar
\sphinxcode{\sphinxupquote{names}}: string array of names for LMI constraints.
\end{quote}
\end{quote}

\subsubsection{Model.setNames()}
\label{\detokenize{javaapi/Model:id113}}\begin{quote}

\sphinxAtStartPar
Set names for given LMI constraints in model.

\sphinxAtStartPar
\sphinxstylestrong{Synopsis}
\begin{quote}

\sphinxAtStartPar
\sphinxcode{\sphinxupquote{void setNames(LmiConstrArray cons, String{[}{]} names)}}
\end{quote}

\sphinxAtStartPar
\sphinxstylestrong{Arguments}
\begin{quote}

\sphinxAtStartPar
\sphinxcode{\sphinxupquote{cons}}: a list of LMI constraints.

\sphinxAtStartPar
\sphinxcode{\sphinxupquote{names}}: string array of names for LMI constraints.
\end{quote}
\end{quote}

\subsubsection{Model.setNames()}
\label{\detokenize{javaapi/Model:id114}}\begin{quote}

\sphinxAtStartPar
Set names for given affine cone constraints in model.

\sphinxAtStartPar
\sphinxstylestrong{Synopsis}
\begin{quote}

\sphinxAtStartPar
\sphinxcode{\sphinxupquote{void setNames(AffineConeArray cones, String{[}{]} names)}}
\end{quote}

\sphinxAtStartPar
\sphinxstylestrong{Arguments}
\begin{quote}

\sphinxAtStartPar
\sphinxcode{\sphinxupquote{cones}}: an array of affine cone constraints.

\sphinxAtStartPar
\sphinxcode{\sphinxupquote{names}}: string array of names for affine cone constraints.
\end{quote}
\end{quote}

\subsubsection{Model.setNlObjective()}
\label{\detokenize{javaapi/Model:model-setnlobjective}}\begin{quote}

\sphinxAtStartPar
Set nonlinear objective for model.

\sphinxAtStartPar
\sphinxstylestrong{Synopsis}
\begin{quote}

\sphinxAtStartPar
\sphinxcode{\sphinxupquote{void setNlObjective(NlExpr expr, int sense)}}
\end{quote}

\sphinxAtStartPar
\sphinxstylestrong{Arguments}
\begin{quote}

\sphinxAtStartPar
\sphinxcode{\sphinxupquote{expr}}: nonlinear expression of the objective.

\sphinxAtStartPar
\sphinxcode{\sphinxupquote{sense}}: optimization sense. optional, default value 0 does not change COPT sense.
\end{quote}
\end{quote}

\subsubsection{Model.setNlPrimalStart()}
\label{\detokenize{javaapi/Model:model-setnlprimalstart}}\begin{quote}

\sphinxAtStartPar
Given count, set initial values for variables of NLP from beginning.

\sphinxAtStartPar
\sphinxstylestrong{Synopsis}
\begin{quote}

\sphinxAtStartPar
\sphinxcode{\sphinxupquote{void setNlPrimalStart(int count, double{[}{]} vals)}}
\end{quote}

\sphinxAtStartPar
\sphinxstylestrong{Arguments}
\begin{quote}

\sphinxAtStartPar
\sphinxcode{\sphinxupquote{count}}: the number of variables to set.

\sphinxAtStartPar
\sphinxcode{\sphinxupquote{vals}}: initial values of variables.
\end{quote}
\end{quote}

\subsubsection{Model.setNlPrimalStart()}
\label{\detokenize{javaapi/Model:id115}}\begin{quote}

\sphinxAtStartPar
Set initial value for the specified variable of NLP.

\sphinxAtStartPar
\sphinxstylestrong{Synopsis}
\begin{quote}

\sphinxAtStartPar
\sphinxcode{\sphinxupquote{void setNlPrimalStart(Var var, double val)}}
\end{quote}

\sphinxAtStartPar
\sphinxstylestrong{Arguments}
\begin{quote}

\sphinxAtStartPar
\sphinxcode{\sphinxupquote{var}}: an interested variable.

\sphinxAtStartPar
\sphinxcode{\sphinxupquote{val}}: initial value of the variable.
\end{quote}
\end{quote}

\subsubsection{Model.setNlPrimalStart()}
\label{\detokenize{javaapi/Model:id116}}\begin{quote}

\sphinxAtStartPar
Set initial values for an array of variables of NLP.

\sphinxAtStartPar
\sphinxstylestrong{Synopsis}
\begin{quote}

\sphinxAtStartPar
\sphinxcode{\sphinxupquote{void setNlPrimalStart(Var{[}{]} vars, double{[}{]} vals)}}
\end{quote}

\sphinxAtStartPar
\sphinxstylestrong{Arguments}
\begin{quote}

\sphinxAtStartPar
\sphinxcode{\sphinxupquote{vars}}: array of interested variables.

\sphinxAtStartPar
\sphinxcode{\sphinxupquote{vals}}: initial values of variables.
\end{quote}
\end{quote}

\subsubsection{Model.setNlPrimalStart()}
\label{\detokenize{javaapi/Model:id117}}\begin{quote}

\sphinxAtStartPar
Set initial values for variable array of NLP.

\sphinxAtStartPar
\sphinxstylestrong{Synopsis}
\begin{quote}

\sphinxAtStartPar
\sphinxcode{\sphinxupquote{void setNlPrimalStart(VarArray vars, double{[}{]} vals)}}
\end{quote}

\sphinxAtStartPar
\sphinxstylestrong{Arguments}
\begin{quote}

\sphinxAtStartPar
\sphinxcode{\sphinxupquote{vars}}: a list of interested variables.

\sphinxAtStartPar
\sphinxcode{\sphinxupquote{vals}}: initial values of variables.
\end{quote}
\end{quote}

\subsubsection{Model.setObjConst()}
\label{\detokenize{javaapi/Model:model-setobjconst}}\begin{quote}

\sphinxAtStartPar
Set objective constant.

\sphinxAtStartPar
\sphinxstylestrong{Synopsis}
\begin{quote}

\sphinxAtStartPar
\sphinxcode{\sphinxupquote{void setObjConst(double constant)}}
\end{quote}

\sphinxAtStartPar
\sphinxstylestrong{Arguments}
\begin{quote}

\sphinxAtStartPar
\sphinxcode{\sphinxupquote{constant}}: constant value to set.
\end{quote}
\end{quote}

\subsubsection{Model.setObjective()}
\label{\detokenize{javaapi/Model:model-setobjective}}\begin{quote}

\sphinxAtStartPar
Set objective for model.

\sphinxAtStartPar
\sphinxstylestrong{Synopsis}
\begin{quote}

\sphinxAtStartPar
\sphinxcode{\sphinxupquote{void setObjective(Expr expr, int sense)}}
\end{quote}

\sphinxAtStartPar
\sphinxstylestrong{Arguments}
\begin{quote}

\sphinxAtStartPar
\sphinxcode{\sphinxupquote{expr}}: expression of the objective.

\sphinxAtStartPar
\sphinxcode{\sphinxupquote{sense}}: optimization sense, which is either Consts.MINIIMIZE or Consts.MAXIMIZE. Set sense to 0 if do not change current sense.
\end{quote}
\end{quote}

\subsubsection{Model.setObjectiveN()}
\label{\detokenize{javaapi/Model:model-setobjectiven}}\begin{quote}

\sphinxAtStartPar
Set a multi\sphinxhyphen{}objective function in model.

\sphinxAtStartPar
\sphinxstylestrong{Synopsis}
\begin{quote}

\sphinxAtStartPar
\sphinxcode{\sphinxupquote{void setObjectiveN(}}
\begin{quote}

\sphinxAtStartPar
\sphinxcode{\sphinxupquote{int idx,}}

\sphinxAtStartPar
\sphinxcode{\sphinxupquote{Expr expr,}}

\sphinxAtStartPar
\sphinxcode{\sphinxupquote{int sense,}}

\sphinxAtStartPar
\sphinxcode{\sphinxupquote{double priority,}}

\sphinxAtStartPar
\sphinxcode{\sphinxupquote{double weight,}}

\sphinxAtStartPar
\sphinxcode{\sphinxupquote{double abstol,}}

\sphinxAtStartPar
\sphinxcode{\sphinxupquote{double reltol)}}
\end{quote}
\end{quote}

\sphinxAtStartPar
\sphinxstylestrong{Arguments}
\begin{quote}

\sphinxAtStartPar
\sphinxcode{\sphinxupquote{idx}}: index of a multi\sphinxhyphen{}objective function.

\sphinxAtStartPar
\sphinxcode{\sphinxupquote{expr}}: linear expression of the multi\sphinxhyphen{}objective function.

\sphinxAtStartPar
\sphinxcode{\sphinxupquote{sense}}: optimization sense with value Consts.MINIIMIZE or Consts.MAXIIMIZE. Set sense to 0 if do not change current sense.

\sphinxAtStartPar
\sphinxcode{\sphinxupquote{priority}}: objective parameter for priority.

\sphinxAtStartPar
\sphinxcode{\sphinxupquote{weight}}: objective parameter for weight.

\sphinxAtStartPar
\sphinxcode{\sphinxupquote{abstol}}: objective parameter for absolute tolerance.

\sphinxAtStartPar
\sphinxcode{\sphinxupquote{reltol}}: objective parameter for relative tolerance.
\end{quote}
\end{quote}

\subsubsection{Model.setObjectiveN()}
\label{\detokenize{javaapi/Model:id118}}\begin{quote}

\sphinxAtStartPar
Set a multi\sphinxhyphen{}objective function in model, with default value of objective parameters.

\sphinxAtStartPar
\sphinxstylestrong{Synopsis}
\begin{quote}

\sphinxAtStartPar
\sphinxcode{\sphinxupquote{void setObjectiveN(}}
\begin{quote}

\sphinxAtStartPar
\sphinxcode{\sphinxupquote{int idx,}}

\sphinxAtStartPar
\sphinxcode{\sphinxupquote{Expr expr,}}

\sphinxAtStartPar
\sphinxcode{\sphinxupquote{int sense)}}
\end{quote}
\end{quote}

\sphinxAtStartPar
\sphinxstylestrong{Arguments}
\begin{quote}

\sphinxAtStartPar
\sphinxcode{\sphinxupquote{idx}}: index of a multi\sphinxhyphen{}objective function.

\sphinxAtStartPar
\sphinxcode{\sphinxupquote{expr}}: linear expression of the multi\sphinxhyphen{}objective function.

\sphinxAtStartPar
\sphinxcode{\sphinxupquote{sense}}: optimization sense with value Consts.MINIIMIZE or Consts.MAXIIMIZE. Set sense to 0 if do not change current sense.
\end{quote}
\end{quote}

\subsubsection{Model.setObjParamN()}
\label{\detokenize{javaapi/Model:model-setobjparamn}}\begin{quote}

\sphinxAtStartPar
Set value of objective parameter of a multi\sphinxhyphen{}objective function.

\sphinxAtStartPar
\sphinxstylestrong{Synopsis}
\begin{quote}

\sphinxAtStartPar
\sphinxcode{\sphinxupquote{void setObjParamN(}}
\begin{quote}

\sphinxAtStartPar
\sphinxcode{\sphinxupquote{int idx,}}

\sphinxAtStartPar
\sphinxcode{\sphinxupquote{String param,}}

\sphinxAtStartPar
\sphinxcode{\sphinxupquote{double val)}}
\end{quote}
\end{quote}

\sphinxAtStartPar
\sphinxstylestrong{Arguments}
\begin{quote}

\sphinxAtStartPar
\sphinxcode{\sphinxupquote{idx}}: index of a multi\sphinxhyphen{}objective function.

\sphinxAtStartPar
\sphinxcode{\sphinxupquote{param}}: name of objective parameter, including priority, weight, abstol and reltol.

\sphinxAtStartPar
\sphinxcode{\sphinxupquote{val}}: new value of objective parameter.
\end{quote}
\end{quote}

\subsubsection{Model.setObjSense()}
\label{\detokenize{javaapi/Model:model-setobjsense}}\begin{quote}

\sphinxAtStartPar
Set objective sense for model.

\sphinxAtStartPar
\sphinxstylestrong{Synopsis}
\begin{quote}

\sphinxAtStartPar
\sphinxcode{\sphinxupquote{void setObjSense(int sense)}}
\end{quote}

\sphinxAtStartPar
\sphinxstylestrong{Arguments}
\begin{quote}

\sphinxAtStartPar
\sphinxcode{\sphinxupquote{sense}}: the objective sense.
\end{quote}
\end{quote}

\subsubsection{Model.setPsdCoeff()}
\label{\detokenize{javaapi/Model:model-setpsdcoeff}}\begin{quote}

\sphinxAtStartPar
Set the coefficient matrix of a PSD variable in a PSD constraint.

\sphinxAtStartPar
\sphinxstylestrong{Synopsis}
\begin{quote}

\sphinxAtStartPar
\sphinxcode{\sphinxupquote{void setPsdCoeff(}}
\begin{quote}

\sphinxAtStartPar
\sphinxcode{\sphinxupquote{PsdConstraint constr,}}

\sphinxAtStartPar
\sphinxcode{\sphinxupquote{PsdVar var,}}

\sphinxAtStartPar
\sphinxcode{\sphinxupquote{SymMatrix mat)}}
\end{quote}
\end{quote}

\sphinxAtStartPar
\sphinxstylestrong{Arguments}
\begin{quote}

\sphinxAtStartPar
\sphinxcode{\sphinxupquote{constr}}: The desired PSD constraint.

\sphinxAtStartPar
\sphinxcode{\sphinxupquote{var}}: The desired PSD variable.

\sphinxAtStartPar
\sphinxcode{\sphinxupquote{mat}}: new coefficient matrix.
\end{quote}
\end{quote}

\subsubsection{Model.setPsdObjective()}
\label{\detokenize{javaapi/Model:model-setpsdobjective}}\begin{quote}

\sphinxAtStartPar
Set PSD objective for model.

\sphinxAtStartPar
\sphinxstylestrong{Synopsis}
\begin{quote}

\sphinxAtStartPar
\sphinxcode{\sphinxupquote{void setPsdObjective(PsdExpr expr, int sense)}}
\end{quote}

\sphinxAtStartPar
\sphinxstylestrong{Arguments}
\begin{quote}

\sphinxAtStartPar
\sphinxcode{\sphinxupquote{expr}}: PSD expression of the objective.

\sphinxAtStartPar
\sphinxcode{\sphinxupquote{sense}}: optimization sense, which is either Consts.MINIIMIZE or Consts.MAXIMIZE. Set sense to 0 if do not change current sense.
\end{quote}
\end{quote}

\subsubsection{Model.setQuadObjective()}
\label{\detokenize{javaapi/Model:model-setquadobjective}}\begin{quote}

\sphinxAtStartPar
Set quadratic objective for model.

\sphinxAtStartPar
\sphinxstylestrong{Synopsis}
\begin{quote}

\sphinxAtStartPar
\sphinxcode{\sphinxupquote{void setQuadObjective(QuadExpr expr, int sense)}}
\end{quote}

\sphinxAtStartPar
\sphinxstylestrong{Arguments}
\begin{quote}

\sphinxAtStartPar
\sphinxcode{\sphinxupquote{expr}}: quadratic expression of the objective.

\sphinxAtStartPar
\sphinxcode{\sphinxupquote{sense}}: optimization sense, which is either Consts.MINIIMIZE or Consts.MAXIMIZE. Set sense to 0 if do not change current sense.
\end{quote}
\end{quote}

\subsubsection{Model.setSlackBasis()}
\label{\detokenize{javaapi/Model:model-setslackbasis}}\begin{quote}

\sphinxAtStartPar
Set slack basis to model.

\sphinxAtStartPar
\sphinxstylestrong{Synopsis}
\begin{quote}

\sphinxAtStartPar
\sphinxcode{\sphinxupquote{void setSlackBasis()}}
\end{quote}
\end{quote}

\subsubsection{Model.setSolverLogFile()}
\label{\detokenize{javaapi/Model:model-setsolverlogfile}}\begin{quote}

\sphinxAtStartPar
Set log file for COPT.

\sphinxAtStartPar
\sphinxstylestrong{Synopsis}
\begin{quote}

\sphinxAtStartPar
\sphinxcode{\sphinxupquote{void setSolverLogFile(String filename)}}
\end{quote}

\sphinxAtStartPar
\sphinxstylestrong{Arguments}
\begin{quote}

\sphinxAtStartPar
\sphinxcode{\sphinxupquote{filename}}: log file name.
\end{quote}
\end{quote}

\subsubsection{Model.solve()}
\label{\detokenize{javaapi/Model:model-solve}}\begin{quote}

\sphinxAtStartPar
Solve the model as MIP.

\sphinxAtStartPar
\sphinxstylestrong{Synopsis}
\begin{quote}

\sphinxAtStartPar
\sphinxcode{\sphinxupquote{void solve()}}
\end{quote}
\end{quote}

\subsubsection{Model.solveLp()}
\label{\detokenize{javaapi/Model:model-solvelp}}\begin{quote}

\sphinxAtStartPar
Solve the model as LP.

\sphinxAtStartPar
\sphinxstylestrong{Synopsis}
\begin{quote}

\sphinxAtStartPar
\sphinxcode{\sphinxupquote{void solveLp()}}
\end{quote}
\end{quote}

\subsubsection{Model.tune()}
\label{\detokenize{javaapi/Model:model-tune}}\begin{quote}

\sphinxAtStartPar
Tune model.

\sphinxAtStartPar
\sphinxstylestrong{Synopsis}
\begin{quote}

\sphinxAtStartPar
\sphinxcode{\sphinxupquote{void tune()}}
\end{quote}
\end{quote}

\subsubsection{Model.write()}
\label{\detokenize{javaapi/Model:model-write}}\begin{quote}

\sphinxAtStartPar
Output problem, solution, basis, MIP start or modified COPT parameters to file.

\sphinxAtStartPar
\sphinxstylestrong{Synopsis}
\begin{quote}

\sphinxAtStartPar
\sphinxcode{\sphinxupquote{void write(String filename)}}
\end{quote}

\sphinxAtStartPar
\sphinxstylestrong{Arguments}
\begin{quote}

\sphinxAtStartPar
\sphinxcode{\sphinxupquote{filename}}: an output file name.
\end{quote}
\end{quote}

\subsubsection{Model.writeBasis()}
\label{\detokenize{javaapi/Model:model-writebasis}}\begin{quote}

\sphinxAtStartPar
Output optimal basis to a file of type ‘.bas’.

\sphinxAtStartPar
\sphinxstylestrong{Synopsis}
\begin{quote}

\sphinxAtStartPar
\sphinxcode{\sphinxupquote{void writeBasis(String filename)}}
\end{quote}

\sphinxAtStartPar
\sphinxstylestrong{Arguments}
\begin{quote}

\sphinxAtStartPar
\sphinxcode{\sphinxupquote{filename}}: an output file name.
\end{quote}
\end{quote}

\subsubsection{Model.writeBin()}
\label{\detokenize{javaapi/Model:model-writebin}}\begin{quote}

\sphinxAtStartPar
Output problem to a file as COPT binary format.

\sphinxAtStartPar
\sphinxstylestrong{Synopsis}
\begin{quote}

\sphinxAtStartPar
\sphinxcode{\sphinxupquote{void writeBin(String filename)}}
\end{quote}

\sphinxAtStartPar
\sphinxstylestrong{Arguments}
\begin{quote}

\sphinxAtStartPar
\sphinxcode{\sphinxupquote{filename}}: an output file name.
\end{quote}
\end{quote}

\subsubsection{Model.writeIIS()}
\label{\detokenize{javaapi/Model:model-writeiis}}\begin{quote}

\sphinxAtStartPar
Output IIS to file.

\sphinxAtStartPar
\sphinxstylestrong{Synopsis}
\begin{quote}

\sphinxAtStartPar
\sphinxcode{\sphinxupquote{void writeIIS(String filename)}}
\end{quote}

\sphinxAtStartPar
\sphinxstylestrong{Arguments}
\begin{quote}

\sphinxAtStartPar
\sphinxcode{\sphinxupquote{filename}}: Output file name.
\end{quote}
\end{quote}

\subsubsection{Model.writeJsonSol()}
\label{\detokenize{javaapi/Model:model-writejsonsol}}\begin{quote}

\sphinxAtStartPar
Output solution to a file of type ‘.json’.

\sphinxAtStartPar
\sphinxstylestrong{Synopsis}
\begin{quote}

\sphinxAtStartPar
\sphinxcode{\sphinxupquote{void writeJsonSol(String filename)}}
\end{quote}

\sphinxAtStartPar
\sphinxstylestrong{Arguments}
\begin{quote}

\sphinxAtStartPar
\sphinxcode{\sphinxupquote{filename}}: an output file name.
\end{quote}
\end{quote}

\subsubsection{Model.writeLp()}
\label{\detokenize{javaapi/Model:model-writelp}}\begin{quote}

\sphinxAtStartPar
Output problem to a file as LP format.

\sphinxAtStartPar
\sphinxstylestrong{Synopsis}
\begin{quote}

\sphinxAtStartPar
\sphinxcode{\sphinxupquote{void writeLp(String filename)}}
\end{quote}

\sphinxAtStartPar
\sphinxstylestrong{Arguments}
\begin{quote}

\sphinxAtStartPar
\sphinxcode{\sphinxupquote{filename}}: an output file name.
\end{quote}
\end{quote}

\subsubsection{Model.writeMps()}
\label{\detokenize{javaapi/Model:model-writemps}}\begin{quote}

\sphinxAtStartPar
Output problem to a file as MPS format.

\sphinxAtStartPar
\sphinxstylestrong{Synopsis}
\begin{quote}

\sphinxAtStartPar
\sphinxcode{\sphinxupquote{void writeMps(String filename)}}
\end{quote}

\sphinxAtStartPar
\sphinxstylestrong{Arguments}
\begin{quote}

\sphinxAtStartPar
\sphinxcode{\sphinxupquote{filename}}: an output file name.
\end{quote}
\end{quote}

\subsubsection{Model.writeMpsStr()}
\label{\detokenize{javaapi/Model:model-writempsstr}}\begin{quote}

\sphinxAtStartPar
Output MPS problem to problem buffer.

\sphinxAtStartPar
\sphinxstylestrong{Synopsis}
\begin{quote}

\sphinxAtStartPar
\sphinxcode{\sphinxupquote{ProbBuffer writeMpsStr()}}
\end{quote}

\sphinxAtStartPar
\sphinxstylestrong{Return}
\begin{quote}

\sphinxAtStartPar
problem buffer for string of MPS problem.
\end{quote}
\end{quote}

\subsubsection{Model.writeMst()}
\label{\detokenize{javaapi/Model:model-writemst}}\begin{quote}

\sphinxAtStartPar
Output MIP start information to a file of type ‘.mst’.

\sphinxAtStartPar
\sphinxstylestrong{Synopsis}
\begin{quote}

\sphinxAtStartPar
\sphinxcode{\sphinxupquote{void writeMst(String filename)}}
\end{quote}

\sphinxAtStartPar
\sphinxstylestrong{Arguments}
\begin{quote}

\sphinxAtStartPar
\sphinxcode{\sphinxupquote{filename}}: an output file name.
\end{quote}
\end{quote}

\subsubsection{Model.writeNL()}
\label{\detokenize{javaapi/Model:model-writenl}}\begin{quote}

\sphinxAtStartPar
Output problem to a file as NL format.

\sphinxAtStartPar
\sphinxstylestrong{Synopsis}
\begin{quote}

\sphinxAtStartPar
\sphinxcode{\sphinxupquote{void writeNL(String filename)}}
\end{quote}

\sphinxAtStartPar
\sphinxstylestrong{Arguments}
\begin{quote}

\sphinxAtStartPar
\sphinxcode{\sphinxupquote{filename}}: an output file name.
\end{quote}
\end{quote}

\subsubsection{Model.writeOrd()}
\label{\detokenize{javaapi/Model:model-writeord}}\begin{quote}

\sphinxAtStartPar
Output branching order to file.

\sphinxAtStartPar
\sphinxstylestrong{Synopsis}
\begin{quote}

\sphinxAtStartPar
\sphinxcode{\sphinxupquote{void writeOrd(String filename)}}
\end{quote}

\sphinxAtStartPar
\sphinxstylestrong{Arguments}
\begin{quote}

\sphinxAtStartPar
\sphinxcode{\sphinxupquote{filename}}: Output file name.
\end{quote}
\end{quote}

\subsubsection{Model.writeParam()}
\label{\detokenize{javaapi/Model:model-writeparam}}\begin{quote}

\sphinxAtStartPar
Output modified COPT parameters to a file of type ‘.par’.

\sphinxAtStartPar
\sphinxstylestrong{Synopsis}
\begin{quote}

\sphinxAtStartPar
\sphinxcode{\sphinxupquote{void writeParam(String filename)}}
\end{quote}

\sphinxAtStartPar
\sphinxstylestrong{Arguments}
\begin{quote}

\sphinxAtStartPar
\sphinxcode{\sphinxupquote{filename}}: an output file name.
\end{quote}
\end{quote}

\subsubsection{Model.writePoolSol()}
\label{\detokenize{javaapi/Model:model-writepoolsol}}\begin{quote}

\sphinxAtStartPar
Output selected pool solution to a file of type ‘.sol’.

\sphinxAtStartPar
\sphinxstylestrong{Synopsis}
\begin{quote}

\sphinxAtStartPar
\sphinxcode{\sphinxupquote{void writePoolSol(int idx, String filename)}}
\end{quote}

\sphinxAtStartPar
\sphinxstylestrong{Arguments}
\begin{quote}

\sphinxAtStartPar
\sphinxcode{\sphinxupquote{idx}}: index of pool solution.

\sphinxAtStartPar
\sphinxcode{\sphinxupquote{filename}}: an output file name.
\end{quote}
\end{quote}

\subsubsection{Model.writeRelax()}
\label{\detokenize{javaapi/Model:model-writerelax}}\begin{quote}

\sphinxAtStartPar
Output feasibility relaxation problem to file.

\sphinxAtStartPar
\sphinxstylestrong{Synopsis}
\begin{quote}

\sphinxAtStartPar
\sphinxcode{\sphinxupquote{void writeRelax(String filename)}}
\end{quote}

\sphinxAtStartPar
\sphinxstylestrong{Arguments}
\begin{quote}

\sphinxAtStartPar
\sphinxcode{\sphinxupquote{filename}}: Output file name.
\end{quote}
\end{quote}

\subsubsection{Model.writeSol()}
\label{\detokenize{javaapi/Model:model-writesol}}\begin{quote}

\sphinxAtStartPar
Output solution to a file of type ‘.sol’.

\sphinxAtStartPar
\sphinxstylestrong{Synopsis}
\begin{quote}

\sphinxAtStartPar
\sphinxcode{\sphinxupquote{void writeSol(String filename)}}
\end{quote}

\sphinxAtStartPar
\sphinxstylestrong{Arguments}
\begin{quote}

\sphinxAtStartPar
\sphinxcode{\sphinxupquote{filename}}: an output file name.
\end{quote}
\end{quote}

\subsubsection{Model.writeTuneParam()}
\label{\detokenize{javaapi/Model:model-writetuneparam}}\begin{quote}

\sphinxAtStartPar
Output specified tuned parameters to a file of type ‘.par’.

\sphinxAtStartPar
\sphinxstylestrong{Synopsis}
\begin{quote}

\sphinxAtStartPar
\sphinxcode{\sphinxupquote{void writeTuneParam(int idx, String filename)}}
\end{quote}

\sphinxAtStartPar
\sphinxstylestrong{Arguments}
\begin{quote}

\sphinxAtStartPar
\sphinxcode{\sphinxupquote{idx}}: Index of tuned parameters.

\sphinxAtStartPar
\sphinxcode{\sphinxupquote{filename}}: Output file name.
\end{quote}
\end{quote}

\subsection{Var}
\label{\detokenize{javaapiref:var}}\label{\detokenize{javaapiref:chapjavaapiref-var}}
\sphinxAtStartPar
COPT variable object. Variables are always associated with a particular model.
User creates a variable object by adding a variable to a model, rather than
by using constructor of Var class.

\sphinxstepscope

\subsubsection{Var.get()}
\label{\detokenize{javaapi/Var:var-get}}\label{\detokenize{javaapi/Var::doc}}\begin{quote}

\sphinxAtStartPar
Get information value of the variable. Support informations of “Value”, “RedCost”, “PrimalRay”, “LB”, “UB”, “Obj” and “BranchFactor”.

\sphinxAtStartPar
\sphinxstylestrong{Synopsis}
\begin{quote}

\sphinxAtStartPar
\sphinxcode{\sphinxupquote{double get(String info)}}
\end{quote}

\sphinxAtStartPar
\sphinxstylestrong{Arguments}
\begin{quote}

\sphinxAtStartPar
\sphinxcode{\sphinxupquote{info}}: information name.
\end{quote}

\sphinxAtStartPar
\sphinxstylestrong{Return}
\begin{quote}

\sphinxAtStartPar
information value.
\end{quote}
\end{quote}

\subsubsection{Var.getBasis()}
\label{\detokenize{javaapi/Var:var-getbasis}}\begin{quote}

\sphinxAtStartPar
Get basis status of the variable.

\sphinxAtStartPar
\sphinxstylestrong{Synopsis}
\begin{quote}

\sphinxAtStartPar
\sphinxcode{\sphinxupquote{int getBasis()}}
\end{quote}

\sphinxAtStartPar
\sphinxstylestrong{Return}
\begin{quote}

\sphinxAtStartPar
Basis status.
\end{quote}
\end{quote}

\subsubsection{Var.getIdx()}
\label{\detokenize{javaapi/Var:var-getidx}}\begin{quote}

\sphinxAtStartPar
Get index of the variable.

\sphinxAtStartPar
\sphinxstylestrong{Synopsis}
\begin{quote}

\sphinxAtStartPar
\sphinxcode{\sphinxupquote{int getIdx()}}
\end{quote}

\sphinxAtStartPar
\sphinxstylestrong{Return}
\begin{quote}

\sphinxAtStartPar
variable index.
\end{quote}
\end{quote}

\subsubsection{Var.getLowerIIS()}
\label{\detokenize{javaapi/Var:var-getloweriis}}\begin{quote}

\sphinxAtStartPar
Get IIS status for lower bound of the variable.

\sphinxAtStartPar
\sphinxstylestrong{Synopsis}
\begin{quote}

\sphinxAtStartPar
\sphinxcode{\sphinxupquote{int getLowerIIS()}}
\end{quote}

\sphinxAtStartPar
\sphinxstylestrong{Return}
\begin{quote}

\sphinxAtStartPar
IIS status.
\end{quote}
\end{quote}

\subsubsection{Var.getName()}
\label{\detokenize{javaapi/Var:var-getname}}\begin{quote}

\sphinxAtStartPar
Get name of the variable.

\sphinxAtStartPar
\sphinxstylestrong{Synopsis}
\begin{quote}

\sphinxAtStartPar
\sphinxcode{\sphinxupquote{String getName()}}
\end{quote}

\sphinxAtStartPar
\sphinxstylestrong{Return}
\begin{quote}

\sphinxAtStartPar
variable name.
\end{quote}
\end{quote}

\subsubsection{Var.getType()}
\label{\detokenize{javaapi/Var:var-gettype}}\begin{quote}

\sphinxAtStartPar
Get type of the variable.

\sphinxAtStartPar
\sphinxstylestrong{Synopsis}
\begin{quote}

\sphinxAtStartPar
\sphinxcode{\sphinxupquote{char getType()}}
\end{quote}

\sphinxAtStartPar
\sphinxstylestrong{Return}
\begin{quote}

\sphinxAtStartPar
variable type.
\end{quote}
\end{quote}

\subsubsection{Var.getUpperIIS()}
\label{\detokenize{javaapi/Var:var-getupperiis}}\begin{quote}

\sphinxAtStartPar
Get IIS status for upper bound of the variable.

\sphinxAtStartPar
\sphinxstylestrong{Synopsis}
\begin{quote}

\sphinxAtStartPar
\sphinxcode{\sphinxupquote{int getUpperIIS()}}
\end{quote}

\sphinxAtStartPar
\sphinxstylestrong{Return}
\begin{quote}

\sphinxAtStartPar
IIS status.
\end{quote}
\end{quote}

\subsubsection{Var.remove()}
\label{\detokenize{javaapi/Var:var-remove}}\begin{quote}

\sphinxAtStartPar
Remove variable from model.

\sphinxAtStartPar
\sphinxstylestrong{Synopsis}
\begin{quote}

\sphinxAtStartPar
\sphinxcode{\sphinxupquote{void remove()}}
\end{quote}
\end{quote}

\subsubsection{Var.set()}
\label{\detokenize{javaapi/Var:var-set}}\begin{quote}

\sphinxAtStartPar
Set information value of the variable. Support informations of “LB”, “UB”, “Obj” and “BranchFactor”.

\sphinxAtStartPar
\sphinxstylestrong{Synopsis}
\begin{quote}

\sphinxAtStartPar
\sphinxcode{\sphinxupquote{void set(String info, double val)}}
\end{quote}

\sphinxAtStartPar
\sphinxstylestrong{Arguments}
\begin{quote}

\sphinxAtStartPar
\sphinxcode{\sphinxupquote{info}}: information name.

\sphinxAtStartPar
\sphinxcode{\sphinxupquote{val}}: new information value.
\end{quote}
\end{quote}

\subsubsection{Var.setName()}
\label{\detokenize{javaapi/Var:var-setname}}\begin{quote}

\sphinxAtStartPar
Set name of the variable.

\sphinxAtStartPar
\sphinxstylestrong{Synopsis}
\begin{quote}

\sphinxAtStartPar
\sphinxcode{\sphinxupquote{void setName(String name)}}
\end{quote}

\sphinxAtStartPar
\sphinxstylestrong{Arguments}
\begin{quote}

\sphinxAtStartPar
\sphinxcode{\sphinxupquote{name}}: variable name.
\end{quote}
\end{quote}

\subsubsection{Var.setType()}
\label{\detokenize{javaapi/Var:var-settype}}\begin{quote}

\sphinxAtStartPar
Set type of the variable.

\sphinxAtStartPar
\sphinxstylestrong{Synopsis}
\begin{quote}

\sphinxAtStartPar
\sphinxcode{\sphinxupquote{void setType(char vtype)}}
\end{quote}

\sphinxAtStartPar
\sphinxstylestrong{Arguments}
\begin{quote}

\sphinxAtStartPar
\sphinxcode{\sphinxupquote{vtype}}: variable type.
\end{quote}
\end{quote}

\subsection{VarArray}
\label{\detokenize{javaapiref:vararray}}\label{\detokenize{javaapiref:chapjavaapiref-vararray}}
\sphinxAtStartPar
COPT variable array object. To store and access a set of Java {\hyperref[\detokenize{javaapiref:chapjavaapiref-var}]{\sphinxcrossref{\DUrole{std,std-ref}{Var}}}} objects,
Cardinal Optimizer provides Java VarArray class, which defines the following methods.

\sphinxstepscope

\subsubsection{VarArray.VarArray()}
\label{\detokenize{javaapi/VarArray:vararray-vararray}}\label{\detokenize{javaapi/VarArray::doc}}\begin{quote}

\sphinxAtStartPar
Constructor of vararray.

\sphinxAtStartPar
\sphinxstylestrong{Synopsis}
\begin{quote}

\sphinxAtStartPar
\sphinxcode{\sphinxupquote{VarArray()}}
\end{quote}
\end{quote}

\subsubsection{VarArray.getVar()}
\label{\detokenize{javaapi/VarArray:vararray-getvar}}\begin{quote}

\sphinxAtStartPar
Get idx\sphinxhyphen{}th variable object.

\sphinxAtStartPar
\sphinxstylestrong{Synopsis}
\begin{quote}

\sphinxAtStartPar
\sphinxcode{\sphinxupquote{Var getVar(int idx)}}
\end{quote}

\sphinxAtStartPar
\sphinxstylestrong{Arguments}
\begin{quote}

\sphinxAtStartPar
\sphinxcode{\sphinxupquote{idx}}: index of the variable.
\end{quote}

\sphinxAtStartPar
\sphinxstylestrong{Return}
\begin{quote}

\sphinxAtStartPar
variable object with index idx.
\end{quote}
\end{quote}

\subsubsection{VarArray.pushBack()}
\label{\detokenize{javaapi/VarArray:vararray-pushback}}\begin{quote}

\sphinxAtStartPar
Add a variable object to variable array.

\sphinxAtStartPar
\sphinxstylestrong{Synopsis}
\begin{quote}

\sphinxAtStartPar
\sphinxcode{\sphinxupquote{void pushBack(Var var)}}
\end{quote}

\sphinxAtStartPar
\sphinxstylestrong{Arguments}
\begin{quote}

\sphinxAtStartPar
\sphinxcode{\sphinxupquote{var}}: a variable object.
\end{quote}
\end{quote}

\subsubsection{VarArray.size()}
\label{\detokenize{javaapi/VarArray:vararray-size}}\begin{quote}

\sphinxAtStartPar
Get the number of variable objects.

\sphinxAtStartPar
\sphinxstylestrong{Synopsis}
\begin{quote}

\sphinxAtStartPar
\sphinxcode{\sphinxupquote{int size()}}
\end{quote}

\sphinxAtStartPar
\sphinxstylestrong{Return}
\begin{quote}

\sphinxAtStartPar
number of variable objects.
\end{quote}
\end{quote}

\subsection{Expr}
\label{\detokenize{javaapiref:expr}}\label{\detokenize{javaapiref:chapjavaapiref-expr}}
\sphinxAtStartPar
COPT linear expression object. A linear expression consists of a constant term, a list of
terms of variables and associated coefficients. Linear expressions are used to build
constraints.

\sphinxstepscope

\subsubsection{Expr.Expr()}
\label{\detokenize{javaapi/Expr:expr-expr}}\label{\detokenize{javaapi/Expr::doc}}\begin{quote}

\sphinxAtStartPar
Constructor of a constant linear expression with constant 0.0

\sphinxAtStartPar
\sphinxstylestrong{Synopsis}
\begin{quote}

\sphinxAtStartPar
\sphinxcode{\sphinxupquote{Expr()}}
\end{quote}
\end{quote}

\subsubsection{Expr.Expr()}
\label{\detokenize{javaapi/Expr:id1}}\begin{quote}

\sphinxAtStartPar
Constructor of a constant linear expression.

\sphinxAtStartPar
\sphinxstylestrong{Synopsis}
\begin{quote}

\sphinxAtStartPar
\sphinxcode{\sphinxupquote{Expr(double constant)}}
\end{quote}

\sphinxAtStartPar
\sphinxstylestrong{Arguments}
\begin{quote}

\sphinxAtStartPar
\sphinxcode{\sphinxupquote{constant}}: constant value in expression object.
\end{quote}
\end{quote}

\subsubsection{Expr.Expr()}
\label{\detokenize{javaapi/Expr:id2}}\begin{quote}

\sphinxAtStartPar
Constructor of a linear expression with one term.

\sphinxAtStartPar
\sphinxstylestrong{Synopsis}
\begin{quote}

\sphinxAtStartPar
\sphinxcode{\sphinxupquote{Expr(Var var)}}
\end{quote}

\sphinxAtStartPar
\sphinxstylestrong{Arguments}
\begin{quote}

\sphinxAtStartPar
\sphinxcode{\sphinxupquote{var}}: variable for the added term.
\end{quote}
\end{quote}

\subsubsection{Expr.Expr()}
\label{\detokenize{javaapi/Expr:id3}}\begin{quote}

\sphinxAtStartPar
Constructor of a linear expression with one term.

\sphinxAtStartPar
\sphinxstylestrong{Synopsis}
\begin{quote}

\sphinxAtStartPar
\sphinxcode{\sphinxupquote{Expr(Var var, double coeff)}}
\end{quote}

\sphinxAtStartPar
\sphinxstylestrong{Arguments}
\begin{quote}

\sphinxAtStartPar
\sphinxcode{\sphinxupquote{var}}: variable for the added term.

\sphinxAtStartPar
\sphinxcode{\sphinxupquote{coeff}}: coefficent for the added term.
\end{quote}
\end{quote}

\subsubsection{Expr.add()}
\label{\detokenize{javaapi/Expr:expr-add}}\begin{quote}

\sphinxAtStartPar
Add itself by a linear expression.

\sphinxAtStartPar
\sphinxstylestrong{Synopsis}
\begin{quote}

\sphinxAtStartPar
\sphinxcode{\sphinxupquote{Expr add(Expr expr, double mult)}}
\end{quote}

\sphinxAtStartPar
\sphinxstylestrong{Arguments}
\begin{quote}

\sphinxAtStartPar
\sphinxcode{\sphinxupquote{expr}}: expression operand, including Expr, Var.

\sphinxAtStartPar
\sphinxcode{\sphinxupquote{mult}}: constant multiplier.
\end{quote}

\sphinxAtStartPar
\sphinxstylestrong{Return}
\begin{quote}

\sphinxAtStartPar
linear expression itself.
\end{quote}
\end{quote}

\subsubsection{Expr.add()}
\label{\detokenize{javaapi/Expr:id4}}\begin{quote}

\sphinxAtStartPar
Add itself by a linear expression.

\sphinxAtStartPar
\sphinxstylestrong{Synopsis}
\begin{quote}

\sphinxAtStartPar
\sphinxcode{\sphinxupquote{Expr add(Expr expr)}}
\end{quote}

\sphinxAtStartPar
\sphinxstylestrong{Arguments}
\begin{quote}

\sphinxAtStartPar
\sphinxcode{\sphinxupquote{expr}}: expression operand, including Expr, Var and constant.
\end{quote}

\sphinxAtStartPar
\sphinxstylestrong{Return}
\begin{quote}

\sphinxAtStartPar
linear expression itself.
\end{quote}
\end{quote}

\subsubsection{Expr.addConstant()}
\label{\detokenize{javaapi/Expr:expr-addconstant}}\begin{quote}

\sphinxAtStartPar
Add extra constant to the expression.

\sphinxAtStartPar
\sphinxstylestrong{Synopsis}
\begin{quote}

\sphinxAtStartPar
\sphinxcode{\sphinxupquote{void addConstant(double constant)}}
\end{quote}

\sphinxAtStartPar
\sphinxstylestrong{Arguments}
\begin{quote}

\sphinxAtStartPar
\sphinxcode{\sphinxupquote{constant}}: delta value to be added to expression constant.
\end{quote}
\end{quote}

\subsubsection{Expr.addExpr()}
\label{\detokenize{javaapi/Expr:expr-addexpr}}\begin{quote}

\sphinxAtStartPar
Add a linear expression to self.

\sphinxAtStartPar
\sphinxstylestrong{Synopsis}
\begin{quote}

\sphinxAtStartPar
\sphinxcode{\sphinxupquote{void addExpr(Expr expr)}}
\end{quote}

\sphinxAtStartPar
\sphinxstylestrong{Arguments}
\begin{quote}

\sphinxAtStartPar
\sphinxcode{\sphinxupquote{expr}}: linear expression to be added.
\end{quote}
\end{quote}

\subsubsection{Expr.addExpr()}
\label{\detokenize{javaapi/Expr:id5}}\begin{quote}

\sphinxAtStartPar
Add a linear expression to self.

\sphinxAtStartPar
\sphinxstylestrong{Synopsis}
\begin{quote}

\sphinxAtStartPar
\sphinxcode{\sphinxupquote{void addExpr(Expr expr, double mult)}}
\end{quote}

\sphinxAtStartPar
\sphinxstylestrong{Arguments}
\begin{quote}

\sphinxAtStartPar
\sphinxcode{\sphinxupquote{expr}}: linear expression to be added.

\sphinxAtStartPar
\sphinxcode{\sphinxupquote{mult}}: multiplier constant.
\end{quote}
\end{quote}

\subsubsection{Expr.addTerm()}
\label{\detokenize{javaapi/Expr:expr-addterm}}\begin{quote}

\sphinxAtStartPar
Add a term to expression object.

\sphinxAtStartPar
\sphinxstylestrong{Synopsis}
\begin{quote}

\sphinxAtStartPar
\sphinxcode{\sphinxupquote{void addTerm(Var var, double coeff)}}
\end{quote}

\sphinxAtStartPar
\sphinxstylestrong{Arguments}
\begin{quote}

\sphinxAtStartPar
\sphinxcode{\sphinxupquote{var}}: a variable for new term.

\sphinxAtStartPar
\sphinxcode{\sphinxupquote{coeff}}: coefficient for new term.
\end{quote}
\end{quote}

\subsubsection{Expr.addTerms()}
\label{\detokenize{javaapi/Expr:expr-addterms}}\begin{quote}

\sphinxAtStartPar
Add terms to expression object.

\sphinxAtStartPar
\sphinxstylestrong{Synopsis}
\begin{quote}

\sphinxAtStartPar
\sphinxcode{\sphinxupquote{void addTerms(Var{[}{]} vars, double coeff)}}
\end{quote}

\sphinxAtStartPar
\sphinxstylestrong{Arguments}
\begin{quote}

\sphinxAtStartPar
\sphinxcode{\sphinxupquote{vars}}: variables for added terms.

\sphinxAtStartPar
\sphinxcode{\sphinxupquote{coeff}}: one coefficient for added terms.
\end{quote}
\end{quote}

\subsubsection{Expr.addTerms()}
\label{\detokenize{javaapi/Expr:id6}}\begin{quote}

\sphinxAtStartPar
Add terms to expression object.

\sphinxAtStartPar
\sphinxstylestrong{Synopsis}
\begin{quote}

\sphinxAtStartPar
\sphinxcode{\sphinxupquote{void addTerms(Var{[}{]} vars, double{[}{]} coeffs)}}
\end{quote}

\sphinxAtStartPar
\sphinxstylestrong{Arguments}
\begin{quote}

\sphinxAtStartPar
\sphinxcode{\sphinxupquote{vars}}: variables for added terms.

\sphinxAtStartPar
\sphinxcode{\sphinxupquote{coeffs}}: coefficients array for added terms.
\end{quote}
\end{quote}

\subsubsection{Expr.addTerms()}
\label{\detokenize{javaapi/Expr:id7}}\begin{quote}

\sphinxAtStartPar
Add terms to expression object.

\sphinxAtStartPar
\sphinxstylestrong{Synopsis}
\begin{quote}

\sphinxAtStartPar
\sphinxcode{\sphinxupquote{void addTerms(VarArray vars, double coeff)}}
\end{quote}

\sphinxAtStartPar
\sphinxstylestrong{Arguments}
\begin{quote}

\sphinxAtStartPar
\sphinxcode{\sphinxupquote{vars}}: variables for added terms.

\sphinxAtStartPar
\sphinxcode{\sphinxupquote{coeff}}: one coefficient for added terms.
\end{quote}
\end{quote}

\subsubsection{Expr.addTerms()}
\label{\detokenize{javaapi/Expr:id8}}\begin{quote}

\sphinxAtStartPar
Add terms to expression object.

\sphinxAtStartPar
\sphinxstylestrong{Synopsis}
\begin{quote}

\sphinxAtStartPar
\sphinxcode{\sphinxupquote{void addTerms(VarArray vars, double{[}{]} coeffs)}}
\end{quote}

\sphinxAtStartPar
\sphinxstylestrong{Arguments}
\begin{quote}

\sphinxAtStartPar
\sphinxcode{\sphinxupquote{vars}}: variables for added terms.

\sphinxAtStartPar
\sphinxcode{\sphinxupquote{coeffs}}: coefficients array for added terms.
\end{quote}
\end{quote}

\subsubsection{Expr.clone()}
\label{\detokenize{javaapi/Expr:expr-clone}}\begin{quote}

\sphinxAtStartPar
Deep copy linear expression object.

\sphinxAtStartPar
\sphinxstylestrong{Synopsis}
\begin{quote}

\sphinxAtStartPar
\sphinxcode{\sphinxupquote{Expr clone()}}
\end{quote}

\sphinxAtStartPar
\sphinxstylestrong{Return}
\begin{quote}

\sphinxAtStartPar
cloned linear expression object.
\end{quote}
\end{quote}

\subsubsection{Expr.divide()}
\label{\detokenize{javaapi/Expr:expr-divide}}\begin{quote}

\sphinxAtStartPar
Divide itself by double constant.

\sphinxAtStartPar
\sphinxstylestrong{Synopsis}
\begin{quote}

\sphinxAtStartPar
\sphinxcode{\sphinxupquote{Expr divide(double c)}}
\end{quote}

\sphinxAtStartPar
\sphinxstylestrong{Arguments}
\begin{quote}

\sphinxAtStartPar
\sphinxcode{\sphinxupquote{c}}: constant operand.
\end{quote}

\sphinxAtStartPar
\sphinxstylestrong{Return}
\begin{quote}

\sphinxAtStartPar
linear expression itself.
\end{quote}
\end{quote}

\subsubsection{Expr.evaluate()}
\label{\detokenize{javaapi/Expr:expr-evaluate}}\begin{quote}

\sphinxAtStartPar
Evaluate linear expression after solving.

\sphinxAtStartPar
\sphinxstylestrong{Synopsis}
\begin{quote}

\sphinxAtStartPar
\sphinxcode{\sphinxupquote{double evaluate()}}
\end{quote}

\sphinxAtStartPar
\sphinxstylestrong{Return}
\begin{quote}

\sphinxAtStartPar
value of linear expression.
\end{quote}
\end{quote}

\subsubsection{Expr.getCoeff()}
\label{\detokenize{javaapi/Expr:expr-getcoeff}}\begin{quote}

\sphinxAtStartPar
Get coefficient from the i\sphinxhyphen{}th term in expression.

\sphinxAtStartPar
\sphinxstylestrong{Synopsis}
\begin{quote}

\sphinxAtStartPar
\sphinxcode{\sphinxupquote{double getCoeff(int i)}}
\end{quote}

\sphinxAtStartPar
\sphinxstylestrong{Arguments}
\begin{quote}

\sphinxAtStartPar
\sphinxcode{\sphinxupquote{i}}: index of the term.
\end{quote}

\sphinxAtStartPar
\sphinxstylestrong{Return}
\begin{quote}

\sphinxAtStartPar
coefficient of the i\sphinxhyphen{}th term in expression object.
\end{quote}
\end{quote}

\subsubsection{Expr.getConstant()}
\label{\detokenize{javaapi/Expr:expr-getconstant}}\begin{quote}

\sphinxAtStartPar
Get constant in expression.

\sphinxAtStartPar
\sphinxstylestrong{Synopsis}
\begin{quote}

\sphinxAtStartPar
\sphinxcode{\sphinxupquote{double getConstant()}}
\end{quote}

\sphinxAtStartPar
\sphinxstylestrong{Return}
\begin{quote}

\sphinxAtStartPar
constant in expression.
\end{quote}
\end{quote}

\subsubsection{Expr.getVar()}
\label{\detokenize{javaapi/Expr:expr-getvar}}\begin{quote}

\sphinxAtStartPar
Get variable from the i\sphinxhyphen{}th term in expression.

\sphinxAtStartPar
\sphinxstylestrong{Synopsis}
\begin{quote}

\sphinxAtStartPar
\sphinxcode{\sphinxupquote{Var getVar(int i)}}
\end{quote}

\sphinxAtStartPar
\sphinxstylestrong{Arguments}
\begin{quote}

\sphinxAtStartPar
\sphinxcode{\sphinxupquote{i}}: index of the term.
\end{quote}

\sphinxAtStartPar
\sphinxstylestrong{Return}
\begin{quote}

\sphinxAtStartPar
variable of the i\sphinxhyphen{}th term in expression object.
\end{quote}
\end{quote}

\subsubsection{Expr.multiply()}
\label{\detokenize{javaapi/Expr:expr-multiply}}\begin{quote}

\sphinxAtStartPar
Multiply itself by double constant.

\sphinxAtStartPar
\sphinxstylestrong{Synopsis}
\begin{quote}

\sphinxAtStartPar
\sphinxcode{\sphinxupquote{Expr multiply(double c)}}
\end{quote}

\sphinxAtStartPar
\sphinxstylestrong{Arguments}
\begin{quote}

\sphinxAtStartPar
\sphinxcode{\sphinxupquote{c}}: constant operand.
\end{quote}

\sphinxAtStartPar
\sphinxstylestrong{Return}
\begin{quote}

\sphinxAtStartPar
linear expression itself.
\end{quote}
\end{quote}

\subsubsection{Expr.remove()}
\label{\detokenize{javaapi/Expr:expr-remove}}\begin{quote}

\sphinxAtStartPar
Remove idx\sphinxhyphen{}th term from expression object.

\sphinxAtStartPar
\sphinxstylestrong{Synopsis}
\begin{quote}

\sphinxAtStartPar
\sphinxcode{\sphinxupquote{void remove(int idx)}}
\end{quote}

\sphinxAtStartPar
\sphinxstylestrong{Arguments}
\begin{quote}

\sphinxAtStartPar
\sphinxcode{\sphinxupquote{idx}}: index of the term to be removed.
\end{quote}
\end{quote}

\subsubsection{Expr.remove()}
\label{\detokenize{javaapi/Expr:id9}}\begin{quote}

\sphinxAtStartPar
Remove the term associated with variable from expression.

\sphinxAtStartPar
\sphinxstylestrong{Synopsis}
\begin{quote}

\sphinxAtStartPar
\sphinxcode{\sphinxupquote{void remove(Var var)}}
\end{quote}

\sphinxAtStartPar
\sphinxstylestrong{Arguments}
\begin{quote}

\sphinxAtStartPar
\sphinxcode{\sphinxupquote{var}}: a variable whose term should be removed.
\end{quote}
\end{quote}

\subsubsection{Expr.setCoeff()}
\label{\detokenize{javaapi/Expr:expr-setcoeff}}\begin{quote}

\sphinxAtStartPar
Set coefficient for the i\sphinxhyphen{}th term in expression.

\sphinxAtStartPar
\sphinxstylestrong{Synopsis}
\begin{quote}

\sphinxAtStartPar
\sphinxcode{\sphinxupquote{void setCoeff(int i, double val)}}
\end{quote}

\sphinxAtStartPar
\sphinxstylestrong{Arguments}
\begin{quote}

\sphinxAtStartPar
\sphinxcode{\sphinxupquote{i}}: index of the term.

\sphinxAtStartPar
\sphinxcode{\sphinxupquote{val}}: coefficient of the term.
\end{quote}
\end{quote}

\subsubsection{Expr.setConstant()}
\label{\detokenize{javaapi/Expr:expr-setconstant}}\begin{quote}

\sphinxAtStartPar
Set constant for the expression.

\sphinxAtStartPar
\sphinxstylestrong{Synopsis}
\begin{quote}

\sphinxAtStartPar
\sphinxcode{\sphinxupquote{void setConstant(double constant)}}
\end{quote}

\sphinxAtStartPar
\sphinxstylestrong{Arguments}
\begin{quote}

\sphinxAtStartPar
\sphinxcode{\sphinxupquote{constant}}: the value of the constant.
\end{quote}
\end{quote}

\subsubsection{Expr.size()}
\label{\detokenize{javaapi/Expr:expr-size}}\begin{quote}

\sphinxAtStartPar
Get number of terms in expression.

\sphinxAtStartPar
\sphinxstylestrong{Synopsis}
\begin{quote}

\sphinxAtStartPar
\sphinxcode{\sphinxupquote{long size()}}
\end{quote}

\sphinxAtStartPar
\sphinxstylestrong{Return}
\begin{quote}

\sphinxAtStartPar
number of terms.
\end{quote}
\end{quote}

\subsection{Constraint}
\label{\detokenize{javaapiref:constraint}}\label{\detokenize{javaapiref:chapjavaapiref-constraint}}
\sphinxAtStartPar
COPT constraint object. Constraints are always associated with a particular model.
User creates a constraint object by adding a constraint to a model,
rather than by using constructor of Constraint class.

\sphinxstepscope

\subsubsection{Constraint.get()}
\label{\detokenize{javaapi/Constraint:constraint-get}}\label{\detokenize{javaapi/Constraint::doc}}\begin{quote}

\sphinxAtStartPar
Get information value of the constraint. Support informations of “Dual”, “Slack”, “LB”, “UB”.

\sphinxAtStartPar
\sphinxstylestrong{Synopsis}
\begin{quote}

\sphinxAtStartPar
\sphinxcode{\sphinxupquote{double get(String info)}}
\end{quote}

\sphinxAtStartPar
\sphinxstylestrong{Arguments}
\begin{quote}

\sphinxAtStartPar
\sphinxcode{\sphinxupquote{info}}: name of the information being queried.
\end{quote}

\sphinxAtStartPar
\sphinxstylestrong{Return}
\begin{quote}

\sphinxAtStartPar
information value.
\end{quote}
\end{quote}

\subsubsection{Constraint.getBasis()}
\label{\detokenize{javaapi/Constraint:constraint-getbasis}}\begin{quote}

\sphinxAtStartPar
Get basis status of this constraint.

\sphinxAtStartPar
\sphinxstylestrong{Synopsis}
\begin{quote}

\sphinxAtStartPar
\sphinxcode{\sphinxupquote{int getBasis()}}
\end{quote}

\sphinxAtStartPar
\sphinxstylestrong{Return}
\begin{quote}

\sphinxAtStartPar
basis status.
\end{quote}
\end{quote}

\subsubsection{Constraint.getIdx()}
\label{\detokenize{javaapi/Constraint:constraint-getidx}}\begin{quote}

\sphinxAtStartPar
Get index of the constraint.

\sphinxAtStartPar
\sphinxstylestrong{Synopsis}
\begin{quote}

\sphinxAtStartPar
\sphinxcode{\sphinxupquote{int getIdx()}}
\end{quote}

\sphinxAtStartPar
\sphinxstylestrong{Return}
\begin{quote}

\sphinxAtStartPar
the index of the constraint.
\end{quote}
\end{quote}

\subsubsection{Constraint.getLowerIIS()}
\label{\detokenize{javaapi/Constraint:constraint-getloweriis}}\begin{quote}

\sphinxAtStartPar
Get IIS status for lower bound of the constraint.

\sphinxAtStartPar
\sphinxstylestrong{Synopsis}
\begin{quote}

\sphinxAtStartPar
\sphinxcode{\sphinxupquote{int getLowerIIS()}}
\end{quote}

\sphinxAtStartPar
\sphinxstylestrong{Return}
\begin{quote}

\sphinxAtStartPar
IIS status.
\end{quote}
\end{quote}

\subsubsection{Constraint.getName()}
\label{\detokenize{javaapi/Constraint:constraint-getname}}\begin{quote}

\sphinxAtStartPar
Get name of the constraint.

\sphinxAtStartPar
\sphinxstylestrong{Synopsis}
\begin{quote}

\sphinxAtStartPar
\sphinxcode{\sphinxupquote{String getName()}}
\end{quote}

\sphinxAtStartPar
\sphinxstylestrong{Return}
\begin{quote}

\sphinxAtStartPar
the name of the constraint.
\end{quote}
\end{quote}

\subsubsection{Constraint.getUpperIIS()}
\label{\detokenize{javaapi/Constraint:constraint-getupperiis}}\begin{quote}

\sphinxAtStartPar
Get IIS status for upper bound of the constraint.

\sphinxAtStartPar
\sphinxstylestrong{Synopsis}
\begin{quote}

\sphinxAtStartPar
\sphinxcode{\sphinxupquote{int getUpperIIS()}}
\end{quote}

\sphinxAtStartPar
\sphinxstylestrong{Return}
\begin{quote}

\sphinxAtStartPar
IIS status.
\end{quote}
\end{quote}

\subsubsection{Constraint.remove()}
\label{\detokenize{javaapi/Constraint:constraint-remove}}\begin{quote}

\sphinxAtStartPar
Remove this constraint from model.

\sphinxAtStartPar
\sphinxstylestrong{Synopsis}
\begin{quote}

\sphinxAtStartPar
\sphinxcode{\sphinxupquote{void remove()}}
\end{quote}
\end{quote}

\subsubsection{Constraint.set()}
\label{\detokenize{javaapi/Constraint:constraint-set}}\begin{quote}

\sphinxAtStartPar
Set information value of the constraint. Support informations of “LB” and “UB”.

\sphinxAtStartPar
\sphinxstylestrong{Synopsis}
\begin{quote}

\sphinxAtStartPar
\sphinxcode{\sphinxupquote{void set(String info, double val)}}
\end{quote}

\sphinxAtStartPar
\sphinxstylestrong{Arguments}
\begin{quote}

\sphinxAtStartPar
\sphinxcode{\sphinxupquote{info}}: name of the information.

\sphinxAtStartPar
\sphinxcode{\sphinxupquote{val}}: new information value.
\end{quote}
\end{quote}

\subsubsection{Constraint.setName()}
\label{\detokenize{javaapi/Constraint:constraint-setname}}\begin{quote}

\sphinxAtStartPar
Set name for the constraint.

\sphinxAtStartPar
\sphinxstylestrong{Synopsis}
\begin{quote}

\sphinxAtStartPar
\sphinxcode{\sphinxupquote{void setName(String name)}}
\end{quote}

\sphinxAtStartPar
\sphinxstylestrong{Arguments}
\begin{quote}

\sphinxAtStartPar
\sphinxcode{\sphinxupquote{name}}: the name to set.
\end{quote}
\end{quote}

\subsection{ConstrArray}
\label{\detokenize{javaapiref:constrarray}}\label{\detokenize{javaapiref:chapjavaapiref-constrarray}}
\sphinxAtStartPar
COPT constraint array object. To store and access a set of Java {\hyperref[\detokenize{javaapiref:chapjavaapiref-constraint}]{\sphinxcrossref{\DUrole{std,std-ref}{Constraint}}}}
objects, Cardinal Optimizer provides Java ConstrArray class, which defines the following methods.

\sphinxstepscope

\subsubsection{ConstrArray.ConstrArray()}
\label{\detokenize{javaapi/ConstrArray:constrarray-constrarray}}\label{\detokenize{javaapi/ConstrArray::doc}}\begin{quote}

\sphinxAtStartPar
Constructor of constrarray object.

\sphinxAtStartPar
\sphinxstylestrong{Synopsis}
\begin{quote}

\sphinxAtStartPar
\sphinxcode{\sphinxupquote{ConstrArray()}}
\end{quote}
\end{quote}

\subsubsection{ConstrArray.getConstr()}
\label{\detokenize{javaapi/ConstrArray:constrarray-getconstr}}\begin{quote}

\sphinxAtStartPar
Get idx\sphinxhyphen{}th constraint object.

\sphinxAtStartPar
\sphinxstylestrong{Synopsis}
\begin{quote}

\sphinxAtStartPar
\sphinxcode{\sphinxupquote{Constraint getConstr(int idx)}}
\end{quote}

\sphinxAtStartPar
\sphinxstylestrong{Arguments}
\begin{quote}

\sphinxAtStartPar
\sphinxcode{\sphinxupquote{idx}}: index of the constraint.
\end{quote}

\sphinxAtStartPar
\sphinxstylestrong{Return}
\begin{quote}

\sphinxAtStartPar
constraint object with index idx.
\end{quote}
\end{quote}

\subsubsection{ConstrArray.pushBack()}
\label{\detokenize{javaapi/ConstrArray:constrarray-pushback}}\begin{quote}

\sphinxAtStartPar
Add a constraint object to constraint array.

\sphinxAtStartPar
\sphinxstylestrong{Synopsis}
\begin{quote}

\sphinxAtStartPar
\sphinxcode{\sphinxupquote{void pushBack(Constraint constr)}}
\end{quote}

\sphinxAtStartPar
\sphinxstylestrong{Arguments}
\begin{quote}

\sphinxAtStartPar
\sphinxcode{\sphinxupquote{constr}}: a constraint object.
\end{quote}
\end{quote}

\subsubsection{ConstrArray.size()}
\label{\detokenize{javaapi/ConstrArray:constrarray-size}}\begin{quote}

\sphinxAtStartPar
Get the number of constraint objects.

\sphinxAtStartPar
\sphinxstylestrong{Synopsis}
\begin{quote}

\sphinxAtStartPar
\sphinxcode{\sphinxupquote{int size()}}
\end{quote}

\sphinxAtStartPar
\sphinxstylestrong{Return}
\begin{quote}

\sphinxAtStartPar
number of constraint objects.
\end{quote}
\end{quote}

\subsection{ConstrBuilder}
\label{\detokenize{javaapiref:constrbuilder}}\label{\detokenize{javaapiref:chapjavaapiref-constrbuilder}}
\sphinxAtStartPar
COPT constraint builder object. To help building a constraint, given a linear expression,
constraint sense and right\sphinxhyphen{}hand side value, Cardinal Optimizer provides Java ConstrBuilder
class, which defines the following methods.

\sphinxstepscope

\subsubsection{ConstrBuilder.ConstrBuilder()}
\label{\detokenize{javaapi/ConstrBuilder:constrbuilder-constrbuilder}}\label{\detokenize{javaapi/ConstrBuilder::doc}}\begin{quote}

\sphinxAtStartPar
Constructor of constrbuilder object.

\sphinxAtStartPar
\sphinxstylestrong{Synopsis}
\begin{quote}

\sphinxAtStartPar
\sphinxcode{\sphinxupquote{ConstrBuilder()}}
\end{quote}
\end{quote}

\subsubsection{ConstrBuilder.getExpr()}
\label{\detokenize{javaapi/ConstrBuilder:constrbuilder-getexpr}}\begin{quote}

\sphinxAtStartPar
Get expression associated with constraint.

\sphinxAtStartPar
\sphinxstylestrong{Synopsis}
\begin{quote}

\sphinxAtStartPar
\sphinxcode{\sphinxupquote{Expr getExpr()}}
\end{quote}

\sphinxAtStartPar
\sphinxstylestrong{Return}
\begin{quote}

\sphinxAtStartPar
expression object.
\end{quote}
\end{quote}

\subsubsection{ConstrBuilder.getRange()}
\label{\detokenize{javaapi/ConstrBuilder:constrbuilder-getrange}}\begin{quote}

\sphinxAtStartPar
Get range from lower bound to upper bound of range constraint.

\sphinxAtStartPar
\sphinxstylestrong{Synopsis}
\begin{quote}

\sphinxAtStartPar
\sphinxcode{\sphinxupquote{double getRange()}}
\end{quote}

\sphinxAtStartPar
\sphinxstylestrong{Return}
\begin{quote}

\sphinxAtStartPar
length from lower bound to upper bound of the constraint.
\end{quote}
\end{quote}

\subsubsection{ConstrBuilder.getSense()}
\label{\detokenize{javaapi/ConstrBuilder:constrbuilder-getsense}}\begin{quote}

\sphinxAtStartPar
Get sense associated with constraint.

\sphinxAtStartPar
\sphinxstylestrong{Synopsis}
\begin{quote}

\sphinxAtStartPar
\sphinxcode{\sphinxupquote{char getSense()}}
\end{quote}

\sphinxAtStartPar
\sphinxstylestrong{Return}
\begin{quote}

\sphinxAtStartPar
constraint sense.
\end{quote}
\end{quote}

\subsubsection{ConstrBuilder.set()}
\label{\detokenize{javaapi/ConstrBuilder:constrbuilder-set}}\begin{quote}

\sphinxAtStartPar
Set detail of a constraint to its builder object.

\sphinxAtStartPar
\sphinxstylestrong{Synopsis}
\begin{quote}

\sphinxAtStartPar
\sphinxcode{\sphinxupquote{void set(}}
\begin{quote}

\sphinxAtStartPar
\sphinxcode{\sphinxupquote{Expr expr,}}

\sphinxAtStartPar
\sphinxcode{\sphinxupquote{char sense,}}

\sphinxAtStartPar
\sphinxcode{\sphinxupquote{double rhs)}}
\end{quote}
\end{quote}

\sphinxAtStartPar
\sphinxstylestrong{Arguments}
\begin{quote}

\sphinxAtStartPar
\sphinxcode{\sphinxupquote{expr}}: expression object at one side of the constraint

\sphinxAtStartPar
\sphinxcode{\sphinxupquote{sense}}: constraint sense other than COPT\_RANGE.

\sphinxAtStartPar
\sphinxcode{\sphinxupquote{rhs}}: constant of right side of the constraint.
\end{quote}
\end{quote}

\subsubsection{ConstrBuilder.setRange()}
\label{\detokenize{javaapi/ConstrBuilder:constrbuilder-setrange}}\begin{quote}

\sphinxAtStartPar
Set a range constraint to its builder.

\sphinxAtStartPar
\sphinxstylestrong{Synopsis}
\begin{quote}

\sphinxAtStartPar
\sphinxcode{\sphinxupquote{void setRange(Expr expr, double range)}}
\end{quote}

\sphinxAtStartPar
\sphinxstylestrong{Arguments}
\begin{quote}

\sphinxAtStartPar
\sphinxcode{\sphinxupquote{expr}}: expression object, whose constant is negative upper bound.

\sphinxAtStartPar
\sphinxcode{\sphinxupquote{range}}: length from lower bound to upper bound of the constraint. Must greater than 0.
\end{quote}
\end{quote}

\subsection{ConstrBuilderArray}
\label{\detokenize{javaapiref:constrbuilderarray}}\label{\detokenize{javaapiref:chapjavaapiref-constrbuilderarray}}
\sphinxAtStartPar
COPT constraint builder array object. To store and access a set of Java {\hyperref[\detokenize{javaapiref:chapjavaapiref-constrbuilder}]{\sphinxcrossref{\DUrole{std,std-ref}{ConstrBuilder}}}}
objects, Cardinal Optimizer provides Java ConstrBuilderArray class, which defines the following methods.

\sphinxstepscope

\subsubsection{ConstrBuilderArray.ConstrBuilderArray()}
\label{\detokenize{javaapi/ConstrBuilderArray:constrbuilderarray-constrbuilderarray}}\label{\detokenize{javaapi/ConstrBuilderArray::doc}}\begin{quote}

\sphinxAtStartPar
Constructor of constrbuilderarray object.

\sphinxAtStartPar
\sphinxstylestrong{Synopsis}
\begin{quote}

\sphinxAtStartPar
\sphinxcode{\sphinxupquote{ConstrBuilderArray()}}
\end{quote}
\end{quote}

\subsubsection{ConstrBuilderArray.getBuilder()}
\label{\detokenize{javaapi/ConstrBuilderArray:constrbuilderarray-getbuilder}}\begin{quote}

\sphinxAtStartPar
Get idx\sphinxhyphen{}th constraint builder object.

\sphinxAtStartPar
\sphinxstylestrong{Synopsis}
\begin{quote}

\sphinxAtStartPar
\sphinxcode{\sphinxupquote{ConstrBuilder getBuilder(int idx)}}
\end{quote}

\sphinxAtStartPar
\sphinxstylestrong{Arguments}
\begin{quote}

\sphinxAtStartPar
\sphinxcode{\sphinxupquote{idx}}: index of the constraint builder.
\end{quote}

\sphinxAtStartPar
\sphinxstylestrong{Return}
\begin{quote}

\sphinxAtStartPar
constraint builder object with index idx.
\end{quote}
\end{quote}

\subsubsection{ConstrBuilderArray.pushBack()}
\label{\detokenize{javaapi/ConstrBuilderArray:constrbuilderarray-pushback}}\begin{quote}

\sphinxAtStartPar
Add a constraint builder object to constraint builder array.

\sphinxAtStartPar
\sphinxstylestrong{Synopsis}
\begin{quote}

\sphinxAtStartPar
\sphinxcode{\sphinxupquote{void pushBack(ConstrBuilder builder)}}
\end{quote}

\sphinxAtStartPar
\sphinxstylestrong{Arguments}
\begin{quote}

\sphinxAtStartPar
\sphinxcode{\sphinxupquote{builder}}: a constraint builder object.
\end{quote}
\end{quote}

\subsubsection{ConstrBuilderArray.size()}
\label{\detokenize{javaapi/ConstrBuilderArray:constrbuilderarray-size}}\begin{quote}

\sphinxAtStartPar
Get the number of constraint builder objects.

\sphinxAtStartPar
\sphinxstylestrong{Synopsis}
\begin{quote}

\sphinxAtStartPar
\sphinxcode{\sphinxupquote{int size()}}
\end{quote}

\sphinxAtStartPar
\sphinxstylestrong{Return}
\begin{quote}

\sphinxAtStartPar
number of constraint builder objects.
\end{quote}
\end{quote}

\subsection{Column}
\label{\detokenize{javaapiref:column}}\label{\detokenize{javaapiref:chapjavaapiref-column}}
\sphinxAtStartPar
COPT column object. A column consists of a list of constraints and associated coefficients.
Columns are used to represent the set of constraints in which a variable participates,
and the asssociated coefficents.

\sphinxstepscope

\subsubsection{Column.Column()}
\label{\detokenize{javaapi/Column:column-column}}\label{\detokenize{javaapi/Column::doc}}\begin{quote}

\sphinxAtStartPar
Constructor of column.

\sphinxAtStartPar
\sphinxstylestrong{Synopsis}
\begin{quote}

\sphinxAtStartPar
\sphinxcode{\sphinxupquote{Column()}}
\end{quote}
\end{quote}

\subsubsection{Column.addColumn()}
\label{\detokenize{javaapi/Column:column-addcolumn}}\begin{quote}

\sphinxAtStartPar
Add a column to self.

\sphinxAtStartPar
\sphinxstylestrong{Synopsis}
\begin{quote}

\sphinxAtStartPar
\sphinxcode{\sphinxupquote{void addColumn(Column col, double mult)}}
\end{quote}

\sphinxAtStartPar
\sphinxstylestrong{Arguments}
\begin{quote}

\sphinxAtStartPar
\sphinxcode{\sphinxupquote{col}}: column object to be added.

\sphinxAtStartPar
\sphinxcode{\sphinxupquote{mult}}: constant multiplier.
\end{quote}
\end{quote}

\subsubsection{Column.addColumn()}
\label{\detokenize{javaapi/Column:id1}}\begin{quote}

\sphinxAtStartPar
Add a column to self.

\sphinxAtStartPar
\sphinxstylestrong{Synopsis}
\begin{quote}

\sphinxAtStartPar
\sphinxcode{\sphinxupquote{void addColumn(Column col)}}
\end{quote}

\sphinxAtStartPar
\sphinxstylestrong{Arguments}
\begin{quote}

\sphinxAtStartPar
\sphinxcode{\sphinxupquote{col}}: column object to be added.
\end{quote}
\end{quote}

\subsubsection{Column.addTerm()}
\label{\detokenize{javaapi/Column:column-addterm}}\begin{quote}

\sphinxAtStartPar
Add a term to column object.

\sphinxAtStartPar
\sphinxstylestrong{Synopsis}
\begin{quote}

\sphinxAtStartPar
\sphinxcode{\sphinxupquote{void addTerm(Constraint constr, double coeff)}}
\end{quote}

\sphinxAtStartPar
\sphinxstylestrong{Arguments}
\begin{quote}

\sphinxAtStartPar
\sphinxcode{\sphinxupquote{constr}}: a constraint for new term.

\sphinxAtStartPar
\sphinxcode{\sphinxupquote{coeff}}: coefficient for new term.
\end{quote}
\end{quote}

\subsubsection{Column.addTerms()}
\label{\detokenize{javaapi/Column:column-addterms}}\begin{quote}

\sphinxAtStartPar
Add terms to column object.

\sphinxAtStartPar
\sphinxstylestrong{Synopsis}
\begin{quote}

\sphinxAtStartPar
\sphinxcode{\sphinxupquote{void addTerms(Constraint{[}{]} constrs, double coeff)}}
\end{quote}

\sphinxAtStartPar
\sphinxstylestrong{Arguments}
\begin{quote}

\sphinxAtStartPar
\sphinxcode{\sphinxupquote{constrs}}: constraints for added terms.

\sphinxAtStartPar
\sphinxcode{\sphinxupquote{coeff}}: coefficient for added terms.
\end{quote}
\end{quote}

\subsubsection{Column.addTerms()}
\label{\detokenize{javaapi/Column:id2}}\begin{quote}

\sphinxAtStartPar
Add terms to column object.

\sphinxAtStartPar
\sphinxstylestrong{Synopsis}
\begin{quote}

\sphinxAtStartPar
\sphinxcode{\sphinxupquote{void addTerms(Constraint{[}{]} constrs, double{[}{]} coeffs)}}
\end{quote}

\sphinxAtStartPar
\sphinxstylestrong{Arguments}
\begin{quote}

\sphinxAtStartPar
\sphinxcode{\sphinxupquote{constrs}}: constraints for added terms.

\sphinxAtStartPar
\sphinxcode{\sphinxupquote{coeffs}}: coefficients for added terms.
\end{quote}
\end{quote}

\subsubsection{Column.addTerms()}
\label{\detokenize{javaapi/Column:id3}}\begin{quote}

\sphinxAtStartPar
Add terms to column object.

\sphinxAtStartPar
\sphinxstylestrong{Synopsis}
\begin{quote}

\sphinxAtStartPar
\sphinxcode{\sphinxupquote{void addTerms(ConstrArray constrs, double coeff)}}
\end{quote}

\sphinxAtStartPar
\sphinxstylestrong{Arguments}
\begin{quote}

\sphinxAtStartPar
\sphinxcode{\sphinxupquote{constrs}}: constraints for added terms.

\sphinxAtStartPar
\sphinxcode{\sphinxupquote{coeff}}: coefficient for added terms.
\end{quote}
\end{quote}

\subsubsection{Column.addTerms()}
\label{\detokenize{javaapi/Column:id4}}\begin{quote}

\sphinxAtStartPar
Add terms to column object.

\sphinxAtStartPar
\sphinxstylestrong{Synopsis}
\begin{quote}

\sphinxAtStartPar
\sphinxcode{\sphinxupquote{void addTerms(ConstrArray constrs, double{[}{]} coeffs)}}
\end{quote}

\sphinxAtStartPar
\sphinxstylestrong{Arguments}
\begin{quote}

\sphinxAtStartPar
\sphinxcode{\sphinxupquote{constrs}}: constraints for added terms.

\sphinxAtStartPar
\sphinxcode{\sphinxupquote{coeffs}}: coefficients for added terms.
\end{quote}
\end{quote}

\subsubsection{Column.clear()}
\label{\detokenize{javaapi/Column:column-clear}}\begin{quote}

\sphinxAtStartPar
Clear all terms.

\sphinxAtStartPar
\sphinxstylestrong{Synopsis}
\begin{quote}

\sphinxAtStartPar
\sphinxcode{\sphinxupquote{void clear()}}
\end{quote}
\end{quote}

\subsubsection{Column.clone()}
\label{\detokenize{javaapi/Column:column-clone}}\begin{quote}

\sphinxAtStartPar
Deep copy column object.

\sphinxAtStartPar
\sphinxstylestrong{Synopsis}
\begin{quote}

\sphinxAtStartPar
\sphinxcode{\sphinxupquote{Column clone()}}
\end{quote}

\sphinxAtStartPar
\sphinxstylestrong{Return}
\begin{quote}

\sphinxAtStartPar
cloned column object.
\end{quote}
\end{quote}

\subsubsection{Column.getCoeff()}
\label{\detokenize{javaapi/Column:column-getcoeff}}\begin{quote}

\sphinxAtStartPar
Get coefficient from the i\sphinxhyphen{}th term in column object.

\sphinxAtStartPar
\sphinxstylestrong{Synopsis}
\begin{quote}

\sphinxAtStartPar
\sphinxcode{\sphinxupquote{double getCoeff(int i)}}
\end{quote}

\sphinxAtStartPar
\sphinxstylestrong{Arguments}
\begin{quote}

\sphinxAtStartPar
\sphinxcode{\sphinxupquote{i}}: index of the term.
\end{quote}

\sphinxAtStartPar
\sphinxstylestrong{Return}
\begin{quote}

\sphinxAtStartPar
coefficient of the i\sphinxhyphen{}th term in column object.
\end{quote}
\end{quote}

\subsubsection{Column.getConstr()}
\label{\detokenize{javaapi/Column:column-getconstr}}\begin{quote}

\sphinxAtStartPar
Get constraint from the i\sphinxhyphen{}th term in column object.

\sphinxAtStartPar
\sphinxstylestrong{Synopsis}
\begin{quote}

\sphinxAtStartPar
\sphinxcode{\sphinxupquote{Constraint getConstr(int i)}}
\end{quote}

\sphinxAtStartPar
\sphinxstylestrong{Arguments}
\begin{quote}

\sphinxAtStartPar
\sphinxcode{\sphinxupquote{i}}: index of the term.
\end{quote}

\sphinxAtStartPar
\sphinxstylestrong{Return}
\begin{quote}

\sphinxAtStartPar
constraint of the i\sphinxhyphen{}th term in column object.
\end{quote}
\end{quote}

\subsubsection{Column.remove()}
\label{\detokenize{javaapi/Column:column-remove}}\begin{quote}

\sphinxAtStartPar
Remove i\sphinxhyphen{}th term from column object.

\sphinxAtStartPar
\sphinxstylestrong{Synopsis}
\begin{quote}

\sphinxAtStartPar
\sphinxcode{\sphinxupquote{void remove(int idx)}}
\end{quote}

\sphinxAtStartPar
\sphinxstylestrong{Arguments}
\begin{quote}

\sphinxAtStartPar
\sphinxcode{\sphinxupquote{idx}}: index of the term to be removed.
\end{quote}
\end{quote}

\subsubsection{Column.remove()}
\label{\detokenize{javaapi/Column:id5}}\begin{quote}

\sphinxAtStartPar
Remove the term associated with constraint from column object.

\sphinxAtStartPar
\sphinxstylestrong{Synopsis}
\begin{quote}

\sphinxAtStartPar
\sphinxcode{\sphinxupquote{void remove(Constraint constr)}}
\end{quote}

\sphinxAtStartPar
\sphinxstylestrong{Arguments}
\begin{quote}

\sphinxAtStartPar
\sphinxcode{\sphinxupquote{constr}}: a constraint whose term should be removed.
\end{quote}
\end{quote}

\subsubsection{Column.size()}
\label{\detokenize{javaapi/Column:column-size}}\begin{quote}

\sphinxAtStartPar
Get number of terms in column object.

\sphinxAtStartPar
\sphinxstylestrong{Synopsis}
\begin{quote}

\sphinxAtStartPar
\sphinxcode{\sphinxupquote{int size()}}
\end{quote}

\sphinxAtStartPar
\sphinxstylestrong{Return}
\begin{quote}

\sphinxAtStartPar
number of terms.
\end{quote}
\end{quote}

\subsection{ColumnArray}
\label{\detokenize{javaapiref:columnarray}}\label{\detokenize{javaapiref:chapjavaapiref-columnarray}}
\sphinxAtStartPar
COPT column array object. To store and access a set of Java
{\hyperref[\detokenize{javaapiref:chapjavaapiref-column}]{\sphinxcrossref{\DUrole{std,std-ref}{Column}}}} objects, Cardinal Optimizer provides Java
ColumnArray class, which defines the following methods.

\sphinxstepscope

\subsubsection{ColumnArray.ColumnArray()}
\label{\detokenize{javaapi/ColumnArray:columnarray-columnarray}}\label{\detokenize{javaapi/ColumnArray::doc}}\begin{quote}

\sphinxAtStartPar
Constructor of columnarray object.

\sphinxAtStartPar
\sphinxstylestrong{Synopsis}
\begin{quote}

\sphinxAtStartPar
\sphinxcode{\sphinxupquote{ColumnArray()}}
\end{quote}
\end{quote}

\subsubsection{ColumnArray.clear()}
\label{\detokenize{javaapi/ColumnArray:columnarray-clear}}\begin{quote}

\sphinxAtStartPar
Clear all column objects.

\sphinxAtStartPar
\sphinxstylestrong{Synopsis}
\begin{quote}

\sphinxAtStartPar
\sphinxcode{\sphinxupquote{void clear()}}
\end{quote}
\end{quote}

\subsubsection{ColumnArray.getColumn()}
\label{\detokenize{javaapi/ColumnArray:columnarray-getcolumn}}\begin{quote}

\sphinxAtStartPar
Get idx\sphinxhyphen{}th column object.

\sphinxAtStartPar
\sphinxstylestrong{Synopsis}
\begin{quote}

\sphinxAtStartPar
\sphinxcode{\sphinxupquote{Column getColumn(int idx)}}
\end{quote}

\sphinxAtStartPar
\sphinxstylestrong{Arguments}
\begin{quote}

\sphinxAtStartPar
\sphinxcode{\sphinxupquote{idx}}: index of the column.
\end{quote}

\sphinxAtStartPar
\sphinxstylestrong{Return}
\begin{quote}

\sphinxAtStartPar
column object with index idx.
\end{quote}
\end{quote}

\subsubsection{ColumnArray.pushBack()}
\label{\detokenize{javaapi/ColumnArray:columnarray-pushback}}\begin{quote}

\sphinxAtStartPar
Add a column object to column array.

\sphinxAtStartPar
\sphinxstylestrong{Synopsis}
\begin{quote}

\sphinxAtStartPar
\sphinxcode{\sphinxupquote{void pushBack(Column col)}}
\end{quote}

\sphinxAtStartPar
\sphinxstylestrong{Arguments}
\begin{quote}

\sphinxAtStartPar
\sphinxcode{\sphinxupquote{col}}: a column object.
\end{quote}
\end{quote}

\subsubsection{ColumnArray.size()}
\label{\detokenize{javaapi/ColumnArray:columnarray-size}}\begin{quote}

\sphinxAtStartPar
Get the number of column objects.

\sphinxAtStartPar
\sphinxstylestrong{Synopsis}
\begin{quote}

\sphinxAtStartPar
\sphinxcode{\sphinxupquote{int size()}}
\end{quote}

\sphinxAtStartPar
\sphinxstylestrong{Return}
\begin{quote}

\sphinxAtStartPar
number of column objects.
\end{quote}
\end{quote}

\subsection{Sos}
\label{\detokenize{javaapiref:sos}}\label{\detokenize{javaapiref:chapjavaapiref-sos}}
\sphinxAtStartPar
COPT SOS constraint object. SOS constraints are always associated with a particular model.
User creates an SOS constraint object by adding an SOS constraint to a model,
rather than by using constructor of Sos class.

\sphinxAtStartPar
An SOS constraint can be type 1 or 2 (\sphinxcode{\sphinxupquote{COPT\_SOS\_TYPE1}} or \sphinxcode{\sphinxupquote{COPT\_SOS\_TYPE2}}).

\sphinxstepscope

\subsubsection{Sos.getIdx()}
\label{\detokenize{javaapi/Sos:sos-getidx}}\label{\detokenize{javaapi/Sos::doc}}\begin{quote}

\sphinxAtStartPar
Get the index of SOS constraint.

\sphinxAtStartPar
\sphinxstylestrong{Synopsis}
\begin{quote}

\sphinxAtStartPar
\sphinxcode{\sphinxupquote{int getIdx()}}
\end{quote}

\sphinxAtStartPar
\sphinxstylestrong{Return}
\begin{quote}

\sphinxAtStartPar
index of SOS constraint.
\end{quote}
\end{quote}

\subsubsection{Sos.getIIS()}
\label{\detokenize{javaapi/Sos:sos-getiis}}\begin{quote}

\sphinxAtStartPar
Get IIS status of the SOS constraint.

\sphinxAtStartPar
\sphinxstylestrong{Synopsis}
\begin{quote}

\sphinxAtStartPar
\sphinxcode{\sphinxupquote{int getIIS()}}
\end{quote}

\sphinxAtStartPar
\sphinxstylestrong{Return}
\begin{quote}

\sphinxAtStartPar
IIS status.
\end{quote}
\end{quote}

\subsubsection{Sos.remove()}
\label{\detokenize{javaapi/Sos:sos-remove}}\begin{quote}

\sphinxAtStartPar
Remove the SOS constraint from model.

\sphinxAtStartPar
\sphinxstylestrong{Synopsis}
\begin{quote}

\sphinxAtStartPar
\sphinxcode{\sphinxupquote{void remove()}}
\end{quote}
\end{quote}

\subsection{SosArray}
\label{\detokenize{javaapiref:sosarray}}\label{\detokenize{javaapiref:chapjavaapiref-sosarray}}
\sphinxAtStartPar
COPT SOS constraint array object. To store and access a set of Java {\hyperref[\detokenize{javaapiref:chapjavaapiref-sos}]{\sphinxcrossref{\DUrole{std,std-ref}{Sos}}}}
objects, Cardinal Optimizer provides Java SosArray class, which defines the following methods.

\sphinxstepscope

\subsubsection{SosArray.SosArray()}
\label{\detokenize{javaapi/SosArray:sosarray-sosarray}}\label{\detokenize{javaapi/SosArray::doc}}\begin{quote}

\sphinxAtStartPar
Constructor of sosarray object.

\sphinxAtStartPar
\sphinxstylestrong{Synopsis}
\begin{quote}

\sphinxAtStartPar
\sphinxcode{\sphinxupquote{SosArray()}}
\end{quote}
\end{quote}

\subsubsection{SosArray.getSos()}
\label{\detokenize{javaapi/SosArray:sosarray-getsos}}\begin{quote}

\sphinxAtStartPar
Get idx\sphinxhyphen{}th SOS object.

\sphinxAtStartPar
\sphinxstylestrong{Synopsis}
\begin{quote}

\sphinxAtStartPar
\sphinxcode{\sphinxupquote{Sos getSos(int idx)}}
\end{quote}

\sphinxAtStartPar
\sphinxstylestrong{Arguments}
\begin{quote}

\sphinxAtStartPar
\sphinxcode{\sphinxupquote{idx}}: index of SOS.
\end{quote}

\sphinxAtStartPar
\sphinxstylestrong{Return}
\begin{quote}

\sphinxAtStartPar
SOS object with index idx.
\end{quote}
\end{quote}

\subsubsection{SosArray.pushBack()}
\label{\detokenize{javaapi/SosArray:sosarray-pushback}}\begin{quote}

\sphinxAtStartPar
Add a SOS constraint object to SOS constraint array.

\sphinxAtStartPar
\sphinxstylestrong{Synopsis}
\begin{quote}

\sphinxAtStartPar
\sphinxcode{\sphinxupquote{void pushBack(Sos sos)}}
\end{quote}

\sphinxAtStartPar
\sphinxstylestrong{Arguments}
\begin{quote}

\sphinxAtStartPar
\sphinxcode{\sphinxupquote{sos}}: a SOS constraint object.
\end{quote}
\end{quote}

\subsubsection{SosArray.size()}
\label{\detokenize{javaapi/SosArray:sosarray-size}}\begin{quote}

\sphinxAtStartPar
Get the number of SOS constraint objects.

\sphinxAtStartPar
\sphinxstylestrong{Synopsis}
\begin{quote}

\sphinxAtStartPar
\sphinxcode{\sphinxupquote{int size()}}
\end{quote}

\sphinxAtStartPar
\sphinxstylestrong{Return}
\begin{quote}

\sphinxAtStartPar
number of SOS constraint objects.
\end{quote}
\end{quote}

\subsection{SosBuilder}
\label{\detokenize{javaapiref:sosbuilder}}\label{\detokenize{javaapiref:chapjavaapiref-sosbuilder}}
\sphinxAtStartPar
COPT SOS constraint builder object. To help building an SOS constraint, given the SOS type,
a set of variables and associated weights, Cardinal Optimizer provides Java SosBuilder
class, which defines the following methods.

\sphinxstepscope

\subsubsection{SosBuilder.SosBuilder()}
\label{\detokenize{javaapi/SosBuilder:sosbuilder-sosbuilder}}\label{\detokenize{javaapi/SosBuilder::doc}}\begin{quote}

\sphinxAtStartPar
Constructor of sosbuilder object.

\sphinxAtStartPar
\sphinxstylestrong{Synopsis}
\begin{quote}

\sphinxAtStartPar
\sphinxcode{\sphinxupquote{SosBuilder()}}
\end{quote}
\end{quote}

\subsubsection{SosBuilder.getSize()}
\label{\detokenize{javaapi/SosBuilder:sosbuilder-getsize}}\begin{quote}

\sphinxAtStartPar
Get number of terms in SOS constraint.

\sphinxAtStartPar
\sphinxstylestrong{Synopsis}
\begin{quote}

\sphinxAtStartPar
\sphinxcode{\sphinxupquote{int getSize()}}
\end{quote}

\sphinxAtStartPar
\sphinxstylestrong{Return}
\begin{quote}

\sphinxAtStartPar
number of terms.
\end{quote}
\end{quote}

\subsubsection{SosBuilder.getType()}
\label{\detokenize{javaapi/SosBuilder:sosbuilder-gettype}}\begin{quote}

\sphinxAtStartPar
Get type of SOS constraint.

\sphinxAtStartPar
\sphinxstylestrong{Synopsis}
\begin{quote}

\sphinxAtStartPar
\sphinxcode{\sphinxupquote{int getType()}}
\end{quote}

\sphinxAtStartPar
\sphinxstylestrong{Return}
\begin{quote}

\sphinxAtStartPar
type of SOS constraint.
\end{quote}
\end{quote}

\subsubsection{SosBuilder.getVar()}
\label{\detokenize{javaapi/SosBuilder:sosbuilder-getvar}}\begin{quote}

\sphinxAtStartPar
Get variable from the idx\sphinxhyphen{}th term in SOS constraint.

\sphinxAtStartPar
\sphinxstylestrong{Synopsis}
\begin{quote}

\sphinxAtStartPar
\sphinxcode{\sphinxupquote{Var getVar(int idx)}}
\end{quote}

\sphinxAtStartPar
\sphinxstylestrong{Arguments}
\begin{quote}

\sphinxAtStartPar
\sphinxcode{\sphinxupquote{idx}}: index of the term.
\end{quote}

\sphinxAtStartPar
\sphinxstylestrong{Return}
\begin{quote}

\sphinxAtStartPar
variable of the idx\sphinxhyphen{}th term in SOS constraint.
\end{quote}
\end{quote}

\subsubsection{SosBuilder.GetVars()}
\label{\detokenize{javaapi/SosBuilder:sosbuilder-getvars}}\begin{quote}

\sphinxAtStartPar
Get all variables in a SOS constraint.

\sphinxAtStartPar
\sphinxstylestrong{Synopsis}
\begin{quote}

\sphinxAtStartPar
\sphinxcode{\sphinxupquote{VarArray GetVars()}}
\end{quote}

\sphinxAtStartPar
\sphinxstylestrong{Return}
\begin{quote}

\sphinxAtStartPar
variables in a SOS constraint.
\end{quote}
\end{quote}

\subsubsection{SosBuilder.getWeight()}
\label{\detokenize{javaapi/SosBuilder:sosbuilder-getweight}}\begin{quote}

\sphinxAtStartPar
Get weight from the idx\sphinxhyphen{}th term in SOS constraint.

\sphinxAtStartPar
\sphinxstylestrong{Synopsis}
\begin{quote}

\sphinxAtStartPar
\sphinxcode{\sphinxupquote{double getWeight(int idx)}}
\end{quote}

\sphinxAtStartPar
\sphinxstylestrong{Arguments}
\begin{quote}

\sphinxAtStartPar
\sphinxcode{\sphinxupquote{idx}}: index of the term.
\end{quote}

\sphinxAtStartPar
\sphinxstylestrong{Return}
\begin{quote}

\sphinxAtStartPar
weight of the idx\sphinxhyphen{}th term in SOS constraint.
\end{quote}
\end{quote}

\subsubsection{SosBuilder.getWeights()}
\label{\detokenize{javaapi/SosBuilder:sosbuilder-getweights}}\begin{quote}

\sphinxAtStartPar
Get weights of all terms in SOS constraint.

\sphinxAtStartPar
\sphinxstylestrong{Synopsis}
\begin{quote}

\sphinxAtStartPar
\sphinxcode{\sphinxupquote{double{[}{]} getWeights()}}
\end{quote}

\sphinxAtStartPar
\sphinxstylestrong{Return}
\begin{quote}

\sphinxAtStartPar
array of weights.
\end{quote}
\end{quote}

\subsubsection{SosBuilder.set()}
\label{\detokenize{javaapi/SosBuilder:sosbuilder-set}}\begin{quote}

\sphinxAtStartPar
Set variables and weights of SOS constraint.

\sphinxAtStartPar
\sphinxstylestrong{Synopsis}
\begin{quote}

\sphinxAtStartPar
\sphinxcode{\sphinxupquote{void set(}}
\begin{quote}

\sphinxAtStartPar
\sphinxcode{\sphinxupquote{VarArray vars,}}

\sphinxAtStartPar
\sphinxcode{\sphinxupquote{double{[}{]} weights,}}

\sphinxAtStartPar
\sphinxcode{\sphinxupquote{int type)}}
\end{quote}
\end{quote}

\sphinxAtStartPar
\sphinxstylestrong{Arguments}
\begin{quote}

\sphinxAtStartPar
\sphinxcode{\sphinxupquote{vars}}: variable array object.

\sphinxAtStartPar
\sphinxcode{\sphinxupquote{weights}}: array of weights.

\sphinxAtStartPar
\sphinxcode{\sphinxupquote{type}}: type of SOS constraint.
\end{quote}
\end{quote}

\subsection{SosBuilderArray}
\label{\detokenize{javaapiref:sosbuilderarray}}\label{\detokenize{javaapiref:chapjavaapiref-sosbuilderarray}}
\sphinxAtStartPar
COPT SOS constraint builder array object. To store and access a set of Java {\hyperref[\detokenize{javaapiref:chapjavaapiref-sosbuilder}]{\sphinxcrossref{\DUrole{std,std-ref}{SosBuilder}}}}
objects, Cardinal Optimizer provides Java SosBuilderArray class, which defines the following methods.

\sphinxstepscope

\subsubsection{SosBuilderArray.SosBuilderArray()}
\label{\detokenize{javaapi/SosBuilderArray:sosbuilderarray-sosbuilderarray}}\label{\detokenize{javaapi/SosBuilderArray::doc}}\begin{quote}

\sphinxAtStartPar
Constructor of sosbuilderarray object.

\sphinxAtStartPar
\sphinxstylestrong{Synopsis}
\begin{quote}

\sphinxAtStartPar
\sphinxcode{\sphinxupquote{SosBuilderArray()}}
\end{quote}
\end{quote}

\subsubsection{SosBuilderArray.getBuilder()}
\label{\detokenize{javaapi/SosBuilderArray:sosbuilderarray-getbuilder}}\begin{quote}

\sphinxAtStartPar
Get idx\sphinxhyphen{}th SOS constraint builder object.

\sphinxAtStartPar
\sphinxstylestrong{Synopsis}
\begin{quote}

\sphinxAtStartPar
\sphinxcode{\sphinxupquote{SosBuilder getBuilder(int idx)}}
\end{quote}

\sphinxAtStartPar
\sphinxstylestrong{Arguments}
\begin{quote}

\sphinxAtStartPar
\sphinxcode{\sphinxupquote{idx}}: index of the SOS constraint builder.
\end{quote}

\sphinxAtStartPar
\sphinxstylestrong{Return}
\begin{quote}

\sphinxAtStartPar
SOS constraint builder object with index idx.
\end{quote}
\end{quote}

\subsubsection{SosBuilderArray.pushBack()}
\label{\detokenize{javaapi/SosBuilderArray:sosbuilderarray-pushback}}\begin{quote}

\sphinxAtStartPar
Add a SOS constraint builder object to SOS constraint builder array.

\sphinxAtStartPar
\sphinxstylestrong{Synopsis}
\begin{quote}

\sphinxAtStartPar
\sphinxcode{\sphinxupquote{void pushBack(SosBuilder builder)}}
\end{quote}

\sphinxAtStartPar
\sphinxstylestrong{Arguments}
\begin{quote}

\sphinxAtStartPar
\sphinxcode{\sphinxupquote{builder}}: a SOS constraint builder object.
\end{quote}
\end{quote}

\subsubsection{SosBuilderArray.size()}
\label{\detokenize{javaapi/SosBuilderArray:sosbuilderarray-size}}\begin{quote}

\sphinxAtStartPar
Get the number of SOS constraint builder objects.

\sphinxAtStartPar
\sphinxstylestrong{Synopsis}
\begin{quote}

\sphinxAtStartPar
\sphinxcode{\sphinxupquote{int size()}}
\end{quote}

\sphinxAtStartPar
\sphinxstylestrong{Return}
\begin{quote}

\sphinxAtStartPar
number of SOS constraint builder objects.
\end{quote}
\end{quote}

\subsection{GenConstr}
\label{\detokenize{javaapiref:genconstr}}\label{\detokenize{javaapiref:chapjavaapiref-genconstr}}
\sphinxAtStartPar
COPT general constraint object. General constraints are always associated with a particular model.
User creates a general constraint object by adding a general constraint to a model,
rather than by using constructor of GenConstr class.

\sphinxstepscope

\subsubsection{GenConstr.getIdx()}
\label{\detokenize{javaapi/GenConstr:genconstr-getidx}}\label{\detokenize{javaapi/GenConstr::doc}}\begin{quote}

\sphinxAtStartPar
Get the index of the general constraint.

\sphinxAtStartPar
\sphinxstylestrong{Synopsis}
\begin{quote}

\sphinxAtStartPar
\sphinxcode{\sphinxupquote{int getIdx()}}
\end{quote}

\sphinxAtStartPar
\sphinxstylestrong{Return}
\begin{quote}

\sphinxAtStartPar
index of the general constraint.
\end{quote}
\end{quote}

\subsubsection{GenConstr.getIIS()}
\label{\detokenize{javaapi/GenConstr:genconstr-getiis}}\begin{quote}

\sphinxAtStartPar
Get IIS status of the general constraint.

\sphinxAtStartPar
\sphinxstylestrong{Synopsis}
\begin{quote}

\sphinxAtStartPar
\sphinxcode{\sphinxupquote{int getIIS()}}
\end{quote}

\sphinxAtStartPar
\sphinxstylestrong{Return}
\begin{quote}

\sphinxAtStartPar
IIS status.
\end{quote}
\end{quote}

\subsubsection{GenConstr.getName()}
\label{\detokenize{javaapi/GenConstr:genconstr-getname}}\begin{quote}

\sphinxAtStartPar
Get name of general constraint.

\sphinxAtStartPar
\sphinxstylestrong{Synopsis}
\begin{quote}

\sphinxAtStartPar
\sphinxcode{\sphinxupquote{String getName()}}
\end{quote}

\sphinxAtStartPar
\sphinxstylestrong{Return}
\begin{quote}

\sphinxAtStartPar
the name of general constraint.
\end{quote}
\end{quote}

\subsubsection{GenConstr.remove()}
\label{\detokenize{javaapi/GenConstr:genconstr-remove}}\begin{quote}

\sphinxAtStartPar
Remove the general constraint from model.

\sphinxAtStartPar
\sphinxstylestrong{Synopsis}
\begin{quote}

\sphinxAtStartPar
\sphinxcode{\sphinxupquote{void remove()}}
\end{quote}
\end{quote}

\subsubsection{GenConstr.setName()}
\label{\detokenize{javaapi/GenConstr:genconstr-setname}}\begin{quote}

\sphinxAtStartPar
Set name for general constraint.

\sphinxAtStartPar
\sphinxstylestrong{Synopsis}
\begin{quote}

\sphinxAtStartPar
\sphinxcode{\sphinxupquote{void setName(String name)}}
\end{quote}

\sphinxAtStartPar
\sphinxstylestrong{Arguments}
\begin{quote}

\sphinxAtStartPar
\sphinxcode{\sphinxupquote{name}}: the name to set.
\end{quote}
\end{quote}

\subsection{GenConstrArray}
\label{\detokenize{javaapiref:genconstrarray}}\label{\detokenize{javaapiref:chapjavaapiref-genconstrarray}}
\sphinxAtStartPar
COPT general constraint array object. To store and access a set of Java {\hyperref[\detokenize{javaapiref:chapjavaapiref-genconstr}]{\sphinxcrossref{\DUrole{std,std-ref}{GenConstr}}}}
objects, Cardinal Optimizer provides Java GenConstrArray class, which defines the following methods.

\sphinxstepscope

\subsubsection{GenConstrArray.GenConstrArray()}
\label{\detokenize{javaapi/GenConstrArray:genconstrarray-genconstrarray}}\label{\detokenize{javaapi/GenConstrArray::doc}}\begin{quote}

\sphinxAtStartPar
Constructor of genconstrarray.

\sphinxAtStartPar
\sphinxstylestrong{Synopsis}
\begin{quote}

\sphinxAtStartPar
\sphinxcode{\sphinxupquote{GenConstrArray()}}
\end{quote}
\end{quote}

\subsubsection{GenConstrArray.getGenConstr()}
\label{\detokenize{javaapi/GenConstrArray:genconstrarray-getgenconstr}}\begin{quote}

\sphinxAtStartPar
Get idx\sphinxhyphen{}th general constraint object.

\sphinxAtStartPar
\sphinxstylestrong{Synopsis}
\begin{quote}

\sphinxAtStartPar
\sphinxcode{\sphinxupquote{GenConstr getGenConstr(int idx)}}
\end{quote}

\sphinxAtStartPar
\sphinxstylestrong{Arguments}
\begin{quote}

\sphinxAtStartPar
\sphinxcode{\sphinxupquote{idx}}: index of the general constraint.
\end{quote}

\sphinxAtStartPar
\sphinxstylestrong{Return}
\begin{quote}

\sphinxAtStartPar
general constraint object with index idx.
\end{quote}
\end{quote}

\subsubsection{GenConstrArray.pushBack()}
\label{\detokenize{javaapi/GenConstrArray:genconstrarray-pushback}}\begin{quote}

\sphinxAtStartPar
Add a general constraint object to general constraint array.

\sphinxAtStartPar
\sphinxstylestrong{Synopsis}
\begin{quote}

\sphinxAtStartPar
\sphinxcode{\sphinxupquote{void pushBack(GenConstr genconstr)}}
\end{quote}

\sphinxAtStartPar
\sphinxstylestrong{Arguments}
\begin{quote}

\sphinxAtStartPar
\sphinxcode{\sphinxupquote{genconstr}}: a general constraint object.
\end{quote}
\end{quote}

\subsubsection{GenConstrArray.reserve()}
\label{\detokenize{javaapi/GenConstrArray:genconstrarray-reserve}}\begin{quote}

\sphinxAtStartPar
Reserve capacity to contain at least n items.

\sphinxAtStartPar
\sphinxstylestrong{Synopsis}
\begin{quote}

\sphinxAtStartPar
\sphinxcode{\sphinxupquote{void reserve(int n)}}
\end{quote}

\sphinxAtStartPar
\sphinxstylestrong{Arguments}
\begin{quote}

\sphinxAtStartPar
\sphinxcode{\sphinxupquote{n}}: capacity number of general constraint object.
\end{quote}
\end{quote}

\subsubsection{GenConstrArray.size()}
\label{\detokenize{javaapi/GenConstrArray:genconstrarray-size}}\begin{quote}

\sphinxAtStartPar
Get the number of general constraint objects.

\sphinxAtStartPar
\sphinxstylestrong{Synopsis}
\begin{quote}

\sphinxAtStartPar
\sphinxcode{\sphinxupquote{int size()}}
\end{quote}

\sphinxAtStartPar
\sphinxstylestrong{Return}
\begin{quote}

\sphinxAtStartPar
number of general constraint objects.
\end{quote}
\end{quote}

\subsection{GenConstrBuilder}
\label{\detokenize{javaapiref:genconstrbuilder}}\label{\detokenize{javaapiref:chapjavaapiref-genconstrbuilder}}
\sphinxAtStartPar
COPT general constraint builder object. To help building a general constraint, given a binary variable
and associated value, a linear expression and constraint sense, Cardinal Optimizer provides Java
GenConstrBuilder class, which defines the following methods.

\sphinxstepscope

\subsubsection{GenConstrBuilder.GenConstrBuilder()}
\label{\detokenize{javaapi/GenConstrBuilder:genconstrbuilder-genconstrbuilder}}\label{\detokenize{javaapi/GenConstrBuilder::doc}}\begin{quote}

\sphinxAtStartPar
Constructor of genconstrbuilder.

\sphinxAtStartPar
\sphinxstylestrong{Synopsis}
\begin{quote}

\sphinxAtStartPar
\sphinxcode{\sphinxupquote{GenConstrBuilder()}}
\end{quote}
\end{quote}

\subsubsection{GenConstrBuilder.getBinVal()}
\label{\detokenize{javaapi/GenConstrBuilder:genconstrbuilder-getbinval}}\begin{quote}

\sphinxAtStartPar
Get binary value associated with general constraint.

\sphinxAtStartPar
\sphinxstylestrong{Synopsis}
\begin{quote}

\sphinxAtStartPar
\sphinxcode{\sphinxupquote{int getBinVal()}}
\end{quote}

\sphinxAtStartPar
\sphinxstylestrong{Return}
\begin{quote}

\sphinxAtStartPar
binary value.
\end{quote}
\end{quote}

\subsubsection{GenConstrBuilder.getBinVar()}
\label{\detokenize{javaapi/GenConstrBuilder:genconstrbuilder-getbinvar}}\begin{quote}

\sphinxAtStartPar
Get binary variable associated with general constraint.

\sphinxAtStartPar
\sphinxstylestrong{Synopsis}
\begin{quote}

\sphinxAtStartPar
\sphinxcode{\sphinxupquote{Var getBinVar()}}
\end{quote}

\sphinxAtStartPar
\sphinxstylestrong{Return}
\begin{quote}

\sphinxAtStartPar
binary vaiable object.
\end{quote}
\end{quote}

\subsubsection{GenConstrBuilder.getExpr()}
\label{\detokenize{javaapi/GenConstrBuilder:genconstrbuilder-getexpr}}\begin{quote}

\sphinxAtStartPar
Get expression associated with general constraint.

\sphinxAtStartPar
\sphinxstylestrong{Synopsis}
\begin{quote}

\sphinxAtStartPar
\sphinxcode{\sphinxupquote{Expr getExpr()}}
\end{quote}

\sphinxAtStartPar
\sphinxstylestrong{Return}
\begin{quote}

\sphinxAtStartPar
expression object.
\end{quote}
\end{quote}

\subsubsection{GenConstrBuilder.getIndType()}
\label{\detokenize{javaapi/GenConstrBuilder:genconstrbuilder-getindtype}}\begin{quote}

\sphinxAtStartPar
Get type of general constraint.

\sphinxAtStartPar
\sphinxstylestrong{Synopsis}
\begin{quote}

\sphinxAtStartPar
\sphinxcode{\sphinxupquote{int getIndType()}}
\end{quote}

\sphinxAtStartPar
\sphinxstylestrong{Return}
\begin{quote}

\sphinxAtStartPar
type of general constraint.
\end{quote}
\end{quote}

\subsubsection{GenConstrBuilder.getSense()}
\label{\detokenize{javaapi/GenConstrBuilder:genconstrbuilder-getsense}}\begin{quote}

\sphinxAtStartPar
Get sense associated with general constraint.

\sphinxAtStartPar
\sphinxstylestrong{Synopsis}
\begin{quote}

\sphinxAtStartPar
\sphinxcode{\sphinxupquote{char getSense()}}
\end{quote}

\sphinxAtStartPar
\sphinxstylestrong{Return}
\begin{quote}

\sphinxAtStartPar
constraint sense.
\end{quote}
\end{quote}

\subsubsection{GenConstrBuilder.set()}
\label{\detokenize{javaapi/GenConstrBuilder:genconstrbuilder-set}}\begin{quote}

\sphinxAtStartPar
Set binary variable, binary value, expression and sense of general constraint.

\sphinxAtStartPar
\sphinxstylestrong{Synopsis}
\begin{quote}

\sphinxAtStartPar
\sphinxcode{\sphinxupquote{void set(}}
\begin{quote}

\sphinxAtStartPar
\sphinxcode{\sphinxupquote{Var binvar,}}

\sphinxAtStartPar
\sphinxcode{\sphinxupquote{int binval,}}

\sphinxAtStartPar
\sphinxcode{\sphinxupquote{Expr expr,}}

\sphinxAtStartPar
\sphinxcode{\sphinxupquote{char sense)}}
\end{quote}
\end{quote}

\sphinxAtStartPar
\sphinxstylestrong{Arguments}
\begin{quote}

\sphinxAtStartPar
\sphinxcode{\sphinxupquote{binvar}}: binary variable.

\sphinxAtStartPar
\sphinxcode{\sphinxupquote{binval}}: binary value.

\sphinxAtStartPar
\sphinxcode{\sphinxupquote{expr}}: expression object.

\sphinxAtStartPar
\sphinxcode{\sphinxupquote{sense}}: general constraint sense.
\end{quote}
\end{quote}

\subsubsection{GenConstrBuilder.set()}
\label{\detokenize{javaapi/GenConstrBuilder:id1}}\begin{quote}

\sphinxAtStartPar
Set binary variable, binary value, expression and sense of general constraint.

\sphinxAtStartPar
\sphinxstylestrong{Synopsis}
\begin{quote}

\sphinxAtStartPar
\sphinxcode{\sphinxupquote{void set(}}
\begin{quote}

\sphinxAtStartPar
\sphinxcode{\sphinxupquote{Var binvar,}}

\sphinxAtStartPar
\sphinxcode{\sphinxupquote{int binval,}}

\sphinxAtStartPar
\sphinxcode{\sphinxupquote{Expr expr,}}

\sphinxAtStartPar
\sphinxcode{\sphinxupquote{char sense,}}

\sphinxAtStartPar
\sphinxcode{\sphinxupquote{int type)}}
\end{quote}
\end{quote}

\sphinxAtStartPar
\sphinxstylestrong{Arguments}
\begin{quote}

\sphinxAtStartPar
\sphinxcode{\sphinxupquote{binvar}}: binary variable.

\sphinxAtStartPar
\sphinxcode{\sphinxupquote{binval}}: binary value.

\sphinxAtStartPar
\sphinxcode{\sphinxupquote{expr}}: expression object.

\sphinxAtStartPar
\sphinxcode{\sphinxupquote{sense}}: general constraint sense.

\sphinxAtStartPar
\sphinxcode{\sphinxupquote{type}}: type of general constraint.
\end{quote}
\end{quote}

\subsection{GenConstrBuilderArray}
\label{\detokenize{javaapiref:genconstrbuilderarray}}\label{\detokenize{javaapiref:chapjavaapiref-genconstrbuilderarray}}
\sphinxAtStartPar
COPT general constraint builder array object. To store and access a set of Java
{\hyperref[\detokenize{javaapiref:chapjavaapiref-genconstrbuilder}]{\sphinxcrossref{\DUrole{std,std-ref}{GenConstrBuilder}}}} objects, Cardinal Optimizer provides Java
GenConstrBuilderArray class, which defines the following methods.

\sphinxstepscope

\subsubsection{GenConstrBuilderArray.GenConstrBuilderArray()}
\label{\detokenize{javaapi/GenConstrBuilderArray:genconstrbuilderarray-genconstrbuilderarray}}\label{\detokenize{javaapi/GenConstrBuilderArray::doc}}\begin{quote}

\sphinxAtStartPar
Constructor of genconstrbuilderarray.

\sphinxAtStartPar
\sphinxstylestrong{Synopsis}
\begin{quote}

\sphinxAtStartPar
\sphinxcode{\sphinxupquote{GenConstrBuilderArray()}}
\end{quote}
\end{quote}

\subsubsection{GenConstrBuilderArray.getBuilder()}
\label{\detokenize{javaapi/GenConstrBuilderArray:genconstrbuilderarray-getbuilder}}\begin{quote}

\sphinxAtStartPar
Get idx\sphinxhyphen{}th general constraint builder object.

\sphinxAtStartPar
\sphinxstylestrong{Synopsis}
\begin{quote}

\sphinxAtStartPar
\sphinxcode{\sphinxupquote{GenConstrBuilder getBuilder(int idx)}}
\end{quote}

\sphinxAtStartPar
\sphinxstylestrong{Arguments}
\begin{quote}

\sphinxAtStartPar
\sphinxcode{\sphinxupquote{idx}}: index of the general constraint builder.
\end{quote}

\sphinxAtStartPar
\sphinxstylestrong{Return}
\begin{quote}

\sphinxAtStartPar
general constraint builder object with index idx.
\end{quote}
\end{quote}

\subsubsection{GenConstrBuilderArray.pushBack()}
\label{\detokenize{javaapi/GenConstrBuilderArray:genconstrbuilderarray-pushback}}\begin{quote}

\sphinxAtStartPar
Add a general constraint builder object to general constraint builder array.

\sphinxAtStartPar
\sphinxstylestrong{Synopsis}
\begin{quote}

\sphinxAtStartPar
\sphinxcode{\sphinxupquote{void pushBack(GenConstrBuilder builder)}}
\end{quote}

\sphinxAtStartPar
\sphinxstylestrong{Arguments}
\begin{quote}

\sphinxAtStartPar
\sphinxcode{\sphinxupquote{builder}}: a general constraint builder object.
\end{quote}
\end{quote}

\subsubsection{GenConstrBuilderArray.size()}
\label{\detokenize{javaapi/GenConstrBuilderArray:genconstrbuilderarray-size}}\begin{quote}

\sphinxAtStartPar
Get the number of general constraint builder objects.

\sphinxAtStartPar
\sphinxstylestrong{Synopsis}
\begin{quote}

\sphinxAtStartPar
\sphinxcode{\sphinxupquote{int size()}}
\end{quote}

\sphinxAtStartPar
\sphinxstylestrong{Return}
\begin{quote}

\sphinxAtStartPar
number of general constraint builder objects.
\end{quote}
\end{quote}

\subsection{Cone}
\label{\detokenize{javaapiref:cone}}\label{\detokenize{javaapiref:chapjavaapiref-cone}}
\sphinxAtStartPar
COPT cone constraint object. Cone constraints are always associated with a particular model.
User creates a cone constraint object by adding a cone constraint to a model,
rather than by using constructor of Cone class.

\sphinxAtStartPar
A cone constraint can be regular or rotated (\sphinxcode{\sphinxupquote{COPT\_CONE\_QUAD}} or \sphinxcode{\sphinxupquote{COPT\_CONE\_RQUAD}}).

\sphinxstepscope

\subsubsection{Cone.getIdx()}
\label{\detokenize{javaapi/Cone:cone-getidx}}\label{\detokenize{javaapi/Cone::doc}}\begin{quote}

\sphinxAtStartPar
Get the index of a cone constraint.

\sphinxAtStartPar
\sphinxstylestrong{Synopsis}
\begin{quote}

\sphinxAtStartPar
\sphinxcode{\sphinxupquote{int getIdx()}}
\end{quote}

\sphinxAtStartPar
\sphinxstylestrong{Return}
\begin{quote}

\sphinxAtStartPar
index of a cone constraint.
\end{quote}
\end{quote}

\subsubsection{Cone.remove()}
\label{\detokenize{javaapi/Cone:cone-remove}}\begin{quote}

\sphinxAtStartPar
Remove the cone constraint from model.

\sphinxAtStartPar
\sphinxstylestrong{Synopsis}
\begin{quote}

\sphinxAtStartPar
\sphinxcode{\sphinxupquote{void remove()}}
\end{quote}
\end{quote}

\subsection{ConeArray}
\label{\detokenize{javaapiref:conearray}}\label{\detokenize{javaapiref:chapjavaapiref-conearray}}
\sphinxAtStartPar
COPT cone constraint array object. To store and access a set of Java {\hyperref[\detokenize{javaapiref:chapjavaapiref-cone}]{\sphinxcrossref{\DUrole{std,std-ref}{Cone}}}}
objects, Cardinal Optimizer provides Java ConeArray class, which defines the following methods.

\sphinxstepscope

\subsubsection{ConeArray.ConeArray()}
\label{\detokenize{javaapi/ConeArray:conearray-conearray}}\label{\detokenize{javaapi/ConeArray::doc}}\begin{quote}

\sphinxAtStartPar
Constructor of ConeArray object.

\sphinxAtStartPar
\sphinxstylestrong{Synopsis}
\begin{quote}

\sphinxAtStartPar
\sphinxcode{\sphinxupquote{ConeArray()}}
\end{quote}
\end{quote}

\subsubsection{ConeArray.getCone()}
\label{\detokenize{javaapi/ConeArray:conearray-getcone}}\begin{quote}

\sphinxAtStartPar
Get idx\sphinxhyphen{}th cone object.

\sphinxAtStartPar
\sphinxstylestrong{Synopsis}
\begin{quote}

\sphinxAtStartPar
\sphinxcode{\sphinxupquote{Cone getCone(int idx)}}
\end{quote}

\sphinxAtStartPar
\sphinxstylestrong{Arguments}
\begin{quote}

\sphinxAtStartPar
\sphinxcode{\sphinxupquote{idx}}: index of cone.
\end{quote}

\sphinxAtStartPar
\sphinxstylestrong{Return}
\begin{quote}

\sphinxAtStartPar
cone object with index idx.
\end{quote}
\end{quote}

\subsubsection{ConeArray.pushBack()}
\label{\detokenize{javaapi/ConeArray:conearray-pushback}}\begin{quote}

\sphinxAtStartPar
Add a cone constraint object to cone constraint array.

\sphinxAtStartPar
\sphinxstylestrong{Synopsis}
\begin{quote}

\sphinxAtStartPar
\sphinxcode{\sphinxupquote{void pushBack(Cone cone)}}
\end{quote}

\sphinxAtStartPar
\sphinxstylestrong{Arguments}
\begin{quote}

\sphinxAtStartPar
\sphinxcode{\sphinxupquote{cone}}: a cone constraint object.
\end{quote}
\end{quote}

\subsubsection{ConeArray.size()}
\label{\detokenize{javaapi/ConeArray:conearray-size}}\begin{quote}

\sphinxAtStartPar
Get the number of cone constraint objects.

\sphinxAtStartPar
\sphinxstylestrong{Synopsis}
\begin{quote}

\sphinxAtStartPar
\sphinxcode{\sphinxupquote{int size()}}
\end{quote}

\sphinxAtStartPar
\sphinxstylestrong{Return}
\begin{quote}

\sphinxAtStartPar
number of cone constraint objects.
\end{quote}
\end{quote}

\subsection{ConeBuilder}
\label{\detokenize{javaapiref:conebuilder}}\label{\detokenize{javaapiref:chapjavaapiref-conebuilder}}
\sphinxAtStartPar
COPT cone constraint builder object. To help building a cone constraint, given the cone type and
a set of variables, Cardinal Optimizer provides Java ConeBuilder class, which defines the following methods.

\sphinxstepscope

\subsubsection{ConeBuilder.ConeBuilder()}
\label{\detokenize{javaapi/ConeBuilder:conebuilder-conebuilder}}\label{\detokenize{javaapi/ConeBuilder::doc}}\begin{quote}

\sphinxAtStartPar
Constructor of ConeBuilder object.

\sphinxAtStartPar
\sphinxstylestrong{Synopsis}
\begin{quote}

\sphinxAtStartPar
\sphinxcode{\sphinxupquote{ConeBuilder()}}
\end{quote}
\end{quote}

\subsubsection{ConeBuilder.getSize()}
\label{\detokenize{javaapi/ConeBuilder:conebuilder-getsize}}\begin{quote}

\sphinxAtStartPar
Get number of variables in a cone constraint.

\sphinxAtStartPar
\sphinxstylestrong{Synopsis}
\begin{quote}

\sphinxAtStartPar
\sphinxcode{\sphinxupquote{int getSize()}}
\end{quote}

\sphinxAtStartPar
\sphinxstylestrong{Return}
\begin{quote}

\sphinxAtStartPar
number of vars.
\end{quote}
\end{quote}

\subsubsection{ConeBuilder.getType()}
\label{\detokenize{javaapi/ConeBuilder:conebuilder-gettype}}\begin{quote}

\sphinxAtStartPar
Get type of a cone constraint.

\sphinxAtStartPar
\sphinxstylestrong{Synopsis}
\begin{quote}

\sphinxAtStartPar
\sphinxcode{\sphinxupquote{int getType()}}
\end{quote}

\sphinxAtStartPar
\sphinxstylestrong{Return}
\begin{quote}

\sphinxAtStartPar
type of a cone constraint.
\end{quote}
\end{quote}

\subsubsection{ConeBuilder.getVar()}
\label{\detokenize{javaapi/ConeBuilder:conebuilder-getvar}}\begin{quote}

\sphinxAtStartPar
Get idx\sphinxhyphen{}th variable in a cone constraint.

\sphinxAtStartPar
\sphinxstylestrong{Synopsis}
\begin{quote}

\sphinxAtStartPar
\sphinxcode{\sphinxupquote{Var getVar(int idx)}}
\end{quote}

\sphinxAtStartPar
\sphinxstylestrong{Arguments}
\begin{quote}

\sphinxAtStartPar
\sphinxcode{\sphinxupquote{idx}}: index of variables.
\end{quote}

\sphinxAtStartPar
\sphinxstylestrong{Return}
\begin{quote}

\sphinxAtStartPar
the idx\sphinxhyphen{}th variable in a cone constraint.
\end{quote}
\end{quote}

\subsubsection{ConeBuilder.getVars()}
\label{\detokenize{javaapi/ConeBuilder:conebuilder-getvars}}\begin{quote}

\sphinxAtStartPar
Get all variables in a cone constraint.

\sphinxAtStartPar
\sphinxstylestrong{Synopsis}
\begin{quote}

\sphinxAtStartPar
\sphinxcode{\sphinxupquote{VarArray getVars()}}
\end{quote}

\sphinxAtStartPar
\sphinxstylestrong{Return}
\begin{quote}

\sphinxAtStartPar
variables in a cone constraint.
\end{quote}
\end{quote}

\subsubsection{ConeBuilder.set()}
\label{\detokenize{javaapi/ConeBuilder:conebuilder-set}}\begin{quote}

\sphinxAtStartPar
Set variables and type of a cone constraint.

\sphinxAtStartPar
\sphinxstylestrong{Synopsis}
\begin{quote}

\sphinxAtStartPar
\sphinxcode{\sphinxupquote{void set(VarArray vars, int type)}}
\end{quote}

\sphinxAtStartPar
\sphinxstylestrong{Arguments}
\begin{quote}

\sphinxAtStartPar
\sphinxcode{\sphinxupquote{vars}}: variable array object.

\sphinxAtStartPar
\sphinxcode{\sphinxupquote{type}}: type of a cone constraint.
\end{quote}
\end{quote}

\subsection{ConeBuilderArray}
\label{\detokenize{javaapiref:conebuilderarray}}\label{\detokenize{javaapiref:chapjavaapiref-conebuilderarray}}
\sphinxAtStartPar
COPT cone constraint builder array object. To store and access a set of Java {\hyperref[\detokenize{javaapiref:chapjavaapiref-conebuilder}]{\sphinxcrossref{\DUrole{std,std-ref}{ConeBuilder}}}}
objects, Cardinal Optimizer provides Java ConeBuilderArray class, which defines the following methods.

\sphinxstepscope

\subsubsection{ConeBuilderArray.ConeBuilderArray()}
\label{\detokenize{javaapi/ConeBuilderArray:conebuilderarray-conebuilderarray}}\label{\detokenize{javaapi/ConeBuilderArray::doc}}\begin{quote}

\sphinxAtStartPar
Constructor of ConeBuilderArray object.

\sphinxAtStartPar
\sphinxstylestrong{Synopsis}
\begin{quote}

\sphinxAtStartPar
\sphinxcode{\sphinxupquote{ConeBuilderArray()}}
\end{quote}
\end{quote}

\subsubsection{ConeBuilderArray.getBuilder()}
\label{\detokenize{javaapi/ConeBuilderArray:conebuilderarray-getbuilder}}\begin{quote}

\sphinxAtStartPar
Get idx\sphinxhyphen{}th cone constraint builder object.

\sphinxAtStartPar
\sphinxstylestrong{Synopsis}
\begin{quote}

\sphinxAtStartPar
\sphinxcode{\sphinxupquote{ConeBuilder getBuilder(int idx)}}
\end{quote}

\sphinxAtStartPar
\sphinxstylestrong{Arguments}
\begin{quote}

\sphinxAtStartPar
\sphinxcode{\sphinxupquote{idx}}: index of the cone constraint builder.
\end{quote}

\sphinxAtStartPar
\sphinxstylestrong{Return}
\begin{quote}

\sphinxAtStartPar
cone constraint builder object with index idx.
\end{quote}
\end{quote}

\subsubsection{ConeBuilderArray.pushBack()}
\label{\detokenize{javaapi/ConeBuilderArray:conebuilderarray-pushback}}\begin{quote}

\sphinxAtStartPar
Add a cone constraint builder object to cone constraint builder array.

\sphinxAtStartPar
\sphinxstylestrong{Synopsis}
\begin{quote}

\sphinxAtStartPar
\sphinxcode{\sphinxupquote{void pushBack(ConeBuilder builder)}}
\end{quote}

\sphinxAtStartPar
\sphinxstylestrong{Arguments}
\begin{quote}

\sphinxAtStartPar
\sphinxcode{\sphinxupquote{builder}}: a cone constraint builder object.
\end{quote}
\end{quote}

\subsubsection{ConeBuilderArray.size()}
\label{\detokenize{javaapi/ConeBuilderArray:conebuilderarray-size}}\begin{quote}

\sphinxAtStartPar
Get the number of cone constraint builder objects.

\sphinxAtStartPar
\sphinxstylestrong{Synopsis}
\begin{quote}

\sphinxAtStartPar
\sphinxcode{\sphinxupquote{int size()}}
\end{quote}

\sphinxAtStartPar
\sphinxstylestrong{Return}
\begin{quote}

\sphinxAtStartPar
number of cone constraint builder objects.
\end{quote}
\end{quote}

\subsection{ExpCone}
\label{\detokenize{javaapiref:expcone}}\label{\detokenize{javaapiref:chapjavaapiref-expcone}}
\sphinxAtStartPar
COPT exponential cone constraint object. ExpCone constraints are always associated with a particular model.
User creates an exponential cone constraint object by adding an exponential cone constraint to a model,
rather than by using constructor of ExpCone class.

\sphinxstepscope

\subsubsection{ExpCone.getIdx()}
\label{\detokenize{javaapi/ExpCone:expcone-getidx}}\label{\detokenize{javaapi/ExpCone::doc}}\begin{quote}

\sphinxAtStartPar
Get the index of an exponential cone constraint.

\sphinxAtStartPar
\sphinxstylestrong{Synopsis}
\begin{quote}

\sphinxAtStartPar
\sphinxcode{\sphinxupquote{int getIdx()}}
\end{quote}

\sphinxAtStartPar
\sphinxstylestrong{Return}
\begin{quote}

\sphinxAtStartPar
index of an exponential cone constraint.
\end{quote}
\end{quote}

\subsubsection{ExpCone.remove()}
\label{\detokenize{javaapi/ExpCone:expcone-remove}}\begin{quote}

\sphinxAtStartPar
Remove the exponential cone constraint from model.

\sphinxAtStartPar
\sphinxstylestrong{Synopsis}
\begin{quote}

\sphinxAtStartPar
\sphinxcode{\sphinxupquote{void remove()}}
\end{quote}
\end{quote}

\subsection{ExpConeArray}
\label{\detokenize{javaapiref:expconearray}}\label{\detokenize{javaapiref:chapjavaapiref-expconearray}}
\sphinxAtStartPar
COPT exponential cone constraint array object. To store and access a set of Java {\hyperref[\detokenize{javaapiref:chapjavaapiref-expcone}]{\sphinxcrossref{\DUrole{std,std-ref}{ExpCone}}}}
objects, Cardinal Optimizer provides Java ExpConeArray class, which defines the following methods.

\sphinxstepscope

\subsubsection{ExpConeArray.ExpConeArray()}
\label{\detokenize{javaapi/ExpConeArray:expconearray-expconearray}}\label{\detokenize{javaapi/ExpConeArray::doc}}\begin{quote}

\sphinxAtStartPar
Constructor of ExpConeArray object.

\sphinxAtStartPar
\sphinxstylestrong{Synopsis}
\begin{quote}

\sphinxAtStartPar
\sphinxcode{\sphinxupquote{ExpConeArray()}}
\end{quote}
\end{quote}

\subsubsection{ExpConeArray.getCone()}
\label{\detokenize{javaapi/ExpConeArray:expconearray-getcone}}\begin{quote}

\sphinxAtStartPar
Get idx\sphinxhyphen{}th exponential cone object.

\sphinxAtStartPar
\sphinxstylestrong{Synopsis}
\begin{quote}

\sphinxAtStartPar
\sphinxcode{\sphinxupquote{ExpCone getCone(int idx)}}
\end{quote}

\sphinxAtStartPar
\sphinxstylestrong{Arguments}
\begin{quote}

\sphinxAtStartPar
\sphinxcode{\sphinxupquote{idx}}: index of exponential cone.
\end{quote}

\sphinxAtStartPar
\sphinxstylestrong{Return}
\begin{quote}

\sphinxAtStartPar
exponential cone object with index idx.
\end{quote}
\end{quote}

\subsubsection{ExpConeArray.pushBack()}
\label{\detokenize{javaapi/ExpConeArray:expconearray-pushback}}\begin{quote}

\sphinxAtStartPar
Add an exponential cone constraint object to exponential cone constraint array.

\sphinxAtStartPar
\sphinxstylestrong{Synopsis}
\begin{quote}

\sphinxAtStartPar
\sphinxcode{\sphinxupquote{void pushBack(ExpCone cone)}}
\end{quote}

\sphinxAtStartPar
\sphinxstylestrong{Arguments}
\begin{quote}

\sphinxAtStartPar
\sphinxcode{\sphinxupquote{cone}}: an exponential cone constraint object.
\end{quote}
\end{quote}

\subsubsection{ExpConeArray.size()}
\label{\detokenize{javaapi/ExpConeArray:expconearray-size}}\begin{quote}

\sphinxAtStartPar
Get the number of exponential cone constraint objects.

\sphinxAtStartPar
\sphinxstylestrong{Synopsis}
\begin{quote}

\sphinxAtStartPar
\sphinxcode{\sphinxupquote{int size()}}
\end{quote}

\sphinxAtStartPar
\sphinxstylestrong{Return}
\begin{quote}

\sphinxAtStartPar
number of exponential cone constraint objects.
\end{quote}
\end{quote}

\subsection{ExpConeBuilder}
\label{\detokenize{javaapiref:expconebuilder}}\label{\detokenize{javaapiref:chapjavaapiref-expconebuilder}}
\sphinxAtStartPar
COPT exponential cone constraint builder object. To help building an exponential cone constraint, given the exponential
cone type and a set of variables, Cardinal Optimizer provides Java ExpConeBuilder class, which defines the following methods.

\sphinxstepscope

\subsubsection{ExpConeBuilder.ExpConeBuilder()}
\label{\detokenize{javaapi/ExpConeBuilder:expconebuilder-expconebuilder}}\label{\detokenize{javaapi/ExpConeBuilder::doc}}\begin{quote}

\sphinxAtStartPar
Constructor of ExpConeBuilder object.

\sphinxAtStartPar
\sphinxstylestrong{Synopsis}
\begin{quote}

\sphinxAtStartPar
\sphinxcode{\sphinxupquote{ExpConeBuilder()}}
\end{quote}
\end{quote}

\subsubsection{ExpConeBuilder.getSize()}
\label{\detokenize{javaapi/ExpConeBuilder:expconebuilder-getsize}}\begin{quote}

\sphinxAtStartPar
Get number of variables in an exponential cone constraint.

\sphinxAtStartPar
\sphinxstylestrong{Synopsis}
\begin{quote}

\sphinxAtStartPar
\sphinxcode{\sphinxupquote{int getSize()}}
\end{quote}

\sphinxAtStartPar
\sphinxstylestrong{Return}
\begin{quote}

\sphinxAtStartPar
number of vars.
\end{quote}
\end{quote}

\subsubsection{ExpConeBuilder.getType()}
\label{\detokenize{javaapi/ExpConeBuilder:expconebuilder-gettype}}\begin{quote}

\sphinxAtStartPar
Get type of an exponential cone constraint.

\sphinxAtStartPar
\sphinxstylestrong{Synopsis}
\begin{quote}

\sphinxAtStartPar
\sphinxcode{\sphinxupquote{int getType()}}
\end{quote}

\sphinxAtStartPar
\sphinxstylestrong{Return}
\begin{quote}

\sphinxAtStartPar
type of an exponential cone constraint.
\end{quote}
\end{quote}

\subsubsection{ExpConeBuilder.getVar()}
\label{\detokenize{javaapi/ExpConeBuilder:expconebuilder-getvar}}\begin{quote}

\sphinxAtStartPar
Get idx\sphinxhyphen{}th variable in an exponential cone constraint.

\sphinxAtStartPar
\sphinxstylestrong{Synopsis}
\begin{quote}

\sphinxAtStartPar
\sphinxcode{\sphinxupquote{Var getVar(int idx)}}
\end{quote}

\sphinxAtStartPar
\sphinxstylestrong{Arguments}
\begin{quote}

\sphinxAtStartPar
\sphinxcode{\sphinxupquote{idx}}: index of variables.
\end{quote}

\sphinxAtStartPar
\sphinxstylestrong{Return}
\begin{quote}

\sphinxAtStartPar
the idx\sphinxhyphen{}th variable in an exponential cone constraint.
\end{quote}
\end{quote}

\subsubsection{ExpConeBuilder.getVars()}
\label{\detokenize{javaapi/ExpConeBuilder:expconebuilder-getvars}}\begin{quote}

\sphinxAtStartPar
Get all variables in an exponential cone constraint.

\sphinxAtStartPar
\sphinxstylestrong{Synopsis}
\begin{quote}

\sphinxAtStartPar
\sphinxcode{\sphinxupquote{VarArray getVars()}}
\end{quote}

\sphinxAtStartPar
\sphinxstylestrong{Return}
\begin{quote}

\sphinxAtStartPar
variables in an exponential cone constraint.
\end{quote}
\end{quote}

\subsubsection{ExpConeBuilder.set()}
\label{\detokenize{javaapi/ExpConeBuilder:expconebuilder-set}}\begin{quote}

\sphinxAtStartPar
Set variables and type of an exponential cone constraint.

\sphinxAtStartPar
\sphinxstylestrong{Synopsis}
\begin{quote}

\sphinxAtStartPar
\sphinxcode{\sphinxupquote{void set(VarArray vars, int type)}}
\end{quote}

\sphinxAtStartPar
\sphinxstylestrong{Arguments}
\begin{quote}

\sphinxAtStartPar
\sphinxcode{\sphinxupquote{vars}}: variable array object.

\sphinxAtStartPar
\sphinxcode{\sphinxupquote{type}}: type of an exponential cone constraint.
\end{quote}
\end{quote}

\subsection{ExpConeBuilderArray}
\label{\detokenize{javaapiref:expconebuilderarray}}\label{\detokenize{javaapiref:chapjavaapiref-expconebuilderarray}}
\sphinxAtStartPar
COPT exponential cone constraint builder array object. To store and access a set of Java {\hyperref[\detokenize{javaapiref:chapjavaapiref-expconebuilder}]{\sphinxcrossref{\DUrole{std,std-ref}{ExpConeBuilder}}}}
objects, Cardinal Optimizer provides Java ExpConeBuilderArray class, which defines the following methods.

\sphinxstepscope

\subsubsection{ExpConeBuilderArray.ExpConeBuilderArray()}
\label{\detokenize{javaapi/ExpConeBuilderArray:expconebuilderarray-expconebuilderarray}}\label{\detokenize{javaapi/ExpConeBuilderArray::doc}}\begin{quote}

\sphinxAtStartPar
Constructor of ExpConeBuilderArray object.

\sphinxAtStartPar
\sphinxstylestrong{Synopsis}
\begin{quote}

\sphinxAtStartPar
\sphinxcode{\sphinxupquote{ExpConeBuilderArray()}}
\end{quote}
\end{quote}

\subsubsection{ExpConeBuilderArray.getBuilder()}
\label{\detokenize{javaapi/ExpConeBuilderArray:expconebuilderarray-getbuilder}}\begin{quote}

\sphinxAtStartPar
Get idx\sphinxhyphen{}th exponential cone constraint builder object.

\sphinxAtStartPar
\sphinxstylestrong{Synopsis}
\begin{quote}

\sphinxAtStartPar
\sphinxcode{\sphinxupquote{ExpConeBuilder getBuilder(int idx)}}
\end{quote}

\sphinxAtStartPar
\sphinxstylestrong{Arguments}
\begin{quote}

\sphinxAtStartPar
\sphinxcode{\sphinxupquote{idx}}: index of the exponential cone constraint builder.
\end{quote}

\sphinxAtStartPar
\sphinxstylestrong{Return}
\begin{quote}

\sphinxAtStartPar
exponential cone constraint builder object with index idx.
\end{quote}
\end{quote}

\subsubsection{ExpConeBuilderArray.pushBack()}
\label{\detokenize{javaapi/ExpConeBuilderArray:expconebuilderarray-pushback}}\begin{quote}

\sphinxAtStartPar
Add an exponential cone constraint builder object to exponential cone constraint builder array.

\sphinxAtStartPar
\sphinxstylestrong{Synopsis}
\begin{quote}

\sphinxAtStartPar
\sphinxcode{\sphinxupquote{void pushBack(ExpConeBuilder builder)}}
\end{quote}

\sphinxAtStartPar
\sphinxstylestrong{Arguments}
\begin{quote}

\sphinxAtStartPar
\sphinxcode{\sphinxupquote{builder}}: an exponential cone constraint builder object.
\end{quote}
\end{quote}

\subsubsection{ExpConeBuilderArray.size()}
\label{\detokenize{javaapi/ExpConeBuilderArray:expconebuilderarray-size}}\begin{quote}

\sphinxAtStartPar
Get the number of exponential cone constraint builder objects.

\sphinxAtStartPar
\sphinxstylestrong{Synopsis}
\begin{quote}

\sphinxAtStartPar
\sphinxcode{\sphinxupquote{int size()}}
\end{quote}

\sphinxAtStartPar
\sphinxstylestrong{Return}
\begin{quote}

\sphinxAtStartPar
number of exponential cone constraint builder objects.
\end{quote}
\end{quote}

\subsection{AffineCone Class}
\label{\detokenize{javaapiref:affinecone-class}}\label{\detokenize{javaapiref:chapjavaapiref-affinecone}}
\sphinxAtStartPar
The \sphinxtitleref{AffineCone} class in COPT encapsulates operations related to affine cones.
The following methods are provided:

\sphinxstepscope

\subsubsection{AffineCone.getIdx()}
\label{\detokenize{javaapi/AffineCone:affinecone-getidx}}\label{\detokenize{javaapi/AffineCone::doc}}\begin{quote}

\sphinxAtStartPar
Get the index of an affine cone constraint.

\sphinxAtStartPar
\sphinxstylestrong{Synopsis}
\begin{quote}

\sphinxAtStartPar
\sphinxcode{\sphinxupquote{int getIdx()}}
\end{quote}

\sphinxAtStartPar
\sphinxstylestrong{Return}
\begin{quote}

\sphinxAtStartPar
index of an affine cone constraint.
\end{quote}
\end{quote}

\subsubsection{AffineCone.getName()}
\label{\detokenize{javaapi/AffineCone:affinecone-getname}}\begin{quote}

\sphinxAtStartPar
Get name of the affine cone.

\sphinxAtStartPar
\sphinxstylestrong{Synopsis}
\begin{quote}

\sphinxAtStartPar
\sphinxcode{\sphinxupquote{String getName()}}
\end{quote}

\sphinxAtStartPar
\sphinxstylestrong{Return}
\begin{quote}

\sphinxAtStartPar
affine cone name.
\end{quote}
\end{quote}

\subsubsection{AffineCone.remove()}
\label{\detokenize{javaapi/AffineCone:affinecone-remove}}\begin{quote}

\sphinxAtStartPar
Remove the affine cone constraint from model.

\sphinxAtStartPar
\sphinxstylestrong{Synopsis}
\begin{quote}

\sphinxAtStartPar
\sphinxcode{\sphinxupquote{void remove()}}
\end{quote}
\end{quote}

\subsubsection{AffineCone.setName()}
\label{\detokenize{javaapi/AffineCone:affinecone-setname}}\begin{quote}

\sphinxAtStartPar
Set name of the affine cone.

\sphinxAtStartPar
\sphinxstylestrong{Synopsis}
\begin{quote}

\sphinxAtStartPar
\sphinxcode{\sphinxupquote{void setName(String name)}}
\end{quote}

\sphinxAtStartPar
\sphinxstylestrong{Arguments}
\begin{quote}

\sphinxAtStartPar
\sphinxcode{\sphinxupquote{name}}: affinecone name.
\end{quote}
\end{quote}

\subsection{AffineConeArray Class}
\label{\detokenize{javaapiref:affineconearray-class}}\label{\detokenize{javaapiref:chapjavaapiref-affineconearray}}
\sphinxAtStartPar
To facilitate user operations on a group of Java {\hyperref[\detokenize{javaapiref:chapjavaapiref-affinecone}]{\sphinxcrossref{\DUrole{std,std-ref}{AffineCone Class}}}} objects,
the Java interface of COPT introduces the \sphinxtitleref{AffineConeArray} class.
The following methods are provided:

\sphinxstepscope

\subsubsection{AffineConeArray.AffineConeArray()}
\label{\detokenize{javaapi/AffineConeArray:affineconearray-affineconearray}}\label{\detokenize{javaapi/AffineConeArray::doc}}\begin{quote}

\sphinxAtStartPar
Constructor of AffineConeArray object.

\sphinxAtStartPar
\sphinxstylestrong{Synopsis}
\begin{quote}

\sphinxAtStartPar
\sphinxcode{\sphinxupquote{AffineConeArray()}}
\end{quote}
\end{quote}

\subsubsection{AffineConeArray.getCone()}
\label{\detokenize{javaapi/AffineConeArray:affineconearray-getcone}}\begin{quote}

\sphinxAtStartPar
Get idx\sphinxhyphen{}th affine cone object.

\sphinxAtStartPar
\sphinxstylestrong{Synopsis}
\begin{quote}

\sphinxAtStartPar
\sphinxcode{\sphinxupquote{AffineCone getCone(int idx)}}
\end{quote}

\sphinxAtStartPar
\sphinxstylestrong{Arguments}
\begin{quote}

\sphinxAtStartPar
\sphinxcode{\sphinxupquote{idx}}: index of affine cone.
\end{quote}

\sphinxAtStartPar
\sphinxstylestrong{Return}
\begin{quote}

\sphinxAtStartPar
affine cone object with index idx.
\end{quote}
\end{quote}

\subsubsection{AffineConeArray.pushBack()}
\label{\detokenize{javaapi/AffineConeArray:affineconearray-pushback}}\begin{quote}

\sphinxAtStartPar
Add an affine cone constraint object to affine cone constraint array.

\sphinxAtStartPar
\sphinxstylestrong{Synopsis}
\begin{quote}

\sphinxAtStartPar
\sphinxcode{\sphinxupquote{void pushBack(AffineCone cone)}}
\end{quote}

\sphinxAtStartPar
\sphinxstylestrong{Arguments}
\begin{quote}

\sphinxAtStartPar
\sphinxcode{\sphinxupquote{cone}}: an affine cone constraint object.
\end{quote}
\end{quote}

\subsubsection{AffineConeArray.size()}
\label{\detokenize{javaapi/AffineConeArray:affineconearray-size}}\begin{quote}

\sphinxAtStartPar
Get the number of affine cone constraint objects.

\sphinxAtStartPar
\sphinxstylestrong{Synopsis}
\begin{quote}

\sphinxAtStartPar
\sphinxcode{\sphinxupquote{int size()}}
\end{quote}

\sphinxAtStartPar
\sphinxstylestrong{Return}
\begin{quote}

\sphinxAtStartPar
number of affine cone constraint objects.
\end{quote}
\end{quote}

\subsection{AffineConeBuilder Class}
\label{\detokenize{javaapiref:affineconebuilder-class}}\label{\detokenize{javaapiref:chapjavaapiref-affineconebuilder}}
\sphinxAtStartPar
The \sphinxtitleref{AffineConeBuilder} class in COPT encapsulates the builder for constructing affine cones.
The following methods are provided:

\sphinxstepscope

\subsubsection{AffineConeBuilder.AffineConeBuilder()}
\label{\detokenize{javaapi/AffineConeBuilder:affineconebuilder-affineconebuilder}}\label{\detokenize{javaapi/AffineConeBuilder::doc}}\begin{quote}

\sphinxAtStartPar
Constructor of AffineConeBuilder object.

\sphinxAtStartPar
\sphinxstylestrong{Synopsis}
\begin{quote}

\sphinxAtStartPar
\sphinxcode{\sphinxupquote{AffineConeBuilder()}}
\end{quote}
\end{quote}

\subsubsection{AffineConeBuilder.getExpr()}
\label{\detokenize{javaapi/AffineConeBuilder:affineconebuilder-getexpr}}\begin{quote}

\sphinxAtStartPar
Get idx\sphinxhyphen{}th linear expression in an affine cone constraint.

\sphinxAtStartPar
\sphinxstylestrong{Synopsis}
\begin{quote}

\sphinxAtStartPar
\sphinxcode{\sphinxupquote{Expr getExpr(int idx)}}
\end{quote}

\sphinxAtStartPar
\sphinxstylestrong{Arguments}
\begin{quote}

\sphinxAtStartPar
\sphinxcode{\sphinxupquote{idx}}: index of linear expression.
\end{quote}

\sphinxAtStartPar
\sphinxstylestrong{Return}
\begin{quote}

\sphinxAtStartPar
the idx\sphinxhyphen{}th linear expression in an affine cone constraint.
\end{quote}
\end{quote}

\subsubsection{AffineConeBuilder.getExprs()}
\label{\detokenize{javaapi/AffineConeBuilder:affineconebuilder-getexprs}}\begin{quote}

\sphinxAtStartPar
Get all linear expressions in an affine cone constraint.

\sphinxAtStartPar
\sphinxstylestrong{Synopsis}
\begin{quote}

\sphinxAtStartPar
\sphinxcode{\sphinxupquote{Expr{[}{]} getExprs()}}
\end{quote}

\sphinxAtStartPar
\sphinxstylestrong{Return}
\begin{quote}

\sphinxAtStartPar
array of linear expressions.
\end{quote}
\end{quote}

\subsubsection{AffineConeBuilder.getPsdExpr()}
\label{\detokenize{javaapi/AffineConeBuilder:affineconebuilder-getpsdexpr}}\begin{quote}

\sphinxAtStartPar
Get idx\sphinxhyphen{}th PSD expression in an affine cone constraint.

\sphinxAtStartPar
\sphinxstylestrong{Synopsis}
\begin{quote}

\sphinxAtStartPar
\sphinxcode{\sphinxupquote{PsdExpr getPsdExpr(int idx)}}
\end{quote}

\sphinxAtStartPar
\sphinxstylestrong{Arguments}
\begin{quote}

\sphinxAtStartPar
\sphinxcode{\sphinxupquote{idx}}: index of PSD expression.
\end{quote}

\sphinxAtStartPar
\sphinxstylestrong{Return}
\begin{quote}

\sphinxAtStartPar
the idx\sphinxhyphen{}th PSD expression in an affine cone constraint.
\end{quote}
\end{quote}

\subsubsection{AffineConeBuilder.getPsdExprs()}
\label{\detokenize{javaapi/AffineConeBuilder:affineconebuilder-getpsdexprs}}\begin{quote}

\sphinxAtStartPar
Get all PSD expressions in an affine cone constraint.

\sphinxAtStartPar
\sphinxstylestrong{Synopsis}
\begin{quote}

\sphinxAtStartPar
\sphinxcode{\sphinxupquote{PsdExpr{[}{]} getPsdExprs()}}
\end{quote}

\sphinxAtStartPar
\sphinxstylestrong{Return}
\begin{quote}

\sphinxAtStartPar
array of PSD expressions.
\end{quote}
\end{quote}

\subsubsection{AffineConeBuilder.getSize()}
\label{\detokenize{javaapi/AffineConeBuilder:affineconebuilder-getsize}}\begin{quote}

\sphinxAtStartPar
Get number of variables in an affine cone constraint.

\sphinxAtStartPar
\sphinxstylestrong{Synopsis}
\begin{quote}

\sphinxAtStartPar
\sphinxcode{\sphinxupquote{int getSize()}}
\end{quote}

\sphinxAtStartPar
\sphinxstylestrong{Return}
\begin{quote}

\sphinxAtStartPar
number of vars.
\end{quote}
\end{quote}

\subsubsection{AffineConeBuilder.getType()}
\label{\detokenize{javaapi/AffineConeBuilder:affineconebuilder-gettype}}\begin{quote}

\sphinxAtStartPar
Get type of an affine cone constraint.

\sphinxAtStartPar
\sphinxstylestrong{Synopsis}
\begin{quote}

\sphinxAtStartPar
\sphinxcode{\sphinxupquote{int getType()}}
\end{quote}

\sphinxAtStartPar
\sphinxstylestrong{Return}
\begin{quote}

\sphinxAtStartPar
type of an affine cone constraint.
\end{quote}
\end{quote}

\subsubsection{AffineConeBuilder.hasPsdTerm()}
\label{\detokenize{javaapi/AffineConeBuilder:affineconebuilder-haspsdterm}}\begin{quote}

\sphinxAtStartPar
Check whether affine cone has PSD terms.

\sphinxAtStartPar
\sphinxstylestrong{Synopsis}
\begin{quote}

\sphinxAtStartPar
\sphinxcode{\sphinxupquote{Boolean hasPsdTerm()}}
\end{quote}

\sphinxAtStartPar
\sphinxstylestrong{Return}
\begin{quote}

\sphinxAtStartPar
flag to indicate whether affine cone has PSD terms.
\end{quote}
\end{quote}

\subsubsection{AffineConeBuilder.set()}
\label{\detokenize{javaapi/AffineConeBuilder:affineconebuilder-set}}\begin{quote}

\sphinxAtStartPar
Set linear expressions and type of an affine cone constraint.

\sphinxAtStartPar
\sphinxstylestrong{Synopsis}
\begin{quote}

\sphinxAtStartPar
\sphinxcode{\sphinxupquote{void set(Expr{[}{]} exprs, int type)}}
\end{quote}

\sphinxAtStartPar
\sphinxstylestrong{Arguments}
\begin{quote}

\sphinxAtStartPar
\sphinxcode{\sphinxupquote{exprs}}: array of linear expressions.

\sphinxAtStartPar
\sphinxcode{\sphinxupquote{type}}: type of an affine cone constraint.
\end{quote}
\end{quote}

\subsubsection{AffineConeBuilder.set()}
\label{\detokenize{javaapi/AffineConeBuilder:id1}}\begin{quote}

\sphinxAtStartPar
Set PSD expressions and type of an affine cone constraint.

\sphinxAtStartPar
\sphinxstylestrong{Synopsis}
\begin{quote}

\sphinxAtStartPar
\sphinxcode{\sphinxupquote{void set(PsdExpr{[}{]} exprs, int type)}}
\end{quote}

\sphinxAtStartPar
\sphinxstylestrong{Arguments}
\begin{quote}

\sphinxAtStartPar
\sphinxcode{\sphinxupquote{exprs}}: array of PSD expressions.

\sphinxAtStartPar
\sphinxcode{\sphinxupquote{type}}: type of an affine cone constraint.
\end{quote}
\end{quote}

\subsection{AffineConeBuilderArray Class}
\label{\detokenize{javaapiref:affineconebuilderarray-class}}\label{\detokenize{javaapiref:chapjavaapiref-affineconebuilderarray}}
\sphinxAtStartPar
To facilitate operations on a group of Java {\hyperref[\detokenize{javaapiref:chapjavaapiref-affineconebuilder}]{\sphinxcrossref{\DUrole{std,std-ref}{AffineConeBuilder Class}}}} objects,
the Java interface of COPT introduces the \sphinxtitleref{AffineConeBuilderArray} class.
The following methods are provided:

\sphinxstepscope

\subsubsection{AffineConeBuilderArray.AffineConeBuilderArray()}
\label{\detokenize{javaapi/AffineConeBuilderArray:affineconebuilderarray-affineconebuilderarray}}\label{\detokenize{javaapi/AffineConeBuilderArray::doc}}\begin{quote}

\sphinxAtStartPar
Constructor of AffineConeBuilderArray object.

\sphinxAtStartPar
\sphinxstylestrong{Synopsis}
\begin{quote}

\sphinxAtStartPar
\sphinxcode{\sphinxupquote{AffineConeBuilderArray()}}
\end{quote}
\end{quote}

\subsubsection{AffineConeBuilderArray.getBuilder()}
\label{\detokenize{javaapi/AffineConeBuilderArray:affineconebuilderarray-getbuilder}}\begin{quote}

\sphinxAtStartPar
Get idx\sphinxhyphen{}th affine cone constraint builder object.

\sphinxAtStartPar
\sphinxstylestrong{Synopsis}
\begin{quote}

\sphinxAtStartPar
\sphinxcode{\sphinxupquote{AffineConeBuilder getBuilder(int idx)}}
\end{quote}

\sphinxAtStartPar
\sphinxstylestrong{Arguments}
\begin{quote}

\sphinxAtStartPar
\sphinxcode{\sphinxupquote{idx}}: index of the affine cone constraint builder.
\end{quote}

\sphinxAtStartPar
\sphinxstylestrong{Return}
\begin{quote}

\sphinxAtStartPar
affine cone constraint builder object with index idx.
\end{quote}
\end{quote}

\subsubsection{AffineConeBuilderArray.pushBack()}
\label{\detokenize{javaapi/AffineConeBuilderArray:affineconebuilderarray-pushback}}\begin{quote}

\sphinxAtStartPar
Add an affine cone constraint builder object to affine cone constraint builder array.

\sphinxAtStartPar
\sphinxstylestrong{Synopsis}
\begin{quote}

\sphinxAtStartPar
\sphinxcode{\sphinxupquote{void pushBack(AffineConeBuilder builder)}}
\end{quote}

\sphinxAtStartPar
\sphinxstylestrong{Arguments}
\begin{quote}

\sphinxAtStartPar
\sphinxcode{\sphinxupquote{builder}}: an affine cone constraint builder object.
\end{quote}
\end{quote}

\subsubsection{AffineConeBuilderArray.size()}
\label{\detokenize{javaapi/AffineConeBuilderArray:affineconebuilderarray-size}}\begin{quote}

\sphinxAtStartPar
Get the number of affine cone constraint builder objects.

\sphinxAtStartPar
\sphinxstylestrong{Synopsis}
\begin{quote}

\sphinxAtStartPar
\sphinxcode{\sphinxupquote{int size()}}
\end{quote}

\sphinxAtStartPar
\sphinxstylestrong{Return}
\begin{quote}

\sphinxAtStartPar
number of affine cone constraint builder objects.
\end{quote}
\end{quote}

\subsection{QuadExpr}
\label{\detokenize{javaapiref:quadexpr}}\label{\detokenize{javaapiref:chapjavaapiref-quadexpr}}
\sphinxAtStartPar
COPT quadratic expression object. A quadratic expression consists of a linear expression,
a list of variable pairs and associated coefficients of quadratic terms. Quadratic expressions
are used to build quadratic constraints and objectives.

\sphinxstepscope

\subsubsection{QuadExpr.QuadExpr()}
\label{\detokenize{javaapi/QuadExpr:quadexpr-quadexpr}}\label{\detokenize{javaapi/QuadExpr::doc}}\begin{quote}

\sphinxAtStartPar
Constructor of a constant quadratic expression with constant 0.0

\sphinxAtStartPar
\sphinxstylestrong{Synopsis}
\begin{quote}

\sphinxAtStartPar
\sphinxcode{\sphinxupquote{QuadExpr()}}
\end{quote}
\end{quote}

\subsubsection{QuadExpr.QuadExpr()}
\label{\detokenize{javaapi/QuadExpr:id1}}\begin{quote}

\sphinxAtStartPar
Constructor of a constant quadratic expression.

\sphinxAtStartPar
\sphinxstylestrong{Synopsis}
\begin{quote}

\sphinxAtStartPar
\sphinxcode{\sphinxupquote{QuadExpr(double constant)}}
\end{quote}

\sphinxAtStartPar
\sphinxstylestrong{Arguments}
\begin{quote}

\sphinxAtStartPar
\sphinxcode{\sphinxupquote{constant}}: constant value in expression object.
\end{quote}
\end{quote}

\subsubsection{QuadExpr.QuadExpr()}
\label{\detokenize{javaapi/QuadExpr:id2}}\begin{quote}

\sphinxAtStartPar
Constructor of a quadratic expression with one term.

\sphinxAtStartPar
\sphinxstylestrong{Synopsis}
\begin{quote}

\sphinxAtStartPar
\sphinxcode{\sphinxupquote{QuadExpr(Var var)}}
\end{quote}

\sphinxAtStartPar
\sphinxstylestrong{Arguments}
\begin{quote}

\sphinxAtStartPar
\sphinxcode{\sphinxupquote{var}}: variable for the added term.
\end{quote}
\end{quote}

\subsubsection{QuadExpr.QuadExpr()}
\label{\detokenize{javaapi/QuadExpr:id3}}\begin{quote}

\sphinxAtStartPar
Constructor of a quadratic expression with one term.

\sphinxAtStartPar
\sphinxstylestrong{Synopsis}
\begin{quote}

\sphinxAtStartPar
\sphinxcode{\sphinxupquote{QuadExpr(Var var, double coeff)}}
\end{quote}

\sphinxAtStartPar
\sphinxstylestrong{Arguments}
\begin{quote}

\sphinxAtStartPar
\sphinxcode{\sphinxupquote{var}}: variable for the added term.

\sphinxAtStartPar
\sphinxcode{\sphinxupquote{coeff}}: coefficent for the added term.
\end{quote}
\end{quote}

\subsubsection{QuadExpr.QuadExpr()}
\label{\detokenize{javaapi/QuadExpr:id4}}\begin{quote}

\sphinxAtStartPar
Constructor of a quadratic expression with a linear expression.

\sphinxAtStartPar
\sphinxstylestrong{Synopsis}
\begin{quote}

\sphinxAtStartPar
\sphinxcode{\sphinxupquote{QuadExpr(Expr expr)}}
\end{quote}

\sphinxAtStartPar
\sphinxstylestrong{Arguments}
\begin{quote}

\sphinxAtStartPar
\sphinxcode{\sphinxupquote{expr}}: linear expression added to the quadratic expression.
\end{quote}
\end{quote}

\subsubsection{QuadExpr.QuadExpr()}
\label{\detokenize{javaapi/QuadExpr:id5}}\begin{quote}

\sphinxAtStartPar
Constructor of a quadratic expression with two linear expression.

\sphinxAtStartPar
\sphinxstylestrong{Synopsis}
\begin{quote}

\sphinxAtStartPar
\sphinxcode{\sphinxupquote{QuadExpr(Var var1, Var var2)}}
\end{quote}

\sphinxAtStartPar
\sphinxstylestrong{Arguments}
\begin{quote}

\sphinxAtStartPar
\sphinxcode{\sphinxupquote{var1}}: one variable.

\sphinxAtStartPar
\sphinxcode{\sphinxupquote{var2}}: another variable.
\end{quote}
\end{quote}

\subsubsection{QuadExpr.QuadExpr()}
\label{\detokenize{javaapi/QuadExpr:id6}}\begin{quote}

\sphinxAtStartPar
Constructor of a quadratic expression with two linear expression.

\sphinxAtStartPar
\sphinxstylestrong{Synopsis}
\begin{quote}

\sphinxAtStartPar
\sphinxcode{\sphinxupquote{QuadExpr(Expr expr, Var var)}}
\end{quote}

\sphinxAtStartPar
\sphinxstylestrong{Arguments}
\begin{quote}

\sphinxAtStartPar
\sphinxcode{\sphinxupquote{expr}}: one linear expression.

\sphinxAtStartPar
\sphinxcode{\sphinxupquote{var}}: another variable.
\end{quote}
\end{quote}

\subsubsection{QuadExpr.QuadExpr()}
\label{\detokenize{javaapi/QuadExpr:id7}}\begin{quote}

\sphinxAtStartPar
Constructor of a quadratic expression with two linear expression.

\sphinxAtStartPar
\sphinxstylestrong{Synopsis}
\begin{quote}

\sphinxAtStartPar
\sphinxcode{\sphinxupquote{QuadExpr(Expr left, Expr right)}}
\end{quote}

\sphinxAtStartPar
\sphinxstylestrong{Arguments}
\begin{quote}

\sphinxAtStartPar
\sphinxcode{\sphinxupquote{left}}: one linear expression.

\sphinxAtStartPar
\sphinxcode{\sphinxupquote{right}}: another linear expression.
\end{quote}
\end{quote}

\subsubsection{QuadExpr.add()}
\label{\detokenize{javaapi/QuadExpr:quadexpr-add}}\begin{quote}

\sphinxAtStartPar
Add itself by a quadratic expression.

\sphinxAtStartPar
\sphinxstylestrong{Synopsis}
\begin{quote}

\sphinxAtStartPar
\sphinxcode{\sphinxupquote{QuadExpr add(QuadExpr expr, double mult)}}
\end{quote}

\sphinxAtStartPar
\sphinxstylestrong{Arguments}
\begin{quote}

\sphinxAtStartPar
\sphinxcode{\sphinxupquote{expr}}: expression operand, including QuadExpr, Expr, Var and constant.

\sphinxAtStartPar
\sphinxcode{\sphinxupquote{mult}}: constant multiplier.
\end{quote}

\sphinxAtStartPar
\sphinxstylestrong{Return}
\begin{quote}

\sphinxAtStartPar
quadratic expression itself.
\end{quote}
\end{quote}

\subsubsection{QuadExpr.add()}
\label{\detokenize{javaapi/QuadExpr:id8}}\begin{quote}

\sphinxAtStartPar
Add itself by a quadratic expression.

\sphinxAtStartPar
\sphinxstylestrong{Synopsis}
\begin{quote}

\sphinxAtStartPar
\sphinxcode{\sphinxupquote{QuadExpr add(QuadExpr expr)}}
\end{quote}

\sphinxAtStartPar
\sphinxstylestrong{Arguments}
\begin{quote}

\sphinxAtStartPar
\sphinxcode{\sphinxupquote{expr}}: expression operand, including QuadExpr, Expr, Var and constant.
\end{quote}

\sphinxAtStartPar
\sphinxstylestrong{Return}
\begin{quote}

\sphinxAtStartPar
quadratic expression itself.
\end{quote}
\end{quote}

\subsubsection{QuadExpr.addConstant()}
\label{\detokenize{javaapi/QuadExpr:quadexpr-addconstant}}\begin{quote}

\sphinxAtStartPar
Add a constant to the quadratic expression.

\sphinxAtStartPar
\sphinxstylestrong{Synopsis}
\begin{quote}

\sphinxAtStartPar
\sphinxcode{\sphinxupquote{void addConstant(double constant)}}
\end{quote}

\sphinxAtStartPar
\sphinxstylestrong{Arguments}
\begin{quote}

\sphinxAtStartPar
\sphinxcode{\sphinxupquote{constant}}: value to be added.
\end{quote}
\end{quote}

\subsubsection{QuadExpr.addLinExpr()}
\label{\detokenize{javaapi/QuadExpr:quadexpr-addlinexpr}}\begin{quote}

\sphinxAtStartPar
Add a linear expression to self.

\sphinxAtStartPar
\sphinxstylestrong{Synopsis}
\begin{quote}

\sphinxAtStartPar
\sphinxcode{\sphinxupquote{void addLinExpr(Expr expr)}}
\end{quote}

\sphinxAtStartPar
\sphinxstylestrong{Arguments}
\begin{quote}

\sphinxAtStartPar
\sphinxcode{\sphinxupquote{expr}}: linear expression to be added.
\end{quote}
\end{quote}

\subsubsection{QuadExpr.addLinExpr()}
\label{\detokenize{javaapi/QuadExpr:id9}}\begin{quote}

\sphinxAtStartPar
Add a linear expression to self.

\sphinxAtStartPar
\sphinxstylestrong{Synopsis}
\begin{quote}

\sphinxAtStartPar
\sphinxcode{\sphinxupquote{void addLinExpr(Expr expr, double mult)}}
\end{quote}

\sphinxAtStartPar
\sphinxstylestrong{Arguments}
\begin{quote}

\sphinxAtStartPar
\sphinxcode{\sphinxupquote{expr}}: linear expression to be added.

\sphinxAtStartPar
\sphinxcode{\sphinxupquote{mult}}: multiplier constant.
\end{quote}
\end{quote}

\subsubsection{QuadExpr.addQuadExpr()}
\label{\detokenize{javaapi/QuadExpr:quadexpr-addquadexpr}}\begin{quote}

\sphinxAtStartPar
Add a quadratic expression to self.

\sphinxAtStartPar
\sphinxstylestrong{Synopsis}
\begin{quote}

\sphinxAtStartPar
\sphinxcode{\sphinxupquote{void addQuadExpr(QuadExpr expr)}}
\end{quote}

\sphinxAtStartPar
\sphinxstylestrong{Arguments}
\begin{quote}

\sphinxAtStartPar
\sphinxcode{\sphinxupquote{expr}}: quadratic expression to be added.
\end{quote}
\end{quote}

\subsubsection{QuadExpr.addQuadExpr()}
\label{\detokenize{javaapi/QuadExpr:id10}}\begin{quote}

\sphinxAtStartPar
Add a quadratic expression to self.

\sphinxAtStartPar
\sphinxstylestrong{Synopsis}
\begin{quote}

\sphinxAtStartPar
\sphinxcode{\sphinxupquote{void addQuadExpr(QuadExpr expr, double mult)}}
\end{quote}

\sphinxAtStartPar
\sphinxstylestrong{Arguments}
\begin{quote}

\sphinxAtStartPar
\sphinxcode{\sphinxupquote{expr}}: quadratic expression to be added.

\sphinxAtStartPar
\sphinxcode{\sphinxupquote{mult}}: multiplier constant.
\end{quote}
\end{quote}

\subsubsection{QuadExpr.addTerm()}
\label{\detokenize{javaapi/QuadExpr:quadexpr-addterm}}\begin{quote}

\sphinxAtStartPar
Add a term to quadratic expression object.

\sphinxAtStartPar
\sphinxstylestrong{Synopsis}
\begin{quote}

\sphinxAtStartPar
\sphinxcode{\sphinxupquote{void addTerm(Var var, double coeff)}}
\end{quote}

\sphinxAtStartPar
\sphinxstylestrong{Arguments}
\begin{quote}

\sphinxAtStartPar
\sphinxcode{\sphinxupquote{var}}: a variable of new term.

\sphinxAtStartPar
\sphinxcode{\sphinxupquote{coeff}}: coefficient of new term.
\end{quote}
\end{quote}

\subsubsection{QuadExpr.addTerm()}
\label{\detokenize{javaapi/QuadExpr:id11}}\begin{quote}

\sphinxAtStartPar
Add a quadratic term to expression object.

\sphinxAtStartPar
\sphinxstylestrong{Synopsis}
\begin{quote}

\sphinxAtStartPar
\sphinxcode{\sphinxupquote{void addTerm(}}
\begin{quote}

\sphinxAtStartPar
\sphinxcode{\sphinxupquote{Var var1,}}

\sphinxAtStartPar
\sphinxcode{\sphinxupquote{Var var2,}}

\sphinxAtStartPar
\sphinxcode{\sphinxupquote{double coeff)}}
\end{quote}
\end{quote}

\sphinxAtStartPar
\sphinxstylestrong{Arguments}
\begin{quote}

\sphinxAtStartPar
\sphinxcode{\sphinxupquote{var1}}: first variable of new quadratic term.

\sphinxAtStartPar
\sphinxcode{\sphinxupquote{var2}}: second variable of new quadratic term.

\sphinxAtStartPar
\sphinxcode{\sphinxupquote{coeff}}: coefficient of new quadratic term.
\end{quote}
\end{quote}

\subsubsection{QuadExpr.addTerms()}
\label{\detokenize{javaapi/QuadExpr:quadexpr-addterms}}\begin{quote}

\sphinxAtStartPar
Add linear terms to quadratic expression object.

\sphinxAtStartPar
\sphinxstylestrong{Synopsis}
\begin{quote}

\sphinxAtStartPar
\sphinxcode{\sphinxupquote{void addTerms(Var{[}{]} vars, double coeff)}}
\end{quote}

\sphinxAtStartPar
\sphinxstylestrong{Arguments}
\begin{quote}

\sphinxAtStartPar
\sphinxcode{\sphinxupquote{vars}}: variables of added linear terms.

\sphinxAtStartPar
\sphinxcode{\sphinxupquote{coeff}}: one coefficient for added linear terms.
\end{quote}
\end{quote}

\subsubsection{QuadExpr.addTerms()}
\label{\detokenize{javaapi/QuadExpr:id12}}\begin{quote}

\sphinxAtStartPar
Add linear terms to quadratic expression object.

\sphinxAtStartPar
\sphinxstylestrong{Synopsis}
\begin{quote}

\sphinxAtStartPar
\sphinxcode{\sphinxupquote{void addTerms(Var{[}{]} vars, double{[}{]} coeffs)}}
\end{quote}

\sphinxAtStartPar
\sphinxstylestrong{Arguments}
\begin{quote}

\sphinxAtStartPar
\sphinxcode{\sphinxupquote{vars}}: variables of added linear terms.

\sphinxAtStartPar
\sphinxcode{\sphinxupquote{coeffs}}: coefficients of added linear terms.
\end{quote}
\end{quote}

\subsubsection{QuadExpr.addTerms()}
\label{\detokenize{javaapi/QuadExpr:id13}}\begin{quote}

\sphinxAtStartPar
Add linear terms to quadratic expression object.

\sphinxAtStartPar
\sphinxstylestrong{Synopsis}
\begin{quote}

\sphinxAtStartPar
\sphinxcode{\sphinxupquote{void addTerms(VarArray vars, double coeff)}}
\end{quote}

\sphinxAtStartPar
\sphinxstylestrong{Arguments}
\begin{quote}

\sphinxAtStartPar
\sphinxcode{\sphinxupquote{vars}}: variables of added linear terms.

\sphinxAtStartPar
\sphinxcode{\sphinxupquote{coeff}}: one coefficient for added linear terms.
\end{quote}
\end{quote}

\subsubsection{QuadExpr.addTerms()}
\label{\detokenize{javaapi/QuadExpr:id14}}\begin{quote}

\sphinxAtStartPar
Add linear terms to quadratic expression object.

\sphinxAtStartPar
\sphinxstylestrong{Synopsis}
\begin{quote}

\sphinxAtStartPar
\sphinxcode{\sphinxupquote{void addTerms(VarArray vars, double{[}{]} coeffs)}}
\end{quote}

\sphinxAtStartPar
\sphinxstylestrong{Arguments}
\begin{quote}

\sphinxAtStartPar
\sphinxcode{\sphinxupquote{vars}}: variables of added terms.

\sphinxAtStartPar
\sphinxcode{\sphinxupquote{coeffs}}: coefficients of added terms.
\end{quote}
\end{quote}

\subsubsection{QuadExpr.addTerms()}
\label{\detokenize{javaapi/QuadExpr:id15}}\begin{quote}

\sphinxAtStartPar
Add quadratic terms to expression object.

\sphinxAtStartPar
\sphinxstylestrong{Synopsis}
\begin{quote}

\sphinxAtStartPar
\sphinxcode{\sphinxupquote{void addTerms(}}
\begin{quote}

\sphinxAtStartPar
\sphinxcode{\sphinxupquote{VarArray vars1,}}

\sphinxAtStartPar
\sphinxcode{\sphinxupquote{VarArray vars2,}}

\sphinxAtStartPar
\sphinxcode{\sphinxupquote{double{[}{]} coeffs)}}
\end{quote}
\end{quote}

\sphinxAtStartPar
\sphinxstylestrong{Arguments}
\begin{quote}

\sphinxAtStartPar
\sphinxcode{\sphinxupquote{vars1}}: first set of variables for added quadratic terms.

\sphinxAtStartPar
\sphinxcode{\sphinxupquote{vars2}}: second set of variables for added quadratic terms.

\sphinxAtStartPar
\sphinxcode{\sphinxupquote{coeffs}}: coefficient array for added quadratic terms.
\end{quote}
\end{quote}

\subsubsection{QuadExpr.addTerms()}
\label{\detokenize{javaapi/QuadExpr:id16}}\begin{quote}

\sphinxAtStartPar
Add quadratic terms to expression object.

\sphinxAtStartPar
\sphinxstylestrong{Synopsis}
\begin{quote}

\sphinxAtStartPar
\sphinxcode{\sphinxupquote{void addTerms(}}
\begin{quote}

\sphinxAtStartPar
\sphinxcode{\sphinxupquote{Var{[}{]} vars1,}}

\sphinxAtStartPar
\sphinxcode{\sphinxupquote{Var{[}{]} vars2,}}

\sphinxAtStartPar
\sphinxcode{\sphinxupquote{double{[}{]} coeffs)}}
\end{quote}
\end{quote}

\sphinxAtStartPar
\sphinxstylestrong{Arguments}
\begin{quote}

\sphinxAtStartPar
\sphinxcode{\sphinxupquote{vars1}}: first set of variables for added quadratic terms.

\sphinxAtStartPar
\sphinxcode{\sphinxupquote{vars2}}: second set of variables for added quadratic terms.

\sphinxAtStartPar
\sphinxcode{\sphinxupquote{coeffs}}: coefficient array for added quadratic terms.
\end{quote}
\end{quote}

\subsubsection{QuadExpr.clone()}
\label{\detokenize{javaapi/QuadExpr:quadexpr-clone}}\begin{quote}

\sphinxAtStartPar
Deep copy quadratic expression object.

\sphinxAtStartPar
\sphinxstylestrong{Synopsis}
\begin{quote}

\sphinxAtStartPar
\sphinxcode{\sphinxupquote{QuadExpr clone()}}
\end{quote}

\sphinxAtStartPar
\sphinxstylestrong{Return}
\begin{quote}

\sphinxAtStartPar
cloned quadratic expression object.
\end{quote}
\end{quote}

\subsubsection{QuadExpr.divide()}
\label{\detokenize{javaapi/QuadExpr:quadexpr-divide}}\begin{quote}

\sphinxAtStartPar
Divide itself by double constant.

\sphinxAtStartPar
\sphinxstylestrong{Synopsis}
\begin{quote}

\sphinxAtStartPar
\sphinxcode{\sphinxupquote{QuadExpr divide(double c)}}
\end{quote}

\sphinxAtStartPar
\sphinxstylestrong{Arguments}
\begin{quote}

\sphinxAtStartPar
\sphinxcode{\sphinxupquote{c}}: constant operand.
\end{quote}

\sphinxAtStartPar
\sphinxstylestrong{Return}
\begin{quote}

\sphinxAtStartPar
quadratic expression itself.
\end{quote}
\end{quote}

\subsubsection{QuadExpr.evaluate()}
\label{\detokenize{javaapi/QuadExpr:quadexpr-evaluate}}\begin{quote}

\sphinxAtStartPar
evaluate quadratic expression after solving.

\sphinxAtStartPar
\sphinxstylestrong{Synopsis}
\begin{quote}

\sphinxAtStartPar
\sphinxcode{\sphinxupquote{double evaluate()}}
\end{quote}

\sphinxAtStartPar
\sphinxstylestrong{Return}
\begin{quote}

\sphinxAtStartPar
value of quadratic expression.
\end{quote}
\end{quote}

\subsubsection{QuadExpr.getCoeff()}
\label{\detokenize{javaapi/QuadExpr:quadexpr-getcoeff}}\begin{quote}

\sphinxAtStartPar
Get coefficient from the i\sphinxhyphen{}th term in quadratic expression.

\sphinxAtStartPar
\sphinxstylestrong{Synopsis}
\begin{quote}

\sphinxAtStartPar
\sphinxcode{\sphinxupquote{double getCoeff(int i)}}
\end{quote}

\sphinxAtStartPar
\sphinxstylestrong{Arguments}
\begin{quote}

\sphinxAtStartPar
\sphinxcode{\sphinxupquote{i}}: index of the term.
\end{quote}

\sphinxAtStartPar
\sphinxstylestrong{Return}
\begin{quote}

\sphinxAtStartPar
coefficient of the i\sphinxhyphen{}th term in quadratic expression object.
\end{quote}
\end{quote}

\subsubsection{QuadExpr.getConstant()}
\label{\detokenize{javaapi/QuadExpr:quadexpr-getconstant}}\begin{quote}

\sphinxAtStartPar
Get constant in quadratic expression.

\sphinxAtStartPar
\sphinxstylestrong{Synopsis}
\begin{quote}

\sphinxAtStartPar
\sphinxcode{\sphinxupquote{double getConstant()}}
\end{quote}

\sphinxAtStartPar
\sphinxstylestrong{Return}
\begin{quote}

\sphinxAtStartPar
constant in quadratic expression.
\end{quote}
\end{quote}

\subsubsection{QuadExpr.getLinExpr()}
\label{\detokenize{javaapi/QuadExpr:quadexpr-getlinexpr}}\begin{quote}

\sphinxAtStartPar
Get linear expression in quadratic expression.

\sphinxAtStartPar
\sphinxstylestrong{Synopsis}
\begin{quote}

\sphinxAtStartPar
\sphinxcode{\sphinxupquote{Expr getLinExpr()}}
\end{quote}

\sphinxAtStartPar
\sphinxstylestrong{Return}
\begin{quote}

\sphinxAtStartPar
linear expression object.
\end{quote}
\end{quote}

\subsubsection{QuadExpr.getVar1()}
\label{\detokenize{javaapi/QuadExpr:quadexpr-getvar1}}\begin{quote}

\sphinxAtStartPar
Get first variable from the i\sphinxhyphen{}th term in quadratic expression.

\sphinxAtStartPar
\sphinxstylestrong{Synopsis}
\begin{quote}

\sphinxAtStartPar
\sphinxcode{\sphinxupquote{Var getVar1(int i)}}
\end{quote}

\sphinxAtStartPar
\sphinxstylestrong{Arguments}
\begin{quote}

\sphinxAtStartPar
\sphinxcode{\sphinxupquote{i}}: index of the term.
\end{quote}

\sphinxAtStartPar
\sphinxstylestrong{Return}
\begin{quote}

\sphinxAtStartPar
first variable of the i\sphinxhyphen{}th term in quadratic expression object.
\end{quote}
\end{quote}

\subsubsection{QuadExpr.getVar2()}
\label{\detokenize{javaapi/QuadExpr:quadexpr-getvar2}}\begin{quote}

\sphinxAtStartPar
Get second variable from the i\sphinxhyphen{}th term in quadratic expression.

\sphinxAtStartPar
\sphinxstylestrong{Synopsis}
\begin{quote}

\sphinxAtStartPar
\sphinxcode{\sphinxupquote{Var getVar2(int i)}}
\end{quote}

\sphinxAtStartPar
\sphinxstylestrong{Arguments}
\begin{quote}

\sphinxAtStartPar
\sphinxcode{\sphinxupquote{i}}: index of the term.
\end{quote}

\sphinxAtStartPar
\sphinxstylestrong{Return}
\begin{quote}

\sphinxAtStartPar
second variable of the i\sphinxhyphen{}th term in quadratic expression object.
\end{quote}
\end{quote}

\subsubsection{QuadExpr.multiply()}
\label{\detokenize{javaapi/QuadExpr:quadexpr-multiply}}\begin{quote}

\sphinxAtStartPar
Multiply itself by double constant.

\sphinxAtStartPar
\sphinxstylestrong{Synopsis}
\begin{quote}

\sphinxAtStartPar
\sphinxcode{\sphinxupquote{QuadExpr multiply(double c)}}
\end{quote}

\sphinxAtStartPar
\sphinxstylestrong{Arguments}
\begin{quote}

\sphinxAtStartPar
\sphinxcode{\sphinxupquote{c}}: constant operand.
\end{quote}

\sphinxAtStartPar
\sphinxstylestrong{Return}
\begin{quote}

\sphinxAtStartPar
quadratic expression itself.
\end{quote}
\end{quote}

\subsubsection{QuadExpr.remove()}
\label{\detokenize{javaapi/QuadExpr:quadexpr-remove}}\begin{quote}

\sphinxAtStartPar
Remove idx\sphinxhyphen{}th term from quadratic expression object.

\sphinxAtStartPar
\sphinxstylestrong{Synopsis}
\begin{quote}

\sphinxAtStartPar
\sphinxcode{\sphinxupquote{void remove(int idx)}}
\end{quote}

\sphinxAtStartPar
\sphinxstylestrong{Arguments}
\begin{quote}

\sphinxAtStartPar
\sphinxcode{\sphinxupquote{idx}}: index of the term to be removed.
\end{quote}
\end{quote}

\subsubsection{QuadExpr.remove()}
\label{\detokenize{javaapi/QuadExpr:id17}}\begin{quote}

\sphinxAtStartPar
Remove the term associated with variable from quadratic expression.

\sphinxAtStartPar
\sphinxstylestrong{Synopsis}
\begin{quote}

\sphinxAtStartPar
\sphinxcode{\sphinxupquote{void remove(Var var)}}
\end{quote}

\sphinxAtStartPar
\sphinxstylestrong{Arguments}
\begin{quote}

\sphinxAtStartPar
\sphinxcode{\sphinxupquote{var}}: a variable whose term should be removed.
\end{quote}
\end{quote}

\subsubsection{QuadExpr.setCoeff()}
\label{\detokenize{javaapi/QuadExpr:quadexpr-setcoeff}}\begin{quote}

\sphinxAtStartPar
Set coefficient of the i\sphinxhyphen{}th term in quadratic expression.

\sphinxAtStartPar
\sphinxstylestrong{Synopsis}
\begin{quote}

\sphinxAtStartPar
\sphinxcode{\sphinxupquote{void setCoeff(int i, double val)}}
\end{quote}

\sphinxAtStartPar
\sphinxstylestrong{Arguments}
\begin{quote}

\sphinxAtStartPar
\sphinxcode{\sphinxupquote{i}}: index of the quadratic term.

\sphinxAtStartPar
\sphinxcode{\sphinxupquote{val}}: coefficient of the term.
\end{quote}
\end{quote}

\subsubsection{QuadExpr.setConstant()}
\label{\detokenize{javaapi/QuadExpr:quadexpr-setconstant}}\begin{quote}

\sphinxAtStartPar
Set constant for the quadratic expression.

\sphinxAtStartPar
\sphinxstylestrong{Synopsis}
\begin{quote}

\sphinxAtStartPar
\sphinxcode{\sphinxupquote{void setConstant(double constant)}}
\end{quote}

\sphinxAtStartPar
\sphinxstylestrong{Arguments}
\begin{quote}

\sphinxAtStartPar
\sphinxcode{\sphinxupquote{constant}}: the value of the constant.
\end{quote}
\end{quote}

\subsubsection{QuadExpr.size()}
\label{\detokenize{javaapi/QuadExpr:quadexpr-size}}\begin{quote}

\sphinxAtStartPar
Get number of terms in quadratic expression.

\sphinxAtStartPar
\sphinxstylestrong{Synopsis}
\begin{quote}

\sphinxAtStartPar
\sphinxcode{\sphinxupquote{long size()}}
\end{quote}

\sphinxAtStartPar
\sphinxstylestrong{Return}
\begin{quote}

\sphinxAtStartPar
number of quadratic terms.
\end{quote}
\end{quote}

\subsection{QConstraint}
\label{\detokenize{javaapiref:qconstraint}}\label{\detokenize{javaapiref:chapjavaapiref-qconstr}}
\sphinxAtStartPar
COPT quadratic constraint object. Quadratic constraints are always associated with a particular model.
User creates a quadratic constraint object by adding a quadratic constraint to a model,
rather than by using constructor of QConstraint class.

\sphinxstepscope

\subsubsection{QConstraint.get()}
\label{\detokenize{javaapi/QConstraint:qconstraint-get}}\label{\detokenize{javaapi/QConstraint::doc}}\begin{quote}

\sphinxAtStartPar
Get information value of the quadratic constraint.

\sphinxAtStartPar
\sphinxstylestrong{Synopsis}
\begin{quote}

\sphinxAtStartPar
\sphinxcode{\sphinxupquote{double get(String info)}}
\end{quote}

\sphinxAtStartPar
\sphinxstylestrong{Arguments}
\begin{quote}

\sphinxAtStartPar
\sphinxcode{\sphinxupquote{info}}: name of the information being queried.
\end{quote}

\sphinxAtStartPar
\sphinxstylestrong{Return}
\begin{quote}

\sphinxAtStartPar
information value.
\end{quote}
\end{quote}

\subsubsection{QConstraint.getIdx()}
\label{\detokenize{javaapi/QConstraint:qconstraint-getidx}}\begin{quote}

\sphinxAtStartPar
Get index of the quadratic constraint.

\sphinxAtStartPar
\sphinxstylestrong{Synopsis}
\begin{quote}

\sphinxAtStartPar
\sphinxcode{\sphinxupquote{int getIdx()}}
\end{quote}

\sphinxAtStartPar
\sphinxstylestrong{Return}
\begin{quote}

\sphinxAtStartPar
the index of the quadratic constraint.
\end{quote}
\end{quote}

\subsubsection{QConstraint.getName()}
\label{\detokenize{javaapi/QConstraint:qconstraint-getname}}\begin{quote}

\sphinxAtStartPar
Get name of the constraint.

\sphinxAtStartPar
\sphinxstylestrong{Synopsis}
\begin{quote}

\sphinxAtStartPar
\sphinxcode{\sphinxupquote{String getName()}}
\end{quote}

\sphinxAtStartPar
\sphinxstylestrong{Return}
\begin{quote}

\sphinxAtStartPar
the name of the constraint.
\end{quote}
\end{quote}

\subsubsection{QConstraint.getRhs()}
\label{\detokenize{javaapi/QConstraint:qconstraint-getrhs}}\begin{quote}

\sphinxAtStartPar
Get rhs of quadratic constraint.

\sphinxAtStartPar
\sphinxstylestrong{Synopsis}
\begin{quote}

\sphinxAtStartPar
\sphinxcode{\sphinxupquote{double getRhs()}}
\end{quote}

\sphinxAtStartPar
\sphinxstylestrong{Return}
\begin{quote}

\sphinxAtStartPar
rhs of quadratic constraint.
\end{quote}
\end{quote}

\subsubsection{QConstraint.getSense()}
\label{\detokenize{javaapi/QConstraint:qconstraint-getsense}}\begin{quote}

\sphinxAtStartPar
Get rhs of quadratic constraint.

\sphinxAtStartPar
\sphinxstylestrong{Synopsis}
\begin{quote}

\sphinxAtStartPar
\sphinxcode{\sphinxupquote{char getSense()}}
\end{quote}

\sphinxAtStartPar
\sphinxstylestrong{Return}
\begin{quote}

\sphinxAtStartPar
rhs of quadratic constraint.
\end{quote}
\end{quote}

\subsubsection{QConstraint.remove()}
\label{\detokenize{javaapi/QConstraint:qconstraint-remove}}\begin{quote}

\sphinxAtStartPar
Remove this constraint from model.

\sphinxAtStartPar
\sphinxstylestrong{Synopsis}
\begin{quote}

\sphinxAtStartPar
\sphinxcode{\sphinxupquote{void remove()}}
\end{quote}
\end{quote}

\subsubsection{QConstraint.set()}
\label{\detokenize{javaapi/QConstraint:qconstraint-set}}\begin{quote}

\sphinxAtStartPar
Set information value of the quadratic constraint.

\sphinxAtStartPar
\sphinxstylestrong{Synopsis}
\begin{quote}

\sphinxAtStartPar
\sphinxcode{\sphinxupquote{void set(String info, double val)}}
\end{quote}

\sphinxAtStartPar
\sphinxstylestrong{Arguments}
\begin{quote}

\sphinxAtStartPar
\sphinxcode{\sphinxupquote{info}}: name of the information.

\sphinxAtStartPar
\sphinxcode{\sphinxupquote{val}}: new information value.
\end{quote}
\end{quote}

\subsubsection{QConstraint.setName()}
\label{\detokenize{javaapi/QConstraint:qconstraint-setname}}\begin{quote}

\sphinxAtStartPar
Set name of quadratic constraint.

\sphinxAtStartPar
\sphinxstylestrong{Synopsis}
\begin{quote}

\sphinxAtStartPar
\sphinxcode{\sphinxupquote{void setName(String name)}}
\end{quote}

\sphinxAtStartPar
\sphinxstylestrong{Arguments}
\begin{quote}

\sphinxAtStartPar
\sphinxcode{\sphinxupquote{name}}: the name to set.
\end{quote}
\end{quote}

\subsubsection{QConstraint.setRhs()}
\label{\detokenize{javaapi/QConstraint:qconstraint-setrhs}}\begin{quote}

\sphinxAtStartPar
Set rhs of quadratic constraint.

\sphinxAtStartPar
\sphinxstylestrong{Synopsis}
\begin{quote}

\sphinxAtStartPar
\sphinxcode{\sphinxupquote{void setRhs(double rhs)}}
\end{quote}

\sphinxAtStartPar
\sphinxstylestrong{Arguments}
\begin{quote}

\sphinxAtStartPar
\sphinxcode{\sphinxupquote{rhs}}: rhs of quadratic constraint.
\end{quote}
\end{quote}

\subsubsection{QConstraint.setSense()}
\label{\detokenize{javaapi/QConstraint:qconstraint-setsense}}\begin{quote}

\sphinxAtStartPar
Set sense of quadratic constraint.

\sphinxAtStartPar
\sphinxstylestrong{Synopsis}
\begin{quote}

\sphinxAtStartPar
\sphinxcode{\sphinxupquote{void setSense(char sense)}}
\end{quote}

\sphinxAtStartPar
\sphinxstylestrong{Arguments}
\begin{quote}

\sphinxAtStartPar
\sphinxcode{\sphinxupquote{sense}}: sense of quadratic constraint.
\end{quote}
\end{quote}

\subsection{QConstrArray}
\label{\detokenize{javaapiref:qconstrarray}}\label{\detokenize{javaapiref:chapjavaapiref-qconstrarray}}
\sphinxAtStartPar
COPT quadratic constraint array object. To store and access a set of Java {\hyperref[\detokenize{javaapiref:chapjavaapiref-qconstr}]{\sphinxcrossref{\DUrole{std,std-ref}{QConstraint}}}}
objects, Cardinal Optimizer provides Java QConstrArray class, which defines the following methods.

\sphinxstepscope

\subsubsection{QConstrArray.QConstrArray()}
\label{\detokenize{javaapi/QConstrArray:qconstrarray-qconstrarray}}\label{\detokenize{javaapi/QConstrArray::doc}}\begin{quote}

\sphinxAtStartPar
QConstructor of constrarray object.

\sphinxAtStartPar
\sphinxstylestrong{Synopsis}
\begin{quote}

\sphinxAtStartPar
\sphinxcode{\sphinxupquote{QConstrArray()}}
\end{quote}
\end{quote}

\subsubsection{QConstrArray.getQConstr()}
\label{\detokenize{javaapi/QConstrArray:qconstrarray-getqconstr}}\begin{quote}

\sphinxAtStartPar
Get idx\sphinxhyphen{}th constraint object.

\sphinxAtStartPar
\sphinxstylestrong{Synopsis}
\begin{quote}

\sphinxAtStartPar
\sphinxcode{\sphinxupquote{QConstraint getQConstr(int idx)}}
\end{quote}

\sphinxAtStartPar
\sphinxstylestrong{Arguments}
\begin{quote}

\sphinxAtStartPar
\sphinxcode{\sphinxupquote{idx}}: index of the constraint.
\end{quote}

\sphinxAtStartPar
\sphinxstylestrong{Return}
\begin{quote}

\sphinxAtStartPar
constraint object with index idx.
\end{quote}
\end{quote}

\subsubsection{QConstrArray.pushBack()}
\label{\detokenize{javaapi/QConstrArray:qconstrarray-pushback}}\begin{quote}

\sphinxAtStartPar
Add a constraint object to constraint array.

\sphinxAtStartPar
\sphinxstylestrong{Synopsis}
\begin{quote}

\sphinxAtStartPar
\sphinxcode{\sphinxupquote{void pushBack(QConstraint constr)}}
\end{quote}

\sphinxAtStartPar
\sphinxstylestrong{Arguments}
\begin{quote}

\sphinxAtStartPar
\sphinxcode{\sphinxupquote{constr}}: a constraint object.
\end{quote}
\end{quote}

\subsubsection{QConstrArray.size()}
\label{\detokenize{javaapi/QConstrArray:qconstrarray-size}}\begin{quote}

\sphinxAtStartPar
Get the number of constraint objects.

\sphinxAtStartPar
\sphinxstylestrong{Synopsis}
\begin{quote}

\sphinxAtStartPar
\sphinxcode{\sphinxupquote{int size()}}
\end{quote}

\sphinxAtStartPar
\sphinxstylestrong{Return}
\begin{quote}

\sphinxAtStartPar
number of constraint objects.
\end{quote}
\end{quote}

\subsection{QConstrBuilder}
\label{\detokenize{javaapiref:qconstrbuilder}}\label{\detokenize{javaapiref:chapjavaapiref-qconstrbuilder}}
\sphinxAtStartPar
COPT quadratic constraint builder object. To help building a quadratic constraint, given a quadratic
expression, constraint sense and right\sphinxhyphen{}hand side value, Cardinal Optimizer provides Java ConeBuilder
class, which defines the following methods.

\sphinxstepscope

\subsubsection{QConstrBuilder.QConstrBuilder()}
\label{\detokenize{javaapi/QConstrBuilder:qconstrbuilder-qconstrbuilder}}\label{\detokenize{javaapi/QConstrBuilder::doc}}\begin{quote}

\sphinxAtStartPar
QConstructor of constrbuilder object.

\sphinxAtStartPar
\sphinxstylestrong{Synopsis}
\begin{quote}

\sphinxAtStartPar
\sphinxcode{\sphinxupquote{QConstrBuilder()}}
\end{quote}
\end{quote}

\subsubsection{QConstrBuilder.getQuadExpr()}
\label{\detokenize{javaapi/QConstrBuilder:qconstrbuilder-getquadexpr}}\begin{quote}

\sphinxAtStartPar
Get expression associated with constraint.

\sphinxAtStartPar
\sphinxstylestrong{Synopsis}
\begin{quote}

\sphinxAtStartPar
\sphinxcode{\sphinxupquote{QuadExpr getQuadExpr()}}
\end{quote}

\sphinxAtStartPar
\sphinxstylestrong{Return}
\begin{quote}

\sphinxAtStartPar
quadratic expression object.
\end{quote}
\end{quote}

\subsubsection{QConstrBuilder.getSense()}
\label{\detokenize{javaapi/QConstrBuilder:qconstrbuilder-getsense}}\begin{quote}

\sphinxAtStartPar
Get sense associated with quadratic constraint.

\sphinxAtStartPar
\sphinxstylestrong{Synopsis}
\begin{quote}

\sphinxAtStartPar
\sphinxcode{\sphinxupquote{char getSense()}}
\end{quote}

\sphinxAtStartPar
\sphinxstylestrong{Return}
\begin{quote}

\sphinxAtStartPar
quadratic constraint sense.
\end{quote}
\end{quote}

\subsubsection{QConstrBuilder.set()}
\label{\detokenize{javaapi/QConstrBuilder:qconstrbuilder-set}}\begin{quote}

\sphinxAtStartPar
Set detail of a quadratic constraint to its builder object.

\sphinxAtStartPar
\sphinxstylestrong{Synopsis}
\begin{quote}

\sphinxAtStartPar
\sphinxcode{\sphinxupquote{void set(}}
\begin{quote}

\sphinxAtStartPar
\sphinxcode{\sphinxupquote{QuadExpr expr,}}

\sphinxAtStartPar
\sphinxcode{\sphinxupquote{char sense,}}

\sphinxAtStartPar
\sphinxcode{\sphinxupquote{double rhs)}}
\end{quote}
\end{quote}

\sphinxAtStartPar
\sphinxstylestrong{Arguments}
\begin{quote}

\sphinxAtStartPar
\sphinxcode{\sphinxupquote{expr}}: expression object at one side of the quadratic constraint.

\sphinxAtStartPar
\sphinxcode{\sphinxupquote{sense}}: quadratic constraint sense.

\sphinxAtStartPar
\sphinxcode{\sphinxupquote{rhs}}: constant of right side of quadratic constraint.
\end{quote}
\end{quote}

\subsection{QConstrBuilderArray}
\label{\detokenize{javaapiref:qconstrbuilderarray}}\label{\detokenize{javaapiref:chapjavaapiref-qconstrbuilderarray}}
\sphinxAtStartPar
COPT quadratic constraint builder array object. To store and access a set of Java {\hyperref[\detokenize{javaapiref:chapjavaapiref-qconstrbuilder}]{\sphinxcrossref{\DUrole{std,std-ref}{QConstrBuilder}}}}
objects, Cardinal Optimizer provides Java QConstrBuilderArray class, which defines the following methods.

\sphinxstepscope

\subsubsection{QConstrBuilderArray.QConstrBuilderArray()}
\label{\detokenize{javaapi/QConstrBuilderArray:qconstrbuilderarray-qconstrbuilderarray}}\label{\detokenize{javaapi/QConstrBuilderArray::doc}}\begin{quote}

\sphinxAtStartPar
QConstructor of constrbuilderarray object.

\sphinxAtStartPar
\sphinxstylestrong{Synopsis}
\begin{quote}

\sphinxAtStartPar
\sphinxcode{\sphinxupquote{QConstrBuilderArray()}}
\end{quote}
\end{quote}

\subsubsection{QConstrBuilderArray.getBuilder()}
\label{\detokenize{javaapi/QConstrBuilderArray:qconstrbuilderarray-getbuilder}}\begin{quote}

\sphinxAtStartPar
Get idx\sphinxhyphen{}th constraint builder object.

\sphinxAtStartPar
\sphinxstylestrong{Synopsis}
\begin{quote}

\sphinxAtStartPar
\sphinxcode{\sphinxupquote{QConstrBuilder getBuilder(int idx)}}
\end{quote}

\sphinxAtStartPar
\sphinxstylestrong{Arguments}
\begin{quote}

\sphinxAtStartPar
\sphinxcode{\sphinxupquote{idx}}: index of the constraint builder.
\end{quote}

\sphinxAtStartPar
\sphinxstylestrong{Return}
\begin{quote}

\sphinxAtStartPar
constraint builder object with index idx.
\end{quote}
\end{quote}

\subsubsection{QConstrBuilderArray.pushBack()}
\label{\detokenize{javaapi/QConstrBuilderArray:qconstrbuilderarray-pushback}}\begin{quote}

\sphinxAtStartPar
Add a constraint builder object to constraint builder array.

\sphinxAtStartPar
\sphinxstylestrong{Synopsis}
\begin{quote}

\sphinxAtStartPar
\sphinxcode{\sphinxupquote{void pushBack(QConstrBuilder builder)}}
\end{quote}

\sphinxAtStartPar
\sphinxstylestrong{Arguments}
\begin{quote}

\sphinxAtStartPar
\sphinxcode{\sphinxupquote{builder}}: a constraint builder object.
\end{quote}
\end{quote}

\subsubsection{QConstrBuilderArray.size()}
\label{\detokenize{javaapi/QConstrBuilderArray:qconstrbuilderarray-size}}\begin{quote}

\sphinxAtStartPar
Get the number of constraint builder objects.

\sphinxAtStartPar
\sphinxstylestrong{Synopsis}
\begin{quote}

\sphinxAtStartPar
\sphinxcode{\sphinxupquote{int size()}}
\end{quote}

\sphinxAtStartPar
\sphinxstylestrong{Return}
\begin{quote}

\sphinxAtStartPar
number of constraint builder objects.
\end{quote}
\end{quote}

\subsection{PsdVar}
\label{\detokenize{javaapiref:psdvar}}\label{\detokenize{javaapiref:chapjavaapiref-psdvar}}
\sphinxAtStartPar
COPT PSD variable object. PSD variables are always associated with a particular model.
User creates a PSD variable object by adding a PSD variable to model, rather than
by constructor of PsdVar class.

\sphinxstepscope

\subsubsection{PsdVar.get()}
\label{\detokenize{javaapi/PsdVar:psdvar-get}}\label{\detokenize{javaapi/PsdVar::doc}}\begin{quote}

\sphinxAtStartPar
Get information values of PSD variable.

\sphinxAtStartPar
\sphinxstylestrong{Synopsis}
\begin{quote}

\sphinxAtStartPar
\sphinxcode{\sphinxupquote{double{[}{]} get(String info)}}
\end{quote}

\sphinxAtStartPar
\sphinxstylestrong{Arguments}
\begin{quote}

\sphinxAtStartPar
\sphinxcode{\sphinxupquote{info}}: information name.
\end{quote}

\sphinxAtStartPar
\sphinxstylestrong{Return}
\begin{quote}

\sphinxAtStartPar
array of information values.
\end{quote}
\end{quote}

\subsubsection{PsdVar.getDim()}
\label{\detokenize{javaapi/PsdVar:psdvar-getdim}}\begin{quote}

\sphinxAtStartPar
Get dimension of PSD variable.

\sphinxAtStartPar
\sphinxstylestrong{Synopsis}
\begin{quote}

\sphinxAtStartPar
\sphinxcode{\sphinxupquote{int getDim()}}
\end{quote}

\sphinxAtStartPar
\sphinxstylestrong{Return}
\begin{quote}

\sphinxAtStartPar
dimension of PSD variable.
\end{quote}
\end{quote}

\subsubsection{PsdVar.getIdx()}
\label{\detokenize{javaapi/PsdVar:psdvar-getidx}}\begin{quote}

\sphinxAtStartPar
Get index of PSD variable.

\sphinxAtStartPar
\sphinxstylestrong{Synopsis}
\begin{quote}

\sphinxAtStartPar
\sphinxcode{\sphinxupquote{int getIdx()}}
\end{quote}

\sphinxAtStartPar
\sphinxstylestrong{Return}
\begin{quote}

\sphinxAtStartPar
index of PSD variable.
\end{quote}
\end{quote}

\subsubsection{PsdVar.getLen()}
\label{\detokenize{javaapi/PsdVar:psdvar-getlen}}\begin{quote}

\sphinxAtStartPar
Get length of PSD variable.

\sphinxAtStartPar
\sphinxstylestrong{Synopsis}
\begin{quote}

\sphinxAtStartPar
\sphinxcode{\sphinxupquote{int getLen()}}
\end{quote}

\sphinxAtStartPar
\sphinxstylestrong{Return}
\begin{quote}

\sphinxAtStartPar
length of PSD variable.
\end{quote}
\end{quote}

\subsubsection{PsdVar.getName()}
\label{\detokenize{javaapi/PsdVar:psdvar-getname}}\begin{quote}

\sphinxAtStartPar
Get name of PSD variable.

\sphinxAtStartPar
\sphinxstylestrong{Synopsis}
\begin{quote}

\sphinxAtStartPar
\sphinxcode{\sphinxupquote{String getName()}}
\end{quote}

\sphinxAtStartPar
\sphinxstylestrong{Return}
\begin{quote}

\sphinxAtStartPar
name of PSD variable.
\end{quote}
\end{quote}

\subsubsection{PsdVar.remove()}
\label{\detokenize{javaapi/PsdVar:psdvar-remove}}\begin{quote}

\sphinxAtStartPar
Remove PSD variable from model.

\sphinxAtStartPar
\sphinxstylestrong{Synopsis}
\begin{quote}

\sphinxAtStartPar
\sphinxcode{\sphinxupquote{void remove()}}
\end{quote}
\end{quote}

\subsection{PsdVarArray}
\label{\detokenize{javaapiref:psdvararray}}\label{\detokenize{javaapiref:chapjavaapiref-psdvararray}}
\sphinxAtStartPar
COPT PSD variable array object. To store and access a set of {\hyperref[\detokenize{javaapiref:chapjavaapiref-psdvar}]{\sphinxcrossref{\DUrole{std,std-ref}{PsdVar}}}} objects,
Cardinal Optimizer provides PsdVarArray class, which defines the following methods.

\sphinxstepscope

\subsubsection{PsdVarArray.PsdVarArray()}
\label{\detokenize{javaapi/PsdVarArray:psdvararray-psdvararray}}\label{\detokenize{javaapi/PsdVarArray::doc}}\begin{quote}

\sphinxAtStartPar
Constructor of PsdVarArray.

\sphinxAtStartPar
\sphinxstylestrong{Synopsis}
\begin{quote}

\sphinxAtStartPar
\sphinxcode{\sphinxupquote{PsdVarArray()}}
\end{quote}
\end{quote}

\subsubsection{PsdVarArray.getPsdVar()}
\label{\detokenize{javaapi/PsdVarArray:psdvararray-getpsdvar}}\begin{quote}

\sphinxAtStartPar
Get idx\sphinxhyphen{}th PSD variable object.

\sphinxAtStartPar
\sphinxstylestrong{Synopsis}
\begin{quote}

\sphinxAtStartPar
\sphinxcode{\sphinxupquote{PsdVar getPsdVar(int idx)}}
\end{quote}

\sphinxAtStartPar
\sphinxstylestrong{Arguments}
\begin{quote}

\sphinxAtStartPar
\sphinxcode{\sphinxupquote{idx}}: index of the PSD variable.
\end{quote}

\sphinxAtStartPar
\sphinxstylestrong{Return}
\begin{quote}

\sphinxAtStartPar
PSD variable object with index idx.
\end{quote}
\end{quote}

\subsubsection{PsdVarArray.pushBack()}
\label{\detokenize{javaapi/PsdVarArray:psdvararray-pushback}}\begin{quote}

\sphinxAtStartPar
Add a PSD variable object to PSD variable array.

\sphinxAtStartPar
\sphinxstylestrong{Synopsis}
\begin{quote}

\sphinxAtStartPar
\sphinxcode{\sphinxupquote{void pushBack(PsdVar var)}}
\end{quote}

\sphinxAtStartPar
\sphinxstylestrong{Arguments}
\begin{quote}

\sphinxAtStartPar
\sphinxcode{\sphinxupquote{var}}: a PSD variable object.
\end{quote}
\end{quote}

\subsubsection{PsdVarArray.reserve()}
\label{\detokenize{javaapi/PsdVarArray:psdvararray-reserve}}\begin{quote}

\sphinxAtStartPar
Reserve capacity to contain at least n items.

\sphinxAtStartPar
\sphinxstylestrong{Synopsis}
\begin{quote}

\sphinxAtStartPar
\sphinxcode{\sphinxupquote{void reserve(int n)}}
\end{quote}

\sphinxAtStartPar
\sphinxstylestrong{Arguments}
\begin{quote}

\sphinxAtStartPar
\sphinxcode{\sphinxupquote{n}}: minimum capacity for PSD variable object.
\end{quote}
\end{quote}

\subsubsection{PsdVarArray.size()}
\label{\detokenize{javaapi/PsdVarArray:psdvararray-size}}\begin{quote}

\sphinxAtStartPar
Get the number of PSD variable objects.

\sphinxAtStartPar
\sphinxstylestrong{Synopsis}
\begin{quote}

\sphinxAtStartPar
\sphinxcode{\sphinxupquote{int size()}}
\end{quote}

\sphinxAtStartPar
\sphinxstylestrong{Return}
\begin{quote}

\sphinxAtStartPar
number of PSD variable objects.
\end{quote}
\end{quote}

\subsection{PsdExpr}
\label{\detokenize{javaapiref:psdexpr}}\label{\detokenize{javaapiref:chapjavaapiref-psdexpr}}
\sphinxAtStartPar
COPT PSD expression object. A PSD expression consists of a linear expression,
a list of PSD variables and associated coefficient matrices of PSD terms. PSD expressions
are used to build PSD constraints and objectives.

\sphinxstepscope

\subsubsection{PsdExpr.PsdExpr()}
\label{\detokenize{javaapi/PsdExpr:psdexpr-psdexpr}}\label{\detokenize{javaapi/PsdExpr::doc}}\begin{quote}

\sphinxAtStartPar
Constructor of a PSD expression with default constant value 0.

\sphinxAtStartPar
\sphinxstylestrong{Synopsis}
\begin{quote}

\sphinxAtStartPar
\sphinxcode{\sphinxupquote{PsdExpr(double constant)}}
\end{quote}

\sphinxAtStartPar
\sphinxstylestrong{Arguments}
\begin{quote}

\sphinxAtStartPar
\sphinxcode{\sphinxupquote{constant}}: constant value in PSD expression object.
\end{quote}
\end{quote}

\subsubsection{PsdExpr.PsdExpr()}
\label{\detokenize{javaapi/PsdExpr:id1}}\begin{quote}

\sphinxAtStartPar
Constructor of a PSD expression with one term.

\sphinxAtStartPar
\sphinxstylestrong{Synopsis}
\begin{quote}

\sphinxAtStartPar
\sphinxcode{\sphinxupquote{PsdExpr(Var var)}}
\end{quote}

\sphinxAtStartPar
\sphinxstylestrong{Arguments}
\begin{quote}

\sphinxAtStartPar
\sphinxcode{\sphinxupquote{var}}: variable for the added term.
\end{quote}
\end{quote}

\subsubsection{PsdExpr.PsdExpr()}
\label{\detokenize{javaapi/PsdExpr:id2}}\begin{quote}

\sphinxAtStartPar
Constructor of a PSD expression with one term.

\sphinxAtStartPar
\sphinxstylestrong{Synopsis}
\begin{quote}

\sphinxAtStartPar
\sphinxcode{\sphinxupquote{PsdExpr(Var var, double coeff)}}
\end{quote}

\sphinxAtStartPar
\sphinxstylestrong{Arguments}
\begin{quote}

\sphinxAtStartPar
\sphinxcode{\sphinxupquote{var}}: variable for the added term.

\sphinxAtStartPar
\sphinxcode{\sphinxupquote{coeff}}: coefficent for the added term.
\end{quote}
\end{quote}

\subsubsection{PsdExpr.PsdExpr()}
\label{\detokenize{javaapi/PsdExpr:id3}}\begin{quote}

\sphinxAtStartPar
Constructor of a PSD expression with a linear expression.

\sphinxAtStartPar
\sphinxstylestrong{Synopsis}
\begin{quote}

\sphinxAtStartPar
\sphinxcode{\sphinxupquote{PsdExpr(Expr expr)}}
\end{quote}

\sphinxAtStartPar
\sphinxstylestrong{Arguments}
\begin{quote}

\sphinxAtStartPar
\sphinxcode{\sphinxupquote{expr}}: input linear expression.
\end{quote}
\end{quote}

\subsubsection{PsdExpr.PsdExpr()}
\label{\detokenize{javaapi/PsdExpr:id4}}\begin{quote}

\sphinxAtStartPar
Constructor of a PSD expression with one term.

\sphinxAtStartPar
\sphinxstylestrong{Synopsis}
\begin{quote}

\sphinxAtStartPar
\sphinxcode{\sphinxupquote{PsdExpr(PsdVar var, SymMatrix mat)}}
\end{quote}

\sphinxAtStartPar
\sphinxstylestrong{Arguments}
\begin{quote}

\sphinxAtStartPar
\sphinxcode{\sphinxupquote{var}}: PSD variable for the added term.

\sphinxAtStartPar
\sphinxcode{\sphinxupquote{mat}}: coefficient matrix for the added term.
\end{quote}
\end{quote}

\subsubsection{PsdExpr.PsdExpr()}
\label{\detokenize{javaapi/PsdExpr:id5}}\begin{quote}

\sphinxAtStartPar
Constructor of a PSD expression with one term.

\sphinxAtStartPar
\sphinxstylestrong{Synopsis}
\begin{quote}

\sphinxAtStartPar
\sphinxcode{\sphinxupquote{PsdExpr(PsdVar var, SymMatExpr expr)}}
\end{quote}

\sphinxAtStartPar
\sphinxstylestrong{Arguments}
\begin{quote}

\sphinxAtStartPar
\sphinxcode{\sphinxupquote{var}}: PSD variable for the added term.

\sphinxAtStartPar
\sphinxcode{\sphinxupquote{expr}}: coefficient expression of symmetric matrices of new PSD term.
\end{quote}
\end{quote}

\subsubsection{PsdExpr.addConstant()}
\label{\detokenize{javaapi/PsdExpr:psdexpr-addconstant}}\begin{quote}

\sphinxAtStartPar
Add constant to the PSD expression.

\sphinxAtStartPar
\sphinxstylestrong{Synopsis}
\begin{quote}

\sphinxAtStartPar
\sphinxcode{\sphinxupquote{void addConstant(double constant)}}
\end{quote}

\sphinxAtStartPar
\sphinxstylestrong{Arguments}
\begin{quote}

\sphinxAtStartPar
\sphinxcode{\sphinxupquote{constant}}: value to be added.
\end{quote}
\end{quote}

\subsubsection{PsdExpr.addLinExpr()}
\label{\detokenize{javaapi/PsdExpr:psdexpr-addlinexpr}}\begin{quote}

\sphinxAtStartPar
Add a linear expression to PSD expression object.

\sphinxAtStartPar
\sphinxstylestrong{Synopsis}
\begin{quote}

\sphinxAtStartPar
\sphinxcode{\sphinxupquote{void addLinExpr(Expr expr)}}
\end{quote}

\sphinxAtStartPar
\sphinxstylestrong{Arguments}
\begin{quote}

\sphinxAtStartPar
\sphinxcode{\sphinxupquote{expr}}: linear expression to be added.
\end{quote}
\end{quote}

\subsubsection{PsdExpr.addLinExpr()}
\label{\detokenize{javaapi/PsdExpr:id6}}\begin{quote}

\sphinxAtStartPar
Add a linear expression to PSD expression object.

\sphinxAtStartPar
\sphinxstylestrong{Synopsis}
\begin{quote}

\sphinxAtStartPar
\sphinxcode{\sphinxupquote{void addLinExpr(Expr expr, double mult)}}
\end{quote}

\sphinxAtStartPar
\sphinxstylestrong{Arguments}
\begin{quote}

\sphinxAtStartPar
\sphinxcode{\sphinxupquote{expr}}: linear expression to be added.

\sphinxAtStartPar
\sphinxcode{\sphinxupquote{mult}}: multiplier constant.
\end{quote}
\end{quote}

\subsubsection{PsdExpr.addPsdExpr()}
\label{\detokenize{javaapi/PsdExpr:psdexpr-addpsdexpr}}\begin{quote}

\sphinxAtStartPar
Add a PSD expression to self.

\sphinxAtStartPar
\sphinxstylestrong{Synopsis}
\begin{quote}

\sphinxAtStartPar
\sphinxcode{\sphinxupquote{void addPsdExpr(PsdExpr expr)}}
\end{quote}

\sphinxAtStartPar
\sphinxstylestrong{Arguments}
\begin{quote}

\sphinxAtStartPar
\sphinxcode{\sphinxupquote{expr}}: PSD expression to be added.
\end{quote}
\end{quote}

\subsubsection{PsdExpr.addPsdExpr()}
\label{\detokenize{javaapi/PsdExpr:id7}}\begin{quote}

\sphinxAtStartPar
Add a PSD expression to self.

\sphinxAtStartPar
\sphinxstylestrong{Synopsis}
\begin{quote}

\sphinxAtStartPar
\sphinxcode{\sphinxupquote{void addPsdExpr(PsdExpr expr, double mult)}}
\end{quote}

\sphinxAtStartPar
\sphinxstylestrong{Arguments}
\begin{quote}

\sphinxAtStartPar
\sphinxcode{\sphinxupquote{expr}}: PSD expression to be added.

\sphinxAtStartPar
\sphinxcode{\sphinxupquote{mult}}: multiplier constant.
\end{quote}
\end{quote}

\subsubsection{PsdExpr.addTerm()}
\label{\detokenize{javaapi/PsdExpr:psdexpr-addterm}}\begin{quote}

\sphinxAtStartPar
Add a linear term to PSD expression object.

\sphinxAtStartPar
\sphinxstylestrong{Synopsis}
\begin{quote}

\sphinxAtStartPar
\sphinxcode{\sphinxupquote{void addTerm(Var var, double coeff)}}
\end{quote}

\sphinxAtStartPar
\sphinxstylestrong{Arguments}
\begin{quote}

\sphinxAtStartPar
\sphinxcode{\sphinxupquote{var}}: variable of new linear term.

\sphinxAtStartPar
\sphinxcode{\sphinxupquote{coeff}}: coefficient of new linear term.
\end{quote}
\end{quote}

\subsubsection{PsdExpr.addTerm()}
\label{\detokenize{javaapi/PsdExpr:id8}}\begin{quote}

\sphinxAtStartPar
Add a PSD term to PSD expression object.

\sphinxAtStartPar
\sphinxstylestrong{Synopsis}
\begin{quote}

\sphinxAtStartPar
\sphinxcode{\sphinxupquote{void addTerm(PsdVar var, SymMatrix mat)}}
\end{quote}

\sphinxAtStartPar
\sphinxstylestrong{Arguments}
\begin{quote}

\sphinxAtStartPar
\sphinxcode{\sphinxupquote{var}}: PSD variable of new PSD term.

\sphinxAtStartPar
\sphinxcode{\sphinxupquote{mat}}: coefficient matrix of new PSD term.
\end{quote}
\end{quote}

\subsubsection{PsdExpr.addTerm()}
\label{\detokenize{javaapi/PsdExpr:id9}}\begin{quote}

\sphinxAtStartPar
Add a PSD term to PSD expression object.

\sphinxAtStartPar
\sphinxstylestrong{Synopsis}
\begin{quote}

\sphinxAtStartPar
\sphinxcode{\sphinxupquote{void addTerm(PsdVar var, SymMatExpr expr)}}
\end{quote}

\sphinxAtStartPar
\sphinxstylestrong{Arguments}
\begin{quote}

\sphinxAtStartPar
\sphinxcode{\sphinxupquote{var}}: PSD variable of new PSD term.

\sphinxAtStartPar
\sphinxcode{\sphinxupquote{expr}}: coefficient expression of symmetric matrices of new PSD term.
\end{quote}
\end{quote}

\subsubsection{PsdExpr.addTerms()}
\label{\detokenize{javaapi/PsdExpr:psdexpr-addterms}}\begin{quote}

\sphinxAtStartPar
Add linear terms to PSD expression object.

\sphinxAtStartPar
\sphinxstylestrong{Synopsis}
\begin{quote}

\sphinxAtStartPar
\sphinxcode{\sphinxupquote{void addTerms(Var{[}{]} vars, double coeff)}}
\end{quote}

\sphinxAtStartPar
\sphinxstylestrong{Arguments}
\begin{quote}

\sphinxAtStartPar
\sphinxcode{\sphinxupquote{vars}}: variables of added linear terms.

\sphinxAtStartPar
\sphinxcode{\sphinxupquote{coeff}}: one coefficient for added linear terms.
\end{quote}
\end{quote}

\subsubsection{PsdExpr.addTerms()}
\label{\detokenize{javaapi/PsdExpr:id10}}\begin{quote}

\sphinxAtStartPar
Add linear terms to PSD expression object.

\sphinxAtStartPar
\sphinxstylestrong{Synopsis}
\begin{quote}

\sphinxAtStartPar
\sphinxcode{\sphinxupquote{void addTerms(Var{[}{]} vars, double{[}{]} coeffs)}}
\end{quote}

\sphinxAtStartPar
\sphinxstylestrong{Arguments}
\begin{quote}

\sphinxAtStartPar
\sphinxcode{\sphinxupquote{vars}}: variables for added linear terms.

\sphinxAtStartPar
\sphinxcode{\sphinxupquote{coeffs}}: coefficient array for added linear terms.
\end{quote}
\end{quote}

\subsubsection{PsdExpr.addTerms()}
\label{\detokenize{javaapi/PsdExpr:id11}}\begin{quote}

\sphinxAtStartPar
Add linear terms to PSD expression object.

\sphinxAtStartPar
\sphinxstylestrong{Synopsis}
\begin{quote}

\sphinxAtStartPar
\sphinxcode{\sphinxupquote{void addTerms(VarArray vars, double coeff)}}
\end{quote}

\sphinxAtStartPar
\sphinxstylestrong{Arguments}
\begin{quote}

\sphinxAtStartPar
\sphinxcode{\sphinxupquote{vars}}: variables of added linear terms.

\sphinxAtStartPar
\sphinxcode{\sphinxupquote{coeff}}: one coefficient for added linear terms.
\end{quote}
\end{quote}

\subsubsection{PsdExpr.addTerms()}
\label{\detokenize{javaapi/PsdExpr:id12}}\begin{quote}

\sphinxAtStartPar
Add linear terms to PSD expression object.

\sphinxAtStartPar
\sphinxstylestrong{Synopsis}
\begin{quote}

\sphinxAtStartPar
\sphinxcode{\sphinxupquote{void addTerms(VarArray vars, double{[}{]} coeffs)}}
\end{quote}

\sphinxAtStartPar
\sphinxstylestrong{Arguments}
\begin{quote}

\sphinxAtStartPar
\sphinxcode{\sphinxupquote{vars}}: variables of added terms.

\sphinxAtStartPar
\sphinxcode{\sphinxupquote{coeffs}}: coefficients of added terms.
\end{quote}
\end{quote}

\subsubsection{PsdExpr.addTerms()}
\label{\detokenize{javaapi/PsdExpr:id13}}\begin{quote}

\sphinxAtStartPar
Add PSD terms to PSD expression object.

\sphinxAtStartPar
\sphinxstylestrong{Synopsis}
\begin{quote}

\sphinxAtStartPar
\sphinxcode{\sphinxupquote{void addTerms(PsdVarArray vars, SymMatrixArray mats)}}
\end{quote}

\sphinxAtStartPar
\sphinxstylestrong{Arguments}
\begin{quote}

\sphinxAtStartPar
\sphinxcode{\sphinxupquote{vars}}: PSD variables for added PSD terms.

\sphinxAtStartPar
\sphinxcode{\sphinxupquote{mats}}: coefficient matrixes for added PSD terms.
\end{quote}
\end{quote}

\subsubsection{PsdExpr.addTerms()}
\label{\detokenize{javaapi/PsdExpr:id14}}\begin{quote}

\sphinxAtStartPar
Add PSD terms to PSD expression object.

\sphinxAtStartPar
\sphinxstylestrong{Synopsis}
\begin{quote}

\sphinxAtStartPar
\sphinxcode{\sphinxupquote{void addTerms(PsdVar{[}{]} vars, SymMatrix{[}{]} mats)}}
\end{quote}

\sphinxAtStartPar
\sphinxstylestrong{Arguments}
\begin{quote}

\sphinxAtStartPar
\sphinxcode{\sphinxupquote{vars}}: PSD variables for added PSD terms.

\sphinxAtStartPar
\sphinxcode{\sphinxupquote{mats}}: coefficient matrixes for added PSD terms.
\end{quote}
\end{quote}

\subsubsection{PsdExpr.clone()}
\label{\detokenize{javaapi/PsdExpr:psdexpr-clone}}\begin{quote}

\sphinxAtStartPar
Deep copy PSD expression object.

\sphinxAtStartPar
\sphinxstylestrong{Synopsis}
\begin{quote}

\sphinxAtStartPar
\sphinxcode{\sphinxupquote{PsdExpr clone()}}
\end{quote}

\sphinxAtStartPar
\sphinxstylestrong{Return}
\begin{quote}

\sphinxAtStartPar
cloned PSD expression object.
\end{quote}
\end{quote}

\subsubsection{PsdExpr.evaluate()}
\label{\detokenize{javaapi/PsdExpr:psdexpr-evaluate}}\begin{quote}

\sphinxAtStartPar
Evaluate PSD expression after solving.

\sphinxAtStartPar
\sphinxstylestrong{Synopsis}
\begin{quote}

\sphinxAtStartPar
\sphinxcode{\sphinxupquote{double evaluate()}}
\end{quote}

\sphinxAtStartPar
\sphinxstylestrong{Return}
\begin{quote}

\sphinxAtStartPar
Value of PSD expression.
\end{quote}
\end{quote}

\subsubsection{PsdExpr.getCoeff()}
\label{\detokenize{javaapi/PsdExpr:psdexpr-getcoeff}}\begin{quote}

\sphinxAtStartPar
Get coefficient from the i\sphinxhyphen{}th term in PSD expression.

\sphinxAtStartPar
\sphinxstylestrong{Synopsis}
\begin{quote}

\sphinxAtStartPar
\sphinxcode{\sphinxupquote{SymMatExpr getCoeff(int i)}}
\end{quote}

\sphinxAtStartPar
\sphinxstylestrong{Arguments}
\begin{quote}

\sphinxAtStartPar
\sphinxcode{\sphinxupquote{i}}: index of the PSD term.
\end{quote}

\sphinxAtStartPar
\sphinxstylestrong{Return}
\begin{quote}

\sphinxAtStartPar
coefficient expression of the i\sphinxhyphen{}th PSD term.
\end{quote}
\end{quote}

\subsubsection{PsdExpr.getConstant()}
\label{\detokenize{javaapi/PsdExpr:psdexpr-getconstant}}\begin{quote}

\sphinxAtStartPar
Get constant in PSD expression.

\sphinxAtStartPar
\sphinxstylestrong{Synopsis}
\begin{quote}

\sphinxAtStartPar
\sphinxcode{\sphinxupquote{double getConstant()}}
\end{quote}

\sphinxAtStartPar
\sphinxstylestrong{Return}
\begin{quote}

\sphinxAtStartPar
constant in PSD expression.
\end{quote}
\end{quote}

\subsubsection{PsdExpr.getLinExpr()}
\label{\detokenize{javaapi/PsdExpr:psdexpr-getlinexpr}}\begin{quote}

\sphinxAtStartPar
Get linear expression in PSD expression.

\sphinxAtStartPar
\sphinxstylestrong{Synopsis}
\begin{quote}

\sphinxAtStartPar
\sphinxcode{\sphinxupquote{Expr getLinExpr()}}
\end{quote}

\sphinxAtStartPar
\sphinxstylestrong{Return}
\begin{quote}

\sphinxAtStartPar
linear expression object.
\end{quote}
\end{quote}

\subsubsection{PsdExpr.getPsdVar()}
\label{\detokenize{javaapi/PsdExpr:psdexpr-getpsdvar}}\begin{quote}

\sphinxAtStartPar
Get the PSD variable from the i\sphinxhyphen{}th term in PSD expression.

\sphinxAtStartPar
\sphinxstylestrong{Synopsis}
\begin{quote}

\sphinxAtStartPar
\sphinxcode{\sphinxupquote{PsdVar getPsdVar(int i)}}
\end{quote}

\sphinxAtStartPar
\sphinxstylestrong{Arguments}
\begin{quote}

\sphinxAtStartPar
\sphinxcode{\sphinxupquote{i}}: index of the term.
\end{quote}

\sphinxAtStartPar
\sphinxstylestrong{Return}
\begin{quote}

\sphinxAtStartPar
the first variable of the i\sphinxhyphen{}th term in PSD expression object.
\end{quote}
\end{quote}

\subsubsection{PsdExpr.multiply()}
\label{\detokenize{javaapi/PsdExpr:psdexpr-multiply}}\begin{quote}

\sphinxAtStartPar
Multiply itself by a constant.

\sphinxAtStartPar
\sphinxstylestrong{Synopsis}
\begin{quote}

\sphinxAtStartPar
\sphinxcode{\sphinxupquote{PsdExpr multiply(double c)}}
\end{quote}

\sphinxAtStartPar
\sphinxstylestrong{Arguments}
\begin{quote}

\sphinxAtStartPar
\sphinxcode{\sphinxupquote{c}}: constant operand.
\end{quote}

\sphinxAtStartPar
\sphinxstylestrong{Return}
\begin{quote}

\sphinxAtStartPar
PSD expression itself.
\end{quote}
\end{quote}

\subsubsection{PsdExpr.remove()}
\label{\detokenize{javaapi/PsdExpr:psdexpr-remove}}\begin{quote}

\sphinxAtStartPar
Remove i\sphinxhyphen{}th term from PSD expression object.

\sphinxAtStartPar
\sphinxstylestrong{Synopsis}
\begin{quote}

\sphinxAtStartPar
\sphinxcode{\sphinxupquote{void remove(int idx)}}
\end{quote}

\sphinxAtStartPar
\sphinxstylestrong{Arguments}
\begin{quote}

\sphinxAtStartPar
\sphinxcode{\sphinxupquote{idx}}: index of the term to be removed.
\end{quote}
\end{quote}

\subsubsection{PsdExpr.remove()}
\label{\detokenize{javaapi/PsdExpr:id15}}\begin{quote}

\sphinxAtStartPar
Remove the term associated with variable from PSD expression.

\sphinxAtStartPar
\sphinxstylestrong{Synopsis}
\begin{quote}

\sphinxAtStartPar
\sphinxcode{\sphinxupquote{void remove(Var var)}}
\end{quote}

\sphinxAtStartPar
\sphinxstylestrong{Arguments}
\begin{quote}

\sphinxAtStartPar
\sphinxcode{\sphinxupquote{var}}: a variable whose term should be removed.
\end{quote}
\end{quote}

\subsubsection{PsdExpr.remove()}
\label{\detokenize{javaapi/PsdExpr:id16}}\begin{quote}

\sphinxAtStartPar
Remove the term associated with PSD variable from PSD expression.

\sphinxAtStartPar
\sphinxstylestrong{Synopsis}
\begin{quote}

\sphinxAtStartPar
\sphinxcode{\sphinxupquote{void remove(PsdVar var)}}
\end{quote}

\sphinxAtStartPar
\sphinxstylestrong{Arguments}
\begin{quote}

\sphinxAtStartPar
\sphinxcode{\sphinxupquote{var}}: a PSD variable whose term should be removed.
\end{quote}
\end{quote}

\subsubsection{PsdExpr.setCoeff()}
\label{\detokenize{javaapi/PsdExpr:psdexpr-setcoeff}}\begin{quote}

\sphinxAtStartPar
Set coefficient matrix of the i\sphinxhyphen{}th term in PSD expression.

\sphinxAtStartPar
\sphinxstylestrong{Synopsis}
\begin{quote}

\sphinxAtStartPar
\sphinxcode{\sphinxupquote{void setCoeff(int i, SymMatrix mat)}}
\end{quote}

\sphinxAtStartPar
\sphinxstylestrong{Arguments}
\begin{quote}

\sphinxAtStartPar
\sphinxcode{\sphinxupquote{i}}: index of the PSD term.

\sphinxAtStartPar
\sphinxcode{\sphinxupquote{mat}}: coefficient matrix of the term.
\end{quote}
\end{quote}

\subsubsection{PsdExpr.setConstant()}
\label{\detokenize{javaapi/PsdExpr:psdexpr-setconstant}}\begin{quote}

\sphinxAtStartPar
Set constant for the PSD expression.

\sphinxAtStartPar
\sphinxstylestrong{Synopsis}
\begin{quote}

\sphinxAtStartPar
\sphinxcode{\sphinxupquote{void setConstant(double constant)}}
\end{quote}

\sphinxAtStartPar
\sphinxstylestrong{Arguments}
\begin{quote}

\sphinxAtStartPar
\sphinxcode{\sphinxupquote{constant}}: the value of the constant.
\end{quote}
\end{quote}

\subsubsection{PsdExpr.size()}
\label{\detokenize{javaapi/PsdExpr:psdexpr-size}}\begin{quote}

\sphinxAtStartPar
Get number of PSD terms in expression.

\sphinxAtStartPar
\sphinxstylestrong{Synopsis}
\begin{quote}

\sphinxAtStartPar
\sphinxcode{\sphinxupquote{long size()}}
\end{quote}

\sphinxAtStartPar
\sphinxstylestrong{Return}
\begin{quote}

\sphinxAtStartPar
number of PSD terms.
\end{quote}
\end{quote}

\subsection{PsdConstraint}
\label{\detokenize{javaapiref:psdconstraint}}\label{\detokenize{javaapiref:chapjavaapiref-psdconstraint}}
\sphinxAtStartPar
COPT PSD constraint object. PSD constraints are always associated with a particular model.
User creates a PSD constraint object by adding a PSD constraint to model,
rather than by constructor of PsdConstraint class.

\sphinxstepscope

\subsubsection{PsdConstraint.get()}
\label{\detokenize{javaapi/PsdConstraint:psdconstraint-get}}\label{\detokenize{javaapi/PsdConstraint::doc}}\begin{quote}

\sphinxAtStartPar
Get information value of the PSD constraint. Support related PSD informations.

\sphinxAtStartPar
\sphinxstylestrong{Synopsis}
\begin{quote}

\sphinxAtStartPar
\sphinxcode{\sphinxupquote{double get(String info)}}
\end{quote}

\sphinxAtStartPar
\sphinxstylestrong{Arguments}
\begin{quote}

\sphinxAtStartPar
\sphinxcode{\sphinxupquote{info}}: name of queried information.
\end{quote}

\sphinxAtStartPar
\sphinxstylestrong{Return}
\begin{quote}

\sphinxAtStartPar
information value.
\end{quote}
\end{quote}

\subsubsection{PsdConstraint.getIdx()}
\label{\detokenize{javaapi/PsdConstraint:psdconstraint-getidx}}\begin{quote}

\sphinxAtStartPar
Get index of the PSD constraint.

\sphinxAtStartPar
\sphinxstylestrong{Synopsis}
\begin{quote}

\sphinxAtStartPar
\sphinxcode{\sphinxupquote{int getIdx()}}
\end{quote}

\sphinxAtStartPar
\sphinxstylestrong{Return}
\begin{quote}

\sphinxAtStartPar
the index of the PSD constraint.
\end{quote}
\end{quote}

\subsubsection{PsdConstraint.getName()}
\label{\detokenize{javaapi/PsdConstraint:psdconstraint-getname}}\begin{quote}

\sphinxAtStartPar
Get name of the PSD constraint.

\sphinxAtStartPar
\sphinxstylestrong{Synopsis}
\begin{quote}

\sphinxAtStartPar
\sphinxcode{\sphinxupquote{String getName()}}
\end{quote}

\sphinxAtStartPar
\sphinxstylestrong{Return}
\begin{quote}

\sphinxAtStartPar
the name of the PSD constraint.
\end{quote}
\end{quote}

\subsubsection{PsdConstraint.remove()}
\label{\detokenize{javaapi/PsdConstraint:psdconstraint-remove}}\begin{quote}

\sphinxAtStartPar
Remove this PSD constraint from model.

\sphinxAtStartPar
\sphinxstylestrong{Synopsis}
\begin{quote}

\sphinxAtStartPar
\sphinxcode{\sphinxupquote{void remove()}}
\end{quote}
\end{quote}

\subsubsection{PsdConstraint.set()}
\label{\detokenize{javaapi/PsdConstraint:psdconstraint-set}}\begin{quote}

\sphinxAtStartPar
Set information value of the PSD constraint. Support related PSD informations.

\sphinxAtStartPar
\sphinxstylestrong{Synopsis}
\begin{quote}

\sphinxAtStartPar
\sphinxcode{\sphinxupquote{void set(String info, double value)}}
\end{quote}

\sphinxAtStartPar
\sphinxstylestrong{Arguments}
\begin{quote}

\sphinxAtStartPar
\sphinxcode{\sphinxupquote{info}}: name of queried information.

\sphinxAtStartPar
\sphinxcode{\sphinxupquote{value}}: new information value.
\end{quote}
\end{quote}

\subsubsection{PsdConstraint.setName()}
\label{\detokenize{javaapi/PsdConstraint:psdconstraint-setname}}\begin{quote}

\sphinxAtStartPar
Set name of a PSD constraint.

\sphinxAtStartPar
\sphinxstylestrong{Synopsis}
\begin{quote}

\sphinxAtStartPar
\sphinxcode{\sphinxupquote{void setName(String name)}}
\end{quote}

\sphinxAtStartPar
\sphinxstylestrong{Arguments}
\begin{quote}

\sphinxAtStartPar
\sphinxcode{\sphinxupquote{name}}: the name to set.
\end{quote}
\end{quote}

\subsection{PsdConstrArray}
\label{\detokenize{javaapiref:psdconstrarray}}\label{\detokenize{javaapiref:chapjavaapiref-psdconstrarray}}
\sphinxAtStartPar
COPT PSD constraint array object. To store and access a set of {\hyperref[\detokenize{javaapiref:chapjavaapiref-psdconstraint}]{\sphinxcrossref{\DUrole{std,std-ref}{PsdConstraint}}}}
objects, Cardinal Optimizer provides PsdConstrArray class, which defines the following methods.

\sphinxstepscope

\subsubsection{PsdConstrArray.PsdConstrArray()}
\label{\detokenize{javaapi/PsdConstrArray:psdconstrarray-psdconstrarray}}\label{\detokenize{javaapi/PsdConstrArray::doc}}\begin{quote}

\sphinxAtStartPar
Constructor of PsdConstrArray object.

\sphinxAtStartPar
\sphinxstylestrong{Synopsis}
\begin{quote}

\sphinxAtStartPar
\sphinxcode{\sphinxupquote{PsdConstrArray()}}
\end{quote}
\end{quote}

\subsubsection{PsdConstrArray.getPsdConstr()}
\label{\detokenize{javaapi/PsdConstrArray:psdconstrarray-getpsdconstr}}\begin{quote}

\sphinxAtStartPar
Get idx\sphinxhyphen{}th PSD constraint object.

\sphinxAtStartPar
\sphinxstylestrong{Synopsis}
\begin{quote}

\sphinxAtStartPar
\sphinxcode{\sphinxupquote{PsdConstraint getPsdConstr(int idx)}}
\end{quote}

\sphinxAtStartPar
\sphinxstylestrong{Arguments}
\begin{quote}

\sphinxAtStartPar
\sphinxcode{\sphinxupquote{idx}}: index of the PSD constraint.
\end{quote}

\sphinxAtStartPar
\sphinxstylestrong{Return}
\begin{quote}

\sphinxAtStartPar
PSD constraint object with index idx.
\end{quote}
\end{quote}

\subsubsection{PsdConstrArray.pushBack()}
\label{\detokenize{javaapi/PsdConstrArray:psdconstrarray-pushback}}\begin{quote}

\sphinxAtStartPar
Add a PSD constraint object to PSD constraint array.

\sphinxAtStartPar
\sphinxstylestrong{Synopsis}
\begin{quote}

\sphinxAtStartPar
\sphinxcode{\sphinxupquote{void pushBack(PsdConstraint constr)}}
\end{quote}

\sphinxAtStartPar
\sphinxstylestrong{Arguments}
\begin{quote}

\sphinxAtStartPar
\sphinxcode{\sphinxupquote{constr}}: a PSD constraint object.
\end{quote}
\end{quote}

\subsubsection{PsdConstrArray.reserve()}
\label{\detokenize{javaapi/PsdConstrArray:psdconstrarray-reserve}}\begin{quote}

\sphinxAtStartPar
Reserve capacity to contain at least n items.

\sphinxAtStartPar
\sphinxstylestrong{Synopsis}
\begin{quote}

\sphinxAtStartPar
\sphinxcode{\sphinxupquote{void reserve(int n)}}
\end{quote}

\sphinxAtStartPar
\sphinxstylestrong{Arguments}
\begin{quote}

\sphinxAtStartPar
\sphinxcode{\sphinxupquote{n}}: minimum capacity for PSD constraint objects.
\end{quote}
\end{quote}

\subsubsection{PsdConstrArray.size()}
\label{\detokenize{javaapi/PsdConstrArray:psdconstrarray-size}}\begin{quote}

\sphinxAtStartPar
Get the number of PSD constraint objects.

\sphinxAtStartPar
\sphinxstylestrong{Synopsis}
\begin{quote}

\sphinxAtStartPar
\sphinxcode{\sphinxupquote{int size()}}
\end{quote}

\sphinxAtStartPar
\sphinxstylestrong{Return}
\begin{quote}

\sphinxAtStartPar
number of PSD constraint objects.
\end{quote}
\end{quote}

\subsection{PsdConstrBuilder}
\label{\detokenize{javaapiref:psdconstrbuilder}}\label{\detokenize{javaapiref:chapjavaapiref-psdconstrbuilder}}
\sphinxAtStartPar
COPT PSD constraint builder object. To help building a PSD constraint, given a PSD
expression, constraint sense and right\sphinxhyphen{}hand side value, Cardinal Optimizer provides PsdConstrBuilder
class, which defines the following methods.

\sphinxstepscope

\subsubsection{PsdConstrBuilder.PsdConstrBuilder()}
\label{\detokenize{javaapi/PsdConstrBuilder:psdconstrbuilder-psdconstrbuilder}}\label{\detokenize{javaapi/PsdConstrBuilder::doc}}\begin{quote}

\sphinxAtStartPar
Constructor of PsdConstrBuilder object.

\sphinxAtStartPar
\sphinxstylestrong{Synopsis}
\begin{quote}

\sphinxAtStartPar
\sphinxcode{\sphinxupquote{PsdConstrBuilder()}}
\end{quote}
\end{quote}

\subsubsection{PsdConstrBuilder.getPsdExpr()}
\label{\detokenize{javaapi/PsdConstrBuilder:psdconstrbuilder-getpsdexpr}}\begin{quote}

\sphinxAtStartPar
Get expression associated with PSD constraint.

\sphinxAtStartPar
\sphinxstylestrong{Synopsis}
\begin{quote}

\sphinxAtStartPar
\sphinxcode{\sphinxupquote{PsdExpr getPsdExpr()}}
\end{quote}

\sphinxAtStartPar
\sphinxstylestrong{Return}
\begin{quote}

\sphinxAtStartPar
PSD expression object.
\end{quote}
\end{quote}

\subsubsection{PsdConstrBuilder.getRange()}
\label{\detokenize{javaapi/PsdConstrBuilder:psdconstrbuilder-getrange}}\begin{quote}

\sphinxAtStartPar
Get range from lower bound to upper bound of range constraint.

\sphinxAtStartPar
\sphinxstylestrong{Synopsis}
\begin{quote}

\sphinxAtStartPar
\sphinxcode{\sphinxupquote{double getRange()}}
\end{quote}

\sphinxAtStartPar
\sphinxstylestrong{Return}
\begin{quote}

\sphinxAtStartPar
length from lower bound to upper bound of the constraint.
\end{quote}
\end{quote}

\subsubsection{PsdConstrBuilder.getSense()}
\label{\detokenize{javaapi/PsdConstrBuilder:psdconstrbuilder-getsense}}\begin{quote}

\sphinxAtStartPar
Get sense associated with PSD constraint.

\sphinxAtStartPar
\sphinxstylestrong{Synopsis}
\begin{quote}

\sphinxAtStartPar
\sphinxcode{\sphinxupquote{char getSense()}}
\end{quote}

\sphinxAtStartPar
\sphinxstylestrong{Return}
\begin{quote}

\sphinxAtStartPar
PSD constraint sense.
\end{quote}
\end{quote}

\subsubsection{PsdConstrBuilder.set()}
\label{\detokenize{javaapi/PsdConstrBuilder:psdconstrbuilder-set}}\begin{quote}

\sphinxAtStartPar
Set detail of a PSD constraint to its builder object.

\sphinxAtStartPar
\sphinxstylestrong{Synopsis}
\begin{quote}

\sphinxAtStartPar
\sphinxcode{\sphinxupquote{void set(}}
\begin{quote}

\sphinxAtStartPar
\sphinxcode{\sphinxupquote{PsdExpr expr,}}

\sphinxAtStartPar
\sphinxcode{\sphinxupquote{char sense,}}

\sphinxAtStartPar
\sphinxcode{\sphinxupquote{double rhs)}}
\end{quote}
\end{quote}

\sphinxAtStartPar
\sphinxstylestrong{Arguments}
\begin{quote}

\sphinxAtStartPar
\sphinxcode{\sphinxupquote{expr}}: expression object at one side of the PSD constraint.

\sphinxAtStartPar
\sphinxcode{\sphinxupquote{sense}}: PSD constraint sense, other than COPT\_RANGE.

\sphinxAtStartPar
\sphinxcode{\sphinxupquote{rhs}}: constant at right side of the PSD constraint.
\end{quote}
\end{quote}

\subsubsection{PsdConstrBuilder.setRange()}
\label{\detokenize{javaapi/PsdConstrBuilder:psdconstrbuilder-setrange}}\begin{quote}

\sphinxAtStartPar
Set a range constraint to its builder.

\sphinxAtStartPar
\sphinxstylestrong{Synopsis}
\begin{quote}

\sphinxAtStartPar
\sphinxcode{\sphinxupquote{void setRange(PsdExpr expr, double range)}}
\end{quote}

\sphinxAtStartPar
\sphinxstylestrong{Arguments}
\begin{quote}

\sphinxAtStartPar
\sphinxcode{\sphinxupquote{expr}}: PSD expression object, whose constant is negative upper bound.

\sphinxAtStartPar
\sphinxcode{\sphinxupquote{range}}: length from lower bound to upper bound of the constraint. Must greater than 0.
\end{quote}
\end{quote}

\subsection{PsdConstrBuilderArray}
\label{\detokenize{javaapiref:psdconstrbuilderarray}}\label{\detokenize{javaapiref:chapjavaapiref-psdconstrbuilderarray}}
\sphinxAtStartPar
COPT PSD constraint builder array object. To store and access a set of {\hyperref[\detokenize{javaapiref:chapjavaapiref-psdconstrbuilder}]{\sphinxcrossref{\DUrole{std,std-ref}{PsdConstrBuilder}}}}
objects, Cardinal Optimizer provides PsdConstrBuilderArray class, which defines the following methods.

\sphinxstepscope

\subsubsection{PsdConstrBuilderArray.PsdConstrBuilderArray()}
\label{\detokenize{javaapi/PsdConstrBuilderArray:psdconstrbuilderarray-psdconstrbuilderarray}}\label{\detokenize{javaapi/PsdConstrBuilderArray::doc}}\begin{quote}

\sphinxAtStartPar
Constructor of PsdConstrBuilderArray object.

\sphinxAtStartPar
\sphinxstylestrong{Synopsis}
\begin{quote}

\sphinxAtStartPar
\sphinxcode{\sphinxupquote{PsdConstrBuilderArray()}}
\end{quote}
\end{quote}

\subsubsection{PsdConstrBuilderArray.getBuilder()}
\label{\detokenize{javaapi/PsdConstrBuilderArray:psdconstrbuilderarray-getbuilder}}\begin{quote}

\sphinxAtStartPar
Get idx\sphinxhyphen{}th PSD constraint builder object.

\sphinxAtStartPar
\sphinxstylestrong{Synopsis}
\begin{quote}

\sphinxAtStartPar
\sphinxcode{\sphinxupquote{PsdConstrBuilder getBuilder(int idx)}}
\end{quote}

\sphinxAtStartPar
\sphinxstylestrong{Arguments}
\begin{quote}

\sphinxAtStartPar
\sphinxcode{\sphinxupquote{idx}}: index of the PSD constraint builder.
\end{quote}

\sphinxAtStartPar
\sphinxstylestrong{Return}
\begin{quote}

\sphinxAtStartPar
PSD constraint builder object with index idx.
\end{quote}
\end{quote}

\subsubsection{PsdConstrBuilderArray.pushBack()}
\label{\detokenize{javaapi/PsdConstrBuilderArray:psdconstrbuilderarray-pushback}}\begin{quote}

\sphinxAtStartPar
Add a PSD constraint builder to PSD constraint builder array.

\sphinxAtStartPar
\sphinxstylestrong{Synopsis}
\begin{quote}

\sphinxAtStartPar
\sphinxcode{\sphinxupquote{void pushBack(PsdConstrBuilder builder)}}
\end{quote}

\sphinxAtStartPar
\sphinxstylestrong{Arguments}
\begin{quote}

\sphinxAtStartPar
\sphinxcode{\sphinxupquote{builder}}: a PSD constraint builder object.
\end{quote}
\end{quote}

\subsubsection{PsdConstrBuilderArray.reserve()}
\label{\detokenize{javaapi/PsdConstrBuilderArray:psdconstrbuilderarray-reserve}}\begin{quote}

\sphinxAtStartPar
Reserve capacity to contain at least n items.

\sphinxAtStartPar
\sphinxstylestrong{Synopsis}
\begin{quote}

\sphinxAtStartPar
\sphinxcode{\sphinxupquote{void reserve(int n)}}
\end{quote}

\sphinxAtStartPar
\sphinxstylestrong{Arguments}
\begin{quote}

\sphinxAtStartPar
\sphinxcode{\sphinxupquote{n}}: minimum capacity for PSD constraint builder object.
\end{quote}
\end{quote}

\subsubsection{PsdConstrBuilderArray.size()}
\label{\detokenize{javaapi/PsdConstrBuilderArray:psdconstrbuilderarray-size}}\begin{quote}

\sphinxAtStartPar
Get the number of PSD constraint builder objects.

\sphinxAtStartPar
\sphinxstylestrong{Synopsis}
\begin{quote}

\sphinxAtStartPar
\sphinxcode{\sphinxupquote{int size()}}
\end{quote}

\sphinxAtStartPar
\sphinxstylestrong{Return}
\begin{quote}

\sphinxAtStartPar
number of PSD constraint builder objects.
\end{quote}
\end{quote}

\subsection{LmiConstraint}
\label{\detokenize{javaapiref:lmiconstraint}}\label{\detokenize{javaapiref:chapjavaapiref-lmiconstraint}}
\sphinxAtStartPar
COPT LMI constraint object. LMI constraints are always associated with a
particular model.  User creates a LMI constraint object by adding a LMI
constraint to model, rather than by constructor of LmiConstraint class.

\sphinxstepscope

\subsubsection{LmiConstraint.get()}
\label{\detokenize{javaapi/LmiConstraint:lmiconstraint-get}}\label{\detokenize{javaapi/LmiConstraint::doc}}\begin{quote}

\sphinxAtStartPar
Get information values of LMI constraint.

\sphinxAtStartPar
\sphinxstylestrong{Synopsis}
\begin{quote}

\sphinxAtStartPar
\sphinxcode{\sphinxupquote{double{[}{]} get(String info)}}
\end{quote}

\sphinxAtStartPar
\sphinxstylestrong{Arguments}
\begin{quote}

\sphinxAtStartPar
\sphinxcode{\sphinxupquote{info}}: information name.
\end{quote}

\sphinxAtStartPar
\sphinxstylestrong{Return}
\begin{quote}

\sphinxAtStartPar
array of information values.
\end{quote}
\end{quote}

\subsubsection{LmiConstraint.getDim()}
\label{\detokenize{javaapi/LmiConstraint:lmiconstraint-getdim}}\begin{quote}

\sphinxAtStartPar
Get dimension of LMI constraint.

\sphinxAtStartPar
\sphinxstylestrong{Synopsis}
\begin{quote}

\sphinxAtStartPar
\sphinxcode{\sphinxupquote{int getDim()}}
\end{quote}

\sphinxAtStartPar
\sphinxstylestrong{Return}
\begin{quote}

\sphinxAtStartPar
dimension of LMI constraint.
\end{quote}
\end{quote}

\subsubsection{LmiConstraint.getIdx()}
\label{\detokenize{javaapi/LmiConstraint:lmiconstraint-getidx}}\begin{quote}

\sphinxAtStartPar
Get index of LMI constraint.

\sphinxAtStartPar
\sphinxstylestrong{Synopsis}
\begin{quote}

\sphinxAtStartPar
\sphinxcode{\sphinxupquote{int getIdx()}}
\end{quote}

\sphinxAtStartPar
\sphinxstylestrong{Return}
\begin{quote}

\sphinxAtStartPar
index of LMI constraint.
\end{quote}
\end{quote}

\subsubsection{LmiConstraint.getLen()}
\label{\detokenize{javaapi/LmiConstraint:lmiconstraint-getlen}}\begin{quote}

\sphinxAtStartPar
Get length of LMI constraint.

\sphinxAtStartPar
\sphinxstylestrong{Synopsis}
\begin{quote}

\sphinxAtStartPar
\sphinxcode{\sphinxupquote{int getLen()}}
\end{quote}

\sphinxAtStartPar
\sphinxstylestrong{Return}
\begin{quote}

\sphinxAtStartPar
length of LMI constraint.
\end{quote}
\end{quote}

\subsubsection{LmiConstraint.getName()}
\label{\detokenize{javaapi/LmiConstraint:lmiconstraint-getname}}\begin{quote}

\sphinxAtStartPar
Get name of LMI constraint.

\sphinxAtStartPar
\sphinxstylestrong{Synopsis}
\begin{quote}

\sphinxAtStartPar
\sphinxcode{\sphinxupquote{String getName()}}
\end{quote}

\sphinxAtStartPar
\sphinxstylestrong{Return}
\begin{quote}

\sphinxAtStartPar
name of LMI constraint.
\end{quote}
\end{quote}

\subsubsection{LmiConstraint.remove()}
\label{\detokenize{javaapi/LmiConstraint:lmiconstraint-remove}}\begin{quote}

\sphinxAtStartPar
Remove this LMI constraint from model.

\sphinxAtStartPar
\sphinxstylestrong{Synopsis}
\begin{quote}

\sphinxAtStartPar
\sphinxcode{\sphinxupquote{void remove()}}
\end{quote}
\end{quote}

\subsubsection{LmiConstraint.setRhs()}
\label{\detokenize{javaapi/LmiConstraint:lmiconstraint-setrhs}}\begin{quote}

\sphinxAtStartPar
Set constant term of LMI constraint.

\sphinxAtStartPar
\sphinxstylestrong{Synopsis}
\begin{quote}

\sphinxAtStartPar
\sphinxcode{\sphinxupquote{void setRhs(SymMatrix mat)}}
\end{quote}

\sphinxAtStartPar
\sphinxstylestrong{Arguments}
\begin{quote}

\sphinxAtStartPar
\sphinxcode{\sphinxupquote{mat}}: new symmetric matrix for constant term.
\end{quote}
\end{quote}

\subsection{LmiConstrArray}
\label{\detokenize{javaapiref:lmiconstrarray}}\label{\detokenize{javaapiref:chapjavaapiref-lmiconstrarray}}
\sphinxAtStartPar
COPT LMI constraint array object. To store and access a set of
{\hyperref[\detokenize{javaapiref:chapjavaapiref-lmiconstraint}]{\sphinxcrossref{\DUrole{std,std-ref}{LmiConstraint}}}} objects, Cardinal Optimizer provides
LmiConstrArray class, which defines the following methods.

\sphinxstepscope

\subsubsection{LmiConstrArray.LmiConstrArray()}
\label{\detokenize{javaapi/LmiConstrArray:lmiconstrarray-lmiconstrarray}}\label{\detokenize{javaapi/LmiConstrArray::doc}}\begin{quote}

\sphinxAtStartPar
Constructor of LmiConstrArray.

\sphinxAtStartPar
\sphinxstylestrong{Synopsis}
\begin{quote}

\sphinxAtStartPar
\sphinxcode{\sphinxupquote{LmiConstrArray()}}
\end{quote}
\end{quote}

\subsubsection{LmiConstrArray.getLmiConstr()}
\label{\detokenize{javaapi/LmiConstrArray:lmiconstrarray-getlmiconstr}}\begin{quote}

\sphinxAtStartPar
Get idx\sphinxhyphen{}th LMI constraint object.

\sphinxAtStartPar
\sphinxstylestrong{Synopsis}
\begin{quote}

\sphinxAtStartPar
\sphinxcode{\sphinxupquote{LmiConstraint getLmiConstr(int idx)}}
\end{quote}

\sphinxAtStartPar
\sphinxstylestrong{Arguments}
\begin{quote}

\sphinxAtStartPar
\sphinxcode{\sphinxupquote{idx}}: index of the LMI constraint.
\end{quote}

\sphinxAtStartPar
\sphinxstylestrong{Return}
\begin{quote}

\sphinxAtStartPar
LMI constraint object with index idx.
\end{quote}
\end{quote}

\subsubsection{LmiConstrArray.pushBack()}
\label{\detokenize{javaapi/LmiConstrArray:lmiconstrarray-pushback}}\begin{quote}

\sphinxAtStartPar
Add an LMI constraint to LMI constraint array.

\sphinxAtStartPar
\sphinxstylestrong{Synopsis}
\begin{quote}

\sphinxAtStartPar
\sphinxcode{\sphinxupquote{void pushBack(LmiConstraint constr)}}
\end{quote}

\sphinxAtStartPar
\sphinxstylestrong{Arguments}
\begin{quote}

\sphinxAtStartPar
\sphinxcode{\sphinxupquote{constr}}: LMI constraint object.
\end{quote}
\end{quote}

\subsubsection{LmiConstrArray.reserve()}
\label{\detokenize{javaapi/LmiConstrArray:lmiconstrarray-reserve}}\begin{quote}

\sphinxAtStartPar
Reserve capacity to contain at least n items.

\sphinxAtStartPar
\sphinxstylestrong{Synopsis}
\begin{quote}

\sphinxAtStartPar
\sphinxcode{\sphinxupquote{void reserve(int n)}}
\end{quote}

\sphinxAtStartPar
\sphinxstylestrong{Arguments}
\begin{quote}

\sphinxAtStartPar
\sphinxcode{\sphinxupquote{n}}: capacity number of LMI constraint objects.
\end{quote}
\end{quote}

\subsubsection{LmiConstrArray.size()}
\label{\detokenize{javaapi/LmiConstrArray:lmiconstrarray-size}}\begin{quote}

\sphinxAtStartPar
Get the number of LMI constraint objects.

\sphinxAtStartPar
\sphinxstylestrong{Synopsis}
\begin{quote}

\sphinxAtStartPar
\sphinxcode{\sphinxupquote{int size()}}
\end{quote}

\sphinxAtStartPar
\sphinxstylestrong{Return}
\begin{quote}

\sphinxAtStartPar
number of LMI constraint objects.
\end{quote}
\end{quote}

\subsection{LmiExpr}
\label{\detokenize{javaapiref:lmiexpr}}\label{\detokenize{javaapiref:chapjavaapiref-lmiexpr}}
\sphinxAtStartPar
COPT LMI expression object. A LMI expression consists of
a list of variables, associated coefficient matrices of LMI term, and constant matrices.
LMI expressions are used to build LMI constraints.

\sphinxstepscope

\subsubsection{LmiExpr.LmiExpr()}
\label{\detokenize{javaapi/LmiExpr:lmiexpr-lmiexpr}}\label{\detokenize{javaapi/LmiExpr::doc}}\begin{quote}

\sphinxAtStartPar
Default constructor of a LMI expression.

\sphinxAtStartPar
\sphinxstylestrong{Synopsis}
\begin{quote}

\sphinxAtStartPar
\sphinxcode{\sphinxupquote{LmiExpr()}}
\end{quote}
\end{quote}

\subsubsection{LmiExpr.LmiExpr()}
\label{\detokenize{javaapi/LmiExpr:id1}}\begin{quote}

\sphinxAtStartPar
Constructor of LMI expression with constant term.

\sphinxAtStartPar
\sphinxstylestrong{Synopsis}
\begin{quote}

\sphinxAtStartPar
\sphinxcode{\sphinxupquote{LmiExpr(SymMatrix mat)}}
\end{quote}

\sphinxAtStartPar
\sphinxstylestrong{Arguments}
\begin{quote}

\sphinxAtStartPar
\sphinxcode{\sphinxupquote{mat}}: symmetric matrix object.
\end{quote}
\end{quote}

\subsubsection{LmiExpr.LmiExpr()}
\label{\detokenize{javaapi/LmiExpr:id2}}\begin{quote}

\sphinxAtStartPar
Constructor of LMI expression with constant term.

\sphinxAtStartPar
\sphinxstylestrong{Synopsis}
\begin{quote}

\sphinxAtStartPar
\sphinxcode{\sphinxupquote{LmiExpr(SymMatExpr expr)}}
\end{quote}

\sphinxAtStartPar
\sphinxstylestrong{Arguments}
\begin{quote}

\sphinxAtStartPar
\sphinxcode{\sphinxupquote{expr}}: matrix expression object.
\end{quote}
\end{quote}

\subsubsection{LmiExpr.LmiExpr()}
\label{\detokenize{javaapi/LmiExpr:id3}}\begin{quote}

\sphinxAtStartPar
Constructor of LMI expression with one term.

\sphinxAtStartPar
\sphinxstylestrong{Synopsis}
\begin{quote}

\sphinxAtStartPar
\sphinxcode{\sphinxupquote{LmiExpr(Var var, SymMatrix mat)}}
\end{quote}

\sphinxAtStartPar
\sphinxstylestrong{Arguments}
\begin{quote}

\sphinxAtStartPar
\sphinxcode{\sphinxupquote{var}}: variable of the added term.

\sphinxAtStartPar
\sphinxcode{\sphinxupquote{mat}}: coefficient matrix of the added term.
\end{quote}
\end{quote}

\subsubsection{LmiExpr.LmiExpr()}
\label{\detokenize{javaapi/LmiExpr:id4}}\begin{quote}

\sphinxAtStartPar
Constructor of LMI expression with one term.

\sphinxAtStartPar
\sphinxstylestrong{Synopsis}
\begin{quote}

\sphinxAtStartPar
\sphinxcode{\sphinxupquote{LmiExpr(Var var, SymMatExpr expr)}}
\end{quote}

\sphinxAtStartPar
\sphinxstylestrong{Arguments}
\begin{quote}

\sphinxAtStartPar
\sphinxcode{\sphinxupquote{var}}: variable of the added term.

\sphinxAtStartPar
\sphinxcode{\sphinxupquote{expr}}: coefficient expression of symmetric matrices of new LMI term.
\end{quote}
\end{quote}

\subsubsection{LmiExpr.addConstant()}
\label{\detokenize{javaapi/LmiExpr:lmiexpr-addconstant}}\begin{quote}

\sphinxAtStartPar
Add to constant term of the LMI expression.

\sphinxAtStartPar
\sphinxstylestrong{Synopsis}
\begin{quote}

\sphinxAtStartPar
\sphinxcode{\sphinxupquote{void addConstant(SymMatExpr expr)}}
\end{quote}

\sphinxAtStartPar
\sphinxstylestrong{Arguments}
\begin{quote}

\sphinxAtStartPar
\sphinxcode{\sphinxupquote{expr}}: matrix expression added to the constant term.
\end{quote}
\end{quote}

\subsubsection{LmiExpr.addLmiExpr()}
\label{\detokenize{javaapi/LmiExpr:lmiexpr-addlmiexpr}}\begin{quote}

\sphinxAtStartPar
Add an LMI expression to self.

\sphinxAtStartPar
\sphinxstylestrong{Synopsis}
\begin{quote}

\sphinxAtStartPar
\sphinxcode{\sphinxupquote{void addLmiExpr(LmiExpr expr)}}
\end{quote}

\sphinxAtStartPar
\sphinxstylestrong{Arguments}
\begin{quote}

\sphinxAtStartPar
\sphinxcode{\sphinxupquote{expr}}: LMI expression to be added.
\end{quote}
\end{quote}

\subsubsection{LmiExpr.addLmiExpr()}
\label{\detokenize{javaapi/LmiExpr:id5}}\begin{quote}

\sphinxAtStartPar
Add an LMI expression to self.

\sphinxAtStartPar
\sphinxstylestrong{Synopsis}
\begin{quote}

\sphinxAtStartPar
\sphinxcode{\sphinxupquote{void addLmiExpr(LmiExpr expr, double mult)}}
\end{quote}

\sphinxAtStartPar
\sphinxstylestrong{Arguments}
\begin{quote}

\sphinxAtStartPar
\sphinxcode{\sphinxupquote{expr}}: LMI expression to be added.

\sphinxAtStartPar
\sphinxcode{\sphinxupquote{mult}}: multiplier constant.
\end{quote}
\end{quote}

\subsubsection{LmiExpr.addTerm()}
\label{\detokenize{javaapi/LmiExpr:lmiexpr-addterm}}\begin{quote}

\sphinxAtStartPar
Add an LMI term to LMI expression object.

\sphinxAtStartPar
\sphinxstylestrong{Synopsis}
\begin{quote}

\sphinxAtStartPar
\sphinxcode{\sphinxupquote{void addTerm(Var var, SymMatrix mat)}}
\end{quote}

\sphinxAtStartPar
\sphinxstylestrong{Arguments}
\begin{quote}

\sphinxAtStartPar
\sphinxcode{\sphinxupquote{var}}: LMI variable of new LMI term.

\sphinxAtStartPar
\sphinxcode{\sphinxupquote{mat}}: coefficient matrix of new LMI term.
\end{quote}
\end{quote}

\subsubsection{LmiExpr.addTerm()}
\label{\detokenize{javaapi/LmiExpr:id6}}\begin{quote}

\sphinxAtStartPar
Add an LMI term to LMI expression object.

\sphinxAtStartPar
\sphinxstylestrong{Synopsis}
\begin{quote}

\sphinxAtStartPar
\sphinxcode{\sphinxupquote{void addTerm(Var var, SymMatExpr expr)}}
\end{quote}

\sphinxAtStartPar
\sphinxstylestrong{Arguments}
\begin{quote}

\sphinxAtStartPar
\sphinxcode{\sphinxupquote{var}}: variable of new LMI term.

\sphinxAtStartPar
\sphinxcode{\sphinxupquote{expr}}: coefficient expression of symmetric matrices of new LMI term.
\end{quote}
\end{quote}

\subsubsection{LmiExpr.addTerms()}
\label{\detokenize{javaapi/LmiExpr:lmiexpr-addterms}}\begin{quote}

\sphinxAtStartPar
Add LMI terms to LMI expression object.

\sphinxAtStartPar
\sphinxstylestrong{Synopsis}
\begin{quote}

\sphinxAtStartPar
\sphinxcode{\sphinxupquote{void addTerms(VarArray vars, SymMatrixArray mats)}}
\end{quote}

\sphinxAtStartPar
\sphinxstylestrong{Arguments}
\begin{quote}

\sphinxAtStartPar
\sphinxcode{\sphinxupquote{vars}}: variables for added LMI terms.

\sphinxAtStartPar
\sphinxcode{\sphinxupquote{mats}}: coefficient matrices for added LMI terms.
\end{quote}
\end{quote}

\subsubsection{LmiExpr.addTerms()}
\label{\detokenize{javaapi/LmiExpr:id7}}\begin{quote}

\sphinxAtStartPar
Add LMI terms to LMI expression object.

\sphinxAtStartPar
\sphinxstylestrong{Synopsis}
\begin{quote}

\sphinxAtStartPar
\sphinxcode{\sphinxupquote{void addTerms(Var{[}{]} vars, SymMatrix{[}{]} mats)}}
\end{quote}

\sphinxAtStartPar
\sphinxstylestrong{Arguments}
\begin{quote}

\sphinxAtStartPar
\sphinxcode{\sphinxupquote{vars}}: variables for added LMI terms.

\sphinxAtStartPar
\sphinxcode{\sphinxupquote{mats}}: coefficient matrices for added LMI terms.
\end{quote}
\end{quote}

\subsubsection{LmiExpr.clone()}
\label{\detokenize{javaapi/LmiExpr:lmiexpr-clone}}\begin{quote}

\sphinxAtStartPar
Deep copy LMI expression.

\sphinxAtStartPar
\sphinxstylestrong{Synopsis}
\begin{quote}

\sphinxAtStartPar
\sphinxcode{\sphinxupquote{LmiExpr clone()}}
\end{quote}

\sphinxAtStartPar
\sphinxstylestrong{Return}
\begin{quote}

\sphinxAtStartPar
cloned LMI expression object.
\end{quote}
\end{quote}

\subsubsection{LmiExpr.getCoeff()}
\label{\detokenize{javaapi/LmiExpr:lmiexpr-getcoeff}}\begin{quote}

\sphinxAtStartPar
Get coefficient from the i\sphinxhyphen{}th term in LMI expression.

\sphinxAtStartPar
\sphinxstylestrong{Synopsis}
\begin{quote}

\sphinxAtStartPar
\sphinxcode{\sphinxupquote{SymMatExpr getCoeff(int i)}}
\end{quote}

\sphinxAtStartPar
\sphinxstylestrong{Arguments}
\begin{quote}

\sphinxAtStartPar
\sphinxcode{\sphinxupquote{i}}: index of the LMI term.
\end{quote}

\sphinxAtStartPar
\sphinxstylestrong{Return}
\begin{quote}

\sphinxAtStartPar
coefficient expression of the i\sphinxhyphen{}th LMI term.
\end{quote}
\end{quote}

\subsubsection{LmiExpr.getConstant()}
\label{\detokenize{javaapi/LmiExpr:lmiexpr-getconstant}}\begin{quote}

\sphinxAtStartPar
Get constant term in LMI expression.

\sphinxAtStartPar
\sphinxstylestrong{Synopsis}
\begin{quote}

\sphinxAtStartPar
\sphinxcode{\sphinxupquote{SymMatExpr getConstant()}}
\end{quote}

\sphinxAtStartPar
\sphinxstylestrong{Return}
\begin{quote}

\sphinxAtStartPar
symmetric matrix expression object.
\end{quote}
\end{quote}

\subsubsection{LmiExpr.getVar()}
\label{\detokenize{javaapi/LmiExpr:lmiexpr-getvar}}\begin{quote}

\sphinxAtStartPar
Get variable from the i\sphinxhyphen{}th term in LMI expression.

\sphinxAtStartPar
\sphinxstylestrong{Synopsis}
\begin{quote}

\sphinxAtStartPar
\sphinxcode{\sphinxupquote{Var getVar(int i)}}
\end{quote}

\sphinxAtStartPar
\sphinxstylestrong{Arguments}
\begin{quote}

\sphinxAtStartPar
\sphinxcode{\sphinxupquote{i}}: index of the term.
\end{quote}

\sphinxAtStartPar
\sphinxstylestrong{Return}
\begin{quote}

\sphinxAtStartPar
variable of the i\sphinxhyphen{}th term in LMI expression object.
\end{quote}
\end{quote}

\subsubsection{LmiExpr.multiply()}
\label{\detokenize{javaapi/LmiExpr:lmiexpr-multiply}}\begin{quote}

\sphinxAtStartPar
Multiply itself by a constant.

\sphinxAtStartPar
\sphinxstylestrong{Synopsis}
\begin{quote}

\sphinxAtStartPar
\sphinxcode{\sphinxupquote{LmiExpr multiply(double c)}}
\end{quote}

\sphinxAtStartPar
\sphinxstylestrong{Arguments}
\begin{quote}

\sphinxAtStartPar
\sphinxcode{\sphinxupquote{c}}: constant operand.
\end{quote}

\sphinxAtStartPar
\sphinxstylestrong{Return}
\begin{quote}

\sphinxAtStartPar
LMI expression itself.
\end{quote}
\end{quote}

\subsubsection{LmiExpr.remove()}
\label{\detokenize{javaapi/LmiExpr:lmiexpr-remove}}\begin{quote}

\sphinxAtStartPar
Remove i\sphinxhyphen{}th term from LMI expression object.

\sphinxAtStartPar
\sphinxstylestrong{Synopsis}
\begin{quote}

\sphinxAtStartPar
\sphinxcode{\sphinxupquote{void remove(int idx)}}
\end{quote}

\sphinxAtStartPar
\sphinxstylestrong{Arguments}
\begin{quote}

\sphinxAtStartPar
\sphinxcode{\sphinxupquote{idx}}: index of the term to be removed.
\end{quote}
\end{quote}

\subsubsection{LmiExpr.remove()}
\label{\detokenize{javaapi/LmiExpr:id8}}\begin{quote}

\sphinxAtStartPar
Remove the term associated with variable from LMI expression.

\sphinxAtStartPar
\sphinxstylestrong{Synopsis}
\begin{quote}

\sphinxAtStartPar
\sphinxcode{\sphinxupquote{void remove(Var var)}}
\end{quote}

\sphinxAtStartPar
\sphinxstylestrong{Arguments}
\begin{quote}

\sphinxAtStartPar
\sphinxcode{\sphinxupquote{var}}: variable whose term should be removed.
\end{quote}
\end{quote}

\subsubsection{LmiExpr.setCoeff()}
\label{\detokenize{javaapi/LmiExpr:lmiexpr-setcoeff}}\begin{quote}

\sphinxAtStartPar
Set coefficient matrix of the i\sphinxhyphen{}th term in LMI expression.

\sphinxAtStartPar
\sphinxstylestrong{Synopsis}
\begin{quote}

\sphinxAtStartPar
\sphinxcode{\sphinxupquote{void setCoeff(int i, SymMatrix mat)}}
\end{quote}

\sphinxAtStartPar
\sphinxstylestrong{Arguments}
\begin{quote}

\sphinxAtStartPar
\sphinxcode{\sphinxupquote{i}}: index of the LMI term.

\sphinxAtStartPar
\sphinxcode{\sphinxupquote{mat}}: coefficient matrix of the term.
\end{quote}
\end{quote}

\subsubsection{LmiExpr.setConstant()}
\label{\detokenize{javaapi/LmiExpr:lmiexpr-setconstant}}\begin{quote}

\sphinxAtStartPar
Set constant term of the LMI expression.

\sphinxAtStartPar
\sphinxstylestrong{Synopsis}
\begin{quote}

\sphinxAtStartPar
\sphinxcode{\sphinxupquote{void setConstant(SymMatrix mat)}}
\end{quote}

\sphinxAtStartPar
\sphinxstylestrong{Arguments}
\begin{quote}

\sphinxAtStartPar
\sphinxcode{\sphinxupquote{mat}}: symmetric matrix of the constant term.
\end{quote}
\end{quote}

\subsubsection{LmiExpr.size()}
\label{\detokenize{javaapi/LmiExpr:lmiexpr-size}}\begin{quote}

\sphinxAtStartPar
Get number of LMI terms in expression.

\sphinxAtStartPar
\sphinxstylestrong{Synopsis}
\begin{quote}

\sphinxAtStartPar
\sphinxcode{\sphinxupquote{long size()}}
\end{quote}

\sphinxAtStartPar
\sphinxstylestrong{Return}
\begin{quote}

\sphinxAtStartPar
number of LMI terms.
\end{quote}
\end{quote}

\subsection{SymMatrix}
\label{\detokenize{javaapiref:symmatrix}}\label{\detokenize{javaapiref:chapjavaapiref-symmatrix}}
\sphinxAtStartPar
COPT symmetric matrix object. Symmetric matrices are always associated with a particular model.
User creates a symmetric matrix object by adding a symmetric matrix to model,
rather than by constructor of SymMatrix class.

\sphinxAtStartPar
Symmetric matrices are used as coefficient matrices of PSD terms in PSD expressions, PSD constraints or PSD objectives.

\sphinxstepscope

\subsubsection{SymMatrix.getDim()}
\label{\detokenize{javaapi/SymMatrix:symmatrix-getdim}}\label{\detokenize{javaapi/SymMatrix::doc}}\begin{quote}

\sphinxAtStartPar
Get the dimension of a symmetric matrix.

\sphinxAtStartPar
\sphinxstylestrong{Synopsis}
\begin{quote}

\sphinxAtStartPar
\sphinxcode{\sphinxupquote{int getDim()}}
\end{quote}

\sphinxAtStartPar
\sphinxstylestrong{Return}
\begin{quote}

\sphinxAtStartPar
Dimension of a symmetric matrix.
\end{quote}
\end{quote}

\subsubsection{SymMatrix.getIdx()}
\label{\detokenize{javaapi/SymMatrix:symmatrix-getidx}}\begin{quote}

\sphinxAtStartPar
Get the index of a symmetric matrix.

\sphinxAtStartPar
\sphinxstylestrong{Synopsis}
\begin{quote}

\sphinxAtStartPar
\sphinxcode{\sphinxupquote{int getIdx()}}
\end{quote}

\sphinxAtStartPar
\sphinxstylestrong{Return}
\begin{quote}

\sphinxAtStartPar
Index of a symmetric matrix.
\end{quote}
\end{quote}

\subsection{SymMatrixArray}
\label{\detokenize{javaapiref:symmatrixarray}}\label{\detokenize{javaapiref:chapjavaapiref-symmatrixarray}}
\sphinxAtStartPar
COPT symmetric matrix object. To store and access a set of {\hyperref[\detokenize{javaapiref:chapjavaapiref-symmatrix}]{\sphinxcrossref{\DUrole{std,std-ref}{SymMatrix}}}}
objects, Cardinal Optimizer provides SymMatrixArray class, which defines the following methods.

\sphinxstepscope

\subsubsection{SymMatrixArray.SymMatrixArray()}
\label{\detokenize{javaapi/SymMatrixArray:symmatrixarray-symmatrixarray}}\label{\detokenize{javaapi/SymMatrixArray::doc}}\begin{quote}

\sphinxAtStartPar
Constructor of SymMatrixAarray.

\sphinxAtStartPar
\sphinxstylestrong{Synopsis}
\begin{quote}

\sphinxAtStartPar
\sphinxcode{\sphinxupquote{SymMatrixArray()}}
\end{quote}
\end{quote}

\subsubsection{SymMatrixArray.getMatrix()}
\label{\detokenize{javaapi/SymMatrixArray:symmatrixarray-getmatrix}}\begin{quote}

\sphinxAtStartPar
Get i\sphinxhyphen{}th SymMatrix object.

\sphinxAtStartPar
\sphinxstylestrong{Synopsis}
\begin{quote}

\sphinxAtStartPar
\sphinxcode{\sphinxupquote{SymMatrix getMatrix(int idx)}}
\end{quote}

\sphinxAtStartPar
\sphinxstylestrong{Arguments}
\begin{quote}

\sphinxAtStartPar
\sphinxcode{\sphinxupquote{idx}}: index of the SymMatrix object.
\end{quote}

\sphinxAtStartPar
\sphinxstylestrong{Return}
\begin{quote}

\sphinxAtStartPar
SymMatrix object with index idx.
\end{quote}
\end{quote}

\subsubsection{SymMatrixArray.pushBack()}
\label{\detokenize{javaapi/SymMatrixArray:symmatrixarray-pushback}}\begin{quote}

\sphinxAtStartPar
Add a SymMatrix object to SymMatrix array.

\sphinxAtStartPar
\sphinxstylestrong{Synopsis}
\begin{quote}

\sphinxAtStartPar
\sphinxcode{\sphinxupquote{void pushBack(SymMatrix mat)}}
\end{quote}

\sphinxAtStartPar
\sphinxstylestrong{Arguments}
\begin{quote}

\sphinxAtStartPar
\sphinxcode{\sphinxupquote{mat}}: a SymMatrix object.
\end{quote}
\end{quote}

\subsubsection{SymMatrixArray.reserve()}
\label{\detokenize{javaapi/SymMatrixArray:symmatrixarray-reserve}}\begin{quote}

\sphinxAtStartPar
Reserve capacity to contain at least n items.

\sphinxAtStartPar
\sphinxstylestrong{Synopsis}
\begin{quote}

\sphinxAtStartPar
\sphinxcode{\sphinxupquote{void reserve(int n)}}
\end{quote}

\sphinxAtStartPar
\sphinxstylestrong{Arguments}
\begin{quote}

\sphinxAtStartPar
\sphinxcode{\sphinxupquote{n}}: minimum capacity for symmetric matrix object.
\end{quote}
\end{quote}

\subsubsection{SymMatrixArray.size()}
\label{\detokenize{javaapi/SymMatrixArray:symmatrixarray-size}}\begin{quote}

\sphinxAtStartPar
Get the number of SymMatrix objects.

\sphinxAtStartPar
\sphinxstylestrong{Synopsis}
\begin{quote}

\sphinxAtStartPar
\sphinxcode{\sphinxupquote{int size()}}
\end{quote}

\sphinxAtStartPar
\sphinxstylestrong{Return}
\begin{quote}

\sphinxAtStartPar
number of SymMatrix objects.
\end{quote}
\end{quote}

\subsection{SymMatExpr}
\label{\detokenize{javaapiref:symmatexpr}}\label{\detokenize{javaapiref:chapjavaapiref-symmatexpr}}
\sphinxAtStartPar
COPT symmetric matrix expression object. A symmetric matrix expression is a
linear combination of symmetric matrices, which is still a symmetric matrix.
However, by doing so, we are able to delay computing the final matrix
until setting PSD constraints or PSD objective.

\sphinxstepscope

\subsubsection{SymMatExpr.SymMatExpr()}
\label{\detokenize{javaapi/SymMatExpr:symmatexpr-symmatexpr}}\label{\detokenize{javaapi/SymMatExpr::doc}}\begin{quote}

\sphinxAtStartPar
Constructor of a symmetric matrix expression.

\sphinxAtStartPar
\sphinxstylestrong{Synopsis}
\begin{quote}

\sphinxAtStartPar
\sphinxcode{\sphinxupquote{SymMatExpr()}}
\end{quote}
\end{quote}

\subsubsection{SymMatExpr.SymMatExpr()}
\label{\detokenize{javaapi/SymMatExpr:id1}}\begin{quote}

\sphinxAtStartPar
Constructor of a symmetric matrix expression with one term.

\sphinxAtStartPar
\sphinxstylestrong{Synopsis}
\begin{quote}

\sphinxAtStartPar
\sphinxcode{\sphinxupquote{SymMatExpr(SymMatrix mat, double coeff)}}
\end{quote}

\sphinxAtStartPar
\sphinxstylestrong{Arguments}
\begin{quote}

\sphinxAtStartPar
\sphinxcode{\sphinxupquote{mat}}: symmetric matrix of the added term.

\sphinxAtStartPar
\sphinxcode{\sphinxupquote{coeff}}: coefficent for the added term.
\end{quote}
\end{quote}

\subsubsection{SymMatExpr.addSymMatExpr()}
\label{\detokenize{javaapi/SymMatExpr:symmatexpr-addsymmatexpr}}\begin{quote}

\sphinxAtStartPar
Add a symmetric matrix expression to self.

\sphinxAtStartPar
\sphinxstylestrong{Synopsis}
\begin{quote}

\sphinxAtStartPar
\sphinxcode{\sphinxupquote{void addSymMatExpr(SymMatExpr expr, double mult)}}
\end{quote}

\sphinxAtStartPar
\sphinxstylestrong{Arguments}
\begin{quote}

\sphinxAtStartPar
\sphinxcode{\sphinxupquote{expr}}: symmetric matrix expression to be added.

\sphinxAtStartPar
\sphinxcode{\sphinxupquote{mult}}: constant multiplier.
\end{quote}
\end{quote}

\subsubsection{SymMatExpr.addTerm()}
\label{\detokenize{javaapi/SymMatExpr:symmatexpr-addterm}}\begin{quote}

\sphinxAtStartPar
Add a term to symmetric matrix expression object.

\sphinxAtStartPar
\sphinxstylestrong{Synopsis}
\begin{quote}

\sphinxAtStartPar
\sphinxcode{\sphinxupquote{Boolean addTerm(SymMatrix mat, double coeff)}}
\end{quote}

\sphinxAtStartPar
\sphinxstylestrong{Arguments}
\begin{quote}

\sphinxAtStartPar
\sphinxcode{\sphinxupquote{mat}}: symmetric matrix of the new term.

\sphinxAtStartPar
\sphinxcode{\sphinxupquote{coeff}}: coefficient of the new term.
\end{quote}

\sphinxAtStartPar
\sphinxstylestrong{Return}
\begin{quote}

\sphinxAtStartPar
True if the term is added successfully.
\end{quote}
\end{quote}

\subsubsection{SymMatExpr.addTerms()}
\label{\detokenize{javaapi/SymMatExpr:symmatexpr-addterms}}\begin{quote}

\sphinxAtStartPar
Add multiple terms to expression object.

\sphinxAtStartPar
\sphinxstylestrong{Synopsis}
\begin{quote}

\sphinxAtStartPar
\sphinxcode{\sphinxupquote{int addTerms(SymMatrix{[}{]} mats, double coeff)}}
\end{quote}

\sphinxAtStartPar
\sphinxstylestrong{Arguments}
\begin{quote}

\sphinxAtStartPar
\sphinxcode{\sphinxupquote{mats}}: symmetric matrix array object for added terms.

\sphinxAtStartPar
\sphinxcode{\sphinxupquote{coeff}}: common coefficient for added terms.
\end{quote}

\sphinxAtStartPar
\sphinxstylestrong{Return}
\begin{quote}

\sphinxAtStartPar
Number of added terms. If negative, fail to add one of terms.
\end{quote}
\end{quote}

\subsubsection{SymMatExpr.addTerms()}
\label{\detokenize{javaapi/SymMatExpr:id2}}\begin{quote}

\sphinxAtStartPar
Add multiple terms to expression object.

\sphinxAtStartPar
\sphinxstylestrong{Synopsis}
\begin{quote}

\sphinxAtStartPar
\sphinxcode{\sphinxupquote{int addTerms(SymMatrixArray mats, double{[}{]} coeffs)}}
\end{quote}

\sphinxAtStartPar
\sphinxstylestrong{Arguments}
\begin{quote}

\sphinxAtStartPar
\sphinxcode{\sphinxupquote{mats}}: symmetric matrix array object for added terms.

\sphinxAtStartPar
\sphinxcode{\sphinxupquote{coeffs}}: coefficient array for added terms.
\end{quote}

\sphinxAtStartPar
\sphinxstylestrong{Return}
\begin{quote}

\sphinxAtStartPar
Number of added terms. If negative, fail to add one of terms.
\end{quote}
\end{quote}

\subsubsection{SymMatExpr.addTerms()}
\label{\detokenize{javaapi/SymMatExpr:id3}}\begin{quote}

\sphinxAtStartPar
Add multiple terms to expression object.

\sphinxAtStartPar
\sphinxstylestrong{Synopsis}
\begin{quote}

\sphinxAtStartPar
\sphinxcode{\sphinxupquote{int addTerms(SymMatrix{[}{]} mats, double{[}{]} coeffs)}}
\end{quote}

\sphinxAtStartPar
\sphinxstylestrong{Arguments}
\begin{quote}

\sphinxAtStartPar
\sphinxcode{\sphinxupquote{mats}}: symmetric matrix array object for added terms.

\sphinxAtStartPar
\sphinxcode{\sphinxupquote{coeffs}}: coefficient array for added terms.
\end{quote}

\sphinxAtStartPar
\sphinxstylestrong{Return}
\begin{quote}

\sphinxAtStartPar
Number of added terms. If negative, fail to add one of terms.
\end{quote}
\end{quote}

\subsubsection{SymMatExpr.clone()}
\label{\detokenize{javaapi/SymMatExpr:symmatexpr-clone}}\begin{quote}

\sphinxAtStartPar
Deep copy symmetric matrix expression object.

\sphinxAtStartPar
\sphinxstylestrong{Synopsis}
\begin{quote}

\sphinxAtStartPar
\sphinxcode{\sphinxupquote{SymMatExpr clone()}}
\end{quote}

\sphinxAtStartPar
\sphinxstylestrong{Return}
\begin{quote}

\sphinxAtStartPar
cloned expression object.
\end{quote}
\end{quote}

\subsubsection{SymMatExpr.getCoeff()}
\label{\detokenize{javaapi/SymMatExpr:symmatexpr-getcoeff}}\begin{quote}

\sphinxAtStartPar
Get coefficient of the i\sphinxhyphen{}th term in expression object.

\sphinxAtStartPar
\sphinxstylestrong{Synopsis}
\begin{quote}

\sphinxAtStartPar
\sphinxcode{\sphinxupquote{double getCoeff(int i)}}
\end{quote}

\sphinxAtStartPar
\sphinxstylestrong{Arguments}
\begin{quote}

\sphinxAtStartPar
\sphinxcode{\sphinxupquote{i}}: index of the term.
\end{quote}

\sphinxAtStartPar
\sphinxstylestrong{Return}
\begin{quote}

\sphinxAtStartPar
coefficient of the i\sphinxhyphen{}th term.
\end{quote}
\end{quote}

\subsubsection{SymMatExpr.getDim()}
\label{\detokenize{javaapi/SymMatExpr:symmatexpr-getdim}}\begin{quote}

\sphinxAtStartPar
Get dimension of symmetric matrix in expression.

\sphinxAtStartPar
\sphinxstylestrong{Synopsis}
\begin{quote}

\sphinxAtStartPar
\sphinxcode{\sphinxupquote{int getDim()}}
\end{quote}

\sphinxAtStartPar
\sphinxstylestrong{Return}
\begin{quote}

\sphinxAtStartPar
dimension of symmetric matrix.
\end{quote}
\end{quote}

\subsubsection{SymMatExpr.getSymMat()}
\label{\detokenize{javaapi/SymMatExpr:symmatexpr-getsymmat}}\begin{quote}

\sphinxAtStartPar
Get symmetric matrix of the i\sphinxhyphen{}th term in expression object.

\sphinxAtStartPar
\sphinxstylestrong{Synopsis}
\begin{quote}

\sphinxAtStartPar
\sphinxcode{\sphinxupquote{SymMatrix getSymMat(int i)}}
\end{quote}

\sphinxAtStartPar
\sphinxstylestrong{Arguments}
\begin{quote}

\sphinxAtStartPar
\sphinxcode{\sphinxupquote{i}}: index of the term.
\end{quote}

\sphinxAtStartPar
\sphinxstylestrong{Return}
\begin{quote}

\sphinxAtStartPar
the symmetric matrix of the i\sphinxhyphen{}th term.
\end{quote}
\end{quote}

\subsubsection{SymMatExpr.multiply()}
\label{\detokenize{javaapi/SymMatExpr:symmatexpr-multiply}}\begin{quote}

\sphinxAtStartPar
Multiply itself by a constant.

\sphinxAtStartPar
\sphinxstylestrong{Synopsis}
\begin{quote}

\sphinxAtStartPar
\sphinxcode{\sphinxupquote{SymMatExpr multiply(double c)}}
\end{quote}

\sphinxAtStartPar
\sphinxstylestrong{Arguments}
\begin{quote}

\sphinxAtStartPar
\sphinxcode{\sphinxupquote{c}}: constant operand.
\end{quote}

\sphinxAtStartPar
\sphinxstylestrong{Return}
\begin{quote}

\sphinxAtStartPar
symmetric matrix expression itself.
\end{quote}
\end{quote}

\subsubsection{SymMatExpr.remove()}
\label{\detokenize{javaapi/SymMatExpr:symmatexpr-remove}}\begin{quote}

\sphinxAtStartPar
Remove i\sphinxhyphen{}th term from expression object.

\sphinxAtStartPar
\sphinxstylestrong{Synopsis}
\begin{quote}

\sphinxAtStartPar
\sphinxcode{\sphinxupquote{void remove(int idx)}}
\end{quote}

\sphinxAtStartPar
\sphinxstylestrong{Arguments}
\begin{quote}

\sphinxAtStartPar
\sphinxcode{\sphinxupquote{idx}}: index of the term to be removed.
\end{quote}
\end{quote}

\subsubsection{SymMatExpr.remove()}
\label{\detokenize{javaapi/SymMatExpr:id4}}\begin{quote}

\sphinxAtStartPar
Remove the term associated with the symmetric matrix.

\sphinxAtStartPar
\sphinxstylestrong{Synopsis}
\begin{quote}

\sphinxAtStartPar
\sphinxcode{\sphinxupquote{void remove(SymMatrix mat)}}
\end{quote}

\sphinxAtStartPar
\sphinxstylestrong{Arguments}
\begin{quote}

\sphinxAtStartPar
\sphinxcode{\sphinxupquote{mat}}: a symmetric matrix whose term should be removed.
\end{quote}
\end{quote}

\subsubsection{SymMatExpr.reserve()}
\label{\detokenize{javaapi/SymMatExpr:symmatexpr-reserve}}\begin{quote}

\sphinxAtStartPar
Reserve capacity to contain at least n items.

\sphinxAtStartPar
\sphinxstylestrong{Synopsis}
\begin{quote}

\sphinxAtStartPar
\sphinxcode{\sphinxupquote{void reserve(int n)}}
\end{quote}

\sphinxAtStartPar
\sphinxstylestrong{Arguments}
\begin{quote}

\sphinxAtStartPar
\sphinxcode{\sphinxupquote{n}}: minimum capacity for expression object.
\end{quote}
\end{quote}

\subsubsection{SymMatExpr.setCoeff()}
\label{\detokenize{javaapi/SymMatExpr:symmatexpr-setcoeff}}\begin{quote}

\sphinxAtStartPar
Set coefficient for the i\sphinxhyphen{}th term in expression object.

\sphinxAtStartPar
\sphinxstylestrong{Synopsis}
\begin{quote}

\sphinxAtStartPar
\sphinxcode{\sphinxupquote{void setCoeff(int i, double val)}}
\end{quote}

\sphinxAtStartPar
\sphinxstylestrong{Arguments}
\begin{quote}

\sphinxAtStartPar
\sphinxcode{\sphinxupquote{i}}: index of the term.

\sphinxAtStartPar
\sphinxcode{\sphinxupquote{val}}: coefficient of the term.
\end{quote}
\end{quote}

\subsubsection{SymMatExpr.size()}
\label{\detokenize{javaapi/SymMatExpr:symmatexpr-size}}\begin{quote}

\sphinxAtStartPar
Get number of terms in expression.

\sphinxAtStartPar
\sphinxstylestrong{Synopsis}
\begin{quote}

\sphinxAtStartPar
\sphinxcode{\sphinxupquote{long size()}}
\end{quote}

\sphinxAtStartPar
\sphinxstylestrong{Return}
\begin{quote}

\sphinxAtStartPar
number of terms.
\end{quote}
\end{quote}

\subsection{NlExpr Class}
\label{\detokenize{javaapiref:nlexpr-class}}\label{\detokenize{javaapiref:chapjavaapiref-nlexpr}}
\sphinxAtStartPar
COPT nonlinear expression object.
The \sphinxcode{\sphinxupquote{NlExpr}} class represents nonlinear expressions in COPT. The nonlinear expressions
are used to build nonlinear constraints. The following methods are provided:

\sphinxstepscope

\subsubsection{NlExpr.NlExpr()}
\label{\detokenize{javaapi/NlExpr:nlexpr-nlexpr}}\label{\detokenize{javaapi/NlExpr::doc}}\begin{quote}

\sphinxAtStartPar
Default constructor of a nonlinear expression.

\sphinxAtStartPar
\sphinxstylestrong{Synopsis}
\begin{quote}

\sphinxAtStartPar
\sphinxcode{\sphinxupquote{NlExpr()}}
\end{quote}
\end{quote}

\subsubsection{NlExpr.NlExpr()}
\label{\detokenize{javaapi/NlExpr:id1}}\begin{quote}

\sphinxAtStartPar
Constructor of a nonlinear expression with a constant.

\sphinxAtStartPar
\sphinxstylestrong{Synopsis}
\begin{quote}

\sphinxAtStartPar
\sphinxcode{\sphinxupquote{NlExpr(double constant)}}
\end{quote}

\sphinxAtStartPar
\sphinxstylestrong{Arguments}
\begin{quote}

\sphinxAtStartPar
\sphinxcode{\sphinxupquote{constant}}: constant value in nonlinear expression object.
\end{quote}
\end{quote}

\subsubsection{NlExpr.NlExpr()}
\label{\detokenize{javaapi/NlExpr:id2}}\begin{quote}

\sphinxAtStartPar
Constructor of a nonlinear expression with a variable.

\sphinxAtStartPar
\sphinxstylestrong{Synopsis}
\begin{quote}

\sphinxAtStartPar
\sphinxcode{\sphinxupquote{NlExpr(Var var)}}
\end{quote}

\sphinxAtStartPar
\sphinxstylestrong{Arguments}
\begin{quote}

\sphinxAtStartPar
\sphinxcode{\sphinxupquote{var}}: the added variable.
\end{quote}
\end{quote}

\subsubsection{NlExpr.NlExpr()}
\label{\detokenize{javaapi/NlExpr:id3}}\begin{quote}

\sphinxAtStartPar
Constructor of a nonlinear expression with one linear term.

\sphinxAtStartPar
\sphinxstylestrong{Synopsis}
\begin{quote}

\sphinxAtStartPar
\sphinxcode{\sphinxupquote{NlExpr(Var var, double coeff)}}
\end{quote}

\sphinxAtStartPar
\sphinxstylestrong{Arguments}
\begin{quote}

\sphinxAtStartPar
\sphinxcode{\sphinxupquote{var}}: variable for the added term.

\sphinxAtStartPar
\sphinxcode{\sphinxupquote{coeff}}: coefficent for the added term.
\end{quote}
\end{quote}

\subsubsection{NlExpr.NlExpr()}
\label{\detokenize{javaapi/NlExpr:id4}}\begin{quote}

\sphinxAtStartPar
Constructor of a nonlinear expression with a linear expression.

\sphinxAtStartPar
\sphinxstylestrong{Synopsis}
\begin{quote}

\sphinxAtStartPar
\sphinxcode{\sphinxupquote{NlExpr(Expr expr)}}
\end{quote}

\sphinxAtStartPar
\sphinxstylestrong{Arguments}
\begin{quote}

\sphinxAtStartPar
\sphinxcode{\sphinxupquote{expr}}: the added linear expression.
\end{quote}
\end{quote}

\subsubsection{NlExpr.NlExpr()}
\label{\detokenize{javaapi/NlExpr:id5}}\begin{quote}

\sphinxAtStartPar
Constructor of a nonlinear expression with a quadratic expression.

\sphinxAtStartPar
\sphinxstylestrong{Synopsis}
\begin{quote}

\sphinxAtStartPar
\sphinxcode{\sphinxupquote{NlExpr(QuadExpr expr)}}
\end{quote}

\sphinxAtStartPar
\sphinxstylestrong{Arguments}
\begin{quote}

\sphinxAtStartPar
\sphinxcode{\sphinxupquote{expr}}: the added quadratic expression.
\end{quote}
\end{quote}

\subsubsection{NlExpr.add()}
\label{\detokenize{javaapi/NlExpr:nlexpr-add}}\begin{quote}

\sphinxAtStartPar
Add itself by an expression.

\sphinxAtStartPar
\sphinxstylestrong{Synopsis}
\begin{quote}

\sphinxAtStartPar
\sphinxcode{\sphinxupquote{NlExpr add(NlExpr expr, double mult)}}
\end{quote}

\sphinxAtStartPar
\sphinxstylestrong{Arguments}
\begin{quote}

\sphinxAtStartPar
\sphinxcode{\sphinxupquote{expr}}: expression operand, including NlExpr, QuadExpr, Expr, Var.

\sphinxAtStartPar
\sphinxcode{\sphinxupquote{mult}}: constant multiplier.
\end{quote}

\sphinxAtStartPar
\sphinxstylestrong{Return}
\begin{quote}

\sphinxAtStartPar
nonlinear expression itself.
\end{quote}
\end{quote}

\subsubsection{NlExpr.add()}
\label{\detokenize{javaapi/NlExpr:id6}}\begin{quote}

\sphinxAtStartPar
Add itself by an expression.

\sphinxAtStartPar
\sphinxstylestrong{Synopsis}
\begin{quote}

\sphinxAtStartPar
\sphinxcode{\sphinxupquote{NlExpr add(NlExpr expr)}}
\end{quote}

\sphinxAtStartPar
\sphinxstylestrong{Arguments}
\begin{quote}

\sphinxAtStartPar
\sphinxcode{\sphinxupquote{expr}}: expression operand, including NlExpr, QuadExpr, Expr, Var and constant.
\end{quote}

\sphinxAtStartPar
\sphinxstylestrong{Return}
\begin{quote}

\sphinxAtStartPar
nonlinear expression itself.
\end{quote}
\end{quote}

\subsubsection{NlExpr.addConstant()}
\label{\detokenize{javaapi/NlExpr:nlexpr-addconstant}}\begin{quote}

\sphinxAtStartPar
Add constant to the nonlinear expression.

\sphinxAtStartPar
\sphinxstylestrong{Synopsis}
\begin{quote}

\sphinxAtStartPar
\sphinxcode{\sphinxupquote{void addConstant(double constant)}}
\end{quote}

\sphinxAtStartPar
\sphinxstylestrong{Arguments}
\begin{quote}

\sphinxAtStartPar
\sphinxcode{\sphinxupquote{constant}}: value to be added.
\end{quote}
\end{quote}

\subsubsection{NlExpr.addLinExpr()}
\label{\detokenize{javaapi/NlExpr:nlexpr-addlinexpr}}\begin{quote}

\sphinxAtStartPar
Add a linear expression to self.

\sphinxAtStartPar
\sphinxstylestrong{Synopsis}
\begin{quote}

\sphinxAtStartPar
\sphinxcode{\sphinxupquote{void addLinExpr(Expr expr, double mult)}}
\end{quote}

\sphinxAtStartPar
\sphinxstylestrong{Arguments}
\begin{quote}

\sphinxAtStartPar
\sphinxcode{\sphinxupquote{expr}}: linear expression to be added.

\sphinxAtStartPar
\sphinxcode{\sphinxupquote{mult}}: constant multiplier.
\end{quote}
\end{quote}

\subsubsection{NlExpr.addNlExpr()}
\label{\detokenize{javaapi/NlExpr:nlexpr-addnlexpr}}\begin{quote}

\sphinxAtStartPar
Add a nonlinear expression to self.

\sphinxAtStartPar
\sphinxstylestrong{Synopsis}
\begin{quote}

\sphinxAtStartPar
\sphinxcode{\sphinxupquote{void addNlExpr(NlExpr expr, double mult)}}
\end{quote}

\sphinxAtStartPar
\sphinxstylestrong{Arguments}
\begin{quote}

\sphinxAtStartPar
\sphinxcode{\sphinxupquote{expr}}: nonlinear expression to be added.

\sphinxAtStartPar
\sphinxcode{\sphinxupquote{mult}}: constant multiplier.
\end{quote}
\end{quote}

\subsubsection{NlExpr.addQuadExpr()}
\label{\detokenize{javaapi/NlExpr:nlexpr-addquadexpr}}\begin{quote}

\sphinxAtStartPar
Add a quadratic expression to self.

\sphinxAtStartPar
\sphinxstylestrong{Synopsis}
\begin{quote}

\sphinxAtStartPar
\sphinxcode{\sphinxupquote{void addQuadExpr(QuadExpr expr, double mult)}}
\end{quote}

\sphinxAtStartPar
\sphinxstylestrong{Arguments}
\begin{quote}

\sphinxAtStartPar
\sphinxcode{\sphinxupquote{expr}}: quadratic expression to be added.

\sphinxAtStartPar
\sphinxcode{\sphinxupquote{mult}}: constant multiplier.
\end{quote}
\end{quote}

\subsubsection{NlExpr.addTerm()}
\label{\detokenize{javaapi/NlExpr:nlexpr-addterm}}\begin{quote}

\sphinxAtStartPar
Add a linear term to nonlinear expression object.

\sphinxAtStartPar
\sphinxstylestrong{Synopsis}
\begin{quote}

\sphinxAtStartPar
\sphinxcode{\sphinxupquote{void addTerm(Var var, double coeff)}}
\end{quote}

\sphinxAtStartPar
\sphinxstylestrong{Arguments}
\begin{quote}

\sphinxAtStartPar
\sphinxcode{\sphinxupquote{var}}: variable of new linear term.

\sphinxAtStartPar
\sphinxcode{\sphinxupquote{coeff}}: coefficient of new linear term.
\end{quote}
\end{quote}

\subsubsection{NlExpr.addTerms()}
\label{\detokenize{javaapi/NlExpr:nlexpr-addterms}}\begin{quote}

\sphinxAtStartPar
Add linear terms to nonlinear expression object.

\sphinxAtStartPar
\sphinxstylestrong{Synopsis}
\begin{quote}

\sphinxAtStartPar
\sphinxcode{\sphinxupquote{void addTerms(Var{[}{]} vars, double{[}{]} coeffs)}}
\end{quote}

\sphinxAtStartPar
\sphinxstylestrong{Arguments}
\begin{quote}

\sphinxAtStartPar
\sphinxcode{\sphinxupquote{vars}}: variable array for added linear terms.

\sphinxAtStartPar
\sphinxcode{\sphinxupquote{coeffs}}: coefficient array for added linear terms.
\end{quote}
\end{quote}

\subsubsection{NlExpr.addTerms()}
\label{\detokenize{javaapi/NlExpr:id7}}\begin{quote}

\sphinxAtStartPar
Add linear terms to nonlinear expression object.

\sphinxAtStartPar
\sphinxstylestrong{Synopsis}
\begin{quote}

\sphinxAtStartPar
\sphinxcode{\sphinxupquote{void addTerms(VarArray vars, double{[}{]} coeffs)}}
\end{quote}

\sphinxAtStartPar
\sphinxstylestrong{Arguments}
\begin{quote}

\sphinxAtStartPar
\sphinxcode{\sphinxupquote{vars}}: variables for added linear terms.

\sphinxAtStartPar
\sphinxcode{\sphinxupquote{coeffs}}: coefficient array for added linear terms.
\end{quote}
\end{quote}

\subsubsection{NlExpr.clear()}
\label{\detokenize{javaapi/NlExpr:nlexpr-clear}}\begin{quote}

\sphinxAtStartPar
Clear nonlinear expression object.

\sphinxAtStartPar
\sphinxstylestrong{Synopsis}
\begin{quote}

\sphinxAtStartPar
\sphinxcode{\sphinxupquote{void clear()}}
\end{quote}
\end{quote}

\subsubsection{NlExpr.clone()}
\label{\detokenize{javaapi/NlExpr:nlexpr-clone}}\begin{quote}

\sphinxAtStartPar
Deep copy nonlinear expression object.

\sphinxAtStartPar
\sphinxstylestrong{Synopsis}
\begin{quote}

\sphinxAtStartPar
\sphinxcode{\sphinxupquote{NlExpr clone()}}
\end{quote}

\sphinxAtStartPar
\sphinxstylestrong{Return}
\begin{quote}

\sphinxAtStartPar
cloned nonlinear expression object.
\end{quote}
\end{quote}

\subsubsection{NlExpr.divide()}
\label{\detokenize{javaapi/NlExpr:nlexpr-divide}}\begin{quote}

\sphinxAtStartPar
Divide itself by an expression.

\sphinxAtStartPar
\sphinxstylestrong{Synopsis}
\begin{quote}

\sphinxAtStartPar
\sphinxcode{\sphinxupquote{NlExpr divide(NlExpr expr)}}
\end{quote}

\sphinxAtStartPar
\sphinxstylestrong{Arguments}
\begin{quote}

\sphinxAtStartPar
\sphinxcode{\sphinxupquote{expr}}: expression operand, including NlExpr, QuadExpr, Expr, Var and constant.
\end{quote}

\sphinxAtStartPar
\sphinxstylestrong{Return}
\begin{quote}

\sphinxAtStartPar
nonlinear expression itself.
\end{quote}
\end{quote}

\subsubsection{NlExpr.evaluate()}
\label{\detokenize{javaapi/NlExpr:nlexpr-evaluate}}\begin{quote}

\sphinxAtStartPar
Evaluate nonlinear expression after solving.

\sphinxAtStartPar
\sphinxstylestrong{Synopsis}
\begin{quote}

\sphinxAtStartPar
\sphinxcode{\sphinxupquote{double evaluate()}}
\end{quote}

\sphinxAtStartPar
\sphinxstylestrong{Return}
\begin{quote}

\sphinxAtStartPar
value of nonlinear expression.
\end{quote}
\end{quote}

\subsubsection{NlExpr.getConstant()}
\label{\detokenize{javaapi/NlExpr:nlexpr-getconstant}}\begin{quote}

\sphinxAtStartPar
Get constant in nonlinear expression.

\sphinxAtStartPar
\sphinxstylestrong{Synopsis}
\begin{quote}

\sphinxAtStartPar
\sphinxcode{\sphinxupquote{double getConstant()}}
\end{quote}

\sphinxAtStartPar
\sphinxstylestrong{Return}
\begin{quote}

\sphinxAtStartPar
constant in nonlinear expression.
\end{quote}
\end{quote}

\subsubsection{NlExpr.getLinExpr()}
\label{\detokenize{javaapi/NlExpr:nlexpr-getlinexpr}}\begin{quote}

\sphinxAtStartPar
Get linear expression of nonlinear expression.

\sphinxAtStartPar
\sphinxstylestrong{Synopsis}
\begin{quote}

\sphinxAtStartPar
\sphinxcode{\sphinxupquote{Expr getLinExpr()}}
\end{quote}

\sphinxAtStartPar
\sphinxstylestrong{Return}
\begin{quote}

\sphinxAtStartPar
linear expression object.
\end{quote}
\end{quote}

\subsubsection{NlExpr.multiply()}
\label{\detokenize{javaapi/NlExpr:nlexpr-multiply}}\begin{quote}

\sphinxAtStartPar
Multiply itself by an expression.

\sphinxAtStartPar
\sphinxstylestrong{Synopsis}
\begin{quote}

\sphinxAtStartPar
\sphinxcode{\sphinxupquote{NlExpr multiply(NlExpr expr)}}
\end{quote}

\sphinxAtStartPar
\sphinxstylestrong{Arguments}
\begin{quote}

\sphinxAtStartPar
\sphinxcode{\sphinxupquote{expr}}: expression operand, including NlExpr, QuadExpr, Expr, Var and constant.
\end{quote}

\sphinxAtStartPar
\sphinxstylestrong{Return}
\begin{quote}

\sphinxAtStartPar
nonlinear expression itself.
\end{quote}
\end{quote}

\subsubsection{NlExpr.negate()}
\label{\detokenize{javaapi/NlExpr:nlexpr-negate}}\begin{quote}

\sphinxAtStartPar
Negate itself.

\sphinxAtStartPar
\sphinxstylestrong{Synopsis}
\begin{quote}

\sphinxAtStartPar
\sphinxcode{\sphinxupquote{NlExpr negate()}}
\end{quote}

\sphinxAtStartPar
\sphinxstylestrong{Return}
\begin{quote}

\sphinxAtStartPar
nonlinear expression itself.
\end{quote}
\end{quote}

\subsubsection{NlExpr.reserve()}
\label{\detokenize{javaapi/NlExpr:nlexpr-reserve}}\begin{quote}

\sphinxAtStartPar
Reserve capacity to contain at least n items.

\sphinxAtStartPar
\sphinxstylestrong{Synopsis}
\begin{quote}

\sphinxAtStartPar
\sphinxcode{\sphinxupquote{void reserve(int n)}}
\end{quote}

\sphinxAtStartPar
\sphinxstylestrong{Arguments}
\begin{quote}

\sphinxAtStartPar
\sphinxcode{\sphinxupquote{n}}: capacity of nonlinear constraint objects.
\end{quote}
\end{quote}

\subsubsection{NlExpr.setConstant()}
\label{\detokenize{javaapi/NlExpr:nlexpr-setconstant}}\begin{quote}

\sphinxAtStartPar
Set constant for the nonlinear expression.

\sphinxAtStartPar
\sphinxstylestrong{Synopsis}
\begin{quote}

\sphinxAtStartPar
\sphinxcode{\sphinxupquote{void setConstant(double constant)}}
\end{quote}

\sphinxAtStartPar
\sphinxstylestrong{Arguments}
\begin{quote}

\sphinxAtStartPar
\sphinxcode{\sphinxupquote{constant}}: the value of the constant.
\end{quote}
\end{quote}

\subsubsection{NlExpr.size()}
\label{\detokenize{javaapi/NlExpr:nlexpr-size}}\begin{quote}

\sphinxAtStartPar
Get size of tokens in nonlinear expression.

\sphinxAtStartPar
\sphinxstylestrong{Synopsis}
\begin{quote}

\sphinxAtStartPar
\sphinxcode{\sphinxupquote{long size()}}
\end{quote}

\sphinxAtStartPar
\sphinxstylestrong{Return}
\begin{quote}

\sphinxAtStartPar
size of none\sphinxhyphen{}linear tokens.
\end{quote}
\end{quote}

\subsection{NlConstraint Class}
\label{\detokenize{javaapiref:nlconstraint-class}}\label{\detokenize{javaapiref:chapjavaapiref-nlconstraint}}
\sphinxAtStartPar
COPT nonlinear constraint object. The \sphinxcode{\sphinxupquote{NlConstraint}} object is always associated with a
particular model. User creates a \sphinxcode{\sphinxupquote{NlConstraint}} object by adding a nonlinear
constraint to model, rather than by constructor of \sphinxcode{\sphinxupquote{NlConstraint}} class.

\sphinxstepscope

\subsubsection{NlConstraint.get()}
\label{\detokenize{javaapi/NlConstraint:nlconstraint-get}}\label{\detokenize{javaapi/NlConstraint::doc}}\begin{quote}

\sphinxAtStartPar
Get information value of the nonlinear constraint. Support informations of “LB”, “UB”, “Slack”.

\sphinxAtStartPar
\sphinxstylestrong{Synopsis}
\begin{quote}

\sphinxAtStartPar
\sphinxcode{\sphinxupquote{double get(String info)}}
\end{quote}

\sphinxAtStartPar
\sphinxstylestrong{Arguments}
\begin{quote}

\sphinxAtStartPar
\sphinxcode{\sphinxupquote{info}}: name of the information being queried.
\end{quote}

\sphinxAtStartPar
\sphinxstylestrong{Return}
\begin{quote}

\sphinxAtStartPar
value of information.
\end{quote}
\end{quote}

\subsubsection{NlConstraint.getIdx()}
\label{\detokenize{javaapi/NlConstraint:nlconstraint-getidx}}\begin{quote}

\sphinxAtStartPar
Get index of nonlinear constraint.

\sphinxAtStartPar
\sphinxstylestrong{Synopsis}
\begin{quote}

\sphinxAtStartPar
\sphinxcode{\sphinxupquote{int getIdx()}}
\end{quote}

\sphinxAtStartPar
\sphinxstylestrong{Return}
\begin{quote}

\sphinxAtStartPar
the index of nonlinear constraint.
\end{quote}
\end{quote}

\subsubsection{NlConstraint.getName()}
\label{\detokenize{javaapi/NlConstraint:nlconstraint-getname}}\begin{quote}

\sphinxAtStartPar
Get name of nonlinear constraint.

\sphinxAtStartPar
\sphinxstylestrong{Synopsis}
\begin{quote}

\sphinxAtStartPar
\sphinxcode{\sphinxupquote{String getName()}}
\end{quote}

\sphinxAtStartPar
\sphinxstylestrong{Return}
\begin{quote}

\sphinxAtStartPar
the name of nonlinear constraint.
\end{quote}
\end{quote}

\subsubsection{NlConstraint.remove()}
\label{\detokenize{javaapi/NlConstraint:nlconstraint-remove}}\begin{quote}

\sphinxAtStartPar
Remove this nonlinear constraint from model.

\sphinxAtStartPar
\sphinxstylestrong{Synopsis}
\begin{quote}

\sphinxAtStartPar
\sphinxcode{\sphinxupquote{void remove()}}
\end{quote}
\end{quote}

\subsubsection{NlConstraint.set()}
\label{\detokenize{javaapi/NlConstraint:nlconstraint-set}}\begin{quote}

\sphinxAtStartPar
Set information value of nonlinear constraint. Support informations of “LB” and “UB”.

\sphinxAtStartPar
\sphinxstylestrong{Synopsis}
\begin{quote}

\sphinxAtStartPar
\sphinxcode{\sphinxupquote{void set(String info, double val)}}
\end{quote}

\sphinxAtStartPar
\sphinxstylestrong{Arguments}
\begin{quote}

\sphinxAtStartPar
\sphinxcode{\sphinxupquote{info}}: name of the information.

\sphinxAtStartPar
\sphinxcode{\sphinxupquote{val}}: new information value.
\end{quote}
\end{quote}

\subsubsection{NlConstraint.setName()}
\label{\detokenize{javaapi/NlConstraint:nlconstraint-setname}}\begin{quote}

\sphinxAtStartPar
Set name for nonlinear constraint.

\sphinxAtStartPar
\sphinxstylestrong{Synopsis}
\begin{quote}

\sphinxAtStartPar
\sphinxcode{\sphinxupquote{void setName(String name)}}
\end{quote}

\sphinxAtStartPar
\sphinxstylestrong{Arguments}
\begin{quote}

\sphinxAtStartPar
\sphinxcode{\sphinxupquote{name}}: the name to set.
\end{quote}
\end{quote}

\subsection{NlConstrArray Class}
\label{\detokenize{javaapiref:nlconstrarray-class}}\label{\detokenize{javaapiref:chapjavaapiref-nlconstrarray}}
\sphinxAtStartPar
COPT nonlinear constraint array object. To store and access a set of
{\hyperref[\detokenize{javaapiref:chapjavaapiref-nlconstraint}]{\sphinxcrossref{\DUrole{std,std-ref}{NlConstraint Class}}}} objects, Cardinal Optimizer provides
NlConstrArray class, which defines the following methods.

\sphinxstepscope

\subsubsection{NlConstrArray.NlConstrArray()}
\label{\detokenize{javaapi/NlConstrArray:nlconstrarray-nlconstrarray}}\label{\detokenize{javaapi/NlConstrArray::doc}}\begin{quote}

\sphinxAtStartPar
Constructor of NlConstrArray object.

\sphinxAtStartPar
\sphinxstylestrong{Synopsis}
\begin{quote}

\sphinxAtStartPar
\sphinxcode{\sphinxupquote{NlConstrArray()}}
\end{quote}
\end{quote}

\subsubsection{NlConstrArray.getNlConstr()}
\label{\detokenize{javaapi/NlConstrArray:nlconstrarray-getnlconstr}}\begin{quote}

\sphinxAtStartPar
Get idx\sphinxhyphen{}th nonlinear constraint object.

\sphinxAtStartPar
\sphinxstylestrong{Synopsis}
\begin{quote}

\sphinxAtStartPar
\sphinxcode{\sphinxupquote{NlConstraint getNlConstr(int idx)}}
\end{quote}

\sphinxAtStartPar
\sphinxstylestrong{Arguments}
\begin{quote}

\sphinxAtStartPar
\sphinxcode{\sphinxupquote{idx}}: index of the nonlinear constraint.
\end{quote}

\sphinxAtStartPar
\sphinxstylestrong{Return}
\begin{quote}

\sphinxAtStartPar
nonlinear constraint object with index value.
\end{quote}
\end{quote}

\subsubsection{NlConstrArray.pushBack()}
\label{\detokenize{javaapi/NlConstrArray:nlconstrarray-pushback}}\begin{quote}

\sphinxAtStartPar
Add a nonlinear constraint to nonlinear constraint array.

\sphinxAtStartPar
\sphinxstylestrong{Synopsis}
\begin{quote}

\sphinxAtStartPar
\sphinxcode{\sphinxupquote{void pushBack(NlConstraint constr)}}
\end{quote}

\sphinxAtStartPar
\sphinxstylestrong{Arguments}
\begin{quote}

\sphinxAtStartPar
\sphinxcode{\sphinxupquote{constr}}: nonlinear constraint object.
\end{quote}
\end{quote}

\subsubsection{NlConstrArray.reserve()}
\label{\detokenize{javaapi/NlConstrArray:nlconstrarray-reserve}}\begin{quote}

\sphinxAtStartPar
Reserve capacity to contain at least n items.

\sphinxAtStartPar
\sphinxstylestrong{Synopsis}
\begin{quote}

\sphinxAtStartPar
\sphinxcode{\sphinxupquote{void reserve(int n)}}
\end{quote}

\sphinxAtStartPar
\sphinxstylestrong{Arguments}
\begin{quote}

\sphinxAtStartPar
\sphinxcode{\sphinxupquote{n}}: capacity of nonlinear constraint objects.
\end{quote}
\end{quote}

\subsubsection{NlConstrArray.size()}
\label{\detokenize{javaapi/NlConstrArray:nlconstrarray-size}}\begin{quote}

\sphinxAtStartPar
Get the number of nonlinear constraint objects.

\sphinxAtStartPar
\sphinxstylestrong{Synopsis}
\begin{quote}

\sphinxAtStartPar
\sphinxcode{\sphinxupquote{int size()}}
\end{quote}

\sphinxAtStartPar
\sphinxstylestrong{Return}
\begin{quote}

\sphinxAtStartPar
number of nonlinear constraint objects.
\end{quote}
\end{quote}

\subsection{NlConstrBuilder Class}
\label{\detokenize{javaapiref:nlconstrbuilder-class}}\label{\detokenize{javaapiref:chapjavaapiref-nlconstrbuilder}}
\sphinxAtStartPar
COPT nonlinear constraint builder object. To help building a nonlinear constraint, given a
nonlinear expression, constraint sense and right\sphinxhyphen{}hand side value, Cardinal Optimizer provides
NlConstrBuilder class, which defines the following methods.

\sphinxstepscope

\subsubsection{NlConstrBuilder.NlConstrBuilder()}
\label{\detokenize{javaapi/NlConstrBuilder:nlconstrbuilder-nlconstrbuilder}}\label{\detokenize{javaapi/NlConstrBuilder::doc}}\begin{quote}

\sphinxAtStartPar
Constructor of NlConstrBuilder object.

\sphinxAtStartPar
\sphinxstylestrong{Synopsis}
\begin{quote}

\sphinxAtStartPar
\sphinxcode{\sphinxupquote{NlConstrBuilder()}}
\end{quote}
\end{quote}

\subsubsection{NlConstrBuilder.getNlExpr()}
\label{\detokenize{javaapi/NlConstrBuilder:nlconstrbuilder-getnlexpr}}\begin{quote}

\sphinxAtStartPar
Get nonlinear expression associated with constraint.

\sphinxAtStartPar
\sphinxstylestrong{Synopsis}
\begin{quote}

\sphinxAtStartPar
\sphinxcode{\sphinxupquote{NlExpr getNlExpr()}}
\end{quote}

\sphinxAtStartPar
\sphinxstylestrong{Return}
\begin{quote}

\sphinxAtStartPar
nonlinear expression object.
\end{quote}
\end{quote}

\subsubsection{NlConstrBuilder.getRange()}
\label{\detokenize{javaapi/NlConstrBuilder:nlconstrbuilder-getrange}}\begin{quote}

\sphinxAtStartPar
Get range from lower bound to upper bound of range constraint.

\sphinxAtStartPar
\sphinxstylestrong{Synopsis}
\begin{quote}

\sphinxAtStartPar
\sphinxcode{\sphinxupquote{double getRange()}}
\end{quote}

\sphinxAtStartPar
\sphinxstylestrong{Return}
\begin{quote}

\sphinxAtStartPar
length from lower bound to upper bound of nonlinear constraint.
\end{quote}
\end{quote}

\subsubsection{NlConstrBuilder.getSense()}
\label{\detokenize{javaapi/NlConstrBuilder:nlconstrbuilder-getsense}}\begin{quote}

\sphinxAtStartPar
Get sense associated with nonlinear constraint.

\sphinxAtStartPar
\sphinxstylestrong{Synopsis}
\begin{quote}

\sphinxAtStartPar
\sphinxcode{\sphinxupquote{char getSense()}}
\end{quote}

\sphinxAtStartPar
\sphinxstylestrong{Return}
\begin{quote}

\sphinxAtStartPar
nonlinear constraint sense.
\end{quote}
\end{quote}

\subsubsection{NlConstrBuilder.set()}
\label{\detokenize{javaapi/NlConstrBuilder:nlconstrbuilder-set}}\begin{quote}

\sphinxAtStartPar
Set detail of a nonlinear constraint to its builder object.

\sphinxAtStartPar
\sphinxstylestrong{Synopsis}
\begin{quote}

\sphinxAtStartPar
\sphinxcode{\sphinxupquote{void set(}}
\begin{quote}

\sphinxAtStartPar
\sphinxcode{\sphinxupquote{NlExpr expr,}}

\sphinxAtStartPar
\sphinxcode{\sphinxupquote{char sense,}}

\sphinxAtStartPar
\sphinxcode{\sphinxupquote{double rhs)}}
\end{quote}
\end{quote}

\sphinxAtStartPar
\sphinxstylestrong{Arguments}
\begin{quote}

\sphinxAtStartPar
\sphinxcode{\sphinxupquote{expr}}: nonlinear expression object at one side of nonlinear constraint

\sphinxAtStartPar
\sphinxcode{\sphinxupquote{sense}}: constraint sense other than COPT\_RANGE.

\sphinxAtStartPar
\sphinxcode{\sphinxupquote{rhs}}: constant of right side of nonlinear constraint.
\end{quote}
\end{quote}

\subsubsection{NlConstrBuilder.setRange()}
\label{\detokenize{javaapi/NlConstrBuilder:nlconstrbuilder-setrange}}\begin{quote}

\sphinxAtStartPar
Set a range constraint to nonlinear constraint builder.

\sphinxAtStartPar
\sphinxstylestrong{Synopsis}
\begin{quote}

\sphinxAtStartPar
\sphinxcode{\sphinxupquote{void setRange(NlExpr expr, double range)}}
\end{quote}

\sphinxAtStartPar
\sphinxstylestrong{Arguments}
\begin{quote}

\sphinxAtStartPar
\sphinxcode{\sphinxupquote{expr}}: nonlinear expression object, whose constant is negative upper bound.

\sphinxAtStartPar
\sphinxcode{\sphinxupquote{range}}: length from lower bound to upper bound of nonlinear constraint. Must greater than 0.
\end{quote}
\end{quote}

\subsection{NlConstrBuilderArray Class}
\label{\detokenize{javaapiref:nlconstrbuilderarray-class}}\label{\detokenize{javaapiref:chapjavaapiref-nlconstrbuilderarray}}
\sphinxAtStartPar
COPT nonlinear constraint builder array object. To store and access a set of
{\hyperref[\detokenize{javaapiref:chapjavaapiref-nlconstrbuilder}]{\sphinxcrossref{\DUrole{std,std-ref}{NlConstrBuilder Class}}}} objects, Cardinal Optimizer provides
\sphinxcode{\sphinxupquote{NlConstrBuilderArray}} class, which defines the following methods.

\sphinxstepscope

\subsubsection{NlConstrBuilderArray.NlConstrBuilderArray()}
\label{\detokenize{javaapi/NlConstrBuilderArray:nlconstrbuilderarray-nlconstrbuilderarray}}\label{\detokenize{javaapi/NlConstrBuilderArray::doc}}\begin{quote}

\sphinxAtStartPar
Constructor of NlConstrBuilderArray object.

\sphinxAtStartPar
\sphinxstylestrong{Synopsis}
\begin{quote}

\sphinxAtStartPar
\sphinxcode{\sphinxupquote{NlConstrBuilderArray()}}
\end{quote}
\end{quote}

\subsubsection{NlConstrBuilderArray.getBuilder()}
\label{\detokenize{javaapi/NlConstrBuilderArray:nlconstrbuilderarray-getbuilder}}\begin{quote}

\sphinxAtStartPar
Get idx\sphinxhyphen{}th nonlinear constraint builder object.

\sphinxAtStartPar
\sphinxstylestrong{Synopsis}
\begin{quote}

\sphinxAtStartPar
\sphinxcode{\sphinxupquote{NlConstrBuilder getBuilder(int idx)}}
\end{quote}

\sphinxAtStartPar
\sphinxstylestrong{Arguments}
\begin{quote}

\sphinxAtStartPar
\sphinxcode{\sphinxupquote{idx}}: index of the nonlinear constraint builder.
\end{quote}

\sphinxAtStartPar
\sphinxstylestrong{Return}
\begin{quote}

\sphinxAtStartPar
nonlinear constraint builder object with index idx.
\end{quote}
\end{quote}

\subsubsection{NlConstrBuilderArray.pushBack()}
\label{\detokenize{javaapi/NlConstrBuilderArray:nlconstrbuilderarray-pushback}}\begin{quote}

\sphinxAtStartPar
Add a nonlinear constraint builder object to nonlinear constraint builder array.

\sphinxAtStartPar
\sphinxstylestrong{Synopsis}
\begin{quote}

\sphinxAtStartPar
\sphinxcode{\sphinxupquote{void pushBack(NlConstrBuilder builder)}}
\end{quote}

\sphinxAtStartPar
\sphinxstylestrong{Arguments}
\begin{quote}

\sphinxAtStartPar
\sphinxcode{\sphinxupquote{builder}}: a nonlinear constraint builder object.
\end{quote}
\end{quote}

\subsubsection{NlConstrBuilderArray.reserve()}
\label{\detokenize{javaapi/NlConstrBuilderArray:nlconstrbuilderarray-reserve}}\begin{quote}

\sphinxAtStartPar
Reserve capacity to contain at least n items.

\sphinxAtStartPar
\sphinxstylestrong{Synopsis}
\begin{quote}

\sphinxAtStartPar
\sphinxcode{\sphinxupquote{void reserve(int n)}}
\end{quote}

\sphinxAtStartPar
\sphinxstylestrong{Arguments}
\begin{quote}

\sphinxAtStartPar
\sphinxcode{\sphinxupquote{n}}: capacity of nonlinear constraint objects.
\end{quote}
\end{quote}

\subsubsection{NlConstrBuilderArray.size()}
\label{\detokenize{javaapi/NlConstrBuilderArray:nlconstrbuilderarray-size}}\begin{quote}

\sphinxAtStartPar
Get the number of nonlinear constraint builder objects.

\sphinxAtStartPar
\sphinxstylestrong{Synopsis}
\begin{quote}

\sphinxAtStartPar
\sphinxcode{\sphinxupquote{int size()}}
\end{quote}

\sphinxAtStartPar
\sphinxstylestrong{Return}
\begin{quote}

\sphinxAtStartPar
number of nonlinear constraint builder objects.
\end{quote}
\end{quote}

\subsection{NL Namespace}
\label{\detokenize{javaapiref:nl-namespace}}\label{\detokenize{javaapiref:chapjavaapiref-nl}}
\sphinxAtStartPar
Common nonlinear functions in the \sphinxcode{\sphinxupquote{NL}} namespace are provided for
constructing nonlinear expressions.  The following methods are provided:

\sphinxstepscope

\subsubsection{NL.abs()}
\label{\detokenize{javaapi/NL:nl-abs}}\label{\detokenize{javaapi/NL::doc}}\begin{quote}

\sphinxAtStartPar
Calculate absolute value of a nonlinear expression.

\sphinxAtStartPar
\sphinxstylestrong{Synopsis}
\begin{quote}

\sphinxAtStartPar
\sphinxcode{\sphinxupquote{static NlExpr abs(NlExpr expr)}}
\end{quote}

\sphinxAtStartPar
\sphinxstylestrong{Arguments}
\begin{quote}

\sphinxAtStartPar
\sphinxcode{\sphinxupquote{expr}}: a nonlinear expression.
\end{quote}

\sphinxAtStartPar
\sphinxstylestrong{Return}
\begin{quote}

\sphinxAtStartPar
result as a nonlinear expression.
\end{quote}
\end{quote}

\subsubsection{NL.acos()}
\label{\detokenize{javaapi/NL:nl-acos}}\begin{quote}

\sphinxAtStartPar
Calculate arccosine of a nonlinear expression.

\sphinxAtStartPar
\sphinxstylestrong{Synopsis}
\begin{quote}

\sphinxAtStartPar
\sphinxcode{\sphinxupquote{static NlExpr acos(NlExpr expr)}}
\end{quote}

\sphinxAtStartPar
\sphinxstylestrong{Arguments}
\begin{quote}

\sphinxAtStartPar
\sphinxcode{\sphinxupquote{expr}}: a nonlinear expression.
\end{quote}

\sphinxAtStartPar
\sphinxstylestrong{Return}
\begin{quote}

\sphinxAtStartPar
result as a nonlinear expression.
\end{quote}
\end{quote}

\subsubsection{NL.acosh()}
\label{\detokenize{javaapi/NL:nl-acosh}}\begin{quote}

\sphinxAtStartPar
Calculate inverse hyperbolic cosine of a nonlinear expression.

\sphinxAtStartPar
\sphinxstylestrong{Synopsis}
\begin{quote}

\sphinxAtStartPar
\sphinxcode{\sphinxupquote{static NlExpr acosh(NlExpr expr)}}
\end{quote}

\sphinxAtStartPar
\sphinxstylestrong{Arguments}
\begin{quote}

\sphinxAtStartPar
\sphinxcode{\sphinxupquote{expr}}: a nonlinear expression.
\end{quote}

\sphinxAtStartPar
\sphinxstylestrong{Return}
\begin{quote}

\sphinxAtStartPar
result as a nonlinear expression.
\end{quote}
\end{quote}

\subsubsection{NL.asin()}
\label{\detokenize{javaapi/NL:nl-asin}}\begin{quote}

\sphinxAtStartPar
Calculate arcsine of a nonlinear expression.

\sphinxAtStartPar
\sphinxstylestrong{Synopsis}
\begin{quote}

\sphinxAtStartPar
\sphinxcode{\sphinxupquote{static NlExpr asin(NlExpr expr)}}
\end{quote}

\sphinxAtStartPar
\sphinxstylestrong{Arguments}
\begin{quote}

\sphinxAtStartPar
\sphinxcode{\sphinxupquote{expr}}: a nonlinear expression.
\end{quote}

\sphinxAtStartPar
\sphinxstylestrong{Return}
\begin{quote}

\sphinxAtStartPar
result as a nonlinear expression.
\end{quote}
\end{quote}

\subsubsection{NL.asinh()}
\label{\detokenize{javaapi/NL:nl-asinh}}\begin{quote}

\sphinxAtStartPar
Calculate inverse hyperbolic sine of a nonlinear expression.

\sphinxAtStartPar
\sphinxstylestrong{Synopsis}
\begin{quote}

\sphinxAtStartPar
\sphinxcode{\sphinxupquote{static NlExpr asinh(NlExpr expr)}}
\end{quote}

\sphinxAtStartPar
\sphinxstylestrong{Arguments}
\begin{quote}

\sphinxAtStartPar
\sphinxcode{\sphinxupquote{expr}}: a nonlinear expression.
\end{quote}

\sphinxAtStartPar
\sphinxstylestrong{Return}
\begin{quote}

\sphinxAtStartPar
result as a nonlinear expression.
\end{quote}
\end{quote}

\subsubsection{NL.atan()}
\label{\detokenize{javaapi/NL:nl-atan}}\begin{quote}

\sphinxAtStartPar
Calculate arctangent of a nonlinear expression.

\sphinxAtStartPar
\sphinxstylestrong{Synopsis}
\begin{quote}

\sphinxAtStartPar
\sphinxcode{\sphinxupquote{static NlExpr atan(NlExpr expr)}}
\end{quote}

\sphinxAtStartPar
\sphinxstylestrong{Arguments}
\begin{quote}

\sphinxAtStartPar
\sphinxcode{\sphinxupquote{expr}}: a nonlinear expression.
\end{quote}

\sphinxAtStartPar
\sphinxstylestrong{Return}
\begin{quote}

\sphinxAtStartPar
result as a nonlinear expression.
\end{quote}
\end{quote}

\subsubsection{NL.atan2()}
\label{\detokenize{javaapi/NL:nl-atan2}}\begin{quote}

\sphinxAtStartPar
Calculate two\sphinxhyphen{}argument arctangent of a nonlinear expression.

\sphinxAtStartPar
\sphinxstylestrong{Synopsis}
\begin{quote}

\sphinxAtStartPar
\sphinxcode{\sphinxupquote{static NlExpr atan2(NlExpr y, NlExpr x)}}
\end{quote}

\sphinxAtStartPar
\sphinxstylestrong{Arguments}
\begin{quote}

\sphinxAtStartPar
\sphinxcode{\sphinxupquote{y}}: y coordinate as a nonlinear expression.

\sphinxAtStartPar
\sphinxcode{\sphinxupquote{x}}: x coordinate as a nonlinear expression.
\end{quote}

\sphinxAtStartPar
\sphinxstylestrong{Return}
\begin{quote}

\sphinxAtStartPar
result as a nonlinear expression.
\end{quote}
\end{quote}

\subsubsection{NL.atanh()}
\label{\detokenize{javaapi/NL:nl-atanh}}\begin{quote}

\sphinxAtStartPar
Calculate inverse hyperbolic tangent of a nonlinear expression.

\sphinxAtStartPar
\sphinxstylestrong{Synopsis}
\begin{quote}

\sphinxAtStartPar
\sphinxcode{\sphinxupquote{static NlExpr atanh(NlExpr expr)}}
\end{quote}

\sphinxAtStartPar
\sphinxstylestrong{Arguments}
\begin{quote}

\sphinxAtStartPar
\sphinxcode{\sphinxupquote{expr}}: a nonlinear expression.
\end{quote}

\sphinxAtStartPar
\sphinxstylestrong{Return}
\begin{quote}

\sphinxAtStartPar
result as a nonlinear expression.
\end{quote}
\end{quote}

\subsubsection{NL.ceil()}
\label{\detokenize{javaapi/NL:nl-ceil}}\begin{quote}

\sphinxAtStartPar
Calculate ceiling value of a nonlinear expression.

\sphinxAtStartPar
\sphinxstylestrong{Synopsis}
\begin{quote}

\sphinxAtStartPar
\sphinxcode{\sphinxupquote{static NlExpr ceil(NlExpr expr)}}
\end{quote}

\sphinxAtStartPar
\sphinxstylestrong{Arguments}
\begin{quote}

\sphinxAtStartPar
\sphinxcode{\sphinxupquote{expr}}: a nonlinear expression.
\end{quote}

\sphinxAtStartPar
\sphinxstylestrong{Return}
\begin{quote}

\sphinxAtStartPar
result as a nonlinear expression.
\end{quote}
\end{quote}

\subsubsection{NL.cos()}
\label{\detokenize{javaapi/NL:nl-cos}}\begin{quote}

\sphinxAtStartPar
Calculate cosine of a nonlinear expression.

\sphinxAtStartPar
\sphinxstylestrong{Synopsis}
\begin{quote}

\sphinxAtStartPar
\sphinxcode{\sphinxupquote{static NlExpr cos(NlExpr expr)}}
\end{quote}

\sphinxAtStartPar
\sphinxstylestrong{Arguments}
\begin{quote}

\sphinxAtStartPar
\sphinxcode{\sphinxupquote{expr}}: a nonlinear expression.
\end{quote}

\sphinxAtStartPar
\sphinxstylestrong{Return}
\begin{quote}

\sphinxAtStartPar
result as a nonlinear expression.
\end{quote}
\end{quote}

\subsubsection{NL.cosh()}
\label{\detokenize{javaapi/NL:nl-cosh}}\begin{quote}

\sphinxAtStartPar
Calculate hyperbolic cosine of a nonlinear expression.

\sphinxAtStartPar
\sphinxstylestrong{Synopsis}
\begin{quote}

\sphinxAtStartPar
\sphinxcode{\sphinxupquote{static NlExpr cosh(NlExpr expr)}}
\end{quote}

\sphinxAtStartPar
\sphinxstylestrong{Arguments}
\begin{quote}

\sphinxAtStartPar
\sphinxcode{\sphinxupquote{expr}}: a nonlinear expression.
\end{quote}

\sphinxAtStartPar
\sphinxstylestrong{Return}
\begin{quote}

\sphinxAtStartPar
result as a nonlinear expression.
\end{quote}
\end{quote}

\subsubsection{NL.exp()}
\label{\detokenize{javaapi/NL:nl-exp}}\begin{quote}

\sphinxAtStartPar
Calculate exponential function of a nonlinear expression.

\sphinxAtStartPar
\sphinxstylestrong{Synopsis}
\begin{quote}

\sphinxAtStartPar
\sphinxcode{\sphinxupquote{static NlExpr exp(NlExpr expo)}}
\end{quote}

\sphinxAtStartPar
\sphinxstylestrong{Arguments}
\begin{quote}

\sphinxAtStartPar
\sphinxcode{\sphinxupquote{expo}}: exponent as a nonlinear expression.
\end{quote}

\sphinxAtStartPar
\sphinxstylestrong{Return}
\begin{quote}

\sphinxAtStartPar
result as a nonlinear expression.
\end{quote}
\end{quote}

\subsubsection{NL.floor()}
\label{\detokenize{javaapi/NL:nl-floor}}\begin{quote}

\sphinxAtStartPar
Calculate floor value of a nonlinear expression.

\sphinxAtStartPar
\sphinxstylestrong{Synopsis}
\begin{quote}

\sphinxAtStartPar
\sphinxcode{\sphinxupquote{static NlExpr floor(NlExpr expr)}}
\end{quote}

\sphinxAtStartPar
\sphinxstylestrong{Arguments}
\begin{quote}

\sphinxAtStartPar
\sphinxcode{\sphinxupquote{expr}}: a nonlinear expression.
\end{quote}

\sphinxAtStartPar
\sphinxstylestrong{Return}
\begin{quote}

\sphinxAtStartPar
result as a nonlinear expression.
\end{quote}
\end{quote}

\subsubsection{NL.log10()}
\label{\detokenize{javaapi/NL:nl-log10}}\begin{quote}

\sphinxAtStartPar
Calculate logarithmic function of a nonlinear expression with base 10.

\sphinxAtStartPar
\sphinxstylestrong{Synopsis}
\begin{quote}

\sphinxAtStartPar
\sphinxcode{\sphinxupquote{static NlExpr log10(NlExpr expr)}}
\end{quote}

\sphinxAtStartPar
\sphinxstylestrong{Arguments}
\begin{quote}

\sphinxAtStartPar
\sphinxcode{\sphinxupquote{expr}}: a nonlinear expression.
\end{quote}

\sphinxAtStartPar
\sphinxstylestrong{Return}
\begin{quote}

\sphinxAtStartPar
result as a nonlinear expression.
\end{quote}
\end{quote}

\subsubsection{NL.log()}
\label{\detokenize{javaapi/NL:nl-log}}\begin{quote}

\sphinxAtStartPar
Calculate nature logarithmic function of a nonlinear expression.

\sphinxAtStartPar
\sphinxstylestrong{Synopsis}
\begin{quote}

\sphinxAtStartPar
\sphinxcode{\sphinxupquote{static NlExpr log(NlExpr expr)}}
\end{quote}

\sphinxAtStartPar
\sphinxstylestrong{Arguments}
\begin{quote}

\sphinxAtStartPar
\sphinxcode{\sphinxupquote{expr}}: a nonlinear expression.
\end{quote}

\sphinxAtStartPar
\sphinxstylestrong{Return}
\begin{quote}

\sphinxAtStartPar
result as a nonlinear expression.
\end{quote}
\end{quote}

\subsubsection{NL.neg()}
\label{\detokenize{javaapi/NL:nl-neg}}\begin{quote}

\sphinxAtStartPar
Calculate negative value of a nonlinear expression.

\sphinxAtStartPar
\sphinxstylestrong{Synopsis}
\begin{quote}

\sphinxAtStartPar
\sphinxcode{\sphinxupquote{static NlExpr neg(NlExpr expr)}}
\end{quote}

\sphinxAtStartPar
\sphinxstylestrong{Arguments}
\begin{quote}

\sphinxAtStartPar
\sphinxcode{\sphinxupquote{expr}}: a nonlinear expression.
\end{quote}

\sphinxAtStartPar
\sphinxstylestrong{Return}
\begin{quote}

\sphinxAtStartPar
result as a nonlinear expression.
\end{quote}
\end{quote}

\subsubsection{NL.pow()}
\label{\detokenize{javaapi/NL:nl-pow}}\begin{quote}

\sphinxAtStartPar
Calculate power function of a nonlinear expression.

\sphinxAtStartPar
\sphinxstylestrong{Synopsis}
\begin{quote}

\sphinxAtStartPar
\sphinxcode{\sphinxupquote{static NlExpr pow(NlExpr base, NlExpr expo)}}
\end{quote}

\sphinxAtStartPar
\sphinxstylestrong{Arguments}
\begin{quote}

\sphinxAtStartPar
\sphinxcode{\sphinxupquote{base}}: base as a nonlinear expression.

\sphinxAtStartPar
\sphinxcode{\sphinxupquote{expo}}: exponent as a nonlinear expression.
\end{quote}

\sphinxAtStartPar
\sphinxstylestrong{Return}
\begin{quote}

\sphinxAtStartPar
result as a nonlinear expression.
\end{quote}
\end{quote}

\subsubsection{NL.sin()}
\label{\detokenize{javaapi/NL:nl-sin}}\begin{quote}

\sphinxAtStartPar
Calculate sine of a nonlinear expression.

\sphinxAtStartPar
\sphinxstylestrong{Synopsis}
\begin{quote}

\sphinxAtStartPar
\sphinxcode{\sphinxupquote{static NlExpr sin(NlExpr expr)}}
\end{quote}

\sphinxAtStartPar
\sphinxstylestrong{Arguments}
\begin{quote}

\sphinxAtStartPar
\sphinxcode{\sphinxupquote{expr}}: a nonlinear expression.
\end{quote}

\sphinxAtStartPar
\sphinxstylestrong{Return}
\begin{quote}

\sphinxAtStartPar
result as a nonlinear expression.
\end{quote}
\end{quote}

\subsubsection{NL.sinh()}
\label{\detokenize{javaapi/NL:nl-sinh}}\begin{quote}

\sphinxAtStartPar
Calculate hyperbolic sine of a nonlinear expression.

\sphinxAtStartPar
\sphinxstylestrong{Synopsis}
\begin{quote}

\sphinxAtStartPar
\sphinxcode{\sphinxupquote{static NlExpr sinh(NlExpr expr)}}
\end{quote}

\sphinxAtStartPar
\sphinxstylestrong{Arguments}
\begin{quote}

\sphinxAtStartPar
\sphinxcode{\sphinxupquote{expr}}: a nonlinear expression.
\end{quote}

\sphinxAtStartPar
\sphinxstylestrong{Return}
\begin{quote}

\sphinxAtStartPar
result as a nonlinear expression.
\end{quote}
\end{quote}

\subsubsection{NL.sqrt()}
\label{\detokenize{javaapi/NL:nl-sqrt}}\begin{quote}

\sphinxAtStartPar
Calculate square root of a nonlinear expression.

\sphinxAtStartPar
\sphinxstylestrong{Synopsis}
\begin{quote}

\sphinxAtStartPar
\sphinxcode{\sphinxupquote{static NlExpr sqrt(NlExpr expr)}}
\end{quote}

\sphinxAtStartPar
\sphinxstylestrong{Arguments}
\begin{quote}

\sphinxAtStartPar
\sphinxcode{\sphinxupquote{expr}}: a nonlinear expression.
\end{quote}

\sphinxAtStartPar
\sphinxstylestrong{Return}
\begin{quote}

\sphinxAtStartPar
result as a nonlinear expression.
\end{quote}
\end{quote}

\subsubsection{NL.sum()}
\label{\detokenize{javaapi/NL:nl-sum}}\begin{quote}

\sphinxAtStartPar
Sum of nonlinear expressions.

\sphinxAtStartPar
\sphinxstylestrong{Synopsis}
\begin{quote}

\sphinxAtStartPar
\sphinxcode{\sphinxupquote{static NlExpr sum(}}
\begin{quote}

\sphinxAtStartPar
\sphinxcode{\sphinxupquote{NlExpr op1,}}

\sphinxAtStartPar
\sphinxcode{\sphinxupquote{NlExpr op2,}}

\sphinxAtStartPar
\sphinxcode{\sphinxupquote{NlExpr op3,}}

\sphinxAtStartPar
\sphinxcode{\sphinxupquote{NlExpr op4)}}
\end{quote}
\end{quote}

\sphinxAtStartPar
\sphinxstylestrong{Arguments}
\begin{quote}

\sphinxAtStartPar
\sphinxcode{\sphinxupquote{op1}}: first nonlinear expression.

\sphinxAtStartPar
\sphinxcode{\sphinxupquote{op2}}: second nonlinear expression.

\sphinxAtStartPar
\sphinxcode{\sphinxupquote{op3}}: third nonlinear expression.

\sphinxAtStartPar
\sphinxcode{\sphinxupquote{op4}}: fourth nonlinear expression.
\end{quote}

\sphinxAtStartPar
\sphinxstylestrong{Return}
\begin{quote}

\sphinxAtStartPar
result as a nonlinear expression.
\end{quote}
\end{quote}

\subsubsection{NL.sum()}
\label{\detokenize{javaapi/NL:id1}}\begin{quote}

\sphinxAtStartPar
Sum of nonlinear expressions.

\sphinxAtStartPar
\sphinxstylestrong{Synopsis}
\begin{quote}

\sphinxAtStartPar
\sphinxcode{\sphinxupquote{static NlExpr sum(}}
\begin{quote}

\sphinxAtStartPar
\sphinxcode{\sphinxupquote{NlExpr op1,}}

\sphinxAtStartPar
\sphinxcode{\sphinxupquote{NlExpr op2,}}

\sphinxAtStartPar
\sphinxcode{\sphinxupquote{NlExpr op3)}}
\end{quote}
\end{quote}

\sphinxAtStartPar
\sphinxstylestrong{Arguments}
\begin{quote}

\sphinxAtStartPar
\sphinxcode{\sphinxupquote{op1}}: first nonlinear expression.

\sphinxAtStartPar
\sphinxcode{\sphinxupquote{op2}}: second nonlinear expression.

\sphinxAtStartPar
\sphinxcode{\sphinxupquote{op3}}: third nonlinear expression.
\end{quote}

\sphinxAtStartPar
\sphinxstylestrong{Return}
\begin{quote}

\sphinxAtStartPar
result as a nonlinear expression.
\end{quote}
\end{quote}

\subsubsection{NL.sum()}
\label{\detokenize{javaapi/NL:id2}}\begin{quote}

\sphinxAtStartPar
Sum of nonlinear expressions.

\sphinxAtStartPar
\sphinxstylestrong{Synopsis}
\begin{quote}

\sphinxAtStartPar
\sphinxcode{\sphinxupquote{static NlExpr sum(NlExpr{[}{]} exprs)}}
\end{quote}

\sphinxAtStartPar
\sphinxstylestrong{Arguments}
\begin{quote}

\sphinxAtStartPar
\sphinxcode{\sphinxupquote{exprs}}: array of nonlinear expressions.
\end{quote}

\sphinxAtStartPar
\sphinxstylestrong{Return}
\begin{quote}

\sphinxAtStartPar
result as a nonlinear expression.
\end{quote}
\end{quote}

\subsubsection{NL.tan()}
\label{\detokenize{javaapi/NL:nl-tan}}\begin{quote}

\sphinxAtStartPar
Calculate tangent of a nonlinear expression.

\sphinxAtStartPar
\sphinxstylestrong{Synopsis}
\begin{quote}

\sphinxAtStartPar
\sphinxcode{\sphinxupquote{static NlExpr tan(NlExpr expr)}}
\end{quote}

\sphinxAtStartPar
\sphinxstylestrong{Arguments}
\begin{quote}

\sphinxAtStartPar
\sphinxcode{\sphinxupquote{expr}}: a nonlinear expression.
\end{quote}

\sphinxAtStartPar
\sphinxstylestrong{Return}
\begin{quote}

\sphinxAtStartPar
result as a nonlinear expression.
\end{quote}
\end{quote}

\subsubsection{NL.tanh()}
\label{\detokenize{javaapi/NL:nl-tanh}}\begin{quote}

\sphinxAtStartPar
Calculat hyperbolic tangent of a nonlinear expression.

\sphinxAtStartPar
\sphinxstylestrong{Synopsis}
\begin{quote}

\sphinxAtStartPar
\sphinxcode{\sphinxupquote{static NlExpr tanh(NlExpr expr)}}
\end{quote}

\sphinxAtStartPar
\sphinxstylestrong{Arguments}
\begin{quote}

\sphinxAtStartPar
\sphinxcode{\sphinxupquote{expr}}: a nonlinear expression.
\end{quote}

\sphinxAtStartPar
\sphinxstylestrong{Return}
\begin{quote}

\sphinxAtStartPar
result as a nonlinear expression.
\end{quote}
\end{quote}

\subsection{CallbackBase}
\label{\detokenize{javaapiref:callbackbase}}\label{\detokenize{javaapiref:chapjavaapiref-callback}}
\sphinxAtStartPar
COPT Callback abstract base object. Users must implment its virtual method
\sphinxcode{\sphinxupquote{virtual void CallbackBase::callback()}} to instantiate an instance, which
pass to \sphinxcode{\sphinxupquote{Model::SetCallback(CallbackBase cb, int cbctx)}} as the first
parameter. Subclass of CallbackBase inherits the following member methods:

\sphinxstepscope

\subsubsection{CallbackBase.CallbackBase()}
\label{\detokenize{javaapi/CallbackBase:callbackbase-callbackbase}}\label{\detokenize{javaapi/CallbackBase::doc}}\begin{quote}

\sphinxAtStartPar
Constructor of CallbackBase, implementing ICallback interface.

\sphinxAtStartPar
\sphinxstylestrong{Synopsis}
\begin{quote}

\sphinxAtStartPar
\sphinxcode{\sphinxupquote{CallbackBase()}}
\end{quote}
\end{quote}

\subsubsection{CallbackBase.addLazyConstr()}
\label{\detokenize{javaapi/CallbackBase:callbackbase-addlazyconstr}}\begin{quote}

\sphinxAtStartPar
Add a lazy constraint to model.

\sphinxAtStartPar
\sphinxstylestrong{Synopsis}
\begin{quote}

\sphinxAtStartPar
\sphinxcode{\sphinxupquote{void addLazyConstr(}}
\begin{quote}

\sphinxAtStartPar
\sphinxcode{\sphinxupquote{Expr lhs,}}

\sphinxAtStartPar
\sphinxcode{\sphinxupquote{char sense,}}

\sphinxAtStartPar
\sphinxcode{\sphinxupquote{double rhs)}}
\end{quote}
\end{quote}

\sphinxAtStartPar
\sphinxstylestrong{Arguments}
\begin{quote}

\sphinxAtStartPar
\sphinxcode{\sphinxupquote{lhs}}: expression for lazy contraint.

\sphinxAtStartPar
\sphinxcode{\sphinxupquote{sense}}: sense for lazy constraint.

\sphinxAtStartPar
\sphinxcode{\sphinxupquote{rhs}}: right hand side value for lazy constraint.
\end{quote}
\end{quote}

\subsubsection{CallbackBase.addLazyConstr()}
\label{\detokenize{javaapi/CallbackBase:id1}}\begin{quote}

\sphinxAtStartPar
Add a lazy constraint to model.

\sphinxAtStartPar
\sphinxstylestrong{Synopsis}
\begin{quote}

\sphinxAtStartPar
\sphinxcode{\sphinxupquote{void addLazyConstr(}}
\begin{quote}

\sphinxAtStartPar
\sphinxcode{\sphinxupquote{Expr lhs,}}

\sphinxAtStartPar
\sphinxcode{\sphinxupquote{char sense,}}

\sphinxAtStartPar
\sphinxcode{\sphinxupquote{Expr rhs)}}
\end{quote}
\end{quote}

\sphinxAtStartPar
\sphinxstylestrong{Arguments}
\begin{quote}

\sphinxAtStartPar
\sphinxcode{\sphinxupquote{lhs}}: left hand side expression for lazy contraint.

\sphinxAtStartPar
\sphinxcode{\sphinxupquote{sense}}: sense for lazy constraint.

\sphinxAtStartPar
\sphinxcode{\sphinxupquote{rhs}}: right hand side expression for lazy contraint.
\end{quote}
\end{quote}

\subsubsection{CallbackBase.addLazyConstr()}
\label{\detokenize{javaapi/CallbackBase:id2}}\begin{quote}

\sphinxAtStartPar
Add a lazy constraint to model.

\sphinxAtStartPar
\sphinxstylestrong{Synopsis}
\begin{quote}

\sphinxAtStartPar
\sphinxcode{\sphinxupquote{void addLazyConstr(ConstrBuilder builder)}}
\end{quote}

\sphinxAtStartPar
\sphinxstylestrong{Arguments}
\begin{quote}

\sphinxAtStartPar
\sphinxcode{\sphinxupquote{builder}}: builder for lazy contraint.
\end{quote}
\end{quote}

\subsubsection{CallbackBase.addLazyConstrs()}
\label{\detokenize{javaapi/CallbackBase:callbackbase-addlazyconstrs}}\begin{quote}

\sphinxAtStartPar
Add lazy constraints to model.

\sphinxAtStartPar
\sphinxstylestrong{Synopsis}
\begin{quote}

\sphinxAtStartPar
\sphinxcode{\sphinxupquote{void addLazyConstrs(ConstrBuilderArray builders)}}
\end{quote}

\sphinxAtStartPar
\sphinxstylestrong{Arguments}
\begin{quote}

\sphinxAtStartPar
\sphinxcode{\sphinxupquote{builders}}: array of builders for lazy contraints.
\end{quote}
\end{quote}

\subsubsection{CallbackBase.addUserCut()}
\label{\detokenize{javaapi/CallbackBase:callbackbase-addusercut}}\begin{quote}

\sphinxAtStartPar
Add a user cut to model.

\sphinxAtStartPar
\sphinxstylestrong{Synopsis}
\begin{quote}

\sphinxAtStartPar
\sphinxcode{\sphinxupquote{void addUserCut(}}
\begin{quote}

\sphinxAtStartPar
\sphinxcode{\sphinxupquote{Expr lhs,}}

\sphinxAtStartPar
\sphinxcode{\sphinxupquote{char sense,}}

\sphinxAtStartPar
\sphinxcode{\sphinxupquote{double rhs)}}
\end{quote}
\end{quote}

\sphinxAtStartPar
\sphinxstylestrong{Arguments}
\begin{quote}

\sphinxAtStartPar
\sphinxcode{\sphinxupquote{lhs}}: expression for user cut.

\sphinxAtStartPar
\sphinxcode{\sphinxupquote{sense}}: sense for user cut.

\sphinxAtStartPar
\sphinxcode{\sphinxupquote{rhs}}: right hand side value for user cut.
\end{quote}
\end{quote}

\subsubsection{CallbackBase.addUserCut()}
\label{\detokenize{javaapi/CallbackBase:id3}}\begin{quote}

\sphinxAtStartPar
Add a user cut to model.

\sphinxAtStartPar
\sphinxstylestrong{Synopsis}
\begin{quote}

\sphinxAtStartPar
\sphinxcode{\sphinxupquote{void addUserCut(}}
\begin{quote}

\sphinxAtStartPar
\sphinxcode{\sphinxupquote{Expr lhs,}}

\sphinxAtStartPar
\sphinxcode{\sphinxupquote{char sense,}}

\sphinxAtStartPar
\sphinxcode{\sphinxupquote{Expr rhs)}}
\end{quote}
\end{quote}

\sphinxAtStartPar
\sphinxstylestrong{Arguments}
\begin{quote}

\sphinxAtStartPar
\sphinxcode{\sphinxupquote{lhs}}: left hand side expression for user cut.

\sphinxAtStartPar
\sphinxcode{\sphinxupquote{sense}}: sense for user cut.

\sphinxAtStartPar
\sphinxcode{\sphinxupquote{rhs}}: right hand side expression for user cut.
\end{quote}
\end{quote}

\subsubsection{CallbackBase.addUserCut()}
\label{\detokenize{javaapi/CallbackBase:id4}}\begin{quote}

\sphinxAtStartPar
Add a user cut to model.

\sphinxAtStartPar
\sphinxstylestrong{Synopsis}
\begin{quote}

\sphinxAtStartPar
\sphinxcode{\sphinxupquote{void addUserCut(ConstrBuilder builder)}}
\end{quote}

\sphinxAtStartPar
\sphinxstylestrong{Arguments}
\begin{quote}

\sphinxAtStartPar
\sphinxcode{\sphinxupquote{builder}}: builder for user cut.
\end{quote}
\end{quote}

\subsubsection{CallbackBase.addUserCuts()}
\label{\detokenize{javaapi/CallbackBase:callbackbase-addusercuts}}\begin{quote}

\sphinxAtStartPar
Add user cuts to model.

\sphinxAtStartPar
\sphinxstylestrong{Synopsis}
\begin{quote}

\sphinxAtStartPar
\sphinxcode{\sphinxupquote{void addUserCuts(ConstrBuilderArray builders)}}
\end{quote}

\sphinxAtStartPar
\sphinxstylestrong{Arguments}
\begin{quote}

\sphinxAtStartPar
\sphinxcode{\sphinxupquote{builders}}: array of builders for user cuts.
\end{quote}
\end{quote}

\subsubsection{CallbackBase.callback()}
\label{\detokenize{javaapi/CallbackBase:callbackbase-callback}}\begin{quote}

\sphinxAtStartPar
Pure virtual function defined in ICallback interface. User must implement it.

\sphinxAtStartPar
\sphinxstylestrong{Synopsis}
\begin{quote}

\sphinxAtStartPar
\sphinxcode{\sphinxupquote{void callback()}}
\end{quote}
\end{quote}

\subsubsection{CallbackBase.getDblInfo()}
\label{\detokenize{javaapi/CallbackBase:callbackbase-getdblinfo}}\begin{quote}

\sphinxAtStartPar
Get double value of given information name in callback.

\sphinxAtStartPar
\sphinxstylestrong{Synopsis}
\begin{quote}

\sphinxAtStartPar
\sphinxcode{\sphinxupquote{double getDblInfo(String cbinfo)}}
\end{quote}

\sphinxAtStartPar
\sphinxstylestrong{Arguments}
\begin{quote}

\sphinxAtStartPar
\sphinxcode{\sphinxupquote{cbinfo}}: name of callback info.
\end{quote}

\sphinxAtStartPar
\sphinxstylestrong{Return}
\begin{quote}

\sphinxAtStartPar
value of desired information.
\end{quote}
\end{quote}

\subsubsection{CallbackBase.getIncumbent()}
\label{\detokenize{javaapi/CallbackBase:callbackbase-getincumbent}}\begin{quote}

\sphinxAtStartPar
Get best feasible solution of given variable in callback.

\sphinxAtStartPar
\sphinxstylestrong{Synopsis}
\begin{quote}

\sphinxAtStartPar
\sphinxcode{\sphinxupquote{double getIncumbent(Var var)}}
\end{quote}

\sphinxAtStartPar
\sphinxstylestrong{Arguments}
\begin{quote}

\sphinxAtStartPar
\sphinxcode{\sphinxupquote{var}}: given variable.
\end{quote}

\sphinxAtStartPar
\sphinxstylestrong{Return}
\begin{quote}

\sphinxAtStartPar
best feasible solution of given variable.
\end{quote}
\end{quote}

\subsubsection{CallbackBase.getIncumbent()}
\label{\detokenize{javaapi/CallbackBase:id5}}\begin{quote}

\sphinxAtStartPar
Get best feasible solution of variables in callback.

\sphinxAtStartPar
\sphinxstylestrong{Synopsis}
\begin{quote}

\sphinxAtStartPar
\sphinxcode{\sphinxupquote{double{[}{]} getIncumbent(VarArray vars)}}
\end{quote}

\sphinxAtStartPar
\sphinxstylestrong{Arguments}
\begin{quote}

\sphinxAtStartPar
\sphinxcode{\sphinxupquote{vars}}: an array of variables.
\end{quote}

\sphinxAtStartPar
\sphinxstylestrong{Return}
\begin{quote}

\sphinxAtStartPar
best feasible solution of desired variables.
\end{quote}
\end{quote}

\subsubsection{CallbackBase.getIncumbent()}
\label{\detokenize{javaapi/CallbackBase:id6}}\begin{quote}

\sphinxAtStartPar
Get best feasible solution of variables in callback.

\sphinxAtStartPar
\sphinxstylestrong{Synopsis}
\begin{quote}

\sphinxAtStartPar
\sphinxcode{\sphinxupquote{double{[}{]} getIncumbent(Var{[}{]} vars)}}
\end{quote}

\sphinxAtStartPar
\sphinxstylestrong{Arguments}
\begin{quote}

\sphinxAtStartPar
\sphinxcode{\sphinxupquote{vars}}: an array of variables.
\end{quote}

\sphinxAtStartPar
\sphinxstylestrong{Return}
\begin{quote}

\sphinxAtStartPar
best feasible solution of desired variables.
\end{quote}
\end{quote}

\subsubsection{CallbackBase.getIncumbent()}
\label{\detokenize{javaapi/CallbackBase:id7}}\begin{quote}

\sphinxAtStartPar
Get best feasible solution of all variables in callback.

\sphinxAtStartPar
\sphinxstylestrong{Synopsis}
\begin{quote}

\sphinxAtStartPar
\sphinxcode{\sphinxupquote{double{[}{]} getIncumbent()}}
\end{quote}

\sphinxAtStartPar
\sphinxstylestrong{Return}
\begin{quote}

\sphinxAtStartPar
best feasible solution of all variables.
\end{quote}
\end{quote}

\subsubsection{CallbackBase.getIntInfo()}
\label{\detokenize{javaapi/CallbackBase:callbackbase-getintinfo}}\begin{quote}

\sphinxAtStartPar
Get integer value of given information name in callback.

\sphinxAtStartPar
\sphinxstylestrong{Synopsis}
\begin{quote}

\sphinxAtStartPar
\sphinxcode{\sphinxupquote{int getIntInfo(String cbinfo)}}
\end{quote}

\sphinxAtStartPar
\sphinxstylestrong{Arguments}
\begin{quote}

\sphinxAtStartPar
\sphinxcode{\sphinxupquote{cbinfo}}: name of callback info.
\end{quote}

\sphinxAtStartPar
\sphinxstylestrong{Return}
\begin{quote}

\sphinxAtStartPar
value of desired information.
\end{quote}
\end{quote}

\subsubsection{CallbackBase.getRelaxSol()}
\label{\detokenize{javaapi/CallbackBase:callbackbase-getrelaxsol}}\begin{quote}

\sphinxAtStartPar
Get LP\sphinxhyphen{}relaxation solution of given variable in callback.

\sphinxAtStartPar
\sphinxstylestrong{Synopsis}
\begin{quote}

\sphinxAtStartPar
\sphinxcode{\sphinxupquote{double getRelaxSol(Var var)}}
\end{quote}

\sphinxAtStartPar
\sphinxstylestrong{Arguments}
\begin{quote}

\sphinxAtStartPar
\sphinxcode{\sphinxupquote{var}}: given variable.
\end{quote}

\sphinxAtStartPar
\sphinxstylestrong{Return}
\begin{quote}

\sphinxAtStartPar
LP\sphinxhyphen{}relaxation solution of given variable.
\end{quote}
\end{quote}

\subsubsection{CallbackBase.getRelaxSol()}
\label{\detokenize{javaapi/CallbackBase:id8}}\begin{quote}

\sphinxAtStartPar
Get LP\sphinxhyphen{}relaxation solution of variables in callback.

\sphinxAtStartPar
\sphinxstylestrong{Synopsis}
\begin{quote}

\sphinxAtStartPar
\sphinxcode{\sphinxupquote{double{[}{]} getRelaxSol(VarArray vars)}}
\end{quote}

\sphinxAtStartPar
\sphinxstylestrong{Arguments}
\begin{quote}

\sphinxAtStartPar
\sphinxcode{\sphinxupquote{vars}}: an array of variables.
\end{quote}

\sphinxAtStartPar
\sphinxstylestrong{Return}
\begin{quote}

\sphinxAtStartPar
LP\sphinxhyphen{}relaxation solution of variables.
\end{quote}
\end{quote}

\subsubsection{CallbackBase.getRelaxSol()}
\label{\detokenize{javaapi/CallbackBase:id9}}\begin{quote}

\sphinxAtStartPar
Get LP\sphinxhyphen{}relaxation solution of variables in callback.

\sphinxAtStartPar
\sphinxstylestrong{Synopsis}
\begin{quote}

\sphinxAtStartPar
\sphinxcode{\sphinxupquote{double{[}{]} getRelaxSol(Var{[}{]} vars)}}
\end{quote}

\sphinxAtStartPar
\sphinxstylestrong{Arguments}
\begin{quote}

\sphinxAtStartPar
\sphinxcode{\sphinxupquote{vars}}: an array of variables.
\end{quote}

\sphinxAtStartPar
\sphinxstylestrong{Return}
\begin{quote}

\sphinxAtStartPar
LP\sphinxhyphen{}relaxation solution of variables.
\end{quote}
\end{quote}

\subsubsection{CallbackBase.getRelaxSol()}
\label{\detokenize{javaapi/CallbackBase:id10}}\begin{quote}

\sphinxAtStartPar
Get LP\sphinxhyphen{}relaxation solution of all variables in callback.

\sphinxAtStartPar
\sphinxstylestrong{Synopsis}
\begin{quote}

\sphinxAtStartPar
\sphinxcode{\sphinxupquote{double{[}{]} getRelaxSol()}}
\end{quote}

\sphinxAtStartPar
\sphinxstylestrong{Return}
\begin{quote}

\sphinxAtStartPar
LP\sphinxhyphen{}relaxation solution of all variables.
\end{quote}
\end{quote}

\subsubsection{CallbackBase.getSolution()}
\label{\detokenize{javaapi/CallbackBase:callbackbase-getsolution}}\begin{quote}

\sphinxAtStartPar
Get solution of given variable in callback.

\sphinxAtStartPar
\sphinxstylestrong{Synopsis}
\begin{quote}

\sphinxAtStartPar
\sphinxcode{\sphinxupquote{double getSolution(Var var)}}
\end{quote}

\sphinxAtStartPar
\sphinxstylestrong{Arguments}
\begin{quote}

\sphinxAtStartPar
\sphinxcode{\sphinxupquote{var}}: given variable.
\end{quote}

\sphinxAtStartPar
\sphinxstylestrong{Return}
\begin{quote}

\sphinxAtStartPar
solution of given variable.
\end{quote}
\end{quote}

\subsubsection{CallbackBase.getSolution()}
\label{\detokenize{javaapi/CallbackBase:id11}}\begin{quote}

\sphinxAtStartPar
Get solution of variables in callback.

\sphinxAtStartPar
\sphinxstylestrong{Synopsis}
\begin{quote}

\sphinxAtStartPar
\sphinxcode{\sphinxupquote{double{[}{]} getSolution(VarArray vars)}}
\end{quote}

\sphinxAtStartPar
\sphinxstylestrong{Arguments}
\begin{quote}

\sphinxAtStartPar
\sphinxcode{\sphinxupquote{vars}}: an array of variables.
\end{quote}

\sphinxAtStartPar
\sphinxstylestrong{Return}
\begin{quote}

\sphinxAtStartPar
solution of variables.
\end{quote}
\end{quote}

\subsubsection{CallbackBase.getSolution()}
\label{\detokenize{javaapi/CallbackBase:id12}}\begin{quote}

\sphinxAtStartPar
Get solution of variables in callback.

\sphinxAtStartPar
\sphinxstylestrong{Synopsis}
\begin{quote}

\sphinxAtStartPar
\sphinxcode{\sphinxupquote{double{[}{]} getSolution(Var{[}{]} vars)}}
\end{quote}

\sphinxAtStartPar
\sphinxstylestrong{Arguments}
\begin{quote}

\sphinxAtStartPar
\sphinxcode{\sphinxupquote{vars}}: an array of variables.
\end{quote}

\sphinxAtStartPar
\sphinxstylestrong{Return}
\begin{quote}

\sphinxAtStartPar
solution of variables.
\end{quote}
\end{quote}

\subsubsection{CallbackBase.getSolution()}
\label{\detokenize{javaapi/CallbackBase:id13}}\begin{quote}

\sphinxAtStartPar
Get solution of all variables in callback.

\sphinxAtStartPar
\sphinxstylestrong{Synopsis}
\begin{quote}

\sphinxAtStartPar
\sphinxcode{\sphinxupquote{double{[}{]} getSolution()}}
\end{quote}

\sphinxAtStartPar
\sphinxstylestrong{Return}
\begin{quote}

\sphinxAtStartPar
solution of all variables.
\end{quote}
\end{quote}

\subsubsection{CallbackBase.interrupt()}
\label{\detokenize{javaapi/CallbackBase:callbackbase-interrupt}}\begin{quote}

\sphinxAtStartPar
Interrupt solving problems in callback

\sphinxAtStartPar
\sphinxstylestrong{Synopsis}
\begin{quote}

\sphinxAtStartPar
\sphinxcode{\sphinxupquote{void interrupt()}}
\end{quote}
\end{quote}

\subsubsection{CallbackBase.loadSolution()}
\label{\detokenize{javaapi/CallbackBase:callbackbase-loadsolution}}\begin{quote}

\sphinxAtStartPar
Load customized solution to model.

\sphinxAtStartPar
\sphinxstylestrong{Synopsis}
\begin{quote}

\sphinxAtStartPar
\sphinxcode{\sphinxupquote{double loadSolution()}}
\end{quote}

\sphinxAtStartPar
\sphinxstylestrong{Return}
\begin{quote}

\sphinxAtStartPar
objective value of given solution.
\end{quote}
\end{quote}

\subsubsection{CallbackBase.setSolution()}
\label{\detokenize{javaapi/CallbackBase:callbackbase-setsolution}}\begin{quote}

\sphinxAtStartPar
Set solution of a given variable in callback.

\sphinxAtStartPar
\sphinxstylestrong{Synopsis}
\begin{quote}

\sphinxAtStartPar
\sphinxcode{\sphinxupquote{void setSolution(Var var, double val)}}
\end{quote}

\sphinxAtStartPar
\sphinxstylestrong{Arguments}
\begin{quote}

\sphinxAtStartPar
\sphinxcode{\sphinxupquote{var}}: a variable object.

\sphinxAtStartPar
\sphinxcode{\sphinxupquote{val}}: double value.
\end{quote}
\end{quote}

\subsubsection{CallbackBase.setSolution()}
\label{\detokenize{javaapi/CallbackBase:id14}}\begin{quote}

\sphinxAtStartPar
Set solution of variables in callback.

\sphinxAtStartPar
\sphinxstylestrong{Synopsis}
\begin{quote}

\sphinxAtStartPar
\sphinxcode{\sphinxupquote{void setSolution(VarArray vars, double{[}{]} vals)}}
\end{quote}

\sphinxAtStartPar
\sphinxstylestrong{Arguments}
\begin{quote}

\sphinxAtStartPar
\sphinxcode{\sphinxupquote{vars}}: an array of variable objects.

\sphinxAtStartPar
\sphinxcode{\sphinxupquote{vals}}: an array of double values.
\end{quote}
\end{quote}

\subsubsection{CallbackBase.setSolution()}
\label{\detokenize{javaapi/CallbackBase:id15}}\begin{quote}

\sphinxAtStartPar
Set solution of variables in callback.

\sphinxAtStartPar
\sphinxstylestrong{Synopsis}
\begin{quote}

\sphinxAtStartPar
\sphinxcode{\sphinxupquote{void setSolution(Var{[}{]} vars, double{[}{]} vals)}}
\end{quote}

\sphinxAtStartPar
\sphinxstylestrong{Arguments}
\begin{quote}

\sphinxAtStartPar
\sphinxcode{\sphinxupquote{vars}}: an array of variable objects.

\sphinxAtStartPar
\sphinxcode{\sphinxupquote{vals}}: an array of double values.
\end{quote}
\end{quote}

\subsubsection{CallbackBase.where()}
\label{\detokenize{javaapi/CallbackBase:callbackbase-where}}\begin{quote}

\sphinxAtStartPar
Get context in callback.

\sphinxAtStartPar
\sphinxstylestrong{Synopsis}
\begin{quote}

\sphinxAtStartPar
\sphinxcode{\sphinxupquote{int where()}}
\end{quote}

\sphinxAtStartPar
\sphinxstylestrong{Return}
\begin{quote}

\sphinxAtStartPar
integer value of context.
\end{quote}
\end{quote}

\subsection{ProbBuffer}
\label{\detokenize{javaapiref:probbuffer}}\label{\detokenize{javaapiref:chapjavaapiref-probbuffer}}
\sphinxAtStartPar
Buffer object for COPT problem. ProbBuffer object holds the (MPS) problem in string format.

\sphinxstepscope

\subsubsection{ProbBuffer.ProbBuffer()}
\label{\detokenize{javaapi/ProbBuffer:probbuffer-probbuffer}}\label{\detokenize{javaapi/ProbBuffer::doc}}\begin{quote}

\sphinxAtStartPar
Constructor of ProbBuffer object.

\sphinxAtStartPar
\sphinxstylestrong{Synopsis}
\begin{quote}

\sphinxAtStartPar
\sphinxcode{\sphinxupquote{ProbBuffer(int sz)}}
\end{quote}

\sphinxAtStartPar
\sphinxstylestrong{Arguments}
\begin{quote}

\sphinxAtStartPar
\sphinxcode{\sphinxupquote{sz}}: initial size of the problem buffer.
\end{quote}
\end{quote}

\subsubsection{ProbBuffer.getData()}
\label{\detokenize{javaapi/ProbBuffer:probbuffer-getdata}}\begin{quote}

\sphinxAtStartPar
Get string of problem in problem buffer.

\sphinxAtStartPar
\sphinxstylestrong{Synopsis}
\begin{quote}

\sphinxAtStartPar
\sphinxcode{\sphinxupquote{String getData()}}
\end{quote}

\sphinxAtStartPar
\sphinxstylestrong{Return}
\begin{quote}

\sphinxAtStartPar
string of problem in problem buffer.
\end{quote}
\end{quote}

\subsubsection{ProbBuffer.resize()}
\label{\detokenize{javaapi/ProbBuffer:probbuffer-resize}}\begin{quote}

\sphinxAtStartPar
Resize buffer to given size, and zero\sphinxhyphen{}ended.

\sphinxAtStartPar
\sphinxstylestrong{Synopsis}
\begin{quote}

\sphinxAtStartPar
\sphinxcode{\sphinxupquote{void resize(int sz)}}
\end{quote}

\sphinxAtStartPar
\sphinxstylestrong{Arguments}
\begin{quote}

\sphinxAtStartPar
\sphinxcode{\sphinxupquote{sz}}: new buffer size.
\end{quote}
\end{quote}

\subsubsection{ProbBuffer.size()}
\label{\detokenize{javaapi/ProbBuffer:probbuffer-size}}\begin{quote}

\sphinxAtStartPar
Get the size of problem buffer.

\sphinxAtStartPar
\sphinxstylestrong{Synopsis}
\begin{quote}

\sphinxAtStartPar
\sphinxcode{\sphinxupquote{int size()}}
\end{quote}

\sphinxAtStartPar
\sphinxstylestrong{Return}
\begin{quote}

\sphinxAtStartPar
size of problem buffer.
\end{quote}
\end{quote}

\subsection{CoptException}
\label{\detokenize{javaapiref:coptexception}}
\sphinxAtStartPar
Copt exception object.

\sphinxstepscope

\subsubsection{CoptException.CoptException()}
\label{\detokenize{javaapi/CoptException:coptexception-coptexception}}\label{\detokenize{javaapi/CoptException::doc}}\begin{quote}

\sphinxAtStartPar
Constructor of coptexception.

\sphinxAtStartPar
\sphinxstylestrong{Synopsis}
\begin{quote}

\sphinxAtStartPar
\sphinxcode{\sphinxupquote{CoptException(int code, String msg)}}
\end{quote}

\sphinxAtStartPar
\sphinxstylestrong{Arguments}
\begin{quote}

\sphinxAtStartPar
\sphinxcode{\sphinxupquote{code}}: error code for exception.

\sphinxAtStartPar
\sphinxcode{\sphinxupquote{msg}}: error message for exception.
\end{quote}
\end{quote}

\subsubsection{CoptException.getCode()}
\label{\detokenize{javaapi/CoptException:coptexception-getcode}}\begin{quote}

\sphinxAtStartPar
Get the error code associated with the exception.

\sphinxAtStartPar
\sphinxstylestrong{Synopsis}
\begin{quote}

\sphinxAtStartPar
\sphinxcode{\sphinxupquote{int getCode()}}
\end{quote}

\sphinxAtStartPar
\sphinxstylestrong{Return}
\begin{quote}

\sphinxAtStartPar
the error code.
\end{quote}
\end{quote}

\renewcommand{\indexname}{Index}
\printindex
\end{document}